\documentclass[numbers=enddot,12pt,final,onecolumn,notitlepage]{scrartcl}%
\usepackage[headsepline,footsepline,manualmark]{scrlayer-scrpage}
\usepackage[all,cmtip]{xy}
\usepackage{amssymb}
\usepackage{amsmath}
\usepackage{amsthm}
\usepackage{framed}
\usepackage{comment}
\usepackage{color}
\usepackage[breaklinks=True]{hyperref}
\usepackage[sc]{mathpazo}
\usepackage[T1]{fontenc}
\usepackage{needspace}
\usepackage{tabls}
\usepackage{ytableau}
\usepackage{tikz}
\usepackage{pgfplots}
\usepackage[type={CC}, modifier={zero}, version={1.0},]{doclicense}
\providecommand{\U}[1]{\protect\rule{.1in}{.1in}}
\pgfplotsset{compat=1.13}
\theoremstyle{definition}
\newtheorem{theo}{Theorem}[subsection]
\newtheorem{exer}{Exercise}[subsection]

\newenvironment{theorem}[1][]
{\begin{theo}[#1]\begin{leftbar}}
{\end{leftbar}\end{theo}}
\newtheorem{lem}[theo]{Lemma}
\newenvironment{lemma}[1][]
{\begin{lem}[#1]\begin{leftbar}}
{\end{leftbar}\end{lem}}
\newtheorem{prop}[theo]{Proposition}
\newenvironment{proposition}[1][]
{\begin{prop}[#1]\begin{leftbar}}
{\end{leftbar}\end{prop}}
\newtheorem{defi}[theo]{Definition}
\newenvironment{definition}[1][]
{\begin{defi}[#1]\begin{leftbar}}
{\end{leftbar}\end{defi}}
\newtheorem{remk}[theo]{Remark}
\newenvironment{remark}[1][]
{\begin{remk}[#1]\begin{leftbar}}
{\end{leftbar}\end{remk}}
\newtheorem{coro}[theo]{Corollary}
\newenvironment{corollary}[1][]
{\begin{coro}[#1]\begin{leftbar}}
{\end{leftbar}\end{coro}}
\newtheorem{conv}[theo]{Convention}
\newenvironment{convention}[1][]
{\begin{conv}[#1]\begin{leftbar}}
{\end{leftbar}\end{conv}}
\newtheorem{quest}[theo]{Question}
\newenvironment{question}[1][]
{\begin{quest}[#1]\begin{leftbar}}
{\end{leftbar}\end{quest}}
\newtheorem{warn}[theo]{Warning}
\newenvironment{warning}[1][]
{\begin{warn}[#1]\begin{leftbar}}
{\end{leftbar}\end{warn}}
\newtheorem{conj}[theo]{Conjecture}
\newenvironment{conjecture}[1][]
{\begin{conj}[#1]\begin{leftbar}}
{\end{leftbar}\end{conj}}
\newtheorem{exam}[theo]{Example}
\newenvironment{example}[1][]
{\begin{exam}[#1]\begin{leftbar}}
{\end{leftbar}\end{exam}}
\newtheorem{exmp}[exer]{Exercise}
\newenvironment{exercise}[1][]
{\begin{exmp}[#1]\begin{leftbar}}
{\end{leftbar}\end{exmp}}
\newenvironment{statement}{\begin{quote}}{\end{quote}}
\newenvironment{fineprint}{\begin{small}}{\end{small}}
\let\sumnonlimits\sum
\let\prodnonlimits\prod
\let\cupnonlimits\bigcup
\let\capnonlimits\bigcap
\renewcommand{\sum}{\sumnonlimits\limits}
\renewcommand{\prod}{\prodnonlimits\limits}
\renewcommand{\bigcup}{\cupnonlimits\limits}
\renewcommand{\bigcap}{\capnonlimits\limits}
\usetikzlibrary{arrows,arrows.meta,decorations.markings}
\newenvironment{noncompile}{}{}
\excludecomment{verlong}
\includecomment{vershort}
\excludecomment{noncompile}

\newcommand{\arinj}{\ar@{_{(}->}}
\newcommand{\arinjrev}{\ar@{^{(}->}}
\newcommand{\arsurj}{\ar@{->>}}
\newcommand{\arelem}{\ar@{|->}}
\newcommand{\arback}{\ar@{<-}}

\newcommand{\MurpF}{\vphantom{\int_g^g}\mathcal{F}}

\definecolor{dbluecolor}{rgb}{0.01,0.02,0.7}
\definecolor{dgreencolor}{rgb}{0.2,0.4,0.0}
\definecolor{darkred}{rgb}{0.7,0,0}
\newtheoremstyle{plainsl}
{8pt plus 2pt minus 4pt}
{8pt plus 2pt minus 4pt}
{\slshape}
{0pt}
{\bfseries}
{.}
{5pt plus 1pt minus 1pt}
{}
\theoremstyle{plainsl}
\ihead{Symmetric group algebras, Spring 2024, version 2026-06-12}
\ohead{page \thepage}
\begin{document}

\title{An introduction to the symmetric group algebra\\{\normalsize [Math 701, Spring 2024 lecture notes]}}
\author{Darij Grinberg}
\date{November 19, 2026}
\maketitle

\begin{abstract}
\textbf{Abstract.} This is an introduction to the group algebras of the
symmetric groups, written for a quarter-long graduate course. After recalling
the definition of group algebras (and monoid algebras) in general, as well as
basic properties of permutations, we introduce several families of elements in
the symmetric group algebras $\mathbf{k}\left[  S_{n}\right]  $ such as the
Young--Jucys--Murphy elements, the (sign-)integrals and the conjugacy class
sums. Then comes a chapter on group actions and representations in general,
followed by the core of this text: a study of the representations of symmetric
groups (i.e., of left $\mathbf{k}\left[  S_{n}\right]  $-modules), including
the classical theory of Young tableaux and Young symmetrizers. We prove in
detail the main facts including the characterization of irreducible
representations (in characteristic $0$), the Garnir relations, the standard
basis theorem, the description of duals of Specht modules, and the hook length
formula, as well as a number of less known results. Finally, we describe
several bases of $\mathbf{k}\left[  S_{n}\right]  $ that arise from the study
of Specht modules, including the Murphy cellular bases.

The methods used are elementary and computational. We aim to assume as little
as possible of the base ring $\mathbf{k}$, and to use as little as possible
from representation theory (nothing more advanced than Maschke and Jordan--H\"{o}lder).

Over 100 exercises (without solutions) are scattered through the text.

\end{abstract}
\tableofcontents

\doclicenseThis

\section{What is this?}

These are the notes for a graduate-level course on the group algebras of the
symmetric groups. In the present state, they tell perhaps half of what I
wanted to tell, but this half is in a finished state. More will be added as I
find the time.

\subsection*{Acknowledgments}

I thank Daniel Harrer, Daniel Kirch, Victor Reiner, Franco Saliola and Jessica
Tomasko for interesting comments and corrections. Travis Scrimshaw has
patiently explained the intricacies of SageMath to me. Minor errors have also
been found by various versions of GPT. Questions, corrections and suggestions
are welcome!

\section{\label{chap.intro2}Before we start...}

\subsection{What is this?}

This is a course on the \emph{group algebras of the symmetric groups}: one of
the most famous and feature-rich objects in algebraic combinatorics.

These group algebras have been studied since the early 20th Century, under the
name \textquotedblleft substitutional analysis\textquotedblright\ (see
\cite{Ruther48} for a textbook treatment on those early results). Here,
\textquotedblleft substitution\textquotedblright\ is an antiquated word for
\textquotedblleft permutation\textquotedblright, whereas \textquotedblleft
analysis\textquotedblright\ is used in the literal sense (analyzing things,
not limits and derivatives).

Originally, this was used as an abstraction for what happens if you permute
the arguments of a multivariate function. Given a function $f:\mathbb{R}%
^{n}\rightarrow\mathbb{R}$, we can obtain $n!$ other functions $f_{\sigma
}:\mathbb{R}^{n}\rightarrow\mathbb{R}$, one for each permutation $\sigma$ of
$\left\{  1,2,\ldots,n\right\}  $, by permuting the arguments of $f$, as
follows:%
\[
f_{\sigma}\left(  x_{1},x_{2},\ldots,x_{n}\right)  =f\left(  x_{\sigma\left(
1\right)  },x_{\sigma\left(  2\right)  },\ldots,x_{\sigma\left(  n\right)
}\right)  .
\]
Then, we can take linear combinations of these $n!$ many functions, for
example the symmetrization $\dfrac{1}{n!}\sum\limits_{\sigma\in S_{n}%
}f_{\sigma}$. The input set $\mathbb{R}^{n}$ could be replaced by $I^{n}$ for
any set $I$; this includes \textquotedblleft tensors\textquotedblright\ in the
vulgar meaning of the word (\textquotedblleft higher-dimensional
matrices\textquotedblright), as well as the multilinear maps known from linear
algebra. In particular, the study of \textquotedblleft symmetries of
tensors\textquotedblright\ in the 1950s (by Hermann Weyl and others) rests on
the group algebras of the symmetric groups.

We will study these group algebras abstractly -- meaning that we won't mainly
be applying them to specific functions $f:\mathbb{R}^{n}\rightarrow\mathbb{R}$
but we will be computing formally with the permutations. This is similar to
studying polynomials without evaluating them at specific points, or studying
permutations as \textquotedblleft abstract\textquotedblright\ maps rather than
using them to permute a given tuple.

Permutations are indeed the backbone of this theory. In a way, symmetric group
algebras are to permutations as polynomial rings are to monomials. An element
of the group algebra of a symmetric group is a \textquotedblleft
formal\textquotedblright\ linear combination of the permutations in this
group, where the coefficients come from a given commutative ring $\mathbf{k}$.
These kinds of combinations can be taken for any group, not just a group of
permutations, but we will say very little about this general case.

This is \textbf{not} a course on representation theory. But we will cover,
among other things, the main results of the representation theory of symmetric
groups over fields of characteristic $0$.

\subsection{Prerequisites}

We will use

\begin{itemize}
\item some abstract algebra: some basic group theory (Lagrange's theorem,
cosets; see \cite{McClea23}, \cite[Chapters 1--4]{DumFoo04} or \cite{Conrad})
and some theory of rings and modules (as covered in \cite{23wa});

\item some combinatorics of permutations (as covered in \cite[Chapter 5]{21s}
and to some extent in \cite{22fco});

\item some theory of determinants (\cite[Chapter 6]{21s}) and characteristic
polynomials (the Cayley--Hamilton theorem; see, e.g., \cite[Chapter
3]{StoLui19}).
\end{itemize}

No representation theory is required to read these notes; what little of it
will be used will be introduced in Chapter \ref{chp.rep}.

We will be covering both classical and fairly recent research, overlapping
with the purview of several textbooks. I can recommend \cite{Wildon18},
\cite{Sagan01}, \cite{Meliot17}, \cite{Fulton97}, \cite{JamKer81},
\cite{Ruther48}, \cite{Kerber99}, \cite{Lorenz18}, \cite{EGHetc11} in
particular (further sources will be cited in the respective chapters and
sections). However, we will not follow any source for too long, and we will
try to approach the material in somewhat novel ways.

\subsection{A note on the exercises}

All exercises in these notes are equipped with a difficulty rating. This is
the number in the square (such as \fbox{3}) placed at the beginning of the
exercise. The harder the exercise, the larger this number is (although the
rating has been artificially lowered for some exercises that are too
tangential to the subject at hand). For example, a \fbox{5} represents a good
graduate-level homework problem that requires thinking and work. A \fbox{3}
usually requires some thinking \textbf{or} work. A \fbox{1} is a warm-up
question that might be obvious after 15 seconds of thought. A \fbox{7} should
be somewhat too hard for homework. Exercises with ratings above \fbox{10} are
hard research results on which I would enjoy hearing your ideas but don't
really expect a solution. The rating \fbox{open} represents an open problem
(conjecture or question). Multi-part exercises sometimes have points split
between the parts.

When I teach from these notes, I assign the exercises as homework. By solving
exercises, students gain points towards their total score. The points gained
from a given exercise equal its difficulty rating (or less when the solution
is incomplete or flawed). In Spring 2024, the requirement for a perfect grade
(A+) was to gain at least 50 points in total. In hindsight, it would have
probably been better to break the text into homework sets and require (say) at
least 10 points from each homework set.

In solving an exercise, you can freely use (without proof) the claims of all
exercises above it.

\subsection{\label{sec.intro.nots}Notations and elementary facts}

Before we properly begin, let me list some notations and basic facts that we
will be using.

\subsubsection{Basics}

First, we introduce some notations and terminology for basic combinatorial concepts:

\begin{itemize}
\item The symbol $\mathbb{N}$ will denote the set $\left\{  0,1,2,3,\ldots
\right\}  $ of nonnegative integers.

\item The size (i.e., cardinality) of a set $A$ will be denoted by $\left\vert
A\right\vert $.

\item The symbol \textquotedblleft\#\textquotedblright\ means
\textquotedblleft number\textquotedblright. For example, the size $\left\vert
A\right\vert $ of a set $A$ is the \# of elements of $A$.

\item For any integer $k$, we let $\left[  k\right]  $ denote the set
$\left\{  1,2,\ldots,k\right\}  $. This is a $k$-element set if $k\geq0$, and
is the empty set if $k\leq0$.
\end{itemize}

\subsubsection{\label{subsec.intro.nots.perms}Permutations}

Permutations are among the main characters of this course, and unsurprisingly
we will use some of their properties. Here we collect a few basic ones; we
refer to \cite[Chapter 5]{21s}, \cite[Chapter 5]{detnotes}, \cite[Chapter
4]{22fco} for proofs and details.

\begin{itemize}
\item A \emph{permutation} of a set $X$ means a bijection (i.e., bijective
map) from $X$ to $X$. Permutations are composed like any functions: If
$\alpha$ and $\beta$ are two permutations of $X$, then
\[
\left(  \alpha\circ\beta\right)  \left(  x\right)  =\alpha\left(  \beta\left(
x\right)  \right)  \ \ \ \ \ \ \ \ \ \ \text{for any }x\in X.
\]
(Beware that some authors -- particularly British ones -- disagree.)

The composition $\alpha\circ\beta$ is again a permutation of $X$.

\item The \emph{inverse} of a permutation $\alpha\in S_{X}$ is the inverse
$\alpha^{-1}$ of the map $\alpha$. It is again a permutation of $X$.

\item The identity map $\operatorname*{id}\nolimits_{X}:X\rightarrow X$ (often
just called $\operatorname*{id}$) is the map sending every $x$ to $x$. It is
also a permutation of $X$.

\item For a set $X$, the \emph{symmetric group} $S_{X}$ is the group of all
permutations of $X$, where the group operation is composition (so that
$\alpha\beta=\alpha\circ\beta$ for any $\alpha,\beta\in S_{X}$). Its size is
$\left\vert S_{X}\right\vert =\left\vert X\right\vert !$ (when $X$ is finite).

(Alternative notations for $S_{X}$ are $\mathfrak{S}_{X}$ and $\Sigma_{X}$ and
$\mathcal{S}_{X}$ and $\operatorname*{Aut}X$. The last one is justified by the
fact that the permutations are the automorphisms in the category of sets. The
other three are justified by the lack of letters that authors often experience
in combinatorics.)

Since $S_{X}$ is a group, we shall use the notation $\alpha\beta$ as a synonym
for $\alpha\circ\beta$, and we shall use the notation $1$ as a synonym for
$\operatorname*{id}\nolimits_{X}$ (the identity map of $X$, aka the neutral
element of $S_{X}$).

\item For a given $n\in\mathbb{N}$, the symmetric group $S_{\left[  n\right]
}=S_{\left\{  1,2,\ldots,n\right\}  }$ is denoted by $S_{n}$. This is called
the $n$\emph{-th symmetric group}, and has size $\left\vert S_{n}\right\vert
=n!$. (Again, some authors call it $\mathfrak{S}_{n}$ or $\Sigma_{n}$ or
$\mathcal{S}_{n}$.)

\item The easiest way to write down a permutation in $S_{n}$ is using its
\emph{one-line notation}.

The \emph{one-line notation} (short: \emph{OLN}) of a permutation $\sigma\in
S_{n}$ is the $n$-tuple $\left(  \sigma\left(  1\right)  ,\sigma\left(
2\right)  ,\ldots,\sigma\left(  n\right)  \right)  $. We denote it by
$\operatorname*{OLN}\sigma$. Thus,%
\[
\operatorname*{OLN}\sigma=\left(  \sigma\left(  1\right)  ,\sigma\left(
2\right)  ,\ldots,\sigma\left(  n\right)  \right)  .
\]
Sometimes we omit the parentheses and the commas in the OLN. For example, if
$\sigma\in S_{3}$ is the permutation that sends $1,2,3$ to $3,1,2$, then
$\operatorname*{OLN}\sigma=\left(  3,1,2\right)  $, and we abbreviate it as
$312$. (This kind of abbreviation becomes ambiguous when $\sigma\in S_{n}$ for
$n>10$, but fortunately such large numbers are rarely used in examples.)

Conversely, if $i_{1},i_{2},\ldots,i_{n}$ are $n$ distinct elements of
$\left[  n\right]  $, then there exists a unique permutation $\sigma\in S_{n}$
with $\operatorname*{OLN}\sigma=\left(  i_{1},i_{2},\ldots,i_{n}\right)  $.
This permutation $\sigma$ is the permutation of $\left[  n\right]  $ that
sends $1,2,\ldots,n$ to $i_{1},i_{2},\ldots,i_{n}$, respectively. It will be
denoted by $\operatorname*{oln}\left(  i_{1},i_{2},\ldots,i_{n}\right)  $.
Sometimes, for the sake of brevity, we omit the commas in this notation, so
that we just write $\operatorname*{oln}\left(  i_{1}i_{2}\ldots i_{n}\right)
$.

\item Some of the most widely used permutations are the \emph{cycles}. If
$i_{1},i_{2},\ldots,i_{k}$ are $k$ distinct elements of a set $X$, then the
\emph{cycle} $\operatorname*{cyc}\nolimits_{i_{1},i_{2},\ldots,i_{k}}$ denotes
the permutation of $X$ that sends $i_{1}$ to $i_{2}$, sends $i_{2}$ to $i_{3}%
$, sends $i_{3}$ to $i_{4}$, and so on, sends $i_{k}$ to $i_{1}$, and leaves
all the remaining elements of $X$ unchanged. In other words, it is defined by%
\[
\operatorname*{cyc}\nolimits_{i_{1},i_{2},\ldots,i_{k}}\left(  x\right)  =%
\begin{cases}
i_{h+1}, & \text{if }x=i_{h}\text{ for some }h\in\left[  k\right]  ;\\
x, & \text{if not}%
\end{cases}
\ \ \ \ \ \ \ \ \ \ \text{for each }x\in X,
\]
where we set $i_{k+1}:=i_{1}$.

(Many authors write $\left(  i_{1},i_{2},\ldots,i_{k}\right)  $ for
$\operatorname*{cyc}\nolimits_{i_{1},i_{2},\ldots,i_{k}}$. Be careful if you
do this, because the notation $\left(  i_{1},i_{2},\ldots,i_{k}\right)  $
already means a $k$-tuple, and thus you are creating an ambiguity. And it can
bite: $\left(  1,2,3\right)  $ and $\left(  2,3,1\right)  $ are not the same
$3$-tuple, but $\operatorname*{cyc}\nolimits_{1,2,3}$ and $\operatorname*{cyc}%
\nolimits_{2,3,1}$ are the same cycle.)

Note that $\operatorname*{cyc}\nolimits_{i}=\operatorname*{id}$ for any $i\in
X$. See \cite[\S 5.2.2 and \S 5.5]{21s} and \cite[\S 4.2 and \S 4.3]{22fco}
for more about cycles.

Strictly speaking, the cycle $\operatorname*{cyc}\nolimits_{i_{1},i_{2}%
,\ldots,i_{k}}$ depends not only on the elements $i_{1},i_{2},\ldots,i_{k}$,
but also on the set $X$, even though the notation $\operatorname*{cyc}%
\nolimits_{i_{1},i_{2},\ldots,i_{k}}$ neglects to mention this. In practice,
this is rarely problematic, since the set $X$ is usually clear from the
context or irrelevant. (It helps that $\operatorname*{cyc}\nolimits_{i_{1}%
,i_{2},\ldots,i_{k}}\left(  x\right)  =x$ for all $x\notin\left\{  i_{1}%
,i_{2},\ldots,i_{k}\right\}  $.)

\item If $i$ and $j$ are two distinct elements of $X$, then we set
\begin{equation}
t_{i,j}:=\operatorname*{cyc}\nolimits_{i,j}\in S_{X}.
\label{eq.intro.perms.cycs.tij=}%
\end{equation}
This is called the \emph{transposition} that swaps $i$ and $j$. It is the
permutation of $X$ that swaps $i$ with $j$ and leaves the other elements unchanged.

Note that $t_{i,j}=t_{j,i}$ for any distinct $i$ and $j$. More generally,
\[
\operatorname*{cyc}\nolimits_{i_{1},i_{2},\ldots,i_{k}}=\operatorname*{cyc}%
\nolimits_{i_{2},i_{3},\ldots,i_{k},i_{1}}=\operatorname*{cyc}\nolimits_{i_{3}%
,i_{4},\ldots,i_{k},i_{1},i_{2}}=\cdots=\operatorname*{cyc}\nolimits_{i_{k}%
,i_{1},i_{2},\ldots,i_{k-1}}%
\]
whenever $i_{1},i_{2},\ldots,i_{k}$ are $k$ distinct elements of $X$.

\item Specifically for $X=\left[  n\right]  $, there is a specific sort of
transpositions known as the \emph{simple transpositions}:

For any $n\in\mathbb{N}$ and any $i\in\left[  n-1\right]  $, we define%
\begin{equation}
s_{i}:=t_{i,i+1}\in S_{n}. \label{eq.intro.perms.cycs.si=}%
\end{equation}
This is the transposition that swaps $i$ with $i+1$. It is called a
\emph{simple transposition}.

An important fact (see \cite[Theorem 5.3.17 \textbf{(a)}]{21s} or most good
courses on group theory) is that every permutation in $S_{n}$ can be written
as a product of simple transpositions.

\item If $\phi:X\rightarrow Y$ is a bijection between two sets $X$ and $Y$,
then there is a group isomorphism%
\begin{align}
S_{\phi}:S_{X}  &  \rightarrow S_{Y},\label{eq.intro.perms.cycs.Sphi}\\
\alpha &  \mapsto\phi\circ\alpha\circ\phi^{-1}.\nonumber
\end{align}
Essentially, it just transforms every permutation of $X$ into a permutation of
$Y$ by \textquotedblleft relabelling\textquotedblright\ each $x\in X$ as
$\phi\left(  x\right)  \in Y$.

In particular, if $X$ is a finite set of size $n\in\mathbb{N}$, then we can
pick any bijection $\phi:X\rightarrow\left[  n\right]  $ and thus obtain a
group isomorphism $S_{\phi}:S_{X}\rightarrow S_{n}$. This is why it is often
sufficient to study $S_{n}$ when you want to understand $S_{X}$ for arbitrary
finite sets $X$.

\item For any set $X$, any cycle $\operatorname*{cyc}\nolimits_{i_{1}%
,i_{2},\ldots,i_{k}}\in S_{X}$ and any permutation $\sigma\in S_{X}$, we have%
\begin{equation}
\sigma\circ\operatorname*{cyc}\nolimits_{i_{1},i_{2},\ldots,i_{k}}\circ
\sigma^{-1}=\operatorname*{cyc}\nolimits_{\sigma\left(  i_{1}\right)
,\sigma\left(  i_{2}\right)  ,\ldots,\sigma\left(  i_{k}\right)  }
\label{eq.intro.perms.cycs.scs-1}%
\end{equation}
(see, e.g., \cite[Exercise 5.17 \textbf{(a)}]{detnotes} for the easy
proof\footnote{Note that \cite[Exercise 5.17 \textbf{(a)}]{detnotes} is only
stated for the case $X=\left[  n\right]  $, but the proof works equally for
any $X$.}). In particular,%
\begin{equation}
\sigma\circ t_{i,j}\circ\sigma^{-1}=t_{\sigma\left(  i\right)  ,\sigma\left(
j\right)  } \label{eq.intro.perms.cycs.sts-1}%
\end{equation}
for any transposition $t_{i,j}\in S_{X}$ and any permutation $\sigma\in S_{X}$.

\item Each cycle is a product (i.e., composition) of transpositions: Indeed,
for any set $X$ and any cycle $\operatorname*{cyc}\nolimits_{i_{1}%
,i_{2},\ldots,i_{k}}\in S_{X}$, we have%
\begin{equation}
\operatorname*{cyc}\nolimits_{i_{1},i_{2},\ldots,i_{k}}=t_{i_{1},i_{2}}\circ
t_{i_{2},i_{3}}\circ\cdots\circ t_{i_{k-1},i_{k}}.
\label{eq.intro.perms.cycs.c=ttt}%
\end{equation}
(The right hand side is understood to be $\operatorname*{id}\nolimits_{X}$ if
$k=1$. Generally, an empty product in a group is understood to be the neutral
element of the group.)

For a proof of (\ref{eq.intro.perms.cycs.c=ttt}), see \cite[Exercise
5.16]{detnotes}\footnote{Again, \cite[Exercise 5.16]{detnotes} is only stated
for $X=\left[  n\right]  $, but the proof applies equally well to all $X$.}.

\item For any set $X$ and any $m$ distinct elements $i_{1},i_{2},\ldots,i_{m}$
of $X$, we have%
\begin{equation}
\operatorname*{cyc}\nolimits_{i_{1},i_{2},\ldots,i_{m}}=\operatorname*{cyc}%
\nolimits_{i_{1},i_{2},\ldots,i_{k}}\circ\operatorname*{cyc}\nolimits_{i_{k}%
,i_{k+1},\ldots,i_{m}} \label{eq.intro.perms.cycs.c=cc}%
\end{equation}
for all $1\leq k\leq m$. This is not hard to check, either by expressing all
three cycles as products of transpositions using
(\ref{eq.intro.perms.cycs.c=ttt}), or by directly comparing how the two sides
act on each $x\in X$.

\item Disjoint cycles commute: If $i_{1},i_{2},\ldots,i_{k}$ and $j_{1}%
,j_{2},\ldots,j_{\ell}$ are $k+\ell$ distinct elements of $X$, then%
\begin{equation}
\operatorname*{cyc}\nolimits_{i_{1},i_{2},\ldots,i_{k}}\circ
\operatorname*{cyc}\nolimits_{j_{1},j_{2},\ldots,j_{\ell}}=\operatorname*{cyc}%
\nolimits_{j_{1},j_{2},\ldots,j_{\ell}}\circ\operatorname*{cyc}%
\nolimits_{i_{1},i_{2},\ldots,i_{k}}. \label{eq.intro.perms.cycs.cc=cc}%
\end{equation}
This is easy to check directly by comparing the action on each $x\in X$.

Note that the converse is not true: Non-disjoint cycles can sometimes also
commute. In other words, (\ref{eq.intro.perms.cycs.cc=cc}) occasionally holds
even if the $k+\ell$ elements $i_{1},i_{2},\ldots,i_{k}$ and $j_{1}%
,j_{2},\ldots,j_{\ell}$ are not all distinct.

\item To invert a cycle, we just read its indices right to left: In other
words, for any cycle $\operatorname*{cyc}\nolimits_{i_{1},i_{2},\ldots,i_{k}}%
$, we have%
\begin{equation}
\left(  \operatorname*{cyc}\nolimits_{i_{1},i_{2},\ldots,i_{k}}\right)
^{-1}=\operatorname*{cyc}\nolimits_{i_{k},i_{k-1},\ldots,i_{1}}.
\label{eq.intro.perms.cycs.inverse}%
\end{equation}

\item Let $X$ be a finite set. The \emph{sign} $\left(  -1\right)  ^{\sigma}$
of a permutation $\sigma\in S_{X}$ is a number that is either $1$ or $-1$, and
can be defined in the following two equivalent ways:

\begin{enumerate}
\item If $\sigma$ is a product of an even number of transpositions, then
$\left(  -1\right)  ^{\sigma}=1$. If $\sigma$ is a product of an odd number of
transpositions, then $\left(  -1\right)  ^{\sigma}=-1$.

\item Alternatively, we can define signs as follows: The sign of a permutation
$\sigma\in S_{n}$ (for some $n\in\mathbb{N}$) is defined to be $\left(
-1\right)  ^{\ell\left(  \sigma\right)  }$, where $\ell\left(  \sigma\right)
$ is the \# of inversions of $\sigma$ (that is, of all pairs $\left(
i,j\right)  \in\left[  n\right]  ^{2}$ such that $i<j$ and $\sigma\left(
i\right)  >\sigma\left(  j\right)  $). The sign of a permutation $\sigma\in
S_{X}$ for an arbitrary finite set $X$ is defined by mapping $S_{X}$ to
$S_{n}$ for $n=\left\vert X\right\vert $ using the above group isomorphism
$S_{\phi}$ defined in (\ref{eq.intro.perms.cycs.Sphi}), and computing the sign
there. The exact choice of $\phi$ does not matter for the sign.
\end{enumerate}

For details and proofs, see \cite[\S 5.3 and \S 5.6]{detnotes} (and
\cite[\S 5.4]{21s} for a sketch). In particular, the second definition is the
one given in \cite[Definition 5.14 and Exercise 5.12]{detnotes}, whereas the
first definition is equivalent to it by \cite[Exercise 5.15 \textbf{(c)}%
]{detnotes}.

For us, the most important properties of signs are%
\begin{align*}
\left(  -1\right)  ^{\sigma\tau}  &  =\left(  -1\right)  ^{\sigma}\cdot\left(
-1\right)  ^{\tau}\ \ \ \ \ \ \ \ \ \ \text{for any }\sigma,\tau\in S_{X};\\
\left(  -1\right)  ^{t_{i,j}}  &  =-1\ \ \ \ \ \ \ \ \ \ \text{for any
distinct }i,j\in X;\\
\left(  -1\right)  ^{\operatorname*{cyc}\nolimits_{i_{1},i_{2},\ldots,i_{k}}}
&  =\left(  -1\right)  ^{k-1}\ \ \ \ \ \ \ \ \ \ \text{for any cycle
}\operatorname*{cyc}\nolimits_{i_{1},i_{2},\ldots,i_{k}}\in S_{X};\\
\left(  -1\right)  ^{\sigma^{-1}}  &  =\left(  -1\right)  ^{\sigma
}\ \ \ \ \ \ \ \ \ \ \text{for any }\sigma\in S_{X}.
\end{align*}
For the proofs, see \cite[Exercise 5.12 \textbf{(c)} and Exercise 5.15
\textbf{(a)}]{detnotes} as well as (\ref{eq.intro.perms.cycs.c=ttt}). The
identity $\left(  -1\right)  ^{\sigma\tau}=\left(  -1\right)  ^{\sigma}%
\cdot\left(  -1\right)  ^{\tau}$ is known as the \emph{multiplicativity of the
sign}, and shows that the \emph{sign homomorphism} (i.e., the map
$S_{X}\rightarrow\left\{  1,-1\right\}  $ that sends each $\sigma\in S_{X}$ to
$\left(  -1\right)  ^{\sigma}$) is a group morphism. An easy consequence of
these properties is that conjugate permutations have equal signs -- i.e., that%
\begin{equation}
\left(  -1\right)  ^{\sigma\tau\sigma^{-1}}=\left(  -1\right)  ^{\tau
}\ \ \ \ \ \ \ \ \ \ \text{for all }\sigma,\tau\in S_{X}.
\label{eq.intro.perms.signs.conj}%
\end{equation}

\item A permutation $\sigma\in S_{X}$ of a finite set $X$ is said to be
\emph{even} if its sign $\left(  -1\right)  ^{\sigma}$ equals $1$, and it is
said to be \emph{odd} if its sign $\left(  -1\right)  ^{\sigma}$ equals $-1$.
\end{itemize}

\subsubsection{Algebraic preliminaries}

Next, we recall some concepts from abstract algebra. We will freely use the
notion of a ring, and the notions of modules and algebras over a commutative
ring. See \cite[Chapters 2 and 3]{23wa} for a detailed discussion of these
notions, and \cite[\S 3.4.2]{21s} for a quick reminder. The word
\textquotedblleft ring\textquotedblright, for us, always means an associative
unital ring (but not necessarily commutative).

We fix a commutative ring $\mathbf{k}$ from now on. It will be called our
\emph{base ring}, and its elements will be called \emph{scalars}. We recall
the basic ideas behind modules and algebras:

\begin{itemize}
\item A $\mathbf{k}$\emph{-module} means an additive abelian group whose
elements can be scaled by elements of $\mathbf{k}$ and satisfy a few simple
axioms. These axioms are precisely the vector space axioms known from linear
algebra, and indeed the $\mathbf{k}$-modules are also known as
\textquotedblleft$\mathbf{k}$-vector spaces\textquotedblright\ when
$\mathbf{k}$ is a field. But we will allow $\mathbf{k}$ to be an arbitrary
commutative ring, in which case it is frowned upon to speak of
\textquotedblleft$\mathbf{k}$-vector spaces\textquotedblright.

The elements of a $\mathbf{k}$-module will nevertheless be called
\emph{vectors}.

Just like in linear algebra, the simplest examples of $\mathbf{k}$-modules are
$\mathbf{k}$ itself, its Cartesian powers $\mathbf{k}^{n}$ (consisting of
$n$-tuples of elements of $\mathbf{k}$, added and scaled entrywise), its
\textquotedblleft matrix spaces\textquotedblright\ $\mathbf{k}^{n\times m}$
(consisting of $n\times m$-matrices with entries in $\mathbf{k}$), as well as
various submodules thereof (such as the submodule of $\mathbf{k}^{n}$ that
consists of all $n$-tuples $\left(  a_{1},a_{2},\ldots,a_{n}\right)  $
satisfying $a_{1}+a_{2}+\cdots+a_{n}=0$).

\item A $\mathbf{k}$\emph{-algebra} means a $\mathbf{k}$-module $A$ that is
also a ring (with the same addition and zero), with an extra compatibility
requirement saying that%
\[
\lambda\left(  ab\right)  =\left(  \lambda a\right)  b=a\left(  \lambda
b\right)  \ \ \ \ \ \ \ \ \ \ \text{for any }\lambda\in\mathbf{k}\text{ and
}a,b\in A.
\]
For example, $\mathbf{k}$ itself is a $\mathbf{k}$-algebra, and so are the
polynomial ring $\mathbf{k}\left[  x\right]  $ and the matrix ring
$\mathbf{k}^{n\times n}$ for any $n\in\mathbb{N}$.

\item The word \textquotedblleft\emph{morphism}\textquotedblright\ is short
for \textquotedblleft homomorphism\textquotedblright.

\item The word \textquotedblleft\emph{ideal}\textquotedblright\ means
\textquotedblleft two-sided ideal\textquotedblright, unless it is part of the
phrasing \textquotedblleft left ideal\textquotedblright\ or \textquotedblleft
right ideal\textquotedblright.
\end{itemize}

\bigskip

\section{Introduction to the symmetric group algebras}

\subsection{\label{sec.intro.monalg}Monoid algebras and group algebras}

\subsubsection{Monoids}

Next, we recall the notion of a monoid.

Informally, a \emph{monoid} is a \textquotedblleft group without
inverses\textquotedblright.

Rigorously, a \emph{monoid} is defined to be a triple $\left(  M,\cdot
,1\right)  $, where $M$ is a set, $\cdot$ is an associative binary operation
on $M$ (that is, a map from $M\times M$ to $M$), and $1$ is an element of $M$
that is neutral for this operation $\cdot$. When $M$ is a monoid, we will
write $mn$ for $m\cdot n$ whenever $m,n\in M$. The element $mn$ will be called
the \emph{product} of $m$ and $n$ in the monoid $M$. We will write $M$ for
$\left(  M,\cdot,1\right)  $ if the operation $\cdot$ and the element $1$ are
clear from the context. The monoid $M$ is said to be \emph{abelian} if all
$m,n\in M$ satisfy $mn=nm$.

The notations $\cdot$ and $1$ are not set in stone; some authors use $\ast$
and $e$ for them. Depending on the monoid, they might have different names.
For example, the most fundamental example of a monoid is the abelian monoid
$\left(  \mathbb{N},+,0\right)  $. Here, the operation is actually the
addition $+$ of the nonnegative integers, so it would be a bad idea to rename
it as $\cdot$. Likewise, its neutral element is the number $0$, not the number
$1$.

When we use the notations $\cdot$ and $1$ for the operation and the neutral
element of a monoid, we say that our monoid is \emph{written multiplicatively}.

Here are some examples of monoids:

\begin{itemize}
\item Any group is a monoid. For instance, $\left(  \mathbb{Z},+,0\right)  $
is a monoid, as is $\left(  \mathbb{Z},\cdot,1\right)  $. Actually, $\left(
\mathbb{Z},+,0\right)  $ is a group, but $\left(  \mathbb{Z},\cdot,1\right)  $
is not.

\item More generally, if $R$ is any ring, then $\left(  R,+,0\right)  $ and
$\left(  R,\cdot,1\right)  $ are monoids (where $+$ and $\cdot$ are the
operations of $R$, and where $0$ and $1$ are the zero and the unity of $R$).
The monoid $\left(  R,+,0\right)  $ is always an abelian group, but the monoid
$\left(  R,\cdot,1\right)  $ is (in general) neither abelian nor a group.

\item As we already mentioned, $\left(  \mathbb{N},+,0\right)  $ is a monoid
as well (but not a group). So is $\left(  \mathbb{N},\cdot,1\right)  $.

\item The Boolean monoid is defined to be the monoid $\left(  \mathbb{B}%
,\vee,\operatorname*{False}\right)  $, where $\mathbb{B}=\left\{
\operatorname*{True},\operatorname*{False}\right\}  $.

\item For any positive integer $n$, the triple $\left(  \left[  n\right]
,\ \max,\ 1\right)  $ is a monoid, where $\max$ is the binary operation
sending each pair $\left(  u,v\right)  \in\left[  n\right]  \times\left[
n\right]  $ to $\max\left\{  u,v\right\}  $. The associativity follows from
the fact that all $a,b,c\in\left[  n\right]  $ satisfy \newline$\max\left\{
a,\max\left\{  b,c\right\}  \right\}  =\max\left\{  \max\left\{  a,b\right\}
,c\right\}  $, whereas the neutrality of $1$ follows from the fact that each
$u\in\left[  n\right]  $ satisfies $\max\left\{  1,u\right\}  =\max\left\{
u,1\right\}  =u$.

\item We can extend the set $\mathbb{Z}$ to a slightly larger set
$\mathbb{Z}^{\ast}:=\mathbb{Z}\cup\left\{  -\infty\right\}  $ by inserting a
new element $-\infty$. We understand this new element $-\infty$ to be smaller
than every integer, so that $\max\left\{  k,-\infty\right\}  =\max\left\{
-\infty,k\right\}  =k$ for all $k\in\mathbb{Z}^{\ast}$. Now, the triple
$\left(  \mathbb{Z}^{\ast},\max,-\infty\right)  $ is a monoid (where $\max$ is
the binary operation sending each pair $\left(  u,v\right)  \in\mathbb{Z}%
^{\ast}\times\mathbb{Z}^{\ast}$ to $\max\left\{  u,v\right\}  $). It is called
the \emph{tropical monoid} on $\mathbb{Z}$. Similarly we can define a tropical
monoid on every totally ordered set.
\end{itemize}

The neutral element of any monoid $M$ will be called $1_{M}$.

A \emph{monoid morphism} is a map $f:M\rightarrow N$ between two monoids $M$
and $N$ that respects the neutral elements (i.e., satisfies $f\left(
1_{M}\right)  =1_{N}$) and respects multiplication (i.e., satisfies $f\left(
mm^{\prime}\right)  =f\left(  m\right)  \cdot f\left(  m^{\prime}\right)  $
for any $m,m^{\prime}\in M$).

\subsubsection{Monoid algebras}

Next, we shall define the \textquotedblleft monoid algebra\textquotedblright%
\ of a monoid. For details, we refer to \cite[\S 4.1]{23wa}.

Informally, the \emph{monoid algebra} of a monoid $M$ over a commutative ring
$\mathbf{k}$ is the $\mathbf{k}$-algebra obtained by \textquotedblleft
adjoining\textquotedblright\ the monoid $M$ to the ring $\mathbf{k}$, which
means \textquotedblleft inserting\textquotedblright\ the elements of $M$ into
$\mathbf{k}$. That is, this algebra consists of \textquotedblleft formal
products\textquotedblright\ $\lambda m$ with $\lambda\in\mathbf{k}$ and $m\in
M$, as well as their \textquotedblleft formal sums\textquotedblright%
\ $\lambda_{1}m_{1}+\lambda_{2}m_{2}+\cdots+\lambda_{i}m_{i}$, which are
called \textquotedblleft formal $\mathbf{k}$-linear combinations of elements
of $M$\textquotedblright. These $\mathbf{k}$-linear combinations are
multiplied using the multiplications of $\mathbf{k}$ and $M$: For example, if
$\lambda_{1}m_{1}+\lambda_{2}m_{2}$ and $\lambda_{3}m_{3}+\lambda_{4}m_{4}$
are two such $\mathbf{k}$-linear combinations, then their product is
\begin{align*}
&  \left(  \lambda_{1}m_{1}+\lambda_{2}m_{2}\right)  \left(  \lambda_{3}%
m_{3}+\lambda_{4}m_{4}\right) \\
&  =\underbrace{\lambda_{1}\lambda_{3}}_{\substack{\text{a product}\\\text{in
}\mathbf{k}}}\ \ \underbrace{m_{1}m_{3}}_{\substack{\text{a product}\\\text{in
}M}}+\,\lambda_{1}\lambda_{4}m_{1}m_{4}+\lambda_{2}\lambda_{3}m_{2}%
m_{3}+\lambda_{2}\lambda_{4}m_{2}m_{4}.
\end{align*}

Let us now define this monoid algebra formally.

First, we define the notion of the \emph{free }$\mathbf{k}$\emph{-module
}$\mathbf{k}^{\left(  S\right)  }$ for a given set $S$:

\begin{definition}
\label{def.monalg.freemod}Let $S$ be an arbitrary set, and (as before) let
$\mathbf{k}$ be a commutative ring.

As usual, $\mathbf{k}^{S}$ denotes the set of all families $\left(
\lambda_{s}\right)  _{s\in S}$ of elements of $\mathbf{k}$ indexed by elements
of $S$. For example, if $S=\left[  n\right]  $ for some $n\in\mathbb{N}$, then
such families are just the $n$-tuples $\left(  \lambda_{1},\lambda_{2}%
,\ldots,\lambda_{n}\right)  $ of elements of $\mathbf{k}$, so that
$\mathbf{k}^{\left[  n\right]  }=\mathbf{k}^{n}$. For another example, if
$S=\mathbb{N}$, then such families are the sequences $\left(  \lambda
_{0},\lambda_{1},\lambda_{2},\ldots\right)  $.

A family $\left(  \lambda_{s}\right)  _{s\in S}\in\mathbf{k}^{S}$ is said to
be \emph{essentially finite} if all but finitely many $s\in S$ satisfy
$\lambda_{s}=0$ (that is, if it has only finitely many nonzero entries).

We let $\mathbf{k}^{\left(  S\right)  }$ denote the set of all essentially
finite families in $\mathbf{k}^{S}$. Thus,%
\[
\mathbf{k}^{\left(  S\right)  }=\left\{  \left(  \lambda_{s}\right)  _{s\in
S}\in\mathbf{k}^{S}\ \mid\ \text{all but finitely many }s\in S\text{ satisfy
}\lambda_{s}=0\right\}  .
\]

The set $\mathbf{k}^{S}$ is a $\mathbf{k}$-module, where addition and scaling
are defined entrywise:%
\begin{align*}
\left(  a_{s}\right)  _{s\in S}+\left(  b_{s}\right)  _{s\in S}  &  :=\left(
a_{s}+b_{s}\right)  _{s\in S}\ \ \ \ \ \ \ \ \ \ \text{for any }\left(
a_{s}\right)  _{s\in S},\left(  b_{s}\right)  _{s\in S}\in\mathbf{k}^{S};\\
\lambda\left(  a_{s}\right)  _{s\in S}  &  :=\left(  \lambda a_{s}\right)
_{s\in S}\ \ \ \ \ \ \ \ \ \ \text{for any }\lambda\in\mathbf{k}\text{ and
}\left(  a_{s}\right)  _{s\in S}\in\mathbf{k}^{S}.
\end{align*}
The subset $\mathbf{k}^{\left(  S\right)  }$ of $\mathbf{k}^{S}$ is a
$\mathbf{k}$-submodule of $\mathbf{k}^{S}$, because the union of two finite
sets is still finite. This submodule $\mathbf{k}^{\left(  S\right)  }$ is
called the \emph{free }$\mathbf{k}$\emph{-module on the set }$S$. (We often
omit the $\mathbf{k}$, and just say \textquotedblleft free
module\textquotedblright.)

This $\mathbf{k}$-module $\mathbf{k}^{\left(  S\right)  }$ is indeed free.
Indeed, we can construct a basis (called the \emph{standard basis}) as
follows: For each $t\in S$, we define the \emph{standard basis vector }$e_{t}$
to be the family $\left(  \delta_{s,t}\right)  _{s\in S}\in\mathbf{k}^{S}$,
where $\delta_{s,t}$ is the Kronecker delta (that is, $1$ if $s=t$ and $0$
otherwise). This standard basis vector $e_{t}$ has only one nonzero entry
(namely, its $t$-th entry, which is $1$), and thus belongs to $\mathbf{k}%
^{\left(  S\right)  }$. The family $\left(  e_{t}\right)  _{t\in S}$ is a
basis of $\mathbf{k}^{\left(  S\right)  }$, called the \emph{standard basis}.
\end{definition}

For $S=\left[  n\right]  $, this all is standard linear algebra: Here,
$e_{t}=\left(  0,0,\ldots,0,1,0,0,\ldots,0\right)  $ with the $1$ in the
$t$-th position, and thus the standard basis is really the standard basis you
know from linear algebra.

Here are a few more remarks about Definition \ref{def.monalg.freemod}:

\begin{itemize}
\item We have $\mathbf{k}^{\left(  S\right)  }=\mathbf{k}^{S}$ when the set
$S$ is finite. When $S$ is infinite, $\mathbf{k}^{\left(  S\right)  }$ is a
free $\mathbf{k}$-module, while $\mathbf{k}^{S}$ is not always (for example,
$\mathbb{Z}^{\mathbb{N}}$ is not a free $\mathbb{Z}$-module).

\item We always have $\mathbf{k}^{S}=\prod\limits_{s\in S}\mathbf{k}$ and
$\mathbf{k}^{\left(  S\right)  }=\bigoplus\limits_{s\in S}\mathbf{k}$ as
$\mathbf{k}$-modules.

\item If $\mathbf{k}=\mathbb{R}$, then an element $\left(  \alpha_{s}\right)
_{s\in S}\in\mathbb{R}^{S}$ with all its entries $\alpha_{s}$ nonnegative and
with $\sum_{s\in S}\alpha_{s}=1$ is called a \emph{probability distribution}
on $S$. This distribution is said to be \emph{finitely supported} if it
belongs to $\mathbb{R}^{\left(  S\right)  }$. Thus, one may think of $\left(
\alpha_{s}\right)  _{s\in S}$ in the general case as generalized probability
distributions (a useful mental picture if you find probabilities intuitive).
\end{itemize}

When $S$ is just a set, $\mathbf{k}^{\left(  S\right)  }$ is just a free
module. But any further structures on $S$ yield corresponding extra structures
on $\mathbf{k}^{\left(  S\right)  }$. In particular, when $S$ is a monoid,
$\mathbf{k}^{\left(  S\right)  }$ inherits a multiplication and becomes a
$\mathbf{k}$-algebra:

\begin{definition}
\label{def.monalg.monalg}Let $M$ be a monoid, written multiplicatively (so
that $M=\left(  M,\cdot,1\right)  $).

The \emph{monoid algebra of }$M$\emph{ over }$\mathbf{k}$ (also known as the
\emph{monoid ring of }$M$\emph{ over }$\mathbf{k}$) is the $\mathbf{k}%
$-algebra $\mathbf{k}\left[  M\right]  $ defined as follows: As a $\mathbf{k}%
$-module, it is just the free $\mathbf{k}$-module $\mathbf{k}^{\left(
M\right)  }$. Its multiplication is defined to be the unique $\mathbf{k}%
$-bilinear map $\mu:\mathbf{k}^{\left(  M\right)  }\times\mathbf{k}^{\left(
M\right)  }\rightarrow\mathbf{k}^{\left(  M\right)  }$ that satisfies%
\begin{equation}
\mu\left(  e_{m},e_{n}\right)  =e_{mn}\ \ \ \ \ \ \ \ \ \ \text{for all
}m,n\in M. \label{eq.def.monalg.monalg.mu}%
\end{equation}
(Here, $\left(  e_{m}\right)  _{m\in M}$ is the standard basis of
$\mathbf{k}^{\left(  M\right)  }$, as defined above. Moreover, a map is said
to be $\mathbf{k}$\emph{-bilinear} if it is $\mathbf{k}$-linear in each of its
two arguments.)

The unity of this $\mathbf{k}$-algebra $\mathbf{k}\left[  M\right]  $ is the
basis element $e_{1}$ corresponding to the neutral element $1$ of $M$.
\end{definition}

\begin{theorem}
\label{thm.monalg.wd}This is indeed a well-defined $\mathbf{k}$-algebra.
\end{theorem}

\begin{proof}
See \cite[Theorem 4.1.2]{23wa}.
\end{proof}

Since the bilinear map $\mu$ in Definition \ref{def.monalg.monalg} is used as
the multiplication of $\mathbf{k}\left[  M\right]  $, we can rewrite
(\ref{eq.def.monalg.monalg.mu}) as follows:%
\[
e_{m}\cdot e_{n}=e_{mn}\ \ \ \ \ \ \ \ \ \ \text{for all }m,n\in M.
\]
Since the multiplication $\mu$ is $\mathbf{k}$-bilinear, we can multiply not
just two basis elements, but any two of their $\mathbf{k}$-linear
combinations. For example, for any $\alpha,\beta,\gamma,\delta\in\mathbf{k}$
and $m,n,p,q\in M$, we have%
\begin{align*}
&  \left(  \alpha e_{m}+\beta e_{n}\right)  \cdot\left(  \gamma e_{p}+\delta
e_{q}\right) \\
&  =\alpha e_{m}\left(  \gamma e_{p}+\delta e_{q}\right)  +\beta e_{n}\left(
\gamma e_{p}+\delta e_{q}\right)  \ \ \ \ \ \ \ \ \ \ \left(
\begin{array}
[c]{c}%
\text{since }\mu\text{ is linear in}\\
\text{the first argument}%
\end{array}
\right) \\
&  =\alpha\gamma\underbrace{e_{m}e_{p}}_{=e_{mp}}+\,\alpha\delta
\underbrace{e_{m}e_{q}}_{=e_{mq}}+\,\beta\gamma\underbrace{e_{n}e_{p}%
}_{=e_{np}}+\,\beta\delta\underbrace{e_{n}e_{q}}_{=e_{nq}}%
\ \ \ \ \ \ \ \ \ \ \left(
\begin{array}
[c]{c}%
\text{since }\mu\text{ is linear in}\\
\text{the second argument}%
\end{array}
\right) \\
&  =\alpha\gamma e_{mp}+\alpha\delta e_{mq}+\beta\gamma e_{np}+\beta\delta
e_{nq}.
\end{align*}

\subsubsection{Group algebras}

\begin{definition}
\label{def.monalg.gralg}When a monoid $M$ is a group, its monoid algebra
$\mathbf{k}\left[  M\right]  $ is called its \emph{group algebra} (or
\emph{group ring}).
\end{definition}

\subsubsection{Some properties}

Some conventions and some simple facts make it easier to work in monoid
algebras. First, we note that the commutativity of a monoid is inherited by
its monoid algebra:

\begin{proposition}
\label{prop.monalg.commut}Let $M$ be an \textbf{abelian} monoid. Then, the
algebra $\mathbf{k}\left[  M\right]  $ is commutative.
\end{proposition}

\begin{proof}
See \cite[Proposition 4.1.9]{23wa}.
\end{proof}

The next two propositions show that the monoid algebra $\mathbf{k}\left[
M\right]  $ contains a \textquotedblleft copy\textquotedblright\ of
$\mathbf{k}$ and (unless $\mathbf{k}$ is trivial) a \textquotedblleft
copy\textquotedblright\ of $M$.

\begin{proposition}
\label{prop.monalg.k-to-kM}Let $M$ be a monoid with neutral element $1$. Then,
the map%
\begin{align*}
\mathbf{k}  &  \rightarrow\mathbf{k}\left[  M\right]  ,\\
\lambda &  \mapsto\lambda e_{1}%
\end{align*}
is an injective $\mathbf{k}$-algebra morphism. Thus, its image is a
$\mathbf{k}$-subalgebra of $\mathbf{k}\left[  M\right]  $ that is isomorphic
to $\mathbf{k}$.
\end{proposition}

\begin{proof}
See \cite[Proposition 4.1.10]{23wa}.
\end{proof}

\begin{proposition}
\label{prop.monalg.M-to-kM}Let $M$ be a monoid. Then, the map%
\begin{align*}
M  &  \rightarrow\mathbf{k}\left[  M\right]  ,\\
m  &  \mapsto e_{m}%
\end{align*}
is a monoid morphism from $M$ to the monoid $\left(  \mathbf{k}\left[
M\right]  ,\cdot,1\right)  $. If the ring $\mathbf{k}$ is nontrivial (i.e., if
$0\neq1$ in $\mathbf{k}$), then this morphism is injective, so that its image
is a submonoid of $\left(  \mathbf{k}\left[  M\right]  ,\cdot,1\right)  $
isomorphic to $M$.
\end{proposition}

\begin{proof}
See \cite[Proposition 4.1.12]{23wa} for the first claim. The second claim
(about injectivity) is clear, since the $e_{m}$ for different $m\in M$ are
distinct elements of the standard basis $\left(  e_{m}\right)  _{m\in M}$ and
therefore distinct (when $\mathbf{k}$ is nontrivial\footnote{Indeed, if $m$
and $n$ are two distinct elements of $M$, then the $m$-th entry of the family
$e_{m}$ is $1$, whereas the $m$-th entry of the family $e_{n}$ is $0$. Thus,
these two entries are distinct, unless $1=0$ in $\mathbf{k}$.}).
\end{proof}

The above two propositions show that $\mathbf{k}\left[  M\right]  $ contains
both a copy of $\mathbf{k}$ and (when $\mathbf{k}$ is nontrivial) a copy of
$M$. Moreover, $\mathbf{k}\left[  M\right]  $ is \textquotedblleft
generated\textquotedblright\ by these two pieces: Every element of
$\mathbf{k}\left[  M\right]  $ has the form
\[
\lambda_{1}e_{m_{1}}+\lambda_{2}e_{m_{2}}+\cdots+\lambda_{k}e_{m_{k}}%
\]
for some scalars $\lambda_{1},\lambda_{2},\ldots,\lambda_{k}\in\mathbf{k}$ and
some monoid elements $m_{1},m_{2},\ldots,m_{k}\in M$. This is just because the
standard basis $\left(  e_{m}\right)  _{m\in M}$ is a basis of $\mathbf{k}%
\left[  M\right]  $.

\subsubsection{\label{subsec.intro.monalg.convs}Some conventions}

It gets tiresome to always write $e_{m}$ for a standard basis vector in
$\mathbf{k}\left[  M\right]  $. Wouldn't it be easier to just write $m$
instead? (After all, Proposition \ref{prop.monalg.M-to-kM} tells us that the
map sending each $m$ to $e_{m}$ is a monoid morphism, so that $e_{m}$ is a
\textquotedblleft copy\textquotedblright\ of $m$ that behaves just like $m$
with respect to multiplication.) This is indeed a good idea, and is
unproblematic as long as it does not get ambiguous, i.e., as long as
expressions like \textquotedblleft$m+n$\textquotedblright\ (for $m,n\in M$)
and \textquotedblleft$\lambda m$\textquotedblright\ (for $\lambda\in
\mathbf{k}$ and $m\in M$) don't already have a different meaning. So let us
agree to do this:

\begin{convention}
\label{conv.monalg.em=m}Let $M$ be a monoid. The standard basis vectors
$e_{m}$ in $\mathbf{k}\left[  M\right]  $ will just be denoted by $m$ (by
abuse of notation). Thus, for example, a linear combination of the form
$\alpha e_{m}+\beta e_{n}+\gamma e_{p}$ will be written $\alpha m+\beta
n+\gamma p$. With this notation, an element of $\mathbf{k}\left[  M\right]  $
really looks like a $\mathbf{k}$-linear combination of elements of $M$
(although in truth, it is not the elements of $M$ but the respective standard
basis vectors that are being combined).

Avoid this notation when the elements of $M$ can be scaled or added in another
way (i.e., when there already is a scaling or an addition defined on the
monoid $M$). For example, if $M$ is the monoid $\left(  \mathbb{N}%
,\cdot,1\right)  $, then this notation must not be used, since there is
already an addition defined on $\mathbb{N}$ which has nothing to do with the
addition in $\mathbf{k}\left[  M\right]  $ (for example, $2+3=5$ but
$e_{2}+e_{3}\neq e_{5}$). For another example in which this notation would be
confusing, consider the quaternion group: a group of size $8$ with elements
$1,i,j,k,-1,-i,-j,-k$. If we wrote $m$ for $e_{m}$ here, then the nonzero sum
$e_{1}+e_{-1}$ would look like the zero sum $1+\left(  -1\right)  $, which is
not a good thing.
\end{convention}

If $M$ is a monoid, then every element $\left(  \alpha_{m}\right)  _{m\in M}$
of the monoid algebra $\mathbf{k}\left[  M\right]  $ can be rewritten as
$\sum\limits_{m\in M}\alpha_{m}e_{m}$ (just like the vector $\left(
\alpha_{1},\alpha_{2},\ldots,\alpha_{n}\right)  \in\mathbb{R}^{n}$ in linear
algebra is rewritten as $\alpha_{1}e_{1}+\alpha_{2}e_{2}+\cdots+\alpha
_{n}e_{n}$). If we use Convention \ref{conv.monalg.em=m}, we can rewrite this
element further as $\sum\limits_{m\in M}\alpha_{m}m$ (since we write $m$ for
$e_{m}$). Thus, using Convention \ref{conv.monalg.em=m}, we have%
\begin{equation}
\left(  \alpha_{m}\right)  _{m\in M}=\sum\limits_{m\in M}\alpha_{m}e_{m}%
=\sum\limits_{m\in M}\alpha_{m}m \label{eq.rmk.monalg.compare.alm=}%
\end{equation}
for any family $\left(  \alpha_{m}\right)  _{m\in M}\in\mathbf{k}\left[
M\right]  $.

We will often use the notation $\sum\limits_{m\in M}\alpha_{m}m$ instead of
the \textquotedblleft native\textquotedblright\ notation $\left(  \alpha
_{m}\right)  _{m\in M}$, because it allows us to multiply two elements of
$\mathbf{k}\left[  M\right]  $ by just expanding the product:%
\[
\left(  \sum\limits_{m\in M}\alpha_{m}m\right)  \left(  \sum\limits_{n\in
M}\beta_{n}n\right)  =\sum\limits_{m\in M}\ \ \sum\limits_{n\in M}%
\ \ \underbrace{\alpha_{m}\beta_{n}}_{\text{a product in }\mathbf{k}%
}\ \ \underbrace{mn}_{\text{a product in }M}.
\]

\begin{remark}
\label{rmk.monalg.compare}Let $M$ be a monoid. Let $\left(  \alpha_{m}\right)
_{m\in M}$ and $\left(  \beta_{m}\right)  _{m\in M}$ be two elements of the
monoid algebra $\mathbf{k}\left[  M\right]  $. Then, the equation%
\begin{equation}
\sum\limits_{m\in M}\alpha_{m}m=\sum\limits_{m\in M}\beta_{m}%
m\ \ \ \ \ \ \ \ \ \ \text{in }\mathbf{k}\left[  M\right]
\label{eq.rmk.monalg.compare.ass}%
\end{equation}
is equivalent to the statement that%
\begin{equation}
\left(  \alpha_{m}=\beta_{m}\ \ \ \ \ \ \ \ \ \ \text{for each }m\in M\right)
. \label{eq.rmk.monalg.compare.goal}%
\end{equation}
This fact is completely trivial (since we can use
(\ref{eq.rmk.monalg.compare.alm=}) to rewrite the equation
(\ref{eq.rmk.monalg.compare.ass}) as $\left(  \alpha_{m}\right)  _{m\in
M}=\left(  \beta_{m}\right)  _{m\in M}$, but this is clearly equivalent to the
statement (\ref{eq.rmk.monalg.compare.goal})). Nevertheless, it is worth
pointing out, since it provides a dictionary between equalities in the monoid
algebra and equalities in the base ring $\mathbf{k}$. We will often use it to
deduce (\ref{eq.rmk.monalg.compare.goal}) from
(\ref{eq.rmk.monalg.compare.ass}) (since the other direction is obvious). This
deduction is called \emph{comparing coefficients}.
\end{remark}

Here is another useful convention:

\begin{convention}
\label{conv.monalg.e1=1}Let $M$ be a monoid. Then, we shall identify each
scalar $\lambda\in\mathbf{k}$ with the corresponding element $\lambda e_{1}%
\in\mathbf{k}\left[  M\right]  $ of the monoid algebra. In other words, we
identify each $\lambda\in\mathbf{k}$ with its image under the injective
$\mathbf{k}$-algebra morphism%
\begin{align*}
\mathbf{k}  &  \rightarrow\mathbf{k}\left[  M\right]  ,\\
\lambda &  \mapsto\lambda e_{1}.
\end{align*}

\end{convention}

Thus, for example, $\alpha e_{1}+\beta e_{m}$ becomes $\alpha1+\beta m$ (by
Convention \ref{conv.monalg.em=m}) or just $\alpha+\beta m$ (by Convention
\ref{conv.monalg.e1=1}).

Note that Convention \ref{conv.monalg.e1=1} is an instance of the standard
convention to embed the base ring $\mathbf{k}$ into each $\mathbf{k}$-algebra
$A$, at least if the canonical morphism%
\begin{align*}
\mathbf{k}  &  \rightarrow A,\\
\lambda &  \mapsto\lambda1_{A}%
\end{align*}
is injective.

\subsubsection{\label{subsec.intro.monalg.exas}Examples of monoid algebras}

Here are some of the most famous examples of monoid algebras:

\begin{enumerate}
\item Let $M$ be the monoid $\left\{  x^{0},x^{1},x^{2},\ldots\right\}  $ with
multiplication given by $x^{i}x^{j}=x^{i+j}$ and with neutral element $x^{0}$.
This monoid is just the familiar monoid $\left(  \mathbb{N},+,0\right)  $,
rewritten multiplicatively (by renaming each $n\in\mathbb{N}$ as $x^{n}$).

The monoid ring $\mathbf{k}\left[  M\right]  $ then has the standard basis
$\left(  e_{x^{0}},e_{x^{1}},e_{x^{2}},\ldots\right)  $, which (by Convention
\ref{conv.monalg.em=m}) we can just write as $\left(  x^{0},x^{1},x^{2}%
,\ldots\right)  $. Thus, each element of this monoid ring is a $\mathbf{k}%
$-linear combination of $x^{0},x^{1},x^{2},\ldots$. Moreover, the
multiplication is given by $e_{x^{i}}e_{x^{j}}=e_{x^{i}x^{j}}=e_{x^{i+j}}$,
that is, by $x^{i}x^{j}=x^{i+j}$ (again using Convention
\ref{conv.monalg.em=m}). This should sound familiar: It shows that this monoid
ring $\mathbf{k}\left[  M\right]  $ is just the univariate polynomial ring
$\mathbf{k}\left[  x\right]  $.

Note that this is not the same as the ring $\mathbf{k}\left[  \left[
x\right]  \right]  $ of formal power series, because we only allow essentially
finite $\mathbf{k}$-linear combinations (i.e., all but finitely many
coefficients must be $0$).

\item Let $M$ be the monoid $\left\{  \ldots,x^{-1},x^{0},x^{1},x^{2}%
,\ldots\right\}  $ with multiplication given by $x^{i}x^{j}=x^{i+j}$ and with
neutral element $x^{0}$. This monoid is just $\left(  \mathbb{Z},+,0\right)
$, rewritten multiplicatively (by renaming each $n\in\mathbb{Z}$ as $x^{n}$).

Then, the monoid ring $\mathbf{k}\left[  M\right]  $ is the Laurent polynomial
ring $\mathbf{k}\left[  x^{\pm1}\right]  $.

\item You can similarly obtain multivariate polynomial rings (or Laurent
polynomial rings) by starting with multivariate monomial monoids (i.e.,
monoids of the form $\left(  \mathbb{N}^{k},+,0\right)  $ or $\left(
\mathbb{Z}^{k},+,0\right)  $, rewritten multiplicatively).
\end{enumerate}

\subsection{The symmetric group algebra}

\begin{convention}
\label{conv.n}From now on, we fix a nonnegative integer $n$.
\end{convention}

\begin{definition}
\label{def.sga.sga}The $n$\emph{-th symmetric group algebra} (over
$\mathbf{k}$) means the group algebra of the $n$-th symmetric group $S_{n}$
over $\mathbf{k}$. In other words, it means the $\mathbf{k}$-algebra
$\mathbf{k}\left[  S_{n}\right]  $. We will be using Convention
\ref{conv.monalg.em=m} when working in this algebra.
\end{definition}

So the elements of $\mathbf{k}\left[  S_{n}\right]  $ are $\mathbf{k}$-linear
combinations of permutations of $\left[  n\right]  $. Note that the monoid
$S_{n}$ is finite, so that we don't need to distinguish between $\mathbf{k}%
^{S_{n}}$ and $\mathbf{k}^{\left(  S_{n}\right)  }$.

Let us play around with $\mathbf{k}\left[  S_{n}\right]  $ for $n=2$ and $n=3$:

\begin{example}
\label{exa.sga.S2}The symmetric group $S_{2}$ has just two elements: the
permutation $\operatorname*{id}$ and the simple transposition $s_{1}$ (recall
that $s_{1}$ is the transposition $t_{1,2}$, swapping $1$ with $2$). Thus, the
standard basis of $\mathbf{k}\left[  S_{2}\right]  $ is $\left(
e_{\operatorname*{id}},e_{s_{1}}\right)  =\left(  \operatorname*{id}%
,s_{1}\right)  $ (using Convention \ref{conv.monalg.em=m}). Hence, each
element of $\mathbf{k}\left[  S_{2}\right]  $ has the form%
\[
\alpha\operatorname*{id}+\,\beta s_{1}\ \ \ \ \ \ \ \ \ \text{\ for some
}\alpha,\beta\in\mathbf{k}%
\]
(of course, $\alpha\operatorname*{id}+\,\beta s_{1}$ really means $\alpha
e_{\operatorname*{id}}+\beta e_{s_{1}}$, since we are using Convention
\ref{conv.monalg.em=m}). Since $\operatorname*{id}$ is the neutral element of
$S_{2}$, we can also write $1$ for it, and we can rewrite $\alpha
\operatorname*{id}$ as $\alpha$ (by Convention \ref{conv.monalg.e1=1}). Hence,
each element of $\mathbf{k}\left[  S_{2}\right]  $ has the form%
\[
\alpha+\beta s_{1}\ \ \ \ \ \ \ \ \ \text{\ for some }\alpha,\beta
\in\mathbf{k}.
\]

Let us multiply two such elements:%
\begin{align*}
\left(  \alpha+\beta s_{1}\right)  \left(  \gamma+\delta s_{1}\right)   &
=\alpha\gamma+\alpha\delta s_{1}+\beta\gamma s_{1}+\beta\delta
\underbrace{s_{1}s_{1}}_{\substack{=\operatorname*{id}\\\text{(because
}t_{i,j}t_{i,j}=\operatorname*{id}\\\text{for any }i,j\text{)}}}\\
&  =\alpha\gamma+\alpha\delta s_{1}+\beta\gamma s_{1}+\underbrace{\beta
\delta\operatorname*{id}}_{=\beta\delta}\\
&  =\alpha\gamma+\alpha\delta s_{1}+\beta\gamma s_{1}+\beta\delta\\
&  =\left(  \alpha\gamma+\beta\delta\right)  +\left(  \alpha\delta+\beta
\gamma\right)  s_{1}.
\end{align*}
The unity of the algebra $\mathbf{k}\left[  S_{2}\right]  $ is the permutation
$\operatorname*{id}$, or more precisely the corresponding basis vector
$e_{\operatorname*{id}}$ (which we just call $\operatorname*{id}$).

The algebra $\mathbf{k}\left[  S_{2}\right]  $ is commutative, since the group
$S_{2}$ is abelian. (Indeed, $S_{2}$ is a cyclic group of order $2$.)
\end{example}

\begin{example}
\label{exa.sga.S3}The symmetric group $S_{3}$ has $6$ elements: the
permutations
\[
\operatorname*{oln}\left(  123\right)  ,\ \operatorname*{oln}\left(
132\right)  ,\ \operatorname*{oln}\left(  213\right)  ,\ \operatorname*{oln}%
\left(  231\right)  ,\ \operatorname*{oln}\left(  312\right)
,\ \operatorname*{oln}\left(  321\right)
\]
(see Subsection \ref{subsec.intro.nots.perms} for the meaning of this
notation). Hence, each element of the symmetric group algebra $\mathbf{k}%
\left[  S_{3}\right]  $ has the form%
\[
\alpha\underbrace{\operatorname*{oln}\left(  123\right)  }%
_{=\operatorname*{id}=1}+\,\beta\underbrace{\operatorname*{oln}\left(
132\right)  }_{=s_{2}}+\,\gamma\underbrace{\operatorname*{oln}\left(
213\right)  }_{=s_{1}}+\,\delta\underbrace{\operatorname*{oln}\left(
231\right)  }_{\substack{=\operatorname*{cyc}\nolimits_{1,2,3}\\=s_{1}s_{2}%
}}+\,\varepsilon\underbrace{\operatorname*{oln}\left(  312\right)
}_{\substack{=\operatorname*{cyc}\nolimits_{1,3,2}\\=s_{2}s_{1}}%
}+\,\zeta\underbrace{\operatorname*{oln}\left(  321\right)  }%
_{\substack{=t_{1,3}\\=s_{1}s_{2}s_{1}\\=s_{2}s_{1}s_{2}}}
\]
for some $\alpha,\beta,\gamma,\delta,\varepsilon,\zeta\in\mathbf{k}$. Just
like in Example \ref{exa.sga.S2}, we could compute the product of any two such
elements. Instead of the general formula, let me show two examples of such
products:%
\begin{align}
\left(  1+s_{1}\right)  \left(  1-s_{1}\right)   &  =1+s_{1}-s_{1}-s_{1}%
s_{1}=1-\underbrace{s_{1}s_{1}}_{=\operatorname*{id}=1}=1-1\nonumber\\
&  =0 \label{eq.exa.sga.S3.f1}%
\end{align}
and%
\begin{align}
\left(  1+s_{2}\right)  \left(  1+s_{1}+s_{1}s_{2}\right)   &  =1+s_{1}%
+s_{1}s_{2}+s_{2}+s_{2}s_{1}+s_{2}s_{1}s_{2}\nonumber\\
&  =\sum_{w\in S_{3}}w. \label{eq.exa.sga.S3.f2}%
\end{align}

\end{example}

The product (\ref{eq.exa.sga.S3.f1}) can be easily generalized: We have%
\[
\left(  1+t_{i,j}\right)  \left(  1-t_{i,j}\right)
=0\ \ \ \ \ \ \ \ \ \ \text{for any transposition }t_{i,j}\in S_{n}.
\]
Even more generally,
\[
\left(  1+w\right)  \left(  1-w\right)  =0\ \ \ \ \ \ \ \ \ \ \text{for any
involution }w\in S_{n}.
\]
(An \emph{involution} means a permutation $w$ such that $w^{2}%
=\operatorname*{id}$, that is, $w=w^{-1}$.)

Soon we will also generalize the product (\ref{eq.exa.sga.S3.f2}).

We note that the $0$-th and the $1$-st symmetric group algebras $\mathbf{k}%
\left[  S_{0}\right]  $ and $\mathbf{k}\left[  S_{1}\right]  $ are both
isomorphic to $\mathbf{k}$, since the groups $S_{0}$ and $S_{1}$ are trivial
(i.e., groups of size $1$). These are not very interesting objects, although
their existence can be helpful in induction proofs.

\subsection{The integral and the sign-integral}

We now define two peculiar elements of any symmetric group algebra.

\begin{definition}
\label{def.sga.integral}The \emph{integral} (aka \emph{symmetrizer}) of the
symmetric group algebra $\mathbf{k}\left[  S_{n}\right]  $ is defined to be
the element%
\[
\nabla:=\sum_{w\in S_{n}}w\in\mathbf{k}\left[  S_{n}\right]  .
\]
The \emph{sign-integral} (aka \emph{antisymmetrizer}) of the symmetric group
algebra $\mathbf{k}\left[  S_{n}\right]  $ is defined to be the element%
\[
\nabla^{-}:=\sum_{w\in S_{n}}\underbrace{\left(  -1\right)  ^{w}}_{\text{the
sign of }w}w\in\mathbf{k}\left[  S_{n}\right]  .
\]

\end{definition}

\begin{example}
For $n=2$, we have%
\[
\nabla=1+s_{1}\ \ \ \ \ \ \ \ \ \ \text{and}\ \ \ \ \ \ \ \ \ \ \nabla
^{-}=1-s_{1}.
\]

For $n=3$, we have%
\begin{align*}
\nabla &  =1+s_{1}+s_{2}+s_{1}s_{2}+s_{2}s_{1}+s_{2}s_{1}s_{2}%
\ \ \ \ \ \ \ \ \ \ \text{and}\\
\nabla^{-}  &  =1-s_{1}-s_{2}+s_{1}s_{2}+s_{2}s_{1}-s_{2}s_{1}s_{2}.
\end{align*}

For $n=0$ or $n=1$, we have $\nabla=\nabla^{-}=1$.
\end{example}

Multiplying the integral or the sign-integral by another permutation has a
rather simple outcome:

\begin{proposition}
\label{prop.integral.fix}Let $w\in S_{n}$. Then,%
\begin{align}
w\nabla &  =\nabla w=\nabla\ \ \ \ \ \ \ \ \ \ \text{and}%
\label{eq.prop.integral.fix.+}\\
w\nabla^{-}  &  =\nabla^{-}w=\left(  -1\right)  ^{w}\nabla^{-}.
\label{eq.prop.integral.fix.-}%
\end{align}

\end{proposition}

\begin{proof}
The symmetric group $S_{n}$ is a group. Thus, the map%
\begin{align*}
S_{n}  &  \rightarrow S_{n},\\
u  &  \mapsto wu
\end{align*}
is a bijection (with inverse given by $u\mapsto w^{-1}u$). But the definition
of $\nabla$ can be rewritten as $\nabla=\sum\limits_{u\in S_{n}}u$. Thus,%
\begin{align*}
w\nabla &  =w\sum\limits_{u\in S_{n}}u=\sum\limits_{u\in S_{n}}wu=\sum_{u\in
S_{n}}u\\
&  \ \ \ \ \ \ \ \ \ \ \ \ \ \ \ \ \ \ \ \ \left(
\begin{array}
[c]{c}%
\text{here, we substituted }u\text{ for }wu\text{ in the sum, since}\\
\text{the map }S_{n}\rightarrow S_{n},\ u\mapsto wu\text{ is a bijection}%
\end{array}
\right) \\
&  =\nabla.
\end{align*}

Similarly, $\nabla w=\nabla$. Combining these two equalities, we obtain
(\ref{eq.prop.integral.fix.+}).

Next, let us prove (\ref{eq.prop.integral.fix.-}). One of the basic properties
of signs of permutations says that $\left(  -1\right)  ^{\sigma\tau}=\left(
-1\right)  ^{\sigma}\cdot\left(  -1\right)  ^{\tau}$ for any $\sigma,\tau\in
S_{n}$. Thus, for any $u\in S_{n}$, we have%
\[
\left(  -1\right)  ^{wu}=\left(  -1\right)  ^{w}\cdot\left(  -1\right)  ^{u},
\]
so that%
\begin{align}
\left(  -1\right)  ^{u}  &  =\underbrace{\dfrac{1}{\left(  -1\right)  ^{w}}%
}_{\substack{=\left(  -1\right)  ^{w}\\\text{(since the number }\left(
-1\right)  ^{w}\text{ is either }1\text{ or }-1\text{,}\\\text{and thus equals
its own reciprocal)}}}\cdot\left(  -1\right)  ^{wu}\nonumber\\
&  =\left(  -1\right)  ^{w}\cdot\left(  -1\right)  ^{wu}.
\label{pf.prop.integral.fix.sign}%
\end{align}

Now, rewrite the definition of $\nabla^{-}$ as $\nabla^{-}=\sum\limits_{u\in
S_{n}}\left(  -1\right)  ^{u}u$. Thus,%
\begin{align*}
w\nabla^{-}  &  =w\sum\limits_{u\in S_{n}}\left(  -1\right)  ^{u}%
u=\sum\limits_{u\in S_{n}}\underbrace{\left(  -1\right)  ^{u}}%
_{\substack{=\left(  -1\right)  ^{w}\cdot\left(  -1\right)  ^{wu}\\\text{(by
(\ref{pf.prop.integral.fix.sign}))}}}wu\\
&  =\sum_{u\in S_{n}}\left(  -1\right)  ^{w}\cdot\left(  -1\right)
^{wu}wu=\sum_{u\in S_{n}}\left(  -1\right)  ^{w}\cdot\left(  -1\right)
^{u}u\\
&  \ \ \ \ \ \ \ \ \ \ \ \ \ \ \ \ \ \ \ \ \left(
\begin{array}
[c]{c}%
\text{here, we substituted }u\text{ for }wu\text{ in the sum, since}\\
\text{the map }S_{n}\rightarrow S_{n},\ u\mapsto wu\text{ is a bijection}%
\end{array}
\right) \\
&  =\left(  -1\right)  ^{w}\cdot\underbrace{\sum_{u\in S_{n}}\left(
-1\right)  ^{u}u}_{=\nabla^{-}}=\left(  -1\right)  ^{w}\nabla^{-}.
\end{align*}
Similarly, $\nabla^{-}w=\left(  -1\right)  ^{w}\nabla^{-}$. Combining these
two equalities, we obtain (\ref{eq.prop.integral.fix.-}). This completes the
proof of Proposition \ref{prop.integral.fix}.
\end{proof}

Note that commutativity in $\mathbf{k}\left[  S_{n}\right]  $ is rare (for
$n\geq3$), so that the equalities $w\nabla=\nabla w$ and $w\nabla^{-}%
=\nabla^{-}w$ in Proposition \ref{prop.integral.fix} should be somewhat
surprising (but not too much, since $\nabla$ and $\nabla^{-}$ are two of the
simplest elements of $\mathbf{k}\left[  S_{n}\right]  $).

The relation (\ref{eq.prop.integral.fix.+}), in a sense, characterizes
$\nabla$ (up to scalar multiplication): Indeed, the only elements
$\mathbf{a}\in\mathbf{k}\left[  S_{n}\right]  $ that satisfy $\mathbf{a}%
w=\mathbf{a}$ for all $w\in S_{n}$ are $\nabla$ and its scalar multiples. This
is part of the following exercise:

\begin{exercise}
\label{exe.integral.fix-converse}Let $\mathbf{a}\in\mathbf{k}\left[
S_{n}\right]  $. Prove the following: \medskip

\textbf{(a)} \fbox{1} If $\mathbf{a}w=\mathbf{a}$ for all $w\in S_{n}$, then
$\mathbf{a}=\lambda\nabla$ for some $\lambda\in\mathbf{k}$. \medskip

\textbf{(b)} \fbox{1} Even better: If $\mathbf{a}s_{i}=\mathbf{a}$ for all
$i\in\left[  n-1\right]  $, then $\mathbf{a}=\lambda\nabla$ for some
$\lambda\in\mathbf{k}$.
\end{exercise}

\begin{corollary}
\label{cor.integral.square}We have $\nabla^{2}=n!\cdot\nabla$.
\end{corollary}

\begin{proof}
We have%
\begin{align*}
\nabla^{2}  &  =\nabla\nabla=\left(  \sum_{w\in S_{n}}w\right)  \nabla
\ \ \ \ \ \ \ \ \ \ \left(  \text{since }\nabla=\sum_{w\in S_{n}}w\right) \\
&  =\sum_{w\in S_{n}}\underbrace{w\nabla}_{\substack{=\nabla\\\text{(by
(\ref{eq.prop.integral.fix.+}))}}}=\sum_{w\in S_{n}}\nabla=n!\cdot\nabla
\end{align*}
(since there are $n!$ many permutations $w\in S_{n}$). This proves Corollary
\ref{cor.integral.square}.
\end{proof}

Corollary \ref{cor.integral.square} can be restated as follows: $\nabla\left(
\nabla-n!\right)  =0$. In other words, $P\left(  \nabla\right)  =0$, where
$P\in\mathbb{Z}\left[  x\right]  $ is the polynomial $x\left(  x-n!\right)  $.
This is the first sound of a recurring theme: Often, you define an element
$\mathbf{a}$ of $\mathbf{k}\left[  S_{n}\right]  $ and find that there is a
\textquotedblleft nice\textquotedblright\ polynomial $P\in\mathbf{k}\left[
x\right]  $ such that $P\left(  \mathbf{a}\right)  =0$. Actually, we will soon
see that there always is a monic polynomial $P\in\mathbf{k}\left[  x\right]  $
of degree $n!$ such that $P\left(  \mathbf{a}\right)  =0$. But in many cases,
there exist such polynomials $P$ of much smaller degree than $n!$, and
sometimes these $P$'s factor into linear factors over $\mathbb{Z}$. Many
interesting results have this form; Corollary \ref{cor.integral.square} is the
simplest of them.

\begin{remark}
\label{rmk.integral.G}We can generalize $\nabla$ to an arbitrary finite group
$G$. Indeed, if $G$ is any finite group, then we can define an element
$\nabla_{G}:=\sum\limits_{g\in G}g\in\mathbf{k}\left[  G\right]  $, and then
every $w\in G$ satisfies%
\begin{equation}
w\nabla_{G}=\nabla_{G}w=\nabla_{G}. \label{eq.rmk.integral.G.+}%
\end{equation}
This can be proved in the same way as we proved (\ref{eq.prop.integral.fix.+})
(and indeed, (\ref{eq.prop.integral.fix.+}) is the particular case of
(\ref{eq.rmk.integral.G.+}) for $G=S_{n}$). As a consequence, $\nabla_{G}%
^{2}=\left\vert G\right\vert \cdot\nabla_{G}$ (which generalizes Corollary
\ref{cor.integral.square}).

For infinite groups $G$, we cannot define $\nabla_{G}$, since linear
combinations are not allowed to have infinitely many nonzero coefficients.

The equality (\ref{eq.prop.integral.fix.-}) can also be generalized with a bit
more work: We let $\mathbf{k}^{\times}$ denote the group of units of
$\mathbf{k}$ (that is, the group of all invertible elements of $\mathbf{k}$,
with respect to multiplication). If $f:G\rightarrow\mathbf{k}^{\times}$ is a
group morphism from $G$ to $\mathbf{k}^{\times}$, then we can define an
element $\nabla_{G}^{f}:=\sum\limits_{g\in G}f\left(  g\right)  g\in
\mathbf{k}\left[  G\right]  $, and then every $w\in G$ satisfies%
\begin{equation}
w\nabla_{G}^{f}=\nabla_{G}^{f}w=f\left(  w^{-1}\right)  \cdot\nabla_{G}^{f}.
\label{eq.rmk.integral.G.f}%
\end{equation}
If we let $G=S_{n}$ and let $f$ be the sign homomorphism (so that $f\left(
w\right)  =\operatorname*{sign}w=\left(  -1\right)  ^{w}$ for each $w\in
S_{n}$), then this recovers (\ref{eq.prop.integral.fix.-}).
\end{remark}

\subsection{\label{sec.intro.YJM}The Young--Jucys--Murphy elements (aka YJM
elements, or jucys--murphies)}

\subsubsection{\label{subsec.intro.YJM.def}Definition and commutativity}

It is time to define the first family of elements of $\mathbf{k}\left[
S_{n}\right]  $ that have nontrivial interesting properties:

\begin{definition}
\label{def.YJM.mk}For each $k\in\left[  n\right]  $, we define the
$k$\emph{-th Young--Jucys--Murphy element} (for short: \emph{YJM element}, or
just \emph{jucys--murphy}) in $\mathbf{k}\left[  S_{n}\right]  $ to be the
element%
\[
\mathbf{m}_{k}:=t_{1,k}+t_{2,k}+\cdots+t_{k-1,k}=\sum_{i=1}^{k-1}t_{i,k}.
\]
This is the sum of all transpositions $t_{i,k}$ with $i<k$.
\end{definition}

In particular,%
\begin{align*}
\mathbf{m}_{1}  &  =0;\\
\mathbf{m}_{2}  &  =t_{1,2};\\
\mathbf{m}_{3}  &  =t_{1,3}+t_{2,3};\\
\mathbf{m}_{4}  &  =t_{1,4}+t_{2,4}+t_{3,4};\\
&  \ldots;\\
\mathbf{m}_{n}  &  =t_{1,n}+t_{2,n}+\cdots+t_{n-1,n}.
\end{align*}
Thus,%
\[
\mathbf{m}_{1}+\mathbf{m}_{2}+\cdots+\mathbf{m}_{n}=\sum_{k=1}^{n}%
\mathbf{m}_{k}=\sum_{k=1}^{n}\ \ \sum_{i=1}^{k-1}t_{i,k}%
\]
is the sum of \textbf{all} transpositions in $S_{n}$.

Here is our first surprise:

\begin{theorem}
\label{thm.YJM.commute}The jucys--murphies $\mathbf{m}_{1},\mathbf{m}%
_{2},\ldots,\mathbf{m}_{n}$ commute. In other words,%
\[
\mathbf{m}_{i}\mathbf{m}_{j}=\mathbf{m}_{j}\mathbf{m}_{i}%
\ \ \ \ \ \ \ \ \ \ \text{for all }i,j\in\left[  n\right]  .
\]

\end{theorem}

\begin{proof}
Let $i,j\in\left[  n\right]  $. We must prove that $\mathbf{m}_{i}%
\mathbf{m}_{j}=\mathbf{m}_{j}\mathbf{m}_{i}$. If $i=j$, then this is obvious.
Thus, we WLOG assume that $i\neq j$. Furthermore, our claim ($\mathbf{m}%
_{i}\mathbf{m}_{j}=\mathbf{m}_{j}\mathbf{m}_{i}$) is symmetric in $i$ and $j$,
so that we can swap $i$ with $j$ at will. Hence, we furthermore WLOG assume
that $i\leq j$. Thus, $i<j$ (since $i\neq j$), so that $i\in\left[
j-1\right]  $.

Recall that
\begin{align*}
\mathbf{m}_{i}  &  =t_{1,i}+t_{2,i}+\cdots+t_{i-1,i}=\sum_{p=1}^{i-1}%
t_{p,i}\ \ \ \ \ \ \ \ \ \ \text{and}\\
\mathbf{m}_{j}  &  =t_{1,j}+t_{2,j}+\cdots+t_{j-1,j}=\sum_{k=1}^{j-1}t_{k,j}.
\end{align*}
Now we shall show that%
\begin{equation}
t_{p,i}\mathbf{m}_{j}=\mathbf{m}_{j}t_{p,i}\text{ for each }p\in\left[
i-1\right]  . \label{pf.thm.YJM.commute.1}%
\end{equation}

[\textit{Proof of (\ref{pf.thm.YJM.commute.1}):} Let $p\in\left[  i-1\right]
$. Then, $p\leq i-1<i<j$. Hence, each $k\in\left[  j-1\right]  \setminus
\left\{  p,i\right\}  $ satisfies $k\notin\left\{  j,p,i\right\}  $ and
therefore
\begin{equation}
t_{p,i}t_{k,j}=t_{k,j}t_{p,i} \label{pf.thm.YJM.commute.1.pf.k}%
\end{equation}
(this is a particular case of (\ref{eq.intro.perms.cycs.cc=cc}), since the
four numbers $p,i,k,j$ are distinct). Moreover, we can easily check that%
\begin{equation}
t_{p,i}t_{p,j}=t_{i,j}t_{p,i} \label{pf.thm.YJM.commute.1.pf.p}%
\end{equation}
(indeed, both sides here equal $\operatorname*{cyc}\nolimits_{i,p,j}$) and%
\begin{equation}
t_{p,i}t_{i,j}=t_{p,j}t_{p,i} \label{pf.thm.YJM.commute.1.pf.i}%
\end{equation}
(indeed, both sides here equal $\operatorname*{cyc}\nolimits_{i,j,p}$).

Now, $p$ and $i$ are two distinct elements of $\left[  j\right]  $ (since
$p<i<j$). Recall that%
\[
\mathbf{m}_{j}=\sum_{k=1}^{j-1}t_{k,j}=t_{p,j}+t_{i,j}+\sum_{k\in\left[
j-1\right]  \setminus\left\{  p,i\right\}  }t_{k,j}%
\]
(here, we have split off the addends for $k=p$ and for $k=i$ from the sum,
since $p$ and $i$ are two distinct elements of $\left[  j-1\right]  $). Hence,%
\begin{align*}
t_{p,i}\mathbf{m}_{j}  &  =t_{p,i}\left(  t_{p,j}+t_{i,j}+\sum_{k\in\left[
j-1\right]  \setminus\left\{  p,i\right\}  }t_{k,j}\right) \\
&  =\underbrace{t_{p,i}t_{p,j}}_{\substack{=t_{i,j}t_{p,i}\\\text{(by
(\ref{pf.thm.YJM.commute.1.pf.p}))}}}+\underbrace{t_{p,i}t_{i,j}%
}_{\substack{=t_{p,j}t_{p,i}\\\text{(by (\ref{pf.thm.YJM.commute.1.pf.i}))}%
}}+\sum_{k\in\left[  j-1\right]  \setminus\left\{  p,i\right\}  }%
\underbrace{t_{p,i}t_{k,j}}_{\substack{=t_{k,j}t_{p,i}\\\text{(by
(\ref{pf.thm.YJM.commute.1.pf.k}))}}}\\
&  =t_{i,j}t_{p,i}+t_{p,j}t_{p,i}+\sum_{k\in\left[  j-1\right]  \setminus
\left\{  p,i\right\}  }t_{k,j}t_{p,i}\\
&  =\underbrace{\left(  t_{i,j}+t_{p,j}+\sum_{k\in\left[  j-1\right]
\setminus\left\{  p,i\right\}  }t_{k,j}\right)  }_{\substack{=t_{p,j}%
+t_{i,j}+\sum_{k\in\left[  j-1\right]  \setminus\left\{  p,i\right\}  }%
t_{k,j}\\=\mathbf{m}_{j}}}t_{p,i}=\mathbf{m}_{j}t_{p,i}.
\end{align*}
This proves (\ref{pf.thm.YJM.commute.1}).] \medskip

Now, from $\mathbf{m}_{i}=\sum\limits_{p=1}^{i-1}t_{p,i}$, we obtain%
\[
\mathbf{m}_{i}\mathbf{m}_{j}=\left(  \sum\limits_{p=1}^{i-1}t_{p,i}\right)
\mathbf{m}_{j}=\sum\limits_{p=1}^{i-1}\underbrace{t_{p,i}\mathbf{m}_{j}%
}_{\substack{=\mathbf{m}_{j}t_{p,i}\\\text{(by (\ref{pf.thm.YJM.commute.1}))}%
}}=\sum\limits_{p=1}^{i-1}\mathbf{m}_{j}t_{p,i}=\mathbf{m}_{j}%
\underbrace{\left(  \sum\limits_{p=1}^{i-1}t_{p,i}\right)  }_{=\mathbf{m}_{i}%
}=\mathbf{m}_{j}\mathbf{m}_{i}.
\]
This completes the proof of Theorem \ref{thm.YJM.commute}.
\end{proof}

\subsubsection{A product}

Now let us compute the product $\left(  1+\mathbf{m}_{1}\right)  \left(
1+\mathbf{m}_{2}\right)  \cdots\left(  1+\mathbf{m}_{n}\right)  $. For $n=3$,
it is%
\begin{align*}
\left(  1+\mathbf{m}_{1}\right)  \left(  1+\mathbf{m}_{2}\right)  \left(
1+\mathbf{m}_{3}\right)   &  =1\left(  1+t_{1,2}\right)  \left(
1+t_{1,3}+t_{2,3}\right) \\
&  =1+t_{1,2}+t_{1,3}+t_{2,3}+t_{1,2}t_{1,3}+t_{1,2}t_{2,3}\\
&  =\sum_{w\in S_{3}}w=\nabla.
\end{align*}
This is not a coincidence:

\begin{theorem}
\label{thm.YJM.prod1+mk}We have%
\[
\left(  1+\mathbf{m}_{1}\right)  \left(  1+\mathbf{m}_{2}\right)
\cdots\left(  1+\mathbf{m}_{n}\right)  =\nabla.
\]

\end{theorem}

\begin{proof}
[First proof of Theorem \ref{thm.YJM.prod1+mk}.]We have so far defined the
transpositions $t_{i,j}$ only for $i\neq j$. Now, let us also set
$t_{k,k}:=\operatorname*{id}$ for each $k\in\left[  n\right]  $. This is not
really a transposition, but it still fits the pattern of the $t_{i,j}$'s,
since it \textquotedblleft swaps $k$ with $k$ and leaves everything else
unchanged\textquotedblright. In other words, we still have $t_{i,j}\left(
i\right)  =j$ and $t_{i,j}\left(  j\right)  =i$ for all $i,j\in\left[
n\right]  $, as well as $t_{i,j}\left(  k\right)  =k$ for all $i,j,k\in\left[
n\right]  $ that satisfy $k\notin\left\{  i,j\right\}  $.

Now, for each $k\in\left[  n\right]  $, we have%
\[
\underbrace{1}_{=\operatorname*{id}=t_{k,k}}+\underbrace{\mathbf{m}_{k}%
}_{=\sum_{j=1}^{k-1}t_{j,k}}=t_{k,k}+\sum_{j=1}^{k-1}t_{j,k}=\sum_{j=1}%
^{k}t_{j,k}.
\]
Multiplying these equalities for all $k\in\left[  n\right]  $, we obtain%
\begin{align*}
&  \left(  1+\mathbf{m}_{1}\right)  \left(  1+\mathbf{m}_{2}\right)
\cdots\left(  1+\mathbf{m}_{n}\right) \\
&  =\left(  \sum_{j=1}^{1}t_{j,1}\right)  \left(  \sum_{j=1}^{2}%
t_{j,2}\right)  \cdots\left(  \sum_{j=1}^{n}t_{j,n}\right) \\
&  =\sum_{\left(  j_{1},j_{2},\ldots,j_{n}\right)  \in\left[  1\right]
\times\left[  2\right]  \times\cdots\times\left[  n\right]  }t_{j_{1}%
,1}t_{j_{2},2}\cdots t_{j_{n},n}%
\end{align*}
(by the product rule, i.e., by expanding the product).\footnote{The condition
$\left(  j_{1},j_{2},\ldots,j_{n}\right)  \in\left[  1\right]  \times\left[
2\right]  \times\cdots\times\left[  n\right]  $ under the summation sign says
precisely that $j_{1},j_{2},\ldots,j_{n}$ are positive integers satisfying
$j_{k}\leq k$ for each $k$.}

Our goal is to prove that this sum equals $\nabla$, i.e., that it contains
each permutation $w\in S_{n}$ exactly once. In other words, we must prove the
following theorem:

\begin{theorem}
\label{thm.perm.transp-code}Each $w\in S_{n}$ can be written as $t_{j_{1}%
,1}t_{j_{2},2}\cdots t_{j_{n},n}$ for exactly one $n$-tuple $\left(
j_{1},j_{2},\ldots,j_{n}\right)  \in\left[  1\right]  \times\left[  2\right]
\times\cdots\times\left[  n\right]  $. In other words, the map%
\begin{align*}
\left[  1\right]  \times\left[  2\right]  \times\cdots\times\left[  n\right]
&  \rightarrow S_{n},\\
\left(  j_{1},j_{2},\ldots,j_{n}\right)   &  \mapsto t_{j_{1},1}t_{j_{2}%
,2}\cdots t_{j_{n},n}%
\end{align*}
is a bijection.
\end{theorem}

\begin{proof}
[Proof of Theorem \ref{thm.perm.transp-code}.]Let $\Phi$ be the map%
\begin{align*}
\left[  1\right]  \times\left[  2\right]  \times\cdots\times\left[  n\right]
&  \rightarrow S_{n},\\
\left(  j_{1},j_{2},\ldots,j_{n}\right)   &  \mapsto t_{j_{1},1}t_{j_{2}%
,2}\cdots t_{j_{n},n}.
\end{align*}
We must show that this map $\Phi$ is a bijection.

Let us first prove that it is injective. In other words, let us prove that any
$n$-tuple $\left(  j_{1},j_{2},\ldots,j_{n}\right)  \in\left[  1\right]
\times\left[  2\right]  \times\cdots\times\left[  n\right]  $ can be
reconstructed from its image $\Phi\left(  j_{1},j_{2},\ldots,j_{n}\right)
=t_{j_{1},1}t_{j_{2},2}\cdots t_{j_{n},n}$.

To do so, we pick an $n$-tuple $\left(  j_{1},j_{2},\ldots,j_{n}\right)
\in\left[  1\right]  \times\left[  2\right]  \times\cdots\times\left[
n\right]  $, and we set
\[
w:=\Phi\left(  j_{1},j_{2},\ldots,j_{n}\right)  =t_{j_{1},1}t_{j_{2},2}\cdots
t_{j_{n},n}.
\]
How can we reconstruct $\left(  j_{1},j_{2},\ldots,j_{n}\right)  $ from this
$w$ ? Let us define the permutation%
\[
w_{k}:=t_{j_{1},1}t_{j_{2},2}\cdots t_{j_{k},k}\in S_{n}%
\ \ \ \ \ \ \ \ \ \ \text{for each }k\in\left\{  0,1,\ldots,n\right\}  .
\]
Thus, $w_{0}=\left(  \text{empty product}\right)  =1=\operatorname*{id}$ and
$w_{n}=t_{j_{1},1}t_{j_{2},2}\cdots t_{j_{n},n}=w$. The definition of $w_{k}$
shows that%
\begin{equation}
w_{k}=w_{k-1}t_{j_{k},k}\ \ \ \ \ \ \ \ \ \ \text{for each }k\in\left[
n\right]  \label{pf.thm.perm.transp-code.wk=wk-1}%
\end{equation}
(since $\underbrace{w_{k-1}}_{=t_{j_{1},1}t_{j_{2},2}\cdots t_{j_{k-1},k-1}%
}t_{j_{k},k}=\left(  t_{j_{1},1}t_{j_{2},2}\cdots t_{j_{k-1},k-1}\right)
t_{j_{k},k}=t_{j_{1},1}t_{j_{2},2}\cdots t_{j_{k},k}=w_{k}$). We can now
reconstruct (in this order) $w_{n},j_{n},w_{n-1},j_{n-1},w_{n-2}%
,j_{n-2},\ldots$ (from $w$) by the following recursive procedure:

\begin{enumerate}
\item We know $w_{n}$ already, since $w_{n}=w$.

\item I claim that $j_{n}=w_{n}^{-1}\left(  n\right)  $. To see this, note
that $w_{n}=t_{j_{1},1}t_{j_{2},2}\cdots t_{j_{n},n}$ entails%
\begin{align*}
w_{n}\left(  j_{n}\right)   &  =\left(  t_{j_{1},1}t_{j_{2},2}\cdots
t_{j_{n},n}\right)  \left(  j_{n}\right)  =\left(  t_{j_{1},1}t_{j_{2}%
,2}\cdots t_{j_{n-1},n-1}\right)  \left(  \underbrace{t_{j_{n},n}\left(
j_{n}\right)  }_{=n}\right) \\
&  =\left(  t_{j_{1},1}t_{j_{2},2}\cdots t_{j_{n-1},n-1}\right)  \left(
n\right)  =n,
\end{align*}
since each of the transpositions $t_{j_{1},1},t_{j_{2},2},\ldots
,t_{j_{n-1},n-1}$ leaves $n$ unchanged (because $j_{1}\in\left[  1\right]  $,
$j_{2}\in\left[  2\right]  $, $\ldots$, $j_{n-1}\in\left[  n-1\right]  $
ensures that all the numbers $j_{1},j_{2},\ldots,j_{n-1}$ are $<n$). Hence,
$j_{n}=w_{n}^{-1}\left(  n\right)  $, as I claimed.

Thus, $j_{n}$ has been reconstructed from $w$.

\item Next, I claim that $w_{n-1}=w_{n}t_{j_{n},n}^{-1}$. (Note that
$t_{j_{n},n}^{-1}$ is the same as $t_{j_{n},n}$, but we write it as
$t_{j_{n},n}^{-1}$ for clarity.) Indeed,
(\ref{pf.thm.perm.transp-code.wk=wk-1}) (applied to $k=n$) yields
$w_{n}=w_{n-1}t_{j_{n},n}$. Solving this for $w_{n-1}$, we obtain
$w_{n-1}=w_{n}t_{j_{n},n}^{-1}$.

Thus, $w_{n-1}$ has been reconstructed (since $j_{n}$ and $w_{n}$ have been reconstructed).

\item Just as we showed that $j_{n}=w_{n}^{-1}\left(  n\right)  $, we can show
that $j_{n-1}=w_{n-1}^{-1}\left(  n-1\right)  $ (since $t_{j_{n-1},n-1}\left(
j_{n-1}\right)  =n-1$, and since each of the transpositions \newline%
$t_{j_{1},1},t_{j_{2},2},\ldots,t_{j_{n-2},n-2}$ leaves $n-1$ unchanged).
Thus, $j_{n-1}$ has been reconstructed from $w$.

\item Just as we showed that $w_{n-1}=w_{n}t_{j_{n},n}^{-1}$, we can show that
$w_{n-2}=w_{n-1}t_{j_{n-1},n-1}^{-1}$. Thus, $w_{n-2}$ has been reconstructed.

\item We continue along the same lines. That is, for each $k$ ranging from $n$
to $1$, we obtain $j_{k}$ from $w_{k}$ by the formula $j_{k}=w_{k}^{-1}\left(
k\right)  $, and then obtain $w_{k-1}$ by the formula $w_{k-1}=w_{k}%
t_{j_{k},k}^{-1}$.
\end{enumerate}

At the end of this procedure, all of $j_{1},j_{2},\ldots,j_{n}$ have been
reconstructed, so that the $n$-tuple $\left(  j_{1},j_{2},\ldots,j_{n}\right)
$ has been reconstructed. This shows that our map $\Phi$ is injective.

As we recall, our goal is to prove that $\Phi$ is a bijection, i.e.,
bijective. A very quick shortcut is available towards this goal: The sets
$\left[  1\right]  \times\left[  2\right]  \times\cdots\times\left[  n\right]
$ and $S_{n}$ are two finite sets of the same size (namely, $n!$), and thus
any injective map from one to the other must be bijective (by the pigeonhole
principle). Thus, $\Phi$ is bijective (since $\Phi$ is injective), and Theorem
\ref{thm.perm.transp-code} is proved.

This trick argument is perfectly legitimate, but let me also show a slower but
more \textquotedblleft natural\textquotedblright\ way to finish this proof of
Theorem \ref{thm.perm.transp-code}. We know that $\Phi$ is injective; we must
prove that $\Phi$ is bijective. So it remains to prove that $\Phi$ is
surjective. In other words, we must show that each $w\in S_{n}$ can be written
as $\Phi\left(  j_{1},j_{2},\ldots,j_{n}\right)  $ for some $\left(
j_{1},j_{2},\ldots,j_{n}\right)  \in\left[  1\right]  \times\left[  2\right]
\times\cdots\times\left[  n\right]  $.

So let $w\in S_{n}$ be arbitrary. Recall the recursive procedure that we used
to reconstruct $j_{1},j_{2},\ldots,j_{n}$ in the above proof of the
injectivity of $\Phi$. We can apply this procedure to $w$ without knowing
a-priori that $j_{1},j_{2},\ldots,j_{n}$ exist. Thus,

\begin{enumerate}
\item we define $w_{n}\in S_{n}$ by $w_{n}:=w$;

\item we define $j_{n}\in\left[  n\right]  $ by $j_{n}:=w_{n}^{-1}\left(
n\right)  $;

\item we define $w_{n-1}\in S_{n}$ by $w_{n-1}:=w_{n}t_{j_{n},n}^{-1}$;

\item we define $j_{n-1}\in\left[  n\right]  $ by $j_{n-1}:=w_{n-1}%
^{-1}\left(  n-1\right)  $;

\item we define $w_{n-2}\in S_{n}$ by $w_{n-2}:=w_{n-1}t_{j_{n-1},n-1}^{-1}$;

\item and so on.
\end{enumerate}

Thus, we have recursively defined $n+1$ permutations $w_{n},w_{n-1}%
,\ldots,w_{0}\in S_{n}$ and $n$ elements $j_{n},j_{n-1},\ldots,j_{1}\in\left[
n\right]  $ by the equalities $w_{n}:=w$ and
\begin{align}
w_{k-1}  &  :=w_{k}t_{j_{k},k}^{-1}\ \ \ \ \ \ \ \ \ \ \text{and}%
\label{pf.thm.perm.transp-code.sur.wk-1=}\\
j_{k}  &  :=w_{k}^{-1}\left(  k\right)
\label{pf.thm.perm.transp-code.sur.jk=}%
\end{align}
for all $k\in\left[  n\right]  $. It remains to show that $\left(  j_{1}%
,j_{2},\ldots,j_{n}\right)  \in\left[  1\right]  \times\left[  2\right]
\times\cdots\times\left[  n\right]  $ and that $\Phi\left(  j_{1},j_{2}%
,\ldots,j_{n}\right)  =w$. Let us now do this. We claim that\footnote{A map
$f$ is said to \emph{fix} an element $x$ if $f\left(  x\right)  =x$.}
\begin{equation}
w_{k}\text{ fixes }k+1,k+2,\ldots,n \label{pf.thm.perm.transp-code.sur.fixes}%
\end{equation}
(i.e., that $w_{k}\left(  i\right)  =i$ for each $i\in\left\{  k+1,k+2,\ldots
,n\right\}  $) for each $k\in\left\{  0,1,\ldots,n\right\}  $. We shall prove
this claim by descending induction on $k$ (that is, by an induction with base
case $k=n$ and with induction step $k\mapsto k-1$). The \textit{base case}
($k=n$) is vacuously true, since in this case the list of numbers
$k+1,k+2,\ldots,n$ is empty. For the \textit{induction step}, we fix some
$k\in\left[  n\right]  $, and assume (as our induction hypothesis) that
$w_{k}$ fixes $k+1,k+2,\ldots,n$; our goal is to show that $w_{k-1}$ fixes
$k,k+1,\ldots,n$. From (\ref{pf.thm.perm.transp-code.sur.jk=}), we see that
$w_{k}\left(  j_{k}\right)  =k$, and thus $j_{k}\leq k$ (because if we had
$j_{k}>k$, then $j_{k}$ would be one of the numbers $k+1,k+2,\ldots,n$, which
(by our induction hypothesis) would imply that $w_{k}$ fixes $j_{k}$, so that
$w_{k}\left(  j_{k}\right)  =j_{k}>k$, which would contradict $w_{k}\left(
j_{k}\right)  =k$). Hence, both $j_{k}$ and $k$ are $\leq k$. Therefore, the
permutation $t_{j_{k},k}$ fixes $k+1,k+2,\ldots,n$. Hence, its inverse
$t_{j_{k},k}^{-1}$ fixes $k+1,k+2,\ldots,n$ as well. So does the permutation
$w_{k}$ (by our induction hypothesis). But
(\ref{pf.thm.perm.transp-code.sur.wk-1=}) yields $w_{k-1}=w_{k}t_{j_{k}%
,k}^{-1}$. Hence, $w_{k-1}$ fixes $k+1,k+2,\ldots,n$ (since both factors
$w_{k}$ and $t_{j_{k},k}^{-1}$ in the product $w_{k}t_{j_{k},k}^{-1}$ fix
$k+1,k+2,\ldots,n$). Moreover, $w_{k-1}$ fixes $k$ as well, since%
\[
\underbrace{w_{k-1}}_{=w_{k}t_{j_{k},k}^{-1}}\left(  k\right)  =\left(
w_{k}t_{j_{k},k}^{-1}\right)  \left(  k\right)  =w_{k}\left(
\underbrace{t_{j_{k},k}^{-1}\left(  k\right)  }_{=j_{k}}\right)  =w_{k}\left(
j_{k}\right)  =k.
\]
Thus, $w_{k-1}$ fixes $k,k+1,\ldots,n$ (since $w_{k-1}$ fixes $k$ and fixes
$k+1,k+2,\ldots,n$). This completes the induction step. Thus, the proof of
(\ref{pf.thm.perm.transp-code.sur.fixes}) is complete.

In the above induction step, we have additionally shown that each $k\in\left[
n\right]  $ satisfies $j_{k}\leq k$ and thus $j_{k}\in\left[  k\right]  $.
Hence, $\left(  j_{1},j_{2},\ldots,j_{n}\right)  \in\left[  1\right]
\times\left[  2\right]  \times\cdots\times\left[  n\right]  $. This proves one
of our two goals. It remains to prove the other goal, which is $\Phi\left(
j_{1},j_{2},\ldots,j_{n}\right)  =w$.

This is not too hard any more. First, we observe that the permutation $w_{0}$
fixes $1,2,\ldots,n$ (by (\ref{pf.thm.perm.transp-code.sur.fixes}), applied to
$k=0$), and thus equals $\operatorname*{id}$. In other words, $w_{0}%
=\operatorname*{id}$. Next, we note that each $k\in\left[  n\right]  $
satisfies%
\begin{equation}
w_{k}=w_{k-1}t_{j_{k},k}\ \ \ \ \ \ \ \ \ \ \left(  \text{by
(\ref{pf.thm.perm.transp-code.sur.wk-1=})}\right)  .
\label{pf.thm.perm.transp-code.sur.wk=}%
\end{equation}
Now,%
\begin{align*}
w  &  =w_{n}=\underbrace{w_{n-1}}_{\substack{=w_{n-2}t_{j_{n-1},n-1}%
\\\text{(by (\ref{pf.thm.perm.transp-code.sur.wk=}))}}}t_{j_{n},n}%
\ \ \ \ \ \ \ \ \ \ \left(  \text{by (\ref{pf.thm.perm.transp-code.sur.wk=}%
)}\right) \\
&  =\underbrace{w_{n-2}}_{\substack{=w_{n-3}t_{j_{n-2},n-2}\\\text{(by
(\ref{pf.thm.perm.transp-code.sur.wk=}))}}}t_{j_{n-1},n-1}t_{j_{n},n}\\
&  =\underbrace{w_{n-3}}_{\substack{=w_{n-4}t_{j_{n-3},n-3}\\\text{(by
(\ref{pf.thm.perm.transp-code.sur.wk=}))}}}t_{j_{n-2},n-2}t_{j_{n-1}%
,n-1}t_{j_{n},n}\\
&  =\cdots\\
&  =\underbrace{w_{0}}_{=\operatorname*{id}}t_{j_{1},1}t_{j_{2},2}\cdots
t_{j_{n},n}=t_{j_{1},1}t_{j_{2},2}\cdots t_{j_{n},n}=\Phi\left(  j_{1}%
,j_{2},\ldots,j_{n}\right)  .
\end{align*}

Thus, we have written our (arbitrary) permutation $w\in S_{n}$ as $\Phi\left(
j_{1},j_{2},\ldots,j_{n}\right)  $ for some $\left(  j_{1},j_{2},\ldots
,j_{n}\right)  \in\left[  1\right]  \times\left[  2\right]  \times\cdots
\times\left[  n\right]  $. This proves that $\Phi$ is surjective. Thus, our
second proof of Theorem \ref{thm.perm.transp-code} is complete.
\end{proof}

(Note that Theorem \ref{thm.perm.transp-code} is equivalent to \cite[Exercise
5.9]{detnotes}, since $t_{j_{k},k}=t_{k,j_{k}}$ for all $k$ and $j_{k}$.)
\medskip

Theorem \ref{thm.YJM.prod1+mk} now easily follows: Indeed, we have
\[
\left(  1+\mathbf{m}_{1}\right)  \left(  1+\mathbf{m}_{2}\right)
\cdots\left(  1+\mathbf{m}_{n}\right)  =\sum_{\left(  j_{1},j_{2},\ldots
,j_{n}\right)  \in\left[  1\right]  \times\left[  2\right]  \times\cdots
\times\left[  n\right]  }t_{j_{1},1}t_{j_{2},2}\cdots t_{j_{n},n}=\sum_{w\in
S_{n}}w
\]
(here, we have substituted $w$ for $t_{j_{1},1}t_{j_{2},2}\cdots t_{j_{n},n}$
in the sum, since Theorem \ref{thm.perm.transp-code} shows that the map
\begin{align*}
\left[  1\right]  \times\left[  2\right]  \times\cdots\times\left[  n\right]
&  \rightarrow S_{n},\\
\left(  j_{1},j_{2},\ldots,j_{n}\right)   &  \mapsto t_{j_{1},1}t_{j_{2}%
,2}\cdots t_{j_{n},n}%
\end{align*}
is a bijection), so that%
\[
\left(  1+\mathbf{m}_{1}\right)  \left(  1+\mathbf{m}_{2}\right)
\cdots\left(  1+\mathbf{m}_{n}\right)  =\sum_{w\in S_{n}}w=\nabla
\]
(by the definition of $\nabla$). Thus, Theorem \ref{thm.YJM.prod1+mk} is proved.
\end{proof}

We will soon give a different (not in substance, but in presentation) proof of
Theorem \ref{thm.YJM.prod1+mk}, which has the advantage of being slicker. But
the above proof has its own claim to importance: It shows that Theorem
\ref{thm.YJM.prod1+mk} is, at heart, a combinatorial result (Theorem
\ref{thm.perm.transp-code}). Indeed, while we derived Theorem
\ref{thm.YJM.prod1+mk} from Theorem \ref{thm.perm.transp-code}, we could just
as well reverse this argument and obtain Theorem \ref{thm.perm.transp-code}
from Theorem \ref{thm.YJM.prod1+mk} (provided that the latter theorem has been
proved in some other way), since a sum of elements of $\mathbf{k}\left[
S_{n}\right]  $ equals $\nabla$ if and only if this sum contains each
permutation $w\in S_{n}$ exactly once. (To be precise, this reasoning requires
$\mathbf{k}$ to be $\mathbb{Z}$ or a ring that contains $\mathbb{Z}$.) The
takeaway from this is that (almost) every identity that we prove in
$\mathbf{k}\left[  S_{n}\right]  $ can (by expanding its two sides) be turned
into a combinatorial statement about permutations. The former usually looks
slicker, but the latter gives it intuition and grounding.

\subsubsection{More about jucys--murphies}

\begin{exercise}
\label{exe.YJM.sm}\fbox{2} Let $i\in\left[  n-1\right]  $. Show the
following:\medskip

\textbf{(a)} We have $s_{i}\mathbf{m}_{j}=\mathbf{m}_{j}s_{i}$ for all
$j\in\left[  n\right]  \setminus\left\{  i,i+1\right\}  $. \medskip

\textbf{(b)} We have $s_{i}\mathbf{m}_{i}+1=\mathbf{m}_{i+1}s_{i}$. \medskip

\textbf{(c)} We have $s_{i}\mathbf{m}_{i+1}-1=\mathbf{m}_{i}s_{i}$.
\end{exercise}

Soon, we will compute $\mathbf{m}_{2}\mathbf{m}_{3}\cdots\mathbf{m}_{n}$ as
well (see Corollary \ref{cor.YJM.prod}). But let us try it out for $n=3$:%
\[
\mathbf{m}_{2}\mathbf{m}_{3}=t_{1,2}\left(  t_{1,3}+t_{2,3}\right)
=t_{1,2}t_{1,3}+t_{1,2}t_{2,3}=\operatorname*{cyc}\nolimits_{1,2,3}%
+\operatorname*{cyc}\nolimits_{1,3,2}.
\]
This is the sum of the two $3$-cycles in $S_{3}$ (that is, the cycles of the
form $\operatorname*{cyc}\nolimits_{i,j,k}$). Is this the beginning of a
pattern? \medskip

Here is another natural question (inspired by Corollary
\ref{cor.integral.square}, which shows that $P\left(  \nabla\right)  =0$ for
the polynomial $P=x\left(  x-n!\right)  $): What polynomials $P\in
\mathbb{Z}\left[  x\right]  $ satisfy $P\left(  \mathbf{m}_{k}\right)  =0$ for
a given $k$ ? Let us see:

\begin{itemize}
\item For $k=1$, this is easy: We have $\mathbf{m}_{1}=0$, so we can pick
$P=x$.

\item For $k=2$, this is still easy: We have $\mathbf{m}_{2}=t_{1,2}=s_{1}$,
so $\mathbf{m}_{2}^{2}=s_{1}^{2}=\operatorname*{id}=1$ and therefore
$\mathbf{m}_{2}^{2}-1=0$. Thus, we can pick $P=x^{2}-1=\left(  x-1\right)
\left(  x+1\right)  $.

\item For $k=3$, we have $\mathbf{m}_{3}=t_{1,3}+t_{2,3}$, and a bit of
computation shows that
\[
\left(  \mathbf{m}_{3}-2\right)  \left(  \mathbf{m}_{3}-1\right)  \left(
\mathbf{m}_{3}+1\right)  \left(  \mathbf{m}_{3}+2\right)  =0,
\]
so that we can pick $P=\left(  x-2\right)  \left(  x-1\right)  \left(
x+1\right)  \left(  x+2\right)  $.
\end{itemize}

In each of these three cases, we observe that $P$ is a polynomial that splits
over $\mathbb{Z}$ (meaning that it is a product of linear factors $x-r$ for
integers $r$). This generalizes to all $k$ and $n$:

\begin{theorem}
\label{thm.YJM.poly=0}Let $k\in\left[  n\right]  $. Then,%
\[
\prod_{i=-\left(  k-1\right)  }^{k-1}\left(  \mathbf{m}_{k}-i\right)  =0.
\]

\end{theorem}

This theorem is not easy at all. The classical proof uses the so-called
\emph{seminormal basis} of $\mathbf{k}\left[  S_{n}\right]  $. One of the
simplest proofs uses some tricky linear algebra (see Igor Makhlin's
MathOverflow post \cite{Makhlin-yjm}). A long but fully elementary proof can
be found in \cite{yjmnote}, along with further references.

Note that the factor $\mathbf{m}_{k}-0$ in Theorem \ref{thm.YJM.poly=0} is
unnecessary when $k=2$ or $k=3$, as we have seen above. However, all $2k-1$
factors in the product in Theorem \ref{thm.YJM.poly=0} are necessary for
$k\geq4$ (as long as $\mathbf{k}=\mathbb{Z}$).

\subsection{Interlude: Subalgebras generated by a set}

\subsubsection{The general case}

Next, we recall/introduce a general concept in algebra:

\begin{definition}
\label{def.AlgGen.AlgGen}Let $A$ be a $\mathbf{k}$-algebra. Let $U$ be a
subset of $A$. Then, the \emph{subalgebra of }$A$ \emph{generated by }$U$ is
defined to be the subset of $A$ consisting of all (finite) $\mathbf{k}$-linear
combinations of products of elements of $U$. In other words, it is the
$\mathbf{k}$-linear span%
\[
\operatorname*{span}\nolimits_{\mathbf{k}}\left\{  u_{1}u_{2}\cdots
u_{i}\ \mid\ i\in\mathbb{N}\text{ and }u_{1},u_{2},\ldots,u_{i}\in U\right\}
.
\]
(Note that $i$ is allowed to be $0$, in which case $u_{1}u_{2}\cdots
u_{i}=\left(  \text{empty product in }A\right)  =1_{A}$.)

This $\mathbf{k}$-linear span is really a $\mathbf{k}$-subalgebra of $A$ that
contains $U$ as a subset. Moreover, it is the smallest such $\mathbf{k}%
$-subalgebra (i.e., it is a subset of any such $\mathbf{k}$-subalgebra). We
denote it by $\operatorname*{AlgGen}U$ or by $\mathbf{k}\left[  U\right]  $.
\end{definition}

\begin{exercise}
\fbox{1} Prove the above claims.
\end{exercise}

\begin{example}
\textbf{(a)} The subalgebra of the $\mathbb{Z}$-algebra $\mathbb{R}$ generated
by the set $\left\{  \sqrt{2}\right\}  $ is the $\mathbb{Z}$-subalgebra%
\[
\left\{  a+b\sqrt{2}\ \mid\ a,b\in\mathbb{Z}\right\}
\]
of $\mathbb{R}$. Indeed, any product of the form $\sqrt{2}\sqrt{2}\cdots
\sqrt{2}$ is either an integer or $\sqrt{2}$ times an integer, and thus
belongs to the set $\left\{  a+b\sqrt{2}\ \mid\ a,b\in\mathbb{Z}\right\}  $;
therefore, any $\mathbb{Z}$-linear combination of such products belongs to
this set as well. \medskip

\textbf{(b)} The subalgebra of the $\mathbb{Z}$-algebra $\mathbb{R}$ generated
by the set $\left\{  \sqrt[3]{2}\right\}  $ is the $\mathbb{Z}$-subalgebra%
\[
\left\{  a+b\sqrt[3]{2}+c\sqrt[3]{4}\ \mid\ a,b,c\in\mathbb{Z}\right\}
\]
of $\mathbb{R}$. \medskip

\textbf{(c)} More generally, if $A$ is any $\mathbf{k}$-algebra, and if $u$ is
any element of $A$, then the subalgebra of $A$ generated by the one-element
subset $\left\{  u\right\}  $ is the $\mathbf{k}$-linear span of all powers
$u^{0},u^{1},u^{2},\ldots$ of $u$. In other words, it is the set of all
polynomials in $u$ (that is, of all values $p\left(  u\right)  $ for
$p\in\mathbf{k}\left[  x\right]  $). \medskip

\textbf{(d)} Now consider the polynomial ring $\mathbf{k}\left[  x\right]  $
as a $\mathbf{k}$-algebra. The subalgebra of $\mathbf{k}\left[  x\right]  $
generated by the set $\left\{  x^{2},x^{3}\right\}  $ is the set of all
polynomials $p\in\mathbf{k}\left[  x\right]  $ whose $x^{1}$-coefficient is
$0$ (that is, whose derivative at $0$ is $0$). Proving this is a quick way to
check your understanding of polynomials and subalgebras. \medskip

\textbf{(e)} The Laurent polynomial ring $\mathbf{k}\left[  x^{\pm1}\right]  $
is its own subalgebra generated by the two-element set $\left\{
x,x^{-1}\right\}  $. This is just saying that every Laurent polynomial in $x$
can be written as a (usual, not Laurent) polynomial in $x$ and $x^{-1}$.
\end{example}

The following exercise gives another example for subalgebras generated by a subset:

\begin{exercise}
\label{exe.AlgGen.laurent}Consider the Laurent polynomial ring $\mathbf{k}%
\left[  x^{\pm1}\right]  $. \medskip

\textbf{(a)} \fbox{2} Prove that%
\[
\left\{  p\in\mathbf{k}\left[  x^{\pm1}\right]  \ \mid\ p\left(
x^{-1}\right)  =p\left(  x\right)  \right\}  =\operatorname*{AlgGen}\left\{
x+x^{-1}\right\}  .
\]
In other words, the Laurent polynomials $p\in\mathbf{k}\left[  x^{\pm
1}\right]  $ that satisfy $p\left(  x^{-1}\right)  =p\left(  x\right)  $ are
precisely the polynomials in $x+x^{-1}$. \medskip

\textbf{(b)} \fbox{1} Assume that $2$ is a regular element of $\mathbf{k}$
(meaning that $2u=0$ implies $u=0$ whenever $u\in\mathbf{k}$). Prove that%
\[
\left\{  p\in\mathbf{k}\left[  x^{\pm1}\right]  \ \mid\ p\left(  -x\right)
=p\left(  x\right)  \right\}  =\operatorname*{AlgGen}\left\{  x^{2}%
,x^{-2}\right\}
\]
and%
\[
\left\{  p\in\mathbf{k}\left[  x\right]  \ \mid\ p\left(  -x\right)  =p\left(
x\right)  \right\}  =\operatorname*{AlgGen}\left\{  x^{2}\right\}  .
\]

\end{exercise}

The following proposition is easy but quite important, as it provides an easy
way to prove the commutativity of an algebra:

\begin{proposition}
\label{prop.AlgGen.comm}Let $A$ and $U$ be as in Definition
\ref{def.AlgGen.AlgGen}. Assume that every two elements of $U$ commute. Then,
the $\mathbf{k}$-algebra $\operatorname*{AlgGen}U$ is commutative.
\end{proposition}

\begin{exercise}
\label{exe.AlgGen.comm}\fbox{2} Prove this.
\end{exercise}

\subsubsection{The Gelfand--Tsetlin algebra}

Thanks to Proposition \ref{prop.AlgGen.comm}, the commutativity of an algebra
needs only to be checked on its generators. This lets us easily show facts
like the following:

\begin{corollary}
\label{cor.GZ.comm}Let $\operatorname*{GZ}\nolimits_{n}$ be the subalgebra of
$\mathbf{k}\left[  S_{n}\right]  $ generated by the set $\left\{
\mathbf{m}_{1},\mathbf{m}_{2},\ldots,\mathbf{m}_{n}\right\}  $ of its
jucys--murphies. (In other words, let $\operatorname*{GZ}\nolimits_{n}%
:=\operatorname*{AlgGen}\left\{  \mathbf{m}_{1},\mathbf{m}_{2},\ldots
,\mathbf{m}_{n}\right\}  $.) Then, $\operatorname*{GZ}\nolimits_{n}$ is a
commutative $\mathbf{k}$-algebra.
\end{corollary}

\begin{proof}
Theorem \ref{thm.YJM.commute} shows that every two elements of $\left\{
\mathbf{m}_{1},\mathbf{m}_{2},\ldots,\mathbf{m}_{n}\right\}  $ commute. Thus,
Proposition \ref{prop.AlgGen.comm} (applied to $A=\mathbf{k}\left[
S_{n}\right]  $ and $U=\left\{  \mathbf{m}_{1},\mathbf{m}_{2},\ldots
,\mathbf{m}_{n}\right\}  $) shows that $\operatorname*{AlgGen}\left\{
\mathbf{m}_{1},\mathbf{m}_{2},\ldots,\mathbf{m}_{n}\right\}  $ is commutative.
This proves Corollary \ref{cor.GZ.comm}.
\end{proof}

The subalgebra $\operatorname*{GZ}\nolimits_{n}$ defined in Corollary
\ref{cor.GZ.comm} is called the \emph{Gelfand--Tsetlin subalgebra} (or
\emph{Gelfand--Cetlin subalgebra}, or \emph{Gelfand--Zetlin subalgebra}) of
$\mathbf{k}\left[  S_{n}\right]  $, and we will eventually learn some more of
its properties.

\subsection{The somewhere-to-below shuffles}

Here is a newer family of elements of $\mathbf{k}\left[  S_{n}\right]  $,
introduced by Lafreni\`{e}re and Grinberg in \cite{GriLaf22}:

\begin{definition}
\label{def.stb.tk}The \emph{somewhere-to-below shuffles} $\mathbf{t}%
_{1},\mathbf{t}_{2},\ldots,\mathbf{t}_{n}$ are the elements of $\mathbf{k}%
\left[  S_{n}\right]  $ defined by%
\[
\mathbf{t}_{k}:=\operatorname*{cyc}\nolimits_{k}+\operatorname*{cyc}%
\nolimits_{k,k+1}+\operatorname*{cyc}\nolimits_{k,k+1,k+2}+\cdots
+\operatorname*{cyc}\nolimits_{k,k+1,\ldots,n}=\sum_{i=k}^{n}%
\operatorname*{cyc}\nolimits_{k,k+1,\ldots,i}\in\mathbf{k}\left[
S_{n}\right]  .
\]

(Recall that $\operatorname*{cyc}\nolimits_{k}=\operatorname*{id}$.)
\end{definition}

In particular, $\mathbf{t}_{n}=\operatorname*{cyc}\nolimits_{n}%
=\operatorname*{id}=1$. The motivation for calling these elements
\textquotedblleft shuffles\textquotedblright\ will eventually be explained.

The following exercise generalizes (\ref{eq.exa.sga.S3.f2}) (why?):

\begin{exercise}
\label{exe.stb.prod=Nabla}\fbox{2} Prove that $\mathbf{t}_{n}\mathbf{t}%
_{n-1}\cdots\mathbf{t}_{1}=\nabla$.
\end{exercise}

The somewhere-to-below shuffles $\mathbf{t}_{1},\mathbf{t}_{2},\ldots
,\mathbf{t}_{n}$ do not commute for $n\geq3$. For instance, for $n=3$, we have%
\[
\mathbf{t}_{1}\mathbf{t}_{2}-\mathbf{t}_{2}\mathbf{t}_{1}=t_{1,2}%
+\operatorname*{cyc}\nolimits_{1,2,3}-\operatorname*{cyc}\nolimits_{1,3,2}%
-t_{1,3}\neq0
\]
(if $\mathbf{k}$ is nontrivial). However, they come fairly close to commuting.
To make sense of this, let us define the \emph{commutator} in a ring:

\begin{definition}
\label{def.ring.commutator}Let $a$ and $b$ be two elements of a ring. Then,
their \emph{commutator} $\left[  a,b\right]  $ is defined to be $ab-ba$.
\end{definition}

Thus, two elements $a$ and $b$ satisfy $\left[  a,b\right]  =0$ if and only if
$a$ and $b$ commute.

Our above computation shows that $\left[  \mathbf{t}_{1},\mathbf{t}%
_{2}\right]  \neq0$ (for $n\geq3$). However, one can also show that $\left[
\mathbf{t}_{1},\mathbf{t}_{2}\right]  ^{2}=0$. Thus, $\left[  \mathbf{t}%
_{1},\mathbf{t}_{2}\right]  $ is a nilpotent element of $\mathbf{k}\left[
S_{n}\right]  $, with nilpotency order\footnote{The \emph{nilpotency order}
(or \emph{order of nilpotency}) of a nilpotent element $a$ of a ring $A$ is
defined to be the smallest $k\in\mathbb{N}$ such that $a^{k}=0$.} $2$. For
comparison, a strictly upper-triangular $k\times k$-matrix is nilpotent with
nilpotency order $\leq k$.

This nilpotency of $\left[  \mathbf{t}_{1},\mathbf{t}_{2}\right]  $
furthermore generalizes: We have%
\[
\left[  \mathbf{t}_{i},\mathbf{t}_{j}\right]  ^{j-i+1}%
=0\ \ \ \ \ \ \ \ \ \ \text{for all }1\leq i<j\leq n.
\]
This is \cite[Corollary 9.11]{Grinbe23}.

\begin{exercise}
\label{exe.stb.ti+1ti=}\fbox{3} Prove that $\mathbf{t}_{i+1}\mathbf{t}%
_{i}=\left(  \mathbf{t}_{i}-1\right)  \mathbf{t}_{i}$ for all $i\in\left[
n-1\right]  $.
\end{exercise}

We can again ask for polynomials $P\in\mathbf{k}\left[  x\right]  $ that
annihilate these $\mathbf{t}_{i}$'s. For example, it turns out (\cite[Theorem
3]{r2t-elem}) that%
\[
P\left(  \mathbf{t}_{1}\right)  =0\ \ \ \ \ \ \ \ \ \ \text{for }P=x\left(
x-1\right)  \left(  x-2\right)  \cdots\left(  x-\left(  n-2\right)  \right)
\left(  x-n\right)  .
\]
More generally,%
\[
P\left(  \mathbf{t}_{i}\right)  =0\ \ \ \ \ \ \ \ \ \ \text{for }P=x\left(
x-1\right)  \left(  x-2\right)  \cdots\left(  x-\left(  n-i-1\right)  \right)
\left(  x-\left(  n-i+1\right)  \right)  .
\]
We might prove this eventually.

\subsection{Partial integrals}

\subsubsection{Definition and properties}

Let us now generalize the integral $\nabla$ and the sign-integral $\nabla^{-}$:

\begin{definition}
\label{def.intX.intX}Let $X$ be a subset of $\left[  n\right]  $. Then, we
define two elements $\nabla_{X}$ and $\nabla_{X}^{-}$ of $\mathbf{k}\left[
S_{n}\right]  $ by%
\[
\nabla_{X}:=\sum_{\substack{w\in S_{n};\\w\left(  i\right)  =i\text{ for all
}i\notin X}}w\ \ \ \ \ \ \ \ \ \ \text{and}\ \ \ \ \ \ \ \ \ \ \nabla_{X}%
^{-}:=\sum_{\substack{w\in S_{n};\\w\left(  i\right)  =i\text{ for all
}i\notin X}}\left(  -1\right)  ^{w}w.
\]
(Here, \textquotedblleft for all $i\notin X$\textquotedblright\ means
\textquotedblleft for all $i\in\left[  n\right]  \setminus X$%
\textquotedblright. Thus, both sums range over all permutations $w\in S_{n}$
that fix all elements $i\in\left[  n\right]  $ that don't belong to $X$.)

We shall call $\nabla_{X}$ the $X$\emph{-integral}, and we shall call
$\nabla_{X}^{-}$ the $X$\emph{-sign-integral}.
\end{definition}

\begin{example}
\label{exa.intX.X=0123}Let $X$ be a subset of $\left[  n\right]  $. \medskip

\textbf{(a)} If $\left\vert X\right\vert =0$, then $X=\varnothing$ and
therefore $\nabla_{X}=\nabla_{X}^{-}=1$ (because if a permutation $w\in S_{n}$
satisfies $w\left(  i\right)  =i$ for all $i\notin\varnothing$, then it simply
satisfies $w\left(  i\right)  =i$ for all $i\in\left[  n\right]  $, which
means that $w=\operatorname*{id}=1$). \medskip

\textbf{(b)} If $\left\vert X\right\vert =1$, then we again have $\nabla
_{X}=\nabla_{X}^{-}=1$ (because if a permutation $w\in S_{n}$ fixes all but
one element of $\left[  n\right]  $, then it must also fix the remaining
element, and thus it equals $\operatorname*{id}=1$). \medskip

\textbf{(c)} If $\left\vert X\right\vert =2$, then $X=\left\{  i,j\right\}  $
for some $i\neq j$, and thus%
\[
\nabla_{X}=1+t_{i,j}\ \ \ \ \ \ \ \ \ \ \text{and}\ \ \ \ \ \ \ \ \ \ \nabla
_{X}^{-}=1-t_{i,j}%
\]
(since the permutations $w\in S_{n}$ that fix all elements of $\left[
n\right]  $ except for $i$ and $j$ are precisely $\operatorname*{id}=1$ and
$t_{i,j}$). \medskip

\textbf{(d)} If $\left\vert X\right\vert =3$, then $X=\left\{  i,j,k\right\}
$ for some distinct $i,j,k$, and thus%
\begin{align*}
\nabla_{X}  &  =1+t_{i,j}+t_{i,k}+t_{j,k}+\operatorname*{cyc}\nolimits_{i,j,k}%
+\operatorname*{cyc}\nolimits_{i,k,j}\ \ \ \ \ \ \ \ \ \ \text{and}\\
\nabla_{X}^{-}  &  =1-t_{i,j}-t_{i,k}-t_{j,k}+\operatorname*{cyc}%
\nolimits_{i,j,k}+\operatorname*{cyc}\nolimits_{i,k,j}.
\end{align*}
(Here, we have used the fact that any transposition $t_{p,q}$ has sign $-1$
whereas any cycle $\operatorname*{cyc}\nolimits_{p,q,r}$ has sign $1$.)
\end{example}

\begin{example}
\label{exa.intX.n}The $\left[  n\right]  $-integral $\nabla_{\left[  n\right]
}$ just equals the usual integral $\nabla$. Indeed, the definition of
$\nabla_{\left[  n\right]  }$ yields%
\[
\nabla_{\left[  n\right]  }=\sum_{\substack{w\in S_{n};\\w\left(  i\right)
=i\text{ for all }i\notin\left[  n\right]  }}w=\sum_{w\in S_{n}}w
\]
(since the condition \textquotedblleft$w\left(  i\right)  =i$ for all
$i\notin\left[  n\right]  $\textquotedblright\ is vacuously true for any $w\in
S_{n}$), so that
\[
\nabla_{\left[  n\right]  }=\sum_{w\in S_{n}}w=\nabla.
\]
Similarly, $\nabla_{\left[  n\right]  }^{-}=\nabla^{-}$.
\end{example}

The following proposition collects a number of basic properties of $\nabla
_{X}$:

\begin{proposition}
\label{prop.intX.basics}Let $X$ be a subset of $\left[  n\right]  $. Let%
\[
S_{n,X}:=\left\{  w\in S_{n}\ \mid\ w\left(  i\right)  =i\text{ for all
}i\notin X\right\}  .
\]
Then: \medskip

\textbf{(a)} The set $S_{n,X}$ is a subgroup of $S_{n}$, and is isomorphic to
$S_{X}$ via the group isomorphism%
\begin{align*}
S_{n,X}  &  \rightarrow S_{X},\\
w  &  \mapsto\left.  w\mid_{X}^{X}\right.  .
\end{align*}
(Here, $w\mid_{X}^{X}$ means the restriction of $w$ to $X$, regarded as a map
from $X$ to $X$.) \medskip

\textbf{(b)} We have
\[
\nabla_{X}=\sum\limits_{w\in S_{n,X}}w\ \ \ \ \ \ \ \ \ \ \text{and}%
\ \ \ \ \ \ \ \ \ \ \nabla_{X}^{-}=\sum\limits_{w\in S_{n,X}}\left(
-1\right)  ^{w}w.
\]

\textbf{(c)} For each $w\in S_{n,X}$, we have%
\begin{align*}
w\nabla_{X}  &  =\nabla_{X}w=\nabla_{X}\ \ \ \ \ \ \ \ \ \ \text{and}\\
w\nabla_{X}^{-}  &  =\nabla_{X}^{-}w=\left(  -1\right)  ^{w}\nabla_{X}^{-}.
\end{align*}

\textbf{(d)} Let $i$ and $j$ be two distinct elements of $X$. Let us set%
\[
\nabla_{X}^{\operatorname*{even}}:=\sum_{\substack{w\in S_{n,X};\\\left(
-1\right)  ^{w}=1}}w.
\]
Then,%
\begin{equation}
\nabla_{X}=\nabla_{X}^{\operatorname*{even}}\cdot\left(  1+t_{i,j}\right)
=\left(  1+t_{i,j}\right)  \cdot\nabla_{X}^{\operatorname*{even}}
\label{prop.intX.basics.d.1}%
\end{equation}
and%
\begin{equation}
\nabla_{X}^{-}=\nabla_{X}^{\operatorname*{even}}\cdot\left(  1-t_{i,j}\right)
=\left(  1-t_{i,j}\right)  \cdot\nabla_{X}^{\operatorname*{even}}.
\label{prop.intX.basics.d.2}%
\end{equation}

\end{proposition}

\begin{example}
Let $X$ be a three-element subset $\left\{  i,j,k\right\}  $ of $\left[
n\right]  $. Then, the element $\nabla_{X}^{\operatorname*{even}}$ in
Proposition \ref{prop.intX.basics} \textbf{(d)} is given by%
\[
\nabla_{X}^{\operatorname*{even}}=1+\operatorname*{cyc}\nolimits_{i,j,k}%
+\operatorname*{cyc}\nolimits_{i,k,j}%
\]
(since the even permutations $w\in S_{n,X}$ are $\operatorname*{id}=1$ and
$\operatorname*{cyc}\nolimits_{i,j,k}$ and $\operatorname*{cyc}%
\nolimits_{i,k,j}$). Hence, for instance, the first equality in
(\ref{prop.intX.basics.d.1}) says that%
\[
\underbrace{1+t_{i,j}+t_{i,k}+t_{j,k}+\operatorname*{cyc}\nolimits_{i,j,k}%
+\operatorname*{cyc}\nolimits_{i,k,j}}_{=\nabla_{X}}=\underbrace{\left(
1+\operatorname*{cyc}\nolimits_{i,j,k}+\operatorname*{cyc}\nolimits_{i,k,j}%
\right)  }_{=\nabla_{X}^{\operatorname*{even}}}\cdot\left(  1+t_{i,j}\right)
.
\]

\end{example}

\begin{proof}
[Proof of Proposition \ref{prop.intX.basics}.]\textbf{(a)} We must prove the
following claims:

\begin{statement}
\textit{Claim 1:} For any $w\in S_{n,X}$, we have $w\left(  X\right)
\subseteq X$ (so that the restriction $w\mid_{X}^{X}$ is really a well-defined
map from $X$ to $X$) and $\left.  w\mid_{X}^{X}\right.  \in S_{X}$.
\end{statement}

\begin{statement}
\textit{Claim 2:} The set $S_{n,X}$ is a subgroup of $S_{n}$.
\end{statement}

\begin{statement}
\textit{Claim 3:} The map
\begin{align*}
S_{n,X}  &  \rightarrow S_{X},\\
w  &  \mapsto\left.  w\mid_{X}^{X}\right.
\end{align*}
is a group morphism.
\end{statement}

\begin{statement}
\textit{Claim 4:} This group morphism is a group isomorphism.
\end{statement}

All these claims are easy, so we just outline their proofs:

\begin{fineprint}
\begin{proof}
[Proof of Claim 1 (sketched).]Let $w\in S_{n,X}$. Let $x\in X$. From $w\in
S_{n,X}$, we obtain $w\left(  i\right)  =i$ for all $i\notin X$. Thus, if we
had $w\left(  x\right)  \notin X$, then we would have $w\left(  w\left(
x\right)  \right)  =w\left(  x\right)  $ (by applying the equality $w\left(
i\right)  =i$ to $i=w\left(  x\right)  $), which would entail $w\left(
x\right)  =x$ (since $w$ is a permutation and thus injective), and therefore
$x=w\left(  x\right)  \notin X$, which would contradict $x\in X$. Hence, we
cannot have $w\left(  x\right)  \notin X$. Thus, $w\left(  x\right)  \in X$.

Forget that we fixed $x$. We thus have shown that $w\left(  x\right)  \in X$
for each $x\in X$. That is, $w\left(  X\right)  \subseteq X$.

This yields that the restriction $w\mid_{X}^{X}$ is well-defined as a map from
$X$ to $X$. It remains to show that it is a permutation of $X$ (that is,
belongs to $S_{X}$). It clearly suffices to show that it is injective and
surjective. The injectivity follows from the injectivity of $w$ (which is
because $w$ is a permutation), since a restriction of an injective map is
again injective. Now, let us prove the surjectivity: Let $y\in X$. Let
$x=w^{-1}\left(  y\right)  $. Thus, $w\left(  x\right)  =y$. From $w\in
S_{n,X}$, we obtain $w\left(  i\right)  =i$ for all $i\notin X$. Thus, if we
had $x\notin X$, then we would have $w\left(  x\right)  =x$ and thus
$x=w\left(  x\right)  =y\in X$, which would contradict $x\notin X$. Hence, we
cannot have $x\notin X$. So we must have $x\in X$. Thus, $\left(  w\mid
_{X}^{X}\right)  \left(  x\right)  =w\left(  x\right)  =y$. This shows that
$y$ is an image under the map $w\mid_{X}^{X}$. Since we have proved this for
each $y\in X$, we have thus shown that each element of $X$ is an image under
$w\mid_{X}^{X}$. In other words, $w\mid_{X}^{X}$ is surjective. This completes
our proof of Claim 1.
\end{proof}

\begin{proof}
[Proof of Claim 2 (sketched).]This is straightforward. For example, if $u,v\in
S_{n,X}$, then $uv\in S_{n,X}$, since each $i\notin X$ satisfies $\left(
uv\right)  \left(  i\right)  =u\left(  \underbrace{v\left(  i\right)  }%
_{=i}\right)  =u\left(  i\right)  =i$.
\end{proof}

\begin{proof}
[Proof of Claim 3 (sketched).]This is again straightforward. (The
well-definedness of the map follows from Claim 1.)
\end{proof}

\begin{proof}
[Proof of Claim 4 (sketched).]Let us denote this group morphism by $\Phi$.

Each permutation $w\in S_{n,X}$ is uniquely determined by its values on the
elements of $X$ (since all its other values are uniquely determined by the
requirement \textquotedblleft$w\left(  i\right)  =i$ for all $i\notin
X$\textquotedblright), that is, by its restriction $\Phi\left(  w\right)
=\left.  w\mid_{X}^{X}\right.  $. Hence, the map $\Phi$ is injective.

Any permutation $u\in S_{X}$ can be extended to a permutation $w\in S_{n,X}$
satisfying $\Phi\left(  w\right)  =u$ (we just need to set $w\left(  x\right)
=u\left(  x\right)  $ for all $x\in X$, and $w\left(  i\right)  =i$ for all
$i\notin X$). Thus, the map $\Phi$ is surjective.

So the map $\Phi$ is both injective and surjective, hence bijective. Since
$\Phi$ is a group morphism, it thus follows that $\Phi$ is a group
isomorphism. This proves Claim 4.
\end{proof}
\end{fineprint}

With these four claims proved, Proposition \ref{prop.intX.basics} \textbf{(a)}
follows. \medskip

\textbf{(b)} This is just a restatement of the definitions of $\nabla_{X}$ and
$\nabla_{X}^{-}$, since the permutations $w\in S_{n}$ that satisfy $w\left(
i\right)  =i$ for all $i\notin X$ are precisely the permutations $w\in
S_{n,X}$. \medskip

\textbf{(c)} This is a generalization of Proposition \ref{prop.integral.fix},
and is proved in the exact same way (using the formulas for $\nabla_{X}$ and
$\nabla_{X}^{-}$ given in Proposition \ref{prop.intX.basics} \textbf{(b)}),
since $S_{n,X}$ is a group. \medskip

\textbf{(d)} Let $i$ and $j$ be two distinct elements of $X$. Then,
Proposition \ref{prop.intX.basics} \textbf{(b)} yields
\begin{equation}
\nabla_{X}=\sum\limits_{w\in S_{n,X}}w=\sum\limits_{\substack{w\in
S_{n,X};\\\left(  -1\right)  ^{w}=1}}w+\sum\limits_{\substack{w\in
S_{n,X};\\\left(  -1\right)  ^{w}=-1}}w \label{pf.prop.intX.basics.d.1}%
\end{equation}
(since each $w\in S_{n,X}$ satisfies either $\left(  -1\right)  ^{w}=1$ or
$\left(  -1\right)  ^{w}=-1$, but not both at the same time). Let us now
simplify the second sum on the right hand side. Indeed, $S_{n,X}$ is a group,
and $t_{i,j}$ is one of its elements (since both $i$ and $j$ belong to $X$,
and thus the transposition $t_{i,j}$ fixes all elements $k\notin X$). Hence,
the map%
\begin{align*}
S_{n,X}  &  \rightarrow S_{n,X},\\
w  &  \mapsto wt_{i,j}%
\end{align*}
is well-defined. This map is furthermore a bijection (with inverse given by
$w\mapsto wt_{i,j}^{-1}$). Thus, we can substitute $wt_{i,j}$ for $w$ in the
sum $\sum\limits_{\substack{w\in S_{n,X};\\\left(  -1\right)  ^{w}=-1}}w$. We
thus obtain%
\begin{equation}
\sum\limits_{\substack{w\in S_{n,X};\\\left(  -1\right)  ^{w}=-1}%
}w=\sum\limits_{\substack{w\in S_{n,X};\\\left(  -1\right)  ^{wt_{i,j}}%
=-1}}wt_{i,j}. \label{pf.prop.intX.basics.d.2}%
\end{equation}
However, the transposition $t_{i,j}$ has sign $\left(  -1\right)  ^{t_{i,j}%
}=-1$ (like any transposition), and thus each $w\in S_{n,X}$ satisfies%
\begin{align*}
\left(  -1\right)  ^{wt_{i,j}}  &  =\left(  -1\right)  ^{w}\underbrace{\left(
-1\right)  ^{t_{i,j}}}_{=-1}\ \ \ \ \ \ \ \ \ \ \left(  \text{by the rule
}\left(  -1\right)  ^{\sigma\tau}=\left(  -1\right)  ^{\sigma}\cdot\left(
-1\right)  ^{\tau}\right) \\
&  =-\left(  -1\right)  ^{w}.
\end{align*}
Hence, the condition \textquotedblleft$\left(  -1\right)  ^{wt_{i,j}}%
=-1$\textquotedblright\ under the summation sign on the right hand side of
(\ref{pf.prop.intX.basics.d.2}) is equivalent to the condition
\textquotedblleft$-\left(  -1\right)  ^{w}=-1$\textquotedblright, which is in
turn equivalent to \textquotedblleft$\left(  -1\right)  ^{w}=1$%
\textquotedblright. Therefore, we can replace the former condition by the
latter. Thus, (\ref{pf.prop.intX.basics.d.2}) rewrites as%
\[
\sum\limits_{\substack{w\in S_{n,X};\\\left(  -1\right)  ^{w}=-1}%
}w=\sum\limits_{\substack{w\in S_{n,X};\\\left(  -1\right)  ^{w}=1}%
}wt_{i,j}=\underbrace{\left(  \sum\limits_{\substack{w\in S_{n,X};\\\left(
-1\right)  ^{w}=1}}w\right)  }_{=\nabla_{X}^{\operatorname*{even}}}%
t_{i,j}=\nabla_{X}^{\operatorname*{even}}t_{i,j}.
\]
Hence, (\ref{pf.prop.intX.basics.d.1}) becomes%
\[
\nabla_{X}=\underbrace{\sum\limits_{\substack{w\in S_{n,X};\\\left(
-1\right)  ^{w}=1}}w}_{=\nabla_{X}^{\operatorname*{even}}}+\underbrace{\sum
\limits_{\substack{w\in S_{n,X};\\\left(  -1\right)  ^{w}=-1}}w}_{=\nabla
_{X}^{\operatorname*{even}}t_{i,j}}=\nabla_{X}^{\operatorname*{even}}%
+\nabla_{X}^{\operatorname*{even}}t_{i,j}=\nabla_{X}^{\operatorname*{even}%
}\cdot\left(  1+t_{i,j}\right)  .
\]
A similar argument (but using the bijection $w\mapsto t_{i,j}w$ instead of
$w\mapsto wt_{i,j}$) shows that $\nabla_{X}=\left(  1+t_{i,j}\right)
\cdot\nabla_{X}^{\operatorname*{even}}$. Combining these two equalities, we
obtain (\ref{prop.intX.basics.d.1}). The proof of (\ref{prop.intX.basics.d.2})
is similar, except that the signs need to be taken into account. Proposition
\ref{prop.intX.basics} \textbf{(d)} is thus completely proved.
\end{proof}

We note that the $t_{i,j}$ in Proposition \ref{prop.intX.basics} \textbf{(d)}
could be replaced by any odd permutation $\tau\in S_{n,X}$:

\begin{proposition}
\label{prop.intX.basics.d-gen}Let $X$ be a subset of $\left[  n\right]  $. Let
$\tau\in S_{n,X}$ be an odd permutation (where $S_{n,X}$ is as in Proposition
\ref{prop.intX.basics}). Let us set%
\[
\nabla_{X}^{\operatorname*{even}}:=\sum_{\substack{w\in S_{n,X};\\\left(
-1\right)  ^{w}=1}}w.
\]
Then,%
\begin{equation}
\nabla_{X}=\nabla_{X}^{\operatorname*{even}}\cdot\left(  1+\tau\right)
=\left(  1+\tau\right)  \cdot\nabla_{X}^{\operatorname*{even}}
\label{eq.prop.intX.basics.d-gen.1}%
\end{equation}
and%
\begin{equation}
\nabla_{X}^{-}=\nabla_{X}^{\operatorname*{even}}\cdot\left(  1-\tau\right)
=\left(  1-\tau\right)  \cdot\nabla_{X}^{\operatorname*{even}}.
\label{eq.prop.intX.basics.d-gen.2}%
\end{equation}

\end{proposition}

\begin{proof}
Exactly the same as for Proposition \ref{prop.intX.basics} \textbf{(d)}.
\end{proof}

\begin{remark}
\label{rmk.intX.birep}The equality (\ref{eq.prop.intX.basics.d-gen.1}) can be
generalized even further to arbitrary finite groups, albeit at the cost of
concreteness. Indeed, let $G$ be a finite group, and let $H$ be a subgroup of
$G$. Then, we can find a subset $P$ of $G$ such that each element of $G$ can
be uniquely expressed as $hp$ with $h\in H$ and $p\in P$. (Such a subset $P$
is called a \emph{right transversal of }$H$ \emph{in }$G$, or a \emph{system
of representatives for the right cosets of }$H$\emph{ in }$G$.) For such a
subset $P$, we then get%
\[
\sum_{g\in G}g=\left(  \sum_{h\in H}h\right)  \left(  \sum_{p\in P}p\right)
.
\]
For $G=S_{n,X}$ and $H=S_{n,X}^{\operatorname*{even}}=\left\{  w\in
S_{n,X}\ \mid\ \left(  -1\right)  ^{w}=1\right\}  $ and $P=\left\{
1,\tau\right\}  $, this recovers the first equality in
(\ref{eq.prop.intX.basics.d-gen.1}).

We can also find a subset $Q$ of $G$ such that each element of $G$ can be
uniquely expressed as $qh$ with $h\in H$ and $q\in Q$. Then we have%
\[
\sum_{g\in G}g=\left(  \sum_{q\in Q}q\right)  \left(  \sum_{h\in H}h\right)
.
\]

More interestingly, we can find a subset $R$ of $G$ that simultaneously
satisfies the requirements on $P$ and on $Q$ (by a result of Miller from 1910
-- see \url{https://math.stackexchange.com/questions/178186} for some
references), so that we get%
\[
\sum_{g\in G}g=\left(  \sum_{h\in H}h\right)  \left(  \sum_{r\in R}r\right)
=\left(  \sum_{r\in R}r\right)  \left(  \sum_{h\in H}h\right)  .
\]
In the case of $G=S_{n,X}$ and $H=S_{n,X}^{\operatorname*{even}}$, we can take
$R=\left\{  1,\tau\right\}  $ for example.
\end{remark}

\subsubsection{A recursion}

We can also use $X$-integrals to give a new proof of Theorem
\ref{thm.YJM.prod1+mk}. For this, we will show the following
lemma:\footnote{Recall that $\left[  i\right]  =\left\{  1,2,\ldots,i\right\}
$ for each $i\in\mathbb{N}$. Also, the Young--Jucys--Murphy elements
$\mathbf{m}_{1},\mathbf{m}_{2},\ldots,\mathbf{m}_{n}$ were defined in
Definition \ref{def.YJM.mk}.}

\begin{lemma}
\label{lem.intX.rec1}For any $k\in\left[  n\right]  $, we have $\nabla
_{\left[  k\right]  }=\nabla_{\left[  k-1\right]  }\left(  1+\mathbf{m}%
_{k}\right)  $.
\end{lemma}

If this lemma is proved, then we get%
\begin{align*}
\nabla &  =\nabla_{\left[  n\right]  }\ \ \ \ \ \ \ \ \ \ \left(  \text{by
Example \ref{exa.intX.n}}\right) \\
&  =\underbrace{\nabla_{\left[  n-1\right]  }}_{\substack{=\nabla_{\left[
n-2\right]  }\left(  1+\mathbf{m}_{n-1}\right)  \\\text{(by Lemma
\ref{lem.intX.rec1})}}}\left(  1+\mathbf{m}_{n}\right)
\ \ \ \ \ \ \ \ \ \ \left(  \text{by Lemma }\ref{lem.intX.rec1}\right) \\
&  =\underbrace{\nabla_{\left[  n-2\right]  }}_{\substack{=\nabla_{\left[
n-3\right]  }\left(  1+\mathbf{m}_{n-2}\right)  \\\text{(by Lemma
\ref{lem.intX.rec1})}}}\left(  1+\mathbf{m}_{n-1}\right)  \left(
1+\mathbf{m}_{n}\right) \\
&  =\cdots\\
&  =\underbrace{\nabla_{\left[  0\right]  }}_{\substack{=1\\\text{(by Example
\ref{exa.intX.X=0123} \textbf{(a)})}}}\left(  1+\mathbf{m}_{1}\right)  \left(
1+\mathbf{m}_{2}\right)  \cdots\left(  1+\mathbf{m}_{n}\right) \\
&  =\left(  1+\mathbf{m}_{1}\right)  \left(  1+\mathbf{m}_{2}\right)
\cdots\left(  1+\mathbf{m}_{n}\right)  ,
\end{align*}
which proves Theorem \ref{thm.YJM.prod1+mk} again. So we only need to prove
Lemma \ref{lem.intX.rec1} in order to obtain a second proof of Theorem
\ref{thm.YJM.prod1+mk}.

We shall not prove Lemma \ref{lem.intX.rec1} immediately, but rather extend it
first, as this allows us to get a more general result at the same price.
First, we extend it by replacing the subset $\left[  k\right]  $ and its
element $k$ with an arbitrary subset $X$ and any element $x\in X$:

\begin{lemma}
\label{lem.intX.rec2}Let $X$ be any subset of $\left[  n\right]  $, and let
$x\in X$. Then,%
\[
\nabla_{X}=\nabla_{X\setminus\left\{  x\right\}  }\left(  1+\sum_{y\in
X\setminus\left\{  x\right\}  }t_{y,x}\right)  .
\]

\end{lemma}

Even more generally, we can replace each transposition $t_{y,x}$ on the right
hand side by an arbitrary permutation $\sigma_{y}\in S_{n,X}$ that sends $y$
to $x$ (where $S_{n,X}$ is as in Proposition \ref{prop.intX.basics}), and also
replace the $1$ by an arbitrary permutation $\sigma_{x}\in S_{n,X}$ that sends
$x$ to $x$. Thus we arrive at the following generalization:

\begin{lemma}
\label{lem.intX.rec3}Let $X$ be any subset of $\left[  n\right]  $, and let
$x\in X$. Let $S_{n,X}$ be as in Proposition \ref{prop.intX.basics}. For each
$y\in X$, fix a permutation $\sigma_{y}\in S_{n,X}$ that sends $y$ to $x$.
Then,%
\[
\nabla_{X}=\nabla_{X\setminus\left\{  x\right\}  }\sum_{y\in X}\sigma_{y}.
\]

\end{lemma}

\begin{proof}
[Proof of Lemma \ref{lem.intX.rec3}.]We first observe that%
\begin{equation}
\left\{  w\in S_{n,X}\ \mid\ w\left(  x\right)  =x\right\}  =S_{n,X\setminus
\left\{  x\right\}  }. \label{pf.lem.intX.rec3.sub}%
\end{equation}

\begin{proof}
[Proof of (\ref{pf.lem.intX.rec3.sub}).]We have $x\in X$ and thus $\left\{
x\right\}  \subseteq X$.

Recall that $S_{n,X}=\left\{  w\in S_{n}\ \mid\ w\left(  i\right)  =i\text{
for all }i\notin X\right\}  $ (by the definition of $S_{n,X}$). Thus,%
\begin{align*}
&  \left\{  w\in S_{n,X}\ \mid\ w\left(  x\right)  =x\right\} \\
&  =\left\{  w\in S_{n}\ \mid\ w\left(  i\right)  =i\text{ for all }i\notin
X\text{, and moreover }w\left(  x\right)  =x\right\} \\
&  =\left\{  w\in S_{n}\ \mid\ w\left(  i\right)  =i\text{ for all }i\notin
X\text{ and also for }i=x\right\} \\
&  =\left\{  w\in S_{n}\ \mid\ w\left(  i\right)  =i\text{ for all }%
i\in\left[  n\right]  \setminus X\text{ and also for }i=x\right\} \\
&  =\left\{  w\in S_{n}\ \mid\ w\left(  i\right)  =i\text{ for all }%
i\in\left(  \left[  n\right]  \setminus X\right)  \cup\left\{  x\right\}
\right\} \\
&  =\left\{  w\in S_{n}\ \mid\ w\left(  i\right)  =i\text{ for all }%
i\in\left[  n\right]  \setminus\left(  X\setminus\left\{  x\right\}  \right)
\right\} \\
&  \ \ \ \ \ \ \ \ \ \ \ \ \ \ \ \ \ \ \ \ \left(
\begin{array}
[c]{c}%
\text{since }\left(  \left[  n\right]  \setminus X\right)  \cup\left\{
x\right\}  =\left[  n\right]  \setminus\left(  X\setminus\left\{  x\right\}
\right) \\
\text{(because }\left\{  x\right\}  \subseteq X\subseteq\left[  n\right]
\text{)}%
\end{array}
\right) \\
&  =\left\{  w\in S_{n}\ \mid\ w\left(  i\right)  =i\text{ for all }i\notin
X\setminus\left\{  x\right\}  \right\} \\
&  =S_{n,X\setminus\left\{  x\right\}  }\ \ \ \ \ \ \ \ \ \ \left(  \text{by
the definition of }S_{n,X\setminus\left\{  x\right\}  }\right)  .
\end{align*}
This proves (\ref{pf.lem.intX.rec3.sub}).
\end{proof}

Proposition \ref{prop.intX.basics} \textbf{(b)} yields%
\begin{align}
\nabla_{X}  &  =\sum\limits_{w\in S_{n,X}}w\ \ \ \ \ \ \ \ \ \ \text{and}%
\label{pf.lem.intX.rec3.1}\\
\nabla_{X\setminus\left\{  x\right\}  }  &  =\sum\limits_{w\in S_{n,X\setminus
\left\{  x\right\}  }}w. \label{pf.lem.intX.rec3.1x}%
\end{align}

If $w\in S_{n,X}$, then there is a unique $y\in X$ such that $w\left(
y\right)  =x$\ \ \ \ \footnote{\textit{Proof.} Let $w\in S_{n,X}$. We must
prove that there is a unique $y\in X$ such that $w\left(  y\right)  =x$.
\par
From $w\in S_{n,X}$, we see that $w\left(  i\right)  =i$ for all $i\notin X$.
Thus, if we had $w^{-1}\left(  x\right)  \notin X$, then we would have
$w\left(  w^{-1}\left(  x\right)  \right)  =w^{-1}\left(  x\right)  $ (by the
equality $w\left(  i\right)  =i$, applied to $i=w^{-1}\left(  x\right)  $), so
that $w^{-1}\left(  x\right)  =w\left(  w^{-1}\left(  x\right)  \right)  =x\in
X$, which would contradict $w^{-1}\left(  x\right)  \notin X$. Thus, we cannot
have $w^{-1}\left(  x\right)  \notin X$. Hence, we must have $w^{-1}\left(
x\right)  \in X$. Thus, there exists a $y\in X$ such that $w\left(  y\right)
=x$ (namely, $y=w^{-1}\left(  x\right)  $). Moreover, this $y$ is clearly
unique (since $w$ is a permutation and thus injective). Hence, we have shown
that there is a unique $y\in X$ such that $w\left(  y\right)  =x$.}. Hence, we
can split the sum $\sum\limits_{w\in S_{n,X}}w$ according to the value of this
$y$. We thus obtain%
\begin{equation}
\sum\limits_{w\in S_{n,X}}w=\sum_{y\in X}\ \ \sum\limits_{\substack{w\in
S_{n,X};\\w\left(  y\right)  =x}}w. \label{pf.lem.intX.rec3.2}%
\end{equation}

Let us next show that each $y\in X$ satisfies%
\begin{equation}
\sum\limits_{\substack{w\in S_{n,X};\\w\left(  y\right)  =x}}w=\nabla
_{X\setminus\left\{  x\right\}  }\sigma_{y}. \label{pf.lem.intX.rec3.3}%
\end{equation}

\begin{proof}
[Proof of (\ref{pf.lem.intX.rec3.3}).]Let $y\in X$. Then, $\sigma_{y}\in
S_{n,X}$ is a permutation that sends $y$ to $x$ (by its definition); thus,
$\sigma_{y}\left(  y\right)  =x$. Moreover, $S_{n,X}$ is a group (by
Proposition \ref{prop.intX.basics} \textbf{(a)}). Thus, $w\sigma_{y}\in
S_{n,X}$ for each $w\in S_{n,X}$ (since $\sigma_{y}\in S_{n,X}$). Hence, the
map
\begin{align*}
S_{n,X}  &  \rightarrow S_{n,X},\\
w  &  \mapsto w\sigma_{y}%
\end{align*}
is well-defined. This map is furthermore bijective (with inverse given by
$w\mapsto w\sigma_{y}^{-1}$). Thus, we can substitute $w\sigma_{y}$ for $w$ in
the sum $\sum\limits_{\substack{w\in S_{n,X};\\w\left(  y\right)  =x}}w$. We
thus find%
\begin{equation}
\sum\limits_{\substack{w\in S_{n,X};\\w\left(  y\right)  =x}}w=\sum
\limits_{\substack{w\in S_{n,X};\\\left(  w\sigma_{y}\right)  \left(
y\right)  =x}}w\sigma_{y}=\sum\limits_{\substack{w\in S_{n,X};\\w\left(
x\right)  =x}}w\sigma_{y} \label{pf.lem.intX.rec3.3.pf.1}%
\end{equation}
(since each $w\in S_{n,X}$ satisfies $\left(  w\sigma_{y}\right)  \left(
y\right)  =w\left(  \underbrace{\sigma_{y}\left(  y\right)  }_{=x}\right)
=w\left(  x\right)  $). But we have%
\[
\left\{  w\in S_{n,X}\ \mid\ w\left(  x\right)  =x\right\}  =S_{n,X\setminus
\left\{  x\right\}  }\ \ \ \ \ \ \ \ \ \ \left(  \text{by
(\ref{pf.lem.intX.rec3.sub})}\right)  .
\]
Thus, the summation sign $\sum\limits_{\substack{w\in S_{n,X};\\w\left(
x\right)  =x}}$ on the right hand side of (\ref{pf.lem.intX.rec3.3.pf.1}) can
be simplified to $\sum\limits_{w\in S_{n,X\setminus\left\{  x\right\}  }}$.
Hence, (\ref{pf.lem.intX.rec3.3.pf.1}) rewrites as%
\[
\sum\limits_{\substack{w\in S_{n,X};\\w\left(  y\right)  =x}}w=\sum
\limits_{w\in S_{n,X\setminus\left\{  x\right\}  }}w\sigma_{y}%
=\underbrace{\left(  \sum\limits_{w\in S_{n,X\setminus\left\{  x\right\}  }%
}w\right)  }_{\substack{=\nabla_{X\setminus\left\{  x\right\}  }\\\text{(by
(\ref{pf.lem.intX.rec3.1x}))}}}\sigma_{y}=\nabla_{X\setminus\left\{
x\right\}  }\sigma_{y}.
\]
This proves (\ref{pf.lem.intX.rec3.3}).
\end{proof}

Now, (\ref{pf.lem.intX.rec3.1}) becomes%
\begin{align*}
\nabla_{X}  &  =\sum\limits_{w\in S_{n,X}}w=\sum_{y\in X}\ \ \underbrace{\sum
\limits_{\substack{w\in S_{n,X};\\w\left(  y\right)  =x}}w}_{\substack{=\nabla
_{X\setminus\left\{  x\right\}  }\sigma_{y}\\\text{(by
(\ref{pf.lem.intX.rec3.3}))}}}\ \ \ \ \ \ \ \ \ \ \left(  \text{by
(\ref{pf.lem.intX.rec3.2})}\right) \\
&  =\sum_{y\in X}\nabla_{X\setminus\left\{  x\right\}  }\sigma_{y}%
=\nabla_{X\setminus\left\{  x\right\}  }\sum_{y\in X}\sigma_{y}.
\end{align*}
This proves Lemma \ref{lem.intX.rec3}.
\end{proof}

\begin{proof}
[Proof of Lemma \ref{lem.intX.rec2}.]Let $S_{n,X}$ be as in Proposition
\ref{prop.intX.basics}. Set $t_{x,x}:=\operatorname*{id}\in S_{n}$. Thus, a
permutation $t_{y,x}$ is defined for all $y\in X$ (including $x$). It is clear
that this permutation $t_{y,x}$ belongs to $S_{n,X}$ (since both $x$ and $y$
belong to $X$) and sends $y$ to $x$. Thus, Lemma \ref{lem.intX.rec3} (applied
to $\sigma_{y}=t_{y,x}$) yields
\begin{align*}
\nabla_{X}  &  =\nabla_{X\setminus\left\{  x\right\}  }\underbrace{\sum_{y\in
X}t_{y,x}}_{\substack{=t_{x,x}+\sum_{y\in X\setminus\left\{  x\right\}
}t_{y,x}\\\text{(here, we have split off the addend for }y=x\text{)}}}\\
&  =\nabla_{X\setminus\left\{  x\right\}  }\left(  \underbrace{t_{x,x}%
}_{=\operatorname*{id}=1}+\sum_{y\in X\setminus\left\{  x\right\}  }%
t_{y,x}\right)  =\nabla_{X\setminus\left\{  x\right\}  }\left(  1+\sum_{y\in
X\setminus\left\{  x\right\}  }t_{y,x}\right)  .
\end{align*}
This proves Lemma \ref{lem.intX.rec2}.
\end{proof}

\begin{proof}
[Proof of Lemma \ref{lem.intX.rec1}.]Lemma \ref{lem.intX.rec2} (applied to
$X=\left[  k\right]  $ and $x=k$) yields
\begin{align*}
\nabla_{\left[  k\right]  }  &  =\nabla_{\left[  k\right]  \setminus\left\{
k\right\}  }\left(  1+\sum_{y\in\left[  k\right]  \setminus\left\{  k\right\}
}t_{y,k}\right) \\
&  =\nabla_{\left[  k-1\right]  }\left(  1+\underbrace{\sum_{y\in\left[
k-1\right]  }t_{y,k}}_{=t_{1,k}+t_{2,k}+\cdots+t_{k-1,k}}\right)
\ \ \ \ \ \ \ \ \ \ \left(  \text{since }\left[  k\right]  \setminus\left\{
k\right\}  =\left[  k-1\right]  \right) \\
&  =\nabla_{\left[  k-1\right]  }\left(  1+\underbrace{t_{1,k}+t_{2,k}%
+\cdots+t_{k-1,k}}_{\substack{=\mathbf{m}_{k}\\\text{(by the definition of
}\mathbf{m}_{k}\text{)}}}\right)  =\nabla_{\left[  k-1\right]  }\left(
1+\mathbf{m}_{k}\right)  .
\end{align*}
This proves Lemma \ref{lem.intX.rec1}.
\end{proof}

\subsubsection{The transposition cancellation theorem}

The following theorem generalizes (\ref{eq.exa.sga.S3.f1}) as well as the
slightly more general formula $\left(  1+t_{i,j}\right)  \left(
1-t_{i,j}\right)  =0$ for all $i\neq j$:

\begin{theorem}
[transposition cancellation theorem]\label{thm.int.+-=0}Let $X$ and $Y$ be two
subsets of $\left[  n\right]  $ such that $\left\vert X\cap Y\right\vert >1$.
Then,%
\begin{align*}
\nabla_{X}\nabla_{Y}^{-}  &  =0\ \ \ \ \ \ \ \ \ \ \text{and}\\
\nabla_{X}^{-}\nabla_{Y}  &  =0.
\end{align*}

\end{theorem}

For example, $\nabla_{\left\{  1,2,3\right\}  }\nabla_{\left\{  1,2,4\right\}
}^{-}=0$ and $\nabla_{\left\{  1,2,3\right\}  }^{-}\nabla_{\left\{
1,2,4\right\}  }=0$.

\begin{proof}
[Proof of Theorem \ref{thm.int.+-=0}.]There exist two distinct elements $i$
and $j$ of $X\cap Y$ (since $\left\vert X\cap Y\right\vert >1$). Pick such $i$
and $j$. Then, $i$ and $j$ are two distinct elements of $X$ (since $i\in X\cap
Y\subseteq X$ and $j\in X\cap Y\subseteq X$). Hence,
(\ref{prop.intX.basics.d.1}) yields
\[
\nabla_{X}=\nabla_{X}^{\operatorname*{even}}\cdot\left(  1+t_{i,j}\right)
=\left(  1+t_{i,j}\right)  \cdot\nabla_{X}^{\operatorname*{even}},
\]
where $\nabla_{X}^{\operatorname*{even}}$ is defined in Proposition
\ref{prop.intX.basics} \textbf{(d)}. However, $i$ and $j$ are also two
distinct elements of $Y$ (since $i\in X\cap Y\subseteq Y$ and $j\in X\cap
Y\subseteq Y$). Hence, (\ref{prop.intX.basics.d.2}) yields%
\[
\nabla_{Y}^{-}=\nabla_{Y}^{\operatorname*{even}}\cdot\left(  1-t_{i,j}\right)
=\left(  1-t_{i,j}\right)  \cdot\nabla_{Y}^{\operatorname*{even}},
\]
where $\nabla_{Y}^{\operatorname*{even}}$ is defined similarly to the
$\nabla_{X}^{\operatorname*{even}}$ in Proposition \ref{prop.intX.basics}
\textbf{(d)}. Thus,%
\[
\underbrace{\nabla_{X}}_{=\nabla_{X}^{\operatorname*{even}}\cdot\left(
1+t_{i,j}\right)  }\ \ \underbrace{\nabla_{Y}^{-}}_{=\left(  1-t_{i,j}\right)
\cdot\nabla_{Y}^{\operatorname*{even}}}=\nabla_{X}^{\operatorname*{even}}%
\cdot\underbrace{\left(  1+t_{i,j}\right)  \cdot\left(  1-t_{i,j}\right)
}_{\substack{=1-t_{i,j}^{2}=0\\\text{(since }t_{i,j}^{2}=\operatorname*{id}%
=1\text{)}}}\cdot\,\nabla_{Y}^{\operatorname*{even}}=0.
\]
Thus, $\nabla_{X}\nabla_{Y}^{-}=0$ is proved. Similarly, we can show
$\nabla_{X}^{-}\nabla_{Y}=0$ (using $\nabla_{X}^{-}=\nabla_{X}%
^{\operatorname*{even}}\cdot\left(  1-t_{i,j}\right)  $ and $\nabla
_{Y}=\left(  1+t_{i,j}\right)  \cdot\nabla_{Y}^{\operatorname*{even}}$). The
proof of Theorem \ref{thm.int.+-=0} is thus complete.
\end{proof}

\subsubsection{Disjoint commutativity}

The following easy fact says that if $X$ and $Y$ are two disjoint subsets of
$\left[  n\right]  $, then the $X$-(sign)-integral commutes with the $Y$-(sign-)integral:

\begin{proposition}
\label{prop.int.commute}Let $X$ and $Y$ be two disjoint subsets of $\left[
n\right]  $. Then: \medskip

\textbf{(a)} If $u\in S_{n}$ is a permutation that satisfies $\left(  u\left(
i\right)  =i\text{ for all }i\notin X\right)  $, and if $v\in S_{n}$ is a
permutation that satisfies $\left(  v\left(  i\right)  =i\text{ for all
}i\notin Y\right)  $, then $uv=vu$. \medskip

\textbf{(b)} We have
\begin{align*}
\nabla_{X}\nabla_{Y}  &  =\nabla_{Y}\nabla_{X};\ \ \ \ \ \ \ \ \ \ \nabla
_{X}\nabla_{Y}^{-}=\nabla_{Y}^{-}\nabla_{X};\\
\nabla_{X}^{-}\nabla_{Y}  &  =\nabla_{Y}\nabla_{X}^{-}%
;\ \ \ \ \ \ \ \ \ \ \nabla_{X}^{-}\nabla_{Y}^{-}=\nabla_{Y}^{-}\nabla_{X}%
^{-}.
\end{align*}

\textbf{(c)} If $v\in S_{n}$ is a permutation that satisfies $\left(  v\left(
i\right)  =i\text{ for all }i\notin Y\right)  $, then%
\[
\nabla_{X}v=v\nabla_{X}\ \ \ \ \ \ \ \ \ \ \text{and}%
\ \ \ \ \ \ \ \ \ \ \nabla_{X}^{-}v=v\nabla_{X}^{-}.
\]

\textbf{(d)} If $p$ and $q$ are two distinct elements of $Y$, then the
transposition $t_{p,q}\in S_{n}$ satisfies%
\[
\nabla_{X}t_{p,q}=t_{p,q}\nabla_{X}\ \ \ \ \ \ \ \ \ \ \text{and}%
\ \ \ \ \ \ \ \ \ \ \nabla_{X}^{-}t_{p,q}=t_{p,q}\nabla_{X}^{-}.
\]

\end{proposition}

\begin{fineprint}
\begin{proof}
\textbf{(a)} Let $u\in S_{n}$ be a permutation that satisfies $\left(
u\left(  i\right)  =i\text{ for all }i\notin X\right)  $. Let $v\in S_{n}$ be
a permutation that satisfies $\left(  v\left(  i\right)  =i\text{ for all
}i\notin Y\right)  $. We must prove that $uv=vu$. In other words, we must
prove that $\left(  uv\right)  \left(  p\right)  =\left(  vu\right)  \left(
p\right)  $ for each $p\in\left[  n\right]  $.

So let $p\in\left[  n\right]  $ be arbitrary. We are in one of the following
three cases:

\textit{Case 1:} We have $p\in X$.

\textit{Case 2:} We have $p\in Y$.

\textit{Case 3:} We have neither $p\in X$ nor $p\in Y$.

Let us first consider Case 1. In this case, we have $p\in X$. Hence, $p\notin
Y$ (since the sets $X$ and $Y$ are disjoint). Thus, $v\left(  p\right)  =p$
(since $v\left(  i\right)  =i$ for all $i\notin Y$).

If we had $u\left(  p\right)  \notin X$, then we would have $u\left(  u\left(
p\right)  \right)  =u\left(  p\right)  $ (since $u\left(  i\right)  =i$ for
all $i\notin X$), which would yield $u\left(  p\right)  =p$ (since the map $u$
is a permutation and thus injective), which would in turn yield $p=u\left(
p\right)  \notin X$, which would contradict $p\in X$. Hence, we cannot have
$u\left(  p\right)  \notin X$. Thus, we have $u\left(  p\right)  \in X$.
Therefore, $u\left(  p\right)  \notin Y$ (since the sets $X$ and $Y$ are
disjoint). Thus, $v\left(  u\left(  p\right)  \right)  =u\left(  p\right)  $
(since $v\left(  i\right)  =i$ for all $i\notin Y$).

Now, $\left(  uv\right)  \left(  p\right)  =u\left(  v\left(  p\right)
\right)  =u\left(  p\right)  $ (since $v\left(  p\right)  =p$). Comparing this
with $\left(  vu\right)  \left(  p\right)  =v\left(  u\left(  p\right)
\right)  =u\left(  p\right)  $, we obtain $\left(  uv\right)  \left(
p\right)  =\left(  vu\right)  \left(  p\right)  $. Hence, we have proved
$\left(  uv\right)  \left(  p\right)  =\left(  vu\right)  \left(  p\right)  $
in Case 1.

Let us now consider Case 2. In this case, we have $p\in Y$. Hence, the same
argument as in Case 1 (but now with the roles of $X$ and $Y$ swapped, and with
the roles of $u$ and $v$ swapped) shows that $\left(  vu\right)  \left(
p\right)  =\left(  uv\right)  \left(  p\right)  $. In other words, $\left(
uv\right)  \left(  p\right)  =\left(  vu\right)  \left(  p\right)  $. Hence,
we have proved $\left(  uv\right)  \left(  p\right)  =\left(  vu\right)
\left(  p\right)  $ in Case 2.

Let us finally consider Case 3. In this case, we have neither $p\in X$ nor
$p\in Y$. Hence, $p\notin X$ and $p\notin Y$. From $p\notin X$, we obtain
$u\left(  p\right)  =p$ (since $u\left(  i\right)  =i$ for all $i\notin X$).
Thus, $v\left(  u\left(  p\right)  \right)  =v\left(  p\right)  $. From
$p\notin Y$, we obtain $v\left(  p\right)  =p$ (since $v\left(  i\right)  =i$
for all $i\notin Y$). Hence, $u\left(  v\left(  p\right)  \right)  =u\left(
p\right)  $. Now, comparing $\left(  uv\right)  \left(  p\right)  =u\left(
v\left(  p\right)  \right)  =u\left(  p\right)  =p$ with $\left(  vu\right)
\left(  p\right)  =v\left(  u\left(  p\right)  \right)  =v\left(  p\right)
=p$, we obtain $\left(  uv\right)  \left(  p\right)  =\left(  vu\right)
\left(  p\right)  $. Hence, we have proved $\left(  uv\right)  \left(
p\right)  =\left(  vu\right)  \left(  p\right)  $ in Case 3.

Now we have proved $\left(  uv\right)  \left(  p\right)  =\left(  vu\right)
\left(  p\right)  $ in all three Cases 1, 2 and 3. Hence, $\left(  uv\right)
\left(  p\right)  =\left(  vu\right)  \left(  p\right)  $ always holds.

Forget that we fixed $p$. We thus have shown that $\left(  uv\right)  \left(
p\right)  =\left(  vu\right)  \left(  p\right)  $ for all $p\in\left[
n\right]  $. In other words, $uv=vu$. This proves Proposition
\ref{prop.int.commute} \textbf{(a)}. \medskip

\textbf{(b)} The definition of $\nabla_{X}$ yields
\begin{equation}
\nabla_{X}=\sum_{\substack{w\in S_{n};\\w\left(  i\right)  =i\text{ for all
}i\notin X}}w=\sum_{\substack{u\in S_{n};\\u\left(  i\right)  =i\text{ for all
}i\notin X}}u. \label{pf.prop.int.commute.b.1}%
\end{equation}
Similarly,%
\begin{equation}
\nabla_{Y}=\sum_{\substack{v\in S_{n};\\v\left(  i\right)  =i\text{ for all
}i\notin Y}}v. \label{pf.prop.int.commute.b.2}%
\end{equation}
Multiplying these two equalities, we obtain%
\begin{align*}
\nabla_{X}\nabla_{Y}  &  =\left(  \sum_{\substack{u\in S_{n};\\u\left(
i\right)  =i\text{ for all }i\notin X}}u\right)  \left(  \sum_{\substack{v\in
S_{n};\\v\left(  i\right)  =i\text{ for all }i\notin Y}}v\right) \\
&  =\sum_{\substack{u\in S_{n};\\u\left(  i\right)  =i\text{ for all }i\notin
X}}\ \ \sum_{\substack{v\in S_{n};\\v\left(  i\right)  =i\text{ for all
}i\notin Y}}\ \ \underbrace{uv}_{\substack{=vu\\\text{(by Proposition
\ref{prop.int.commute} \textbf{(a)})}}}\\
&  =\sum_{\substack{u\in S_{n};\\u\left(  i\right)  =i\text{ for all }i\notin
X}}\ \ \sum_{\substack{v\in S_{n};\\v\left(  i\right)  =i\text{ for all
}i\notin Y}}vu\\
&  =\underbrace{\left(  \sum_{\substack{v\in S_{n};\\v\left(  i\right)
=i\text{ for all }i\notin Y}}v\right)  }_{\substack{=\nabla_{Y}\\\text{(by
(\ref{pf.prop.int.commute.b.2}))}}}\underbrace{\left(  \sum_{\substack{u\in
S_{n};\\u\left(  i\right)  =i\text{ for all }i\notin X}}u\right)
}_{\substack{=\nabla_{X}\\\text{(by (\ref{pf.prop.int.commute.b.1}))}}%
}=\nabla_{Y}\nabla_{X}.
\end{align*}
Analogous computations (with some signs injected in) prove the equalities
$\nabla_{X}\nabla_{Y}^{-}=\nabla_{Y}^{-}\nabla_{X}$ and $\nabla_{X}^{-}%
\nabla_{Y}=\nabla_{Y}\nabla_{X}^{-}$ and $\nabla_{X}^{-}\nabla_{Y}^{-}%
=\nabla_{Y}^{-}\nabla_{X}^{-}$. This completes the proof of Proposition
\ref{prop.int.commute} \textbf{(b)}. \medskip

\textbf{(c)} Let $v\in S_{n}$ be a permutation that satisfies $\left(
v\left(  i\right)  =i\text{ for all }i\notin Y\right)  $. Then, from
(\ref{pf.prop.int.commute.b.1}), we obtain%
\begin{align*}
\nabla_{X}v  &  =\left(  \sum_{\substack{u\in S_{n};\\u\left(  i\right)
=i\text{ for all }i\notin X}}u\right)  v=\sum_{\substack{u\in S_{n};\\u\left(
i\right)  =i\text{ for all }i\notin X}}\ \ \underbrace{uv}%
_{\substack{=vu\\\text{(by Proposition \ref{prop.int.commute} \textbf{(a)})}%
}}\\
&  =\sum_{\substack{u\in S_{n};\\u\left(  i\right)  =i\text{ for all }i\notin
X}}\ \ vu=v\underbrace{\left(  \sum_{\substack{u\in S_{n};\\u\left(  i\right)
=i\text{ for all }i\notin X}}u\right)  }_{\substack{=\nabla_{X}\\\text{(by
(\ref{pf.prop.int.commute.b.1}))}}}=v\nabla_{X}.
\end{align*}
An analogous computation (with some signs injected in) proves the equality
$\nabla_{X}^{-}v=v\nabla_{X}^{-}$. Thus, Proposition \ref{prop.int.commute}
\textbf{(c)} is proved. \medskip

\textbf{(d)} Let $p$ and $q$ be two distinct elements of $Y$. Then, $\left\{
p,q\right\}  \subseteq Y$. Hence, each $i\notin Y$ satisfies $i\notin\left\{
p,q\right\}  $.

But the transposition $t_{p,q}$ satisfies $t_{p,q}\left(  i\right)  =i$ for
all $i\notin\left\{  p,q\right\}  $ (by its definition). Therefore, it
satisfies $t_{p,q}\left(  i\right)  =i$ for all $i\notin Y$ (since each
$i\notin Y$ satisfies $i\notin\left\{  p,q\right\}  $). Hence, Proposition
\ref{prop.int.commute} \textbf{(c)} (applied to $v=t_{p,q}$) yields
\[
\nabla_{X}t_{p,q}=t_{p,q}\nabla_{X}\ \ \ \ \ \ \ \ \ \ \text{and}%
\ \ \ \ \ \ \ \ \ \ \nabla_{X}^{-}t_{p,q}=t_{p,q}\nabla_{X}^{-}.
\]
This proves Proposition \ref{prop.int.commute} \textbf{(d)}.
\end{proof}
\end{fineprint}

\subsubsection{Products of $X$-integrals}

The next proposition will prove itself useful in the study of Specht modules.

\begin{proposition}
\label{prop.int.wXi=Xi}Let $X_{1},X_{2},\ldots,X_{k}$ be $k$ disjoint subsets
of $\left[  n\right]  $ such that $X_{1}\cup X_{2}\cup\cdots\cup X_{k}=\left[
n\right]  $. Then: \medskip

\textbf{(a)} If $w\in S_{n}$ is any permutation, then the four statements%
\begin{align}
&  \left(  w\left(  X_{i}\right)  \subseteq X_{i}\text{ for all }i\in\left[
k\right]  \right)  ,\label{eq.prop.int.wXi=Xi.a.1}\\
&  \left(  w\left(  X_{i}\right)  =X_{i}\text{ for all }i\in\left[  k\right]
\right)  ,\label{eq.prop.int.wXi=Xi.a.2}\\
&  \left(  w\left(  X_{i}\right)  \subseteq X_{i}\text{ for all }i\in\left[
k-1\right]  \right)  \ \ \ \ \ \ \ \ \ \ \text{and}%
\label{eq.prop.int.wXi=Xi.a.3}\\
&  \left(  w\left(  X_{i}\right)  =X_{i}\text{ for all }i\in\left[
k-1\right]  \right)  \label{eq.prop.int.wXi=Xi.a.4}%
\end{align}
are equivalent. \medskip

\textbf{(b)} We have%
\begin{equation}
\sum_{\substack{w\in S_{n};\\w\left(  X_{i}\right)  \subseteq X_{i}\text{ for
all }i\in\left[  k\right]  }}w=\nabla_{X_{1}}\nabla_{X_{2}}\cdots\nabla
_{X_{k}} \label{eq.prop.int.wXi=Xi.b.1}%
\end{equation}
and%
\begin{equation}
\sum_{\substack{w\in S_{n};\\w\left(  X_{i}\right)  \subseteq X_{i}\text{ for
all }i\in\left[  k\right]  }}\left(  -1\right)  ^{w}w=\nabla_{X_{1}}^{-}%
\nabla_{X_{2}}^{-}\cdots\nabla_{X_{k}}^{-}. \label{eq.prop.int.wXi=Xi.b.2}%
\end{equation}

\textbf{(c)} The factors $\nabla_{X_{1}},\nabla_{X_{2}},\ldots,\nabla_{X_{k}}$
on the right hand side of (\ref{eq.prop.int.wXi=Xi.b.1}) pairwise commute. So
do the factors $\nabla_{X_{1}}^{-},\nabla_{X_{2}}^{-},\ldots,\nabla_{X_{k}%
}^{-}$ on the right hand side of (\ref{eq.prop.int.wXi=Xi.b.2}).
\end{proposition}

\begin{example}
Let $n=8$ and $k=4$ and $X_{1}=\left\{  2,5,6\right\}  $ and $X_{2}=\left\{
1,3\right\}  $ and $X_{3}=\left\{  4,8\right\}  $ and $X_{4}=\left\{
7\right\}  $ (so that the four sets $X_{1},X_{2},X_{3},X_{4}$ are disjoint and
their union is $\left[  n\right]  $). Then, (\ref{eq.prop.int.wXi=Xi.b.1})
says that%
\[
\sum_{\substack{w\in S_{8};\\w\left(  \left\{  2,5,6\right\}  \right)
\subseteq\left\{  2,5,6\right\}  ;\\w\left(  \left\{  1,3\right\}  \right)
\subseteq\left\{  1,3\right\}  ;\\w\left(  \left\{  4,8\right\}  \right)
\subseteq\left\{  4,8\right\}  ;\\w\left(  \left\{  7\right\}  \right)
\subseteq\left\{  7\right\}  }}w=\nabla_{\left\{  2,5,6\right\}  }%
\nabla_{\left\{  1,3\right\}  }\nabla_{\left\{  4,8\right\}  }\nabla_{\left\{
7\right\}  }.
\]
Note that the factor $\nabla_{\left\{  7\right\}  }$ equals $1$ (by Example
\ref{exa.intX.X=0123} \textbf{(b)}) and thus can be omitted.
\end{example}

\begin{proof}
[Proof of Proposition \ref{prop.int.wXi=Xi}.]\textbf{(a)} Let $w\in S_{n}$ be
any permutation. We must prove that the four statements
(\ref{eq.prop.int.wXi=Xi.a.1}), (\ref{eq.prop.int.wXi=Xi.a.2}),
(\ref{eq.prop.int.wXi=Xi.a.3}) and (\ref{eq.prop.int.wXi=Xi.a.4}) are
equivalent. Clearly, (\ref{eq.prop.int.wXi=Xi.a.2}) implies
(\ref{eq.prop.int.wXi=Xi.a.1}). Let us now show the converse implication:

\begin{proof}
[Proof that (\ref{eq.prop.int.wXi=Xi.a.1}) implies
(\ref{eq.prop.int.wXi=Xi.a.2}):]Assume that (\ref{eq.prop.int.wXi=Xi.a.1})
holds. We must prove that (\ref{eq.prop.int.wXi=Xi.a.2}) holds.

Let $i\in\left[  k\right]  $. Then, $w\left(  X_{i}\right)  \subseteq X_{i}$
(since we assumed that (\ref{eq.prop.int.wXi=Xi.a.1}) holds). But $w$ is a
permutation, and thus is injective. Hence, $\left\vert w\left(  P\right)
\right\vert =\left\vert P\right\vert $ for any subset $P$ of $\left[
n\right]  $. In particular, this shows that $\left\vert w\left(  X_{i}\right)
\right\vert =\left\vert X_{i}\right\vert $.

But $X_{i}$ is a finite set. Hence, the only subset of $X_{i}$ that has the
same size as $X_{i}$ is $X_{i}$ itself. In other words, if $P$ is a subset of
$X_{i}$ such that $\left\vert P\right\vert =\left\vert X_{i}\right\vert $,
then $P=X_{i}$. Applying this to $P=w\left(  X_{i}\right)  $, we obtain
$w\left(  X_{i}\right)  =X_{i}$ (since $w\left(  X_{i}\right)  \subseteq
X_{i}$ and $\left\vert w\left(  X_{i}\right)  \right\vert =\left\vert
X_{i}\right\vert $).

Forget that we fixed $i$. We thus have shown that $w\left(  X_{i}\right)
=X_{i}$ for all $i\in\left[  k\right]  $. In other words,
(\ref{eq.prop.int.wXi=Xi.a.2}) holds. This shows that
(\ref{eq.prop.int.wXi=Xi.a.1}) implies (\ref{eq.prop.int.wXi=Xi.a.2}).
\end{proof}

We have now shown that (\ref{eq.prop.int.wXi=Xi.a.1}) implies
(\ref{eq.prop.int.wXi=Xi.a.2}), whereas we know that
(\ref{eq.prop.int.wXi=Xi.a.2}) implies (\ref{eq.prop.int.wXi=Xi.a.1}). Hence,
the statements (\ref{eq.prop.int.wXi=Xi.a.1}) and
(\ref{eq.prop.int.wXi=Xi.a.2}) are equivalent. The same argument (with $k-1$
instead of $k$) shows that the statements (\ref{eq.prop.int.wXi=Xi.a.3}) and
(\ref{eq.prop.int.wXi=Xi.a.4}) are equivalent.

It remains to prove that the statements (\ref{eq.prop.int.wXi=Xi.a.1}) and
(\ref{eq.prop.int.wXi=Xi.a.3}) are equivalent. Clearly,
(\ref{eq.prop.int.wXi=Xi.a.1}) implies (\ref{eq.prop.int.wXi=Xi.a.3}), so we
only need to check that (\ref{eq.prop.int.wXi=Xi.a.3}) implies
(\ref{eq.prop.int.wXi=Xi.a.1}).

\begin{proof}
[Proof that (\ref{eq.prop.int.wXi=Xi.a.3}) implies
(\ref{eq.prop.int.wXi=Xi.a.1}):]Assume that (\ref{eq.prop.int.wXi=Xi.a.3})
holds. We must prove that (\ref{eq.prop.int.wXi=Xi.a.1}) holds. In other
words, we must prove that $w\left(  X_{i}\right)  \subseteq X_{i}$ for all
$i\in\left[  k\right]  $.

First let us prove that $w\left(  X_{k}\right)  \subseteq X_{k}$. Indeed, let
$x\in X_{k}$ be arbitrary. We shall show that $w\left(  x\right)  \in X_{k}$.
Indeed, assume the contrary. Then, $w\left(  x\right)  \notin X_{k}$. But
$w\left(  x\right)  \in\left[  n\right]  =X_{1}\cup X_{2}\cup\cdots\cup X_{k}%
$, so that $w\left(  x\right)  \in X_{i}$ for some $i\in\left[  k\right]  $.
Consider this $i$. If we had $i=k$, then we would obtain $w\left(  x\right)
\in X_{k}$ (since $w\left(  x\right)  \in X_{i}$), which would contradict
$w\left(  x\right)  \notin X_{k}$. Thus, we cannot have $i=k$. Hence, $i\neq
k$, so that $i\in\left[  k\right]  \setminus\left\{  k\right\}  =\left[
k-1\right]  $.

But we have assumed that (\ref{eq.prop.int.wXi=Xi.a.3}) holds. Thus,
(\ref{eq.prop.int.wXi=Xi.a.4}) also holds (since the statements
(\ref{eq.prop.int.wXi=Xi.a.3}) and (\ref{eq.prop.int.wXi=Xi.a.4}) are
equivalent). Therefore, $w\left(  X_{i}\right)  =X_{i}$ (since $i\in\left[
k-1\right]  $). Hence, $w\left(  x\right)  \in X_{i}=w\left(  X_{i}\right)  $,
so that $w\left(  x\right)  =w\left(  y\right)  $ for some $y\in X_{i}$.
Consider this $y$.

However, $w$ is a permutation and thus injective. Hence, from $w\left(
x\right)  =w\left(  y\right)  $, we obtain $x=y$, so that $x=y\in X_{i}$.
Combining this with $x\in X_{k}$, we obtain $x\in X_{i}\cap X_{k}$.

But the sets $X_{1},X_{2},\ldots,X_{k}$ are disjoint. Hence, $X_{i}\cap
X_{k}=\varnothing$ (since $i\neq k$). Hence, $x\in X_{i}\cap X_{k}%
=\varnothing$, which is absurd. This contradiction shows that our assumption
was false. Thus, $w\left(  x\right)  \in X_{k}$ is proved.

Forget that we fixed $x$. We thus have shown that $w\left(  x\right)  \in
X_{k}$ for each $x\in X_{k}$. In other words, $w\left(  X_{k}\right)
\subseteq X_{k}$.

Now, we know that the relation $w\left(  X_{i}\right)  \subseteq X_{i}$ holds
for all $i\in\left[  k-1\right]  $ (by our assumption
(\ref{eq.prop.int.wXi=Xi.a.3})), but we also know that the same relation holds
for $i=k$ (since $w\left(  X_{k}\right)  \subseteq X_{k}$). Hence, this
relation holds for all $i\in\left[  k-1\right]  \cup\left\{  k\right\}  $,
that is, for all $i\in\left[  k\right]  $. In other words,
(\ref{eq.prop.int.wXi=Xi.a.1}) holds. This shows that
(\ref{eq.prop.int.wXi=Xi.a.3}) implies (\ref{eq.prop.int.wXi=Xi.a.1}).
\end{proof}

We thus conclude that the statements (\ref{eq.prop.int.wXi=Xi.a.3}) and
(\ref{eq.prop.int.wXi=Xi.a.1}) are equivalent (since
(\ref{eq.prop.int.wXi=Xi.a.3}) implies (\ref{eq.prop.int.wXi=Xi.a.1}), whereas
(\ref{eq.prop.int.wXi=Xi.a.1}) implies (\ref{eq.prop.int.wXi=Xi.a.3})).
Combined with the other equivalences we have shown above, this entails that
all four statements (\ref{eq.prop.int.wXi=Xi.a.1}),
(\ref{eq.prop.int.wXi=Xi.a.2}), (\ref{eq.prop.int.wXi=Xi.a.3}) and
(\ref{eq.prop.int.wXi=Xi.a.4}) are equivalent. Thus, Proposition
\ref{prop.int.wXi=Xi} \textbf{(a)} is proved. \medskip

\textbf{(b)} The following proof is long and yet free of any substance. The
underlying idea is that a permutation $w\in S_{n}$ that satisfies $w\left(
X_{i}\right)  \subseteq X_{i}$ for all $i\in\left[  k\right]  $ is a
permutation that permutes each of the $k$ subsets $X_{1},X_{2},\ldots,X_{k}$
individually, and thus can be split up into a permutation of $X_{1}$, a
permutation of $X_{2}$, and so on. This splitting-up can be used to write $w$
as a product $w_{1}w_{2}\cdots w_{k}$ with $w_{i}\in S_{n,X_{i}}$ for each
$i\in\left[  k\right]  $ (where $S_{n,X_{i}}$ is defined as in Proposition
\ref{prop.intX.basics}), and thus the sums on the left hand sides of
(\ref{eq.prop.int.wXi=Xi.b.1}) and (\ref{eq.prop.int.wXi=Xi.b.2}) can be
reindexed as sums over $k$-tuples $\left(  w_{1},w_{2},\ldots,w_{k}\right)  $,
which then (by the product rule) become the right hand sides.

To make this proof rigorous, we need to formally define the \textquotedblleft
splitting-up\textquotedblright\ procedure and prove all its intuitively
obvious properties. Thus we arrive at the following proof:

For any subset $X$ of $\left[  n\right]  $, define a subset $S_{n,X}$ of
$S_{n}$ as in Proposition \ref{prop.intX.basics}.

Define a subset $K$ of $S_{n}$ by%
\[
K:=\left\{  w\in S_{n}\ \mid\ w\left(  X_{i}\right)  \subseteq X_{i}\text{ for
all }i\in\left[  k\right]  \right\}  .
\]

Clearly, this set $K$ is a subgroup of $S_{n}$ (for example, it is closed
under multiplication because if two permutations $u,v\in S_{n}$ satisfy
$u\left(  X_{i}\right)  \subseteq X_{i}$ and $v\left(  X_{i}\right)  \subseteq
X_{i}$, then their product $uv$ satisfies $\left(  uv\right)  \left(
X_{i}\right)  =u\left(  \underbrace{v\left(  X_{i}\right)  }_{\subseteq X_{i}%
}\right)  \subseteq u\left(  X_{i}\right)  \subseteq X_{i}$).

We note the following simple fact first:

\begin{statement}
\textit{Claim 1:} For any $j\in\left[  k\right]  $, the set $S_{n,X_{j}}$ is a
subgroup of $K$.
\end{statement}

\begin{proof}
[Proof of Claim 1.]Let $j\in\left[  k\right]  $. We must prove that the set
$S_{n,X_{j}}$ is a subgroup of $K$. We already know (from Proposition
\ref{prop.intX.basics} \textbf{(a)}, applied to $X=X_{j}$) that $S_{n,X_{j}}$
is a subgroup of $S_{n}$, and thus is closed under multiplication and taking
inverses. Hence, it remains to prove that $S_{n,X_{j}}$ is a subset of $K$.

Let $u\in S_{n,X_{j}}$ be arbitrary. Thus,
\begin{equation}
u\left(  i\right)  =i\ \ \ \ \ \ \ \ \ \ \text{for all }i\notin X_{j}
\label{pf.prop.int.wXi=Xi.b.1}%
\end{equation}
(by the definition of $S_{n,X_{j}}$). Moreover, $u$ is a permutation of
$\left[  n\right]  $, thus injective.

We shall now prove that $u\left(  X_{j}\right)  \subseteq X_{j}$.

Indeed, let $x\in X_{j}$. If we had $u\left(  x\right)  \notin X_{j}$, then we
would have $u\left(  u\left(  x\right)  \right)  =u\left(  x\right)  $ (by
(\ref{pf.prop.int.wXi=Xi.b.1}), applied to $i=u\left(  x\right)  $), which
would yield $u\left(  x\right)  =x$ (since $u$ is injective) and thus
$x=u\left(  x\right)  \notin X_{j}$, which would contradict $x\in X_{j}$.
Thus, we cannot have $u\left(  x\right)  \notin X_{j}$. Hence, $u\left(
x\right)  \in X_{j}$.

Forget that we fixed $x$. We thus have shown that $u\left(  x\right)  \in
X_{j}$ for each $x\in X_{j}$. In other words, $u\left(  X_{j}\right)
\subseteq X_{j}$.

Next, we shall show that
\begin{equation}
u\left(  X_{i}\right)  \subseteq X_{i}\ \ \ \ \ \ \ \ \ \ \text{for each }%
i\in\left[  k\right]  . \label{pf.prop.int.wXi=Xi.b.2}%
\end{equation}
Indeed, let $i\in\left[  k\right]  $ be arbitrary. We must prove that
$u\left(  X_{i}\right)  \subseteq X_{i}$. If $i=j$, then this follows from
$u\left(  X_{j}\right)  \subseteq X_{j}$. Thus, we WLOG assume that $i\neq j$.
Hence, the sets $X_{i}$ and $X_{j}$ are disjoint (since the sets $X_{1}%
,X_{2},\ldots,X_{k}$ are disjoint). Now, let $x\in X_{i}$. Therefore, $x\notin
X_{j}$ (since the sets $X_{i}$ and $X_{j}$ are disjoint, so that an element of
$X_{i}$ cannot belong to $X_{j}$). Hence, (\ref{pf.prop.int.wXi=Xi.b.1})
(applied to $x$ instead of $i$) yields $u\left(  x\right)  =x\in X_{i}$.

Forget that we fixed $x$. We thus have shown that $u\left(  x\right)  \in
X_{i}$ for each $x\in X_{i}$. In other words, $u\left(  X_{i}\right)
\subseteq X_{i}$.

Forget that we fixed $i$. Thus we have proved that $u\left(  X_{i}\right)
\subseteq X_{i}$ for each $i\in\left[  k\right]  $. In other words, $u\in K$
(by the definition of $K$).

Forget that we fixed $u$. We thus have shown that $u\in K$ for each $u\in
S_{n,X_{j}}$. In other words, $S_{n,X_{j}}\subseteq K$. Hence, $S_{n,X_{j}}$
is a subgroup of $K$ (since $S_{n,X_{j}}$ is a subgroup of $S_{n}$). This
proves Claim 1.
\end{proof}

Now, each $k$-tuple $\left(  w_{1},w_{2},\ldots,w_{k}\right)  \in S_{n,X_{1}%
}\times S_{n,X_{2}}\times\cdots\times S_{n,X_{k}}$ satisfies $w_{j}\in
S_{n,X_{j}}\subseteq K$ for each $j\in\left[  k\right]  $ (by Claim 1) and
therefore $w_{1}w_{2}\cdots w_{k}\in K$ (since $K$ is a subgroup of $S_{n}$).
Thus, we can define a map%
\begin{align*}
\Psi:S_{n,X_{1}}\times S_{n,X_{2}}\times\cdots\times S_{n,X_{k}}  &
\rightarrow K,\\
\left(  w_{1},w_{2},\ldots,w_{k}\right)   &  \mapsto w_{1}w_{2}\cdots w_{k}.
\end{align*}
Consider this map $\Psi$.

For each $j\in\left[  k\right]  $ and each $w\in K$, we define a map
$w^{\left(  j\right)  }:\left[  n\right]  \rightarrow\left[  n\right]  $ by
setting%
\[
w^{\left(  j\right)  }\left(  x\right)  =%
\begin{cases}
w\left(  x\right)  , & \text{if }x\in X_{j};\\
x, & \text{if }x\notin X_{j}%
\end{cases}
\ \ \ \ \ \ \ \ \ \ \text{for each }x\in\left[  n\right]  .
\]
We claim the following:

\begin{statement}
\textit{Claim 2:} Let $j\in\left[  k\right]  $ and $w\in K$. Then, $w^{\left(
j\right)  }\in S_{n,X_{j}}$.
\end{statement}

\begin{proof}
[Proof of Claim 2.]The permutation $w$ belongs to $K$, and thus satisfies
\newline$\left(  w\left(  X_{i}\right)  \subseteq X_{i}\text{ for all }%
i\in\left[  k\right]  \right)  $ (by the definition of $K$). Hence, it
satisfies \newline$\left(  w\left(  X_{i}\right)  =X_{i}\text{ for all }%
i\in\left[  k\right]  \right)  $ (since Proposition \ref{prop.int.wXi=Xi}
\textbf{(a)} shows that the two statements $\left(  w\left(  X_{i}\right)
\subseteq X_{i}\text{ for all }i\in\left[  k\right]  \right)  $ and $\left(
w\left(  X_{i}\right)  =X_{i}\text{ for all }i\in\left[  k\right]  \right)  $
are equivalent). Applying this to $i=j$, we obtain $w\left(  X_{j}\right)
=X_{j}$.

Now, we shall prove that the map $w^{\left(  j\right)  }$ is surjective. In
other words, we shall prove that $y\in w^{\left(  j\right)  }\left(  \left[
n\right]  \right)  $ for each $y\in\left[  n\right]  $.

Indeed, fix $y\in\left[  n\right]  $. We must prove that $y\in w^{\left(
j\right)  }\left(  \left[  n\right]  \right)  $.

If $y\notin X_{j}$, then $w^{\left(  j\right)  }\left(  y\right)  =y$ (since
the definition of $w^{\left(  j\right)  }$ yields that $w^{\left(  j\right)
}\left(  x\right)  =x$ whenever $x\notin X_{j}$) and thus $y=w^{\left(
j\right)  }\left(  y\right)  \in w^{\left(  j\right)  }\left(  \left[
n\right]  \right)  $. Hence, we can WLOG assume that we don't have $y\notin
X_{j}$. Thus, $y\in X_{j}$. In other words, $y\in w\left(  X_{j}\right)  $
(since $w\left(  X_{j}\right)  =X_{j}$). In other words, $y=w\left(  x\right)
$ for some $x\in X_{j}$. Consider this $x$. The definition of $w^{\left(
j\right)  }$ yields $w^{\left(  j\right)  }\left(  x\right)  =w\left(
x\right)  $ (since $x\in X_{j}$). Hence, $y=w\left(  x\right)  =w^{\left(
j\right)  }\left(  x\right)  \in w^{\left(  j\right)  }\left(  \left[
n\right]  \right)  $.

Forget that we fixed $y$. We thus have proved that $y\in w^{\left(  j\right)
}\left(  \left[  n\right]  \right)  $ for each $y\in\left[  n\right]  $. In
other words, $\left[  n\right]  \subseteq w^{\left(  j\right)  }\left(
\left[  n\right]  \right)  $. Thus, the map $w^{\left(  j\right)  }$ is surjective.

But the pigeonhole principle says that any surjective map from a finite set to
itself must be bijective. Hence, the map $w^{\left(  j\right)  }$ is bijective
(since $w^{\left(  j\right)  }$ is a surjective map from the finite set
$\left[  n\right]  $ to itself). Thus, $w^{\left(  j\right)  }$ is a
permutation of $\left[  n\right]  $. In other words, $w^{\left(  j\right)
}\in S_{n}$. Moreover, $w^{\left(  j\right)  }\left(  x\right)  =x$ for all
$x\notin X_{j}$ (by the definition of $w^{\left(  j\right)  }$). In other
words, $w^{\left(  j\right)  }\left(  i\right)  =i$ for all $i\notin X_{j}$.
This entails $w^{\left(  j\right)  }\in S_{n,X_{j}}$ (by the definition of
$S_{n,X_{j}}$, since $w^{\left(  j\right)  }\in S_{n}$). This proves Claim 2.
\end{proof}

Claim 2 shows that each $w\in K$ satisfies $w^{\left(  j\right)  }\in
S_{n,X_{j}}$ for all $j\in\left[  k\right]  $, and therefore $\left(
w^{\left(  1\right)  },w^{\left(  2\right)  },\ldots,w^{\left(  k\right)
}\right)  \in S_{n,X_{1}}\times S_{n,X_{2}}\times\cdots\times S_{n,X_{k}}$.
Hence, we can define a map
\begin{align*}
\Phi:K  &  \rightarrow S_{n,X_{1}}\times S_{n,X_{2}}\times\cdots\times
S_{n,X_{k}},\\
w  &  \mapsto\left(  w^{\left(  1\right)  },w^{\left(  2\right)  }%
,\ldots,w^{\left(  k\right)  }\right)  .
\end{align*}
Consider this map $\Phi$. We shall now show that the maps $\Phi$ and $\Psi$
are mutually inverse. First, we shall show the following claim:

\begin{statement}
\textit{Claim 3:} Let $w\in K$. Then, $w=w^{\left(  1\right)  }w^{\left(
2\right)  }\cdots w^{\left(  k\right)  }$.
\end{statement}

\begin{proof}
[Proof of Claim 3.]Let $x\in\left[  n\right]  $. Then, $x\in\left[  n\right]
=X_{1}\cup X_{2}\cup\cdots\cup X_{k}$. Hence, $x\in X_{p}$ for some
$p\in\left[  k\right]  $. Consider this $p$. By the definition of $w^{\left(
p\right)  }$, we have $w^{\left(  p\right)  }\left(  x\right)  =%
\begin{cases}
w\left(  x\right)  , & \text{if }x\in X_{p};\\
x, & \text{if }x\notin X_{p}%
\end{cases}
\ \ =w\left(  x\right)  $ (since $x\in X_{p}$). But $w\in K$, so that $\left(
w\left(  X_{i}\right)  \subseteq X_{i}\text{ for all }i\in\left[  k\right]
\right)  $ (by the definition of $K$). Applying this to $i=p$, we find
$w\left(  X_{p}\right)  \subseteq X_{p}$.

Set $y:=w\left(  x\right)  $. Then, $y=w\left(  x\right)  \in w\left(
X_{p}\right)  $ (since $x\in X_{p}$), so that $y\in w\left(  X_{p}\right)
\subseteq X_{p}$.

The $k$ sets $X_{1},X_{2},\ldots,X_{k}$ are disjoint. Thus, for each
$i\in\left\{  p+1,p+2,\ldots,k\right\}  $, we have $w^{\left(  i\right)
}\left(  x\right)  =x$\ \ \ \ \footnote{\textit{Proof.} Let $i\in\left\{
p+1,p+2,\ldots,k\right\}  $. Thus, $i\neq p$. Hence, the sets $X_{i}$ and
$X_{p}$ are disjoint (since $X_{1},X_{2},\ldots,X_{k}$ are disjoint). Thus, an
element of $X_{p}$ cannot belong to $X_{i}$. From $x\in X_{p}$, we thus obtain
$x\notin X_{i}$. Hence, $w^{\left(  i\right)  }\left(  x\right)  =x$ (since
the definition of $w^{\left(  i\right)  }$ yields that $w^{\left(  i\right)
}\left(  z\right)  =z$ for each $z\notin X_{i}$).}. In other words, all $k-p$
permutations $w^{\left(  p+1\right)  },w^{\left(  p+2\right)  },\ldots
,w^{\left(  k\right)  }$ leave the element $x$ unchanged. Hence, so does their
product $w^{\left(  p+1\right)  }w^{\left(  p+2\right)  }\cdots w^{\left(
k\right)  }$. In other words,%
\[
\left(  w^{\left(  p+1\right)  }w^{\left(  p+2\right)  }\cdots w^{\left(
k\right)  }\right)  \left(  x\right)  =x.
\]
Hence,%
\begin{align*}
\underbrace{\left(  w^{\left(  p\right)  }w^{\left(  p+1\right)  }\cdots
w^{\left(  k\right)  }\right)  }_{=w^{\left(  p\right)  }\circ\left(
w^{\left(  p+1\right)  }w^{\left(  p+2\right)  }\cdots w^{\left(  k\right)
}\right)  }\left(  x\right)   &  =w^{\left(  p\right)  }\left(
\underbrace{\left(  w^{\left(  p+1\right)  }w^{\left(  p+2\right)  }\cdots
w^{\left(  k\right)  }\right)  \left(  x\right)  }_{=x}\right) \\
&  =w^{\left(  p\right)  }\left(  x\right)  =w\left(  x\right)  =y.
\end{align*}

The $k$ sets $X_{1},X_{2},\ldots,X_{k}$ are disjoint. Thus, for each
$i\in\left\{  1,2,\ldots,p-1\right\}  $, we have $w^{\left(  i\right)
}\left(  y\right)  =y$\ \ \ \ \footnote{\textit{Proof.} Let $i\in\left\{
1,2,\ldots,p-1\right\}  $. Thus, $i\neq p$. Hence, the sets $X_{i}$ and
$X_{p}$ are disjoint (since $X_{1},X_{2},\ldots,X_{k}$ are disjoint). Thus, an
element of $X_{p}$ cannot belong to $X_{i}$. From $y\in X_{p}$, we thus obtain
$y\notin X_{i}$. Hence, $w^{\left(  i\right)  }\left(  y\right)  =y$ (since
the definition of $w^{\left(  i\right)  }$ yields that $w^{\left(  i\right)
}\left(  z\right)  =z$ for each $z\notin X_{i}$).}. In other words, all $p-1$
permutations $w^{\left(  1\right)  },w^{\left(  2\right)  },\ldots,w^{\left(
p-1\right)  }$ leave the element $y$ unchanged. Hence, so does their product
$w^{\left(  1\right)  }w^{\left(  2\right)  }\cdots w^{\left(  p-1\right)  }$.
In other words,%
\[
\left(  w^{\left(  1\right)  }w^{\left(  2\right)  }\cdots w^{\left(
p-1\right)  }\right)  \left(  y\right)  =y.
\]

Now,%
\begin{align*}
&  \underbrace{\left(  w^{\left(  1\right)  }w^{\left(  2\right)  }\cdots
w^{\left(  k\right)  }\right)  }_{=\left(  w^{\left(  1\right)  }w^{\left(
2\right)  }\cdots w^{\left(  p-1\right)  }\right)  \circ\left(  w^{\left(
p\right)  }w^{\left(  p+1\right)  }\cdots w^{\left(  k\right)  }\right)
}\left(  x\right) \\
&  =\left(  w^{\left(  1\right)  }w^{\left(  2\right)  }\cdots w^{\left(
p-1\right)  }\right)  \left(  \underbrace{\left(  w^{\left(  p\right)
}w^{\left(  p+1\right)  }\cdots w^{\left(  k\right)  }\right)  \left(
x\right)  }_{=y}\right) \\
&  =\left(  w^{\left(  1\right)  }w^{\left(  2\right)  }\cdots w^{\left(
p-1\right)  }\right)  \left(  y\right)  =y=w\left(  x\right)  .
\end{align*}

Forget that we fixed $x$. We thus have shown that $\left(  w^{\left(
1\right)  }w^{\left(  2\right)  }\cdots w^{\left(  k\right)  }\right)  \left(
x\right)  =w\left(  x\right)  $ for each $x\in\left[  n\right]  $. In other
words, $w^{\left(  1\right)  }w^{\left(  2\right)  }\cdots w^{\left(
k\right)  }=w$. This proves Claim 3.
\end{proof}

Claim 3 easily yields the following:

\begin{statement}
\textit{Claim 4:} We have $\Psi\circ\Phi=\operatorname*{id}$.
\end{statement}

\begin{proof}
[Proof of Claim 4.]Let $w\in K$. Then, the definition of $\Phi$ yields
\[
\Phi\left(  w\right)  =\left(  w^{\left(  1\right)  },w^{\left(  2\right)
},\ldots,w^{\left(  k\right)  }\right)  .
\]
Applying the map $\Psi$ to this equality, we obtain%
\begin{align*}
\Psi\left(  \Phi\left(  w\right)  \right)   &  =\Psi\left(  w^{\left(
1\right)  },w^{\left(  2\right)  },\ldots,w^{\left(  k\right)  }\right) \\
&  =w^{\left(  1\right)  }w^{\left(  2\right)  }\cdots w^{\left(  k\right)
}\ \ \ \ \ \ \ \ \ \ \left(  \text{by the definition of }\Psi\right) \\
&  =w\ \ \ \ \ \ \ \ \ \ \left(  \text{by Claim 3}\right)  .
\end{align*}
Hence, $\left(  \Psi\circ\Phi\right)  \left(  w\right)  =\Psi\left(
\Phi\left(  w\right)  \right)  =w=\operatorname*{id}\left(  w\right)  $.

Forget that we fixed $w$. We thus have shown that $\left(  \Psi\circ
\Phi\right)  \left(  w\right)  =\operatorname*{id}\left(  w\right)  $ for each
$w\in K$. In other words, $\Psi\circ\Phi=\operatorname*{id}$. This proves
Claim 4.
\end{proof}

On the other hand, the following holds:

\begin{statement}
\textit{Claim 5:} Let $\left(  w_{1},w_{2},\ldots,w_{k}\right)  \in
S_{n,X_{1}}\times S_{n,X_{2}}\times\cdots\times S_{n,X_{k}}$. Set
$w=w_{1}w_{2}\cdots w_{k}$. Then, $w^{\left(  j\right)  }=w_{j}$ for each
$j\in\left[  k\right]  $.
\end{statement}

\begin{proof}
[Proof of Claim 5.]We have $\left(  w_{1},w_{2},\ldots,w_{k}\right)  \in
S_{n,X_{1}}\times S_{n,X_{2}}\times\cdots\times S_{n,X_{k}}$. Thus, for each
$j\in\left[  k\right]  $, we have%
\begin{align}
w_{j}  &  \in S_{n,X_{j}}\label{pf.prop.int.wXi=Xi.b.4}\\
&  =\left\{  w\in S_{n}\ \mid\ w\left(  i\right)  =i\text{ for all }i\notin
X_{j}\right\} \nonumber
\end{align}
(by the definition of $S_{n,X_{j}}$) and therefore%
\begin{equation}
w_{j}\left(  i\right)  =i\ \ \ \ \ \ \ \ \ \ \text{for all }i\notin X_{j}.
\label{pf.prop.int.wXi=Xi.b.5}%
\end{equation}

Let $j\in\left[  k\right]  $. Let $x\in\left[  n\right]  $. We shall show that
$w^{\left(  j\right)  }\left(  x\right)  =w_{j}\left(  x\right)  $.

We are in one of the following two cases:

\textit{Case 1:} We have $x\in X_{j}$.

\textit{Case 2:} We have $x\notin X_{j}$.

Let us first consider Case 1. In this case, we have $x\in X_{j}$. Hence,
$w^{\left(  j\right)  }\left(  x\right)  =w\left(  x\right)  $ (by the
definition of $w^{\left(  j\right)  }$).

Set $y:=w_{j}\left(  x\right)  $. Note that $S_{n,X_{j}}\subseteq K$ (by Claim
1), so that (\ref{pf.prop.int.wXi=Xi.b.4}) becomes $w_{j}\in S_{n,X_{j}%
}\subseteq K$. Hence, we have $w_{j}\left(  X_{i}\right)  \subseteq X_{i}$ for
all $i\in\left[  k\right]  $ (by the definition of $K$). Applying this to
$i=j$, we find $w_{j}\left(  X_{j}\right)  \subseteq X_{j}$. Hence, from $x\in
X_{j}$, we obtain $w_{j}\left(  x\right)  \in w_{j}\left(  X_{j}\right)
\subseteq X_{j}$. In other words, $y\in X_{j}$ (since $y=w_{j}\left(
x\right)  $).

The $k$ sets $X_{1},X_{2},\ldots,X_{k}$ are disjoint. Thus, for each
$i\in\left\{  j+1,j+2,\ldots,k\right\}  $, we have $w_{i}\left(  x\right)
=x$\ \ \ \ \footnote{\textit{Proof.} Let $i\in\left\{  j+1,j+2,\ldots
,k\right\}  $. Thus, $i\neq j$. Hence, the sets $X_{i}$ and $X_{j}$ are
disjoint (since $X_{1},X_{2},\ldots,X_{k}$ are disjoint). Thus, an element of
$X_{j}$ cannot belong to $X_{i}$. From $x\in X_{j}$, we thus obtain $x\notin
X_{i}$. Hence, (\ref{pf.prop.int.wXi=Xi.b.5}) (applied to $i$ and $x$ instead
of $j$ and $i$) yields $w_{i}\left(  x\right)  =x$.}. In other words, all
$k-j$ permutations $w_{j+1},w_{j+2},\ldots,w_{k}$ leave the element $x$
unchanged. Hence, so does their product $w_{j+1}w_{j+2}\cdots w_{k}$. In other
words,%
\[
\left(  w_{j+1}w_{j+2}\cdots w_{k}\right)  \left(  x\right)  =x.
\]
Hence,%
\[
\underbrace{\left(  w_{j}w_{j+1}\cdots w_{k}\right)  }_{=w_{j}\circ\left(
w_{j+1}w_{j+2}\cdots w_{k}\right)  }\left(  x\right)  =w_{j}\left(
\underbrace{\left(  w_{j+1}w_{j+2}\cdots w_{k}\right)  \left(  x\right)
}_{=x}\right)  =w_{j}\left(  x\right)  =y.
\]

The $k$ sets $X_{1},X_{2},\ldots,X_{k}$ are disjoint. Thus, for each
$i\in\left\{  1,2,\ldots,j-1\right\}  $, we have $w_{i}\left(  y\right)
=y$\ \ \ \ \footnote{\textit{Proof.} Let $i\in\left\{  1,2,\ldots,j-1\right\}
$. Thus, $i\neq j$. Hence, the sets $X_{i}$ and $X_{j}$ are disjoint (since
$X_{1},X_{2},\ldots,X_{k}$ are disjoint). Thus, an element of $X_{j}$ cannot
belong to $X_{i}$. From $y\in X_{j}$, we thus obtain $y\notin X_{i}$. Hence,
(\ref{pf.prop.int.wXi=Xi.b.5}) (applied to $i$ and $y$ instead of $j$ and $i$)
yields $w_{i}\left(  y\right)  =y$.}. In other words, all $j-1$ permutations
$w_{1},w_{2},\ldots,w_{j-1}$ leave the element $y$ unchanged. Hence, so does
their product $w_{1}w_{2}\cdots w_{j-1}$. In other words,%
\[
\left(  w_{1}w_{2}\cdots w_{j-1}\right)  \left(  y\right)  =y.
\]

Now,%
\begin{align*}
w^{\left(  j\right)  }\left(  x\right)   &  =\underbrace{w}_{\substack{=w_{1}%
w_{2}\cdots w_{k}\\=\left(  w_{1}w_{2}\cdots w_{j-1}\right)  \circ\left(
w_{j}w_{j+1}\cdots w_{k}\right)  }}\left(  x\right)  =\left(  w_{1}w_{2}\cdots
w_{j-1}\right)  \left(  \underbrace{\left(  w_{j}w_{j+1}\cdots w_{k}\right)
\left(  x\right)  }_{=y}\right) \\
&  =\left(  w_{1}w_{2}\cdots w_{j-1}\right)  \left(  y\right)  =y=w_{j}\left(
x\right)  .
\end{align*}

Thus, we have proved $w^{\left(  j\right)  }\left(  x\right)  =w_{j}\left(
x\right)  $ in Case 1.

Let us now consider Case 2. In this case, we have $x\notin X_{j}$. Hence,
$w^{\left(  j\right)  }\left(  x\right)  =x$ (by the definition of $w^{\left(
j\right)  }$) and $w_{j}\left(  x\right)  =x$ (by
(\ref{pf.prop.int.wXi=Xi.b.5}), applied to $i=x$). Comparing these two
equalities, we obtain $w^{\left(  j\right)  }\left(  x\right)  =w_{j}\left(
x\right)  $. Thus, we have proved $w^{\left(  j\right)  }\left(  x\right)
=w_{j}\left(  x\right)  $ in Case 2.

We have now proved $w^{\left(  j\right)  }\left(  x\right)  =w_{j}\left(
x\right)  $ in both Cases 1 and 2. Hence, $w^{\left(  j\right)  }\left(
x\right)  =w_{j}\left(  x\right)  $ always holds.

Forget that we fixed $x$. We have now shown that $w^{\left(  j\right)
}\left(  x\right)  =w_{j}\left(  x\right)  $ for all $x\in\left[  n\right]  $.
In other words, $w^{\left(  j\right)  }=w_{j}$. This proves Claim 5.
\end{proof}

Claim 5, in turn, easily yields the following:

\begin{statement}
\textit{Claim 6:} We have $\Phi\circ\Psi=\operatorname*{id}$.
\end{statement}

\begin{proof}
[Proof of Claim 6.]Let $\mathbf{w}\in S_{n,X_{1}}\times S_{n,X_{2}}%
\times\cdots\times S_{n,X_{k}}$. Write this $k$-tuple $\mathbf{w}$ in the form
$\mathbf{w}=\left(  w_{1},w_{2},\ldots,w_{k}\right)  $. Thus,%
\[
\Psi\left(  \mathbf{w}\right)  =\Psi\left(  w_{1},w_{2},\ldots,w_{k}\right)
=w_{1}w_{2}\cdots w_{k}\ \ \ \ \ \ \ \ \ \ \left(  \text{by the definition of
}\Psi\right)  .
\]

Set $w:=w_{1}w_{2}\cdots w_{k}$. Thus, $\Psi\left(  \mathbf{w}\right)
=w_{1}w_{2}\cdots w_{k}=w$. Applying the map $\Phi$ to this equality, we
obtain%
\[
\Phi\left(  \Psi\left(  \mathbf{w}\right)  \right)  =\Phi\left(  w\right)
=\left(  w^{\left(  1\right)  },w^{\left(  2\right)  },\ldots,w^{\left(
k\right)  }\right)  \ \ \ \ \ \ \ \ \ \ \left(  \text{by the definition of
}\Phi\right)  .
\]
However, Claim 5 yields that $w^{\left(  j\right)  }=w_{j}$ for each
$j\in\left[  k\right]  $. In other words,%
\[
\left(  w^{\left(  1\right)  },w^{\left(  2\right)  },\ldots,w^{\left(
k\right)  }\right)  =\left(  w_{1},w_{2},\ldots,w_{k}\right)  .
\]
Hence,%
\begin{align*}
\left(  \Phi\circ\Psi\right)  \left(  \mathbf{w}\right)   &  =\Phi\left(
\Psi\left(  \mathbf{w}\right)  \right)  =\left(  w^{\left(  1\right)
},w^{\left(  2\right)  },\ldots,w^{\left(  k\right)  }\right) \\
&  =\left(  w_{1},w_{2},\ldots,w_{k}\right)  =\mathbf{w}=\operatorname*{id}%
\left(  \mathbf{w}\right)  .
\end{align*}

Since we have proved this for each $\mathbf{w}\in S_{n,X_{1}}\times
S_{n,X_{2}}\times\cdots\times S_{n,X_{k}}$, we thus conclude that $\Phi
\circ\Psi=\operatorname*{id}$. Thus, Claim 6 is proved.
\end{proof}

Combining Claim 4 with Claim 6, we see that the maps $\Phi$ and $\Psi$ are
mutually inverse. Hence, the maps $\Phi$ and $\Psi$ are invertible, i.e., are
bijections. Now, recall that $K$ was defined as the set of all permutations
$w\in S_{n}$ that satisfy $\left(  w\left(  X_{i}\right)  \subseteq
X_{i}\text{ for all }i\in\left[  k\right]  \right)  $. Hence, the summation
sign $\sum_{\substack{w\in S_{n};\\w\left(  X_{i}\right)  \subseteq
X_{i}\text{ for all }i\in\left[  k\right]  }}$ can be rewritten as $\sum_{w\in
K}$. Thus,%
\begin{align}
&  \sum_{\substack{w\in S_{n};\\w\left(  X_{i}\right)  \subseteq X_{i}\text{
for all }i\in\left[  k\right]  }}\left(  -1\right)  ^{w}w\nonumber\\
&  =\sum_{w\in K}\left(  -1\right)  ^{w}w\nonumber\\
&  =\sum_{\left(  w_{1},w_{2},\ldots,w_{k}\right)  \in S_{n,X_{1}}\times
S_{n,X_{2}}\times\cdots\times S_{n,X_{k}}}\underbrace{\left(  -1\right)
^{\Psi\left(  w_{1},w_{2},\ldots,w_{k}\right)  }\Psi\left(  w_{1},w_{2}%
,\ldots,w_{k}\right)  }_{\substack{=\left(  -1\right)  ^{w_{1}w_{2}\cdots
w_{k}}w_{1}w_{2}\cdots w_{k}\\\text{(since }\Psi\left(  w_{1},w_{2}%
,\ldots,w_{k}\right)  =w_{1}w_{2}\cdots w_{k}\\\text{(by the definition of
}\Psi\text{))}}}\nonumber\\
&  \ \ \ \ \ \ \ \ \ \ \ \ \ \ \ \ \ \ \ \ \left(
\begin{array}
[c]{c}%
\text{here, we have substituted }\Psi\left(  w_{1},w_{2},\ldots,w_{k}\right)
\\
\text{for }w\text{ in the sum, since the}\\
\text{map }\Psi:S_{n,X_{1}}\times S_{n,X_{2}}\times\cdots\times S_{n,X_{k}%
}\rightarrow K\\
\text{is a bijection}%
\end{array}
\right) \nonumber\\
&  =\sum_{\left(  w_{1},w_{2},\ldots,w_{k}\right)  \in S_{n,X_{1}}\times
S_{n,X_{2}}\times\cdots\times S_{n,X_{k}}}\underbrace{\left(  -1\right)
^{w_{1}w_{2}\cdots w_{k}}}_{\substack{=\left(  -1\right)  ^{w_{1}}\left(
-1\right)  ^{w_{2}}\cdots\left(  -1\right)  ^{w_{k}}\\\text{(by the
multiplicativity of the sign)}}}w_{1}w_{2}\cdots w_{k}\nonumber\\
&  =\sum_{\left(  w_{1},w_{2},\ldots,w_{k}\right)  \in S_{n,X_{1}}\times
S_{n,X_{2}}\times\cdots\times S_{n,X_{k}}}\underbrace{\left(  -1\right)
^{w_{1}}\left(  -1\right)  ^{w_{2}}\cdots\left(  -1\right)  ^{w_{k}}w_{1}%
w_{2}\cdots w_{k}}_{=\left(  \left(  -1\right)  ^{w_{1}}w_{1}\right)  \left(
\left(  -1\right)  ^{w_{2}}w_{2}\right)  \cdots\left(  \left(  -1\right)
^{w_{k}}w_{k}\right)  }\nonumber\\
&  =\sum_{\left(  w_{1},w_{2},\ldots,w_{k}\right)  \in S_{n,X_{1}}\times
S_{n,X_{2}}\times\cdots\times S_{n,X_{k}}}\left(  \left(  -1\right)  ^{w_{1}%
}w_{1}\right)  \left(  \left(  -1\right)  ^{w_{2}}w_{2}\right)  \cdots\left(
\left(  -1\right)  ^{w_{k}}w_{k}\right)  . \label{pf.prop.int.wXi=Xi.b.-3}%
\end{align}

However, each $i\in\left[  k\right]  $ satisfies%
\begin{equation}
\nabla_{X_{i}}^{-}=\sum\limits_{w\in S_{n,X_{i}}}\left(  -1\right)  ^{w}w
\label{pf.prop.int.wXi=Xi.b.-2}%
\end{equation}
(by the second equality in Proposition \ref{prop.intX.basics} \textbf{(b)},
applied to $X=X_{i}$). Multiplying the equalities
(\ref{pf.prop.int.wXi=Xi.b.-2}) over all $i\in\left\{  1,2,\ldots,k\right\}
$, we obtain%
\begin{align*}
&  \nabla_{X_{1}}^{-}\nabla_{X_{2}}^{-}\cdots\nabla_{X_{k}}^{-}\\
&  =\left(  \sum\limits_{w\in S_{n,X_{1}}}\left(  -1\right)  ^{w}w\right)
\left(  \sum\limits_{w\in S_{n,X_{2}}}\left(  -1\right)  ^{w}w\right)
\cdots\left(  \sum\limits_{w\in S_{n,X_{k}}}\left(  -1\right)  ^{w}w\right) \\
&  =\sum_{\left(  w_{1},w_{2},\ldots,w_{k}\right)  \in S_{n,X_{1}}\times
S_{n,X_{2}}\times\cdots\times S_{n,X_{k}}}\left(  \left(  -1\right)  ^{w_{1}%
}w_{1}\right)  \left(  \left(  -1\right)  ^{w_{2}}w_{2}\right)  \cdots\left(
\left(  -1\right)  ^{w_{k}}w_{k}\right)
\end{align*}
(by the product rule, i.e., by expanding the product). Comparing this with
(\ref{pf.prop.int.wXi=Xi.b.-3}), we obtain%
\[
\sum_{\substack{w\in S_{n};\\w\left(  X_{i}\right)  \subseteq X_{i}\text{ for
all }i\in\left[  k\right]  }}\left(  -1\right)  ^{w}w=\nabla_{X_{1}}^{-}%
\nabla_{X_{2}}^{-}\cdots\nabla_{X_{k}}^{-}.
\]
An analogous argument (but with all the signs removed) shows that%
\[
\sum_{\substack{w\in S_{n};\\w\left(  X_{i}\right)  \subseteq X_{i}\text{ for
all }i\in\left[  k\right]  }}w=\nabla_{X_{1}}\nabla_{X_{2}}\cdots\nabla
_{X_{k}}.
\]
Thus, Proposition \ref{prop.int.wXi=Xi} \textbf{(b)} is proved. \medskip

\textbf{(c)} If $i$ and $j$ are two distinct elements of $\left[  k\right]  $,
then the sets $X_{i}$ and $X_{j}$ are disjoint (since all $k$ sets
$X_{1},X_{2},\ldots,X_{k}$ are disjoint) and thus satisfy $\nabla_{X_{i}%
}\nabla_{X_{j}}=\nabla_{X_{j}}\nabla_{X_{i}}$ (by the equality $\nabla
_{X}\nabla_{Y}=\nabla_{Y}\nabla_{X}$ from Proposition \ref{prop.int.commute}
\textbf{(b)}, applied to $X=X_{i}$ and $Y=X_{j}$). In other words, the $k$
elements $\nabla_{X_{1}},\nabla_{X_{2}},\ldots,\nabla_{X_{k}}$ pairwise
commute. Similarly, the $k$ elements $\nabla_{X_{1}}^{-},\nabla_{X_{2}}%
^{-},\ldots,\nabla_{X_{k}}^{-}$ pairwise commute (here we need to use the
equality $\nabla_{X}^{-}\nabla_{Y}^{-}=\nabla_{Y}^{-}\nabla_{X}^{-}$ from
Proposition \ref{prop.int.commute} \textbf{(b)}). Thus, Proposition
\ref{prop.int.wXi=Xi} \textbf{(c)} follows.
\end{proof}

\subsection{Orbits and cycles}

\subsubsection{Orbits, cycles and reflection lengths}

Let us now recall some elementary combinatorics of permutations (see, e.g.,
\cite[Lectures 27--28, \S 4.3]{22fco}, \cite[\S 5.5]{21s}, \cite[\S 3.3]{26fa}):

\begin{definition}
\label{def.orbits.orbits}Let $X$ be a finite set. Let $\sigma$ be a
permutation of $X$.

Consider the binary relation $\overset{\sigma}{\sim}$ on the set $X$ defined
as follows: For any $i,j\in X$, we set%
\[
\left(  i\overset{\sigma}{\sim}j\right)  \ \Longleftrightarrow\ \left(
i=\sigma^{k}\left(  j\right)  \text{ for some }k\in\mathbb{N}\right)  .
\]

It is not hard to see that this relation $\overset{\sigma}{\sim}$ is an
equivalence relation (see \cite[Lecture 27, Proposition 4.3.5]{22fco} or
\cite[Theorem 3.3.2]{26fa}). Its equivalence classes are called the
\emph{orbits} of $\sigma$. Each element of $X$ belongs to exactly one orbit of
$\sigma$.
\end{definition}

\begin{example}
\label{exa.orbits.orbits-1}The orbits of $\operatorname*{oln}\left(
432651\right)  \in S_{6}$ are $\left\{  1,4,6\right\}  $, $\left\{
2,3\right\}  $ and $\left\{  5\right\}  $. This is easiest to see using the
\emph{cycle digraph} of this permutation, which is a convenient way of
visualizing it:

Generally, if $\sigma$ is a permutation of a finite set $X$, then the
\emph{cycle digraph} of $\sigma$ is the digraph (= directed graph) whose
vertices are the elements of $X$, and whose arcs (= directed edges) are the
pairs $\left(  i,\sigma\left(  i\right)  \right)  $ for each vertex $i$. Here
is how this digraph looks like for $\sigma=\operatorname*{oln}\left(
432651\right)  \in S_{6}$:%
\[%
%TCIMACRO{\TeXButton{tikz cycle digraph}{\begin{tikzpicture}%
%[->,shorten >=1pt,auto,node distance=3cm, thick,main node/.style={circle,fill=blue!20,draw}%
%]
%\node[main node] (1) {1};
%\node[main node] (4) [above right of=1] {4};
%\node[main node] (6) [below right of=1] {6};
%\node[main node] (2) [right of=4] {2};
%\node[main node] (3) [right of=6] {3};
%\node[main node] (5) [above of=1] {5};
%\path[-{Stealth[length=4mm]}]
%(1) edge [bend left] (4)
%(4) edge [bend left] (6)
%(6) edge [bend left] (1)
%(2) edge [bend right] (3)
%(3) edge [bend right] (2)
%(5) edge [loop right] (5);
%\end{tikzpicture}
%}}%
%BeginExpansion
\begin{tikzpicture}%
[->,shorten >=1pt,auto,node distance=3cm, thick,main node/.style={circle,fill=blue!20,draw}%
]
\node[main node] (1) {1};
\node[main node] (4) [above right of=1] {4};
\node[main node] (6) [below right of=1] {6};
\node[main node] (2) [right of=4] {2};
\node[main node] (3) [right of=6] {3};
\node[main node] (5) [above of=1] {5};
\path[-{Stealth[length=4mm]}]
(1) edge [bend left] (4)
(4) edge [bend left] (6)
(6) edge [bend left] (1)
(2) edge [bend right] (3)
(3) edge [bend right] (2)
(5) edge [loop right] (5);
\end{tikzpicture}
%EndExpansion
\ \ .
\]
The digraph visibly consists of three vertex-disjoint cycles (i.e., three
cycles that pairwise have no vertices in common). The vertex sets of these
three cycles are exactly the three orbits of the permutation.
\end{example}

The orbits of a permutation have a rather simple structure (see \cite[Lecture
27, Proposition 4.3.4]{22fco} for a proof):

\begin{remark}
\label{rmk.orbits.O=}Let $X$ be a finite set. Let $\sigma$ be a permutation of
$X$. Let $O$ be an orbit of $\sigma$. Let $a\in O$. Then,%
\begin{align*}
O  &  =\left\{  \sigma^{k}\left(  a\right)  \ \mid\ k\in\mathbb{N}\right\} \\
&  =\left\{  \sigma^{0}\left(  a\right)  ,\ \sigma^{1}\left(  a\right)
,\ \sigma^{2}\left(  a\right)  ,\ \ldots\right\}  .
\end{align*}
Moreover, if we set $m=\left\vert O\right\vert $, then $\sigma^{m}\left(
a\right)  =a$ and%
\[
O=\left\{  \sigma^{0}\left(  a\right)  ,\ \sigma^{1}\left(  a\right)
,\ \ldots,\ \sigma^{m-1}\left(  a\right)  \right\}  .
\]

\end{remark}

The orbits of a permutation do not uniquely determine it. For example, the
permutations $\operatorname*{cyc}\nolimits_{1,2,3}$ and $\operatorname*{cyc}%
\nolimits_{1,3,2}$ in $S_{3}$ have the same orbits (namely, the single orbit
$\left\{  1,2,3\right\}  $). However, there is a way to \textquotedblleft
refine\textquotedblright\ the structure of the orbits to get something that
does determine the permutation. This refinement is the \emph{disjoint cycle
decomposition}\ of the permutation:

\begin{definition}
[disjoint cycle decomposition]\label{def.dcd.dcd}Let $X$ be a finite set. Let
$\sigma$ be a permutation of $X$. Then, there is a list%
\begin{align}
&  \Big(\left(  a_{1,1},a_{1,2},\ldots,a_{1,n_{1}}\right)  ,\nonumber\\
&  \ \ \ \left(  a_{2,1},a_{2,2},\ldots,a_{2,n_{2}}\right)  ,\nonumber\\
&  \ \ \ \ldots,\nonumber\\
&  \ \ \ \left(  a_{k,1},a_{k,2},\ldots,a_{k,n_{k}}\right)
\Big) \label{eq.def.dcd.dcd.dcd}%
\end{align}
of nonempty lists of elements of $X$ such that:

\begin{itemize}
\item each element of $X$ appears exactly once in the composite list%
\begin{align}
&  (a_{1,1},a_{1,2},\ldots,a_{1,n_{1}},\nonumber\\
&  \ \ \ a_{2,1},a_{2,2},\ldots,a_{2,n_{2}},\nonumber\\
&  \ \ \ \ldots,\nonumber\\
&  \ \ \ a_{k,1},a_{k,2},\ldots,a_{k,n_{k}}), \label{eq.def.dcd.dcd.comp}%
\end{align}
and

\item we have%
\[
\sigma=\operatorname*{cyc}\nolimits_{a_{1,1},a_{1,2},\ldots,a_{1,n_{1}}}%
\circ\operatorname*{cyc}\nolimits_{a_{2,1},a_{2,2},\ldots,a_{2,n_{2}}}%
\circ\cdots\circ\operatorname*{cyc}\nolimits_{a_{k,1},a_{k,2},\ldots
,a_{k,n_{k}}}%
\]

\end{itemize}

\noindent(see \cite[Lecture 27, Theorem 4.3.11]{22fco} or \cite[Theorem
5.5.2]{21s} for the proof). Such a list is called a \emph{disjoint cycle
decomposition} (or, for short, \emph{DCD}) of $\sigma$. Its entries (which
themselves are lists of elements of $X$) are called the \emph{cycles} of
$\sigma$.
\end{definition}

The DCD of $\sigma$ is unique up to permuting the cycles and rotating each
cycle. See \cite[Lecture 27, Theorem 4.3.11]{22fco} or \cite[Theorem
5.5.2]{21s} for the details of this statement and for a proof outline.

\begin{example}
The permutation $\operatorname*{oln}\left(  432651\right)  \in S_{6}$ has
disjoint cycle decomposition%
\[
\left(  \left(  1,4,6\right)  ,\ \left(  2,3\right)  ,\ \left(  5\right)
\right)  ,
\]
since each element of $\left[  6\right]  $ appears exactly once in the
composite list $\left(  1,4,6,2,3,5\right)  $, and because $\sigma
=\operatorname*{cyc}\nolimits_{1,4,6}\circ\operatorname*{cyc}\nolimits_{2,3}%
\circ\operatorname*{cyc}\nolimits_{5}$. Thus, the cycles of $\sigma$ are
$\left(  1,4,6\right)  $, $\left(  2,3\right)  $ and $\left(  5\right)  $. (Do
not forget the length-$1$ cycle $\left(  5\right)  $. Even though the
permutation $\operatorname*{cyc}\nolimits_{5}$ is just the identity, the cycle
$\left(  5\right)  $ needs to be listed, since otherwise $5$ would not appear
in the composite list.)

As we said, the DCD of $\sigma$ is unique up to permuting the cycles and
rotating each cycle. Thus, another DCD of $\sigma$ is%
\[
\left(  \left(  2,3\right)  ,\ \left(  5\right)  ,\ \left(  6,1,4\right)
\right)  .
\]

\end{example}

We stress that the cycles of $\sigma$ are determined only up to cyclic
rotation (e.g., we cannot distinguish between $\left(  1,4,6\right)  $ and
$\left(  6,1,4\right)  $, although we can distinguish them both from $\left(
1,6,4\right)  $); but their lengths are uniquely determined by $\sigma$ (as a
multiset, since the order in which the cycles appear in the DCD is also not fixed).

Obviously, a permutation $\sigma$ is uniquely determined by its cycles (unlike
by its orbits). The orbits of $\sigma$ are just a shadow of the cycles:
Indeed, they are just the cycles of $\sigma$ transformed into sets (i.e., the
sets of entries of all cycles of $\sigma$). Let us state this rigorously:

\begin{proposition}
\label{prop.dcd.orbs-from-cycs}Let $X$ be a finite set. Let $\sigma\in S_{X}$
be any permutation. Let%
\begin{align}
&  \Big(\left(  a_{1,1},a_{1,2},\ldots,a_{1,n_{1}}\right)  ,\nonumber\\
&  \ \ \ \left(  a_{2,1},a_{2,2},\ldots,a_{2,n_{2}}\right)  ,\nonumber\\
&  \ \ \ \ldots,\nonumber\\
&  \ \ \ \left(  a_{k,1},a_{k,2},\ldots,a_{k,n_{k}}\right)
\Big) \label{eq.prop.dcd.orbs-from-cycs.dcd}%
\end{align}
be a DCD of $\sigma$. Then, the orbits of $\sigma$ are%
\begin{align*}
&  \left\{  a_{1,1},a_{1,2},\ldots,a_{1,n_{1}}\right\}  ,\\
&  \left\{  a_{2,1},a_{2,2},\ldots,a_{2,n_{2}}\right\}  ,\\
&  \ldots,\\
&  \left\{  a_{k,1},a_{k,2},\ldots,a_{k,n_{k}}\right\}  .
\end{align*}

\end{proposition}

\begin{proof}
[Proof of Proposition \ref{prop.dcd.orbs-from-cycs}.]We shall use the notation
\textquotedblleft$x\overset{f}{\mapsto}y$\textquotedblright\ for
\textquotedblleft$f\left(  x\right)  =y$\textquotedblright, where $f$ is a map
and where $x$ is some value that $f$ can be applied to.

All the elements $a_{u,v}$ (for $u\in\left[  k\right]  $ and $v\in\left[
n_{u}\right]  $) are distinct (because of the words \textquotedblleft exactly
once\textquotedblright\ in Definition \ref{def.dcd.dcd}). Moreover, each
element of $X$ has the form $a_{u,v}$ for some $u\in\left[  k\right]  $ and
$v\in\left[  n_{u}\right]  $ (for the same reason).

For each $i\in\left[  k\right]  $, let us set%
\[
c_{i}:=\operatorname*{cyc}\nolimits_{a_{i,1},a_{i,2},\ldots,a_{i,n_{i}}%
}\ \ \ \ \ \ \ \ \ \ \text{and}\ \ \ \ \ \ \ \ \ \ a_{i,n_{i}+1}:=a_{i,1}.
\]
Thus, for each $i\in\left[  k\right]  $ and each $j\in\left[  n_{i}\right]  $,
we have%
\begin{equation}
c_{i}\left(  a_{i,j}\right)  =a_{i,j+1}. \label{pf.prop.dcd.orbs-from-cycs.1}%
\end{equation}

Moreover, for each $i\in\left[  k\right]  $ and $r\in\left[  k\right]
\setminus\left\{  i\right\}  $ and $j\in\left[  n_{r}+1\right]  $, we have%
\begin{equation}
c_{i}\left(  a_{r,j}\right)  =a_{r,j} \label{pf.prop.dcd.orbs-from-cycs.2}%
\end{equation}
(since $a_{r,j}$ is one of the elements $a_{r,1},a_{r,2},\ldots,a_{r,n_{r}}$
(because $a_{r,n_{r}+1}=a_{r,1}$), and thus is \textbf{not} one of the
elements $a_{i,1},a_{i,2},\ldots,a_{i,n_{i}}$ (because all the elements
$a_{u,v}$ are distinct), and therefore is fixed by the cycle $c_{i}%
=\operatorname*{cyc}\nolimits_{a_{i,1},a_{i,2},\ldots,a_{i,n_{i}}}$).

We assumed that the list-of-lists (\ref{eq.prop.dcd.orbs-from-cycs.dcd}) is a
DCD of $\sigma$. Thus,%
\begin{align*}
\sigma &  =\underbrace{\operatorname*{cyc}\nolimits_{a_{1,1},a_{1,2}%
,\ldots,a_{1,n_{1}}}}_{\substack{=c_{1}\\\text{(by the definition of }%
c_{1}\text{)}}}\circ\underbrace{\operatorname*{cyc}\nolimits_{a_{2,1}%
,a_{2,2},\ldots,a_{2,n_{2}}}}_{\substack{=c_{2}\\\text{(by the definition of
}c_{2}\text{)}}}\circ\cdots\circ\underbrace{\operatorname*{cyc}%
\nolimits_{a_{k,1},a_{k,2},\ldots,a_{k,n_{k}}}}_{\substack{=c_{k}\\\text{(by
the definition of }c_{k}\text{)}}}\\
&  \ \ \ \ \ \ \ \ \ \ \ \ \ \ \ \ \ \ \ \ \left(  \text{by the definition of
a DCD}\right) \\
&  =c_{1}\circ c_{2}\circ\cdots\circ c_{k}.
\end{align*}

Now, let $i\in\left[  k\right]  $. Furthermore, let $j\in\left[  n_{i}\right]
$. Recall that $\sigma=c_{1}\circ c_{2}\circ\cdots\circ c_{k}$. Thus, in order
to compute $\sigma\left(  a_{i,j}\right)  $, we just need to apply the
permutations $c_{k},c_{k-1},\ldots,c_{1}$ to $a_{i,j}$ successively in this
order (i.e., to first apply $c_{k}$ to $a_{i,j}$, then apply $c_{k-1}$ to the
result, then apply $c_{k-2}$ to the result, and so on). Let us do this:%
\begin{align*}
a_{i,j}  &  \overset{c_{k}}{\mapsto}a_{i,j}\ \ \ \ \ \ \ \ \ \ \left(
\text{by (\ref{pf.prop.dcd.orbs-from-cycs.2})}\right) \\
&  \overset{c_{k-1}}{\mapsto}a_{i,j}\ \ \ \ \ \ \ \ \ \ \left(  \text{by
(\ref{pf.prop.dcd.orbs-from-cycs.2})}\right) \\
&  \overset{c_{k-2}}{\mapsto}\cdots\\
&  \overset{c_{i+1}}{\mapsto}a_{i,j}\ \ \ \ \ \ \ \ \ \ \left(  \text{by
(\ref{pf.prop.dcd.orbs-from-cycs.2})}\right) \\
&  \overset{c_{i}}{\mapsto}a_{i,j+1}\ \ \ \ \ \ \ \ \ \ \left(  \text{by
(\ref{pf.prop.dcd.orbs-from-cycs.1})}\right) \\
&  \overset{c_{i-1}}{\mapsto}a_{i,j+1}\ \ \ \ \ \ \ \ \ \ \left(  \text{by
(\ref{pf.prop.dcd.orbs-from-cycs.2})}\right) \\
&  \overset{c_{i-2}}{\mapsto}a_{i,j+1}\ \ \ \ \ \ \ \ \ \ \left(  \text{by
(\ref{pf.prop.dcd.orbs-from-cycs.2})}\right) \\
&  \overset{c_{i-3}}{\mapsto}\cdots\\
&  \overset{c_{1}}{\mapsto}a_{i,j+1}\ \ \ \ \ \ \ \ \ \ \left(  \text{by
(\ref{pf.prop.dcd.orbs-from-cycs.2})}\right)  .
\end{align*}
The final result of this procedure is $a_{i,j+1}$. Thus, $\sigma\left(
a_{i,j}\right)  =a_{i,j+1}$. In other words, $a_{i,j}\overset{\sigma}{\mapsto
}a_{i,j+1}$.

Forget that we fixed $j$. We thus have shown that $a_{i,j}\overset{\sigma
}{\mapsto}a_{i,j+1}$ for each $j\in\left[  n_{i}\right]  $. Hence,%
\[
a_{i,1}\overset{\sigma}{\mapsto}a_{i,2}\overset{\sigma}{\mapsto}%
a_{i,3}\overset{\sigma}{\mapsto}\cdots\overset{\sigma}{\mapsto}a_{i,n_{i}%
}\overset{\sigma}{\mapsto}a_{i,n_{i}+1}=a_{i,1}.
\]
This shows that the elements obtained from $a_{i,1}$ by repeated application
of $\sigma$ (that is, the elements of the form $\sigma^{k}\left(
a_{i,1}\right)  $ for $k\in\mathbb{N}$) are precisely $a_{i,1},a_{i,2}%
,\ldots,a_{i,n_{i}}$. Thus, the set $\left\{  a_{i,1},a_{i,2},\ldots
,a_{i,n_{i}}\right\}  $ is an orbit of $\sigma$.

Forget that we fixed $i$. We thus have shown that $\left\{  a_{i,1}%
,a_{i,2},\ldots,a_{i,n_{i}}\right\}  $ is an orbit of $\sigma$ for each
$i\in\left[  k\right]  $. In other words, all the $k$ sets%
\begin{align*}
&  \left\{  a_{1,1},a_{1,2},\ldots,a_{1,n_{1}}\right\}  ,\\
&  \left\{  a_{2,1},a_{2,2},\ldots,a_{2,n_{2}}\right\}  ,\\
&  \ldots,\\
&  \left\{  a_{k,1},a_{k,2},\ldots,a_{k,n_{k}}\right\}
\end{align*}
are orbits of $\sigma$. Since these $k$ sets include all elements of $X$
(because each element of $X$ has the form $a_{u,v}$ for some $u\in\left[
k\right]  $ and $v\in\left[  n_{u}\right]  $), it thus follows that $\sigma$
cannot have any further orbits; thus, the orbits of $\sigma$ are precisely the
$k$ sets that we found. Proposition \ref{prop.dcd.orbs-from-cycs} is now proved.
\end{proof}

\begin{noncompile}
\begin{proof}
[Old proof sketch for Proposition \ref{prop.dcd.orbs-from-cycs}.]Let
$i\in\left[  k\right]  $. We have
\[
\sigma=\operatorname*{cyc}\nolimits_{a_{1,1},a_{1,2},\ldots,a_{1,n_{1}}}%
\circ\operatorname*{cyc}\nolimits_{a_{2,1},a_{2,2},\ldots,a_{2,n_{2}}}%
\circ\cdots\circ\operatorname*{cyc}\nolimits_{a_{k,1},a_{k,2},\ldots
,a_{k,n_{k}}}%
\]
(by the definition of a DCD). Thus, applying $\sigma$ to an element $x\in X$
is tantamount to the following $k$-step procedure: First apply
$\operatorname*{cyc}\nolimits_{a_{k,1},a_{k,2},\ldots,a_{k,n_{k}}}$ to $x$;
then apply $\operatorname*{cyc}\nolimits_{a_{k-1,1},a_{k-1,2},\ldots
,a_{k-1,n_{k-1}}}$ to the result; then apply $\operatorname*{cyc}%
\nolimits_{a_{k-2,1},a_{k-2,2},\ldots,a_{k-2,n_{k-2}}}$ to the result; and so
on. When $x=a_{i,1}$, the first $k-i$ steps of this $k$-step procedure leave
$x$ unchanged (since $x=a_{i,1}\neq a_{u,v}$ for all $u>i$). The next step
(i.e., the $\left(  k-i+1\right)  $-th step) applies $\operatorname*{cyc}%
\nolimits_{a_{i,1},a_{i,2},\ldots,a_{i,n_{i}}}$, and thus sends $x$ to
$\operatorname*{cyc}\nolimits_{a_{i,1},a_{i,2},\ldots,a_{i,n_{i}}}\left(
x\right)  =a_{i,2}$ (since $x=a_{i,1}$). The remaining $i-1$ steps again leave
$a_{i,2}$ unchanged (since $a_{i,2}\neq a_{u,v}$ for all $u<i$). Thus,
$\sigma$ sends $a_{i,1}$ to $a_{i,2}$. In other words, $a_{i,1}\overset{\sigma
}{\mapsto}a_{i,2}$ (where the notation $x\overset{\sigma}{\mapsto}y$ means
\textquotedblleft$\sigma\left(  x\right)  =y$\textquotedblright). Similarly,
$a_{i,2}\overset{\sigma}{\mapsto}a_{i,3}$ and $a_{i,3}\overset{\sigma
}{\mapsto}a_{i,4}$ and so on, and finally $a_{i,n_{i}}\overset{\sigma
}{\mapsto}a_{i,1}$. Hence,%
\[
a_{i,1}\overset{\sigma}{\mapsto}a_{i,2}\overset{\sigma}{\mapsto}%
a_{i,3}\overset{\sigma}{\mapsto}\cdots\overset{\sigma}{\mapsto}a_{i,n_{i}%
}\overset{\sigma}{\mapsto}a_{i,1}.
\]
This shows that the elements of the form $\sigma^{k}\left(  a_{i,1}\right)  $
for $k\in\mathbb{N}$ are precisely $a_{i,1},a_{i,2},\ldots,a_{i,n_{i}}$. Thus,
the set $\left\{  a_{i,1},a_{i,2},\ldots,a_{i,n_{i}}\right\}  $ is an orbit of
$\sigma$.

Forget that we fixed $i$. We thus have shown that $\left\{  a_{i,1}%
,a_{i,2},\ldots,a_{i,n_{i}}\right\}  $ is an orbit of $\sigma$ for each
$i\in\left[  k\right]  $. In other words, all the $k$ sets%
\begin{align*}
&  \left\{  a_{1,1},a_{1,2},\ldots,a_{1,n_{1}}\right\}  ,\\
&  \left\{  a_{2,1},a_{2,2},\ldots,a_{2,n_{2}}\right\}  ,\\
&  \ldots,\\
&  \left\{  a_{k,1},a_{k,2},\ldots,a_{k,n_{k}}\right\}
\end{align*}
are orbits of $\sigma$. Since these $k$ sets include all elements of $X$ (by
the definition of a DCD), it thus follows that $\sigma$ cannot have any
further orbits; thus, the orbits of $\sigma$ are precisely the $k$ sets that
we found. Proposition \ref{prop.dcd.orbs-from-cycs} is now proved.
\end{proof}
\end{noncompile}

\begin{proposition}
\label{prop.dcd.inverse}Let $X$ be a finite set. Let $\sigma\in S_{X}$ be any
permutation. Then: \medskip

\textbf{(a)} The cycles of $\sigma^{-1}$ are the cycles of $\sigma$, read
backwards. (In other words, if $\left(  a_{1},a_{2},\ldots,a_{k}\right)  $ is
a cycle of $\sigma$, then $\left(  a_{k},a_{k-1},\ldots,a_{1}\right)  $ is a
cycle of $\sigma^{-1}$, and vice versa.) \medskip

\textbf{(b)} The orbits of $\sigma^{-1}$ are the orbits of $\sigma$.
\end{proposition}

\begin{proof}
\textbf{(a)} Let (\ref{eq.def.dcd.dcd.dcd}) be a DCD of $\sigma$. Then, each
element of $X$ appears exactly once in the composite list
(\ref{eq.def.dcd.dcd.comp}), and the equality
\[
\sigma=\operatorname*{cyc}\nolimits_{a_{1,1},a_{1,2},\ldots,a_{1,n_{1}}}%
\circ\operatorname*{cyc}\nolimits_{a_{2,1},a_{2,2},\ldots,a_{2,n_{2}}}%
\circ\cdots\circ\operatorname*{cyc}\nolimits_{a_{k,1},a_{k,2},\ldots
,a_{k,n_{k}}}%
\]
holds (by the definition of a DCD). From the latter equality, we obtain%
\begin{align*}
\sigma^{-1}  &  =\left(  \operatorname*{cyc}\nolimits_{a_{1,1},a_{1,2}%
,\ldots,a_{1,n_{1}}}\circ\operatorname*{cyc}\nolimits_{a_{2,1},a_{2,2}%
,\ldots,a_{2,n_{2}}}\circ\cdots\circ\operatorname*{cyc}\nolimits_{a_{k,1}%
,a_{k,2},\ldots,a_{k,n_{k}}}\right)  ^{-1}\\
&  =\underbrace{\operatorname*{cyc}\nolimits_{a_{k,1},a_{k,2},\ldots
,a_{k,n_{k}}}^{-1}}_{\substack{=\operatorname*{cyc}\nolimits_{a_{k,n_{k}%
},a_{k,n_{k}-1},\ldots,a_{k,1}}\\\text{(by (\ref{eq.intro.perms.cycs.inverse}%
))}}}\circ\underbrace{\operatorname*{cyc}\nolimits_{a_{k-1,1},a_{k-1,2}%
,\ldots,a_{k-1,n_{k-1}}}^{-1}}_{\substack{=\operatorname*{cyc}%
\nolimits_{a_{k-1,n_{k-1}},a_{k-1,n_{k-1}-1},\ldots,a_{k-1,1}}\\\text{(by
(\ref{eq.intro.perms.cycs.inverse}))}}}\circ\cdots\circ
\underbrace{\operatorname*{cyc}\nolimits_{a_{1,1},a_{1,2},\ldots,a_{1,n_{1}}%
}^{-1}}_{\substack{=\operatorname*{cyc}\nolimits_{a_{1,n_{1}},a_{1,n_{1}%
-1},\ldots,a_{1,1}}\\\text{(by (\ref{eq.intro.perms.cycs.inverse}))}}}\\
&  =\operatorname*{cyc}\nolimits_{a_{k,n_{k}},a_{k,n_{k}-1},\ldots,a_{k,1}%
}\circ\operatorname*{cyc}\nolimits_{a_{k-1,n_{k-1}},a_{k-1,n_{k-1}-1}%
,\ldots,a_{k-1,1}}\circ\cdots\circ\operatorname*{cyc}\nolimits_{a_{1,n_{1}%
},a_{1,n_{1}-1},\ldots,a_{1,1}}.
\end{align*}
This shows that the nested list
\begin{align*}
&  \Big(\left(  a_{k,n_{k}},a_{k,n_{k}-1},\ldots,a_{k,1}\right)  ,\\
&  \ \ \ \left(  a_{k-1,n_{k-1}},a_{k-1,n_{k-1}-1},\ldots,a_{k-1,1}\right)
,\\
&  \ \ \ \ldots,\\
&  \ \ \ \left(  a_{1,n_{1}},a_{1,n_{1}-1},\ldots,a_{1,1}\right)  \Big)
\end{align*}
is a DCD of $\sigma^{-1}$ (since each element of $X$ appears exactly once in
the composite list%
\begin{align*}
&  (a_{k,n_{k}},a_{k,n_{k}-1},\ldots,a_{k,1},\\
&  \ \ \ a_{k-1,n_{k-1}},a_{k-1,n_{k-1}-1},\ldots,a_{k-1,1},\\
&  \ \ \ \ldots,\\
&  \ \ \ a_{1,n_{1}},a_{1,n_{1}-1},\ldots,a_{1,1})
\end{align*}
(because each element of $X$ appears exactly once in the composite list
(\ref{eq.def.dcd.dcd.comp}))). Thus, the cycles of $\sigma^{-1}$ are%
\begin{align*}
&  \left(  a_{k,n_{k}},a_{k,n_{k}-1},\ldots,a_{k,1}\right)  ,\\
&  \left(  a_{k-1,n_{k-1}},a_{k-1,n_{k-1}-1},\ldots,a_{k-1,1}\right)  ,\\
&  \ldots,\\
&  \left(  a_{1,n_{1}},a_{1,n_{1}-1},\ldots,a_{1,1}\right)  .
\end{align*}
On the other hand, the cycles of $\sigma$ are%
\begin{align*}
&  \left(  a_{1,1},a_{1,2},\ldots,a_{1,n_{1}}\right)  ,\\
&  \left(  a_{2,1},a_{2,2},\ldots,a_{2,n_{2}}\right)  ,\\
&  \ldots,\\
&  \left(  a_{k,1},a_{k,2},\ldots,a_{k,n_{k}}\right)
\end{align*}
(since (\ref{eq.def.dcd.dcd.dcd}) is a DCD of $\sigma$). Comparing these, we
see that the cycles of $\sigma^{-1}$ are the cycles of $\sigma$, read
backwards (and listed in reverse order, but this does not matter, since the
different cycles can be permuted at will anyway). This proves Proposition
\ref{prop.dcd.inverse} \textbf{(a)}. \medskip

\textbf{(b)} This follows easily from part \textbf{(a)} using Proposition
\ref{prop.dcd.orbs-from-cycs}. But let us give a more direct proof:

Recall the binary relation $\overset{\sigma}{\sim}$ (from Definition
\ref{def.orbits.orbits}) and the analogously defined relation $\overset{\sigma
^{-1}}{\sim}$. As we know, the relation $\overset{\sigma}{\sim}$ is an
equivalence relation, thus is symmetric. Hence, for any $i,j\in X$, we have
the following chain of equivalences:%
\begin{align*}
\left(  i\overset{\sigma}{\sim}j\right)  \  &  \Longleftrightarrow\ \left(
j\overset{\sigma}{\sim}i\right) \\
&  \Longleftrightarrow\ \left(  j=\sigma^{k}\left(  i\right)  \text{ for some
}k\in\mathbb{N}\right) \\
&  \ \ \ \ \ \ \ \ \ \ \ \ \ \ \ \ \ \ \ \ \left(  \text{by the definition of
the relation }\overset{\sigma}{\sim}\right) \\
&  \Longleftrightarrow\ \left(  i=\left(  \sigma^{k}\right)  ^{-1}\left(
j\right)  \text{ for some }k\in\mathbb{N}\right) \\
&  \ \ \ \ \ \ \ \ \ \ \ \ \ \ \ \ \ \ \ \ \left(
\begin{array}
[c]{c}%
\text{because the equality }j=\sigma^{k}\left(  i\right)  \text{ can be
equivalently}\\
\text{rewritten as }i=\left(  \sigma^{k}\right)  ^{-1}\left(  j\right)
\end{array}
\right) \\
&  \Longleftrightarrow\ \left(  i=\left(  \sigma^{-1}\right)  ^{k}\left(
j\right)  \text{ for some }k\in\mathbb{N}\right) \\
&  \ \ \ \ \ \ \ \ \ \ \ \ \ \ \ \ \ \ \ \ \left(  \text{since }\left(
\sigma^{k}\right)  ^{-1}=\left(  \sigma^{-1}\right)  ^{k}\text{ for any }%
k\in\mathbb{N}\right) \\
&  \Longleftrightarrow\ \left(  i\overset{\sigma^{-1}}{\sim}j\right)
\ \ \ \ \ \ \ \ \ \ \left(  \text{by the definition of the relation
}\overset{\sigma^{-1}}{\sim}\right)  .
\end{align*}
In other words, the relation $\overset{\sigma}{\sim}$ and the relation
$\overset{\sigma^{-1}}{\sim}$ are the same. Hence, the equivalence classes of
the latter relation are the equivalence classes of the former. In other words,
the orbits of $\sigma^{-1}$ are the orbits of $\sigma$ (since the orbits of
$\sigma^{-1}$ are defined as the equivalence classes of $\overset{\sigma
^{-1}}{\sim}$, while the orbits of $\sigma$ are defined as the equivalence
classes of $\overset{\sigma}{\sim}$). This proves Proposition
\ref{prop.dcd.inverse} \textbf{(b)}.
\end{proof}

\begin{definition}
\label{def.rl.rl}Let $X$ be a finite set. Let $\sigma\in S_{X}$ be any
permutation. The \emph{reflection length} of $\sigma$ is defined to be the
number%
\[
\operatorname*{rl}\sigma:=\left\vert X\right\vert -\left(  \text{\# of orbits
of }\sigma\right)  .
\]
This is a nonnegative integer, since $\left(  \text{\# of orbits of }%
\sigma\right)  \leq\left\vert X\right\vert $. It is no larger than $\left\vert
X\right\vert -1$ (unless $X=\varnothing$), since $\sigma$ has at least one orbit.
\end{definition}

The reflection length is closely connected to transpositions, as the following
theorem (\cite[Lecture 28, Theorem 4.3.23 \textbf{(b)}]{22fco},
\cite[\textquotedblleft The sign of a permutation\textquotedblright, Theorem
6.1]{Conrad}) shows:

\begin{theorem}
\label{thm.rl.transps}Let $X$ be a finite set. Let $\sigma\in S_{X}$ be any
permutation. Then, $\sigma$ can be written as a product of $\operatorname*{rl}%
\sigma$ many transpositions in $S_{X}$. Moreover, $\operatorname*{rl}\sigma$
is the smallest number $p\in\mathbb{N}$ such that $\sigma$ can be written as a
product of $p$ transpositions in $S_{X}$.
\end{theorem}

Hence the term \textquotedblleft reflection length\textquotedblright.
(\textquotedblleft Reflection\textquotedblright\ is how geometers refer to transpositions.)

The proof of Theorem \ref{thm.rl.transps} boils down to the following
important lemma:

\begin{lemma}
\label{lem.rl.change}Let $X$ be a finite set. Let $\sigma\in S_{X}$ be any
permutation. Let $i$ and $j$ be two distinct elements of $X$. Let
$\tau=t_{i,j}\circ\sigma$. Then: \medskip

\textbf{(a)} If $i\overset{\sigma}{\sim}j$, then the permutation $\tau$ has
$1$ more orbit than $\sigma$. \medskip

\textbf{(b)} If $i\overset{\sigma}{\sim}j$, then we don't have $i\overset{\tau
}{\sim}j$. \medskip

\textbf{(c)} If we don't have $i\overset{\sigma}{\sim}j$, then the permutation
$\tau$ has $1$ fewer orbits than $\sigma$. \medskip

\textbf{(d)} If we don't have $i\overset{\sigma}{\sim}j$, then we do have
$i\overset{\tau}{\sim}j$.
\end{lemma}

\begin{proof}
See \cite[Lecture 27, Lemma 4.3.13]{22fco}.
\end{proof}

\subsubsection{Nonstarters and starters}

Next, we introduce some nonstandard but useful notation:

\begin{definition}
\label{def.orbits.nost}Let $\sigma\in S_{n}$.

An element $i\in\left[  n\right]  $ is called a \emph{nonstarter} of $\sigma$
if the orbit of $\sigma$ that contains $i$ also contains some number smaller
than $i$ (that is, if $i$ is \textbf{not} the smallest element of the orbit of
$\sigma$ that contains $i$). Otherwise, $i$ is called a \emph{starter} of
$\sigma$. (Thus, $i$ is called a starter of $\sigma$ if and only if $i$ is the
smallest element of the orbit of $\sigma$ that contains $i$.)

We let $\operatorname*{NoSt}\sigma$ denote the set of all nonstarters of
$\sigma$.
\end{definition}

\begin{example}
The nonstarters of $\operatorname*{oln}\left(  432651\right)  \in S_{6}$ are
$3$, $4$ and $6$. The starters of this permutation are thus $1$, $2$ and $5$.
\end{example}

The following is easy to see:

\begin{remark}
\label{rmk.orbits.nost.easy}Let $\sigma\in S_{n}$ be a permutation. Then:
\medskip

\textbf{(a)} Any fixed point of $\sigma$ is a starter of $\sigma$. \medskip

\textbf{(b)} The number $1$ is always a starter of $\sigma$, as long as
$n\geq1$. \medskip

\textbf{(c)} If $\sigma$ has no nonstarters, then $\sigma=\operatorname*{id}$.
\medskip

\textbf{(d)} We have $\left(  \text{\# of starters of }\sigma\right)  =\left(
\text{\# of orbits of }\sigma\right)  $. \medskip

\textbf{(e)} We have $\left\vert \operatorname*{NoSt}\sigma\right\vert
=\operatorname*{rl}\sigma$.
\end{remark}

\begin{proof}
[Proof of Remark \ref{rmk.orbits.nost.easy}.]\textbf{(a)} Let $i\in\left[
n\right]  $ be a fixed point of $\sigma$. Then, $\sigma\left(  i\right)  =i$,
so that $\sigma^{k}\left(  i\right)  =i$ for each $k\in\mathbb{N}$. Thus, the
orbit of $\sigma$ that contains $i$ is just the one-element set $\left\{
i\right\}  $. Hence, $i$ is the smallest element of this orbit, and thus a
starter of $\sigma$ (by the definition of a starter). So we have shown that
any fixed point of $\sigma$ is a starter of $\sigma$. This proves Remark
\ref{rmk.orbits.nost.easy} \textbf{(a)}. \medskip

\textbf{(b)} Assume that $n\geq1$. Then, $1$ is the smallest element of the
orbit of $\sigma$ that contains $1$ (since $1$ is the smallest element of
$\left[  n\right]  $). Hence, $1$ is a starter of $\sigma$. This proves Remark
\ref{rmk.orbits.nost.easy} \textbf{(b)}. \medskip

\textbf{(c)} Assume that $\sigma$ has no nonstarters. Let $i\in\left[
n\right]  $. Then, $\sigma\left(  i\right)  \overset{\sigma}{\sim}i$ (by the
definition of $\overset{\sigma}{\sim}$), and thus the orbit of $\sigma$ that
contains $i$ must also contain $\sigma\left(  i\right)  $. If we had
$\sigma\left(  i\right)  \neq i$, then the orbit of $\sigma$ that contains $i$
would contain at least two elements (since it would contain the two distinct
elements $i$ and $\sigma\left(  i\right)  $), and thus it would contain at
least one element that is \textbf{not} the smallest element of this orbit.
This element would then be a nonstarter of $\sigma$ (by the definition of a
nonstarter), which would contradict the assumption that $\sigma$ has no
nonstarters. Thus, $\sigma\left(  i\right)  \neq i$ is impossible. Hence,
$\sigma\left(  i\right)  =i$. Since we have proved this for each $i\in\left[
n\right]  $, we thus conclude that $\sigma=\operatorname*{id}$. This proves
Remark \ref{rmk.orbits.nost.easy} \textbf{(c)}. \medskip

\textbf{(d)} A starter of $\sigma$ is defined to be an element of $\left[
n\right]  $ that is the smallest element of the orbit of $\sigma$ that
contains it. Thus, each orbit of $\sigma$ contains exactly one starter of
$\sigma$ (namely, its smallest element). Hence, the permutation $\sigma$ has
as many starters as it has orbits. In other words, $\left(  \text{\# of
starters of }\sigma\right)  =\left(  \text{\# of orbits of }\sigma\right)  $.
This proves Remark \ref{rmk.orbits.nost.easy} \textbf{(d)}. \medskip

\textbf{(e)} Recall that $\operatorname*{NoSt}\sigma$ is defined as the set of
all nonstarters of $\sigma$. Thus,%
\begin{align*}
\left\vert \operatorname*{NoSt}\sigma\right\vert  &  =\left(  \text{\# of
nonstarters of }\sigma\right)  =\left\vert \left[  n\right]  \right\vert
-\underbrace{\left(  \text{\# of starters of }\sigma\right)  }%
_{\substack{=\left(  \text{\# of orbits of }\sigma\right)  \\\text{(by Remark
\ref{rmk.orbits.nost.easy} \textbf{(d)})}}}\\
&  \ \ \ \ \ \ \ \ \ \ \ \ \ \ \ \ \ \ \ \ \left(
\begin{array}
[c]{c}%
\text{since each }i\in\left[  n\right]  \text{ is either a starter of }%
\sigma\\
\text{or a nonstarter of }\sigma\text{, but not both}%
\end{array}
\right) \\
&  =\left\vert \left[  n\right]  \right\vert -\left(  \text{\# of orbits of
}\sigma\right)  =\operatorname*{rl}\sigma
\end{align*}
(by Definition \ref{def.rl.rl}). This proves Remark \ref{rmk.orbits.nost.easy}
\textbf{(e)}.
\end{proof}

\subsubsection{Products of distinct jucys--murphies}

We can now give an explicit formula for an arbitrary product of distinct jucys--murphies:

\begin{theorem}
\label{thm.YJM.some-prod}Let $i_{1},i_{2},\ldots,i_{k}$ be $k$ distinct
elements of $\left[  n\right]  $. Then,%
\[
\mathbf{m}_{i_{1}}\mathbf{m}_{i_{2}}\cdots\mathbf{m}_{i_{k}}=\sum
_{\substack{w\in S_{n};\\\operatorname*{NoSt}w=\left\{  i_{1},i_{2}%
,\ldots,i_{k}\right\}  }}w.
\]

\end{theorem}

\begin{remark}
If $i_{1}=1$, then both sides of this equality are $0$, since $\mathbf{m}%
_{1}=0$ and since $1$ is never a nonstarter of a permutation (by Remark
\ref{rmk.orbits.nost.easy} \textbf{(b)}).
\end{remark}

\begin{example}
Let $n=5$. Then, Theorem \ref{thm.YJM.some-prod} yields%
\[
\mathbf{m}_{2}\mathbf{m}_{4}=\sum_{\substack{w\in S_{5};\\\operatorname*{NoSt}%
w=\left\{  2,4\right\}  }}w.
\]
And indeed, it is not hard to see that both sides equal $\operatorname*{cyc}%
\nolimits_{1,4,2}+\operatorname*{cyc}\nolimits_{1,2,4}+\,t_{1,2}t_{3,4}$.
\end{example}

To prove Theorem \ref{thm.YJM.some-prod}, we need the following two lemmas:

\begin{lemma}
\label{lem.YJM.some-prod.1}Let $w\in S_{n}$. Then: \medskip

\textbf{(a)} If $p\in\left[  n\right]  $ is larger than each nonstarter of
$w$, then $w\left(  p\right)  =p$. \medskip

\textbf{(b)} If $q$ is the largest nonstarter of $w$, then $w\left(  q\right)
\in\left[  q-1\right]  $.
\end{lemma}

\begin{proof}
[Proof of Lemma \ref{lem.YJM.some-prod.1}.]In the following, the word
\textquotedblleft nonstarter\textquotedblright\ will always mean
\textquotedblleft nonstarter of $w$\textquotedblright. Likewise, the words
\textquotedblleft starter\textquotedblright\ and \textquotedblleft
orbit\textquotedblright\ will mean \textquotedblleft starter of $w$%
\textquotedblright\ or \textquotedblleft orbit of $w$\textquotedblright.
\medskip

\textbf{(a)} Let $p\in\left[  n\right]  $ be larger than each nonstarter of
$w$. Then, $p$ cannot be a nonstarter. In other words, $p$ is a starter. In
yet other words, $p$ is the smallest element of the orbit that contains $p$
(by the definition of \textquotedblleft starter\textquotedblright).

But $w\left(  p\right)  \overset{w}{\sim}p$ (by the definition of the relation
$\overset{w}{\sim}$), and therefore $w\left(  p\right)  $ belongs to the orbit
that contains $p$. Since $p$ is the smallest element of this orbit, we thus
obtain $w\left(  p\right)  \geq p$. In other words, $p$ is not larger than
$w\left(  p\right)  $.

Hence, $w\left(  p\right)  $ is not a nonstarter (since $p$ is larger than
each nonstarter, but not larger than $w\left(  p\right)  $). In other words,
$w\left(  p\right)  $ is a starter.

But recall that a starter is an element that is the smallest element of the
orbit that contains it. Hence, each orbit contains exactly one starter.
Therefore, any two starters belonging to the same orbit must be equal. Since
$w\left(  p\right)  $ and $p$ are two starters belonging to the same orbit
(because $w\left(  p\right)  \overset{w}{\sim}p$), we thus conclude that
$w\left(  p\right)  =p$. This proves Lemma \ref{lem.YJM.some-prod.1}
\textbf{(a)}. \medskip

\textbf{(b)} First, we note that $w$ is a permutation and thus injective.

Let $q$ be the largest nonstarter of $w$. If we had $w\left(  q\right)  =q$,
then $q$ would be a fixed point of $w$ and therefore a starter (since Remark
\ref{rmk.orbits.nost.easy} \textbf{(a)} yields that any fixed point of $w$ is
a starter of $w$), which would contradict the fact that $q$ is a nonstarter.
Hence, we cannot have $w\left(  q\right)  =q$. Thus, $w\left(  q\right)  \neq
q$.

Let $p=w\left(  q\right)  $. If we had $p>q$, then $p$ would be larger than
each nonstarter of $w$ (since $q$ is the largest nonstarter of $w$), and thus
we would have $w\left(  p\right)  =p$ (by Lemma \ref{lem.YJM.some-prod.1}
\textbf{(a)}); but this would yield $w\left(  p\right)  =p=w\left(  q\right)
$ and therefore $p=q$ (since $w$ is injective), which would contradict $p>q$.
Hence, we cannot have $p>q$. Thus, we have $p\leq q$. In other words,
$w\left(  q\right)  \leq q$ (since $p=w\left(  q\right)  $). Combined with
$w\left(  q\right)  \neq q$, this yields $w\left(  q\right)  <q$. Hence,
$w\left(  q\right)  \leq q-1$, so that $w\left(  q\right)  \in\left[
q-1\right]  $. This proves Lemma \ref{lem.YJM.some-prod.1} \textbf{(b)}.
\end{proof}

\begin{lemma}
\label{lem.YJM.some-prod.2}Let $w\in S_{n}$. Let $p,q\in\left[  n\right]  $ be
two elements with $p<q$. Then: \medskip

\textbf{(a)} If $w\left(  q\right)  =q$, then $\operatorname*{NoSt}\left(
t_{p,q}w\right)  =\left(  \operatorname*{NoSt}w\right)  \cup\left\{
q\right\}  $. \medskip

\textbf{(b)} If $w\left(  q\right)  =p$, then $\operatorname*{NoSt}\left(
t_{p,q}w\right)  =\left(  \operatorname*{NoSt}w\right)  \setminus\left\{
q\right\}  $.
\end{lemma}

\begin{proof}
[Proof of Lemma \ref{lem.YJM.some-prod.2}.]\textbf{(a)} Assume that $w\left(
q\right)  =q$. Thus, $w^{k}\left(  q\right)  =q$ for each $k\in\mathbb{N}$.
Hence, the one-element set $\left\{  q\right\}  $ is an orbit of $w$. Let us
denote this orbit by $Q$. The permutation $w$ acts on this orbit $Q$ by
\begin{equation}
q\overset{w}{\mapsto}q \label{pf.lem.YJM.some-prod.2.a.orbQ}%
\end{equation}
(where the notation \textquotedblleft$x\overset{w}{\mapsto}y$%
\textquotedblright\ means \textquotedblleft$w\left(  x\right)  =y$%
\textquotedblright).

Let $P$ be the orbit of $w$ that contains $p$. Set $m=\left\vert P\right\vert
$. Then, by Remark \ref{rmk.orbits.O=} (applied to $X=\left[  n\right]  $ and
$\sigma=w$ and $O=P$ and $a=p$), we have $w^{m}\left(  p\right)  =p$ and%
\begin{equation}
P=\left\{  w^{0}\left(  p\right)  ,\ w^{1}\left(  p\right)  ,\ \ldots
,\ w^{m-1}\left(  p\right)  \right\}  . \label{pf.lem.YJM.some-prod.2.a.P=}%
\end{equation}
The permutation $w$ acts on the orbit $P$ by%
\[
w^{0}\left(  p\right)  \overset{w}{\mapsto}w^{1}\left(  p\right)
\overset{w}{\mapsto}w^{2}\left(  p\right)  \overset{w}{\mapsto}\cdots
\overset{w}{\mapsto}w^{m-1}\left(  p\right)  \overset{w}{\mapsto}w^{m}\left(
p\right)  =p=w^{0}\left(  p\right)  .
\]

Recall that $\left\vert P\right\vert =m$. In view of
(\ref{pf.lem.YJM.some-prod.2.a.P=}), we can rewrite this as
\[
\left\vert \left\{  w^{0}\left(  p\right)  ,\ w^{1}\left(  p\right)
,\ \ldots,\ w^{m-1}\left(  p\right)  \right\}  \right\vert =m.
\]
Hence, the $m$ elements $w^{0}\left(  p\right)  ,\ w^{1}\left(  p\right)
,\ \ldots,\ w^{m-1}\left(  p\right)  $ are distinct.

The two orbits $P$ and $Q$ are distinct\footnote{\textit{Proof.} From $p<q$,
we obtain $p\neq q$, so that $p\notin\left\{  q\right\}  =Q$. Consider the two
orbits $P$ and $Q$ of $w$. The orbit $P$ contains $p$ (by its definition),
whereas the orbit $Q$ does not (since $p\notin Q$). Thus, the two orbits $P$
and $Q$ are distinct, qed.}. Hence, none of the elements \newline$w^{0}\left(
p\right)  ,\ w^{1}\left(  p\right)  ,\ \ldots,\ w^{m-1}\left(  p\right)  $
equals $q$ (since they belong to a different orbit of $w$ than $q$
does\footnote{\textit{Proof.} The two orbits $P$ and $Q$ are distinct, and
therefore disjoint (since any two distinct orbits of $w$ are disjoint). Thus,
none of the elements $w^{0}\left(  p\right)  ,\ w^{1}\left(  p\right)
,\ \ldots,\ w^{m-1}\left(  p\right)  $ of $P$ belongs to $Q$. Hence, none of
the elements $w^{0}\left(  p\right)  ,\ w^{1}\left(  p\right)  ,\ \ldots
,\ w^{m-1}\left(  p\right)  $ of $P$ equals $q$ (since $q$ belongs to $Q$).}).
Thus, in particular, none of the elements $w^{0}\left(  p\right)
,\ w^{1}\left(  p\right)  ,\ \ldots,\ w^{m-2}\left(  p\right)  $ equals $q$.
Moreover, none of these elements $w^{0}\left(  p\right)  ,\ w^{1}\left(
p\right)  ,\ \ldots,\ w^{m-2}\left(  p\right)  $ equals $w^{m-1}\left(
p\right)  $ either (since the $m$ elements $w^{0}\left(  p\right)
,\ w^{1}\left(  p\right)  ,\ \ldots,\ w^{m-1}\left(  p\right)  $ are distinct).

Now, let $v$ be the permutation $t_{p,q}w$. This permutation $v$ is obtained
from $w$ by swapping the two values $p$ and $q$ (since it is the composition
of $t_{p,q}$ with $w$). Thus, $v$ sends $q$ not to $w\left(  q\right)  =q$ but
to $p$ instead, while sending $w^{m-1}\left(  p\right)  $ not to $w^{m}\left(
p\right)  =p$ but to $q$ instead; on all other inputs, $v$ acts in the exact
same way as $w$. Hence, the relations%
\[
w^{0}\left(  p\right)  \overset{w}{\mapsto}w^{1}\left(  p\right)
\overset{w}{\mapsto}w^{2}\left(  p\right)  \overset{w}{\mapsto}\cdots
\overset{w}{\mapsto}w^{m-1}\left(  p\right)
\]
(which we have proved above) become%
\[
w^{0}\left(  p\right)  \overset{v}{\mapsto}w^{1}\left(  p\right)
\overset{v}{\mapsto}w^{2}\left(  p\right)  \overset{v}{\mapsto}\cdots
\overset{v}{\mapsto}w^{m-1}\left(  p\right)
\]
(since none of the elements $w^{0}\left(  p\right)  ,\ w^{1}\left(  p\right)
,\ \ldots,\ w^{m-2}\left(  p\right)  $ equals $q$, and none of them equals
$w^{m-1}\left(  p\right)  $ either). Combining this chain of relations with
$w^{m-1}\left(  p\right)  \overset{v}{\mapsto}q$ (since $v$ sends
$w^{m-1}\left(  p\right)  $ to $q$) and with $q\overset{v}{\mapsto}%
w^{0}\left(  p\right)  $ (since $v$ sends $q$ to $p=w^{0}\left(  p\right)  $),
we obtain a single circular chain%
\[
q\overset{v}{\mapsto}w^{0}\left(  p\right)  \overset{v}{\mapsto}w^{1}\left(
p\right)  \overset{v}{\mapsto}w^{2}\left(  p\right)  \overset{v}{\mapsto
}\cdots\overset{v}{\mapsto}w^{m-1}\left(  p\right)  \overset{v}{\mapsto}q.
\]
Hence, the set $\left\{  q,\ w^{0}\left(  p\right)  ,\ w^{1}\left(  p\right)
,\ \ldots,\ w^{m-1}\left(  p\right)  \right\}  $ is an orbit of $v$. In other
words, $P\cup Q$ is an orbit of $v$ (since
\begin{align*}
\left\{  q,\ w^{0}\left(  p\right)  ,\ w^{1}\left(  p\right)  ,\ \ldots
,\ w^{m-1}\left(  p\right)  \right\}   &  =\underbrace{\left\{  q\right\}
}_{=Q}\cup\underbrace{\left\{  w^{0}\left(  p\right)  ,\ w^{1}\left(
p\right)  ,\ \ldots,\ w^{m-1}\left(  p\right)  \right\}  }_{=P}\\
&  =Q\cup P=P\cup Q
\end{align*}
). Thus, when we pass from the permutation $w$ to $v$, the two orbits $P$ and
$Q$ are merged into a single orbit $P\cup Q$. All the other orbits of $w$
remain orbits of $v$ (since $v$ differs from $w$ only in that the values $p$
and $q$ are swapped, and this swap clearly does not affect any orbits except
for the ones that contain $p$ or $q$). In other words, the orbits of $v$
differ from the orbits of $w$ only in that the two orbits $P$ and $Q$ are
merged into a single orbit $P\cup Q$.

Now recall that the starters of a permutation are the smallest elements of its
orbits (by the definition of a starter). Each orbit $O$ thus produces exactly
one starter, namely its smallest element $\min O$. In particular, the two
distinct orbits $P$ and $Q$ of $w$ thus produce two distinct starters $\min P$
and $\min Q$ of $w$. Hence, $\min P$ and $\min Q$ are two distinct starters of
$w$, that is, two distinct elements of the set $\left\{  \text{starters of
}w\right\}  $. Therefore,%
\[
\min P\in\left\{  \text{starters of }w\right\}  \setminus\left\{  \min
Q\right\}  .
\]

Recall that when we pass from the permutation $w$ to $v$, the two orbits $P$
and $Q$ are merged into a single orbit $P\cup Q$, while all other orbits stay
unchanged. Thus, when we pass from $w$ to $v$, all the starters stay
unchanged, except that the two starters $\min P$ and $\min Q$ coming from the
orbits $P$ and $Q$ get replaced by a single starter $\min\left(  P\cup
Q\right)  $ coming from the new orbit $P\cup Q$ (since the starters of a
permutation are the smallest elements of its orbits). Hence,%
\[
\left\{  \text{starters of }v\right\}  =\left(  \left\{  \text{starters of
}w\right\}  \setminus\left\{  \min P,\ \min Q\right\}  \right)  \cup\left\{
\min\left(  P\cup Q\right)  \right\}  .
\]

Now, recall that $Q=\left\{  q\right\}  $, so that $\min Q=q$. This shows that
$q$ is a starter of $w$ (since $\min Q$ is a starter of $w$). In other words,
$q\in\left\{  \text{starters of }w\right\}  $.

However, $p\in P$ (by the definition of $P$) and thus $\min P\leq p$ (since
the smallest element of a set is $\leq$ to any element of this set). Hence,
$\min P\leq p<q=\min Q$. Furthermore, recall that the minimum of the union of
two sets is the smaller one among their individual minima. Thus,
\[
\min\left(  P\cup Q\right)  =\min\left\{  \min P,\ \min Q\right\}  =\min P
\]
(since $\min P<\min Q$). Thus, our above computation becomes%
\begin{align*}
\left\{  \text{starters of }v\right\}   &  =\underbrace{\left(  \left\{
\text{starters of }w\right\}  \setminus\left\{  \min P,\ \min Q\right\}
\right)  }_{=\left(  \left\{  \text{starters of }w\right\}  \setminus\left\{
\min Q\right\}  \right)  \setminus\left\{  \min P\right\}  }\cup\left\{
\underbrace{\min\left(  P\cup Q\right)  }_{=\min P}\right\} \\
&  =\left(  \left(  \left\{  \text{starters of }w\right\}  \setminus\left\{
\min Q\right\}  \right)  \setminus\left\{  \min P\right\}  \right)
\cup\left\{  \min P\right\} \\
&  =\left\{  \text{starters of }w\right\}  \setminus\left\{  \min Q\right\}
\end{align*}
(since $\min P\in\left\{  \text{starters of }w\right\}  \setminus\left\{  \min
Q\right\}  $). Since $\min Q=q$, we can rewrite this as%
\[
\left\{  \text{starters of }v\right\}  =\left\{  \text{starters of }w\right\}
\setminus\left\{  q\right\}  .
\]

Now, the definition of $\operatorname*{NoSt}w$ yields%
\[
\operatorname*{NoSt}w=\left\{  \text{nonstarters of }w\right\}  =\left[
n\right]  \setminus\left\{  \text{starters of }w\right\}
\]
(since the nonstarters of $w$ are the elements $i\in\left[  n\right]  $ that
are not starters of $w$). Similarly,%
\begin{align*}
\operatorname*{NoSt}v  &  =\left[  n\right]  \setminus\underbrace{\left\{
\text{starters of }v\right\}  }_{=\left\{  \text{starters of }w\right\}
\setminus\left\{  q\right\}  }\\
&  =\left[  n\right]  \setminus\left(  \left\{  \text{starters of }w\right\}
\setminus\left\{  q\right\}  \right) \\
&  =\underbrace{\left(  \left[  n\right]  \setminus\left\{  \text{starters of
}w\right\}  \right)  }_{=\operatorname*{NoSt}w}\cup\left\{  q\right\}
\ \ \ \ \ \ \ \ \ \ \left(
\begin{array}
[c]{c}%
\text{since }q\in\left\{  \text{starters of }w\right\} \\
\text{and }\left\{  \text{starters of }w\right\}  \subseteq\left[  n\right]
\end{array}
\right) \\
&  =\left(  \operatorname*{NoSt}w\right)  \cup\left\{  q\right\}  .
\end{align*}
In other words, $\operatorname*{NoSt}\left(  t_{p,q}w\right)  =\left(
\operatorname*{NoSt}w\right)  \cup\left\{  q\right\}  $ (since $v=t_{p,q}w$).
This proves Lemma \ref{lem.YJM.some-prod.2} \textbf{(a)}. \medskip

\textbf{(b)} Assume that $w\left(  q\right)  =p$. Hence, $\left(
t_{p,q}w\right)  \left(  q\right)  =t_{p,q}\left(  \underbrace{w\left(
q\right)  }_{=p}\right)  =t_{p,q}\left(  p\right)  =q$. Hence, $q$ is a fixed
point of $t_{p,q}w$, and thus is a starter of $t_{p,q}w$ (since Remark
\ref{rmk.orbits.nost.easy} \textbf{(a)} yields that any fixed point of
$t_{p,q}w$ is a starter of $t_{p,q}w$). Thus, $q$ is not a nonstarter of
$t_{p,q}w$. In other words, $q\notin\operatorname*{NoSt}\left(  t_{p,q}%
w\right)  $.

We have $\left(  t_{p,q}w\right)  \left(  q\right)  =q$. Thus, Lemma
\ref{lem.YJM.some-prod.2} \textbf{(a)} (applied to $t_{p,q}w$ instead of $w$)
yields $\operatorname*{NoSt}\left(  t_{p,q}t_{p,q}w\right)  =\left(
\operatorname*{NoSt}\left(  t_{p,q}w\right)  \right)  \cup\left\{  q\right\}
$. In view of $\underbrace{t_{p,q}t_{p,q}}_{=t_{p,q}^{2}=\operatorname*{id}%
}w=w$, we can rewrite this as%
\[
\operatorname*{NoSt}w=\left(  \operatorname*{NoSt}\left(  t_{p,q}w\right)
\right)  \cup\left\{  q\right\}  .
\]
Hence,%
\[
\left(  \operatorname*{NoSt}w\right)  \setminus\left\{  q\right\}  =\left(
\left(  \operatorname*{NoSt}\left(  t_{p,q}w\right)  \right)  \cup\left\{
q\right\}  \right)  \setminus\left\{  q\right\}  =\operatorname*{NoSt}\left(
t_{p,q}w\right)
\]
(since $q\notin\operatorname*{NoSt}\left(  t_{p,q}w\right)  $). In other
words, $\operatorname*{NoSt}\left(  t_{p,q}w\right)  =\left(
\operatorname*{NoSt}w\right)  \setminus\left\{  q\right\}  $. This proves
Lemma \ref{lem.YJM.some-prod.2} \textbf{(b)}.
\end{proof}

\begin{proof}
[Proof of Theorem \ref{thm.YJM.some-prod}.]We induct on $k$.

\textit{Base case:} We must prove Theorem \ref{thm.YJM.some-prod} in the case
when $k=0$. In this case, the product $\mathbf{m}_{i_{1}}\mathbf{m}_{i_{2}%
}\cdots\mathbf{m}_{i_{k}}$ is an empty product, and thus equals $1$. On the
other hand, the sum $\sum_{\substack{w\in S_{n};\\\operatorname*{NoSt}%
w=\left\{  i_{1},i_{2},\ldots,i_{k}\right\}  }}w$ equals $\operatorname*{id}$
in this case, because its only addend is $\operatorname*{id}$ (since the only
permutation $w\in S_{n}$ that has no nonstarters is $\operatorname*{id}$ (by
Remark \ref{rmk.orbits.nost.easy} \textbf{(c)})). Thus, in the case when
$k=0$, the claim of Theorem \ref{thm.YJM.some-prod} boils down to
$1=\operatorname*{id}$, which is true. This completes the base case.

\textit{Induction step:} We proceed from $k-1$ to $k$. Thus, we fix a positive
integer $k$, and we set out to prove Theorem \ref{thm.YJM.some-prod} for $k$,
assuming (as the induction hypothesis) that it has already been proved for
$k-1$.

So let $i_{1},i_{2},\ldots,i_{k}$ be $k$ distinct elements of $\left[
n\right]  $. We must prove that%
\begin{equation}
\mathbf{m}_{i_{1}}\mathbf{m}_{i_{2}}\cdots\mathbf{m}_{i_{k}}=\sum
_{\substack{w\in S_{n};\\\operatorname*{NoSt}w=\left\{  i_{1},i_{2}%
,\ldots,i_{k}\right\}  }}w. \label{pf.thm.YJM.some-prod.IG}%
\end{equation}
This claim does not change if we permute the elements $i_{1},i_{2}%
,\ldots,i_{k}$ (since the jucys--murphies commute). Thus, we can WLOG assume
that $i_{1}\leq i_{2}\leq\cdots\leq i_{k}$. Assume this. Hence, $i_{1}%
<i_{2}<\cdots<i_{k}$ (since $i_{1},i_{2},\ldots,i_{k}$ are distinct).

Our induction hypothesis yields%
\begin{equation}
\mathbf{m}_{i_{1}}\mathbf{m}_{i_{2}}\cdots\mathbf{m}_{i_{k-1}}=\sum
_{\substack{w\in S_{n};\\\operatorname*{NoSt}w=\left\{  i_{1},i_{2}%
,\ldots,i_{k-1}\right\}  }}w. \label{pf.thm.YJM.some-prod.IH}%
\end{equation}

Now,%
\begin{align}
\mathbf{m}_{i_{1}}\mathbf{m}_{i_{2}}\cdots\mathbf{m}_{i_{k}}  &  =\left(
\mathbf{m}_{i_{1}}\mathbf{m}_{i_{2}}\cdots\mathbf{m}_{i_{k-1}}\right)
\mathbf{m}_{i_{k}}\nonumber\\
&  =\mathbf{m}_{i_{k}}\left(  \mathbf{m}_{i_{1}}\mathbf{m}_{i_{2}}%
\cdots\mathbf{m}_{i_{k-1}}\right)  \ \ \ \ \ \ \ \ \ \ \left(  \text{since the
jucys--murphies commute}\right) \nonumber\\
&  =\mathbf{m}_{i_{k}}\sum_{\substack{w\in S_{n};\\\operatorname*{NoSt}%
w=\left\{  i_{1},i_{2},\ldots,i_{k-1}\right\}  }}w\ \ \ \ \ \ \ \ \ \ \left(
\text{by (\ref{pf.thm.YJM.some-prod.IH})}\right) \nonumber\\
&  =\sum_{p=1}^{i_{k}-1}t_{p,i_{k}}\sum_{\substack{w\in S_{n}%
;\\\operatorname*{NoSt}w=\left\{  i_{1},i_{2},\ldots,i_{k-1}\right\}
}}w\nonumber\\
&  \ \ \ \ \ \ \ \ \ \ \ \ \ \ \ \ \ \ \ \ \left(
\begin{array}
[c]{c}%
\text{since the definition of }\mathbf{m}_{i_{k}}\\
\text{says }\mathbf{m}_{i_{k}}=t_{1,i_{k}}+t_{2,i_{k}}+\cdots+t_{i_{k}%
-1,i_{k}}=\sum_{p=1}^{i_{k}-1}t_{p,i_{k}}%
\end{array}
\right) \nonumber\\
&  =\sum_{p=1}^{i_{k}-1}\ \ \sum_{\substack{w\in S_{n};\\\operatorname*{NoSt}%
w=\left\{  i_{1},i_{2},\ldots,i_{k-1}\right\}  }}t_{p,i_{k}}w.
\label{pf.thm.YJM.some-prod.4}%
\end{align}

Let us now show that each $p\in\left[  i_{k}-1\right]  $ satisfies%
\begin{equation}
\sum_{\substack{w\in S_{n};\\\operatorname*{NoSt}w=\left\{  i_{1},i_{2}%
,\ldots,i_{k-1}\right\}  }}t_{p,i_{k}}w=\sum_{\substack{w\in S_{n}%
;\\\operatorname*{NoSt}w=\left\{  i_{1},i_{2},\ldots,i_{k}\right\}
;\\w\left(  i_{k}\right)  =p}}w. \label{pf.thm.YJM.some-prod.5}%
\end{equation}

\begin{proof}
[Proof of (\ref{pf.thm.YJM.some-prod.5}).]Let $p\in\left[  i_{k}-1\right]  $.
Thus, $p\leq i_{k}-1<i_{k}$. Define the two sets%
\begin{align*}
U  &  =\left\{  w\in S_{n}\ \mid\ \operatorname*{NoSt}w=\left\{  i_{1}%
,i_{2},\ldots,i_{k}\right\}  \text{ and }w\left(  i_{k}\right)  =p\right\}
\ \ \ \ \ \ \ \ \ \ \text{and}\\
V  &  =\left\{  w\in S_{n}\ \mid\ \operatorname*{NoSt}w=\left\{  i_{1}%
,i_{2},\ldots,i_{k-1}\right\}  \right\}  .
\end{align*}
We shall now show that the maps%
\begin{align*}
\Phi:V  &  \rightarrow U,\\
w  &  \mapsto t_{p,i_{k}}w
\end{align*}
and%
\begin{align*}
\Psi:U  &  \rightarrow V,\\
w  &  \mapsto t_{p,i_{k}}w
\end{align*}
are well-defined.

\textit{Well-definedness of }$\Phi$\textit{:} In order to show that $\Phi$ is
well-defined, we must prove that $t_{p,i_{k}}w\in U$ for each $w\in V$. So let
$w\in V$ be arbitrary. Thus, $w\in S_{n}$ and $\operatorname*{NoSt}w=\left\{
i_{1},i_{2},\ldots,i_{k-1}\right\}  $ (by the definition of $V$). The latter
equality shows that the nonstarters of $w$ are $i_{1},i_{2},\ldots,i_{k-1}$.
Hence, $i_{k}$ is larger than each nonstarter of $w$ (since $i_{1}%
<i_{2}<\cdots<i_{k}$). Thus, Lemma \ref{lem.YJM.some-prod.1} \textbf{(a)}
(applied to $i_{k}$ instead of $p$) shows that $w\left(  i_{k}\right)  =i_{k}%
$. Hence, Lemma \ref{lem.YJM.some-prod.2} \textbf{(a)} (applied to $q=i_{k}$)
yields that
\[
\operatorname*{NoSt}\left(  t_{p,i_{k}}w\right)  =\underbrace{\left(
\operatorname*{NoSt}w\right)  }_{=\left\{  i_{1},i_{2},\ldots,i_{k-1}\right\}
}\cup\left\{  i_{k}\right\}  =\left\{  i_{1},i_{2},\ldots,i_{k-1}\right\}
\cup\left\{  i_{k}\right\}  =\left\{  i_{1},i_{2},\ldots,i_{k}\right\}  .
\]
Furthermore, $\left(  t_{p,i_{k}}w\right)  \left(  i_{k}\right)  =t_{p,i_{k}%
}\left(  \underbrace{w\left(  i_{k}\right)  }_{=i_{k}}\right)  =t_{p,i_{k}%
}\left(  i_{k}\right)  =p$. These two equalities, combined, show that
$t_{p,i_{k}}w\in U$ (by the definition of $U$). Thus, we have shown that
$\Phi$ is well-defined.

\textit{Well-definedness of }$\Psi$\textit{:} In order to show that $\Psi$ is
well-defined, we must prove that $t_{p,i_{k}}w\in V$ for each $w\in U$. So let
$w\in U$ be arbitrary. Thus, $w\in S_{n}$ and $\operatorname*{NoSt}w=\left\{
i_{1},i_{2},\ldots,i_{k}\right\}  $ and $w\left(  i_{k}\right)  =p$ (by the
definition of $U$). The equality $\operatorname*{NoSt}w=\left\{  i_{1}%
,i_{2},\ldots,i_{k}\right\}  $ shows that the nonstarters of $w$ are
$i_{1},i_{2},\ldots,i_{k}$. Hence, $i_{k}$ is the largest nonstarter of $w$
(since $i_{1}<i_{2}<\cdots<i_{k}$). Thus, Lemma \ref{lem.YJM.some-prod.2}
\textbf{(b)} (applied to $i_{k}$ instead of $q$) shows that
\begin{align*}
\operatorname*{NoSt}\left(  t_{p,i_{k}}w\right)   &  =\underbrace{\left(
\operatorname*{NoSt}w\right)  }_{=\left\{  i_{1},i_{2},\ldots,i_{k}\right\}
}\setminus\,\left\{  i_{k}\right\}  \ \ \ \ \ \ \ \ \ \ \left(  \text{since
}w\left(  i_{k}\right)  =p\right) \\
&  =\left\{  i_{1},i_{2},\ldots,i_{k}\right\}  \setminus\left\{
i_{k}\right\}  =\left\{  i_{1},i_{2},\ldots,i_{k-1}\right\}
\end{align*}
(since $i_{1},i_{2},\ldots,i_{k}$ are distinct). Hence, $t_{p,i_{k}}w\in V$
(by the definition of $V$). Thus, we have shown that $\Psi$ is well-defined.

The maps $\Phi$ and $\Psi$ satisfy $\Phi\circ\Psi=\operatorname*{id}$ (since
the definitions of $\Phi$ and $\Psi$ yield $\Phi\left(  \Psi\left(  w\right)
\right)  =t_{p,i_{k}}\left(  t_{p,i_{k}}w\right)  =\underbrace{t_{p,i_{k}%
}t_{p,i_{k}}}_{=t_{p,i_{k}}^{2}=\operatorname*{id}}w=w$ for each $w\in U$) and
$\Psi\circ\Phi=\operatorname*{id}$ (for similar reasons). Thus, they are
mutually inverse. Hence, the map%
\begin{align*}
\Phi:V  &  \rightarrow U,\\
w  &  \mapsto t_{p,i_{k}}w
\end{align*}
is invertible (with inverse $\Psi$), and therefore is a bijection. Hence, we
can substitute $t_{p,i_{k}}w$ for $w$ in the sum $\sum_{w\in U}w$, thus
obtaining the equality%
\[
\sum_{w\in U}w=\sum_{w\in V}t_{p,i_{k}}w.
\]
Recalling the definitions of $U$ and $V$, we can rewrite the two summation
signs $\sum_{w\in U}$ and $\sum_{w\in V}$ in this equality as $\sum
_{\substack{w\in S_{n};\\\operatorname*{NoSt}w=\left\{  i_{1},i_{2}%
,\ldots,i_{k}\right\}  ;\\w\left(  i_{k}\right)  =p}}$ and $\sum
_{\substack{w\in S_{n};\\\operatorname*{NoSt}w=\left\{  i_{1},i_{2}%
,\ldots,i_{k-1}\right\}  }}$, respectively. Thus, this equality rewrites as%
\[
\sum_{\substack{w\in S_{n};\\\operatorname*{NoSt}w=\left\{  i_{1},i_{2}%
,\ldots,i_{k}\right\}  ;\\w\left(  i_{k}\right)  =p}}w=\sum_{\substack{w\in
S_{n};\\\operatorname*{NoSt}w=\left\{  i_{1},i_{2},\ldots,i_{k-1}\right\}
}}t_{p,i_{k}}w.
\]
Hence, (\ref{pf.thm.YJM.some-prod.5}) is proved.
\end{proof}

However, it is easy to see that if a permutation $w\in S_{n}$ satisfies
$\operatorname*{NoSt}w=\left\{  i_{1},i_{2},\ldots,i_{k}\right\}  $, then
$w\left(  i_{k}\right)  \in\left[  i_{k}-1\right]  $%
\ \ \ \ \footnote{\textit{Proof.} Let $w\in S_{n}$ be a permutation that
satisfies $\operatorname*{NoSt}w=\left\{  i_{1},i_{2},\ldots,i_{k}\right\}  $.
The equality $\operatorname*{NoSt}w=\left\{  i_{1},i_{2},\ldots,i_{k}\right\}
$ shows that the nonstarters of $w$ are $i_{1},i_{2},\ldots,i_{k}$. Hence,
$i_{k}$ is the largest nonstarter of $w$ (since $i_{1}<i_{2}<\cdots<i_{k}$).
Hence, Lemma \ref{lem.YJM.some-prod.1} \textbf{(b)} (applied to $q=i_{k}$)
yields that $w\left(  i_{k}\right)  \in\left[  i_{k}-1\right]  $, qed.}.
Hence, we can split the sum $\sum_{\substack{w\in S_{n};\\\operatorname*{NoSt}%
w=\left\{  i_{1},i_{2},\ldots,i_{k}\right\}  }}w$ according to the value of
$w\left(  i_{k}\right)  $ as follows:%
\[
\sum_{\substack{w\in S_{n};\\\operatorname*{NoSt}w=\left\{  i_{1},i_{2}%
,\ldots,i_{k}\right\}  }}w=\underbrace{\sum_{p\in\left[  i_{k}-1\right]  }%
}_{=\sum_{p=1}^{i_{k}-1}}\ \ \underbrace{\sum_{\substack{w\in S_{n}%
;\\\operatorname*{NoSt}w=\left\{  i_{1},i_{2},\ldots,i_{k}\right\}
;\\w\left(  i_{k}\right)  =p}}w}_{\substack{=\sum_{\substack{w\in
S_{n};\\\operatorname*{NoSt}w=\left\{  i_{1},i_{2},\ldots,i_{k-1}\right\}
}}t_{p,i_{k}}w\\\text{(by (\ref{pf.thm.YJM.some-prod.5}))}}}=\sum_{p=1}%
^{i_{k}-1}\ \ \sum_{\substack{w\in S_{n};\\\operatorname*{NoSt}w=\left\{
i_{1},i_{2},\ldots,i_{k-1}\right\}  }}t_{p,i_{k}}w.
\]
Comparing this with (\ref{pf.thm.YJM.some-prod.4}), we find
\[
\mathbf{m}_{i_{1}}\mathbf{m}_{i_{2}}\cdots\mathbf{m}_{i_{k}}=\sum
_{\substack{w\in S_{n};\\\operatorname*{NoSt}w=\left\{  i_{1},i_{2}%
,\ldots,i_{k}\right\}  }}w.
\]
Thus, (\ref{pf.thm.YJM.some-prod.IG}) is proved. This completes the induction
step. Thus, Theorem \ref{thm.YJM.some-prod} is proved.
\end{proof}

\begin{exercise}
\fbox{1} Reprove Theorem \ref{thm.YJM.prod1+mk} using Theorem
\ref{thm.YJM.some-prod}.
\end{exercise}

But we can get more mileage out of Theorem \ref{thm.YJM.some-prod}. Indeed,
recall the following definition (a consequence of \cite[Definition 7.1.9
\textbf{(a)}]{21s}):

\begin{definition}
\label{def.ek}Let $u_{1},u_{2},\ldots,u_{n}$ be any $n$ commuting elements of
a ring. Then, for each $k\in\mathbb{Z}$, we define
\[
e_{k}\left(  u_{1},u_{2},\ldots,u_{n}\right)  :=\sum_{1\leq i_{1}<i_{2}%
<\cdots<i_{k}\leq n}u_{i_{1}}u_{i_{2}}\cdots u_{i_{k}}%
\]
(the sum ranges over all $k$-tuples $\left(  i_{1},i_{2},\ldots,i_{k}\right)
\in\left[  n\right]  ^{k}$ satisfying $i_{1}<i_{2}<\cdots<i_{k}$). In other
words, we define $e_{k}\left(  u_{1},u_{2},\ldots,u_{n}\right)  $ to be the
sum of all products of $k$ among the $n$ inputs $u_{1},u_{2},\ldots,u_{n}$.

This $e_{k}\left(  u_{1},u_{2},\ldots,u_{n}\right)  $ is called the
$k$\emph{-th elementary symmetric polynomial} in the inputs $u_{1}%
,u_{2},\ldots,u_{n}$.
\end{definition}

For example,%
\begin{align*}
e_{0}\left(  u_{1},u_{2},\ldots,u_{n}\right)   &  =1;\\
e_{1}\left(  u_{1},u_{2},\ldots,u_{n}\right)   &  =u_{1}+u_{2}+\cdots+u_{n};\\
e_{2}\left(  u_{1},u_{2},\ldots,u_{n}\right)   &  =u_{1}u_{2}+u_{1}%
u_{3}+\cdots+u_{1}u_{n}\\
&  \ \ \ \ \ \ \ \ \ \ \ \ +u_{2}u_{3}+\cdots+u_{2}u_{n}\\
&  \ \ \ \ \ \ \ \ \ \ \ \ \ \ \ \ \ \ \ \ \ \ \ +\cdots\\
&  \ \ \ \ \ \ \ \ \ \ \ \ \ \ \ \ \ \ \ \ \ \ \ \ \ \ \ \ \ \ \ \ +u_{n-1}%
u_{n};\\
e_{n}\left(  u_{1},u_{2},\ldots,u_{n}\right)   &  =u_{1}u_{2}\cdots u_{n};\\
e_{k}\left(  u_{1},u_{2},\ldots,u_{n}\right)   &
=0\ \ \ \ \ \ \ \ \ \ \text{for }k>n.
\end{align*}
The best-known example of $e_{k}\left(  u_{1},u_{2},\ldots,u_{n}\right)  $ is
obtained when the inputs $u_{1},u_{2},\ldots,u_{n}$ are the $n$ indeterminates
$x_{1},x_{2},\ldots,x_{n}$ of the polynomial ring $\mathbf{k}\left[
x_{1},x_{2},\ldots,x_{n}\right]  $. In this case, we obtain the polynomial
$e_{k}\left(  x_{1},x_{2},\ldots,x_{n}\right)  $, which is known as the
$k$\emph{-th elementary symmetric polynomial} in the $n$ indeterminates
$x_{1},x_{2},\ldots,x_{n}$. This polynomial is the $e_{k}$ defined in
\cite[Definition 7.1.9 \textbf{(a)}]{21s}, and has many important properties
(figuring, e.g., in
\href{https://en.wikipedia.org/wiki/Vieta's_formulas}{Vieta's formulas}). In
the general case, $e_{k}\left(  u_{1},u_{2},\ldots,u_{n}\right)  $ can be
obtained from this polynomial $e_{k}$ by substituting $u_{i}$ for each $x_{i}$.

The jucys--murphies $\mathbf{m}_{1},\mathbf{m}_{2},\ldots,\mathbf{m}_{n}$
commute (by Theorem \ref{thm.YJM.commute}); thus, we can apply the elementary
symmetric polynomials to them. We obtain something rather nice (\cite[(8)]%
{Jucys74}):

\begin{corollary}
\label{cor.YJM.ek}For each $k\in\mathbb{N}$, we have%
\[
e_{k}\left(  \mathbf{m}_{1},\mathbf{m}_{2},\ldots,\mathbf{m}_{n}\right)
=\sum_{\substack{w\in S_{n};\\w\text{ has exactly }n-k\text{ orbits}}%
}w=\sum_{\substack{w\in S_{n};\\\operatorname*{rl}w=k}}w.
\]

\end{corollary}

\begin{proof}
For each $w\in S_{n}$, we have%
\begin{align}
\left\vert \operatorname*{NoSt}w\right\vert  &  =\left(  \text{\# of
nonstarters of }w\right)  \ \ \ \ \ \ \ \ \ \ \left(  \text{by the definition
of }\operatorname*{NoSt}w\right) \nonumber\\
&  =n-\underbrace{\left(  \text{\# of starters of }w\right)  }%
_{\substack{=\left(  \text{\# of orbits of }w\right)  \\\text{(by Remark
\ref{rmk.orbits.nost.easy} \textbf{(d)})}}}\nonumber\\
&  \ \ \ \ \ \ \ \ \ \ \ \ \ \ \ \ \ \ \ \ \left(
\begin{array}
[c]{c}%
\text{since each }i\in\left[  n\right]  \text{ is either a starter}\\
\text{or a nonstarter, but not both}%
\end{array}
\right) \nonumber\\
&  =\underbrace{n}_{=\left\vert \left[  n\right]  \right\vert }-\left(
\text{\# of orbits of }w\right) \label{pf.cor.YJM.ek.0}\\
&  =\left\vert \left[  n\right]  \right\vert -\left(  \text{\# of orbits of
}w\right) \nonumber\\
&  =\operatorname*{rl}w \label{pf.cor.YJM.ek.1}%
\end{align}
(by Definition \ref{def.rl.rl}).

By Definition \ref{def.ek}, we have%
\begin{align*}
&  e_{k}\left(  \mathbf{m}_{1},\mathbf{m}_{2},\ldots,\mathbf{m}_{n}\right) \\
&  =\sum_{1\leq i_{1}<i_{2}<\cdots<i_{k}\leq n}\mathbf{m}_{i_{1}}%
\mathbf{m}_{i_{2}}\cdots\mathbf{m}_{i_{k}}\\
&  =\sum_{1\leq i_{1}<i_{2}<\cdots<i_{k}\leq n}\ \ \sum_{\substack{w\in
S_{n};\\\operatorname*{NoSt}w=\left\{  i_{1},i_{2},\ldots,i_{k}\right\}
}}w\ \ \ \ \ \ \ \ \ \ \left(  \text{by Theorem \ref{thm.YJM.some-prod}%
}\right) \\
&  =\underbrace{\sum_{\substack{I\subseteq\left[  n\right]  ;\\\left\vert
I\right\vert =k}}\ \ \sum_{\substack{w\in S_{n};\\\operatorname*{NoSt}w=I}%
}}_{=\sum_{\substack{w\in S_{n};\\\left\vert \operatorname*{NoSt}w\right\vert
=k}}}w\ \ \ \ \ \ \ \ \ \ \left(
\begin{array}
[c]{c}%
\text{here, we have substituted }I\text{ for }\left\{  i_{1},i_{2}%
,\ldots,i_{k}\right\} \\
\text{in the sum, since each }k\text{-element subset }I\\
\text{of }\left[  n\right]  \text{ can be written as }\left\{  i_{1}%
,i_{2},\ldots,i_{k}\right\} \\
\text{for a unique }k\text{-tuple }\left(  i_{1},i_{2},\ldots,i_{k}\right)
\in\left[  n\right]  ^{k}\\
\text{satisfying }i_{1}<i_{2}<\cdots<i_{k}%
\end{array}
\right) \\
&  =\sum_{\substack{w\in S_{n};\\\left\vert \operatorname*{NoSt}w\right\vert
=k}}w=\sum_{\substack{w\in S_{n};\\\operatorname*{rl}w=k}%
}w\ \ \ \ \ \ \ \ \ \ \left(  \text{by (\ref{pf.cor.YJM.ek.1})}\right) \\
&  =\sum_{\substack{w\in S_{n};\\\left\vert \operatorname*{NoSt}w\right\vert
=k}}w=\sum_{\substack{w\in S_{n};\\n-\left(  \text{\# of orbits of }w\right)
=k}}w\ \ \ \ \ \ \ \ \ \ \left(  \text{by (\ref{pf.cor.YJM.ek.0})}\right) \\
&  =\sum_{\substack{w\in S_{n};\\\left(  \text{\# of orbits of }w\right)
=n-k}}w\ \ \ \ \ \ \ \ \ \ \left(
\begin{array}
[c]{c}%
\text{since the condition \textquotedblleft}n-\left(  \text{\# of orbits of
}w\right)  =k\text{\textquotedblright}\\
\text{is equivalent to \textquotedblleft}\left(  \text{\# of orbits of
}w\right)  =n-k\text{\textquotedblright}%
\end{array}
\right) \\
&  =\sum_{\substack{w\in S_{n};\\w\text{ has exactly }n-k\text{ orbits}}}w.
\end{align*}
Thus, Corollary \ref{cor.YJM.ek} is proved.
\end{proof}

\begin{corollary}
\label{cor.YJM.prod}If $n\geq1$, then%
\[
\mathbf{m}_{2}\mathbf{m}_{3}\cdots\mathbf{m}_{n}=\sum_{\substack{w\in
S_{n};\\w\text{ has exactly }1\text{ cycle}}}w.
\]

\end{corollary}

\begin{proof}
Assume that $n\geq1$. Then, Corollary \ref{cor.YJM.ek} (applied to $k=n-1$)
yields%
\begin{equation}
e_{n-1}\left(  \mathbf{m}_{1},\mathbf{m}_{2},\ldots,\mathbf{m}_{n}\right)
=\sum_{\substack{w\in S_{n};\\w\text{ has exactly }1\text{ cycle}}}w.
\label{pf.cor.YJM.prod.1}%
\end{equation}
Now, recall that $e_{n-1}\left(  \mathbf{m}_{1},\mathbf{m}_{2},\ldots
,\mathbf{m}_{n}\right)  $ is the sum of all products of $n-1$ among the $n$
jucys--murphies $\mathbf{m}_{1},\mathbf{m}_{2},\ldots,\mathbf{m}_{n}$ (by
Definition \ref{def.ek}). However, all such products except for $\mathbf{m}%
_{2}\mathbf{m}_{3}\cdots\mathbf{m}_{n}$ are $0$ (since they contain the factor
$\mathbf{m}_{1}$, which is $0$). Thus, $e_{n-1}\left(  \mathbf{m}%
_{1},\mathbf{m}_{2},\ldots,\mathbf{m}_{n}\right)  =\mathbf{m}_{2}%
\mathbf{m}_{3}\cdots\mathbf{m}_{n}$. Comparing this with
(\ref{pf.cor.YJM.prod.1}), we obtain%
\[
\mathbf{m}_{2}\mathbf{m}_{3}\cdots\mathbf{m}_{n}=\sum_{\substack{w\in
S_{n};\\w\text{ has exactly }1\text{ cycle}}}w.
\]
This proves Corollary \ref{cor.YJM.prod}.
\end{proof}

\begin{question}
Can we generalize Theorem \ref{thm.YJM.some-prod} to non-distinct $i_{1}%
,i_{2},\ldots,i_{k}$ in some way? In other words, is there an explicit formula
for a product $\mathbf{m}_{i_{1}}\mathbf{m}_{i_{2}}\cdots\mathbf{m}_{i_{k}}$
of non-distinct jucys--murphies? In particular, is there a formula for
$\mathbf{m}_{i}^{k}$ ?

Be ready to see the same permutation multiple times in such a formula. For
instance,%
\[
\mathbf{m}_{3}^{2}=\left(  t_{1,3}+t_{2,3}\right)  ^{2}=\underbrace{2}%
_{=2\operatorname*{id}}+\operatorname*{cyc}\nolimits_{1,2,3}%
+\operatorname*{cyc}\nolimits_{1,3,2}.
\]

\end{question}

The following corollary is just a neat way to repack Corollary
\ref{cor.YJM.ek}:

\begin{corollary}
\label{cor.YJM.1+xm}Let $x\in\mathbf{k}$. Then,%
\begin{align*}
&  \left(  1+x\mathbf{m}_{1}\right)  \left(  1+x\mathbf{m}_{2}\right)
\cdots\left(  1+x\mathbf{m}_{n}\right) \\
&  =\sum_{k\in\mathbb{N}}x^{k}e_{k}\left(  \mathbf{m}_{1},\mathbf{m}%
_{2},\ldots,\mathbf{m}_{n}\right)  =\sum_{k\in\mathbb{N}}x^{k}\sum
_{\substack{w\in S_{n};\\\operatorname*{rl}w=k}}w=\sum_{w\in S_{n}%
}x^{\operatorname*{rl}w}w.
\end{align*}

\end{corollary}

For instance, for $n=3$, this is saying that%
\begin{align*}
&  \left(  1+x\mathbf{m}_{1}\right)  \left(  1+x\mathbf{m}_{2}\right)  \left(
1+x\mathbf{m}_{3}\right) \\
&  =1+xt_{1,2}+xt_{1,3}+xt_{2,3}+x^{2}\operatorname*{cyc}\nolimits_{1,2,3}%
+\,x^{2}\operatorname*{cyc}\nolimits_{1,3,2}.
\end{align*}

\begin{proof}
[Proof of Corollary \ref{cor.YJM.1+xm}.]It is well-known that any $n$ elements
$a_{1},a_{2},\ldots,a_{n}$ of a ring $R$ satisfy the equality%
\[
\left(  1+a_{1}\right)  \left(  1+a_{2}\right)  \cdots\left(  1+a_{n}\right)
=\sum_{1\leq i_{1}<i_{2}<\cdots<i_{k}\leq n}a_{i_{1}}a_{i_{2}}\cdots a_{i_{k}}%
\]
(where the sum ranges over all strictly increasing tuples $\left(  i_{1}%
<i_{2}<\cdots<i_{k}\right)  $ of elements of $\left[  n\right]  $, with
arbitrary $k$). (Indeed, this equality can be obtained by expanding the left
hand side.) Applying this equality to $R=\mathbf{k}\left[  S_{n}\right]  $ and
$a_{i}=x\mathbf{m}_{i}$, we obtain%
\begin{align*}
&  \left(  1+x\mathbf{m}_{1}\right)  \left(  1+x\mathbf{m}_{2}\right)
\cdots\left(  1+x\mathbf{m}_{n}\right) \\
&  =\underbrace{\sum_{1\leq i_{1}<i_{2}<\cdots<i_{k}\leq n}}_{=\sum
_{k\in\mathbb{N}}\ \ \sum_{1\leq i_{1}<i_{2}<\cdots<i_{k}\leq n}%
}\underbrace{\left(  x\mathbf{m}_{i_{1}}\right)  \left(  x\mathbf{m}_{i_{2}%
}\right)  \cdots\left(  x\mathbf{m}_{i_{k}}\right)  }_{=x^{k}\mathbf{m}%
_{i_{1}}\mathbf{m}_{i_{2}}\cdots\mathbf{m}_{i_{k}}}\\
&  =\sum_{k\in\mathbb{N}}\ \ \sum_{1\leq i_{1}<i_{2}<\cdots<i_{k}\leq n}%
x^{k}\mathbf{m}_{i_{1}}\mathbf{m}_{i_{2}}\cdots\mathbf{m}_{i_{k}}\\
&  =\sum_{k\in\mathbb{N}}x^{k}\underbrace{\sum_{1\leq i_{1}<i_{2}<\cdots
<i_{k}\leq n}\mathbf{m}_{i_{1}}\mathbf{m}_{i_{2}}\cdots\mathbf{m}_{i_{k}}%
}_{\substack{=e_{k}\left(  \mathbf{m}_{1},\mathbf{m}_{2},\ldots,\mathbf{m}%
_{n}\right)  \\\text{(by the definition of }e_{k}\left(  \mathbf{m}%
_{1},\mathbf{m}_{2},\ldots,\mathbf{m}_{n}\right)  \text{)}}}\\
&  =\sum_{k\in\mathbb{N}}x^{k}\underbrace{e_{k}\left(  \mathbf{m}%
_{1},\mathbf{m}_{2},\ldots,\mathbf{m}_{n}\right)  }_{\substack{=\sum
_{\substack{w\in S_{n};\\\operatorname*{rl}w=k}}w\\\text{(by Corollary
\ref{cor.YJM.ek})}}}=\sum_{k\in\mathbb{N}}x^{k}\sum_{\substack{w\in
S_{n};\\\operatorname*{rl}w=k}}w\\
&  =\sum_{k\in\mathbb{N}}\ \ \sum_{\substack{w\in S_{n};\\\operatorname*{rl}%
w=k}}\ \ \underbrace{x^{k}}_{\substack{=x^{\operatorname*{rl}w}\\\text{(since
}k=\operatorname*{rl}w\text{)}}}w=\underbrace{\sum_{k\in\mathbb{N}}%
\ \ \sum_{\substack{w\in S_{n};\\\operatorname*{rl}w=k}}}_{=\sum_{w\in S_{n}}%
}x^{\operatorname*{rl}w}w=\sum_{w\in S_{n}}x^{\operatorname*{rl}w}w.
\end{align*}
Thus, Corollary \ref{cor.YJM.1+xm} is proved.
\end{proof}

\subsection{Conjugacy classes and the center}

\subsubsection{The general case}

Let us now recall some more abstract algebra, namely the definition of a
\emph{center} for groups and for rings:

\begin{definition}
\label{def.center.center}\textbf{(a)} The \emph{center} of a group $G$ is
defined to be the subset%
\[
Z\left(  G\right)  :=\left\{  x\in G\ \mid\ xy=yx\text{ for all }y\in
G\right\}  \text{ of }G.
\]
This is the set of all elements of $G$ that commute with all elements of $G$.
It is always an abelian subgroup of $G$. If $G$ itself is abelian, then
$Z\left(  G\right)  =G$. \medskip

\textbf{(b)} The \emph{center} of a ring $R$ is defined to be the subset%
\[
Z\left(  R\right)  :=\left\{  x\in R\ \mid\ xy=yx\text{ for all }y\in
R\right\}  \text{ of }R.
\]
This is the set of all elements of $R$ that commute with all elements of $R$.
It is always a commutative subring of $R$. Moreover, if $R$ is a $\mathbf{k}%
$-algebra, then this center is also a $\mathbf{k}$-subalgebra of $R$. If $R$
itself is commutative, then $Z\left(  R\right)  =R$.
\end{definition}

The two kinds of centers are clearly analogous, and one might wonder how they
interact in the case of a group algebra. It is easy to see that if $g$ is an
element of the center $Z\left(  G\right)  $ of a group $G$, then the
corresponding standard basis vector $g=e_{g}\in\mathbf{k}\left[  G\right]  $
belongs to the center $Z\left(  \mathbf{k}\left[  G\right]  \right)  $ of the
group algebra $\mathbf{k}\left[  G\right]  $ (this is by linearity). Taking
linear combinations of such basis vectors yields further elements of $Z\left(
\mathbf{k}\left[  G\right]  \right)  $ (since $Z\left(  \mathbf{k}\left[
G\right]  \right)  $ is a $\mathbf{k}$-subalgebra). However, not every element
of $Z\left(  \mathbf{k}\left[  G\right]  \right)  $ is obtained this way (in
general). For example:

\begin{proposition}
\label{prop.center.sum-mi}The sum%
\[
\mathbf{m}_{1}+\mathbf{m}_{2}+\cdots+\mathbf{m}_{n}=\sum_{i<j}t_{i,j}%
\]
(that is, the sum of all transpositions in $S_{n}$) belongs to the center of
$\mathbf{k}\left[  S_{n}\right]  $.
\end{proposition}

We will soon prove this. Note that none of the transpositions $t_{i,j}$
themselves belongs to $Z\left(  S_{n}\right)  $ (unless $n\leq2$).

Let us first study the center of $\mathbf{k}\left[  G\right]  $ for arbitrary
$G$:

\begin{lemma}
\label{lem.center.waw}Let $G$ be a group. Let $\mathbf{a}\in\mathbf{k}\left[
G\right]  $. Then, $\mathbf{a}\in Z\left(  \mathbf{k}\left[  G\right]
\right)  $ if and only if each $w\in G$ satisfies $w\mathbf{a}w^{-1}%
=\mathbf{a}$.
\end{lemma}

\begin{proof}
$\Longrightarrow:$ Assume that $\mathbf{a}\in Z\left(  \mathbf{k}\left[
G\right]  \right)  $. Then, each $w\in G$ satisfies $w\mathbf{a}=\mathbf{a}w$
and therefore $\underbrace{w\mathbf{a}}_{=\mathbf{a}w}w^{-1}=\mathbf{a}%
\underbrace{ww^{-1}}_{=1}=\mathbf{a}$. In other words, each $w\in G$ satisfies
$w\mathbf{a}w^{-1}=\mathbf{a}$. This proves the \textquotedblleft%
$\Longrightarrow$\textquotedblright\ direction of Lemma \ref{lem.center.waw}.

$\Longleftarrow:$ Assume that each $w\in G$ satisfies $w\mathbf{a}%
w^{-1}=\mathbf{a}$. Thus, each $w\in G$ satisfies
\begin{equation}
w\mathbf{a}=\mathbf{a}w \label{pf.lem.center.waw.2}%
\end{equation}
(since multiplying the equality $w\mathbf{a}w^{-1}=\mathbf{a}$ by $w$ on the
right yields $w\mathbf{a}w^{-1}w=\mathbf{a}w$, thus $\mathbf{a}w=w\mathbf{a}%
\underbrace{w^{-1}w}_{=1}=w\mathbf{a}$). By linearity, this entails that any
$\mathbf{k}$-linear combination $\sum\limits_{w\in G}\lambda_{w}w$ of elements
of $G$ (with $\lambda_{w}\in\mathbf{k}$ for all $w\in G$) satisfies%
\[
\left(  \sum\limits_{w\in G}\lambda_{w}w\right)  \mathbf{a}=\sum\limits_{w\in
G}\lambda_{w}\underbrace{w\mathbf{a}}_{\substack{=\mathbf{a}w\\\text{(by
(\ref{pf.lem.center.waw.2}))}}}=\sum\limits_{w\in G}\lambda_{w}\mathbf{a}%
w=\mathbf{a}\left(  \sum\limits_{w\in G}\lambda_{w}w\right)  .
\]
In other words, $\mathbf{a}$ commutes with any $\mathbf{k}$-linear combination
$\sum\limits_{w\in G}\lambda_{w}w$ of elements of $G$. But this is just saying
that $\mathbf{a}$ commutes with every element of $\mathbf{k}\left[  G\right]
$ (since every element of $\mathbf{k}\left[  G\right]  $ is such a
$\mathbf{k}$-linear combination). In other words, $\mathbf{a}\in Z\left(
\mathbf{k}\left[  G\right]  \right)  $. This proves the \textquotedblleft%
$\Longleftarrow$\textquotedblright\ direction of Lemma \ref{lem.center.waw}.
\end{proof}

Using Lemma \ref{lem.center.waw}, we can now prove Proposition
\ref{prop.center.sum-mi}:

\begin{proof}
[Proof of Proposition \ref{prop.center.sum-mi}.]Let $T=\left\{  t_{i,j}%
\ \mid\ i,j\in\left[  n\right]  \text{ with }i<j\right\}  $ be the set of all
transpositions in $S_{n}$. Then,
\[
\sum_{\sigma\in T}\sigma=\sum_{i<j}t_{i,j}%
\]
(since each transposition in $S_{n}$ can be written as $t_{i,j}$ for a unique
pair of integers $i,j\in\left[  n\right]  $ with $i<j$).

Set $\mathbf{a}=\sum_{\sigma\in T}\sigma$. Then, $\mathbf{a}$ is the sum of
all transpositions in $S_{n}$. Moreover,%
\begin{align*}
\mathbf{a}  &  =\sum_{\sigma\in T}\sigma=\sum_{i<j}t_{i,j}=\sum_{j=1}%
^{n}\ \ \underbrace{\sum_{i=1}^{j-1}t_{i,j}}_{\substack{=\mathbf{m}%
_{j}\\\text{(by the definition of }\mathbf{m}_{j}\text{)}}}=\sum_{j=1}%
^{n}\mathbf{m}_{j}\\
&  =\mathbf{m}_{1}+\mathbf{m}_{2}+\cdots+\mathbf{m}_{n}.
\end{align*}

We shall now show that $\mathbf{a}\in Z\left(  \mathbf{k}\left[  S_{n}\right]
\right)  $.

Indeed, let $w\in S_{n}$. For any transposition $\sigma\in T$, the permutation
$w\sigma w^{-1}$ is again a transposition (since
(\ref{eq.intro.perms.cycs.sts-1}) shows that $wt_{i,j}w^{-1}=t_{w\left(
i\right)  ,w\left(  j\right)  }$ for any distinct elements $i,j\in\left[
n\right]  $). In other words, $w\sigma w^{-1}\in T$ for any $\sigma\in T$.
Hence, the map%
\begin{align*}
T  &  \rightarrow T,\\
\sigma &  \mapsto w\sigma w^{-1}%
\end{align*}
is well-defined. Similarly, the map%
\begin{align*}
T  &  \rightarrow T,\\
\sigma &  \mapsto w^{-1}\sigma w
\end{align*}
is well-defined. These two maps are easily seen to be mutually inverse
(because for each $\sigma\in T$, we have $w^{-1}\left(  w\sigma w^{-1}\right)
w=\underbrace{w^{-1}w}_{=1}\sigma\underbrace{w^{-1}w}_{=1}=\sigma$ and
similarly $w\left(  w^{-1}\sigma w\right)  w^{-1}=\sigma$), and thus are
bijections. Hence, in particular, the map
\begin{align*}
T  &  \rightarrow T,\\
\sigma &  \mapsto w\sigma w^{-1}%
\end{align*}
is a bijection. Thus, we can substitute $w\sigma w^{-1}$ for $\sigma$ in the
sum $\sum_{\sigma\in T}\sigma$. We thus obtain%
\[
\sum_{\sigma\in T}\sigma=\sum_{\sigma\in T}w\sigma w^{-1}=w\left(
\sum_{\sigma\in T}\sigma\right)  w^{-1}.
\]
In view of $\mathbf{a}=\sum_{\sigma\in T}\sigma$, we can rewrite this as
$\mathbf{a}=w\mathbf{a}w^{-1}$. In other words, $w\mathbf{a}w^{-1}=\mathbf{a}$.

Forget now that we fixed $w$. We thus have shown that each $w\in S_{n}$
satisfies $w\mathbf{a}w^{-1}=\mathbf{a}$. Hence, by Lemma \ref{lem.center.waw}
(applied to $G=S_{n}$), we conclude that $\mathbf{a}\in Z\left(
\mathbf{k}\left[  S_{n}\right]  \right)  $. In other words, $\mathbf{a}$
belongs to the center of $\mathbf{k}\left[  S_{n}\right]  $. In other words,%
\[
\mathbf{m}_{1}+\mathbf{m}_{2}+\cdots+\mathbf{m}_{n}=\sum_{i<j}t_{i,j}%
\]
belongs to the center of $\mathbf{k}\left[  S_{n}\right]  $ (since
$\mathbf{a}=\mathbf{m}_{1}+\mathbf{m}_{2}+\cdots+\mathbf{m}_{n}=\sum
_{i<j}t_{i,j}$). This proves Proposition \ref{prop.center.sum-mi}.
\end{proof}

Let us try to describe $Z\left(  \mathbf{k}\left[  G\right]  \right)  $ more
explicitly. For this we recall some classical notations from group theory:

\begin{definition}
\label{def.groups.conj}Let $G$ be a group. \medskip

\textbf{(a)} Two elements $x$ and $y$ of $G$ are said to be \emph{conjugate}
if there exists a $g\in G$ such that $y=gxg^{-1}$. In this case, we write
$x\sim y$. \medskip

\textbf{(b)} The binary relation $\sim$ on $G$ is an equivalence relation, and
is known as \emph{conjugacy}. Its equivalence classes are called the
\emph{conjugacy classes} of $G$. Thus, each element of $G$ belongs to exactly
one conjugacy class. \medskip

\textbf{(c)} If $C$ is a finite conjugacy class of $G$, then we define
$\mathbf{z}_{C}$ to be the sum $\sum_{c\in C}c$ in $\mathbf{k}\left[
G\right]  $. This is called a \emph{conjugacy class sum} of $G$. (We require
$C$ to be finite, since infinite sums do not make sense in $\mathbf{k}\left[
G\right]  $. However, this will not be important for us in the case $G=S_{n}$,
since $S_{n}$ is finite.)
\end{definition}

Here are some examples:

\begin{itemize}
\item If $G$ is an abelian group, then conjugacy in $G$ is just equality
(i.e., we have $x\sim y$ if and only if $x=y$), so all the conjugacy classes
are singletons $\left\{  x\right\}  $ (with $x\in G$), and thus the conjugacy
class sums are just the standard basis vectors $x=e_{x}$ of $\mathbf{k}\left[
G\right]  $.

More generally, if $x$ is any element of the center of any group $G$, then its
conjugacy class is the singleton $\left\{  x\right\}  $ (since $gxg^{-1}=x$
for all $g\in G$).

\item The conjugacy classes of the symmetric group $S_{3}$ are%
\[
\left\{  \operatorname*{id}\right\}  ,\ \ \ \ \ \ \ \ \ \ \left\{
t_{1,2},\ t_{1,3},\ t_{2,3}\right\}  ,\ \ \ \ \ \ \ \ \ \ \left\{
\operatorname*{cyc}\nolimits_{1,2,3},\ \operatorname*{cyc}\nolimits_{1,3,2}%
\right\}  .
\]
Thus, the conjugacy class sums of $S_{3}$ are%
\[
\operatorname*{id}=1,\ \ \ \ \ \ \ \ \ \ t_{1,2}+t_{1,3}+t_{2,3}%
,\ \ \ \ \ \ \ \ \ \ \operatorname*{cyc}\nolimits_{1,2,3}+\operatorname*{cyc}%
\nolimits_{1,3,2}.
\]

\end{itemize}

See \cite[\textquotedblleft Conjugation in a group\textquotedblright]{Conrad}
for more examples. You might notice that in all these examples, the conjugacy
class sums all belong to $Z\left(  \mathbf{k}\left[  G\right]  \right)  $.
This does generalize, and even more is true:

\begin{theorem}
\label{thm.center.conjsums}Let $G$ be a group. Then, the center $Z\left(
\mathbf{k}\left[  G\right]  \right)  $ of the group algebra $\mathbf{k}\left[
G\right]  $ is the $\mathbf{k}$-linear span of the conjugacy class sums.
Moreover, these sums are $\mathbf{k}$-linearly independent, so they form a
basis of $Z\left(  \mathbf{k}\left[  G\right]  \right)  $.
\end{theorem}

\begin{proof}
[Proof of Theorem \ref{thm.center.conjsums}.]We need to prove the following
three claims:

\begin{statement}
\textit{Claim 1:} Each conjugacy class sum $\mathbf{z}_{C}$ belongs to
$Z\left(  \mathbf{k}\left[  G\right]  \right)  $.
\end{statement}

\begin{statement}
\textit{Claim 2:} The conjugacy class sums span the $\mathbf{k}$-module
$Z\left(  \mathbf{k}\left[  G\right]  \right)  $.
\end{statement}

\begin{statement}
\textit{Claim 3:} The conjugacy class sums are $\mathbf{k}$-linearly independent.
\end{statement}

\begin{proof}
[Proof of Claim 1.]Let $C$ be a finite conjugacy class of $G$. We must prove
that $\mathbf{z}_{C}\in Z\left(  \mathbf{k}\left[  G\right]  \right)  $.

We will prove this using Lemma \ref{lem.center.waw}. Recall that
$\mathbf{z}_{C}=\sum_{c\in C}c$ (by the definition of $\mathbf{z}_{C}$).

Let $w\in G$. Then, $wcw^{-1}\in C$ for each $c\in C$ (since the element
$wcw^{-1}$ is clearly conjugate to $c$, and thus belongs to the same conjugacy
class as $c$ does, but the latter class is $C$ of course). Hence, the map%
\begin{align*}
C  &  \rightarrow C,\\
c  &  \mapsto wcw^{-1}%
\end{align*}
is well-defined. Similarly, the map%
\begin{align*}
C  &  \rightarrow C,\\
c  &  \mapsto w^{-1}cw
\end{align*}
is well-defined. These two maps are mutually inverse (since each $c\in C$
satisfies $w^{-1}\left(  wcw^{-1}\right)  w=\underbrace{w^{-1}w}%
_{=1}c\underbrace{w^{-1}w}_{=1}=c$ and similarly $w\left(  w^{-1}cw\right)
w^{-1}=c$), and thus are bijections. In particular, the map%
\begin{align*}
C  &  \rightarrow C,\\
c  &  \mapsto wcw^{-1}%
\end{align*}
is a bijection. Hence, we can substitute $wcw^{-1}$ for $c$ in the sum
$\sum_{c\in C}c$. We thus obtain
\[
\sum_{c\in C}c=\sum_{c\in C}wcw^{-1}=w\left(  \sum_{c\in C}c\right)  w^{-1}.
\]
In view of $\mathbf{z}_{C}=\sum_{c\in C}c$, we can rewrite this as
$\mathbf{z}_{C}=w\mathbf{z}_{C}w^{-1}$. In other words, $w\mathbf{z}_{C}%
w^{-1}=\mathbf{z}_{C}$.

Forget that we fixed $w$. We thus have shown that each $w\in G$ satisfies
$w\mathbf{z}_{C}w^{-1}=\mathbf{z}_{C}$. According to Lemma
\ref{lem.center.waw} (applied to $\mathbf{a}=\mathbf{z}_{C}$), we thus
conclude that $\mathbf{z}_{C}\in Z\left(  \mathbf{k}\left[  G\right]  \right)
$. This proves Claim 1.
\end{proof}

\begin{proof}
[Proof of Claim 2.]Claim 1 shows that these conjugacy class sums
$\mathbf{z}_{C}$ all belong to $Z\left(  \mathbf{k}\left[  G\right]  \right)
$. It remains to show that each $\mathbf{a}\in Z\left(  \mathbf{k}\left[
G\right]  \right)  $ is a $\mathbf{k}$-linear combination of these sums. So
let us do this.

Let $\mathbf{a}\in Z\left(  \mathbf{k}\left[  G\right]  \right)  $ be
arbitrary. Recall that $\mathbf{k}\left[  G\right]  =\mathbf{k}^{\left(
G\right)  }$ as a $\mathbf{k}$-module (by the definition of a monoid algebra).
Hence, $\mathbf{a}\in Z\left(  \mathbf{k}\left[  G\right]  \right)
\subseteq\mathbf{k}\left[  G\right]  =\mathbf{k}^{\left(  G\right)  }%
\subseteq\mathbf{k}^{G}$. Hence, $\mathbf{a}$ is a family $\left(  \alpha
_{w}\right)  _{w\in G}$ of elements of $\mathbf{k}$. Consider these elements.
Thus, $\mathbf{a}=\left(  \alpha_{w}\right)  _{w\in G}=\sum_{w\in G}\alpha
_{w}w$ (as we know from Subsection \ref{subsec.intro.monalg.convs}). Moreover,
all but finitely many $w\in G$ satisfy $\alpha_{w}=0$ (since $\left(
\alpha_{w}\right)  _{w\in G}=\mathbf{a}\in\mathbf{k}^{\left(  G\right)  }$).
In other words, only finitely many $w\in G$ satisfy $\alpha_{w}\neq0$. We now
make the following observations:

\begin{itemize}
\item Any two conjugate elements $x$ and $y$ of $G$ satisfy
\begin{equation}
\alpha_{x}=\alpha_{y}. \label{pf.thm.center.conjsums.c2.pf.1}%
\end{equation}

[\textit{Proof:} Let $x$ and $y$ be two conjugate elements of $G$. Thus,
$y=gxg^{-1}$ for some $g\in G$. Consider this $g$. Now, $\mathbf{a}\in
Z\left(  \mathbf{k}\left[  G\right]  \right)  $, so that each $w\in G$
satisfies $w\mathbf{a}w^{-1}=\mathbf{a}$ (by Lemma \ref{lem.center.waw}).
Applying this to $w=g$, we obtain $g\mathbf{a}g^{-1}=\mathbf{a}=\sum
\limits_{w\in G}\alpha_{w}w$. Comparing this with%
\begin{align*}
g\mathbf{a}g^{-1}  &  =g\left(  \sum\limits_{w\in G}\alpha_{w}w\right)
g^{-1}\ \ \ \ \ \ \ \ \ \ \left(  \text{since }\mathbf{a}=\sum\limits_{w\in
G}\alpha_{w}w\right) \\
&  =\sum\limits_{w\in G}\alpha_{w}gwg^{-1}\\
&  =\sum\limits_{w\in G}\alpha_{g^{-1}wg}\underbrace{gg^{-1}}_{=1}%
w\underbrace{gg^{-1}}_{=1}\ \ \ \ \ \ \ \ \ \ \left(
\begin{array}
[c]{c}%
\text{here, we have substituted }gwg^{-1}\\
\text{for }w\text{ in the sum, since the}\\
\text{map }G\rightarrow G,\ w\mapsto gwg^{-1}\\
\text{is a bijection}%
\end{array}
\right) \\
&  =\sum\limits_{w\in G}\alpha_{g^{-1}wg}w,
\end{align*}
we obtain%
\[
\sum\limits_{w\in G}\alpha_{w}w=\sum\limits_{w\in G}\alpha_{g^{-1}wg}w.
\]
Now, the left hand side of this equality can be rewritten as the family
$\left(  \alpha_{w}\right)  _{w\in G}$ (since $\sum\limits_{w\in G}\alpha
_{w}w=\left(  \alpha_{w}\right)  _{w\in G}$), whereas the right hand side can
be rewritten as the family $\left(  \alpha_{g^{-1}wg}\right)  _{w\in G}$ (for
similar reasons). Thus, this equality rewrites as%
\[
\left(  \alpha_{w}\right)  _{w\in G}=\left(  \alpha_{g^{-1}wg}\right)  _{w\in
G}.
\]
In other words, $\alpha_{w}=\alpha_{g^{-1}wg}$ for each $w\in G$. Applying
this to $w=y$, we obtain $\alpha_{y}=\alpha_{g^{-1}yg}$. Since $g^{-1}%
\underbrace{y}_{=gxg^{-1}}g=\underbrace{g^{-1}g}_{=1}x\underbrace{g^{-1}%
g}_{=1}=x$, we can rewrite this as $\alpha_{y}=\alpha_{x}$. This proves
(\ref{pf.thm.center.conjsums.c2.pf.1}).]

\item For each conjugacy class $C$ of $G$, there exists an element $\beta_{C}$
of $\mathbf{k}$ such that%
\begin{equation}
\text{all }w\in C\text{ satisfy }\alpha_{w}=\beta_{C}.
\label{pf.thm.center.conjsums.c2.pf.2}%
\end{equation}

[\textit{Proof of (\ref{pf.thm.center.conjsums.c2.pf.2}):} Let $C$ be a
conjugacy class of $G$. Pick any element $u$ of $C$, and set $\beta
_{C}:=\alpha_{u}$. We shall now show that all $w\in C$ satisfy $\alpha
_{w}=\beta_{C}$. Indeed, for any $w\in C$, the elements $w$ and $u$ of $G$ are
conjugate (since they both belong to the same conjugacy class, namely $C$),
and thus satisfy $\alpha_{w}=\alpha_{u}$ (by
(\ref{pf.thm.center.conjsums.c2.pf.1}), applied to $x=w$ and $y=u$) and
therefore $\alpha_{w}=\alpha_{u}=\beta_{C}$. Thus, we have found an element
$\beta_{C}$ of $\mathbf{k}$ such that all $w\in C$ satisfy $\alpha_{w}%
=\beta_{C}$. This proves (\ref{pf.thm.center.conjsums.c2.pf.2}).]

\item If $C$ is an infinite conjugacy class of $G$, then the element
$\beta_{C}$ of $\mathbf{k}$ that satisfies
(\ref{pf.thm.center.conjsums.c2.pf.2}) is given by%
\begin{equation}
\beta_{C}=0. \label{pf.thm.center.conjsums.c2.pf.3}%
\end{equation}

[\textit{Proof:} Let $C$ be an infinite conjugacy class of $G$. Then, there
are infinitely many $w\in C$. If we had $\beta_{C}\neq0$, then all these
infinitely many $w\in C$ would satisfy $\alpha_{w}\neq0$ (since
(\ref{pf.thm.center.conjsums.c2.pf.2}) would yield $\alpha_{w}=\beta_{C}\neq
0$), which would contradict the fact that only finitely many $w\in G$ satisfy
$\alpha_{w}\neq0$. Hence, we cannot have $\beta_{C}\neq0$. Thus, $\beta_{C}%
=0$. This proves (\ref{pf.thm.center.conjsums.c2.pf.3}).]
\end{itemize}

For each conjugacy class $C$ of $G$, let us pick an element $\beta_{C}$ of
$\mathbf{k}$ that satisfies (\ref{pf.thm.center.conjsums.c2.pf.2}). (We have
already shown that such a $\beta_{C}$ exists.)

Now, recall that the conjugacy classes of $G$ are the equivalence classes of
an equivalence relation, and thus are disjoint nonempty subsets of $G$ whose
union is $G$. Thus, we can split the sum $\sum\limits_{w\in G}\alpha_{w}w$ as
follows:%
\begin{align*}
\sum\limits_{w\in G}\alpha_{w}w  &  =\sum_{\substack{C\text{ is a
conjugacy}\\\text{class of }G}}\ \ \sum\limits_{w\in C}\underbrace{\alpha_{w}%
}_{\substack{=\beta_{C}\\\text{(by (\ref{pf.thm.center.conjsums.c2.pf.2}))}%
}}w=\sum_{\substack{C\text{ is a conjugacy}\\\text{class of }G}}\ \ \sum
\limits_{w\in C}\beta_{C}w\\
&  =\sum_{\substack{C\text{ is a finite conjugacy}\\\text{class of }%
G}}\ \ \underbrace{\sum\limits_{w\in C}\beta_{C}w}_{=\beta_{C}\sum
\limits_{w\in C}w}+\sum_{\substack{C\text{ is an infinite conjugacy}%
\\\text{class of }G}}\ \ \sum\limits_{w\in C}\underbrace{\beta_{C}%
}_{\substack{=0\\\text{(by (\ref{pf.thm.center.conjsums.c2.pf.3}))}}}w\\
&  =\sum_{\substack{C\text{ is a finite conjugacy}\\\text{class of }G}%
}\beta_{C}\underbrace{\sum\limits_{w\in C}w}_{\substack{=\sum_{c\in
C}c=\mathbf{z}_{C}\\\text{(by the definition}\\\text{of }\mathbf{z}%
_{C}\text{)}}}+\underbrace{\sum_{\substack{C\text{ is an infinite
conjugacy}\\\text{class of }G}}\ \ \sum\limits_{w\in C}0w}_{=0}\\
&  =\sum_{\substack{C\text{ is a finite conjugacy}\\\text{class of }G}%
}\beta_{C}\mathbf{z}_{C}.
\end{align*}
This is clearly a $\mathbf{k}$-linear combination of the conjugacy class sums
$\mathbf{z}_{C}$. Hence, we have shown that $\sum\limits_{w\in G}\alpha_{w}w$
is a $\mathbf{k}$-linear combination of the conjugacy class sums
$\mathbf{z}_{C}$. In other words, $\mathbf{a}$ is a $\mathbf{k}$-linear
combination of the conjugacy class sums $\mathbf{z}_{C}$ (since $\mathbf{a}%
=\sum\limits_{w\in G}\alpha_{w}w$). This completes our proof of Claim 2.
\end{proof}

\begin{proof}
[Proof of Claim 3.]This follows easily from the fact that the conjugacy
classes $C$ are nonempty and disjoint. Here are the details:

We must prove that the family $\left(  \mathbf{z}_{C}\right)  _{C\in\left\{
\text{finite conjugacy classes}\right\}  }$ is $\mathbf{k}$-linearly
independent. In other words, we must prove that if we have an essentially
finite family of scalars $\left(  \beta_{C}\right)  _{C\in\left\{
\text{finite conjugacy classes}\right\}  }$ that satisfies%
\begin{equation}
\sum_{\substack{C\text{ is a finite}\\\text{conjugacy class}}}\beta
_{C}\mathbf{z}_{C}=0, \label{pf.thm.center.conjsums.c3.pf.ass}%
\end{equation}
then all these scalars $\beta_{C}$ are $0$. So let us assume that we have an
essentially finite family of scalars $\left(  \beta_{C}\right)  _{C\in\left\{
\text{finite conjugacy classes}\right\}  }$ that satisfies
(\ref{pf.thm.center.conjsums.c3.pf.ass}). We must prove that all these scalars
$\beta_{C}$ are $0$.

Let us extend the definition of $\beta_{C}$ to \textbf{all} conjugacy classes
$C$ (not just the finite ones), by setting
\begin{equation}
\beta_{C}:=0\ \ \ \ \ \ \ \ \ \ \text{for all infinite conjugacy classes }C.
\label{pf.thm.center.conjsums.c3.pf.inf}%
\end{equation}
For any $g\in G$, let $K\left(  g\right)  $ denote the conjugacy class of $G$
that contains $g$. Thus, if $g$ is an element of a conjugacy class $C$, then
\begin{equation}
C=K\left(  g\right)  . \label{pf.thm.center.conjsums.c3.pf.C=Kg}%
\end{equation}

From (\ref{pf.thm.center.conjsums.c3.pf.ass}), we obtain%
\begin{align*}
0  &  =\sum_{\substack{C\text{ is a finite}\\\text{conjugacy class}}}\beta
_{C}\underbrace{\mathbf{z}_{C}}_{\substack{=\sum_{c\in C}c\\\text{(by the
definition of }\mathbf{z}_{C}\text{)}}}=\sum_{\substack{C\text{ is a
finite}\\\text{conjugacy class}}}\beta_{C}\sum_{c\in C}c\\
&  =\sum_{\substack{C\text{ is a finite}\\\text{conjugacy class}}}\beta
_{C}\sum_{g\in C}g=\sum_{\substack{C\text{ is a finite}\\\text{conjugacy
class}}}\ \ \sum_{g\in C}\beta_{C}g\\
&  =\sum_{\substack{C\text{ is a}\\\text{conjugacy class}}}\ \ \sum_{g\in
C}\underbrace{\beta_{C}}_{\substack{=\beta_{K\left(  g\right)  }\\\text{(since
}C=K\left(  g\right)  \\\text{(by (\ref{pf.thm.center.conjsums.c3.pf.C=Kg}%
)))}}}g\\
&  \ \ \ \ \ \ \ \ \ \ \ \ \ \ \ \ \ \ \ \ \left(
\begin{array}
[c]{c}%
\text{here, we have extended the range of the outer sum}\\
\text{from }\left\{  \text{finite conjugacy classes}\right\}  \text{ to
}\left\{  \text{all conjugacy classes}\right\}  \text{,}\\
\text{which does not change the sum}\\
\text{because (\ref{pf.thm.center.conjsums.c3.pf.inf}) shows that all the new
addends are }0
\end{array}
\right) \\
&  =\underbrace{\sum_{\substack{C\text{ is a}\\\text{conjugacy class}%
}}\ \ \sum_{g\in C}}_{\substack{=\sum_{\substack{g\in G}}\\\text{(since each
}g\in G\\\text{belongs to exactly one}\\\text{conjugacy class)}}%
}\beta_{K\left(  g\right)  }g=\sum_{\substack{g\in G}}\beta_{K\left(
g\right)  }g=\left(  \beta_{K\left(  g\right)  }\right)  _{g\in G}%
\end{align*}
(since $\sum_{g\in G}\alpha_{g}g=\left(  \alpha_{g}\right)  _{g\in G}$ for any
essentially finite family $\left(  \alpha_{g}\right)  _{g\in G}$ of scalars).
In other words,%
\[
\left(  \beta_{K\left(  g\right)  }\right)  _{g\in G}=0=\left(  0\right)
_{g\in G}.
\]
In other words, $\beta_{K\left(  g\right)  }=0$ for each $g\in G$.

Now, let $C$ be any finite conjugacy class. Then, $C$ is nonempty, so that
there exists some $w\in C$. Consider this $w$. Then, $C=K\left(  w\right)  $
(by (\ref{pf.thm.center.conjsums.c3.pf.C=Kg}), applied to $g=w$). But we
showed that $\beta_{K\left(  g\right)  }=0$ for each $g\in G$. Applying this
to $g=w$, we obtain $\beta_{K\left(  w\right)  }=0$. In other words,
$\beta_{C}=0$ (since $C=K\left(  w\right)  $).

Forget that we fixed $C$. We thus have shown that $\beta_{C}=0$ for all finite
conjugacy classes $C$. In other words, all the scalars $\beta_{C}$ are $0$.
This completes the proof of Claim 3.
\end{proof}

Having proved all three claims, we are now done proving Theorem
\ref{thm.center.conjsums}.
\end{proof}

\begin{exercise}
\label{exe.center.waw-sums}Let $G$ be a finite group. Prove the following:
\medskip

\textbf{(a)} \fbox{1} For each $\mathbf{a}\in\mathbf{k}\left[  G\right]  $, we
have $\sum_{w\in G}w\mathbf{a}w^{-1}\in Z\left(  \mathbf{k}\left[  G\right]
\right)  $. \medskip

\textbf{(b)} \fbox{2} For each $\mathbf{a},\mathbf{b}\in\mathbf{k}\left[
G\right]  $, we have $\sum_{w\in G}w\mathbf{ab}w^{-1}=\sum_{w\in
G}w\mathbf{ba}w^{-1}$.
\end{exercise}

\begin{exercise}
Let $G$ be a group. For any $\mathbf{a}=\left(  \alpha_{g}\right)  _{g\in
G}\in\mathbf{k}\left[  G\right]  $ and any element $h\in G$, we let $\left[
h\right]  \mathbf{a}$ denote the scalar $\alpha_{h}$ (that is, the $h$-th
entry of the family $\mathbf{a}=\left(  \alpha_{g}\right)  _{g\in G}$); this
is also known as the $h$\emph{-coefficient} of $\mathbf{a}$. \medskip

\textbf{(a)} \fbox{1} Prove that every $\mathbf{a},\mathbf{b}\in
\mathbf{k}\left[  G\right]  $ satisfy $\left[  1\right]  \left(
\mathbf{ab}\right)  =\left[  1\right]  \left(  \mathbf{ba}\right)  $ (even if
$\mathbf{ab}$ and $\mathbf{ba}$ are not necessarily equal). Here, as usual,
$1$ denotes the neutral element of $G$. \medskip

\textbf{(b)} \fbox{1} Let $h\in G$. Prove that the equality $\left[  h\right]
\left(  \mathbf{ab}\right)  =\left[  h\right]  \left(  \mathbf{ba}\right)  $
holds for all $\mathbf{a},\mathbf{b}\in\mathbf{k}\left[  G\right]  $ if and
only if $h\in Z\left(  G\right)  $ or the ring $\mathbf{k}$ is trivial (i.e.,
we have $0=1$ in $\mathbf{k}$).
\end{exercise}

\subsubsection{Symmetric groups}

So much for conjugacy in arbitrary groups. Now, let us find out what conjugacy
means in the symmetric group $S_{n}$ and, more generally, in $S_{X}$ for $X$ finite.

Recall that a \emph{partition} (or, to be more precise, \emph{integer
partition}) means a weakly decreasing finite tuple of positive integers. For
instance, $\left(  5,2,2,1\right)  $ is a partition, whereas $\left(
5,1,2\right)  $ and $\left(  5,2,0\right)  $ are not. The \emph{size} of a
partition $\lambda=\left(  \lambda_{1},\lambda_{2},\ldots,\lambda_{k}\right)
$ is defined to be the sum $\lambda_{1}+\lambda_{2}+\cdots+\lambda_{k}%
\in\mathbb{N}$. A \emph{partition of }$n$ (for a given number $n\in\mathbb{Z}%
$) means a partition having size $n$. The empty tuple $\left(  {}\right)  $ is
a partition of $0$. See \cite[Chapter 4]{21s} for more about partitions.

\begin{definition}
\label{def.type.type}Let $X$ be a finite set. Let $w\in S_{X}$. Then, the
\emph{cycle type} of $w$ is the partition whose entries are the sizes of the
orbits of $w$ (from largest to smallest). Of course, you can replace
\textquotedblleft sizes of the orbits\textquotedblright\ by \textquotedblleft
lengths of the cycles\textquotedblright\ in this definition without changing
the outcome (because of Proposition \ref{prop.dcd.orbs-from-cycs}).

We denote the cycle type of $w$ by $\operatorname*{type}w$.
\end{definition}

\begin{example}
Let $w=\operatorname*{oln}\left(  432651\right)  \in S_{6}$. Then,
$\operatorname*{type}w=\left(  3,2,1\right)  $, since the orbits of $w$ are
the sets $\left\{  1,4,6\right\}  $, $\left\{  2,3\right\}  $ and $\left\{
5\right\}  $ with respective sizes $3$, $2$ and $1$. (This $w$ is the
permutation we analyzed in Example \ref{exa.orbits.orbits-1}.)
\end{example}

\begin{example}
Let $w=\operatorname*{oln}\left(  426153\right)  \in S_{6}$. Then, the cycle
digraph of $w$ looks as follows:%
\[%
%TCIMACRO{\TeXButton{tikz cycle digraph}{\begin{tikzpicture}%
%[->,shorten >=1pt,auto,node distance=3cm, thick,main node/.style={circle,fill=blue!20,draw}%
%]
%\node[main node] (1) {1};
%\node[main node] (4) [below of=1] {4};
%\node[main node] (3) [right of=1] {3};
%\node[main node] (6) [right of=4] {6};
%\node[main node] (2) [right of=3] {2};
%\node[main node] (5) [right of=6] {5};
%\path[-{Stealth[length=4mm]}]
%(1) edge [bend right] (4)
%(4) edge [bend right] (1)
%(3) edge [bend right] (6)
%(6) edge [bend right] (3)
%(2) edge [loop right] (2)
%(5) edge [loop right] (5);
%\end{tikzpicture}}}%
%BeginExpansion
\begin{tikzpicture}%
[->,shorten >=1pt,auto,node distance=3cm, thick,main node/.style={circle,fill=blue!20,draw}%
]
\node[main node] (1) {1};
\node[main node] (4) [below of=1] {4};
\node[main node] (3) [right of=1] {3};
\node[main node] (6) [right of=4] {6};
\node[main node] (2) [right of=3] {2};
\node[main node] (5) [right of=6] {5};
\path[-{Stealth[length=4mm]}]
(1) edge [bend right] (4)
(4) edge [bend right] (1)
(3) edge [bend right] (6)
(6) edge [bend right] (3)
(2) edge [loop right] (2)
(5) edge [loop right] (5);
\end{tikzpicture}%
%EndExpansion
\ \ .
\]
Thus, the orbits of $w$ are the sets $\left\{  1,4\right\}  ,\ \left\{
2\right\}  ,\ \left\{  3,6\right\}  ,\ \left\{  5\right\}  $ with respective
sizes $2$, $1$, $2$ and $1$. Hence, the cycle type of $w$ is
$\operatorname*{type}w=\left(  2,2,1,1\right)  $ (don't forget to arrange the
entries in decreasing order).
\end{example}

Note that the cycle type of a permutation $w\in S_{X}$ is always a partition
of $\left\vert X\right\vert $.

Using cycle types, we can now characterize conjugacy in a symmetric group
$S_{X}$ (for finite $X$): Two permutations $u,v\in S_{X}$ are conjugate if and
only if they have the same cycle type. In other words:

\begin{theorem}
\label{thm.type.conj=type}Let $X$ be a finite set. Let $u,v\in S_{X}$. Then,
$u\sim v$ in $S_{X}$ if and only if $\operatorname*{type}%
u=\operatorname*{type}v$.
\end{theorem}

\begin{proof}
[Proof sketch.]$\Longrightarrow:$ Assume that $u\sim v$ in $S_{X}$. Then,
$u=pvp^{-1}$ for some $p\in S_{X}$ (by the definition of conjugacy). Consider
this $p$. We must show that $\operatorname*{type}u=\operatorname*{type}v$.

Intuitively, this is pretty clear: The permutation $u=pvp^{-1}$ is obtained
from $v$ by \textquotedblleft relabelling\textquotedblright\ the elements $x$
of $X$ as $p\left(  x\right)  $, and it is obvious that the sizes of the
orbits do not change under this \textquotedblleft
relabelling\textquotedblright\ operation. Thus, $\operatorname*{type}%
u=\operatorname*{type}v$.

To formalize this argument, we can argue as follows: Note that $p$ is a
permutation, and thus an injective map; hence, any subset $U$ of $X$ satisfies%
\begin{equation}
\left\vert p\left(  U\right)  \right\vert =\left\vert U\right\vert
\label{pf.thm.type.conj=type.1}%
\end{equation}
(where $p\left(  U\right)  $ means the set $\left\{  p\left(  x\right)
\ \mid\ x\in U\right\}  $).

Now, if $O$ is an orbit of $v$, then $p\left(  O\right)  $ is an orbit of
$pvp^{-1}$ (this is straightforward to check: just show that two elements
$x,y\in X$ satisfy $x\overset{v}{\sim}y$ if and only if $p\left(  x\right)
\overset{pvp^{-1}}{\sim}p\left(  y\right)  $). Since $p$ is bijective, this
construction yields all orbits of $pvp^{-1}$. Thus, the orbits of $pvp^{-1}$
are just the images of the orbits of $v$ under application of the map $p$.
Since this application does not change the size of an orbit (by
(\ref{pf.thm.type.conj=type.1})), we thus conclude that the orbits of
$pvp^{-1}$ have the same sizes as the orbits of $v$. In other words,
$\operatorname*{type}\left(  pvp^{-1}\right)  =\operatorname*{type}v$. Since
$pvp^{-1}=u$, we can rewrite this as $\operatorname*{type}%
u=\operatorname*{type}v$. This proves the \textquotedblleft$\Longrightarrow
$\textquotedblright\ direction of Theorem \ref{thm.type.conj=type}.\medskip

$\Longleftarrow:$ Assume that $\operatorname*{type}u=\operatorname*{type}v$.
Let us write the partition $\operatorname*{type}u=\operatorname*{type}v$ as
$\left(  n_{1},n_{2},\ldots,n_{k}\right)  $. Thus, the cycle type of $u$ is
$\operatorname*{type}u=\left(  n_{1},n_{2},\ldots,n_{k}\right)  $. In other
words, the cycles of $u$ have lengths $n_{1},n_{2},\ldots,n_{k}$ (since the
cycle type of a permutation consists of the lengths of its cycles). Hence, the
permutation $u$ has a DCD%
\[
u=\operatorname*{cyc}\nolimits_{a_{1,1},a_{1,2},\ldots,a_{1,n_{1}}}%
\circ\operatorname*{cyc}\nolimits_{a_{2,1},a_{2,2},\ldots,a_{2,n_{2}}}%
\circ\cdots\circ\operatorname*{cyc}\nolimits_{a_{k,1},a_{k,2},\ldots
,a_{k,n_{k}}}%
\]
(since disjoint cycles commute and thus can be arranged in any order).
Similarly, the permutation $v$ has a DCD%
\[
v=\operatorname*{cyc}\nolimits_{b_{1,1},b_{1,2},\ldots,b_{1,n_{1}}}%
\circ\operatorname*{cyc}\nolimits_{b_{2,1},b_{2,2},\ldots,b_{2,n_{2}}}%
\circ\cdots\circ\operatorname*{cyc}\nolimits_{b_{k,1},b_{k,2},\ldots
,b_{k,n_{k}}}.
\]
Now, let $p$ be the permutation in $S_{X}$ that sends each $a_{i,j}$ to the
corresponding $b_{i,j}$. (This permutation exists because the definition of a
DCD guarantees that each element of $X$ appears exactly once in the composite
list
\begin{align*}
&  (a_{1,1},a_{1,2},\ldots,a_{1,n_{1}},\\
&  \ \ \ a_{2,1},a_{2,2},\ldots,a_{2,n_{2}},\\
&  \ \ \ \ldots,\\
&  \ \ \ a_{k,1},a_{k,2},\ldots,a_{k,n_{k}})
\end{align*}
and exactly once in the composite list%
\begin{align*}
&  (b_{1,1},b_{1,2},\ldots,b_{1,n_{1}},\\
&  \ \ \ b_{2,1},b_{2,2},\ldots,b_{2,n_{2}},\\
&  \ \ \ \ldots,\\
&  \ \ \ b_{k,1},b_{k,2},\ldots,b_{k,n_{k}}).
\end{align*}
) Then, it is easy to see that $v=pup^{-1}$ (indeed, just check that both maps
$v$ and $pup^{-1}$ send each $b_{i,j}$ to $b_{i,j+1}$, where $b_{i,n_{i}+1}$
is defined to be $b_{i,1}$\ \ \ \ \footnote{Here are the details: Set
$b_{i,n_{i}+1}:=b_{i,1}$ and $a_{i,n_{i}+1}:=a_{i,1}$ for each $i\in\left[
k\right]  $. Then, for each $i\in\left[  k\right]  $ and $j\in\left[
n_{i}\right]  $, we have $u\left(  a_{i,j}\right)  =a_{i,j+1}$ (since
$u=\operatorname*{cyc}\nolimits_{a_{1,1},a_{1,2},\ldots,a_{1,n_{1}}}%
\circ\operatorname*{cyc}\nolimits_{a_{2,1},a_{2,2},\ldots,a_{2,n_{2}}}%
\circ\cdots\circ\operatorname*{cyc}\nolimits_{a_{k,1},a_{k,2},\ldots
,a_{k,n_{k}}}$) and $v\left(  b_{i,j}\right)  =b_{i,j+1}$ (for similar
reasons) and $p\left(  a_{i,j}\right)  =b_{i,j}$ (by the definition of $p$)
and $p\left(  a_{i,j+1}\right)  =b_{i,j+1}$ (by the definition of $p$, since
$b_{i,n_{i}+1}=b_{i,1}$ and $a_{i,n_{i}+1}=a_{i,1}$), so that%
\[
\left(  pup^{-1}\right)  \left(  b_{i,j}\right)  =p\left(  u\left(
\underbrace{p^{-1}\left(  b_{i,j}\right)  }_{\substack{=a_{i,j}\\\text{(since
}p\left(  a_{i,j}\right)  =b_{i,j}\text{)}}}\right)  \right)  =p\left(
\underbrace{u\left(  a_{i,j}\right)  }_{=a_{i,j+1}}\right)  =p\left(
a_{i,j+1}\right)  =b_{i,j+1}=v\left(  b_{i,j}\right)  .
\]
Since each element of $X$ has the form $b_{i,j}$ for some $i$ and $j$, we thus
have shown that $\left(  pup^{-1}\right)  \left(  x\right)  =v\left(
x\right)  $ for each $x\in X$. In other words, $pup^{-1}=v$, qed.}). Thus,
$v\sim u$, so that $u\sim v$. This proves the \textquotedblleft$\Longleftarrow
$\textquotedblright\ direction of Theorem \ref{thm.type.conj=type}.
\end{proof}

\begin{corollary}
\label{cor.type.conjclasses}Let $X$ be a finite set. \medskip

\textbf{(a)} For each partition $\lambda$ of $\left\vert X\right\vert $, the
set%
\[
C_{\lambda}\left(  X\right)  :=\left\{  w\in S_{X}\ \mid\ \operatorname*{type}%
w=\lambda\right\}
\]
(consisting of all permutations in $S_{X}$ having cycle type $\lambda$) is a
conjugacy class of $S_{X}$. \medskip

\textbf{(b)} Each conjugacy class of $S_{X}$ has the form $C_{\lambda}\left(
X\right)  $ for some partition $\lambda$ of $\left\vert X\right\vert $.
\medskip

\textbf{(c)} The \# of conjugacy classes of $S_{X}$ equals the \# of
partitions of $\left\vert X\right\vert $.
\end{corollary}

\begin{proof}
[Proof sketch.]\textbf{(a)} Let $\lambda$ be a partition of $\left\vert
X\right\vert $. Write $\lambda$ as $\lambda=\left(  \lambda_{1},\lambda
_{2},\ldots,\lambda_{k}\right)  $, so that $\lambda_{1}+\lambda_{2}%
+\cdots+\lambda_{k}=\left(  \text{size of }\lambda\right)  =\left\vert
X\right\vert $. Hence, we can subdivide the set $X$ into $k$ disjoint subsets
\begin{align*}
L_{1}  &  =\left\{  x_{1,1},x_{1,2},\ldots,x_{1,\lambda_{1}}\right\}  ,\\
L_{2}  &  =\left\{  x_{2,1},x_{2,2},\ldots,x_{2,\lambda_{2}}\right\}  ,\\
&  \ldots,\\
L_{k}  &  =\left\{  x_{k,1},x_{k,2},\ldots,x_{k,\lambda_{k}}\right\}
\end{align*}
of sizes $\lambda_{1},\lambda_{2},\ldots,\lambda_{k}$ (respectively).
Obviously, neither these $k$ subsets nor the specific elements $x_{i,j}$ are
uniquely determined, but we just make some choice and stick with it. Now, the
permutation%
\[
u:=\operatorname*{cyc}\nolimits_{x_{1,1},x_{1,2},\ldots,x_{1,\lambda_{1}}%
}\circ\operatorname*{cyc}\nolimits_{x_{2,1},x_{2,2},\ldots,x_{2,\lambda_{2}}%
}\circ\cdots\circ\operatorname*{cyc}\nolimits_{x_{k,1},x_{k,2},\ldots
,x_{k,\lambda_{k}}}\in S_{X}%
\]
has cycle type $\lambda$, and thus belongs to $C_{\lambda}\left(  X\right)  $.
Hence, the set $C_{\lambda}\left(  X\right)  $ is nonempty.

Moreover, $\operatorname*{type}u=\lambda$ (since $u$ has cycle type $\lambda
$). The conjugacy class of the group $S_{X}$ that contains $u$ is%
\begin{align*}
&  \left\{  w\in S_{X}\ \mid\ w\sim u\right\} \\
&  =\left\{  w\in S_{X}\ \mid\ \operatorname*{type}w=\operatorname*{type}%
u\right\} \\
&  \ \ \ \ \ \ \ \ \ \ \ \ \ \ \ \ \ \ \ \ \left(
\begin{array}
[c]{c}%
\text{since Theorem \ref{thm.type.conj=type} shows that}\\
\text{we have }w\sim u\text{ if and only if }\operatorname*{type}%
w=\operatorname*{type}u
\end{array}
\right) \\
&  =\left\{  w\in S_{X}\ \mid\ \operatorname*{type}w=\lambda\right\}
\ \ \ \ \ \ \ \ \ \ \left(  \text{since }\operatorname*{type}u=\lambda\right)
\\
&  =C_{\lambda}\left(  X\right)  \ \ \ \ \ \ \ \ \ \ \left(  \text{by the
definition of }C_{\lambda}\left(  X\right)  \right)  .
\end{align*}
This shows that $C_{\lambda}\left(  X\right)  $ is a conjugacy class of
$S_{X}$. This proves Corollary \ref{cor.type.conjclasses} \textbf{(a)}.
\medskip

\textbf{(b)} Corollary \ref{cor.type.conjclasses} \textbf{(a)} shows that
$C_{\lambda}\left(  X\right)  $ is a conjugacy class of $S_{X}$ whenever
$\lambda$ is a partition of $\left\vert X\right\vert $. Clearly, these
conjugacy classes $C_{\lambda}\left(  X\right)  $ for different partitions
$\lambda$ are disjoint (since a permutation $w\in S_{X}$ cannot have two
different cycle types at the same time) and thus distinct. Moreover, these
classes $C_{\lambda}\left(  X\right)  $ cover the entire set $S_{X}$ (since
each permutation $w\in S_{X}$ has cycle type $\lambda$ for some partition
$\lambda$ of $\left\vert X\right\vert $, and thus belongs to the corresponding
class $C_{\lambda}\left(  X\right)  $). Hence, they are all the conjugacy
classes of $S_{X}$. In particular, each conjugacy class of $S_{X}$ has the
form $C_{\lambda}\left(  X\right)  $ for some partition $\lambda$ of
$\left\vert X\right\vert $. This proves Corollary \ref{cor.type.conjclasses}
\textbf{(b)}. \medskip

\textbf{(c)} Combining parts \textbf{(a)} and \textbf{(b)} of Corollary
\ref{cor.type.conjclasses}, we see that the conjugacy classes of $S_{X}$ are
precisely the sets $C_{\lambda}\left(  X\right)  $ for all partitions
$\lambda$ of $\left\vert X\right\vert $. But the latter sets $C_{\lambda
}\left(  X\right)  $ are distinct (as we saw in the above proof of part
\textbf{(b)}), and thus their \# equals the \# of partitions of $\left\vert
X\right\vert $. Hence, we conclude that the \# of conjugacy classes of $S_{X}$
equals the \# of partitions of $\left\vert X\right\vert $. This proves
Corollary \ref{cor.type.conjclasses} \textbf{(c)}.
\end{proof}

For example:

\begin{itemize}
\item The partitions of $3$ are $\left(  1,1,1\right)  $, $\left(  2,1\right)
$ and $\left(  3\right)  $. Hence, Corollary \ref{cor.type.conjclasses}
\textbf{(b)} shows that the conjugacy classes of $S_{3}=S_{\left[  3\right]
}$ are%
\begin{align*}
C_{\left(  1,1,1\right)  }\left(  \left[  3\right]  \right)   &  =\left\{
w\in S_{3}\ \mid\ \operatorname*{type}w=\left(  1,1,1\right)  \right\}
=\left\{  \operatorname*{id}\right\}  ,\\
C_{\left(  2,1\right)  }\left(  \left[  3\right]  \right)   &  =\left\{  w\in
S_{3}\ \mid\ \operatorname*{type}w=\left(  2,1\right)  \right\}  =\left\{
t_{1,2},\ t_{1,3},\ t_{2,3}\right\}  ,\\
C_{\left(  3\right)  }\left(  \left[  3\right]  \right)   &  =\left\{  w\in
S_{3}\ \mid\ \operatorname*{type}w=\left(  3\right)  \right\}  =\left\{
\operatorname*{cyc}\nolimits_{1,2,3},\ \operatorname*{cyc}\nolimits_{1,3,2}%
\right\}  .
\end{align*}

\item The partitions of $4$ are $\left(  1,1,1,1\right)  $, $\left(
2,1,1\right)  $, $\left(  2,2\right)  $, $\left(  3,1\right)  $ and $\left(
4\right)  $. Hence, Corollary \ref{cor.type.conjclasses} \textbf{(b)} shows
that the conjugacy classes of $S_{4}=S_{\left[  4\right]  }$ are%
\begin{align*}
C_{\left(  1,1,1,1\right)  }\left(  \left[  4\right]  \right)   &  =\left\{
w\in S_{4}\ \mid\ \operatorname*{type}w=\left(  1,1,1,1\right)  \right\}
=\left\{  \operatorname*{id}\right\}  ,\\
C_{\left(  2,1,1\right)  }\left(  \left[  4\right]  \right)   &  =\left\{
w\in S_{4}\ \mid\ \operatorname*{type}w=\left(  2,1,1\right)  \right\}
=\left\{  t_{i,j}\ \mid\ i<j\right\} \\
&  =\left\{  \text{transpositions in }S_{4}\right\}  ,\\
C_{\left(  2,2\right)  }\left(  \left[  4\right]  \right)   &  =\left\{  w\in
S_{4}\ \mid\ \operatorname*{type}w=\left(  2,2\right)  \right\}  =\left\{
t_{i,j}t_{u,v}\ \mid\ i,j,u,v\text{ distinct}\right\} \\
&  =\left\{  w\in S_{4}\ \mid\ w\text{ is an involution and has no fixed
points}\right\}  ,\\
C_{\left(  3,1\right)  }\left(  \left[  4\right]  \right)   &  =\left\{  w\in
S_{4}\ \mid\ \operatorname*{type}w=\left(  3,1\right)  \right\}  =\left\{
\operatorname*{cyc}\nolimits_{i,j,k}\ \mid\ i,j,k\text{ distinct}\right\}  ,\\
C_{\left(  4\right)  }\left(  \left[  4\right]  \right)   &  =\left\{  w\in
S_{4}\ \mid\ \operatorname*{type}w=\left(  4\right)  \right\}  =\left\{
\operatorname*{cyc}\nolimits_{i,j,k,\ell}\ \mid\ i,j,k,\ell\text{
distinct}\right\}  .
\end{align*}

\end{itemize}

In general, for each $n\in\mathbb{N}$, Corollary \ref{cor.type.conjclasses}
\textbf{(b)} shows that there is one conjugacy class of $S_{n}$ for each
partition of $n$, and it consists of all permutations $w\in S_{n}$ whose cycle
type is this partition. Thus, the center $Z\left(  \mathbf{k}\left[
S_{n}\right]  \right)  $ is a free $\mathbf{k}$-module with a basis indexed by
the partitions of $n$ (by Theorem \ref{thm.center.conjsums}).

The set of all transpositions in $S_{n}$ is precisely the conjugacy class%
\[
\left\{  w\in S_{n}\ \mid\ \operatorname*{type}w=\left(
2,\underbrace{1,1,\ldots,1}_{n-2\text{ ones}}\right)  \right\}
\ \ \ \ \ \ \ \ \ \ \left(  \text{when }n\geq2\right)  .
\]
Thus, we obtain Proposition \ref{prop.center.sum-mi} as a particular case.

\subsubsection{Centers and jucys--murphies}

Having described the center of $\mathbf{k}\left[  S_{n}\right]  $ explicitly
in terms of the conjugacy class sums, we can consider it well-understood, but
some questions still remain: For instance, what is it \textbf{as an algebra}
(i.e., how does its multiplication work)? Even the question \textquotedblleft
does a given element belong to the center?\textquotedblright\ is not
completely straightforward if the coefficients of this element are not known
explicitly. Thus, we have some more to say.

First, we observe something simple: The integral $\nabla$ and the
sign-integral $\nabla^{-}$ belong to the center of $\mathbf{k}\left[
S_{n}\right]  $. This follows directly from Proposition
\ref{prop.integral.fix}, although it can also be derived from Lemma
\ref{lem.center.waw}. The jucys--murphies $\mathbf{m}_{1},\mathbf{m}%
_{2},\ldots,\mathbf{m}_{n}$ commute with each other, but (except for
$\mathbf{m}_{1}=0$ and for $\mathbf{m}_{2}$ when $n=2$) do not belong to the
center of $\mathbf{k}\left[  S_{n}\right]  $. However, their sum
$\mathbf{m}_{1}+\mathbf{m}_{2}+\cdots+\mathbf{m}_{n}$ does (as we saw in
Proposition \ref{prop.center.sum-mi}), and so do their elementary symmetric
polynomials $e_{k}\left(  \mathbf{m}_{1},\mathbf{m}_{2},\ldots,\mathbf{m}%
_{n}\right)  $ that we computed in Corollary \ref{cor.YJM.ek} (since the
number of orbits of a permutation $w$ depends only on the conjugacy class of
$w$).

We can generalize this latter result somewhat. Recall that a polynomial
$f\in\mathbf{k}\left[  x_{1},x_{2},\ldots,x_{n}\right]  $ in $n$
indeterminates $x_{1},x_{2},\ldots,x_{n}$ is said to be \emph{symmetric} if
each permutation $\sigma\in S_{n}$ satisfies%
\[
f\left(  x_{\sigma\left(  1\right)  },x_{\sigma\left(  2\right)  }%
,\ldots,x_{\sigma\left(  n\right)  }\right)  =f.
\]
The elementary symmetric polynomials $e_{k}\left(  x_{1},x_{2},\ldots
,x_{n}\right)  $ are examples of symmetric polynomials, but there are others,
such as the \emph{power sum symmetric polynomials} $x_{1}^{k}+x_{2}^{k}%
+\cdots+x_{n}^{k}$ for all $k\in\mathbb{N}$. See \cite[Chapter 7]{21s} for an
introduction to symmetric polynomials. Now we can give another
characterization of $Z\left(  \mathbf{k}\left[  S_{n}\right]  \right)  $
(\cite[Theorem 1.9]{Murphy83}, \cite[Theorem 1']{Moran92}, \cite[Theorem
4.4.5]{CSScTo10}):

\begin{theorem}
[Jucys--Murphy theorem]\label{thm.center.murphies}The center $Z\left(
\mathbf{k}\left[  S_{n}\right]  \right)  $ is the set of all symmetric
polynomials in $n$ variables, applied to the jucys--murphies $\mathbf{m}%
_{1},\mathbf{m}_{2},\ldots,\mathbf{m}_{n}$. That is,%
\[
Z\left(  \mathbf{k}\left[  S_{n}\right]  \right)  =\left\{  f\left(
\mathbf{m}_{1},\mathbf{m}_{2},\ldots,\mathbf{m}_{n}\right)  \ \mid
\ f\in\mathbf{k}\left[  x_{1},x_{2},\ldots,x_{n}\right]  \text{ is
symmetric}\right\}  .
\]

\end{theorem}

\begin{proof}
[Partial proof of Theorem \ref{thm.center.murphies} (sketched).]Let us show
the \textquotedblleft$\supseteq$\textquotedblright\ part (the
\textquotedblleft$\subseteq$\textquotedblright\ part is much harder). So we
need to prove that
\[
f\left(  \mathbf{m}_{1},\mathbf{m}_{2},\ldots,\mathbf{m}_{n}\right)  \in
Z\left(  \mathbf{k}\left[  S_{n}\right]  \right)
\]
for any symmetric polynomial $f\in\mathbf{k}\left[  x_{1},x_{2},\ldots
,x_{n}\right]  $.

A theorem of Gauss (known as the \textquotedblleft fundamental theorem of
symmetric polynomials\textquotedblright)\footnote{Proofs of this theorem can
be found in \cite[proof of Theorem 1]{BluCos16}, in \cite[Theorem
1.2.1]{Dumas08}, in \cite[Theorem 9.6.6]{Goodman}, in \cite[\S 1.1]{Smith95}
and in many other sources. (Some of these sources only state the result in the
case when $\mathbf{k}$ is a field, or when $\mathbf{k}=\mathbb{C}$; but the
same proof applies more generally for any commutative ring $\mathbf{k}$.)}
says that each symmetric polynomial in $\mathbf{k}\left[  x_{1},x_{2}%
,\ldots,x_{n}\right]  $ can be written as a polynomial in the elementary
symmetric polynomials $e_{1},e_{2},\ldots,e_{n}$. For instance,%
\[
x_{1}^{2}+x_{2}^{2}+\cdots+x_{n}^{2}=\left(  \underbrace{x_{1}+x_{2}%
+\cdots+x_{n}}_{=e_{1}}\right)  ^{2}-2\underbrace{\sum_{i<j}x_{i}x_{j}%
}_{=e_{2}}=e_{1}^{2}-2e_{2}%
\]
and%
\begin{align*}
&  x_{1}^{3}+x_{2}^{3}+\cdots+x_{n}^{3}\\
&  =\left(  \underbrace{x_{1}+x_{2}+\cdots+x_{n}}_{=e_{1}}\right)
^{3}-3\underbrace{\sum_{i<j}\left(  x_{i}^{2}x_{j}+x_{i}x_{j}^{2}\right)
}_{=e_{1}e_{2}-3e_{3}}-\,6\underbrace{\sum_{i<j<k}x_{i}x_{j}x_{k}}_{=e_{3}}\\
&  =e_{1}^{3}-3\left(  e_{1}e_{2}-3e_{3}\right)  -6e_{3}=e_{1}^{3}-3e_{1}%
e_{2}+3e_{3}.
\end{align*}

Thanks to this theorem, it suffices to show that
\[
e_{k}\left(  \mathbf{m}_{1},\mathbf{m}_{2},\ldots,\mathbf{m}_{n}\right)  \in
Z\left(  \mathbf{k}\left[  S_{n}\right]  \right)
\ \ \ \ \ \ \ \ \ \ \text{for each }k\in\mathbb{N}%
\]
(since $Z\left(  \mathbf{k}\left[  S_{n}\right]  \right)  $ is a $\mathbf{k}%
$-subalgebra of $\mathbf{k}\left[  S_{n}\right]  $). But this follows from
\[
e_{k}\left(  \mathbf{m}_{1},\mathbf{m}_{2},\ldots,\mathbf{m}_{n}\right)
=\sum_{\substack{w\in S_{n};\\w\text{ has exactly }n-k\text{ orbits}}}w
\]
(which is clearly preserved under conjugation).
\end{proof}

As a consequence, we have the following chain of subalgebras:%
\begin{equation}
Z\left(  \mathbf{k}\left[  S_{n}\right]  \right)  \subseteq\operatorname*{GZ}%
\nolimits_{n}\subseteq\mathbf{k}\left[  S_{n}\right]  ,
\label{eq.thm.center.murphies.incs}%
\end{equation}
where $\operatorname*{GZ}\nolimits_{n}$ is the Gelfand--Tsetlin subalgebra
defined in Corollary \ref{cor.GZ.comm}.

\subsection{A potpourri of other elements}

In this section, I will showcase several other families of elements of
$\mathbf{k}\left[  S_{n}\right]  $, stating a few properties of each. Some of
these properties will be proved in later chapters.

\subsubsection{Descent numbers and the Eulerian subalgebra}

Descents are one of the simplest features of permutations:

\begin{definition}
\label{def.Des.des}Let $w\in S_{n}$ be a permutation. \medskip

\textbf{(a)} A \emph{descent} of $w$ means a number $i\in\left[  n-1\right]  $
such that $w\left(  i\right)  >w\left(  i+1\right)  $. \medskip

\textbf{(b)} The \emph{descent set} of $w$ means the set of all descents of
$w$. It is called $\operatorname*{Des}w$. \medskip

\textbf{(c)} The \emph{descent number} of $w$ means the \# of all descents of
$w$. It is called $\operatorname*{des}w$.
\end{definition}

\begin{example}
If $w=\operatorname*{oln}\left(  516324\right)  \in S_{6}$, then the descents
of $w$ are $1,3,4$, so that $\operatorname*{Des}w=\left\{  1,3,4\right\}  $
and $\operatorname*{des}w=3$.
\end{example}

See \cite[Lecture 28, \S 4.6]{22fco} and \cite[\S A.4.11]{21s} for several
properties of descents, and particularly for the \emph{Eulerian numbers},
which count the permutations $w\in S_{n}$ with a given number of descents. We
will instead take the sum of these permutations in $\mathbf{k}\left[
S_{n}\right]  $:

\begin{definition}
\label{def.Des.dk}For any $k\in\mathbb{N}$, let us set%
\[
\mathbf{d}_{k}:=\sum_{\substack{w\in S_{n};\\\operatorname*{des}w=k}%
}w\in\mathbf{k}\left[  S_{n}\right]  .
\]

\end{definition}

\begin{example}
For $n=3$, we have%
\begin{align*}
\mathbf{d}_{0}  &  =\operatorname*{id}=1;\\
\mathbf{d}_{1}  &  =s_{1}+s_{2}+\operatorname*{cyc}\nolimits_{1,2,3}%
+\operatorname*{cyc}\nolimits_{1,3,2};\\
\mathbf{d}_{2}  &  =t_{1,3};\\
\mathbf{d}_{k}  &  =0\ \ \ \ \ \ \ \ \ \ \text{for all }k>2.
\end{align*}

\end{example}

More generally:

\begin{proposition}
\textbf{(a)} We have $\mathbf{d}_{k}=0$ for all $k\geq n>0$, since a
permutation $w\in S_{n}$ has no more than $n-1$ descents. \medskip

\textbf{(b)} We have $\mathbf{d}_{0}=\operatorname*{id}=1$. \medskip

\textbf{(c)} Assume that $n>0$. Then, $\mathbf{d}_{n-1}=w_{0}$, where
$w_{0}\in S_{n}$ is the permutation with OLN $\left(  n,n-1,\ldots,2,1\right)
$.
\end{proposition}

Much less obviously, the following holds:

\begin{theorem}
\label{thm.Des.dk-subalg}Assume that $n>0$. Then, the $\mathbf{k}$-linear span
of $\left\{  \mathbf{d}_{0},\mathbf{d}_{1},\ldots,\mathbf{d}_{n-1}\right\}  $
is a commutative $\mathbf{k}$-subalgebra of $\mathbf{k}\left[  S_{n}\right]
$. In other words, there exist coefficients $c_{i,j,k}\in\mathbf{k}$
(actually, $\in\mathbb{N}$) such that all $i,j\in\left\{  0,1,\ldots
,n-1\right\}  $ satisfy%
\[
\mathbf{d}_{i}\mathbf{d}_{j}=\mathbf{d}_{j}\mathbf{d}_{i}=c_{i,j,0}%
\mathbf{d}_{0}+c_{i,j,1}\mathbf{d}_{1}+\cdots+c_{i,j,n-1}\mathbf{d}_{n-1}.
\]

\end{theorem}

This subalgebra is called the \emph{Eulerian subalgebra} of $\mathbf{k}\left[
S_{n}\right]  $. The name harkens back to the \emph{Eulerian numbers}, which
count the permutations $w\in S_{n}$ having $k$ descents (\cite[Lecture 28,
Definition 4.6.3]{22fco}).

\begin{question}
What are the coefficients $c_{i,j,k}$ in Theorem \ref{thm.Des.dk-subalg}? Can
they be computed explicitly?
\end{question}

Here is another curious property of descent-related sums (\cite[part of
Corollary 2]{VerTsi16}, \cite[(1)]{Rentel23}):

\begin{exercise}
\fbox{10} Assume that $n\geq1$. Let
\[
\widetilde{\mathbf{d}}:=\sum_{w\in S_{n}}\operatorname*{des}w\cdot
w=\sum_{k=0}^{n}k\mathbf{d}_{k}.
\]
Prove that
\[
\widetilde{\mathbf{d}}^{2}=-\left(  n-1\right)  !\cdot\widetilde{\mathbf{d}%
}+\dfrac{1}{4}\left(  n-1\right)  !\left(  n-1\right)  \left(  n^{2}%
-n+2\right)  \nabla.
\]

\end{exercise}

\subsubsection{V-permutations}

Let us next define a special type of permutations:

\begin{definition}
\label{def.Vperm.Vperm}A permutation $w\in S_{n}$ is said to be a
\emph{V-permutation} if there exists some $k\in\left[  n\right]  $ such that%
\[
w\left(  1\right)  >w\left(  2\right)  >\cdots>w\left(  k\right)  <w\left(
k+1\right)  <w\left(  k+2\right)  <\cdots<w\left(  n\right)
\]
(that is, $w\left(  i\right)  >w\left(  i+1\right)  $ for each $i<k$, and
$w\left(  i\right)  <w\left(  i+1\right)  $ for each $i\geq k$).
\end{definition}

Note that the $k$ here is uniquely determined by $w$, and in fact satisfies
$k=w^{-1}\left(  1\right)  $ (why?).

\begin{example}
The permutation $\operatorname*{oln}\left(  643125\right)  $ is a V-permutation.
\end{example}

V-permutations have got their name from the fact that their \textquotedblleft
plot\textquotedblright\ (obtained by marking the points $\left(  i,w\left(
i\right)  \right)  $ for all $i\in\left[  n\right]  $ in the Cartesian plane,
and connecting adjacent points by straight lines) looks \textquotedblleft
V-shaped\textquotedblright\ (i.e., first decreasing, then increasing). For
example, the \textquotedblleft plot\textquotedblright\ of the permutation
$\operatorname*{oln}\left(  643125\right)  $ looks as follows:%
\[
\begin{tikzpicture}
\clip(-0.3, -0.3) rectangle (7.5, 7.5);
\draw[very thin,color=gray] (0, 0) grid (7, 7);
\draw[->] (0,0) -- (7,0) node[right] {$x$};
\draw[->] (0,0) -- (0,7) node[above] {$y$};
\draw[very thick, color=blue] (1, 6) -- (2, 4) -- (3, 3) -- (4, 1) -- (5, 2) -- (6, 5);
\filldraw[blue] (1, 6) circle (0.1);
\filldraw[blue] (2, 4) circle (0.1);
\filldraw[blue] (3, 3) circle (0.1);
\filldraw[blue] (4, 1) circle (0.1);
\filldraw[blue] (5, 2) circle (0.1);
\filldraw[blue] (6, 5) circle (0.1);
\end{tikzpicture}
\ \ .
\]

Some properties of V-permutations are discussed in \cite[\S A.4.4]{21s}. We
shall now encounter them when multiplying a certain product in $\mathbf{k}%
\left[  S_{n}\right]  $:

\begin{theorem}
\label{thm.Vperm.Vn}Assume that $n\geq1$. Let
\begin{align*}
\mathbf{V}_{n}:=  &  \left(  1-\operatorname*{cyc}\nolimits_{2,1}\right)
\left(  1-\operatorname*{cyc}\nolimits_{3,2,1}\right)  \left(
1-\operatorname*{cyc}\nolimits_{4,3,2,1}\right)  \cdots\left(
1-\operatorname*{cyc}\nolimits_{n,n-1,\ldots,1}\right) \\
=  &  \prod_{i=2}^{n}\left(  1-\operatorname*{cyc}\nolimits_{i,i-1,\ldots
,1}\right)  \in\mathbf{k}\left[  S_{n}\right]  ,
\end{align*}
where we agree to read a product $\prod_{i=2}^{n}a_{i}$ of noncommuting
elements $a_{i}$ as $a_{2}a_{3}\cdots a_{n}$. Then: \medskip

\textbf{(a)} We have%
\[
\mathbf{V}_{n}=\sum_{\substack{w\in S_{n}\text{ is a}\\\text{V-permutation}%
}}\left(  -1\right)  ^{w^{-1}\left(  1\right)  -1}w.
\]

\textbf{(b)} Moreover,%
\[
\mathbf{V}_{n}^{2}=n\mathbf{V}_{n}.
\]
Hence, if $n$ is invertible in $\mathbf{k}$, then $\dfrac{1}{n}\mathbf{V}_{n}$
is an idempotent, known as the \emph{Dynkin idempotent}.
\end{theorem}

Theorem \ref{thm.Vperm.Vn} \textbf{(b)} is a result of Wever (\cite[page
571]{Wever47}, see also \cite[Theorem 8.16]{Reuten93} and \cite[Remark
2.1]{Garsia90}).

\begin{remark}
\label{rmk.Vperm.Vn-comb}As an illustration for the combinatorial nature of
these results, let me show how Theorem \ref{thm.Vperm.Vn} \textbf{(b)} can be
restated as a purely combinatorial equality between two integers.

Indeed, set $\mathbf{k}=\mathbb{Z}$, and assume that $n\geq1$. Then, Theorem
\ref{thm.Vperm.Vn} \textbf{(a)} yields%
\begin{equation}
\mathbf{V}_{n}=\sum_{\substack{w\in S_{n}\text{ is a}\\\text{V-permutation}%
}}\left(  -1\right)  ^{w^{-1}\left(  1\right)  -1}w.
\label{eq.rmk.Vperm.Vn-comb.a}%
\end{equation}
Multiplying this equality by itself, we obtain%
\begin{align*}
\mathbf{V}_{n}^{2}  &  =\left(  \sum_{\substack{w\in S_{n}\text{ is
a}\\\text{V-permutation}}}\left(  -1\right)  ^{w^{-1}\left(  1\right)
-1}w\right)  \left(  \sum_{\substack{w\in S_{n}\text{ is a}%
\\\text{V-permutation}}}\left(  -1\right)  ^{w^{-1}\left(  1\right)
-1}w\right) \\
&  =\left(  \sum_{\substack{u\in S_{n}\text{ is a}\\\text{V-permutation}%
}}\left(  -1\right)  ^{u^{-1}\left(  1\right)  -1}u\right)  \left(
\sum_{\substack{v\in S_{n}\text{ is a}\\\text{V-permutation}}}\left(
-1\right)  ^{v^{-1}\left(  1\right)  -1}v\right) \\
&  \ \ \ \ \ \ \ \ \ \ \ \ \ \ \ \ \ \ \ \ \left(
\begin{array}
[c]{c}%
\text{here, we have renamed the summation indices }w\text{ and }w\\
\text{in the two sums as }u\text{ and }v
\end{array}
\right) \\
&  =\underbrace{\sum_{\substack{u\in S_{n}\text{ is a}\\\text{V-permutation}%
}}\ \ \sum_{\substack{v\in S_{n}\text{ is a}\\\text{V-permutation}}}}%
_{=\sum_{\substack{u,v\in S_{n}\text{ are two}\\\text{V-permutations}}%
}}\underbrace{\left(  -1\right)  ^{u^{-1}\left(  1\right)  -1}u\left(
-1\right)  ^{v^{-1}\left(  1\right)  -1}v}_{\substack{=\left(  -1\right)
^{\left(  u^{-1}\left(  1\right)  -1\right)  +\left(  v^{-1}\left(  1\right)
-1\right)  }uv\\=\left(  -1\right)  ^{u^{-1}\left(  1\right)  +v^{-1}\left(
1\right)  }uv\\\text{(since }\left(  u^{-1}\left(  1\right)  -1\right)
+\left(  v^{-1}\left(  1\right)  -1\right)  \equiv u^{-1}\left(  1\right)
+v^{-1}\left(  1\right)  \operatorname{mod}2\text{)}}}\\
&  =\sum_{\substack{u,v\in S_{n}\text{ are two}\\\text{V-permutations}%
}}\left(  -1\right)  ^{u^{-1}\left(  1\right)  +v^{-1}\left(  1\right)  }uv\\
&  =\sum_{w\in S_{n}}\ \ \sum_{\substack{u,v\in S_{n}\text{ are two}%
\\\text{V-permutations;}\\uv=w}}\left(  -1\right)  ^{u^{-1}\left(  1\right)
+v^{-1}\left(  1\right)  }w.
\end{align*}
Meanwhile, (\ref{eq.rmk.Vperm.Vn-comb.a}) also yields%
\begin{align*}
n\mathbf{V}_{n}  &  =n\sum_{\substack{w\in S_{n}\text{ is a}%
\\\text{V-permutation}}}\left(  -1\right)  ^{w^{-1}\left(  1\right)  -1}%
w=\sum_{\substack{w\in S_{n}\text{ is a}\\\text{V-permutation}}}n\left(
-1\right)  ^{w^{-1}\left(  1\right)  -1}w\\
&  =\sum_{w\in S_{n}}%
\begin{cases}
n\left(  -1\right)  ^{w^{-1}\left(  1\right)  -1}, & \text{if }w\text{ is a
V-permutation};\\
0, & \text{if not}%
\end{cases}
\cdot w.
\end{align*}

Now, Theorem \ref{thm.Vperm.Vn} \textbf{(b)} says that the left hand sides of
these two equalities are equal. Equivalently, it thus says that their right
hand sides are equal. In other words, it says that the equality%
\begin{align*}
&  \sum_{w\in S_{n}}\ \ \sum_{\substack{u,v\in S_{n}\text{ are two}%
\\\text{V-permutations;}\\uv=w}}\left(  -1\right)  ^{u^{-1}\left(  1\right)
+v^{-1}\left(  1\right)  }w\\
&  =\sum_{w\in S_{n}}%
\begin{cases}
n\left(  -1\right)  ^{w^{-1}\left(  1\right)  -1}, & \text{if }w\text{ is a
V-permutation};\\
0, & \text{if not}%
\end{cases}
\cdot w
\end{align*}
holds. By comparing coefficients (see Remark \ref{rmk.monalg.compare}), this
equality is in turn equivalent to the statement that%
\[
\sum_{\substack{u,v\in S_{n}\text{ are two}\\\text{V-permutations;}%
\\uv=w}}\left(  -1\right)  ^{u^{-1}\left(  1\right)  +v^{-1}\left(  1\right)
}=%
\begin{cases}
n\left(  -1\right)  ^{w^{-1}\left(  1\right)  -1}, & \text{if }w\text{ is a
V-permutation};\\
0, & \text{if not}%
\end{cases}
\]
for each $w\in S_{n}$. So this latter statement (an equality between two
integers, explicitly given as sums) is a combinatorial restatement of Theorem
\ref{thm.Vperm.Vn} \textbf{(b)}. Does this restatement make the theorem easier
to prove? I don't know!
\end{remark}

\begin{exercise}
\fbox{2} Prove Theorem \ref{thm.Vperm.Vn} \textbf{(a)}.
\end{exercise}

\begin{exercise}
\fbox{1} Let $w_{0}\in S_{n}$ be the permutation $\operatorname*{oln}\left(
n,n-1,\ldots,2,1\right)  $, which sends each $i\in\left[  n\right]  $ to
$n+1-i$. Prove that $\mathbf{V}_{n}=\left(  -1\right)  ^{n-1}\mathbf{V}%
_{n}w_{0}$, where $\mathbf{V}_{n}$ is as in Theorem \ref{thm.Vperm.Vn}.
\end{exercise}

\begin{exercise}
\fbox{2} For each $k\in\left[  n\right]  $, let $\alpha_{k}\in S_{n}$ be the
permutation that sends the numbers $1,2,\ldots,k$ to $k,k-1,\ldots,1$
(respectively) while leaving all remaining numbers $k+1,k+2,\ldots,n$
unchanged. Let $\mathbf{V}_{n}$ be as in Theorem \ref{thm.Vperm.Vn}. Prove
that%
\begin{align*}
\mathbf{V}_{n}  &  =\left(  1-\alpha_{2}\right)  \left(  1+\alpha_{3}\right)
\left(  1-\alpha_{4}\right)  \left(  1+\alpha_{5}\right)  \cdots\left(
1-\left(  -1\right)  ^{n}\alpha_{n}\right) \\
&  =\prod_{i=2}^{n}\left(  1-\left(  -1\right)  ^{i}\alpha_{i}\right)
\end{align*}
(where the product is understood as in Theorem \ref{thm.Vperm.Vn}).
\end{exercise}

\begin{exercise}
\fbox{15} Prove Theorem \ref{thm.Vperm.Vn} \textbf{(b)}.

[\textbf{Remark:} A direct or combinatorial proof would be particularly welcome.]
\end{exercise}

\subsubsection{Random-to-random shuffles}

We shall next define another kind of peculiar elements in $\mathbf{k}\left[
S_{n}\right]  $. They rely on the combinatorial notion of \emph{noninversions}:

\begin{definition}
\label{def.R2R.noninv}Let $i\in\mathbb{Z}$. Let $w\in S_{n}$ be a permutation.
Then, $\operatorname*{noninv}\nolimits_{i}\left(  w\right)  $ shall denote the
\# of $i$-element subsets $I$ of $\left[  n\right]  $ on which $w$ is
increasing (i.e., for which the restriction $\left.  w\mid_{I}\right.
:I\rightarrow\left[  n\right]  $ is increasing). In other words, we let
$\operatorname*{noninv}\nolimits_{i}\left(  w\right)  $ denote the \# of ways
to pick $i$ numbers $j_{1}<j_{2}<\cdots<j_{i}$ in $\left[  n\right]  $ such
that $w\left(  j_{1}\right)  <w\left(  j_{2}\right)  <\cdots<w\left(
j_{i}\right)  $. This number $\operatorname*{noninv}\nolimits_{i}\left(
w\right)  $ is called the $i$\emph{-noninversion number} of $w$.
\end{definition}

\begin{example}
Let $n=6$ and $w=\operatorname*{oln}\left(  642135\right)  \in S_{6}$. Then,
$\operatorname*{noninv}\nolimits_{0}\left(  w\right)  =1$ (since the only
$0$-element subset $I$ of $\left[  n\right]  $ is $\varnothing$, and our
permutation $w$ is increasing on $\varnothing$ for trivial reasons) and
$\operatorname*{noninv}\nolimits_{1}\left(  w\right)  =6$ (since there are six
$1$-element subsets $I$ of $\left[  n\right]  $, and again $w$ is increasing
on each of them for trivial reasons) and $\operatorname*{noninv}%
\nolimits_{2}\left(  w\right)  =6$ (since the $2$-element subsets $I$ of
$\left[  n\right]  $ on which $w$ is increasing are $\left\{  2,6\right\}  $,
$\left\{  3,5\right\}  $, $\left\{  3,6\right\}  $, $\left\{  4,5\right\}  $,
$\left\{  4,6\right\}  $ and $\left\{  5,6\right\}  $) and
$\operatorname*{noninv}\nolimits_{3}\left(  w\right)  =2$ (since the
$3$-element subsets $I$ of $\left[  n\right]  $ on which $w$ is increasing are
$\left\{  3,5,6\right\}  $ and $\left\{  4,5,6\right\}  $) and
$\operatorname*{noninv}\nolimits_{k}\left(  w\right)  =0$ for all $k>3$. Of
course, $\operatorname*{noninv}\nolimits_{k}\left(  w\right)  =0$ for all
$k<0$ as well.
\end{example}

We can now define a new family of elements of $\mathbf{k}\left[  S_{n}\right]
$:

\begin{definition}
\label{def.R2R.Rk}For each $k\in\mathbb{N}$, we define%
\[
\mathbf{R}_{k}:=\sum_{w\in S_{n}}\operatorname*{noninv}\nolimits_{n-k}\left(
w\right)  \cdot w\in\mathbf{k}\left[  S_{n}\right]  .
\]
This element $\mathbf{R}_{k}$ is called the $k$\emph{-random-to-random
shuffle} (for reasons we will later explain).
\end{definition}

For example, for $n=3$, we have%
\begin{align*}
\mathbf{R}_{0}  &  =\sum_{w\in S_{3}}\operatorname*{noninv}\nolimits_{3}%
\left(  w\right)  \cdot w=\operatorname*{id}=1;\\
\mathbf{R}_{1}  &  =\sum_{w\in S_{3}}\operatorname*{noninv}\nolimits_{2}%
\left(  w\right)  \cdot w=3\operatorname*{id}+\,2s_{1}+2s_{2}%
+\operatorname*{cyc}\nolimits_{1,2,3}+\operatorname*{cyc}\nolimits_{1,3,2};\\
\mathbf{R}_{2}  &  =\sum_{w\in S_{3}}\underbrace{\operatorname*{noninv}%
\nolimits_{1}\left(  w\right)  }_{=3}\cdot\,w=3\sum_{w\in S_{3}}w=3\nabla;\\
\mathbf{R}_{3}  &  =\sum_{w\in S_{3}}\underbrace{\operatorname*{noninv}%
\nolimits_{0}\left(  w\right)  }_{=1}\cdot\,w=\sum_{w\in S_{3}}w=\nabla.
\end{align*}
Some of this obviously generalizes:

\begin{itemize}
\item We always have $\mathbf{R}_{0}=1$ (since the only permutation $w\in
S_{n}$ with $\operatorname*{noninv}\nolimits_{n}w\neq0$ is $\operatorname*{id}%
=1$).

\item We have $\mathbf{R}_{n}=\nabla$ (since $\operatorname*{noninv}%
\nolimits_{0}w=1$ for each $w\in S_{n}$).

\item If $n\geq1$, then $\mathbf{R}_{n-1}=n\nabla$ (since
$\operatorname*{noninv}\nolimits_{1}w=n$ for each $w\in S_{n}$).

\item If $k>n$, then $\mathbf{R}_{k}=0$.
\end{itemize}

But the other elements $\mathbf{R}_{k}$ with $0<k<n-1$ are interesting. The
following fact was discovered by Reiner, Saliola and Welker around 2011
(\cite[Theorem 1.1]{RSW}), and still retains its mystery:

\begin{theorem}
\label{thm.R2R.commute}The $n+1$ elements $\mathbf{R}_{0},\mathbf{R}%
_{1},\ldots,\mathbf{R}_{n}$ commute.
\end{theorem}

Three proofs of Theorem \ref{thm.R2R.commute} are known by now: \cite[\S 6.1]%
{RSW} gives a terse and handwavy combinatorial argument; \cite[Chapitre
IV]{Lafren19} a long and computational proof; \cite[Theorem 4.1]{BCGS25} a
simpler but still fairly long proof (which applies to a more general setting,
in which the symmetric group algebra $\mathbf{k}\left[  S_{n}\right]  $ is
replaced by the Hecke algebra $\mathcal{H}_{n}$).

\begin{question}
Is there a proof of Theorem \ref{thm.R2R.commute} that is as simple as the
theorem would suggest?
\end{question}

Another recent theorem (by Dieker, Saliola and Lafreni\`{e}re) says:

\begin{theorem}
\label{thm.R2R.mipo-generic}Let $k\in\left\{  0,1,\ldots,n\right\}  $. Then,
there exists a polynomial $P\in\mathbb{Z}\left[  x\right]  $ that splits over
$\mathbb{Z}$ (that is, that can be written in the form $P=\prod_{r\in
S}\left(  x-r\right)  $ for some finite multiset $S$ of integers) such that
$P\left(  \mathbf{R}_{k}\right)  =0$.

For $k=1$, we can be more explicit:%
\[
\prod_{i=0}^{n^{2}}\left(  \mathbf{R}_{1}-i\right)  =0.
\]

\end{theorem}

The original proof (in \cite{Lafren19}, building upon \cite{DieSal18}) takes
about 100 pages. Recently, a shorter proof was found by Axelrod-Freed,
Brauner, Chiang, Commins, Lang and the present author (see \cite{AFBCCL24} for
the $k=1$ case, and \cite[Proposition D.1]{BCGS25} for all other $k$). Still,
the proof is technical and uses a significant amount of representation theory.

\subsubsection{Excedance and antiexcedance elements}

Here is something even stranger: an open problem (albeit one that has not been
advertised very widely). We again begin with a combinatorial definition:

\begin{definition}
\label{def.exc-anxc}Let $w\in S_{n}$ be a permutation. Then, we define the two
numbers%
\begin{align*}
\operatorname*{exc}w  &  :=\left(  \text{\# of }i\in\left[  n\right]  \text{
such that }w\left(  i\right)  >i\right)  \ \ \ \ \ \ \ \ \ \ \text{and}\\
\operatorname*{anxc}w  &  :=\left(  \text{\# of }i\in\left[  n\right]  \text{
such that }w\left(  i\right)  <i\right)  .
\end{align*}
These numbers $\operatorname*{exc}w$ and $\operatorname*{anxc}w$ are called
the \emph{excedance number} and the \emph{antiexcedance number} of $w$, respectively.
\end{definition}

\begin{example}
Let $w=\operatorname*{oln}\left(  3524176\right)  $. Then,
$\operatorname*{exc}w=3$ and $\operatorname*{anxc}w=3$.
\end{example}

\begin{definition}
\label{def.Xab.Xab}For any $a,b\in\mathbb{N}$, define%
\[
\mathbf{X}_{a,b}:=\sum_{\substack{w\in S_{n};\\\operatorname*{exc}%
w=a;\\\operatorname*{anxc}w=b}}w\in\mathbf{k}\left[  S_{n}\right]  .
\]

\end{definition}

\begin{example}
For $n=3$, we have%
\begin{align*}
\mathbf{X}_{0,0}  &  =\operatorname*{id}=1;\ \ \ \ \ \ \ \ \ \ \mathbf{X}%
_{1,2}=\operatorname*{cyc}\nolimits_{1,3,2};\ \ \ \ \ \ \ \ \ \ \mathbf{X}%
_{2,1}=\operatorname*{cyc}\nolimits_{1,2,3};\\
\mathbf{X}_{1,1}  &  =t_{1,2}+t_{1,3}+t_{2,3};
\end{align*}
and $\mathbf{X}_{a,b}=0$ for all remaining choices of $a,b\in\mathbb{N}$.
\end{example}

\begin{conjecture}
\label{conj.exc-anxc.comm}These elements $\mathbf{X}_{a,b}$ for all
$a,b\in\mathbb{N}$ commute (for fixed $n$). In other words, $\mathbf{X}%
_{a,b}\mathbf{X}_{c,d}=\mathbf{X}_{c,d}\mathbf{X}_{a,b}$ for all
$a,b,c,d\in\mathbb{N}$.
\end{conjecture}

This conjecture was inspired by a conversation with Theo Douvropoulos; I have
checked it for all $n\leq8$.

We can restate this conjecture explicitly (by comparing coefficients, just as
in Remark \ref{rmk.Vperm.Vn-comb}): It is saying that for any permutation
$w\in S_{n}$ and any numbers $a,b,c,d\in\mathbb{N}$, we have%
\begin{align*}
&  \left(  \text{\# of pairs }\left(  u,v\right)  \in S_{n}\times S_{n}\text{
with }w=uv\right. \\
&  \ \ \ \ \ \ \ \ \ \ \left.  \text{and }\operatorname*{exc}u=a\text{ and
}\operatorname*{anxc}u=b\text{ and }\operatorname*{exc}v=c\text{ and
}\operatorname*{anxc}v=d\right) \\
&  =\left(  \text{\# of pairs }\left(  u,v\right)  \in S_{n}\times S_{n}\text{
with }w=vu\right. \\
&  \ \ \ \ \ \ \ \ \ \ \left.  \text{and }\operatorname*{exc}u=a\text{ and
}\operatorname*{anxc}u=b\text{ and }\operatorname*{exc}v=c\text{ and
}\operatorname*{anxc}v=d\right)  .
\end{align*}
(Indeed, this equality is what we obtain if we compare the coefficients of $w$
on both sides of the equality $\mathbf{X}_{a,b}\mathbf{X}_{c,d}=\mathbf{X}%
_{c,d}\mathbf{X}_{a,b}$.)

For a change, the elements $\mathbf{X}_{a,b}$ are \textbf{not} (in general)
annihilated by polynomials that split over $\mathbb{Z}$.

\subsection{The antipode and the sign-twist}

Now we shall define some maps on $\mathbf{k}\left[  S_{n}\right]  $. Such maps
can be used to transform the above-defined elements into others, thus doubling
our repository of interesting elements. Moreover, if these maps are
structure-preserving (e.g., $\mathbf{k}$-linear maps or $\mathbf{k}$-algebra
homomorphisms), then they can be used to carry some of the properties of the
elements along, hence saving us a lot of proof work.

\subsubsection{Morphisms}

We begin by recalling the concept of a $\mathbf{k}$-algebra automorphism, and
related concepts:

\begin{definition}
\label{def.alghom.alghom}\textbf{(a)} A $\mathbf{k}$\emph{-algebra
homomorphism} (for short: $\mathbf{k}$\emph{-algebra morphism}) is a
$\mathbf{k}$-linear map $f:A\rightarrow B$ between two $\mathbf{k}$-algebras
$A$ and $B$ that respects the unity (that is, satisfies $f\left(
1_{A}\right)  =1_{B}$) and respects multiplication (that is, satisfies
$f\left(  aa^{\prime}\right)  =f\left(  a\right)  f\left(  a^{\prime}\right)
$ for all $a,a^{\prime}\in A$). Equivalently, it is a map that is
simultaneously a $\mathbf{k}$-module morphism and a ring morphism. \medskip

\textbf{(b)} A $\mathbf{k}$\emph{-algebra endomorphism} is a $\mathbf{k}%
$-algebra morphism from a $\mathbf{k}$-algebra to itself. \medskip

\textbf{(c)} A $\mathbf{k}$\emph{-algebra isomorphism} is a $\mathbf{k}%
$-algebra morphism that has an inverse, which is also a $\mathbf{k}$-algebra
morphism. (Actually, the \textquotedblleft which is also a $\mathbf{k}%
$-algebra morphism\textquotedblright\ requirement comes for free; see
\cite[Proposition 3.11.8]{23wa} for a proof.) \medskip

\textbf{(d)} A $\mathbf{k}$\emph{-algebra automorphism} is a $\mathbf{k}%
$-algebra isomorphism from a $\mathbf{k}$-algebra to itself.
\end{definition}

A ubiquitous source of $\mathbf{k}$-algebra automorphisms is conjugation: If
$u$ is an invertible element of some $\mathbf{k}$-algebra $A$, then the map%
\begin{align*}
A  &  \rightarrow A,\\
a  &  \mapsto uau^{-1}%
\end{align*}
is a $\mathbf{k}$-algebra automorphism of $A$, known as \emph{conjugation by
}$u$. (This is easy to check -- do it if you have not seen it before!)
Automorphisms of $A$ obtained in this way (i.e., obtained as conjugations by
invertible elements $u\in A$) are called \emph{inner automorphisms}.

Note that any element $g$ of a group $G$ is automatically invertible in its
group algebra $\mathbf{k}\left[  G\right]  $ (or rather, to be precise: the
corresponding standard basis vector $e_{g}$ is invertible in $\mathbf{k}%
\left[  G\right]  $), and thus gives rise to an inner automorphism of
$\mathbf{k}\left[  G\right]  $. There are typically many invertible elements
of $\mathbf{k}\left[  G\right]  $ besides the standard basis
vectors\footnote{The simplest example is probably $-1$, but you can easily
find more interesting ones. For example, if $n\geq3$, then the element
$1+t_{1,2}+\operatorname*{cyc}\nolimits_{1,2,3}-\operatorname*{cyc}%
\nolimits_{1,3,2}-t_{1,3}\in\mathbf{k}\left[  S_{n}\right]  $ is invertible,
with inverse $1-t_{1,2}-\operatorname*{cyc}\nolimits_{1,2,3}%
+\operatorname*{cyc}\nolimits_{1,3,2}+t_{1,3}$.}, and thus many more inner automorphisms.

But even these are often not the only automorphisms of $\mathbf{k}\left[
G\right]  $. For $\mathbf{k}\left[  S_{n}\right]  $ in particular, we have an
important automorphism that is (almost) never an inner automorphism:

\subsubsection{The sign-twist}

\begin{definition}
\label{def.Tsign.Tsign}The \emph{sign-twist} is the $\mathbf{k}$-linear map%
\begin{align*}
T_{\operatorname*{sign}}:\mathbf{k}\left[  S_{n}\right]   &  \rightarrow
\mathbf{k}\left[  S_{n}\right]  ,\\
w  &  \mapsto\left(  -1\right)  ^{w}w\ \ \ \ \ \ \ \ \ \ \text{for all }w\in
S_{n}.
\end{align*}
This definition should be read as follows: \textquotedblleft The
\emph{sign-twist} is the $\mathbf{k}$-linear map $T_{\operatorname*{sign}%
}:\mathbf{k}\left[  S_{n}\right]  \rightarrow\mathbf{k}\left[  S_{n}\right]  $
that sends each standard basis vector $w$ (with $w\in S_{n}$) to $\left(
-1\right)  ^{w}w$\textquotedblright. This requirement uniquely determines this
map, since a $\mathbf{k}$-linear map from $\mathbf{k}\left[  S_{n}\right]  $
is uniquely determined by its values on the standard basis vectors.
Explicitly, it thus is given by%
\[
T_{\operatorname*{sign}}\left(  \sum_{w\in S_{n}}\alpha_{w}w\right)
=\sum_{w\in S_{n}}\alpha_{w}\left(  -1\right)  ^{w}%
w\ \ \ \ \ \ \ \ \ \ \text{for any scalars }\alpha_{w}\in\mathbf{k}.
\]

\end{definition}

\begin{example}
\label{exa.Tsign.exas1}\textbf{(a)} For $n=3$, we have%
\begin{align*}
&  T_{\operatorname*{sign}}\left(  2\operatorname*{id}+\,7s_{1}%
-9\operatorname*{cyc}\nolimits_{1,2,3}\right) \\
&  =2\cdot\underbrace{\left(  -1\right)  ^{\operatorname*{id}}}_{=1}%
\operatorname*{id}+\,7\cdot\underbrace{\left(  -1\right)  ^{s_{1}}}_{=-1}%
s_{1}-9\cdot\underbrace{\left(  -1\right)  ^{\operatorname*{cyc}%
\nolimits_{1,2,3}}}_{=1}\operatorname*{cyc}\nolimits_{1,2,3}\\
&  =2\operatorname*{id}-\,7s_{1}-9\operatorname*{cyc}\nolimits_{1,2,3}.
\end{align*}

\textbf{(b)} For each transposition $t_{i,j}\in S_{n}$, we have%
\begin{equation}
T_{\operatorname*{sign}}\left(  t_{i,j}\right)  =\underbrace{\left(
-1\right)  ^{t_{i,j}}}_{=-1}t_{i,j}=-t_{i,j}. \label{eq.exa.Tsign.exas1.b.1}%
\end{equation}
Thus, for any $k\in\left[  n\right]  $, the jucys--murphy $\mathbf{m}_{k}$
satisfies%
\[
T_{\operatorname*{sign}}\left(  \mathbf{m}_{k}\right)  =-\mathbf{m}_{k}%
\]
(since $\mathbf{m}_{k}$ is a sum of transpositions $t_{i,j}$). \medskip

\textbf{(c)} We have%
\begin{equation}
T_{\operatorname*{sign}}\left(  \nabla\right)  =\nabla^{-}.
\label{eq.exa.Tsign.exas1.c.1}%
\end{equation}
More generally, for any subset $X$ of $\left[  n\right]  $, we have%
\begin{equation}
T_{\operatorname*{sign}}\left(  \nabla_{X}\right)  =\nabla_{X}^{-}.
\label{eq.exa.Tsign.exas1.c.2}%
\end{equation}
This is easy to see directly from the definitions of $\nabla_{X}$ and
$\nabla_{X}^{-}$. For example, to prove (\ref{eq.exa.Tsign.exas1.c.2}), we can
argue that every subset $X$ of $\left[  n\right]  $ satisfies%
\begin{align*}
T_{\operatorname*{sign}}\left(  \nabla_{X}\right)   &
=T_{\operatorname*{sign}}\left(  \sum_{\substack{w\in S_{n};\\w\left(
i\right)  =i\text{ for all }i\notin X}}w\right)  \ \ \ \ \ \ \ \ \ \ \left(
\text{since }\nabla_{X}=\sum_{\substack{w\in S_{n};\\w\left(  i\right)
=i\text{ for all }i\notin X}}w\right) \\
&  =\sum_{\substack{w\in S_{n};\\w\left(  i\right)  =i\text{ for all }i\notin
X}}\underbrace{T_{\operatorname*{sign}}\left(  w\right)  }_{\substack{=\left(
-1\right)  ^{w}w\\\text{(by the definition of }T_{\operatorname*{sign}%
}\text{)}}}\ \ \ \ \ \ \ \ \ \ \left(
\begin{array}
[c]{c}%
\text{since }T_{\operatorname*{sign}}\text{ is}\\
\text{a }\mathbf{k}\text{-linear map}%
\end{array}
\right) \\
&  =\sum_{\substack{w\in S_{n};\\w\left(  i\right)  =i\text{ for all }i\notin
X}}\left(  -1\right)  ^{w}w=\nabla_{X}^{-}\ \ \ \ \ \ \ \ \ \ \left(  \text{by
the definition of }\nabla_{X}^{-}\right)  .
\end{align*}
Thus, (\ref{eq.exa.Tsign.exas1.c.2}) is proved. Of course,
(\ref{eq.exa.Tsign.exas1.c.1}) follows by applying
(\ref{eq.exa.Tsign.exas1.c.2}) to $X=\left[  n\right]  $. \medskip

\textbf{(d)} For the somewhere-to-below shuffles $\mathbf{t}_{k}$ (defined in
Definition \ref{def.stb.tk}), we have%
\begin{align*}
T_{\operatorname*{sign}}\left(  \mathbf{t}_{k}\right)   &
=T_{\operatorname*{sign}}\left(  \operatorname*{cyc}\nolimits_{k}%
+\operatorname*{cyc}\nolimits_{k,k+1}+\operatorname*{cyc}\nolimits_{k,k+1,k+2}%
+\cdots+\operatorname*{cyc}\nolimits_{k,k+1,\ldots,n}\right) \\
&  =\operatorname*{cyc}\nolimits_{k}-\operatorname*{cyc}\nolimits_{k,k+1}%
+\operatorname*{cyc}\nolimits_{k,k+1,k+2}\pm\cdots+\left(  -1\right)
^{n-k}\operatorname*{cyc}\nolimits_{k,k+1,\ldots,n}%
\end{align*}
(an alternating sum, in which the signs alternate between $+$ and $-$). This
follows easily from the formula $\left(  -1\right)  ^{\operatorname*{cyc}%
\nolimits_{i_{1},i_{2},\ldots,i_{k}}}=\left(  -1\right)  ^{k-1}$ for the sign
of a cycle.
\end{example}

The sign-twist $T_{\operatorname*{sign}}$ has some nice but simple properties.
To state them, we need one more notion:

\begin{definition}
\label{def.invol}An \emph{involution} means a map $f:X\rightarrow X$ that
satisfies $f\circ f=\operatorname*{id}$.
\end{definition}

For example, all transpositions $t_{i,j}\in S_{X}$ for any set $X$ are
involutions. But so is the sign-twist $T_{\operatorname*{sign}}$, as we shall
now show.

\begin{theorem}
\label{thm.Tsign.auto}\textbf{(a)} The map $T_{\operatorname*{sign}%
}:\mathbf{k}\left[  S_{n}\right]  \rightarrow\mathbf{k}\left[  S_{n}\right]  $
is a $\mathbf{k}$-algebra automorphism. \medskip

\textbf{(b)} It is furthermore an involution (i.e., it satisfies
$T_{\operatorname*{sign}}\circ T_{\operatorname*{sign}}=\operatorname*{id}$).
\end{theorem}

\begin{proof}
[Proof of Theorem \ref{thm.Tsign.auto}.]This is a fairly straightforward use
of linearity. To illustrate this strategy, let me give a detailed proof. We
begin with part \textbf{(b)}. \medskip

\textbf{(b)} We must show that $T_{\operatorname*{sign}}\circ
T_{\operatorname*{sign}}=\operatorname*{id}$. In other words, we must show
that $T_{\operatorname*{sign}}\left(  T_{\operatorname*{sign}}\left(
\mathbf{a}\right)  \right)  =\mathbf{a}$ for any $\mathbf{a}\in\mathbf{k}%
\left[  S_{n}\right]  $.

So let $\mathbf{a}\in\mathbf{k}\left[  S_{n}\right]  $ be arbitrary. We must
prove the equality $T_{\operatorname*{sign}}\left(  T_{\operatorname*{sign}%
}\left(  \mathbf{a}\right)  \right)  =\mathbf{a}$. This equality is
$\mathbf{k}$-linear in $\mathbf{a}$ (since the map $T_{\operatorname*{sign}}$
is $\mathbf{k}$-linear). Thus, we can WLOG assume that $\mathbf{a}$ is a
standard basis vector, i.e., has the form $\mathbf{a}=u$ for some $u\in S_{n}%
$. Assume this, and consider this $u$. From $\mathbf{a}=u$, we obtain
$T_{\operatorname*{sign}}\left(  \mathbf{a}\right)  =T_{\operatorname*{sign}%
}\left(  u\right)  =\left(  -1\right)  ^{u}u$ (by the definition of
$T_{\operatorname*{sign}}$, since $u\in S_{n}$). Applying the map
$T_{\operatorname*{sign}}$ to both sides of this equality, we find%
\begin{align*}
T_{\operatorname*{sign}}\left(  T_{\operatorname*{sign}}\left(  \mathbf{a}%
\right)  \right)   &  =T_{\operatorname*{sign}}\left(  \left(  -1\right)
^{u}u\right) \\
&  =\left(  -1\right)  ^{u}\underbrace{T_{\operatorname*{sign}}\left(
u\right)  }_{=\left(  -1\right)  ^{u}u}\ \ \ \ \ \ \ \ \ \ \left(  \text{since
the map }T_{\operatorname*{sign}}\text{ is }\mathbf{k}\text{-linear}\right) \\
&  =\underbrace{\left(  -1\right)  ^{u}\left(  -1\right)  ^{u}}%
_{\substack{=\left(  \left(  -1\right)  ^{u}\right)  ^{2}=1\\\text{(since
}\left(  -1\right)  ^{u}\text{ is either }1\text{ or }-1\text{,}\\\text{and
thus squares to }1\text{)}}}u=u=\mathbf{a}.
\end{align*}

We thus have shown that $T_{\operatorname*{sign}}\left(
T_{\operatorname*{sign}}\left(  \mathbf{a}\right)  \right)  =\mathbf{a}$. As
explained above, this completes the proof of Theorem \ref{thm.Tsign.auto}
\textbf{(b)}. \medskip

\textbf{(a)} The map $T_{\operatorname*{sign}}$ is $\mathbf{k}$-linear (by
definition) and respects the unity (since it sends the identity permutation
$\operatorname*{id}\in S_{n}$ to $T_{\operatorname*{sign}}\left(
\operatorname*{id}\right)  =\underbrace{\left(  -1\right)
^{\operatorname*{id}}}_{=1}\operatorname*{id}=\operatorname*{id}$). We shall
now show that it respects multiplication. In other words, we shall show that
$T_{\operatorname*{sign}}\left(  \mathbf{ab}\right)  =T_{\operatorname*{sign}%
}\left(  \mathbf{a}\right)  \cdot T_{\operatorname*{sign}}\left(
\mathbf{b}\right)  $ for any $\mathbf{a},\mathbf{b}\in\mathbf{k}\left[
S_{n}\right]  $.

So let $\mathbf{a},\mathbf{b}\in\mathbf{k}\left[  S_{n}\right]  $ be
arbitrary. We must prove the equality $T_{\operatorname*{sign}}\left(
\mathbf{ab}\right)  =T_{\operatorname*{sign}}\left(  \mathbf{a}\right)  \cdot
T_{\operatorname*{sign}}\left(  \mathbf{b}\right)  $. This equality is
$\mathbf{k}$-linear in $\mathbf{a}$ (since the map $T_{\operatorname*{sign}}$
is $\mathbf{k}$-linear). Thus, we can WLOG assume that $\mathbf{a}$ is a
standard basis vector (i.e., has the form $\mathbf{a}=u$ for some $u\in S_{n}%
$). Similarly, we can WLOG assume that $\mathbf{b}$ is a standard basis
vector. Let us make both of these assumptions. Hence, $\mathbf{a}=u$ and
$\mathbf{b}=v$ for some permutations $u,v\in S_{n}$. Consider these $u,v$.
From $\mathbf{a}=u$, we obtain $T_{\operatorname*{sign}}\left(  \mathbf{a}%
\right)  =T_{\operatorname*{sign}}\left(  u\right)  =\left(  -1\right)  ^{u}u$
(by the definition of $T_{\operatorname*{sign}}$, since $u\in S_{n}$).
Similarly, $T_{\operatorname*{sign}}\left(  \mathbf{b}\right)  =\left(
-1\right)  ^{v}v$.

One of the basic properties of signs shows that $\left(  -1\right)
^{uv}=\left(  -1\right)  ^{u}\left(  -1\right)  ^{v}$.

From $\mathbf{a}=u$ and $\mathbf{b}=v$, we obtain%
\begin{align*}
T_{\operatorname*{sign}}\left(  \mathbf{ab}\right)   &
=T_{\operatorname*{sign}}\left(  uv\right) \\
&  =\underbrace{\left(  -1\right)  ^{uv}}_{=\left(  -1\right)  ^{u}\left(
-1\right)  ^{v}}uv\ \ \ \ \ \ \ \ \ \ \left(  \text{by the definition of
}T_{\operatorname*{sign}}\text{, since }uv\in S_{n}\right) \\
&  =\left(  -1\right)  ^{u}\left(  -1\right)  ^{v}uv=\underbrace{\left(
-1\right)  ^{u}u}_{=T_{\operatorname*{sign}}\left(  \mathbf{a}\right)  }%
\cdot\underbrace{\left(  -1\right)  ^{v}v}_{=T_{\operatorname*{sign}}\left(
\mathbf{b}\right)  }=T_{\operatorname*{sign}}\left(  \mathbf{a}\right)  \cdot
T_{\operatorname*{sign}}\left(  \mathbf{b}\right)  .
\end{align*}

Forget that we fixed $\mathbf{a},\mathbf{b}$. We thus have shown that
$T_{\operatorname*{sign}}\left(  \mathbf{ab}\right)  =T_{\operatorname*{sign}%
}\left(  \mathbf{a}\right)  \cdot T_{\operatorname*{sign}}\left(
\mathbf{b}\right)  $ for any $\mathbf{a},\mathbf{b}\in\mathbf{k}\left[
S_{n}\right]  $. In other words, the map $T_{\operatorname*{sign}}$ respects
multiplication. Since we already know that $T_{\operatorname*{sign}}$ is
$\mathbf{k}$-linear and respects the unity, we conclude that
$T_{\operatorname*{sign}}$ is a $\mathbf{k}$-algebra morphism.

From Theorem \ref{thm.Tsign.auto} \textbf{(b)}, we know that
$T_{\operatorname*{sign}}\circ T_{\operatorname*{sign}}=\operatorname*{id}$.
In other words, the map $T_{\operatorname*{sign}}$ is inverse to itself. Thus,
this map $T_{\operatorname*{sign}}$ is invertible, and its inverse is
$T_{\operatorname*{sign}}$ itself. Since $T_{\operatorname*{sign}}$ is a
$\mathbf{k}$-algebra morphism, we thus conclude that $T_{\operatorname*{sign}%
}$ is a $\mathbf{k}$-algebra isomorphism. Hence, $T_{\operatorname*{sign}}$ is
a $\mathbf{k}$-algebra automorphism (since it is a map from $\mathbf{k}\left[
S_{n}\right]  $ to itself). This proves Theorem \ref{thm.Tsign.auto}
\textbf{(a)}.
\end{proof}

\begin{exercise}
\fbox{1} Assume that $n\geq2$ and $\mathbb{Z}\subseteq\mathbf{k}$. Prove that
the sign-twist $T_{\operatorname*{sign}}$ is not an inner automorphism of
$\mathbf{k}\left[  S_{n}\right]  $ (that is, there exists no invertible
$\mathbf{u}\in\mathbf{k}\left[  S_{n}\right]  $ such that all $\mathbf{a}%
\in\mathbf{k}\left[  S_{n}\right]  $ satisfy $T_{\operatorname*{sign}}\left(
\mathbf{a}\right)  =\mathbf{uau}^{-1}$).
\end{exercise}

We can use the sign-twist $T_{\operatorname*{sign}}$ to turn identities in
$\mathbf{k}\left[  S_{n}\right]  $ into new identities. For instance, by
applying $T_{\operatorname*{sign}}$ to Lemma \ref{lem.intX.rec2}, we obtain
the following \textquotedblleft negative version\textquotedblright\ of Lemma
\ref{lem.intX.rec2}:

\begin{lemma}
\label{lem.intX.rec2sign}Let $X$ be any subset of $\left[  n\right]  $, and
let $x\in X$. Then,%
\[
\nabla_{X}^{-}=\nabla_{X\setminus\left\{  x\right\}  }^{-}\left(  1-\sum_{y\in
X\setminus\left\{  x\right\}  }t_{y,x}\right)  .
\]

\end{lemma}

\begin{proof}
Lemma \ref{lem.intX.rec2} yields
\[
\nabla_{X}=\nabla_{X\setminus\left\{  x\right\}  }\left(  1+\sum_{y\in
X\setminus\left\{  x\right\}  }t_{y,x}\right)  .
\]
Applying the map $T_{\operatorname*{sign}}$ to both sides of this equality, we
find%
\begin{align*}
T_{\operatorname*{sign}}\left(  \nabla_{X}\right)   &
=T_{\operatorname*{sign}}\left(  \nabla_{X\setminus\left\{  x\right\}
}\left(  1+\sum_{y\in X\setminus\left\{  x\right\}  }t_{y,x}\right)  \right)
\\
&  =\underbrace{T_{\operatorname*{sign}}\left(  \nabla_{X\setminus\left\{
x\right\}  }\right)  }_{\substack{=\nabla_{X\setminus\left\{  x\right\}  }%
^{-}\\\text{(by (\ref{eq.exa.Tsign.exas1.c.2}), applied to }X\setminus\left\{
x\right\}  \\\text{instead of }X\text{)}}}\left(  1+\sum_{y\in X\setminus
\left\{  x\right\}  }\underbrace{T_{\operatorname*{sign}}\left(
t_{y,x}\right)  }_{\substack{=-t_{y,x}\\\text{(by
(\ref{eq.exa.Tsign.exas1.b.1}))}}}\right) \\
&  \ \ \ \ \ \ \ \ \ \ \ \ \ \ \ \ \ \ \ \ \left(
\begin{array}
[c]{c}%
\text{since }T_{\operatorname*{sign}}\text{ is a }\mathbf{k}\text{-algebra
morphism}\\
\text{(by Theorem \ref{thm.Tsign.auto} \textbf{(a)})}%
\end{array}
\right) \\
&  =\nabla_{X\setminus\left\{  x\right\}  }^{-}\underbrace{\left(
1+\sum_{y\in X\setminus\left\{  x\right\}  }\left(  -t_{y,x}\right)  \right)
}_{=1-\sum_{y\in X\setminus\left\{  x\right\}  }t_{y,x}}=\nabla_{X\setminus
\left\{  x\right\}  }^{-}\left(  1-\sum_{y\in X\setminus\left\{  x\right\}
}t_{y,x}\right)  .
\end{align*}
Comparing this with (\ref{eq.exa.Tsign.exas1.c.2}), we find%
\[
\nabla_{X}^{-}=\nabla_{X\setminus\left\{  x\right\}  }^{-}\left(  1-\sum_{y\in
X\setminus\left\{  x\right\}  }t_{y,x}\right)  .
\]
This proves Lemma \ref{lem.intX.rec2sign}.
\end{proof}

\subsubsection{Anti-morphisms}

There is yet another important map on $\mathbf{k}\left[  S_{n}\right]  $,
called the \emph{antipode}, which is not quite an automorphism but fairly
close: It is an \emph{algebra anti-automorphism}. Let us define this general concept:

\begin{definition}
\textbf{(a)} A $\mathbf{k}$\emph{-algebra anti-homomorphism} (for short:
$\mathbf{k}$\emph{-algebra anti-morphism}) is a $\mathbf{k}$-linear map
$f:A\rightarrow B$ between two $\mathbf{k}$-algebras $A$ and $B$ that respects
the unity (that is, satisfies $f\left(  1_{A}\right)  =1_{B}$) and reverses
multiplication (that is, satisfies $f\left(  aa^{\prime}\right)  =f\left(
a^{\prime}\right)  f\left(  a\right)  $ for all $a,a^{\prime}\in A$). \medskip

\textbf{(b)} A $\mathbf{k}$\emph{-algebra anti-endomorphism} is a $\mathbf{k}%
$-algebra anti-morphism from a $\mathbf{k}$-algebra to itself. \medskip

\textbf{(c)} A $\mathbf{k}$\emph{-algebra anti-isomorphism} is a $\mathbf{k}%
$-algebra anti-morphism that has an inverse, which is also a $\mathbf{k}%
$-algebra anti-morphism. (Actually, the \textquotedblleft which is also a
$\mathbf{k}$-algebra anti-morphism\textquotedblright\ requirement comes for
free again.) \medskip

\textbf{(d)} A $\mathbf{k}$\emph{-algebra anti-automorphism} is a $\mathbf{k}%
$-algebra anti-isomorphism from a $\mathbf{k}$-algebra to itself.
\end{definition}

Anti-morphisms are not as frequently found as morphisms, but here is a very
popular example:

\begin{example}
Let $n\in\mathbb{N}$, and consider the $\mathbf{k}$-algebra $\mathbf{k}%
^{n\times n}$ of all $n\times n$-matrices with entries in $\mathbf{k}$. Then,
the map%
\begin{align*}
\mathbf{k}^{n\times n}  &  \rightarrow\mathbf{k}^{n\times n},\\
A  &  \mapsto A^{T}%
\end{align*}
(which transposes every matrix) is a $\mathbf{k}$-algebra anti-automorphism of
$\mathbf{k}^{n\times n}$ (and, by the way, an involution). In particular, it
reverses multiplication, since $\left(  AB\right)  ^{T}=B^{T}A^{T}$ for any
$n\times n$-matrices $A$ and $B$.
\end{example}

\begin{remark}
If any one (or both) of two $\mathbf{k}$-algebras $A$ and $B$ is commutative,
then there is no difference between \textquotedblleft morphisms from $A$ to
$B$\textquotedblright\ and \textquotedblleft anti-morphisms from $A$ to
$B$\textquotedblright. Thus, the notion of an anti-morphism is interesting
only in the noncommutative setting.
\end{remark}

\begin{remark}
Anti-morphisms can be viewed as (regular) morphisms to a slightly different
target. Indeed, if you are familiar with the \textquotedblleft opposite
ring\textquotedblright\ $B^{\operatorname*{op}}$ of a ring $B$ (see, e.g.,
\cite[Exercise 2.7.6]{23wa}), then you can probably guess what the opposite
$\mathbf{k}$-algebra of a $\mathbf{k}$-algebra is (essentially: it is the same
$\mathbf{k}$-algebra, but the order of factors in a product is reversed,
meaning that the product $b_{1}b_{2}$ has been renamed as $b_{2}b_{1}$). Now,
a $\mathbf{k}$-algebra anti-morphism from $A$ to $B$ is nothing but a
$\mathbf{k}$-algebra morphism from $A$ to the opposite $\mathbf{k}$-algebra of
$B$.
\end{remark}

\begin{remark}
Let $A$, $B$ and $C$ be three $\mathbf{k}$-algebras. If $f:B\rightarrow C$ and
$g:A\rightarrow B$ are two $\mathbf{k}$-algebra anti-morphisms, then their
composition $f\circ g:A\rightarrow C$ is a $\mathbf{k}$-algebra morphism (not
anti-morphism). This might sound surprising, but is completely natural, since
a product reversed twice reverts to its original order.
\end{remark}

\subsubsection{The antipode}

We now define the antipode of $\mathbf{k}\left[  S_{n}\right]  $:

\begin{definition}
\label{def.S.S}The \emph{antipode} of $\mathbf{k}\left[  S_{n}\right]  $ is
defined to be the $\mathbf{k}$-linear map%
\begin{align*}
S:\mathbf{k}\left[  S_{n}\right]   &  \rightarrow\mathbf{k}\left[
S_{n}\right]  ,\\
w  &  \mapsto w^{-1}\ \ \ \ \ \ \ \ \ \ \text{for all }w\in S_{n}.
\end{align*}
This is understood as in Definition \ref{def.Tsign.Tsign}. Thus, the map $S$
does \textbf{not} send every element of $\mathbf{k}\left[  S_{n}\right]  $ to
its inverse, but only the standard basis vectors. Any other element of
$\mathbf{k}\left[  S_{n}\right]  $ is written as a $\mathbf{k}$-linear
combination of standard basis vectors, and then the standard basis vectors in
this combination are inverted. Explicitly,%
\[
S\left(  \sum_{w\in S_{n}}\alpha_{w}w\right)  =\sum_{w\in S_{n}}\alpha
_{w}w^{-1}\ \ \ \ \ \ \ \ \ \ \text{for any scalars }\alpha_{w}\in\mathbf{k}.
\]

\end{definition}

\begin{example}
\label{exa.S.exas1}\textbf{(a)} For $n=3$, we have%
\begin{align*}
S\left(  2\operatorname*{id}+\,7s_{1}-9\operatorname*{cyc}\nolimits_{1,2,3}%
\right)   &  =2\underbrace{\operatorname*{id}\nolimits^{-1}}%
_{=\operatorname*{id}}+\,7\underbrace{s_{1}^{-1}}_{=s_{1}}%
-\,9\underbrace{\left(  \operatorname*{cyc}\nolimits_{1,2,3}\right)  ^{-1}%
}_{\substack{=\operatorname*{cyc}\nolimits_{1,3,2}\\\text{(by
(\ref{eq.intro.perms.cycs.inverse}))}}}\\
&  =2\operatorname*{id}+\,7s_{1}-9\operatorname*{cyc}\nolimits_{1,3,2}.
\end{align*}

\textbf{(b)} For each transposition $t_{i,j}\in S_{n}$, we have%
\[
S\left(  t_{i,j}\right)  =t_{i,j}^{-1}=t_{i,j}.
\]
Thus, for any $k\in\left[  n\right]  $, the jucys--murphy $\mathbf{m}_{k}$
satisfies%
\[
S\left(  \mathbf{m}_{k}\right)  =\mathbf{m}_{k}%
\]
(since $\mathbf{m}_{k}$ is a sum of transpositions $t_{i,j}$). \medskip

\textbf{(c)} We have%
\[
S\left(  \nabla\right)  =\nabla\ \ \ \ \ \ \ \ \ \ \text{and}%
\ \ \ \ \ \ \ \ \ \ S\left(  \nabla^{-}\right)  =\nabla^{-}.
\]
Indeed, we have%
\begin{align*}
S\left(  \nabla^{-}\right)   &  =S\left(  \sum_{w\in S_{n}}\left(  -1\right)
^{w}w\right)  \ \ \ \ \ \ \ \ \ \ \left(  \text{since }\nabla^{-}=\sum_{w\in
S_{n}}\left(  -1\right)  ^{w}w\right) \\
&  =\sum_{w\in S_{n}}\left(  -1\right)  ^{w}w^{-1}\ \ \ \ \ \ \ \ \ \ \left(
\text{by the definition of }S\right) \\
&  =\sum_{w\in S_{n}}\underbrace{\left(  -1\right)  ^{w^{-1}}}_{=\left(
-1\right)  ^{w}}\underbrace{\left(  w^{-1}\right)  ^{-1}}_{=w}%
\ \ \ \ \ \ \ \ \ \ \left(
\begin{array}
[c]{c}%
\text{here, we have substituted }w^{-1}\\
\text{for }w\text{ in the sum, because}\\
\text{the map }S_{n}\rightarrow S_{n},\ w\mapsto w^{-1}\\
\text{is a bijection (since }S_{n}\text{ is a group)}%
\end{array}
\right) \\
&  =\sum_{w\in S_{n}}\left(  -1\right)  ^{w}w=\nabla^{-}%
\end{align*}
and similarly $S\left(  \nabla\right)  =\nabla$.

More generally, for any subset $X$ of $\left[  n\right]  $, we have%
\[
S\left(  \nabla_{X}\right)  =\nabla_{X}\ \ \ \ \ \ \ \ \ \ \text{and}%
\ \ \ \ \ \ \ \ \ \ S\left(  \nabla_{X}^{-}\right)  =\nabla_{X}^{-}.
\]
(The proofs are analogous to the proofs of $S\left(  \nabla\right)  =\nabla$
and $S\left(  \nabla^{-}\right)  =\nabla^{-}$.) \medskip

\textbf{(d)} For the somewhere-to-below shuffles $\mathbf{t}_{k}$ (defined in
Definition \ref{def.stb.tk}), we have%
\begin{align*}
S\left(  \mathbf{t}_{k}\right)   &  =S\left(  \operatorname*{cyc}%
\nolimits_{k}+\operatorname*{cyc}\nolimits_{k,k+1}+\operatorname*{cyc}%
\nolimits_{k,k+1,k+2}+\cdots+\operatorname*{cyc}\nolimits_{k,k+1,\ldots
,n}\right) \\
&  =\left(  \operatorname*{cyc}\nolimits_{k}\right)  ^{-1}+\left(
\operatorname*{cyc}\nolimits_{k,k+1}\right)  ^{-1}+\left(  \operatorname*{cyc}%
\nolimits_{k,k+1,k+2}\right)  ^{-1}+\cdots+\left(  \operatorname*{cyc}%
\nolimits_{k,k+1,\ldots,n}\right)  ^{-1}\\
&  =\operatorname*{cyc}\nolimits_{k}+\operatorname*{cyc}\nolimits_{k+1,k}%
+\operatorname*{cyc}\nolimits_{k+2,k+1,k}+\cdots+\operatorname*{cyc}%
\nolimits_{n,n-1,\ldots,k}%
\end{align*}
(by (\ref{eq.intro.perms.cycs.inverse})). I call these images $S\left(
\mathbf{t}_{k}\right)  $ the \emph{below-to-somewhere shuffles}.
\end{example}

The map $S$ has almost the same nice properties as $T_{\operatorname*{sign}}$
(see Theorem \ref{thm.Tsign.auto}), with the exception that it is an
anti-automorphism rather than an automorphism:

\begin{theorem}
\label{thm.S.auto}\textbf{(a)} The map $S:\mathbf{k}\left[  S_{n}\right]
\rightarrow\mathbf{k}\left[  S_{n}\right]  $ is a $\mathbf{k}$-algebra
anti-automorphism. \medskip

\textbf{(b)} It is furthermore an involution (i.e., it satisfies $S\circ
S=\operatorname*{id}$).
\end{theorem}

\begin{proof}
[Proof of Theorem \ref{thm.S.auto}.]This is analogous to Theorem
\ref{thm.Tsign.auto}. (Part \textbf{(b)} boils down to the fact that $\left(
u^{-1}\right)  ^{-1}=u$ for each $u\in S_{n}$, whereas part \textbf{(a)}
relies on the fact that $\left(  uv\right)  ^{-1}=v^{-1}u^{-1}$ for any
$u,v\in S_{n}$.)
\end{proof}

\begin{exercise}
\fbox{1} Recall the $\mathbf{R}_{k}$ from Definition \ref{def.R2R.Rk}. Prove
that $S\left(  \mathbf{R}_{k}\right)  =\mathbf{R}_{k}$ for all $k\in
\mathbb{N}$.
\end{exercise}

\begin{exercise}
\fbox{1} Recall the $\mathbf{X}_{a,b}$ from Definition \ref{def.Xab.Xab}.
Prove that $S\left(  \mathbf{X}_{a,b}\right)  =\mathbf{X}_{b,a}$ for all
$a,b\in\mathbb{N}$.
\end{exercise}

\begin{exercise}
\label{exe.S.STsign}\fbox{1} \textbf{(a)} Prove that the maps $S$ and
$T_{\operatorname*{sign}}$ commute, i.e., we have $S\circ
T_{\operatorname*{sign}}=T_{\operatorname*{sign}}\circ S$. \medskip

\textbf{(b)} Prove that the map $S\circ T_{\operatorname*{sign}}$ is an involution.
\end{exercise}

\begin{exercise}
\label{exe.S.GZ}\fbox{2} Prove that $S\left(  \mathbf{a}\right)  =\mathbf{a}$
for all $\mathbf{a}\in\operatorname*{GZ}\nolimits_{n}$ (where
$\operatorname*{GZ}\nolimits_{n}$ is as in Corollary \ref{cor.GZ.comm}).
\end{exercise}

\begin{remark}
\label{rmk.S.all-G}The antipode $S$ exists not just for $\mathbf{k}\left[
S_{n}\right]  $, but more generally for the group algebra $\mathbf{k}\left[
G\right]  $ of any group $G$. Indeed, Definition \ref{def.S.S} extends
verbatim to this general case (we just need to replace $S_{n}$ by $G$).
Theorem \ref{thm.S.auto} is still valid in this general case as well.
\end{remark}

Even more generally, antipode-like maps exist for a wider class of algebraic
structures known as \emph{Hopf algebras} (see, e.g., \cite[Chapter
1]{GriRei14}). Their theory is actually the source of the name
\textquotedblleft antipode\textquotedblright.

\subsection{Other maps around group algebras}

Next, we will introduce some very simple but useful maps from group algebras
(and, more generally, monoid algebras).

\subsubsection{The counit}

We begin with a morphism from $\mathbf{k}\left[  S_{n}\right]  $ to the base
ring $\mathbf{k}$:

\begin{definition}
\label{def.eps.eps}The \emph{counit} of $\mathbf{k}\left[  S_{n}\right]  $ is
defined to be the $\mathbf{k}$-linear map%
\begin{align*}
\varepsilon:\mathbf{k}\left[  S_{n}\right]   &  \rightarrow\mathbf{k},\\
w  &  \mapsto1\ \ \ \ \ \ \ \ \ \ \text{for all }w\in S_{n}.
\end{align*}
Thus, explicitly,%
\[
\varepsilon\left(  \sum_{w\in S_{n}}\alpha_{w}w\right)  =\sum_{w\in S_{n}%
}\alpha_{w}\ \ \ \ \ \ \ \ \ \ \text{for any scalars }\alpha_{w}\in
\mathbf{k}.
\]

\end{definition}

\begin{example}
\label{exa.eps.mk}For any $k\in\left[  n\right]  $, the counit $\varepsilon$
sends the jucys--murphy $\mathbf{m}_{k}=\sum_{i=1}^{k-1}t_{i,k}$ to%
\begin{align*}
\varepsilon\left(  \mathbf{m}_{k}\right)   &  =\varepsilon\left(  \sum
_{i=1}^{k-1}t_{i,k}\right)  =\sum_{i=1}^{k-1}\underbrace{\varepsilon\left(
t_{i,k}\right)  }_{\substack{=1\\\text{(since }t_{i,k}\in S_{n}\text{)}%
}}\ \ \ \ \ \ \ \ \ \ \left(  \text{since }\varepsilon\text{ is }%
\mathbf{k}\text{-linear}\right) \\
&  =\sum_{i=1}^{k-1}1=k-1.
\end{align*}

\end{example}

\begin{theorem}
\label{thm.eps.mor}The map $\varepsilon$ is a $\mathbf{k}$-algebra morphism.
\end{theorem}

\begin{proof}
This is analogous to Theorem \ref{thm.Tsign.auto} \textbf{(a)} but even easier.
\end{proof}

\begin{remark}
\label{rmk.eps.all-G}The counit $\varepsilon$ exists not just for
$\mathbf{k}\left[  S_{n}\right]  $, but more generally for the monoid algebra
$\mathbf{k}\left[  M\right]  $ of any monoid $M$. Indeed, Definition
\ref{def.eps.eps} extends verbatim to this general case (we just need to
replace $S_{n}$ by $M$). Theorem \ref{thm.eps.mor} is still valid in this
general case as well.
\end{remark}

Simple as it is, the counit $\varepsilon$ can be quite useful in drawing
combinatorial conclusions from identities in $\mathbf{k}\left[  S_{n}\right]
$. For example:

\begin{example}
\label{exa.eps.stirl}Let $\mathbf{k}=\mathbb{Z}$ and $k\in\mathbb{N}$. We know
from Corollary \ref{cor.YJM.ek} that the $k$-th elementary symmetric
polynomial $e_{k}$ applied to the jucys--murphies $\mathbf{m}_{1}%
,\mathbf{m}_{2},\ldots,\mathbf{m}_{n}$ yields%
\[
e_{k}\left(  \mathbf{m}_{1},\mathbf{m}_{2},\ldots,\mathbf{m}_{n}\right)
=\sum_{\substack{w\in S_{n};\\w\text{ has exactly }n-k\text{ orbits}}}w.
\]
Applying the counit $\varepsilon$ to both sides of this equality, we find%
\[
\varepsilon\left(  e_{k}\left(  \mathbf{m}_{1},\mathbf{m}_{2},\ldots
,\mathbf{m}_{n}\right)  \right)  =\varepsilon\left(  \sum_{\substack{w\in
S_{n};\\w\text{ has exactly }n-k\text{ orbits}}}w\right)  .
\]
Since $\varepsilon$ is a $\mathbf{k}$-algebra morphism, we can rewrite this as%
\[
e_{k}\left(  \varepsilon\left(  \mathbf{m}_{1}\right)  ,\varepsilon\left(
\mathbf{m}_{2}\right)  ,\ldots,\varepsilon\left(  \mathbf{m}_{n}\right)
\right)  =\sum_{\substack{w\in S_{n};\\w\text{ has exactly }n-k\text{ orbits}%
}}\varepsilon\left(  w\right)
\]
(here, we were able to rewrite $\varepsilon\left(  e_{k}\left(  \mathbf{m}%
_{1},\mathbf{m}_{2},\ldots,\mathbf{m}_{n}\right)  \right)  $ as $e_{k}\left(
\varepsilon\left(  \mathbf{m}_{1}\right)  ,\varepsilon\left(  \mathbf{m}%
_{2}\right)  ,\ldots,\varepsilon\left(  \mathbf{m}_{n}\right)  \right)  $
because $\mathbf{k}$-algebra morphisms \textquotedblleft
commute\textquotedblright\ with polynomials). Since $\varepsilon\left(
\mathbf{m}_{i}\right)  =i-1$ (by Example \ref{exa.eps.mk}, applied to $i$
instead of $k$), this can be further rewritten as%
\begin{align*}
&  e_{k}\left(  0,1,\ldots,n-1\right) \\
&  =\sum_{\substack{w\in S_{n};\\w\text{ has exactly }n-k\text{ orbits}%
}}\ \ \underbrace{\varepsilon\left(  w\right)  }_{\substack{=1\\\text{(by the
definition of }\varepsilon\text{)}}}\\
&  =\sum_{\substack{w\in S_{n};\\w\text{ has exactly }n-k\text{ orbits}}}1\\
&  =\left(  \text{\# of permutations }w\in S_{n}\text{ that have exactly
}n-k\text{ orbits}\right)  .
\end{align*}
Renaming $k$ as $n-i$, we can rewrite this as%
\[
e_{n-i}\left(  0,1,\ldots,n-1\right)  =\left(  \text{\# of permutations }w\in
S_{n}\text{ that have exactly }i\text{ orbits}\right)  .
\]
This formula (which is true for each $i\in\left\{  0,1,\ldots,n\right\}  $) is
one of the simplest expressions for the \# of permutations $w\in S_{n}$ that
have exactly $i$ orbits. This \# is known as the \emph{unsigned Stirling
number of the first kind} $c\left(  n,i\right)  $, and is commonly discussed
in textbooks on enumerative combinatorics (see, e.g., \cite[Lecture 27,
\S 4.3.4]{22fco} or \cite[\S 1.3]{Stanley-EC1}).
\end{example}

\begin{exercise}
\label{exe.eps.Nab}\fbox{1} Let $\mathbf{a}\in\mathbf{k}\left[  S_{n}\right]
$. Prove that%
\[
\nabla\mathbf{a}=\mathbf{a}\nabla=\varepsilon\left(  \mathbf{a}\right)
\cdot\nabla.
\]

\end{exercise}

For later use, let us also observe how the counit $\varepsilon$ (introduced in
Definition \ref{def.eps.eps}) acts on an $X$-integral:

\begin{proposition}
\label{prop.eps.Xint}Let $X$ be a subset of $\left[  n\right]  $. Then,
$\varepsilon\left(  \nabla_{X}\right)  =\left\vert X\right\vert !$.
\end{proposition}

\begin{proof}
Define the set $S_{n,X}$ as in Proposition \ref{prop.intX.basics}. Then,
Proposition \ref{prop.intX.basics} \textbf{(a)} shows that $S_{n,X}$ is a
subgroup of $S_{n}$, and is isomorphic to $S_{X}$. Hence, $\left\vert
S_{n,X}\right\vert =\left\vert S_{X}\right\vert $ (since two isomorphic groups
always have the same size). But we know that $\left\vert S_{X}\right\vert
=\left\vert X\right\vert !$.

However, Proposition \ref{prop.intX.basics} \textbf{(b)} yields $\nabla
_{X}=\sum\limits_{w\in S_{n,X}}w$. Thus,%
\begin{align*}
\varepsilon\left(  \nabla_{X}\right)   &  =\varepsilon\left(  \sum
\limits_{w\in S_{n,X}}w\right)  =\sum\limits_{w\in S_{n,X}}%
\underbrace{\varepsilon\left(  w\right)  }_{\substack{=1\\\text{(by the
definition of }\varepsilon\text{,}\\\text{since }w\in S_{n,X}\subseteq
S_{n}\text{)}}}\ \ \ \ \ \ \ \ \ \ \left(
\begin{array}
[c]{c}%
\text{since the map }\varepsilon\\
\text{is }\mathbf{k}\text{-linear}%
\end{array}
\right) \\
&  =\sum\limits_{w\in S_{n,X}}1=\left\vert S_{n,X}\right\vert \cdot
1=\left\vert S_{n,X}\right\vert =\left\vert S_{X}\right\vert =\left\vert
X\right\vert !.
\end{align*}
This proves Proposition \ref{prop.eps.Xint}.
\end{proof}

\subsubsection{Morphisms of monoids yield morphisms of monoid algebras}

Let us now explore some morphisms between different monoid algebras.

We begin with a simple construction in linear algebra:\footnote{Recall the
free $\mathbf{k}$-module $\mathbf{k}^{\left(  S\right)  }$ on a set $S$
(defined in Definition \ref{def.monalg.freemod}).}

\begin{definition}
\label{def.linearization.gen}Let $M$ and $N$ be two sets, and let
$f:M\rightarrow N$ be any map. Then,
\[
f_{\ast}:\mathbf{k}^{\left(  M\right)  }\rightarrow\mathbf{k}^{\left(
N\right)  }%
\]
shall denote the $\mathbf{k}$-linear map that sends each standard basis vector
$e_{m}\in\mathbf{k}^{\left(  M\right)  }$ to the standard basis vector
$e_{f\left(  m\right)  }\in\mathbf{k}^{\left(  N\right)  }$. This map
$f_{\ast}$ is called the \emph{linearization} of $f$. Explicitly, it is given
by%
\begin{equation}
f_{\ast}\left(  \sum_{m\in M}\alpha_{m}e_{m}\right)  =\sum_{m\in M}\alpha
_{m}e_{f\left(  m\right)  } \label{eq.def.linearization.gen.expl}%
\end{equation}
for any scalars $\alpha_{m}\in\mathbf{k}$.
\end{definition}

\begin{proposition}
\label{prop.linearization.ject}If the map $f$ in Definition
\ref{def.linearization.gen} is injective, then so is $f_{\ast}$. The same
holds if \textquotedblleft injective\textquotedblright\ is replaced by
\textquotedblleft surjective\textquotedblright\ or \textquotedblleft
bijective\textquotedblright.
\end{proposition}

\begin{proof}
Straightforward and left to the reader.
\end{proof}

\begin{proposition}
\label{prop.linearization.mon}Let $M$ and $N$ be two monoids, and let
$f:M\rightarrow N$ be a monoid morphism. Then the linearization $f_{\ast
}:\mathbf{k}^{\left(  M\right)  }\rightarrow\mathbf{k}^{\left(  N\right)  }$
is a $\mathbf{k}$-algebra morphism from $\mathbf{k}\left[  M\right]  $ to
$\mathbf{k}\left[  N\right]  $.

Using Convention \ref{conv.monalg.em=m}, we can rewrite the above explicit
formula (\ref{eq.def.linearization.gen.expl}) as%
\[
f_{\ast}\left(  \sum_{m\in M}\alpha_{m}m\right)  =\sum_{m\in M}\alpha
_{m}f\left(  m\right)  .
\]
Thus, $f_{\ast}$ is the $\mathbf{k}$-linear map from $\mathbf{k}\left[
M\right]  $ to $\mathbf{k}\left[  N\right]  $ that sends each standard basis
vector $m\in M$ to $f\left(  m\right)  \in N$.
\end{proposition}

\begin{proof}
Straightforward and left to the reader. (Again, use linearity as in the proof
of Theorem \ref{thm.Tsign.auto} \textbf{(a)}.)
\end{proof}

Thus, for any monoid morphism $f:M\rightarrow N$, we have constructed a
$\mathbf{k}$-algebra morphism $f_{\ast}:\mathbf{k}\left[  M\right]
\rightarrow\mathbf{k}\left[  N\right]  $ between the respective monoid
algebras. Not all $\mathbf{k}$-algebra morphisms between monoid algebras arise
in this way, but the ones that do are particularly easy and well-behaved.

Monoid morphisms between groups are group morphisms. Here is a simple one
between symmetric groups:

\begin{definition}
\label{def.default-ex.IXn}Let $X$ be a subset of $\left[  n\right]  $. Then,
for each permutation $\sigma\in S_{X}$, we can define a permutation
$\sigma\uparrow^{n}\ \in S_{n}$ by setting%
\[
\left(  \sigma\uparrow^{n}\right)  \left(  i\right)  =%
\begin{cases}
\sigma\left(  i\right)  , & \text{if }i\in X;\\
i, & \text{if }i\notin X
\end{cases}
\ \ \ \ \ \ \ \ \ \ \text{for each }i\in\left[  n\right]  .
\]
This permutation $\sigma\uparrow^{n}$ is called the \emph{default extension}
of $\sigma$ to $\left[  n\right]  $. In words, it can be described as
\textquotedblleft permuting the elements of $X$ like $\sigma$ does, and
leaving the elements of $\left[  n\right]  \setminus X$
unchanged\textquotedblright. (It is easy to check that it really is a
permutation of $\left[  n\right]  $.)

We define the map%
\begin{align*}
I_{X\rightarrow n}:S_{X}  &  \rightarrow S_{n},\\
\sigma &  \mapsto\sigma\uparrow^{n},
\end{align*}
which sends each permutation $\sigma\in S_{X}$ to its default extension
$\sigma\uparrow^{n}$. This map $I_{X\rightarrow n}$ will be called the
\emph{default embedding} of $S_{X}$ into $S_{n}$.
\end{definition}

\begin{proposition}
\label{prop.default-ex.mor}Let $X$ be a subset of $\left[  n\right]  $. Then:
\medskip

\textbf{(a)} The default embedding $I_{X\rightarrow n}$ is an injective group
morphism. \medskip

\textbf{(b)} Its image is the subgroup
\[
S_{n,X}=\left\{  w\in S_{n}\ \mid\ w\left(  i\right)  =i\text{ for all
}i\notin X\right\}
\]
of $S_{n}$ (already defined in Proposition \ref{prop.intX.basics}). \medskip

\textbf{(c)} The default embedding $I_{X\rightarrow n}$ sends cycles in
$S_{X}$ to the \textquotedblleft same\textquotedblright\ cycles in $S_{n}$, in
the sense that it sends each cycle $\operatorname*{cyc}\nolimits_{i_{1}%
,i_{2},\ldots,i_{k}}\in S_{X}$ to $\operatorname*{cyc}\nolimits_{i_{1}%
,i_{2},\ldots,i_{k}}\in S_{n}$. \medskip

\textbf{(d)} The default embedding $I_{X\rightarrow n}$ preserves signs; i.e.,
we have
\[
\left(  -1\right)  ^{\sigma}=\left(  -1\right)  ^{I_{X\rightarrow n}\left(
\sigma\right)  }\ \ \ \ \ \ \ \ \ \ \text{for any }\sigma\in S_{X}.
\]

\end{proposition}

\begin{proof}
[Proof sketch.]\textbf{(a)} This is straightforward. \medskip

\textbf{(b)} This is also pretty easy. Indeed, it is clear that the image of
$I_{X\rightarrow n}$ is a subset of $S_{n,X}$; it thus remains to show that
all $w\in S_{n,X}$ belong to the image of $I_{X\rightarrow n}$. So let $w\in
S_{n,X}$ be arbitrary. Then, Claim 1 in the proof of Proposition
\ref{prop.intX.basics} \textbf{(a)} shows that $w\left(  X\right)  \subseteq
X$ (so that the restriction $w\mid_{X}^{X}$ is really a well-defined map from
$X$ to $X$) and that $\left.  w\mid_{X}^{X}\right.  \in S_{X}$. Now, it is
easy to see that $w=I_{X\rightarrow n}\left(  w\mid_{X}^{X}\right)  $ (since
$w\in S_{n,X}$ entails that $w\left(  i\right)  =i$ for all $i\notin X$,
whereas the definition of $w\mid_{X}^{X}$ entails that $w\left(  i\right)
=\left(  w\mid_{X}^{X}\right)  \left(  i\right)  $ for all $i\in X$). Thus,
$w$ belongs to the image of $I_{X\rightarrow n}$. This completes the proof of
Proposition \ref{prop.default-ex.mor} \textbf{(b)}. \medskip

\textbf{(c)} Let $i_{1},i_{2},\ldots,i_{k}$ be $k$ distinct elements of $X$.
By its definition, the cycle $\operatorname*{cyc}\nolimits_{i_{1},i_{2}%
,\ldots,i_{k}}\in S_{X}$ sends the elements $i_{1},i_{2},\ldots,i_{k}$ to
$i_{2},i_{3},\ldots,i_{k},i_{1}$, while leaving all remaining elements of $X$
unchanged. Its default extension $\operatorname*{cyc}\nolimits_{i_{1}%
,i_{2},\ldots,i_{k}}\uparrow^{n}$ thus does the same thing, but also leaves
the elements of $\left[  n\right]  \setminus X$ unchanged. In other words, it
sends the elements $i_{1},i_{2},\ldots,i_{k}$ to $i_{2},i_{3},\ldots
,i_{k},i_{1}$, while leaving all remaining elements of $\left[  n\right]  $
unchanged. But this is exactly what the cycle $\operatorname*{cyc}%
\nolimits_{i_{1},i_{2},\ldots,i_{k}}\in S_{n}$ does. Therefore,%
\[
\underbrace{\operatorname*{cyc}\nolimits_{i_{1},i_{2},\ldots,i_{k}}}_{\text{in
}S_{X}}\uparrow^{n}\ =\underbrace{\operatorname*{cyc}\nolimits_{i_{1}%
,i_{2},\ldots,i_{k}}}_{\text{in }S_{n}}.
\]
Now, the definition of $I_{X\rightarrow n}$ yields%
\[
I_{X\rightarrow n}\left(  \underbrace{\operatorname*{cyc}\nolimits_{i_{1}%
,i_{2},\ldots,i_{k}}}_{\text{in }S_{X}}\right)
=\underbrace{\operatorname*{cyc}\nolimits_{i_{1},i_{2},\ldots,i_{k}}%
}_{\text{in }S_{X}}\uparrow^{n}\ =\underbrace{\operatorname*{cyc}%
\nolimits_{i_{1},i_{2},\ldots,i_{k}}}_{\text{in }S_{n}}.
\]
This proves Proposition \ref{prop.default-ex.mor} \textbf{(c)}. \medskip

\textbf{(d)} The easiest way to prove this is to recall the following two facts:

\begin{statement}
\textit{Claim 1:} Each permutation in $S_{X}$ can be written as a product of
(finitely many) transpositions.
\end{statement}

\begin{statement}
\textit{Claim 2:} If a permutation $\sigma\in S_{X}$ is a product of $k$
transpositions, then $\left(  -1\right)  ^{\sigma}=\left(  -1\right)  ^{k}$.
\end{statement}

Claim 1 is \cite[Exercise 5.15 \textbf{(b)}]{detnotes}, whereas Claim 2 is
\cite[Exercise 5.15 \textbf{(c)}]{detnotes}.

Now, let $\sigma\in S_{X}$. Then, Claim 1 shows that $\sigma$ can be written
as a product of (finitely many) transpositions. In other words, there exists
some transpositions $t_{i_{1},j_{1}},t_{i_{2},j_{2}},\ldots,t_{i_{k},j_{k}}\in
S_{X}$ such that%
\begin{equation}
\sigma=t_{i_{1},j_{1}}t_{i_{2},j_{2}}\cdots t_{i_{k},j_{k}}.
\label{pf.prop.default-ex.mor.d.0}%
\end{equation}
Consider these transpositions. Thus, $\sigma$ is a product of $k$
transpositions. Hence, Claim 2 yields $\left(  -1\right)  ^{\sigma}=\left(
-1\right)  ^{k}$. However, for any two distinct elements $u,v\in X$, the map
$I_{X\rightarrow n}$ sends the transposition $t_{u,v}\in S_{X}$ to the
transposition $t_{u,v}\in S_{n}$ (since Proposition \ref{prop.default-ex.mor}
\textbf{(c)} yields that $I_{X\rightarrow n}$ maps $\operatorname*{cyc}%
\nolimits_{u,v}\in S_{X}$ to $\operatorname*{cyc}\nolimits_{u,v}\in S_{n}$,
but the transposition $t_{u,v}$ is just the cycle $\operatorname*{cyc}%
\nolimits_{u,v}$). In other words,
\begin{equation}
I_{X\rightarrow n}\left(  t_{u,v}\right)  =t_{u,v}%
\ \ \ \ \ \ \ \ \ \ \text{for any distinct }u,v\in X.
\label{pf.prop.default-ex.mor.d.t}%
\end{equation}
Now, applying the map $I_{X\rightarrow n}$ to the equality
(\ref{pf.prop.default-ex.mor.d.0}), we obtain%
\begin{align*}
I_{X\rightarrow n}\left(  \sigma\right)   &  =I_{X\rightarrow n}\left(
t_{i_{1},j_{1}}t_{i_{2},j_{2}}\cdots t_{i_{k},j_{k}}\right) \\
&  =\underbrace{I_{X\rightarrow n}\left(  t_{i_{1},j_{1}}\right)
}_{\substack{=t_{i_{1},j_{1}}\\\text{(by (\ref{pf.prop.default-ex.mor.d.t}))}%
}}\underbrace{I_{X\rightarrow n}\left(  t_{i_{2},j_{2}}\right)  }%
_{\substack{=t_{i_{2},j_{2}}\\\text{(by (\ref{pf.prop.default-ex.mor.d.t}))}%
}}\cdots\underbrace{I_{X\rightarrow n}\left(  t_{i_{k},j_{k}}\right)
}_{\substack{=t_{i_{k},j_{k}}\\\text{(by (\ref{pf.prop.default-ex.mor.d.t}))}%
}}\\
&  \ \ \ \ \ \ \ \ \ \ \ \ \ \ \ \ \ \ \ \ \left(  \text{since }%
I_{X\rightarrow n}\text{ is a group morphism}\right) \\
&  =t_{i_{1},j_{1}}t_{i_{2},j_{2}}\cdots t_{i_{k},j_{k}}.
\end{align*}
This shows that $I_{X\rightarrow n}\left(  \sigma\right)  $ is a product of
$k$ transpositions in $S_{n}$. Hence, we have $\left(  -1\right)
^{I_{X\rightarrow n}\left(  \sigma\right)  }=\left(  -1\right)  ^{k}$ (by
Claim 2, applied to $\left[  n\right]  $ instead of $X$). Comparing this with
$\left(  -1\right)  ^{\sigma}=\left(  -1\right)  ^{k}$, we obtain $\left(
-1\right)  ^{\sigma}=\left(  -1\right)  ^{I_{X\rightarrow n}\left(
\sigma\right)  }$. This proves Proposition \ref{prop.default-ex.mor}
\textbf{(d)}.
\end{proof}

\begin{remark}
\label{rmk.default-ex.abuse1}Default embeddings facilitate a rather popular
abuse of notation: Namely, if $X$ is a subset of $\left[  n\right]  $, then
each permutation $\sigma\in S_{X}$ can be identified with its default
extension $I_{X\rightarrow n}\left(  \sigma\right)  =\sigma\uparrow^{n}\ \in
S_{n}$ (by abuse of notation). This identification is usually harmless, since
the default embedding $I_{X\rightarrow n}$ is an injective group morphism from
$S_{X}$ to $S_{n}$ (by Proposition \ref{prop.default-ex.mor} \textbf{(a)}) and
thus respects multiplication (i.e., it does not matter if we compute a product
$\sigma\tau$ of two permutations in $S_{X}$ or instead compute the product
$I_{X\rightarrow n}\left(  \sigma\right)  I_{X\rightarrow n}\left(
\tau\right)  $ of their default extensions).\footnotemark

As a consequence of this identification, the group $S_{X}$ becomes identified
with its image under the morphism $I_{X\rightarrow n}$; this image is
$S_{n,X}$ (by Proposition \ref{prop.default-ex.mor} \textbf{(b)}). We shall
mostly avoid this identification in the following, but many authors use it, so
it is worth knowing.
\end{remark}

\footnotetext{Some care needs to be taken with this identification. For
example, the \# of fixed points of a permutation $\sigma\in S_{X}$ is
\textbf{not} the \# of fixed points of its default extension $I_{X\rightarrow
n}\left(  \sigma\right)  $ (unless $X=\left[  n\right]  $), since the latter
has all the elements of $\left[  n\right]  \setminus X$ as additional fixed
points. Thus, when counting (or listing) the fixed points, we must distinguish
between $\sigma$ and $I_{X\rightarrow n}\left(  \sigma\right)  $. Similarly
for descents and various other features.}

\begin{corollary}
\label{cor.default-ex.injmor}Let $X$ be a subset of $\left[  n\right]  $.
Then, the linearization of the default embedding $I_{X\rightarrow n}%
:S_{X}\rightarrow S_{n}$ is an injective $\mathbf{k}$-algebra morphism%
\[
\left(  I_{X\rightarrow n}\right)  _{\ast}:\mathbf{k}\left[  S_{X}\right]
\rightarrow\mathbf{k}\left[  S_{n}\right]  .
\]

\end{corollary}

\begin{proof}
Proposition \ref{prop.linearization.ject} shows that $\left(  I_{X\rightarrow
n}\right)  _{\ast}$ is injective (since $I_{X\rightarrow n}$ is injective).
Proposition \ref{prop.linearization.mon} shows that $\left(  I_{X\rightarrow
n}\right)  _{\ast}$ is a $\mathbf{k}$-algebra morphism.
\end{proof}

\begin{definition}
We will refer to the injective $\mathbf{k}$-algebra morphism $\left(
I_{X\rightarrow n}\right)  _{\ast}$ in Corollary \ref{cor.default-ex.injmor}
as a \emph{default embedding of group algebras}.
\end{definition}

\begin{example}
\label{exa.default-ex.Nabla}Let $X$ be a subset of $\left[  n\right]  $. Then,
the $X$-integral $\nabla_{X}$ is%
\begin{align}
\nabla_{X}  &  =\sum_{w\in S_{n,X}}w\ \ \ \ \ \ \ \ \ \ \left(  \text{by
Proposition \ref{prop.intX.basics} \textbf{(b)}}\right) \nonumber\\
&  =\sum_{u\in S_{X}}I_{X\rightarrow n}\left(  u\right)
\ \ \ \ \ \ \ \ \ \ \left(
\begin{array}
[c]{c}%
\text{here, we have substituted }I_{X\rightarrow n}\left(  u\right) \\
\text{for }w\text{ in the sum, since the}\\
\text{map }S_{X}\rightarrow S_{n,X},\ u\mapsto I_{X\rightarrow n}\left(
u\right) \\
\text{is a bijection (by parts \textbf{(a)} and \textbf{(b)}}\\
\text{of Proposition \ref{prop.default-ex.mor})}%
\end{array}
\right) \nonumber\\
&  =\left(  I_{X\rightarrow n}\right)  _{\ast}\left(  \sum_{u\in S_{X}%
}u\right)  . \label{eq.exa.default-ex.Nabla.1}%
\end{align}
Similarly, using Proposition \ref{prop.default-ex.mor} \textbf{(d)}, we can
show that%
\[
\nabla_{X}^{-}=\left(  I_{X\rightarrow n}\right)  _{\ast}\left(  \sum_{u\in
S_{X}}\left(  -1\right)  ^{u}u\right)  .
\]

\end{example}

\begin{remark}
\label{rmk.default-ex.abuse2}Let $X$ be a subset of $\left[  n\right]  $. If
you are willing to identify each permutation $\sigma\in S_{X}$ with its
default extension $I_{X\rightarrow n}\left(  \sigma\right)  $ (by the abuse of
notation explained in Remark \ref{rmk.default-ex.abuse1}), then you can also
identify each element $\mathbf{a}\in\mathbf{k}\left[  S_{X}\right]  $ with its
image $\left(  I_{X\rightarrow n}\right)  _{\ast}\left(  \mathbf{a}\right)
\in\mathbf{k}\left[  S_{n}\right]  $ under the default embedding of group
algebras. Thus, $\mathbf{k}\left[  S_{X}\right]  $ becomes a $\mathbf{k}%
$-subalgebra of $\mathbf{k}\left[  S_{n}\right]  $, and the equality
(\ref{eq.exa.default-ex.Nabla.1}) can be rewritten as $\nabla_{X}=\sum_{u\in
S_{X}}u$.
\end{remark}

Default embeddings $I_{X\rightarrow n}$ are used particularly often when $X$
is itself a set of the form $\left[  m\right]  =\left\{  1,2,\ldots,m\right\}
$ for some $m\leq n$:

\begin{definition}
\label{def.default-ex.Imn}Let $m$ be a nonnegative integer with $m\leq n$.
Then, $\left[  m\right]  \subseteq\left[  n\right]  $. The default embedding
$I_{\left[  m\right]  \rightarrow n}:S_{\left[  m\right]  }\rightarrow S_{n}$
will be called $I_{m\rightarrow n}$. This is an injective group morphism from
$S_{\left[  m\right]  }=S_{m}$ to $S_{n}$.
\end{definition}

Thus, for any nonnegative integer $m\leq n$, we have a default embedding%
\[
I_{m\rightarrow n}=I_{\left[  m\right]  \rightarrow n}:S_{m}\rightarrow
S_{n},
\]
which is an injective group morphism, and thus (by Corollary
\ref{cor.default-ex.injmor}) produces an injective $\mathbf{k}$-algebra
morphism%
\[
\left(  I_{m\rightarrow n}\right)  _{\ast}:\mathbf{k}\left[  S_{m}\right]
\rightarrow\mathbf{k}\left[  S_{n}\right]  .
\]
In terms of one-line notations, $I_{m\rightarrow n}$ is given by
\[
I_{m\rightarrow n}\left(  \operatorname*{oln}\left(  i_{1},i_{2},\ldots
,i_{m}\right)  \right)  =\operatorname*{oln}\left(  i_{1},i_{2},\ldots
,i_{m},m+1,m+2,\ldots,n\right)
\]
for every permutation $\operatorname*{oln}\left(  i_{1},i_{2},\ldots
,i_{m}\right)  \in S_{m}$.

Composing two of these default embeddings yields another:

\begin{proposition}
For any nonnegative integers $m\leq n\leq p$, we have%
\[
I_{n\rightarrow p}\circ I_{m\rightarrow n}=I_{m\rightarrow p}%
\ \ \ \ \ \ \ \ \ \ \text{and}\ \ \ \ \ \ \ \ \ \ \left(  I_{n\rightarrow
p}\right)  _{\ast}\circ\left(  I_{m\rightarrow n}\right)  _{\ast}=\left(
I_{m\rightarrow p}\right)  _{\ast}.
\]

\end{proposition}

\begin{proof}
Straightforward.
\end{proof}

\begin{remark}
\label{rmk.default-ex.abuse3}Once again, things get easier if you are a bit
sloppy. By Remark \ref{rmk.default-ex.abuse1} (applied to $X=\left[  m\right]
$), we can identify each permutation $\sigma\in S_{m}$ with its default
extension $I_{m\rightarrow n}\left(  \sigma\right)  =\sigma\uparrow^{n}\ \in
S_{n}$ for any $m\leq n$. Likewise, by Remark \ref{rmk.default-ex.abuse2}, we
can identify each element $\mathbf{a}\in\mathbf{k}\left[  S_{m}\right]  $ with
its image $\left(  I_{m\rightarrow n}\right)  _{\ast}\left(  \mathbf{a}%
\right)  \in\mathbf{k}\left[  S_{n}\right]  $ under the default embedding of
group algebras. This is an abuse of notation, but it is useful. It allows us
to pretend that the default embeddings $I_{m\rightarrow n}$ and $\left(
I_{m\rightarrow n}\right)  _{\ast}$ are actual embeddings, i.e., that we have
a chain of subgroups%
\[
S_{0}\subseteq S_{1}\subseteq S_{2}\subseteq S_{3}\subseteq\cdots
\ \ \ \ \ \ \ \ \ \ \left(  \text{note that the first }\subseteq\text{ sign is
an }=\right)
\]
and a chain of subalgebras%
\[
\mathbf{k}\left[  S_{0}\right]  \subseteq\mathbf{k}\left[  S_{1}\right]
\subseteq\mathbf{k}\left[  S_{2}\right]  \subseteq\mathbf{k}\left[
S_{3}\right]  \subseteq\cdots\ \ \ \ \ \ \ \ \ \ \left(  \text{again, the
first }\subseteq\text{ sign is an }=\right)  .
\]
Keep in mind that this is just a figure of speech, and beware of the confusion
that can arise when you compare elements of $\mathbf{k}\left[  S_{n}\right]  $
for different $n$'s. For example, the $3$rd jucys--murphy $\mathbf{m}%
_{3}=t_{1,3}+t_{2,3}$ of $\mathbf{k}\left[  S_{3}\right]  $ agrees with the
$3$rd jucys--murphy $\mathbf{m}_{3}=t_{1,3}+t_{2,3}$ of $\mathbf{k}\left[
S_{4}\right]  $ (since the default embedding $I_{3\rightarrow4}$ sends each
$t_{i,j}\in S_{3}$ to $t_{i,j}\in S_{4}$), but the $3$rd somewhere-to-below
shuffle $\mathbf{t}_{3}=\operatorname*{cyc}\nolimits_{3}$ of $\mathbf{k}%
\left[  S_{3}\right]  $ does not agree with the $3$rd somewhere-to-below
shuffle $\mathbf{t}_{3}=\operatorname*{cyc}\nolimits_{3}+\operatorname*{cyc}%
\nolimits_{3,4}$ of $\mathbf{k}\left[  S_{4}\right]  $. When necessary, be
explicit about what $n$ is.

Beware also that the center $Z\left(  \mathbf{k}\left[  S_{2}\right]  \right)
$ is not a subset of $Z\left(  \mathbf{k}\left[  S_{3}\right]  \right)  $, no
matter how hard you abuse notation. (More generally, if $A$ is a subring of a
ring $B$, then $Z\left(  A\right)  $ need not be a subring of $Z\left(
B\right)  $.)
\end{remark}

Default embeddings (and Proposition \ref{prop.default-ex.mor} \textbf{(d)}
specifically) help us prove the following purely combinatorial lemma, which is
a close relative of Proposition \ref{prop.int.wXi=Xi}:

\begin{lemma}
\label{lem.Xk.bij1}Let $X_{1},X_{2},\ldots,X_{k}$ be $k$ disjoint subsets of
$\left[  n\right]  $ such that $X_{1}\cup X_{2}\cup\cdots\cup X_{k}=\left[
n\right]  $.

Define a subset $K$ of $S_{n}$ by%
\[
K:=\left\{  w\in S_{n}\ \mid\ w\left(  X_{i}\right)  \subseteq X_{i}\text{ for
all }i\in\left[  k\right]  \right\}  .
\]

For each permutation $w\in K$ and each $j\in\left[  k\right]  $, we define a
map $w^{\left\langle j\right\rangle }:X_{j}\rightarrow X_{j}$ by setting%
\[
w^{\left\langle j\right\rangle }\left(  x\right)  =w\left(  x\right)
\ \ \ \ \ \ \ \ \ \ \text{for each }x\in X_{j}.
\]
Then: \medskip

\textbf{(a)} This map $w^{\left\langle j\right\rangle }$ is well-defined and
belongs to the symmetric group $S_{X_{j}}$ for each $w\in K$ and each
$j\in\left[  k\right]  $. \medskip

\textbf{(b)} The map%
\begin{align*}
K  &  \rightarrow S_{X_{1}}\times S_{X_{2}}\times\cdots\times S_{X_{k}},\\
w  &  \mapsto\left(  w^{\left\langle 1\right\rangle },w^{\left\langle
2\right\rangle },\ldots,w^{\left\langle k\right\rangle }\right)
\end{align*}
is a bijection. \medskip

\textbf{(c)} For each $w\in K$, we have
\[
\left(  -1\right)  ^{w}=\left(  -1\right)  ^{w^{\left\langle 1\right\rangle }%
}\left(  -1\right)  ^{w^{\left\langle 2\right\rangle }}\cdots\left(
-1\right)  ^{w^{\left\langle k\right\rangle }}.
\]

\end{lemma}

\begin{fineprint}
\begin{proof}
For any subset $X$ of $\left[  n\right]  $, define a subset $S_{n,X}$ of
$S_{n}$ as in Proposition \ref{prop.intX.basics}. If $X$ is any subset of
$\left[  n\right]  $, then Proposition \ref{prop.intX.basics} \textbf{(a)}
shows that the set $S_{n,X}$ is a subgroup of $S_{n}$, and is isomorphic to
$S_{X}$ via the group isomorphism%
\begin{align*}
S_{n,X}  &  \rightarrow S_{X},\\
w  &  \mapsto\left.  w\mid_{X}^{X}\right.
\end{align*}
(where $w\mid_{X}^{X}$ means the restriction of $w$ to $X$, regarded as a map
from $X$ to $X$). Let us denote this group isomorphism by $\varphi_{X}$.

Thus, in particular, we obtain a group isomorphism $\varphi_{X_{j}}%
:S_{n,X_{j}}\rightarrow S_{X_{j}}$ for each $j\in\left[  k\right]  $. The
direct product of these $k$ isomorphisms $\varphi_{X_{j}}$ is a group
isomorphism%
\begin{align*}
\varphi_{X_{1}}\times\varphi_{X_{2}}\times\cdots\times\varphi_{X_{k}%
}:S_{n,X_{1}}\times S_{n,X_{2}}\times\cdots\times S_{n,X_{k}}  &  \rightarrow
S_{X_{1}}\times S_{X_{2}}\times\cdots\times S_{X_{k}},\\
\left(  w_{1},w_{2},\ldots,w_{k}\right)   &  \mapsto\left(  \varphi_{X_{1}%
}\left(  w_{1}\right)  ,\ \varphi_{X_{2}}\left(  w_{2}\right)  ,\ \ldots
,\ \varphi_{X_{k}}\left(  w_{k}\right)  \right)  .
\end{align*}
Let us denote this group isomorphism $\varphi_{X_{1}}\times\varphi_{X_{2}%
}\times\cdots\times\varphi_{X_{k}}$ by $\varphi_{\Pi}$. In particular,
$\varphi_{\Pi}$ is a bijection (since $\varphi_{\Pi}$ is a group isomorphism).

For each $j\in\left[  k\right]  $ and each $w\in K$, we define a map
$w^{\left(  j\right)  }:\left[  n\right]  \rightarrow\left[  n\right]  $ by
setting%
\begin{equation}
w^{\left(  j\right)  }\left(  x\right)  =%
\begin{cases}
w\left(  x\right)  , & \text{if }x\in X_{j};\\
x, & \text{if }x\notin X_{j}%
\end{cases}
\ \ \ \ \ \ \ \ \ \ \text{for each }x\in\left[  n\right]  .
\label{pf.lem.Xk.bij1.2}%
\end{equation}
Then, Claim 2 from the proof of Proposition \ref{prop.int.wXi=Xi} \textbf{(b)}
shows that%
\begin{equation}
w^{\left(  j\right)  }\in S_{n,X_{j}}\ \ \ \ \ \ \ \ \ \ \text{for each }%
j\in\left[  k\right]  \text{ and }w\in K. \label{pf.lem.Xk.bij1.3}%
\end{equation}
We now claim the following:

\begin{statement}
\textit{Claim 1:} Let $j\in\left[  k\right]  $ and $w\in K$. Then, the map
$w^{\left\langle j\right\rangle }$ is well-defined and equals $\varphi_{X_{j}%
}\left(  w^{\left(  j\right)  }\right)  $.
\end{statement}

\begin{proof}
[Proof of Claim 1.]We have $w\in K$. In other words, $w\in S_{n}$ is a
permutation with the property that $w\left(  X_{i}\right)  \subseteq X_{i}$
for all $i\in\left[  k\right]  $ (by the definition of $K$). Applying the
latter property to $i=j$, we obtain $w\left(  X_{j}\right)  \subseteq X_{j}$.
Hence, for each $x\in X_{j}$, we have $w\left(  \underbrace{x}_{\in X_{j}%
}\right)  \in w\left(  X_{j}\right)  \subseteq X_{j}$. Thus, the map
$w^{\left\langle j\right\rangle }:X_{j}\rightarrow X_{j}$ is well-defined
(since the definition of this map says that $w^{\left\langle j\right\rangle
}\left(  x\right)  =w\left(  x\right)  $ for each $x\in X_{j}$). It remains to
show that it equals $\varphi_{X_{j}}\left(  w^{\left(  j\right)  }\right)  $.

Indeed, we first observe that $w^{\left(  j\right)  }\in S_{n,X_{j}}$ (by
(\ref{pf.lem.Xk.bij1.3})). Hence, $\varphi_{X_{j}}\left(  w^{\left(  j\right)
}\right)  $ is well-defined and is a map from $X_{j}$ to $X_{j}$. The
definition of $\varphi_{X_{j}}$ yields $\varphi_{X_{j}}\left(  w^{\left(
j\right)  }\right)  =\left.  w^{\left(  j\right)  }\mid_{X_{j}}^{X_{j}%
}\right.  $ (where $w^{\left(  j\right)  }\mid_{X_{j}}^{X_{j}}$ means the
restriction of $w^{\left(  j\right)  }$ to $X_{j}$, regarded as a map from
$X_{j}$ to $X_{j}$). Hence, for each $x\in X_{j}$, we have%
\begin{align*}
\underbrace{\left(  \varphi_{X_{j}}\left(  w^{\left(  j\right)  }\right)
\right)  }_{=\left.  w^{\left(  j\right)  }\mid_{X_{j}}^{X_{j}}\right.
}\left(  x\right)   &  =\left(  w^{\left(  j\right)  }\mid_{X_{j}}^{X_{j}%
}\right)  \left(  x\right)  =w^{\left(  j\right)  }\left(  x\right) \\
&  =%
\begin{cases}
w\left(  x\right)  , & \text{if }x\in X_{j};\\
x, & \text{if }x\notin X_{j}%
\end{cases}
\ \ \ \ \ \ \ \ \ \ \left(  \text{by the definition of }w^{\left(  j\right)
}\right) \\
&  =w\left(  x\right)  \ \ \ \ \ \ \ \ \ \ \left(  \text{since }x\in
X_{j}\right) \\
&  =w^{\left\langle j\right\rangle }\left(  x\right)
\ \ \ \ \ \ \ \ \ \ \left(  \text{since the definition of }w^{\left\langle
j\right\rangle }\text{ yields }w^{\left\langle j\right\rangle }\left(
x\right)  =w\left(  x\right)  \right)  .
\end{align*}
Thus, $\varphi_{X_{j}}\left(  w^{\left(  j\right)  }\right)  =w^{\left\langle
j\right\rangle }$. In other words, $w^{\left\langle j\right\rangle }$ equals
$\varphi_{X_{j}}\left(  w^{\left(  j\right)  }\right)  $. Thus, the proof of
Claim 1 is complete.
\end{proof}

\medskip

\textbf{(a)} Let $w\in K$ and $j\in\left[  k\right]  $. We must prove that the
map $w^{\left\langle j\right\rangle }$ is well-defined and belongs to the
symmetric group $S_{X_{j}}$.

Claim 1 shows that the map $w^{\left\langle j\right\rangle }$ is well-defined
and equals $\varphi_{X_{j}}\left(  w^{\left(  j\right)  }\right)  $. Hence,
$w^{\left\langle j\right\rangle }=\varphi_{X_{j}}\left(  w^{\left(  j\right)
}\right)  \in S_{X_{j}}$ (since $\varphi_{X_{j}}$ is a map from $S_{n,X_{j}}$
to $S_{X_{j}}$). In other words, $w^{\left\langle j\right\rangle }$ belongs to
the symmetric group $S_{X_{j}}$. The proof of Lemma \ref{lem.Xk.bij1}
\textbf{(a)} is thus complete. \medskip

\textbf{(b)} In the proof of Proposition \ref{prop.int.wXi=Xi} \textbf{(b)},
we defined a map%
\begin{align*}
\Phi:K  &  \rightarrow S_{n,X_{1}}\times S_{n,X_{2}}\times\cdots\times
S_{n,X_{k}},\\
w  &  \mapsto\left(  w^{\left(  1\right)  },w^{\left(  2\right)  }%
,\ldots,w^{\left(  k\right)  }\right)  .
\end{align*}
Consider this map $\Phi$. In the proof of Proposition \ref{prop.int.wXi=Xi}
\textbf{(b)}, we have shown that this map $\Phi$ is a bijection.

Each $j\in\left[  k\right]  $ and $w\in K$ satisfy%
\begin{equation}
\varphi_{X_{j}}\left(  w^{\left(  j\right)  }\right)  =w^{\left\langle
j\right\rangle } \label{pf.lem.Xk.bij1.5}%
\end{equation}
(since Claim 1 shows that $w^{\left\langle j\right\rangle }$ equals
$\varphi_{X_{j}}\left(  w^{\left(  j\right)  }\right)  $).

The map $\varphi_{\Pi}\circ\Phi:K\rightarrow S_{X_{1}}\times S_{X_{2}}%
\times\cdots\times S_{X_{k}}$ is a bijection (since it is the composition of
the two bijections $\varphi_{\Pi}$ and $\Phi$). But each $w\in K$ satisfies%
\begin{align*}
\left(  \varphi_{\Pi}\circ\Phi\right)  \left(  w\right)   &  =\varphi_{\Pi
}\left(  \Phi\left(  w\right)  \right)  =\varphi_{\Pi}\left(  w^{\left(
1\right)  },w^{\left(  2\right)  },\ldots,w^{\left(  k\right)  }\right) \\
&  \ \ \ \ \ \ \ \ \ \ \ \ \ \ \ \ \ \ \ \ \left(
\begin{array}
[c]{c}%
\text{since the definition of }\Phi\\
\text{yields }\Phi\left(  w\right)  =\left(  w^{\left(  1\right)  },w^{\left(
2\right)  },\ldots,w^{\left(  k\right)  }\right)
\end{array}
\right) \\
&  =\left(  \varphi_{X_{1}}\times\varphi_{X_{2}}\times\cdots\times
\varphi_{X_{k}}\right)  \left(  w^{\left(  1\right)  },w^{\left(  2\right)
},\ldots,w^{\left(  k\right)  }\right) \\
&  \ \ \ \ \ \ \ \ \ \ \ \ \ \ \ \ \ \ \ \ \left(  \text{since }\varphi_{\Pi
}\text{ was defined to be }\varphi_{X_{1}}\times\varphi_{X_{2}}\times
\cdots\times\varphi_{X_{k}}\right) \\
&  =\left(  \varphi_{X_{1}}\left(  w^{\left(  1\right)  }\right)
,\ \varphi_{X_{2}}\left(  w^{\left(  2\right)  }\right)  ,\ \ldots
,\ \varphi_{X_{k}}\left(  w^{\left(  k\right)  }\right)  \right) \\
&  \ \ \ \ \ \ \ \ \ \ \ \ \ \ \ \ \ \ \ \ \left(  \text{by the definition of
}\varphi_{X_{1}}\times\varphi_{X_{2}}\times\cdots\times\varphi_{X_{k}}\right)
\\
&  =\left(  w^{\left\langle 1\right\rangle },w^{\left\langle 2\right\rangle
},\ldots,w^{\left\langle k\right\rangle }\right) \\
&  \ \ \ \ \ \ \ \ \ \ \ \ \ \ \ \ \ \ \ \ \left(  \text{since
(\ref{pf.lem.Xk.bij1.5}) says that }\varphi_{X_{j}}\left(  w^{\left(
j\right)  }\right)  =w^{\left\langle j\right\rangle }\text{ for each }%
j\in\left[  k\right]  \right)  .
\end{align*}
Hence, the map $\varphi_{\Pi}\circ\Phi$ can be rewritten as the map%
\begin{align*}
K  &  \rightarrow S_{X_{1}}\times S_{X_{2}}\times\cdots\times S_{X_{k}},\\
w  &  \mapsto\left(  w^{\left\langle 1\right\rangle },w^{\left\langle
2\right\rangle },\ldots,w^{\left\langle k\right\rangle }\right)  .
\end{align*}
Thus, the map%
\begin{align*}
K  &  \rightarrow S_{X_{1}}\times S_{X_{2}}\times\cdots\times S_{X_{k}},\\
w  &  \mapsto\left(  w^{\left\langle 1\right\rangle },w^{\left\langle
2\right\rangle },\ldots,w^{\left\langle k\right\rangle }\right)
\end{align*}
is a bijection (since $\varphi_{\Pi}\circ\Phi$ is a bijection). This proves
Lemma \ref{lem.Xk.bij1} \textbf{(b)}. \medskip

\textbf{(c)} Let $w\in K$. For each subset $X$ of $\left[  n\right]  $, recall
the default embedding $I_{X\rightarrow n}:S_{X}\rightarrow S_{n}$ defined in
Definition \ref{def.default-ex.IXn}. We claim the following:

\begin{statement}
\textit{Claim 2:} Let $j\in\left[  k\right]  $. Then, $I_{X_{j}\rightarrow
n}\left(  w^{\left\langle j\right\rangle }\right)  =w^{\left(  j\right)  }$.
\end{statement}

\begin{proof}
[Proof of Claim 2.]The definition of $I_{X_{j}\rightarrow n}$ yields
$I_{X_{j}\rightarrow n}\left(  w^{\left\langle j\right\rangle }\right)
=w^{\left\langle j\right\rangle }\uparrow^{n}$. (See Definition
\ref{def.default-ex.IXn} for the meaning of the \textquotedblleft$\uparrow
$\textquotedblright\ symbol.) The definition of $w^{\left\langle
j\right\rangle }\uparrow^{n}$ (see Definition \ref{def.default-ex.IXn}) shows
that%
\begin{equation}
\left(  w^{\left\langle j\right\rangle }\uparrow^{n}\right)  \left(  i\right)
=%
\begin{cases}
w^{\left\langle j\right\rangle }\left(  i\right)  , & \text{if }i\in X_{j};\\
i, & \text{if }i\notin X_{j}%
\end{cases}
\label{pf.lem.Xk.bij1.c.1}%
\end{equation}
for each $i\in\left[  n\right]  $. Thus, for each $x\in\left[  n\right]  $, we
have%
\begin{align*}
\left(  w^{\left\langle j\right\rangle }\uparrow^{n}\right)  \left(  x\right)
&  =%
\begin{cases}
w^{\left\langle j\right\rangle }\left(  x\right)  , & \text{if }x\in X_{j};\\
x, & \text{if }x\notin X_{j}%
\end{cases}
\ \ \ \ \ \ \ \ \ \ \left(  \text{by (\ref{pf.lem.Xk.bij1.c.1}), applied to
}i=x\right) \\
&  =%
\begin{cases}
w\left(  x\right)  , & \text{if }x\in X_{j};\\
x, & \text{if }x\notin X_{j}%
\end{cases}
\ \ \ \ \ \ \ \ \ \ \left(
\begin{array}
[c]{c}%
\text{since }w^{\left\langle j\right\rangle }\left(  x\right)  =w\left(
x\right)  \text{ if }x\in X_{j}\\
\text{(by the definition of }w^{\left\langle j\right\rangle }\text{)}%
\end{array}
\right) \\
&  =w^{\left(  j\right)  }\left(  x\right)
\end{align*}
(since the definition of $w^{\left(  j\right)  }$ says that $w^{\left(
j\right)  }\left(  x\right)  =%
\begin{cases}
w\left(  x\right)  , & \text{if }x\in X_{j};\\
x, & \text{if }x\notin X_{j}%
\end{cases}
\ $). Therefore, $\left.  w^{\left\langle j\right\rangle }\uparrow^{n}\right.
=w^{\left(  j\right)  }$. Thus, $I_{X_{j}\rightarrow n}\left(  w^{\left\langle
j\right\rangle }\right)  =\left.  w^{\left\langle j\right\rangle }\uparrow
^{n}\right.  =w^{\left(  j\right)  }$. This proves Claim 2.
\end{proof}

\begin{statement}
\textit{Claim 3:} Let $j\in\left[  k\right]  $. Then, $\left(  -1\right)
^{w^{\left(  j\right)  }}=\left(  -1\right)  ^{w^{\left\langle j\right\rangle
}}$.
\end{statement}

\begin{proof}
[Proof of Claim 3.]Proposition \ref{prop.default-ex.mor} \textbf{(d)} (applied
to $X=X_{j}$) shows that the default embedding $I_{X_{j}\rightarrow n}$
preserves signs; i.e., that we have
\[
\left(  -1\right)  ^{\sigma}=\left(  -1\right)  ^{I_{X_{j}\rightarrow
n}\left(  \sigma\right)  }\ \ \ \ \ \ \ \ \ \ \text{for any }\sigma\in
S_{X_{j}}.
\]
Applying the latter equality to $\sigma=w^{\left\langle j\right\rangle }$, we
obtain
\[
\left(  -1\right)  ^{w^{\left\langle j\right\rangle }}=\left(  -1\right)
^{I_{X_{j}\rightarrow n}\left(  w^{\left\langle j\right\rangle }\right)
}=\left(  -1\right)  ^{w^{\left(  j\right)  }}%
\]
(since Claim 2 yields $I_{X_{j}\rightarrow n}\left(  w^{\left\langle
j\right\rangle }\right)  =w^{\left(  j\right)  }$). This proves Claim 3.
\end{proof}

Now, Claim 3 from the proof of Proposition \ref{prop.int.wXi=Xi} \textbf{(b)}
shows that $w=w^{\left(  1\right)  }w^{\left(  2\right)  }\cdots w^{\left(
k\right)  }$. Hence,%
\begin{align*}
\left(  -1\right)  ^{w}  &  =\left(  -1\right)  ^{w^{\left(  1\right)
}w^{\left(  2\right)  }\cdots w^{\left(  k\right)  }}\\
&  =\underbrace{\left(  -1\right)  ^{w^{\left(  1\right)  }}}%
_{\substack{=\left(  -1\right)  ^{w^{\left\langle 1\right\rangle }}\\\text{(by
Claim 3)}}}\underbrace{\left(  -1\right)  ^{w^{\left(  2\right)  }}%
}_{\substack{=\left(  -1\right)  ^{w^{\left\langle 2\right\rangle }%
}\\\text{(by Claim 3)}}}\cdots\underbrace{\left(  -1\right)  ^{w^{\left(
k\right)  }}}_{\substack{=\left(  -1\right)  ^{w^{\left\langle k\right\rangle
}}\\\text{(by Claim 3)}}}\ \ \ \ \ \ \ \ \ \ \left(  \text{by the
multiplicativity of the sign}\right) \\
&  =\left(  -1\right)  ^{w^{\left\langle 1\right\rangle }}\left(  -1\right)
^{w^{\left\langle 2\right\rangle }}\cdots\left(  -1\right)  ^{w^{\left\langle
k\right\rangle }}.
\end{align*}
This proves Lemma \ref{lem.Xk.bij1} \textbf{(c)}.
\end{proof}
\end{fineprint}

\begin{exercise}
Let $n$ be a positive integer. Consider the map%
\begin{align*}
\rho_{n}:S_{n}  &  \rightarrow S_{n-1},\\
w  &  \mapsto I_{n-1\rightarrow n}^{-1}\left(  t_{n,w\left(  n\right)
}w\right)
\end{align*}
(where $t_{n,n}$ is understood to mean $\operatorname*{id}$). \medskip

\textbf{(a)} \fbox{1} Show that this is well-defined (i.e., that
$t_{n,w\left(  n\right)  }w$ really lies in the image of the default embedding
$I_{n-1\rightarrow n}$ for each $w\in S_{n}$). \medskip

\textbf{(b)} \fbox{1} Describe what $\rho_{n}$ does to the disjoint cycle
decomposition of a permutation $w$. \medskip

\textbf{(c)} \fbox{1} Show that $\rho_{n}$ is \textbf{not} a group morphism
whenever $n\geq3$. \medskip

\textbf{(d)} \fbox{1} Show that nevertheless, we have $\rho_{n}\left(
uv\right)  =\rho_{n}\left(  u\right)  \rho_{n}\left(  v\right)  $ and
$\rho_{n}\left(  vu\right)  =\rho_{n}\left(  v\right)  \rho_{n}\left(
u\right)  $ whenever $u\in S_{n}$ and $v\in I_{n-1\rightarrow n}\left(
S_{n-1}\right)  $.
\end{exercise}

\subsubsection{Base change}

So far, we have constructed morphisms between monoid algebras over the same
base ring $\mathbf{k}$. But there are also ring morphisms between monoid
algebras over different base rings. The simplest ones are constructed as follows:

\begin{proposition}
\label{prop.monalg.base-change}Let $f:\mathbf{k}\rightarrow\mathbf{l}$ be a
ring morphism between two commutative rings $\mathbf{k}$ and $\mathbf{l}$. Let
$M$ be any monoid. Then, the map%
\begin{align*}
f_{\ast}:\mathbf{k}\left[  M\right]   &  \rightarrow\mathbf{l}\left[
M\right]  ,\\
\left(  \alpha_{m}\right)  _{m\in M}  &  \mapsto\left(  f\left(  \alpha
_{m}\right)  \right)  _{m\in M}%
\end{align*}
(sending each $\sum_{m\in M}\alpha_{m}m$ to $\sum_{m\in M}f\left(  \alpha
_{m}\right)  m$, where $\alpha_{m}\in\mathbf{k}$ are scalars) is a ring
morphism. It is furthermore a $\mathbf{k}$-algebra morphism if we equip
$\mathbf{l}\left[  M\right]  $ with the $\mathbf{k}$-algebra structure that
comes from restriction of scalars\footnotemark\ via $f:\mathbf{k}%
\rightarrow\mathbf{l}$.
\end{proposition}

\footnotetext{A reminder from abstract algebra: \emph{Restriction of scalars}
is the easiest way to turn an $\mathbf{l}$-module or $\mathbf{l}$-algebra into
an $\mathbf{k}$-module or $\mathbf{k}$-algebra (respectively) using a ring
morphism $f:\mathbf{k}\rightarrow\mathbf{l}$. Namely, if $V$ is an
$\mathbf{l}$-module, then we make $V$ into a $\mathbf{k}$-module by setting%
\[
\kappa v:=f\left(  \kappa\right)  v\ \ \ \ \ \ \ \ \ \ \text{for all }%
\kappa\in\mathbf{k}\text{ and }v\in V.
\]
(This defines a scaling; an addition is already defined on $V$.) By the same
rule, we can make $V$ into a $\mathbf{k}$-algebra if $V$ is an $\mathbf{l}%
$-algebra. This is what we do here (for $V=\mathbf{l}\left[  M\right]  $).}

\begin{proof}
[Proof sketch.]This is just a matter of juggling definitions around and using
linearity, but let me outline the process for the sake of completeness:

It suffices to prove the second claim (since a $\mathbf{k}$-algebra morphism
is automatically a ring morphism).

First, it is easy to see that the map $f_{\ast}$ is $\mathbf{k}$-linear (since
the map $f:\mathbf{k}\rightarrow\mathbf{l}$ itself is $\mathbf{k}$-linear if
we equip $\mathbf{l}$ with the $\mathbf{k}$-algebra structure coming from
restriction of scalars). Moreover, it is easy to see that $f_{\ast}$ sends
each standard basis vector $e_{m}$ of $\mathbf{k}\left[  M\right]  $ to the
corresponding standard basis vector $e_{m}$ of $\mathbf{l}\left[  M\right]  $:
That is, we have%
\begin{equation}
f_{\ast}\left(  e_{m}\right)  =e_{m}\ \ \ \ \ \ \ \ \ \ \text{for each }m\in
M. \label{pf.prop.monalg.base-change.em}%
\end{equation}

\begin{proof}
[Proof of (\ref{pf.prop.monalg.base-change.em}).]Let $m\in M$. Recall that $f$
is a ring homomorphism, and thus satisfies $f\left(  1_{\mathbf{k}}\right)
=1_{\mathbf{l}}$ and $f\left(  0_{\mathbf{k}}\right)  =0_{\mathbf{l}}$.

By its definition, the map $f_{\ast}$ acts on a family $\left(  \alpha
_{m}\right)  _{m\in M}\in\mathbf{k}\left[  M\right]  $ by applying $f$ to each
entry of this family. But the standard basis vector $e_{m}$ of $\mathbf{k}%
\left[  M\right]  $ is defined to be the family whose $m$-th entry is
$1_{\mathbf{k}}$ and whose all other entries are $0_{\mathbf{k}}$. Thus, after
we apply $f_{\ast}$ to this vector $e_{m}$, the $m$-th entry becomes $f\left(
1_{\mathbf{k}}\right)  =1_{\mathbf{l}}$, while all the other entries become
$f\left(  0_{\mathbf{k}}\right)  =0_{\mathbf{l}}$. Hence, we obtain the family
whose $m$-th entry is $1_{\mathbf{l}}$ and whose all other entries are
$0_{\mathbf{l}}$. But this family is precisely the standard basis vector
$e_{m}$ of $\mathbf{l}\left[  M\right]  $. So we have shown that $f_{\ast}$
sends the standard basis vector $e_{m}$ of $\mathbf{k}\left[  M\right]  $ to
the corresponding standard basis vector $e_{m}$ of $\mathbf{l}\left[
M\right]  $. This proves (\ref{pf.prop.monalg.base-change.em}).
\end{proof}

In particular, (\ref{pf.prop.monalg.base-change.em}) yields $f_{\ast}\left(
e_{1}\right)  =e_{1}$ (where $1$ is the neutral element of $M$). In other
words, $f_{\ast}\left(  1_{\mathbf{k}\left[  M\right]  }\right)
=1_{\mathbf{l}\left[  M\right]  }$ (since $1_{\mathbf{k}\left[  M\right]
}=e_{1}$ and $1_{\mathbf{l}\left[  M\right]  }=e_{1}$). In other words, the
map $f_{\ast}$ respects the unity.

Next, we shall prove that $f_{\ast}$ respects multiplication. In other words,
we shall prove that $f_{\ast}\left(  \mathbf{ab}\right)  =f_{\ast}\left(
\mathbf{a}\right)  \cdot f_{\ast}\left(  \mathbf{b}\right)  $ for any
$\mathbf{a},\mathbf{b}\in\mathbf{k}\left[  M\right]  $. Since both sides of
this equality are $\mathbf{k}$-linear in $\mathbf{a}$ and $\mathbf{b}$, we can
reduce this to the case when both $\mathbf{a}$ and $\mathbf{b}$ are standard
basis vectors (just like we did in the proof of Theorem \ref{thm.Tsign.auto}
\textbf{(a)}). In other words, it suffices to prove that $f_{\ast}\left(
e_{m}e_{n}\right)  =f_{\ast}\left(  e_{m}\right)  \cdot f_{\ast}\left(
e_{n}\right)  $ for any $m,n\in M$. But this is easy: For any $m,n\in M$, we
have $e_{m}e_{n}=e_{mn}$ and thus%
\begin{align*}
f_{\ast}\left(  e_{m}e_{n}\right)   &  =f_{\ast}\left(  e_{mn}\right)
=e_{mn}\ \ \ \ \ \ \ \ \ \ \left(  \text{by
(\ref{pf.prop.monalg.base-change.em})}\right) \\
&  =\underbrace{e_{m}}_{\substack{=f_{\ast}\left(  e_{m}\right)  \\\text{(by
(\ref{pf.prop.monalg.base-change.em}))}}}\cdot\underbrace{e_{n}}%
_{\substack{=f_{\ast}\left(  e_{n}\right)  \\\text{(by
(\ref{pf.prop.monalg.base-change.em}))}}}=f_{\ast}\left(  e_{m}\right)  \cdot
f_{\ast}\left(  e_{n}\right)  ,
\end{align*}
qed. So we have shown that the map $f_{\ast}$ respects multiplication. Since
$f_{\ast}$ is $\mathbf{k}$-linear and respects the unity, we conclude that
$f_{\ast}$ is a $\mathbf{k}$-algebra morphism, and therefore a ring morphism.
Proposition \ref{prop.monalg.base-change} is thus proved.
\end{proof}

\begin{definition}
\label{def.monalg.base-change}The ring morphism $f_{\ast}$ constructed in
Proposition \ref{prop.monalg.base-change} is called the \emph{base change
morphism} for $M$ from $\mathbf{k}$ to $\mathbf{l}$.
\end{definition}

\begin{example}
Let $M$ be the monoid $\left\{  x^{0},x^{1},x^{2},\ldots\right\}  $ with
multiplication given by $x^{i}x^{j}=x^{i+j}$ and with neutral element $x^{0}$.
As we have seen in Subsection \ref{subsec.intro.monalg.exas}, its monoid
algebra $\mathbf{k}\left[  M\right]  $ is just the polynomial ring
$\mathbf{k}\left[  x\right]  $.

Let $f:\mathbb{C}\rightarrow\mathbb{C}$ be complex conjugation (i.e., the map
$\mathbb{C}\rightarrow\mathbb{C},\ z\mapsto\overline{z}$); this is a ring morphism.

Now, the base change morphism $f_{\ast}:\mathbb{C}\left[  x\right]
\rightarrow\mathbb{C}\left[  x\right]  $ defined in Proposition
\ref{prop.monalg.base-change} is the automorphism of the polynomial ring
$\mathbb{C}\left[  x\right]  $ that conjugates every coefficient of the polynomial.
\end{example}

In particular, we can apply Definition \ref{def.monalg.base-change} to
$M=S_{n}$, obtaining a base change morphism $f_{\ast}:\mathbf{k}\left[
S_{n}\right]  \rightarrow\mathbf{l}\left[  S_{n}\right]  $. \medskip

For any commutative ring $\mathbf{k}$, there is a unique ring morphism
$g:\mathbb{Z}\rightarrow\mathbf{k}$ (sending $1_{\mathbb{Z}}$ to
$1_{\mathbf{k}}$, and thus sending each integer $m\in\mathbb{Z}$ to
$m\cdot1_{\mathbf{k}}$). Thus, by Proposition \ref{prop.monalg.base-change},
this ring morphism yields a ring morphism%
\[
g_{\ast}:\mathbb{Z}\left[  S_{n}\right]  \rightarrow\mathbf{k}\left[
S_{n}\right]  ,
\]
which sends every $\sum_{w\in S_{n}}\alpha_{w}w$ (with $\alpha_{w}%
\in\mathbb{Z}$) to $\sum_{w\in S_{n}}\left(  \alpha_{w}\cdot1_{\mathbf{k}%
}\right)  w$. This morphism $g_{\ast}$ (and, more generally, any base change
morphism for $S_{n}$) sends every standard basis vector $w\in S_{n}$ of
$\mathbb{Z}\left[  S_{n}\right]  $ to the corresponding standard basis vector
$w$ of $\mathbf{k}\left[  S_{n}\right]  $. Thus, it sends $\nabla$ to $\nabla
$, sends each jucys--murphy $\mathbf{m}_{k}$ to $\mathbf{m}_{k}$, sends every
$\mathbf{t}_{k}$ to $\mathbf{t}_{k}$, and so on.

This has a nice application: If you want to prove an equality in
$\mathbf{k}\left[  S_{n}\right]  $, such as $\mathbf{m}_{2}\mathbf{m}%
_{5}=\mathbf{m}_{5}\mathbf{m}_{2}$, then it suffices to prove this equality in
$\mathbb{Z}\left[  S_{n}\right]  $ (that is, it suffices to prove it for
$\mathbf{k}=\mathbb{Z}$), because you can then apply the ring morphism
$g_{\ast}$ to both sides and obtain the corresponding equality in
$\mathbf{k}\left[  S_{n}\right]  $. This is a useful strategy for proving
equalities (usually more complicated ones than $\mathbf{m}_{2}\mathbf{m}%
_{5}=\mathbf{m}_{5}\mathbf{m}_{2}$), since $\mathbb{Z}\left[  S_{n}\right]  $
is a better-behaved ring than the generic $\mathbf{k}\left[  S_{n}\right]  $.
This strategy is known as \emph{base change}, and the slogan for it is
\textquotedblleft working over $\mathbb{Z}$ is enough\textquotedblright. (Of
course, this holds only if the equality you are proving can be stated in
$\mathbb{Z}$, meaning that it contains no fractions or irrationals or other
non-integer scalars. But it can be easily adapted to other situations. For
instance, if your equality involves fractions and is stated for $\mathbb{Q}%
$-algebras $\mathbf{k}$, then you only need to prove it in $\mathbb{Q}\left[
S_{n}\right]  $. If your equality involves an arbitrary (unspecified) element
$x\in\mathbf{k}$ (for instance, this applies to Corollary \ref{cor.YJM.1+xm}),
then you only need to prove it for $x=X$ in the ring $\left(  \mathbb{Z}%
\left[  X\right]  \right)  \left[  S_{n}\right]  $ (where $\mathbb{Z}\left[
X\right]  $ is the usual polynomial ring in one indeterminate $X$ over
$\mathbb{Z}$), because then you can apply $g_{\ast}$ where $g:\mathbb{Z}%
\left[  X\right]  \rightarrow\mathbf{k}$ is the unique ring morphism that
sends $X$ to $x$.)

\begin{remark}
\label{rmk.monalg.base-change.incluses}A particularly simple kind of ring
morphisms are the \emph{canonical inclusions} (i.e., the obvious maps from a
subring to a ring that send each element of the subring to itself). If
$\mathbf{k}$ is a subring of a commutative ring $\mathbf{l}$, then the
canonical inclusion $g:\mathbf{k}\rightarrow\mathbf{l}$ is a ring morphism,
and thus yields a base change morphism $g_{\ast}:\mathbf{k}\left[
S_{n}\right]  \rightarrow\mathbf{l}\left[  S_{n}\right]  $. This morphism
$g_{\ast}$ sends each $\sum_{w\in S_{n}}\alpha_{w}w$ to $\sum_{w\in S_{n}%
}\alpha_{w}w$, and thus we can treat it as being a canonical inclusion itself,
so that $\mathbf{k}\left[  S_{n}\right]  $ becomes a subring of $\mathbf{l}%
\left[  S_{n}\right]  $. In particular,%
\[
\mathbb{Z}\left[  S_{n}\right]  \subseteq\mathbb{Q}\left[  S_{n}\right]
\subseteq\mathbb{R}\left[  S_{n}\right]  \subseteq\mathbb{C}\left[
S_{n}\right]
\]
(and all these \textquotedblleft$\subseteq$\textquotedblright\ signs are
\textquotedblleft subring\textquotedblright\ signs).
\end{remark}

\subsection{Antipodal conjugacy}

Let us briefly return to the combinatorics of permutations:

\begin{theorem}
\label{thm.antip-conj.perm}Each permutation $w\in S_{n}$ is conjugate to its
inverse $w^{-1}$ in the group $S_{n}$.
\end{theorem}

\begin{proof}
Let $w\in S_{n}$ be a permutation. Then, Proposition \ref{prop.dcd.inverse}
\textbf{(b)} shows that the orbits of $w^{-1}$ are the orbits of $w$. In other
words, the permutations $w$ and $w^{-1}$ have the same orbits, and therefore
the same cycle type (since the cycle type is defined in terms of the orbits).
But any two permutations having the same cycle type are conjugate (by Theorem
\ref{thm.type.conj=type}). Thus, the permutations $w$ and $w^{-1}$ are
conjugate. This proves Theorem \ref{thm.antip-conj.perm}.
\end{proof}

Note that Theorem \ref{thm.antip-conj.perm} is specific to symmetric groups;
it is false (e.g.) for all cyclic groups of order $>2$, and for many other groups.

\begin{exercise}
\label{exe.antip-conj.perm.invol}\fbox{2} Prove the following stronger version
of Theorem \ref{thm.antip-conj.perm}: Let $w\in S_{n}$. Then, there exists an
involution $x\in S_{n}$ such that $w=xw^{-1}x^{-1}$ (or, equivalently,
$w=xw^{-1}x$).
\end{exercise}

Theorem \ref{thm.antip-conj.perm} has a linear analogue, which is surprisingly nontrivial:

\begin{theorem}
\label{thm.antip-conj.linear}Assume that $\mathbf{k}$ is a field of
characteristic $0$ (such as $\mathbb{Q}$, $\mathbb{R}$ or $\mathbb{C}$). Let
$\mathbf{a}\in\mathbf{k}\left[  S_{n}\right]  $. Then, $\mathbf{a}$ is
conjugate to $S\left(  \mathbf{a}\right)  $ in $\mathbf{k}\left[
S_{n}\right]  $ (meaning that there exists some invertible element
$\mathbf{u}\in\mathbf{k}\left[  S_{n}\right]  $ such that $\mathbf{a}%
=\mathbf{u}S\left(  \mathbf{a}\right)  \mathbf{u}^{-1}$).
\end{theorem}

I will prove this theorem later (in Section \ref{sec.rep.antip-conj}). Note
that the \textquotedblleft characteristic $0$\textquotedblright\ requirement
is really needed here. \medskip

The following corollary of Theorem \ref{thm.antip-conj.perm} will be used
later on:

\begin{proposition}
\label{prop.ZkSn.S=id}We have $S\left(  \mathbf{a}\right)  =\mathbf{a}$ for
all $\mathbf{a}\in Z\left(  \mathbf{k}\left[  S_{n}\right]  \right)  $.
\end{proposition}

\begin{proof}
Theorem \ref{thm.center.conjsums} (applied to $G=S_{n}$) shows that the center
$Z\left(  \mathbf{k}\left[  S_{n}\right]  \right)  $ of the group algebra
$\mathbf{k}\left[  S_{n}\right]  $ is the $\mathbf{k}$-linear span of the
conjugacy class sums. In other words, the conjugacy class sums span $Z\left(
\mathbf{k}\left[  S_{n}\right]  \right)  $.

Let $\mathbf{a}\in Z\left(  \mathbf{k}\left[  S_{n}\right]  \right)  $. We
must prove the equality $S\left(  \mathbf{a}\right)  =\mathbf{a}$. Both sides
of this equality are $\mathbf{k}$-linear in $\mathbf{a}$. Hence, in proving
this equality, we can WLOG assume that $\mathbf{a}$ is a conjugacy class sum
(since the conjugacy class sums span $Z\left(  \mathbf{k}\left[  S_{n}\right]
\right)  $). Assume this. Thus, $\mathbf{a}=\mathbf{z}_{C}$ for some conjugacy
class $C$ of $S_{n}$ (where $\mathbf{z}_{C}$ is as defined in Definition
\ref{def.groups.conj} \textbf{(c)}). Consider this $C$. Thus, $\mathbf{z}%
_{C}=\sum_{c\in C}c$ (by the definition of $\mathbf{z}_{C}$).

However, Theorem \ref{thm.antip-conj.perm} shows that each $c\in C$ satisfies
$c^{-1}\in C$\ \ \ \ \footnote{\textit{Proof.} Let $c\in C$. Then, $c\in
C\subseteq S_{n}$. Hence, Theorem \ref{thm.antip-conj.perm} (applied to $w=c$)
shows that the permutation $c$ is conjugate to its inverse $c^{-1}$ in the
group $S_{n}$. In other words, $c$ and $c^{-1}$ belong to the same conjugacy
class of $S_{n}$. This conjugacy class must be $C$ (since $c\in C$). Thus, we
conclude that $c^{-1}$ belongs to $C$ as well. In other words, $c^{-1}\in C$,
qed.}. Thus, the map%
\begin{align*}
C  &  \rightarrow C,\\
c  &  \mapsto c^{-1}%
\end{align*}
is well-defined. This map is clearly inverse to itself (since $\left(
c^{-1}\right)  ^{-1}=c$ for each $c\in C$), and thus is invertible, i.e., a bijection.

Applying the map $S$ to the equality $\mathbf{z}_{C}=\sum_{c\in C}c$, we
obtain%
\begin{align*}
S\left(  \mathbf{z}_{C}\right)   &  =S\left(  \sum_{c\in C}c\right)
=\sum_{c\in C}\underbrace{S\left(  c\right)  }_{\substack{=c^{-1}\\\text{(by
the definition of }S\text{)}}}\ \ \ \ \ \ \ \ \ \ \left(  \text{since the map
}S\text{ is }\mathbf{k}\text{-linear}\right) \\
&  =\sum_{c\in C}c^{-1}=\sum_{c\in C}c\ \ \ \ \ \ \ \ \ \ \left(
\begin{array}
[c]{c}%
\text{here, we have substituted }c\text{ for }c^{-1}\text{ in}\\
\text{the sum, since the map }C\rightarrow C,\ c\mapsto c^{-1}\\
\text{is a bijection}%
\end{array}
\right) \\
&  =\mathbf{z}_{C}.
\end{align*}
In other words, $S\left(  \mathbf{a}\right)  =\mathbf{a}$ (since
$\mathbf{a}=\mathbf{z}_{C}$). Proposition \ref{prop.ZkSn.S=id} is thus proved.
\end{proof}

Another proof of Proposition \ref{prop.ZkSn.S=id} can be obtained by combining
Exercise \ref{exe.S.GZ} with (\ref{eq.thm.center.murphies.incs}%
).\footnote{Also, Theorem \ref{thm.antip-conj.linear} yields Proposition
\ref{prop.ZkSn.S=id} when $\mathbf{k}$ is a field of characteristic $0$.}

\subsection{Computer algebra}

\subsubsection{SageMath}

A symmetric group algebra $\mathbf{k}\left[  S_{n}\right]  $ is a rather
tangible object, perfectly suited for computation and experiments. However,
computing inside $\mathbf{k}\left[  S_{n}\right]  $ by hand quickly gets
tiresome as $n$ gets large, since $\left\vert S_{n}\right\vert =n!$ grows very
fast. Thus, it is advantageous to use computer algebra when working in
$\mathbf{k}\left[  S_{n}\right]  $.

Alas, most CASs (computer algebra systems) do not know $\mathbf{k}\left[
S_{n}\right]  $. You can work around this by representing elements of
$\mathbf{k}\left[  S_{n}\right]  $ as matrices (we will see some ways to do
this in the next two chapters), or you can write your own code for
programmable CASs like Mathematica or Maple. In this section, I shall instead
show how to use SageMath: one of the few CASs that \textquotedblleft
natively\textquotedblright\ know $\mathbf{k}\left[  S_{n}\right]  $.

\href{https://www.sagemath.org/}{\emph{SageMath} (for short: \emph{Sage})} is
a CAS built atop the Python programming language; you can view it as a Python
dialect with a huge mathematical library attached to it (particularly
well-suited for algebra and combinatorics). If you are familiar with Python,
you will find SageMath fairly intuitive (although it redefines some Python
syntax).\footnote{Another CAS that knows group rings is \emph{GAP} (see
\url{https://www.gap-system.org/} ).}

You can install SageMath on your computer (it runs natively on Linux and Mac;
on Windows it runs via
\href{https://learn.microsoft.com/en-us/windows/wsl/install}{the WSL
subsystem}). But the easiest way to use SageMath is through the Sage cell
server:%
\[
\text{\text{\texttt{\href{https://sagecell.sagemath.org/}{\text{\texttt{https://sagecell.sagemath.org/}%
}}}}}%
\]
This is an anonymous cloud server: You type in your code and it runs the
computation for you (as long as it takes no longer than 30 seconds).

\subsubsection{Permutations in SageMath}

Let us see how SageMath is used. First, we play around with permutations in
the symmetric groups $S_{n}$; the group algebras will come later.

I highly recommend pasting (or typing) all the examples into the Sage cell
server (
\text{\texttt{\href{https://sagecell.sagemath.org/}{\text{\texttt{https://sagecell.sagemath.org/}%
}}} ) as you read this section.}

\paragraph{Entering permutations.}

Permutations in Sage are usually expressed in OLN. For example, the
permutation $\operatorname*{oln}\left(  512364\right)  $ can be entered as
follows:\footnote{Everything printed in \texttt{typewritten font} in this
section is either an input for Sage or an output of Sage.}

\qquad\texttt{w = Permutation([5, 1, 2, 3, 6, 4])}

Sage can then compute values of \texttt{w} in the obvious way: Typing in
\texttt{w(1)} will yield the output \texttt{5}, whereas typing in
\texttt{w(2)} will yield \texttt{1}, and so on. \medskip

There are other ways to construct permutations in SageMath. For example,
permutations can be entered in cycle notation (i.e., via their DCD) as follows:

\qquad\texttt{u = Permutation(((1,3),(2,4,6),(5,)))}

This is the permutation with DCD $\left(  \left(  1,3\right)  ,\ \left(
2,4,6\right)  ,\ \left(  5\right)  \right)  $, that is, the permutation
$\operatorname*{cyc}\nolimits_{1,3}\operatorname*{cyc}\nolimits_{2,4,6}%
\operatorname*{cyc}\nolimits_{5}\in S_{6}$. (The extraneous-looking comma in
\textquotedblleft\texttt{(5,)}\textquotedblright\ is a programming quirk:
Without the comma, Sage would understand \textquotedblleft\texttt{(5)}%
\textquotedblright\ as meaning the \textbf{number }$5$, rather than the
$1$-tuple $\left(  5\right)  $ that we want to enter.)

You can skip length-$1$ cycles when entering a permutation in DCD, as long as
Sage understands correctly what the $n$ is. For example, the above permutation
\texttt{u} can also be entered as \texttt{Permutation(((1,3),(2,4,6)))}; Sage
will automatically understand that $5$ is meant to be a fixed
point.\footnote{Of course, if you want $\operatorname*{cyc}\nolimits_{1,2}\in
S_{3}$, then you need to write \texttt{Permutation(((1,2),(3,)))}, since
\texttt{Permutation(((1,2),))} will yield $\operatorname*{cyc}\nolimits_{1,2}%
\in S_{2}$. Although this is not that bad -- see the handling of default
extensions below.} \medskip

In case of doubt, always consult the documentation. For instance,
\href{https://doc.sagemath.org/html/en/reference/combinat/sage/combinat/permutation.html}{the
documentation for permutations in Sage} explains all the above and much more.
(On the cell server, you can also type \texttt{Permutation??} for a quick
overview of the class-level documentation, although this is often not sufficient.)

I also recommend that you acquaint yourself with some of Sage/Python's
syntactic quirks and pitfalls (e.g., its load-bearing
whitespace/indentation\footnote{See, e.g.,
\url{https://note.nkmk.me/en/python-indentation-usage/} .
\par
This has the unfortunate consequence that if you copypaste my sample code from
these notes, you will probably get the indentation wrong and you will have to
fix it manually.}; its use of square brackets for lists\footnote{So
\texttt{[4, 2, 5, 7]} in Sage/Python is what we would call $\left(
4,2,5,7\right)  $.}; its use of \texttt{=} for assignment and \texttt{==} for
equality\footnote{Thus, saying \texttt{a = 15} in Sage/Python means
\textquotedblleft let $a$ be $15$\textquotedblright, whereas saying \texttt{a
== 15} means \textquotedblleft$a$ equals $15$\textquotedblright.}).

\paragraph{Composing permutations.}

Let us now multiply (i.e., compose) the two permutations
$u=\operatorname*{oln}\left(  312\right)  $ and $v=\operatorname*{oln}\left(
132\right)  $ in Sage. The symbol for multiplication in Sage is \texttt{*}
(you cannot just write \texttt{uv} or \texttt{u v} and hope that Sage will
understand it as a product); so we try the following:

\qquad\texttt{u = Permutation([3,1,2])}

\qquad\texttt{v = Permutation([1,3,2])}

\qquad\texttt{u*v}

We get \texttt{[2, 1, 3]}, which (like any permutation that Sage outputs)
should be read in one-line notation, and thus means $\operatorname*{oln}%
\left(  213\right)  $. But wait -- the product $uv$ is actually
$\operatorname*{oln}\left(  321\right)  $, not $\operatorname*{oln}\left(
213\right)  $. What has gone wrong?

The answer is that we got bitten by Sage's nonstandard notations. Sage defines
the product of permutations in the \textbf{opposite} way from us; i.e., it
defines a product $uv$ of two permutations to be \textquotedblleft first $u$
then $v$\textquotedblright\ (so that $\left(  uv\right)  \left(  i\right)
=v\left(  u\left(  i\right)  \right)  $ for each $i$). Thus, what we call $uv$
is rather called $vu$ (or, more precisely, \texttt{v*u}) in Sage. Thus, if you
are using our convention for multiplying permutations, then you should write
all products in reverse when you feed them into Sage. This is an annoying
source of mistakes and anxiety which it is, alas, too late to
fix.\footnote{Sage's convention is in fact inspired by the GAP CAS, which was
written by group theorists (the \textquotedblleft G\textquotedblright\ in
\textquotedblleft GAP\textquotedblright\ stands for \textquotedblleft
groups\textquotedblright), who often write functions to the right of the
elements they act on: i.e., they write $x.f$ instead of $f\left(  x\right)  $.
From this point of view, the convention makes sense, because a product $uv$ of
two permutations should send any $i$ to $\left(  i.u\right)  .v$. But Sage
writes $f\left(  x\right)  $, not $x.f$, so this product sends any $i$ to
$v\left(  u\left(  i\right)  \right)  $, which looks awful.}

(If you don't want to write products from right to left, you can also get
\textbf{our} product $uv$ by \texttt{u.left\_action\_product(v)} .)

\paragraph{Other features of permutations.}

Sage can do a lot more with permutations. It can take their inverses:

\qquad\texttt{u.inverse()}

It can compute their signs:

\qquad\texttt{u.sign()} \medskip

It can even multiply two permutations in different symmetric groups:

\qquad\texttt{Permutation([2,1]) * Permutation([3,2,1])}

This may appear baroque, but there is nothing wild going on here; Sage is
simply replacing the first permutation $\operatorname*{oln}\left(  21\right)
$ with its default extension $I_{2\rightarrow3}\left(  \operatorname*{oln}%
\left(  21\right)  \right)  =\operatorname*{oln}\left(  213\right)  $ in order
to get it into the same group as the second permutation. In other words, Sage
is automatically using the default embeddings $I_{m\rightarrow n}%
:S_{m}\rightarrow S_{n}$ as \textquotedblleft coercions\textquotedblright%
\ (i.e., replacing each $\sigma\in S_{m}$ by $I_{m\rightarrow n}\left(
\sigma\right)  $ whenever necessary).

\paragraph{Listing permutations.}

Now let's list all the $24$ permutations of $\left[  4\right]  $:

\qquad\texttt{Permutations(4)}

Sage dutifully answers \textquotedblleft\texttt{Standard permutations of
4}\textquotedblright, which is correct but not very useful. So we need to be
explicit and ask for a list if we want a list:

\qquad\texttt{list(Permutations(4))}

Now Sage indeed produces a list of all permutations $w\in S_{4}$, in
lexicographically increasing order (i.e., sorted by $w\left(  1\right)  $,
then by $w\left(  2\right)  $, then by $w\left(  3\right)  $, then by
$w\left(  4\right)  $). \medskip

We can quickly access specific entries of this list using Python's
\textquotedblleft bracket\textquotedblright\ notation for the entries of a list:

\qquad\texttt{list(Permutations(4))[5]}

This gives us the $6$-th permutation in the list -- yes, the $6$-th, not the
$5$-th, because Sage (like Python) numbers things starting at $0$ rather than
$1$. Thus, the first permutation in the list can be obtained by

\qquad\texttt{list(Permutations(4))[0]}

There is also a shortcut to the last entry of a list:

\qquad\texttt{list(Permutations(4))[-1]}

Generally, if \texttt{a} is a list, then \texttt{a[-1]} means its last entry
in Python and in Sage. \medskip

By the way, the last permutation \texttt{list(Permutations(4))[-1]} is
$\operatorname*{oln}\left(  4321\right)  $, and this is not a coincidence: The
last permutation in $S_{n}$ (in lexicographically increasing order) is always
the permutation $\operatorname*{oln}\left(  n,n-1,\ldots,2,1\right)  $. In
Sage, you can also construct it as follows:

\qquad\texttt{Permutation(list(range(n, 0, -1)))}

How does this work? Python's \texttt{list(range(a, b, c))} creates an
arithmetic sequence starting at $a$ and increasing by a step size of $c$ until
reaching or surpassing $b$ (but $b$ itself is not included)\footnote{This
works even if $c$ is negative: For instance, \texttt{list(range(7, 1, -2))}
yields \texttt{[7, 5, 3]}.}. Thus, \texttt{list(range(n, 0, -1))} gives the
sequence $\left(  n,n-1,\ldots,2,1\right)  $. Finally, by calling
\texttt{Permutation} on this sequence, we transform it into a permutation.

Likewise, the identity permutation $\operatorname*{id}=\operatorname*{oln}%
\left(  1,2,\ldots,n\right)  \in S_{n}$ can be obtained by

\qquad\texttt{Permutation(list(range(1, n+1, 1)))}

Here you can leave out the \texttt{1} at the end (i.e., you can replace
\texttt{range(1, n+1, 1)} by \texttt{range(1, n+1)}), since $1$ is the default
value for the step size of an arithmetic sequence in Sage.

\paragraph{Permutations with a given property.}

Let us now list all involutions in $S_{4}$, that is, all permutations $w\in
S_{4}$ that satisfy $w^{-1}=w$. Sage has no inbuilt function for this (I
believe), but we can just use the definition:

\qquad\texttt{[w for w in Permutations(4) if w.inverse() == w]}

This syntax is fairly self-explanatory, once you recall that square brackets
in Sage stand for lists (so that we are asking for a list), and that the
\textquotedblleft double equality sign\textquotedblright\ \texttt{==} just
means regular equality (whereas the \textquotedblleft single equality
sign\textquotedblright\ \texttt{=} means assignment). Thus, the expression
above means \textquotedblleft the list of all $w$ where $w$ ranges over the
permutations of $\left[  4\right]  $ that satisfy $w^{-1}=w$\textquotedblright%
. And Sage responds with the actual list of these $w$'s. \medskip

If we want to know the number of these $w$'s (rather than a list of them),
then we can ask Sage for the length of this list:

\qquad\texttt{len([w for w in Permutations(4) if w.inverse() == w])}

Or, even better (certainly faster!), we can let Sage count them:

\qquad\texttt{sum(1 for w in Permutations(4) if w.inverse() == w)}

Yes, we have just asked Sage for the sum $\sum_{\substack{w\in S_{4}%
;\\w^{-1}=w}}1$. This is the most natural way to count these $w$'s, and unlike
the previous way, it does not make Sage waste its time by listing them first.
\medskip

Similar syntax can be used for other classes of permutations. For instance,
let's say you want to count the permutations $w\in S_{7}$ that satisfy
$w\left(  2\right)  >w\left(  3\right)  >\cdots>w\left(  7\right)  $. You can
do this via

\qquad\texttt{sum(1 for w in Permutations(7) if w(2)
%TCIMACRO{\TEXTsymbol{>} }%
%BeginExpansion
$>$
%EndExpansion
w(3)}

\qquad\qquad\qquad\qquad\qquad\texttt{
%TCIMACRO{\TEXTsymbol{>} }%
%BeginExpansion
$>$
%EndExpansion
w(4)
%TCIMACRO{\TEXTsymbol{>} }%
%BeginExpansion
$>$
%EndExpansion
w(5)
%TCIMACRO{\TEXTsymbol{>} }%
%BeginExpansion
$>$
%EndExpansion
w(6)
%TCIMACRO{\TEXTsymbol{>} }%
%BeginExpansion
$>$
%EndExpansion
w(7))}

A more conceptual way would be

\qquad\texttt{sum(1 for w in Permutations(7) if}

\qquad\qquad\qquad\qquad\qquad\texttt{all(w(i)
%TCIMACRO{\TEXTsymbol{>} }%
%BeginExpansion
$>$
%EndExpansion
w(i+1) for i in range(2, 7)))}

In fact, \texttt{all(w(i)
%TCIMACRO{\TEXTsymbol{>} }%
%BeginExpansion
$>$
%EndExpansion
w(i+1) for i in range(2, 7))} is Sage's (and Python's) syntax for
\textquotedblleft$w\left(  i\right)  >w\left(  i+1\right)  $ holds for all
$i\in\left\{  2,3,\ldots,6\right\}  $\textquotedblright, which is the
mathematically rigorous way of saying $w\left(  2\right)  >w\left(  3\right)
>\cdots>w\left(  7\right)  $.

\subsubsection{Symmetric group algebras in SageMath}

\paragraph{Invoking the symmetric group algebra.}

Now let us work in a symmetric group algebra. We will work in the $\mathbb{Q}%
$-algebra $\mathbb{Q}\left[  S_{n}\right]  $, since it includes $\mathbb{Z}%
\left[  S_{n}\right]  $ as a subring and contains almost everything we will
ever need.

Let us start with $\mathbb{Q}\left[  S_{3}\right]  $. Here is how to generate
this algebra in Sage:

\qquad\texttt{A = SymmetricGroupAlgebra(QQ, 3)}

Let us now compute the somewhere-to-below shuffle
\[
\mathbf{t}_{1}=\operatorname*{cyc}\nolimits_{1}+\operatorname*{cyc}%
\nolimits_{1,2}+\operatorname*{cyc}\nolimits_{1,2,3}=\operatorname*{oln}%
\left(  123\right)  +\operatorname*{oln}\left(  213\right)
+\operatorname*{oln}\left(  231\right)
\]
from Definition \ref{def.stb.tk}. We try the obvious thing:

\qquad\texttt{Permutation([1,2,3]) + Permutation([2,1,3]) +
Permutation([2,3,1])}

Sadly, this does not work. Indeed, it cannot work: The three permutations we
are trying to add are just permutations, but we want to add the corresponding
standard basis vectors, not the permutations themselves. By abuse of notation,
we might write the former while meaning the latter, but this does not make
them literally equal, and Sage cannot be fooled. Nor does Sage have the
context needed to understand our abuse correctly: Are we asking for the sum
$\operatorname*{oln}\left(  123\right)  +\operatorname*{oln}\left(
213\right)  +\operatorname*{oln}\left(  231\right)  $ in $\mathbb{Z}\left[
S_{3}\right]  $ or in $\mathbb{Q}\left[  S_{3}\right]  $ or in $\mathbf{k}%
\left[  S_{3}\right]  $ for some completely different ring $\mathbf{k}$ ?

So we need to ask Sage to add the standard basis vectors, not the
permutations. The surefire way to do this is by wrapping the permutations in
\texttt{A(...)}'s (remember that \texttt{A} is $\mathbb{Q}\left[
S_{3}\right]  $). So we write

\qquad\texttt{A(Permutation([1,2,3])) + A(Permutation([2,1,3]))}

\qquad\qquad\qquad\qquad\texttt{ + A(Permutation([2,3,1]))}

\noindent(this should be a single line in your Sage input; if you really want
to break it into two, then put a backslash \texttt{%
%TCIMACRO{\TEXTsymbol{\backslash}}%
%BeginExpansion
$\backslash$%
%EndExpansion
} at the end of the first line). Or, even simpler, we can just write

\qquad\texttt{A([1,2,3]) + A([2,1,3]) + A([2,3,1])}

Sage understands \texttt{A([...])} to mean \texttt{A(Permutation([...]))} .

\paragraph{Computing in the symmetric group algebra.}

Let us call this sum \texttt{t1} (as it is our $\mathbf{t}_{1}$):

\qquad\texttt{A = SymmetricGroupAlgebra(QQ, 3)}

\qquad\texttt{t1 = A([1,2,3]) + A([2,1,3]) + A([2,3,1])}

Let us similarly define \texttt{t2} to be our $\mathbf{t}_{2}%
=\operatorname*{cyc}\nolimits_{2}+\operatorname*{cyc}\nolimits_{2,3}$:

\qquad\texttt{t2 = A([1,2,3]) + A([1,3,2])}

Now let us ask Sage whether we have $\mathbf{t}_{1}\mathbf{t}_{2}%
=\mathbf{t}_{2}\mathbf{t}_{1}$ (that is, whether $\mathbf{t}_{1}$ and
$\mathbf{t}_{2}$ commute):

\qquad\texttt{t1*t2 == t2*t1}

Sage answers \textquotedblleft\texttt{False}\textquotedblright, so we see once
again that $\mathbf{t}_{1}$ and $\mathbf{t}_{2}$ do not commute. Keep in mind
that Sage's convention for multiplying elements in $\mathbf{k}\left[
S_{n}\right]  $ is opposite from ours (just like for multiplying
permutations), but fortunately this does not matter for checking whether two
elements commute.

We can also check that the commutator $\mathbf{t}_{1}\mathbf{t}_{2}%
-\mathbf{t}_{2}\mathbf{t}_{1}$ is nilpotent:

\qquad\texttt{(t1*t2 - t2*t1) ** 2}

(Again, \texttt{t1*t2 - t2*t1} is our $\mathbf{t}_{2}\mathbf{t}_{1}%
-\mathbf{t}_{1}\mathbf{t}_{2}$, not our $\mathbf{t}_{1}\mathbf{t}%
_{2}-\mathbf{t}_{2}\mathbf{t}_{1}$, but this does not matter as we square the
element.) \medskip

Sage knows quite a few more things about symmetric group algebras. For
example, we can compute $S\left(  \mathbf{t}_{1}\right)  $ by

\qquad\texttt{t1.antipode()}

And we can compute the jucys--murphy $\mathbf{m}_{3}$ by

\qquad\texttt{A.jucys\_murphy(3)}

Sage can also compute inverses of invertible elements of $\mathbf{k}\left[
S_{n}\right]  $ (at least for $\mathbf{k}=\mathbb{Q}$):

\qquad\texttt{(t1 + t2).inverse()}

\paragraph{A bit of programming.}

Let's now be more systematic. We are usually interested in families of
elements, not single elements, and it gets tiresome to manually define them
one by one. Sage is a programming language (a dialect of Python), so we can
\textquotedblleft teach it how to fish\textquotedblright. For instance, let us
code a function that computes the somewhere-to-below shuffles%
\[
\mathbf{t}_{k}=\sum_{\ell=k}^{n}\operatorname*{cyc}\nolimits_{k,k+1,\ldots
,\ell}\in\mathbf{k}\left[  S_{n}\right]
\]
for arbitrary $n\in\mathbb{N}$ and $k\in\left[  n\right]  $. First, we observe
that the cycle $\operatorname*{cyc}\nolimits_{k,k+1,\ldots,\ell}$ can be
written in OLN as follows:%
\[
\operatorname*{cyc}\nolimits_{k,k+1,\ldots,\ell}=\operatorname*{oln}\left(
\underbrace{1,2,\ldots,k-1}_{\text{\texttt{list(range(1, k))}}}%
,\underbrace{k+1,k+2,\ldots,\ell}_{\text{\texttt{list(range(k+1, l+1))}}%
},k,\underbrace{\ell+1,\ell+2,\ldots,n}_{\text{\texttt{list(range(l+1, n+1))}%
}}\right)  ,
\]
and thus can be obtained in Sage as follows:

\texttt{Permutation(list(range(1, k)) + list(range(k+1, l+1))}

\texttt{\qquad\qquad\qquad\qquad\qquad\qquad\qquad\qquad+ [k] +
list(range(l+1, n+1)))}

(The \textquotedblleft\texttt{+}\textquotedblright\ operator in Python/Sage
does not only add numbers; it also concatenates lists: For instance, $\left[
u,v,w\right]  +\left[  x,y\right]  =\left[  u,v,w,x,y\right]  $.)

To obtain $\mathbf{t}_{k}$, we need to sum these cycles $\operatorname*{cyc}%
\nolimits_{k,k+1,\ldots,\ell}$ for all $\ell\in\left\{  k,k+1,\ldots
,n\right\}  $ (or, as Sage would say, for all \texttt{l in range(k, n+1)}).
Thus, we end up with the following code:

\texttt{\qquad def QS(n):} \# this is the group algebra $\mathbb{Q}\left[
S_{n}\right]  $

\texttt{\qquad\qquad\qquad return SymmetricGroupAlgebra(QQ, n)}

\texttt{\qquad def t(n, k):} \# this is our $\mathbf{t}_{k}$ in $\mathbb{Q}%
\left[  S_{n}\right]  $

\texttt{\qquad\qquad\qquad A = QS(n)}

\texttt{\qquad\qquad\qquad return A.sum(A(list(range(1, k)) + list(range(k+1,
l+1))}

\texttt{\qquad\qquad\qquad\qquad\qquad\qquad\qquad\qquad+ [k] +
list(range(l+1, n+1)))}

\texttt{\qquad\qquad\qquad\qquad\qquad\qquad\qquad for l in range(k, n+1))}

You can check this by letting Sage compute \texttt{t(3, 1)} or \texttt{t(5,
2)} (for example). \medskip

Let us play around a bit more with these elements. For instance, Exercise
\ref{exe.stb.prod=Nabla} for $n=5$ claims that $\mathbf{t}_{5}\mathbf{t}%
_{4}\mathbf{t}_{3}\mathbf{t}_{2}\mathbf{t}_{1}=\nabla$ in $\mathbf{k}\left[
S_{5}\right]  $. Let us check this with Sage (for $\mathbf{k}=\mathbb{Q}$):

\texttt{\qquad LHS = QS(5).prod(t(5, k) for k in range(1, 6))}

\texttt{\qquad RHS = QS(5).sum(QS(5)(w) for w in Permutations(5))}

\texttt{\qquad LHS == RHS}

The \texttt{LHS} here is saying \textquotedblleft the product of \texttt{t(5,
k)} for $k\in\left\{  1,2,3,4,5\right\}  $ in the ring \texttt{QS(5)}%
\textquotedblright, which in our language is $\mathbf{t}_{5}\mathbf{t}%
_{4}\mathbf{t}_{3}\mathbf{t}_{2}\mathbf{t}_{1}$. (Yes, Sage multiplies the
\texttt{t(5, k)} in the order of increasing $k$, but because Sage's order of
multiplication in $\mathbf{k}\left[  S_{n}\right]  $ is the opposite of ours,
this results in $\mathbf{t}_{5}\mathbf{t}_{4}\mathbf{t}_{3}\mathbf{t}%
_{2}\mathbf{t}_{1}$ rather than $\mathbf{t}_{1}\mathbf{t}_{2}\mathbf{t}%
_{3}\mathbf{t}_{4}\mathbf{t}_{5}$.) The \texttt{RHS} is our $\nabla$. The last
line (\texttt{LHS == RHS}) thus asks Sage whether $\mathbf{t}_{5}%
\mathbf{t}_{4}\mathbf{t}_{3}\mathbf{t}_{2}\mathbf{t}_{1}=\nabla$. And Sage
quickly answers \texttt{True}.

\medskip

Incidentally, what happens if we multiply the $\mathbf{t}_{k}$'s in the
opposite order? In other words, what is $\mathbf{t}_{1}\mathbf{t}%
_{2}\mathbf{t}_{3}\mathbf{t}_{4}\mathbf{t}_{5}$ ? Sage can easily answer this:

\texttt{\qquad QS(5).prod(t(5, k) for k in range(5, 0, -1))}

Alas, the result does not look like there is a simple pattern to its coefficients.

\begin{question}
Is there? Can we at least describe the permutations that appear in the result
with nonzero coefficient? In other words, what permutations $w\in S_{n}$ can
be written in the form
\[
w=\operatorname*{cyc}\nolimits_{1,2,\ldots,i_{1}}\operatorname*{cyc}%
\nolimits_{2,3,\ldots,i_{2}}\operatorname*{cyc}\nolimits_{3,4,\ldots,i_{3}%
}\cdots\operatorname*{cyc}\nolimits_{n}%
\]
for some $i_{1},i_{2},\ldots,i_{n}\in\left[  n\right]  $ with $i_{k}\geq k$
for each $k\in\left[  n\right]  $ ? Apparently, the number of such
permutations $w$ is given by \href{https://oeis.org/A229046}{OEIS Sequence
A229046}; can this be proved?
\end{question}

\subsubsection{Exercises}

In the following exercises, you are supposed to give both the answer and the
Sage code used to obtain it.

\begin{exercise}
\fbox{1} Write Sage code to compute the $\mathbf{X}_{a,b}$ from Definition
\ref{def.Xab.Xab}.
\end{exercise}

\begin{exercise}
\fbox{1} Write Sage code to compute the random-to-random shuffles
$\mathbf{R}_{k}$ (Definition \ref{def.R2R.Rk}) or find the inbuilt function
that computes them.

[\textbf{Hint:} They appear in Sage's \texttt{SymmetricGroupAlgebra} class
under a different name.]
\end{exercise}

\begin{exercise}
\fbox{2} Playing with Exercise \ref{exe.antip-conj.perm.invol}, you may have
made the following conjecture: \textquotedblleft If two permutations $u,v\in
S_{n}$ are conjugate, then there exists an involution $w\in S_{n}$ such that
$u=wvw^{-1}$ (or, equivalently, $u=wvw$).\textquotedblright

Use Sage to check whether this conjecture holds for all $n\leq6$.

[\textbf{Hint:} This requires \texttt{for}-loops.]
\end{exercise}

\begin{exercise}
An \emph{inversion} of a permutation $w\in S_{n}$ means a pair $\left(
i,j\right)  \in\left[  n\right]  ^{2}$ satisfying $i<j$ and $w\left(
i\right)  >w\left(  j\right)  $. The \# of inversions of a given permutation
$w\in S_{n}$ will be denoted $\ell\left(  w\right)  $.

For each $t\in\mathbf{k}$, define the element
\[
\mathbf{L}_{t}:=\sum_{w\in S_{n}}t^{\ell\left(  w\right)  }w\in\mathbf{k}%
\left[  S_{n}\right]  .
\]

(In particular, $\mathbf{L}_{1}=\nabla$ and $\mathbf{L}_{-1}=\nabla^{-}$ and
$\mathbf{L}_{0}=1$.) \medskip

\textbf{(a)} \fbox{1} Is it true that $\mathbf{L}_{t}\mathbf{L}_{s}%
=\mathbf{L}_{s}\mathbf{L}_{t}$ for all $t,s\in\mathbf{k}$ ? \medskip

\textbf{(b)} \fbox{1} Is it true that $\operatorname*{span}%
\nolimits_{\mathbf{k}}\left\{  \mathbf{L}_{t}\ \mid\ t\in\mathbf{k}\right\}  $
is a $\mathbf{k}$-subalgebra of $\mathbf{k}\left[  S_{n}\right]  $ ? \medskip

[\textbf{Hint:} The cases $n\leq3$ are not representative of the situation.
See \cite{AdeGri26} for more about these elements $\mathbf{L}_{t}$.]
\end{exercise}

\subsection{\label{sec.intro.struct}A look at the structure of symmetric group
algebras}

Now let us determine the structure of the $\mathbf{k}$-algebra $\mathbf{k}%
\left[  S_{n}\right]  $ for small values of $n$ (specifically, $2$ and $3$).
The word \textquotedblleft structure\textquotedblright\ here refers to the
question \textquotedblleft is it isomorphic to an algebra we
know?\textquotedblright\ and related questions. While at that, let us also try
to understand the subalgebras $Z\left(  \mathbf{k}\left[  S_{n}\right]
\right)  $ and $\operatorname*{GZ}\nolimits_{n}$ (defined in Corollary
\ref{cor.GZ.comm}).

\subsubsection{$\mathbf{k}\left[  S_{2}\right]  $}

We begin with the case $n=2$.

The group $S_{2}$ is abelian, so its group algebra $\mathbf{k}\left[
S_{2}\right]  $ is commutative. This algebra has basis $\left(  1,s_{1}%
\right)  $ (as a $\mathbf{k}$-module) with multiplication table%
\[
1\cdot1=1,\ \ \ \ \ \ \ \ \ \ 1s_{1}=s_{1},\ \ \ \ \ \ \ \ \ \ s_{1}%
1=s_{1},\ \ \ \ \ \ \ \ \ \ s_{1}s_{1}=1.
\]
Have you seen this $\mathbf{k}$-algebra before?

For comparison, consider the quotient ring $\mathbf{k}\left[  x\right]
/\left(  x^{2}-1\right)  $ of the polynomial ring $\mathbf{k}\left[  x\right]
$ modulo the principal ideal generated by $x^{2}-1$. This quotient ring is
known as the ring of
\href{https://en.wikipedia.org/wiki/Split-complex_number}{\emph{split-complex
numbers}}\footnote{The name is a reference to the ring of \emph{complex
numbers} over $\mathbf{k}$, which is defined as $\mathbf{k}\left[  x\right]
/\left(  x^{2}+1\right)  $. (Of course, for $\mathbf{k}=\mathbb{R}$, this is
the standard ring $\mathbb{C}$ of complex numbers.) Thus, the split-complex
numbers differ from the complex numbers only in that the polynomial $x^{2}-1$
is used instead of $x^{2}+1$.} over $\mathbf{k}$, and has a basis $\left(
\overline{1},\overline{x}\right)  $ (where $\overline{p}$ means the residue
class of $p\in\mathbf{k}\left[  x\right]  $ in the quotient ring) with
multiplication table%
\[
\overline{1}\cdot\overline{1}=\overline{1},\ \ \ \ \ \ \ \ \ \ \overline
{1}\cdot\overline{x}=\overline{x},\ \ \ \ \ \ \ \ \ \ \overline{x}%
\cdot\overline{1}=\overline{x},\ \ \ \ \ \ \ \ \ \ \overline{x}\cdot
\overline{x}=\overline{1}.
\]
This is exactly the multiplication table we found for $\mathbf{k}\left[
S_{2}\right]  $ above, except that $1$ and $s_{1}$ have been renamed as
$\overline{1}$ and $\overline{x}$.

Thus, we obtain a $\mathbf{k}$-algebra isomorphism%
\begin{align*}
\mathbf{k}\left[  x\right]  /\left(  x^{2}-1\right)   &  \rightarrow
\mathbf{k}\left[  S_{2}\right]  ,\\
\overline{x}  &  \mapsto s_{1}.
\end{align*}
(We need not say that it sends $\overline{1}\mapsto1$, because this is
automatically true for $\mathbf{k}$-algebra morphisms.)

Thus, $\mathbf{k}\left[  S_{2}\right]  $ is isomorphic to $\mathbf{k}\left[
x\right]  /\left(  x^{2}-1\right)  $. We can say more if we know a bit about
$\mathbf{k}$. For instance, for $\mathbf{k}=\mathbb{Q}$, we have%
\begin{align*}
\mathbb{Q}\left[  x\right]  /\left(  x^{2}-1\right)   &  =\mathbb{Q}\left[
x\right]  /\left(  \left(  x-1\right)  \left(  x+1\right)  \right) \\
&  \cong\underbrace{\left(  \mathbb{Q}\left[  x\right]  /\left(  x-1\right)
\right)  }_{\cong\mathbb{Q}}\times\underbrace{\left(  \mathbb{Q}\left[
x\right]  /\left(  x+1\right)  \right)  }_{\cong\mathbb{Q}}\\
&  \ \ \ \ \ \ \ \ \ \ \ \ \ \ \ \ \ \ \ \ \left(
\begin{array}
[c]{c}%
\text{by the Chinese Remainder Theorem}\\
\text{(for }\mathbb{Q}\text{-algebras), since the polynomials }x-1\\
\text{and }x+1\text{ are coprime}%
\end{array}
\right) \\
&  \cong\mathbb{Q}\times\mathbb{Q}\ \ \ \ \ \ \ \ \ \ \left(
\begin{array}
[c]{c}%
\text{a direct product of }\mathbb{Q}\text{-algebras, with}\\
\text{entrywise addition and multiplication}%
\end{array}
\right)  .
\end{align*}
This can also be shown more explicitly: There is a $\mathbb{Q}$-algebra
isomorphism%
\begin{align*}
\mathbb{Q}\left[  x\right]  /\left(  x^{2}-1\right)   &  \rightarrow
\mathbb{Q}\times\mathbb{Q},\\
\overline{a+bx}  &  \mapsto\left(  a+b,\ a-b\right)
\ \ \ \ \ \ \ \ \ \ \text{for all }a,b\in\mathbb{Q}.
\end{align*}

More generally, for any commutative ring $\mathbf{k}$, there is a $\mathbf{k}%
$-algebra morphism%
\begin{align*}
\mathbf{k}\left[  x\right]  /\left(  x^{2}-1\right)   &  \rightarrow
\mathbf{k}\times\mathbf{k},\\
\overline{a+bx}  &  \mapsto\left(  a+b,\ a-b\right)
\ \ \ \ \ \ \ \ \ \ \text{for all }a,b\in\mathbf{k}.
\end{align*}
When $2$ is invertible in $\mathbf{k}$, this is a $\mathbf{k}$-algebra
\textbf{iso}morphism, so we conclude that%
\[
\mathbf{k}\left[  S_{2}\right]  \cong\mathbf{k}\left[  x\right]  /\left(
x^{2}-1\right)  \cong\mathbf{k}\times\mathbf{k}\ \ \ \ \ \ \ \ \ \ \text{in
this case.}%
\]
If $\mathbf{k}$ is a field of characteristic $2$, then the ring $\mathbf{k}%
\left[  S_{2}\right]  $ has a nonzero nilpotent element $\nabla$ (since
Corollary \ref{cor.integral.square} yields $\nabla^{2}=2!\cdot\nabla=0$), so
it cannot be isomorphic to $\mathbf{k}\times\mathbf{k}$. Instead, in this
case, we have $\mathbf{k}\left[  S_{2}\right]  \cong\mathbf{k}\left[
x\right]  /\left(  x^{2}\right)  $ (the ring of
\href{https://en.wikipedia.org/wiki/Dual_number}{\emph{dual numbers}} over
$\mathbf{k}$), since $x^{2}-1=x^{2}-1^{2}=\left(  x-1\right)  ^{2}$ can be
turned into $x^{2}$ by substituting $x$ for $x-1$.

We need not analyze $Z\left(  \mathbf{k}\left[  S_{2}\right]  \right)  $ or
$\operatorname*{GZ}\nolimits_{2}$ separately, since we have $Z\left(
\mathbf{k}\left[  S_{2}\right]  \right)  =\operatorname*{GZ}\nolimits_{2}%
=\mathbf{k}\left[  S_{2}\right]  $.

\subsubsection{\label{subsec.intro.struct.3}$\mathbf{k}\left[  S_{3}\right]
$}

Now, we move on to the case $n=3$.

We have $\left\vert S_{3}\right\vert =6$ and $S_{3}=\left\{  1,\ s_{1}%
,\ s_{2},\ t_{1,3},\ z,\ z^{-1}\right\}  $, where%
\[
z:=\operatorname*{cyc}\nolimits_{1,2,3}\ \ \ \ \ \ \ \ \ \ \text{and
thus}\ \ \ \ \ \ \ \ \ \ z^{-1}=\operatorname*{cyc}\nolimits_{1,3,2}=z^{2}.
\]
So the $\mathbf{k}$-algebra $\mathbf{k}\left[  S_{3}\right]  $ has basis
$\left(  1,\ s_{1},\ s_{2},\ t_{1,3},\ z,\ z^{-1}\right)  $ as a $\mathbf{k}%
$-module. What does it look like as a $\mathbf{k}$-algebra?

We begin with its center $Z\left(  \mathbf{k}\left[  S_{3}\right]  \right)  $,
since it is smaller. By our study of centers of group algebras (specifically,
Theorem \ref{thm.center.conjsums} and Corollary \ref{cor.type.conjclasses}),
we know that $Z\left(  \mathbf{k}\left[  S_{3}\right]  \right)  $ has basis
$\left(  1,\ \mathbf{t},\ \mathbf{c}\right)  $, where
\begin{align*}
\mathbf{t}  &  =s_{1}+s_{2}+t_{1,3}=\left(  \text{sum of all transpositions}%
\right)  \ \ \ \ \ \ \ \ \ \ \text{and}\\
\mathbf{c}  &  =z+z^{-1}=\left(  \text{sum of all cycles }\operatorname*{cyc}%
\nolimits_{i,j,k}\right)  .
\end{align*}
The multiplication table for this basis is%
\[
\mathbf{t}^{2}=3+3\mathbf{c},\ \ \ \ \ \ \ \ \ \ \mathbf{tc}=2\mathbf{t}%
,\ \ \ \ \ \ \ \ \ \ \mathbf{c}^{2}=2+\mathbf{c}%
\]
(the unitality and the commutativity of $Z\left(  \mathbf{k}\left[
S_{3}\right]  \right)  $ take care of the rest). For $\mathbf{k}=\mathbb{Q}$,
the first equality can be solved for $\mathbf{c}$, leading to $\mathbf{c}%
=\dfrac{\mathbf{t}^{2}-3}{3}$, and then the other two equalities become%
\[
\mathbf{t}\left(  \mathbf{t}-3\right)  \left(  \mathbf{t}+3\right)
=0\ \ \ \ \ \ \ \ \ \ \text{and}\ \ \ \ \ \ \ \ \ \ \mathbf{t}^{2}\left(
\mathbf{t}-3\right)  \left(  \mathbf{t}+3\right)  =0
\]
(after a bit of computation). Of course, the last equality is redundant, so we
only need to deal with $\mathbf{t}\left(  \mathbf{t}-3\right)  \left(
\mathbf{t}+3\right)  =0$. The three elements $1,\ \mathbf{t},\ \mathbf{t}^{2}$
are $\mathbb{Q}$-linearly independent, and the equality $\mathbf{t}\left(
\mathbf{t}-3\right)  \left(  \mathbf{t}+3\right)  =0$ lets us express any
higher power of $\mathbf{t}$ in terms of them. Thus, $\left(  1,\ \mathbf{t}%
,\ \mathbf{t}^{2}\right)  $ is a basis of the $\mathbb{Q}$-vector space
$Z\left(  \mathbb{Q}\left[  S_{3}\right]  \right)  $, and we get%
\begin{align*}
Z\left(  \mathbb{Q}\left[  S_{3}\right]  \right)   &  \cong\mathbb{Q}\left[
x\right]  /\left(  x\left(  x-3\right)  \left(  x+3\right)  \right) \\
&  \ \ \ \ \ \ \ \ \ \ \ \ \ \ \ \ \ \ \ \ \left(
\begin{array}
[c]{c}%
\text{via the }\mathbb{Q}\text{-algebra isomorphism}\\
\text{from }\mathbb{Q}\left[  x\right]  /\left(  x\left(  x-3\right)  \left(
x+3\right)  \right)  \text{ to }Z\left(  \mathbb{Q}\left[  S_{3}\right]
\right) \\
\text{that sends }\overline{x}\text{ to }\mathbf{t}%
\end{array}
\right) \\
&  \cong\mathbb{Q}\times\mathbb{Q}\times\mathbb{Q}\ \ \ \ \ \ \ \ \ \ \left(
\text{again by the Chinese Remainder Theorem}\right)  .
\end{align*}
More generally, we get a $\mathbf{k}$-algebra isomorphism%
\[
Z\left(  \mathbf{k}\left[  S_{3}\right]  \right)  \cong\mathbf{k}%
\times\mathbf{k}\times\mathbf{k}%
\]
whenever $2$ and $3$ are invertible in $\mathbf{k}$. When $\mathbf{k}$ is a
field satisfying $3=0$ in $\mathbf{k}$, this is false, since $\mathbf{t}%
^{2}=0$ in this case. This is also false when $\mathbf{k}$ is a field with
$2=0$ in $\mathbf{k}$, since $\left(  \mathbf{c}+\mathbf{t}+1\right)  ^{2}=0$
in this case.

\begin{exercise}
Let $\mathbf{k}$ be the field with $3$ elements (that is, $\mathbb{F}%
_{3}=\mathbb{Z}/3$). \medskip

\textbf{(a)} \fbox{2} Prove that the center $Z\left(  \mathbf{k}\left[
S_{3}\right]  \right)  $ is neither isomorphic to $\mathbf{k}\times
\mathbf{k}\times\mathbf{k}$, nor isomorphic to a quotient ring of
$\mathbf{k}\left[  x\right]  $. \medskip

\textbf{(b)} \fbox{2} Prove that $Z\left(  \mathbf{k}\left[  S_{3}\right]
\right)  $ is isomorphic to the quotient ring $\mathbf{k}\left[  x,y\right]
/\left(  x^{2},xy,y^{2}\right)  $. \medskip

[\textbf{Hint:} For part \textbf{(a)}, compute $\left(  \kappa+\lambda
\mathbf{t}+\mu\mathbf{c}\right)  ^{2}$ for $\lambda,\mu,\kappa\in\mathbf{k}$,
and show that this is always a $\mathbf{k}$-linear combination of $1$ and
$\lambda\mathbf{t}+\mu\mathbf{c}$.]
\end{exercise}

Next, we consider the Gelfand--Tsetlin subalgebra $\operatorname*{GZ}%
\nolimits_{3}$. We can easily see that $\operatorname*{GZ}\nolimits_{3}$ has
basis $\left(  1,\ \mathbf{t},\ \mathbf{c},\ s_{1}\right)  $ (note that
$\mathbf{m}_{3}=\mathbf{t}-s_{1}$). The multiplication table consists of the
formulas for $Z\left(  \mathbf{k}\left[  S_{3}\right]  \right)  $ we found
above, plus%
\[
s_{1}^{2}=1,\ \ \ \ \ \ \ \ \ \ \mathbf{t}s_{1}=\mathbf{c}%
+1,\ \ \ \ \ \ \ \ \ \ \mathbf{c}s_{1}=\mathbf{t}-s_{1}.
\]

Can we identify this algebra? At least for $\mathbf{k}=\mathbb{Q}$, we can try
to find an element with $4$ linearly independent powers, since such an element
will generate $\operatorname*{GZ}\nolimits_{3}$ and thus allow us to write
$\operatorname*{GZ}\nolimits_{3}$ as a quotient of $\mathbf{k}\left[
x\right]  $.

We know such an element: the jucys--murphy $\mathbf{m}_{3}$, whose first $4$
powers are linearly independent (at least when $2$ is invertible in
$\mathbf{k}$), and which satisfies%
\[
\left(  \mathbf{m}_{3}-2\right)  \left(  \mathbf{m}_{3}-1\right)  \left(
\mathbf{m}_{3}+1\right)  \left(  \mathbf{m}_{3}+2\right)  =0.
\]
Thus, when $2$ is invertible in $\mathbf{k}$, we obtain%
\[
\operatorname*{GZ}\nolimits_{3}\cong\mathbf{k}\left[  x\right]  /\left(
\left(  x-2\right)  \left(  x-1\right)  \left(  x+1\right)  \left(
x+2\right)  \right)  .
\]
When $2$ and $3$ are invertible in $\mathbf{k}$, this is furthermore
isomorphic to $\mathbf{k}\times\mathbf{k}\times\mathbf{k}\times\mathbf{k}$
using the Chinese Remainder Theorem. \medskip

What about the whole group algebra $\mathbf{k}\left[  S_{3}\right]  $ ? Let us
again start with the case $\mathbf{k}=\mathbb{Q}$. We begin by identifying
central idempotents. A \emph{central idempotent} in a $\mathbf{k}$-algebra $A$
means an idempotent in $Z\left(  A\right)  $. When $z$ is a central idempotent
of $A$, we can decompose $A$ as a direct product $\left(  zA\right)
\times\left(  \left(  1-z\right)  A\right)  $ of two $\mathbf{k}%
$-algebras\footnote{Here we are using the standard notation%
\[
zA=\left\{  za\ \mid\ a\in A\right\}  .
\]
For an arbitrary $z\in Z\left(  A\right)  $, this is an ideal of $A$, and thus
closed under multiplication. When $z$ is furthermore idempotent, this ideal
$zA$ furthermore is a $\mathbf{k}$-algebra with unity $z$ (and multiplication
inherited from $A$). Moreover, when $z$ is a central idempotent, $1-z$ is a
central idempotent as well, and the two ideals $zA$ and $\left(  1-z\right)
A$ add up to the whole $A$ (that is, $A=\left(  zA\right)  \oplus\left(
\left(  1-z\right)  A\right)  $ as an internal direct sum).} (see
\cite[Exercise 2.10.4 and the following few paragraphs]{23wa}). We already
know a central idempotent of $\mathbb{Q}\left[  S_{3}\right]  $: namely,
$\dfrac{1}{6}\nabla$ (since Corollary \ref{cor.integral.square} shows that
$\dfrac{1}{n!}\nabla$ is idempotent for any $n\in\mathbb{N}$). Thus, we can
write $\mathbb{Q}\left[  S_{3}\right]  $ as a direct product%
\[
\mathbb{Q}\left[  S_{3}\right]  \cong\left(  \left(  \dfrac{1}{6}%
\nabla\right)  \mathbb{Q}\left[  S_{3}\right]  \right)  \times\left(  \left(
1-\dfrac{1}{6}\nabla\right)  \mathbb{Q}\left[  S_{3}\right]  \right)  .
\]
The first factor here is the $\mathbb{Q}$-algebra $\left(  \dfrac{1}{6}%
\nabla\right)  \mathbb{Q}\left[  S_{3}\right]  $, which has dimension $1$ as a
$\mathbb{Q}$-vector space (since any multiple of $\dfrac{1}{6}\nabla$ in
$\mathbb{Q}\left[  S_{3}\right]  $ is a \textbf{scalar} multiple of $\nabla$,
due to (\ref{eq.prop.integral.fix.+})), and thus is $\cong\mathbb{Q}$ as a
$\mathbb{Q}$-algebra (why?). The second factor is the $\mathbb{Q}$-algebra
$\left(  1-\dfrac{1}{6}\nabla\right)  \mathbb{Q}\left[  S_{3}\right]  $, which
has dimension $5$.

This second factor can be decomposed further, since it contains another
central idempotent: namely, $\dfrac{1}{6}\nabla^{-}$. (The centrality and the
idempotency of this element are shown just as easily as for $\dfrac{1}%
{6}\nabla$. To see that it belongs to $\left(  1-\dfrac{1}{6}\nabla\right)
\mathbb{Q}\left[  S_{3}\right]  $, we must check that $\nabla\nabla^{-}=0$,
but this follows easily from Theorem \ref{thm.int.+-=0}.) Thus, we obtain%
\[
\left(  1-\dfrac{1}{6}\nabla\right)  \mathbb{Q}\left[  S_{3}\right]
\cong\left(  \left(  \dfrac{1}{6}\nabla^{-}\right)  \mathbb{Q}\left[
S_{3}\right]  \right)  \times\left(  \left(  1-\dfrac{1}{6}\nabla-\dfrac{1}%
{6}\nabla^{-}\right)  \mathbb{Q}\left[  S_{3}\right]  \right)  .
\]

Again, the first factor is $1$-dimensional and isomorphic to $\mathbb{Q}$. It
remains to decompose the $\mathbb{Q}$-algebra $\left(  1-\dfrac{1}{6}%
\nabla-\dfrac{1}{6}\nabla^{-}\right)  \mathbb{Q}\left[  S_{3}\right]  $ here,
which has dimension $4$ as a $\mathbb{Q}$-vector space. At this point, we have
run out of central idempotents (except the obvious ones -- $0$ and $1$ --
which are useless, as they yield trivial decompositions). What can we say
about this $4$-dimensional noncommutative algebra?

What follows is a deus ex machina; I don't know how to motivate it at this
point (although eventually I will). Consider the four elements%
\begin{align*}
\mathbf{x}  &  :=\dfrac{1}{3}\operatorname*{id}+\,\dfrac{1}{6}\mathbf{m}%
_{3}-\dfrac{1}{3}s_{1}-\dfrac{1}{6}\left(  z+z^{-1}\right)  ,\\
\mathbf{y}  &  :=\dfrac{1}{3}\left(  s_{2}-z+z^{-1}-t_{1,3}\right)  ,\\
\mathbf{z}  &  :=\dfrac{3}{4}S\left(  \mathbf{y}\right)  =\dfrac{-3}%
{4}T_{\operatorname*{sign}}\left(  \mathbf{y}\right)  ,\\
\mathbf{w}  &  :=T_{\operatorname*{sign}}\left(  \mathbf{x}\right)
\end{align*}
of $\mathbb{Q}\left[  S_{3}\right]  $. Then, straightforward (if tiresome)
computations show that $\left(  \mathbf{x},\ \mathbf{y},\ \mathbf{z}%
,\ \mathbf{w}\right)  $ is a basis of the very subspace $\left(  1-\dfrac
{1}{6}\nabla-\dfrac{1}{6}\nabla^{-}\right)  \mathbb{Q}\left[  S_{3}\right]  $
that we are trying to analyze. Moreover, this basis has a nice multiplication
table:%
\begin{align*}
\mathbf{xx}  &  =\mathbf{x},\ \ \ \ \ \ \ \ \ \ \mathbf{xy}=\mathbf{y}%
,\ \ \ \ \ \ \ \ \ \ \mathbf{xz}=0,\ \ \ \ \ \ \ \ \ \ \mathbf{xw}=0,\\
\mathbf{yx}  &  =0,\ \ \ \ \ \ \ \ \ \ \mathbf{yy}%
=0,\ \ \ \ \ \ \ \ \ \ \mathbf{yz}=\mathbf{x},\ \ \ \ \ \ \ \ \ \ \mathbf{yw}%
=\mathbf{y},\\
\mathbf{zx}  &  =\mathbf{z},\ \ \ \ \ \ \ \ \ \ \mathbf{zy}=\mathbf{w}%
,\ \ \ \ \ \ \ \ \ \ \mathbf{zz}=0,\ \ \ \ \ \ \ \ \ \ \mathbf{zw}=0,\\
\mathbf{wx}  &  =0,\ \ \ \ \ \ \ \ \ \ \mathbf{wy}%
=0,\ \ \ \ \ \ \ \ \ \ \mathbf{wz}=\mathbf{z},\ \ \ \ \ \ \ \ \ \ \mathbf{ww}%
=\mathbf{w}.
\end{align*}
This is exactly the multiplication table for the elementary matrices%
\[
E_{1,1}=\left(
\begin{array}
[c]{cc}%
1 & 0\\
0 & 0
\end{array}
\right)  ,\ \ \ \ \ \ \ \ \ \ E_{1,2}=\left(
\begin{array}
[c]{cc}%
0 & 1\\
0 & 0
\end{array}
\right)  ,\ \ \ \ \ \ \ \ \ \ E_{2,1}=\left(
\begin{array}
[c]{cc}%
0 & 0\\
1 & 0
\end{array}
\right)  ,\ \ \ \ \ \ \ \ \ \ E_{2,2}=\left(
\begin{array}
[c]{cc}%
0 & 0\\
0 & 1
\end{array}
\right)
\]
in the matrix ring $\mathbb{Q}^{2\times2}$, with $E_{1,1},\ E_{1,2}%
,\ E_{2,1},\ E_{2,2}$ renamed as $\mathbf{x},\ \mathbf{y},\ \mathbf{z}%
,\ \mathbf{w}$. Thus, there is a $\mathbb{Q}$-algebra isomorphism%
\[
\mathbb{Q}^{2\times2}\rightarrow\left(  1-\dfrac{1}{6}\nabla-\dfrac{1}%
{6}\nabla^{-}\right)  \mathbb{Q}\left[  S_{3}\right]
\]
that sends $E_{1,1},\ E_{1,2},\ E_{2,1},\ E_{2,2}$ to $\mathbf{x}%
,\ \mathbf{y},\ \mathbf{z},\ \mathbf{w}$, respectively. Consequently,%
\[
\left(  1-\dfrac{1}{6}\nabla-\dfrac{1}{6}\nabla^{-}\right)  \mathbb{Q}\left[
S_{3}\right]  \cong\mathbb{Q}^{2\times2}.
\]

Altogether, we have obtained%
\begin{align*}
\mathbb{Q}\left[  S_{3}\right]   &  \cong\left(  \left(  \dfrac{1}{6}%
\nabla\right)  \mathbb{Q}\left[  S_{3}\right]  \right)  \times
\underbrace{\left(  \left(  1-\dfrac{1}{6}\nabla\right)  \mathbb{Q}\left[
S_{3}\right]  \right)  }_{\cong\left(  \left(  \dfrac{1}{6}\nabla^{-}\right)
\mathbb{Q}\left[  S_{3}\right]  \right)  \times\left(  \left(  1-\dfrac{1}%
{6}\nabla-\dfrac{1}{6}\nabla^{-}\right)  \mathbb{Q}\left[  S_{3}\right]
\right)  }\\
&  \cong\underbrace{\left(  \left(  \dfrac{1}{6}\nabla\right)  \mathbb{Q}%
\left[  S_{3}\right]  \right)  }_{\cong\mathbb{Q}}\times\underbrace{\left(
\left(  \dfrac{1}{6}\nabla^{-}\right)  \mathbb{Q}\left[  S_{3}\right]
\right)  }_{\cong\mathbb{Q}}\times\underbrace{\left(  \left(  1-\dfrac{1}%
{6}\nabla-\dfrac{1}{6}\nabla^{-}\right)  \mathbb{Q}\left[  S_{3}\right]
\right)  }_{\cong\mathbb{Q}^{2\times2}}\\
&  \cong\mathbb{Q}\times\mathbb{Q}\times\mathbb{Q}^{2\times2}.
\end{align*}
Explicitly, we have a $\mathbb{Q}$-algebra isomorphism%
\begin{align*}
\mathbb{Q}\times\mathbb{Q}\times\mathbb{Q}^{2\times2}  &  \rightarrow
\mathbb{Q}\left[  S_{3}\right]  ,\\
\left(  1,0,0\right)   &  \mapsto\dfrac{1}{6}\nabla,\\
\left(  0,1,0\right)   &  \mapsto\dfrac{1}{6}\nabla^{-},\\
\left(  0,0,E_{1,1}\right)   &  \mapsto\mathbf{x},\\
\left(  0,0,E_{1,2}\right)   &  \mapsto\mathbf{y},\\
\left(  0,0,E_{2,1}\right)   &  \mapsto\mathbf{z},\\
\left(  0,0,E_{2,2}\right)   &  \mapsto\mathbf{w}.
\end{align*}

More generally, the same construction defines a $\mathbf{k}$-algebra
isomorphism%
\[
\mathbf{k}\times\mathbf{k}\times\mathbf{k}^{2\times2}\rightarrow
\mathbf{k}\left[  S_{3}\right]
\]
for any commutative ring $\mathbf{k}$ as long as $2$ and $3$ are invertible in
$\mathbf{k}$.

Now how on earth did I find the elements $\mathbf{x},\ \mathbf{y}%
,\ \mathbf{z},\ \mathbf{w}$ that came to our rescue here? The next chapter
will give us an idea.

\subsubsection{About the general case}

We cannot go much farther than $n=3$ with such ad-hoc methods. But some
patterns are already visible in these cases. Not all these patterns persist.
For example, we observed above that $\operatorname*{GZ}\nolimits_{3}%
=\operatorname*{AlgGen}\left\{  \mathbf{m}_{3}\right\}  $ for $\mathbf{k}%
=\mathbb{Q}$, but not all $n$ satisfy $\operatorname*{GZ}\nolimits_{n}%
=\operatorname*{AlgGen}\left\{  \mathbf{m}_{n}\right\}  $. Likewise, $Z\left(
\mathbb{Q}\left[  S_{n}\right]  \right)  =\operatorname*{AlgGen}\left\{
\mathbf{t}\right\}  $ (where $\mathbf{t}$ is the sum of all transpositions in
$S_{n}$) holds for $n=3$ but not for $n=6$ (see Exercise
\ref{exe.spechtmod.sum-transp.gen} \textbf{(b)} for the reason why).

But what is true is the following:

\begin{theorem}
\label{thm.AWS.demo}Assume that $n!$ is invertible in $\mathbf{k}$. Then:
\medskip

\textbf{(a)} The $\mathbf{k}$-algebra $\mathbf{k}\left[  S_{n}\right]  $ is
isomorphic to a direct product of matrix rings over $\mathbf{k}$ (including
$\mathbf{k}^{1\times1}\cong\mathbf{k}$). That is,%
\[
\mathbf{k}\left[  S_{n}\right]  \cong\prod_{\lambda\in\Lambda}\mathbf{k}%
^{f_{\lambda}\times f_{\lambda}}%
\]
for some finite set $\Lambda$ and some positive integers $f_{\lambda}$ for all
$\lambda\in\Lambda$. Moreover, the set $\Lambda$ is the set of all partitions
of $n$, whereas the integer $f_{\lambda}$ is the \# of so-called standard
tableaux of shape $\lambda$ (a combinatorial object we will define later).
\medskip

\textbf{(b)} The center $Z\left(  \mathbf{k}\left[  S_{n}\right]  \right)  $
is isomorphic to the direct product $\mathbf{k}^{\left\vert \Lambda\right\vert
}=\prod_{\lambda\in\Lambda}\mathbf{k}$. In other words, it is isomorphic to
$\mathbf{k}^{p\left(  n\right)  }$, where $p\left(  n\right)  $ is the \# of
partitions of $n$. \medskip

\textbf{(c)} The Gelfand--Tsetlin subalgebra $\operatorname*{GZ}\nolimits_{n}$
is isomorphic to the direct product $\mathbf{k}^{t\left(  n\right)  }$, where
$t\left(  n\right)  $ is the \# of all standard tableaux of all size-$n$
partition shapes, or, equivalently, the \# of all involutions in $S_{n}$.
\end{theorem}

This is a deep and surprising theorem. We will prove its parts \textbf{(a)}
and \textbf{(b)} later. Part \textbf{(c)} requires the so-called seminormal
basis of $\mathbf{k}\left[  S_{n}\right]  $, which is beyond the scope of this course.

\begin{remark}
\label{rmk.AWS.OEIS}When $\mathbf{k}$ is a field of characteristic $0$, the
numbers $p\left(  n\right)  $ and $t\left(  n\right)  $ in Theorem
\ref{thm.AWS.demo} are the dimensions of the $\mathbf{k}$-vector spaces
$Z\left(  \mathbf{k}\left[  S_{n}\right]  \right)  $ and $\operatorname*{GZ}%
\nolimits_{n}$. Here is a little table of these numbers in Theorem
\ref{thm.AWS.demo}:%
\[%
\begin{tabular}
[c]{|c||c|c|c|c|c|c|c|c|c|c|c|c|}\hline
$n$ & $0$ & $1$ & $2$ & $3$ & $4$ & $5$ & $6$ & $7$ & $8$ & $9$ & $10$ &
$11$\\\hline
$p\left(  n\right)  $ & $1$ & $1$ & $2$ & $3$ & $5$ & $7$ & $11$ & $15$ & $22$
& $30$ & $42$ & $56$\\\hline
$t\left(  n\right)  $ & $1$ & $1$ & $2$ & $4$ & $10$ & $26$ & $76$ & $232$ &
$764$ & $2620$ & $9496$ & $35696$\\\hline
\end{tabular}
\ \ \ \ \ .
\]
Note that both numbers $p\left(  n\right)  $ and $t\left(  n\right)  $ are
much smaller than $n!=\left\vert S_{n}\right\vert =\dim\left(  \mathbf{k}%
\left[  S_{n}\right]  \right)  $ for $n\geq3$.

In the OEIS (Online Encyclopedia of Integer Sequences), the sequences $\left(
p\left(  n\right)  \right)  _{n\geq0}$ and $\left(  t\left(  n\right)
\right)  _{n\geq0}$ appear as \href{https://oeis.org/A000041}{A000041} and
\href{https://oeis.org/A000085}{A000085}, respectively.
\end{remark}

The case when $n!$ is \textbf{not} invertible in $\mathbf{k}$ (this is known
as the \emph{modular case}) is much less well-behaved; in this case,
$\mathbf{k}\left[  S_{n}\right]  $ is not a direct product of known rings. In
particular, if $\mathbf{k}$ is a field of positive characteristic $\leq n$,
then $\nabla\in\mathbf{k}\left[  S_{n}\right]  $ is nilpotent (with
$\nabla^{2}=0$), so any description of $\mathbf{k}\left[  S_{n}\right]  $ will
have to deal with both idempotents and nilpotents, which dooms it to be more
complicated. No good descriptions of $\mathbf{k}\left[  S_{n}\right]  $ in the
general case are known, although some partial results have been revealing
themselves slowly over the last 80 years (e.g., see \cite{Kuenze15} for some work).

\subsection{More exercises}

The following is an assortment of further exercises on $\mathbf{k}\left[
S_{n}\right]  $.

\begin{exercise}
\label{exe.Nabab.1}For any $a,b\in\left[  n\right]  $, we define the element%
\[
\nabla_{b,a}:=\sum_{\substack{w\in S_{n};\\w\left(  a\right)  =b}%
}w\in\mathbf{k}\left[  S_{n}\right]  .
\]
(This is a particular case of the rook sums $\nabla_{B,A}$ we will define
later on.) \medskip

\textbf{(a)} \fbox{3} Prove that every $a,b,c,d\in\left[  n\right]  $ satisfy%
\[
\nabla_{d,c}\nabla_{b,a}=%
\begin{cases}
\left(  n-2\right)  !\left(  \nabla-\nabla_{d,a}\right)  , & \text{if }b\neq
c;\\
\left(  n-1\right)  !\nabla_{d,a}, & \text{if }b=c.
\end{cases}
\]

\textbf{(b)} \fbox{2} Let $a,b\in\left[  n\right]  $, and let $X$ be a subset
of $\left[  n\right]  $ having size $\left\vert X\right\vert >2$. Prove that%
\[
\nabla_{X}^{-}\nabla_{b,a}=\nabla_{b,a}\nabla_{X}^{-}=0.
\]

\end{exercise}

The next exercise introduces another family of commuting elements of
$\mathbf{k}\left[  S_{n}\right]  $:

\begin{exercise}
\fbox{3} Let $z_{1},z_{2},\ldots,z_{n}\in\mathbf{k}$ be some scalars such that
each difference $z_{i}-z_{j}$ with $i\neq j$ is invertible in $\mathbf{k}$.
(For instance, if $\mathbf{k}$ is a field, then $z_{1},z_{2},\ldots,z_{n}$ can
be any $n$ distinct elements of $\mathbf{k}$.) Set%
\[
\mathbf{H}_{i}:=\sum_{\substack{j\in\left[  n\right]  ;\\j\neq i}}\dfrac
{1}{z_{j}-z_{i}}t_{i,j}\in\mathbf{k}\left[  S_{n}\right]
\ \ \ \ \ \ \ \ \ \ \text{for each }i\in\left[  n\right]  .
\]
Prove that the elements $\mathbf{H}_{1},\ \mathbf{H}_{2},\ \ldots
,\ \mathbf{H}_{n}$ commute. (These elements are known as the
\emph{Knizhnik--Zamolodchikov sums} in $\mathbf{k}\left[  S_{n}\right]  $.)
\end{exercise}

The next exercise continues the hunt for polynomials that annihilate elements
of $\mathbf{k}\left[  S_{n}\right]  $:

\begin{exercise}
\textbf{(a)} \fbox{2} Let $P\in\mathbb{Z}\left[  x\right]  $ be the polynomial
$\left(  x-n!+1\right)  \left(  x-1\right)  \left(  x+1\right)  $. Let $w_{0}$
be the permutation $\operatorname*{oln}\left(  n,n-1,\ldots,2,1\right)  \in
S_{n}$. Prove that $P\left(  \nabla-w_{0}\right)  =0$. \medskip

\textbf{(b)} \fbox{1} For each $\mathbf{a}\in\mathbf{k}\left[  S_{n}\right]  $
and each $m\in\mathbb{N}$, prove that%
\[
\left(  \nabla+\mathbf{a}\right)  ^{m}=\mathbf{a}^{m}+\lambda\nabla
\ \ \ \ \ \ \ \ \ \ \text{for some scalar }\lambda\in\mathbf{k}.
\]

\textbf{(c)} \fbox{1} Let $w\in S_{n}$, and let $m$ be the order of $w$ in
$S_{n}$ (that is, the smallest positive integer satisfying $w^{m}=1$). Let
$Q\in\mathbb{Z}\left[  x\right]  $ be the polynomial $\left(  x-n!+1\right)
\left(  x^{m}-\left(  -1\right)  ^{m}\right)  $. Prove that $Q\left(
\nabla-w\right)  =0$.
\end{exercise}

Next, we present a celebrated result of Olshanskii as a (difficult) exercise.
First, we recall a notion from ring theory:

\begin{definition}
Let $B$ be a subset of a ring $A$. Then, the \emph{centralizer} of $B$ in $A$
denotes the subset%
\[
\left\{  a\in A\ \mid\ ab=ba\text{ for all }b\in B\right\}  \text{ of }A
\]
(consisting of all elements of $A$ that commute with all elements of $B$).
This is itself a subring of $A$, and is denoted by $Z_{A}\left(  B\right)  $.
(See \cite[\S 2.3.5]{23wa} for some examples.)
\end{definition}

\begin{exercise}
\label{exe.YJM.olshanski}Assume that $n\geq1$. As in Remark
\ref{rmk.default-ex.abuse2}, let us regard $\mathbf{k}\left[  S_{n-1}\right]
$ as a subring of $\mathbf{k}\left[  S_{n}\right]  $ by identifying each
$\mathbf{a}\in\mathbf{k}\left[  S_{n-1}\right]  $ with its image $\left(
I_{n-1\rightarrow n}\right)  _{\ast}\left(  \mathbf{a}\right)  $ under the
default embedding of group algebras. \medskip

\textbf{(a)} \fbox{1} Prove that the jucys--murphy $\mathbf{m}_{n}%
\in\mathbf{k}\left[  S_{n}\right]  $ commutes with all elements of
$\mathbf{k}\left[  S_{n-1}\right]  $. In other words, prove that
$\mathbf{m}_{n}\in Z_{\mathbf{k}\left[  S_{n}\right]  }\left(  \mathbf{k}%
\left[  S_{n-1}\right]  \right)  $. \medskip

\textbf{(b)} \fbox{7} Prove that%
\[
Z_{\mathbf{k}\left[  S_{n}\right]  }\left(  \mathbf{k}\left[  S_{n-1}\right]
\right)  =\operatorname*{AlgGen}\left(  Z\left(  \mathbf{k}\left[
S_{n-1}\right]  \right)  \cup\left\{  \mathbf{m}_{n}\right\}  \right)  .
\]
In other words, prove that an element $\mathbf{a}\in\mathbf{k}\left[
S_{n}\right]  $ commutes with all elements of $\mathbf{k}\left[
S_{n-1}\right]  $ if and only if $\mathbf{a}$ lies in the $\mathbf{k}%
$-subalgebra of $\mathbf{k}\left[  S_{n}\right]  $ generated by the elements
of the center of $\mathbf{k}\left[  S_{n-1}\right]  $ and the $n$-th
jucys--murphy $\mathbf{m}_{n}$.
\end{exercise}

The following exercise is about the group algebras of cyclic groups, not
symmetric groups. It offers some contrast to the analysis of $\mathbf{k}%
\left[  S_{n}\right]  $ we have seen above.

\begin{exercise}
\label{exe.Cn.DFT}Let $n$ be a positive integer. Consider the cyclic group
$C_{n}=\left\{  g^{0},g^{1},\ldots,g^{n-1}\right\}  $ of size $n$. Let
$\zeta:=e^{2\pi i/n}\in\mathbb{C}$ (so that the complex numbers $\zeta
^{0},\zeta^{1},\ldots,\zeta^{n-1}$ are the vertices of a regular $n$-gon
inscribed in the unit circle).

For each $i\in\mathbb{Z}$, define the element%
\[
\Delta_{i}:=\sum_{k=0}^{n-1}\zeta^{ik}g^{k}\in\mathbb{C}\left[  C_{n}\right]
.
\]
(Note that this depends only on the residue class of $i$ modulo $n$, since
$\zeta^{n}=1$.) \medskip

\textbf{(a)} \fbox{1} Prove that $\Delta_{i}^{2}=n\Delta_{i}$ for each
$i\in\mathbb{Z}$. \medskip

\textbf{(b)} \fbox{2} Prove that $\Delta_{i}\Delta_{j}=0$ for any distinct
$i,j\in\left\{  0,1,\ldots,n-1\right\}  $. \medskip

\textbf{(c)} \fbox{2} Consider the direct product $\mathbb{C}^{n}$ of the
$\mathbb{C}$-algebras $\underbrace{\mathbb{C},\mathbb{C},\ldots,\mathbb{C}%
}_{n\text{ times}}$. Prove that the map%
\begin{align*}
\mathbb{C}^{n}  &  \rightarrow\mathbb{C}\left[  C_{n}\right]  ,\\
\left(  a_{0},a_{1},\ldots,a_{n-1}\right)   &  \mapsto\dfrac{1}{n}\sum
_{i=0}^{n-1}a_{i}\Delta_{i}%
\end{align*}
is a $\mathbb{C}$-algebra isomorphism. (This isomorphism is -- up to scalar
multiple -- known as the \emph{discrete Fourier transform}.) \medskip

\textbf{(d)} \fbox{2} Prove that the $\mathbb{R}$-algebra $\mathbb{R}\left[
C_{3}\right]  $ is \textbf{not} isomorphic to $\mathbb{R}^{3}$. (It is not
enough to show that $\zeta\notin\mathbb{R}$ for $n=3$; you should also show
that no other isomorphism can exist.) \medskip

\textbf{(e)} \fbox{2} Prove that $\mathbb{R}\left[  C_{4}\right]
\cong\mathbb{R}\times\mathbb{R}\times\mathbb{C}$ as $\mathbb{R}$-algebras.
\end{exercise}

Note that this exercise shows that the structure of $\mathbf{k}\left[
C_{n}\right]  $ can change from $\mathbf{k}=\mathbb{R}$ to $\mathbf{k}%
=\mathbb{C}$. But this does not happen for $\mathbf{k}\left[  S_{n}\right]  $,
since Theorem \ref{thm.AWS.demo} \textbf{(a)} describes the structure
independently of $\mathbf{k}$ when $\mathbf{k}$ is a $\mathbb{Q}$-algebra.

\bigskip

\section{\label{chp.rep}Actions and representations}

In this chapter, we will introduce group actions and representations of groups
and algebras. This might feel like a sideshow at first, but I hope to show
that it is actually one of the best ways to study groups and algebras. Viewing
groups or algebras through the prism of their actions is how many of the above
results were first discovered, and in some cases (e.g., for Theorem
\ref{thm.Vperm.Vn} \textbf{(b)}) it is the only way I know to prove them. It
will also explain to some extent how the elements $\mathbf{x},\mathbf{y}%
,\mathbf{z},\mathbf{w}$ of $\mathbb{Q}\left[  S_{3}\right]  $ were found in
Subsection \ref{subsec.intro.struct.3}.

\subsection{\label{sec.rep.G-sets}Group actions and $G$-sets}

\begin{convention}
\label{conv.rep.G-sets.G}For the entirety of Section \ref{sec.rep.G-sets}, we
let $G$ be a group (not necessarily finite).
\end{convention}

\subsubsection{Motivation}

Here is a slightly provocative slogan: Group elements want to be permutations.
I.e., a group wants to act.

This might sound strange. After all, we nowadays define groups as abstract
structures consisting of disembodied objects that multiply via a
multiplication table. But the elements of the (historically) first groups ever
defined (by Galois in the early 19th century) were permutations of sets. That
is, the first groups ever defined were subgroups of symmetric groups.

Even groups that are not defined in this way can be forced into this mold:

\begin{example}
\label{exa.rep.G-sets.taun}Recall the additive group $\left(  \mathbb{Z}%
,+,0\right)  $. Its elements are integers, not permutations, but you can
replace them by permutations that get multiplied (i.e., composed) just like
the corresponding integers are added. To be specific: For each $n\in
\mathbb{Z}$, we consider the map%
\begin{align*}
\tau_{n}:\mathbb{Z}  &  \rightarrow\mathbb{Z},\\
k  &  \mapsto n+k,
\end{align*}
which we call the \emph{translation} by $n$ (since it is geometrically a
translation along the number line). These translations $\tau_{n}$ are
permutations of $\mathbb{Z}$ (that is, they belong to the symmetric group
$S_{\mathbb{Z}}$), and they multiply like the $n$'s add: That is, they satisfy
$\tau_{n}\tau_{m}=\tau_{n+m}$ for all $n,m\in\mathbb{Z}$, as well as $\tau
_{0}=\operatorname*{id}$. Hence, there is a group morphism%
\begin{align*}
\mathbb{Z}  &  \rightarrow S_{\mathbb{Z}},\\
n  &  \mapsto\tau_{n}.
\end{align*}
Moreover, this morphism is injective (meaning that $\tau_{n}\neq\tau_{m}$ for
all $n\neq m$), so its image is a subgroup of $S_{\mathbb{Z}}$ isomorphic to
$\mathbb{Z}$. Thus, $\mathbb{Z}$ is isomorphic to a group of actual
permutations (namely, the image we just found).
\end{example}

We can do this trick for any group, not just for $\mathbb{Z}$. Thus we obtain:

\begin{theorem}
[Cayley's theorem]\label{thm.rep.G-sets.cay}Let $G$ be any group. For each
$g\in G$, we define the permutation%
\begin{align*}
\tau_{g}:G  &  \rightarrow G,\\
k  &  \mapsto gk
\end{align*}
of $G$. Then, the map%
\begin{align*}
G  &  \rightarrow S_{G},\\
g  &  \mapsto\tau_{g}%
\end{align*}
is an injective group morphism. Thus, its image is a subgroup of $S_{G}$
isomorphic to $G$.
\end{theorem}

And therefore $G$ is isomorphic to an actual group of permutations. This can
be useful, since permutations have some features that elements of an arbitrary
group might not have.\footnote{A well-known application of Theorem
\ref{thm.rep.G-sets.cay} is the fact that any finite group of order $4n+2$
(with $n\in\mathbb{N}$) must have a subgroup of index $2$ (that is, of order
$2n+1$). See \url{https://math.stackexchange.com/questions/225987/} for a
proof of this fact.}

Yet this is just the tip of the iceberg. We can turn the elements of a group
$G$ not just into the permutations $\tau_{g}$, but also into other
permutations of other sets, and we will soon see that this too can be rather useful.

\subsubsection{Definition}

Let us give this a name:

\begin{definition}
\label{def.rep.G-sets.left}\textbf{(a)} Let $G$ be a group, and let $X$ be a
set. A \emph{left }$G$\emph{-action on }$X$ (aka a \emph{left action of }$G$
\emph{on }$X$) means a map $\rightharpoonup\ :G\times X\rightarrow X$ that
satisfies the following two axioms:

\begin{enumerate}
\item \textit{Associativity:} We have $g\rightharpoonup\left(
h\rightharpoonup x\right)  =\left(  gh\right)  \rightharpoonup x$ for all
$g,h\in G$ and $x\in X$. Here and in the following, we use infix notation for
the map $\rightharpoonup$, meaning that we write $g\rightharpoonup x$ for the
image of $\left(  g,x\right)  \in G\times X$ under this map.

\item \textit{Unitality:} We have $1_{G}\rightharpoonup x=x$ for all $x\in X$.
(Here, $1_{G}$ means the neutral element of $G$.)
\end{enumerate}

We often omit the word \textquotedblleft left\textquotedblright\ in
\textquotedblleft left $G$-action\textquotedblright\ (or \textquotedblleft
left action\textquotedblright). We also often write $g\cdot x$ or $gx$ for
$g\rightharpoonup x$, as long as this creates no ambiguity (i.e., as long as
it cannot be mistaken for another kind of product, as it could be in Example
\ref{exa.rep.G-sets.conj} below). Thus, the axioms
\[
g\rightharpoonup\left(  h\rightharpoonup x\right)  =\left(  gh\right)
\rightharpoonup x\ \ \ \ \ \ \ \ \ \ \text{and}\ \ \ \ \ \ \ \ \ \ 1_{G}%
\rightharpoonup x=x
\]
become%
\[
g\left(  hx\right)  =\left(  gh\right)  x\ \ \ \ \ \ \ \ \ \ \text{and}%
\ \ \ \ \ \ \ \ \ \ 1_{G}x=x,
\]
which look exactly like the associativity and the unitality axioms for a group
(or monoid), but of course are different in that we are not composing two
elements of the same set but rather one element of $G$ with one element of
$X$. \medskip

\textbf{(b)} A \emph{left }$G$\emph{-set} means a set $X$ equipped with a left
$G$-action on $X$. We furthermore say that \textquotedblleft$G$ \emph{acts on
}$X$ \emph{from the left}\textquotedblright\ to mean that $X$ is a left
$G$-set. The words \textquotedblleft from the left\textquotedblright\ are
frequently omitted here.
\end{definition}

Analogously, we can define right $G$-actions and right $G$-sets. Here, the map
$\rightharpoonup\ :G\times X\rightarrow X$ is replaced by a map%
\[
\leftharpoonup\ :X\times G\rightarrow X,
\]
and the axioms are replaced by%
\[
\left(  x\leftharpoonup h\right)  \leftharpoonup g=x\leftharpoonup\left(
hg\right)  \ \ \ \ \ \ \ \ \ \ \text{and}\ \ \ \ \ \ \ \ \ \ x\leftharpoonup
1_{G}=x.
\]
We prefer to use left $G$-actions, but sometimes right $G$-actions are more
natural or convenient. The two concepts are easily translated into one another
anyway (see Exercise \ref{exe.rep.G-sets.l2r} below).

\begin{convention}
The words \textquotedblleft action\textquotedblright, \textquotedblleft%
$G$-set\textquotedblright\ and \textquotedblleft act\textquotedblright\ shall
be understood by default to mean \textquotedblleft left
action\textquotedblright, \textquotedblleft left $G$-set\textquotedblright%
\ and \textquotedblleft act from the left\textquotedblright, unless we qualify
them by the word \textquotedblleft right\textquotedblright.
\end{convention}

Definition \ref{def.rep.G-sets.left} can be generalized by replacing
\textquotedblleft group\textquotedblright\ by \textquotedblleft
monoid\textquotedblright. No other changes are necessary, and the notion of a
left action of a monoid is indeed useful, although not as mainstream and
well-behaved as that of a left action of a group.

\subsubsection{Examples}

There are lots of $G$-actions in nature. Let us list just a few examples. We
begin with three $G$-actions on the set $G$ itself (i.e., the set $X$ in
Definition \ref{def.rep.G-sets.left} is $G$ itself):

\begin{example}
\label{exa.rep.G-sets.lreg}Let $G$ be any group. The \emph{left regular }%
$G$\emph{-action} is the left $G$-action on $G$ given by%
\[
g\rightharpoonup h=gh\ \ \ \ \ \ \ \ \ \ \text{for all }g\in G\text{ and }h\in
G.
\]
In other words, this action is just the multiplication map of $G$.
\end{example}

\begin{example}
\label{exa.rep.G-sets.irreg}Let $G$ be any group. The \emph{inverse right
regular }$G$\emph{-action} is the left $G$-action on $G$ given by
\[
g\rightharpoonup h=hg^{-1}\ \ \ \ \ \ \ \ \ \ \text{for all }g\in G\text{ and
}h\in G.
\]
Note that the formula $g\rightharpoonup h=hg$ would not define a left
$G$-action, since it would fail the associativity axiom.
\end{example}

\begin{example}
\label{exa.rep.G-sets.conj}Let $G$ be any group. The \emph{conjugation }%
$G$\emph{-action} is the left $G$-action on $G$ given by
\[
g\rightharpoonup h=ghg^{-1}\ \ \ \ \ \ \ \ \ \ \text{for all }g\in G\text{ and
}h\in G.
\]

\end{example}

\begin{remark}
In Definition \ref{def.rep.G-sets.left}, it was suggested to abbreviate
$g\rightharpoonup x$ as $gx$ when $g$ is an element of the group $G$ and $x$
is an element of a left $G$-set $X$. When the $G$-set $X$ is $G$ itself, this
abbreviation is dangerous, since $gx$ could also mean the product of $g$ and
$x$ in the group $G$. In Example \ref{exa.rep.G-sets.lreg}, this would be
unproblematic, since $g\rightharpoonup x$ agrees with the product of $g$ and
$x$ by definition. However, doing the same in Example
\ref{exa.rep.G-sets.irreg} or Example \ref{exa.rep.G-sets.conj} would be a
recipe for disaster.
\end{remark}

Our next example is a $G$-action that can be defined on \textbf{any} set $X$,
although not a very interesting one:

\begin{example}
\label{exa.rep.G-sets.triv}Let $G$ be any group, and let $X$ be any set. The
\emph{trivial }$G$\emph{-action on }$X$ is the left $G$-action on $X$ given
by
\[
g\rightharpoonup x=x\ \ \ \ \ \ \ \ \ \ \text{for all }g\in G\text{ and }x\in
X.
\]

\end{example}

The next example is one of the most important in group theory. It answers the
rather natural question \textquotedblleft if a quotient $G/H$ of a group $G$
is not a group, then what is it?\textquotedblright:

\begin{example}
\label{exa.rep.G-sets.G/H}Let $G$ be any group. Let $H$ be a subgroup of $G$.
Let%
\[
G/H=\left\{  \text{left cosets of }H\text{ in }G\right\}  =\left\{
uH\ \mid\ u\in G\right\}  .
\]
This set $G/H$ is (in general) not a group (unless $H$ is a normal subgroup of
$G$), but it becomes a left $G$-set if we define the left $G$-action on $G/H$
by%
\[
g\rightharpoonup\left(  uH\right)  =\left(  gu\right)
H\ \ \ \ \ \ \ \ \ \ \text{for all }g\in G\text{ and }u\in G.
\]
Here, of course, you can abbreviate $g\rightharpoonup x$ as $gx$, so this just
takes the form%
\[
g\left(  uH\right)  =\left(  gu\right)  H\ \ \ \ \ \ \ \ \ \ \text{for all
}g\in G\text{ and }u\in G.
\]

\end{example}

So much for actions of arbitrary groups $G$. But there are many more examples
for specific groups:

\begin{example}
\label{exa.rep.G-sets.nat}Let $A$ be any set. The symmetric group $S_{A}$ acts
on the set $A$ by%
\[
g\rightharpoonup x=g\left(  x\right)  \ \ \ \ \ \ \ \ \ \ \text{for all }g\in
S_{A}\text{ and }x\in A.
\]
This is a left $S_{A}$-action on $A$, and is called the \emph{natural action}
of $S_{A}$.
\end{example}

\begin{example}
\label{exa.rep.G-sets.PA}Let $A$ be any set. Consider the set $\mathcal{P}%
\left(  A\right)  $ of all subsets of $A$. (This is the power set of $A$.) The
symmetric group $S_{A}$ acts on $\mathcal{P}\left(  A\right)  $ by%
\[
g\rightharpoonup U=g\left(  U\right)  \ \ \ \ \ \ \ \ \ \ \text{for all }g\in
S_{A}\text{ and }U\in\mathcal{P}\left(  A\right)  .
\]
(Here, as usual, $g\left(  U\right)  $ means the set $\left\{  g\left(
u\right)  \ \mid\ u\in U\right\}  $.)
\end{example}

For example, for $A=\left[  5\right]  $ and $g=\operatorname*{cyc}%
\nolimits_{1,2,3}\in S_{A}$, we have
\[
g\rightharpoonup\left\{  2,3,5\right\}  =g\left(  \left\{  2,3,5\right\}
\right)  =\left\{  g\left(  2\right)  ,g\left(  3\right)  ,g\left(  5\right)
\right\}  =\left\{  3,1,5\right\}  =\left\{  1,3,5\right\}  .
\]

Note that if $g\in S_{A}$ is a permutation of $A$, then $\left\vert g\left(
U\right)  \right\vert =\left\vert U\right\vert $ for any subset $U$ of $A$.
Thus, we can restrict the action in Example \ref{exa.rep.G-sets.PA} to the
subsets of a given size $k$:

\begin{example}
\label{exa.rep.G-sets.PkA}Let $A$ be any set. Let $k\in\mathbb{N}$. Consider
the set $\mathcal{P}_{k}\left(  A\right)  $ of all $k$-element subsets of $A$.
The symmetric group $S_{A}$ acts on $\mathcal{P}_{k}\left(  A\right)  $ by%
\[
g\rightharpoonup U=g\left(  U\right)  \ \ \ \ \ \ \ \ \ \ \text{for all }g\in
S_{A}\text{ and }U\in\mathcal{P}_{k}\left(  A\right)  .
\]

\end{example}

Here is a slightly more exotic action of a symmetric group:

\begin{example}
\label{exa.rep.G-sets.sign}Let $A$ be a finite set. Then, the symmetric group
$S_{A}$ acts on the two-element set $\left\{  1,-1\right\}  $ by%
\[
g\rightharpoonup x=\left(  -1\right)  ^{g}x\ \ \ \ \ \ \ \ \ \ \text{for all
}g\in S_{A}\text{ and }x\in\left\{  1,-1\right\}  .
\]
In other words, $g\rightharpoonup x=x$ if $g$ is even, and $g\rightharpoonup
x=-x$ if $g$ is odd.
\end{example}

\begin{fineprint}
\begin{proof}
[Proof of Example \ref{exa.rep.G-sets.sign} (sketched).]Consider the map
$\rightharpoonup\ :S_{A}\times\left\{  1,-1\right\}  \rightarrow\left\{
1,-1\right\}  $ defined by%
\[
g\rightharpoonup x=\left(  -1\right)  ^{g}x\ \ \ \ \ \ \ \ \ \ \text{for all
}g\in S_{A}\text{ and }x\in\left\{  1,-1\right\}  .
\]
We must prove that this map is a left $S_{A}$-action on $\left\{
1,-1\right\}  $. In other words, we must prove that it satisfies the two
axioms in Definition \ref{def.rep.G-sets.left} \textbf{(a)}.

We begin with the first axiom. This axiom claims that $g\rightharpoonup\left(
h\rightharpoonup x\right)  =\left(  gh\right)  \rightharpoonup x$ for all
$g,h\in S_{A}$ and $x\in\left\{  1,-1\right\}  $. To prove this, we fix
$g,h\in S_{A}$ and $x\in\left\{  1,-1\right\}  $. Then, one of the basic
properties of signs of permutations shows that $\left(  -1\right)
^{gh}=\left(  -1\right)  ^{g}\cdot\left(  -1\right)  ^{h}$. Now, the
definition of $\rightharpoonup$ yields%
\[
g\rightharpoonup\left(  h\rightharpoonup x\right)  =\left(  -1\right)
^{g}\underbrace{\left(  h\rightharpoonup x\right)  }_{\substack{=\left(
-1\right)  ^{h}x\\\text{(by the definition of }\rightharpoonup\text{)}%
}}=\left(  -1\right)  ^{g}\left(  -1\right)  ^{h}x
\]
and%
\[
\left(  gh\right)  \rightharpoonup x=\underbrace{\left(  -1\right)  ^{gh}%
}_{=\left(  -1\right)  ^{g}\cdot\left(  -1\right)  ^{h}}x=\left(  -1\right)
^{g}\cdot\left(  -1\right)  ^{h}x=\left(  -1\right)  ^{g}\left(  -1\right)
^{h}x.
\]
Comparing these two equalities, we obtain $g\rightharpoonup\left(
h\rightharpoonup x\right)  =\left(  gh\right)  \rightharpoonup x$. Thus, we
have shown that $\rightharpoonup$ satisfies the first axiom in Definition
\ref{def.rep.G-sets.left} \textbf{(a)}.

The proof of the second axiom is similar but even more trivial (we need the
fact that $\left(  -1\right)  ^{\operatorname*{id}}=1$). Thus, both axioms are
satisfied, so that $\rightharpoonup$ is really a left $S_{A}$-action. This
proves Example \ref{exa.rep.G-sets.sign}.
\end{proof}
\end{fineprint}

\begin{example}
\label{exa.rep.G-sets.An}Let $A$ be any set. Let $n\in\mathbb{N}$. Then, the
symmetric group $S_{A}$ acts on the set%
\[
A^{n}=\left\{  \text{all }n\text{-tuples of elements of }A\right\}
\]
by%
\begin{align*}
g  &  \rightharpoonup\left(  a_{1},a_{2},\ldots,a_{n}\right)  =\left(
g\left(  a_{1}\right)  ,g\left(  a_{2}\right)  ,\ldots,g\left(  a_{n}\right)
\right) \\
&  \ \ \ \ \ \ \ \ \ \ \ \ \ \ \ \ \ \ \ \ \text{for all }g\in S_{A}\text{ and
}\left(  a_{1},a_{2},\ldots,a_{n}\right)  \in A^{n}.
\end{align*}
This is called \emph{entrywise action} of $S_{A}$, or the \emph{action on
entries}.
\end{example}

For example, for $A=\left\{  a,b,c,d,e\right\}  $ and $n=3$, we have%
\[
\operatorname*{cyc}\nolimits_{a,b,c}\rightharpoonup\left(  d,c,a\right)
=\left(  d,a,b\right)  .
\]

\begin{example}
\label{exa.rep.G-sets.place}Let $A$ be any set. Let $n\in\mathbb{N}$. Then,
the symmetric group $S_{n}$ acts on the set%
\[
A^{n}=\left\{  \text{all }n\text{-tuples of elements of }A\right\}
\]
by%
\begin{align*}
g  &  \rightharpoonup\left(  a_{1},a_{2},\ldots,a_{n}\right)  =\left(
a_{g^{-1}\left(  1\right)  },a_{g^{-1}\left(  2\right)  },\ldots
,a_{g^{-1}\left(  n\right)  }\right) \\
&  \ \ \ \ \ \ \ \ \ \ \ \ \ \ \ \ \ \ \ \ \text{for all }g\in S_{n}\text{ and
}\left(  a_{1},a_{2},\ldots,a_{n}\right)  \in A^{n}.
\end{align*}
This is called the \emph{place permutation action} of $S_{n}$, or the
\emph{action on places}.
\end{example}

For example, for $A=\left\{  a,b,c,d,e\right\}  $ and $n=5$, we have%
\begin{align*}
\operatorname*{cyc}\nolimits_{1,2,3}  &  \rightharpoonup\left(
a,b,c,d,e\right)  =\left(  c,a,b,d,e\right)  \ \ \ \ \ \ \ \ \ \ \text{and}\\
t_{2,5}  &  \rightharpoonup\left(  a,b,c,d,e\right)  =\left(
a,e,c,d,b\right)  .
\end{align*}

\begin{fineprint}
\begin{proof}
[Proof of Example \ref{exa.rep.G-sets.place}.]Let us prove the claim of
Example \ref{exa.rep.G-sets.place}, to illustrate how such claims can be
shown. The proof is simple and straightforward, but the road is somewhat
slippery due to the many subscripts involved.

We must prove that the map $\rightharpoonup$ introduced in Example
\ref{exa.rep.G-sets.place} really is a left $S_{n}$-action on $A^{n}$. In
other words, we must show that it satisfies the two axioms in Definition
\ref{def.rep.G-sets.left} \textbf{(a)}. The second axiom (unitality) is
satisfied for pretty obvious reasons\footnote{Indeed, $\operatorname*{id}%
\rightharpoonup\left(  a_{1},a_{2},\ldots,a_{n}\right)  =\left(
a_{\operatorname*{id}\nolimits^{-1}\left(  1\right)  },a_{\operatorname*{id}%
\nolimits^{-1}\left(  2\right)  },\ldots,a_{\operatorname*{id}\nolimits^{-1}%
\left(  n\right)  }\right)  =\left(  a_{1},a_{2},\ldots,a_{n}\right)  $ for
all $\left(  a_{1},a_{2},\ldots,a_{n}\right)  \in A^{n}$.}, so we focus on
verifying the first (associativity). Thus, we need to prove that
\[
g\rightharpoonup\left(  h\rightharpoonup x\right)  =\left(  gh\right)
\rightharpoonup x\ \ \ \ \ \ \ \ \ \ \text{for all }g,h\in S_{n}\text{ and
}x\in A^{n}.
\]

For this purpose, we fix $g,h\in S_{n}$ and $x\in A^{n}$. Write $x$ as
$x=\left(  x_{1},x_{2},\ldots,x_{n}\right)  $. Thus,%
\[
h\rightharpoonup x=h\rightharpoonup\left(  x_{1},x_{2},\ldots,x_{n}\right)
=\left(  x_{h^{-1}\left(  1\right)  },x_{h^{-1}\left(  2\right)  }%
,\ldots,x_{h^{-1}\left(  n\right)  }\right)
\]
(by the definition of $\rightharpoonup$). Let us denote this $n$-tuple
$h\rightharpoonup x$ by $\left(  y_{1},y_{2},\ldots,y_{n}\right)  $. Thus,%
\[
\left(  y_{1},y_{2},\ldots,y_{n}\right)  =h\rightharpoonup x=\left(
x_{h^{-1}\left(  1\right)  },x_{h^{-1}\left(  2\right)  },\ldots
,x_{h^{-1}\left(  n\right)  }\right)  .
\]
In other words,%
\begin{equation}
y_{i}=x_{h^{-1}\left(  i\right)  }\ \ \ \ \ \ \ \ \ \ \text{for each }%
i\in\left[  n\right]  . \label{pf.exa.rep.G-sets.place.2}%
\end{equation}

Now, from $h\rightharpoonup x=\left(  y_{1},y_{2},\ldots,y_{n}\right)  $, we
obtain%
\begin{align}
g\rightharpoonup\left(  h\rightharpoonup x\right)   &  =g\rightharpoonup
\left(  y_{1},y_{2},\ldots,y_{n}\right) \nonumber\\
&  =\left(  y_{g^{-1}\left(  1\right)  },y_{g^{-1}\left(  2\right)  }%
,\ldots,y_{g^{-1}\left(  n\right)  }\right)  \label{pf.exa.rep.G-sets.place.3}%
\end{align}
(by the definition of $\rightharpoonup$). However, each $k\in\left[  n\right]
$ satisfies%
\begin{align*}
y_{g^{-1}\left(  k\right)  }  &  =x_{h^{-1}\left(  g^{-1}\left(  k\right)
\right)  }\ \ \ \ \ \ \ \ \ \ \left(  \text{by
(\ref{pf.exa.rep.G-sets.place.2}), applied to }i=g^{-1}\left(  k\right)
\right) \\
&  =x_{\left(  gh\right)  ^{-1}\left(  k\right)  }\ \ \ \ \ \ \ \ \ \ \left(
\text{since }h^{-1}\left(  g^{-1}\left(  k\right)  \right)
=\underbrace{\left(  h^{-1}g^{-1}\right)  }_{=\left(  gh\right)  ^{-1}}\left(
k\right)  =\left(  gh\right)  ^{-1}\left(  k\right)  \right)  .
\end{align*}
Hence,
\[
\left(  y_{g^{-1}\left(  1\right)  },y_{g^{-1}\left(  2\right)  }%
,\ldots,y_{g^{-1}\left(  n\right)  }\right)  =\left(  x_{\left(  gh\right)
^{-1}\left(  1\right)  },x_{\left(  gh\right)  ^{-1}\left(  2\right)  }%
,\ldots,x_{\left(  gh\right)  ^{-1}\left(  n\right)  }\right)  .
\]
Thus, (\ref{pf.exa.rep.G-sets.place.3}) can be rewritten as%
\[
g\rightharpoonup\left(  h\rightharpoonup x\right)  =\left(  x_{\left(
gh\right)  ^{-1}\left(  1\right)  },x_{\left(  gh\right)  ^{-1}\left(
2\right)  },\ldots,x_{\left(  gh\right)  ^{-1}\left(  n\right)  }\right)  .
\]
Comparing this with%
\begin{align*}
\left(  gh\right)  \rightharpoonup x  &  =\left(  gh\right)  \rightharpoonup
\left(  x_{1},x_{2},\ldots,x_{n}\right)  \ \ \ \ \ \ \ \ \ \ \left(
\text{since }x=\left(  x_{1},x_{2},\ldots,x_{n}\right)  \right) \\
&  =\left(  x_{\left(  gh\right)  ^{-1}\left(  1\right)  },x_{\left(
gh\right)  ^{-1}\left(  2\right)  },\ldots,x_{\left(  gh\right)  ^{-1}\left(
n\right)  }\right)  \ \ \ \ \ \ \ \ \ \ \left(  \text{by the definition of
}\rightharpoonup\right)  ,
\end{align*}
we obtain $g\rightharpoonup\left(  h\rightharpoonup x\right)  =\left(
gh\right)  \rightharpoonup x$. Thus, we have shown that the map
$\rightharpoonup$ satisfies the first axiom (associativity) in Definition
\ref{def.rep.G-sets.left} \textbf{(a)}. As we said, this completes the proof
of Example \ref{exa.rep.G-sets.place}.
\end{proof}
\end{fineprint}

Note that Example \ref{exa.rep.G-sets.An} and Example
\ref{exa.rep.G-sets.place} define actions of two different symmetric groups
($S_{A}$ and $S_{n}$) on the same set $A^{n}$. When $A=\left[  n\right]  $,
the two groups agree, but the actions do not.

\begin{example}
\label{exa.rep.G-sets.knm}Let $n,m\in\mathbb{N}$. Consider the set
$\mathbf{k}^{n\times m}$ of $n\times m$-matrices with entries in $\mathbf{k}$.
We use the notation $\left(  a_{i,j}\right)  _{i\in\left[  n\right]
,\ j\in\left[  m\right]  }$ for the $n\times m$-matrix whose $\left(
i,j\right)  $-th entry is $a_{i,j}$ for all $i\in\left[  n\right]  $ and
$j\in\left[  m\right]  $.

Then, the Cartesian product $S_{n}\times S_{m}$ of the symmetric groups
$S_{n}$ and $S_{m}$ acts on $\mathbf{k}^{n\times m}$ by permuting rows and
columns:%
\begin{align*}
&  \left(  g,h\right)  \rightharpoonup\left(  a_{i,j}\right)  _{i\in\left[
n\right]  ,\ j\in\left[  m\right]  }=\left(  a_{g^{-1}\left(  i\right)
,h^{-1}\left(  j\right)  }\right)  _{i\in\left[  n\right]  ,\ j\in\left[
m\right]  }\\
&  \ \ \ \ \ \ \ \ \ \ \text{for all }g\in S_{n}\text{ and }h\in S_{m}\text{
and }\left(  a_{i,j}\right)  _{i\in\left[  n\right]  ,\ j\in\left[  m\right]
}\in\mathbf{k}^{n\times m}\text{.}%
\end{align*}
In other words, for any $\left(  g,h\right)  \in S_{n}\times S_{m}$ and any
matrix $A\in\mathbf{k}^{n\times m}$, the matrix $\left(  g,h\right)
\rightharpoonup A$ is obtained from $A$ by permuting the rows using $g$ and
permuting the columns using $h$.
\end{example}

For instance, for $n=3$ and $m=3$, we have%
\[
\left(  s_{1},\ \operatorname*{cyc}\nolimits_{1,2,3}\right)  \rightharpoonup
\left(
\begin{array}
[c]{ccc}%
a_{1,1} & a_{1,2} & a_{1,3}\\
a_{2,1} & a_{2,2} & a_{2,3}\\
a_{3,1} & a_{3,2} & a_{3,3}%
\end{array}
\right)  =\left(
\begin{array}
[c]{ccc}%
a_{2,3} & a_{2,1} & a_{2,2}\\
a_{1,3} & a_{1,1} & a_{1,2}\\
a_{3,3} & a_{3,1} & a_{3,2}%
\end{array}
\right)  .
\]

Other groups have interesting actions, too:

\begin{example}
\label{exa.rep.G-sets.GL}Let $n\in\mathbb{N}$. The group%
\[
\operatorname*{GL}\nolimits_{n}\left(  \mathbb{R}\right)  =\left\{
\text{invertible }n\times n\text{-matrices with real entries}\right\}
\]
acts on the set $\mathbb{R}^{n\times1}=\left\{  \text{column vectors of size
}n\text{ over }\mathbb{R}\right\}  $ by%
\[
A\rightharpoonup v=Av\ \ \ \ \ \ \ \ \ \ \text{for all }A\in\operatorname*{GL}%
\nolimits_{n}\left(  \mathbb{R}\right)  \text{ and }v\in\mathbb{R}^{n\times
1}.
\]

The same applies to any ring instead of $\mathbb{R}$.
\end{example}

And there are many more (see, e.g., \cite[\textquotedblleft Group
actions\textquotedblright, \S 2]{Conrad}).

\subsubsection{$G$-sets vs. morphisms $G\rightarrow S_{X}$}

You can think of $G$-sets as a non-additive analogue of $\mathbf{k}$-modules
(\textquotedblleft non-additive\textquotedblright\ because neither $G$ nor the
$G$-set are equipped with an addition). The left $G$-action $\rightharpoonup
\ :G\times X\rightarrow X$ for a $G$-set is thus an analogue of the scaling
map $\mathbf{k}\times V\rightarrow V,\ \left(  \lambda,v\right)
\mapsto\lambda v$ that scales vectors in a $\mathbf{k}$-module $V$. The group
elements $g\in G$ here correspond to the scalars $\lambda\in\mathbf{k}$.

Just like for $\mathbf{k}$-modules, we can fix a \textquotedblleft
scalar\textquotedblright\ (really a group element) $g\in G$ and consider the
\textquotedblleft scaling\textquotedblright\ by $g$ as a map from $X$ to $X$:

\begin{definition}
\label{def.rep.G-sets.act}Let $G$ be a group, and let $X$ be a left $G$-set.
Let $g\in G$. Then, the map%
\begin{align*}
\tau_{g}:X  &  \rightarrow X,\\
x  &  \mapsto g\rightharpoonup x
\end{align*}
is called the \emph{action} of $g$ on $X$.
\end{definition}

This map $\tau_{g}$ generalizes the translations $\tau_{n}:\mathbb{Z}%
\rightarrow\mathbb{Z}$ from Example \ref{exa.rep.G-sets.taun} (which are
obtained as a special case for $G=\mathbb{Z}$ and $X=\mathbb{Z}$ and $g=n$,
where $G$ acts on $X$ by the left regular $G$-action defined in Example
\ref{exa.rep.G-sets.lreg}), and the translations $\tau_{g}:G\rightarrow G$
from Theorem \ref{thm.rep.G-sets.cay} (which, again, arise by letting $X$ be
$G$ with the left regular $G$-action). This fulfills the \textquotedblleft
dream\textquotedblright\ of $g$ to act on something.

Let $G$ be a group, and let $X$ be a left $G$-set. The actions of the elements
of $G$ on $X$ behave very much like the elements themselves: For example, the
action of $1_{G}$ is
\begin{equation}
\tau_{1_{G}}=\operatorname*{id}\nolimits_{X} \label{eq.rep.G-sets.curr.0}%
\end{equation}
(since $\tau_{1_{G}}\left(  x\right)  =1\rightharpoonup x=x$ for each $x\in
X$). Moreover, for any two elements $g,h\in G$, the action of the product $gh$
is the composition of the actions of $g$ and $h$: that is, we have%
\begin{equation}
\tau_{gh}=\tau_{g}\tau_{h} \label{eq.rep.G-sets.curr.1}%
\end{equation}
(this is just restating the axiom $\left(  gh\right)  \rightharpoonup
x=g\rightharpoonup\left(  h\rightharpoonup x\right)  $ in Definition
\ref{def.rep.G-sets.left} \textbf{(a)}). This easily yields that the action
$\tau_{g}$ of any element $g\in G$ is invertible, with inverse $\tau_{g^{-1}}$
(since (\ref{eq.rep.G-sets.curr.1}) yields $\tau_{gg^{-1}}=\tau_{g}%
\tau_{g^{-1}}$ and thus $\tau_{g}\tau_{g^{-1}}=\tau_{gg^{-1}}=\tau_{1_{G}%
}=\operatorname*{id}\nolimits_{X}$ (by (\ref{eq.rep.G-sets.curr.0})), and
similarly $\tau_{g^{-1}}\tau_{g}=\operatorname*{id}$). Hence, this action
$\tau_{g}$ is a permutation of $X$, and therefore belongs to the symmetric
group $S_{X}$. We thus obtain a map%
\begin{align*}
\rho:G  &  \rightarrow S_{X},\\
g  &  \mapsto\tau_{g}%
\end{align*}
(which sends each $g\in G$ to its action on $X$). This map $\rho$ is
furthermore a group morphism (by (\ref{eq.rep.G-sets.curr.1}) and
(\ref{eq.rep.G-sets.curr.0})).

\begin{definition}
\label{def.rep.G-sets.curr}We shall call this morphism $\rho$ the
\emph{curried form} of the left $G$-action $\rightharpoonup$ on $X$.
\end{definition}

We note that we can reconstruct the left $G$-action $\rightharpoonup$ on $X$
from this morphism $\rho$, since all $g\in G$ and $x\in X$ satisfy%
\[
g\rightharpoonup x=\underbrace{\tau_{g}}_{=\rho\left(  g\right)  }\left(
x\right)  =\left(  \rho\left(  g\right)  \right)  \left(  x\right)  .
\]
Thus, the morphism $\rho$ encodes the left $G$-action $\rightharpoonup$.

Note that $\rho$ is not always injective. For example, if $X$ has a trivial
$G$-action (as defined in Example \ref{exa.rep.G-sets.triv}), then
$\rho\left(  g\right)  =\operatorname*{id}\nolimits_{X}$ for each $g\in G$.
\medskip

Thus, we have transformed each left $G$-set $X$ into a morphism $\rho
:G\rightarrow S_{X}$.

Let us now forget that we fixed $X$, and work in the opposite direction: If we
are given a set $X$ (without a $G$-action provided) and a group morphism
$\rho:G\rightarrow S_{X}$, then we can define a left $G$-action on $X$ by%
\[
g\rightharpoonup x=\left(  \rho\left(  g\right)  \right)  \left(  x\right)
\ \ \ \ \ \ \ \ \ \ \text{for all }g\in G\text{ and }x\in X.
\]

\begin{definition}
\label{def.rep.G-sets.uncurr}We shall call this $G$-action the \emph{uncurried
form} of $\rho$.
\end{definition}

Altogether, we thus obtain the following theorem:

\begin{theorem}
\label{thm.rep.G-sets.curr}Let $G$ be a group. Let $X$ be a set. Then, a left
$G$-action on $X$ is \textquotedblleft the same as\textquotedblright\ a group
morphism from $G$ to $S_{X}$. To be more precise: \medskip

\textbf{(a)} If $\rightharpoonup$ is a left $G$-action on $X$, then the
curried form of $\rightharpoonup$ is a group morphism from $G$ to $S_{X}$.
\medskip

\textbf{(b)} If $\rho$ is a group morphism from $G$ to $S_{X}$, then the
uncurried form of $\rho$ is a left $G$-action on $X$. \medskip

\textbf{(c)} Currying and uncurrying are mutually inverse: I.e., the uncurried
form of the curried form of a left $G$-action $\rightharpoonup$ is
$\rightharpoonup$ itself, and conversely, the curried form of the uncurried
form of a group morphism $\rho:G\rightarrow S_{X}$ is $\rho$ itself.
\end{theorem}

\begin{proof}
[Proof sketch.]Straightforward. (Part \textbf{(a)} was essentially proved above.)
\end{proof}

\begin{remark}
The words \textquotedblleft curried\textquotedblright\ and \textquotedblleft
uncurried\textquotedblright\ refer to the general notions of \emph{currying}
and \emph{uncurrying} in computer science. Namely, \emph{currying} is the
operation that transforms an arbitrary map $f:A\times B\rightarrow C$ with two
inputs into a map%
\begin{align*}
\widetilde{f}:A  &  \rightarrow\left\{  \text{maps from }B\text{ to
}C\right\}  ,\\
a  &  \mapsto\left(  \text{the map from }B\text{ to }C\text{ that sends each
}b\text{ to }f\left(  a,b\right)  \right)  .
\end{align*}
(In essence, $\widetilde{f}$ is a version of $f$ that takes the two inputs of
$f$ \textquotedblleft one at a time\textquotedblright.) Conversely,
\emph{uncurrying} is the operation that transforms a map $g:A\rightarrow
\left\{  \text{maps from }B\text{ to }C\right\}  $ (whose outputs themselves
are maps) into a two-input map%
\begin{align*}
g^{\prime}:A\times B  &  \rightarrow C,\\
\left(  a,b\right)   &  \mapsto\left(  g\left(  a\right)  \right)  \left(
b\right)  .
\end{align*}
It is easy to see that currying and uncurrying are mutually inverse
operations. The curried form of a left $G$-action $\rightharpoonup$ is exactly
the result of currying $\rightharpoonup$, except that its target is not
$\left\{  \text{maps from }X\text{ to }X\right\}  $ but rather the symmetric
group $S_{X}$. Likewise, the uncurried form of a group morphism $\rho
:G\rightarrow S_{X}$ is the result of uncurrying $\rho$, where we treat $\rho$
as a map from $G$ to $\left\{  \text{maps from }X\text{ to }X\right\}  $.
\end{remark}

Theorem \ref{thm.rep.G-sets.curr} can be extended to the case when $G$ is just
a monoid (rather than a group), provided that we replace the symmetric group
$S_{X}$ by the \textquotedblleft symmetric monoid\textquotedblright\ $\left\{
\text{all maps from }X\text{ to }X\right\}  $.

\begin{example}
\label{exa.rep.G-sets.curr.lreg}Consider the left regular $G$-action (as
defined in Example \ref{exa.rep.G-sets.lreg}). Its curried form is just the
morphism $G\rightarrow S_{G}$ from Cayley's theorem (Theorem
\ref{thm.rep.G-sets.cay}).
\end{example}

\begin{example}
Let $A$ be any set. Consider the natural action of $S_{A}$ on $A$ (as defined
in Example \ref{exa.rep.G-sets.nat}). Its curried form is the identity
morphism $S_{A}\rightarrow S_{A}$.
\end{example}

\begin{example}
\label{exa.rep.G-sets.curr.sign}Let $A$ be a finite set. Consider the left
$S_{A}$-action on $\left\{  1,-1\right\}  $ given by%
\[
g\rightharpoonup x=\left(  -1\right)  ^{g}x\ \ \ \ \ \ \ \ \ \ \text{for all
}g\in S_{A}\text{ and }x\in\left\{  1,-1\right\}  .
\]
The curried form of this left $S_{A}$-action $\rightharpoonup$ is the map%
\begin{align*}
S_{A}  &  \rightarrow S_{\left\{  1,-1\right\}  },\\
g  &  \mapsto\left(  \text{multiplication by }\left(  -1\right)  ^{g}\right)
.
\end{align*}
This is \textquotedblleft more or less\textquotedblright\ the sign
homomorphism%
\begin{align*}
S_{A}  &  \rightarrow\left\{  1,-1\right\}  ,\\
g  &  \mapsto\left(  -1\right)  ^{g},
\end{align*}
except that instead of the actual numbers $1$ and $-1$ it produces the
multiplications by these numbers (but this is the same thing up to isomorphism).
\end{example}

\subsubsection{$G$-equivariant maps aka $G$-set morphisms}

Whenever you define an algebraic object, the next logical step is to define
morphisms between these objects. For $G$-sets, it is pretty clear how to do this:

\begin{definition}
\label{def.rep.G-sets.mor}Let $X$ and $Y$ be two left $G$-sets. Then, a
$G$\emph{-set morphism} (aka a \emph{morphism of }$G$\emph{-sets}, aka a
$G$\emph{-equivariant map}) means a map $f:X\rightarrow Y$ with the property
that all $g\in G$ and $x\in X$ satisfy%
\begin{equation}
f\left(  g\rightharpoonup x\right)  =g\rightharpoonup f\left(  x\right)  .
\label{eq.def.rep.G-sets.mor.req}%
\end{equation}

\end{definition}

This requirement is similar to the axiom $f\left(  \lambda x\right)  =\lambda
f\left(  x\right)  $ in the definition of a $\mathbf{k}$-module morphism (aka
$\mathbf{k}$-linear map).

\begin{example}
\label{exa.rep.G-sets.mor-Sn}Let $A$ be any set. Let $n\in\mathbb{N}$.
Consider $A^{n}$ as a left $S_{n}$-set via the place permutation action of
$S_{n}$ on $A^{n}$ (defined in Example \ref{exa.rep.G-sets.place}). Also
consider $A$ as a left $S_{n}$-set via the trivial $S_{n}$-action on $A$ (as
defined in Example \ref{exa.rep.G-sets.triv}, i.e., defined by
$g\rightharpoonup x=x$ for all $g\in S_{n}$ and $x\in A$). Then: \medskip

\textbf{(a)} The map
\begin{align*}
f_{1}:A  &  \rightarrow A^{n},\\
a  &  \mapsto\left(  \underbrace{a,a,\ldots,a}_{n\text{ times}}\right)
\end{align*}
is an $S_{n}$-set morphism. Indeed, all $g\in S_{n}$ and $x\in A$ satisfy%
\begin{align*}
f_{1}\left(  \underbrace{g\rightharpoonup x}_{=x}\right)   &  =f_{1}\left(
x\right)  =\left(  \underbrace{x,x,\ldots,x}_{n\text{ times}}\right)
\ \ \ \ \ \ \ \ \ \ \text{and}\\
g\rightharpoonup f_{1}\left(  x\right)   &  =g\rightharpoonup\left(
\underbrace{x,x,\ldots,x}_{n\text{ times}}\right)  =\left(
\underbrace{x,x,\ldots,x}_{n\text{ times}}\right)
\end{align*}
(by the definition of $\rightharpoonup$) and therefore (by comparing these two
equalities) $f_{1}\left(  g\rightharpoonup x\right)  =g\rightharpoonup
f_{1}\left(  x\right)  $, so that the requirement
(\ref{eq.def.rep.G-sets.mor.req}) is satisfied for $f=f_{1}$. \medskip

\textbf{(b)} Assume that $A$ is an additive abelian monoid (i.e., a set
equipped with a binary operation $+$ that is associative and commutative and
has a neutral element). For instance, $A$ can be $\mathbb{N}$ or any ring like
$\mathbb{Z}$ or $\mathbb{Q}$.

Then, the map%
\begin{align*}
f_{2}:A^{n}  &  \rightarrow A,\\
\left(  a_{1},a_{2},\ldots,a_{n}\right)   &  \mapsto a_{1}+a_{2}+\cdots+a_{n}%
\end{align*}
is also an $S_{n}$-set morphism. Indeed, the requirement
(\ref{eq.def.rep.G-sets.mor.req}) (applied to $x=\left(  a_{1},a_{2}%
,\ldots,a_{n}\right)  $) is here just saying that%
\[
a_{g^{-1}\left(  1\right)  }+a_{g^{-1}\left(  2\right)  }+\cdots
+a_{g^{-1}\left(  n\right)  }=a_{1}+a_{2}+\cdots+a_{n}%
\]
for all $g\in S_{n}$ and $\left(  a_{1},a_{2},\ldots,a_{n}\right)  \in A^{n}$;
but this is a known property of abelian monoids (\textquotedblleft generalized
commutativity\textquotedblright). \medskip

\textbf{(c)} The map%
\begin{align*}
f_{3}:A^{n}  &  \rightarrow A,\\
\left(  a_{1},a_{2},\ldots,a_{n}\right)   &  \mapsto a_{1}%
\end{align*}
is \textbf{not} an $S_{n}$-set morphism (unless $n\leq1$ or $\left\vert
A\right\vert \leq1$). Indeed, if we pick any $x=\left(  x_{1},x_{2}%
,\ldots,x_{n}\right)  \in A^{n}$ with $x_{1}\neq x_{2}$, then we have
$f_{3}\left(  s_{1}\rightharpoonup x\right)  =f_{3}\left(  x_{2},x_{1}%
,x_{3},x_{4},\ldots,x_{n}\right)  =x_{2}$ but $s_{1}\rightharpoonup
f_{3}\left(  x\right)  =f_{3}\left(  x\right)  =x_{1}\neq x_{2}$. \medskip

\textbf{(d)} The map%
\begin{align*}
f_{4}:A^{n}  &  \rightarrow A^{n},\\
\left(  a_{1},a_{2},\ldots,a_{n}\right)   &  \mapsto\left(  a_{n}%
,a_{n-1},\ldots,a_{1}\right)
\end{align*}
(which reverses each $n$-tuple given to it as an input) is not an $S_{n}$-set
morphism (unless $n\leq2$ or $\left\vert A\right\vert \leq1$). Indeed, if we
pick any $x=\left(  x_{1},x_{2},\ldots,x_{n}\right)  \in A^{n}$ with
$x_{1}\neq x_{2}$, then we have
\begin{align*}
f_{4}\left(  s_{1}\rightharpoonup x\right)   &  =f_{4}\left(  x_{2}%
,x_{1},x_{3},x_{4},\ldots,x_{n}\right)  =\left(  x_{n},x_{n-1},\ldots
,x_{3},x_{1},x_{2}\right)  \ \ \ \ \ \ \ \ \ \ \text{but}\\
s_{1}\rightharpoonup f_{4}\left(  x\right)   &  =s_{1}\rightharpoonup\left(
x_{n},x_{n-1},\ldots,x_{1}\right)  =\left(  x_{n-1},x_{n},x_{n-2}%
,x_{n-3},\ldots,x_{1}\right)  ,
\end{align*}
which are different.
\end{example}

\begin{example}
\label{exa.rep.G-sets.mor-triv}Let $X$ and $Y$ be two left $G$-sets with
trivial $G$-actions (i.e., we assume that $g\rightharpoonup x=x$ for all $g\in
G$ and $x\in X$, and that $g\rightharpoonup y=y$ for all $g\in G$ and $y\in
Y$). Then, \textbf{any} map $f:X\rightarrow Y$ is $G$-equivariant.
\end{example}

\begin{fineprint}
\begin{proof}
Let $f:X\rightarrow Y$ be any map. For each $g\in G$ and each $x\in X$, we
then have $g\rightharpoonup f\left(  x\right)  =f\left(  x\right)  $ (since
the $G$-action on $Y$ is trivial) and $g\rightharpoonup x=x$ (since the
$G$-action on $X$ is trivial) and thus
\[
f\left(  g\rightharpoonup x\right)  =f\left(  x\right)  =g\rightharpoonup
f\left(  x\right)  .
\]
But this means precisely that $f$ is $G$-equivariant. This proves Example
\ref{exa.rep.G-sets.mor-triv}.
\end{proof}
\end{fineprint}

As usual, the nicest morphisms are isomorphisms:

\begin{definition}
\label{def.rep.G-sets.iso}Let $X$ and $Y$ be two left $G$-sets. Then, a
$G$\emph{-set isomorphism} (aka \emph{isomorphism of }$G$\emph{-sets}) means a
$G$-set morphism $f:X\rightarrow Y$ whose inverse exists and is again a
$G$-set morphism.
\end{definition}

\begin{proposition}
\label{prop.rep.G-sets.iso=bij}The \textquotedblleft is again a $G$-set
morphism\textquotedblright\ part in Definition \ref{def.rep.G-sets.iso} is
redundant (i.e., follows from the other requirements). In other words, if a
$G$-set morphism $f:X\rightarrow Y$ is invertible, then its inverse is again a
$G$-set morphism.
\end{proposition}

\begin{proof}
This is proved in the same straightforward way as the analogous claims about
groups, rings, $\mathbf{k}$-modules and many other algebraic structures.
\end{proof}

\begin{definition}
\label{def.rep.G-sets.isomorphic}Let $X$ and $Y$ be two left $G$-sets. Then,
we say that $X$ and $Y$ are \emph{isomorphic} (this is written $X\cong Y$) if
there exists a $G$-set isomorphism from $X$ to $Y$.
\end{definition}

Here are some examples of $G$-set isomorphisms:

\begin{example}
We have $G/G\cong\left\{  1\right\}  $ as left $G$-sets (where $G/G$ is a left
$G$-set as defined in Example \ref{exa.rep.G-sets.G/H}, whereas $\left\{
1\right\}  $ has a trivial $G$-action as defined in Example
\ref{exa.rep.G-sets.triv}). Indeed, there is a $G$-set isomorphism%
\begin{align*}
G/G  &  \rightarrow\left\{  1\right\}  ,\\
gG  &  \mapsto1.
\end{align*}

\end{example}

\begin{example}
We have $G/\left\{  1_{G}\right\}  \cong G$ as left $G$-sets (where
$G/\left\{  1_{G}\right\}  $ is a left $G$-set as defined in Example
\ref{exa.rep.G-sets.G/H}, whereas $G$ has the left regular $G$-action as
defined in Example \ref{exa.rep.G-sets.lreg}). Indeed, there is a $G$-set
isomorphism%
\begin{align*}
G/\left\{  1_{G}\right\}   &  \rightarrow G,\\
g\left\{  1_{G}\right\}   &  \mapsto g.
\end{align*}

\end{example}

\begin{example}
Let $A$ be any set. Consider the left $S_{A}$-set $\mathcal{P}\left(
A\right)  $ defined in Example \ref{exa.rep.G-sets.PA}. Then, the map%
\begin{align*}
\mathcal{P}\left(  A\right)   &  \rightarrow\mathcal{P}\left(  A\right)  ,\\
U  &  \mapsto A\setminus U
\end{align*}
(which sends each subset of $A$ to its complement) is an $S_{A}$-set isomorphism.
\end{example}

\begin{example}
Let $A$ be any set. Consider the left $S_{A}$-set $\mathcal{P}_{k}\left(
A\right)  $ defined in Example \ref{exa.rep.G-sets.PkA} for each
$k\in\mathbb{N}$. Also, consider the left $S_{A}$-set $A$ defined in Example
\ref{exa.rep.G-sets.nat}. Then, we have $\mathcal{P}_{1}\left(  A\right)
\cong A$ as $S_{A}$-sets. To be more specific: the map%
\begin{align*}
A  &  \rightarrow\mathcal{P}_{1}\left(  A\right)  ,\\
a  &  \mapsto\left\{  a\right\}
\end{align*}
is an $S_{A}$-set isomorphism.
\end{example}

\begin{example}
Consider the left $S_{n}$-action on $\mathcal{P}\left(  \left[  n\right]
\right)  $ defined in Example \ref{exa.rep.G-sets.PA} (for $A=\left[
n\right]  $), and consider the left $S_{n}$-action on $\left\{  0,1\right\}
^{n}$ defined as in Example \ref{exa.rep.G-sets.place}. Then, we have
$\mathcal{P}\left(  \left[  n\right]  \right)  \cong\left\{  0,1\right\}
^{n}$ as left $S_{n}$-sets. Indeed, there is an $S_{n}$-set isomorphism%
\begin{align*}
\left\{  0,1\right\}  ^{n}  &  \rightarrow\mathcal{P}\left(  \left[  n\right]
\right)  ,\\
\left(  a_{1},a_{2},\ldots,a_{n}\right)   &  \mapsto\left\{  i\in\left[
n\right]  \ \mid\ a_{i}=1\right\}  .
\end{align*}

\end{example}

\begin{exercise}
\fbox{1} Prove the latter claim about the map $\left\{  0,1\right\}
^{n}\rightarrow\mathcal{P}\left(  \left[  n\right]  \right)  $. You can take
for granted the well-known fact that this map is well-defined and bijective.
\end{exercise}

\subsubsection{Combining $G$-sets}

We can think of $G$-sets as a set-theoretical analogue of vector spaces (with
the map $\rightharpoonup\ :G\times X\rightarrow X$ being analogous to the
scaling map $\mathbf{k}\times V\rightarrow V$ on a $\mathbf{k}$-vector space).
Two vector spaces can be combined in two important ways: tensor products and
direct products. The analogues of these two combinations for $G$-sets are
Cartesian products and disjoint unions. \medskip

We begin with the Cartesian products, since they are the easiest to define
(even though the disjoint unions are easier to visualize):

\begin{definition}
\label{def.rep.G-sets.cartprod}\textbf{(a)} If $X$ and $Y$ are two left
$G$-sets, then the Cartesian product $X\times Y$ also becomes a left $G$-set
by setting%
\[
g\rightharpoonup\left(  x,y\right)  =\left(  g\rightharpoonup
x,\ g\rightharpoonup y\right)  \ \ \ \ \ \ \ \ \ \ \text{for all }g\in G\text{
and }x\in X\text{ and }y\in Y.
\]

\textbf{(b)} More generally, for any left $G$-sets $X_{1},X_{2},\ldots,X_{n}$,
we make the Cartesian product $X_{1}\times X_{2}\times\cdots\times X_{n}$ into
a left $G$-set by setting%
\begin{align*}
&  g\rightharpoonup\left(  x_{1},x_{2},\ldots,x_{n}\right)  =\left(
g\rightharpoonup x_{1},\ g\rightharpoonup x_{2},\ \ldots,\ g\rightharpoonup
x_{n}\right) \\
&  \ \ \ \ \ \ \ \ \ \ \text{for all }g\in G\text{ and }\left(  x_{1}%
,x_{2},\ldots,x_{n}\right)  \in X_{1}\times X_{2}\times\cdots\times X_{n}.
\end{align*}
This $G$-set is called the \emph{Cartesian product} (or \emph{direct product})
of these left $G$-sets $X_{1},X_{2},\ldots,X_{n}$. The $G$-action on this
$G$-set is called the \emph{entrywise action} (since $g$ acts on each entry of
a $n$-tuple independently). \medskip

\textbf{(c)} For any left $G$-set $X$ and any $n\in\mathbb{N}$, we let $X^{n}$
denote the Cartesian product $\underbrace{X\times X\times\cdots\times
X}_{n\text{ times}}$. This is called the $n$\emph{-th Cartesian power} of $X$.
\end{definition}

\begin{example}
Let $A$ be any set, and let $n\in\mathbb{N}$. As we saw in Example
\ref{exa.rep.G-sets.nat}, the set $A$ becomes a left $S_{A}$-set using the
natural action (given by $g\rightharpoonup a=g\left(  a\right)  $ for all
$g\in S_{A}$ and $a\in A$). The $n$-th Cartesian power $A^{n}$ of this $S_{A}%
$-set $A$ (as defined in Definition \ref{def.rep.G-sets.cartprod}
\textbf{(c)}) is precisely the $S_{A}$-set $A^{n}$ introduced in Example
\ref{exa.rep.G-sets.An} (i.e., the one where a permutation $g\in S_{A}$ is
applied to each entry of the $n$-tuple independently).
\end{example}

Now let us move on to the disjoint unions.

Recall that the \emph{disjoint union} $X\sqcup Y$ of two (not necessarily
disjoint) sets $X$ and $Y$ is defined by%
\[
X\sqcup Y:=\left(  X\times\left\{  1\right\}  \right)  \cup\left(
Y\times\left\{  2\right\}  \right)  .
\]
Thus, it contains a \textquotedblleft clone\textquotedblright\ $\left(
x,1\right)  $ of each element $x\in X$ as well as a \textquotedblleft
clone\textquotedblright\ $\left(  y,2\right)  $ of each element $y\in Y$. The
\textquotedblleft clones\textquotedblright\ are pairs, whose second entries
$1$ and $2$ have been chosen to disambiguate the \textquotedblleft
clones\textquotedblright\ coming from $X$ from the \textquotedblleft
clones\textquotedblright\ coming from $Y$ (so that, even if the sets $X$ and
$Y$ have an element $z$ in common, the two \textquotedblleft
clones\textquotedblright\ $\left(  z,1\right)  $ and $\left(  z,2\right)  $
will be distinguishable in $X\sqcup Y$). When the sets $X$ and $Y$ are
themselves disjoint, we can (by abuse of notation) identify the clones
$\left(  x,1\right)  $ and $\left(  y,2\right)  $ with the elements $x$ and
$y$ themselves, which allows us to pretend that $X\sqcup Y$ is the literal
union $X\cup Y$ in this case.

More generally, the \emph{disjoint union} $X_{1}\sqcup X_{2}\sqcup\cdots\sqcup
X_{n}$ of $n$ (not necessarily disjoint) sets $X_{1},X_{2},\ldots,X_{n}$ is
defined by%
\begin{align*}
X_{1}\sqcup X_{2}\sqcup\cdots\sqcup X_{n}:=  &  \ \left(  X_{1}\times\left\{
1\right\}  \right)  \cup\left(  X_{2}\times\left\{  2\right\}  \right)
\cup\cdots\cup\left(  X_{n}\times\left\{  n\right\}  \right) \\
=  &  \ \left\{  \left(  x,i\right)  \ \mid\ i\in\left[  n\right]  \text{ and
}x\in X_{i}\right\}  .
\end{align*}
It contains a \textquotedblleft clone\textquotedblright\ $\left(  x,i\right)
$ of each element $x\in X_{i}$ for each $i\in\left[  n\right]  $, and again
the second entries of these \textquotedblleft clones\textquotedblright\ allow
us to distinguish \textquotedblleft clones\textquotedblright\ coming from
different sets.

We can equip disjoint unions of $G$-sets with a canonical $G$-set structure:

\begin{definition}
\label{def.rep.G-sets.djun}\textbf{(a)} If $X$ and $Y$ are two left $G$-sets,
then the disjoint union%
\[
X\sqcup Y=\left(  X\times\left\{  1\right\}  \right)  \cup\left(
Y\times\left\{  2\right\}  \right)
\]
also becomes a left $G$-set by setting%
\begin{align*}
g  &  \rightharpoonup\left(  x,1\right)  =\left(  \underbrace{g\rightharpoonup
x}_{G\text{-action on }X},\ 1\right)  \ \ \ \ \ \ \ \ \ \ \text{for all }g\in
G\text{ and }x\in X,\ \ \ \ \ \ \ \ \ \ \text{and}\\
g  &  \rightharpoonup\left(  y,2\right)  =\left(  \underbrace{g\rightharpoonup
y}_{G\text{-action on }Y},\ 2\right)  \ \ \ \ \ \ \ \ \ \ \text{for all }g\in
G\text{ and }y\in Y.
\end{align*}

\textbf{(b)} More generally, for any left $G$-sets $X_{1},X_{2},\ldots,X_{n}$,
we make the disjoint union $X_{1}\sqcup X_{2}\sqcup\cdots\sqcup X_{n}$ into a
left $G$-set by setting%
\[
g\rightharpoonup\left(  x,i\right)  =\left(  \underbrace{g\rightharpoonup
x}_{G\text{-action on }X_{i}},\ i\right)  \ \ \ \ \ \ \ \ \ \ \text{for all
}g\in G\text{, }i\in\left[  n\right]  \text{ and }x\in X_{i}.
\]
This $G$-set is called the \emph{disjoint union} (or \emph{direct sum}) of
these left $G$-sets $X_{1},X_{2},\ldots,X_{n}$. \medskip

\textbf{(c)} For any left $G$-set $X$ and any $n\in\mathbb{N}$, we let
$X^{\sqcup n}$ denote the disjoint union $\underbrace{X\sqcup X\sqcup
\cdots\sqcup X}_{n\text{ times}}$. This is called the $n$\emph{-th multiple}
of $X$.
\end{definition}

We note a theorem that we will not need, but is worth knowing:

\begin{theorem}
\label{thm.rep.G-sets.decomp}Any finite left $G$-set is isomorphic to a
disjoint union of $G/H$'s corresponding to certain subgroups $H$ of $G$. (Any
such $G/H$ can be used multiple times in the union.)
\end{theorem}

\begin{proof}
[Proof sketch.]Combine \S 1 and Theorem 3.6 of \cite[\textquotedblleft
Transitive group actions\textquotedblright]{Conrad}. Specifically, the former
shows that any finite left $G$-set is isomorphic to a disjoint union of
transitive $G$-sets, whereas the latter reveals that each transitive $G$-set
is isomorphic to $G/H$ for some subgroup $H$ of $G$.
\end{proof}

This is also true for infinite $G$-sets, if we define infinite disjoint unions
in the obvious way.

\subsubsection{$G$-subsets and quotients}

Just like $\mathbf{k}$-modules have $\mathbf{k}$-submodules (and likewise for
groups, rings, etc.), $G$-sets have subobjects called \textquotedblleft%
$G$-subsets\textquotedblright:

\begin{definition}
\label{def.rep.G-sets.sub}Let $X$ be a left $G$-set. \medskip

\textbf{(a)} A $G$\emph{-subset} (aka \emph{sub-}$G$\emph{-set}) of $X$ means
a subset $Y$ of $X$ such that%
\[
\text{all }g\in G\text{ and }y\in Y\text{ satisfy }g\rightharpoonup y\in Y.
\]

\textbf{(b)} If $Y$ is a $G$-subset of $X$, then $Y$ becomes a left $G$-set in
its own right, by simply restricting the left $G$-action on $X$ to $Y$. (That
is, the left $G$-action on $Y$ is just the restriction of the map
$\rightharpoonup\ :G\times X\rightarrow X$ to $G\times Y$.) Clearly, the
canonical inclusion map $i:Y\rightarrow X$ (sending each $y\in Y$ to $y$) is a
$G$-set morphism.
\end{definition}

\begin{example}
Let $A$ be any set. Consider $A^{3}$ as a left $S_{3}$-set via the place
permutation action of $S_{3}$ on $A^{3}$ (defined in Example
\ref{exa.rep.G-sets.place}, for $n=3$). Then, the subset%
\[
\left\{  \left(  a,b,c\right)  \in A^{3}\ \mid\ a=b=c\right\}  =\left\{
\left(  a,a,a\right)  \ \mid\ a\in A\right\}
\]
of $A^{3}$ is an $S_{3}$-subset of $A^{3}$. Moreover, if $A$ is an additive
abelian group, then the subset%
\[
\left\{  \left(  a,b,c\right)  \in A^{3}\ \mid\ a+b+c=0\right\}
\]
of $A^{3}$ is an $S_{3}$-subset of $A^{3}$ as well (since the sum $a+b+c$ does
not change when we permute the entries of the triple $\left(  a,b,c\right)
$). But the subset%
\[
\left\{  \left(  a,b,c\right)  \in A^{3}\ \mid\ a\neq b\text{ and }b\neq
c\right\}
\]
is not an $S_{3}$-subset of $A^{3}$ (unless $\left\vert A\right\vert \leq1$),
since (for any two distinct elements $u$ and $v$ of $A$) the triple $\left(
u,v,u\right)  $ belongs to it but the triple $s_{1}\rightharpoonup\left(
u,v,u\right)  =\left(  v,u,u\right)  $ does not. The reader can easily
generalize these examples to $A^{n}$ instead of $A^{3}$.
\end{example}

Just like $\mathbf{k}$-modules have quotient modules, $G$-sets have quotient
$G$-sets. However, the quotient is not a quotient by a $G$-subset, but rather
a quotient by a \textquotedblleft$G$-invariant equivalence
relation\textquotedblright. To explain this, let me first recall how quotients
by equivalence relations are defined for usual sets (not $G$-sets): If $\sim$
is an equivalence relation on a given set $X$, then the \emph{quotient set} of
$X$ by this relation $\sim$ is defined to be the set of all $\sim$-equivalence
classes on $X$. This quotient set is denoted by $X/\left.  \sim\right.  $. For
each $x\in X$, we denote the $\sim$-equivalence class that contains $x$ as
$\left[  x\right]  _{\sim}$. (Thus, $X/\left.  \sim\right.  =\left\{  \left[
x\right]  _{\sim}\ \mid\ x\in X\right\}  $.) For example, if $n$ is an
integer, then the residue classes of integers modulo $n$ are the equivalence
classes of the relation \textquotedblleft congruent modulo $n$%
\textquotedblright, and thus their set $\mathbb{Z}/n$ is the quotient set
$\mathbb{Z}/\left.  \overset{n}{\equiv}\right.  $, where $\overset{n}{\equiv}$
is this relation.

Now, we can define quotient $G$-sets:

\begin{definition}
\label{def.rep.G-sets.eqrel}Let $\sim$ be an equivalence relation on a left
$G$-set $X$. \medskip

\textbf{(a)} We say that this relation $\sim$ is $G$\emph{-invariant} if it
has the property that for every $x,y\in X$ and every $g\in G$, the implication%
\[
\left(  x\sim y\right)  \ \Longrightarrow\ \left(  g\rightharpoonup x\sim
g\rightharpoonup y\right)  \ \ \ \ \ \ \ \ \ \ \text{holds.}%
\]
In words: The relation $\sim$ is $G$-invariant if and only if any two $\sim
$-equivalent elements of $X$ remain $\sim$-equivalent to each other under the
action of any $g\in G$. (This does \textbf{not} mean that $x\sim
g\rightharpoonup x$ for all $x\in X$.) \medskip

\textbf{(b)} If $\sim$ is a $G$-invariant equivalence relation on $X$, then
the quotient set%
\[
X/\left.  \sim\right.  \ =\left\{  \text{all }\sim\text{-equivalence classes
on }X\right\}
\]
becomes a left $G$-set, with the $G$-action defined by%
\[
g\rightharpoonup U=g\left(  U\right)  \ \ \ \ \ \ \ \ \ \ \text{for any }%
\sim\text{-equivalence class }U.
\]
Equivalently, we can rewrite this definition of $\rightharpoonup$ as%
\[
g\rightharpoonup\left[  x\right]  _{\sim}=\left[  g\rightharpoonup x\right]
_{\sim}\ \ \ \ \ \ \ \ \ \ \text{for any }\sim\text{-equivalence class
}\left[  x\right]  _{\sim}.
\]
The map%
\begin{align*}
X  &  \rightarrow X/\left.  \sim\right.  ,\\
x  &  \mapsto\left[  x\right]  _{\sim}%
\end{align*}
(which sends each element of $X$ to its $\sim$-equivalence class) is thus a
$G$-set morphism.
\end{definition}

\begin{exercise}
\fbox{1} Prove all claims implicitly made in Definition
\ref{def.rep.G-sets.eqrel} (i.e., that $X/\left.  \sim\right.  $ really
becomes a left $G$-set, and that the map is a $G$-set morphism).
\end{exercise}

\begin{example}
Let $A$ be any set. Let $n\in\mathbb{N}$. Consider $A^{n}$ as a left $S_{A}%
$-set via the entrywise action (defined in Example \ref{exa.rep.G-sets.An}).
We can define an equivalence relation $\sim$ on $A^{n}$ by setting%
\begin{align*}
\left(  a_{1},a_{2},\ldots,a_{n}\right)   &  \sim\left(  a_{\sigma\left(
1\right)  },a_{\sigma\left(  2\right)  },\ldots,a_{\sigma\left(  n\right)
}\right) \\
&  \ \ \ \ \ \ \ \ \ \ \text{for all }\left(  a_{1},a_{2},\ldots,a_{n}\right)
\in A^{n}\text{ and }\sigma\in S_{n}.
\end{align*}
In other words, we declare two $n$-tuples in $A^{n}$ to be equivalent if they
can be obtained from each other by permuting the entries. (For example,
$\left(  2,4,2,5\right)  \sim\left(  4,5,2,2\right)  $ but not $\sim\left(
1,4,1,5\right)  $.)

It is easy to see that this relation $\sim$ is $S_{A}$-invariant, so that the
quotient set $A^{n}/\left.  \sim\right.  $ becomes an $S_{A}$-set. The
elements of this quotient set $A^{n}/\left.  \sim\right.  $ are $\sim
$-equivalence classes, i.e., $n$-tuples \textquotedblleft up to
permutation\textquotedblright. They are known as \emph{unordered }%
$n$\emph{-tuples} of elements of $A$. Thus, we have obtained an $S_{A}$-action
on the set of unordered $n$-tuples of elements of $A$.
\end{example}

\subsubsection{Restriction}

Any $\mathbb{C}$-vector space is automatically an $\mathbb{R}$-vector space,
since $\mathbb{R}$ is a subfield of $\mathbb{C}$. More generally, if
$\mathbf{l}$ is a subring of $\mathbf{k}$, then any $\mathbf{k}$-module
automatically becomes an $\mathbf{l}$-module. Something analogous holds for
$G$-sets and subgroups:

\begin{definition}
\label{def.rep.G-sets.restr}Let $H$ be a subgroup of the group $G$. Then,
every left $G$-set $X$ automatically becomes a left $H$-set by setting%
\[
h\rightharpoonup x=h\rightharpoonup x\ \ \ \ \ \ \ \ \ \ \text{for each }h\in
H\text{ and }x\in X
\]
(where the $\rightharpoonup$ on the right hand side is the left $G$-action
that we are given, whereas the $\rightharpoonup$ on the left hand side is the
left $H$-action that we are defining). This left $H$-set is called the
\emph{restriction} of the original $G$-set to the subgroup $H$.
\end{definition}

This notion of restriction can be generalized:

\begin{definition}
\label{def.rep.G-sets.restr-mor}Let $G$ and $H$ be any two groups. Let
$\alpha:H\rightarrow G$ be a group morphism. Then, every left $G$-set $X$
automatically becomes a left $H$-set by setting%
\[
h\rightharpoonup x=\alpha\left(  h\right)  \rightharpoonup
x\ \ \ \ \ \ \ \ \ \ \text{for each }h\in H\text{ and }x\in X.
\]
This left $H$-set is called the \emph{restriction} of the original $G$-set
along the morphism $\alpha$.
\end{definition}

When $H$ is a subgroup of $G$, the canonical inclusion map $i:H\rightarrow
G,\ h\mapsto h$ is a group morphism, and then the restriction of a left
$G$-set along this morphism $i$ is precisely its restriction to the subgroup
$H$ (as defined in Definition \ref{def.rep.G-sets.restr}). Thus, Definition
\ref{def.rep.G-sets.restr-mor} generalizes Definition
\ref{def.rep.G-sets.restr}.

\begin{example}
Let $X$ be any set. The trivial $G$-action on $X$ (that is, the left
$G$-action defined by $g\rightharpoonup x=x$, as in Example
\ref{exa.rep.G-sets.triv}) is the restriction of the trivial $\left\{
1\right\}  $-action on $X$ along the morphism%
\begin{align*}
\pi:G  &  \rightarrow\left\{  1\right\}  ,\\
g  &  \mapsto1.
\end{align*}

\end{example}

\begin{exercise}
\label{exe.rep.G-sets.l2r}\fbox{2} Let $G$ be a group. Let $X$ be a set. Prove
the following: \medskip

\textbf{(a)} If $\rightharpoonup\ :G\times X\rightarrow X$ is a left
$G$-action on $X$, then the map $\leftharpoonup\ :X\times G\rightarrow X$
defined by%
\[
x\leftharpoonup g=g^{-1}\rightharpoonup x\ \ \ \ \ \ \ \ \ \ \text{for all
}g\in G\text{ and }x\in X
\]
is a right $G$-action on $X$. Thus, any left $G$-action on $X$ can be
converted into a right $G$-action on $X$. \medskip

\textbf{(b)} If $\leftharpoonup\ :X\times G\rightarrow X$ is a right
$G$-action on $X$, then the map $\rightharpoonup\ :G\times X\rightarrow X$
defined by%
\[
g\rightharpoonup x=x\leftharpoonup g^{-1}\ \ \ \ \ \ \ \ \ \ \text{for all
}g\in G\text{ and }x\in X
\]
is a left $G$-action on $X$. Thus, any right $G$-action on $X$ can be
converted into a left $G$-action on $X$. \medskip

\textbf{(c)} These two conversions are mutually inverse: I.e.,

\begin{itemize}
\item if we convert a left $G$-action into a right $G$-action and then convert
the resulting right $G$-action back into a left $G$-action, then we recover
the original left $G$-action,

\item and likewise for right $G$-actions.
\end{itemize}
\end{exercise}

\subsubsection{Understanding $G$-actions for cyclic groups $G$}

Recall that a group is called \emph{cyclic} if it can be generated by a single
element. From basic abstract algebra, we know what the cyclic groups are (up
to isomorphism):

\begin{itemize}
\item For each positive integer $n$, there is a cyclic group $C_{n}$ of size
$n$. This group is generated by an element $g$ that satisfies $g^{n}=1$; thus
it consists of $n$ distinct elements $g^{0},g^{1},\ldots,g^{n-1}$ with
$g^{n}=g^{0}=1$. Of course, this group is nothing but the additive abelian
group $\left(  \mathbb{Z}/n,+,0\right)  $, rewritten multiplicatively (so we
write $g^{i}$ for the residue class $\overline{i}\in\mathbb{Z}/n$).

\item There is also an infinite cyclic group $C_{\infty}$, which is generated
by an element $g$ that satisfies no relations. This is just the additive
abelian group $\left(  \mathbb{Z},+,0\right)  $, rewritten multiplicatively
(so we write $g^{i}$ for the integer $i$).
\end{itemize}

Any cyclic group is isomorphic to one of the $C_{n}$'s or $C_{\infty}$. The
cyclic groups have the reputation of being the simplest groups to study, and
so we shall now try to better understand their actions.

We ask ourselves: What \textquotedblleft is\textquotedblright\ a $G$-set
\textquotedblleft really\textquotedblright? In other words, what is the
easiest way to view a $G$-action? For general $G$, we cannot improve on the
definition of a $G$-action (or on the equivalent description in Theorem
\ref{thm.rep.G-sets.curr}). However, when $G$ is a cyclic group, we can
reframe this definition in a much more compact form.

We begin with the most trivial example: that of the cyclic group $C_{1}$, also
known as the trivial group $\left\{  1\right\}  $, also isomorphic to the
symmetric groups $S_{0}$ and $S_{1}$.

\begin{example}
\label{exa.rep.G-sets.C1}Consider the trivial group $\left\{  1\right\}  $.

What is a left $\left\{  1\right\}  $-set? In other words, what does a left
$\left\{  1\right\}  $-action on a set $X$ mean? It means a map
$\rightharpoonup\ :\left\{  1\right\}  \times X\rightarrow X$ that satisfies
the two axioms from Definition \ref{def.rep.G-sets.left} \textbf{(a)}. But the
second axiom is saying that $1\rightharpoonup x=x$ for all $x\in X$ (which
uniquely determines the map $\rightharpoonup$), whereas the first axiom is
just saying that $1\rightharpoonup\left(  1\rightharpoonup x\right)  =\left(
1\cdot1\right)  \rightharpoonup x$ (since the only $g\in\left\{  1\right\}  $
is $1$), which follows automatically from the second axiom. Thus, a left
$\left\{  1\right\}  $-action $\rightharpoonup$ on $X$ is uniquely determined
by the rule%
\[
1\rightharpoonup x=x\ \ \ \ \ \ \ \ \ \ \text{for all }x\in X,
\]
and moreover, the map $\rightharpoonup$ defined by this rule really is a left
$\left\{  1\right\}  $-action on $X$. Hence, any set $X$ has a unique left
$\left\{  1\right\}  $-action on it. Consequently, a left $\left\{  1\right\}
$-set is just a set equipped with this unique left $\left\{  1\right\}
$-action. Roughly speaking, this means that a left $\left\{  1\right\}  $-set
is just a set (the $\left\{  1\right\}  $-action \textquotedblleft comes for
free\textquotedblright).
\end{example}

This shows that $G$-sets are a generalization of sets. Nice, but not very
exciting. But this is what we get by looking at the trivial group. Let us now
move onward to slightly more interesting groups.

\begin{example}
\label{exa.rep.G-sets.C2}Consider the cyclic group $C_{2}=\left\{
1,g\right\}  $ (also known as the symmetric group $S_{2}=\left\{
1,s_{1}\right\}  $ if we identify the generator $g$ with $s_{1}$).

What is a left $C_{2}$-set? In other words, what does a left $C_{2}$-action on
a set $X$ mean? To define such an action, we must specify $1\rightharpoonup x$
and $g\rightharpoonup x$ for each $x\in X$. Again, $1\rightharpoonup x$ is
uniquely determined by the second axiom in Definition
\ref{def.rep.G-sets.left} \textbf{(a)}, which requires $1\rightharpoonup x$ to
be $x$. As for $g\rightharpoonup x$, we know by the first axiom that%
\[
g\rightharpoonup\left(  g\rightharpoonup x\right)  =\underbrace{\left(
gg\right)  }_{=g^{2}=1}\rightharpoonup x=1\rightharpoonup x=x.
\]
Thus, the map%
\begin{align*}
\tau_{g}:X  &  \rightarrow X,\\
x  &  \mapsto g\rightharpoonup x
\end{align*}
must be an involution. Conversely, any involution of $X$ can be used as
$\tau_{g}$ and thus defines a $C_{2}$-action on $X$ (why?). So a $C_{2}%
$-action on $X$ is \textquotedblleft the same as\textquotedblright\ an
involution of $X$. And thus a left $C_{2}$-set is just a set with an
involution on it.

This can also be seen using the curried form (i.e., using Theorem
\ref{thm.rep.G-sets.curr}): In its curried form, a $C_{2}$-action on $X$ is a
group morphism $\rho:C_{2}\rightarrow S_{X}$. Such a group morphism is
uniquely determined by its value $\rho\left(  g\right)  $ (since $\rho\left(
1\right)  $ must be $\operatorname*{id}\in S_{X}$), and this value must be an
involution since $\left(  \rho\left(  g\right)  \right)  ^{2}=\rho\left(
g^{2}\right)  =\rho\left(  1\right)  =\operatorname*{id}$. Conversely, for any
involution $w\in S_{X}$, we can define a group morphism $\rho:C_{2}\rightarrow
S_{X}$ by setting $\rho\left(  1\right)  =\operatorname*{id}$ and $\rho\left(
g\right)  =w$, and this defines a left $C_{2}$-action on $X$.
\end{example}

So a left $C_{2}$-set is just a set $X$ with an involution, i.e., with a
permutation $w\in S_{X}$ satisfying $w^{2}=\operatorname*{id}$. This
generalizes to arbitrary cyclic groups $C_{n}$:

\begin{theorem}
\label{thm.rep.G-sets.Cn}Let $n$ be a positive integer. Let $C_{n}$ be the
cyclic group of size $n$, with generator $g$ (so that $C_{n}=\left\{
g^{0},g^{1},\ldots,g^{n-1}\right\}  $ and $g^{n}=g^{0}=1$). Let $X$ be a set.
Then, a left $C_{n}$-action on $X$ is \textquotedblleft the same
as\textquotedblright\ a permutation $w\in S_{X}$ satisfying $w^{n}%
=\operatorname*{id}$. To be more precise: \medskip

\textbf{(a)} If $\rightharpoonup$ is a left $C_{n}$-action on $X$, then the
map%
\begin{align*}
\tau_{g}:X  &  \rightarrow X,\\
x  &  \mapsto g\rightharpoonup x
\end{align*}
is a permutation $w\in S_{X}$ satisfying $w^{n}=\operatorname*{id}$. \medskip

\textbf{(b)} Conversely, if $w\in S_{X}$ is a permutation satisfying
$w^{n}=\operatorname*{id}$, then we can define a left $C_{n}$-action on $X$ by
setting%
\[
g^{i}\rightharpoonup x=w^{i}\left(  x\right)  \ \ \ \ \ \ \ \ \ \ \text{for
all }i\in\mathbb{Z}\text{ and }x\in X.
\]

\textbf{(c)} These two conversions (from a left $C_{n}$-action to a
permutation $w$, and vice versa) are mutually inverse.
\end{theorem}

\begin{fineprint}
\begin{proof}
[Proof sketch.]Much of this is straightforward, so I only show some
milestones. \medskip

\textbf{(a)} Let $\rightharpoonup$ be a left $C_{n}$-action on $X$. It
suffices to show that $\tau_{g}^{n}=\operatorname*{id}$.

Consider the curried form $\rho:C_{n}\rightarrow S_{X}$ of the left $C_{n}%
$-action $\rightharpoonup$. Then, $\rho$ is a group morphism. Hence,
$\rho\left(  g^{n}\right)  =\left(  \rho\left(  g\right)  \right)  ^{n}$. In
view of $\rho\left(  g\right)  =\tau_{g}$ (by the definition of $\rho$) and
$g^{n}=1$, we can rewrite this as $\rho\left(  1\right)  =\tau_{g}^{n}$.
However, $\rho\left(  1\right)  =\operatorname*{id}$ (since $\rho$ is a group
morphism). Comparing these two equalities, we obtain $\tau_{g}^{n}%
=\operatorname*{id}$. Thus, Theorem \ref{thm.rep.G-sets.Cn} \textbf{(a)} is
proved. \medskip

\textbf{(b)} Let $w\in S_{X}$ be a permutation satisfying $w^{n}%
=\operatorname*{id}$. We want to define a left $C_{n}$-action on $X$ by
setting%
\[
g^{i}\rightharpoonup x=w^{i}\left(  x\right)  \ \ \ \ \ \ \ \ \ \ \text{for
all }i\in\mathbb{Z}\text{ and }x\in X.
\]
First, we need to show that this map $\rightharpoonup$ is well-defined, i.e.,
that the value $w^{i}\left(  x\right)  $ that we are assigning to
$g^{i}\rightharpoonup x$ depends only on $g^{i}$ and on $x$ (but not on the
specific choice of $i$). In other words, we must prove that if $i$ and $j$ are
two integers satisfying $g^{i}=g^{j}$, then $w^{i}\left(  x\right)
=w^{j}\left(  x\right)  $ for each $x\in X$. But this is not hard: If two
integers $i$ and $j$ satisfy $g^{i}=g^{j}$, then $i\equiv j\operatorname{mod}%
n$ (since $g\in C_{n}$ has order $n$), and thus $i=j+nk$ for some
$k\in\mathbb{Z}$; but this entails $w^{i}=w^{j+nk}=w^{j}\left(
\underbrace{w^{n}}_{=\operatorname*{id}}\right)  ^{k}=w^{j}\operatorname*{id}%
\nolimits^{k}=w^{j}$, so that $w^{i}\left(  x\right)  =w^{j}\left(  x\right)
$ for each $x\in X$.

Thus, we can define a map $\rightharpoonup\ :C_{n}\times X\rightarrow X$ by
setting
\[
g^{i}\rightharpoonup x=w^{i}\left(  x\right)  \ \ \ \ \ \ \ \ \ \ \text{for
all }i\in\mathbb{Z}\text{ and }x\in X.
\]
It remains to show that this map $\rightharpoonup$ is an action of $C_{n}$ on
$X$. This requires verifying the two axioms in Definition
\ref{def.rep.G-sets.left} \textbf{(a)}. But this is straightforward (e.g., the
first axiom boils down to $w^{i}\left(  w^{j}\left(  x\right)  \right)
=w^{i+j}\left(  x\right)  $ for all $i,j\in\mathbb{Z}$ and $x\in X$). This
proves Theorem \ref{thm.rep.G-sets.Cn} \textbf{(b)}. \medskip

\textbf{(c)} We must prove the following claims:

\begin{statement}
\textit{Claim 1:} If we convert a left $C_{n}$-action $\rightharpoonup$ into a
permutation $w\in S_{X}$, and then convert the latter permutation $w$ back
into a left $C_{n}$-action $\rightharpoonup^{\prime}$, then this new left
$C_{n}$-action $\rightharpoonup^{\prime}$ is the original $C_{n}$-action
$\rightharpoonup$.
\end{statement}

\begin{statement}
\textit{Claim 2:} If we convert a permutation $w\in S_{X}$ satisfying
$w^{n}=\operatorname*{id}$ into a left $C_{n}$-action $\rightharpoonup$, and
then convert the latter $C_{n}$-action $\rightharpoonup$ back into a
permutation $w^{\prime}\in S_{X}$, then this $w^{\prime}$ is the original $w$.
\end{statement}

\begin{proof}
[Proof of Claim 1.]Let $\rightharpoonup$ be a left $C_{n}$-action on $X$.
Convert this left $C_{n}$-action into a permutation $w\in S_{X}$, and then
convert the latter permutation $w$ back into a left $C_{n}$-action
$\rightharpoonup^{\prime}$. We must prove that this $\rightharpoonup^{\prime}$
is the original $C_{n}$-action $\rightharpoonup$.

Consider the curried form $\rho:C_{n}\rightarrow S_{X}$ of the original left
$C_{n}$-action $\rightharpoonup$. Then, $\rho$ is a group morphism. The
definition of $\rho$ yields $\rho\left(  g\right)  =\tau_{g}$ (where $\tau
_{g}$ is defined as in Theorem \ref{thm.rep.G-sets.Cn} \textbf{(a)}).

By the definition of $w$, we have $w=\tau_{g}$. Comparing this with
$\rho\left(  g\right)  =\tau_{g}$, we find $w=\rho\left(  g\right)  $. Thus,
for each $i\in\mathbb{Z}$, we have%
\begin{equation}
w^{i}=\left(  \rho\left(  g\right)  \right)  ^{i}=\rho\left(  g^{i}\right)
\label{pf.thm.rep.G-sets.Cn.c.i}%
\end{equation}
(since $\rho$ is a group morphism).

By the definition of $\rightharpoonup^{\prime}$, we have%
\[
g^{i}\rightharpoonup^{\prime}x=w^{i}\left(  x\right)
\ \ \ \ \ \ \ \ \ \ \text{for all }i\in\mathbb{Z}\text{ and }x\in X.
\]
Hence, for all $i\in\mathbb{Z}$ and $x\in X$, we have%
\begin{align*}
g^{i}\rightharpoonup^{\prime}x  &  =w^{i}\left(  x\right) \\
&  =\left(  \rho\left(  g^{i}\right)  \right)  \left(  x\right)
\ \ \ \ \ \ \ \ \ \ \left(  \text{since (\ref{pf.thm.rep.G-sets.Cn.c.i})
yields }w^{i}=\rho\left(  g^{i}\right)  \right) \\
&  =g^{i}\rightharpoonup x\ \ \ \ \ \ \ \ \ \ \left(
\begin{array}
[c]{c}%
\text{since }\rho\left(  g^{i}\right)  \text{ is defined to be the map }%
\tau_{g^{i}}:X\rightarrow X\\
\text{that sends each }y\in X\text{ to }g^{i}\rightharpoonup y
\end{array}
\right)  .
\end{align*}
In other words, for all $h\in C_{n}$ and $x\in X$, we have $h\rightharpoonup
^{\prime}x=h\rightharpoonup x$ (since each $h\in C_{n}$ can be written as
$g^{i}$ for some $i\in\mathbb{Z}$). In other words, the map $\rightharpoonup
^{\prime}$ equals the map $\rightharpoonup$. This proves Claim 1.
\end{proof}

\begin{proof}
[Proof of Claim 2.]Let $w\in S_{X}$ be a permutation satisfying $w^{n}%
=\operatorname*{id}$. Convert this permutation $w$ into a left $C_{n}$-action
$\rightharpoonup$, and then convert the latter $C_{n}$-action $\rightharpoonup
$ back into a permutation $w^{\prime}\in S_{X}$. We must prove that this
$w^{\prime}$ is the original $w$.

The definition of $w^{\prime}$ yields $w^{\prime}=\tau_{g}$ (where $\tau_{g}$
is defined as in Theorem \ref{thm.rep.G-sets.Cn} \textbf{(a)}). Thus, for each
$x\in X$, we have%
\begin{align*}
w^{\prime}\left(  x\right)   &  =\tau_{g}\left(  x\right)  =g\rightharpoonup
x\ \ \ \ \ \ \ \ \ \ \left(  \text{by the definition of }\tau_{g}\right) \\
&  =g^{1}\rightharpoonup x\ \ \ \ \ \ \ \ \ \ \left(  \text{since }%
g=g^{1}\right) \\
&  =w^{1}\left(  x\right)  \ \ \ \ \ \ \ \ \ \ \left(  \text{by the definition
of }\rightharpoonup\right)  .
\end{align*}
Hence, $w^{\prime}=w^{1}=w$. This proves Claim 2.
\end{proof}

Thus, the proof of Theorem \ref{thm.rep.G-sets.Cn} \textbf{(c)} is complete.
\end{proof}
\end{fineprint}

Something similar (but simpler) holds for the infinite cyclic group
$C_{\infty}$:

\begin{theorem}
Let $C_{\infty}$ be the infinite cyclic group, with generator $g$. Let $X$ be
a set. Then, a left $C_{\infty}$-action on $X$ is \textquotedblleft the same
as\textquotedblright\ a permutation $w\in S_{X}$. To be more precise: \medskip

\textbf{(a)} If $\rightharpoonup$ is a left $C_{\infty}$-action on $X$, then
the map%
\begin{align*}
\tau_{g}:X  &  \rightarrow X,\\
x  &  \mapsto g\rightharpoonup x
\end{align*}
is a permutation $w\in S_{X}$. \medskip

\textbf{(b)} Conversely, if $w\in S_{X}$ is a permutation, then we can define
a left $C_{\infty}$-action on $X$ by setting%
\[
g^{i}\rightharpoonup x=w^{i}\left(  x\right)  \ \ \ \ \ \ \ \ \ \ \text{for
all }i\in\mathbb{Z}\text{ and }x\in X.
\]

\textbf{(c)} These two conversions (from a left $C_{\infty}$-action to a
permutation $w$, and vice versa) are mutually inverse.
\end{theorem}

\begin{fineprint}
\begin{proof}
Similar to Theorem \ref{thm.rep.G-sets.Cn}.
\end{proof}
\end{fineprint}

Thus, a left $C_{\infty}$-set is nothing but a set with a permutation.

\begin{remark}
There is a similar result about monoid actions: Let $M$ be the monoid $\left(
\mathbb{N},+,0\right)  $. Then, a left $M$-set is just a set with a function
from this set to itself.
\end{remark}

Alas, $S_{n}$-actions for $n\geq3$ cannot be described as just
\textquotedblleft sets with a map on them\textquotedblright, since the group
$S_{n}$ cannot be expressed through a single generator.\footnote{They can be
described in a similar, but more complicated way. For instance, a known fact
(the \textquotedblleft Moore presentation\textquotedblright\ of $S_{n}$) says
that the symmetric group $S_{n}$ can be described by the generators
$s_{1},s_{2},\ldots,s_{n-1}$ and the relations%
\begin{align*}
s_{i}^{2}  &  =1\ \ \ \ \ \ \ \ \ \ \text{for all }i\in\left[  n-1\right]  ;\\
s_{i}s_{j}  &  =s_{j}s_{i}\ \ \ \ \ \ \ \ \ \ \text{for all }i,j\in\left[
n-1\right]  \text{ satisfying }\left\vert i-j\right\vert >1;\\
s_{i}s_{i+1}s_{i}  &  =s_{i+1}s_{i}s_{i+1}\ \ \ \ \ \ \ \ \ \ \text{for all
}i\in\left[  n-2\right]
\end{align*}
(see \cite[Theorem 1.2.4]{Willia03} for a proof). Using this fact, we can then
describe a left $S_{n}$-set as being \textquotedblleft just\textquotedblright%
\ a set $X$ equipped with $n-1$ permutations $w_{1},w_{2},\ldots,w_{n-1}\in
S_{X}$ that satisfy%
\begin{align*}
w_{i}^{2}  &  =1\ \ \ \ \ \ \ \ \ \ \text{for all }i\in\left[  n-1\right]  ;\\
w_{i}w_{j}  &  =w_{j}w_{i}\ \ \ \ \ \ \ \ \ \ \text{for all }i,j\in\left[
n-1\right]  \text{ satisfying }\left\vert i-j\right\vert >1;\\
w_{i}w_{i+1}w_{i}  &  =w_{i+1}w_{i}w_{i+1}\ \ \ \ \ \ \ \ \ \ \text{for all
}i\in\left[  n-2\right]  .
\end{align*}
(These permutations $w_{i}$ are just the actions of the corresponding
$s_{i}\in S_{n}$ on $X$.) Whether this description is any easier (or more
intuitive) than the original definition of a left $S_{n}$-set, however, is
another question.} This makes them all the more interesting.

\subsubsection{Faithfulness}

Here is another useful property of certain actions:

\begin{definition}
\label{def.rep.G-sets.faith}\textbf{(a)} A left $G$-set $X$ is said to be
\emph{faithful} if it has the following property: If $g,h\in G$ are two
elements that satisfy
\[
\left(  g\rightharpoonup x=h\rightharpoonup x\text{ for all }x\in X\right)  ,
\]
then $g=h$.

In other words, a left $G$-set $X$ is \emph{faithful} if and only if the
curried form $\rho:G\rightarrow S_{X}$ of its $G$-action is injective.
\medskip

\textbf{(b)} Instead of saying that a left $G$-set $X$ is faithful, we also
say that the action of $G$ on $X$ is faithful.
\end{definition}

So a left $G$-set is faithful if different elements of $G$ can be told apart
based on how they act on it. Thus, a faithful $G$-set allows us to replace the
elements $g$ of $G$ by the corresponding permutations $\tau_{g}:X\rightarrow
X$ without accidentally equating two different elements.

The following are some examples of faithfulness and its lack:

\begin{itemize}
\item The left regular $G$-action on $G$ (defined in Example
\ref{exa.rep.G-sets.lreg}) is faithful. This is the \textquotedblleft
injective\textquotedblright\ part of Cayley's Theorem (Theorem
\ref{thm.rep.G-sets.cay}).

\item The inverse right regular $G$-action on $G$ (defined in Example
\ref{exa.rep.G-sets.irreg}) is also faithful.

\item The conjugation action of $G$ on $G$ (defined in Example
\ref{exa.rep.G-sets.conj}), given by%
\[
g\rightharpoonup x=gxg^{-1}\ \ \ \ \ \ \ \ \ \ \text{for all }g\in G\text{ and
}x\in G,
\]
is often not faithful. For example, if $g$ is any element of the center of
$G$, then $g\rightharpoonup x=\underbrace{gx}_{=xg}g^{-1}=x\underbrace{gg^{-1}%
}_{=1}=x$ for each $x\in G$; thus, all elements of the center of $G$ act in
the same way (as the identity map). Hence, the conjugation action of $G$ is
not faithful unless the center of $G$ is trivial. The converse is also true
(see Exercise \ref{exe.rep.G-sets.faith.conj} for a proof).

\item The trivial action of $G$ on a set $X$ is only faithful if $G$ is trivial.
\end{itemize}

\begin{exercise}
\fbox{1} Let $X$ be a left $G$-set. Prove that $X$ is faithful if and only if
the only $g\in G$ satisfying $\left(  g\rightharpoonup x=x\text{ for all }x\in
X\right)  $ is $1_{G}$.
\end{exercise}

\begin{exercise}
\label{exe.rep.G-sets.faith.conj}\fbox{1} Prove that the conjugation action of
$G$ on $G$ is faithful if and only if $Z\left(  G\right)  =\left\{  1\right\}
$.
\end{exercise}

\begin{exercise}
\fbox{1} Let $A$ be any set. Let $k\in\mathbb{N}$. Consider the left $S_{A}%
$-set $\mathcal{P}_{k}\left(  A\right)  $ from Example
\ref{exa.rep.G-sets.PkA}. What conditions must $k$ satisfy in order for this
$S_{A}$-set to be faithful?
\end{exercise}

\begin{exercise}
\fbox{1} Let $n$ and $m$ be positive integers, and assume that the ring
$\mathbf{k}$ is not trivial. Prove that the action of $S_{n}\times S_{m}$ on
$\mathbf{k}^{n\times m}$ is faithful.
\end{exercise}

\subsubsection{Further reading}

Group actions have lots of applications; in particular they are helpful in
proving group-theoretical results such as the Sylow theorems. See
\cite[\textquotedblleft Group actions\textquotedblright, \textquotedblleft
Transitive group actions\textquotedblright, \textquotedblleft The Sylow
theorems\textquotedblright]{Conrad}, \cite[Chapter 5]{Goodman}, \cite[Chapter
2]{McClea23}, \cite[Chapter 9]{Carter09} and other texts for this and much
more about group actions.

Group actions can also be used to solve enumerative problems, like counting
\emph{necklaces} (equivalence classes of $n$-tuples under cyclic rotations).
For an introduction to this, see \cite{Davis-Polya} or \cite[Chapter
9]{Loehr-BC}. A book-length treatment can be found in \cite{Kerber99}.

\subsection{\label{sec.rep.G-rep}Linear group actions and representations of
$G$}

Let us continue the study of $G$-sets, but now combining it with the
well-known concept of a $\mathbf{k}$-module. (You can think of $\mathbf{k}$ as
being your favorite field -- such as $\mathbb{Q}$, $\mathbb{R}$, $\mathbb{C}$
or $\mathbb{Z}/p=\mathbb{F}_{p}$ for some prime $p$. Many authors restrict
themselves to the case of $\mathbf{k}=\mathbb{C}$ and find an interesting
theory already at that level. When $\mathbf{k}$ is a field, the $\mathbf{k}%
$-modules are also known as $\mathbf{k}$-vector spaces, so you have linear
algebra at your disposal.)

\begin{convention}
\label{conv.rep.G-rep.G}For the entirety of Section \ref{sec.rep.G-rep}, we
let $G$ be a group (not necessarily finite).
\end{convention}

\subsubsection{Linear actions, aka representations}

\begin{definition}
\label{def.rep.G-rep.G-rep}Let $G$ be a group. \medskip

\textbf{(a)} Let $V$ be a $\mathbf{k}$-module. A left action $\rightharpoonup
\ :G\times V\rightarrow V$ of $G$ on $V$ is said to be $\mathbf{k}%
$\emph{-linear} (or, for short, \emph{linear}) if for each $g\in G$, the map%
\begin{align*}
\tau_{g}:V  &  \rightarrow V,\\
x  &  \mapsto g\rightharpoonup x
\end{align*}
is $\mathbf{k}$-linear. Explicitly, this condition is saying that%
\begin{align*}
g  &  \rightharpoonup\left(  x+y\right)  =\left(  g\rightharpoonup x\right)
+\left(  g\rightharpoonup y\right)  \ \ \ \ \ \ \ \ \ \ \text{for all }g\in
G\text{ and }x,y\in V;\\
g  &  \rightharpoonup\left(  \lambda x\right)  =\lambda\left(
g\rightharpoonup x\right)  \ \ \ \ \ \ \ \ \ \ \text{for all }g\in G\text{ and
}\lambda\in\mathbf{k}\text{ and }x\in V;\\
g  &  \rightharpoonup0_{V}=0_{V}\ \ \ \ \ \ \ \ \ \ \text{for all }g\in G.
\end{align*}

\textbf{(b)} A \emph{left representation of }$G$\emph{ over }$\mathbf{k}$ (aka
a \emph{left }$G$\emph{-representation over }$\mathbf{k}$) means a
$\mathbf{k}$-module $V$ equipped with a $\mathbf{k}$-linear left $G$-action on
$V$.

We usually omit the words \textquotedblleft left\textquotedblright\ and
\textquotedblleft over $\mathbf{k}$\textquotedblright, unless the context
suggests something different. Thus, we just say \textquotedblleft
representation of $G$\textquotedblright\ instead of \textquotedblleft left
representation of $G$ over $\mathbf{k}$\textquotedblright. Some authors also
refer to representations of $G$ as \textquotedblleft$G$\emph{-modules}%
\textquotedblright, but we will avoid this notation in favor of a better one
that we will soon introduce.
\end{definition}

Let us revisit the above examples of $G$-actions, to see which of them are
$\mathbf{k}$-linear. Obviously, we only need to consider the ones where the
$G$-set is a $\mathbf{k}$-module.

\begin{example}
If $X$ is a $\mathbf{k}$-module, then the trivial $G$-action on $X$ (as
defined in Example \ref{exa.rep.G-sets.triv}) is linear (since the map
$\operatorname*{id}\nolimits_{X}:X\rightarrow X$ is $\mathbf{k}$-linear).
\end{example}

\begin{example}
If $A$ is a nontrivial $\mathbf{k}$-module, then the natural action of $S_{A}$
on $A$ (as defined in Example \ref{exa.rep.G-sets.nat}) is not linear, since a
permutation of $A$ is usually not $\mathbf{k}$-linear.
\end{example}

\begin{example}
\label{exa.rep.G-rep.lin.place}Let $A$ be a $\mathbf{k}$-module, and
$n\in\mathbb{N}$. Then, the place permutation action of $S_{n}$ on $A^{n}$ (as
defined in Example \ref{exa.rep.G-sets.place}) is $\mathbf{k}$-linear. To
check this, we need to verify the three equalities in Definition
\ref{def.rep.G-rep.G-rep} \textbf{(a)}; let us check the first of them: We
must prove that
\[
g\rightharpoonup\left(  x+y\right)  =\left(  g\rightharpoonup x\right)
+\left(  g\rightharpoonup y\right)  \ \ \ \ \ \ \ \ \ \ \text{for all }g\in
S_{n}\text{ and }x,y\in A^{n}.
\]
To prove this, we write $x$ and $y$ as $x=\left(  x_{1},x_{2},\ldots
,x_{n}\right)  $ and $y=\left(  y_{1},y_{2},\ldots,y_{n}\right)  $, and
observe that%
\begin{align*}
g\rightharpoonup\left(  x+y\right)   &  =g\rightharpoonup\underbrace{\left(
\left(  x_{1},x_{2},\ldots,x_{n}\right)  +\left(  y_{1},y_{2},\ldots
,y_{n}\right)  \right)  }_{=\left(  x_{1}+y_{1},\ x_{2}+y_{2},\ \ldots
,\ x_{n}+y_{n}\right)  }\\
&  =g\rightharpoonup\left(  x_{1}+y_{1},\ x_{2}+y_{2},\ \ldots,\ x_{n}%
+y_{n}\right) \\
&  =\left(  x_{g^{-1}\left(  1\right)  }+y_{g^{-1}\left(  1\right)
},\ x_{g^{-1}\left(  2\right)  }+y_{g^{-1}\left(  2\right)  },\ \ldots
,\ x_{g^{-1}\left(  n\right)  }+y_{g^{-1}\left(  n\right)  }\right) \\
&  \ \ \ \ \ \ \ \ \ \ \ \ \ \ \ \ \ \ \ \ \left(  \text{by the definition of
}\rightharpoonup\right) \\
&  =\underbrace{\left(  x_{g^{-1}\left(  1\right)  },\ x_{g^{-1}\left(
2\right)  },\ \ldots,\ x_{g^{-1}\left(  n\right)  }\right)  }%
_{\substack{=g\rightharpoonup\left(  x_{1},x_{2},\ldots,x_{n}\right)
\\\text{(by the definition of }\rightharpoonup\text{)}}}+\underbrace{\left(
y_{g^{-1}\left(  1\right)  },\ y_{g^{-1}\left(  2\right)  },\ \ldots
,\ y_{g^{-1}\left(  n\right)  }\right)  }_{\substack{=g\rightharpoonup\left(
y_{1},y_{2},\ldots,y_{n}\right)  \\\text{(by the definition of }%
\rightharpoonup\text{)}}}\\
&  =\left(  g\rightharpoonup\underbrace{\left(  x_{1},x_{2},\ldots
,x_{n}\right)  }_{=x}\right)  +\left(  g\rightharpoonup\underbrace{\left(
y_{1},y_{2},\ldots,y_{n}\right)  }_{=y}\right) \\
&  =\left(  g\rightharpoonup x\right)  +\left(  g\rightharpoonup y\right)  .
\end{align*}
The other two equalities are proved similarly.

Thus, $A^{n}$ is a representation of $S_{n}$ over $\mathbf{k}$.
\end{example}

\begin{example}
The action of $S_{n}\times S_{m}$ on $\mathbf{k}^{n\times m}$ defined in
Example \ref{exa.rep.G-sets.knm} is $\mathbf{k}$-linear. This is
straightforward to check.
\end{example}

\begin{example}
\label{exa.rep.G-rep.lin.1}Consider the trivial group $\left\{  1\right\}  $.
For any $\mathbf{k}$-module $X$, there is a unique left $\left\{  1\right\}
$-action on $X$ (as we learned in Example \ref{exa.rep.G-sets.C1}), and this
$\left\{  1\right\}  $-action is $\mathbf{k}$-linear (since the action
$\tau_{1}$ of $1$ is $\operatorname*{id}\nolimits_{X}$, which is a
$\mathbf{k}$-linear map). Thus, a representation of $\left\{  1\right\}  $ is
just a $\mathbf{k}$-module.
\end{example}

\begin{example}
\label{exa.rep.G-rep.lin.C2}As we saw in Example \ref{exa.rep.G-sets.C2}, an
action of the size-$2$ cyclic group $C_{2}$ on a set $X$ is just an involution
of $X$. Assuming that $X$ is a $\mathbf{k}$-module, this action is
$\mathbf{k}$-linear if and only if the involution is $\mathbf{k}$-linear. For
example, the involution
\begin{align*}
\mathbf{k}  &  \rightarrow\mathbf{k},\\
x  &  \mapsto-x
\end{align*}
is $\mathbf{k}$-linear (and thus gives a $\mathbf{k}$-linear action of $C_{2}%
$), but the involution%
\begin{align*}
\mathbf{k}  &  \rightarrow\mathbf{k},\\
x  &  \mapsto1-x
\end{align*}
is not.

Thus, a representation of $C_{2}$ is just a $\mathbf{k}$-module equipped with
a $\mathbf{k}$-linear involution.
\end{example}

\begin{example}
\label{exa.rep.G-rep.lin.Cn}Let $n$ be a positive integer. As we saw in
Theorem \ref{thm.rep.G-sets.Cn}, an action of the size-$n$ cyclic group
$C_{n}$ on a set $X$ is just a permutation $w$ of $X$ satisfying
$w^{n}=\operatorname*{id}$. Assuming that $X$ is a $\mathbf{k}$-module, this
action is $\mathbf{k}$-linear if and only if the map $w$ is $\mathbf{k}%
$-linear. For example, for $\mathbf{k}=\mathbb{R}$ and $X=\mathbb{R}^{2}$, we
can take $w$ to be the left multiplication by the matrix%
\[
\left(
\begin{array}
[c]{cc}%
\cos\dfrac{2\pi}{n} & -\sin\dfrac{2\pi}{n}\\
\sin\dfrac{2\pi}{n} & \cos\dfrac{2\pi}{n}%
\end{array}
\right)  \in\mathbb{R}^{2\times2},
\]
aka a rotation by $2\pi/n$ around the origin (viewing $X=\mathbb{R}^{2}$ as
the Cartesian plane). If we take $w$ to be this rotation, then $\mathbb{R}%
^{2}$ thus becomes a representation of $C_{n}$ over $\mathbb{R}$, and the
action of $C_{n}$ is given explicitly by%
\begin{align*}
g^{k}  &  \rightharpoonup v=\left(  v\text{ rotated by }2\pi k/n\text{ around
the origin}\right) \\
&  \ \ \ \ \ \ \ \ \ \ \text{for each }k\in\mathbb{Z}\text{ and }%
v\in\mathbb{R}^{2}.
\end{align*}

\end{example}

\begin{example}
\label{exa.rep.G-rep.GL-nat}The natural action of $\operatorname*{GL}%
\nolimits_{n}\left(  \mathbb{R}\right)  $ on $\mathbb{R}^{n\times1}$ (that is,
the action defined in Example \ref{exa.rep.G-sets.GL}) is $\mathbb{R}$-linear.
More generally, the natural action of the group%
\[
\operatorname*{GL}\nolimits_{n}\left(  \mathbf{k}\right)  :=\left\{  \text{all
invertible }n\times n\text{-matrices over }\mathbf{k}\right\}
\]
on $\mathbf{k}^{n\times1}$ is $\mathbf{k}$-linear. This is called the
\emph{natural representation} of $\operatorname*{GL}\nolimits_{n}\left(
\mathbf{k}\right)  $.
\end{example}

We can generalize the group $\operatorname*{GL}\nolimits_{n}\left(
\mathbf{k}\right)  $ as follows:

\begin{definition}
\label{def.GLV}Let $V$ be any $\mathbf{k}$-module. Then, we define%
\begin{align*}
\operatorname*{GL}\left(  V\right)  :=  &  \ \left\{  \text{all invertible
}\mathbf{k}\text{-linear maps from }V\text{ to }V\right\} \\
=  &  \ \left\{  \text{all automorphisms of the }\mathbf{k}\text{-module
}V\right\} \\
=  &  \ \left\{  \text{all }\mathbf{k}\text{-linear permutations of
}V\right\}  .
\end{align*}
This is a group under composition, and is called the \emph{general linear
group} of $V$. It is a subgroup of the symmetric group $S_{V}$.
\end{definition}

In particular, if $V=\mathbf{k}^{n\times1}$ (the $\mathbf{k}$-module of column
vectors of size $n$), then $\operatorname*{GL}\left(  V\right)  \cong%
\operatorname*{GL}\nolimits_{n}\left(  \mathbf{k}\right)  $ by the canonical
isomorphism that sends each linear map $f:V\rightarrow V$ to the $n\times
n$-matrix that represents this map in the standard basis of $\mathbf{k}%
^{n\times1}$. Now, Example \ref{exa.rep.G-rep.GL-nat} generalizes as follows:

\begin{example}
\label{exa.rep.G-rep.GLV-nat}For any $\mathbf{k}$-module $V$, we have a left
$\operatorname*{GL}\left(  V\right)  $-action on $V$ given by%
\[
f\rightharpoonup v=f\left(  v\right)  \ \ \ \ \ \ \ \ \ \ \text{for all }%
f\in\operatorname*{GL}\left(  V\right)  \text{ and }v\in V.
\]
This action is $\mathbf{k}$-linear, and is called the \emph{natural action} of
$\operatorname*{GL}\left(  V\right)  $. The $\mathbf{k}$-module $V$, equipped
with this action, is called the \emph{natural representation} of
$\operatorname*{GL}\left(  V\right)  $.
\end{example}

\subsubsection{The curried form}

As we have seen in Definition \ref{def.rep.G-sets.curr}, we can curry an
arbitrary action of $G$ on any set $X$ to obtain a group morphism
$\rho:G\rightarrow S_{X}$. Similarly, we can curry a linear action of $G$ on a
$\mathbf{k}$-module $X$ to obtain a group morphism $\rho:G\rightarrow
\operatorname*{GL}\left(  X\right)  $. In fact, this is the same morphism
$\rho$ as in Definition \ref{def.rep.G-sets.curr}, just with its target
$S_{X}$ replaced by $\operatorname*{GL}\left(  X\right)  $:

\begin{theorem}
\label{thm.rep.G-rep.curry}Let $X$ be a left $G$-set that is also a
$\mathbf{k}$-module. Then, the left action $\rightharpoonup$ of $G$ on $X$ is
$\mathbf{k}$-linear if and only if its curried form $\rho:G\rightarrow S_{X}$
satisfies $\rho\left(  G\right)  \subseteq\operatorname*{GL}\left(  X\right)
$.
\end{theorem}

\begin{proof}
[Proof sketch.]This is more or less a restatement of Definition
\ref{def.rep.G-rep.G-rep} \textbf{(a)}. To wit, consider the map%
\begin{align*}
\tau_{g}:X  &  \rightarrow X,\\
x  &  \mapsto g\rightharpoonup x
\end{align*}
for each $g\in G$. Then, the definition of the curried form $\rho$ shows that%
\begin{equation}
\rho\left(  g\right)  =\tau_{g}\ \ \ \ \ \ \ \ \ \ \text{for each }g\in G.
\label{pf.thm.rep.G-rep.curry.rho}%
\end{equation}
But we have the following chain of equivalences:%
\begin{align*}
&  \ \left(  \text{the left action }\rightharpoonup\text{ of }G\text{ on
}X\text{ is }\mathbf{k}\text{-linear}\right) \\
&  \Longleftrightarrow\ \left(  \text{the map }\tau_{g}\text{ is }%
\mathbf{k}\text{-linear for each }g\in G\right) \\
&  \ \ \ \ \ \ \ \ \ \ \ \ \ \ \ \ \ \ \ \ \left(  \text{by Definition
\ref{def.rep.G-rep.G-rep} \textbf{(a)}}\right) \\
&  \Longleftrightarrow\ \left(  \text{the map }\rho\left(  g\right)  \text{ is
}\mathbf{k}\text{-linear for each }g\in G\right) \\
&  \ \ \ \ \ \ \ \ \ \ \ \ \ \ \ \ \ \ \ \ \left(  \text{since
(\ref{pf.thm.rep.G-rep.curry.rho}) yields }\tau_{g}=\rho\left(  g\right)
\text{ for each }g\in G\right) \\
&  \Longleftrightarrow\ \left(  \text{we have }\rho\left(  g\right)
\in\operatorname*{GL}\left(  X\right)  \text{ for each }g\in G\right) \\
&  \ \ \ \ \ \ \ \ \ \ \ \ \ \ \ \ \ \ \ \ \left(
\begin{array}
[c]{c}%
\text{because for each }g\in G\text{, we have }\rho\left(  g\right)  \in
S_{X}\\
\text{(since }\rho\text{ is a map from }G\text{ to }S_{X}\text{), which shows
that}\\
\text{the map }\rho\left(  g\right)  \text{ is invertible, and thus }%
\rho\left(  g\right)  \text{ is }\mathbf{k}\text{-linear}\\
\text{if and only if }\rho\left(  g\right)  \in\operatorname*{GL}\left(
X\right)
\end{array}
\right) \\
&  \Longleftrightarrow\ \left(  \rho\left(  G\right)  \subseteq
\operatorname*{GL}\left(  X\right)  \right)  .
\end{align*}
This proves Theorem \ref{thm.rep.G-rep.curry}.
\end{proof}

Theorem \ref{thm.rep.G-rep.curry} yields the following $\mathbf{k}$-linear
analogue of Theorem \ref{thm.rep.G-sets.curr}:

\begin{theorem}
\label{thm.rep.G-rep.curr}Let $G$ be a group. Let $X$ be a $\mathbf{k}%
$-module. Then, a $\mathbf{k}$-linear left $G$-action on $X$ is
\textquotedblleft the same as\textquotedblright\ a group morphism from $G$ to
$\operatorname*{GL}\left(  X\right)  $. To be more precise: \medskip

\textbf{(a)} If $\rightharpoonup$ is a $\mathbf{k}$-linear left $G$-action on
$X$, then the curried form of $\rightharpoonup$ is a group morphism from $G$
to $\operatorname*{GL}\left(  X\right)  $. \medskip

\textbf{(b)} If $\rho$ is a group morphism from $G$ to $\operatorname*{GL}%
\left(  X\right)  $, then the uncurried form of $\rho$ is a $\mathbf{k}%
$-linear left $G$-action on $X$. \medskip

\textbf{(c)} Currying and uncurrying are mutually inverse.
\end{theorem}

\begin{proof}
[Proof sketch.]\textbf{(a)} Let $\rightharpoonup$ be a $\mathbf{k}$-linear
left $G$-action on $X$. Let $\rho$ be the curried form of $\rightharpoonup$.
Then, $\rho$ is a group morphism from $G$ to $S_{X}$ (by the definition of a
curried form), but also satisfies $\rho\left(  G\right)  \subseteq
\operatorname*{GL}\left(  X\right)  $ (by Theorem \ref{thm.rep.G-rep.curry}).
Hence, $\rho$ is a group morphism from $G$ to $\operatorname*{GL}\left(
X\right)  $. This proves Theorem \ref{thm.rep.G-rep.curr} \textbf{(a)}.
\medskip

\textbf{(b)} Let $\rho$ be a group morphism from $G$ to $\operatorname*{GL}%
\left(  X\right)  $. We must prove that the uncurried form of $\rho$ is a
$\mathbf{k}$-linear left $G$-action on $X$.

Let $\rightharpoonup$ be this uncurried form. Then, $\rightharpoonup$ is a
$G$-action on $X$, and the curried form of $\rightharpoonup$ is $\rho$ (by
Theorem \ref{thm.rep.G-sets.curr} \textbf{(c)}). Hence, Theorem
\ref{thm.rep.G-rep.curry} yields that the left action $\rightharpoonup$ of $G$
on $X$ is $\mathbf{k}$-linear (since $\rho\left(  G\right)  \subseteq
\operatorname*{GL}\left(  X\right)  $). In other words, the uncurried form of
$\rho$ is a $\mathbf{k}$-linear left $G$-action on $X$ (since this uncurried
form is exactly $\rightharpoonup$). This proves Theorem
\ref{thm.rep.G-rep.curr} \textbf{(b)}. \medskip

\textbf{(c)} This is just a particular case of Theorem
\ref{thm.rep.G-sets.curr} \textbf{(c)}.
\end{proof}

\begin{corollary}
\label{cor.rep.G-rep.curr2}A representation of $G$ is \textquotedblleft the
same as\textquotedblright\ a $\mathbf{k}$-module $X$ equipped with a group
morphism $\rho:G\rightarrow\operatorname*{GL}\left(  X\right)  $.
\end{corollary}

This description is used as the definition of a representation of $G$ in most
textbooks (e.g., \cite[Definition 1.1.1]{Prasad-rep}, \cite[\S 1.1]{Serre77},
\cite[Definition 4.1]{McClea23}, \cite[\S 1.1.1]{CSScTo10}, \cite[\S 109]%
{Elman22}, \cite[\S 1.5]{Roseng15}, \cite[Definition 1.1]{GruSer18}).

\subsubsection{Outlook on representation theory}

The above definitions are the starting point for the study of representations
of groups -- a wide field of mathematical research known as
\emph{representation theory of groups}. Even the most classical setting, in
which the group $G$ is finite and the ring $\mathbf{k}$ is a field of
characteristic $0$, is full of surprises and nontrivial results; by now, this
part of the theory is well-exposed in several textbooks (\cite{EGHetc11},
\cite{Serre77}, \cite{Prasad-rep}, \cite{CurRei62}, \cite{AlFais25},
\cite[Chapters 1--2]{GruSer18}, etc.). The case of finite groups $G$ and
positive-characteristic fields $\mathbf{k}$ is much wilder (in the sense that
there are much fewer general results and much more variety of behavior), and
is a topic of ongoing research known as \textquotedblleft\emph{modular
representation theory}\textquotedblright\ (see \cite[Part III]{Serre77},
\cite{Lassue23}, \cite[Chapter III]{Morel19}, \cite[Chapter XII]{CurRei62},
\cite{Schnei13} or \cite{Webb16} for introductions, and -- e.g. --
\cite[Chapter 2]{Zimmer14} for a deeper dive). Infinite groups $G$ have also
been studied extensively, often based on their topological and geometrical
structures (see \cite[Chapters 4--5]{Sternb95}, \cite{Proces07}, \cite[Chapter
V]{Morel19} and many other places\footnote{For an infinite group $G$, there is
usually not much that can be said about all $G$-representations without any
restrictions; they are \textquotedblleft too wild\textquotedblright\ (e.g.,
there are uncountably many $1$-dimensional $\mathbb{Q}$-linear representations
of the multiplicative group $\mathbb{Q}^{\times}=\left(  \mathbb{Q}%
\setminus\left\{  0\right\}  ,\ \cdot,\ 1\right)  $ over $\mathbb{Q}$; and
this is far from the worst that can happen). Thus, much of the actual science
is devoted to specific classes of groups (topological groups, Lie groups,
pro-finite groups, etc.) and specific classes of representations (polynomial,
rational, holomorphic, etc.). Among all these settings, two are particularly
related to ours:
\par
\begin{enumerate}
\item The polynomial representations of the general linear group
$\operatorname*{GL}\nolimits_{n}\mathbf{k}$ have a particularly nice theory
that mirrors and is closely related to the representation theory of symmetric
groups; see \cite[Chapter 8]{Fulton97} or \cite[Chapter 9]{Proces07} for
fairly elementary treatments.
\par
\item There is also a representation theory of the infinite symmetric group
(which, in this context, means the group of all permutations $\sigma$ of
$\left\{  1,2,3,\ldots\right\}  $ that fix all but finitely many elements);
see \cite[Chapter 11]{Meliot17} or \cite{BorOls16}.
\end{enumerate}
\par
However, we will not veer into any of this in the present text.}).

Our focus in this course is different: We are concerned almost entirely with
the case when $G$ is the symmetric group $S_{n}$, but we try to avoid any
unnecessary assumptions on $\mathbf{k}$; we will even allow $\mathbf{k}$ to be
an arbitrary commutative ring much of the time. We will classify the
irreducible\footnote{Roughly speaking, a representation is said to be
\textquotedblleft irreducible\textquotedblright\ if it has no smaller
representations \textquotedblleft inside it\textquotedblright. We will make
this precise soon.} representations of $S_{n}$ over characteristic-$0$ fields
$\mathbf{k}$, but we will also (try to) analyze some non-irreducible ones in
some detail, in particular decomposing them explicitly into
subrepresentations. Thus, compared to most of the literature, we restrict
ourselves to a much smaller set of groups $G$ (viz., the symmetric groups),
but we allow ourselves to require less of the field $\mathbf{k}$ and to ask
some deeper questions that make no sense for arbitrary $G$.

This said, there already exist several texts that are (fully or in great part)
devoted to the representation theory of the symmetric groups specifically:
\cite{JamKer81}, \cite{Sagan01}, \cite{Feray16}, \cite{Kerber99},
\cite{James78}, \cite{Wildon18}, \cite{Howe22}, \cite{BleSch05},
\cite{Zelevi81b}, \cite{Knutso73}, \cite{Meliot17}, \cite{Prasad-rep} and many
others. We will see a lot of overlap with many of these texts, but we will
often make different choices in what results we cover and how we prove them.

\subsubsection{Representations of $S_{n}$: first examples}

We have already seen a few examples of representations of $S_{n}$. Let us
recall them again:

\begin{example}
\label{exa.rep.Sn-rep.triv}For any $\mathbf{k}$-module $V$, the trivial
$S_{n}$-action on $V$ is $\mathbf{k}$-linear, and thus gives a representation
of $S_{n}$. This is called a \emph{trivial representation} of $S_{n}$. Recall
that its action is given by%
\[
g\rightharpoonup v=v\ \ \ \ \ \ \ \ \ \ \text{for all }g\in S_{n}\text{ and
}v\in V.
\]
When $V=\mathbf{k}$, we refer to this trivial representation as
\textquotedblleft\textbf{the} trivial representation\textquotedblright\ of
$S_{n}$ (with a definite article), and we denote it by $\mathbf{k}%
_{\operatorname*{triv}}$.
\end{example}

\begin{example}
\label{exa.rep.Sn-rep.nat}The symmetric group $S_{n}$ acts on the $\mathbf{k}%
$-module $\mathbf{k}^{n}$ by%
\begin{align*}
g  &  \rightharpoonup\left(  a_{1},a_{2},\ldots,a_{n}\right)  =\left(
a_{g^{-1}\left(  1\right)  },a_{g^{-1}\left(  2\right)  },\ldots
,a_{g^{-1}\left(  n\right)  }\right) \\
&  \ \ \ \ \ \ \ \ \ \ \ \ \ \ \ \ \ \ \ \ \text{for all }g\in S_{n}\text{ and
}\left(  a_{1},a_{2},\ldots,a_{n}\right)  \in\mathbf{k}^{n}.
\end{align*}
This is the place permutation action on $\mathbf{k}^{n}$ we introduced in
Example \ref{exa.rep.G-sets.place} (for $A=\mathbf{k}$). As we know from
Example \ref{exa.rep.G-rep.lin.place}, this left $S_{n}$-action is
$\mathbf{k}$-linear, and therefore $\mathbf{k}^{n}$ becomes a representation
of $S_{n}$. This representation is called the \emph{natural representation} of
$S_{n}$.
\end{example}

The trivial representation has a slightly less trivial sibling, which is also
defined for any $\mathbf{k}$-module $V$:

\begin{example}
\label{exa.rep.Sn-rep.sign}For any $\mathbf{k}$-module $V$, there is a
$\mathbf{k}$-linear left $S_{n}$-action on $V$ given by%
\[
g\rightharpoonup v=\left(  -1\right)  ^{g}v\ \ \ \ \ \ \ \ \ \ \text{for all
}g\in S_{n}\text{ and }v\in V.
\]
This gives a representation of $S_{n}$, which is called a \emph{sign
representation} of $S_{n}$.

When $V=\mathbf{k}$, we refer to this sign representation as \textquotedblleft%
\textbf{the} sign representation\textquotedblright\ of $S_{n}$ (with a
definite article), and we denote it by $\mathbf{k}_{\operatorname*{sign}}$.
\end{example}

\begin{fineprint}
\begin{proof}
[Proof of Example \ref{exa.rep.Sn-rep.sign} (sketched).]Similarly to Example
\ref{exa.rep.G-sets.sign}, we can show that the formula $g\rightharpoonup
v=\left(  -1\right)  ^{g}v$ really defines a left $S_{n}$-action on $V$. That
this $S_{n}$-action is $\mathbf{k}$-linear is obvious (since the action of a
$g\in S_{n}$ on $V$ is just scaling by the number $\left(  -1\right)  ^{g}$).
\end{proof}
\end{fineprint}

Inspired by these examples, we can easily come up with another (combining the
sign representation with the natural one):

\begin{example}
\label{exa.rep.Sn-rep.sign-nat}The \emph{sign-twisted natural representation}
of $S_{n}$ is the $\mathbf{k}$-module $\mathbf{k}^{n}$ equipped with the left
$S_{n}$-action given by%
\begin{align*}
g  &  \rightharpoonup\left(  a_{1},a_{2},\ldots,a_{n}\right)  =\left(
-1\right)  ^{g}\left(  a_{g^{-1}\left(  1\right)  },a_{g^{-1}\left(  2\right)
},\ldots,a_{g^{-1}\left(  n\right)  }\right) \\
&  \ \ \ \ \ \ \ \ \ \ \ \ \ \ \ \ \ \ \ \ \text{for all }g\in S_{n}\text{ and
}\left(  a_{1},a_{2},\ldots,a_{n}\right)  \in\mathbf{k}^{n}.
\end{align*}

\end{example}

Here are some further representations of $S_{n}$:

\begin{example}
\label{exa.rep.Sn-rep.polyring}\textbf{(a)} Consider the polynomial ring
$\mathbf{k}\left[  x_{1},x_{2},\ldots,x_{n}\right]  $ in $n$ indeterminates
$x_{1},x_{2},\ldots,x_{n}$ over $\mathbf{k}$. We define a left $S_{n}$-action
on this polynomial ring $\mathbf{k}\left[  x_{1},x_{2},\ldots,x_{n}\right]  $
by setting%
\begin{align*}
&  g\rightharpoonup f=f\left(  x_{g\left(  1\right)  },x_{g\left(  2\right)
},\ldots,x_{g\left(  n\right)  }\right) \\
&  \ \ \ \ \ \ \ \ \ \ \text{for all }g\in S_{n}\text{ and }f\left(
x_{1},x_{2},\ldots,x_{n}\right)  \in\mathbf{k}\left[  x_{1},x_{2},\ldots
,x_{n}\right]  .
\end{align*}
(Thus, a permutation $g\in S_{n}$ acts on a polynomial $f$ by permuting the
indeterminates.) This is indeed a left $S_{n}$-action (see \cite[Proposition
7.1.4]{21s} for a proof), and is easily seen to be $\mathbf{k}$-linear, so
that the polynomial ring $\mathbf{k}\left[  x_{1},x_{2},\ldots,x_{n}\right]  $
becomes a representation of $S_{n}$. \medskip

\textbf{(b)} Let $d\in\mathbb{N}$. Define $\mathbf{k}\left[  x_{1}%
,x_{2},\ldots,x_{n}\right]  _{d}$ to be the set of all homogeneous polynomials
of degree $d$ in $\mathbf{k}\left[  x_{1},x_{2},\ldots,x_{n}\right]  $. That
is,%
\[
\mathbf{k}\left[  x_{1},x_{2},\ldots,x_{n}\right]  _{d}=\operatorname*{span}%
\nolimits_{\mathbf{k}}\left\{  x_{1}^{a_{1}}x_{2}^{a_{2}}\cdots x_{n}^{a_{n}%
}\ \mid\ a_{1}+a_{2}+\cdots+a_{n}=d\right\}  .
\]
Thus, $\mathbf{k}\left[  x_{1},x_{2},\ldots,x_{n}\right]  _{d}$ is a
$\mathbf{k}$-submodule of $\mathbf{k}\left[  x_{1},x_{2},\ldots,x_{n}\right]
$. We define a left $S_{n}$-action on this $\mathbf{k}$-submodule
$\mathbf{k}\left[  x_{1},x_{2},\ldots,x_{n}\right]  _{d}$ in the same way as
we defined it on the whole $\mathbf{k}\left[  x_{1},x_{2},\ldots,x_{n}\right]
$ in part \textbf{(a)}. Thus, $\mathbf{k}\left[  x_{1},x_{2},\ldots
,x_{n}\right]  _{d}$ becomes a representation of $S_{n}$ as well.
\end{example}

What other representations of $S_{n}$ are there? We will soon see.

\subsubsection{Morphisms of representations}

Let us now define morphisms of representations. As with any algebraic
structure, this is done in the obvious way: The morphisms are just the maps
that preserve all the structures. In detail:

\begin{definition}
\label{def.rep.G-rep.mor}Let $G$ be a group. Let $V$ and $W$ be two
representations of $G$ over $\mathbf{k}$. Let $f:V\rightarrow W$ be any map.
Then, we say that $f$ is a \emph{morphism of representations} (aka a
$G$\emph{-representation morphism}) if $f$ is both $\mathbf{k}$-linear (i.e.,
a morphism of $\mathbf{k}$-modules) and $G$-equivariant (i.e., a morphism of
$G$-sets).
\end{definition}

More classical synonyms for \textquotedblleft morphism of
representations\textquotedblright\ are \textquotedblleft\emph{intertwining
map}\textquotedblright\ or \textquotedblleft\emph{intertwiner}%
\textquotedblright, although I am not sure how standard the meanings of these
words are (and I believe that authors are avoiding them nowadays).

\begin{example}
\label{exa.rep.G-rep.mor.Sn}Consider the natural representation $\mathbf{k}%
^{n}$ of the symmetric group $S_{n}$ (defined in Example
\ref{exa.rep.Sn-rep.nat}). Consider furthermore the trivial representation
$\mathbf{k}$ of the symmetric group $S_{n}$ (defined in Example
\ref{exa.rep.Sn-rep.triv}). Then: \medskip

\textbf{(a)} The map%
\begin{align*}
\mathbf{k}  &  \rightarrow\mathbf{k}^{n},\\
a  &  \mapsto\left(  \underbrace{a,a,\ldots,a}_{n\text{ times}}\right)
\end{align*}
is a morphism of representations (since it is clearly $\mathbf{k}$-linear, and
since we know from Example \ref{exa.rep.G-sets.mor-Sn} \textbf{(a)} that it is
$S_{n}$-equivariant). \medskip

\textbf{(b)} The map%
\begin{align*}
\mathbf{k}^{n}  &  \rightarrow\mathbf{k},\\
\left(  a_{1},a_{2},\ldots,a_{n}\right)   &  \mapsto a_{1}+a_{2}+\cdots+a_{n}%
\end{align*}
is a morphism of representations (since it is clearly $\mathbf{k}$-linear, and
since we know from Example \ref{exa.rep.G-sets.mor-Sn} \textbf{(b)} that it is
$S_{n}$-equivariant). \medskip

\textbf{(c)} The map%
\begin{align*}
\mathbf{k}^{n}  &  \rightarrow\mathbf{k},\\
\left(  a_{1},a_{2},\ldots,a_{n}\right)   &  \mapsto a_{1}a_{2}\cdots a_{n}%
\end{align*}
is $S_{n}$-equivariant, but not a morphism of representations, since it is not
$\mathbf{k}$-linear (unless $n=1$ or $\mathbf{k}$ is trivial).
\end{example}

With morphisms come isomorphisms, which again are defined in the usual way:

\begin{definition}
\label{def.rep.G-rep.iso}An \emph{isomorphism of representations} is an
invertible morphism of representations whose inverse is also a morphism of representations.
\end{definition}

\begin{proposition}
\label{prop.rep.G-rep.iso=bij}The \textquotedblleft is also a morphism of
representations\textquotedblright\ part in Definition \ref{def.rep.G-rep.iso}
is redundant (i.e., follows from the other requirements). In other words, if a
morphism of representations is invertible, then its inverse is again a
morphism of representations.
\end{proposition}

\begin{proof}
Straightforward and standard, just like Proposition
\ref{prop.rep.G-sets.iso=bij}.
\end{proof}

\begin{definition}
\label{def.rep.G-rep.isomorphic}Let $X$ and $Y$ be two representations of a
group $G$. Then, we say that $X$ and $Y$ are \emph{isomorphic} (this is
written $X\cong Y$) if there exists an isomorphism of representations from $X$
to $Y$.
\end{definition}

\subsubsection{\label{subsec.rep.G-rep.subreps}Subrepresentations}

Subrepresentations are defined in the same way as any kinds of sub-objects in algebra:

\begin{definition}
\label{def.rep.G-rep.sub}Let $V$ be a representation of a group $G$ over
$\mathbf{k}$. A \emph{subrepresentation} of $V$ means a $\mathbf{k}$-submodule
of $V$ that is simultaneously a $G$-subset of $V$.

Any such subrepresentation automatically is a representation of $G$ itself (by
restricting the addition, the scaling and the $G$-action from $V$ to it).
\end{definition}

Recall the natural representation $\mathbf{k}^{n}$ of the symmetric group
$S_{n}$ (defined in Example \ref{exa.rep.Sn-rep.nat}). Let us study its
subrepresentations. For convenience's sake, let us define a notation: Let
$\left(  e_{1},e_{2},\ldots,e_{n}\right)  $ be the standard basis of
$\mathbf{k}^{n}$ (so that $e_{i}\in\mathbf{k}^{n}$ is the $n$-tuple $\left(
0,0,\ldots,0,1,0,0,\ldots,0\right)  $ whose $i$-th entry is $1$ and whose all
other entries are $0$). The action of the symmetric group $S_{n}$ then
satisfies%
\begin{equation}
g\rightharpoonup e_{i}=e_{g\left(  i\right)  }\ \ \ \ \ \ \ \ \ \ \text{for
all }i\in\left[  n\right]  . \label{eq.rep.Sn-rep.nat.gei}%
\end{equation}

Now, what subrepresentations does the natural representation $\mathbf{k}^{n}$
of $S_{n}$ have? Here are four:

\begin{itemize}
\item One subrepresentation of $\mathbf{k}^{n}$ is the \emph{diagonal
subrepresentation}%
\begin{align*}
D\left(  \mathbf{k}^{n}\right)  :=  &  \ \left\{  \left(  a_{1},a_{2}%
,\ldots,a_{n}\right)  \in\mathbf{k}^{n}\ \mid\ a_{1}=a_{2}=\cdots
=a_{n}\right\} \\
=  &  \ \left\{  \left(  a,a,\ldots,a\right)  \ \mid\ a\in\mathbf{k}\right\}
\end{align*}
(we understand $\left(  a,a,\ldots,a\right)  $ to be an $n$-tuple here). If
$n>0$, then this is a free $\mathbf{k}$-module of rank\footnote{The
\emph{rank} of a free $\mathbf{k}$-module is the number of elements in its
basis. When $\mathbf{k}$ is a field, this is the classical concept of the
dimension of a $\mathbf{k}$-vector space. I don't know why the word
\textquotedblleft dimension\textquotedblright\ is not commonly used for
arbitrary $\mathbf{k}$, but the custom is to call it \textquotedblleft
rank\textquotedblright.} $1$, with basis
\[
\left(  \left(  1,1,\ldots,1\right)  \right)  =\left(  e_{1}+e_{2}%
+\cdots+e_{n}\right)  .
\]
In geometric terms, $D\left(  \mathbf{k}^{n}\right)  $ is a line through the
origin.\footnote{These \textquotedblleft geometric terms\textquotedblright%
\ are, of course, not to be taken too literally; after all, $\mathbf{k}$ is
just a commutative ring here.} As an $S_{n}$-representation, $D\left(
\mathbf{k}^{n}\right)  $ is isomorphic to the trivial representation
$\mathbf{k}_{\operatorname*{triv}}$ of $S_{n}$ (defined in Example
\ref{exa.rep.Sn-rep.triv})\footnote{again assuming that $n>0$}, since the map%
\begin{align*}
\mathbf{k}_{\operatorname*{triv}}  &  \rightarrow D\left(  \mathbf{k}%
^{n}\right)  ,\\
a  &  \mapsto\left(  a,a,\ldots,a\right)
\end{align*}
is an isomorphism of representations (why?).

\item Another subrepresentation of $\mathbf{k}^{n}$ is the \emph{reflection
subrepresentation} (aka \emph{zero-sum subrepresentation})%
\[
R\left(  \mathbf{k}^{n}\right)  :=\left\{  \left(  a_{1},a_{2},\ldots
,a_{n}\right)  \in\mathbf{k}^{n}\ \mid\ a_{1}+a_{2}+\cdots+a_{n}=0\right\}  .
\]
It consists of all $n$-tuples in $\mathbf{k}^{n}$ whose entries sum to $0$.
This is a free $\mathbf{k}$-module of rank $n-1$ (unless $n=0$), with basis%
\[
\left(  e_{1}-e_{n},\ e_{2}-e_{n},\ \ldots,\ e_{n-1}-e_{n}\right)
\]
(prove this, if you have not seen this before!). Geometrically, $R\left(
\mathbf{k}^{n}\right)  $ is the hyperplane through the origin perpendicular to
the line $D\left(  \mathbf{k}^{n}\right)  $ (with respect to the dot product).

\item Another subrepresentation of $\mathbf{k}^{n}$ is the zero submodule
$\left\{  0\right\}  $. (In fact, for any representation $V$ of any group $G$,
its zero submodule $\left\{  0_{V}\right\}  $ is a subrepresentation of $V$.)

\item Another subrepresentation of $\mathbf{k}^{n}$ is $\mathbf{k}^{n}$
itself. (In fact, any representation $V$ has itself as a subrepresentation.)
\end{itemize}

Can we find any others?

For $\mathbf{k}=\mathbb{Z}$, we can find a lot. For instance, for any integers
$u,v,w$, the subset%
\begin{align*}
R_{u,v,w}\left(  \mathbb{Z}^{n}\right)   &  :=\Big\{\ \left(  a_{1}%
,a_{2},\ldots,a_{n}\right)  \in\mathbb{Z}^{n}\ \mid\ a_{1}\equiv a_{2}%
\equiv\cdots\equiv a_{n}\equiv0\operatorname{mod}u\\
&  \ \ \ \ \ \ \ \ \ \ \ \ \ \ \ \ \ \ \ \ \text{and }a_{1}\equiv a_{2}%
\equiv\cdots\equiv a_{n}\operatorname{mod}v\\
&  \ \ \ \ \ \ \ \ \ \ \ \ \ \ \ \ \ \ \ \ \text{and }a_{1}+a_{2}+\cdots
+a_{n}\equiv0\operatorname{mod}w\ \Big\}
\end{align*}
of $\mathbb{Z}^{n}$ is a subrepresentation of the natural $S_{n}%
$-representation $\mathbb{Z}^{n}$. (Keep in mind that \textquotedblleft%
$a\equiv b\operatorname{mod}0$\textquotedblright\ means \textquotedblleft%
$a=b$\textquotedblright, so that we can encode any equalities as congruences.)
Are there even more? I don't know.

When $\mathbf{k}$ is a field, however, the answer is much simpler:

\begin{theorem}
\label{thm.rep.G-rep.Sn-nat.subreps}Assume that $\mathbf{k}$ is a field. Then,
the only subrepresentations of the natural $S_{n}$-representation
$\mathbf{k}^{n}$ are $D\left(  \mathbf{k}^{n}\right)  $, $R\left(
\mathbf{k}^{n}\right)  $, $\left\{  0\right\}  $ and $\mathbf{k}^{n}$.
\end{theorem}

\begin{proof}
[Proof idea.]Let $W$ be a subrepresentation of $\mathbf{k}^{n}$. If
$W\subseteq D\left(  \mathbf{k}^{n}\right)  $, then $W$ is either $D\left(
\mathbf{k}^{n}\right)  $ or $\left\{  0\right\}  $ (why?). If not, then $W$
contains some vector
\[
w=\left(  w_{1},w_{2},\ldots,w_{n}\right)
\]
with two distinct consecutive entries (why?). Pick such a $w=\left(
w_{1},w_{2},\ldots,w_{n}\right)  $, and pick some $i\in\left[  n-1\right]  $
such that $w_{i}\neq w_{i+1}$. Then, since $W$ is a subrepresentation, we must
have $s_{i}\rightharpoonup w\in W$, and thus $w-s_{i}\rightharpoonup w\in W$.
But%
\begin{align*}
w-s_{i}\rightharpoonup w  &  =\left(  w_{1},w_{2},\ldots,w_{n}\right)
-\left(  w_{1},w_{2},\ldots,w_{i-1},w_{i+1},w_{i},w_{i+2},\ldots,w_{n}\right)
\\
&  =\left(  0,0,\ldots,0,w_{i}-w_{i+1},w_{i+1}-w_{i},0,0,\ldots,0\right) \\
&  =\left(  w_{i}-w_{i+1}\right)  \left(  e_{i}-e_{i+1}\right)
\end{align*}
(where $e_{1},e_{2},\ldots,e_{n}$ are the standard basis vectors again).
Thus,
\[
\left(  w_{i}-w_{i+1}\right)  \left(  e_{i}-e_{i+1}\right)  =w-s_{i}%
\rightharpoonup w\in W.
\]
Hence, $e_{i}-e_{i+1}\in W$ (since $w_{i}-w_{i+1}$ is nonzero and therefore
invertible in the field $\mathbf{k}$). From this, we can easily derive that
$e_{i}-e_{n}\in W$ (why?). This, in turn, easily entails that $e_{k}-e_{n}\in
W$ for each $k\in\left[  n-1\right]  $ (why?). From this, we obtain $R\left(
\mathbf{k}^{n}\right)  \subseteq W$ (why?). But this entails that $W$ is
either $R\left(  \mathbf{k}^{n}\right)  $ or $\mathbf{k}^{n}$ (why?). Thus,
Theorem \ref{thm.rep.G-rep.Sn-nat.subreps} is proved.
\end{proof}

\begin{exercise}
\fbox{3} Fill in the \textquotedblleft why?\textquotedblright s in the above proof.
\end{exercise}

As mentioned above, when $\mathbf{k}=\mathbb{R}$, we can visualize
$\mathbf{k}^{n}=\mathbb{R}^{n}$ geometrically as $n$-dimensional space, and
then the subrepresentations $D\left(  \mathbf{k}^{n}\right)  $ and $R\left(
\mathbf{k}^{n}\right)  $ become (respectively) a line and a hyperplane
orthogonal to it. This geometric picture is a good guide for the case
$\mathbf{k}=\mathbb{R}$, but it can mislead when $\mathbf{k}$ is (e.g.) a
finite field. Indeed, you might expect the orthogonality to ensure that
$R\left(  \mathbf{k}^{n}\right)  $ and $D\left(  \mathbf{k}^{n}\right)  $ are
complementary subspaces (or $\mathbf{k}$-submodules, to be more precise); but
this reasoning only works when $\mathbf{k}$ is an ordered field like
$\mathbb{R}$ or $\mathbb{Q}$. Thus, the question of how $R\left(
\mathbf{k}^{n}\right)  $ and $D\left(  \mathbf{k}^{n}\right)  $ interact in
general needs to be answered algebraically. Here is a partial answer:

\begin{proposition}
\label{prop.rep.G-rep.Sn-nat.subreps.dirsum}\textbf{(a)} If $n$ is invertible
in $\mathbf{k}$, then
\[
\mathbf{k}^{n}=R\left(  \mathbf{k}^{n}\right)  \oplus D\left(  \mathbf{k}%
^{n}\right)  \ \ \ \ \ \ \ \ \ \ \left(  \text{as }\mathbf{k}\text{-modules}%
\right)  .
\]
That is, $R\left(  \mathbf{k}^{n}\right)  $ and $D\left(  \mathbf{k}%
^{n}\right)  $ are complementary $\mathbf{k}$-submodules of $\mathbf{k}^{n}$
in this case. \medskip

\textbf{(b)} If $n=0$ in $\mathbf{k}$, then $D\left(  \mathbf{k}^{n}\right)
\subseteq R\left(  \mathbf{k}^{n}\right)  $.
\end{proposition}

\begin{exercise}
\fbox{2} Prove this proposition.
\end{exercise}

Note that Proposition \ref{prop.rep.G-rep.Sn-nat.subreps.dirsum} does not
cover all possibilities, since there can be commutative rings $\mathbf{k}$ in
which $n$ is neither invertible nor $0$ (for example, $2$ is neither
invertible nor zero in $\mathbb{Z}/4$ or $\mathbb{Z}/6$). But any field
$\mathbf{k}$ falls under either part \textbf{(a)} or part \textbf{(b)} of
Proposition \ref{prop.rep.G-rep.Sn-nat.subreps.dirsum}, since any element of a
field is either invertible or zero.

\subsubsection{Quotients}

We know what quotients of $\mathbf{k}$-modules are, and we know (from
Definition \ref{def.rep.G-sets.eqrel} \textbf{(b)}) what quotients of $G$-sets
are. Thus we can easily define quotients of representations:

\begin{definition}
\label{def.rep.G-rep.quot}Let $G$ be a group. Let $V$ be a representation of
$G$ over $\mathbf{k}$. Let $W$ be a subrepresentation of $V$. Consider the
quotient $\mathbf{k}$-module%
\begin{align*}
V/W  &  =\left\{  \text{cosets of }W\text{ in }V\right\} \\
&  =\left\{  \text{residue classes of vectors in }V\text{ modulo }W\right\}  .
\end{align*}
(\textquotedblleft Residue class\textquotedblright\ is just a synonym for
\textquotedblleft coset\textquotedblright\ here.) Let us write $\overline{v}$
for the coset $v+W\in V/W$ whenever $v\in V$.

Now, the $\mathbf{k}$-module $V/W$ becomes a representation of $G$ in a
natural way, by defining the left $G$-action on $V/W$ by%
\[
g\rightharpoonup\overline{v}=\overline{g\rightharpoonup v}%
\ \ \ \ \ \ \ \ \ \ \text{for all }g\in G\text{ and }v\in V.
\]
Alternatively, this left $G$-action can be obtained as a particular case of
Definition \ref{def.rep.G-sets.eqrel} \textbf{(b)}: Namely, let
$\overset{W}{\equiv}$ be the binary relation \textquotedblleft congruent
modulo $W$\textquotedblright\ on $V$ (that is, the relation defined by
\[
\left(  v_{1}\overset{W}{\equiv}v_{2}\right)  \ \Longleftrightarrow\ \left(
v_{1}-v_{2}\in W\right)
\]
for all $v_{1},v_{2}\in V$). This relation $\overset{W}{\equiv}$ is an
equivalence relation and is $G$-invariant (since $v_{1}\overset{W}{\equiv
}v_{2}$ yields $v_{1}-v_{2}\in W$ and thus $g\rightharpoonup\left(
v_{1}-v_{2}\right)  \in W$ for each $g\in G$, so that
\begin{align*}
g\rightharpoonup v_{1}-g\rightharpoonup v_{2}  &  =g\rightharpoonup\left(
v_{1}-v_{2}\right)  \ \ \ \ \ \ \ \ \ \ \left(  \text{by the }\mathbf{k}%
\text{-linearity of }\rightharpoonup\right) \\
&  \in W
\end{align*}
and therefore $g\rightharpoonup v_{1}\overset{W}{\equiv}g\rightharpoonup
v_{2}$). Thus, the quotient set $V/\left.  \overset{W}{\equiv}\right.  $
itself acquires a left $G$-action, according to Definition
\ref{def.rep.G-sets.eqrel} \textbf{(b)}. This is precisely the left $G$-action
on $V/W$ defined above (since $V/W$ is precisely $V/\left.  \overset{W}{\equiv
}\right.  $).
\end{definition}

Now, let us study the quotients of the natural representation $\mathbf{k}^{n}$
of $S_{n}$ by its two subrepresentations $D\left(  \mathbf{k}^{n}\right)  $
and $R\left(  \mathbf{k}^{n}\right)  $. Here, we are using the notations
introduced in Subsection \ref{subsec.rep.G-rep.subreps}.

\begin{theorem}
\label{thm.rep.G-rep.Sn-nat.quots}\textbf{(a)} We have $\mathbf{k}%
^{n}/D\left(  \mathbf{k}^{n}\right)  \cong R\left(  \mathbf{k}^{n}\right)  $
as $S_{n}$-representations if $n$ is invertible in $\mathbf{k}$. \medskip

\textbf{(b)} We always have $\mathbf{k}^{n}/D\left(  \mathbf{k}^{n}\right)
\cong R\left(  \mathbf{k}^{n}\right)  $ as $\mathbf{k}$-modules. \medskip

\textbf{(c)} For $n\geq3$, we don't always have $\mathbf{k}^{n}/D\left(
\mathbf{k}^{n}\right)  \cong R\left(  \mathbf{k}^{n}\right)  $ as $S_{n}%
$-representations. \medskip

\textbf{(d)} We always have $\mathbf{k}^{n}/R\left(  \mathbf{k}^{n}\right)
\cong D\left(  \mathbf{k}^{n}\right)  $ as $S_{n}$-representations.
\end{theorem}

\begin{proof}
[Proof of Theorem \ref{thm.rep.G-rep.Sn-nat.quots}.]WLOG assume that $n\neq0$
(since the $n=0$ case is obvious). Let $\left(  e_{1},e_{2},\ldots
,e_{n}\right)  $ be the standard basis of $\mathbf{k}^{n}$.

Some intuition first: We have%
\[
D\left(  \mathbf{k}^{n}\right)  =\operatorname*{span}\nolimits_{\mathbf{k}%
}\left\{  \left(  1,1,\ldots,1\right)  \right\}  =\left\{  \left(
a,a,\ldots,a\right)  \ \mid\ a\in\mathbf{k}\right\}
\]
(where the tuples are $n$-tuples). The elements of $\mathbf{k}^{n}/D\left(
\mathbf{k}^{n}\right)  $ are residue classes of vectors in $\mathbf{k}^{n}$
modulo this submodule $D\left(  \mathbf{k}^{n}\right)  $. Thus, roughly
speaking, they are vectors in $\mathbf{k}^{n}$ defined up to adding a vector
of the form $\left(  a,a,\ldots,a\right)  $ with $a\in\mathbf{k}$. In other
words, they are vectors in $\mathbf{k}^{n}$ defined up to adding the same
scalar to each entry. \medskip

In the following, for any $v\in\mathbf{k}^{n}$, we shall let $\overline{v}$
denote the residue class of $v$ in $\mathbf{k}^{n}/D\left(  \mathbf{k}%
^{n}\right)  $ or in $\mathbf{k}^{n}/R\left(  \mathbf{k}^{n}\right)  $
depending on the context. (Specifically, in the proofs of part \textbf{(a)},
\textbf{(b)} and \textbf{(c)}, it will mean the residue class in
$\mathbf{k}^{n}/D\left(  \mathbf{k}^{n}\right)  $, whereas in the proof of
part \textbf{(d)} it will mean the residue class in $\mathbf{k}^{n}/R\left(
\mathbf{k}^{n}\right)  $.) \medskip

\textbf{(a)} Assume that $n$ is invertible in $\mathbf{k}$. Consider the map%
\begin{align*}
f:R\left(  \mathbf{k}^{n}\right)   &  \rightarrow\mathbf{k}^{n}/D\left(
\mathbf{k}^{n}\right)  ,\\
v  &  \mapsto\overline{v}.
\end{align*}
This map $f$ is $\mathbf{k}$-linear (obviously) and $S_{n}$-equivariant (since
$g\rightharpoonup\overline{v}=\overline{g\rightharpoonup v}$ for all $g\in
S_{n}$ and $v\in\mathbf{k}^{n}$), so it is a morphism of $S_{n}$%
-representations. Moreover:

\begin{itemize}
\item The map $f$ is injective.

[\textit{Proof:} Let $v$ lie in the kernel of the $\mathbf{k}$-linear map $f$.
Thus, $v\in R\left(  \mathbf{k}^{n}\right)  $ and $\overline{v}=\overline{0}$.
From $v\in D\left(  \mathbf{k}^{n}\right)  $, we obtain $v=\left(
a,a,\ldots,a\right)  $ for some $a\in\mathbf{k}$. Consider this $a$. Then,
$\left(  a,a,\ldots,a\right)  =v\in R\left(  \mathbf{k}^{n}\right)  $, so that
$\underbrace{a+a+\cdots+a}_{n\text{ times}}=0$ (by the definition of $R\left(
\mathbf{k}^{n}\right)  $). Hence, $na=\underbrace{a+a+\cdots+a}_{n\text{
times}}=0$. Since $n$ is invertible in $\mathbf{k}$, we can multiply this
equality by $n^{-1}$ and obtain $a=0$. Hence, we can rewrite $v=\left(
a,a,\ldots,a\right)  $ as $v=\left(  0,0,\ldots,0\right)  =0$.

Forget that we fixed $v$. We thus have shown that each $v$ in the kernel of
the $\mathbf{k}$-linear map $f$ is $0$. In other words, the kernel of the
$\mathbf{k}$-linear map $f$ is $\left\{  0\right\}  $. Hence, $f$ is injective.]

\item The map $f$ is surjective.

[\textit{Proof:} Let $x\in\mathbf{k}^{n}/D\left(  \mathbf{k}^{n}\right)  $. We
want to find a $y\in R\left(  \mathbf{k}^{n}\right)  $ such that $f\left(
y\right)  =x$.

We have $x\in\mathbf{k}^{n}/D\left(  \mathbf{k}^{n}\right)  $, so that
$x=\overline{v}$ for some $v\in\mathbf{k}^{n}$. Consider this $v$. This vector
$v$ may or may not lie in $R\left(  \mathbf{k}^{n}\right)  $. However, we can
find a vector $y\in R\left(  \mathbf{k}^{n}\right)  $ such that $\overline
{y}=\overline{v}$. To this purpose, we write $v$ as $v=\left(  v_{1}%
,v_{2},\ldots,v_{n}\right)  $, and we let $v_{\varnothing}:=\dfrac{1}%
{n}\left(  v_{1}+v_{2}+\cdots+v_{n}\right)  $ denote the average entry of $v$
(this is well-defined, since $n$ is invertible in $\mathbf{k}$). Then, we set
$y:=v-\left(  v_{\varnothing},v_{\varnothing},\ldots,v_{\varnothing}\right)
$. Then, $v-y=\left(  v_{\varnothing},v_{\varnothing},\ldots,v_{\varnothing
}\right)  \in D\left(  \mathbf{k}^{n}\right)  $, so that $\overline
{v}=\overline{y}$ and therefore $\overline{y}=\overline{v}=x$. Moreover, from%
\begin{align*}
y  &  =\underbrace{v}_{=\left(  v_{1},v_{2},\ldots,v_{n}\right)  }-\left(
v_{\varnothing},v_{\varnothing},\ldots,v_{\varnothing}\right)  =\left(
v_{1},v_{2},\ldots,v_{n}\right)  -\left(  v_{\varnothing},v_{\varnothing
},\ldots,v_{\varnothing}\right) \\
&  =\left(  v_{1}-v_{\varnothing},\ v_{2}-v_{\varnothing},\ \ldots
,\ v_{n}-v_{\varnothing}\right)  ,
\end{align*}
we see that the sum of all entries of $y$ is%
\begin{align*}
&  \left(  v_{1}-v_{\varnothing}\right)  +\left(  v_{2}-v_{\varnothing
}\right)  +\cdots+\left(  v_{n}-v_{\varnothing}\right) \\
&  =\left(  v_{1}+v_{2}+\cdots+v_{n}\right)  -n\underbrace{v_{\varnothing}%
}_{=\dfrac{1}{n}\left(  v_{1}+v_{2}+\cdots+v_{n}\right)  }\\
&  =\left(  v_{1}+v_{2}+\cdots+v_{n}\right)  -n\cdot\dfrac{1}{n}\left(
v_{1}+v_{2}+\cdots+v_{n}\right)  =0.
\end{align*}
Thus, $y\in R\left(  \mathbf{k}^{n}\right)  $. Hence, $f\left(  y\right)  $ is
well-defined and equals $f\left(  y\right)  =\overline{y}=x$.

Forget that we fixed $x$. Thus, for each $x\in\mathbf{k}^{n}/D\left(
\mathbf{k}^{n}\right)  $, we have found a $y\in R\left(  \mathbf{k}%
^{n}\right)  $ such that $f\left(  y\right)  =x$. In other words, we have
shown that the map $f$ is surjective.]
\end{itemize}

Now we know that the map $f$ is both injective and surjective, thus bijective,
i.e., invertible. Hence, $f$ is an invertible morphism of $S_{n}%
$-representations, and therefore an isomorphism of $S_{n}$-representations (by
Proposition \ref{prop.rep.G-rep.iso=bij}). This proves Theorem
\ref{thm.rep.G-rep.Sn-nat.quots} \textbf{(a)}. \medskip

\textbf{(b)} This is an exercise in linear algebra. We shall show that the two
$\mathbf{k}$-modules $\mathbf{k}^{n}/D\left(  \mathbf{k}^{n}\right)  $ and
$R\left(  \mathbf{k}^{n}\right)  $ are free of the same rank. This will
automatically entail their isomorphism (as $\mathbf{k}$-modules, not as
$S_{n}$-representations, of course), because any two free $\mathbf{k}$-modules
of the same rank are isomorphic\footnote{\textit{Proof.} Let $M$ and $N$ be
two free $\mathbf{k}$-modules of the same rank $k$. Then, $M$ and $N$ have
bases $\left(  m_{1},m_{2},\ldots,m_{k}\right)  $ and $\left(  n_{1}%
,n_{2},\ldots,n_{k}\right)  $, respectively. Thus, we can define a
$\mathbf{k}$-linear map $f:M\rightarrow N$ that sends each $m_{i}$ to the
corresponding $n_{i}$. This $\mathbf{k}$-linear map $f$ is easily seen to be
an isomorphism. Hence, $M$ and $N$ are isomorphic.}.

We already know that the $\mathbf{k}$-module $R\left(  \mathbf{k}^{n}\right)
$ is free of rank $n-1$, with basis $\left(  e_{1}-e_{n},\ e_{2}%
-e_{n},\ \ldots,\ e_{n-1}-e_{n}\right)  $.

Now, we claim that the $\mathbf{k}$-module $\mathbf{k}^{n}/D\left(
\mathbf{k}^{n}\right)  $ is free of rank $n-1$, with basis $\left(
\overline{e_{1}},\overline{e_{2}},\ldots,\overline{e_{n-1}}\right)  $. Indeed:

\begin{itemize}
\item The residue classes $\overline{e_{1}},\overline{e_{2}},\ldots
,\overline{e_{n-1}}$ in $\mathbf{k}^{n}/D\left(  \mathbf{k}^{n}\right)  $ are
$\mathbf{k}$-linearly independent.

[\textit{Proof:} Let $a_{1},a_{2},\ldots,a_{n-1}\in\mathbf{k}$ be scalars such
that $a_{1}\overline{e_{1}}+a_{2}\overline{e_{2}}+\cdots+a_{n-1}%
\overline{e_{n-1}}=\overline{0}$. We must show that $a_{1}=a_{2}%
=\cdots=a_{n-1}=0$.

By assumption, we have $a_{1}\overline{e_{1}}+a_{2}\overline{e_{2}}%
+\cdots+a_{n-1}\overline{e_{n-1}}=\overline{0}$. In view of%
\begin{align*}
a_{1}\overline{e_{1}}+a_{2}\overline{e_{2}}+\cdots+a_{n-1}\overline{e_{n-1}}
&  =\overline{a_{1}e_{1}+a_{2}e_{2}+\cdots+a_{n-1}e_{n-1}}\\
&  =\overline{\left(  a_{1},a_{2},\ldots,a_{n-1},0\right)  },
\end{align*}
we can rewrite this as $\overline{\left(  a_{1},a_{2},\ldots,a_{n-1},0\right)
}=\overline{0}$. In other words, \newline$\left(  a_{1},a_{2},\ldots
,a_{n-1},0\right)  \in D\left(  \mathbf{k}^{n}\right)  $. By the definition of
$D\left(  \mathbf{k}^{n}\right)  $, this means that $a_{1}=a_{2}%
=\cdots=a_{n-1}=0$. But this is precisely what we needed to prove. Thus, the
$\mathbf{k}$-linear independence of $\overline{e_{1}},\overline{e_{2}}%
,\ldots,\overline{e_{n-1}}$ in $\mathbf{k}^{n}/D\left(  \mathbf{k}^{n}\right)
$ is proved.]

\item The residue classes $\overline{e_{1}},\overline{e_{2}},\ldots
,\overline{e_{n-1}}$ span the $\mathbf{k}$-module $\mathbf{k}^{n}/D\left(
\mathbf{k}^{n}\right)  $.

[\textit{Proof:} We have%
\[
\left(  \overline{e_{1}}+\overline{e_{2}}+\cdots+\overline{e_{n-1}}\right)
+\overline{e_{n}}=\overline{e_{1}}+\overline{e_{2}}+\cdots+\overline{e_{n}%
}=\overline{e_{1}+e_{2}+\cdots+e_{n}}=\overline{0}%
\]
(since $e_{1}+e_{2}+\cdots+e_{n}=\left(  1,1,\ldots,1\right)  \in D\left(
\mathbf{k}^{n}\right)  $), so that%
\begin{equation}
\overline{e_{n}}=-\left(  \overline{e_{1}}+\overline{e_{2}}+\cdots
+\overline{e_{n-1}}\right)  . \label{pf.thm.rep.G-rep.Sn-nat.quots.b.4}%
\end{equation}
However, the standard basis vectors $e_{1},e_{2},\ldots,e_{n}$ clearly span
the $\mathbf{k}$-module $\mathbf{k}^{n}$. Thus, their residue classes
$\overline{e_{1}},\overline{e_{2}},\ldots,\overline{e_{n}}$ span the quotient
$\mathbf{k}$-module $\mathbf{k}^{n}/D\left(  \mathbf{k}^{n}\right)  $.
Moreover, the last of the latter $n$ classes (that is, $\overline{e_{n}}$) is
actually redundant (since the equality
(\ref{pf.thm.rep.G-rep.Sn-nat.quots.b.4}) expresses it as a $\mathbf{k}%
$-linear combination of the remaining $n-1$ classes $\overline{e_{1}%
},\overline{e_{2}},\ldots,\overline{e_{n-1}}$). Hence, the other residue
classes $\overline{e_{1}},\overline{e_{2}},\ldots,\overline{e_{n-1}}$ already
span the $\mathbf{k}$-module $\mathbf{k}^{n}/D\left(  \mathbf{k}^{n}\right)  $.]
\end{itemize}

Thus, we know that the residue classes $\overline{e_{1}},\overline{e_{2}%
},\ldots,\overline{e_{n-1}}$ are $\mathbf{k}$-linearly independent and span
the $\mathbf{k}$-module $\mathbf{k}^{n}/D\left(  \mathbf{k}^{n}\right)  $.
Hence, they form a basis of $\mathbf{k}^{n}/D\left(  \mathbf{k}^{n}\right)  $.
Thus, the $\mathbf{k}$-module $\mathbf{k}^{n}/D\left(  \mathbf{k}^{n}\right)
$ is free of rank $n-1$.

Now, we know that the two $\mathbf{k}$-modules $\mathbf{k}^{n}/D\left(
\mathbf{k}^{n}\right)  $ and $R\left(  \mathbf{k}^{n}\right)  $ are both free
of rank $n-1$. Hence, they are free of the same rank, and thus isomorphic.
This proves Theorem \ref{thm.rep.G-rep.Sn-nat.quots} \textbf{(b)}. \medskip

\textbf{(c)} Assume that $n\geq3$, and that $\mathbf{k}$ is a field of
characteristic $p$ with $p\mid n$. We shall show that the $S_{n}%
$-representations $\mathbf{k}^{n}/D\left(  \mathbf{k}^{n}\right)  $ and
$R\left(  \mathbf{k}^{n}\right)  $ are not isomorphic.

Note that it does \textbf{not} suffice that the map $f$ from the above proof
of Theorem \ref{thm.rep.G-rep.Sn-nat.quots} \textbf{(a)} fails to be an
isomorphism, or that the $\mathbf{k}$-module isomorphism from Theorem
\ref{thm.rep.G-rep.Sn-nat.quots} \textbf{(b)} fails to be $S_{n}$-equivariant.
After all, perhaps these two maps are just the wrong ones for the job, but
some other map would qualify?

We must instead show that \textbf{no} $S_{n}$-representation isomorphism from
$\mathbf{k}^{n}/D\left(  \mathbf{k}^{n}\right)  $ to $R\left(  \mathbf{k}%
^{n}\right)  $ exists.

Here is a general trick for proving non-isomorphy:

If $G$ is any group and $V$ is any representation of $G$, then we set%
\[
V^{G-\operatorname*{fix}}:=\left\{  v\in V\ \mid\ g\rightharpoonup v=v\text{
for all }g\in G\right\}  .
\]
This is called the \emph{space of }$G$\emph{-fixed points in }$V$. It is a
$\mathbf{k}$-submodule of $V$ and even a $G$-subrepresentation of $V$ (a
trivial one, by its definition). Its elements are called the $G$\emph{-fixed
points in }$V$. (Classically, it is denoted by $V^{G}$, but this notation can
also mean the set of all maps from $G$ to $V$, which is different.)

Obviously, if $V$ and $W$ are two isomorphic representations of $G$, then
$V^{G-\operatorname*{fix}}\cong W^{G-\operatorname*{fix}}$ as well (as
$\mathbf{k}$-modules or as $G$-representations, as you prefer). Thus, if two
$G$-representations have non-isomorphic spaces of $G$-fixed points, then they
themselves are not isomorphic.

We shall now show that the $S_{n}$-representations $\mathbf{k}^{n}/D\left(
\mathbf{k}^{n}\right)  $ and $R\left(  \mathbf{k}^{n}\right)  $ have
non-isomorphic spaces of $S_{n}$-fixed points:

\begin{itemize}
\item The only $S_{n}$-fixed point in $\mathbf{k}^{n}/D\left(  \mathbf{k}%
^{n}\right)  $ is the zero vector $\overline{0}$.

[\textit{Proof:} Clearly, the zero vector $\overline{0}$ is an $S_{n}$-fixed
point. Now let us prove the converse. Let $\overline{\left(  a_{1}%
,a_{2},\ldots,a_{n}\right)  }\in\mathbf{k}^{n}/D\left(  \mathbf{k}^{n}\right)
$ be an $S_{n}$-fixed point. Then, by the definition of an $S_{n}$-fixed
point, we have%
\[
s_{1}\rightharpoonup\overline{\left(  a_{1},a_{2},\ldots,a_{n}\right)
}=\overline{\left(  a_{1},a_{2},\ldots,a_{n}\right)  }.
\]
In view of%
\[
s_{1}\rightharpoonup\overline{\left(  a_{1},a_{2},\ldots,a_{n}\right)
}=\overline{s_{1}\rightharpoonup\left(  a_{1},a_{2},\ldots,a_{n}\right)
}=\overline{\left(  a_{2},a_{1},a_{3},a_{4},\ldots,a_{n}\right)  },
\]
this rewrites as
\[
\overline{\left(  a_{2},a_{1},a_{3},a_{4},\ldots,a_{n}\right)  }%
=\overline{\left(  a_{1},a_{2},\ldots,a_{n}\right)  }.
\]
In other words,%
\[
\left(  a_{2},a_{1},a_{3},a_{4},\ldots,a_{n}\right)  -\left(  a_{1}%
,a_{2},\ldots,a_{n}\right)  \in D\left(  \mathbf{k}^{n}\right)  .
\]
In view of%
\begin{align*}
&  \left(  a_{2},a_{1},a_{3},a_{4},\ldots,a_{n}\right)  -\left(  a_{1}%
,a_{2},\ldots,a_{n}\right) \\
&  =\left(  a_{2}-a_{1},\ a_{1}-a_{2},\ 0,\ 0,\ \ldots,\ 0\right)  ,
\end{align*}
this rewrites further as%
\[
\left(  a_{2}-a_{1},\ a_{1}-a_{2},\ 0,\ 0,\ \ldots,\ 0\right)  \in D\left(
\mathbf{k}^{n}\right)  .
\]
By the definition of $D\left(  \mathbf{k}^{n}\right)  $, we can rewrite this
as%
\[
a_{2}-a_{1}=a_{1}-a_{2}=0=0=\cdots=0.
\]
Since $n\geq3$, there is at least one $0$ in this equality, so we obtain
$a_{1}-a_{2}=0$ and therefore $a_{1}=a_{2}$. A similar argument (using $s_{2}$
instead of $s_{1}$) yields $a_{2}=a_{3}$, and likewise we find $a_{3}=a_{4}$
and $a_{4}=a_{5}$ and so on. Combining all these equalities, we obtain%
\[
a_{1}=a_{2}=\cdots=a_{n},
\]
so that $\left(  a_{1},a_{2},\ldots,a_{n}\right)  \in D\left(  \mathbf{k}%
^{n}\right)  $ and thus $\overline{\left(  a_{1},a_{2},\ldots,a_{n}\right)
}=\overline{0}$.

Forget that we fixed $\overline{\left(  a_{1},a_{2},\ldots,a_{n}\right)  }$.
We thus have shown that any $S_{n}$-fixed point $\overline{\left(  a_{1}%
,a_{2},\ldots,a_{n}\right)  }\in\mathbf{k}^{n}/D\left(  \mathbf{k}^{n}\right)
$ equals $\overline{0}$. In other words, the only $S_{n}$-fixed point in
$\mathbf{k}^{n}/D\left(  \mathbf{k}^{n}\right)  $ is the zero vector
$\overline{0}$.]

\item I claim that the vector $\left(  1,1,\ldots,1\right)  \in\mathbf{k}^{n}$
is a nonzero $S_{n}$-fixed point in $R\left(  \mathbf{k}^{n}\right)  $.

[\textit{Proof:} The sum of all entries of the vector $\left(  1,1,\ldots
,1\right)  \in\mathbf{k}^{n}$ is%
\[
\underbrace{1+1+\cdots+1}_{n\text{ times}}=n\cdot
1=0\ \ \ \ \ \ \ \ \ \ \text{in }\mathbf{k}%
\]
(since $\operatorname*{char}\mathbf{k}=p\mid n$). Thus, this vector $\left(
1,1,\ldots,1\right)  $ belongs to $R\left(  \mathbf{k}^{n}\right)  $ (by the
definition of $R\left(  \mathbf{k}^{n}\right)  $). Moreover, this vector
$\left(  1,1,\ldots,1\right)  $ is nonzero (since $\mathbf{k}$ is a field, so
that $1\neq0$ in $\mathbf{k}$), and is an $S_{n}$-fixed point (since the
action of a permutation $g\in S_{n}$ just permutes its entries, but all its
entries are equal and therefore do not change under permutation). Hence,
$\left(  1,1,\ldots,1\right)  \in\mathbf{k}^{n}$ is a nonzero $S_{n}$-fixed
point in $R\left(  \mathbf{k}^{n}\right)  $.]
\end{itemize}

Thus, the $S_{n}$-representation $R\left(  \mathbf{k}^{n}\right)  $ has a
nonzero $S_{n}$-fixed point (namely, $\left(  1,1,\ldots,1\right)  $), whereas
the $S_{n}$-representation $\mathbf{k}^{n}/D\left(  \mathbf{k}^{n}\right)  $
does not. Therefore, these two representations cannot be isomorphic. This
proves Theorem \ref{thm.rep.G-rep.Sn-nat.quots} \textbf{(c)}. \medskip

\textbf{(d)} Consider the $n$ residue classes $\overline{e_{1}},\overline
{e_{2}},\ldots,\overline{e_{n}}$ in the quotient $\mathbf{k}$-module
$\mathbf{k}^{n}/R\left(  \mathbf{k}^{n}\right)  $. We claim that all these
classes are equal, i.e., that we have%
\begin{equation}
\overline{e_{1}}=\overline{e_{2}}=\cdots=\overline{e_{n}}.
\label{pf.thm.rep.G-rep.Sn-nat.quots.d.4}%
\end{equation}

[\textit{Proof:} For each $i\in\left\{  2,3,\ldots,n\right\}  $, we have
\[
e_{i}-e_{1}=\left(  -1,0,0,\ldots,0,1,0,0,\ldots,0\right)  \in R\left(
\mathbf{k}^{n}\right)
\]
(since the sum of the entries of the vector $\left(  -1,0,0,\ldots
,0,1,0,0,\ldots,0\right)  $ is $\left(  -1\right)  +1=0$) and thus
$\overline{e_{i}}=\overline{e_{1}}$. In other words, $\overline{e_{1}%
}=\overline{e_{2}}=\cdots=\overline{e_{n}}$. This proves
(\ref{pf.thm.rep.G-rep.Sn-nat.quots.d.4}).] \medskip

We shall now show that the $\mathbf{k}$-module $\mathbf{k}^{n}/R\left(
\mathbf{k}^{n}\right)  $ is free of rank $1$, with basis $\left(
\overline{e_{1}}\right)  $. Indeed:

\begin{itemize}
\item The residue class $\overline{e_{1}}$ in $\mathbf{k}^{n}/R\left(
\mathbf{k}^{n}\right)  $ is $\mathbf{k}$-linearly independent.

[\textit{Proof:} Let $a_{1}\in\mathbf{k}$ be a scalar such that $a_{1}%
\overline{e_{1}}=\overline{0}$. We must show that $a_{1}=0$.

By assumption, we have $a_{1}\overline{e_{1}}=\overline{0}$. In view of%
\[
a_{1}\overline{e_{1}}=\overline{a_{1}e_{1}}=\overline{\left(  a_{1}%
,0,0,\ldots,0\right)  },
\]
we can rewrite this as $\overline{\left(  a_{1},0,0,\ldots,0\right)
}=\overline{0}$. In other words, $\left(  a_{1},0,0,\ldots,0\right)  \in
R\left(  \mathbf{k}^{n}\right)  $. By the definition of $R\left(
\mathbf{k}^{n}\right)  $, this means that $a_{1}+0+0+\cdots+0=0$. This
obviously simplifies to $a_{1}=0$. But this is precisely what we needed to
prove. Thus, the $\mathbf{k}$-linear independence of $\overline{e_{1}}$ in
$\mathbf{k}^{n}/R\left(  \mathbf{k}^{n}\right)  $ is proved.]

\item The residue class $\overline{e_{1}}$ spans the $\mathbf{k}$-module
$\mathbf{k}^{n}/R\left(  \mathbf{k}^{n}\right)  $.

[\textit{Proof:} The standard basis vectors $e_{1},e_{2},\ldots,e_{n}$ clearly
span the $\mathbf{k}$-module $\mathbf{k}^{n}$. Thus, their residue classes
$\overline{e_{1}},\overline{e_{2}},\ldots,\overline{e_{n}}$ span the quotient
$\mathbf{k}$-module $\mathbf{k}^{n}/R\left(  \mathbf{k}^{n}\right)  $. But all
these $n$ residue classes are equal (by
(\ref{pf.thm.rep.G-rep.Sn-nat.quots.d.4})), so that we need only one of them
to span $\mathbf{k}^{n}/R\left(  \mathbf{k}^{n}\right)  $. Thus, the single
residue class $\overline{e_{1}}$ already spans the $\mathbf{k}$-module
$\mathbf{k}^{n}/R\left(  \mathbf{k}^{n}\right)  $.]
\end{itemize}

Thus, we have shown that the residue class $\overline{e_{1}}$ is $\mathbf{k}%
$-linearly independent and spans the $\mathbf{k}$-module $\mathbf{k}%
^{n}/R\left(  \mathbf{k}^{n}\right)  $. Hence, it forms a basis of
$\mathbf{k}^{n}/R\left(  \mathbf{k}^{n}\right)  $. Thus, the $\mathbf{k}%
$-module $\mathbf{k}^{n}/R\left(  \mathbf{k}^{n}\right)  $ is free of rank $1$.

Now, we know that the two $\mathbf{k}$-modules $\mathbf{k}^{n}/R\left(
\mathbf{k}^{n}\right)  $ and $D\left(  \mathbf{k}^{n}\right)  $ are both free
of rank $1$. Hence, they are free of the same rank, and thus isomorphic as
$\mathbf{k}$-modules (for the same reason as in the proof of Theorem
\ref{thm.rep.G-rep.Sn-nat.quots} \textbf{(b)}).

Are they also isomorphic as $S_{n}$-representations? To find out, we analyze
their $S_{n}$-actions:

\begin{itemize}
\item The $S_{n}$-action on $\mathbf{k}^{n}/R\left(  \mathbf{k}^{n}\right)  $
is trivial (i.e., we have $g\rightharpoonup v=v$ for all $g\in S_{n}$ and
$v\in\mathbf{k}^{n}/R\left(  \mathbf{k}^{n}\right)  $).

[\textit{Proof:} Let $g\in S_{n}$ and $v\in\mathbf{k}^{n}/R\left(
\mathbf{k}^{n}\right)  $. We must prove that $g\rightharpoonup v=v$.

Since both sides of this equality are $\mathbf{k}$-linear in $v$, we can WLOG
assume that $v=\overline{e_{1}}$ (since $\overline{e_{1}}$ spans the
$\mathbf{k}$-module $\mathbf{k}^{n}/R\left(  \mathbf{k}^{n}\right)  $). Assume
this. Then, $g\rightharpoonup v=g\rightharpoonup\overline{e_{1}}%
=\overline{g\rightharpoonup e_{1}}=\overline{e_{g\left(  1\right)  }}$ (since
(\ref{eq.rep.Sn-rep.nat.gei}) yields $g\rightharpoonup e_{1}=e_{g\left(
1\right)  }$). However, (\ref{pf.thm.rep.G-rep.Sn-nat.quots.d.4}) yields
$\overline{e_{g\left(  1\right)  }}=\overline{e_{1}}$, so that
$g\rightharpoonup v=\overline{e_{g\left(  1\right)  }}=\overline{e_{1}}=v$.
This completes the proof that the $S_{n}$-action on $\mathbf{k}^{n}/R\left(
\mathbf{k}^{n}\right)  $ is trivial.]

\item The $S_{n}$-action on $D\left(  \mathbf{k}^{n}\right)  $ is trivial
(i.e., we have $g\rightharpoonup v=v$ for all $g\in S_{n}$ and $v\in D\left(
\mathbf{k}^{n}\right)  $).

[\textit{Proof:} Let $g\in S_{n}$ and $v\in D\left(  \mathbf{k}^{n}\right)  $.
We must prove that $g\rightharpoonup v=v$.

We have $v\in D\left(  \mathbf{k}^{n}\right)  $, so that $v=\left(
a,a,\ldots,a\right)  $ for some $a\in\mathbf{k}$. Thus, all entries of $v$ are
equal. Hence, permuting the entries of $v$ does not change $v$. Therefore,
$g\rightharpoonup v=v$ (since $g$ acts on $\mathbf{k}^{n}$ by permuting the
entries). This completes the proof that the $S_{n}$-action on $D\left(
\mathbf{k}^{n}\right)  $ is trivial.]
\end{itemize}

Now, recall that the two $\mathbf{k}$-modules $\mathbf{k}^{n}/R\left(
\mathbf{k}^{n}\right)  $ and $D\left(  \mathbf{k}^{n}\right)  $ are isomorphic
as $\mathbf{k}$-modules. In other words, there exists a $\mathbf{k}$-module
isomorphism $f:\mathbf{k}^{n}/R\left(  \mathbf{k}^{n}\right)  \rightarrow
D\left(  \mathbf{k}^{n}\right)  $. Since the $S_{n}$-actions on both
$\mathbf{k}^{n}/R\left(  \mathbf{k}^{n}\right)  $ and $D\left(  \mathbf{k}%
^{n}\right)  $ are trivial, this $f$ is automatically $S_{n}$-equivariant (by
Example \ref{exa.rep.G-sets.mor-triv}, applied to $G=S_{n}$, $X=\mathbf{k}%
^{n}/R\left(  \mathbf{k}^{n}\right)  $ and $Y=D\left(  \mathbf{k}^{n}\right)
$), hence is a morphism of $S_{n}$-representations, and thus is an isomorphism
of $S_{n}$-representations (by Proposition \ref{prop.rep.G-rep.iso=bij}). This
proves Theorem \ref{thm.rep.G-rep.Sn-nat.quots} \textbf{(d)}.
\end{proof}

Theorem \ref{thm.rep.G-rep.Sn-nat.quots} describes $\mathbf{k}^{n}/R\left(
\mathbf{k}^{n}\right)  $ always and $\mathbf{k}^{n}/D\left(  \mathbf{k}%
^{n}\right)  $ in most cases. (It is somewhat surprising that $\mathbf{k}%
^{n}/R\left(  \mathbf{k}^{n}\right)  \cong D\left(  \mathbf{k}^{n}\right)  $
is always true but $\mathbf{k}^{n}/D\left(  \mathbf{k}^{n}\right)  \cong
R\left(  \mathbf{k}^{n}\right)  $ is not!)

\begin{exercise}
Consider the left action of the symmetric group $S_{n}$ on the matrix ring
$\mathbf{k}^{n\times n}$ defined by%
\begin{align*}
&  g\rightharpoonup\left(  a_{i,j}\right)  _{i\in\left[  n\right]
,\ j\in\left[  n\right]  }=\left(  a_{g^{-1}\left(  i\right)  ,g^{-1}\left(
j\right)  }\right)  _{i\in\left[  n\right]  ,\ j\in\left[  n\right]  }\\
&  \ \ \ \ \ \ \ \ \ \ \text{for all }g\in S_{n}\text{ and }\left(
a_{i,j}\right)  _{i\in\left[  n\right]  ,\ j\in\left[  n\right]  }\text{.}%
\end{align*}
(This is the restriction of the action of $S_{n}\times S_{n}$ on
$\mathbf{k}^{n\times n}$ defined in Example \ref{exa.rep.G-sets.knm} along the
morphism $\alpha:S_{n}\rightarrow S_{n}\times S_{n},\ g\mapsto\left(
g,g\right)  $. Unlike the latter action, which permutes the rows and the
columns by two different permutations, this one permutes the rows and the
columns by the \textbf{same} permutation.) This left action is $\mathbf{k}%
$-linear, thus making $\mathbf{k}^{n\times n}$ into a representation of
$S_{n}$. \medskip

\textbf{(a)} \fbox{1} Prove that the set of all symmetric matrices in
$\mathbf{k}^{n\times n}$ is a subrepresentation of $\mathbf{k}^{n\times n}$.
\medskip

\textbf{(b)} \fbox{1} Prove that the set of all matrices in $\mathbf{k}%
^{n\times n}$ whose columns sum to the zero vector is a subrepresentation of
$\mathbf{k}^{n\times n}$. \medskip

\textbf{(c)} \fbox{2} Find four more subrepresentations of $\mathbf{k}%
^{n\times n}$, not counting $\left\{  0\right\}  $ and $\mathbf{k}^{n\times
n}$ itself. (The representations should be different for $\mathbf{k}%
=\mathbb{Q}$ and sufficiently large $n$. Proofs are not required in this part.)
\end{exercise}

\begin{exercise}
Let $V$ be the representation $\mathbf{k}\left[  x_{1},x_{2},\ldots
,x_{n}\right]  _{2}$ of $S_{n}$ constructed in Example
\ref{exa.rep.Sn-rep.polyring} \textbf{(b)} (applied to $d=2$). This
representation $V$ has basis $\left(  x_{i}x_{j}\right)  _{1\leq i\leq j\leq
n}$, with $S_{n}$-action given by%
\[
g\rightharpoonup\left(  x_{i}x_{j}\right)  =x_{g\left(  i\right)  }x_{g\left(
j\right)  }\ \ \ \ \ \ \ \ \ \ \text{for all }i,j\in\left[  n\right]  .
\]
Consider the two $\mathbf{k}$-submodules%
\begin{align*}
W_{1}  &  :=\operatorname*{span}\nolimits_{\mathbf{k}}\left\{  x_{1}^{2}%
,x_{2}^{2},\ldots,x_{n}^{2}\right\}  \ \ \ \ \ \ \ \ \ \ \text{and}\\
W_{2}  &  :=\operatorname*{span}\nolimits_{\mathbf{k}}\left\{  x_{i}%
x_{j}\ \mid\ 1\leq i<j\leq n\right\}
\end{align*}
of $V$. \medskip

\textbf{(a)} \fbox{1} Show that $W_{1}$ and $W_{2}$ are subrepresentations of
$V$ such that $V=W_{1}\oplus W_{2}$ (as a $\mathbf{k}$-module). \medskip

\textbf{(b)} \fbox{1} One of $W_{1}$ and $W_{2}$ is isomorphic to an $S_{n}%
$-representation we have thoroughly studied. Which one? \medskip

\textbf{(c)} \fbox{1} Prove that $W_{1}\cong W_{2}$ (as $S_{n}$%
-representations) when $n=3$. \medskip

\textbf{(d)} \fbox{1} Show that each of $W_{1}$ and $W_{2}$ contains a nonzero
trivial subrepresentation if $n\geq2$ (and if $\mathbf{k}$ is nontrivial).
\end{exercise}

\subsection{Modules over rings, and a different way of thinking about
representations}

\subsubsection{Definitions}

In the preceding section (Corollary \ref{cor.rep.G-rep.curr2} in particular),
we have seen two equivalent ways to view a representation of a group $G$ over
$\mathbf{k}$: as a $\mathbf{k}$-linear left $G$-action on a $\mathbf{k}%
$-module $V$, or as a group morphism $\rho:G\rightarrow\operatorname*{GL}%
\left(  V\right)  $ (its curried form). Soon we will see another: as a left
$\mathbf{k}\left[  G\right]  $-module $V$. This is where the group algebra
$\mathbf{k}\left[  G\right]  $ finally shows its use.

First, let us recall the notion of a left module (copied from \cite[\S 3.1]%
{23wa}):

\begin{definition}
\label{def.mod.leftmod}Let $R$ be a ring. A \emph{left }$R$\emph{-module} (or
a \emph{left module over }$R$) means a set $M$ equipped with

\begin{itemize}
\item a binary operation $+$ (that is, a map from $M\times M$ to $M$) that is
called \emph{addition};

\item an element $0_{M}\in M$ that is called the \emph{zero element} or the
\emph{zero vector} or just the \emph{zero}, and is just denoted by $0$ when no
confusion is likely;

\item a map from $R\times M$ to $M$ that is called the \emph{action of }$R$
\emph{on }$M$ (or \emph{scaling of }$M$ \emph{by }$R$), and is written as
multiplication (i.e., we denote the image of a pair $\left(  r,m\right)  \in
R\times M$ under this map by $rm$ or $r\cdot m$)
\end{itemize}

\noindent such that the following properties (the \textquotedblleft%
\emph{module axioms}\textquotedblright) hold:

\begin{itemize}
\item $\left(  M,+,0\right)  $ is an abelian group.

\item The \emph{right distributivity law} holds: We have $\left(  r+s\right)
m=rm+sm$ for all $r,s\in R$ and $m\in M$.

\item The \emph{left distributivity law} holds: We have $r\left(  m+n\right)
=rm+rn$ for all $r\in R$ and $m,n\in M$.

\item The \emph{associativity law} holds: We have $\left(  rs\right)
m=r\left(  sm\right)  $ for all $r,s\in R$ and $m\in M$.

\item We have $0_{R}m=0_{M}$ for every $m\in M$.

\item We have $r\cdot0_{M}=0_{M}$ for every $r\in R$.

\item We have $1m=m$ for every $m\in M$. (Here, \textquotedblleft%
$1$\textquotedblright\ means the unity of $R$.)
\end{itemize}

When $M$ is a left $R$-module, the elements of $M$ are called \emph{vectors},
and the elements of $R$ are called \emph{scalars}. (Note that this sometimes
flies in the face of standard intuition for scalars and vectors; see Example
\ref{exa.mod.knx1} below.)
\end{definition}

As the name \textquotedblleft left $R$-module\textquotedblright\ suggests,
there is an analogous notion of a \emph{right }$R$\emph{-module}:

\begin{definition}
\label{def.mod.rimod}Let $R$ be a ring. A \emph{right }$R$\emph{-module} is
defined just as a left $R$-module was defined in Definition
\ref{def.mod.leftmod}, but with the following changes:

\begin{itemize}
\item For a right $R$-module $M$, the action is not a map from $R\times M$ to
$M$, but rather a map from $M\times R$ to $M$.

\item Accordingly, we use the notation $mr$ (rather than $rm$) for the image
of a pair $\left(  m,r\right)  $ under this map.

\item The axioms for a right $R$-module are similar to the module axioms for a
left $R$-module, accounting for the different form of the action. For example,
the associativity law for a right $R$-module is saying that $m\left(
rs\right)  =\left(  mr\right)  s$ for all $r,s\in R$ and $m\in M$.
\end{itemize}
\end{definition}

When $R$ is commutative, there is a natural way to turn any left $R$-module
into a right $R$-module and vice versa:

\begin{proposition}
\label{prop.mod.leftmod-is-rightmod}Let $R$ be a commutative ring. Then, we
can make any left $R$-module $M$ into a right $R$-module by setting%
\begin{equation}
mr=rm\ \ \ \ \ \ \ \ \ \ \text{for all }r\in R\text{ and }m\in M.
\label{eq.prop.mod.leftmod-is-rightmod.mr=rm}%
\end{equation}
Similarly, we can make any right $R$-module $M$ into a left $R$-module by
setting%
\begin{equation}
rm=mr\ \ \ \ \ \ \ \ \ \ \text{for all }r\in R\text{ and }m\in M.
\label{eq.prop.mod.leftmod-is-rightmod.rm=mr}%
\end{equation}
These two transformations are mutually inverse, so we shall use them to
identify left $R$-modules with right $R$-modules. Thus, we use the words
\textquotedblleft left $R$-module\textquotedblright\ and \textquotedblleft
right $R$-module\textquotedblright\ interchangeably, and just speak of
\textquotedblleft$R$\emph{-modules}\textquotedblright\ instead (without
specifying whether they are left or right).
\end{proposition}

Of course, this is precisely the notion of $\mathbf{k}$-modules that we have
always been using, since $\mathbf{k}$ is assumed to be commutative.

\begin{remark}
\label{rmk.mod.to-op}Even when the ring $R$ is noncommutative, we can
translate left modules into right modules, but not over the same ring: A left
$R$-module naturally becomes a right $R^{\operatorname*{op}}$-module, and vice
versa, where $R^{\operatorname*{op}}$ is the \emph{opposite ring} of $R$ (that
is, the ring $R$ with its multiplication reversed -- i.e., the ring $R$ in
which all products $rs$ have been renamed as $sr$).
\end{remark}

When $R$ is a field, the $R$-modules are also known as $R$\emph{-vector
spaces}. \medskip

Most of the rings that we shall consider in this course are not just rings,
but $\mathbf{k}$-algebras. If $R$ is a $\mathbf{k}$-algebra, then any
$R$-module becomes a $\mathbf{k}$-module in a natural way:

\begin{definition}
\label{def.mod.R-to-k}Let $R$ be a $\mathbf{k}$-algebra. Then, any left
$R$-module $M$ automatically becomes a $\mathbf{k}$-module by setting%
\[
\lambda v:=\underbrace{\left(  \lambda\cdot1_{R}\right)  }_{\in R}%
v\ \ \ \ \ \ \ \ \ \ \text{for all }\lambda\in\mathbf{k}\text{ and }v\in M.
\]
This way of transforming $R$-modules into $\mathbf{k}$-modules is called
\emph{restriction of scalars}. The original action of $R$ on $M$ is easily
seen to be $\mathbf{k}$-bilinear as a map from $R\times M$ to $M$.
\end{definition}

Of course, the same holds for right $R$-modules. In the following, we shall
rarely talk about right $R$-modules, since all of their theory is a carbon
copy of the theory of left $R$-modules. Nevertheless, right $R$-modules
sometimes appear in nature, and you should be ready to formulate and prove
their properties as the need for them arises.

\subsubsection{Examples}

Let us see some examples of $R$-modules for various rings $R$. We do not
linger on the case of commutative $R$, since we have already been discussing
$\mathbf{k}$-modules for a long time. We recall that $\mathbf{k}$-modules are
the same as left $\mathbf{k}$-modules (and also the same as right $\mathbf{k}%
$-modules, if we write the scalars on the right). Let us now see some examples
of $R$-modules where $R$ is not necessarily commutative:

\begin{example}
\label{exa.mod.lreg}Any ring $R$ is automatically a left $R$-module, if we
define the action of $R$ on $R$ to be just the multiplication map of $R$ (that
is, $rm$ will just be the product of $r$ and $m$ in $R$). This is called the
\emph{left regular }$R$\emph{-module}, and its action is called the \emph{left
regular }$R$\emph{-action}.
\end{example}

\begin{example}
\label{exa.mod.knx1}The space $\mathbf{k}^{n\times1}$ of column vectors is not
only a $\mathbf{k}$-module, but also a left $\mathbf{k}^{n\times n}$-module.
Here, the action of $\mathbf{k}^{n\times n}$ on $\mathbf{k}^{n\times1}$ is
matrix-vector multiplication.

(Thus, Definition \ref{def.mod.leftmod} tells us here that we can refer to the
elements of $\mathbf{k}^{n\times n}$ as \textquotedblleft
scalars\textquotedblright\ and to the elements of $\mathbf{k}^{n\times1}$ as
\textquotedblleft vectors\textquotedblright. This is rather baroque, since the
\textquotedblleft scalars\textquotedblright\ have $n^{2}$ degrees of freedom
while the \textquotedblleft vectors\textquotedblright\ have only $n$. In other
words, our \textquotedblleft scalars\textquotedblright\ are \textquotedblleft
large\textquotedblright\ while our \textquotedblleft vectors\textquotedblright%
\ are \textquotedblleft small\textquotedblright, quite contrary to the
intuition we might have from linear algebra. But such is the price of working
with arbitrary rings instead of fields. For another example, $\mathbb{Z}/2$ is
a $\mathbb{Z}$-module, so that the elements of $\mathbb{Z}$ can be
\textquotedblleft scalars\textquotedblright\ while those of $\mathbb{Z}/2$ are
\textquotedblleft vectors\textquotedblright.)
\end{example}

\begin{example}
\label{exa.mod.knxm}More generally, for any $n,m\in\mathbb{N}$, the matrix
space $\mathbf{k}^{n\times m}$ is a left $\mathbf{k}^{n\times n}$-module. The
action is again given by matrix multiplication:%
\[
\left(  n\times n\text{-matrix}\right)  \cdot\left(  n\times m\text{-matrix}%
\right)  =\left(  n\times m\text{-matrix}\right)  .
\]

\end{example}

The next example is particularly important for us, as it defines one of the
simplest types of $\mathbf{k}\left[  G\right]  $-modules for a group algebra
$\mathbf{k}\left[  G\right]  $:

\begin{example}
\label{exa.mod.kG-on-kX}Let $G$ be a group, and let $X$ be a left $G$-set.
Then, the free $\mathbf{k}$-module $\mathbf{k}^{\left(  X\right)  }$ becomes a
left $\mathbf{k}\left[  G\right]  $-module by defining the action of
$\mathbf{k}\left[  G\right]  $ on $\mathbf{k}^{\left(  X\right)  }$ as
follows:%
\begin{align}
&  \left(  \sum_{g\in G}\alpha_{g}e_{g}\right)  \left(  \sum_{x\in X}\beta
_{x}e_{x}\right)  =\sum_{g\in G}\ \ \sum_{x\in X}\alpha_{g}\beta
_{x}e_{g\rightharpoonup x}\label{eq.exa.mod.kG-on-kX.1}\\
&  \ \ \ \ \ \ \ \ \ \ \text{for any }\left(  \alpha_{g}\right)  _{g\in G}%
\in\mathbf{k}\left[  G\right]  \text{ and }\left(  \beta_{x}\right)  _{x\in
X}\in\mathbf{k}^{\left(  X\right)  }.\nonumber
\end{align}
Using our convention to write $e_{g}$ as $g$ for each $g\in G$ (Convention
\ref{conv.monalg.em=m}), and a similar convention to write $e_{x}$ as $x$ for
each $x\in X$, we can rewrite this equality as%
\begin{equation}
\left(  \sum_{g\in G}\alpha_{g}g\right)  \left(  \sum_{x\in X}\beta
_{x}x\right)  =\sum_{g\in G}\ \ \sum_{x\in X}\alpha_{g}\beta_{x}\left(
g\rightharpoonup x\right)  . \label{eq.exa.mod.kG-on-kX.2}%
\end{equation}
This left $\mathbf{k}\left[  G\right]  $-module $\mathbf{k}^{\left(  X\right)
}$ is called the \emph{permutation module} corresponding to the left $G$-set
$X$. (It is named so because the elements of $G$ \textquotedblleft
permute\textquotedblright\ the basis vectors of $\mathbf{k}^{\left(  X\right)
}$.)
\end{example}

\begin{fineprint}
\begin{proof}
[Proof of Example \ref{exa.mod.kG-on-kX} (sketched).]We must show that the
module axioms hold. All of them are straightforward, and most are obvious. The
only one that takes a bit of effort is the associativity law. So let us prove
it. We must show that $\left(  rs\right)  m=r\left(  sm\right)  $ for all
$r,s\in\mathbf{k}\left[  G\right]  $ and $m\in\mathbf{k}^{\left(  X\right)  }$.

Fix $r,s\in\mathbf{k}\left[  G\right]  $ and $m\in\mathbf{k}^{\left(
X\right)  }$. We must prove the equality $\left(  rs\right)  m=r\left(
sm\right)  $. However, the action of $\mathbf{k}\left[  G\right]  $ on
$\mathbf{k}^{\left(  X\right)  }$ (more precisely, the map from $\mathbf{k}%
\left[  G\right]  \times\mathbf{k}^{\left(  X\right)  }$ to $\mathbf{k}%
^{\left(  X\right)  }$ that is defined by (\ref{eq.exa.mod.kG-on-kX.1}) and
that we are claiming to be an action of $\mathbf{k}\left[  G\right]  $ on
$\mathbf{k}^{\left(  X\right)  }$, even though we have not proved this yet) is
easily seen to be $\mathbf{k}$-bilinear, because each addend $\alpha_{g}%
\beta_{x}e_{g\rightharpoonup x}$ on the right hand side of
(\ref{eq.exa.mod.kG-on-kX.1}) depends $\mathbf{k}$-linearly on $\left(
\alpha_{g}\right)  _{g\in G}$ and depends $\mathbf{k}$-linearly on $\left(
\beta_{x}\right)  _{x\in X}\in\mathbf{k}^{\left(  X\right)  }$. Thus, both
sides of the equality $\left(  rs\right)  m=r\left(  sm\right)  $ depend
$\mathbf{k}$-linearly on $r$. Hence, in proving this equality, we can WLOG
assume that $r$ is an element of the standard basis $\left(  e_{g}\right)
_{g\in G}$ of $\mathbf{k}\left[  G\right]  $. For similar reasons, we can WLOG
assume that $s$ is an element of this standard basis as well, and that $m$ is
an element of the standard basis $\left(  e_{x}\right)  _{x\in X}$ of
$\mathbf{k}^{\left(  X\right)  }$. Let us make all these three assumptions.
Hence,%
\begin{equation}
r=e_{g},\ \ \ \ \ \ \ \ \ \ s=e_{h},\ \ \ \ \ \ \ \ \ \ \text{and}%
\ \ \ \ \ \ \ \ \ \ m=e_{x} \label{pf.exa.mod.kG-on-kX.rsm}%
\end{equation}
for some $g\in G$, $h\in G$ and $x\in X$. Consider these $g$, $h$ and $x$.
Now, (\ref{eq.exa.mod.kG-on-kX.1}) easily yields
\begin{equation}
e_{k}e_{y}=e_{k\rightharpoonup y}\ \ \ \ \ \ \ \ \ \ \text{for all }k\in
G\text{ and }y\in X \label{pf.exa.mod.kG-on-kX.ekey}%
\end{equation}
(indeed, this follows by applying (\ref{eq.exa.mod.kG-on-kX.1}) to $\alpha
_{g}=\delta_{g,k}$ and $\beta_{x}=\delta_{x,y}$, where $\delta$ is the
Kronecker delta). Now, from (\ref{pf.exa.mod.kG-on-kX.rsm}), we obtain%
\[
r\left(  sm\right)  =e_{g}\underbrace{\left(  e_{h}e_{x}\right)
}_{\substack{=e_{h\rightharpoonup x}\\\text{(by
(\ref{pf.exa.mod.kG-on-kX.ekey}))}}}=e_{g}e_{h\rightharpoonup x}%
=e_{g\rightharpoonup\left(  h\rightharpoonup x\right)  }%
\ \ \ \ \ \ \ \ \ \ \left(  \text{by (\ref{pf.exa.mod.kG-on-kX.ekey})}\right)
\]
and%
\begin{align*}
\left(  rs\right)  m  &  =\underbrace{\left(  e_{g}e_{h}\right)
}_{\substack{=e_{gh}\\\text{(by the definition of}\\\text{multiplication in
}\mathbf{k}\left[  G\right]  \text{)}}}e_{x}=e_{gh}e_{x}=e_{\left(  gh\right)
\rightharpoonup x}\ \ \ \ \ \ \ \ \ \ \left(  \text{by
(\ref{pf.exa.mod.kG-on-kX.ekey})}\right) \\
&  =e_{g\rightharpoonup\left(  h\rightharpoonup x\right)  }%
\ \ \ \ \ \ \ \ \ \ \left(
\begin{array}
[c]{c}%
\text{since }\left(  gh\right)  \rightharpoonup x=g\rightharpoonup\left(
h\rightharpoonup x\right) \\
\text{(because }X\text{ is a left }G\text{-set)}%
\end{array}
\right)  .
\end{align*}
Comparing these two equalities, we obtain $r\left(  sm\right)  =\left(
rs\right)  m$, which is precisely what we wanted to show. Thus, the
associativity law has been verified, and the proof of Example
\ref{exa.mod.kG-on-kX} is complete.
\end{proof}
\end{fineprint}

\begin{example}
In Example \ref{exa.mod.kG-on-kX}, we have defined permutation modules in
general; let us now show a specific hands-on instance of this construction.

Let $A$ be any set. Let $k\in\mathbb{N}$. As we saw in Example
\ref{exa.rep.G-sets.PkA}, the symmetric group $S_{A}$ acts on the set
$\mathcal{P}_{k}\left(  A\right)  $ of all $k$-element subsets of $A$. Thus,
$\mathcal{P}_{k}\left(  A\right)  $ is a left $S_{A}$-set. For instance,
taking $A=\left[  4\right]  $ and $k=2$, we conclude that
\[
\mathcal{P}_{2}\left(  \left[  4\right]  \right)  =\left\{  \left\{
1,2\right\}  ,\ \left\{  1,3\right\}  ,\ \left\{  1,4\right\}  ,\ \left\{
2,3\right\}  ,\ \left\{  2,4\right\}  ,\ \left\{  3,4\right\}  \right\}
\]
is a left $S_{\left[  4\right]  }$-set, i.e., a left $S_{4}$-set, with action
given by $g\rightharpoonup U=g\left(  U\right)  $ for all $g\in S_{4}$ and
$U\in\mathcal{P}_{2}\left(  \left[  4\right]  \right)  $ (that is,
$g\rightharpoonup\left\{  u,v\right\}  =\left\{  g\left(  u\right)  ,g\left(
v\right)  \right\}  $ for all $g\in S_{4}$ and $u\neq v$ in $\left[  4\right]
$).

The permutation module corresponding to this left $S_{4}$-set is a left
$\mathbf{k}\left[  S_{4}\right]  $-module $\mathbf{k}^{\left(  \mathcal{P}%
_{2}\left(  \left[  4\right]  \right)  \right)  }$. As a $\mathbf{k}$-module,
it is free with basis $\left(  e_{\left\{  1,2\right\}  },\ e_{\left\{
1,3\right\}  },\ e_{\left\{  1,4\right\}  },\ e_{\left\{  2,3\right\}
},\ e_{\left\{  2,4\right\}  },\ e_{\left\{  3,4\right\}  }\right)  $. The
action of $\mathbf{k}\left[  S_{4}\right]  $ is defined by
(\ref{eq.exa.mod.kG-on-kX.1}). Thus, for instance, for any two permutations
$g,h\in S_{4}$ and any two $2$-element subsets $\left\{  u,v\right\}  $ and
$\left\{  x,y\right\}  $ of $\left[  4\right]  $, we have%
\begin{align*}
&  \left(  e_{g}+e_{h}\right)  \left(  e_{\left\{  u,v\right\}  }+e_{\left\{
x,y\right\}  }\right) \\
&  =e_{g\rightharpoonup\left\{  u,v\right\}  }+e_{g\rightharpoonup\left\{
x,y\right\}  }+e_{h\rightharpoonup\left\{  u,v\right\}  }+e_{h\rightharpoonup
\left\{  x,y\right\}  }\\
&  =e_{\left\{  g\left(  u\right)  ,g\left(  v\right)  \right\}  }+e_{\left\{
g\left(  x\right)  ,g\left(  y\right)  \right\}  }+e_{\left\{  h\left(
u\right)  ,h\left(  v\right)  \right\}  }+e_{\left\{  h\left(  x\right)
,h\left(  y\right)  \right\}  }.
\end{align*}
Using our convention to write $e_{g}$ as $g$ for each $g\in G$ (Convention
\ref{conv.monalg.em=m}), and a similar convention to write $e_{x}$ as $x$ for
each $x\in X$, we can express this in the more pleasant form%
\begin{align*}
&  \left(  g+h\right)  \left(  \left\{  u,v\right\}  +\left\{  x,y\right\}
\right) \\
&  =\left\{  g\left(  u\right)  ,g\left(  v\right)  \right\}  +\left\{
g\left(  x\right)  ,g\left(  y\right)  \right\}  +\left\{  h\left(  u\right)
,h\left(  v\right)  \right\}  +\left\{  h\left(  x\right)  ,h\left(  y\right)
\right\}  ;
\end{align*}
of course, it needs to be kept in mind that we are not adding the literal
permutations $g,h$ or the literal subsets $\left\{  u,v\right\}  ,\left\{
x,y\right\}  $ but rather the corresponding basis vectors of $\mathbf{k}%
\left[  S_{4}\right]  $ or $\mathbf{k}^{\left(  \mathcal{P}_{2}\left(  \left[
4\right]  \right)  \right)  }$.
\end{example}

Here is a simpler example of a left $\mathbf{k}\left[  G\right]  $-module:

\begin{example}
\label{exa.mod.kG-eps}Let $G$ be a group. Any $\mathbf{k}$-module $V$ becomes
a left $\mathbf{k}\left[  G\right]  $-module if we define the action of
$\mathbf{k}\left[  G\right]  $ on $V$ by%
\[
\mathbf{a}v:=\varepsilon\left(  \mathbf{a}\right)
v\ \ \ \ \ \ \ \ \ \ \text{for all }\mathbf{a}\in\mathbf{k}\left[  G\right]
\text{ and }v\in V
\]
(where $\varepsilon:\mathbf{k}\left[  G\right]  \rightarrow\mathbf{k}$ is the
$\mathbf{k}$-algebra morphism that sends each $g\in G$ to $1$, as suggested in
Remark \ref{rmk.eps.all-G}).
\end{example}

\subsubsection{Properties and features}

Many of the basic concepts related to left $R$-modules are not much different
from their analogues for $\mathbf{k}$-modules. In particular, the following
concepts can be defined just as in the case of $\mathbf{k}$-modules (with the
obvious insertion of the word \textquotedblleft left\textquotedblright):

\begin{itemize}
\item $R$-submodules of an $R$-module (\cite[Definition 3.1.4]{23wa});

\item quotients of $R$-modules by $R$-submodules (\cite[Definition
3.6.1]{23wa});

\item direct sums and direct products of $R$-modules (\cite[\S 3.3.1]{23wa});

\item morphisms of $R$-modules, also known as \emph{left }$R$\emph{-linear
maps} (\cite[Definition 3.5.1]{23wa});

\item isomorphisms of $R$-modules (\cite[Definition 3.5.1]{23wa});

\item kernels and images of morphisms of $R$-modules (\cite[\S 3.5.6]{23wa}).
\end{itemize}

We refer to \cite[Chapter 3]{23wa} for more about these concepts. We shall
also use some basic properties of these constructions, such as:

\begin{itemize}
\item The kernel and the image of an $R$-module morphism $f:V\rightarrow W$
are always $R$-submodules of $V$ and $W$, respectively (\cite[\S 3.5.6]%
{23wa}). We denote them by $\operatorname*{Ker}f$ and $\operatorname{Im}f$,
respectively (but usually we just write $f\left(  V\right)  $ for
$\operatorname{Im}f$).

\item The \emph{first isomorphism theorem}\footnote{also known as the
\emph{first homomorphism theorem}}: If $f:V\rightarrow W$ is an $R$-module
morphism, then $V/\operatorname*{Ker}f\cong f\left(  V\right)  $ as
$R$-modules (\cite[Theorem 3.6.8]{23wa}).
\end{itemize}

\subsubsection{Representations of $G$ as left $\mathbf{k}\left[  G\right]
$-modules}

Now we take a closer look at left $\mathbf{k}\left[  G\right]  $-modules when
$G$ is a group. It turns out that they are just representations of $G$ over
$\mathbf{k}$ in disguise:

\begin{theorem}
\label{thm.rep.G-rep.mod}Let $G$ be a group. \medskip

\textbf{(a)} A left $\mathbf{k}\left[  G\right]  $-module $V$ is always a
$\mathbf{k}$-module (by restriction of scalars, as defined in Definition
\ref{def.mod.R-to-k}) and (at the same time) a left $G$-set (with the left
$G$-action defined by $g\rightharpoonup v=e_{g}v$ for all $g\in G$ and $v\in
V$). Moreover, the left $G$-action on this left $G$-set is $\mathbf{k}%
$-linear, so that $V$ becomes a representation of $G$ over $\mathbf{k}$.
\medskip

\textbf{(b)} Conversely: If $V$ is a representation of $G$ over $\mathbf{k}$,
then we can turn $V$ into a left $\mathbf{k}\left[  G\right]  $-module by
defining the action of $\mathbf{k}\left[  G\right]  $ on $V$ by
\begin{align}
&  \left(  \sum_{g\in G}\alpha_{g}g\right)  v:=\sum_{g\in G}\alpha_{g}\left(
g\rightharpoonup v\right) \label{eq.thm.rep.G-rep.mod.b.act}\\
&  \ \ \ \ \ \ \ \ \ \ \ \ \ \ \ \ \ \ \ \ \text{for all }\left(  \alpha
_{g}\right)  _{g\in G}\in\mathbf{k}\left[  G\right]  \text{ and }v\in
V.\nonumber
\end{align}

\textbf{(c)} These two transformations (from left $\mathbf{k}\left[  G\right]
$-modules to representations of $G$ and backwards) are mutually inverse. That is:

\begin{itemize}
\item If we convert a left $\mathbf{k}\left[  G\right]  $-module $V$ into a
representation of $G$ (according to Theorem \ref{thm.rep.G-rep.mod}
\textbf{(a)}), and then convert the resulting representation back into a left
$\mathbf{k}\left[  G\right]  $-module (according to Theorem
\ref{thm.rep.G-rep.mod} \textbf{(b)}), then we recover the original left
$\mathbf{k}\left[  G\right]  $-module $V$.

\item If we convert a representation $V$ of $G$ into a left $\mathbf{k}\left[
G\right]  $-module (according to Theorem \ref{thm.rep.G-rep.mod}
\textbf{(b)}), and then convert the resulting left $\mathbf{k}\left[
G\right]  $-module back into a representation of $G$ (according to Theorem
\ref{thm.rep.G-rep.mod} \textbf{(a)}), then we recover the original
representation $V$.
\end{itemize}
\end{theorem}

\begin{proof}
[Proof sketch.]\textbf{(a)} This is straightforward. For instance, the
associativity axiom $g\rightharpoonup\left(  h\rightharpoonup x\right)
=\left(  gh\right)  \rightharpoonup x$ follows from $e_{g}\left(
e_{h}x\right)  =\underbrace{\left(  e_{g}e_{h}\right)  }_{\substack{=e_{gh}%
\\\text{(by the definition of}\\\text{multiplication in }\mathbf{k}\left[
G\right]  \text{)}}}x=e_{gh}x$. \medskip

\textbf{(b)} This is a fairly easy linearity argument, similar to the proof of
Example \ref{exa.mod.kG-on-kX}. Here are the details:

\begin{fineprint}
Let $V$ be a representation of $G$ over $\mathbf{k}$. Define an action of
$\mathbf{k}\left[  G\right]  $ on $V$ by (\ref{eq.thm.rep.G-rep.mod.b.act}).
We must prove that it satisfies the module axioms. As usual, the hardest of
them is the associativity law. So let us check it.

We must show that $\left(  rs\right)  m=r\left(  sm\right)  $ for all
$r,s\in\mathbf{k}\left[  G\right]  $ and $m\in V$. The action of
$\mathbf{k}\left[  G\right]  $ on $V$ (more precisely, the map from
$\mathbf{k}\left[  G\right]  \times V$ to $V$ that is defined by
(\ref{eq.thm.rep.G-rep.mod.b.act}) and that we are claiming to be an action of
$\mathbf{k}\left[  G\right]  $ on $V$, even though we have not proved this
yet) is easily seen to be $\mathbf{k}$-bilinear (indeed,
(\ref{eq.thm.rep.G-rep.mod.b.act}) shows that it is $\mathbf{k}$-linear in
$\left(  \alpha_{g}\right)  _{g\in G}$, whereas the $\mathbf{k}$-linearity of
the $G$-action on $V$ entails that it is $\mathbf{k}$-linear in $v$). Thus,
both sides of the equality $\left(  rs\right)  m=r\left(  sm\right)  $ (which
we want to prove) depend $\mathbf{k}$-linearly on $r$. Hence, in proving this
equality, we can WLOG assume that $r$ is an element of the standard basis
$\left(  e_{g}\right)  _{g\in G}$ of $\mathbf{k}\left[  G\right]  $. For
similar reasons, we can WLOG assume that $s$ is an element of this standard
basis as well. Let us make both assumptions. Hence, $r=e_{g}$ and $s=e_{h}$
for some $g\in G$ and $h\in G$. Consider these $g$ and $h$. Following
Convention \ref{conv.monalg.em=m}, we rewrite $r=e_{g}$ and $s=e_{h}$ as $r=g$
and $s=h$.

Now, (\ref{eq.thm.rep.G-rep.mod.b.act}) easily yields
\begin{equation}
kv=k\rightharpoonup v\ \ \ \ \ \ \ \ \ \ \text{for all }k\in G\text{ and }v\in
V \label{pf.thm.rep.G-rep.mod.b.kv=}%
\end{equation}
(indeed, this follows by applying (\ref{eq.thm.rep.G-rep.mod.b.act}) to
$\alpha_{g}=\delta_{g,k}$, where $\delta$ is the Kronecker delta). Now, from
$r=g$ and $s=h$, we obtain%
\begin{align*}
\left(  rs\right)  m  &  =\left(  gh\right)  m=\left(  gh\right)
\rightharpoonup m\ \ \ \ \ \ \ \ \ \ \left(  \text{by
(\ref{pf.thm.rep.G-rep.mod.b.kv=})}\right) \\
&  =g\rightharpoonup\left(  h\rightharpoonup m\right)
\ \ \ \ \ \ \ \ \ \ \left(  \text{since }\rightharpoonup\text{ is a left
}G\text{-action on }V\right)
\end{align*}
and%
\begin{align*}
r\left(  sm\right)   &  =g\left(  hm\right)  =g\rightharpoonup
\underbrace{\left(  hm\right)  }_{\substack{=h\rightharpoonup m\\\text{(by
(\ref{pf.thm.rep.G-rep.mod.b.kv=}))}}}\ \ \ \ \ \ \ \ \ \ \left(  \text{by
(\ref{pf.thm.rep.G-rep.mod.b.kv=})}\right) \\
&  =g\rightharpoonup\left(  h\rightharpoonup m\right)  .
\end{align*}
Comparing these two equalities, we obtain $\left(  rs\right)  m=r\left(
sm\right)  $, as desired. This completes the proof of Theorem
\ref{thm.rep.G-rep.mod} \textbf{(b)}. \medskip
\end{fineprint}

\textbf{(c)} This is again easy. The \textquotedblleft representation
$\rightarrow$ $\mathbf{k}\left[  G\right]  $-module $\rightarrow$
representation\textquotedblright\ part is obvious. The \textquotedblleft%
$\mathbf{k}\left[  G\right]  $-module $\rightarrow$ representation
$\rightarrow$ $\mathbf{k}\left[  G\right]  $-module\textquotedblright%
\ direction requires a little linearity argument: Let $V$ be a left
$\mathbf{k}\left[  G\right]  $-module. Convert it into a representation of $G$
(according to Theorem \ref{thm.rep.G-rep.mod} \textbf{(a)}), denoting the
$G$-action by $\rightharpoonup$ as usual. Then, convert the resulting
representation back into a left $\mathbf{k}\left[  G\right]  $-module
(according to Theorem \ref{thm.rep.G-rep.mod} \textbf{(b)}). The action of
$\mathbf{k}\left[  G\right]  $ on this new $\mathbf{k}\left[  G\right]
$-module is then given by (\ref{eq.thm.rep.G-rep.mod.b.act}). But the action
of $\mathbf{k}\left[  G\right]  $ on the \textbf{original} $\mathbf{k}\left[
G\right]  $-module $V$ is also given by (\ref{eq.thm.rep.G-rep.mod.b.act}),
since it satisfies%
\begin{align*}
\left(  \sum_{g\in G}\alpha_{g}g\right)  v  &  =\sum_{g\in G}\alpha
_{g}\underbrace{gv}_{\substack{=g\rightharpoonup v\\\text{(by the definition
of}\\\text{the }G\text{-action }\rightharpoonup\text{)}}%
}\ \ \ \ \ \ \ \ \ \ \left(  \text{by }\mathbf{k}\text{-bilinearity}\right) \\
&  =\sum_{g\in G}\alpha_{g}\left(  g\rightharpoonup v\right)
\ \ \ \ \ \ \ \ \ \ \text{for all }\left(  \alpha_{g}\right)  _{g\in G}%
\in\mathbf{k}\left[  G\right]  \text{ and }v\in V\text{.}%
\end{align*}
Hence, the two actions of $\mathbf{k}\left[  G\right]  $ are the same. But
this means that the new left $\mathbf{k}\left[  G\right]  $-module is
precisely the original $\mathbf{k}\left[  G\right]  $-module $V$. This proves
Theorem \ref{thm.rep.G-rep.mod} \textbf{(c)}.
\end{proof}

A good way to think of Theorem \ref{thm.rep.G-rep.mod} \textbf{(b)} is as
\textquotedblleft packing\textquotedblright\ a $\mathbf{k}$-module structure
and a ($\mathbf{k}$-linear) left $G$-action into a left $\mathbf{k}\left[
G\right]  $-module structure (which combines them both, since $\mathbf{k}%
\left[  G\right]  $ includes both a \textquotedblleft copy\textquotedblright%
\ of $\mathbf{k}$ and a \textquotedblleft copy\textquotedblright\ of $G$).
Likewise, Theorem \ref{thm.rep.G-rep.mod} \textbf{(a)} \textquotedblleft
unpacks\textquotedblright\ a left $\mathbf{k}\left[  G\right]  $-module
structure into a $\mathbf{k}$-module structure and a left $G$-action. Thus, in
short, we have proved the following:

\begin{corollary}
\label{cor.rep.G-rep.mod}Let $G$ be a group. A representation of $G$ over
$\mathbf{k}$ is \textquotedblleft the same as\textquotedblright\ a left
$\mathbf{k}\left[  G\right]  $-module (since Theorem \ref{thm.rep.G-rep.mod}
allows us to convert either one into the other).
\end{corollary}

\begin{convention}
\label{conv.rep.G-rep.rep=mod}We shall thus use the notions of
\textquotedblleft representation of $G$ over $\mathbf{k}$\textquotedblright%
\ and \textquotedblleft left $\mathbf{k}\left[  G\right]  $%
-module\textquotedblright\ interchangeably, trusting the reader to convert
from one to the other at the slightest need. As a consequence, the notations
\textquotedblleft$g\rightharpoonup v$\textquotedblright\ and \textquotedblleft%
$gv$\textquotedblright\ and \textquotedblleft$e_{g}v$\textquotedblright%
\ (where $g$ belongs to the group $G$, and where $v$ belongs to the left
$\mathbf{k}\left[  G\right]  $-module $V$) all mean the same thing and will be
used as synonyms.
\end{convention}

\begin{example}
\label{exa.rep.G-rep.triv=mod}Let $G$ be a group. Let $V$ be a trivial
representation of $G$ (that is, a left $\mathbf{k}$-module equipped with a
trivial $G$-action, as defined in Example \ref{exa.rep.G-sets.triv}). Thus,
$V$ is a left $\mathbf{k}\left[  G\right]  $-module. Let us see explicitly how
an element $\mathbf{a}=\sum\limits_{g\in G}\alpha_{g}g\in\mathbf{k}\left[
G\right]  $ (with $\left(  \alpha_{g}\right)  _{g\in G}\in\mathbf{k}^{\left(
G\right)  }$) acts on a vector $v\in V$: Namely, from $\mathbf{a}%
=\sum\limits_{g\in G}\alpha_{g}g$, we obtain%
\begin{align}
\mathbf{a}v  &  =\left(  \sum\limits_{g\in G}\alpha_{g}g\right)
v=\sum\limits_{g\in G}\alpha_{g}\underbrace{\left(  g\rightharpoonup v\right)
}_{\substack{=v\\\text{(since the }G\text{-action}\\\text{on }V\text{ is
trivial)}}}\ \ \ \ \ \ \ \ \ \ \left(  \text{by
(\ref{eq.thm.rep.G-rep.mod.b.act})}\right) \nonumber\\
&  =\sum\limits_{g\in G}\alpha_{g}v=\left(  \sum\limits_{g\in G}\alpha
_{g}\right)  v. \label{eq.exa.rep.G-rep.triv=mod.1}%
\end{align}
We can rewrite this in a nicer form: Indeed, consider the $\mathbf{k}$-algebra
morphism $\varepsilon:\mathbf{k}\left[  G\right]  \rightarrow\mathbf{k}$ that
sends each $g\in G$ to $1$ (as suggested in Remark \ref{rmk.eps.all-G}). Then,
for each $\mathbf{a}=\sum\limits_{g\in G}\alpha_{g}g\in\mathbf{k}\left[
G\right]  $ (with $\left(  \alpha_{g}\right)  _{g\in G}\in\mathbf{k}^{\left(
G\right)  }$), we have%
\[
\varepsilon\left(  \mathbf{a}\right)  =\varepsilon\left(  \sum\limits_{g\in
G}\alpha_{g}g\right)  =\sum\limits_{g\in G}\alpha_{g}\underbrace{\varepsilon
\left(  g\right)  }_{=1}=\sum\limits_{g\in G}\alpha_{g}.
\]
Hence, we can rewrite (\ref{eq.exa.rep.G-rep.triv=mod.1}) as%
\[
\mathbf{a}v=\varepsilon\left(  \mathbf{a}\right)  v.
\]
Thus, our action of $\mathbf{k}\left[  G\right]  $ on $V$ is the same action
that we defined in Example \ref{exa.mod.kG-eps}. Thus, the latter action is
just a trivial representation of $G$ in disguise.
\end{example}

Let us now define a more useful representation of any group $G$:

\begin{definition}
\label{def.rep.G-rep.lreg}Let $G$ be a group. The \emph{left regular
representation} of $G$ is the left regular $\mathbf{k}\left[  G\right]
$-module $\mathbf{k}\left[  G\right]  $ (where the action is just the
multiplication in $\mathbf{k}\left[  G\right]  $, as in Example
\ref{exa.mod.lreg}). Thus, as a representation of $G$, it has its left
$G$-action given by%
\[
g\rightharpoonup\mathbf{a}=g\mathbf{a}=e_{g}\mathbf{a}%
\ \ \ \ \ \ \ \ \ \ \text{for all }g\in G\text{ and }\mathbf{a}\in
\mathbf{k}\left[  G\right]  .
\]

\end{definition}

For example, the left regular representation of $S_{n}$ is $\mathbf{k}\left[
S_{n}\right]  $, with action by left multiplication:%
\[
g\rightharpoonup\mathbf{a}=g\mathbf{a}=e_{g}\mathbf{a}%
\ \ \ \ \ \ \ \ \ \ \text{for all }g\in S_{n}\text{ and }\mathbf{a}%
\in\mathbf{k}\left[  S_{n}\right]  .
\]

Convention \ref{conv.rep.G-rep.rep=mod} gives us a new way to think about
representations of a group $G$. In particular, all concepts that are defined
for left $\mathbf{k}\left[  G\right]  $-modules (e.g., submodules, quotients,
morphisms, direct sums, direct products) can thus be automatically transferred
to representations of $G$. Thus, we obtain concepts of subrepresentation,
quotient of a representation, morphism of representations, direct sum, direct
product, etc. of representations of $G$ (by just viewing these representations
as left $\mathbf{k}\left[  G\right]  $-modules).

However, some of these concepts have already been defined for representations
before. To be safe, we need to verify that the concepts transferred from left
$\mathbf{k}\left[  G\right]  $-modules agree with these formerly defined
concepts for representations of $G$. But this is easy. For instance, we can
easily convince ourselves that subrepresentations are the same thing as
$\mathbf{k}\left[  G\right]  $-submodules:

\begin{proposition}
\label{prop.rep.G-rep.sub=sub}Let $G$ be a group. Let $V$ be a left
$\mathbf{k}\left[  G\right]  $-module (i.e., a representation of $G$ over
$\mathbf{k}$). Then, the left $\mathbf{k}\left[  G\right]  $-submodules of $V$
are precisely the subrepresentations of $V$.
\end{proposition}

\begin{fineprint}
\begin{proof}
[Proof sketch.]It is clear that any left $\mathbf{k}\left[  G\right]
$-submodule of $V$ is a subrepresentation of $V$. It remains to prove the
converse, i.e., that any subrepresentation of $V$ is a left $\mathbf{k}\left[
G\right]  $-submodule of $V$.

Let $W$ be a subrepresentation of $V$. We must show that $W$ is a left
$\mathbf{k}\left[  G\right]  $-submodule of $V$.

We know that $W$ is a subrepresentation of $V$. In other words, $W$ is a
$\mathbf{k}$-submodule of $V$ that is simultaneously a $G$-subset of $V$. In
particular, $W$ is a $G$-subset of $V$; hence,%
\begin{equation}
g\rightharpoonup w\in W\ \ \ \ \ \ \ \ \ \ \text{for each }g\in G\text{ and
}w\in W. \label{pf.prop.rep.G-rep.sub=sub.1}%
\end{equation}

Let $\mathbf{a}\in\mathbf{k}\left[  G\right]  $ and $w\in W$. Write
$\mathbf{a}$ as $\mathbf{a}=\sum_{g\in G}\alpha_{g}g$, where $\left(
\alpha_{g}\right)  _{g\in G}\in\mathbf{k}\left[  G\right]  $. Then,%
\begin{align*}
\mathbf{a}w  &  =\left(  \sum_{g\in G}\alpha_{g}g\right)  w=\sum_{g\in
G}\alpha_{g}\underbrace{\left(  g\rightharpoonup w\right)  }_{\substack{\in
W\\\text{(by (\ref{pf.prop.rep.G-rep.sub=sub.1}))}}%
}\ \ \ \ \ \ \ \ \ \ \left(  \text{by (\ref{eq.thm.rep.G-rep.mod.b.act}%
)}\right) \\
&  =\left(  \text{a }\mathbf{k}\text{-linear combination of elements of
}W\right)  \in W
\end{align*}
(since $W$ is a $\mathbf{k}$-submodule of $V$).

Forget that we fixed $\mathbf{a}$ and $w$. We thus have shown that
$\mathbf{a}w\in W$ for each $\mathbf{a}\in\mathbf{k}\left[  G\right]  $ and
$w\in W$. In other words, $W$ is a left $\mathbf{k}\left[  G\right]
$-submodule of $V$. This completes the proof of Proposition
\ref{prop.rep.G-rep.sub=sub}.
\end{proof}
\end{fineprint}

Similarly, morphisms of representations of $G$ are the same thing as morphisms
of left $\mathbf{k}\left[  G\right]  $-modules:

\begin{proposition}
\label{prop.rep.G-rep.mor=mor}Let $G$ be a group. Let $V$ and $W$ be two left
$\mathbf{k}\left[  G\right]  $-modules (i.e., representations of $G$ over
$\mathbf{k}$). Let $f:V\rightarrow W$ be a map. Then, $f$ is a morphism of
left $\mathbf{k}\left[  G\right]  $-modules if and only if $f$ is a morphism
of representations.
\end{proposition}

\begin{fineprint}
\begin{proof}
[Proof sketch.]$\Longrightarrow:$ This is easy.

$\Longleftarrow:$ Assume that $f$ is a morphism of representations. We must
show that $f$ is a morphism of left $\mathbf{k}\left[  G\right]  $-modules.

We have assumed that $f$ is a morphism of representations. In other words, $f$
is both $\mathbf{k}$-linear (i.e., a morphism of $\mathbf{k}$-modules) and
$G$-equivariant (i.e., a morphism of $G$-sets).

Now, let $\mathbf{a}\in\mathbf{k}\left[  G\right]  $ and $v\in V$. We shall
show that $f\left(  \mathbf{a}v\right)  =\mathbf{a}f\left(  v\right)  $.

Write $\mathbf{a}$ as $\mathbf{a}=\sum_{g\in G}\alpha_{g}g$, where $\left(
\alpha_{g}\right)  _{g\in G}\in\mathbf{k}\left[  G\right]  $. Then,%
\[
\mathbf{a}v=\left(  \sum_{g\in G}\alpha_{g}g\right)  v=\sum_{g\in G}\alpha
_{g}\left(  g\rightharpoonup v\right)  \ \ \ \ \ \ \ \ \ \ \left(  \text{by
(\ref{eq.thm.rep.G-rep.mod.b.act})}\right)  .
\]
Applying the map $f$ to both sides of this equality, we obtain%
\begin{align*}
f\left(  \mathbf{a}v\right)   &  =f\left(  \sum_{g\in G}\alpha_{g}\left(
g\rightharpoonup v\right)  \right)  =\sum_{g\in G}\alpha_{g}%
\underbrace{f\left(  g\rightharpoonup v\right)  }_{\substack{=g\rightharpoonup
f\left(  v\right)  \\\text{(since }f\text{ is }G\text{-equivariant)}%
}}\ \ \ \ \ \ \ \ \ \ \left(  \text{since }f\text{ is }\mathbf{k}%
\text{-linear}\right) \\
&  =\sum_{g\in G}\alpha_{g}\left(  g\rightharpoonup f\left(  v\right)
\right)  .
\end{align*}
Comparing this with%
\begin{align*}
\mathbf{a}f\left(  v\right)   &  =\left(  \sum_{g\in G}\alpha_{g}g\right)
f\left(  v\right)  \ \ \ \ \ \ \ \ \ \ \left(  \text{since }\mathbf{a}%
=\sum_{g\in G}\alpha_{g}g\right) \\
&  =\sum_{g\in G}\alpha_{g}\left(  g\rightharpoonup f\left(  v\right)
\right)  \ \ \ \ \ \ \ \ \ \ \left(
\begin{array}
[c]{c}%
\text{by (\ref{eq.thm.rep.G-rep.mod.b.act}), applied to }W\text{ and }f\left(
v\right) \\
\text{instead of }V\text{ and }v
\end{array}
\right)  ,
\end{align*}
we obtain $f\left(  \mathbf{a}v\right)  =\mathbf{a}f\left(  v\right)  $.

Forget that we fixed $\mathbf{a}$ and $v$. We have shown that every
$\mathbf{a}\in\mathbf{k}\left[  G\right]  $ and $v\in V$ satisfy $f\left(
\mathbf{a}v\right)  =\mathbf{a}f\left(  v\right)  $. In other words, the map
$f$ is a morphism of left $\mathbf{k}\left[  G\right]  $-modules (since $f$ is
$\mathbf{k}$-linear and thus respects addition). Thus, the proof of
Proposition \ref{prop.rep.G-rep.mor=mor} is complete.
\end{proof}
\end{fineprint}

\begin{corollary}
\label{cor.rep.G-rep.iso=iso}Let $G$ be a group. Let $V$ and $W$ be two left
$\mathbf{k}\left[  G\right]  $-modules (i.e., representations of $G$ over
$\mathbf{k}$). Let $f:V\rightarrow W$ be a map. Then, $f$ is an isomorphism of
left $\mathbf{k}\left[  G\right]  $-modules if and only if $f$ is an
isomorphism of representations.
\end{corollary}

\begin{fineprint}
\begin{proof}
[Proof sketch.]This follows easily from Proposition
\ref{prop.rep.G-rep.mor=mor}.
\end{proof}
\end{fineprint}

The reader is invited to convince himself that quotients of a representation
of $G$ are the same thing as quotients of a left $\mathbf{k}\left[  G\right]
$-module.

\begin{exercise}
\label{exe.mod.kX-on-kX.trivsub}\fbox{1} Let $G$ be a group, and let $X$ be a
finite left $G$-set. Consider the permutation module $\mathbf{k}^{\left(
X\right)  }$ defined in Example \ref{exa.mod.kG-on-kX}. Let $s:=\sum_{x\in
X}e_{x}\in\mathbf{k}^{\left(  X\right)  }$. Prove that $\operatorname*{span}%
\nolimits_{\mathbf{k}}\left\{  s\right\}  $ is a trivial subrepresentation of
the $G$-representation $\mathbf{k}^{\left(  X\right)  }$. This
subrepresentation is $1$-dimensional (i.e., a free $\mathbf{k}$-module of rank
$1$) unless $X$ is empty.
\end{exercise}

\subsubsection{Irreducible representations}

Clearly, any representation $V$ of any group $G$ has two subrepresentations:
the subset $\left\{  0\right\}  $ and $V$ itself. (These are distinct unless
$V$ itself is $\left\{  0\right\}  $.) Sometimes, these are the only
subrepresentations of $V$. In this case, we call $V$ \emph{simple}:

\begin{definition}
\label{def.rep.G-rep.irrep}Let $G$ be a group. A representation $V$ of $G$
(or, equivalently, a left $\mathbf{k}\left[  G\right]  $-module $V$) is said
to be \emph{simple} (aka \emph{irreducible}) if it is not $\left\{  0\right\}
$, but its only subrepresentations are $V$ and $\left\{  0\right\}  $.
\end{definition}

More generally, we can make the same definition for modules over a ring:

\begin{definition}
\label{def.mod.irrep}Let $R$ be a ring. A left $R$-module $V$ is said to be
\emph{simple} (aka \emph{irreducible}) if it is not $\left\{  0\right\}  $,
but its only $R$-submodules are $V$ and $\left\{  0\right\}  $.
\end{definition}

Clearly, Definition \ref{def.rep.G-rep.irrep} is the particular case of
Definition \ref{def.mod.irrep} for $R=\mathbf{k}\left[  G\right]  $.
Traditionally, the word \textquotedblleft simple\textquotedblright\ is used
for left $\mathbf{k}\left[  G\right]  $-modules, whereas the word
\textquotedblleft irreducible\textquotedblright\ applies to representations.
But the two concepts are equivalent, and there is no reason other than habit
to distinguish between them.

The following abbreviations are particularly popular with physicists:
Irreducible representations are called \emph{irreps}, whereas representations
in general are called \emph{reps}.

The irreducibility of a representation depends heavily on the base ring
$\mathbf{k}$. The most well-behaved case is when $\mathbf{k}$ is a field:

\begin{example}
\label{exa.rep.G-rep.irrep.1}Let $\mathbf{k}$ be a field. The representations
of the trivial group $\left\{  1\right\}  $ over $\mathbf{k}$ are just the
$\mathbf{k}$-vector spaces (by Example \ref{exa.rep.G-rep.lin.1}). Their
subrepresentations are simply their $\mathbf{k}$-vector subspaces. Thus, an
irreducible representation of $\left\{  1\right\}  $ is simply a $\mathbf{k}%
$-vector space $V$ that is not $\left\{  0\right\}  $ but whose only vector
subspaces are $V$ and $\left\{  0\right\}  $. What are these spaces?

If a $\mathbf{k}$-vector space $V$ has dimension $1$, then its only vector
subspaces are $V$ and $\left\{  0\right\}  $, and thus $V$ is an irreducible
representation of $\left\{  1\right\}  $. But any $\mathbf{k}$-vector space
$V$ of dimension $>1$ has a subspace of dimension $1$, and thus cannot be
irreducible (since such a subspace is neither $V$ nor $\left\{  0\right\}  $).
Thus, the irreducible representations of $\left\{  1\right\}  $ are exactly
the $1$-dimensional $\mathbf{k}$-vector spaces. Up to isomorphism, there is
only one such ($\mathbf{k}$ itself).
\end{example}

\begin{example}
\label{exa.rep.G-rep.irrep.Sn}Let $\mathbf{k}$ be a field. Then, the natural
representation $\mathbf{k}^{n}$ of $S_{n}$ is not irreducible (unless $n=1$),
but (for $n\geq1$) its subrepresentation $D\left(  \mathbf{k}^{n}\right)  $
(defined in Subsection \ref{subsec.rep.G-rep.subreps}) is (since it is
$1$-dimensional as a $\mathbf{k}$-vector space, and thus has no subspaces
other than itself and $\left\{  0\right\}  $). The subrepresentation $R\left(
\mathbf{k}^{n}\right)  $ is also irreducible if $n\geq2$ and
$\operatorname*{char}\mathbf{k}\nmid n$ (this follows easily from Theorem
\ref{thm.rep.G-rep.Sn-nat.subreps}), but not if $n>2$ and
$\operatorname*{char}\mathbf{k}\mid n$ (since Proposition
\ref{prop.rep.G-rep.Sn-nat.subreps.dirsum} \textbf{(b)} yields $D\left(
\mathbf{k}^{n}\right)  \subseteq R\left(  \mathbf{k}^{n}\right)  $ in this case).

Another irreducible representation of $S_{n}$ is the sign representation
$\mathbf{k}_{\operatorname*{sign}}$ (as defined in Example
\ref{exa.rep.Sn-rep.sign}). This is again because this representation is
$1$-dimensional as a $\mathbf{k}$-vector space. (Note, however, that the sign
representation is isomorphic to $R\left(  \mathbf{k}^{n}\right)  $ if $n=2$.)

Are there more irreducible representations of $S_{n}$ ? Yes, and we will get
to meet them.
\end{example}

\begin{example}
\label{exa.rep.G-rep.irrep.C2}Let $\mathbf{k}$ be a field. The representations
of the cyclic group $C_{2}=\left\{  1,g\right\}  $ over $\mathbf{k}$ are
$\mathbf{k}$-vector spaces equipped with a $\mathbf{k}$-linear involution $w$
(as we saw in Example \ref{exa.rep.G-rep.lin.C2}). A subrepresentation of such
a representation is a vector subspace $W$ that is preserved by this involution
(i.e., that satisfies $w\left(  W\right)  \subseteq W$). Thus, the irreducible
representations of $C_{2}$ are the nonzero $\mathbf{k}$-vector spaces $V$
equipped with a $\mathbf{k}$-linear involution $w:V\rightarrow V$ that
preserves no vector subspaces of $V$ other than $V$ and $\left\{  0\right\}
$. What are these?

We claim the following: \medskip

\textbf{(a)} All irreducible representations of $C_{2}$ are $1$-dimensional.
Moreover, the $\mathbf{k}$-linear involution $w$ (that is, the action of the
generator $g\in C_{2}$) on any such representation is either the identity map
$\operatorname*{id}$ or its negative $-\operatorname*{id}$. \medskip

\textbf{(b)} If $\operatorname*{char}\mathbf{k}\neq2$, then $C_{2}$ has (up to
isomorphism) exactly two different irreducible representations: the trivial
representation (the vector space $\mathbf{k}$, with $g$ acting as
$\operatorname*{id}$) and the \textquotedblleft sign
representation\textquotedblright\ (the vector space $\mathbf{k}$, with $g$
acting as $-\operatorname*{id}$). \medskip

\textbf{(c)} If $\operatorname*{char}\mathbf{k}=2$, then $C_{2}$ has (up to
isomorphism) only one irreducible representation: the trivial representation
(the vector space $\mathbf{k}$, with $g$ acting as $\operatorname*{id}$).
\end{example}

\begin{fineprint}
\begin{proof}
[Proof of Example \ref{exa.rep.G-rep.irrep.C2} (sketched).]\textbf{(a)} Let
$V$ be an irreducible representation of $C_{2}$, and let $w$ be the
corresponding $\mathbf{k}$-linear involution of $V$ (that is, the action of
$g$). We must prove that $\dim V=1$ and that $w$ is either $\operatorname*{id}%
$ or $-\operatorname*{id}$.

Since $w$ is an involution, we have $w^{2}=\operatorname*{id}=1$, so that
$w^{2}-1=0$. But $w^{2}-1=\left(  w-1\right)  \left(  w+1\right)  $ (since $w$
is a $\mathbf{k}$-linear map). Thus, $\left(  w-1\right)  \left(  w+1\right)
=0$.

We are in one of the two cases:

\textit{Case 1:} We have $\operatorname*{Ker}\left(  w-1\right)  =\left\{
0\right\}  $.

\textit{Case 2:} We have $\operatorname*{Ker}\left(  w-1\right)  \neq\left\{
0\right\}  $.

Let us first consider Case 1. In this case, we have $\operatorname*{Ker}%
\left(  w-1\right)  =\left\{  0\right\}  $. In other words, the $\mathbf{k}%
$-linear map $w-1$ is injective. Hence, we can cancel $w-1$ from $\left(
w-1\right)  \left(  w+1\right)  =0$, obtaining $w+1=0$. Thus,
$w=-1=-\operatorname*{id}$. Hence, any vector subspace of $V$ is automatically
preserved by $w$ (since $-\operatorname*{id}$ preserves any subspace), and
thus is a subrepresentation of $V$. If $\dim V>1$, then $V$ has at least one
vector subspace that is neither $V$ nor $\left\{  0\right\}  $, and thus (by
the preceding sentence) at least one subrepresentation that is neither $V$ nor
$\left\{  0\right\}  $; but this contradicts the irreducibility of $V$. Hence,
we cannot have $\dim V>1$. Thus, $\dim V\leq1$. But the irreducibility of $V$
also yields $V\neq\left\{  0\right\}  $ and therefore $\dim V\geq1$.
Combining, we find $\dim V=1$. Hence, we have shown that $\dim V=1$ and
$w=-\operatorname*{id}$. Thus, Example \ref{exa.rep.G-rep.irrep.C2}
\textbf{(a)} is proved in Case 1.

Let us now consider Case 2. In this case, we have $\operatorname*{Ker}\left(
w-1\right)  \neq\left\{  0\right\}  $. Hence, there exists a nonzero vector
$v\in\operatorname*{Ker}\left(  w-1\right)  $. Consider this $v$. From
$v\in\operatorname*{Ker}\left(  w-1\right)  $, we obtain $\left(  w-1\right)
\left(  v\right)  =0$, so that $w\left(  v\right)  -v=0$ and thus $w\left(
v\right)  =v$. Therefore, the $1$-dimensional vector subspace
$\operatorname*{span}\nolimits_{\mathbf{k}}\left\{  v\right\}  $ of $V$ is
preserved by $w$. In other words, $\operatorname*{span}\nolimits_{\mathbf{k}%
}\left\{  v\right\}  $ is a subrepresentation of $V$. By the irreducibility of
$V$, this subrepresentation must be either $V$ or $\left\{  0\right\}  $.
Hence, $\operatorname*{span}\nolimits_{\mathbf{k}}\left\{  v\right\}  $ must
be $V$ (since it cannot be $\left\{  0\right\}  $). But this necessarily
entails that $\dim V=1$ and $w=\operatorname*{id}$ (since $w\left(  v\right)
=v$). Thus, Example \ref{exa.rep.G-rep.irrep.C2} \textbf{(a)} is proved in
Case 2.

We have now proved Example \ref{exa.rep.G-rep.irrep.C2} \textbf{(a)} in both
Cases 1 and 2. \medskip

\textbf{(b)} From $\operatorname*{char}\mathbf{k}\neq2$, we obtain $1\neq-1$
in $\mathbf{k}$. Thus, the trivial representation and the \textquotedblleft
sign representation\textquotedblright\ of $C_{2}$ are not isomorphic (since
$w$ acts as the identity on the former but not on the latter). Moreover,
Example \ref{exa.rep.G-rep.irrep.C2} \textbf{(a)} shows that any irreducible
representation of $C_{2}$ is isomorphic to one of the two. This proves Example
\ref{exa.rep.G-rep.irrep.C2} \textbf{(b)}. \medskip

\textbf{(c)} This is just like part \textbf{(b)}, but with one difference:
Since $\operatorname*{char}\mathbf{k}=2$, we now have $1=-1$ in $\mathbf{k}$,
so that the trivial representation and the \textquotedblleft sign
representation\textquotedblright\ are the same thing.
\end{proof}
\end{fineprint}

This last example shows that the description of the irreps of a group $G$ over
$\mathbf{k}$ can depend on $\operatorname*{char}\mathbf{k}$. The situation
where $\mathbf{k}$ is not a field is even weirder:

\begin{example}
\label{exa.rep.G-rep.irrep.1Z}The irreducible representations of the trivial
group $\left\{  1\right\}  $ over $\mathbb{Z}$ are the simple $\mathbb{Z}%
\left[  \left\{  1\right\}  \right]  $-modules, i.e., the simple $\mathbb{Z}%
$-modules (since $\mathbb{Z}\left[  \left\{  1\right\}  \right]
\cong\mathbb{Z}$). They are (up to isomorphism) the $\mathbb{Z}$-modules
$\mathbb{Z}/p$ for all primes $p$.
\end{example}

\begin{exercise}
Consider the cyclic group $C_{3}$. \medskip

\textbf{(a)} \fbox{2} Prove that $C_{3}$ has (up to isomorphism) three
different irreducible representations over $\mathbb{C}$, and each of them has
dimension $1$. \medskip

\textbf{(b)} \fbox{2} Prove that $C_{3}$ has (up to isomorphism) two different
irreducible representations over $\mathbb{R}$: one of dimension $1$ and one of
dimension $2$. \medskip

[\textbf{Hint:} What eigenvalues can a $\mathbb{C}$-linear map (i.e., complex
matrix) $w$ satisfying $w^{3}=\operatorname*{id}$ have?]
\end{exercise}

\subsubsection{Some further properties of representations/modules}

We recall a basic concept from linear algebra (\cite[Exercise 3.11.3]{23wa}):

\begin{definition}
\label{def.EndV}Let $V$ be a $\mathbf{k}$-module. Then, we set%
\[
\operatorname*{End}\nolimits_{\mathbf{k}}V:=\left\{  \mathbf{k}\text{-module
endomorphisms of }V\right\}  =\left\{  \mathbf{k}\text{-linear maps
}V\rightarrow V\right\}  .
\]
This is a $\mathbf{k}$-algebra (with addition and scaling being
pointwise\footnotemark, and with multiplication being composition of maps). It
is called the \emph{endomorphism ring} of $V$. Its group of units (i.e.,
invertible elements) is the general linear group $\operatorname*{GL}\left(
V\right)  $.
\end{definition}

\footnotetext{That is:
\par
\begin{itemize}
\item The sum $f+g$ of two endomorphisms $f,g\in\operatorname*{End}%
\nolimits_{\mathbf{k}}V$ is defined by%
\[
\left(  f+g\right)  \left(  v\right)  =f\left(  v\right)  +g\left(  v\right)
\ \ \ \ \ \ \ \ \ \ \text{for all }v\in V.
\]
\par
\item Scaling an endomorphism $f\in\operatorname*{End}\nolimits_{\mathbf{k}}V$
by a scalar $\lambda\in\mathbf{k}$ results in an endomorphism $\lambda
f\in\operatorname*{End}\nolimits_{\mathbf{k}}V$ defined by%
\[
\left(  \lambda f\right)  \left(  v\right)  =\lambda f\left(  v\right)
\ \ \ \ \ \ \ \ \ \ \text{for all }v\in V.
\]
\end{itemize}
}Note that when $V$ is the free $\mathbf{k}$-module $\mathbf{k}^{n}$ for some
$n\in\mathbb{N}$, the endomorphism ring $\operatorname*{End}%
\nolimits_{\mathbf{k}}V$ is isomorphic to the matrix ring $\mathbf{k}^{n\times
n}$ (as a $\mathbf{k}$-algebra). The isomorphism simply sends a $\mathbf{k}%
$-linear map $V\rightarrow V$ to the matrix that represents this map in the
standard basis of $V$.

We can now define an analogue of currying for modules over any $\mathbf{k}$-algebra:

\begin{theorem}
\label{thm.mod.curry-into-EndV}Let $R$ be a $\mathbf{k}$-algebra. Let $V$ be a
left $R$-module. Then, for each $r\in R$, the map%
\begin{align*}
\tau_{r}:V  &  \rightarrow V,\\
v  &  \mapsto rv
\end{align*}
is $\mathbf{k}$-linear and thus belongs to $\operatorname*{End}%
\nolimits_{\mathbf{k}}V$. Moreover, the map%
\begin{align*}
\rho:R  &  \rightarrow\operatorname*{End}\nolimits_{\mathbf{k}}V,\\
r  &  \mapsto\tau_{r}%
\end{align*}
is a $\mathbf{k}$-algebra morphism.
\end{theorem}

\begin{proof}
[Proof sketch.]Straightforward and left to the reader. (For instance,
$\tau_{rs}=\tau_{r}\tau_{s}$ follows from the associativity law $\left(
rs\right)  m=r\left(  sm\right)  $ for the left $R$-module $V$. Most other
properties follow from the $\mathbf{k}$-bilinearity of the action $R\times
V\rightarrow V,\ \left(  r,v\right)  \mapsto rv$.)
\end{proof}

\begin{definition}
\label{def.mod.curry-into-EndV}The morphism $\rho$ in Theorem
\ref{thm.mod.curry-into-EndV} is called the \emph{curried form} of the left
$R$-action on $V$ (as a $\mathbf{k}$-module).

The map $\tau_{r}$ in Theorem \ref{thm.mod.curry-into-EndV} is called the
\emph{action} of $r$ on $V$.
\end{definition}

This notion of \textquotedblleft action\textquotedblright\ is perfectly
aligned with the notion of \textquotedblleft action\textquotedblright\ for a
$G$-set (see Definition \ref{def.rep.G-sets.act}): Indeed, if $G$ is a group
and if $V$ is a left $\mathbf{k}\left[  G\right]  $-module (i.e., a
representation of $G$), then the action of a group element $g\in G$ on $V$ in
the sense of Definition \ref{def.rep.G-sets.act} is precisely the action of
the corresponding standard basis vector $g=e_{g}\in\mathbf{k}\left[  G\right]
$ on $V$ in the sense of Definition \ref{def.mod.curry-into-EndV} (since
$g\rightharpoonup v=e_{g}v$ for each $v\in V$). \medskip

Currying lets us define the notion of a faithful module:

\begin{definition}
\label{def.mod.faith}Let $R$ be a $\mathbf{k}$-algebra. A left $R$-module $V$
is said to be \emph{faithful} if it has the following property: If $r,s\in R$
are two elements that satisfy%
\[
\left(  rv=sv\ \ \ \ \ \ \ \ \ \ \text{for all }v\in V\right)  ,
\]
then $r=s$.

In other words, a left $R$-module $V$ is said to be \emph{faithful} if the
curried form $\rho:R\rightarrow\operatorname*{End}\nolimits_{\mathbf{k}}V$ of
its action is injective.
\end{definition}

Roughly speaking, a left $R$-module $V$ is faithful if and only if you can
tell elements of $R$ apart by their action on $V$ (that is, two different
elements of $R$ never act the same way on $V$).

\begin{warning}
\label{warn.rep.G-rep.faith}Let $G$ be a group. A faithful left $\mathbf{k}%
\left[  G\right]  $-module is not the same as a faithful representation of
$G$. Indeed, for a representation of $G$ to be faithful, we only need
different elements of $G$ to act differently. However, for a left
$\mathbf{k}\left[  G\right]  $-module to be faithful, we need different
\textbf{linear combinations} of elements of $G$ to act differently. This is a
much stronger requirement! Here is an example for this difference:
\end{warning}

\begin{example}
\label{exa.rep.G-rep.faith-Sn}Consider the natural representation
$\mathbf{k}^{n}$ of the symmetric group $S_{n}$.

This representation is faithful as a representation of $S_{n}$ (if
$\mathbf{k}$ is nontrivial), since different permutations $g\in S_{n}$ can be
told apart by their action on the standard basis vectors $e_{1},e_{2}%
,\ldots,e_{n}$ of $\mathbf{k}^{n}$.

However, the natural representation $\mathbf{k}^{n}$ is \textbf{not} faithful
as a $\mathbf{k}\left[  S_{n}\right]  $-module, unless $n\in\left\{
1,2\right\}  $. For example, let $n\geq3$, and consider the element%
\[
\nabla_{\left[  3\right]  }^{-}=1-s_{1}-s_{2}-t_{1,3}+\operatorname*{cyc}%
\nolimits_{1,2,3}+\operatorname*{cyc}\nolimits_{1,3,2}\in\mathbf{k}\left[
S_{n}\right]  .
\]
We claim that this element $\nabla_{\left[  3\right]  }^{-}$ (which is nonzero
unless $\mathbf{k}$ is trivial) acts as $0$ on $\mathbf{k}^{n}$, meaning that
$\nabla_{\left[  3\right]  }^{-}v=0$ for all $v\in\mathbf{k}^{n}$. In order to
prove this, it suffices to show that $\nabla_{\left[  3\right]  }^{-}e_{k}=0$
for each $k\in\left[  n\right]  $ (since the standard basis vectors
$e_{1},e_{2},\ldots,e_{n}$ span $\mathbf{k}^{n}$). We can check this directly:
We have%
\begin{align*}
\nabla_{\left[  3\right]  }^{-}e_{1}  &  =\underbrace{1e_{1}}_{=e_{1}%
}-\underbrace{s_{1}e_{1}}_{=e_{2}}-\underbrace{s_{2}e_{1}}_{=e_{1}%
}-\underbrace{t_{1,3}e_{1}}_{=e_{3}}+\underbrace{\operatorname*{cyc}%
\nolimits_{1,2,3}e_{1}}_{=e_{2}}+\underbrace{\operatorname*{cyc}%
\nolimits_{1,3,2}e_{1}}_{=e_{3}}\\
&  =e_{1}-e_{2}-e_{1}-e_{3}+e_{2}+e_{3}=0
\end{align*}
and similarly $\nabla_{\left[  3\right]  }^{-}e_{2}=0$ and $\nabla_{\left[
3\right]  }^{-}e_{3}=0$, and it is even easier to see that $\nabla_{\left[
3\right]  }^{-}e_{k}=0$ for all $k\geq4$.

More generally, if $X$ is any subset of $\left[  n\right]  $ satisfying
$\left\vert X\right\vert >2$, then%
\[
\nabla_{X}^{-}v=0\ \ \ \ \ \ \ \ \ \ \text{for any }v\in\mathbf{k}^{n}.
\]

\end{example}

The simplest example of a faithful $R$-module (for any ring $R$) is $R$ itself:

\begin{example}
Let $R$ be any ring. Then, the left regular $R$-module $R$ is faithful.
Indeed, if $r,s\in R$ are two elements that satisfy%
\[
\left(  rv=sv\ \ \ \ \ \ \ \ \ \ \text{for all }v\in R\right)  ,
\]
then $r=s$ (because applying $rv=sv$ to $v=1$ yields $r1=s1$, that is, $r=s$).
\end{example}

The following exercise shows that a faithful left $\mathbf{k}\left[  G\right]
$-module (for a group $G$) has to be \textquotedblleft not too
small\textquotedblright\ (whereas a faithful representation of a group $G$ can
often be $1$-dimensional\footnote{More precisely: Any cyclic group has a
$1$-dimensional faithful representation over $\mathbb{C}$ (why?). Non-cyclic
groups often don't.}):

\begin{exercise}
\fbox{1} Let $\mathbf{k}$ be a field. Let $G$ be a group. Let $V$ be a
faithful left $\mathbf{k}\left[  G\right]  $-module. Prove that $\dim
V\geq\sqrt{\left\vert G\right\vert }$, where $\dim V$ denotes the dimension of
$V$ as a $\mathbf{k}$-vector space. (If $G$ is infinite, then this is saying
that $V$ is infinite-dimensional.)
\end{exercise}

\subsection{\label{sec.rep.maschke}Direct addends, projections, filtrations,
Maschke, Jordan--H\"{o}lder}

Direct addends of an $R$-module are a particularly well-behaved kind of
submodules. We will first briefly discuss them for arbitrary rings $R$, then
focus on the case when $R$ is the group algebra $\mathbf{k}\left[  G\right]  $
of a group $G$.

\subsubsection{\label{subsec.rep.maschke.diradd}Direct addends and
projections}

First, let us do some elementary linear algebra, albeit with $R$-modules
instead of vector spaces.

Recall the following: If $R$ is a ring, and if $W_{1},W_{2},\ldots,W_{k}$ are
finitely many left $R$-modules, then the direct sum $W_{1}\oplus W_{2}%
\oplus\cdots\oplus W_{k}$ of these left $R$-modules is the left $R$-module%
\begin{align*}
W_{1}\oplus W_{2}\oplus\cdots\oplus W_{k}  &  =W_{1}\times W_{2}\times
\cdots\times W_{k}\\
&  =\left\{  \left(  w_{1},w_{2},\ldots,w_{k}\right)  \ \mid\ w_{i}\in
W_{i}\text{ for each }i\right\}
\end{align*}
with entrywise addition and entrywise action of $R$. (Direct sums of
infinitely many addends are a bit more complicated, since -- unlike direct
products -- they only allow essentially finite families of vectors. See
\cite[Definition 3.3.5]{23wa} for the proper definition of a direct sum in
this case.)

This kind of direct sum is called an \emph{external direct sum}, since it
constructs a new $R$-module ($W_{1}\oplus W_{2}\oplus\cdots\oplus W_{k}$) from
existing ones ($W_{1},W_{2},\ldots,W_{k}$). Sometimes, this feels unnecessary,
since $W_{1},W_{2},\ldots,W_{k}$ are already contained in a single $R$-module
$V$ and are \textquotedblleft sufficiently spread out\textquotedblright\ that
their direct sum can be identified with an $R$-submodule of $V$ as well. In
this situation, the latter $R$-submodule is called the \emph{internal direct
sum} of $W_{1},W_{2},\ldots,W_{k}$. Here is the actual definition:

\begin{definition}
\label{def.mod.sum}Let $R$ be a ring. Let $V$ be a left $R$-module. Let
$W_{1},W_{2},\ldots,W_{k}$ be $R$-submodules of $V$. \medskip

\textbf{(a)} The \emph{sum} $W_{1}+W_{2}+\cdots+W_{k}$ of these $R$-submodules
is defined to be the set%
\[
\left\{  w_{1}+w_{2}+\cdots+w_{k}\ \mid\ w_{i}\in W_{i}\text{ for each }%
i\in\left[  k\right]  \right\}  .
\]
This is an $R$-submodule of $V$ that contains $W_{1},W_{2},\ldots,W_{k}$ as
$R$-submodules. It is furthermore the smallest $R$-submodule of $V$ with this
property. \medskip

\textbf{(b)} There is a canonical $R$-module morphism
\begin{align*}
\operatorname*{sum}:W_{1}\oplus W_{2}\oplus\cdots\oplus W_{k}  &  \rightarrow
W_{1}+W_{2}+\cdots+W_{k},\\
\left(  w_{1},w_{2},\ldots,w_{k}\right)   &  \mapsto w_{1}+w_{2}+\cdots+w_{k}.
\end{align*}
This morphism is always surjective. \medskip

\textbf{(c)} If this $R$-module morphism $\operatorname*{sum}$ is also
injective (and thus an isomorphism), then we call the sum $W_{1}+W_{2}%
+\cdots+W_{k}$ an \emph{internal direct sum}, and we denote it by $W_{1}\oplus
W_{2}\oplus\cdots\oplus W_{k}$ as well (although, strictly speaking, it is
merely isomorphic to $W_{1}\oplus W_{2}\oplus\cdots\oplus W_{k}$), thus
implicitly identifying each element $\left(  w_{1},w_{2},\ldots,w_{k}\right)
$ of the external direct sum $W_{1}\oplus W_{2}\oplus\cdots\oplus W_{k}$ with
the element $\operatorname*{sum}\left(  w_{1},w_{2},\ldots,w_{k}\right)
=w_{1}+w_{2}+\cdots+w_{k}$ of the internal direct sum $W_{1}+W_{2}%
+\cdots+W_{k}$. In this situation, we also say that \emph{the sum }%
$W_{1}+W_{2}+\cdots+W_{k}$ \emph{is direct}, or that the $R$-submodules
$W_{1},W_{2},\ldots,W_{k}$ are \emph{linearly disjoint}.
\end{definition}

This definition is frequently used in abstract linear algebra, e.g., when
saying that the eigenspaces of a linear map have a direct sum (which almost
always means the internal direct sum). For example, the two $\mathbb{Q}%
$-vector subspaces $\operatorname*{span}\nolimits_{\mathbb{Q}}\left\{
e_{1},\ e_{2}\right\}  $ and $\operatorname*{span}\nolimits_{\mathbb{Q}%
}\left\{  e_{3}+e_{4}\right\}  $ of $\mathbb{Q}^{4}$ are linearly disjoint and
have the internal direct sum $\operatorname*{span}\nolimits_{\mathbb{Q}%
}\left\{  e_{1},\ e_{2},\ e_{3}+e_{4}\right\}  $, whereas the two $\mathbb{Q}%
$-vector subspaces $\operatorname*{span}\nolimits_{\mathbb{Q}}\left\{
e_{1},\ e_{2}\right\}  $ and $\operatorname*{span}\nolimits_{\mathbb{Q}%
}\left\{  e_{1}+e_{2},\ e_{3}+e_{4}\right\}  $ of $\mathbb{Q}^{4}$ are not
linearly disjoint and thus don't have an internal direct sum (their sum
$\operatorname*{span}\nolimits_{\mathbb{Q}}\left\{  e_{1},\ e_{2}%
,\ e_{1}+e_{2},\ e_{3}+e_{4}\right\}  $ is only $3$-dimensional, since the
vector $e_{1}+e_{2}$ is redundant).

You can think of linear disjointness as a version of linear dependence for
submodules instead of single vectors. It is different in one aspect, however:
It does not depend on $R$. In fact, the definition of linear disjointness
involves only the addition on $V$, not the action of $R$ on $V$. Thus, whether
or not some $R$-submodules $W_{1},W_{2},\ldots,W_{k}$ are linearly disjoint
does not change if we replace $R$ by a subring of $R$. We could just as well
have replaced \textquotedblleft$R$-module\textquotedblright\ by
\textquotedblleft additive abelian group\textquotedblright\ throughout
Definition \ref{def.mod.sum} (though this would not have increased its
generality, since additive abelian groups are just $\mathbb{Z}$-modules).
\medskip

Here is a well-known alternative characterization of linear disjointness:

\begin{proposition}
\label{prop.mod.dirsum=uniq}Let $R$ be a ring. Let $V$ be a left $R$-module.
Let $W_{1},W_{2},\ldots,W_{k}$ be $R$-submodules of $V$. Then, the following
statements are equivalent:

\begin{itemize}
\item The sum $W_{1}+W_{2}+\cdots+W_{k}$ is direct (i.e., the $R$-submodules
$W_{1},W_{2},\ldots,W_{k}$ are linearly disjoint).

\item The only $k$-tuple $\left(  w_{1},w_{2},\ldots,w_{k}\right)  \in
W_{1}\times W_{2}\times\cdots\times W_{k}$ that satisfies $w_{1}+w_{2}%
+\cdots+w_{k}=0$ is the $k$-tuple $\left(  0,0,\ldots,0\right)  $.

\item For each $v\in V$, there exists \textbf{at most one} $k$-tuple $\left(
w_{1},w_{2},\ldots,w_{k}\right)  \in W_{1}\times W_{2}\times\cdots\times
W_{k}$ that satisfies $w_{1}+w_{2}+\cdots+w_{k}=v$.
\end{itemize}
\end{proposition}

Linear disjointness is particularly simple for two submodules:

\begin{proposition}
\label{prop.mod.dirsumW2}Let $R$ be a ring. Let $V$ be a left $R$-module. Let
$W_{1}$ and $W_{2}$ be two $R$-submodules of $V$. Then, the sum $W_{1}+W_{2}$
is direct if and only if $W_{1}\cap W_{2}=\left\{  0\right\}  $.
\end{proposition}

However, the case of several submodules is less simple:

\begin{warning}
Let $R$ be a ring. Let $V$ be a left $R$-module. Let $W_{1},W_{2},\ldots
,W_{k}$ be $R$-submodules of $V$. In order for the sum $W_{1}+W_{2}%
+\cdots+W_{k}$ to be direct, it is necessary \textbf{but not sufficient} that
$W_{i}\cap W_{j}=\left\{  0\right\}  $ for all $i\neq j$. In other words,
linear disjointness of $W_{1},W_{2},\ldots,W_{k}$ is \textbf{not} a
\textquotedblleft pairwise\textquotedblright\ property.

For an example, the three $\mathbb{Q}$-vector subspaces $\operatorname*{span}%
\nolimits_{\mathbb{Q}}\left\{  e_{1}\right\}  $ and $\operatorname*{span}%
\nolimits_{\mathbb{Q}}\left\{  e_{2}\right\}  $ and $\operatorname*{span}%
\nolimits_{\mathbb{Q}}\left\{  e_{1}+e_{2}\right\}  $ of $\mathbb{Q}^{2}$ are
not linearly disjoint, even though any two of them are.
\end{warning}

\begin{exercise}
\fbox{2} Let $R$ be a ring. Let $V$ be a left $R$-module. Let $X_{1}%
,X_{2},\ldots,X_{i},Y_{1},Y_{2},\ldots,Y_{j}$ be $R$-submodules of $V$. Prove
that the sum $X_{1}+X_{2}+\cdots+X_{i}+Y_{1}+Y_{2}+\cdots+Y_{j}$ is direct if
and only if the three sums%
\begin{align*}
&  X_{1}+X_{2}+\cdots+X_{i},\\
&  Y_{1}+Y_{2}+\cdots+Y_{j},\ \ \ \ \ \ \ \ \ \ \text{and}\\
&  \left(  X_{1}+X_{2}+\cdots+X_{i}\right)  +\left(  Y_{1}+Y_{2}+\cdots
+Y_{j}\right)
\end{align*}
are direct. (Note that the last of these sums is a sum of two addends only!
And note that this is not obvious.)
\end{exercise}

Let us now introduce a notation for the sort of $R$-submodules that are part
of a direct sum:

\begin{definition}
\label{def.mod.diradd}Let $R$ be a ring. Let $W$ be an $R$-submodule of a left
$R$-module $V$. We say that $W$ is a \emph{direct addend} of $V$ if there
exists a further $R$-submodule $U$ of $V$ such that $V=W\oplus U$ (internal
direct sum). In this case, the latter $R$-submodule $U$ is called a
\emph{complement} of $W$ in $V$. (Not \textquotedblleft the
complement\textquotedblright\ because $U$ is not uniquely determined by $V$
and $W$.)
\end{definition}

\begin{example}
\label{exa.mod.diradd.11}The $\mathbb{Z}$-submodule $\operatorname*{span}%
\nolimits_{\mathbb{Z}}\left\{  \left(  1,1\right)  \right\}  $ of
$\mathbb{Z}^{2}$ is a direct addend of $\mathbb{Z}^{2}$, since
$\operatorname*{span}\nolimits_{\mathbb{Z}}\left\{  \left(  0,1\right)
\right\}  $ is a complement of it. It has many other complements as well. In
fact, it is not hard to see that two \textquotedblleft
one-dimensional\textquotedblright\ $\mathbb{Z}$-submodules
$\operatorname*{span}\nolimits_{\mathbb{Z}}\left\{  u\right\}  $ and
$\operatorname*{span}\nolimits_{\mathbb{Z}}\left\{  v\right\}  $ of
$\mathbb{Z}^{2}$ are complements of each other if and only if $\left(
u,v\right)  $ is a basis of $\mathbb{Z}^{2}$.
\end{example}

\begin{example}
\label{exa.mod.diradd.2Z}The $\mathbb{Z}$-submodule $2\mathbb{Z}=\left\{
\text{even integers}\right\}  $ of $\mathbb{Z}$ is not a direct addend, since
it has no complement. (The easiest way to see this is to show that no nonzero
submodule of $\mathbb{Z}$ is linearly disjoint from $2\mathbb{Z}$, but
$2\mathbb{Z}\oplus0$ is not $\mathbb{Z}$.) That is, we cannot write
$\mathbb{Z}$ as an internal direct sum $2\mathbb{Z}\oplus\left(  \text{another
submodule of }\mathbb{Z}\right)  $.
\end{example}

\begin{remark}
As we have seen above, the linear disjointness of two given $R$-submodules of
a given $R$-module $V$ does not depend on $R$. However, whether a given
$R$-submodule is a direct addend does depend on $R$: It can have a complement
for some choices of $R$ but not for others. In Example
\ref{exa.rep.G-rep.maschke.kn0}, we will see an instance of this.
\end{remark}

As we see from Example \ref{exa.mod.diradd.2Z}, not every submodule is a
direct addend. But over a field, it is:

\begin{proposition}
\label{prop.mod.diradd.field}Let $\mathbf{k}$ be a field. Then, any vector
subspace of a $\mathbf{k}$-vector space $V$ is a direct addend of $V$.
\end{proposition}

\begin{fineprint}
\begin{proof}
[Proof idea.]Let $W$ be a vector subspace of $V$. We must show that $W$ is a
direct addend of $V$.

Pick a basis $\left(  e_{i}\right)  _{i\in I}$ of $W$, and extend it to a
basis $\left(  e_{i}\right)  _{i\in J}$ of $V$. (Both of these can be done,
since $\mathbf{k}$ is a field.)

Let $U$ be the span of the \textquotedblleft new\textquotedblright\ basis
elements -- i.e., let $U=\operatorname*{span}\nolimits_{\mathbf{k}}\left\{
e_{i}\ \mid\ i\in J\setminus I\right\}  $. Then, $U$ is a $\mathbf{k}%
$-submodule of $V$. Moreover, it is easy to see that $W\cap U=\left\{
0\right\}  $ (since $W$ and $U$ are spanned by two disjoint subsets of the
basis $\left(  e_{i}\right)  _{i\in J}$). Hence, the sum $W+U$ is direct (by
Proposition \ref{prop.mod.dirsumW2}). Thus, $W+U=W\oplus U$. However,
$W=\operatorname*{span}\nolimits_{\mathbf{k}}\left\{  e_{i}\ \mid\ i\in
I\right\}  $ (since $\left(  e_{i}\right)  _{i\in I}$ is a basis of $W$) and
$U=\operatorname*{span}\nolimits_{\mathbf{k}}\left\{  e_{i}\ \mid\ i\in
J\setminus I\right\}  $, so that%
\begin{align*}
W+U  &  =\operatorname*{span}\nolimits_{\mathbf{k}}\left\{  e_{i}\ \mid\ i\in
I\right\}  +\operatorname*{span}\nolimits_{\mathbf{k}}\left\{  e_{i}%
\ \mid\ i\in J\setminus I\right\} \\
&  =\operatorname*{span}\nolimits_{\mathbf{k}}\underbrace{\left(  \left\{
e_{i}\ \mid\ i\in I\right\}  \cup\left\{  e_{i}\ \mid\ i\in J\setminus
I\right\}  \right)  }_{=\left\{  e_{i}\ \mid\ i\in J\right\}  }\\
&  =\operatorname*{span}\nolimits_{\mathbf{k}}\left\{  e_{i}\ \mid\ i\in
J\right\}  =V
\end{align*}
(since $\left(  e_{i}\right)  _{i\in J}$ is a basis of $V$). Thus,
$V=W+U=W\oplus U$. This shows that $W$ is a direct addend of $V$. Proposition
\ref{prop.mod.diradd.field} is thus proved.
\end{proof}
\end{fineprint}

The direct addends of an $R$-module can also be characterized in terms of
projections. We recall what a projection is:

\begin{definition}
\label{def.proj}Let $W$ be a subset of a set $V$. A \emph{projection} from $V$
to $W$ means a map $\pi:V\rightarrow W$ such that $\left.  \pi\mid_{W}\right.
=\operatorname*{id}$ (that is, $\pi\left(  w\right)  =w$ for each $w\in W$).
\end{definition}

Clearly, any projection is surjective. An example of a projection is the
\textquotedblleft floor function\textquotedblright\ -- i.e., the map
$\mathbb{R}\rightarrow\mathbb{Z},\ x\mapsto\left\lfloor x\right\rfloor $ which
sends each real number $x$ to its integer part $\left\lfloor x\right\rfloor $
(that is, the largest integer $\leq x$). But we are more interested in
$R$-linear projections.

\begin{proposition}
\label{prop.mod.diradd.proj}Let $R$ be a ring. Let $V$ be a left $R$-module.
Let $W$ be an $R$-submodule of $V$. Then, $W$ is a direct addend of $V$ if and
only if there exists a left $R$-linear projection $\pi:V\rightarrow W$.
\end{proposition}

\begin{proof}
$\Longrightarrow:$ Assume that $W$ is a direct addend. Thus, $V=W\oplus U$
(internal direct sum) for some $R$-submodule $U$ of $V$. Consider this $U$.

Let us agree to use the symbol \textquotedblleft$\oplus$\textquotedblright%
\ only for internal direct sums in this proof, while using the symbol
\textquotedblleft$\times$\textquotedblright\ for external direct sums. (This
is legitimate, since we are only dealing with direct sums of finitely many
$R$-modules, and these direct sums are the same as the corresponding direct products.)

We have $V=W\oplus U=W+U$, and moreover, the sum $W+U$ is direct (since we are
writing it as $W\oplus U$). The latter fact means that the map%
\begin{align*}
\operatorname*{sum}:W\times U  &  \rightarrow W+U,\\
\left(  w,u\right)   &  \mapsto w+u
\end{align*}
(recall that $W\times U$ is the external direct sum of $W$ and $U$) is
injective. Hence, this map $\operatorname*{sum}$ is bijective (since it is
clearly surjective), and therefore an $R$-module isomorphism (since it is
obviously left $R$-linear). Thus, its inverse map $\operatorname*{sum}%
\nolimits^{-1}:W+U\rightarrow W\times U$ is an $R$-module morphism as well. In
view of $W+U=V$, we can rewrite the latter map as $\operatorname*{sum}%
\nolimits^{-1}:V\rightarrow W\times U$. Explicitly, this map
$\operatorname*{sum}\nolimits^{-1}$ sends each $v\in V$ to the unique pair
$\left(  w,u\right)  \in W\times U$ that satisfies $v=w+u$.

Let $\rho:W\times U\rightarrow W$ be the map that sends each pair $\left(
w,u\right)  $ to its first entry $w$. This is also a left $R$-linear map.

Consider the composition $\rho\circ\operatorname*{sum}\nolimits^{-1}%
:V\rightarrow W$ of these two maps $\rho:W\times U\rightarrow W$ and
$\operatorname*{sum}\nolimits^{-1}:V\rightarrow W\times U$. This composition
is left $R$-linear (since both $\rho$ and $\operatorname*{sum}\nolimits^{-1}$
are left $R$-linear). Moreover, each $w\in W$ satisfies $\operatorname*{sum}%
\left(  w,0\right)  =w+0=w$ and thus $\operatorname*{sum}\nolimits^{-1}\left(
w\right)  =\left(  w,0\right)  $, so that
\[
\left(  \rho\circ\operatorname*{sum}\nolimits^{-1}\right)  \left(  w\right)
=\rho\left(  \underbrace{\operatorname*{sum}\nolimits^{-1}\left(  w\right)
}_{=\left(  w,0\right)  }\right)  =\rho\left(  w,0\right)  =w
\]
(by the definition of $\rho$). In other words, $\left.  \left(  \rho
\circ\operatorname*{sum}\nolimits^{-1}\right)  \mid_{W}\right.
=\operatorname*{id}$. Hence, the map $\rho\circ\operatorname*{sum}%
\nolimits^{-1}:V\rightarrow W$ is a projection from $V$ to $W$. Thus, there
exists a left $R$-linear projection $\pi:V\rightarrow W$ (namely, $\pi
=\rho\circ\operatorname*{sum}\nolimits^{-1}$). This proves the
\textquotedblleft$\Longrightarrow$\textquotedblright\ direction of Proposition
\ref{prop.mod.diradd.proj}. \medskip

$\Longleftarrow:$ Assume that there is a left $R$-linear projection
$\pi:V\rightarrow W$. Consider this $\pi$.

Now define $U=\operatorname*{Ker}\pi$. This is clearly an $R$-submodule of
$V$. We shall now show that $V=W\oplus U$ (internal direct sum). Indeed:

\begin{itemize}
\item The sum $W+U$ is direct. [\textit{Proof:} Let $x\in W\cap U$. Thus,
$x\in W\cap U\subseteq W$, so that $\pi\left(  x\right)  =x$ (since
$\pi:V\rightarrow W$ is a projection and thus satisfies $\pi\left(  w\right)
=w$ for all $w\in W$). But we also have $x\in W\cap U\subseteq
U=\operatorname*{Ker}\pi$ and thus $\pi\left(  x\right)  =0$. Comparing this
with $\pi\left(  x\right)  =x$, we find $x=0$. Forget that we fixed $x$. We
thus have shown that $x=0$ for each $x\in W\cap U$. In other words, $W\cap
U\subseteq\left\{  0\right\}  $. Hence, $W\cap U=\left\{  0\right\}  $ (since
$W\cap U$ is an $R$-submodule of $V$). Thus, the sum $W+U$ is direct (by
Proposition \ref{prop.mod.dirsumW2}).]

\item We have $V=W+U$. [\textit{Proof:} Let $v\in V$. Then, $\pi\left(
v\right)  \in W$ (since $\pi$ is a map from $V$ to $W$). But $\pi:V\rightarrow
W$ is a projection and thus satisfies $\pi\left(  w\right)  =w$ for all $w\in
W$. Applying this to $w=\pi\left(  v\right)  $, we obtain $\pi\left(
\pi\left(  v\right)  \right)  =\pi\left(  v\right)  $ (since $\pi\left(
v\right)  \in W$). Since $\pi$ is a linear map, we have $\pi\left(
v-\pi\left(  v\right)  \right)  =\pi\left(  v\right)  -\underbrace{\pi\left(
\pi\left(  v\right)  \right)  }_{=\pi\left(  v\right)  }=\pi\left(  v\right)
-\pi\left(  v\right)  =0$. Thus, $v-\pi\left(  v\right)  \in
\operatorname*{Ker}\pi=U$. Hence, $v=\underbrace{\pi\left(  v\right)  }_{\in
W}+\underbrace{v-\pi\left(  v\right)  }_{\in U}\in W+U$.

Forget that we fixed $v$. We thus have shown that $v\in W+U$ for each $v\in
V$. In other words, $V\subseteq W+U$. Hence, $V=W+U$ (since we clearly have
$W+U\subseteq V$).]
\end{itemize}

Now, $V=W+U=W\oplus U$ (since the sum $W+U$ is direct). Hence, $W$ is a direct
addend of $V$. This proves the \textquotedblleft$\Longleftarrow$%
\textquotedblright\ direction of Proposition \ref{prop.mod.diradd.proj}.
\end{proof}

\begin{example}
\label{exa.rep.G-rep.maschke.kn0}Assume that $n>0$. Consider again the natural
representation $\mathbf{k}^{n}$ of $S_{n}$. This is a left $\mathbf{k}\left[
S_{n}\right]  $-module. Define its $\mathbf{k}\left[  S_{n}\right]
$-submodule $D\left(  \mathbf{k}^{n}\right)  $ as in Subsection
\ref{subsec.rep.G-rep.subreps}.

For any $n$-tuple $\left(  a_{1},a_{2},\ldots,a_{n}\right)  \in\mathbf{k}^{n}%
$, set $a_{\Sigma}:=a_{1}+a_{2}+\cdots+a_{n}$. \medskip

\textbf{(a)} Consider the $\mathbf{k}$-linear map%
\begin{align*}
f_{1}:\mathbf{k}^{n}  &  \rightarrow D\left(  \mathbf{k}^{n}\right)  ,\\
\left(  a_{1},a_{2},\ldots,a_{n}\right)   &  \mapsto\left(  a_{\Sigma
},a_{\Sigma},\ldots,a_{\Sigma}\right)  .
\end{align*}
This map $f_{1}$ is easily seen to be $S_{n}$-equivariant, and thus is a
morphism of representations of $S_{n}$ (since it is $\mathbf{k}$-linear).
Hence, by Proposition \ref{prop.rep.G-rep.mor=mor}, this map $f_{1}$ is left
$\mathbf{k}\left[  S_{n}\right]  $-linear.

However, the map $f_{1}$ is not a projection (unless $n=1$ in $\mathbf{k}$),
since it sends a vector $\left(  a,a,\ldots,a\right)  \in D\left(
\mathbf{k}^{n}\right)  $ to $\left(  na,na,\ldots,na\right)  $, which is
generally not the same vector. \medskip

\textbf{(b)} Consider the $\mathbf{k}$-linear map%
\begin{align*}
f_{2}:\mathbf{k}^{n}  &  \rightarrow D\left(  \mathbf{k}^{n}\right)  ,\\
\left(  a_{1},a_{2},\ldots,a_{n}\right)   &  \mapsto\left(  a_{1},a_{1}%
,\ldots,a_{1}\right)  .
\end{align*}
This map $f_{2}$ is a projection (since it sends each vector $\left(
a,a,\ldots,a\right)  \in D\left(  \mathbf{k}^{n}\right)  $ to itself). But it
is not $S_{n}$-equivariant, and thus not left $\mathbf{k}\left[  S_{n}\right]
$-linear (unless $n=1$). \medskip

\textbf{(c)} Let us now assume that $n$ is invertible in $\mathbf{k}$.
Consider the $\mathbf{k}$-linear map%
\begin{align*}
f_{3}:\mathbf{k}^{n}  &  \rightarrow D\left(  \mathbf{k}^{n}\right)  ,\\
\left(  a_{1},a_{2},\ldots,a_{n}\right)   &  \mapsto\left(  \dfrac{a_{\Sigma}%
}{n},\dfrac{a_{\Sigma}}{n},\ldots,\dfrac{a_{\Sigma}}{n}\right)  .
\end{align*}
(In other words, this map $f_{3}$ replaces all entries of an $n$-tuple by
their average.)

It is easy to see that this map $f_{3}$ is both left $\mathbf{k}\left[
S_{n}\right]  $-linear and a projection. It is thus a left $\mathbf{k}\left[
S_{n}\right]  $-linear projection.

Had we not required $n$ to be invertible in $\mathbf{k}$, such a map would not
exist. Thus, the existence of an $S_{n}$-equivariant map and the existence of
a $\mathbf{k}$-linear projection do not guarantee the existence of an $S_{n}%
$-equivariant projection.
\end{example}

We note one last general property of direct addends before we resume studying
modules over group algebras:

\begin{proposition}
\label{prop.mod.diradd.V/W}Let $R$ be a ring. Let $V$ be a left $R$-module.
Let $W$ be a direct addend of $V$. Then,%
\[
V\cong W\oplus\left(  V/W\right)  \ \ \ \ \ \ \ \ \ \ \text{as left
}R\text{-modules.}%
\]
(Here, $W\oplus\left(  V/W\right)  $ means an external direct sum, since $V/W$
is not a submodule of $V$.)
\end{proposition}

\begin{proof}
Proceed as in the proof of the \textquotedblleft$\Longrightarrow
$\textquotedblright\ direction of Proposition \ref{prop.mod.diradd.proj}.
Thus, write $V=W\oplus U$ (internal direct sum) for some $R$-submodule $U$ of
$V$, and observe that the map%
\begin{align*}
\operatorname*{sum}:W\times U  &  \rightarrow W+U,\\
\left(  w,u\right)   &  \mapsto w+u
\end{align*}
is an $R$-module isomorphism. Thus, this map $\operatorname*{sum}$ has an
inverse map $\operatorname*{sum}\nolimits^{-1}:W+U\rightarrow W\times U$,
which is also an $R$-module isomorphism. Since $V=W\oplus U=W+U$, we can
rewrite this as follows: The map $\operatorname*{sum}$ has an inverse map
$\operatorname*{sum}\nolimits^{-1}:V\rightarrow W\times U$, which is also an
$R$-module isomorphism. Hence, in particular, the map $\operatorname*{sum}%
\nolimits^{-1}$ is $R$-linear and bijective, hence surjective.

Furthermore, consider the $R$-linear map $\tau:W\times U\rightarrow U$ that
sends each pair $\left(  w,u\right)  $ to its second entry $u$. This map
$\tau$ is clearly surjective (since each $u\in U$ is the image of the pair
$\left(  0,u\right)  $ under $\tau$).

Let $\omega:V\rightarrow U$ be the composition $\tau\circ\operatorname*{sum}%
\nolimits^{-1}$ of the two maps $\tau:W\times U\rightarrow U$ and
$\operatorname*{sum}\nolimits^{-1}:V\rightarrow W\times U$. This composition
$\omega$ is surjective (since both $\tau$ and $\operatorname*{sum}%
\nolimits^{-1}$ are surjective) and left $R$-linear (since both $\tau$ and
$\operatorname*{sum}\nolimits^{-1}$ are $R$-linear).

Moreover, it is easy to see that the kernel $\operatorname*{Ker}\omega$ of
this composition is $W$. [\textit{Proof:} Let $v\in V$. Let $\left(
w,u\right)  =\operatorname*{sum}\nolimits^{-1}\left(  v\right)  $. Thus,
$\left(  w,u\right)  \in W\times U$ and $\operatorname*{sum}\left(
w,u\right)  =v$, so that $\operatorname*{sum}\nolimits^{-1}\left(  v\right)
=\left(  w,u\right)  $. Now, from $\omega=\tau\circ\operatorname*{sum}%
\nolimits^{-1}$, we obtain%
\begin{align*}
\omega\left(  v\right)   &  =\left(  \tau\circ\operatorname*{sum}%
\nolimits^{-1}\right)  \left(  v\right)  =\tau\left(
\underbrace{\operatorname*{sum}\nolimits^{-1}\left(  v\right)  }_{=\left(
w,u\right)  }\right)  =\tau\left(  w,u\right) \\
&  =u\ \ \ \ \ \ \ \ \ \ \left(  \text{by the definition of }\tau\right)  .
\end{align*}

Now, if $u=0$, then $v=\operatorname*{sum}\left(  w,u\right)
=w+\underbrace{u}_{=0}=w\in W$. Conversely, if $v\in W$, then $\left(
v,0\right)  \in W\times U$ and therefore $\operatorname*{sum}\left(
v,0\right)  $ is well-defined and satisfies $\operatorname*{sum}\left(
v,0\right)  =v+0=v$, so that $\left(  v,0\right)  =\operatorname*{sum}%
\nolimits^{-1}\left(  v\right)  =\left(  w,u\right)  $ and therefore $0=u$, so
that $u=0$. Thus, each of the two statements $u=0$ and $v\in W$ implies the
other. In other words, we have the equivalence $\left(  u=0\right)
\ \Longleftrightarrow\ \left(  v\in W\right)  $. Now, we have the chain of
equivalences
\[
\left(  v\in\operatorname*{Ker}\omega\right)  \ \Longleftrightarrow\ \left(
\underbrace{\omega\left(  v\right)  }_{=u}=0\right)  \ \Longleftrightarrow
\ \left(  u=0\right)  \ \Longleftrightarrow\ \left(  v\in W\right)  .
\]

Forget that we fixed $v$. We thus have shown that the equivalence $\left(
v\in\operatorname*{Ker}\omega\right)  \ \Longleftrightarrow\ \left(  v\in
W\right)  $ holds for each $v\in V$. Therefore, $\operatorname*{Ker}\omega=W$
(since both $\operatorname*{Ker}\omega$ and $W$ are subsets of $V$).]

Now, the first isomorphism theorem (applied to the $R$-linear map
$\omega:V\rightarrow U$) yields%
\[
V/\operatorname*{Ker}\omega\cong\omega\left(  V\right)
=U\ \ \ \ \ \ \ \ \ \ \left(  \text{since }\omega\text{ is surjective}\right)
.
\]
In view of $\operatorname*{Ker}\omega=W$, we can rewrite this as%
\[
V/W\cong U.
\]
Thus, $W\times\left(  V/W\right)  \cong W\times U$ (because of a general fact
that if $X$, $Y$ and $Z$ are three $R$-modules satisfying $Y\cong Z$, then
$X\times Y\cong X\times Z$). However, we also have $V\cong W\times U$ (since
$\operatorname*{sum}\nolimits^{-1}:V\rightarrow W\times U$ is an $R$-module
isomorphism). Altogether, we thus find $V\cong W\times U\cong W\times\left(
V/W\right)  =W\oplus\left(  V/W\right)  $ (since a direct product of finitely
many $R$-modules is the same as their direct sum). This proves Proposition
\ref{prop.mod.diradd.V/W}.
\end{proof}

\subsubsection{\label{subsec.rep.maschke.maschke}Maschke's theorem}

Recall that if $G$ is a group, then the left $\mathbf{k}\left[  G\right]
$-modules are just the representations of $G$ over $\mathbf{k}$. Thus, we
obtain a notion of a direct addend for representations of $G$.

So when is a subrepresentation of a representation $V$ of $G$ a direct addend?
i.e., when is there a $\mathbf{k}\left[  G\right]  $-linear projection from
$V$ to the subrepresentation? The following theorem (known as \emph{Maschke's
theorem}) yields a surprisingly general answer to this question, which applies
not quite always but in many situations:

\begin{theorem}
[Maschke's theorem]\label{thm.rep.G-rep.maschke}Let $G$ be a finite group. Let
$V$ be a left $\mathbf{k}\left[  G\right]  $-module (i.e., a representation of
$G$ over $\mathbf{k}$). Let $W$ be a left $\mathbf{k}\left[  G\right]
$-submodule of $V$ (that is, a subrepresentation of $V$). Assume that

\begin{enumerate}
\item the number $\left\vert G\right\vert $ is invertible in $\mathbf{k}$;

\item the $\mathbf{k}$-submodule $W$ is a direct addend of $V$ as a
$\mathbf{k}$-module.
\end{enumerate}

Then, the left $\mathbf{k}\left[  G\right]  $-submodule $W$ is also a direct
addend of $V$ as a left $\mathbf{k}\left[  G\right]  $-module.
\end{theorem}

\begin{remark}
\label{rmk.rep.G-rep.maschke-k}\textbf{(a)} If $\mathbf{k}$ is a field of
characteristic $0$, or (more generally) a commutative $\mathbb{Q}$-algebra,
then the first assumption of Theorem \ref{thm.rep.G-rep.maschke} is
automatically satisfied. \medskip

\textbf{(b)} If $\mathbf{k}$ is a field, then the second assumption of Theorem
\ref{thm.rep.G-rep.maschke} is automatically satisfied (by Proposition
\ref{prop.mod.diradd.field}). \medskip

\textbf{(c)} Thus, if $\mathbf{k}$ is a field of characteristic $0$, then both
assumptions of Theorem \ref{thm.rep.G-rep.maschke} are automatically
satisfied, and therefore the theorem simply says that any left $\mathbf{k}%
\left[  G\right]  $-submodule $W$ of any left $\mathbf{k}\left[  G\right]
$-module $V$ is a direct addend. This is the best-known case of Maschke's theorem.
\end{remark}

\begin{example}
\label{exa.rep.G-rep.maschke.kn}Recall the natural representation
$\mathbf{k}^{n}$ of $S_{n}$ and its subrepresentation $D\left(  \mathbf{k}%
^{n}\right)  $. Then, $D\left(  \mathbf{k}^{n}\right)  $ is always a direct
addend of $\mathbf{k}^{n}$ as a $\mathbf{k}$-module (e.g., with complement
$\operatorname*{span}\nolimits_{\mathbf{k}}\left\{  e_{1},e_{2},\ldots
,e_{n-1}\right\}  $), but not always a direct addend of $\mathbf{k}^{n}$ as a
left $\mathbf{k}\left[  S_{n}\right]  $-module. Maschke's theorem (Theorem
\ref{thm.rep.G-rep.maschke}) guarantees that $D\left(  \mathbf{k}^{n}\right)
$ is a direct addend of $\mathbf{k}^{n}$ as a left $\mathbf{k}\left[
S_{n}\right]  $-module whenever $\left\vert S_{n}\right\vert =n!$ is
invertible in $\mathbf{k}$, but we already know from Proposition
\ref{prop.rep.G-rep.Sn-nat.subreps.dirsum} \textbf{(a)} that this is true
under a weaker assumption (namely, whenever $n$ is invertible in $\mathbf{k}%
$). Thus, Maschke's theorem does not tell us anything new here. But due to its
generality, Maschke's theorem can be used for many other representations,
which are much harder to analyze \textquotedblleft by hand\textquotedblright%
\ than $\mathbf{k}^{n}$.
\end{example}

\begin{proof}
[Proof of Theorem \ref{thm.rep.G-rep.maschke}.]We have assumed that $W$ is a
direct addend of the $\mathbf{k}$-module $V$. Hence, there exists a
$\mathbf{k}$-linear projection $\pi:V\rightarrow W$ (by Proposition
\ref{prop.mod.diradd.proj}, applied to $R=\mathbf{k}$). Consider this $\pi$.
Note that $\pi$ is only $\mathbf{k}$-linear, but usually not left
$\mathbf{k}\left[  G\right]  $-linear.

Define a new map $\pi^{\prime}:V\rightarrow W$ by setting%
\begin{equation}
\pi^{\prime}\left(  v\right)  :=\dfrac{1}{\left\vert G\right\vert }\sum_{g\in
G}g^{-1}\pi\left(  gv\right)  \ \ \ \ \ \ \ \ \ \ \text{for each }v\in V.
\label{pf.thm.rep.G-rep.maschke.pi}%
\end{equation}
We claim that this is a left $\mathbf{k}\left[  G\right]  $-linear projection
from $V$ to $W$. Indeed, let us prove the $\mathbf{k}\left[  G\right]
$-linearity and the \textquotedblleft projection-ness\textquotedblright\ separately:

\begin{itemize}
\item \textit{Proof of }$\mathbf{k}\left[  G\right]  $\textit{-linearity:} It
is clear that the map $\pi^{\prime}$ is $\mathbf{k}$-linear (since $\pi$ is
$\mathbf{k}$-linear, whereas the actions of each $g\in G$ and each $g^{-1}$
are $\mathbf{k}$-linear). Now, we shall show that $\pi^{\prime}$ is
$G$-equivariant, i.e., that we have $\pi^{\prime}\left(  hv\right)
=h\pi^{\prime}\left(  v\right)  $ for all $h\in G$ and $v\in V$.

Indeed, let $h\in G$ and $v\in V$. Recall that $G$ is a group. Thus, the map
$G\rightarrow G,\ g\mapsto gh$ is a bijection (with inverse $g\mapsto gh^{-1}%
$). The definition of $\pi^{\prime}$ yields%
\[
\pi^{\prime}\left(  hv\right)  =\dfrac{1}{\left\vert G\right\vert }\sum_{g\in
G}g^{-1}\pi\left(  ghv\right)  .
\]
Comparing this with%
\begin{align*}
h\pi^{\prime}\left(  v\right)   &  =h\cdot\dfrac{1}{\left\vert G\right\vert
}\sum_{g\in G}g^{-1}\pi\left(  gv\right)  \ \ \ \ \ \ \ \ \ \ \left(  \text{by
(\ref{pf.thm.rep.G-rep.maschke.pi})}\right) \\
&  =\dfrac{1}{\left\vert G\right\vert }\sum_{g\in G}hg^{-1}\pi\left(
gv\right)  =\dfrac{1}{\left\vert G\right\vert }\sum_{g\in G}%
\underbrace{h\left(  gh\right)  ^{-1}}_{=hh^{-1}g^{-1}=g^{-1}}\pi\left(
ghv\right) \\
&  \ \ \ \ \ \ \ \ \ \ \ \ \ \ \ \ \ \ \ \ \left(
\begin{array}
[c]{c}%
\text{here, we have substituted }gh\text{ for }g\text{ in the sum,}\\
\text{since the map }G\rightarrow G,\ g\mapsto gh\text{ is a bijection}%
\end{array}
\right) \\
&  =\dfrac{1}{\left\vert G\right\vert }\sum_{g\in G}g^{-1}\pi\left(
ghv\right)  ,
\end{align*}
we find $\pi^{\prime}\left(  hv\right)  =h\pi^{\prime}\left(  v\right)  $.

Forget that we fixed $h$ and $v$. We thus have shown that $\pi^{\prime}\left(
hv\right)  =h\pi^{\prime}\left(  v\right)  $ for all $h\in G$ and $v\in V$. In
other words, the map $\pi^{\prime}$ is $G$-equivariant. Hence, $\pi^{\prime}$
is a morphism of $G$-representations (since $\pi^{\prime}$ is $\mathbf{k}%
$-linear). By Proposition \ref{prop.rep.G-rep.mor=mor}, this entails that
$\pi^{\prime}$ is a $\mathbf{k}\left[  G\right]  $-module morphism, i.e., a
left $\mathbf{k}\left[  G\right]  $-linear map.

\item \textit{Proof of \textquotedblleft projection-ness\textquotedblright:}
We must now show that $\pi^{\prime}:V\rightarrow W$ is a projection. In other
words, we must show that $\pi^{\prime}\left(  w\right)  =w$ for each $w\in W$.

Let $w\in W$. Then, each $g\in G$ satisfies $gw\in W$ (since $W$ is a
subrepresentation of $V$ and thus a left $\mathbf{k}\left[  G\right]
$-submodule of $V$) and thus%
\begin{equation}
\pi\left(  gw\right)  =\left(  \pi\mid_{W}\right)  \left(  gw\right)  =gw
\label{pf.thm.rep.G-rep.maschke.pigw}%
\end{equation}
(since $\pi:V\rightarrow W$ is a projection and thus satisfies $\left.
\pi\mid_{W}\right.  =\operatorname*{id}$). Now, the definition of $\pi
^{\prime}$ yields%
\begin{align*}
\pi^{\prime}\left(  w\right)   &  =\dfrac{1}{\left\vert G\right\vert }%
\sum_{g\in G}g^{-1}\underbrace{\pi\left(  gw\right)  }%
_{\substack{=gw\\\text{(by (\ref{pf.thm.rep.G-rep.maschke.pigw}))}}}=\dfrac
{1}{\left\vert G\right\vert }\sum_{g\in G}\underbrace{g^{-1}g}_{=1}w=\dfrac
{1}{\left\vert G\right\vert }\underbrace{\sum_{g\in G}w}_{=\left\vert
G\right\vert w}=\dfrac{1}{\left\vert G\right\vert }\left\vert G\right\vert w\\
&  =w.
\end{align*}

Forget that we fixed $w$. We thus have shown that $\pi^{\prime}\left(
w\right)  =w$ for each $w\in W$. Hence, the map $\pi^{\prime}:V\rightarrow W$
is a projection.
\end{itemize}

Altogether, we have now shown that $\pi^{\prime}$ is a left $\mathbf{k}\left[
G\right]  $-linear projection from $V$ to $W$. Hence, such a projection
exists. It thus follows that $W$ is a direct addend of $V$ as a left
$\mathbf{k}\left[  G\right]  $-module (by Proposition
\ref{prop.mod.diradd.proj}, applied to $R=\mathbf{k}\left[  G\right]  $). This
proves Theorem \ref{thm.rep.G-rep.maschke}.
\end{proof}

\begin{remark}
What is the motivation behind the definition of the map $\pi^{\prime}$ in the
above proof?

In order to apply Proposition \ref{prop.mod.diradd.proj}, we need a left
$\mathbf{k}\left[  G\right]  $-linear (that is, $\mathbf{k}$-linear and
$G$-equivariant) projection from $V$ to $W$. We are given a $\mathbf{k}%
$-linear projection $\pi:V\rightarrow W$, but it is (usually) not
$G$-equivariant. This lack of $G$-equivariance can be used to define a family
of \textquotedblleft doppelg\"{a}ngers\textquotedblright\ of $\pi$: Namely,
for each $g\in G$, we can consider the \textquotedblleft$g$-tilted
version\textquotedblright\ $\pi_{g}:V\rightarrow W$ of $\pi$ given by%
\[
\pi_{g}\left(  v\right)  :=g^{-1}\pi\left(  gv\right)
\ \ \ \ \ \ \ \ \ \ \text{for each }v\in V.
\]
It is easy to see that each such $\pi_{g}$ is again a $\mathbf{k}$-linear
projection $V\rightarrow W$ (indeed, $\pi_{g}$ is the conjugate of $\pi$ by
the action of $g^{-1}$, in the appropriate sense).

Thus, we now have not one but $\left\vert G\right\vert $ many $\mathbf{k}%
$-linear projections from $V$ to $W$. None of them is (usually) $G$%
-equivariant, but something similar holds: For any $g\in G$, $h\in G$ and
$v\in V$, we have%
\begin{equation}
\pi_{g}\left(  hv\right)  =h\pi_{gh}\left(  v\right)  .
\label{pf.thm.rep.G-rep.maschke.pig}%
\end{equation}
If the maps $\pi_{g}$ and $\pi_{gh}$ were the same, then this would show that
this map is $G$-equivariant. But they are different. How do we salvage them?
Simply: We take the average of the $\pi_{g}$ over all group elements $g\in G$.
Upon averaging over all $g\in G$, the $\pi_{g}$ and the $\pi_{gh}$ in the
equality (\ref{pf.thm.rep.G-rep.maschke.pig}) become the same, since
$\sum_{g\in G}\pi_{g}=\sum_{g\in G}\pi_{gh}$. Thus, the average of the
$\pi_{g}$ over all $g\in G$ is $G$-equivariant. Moreover, this average is
still a projection (since any average of $\mathbf{k}$-linear projections from
$V$ to $W$ is again a projection from $V$ to $W$). This average is precisely
the map $\pi^{\prime}$ in the above proof.
\end{remark}

\addtocounter{theo}{-1}

\begin{example}
\label{exa.rep.G-rep.maschke.natD}Let us see how the above proof of Theorem
\ref{thm.rep.G-rep.maschke} works in a specific situation. Namely, we let $n$
be a positive integer, let $G$ be the group $S_{n}$, and let $V$ be its
natural representation $\mathbf{k}^{n}$. Let $W$ be the subrepresentation
$D\left(  \mathbf{k}^{n}\right)  =\left\{  \left(  a,a,\ldots,a\right)
\ \mid\ a\in\mathbf{k}\right\}  $ of $V$. Assume that $\left\vert G\right\vert
=\left\vert S_{n}\right\vert =n!$ is invertible in $\mathbf{k}$. We pick the
$\mathbf{k}$-linear projection%
\begin{align*}
\pi:V  &  \rightarrow W,\\
\left(  a_{1},a_{2},\ldots,a_{n}\right)   &  \mapsto\left(  a_{1},a_{1}%
,\ldots,a_{1}\right)
\end{align*}
(which replaces all coordinates of a vector by copies of its first
coordinate). This projection is not $G$-equivariant (unless $n=1$), that is,
not left $\mathbf{k}\left[  G\right]  $-linear. In our above proof of Theorem
\ref{thm.rep.G-rep.maschke}, we constructed a left $\mathbf{k}\left[
G\right]  $-linear projection $\pi^{\prime}$ from $V$ to $W$ by setting%
\[
\pi^{\prime}\left(  v\right)  :=\dfrac{1}{\left\vert G\right\vert }\sum_{g\in
G}g^{-1}\pi\left(  gv\right)  \ \ \ \ \ \ \ \ \ \ \text{for each }v\in V.
\]
Let us make this formula more explicit. For any vector $v=\left(  a_{1}%
,a_{2},\ldots,a_{n}\right)  \in\mathbf{k}^{n}=V$ and any permutation $g\in
S_{n}=G$, we have%
\[
gv=g\rightharpoonup v=g\rightharpoonup\left(  a_{1},a_{2},\ldots,a_{n}\right)
=\left(  a_{g^{-1}\left(  1\right)  },a_{g^{-1}\left(  2\right)  }%
,\ldots,a_{g^{-1}\left(  n\right)  }\right)
\]
and therefore%
\[
\pi\left(  gv\right)  =\pi\left(  a_{g^{-1}\left(  1\right)  },a_{g^{-1}%
\left(  2\right)  },\ldots,a_{g^{-1}\left(  n\right)  }\right)  =\left(
a_{g^{-1}\left(  1\right)  },a_{g^{-1}\left(  1\right)  },\ldots
,a_{g^{-1}\left(  1\right)  }\right)
\]
(by the definition of $\pi$), so that%
\begin{align*}
g^{-1}\pi\left(  gv\right)   &  =g^{-1}\left(  a_{g^{-1}\left(  1\right)
},a_{g^{-1}\left(  1\right)  },\ldots,a_{g^{-1}\left(  1\right)  }\right) \\
&  =\left(  a_{g^{-1}\left(  1\right)  },a_{g^{-1}\left(  1\right)  }%
,\ldots,a_{g^{-1}\left(  1\right)  }\right)
\end{align*}
(since a vector whose entries are all equal does not change when its entries
are permuted). Thus, for any vector $v=\left(  a_{1},a_{2},\ldots
,a_{n}\right)  \in\mathbf{k}^{n}=V$, we have%
\begin{align*}
\pi^{\prime}\left(  v\right)   &  =\dfrac{1}{\left\vert G\right\vert }%
\sum_{g\in G}\underbrace{g^{-1}\pi\left(  gv\right)  }_{=\left(
a_{g^{-1}\left(  1\right)  },a_{g^{-1}\left(  1\right)  },\ldots
,a_{g^{-1}\left(  1\right)  }\right)  }\\
&  =\dfrac{1}{\left\vert G\right\vert }\sum_{g\in G}\left(  a_{g^{-1}\left(
1\right)  },a_{g^{-1}\left(  1\right)  },\ldots,a_{g^{-1}\left(  1\right)
}\right) \\
&  =\left(  \dfrac{1}{\left\vert G\right\vert }\sum_{g\in G}a_{g^{-1}\left(
1\right)  },\ \dfrac{1}{\left\vert G\right\vert }\sum_{g\in G}a_{g^{-1}\left(
1\right)  },\ \ldots,\ \dfrac{1}{\left\vert G\right\vert }\sum_{g\in
G}a_{g^{-1}\left(  1\right)  }\right)  .
\end{align*}
But the average $\dfrac{1}{\left\vert G\right\vert }\sum_{g\in G}%
a_{g^{-1}\left(  1\right)  }$ in this formula is simply the average $\dfrac
{1}{n}\sum_{i=1}^{n}a_{i}$ (because when $g$ ranges over all the $n!$
permutations in $S_{n}=G$, the number $g^{-1}\left(  1\right)  $ takes each of
the values $1,2,\ldots,n$ equally often). Thus, the above formula rewrites as%
\begin{equation}
\pi^{\prime}\left(  v\right)  =\left(  \dfrac{1}{n}\sum_{i=1}^{n}%
a_{i},\ \dfrac{1}{n}\sum_{i=1}^{n}a_{i},\ \ldots,\ \dfrac{1}{n}\sum_{i=1}%
^{n}a_{i}\right)  . \label{eq.exa.rep.G-rep.maschke.natD.5}%
\end{equation}
In other words, $\pi^{\prime}$ replaces each coordinate of the input vector by
the average of all coordinates.

The proof of Proposition \ref{prop.mod.diradd.proj} now says that
$U:=\operatorname*{Ker}\pi^{\prime}$ is a $\mathbf{k}\left[  G\right]
$-submodule (i.e., a $G$-subrepresentation) of $V$ satisfying $V=W\oplus U$.
As we can see from (\ref{eq.exa.rep.G-rep.maschke.natD.5}), this submodule
$U=\operatorname*{Ker}\pi^{\prime}$ is the set of all vectors whose
coordinates sum to $0$, that is, the zero-sum subrepresentation $R\left(
\mathbf{k}^{n}\right)  $ of $\mathbf{k}^{n}$.
\end{example}

Maschke's theorem can be generalized from the \textquotedblleft single group
$G$\textquotedblright\ setting to the \textquotedblleft group $G$ and subgroup
$H$\textquotedblright\ setting:

\begin{exercise}
\fbox{3} Let $G$ be a group, and $H$ a subgroup of $G$. Assume that the index
$\left[  G:H\right]  $ (which is defined as the size of the set $G/H$) is
finite. Let $V$ be a left $\mathbf{k}\left[  G\right]  $-module (i.e., a
representation of $G$ over $\mathbf{k}$). Let $W$ be a left $\mathbf{k}\left[
G\right]  $-submodule of $V$ (that is, a subrepresentation of $V$). Note that
every left $\mathbf{k}\left[  G\right]  $-module automatically becomes a left
$\mathbf{k}\left[  H\right]  $-module (by restriction of scalars, since the
inclusion map $H\rightarrow G$ gives rise to a $\mathbf{k}$-algebra morphism
$\mathbf{k}\left[  H\right]  \rightarrow\mathbf{k}\left[  G\right]  $). Assume that

\begin{enumerate}
\item the number $\left[  G:H\right]  $ is invertible in $\mathbf{k}$;

\item the $\mathbf{k}\left[  H\right]  $-submodule $W$ is a direct addend of
$V$ as a $\mathbf{k}\left[  H\right]  $-module.
\end{enumerate}

Prove that the left $\mathbf{k}\left[  G\right]  $-submodule $W$ is also a
direct addend of $V$ as a left $\mathbf{k}\left[  G\right]  $-module. \medskip

[\textbf{Hint:} Cosets.]
\end{exercise}

Maschke's theorem also has the following consequence:

\begin{corollary}
\label{cor.rep.G-rep.maschke-quot}Let $G$, $V$ and $W$ be as in Theorem
\ref{thm.rep.G-rep.maschke}. Then,%
\[
V\cong W\oplus\left(  V/W\right)  \ \ \ \ \ \ \ \ \ \ \text{as left
}\mathbf{k}\left[  G\right]  \text{-modules.}%
\]
(Here, $W\oplus\left(  V/W\right)  $ is of course an external direct sum,
since $V/W$ is not a subset of $V$.)
\end{corollary}

\begin{proof}
Theorem \ref{thm.rep.G-rep.maschke} yields that the left $\mathbf{k}\left[
G\right]  $-submodule $W$ is a direct addend of $V$ as a left $\mathbf{k}%
\left[  G\right]  $-module. Hence, Proposition \ref{prop.mod.diradd.V/W}
(applied to $R=\mathbf{k}\left[  G\right]  $) yields that
\[
V\cong W\oplus\left(  V/W\right)  \ \ \ \ \ \ \ \ \ \ \text{as left
}\mathbf{k}\left[  G\right]  \text{-modules.}%
\]
This proves Corollary \ref{cor.rep.G-rep.maschke-quot}.
\end{proof}

\begin{exercise}
Let $n$ be a positive integer. Let $G=C_{n}$ (the cyclic group with $n$
elements), and let $\mathbf{k}=\mathbb{Z}/n$. Let $V$ be the representation
$\mathbf{k}^{2}$ of $G$, where the generator $g\in C_{n}$ acts as the
$\mathbf{k}$-linear map
\begin{align*}
w:\mathbf{k}^{2}  &  \rightarrow\mathbf{k}^{2},\\
\left(  x,y\right)   &  \mapsto\left(  x,\ x+y\right)  .
\end{align*}

\textbf{(a)} \fbox{1} Prove that such a representation really exists (i.e.,
that the above $\mathbf{k}$-linear map $w$ satisfies $w^{n}=\operatorname*{id}%
$, so that Example \ref{exa.rep.G-rep.lin.Cn} constructs a representation of
$C_{n}$ from it). \medskip

\textbf{(b)} \fbox{2} Assume that $n$ is prime (so that $\mathbf{k}%
=\mathbb{Z}/n$ is a finite field). Prove that $\operatorname*{span}%
\nolimits_{\mathbf{k}}\left\{  \left(  0,1\right)  \right\}  $ is a
subrepresentation of $V$, but is not a direct addend. More generally, prove
that $V$ cannot be written as a direct sum of two nonzero subrepresentations.
(Unsurprisingly, Theorem \ref{thm.rep.G-rep.maschke} does not apply here,
since $\left\vert G\right\vert =\left\vert C_{n}\right\vert =n$ is not
invertible in $\mathbf{k}$.)
\end{exercise}

\subsubsection{\label{subsec.rep.maschke.filt}Filtrations}

Throughout mathematics, a nested chain of subobjects of a given object
$\mathcal{O}$ is called a \emph{filtration} of $\mathcal{O}$ if it starts with
the \textquotedblleft trivial\textquotedblright\ subobject (whatever this
means in the context) and ends with $\mathcal{O}$ itself. Let us be more
specific and define this notion for left $R$-modules in particular:

\begin{definition}
\label{def.mod.filtr}Let $R$ be a ring. Let $V$ be a left $R$-module. \medskip

\textbf{(a)} A \emph{filtration} (aka \emph{complete filtration}) of $V$ means
a tuple $\left(  W_{0},W_{1},\ldots,W_{k}\right)  $ of left $R$-submodules of
$V$ satisfying%
\[
0=W_{0}\subseteq W_{1}\subseteq W_{2}\subseteq\cdots\subseteq W_{k}=V,
\]
where the $0$ at the front means the $R$-submodule $\left\{  0\right\}  $.
Note that we allow the \textquotedblleft$\subseteq$\textquotedblright\ signs
to be equalities.

When the ring $R$ might not be clear from the context, we can disambiguate it
by saying \textquotedblleft\emph{left }$R$\emph{-module filtration}%
\textquotedblright\ or \textquotedblleft\emph{filtration by left }%
$R$\emph{-submodules}\textquotedblright\ instead of just \textquotedblleft
filtration\textquotedblright. \medskip

\textbf{(b)} If we drop the $W_{k}=V$ requirement in Definition
\ref{def.mod.filtr} \textbf{(a)}, then $\left(  W_{0},W_{1},\ldots
,W_{k}\right)  $ is called an \emph{incomplete filtration} of $V$. \medskip

\textbf{(c)} If $\left(  W_{0},W_{1},\ldots,W_{k}\right)  $ is a filtration of
$V$, then the quotient $R$-modules $W_{i}/W_{i-1}$ for all $i\in\left[
k\right]  $ are called the \emph{subquotients} of this filtration. \medskip

\textbf{(d)} If $\left(  W_{0},W_{1},\ldots,W_{k}\right)  $ is a filtration of
$V$, then the number $k$ is called the \emph{length} of this filtration.
\end{definition}

Some other variants of filtrations can be found in the literature: e.g.,
decreasing filtrations (where the \textquotedblleft$\subseteq$%
\textquotedblright\ signs are replaced by \textquotedblleft$\supseteq
$\textquotedblright\ signs), infinite filtrations (where the tuple is replaced
by an infinite sequence), poset filtrations (where the tuple is replaced by a
partially ordered family) and others.

Intuitively, a filtration of $V$ can be regarded as a \textquotedblleft
staircase\textquotedblright\ going from $0$ to $V$, with each stair being a
left $R$-submodule of $V$. Its subquotients are the \textquotedblleft
steps\textquotedblright\ from a stair to the next.

\begin{example}
Proposition \ref{prop.rep.G-rep.Sn-nat.subreps.dirsum} \textbf{(b)} shows that
if $n=0$ in $\mathbf{k}$, then $\left(  0,\ D\left(  \mathbf{k}^{n}\right)
,\ R\left(  \mathbf{k}^{n}\right)  ,\ \mathbf{k}^{n}\right)  $ is a filtration
of the left $\mathbf{k}\left[  S_{n}\right]  $-module $\mathbf{k}^{n}$.
\end{example}

\begin{example}
Let $R$ be a ring. If $A$ and $B$ are two $R$-submodules of a left $R$-module
$V$, then $\left(  0,\ A\cap B,\ A,\ A+B,\ V\right)  $ is a filtration of $V$.
\end{example}

Maschke's theorem (Corollary \ref{cor.rep.G-rep.maschke-quot} specifically) yields:

\begin{corollary}
\label{cor.rep.G-rep.maschke-filt}Let $\mathbf{k}$ be a field of
characteristic $0$. Let $G$ be a finite group. Let $V$ be a left
$\mathbf{k}\left[  G\right]  $-module. If $\left(  W_{0},W_{1},\ldots
,W_{k}\right)  $ is a filtration of $V$, then%
\[
V\cong\left(  W_{1}/W_{0}\right)  \oplus\left(  W_{2}/W_{1}\right)
\oplus\cdots\oplus\left(  W_{k}/W_{k-1}\right)  .
\]

\end{corollary}

\begin{proof}
[Proof of Corollary \ref{cor.rep.G-rep.maschke-filt} (sketched).]Induct on $k$.

The \textit{base case} $k=0$ is obvious (here, $V$ must be $\left\{
0\right\}  $, since the definition of a filtration yields $0=W_{0}=W_{k}=V$).

For the \textit{induction step} (from $k-1$ to $k$), we consider a filtration
$\left(  W_{0},W_{1},\ldots,W_{k}\right)  $ of $V$. Then, $0=W_{0}\subseteq
W_{1}\subseteq W_{2}\subseteq\cdots\subseteq W_{k-1}\subseteq W_{k}=V$. Now,
$W_{k-1}$ is a left $\mathbf{k}\left[  G\right]  $-submodule of $V$, and thus
we can apply Corollary \ref{cor.rep.G-rep.maschke-quot} to $W=W_{k-1}$
(indeed, both assumptions of Theorem \ref{thm.rep.G-rep.maschke} are satisfied
by Remark \ref{rmk.rep.G-rep.maschke-k}). Thus, we obtain%
\begin{equation}
V\cong W_{k-1}\oplus\left(  V/W_{k-1}\right)  =W_{k-1}\oplus\left(
W_{k}/W_{k-1}\right)  \label{pf.cor.rep.G-rep.maschke-filt.2}%
\end{equation}
(since $V=W_{k}$). But the tuple $\left(  W_{0},W_{1},\ldots,W_{k-1}\right)  $
is a filtration of $W_{k-1}$ (since $0=W_{0}\subseteq W_{1}\subseteq
\cdots\subseteq W_{k-1}=W_{k-1}$), and thus our induction hypothesis (applied
to $W_{k-1}$ and $\left(  W_{0},W_{1},\ldots,W_{k-1}\right)  $ instead of $V$
and $\left(  W_{0},W_{1},\ldots,W_{k}\right)  $) yields that%
\[
W_{k-1}\cong\left(  W_{1}/W_{0}\right)  \oplus\left(  W_{2}/W_{1}\right)
\oplus\cdots\oplus\left(  W_{k-1}/W_{k-2}\right)  .
\]
Thus, we can rewrite (\ref{pf.cor.rep.G-rep.maschke-filt.2}) as%
\begin{align*}
V  &  \cong\left(  \left(  W_{1}/W_{0}\right)  \oplus\left(  W_{2}%
/W_{1}\right)  \oplus\cdots\oplus\left(  W_{k-1}/W_{k-2}\right)  \right)
\oplus\left(  W_{k}/W_{k-1}\right) \\
&  \cong\left(  W_{1}/W_{0}\right)  \oplus\left(  W_{2}/W_{1}\right)
\oplus\cdots\oplus\left(  W_{k}/W_{k-1}\right)  .
\end{align*}
This proves Corollary \ref{cor.rep.G-rep.maschke-filt} for our $k$, and thus
completes the induction step.
\end{proof}

Recall that representations of a group $G$ over $\mathbf{k}$ are just left
$\mathbf{k}\left[  G\right]  $-modules. Thus, the concept of a filtration can
be transferred automatically from left $\mathbf{k}\left[  G\right]  $-modules
to representations of $G$.

If $\mathbf{k}$ is a field of characteristic $0$, and if $G$ is a finite
group, then Corollary \ref{cor.rep.G-rep.maschke-filt} shows that any
filtration of a representation of $G$ can be used to decompose this
representation into a direct sum. But in general (when $\mathbf{k}$ is not a
field of characteristic $0$), filtrations are more frequent (and more easily
found) than direct sum decompositions.

What is the longest filtration that we can find for a given left $R$-module
$V$ ? Taken literally, this question is stupid (we can take an arbitrarily
long tuple of the form $\left(  0,V,V,\ldots,V\right)  $), but we can make it
more sensible by restricting ourselves to \textquotedblleft
nondegenerate\textquotedblright\ filtrations $\left(  W_{0},W_{1},\ldots
,W_{k}\right)  $ (that is, ones in which each $W_{i-1}$ is a \textbf{proper}
subset of $W_{i}$). Even then, the question does not always have an answer,
since (e.g.) the $\mathbb{Z}$-module $\mathbb{Z}$ has arbitrarily long
nondegenerate filtrations:
\[
\left(  0\mathbb{Z},\ 2^{k}\mathbb{Z},\ 2^{k-1}\mathbb{Z},\ \ldots
,\ 2^{0}\mathbb{Z}\right)  \text{ is a filtration of }\mathbb{Z}\text{ for
each }k\in\mathbb{N}.
\]
However, if $\mathbf{k}$ is a field and $R$ is a $\mathbf{k}$-algebra, and if
$V$ is a left $R$-module that is finite-dimensional as a $\mathbf{k}$-vector
space, then a nondegenerate filtration of $V$ will always have length
$\leq\dim V$, where $\dim V$ is the dimension of $V$ as a $\mathbf{k}$-vector
space (why?). Thus, longest nondegenerate filtrations exist in this case. Such
filtrations are known as \emph{composition series} of $V$, and can be
equivalently described as those filtrations whose subquotients are irreducible
$R$-modules (why?). A classical result known as the \emph{Jordan--H\"{o}lder
theorem} (see, e.g., \cite[Theorem 3.7.1]{EGHetc11}) says that the
subquotients of such filtrations are uniquely determined (up to reordering and isomorphism):

\begin{theorem}
[Jordan--H\"{o}lder theorem]\label{thm.rep.G-rep.JH}Let $\mathbf{k}$ be a
field. Let $R$ be a $\mathbf{k}$-algebra. Let $V$ be a left $R$-module that is
finite-dimensional as a $\mathbf{k}$-vector space. Then: \medskip

\textbf{(a)} There exists a filtration $\left(  W_{0},W_{1},\ldots
,W_{k}\right)  $ of $V$ such that all its subquotients $W_{i}/W_{i-1}$ are
irreducible left $R$-modules. Such a filtration is called a \emph{composition
series} (or \emph{Jordan--H\"{o}lder series}) of $V$. \medskip

\textbf{(b)} Any two composition series of $V$ have the same subquotients (up
to reordering and isomorphism).
\end{theorem}

\begin{proof}
[Proof idea.]\textbf{(a)} As we said above, any nondegenerate filtration of
$V$ must have length $\leq\dim V$, so that there exists a longest
nondegenerate filtration of $V$, and any such longest filtration must be a
composition series.

Alternatively, we can construct a composition series by hand: Start by setting
$W_{0}=0$; then let $W_{1}$ be a minimum-dimensional $R$-submodule of $V$ that
properly contains $W_{0}$; then let $W_{2}$ be a minimum-dimensional
$R$-submodule of $V$ that properly contains $W_{1}$; and so on. The sequence
must come to an end, since $\dim W_{0}<\dim W_{1}<\dim W_{2}<\cdots\leq\dim
V$. \medskip

\textbf{(b)} This is not that easy. See \cite[Theorem 3.7.1]{EGHetc11} or any
text on representation theory. (There is also an analogous \textquotedblleft
Jordan--H\"{o}lder theorem\textquotedblright\ in group theory -- see
\cite[\textquotedblleft Subgroup Series I\textquotedblright, Theorem
2.8]{Conrad} -- whose proof is rather analogous to the proof of Theorem
\ref{thm.rep.G-rep.JH} \textbf{(b)}.)
\end{proof}

Composition series can be useful in some obvious and some less obvious ways:

\begin{itemize}
\item Their most obvious use is to help break an $R$-module $V$ into
irreducible $R$-modules (the subquotients of a composition series). This works
particularly well when $R$ is a group algebra over a field of characteristic
$0$, since in this case $V$ is really decomposed into a direct sum (by
Corollary \ref{cor.rep.G-rep.maschke-filt}). In other cases, the
\textquotedblleft pieces\textquotedblright\ might not determine $V$
completely, but they still encode a lot of information about $V$.

\item A less obvious use: They impose a \textquotedblleft
block-triangularity\textquotedblright\ on the action of any element of
$\mathbf{k}\left[  G\right]  $. Indeed, if $R$ is a $\mathbf{k}$-algebra, and
if $\left(  W_{0},W_{1},\ldots,W_{k}\right)  $ is a filtration of a left
$R$-module $V$, then the action $\tau_{r}$ of an element $r\in R$ on $V$ is
not just a $\mathbf{k}$-module endomorphism of $V$, but a $\mathbf{k}$-module
endomorphism that preserves each of the submodules $W_{0},W_{1},\ldots,W_{k}$.
Translated into the language of matrices (where we assume that $\mathbf{k}$ is
a field and $V$ is a finite-dimensional $\mathbf{k}$-vector space), such an
endomorphism becomes a block-upper-triangular matrix in an appropriate basis
of $V$. We will see some instances of this in Section \ref{sec.bas.mur}.

\item Even less obviously: They can be used to find all the irreps (i.e.,
irreducible representations) of a finite group $G$, or more generally, all
irreducible modules over a $\mathbf{k}$-algebra $R$. This relies on the
following theorem:
\end{itemize}

\begin{theorem}
\label{thm.mod.filt.irrep-in-filt}Let $\mathbf{k}$ be a field. Let $R$ be a
$\mathbf{k}$-algebra that is finite-dimensional as a $\mathbf{k}$-vector
space. Let $\left(  W_{0},W_{1},\ldots,W_{k}\right)  $ be a composition series
of the left regular $R$-module $R$. Then, each irreducible left $R$-module is
isomorphic to $W_{i}/W_{i-1}$ for some $i\in\left[  k\right]  $. (Note that
this $i$ does not have to be unique.) In other words, each irreducible left
$R$-module appears (up to isomorphism) as a subquotient of our composition series.
\end{theorem}

\begin{fineprint}
\begin{proof}
[Proof of Theorem \ref{thm.mod.filt.irrep-in-filt} (sketched).]Let $J$ be an
irreducible left $R$-module. We must show that $J\cong W_{i}/W_{i-1}$ for some
$i\in\left[  k\right]  $.

The $R$-module $J$ is nonzero (since $J$ is irreducible). Hence, there exists
a nonzero vector $v\in J$. Consider this $v$.

For each left $R$-submodule $M$ of $R$, we let%
\[
Mv:=\left\{  mv\ \mid\ m\in M\right\}  .
\]
This $Mv$ is a left $R$-submodule of $J$, since $M$ is a left $R$-submodule of
$R$ (thus fixed under addition and left scaling by $R$).

The tuple $\left(  W_{0},W_{1},\ldots,W_{k}\right)  $ is a composition series
and thus a filtration of the left $R$-module $R$. Hence, $0=W_{0}\subseteq
W_{1}\subseteq W_{2}\subseteq\cdots\subseteq W_{k}=R$. In particular,
$W_{k}=R$. But $1\in R$ and thus $1v\in Rv$, so that $v=1v\in Rv=W_{k}v$
(since $R=W_{k}$). Hence, $W_{k}v\neq0$ (since $v$ is nonzero). On the other
hand, from $W_{0}=0$, we obtain $W_{0}v=0v=0$.

We have $W_{k}v\neq0$. Hence, there exists some $i\in\left\{  0,1,\ldots
,k\right\}  $ such that $W_{i}v\neq0$ (namely, $i=k$). Pick the
\textbf{smallest} such $i$. Thus, $i\neq0$ (since $W_{0}v=0$), so that
$i-1\in\left\{  0,1,\ldots,k\right\}  $ as well. Moreover, $W_{i-1}v=0$ (since
$i$ is the \textbf{smallest} number to satisfy $W_{i}v\neq0$).

Consider the map%
\begin{align*}
f:W_{i}/W_{i-1}  &  \rightarrow J,\\
\overline{r}  &  \mapsto rv,
\end{align*}
where $\overline{r}$ denotes the residue class of a vector $r\in W_{i}$ modulo
$W_{i-1}$. This map $f$ is well-defined, since two vectors $r,s\in W_{i}$ with
equal residue classes $\overline{r}=\overline{s}\in W_{i}/W_{i-1}$ will
necessarily satisfy $rv=sv$ (because $\overline{r}=\overline{s}$ yields
$r-s\in W_{i-1}$ and thus $\left(  r-s\right)  v\in W_{i-1}v=0$, so that
$\left(  r-s\right)  v=0$, so that $0=\left(  r-s\right)  v=rv-sv$ and
therefore $rv=sv$). Moreover, the map $f$ is left $R$-linear (by a trivial computation).

Therefore, the image $f\left(  W_{i}/W_{i-1}\right)  $ of $f$ is a left
$R$-submodule of $J$. Since $J$ is irreducible, this image must therefore be
either $J$ or $0$. Since this image cannot be $0$ (because the definition of
$f$ yields that $f\left(  W_{i}/W_{i-1}\right)  =\left\{  rv\ \mid\ r\in
W_{i}\right\}  =W_{i}v\neq0$), we thus conclude that it is $J$. In other
words, the map $f$ is surjective.

The kernel $\operatorname*{Ker}f$ of $f$ is a left $R$-submodule of
$W_{i}/W_{i-1}$ (since $f$ is left $R$-linear). But the left $R$-module
$W_{i}/W_{i-1}$ is irreducible (since $\left(  W_{0},W_{1},\ldots
,W_{k}\right)  $ is a composition series). Hence, its only $R$-submodules are
$W_{i}/W_{i-1}$ and $0$. Consequently, $\operatorname*{Ker}f$ (being its left
$R$-submodule) must be either $W_{i}/W_{i-1}$ or $0$. Since
$\operatorname*{Ker}f$ cannot be $W_{i}/W_{i-1}$ (because this would entail
$f=0$, and thus the image of $f$ would be $0$, but we have already shown that
this image cannot be $0$), we thus conclude that $\operatorname*{Ker}f$ must
be $0$. In other words, the map $f$ is injective (since $f$ is $R$-linear).

Now we know that the map $f$ is injective and surjective, hence bijective.
Since $f$ is furthermore left $R$-linear, we thus conclude that $f$ is an
$R$-module isomorphism. Therefore, $J\cong W_{i}/W_{i-1}$ as left $R$-modules.

Thus we have found an $i\in\left[  k\right]  $ such that $J\cong W_{i}%
/W_{i-1}$. Hence, such an $i$ exists. This proves Theorem
\ref{thm.mod.filt.irrep-in-filt}.
\end{proof}
\end{fineprint}

\begin{noncompile}
The above proof is essentially taken from
https://mathoverflow.net/questions/14514/how-big-can-the-irreps-of-a-finite-group-be-over-an-arbitrary-field/14516\#14516
second proof of Theorem 1. (Note that I call left $R$-modules
\textquotedblleft representations of $R$\textquotedblright\ there, which is a
widespread terminology in abstract algebra.)
\end{noncompile}

Since we are talking about filtrations of the left regular $R$-module $R$, it
is worth mentioning that the left $R$-submodules of this $R$-module are better
known under a different name:

\begin{remark}
\label{rmk.mod.left-R-sub-R}Let $R$ be a ring. A left $R$-submodule of the
left regular $R$-module $R$ is the same thing as a left ideal of $R$, that is,
a subgroup $J$ of the additive group $\left(  R,+,0\right)  $ with the
property that%
\[
\text{all }r\in R\text{ and }j\in J\text{ satisfy }rj\in J
\]
(but not necessarily $jr\in J$).
\end{remark}

\begin{proof}
Just compare the definitions of \textquotedblleft left ideal\textquotedblright%
\ and \textquotedblleft left $R$-submodule\textquotedblright. (Note that
existence of negative inverses follows from scaling by $-1_{R}$.)
\end{proof}

We note the following simple consequence of Theorem
\ref{thm.mod.filt.irrep-in-filt}:

\begin{corollary}
\label{cor.mod.irrep.fin}Let $\mathbf{k}$ be a field. Let $R$ be a
$\mathbf{k}$-algebra that is finite-dimensional as a $\mathbf{k}$-vector
space. Then, there are only finitely many irreducible left $R$-modules (up to isomorphism).
\end{corollary}

In particular, if $\mathbf{k}$ is a field and $G$ is a finite group, then
there are (up to isomorphism) only finitely many irreducible left
$\mathbf{k}\left[  G\right]  $-modules, i.e., only finitely many
irreps\footnote{Recall that \textquotedblleft irrep\textquotedblright\ is
short for \textquotedblleft irreducible representation\textquotedblright.} of
$G$. If we want to find these irreps, then (by Theorem
\ref{thm.mod.filt.irrep-in-filt}) we only need to look at a composition series
of its left regular representation $\mathbf{k}\left[  G\right]  $. The
subquotients of this series will then be the irreps (though some might appear
multiple times). If $\mathbf{k}$ is a field of characteristic $0$, then we
furthermore conclude (by Corollary \ref{cor.rep.G-rep.maschke-filt}) that the
left regular representation $\mathbf{k}\left[  G\right]  $ can be decomposed
as a direct sum of these irreps (with some irreps appearing multiple times as
addends), and therefore we see that every irrep of $G$ is isomorphic to a
subrepresentation of $\mathbf{k}\left[  G\right]  $. The latter part of this
claim is true even when $\mathbf{k}$ does not have characteristic $0$:

\begin{theorem}
\label{thm.mod.irrep.sub-kG}Let $\mathbf{k}$ be any field, and let $G$ be a
finite group. Then, any irreducible representation of $G$ is isomorphic to a
subrepresentation of the left regular representation $\mathbf{k}\left[
G\right]  $.
\end{theorem}

\begin{fineprint}
\begin{proof}
See \url{https://mathoverflow.net/questions/72844} or Exercise
\ref{exe.mod.irrep.sub-kG} \textbf{(c)}.
\end{proof}
\end{fineprint}

\begin{exercise}
\label{exe.mod.irrep.sub-kG}Let $\mathbf{k}$ be any field, and let $G$ be a
finite group. Let $V$ be an irreducible representation of $G$. Prove the
following: \medskip

\textbf{(a)} \fbox{2} If $v\in V$ is a nonzero vector, then
$V=\operatorname*{span}\nolimits_{\mathbf{k}}\left\{  gv\ \mid\ g\in
G\right\}  $. \medskip

\textbf{(b)} \fbox{2} If $f:V\rightarrow\mathbf{k}$ is a nonzero $\mathbf{k}%
$-linear map, and if we define a map $\alpha:V\rightarrow\mathbf{k}\left[
G\right]  $ by setting%
\[
\alpha\left(  v\right)  =\sum_{g\in G}f\left(  gv\right)  g^{-1}%
\ \ \ \ \ \ \ \ \ \ \text{for all }v\in V,
\]
then $\alpha$ is an injective morphism of $G$-representations (from $V$ to the
left regular representation of $G$). \medskip

\textbf{(c)} \fbox{1} The representation $V$ is isomorphic to a
subrepresentation of the left regular representation $\mathbf{k}\left[
G\right]  $.
\end{exercise}

\begin{exercise}
\fbox{2} Prove the following variant of Theorem
\ref{thm.mod.filt.irrep-in-filt}:

Let $\mathbf{k}$ be a field. Let $R$ be a $\mathbf{k}$-algebra. Let $V$ be a
left $R$-module that is finite-dimensional as a $\mathbf{k}$-vector space. Let
$\left(  W_{0},W_{1},\ldots,W_{k}\right)  $ be a composition series of the
left $R$-module $V$. Then, each irreducible left $R$-submodule of $V$ is
isomorphic to $W_{i}/W_{i-1}$ for some $i\in\left[  k\right]  $.
\end{exercise}

\subsubsection{\label{subsec.rep.maschke.kS3}An example: $\mathbf{k}\left[
S_{3}\right]  $}

Let us now get our hands dirty and find a composition series in the wild.

We let $\mathbf{k}$ be a field of characteristic $\neq3$ (from here on until
the end of Subsection \ref{subsec.rep.maschke.kS3}). We want to find a
composition series of the left regular representation of $S_{3}$.

The left regular representation of $S_{3}$ is the left regular $\mathbf{k}%
\left[  S_{3}\right]  $-module $\mathbf{k}\left[  S_{3}\right]  $. As a
$\mathbf{k}$-vector space, it is $6$-dimensional, with basis $\left(
1,s_{1},s_{2},s_{1}s_{2},s_{2}s_{1},s_{1}s_{2}s_{1}\right)  $ (or any other
way to list the six permutations in $S_{3}$). We define the following left
$\mathbf{k}\left[  S_{3}\right]  $-submodules of the left regular
representation $\mathbf{k}\left[  S_{3}\right]  $:%
\begin{align*}
\left\langle \nabla\right\rangle  &  :=\operatorname*{span}%
\nolimits_{\mathbf{k}}\left\{  \nabla\right\}  \ \ \ \ \ \ \ \ \ \ \left(
\begin{array}
[c]{c}%
\text{this is a left }\mathbf{k}\left[  S_{3}\right]  \text{-submodule,}\\
\text{by Exercise \ref{exe.eps.Nab}}%
\end{array}
\right)  ;\\
A_{1}  &  :=\left\{  \mathbf{a}\in\mathbf{k}\left[  S_{3}\right]
\ \mid\ \mathbf{a}s_{1}=\mathbf{a}\right\}  \ \ \ \ \ \ \ \ \ \ \left(
\begin{array}
[c]{c}%
\text{this is a left }\mathbf{k}\left[  S_{3}\right]  \text{-submodule,}\\
\text{since }\mathbf{a}s_{1}=\mathbf{a}\text{ implies }\mathbf{ba}%
s_{1}=\mathbf{ba}\\
\text{for any }\mathbf{b}\in\mathbf{k}\left[  S_{3}\right]
\end{array}
\right)  ;\\
A_{2}  &  :=\left\{  \mathbf{a}\in\mathbf{k}\left[  S_{3}\right]
\ \mid\ \mathbf{a}s_{2}=\mathbf{a}\right\}  \ \ \ \ \ \ \ \ \ \ \left(
\begin{array}
[c]{c}%
\text{this is a left }\mathbf{k}\left[  S_{3}\right]  \text{-submodule,}\\
\text{for similar reasons as for }A_{1}%
\end{array}
\right)  .
\end{align*}
It is easy to see that $A_{1}\cap A_{2}=\left\langle \nabla\right\rangle $ (in
fact, this is just Exercise \ref{exe.integral.fix-converse} \textbf{(b)} for
$n=3$). Thus, the tuple%
\begin{equation}
\left(  0,\ \left\langle \nabla\right\rangle ,\ A_{1},\ A_{1}+A_{2}%
,\ \mathbf{k}\left[  S_{3}\right]  \right)
\label{eq.subsec.rep.maschke.kS3.filt}%
\end{equation}
(where $0$ denotes $\left\{  0\right\}  $ as usual in linear algebra) is a
filtration of the left $\mathbf{k}\left[  S_{3}\right]  $-module
$\mathbf{k}\left[  S_{3}\right]  $. Now we shall show that it is a composition
series of $\mathbf{k}\left[  S_{3}\right]  $.

The easiest way to prove this is to find bases of all the submodules and the
subquotients. This is best achieved by finding a basis of $\mathbf{k}\left[
S_{3}\right]  $ that is \textquotedblleft adapted\textquotedblright\ to this
filtration, meaning that each of the submodules $0,\ \left\langle
\nabla\right\rangle ,\ A_{1},\ A_{1}+A_{2},\ \mathbf{k}\left[  S_{3}\right]  $
in the filtration is spanned by some part of this basis. Indeed, such a basis
will automatically provide bases for all subquotients, because of the
following simple fact from linear algebra:

\begin{lemma}
\label{lem.mod.quot.basis}Let $\mathbf{k}$ be any commutative ring. Let $V$ be
a $\mathbf{k}$-module, and let $W$ be a $\mathbf{k}$-submodule of $V$. For
each $v\in V$, we let $\overline{v}$ denote the residue class of $v$ in $V/W$
(that is, the coset $v+W$).

Assume that we have a basis $\left(  b_{1},b_{2},\ldots,b_{k}\right)  $ of $V$
and a number $i\in\left\{  0,1,\ldots,k\right\}  $ such that the $i$-tuple
$\left(  b_{1},b_{2},\ldots,b_{i}\right)  $ is a basis of $W$. Then, the
$\left(  k-i\right)  $-tuple $\left(  \overline{b_{i+1}},\ \overline{b_{i+2}%
},\ \ldots,\ \overline{b_{k}}\right)  $ is a basis of the quotient
$\mathbf{k}$-module $V/W$.
\end{lemma}

Of course, a basis $\left(  b_{1},b_{2},\ldots,b_{k}\right)  $ and a number
$i$ as in Lemma \ref{lem.mod.quot.basis} do not always exist for general
$\mathbf{k}$ and $V$. But we are dealing with a rather nice situation, in
which $\mathbf{k}$ is a field and $V$ is a finite-dimensional vector space. In
this situation, $\left(  b_{1},b_{2},\ldots,b_{k}\right)  $ and $i$ always
exist, since we can choose any basis $\left(  b_{1},b_{2},\ldots,b_{i}\right)
$ of $W$ and extend it to a basis $\left(  b_{1},b_{2},\ldots,b_{k}\right)  $
of $V$.

This is also how we can find a basis of $\mathbf{k}\left[  S_{3}\right]  $
that is \textquotedblleft adapted\textquotedblright\ to our filtration
(\ref{eq.subsec.rep.maschke.kS3.filt}): We start with the (empty) basis of the
zero subspace $0$; then we extend it to a basis of $\left\langle
\nabla\right\rangle $; then we extend the latter to a basis of $A_{1}$; and so
on. There are many possible options for these steps, but let me choose a
particularly nice one.

To define this basis, recall the elements%
\[
\nabla_{b,a}:=\sum_{\substack{w\in S_{n};\\w\left(  a\right)  =b}%
}w\ \ \ \ \ \ \ \ \ \ \text{for all }a,b\in\left[  n\right]
\]
of $\mathbf{k}\left[  S_{n}\right]  $ (these were defined in Exercise
\ref{exe.Nabab.1}). Now, I claim that the six vectors%
\begin{align*}
\nabla &  =1+s_{1}+s_{2}+s_{1}s_{2}+s_{2}s_{1}+s_{1}s_{2}s_{1},\\
\nabla_{3,3}  &  =1+s_{1},\\
\nabla_{2,3}  &  =s_{2}\left(  1+s_{1}\right)  =s_{2}+s_{2}s_{1},\\
\nabla_{3,2}  &  =\left(  1+s_{1}\right)  s_{2}=s_{2}+s_{1}s_{2},\\
\nabla_{2,2}  &  =1+t_{1,3}=s_{2}\left(  1+s_{1}\right)  s_{2}=1+s_{2}%
s_{1}s_{2}=1+s_{1}s_{2}s_{1},\\
1  &  =1
\end{align*}
form a basis of $\mathbf{k}\left[  S_{3}\right]  $ that is \textquotedblleft
adapted\textquotedblright\ to our filtration
(\ref{eq.subsec.rep.maschke.kS3.filt}). Indeed, it is not hard to see that%
\begin{align}
0  &  =\operatorname*{span}\nolimits_{\mathbf{k}}\left\{  {}\right\}
;\label{eq.subsec.rep.maschke.kS3.filt0}\\
\left\langle \nabla\right\rangle  &  =\operatorname*{span}%
\nolimits_{\mathbf{k}}\left\{  \nabla\right\}
;\label{eq.subsec.rep.maschke.kS3.filt1}\\
A_{1}  &  =\operatorname*{span}\nolimits_{\mathbf{k}}\left\{  \nabla
,\ \nabla_{3,3},\ \nabla_{2,3}\right\}
;\label{eq.subsec.rep.maschke.kS3.filt2}\\
A_{1}+A_{2}  &  =\operatorname*{span}\nolimits_{\mathbf{k}}\left\{
\nabla,\ \nabla_{3,3},\ \nabla_{2,3},\ \nabla_{3,2},\ \nabla_{2,2}\right\}
,\label{eq.subsec.rep.maschke.kS3.filt3}\\
\mathbf{k}\left[  S_{3}\right]   &  =\operatorname*{span}\nolimits_{\mathbf{k}%
}\left\{  \nabla,\ \nabla_{3,3},\ \nabla_{2,3},\ \nabla_{3,2},\ \nabla
_{2,2},\ 1\right\}  . \label{eq.subsec.rep.maschke.kS3.filt4}%
\end{align}
(The first two of these equalities here are obvious. The third is easiest to
prove in the form $A_{1}=\operatorname*{span}\nolimits_{\mathbf{k}}\left\{
\nabla_{3,3},\ \nabla_{2,3},\ \nabla_{1,3}\right\}  $, which is equivalent
because $\nabla=\nabla_{3,3}+\nabla_{2,3}+\nabla_{1,3}$. This can be shown by
writing an arbitrary element $\mathbf{a}\in\mathbf{k}\left[  S_{3}\right]  $
as $\mathbf{a}=\alpha+\beta s_{1}+\gamma s_{2}+\delta t_{1,3}+\varphi
s_{1}s_{2}+\zeta s_{2}s_{1}$ for $\alpha,\beta,\gamma,\delta,\varphi,\zeta
\in\mathbf{k}$, and observing that the equality $\mathbf{a}s_{1}=\mathbf{a}$
is equivalent to $\left(  \alpha=\beta\text{ and }\gamma=\zeta\text{ and
}\delta=\varphi\right)  $, which allows us to write $\mathbf{a}$ as
$\mathbf{a}=\alpha\left(  1+s_{1}\right)  +\gamma s_{2}\left(  1+s_{1}\right)
+\delta t_{1,3}\left(  1+s_{1}\right)  =\alpha\nabla_{3,3}+\gamma\nabla
_{2,3}+\delta\nabla_{1,3}$. The fourth equation follows by combining the third
with the analogous equation $A_{2}=\operatorname*{span}\nolimits_{\mathbf{k}%
}\left\{  \nabla,\ \nabla_{3,1},\ \nabla_{2,1}\right\}  $ and observing that
$\nabla_{3,1}=\nabla-\nabla_{3,3}-\nabla_{3,2}$ and $\nabla_{2,1}%
=\nabla-\nabla_{2,3}-\nabla_{2,2}$. The least obvious of these equalities is
perhaps the fifth and last one, but even that is easy\footnote{In order to see
that our list $\left(  \nabla,\ \nabla_{3,3},\ \nabla_{2,3},\ \nabla
_{3,2},\ \nabla_{2,2},\ 1\right)  $ spans $\mathbf{k}\left[  S_{3}\right]  $,
we just observe that%
\begin{align*}
1  &  =1,\ \ \ \ \ \ \ \ \ \ s_{1}=\nabla_{3,3}-1,\ \ \ \ \ \ \ \ \ \ s_{2}%
=-\nabla+\nabla_{3,3}+\nabla_{2,3}+\nabla_{3,2}+\nabla_{2,2}-1,\\
s_{1}s_{2}  &  =\nabla-\nabla_{2,2}-\nabla_{3,3}-\nabla_{2,3}%
+1,\ \ \ \ \ \ \ \ \ \ s_{2}s_{1}=\nabla-\nabla_{2,2}-\nabla_{3,3}%
-\nabla_{3,2}+1,\\
s_{1}s_{2}s_{1}  &  =\nabla_{2,2}-1.
\end{align*}
Its linear independence can be proved by solving a system of linear equations,
or by recalling the fact (see, e.g., \cite[Exercise 2.5.18 \textbf{(a)}%
]{GriRei14}) that any list of $d$ vectors that spans a free $\mathbf{k}%
$-module of rank $d$ must be a basis and therefore linearly independent.}.)
This basis is an instance of what we will later call the \emph{(row) Murphy
basis} of $\mathbf{k}\left[  S_{n}\right]  $ (see Definition
\ref{def.bas.mur.bas}).

Now let us use this basis to compute the subquotients of our filtration
(\ref{eq.subsec.rep.maschke.kS3.filt}):

\begin{itemize}
\item The subquotient $\left\langle \nabla\right\rangle /0$ is a free
$\mathbf{k}$-module with basis $\left(  \overline{\nabla}\right)  $ (by
(\ref{eq.subsec.rep.maschke.kS3.filt0}) and
(\ref{eq.subsec.rep.maschke.kS3.filt1}) and Lemma \ref{lem.mod.quot.basis}).
As a representation of $S_{3}$, it is isomorphic to the trivial representation
$\mathbf{k}_{\operatorname*{triv}}$ of $S_{3}$ (as defined in Example
\ref{exa.rep.Sn-rep.triv}). Indeed, the $\mathbf{k}$-linear map%
\begin{align*}
\mathbf{k}_{\operatorname*{triv}}  &  \rightarrow\left\langle \nabla
\right\rangle ,\\
1  &  \mapsto\nabla
\end{align*}
is $S_{3}$-equivariant (since $g\nabla=\nabla$ for any $g\in S_{3}$) and
bijective, thus an isomorphism of representations of $S_{3}$. Hence, as
representations of $S_{3}$, we have $\mathbf{k}_{\operatorname*{triv}}%
\cong\left\langle \nabla\right\rangle \cong\left\langle \nabla\right\rangle
/0$. Thus,%
\[
\left\langle \nabla\right\rangle /0\cong\mathbf{k}_{\operatorname*{triv}%
}\ \ \ \ \ \ \ \ \ \ \text{as }S_{3}\text{-representations.}%
\]
Since $\mathbf{k}_{\operatorname*{triv}}$ is irreducible, we conclude that
$\left\langle \nabla\right\rangle /0$ is irreducible.

\item The subquotient $A_{1}/\left\langle \nabla\right\rangle $ is a free
$\mathbf{k}$-module with basis $\left(  \overline{\nabla_{3,3}},\ \overline
{\nabla_{2,3}}\right)  $ (by (\ref{eq.subsec.rep.maschke.kS3.filt1}) and
(\ref{eq.subsec.rep.maschke.kS3.filt2}) and Lemma \ref{lem.mod.quot.basis}). I
claim that, as a representation of $S_{3}$, it is isomorphic to the quotient
representation $\mathbf{k}^{3}/D\left(  \mathbf{k}^{3}\right)  $ (where
$\mathbf{k}^{3}$ is the natural representation of $S_{3}$), or, equivalently,
to the zero-sum subrepresentation $R\left(  \mathbf{k}^{3}\right)  $ (where we
are using the notations of Subsection \ref{subsec.rep.G-rep.subreps}). The
best way to see this is the following: Consider the standard basis $\left(
e_{1},e_{2},e_{3}\right)  $ of the natural representation $\mathbf{k}^{3}$,
and define a $\mathbf{k}$-linear map%
\begin{align*}
\mathbf{k}^{3}  &  \rightarrow A_{1},\\
e_{i}  &  \mapsto\nabla_{i,3}\ \ \ \ \ \ \ \ \ \ \text{for each }i\in\left[
3\right]  .
\end{align*}
This map is bijective (since $\left(  \nabla_{3,3},\ \nabla_{2,3}%
,\ \nabla_{1,3}\right)  $ is a basis of $A_{1}$) and $S_{3}$-equivariant
(since all $i\in\left[  3\right]  $ and $g\in S_{3}$ satisfy $g\rightharpoonup
e_{i}=e_{g\left(  i\right)  }$ and $g\nabla_{i,3}=\nabla_{g\left(  i\right)
,3}$). Thus, it is an isomorphism of $S_{3}$-representations. Moreover, this
isomorphism sends the submodule $D\left(  \mathbf{k}^{3}\right)  $ of
$\mathbf{k}^{3}$ to the submodule $\left\langle \nabla\right\rangle $ of
$A_{1}$ (since it sends $\left(  1,1,1\right)  \in\mathbf{k}^{3}$ to
$\nabla_{3,3}+\nabla_{2,3}+\nabla_{1,3}=\nabla$), and thus it gives rise to a
quotient isomorphism $\mathbf{k}^{3}/D\left(  \mathbf{k}^{3}\right)
\rightarrow A_{1}/\left\langle \nabla\right\rangle $ of $S_{3}$%
-representations. Hence,%
\[
A_{1}/\left\langle \nabla\right\rangle \cong\mathbf{k}^{3}/D\left(
\mathbf{k}^{3}\right)  \ \ \ \ \ \ \ \ \ \ \text{as }S_{3}%
\text{-representations.}%
\]
Furthermore, Theorem \ref{thm.rep.G-rep.Sn-nat.quots} \textbf{(a)} yields that
$\mathbf{k}^{3}/D\left(  \mathbf{k}^{3}\right)  \cong R\left(  \mathbf{k}%
^{3}\right)  $ as $S_{3}$-representations (since $\operatorname*{char}%
\mathbf{k}\neq3$ shows that $3$ is invertible in $\mathbf{k}$). Thus,
\[
A_{1}/\left\langle \nabla\right\rangle \cong\mathbf{k}^{3}/D\left(
\mathbf{k}^{3}\right)  \cong R\left(  \mathbf{k}^{3}\right)
\ \ \ \ \ \ \ \ \ \ \text{as }S_{3}\text{-representations.}%
\]
Finally, Theorem \ref{thm.rep.G-rep.Sn-nat.subreps} easily yields that
$R\left(  \mathbf{k}^{3}\right)  $ is irreducible (since $\operatorname*{char}%
\mathbf{k}\neq3$ and thus $D\left(  \mathbf{k}^{3}\right)  \not \subseteq
R\left(  \mathbf{k}^{3}\right)  $).

\item What about the subquotient $\left(  A_{1}+A_{2}\right)  /A_{1}$ ? We
could analyze this by hand, but it is easier to use a trick:

A general fact about modules (known as the \emph{second isomorphism theorem},
true for left $R$-modules over an arbitrary ring $R$) says that%
\[
\left(  A_{1}+A_{2}\right)  /A_{1}\cong A_{2}/\left(  A_{1}\cap A_{2}\right)
.
\]
But recall that $A_{1}\cap A_{2}=\left\langle \nabla\right\rangle $. Thus, we
find%
\[
\left(  A_{1}+A_{2}\right)  /A_{1}\cong A_{2}/\underbrace{\left(  A_{1}\cap
A_{2}\right)  }_{=\left\langle \nabla\right\rangle }=A_{2}/\left\langle
\nabla\right\rangle .
\]
Just as we showed $A_{1}/\left\langle \nabla\right\rangle \cong R\left(
\mathbf{k}^{3}\right)  $ above, we can show that $A_{2}/\left\langle
\nabla\right\rangle \cong R\left(  \mathbf{k}^{3}\right)  $ (as $S_{3}%
$-representations). Thus, we find that
\[
\left(  A_{1}+A_{2}\right)  /A_{1}\cong A_{2}/\left\langle \nabla\right\rangle
\cong R\left(  \mathbf{k}^{3}\right)  \ \ \ \ \ \ \ \ \ \ \text{as }%
S_{3}\text{-representations.}%
\]

\item Finally, what is the last subquotient $\mathbf{k}\left[  S_{3}\right]
/\left(  A_{1}+A_{2}\right)  $ ? As a $\mathbf{k}$-module, it has basis
$\left(  \overline{1}\right)  $ (again by a use of Lemma
\ref{lem.mod.quot.basis}). To see what $S_{3}$-action it bears, we probe its
basis element $\overline{1}$ with the two group elements $s_{1}$ and $t_{1,3}$
(which generate the group $S_{3}$). To wit, we have%
\[
s_{1}\rightharpoonup\overline{1}=\overline{s_{1}\rightharpoonup1}%
=\overline{s_{1}}=-\overline{1}\ \ \ \ \ \ \ \ \ \ \left(  \text{since
}1+s_{1}\in A_{1}+A_{2}\right)
\]
and%
\[
t_{1,3}\rightharpoonup\overline{1}=\overline{t_{1,3}\rightharpoonup
1}=\overline{t_{1,3}}=-\overline{1}\ \ \ \ \ \ \ \ \ \ \left(  \text{since
}1+t_{1,3}\in A_{1}+A_{2}\right)  .
\]
From these two equalities, it easily follows that%
\[
g\rightharpoonup\overline{1}=\left(  -1\right)  ^{g}\overline{1}%
\ \ \ \ \ \ \ \ \ \ \text{for each }g\in S_{3}%
\]
(since each $g\in S_{3}$ is a product of several $s_{1}$'s and $t_{1,3}$'s).
But this is precisely the rule for the $S_{3}$-action on the sign
representation $\mathbf{k}_{\operatorname*{sign}}$. Hence, the $\mathbf{k}%
$-linear map%
\begin{align*}
\mathbf{k}_{\operatorname*{sign}}  &  \rightarrow\mathbf{k}\left[
S_{3}\right]  /\left(  A_{1}+A_{2}\right)  ,\\
1  &  \mapsto\overline{1}%
\end{align*}
is $S_{3}$-equivariant. Since this map is also bijective\footnote{because
$\mathbf{k}\left[  S_{3}\right]  /\left(  A_{1}+A_{2}\right)  $ is a free
$\mathbf{k}$-module with basis $\left(  \overline{1}\right)  $}, we conclude
that it is an isomorphism of $S_{3}$-modules. Thus,%
\[
\mathbf{k}\left[  S_{3}\right]  /\left(  A_{1}+A_{2}\right)  \cong%
\mathbf{k}_{\operatorname*{sign}}\ \ \ \ \ \ \ \ \ \ \text{as }S_{3}%
\text{-representations.}%
\]
Of course, this shows that $\mathbf{k}\left[  S_{3}\right]  /\left(
A_{1}+A_{2}\right)  $ is irreducible.
\end{itemize}

Thus, all four subquotients of our filtration
(\ref{eq.subsec.rep.maschke.kS3.filt}) are irreducible, so that this
filtration is a composition series. Let us combine this and our other results
in a proposition:

\begin{proposition}
Let $\mathbf{k}$ be a field of characteristic $\neq3$. Then, the filtration%
\[
\left(  0,\ \left\langle \nabla\right\rangle ,\ A_{1},\ A_{1}+A_{2}%
,\ \mathbf{k}\left[  S_{3}\right]  \right)
\]
is a composition series of the left regular $\mathbf{k}\left[  S_{3}\right]
$-module $\mathbf{k}\left[  S_{3}\right]  $. Its subquotients (from left to
right) are isomorphic to $\mathbf{k}_{\operatorname*{triv}}$, $R\left(
\mathbf{k}^{3}\right)  $, $R\left(  \mathbf{k}^{3}\right)  $ and
$\mathbf{k}_{\operatorname*{sign}}$. Thus, $\mathbf{k}_{\operatorname*{triv}}%
$, $R\left(  \mathbf{k}^{3}\right)  $ and $\mathbf{k}_{\operatorname*{sign}}$
are all the irreducible representations of $S_{3}$.

Note that $\mathbf{k}_{\operatorname*{triv}}=\mathbf{k}_{\operatorname*{sign}%
}$ if $\operatorname*{char}\mathbf{k}=2$.
\end{proposition}

\begin{exercise}
\fbox{2} Find a composition series of $\mathbf{k}\left[  S_{3}\right]  $ when
$\mathbf{k}$ is a field of characteristic $3$. \medskip

[\textbf{Hint:} In this case, all irreps of $S_{3}$ are $1$-dimensional, so
the composition series will be a filtration of length $6$.]
\end{exercise}

So much for $\mathbf{k}\left[  S_{3}\right]  $. Analyzing $\mathbf{k}\left[
S_{n}\right]  $ with the same primitive tools is a lot harder. But I hope it
is clear that finding a composition series of $\mathbf{k}\left[  S_{n}\right]
$ is a worthwhile endeavor. Thus, we set ourselves the following goal:

\begin{question}
\label{quest.rep.kSn-compser}Find a composition series of $\mathbf{k}\left[
S_{n}\right]  $, at least when $\mathbf{k}$ is a field of characteristic $0$.
\end{question}

We will answer this question in Subsection \ref{subsec.bas.mur.triang2}, by
constructing such a composition series (actually, we will construct two, one
of which generalizes (\ref{eq.subsec.rep.maschke.kS3.filt})). This will also
yield a complete description of the irreps of $S_{n}$, although we will have
already obtained such a characterization by then (Corollary
\ref{cor.spechtmod.complete}).

\begin{exercise}
\textbf{(a)} \fbox{3} Set $A_{i}:=\left\{  \mathbf{a}\in\mathbf{k}\left[
S_{n}\right]  \ \mid\ \mathbf{a}s_{i}=\mathbf{a}\right\}  $ for each
$i\in\left[  n-1\right]  $ (where $n\in\mathbb{N}$ and $\mathbf{k}$ are again
arbitrary). Set $A:=A_{1}+A_{2}+\cdots+A_{n-1}$. Show that
\[
A=\operatorname*{Ker}\left(  \varepsilon\circ T_{\operatorname*{sign}}\right)
\ \ \ \ \ \ \ \ \ \ \text{and}\ \ \ \ \ \ \ \ \ \ \mathbf{k}\left[
S_{n}\right]  =A\oplus\operatorname*{span}\nolimits_{\mathbf{k}}\left\{
1\right\}  .
\]
(See Definition \ref{def.Tsign.Tsign} and Definition \ref{def.eps.eps} for the
meanings of $T_{\operatorname*{sign}}$ and $\varepsilon$.) \medskip

\textbf{(b)} \fbox{2} Assume that $\mathbf{k}$ is a field, and that $n\geq2$.
Prove that the left regular $\mathbf{k}\left[  S_{n}\right]  $-module
$\mathbf{k}\left[  S_{n}\right]  $ has a composition series of the form%
\[
\left(  0,\ \left\langle \nabla\right\rangle ,\ \ldots,\ A,\ \mathbf{k}\left[
S_{n}\right]  \right)  ,
\]
where $\left\langle \nabla\right\rangle =\operatorname*{span}%
\nolimits_{\mathbf{k}}\left\{  \nabla\right\}  $ and where $A$ is as in part
\textbf{(a)}. (The \textquotedblleft$\ldots$\textquotedblright\ part is far
less clear at this point. For $n=2$, the \textquotedblleft$\left\langle
\nabla\right\rangle ,\ \ldots,\ A$\textquotedblright\ segment should be read
as a single \textquotedblleft$\left\langle \nabla\right\rangle $%
\textquotedblright.) Moreover, prove that the first and the last subquotients
of this composition series are $\left\langle \nabla\right\rangle
/0\cong\mathbf{k}_{\operatorname*{triv}}$ and $\mathbf{k}\left[  S_{n}\right]
/A\cong\mathbf{k}_{\operatorname*{sign}}$ as representations of $S_{n}$.
\medskip

(Note that applying the automorphism $T_{\operatorname*{sign}}$ to the above
composition series yields another composition series, which starts with
$0,\ \left\langle \nabla^{-}\right\rangle $ instead of $0,\ \left\langle
\nabla\right\rangle $. Composition series are rarely unique!)
\end{exercise}

\subsubsection{Artin--Wedderburn}

The theorems of Maschke (Theorem \ref{thm.rep.G-rep.maschke}) and
Jordan--H\"{o}lder (Theorem \ref{thm.rep.G-rep.JH}) are two of the most
fundamental results in the representation theory of finite groups. They are
among the best results of their kind that hold at their level of generality.

Stronger results exist when we assume some things about $\mathbf{k}$ (such as
$\mathbf{k}$ being algebraically closed) or about $G$ (such as $G$ being a
$p$-group for a prime $p$). These results are commonly found in textbooks on
representation theory, such as \cite{EGHetc11} or \cite{Webb16}. We shall
abstain from discussing them here, except for briefly presenting one known as
the \emph{Artin--Wedderburn theorem}. We will not prove nor use this theorem,
but it helps us put Theorem \ref{thm.AWS.demo} \textbf{(a)} in context.

\begin{theorem}
[Artin--Wedderburn theorem]\label{thm.rep.AW}Let $\mathbf{k}$ be an
algebraically closed field of characteristic $0$. Let $G$ be a finite group.
Let $V_{1},V_{2},\ldots,V_{k}$ be all irreducible representations of $G$ over
$\mathbf{k}$, each listed exactly once (up to isomorphism). Consider the
curried forms $\rho_{i}:\mathbf{k}\left[  G\right]  \rightarrow
\operatorname*{End}\nolimits_{\mathbf{k}}\left(  V_{i}\right)  $ of these
representations $V_{i}$. Then, we can combine these curried forms into a
single $\mathbf{k}$-algebra morphism%
\begin{align*}
\rho:\mathbf{k}\left[  G\right]   &  \rightarrow\prod_{i=1}^{k}%
\operatorname*{End}\nolimits_{\mathbf{k}}\left(  V_{i}\right)  ,\\
\mathbf{a}  &  \mapsto\left(  \rho_{1}\left(  \mathbf{a}\right)  ,\ \rho
_{2}\left(  \mathbf{a}\right)  ,\ \ldots,\ \rho_{k}\left(  \mathbf{a}\right)
\right)  .
\end{align*}

This morphism is a $\mathbf{k}$-algebra isomorphism! Thus, in particular,%
\[
\mathbf{k}\left[  G\right]  \cong\prod_{i=1}^{k}\operatorname*{End}%
\nolimits_{\mathbf{k}}\left(  V_{i}\right)  \ \ \ \ \ \ \ \ \ \ \text{as
}\mathbf{k}\text{-algebras.}%
\]
Since each $V_{i}$ is a finite-dimensional $\mathbf{k}$-vector space, we can
identify its endomorphism ring $\operatorname*{End}\nolimits_{\mathbf{k}%
}\left(  V_{i}\right)  $ with a matrix ring $\mathbf{k}^{f_{i}\times f_{i}}$,
where $f_{i}=\dim V_{i}$. Thus,
\[
\mathbf{k}\left[  G\right]  \cong\prod_{i=1}^{k}\operatorname*{End}%
\nolimits_{\mathbf{k}}\left(  V_{i}\right)  \cong\prod_{i=1}^{k}%
\mathbf{k}^{f_{i}\times f_{i}}\ \ \ \ \ \ \ \ \ \ \text{as }\mathbf{k}%
\text{-algebras}%
\]
(a direct product of matrix rings). By comparing dimensions, this entails%
\begin{equation}
\left\vert G\right\vert =\sum_{i=1}^{k}f_{i}^{2}. \label{eq.thm.rep.AW.size}%
\end{equation}

\end{theorem}

This generalizes the $\mathbf{k}\left[  S_{3}\right]  \cong\mathbf{k}%
\times\mathbf{k}\times\mathbf{k}^{2\times2}$ isomorphism that we found in
Subsection \ref{subsec.intro.struct.3}, and to some extent generalizes Theorem
\ref{thm.AWS.demo} \textbf{(a)}. However, the generality comes at the cost of
stronger requirements (since the Artin--Wedderburn theorem requires
$\mathbf{k}$ to be algebraically closed) and weaker claims (the
Artin--Wedderburn theorem says nothing about what the $f_{i}$ are for any
given $G$). We will thus have to prove Theorem \ref{thm.AWS.demo} \textbf{(a)}
independently of the Artin--Wedderburn theorem.

Theorem \ref{thm.AWS.demo} \textbf{(a)} shows that the Artin--Wedderburn
theorem holds for $G=S_{n}$ even if $\mathbf{k}$ is not algebraically
closed.\footnote{Caveat: Theorem \ref{thm.AWS.demo} \textbf{(a)}, as stated,
says that $\mathbf{k}\left[  S_{n}\right]  $ is isomorphic to a direct product
of matrix rings, but it does not give any specific isomorphism. Thus, it does
not properly yield the claim of the Artin--Wedderburn theorem for $G=S_{n}$.
But the proof does.} For other finite groups $G$, the \textquotedblleft
algebraically closed\textquotedblright\ requirement in the Artin--Wedderburn
theorem can be necessary. For example, the cyclic group $C_{n}$ has $n$
mutually non-isomorphic $1$-dimensional representations over the field
$\mathbb{C}$ (all irreducible, of course), and no other irreps; thus, the
Artin--Wedderburn theorem yields a $\mathbb{C}$-algebra isomorphism%
\[
\mathbb{C}\left[  C_{n}\right]  \cong\underbrace{\mathbb{C}^{1\times1}%
\times\mathbb{C}^{1\times1}\times\cdots\times\mathbb{C}^{1\times1}}_{n\text{
times}}=\left(  \mathbb{C}^{1\times1}\right)  ^{n}\cong\mathbb{C}^{n}.
\]
This is precisely the isomorphism we saw in Exercise \ref{exe.Cn.DFT} (and
called the discrete Fourier transform). But no such isomorphism exists for
$\mathbf{k}=\mathbb{R}$ (unless $n\leq2$). \medskip

The Artin--Wedderburn theorem is very useful for computations:

\begin{itemize}
\item It allows you to work in the $\mathbf{k}$-algebra $\prod_{i=1}%
^{k}\mathbf{k}^{f_{i}\times f_{i}}$ instead of $\mathbf{k}\left[  G\right]  $.
In other words, instead of calculating with formal linear combinations of
elements of $G$, you can calculate with $k$-tuples of matrices. The latter are
better supported by computer algebra software (and the entrywise nature of
operations in $\prod_{i=1}^{k}\mathbf{k}^{f_{i}\times f_{i}}$ is conducive to parallelization).

\item If you have found a bunch of mutually non-isomorphic irreps $V_{1}%
,V_{2},\ldots,V_{m}$ of a finite group $G$ over an algebraically closed field
$\mathbf{k}$, then the formula (\ref{eq.thm.rep.AW.size}) helps you determine
whether you have found all irreps of $G$ or there are more to be found:
Indeed, if $\left\vert G\right\vert =\sum_{i=1}^{m}f_{i}^{2}$, then you have
found them all. If $\left\vert G\right\vert >\sum_{i=1}^{m}f_{i}^{2}$, then
there are more. (Why?)

\item The isomorphism $\rho$ in Theorem \ref{thm.rep.AW} is how the magic
elements $\mathbf{x},\mathbf{y},\mathbf{z},\mathbf{w}$ were found in
Subsection \ref{subsec.intro.struct.3}. Indeed, they should be the preimages
of $\left(  0,0,E_{i,j}\right)  $ under $\rho$, so we can compute them using
linear algebra (viewing $\rho$ as a matrix, and inverting this matrix).
\end{itemize}

\begin{exercise}
\fbox{1} Show that the Artin--Wedderburn theorem fails for $G=C_{3}$ and
$\mathbf{k}=\mathbb{R}$. In other words, show that $\mathbb{R}\left[
C_{3}\right]  $ is not isomorphic (as an $\mathbb{R}$-algebra) to the direct
product of the endomorphism rings of the irreps of $C_{3}$ over $\mathbb{R}$.
\end{exercise}

\subsection{\label{sec.rep.ncp}Noncommutative polynomials}

In this section, we will study yet another way to obtain representations of
symmetric groups, or two such ways to be precise, since two different
symmetric groups will be acting on the same space here.

The space in question is a space of \emph{noncommutative polynomials}. We
begin by defining them first. Then we will define the two actions of symmetric
groups on them, and we will use them to half-prove Theorem \ref{thm.Vperm.Vn}
\textbf{(b)} (\textquotedblleft half-prove\textquotedblright\ because we will
only reduce it to a result from noncommutative algebra that we won't prove).
Note that noncommutative polynomials are useful far beyond representation theory.

\subsubsection{Definition and symmetric group actions}

Let $k\in\mathbb{N}$. Recall that the usual (i.e., commutative) polynomial
ring $\mathbf{k}\left[  x_{1},x_{2},\ldots,x_{k}\right]  $ was defined as the
monoid algebra of the monoid of monomials over $\mathbf{k}$. Monomials are, in
essence, expressions of the form $x_{1}^{a_{1}}x_{2}^{a_{2}}\cdots
x_{k}^{a_{k}}$ with $a_{1},a_{2},\ldots,a_{k}\in\mathbb{N}$. They are
multiplied by adding the respective exponents:%
\[
\left(  x_{1}^{a_{1}}x_{2}^{a_{2}}\cdots x_{k}^{a_{k}}\right)  \left(
x_{1}^{b_{1}}x_{2}^{b_{2}}\cdots x_{k}^{b_{k}}\right)  =x_{1}^{a_{1}+b_{1}%
}x_{2}^{a_{2}+b_{2}}\cdots x_{k}^{a_{k}+b_{k}}.
\]

Let us try to define a noncommutative polynomial ring $\mathbf{k}\left\langle
x_{1},x_{2},\ldots,x_{k}\right\rangle $ similarly: It shall be the monoid
algebra of the monoid of \textquotedblleft noncommutative
monomials\textquotedblright\ over $\mathbf{k}$. Noncommutative monomials are,
in essence, expressions of the form $x_{i_{1}}x_{i_{2}}\cdots x_{i_{d}}$ with
$i_{1},i_{2},\ldots,i_{d}\in\left[  k\right]  $. (We want noncommutativity, so
we must not equate $x_{1}x_{2}x_{1}$ with $x_{1}x_{1}x_{2}$. Thus, we cannot
rewrite such monomials as $x_{1}^{a_{1}}x_{2}^{a_{2}}\cdots x_{k}^{a_{k}}$ in
general!) Such noncommutative monomials are multiplied by concatenation:%
\[
\left(  x_{i_{1}}x_{i_{2}}\cdots x_{i_{d}}\right)  \left(  x_{j_{1}}x_{j_{2}%
}\cdots x_{j_{e}}\right)  =x_{i_{1}}x_{i_{2}}\cdots x_{i_{d}}x_{j_{1}}%
x_{j_{2}}\cdots x_{j_{e}}.
\]
To make this rigorous, we encode these (so far undefined) noncommutative
monomials $x_{i_{1}}x_{i_{2}}\cdots x_{i_{d}}$ as the tuples $\left(
i_{1},i_{2},\ldots,i_{d}\right)  $ consisting of their subscripts. Thus, we
arrive at the following rigorous definition:

\begin{definition}
\label{def.ncpol.ncpol}Let $k\in\mathbb{N}$. \medskip

\textbf{(a)} A \emph{noncommutative monomial} in $k$ indeterminates
$x_{1},x_{2},\ldots,x_{k}$ is just a tuple (i.e., finite list) of elements of
$\left[  k\right]  $. The \emph{concatenation} $ij$ of two such tuples
$i=\left(  i_{1},i_{2},\ldots,i_{d}\right)  $ and $j=\left(  j_{1}%
,j_{2},\ldots,j_{e}\right)  $ is defined to be the tuple%
\[
ij:=\left(  i_{1},i_{2},\ldots,i_{d},j_{1},j_{2},\ldots,j_{e}\right)  .
\]
We define $\mathfrak{NM}_{k}$ to be the monoid of all noncommutative monomials
(in $k$ indeterminates $x_{1},x_{2},\ldots,x_{k}$), where the multiplication
is defined to be concatenation. The neutral element of $\mathfrak{NM}_{k}$ is
the empty tuple $\left(  {}\right)  $. \medskip

\textbf{(b)} For each $i\in\left[  k\right]  $, we let $x_{i}$ denote the
$1$-tuple $\left(  i\right)  \in\mathfrak{NM}_{k}$. Thus, each noncommutative
monomial $\left(  i_{1},i_{2},\ldots,i_{d}\right)  \in\mathfrak{NM}_{k}$ can
be rewritten as $x_{i_{1}}x_{i_{2}}\cdots x_{i_{d}}$. \medskip

\textbf{(c)} The \emph{noncommutative polynomial ring} $\mathbf{k}\left\langle
x_{1},x_{2},\ldots,x_{k}\right\rangle $ (also denoted $\mathbf{NP}_{k}$ in
this section) is the monoid algebra of the monoid $\mathfrak{NM}_{k}$ over
$\mathbf{k}$. Its elements are called \emph{noncommutative polynomials} in $k$
indeterminates $x_{1},x_{2},\ldots,x_{k}$ over $\mathbf{k}$.
\end{definition}

For example, $x_{3}x_{1}-5x_{1}x_{2}^{4}x_{3}x_{2}+7$ is a noncommutative
polynomial over $\mathbb{Z}$. (Note that $x_{2}^{4}=x_{2}x_{2}x_{2}x_{2}$.)
\medskip

Next, we define some basic features of noncommutative polynomials, analogous
to corresponding notions for usual (commutative) polynomials:

\begin{definition}
\label{def.ncpol.deg}\textbf{(a)} The \emph{degree} of a noncommutative
monomial $\mathfrak{m}=x_{i_{1}}x_{i_{2}}\cdots x_{i_{d}}$ is defined to be
$d$. (This is just the length of the tuple $\mathfrak{m}$.) This degree is
denoted $\deg\mathfrak{m}$. \medskip

\textbf{(b)} Let $d\in\mathbb{N}$. A noncommutative polynomial $p\in
\mathbf{NP}_{k}$ is said to be \emph{homogeneous of degree }$d$ if it is a
$\mathbf{k}$-linear combination of noncommutative monomials of degree $d$.
(Note that the zero of $\mathbf{NP}_{k}$ is homogeneous of any degree.)
\medskip

\textbf{(c)} Let $d\in\mathbb{N}$. We let $\mathbf{NP}_{k,d}$ be the
$\mathbf{k}$-submodule of $\mathbf{NP}_{k}$ consisting of all noncommutative
polynomials homogeneous of degree $d$. Equivalently, this is the $\mathbf{k}%
$-linear span of all noncommutative monomials of degree $d$. We call this the
\emph{degree-}$d$ \emph{component} of $\mathbf{NP}_{k}$.
\end{definition}

Here are some properties of noncommutative polynomials:

\begin{itemize}
\item The degree of a product is the sum of the degrees: If $p\in
\mathbf{NP}_{k,d}$ and $q\in\mathbf{NP}_{k,e}$, then $pq\in\mathbf{NP}%
_{k,d+e}$.

\item Each noncommutative polynomial can be written as the sum of its
homogeneous components (almost all of which are zero). Thus, $\mathbf{NP}%
_{k}=\bigoplus\limits_{d\in\mathbb{N}}\mathbf{NP}_{k,d}$ as a $\mathbf{k}$-module.

\item We can \emph{evaluate} a noncommutative polynomial $p\in\mathbf{NP}_{k}$
at any $k$-tuple $\left(  a_{1},a_{2},\ldots,a_{k}\right)  $ of elements of
any $\mathbf{k}$-algebra $A$. These elements $a_{1},a_{2},\ldots,a_{k}$ don't
even have to commute, unlike for usual (commutative) polynomials. This
evaluation is defined in the obvious way: We \textquotedblleft
plug\textquotedblright\ $a_{i}$ for each $x_{i}$ (that is, we replace each
noncommutative monomial $x_{i_{1}}x_{i_{2}}\cdots x_{i_{d}}$ by the product
$a_{i_{1}}a_{i_{2}}\cdots a_{i_{d}}$). The result of this evaluation is
denoted by $p\left(  a_{1},a_{2},\ldots,a_{k}\right)  $.

If $a_{1},a_{2},\ldots,a_{k}$ are $k$ elements of a $\mathbf{k}$-algebra $A$,
then the map%
\begin{align*}
\mathbf{NP}_{k}  &  \rightarrow A,\\
p  &  \mapsto p\left(  a_{1},a_{2},\ldots,a_{k}\right)
\end{align*}
is a $\mathbf{k}$-algebra morphism (called the \emph{evaluation morphism} at
$a_{1},a_{2},\ldots,a_{k}$). This yields a universal property of the
noncommutative polynomial ring $\mathbf{NP}_{k}$ (similar to the one of a
usual polynomial ring, but no longer requiring commutativity).

\item If $\mathbf{k}$ is an integral domain, then $\mathbf{NP}_{k}$ is a
\textquotedblleft noncommutative integral domain\textquotedblright\ (i.e., a
ring in which a product of nonzero elements is always nonzero). This can be
proved using leading terms with respect to an appropriately defined
lexicographic order on monomials.

\item For $k=1$, there is no difference between noncommutative polynomials and
usual polynomials. That is, we have $\mathbf{NP}_{1}\cong\mathbf{k}\left[
x\right]  $ as $\mathbf{k}$-algebras. (The isomorphism sends each
noncommutative monomial $\underbrace{x_{1}x_{1}\cdots x_{1}}_{i\text{ times}}$
to the commutative monomial $x^{i}$.)
\end{itemize}

\begin{remark}
\label{rmk.ncpol.tensor}If you know tensor algebra, then you will recognize
the noncommutative polynomial ring $\mathbf{NP}_{k}$ as (an isomorphic copy
of) the tensor algebra of the free $\mathbf{k}$-module $\mathbf{k}^{k}$.

More concretely: For any $d\in\mathbb{N}$, we have a $\mathbf{k}$-module
isomorphism%
\[
\mathbf{NP}_{k,d}\rightarrow\underbrace{\mathbf{k}^{k}\otimes\mathbf{k}%
^{k}\otimes\cdots\otimes\mathbf{k}^{k}}_{d\text{ times}},
\]
which sends each noncommutative monomial $x_{i_{1}}x_{i_{2}}\cdots x_{i_{d}}$
to the pure tensor $e_{i_{1}}\otimes e_{i_{2}}\otimes\cdots\otimes e_{i_{d}}$.
Taking the direct sum of these isomorphisms over all $d\in\mathbb{N}$, we
obtain a $\mathbf{k}$-algebra isomorphism from $\mathbf{NP}_{k}$ to the tensor
algebra of $\mathbf{k}^{k}$.
\end{remark}

Now we define actions of the symmetric groups $S_{k}$ and $S_{d}$ on
$\mathbf{NP}_{k,d}$:

\begin{definition}
\label{def.rep.NP-acts}Let $k,d\in\mathbb{N}$. Let $\mathfrak{NM}_{k,d}$ be
the set of all noncommutative monomials of degree $d$ in the $k$
indeterminates $x_{1},x_{2},\ldots,x_{k}$. These monomials form a basis of the
$\mathbf{k}$-module $\mathbf{NP}_{k,d}$, and thus we can view $\mathbf{NP}%
_{k,d}$ as the free $\mathbf{k}$-module $\mathbf{k}^{\left(  \mathfrak{NM}%
_{k,d}\right)  }$. \medskip

\textbf{(a)} We define a left action of the symmetric group $S_{k}$ on
$\mathfrak{NM}_{k,d}$ by setting%
\begin{align*}
g  &  \rightharpoonup\left(  x_{i_{1}}x_{i_{2}}\cdots x_{i_{d}}\right)
=x_{g\left(  i_{1}\right)  }x_{g\left(  i_{2}\right)  }\cdots x_{g\left(
i_{d}\right)  }\\
&  \ \ \ \ \ \ \ \ \ \ \text{for all }g\in S_{k}\text{ and }x_{i_{1}}x_{i_{2}%
}\cdots x_{i_{d}}\in\mathfrak{NM}_{k,d}.
\end{align*}
This is called the \emph{action on entries}. The permutation module
corresponding to this action (as explained in Example \ref{exa.mod.kG-on-kX})
is thus a left $\mathbf{k}\left[  S_{k}\right]  $-module $\mathbf{k}^{\left(
\mathfrak{NM}_{k,d}\right)  }=\mathbf{NP}_{k,d}$. Hence, $\mathbf{NP}_{k,d}$
becomes a representation of $S_{k}$. \medskip

\textbf{(b)} We define a left action of the symmetric group $S_{d}$ on
$\mathfrak{NM}_{k,d}$ by setting%
\begin{align*}
g  &  \rightharpoonup\left(  x_{i_{1}}x_{i_{2}}\cdots x_{i_{d}}\right)
=x_{i_{g^{-1}\left(  1\right)  }}x_{i_{g^{-1}\left(  2\right)  }}\cdots
x_{i_{g^{-1}\left(  d\right)  }}\\
&  \ \ \ \ \ \ \ \ \ \ \text{for all }g\in S_{d}\text{ and }x_{i_{1}}x_{i_{2}%
}\cdots x_{i_{d}}\in\mathfrak{NM}_{k,d}.
\end{align*}
This is called the \emph{action on places}. The permutation module
corresponding to this action (as explained in Example \ref{exa.mod.kG-on-kX})
is thus a left $\mathbf{k}\left[  S_{d}\right]  $-module $\mathbf{k}^{\left(
\mathfrak{NM}_{k,d}\right)  }=\mathbf{NP}_{k,d}$. Hence, $\mathbf{NP}_{k,d}$
becomes a representation of $S_{d}$. \medskip

(These two representations are not the same even when $k=d$.)
\end{definition}

Strictly speaking, much of this definition is redundant: After all, our
noncommutative monomials $x_{i_{1}}x_{i_{2}}\cdots x_{i_{d}}$ of degree $d$
are just the $d$-tuples $\left(  i_{1},i_{2},\ldots,i_{d}\right)  \in\left[
k\right]  ^{d}$ disguised by a fancy notation, and thus the set $\mathfrak{NM}%
_{k,d}$ is just $\left[  k\right]  ^{d}$. The two actions of $S_{k}$ and
$S_{d}$ that we have just defined on this set are precisely the action on
entries defined in Example \ref{exa.rep.G-sets.An} (for $A=\left[  k\right]
$) and the action on places defined in Example \ref{exa.rep.G-sets.place}
(again for $A=\left[  k\right]  $). The only new thing we are doing now is
extending these actions (by linearity) to the noncommutative \textbf{poly}%
nomials, i.e., to the $\mathbf{k}$-linear combinations of the noncommutative
monomials. For instance, for $k=4$ and $d=4$, we have%
\[
\operatorname*{cyc}\nolimits_{1,2,3}\rightharpoonup\left(  x_{2}x_{4}%
x_{1}x_{2}+5x_{1}^{3}x_{2}\right)  =x_{3}x_{4}x_{2}x_{3}+5x_{2}^{3}%
x_{3}\ \ \ \ \ \ \ \ \ \ \left(  \text{using the action on entries}\right)
\]
and%
\[
\operatorname*{cyc}\nolimits_{1,2,3}\rightharpoonup\left(  x_{2}x_{4}%
x_{1}x_{2}+5x_{1}^{3}x_{2}\right)  =x_{1}x_{2}x_{4}x_{2}+5x_{1}^{3}%
x_{2}\ \ \ \ \ \ \ \ \ \ \left(  \text{using the action on places}\right)  .
\]

The action on entries can be extended to the whole $\mathbf{k}$-algebra
$\mathbf{NP}_{k}$ (not just the degree-$d$ component $\mathbf{NP}_{k,d}$), but
the action on places cannot (since different groups would be acting on
different degrees).

\begin{remark}
The two actions commute, i.e., we have $g\rightharpoonup\left(
h\rightharpoonup p\right)  =h\rightharpoonup\left(  g\rightharpoonup p\right)
$ for any $g\in S_{k}$, any $h\in S_{d}$ and any $p\in\mathfrak{NM}_{k,d}$.
Thus, they can be combined into a single action of $S_{k}\times S_{d}$ on
$\mathfrak{NM}_{k,d}$ given by%
\[
\left(  g,h\right)  \rightharpoonup\left(  x_{i_{1}}x_{i_{2}}\cdots x_{i_{d}%
}\right)  =x_{g\left(  i_{h^{-1}\left(  1\right)  }\right)  }x_{g\left(
i_{h^{-1}\left(  2\right)  }\right)  }\cdots x_{g\left(  i_{h^{-1}\left(
d\right)  }\right)  }.
\]
This will not be used in what follows.
\end{remark}

Let us now explore the faithfulness of above-defined module structures on
$\mathbf{NP}_{k,d}$:

\begin{proposition}
\label{prop.rep.NP-acts-faith}\textbf{(a)} The left $\mathbf{k}\left[
S_{d}\right]  $-module $\mathbf{NP}_{k,d}$ (coming from the action on places)
is faithful if $k\geq d$. \medskip

\textbf{(b)} If $k<d$, then the element $\nabla^{-}\in\mathbf{k}\left[
S_{d}\right]  $ acts as $0$ on the left $\mathbf{k}\left[  S_{d}\right]
$-module $\mathbf{NP}_{k,d}$ (coming from the action on places). \medskip

\textbf{(c)} The left $\mathbf{k}\left[  S_{k}\right]  $-module $\mathbf{NP}%
_{k,d}$ (coming from the action on entries) is faithful if $d\geq k-1$.
\medskip

\textbf{(d)} If $d<k-1$, then the element $\nabla^{-}\in\mathbf{k}\left[
S_{k}\right]  $ acts as $0$ on the left $\mathbf{k}\left[  S_{k}\right]
$-module $\mathbf{NP}_{k,d}$ (coming from the action on entries).
\end{proposition}

\begin{proof}
[Proof of Proposition \ref{prop.rep.NP-acts-faith} \textbf{(a)}.]Assume that
$k\geq d$. Consider the noncommutative monomial $\mathfrak{m}=x_{1}x_{2}\cdots
x_{d}\in\mathfrak{NM}_{k,d}$ (this is well-defined, since $k\geq d$; otherwise
there would be no $x_{d}$). Then, any $g\in S_{d}$ satisfies%
\begin{equation}
g\rightharpoonup\mathfrak{m}=x_{g^{-1}\left(  1\right)  }x_{g^{-1}\left(
2\right)  }\cdots x_{g^{-1}\left(  d\right)  }
\label{pf.prop.rep.NP-acts-faith.1}%
\end{equation}
(by the definition of the action on places).

Now, let $\mathbf{a}$ and $\mathbf{b}$ be two elements of $\mathbf{k}\left[
S_{d}\right]  $ such that $\mathbf{a}\mathfrak{m}=\mathbf{b}\mathfrak{m}$. We
shall show that $\mathbf{a}=\mathbf{b}$.

Indeed, write the difference $\mathbf{a}-\mathbf{b}\in\mathbf{k}\left[
S_{d}\right]  $ in the form
\begin{equation}
\mathbf{a}-\mathbf{b}=\sum_{g\in S_{d}}\lambda_{g}g
\label{pf.prop.rep.NP-acts-faith.3}%
\end{equation}
for some scalars $\lambda_{g}\in\mathbf{k}$. Then,%
\begin{align*}
\left(  \mathbf{a}-\mathbf{b}\right)  \mathfrak{m}  &  =\left(  \sum_{g\in
S_{d}}\lambda_{g}g\right)  \mathfrak{m}=\sum_{g\in S_{d}}\lambda
_{g}\underbrace{g\rightharpoonup\mathfrak{m}}_{\substack{=x_{g^{-1}\left(
1\right)  }x_{g^{-1}\left(  2\right)  }\cdots x_{g^{-1}\left(  d\right)
}\\\text{(by (\ref{pf.prop.rep.NP-acts-faith.1}))}}}\\
&  =\sum_{g\in S_{d}}\lambda_{g}x_{g^{-1}\left(  1\right)  }x_{g^{-1}\left(
2\right)  }\cdots x_{g^{-1}\left(  d\right)  }.
\end{align*}
Hence,%
\begin{align}
\sum_{g\in S_{d}}\lambda_{g}x_{g^{-1}\left(  1\right)  }x_{g^{-1}\left(
2\right)  }\cdots x_{g^{-1}\left(  d\right)  }  &  =\left(  \mathbf{a}%
-\mathbf{b}\right)  \mathfrak{m}=\mathbf{a}\mathfrak{m}-\mathbf{b}%
\mathfrak{m}\nonumber\\
&  =0 \label{pf.prop.rep.NP-acts-faith.4}%
\end{align}
(since we assumed that $\mathbf{a}\mathfrak{m}=\mathbf{b}\mathfrak{m}$).

But the noncommutative monomials $x_{g^{-1}\left(  1\right)  }x_{g^{-1}\left(
2\right)  }\cdots x_{g^{-1}\left(  d\right)  }$ for all $g\in S_{d}$ are
distinct (since they are just the $d$-tuples $\left(  g^{-1}\left(  1\right)
,g^{-1}\left(  2\right)  ,\ldots,g^{-1}\left(  d\right)  \right)  $, and
clearly you can reconstruct a permutation $g\in S_{d}$ from this $d$-tuple),
and thus are $\mathbf{k}$-linearly independent in $\mathbf{NP}_{k,d}$. In
other words, a $\mathbf{k}$-linear combination of these monomials
$x_{g^{-1}\left(  1\right)  }x_{g^{-1}\left(  2\right)  }\cdots x_{g^{-1}%
\left(  d\right)  }$ can only be $0$ if all coefficients in it are $0$.

Hence, the equality (\ref{pf.prop.rep.NP-acts-faith.4}) shows that all
coefficients $\lambda_{g}$ are $0$. Thus, we can rewrite
(\ref{pf.prop.rep.NP-acts-faith.3}) as $\mathbf{a}-\mathbf{b}=\sum_{g\in
S_{d}}0g=0$, so that $\mathbf{a}=\mathbf{b}$.

Forget that we fixed $\mathbf{a}$ and $\mathbf{b}$. We thus have shown that if
$\mathbf{a}$ and $\mathbf{b}$ are two elements of $\mathbf{k}\left[
S_{d}\right]  $ such that $\mathbf{a}\mathfrak{m}=\mathbf{b}\mathfrak{m}$,
then $\mathbf{a}=\mathbf{b}$. In other words, if two elements of
$\mathbf{k}\left[  S_{d}\right]  $ act in the same way on $\mathfrak{m}$, then
they are equal. Thus, we can tell distinct elements of $\mathbf{k}\left[
S_{d}\right]  $ apart by looking at how they act on $\mathfrak{m}$. Therefore,
the action of $\mathbf{k}\left[  S_{d}\right]  $ is faithful (and moreover,
just looking at the action on the single element $\mathfrak{m}$ is already
enough for this faithfulness). This proves Proposition
\ref{prop.rep.NP-acts-faith} \textbf{(a)}.
\end{proof}

\begin{exercise}
\fbox{5 (2, 1, 2)} Prove the other three parts of Proposition
\ref{prop.rep.NP-acts-faith}.
\end{exercise}

More generally, the elements of $\mathbf{k}\left[  S_{d}\right]  $ or of
$\mathbf{k}\left[  S_{k}\right]  $ that act as $0$ on $\mathbf{NP}_{k,d}$ can
be characterized -- see \cite[Theorems 2.8.1 and 2.8.3]{rooksn} (where
$\mathbf{NP}_{k,d}$ appears in the equivalent guise of the tensor power
$\left(  \mathbf{k}^{k}\right)  ^{\otimes d}$).

\subsubsection{\label{subsec.rep.ncp.Vn}Application: the Dynkin element
$\mathbf{V}_{n}$}

Using Proposition \ref{prop.rep.NP-acts-faith} \textbf{(a)}, we can prove the
equality of two elements of $\mathbf{k}\left[  S_{d}\right]  $ by showing that
they act equally on $\mathbf{NP}_{k,d}$ for some $k\geq d$. This technique can
be quite useful; let me show an example of its use. Sadly, the proof will be
incomplete, since it boils down to an algebraic result that cannot be shown
very easily.

Assume that $n\geq1$ for the rest of Subsection \ref{subsec.rep.ncp.Vn}.

Recall from Theorem \ref{thm.Vperm.Vn} that%
\begin{align*}
\mathbf{V}_{n}:=  &  \ \left(  1-\operatorname*{cyc}\nolimits_{2,1}\right)
\left(  1-\operatorname*{cyc}\nolimits_{3,2,1}\right)  \left(
1-\operatorname*{cyc}\nolimits_{4,3,2,1}\right)  \cdots\left(
1-\operatorname*{cyc}\nolimits_{n,n-1,\ldots,1}\right) \\
=  &  \ \prod_{i=2}^{n}\left(  1-\operatorname*{cyc}\nolimits_{i,i-1,\ldots
,1}\right)  \in\mathbf{k}\left[  S_{n}\right]  .
\end{align*}
Theorem \ref{thm.Vperm.Vn} \textbf{(b)} claims that $\mathbf{V}_{n}%
^{2}=n\mathbf{V}_{n}$. We have not proved this so far.

How can we prove this? The first step will be a simple transformation inside
$\mathbf{k}\left[  S_{n}\right]  $. Indeed, recall the antipode $S$ of
$\mathbf{k}\left[  S_{n}\right]  $ (see Definition \ref{def.S.S}). We know
from Theorem \ref{thm.S.auto} \textbf{(a)} that the antipode $S$ is an algebra
anti-automorphism. In particular, $S$ is injective. Thus, in order to prove
$\mathbf{V}_{n}^{2}=n\mathbf{V}_{n}$, it will suffice to show that $S\left(
\mathbf{V}_{n}^{2}\right)  =S\left(  n\mathbf{V}_{n}\right)  $. But we have
$S\left(  \mathbf{V}_{n}^{2}\right)  =S\left(  \mathbf{V}_{n}\right)  ^{2}$
(since $S$ is an algebra anti-automorphism) and $S\left(  n\mathbf{V}%
_{n}\right)  =nS\left(  \mathbf{V}_{n}\right)  $ (likewise). Hence, it will
suffice to show that%
\begin{equation}
S\left(  \mathbf{V}_{n}\right)  ^{2}=nS\left(  \mathbf{V}_{n}\right)  .
\label{pf.thm.Vperm.Vn.suff2}%
\end{equation}

Note that the definition of $\mathbf{V}_{n}$ yields%
\begin{align}
S\left(  \mathbf{V}_{n}\right)   &  =S\left(  \left(  1-\operatorname*{cyc}%
\nolimits_{2,1}\right)  \left(  1-\operatorname*{cyc}\nolimits_{3,2,1}\right)
\left(  1-\operatorname*{cyc}\nolimits_{4,3,2,1}\right)  \cdots\left(
1-\operatorname*{cyc}\nolimits_{n,n-1,\ldots,1}\right)  \right) \nonumber\\
&  =S\left(  1-\operatorname*{cyc}\nolimits_{n,n-1,\ldots,1}\right)  \cdot
S\left(  1-\operatorname*{cyc}\nolimits_{n-1,n-2,\ldots,1}\right)  \cdot
\cdots\cdot S\left(  1-\operatorname*{cyc}\nolimits_{2,1}\right) \nonumber\\
&  \ \ \ \ \ \ \ \ \ \ \ \ \ \ \ \ \ \ \ \ \left(  \text{since }S\text{ is an
algebra anti-automorphism}\right) \nonumber\\
&  =\left(  1-\operatorname*{cyc}\nolimits_{1,2,\ldots,n}\right)  \left(
1-\operatorname*{cyc}\nolimits_{1,2,\ldots,n-1}\right)  \cdots\left(
1-\operatorname*{cyc}\nolimits_{1,2}\right)  \label{pf.thm.Vperm.Vn.Svn=}%
\end{align}
(since each $k\in\left\{  2,3,\ldots,n\right\}  $ satisfies%
\begin{align*}
S\left(  1-\operatorname*{cyc}\nolimits_{k,k-1,\ldots,1}\right)   &
=\underbrace{1^{-1}}_{=1}-\underbrace{\operatorname*{cyc}%
\nolimits_{k,k-1,\ldots,1}^{-1}}_{=\operatorname*{cyc}\nolimits_{1,2,\ldots
,k}}\ \ \ \ \ \ \ \ \ \ \left(  \text{by the definition of }S\right) \\
&  =1-\operatorname*{cyc}\nolimits_{1,2,\ldots,k}%
\end{align*}
).

Now is the time to use faithful modules. Recall that the left $\mathbf{k}%
\left[  S_{n}\right]  $-module $\mathbf{NP}_{k,n}$ (coming from the action on
places) is faithful for any integer $k\geq n$ (by Proposition
\ref{prop.rep.NP-acts-faith} \textbf{(a)}, applied to $d=n$). Hence, in order
to show that $S\left(  \mathbf{V}_{n}\right)  ^{2}=nS\left(  \mathbf{V}%
_{n}\right)  $, it will suffice to prove that%
\begin{equation}
S\left(  \mathbf{V}_{n}\right)  ^{2}\cdot p=nS\left(  \mathbf{V}_{n}\right)
\cdot p \label{pf.thm.Vperm.Vn.suff3}%
\end{equation}
for all $k\in\mathbb{N}$ and all noncommutative polynomials $p\in
\mathbf{NP}_{k,n}$, where we are using the action on places. (Actually, it
will even suffice to prove (\ref{pf.thm.Vperm.Vn.suff3}) for any one integer
$k\geq n$; but this is no easier than proving (\ref{pf.thm.Vperm.Vn.suff3})
for all $k\in\mathbb{N}$.)

We fix an integer $k\in\mathbb{N}$ for the rest of Subsection
\ref{subsec.rep.ncp.Vn}.

In order to prove (\ref{pf.thm.Vperm.Vn.suff3}), we must figure out how
$S\left(  \mathbf{V}_{n}\right)  $ acts on $\mathbf{NP}_{k,n}$. There is a
really nice formula for this, which is in fact the reason why $\mathbf{V}_{n}$
was first introduced.

Before we state it, we recall that the set $\mathbf{NP}_{k,1}$ consists of
homogeneous noncommutative polynomials of degree $1$; these are $\mathbf{k}%
$-linear combinations $\lambda_{1}x_{1}+\lambda_{2}x_{2}+\cdots+\lambda
_{k}x_{k}$ of indeterminates. Any product $a_{1}a_{2}\cdots a_{n}$ of $n$ such
combinations $a_{1},a_{2},\ldots,a_{n}\in\mathbf{NP}_{k,1}$ is a homogeneous
polynomial of degree $n$, that is, an element of $\mathbf{NP}_{k,n}$.
Moreover, all noncommutative monomials $x_{i_{1}}x_{i_{2}}\cdots x_{i_{n}}$ of
degree $n$ are such products (since the indeterminates $x_{i_{1}},x_{i_{2}%
},\ldots,x_{i_{n}}$ themselves belong to $\mathbf{NP}_{k,1}$). Hence, such
products $a_{1}a_{2}\cdots a_{n}$ span the $\mathbf{k}$-module $\mathbf{NP}%
_{k,n}$. Thus, a $\mathbf{k}$-linear map on $\mathbf{NP}_{k,n}$ is uniquely
determined by its values on such products $a_{1}a_{2}\cdots a_{n}$.

We shall describe the action of $S\left(  \mathbf{V}_{n}\right)  $ on
$\mathbf{NP}_{k,n}$ by how it acts on such products. To do so, we introduce
another notation:

Recall that the commutator $ab-ba$ of two elements $a$ and $b$ of a ring is
denoted by $\left[  a,b\right]  $. More generally, if $a_{1},a_{2}%
,\ldots,a_{m}$ are $m$ elements of a ring (where $m\geq1$), then their
\emph{left-nested commutator }$\left[  \cdots\left[  \left[  a_{1}%
,a_{2}\right]  ,a_{3}\right]  ,\ldots,a_{m}\right]  $ is defined recursively
as follows:

\begin{itemize}
\item For $m=1$, it is defined to be $a_{1}$.

\item For $m>1$, it is defined to be the (usual) commutator $\left[
u,a_{m}\right]  $, where $u$ is the left-nested commutator $\left[
\cdots\left[  \left[  a_{1},a_{2}\right]  ,a_{3}\right]  ,\ldots
,a_{m-1}\right]  $ of the first $m-1$ elements $a_{1},a_{2},\ldots,a_{m-1}$.
\end{itemize}

\noindent Thus, the left-nested commutator $\left[  \cdots\left[  \left[
a_{1},a_{2}\right]  ,a_{3}\right]  ,\ldots,a_{m}\right]  $ is the usual
commutator $\left[  a_{1},a_{2}\right]  $ when $m=2$, the nested commutator
$\left[  \left[  a_{1},a_{2}\right]  ,a_{3}\right]  $ when $m=3$, and the
nested commutator $\left[  \left[  \left[  a_{1},a_{2}\right]  ,a_{3}\right]
,a_{4}\right]  $ when $m=4$.

Now we can finally describe the action of $S\left(  \mathbf{V}_{n}\right)  $
on $\mathbf{NP}_{k,n}$:

\begin{proposition}
\label{prop.Vperm.Vn.LNC}Let $n\geq1$. Let $a_{1},a_{2},\ldots,a_{n}%
\in\mathbf{NP}_{k,1}$. Then,%
\[
S\left(  \mathbf{V}_{n}\right)  \cdot\left(  a_{1}a_{2}\cdots a_{n}\right)
=\left[  \cdots\left[  \left[  a_{1},a_{2}\right]  ,a_{3}\right]
,\ldots,a_{n}\right]  .
\]

\end{proposition}

For instance, for $n=3$, this is saying that%
\[
S\left(  \mathbf{V}_{3}\right)  \cdot\left(  abc\right)  =\left[  \left[
a,b\right]  ,c\right]  \ \ \ \ \ \ \ \ \ \ \text{for all }a,b,c\in
\mathbf{NP}_{k,1}.
\]
Let us check this: By (\ref{pf.thm.Vperm.Vn.Svn=}), we have%
\[
S\left(  \mathbf{V}_{3}\right)  =\left(  1-\operatorname*{cyc}%
\nolimits_{1,2,3}\right)  \left(  1-\operatorname*{cyc}\nolimits_{1,2}\right)
=1-\operatorname*{cyc}\nolimits_{1,2,3}-\,s_{1}+t_{1,3}%
\]
and thus%
\begin{align*}
S\left(  \mathbf{V}_{3}\right)  \cdot\left(  abc\right)   &  =\left(
1-\operatorname*{cyc}\nolimits_{1,2,3}-\,s_{1}+t_{1,3}\right)  \cdot abc\\
&  =abc-cab-bac+cba\ \ \ \ \ \ \ \ \ \ \left(  \text{action on places!}%
\right)  ,
\end{align*}
which indeed agrees with%
\begin{align*}
\left[  \left[  a,b\right]  ,c\right]   &  =\left[  ab-ba,\ c\right]  =\left(
ab-ba\right)  c-c\left(  ab-ba\right) \\
&  =abc-bac-cab+cba.
\end{align*}

The proof of Proposition \ref{prop.Vperm.Vn.LNC} is not hard. The main step is
the following lemma:

\begin{lemma}
\label{lem.Vperm.Vn.LNC-k}Let $j\in\left[  n\right]  $. Let $u,g,v\in
\mathbf{NP}_{k}$ be homogeneous noncommutative polynomials of respective
degrees $j-1$, $1$ and $n-j$. Then,%
\begin{equation}
\operatorname*{cyc}\nolimits_{1,2,\ldots,j}\cdot\left(  ugv\right)  =guv
\label{eq.lem.Vperm.Vn.LNC-k.c}%
\end{equation}
and%
\begin{equation}
\left(  1-\operatorname*{cyc}\nolimits_{1,2,\ldots,j}\right)  \cdot\left(
ugv\right)  =\left[  u,g\right]  v. \label{eq.lem.Vperm.Vn.LNC-k.1-c}%
\end{equation}

\end{lemma}

\begin{proof}
Both equalities (\ref{eq.lem.Vperm.Vn.LNC-k.c}) and
(\ref{eq.lem.Vperm.Vn.LNC-k.1-c}) depend $\mathbf{k}$-linearly on each of $u$,
$g$ and $v$. Thus, we can WLOG assume that the homogeneous noncommutative
polynomials $u$, $g$ and $v$ are actually noncommutative monomials of the
respective degrees (since the noncommutative monomials of a given degree $d$
span $\mathbf{NP}_{k,d}$). In other words, we can WLOG assume that%
\[
u=x_{i_{1}}x_{i_{2}}\cdots x_{i_{j-1}},\ \ \ \ \ \ \ \ \ \ g=x_{i_{j}%
},\ \ \ \ \ \ \ \ \ \ \text{and}\ \ \ \ \ \ \ \ \ \ v=x_{i_{j+1}}x_{i_{j+2}%
}\cdots x_{i_{n}}%
\]
for some $i_{1},i_{2},\ldots,i_{n}\in\left[  k\right]  $. Assume this, and
consider these $i_{1},i_{2},\ldots,i_{n}$. Multiplying the three equalities
\[
u=x_{i_{1}}x_{i_{2}}\cdots x_{i_{j-1}},\ \ \ \ \ \ \ \ \ \ g=x_{i_{j}%
},\ \ \ \ \ \ \ \ \ \ \text{and}\ \ \ \ \ \ \ \ \ \ v=x_{i_{j+1}}x_{i_{j+2}%
}\cdots x_{i_{n}}%
\]
together, we obtain%
\[
ugv=\left(  x_{i_{1}}x_{i_{2}}\cdots x_{i_{j-1}}\right)  x_{i_{j}}\left(
x_{i_{j+1}}x_{i_{j+2}}\cdots x_{i_{n}}\right)  =x_{i_{1}}x_{i_{2}}\cdots
x_{i_{n}}.
\]
Set $\sigma:=\operatorname*{cyc}\nolimits_{1,2,\ldots,j}$. Then, $\sigma
^{-1}=\operatorname*{cyc}\nolimits_{1,2,\ldots,j}^{-1}=\operatorname*{cyc}%
\nolimits_{j,j-1,\ldots,1}$ and therefore
\begin{align}
&  \left(  \sigma^{-1}\left(  1\right)  ,\sigma^{-1}\left(  2\right)
,\ldots,\sigma^{-1}\left(  n\right)  \right) \nonumber\\
&  =\left(  j,1,2,\ldots,j-1,j+1,j+2,\ldots,n\right)  .
\label{pf.lem.Vperm.Vn.LNC-k.3}%
\end{align}
But%
\begin{align*}
\sigma\cdot\underbrace{\left(  ugv\right)  }_{=x_{i_{1}}x_{i_{2}}\cdots
x_{i_{n}}}  &  =\sigma\cdot\left(  x_{i_{1}}x_{i_{2}}\cdots x_{i_{n}}\right)
\\
&  =x_{i_{\sigma^{-1}\left(  1\right)  }}x_{i_{\sigma^{-1}\left(  2\right)  }%
}\cdots x_{i_{\sigma^{-1}\left(  n\right)  }}\ \ \ \ \ \ \ \ \ \ \left(
\begin{array}
[c]{c}%
\text{by the definition of}\\
\text{the action on places}%
\end{array}
\right) \\
&  =\underbrace{x_{i_{j}}}_{=g}\underbrace{x_{i_{1}}x_{i_{2}}\cdots
x_{i_{j-1}}}_{=u}\underbrace{x_{i_{j+1}}x_{i_{j+2}}\cdots x_{i_{n}}}%
_{=v}\ \ \ \ \ \ \ \ \ \ \left(  \text{by (\ref{pf.lem.Vperm.Vn.LNC-k.3}%
)}\right) \\
&  =guv
\end{align*}
and thus%
\[
\left(  1-\sigma\right)  \cdot\left(  ugv\right)  =ugv-\underbrace{\sigma
\cdot\left(  ugv\right)  }_{=guv}=ugv-guv=\underbrace{\left(  ug-gu\right)
}_{=\left[  u,g\right]  }v=\left[  u,g\right]  v.
\]
In view of $\sigma=\operatorname*{cyc}\nolimits_{1,2,\ldots,j}$, we can
rewrite these two equalities as%
\[
\operatorname*{cyc}\nolimits_{1,2,\ldots,j}\cdot\left(  ugv\right)  =guv
\]
and%
\[
\left(  1-\operatorname*{cyc}\nolimits_{1,2,\ldots,j}\right)  \cdot\left(
ugv\right)  =\left[  u,g\right]  v.
\]
Thus, Lemma \ref{lem.Vperm.Vn.LNC-k} is proved.
\end{proof}

\begin{proof}
[Proof of Proposition \ref{prop.Vperm.Vn.LNC}.]We start with an example: For
$n=4$, we have%
\[
S\left(  \mathbf{V}_{4}\right)  =\left(  1-\operatorname*{cyc}%
\nolimits_{1,2,3,4}\right)  \left(  1-\operatorname*{cyc}\nolimits_{1,2,3}%
\right)  \left(  1-\operatorname*{cyc}\nolimits_{1,2}\right)
\ \ \ \ \ \ \ \ \ \ \left(  \text{by (\ref{pf.thm.Vperm.Vn.Svn=})}\right)
\]
and thus%
\begin{align*}
&  S\left(  \mathbf{V}_{4}\right)  \cdot\left(  abcd\right) \\
&  =\left(  1-\operatorname*{cyc}\nolimits_{1,2,3,4}\right)  \left(
1-\operatorname*{cyc}\nolimits_{1,2,3}\right)  \underbrace{\left(
1-\operatorname*{cyc}\nolimits_{1,2}\right)  \cdot\left(  abcd\right)
}_{\substack{=\left[  a,b\right]  cd\\\text{(by
(\ref{eq.lem.Vperm.Vn.LNC-k.1-c}), applied to }j=2\text{, }u=a\text{,
}g=b\text{ and }v=cd\text{)}}}\\
&  =\left(  1-\operatorname*{cyc}\nolimits_{1,2,3,4}\right)
\underbrace{\left(  1-\operatorname*{cyc}\nolimits_{1,2,3}\right)
\cdot\left(  \left[  a,b\right]  cd\right)  }_{\substack{=\left[  \left[
a,b\right]  ,c\right]  d\\\text{(by (\ref{eq.lem.Vperm.Vn.LNC-k.1-c}), applied
to }j=3\text{, }u=\left[  a,b\right]  \text{, }g=c\text{ and }v=d\text{)}}}\\
&  =\left(  1-\operatorname*{cyc}\nolimits_{1,2,3,4}\right)  \cdot\left(
\left[  \left[  a,b\right]  ,c\right]  d\right) \\
&  =\left[  \left[  \left[  a,b\right]  ,c\right]  ,d\right]
\ \ \ \ \ \ \ \ \ \ \left(  \text{by (\ref{eq.lem.Vperm.Vn.LNC-k.1-c}),
applied to }j=4\text{, }u=\left[  \left[  a,b\right]  ,c\right]  \text{,
}g=d\text{ and }v=1\right)  .
\end{align*}
The same type of computation works for any $n$: Repeated use of
(\ref{eq.lem.Vperm.Vn.LNC-k.1-c}) gradually turns the product $a_{1}%
a_{2}\cdots a_{n}$ into the left-nested commutator $\left[  \cdots\left[
\left[  a_{1},a_{2}\right]  ,a_{3}\right]  ,\ldots,a_{n}\right]  $ by
progressively \textquotedblleft enclosing\textquotedblright\ more and more
factors into commutator brackets. This can be made rigorous as an induction
argument as follows: \medskip

\begin{fineprint}
We shall show that
\begin{align}
&  \left(  1-\operatorname*{cyc}\nolimits_{1,2,\ldots,j}\right)  \left(
1-\operatorname*{cyc}\nolimits_{1,2,\ldots,j-1}\right)  \cdots\left(
1-\operatorname*{cyc}\nolimits_{1,2}\right)  \cdot\left(  a_{1}a_{2}\cdots
a_{n}\right) \nonumber\\
&  =\left[  \cdots\left[  \left[  a_{1},a_{2}\right]  ,a_{3}\right]
,\ldots,a_{j}\right]  a_{j+1}a_{j+2}\cdots a_{n}
\label{pf.prop.Vperm.Vn.LNC.goal}%
\end{align}
for each $j\in\left[  n\right]  $.

Indeed, we shall prove this by induction on $j$:

\textit{Base case:} For $j=1$, both sides of (\ref{pf.prop.Vperm.Vn.LNC.goal})
equal $a_{1}a_{2}\cdots a_{n}$. (For the left hand side, this is because the
product $\left(  1-\operatorname*{cyc}\nolimits_{1,2,\ldots,j}\right)  \left(
1-\operatorname*{cyc}\nolimits_{1,2,\ldots,j-1}\right)  \cdots\left(
1-\operatorname*{cyc}\nolimits_{1,2}\right)  $ is empty and thus equals $1$.
For the right hand side, this is because the left-nested commutator $\left[
\cdots\left[  \left[  a_{1},a_{2}\right]  ,a_{3}\right]  ,\ldots,a_{j}\right]
$ is defined to be $a_{1}$ in this case, and therefore the right hand side is
$\underbrace{\left[  \cdots\left[  \left[  a_{1},a_{2}\right]  ,a_{3}\right]
,\ldots,a_{j}\right]  }_{=a_{1}}\underbrace{a_{j+1}a_{j+2}\cdots a_{n}%
}_{=a_{2}a_{3}\cdots a_{n}}=a_{1}\left(  a_{2}a_{3}\cdots a_{n}\right)
=a_{1}a_{2}\cdots a_{n}$.) Hence, the equality
(\ref{pf.prop.Vperm.Vn.LNC.goal}) holds for $j=1$.

\textit{Induction step:} Let $j\in\left[  n\right]  $ be such that $j>1$.
Assume as the induction hypothesis that (\ref{pf.prop.Vperm.Vn.LNC.goal})
holds for $j-1$ instead of $j$. In other words, assume that%
\begin{align}
&  \left(  1-\operatorname*{cyc}\nolimits_{1,2,\ldots,j-1}\right)  \left(
1-\operatorname*{cyc}\nolimits_{1,2,\ldots,j-2}\right)  \cdots\left(
1-\operatorname*{cyc}\nolimits_{1,2}\right)  \cdot\left(  a_{1}a_{2}\cdots
a_{n}\right) \nonumber\\
&  =\left[  \cdots\left[  \left[  a_{1},a_{2}\right]  ,a_{3}\right]
,\ldots,a_{j-1}\right]  a_{j}a_{j+1}\cdots a_{n}.
\label{pf.prop.Vperm.Vn.LNC.IH}%
\end{align}
We must now prove that (\ref{pf.prop.Vperm.Vn.LNC.goal}) holds for $j$ as
well. But this follows from%
\begin{align*}
&  \left(  1-\operatorname*{cyc}\nolimits_{1,2,\ldots,j}\right)  \left(
1-\operatorname*{cyc}\nolimits_{1,2,\ldots,j-1}\right)  \cdots\left(
1-\operatorname*{cyc}\nolimits_{1,2}\right)  \cdot\left(  a_{1}a_{2}\cdots
a_{n}\right) \\
&  =\left(  1-\operatorname*{cyc}\nolimits_{1,2,\ldots,j}\right)
\cdot\underbrace{\left(  1-\operatorname*{cyc}\nolimits_{1,2,\ldots
,j-1}\right)  \left(  1-\operatorname*{cyc}\nolimits_{1,2,\ldots,j-2}\right)
\cdots\left(  1-\operatorname*{cyc}\nolimits_{1,2}\right)  \cdot\left(
a_{1}a_{2}\cdots a_{n}\right)  }_{\substack{=\left[  \cdots\left[  \left[
a_{1},a_{2}\right]  ,a_{3}\right]  ,\ldots,a_{j-1}\right]  a_{j}a_{j+1}\cdots
a_{n}\\\text{(by (\ref{pf.prop.Vperm.Vn.LNC.IH}))}}}\\
&  =\left(  1-\operatorname*{cyc}\nolimits_{1,2,\ldots,j}\right)  \cdot\left(
\left[  \cdots\left[  \left[  a_{1},a_{2}\right]  ,a_{3}\right]
,\ldots,a_{j-1}\right]  \underbrace{a_{j}a_{j+1}\cdots a_{n}}_{=a_{j}\left(
a_{j+1}a_{j+2}\cdots a_{n}\right)  }\right) \\
&  =\left(  1-\operatorname*{cyc}\nolimits_{1,2,\ldots,j}\right)  \cdot\left(
\left[  \cdots\left[  \left[  a_{1},a_{2}\right]  ,a_{3}\right]
,\ldots,a_{j-1}\right]  a_{j}\left(  a_{j+1}a_{j+2}\cdots a_{n}\right)
\right) \\
&  =\underbrace{\left[  \left[  \cdots\left[  \left[  a_{1},a_{2}\right]
,a_{3}\right]  ,\ldots,a_{j-1}\right]  ,a_{j}\right]  }_{=\left[
\cdots\left[  \left[  a_{1},a_{2}\right]  ,a_{3}\right]  ,\ldots,a_{j}\right]
}a_{j+1}a_{j+2}\cdots a_{n}\\
&  \ \ \ \ \ \ \ \ \ \ \ \ \ \ \ \ \ \ \ \ \left(
\begin{array}
[c]{c}%
\text{by (\ref{eq.lem.Vperm.Vn.LNC-k.1-c}), applied to }u=\left[
\cdots\left[  \left[  a_{1},a_{2}\right]  ,a_{3}\right]  ,\ldots
,a_{j-1}\right]  \text{, }g=a_{j}\\
\text{and }v=a_{j+1}a_{j+2}\cdots a_{n}\text{ (since a straightforward
induction shows}\\
\text{that }\left[  \cdots\left[  \left[  a_{1},a_{2}\right]  ,a_{3}\right]
,\ldots,a_{j-1}\right]  \text{ is homogeneous of degree }j-1\text{)}%
\end{array}
\right) \\
&  =\left[  \cdots\left[  \left[  a_{1},a_{2}\right]  ,a_{3}\right]
,\ldots,a_{j}\right]  a_{j+1}a_{j+2}\cdots a_{n}.
\end{align*}
Thus, the induction step is complete. Hence, (\ref{pf.prop.Vperm.Vn.LNC.goal})
is proved by induction.

Now, (\ref{pf.thm.Vperm.Vn.Svn=}) yields%
\begin{align*}
&  S\left(  \mathbf{V}_{n}\right)  \cdot\left(  a_{1}a_{2}\cdots a_{n}\right)
\\
&  =\left(  1-\operatorname*{cyc}\nolimits_{1,2,\ldots,n}\right)  \left(
1-\operatorname*{cyc}\nolimits_{1,2,\ldots,n-1}\right)  \cdots\left(
1-\operatorname*{cyc}\nolimits_{1,2}\right)  \cdot\left(  a_{1}a_{2}\cdots
a_{n}\right) \\
&  =\left[  \cdots\left[  \left[  a_{1},a_{2}\right]  ,a_{3}\right]
,\ldots,a_{n}\right]  \underbrace{a_{n+1}a_{n+2}\cdots a_{n}}_{=\left(
\text{empty product}\right)  =1}\ \ \ \ \ \ \ \ \ \ \left(  \text{by
(\ref{pf.prop.Vperm.Vn.LNC.goal}), applied to }j=n\right) \\
&  =\left[  \cdots\left[  \left[  a_{1},a_{2}\right]  ,a_{3}\right]
,\ldots,a_{n}\right]  .
\end{align*}
This proves Proposition \ref{prop.Vperm.Vn.LNC}.
\end{fineprint}
\end{proof}

Let us denote the action of $S\left(  \mathbf{V}_{n}\right)  $ on
$\mathbf{NP}_{k,n}$ by $\mathcal{L}_{n}$. Thus, $\mathcal{L}_{n}$ is the
$\mathbf{k}$-linear map%
\begin{align*}
\mathbf{NP}_{k,n}  &  \rightarrow\mathbf{NP}_{k,n},\\
p  &  \mapsto S\left(  \mathbf{V}_{n}\right)  \cdot p.
\end{align*}
Proposition \ref{prop.Vperm.Vn.LNC} tells us that this map $\mathcal{L}_{n}$
sends each product $a_{1}a_{2}\cdots a_{n}$ (with $a_{1},a_{2},\ldots,a_{n}%
\in\mathbf{NP}_{k,1}$) to $\left[  \cdots\left[  \left[  a_{1},a_{2}\right]
,a_{3}\right]  ,\ldots,a_{n}\right]  $. Thus, this map is sometimes called the
\emph{left-bracketing operator}. And the claimed equality
(\ref{pf.thm.Vperm.Vn.suff3}) -- which we are trying to prove -- becomes%
\begin{equation}
\mathcal{L}_{n}^{2}=n\mathcal{L}_{n}. \label{eq.prop.Vperm.Vn.LNC.LL}%
\end{equation}

This equality was first proved by Franz Wever in 1947 (\cite{Wever47}) by a
rather messy semi-combinatorial argument. Nowadays, the slickest proof of
(\ref{eq.prop.Vperm.Vn.LNC.LL}) uses Hopf algebras, specifically
\cite[Exercise 1.5.14]{GriRei14}:

\begin{proof}
[Proof of (\ref{eq.prop.Vperm.Vn.LNC.LL}).]Let $V$ be the free $\mathbf{k}%
$-module $\mathbf{k}^{k}$. Let $A$ be its tensor algebra $T\left(  V\right)
=\bigoplus\limits_{d\in\mathbb{N}}V^{\otimes d}$. As we know from Remark
\ref{rmk.ncpol.tensor}, we can identify this tensor algebra $A$ with the
noncommutative polynomial ring $\mathbf{NP}_{k}$. In particular, its
degree-$1$ homogeneous component $V^{\otimes1}=V$ is thus $\mathbf{NP}_{k,1}$.

But as a tensor algebra, $A$ also has the structure of a cocommutative graded
Hopf algebra (see \cite[Examples 1.3.11, 1.3.17 and 1.5.3]{GriRei14}). Let us
denote its antipode by $S_{T}$ (note that it is called $S$ in \cite{GriRei14},
but is unrelated to our antipode $S$ of $\mathbf{k}\left[  S_{n}\right]  $).
Let $E:A\rightarrow A$ be the $\mathbf{k}$-linear map that sends each
homogeneous element $a\in A$ to $\left(  \deg a\right)  \cdot a$ (that is,
sends each $a\in\mathbf{NP}_{k,d}$ to $da$)\ \ \ \ \footnote{A more formal
definition of $E$ is as follows: Each element of $A$ can be written uniquely
as $\sum_{d\in\mathbb{N}}a_{d}$, where each $a_{d}$ belongs to $V^{\otimes d}%
$. The map $E$ sends this element to $\sum_{d\in\mathbb{N}}da_{d}$.}. Consider
the $\mathbf{k}$-linear map $S_{T}\star E:A\rightarrow A$, where the symbol
$\star$ denotes convolution (defined in \cite[Definition 1.4.1]{GriRei14}).

\begin{noncompile}
It is easy to see that this map $S_{T}\star E$ is graded (i.e., sends
homogeneous elements of any degree $d$ to homogeneous elements of degree $d$),
since it is the convolution of the two graded maps $S_{T}$ and $E$.
\end{noncompile}

Now, \cite[Exercise 1.5.14 (e)]{GriRei14} shows that%
\begin{equation}
\left(  S_{T}\star E\right)  \left(  a_{1}a_{2}\cdots a_{n}\right)  =\left[
\ldots\left[  \left[  a_{1},a_{2}\right]  ,a_{3}\right]  ,\ldots,a_{n}\right]
\label{pf.eq.prop.Vperm.Vn.LNC.LL.LB}%
\end{equation}
for any $a_{1},a_{2},\ldots,a_{n}\in V=\mathbf{NP}_{k,1}$. Thus, each
$p\in\mathbf{NP}_{k,n}$ satisfies%
\begin{equation}
\left(  S_{T}\star E\right)  \left(  p\right)  =\mathcal{L}_{n}\left(
p\right)  . \label{pf.eq.prop.Vperm.Vn.LNC.LL.STEP}%
\end{equation}

\begin{proof}
[Proof of (\ref{pf.eq.prop.Vperm.Vn.LNC.LL.STEP}).]Both sides of the equality
(\ref{pf.eq.prop.Vperm.Vn.LNC.LL.STEP}) are $\mathbf{k}$-linear in $p$. Thus,
in order to prove it, we can WLOG assume that $p$ is a product $a_{1}%
a_{2}\cdots a_{n}$ of $n$ elements $a_{1},a_{2},\ldots,a_{n}\in\mathbf{NP}%
_{k,1}$ (since such products span the $\mathbf{k}$-module $\mathbf{NP}_{k,n}%
$). Assume this. Then, $p=a_{1}a_{2}\cdots a_{n}$, so that%
\begin{align*}
\left(  S_{T}\star E\right)  \left(  p\right)   &  =\left(  S_{T}\star
E\right)  \left(  a_{1}a_{2}\cdots a_{n}\right) \\
&  =\left[  \ldots\left[  \left[  a_{1},a_{2}\right]  ,a_{3}\right]
,\ldots,a_{n}\right]  \ \ \ \ \ \ \ \ \ \ \left(  \text{by
(\ref{pf.eq.prop.Vperm.Vn.LNC.LL.LB})}\right) \\
&  =S\left(  \mathbf{V}_{n}\right)  \cdot\underbrace{\left(  a_{1}a_{2}\cdots
a_{n}\right)  }_{=p}\ \ \ \ \ \ \ \ \ \ \left(  \text{by Proposition
\ref{prop.Vperm.Vn.LNC}}\right) \\
&  =S\left(  \mathbf{V}_{n}\right)  \cdot p=\mathcal{L}_{n}\left(  p\right)
\ \ \ \ \ \ \ \ \ \ \left(  \text{by the definition of }\mathcal{L}%
_{n}\right)  .
\end{align*}
This proves (\ref{pf.eq.prop.Vperm.Vn.LNC.LL.STEP}).
\end{proof}

Now, we recall some properties of the map $S_{T}\star E$ proved in
\cite{GriRei14}. From \cite[Exercise 1.5.14 (a)]{GriRei14}, it follows that
for every $a\in A$,
\begin{equation}
\text{the element }\left(  S_{T}\star E\right)  \left(  a\right)  \text{ is
primitive.} \label{pf.eq.prop.Vperm.Vn.LNC.LL.prim}%
\end{equation}
But \cite[Exercise 1.5.14 (b)]{GriRei14} yields that
\begin{equation}
\left(  S_{T}\star E\right)  \left(  p\right)  =E\left(  p\right)
\label{pf.eq.prop.Vperm.Vn.LNC.LL.Ep}%
\end{equation}
for each primitive element $p$ of $A$.

Now, let $a\in\mathbf{NP}_{k,n}$. Hence,
(\ref{pf.eq.prop.Vperm.Vn.LNC.LL.STEP}) yields $\left(  S_{T}\star E\right)
\left(  a\right)  =\mathcal{L}_{n}\left(  a\right)  $. But
(\ref{pf.eq.prop.Vperm.Vn.LNC.LL.prim}) shows that the element $\left(
S_{T}\star E\right)  \left(  a\right)  $ is primitive. In other words,
$\mathcal{L}_{n}\left(  a\right)  $ is primitive (since $\left(  S_{T}\star
E\right)  \left(  a\right)  =\mathcal{L}_{n}\left(  a\right)  $). Hence,
(\ref{pf.eq.prop.Vperm.Vn.LNC.LL.Ep}) (applied to $p=\mathcal{L}_{n}\left(
a\right)  $) yields $\left(  S_{T}\star E\right)  \left(  \mathcal{L}%
_{n}\left(  a\right)  \right)  =E\left(  \mathcal{L}_{n}\left(  a\right)
\right)  $. But $\mathcal{L}_{n}\left(  a\right)  $ belongs to $\mathbf{NP}%
_{k,n}$ and thus is a homogeneous element of $A$ of degree $n$. Hence,
$E\left(  \mathcal{L}_{n}\left(  a\right)  \right)  =n\mathcal{L}_{n}\left(
a\right)  $. Furthermore, (\ref{pf.eq.prop.Vperm.Vn.LNC.LL.STEP}) (applied to
$p=\mathcal{L}_{n}\left(  a\right)  $) yields $\left(  S_{T}\star E\right)
\left(  \mathcal{L}_{n}\left(  a\right)  \right)  =\mathcal{L}_{n}\left(
\mathcal{L}_{n}\left(  a\right)  \right)  =\mathcal{L}_{n}^{2}\left(
a\right)  $. Thus,%
\[
\mathcal{L}_{n}^{2}\left(  a\right)  =\left(  S_{T}\star E\right)  \left(
\mathcal{L}_{n}\left(  a\right)  \right)  =E\left(  \mathcal{L}_{n}\left(
a\right)  \right)  =n\mathcal{L}_{n}\left(  a\right)  .
\]

Forget that we fixed $a$. We thus have proved that $\mathcal{L}_{n}^{2}\left(
a\right)  =n\mathcal{L}_{n}\left(  a\right)  $ for each $a\in\mathbf{NP}%
_{k,n}$. Thus, $\mathcal{L}_{n}^{2}=n\mathcal{L}_{n}$. This proves
(\ref{eq.prop.Vperm.Vn.LNC.LL}).
\end{proof}

As we said, we now easily obtain (\ref{pf.thm.Vperm.Vn.suff3}), since each
$p\in\mathbf{NP}_{k,n}$ satisfies%
\begin{align*}
S\left(  \mathbf{V}_{n}\right)  ^{2}\cdot p  &  =S\left(  \mathbf{V}%
_{n}\right)  \cdot\underbrace{S\left(  \mathbf{V}_{n}\right)  \cdot
p}_{\substack{=\mathcal{L}_{n}\left(  p\right)  \\\text{(by the definition of
}\mathcal{L}_{n}\text{)}}}=S\left(  \mathbf{V}_{n}\right)  \cdot
\mathcal{L}_{n}\left(  p\right) \\
&  =\mathcal{L}_{n}\left(  \mathcal{L}_{n}\left(  p\right)  \right)
\ \ \ \ \ \ \ \ \ \ \left(  \text{by the definition of }\mathcal{L}_{n}\right)
\\
&  =\underbrace{\mathcal{L}_{n}^{2}}_{\substack{=n\mathcal{L}_{n}\\\text{(by
(\ref{eq.prop.Vperm.Vn.LNC.LL}))}}}\left(  p\right)  =n\underbrace{\mathcal{L}%
_{n}\left(  p\right)  }_{\substack{=S\left(  \mathbf{V}_{n}\right)  \cdot
p\\\text{(by the definition of }\mathcal{L}_{n}\text{)}}}=nS\left(
\mathbf{V}_{n}\right)  \cdot p.
\end{align*}
So we have proved that (\ref{pf.thm.Vperm.Vn.suff3}) holds for each
$k\in\mathbb{N}$. By picking $k$ to be $\geq n$, we can thus use the
faithfulness of the left $\mathbf{k}\left[  S_{n}\right]  $-module
$\mathbf{NP}_{k,n}$ (Proposition \ref{prop.rep.NP-acts-faith} \textbf{(a)},
applied to $d=n$) to conclude that $S\left(  \mathbf{V}_{n}\right)
^{2}=nS\left(  \mathbf{V}_{n}\right)  $. That is, (\ref{pf.thm.Vperm.Vn.suff2}%
) holds. As we explained, this proves $\mathbf{V}_{n}^{2}=n\mathbf{V}_{n}$,
and thus finishes our proof of Theorem \ref{thm.Vperm.Vn} \textbf{(b)}.

\bigskip

\section{\label{chp.specht}Young and Specht modules}

Now we aim to introduce and study a wider class of representations of $S_{n}$,
which will include the irreducible representations over fields of
characteristic $0$. These are the \emph{Specht modules}, introduced by Wilhelm
Specht in 1935 after Alfred Young planted the seeds in his \textquotedblleft
substitutional analysis\textquotedblright\ (1901--1931)\footnote{Young
published his work in a sequence of nine papers, titled \textquotedblleft%
\textit{On Quantitative Substitutional Analysis}\textquotedblright\ (all
reprinted in \cite{Young77}). These papers contain a phenomenal amount of
ideas, often concealed behind antiquated notations and a terse Victorian
writing style. Many important concepts in this field that are commonly
associated with more recent names (Specht modules, von Neumann sandwich lemma,
Garnir relations) actually appear in some inchoate form in Young's papers.}.
There are several equivalent (up to isomorphism) ways to define Specht modules
(\textquotedblleft avatars\textquotedblright\ as I will call them), and
several levels of generality at which they can be studied (Young diagrams,
skew Young diagrams and arbitrary diagrams). I aim to give a maximally broad
and general approach, although most of the substantial results hold only for
more specialized situations.

\subsection{\label{sec.specht.pars}Partitions, compositions, diagrams}

\subsubsection{Definitions and examples}

To define Specht modules, we will have to define Young modules beforehand.
This, in turn, relies on some basic combinatorics that I will now partly
recall, partly introduce.

We begin by recalling some combinatorial concepts (see, e.g., \cite[\S 3.9 and
Chapter 4]{21s} for more details):

\begin{definition}
\label{def.partitions.partitions}\textbf{(a)} An \emph{integer composition}
(or, for short, \emph{composition}) means a tuple $\left(  a_{1},a_{2}%
,\ldots,a_{k}\right)  $ of positive integers. \medskip

\textbf{(b)} The \emph{size} $\left\vert a\right\vert $ of a composition
$a=\left(  a_{1},a_{2},\ldots,a_{k}\right)  $ is defined to be $a_{1}%
+a_{2}+\cdots+a_{k}\in\mathbb{N}$. \medskip

\textbf{(c)} The \emph{length} $\ell\left(  a\right)  $ of a composition
$a=\left(  a_{1},a_{2},\ldots,a_{k}\right)  $ is defined to be $k$. \medskip

\textbf{(d)} Given $n\in\mathbb{Z}$, a \emph{composition of }$n$ means a
composition with size $n$. \medskip

\textbf{(e)} An \emph{integer partition} (or, for short, \emph{partition})
means an integer composition $\left(  a_{1},a_{2},\ldots,a_{k}\right)  $ that
is weakly decreasing (i.e., satisfies $a_{1}\geq a_{2}\geq\cdots\geq a_{k}$).
\medskip

\textbf{(f)} Given $n\in\mathbb{Z}$, a \emph{partition of }$n$ means a
partition with size $n$. \medskip

\textbf{(g)} We will use the notation $a_{i}$ for the $i$-th entry of a
composition $a$. (This applies, in particular, to all partitions $a$.)
\medskip

\textbf{(h)} If $a$ is a partition, then we furthermore set $a_{i}:=0$ for all
$i>\ell\left(  a\right)  $. That is, we pretend that a partition consists not
only of its original entries, but also contains an infinite row of zeroes
following them. Thus, a partition $a$ satisfies not only the finite chain of
inequalities $a_{1}\geq a_{2}\geq\cdots\geq a_{\ell\left(  a\right)  }$, but
also the infinite chain $a_{1}\geq a_{2}\geq a_{3}\geq\cdots$.

We thus identify each partition $a$ with the infinite sequence $\left(
a_{1},a_{2},a_{3},\ldots\right)  $. From this point of view, a partition can
also be defined as an infinite weakly decreasing sequence of nonnegative
integers that contains only finitely many positive integers. \medskip

\textbf{(i)} Two partitions $a$ and $b$ are said to satisfy $a\supseteq b$
(or, equivalently, $b\subseteq a$) if every positive integer $i$ satisfies
$a_{i}\geq b_{i}$. (It suffices to check this for all $i\leq\ell\left(
b\right)  $, because all $i>\ell\left(  b\right)  $ satisfy $b_{i}=0\leq
a_{i}$ and thus automatically $a_{i}\geq b_{i}$. But it does not suffice to
check this for all $i\leq\ell\left(  a\right)  $.)

If two partitions $a$ and $b$ satisfy $a\supseteq b$, then the pair $\left(
b,a\right)  $ is called a \emph{skew partition} and denoted by $a/b$. \medskip

\textbf{(j)} We denote the empty partition $\left(  {}\right)  =\left(
0,0,0,\ldots\right)  $ by $\varnothing$. \medskip

\textbf{(k)} The notation \textquotedblleft$\lambda\vdash n$\textquotedblright%
\ means \textquotedblleft$\lambda$ is a partition of $n$\textquotedblright.
\end{definition}

For instance:

\begin{itemize}
\item The composition $\left(  2,1,2\right)  $ has size $5$ and length $3$.

\item There are four compositions of $3$, namely $\left(  3\right)  $,
$\left(  2,1\right)  $, $\left(  1,2\right)  $ and $\left(  1,1,1\right)  $.
Only three of them are partitions (indeed, $\left(  1,2\right)  $ is not).

\item The composition $\left(  2,2,1\right)  $ is a partition of $5$. Thus, we
can write $\left(  2,2,1\right)  \vdash5$.

\item The partition $\left(  4,2,2,1\right)  $ is identified with the infinite
sequence $\left(  4,2,2,1,0,0,0,\ldots\right)  $.

\item The partitions $\left(  3,1\right)  $ and $\left(  4,2,1\right)  $
satisfy $\left(  3,1\right)  \subseteq\left(  4,2,1\right)  $ (since $3\leq4$
and $1\leq2$). Thus, the pair $\left(  \left(  3,1\right)  ,\left(
4,2,1\right)  \right)  $ is a skew partition, and is denoted by $\left(
4,2,1\right)  /\left(  3,1\right)  $.

The partitions $\left(  2,2,1\right)  $ and $\left(  3,1,1\right)  $ do not
satisfy $\left(  2,2,1\right)  \subseteq\left(  3,1,1\right)  $, since their
second entries do not satisfy $2\leq1$.

The partitions $\left(  3,1\right)  $ and $\left(  4\right)  $ do not satisfy
$\left(  3,1\right)  \subseteq\left(  4\right)  $, since their second entries
do not satisfy $1\leq0$ (remember that the second entry of $\left(  4\right)
$ is $0$).

\item The empty tuple $\left(  {}\right)  $ is a composition of $0$.
\end{itemize}

Next, we define some \textquotedblleft two-dimensional\textquotedblright%
\ combinatorial objects:

\begin{definition}
\label{def.diagrams.diagrams}\textbf{(a)} The \emph{English coordinate system}
is a Cartesian coordinate system on the plane, in which the x-axis goes
top-down and the y-axis goes left-to-right. (This is exactly the way the
entries of a matrix are indexed: The rows are counted from the top; the
columns are counted from the left.)

(In contrast, the usual Cartesian coordinate system has its x-axis go
left-to-right and its y-axis go bottom-to-top.)

We will be using the English system throughout this topic. \medskip

\textbf{(b)} A \emph{cell} (or \emph{box} or \emph{square}) shall mean a pair
$\left(  i,j\right)  $ of integers. Thus, the set of all cells is
$\mathbb{Z}^{2}=\mathbb{Z}\times\mathbb{Z}$.

We draw each cell $\left(  i,j\right)  $ as a $1\times1$-square in the plane,
centered at the point $\left(  i,j\right)  $ in the English coordinate system.
For now we leave these cells empty; soon we will put numbers into them.
\medskip

\textbf{(c)} We will use the obvious geographical directions (left, right,
above, below, north, east, south, west, northwest, etc.) in reference to the
English coordinate system. For example, a cell $\left(  i,j\right)  $ has the
westward neighbor $\left(  i,j-1\right)  $, the northward neighbor $\left(
i-1,j\right)  $, the eastward neighbor $\left(  i,j+1\right)  $ and the
southward neighbor $\left(  i+1,j\right)  $, and we can also refer to these
four neighbors as the \textquotedblleft left neighbor\textquotedblright, the
\textquotedblleft upward neighbor\textquotedblright, the \textquotedblleft
right neighbor\textquotedblright\ and the \textquotedblleft downward
neighbor\textquotedblright, respectively.

We say that a cell $\left(  i,j\right)  $ lies \emph{weakly northwest} of a
cell $\left(  u,v\right)  $ if $i\leq u$ and $j\leq v$. In this case, we write
$\left(  i,j\right)  \leq\left(  u,v\right)  $. Similarly, we say that a cell
$\left(  i,j\right)  $ lies \emph{weakly northeast} of a cell $\left(
u,v\right)  $ if $i\leq u$ and $j\geq v$. (But we do not have a symbol for
this.) Generally, the words \textquotedblleft weakly\textquotedblright\ and
\textquotedblleft strictly\textquotedblright\ will typically refer to weak and
strict inequalities (so that a cell lies weakly north of itself but not
strictly north of itself, etc.). \medskip

\textbf{(d)} A \emph{diagram} (or \emph{shape}) will mean a finite set of
cells. For instance, the set%
\[
\left\{  \left(  1,2\right)  ,\ \left(  2,1\right)  ,\ \left(  2,3\right)
,\ \left(  1,4\right)  ,\ \left(  1,6\right)  ,\ \left(  2,6\right)
,\ \left(  3,6\right)  \right\}
\]
is a diagram. In the English coordinate system, it is drawn as follows:%
\[%
%TCIMACRO{\TeXButton{tikz turtle-and-wall diagram}{\begin{tikzpicture}%
%[scale=0.7]
%\draw[fill=red!50] (0, 0) rectangle (1, 1);
%\draw[fill=red!50] (1, 1) rectangle (2, 2);
%\draw[fill=red!50] (2, 0) rectangle (3, 1);
%\draw[fill=red!50] (3, 1) rectangle (4, 2);
%\draw[fill=red!50] (5, 0) rectangle (6, 1);
%\draw[fill=red!50] (5, 1) rectangle (6, 2);
%\draw[fill=red!50] (5, -1) rectangle (6, 0);
%\end{tikzpicture}}}%
%BeginExpansion
\begin{tikzpicture}[scale=0.7]
\draw[fill=red!50] (0, 0) rectangle (1, 1);
\draw[fill=red!50] (1, 1) rectangle (2, 2);
\draw[fill=red!50] (2, 0) rectangle (3, 1);
\draw[fill=red!50] (3, 1) rectangle (4, 2);
\draw[fill=red!50] (5, 0) rectangle (6, 1);
\draw[fill=red!50] (5, 1) rectangle (6, 2);
\draw[fill=red!50] (5, -1) rectangle (6, 0);
\end{tikzpicture}%
%EndExpansion
\ \ .
\]
(Beware that, if the axes are not shown, such a drawing determines the diagram
only up to parallel translation. But this will usually suffice for us.)

We always draw cells and diagrams in the English coordinate system. This is
called the \emph{English notation} for diagrams. (There is an alternative
notation used in parts of the literature, known as \emph{French notation},
which has the x-axis go bottom-to-top rather than top-to-bottom.)

The \emph{cells} of a diagram $D$ mean the cells that belong to $D$. \medskip

\textbf{(e)} The \emph{rows} and the \emph{columns} are defined in the way you
would expect:

\begin{itemize}
\item For any given $i\in\mathbb{Z}$, we define the $i$\emph{-th row} (aka
\emph{row }$i$) to be the set of all cells of the form $\left(  i,j\right)  $
for $j\in\mathbb{Z}$.

\item For any diagram $D$, the $i$\emph{-th row of }$D$ (aka \emph{row }%
$i$\emph{ of }$D$) is the set of all such cells in $D$ (that is, the
intersection of $D$ with the $i$-th row).

\item For any given $j\in\mathbb{Z}$, we define the $j$\emph{-th column} (aka
\emph{column }$j$) to be the set of all cells of the form $\left(  i,j\right)
$ for $i\in\mathbb{Z}$.

\item For any diagram $D$, the $j$\emph{-th column of }$D$ (aka \emph{column
}$j$\emph{ of }$D$) is the set of all such cells in $D$ (that is, the
intersection of $D$ with the $j$-th column).
\end{itemize}

\textbf{(f)} The \emph{Young diagram} $Y\left(  a\right)  $ of a partition
$a=\left(  a_{1},a_{2},\ldots,a_{k}\right)  $ is defined to be the diagram%
\begin{align*}
&  \left\{  \left(  i,j\right)  \in\left\{  1,2,3,\ldots\right\}  ^{2}%
\ \mid\ j\leq a_{i}\right\} \\
&  =\left\{  \left(  i,j\right)  \in\left\{  1,2,\ldots,k\right\}
\times\left\{  1,2,3,\ldots\right\}  \ \mid\ j\leq a_{i}\right\} \\
&  \ \ \ \ \ \ \ \ \ \ \ \ \ \ \ \ \ \ \ \ \left(
\begin{array}
[c]{c}%
\text{since }j\leq a_{i}\text{ can happen only if }a_{i}>0\text{,}\\
\text{which entails that }i\in\left\{  1,2,\ldots,k\right\}
\end{array}
\right) \\
&  =\bigcup\limits_{i=1}^{k}\left(  \left\{  i\right\}  \times\left[
a_{i}\right]  \right)  \ \ \ \ \ \ \ \ \ \ \left(
\begin{array}
[c]{c}%
\text{since a positive integer }j\in\left\{  1,2,3,\ldots\right\} \\
\text{satisfies }j\leq a_{i}\text{ if and only if }j\in\left[  a_{i}\right]
\end{array}
\right) \\
&  =\bigcup\limits_{i\geq1}\left(  \left\{  i\right\}  \times\left[
a_{i}\right]  \right)  \ \ \ \ \ \ \ \ \ \ \left(  \text{since }\left\{
i\right\}  \times\left[  a_{i}\right]  \text{ is empty whenever }i>k\right)  .
\end{align*}
Visually, this is a collection of cells that forms a \textquotedblleft table
with $k$ left-aligned rows\textquotedblright, where row $i$ has $a_{i}$ cells
for each $i\in\left[  k\right]  $.

For example, the Young diagram $Y\left(  3,2\right)  $ of the partition
$\left(  3,2\right)  $ is the set
\[
\left\{  \left(  1,1\right)  ,\ \left(  1,2\right)  ,\ \left(  1,3\right)
,\ \left(  2,1\right)  ,\ \left(  2,2\right)  \right\}
\]
and is drawn as $%
%TCIMACRO{\TeXButton{tikz Y522}{\begin{tikzpicture}[scale=0.7]
%\draw[fill=red!50] (0, 0) rectangle (1, 1);
%\draw[fill=red!50] (1, 0) rectangle (2, 1);
%\draw[fill=red!50] (0, 1) rectangle (1, 2);
%\draw[fill=red!50] (1, 1) rectangle (2, 2);
%\draw[fill=red!50] (2, 1) rectangle (3, 2);
%\end{tikzpicture}}}%
%BeginExpansion
\begin{tikzpicture}[scale=0.7]
\draw[fill=red!50] (0, 0) rectangle (1, 1);
\draw[fill=red!50] (1, 0) rectangle (2, 1);
\draw[fill=red!50] (0, 1) rectangle (1, 2);
\draw[fill=red!50] (1, 1) rectangle (2, 2);
\draw[fill=red!50] (2, 1) rectangle (3, 2);
\end{tikzpicture}%
%EndExpansion
$ .

For some more examples, the Young diagrams $Y\left(  5,2,2\right)  $ and
$Y\left(  3,3\right)  $ of the partitions $\left(  5,2,2\right)  $ and
$\left(  3,3\right)  $ are drawn below:%
\[%
\begin{tabular}
[c]{|cc|}\hline
\multicolumn{1}{|c|}{$Y\left(  5,2,2\right)  \vphantom{\int_g^g}$} & $Y\left(
3,3\right)  $\\
\multicolumn{1}{|c|}{$\phantom{\int_g^g}%
%TCIMACRO{\TeXButton{tikz Y522}{\begin{tikzpicture}[scale=0.7]
%\draw[fill=red!50] (0, 0) rectangle (1, 1);
%\draw[fill=red!50] (1, 0) rectangle (2, 1);
%\draw[fill=red!50] (1, 1) rectangle (2, 2);
%\draw[fill=red!50] (0, 1) rectangle (1, 2);
%\draw[fill=red!50] (0, 2) rectangle (1, 3);
%\draw[fill=red!50] (1, 2) rectangle (2, 3);
%\draw[fill=red!50] (2, 2) rectangle (3, 3);
%\draw[fill=red!50] (3, 2) rectangle (4, 3);
%\draw[fill=red!50] (4, 2) rectangle (5, 3);
%\end{tikzpicture}}}%
%BeginExpansion
\begin{tikzpicture}[scale=0.7]
\draw[fill=red!50] (0, 0) rectangle (1, 1);
\draw[fill=red!50] (1, 0) rectangle (2, 1);
\draw[fill=red!50] (1, 1) rectangle (2, 2);
\draw[fill=red!50] (0, 1) rectangle (1, 2);
\draw[fill=red!50] (0, 2) rectangle (1, 3);
\draw[fill=red!50] (1, 2) rectangle (2, 3);
\draw[fill=red!50] (2, 2) rectangle (3, 3);
\draw[fill=red!50] (3, 2) rectangle (4, 3);
\draw[fill=red!50] (4, 2) rectangle (5, 3);
\end{tikzpicture}%
%EndExpansion
\phantom{\int_g^g}$} & $\phantom{\int_g^g}%
%TCIMACRO{\TeXButton{tikz Y33}{\begin{tikzpicture}[scale=0.7]
%\draw[fill=red!50] (0, 0) rectangle (1, 1);
%\draw[fill=red!50] (1, 0) rectangle (3, 1);
%\draw[fill=red!50] (2, 0) rectangle (3, 1);
%\draw[fill=red!50] (1, 0) rectangle (2, 1);
%\draw[fill=red!50] (0, 1) rectangle (1, 2);
%\draw[fill=red!50] (1, 1) rectangle (2, 2);
%\draw[fill=red!50] (2, 1) rectangle (3, 2);
%\end{tikzpicture}}}%
%BeginExpansion
\begin{tikzpicture}[scale=0.7]
\draw[fill=red!50] (0, 0) rectangle (1, 1);
\draw[fill=red!50] (1, 0) rectangle (3, 1);
\draw[fill=red!50] (2, 0) rectangle (3, 1);
\draw[fill=red!50] (1, 0) rectangle (2, 1);
\draw[fill=red!50] (0, 1) rectangle (1, 2);
\draw[fill=red!50] (1, 1) rectangle (2, 2);
\draw[fill=red!50] (2, 1) rectangle (3, 2);
\end{tikzpicture}%
%EndExpansion
\phantom{\int_g^g}$\\\hline
\end{tabular}
\ \ \ \ .
\]

\textbf{(g)} The \emph{skew Young diagram} $Y\left(  a/b\right)  $ of a skew
partition $a/b$ is defined to be the diagram%
\[
Y\left(  a\right)  \setminus Y\left(  b\right)  =\left\{  \left(  i,j\right)
\in\left\{  1,2,3,\ldots\right\}  ^{2}\ \mid\ b_{i}<j\leq a_{i}\right\}  .
\]

For example, the skew Young diagram $Y\left(  \left(  4,3,3,2,1\right)
/\left(  2,2,1,1\right)  \right)  $ is drawn below:%
\[%
%TCIMACRO{\TeXButton{tikz Yskew}{\begin{tikzpicture}[scale=0.7]
%\draw[fill=red!50] (0, 0) rectangle (1, 1);
%\draw[fill=red!50] (1, 1) rectangle (2, 2);
%\draw[fill=red!50] (1, 2) rectangle (2, 3);
%\draw[fill=red!50] (2, 2) rectangle (3, 3);
%\draw[fill=red!50] (2, 3) rectangle (3, 4);
%\draw[fill=red!50] (2, 4) rectangle (3, 5);
%\draw[fill=red!50] (3, 4) rectangle (4, 5);
%\end{tikzpicture}}}%
%BeginExpansion
\begin{tikzpicture}[scale=0.7]
\draw[fill=red!50] (0, 0) rectangle (1, 1);
\draw[fill=red!50] (1, 1) rectangle (2, 2);
\draw[fill=red!50] (1, 2) rectangle (2, 3);
\draw[fill=red!50] (2, 2) rectangle (3, 3);
\draw[fill=red!50] (2, 3) rectangle (3, 4);
\draw[fill=red!50] (2, 4) rectangle (3, 5);
\draw[fill=red!50] (3, 4) rectangle (4, 5);
\end{tikzpicture}%
%EndExpansion
\ \ .
\]

\end{definition}

As the above definitions of $Y\left(  a\right)  $ and $Y\left(  a/b\right)  $
suggest, we will mostly use cells $\left(  i,j\right)  $ whose coordinates $i$
and $j$ are positive, although we do not formally require this.

We note that whenever $a$ is a partition, the pair $a/\varnothing=\left(
\varnothing,a\right)  $ is a skew partition (i.e., we have $\varnothing
\subseteq a$), and its skew Young diagram $Y\left(  a/\varnothing\right)  $ is
just the usual Young diagram $Y\left(  a\right)  $. Thus, skew Young diagrams
generalize Young diagrams.

Here is a little table of Young diagrams:

\begin{itemize}
\item The only partition of $0$ is the empty partition $\varnothing=\left(
{}\right)  $, and its Young diagram $Y\left(  \varnothing\right)  $ is empty
(i.e., has no cells).

\item The only partition of $1$ is $\left(  1\right)  $, and its Young diagram
$Y\left(  1\right)  $ is a single cell: $\begin{tikzpicture}[scale=0.7]
\draw[fill=red!50] (0, 0) rectangle (1, 1);
\end{tikzpicture}$ .

\item The partitions of $2$ are $\left(  2\right)  $ and $\left(  1,1\right)
$, and their Young diagrams are
\[
Y\left(  2\right)  =\begin{tikzpicture}[scale=0.7]
\draw[fill=red!50] (0, 0) rectangle (1, 1);
\draw[fill=red!50] (1, 0) rectangle (2, 1);
\end{tikzpicture}\ \ \ \ \ \ \ \ \ \ \text{and}\ \ \ \ \ \ \ \ \ \ Y\left(
1,1\right)  =\begin{tikzpicture}[scale=0.7]
\draw[fill=red!50] (0, 0) rectangle (1, 1);
\draw[fill=red!50] (0, 1) rectangle (1, 2);
\end{tikzpicture}.
\]

\item Here are the partitions of $3$ and their Young diagrams:%
\[%
\begin{tabular}
[c]{|c||c|c|c|}\hline
$\lambda$ & $\left(  3\right)  $ & $\left(  2,1\right)  $ & $\left(
1,1,1\right)  $\\\hline
$Y\left(  \lambda\right)  $ &
$\begin{tikzpicture}[scale=0.7] \draw[fill=red!50] (0, 0) rectangle (1, 1); \draw[fill=red!50] (1, 0) rectangle (2, 1); \draw[fill=red!50] (2, 0) rectangle (3, 1); \end{tikzpicture}$
&
$\begin{tikzpicture}[scale=0.7] \draw[fill=red!50] (0, 0) rectangle (1, 1); \draw[fill=red!50] (0, 1) rectangle (1, 2); \draw[fill=red!50] (1, 1) rectangle (2, 2); \end{tikzpicture}$
&
$\begin{tikzpicture}[scale=0.7] \draw[fill=red!50] (0, 0) rectangle (1, 1); \draw[fill=red!50] (0, 1) rectangle (1, 2); \draw[fill=red!50] (0, 2) rectangle (1, 3); \end{tikzpicture}$%
\\\hline
\end{tabular}
\]

\item Here are the partitions of $4$ and their Young diagrams:%
\[%
\begin{tabular}
[c]{|c||c|c|c|c|c|}\hline
$\lambda$ & $\left(  4\right)  $ & $\left(  3,1\right)  $ & $\left(
2,2\right)  $ & $\left(  2,1,1\right)  $ & $\left(  1,1,1,1\right)  $\\\hline
$Y\left(  \lambda\right)  $ &
$\begin{tikzpicture}[scale=0.7] \draw[fill=red!50] (0, 0) rectangle (1, 1); \draw[fill=red!50] (1, 0) rectangle (2, 1); \draw[fill=red!50] (2, 0) rectangle (3, 1); \draw[fill=red!50] (3, 0) rectangle (4, 1); \end{tikzpicture}$
&
$\begin{tikzpicture}[scale=0.7] \draw[fill=red!50] (0, 0) rectangle (1, 1); \draw[fill=red!50] (0, 1) rectangle (1, 2); \draw[fill=red!50] (1, 1) rectangle (2, 2); \draw[fill=red!50] (2, 1) rectangle (3, 2); \end{tikzpicture}$
&
$\begin{tikzpicture}[scale=0.7] \draw[fill=red!50] (0, 0) rectangle (1, 1); \draw[fill=red!50] (0, 1) rectangle (1, 2); \draw[fill=red!50] (1, 0) rectangle (2, 1); \draw[fill=red!50] (1, 1) rectangle (2, 2); \end{tikzpicture}$
&
$\begin{tikzpicture}[scale=0.7] \draw[fill=red!50] (0, 0) rectangle (1, 1); \draw[fill=red!50] (0, 1) rectangle (1, 2); \draw[fill=red!50] (0, 2) rectangle (1, 3); \draw[fill=red!50] (1, 2) rectangle (2, 3); \end{tikzpicture}$
&
$\begin{tikzpicture}[scale=0.7] \draw[fill=red!50] (0, 0) rectangle (1, 1); \draw[fill=red!50] (0, 1) rectangle (1, 2); \draw[fill=red!50] (0, 2) rectangle (1, 3); \draw[fill=red!50] (0, 3) rectangle (1, 4); \end{tikzpicture}$%
\\\hline
\end{tabular}
\ \
\]

\end{itemize}

\subsubsection{Basic properties}

Let us prove some basic properties of Young diagrams as a warmup. First, let
us convince ourselves that the relation $\subseteq$ on partitions deserves its
\textquotedblleft subset\textquotedblright\ symbol:

\begin{proposition}
\label{prop.young.b-in-a}Let $a$ and $b$ be two partitions. Then, $b\subseteq
a$ holds if and only if $Y\left(  b\right)  \subseteq Y\left(  a\right)  $.
\end{proposition}

\begin{fineprint}
\begin{proof}
Recall that the definition of $Y\left(  a\right)  $ says that%
\begin{equation}
Y\left(  a\right)  =\left\{  \left(  i,j\right)  \in\left\{  1,2,3,\ldots
\right\}  ^{2}\ \mid\ j\leq a_{i}\right\}  . \label{pf.prop.young.b-in-a.Ya=}%
\end{equation}
Similarly,%
\begin{equation}
Y\left(  b\right)  =\left\{  \left(  i,j\right)  \in\left\{  1,2,3,\ldots
\right\}  ^{2}\ \mid\ j\leq b_{i}\right\}  . \label{pf.prop.young.b-in-a.Yb=}%
\end{equation}

We need to prove the equivalence $\left(  b\subseteq a\right)
\ \Longleftrightarrow\ \left(  Y\left(  b\right)  \subseteq Y\left(  a\right)
\right)  $. We shall prove its \textquotedblleft$\Longrightarrow
$\textquotedblright\ and \textquotedblleft$\Longleftarrow$\textquotedblright%
\ implications separately: \medskip

$\Longrightarrow:$ Assume that $b\subseteq a$. We must prove that $Y\left(
b\right)  \subseteq Y\left(  a\right)  $.

Let $\left(  i,j\right)  \in Y\left(  b\right)  $ be arbitrary. By
(\ref{pf.prop.young.b-in-a.Yb=}), this entails $i\geq1$ and $j\leq b_{i}$, so
that $j\leq b_{i}\leq a_{i}$ (since $b\subseteq a$ means that $b_{i}\leq
a_{i}$ for each $i\geq1$). Hence, $i\geq1$ and $j\leq a_{i}$. This, in turn,
shows that $\left(  i,j\right)  \in Y\left(  a\right)  $ (by
(\ref{pf.prop.young.b-in-a.Ya=})).

Forget that we fixed $\left(  i,j\right)  $. We thus have shown that $\left(
i,j\right)  \in Y\left(  a\right)  $ for each $\left(  i,j\right)  \in
Y\left(  b\right)  $. In other words, $Y\left(  b\right)  \subseteq Y\left(
a\right)  $. This proves the \textquotedblleft$\Longrightarrow$%
\textquotedblright\ direction of the equivalence $\left(  b\subseteq a\right)
\ \Longleftrightarrow\ \left(  Y\left(  b\right)  \subseteq Y\left(  a\right)
\right)  $. \medskip

$\Longleftarrow:$ Assume that $Y\left(  b\right)  \subseteq Y\left(  a\right)
$. We must prove that $b\subseteq a$. In other words, we must prove that
$b_{i}\leq a_{i}$ for each $i\geq1$.

Let $i\geq1$. We must prove that $b_{i}\leq a_{i}$. Assume the contrary. Thus,
$b_{i}>a_{i}\geq0$, so that $b_{i}\in\left\{  1,2,3,\ldots\right\}  $. By
(\ref{pf.prop.young.b-in-a.Yb=}), this entails that $\left(  i,b_{i}\right)
\in Y\left(  b\right)  $ (since $b_{i}\leq b_{i}$). Hence, $\left(
i,b_{i}\right)  \in Y\left(  b\right)  \subseteq Y\left(  a\right)  $. By
(\ref{pf.prop.young.b-in-a.Ya=}), this shows that $b_{i}\leq a_{i}$, which
contradicts $b_{i}>a_{i}$. This contradiction shows that our assumption was
false. Hence, the proof of the \textquotedblleft$\Longleftarrow$%
\textquotedblright\ direction of the equivalence $\left(  b\subseteq a\right)
\ \Longleftrightarrow\ \left(  Y\left(  b\right)  \subseteq Y\left(  a\right)
\right)  $ is complete.

We have now proved both directions of the equivalence; thus, the equivalence
(i.e., Proposition \ref{prop.young.b-in-a}) holds.
\end{proof}
\end{fineprint}

Next, we show that the size (i.e., number of cells) of the Young diagram
$Y\left(  \lambda\right)  $ of a partition $\lambda$ is just the size
$\left\vert \lambda\right\vert $ of $\lambda$:

\begin{proposition}
\label{prop.young.n}Let $\lambda$ be any partition. Then, its Young diagram
$Y\left(  \lambda\right)  $ has size $\left\vert Y\left(  \lambda\right)
\right\vert =\left\vert \lambda\right\vert $.
\end{proposition}

\begin{fineprint}
\begin{proof}
Write $\lambda$ as $\lambda=\left(  \lambda_{1},\lambda_{2},\ldots,\lambda
_{k}\right)  $. Then, the definition of the size $\left\vert \lambda
\right\vert $ shows that $\left\vert \lambda\right\vert =\lambda_{1}%
+\lambda_{2}+\cdots+\lambda_{k}$.

However, the definition of the Young diagram $Y\left(  \lambda\right)  $
yields $Y\left(  \lambda\right)  =\bigcup\limits_{i=1}^{k}\left(  \left\{
i\right\}  \times\left[  \lambda_{i}\right]  \right)  $. Thus,%
\begin{align*}
\left\vert Y\left(  \lambda\right)  \right\vert  &  =\left\vert \bigcup
\limits_{i=1}^{k}\left(  \left\{  i\right\}  \times\left[  \lambda_{i}\right]
\right)  \right\vert \\
&  =\sum_{i=1}^{k}\underbrace{\left\vert \left(  \left\{  i\right\}
\times\left[  \lambda_{i}\right]  \right)  \right\vert }_{=\left\vert \left\{
i\right\}  \right\vert \cdot\left\vert \left[  \lambda_{i}\right]  \right\vert
}\ \ \ \ \ \ \ \ \ \ \left(
\begin{array}
[c]{c}%
\text{by the sum rule, since the sets }\left\{  i\right\}  \times\left[
\lambda_{i}\right]  \text{ for}\\
\text{all }i\in\left[  k\right]  \text{ are distinct (because any element}\\
\text{of }\left\{  i\right\}  \times\left[  \lambda_{i}\right]  \text{ is a
pair whose first entry is }i\text{)}%
\end{array}
\right) \\
&  =\sum_{i=1}^{k}\underbrace{\left\vert \left\{  i\right\}  \right\vert
}_{=1}\cdot\underbrace{\left\vert \left[  \lambda_{i}\right]  \right\vert
}_{=\lambda_{i}}=\sum_{i=1}^{k}\lambda_{i}=\lambda_{1}+\lambda_{2}%
+\cdots+\lambda_{k}=\left\vert \lambda\right\vert
\end{align*}
(since $\left\vert \lambda\right\vert =\lambda_{1}+\lambda_{2}+\cdots
+\lambda_{k}$). This proves Proposition \ref{prop.young.n}.
\end{proof}
\end{fineprint}

\begin{proposition}
\label{prop.young.NE}Let $\lambda$ be any partition. Then, its Young diagram
$Y\left(  \lambda\right)  $ has the following property: If a cell
$c\in\left\{  1,2,3,\ldots\right\}  ^{2}$ lies weakly northwest of some cell
$d\in Y\left(  \lambda\right)  $, then we have $c\in Y\left(  \lambda\right)
$ as well.
\end{proposition}

\begin{fineprint}
\begin{proof}
Let $c\in\left\{  1,2,3,\ldots\right\}  ^{2}$ be a cell that lies weakly
northwest of some cell $d\in Y\left(  \lambda\right)  $. We must show that
$c\in Y\left(  \lambda\right)  $.

The definition of $Y\left(  \lambda\right)  $ yields that%
\begin{equation}
Y\left(  \lambda\right)  =\left\{  \left(  i,j\right)  \in\left\{
1,2,3,\ldots\right\}  ^{2}\ \mid\ j\leq\lambda_{i}\right\}  .
\label{pf.prop.young.NE.1}%
\end{equation}

Write the cells $c$ and $d$ as $c=\left(  i,j\right)  $ and $d=\left(
u,v\right)  $. Since $c$ lies weakly northwest of $d$, we thus have $i\leq u$
and $j\leq v$. Moreover, $i$ and $j$ are positive integers (since $\left(
i,j\right)  =c\in\left\{  1,2,3,\ldots\right\}  ^{2}$).

But $\left(  u,v\right)  =d\in Y\left(  \lambda\right)  $. Hence,
$u,v\in\left\{  1,2,3,\ldots\right\}  $ and $v\leq\lambda_{u}$ (by
(\ref{pf.prop.young.NE.1})).

However, the sequence $\left(  \lambda_{1},\lambda_{2},\lambda_{3}%
,\ldots\right)  $ is weakly decreasing (since $\lambda$ is a partition). Thus,
from $i\leq u$, we obtain $\lambda_{i}\geq\lambda_{u}$. Now, $j\leq
v\leq\lambda_{u}\leq\lambda_{i}$ (since $\lambda_{i}\geq\lambda_{u}$). Since
$i$ and $j$ are positive integers, this entails that $\left(  i,j\right)  \in
Y\left(  \lambda\right)  $ (by (\ref{pf.prop.young.NE.1})). In other words,
$c\in Y\left(  \lambda\right)  $ (since $c=\left(  i,j\right)  $). This proves
Proposition \ref{prop.young.NE}.
\end{proof}
\end{fineprint}

\begin{proposition}
\label{prop.young.convexity}Let $\lambda/\mu$ be any skew partition. Then, its
skew Young diagram $Y\left(  \lambda/\mu\right)  $ has the following property:
If $c,d,e\in\mathbb{Z}^{2}$ are three cells satisfying $c\leq d\leq e$ and
$c,e\in Y\left(  \lambda/\mu\right)  $, then $d\in Y\left(  \lambda
/\mu\right)  $.
\end{proposition}

\begin{fineprint}
\begin{proof}
Let $c,d,e\in\mathbb{Z}^{2}$ be three cells satisfying $c\leq d\leq e$ and
$c,e\in Y\left(  \lambda/\mu\right)  $. We must show that $d\in Y\left(
\lambda/\mu\right)  $.

We have $Y\left(  \lambda/\mu\right)  =Y\left(  \lambda\right)  \setminus
Y\left(  \mu\right)  $ (by the definition of a skew Young diagram). Now, $c\in
Y\left(  \lambda/\mu\right)  =Y\left(  \lambda\right)  \setminus Y\left(
\mu\right)  $, so that $c\notin Y\left(  \mu\right)  $. Furthermore, $e\in
Y\left(  \lambda/\mu\right)  =Y\left(  \lambda\right)  \setminus Y\left(
\mu\right)  $, so that $e\in Y\left(  \lambda\right)  $.

Moreover, $c\in Y\left(  \lambda\right)  \setminus Y\left(  \mu\right)
\subseteq Y\left(  \lambda\right)  \subseteq\left\{  1,2,3,\ldots\right\}
^{2}$. Thus, because of $c\leq d$, it is easy to see that $d\in\left\{
1,2,3,\ldots\right\}  ^{2}$ as well\footnote{\textit{Proof.} Write the cells
$c$ and $d$ as $c=\left(  i,j\right)  $ and $d=\left(  i^{\prime},j^{\prime
}\right)  $. Then, $\left(  i,j\right)  =c\leq d=\left(  i^{\prime},j^{\prime
}\right)  $; in other words, $i\leq i^{\prime}$ and $j\leq j^{\prime}$. But
$\left(  i,j\right)  =c\in\left\{  1,2,3,\ldots\right\}  ^{2}$, so that
$i\geq1$ and $j\geq1$. Thus, from $i\leq i^{\prime}$, we obtain $i^{\prime
}\geq i\geq1$ and thus $i^{\prime}\in\left\{  1,2,3,\ldots\right\}  $, and
similarly $j^{\prime}\in\left\{  1,2,3,\ldots\right\}  $. Hence, $\left(
i^{\prime},j^{\prime}\right)  \in\left\{  1,2,3,\ldots\right\}  ^{2}$, so that
$d=\left(  i^{\prime},j^{\prime}\right)  \in\left\{  1,2,3,\ldots\right\}
^{2}$.}.

The cell $c$ lies weakly northwest of $d$ (since $c\leq d$). Hence, if we had
$d\in Y\left(  \mu\right)  $, then we would have $c\in Y\left(  \mu\right)  $
(by Proposition \ref{prop.young.NE}, applied to $\mu$ instead of $\lambda$),
which would contradict $c\notin Y\left(  \mu\right)  $. Thus, we cannot have
$d\in Y\left(  \mu\right)  $. In other words, $d\notin Y\left(  \mu\right)  $.

The cell $d$ lies weakly northwest of $e$ (since $d\leq e$). Hence, from $e\in
Y\left(  \lambda\right)  $, we obtain $d\in Y\left(  \lambda\right)  $ (by
Proposition \ref{prop.young.NE}, applied to $d$ and $e$ instead of $c$ and $d$).

Combining $d\in Y\left(  \lambda\right)  $ with $d\notin Y\left(  \mu\right)
$, we obtain $d\in Y\left(  \lambda\right)  \setminus Y\left(  \mu\right)
=Y\left(  \lambda/\mu\right)  $. This proves Proposition
\ref{prop.young.convexity}.
\end{proof}
\end{fineprint}

\begin{proposition}
\label{prop.young.convex-etc}Let $D\subseteq\left\{  1,2,3,\ldots\right\}
^{2}$ be a diagram. Then: \medskip

\textbf{(a)} The diagram $D$ is the Young diagram $Y\left(  \lambda\right)  $
of some partition $\lambda$ if and only if it has the following property: If a
cell $c\in\left\{  1,2,3,\ldots\right\}  ^{2}$ lies weakly northwest of a cell
$d\in D$, then we have $c\in D$ as well. \medskip

\textbf{(b)} The diagram $D$ is the skew Young diagram $Y\left(  \lambda
/\mu\right)  $ of some skew partition $\lambda/\mu$ if and only if it has the
following property: If $c,d,e\in\left\{  1,2,3,\ldots\right\}  ^{2}$ are three
cells satisfying $c\leq d\leq e$ and $c,e\in D$, then $d\in D$. \medskip

\textbf{(c)} A partition $\lambda$ is uniquely determined by its Young diagram
$Y\left(  \lambda\right)  $. \medskip

\textbf{(d)} A skew partition $\lambda/\mu$, however, is not always uniquely
determined by its skew Young diagram $Y\left(  \lambda/\mu\right)  $.
\end{proposition}

\begin{fineprint}
\begin{proof}
[Partial proof.]\textbf{(a)} $\Longrightarrow:$ This follows from Proposition
\ref{prop.young.NE}.

$\Longleftarrow:$ Assume that $D$ has the following property: If a cell
$c\in\left\{  1,2,3,\ldots\right\}  ^{2}$ lies weakly northwest of a cell
$d\in D$, then we have $c\in D$ as well. We shall refer to this property as
the \textquotedblleft northwest-closure property\textquotedblright.

We must show that $D$ is the Young diagram $Y\left(  \lambda\right)  $ of some
partition $\lambda$.

If $D$ is empty, then this is clear (since the empty diagram is $Y\left(
\varnothing\right)  $ for the empty partition $\varnothing=\left(  {}\right)
$). Thus, we WLOG assume that $D$ is nonempty.

Let $\left(  p,q\right)  $ be a southernmost cell in $D$ (that is, a cell with
maximum $p$). This is well-defined, since $D$ is finite and nonempty. By the
northwest-closure property, it follows that all $p$ cells $\left(  1,q\right)
,\ \left(  2,q\right)  ,\ \ldots,\ \left(  p,q\right)  $ belong to $D$ as well
(since they lie weakly northwest of the cell $\left(  p,q\right)  \in D$).
Hence, the diagram $D$ has cells in each of rows $1,2,\ldots,p$ (namely, the
$p$ cells we just mentioned). Furthermore, $D$ has no cells south of row $p$
(since $\left(  p,q\right)  $ is a southernmost cell in $D$).

For each $i\in\left[  p\right]  $, let $\left(  i,m_{i}\right)  $ be the
easternmost cell in row $i$ of $D$ (that is, the cell in row $i$ of $D$ having
maximum $m_{i}$). This is well-defined (since $D$ is finite and has at least
one cell in row $i$ (because the diagram $D$ has cells in each of rows
$1,2,\ldots,p$)). Thus, we have defined $p$ positive integers $m_{1}%
,m_{2},\ldots,m_{p}$, which form a composition $\left(  m_{1},m_{2}%
,\ldots,m_{p}\right)  $.

We next claim that this composition $\left(  m_{1},m_{2},\ldots,m_{p}\right)
$ is a partition.

[\textit{Proof:} Let $i\in\left\{  2,3,\ldots,p\right\}  $. We shall prove
that $m_{i-1}\geq m_{i}$. Indeed, $\left(  i,m_{i}\right)  $ is a cell in $D$
(by the definition of $m_{i}$). Hence, its northward neighbor $\left(
i-1,m_{i}\right)  $ also belongs to $D$ (by the northwest-closure property,
since $i-1\geq1$). But $\left(  i-1,m_{i-1}\right)  $ is defined as the
easternmost cell in row $i-1$ of $D$. Thus, each cell in row $i-1$ of $D$ must
lie weakly west of this cell $\left(  i-1,m_{i-1}\right)  $. Applying this to
the cell $\left(  i-1,m_{i}\right)  $ (which clearly lies in row $i-1$ and
belongs to $D$), we conclude that $\left(  i-1,m_{i}\right)  $ lies weakly
west of the cell $\left(  i-1,m_{i-1}\right)  $. In other words, $m_{i}\leq
m_{i-1}$. In other words, $m_{i-1}\geq m_{i}$.

Forget that we fixed $i$. We thus have shown that $m_{i-1}\geq m_{i}$ for each
$i\in\left\{  2,3,\ldots,p\right\}  $. Combining these inequalities, we arrive
at $m_{1}\geq m_{2}\geq\cdots\geq m_{p}$. In other words, the composition
$\left(  m_{1},m_{2},\ldots,m_{p}\right)  $ is a partition.]

The definition of the Young diagram $Y\left(  m_{1},m_{2},\ldots,m_{p}\right)
$ yields%
\begin{equation}
Y\left(  m_{1},m_{2},\ldots,m_{p}\right)  =\bigcup\limits_{i=1}^{p}\left(
\left\{  i\right\}  \times\left[  m_{i}\right]  \right)  .
\label{pf.prop.young.convex.a.5}%
\end{equation}

Let us compare this with the diagram $D$. Let $i\in\left[  p\right]  $. Then,
$\left(  i,m_{i}\right)  $ is the easternmost cell in row $i$ of $D$. Hence,
$\left(  i,m_{i}\right)  \in D$. Therefore, all $m_{i}$ cells $\left(
i,1\right)  ,\ \left(  i,2\right)  ,\ \ldots,\ \left(  i,m_{i}\right)  $ must
belong to $D$ (by the northwest-closure property, since they lie weakly west
of $\left(  i,m_{i}\right)  $). No other cells in the $i$-th row can belong to
$D$, since $\left(  i,m_{i}\right)  $ is the easternmost cell in this row.
Thus, in total, the $i$-th row of $D$ consists of the cells $\left(
i,1\right)  ,\ \left(  i,2\right)  ,\ \ldots,\ \left(  i,m_{i}\right)  $ and
no others. In other words,%
\begin{align*}
\left(  \text{the }i\text{-th row of }D\right)   &  =\left\{  \left(
i,1\right)  ,\ \left(  i,2\right)  ,\ \ldots,\ \left(  i,m_{i}\right)
\right\} \\
&  =\left\{  \left(  i,j\right)  \ \mid\ j\in\left[  m_{i}\right]  \right\}
=\left\{  i\right\}  \times\left[  m_{i}\right]  .
\end{align*}

Forget that we fixed $i$. We thus have shown that%
\begin{equation}
\left(  \text{the }i\text{-th row of }D\right)  =\left\{  i\right\}
\times\left[  m_{i}\right]  \label{pf.prop.young.convex.a.7}%
\end{equation}
for each $i\in\left[  p\right]  $. Moreover, recall that $D$ has no cells
south of row $p$. Hence,%
\[
D=\bigcup\limits_{i=1}^{p}\underbrace{\left(  \text{the }i\text{-th row of
}D\right)  }_{\substack{=\left\{  i\right\}  \times\left[  m_{i}\right]
\\\text{(by (\ref{pf.prop.young.convex.a.7}))}}}=\bigcup\limits_{i=1}%
^{p}\left(  \left\{  i\right\}  \times\left[  m_{i}\right]  \right)  =Y\left(
m_{1},m_{2},\ldots,m_{p}\right)
\]
(by (\ref{pf.prop.young.convex.a.5})). Hence, $D$ is the Young diagram
$Y\left(  \lambda\right)  $ of some partition $\lambda$ (namely, of
$\lambda=\left(  m_{1},m_{2},\ldots,m_{p}\right)  $). This proves the
\textquotedblleft$\Longleftarrow$\textquotedblright\ direction of Proposition
\ref{prop.young.convex-etc} \textbf{(a)}. \medskip

\textbf{(b)} $\Longrightarrow:$ This follows from Proposition
\ref{prop.young.convexity}.

$\Longleftarrow:$ See Exercise \ref{exe.young.convexity-converse} below.
\medskip

\textbf{(c)} It is easy to see that for any partition $\lambda$ and any
$i\geq1$, we have%
\[
\lambda_{i}=\left(  \text{\# of cells in the }i\text{-th row of }Y\left(
\lambda\right)  \right)
\]
(indeed, the cells in the $i$-th row of $Y\left(  \lambda\right)  $ are
$\left(  i,1\right)  ,\ \left(  i,2\right)  ,\ \ldots,\ \left(  i,\lambda
_{i}\right)  $). Thus, each $\lambda_{i}$ and therefore the whole partition
$\lambda$ can be uniquely reconstructed from $Y\left(  \lambda\right)  $. This
proves Proposition \ref{prop.young.convex-etc} \textbf{(c)}. \medskip

\textbf{(d)} For a counterexample, observe that the skew partitions $\left(
1\right)  /\left(  1\right)  $ and $\varnothing/\varnothing$ have the same
skew Young diagram (namely, the empty set). More generally, for any partition
$\lambda$, the skew partition $\lambda/\lambda$ is well-defined and its skew
Young diagram $Y\left(  \lambda/\lambda\right)  $ is the empty set.

Nonempty counterexamples exist as well: For instance,%
\[
Y\left(  \left(  3,1,1\right)  /\left(  2,1\right)  \right)  =Y\left(  \left(
3,2,1\right)  /\left(  2,2\right)  \right)  =\left\{  \left(  1,3\right)
,\ \left(  3,1\right)  \right\}  .
\]

\end{proof}
\end{fineprint}

\begin{exercise}
\label{exe.young.convexity-converse}\fbox{3} Prove the \textquotedblleft%
$\Longleftarrow$\textquotedblright\ direction of Proposition
\ref{prop.young.convex-etc} \textbf{(b)}.
\end{exercise}

\subsubsection{Some notations regarding partitions}

\begin{definition}
\label{def.partitions.expnot}When writing partitions or compositions, we will
often use \emph{exponential notation}, which means that we will abbreviate a
sequence $\underbrace{i,i,\ldots,i}_{k\text{ times}}$ of equal entries by
writing $i^{k}$. To avoid confusion between this notation and an actual power,
we shall never put powers in a partition or composition without warning.
\end{definition}

For example, $\left(  4,2^{3},1\right)  $ means the partition $\left(
4,2,2,2,1\right)  $.

\begin{definition}
\label{def.partitions.hook}A \emph{hook partition} means a partition of the
form $\left(  i,1^{j}\right)  $ for some positive integer $i$ and some
nonnegative integer $j$. (If $j=0$, then this just means $\left(  i\right)  $.
If $i=1$, then this just means $\left(  1^{j+1}\right)  $.)
\end{definition}

Hook partitions owe their name to the fact that their Young diagrams look like
hooks: For example, the Young diagram $Y\left(  4,1^{2}\right)  =Y\left(
4,1,1\right)  $ is%
\[%
%TCIMACRO{\TeXButton{tikz Y411}{\begin{tikzpicture}[scale=0.7]
%\draw[fill=red!50] (0, 0) rectangle (1, 1);
%\draw[fill=red!50] (0, 1) rectangle (1, 2);
%\draw[fill=red!50] (0, 2) rectangle (1, 3);
%\draw[fill=red!50] (1, 2) rectangle (2, 3);
%\draw[fill=red!50] (2, 2) rectangle (3, 3);
%\draw[fill=red!50] (3, 2) rectangle (4, 3);
%\end{tikzpicture}}}%
%BeginExpansion
\begin{tikzpicture}[scale=0.7]
\draw[fill=red!50] (0, 0) rectangle (1, 1);
\draw[fill=red!50] (0, 1) rectangle (1, 2);
\draw[fill=red!50] (0, 2) rectangle (1, 3);
\draw[fill=red!50] (1, 2) rectangle (2, 3);
\draw[fill=red!50] (2, 2) rectangle (3, 3);
\draw[fill=red!50] (3, 2) rectangle (4, 3);
\end{tikzpicture}%
%EndExpansion
\ \ .
\]

\subsubsection{Conjugate partitions}

We next introduce a fundamental transformation on the set of partitions,
called \emph{conjugation} or \emph{transposition}:

\begin{theorem}
\label{thm.partitions.conj}Let $\mathbf{r}:\mathbb{Z}^{2}\rightarrow
\mathbb{Z}^{2}$ be the map that sends each cell $\left(  i,j\right)
\in\mathbb{Z}^{2}$ to $\left(  j,i\right)  $. (In the English coordinate
system, this map $\mathbf{r}$ is just the reflection across the
northwest-to-southeast diagonal.)

Let $\lambda$ be a partition. Then: \medskip

\textbf{(a)} There exists a unique partition $\mu$ such that $Y\left(
\mu\right)  =\mathbf{r}\left(  Y\left(  \lambda\right)  \right)  $. \medskip

This partition $\mu$ is called the \emph{conjugate} (or the \emph{transpose},
or the \emph{conjugate partition}) of $\lambda$. It is denoted by $\lambda
^{t}$ (or, by many authors, by $\lambda^{\prime}$). \medskip

\textbf{(b)} This partition $\lambda^{t}$ is explicitly given by the formula%
\begin{align*}
\lambda_{j}^{t}  &  =\left(  \text{\# of positive integers }i\text{ such that
}\lambda_{i}\geq j\right) \\
&  =\left(  \text{\# of cells in the }j\text{-th column of }Y\left(
\lambda\right)  \right)  \ \ \ \ \ \ \ \ \ \ \text{for all }j\geq1.
\end{align*}

\textbf{(c)} This partition $\lambda^{t}$ satisfies $\left\vert \lambda
^{t}\right\vert =\left\vert \lambda\right\vert $ and $\lambda_{1}^{t}%
=\ell\left(  \lambda\right)  $ and $\ell\left(  \lambda^{t}\right)
=\lambda_{1}$. \medskip

\textbf{(d)} The conjugate $\left(  \lambda^{t}\right)  ^{t}$ of this
partition $\lambda^{t}$ is $\lambda$. \medskip

\textbf{(e)} Let $j$ be a positive integer. Then, the cells in the $j$-th
column of $Y\left(  \lambda\right)  $ are%
\[
\left(  1,j\right)  ,\ \left(  2,j\right)  ,\ \ldots,\ \left(  \lambda_{j}%
^{t},j\right)  .
\]

\end{theorem}

\begin{example}
If $\lambda=\left(  4,2,1\right)  $, then the conjugate of $\lambda$ is
$\lambda^{t}=\left(  3,2,1,1\right)  $. Here are the Young diagrams of these
two partitions:%
\[
Y\left(  \lambda\right)  =%
%TCIMACRO{\TeXButton{tikz Y411}{\begin{tikzpicture}[scale=0.7]
%\draw[fill=red!50] (0, 0) rectangle (1, 1);
%\draw[fill=red!50] (0, 1) rectangle (1, 2);
%\draw[fill=red!50] (0, 2) rectangle (1, 3);
%\draw[fill=red!50] (1, 1) rectangle (2, 2);
%\draw[fill=red!50] (1, 2) rectangle (2, 3);
%\draw[fill=red!50] (2, 2) rectangle (3, 3);
%\draw[fill=red!50] (3, 2) rectangle (4, 3);
%\end{tikzpicture}}}%
%BeginExpansion
\begin{tikzpicture}[scale=0.7]
\draw[fill=red!50] (0, 0) rectangle (1, 1);
\draw[fill=red!50] (0, 1) rectangle (1, 2);
\draw[fill=red!50] (0, 2) rectangle (1, 3);
\draw[fill=red!50] (1, 1) rectangle (2, 2);
\draw[fill=red!50] (1, 2) rectangle (2, 3);
\draw[fill=red!50] (2, 2) rectangle (3, 3);
\draw[fill=red!50] (3, 2) rectangle (4, 3);
\end{tikzpicture}%
%EndExpansion
\ \ \ \ \ \ \ \ \ \ \text{and}\ \ \ \ \ \ \ \ \ \ Y\left(  \lambda^{t}\right)
=%
%TCIMACRO{\TeXButton{tikz Y3211}{\begin{tikzpicture}[scale=0.7]
%\draw[fill=red!50] (0, 0) rectangle (-1, -1);
%\draw[fill=red!50] (-1, 0) rectangle (-2, -1);
%\draw[fill=red!50] (-2, 0) rectangle (-3, -1);
%\draw[fill=red!50] (-1, -1) rectangle (-2, -2);
%\draw[fill=red!50] (-2, -1) rectangle (-3, -2);
%\draw[fill=red!50] (-2, -2) rectangle (-3, -3);
%\draw[fill=red!50] (-2, -3) rectangle (-3, -4);
%\end{tikzpicture}}}%
%BeginExpansion
\begin{tikzpicture}[scale=0.7]
\draw[fill=red!50] (0, 0) rectangle (-1, -1);
\draw[fill=red!50] (-1, 0) rectangle (-2, -1);
\draw[fill=red!50] (-2, 0) rectangle (-3, -1);
\draw[fill=red!50] (-1, -1) rectangle (-2, -2);
\draw[fill=red!50] (-2, -1) rectangle (-3, -2);
\draw[fill=red!50] (-2, -2) rectangle (-3, -3);
\draw[fill=red!50] (-2, -3) rectangle (-3, -4);
\end{tikzpicture}%
%EndExpansion
\ \ .
\]
(In the future, we will often be sloppy and just write $\lambda$ for $Y\left(
\lambda\right)  $.)
\end{example}

\begin{fineprint}
\begin{proof}
[Proof of Theorem \ref{thm.partitions.conj} (sketched).]\textbf{(a)} The
uniqueness of the partition $\mu$ follows from Proposition
\ref{prop.young.convex-etc} \textbf{(c)}. It remains to prove its existence.
In other words, it remains to prove that $\mathbf{r}\left(  Y\left(
\lambda\right)  \right)  $ is the Young diagram of some partition.

By Proposition \ref{prop.young.convex-etc} \textbf{(a)}, this will follow if
we can show that $\mathbf{r}\left(  Y\left(  \lambda\right)  \right)  $ has
the following property: If a cell $c\in\left\{  1,2,3,\ldots\right\}  ^{2}$
lies weakly northwest of a cell $d\in\mathbf{r}\left(  Y\left(  \lambda
\right)  \right)  $, then we have $c\in\mathbf{r}\left(  Y\left(
\lambda\right)  \right)  $ as well.

But this property of $\mathbf{r}\left(  Y\left(  \lambda\right)  \right)  $
follows easily from the analogous property of $Y\left(  \lambda\right)  $
(which is true by Proposition \ref{prop.young.NE}), since the reflection
$\mathbf{r}$ preserves the \textquotedblleft lying weakly northwest
of\textquotedblright\ relation. Thus, the existence of $\mu$ is proved, and
the proof of Theorem \ref{thm.partitions.conj} \textbf{(a)} is complete.
\medskip

\textbf{(b)} The definition of $\lambda^{t}$ yields $Y\left(  \lambda
^{t}\right)  =\mathbf{r}\left(  Y\left(  \lambda\right)  \right)  $. Thus,
each row of $Y\left(  \lambda^{t}\right)  $ is the reflection of the
corresponding column of $Y\left(  \lambda\right)  $ across the
northwest-to-southeast diagonal (since the reflection $\mathbf{r}$ turns
columns into rows). In particular, the length of each row of $Y\left(
\lambda^{t}\right)  $ equals the length of the corresponding column of
$Y\left(  \lambda\right)  $. In other words, for each $j\geq1$, we have%
\begin{align*}
&  \left(  \text{the length of the }j\text{-th row of }Y\left(  \lambda
^{t}\right)  \right) \\
&  =\left(  \text{the length of the }j\text{-th column of }Y\left(
\lambda\right)  \right) \\
&  =\left(  \text{\# of cells in the }j\text{-th column of }Y\left(
\lambda\right)  \right)  .
\end{align*}
But the left hand side of this equality is clearly $\lambda_{j}^{t}$, since
the length of the $j$-th row of any Young diagram $Y\left(  \mu\right)  $ is
$\mu_{j}$. Thus, this equality can be rewritten as%
\begin{align*}
\lambda_{j}^{t}  &  =\left(  \text{\# of cells in the }j\text{-th column of
}Y\left(  \lambda\right)  \right) \\
&  =\left(  \text{\# of positive integers }i\text{ such that }\lambda_{i}\geq
j\right)
\end{align*}
(since the cells in the $j$-th column of $Y\left(  \lambda\right)  $ are the
cells $\left(  i,j\right)  $ for all positive integers $i$ satisfying
$\lambda_{i}\geq j$). This proves Theorem \ref{thm.partitions.conj}
\textbf{(b)}. \medskip

\textbf{(d)} This follows from the fact that the reflection $\mathbf{r}$ is an
involution (i.e., satisfies $\mathbf{r}\circ\mathbf{r}=\operatorname*{id}$).
\medskip

\textbf{(c)} Proposition \ref{prop.young.n} yields $\left\vert Y\left(
\lambda\right)  \right\vert =\left\vert \lambda\right\vert $ and $\left\vert
Y\left(  \lambda^{t}\right)  \right\vert =\left\vert \lambda^{t}\right\vert $.
However, the definition of $\lambda^{t}$ yields $Y\left(  \lambda^{t}\right)
=\mathbf{r}\left(  Y\left(  \lambda\right)  \right)  $, so that $\left\vert
Y\left(  \lambda^{t}\right)  \right\vert =\left\vert Y\left(  \lambda\right)
\right\vert $ (since the reflection $\mathbf{r}$ is a bijection and thus
preserves the size of a set). This rewrites as $\left\vert \lambda
^{t}\right\vert =\left\vert \lambda\right\vert $ (since $\left\vert Y\left(
\lambda\right)  \right\vert =\left\vert \lambda\right\vert $ and $\left\vert
Y\left(  \lambda^{t}\right)  \right\vert =\left\vert \lambda^{t}\right\vert $).

Theorem \ref{thm.partitions.conj} \textbf{(b)} yields
\[
\lambda_{1}^{t}=\left(  \text{\# of positive integers }i\text{ such that
}\lambda_{i}\geq1\right)  =\ell\left(  \lambda\right)
\]
(because if $\ell\left(  \lambda\right)  =k$, then $\lambda=\left(
\lambda_{1},\lambda_{2},\ldots,\lambda_{k}\right)  $, and the positive
integers $i$ such that $\lambda_{i}\geq1$ are precisely the numbers
$1,2,\ldots,k$).

Thus we have proved the equality $\lambda_{1}^{t}=\ell\left(  \lambda\right)
$. Applying this equality to $\lambda^{t}$ instead of $\lambda$, we obtain
$\left(  \lambda^{t}\right)  _{1}^{t}=\ell\left(  \lambda^{t}\right)  $. Thus,
$\ell\left(  \lambda^{t}\right)  =\left(  \lambda^{t}\right)  _{1}^{t}%
=\lambda_{1}$ (since Theorem \ref{thm.partitions.conj} \textbf{(d)} yields
$\left(  \lambda^{t}\right)  ^{t}=\lambda$). With this equality, the proof of
Theorem \ref{thm.partitions.conj} \textbf{(c)} is complete. \medskip

\textbf{(e)} The reflection $\mathbf{r}$ is an involution, i.e., satisfies
$\mathbf{r}\circ\mathbf{r}=\operatorname*{id}$. Hence, $\mathbf{r}$ is
invertible, with inverse $\mathbf{r}^{-1}=\mathbf{r}$.

As we have seen in the proof of Theorem \ref{thm.partitions.conj} \textbf{(b)}
above, each row of $Y\left(  \lambda^{t}\right)  $ is the reflection of the
corresponding column of $Y\left(  \lambda\right)  $ across the
northwest-to-southeast diagonal. Thus, the $j$-th row of $Y\left(  \lambda
^{t}\right)  $ is the reflection of the $j$-th column of $Y\left(
\lambda\right)  $ across this diagonal. In other words,%
\[
\left(  \text{the }j\text{-th row of }Y\left(  \lambda^{t}\right)  \right)
=\mathbf{r}\left(  \text{the }j\text{-th column of }Y\left(  \lambda\right)
\right)  .
\]
Since $\mathbf{r}$ is invertible, this entails that%
\begin{align*}
\left(  \text{the }j\text{-th column of }Y\left(  \lambda\right)  \right)   &
=\underbrace{\mathbf{r}^{-1}}_{=\mathbf{r}}\underbrace{\left(  \text{the
}j\text{-th row of }Y\left(  \lambda^{t}\right)  \right)  }_{=\left\{  \left(
j,1\right)  ,\ \left(  j,2\right)  ,\ \ldots,\ \left(  j,\lambda_{j}%
^{t}\right)  \right\}  }\\
&  =\mathbf{r}\left(  \left\{  \left(  j,1\right)  ,\ \left(  j,2\right)
,\ \ldots,\ \left(  j,\lambda_{j}^{t}\right)  \right\}  \right) \\
&  =\left\{  \mathbf{r}\left(  j,1\right)  ,\ \mathbf{r}\left(  j,2\right)
,\ \ldots,\ \mathbf{r}\left(  j,\lambda_{j}^{t}\right)  \right\} \\
&  =\left\{  \left(  1,j\right)  ,\ \left(  2,j\right)  ,\ \ldots,\ \left(
\lambda_{j}^{t},j\right)  \right\}
\end{align*}
(since $\mathbf{r}\left(  u,v\right)  =\left(  v,u\right)  $ for any cell
$\left(  u,v\right)  \in\mathbb{Z}^{2}$). In other words, the cells in the
$j$-th column of $Y\left(  \lambda\right)  $ are $\left(  1,j\right)
,\ \left(  2,j\right)  ,\ \ldots,\ \left(  \lambda_{j}^{t},j\right)  $. This
proves Theorem \ref{thm.partitions.conj} \textbf{(e)}.
\end{proof}
\end{fineprint}

See \cite[\S 4.1.5]{21s} for a basic application of conjugate partitions.

\begin{remark}
\label{rmk.partitions.conj.Ylm}Let $\lambda/\mu$ be a skew partition. Then,
$\lambda^{t}/\mu^{t}$ is again a skew partition, and we have $Y\left(
\lambda^{t}/\mu^{t}\right)  =\mathbf{r}\left(  Y\left(  \lambda/\mu\right)
\right)  $.
\end{remark}

\begin{fineprint}
\begin{proof}
We know that $\lambda/\mu$ is a skew partition, i.e., a pair $\left(
\mu,\lambda\right)  $ of two partitions satisfying $\lambda\supseteq\mu$ (by
the definition of a skew partition). We have $\lambda\supseteq\mu$; in other
words, $\mu\subseteq\lambda$. Hence, Proposition \ref{prop.young.b-in-a}
(applied to $a=\lambda$ and $b=\mu$) shows that $Y\left(  \mu\right)
\subseteq Y\left(  \lambda\right)  $. Hence, $\mathbf{r}\left(  Y\left(
\mu\right)  \right)  \subseteq\mathbf{r}\left(  Y\left(  \lambda\right)
\right)  $.

But the definition of $\lambda^{t}$ shows that $Y\left(  \lambda^{t}\right)
=\mathbf{r}\left(  Y\left(  \lambda\right)  \right)  $. Similarly, $Y\left(
\mu^{t}\right)  =\mathbf{r}\left(  Y\left(  \mu\right)  \right)  $. Thus,
$Y\left(  \mu^{t}\right)  =\mathbf{r}\left(  Y\left(  \mu\right)  \right)
\subseteq\mathbf{r}\left(  Y\left(  \lambda\right)  \right)  =Y\left(
\lambda^{t}\right)  $. By Proposition \ref{prop.young.b-in-a} (applied to
$a=\lambda^{t}$ and $b=\mu^{t}$), we thus conclude that $\mu^{t}%
\subseteq\lambda^{t}$. Thus, $\lambda^{t}/\mu^{t}$ is a skew partition.

We have $Y\left(  \lambda/\mu\right)  =Y\left(  \lambda\right)  \setminus
Y\left(  \mu\right)  $ (by the definition of a skew Young diagram). Hence,
\[
\mathbf{r}\left(  Y\left(  \lambda/\mu\right)  \right)  =\mathbf{r}\left(
Y\left(  \lambda\right)  \setminus Y\left(  \mu\right)  \right)
=\mathbf{r}\left(  Y\left(  \lambda\right)  \right)  \setminus\mathbf{r}%
\left(  Y\left(  \mu\right)  \right)
\]
(since $\mathbf{r}$ is a bijection, and thus respects set differences).
Comparing this with%
\begin{align*}
Y\left(  \lambda^{t}/\mu^{t}\right)   &  =\underbrace{Y\left(  \lambda
^{t}\right)  }_{=\mathbf{r}\left(  Y\left(  \lambda\right)  \right)
}\setminus\underbrace{Y\left(  \mu^{t}\right)  }_{=\mathbf{r}\left(  Y\left(
\mu\right)  \right)  }\ \ \ \ \ \ \ \ \ \ \left(
\begin{array}
[c]{c}%
\text{by the definition of}\\
\text{a skew Young diagram}%
\end{array}
\right) \\
&  =\mathbf{r}\left(  Y\left(  \lambda\right)  \right)  \setminus
\mathbf{r}\left(  Y\left(  \mu\right)  \right)  ,
\end{align*}
we obtain $Y\left(  \lambda^{t}/\mu^{t}\right)  =\mathbf{r}\left(  Y\left(
\lambda/\mu\right)  \right)  $. The proof of Remark
\ref{rmk.partitions.conj.Ylm} is thus complete.
\end{proof}
\end{fineprint}

\subsection{\label{sec.specht.tabs}Tableaux}

\subsubsection{Definitions and examples}

While Young diagrams have some uses visualizing partitions (see, e.g.,
\cite[Proof of Proposition 4.1.15]{21s} for a nice application), to us they
are mainly important as tables to be filled with numbers. Once filled with
numbers, they are called \emph{tableaux}. (This is a loanword from French; its
singular form is \textquotedblleft tableau\textquotedblright.) Let us define
tableaux accurately and rigorously:

\begin{definition}
Let $D$ be a diagram. \medskip

\textbf{(a)} A \emph{tableau} of shape $D$ (aka a \emph{filling} of $D$) means
a map $T:D\rightarrow\left\{  1,2,3,\ldots\right\}  $ that assigns a positive
integer $T\left(  c\right)  $ to each cell $c\in D$. In this case, $T\left(
c\right)  $ is called the \emph{entry} of $T$ in cell $c$, and we say that the
cell $c$ of $T$ \emph{contains} the number $T\left(  c\right)  $. We visualize
a tableau $T$ by placing the entry $T\left(  c\right)  $ into the cell $c$ for
each cell $c$.

For example, here is a tableau of shape $Y\left(  \left(  5,4,4\right)
/\left(  2,1\right)  \right)  $ visualized in this way:%
\[
\ytableaushort{\none\none231,\none131,5351}\ \ .
\]
Rigorously speaking, it is a map $T:Y\left(  \left(  5,4,4\right)  /\left(
2,1\right)  \right)  \rightarrow\left\{  1,2,3,\ldots\right\}  $ that sends
the cell $\left(  1,3\right)  $ to $2$, the cell $\left(  1,4\right)  $ to $3$
and so on. \medskip

\textbf{(b)} A tableau $T$ is said to be an $m$\emph{-tableau} for some
$m\in\mathbb{N}$ if it is injective (i.e., has no two equal entries) and its
image is $\left[  m\right]  $. (This means that the entries of $T$ are the
numbers $1,2,\ldots,m$, with each of these numbers appearing exactly once.
Thus, an $m$-tableau of shape $D$ can be viewed as a bijection from $D$ to
$\left[  m\right]  $.)

Obviously, there are $m!$ many $m$-tableaux of shape $D$ if $\left\vert
D\right\vert =m$, and otherwise there are none. \medskip

\textbf{(c)} If $T$ is a tableau of shape $D$, and $i$ is a positive integer,
then the $i$\emph{-th row} of $T$ shall mean the restriction of $T$ to the
$i$-th row of $D$. This restriction contains the entries of $T$ that lie in
the $i$-th row. Similarly, the $j$\emph{-th column} of $T$ is defined. The
\emph{cells} of a tableau $T$ are just defined to be the cells of its shape
$D$. \medskip

\textbf{(d)} A tableau $T$ is said to be

\begin{itemize}
\item \emph{row-standard} if its entries strictly increase (left to right)
along each row.

\item \emph{row-semistandard} if its entries weakly increase (left to right)
along each row.

\item \emph{column-standard} if its entries strictly increase (top to bottom)
down each column.

\item \emph{column-semistandard} if its entries weakly increase (top to
bottom) down each column.

\item \emph{standard} if it is row-standard and column-standard and an
$m$-tableau for some $m$.

\item \emph{semistandard} if it is row-semistandard and column-standard (not column-semistandard!).

\item \emph{straight-shaped} (or of \emph{straight shape}, or just
\emph{straight}) if it has shape $Y\left(  \lambda\right)  $ for a partition
$\lambda$;

\item \emph{skew-shaped} (or of \emph{skew shape}, or just \emph{skew}) if it
has shape $Y\left(  \lambda/\mu\right)  $ for a skew partition $\lambda/\mu$
(this includes the straight-shaped case, since $Y\left(  \lambda
/\varnothing\right)  =Y\left(  \lambda\right)  $).

\item \emph{bad-shaped} if it is not skew-shaped.
\end{itemize}

We often abbreviate \textquotedblleft tableau of shape $Y\left(
\lambda\right)  $\textquotedblright\ (where $\lambda$ is a partition) as
\textquotedblleft tableau of shape $\lambda$\textquotedblright, and likewise
we abbreviate \textquotedblleft tableau of shape $Y\left(  \lambda/\mu\right)
$\textquotedblright\ (where $\lambda/\mu$ is a skew partition) as
\textquotedblleft tableau of shape $\lambda/\mu$\textquotedblright.
\end{definition}

\begin{remark}
Tableaux are also known as \emph{Young tableaux}. Depending on the author,
this name can come with extra requirements, such as standard or semistandard.
Different books do not agree on what exactly they require! (For instance,
Fulton's \cite{Fulton97} defines \textquotedblleft tableaux\textquotedblright%
\ and \textquotedblleft Young tableaux\textquotedblright\ to mean what we call
straight semistandard tableaux, whereas all other tableaux are just called
\textquotedblleft fillings\textquotedblright.)
\end{remark}

\begin{example}
\label{exa.tableau.tableau1}There are only $2$ standard tableaux of shape
$Y\left(  2,2\right)  $, namely%
\[
\ytableaushort{12,34}\ \ \ \ \ \ \ \ \ \ \text{and}%
\ \ \ \ \ \ \ \ \ \ \ytableaushort{13,24}\ \ .
\]

There are $5$ standard tableaux of shape $Y\left(  3,2\right)  $, namely%
\[
\ytableaushort{123,45}\qquad\ytableaushort{124,35}\qquad
\ytableaushort{125,34}\qquad\ytableaushort{134,25}\qquad
\ytableaushort{135,24}\ \ .
\]

We cannot list the other types of tableaux, since there are usually infinitely
many of them. But here are some examples:

\begin{itemize}
\item The tableau $\ytableaushort{\none23,133}$ is row-semistandard,
column-semistandard and skew (having shape $Y\left(  \left(  3,3\right)
/\left(  1\right)  \right)  $), but neither row-standard (since $3$ appears
twice in the second row) nor column-standard (since $3$ appears twice in the
third column) nor straight nor semistandard.

\item The tableau $\ytableaushort{\none\none12,1\none3,\none\none46}$ is
row-standard, column-standard and bad-shaped. It is nevertheless not standard,
since it is not an $m$-tableau for any $m\in\mathbb{N}$ (indeed, it contains
the entry $1$ twice and the entry $5$ never).

\item The tableau $\ytableaushort{\none13,224,355}$ is semistandard and skew.
\end{itemize}
\end{example}

\begin{fineprint}
For users of LaTeX: \href{https://ctan.org/pkg/ytableau?lang=en}{the
\texttt{ytableau} package} provides an easy way to typeset tableaux of various
kinds. See
\href{https://mirrors.mit.edu/CTAN/macros/latex/contrib/ytableau/ytableau.pdf}{its
documentation} for several ways it can be used.

For users of SageMath: There is a lot of infrastructure for working with
\href{https://doc.sagemath.org/html/en/reference/combinat/sage/combinat/tableau.html}{straight-shaped}
and
\href{https://doc.sagemath.org/html/en/reference/combinat/sage/combinat/skew_tableau.html}{skew-shaped}
tableaux in SageMath (not so much for bad-shaped tableaux). For instance,
typing in \texttt{list(StandardTableaux([4,2,1]))} will return all standard
tableaux of shape $Y\left(  4,2,1\right)  $. For another instance,
\texttt{SkewTableau([[None,2,3],[1,4]])} returns the skew (standard) tableau
$\ytableaushort{\none23,14}$ .
\end{fineprint}

\Needspace{5cm}

\begin{convention}
\label{conv.tableau.poetic}Sometimes we will use a shorthand notation for
tableaux (particularly when tableaux will appear as subscripts), in which we
write all rows of a tableau in a single line (from top to bottom), separated
by \textquotedblleft$\backslash\backslash$\textquotedblright\ signs, and with
each row being written as the list of its entries (with the \textquotedblleft%
$\circ$\textquotedblright\ sign signifying a missing cell). For example, the
tableau%
\[
\ytableaushort{\none13,224,355}
\]
will thus be written as%
\[
\circ,1,3\backslash\backslash2,2,4\backslash\backslash3,5,5.
\]
Often we will also omit the commas between the entries (as long as the entries
are single-digit numbers). Thus, the above tableau will be written%
\[
\circ13\backslash\backslash224\backslash\backslash355.
\]

\end{convention}

\subsubsection{First properties}

Obviously, all straight tableaux are skew as well (since $Y\left(
\lambda\right)  =Y\left(  \lambda/\varnothing\right)  $ for any partition
$\lambda$). Thus, the following basic property of skew tableaux applies to
straight ones as well:

\begin{proposition}
\label{prop.tableau.NW}Let $T$ be a skew tableau of some shape $D$. Let $c$
and $e$ be two cells of $D$ such that $c$ lies weakly northwest of $e$. Then:
\medskip

\textbf{(a)} If $T$ is row-semistandard and column-semistandard, then
$T\left(  c\right)  \leq T\left(  e\right)  $. \medskip

\textbf{(b)} If $T$ is semistandard, and if $c$ and $e$ lie in different rows,
then $T\left(  c\right)  <T\left(  e\right)  $.
\end{proposition}

\begin{proof}
The tableau $T$ is skew; thus, its shape $D$ has the form $Y\left(
\lambda/\mu\right)  $ for some skew partition $\lambda/\mu$. Consider this
$\lambda/\mu$.

Write the cells $c$ and $e$ in the forms $c=\left(  i,j\right)  $ and
$e=\left(  u,v\right)  $. Since $c$ lies weakly northwest of $e$, we thus have
$i\leq u$ and $j\leq v$.

Let $d$ be the cell $\left(  i,v\right)  $. This cell $d$ lies in the same row
as $c=\left(  i,j\right)  $ but weakly east of $c$ (since $v\geq j$). It
furthermore lies in the same column as $e=\left(  u,v\right)  $ but weakly
north of $e$ (since $i\leq u$). In particular, we have $c\leq d\leq e$. Hence,
Proposition \ref{prop.young.convexity} yields $d\in Y\left(  \lambda
/\mu\right)  $ (since $c\in D=Y\left(  \lambda/\mu\right)  $ and $e\in
D=Y\left(  \lambda/\mu\right)  $). In other words, $d\in D$ (since $D=Y\left(
\lambda/\mu\right)  $). Hence, the entry $T\left(  d\right)  $ exists.
\medskip

\textbf{(a)} Assume that $T$ is row-semistandard and column-semistandard.

Since $T$ is row-semistandard, the entries in the $i$-th row of $T$ weakly
increase left to right. Hence, $T\left(  c\right)  \leq T\left(  d\right)  $
(since the cells $c=\left(  i,j\right)  $ and $d=\left(  i,v\right)  $ both
lie in the $i$-th row, and $d$ lies weakly east of $c$).

Since $T$ is column-semistandard, the entries in the $v$-th column of $T$
weakly increase top to bottom. Hence, $T\left(  d\right)  \leq T\left(
e\right)  $ (since the cells $d=\left(  i,v\right)  $ and $e=\left(
u,v\right)  $ both lie in the $v$-th column, and $d$ lies weakly north of $e$).

Altogether, we conclude that $T\left(  c\right)  \leq T\left(  d\right)  \leq
T\left(  e\right)  $. This proves Proposition \ref{prop.tableau.NW}
\textbf{(a)}. \medskip

\textbf{(b)} Assume that $T$ is semistandard. Thus, $T$ is row-semistandard
and column-standard. Hence, as in the proof of part \textbf{(a)}, we can show
that $T\left(  c\right)  \leq T\left(  d\right)  $.

Now, let us additionally assume that $c$ and $e$ lie in different rows. In
other words, $i\neq u$ (since $c=\left(  i,j\right)  $ lies in row $i$,
whereas $e=\left(  u,v\right)  $ lies in row $u$). Hence, $i<u$ (since $i\leq
u$). Thus, the cell $d=\left(  i,v\right)  $ lies strictly (not just weakly)
north of $e=\left(  u,v\right)  $.

Since $T$ is column-standard, the entries in the $v$-th column of $T$ strictly
increase top to bottom. Hence, $T\left(  d\right)  <T\left(  e\right)  $
(since the cells $d=\left(  i,v\right)  $ and $e=\left(  u,v\right)  $ both
lie in the $v$-th column, and $d$ lies strictly north of $e$).

Altogether, we conclude that $T\left(  c\right)  \leq T\left(  d\right)
<T\left(  e\right)  $. This proves Proposition \ref{prop.tableau.NW}
\textbf{(b)}.
\end{proof}

Even more obvious is the following property:

\begin{proposition}
\label{prop.tableau.std-n}Let $D$ be any diagram with $\left\vert D\right\vert
=n$. Then, the standard tableaux of shape $D$ are precisely the standard
$n$-tableaux of shape $D$.
\end{proposition}

\begin{fineprint}
\begin{proof}
Let $T$ be a standard tableau of shape $D$. We shall show that $T$ is an $n$-tableau.

Indeed, $T$ is standard, and thus is an $m$-tableau for some $m\in\mathbb{N}$.
Consider this $m$. Then, $T$ is a bijection from $D$ to $\left[  m\right]  $
(since any $m$-tableau of shape $D$ can be viewed as a bijection from $D$ to
$\left[  m\right]  $). Hence, $\left\vert D\right\vert =\left\vert \left[
m\right]  \right\vert =m$. Thus, $m=\left\vert D\right\vert =n$. But $T$ is an
$m$-tableau. In other words, $T$ is an $n$-tableau (since $m=n$). Thus, $T$ is
a standard $n$-tableau of shape $D$.

Forget that we fixed $T$. We thus have shown that each standard tableau $T$ of
shape $D$ is a standard $n$-tableau of shape $D$. Conversely, of course, any
standard $n$-tableau of shape $D$ is a standard tableau of shape $D$.
Combining these two facts, we conclude that the standard tableaux of shape $D$
are precisely the standard $n$-tableaux of shape $D$. This proves Proposition
\ref{prop.tableau.std-n}.
\end{proof}
\end{fineprint}

The following \textquotedblleft local\textquotedblright\ criterion for row-
and column-standardness is an important advantage of skew-shaped tableaux (it
does not hold for bad-shaped tableaux):

\begin{proposition}
\label{prop.tableau.std-loc}Let $D$ be a skew Young diagram. Let $T$ be a
tableau of shape $D$. Then: \medskip

\textbf{(a)} If every $\left(  i,j\right)  \in D$ satisfying $\left(
i,j+1\right)  \in D$ satisfies $T\left(  i,j\right)  <T\left(  i,j+1\right)
$, then $T$ is row-standard. \medskip

\textbf{(b)} If every $\left(  i,j\right)  \in D$ satisfying $\left(
i+1,j\right)  \in D$ satisfies $T\left(  i,j\right)  <T\left(  i+1,j\right)
$, then $T$ is column-standard.
\end{proposition}

In words, Proposition \ref{prop.tableau.std-loc} \textbf{(a)} is saying that
in order to check that a skew-shaped tableau is row-standard, we only need to
check that its entries in two \textbf{horizontally adjacent} cells (i.e., two
cells $\left(  i,j\right)  $ and $\left(  i,j+1\right)  $) are always in the
right order (i.e., the left one is always smaller than the right one).
Intuitively speaking, this is because the rows of $D$ are contiguous blocks of
cells, so you can walk from every cell to any other (in the same row) without
having to jump over any missing cells. Proposition \ref{prop.tableau.std-loc}
\textbf{(b)} is analogous, just replacing rows by columns. For the sake of
completeness, we give a rigorous proof of both parts of Proposition
\ref{prop.tableau.std-loc} in the Appendix (Section
\ref{sec.details.specht.tabs}).

\subsubsection{The hook length formula}

It is easy to count $n$-tableaux that are row-standard or column-standard:

\begin{exercise}
Let $\lambda=\left(  \lambda_{1},\lambda_{2},\ldots,\lambda_{k}\right)  $ be
any partition with $\left\vert \lambda\right\vert =n$. Show the following:
\medskip

\textbf{(a)} \fbox{1} The \# of all $n$-tableaux of shape $Y\left(
\lambda\right)  $ is $n!$. \medskip

\textbf{(b)} \fbox{2} The \# of all row-standard $n$-tableaux of shape
$Y\left(  \lambda\right)  $ is $\dfrac{n!}{\lambda_{1}!\lambda_{2}%
!\cdots\lambda_{k}!}$. \medskip

\textbf{(c)} \fbox{1} The \# of all column-standard $n$-tableaux of shape
$Y\left(  \lambda\right)  $ is $\dfrac{n!}{\lambda_{1}^{t}!\lambda_{2}%
^{t}!\cdots\lambda_{m}^{t}!}$, where $\lambda^{t}=\left(  \lambda_{1}%
^{t},\lambda_{2}^{t},\ldots,\lambda_{m}^{t}\right)  $ is the conjugate of the
partition $\lambda$. \medskip

\textbf{(d)} \fbox{1} The \# of all column-standard $n$-tableaux of shape
$Y\left(  \lambda\right)  $ is also equal to $\dfrac{n!}{1^{\lambda_{1}%
}2^{\lambda_{2}}\cdots k^{\lambda_{k}}}$.
\end{exercise}

Now here is something a lot less easy:

\begin{theorem}
[Hook length formula]\label{thm.tableau.hlf}Let $\lambda$ be any partition
with $\left\vert \lambda\right\vert =n$. Then, the \# of standard tableaux of
shape $Y\left(  \lambda\right)  $ is%
\[
\dfrac{n!}{%
%TCIMACRO{\dprod \limits_{c\in Y\left(  \lambda\right)  }}%
%BeginExpansion
{\displaystyle\prod\limits_{c\in Y\left(  \lambda\right)  }}
%EndExpansion
\left\vert H_{\lambda}\left(  c\right)  \right\vert }.
\]
Here, for each cell $c\in Y\left(  \lambda\right)  $, we define the
\emph{hook} $H_{\lambda}\left(  c\right)  $ of this cell in $\lambda$ to be
\begin{align*}
H_{\lambda}\left(  c\right)   &  :=\left\{  c\right\}  \cup\left\{  \text{all
cells in }Y\left(  \lambda\right)  \text{ that lie east of }c\text{ (in the
same row as }c\text{)}\right\} \\
&  \ \ \ \ \ \ \ \ \ \ \cup\left\{  \text{all cells in }Y\left(
\lambda\right)  \text{ that lie south of }c\text{ (in the same column as
}c\text{)}\right\}  .
\end{align*}

\end{theorem}

\begin{example}
Let $n=5$ and $\lambda=\left(  3,2\right)  $. This partition $\lambda$ has
Young diagram%
\[
Y\left(  \lambda\right)  =%
%TCIMACRO{\TeXButton{tikz Y32}{\begin{tikzpicture}[scale=0.7]
%\draw[fill=red!50] (0, 0) rectangle (1, 1);
%\draw[fill=red!50] (1, 0) rectangle (2, 1);
%\draw[fill=red!50] (2, 1) rectangle (3, 2);
%\draw[fill=red!50] (0, 1) rectangle (1, 2);
%\draw[fill=red!50] (1, 1) rectangle (2, 2);
%\end{tikzpicture}}}%
%BeginExpansion
\begin{tikzpicture}[scale=0.7]
\draw[fill=red!50] (0, 0) rectangle (1, 1);
\draw[fill=red!50] (1, 0) rectangle (2, 1);
\draw[fill=red!50] (2, 1) rectangle (3, 2);
\draw[fill=red!50] (0, 1) rectangle (1, 2);
\draw[fill=red!50] (1, 1) rectangle (2, 2);
\end{tikzpicture}%
%EndExpansion
\ =\left\{  \left(  1,1\right)  ,\ \left(  1,2\right)  ,\ \left(  1,3\right)
,\ \left(  2,1\right)  ,\ \left(  2,2\right)  \right\}  .
\]
The hooks of its cells are
\begin{align*}
H_{\lambda}\left(  1,1\right)   &  =\left\{  \left(  1,1\right)  \right\}
\cup\left\{  \left(  1,2\right)  ,\ \left(  1,3\right)  \right\}  \cup\left\{
\left(  2,1\right)  \right\}  ,\\
H_{\lambda}\left(  1,2\right)   &  =\left\{  \left(  1,2\right)  \right\}
\cup\left\{  \left(  1,3\right)  \right\}  \cup\left\{  \left(  2,2\right)
\right\}  ,\\
H_{\lambda}\left(  1,3\right)   &  =\left\{  \left(  1,3\right)  \right\}
\cup\varnothing\cup\varnothing,\\
H_{\lambda}\left(  2,1\right)   &  =\left\{  \left(  2,1\right)  \right\}
\cup\left\{  \left(  2,2\right)  \right\}  \cup\varnothing,\\
H_{\lambda}\left(  2,2\right)   &  =\left\{  \left(  2,2\right)  \right\}
\cup\varnothing\cup\varnothing,
\end{align*}
with respective sizes $4$, $3$, $1$, $2$ and $1$. Thus, the hook length
formula (Theorem \ref{thm.tableau.hlf}) yields that the \# of standard
tableaux of shape $Y\left(  \lambda\right)  $ is%
\[
\dfrac{n!}{%
%TCIMACRO{\dprod \limits_{c\in Y\left(  \lambda\right)  }}%
%BeginExpansion
{\displaystyle\prod\limits_{c\in Y\left(  \lambda\right)  }}
%EndExpansion
\left\vert H_{\lambda}\left(  c\right)  \right\vert }=\dfrac{5!}{4\cdot
3\cdot1\cdot2\cdot1}=5,
\]
just as we found in Example \ref{exa.tableau.tableau1}.
\end{example}

The hook length formula is a surprising and nontrivial result (in general, it
is far from obvious that $\dfrac{n!}{%
%TCIMACRO{\dprod \limits_{c\in Y\left(  \lambda\right)  }}%
%BeginExpansion
{\displaystyle\prod\limits_{c\in Y\left(  \lambda\right)  }}
%EndExpansion
\left\vert H_{\lambda}\left(  c\right)  \right\vert }$ is an integer!), and
provides some first evidence that standard tableaux have a deeper meaning. We
will prove it in Section \ref{sec.specht.hlf-young}. Other (usually simpler
and shorter) proofs of Theorem \ref{thm.tableau.hlf} can also be found in
\cite[Chapter 8, Highlight]{Aigner07}, \cite[Theorem 6.5]{Bona22},
\cite[Theorem 4.2.14]{CSScTo10}, \cite[Theorem 5.17.1]{EGHetc11},
\cite[\S 4.3, Exercise 10]{Fulton97}, \cite[Remark 19]{Hopkin22},
\cite[\S 20]{James78}, \cite[2.3.21]{JamKer81}, \cite[6.3.8]{Kerber99},
\cite[\S 5.1.4, Theorem H]{Knuth-TAoCP3}, \cite[\S V.3]{Krishn86},
\cite[Theorem 5.4.10]{LakBro18}, \cite[\S 12.10, Theorem 12.50]{Loehr-BC},
\cite[\S 4.3.5]{Lorenz18}, \cite[\S I.5, Example 2]{Macdon95}, \cite[Theorem
3.41]{Meliot17}, \cite[(5.7)]{MenRem15}, \cite[Theorem 5.8.3]{Prasad-rep},
\cite[(5.2.1)]{Proces07}, \cite[Theorem 1.12]{Romik15}, \cite[Theorem
3.10.2]{Sagan01}, \cite[Theorem 7.3.1]{Sagan19}, \cite[Corollary
7.21.6]{Stanley-EC2}, \cite[Theorem 8.1]{Stanley-AC}, \cite[\S 6.5,
Proposition (b)]{Zelevi81b}, and myriad other sources. See \cite[\S 11.2]%
{Pak22} for a taxonomy of different proofs and references to each type.

\begin{exercise}
\fbox{2} Prove Theorem \ref{thm.tableau.hlf} by hand in the case when
$\lambda$ is a hook partition.
\end{exercise}

The following exercise applies Theorem \ref{thm.tableau.hlf} to Young diagrams
with two rows (i.e., diagrams of the form $Y\left(  \lambda\right)  $ where
$\lambda=\left(  \lambda_{1},\lambda_{2}\right)  $):

\begin{exercise}
\textbf{(a)} \fbox{1} Use Theorem \ref{thm.tableau.hlf} to show that for each
$m\in\mathbb{N}$, the \# of standard tableaux of shape $Y\left(  m,m\right)  $
is $\dfrac{1}{m+1}\dbinom{2m}{m}$. This number is known as the $m$-th Catalan
number (see, e.g., \cite[\S 3.1.2]{21s} or \cite[Lecture 29, \S 5.3]{22fco}).
\medskip

\textbf{(b)} \fbox{2} Prove the same claim without using Theorem
\ref{thm.tableau.hlf}, by instead constructing a bijection between the
standard tableaux of shape $Y\left(  m,m\right)  $ and the Catalan lattice
paths from $\left(  0,0\right)  $ to $\left(  m,m\right)  $ (see \cite[Lecture
29, \S 5.2]{22fco}). \medskip

\textbf{(c)} \fbox{3} More generally, compute the \# of standard tableaux of
shape $Y\left(  a,b\right)  $ for all $a,b\in\mathbb{N}$ with $a\geq b$ (that
is, of shape $Y\left(  \lambda\right)  $ for any partition $\lambda$ of length
$\leq2$) both using Theorem \ref{thm.tableau.hlf} and without Theorem
\ref{thm.tableau.hlf}. \medskip

[\textbf{Hint:} For part \textbf{(c)}, use \cite[Theorem 5.2.4 \textbf{(d)}%
]{22fco}.]
\end{exercise}

In the form we presented it, Theorem \ref{thm.tableau.hlf} was first stated by
Frame, Robinson and Thrall in 1954 (\cite[Theorem 1]{FrRoTh54}), who proved it
by deriving it from an older (1902) formula due to Frobenius and Young:

\begin{theorem}
[Frobenius--Young formula]\label{thm.tableau.hlf-young}Let $\lambda=\left(
\lambda_{1},\lambda_{2},\ldots,\lambda_{k}\right)  $ be any partition with
$\left\vert \lambda\right\vert =n$. Let $\ell_{i}:=\lambda_{i}+k-i$ for each
$i\in\left[  k\right]  $. Then, the \# of standard tableaux of shape $Y\left(
\lambda\right)  $ is%
\[
\dfrac{n!}{\ell_{1}!\cdot\ell_{2}!\cdot\cdots\cdot\ell_{k}!}\prod_{1\leq
i<j\leq k}\left(  \ell_{i}-\ell_{j}\right)  .
\]

\end{theorem}

\begin{exercise}
\fbox{3} Prove that Theorem \ref{thm.tableau.hlf} and Theorem
\ref{thm.tableau.hlf-young} are equivalent, by showing that%
\[
\dfrac{1}{%
%TCIMACRO{\dprod \limits_{c\in Y\left(  \lambda\right)  }}%
%BeginExpansion
{\displaystyle\prod\limits_{c\in Y\left(  \lambda\right)  }}
%EndExpansion
\left\vert H_{\lambda}\left(  c\right)  \right\vert }=\dfrac{1}{\ell_{1}%
!\cdot\ell_{2}!\cdot\cdots\cdot\ell_{k}!}\prod_{1\leq i<j\leq k}\left(
\ell_{i}-\ell_{j}\right)
\]
(where all the notations are as in Theorem \ref{thm.tableau.hlf-young}).
\end{exercise}

The hook length formula (Theorem \ref{thm.tableau.hlf}) can be generalized to
skew shapes in two ways, yielding respectively the Okounkov--Olshanski formula
(\cite[Theorem 8.1]{OkoOls98}, \cite{MorZhu20}, \cite[Proposition
4.12]{PanPet24}) and the Naruse formula (\cite[Theorem 2.11]{GKMK23},
\cite[(NHLF)]{PanPet24}, \cite{MPP1}, \cite{MPP2}, \cite{Konva1},
\cite{Konva2}). Since they are significantly more complex than Theorem
\ref{thm.tableau.hlf} (in particular, they count the \# of standard tableaux
of shape $Y\left(  \lambda/\mu\right)  $ not by a single fraction, but by a
sum of several such, which are usually not integers), we shall not state them here.

Another nice and useful property of skew tableaux is the \emph{non-messing-up
theorem}, which we give in two variants as an exercise:

\begin{exercise}
\label{exe.tableau.nmu}\textbf{(a)} \fbox{2} Let $T$ be a skew
row-semistandard tableau.

Let $S$ be the tableau obtained from $T$ by sorting each column in weakly
increasing order. (This is supposed to mean that for each $j\in\mathbb{Z}$,
the $j$-th column of $S$ contains the same entries as the $j$-th column of $T$
but rearranged in weakly increasing order from top to bottom. The shapes of
$S$ and $T$ are the same, and no entries get moved out of their column.)

Prove that $S$ is row-semistandard (in addition to obviously being
column-semistandard). \medskip

\textbf{(b)} \fbox{1} Prove that the same is true if we replace
\textquotedblleft row-semistandard\textquotedblright\ by \textquotedblleft
row-standard\textquotedblright\ (both for $T$ and for $S$). \medskip

\textbf{(c)} \fbox{1} Prove that neither part \textbf{(a)} nor part
\textbf{(b)} remains true if the tableaux are allowed to be bad-shape.
\end{exercise}

However, don't get too encouraged:

\begin{exercise}
\fbox{1} If $T$ is an $n$-tableau, then $T^{\operatorname*{row}}$ shall denote
the tableau obtained from $T$ by sorting each row in increasing order, whereas
$T^{\operatorname{col}}$ shall denote the tableau obtained from $T$ by sorting
each column in increasing order.

Find an example of a straight-shaped $n$-tableau $T$ such that $\left(
T^{\operatorname{row}}\right)  ^{\operatorname{col}}\neq\left(
T^{\operatorname{col}}\right)  ^{\operatorname{row}}$.
\end{exercise}

\subsection{Young modules}

\subsubsection{Horizontal and vertical permutations}

Let us now get symmetric groups into the action.

\begin{noncompile}
Let $D$ be a diagram. Then: \medskip

\textbf{(a)} Its \emph{row group} $\mathcal{R}\left(  D\right)  $ is the
subgroup of $S_{D}$ consisting of all permutations $w\in S_{D}$ that send each
cell to a cell in the same row. These permutations are also known as
\emph{horizontal permutations}. \medskip

\textbf{(b)} Its \emph{column group} $\mathcal{C}\left(  D\right)  $ is the
subgroup of $S_{D}$ consisting of all permutations $w\in S_{D}$ that send each
cell to a cell in the same column. These permutations are also known as
\emph{vertical permutations}.
\end{noncompile}

\begin{definition}
\label{def.tableau.cell-cts}Let $T$ be an $n$-tableau (of any shape $D$).
Thus, each $i\in\left[  n\right]  $ lies in exactly one cell of $T$, so that
we can talk about \textquotedblleft the cell of $T$ that contains
$i$\textquotedblright. (Rigorously speaking, this cell is simply
$T^{-1}\left(  i\right)  $.)\medskip

\textbf{(a)} A permutation $w\in S_{n}$ is said to be \emph{horizontal} for
$T$ if it has the property that for each $i\in\left[  n\right]  $, the number
$w\left(  i\right)  $ lies in the same row of $T$ as $i$ does.

The \emph{row group} $\mathcal{R}\left(  T\right)  $ of $T$ is the subgroup of
$S_{n}$ consisting of all permutations $w\in S_{n}$ that are horizontal for
$T$. \medskip

\textbf{(b)} A permutation $w\in S_{n}$ is said to be \emph{vertical} for $T$
if it has the property that for each $i\in\left[  n\right]  $, the number
$w\left(  i\right)  $ lies in the same column of $T$ as $i$ does.

The \emph{column group} $\mathcal{C}\left(  T\right)  $ of $T$ is the subgroup
of $S_{n}$ consisting of all permutations $w\in S_{n}$ that are vertical for
$T$.
\end{definition}

In other words, a permutation $w\in S_{n}$ is \emph{horizontal} for an
$n$-tableau $T$ if it permutes the set of entries of each row\footnote{A
permutation $w\in S_{n}$ is said to \emph{permute} some subset $I$ of $\left[
n\right]  $ if it sends $I$ to itself (i.e., if $w\left(  i\right)  \in I$ for
each $i\in I$).}. Likewise, a permutation $w\in S_{n}$ is \emph{vertical} for
an $n$-tableau $T$ if it permutes the set of entries of each column.

For example, for $T=\ytableaushort{134,25}$, we have%
\[
\mathcal{R}\left(  T\right)  =\left\{  \text{all permutations in }S_{5}\text{
permuting }\left\{  1,3,4\right\}  \text{ and permuting }\left\{  2,5\right\}
\right\}
\]
(a subgroup of $S_{5}$ having size $3!\cdot2!=12$) and%
\begin{align*}
\mathcal{C}\left(  T\right)   &  =\left\{  \text{all permutations in }%
S_{5}\text{ permuting }\left\{  1,2\right\}  \text{ and permuting }\left\{
3,5\right\}  \right. \\
&  \ \ \ \ \ \ \ \ \ \ \left.  \text{and permuting }\left\{  4\right\}
\right\}
\end{align*}
(a subgroup of $S_{5}$ having size $2!\cdot2!\cdot1!=4$; its elements are
$\operatorname*{id}$, $t_{1,2}$, $t_{3,5}$ and $t_{1,2}t_{3,5}$).

\subsubsection{Row-equivalence and column-equivalence}

\begin{definition}
\label{def.tableau.roweq}\textbf{(a)} Two tableaux $S$ and $T$ are said to be
\emph{row-equivalent} if they have the same entries in corresponding rows
(i.e., if for every $i\in\mathbb{Z}$, the multiset of entries in the $i$-th
row of $S$ is the multiset of entries in the $i$-th row of $T$). \medskip

\textbf{(b)} Two tableaux $S$ and $T$ are said to be \emph{column-equivalent}
if they have the same entries in corresponding columns (i.e., if for every
$j\in\mathbb{Z}$, the multiset of entries in the $j$-th column of $S$ is the
multiset of entries in the $j$-th column of $T$).
\end{definition}

\begin{example}
\textbf{(a)} The tableaux $\ytableaushort{114,23}$ and
$\ytableaushort{141,32}$ are row-equivalent, but they are not row-equivalent
to $\ytableaushort{144,23}$ (since multisets are not equal if they have
different multiplicities). \medskip

\textbf{(b)} The tableaux $\ytableaushort{114,23}$ and
$\ytableaushort{134,21}$ are column-equivalent, but they are not
column-equivalent to $\ytableaushort{114,32}$\ .
\end{example}

\begin{fineprint}
Note that tableaux of different shapes can be row-equivalent! For instance,
the tableaux $\ytableaushort{12,3}$ and $\ytableaushort{12,\none3}$ are
row-equivalent. That said, we will usually not compare tableaux of different
shapes in what follows. \medskip
\end{fineprint}

Obviously, row-equivalence is an equivalence relation, and so is
column-equivalence. The following is obvious as well:

\begin{proposition}
\label{prop.tableau.roweq.H=H}\textbf{(a)} If two $n$-tableaux $S$ and $T$ are
row-equivalent, then $\mathcal{R}\left(  S\right)  =\mathcal{R}\left(
T\right)  $. \medskip

\textbf{(b)} If two $n$-tableaux $S$ and $T$ are column-equivalent, then
$\mathcal{C}\left(  S\right)  =\mathcal{C}\left(  T\right)  $.
\end{proposition}

Note that the converse of either statement is false; for example, the two
$4$-tableaux $S=\ytableaushort{12,34}$ and $T=\ytableaushort{34,12}$ are not
row-equivalent but satisfy $\mathcal{R}\left(  S\right)  =\mathcal{R}\left(
T\right)  $. \medskip

Here is a property of row and column groups that is only slightly less obvious
than Proposition \ref{prop.tableau.roweq.H=H}:

\begin{proposition}
\label{prop.tableau.RnC}Let $T$ be an $n$-tableau of any shape $D$. Then:
\medskip

\textbf{(a)} We have $\mathcal{R}\left(  T\right)  \cap\mathcal{C}\left(
T\right)  =\left\{  \operatorname*{id}\right\}  $. \medskip

\textbf{(b)} The only pair $\left(  c,r\right)  \in\mathcal{C}\left(
T\right)  \times\mathcal{R}\left(  T\right)  $ satisfying
$cr=\operatorname*{id}$ is the pair $\left(  \operatorname*{id}%
,\operatorname*{id}\right)  $.
\end{proposition}

\begin{fineprint}
\begin{proof}
\textbf{(a)} Let $w\in\mathcal{R}\left(  T\right)  \cap\mathcal{C}\left(
T\right)  $. Thus, $w\in\mathcal{R}\left(  T\right)  $ and $w\in
\mathcal{C}\left(  T\right)  $.

Let $i\in\left[  n\right]  $. The tableau $T$ is an $n$-tableau, and thus
contains each of the integers $1,2,\ldots,n$ exactly once. Hence, there is a
unique cell $c$ of $T$ that contains the number $i$, and there is a unique
cell $d$ of $T$ that contains the number $w\left(  i\right)  $. Consider these
two cells $c$ and $d$. Thus, $T\left(  c\right)  =i$ (since the cell $c$ of
$T$ contains the number $i$) and $T\left(  d\right)  =w\left(  i\right)  $ (similarly).

From $w\in\mathcal{R}\left(  T\right)  $, we see that the permutation $w$ is
horizontal for $T$ (by the definition of $\mathcal{R}\left(  T\right)  $).
Hence, the number $w\left(  i\right)  $ lies in the same row of $T$ as $i$
does (by the definition of \textquotedblleft horizontal\textquotedblright). In
other words, the cell of $T$ that contains the number $i$ lies in the same row
as the cell of $T$ that contains the number $w\left(  i\right)  $. Since the
former cell is $c$, and since the latter cell is $d$, we can rewrite this as
follows: The cell $c$ lies in the same row as the cell $d$.

Similarly, from $w\in\mathcal{C}\left(  T\right)  $, we can see that the cell
$c$ lies in the same column as the cell $d$.

However, if two cells lie in the same row and in the same column, then they
are simply the same cell. Thus, $c$ and $d$ are the same cell (since $c$ lies
in the same row and in the same column as $d$). In other words, $c=d$. Hence,
$T\left(  c\right)  =T\left(  d\right)  =w\left(  i\right)  $, so that
$w\left(  i\right)  =T\left(  c\right)  =i$.

Forget that we fixed $i$. We thus have proved that $w\left(  i\right)  =i$ for
each $i\in\left[  n\right]  $. Thus, $w=\operatorname*{id}$.

Forget that we fixed $w$. We have now shown that $w=\operatorname*{id}$ for
each $w\in\mathcal{R}\left(  T\right)  \cap\mathcal{C}\left(  T\right)  $. In
other words, $\mathcal{R}\left(  T\right)  \cap\mathcal{C}\left(  T\right)
\subseteq\left\{  \operatorname*{id}\right\}  $. Combining this with $\left\{
\operatorname*{id}\right\}  \subseteq\mathcal{R}\left(  T\right)
\cap\mathcal{C}\left(  T\right)  $ (which is obvious, since the permutation
$\operatorname*{id}$ is both horizontal and vertical for $T$), we obtain
$\mathcal{R}\left(  T\right)  \cap\mathcal{C}\left(  T\right)  =\left\{
\operatorname*{id}\right\}  $. Proposition \ref{prop.tableau.RnC} \textbf{(a)}
is now proved. \medskip

\textbf{(b)} Clearly, $\left(  \operatorname*{id},\operatorname*{id}\right)  $
is a pair $\left(  c,r\right)  \in\mathcal{C}\left(  T\right)  \times
\mathcal{R}\left(  T\right)  $ satisfying $cr=\operatorname*{id}$. It remains
to show that there are no other such pairs. So let $\left(  c,r\right)
\in\mathcal{C}\left(  T\right)  \times\mathcal{R}\left(  T\right)  $ be a pair
that satisfies $cr=\operatorname*{id}$. We must prove that $\left(
c,r\right)  =\left(  \operatorname*{id},\operatorname*{id}\right)  $.

From $\left(  c,r\right)  \in\mathcal{C}\left(  T\right)  \times
\mathcal{R}\left(  T\right)  $, we obtain $c\in\mathcal{C}\left(  T\right)  $
and $r\in\mathcal{R}\left(  T\right)  $. From $cr=\operatorname*{id}$, we
obtain $c=r^{-1}\in\mathcal{R}\left(  T\right)  $ (since $r\in\mathcal{R}%
\left(  T\right)  $, but $\mathcal{R}\left(  T\right)  $ is a group).
Combining this with $c\in\mathcal{C}\left(  T\right)  $, we obtain
$c\in\mathcal{R}\left(  T\right)  \cap\mathcal{C}\left(  T\right)  =\left\{
\operatorname*{id}\right\}  $ (by Proposition \ref{prop.tableau.RnC}
\textbf{(a)}), so that $c=\operatorname*{id}$. Hence, the equality
$cr=\operatorname*{id}$ rewrites as $\operatorname*{id}r=\operatorname*{id}$,
thus $r=\operatorname*{id}$. Hence, $\left(  c,r\right)  =\left(
\operatorname*{id},\operatorname*{id}\right)  $ (since $c=\operatorname*{id}$
and $r=\operatorname*{id}$).

Forget that we fixed $\left(  c,r\right)  $. We thus have shown that if
$\left(  c,r\right)  \in\mathcal{C}\left(  T\right)  \times\mathcal{R}\left(
T\right)  $ is a pair that satisfies $cr=\operatorname*{id}$, then $\left(
c,r\right)  =\left(  \operatorname*{id},\operatorname*{id}\right)  $. Hence,
the only such pair is $\left(  \operatorname*{id},\operatorname*{id}\right)  $
(since we already know that $\left(  \operatorname*{id},\operatorname*{id}%
\right)  $ is such a pair). This proves Proposition \ref{prop.tableau.RnC}
\textbf{(b)}.
\end{proof}
\end{fineprint}

\subsubsection{Tabloids}

\begin{definition}
\label{def.tabloid.tabloid}\textbf{(a)} Let $D$ be a diagram. An
$n$\emph{-tabloid} of shape $D$ is an equivalence class of $n$-tableaux of
shape $D$ with respect to row-equivalence. For example, for $D=Y\left(
2,1\right)  $ and $n=3$, the set $\left\{  12\backslash\backslash
3,\ 21\backslash\backslash3\right\}  $ (where we are using the shorthand
notation from Convention \ref{conv.tableau.poetic}) is an $n$-tabloid of shape
$D$. \medskip

\textbf{(b)} We use the notation $\overline{T}$ for the $n$-tabloid
corresponding to an $n$-tableau $T$ (that is, the equivalence class containing
$T$). For instance, the $3$-tabloid $\left\{  12\backslash\backslash
3,\ 21\backslash\backslash3\right\}  $ can thus be written as $\overline
{12\backslash\backslash3}$ or (equivalently) as $\overline{21\backslash
\backslash3}$.

To draw an $n$-tabloid, we just draw one of the $n$-tableaux in it (often we
pick the row-standard one), and we drop the vertical subdividers between the
entries of a row. For instance, the $3$-tabloid $\overline{12\backslash
\backslash3}$ is thus drawn as
$\ytableausetup{tabloids}\ytableaushort{12,3}\ytableausetup{notabloids}$ .
\medskip

\textbf{(c)} If $\overline{T}$ is an $n$-tabloid, and if $i\in\mathbb{Z}$,
then the $i$\emph{-th row of }$\overline{T}$ shall mean the set of the entries
of the $i$-th row of $T$. This is well-defined, since two row-equivalent
tableaux will have the same entries in their $i$-th row. For example, the
$2$-nd row of the $5$-tabloid $\ytableausetup{tabloids}
\ytableaushort{24,15,3} \ytableausetup{notabloids}$ is the set $\left\{
1,5\right\}  $.
\end{definition}

\begin{example}
\label{exa.tabloid.n-1-1}\textbf{(a)} There are four $4$-tabloids of shape
$Y\left(  3,1\right)  $, namely%
\[
\ytableausetup{tabloids} \ytableaushort{123,4}\ ,\qquad
\ytableaushort{124,3}\ ,\qquad\ytableaushort{134,2}\ ,\qquad
\ytableaushort{234,1}\ . \ytableausetup{notabloids}
\]

\textbf{(b)} More generally: Assume that $n\geq2$. What are the $n$-tabloids
of shape $Y\left(  n-1,1\right)  $ ?

Such an $n$-tabloid is uniquely determined by its only entry in the $2$-nd
row. Indeed, if its entry in the $2$-nd row is $i$, then its $1$-st row must
contain the remaining numbers $1,2,\ldots,i-1,i+1,i+2,\ldots,n$, and the order
in which they appear is immaterial (since row-equivalent tableaux produce the
same tabloid). The entry $i$ in the $2$-nd row can be any $i\in\left[
n\right]  $, and we get a different $n$-tabloid for each $i$. Thus, in total,
there are $n$ many $n$-tabloids of shape $Y\left(  n-1,1\right)  $.
\end{example}

\begin{exercise}
Let $k\in\left[  n-1\right]  $. \medskip

\textbf{(a)} \fbox{1} What are the $n$-tabloids of shape $Y\left(
n-k,1^{k}\right)  $ ? How many are there? \medskip

\textbf{(b)} \fbox{1} Assume that $k\leq n/2$. What are the $n$-tabloids of
shape $Y\left(  n-k,k\right)  $ ? How many are there?
\end{exercise}

\begin{fineprint}
We note that many authors denote an $n$-tabloid $\overline{T}$ by $\left\{
T\right\}  $ instead, but I find this notation misleading: An $n$-tabloid
surely is a set and contains $T$ (being the equivalence class of $T$), but
usually it does not consist of $T$ alone! \medskip
\end{fineprint}

While $n$-tabloids are equivalence classes, they also can be represented in a
canonical way, namely by row-standard $n$-tableaux:

\begin{proposition}
\label{prop.tabloid.row-st}Let $D$ be any diagram. Then, there is a bijection%
\begin{align*}
&  \text{from }\left\{  \text{row-standard }n\text{-tableaux of shape
}D\right\} \\
&  \text{to }\left\{  n\text{-tabloids of shape }D\right\}
\end{align*}
that sends each row-standard $n$-tableau $T$ to its tabloid $\overline{T}$.
\end{proposition}

\begin{proof}
The map from $\left\{  \text{row-standard }n\text{-tableaux of shape
}D\right\}  $ to $\left\{  n\text{-tabloids of shape }D\right\}  $ that sends
each $T$ to $\overline{T}$ is clearly well-defined. It remains to show that
this map is a bijection. In other words, it remains to show that each
$n$-tabloid of shape $D$ can be written as $\overline{T}$ for a unique
row-standard $n$-tableau $T$ of shape $D$. In other words, we need to show
that each $n$-tabloid of shape $D$ has a unique row-standard
representative\footnote{A \emph{representative} of an equivalence class just
means an element of this class.}.

But this is easy: We can pick any representative of our $n$-tabloid, and then
sort each row in increasing order (from left to right) so that the
representative becomes row-standard. This row-standard representative is
furthermore unique, since the row-standardness requirement forces the entries
in each row to appear in increasing order and therefore binds each entry to a
unique position. Thus, we have shown that our $n$-tabloid has a unique
row-standard representative. As explained above, this completes the proof of
Proposition \ref{prop.tabloid.row-st}.
\end{proof}

\subsubsection{The symmetric group acts on tableaux and tabloids}

Now, let us define an action of the symmetric group $S_{n}$ on sets of tableaux:

\begin{definition}
\label{def.tableau.Sn-act}Let $D$ be any diagram. Recall that an $n$-tableau
of shape $D$ is defined as a map from $D$ to $\left\{  1,2,3,\ldots\right\}  $
that is injective and whose image is $\left[  n\right]  $. Thus, an
$n$-tableau of shape $D$ can be viewed as a bijection from $D$ to $\left[
n\right]  $. Hence, if $T$ is an $n$-tableau of shape $D$, and if $w\in S_{n}$
is a permutation of $\left[  n\right]  $, then the composition $w\circ T$ is a
bijection from $D$ to $\left[  n\right]  $, that is, an $n$-tableau of shape
$D$ again.

Now, we define a left $S_{n}$-action on the set $\left\{  n\text{-tableaux of
shape }D\right\}  $ by setting%
\[
w\rightharpoonup T:=w\circ T\ \ \ \ \ \ \ \ \ \ \text{for all }w\in
S_{n}\text{ and all }n\text{-tableaux }T
\]
(recalling that $T$ is a map from $D$ to $\left[  n\right]  $). This really is
a well-defined action of $S_{n}$, since composition of maps is associative.
Explicitly, this action can be described as follows: The tableau
$w\rightharpoonup T$ is obtained from $T$ by applying $w$ to each entry (i.e.,
replacing each entry $i$ by $w\left(  i\right)  $). Thus, this action is
called the \emph{action on entries}.
\end{definition}

For example, for $D=Y\left(  2,1\right)  $ and $w\in S_{3}$, we have%
\[
w\rightharpoonup
\ytableaushort{12,3}=\ytableausetup{boxsize=3em}\begin{ytableau}
w\left(1\right) & w\left(2\right) \\
w\left(3\right)
\end{ytableau}\ytableausetup{boxsize=normal}\ \ \ \ \ \ \ \ \ \ \text{and}%
\ \ \ \ \ \ \ \ \ \ w\rightharpoonup
\ytableaushort{21,3}=\ytableausetup{notabloids, boxsize=3em}\begin{ytableau}
w\left(2\right) & w\left(1\right) \\
w\left(3\right)
\end{ytableau}\ytableausetup{boxsize=normal}\ \ .
\]
(The cells have been enlarged merely to fit the entries into them; they are
still meant to be $1\times1$-squares.) \medskip

The $S_{n}$-action on $\left\{  n\text{-tableaux of shape }D\right\}  $
defined in Definition \ref{def.tableau.Sn-act} can be furthermore
\textquotedblleft pulled down\textquotedblright\ to the set $\left\{
n\text{-tabloids of shape }D\right\}  $, in a rather obvious way:

\begin{definition}
\label{def.tabloid.Sn-act}Let $D$ be any diagram. Then, we define a left
$S_{n}$-action on the set $\left\{  n\text{-tabloids of shape }D\right\}  $ by
setting%
\[
w\rightharpoonup\overline{T}:=\overline{w\rightharpoonup T}=\overline{w\circ
T}\ \ \ \ \ \ \ \ \ \ \text{for all }w\in S_{n}\text{ and all }%
n\text{-tableaux }T.
\]
This action is well-defined (we will justify this in a moment), and is again
called the \emph{action on entries} (but this time on tabloids instead of
tableaux). This action makes the set $\left\{  n\text{-tabloids of shape
}D\right\}  $ into a left $S_{n}$-set.
\end{definition}

Why is this action well-defined? The haughty answer is \textquotedblleft by
Definition \ref{def.rep.G-sets.eqrel} \textbf{(b)}\textquotedblright. Indeed,
as we know, the $n$-tabloids of shape $D$ are the equivalence classes of the
\textquotedblleft row-equivalence\textquotedblright\ relation on the set
$\left\{  n\text{-tableaux of shape }D\right\}  $. The latter relation is
easily seen to be $S_{n}$-invariant\footnote{This is just saying that if $w\in
S_{n}$ is any permutation, and if $T$ and $S$ are two row-equivalent
$n$-tableaux, then the compositions $w\circ T$ and $w\circ S$ are
row-equivalent as well. Convince yourself that this is indeed obvious.}. Thus,
Definition \ref{def.rep.G-sets.eqrel} \textbf{(b)} yields that the quotient
set
\[
\left\{  n\text{-tableaux of shape }D\right\}  /\left(  \text{row-equivalence}%
\right)
\]
becomes a left $S_{n}$-set. But this quotient set is precisely the set
$\left\{  n\text{-tabloids of shape }D\right\}  $, and the left $S_{n}$-action
we get is precisely the one defined in Definition \ref{def.tabloid.Sn-act}.
Thus, the $S_{n}$-action in Definition \ref{def.tabloid.Sn-act} is
well-defined. Of course, we can obtain the same well-definedness
\textquotedblleft by foot\textquotedblright: Observe that if $w\in S_{n}$ is
any permutation, and if $T$ and $S$ are two row-equivalent $n$-tableaux, then
the compositions $w\circ T$ and $w\circ S$ are row-equivalent as well. Thus,
the $n$-tabloid $\overline{w\circ T}$ depends only on $w$ and on $\overline
{T}$ but not on the representative $T$. This shows that $w\rightharpoonup
\overline{T}$ is well-defined. Verifying the associativity and unitality
axioms is completely straightforward. \medskip

As usual with group actions, we abbreviate $w\rightharpoonup T$ by $wT$ (and
likewise $w\rightharpoonup\overline{T}$ by $w\overline{T}$) whenever we feel
like doing so. Thus, for example, the equality in Definition
\ref{def.tabloid.Sn-act} can be rewritten as%
\begin{equation}
w\overline{T}:=\overline{wT}=\overline{w\circ T}.
\label{eq.def.tabloid.Sn-act.rewr}%
\end{equation}

\begin{fineprint}
\begin{remark}
There is also an action on places: Let $D$ be any diagram. Then, the symmetric
group $S_{D}$ (that is, the group of all permutations of the cells of $D$)
acts on the sets%
\[
\left\{  n\text{-tableaux of shape }D\right\}  \ \ \ \ \ \ \ \ \ \ \text{and}%
\ \ \ \ \ \ \ \ \ \ \left\{  \text{tableaux of shape }D\right\}
\]
from the \textbf{right}. In both cases, the action is defined by%
\[
T\leftharpoonup u:=T\circ u\ \ \ \ \ \ \ \ \ \ \text{for all tableaux }T\text{
and all }u\in S_{D}.
\]
These actions are called the \emph{action on places}. Note that the group
acting here is $S_{D}$, not $S_{n}$.

For instance, for $D=Y\left(  2,1\right)  $ and $u=t_{\left(  1,1\right)
,\left(  1,2\right)  }\in S_{D}$, we have%
\[
\ytableaushort{ij,k}\leftharpoonup u=\ytableaushort{ji,k}.
\]

No such action can be defined on any set of $n$-tabloids, since the definition
of an $n$-tabloid \textquotedblleft smears an entry across the entire
row\textquotedblright\ and thus makes it unclear what entries to move where.
\end{remark}
\end{fineprint}

We shall next prove a few basic properties of our $S_{n}$-actions:

\begin{proposition}
\label{prop.tableau.Sn-act.0}Let $D$ be any diagram. Let $T$ be an $n$-tableau
of shape $D$. Let $w\in S_{n}$. Then: \medskip

\textbf{(a)} The tableau $w\rightharpoonup T$ is row-equivalent to $T$ if and
only if $w\in\mathcal{R}\left(  T\right)  $. \medskip

\textbf{(b)} The tableau $w\rightharpoonup T$ is column-equivalent to $T$ if
and only if $w\in\mathcal{C}\left(  T\right)  $. \medskip

\textbf{(c)} We have%
\[
\mathcal{R}\left(  w\rightharpoonup T\right)  =w\mathcal{R}\left(  T\right)
w^{-1}%
\]
and%
\[
\mathcal{C}\left(  w\rightharpoonup T\right)  =w\mathcal{C}\left(  T\right)
w^{-1}.
\]

\end{proposition}

\begin{fineprint}
\begin{proof}
\textbf{(a)} We must prove the equivalence
\begin{equation}
\left(  w\rightharpoonup T\text{ is row-equivalent to }T\right)
\ \Longleftrightarrow\ \left(  w\in\mathcal{R}\left(  T\right)  \right)  .
\label{pf.prop.tableau.Sn-act.0.a.goal}%
\end{equation}

We will prove the \textquotedblleft$\Longleftarrow$\textquotedblright\ and the
\textquotedblleft$\Longrightarrow$\textquotedblright\ directions of this
equivalence separately:

$\Longleftarrow:$ Assume that $w\in\mathcal{R}\left(  T\right)  $. Thus, the
permutation $w$ is horizontal for $T$ (by the definition of $\mathcal{R}%
\left(  T\right)  $). In other words, for each $i\in\left[  n\right]  $, the
number $w\left(  i\right)  $ lies in the same row of $T$ as $i$ does (by the
definition of \textquotedblleft horizontal\textquotedblright).

Recall that $T$ is an $n$-tableau, thus injective (by the definition of an
$n$-tableau). Hence, $T$ has no two equal entries. Similarly,
$w\rightharpoonup T$ has no two equal entries (since $w\rightharpoonup T$ is
an $n$-tableau).

Let $r\in\mathbb{Z}$. Let $t_{1},t_{2},\ldots,t_{k}$ be all the entries in the
$r$-th row of $T$.

These entries $t_{1},t_{2},\ldots,t_{k}$ are distinct (since $T$ has no two
equal entries). Hence, their images $w\left(  t_{1}\right)  ,w\left(
t_{2}\right)  ,\ldots,w\left(  t_{k}\right)  $ are distinct as well (since $w$
is a permutation and thus injective).

However, we know that for each $i\in\left[  n\right]  $, the number $w\left(
i\right)  $ lies in the same row of $T$ as $i$ does. Thus, the numbers
$w\left(  t_{1}\right)  ,w\left(  t_{2}\right)  ,\ldots,w\left(  t_{k}\right)
$ lie in the same rows of $T$ as the numbers $t_{1},t_{2},\ldots,t_{k}$ do.
Since the latter numbers lie in the $r$-th row of $T$, it follows that the
former do as well. In other words,%
\begin{align*}
\left\{  w\left(  t_{1}\right)  ,w\left(  t_{2}\right)  ,\ldots,w\left(
t_{k}\right)  \right\}   &  \subseteq\left\{  \text{the entries in the
}r\text{-th row of }T\right\} \\
&  =\left\{  t_{1},t_{2},\ldots,t_{k}\right\}
\end{align*}
(since $t_{1},t_{2},\ldots,t_{k}$ are all the entries in the $r$-th row of $T$).

The set $\left\{  w\left(  t_{1}\right)  ,w\left(  t_{2}\right)
,\ldots,w\left(  t_{k}\right)  \right\}  $ has size $k$ (since $w\left(
t_{1}\right)  ,w\left(  t_{2}\right)  ,\ldots,w\left(  t_{k}\right)  $ are
distinct), and the set $\left\{  t_{1},t_{2},\ldots,t_{k}\right\}  $ has size
$k$ as well (since $t_{1},t_{2},\ldots,t_{k}$ are distinct). Thus, the two
finite sets $\left\{  w\left(  t_{1}\right)  ,w\left(  t_{2}\right)
,\ldots,w\left(  t_{k}\right)  \right\}  $ and $\left\{  t_{1},t_{2}%
,\ldots,t_{k}\right\}  $ have equal sizes (namely, $k$). It is well-known that
if two finite sets $U$ and $V$ have equal sizes and satisfy $U\subseteq V$,
then $U=V$. Applying this to $U=\left\{  w\left(  t_{1}\right)  ,w\left(
t_{2}\right)  ,\ldots,w\left(  t_{k}\right)  \right\}  $ and $V=\left\{
t_{1},t_{2},\ldots,t_{k}\right\}  $, we thus obtain
\begin{equation}
\left\{  w\left(  t_{1}\right)  ,w\left(  t_{2}\right)  ,\ldots,w\left(
t_{k}\right)  \right\}  =\left\{  t_{1},t_{2},\ldots,t_{k}\right\}
\label{pf.prop.tableau.Sn-act.1.a.2}%
\end{equation}
(since $\left\{  w\left(  t_{1}\right)  ,w\left(  t_{2}\right)  ,\ldots
,w\left(  t_{k}\right)  \right\}  \subseteq\left\{  t_{1},t_{2},\ldots
,t_{k}\right\}  $).

But the entries in the $r$-th row of $T$ are $t_{1},t_{2},\ldots,t_{k}$ (by
the definition of $t_{1},t_{2},\ldots,t_{k}$). Hence, the entries in the
$r$-th row of $w\rightharpoonup T$ are $w\left(  t_{1}\right)  ,w\left(
t_{2}\right)  ,\ldots,w\left(  t_{k}\right)  $ (since the tableau
$w\rightharpoonup T$ is obtained from $T$ by replacing each entry $i$ by
$w\left(  i\right)  $). In other words,%
\begin{align*}
\left\{  \text{the entries in the }r\text{-th row of }w\rightharpoonup
T\right\}   &  =\left\{  w\left(  t_{1}\right)  ,w\left(  t_{2}\right)
,\ldots,w\left(  t_{k}\right)  \right\} \\
&  =\left\{  t_{1},t_{2},\ldots,t_{k}\right\}  \ \ \ \ \ \ \ \ \ \ \left(
\text{by (\ref{pf.prop.tableau.Sn-act.1.a.2})}\right) \\
&  =\left\{  \text{the entries in the }r\text{-th row of }T\right\}
\end{align*}
(by the definition of $t_{1},t_{2},\ldots,t_{k}$). In other words, the set of
entries in the $r$-th row of $w\rightharpoonup T$ is the set of entries in the
$r$-th row of $T$. Moreover, we can replace both words \textquotedblleft
set\textquotedblright\ in this sentence by \textquotedblleft
multiset\textquotedblright, since the multisets in question contain no element
more than once (because neither $w\rightharpoonup T$ nor $T$ has two equal
entries). Thus, the multiset of entries in the $r$-th row of $w\rightharpoonup
T$ is the multiset of entries in the $r$-th row of $T$.

Forget that we fixed $r$. We thus have shown that for every $r\in\mathbb{Z}$,
the multiset of entries in the $r$-th row of $w\rightharpoonup T$ is the
multiset of entries in the $r$-th row of $T$. In other words, the tableaux
$w\rightharpoonup T$ and $T$ have the same entries in corresponding rows. In
other words, these tableaux $w\rightharpoonup T$ and $T$ are row-equivalent.
That is, $w\rightharpoonup T$ is row-equivalent to $T$. This proves the
\textquotedblleft$\Longleftarrow$\textquotedblright\ direction of the
equivalence (\ref{pf.prop.tableau.Sn-act.0.a.goal}).

$\Longrightarrow:$ Assume that $w\rightharpoonup T$ is row-equivalent to $T$.

Now, let $i\in\left[  n\right]  $. Then, $i$ appears in some cell of $T$
(since $T$ is an $n$-tableau). Let $r$ be the row in which this cell lies.
Hence, $i$ appears in the $r$-th row of $T$.

Recall that the tableau $w\rightharpoonup T$ is obtained from $T$ by applying
$w$ to each entry. Hence, whatever cell of $T$ contains the entry $i$ will
contain the entry $w\left(  i\right)  $ in $w\rightharpoonup T$. In
particular, this shows that $w\left(  i\right)  $ appears in the $r$-th row of
$w\rightharpoonup T$ (since $i$ appears in the $r$-th row of $T$).

But the tableaux $w\rightharpoonup T$ and $T$ are row-equivalent, i.e., have
the same entries in corresponding rows (by the definition of \textquotedblleft
row-equivalent\textquotedblright). Thus, any entry that appears in the $r$-th
row of $w\rightharpoonup T$ must also appear in the $r$-th row of $T$. Thus,
$w\left(  i\right)  $ appears in the $r$-th row of $T$ (since $w\left(
i\right)  $ appears in the $r$-th row of $w\rightharpoonup T$). Therefore, the
number $w\left(  i\right)  $ lies in the same row of $T$ as $i$ does (namely,
in the $r$-th row).

Forget that we fixed $i$. We thus have shown that for each $i\in\left[
n\right]  $, the number $w\left(  i\right)  $ lies in the same row of $T$ as
$i$ does. In other words, the permutation $w$ is horizontal for $T$ (by the
definition of \textquotedblleft horizontal\textquotedblright). In other words,
$w\in\mathcal{R}\left(  T\right)  $ (by the definition of $\mathcal{R}\left(
T\right)  $). This proves the \textquotedblleft$\Longrightarrow$%
\textquotedblright\ direction of the equivalence
(\ref{pf.prop.tableau.Sn-act.0.a.goal}).

Now, we have proved both \textquotedblleft$\Longleftarrow$\textquotedblright%
\ and \textquotedblleft$\Longrightarrow$\textquotedblright\ directions of the
equivalence (\ref{pf.prop.tableau.Sn-act.0.a.goal}). Thus, this equivalence
holds. This proves Proposition \ref{prop.tableau.Sn-act.0} \textbf{(a)}.
\medskip

\textbf{(b)} This is analogous to part \textbf{(a)}. (We just need to replace
rows by columns, and the symbol $\mathcal{R}$ by the symbol $\mathcal{C}$.)
\medskip

\textbf{(c)} We have $w\in S_{n}$. Thus, $w$ is a bijection from $\left[
n\right]  $ to $\left[  n\right]  $.

If $u\in S_{n}$ is any permutation, then we have the following chain of
equivalences:%
\begin{align}
&  \ \left(  u\in\mathcal{R}\left(  w\rightharpoonup T\right)  \right)
\nonumber\\
&  \Longleftrightarrow\ \left(  \text{the permutation }u\text{ is horizontal
for }w\rightharpoonup T\right) \nonumber\\
&  \ \ \ \ \ \ \ \ \ \ \ \ \ \ \ \ \ \ \ \ \left(  \text{by the definition of
}\mathcal{R}\left(  w\rightharpoonup T\right)  \right) \nonumber\\
&  \Longleftrightarrow\ \left(  \text{for each }i\in\left[  n\right]  \text{,
the number }u\left(  i\right)  \text{ lies in the same row of }%
w\rightharpoonup T\text{ as }i\right) \nonumber\\
&  \ \ \ \ \ \ \ \ \ \ \ \ \ \ \ \ \ \ \ \ \left(  \text{by the definition of
a horizontal permutation}\right) \nonumber\\
&  \Longleftrightarrow\ \left(  \text{for each }j\in\left[  n\right]  \text{,
the number }u\left(  w\left(  j\right)  \right)  \text{ lies in the same row
of }w\rightharpoonup T\text{ as }w\left(  j\right)  \right) \nonumber\\
&  \ \ \ \ \ \ \ \ \ \ \ \ \ \ \ \ \ \ \ \ \left(
\begin{array}
[c]{c}%
\text{here, we have substituted }w\left(  j\right)  \text{ for }i\text{ in the
\textquotedblleft for all\textquotedblright-sentence,}\\
\text{since }w:\left[  n\right]  \rightarrow\left[  n\right]  \text{ is a
bijection}%
\end{array}
\right) \nonumber\\
&  \Longleftrightarrow\ \left(  \text{for each }j\in\left[  n\right]  \text{,
the number }w^{-1}\left(  u\left(  w\left(  j\right)  \right)  \right)  \text{
lies in the same row of }T\text{ as }j\right) \nonumber\\
&  \ \ \ \ \ \ \ \ \ \ \ \ \ \ \ \ \ \ \ \ \left(
\begin{array}
[c]{c}%
\text{indeed, }w\rightharpoonup T\text{ is obtained from }T\text{ by applying
}w\text{ to each entry;}\\
\text{thus, the row of }w\rightharpoonup T\text{ that contains }w\left(
j\right)  \text{ is precisely the}\\
\text{row of }T\text{ that contains }j\text{, whereas the row of
}w\rightharpoonup T\text{ that}\\
\text{contains }u\left(  w\left(  j\right)  \right)  \text{ is precisely the
row of }T\text{ that contains }w^{-1}\left(  u\left(  w\left(  j\right)
\right)  \right)
\end{array}
\right) \nonumber\\
&  \Longleftrightarrow\ \left(  \text{for each }i\in\left[  n\right]  \text{,
the number }w^{-1}\left(  u\left(  w\left(  i\right)  \right)  \right)  \text{
lies in the same row of }T\text{ as }i\right) \nonumber\\
&  \ \ \ \ \ \ \ \ \ \ \ \ \ \ \ \ \ \ \ \ \left(  \text{here, we renamed the
index }j\text{ as }i\right) \nonumber\\
&  \Longleftrightarrow\ \left(  \text{for each }i\in\left[  n\right]  \text{,
the number }\left(  w^{-1}uw\right)  \left(  i\right)  \text{ lies in the same
row of }T\text{ as }i\right) \nonumber\\
&  \ \ \ \ \ \ \ \ \ \ \ \ \ \ \ \ \ \ \ \ \left(  \text{since }w^{-1}\left(
u\left(  w\left(  i\right)  \right)  \right)  =\left(  w^{-1}uw\right)
\left(  i\right)  \text{ for each }i\in\left[  n\right]  \right) \nonumber\\
&  \Longleftrightarrow\ \left(  \text{the permutation }w^{-1}uw\text{ is
horizontal for }T\right) \nonumber\\
&  \ \ \ \ \ \ \ \ \ \ \ \ \ \ \ \ \ \ \ \ \left(  \text{by the definition of
a horizontal permutation}\right) \nonumber\\
&  \Longleftrightarrow\ \left(  w^{-1}uw\in\mathcal{R}\left(  T\right)
\right)  \ \ \ \ \ \ \ \ \ \ \left(  \text{by the definition of }%
\mathcal{R}\left(  T\right)  \right) \nonumber\\
&  \Longleftrightarrow\ \left(  u\in w\mathcal{R}\left(  T\right)
w^{-1}\right)  . \label{pf.prop.tableau.Sn-act.1.b.1}%
\end{align}

But $\mathcal{R}\left(  w\rightharpoonup T\right)  $ is a subset of $S_{n}$.
Hence,%
\begin{align*}
\mathcal{R}\left(  w\rightharpoonup T\right)   &  =\left\{  u\in S_{n}%
\ \mid\ u\in\mathcal{R}\left(  w\rightharpoonup T\right)  \right\} \\
&  =\left\{  u\in S_{n}\ \mid\ u\in w\mathcal{R}\left(  T\right)
w^{-1}\right\}  \ \ \ \ \ \ \ \ \ \ \left(  \text{by the equivalence
(\ref{pf.prop.tableau.Sn-act.1.b.1})}\right) \\
&  =w\mathcal{R}\left(  T\right)  w^{-1}\ \ \ \ \ \ \ \ \ \ \left(
\text{since }w\mathcal{R}\left(  T\right)  w^{-1}\text{ is a subset of }%
S_{n}\right)  .
\end{align*}
An analogous argument (using columns instead of rows) shows that
$\mathcal{C}\left(  w\rightharpoonup T\right)  =w\mathcal{C}\left(  T\right)
w^{-1}$. Thus, Proposition \ref{prop.tableau.Sn-act.0} \textbf{(c)} is proved.
\end{proof}
\end{fineprint}

\begin{proposition}
\label{prop.tableau.Sn-act.1}Let $D$ be any diagram. Let $T$ be an $n$-tableau
of shape $D$. Then: \medskip

\textbf{(a)} We have $\mathcal{R}\left(  T\right)  =\left\{  w\in S_{n}%
\ \mid\ \overline{w\rightharpoonup T}=\overline{T}\right\}  $. \medskip

\textbf{(b)} Any $n$-tableau of shape $D$ can be written as $w\rightharpoonup
T$ for some $w\in S_{n}$. \medskip

\textbf{(c)} Any $n$-tabloid of shape $D$ can be written as $w\rightharpoonup
\overline{T}$ for some $w\in S_{n}$. \medskip

\textbf{(d)} The map%
\begin{align*}
f:S_{n}  &  \rightarrow\left\{  n\text{-tableaux of shape }D\right\}  ,\\
w  &  \mapsto w\rightharpoonup T
\end{align*}
is a bijection.
\end{proposition}

\begin{fineprint}
\begin{proof}
\textbf{(a)} Let $w\in S_{n}$. Then, we have the equivalence%
\[
\left(  w\in\mathcal{R}\left(  T\right)  \right)  \ \Longleftrightarrow
\ \left(  w\rightharpoonup T\text{ is row-equivalent to }T\right)
\ \ \ \ \ \ \ \ \ \ \left(  \text{by Proposition \ref{prop.tableau.Sn-act.0}
\textbf{(a)}}\right)  .
\]
But an $n$-tabloid is defined as an equivalence class of $n$-tableaux with
respect to row-equivalence. Thus, the $n$-tableaux $w\rightharpoonup T$ and
$T$ are row-equivalent if and only if their equivalence classes $\overline
{w\rightharpoonup T}$ and $\overline{T}$ are the same. In other words, the
equivalence
\[
\left(  w\rightharpoonup T\text{ is row-equivalent to }T\right)
\ \Longleftrightarrow\ \left(  \overline{w\rightharpoonup T}=\overline
{T}\right)
\]
holds. Altogether, we thus have proved the chain of equivalences%
\begin{align}
\left(  w\in\mathcal{R}\left(  T\right)  \right)  \  &  \Longleftrightarrow
\ \left(  w\rightharpoonup T\text{ is row-equivalent to }T\right) \nonumber\\
&  \Longleftrightarrow\ \left(  \overline{w\rightharpoonup T}=\overline
{T}\right)  . \label{pf.prop.tableau.Sn-act.1.a.1}%
\end{align}

Now, forget that we fixed $w$. We thus have proved the equivalence
(\ref{pf.prop.tableau.Sn-act.1.a.1}) for each $w\in S_{n}$.

But $\mathcal{R}\left(  T\right)  $ is a subset (and in fact a subgroup) of
$S_{n}$. Thus,%
\begin{align*}
\mathcal{R}\left(  T\right)   &  =\left\{  w\in S_{n}\ \mid\ w\in
\mathcal{R}\left(  T\right)  \right\} \\
&  =\left\{  w\in S_{n}\ \mid\ \overline{w\rightharpoonup T}=\overline
{T}\right\}  \ \ \ \ \ \ \ \ \ \ \left(  \text{by the equivalence
(\ref{pf.prop.tableau.Sn-act.1.a.1})}\right)  .
\end{align*}

This proves Proposition \ref{prop.tableau.Sn-act.1} \textbf{(a)}. \medskip

\textbf{(b)} Let $U$ be any $n$-tableau of shape $D$. We must show that
$U=w\rightharpoonup T$ for some $w\in S_{n}$.

Intuitively, this is clear: Both tableaux $T$ and $U$ contain the entries
$1,2,\ldots,n$ (each exactly once); thus, we can obtain $U$ from $T$ by
permuting these entries.

The rigorous proof is not much harder: We know that $T$ is an $n$-tableau of
shape $D$, that is, an injective map from $D$ to $\left[  n\right]  $ whose
image is $\left[  n\right]  $. Thus, we can view $T$ as a bijection from $D$
to $\left[  n\right]  $. Likewise, we can view $U$ as a bijection from $D$ to
$\left[  n\right]  $. Bijections are always invertible, so that $T^{-1}$ is a
bijection from $\left[  n\right]  $ to $D$ (since $T$ is a bijection from $D$
to $\left[  n\right]  $). Hence, $U\circ T^{-1}$ is a bijection from $\left[
n\right]  $ to $\left[  n\right]  $ (since a composition of two bijections is
a bijection). In other words, $U\circ T^{-1}$ is a permutation of $\left[
n\right]  $. That is, $U\circ T^{-1}\in S_{n}$.

Now, the definition of the action on entries yields $\left(  U\circ
T^{-1}\right)  \rightharpoonup T=U\circ\underbrace{T^{-1}\circ T}%
_{=\operatorname*{id}}=U$. Thus, $U=w\rightharpoonup T$ for some $w\in S_{n}$
(namely, for $w=U\circ T^{-1}$). This proves Proposition
\ref{prop.tableau.Sn-act.1} \textbf{(b)}. \medskip

\textbf{(c)} Let $\overline{U}$ be an $n$-tabloid of shape $D$. We must show
that $\overline{U}=w\rightharpoonup\overline{T}$ for some $w\in S_{n}$.

Proposition \ref{prop.tableau.Sn-act.1} \textbf{(b)} shows that
$U=w\rightharpoonup T$ for some $w\in S_{n}$. Using this $w$, we then have%
\begin{align*}
\overline{U}  &  =\overline{w\rightharpoonup T}\ \ \ \ \ \ \ \ \ \ \left(
\text{since }U=w\rightharpoonup T\right) \\
&  =w\rightharpoonup\overline{T}\ \ \ \ \ \ \ \ \ \ \left(  \text{by
Definition \ref{def.tabloid.Sn-act}}\right)  .
\end{align*}

Thus, Proposition \ref{prop.tableau.Sn-act.1} \textbf{(c)} is proved. \medskip

\textbf{(d)} Consider the map%
\begin{align*}
f:S_{n}  &  \rightarrow\left\{  n\text{-tableaux of shape }D\right\}  ,\\
w  &  \mapsto w\rightharpoonup T.
\end{align*}
This map is clearly well-defined. Moreover:

\begin{itemize}
\item The map $f$ is injective.

[\textit{Proof:} Let $u,v\in S_{n}$ be such that $f\left(  u\right)  =f\left(
v\right)  $. We must prove that $u=v$.

We have $f\left(  u\right)  =u\rightharpoonup T$ (by the definition of $f$).
But $u\rightharpoonup T=u\circ T$ (by Definition \ref{def.tableau.Sn-act}).
Hence, $f\left(  u\right)  =u\rightharpoonup T=u\circ T$. Similarly, $f\left(
v\right)  =v\circ T$. Thus, $u\circ T=f\left(  u\right)  =f\left(  v\right)
=v\circ T$.

But $T$ is an $n$-tableau of shape $D$. In other words, $T$ is a bijection
from $D$ to $\left[  n\right]  $. Hence, the map $T$ has an inverse
$T^{-1}:\left[  n\right]  \rightarrow D$. Clearly, $u\circ\underbrace{T\circ
T^{-1}}_{=\operatorname*{id}}=u$, so that $u=\underbrace{u\circ T}_{=v\circ
T}\circ\,T^{-1}=v\circ\underbrace{T\circ T^{-1}}_{=\operatorname*{id}}=v$.

Forget that we fixed $u,v$. We thus have shown that if $u,v\in S_{n}$ are such
that $f\left(  u\right)  =f\left(  v\right)  $, then $u=v$. In other words,
the map $f$ is injective.]

\item The map $f$ is surjective.

[\textit{Proof:} Let $Q\in\left\{  n\text{-tableaux of shape }D\right\}  $.
Thus, $Q$ is an $n$-tableau of shape $D$. But Proposition
\ref{prop.tableau.Sn-act.1} \textbf{(b)} shows that any $n$-tableau of shape
$D$ can be written as $w\rightharpoonup T$ for some $w\in S_{n}$. In
particular, $Q$ can be written as $w\rightharpoonup T$ for some $w\in S_{n}$
(since $Q$ is an $n$-tableau of shape $D$). In other words,
$Q=w\rightharpoonup T$ for some $w\in S_{n}$. Consider this $w$.

Now, the definition of $f$ yields $f\left(  w\right)  =w\rightharpoonup T=Q$.
Hence, $Q=f\left(  w\right)  \in f\left(  S_{n}\right)  $.

Forget that we fixed $Q$. We thus have shown that $Q\in f\left(  S_{n}\right)
$ for every $Q\in\left\{  n\text{-tableaux of shape }D\right\}  $. In other
words, $\left\{  n\text{-tableaux of shape }D\right\}  \subseteq f\left(
S_{n}\right)  $. This shows that $f$ is surjective.]
\end{itemize}

Now, we have shown that $f$ is both injective and surjective. Hence, $f$ is
bijective, i.e., is a bijection. This proves Proposition
\ref{prop.tableau.Sn-act.1} \textbf{(d)}.
\end{proof}
\end{fineprint}

For $n$-tableaux, we can also describe row-equivalence in terms of horizontal
permutations, and likewise for column-equivalence:

\begin{proposition}
\label{prop.tableau.req-R}Let $D$ be a diagram. Let $T$ and $S$ be two
$n$-tableaux of shape $D$. Then: \medskip

\textbf{(a)} The $n$-tableaux $T$ and $S$ are row-equivalent if and only if
there exists some $w\in\mathcal{R}\left(  T\right)  $ such that
$S=w\rightharpoonup T$. \medskip

\textbf{(b)} The $n$-tableaux $T$ and $S$ are column-equivalent if and only if
there exists some $w\in\mathcal{C}\left(  T\right)  $ such that
$S=w\rightharpoonup T$.
\end{proposition}

\begin{fineprint}
\begin{proof}
\textbf{(a)} We prove the \textquotedblleft$\Longleftarrow$\textquotedblright%
\ and \textquotedblleft$\Longrightarrow$\textquotedblright\ directions separately:

$\Longleftarrow:$ Assume that there exists some $w\in\mathcal{R}\left(
T\right)  $ such that $S=w\rightharpoonup T$. Consider this $w$. The tableau
$w\rightharpoonup T$ is row-equivalent to $T$ (by Proposition
\ref{prop.tableau.Sn-act.0} \textbf{(a)}, since $w\in\mathcal{R}\left(
T\right)  $). In other words, $S$ is row-equivalent to $T$ (since
$S=w\rightharpoonup T$). In other words, $T$ and $S$ are row-equivalent. This
proves the \textquotedblleft$\Longleftarrow$\textquotedblright\ direction of
Proposition \ref{prop.tableau.req-R} \textbf{(a)}.

$\Longrightarrow:$ Assume that the tableaux $T$ and $S$ are row-equivalent.
Proposition \ref{prop.tableau.Sn-act.1} \textbf{(b)} shows that the tableau
$S$ can be written as $w\rightharpoonup T$ for some $w\in S_{n}$. Consider
this $w$.

We assumed that the tableaux $T$ and $S$ are row-equivalent. In other words,
the tableaux $T$ and $w\rightharpoonup T$ are row-equivalent (since
$S=w\rightharpoonup T$). By Proposition \ref{prop.tableau.Sn-act.0}
\textbf{(a)}, this entails that $w\in\mathcal{R}\left(  T\right)  $. Thus, we
have found a $w\in\mathcal{R}\left(  T\right)  $ such that $S=w\rightharpoonup
T$. Hence, such a $w$ exists. This proves the \textquotedblleft%
$\Longrightarrow$\textquotedblright\ direction of Proposition
\ref{prop.tableau.req-R} \textbf{(a)}. \medskip

\textbf{(b)} This is analogous to part \textbf{(a)}.
\end{proof}
\end{fineprint}

\begin{exercise}
\fbox{2} Let $D$ be a diagram with $\left\vert D\right\vert =n$. Fix an
$n$-tableau $T$ of shape $D$. Then, its row group $\mathcal{R}\left(
T\right)  $ is a subgroup of $S_{n}$, and thus there is a left $S_{n}$-set
$S_{n}/\mathcal{R}\left(  T\right)  $ consisting of all left cosets of
$\mathcal{R}\left(  T\right)  $ in $S_{n}$ (see Example
\ref{exa.rep.G-sets.G/H} for its definition). On the other hand, the set
$\left\{  n\text{-tabloids of shape }D\right\}  $ is a left $S_{n}$-set as
well. Prove that these two left $S_{n}$-sets are isomorphic. To be more
specific: Prove that the map%
\begin{align*}
S_{n}/\mathcal{R}\left(  T\right)   &  \rightarrow\left\{  n\text{-tabloids of
shape }D\right\}  ,\\
w\mathcal{R}\left(  T\right)   &  \mapsto\overline{w\rightharpoonup T}%
\end{align*}
is a left $S_{n}$-set isomorphism.
\end{exercise}

\subsubsection{Young modules}

Now, recall the concept of a permutation module (defined in Example
\ref{exa.mod.kG-on-kX}), which automatically turns any left $G$-set $X$ (where
$G$ is a group) into a left $\mathbf{k}\left[  G\right]  $-module
$\mathbf{k}^{\left(  X\right)  }$ by linearization (i.e., by extending the
left $G$-action $G\times X\rightarrow X$ to a $\mathbf{k}$-bilinear map
$\mathbf{k}\left[  G\right]  \times\mathbf{k}^{\left(  X\right)  }%
\rightarrow\mathbf{k}^{\left(  X\right)  }$ in the only reasonable way). We
shall now use this technique to turn the left $S_{n}$-set $\left\{
n\text{-tabloids of shape }D\right\}  $ (as defined in Definition
\ref{def.tabloid.Sn-act}) into a left $\mathbf{k}\left[  S_{n}\right]
$-module, which is known as a \emph{Young module}:

\begin{definition}
\label{def.youngmod.youngmod}Let $D$ be a diagram with $\left\vert
D\right\vert =n$. The \emph{Young module} $\mathcal{M}^{D}$ (sometimes also
denoted by $\mathcal{Y}^{D}$) is the permutation module corresponding to the
left $S_{n}$-set%
\[
\left\{  n\text{-tabloids of shape }D\right\}
\]
(which was introduced in Definition \ref{def.tabloid.Sn-act}). As a
$\mathbf{k}$-module, it is the free $\mathbf{k}$-module $\mathbf{k}^{\left(
\left\{  n\text{-tabloids of shape }D\right\}  \right)  }$ on the set
$\left\{  n\text{-tabloids of shape }D\right\}  $. It has a standard basis
vector $e_{\overline{T}}$ for each $n$-tabloid $\overline{T}$ of shape $D$,
but we shall just denote this basis vector by $\overline{T}$ by the usual
abuse of notation. An arbitrary element of $\mathcal{M}^{D}$ is thus a formal
$\mathbf{k}$-linear combination of $n$-tabloids of shape $D$. The symmetric
group $S_{n}$ acts on $\mathcal{M}^{D}$ by permuting these $n$-tabloids (i.e.,
the respective standard basis vectors).
\end{definition}

\begin{example}
\label{exa.youngmod.n-1-1}Assume that $n\geq2$. Let $D$ be the Young diagram
$Y\left(  n-1,1\right)  $. This Young diagram has two rows of lengths $n-1$
and $1$, respectively.

As we know from Example \ref{exa.tabloid.n-1-1} \textbf{(b)}, there are $n$
many $n$-tabloids of shape $D$. We can label these $n$-tabloids as
$\underline{\overline{1}},\ \underline{\overline{2}},\ \ldots
,\ \underline{\overline{n}}$, where $\underline{\overline{i}}$ means the
$n$-tabloid whose second row has entry $i$ and whose first row contains all
other elements of $\left[  n\right]  $.

Thus, the Young module $\mathcal{M}^{D}$ has basis $\left(
\underline{\overline{1}},\ \underline{\overline{2}},\ \ldots
,\ \underline{\overline{n}}\right)  $ as a $\mathbf{k}$-module. The action of
$S_{n}$ on this module is given by%
\[
w\ \underline{\overline{i}}=\underline{\overline{w\left(  i\right)  }%
}\ \ \ \ \ \ \ \ \ \ \text{for each }w\in S_{n}\text{ and each }i\in\left[
n\right]  .
\]
For instance, for $n=4$, we have%
\[
\operatorname*{cyc}\nolimits_{1,2,3}\ \underline{\overline{2}}%
=\operatorname*{cyc}\nolimits_{1,2,3}%
\ytableausetup{tabloids}\ytableaushort{134,2}\ytableausetup{notabloids}=\ytableausetup{tabloids}\ytableaushort{214,3}\ytableausetup{notabloids}=\ytableausetup{tabloids}\ytableaushort{124,3}\ytableausetup{notabloids}=\underline{\overline
{3}}=\underline{\overline{\operatorname*{cyc}\nolimits_{1,2,3}\left(
2\right)  }}.
\]

The Young module $\mathcal{M}^{D}$ is just the natural $S_{n}$-representation
$\mathbf{k}^{n}$ in disguise. Indeed, the $\mathbf{k}$-linear map%
\begin{align*}
\mathbf{k}^{n}  &  \rightarrow\mathcal{M}^{D},\\
e_{i}  &  \mapsto\underline{\overline{i}}%
\end{align*}
is an isomorphism of $S_{n}$-representations (since any $w\in S_{n}$ sends
$e_{i}$ to $e_{w\left(  i\right)  }$ and sends $\underline{\overline{i}}$ to
$\underline{\overline{w\left(  i\right)  }}$).
\end{example}

\begin{example}
\label{exa.youngmod.trivrows}Let $D$ be a diagram that has no two cells in the
same row and satisfies $\left\vert D\right\vert =n$. (For example, we can take
$D$ to be the Young diagram $Y\left(  \underbrace{1,1,\ldots,1}_{n\text{
times}}\right)  $, or the skew Young diagram $Y\left(  \left(  n,n-1,\ldots
,1\right)  /\left(  n-1,n-2,\ldots,1\right)  \right)  $. For $n=4$, these two
diagrams are respectively $\ydiagram{1,1,1,1}$ and
$\ydiagram{3+1,2+1,1+1,0+1}$ .)

Then, the $n$-tabloids of shape $D$ are \textquotedblleft
essentially\textquotedblright\ the permutations $w\in S_{n}$. Indeed, no two
distinct $n$-tableaux are row-equivalent (since no row has more than one
entry). Thus, each $n$-tabloid consists of just one $n$-tableau, and can be
identified with this one $n$-tableau that it contains. But there are $n!$ many
$n$-tableaux, since we can easily find a bijection from $S_{n}$ to the set of
all $n$-tableaux (by choosing an arbitrary $n$-tableau $T$, and sending each
$w\in S_{n}$ to the $n$-tableau $w\rightharpoonup T$).

The Young module $\mathcal{M}^{D}$ (for such a diagram $D$) is thus just the
left regular representation $\mathbf{k}\left[  S_{n}\right]  $ of $S_{n}$.
More precisely, if we fix any $n$-tableau $T$ of shape $D$, then there is a
left $\mathbf{k}\left[  S_{n}\right]  $-module isomorphism%
\begin{align*}
\mathbf{k}\left[  S_{n}\right]   &  \rightarrow\mathcal{M}^{D},\\
w  &  \mapsto\overline{w\rightharpoonup T}=w\rightharpoonup\overline{T}.
\end{align*}

\end{example}

\begin{exercise}
\label{exe.youngmod.1dim-sub}\fbox{2} Assume that $\mathbf{k}$ is a field. Let
$D$ be a diagram that contains at least one row having more than one cell. Let
$s\in\mathcal{M}^{D}$ be the sum of all $n$-tabloids of shape $D$ (or, rather,
of the corresponding standard basis vectors). Recall from Exercise
\ref{exe.mod.kX-on-kX.trivsub} that the $1$-dimensional submodule
$\operatorname*{span}\nolimits_{\mathbf{k}}\left\{  s\right\}  $ is a trivial
subrepresentation of the Young module $\mathcal{M}^{D}$.

Prove that this submodule is the only $1$-dimensional subrepresentation of
$\mathcal{M}^{D}$. \medskip

[\textbf{Hint:} Any $1$-dimensional subrepresentation of $\mathcal{M}^{D}$ is
spanned by a single element, which we can write as $\sum_{\overline{T}}%
\alpha_{\overline{T}}\overline{T}$. Each $w\in S_{n}$ must multiply this
element by a scalar $\kappa_{w}$. This gives $\alpha_{\overline{T}}=\kappa
_{w}\alpha_{\overline{wT}}$ for all $T$ and all $w$. This shows that each
$\alpha_{\overline{T}}$ must be nonzero (why?), but what else does it show?]
\end{exercise}

\subsection{Specht modules}

\subsubsection{Definition and examples}

Next, we will define a certain subrepresentation of each Young module, which
is called a \emph{Specht module} and will later turn out to be a basic
building block of the representation theory of $S_{n}$.

\begin{definition}
\label{def.spechtmod.spechtmod}Let $D$ be a diagram with $\left\vert
D\right\vert =n$. \medskip

\textbf{(a)} If $T$ is an $n$-tableau of shape $D$, then we define the
so-called \emph{polytabloid}%
\[
\mathbf{e}_{T}:=\sum_{w\in\mathcal{C}\left(  T\right)  }\left(  -1\right)
^{w}\overline{wT}\in\mathcal{M}^{D}.
\]

\textbf{(b)} The span of these elements $\mathbf{e}_{T}$, where $T$ ranges
over all $n$-tableaux of shape $D$, is called the \emph{Specht module}
$\mathcal{S}^{D}$. Thus,%
\[
\mathcal{S}^{D}=\operatorname*{span}\nolimits_{\mathbf{k}}\left\{
\mathbf{e}_{T}\ \mid\ T\text{ is an }n\text{-tableau of shape }D\right\}  .
\]
As we will soon see (in Lemma \ref{lem.spechtmod.submod} \textbf{(b)}), this
Specht module is a left $\mathbf{k}\left[  S_{n}\right]  $-submodule (i.e., an
$S_{n}$-subrepresentation) of the Young module $\mathcal{M}^{D}$. \medskip

\textbf{Warning:} Do not mistake the polytabloid $\mathbf{e}_{T}$ for the
standard basis vector corresponding to the tabloid $\overline{T}$ (which we
just call $\overline{T}$). The letter \textquotedblleft$\mathbf{e}%
$\textquotedblright\ in \textquotedblleft$\mathbf{e}_{T}$\textquotedblright%
\ has nothing to do with the notation $e_{i}$ for a standard basis vector.

Also, keep in mind that tabloids are combinatorial objects (sets of tableaux),
while polytabloids are algebraic ones (elements of the $S_{n}$-representation
$\mathcal{M}^{D}$).
\end{definition}

\begin{example}
\label{exa.spechtmod.Yn-1-1}Let $n=3$ and $D=Y\left(  2,1\right)  $. Then, the
$n$-tableaux of shape $D$ are%
\[
\ytableaushort{12,3}\ \ \ \ \ \ \ \ \ \ \ytableaushort{13,2}\ \ \ \ \ \ \ \ \ \ \ytableaushort{21,3}\ \ \ \ \ \ \ \ \ \ \ytableaushort{23,1}\ \ \ \ \ \ \ \ \ \ \ytableaushort{31,2}\ \ \ \ \ \ \ \ \ \ \ytableaushort{32,1}.
\]
Following Convention \ref{conv.tableau.poetic}, we can just call them
$12\backslash\backslash3$ and $13\backslash\backslash2$ and so on.

The $n$-tabloids of shape $D$ are%
\[
\ytableausetup{tabloids}\ytableaushort{12,3}\ytableausetup{notabloids}\ \ \ \ \ \ \ \ \ \ \ \ytableausetup{tabloids}\ytableaushort{13,2}\ytableausetup{notabloids}\ \ \ \ \ \ \ \ \ \ \ \ytableausetup{tabloids}\ytableaushort{23,1}\ytableausetup{notabloids}\ \ \ ,
\]
aka $\overline{12\backslash\backslash3}$ and $\overline{13\backslash
\backslash2}$ and $\overline{23\backslash\backslash1}$. Each $n$-tabloid
contains two $n$-tableaux.

Each of the six $n$-tableaux $T$ of shape $D$ gives rise to a polytabloid.
Here are two of them:
\begin{align*}
\mathbf{e}_{12\backslash\backslash3}  &  =\sum_{w\in\mathcal{C}\left(
12\backslash\backslash3\right)  }\left(  -1\right)  ^{w}\overline{w\left(
12\backslash\backslash3\right)  }\\
&  =\underbrace{\left(  -1\right)  ^{\operatorname*{id}}}_{=1}%
\underbrace{\overline{\operatorname*{id}\left(  12\backslash\backslash
3\right)  }}_{=\overline{12\backslash\backslash3}}+\underbrace{\left(
-1\right)  ^{t_{1,3}}}_{=-1}\underbrace{\overline{t_{1,3}\left(
12\backslash\backslash3\right)  }}_{=\overline{32\backslash\backslash1}}\\
&  \ \ \ \ \ \ \ \ \ \ \ \ \ \ \ \ \ \ \ \ \left(  \text{since }%
\mathcal{C}\left(  12\backslash\backslash3\right)  =\left\{
\operatorname*{id},t_{1,3}\right\}  \right) \\
&  =\overline{12\backslash\backslash3}-\underbrace{\overline{32\backslash
\backslash1}}_{=\overline{23\backslash\backslash1}}=\overline{12\backslash
\backslash3}-\overline{23\backslash\backslash1}%
\end{align*}
and%
\begin{align*}
\mathbf{e}_{21\backslash\backslash3}  &  =\sum_{w\in\mathcal{C}\left(
21\backslash\backslash3\right)  }\left(  -1\right)  ^{w}\overline{w\left(
21\backslash\backslash3\right)  }\\
&  =\left(  -1\right)  ^{\operatorname*{id}}\overline{\operatorname*{id}%
\left(  21\backslash\backslash3\right)  }+\left(  -1\right)  ^{t_{2,3}%
}\overline{t_{2,3}\left(  21\backslash\backslash3\right)  }\\
&  \ \ \ \ \ \ \ \ \ \ \ \ \ \ \ \ \ \ \ \ \left(  \text{since }%
\mathcal{C}\left(  21\backslash\backslash3\right)  =\left\{
\operatorname*{id},t_{2,3}\right\}  \right) \\
&  =\overline{21\backslash\backslash3}-\overline{31\backslash\backslash
2}=\overline{12\backslash\backslash3}-\overline{13\backslash\backslash2};
\end{align*}
and there are four more. Rather than compute them one by one, we can just
compute them all at once: For any distinct elements $i,j,k\in\left[  3\right]
$, we have%
\begin{align*}
\mathbf{e}_{ij\backslash\backslash k}  &  =\sum_{w\in\mathcal{C}\left(
ij\backslash\backslash k\right)  }\left(  -1\right)  ^{w}\overline{w\left(
ij\backslash\backslash k\right)  }\\
&  =\underbrace{\left(  -1\right)  ^{\operatorname*{id}}}_{=1}%
\underbrace{\overline{\operatorname*{id}\left(  ij\backslash\backslash
k\right)  }}_{=\overline{ij\backslash\backslash k}}+\underbrace{\left(
-1\right)  ^{t_{i,k}}}_{=-1}\underbrace{\overline{t_{i,k}\left(
ij\backslash\backslash k\right)  }}_{=\overline{kj\backslash\backslash i}}\\
&  \ \ \ \ \ \ \ \ \ \ \ \ \ \ \ \ \ \ \ \ \left(  \text{since }%
\mathcal{C}\left(  ij\backslash\backslash k\right)  =\left\{
\operatorname*{id},t_{i,k}\right\}  \right) \\
&  =\overline{ij\backslash\backslash k}-\overline{kj\backslash\backslash i}.
\end{align*}
Using the notation $\underline{\overline{i}}$ from Example
\ref{exa.youngmod.n-1-1}, we can rewrite this as%
\[
\mathbf{e}_{ij\backslash\backslash k}=\underline{\overline{k}}%
-\underline{\overline{i}}%
\]
(since $\underline{\overline{k}}=\overline{ij\backslash\backslash k}$ and
$\underline{\overline{i}}=\overline{kj\backslash\backslash i}$).

Thus, the Specht module $\mathcal{S}^{D}$ is spanned by the differences
$\underline{\overline{k}}-\underline{\overline{i}}$ for all triples of
distinct elements $i,j,k\in\left[  3\right]  $. In other words, it is spanned
by the differences $\underline{\overline{k}}-\underline{\overline{i}}$ for all
pairs of distinct elements $k\neq i$ in $\left[  3\right]  $. Identifying the
Young module $\mathcal{M}^{D}$ with the natural representation $\mathbf{k}%
^{n}$ of $S_{n}$ for $n=3$ (as explained in Example \ref{exa.youngmod.n-1-1}),
we can thus view $\mathcal{S}^{D}$ as the span of the differences $e_{k}%
-e_{i}$ of distinct standard basis vectors. But, as we know, this span is the
zero-sum subrepresentation $R\left(  \mathbf{k}^{n}\right)  $.

The same reasoning works not just for $n=3$ but for any $n\geq2$. Thus, we can
see that the Specht module $\mathcal{S}^{Y\left(  n-1,1\right)  }$ is
isomorphic to the zero-sum subrepresentation $R\left(  \mathbf{k}^{n}\right)
$ of the natural $S_{n}$-representation $\mathbf{k}^{n}$.
\end{example}

\begin{example}
\label{exa.spechtmod.Yn}Now consider the partition $\left(  n\right)  $ (by
which we understand the empty partition if $n=0$). Let $D=Y\left(  n\right)
$. Then, all cells of $D$ lie in the same row. Hence, any two $n$-tableaux of
shape $D$ are row-equivalent. Thus, there is only one $n$-tabloid, namely
$\ytableausetup{tabloids} \ytableaushort{12\cdots n}
\ytableausetup{notabloids}$. The left $S_{n}$-action on the set of this one
$n$-tabloid is trivial, and thus the Young module $\mathcal{M}^{Y\left(
n\right)  }$ is the trivial $S_{n}$-representation. The Specht module
$\mathcal{S}^{Y\left(  n\right)  }$ is also the trivial $S_{n}$%
-representation, since $\mathcal{C}\left(  T\right)  =\left\{
\operatorname*{id}\right\}  $ for any $n$-tableau $T$ of shape $D$.
\end{example}

\begin{example}
\label{exa.spechtmod.Y111}Now consider the partition $\left(  1^{n}\right)  $
(this is exponential notation for $\underbrace{\left(  1,1,\ldots,1\right)
}_{n\text{ times}}$, as explained in Definition \ref{def.partitions.expnot}).
Let $D=Y\left(  1^{n}\right)  $. Then, the diagram $D$ consists of a single
column, so it has no two cells in the same row. As we saw in Example
\ref{exa.youngmod.trivrows}, this entails that the $n$-tabloids of shape $D$
are all the permutations $w\in S_{n}$ (encoded as the respective $n$-tabloids
$\ytableausetup{tabloids, boxsize=2.5em} \begin{ytableau}
w\left(1\right) \\
w\left(2\right) \\
\vdots \\
w\left(n\right)
\end{ytableau}
\ytableausetup{notabloids, boxsize=normal}$), and the Young module
$\mathcal{M}^{Y\left(  1^{n}\right)  }$ is (isomorphic to) the left regular
representation $\mathbf{k}\left[  S_{n}\right]  $ of $S_{n}$.

What about the Specht module $\mathcal{S}^{Y\left(  1^{n}\right)  }$ ? It is
spanned by the polytabloids%
\[
\mathbf{e}_{T}=\sum_{w\in\mathcal{C}\left(  T\right)  }\left(  -1\right)
^{w}\overline{wT}=\sum_{w\in S_{n}}\left(  -1\right)  ^{w}\overline
{wT}\ \ \ \ \ \ \ \ \ \ \left(  \text{since }\mathcal{C}\left(  T\right)
=S_{n}\text{ here}\right)
\]
for all $n$-tableaux $T$ of shape $Y\left(  1^{n}\right)  $. These look like a
lot of polytabloids ($n!$ to be precise), but I claim that they are all equal
up to sign: Indeed, if we identify $\mathcal{M}^{Y\left(  1^{n}\right)  }$
with the left regular representation $\mathbf{k}\left[  S_{n}\right]  $, then
each of the polytabloids $\mathbf{e}_{T}$ is either the sign-integral
$\nabla^{-}$ or $-\nabla^{-}$.

The easiest way to see this is by identifying each permutation $w\in S_{n}$
with the $n$-tableau $T_{w}=\ytableausetup{boxsize=3em}\begin{ytableau}
w\left(1\right) \\
w\left(2\right) \\
\vdots \\
w\left(n\right)
\end{ytableau}\ytableausetup{boxsize=normal}$ and the corresponding
$n$-tabloid $\overline{T_{w}}$. Then, conversely, any $n$-tableau $T$ can be
viewed as a permutation in $S_{n}$ (namely, as the permutation $w\in S_{n}$
for which $T$ is $T_{w}$). Hence, for any $n$-tableau $T$, we have
\begin{align*}
\mathbf{e}_{T}  &  =\sum_{w\in S_{n}}\left(  -1\right)  ^{w}%
\underbrace{\overline{wT}}_{\substack{=wT\\\text{(where we regard }T\\\text{as
a permutation in }S_{n}\text{)}}}=\underbrace{\sum_{w\in S_{n}}\left(
-1\right)  ^{w}w}_{\substack{=\nabla^{-}\\\text{(by the definition of }%
\nabla^{-}\text{)}}}T=\nabla^{-}T\\
&  =\left(  -1\right)  ^{T}\nabla^{-}\ \ \ \ \ \ \ \ \ \ \left(  \text{by
(\ref{eq.prop.integral.fix.-}), applied to }T\text{ instead of }w\right)  .
\end{align*}
Therefore, $\mathbf{e}_{T}$ is $\nabla^{-}$ or $-\nabla^{-}$.

Thus, as a left $\mathbf{k}\left[  S_{n}\right]  $-module, the Specht module
$\mathcal{S}^{Y\left(  1^{n}\right)  }$ is isomorphic to the ($\mathbf{k}%
$-linear) span of $\nabla^{-}$. It is furthermore isomorphic to the sign
representation of $S_{n}$, since every $w\in S_{n}$ satisfies $w\nabla
^{-}=\left(  -1\right)  ^{w}\nabla^{-}$.
\end{example}

\begin{example}
\label{exa.spechtmod.lrr}Let $D$ be a diagram that has no two cells in the
same row and has no two cells in the same column and satisfies $\left\vert
D\right\vert =n$. (For example, $D$ can be the skew Young diagram $Y\left(
\left(  n,n-1,\ldots,1\right)  /\left(  n-1,n-2,\ldots,1\right)  \right)  $,
as in Example \ref{exa.youngmod.trivrows}. But $D$ cannot be a straight Young
diagram $Y\left(  \lambda\right)  $, unless $n\leq1$.)

As we know from Example \ref{exa.youngmod.trivrows}, the Young module
$\mathcal{M}^{D}$ is (an isomorphic copy of) the left regular representation
$\mathbf{k}\left[  S_{n}\right]  $ of $S_{n}$, since each $n$-tabloid of shape
$D$ contains just one $n$-tableau. Furthermore, since $D$ has no two cells in
the same column, each $n$-tableau $T$ of shape $D$ satisfies $\mathcal{C}%
\left(  T\right)  =\left\{  \operatorname*{id}\right\}  $, and thus the
polytabloids $\mathbf{e}_{T}$ are just the corresponding tabloids
$\overline{T}$. Hence, the Specht module $\mathcal{S}^{D}$ is the span of all
$n$-tabloids $\overline{T}$, which is of course the whole Young module
$\mathcal{M}^{D}$. Thus, we have%
\[
\mathcal{S}^{D}=\mathcal{M}^{D}\cong\mathbf{k}\left[  S_{n}\right]
\]
in this case.
\end{example}

\begin{fineprint}
Some remarks on software implementations are in order. The SageMath CAS has
\href{https://doc.sagemath.org/html/en/reference/combinat/sage/combinat/specht_module.html}{some
support for working in Specht modules} $\mathcal{S}^{D}$ (at least when
$\mathbf{k}$ is a field). This can be accessed as follows (shown here for the
diagram $D$ from Exercise \ref{exe.spechtmod.bad.6}):

\qquad\texttt{SGA =\ SymmetricGroupAlgebra(QQ, 6)}

\qquad\texttt{D = [(1,2), (1,3), (2,1), (2,3), (3,1), (3,2)]}

\qquad\texttt{SD = SGA.specht\_module(D)}

When $D$ is a Young diagram, it is better to access $\mathcal{S}^{D}$ as follows:

\qquad\texttt{SGA =\ SymmetricGroupAlgebra(QQ, 6)}

\qquad\texttt{SD = SGA.specht\_module(Partition([3,2,1]))}

\noindent as this constructs the basis consisting of the standard polytabloids
(see Theorem \ref{thm.spechtmod.basis}) rather than just a nondescript basis
computed using linear algebra. When $\mathcal{S}^{D}$ is constructed in this
way, the standard polytabloids $\mathbf{e}_{T}$ can be called in an obvious
way: e.g., \texttt{SD(Tableau([[1,3,4],[2,5],[6]]))} yields $\mathbf{e}%
_{134\backslash\backslash25\backslash\backslash6}$. At the present moment,
non-standard polytabloids cannot be constructed this way; instead, write $T$
as $wP$ for some standard $P$ and then compute $\mathbf{e}_{T}$ as
$w\mathbf{e}_{P}$.

However, at the present state (mid-2025), all this functionality is in
development and far from optimized. For the Young diagram case (that is,
$D=Y\left(  \lambda\right)  $) over the field $\mathbb{Q}$, I thus recommend
Macaulay2 -- a different CAS -- and
\href{https://macaulay2.com/doc/Macaulay2/share/doc/Macaulay2/SpechtModule/html/index.html}{its
SpechtModule package} instead. The Specht--Vandermonde avatar (Theorem
\ref{thm.spechtmod.vdm}) can also be used to reduce computations in Young or
Specht modules to computations in a polynomial ring, which are supported by
many more software packages.
\end{fineprint}

\subsubsection{Properties}

It is time to state some general properties of polytabloids and Specht modules:

\begin{lemma}
\label{lem.spechtmod.submod}Let $D$ be any diagram with $\left\vert
D\right\vert =n$. Then: \medskip

\textbf{(a)} If $T$ is any $n$-tableau of shape $D$, and if $u\in S_{n}$, then%
\[
\mathbf{e}_{uT}=u\mathbf{e}_{T}.
\]

\textbf{(b)} The set $\mathcal{S}^{D}$ is a left $\mathbf{k}\left[
S_{n}\right]  $-submodule of $\mathcal{M}^{D}$. \medskip

\textbf{(c)} If $T$ is any $n$-tableau of shape $D$, and if $u\in
\mathcal{C}\left(  T\right)  $, then $\mathbf{e}_{uT}=\left(  -1\right)
^{u}\mathbf{e}_{T}$. \medskip

\textbf{(d)} If $T$ and $S$ are two column-equivalent $n$-tableaux of shape
$D$, then $\mathbf{e}_{T}=\pm\mathbf{e}_{S}$.
\end{lemma}

\begin{proof}
\textbf{(a)} Let $T$ be any $n$-tableau of shape $D$. Let $u\in S_{n}$. Then,
Proposition \ref{prop.tableau.Sn-act.0} \textbf{(c)} (applied to $w=u$) yields
$\mathcal{R}\left(  u\rightharpoonup T\right)  =u\mathcal{R}\left(  T\right)
u^{-1}$ and $\mathcal{C}\left(  u\rightharpoonup T\right)  =u\mathcal{C}%
\left(  T\right)  u^{-1}$. Since $uT$ is a shorthand for $u\rightharpoonup T$,
we can rewrite these two equalities as $\mathcal{R}\left(  uT\right)
=u\mathcal{R}\left(  T\right)  u^{-1}$ and $\mathcal{C}\left(  uT\right)
=u\mathcal{C}\left(  T\right)  u^{-1}$.

From $\mathcal{C}\left(  uT\right)  =u\mathcal{C}\left(  T\right)  u^{-1}$, we
see that the map
\begin{align*}
\mathcal{C}\left(  T\right)   &  \rightarrow\mathcal{C}\left(  uT\right)  ,\\
w  &  \mapsto uwu^{-1}%
\end{align*}
is well-defined and surjective. This map is injective as well (since $S_{n}$
is a group), and thus is a bijection.

But Definition \ref{def.spechtmod.spechtmod} \textbf{(a)} yields%
\begin{equation}
\mathbf{e}_{T}=\sum_{w\in\mathcal{C}\left(  T\right)  }\left(  -1\right)
^{w}\overline{wT} \label{pf.lem.spechtmod.submod.a.1}%
\end{equation}
and%
\begin{align*}
\mathbf{e}_{uT}  &  =\sum_{w\in\mathcal{C}\left(  uT\right)  }\left(
-1\right)  ^{w}\overline{wuT}=\sum_{w\in\mathcal{C}\left(  T\right)
}\underbrace{\left(  -1\right)  ^{uwu^{-1}}}_{\substack{=\left(  -1\right)
^{w}\\\text{(by (\ref{eq.intro.perms.signs.conj}))}}}\overline
{u\underbrace{wu^{-1}u}_{=w}T}\\
&  \ \ \ \ \ \ \ \ \ \ \ \ \ \ \ \ \ \ \ \ \left(
\begin{array}
[c]{c}%
\text{here, we have substituted }uwu^{-1}\text{ for }w\text{ in the sum,}\\
\text{since the map }\mathcal{C}\left(  T\right)  \rightarrow\mathcal{C}%
\left(  uT\right)  ,\ w\mapsto uwu^{-1}\\
\text{is a bijection}%
\end{array}
\right) \\
&  =\sum_{w\in\mathcal{C}\left(  T\right)  }\left(  -1\right)  ^{w}%
\underbrace{\overline{uwT}}_{\substack{=u\overline{wT}\\\text{(by
(\ref{eq.def.tabloid.Sn-act.rewr}))}}}=\sum_{w\in\mathcal{C}\left(  T\right)
}\left(  -1\right)  ^{w}u\overline{wT}=u\underbrace{\sum_{w\in\mathcal{C}%
\left(  T\right)  }\left(  -1\right)  ^{w}\overline{wT}}%
_{\substack{=\mathbf{e}_{T}\\\text{(by (\ref{pf.lem.spechtmod.submod.a.1}))}%
}}=u\mathbf{e}_{T}.
\end{align*}
This proves Lemma \ref{lem.spechtmod.submod} \textbf{(a)}. \medskip

\textbf{(b)} Clearly, $\mathcal{S}^{D}$ is a $\mathbf{k}$-submodule of
$\mathcal{M}^{D}$ (being defined as a span). Next, we shall show that
$\mathcal{S}^{D}$ is an $S_{n}$-subset of $\mathcal{M}^{D}$.

To show this, we must prove that $ux\in\mathcal{S}^{D}$ for each $u\in S_{n}$
and each $x\in\mathcal{S}^{D}$.

So let us prove this. Let $u\in S_{n}$ and $x\in\mathcal{S}^{D}$. We must show
that $ux\in\mathcal{S}^{D}$. Since $\mathcal{S}^{D}$ is spanned by the
polytabloids $\mathbf{e}_{T}$, we can WLOG assume (by linearity) that $x$ is
one of these polytabloids $\mathbf{e}_{T}$. Assume this. Thus, $x=\mathbf{e}%
_{T}$, so that $ux=u\mathbf{e}_{T}=\mathbf{e}_{uT}$ (by Lemma
\ref{lem.spechtmod.submod} \textbf{(a)}). Thus, $ux$ is a polytabloid, and
therefore belongs to $\mathcal{S}^{D}$ (since $\mathcal{S}^{D}$ is the span of
all polytabloids).

Forget that we fixed $u$ and $x$. We thus have shown that $ux\in
\mathcal{S}^{D}$ for each $u\in S_{n}$ and each $x\in\mathcal{S}^{D}$. In
other words, $\mathcal{S}^{D}$ is an $S_{n}$-subset of $\mathcal{M}^{D}$.
Since $\mathcal{S}^{D}$ is also a $\mathbf{k}$-submodule of $\mathcal{M}^{D}$,
we thus conclude that $\mathcal{S}^{D}$ is a subrepresentation of
$\mathcal{M}^{D}$. Hence, by Proposition \ref{prop.rep.G-rep.sub=sub}, we
conclude that $\mathcal{S}^{D}$ is a left $\mathbf{k}\left[  S_{n}\right]
$-submodule of $\mathcal{M}^{D}$. This proves Lemma \ref{lem.spechtmod.submod}
\textbf{(b)}. \medskip

\textbf{(c)} Let $T$ be any $n$-tableau of shape $D$. Let $u\in\mathcal{C}%
\left(  T\right)  $.

As in the proof of part \textbf{(a)}, we can see that $\mathcal{C}\left(
uT\right)  =u\mathcal{C}\left(  T\right)  u^{-1}$. However, $\mathcal{C}%
\left(  T\right)  $ is a group, and $u$ is one of its elements (since
$u\in\mathcal{C}\left(  T\right)  $). Thus, $u\mathcal{C}\left(  T\right)
u^{-1}$ is just $\mathcal{C}\left(  T\right)  $. Hence, we can rewrite the
equality $\mathcal{C}\left(  uT\right)  =u\mathcal{C}\left(  T\right)  u^{-1}$
as $\mathcal{C}\left(  uT\right)  =\mathcal{C}\left(  T\right)  $.

Moreover, $u$ is an element of the group $\mathcal{C}\left(  T\right)  $.
Thus, so is its inverse $u^{-1}$. Hence, the map $\mathcal{C}\left(  T\right)
\rightarrow\mathcal{C}\left(  T\right)  ,\ w\mapsto wu^{-1}$ is a bijection
(with inverse map given by $w\mapsto wu$).

Now, Definition \ref{def.spechtmod.spechtmod} \textbf{(a)} yields%
\begin{align*}
\mathbf{e}_{uT}  &  =\sum_{w\in\mathcal{C}\left(  uT\right)  }\left(
-1\right)  ^{w}\overline{wuT}=\sum_{w\in\mathcal{C}\left(  T\right)  }\left(
-1\right)  ^{w}\overline{wuT}\ \ \ \ \ \ \ \ \ \ \left(  \text{since
}\mathcal{C}\left(  uT\right)  =\mathcal{C}\left(  T\right)  \right) \\
&  =\sum_{w\in\mathcal{C}\left(  T\right)  }\underbrace{\left(  -1\right)
^{wu^{-1}}}_{\substack{=\left(  -1\right)  ^{w}\left(  -1\right)  ^{u^{-1}%
}\\\text{(by the multiplicativity}\\\text{of the sign)}}}\overline
{\underbrace{wu^{-1}u}_{=w}T}\\
&  \ \ \ \ \ \ \ \ \ \ \ \ \ \ \ \ \ \ \ \ \left(
\begin{array}
[c]{c}%
\text{here, we have substituted }wu^{-1}\text{ for }w\text{ in the sum,}\\
\text{since the map }\mathcal{C}\left(  T\right)  \rightarrow\mathcal{C}%
\left(  T\right)  ,\ w\mapsto wu^{-1}\\
\text{is a bijection}%
\end{array}
\right) \\
&  =\sum_{w\in\mathcal{C}\left(  T\right)  }\left(  -1\right)  ^{w}\left(
-1\right)  ^{u^{-1}}\overline{wT}=\left(  -1\right)  ^{u^{-1}}\underbrace{\sum
_{w\in\mathcal{C}\left(  T\right)  }\left(  -1\right)  ^{w}\overline{wT}%
}_{\substack{=\mathbf{e}_{T}\\\text{(by (\ref{pf.lem.spechtmod.submod.a.1}))}%
}}\\
&  =\underbrace{\left(  -1\right)  ^{u^{-1}}}_{\substack{=\left(  -1\right)
^{u}\\\text{(by the basic properties}\\\text{of signs of permutations)}%
}}\mathbf{e}_{T}=\left(  -1\right)  ^{u}\mathbf{e}_{T}.
\end{align*}
This proves Lemma \ref{lem.spechtmod.submod} \textbf{(c)}. \medskip

\textbf{(d)} Let $T$ and $S$ be two column-equivalent $n$-tableaux of shape
$D$. Then, Proposition \ref{prop.tableau.req-R} \textbf{(b)} shows that there
exists some $w\in\mathcal{C}\left(  T\right)  $ such that $S=w\rightharpoonup
T$. Consider this $w$. Then, $S=w\rightharpoonup T=wT$. But Lemma
\ref{lem.spechtmod.submod} \textbf{(c)} (applied to $u=w$) yields
$\mathbf{e}_{wT}=\left(  -1\right)  ^{w}\mathbf{e}_{T}$. In other words,
$\mathbf{e}_{S}=\left(  -1\right)  ^{w}\mathbf{e}_{T}$ (since $S=wT$). Hence,
$\mathbf{e}_{S}=\pm\mathbf{e}_{T}$, so that $\mathbf{e}_{T}=\pm\mathbf{e}_{S}%
$. This proves Lemma \ref{lem.spechtmod.submod} \textbf{(d)}.
\end{proof}

\subsubsection{More examples}

Our above examples show that the Specht modules $\mathcal{S}^{Y\left(
\lambda\right)  }$ for various partitions $\lambda$ include

\begin{itemize}
\item the trivial representation of $S_{n}$ (obtained for $\lambda=\left(
n\right)  $, as seen in Example \ref{exa.spechtmod.Yn}),

\item the sign representation of $S_{n}$ (obtained for $\lambda=\left(
1^{n}\right)  $, as seen in Example \ref{exa.spechtmod.Y111}), and

\item the zero-sum representation $R\left(  \mathbf{k}^{n}\right)  $ (obtained
for $\lambda=\left(  n-1,1\right)  $, as seen in Example
\ref{exa.spechtmod.Yn-1-1}).
\end{itemize}

For $n=3$, all the irreducible representations of $S_{n}$ (over a field of
characteristic $0$) are contained in this list. For $n=4$, there are others.
Can we also identify them as Specht modules?

\begin{example}
\label{exa.spechtmod.Y22}Let $n=4$ and $D=Y\left(  2,2\right)  $. This diagram
$D$ has the form $\ydiagram{2,2}\ \ $.

There are six $n$-tabloids of shape $D$, corresponding bijectively to the
$2$-element subsets of $\left[  4\right]  $. Indeed, any $n$-tabloid of shape
$D$ is uniquely determined by the set of the two entries in its second row,
and will be denoted $\underline{\overline{i\ j}}$, where $i$ and $j$ are these
two entries. For instance,%
\[
\underline{\overline{1\ 4}}=\ytableausetup{tabloids} \ytableaushort{23,14}
\ytableausetup{notabloids}\ \ =\overline{23\backslash\backslash14}.
\]
(This can also be written as $\underline{\overline{4\ 1}}$, since the entries
of a row in an $n$-tabloid do not come ordered. But this is not the same
$n$-tabloid as $\underline{\overline{2\ 3}}$.)

Thus, the Young module $\mathcal{M}^{D}$ is a free $\mathbf{k}$-module of rank
$6$, with basis $\left(  \underline{\overline{1\ 2}},\ \underline{\overline
{1\ 3}},\ \underline{\overline{1\ 4}},\ \underline{\overline{2\ 3}%
},\ \underline{\overline{2\ 4}},\ \underline{\overline{3\ 4}}\right)  $. What
is the Specht module $\mathcal{S}^{D}$ ? For any $n$-tableau $ij\backslash
\backslash k\ell=\ytableaushort{ij,k\ell}\ $, the corresponding polytabloid
$\mathbf{e}_{ij\backslash\backslash k\ell}$ is%
\begin{align*}
\mathbf{e}_{ij\backslash\backslash k\ell}  &  =\sum_{w\in\mathcal{C}\left(
ij\backslash\backslash k\ell\right)  }\left(  -1\right)  ^{w}\overline
{w\left(  ij\backslash\backslash k\ell\right)  }\\
&  =\left(  -1\right)  ^{\operatorname*{id}}\overline{\operatorname*{id}%
\left(  ij\backslash\backslash k\ell\right)  }+\left(  -1\right)  ^{t_{i,k}%
}\overline{t_{i,k}\left(  ij\backslash\backslash k\ell\right)  }\\
&  \ \ \ \ \ \ \ \ \ \ +\left(  -1\right)  ^{t_{j,\ell}}\overline{t_{j,\ell
}\left(  ij\backslash\backslash k\ell\right)  }+\left(  -1\right)
^{t_{i,k}t_{j,\ell}}\overline{t_{i,k}t_{j,\ell}\left(  ij\backslash\backslash
k\ell\right)  }\\
&  \ \ \ \ \ \ \ \ \ \ \ \ \ \ \ \ \ \ \ \ \left(  \text{since }%
\mathcal{C}\left(  ij\backslash\backslash k\ell\right)  =\left\{
\operatorname*{id},\ t_{i,k},\ t_{j,\ell},\ t_{i,k}t_{j,\ell}\right\}  \right)
\\
&  =\overline{ij\backslash\backslash k\ell}-\overline{kj\backslash\backslash
i\ell}-\overline{i\ell\backslash\backslash kj}+\overline{k\ell\backslash
\backslash ij}\\
&  =\underline{\overline{k\ \ell}}-\underline{\overline{i\ \ell}%
}-\underline{\overline{k\ j}}+\underline{\overline{i\ j}}.
\end{align*}
Thus, the Specht module $\mathcal{S}^{D}$ is the span of all these elements
$\underline{\overline{k\ \ell}}-\underline{\overline{i\ \ell}}%
-\underline{\overline{k\ j}}+\underline{\overline{i\ j}}$, where $\left(
i,j,k,\ell\right)  $ ranges over all $4$-tuples of distinct elements of
$\left[  4\right]  $. Many of these elements $\mathbf{e}_{ij\backslash
\backslash k\ell}$ are equal (up to sign), but even the distinct ones are
still linearly dependent. Their span $\mathcal{S}^{D}$ turns out to be only
$2$-dimensional, i.e., a free $\mathbf{k}$-module of rank $2$. One basis of it
is%
\[
\left(  \mathbf{e}_{12\backslash\backslash34},\ \mathbf{e}_{13\backslash
\backslash24}\right)  .
\]
Do the two tableaux in the subscripts remind you of something? For the answer,
look at Example \ref{exa.tableau.tableau1}.
\end{example}

\begin{exercise}
\fbox{2} Assume that $\mathbf{k}$ is a field. Show that the Specht module
$\mathcal{S}^{D}$ in Example \ref{exa.spechtmod.Y22} is an irreducible
representation of $S_{4}$ if and only if $\operatorname*{char}\mathbf{k}\neq
3$. \medskip

[\textbf{Hint:} Use Exercise \ref{exe.youngmod.1dim-sub}.]
\end{exercise}

We have not seen this representation before.

We could analyze the Specht module $\mathcal{S}^{Y\left(  2,1,1\right)  }$
similarly, although this would be more laborious (the Young module
$\mathcal{M}^{Y\left(  2,1,1\right)  }$ is $12$-dimensional). It is yet
another irreducible representation of $S_{4}$ in characteristic $0$. Do we
have all of them now? We will see -- but not very soon (Theorem
\ref{thm.spechtmod.irred}). We have quite a way to go until we understand
Specht modules well enough to answer this. \medskip

The partition $\left(  2,2\right)  $ in Example \ref{exa.spechtmod.Y22} is a
particular case of the partition $\left(  n-2,2\right)  $, which exists for
each $n\geq4$. The Specht module corresponding to the latter partition is the
subject of the following exercise:

\begin{exercise}
\label{exe.spechtmod.Yn-22}Let $n\geq4$, and let $\lambda$ be the partition
$\left(  n-2,2\right)  $. Let $D=Y\left(  \lambda\right)  $. As in Example
\ref{exa.spechtmod.Y22}, let us write $\underline{\overline{i\ j}}$ for the
$n$-tabloid of shape $D$ whose second row has entries $i$ and $j$. \medskip

\textbf{(a)} \fbox{1} Prove that%
\[
\mathcal{S}^{D}=\operatorname*{span}\nolimits_{\mathbf{k}}\left\{
\underline{\overline{k\ \ell}}-\underline{\overline{i\ \ell}}%
-\underline{\overline{k\ j}}+\underline{\overline{i\ j}}\ \mid\ i,j,k,\ell
\in\left[  n\right]  \text{ distinct}\right\}  .
\]

\textbf{(b)} \fbox{1} Let $E=Y\left(  n-1,1\right)  $. As in Example
\ref{exa.youngmod.n-1-1}, let us write $\underline{\overline{i}}$ for the
$n$-tabloid of shape $E$ whose second row has entry $i$. Consider the
$\mathbf{k}$-linear map
\begin{align*}
\psi_{1}:\mathcal{M}^{D}  &  \rightarrow\mathcal{M}^{E},\\
\underline{\overline{i\ j}}  &  \mapsto\underline{\overline{i}}%
+\underline{\overline{j}}%
\end{align*}
and the $\mathbf{k}$-linear map
\begin{align*}
\psi_{0}:\mathcal{M}^{D}  &  \rightarrow\mathbf{k},\\
\underline{\overline{i\ j}}  &  \mapsto1.
\end{align*}
Prove that both these maps $\psi_{1}$ and $\psi_{0}$ are left $\mathbf{k}%
\left[  S_{n}\right]  $-linear (where $\mathbf{k}$ is equipped with the
trivial $S_{n}$-action). \medskip

\textbf{(c)} \fbox{4} Prove that $\mathcal{S}^{D}=\left(  \operatorname*{Ker}%
\psi_{1}\right)  \cap\left(  \operatorname*{Ker}\psi_{0}\right)  $. \medskip

\textbf{(d)} \fbox{1} Prove that $\mathcal{S}^{D}=\operatorname*{Ker}\psi_{1}$
if $2$ is invertible in $\mathbf{k}$. \medskip

\textbf{(e)} \fbox{1} Assume that both $n-1$ and $\dfrac{n\left(  n-1\right)
}{2}$ equal $0$ in $\mathbf{k}$. Prove that $\mathcal{S}^{D}$ contains the
element $\sum_{1\leq i<j\leq n}\underline{\overline{i\ j}}$ (this is the sum
of all $n$-tabloids of shape $D$), and thus contains a trivial
subrepresentation spanned by this element. \medskip

[\textbf{Hint:} For part \textbf{(c)}, first show that each element of
$\mathcal{M}^{D}$ can be written as some element of $\mathcal{S}^{D}$ plus
some element of the form $\sum_{i=3}^{n}\alpha_{i}\underline{\overline{1\ i}%
}+\sum_{i=3}^{n}\beta_{i}\underline{\overline{2\ i}}+\gamma
\underline{\overline{1\ 2}}$ with $\alpha_{i},\beta_{i},\gamma\in\mathbf{k}$;
then analyze the conditions under which an element of the latter form belongs
to $\operatorname*{Ker}\psi_{1}$ and to $\operatorname*{Ker}\psi_{0}$.]
\end{exercise}

\begin{exercise}
\fbox{2} Let $D$ be any diagram. Prove that $\mathcal{S}^{D}=\mathcal{M}^{D}$
holds if and only if the ring $\mathbf{k}$ is trivial or the diagram $D$ has
no two cells in the same column. \medskip

[\textbf{Hint:} In all other cases, find a nontrivial $\mathbf{k}$-linear map
$f:\mathcal{M}^{D}\rightarrow\mathbf{k}$ that sends every polytabloid
$\mathbf{e}_{T}$ to $0$.]
\end{exercise}

\subsection{Symmetrizers, antisymmetrizers and the left ideal avatar}

We will now take a look at the Specht modules from a different direction. This
relies on a few preliminary results about $\mathbf{k}\left[  S_{n}\right]  $.

\subsubsection{Row symmetrizers and column antisymmetrizers}

We begin by introducing some more elements of the group algebra $\mathbf{k}%
\left[  S_{n}\right]  $:

\begin{definition}
\label{def.symmetrizers.symmetrizers}Let $T$ be an $n$-tableau of any shape
$D$. We define its \emph{row symmetrizer}%
\[
\nabla_{\operatorname*{Row}T}:=\sum_{w\in\mathcal{R}\left(  T\right)  }%
w\in\mathbf{k}\left[  S_{n}\right]
\]
and its \emph{column antisymmetrizer}%
\[
\nabla_{\operatorname*{Col}T}^{-}:=\sum_{w\in\mathcal{C}\left(  T\right)
}\left(  -1\right)  ^{w}w\in\mathbf{k}\left[  S_{n}\right]  .
\]
(We could similarly define the column symmetrizer $\nabla_{\operatorname*{Col}%
T}$ and the row antisymmetrizer $\nabla_{\operatorname*{Row}T}^{-}$, but we
probably won't need them.)
\end{definition}

\begin{example}
For the $5$-tableau $12\backslash\backslash34\backslash\backslash
5=\ytableaushort{12,34,5}$ of shape $Y\left(  2,2,1\right)  $, we have%
\begin{align*}
\nabla_{\operatorname*{Row}\left(  12\backslash\backslash34\backslash
\backslash5\right)  }  &  =1+t_{1,2}+t_{3,4}+t_{1,2}t_{3,4}%
\ \ \ \ \ \ \ \ \ \ \text{and}\\
\nabla_{\operatorname*{Col}\left(  12\backslash\backslash34\backslash
\backslash5\right)  }^{-}  &  =1-t_{1,3}-t_{1,5}-t_{3,5}+\operatorname*{cyc}%
\nolimits_{1,3,5}+\operatorname*{cyc}\nolimits_{1,5,3}-\,t_{2,4}%
+t_{1,3}t_{2,4}\\
&  \ \ \ \ \ \ \ \ \ \ +t_{1,5}t_{2,4}+t_{3,5}t_{2,4}-\operatorname*{cyc}%
\nolimits_{1,3,5}t_{2,4}-\operatorname*{cyc}\nolimits_{1,5,3}t_{2,4}.
\end{align*}

\end{example}

The letter $\nabla$ was chosen because the definitions of $\nabla
_{\operatorname*{Row}T}$ and $\nabla_{\operatorname*{Col}T}^{-}$ are rather
similar to those of the $X$-integrals $\nabla_{X}$ and the $X$-sign-integrals
$\nabla_{X}^{-}$ (see Definition \ref{def.intX.intX}). There is in fact a
direct relation between the former and the latter:

\begin{proposition}
\label{prop.symmetrizers.int}Let $T$ be an $n$-tableau of any shape $D$. Then:
\medskip

\textbf{(a)} For any $i\in\mathbb{Z}$, let $\operatorname*{Row}\left(
i,T\right)  $ denote the set of all entries in the $i$-th row of $T$. If
$\ell\in\mathbb{N}$ is chosen such that all cells of $D$ lie in rows
$1,2,\ldots,\ell$, then%
\[
\nabla_{\operatorname*{Row}T}=\prod_{i=1}^{\ell}\nabla_{\operatorname*{Row}%
\left(  i,T\right)  }.
\]
Here we are not specifying the order of the factors in the product, since it
does not matter (the factors commute). \medskip

\textbf{(b)} For any $j\in\mathbb{Z}$, let $\operatorname*{Col}\left(
j,T\right)  $ denote the set of all entries in the $j$-th column of $T$. If
$k\in\mathbb{N}$ is chosen such that all cells of $D$ lie in columns
$1,2,\ldots,k$, then%
\[
\nabla_{\operatorname*{Col}T}^{-}=\prod_{j=1}^{k}\nabla_{\operatorname*{Col}%
\left(  j,T\right)  }^{-}.
\]
Here, again, we are not specifying the order of the factors, since it does not matter.
\end{proposition}

\begin{example}
Let $n=8$ and $T=123\backslash\backslash45\backslash\backslash67\backslash
\backslash8=\ytableaushort{123,45,67,8}$\ . Then, Proposition
\ref{prop.symmetrizers.int} \textbf{(a)} (for $\ell=4$) yields%
\[
\nabla_{\operatorname*{Row}\left(  123\backslash\backslash45\backslash
\backslash67\backslash\backslash8\right)  }=\nabla_{\left\{  1,2,3\right\}
}\nabla_{\left\{  4,5\right\}  }\nabla_{\left\{  6,7\right\}  }%
\underbrace{\nabla_{\left\{  8\right\}  }}_{=1}=\nabla_{\left\{
1,2,3\right\}  }\nabla_{\left\{  4,5\right\}  }\nabla_{\left\{  6,7\right\}
}.
\]
(We could also take $\ell$ larger than $4$, but this would merely introduce
extra factors $\nabla_{\operatorname*{Row}\left(  i,T\right)  }$ corresponding
to empty rows, and all these factors are $\nabla_{\varnothing}=1$.)

Furthermore, Proposition \ref{prop.symmetrizers.int} \textbf{(b)} (for $k=3$)
yields%
\[
\nabla_{\operatorname*{Col}\left(  123\backslash\backslash45\backslash
\backslash67\backslash\backslash8\right)  }^{-}=\nabla_{\left\{
1,4,6,8\right\}  }^{-}\nabla_{\left\{  2,5,7\right\}  }^{-}\underbrace{\nabla
_{\left\{  3\right\}  }^{-}}_{=1}.
\]
(Here, again, we could take $k$ larger than $3$, but this would merely
introduce factors of the form $\nabla_{\varnothing}^{-}=1$ corresponding to
empty columns.)
\end{example}

\begin{proof}
[Proof of Proposition \ref{prop.symmetrizers.int}.]\textbf{(b)} We know that
$T$ is an $n$-tableau. Hence, the entries of $T$ are the numbers
$1,2,\ldots,n$, and each of these numbers appears exactly once (i.e., the
tableau $T$ has no two equal entries).

Let $k\in\mathbb{N}$ be such that all cells of $D$ lie in columns
$1,2,\ldots,k$. Then, each entry of $T$ lies in exactly one of the columns
$1,2,\ldots,k$ (indeed, not more than one, since $T$ has no two equal
entries). In other words, each entry of $T$ lies in exactly one of the sets
$\operatorname*{Col}\left(  1,T\right)  ,\ \operatorname*{Col}\left(
2,T\right)  ,\ \ldots,\ \operatorname*{Col}\left(  k,T\right)  $.

In other words, each of the numbers $1,2,\ldots,n$ lies in exactly one of the
sets $\operatorname*{Col}\left(  1,T\right)  ,\ \operatorname*{Col}\left(
2,T\right)  ,\ \ldots,\ \operatorname*{Col}\left(  k,T\right)  $ (since the
entries of $T$ are the numbers $1,2,\ldots,n$). Thus, $\operatorname*{Col}%
\left(  1,T\right)  ,\ \operatorname*{Col}\left(  2,T\right)  ,\ \ldots
,\ \operatorname*{Col}\left(  k,T\right)  $ are $k$ disjoint subsets of
$\left[  n\right]  $ and satisfy $\operatorname*{Col}\left(  1,T\right)
\cup\operatorname*{Col}\left(  2,T\right)  \cup\cdots\cup\operatorname*{Col}%
\left(  k,T\right)  =\left[  n\right]  $. Hence, the equality
(\ref{eq.prop.int.wXi=Xi.b.2}) from Proposition \ref{prop.int.wXi=Xi}
\textbf{(b)} (applied to $X_{j}=\operatorname*{Col}\left(  j,T\right)  $)
yields
\[
\sum_{\substack{w\in S_{n};\\w\left(  \operatorname*{Col}\left(  i,T\right)
\right)  \subseteq\operatorname*{Col}\left(  i,T\right)  \text{ for all }%
i\in\left[  k\right]  }}\left(  -1\right)  ^{w}w=\nabla_{\operatorname*{Col}%
\left(  1,T\right)  }^{-}\nabla_{\operatorname*{Col}\left(  2,T\right)  }%
^{-}\cdots\nabla_{\operatorname*{Col}\left(  k,T\right)  }^{-}.
\]
Moreover, Proposition \ref{prop.int.wXi=Xi} \textbf{(c)} (applied to
$X_{j}=\operatorname*{Col}\left(  j,T\right)  $) yields that the factors
$\nabla_{\operatorname*{Col}\left(  1,T\right)  }^{-},\nabla
_{\operatorname*{Col}\left(  2,T\right)  }^{-},\ldots,\nabla
_{\operatorname*{Col}\left(  k,T\right)  }^{-}$ commute. Thus, we can write
their product as $\prod_{j=1}^{k}\nabla_{\operatorname*{Col}\left(
j,T\right)  }^{-}$ without specifying the order of multiplication.

However, for any permutation $w\in S_{n}$, we have the logical equivalence%
\[
\left(  w\in\mathcal{C}\left(  T\right)  \right)  \ \Longleftrightarrow
\ \left(  w\left(  \operatorname*{Col}\left(  i,T\right)  \right)
\subseteq\operatorname*{Col}\left(  i,T\right)  \text{ for all }i\in\left[
k\right]  \right)
\]
\footnote{\textit{Proof.} Let $w\in S_{n}$. Then, we have the following chain
of logical equivalences:%
\begin{align*}
&  \ \left(  w\in\mathcal{C}\left(  T\right)  \right) \\
&  \Longleftrightarrow\ \left(  w\text{ is a vertical permutation for
}T\right)  \ \ \ \ \ \ \ \ \ \ \left(  \text{by the definition of }%
\mathcal{C}\left(  T\right)  \right) \\
&  \Longleftrightarrow\ \left(  \text{for each }i\in\left[  n\right]  \text{,
the number }w\left(  i\right)  \text{ lies in the same column of }T\text{ as
}i\right) \\
&  \ \ \ \ \ \ \ \ \ \ \ \ \ \ \ \ \ \ \ \ \left(  \text{by the definition of
vertical permutations}\right) \\
&  \Longleftrightarrow\ \left(  \text{for each }i\in\left[  n\right]  \text{,
for each }c\in\mathbb{Z}\text{, if }i\text{ lies in the }c\text{-th column of
}T\text{, then so does }w\left(  i\right)  \right) \\
&  \ \ \ \ \ \ \ \ \ \ \ \ \ \ \ \ \ \ \ \ \left(  \text{here, we just
restated the \textquotedblleft lies in the same column\textquotedblright\ part
in words}\right) \\
&  \Longleftrightarrow\ \left(  \text{for each }i\in\left[  n\right]  \text{,
for each }c\in\left[  k\right]  \text{, if }i\text{ lies in the }c\text{-th
column of }T\text{, then so does }w\left(  i\right)  \right) \\
&  \ \ \ \ \ \ \ \ \ \ \ \ \ \ \ \ \ \ \ \ \left(
\begin{array}
[c]{c}%
\text{here, we have restricted the index }c\text{ to the set }\left[
k\right]  \text{,}\\
\text{since the only columns of }T\text{ that can possibly}\\
\text{contain }i\text{ are the columns }1,2,\ldots,k\\
\text{(because all cells of }D\text{ lie in columns }1,2,\ldots,k\text{)}%
\end{array}
\right) \\
&  \Longleftrightarrow\ \left(  \text{for each }c\in\left[  k\right]  \text{,
for each }i\in\left[  n\right]  \text{, if }i\text{ lies in the }c\text{-th
column of }T\text{, then so does }w\left(  i\right)  \right) \\
&  \ \ \ \ \ \ \ \ \ \ \ \ \ \ \ \ \ \ \ \ \left(  \text{here, we swapped the
two \textquotedblleft for all\textquotedblright-quantifiers}\right) \\
&  \Longleftrightarrow\ \left(  \text{for each }c\in\left[  k\right]  \text{,
for each }i\in\left[  n\right]  \text{, if }i\in\operatorname*{Col}\left(
c,T\right)  \text{, then }w\left(  i\right)  \in\operatorname*{Col}\left(
c,T\right)  \right) \\
&  \ \ \ \ \ \ \ \ \ \ \ \ \ \ \ \ \ \ \ \ \left(
\begin{array}
[c]{c}%
\text{since \textquotedblleft}i\text{ lies in the }c\text{-th column of
}T\text{\textquotedblright}\\
\text{is equivalent to \textquotedblleft}i\in\operatorname*{Col}\left(
c,T\right)  \text{\textquotedblright,}\\
\text{and \textquotedblleft}w\left(  i\right)  \text{ lies in the }c\text{-th
column of }T\text{\textquotedblright}\\
\text{is equivalent to \textquotedblleft}w\left(  i\right)  \in
\operatorname*{Col}\left(  c,T\right)  \text{\textquotedblright}%
\end{array}
\right) \\
&  \Longleftrightarrow\ \left(  \text{for each }c\in\left[  k\right]  \text{,
for each }i\in\operatorname*{Col}\left(  c,T\right)  \text{, we have }w\left(
i\right)  \in\operatorname*{Col}\left(  c,T\right)  \right) \\
&  \ \ \ \ \ \ \ \ \ \ \ \ \ \ \ \ \ \ \ \ \left(
\begin{array}
[c]{c}%
\text{since the elements }i\in\left[  n\right]  \text{ that satisfy }%
i\in\operatorname*{Col}\left(  c,T\right)  \text{ are}\\
\text{precisely the elements of }\operatorname*{Col}\left(  c,T\right)  \text{
(because }\operatorname*{Col}\left(  c,T\right)  \subseteq\left[  n\right]
\text{)}%
\end{array}
\right) \\
&  \Longleftrightarrow\ \left(  \text{for each }c\in\left[  k\right]  \text{,
we have }w\left(  \operatorname*{Col}\left(  c,T\right)  \right)
\subseteq\operatorname*{Col}\left(  c,T\right)  \right) \\
&  \ \ \ \ \ \ \ \ \ \ \ \ \ \ \ \ \ \ \ \ \left(
\begin{array}
[c]{c}%
\text{since the statement \textquotedblleft for each }i\in\operatorname*{Col}%
\left(  c,T\right)  \text{, we have }w\left(  i\right)  \in\operatorname*{Col}%
\left(  c,T\right)  \text{\textquotedblright}\\
\text{is equivalent to \textquotedblleft}w\left(  \operatorname*{Col}\left(
c,T\right)  \right)  \subseteq\operatorname*{Col}\left(  c,T\right)
\text{\textquotedblright}%
\end{array}
\right) \\
&  \Longleftrightarrow\ \left(  w\left(  \operatorname*{Col}\left(
c,T\right)  \right)  \subseteq\operatorname*{Col}\left(  c,T\right)  \text{
for all }c\in\left[  k\right]  \right) \\
&  \Longleftrightarrow\ \left(  w\left(  \operatorname*{Col}\left(
i,T\right)  \right)  \subseteq\operatorname*{Col}\left(  i,T\right)  \text{
for all }i\in\left[  k\right]  \right)  \ \ \ \ \ \ \ \ \ \ \left(
\text{here, we renamed the index }c\text{ as }i\right)  .
\end{align*}
}. Hence,%
\[
\sum_{w\in\mathcal{C}\left(  T\right)  }\left(  -1\right)  ^{w}w=\sum
_{\substack{w\in S_{n};\\w\left(  \operatorname*{Col}\left(  i,T\right)
\right)  \subseteq\operatorname*{Col}\left(  i,T\right)  \text{ for all }%
i\in\left[  k\right]  }}\left(  -1\right)  ^{w}w.
\]

But the definition of $\nabla_{\operatorname*{Col}T}^{-}$ yields%
\begin{align*}
\nabla_{\operatorname*{Col}T}^{-}  &  =\sum_{w\in\mathcal{C}\left(  T\right)
}\left(  -1\right)  ^{w}w=\sum_{\substack{w\in S_{n};\\w\left(
\operatorname*{Col}\left(  i,T\right)  \right)  \subseteq\operatorname*{Col}%
\left(  i,T\right)  \text{ for all }i\in\left[  k\right]  }}\left(  -1\right)
^{w}w\\
&  =\nabla_{\operatorname*{Col}\left(  1,T\right)  }^{-}\nabla
_{\operatorname*{Col}\left(  2,T\right)  }^{-}\cdots\nabla
_{\operatorname*{Col}\left(  k,T\right)  }^{-}=\prod_{j=1}^{k}\nabla
_{\operatorname*{Col}\left(  j,T\right)  }^{-}.
\end{align*}
This proves Proposition \ref{prop.symmetrizers.int} \textbf{(b)}. \medskip

\textbf{(a)} Analogous to part \textbf{(b)}, but without the signs.
\end{proof}

Some other properties of row symmetrizers and column antisymmetrizers also
resemble analogous properties of integrals and sign-integrals:

\begin{proposition}
\label{prop.symmetrizers.fix}Let $T$ be an $n$-tableau of any shape $D$. Then:
\medskip

\textbf{(a)} We have $w\nabla_{\operatorname*{Row}T}=\nabla
_{\operatorname*{Row}T}w=\nabla_{\operatorname*{Row}T}$ for each
$w\in\mathcal{R}\left(  T\right)  $. \medskip

\textbf{(b)} We have $w\nabla_{\operatorname*{Col}T}^{-}=\nabla
_{\operatorname*{Col}T}^{-}w=\left(  -1\right)  ^{w}\nabla
_{\operatorname*{Col}T}^{-}$ for each $w\in\mathcal{C}\left(  T\right)  $.
\end{proposition}

\begin{proof}
Analogous to Proposition \ref{prop.integral.fix} (but using the group
$\mathcal{R}\left(  T\right)  $ or $\mathcal{C}\left(  T\right)  $ instead of
$S_{n}$).
\end{proof}

\begin{corollary}
\label{cor.symmetrizers.square}Let $T$ be an $n$-tableau of any shape $D$.
Then: \medskip

\textbf{(a)} We have $\left(  \nabla_{\operatorname*{Row}T}\right)
^{2}=\left\vert \mathcal{R}\left(  T\right)  \right\vert \cdot\nabla
_{\operatorname*{Row}T}$. \medskip

\textbf{(b)} We have $\left(  \nabla_{\operatorname*{Col}T}^{-}\right)
^{2}=\left\vert \mathcal{C}\left(  T\right)  \right\vert \cdot\nabla
_{\operatorname*{Col}T}^{-}$.
\end{corollary}

\begin{proof}
Each of the two parts is proved similarly to Corollary
\ref{cor.integral.square}, but using Proposition \ref{prop.symmetrizers.fix}
instead of Proposition \ref{prop.integral.fix}. (Part \textbf{(b)} involves
signs as well, but they do not complicate the computation.)
\end{proof}

We note that the factors $\left\vert \mathcal{R}\left(  T\right)  \right\vert
$ and $\left\vert \mathcal{C}\left(  T\right)  \right\vert $ in Corollary
\ref{cor.symmetrizers.square} can be computed rather explicitly:

\begin{proposition}
\label{prop.tableau.RTsize}Let $T$ be an $n$-tableau of any shape $D$. Then:
\medskip

\textbf{(a)} For any $i\in\mathbb{Z}$, let $\operatorname*{Row}\left(
i,T\right)  $ denote the set of all entries in the $i$-th row of $T$. If
$\ell\in\mathbb{N}$ is chosen such that all cells of $D$ lie in rows
$1,2,\ldots,\ell$, then%
\[
\left\vert \mathcal{R}\left(  T\right)  \right\vert =\prod_{i=1}^{\ell
}\left\vert \operatorname*{Row}\left(  i,T\right)  \right\vert !.
\]

\textbf{(b)} For any $j\in\mathbb{Z}$, let $\operatorname*{Col}\left(
j,T\right)  $ denote the set of all entries in the $j$-th column of $T$. If
$k\in\mathbb{N}$ is chosen such that all cells of $D$ lie in columns
$1,2,\ldots,k$, then%
\[
\left\vert \mathcal{C}\left(  T\right)  \right\vert =\prod_{j=1}^{k}\left\vert
\operatorname*{Col}\left(  j,T\right)  \right\vert !.
\]

\end{proposition}

\begin{fineprint}
\begin{proof}
\textbf{(a)} This is intuitively obvious, but the slickest proof uses a bit of
algebra: Let $\mathbf{k}$ be the ring $\mathbb{Z}$. Recall the counit
$\varepsilon:\mathbf{k}\left[  S_{n}\right]  \rightarrow\mathbf{k}$ (which was
defined in Definition \ref{def.eps.eps}). This counit $\varepsilon$ is a
$\mathbf{k}$-algebra morphism (by Theorem \ref{thm.eps.mor}). Now, let
$\ell\in\mathbb{N}$ be chosen such that all cells of $D$ lie in rows
$1,2,\ldots,\ell$. Then, Proposition \ref{prop.symmetrizers.int} \textbf{(a)}
yields%
\[
\nabla_{\operatorname*{Row}T}=\prod_{i=1}^{\ell}\nabla_{\operatorname*{Row}%
\left(  i,T\right)  }.
\]
Applying the map $\varepsilon$ to both sides of this equality, we find%
\begin{align*}
\varepsilon\left(  \nabla_{\operatorname*{Row}T}\right)   &  =\varepsilon
\left(  \prod_{i=1}^{\ell}\nabla_{\operatorname*{Row}\left(  i,T\right)
}\right) \\
&  =\prod_{i=1}^{\ell}\underbrace{\varepsilon\left(  \nabla
_{\operatorname*{Row}\left(  i,T\right)  }\right)  }_{\substack{=\left\vert
\operatorname*{Row}\left(  i,T\right)  \right\vert !\\\text{(by Proposition
\ref{prop.eps.Xint},}\\\text{applied to }X=\operatorname*{Row}\left(
i,T\right)  \text{)}}}\ \ \ \ \ \ \ \ \ \ \left(  \text{since }\varepsilon
\text{ is a }\mathbf{k}\text{-algebra morphism}\right) \\
&  =\prod_{i=1}^{\ell}\left\vert \operatorname*{Row}\left(  i,T\right)
\right\vert !.
\end{align*}
Comparing this with
\begin{align*}
\varepsilon\left(  \nabla_{\operatorname*{Row}T}\right)   &  =\varepsilon
\left(  \sum_{w\in\mathcal{R}\left(  T\right)  }w\right)
\ \ \ \ \ \ \ \ \ \ \left(  \text{since }\nabla_{\operatorname*{Row}T}%
=\sum_{w\in\mathcal{R}\left(  T\right)  }w\right) \\
&  =\sum_{w\in\mathcal{R}\left(  T\right)  }\underbrace{\varepsilon\left(
w\right)  }_{\substack{=1\\\text{(by the definition of }\varepsilon
\text{,}\\\text{since }w\in\mathcal{R}\left(  T\right)  \subseteq
S_{n}\text{)}}}\ \ \ \ \ \ \ \ \ \ \left(  \text{since the map }%
\varepsilon\text{ is }\mathbf{k}\text{-linear}\right) \\
&  =\sum_{w\in\mathcal{R}\left(  T\right)  }1=\left\vert \mathcal{R}\left(
T\right)  \right\vert \cdot1=\left\vert \mathcal{R}\left(  T\right)
\right\vert ,
\end{align*}
we obtain $\left\vert \mathcal{R}\left(  T\right)  \right\vert =\prod
_{i=1}^{\ell}\left\vert \operatorname*{Row}\left(  i,T\right)  \right\vert !$.
This proves Proposition \ref{prop.tableau.RTsize} \textbf{(a)}. \medskip

\textbf{(b)} This is analogous to part \textbf{(a)}. Namely, we have to define
an element $\nabla_{\operatorname*{Col}T}:=\sum_{w\in\mathcal{C}\left(
T\right)  }w$ of $\mathbf{k}\left[  S_{n}\right]  $, and we have to show that
this element satisfies $\nabla_{\operatorname*{Col}T}=\prod_{j=1}^{k}%
\nabla_{\operatorname*{Col}\left(  j,T\right)  }$ (this is an analogue of
Proposition \ref{prop.symmetrizers.int} \textbf{(a)}, and can be proved in the
same way). Then, we have to apply $\varepsilon$ as in our above proof of
Proposition \ref{prop.tableau.RTsize} \textbf{(a)}.
\end{proof}
\end{fineprint}

Next, let us show some ways to factor out simple factors from $\nabla
_{\operatorname*{Row}T}$ and $\nabla_{\operatorname*{Col}T}^{-}$, similarly to
Proposition \ref{prop.intX.basics} \textbf{(d)}:

\begin{proposition}
\label{prop.symmetrizers.factor-out-row}Let $T$ be an $n$-tableau of any shape
$D$. Let $i$ and $j$ be two distinct elements of $\left[  n\right]  $ that lie
in the same row of $T$. Then: \medskip

\textbf{(a)} The transposition $t_{i,j}$ belongs to $\mathcal{R}\left(
T\right)  $. \medskip

\textbf{(b)} We have
\[
\nabla_{\operatorname*{Row}T}=\nabla_{\operatorname*{Row}T}%
^{\operatorname*{even}}\cdot\left(  1+t_{i,j}\right)  =\left(  1+t_{i,j}%
\right)  \cdot\nabla_{\operatorname*{Row}T}^{\operatorname*{even}},
\]
where we set%
\[
\nabla_{\operatorname*{Row}T}^{\operatorname*{even}}:=\sum_{\substack{w\in
\mathcal{R}\left(  T\right)  ;\\\left(  -1\right)  ^{w}=1}}w.
\]

\textbf{(c)} We have $t_{i,j}\nabla_{\operatorname*{Row}T}=\nabla
_{\operatorname*{Row}T}t_{i,j}=\nabla_{\operatorname*{Row}T}$.
\end{proposition}

\begin{proof}
\textbf{(a)} It is easy to see that the transposition $t_{i,j}\in S_{n}$ is a
horizontal permutation for $T$. (\textit{Proof:} We must show that for each
$k\in\left[  n\right]  $, the number $t_{i,j}\left(  k\right)  $ lies in the
same row of $T$ as $k$ does. This is true for $k=i$ (since $t_{i,j}\left(
i\right)  =j$ lies in the same row of $T$ as $i$ does) and also true for $k=j$
(for similar reasons), but it is even more obvious for all remaining values of
$k$ (since $t_{i,j}\left(  k\right)  =k$ for all those values). Thus, this is
true for all $k\in\left[  n\right]  $. It follows that the permutation
$t_{i,j}$ is horizontal for $T$.)

So we have shown that $t_{i,j}$ is a horizontal permutation for $T$. In other
words, $t_{i,j}$ belongs to $\mathcal{R}\left(  T\right)  $ (since
$\mathcal{R}\left(  T\right)  $ is the set of all permutations $w\in S_{n}$
that are horizontal for $T$). This proves Proposition
\ref{prop.symmetrizers.factor-out-row} \textbf{(a)}. \medskip

\textbf{(b)} We know that $\mathcal{R}\left(  T\right)  $ is a group, and
$t_{i,j}$ is one of its elements (by Proposition
\ref{prop.symmetrizers.factor-out-row} \textbf{(a)}). From here on, the proof
of Proposition \ref{prop.symmetrizers.factor-out-row} \textbf{(b)} proceeds
exactly like the proof of Proposition \ref{prop.intX.basics} \textbf{(d)}
(more precisely, the proof of (\ref{prop.intX.basics.d.1})), except that
$S_{n,X}$, $\nabla_{X}$ and $\nabla_{X}^{\operatorname*{even}}$ must be
replaced by $\mathcal{R}\left(  T\right)  $, $\nabla_{\operatorname*{Row}T}$
and $\nabla_{\operatorname*{Row}T}^{\operatorname*{even}}$. Thus, Proposition
\ref{prop.symmetrizers.factor-out-row} \textbf{(b)} is proved. \medskip

\textbf{(c)} Proposition \ref{prop.symmetrizers.factor-out-row} \textbf{(a)}
shows that $t_{i,j}\in\mathcal{R}\left(  T\right)  $. Thus, Proposition
\ref{prop.symmetrizers.fix} \textbf{(a)} (applied to $w=t_{i,j}$) yields
$t_{i,j}\nabla_{\operatorname*{Row}T}=\nabla_{\operatorname*{Row}T}%
t_{i,j}=\nabla_{\operatorname*{Row}T}$. This proves Proposition
\ref{prop.symmetrizers.factor-out-row} \textbf{(c)}.
\end{proof}

\begin{proposition}
\label{prop.symmetrizers.factor-out-col}Let $T$ be an $n$-tableau of any shape
$D$. Let $i$ and $j$ be two distinct elements of $\left[  n\right]  $ that lie
in the same column of $T$. Then: \medskip

\textbf{(a)} The transposition $t_{i,j}$ belongs to $\mathcal{C}\left(
T\right)  $. \medskip

\textbf{(b)} We have%
\[
\nabla_{\operatorname*{Col}T}^{-}=\nabla_{\operatorname*{Col}T}%
^{\operatorname*{even}}\cdot\left(  1-t_{i,j}\right)  =\left(  1-t_{i,j}%
\right)  \cdot\nabla_{\operatorname*{Col}T}^{\operatorname*{even}},
\]
where we set%
\[
\nabla_{\operatorname*{Col}T}^{\operatorname*{even}}:=\sum_{\substack{w\in
\mathcal{C}\left(  T\right)  ;\\\left(  -1\right)  ^{w}=1}}w.
\]

\textbf{(c)} We have $t_{i,j}\nabla^{-}_{\operatorname*{Col}T} =\nabla
^{-}_{\operatorname*{Col}T}t_{i,j}=-\nabla^{-}_{\operatorname*{Col}T}$.
\end{proposition}

\begin{proof}
This is analogous to Proposition \ref{prop.symmetrizers.factor-out-row},
except that the signs need to be taken into account (which is easy, since we
know that $\left(  -1\right)  ^{t_{i,j}}=-1$).
\end{proof}

The antipode $S$ preserves both row symmetrizers and column antisymmetrizers:

\begin{proposition}
\label{prop.symmetrizers.antipode}Let $T$ be an $n$-tableau of any shape $D$.
Then, the antipode $S$ of $\mathbf{k}\left[  S_{n}\right]  $ satisfies%
\[
S\left(  \nabla_{\operatorname*{Row}T}\right)  =\nabla_{\operatorname*{Row}%
T}\ \ \ \ \ \ \ \ \ \ \text{and}\ \ \ \ \ \ \ \ \ \ S\left(  \nabla
_{\operatorname*{Col}T}^{-}\right)  =\nabla_{\operatorname*{Col}T}^{-}.
\]

\end{proposition}

\begin{proof}
The proof of $S\left(  \nabla_{\operatorname*{Col}T}^{-}\right)
=\nabla_{\operatorname*{Col}T}^{-}$ is analogous to the proof of $S\left(
\nabla^{-}\right)  =\nabla^{-}$ in Example \ref{exa.S.exas1} \textbf{(c)}
(indeed, $\mathcal{C}\left(  T\right)  $ is a group, so that the map
$\mathcal{C}\left(  T\right)  \rightarrow\mathcal{C}\left(  T\right)
,\ w\mapsto w^{-1}$ is a bijection). The proof of $S\left(  \nabla
_{\operatorname*{Row}T}\right)  =\nabla_{\operatorname*{Row}T}$ is analogous
but easier.
\end{proof}

Next, let us see how row symmetrizers and column antisymmetrizers are affected
when a permutation $w\in S_{n}$ acts on an $n$-tableau $T$:

\begin{proposition}
\label{prop.symmetrizers.conj}Let $T$ be an $n$-tableau of any shape $D$. Let
$w\in S_{n}$. Then,%
\begin{align*}
\nabla_{\operatorname*{Row}\left(  w\rightharpoonup T\right)  }  &
=w\nabla_{\operatorname*{Row}T}w^{-1}\ \ \ \ \ \ \ \ \ \ \text{and}\\
\nabla_{\operatorname*{Col}\left(  w\rightharpoonup T\right)  }^{-}  &
=w\nabla_{\operatorname*{Col}T}^{-}w^{-1}.
\end{align*}

\end{proposition}

\begin{proof}
Proposition \ref{prop.tableau.Sn-act.0} \textbf{(c)} yields $\mathcal{R}%
\left(  w\rightharpoonup T\right)  =w\mathcal{R}\left(  T\right)  w^{-1}$ and
$\mathcal{C}\left(  w\rightharpoonup T\right)  =w\mathcal{C}\left(  T\right)
w^{-1}$. The latter equality shows that the map
\begin{align*}
\mathcal{C}\left(  T\right)   &  \rightarrow\mathcal{C}\left(
w\rightharpoonup T\right)  ,\\
u  &  \mapsto wuw^{-1}%
\end{align*}
is well-defined and surjective. This map is also injective\footnote{Indeed, if
two elements $u,v\in\mathcal{C}\left(  T\right)  $ satisfy $wuw^{-1}=wvw^{-1}%
$, then $u=v$ (because $S_{n}$ is a group, and thus we can cancel both $w$ and
$w^{-1}$ from the equality $wuw^{-1}=wvw^{-1}$).}, and thus is a bijection.

However, the definition of $\nabla_{\operatorname*{Col}T}^{-}$ says that
$\nabla_{\operatorname*{Col}T}^{-}=\sum_{u\in\mathcal{C}\left(  T\right)
}\left(  -1\right)  ^{u}u$. For the same reason, we have%
\begin{align*}
\nabla_{\operatorname*{Col}\left(  w\rightharpoonup T\right)  }^{-}  &
=\sum_{u\in\mathcal{C}\left(  w\rightharpoonup T\right)  }\left(  -1\right)
^{u}u=\sum_{u\in\mathcal{C}\left(  T\right)  }\underbrace{\left(  -1\right)
^{wuw^{-1}}}_{\substack{=\left(  -1\right)  ^{u}\\\text{(by
(\ref{eq.intro.perms.signs.conj}))}}}wuw^{-1}\\
&  \ \ \ \ \ \ \ \ \ \ \ \ \ \ \ \ \ \ \ \ \left(
\begin{array}
[c]{c}%
\text{here, we have substituted }wuw^{-1}\text{ for }u\text{ in the sum,}\\
\text{since the map }\mathcal{C}\left(  T\right)  \rightarrow\mathcal{C}%
\left(  w\rightharpoonup T\right)  ,\ u\mapsto wuw^{-1}\\
\text{is a bijection}%
\end{array}
\right) \\
&  =\sum_{u\in\mathcal{C}\left(  T\right)  }\left(  -1\right)  ^{u}%
wuw^{-1}=w\underbrace{\left(  \sum_{u\in\mathcal{C}\left(  T\right)  }\left(
-1\right)  ^{u}u\right)  }_{=\nabla_{\operatorname*{Col}T}^{-}}w^{-1}%
=w\nabla_{\operatorname*{Col}T}^{-}w^{-1}.
\end{align*}
A similar computation (but simpler due to the absence of signs) shows that
$\nabla_{\operatorname*{Row}\left(  w\rightharpoonup T\right)  }%
=w\nabla_{\operatorname*{Row}T}w^{-1}$. Proposition
\ref{prop.symmetrizers.conj} is thus proved.
\end{proof}

We can use column antisymmetrizers to rewrite the definition of polytabloids
in the following compact form:

\begin{proposition}
\label{prop.symmetrizers.eT}Let $T$ be an $n$-tableau of any shape $D$. Then,
in the Young module $\mathcal{M}^{D}$, we have%
\[
\mathbf{e}_{T}=\nabla_{\operatorname*{Col}T}^{-}\overline{T}.
\]

\end{proposition}

\begin{proof}
The definition of $\nabla_{\operatorname*{Col}T}^{-}$ yields $\nabla
_{\operatorname*{Col}T}^{-}=\sum_{w\in\mathcal{C}\left(  T\right)  }\left(
-1\right)  ^{w}w$. Thus,%
\[
\nabla_{\operatorname*{Col}T}^{-}\overline{T}=\sum_{w\in\mathcal{C}\left(
T\right)  }\left(  -1\right)  ^{w}\underbrace{w\overline{T}}%
_{\substack{=\overline{wT}\\\text{(by (\ref{eq.def.tabloid.Sn-act.rewr}))}%
}}=\sum_{w\in\mathcal{C}\left(  T\right)  }\left(  -1\right)  ^{w}%
\overline{wT}=\mathbf{e}_{T}%
\]
(by the definition of $\mathbf{e}_{T}$). This proves Proposition
\ref{prop.symmetrizers.eT}.
\end{proof}

\subsubsection{The left ideal avatars of Young and Specht modules}

Using row symmetrizers and column antisymmetrizers, we can find new avatars
(i.e., isomorphic copies) of the Young and Specht modules as
subrepresentations of the left regular representation $\mathbf{k}\left[
S_{n}\right]  $.

We recall that any element $a$ of any ring $R$ gives rise to a left ideal
$Ra=\left\{  ra\ \mid\ r\in R\right\}  $ of $R$, called the \emph{principal
left ideal} of $R$ generated by $a$. But a left ideal of $R$ is the same thing
as a left $R$-submodule of the left regular $R$-module $R$. Applying this to
$R=\mathbf{k}\left[  S_{n}\right]  $, we see that any element $\mathbf{a}%
\in\mathbf{k}\left[  S_{n}\right]  $ gives rise to a principal left ideal
$\mathbf{k}\left[  S_{n}\right]  \cdot\mathbf{a}$ of $\mathbf{k}\left[
S_{n}\right]  $, that is, to a left $\mathbf{k}\left[  S_{n}\right]  $-module
of $\mathbf{k}\left[  S_{n}\right]  $, that is, to a subrepresentation of the
left regular representation $\mathbf{k}\left[  S_{n}\right]  $ of $S_{n}$.
This gives us an abundant supply of representations of $S_{n}$. Some of these,
as we will now see, are isomorphic to Young modules $\mathcal{M}^{D}$ or to
Specht modules $\mathcal{S}^{D}$:

\begin{theorem}
[left ideal avatars of Young and Specht modules]%
\label{thm.spechtmod.leftideal}Let $D$ be any diagram with $\left\vert
D\right\vert =n$. Let $T$ be an $n$-tableau of shape $D$. Then: \medskip

\textbf{(a)} The Young module $\mathcal{M}^{D}$ is isomorphic (as an $S_{n}%
$-representation) to the left ideal
\[
\mathbf{k}\left[  S_{n}\right]  \cdot\nabla_{\operatorname*{Row}T}=\left\{
\mathbf{a}\nabla_{\operatorname*{Row}T}\ \mid\ \mathbf{a}\in\mathbf{k}\left[
S_{n}\right]  \right\}
\]
of $\mathbf{k}\left[  S_{n}\right]  $. To be more specific: There is a unique
left $\mathbf{k}\left[  S_{n}\right]  $-linear map%
\[
\alpha:\mathbf{k}\left[  S_{n}\right]  \cdot\nabla_{\operatorname*{Row}%
T}\rightarrow\mathcal{M}^{D}%
\]
that sends $\nabla_{\operatorname*{Row}T}$ to $\overline{T}$, and this map is
an isomorphism of $S_{n}$-representations (i.e., of left $\mathbf{k}\left[
S_{n}\right]  $-modules). \medskip

\textbf{(b)} This isomorphism $\alpha$ sends the submodule%
\[
\mathbf{k}\left[  S_{n}\right]  \cdot\nabla_{\operatorname*{Col}T}^{-}%
\nabla_{\operatorname*{Row}T}=\left\{  \mathbf{a}\nabla_{\operatorname*{Col}%
T}^{-}\nabla_{\operatorname*{Row}T}\ \mid\ \mathbf{a}\in\mathbf{k}\left[
S_{n}\right]  \right\}
\]
of $\mathbf{k}\left[  S_{n}\right]  \cdot\nabla_{\operatorname*{Row}T}$ to the
Specht module $\mathcal{S}^{D}$. Thus, the Specht module $\mathcal{S}^{D}$ is
isomorphic to the left ideal $\mathbf{k}\left[  S_{n}\right]  \cdot
\nabla_{\operatorname*{Col}T}^{-}\nabla_{\operatorname*{Row}T}$ of
$\mathbf{k}\left[  S_{n}\right]  $.
\end{theorem}

The product $\nabla_{\operatorname*{Col}T}^{-}\nabla_{\operatorname*{Row}T}$
in Theorem \ref{thm.spechtmod.leftideal} \textbf{(b)} is often called the
\emph{Young symmetrizer} of $T$ (although sometimes the word is used instead
for a certain scalar multiple of it, which we will see later). We will denote
it by $\mathbf{E}_{T}$ in a later section.

\begin{example}
\label{exa.spechtmod.leftideal.221}Let $n=5$ and $D=Y\left(  2,2,1\right)  $
and $T=12\backslash\backslash34\backslash\backslash5$. Then,%
\begin{align*}
\nabla_{\operatorname*{Col}T}^{-}  &  =\nabla_{\left\{  1,3,5\right\}  }%
^{-}\nabla_{\left\{  2,4\right\}  }^{-}\ \ \ \ \ \ \ \ \ \ \text{and}\\
\nabla_{\operatorname*{Row}T}  &  =\nabla_{\left\{  1,2\right\}  }%
\nabla_{\left\{  3,4\right\}  }\underbrace{\nabla_{\left\{  5\right\}  }}%
_{=1}=\nabla_{\left\{  1,2\right\}  }\nabla_{\left\{  3,4\right\}  },
\end{align*}
so we obtain the left ideal%
\[
\mathbf{k}\left[  S_{n}\right]  \cdot\nabla_{\operatorname*{Col}T}^{-}%
\nabla_{\operatorname*{Row}T}=\mathbf{k}\left[  S_{5}\right]  \cdot
\nabla_{\left\{  1,3,5\right\}  }^{-}\nabla_{\left\{  2,4\right\}  }^{-}%
\nabla_{\left\{  1,2\right\}  }\nabla_{\left\{  3,4\right\}  }.
\]
Theorem \ref{thm.spechtmod.leftideal} \textbf{(b)} claims that this left ideal
of $\mathbf{k}\left[  S_{5}\right]  $ is isomorphic to the Specht module
$\mathcal{S}^{D}=\mathcal{S}^{Y\left(  2,2,1\right)  }$. Note that taking a
different $n$-tableau $T$ of shape $D$ would yield a different but isomorphic
left ideal. Thus, Theorem \ref{thm.spechtmod.leftideal} \textbf{(b)} produces
not one, but many different isomorphic copies of a Specht module inside
$\mathbf{k}\left[  S_{n}\right]  $. And Theorem \ref{thm.spechtmod.leftideal}
\textbf{(a)} does likewise for Young modules.
\end{example}

There are many ways to prove Theorem \ref{thm.spechtmod.leftideal} (note that
part \textbf{(b)} follows easily from part \textbf{(a)}, but part \textbf{(a)}
presents some -- at least notational -- difficulties). In particular, we could
generalize it to arbitrary groups, since the Young module $\mathcal{M}^{D}$ is
the permutation module corresponding to a left $S_{n}$-set that can be
identified with a quotient of $S_{n}$. However, we shall take a different
road, which will net us some helpful results along the way that will be used
later on.

\subsubsection{The row-to-row sums}

This road leads us through what we call the \emph{row-to-row sums}:

\begin{definition}
\label{def.row-to-row.row-to-row}Let $D$ be any diagram. Let $\overline{S}$
and $\overline{T}$ be two $n$-tabloids of shape $D$. Then, we define the
\emph{row-to-row sum} $\nabla_{\overline{S},\overline{T}}\in\mathbf{k}\left[
S_{n}\right]  $ by%
\[
\nabla_{\overline{S},\overline{T}}:=\sum_{\substack{u\in S_{n};\\u\overline
{T}=\overline{S}}}u.
\]
(Of course, we could just as well rewrite $u\overline{T}$ as $u\rightharpoonup
\overline{T}$ or as $\overline{uT}$ here.)
\end{definition}

\begin{example}
Let $n=4$ and $D=Y\left(  2,2\right)  $ and $\overline{S}%
=\ytableausetup{tabloids}\ytableaushort{24,13}\ytableausetup{notabloids}$ and
$\overline{T}%
=\ytableausetup{tabloids}\ytableaushort{12,34}\ytableausetup{notabloids}$.
Then,%
\begin{align*}
\nabla_{\overline{S},\overline{T}}  &  =\sum_{\substack{u\in S_{n}%
;\\u\overline{T}=\overline{S}}}u=\sum_{\substack{u\in S_{n};\\u\left(
\left\{  1,2\right\}  \right)  =\left\{  2,4\right\}  ;\\u\left(  \left\{
3,4\right\}  \right)  =\left\{  1,3\right\}  }}u\\
&  =\operatorname*{oln}\left(  2413\right)  +\operatorname*{oln}\left(
2431\right)  +\operatorname*{oln}\left(  4213\right)  +\operatorname*{oln}%
\left(  4231\right)  .
\end{align*}
(See Subsection \ref{subsec.intro.nots.perms} for the meaning of
$\operatorname*{oln}$.)
\end{example}

We begin by proving some easy properties of these elements. First, we observe
that row-to-row sums include row symmetrizers as a particular case:

\begin{proposition}
\label{prop.row-to-row.rsym}Let $D$ be any diagram. Let $T$ be any $n$-tableau
of shape $D$. Then,%
\[
\nabla_{\overline{T},\overline{T}}=\nabla_{\operatorname*{Row}T}.
\]

\end{proposition}

\begin{proof}
Definition \ref{def.symmetrizers.symmetrizers} yields $\nabla
_{\operatorname*{Row}T}=\sum_{w\in\mathcal{R}\left(  T\right)  }w$. But
Proposition \ref{prop.tableau.Sn-act.1} \textbf{(a)} shows that
\[
\mathcal{R}\left(  T\right)  =\left\{  w\in S_{n}\ \mid\ \overline
{w\rightharpoonup T}=\overline{T}\right\}  =\left\{  w\in S_{n}\ \mid
\ w\overline{T}=\overline{T}\right\}
\]
(since $\overline{w\rightharpoonup T}=w\rightharpoonup\overline{T}%
=w\overline{T}$ for each $w\in S_{n}$). Thus, the summation sign $\sum
_{w\in\mathcal{R}\left(  T\right)  }$ can be rewritten as $\sum
_{\substack{w\in S_{n};\\w\overline{T}=\overline{T}}}$. Hence,%
\[
\sum_{w\in\mathcal{R}\left(  T\right)  }w=\sum_{\substack{w\in S_{n}%
;\\w\overline{T}=\overline{T}}}w=\sum_{\substack{u\in S_{n};\\u\overline
{T}=\overline{T}}}u
\]
(here, we have renamed the summation index $w$ as $u$). Therefore,%
\[
\nabla_{\operatorname*{Row}T}=\sum_{w\in\mathcal{R}\left(  T\right)  }%
w=\sum_{\substack{u\in S_{n};\\u\overline{T}=\overline{T}}}u=\nabla
_{\overline{T},\overline{T}}%
\]
(since $\nabla_{\overline{T},\overline{T}}$ was defined as $\sum
_{\substack{u\in S_{n};\\u\overline{T}=\overline{T}}}u$). This proves
Proposition \ref{prop.row-to-row.rsym}.
\end{proof}

The next property is a kind of equivariance:

\begin{proposition}
\label{prop.row-to-row.Sn-act}Let $D$ be any diagram. Let $\overline{S}$ and
$\overline{T}$ be two $n$-tabloids of shape $D$. Let $g\in S_{n}$ and $h\in
S_{n}$ be two permutations. Then,%
\[
\nabla_{g\overline{S},h\overline{T}}=g\nabla_{\overline{S},\overline{T}}%
h^{-1}.
\]

\end{proposition}

\begin{proof}
Recall that $S_{n}$ is a group. Thus, the map%
\begin{align*}
S_{n}  &  \rightarrow S_{n},\\
u  &  \mapsto guh^{-1}%
\end{align*}
is a bijection (with inverse given by $u\mapsto g^{-1}uh$). The definition of
$\nabla_{\overline{S},\overline{T}}$ yields%
\begin{equation}
\nabla_{\overline{S},\overline{T}}=\sum_{\substack{u\in S_{n};\\u\overline
{T}=\overline{S}}}u. \label{pf.prop.row-to-row.Sn-act.1}%
\end{equation}
The same definition (but applied to $g\overline{S}$ and $h\overline{T}$
instead of $\overline{S}$ and $\overline{T}$) yields%
\begin{equation}
\nabla_{g\overline{S},h\overline{T}}=\sum_{\substack{u\in S_{n};\\uh\overline
{T}=g\overline{S}}}u=\sum_{\substack{u\in S_{n};\\guh^{-1}h\overline
{T}=g\overline{S}}}guh^{-1} \label{pf.prop.row-to-row.Sn-act.2}%
\end{equation}
(here, we have substituted $guh^{-1}$ for $u$ in the sum, since the map
$S_{n}\rightarrow S_{n},\ u\mapsto guh^{-1}$ is a bijection). However, for any
permutation $u\in S_{n}$, we have the chain of equivalences%
\[
\left(  g\underbrace{uh^{-1}h}_{=u}\overline{T}=g\overline{S}\right)
\ \Longleftrightarrow\ \left(  gu\overline{T}=g\overline{S}\right)
\ \Longleftrightarrow\ \left(  u\overline{T}=\overline{S}\right)
\]
(since the statement $\left(  u\overline{T}=\overline{S}\right)  $ implies
$\left(  gu\overline{T}=g\overline{S}\right)  $ by left-multiplying it by $g$,
whereas conversely, the statement $\left(  gu\overline{T}=g\overline
{S}\right)  $ implies $\left(  u\overline{T}=\overline{S}\right)  $ by
left-multiplying it by $g^{-1}$). Hence, we can replace the condition
\textquotedblleft$guh^{-1}h\overline{T}=g\overline{S}$\textquotedblright%
\ under the summation sign in (\ref{pf.prop.row-to-row.Sn-act.2}) by
\textquotedblleft$u\overline{T}=\overline{S}$\textquotedblright. Thus,
(\ref{pf.prop.row-to-row.Sn-act.2}) rewrites as%
\[
\nabla_{g\overline{S},h\overline{T}}=\sum_{\substack{u\in S_{n};\\u\overline
{T}=\overline{S}}}guh^{-1}=g\underbrace{\left(  \sum_{\substack{u\in
S_{n};\\u\overline{T}=\overline{S}}}u\right)  }_{\substack{=\nabla
_{\overline{S},\overline{T}}\\\text{(by (\ref{pf.prop.row-to-row.Sn-act.1}))}%
}}h^{-1}=g\nabla_{\overline{S},\overline{T}}h^{-1}.
\]
This proves Proposition \ref{prop.row-to-row.Sn-act}.
\end{proof}

\begin{corollary}
\label{cor.row-to-row.surj}Let $D$ be any diagram. Let $\overline{T}$ be an
$n$-tabloid of shape $D$. Let $w\in S_{n}$. Then, $\nabla_{w\overline
{T},\overline{T}}=w\nabla_{\operatorname*{Row}T}$.
\end{corollary}

\begin{proof}
Proposition \ref{prop.row-to-row.Sn-act} (applied to $\overline{S}%
=\overline{T}$, $g=w$ and $h=\operatorname*{id}$) yields $\nabla
_{w\overline{T},\operatorname*{id}\overline{T}}=w\nabla_{\overline
{T},\overline{T}}\underbrace{\operatorname*{id}\nolimits^{-1}}%
_{=\operatorname*{id}}=w\underbrace{\nabla_{\overline{T},\overline{T}}%
}_{\substack{=\nabla_{\operatorname*{Row}T}\\\text{(by Proposition
\ref{prop.row-to-row.rsym})}}}=w\nabla_{\operatorname*{Row}T}$. In view of
$\operatorname*{id}\overline{T}=\overline{T}$, we can rewrite this as
$\nabla_{w\overline{T},\overline{T}}=w\nabla_{\operatorname*{Row}T}$. This
proves Corollary \ref{cor.row-to-row.surj}.
\end{proof}

\begin{corollary}
\label{cor.row-to-row.in-left-ideal}Let $D$ be any diagram. Let $\overline{S}$
and $\overline{T}$ be two $n$-tabloids of shape $D$. Then, $\nabla
_{\overline{S},\overline{T}}\in\mathbf{k}\left[  S_{n}\right]  \cdot
\nabla_{\operatorname*{Row}T}$.
\end{corollary}

\begin{proof}
Proposition \ref{prop.tableau.Sn-act.1} \textbf{(c)} shows that the
$n$-tabloid $\overline{S}$ can be written as $w\rightharpoonup\overline{T}$
for some $w\in S_{n}$. Consider this $w$. Thus, $\overline{S}=w\rightharpoonup
\overline{T}=w\overline{T}$, so that%
\begin{align*}
\nabla_{\overline{S},\overline{T}}  &  =\nabla_{w\overline{T},\overline{T}%
}=\underbrace{w}_{\in\mathbf{k}\left[  S_{n}\right]  }\nabla
_{\operatorname*{Row}T}\ \ \ \ \ \ \ \ \ \ \left(  \text{by Corollary
\ref{cor.row-to-row.surj}}\right) \\
&  \in\mathbf{k}\left[  S_{n}\right]  \cdot\nabla_{\operatorname*{Row}T}.
\end{align*}
This proves Corollary \ref{cor.row-to-row.in-left-ideal}.
\end{proof}

\subsubsection{Proving the left ideal avatars}

Now, we shall define a map $\zeta:\mathcal{M}^{D}\rightarrow\mathbf{k}\left[
S_{n}\right]  \cdot\nabla_{\operatorname*{Row}T}$ that will later turn out to
be inverse to the map $\alpha$ in Theorem \ref{thm.spechtmod.leftideal}
\textbf{(a)}:

\begin{proposition}
\label{prop.row-to-row.zeta}Let $D$ be any diagram with $\left\vert
D\right\vert =n$. Let $T$ be an $n$-tableau of shape $D$. Consider the
$\mathbf{k}$-linear map%
\[
\zeta:\mathcal{M}^{D}\rightarrow\mathbf{k}\left[  S_{n}\right]  \cdot
\nabla_{\operatorname*{Row}T}%
\]
that sends each $n$-tabloid $\overline{S}$ of shape $D$ to $\nabla
_{\overline{S},\overline{T}}$. This map $\zeta$ is well-defined and is a left
$\mathbf{k}\left[  S_{n}\right]  $-module isomorphism.
\end{proposition}

\begin{proof}
Let us first show that the map $\zeta$ is well-defined. First of all, we
recall that the $n$-tabloids form a basis of the $\mathbf{k}$-module
$\mathcal{M}^{D}$ (by the definition of $\mathcal{M}^{D}$). Hence, we can
uniquely define a $\mathbf{k}$-linear map from $\mathcal{M}^{D}$ by specifying
its values on the $n$-tabloids. Thus, in particular, the definition of $\zeta$
given in Proposition \ref{prop.row-to-row.zeta} is legitimate, once we show
that the value $\nabla_{\overline{S},\overline{T}}$ really belongs to
$\mathbf{k}\left[  S_{n}\right]  \cdot\nabla_{\operatorname*{Row}T}$ for each
$n$-tabloid $\overline{S}$ of shape $D$. But the latter claim is just
Corollary \ref{cor.row-to-row.in-left-ideal}. Thus, the definition of $\zeta$
is legitimate -- i.e., the map $\zeta$ is well-defined.

Next, let us show some further properties of $\zeta$:

\begin{statement}
\textit{Claim 1:} For each $w\in S_{n}$, we have $\zeta\left(  \overline
{wT}\right)  =w\nabla_{\operatorname*{Row}T}$.
\end{statement}

\begin{proof}
[Proof of Claim 1.]Let $w\in S_{n}$. The definition of $\zeta$ yields
\begin{align*}
\zeta\left(  \overline{wT}\right)   &  =\nabla_{\overline{wT},\overline{T}%
}=\nabla_{w\overline{T},\overline{T}}\ \ \ \ \ \ \ \ \ \ \left(  \text{since
}\overline{wT}=w\overline{T}\right) \\
&  =w\nabla_{\operatorname*{Row}T}\ \ \ \ \ \ \ \ \ \ \left(  \text{by
Corollary \ref{cor.row-to-row.surj}}\right)  .
\end{align*}
This proves Claim 1.
\end{proof}

\begin{statement}
\textit{Claim 2:} The map $\zeta$ is left $\mathbf{k}\left[  S_{n}\right]  $-linear.
\end{statement}

\begin{proof}
[Proof of Claim 2.]Let $w\in S_{n}$, and let $\mathbf{a}\in\mathcal{M}^{D}$.
We shall show that $\zeta\left(  w\mathbf{a}\right)  =w\zeta\left(
\mathbf{a}\right)  $.

Indeed, both sides of this equality are $\mathbf{k}$-linear in $\mathbf{a}$.
Hence, we can WLOG assume that $\mathbf{a}$ is an $n$-tabloid of shape $D$
(since the $\mathbf{k}$-module $\mathcal{M}^{D}$ is spanned by such
$n$-tabloids). Assume this. Thus, $\mathbf{a}=\overline{S}$ for some
$n$-tabloid $\overline{S}$ of shape $D$. Consider this $\overline{S}$. From
$\mathbf{a}=\overline{S}$, we obtain $\zeta\left(  \mathbf{a}\right)
=\zeta\left(  \overline{S}\right)  =\nabla_{\overline{S},\overline{T}}$ (by
the definition of $\zeta$). Thus,
\begin{equation}
w\zeta\left(  \mathbf{a}\right)  =w\nabla_{\overline{S},\overline{T}}.
\label{pf.prop.row-to-row.zeta.1}%
\end{equation}
However, from $\mathbf{a}=\overline{S}$, we also obtain $w\mathbf{a}%
=w\overline{S}=\overline{wS}$ and therefore
\begin{align*}
\zeta\left(  w\mathbf{a}\right)   &  =\zeta\left(  \overline{wS}\right)
=\nabla_{\overline{wS},\overline{T}}\ \ \ \ \ \ \ \ \ \ \left(  \text{by the
definition of }\zeta\right) \\
&  =\nabla_{w\overline{S},\operatorname*{id}\overline{T}}%
\ \ \ \ \ \ \ \ \ \ \left(  \text{since }\overline{wS}=w\overline{S}\text{ and
}\overline{T}=\operatorname*{id}\overline{T}\right) \\
&  =w\nabla_{\overline{S},\overline{T}}\underbrace{\operatorname*{id}%
\nolimits^{-1}}_{=\operatorname*{id}=1}\ \ \ \ \ \ \ \ \ \ \left(
\begin{array}
[c]{c}%
\text{by Proposition \ref{prop.row-to-row.Sn-act},}\\
\text{applied to }g=w\text{ and }h=\operatorname*{id}%
\end{array}
\right) \\
&  =w\nabla_{\overline{S},\overline{T}}=w\zeta\left(  \mathbf{a}\right)
\ \ \ \ \ \ \ \ \ \ \left(  \text{by (\ref{pf.prop.row-to-row.zeta.1}%
)}\right)  .
\end{align*}

Forget that we fixed $w$ and $\mathbf{a}$. We thus have proved that
$\zeta\left(  w\mathbf{a}\right)  =w\zeta\left(  \mathbf{a}\right)  $ for all
$w\in S_{n}$ and all $\mathbf{a}\in\mathcal{M}^{D}$. In other words, the map
$\zeta$ is $S_{n}$-equivariant. Hence, $\zeta$ is a morphism of $S_{n}%
$-representations (since $\zeta$ is $\mathbf{k}$-linear). By Proposition
\ref{prop.rep.G-rep.mor=mor}, this entails that $\zeta$ is a left
$\mathbf{k}\left[  S_{n}\right]  $-module morphism, i.e., a left
$\mathbf{k}\left[  S_{n}\right]  $-linear map. This proves Claim 2.
\end{proof}

\begin{statement}
\textit{Claim 3:} The map $\zeta$ is surjective.
\end{statement}

\begin{proof}
[Proof of Claim 3.]Let $\mathbf{a}\in\mathbf{k}\left[  S_{n}\right]  $. We
shall show that $\mathbf{a}\nabla_{\operatorname*{Row}T}\in\zeta\left(
\mathcal{M}^{D}\right)  $.

Indeed, the definition of $\zeta$ yields $\zeta\left(  \overline{T}\right)
=\nabla_{\overline{T},\overline{T}}=\nabla_{\operatorname*{Row}T}$ (by
Proposition \ref{prop.row-to-row.rsym}). However, the map $\zeta$ is left
$\mathbf{k}\left[  S_{n}\right]  $-linear (by Claim 2), and thus we have
$\zeta\left(  \mathbf{a}\overline{T}\right)  =\mathbf{a}\underbrace{\zeta
\left(  \overline{T}\right)  }_{=\nabla_{\operatorname*{Row}T}}=\mathbf{a}%
\nabla_{\operatorname*{Row}T}$. Thus, $\mathbf{a}\nabla_{\operatorname*{Row}%
T}=\zeta\left(  \mathbf{a}\overline{T}\right)  \in\zeta\left(  \mathcal{M}%
^{D}\right)  $.

Forget that we fixed $\mathbf{a}$. We thus have shown that $\mathbf{a}%
\nabla_{\operatorname*{Row}T}\in\zeta\left(  \mathcal{M}^{D}\right)  $ for
each $\mathbf{a}\in\mathbf{k}\left[  S_{n}\right]  $. In other words,
$\left\{  \mathbf{a}\nabla_{\operatorname*{Row}T}\ \mid\ \mathbf{a}%
\in\mathbf{k}\left[  S_{n}\right]  \right\}  \subseteq\zeta\left(
\mathcal{M}^{D}\right)  $. Now,%
\[
\mathbf{k}\left[  S_{n}\right]  \cdot\nabla_{\operatorname*{Row}T}=\left\{
\mathbf{a}\nabla_{\operatorname*{Row}T}\ \mid\ \mathbf{a}\in\mathbf{k}\left[
S_{n}\right]  \right\}  \subseteq\zeta\left(  \mathcal{M}^{D}\right)  .
\]
In other words, the map $\zeta$ is surjective. This proves Claim 3.
\end{proof}

\begin{statement}
\textit{Claim 4:} Let $a_{\overline{S}}\in\mathbf{k}$ be a scalar for each
$n$-tabloid $\overline{S}$ of shape $D$. Then,%
\[
\zeta\left(  \sum_{\substack{\overline{S}\text{ is an }n\text{-tabloid}%
\\\text{of shape }D}}a_{\overline{S}}\overline{S}\right)  =\sum_{u\in S_{n}%
}a_{u\overline{T}}u.
\]

\end{statement}

\begin{proof}
[Proof of Claim 4.]The map $\zeta$ is $\mathbf{k}$-linear. Thus,
\begin{align*}
&  \zeta\left(  \sum_{\substack{\overline{S}\text{ is an }n\text{-tabloid}%
\\\text{of shape }D}}a_{\overline{S}}\overline{S}\right) \\
&  =\sum_{\substack{\overline{S}\text{ is an }n\text{-tabloid}\\\text{of shape
}D}}a_{\overline{S}}\underbrace{\zeta\left(  \overline{S}\right)
}_{\substack{=\nabla_{\overline{S},\overline{T}}\\\text{(by the definition of
}\zeta\text{)}}}\\
&  =\sum_{\substack{\overline{S}\text{ is an }n\text{-tabloid}\\\text{of shape
}D}}a_{\overline{S}}\underbrace{\nabla_{\overline{S},\overline{T}}%
}_{\substack{=\sum_{\substack{u\in S_{n};\\u\overline{T}=\overline{S}%
}}u\\\text{(by the definition of }\nabla_{\overline{S},\overline{T}}\text{)}%
}}=\sum_{\substack{\overline{S}\text{ is an }n\text{-tabloid}\\\text{of shape
}D}}a_{\overline{S}}\sum_{\substack{u\in S_{n};\\u\overline{T}=\overline{S}%
}}u\\
&  =\sum_{\substack{\overline{S}\text{ is an }n\text{-tabloid}\\\text{of shape
}D}}\ \ \sum_{\substack{u\in S_{n};\\u\overline{T}=\overline{S}}%
}\underbrace{a_{\overline{S}}}_{\substack{=a_{u\overline{T}}\\\text{(since
}\overline{S}=u\overline{T}\text{)}}}u=\underbrace{\sum_{\substack{\overline
{S}\text{ is an }n\text{-tabloid}\\\text{of shape }D}}\ \ \sum_{\substack{u\in
S_{n};\\u\overline{T}=\overline{S}}}}_{\substack{=\sum_{u\in S_{n}%
}\\\text{(since each }u\in S_{n}\text{ appears in the}\\\text{inner sum for
exactly one }\overline{S}\text{)}}}a_{u\overline{T}}u\\
&  =\sum_{u\in S_{n}}a_{u\overline{T}}u.
\end{align*}
This proves Claim 4.
\end{proof}

\begin{statement}
\textit{Claim 5:} The map $\zeta$ is injective.
\end{statement}

\begin{proof}
[Proof of Claim 5.]Let $\mathbf{s}\in\operatorname*{Ker}\zeta$. We shall show
that $\mathbf{s}=0$.

We have $\mathbf{s}\in\operatorname*{Ker}\zeta\subseteq\mathcal{M}^{D}$.
Hence, we can write $\mathbf{s}$ as a $\mathbf{k}$-linear combination of
$n$-tabloids of shape $D$ (since $\mathcal{M}^{D}$ is spanned by such
$n$-tabloids). In other words, we can write $\mathbf{s}$ as
\begin{equation}
\mathbf{s}=\sum_{\substack{\overline{S}\text{ is an }n\text{-tabloid}%
\\\text{of shape }D}}a_{\overline{S}}\overline{S}
\label{pf.prop.row-to-row.zeta.5}%
\end{equation}
for some coefficients $a_{\overline{S}}\in\mathbf{k}$. Consider these
$a_{\overline{S}}$.

Applying $\zeta$ to the equality (\ref{pf.prop.row-to-row.zeta.5}), we find%
\begin{align*}
\zeta\left(  \mathbf{s}\right)   &  =\zeta\left(  \sum_{\substack{\overline
{S}\text{ is an }n\text{-tabloid}\\\text{of shape }D}}a_{\overline{S}%
}\overline{S}\right)  =\sum_{u\in S_{n}}a_{u\overline{T}}%
u\ \ \ \ \ \ \ \ \ \ \left(  \text{by Claim 4}\right) \\
&  =\left(  a_{u\overline{T}}\right)  _{u\in S_{n}}\ \ \ \ \ \ \ \ \ \ \left(
\text{by (\ref{eq.rmk.monalg.compare.alm=})}\right)  .
\end{align*}
Thus, $\left(  a_{u\overline{T}}\right)  _{u\in S_{n}}=\zeta\left(
\mathbf{s}\right)  =0$ (since $\mathbf{s}\in\operatorname*{Ker}\zeta$). In
other words,
\begin{equation}
a_{u\overline{T}}=0\ \ \ \ \ \ \ \ \ \ \text{for each }u\in S_{n}.
\label{pf.prop.row-to-row.zeta.6}%
\end{equation}

Now, let $\overline{S}$ be any $n$-tabloid of shape $D$. Then, Proposition
\ref{prop.tableau.Sn-act.1} \textbf{(c)} shows that the $n$-tabloid
$\overline{S}$ can be written as $w\rightharpoonup\overline{T}$ for some $w\in
S_{n}$. Consider this $w$. Thus, $\overline{S}=w\rightharpoonup\overline
{T}=w\overline{T}$, so that $a_{\overline{S}}=a_{w\overline{T}}=0$ (by
(\ref{pf.prop.row-to-row.zeta.6}), applied to $u=w$).

Forget that we fixed $\overline{S}$. We thus have shown that $a_{\overline{S}%
}=0$ for each $n$-tabloid $\overline{S}$ of shape $D$. Hence, we can rewrite
(\ref{pf.prop.row-to-row.zeta.5}) as $\mathbf{s}=\sum_{\substack{\overline
{S}\text{ is an }n\text{-tabloid}\\\text{of shape }D}}0\overline{S}=0$.

Forget that we fixed $\mathbf{s}$. We thus have shown that $\mathbf{s}=0$ for
each $\mathbf{s}\in\operatorname*{Ker}\zeta$. In other words,
$\operatorname*{Ker}\zeta=0$. Since $\zeta$ is a $\mathbf{k}$-linear map, this
entails that $\zeta$ is injective. This proves Claim 5.
\end{proof}

The map $\zeta$ is injective (by Claim 5) and surjective (by Claim 3); thus,
it is bijective, i.e., invertible. Since it is furthermore left $\mathbf{k}%
\left[  S_{n}\right]  $-linear (by Claim 2), we thus conclude that $\zeta$ is
a left $\mathbf{k}\left[  S_{n}\right]  $-module isomorphism. This proves
Proposition \ref{prop.row-to-row.zeta}.
\end{proof}

We can now easily prove Theorem \ref{thm.spechtmod.leftideal}:

\begin{proof}
[Proof of Theorem \ref{thm.spechtmod.leftideal}.]\textbf{(a)} We must prove
the following two claims:

\begin{statement}
\textit{Claim 1:} There is a unique left $\mathbf{k}\left[  S_{n}\right]
$-linear map%
\[
\alpha:\mathbf{k}\left[  S_{n}\right]  \cdot\nabla_{\operatorname*{Row}%
T}\rightarrow\mathcal{M}^{D}%
\]
that sends $\nabla_{\operatorname*{Row}T}$ to $\overline{T}$.
\end{statement}

\begin{statement}
\textit{Claim 2:} This map $\alpha$ is an isomorphism of $S_{n}$-representations.
\end{statement}

\begin{proof}
[Proof of Claim 1.]The uniqueness of $\alpha$ is easy to see: If
$\alpha:\mathbf{k}\left[  S_{n}\right]  \cdot\nabla_{\operatorname*{Row}%
T}\rightarrow\mathcal{M}^{D}$ is a left $\mathbf{k}\left[  S_{n}\right]
$-linear map that sends $\nabla_{\operatorname*{Row}T}$ to $\overline{T}$,
then each $\mathbf{a}\in\mathbf{k}\left[  S_{n}\right]  $ satisfies
\begin{align*}
\alpha\left(  \mathbf{a}\nabla_{\operatorname*{Row}T}\right)   &
=\mathbf{a}\underbrace{\alpha\left(  \nabla_{\operatorname*{Row}T}\right)
}_{\substack{=\overline{T}\\\text{(since }\alpha\text{ sends }\nabla
_{\operatorname*{Row}T}\text{ to }\overline{T}\text{)}}%
}\ \ \ \ \ \ \ \ \ \ \left(  \text{since }\alpha\text{ is left }%
\mathbf{k}\left[  S_{n}\right]  \text{-linear}\right) \\
&  =\mathbf{a}\overline{T},
\end{align*}
and this determines all values of $\alpha$ (since each element of
$\mathbf{k}\left[  S_{n}\right]  \cdot\nabla_{\operatorname*{Row}T}$ has the
form $\mathbf{a}\nabla_{\operatorname*{Row}T}$ for some $\mathbf{a}%
\in\mathbf{k}\left[  S_{n}\right]  $). Thus, the map $\alpha$ is unique. It
remains to prove that such a map exists.

To do so, we construct it explicitly: Consider the left $\mathbf{k}\left[
S_{n}\right]  $-module isomorphism $\zeta:\mathcal{M}^{D}\rightarrow
\mathbf{k}\left[  S_{n}\right]  \cdot\nabla_{\operatorname*{Row}T}$
constructed in Proposition \ref{prop.row-to-row.zeta}. The inverse $\zeta
^{-1}$ of $\zeta$ is then a left $\mathbf{k}\left[  S_{n}\right]  $-module
isomorphism from $\mathbf{k}\left[  S_{n}\right]  \cdot\nabla
_{\operatorname*{Row}T}$ to $\mathcal{M}^{D}$. Moreover, the definition of
$\zeta$ yields $\zeta\left(  \overline{T}\right)  =\nabla_{\overline
{T},\overline{T}}=\nabla_{\operatorname*{Row}T}$ (by Proposition
\ref{prop.row-to-row.rsym}). Thus, $\zeta^{-1}\left(  \nabla
_{\operatorname*{Row}T}\right)  =\overline{T}$. In other words, $\zeta^{-1}$
sends $\nabla_{\operatorname*{Row}T}$ to $\overline{T}$. Hence, there exists a
left $\mathbf{k}\left[  S_{n}\right]  $-linear map $\alpha:\mathbf{k}\left[
S_{n}\right]  \cdot\nabla_{\operatorname*{Row}T}\rightarrow\mathcal{M}^{D}$
that sends $\nabla_{\operatorname*{Row}T}$ to $\overline{T}$: namely,
$\zeta^{-1}$ qualifies as such a map $\alpha$. This proves the existence of
$\alpha$, and thus Claim 1 is fully proved.
\end{proof}

\begin{proof}
[Proof of Claim 2.]Consider the left $\mathbf{k}\left[  S_{n}\right]  $-module
isomorphism $\zeta:\mathcal{M}^{D}\rightarrow\mathbf{k}\left[  S_{n}\right]
\cdot\nabla_{\operatorname*{Row}T}$ constructed in Proposition
\ref{prop.row-to-row.zeta}. Then, its inverse $\zeta^{-1}$ is also a left
$\mathbf{k}\left[  S_{n}\right]  $-module isomorphism, i.e., an isomorphism of
$S_{n}$-representations.

However, this inverse $\zeta^{-1}$ is precisely the map $\alpha$ whose
existence was claimed in Claim 1 (because we saw in the proof of Claim 1 that
$\zeta^{-1}$ qualifies as such a map $\alpha$, but since $\alpha$ is unique,
this means that $\zeta^{-1}$ \textbf{is} $\alpha$). Thus, the map $\alpha$ is
an isomorphism of $S_{n}$-representations (since $\zeta^{-1}$ is an
isomorphism of $S_{n}$-representations). This proves Claim 2.
\end{proof}

The proof of Theorem \ref{thm.spechtmod.leftideal} \textbf{(a)} is now
complete. \medskip

\textbf{(b)} We must prove that $\alpha$ sends $\mathbf{k}\left[
S_{n}\right]  \cdot\nabla_{\operatorname*{Col}T}^{-}\nabla
_{\operatorname*{Row}T}$ to $\mathcal{S}^{D}$.

The Specht module $\mathcal{S}^{D}$ is defined as the span of all polytabloids
of shape $D$. Thus,
\begin{equation}
\mathcal{S}^{D}=\operatorname*{span}\nolimits_{\mathbf{k}}\left\{
\mathbf{e}_{S}\ \mid\ S\text{ is an }n\text{-tableau of shape }D\right\}  .
\label{pf.thm.spechtmod.leftideal.b.1}%
\end{equation}

However, each $n$-tableau of shape $D$ can be written as $w\rightharpoonup T$
for some $w\in S_{n}$ (by Proposition \ref{prop.tableau.Sn-act.1}
\textbf{(b)}). Thus,%
\[
\left\{  n\text{-tableaux of shape }D\right\}  =\left\{  w\rightharpoonup
T\ \mid\ w\in S_{n}\right\}  =\left\{  wT\ \mid\ w\in S_{n}\right\}
\]
(since $wT$ is a shorthand for $w\rightharpoonup T$). Therefore,%
\[
\left\{  \mathbf{e}_{S}\ \mid\ S\text{ is an }n\text{-tableau of shape
}D\right\}  =\left\{  \mathbf{e}_{wT}\ \mid\ w\in S_{n}\right\}  .
\]
Thus, we can rewrite (\ref{pf.thm.spechtmod.leftideal.b.1}) as%
\begin{equation}
\mathcal{S}^{D}=\operatorname*{span}\nolimits_{\mathbf{k}}\left\{
\mathbf{e}_{wT}\ \mid\ w\in S_{n}\right\}  .
\label{pf.thm.spechtmod.leftideal.b.2}%
\end{equation}

On the other hand, $\mathbf{k}\left[  S_{n}\right]  =\operatorname*{span}%
\nolimits_{\mathbf{k}}\left\{  w\ \mid\ w\in S_{n}\right\}  $. Multiplying
both sides of this equality with $\nabla_{\operatorname*{Col}T}^{-}%
\nabla_{\operatorname*{Row}T}$ on the right, we find%
\begin{align*}
\mathbf{k}\left[  S_{n}\right]  \cdot\nabla_{\operatorname*{Col}T}^{-}%
\nabla_{\operatorname*{Row}T}  &  =\operatorname*{span}\nolimits_{\mathbf{k}%
}\left\{  w\ \mid\ w\in S_{n}\right\}  \cdot\nabla_{\operatorname*{Col}T}%
^{-}\nabla_{\operatorname*{Row}T}\\
&  =\operatorname*{span}\nolimits_{\mathbf{k}}\left\{  w\nabla
_{\operatorname*{Col}T}^{-}\nabla_{\operatorname*{Row}T}\ \mid\ w\in
S_{n}\right\}
\end{align*}
(since $\operatorname*{span}\nolimits_{\mathbf{k}}\left\{  a_{1},a_{2}%
,\ldots,a_{k}\right\}  \cdot b=\operatorname*{span}\nolimits_{\mathbf{k}%
}\left\{  a_{1}b,\ a_{2}b,\ \ldots,\ a_{k}b\right\}  $ for any elements
$a_{1},a_{2},\ldots,a_{k},b$ of any $\mathbf{k}$-algebra). Applying the map
$\alpha$ to this equality, we find%
\begin{align}
&  \alpha\left(  \mathbf{k}\left[  S_{n}\right]  \cdot\nabla
_{\operatorname*{Col}T}^{-}\nabla_{\operatorname*{Row}T}\right) \nonumber\\
&  =\alpha\left(  \operatorname*{span}\nolimits_{\mathbf{k}}\left\{
w\nabla_{\operatorname*{Col}T}^{-}\nabla_{\operatorname*{Row}T}\ \mid\ w\in
S_{n}\right\}  \right) \nonumber\\
&  =\operatorname*{span}\nolimits_{\mathbf{k}}\left\{  \alpha\left(
w\nabla_{\operatorname*{Col}T}^{-}\nabla_{\operatorname*{Row}T}\right)
\ \mid\ w\in S_{n}\right\}  \label{pf.thm.spechtmod.leftideal.b.3}%
\end{align}
(since the map $\alpha$ is $\mathbf{k}$-linear and thus respects $\mathbf{k}%
$-linear combinations). But the map $\alpha$ is left $\mathbf{k}\left[
S_{n}\right]  $-linear. Hence,
\[
\alpha\left(  \nabla_{\operatorname*{Col}T}^{-}\nabla_{\operatorname*{Row}%
T}\right)  =\nabla_{\operatorname*{Col}T}^{-}\underbrace{\alpha\left(
\nabla_{\operatorname*{Row}T}\right)  }_{\substack{=\overline{T}\\\text{(since
}\alpha\text{ sends }\nabla_{\operatorname*{Row}T}\text{ to }\overline
{T}\text{)}}}=\nabla_{\operatorname*{Col}T}^{-}\overline{T}=\mathbf{e}_{T}%
\]
(by Proposition \ref{prop.symmetrizers.eT}). Moreover, let us recall again
that the map $\alpha$ is left $\mathbf{k}\left[  S_{n}\right]  $-linear.
Hence, each $w\in S_{n}$ satisfies%
\[
\alpha\left(  w\nabla_{\operatorname*{Col}T}^{-}\nabla_{\operatorname*{Row}%
T}\right)  =w\underbrace{\alpha\left(  \nabla_{\operatorname*{Col}T}^{-}%
\nabla_{\operatorname*{Row}T}\right)  }_{=\mathbf{e}_{T}}=w\mathbf{e}%
_{T}=\mathbf{e}_{wT}%
\]
(since Lemma \ref{lem.spechtmod.submod} \textbf{(a)} yields $\mathbf{e}%
_{wT}=w\mathbf{e}_{T}$). Thus, we can rewrite
(\ref{pf.thm.spechtmod.leftideal.b.3}) as%
\[
\alpha\left(  \mathbf{k}\left[  S_{n}\right]  \cdot\nabla_{\operatorname*{Col}%
T}^{-}\nabla_{\operatorname*{Row}T}\right)  =\operatorname*{span}%
\nolimits_{\mathbf{k}}\left\{  \mathbf{e}_{wT}\ \mid\ w\in S_{n}\right\}  .
\]
Comparing this with (\ref{pf.thm.spechtmod.leftideal.b.2}), we obtain%
\[
\alpha\left(  \mathbf{k}\left[  S_{n}\right]  \cdot\nabla_{\operatorname*{Col}%
T}^{-}\nabla_{\operatorname*{Row}T}\right)  =\mathcal{S}^{D}.
\]
In other words, $\alpha$ sends $\mathbf{k}\left[  S_{n}\right]  \cdot
\nabla_{\operatorname*{Col}T}^{-}\nabla_{\operatorname*{Row}T}$ to
$\mathcal{S}^{D}$. Since $\alpha$ is an isomorphism of $S_{n}$-representations
(by part \textbf{(a)}), this entails that $\mathcal{S}^{D}$ is isomorphic to
$\mathbf{k}\left[  S_{n}\right]  \cdot\nabla_{\operatorname*{Col}T}^{-}%
\nabla_{\operatorname*{Row}T}$ as $S_{n}$-representations. Thus, the proof of
Theorem \ref{thm.spechtmod.leftideal} \textbf{(b)} is complete.
\end{proof}

\subsubsection{Explicit formulas for $\alpha^{-1}$}

From our above proof of Theorem \ref{thm.spechtmod.leftideal}, we can recover
some explicit formulas for the inverse $\alpha^{-1}$ of the isomorphism
$\alpha$:

\begin{proposition}
\label{prop.spechtmod.leftideal.al-1}Let $D$ be any diagram with $\left\vert
D\right\vert =n$. Let $T$ be an $n$-tableau of shape $D$. Consider the
isomorphism%
\[
\alpha:\mathbf{k}\left[  S_{n}\right]  \cdot\nabla_{\operatorname*{Row}%
T}\rightarrow\mathcal{M}^{D}%
\]
from Theorem \ref{thm.spechtmod.leftideal} \textbf{(a)}. Being an isomorphism,
it has an inverse $\alpha^{-1}:\mathcal{M}^{D}\rightarrow\mathbf{k}\left[
S_{n}\right]  \cdot\nabla_{\operatorname*{Row}T}$. Then: \medskip

\textbf{(a)} This $\alpha^{-1}$ is the isomorphism $\zeta$ from Proposition
\ref{prop.row-to-row.zeta}. \medskip

\textbf{(b)} We have%
\[
\alpha^{-1}\left(  \overline{S}\right)  =\nabla_{\overline{S},\overline{T}%
}\ \ \ \ \ \ \ \ \ \ \text{for each }n\text{-tabloid }\overline{S}\text{ of
shape }D.
\]

\textbf{(c)} Let $a_{\overline{S}}\in\mathbf{k}$ be a scalar for each
$n$-tabloid $\overline{S}$ of shape $D$. Then,%
\[
\alpha^{-1}\left(  \sum_{\substack{\overline{S}\text{ is an }n\text{-tabloid}%
\\\text{of shape }D}}a_{\overline{S}}\overline{S}\right)  =\sum_{u\in S_{n}%
}a_{u\overline{T}}u.
\]

\end{proposition}

\begin{proof}
Proposition \ref{prop.row-to-row.zeta} shows that the map $\zeta
:\mathcal{M}^{D}\rightarrow\mathbf{k}\left[  S_{n}\right]  \cdot
\nabla_{\operatorname*{Row}T}$ is a left $\mathbf{k}\left[  S_{n}\right]
$-module isomorphism. Hence, its inverse $\zeta^{-1}:\mathbf{k}\left[
S_{n}\right]  \cdot\nabla_{\operatorname*{Row}T}\rightarrow\mathcal{M}^{D}$ is
left $\mathbf{k}\left[  S_{n}\right]  $-linear as well. Moreover, the
definition of $\zeta$ yields $\zeta\left(  \overline{T}\right)  =\nabla
_{\overline{T},\overline{T}}=\nabla_{\operatorname*{Row}T}$ (by Proposition
\ref{prop.row-to-row.rsym}). Thus, $\zeta^{-1}\left(  \nabla
_{\operatorname*{Row}T}\right)  =\overline{T}$. In other words, $\zeta^{-1}$
sends $\nabla_{\operatorname*{Row}T}$ to $\overline{T}$.

But $\alpha$ was defined as the unique left $\mathbf{k}\left[  S_{n}\right]
$-linear map%
\[
\alpha:\mathbf{k}\left[  S_{n}\right]  \cdot\nabla_{\operatorname*{Row}%
T}\rightarrow\mathcal{M}^{D}%
\]
that sends $\nabla_{\operatorname*{Row}T}$ to $\overline{T}$. Hence, if
$\omega:\mathbf{k}\left[  S_{n}\right]  \cdot\nabla_{\operatorname*{Row}%
T}\rightarrow\mathcal{M}^{D}$ is any left $\mathbf{k}\left[  S_{n}\right]
$-linear map that sends $\nabla_{\operatorname*{Row}T}$ to $\overline{T}$,
then $\omega=\alpha$ (by the uniqueness part of the preceding sentence).
Applying this to $\omega=\zeta^{-1}$, we conclude that $\zeta^{-1}=\alpha$
(since $\zeta^{-1}:\mathbf{k}\left[  S_{n}\right]  \cdot\nabla
_{\operatorname*{Row}T}\rightarrow\mathcal{M}^{D}$ is a left $\mathbf{k}%
\left[  S_{n}\right]  $-linear map that sends $\nabla_{\operatorname*{Row}T}$
to $\overline{T}$). In other words, $\alpha$ is $\zeta^{-1}$. This proves
Proposition \ref{prop.spechtmod.leftideal.al-1} \textbf{(a)}. \medskip

\textbf{(b)} Let $\overline{S}$ be an $n$-tabloid of shape $D$. The definition
of $\zeta$ shows that $\zeta\left(  \overline{S}\right)  =\nabla_{\overline
{S},\overline{T}}$. Since $\zeta=\alpha^{-1}$ (because $\zeta^{-1}=\alpha$),
we can rewrite this as $\alpha^{-1}\left(  \overline{S}\right)  =\nabla
_{\overline{S},\overline{T}}$. This proves Proposition
\ref{prop.spechtmod.leftideal.al-1} \textbf{(b)}. \medskip

\textbf{(c)} Claim 4 in the above proof of Proposition
\ref{prop.row-to-row.zeta} shows that
\[
\zeta\left(  \sum_{\substack{\overline{S}\text{ is an }n\text{-tabloid}%
\\\text{of shape }D}}a_{\overline{S}}\overline{S}\right)  =\sum_{u\in S_{n}%
}a_{u\overline{T}}u.
\]
Since $\zeta=\alpha^{-1}$ (because $\zeta^{-1}=\alpha$), we can rewrite this
as%
\[
\alpha^{-1}\left(  \sum_{\substack{\overline{S}\text{ is an }n\text{-tabloid}%
\\\text{of shape }D}}a_{\overline{S}}\overline{S}\right)  =\sum_{u\in S_{n}%
}a_{u\overline{T}}u.
\]
This proves Proposition \ref{prop.spechtmod.leftideal.al-1} \textbf{(c)}.
\end{proof}

\subsubsection{A restatement of the left ideal avatar of the Specht module}

Let us restate part of Theorem \ref{thm.spechtmod.leftideal} in a form that
involves only the Specht module $\mathcal{S}^{D}$ (not the Young module
$\mathcal{M}^{D}$); this will be convenient to us later.

\begin{corollary}
\label{cor.spechtmod.leftideal.S}Let $D$ be any diagram with $\left\vert
D\right\vert =n$. Let $T$ be an $n$-tableau of shape $D$. Then, there is a
left $\mathbf{k}\left[  S_{n}\right]  $-module isomorphism%
\[
\overline{\alpha}:\mathbf{k}\left[  S_{n}\right]  \cdot\nabla
_{\operatorname*{Col}T}^{-}\nabla_{\operatorname*{Row}T}\rightarrow
\mathcal{S}^{D}%
\]
that sends $\nabla_{\operatorname*{Col}T}^{-}\nabla_{\operatorname*{Row}T}$ to
$\mathbf{e}_{T}$.
\end{corollary}

\begin{proof}
Theorem \ref{thm.spechtmod.leftideal} \textbf{(a)} shows that there is a
unique left $\mathbf{k}\left[  S_{n}\right]  $-linear map%
\[
\alpha:\mathbf{k}\left[  S_{n}\right]  \cdot\nabla_{\operatorname*{Row}%
T}\rightarrow\mathcal{M}^{D}%
\]
that sends $\nabla_{\operatorname*{Row}T}$ to $\overline{T}$, and that this
map is an isomorphism of $S_{n}$-representations (i.e., of left $\mathbf{k}%
\left[  S_{n}\right]  $-modules). Consider this isomorphism $\alpha$.

Theorem \ref{thm.spechtmod.leftideal} \textbf{(b)} shows that this isomorphism
$\alpha$ sends the submodule $\mathbf{k}\left[  S_{n}\right]  \cdot
\nabla_{\operatorname*{Col}T}^{-}\nabla_{\operatorname*{Row}T}$ of
$\mathbf{k}\left[  S_{n}\right]  \cdot\nabla_{\operatorname*{Row}T}$ to the
Specht module $\mathcal{S}^{D}$. Thus, by restricting $\alpha$ to the former
submodule $\mathbf{k}\left[  S_{n}\right]  \cdot\nabla_{\operatorname*{Col}%
T}^{-}\nabla_{\operatorname*{Row}T}$, we obtain an isomorphism%
\begin{align*}
\overline{\alpha}:\mathbf{k}\left[  S_{n}\right]  \cdot\nabla
_{\operatorname*{Col}T}^{-}\nabla_{\operatorname*{Row}T}  &  \rightarrow
\mathcal{S}^{D},\\
\mathbf{x}  &  \mapsto\alpha\left(  \mathbf{x}\right)  .
\end{align*}
Consider this isomorphism $\overline{\alpha}$.

The map $\alpha:\mathbf{k}\left[  S_{n}\right]  \cdot\nabla
_{\operatorname*{Row}T}\rightarrow\mathcal{M}^{D}$ is a morphism of $S_{n}%
$-representations, and thus is a morphism of left $\mathbf{k}\left[
S_{n}\right]  $-modules (since Proposition \ref{prop.rep.G-rep.mor=mor} shows
that any morphism of $S_{n}$-representations is a morphism of left
$\mathbf{k}\left[  S_{n}\right]  $-modules). Hence, its restriction
$\overline{\alpha}:\mathbf{k}\left[  S_{n}\right]  \cdot\nabla
_{\operatorname*{Col}T}^{-}\nabla_{\operatorname*{Row}T}\rightarrow
\mathcal{S}^{D}$ is a morphism of left $\mathbf{k}\left[  S_{n}\right]
$-modules as well. Since $\overline{\alpha}$ is furthermore invertible
(because it is an isomorphism of $S_{n}$-representations), we thus conclude
that $\overline{\alpha}$ is an isomorphism of left $\mathbf{k}\left[
S_{n}\right]  $-modules.

It remains to show that $\overline{\alpha}$ sends $\nabla_{\operatorname*{Col}%
T}^{-}\nabla_{\operatorname*{Row}T}$ to $\mathbf{e}_{T}$. For this purpose, we
observe that $\alpha\left(  \nabla_{\operatorname*{Row}T}\right)
=\overline{T}$ (since $\alpha$ sends $\nabla_{\operatorname*{Row}T}$ to
$\overline{T}$). But $\overline{\alpha}$ was defined as a restriction of
$\alpha$; therefore,%
\begin{align*}
&  \overline{\alpha}\left(  \nabla_{\operatorname*{Col}T}^{-}\nabla
_{\operatorname*{Row}T}\right) \\
&  =\alpha\left(  \nabla_{\operatorname*{Col}T}^{-}\nabla_{\operatorname*{Row}%
T}\right) \\
&  =\nabla_{\operatorname*{Col}T}^{-}\underbrace{\alpha\left(  \nabla
_{\operatorname*{Row}T}\right)  }_{=\overline{T}}\ \ \ \ \ \ \ \ \ \ \left(
\text{since }\alpha\text{ is a morphism of left }\mathbf{k}\left[
S_{n}\right]  \text{-modules}\right) \\
&  =\nabla_{\operatorname*{Col}T}^{-}\overline{T}=\mathbf{e}_{T}%
\ \ \ \ \ \ \ \ \ \ \left(  \text{by Proposition \ref{prop.symmetrizers.eT}%
}\right)  .
\end{align*}
In other words, $\overline{\alpha}$ sends $\nabla_{\operatorname*{Col}T}%
^{-}\nabla_{\operatorname*{Row}T}$ to $\mathbf{e}_{T}$. This completes the
proof of Corollary \ref{cor.spechtmod.leftideal.S}.
\end{proof}

\subsection{\label{sec.spechtmod.detava}Determinantal avatars of Young and
Specht modules}

There are several other avatars (i.e., isomorphic versions) of the Specht
modules $\mathcal{S}^{D}$. Two such avatars involve matrix algebra. They are
constructed by taking certain matrices (which consist of indeterminates or
polynomials) and assigning a certain product of determinants (specifically, of
minors -- i.e., determinants of submatrices) to each $n$-tableau of shape $D$.

Both avatars cover not only the Specht module $\mathcal{S}^{D}$ but also the
Young module $\mathcal{M}^{D}$, and indeed are easiest to define starting with
$\mathcal{M}^{D}$. We shall do so now.

First, however, we recall a notion of matrices that is slightly more general
than the one commonly studied in linear algebra. Namely, if $P$ and $Q$ are
two finite sets, then a $P\times Q$\emph{-matrix} means a family $\left(
a_{p,q}\right)  _{\left(  p,q\right)  \in P\times Q}$ of elements $a_{p,q}$ of
a ring indexed by all pairs $\left(  p,q\right)  \in P\times Q$. When
$P=\left[  k\right]  $ and $Q=\left[  \ell\right]  $ for some nonnegative
integers $k$ and $\ell$, this notion coincides with the classical notion of a
$k\times\ell$-matrix. A $P\times Q$-matrix is said to be \emph{square} if
$P=Q$. For a square matrix $\left(  a_{p,q}\right)  _{\left(  p,q\right)  \in
P\times P}$, we can also use the shorthand notation $\left(  a_{p,q}\right)
_{p,q\in P}$. The determinant $\det A$ of a square matrix $\left(
a_{p,q}\right)  _{\left(  p,q\right)  \in P\times P}$ (over a commutative
ring) is defined to be the sum%
\[
\sum_{\sigma\in S_{P}}\left(  -1\right)  ^{\sigma}\prod_{p\in P}%
a_{p,\sigma\left(  p\right)  };
\]
this is the same sum that would be obtained if we relabeled the elements of
$P$ as $1,2,\ldots,k$ (for $k=\left\vert P\right\vert $) and then took the
(usual) determinant of the resulting $k\times k$-matrix.

\subsubsection{The Specht--Vandermonde avatar}

Now we come to the two matrix avatars. We begin with the more classical one,
which we call the \emph{Specht--Vandermonde avatar}, since it involves
something similar to Vandermonde determinants and was introduced by Specht in
\cite{Specht35}:

\begin{theorem}
[Specht--Vandermonde avatar of Young and Specht modules]%
\label{thm.spechtmod.vdm}Let $D$ be a diagram with $\left\vert D\right\vert
=n$. Assume that all $\left(  i,j\right)  \in D$ satisfy $i>0$.

Consider the polynomial ring $\mathbf{k}\left[  x_{1},x_{2},\ldots
,x_{n}\right]  $. The symmetric group $S_{n}$ acts on this ring $\mathbf{k}%
\left[  x_{1},x_{2},\ldots,x_{n}\right]  $ by permuting the variables (that
is, $g\rightharpoonup x_{i}=x_{g\left(  i\right)  }$ for all $g\in S_{n}$ and
$i\in\left[  n\right]  $).

We define a $\mathbf{k}$-linear map%
\begin{align*}
\beta:\mathcal{M}^{D}  &  \rightarrow\mathbf{k}\left[  x_{1},x_{2}%
,\ldots,x_{n}\right]  ,\\
\overline{T}  &  \mapsto\prod_{\left(  i,j\right)  \in D}x_{T\left(
i,j\right)  }^{i-1}=\prod_{i\geq1}\ \ \prod_{m\in\operatorname*{Row}\left(
i,T\right)  }x_{m}^{i-1},
\end{align*}
where $\operatorname*{Row}\left(  i,T\right)  $ denotes the set of all entries
in the $i$-th row of $T$. Then: \medskip

\textbf{(a)} This map $\beta$ is an injective left $\mathbf{k}\left[
S_{n}\right]  $-module morphism. \medskip

\textbf{(b)} This map $\beta$ sends any polytabloid $\mathbf{e}_{T}%
\in\mathcal{M}^{D}$ to the product%
\begin{equation}
\prod_{j\in\mathbb{Z}}\det\left(  x_{T\left(  i,j\right)  }^{p-1}\right)
_{i,p\in D\left\langle j\right\rangle }, \label{eq.thm.spechtmod.vdm.b.prod}%
\end{equation}
where $D\left\langle j\right\rangle $ is the set of all integers $i$ that
satisfy $\left(  i,j\right)  \in D$. (Note that the matrix $\left(
x_{T\left(  i,j\right)  }^{p-1}\right)  _{i,p\in D\left\langle j\right\rangle
}$ is thus a square matrix, whose rows and columns are indexed by the elements
of $D\left\langle j\right\rangle $. Such square matrices have well-defined
determinants, since $D\left\langle j\right\rangle $ is a finite set. When
$D\left\langle j\right\rangle $ is empty, this matrix is a $0\times0$-matrix
and thus has determinant $1$, by definition. Hence, all but finitely many
factors of the product in (\ref{eq.thm.spechtmod.vdm.b.prod}) equal $1$.)

Thus, the map $\beta$ sends the Specht module $\mathcal{S}^{D}$ to the span of
the products shown in (\ref{eq.thm.spechtmod.vdm.b.prod}). \medskip

\textbf{(c)} Assume that $D$ is the Young diagram $Y\left(  \lambda\right)  $
of a partition $\lambda$. Let $\lambda^{t}$ be the conjugate of the partition
$\lambda$ (as defined in Theorem \ref{thm.partitions.conj}). Then, the map
$\beta$ sends any polytabloid $\mathbf{e}_{T}$ to%
\[
\prod_{j\geq1}\ \ \prod_{1\leq i_{1}<i_{2}\leq\lambda_{j}^{t}}\left(
x_{T\left(  i_{2},j\right)  }-x_{T\left(  i_{1},j\right)  }\right)  .
\]

\end{theorem}

\begin{example}
\label{exa.spechtmod.vdm.8}Let $n=8$ and%
\[
D=\left\{  \left(  1,1\right)  ,\ \left(  1,3\right)  ,\ \left(  1,4\right)
,\ \left(  2,1\right)  ,\ \left(  2,2\right)  ,\ \left(  3,1\right)
,\ \left(  3,2\right)  ,\ \left(  3,4\right)  \right\}  .
\]
Visually, $D$ looks as follows:%
\[%
%TCIMACRO{\TeXButton{tikz weird bad shape}{\begin{tikzpicture}[scale=0.7]
%\draw[fill=red!50] (0, 0) rectangle (1, 1);
%\draw[fill=red!50] (0, 1) rectangle (1, 2);
%\draw[fill=red!50] (0, 2) rectangle (1, 3);
%\draw[fill=red!50] (1, 0) rectangle (2, 1);
%\draw[fill=red!50] (1, 1) rectangle (2, 2);
%\draw[fill=red!50] (2, 2) rectangle (3, 3);
%\draw[fill=red!50] (3, 0) rectangle (4, 1);
%\draw[fill=red!50] (3, 2) rectangle (4, 3);
%\end{tikzpicture}}}%
%BeginExpansion
\begin{tikzpicture}[scale=0.7]
\draw[fill=red!50] (0, 0) rectangle (1, 1);
\draw[fill=red!50] (0, 1) rectangle (1, 2);
\draw[fill=red!50] (0, 2) rectangle (1, 3);
\draw[fill=red!50] (1, 0) rectangle (2, 1);
\draw[fill=red!50] (1, 1) rectangle (2, 2);
\draw[fill=red!50] (2, 2) rectangle (3, 3);
\draw[fill=red!50] (3, 0) rectangle (4, 1);
\draw[fill=red!50] (3, 2) rectangle (4, 3);
\end{tikzpicture}%
%EndExpansion
\ \ .
\]
Define an $n$-tableau $T$ of shape $D$ by%
\[
T=\ytableaushort{4\none13,28,76\none5}\ \ .
\]
Then, the map $\beta$ from Theorem \ref{thm.spechtmod.vdm} sends the
$n$-tabloid $\overline{T}\in\mathcal{M}^{D}$ to the monomial
\begin{align*}
\prod_{\left(  i,j\right)  \in D}x_{T\left(  i,j\right)  }^{i-1}  &
=x_{T\left(  1,1\right)  }^{0}x_{T\left(  1,3\right)  }^{0}x_{T\left(
1,4\right)  }^{0}x_{T\left(  2,1\right)  }^{1}x_{T\left(  2,2\right)  }%
^{1}x_{T\left(  3,1\right)  }^{2}x_{T\left(  3,2\right)  }^{2}x_{T\left(
3,4\right)  }^{2}\\
&  =x_{4}^{0}x_{1}^{0}x_{3}^{0}x_{2}^{1}x_{8}^{1}x_{7}^{2}x_{6}^{2}x_{5}%
^{2}=x_{2}x_{8}x_{7}^{2}x_{6}^{2}x_{5}^{2},
\end{align*}
and (by Theorem \ref{thm.spechtmod.vdm} \textbf{(b)}) sends the polytabloid
$\mathbf{e}_{T}$ to the product%
\begin{align*}
&  \prod_{j\in\mathbb{Z}}\det\left(  x_{T\left(  i,j\right)  }^{p-1}\right)
_{i,p\in D\left\langle j\right\rangle }\\
&  =\det\left(  x_{T\left(  i,1\right)  }^{p-1}\right)  _{i,p\in D\left\langle
1\right\rangle }\cdot\det\left(  x_{T\left(  i,2\right)  }^{p-1}\right)
_{i,p\in D\left\langle 2\right\rangle }\\
&  \ \ \ \ \ \ \ \ \ \ \cdot\det\left(  x_{T\left(  i,3\right)  }%
^{p-1}\right)  _{i,p\in D\left\langle 3\right\rangle }\cdot\det\left(
x_{T\left(  i,4\right)  }^{p-1}\right)  _{i,p\in D\left\langle 4\right\rangle
}\\
&  \ \ \ \ \ \ \ \ \ \ \ \ \ \ \ \ \ \ \ \ \left(
\begin{array}
[c]{c}%
\text{because for all }j\notin\left\{  1,2,3,4\right\}  \text{, the set
}D\left\langle j\right\rangle \text{ is empty,}\\
\text{and thus }\det\left(  x_{T\left(  i,j\right)  }^{p-1}\right)  _{i,p\in
D\left\langle j\right\rangle }=1
\end{array}
\right) \\
&  =\det\left(
\begin{array}
[c]{ccc}%
x_{4}^{0} & x_{4}^{1} & x_{4}^{2}\\
x_{2}^{0} & x_{2}^{1} & x_{2}^{2}\\
x_{7}^{0} & x_{7}^{1} & x_{7}^{2}%
\end{array}
\right)  \cdot\det\left(
\begin{array}
[c]{cc}%
x_{8}^{1} & x_{8}^{2}\\
x_{6}^{1} & x_{6}^{2}%
\end{array}
\right)  \cdot\det\left(
\begin{array}
[c]{c}%
x_{1}^{0}%
\end{array}
\right)  \cdot\det\left(
\begin{array}
[c]{cc}%
x_{3}^{0} & x_{3}^{2}\\
x_{5}^{0} & x_{5}^{2}%
\end{array}
\right) \\
&  \ \ \ \ \ \ \ \ \ \ \ \ \ \ \ \ \ \ \ \ \left(
\begin{array}
[c]{c}%
\text{since }D\left\langle 1\right\rangle =\left\{  1,2,3\right\}  \text{ and
}D\left\langle 2\right\rangle =\left\{  2,3\right\} \\
\text{and }D\left\langle 3\right\rangle =\left\{  1\right\}  \text{ and
}D\left\langle 4\right\rangle =\left\{  1,3\right\}
\end{array}
\right)  .
\end{align*}

\end{example}

\begin{example}
Let $n=6$ and $D=Y\left(  3,2,1\right)  $, and let $T$ be the tableau
$\ytableaushort{126,34,5}$ of shape $D$. Then, Theorem \ref{thm.spechtmod.vdm}
\textbf{(c)} says that
\[
\beta\left(  \mathbf{e}_{T}\right)  =\underbrace{\left(  x_{3}-x_{1}\right)
\left(  x_{5}-x_{1}\right)  \left(  x_{5}-x_{3}\right)  }_{\text{from column
}1}\cdot\underbrace{\left(  x_{4}-x_{2}\right)  }_{\text{from column }2}.
\]
Note that the third column of $T$ contributes nothing (or, to be more precise,
an empty product) to this product, since a product of the form $\prod_{1\leq
i_{1}<i_{2}\leq\lambda_{j}^{t}}$ is empty when $\lambda_{j}^{t}\leq1$. More
generally, any column of $Y\left(  \lambda\right)  $ that contains only one
cell (or no cells at all) contributes an empty product to the product in
Theorem \ref{thm.spechtmod.vdm} \textbf{(c)}.
\end{example}

\begin{proof}
[Proof of Theorem \ref{thm.spechtmod.vdm}.]\textbf{(a)} First, we need to show
that the map $\beta$ is well-defined. This requires proving the following two claims:

\begin{statement}
\textit{Claim 1:} If $T$ is an $n$-tableau of shape $D$, then the product
$\prod_{i\geq1}\ \ \prod_{m\in\operatorname*{Row}\left(  i,T\right)  }%
x_{m}^{i-1}$ depends only on the $n$-tabloid $\overline{T}$, not on the
$n$-tableau $T$.
\end{statement}

\begin{statement}
\textit{Claim 2:} If $T$ is an $n$-tableau of shape $D$, then $\prod_{\left(
i,j\right)  \in D}x_{T\left(  i,j\right)  }^{i-1}=\prod_{i\geq1}%
\ \ \prod_{m\in\operatorname*{Row}\left(  i,T\right)  }x_{m}^{i-1}$.
\end{statement}

\begin{proof}
[Proof of Claim 1.]Let $T$ be an $n$-tableau of shape $D$. Then, for each
$i\in\mathbb{Z}$, the set $\operatorname*{Row}\left(  i,T\right)  $ is the set
of all entries in the $i$-th row of $T$, and thus depends only on the
$n$-tabloid $\overline{T}$, not on the $n$-tableau $T$ (because row-equivalent
tableaux have the same set of entries in their $i$-th row). Therefore, the
product $\prod_{i\geq1}\ \ \prod_{m\in\operatorname*{Row}\left(  i,T\right)
}x_{m}^{i-1}$ also depends only on the $n$-tabloid $\overline{T}$, not on the
$n$-tableau $T$. This proves Claim 1.
\end{proof}

\begin{proof}
[Proof of Claim 2.]Let $T$ be an $n$-tableau of shape $D$. Then, $T$ is an
injective tableau, i.e., all the entries of $T$ are distinct. Moreover, all
$\left(  i,j\right)  \in D$ satisfy $i>0$ and thus $i\geq1$ (since $i$ is an
integer). Hence, we can rewrite the product sign $\prod_{\left(  i,j\right)
\in D}$ as $\prod_{i\geq1}\ \ \prod_{\substack{j\in\mathbb{Z};\\\left(
i,j\right)  \in D}}$. Hence,%
\begin{equation}
\prod_{\left(  i,j\right)  \in D}x_{T\left(  i,j\right)  }^{i-1}=\prod
_{i\geq1}\ \ \prod_{\substack{j\in\mathbb{Z};\\\left(  i,j\right)  \in
D}}x_{T\left(  i,j\right)  }^{i-1}. \label{pf.thm.spechtmod.vdm.a.1}%
\end{equation}
Now, let $i$ be any integer. We shall show that
\begin{equation}
\prod_{\substack{j\in\mathbb{Z};\\\left(  i,j\right)  \in D}}x_{T\left(
i,j\right)  }^{i-1}=\prod_{m\in\operatorname*{Row}\left(  i,T\right)  }%
x_{m}^{i-1}. \label{pf.thm.spechtmod.vdm.a.2}%
\end{equation}
Indeed, $\operatorname*{Row}\left(  i,T\right)  $ is the set of all entries in
the $i$-th row of $T$ (by the definition of $\operatorname*{Row}\left(
i,T\right)  $), and these entries are all distinct (since all the entries of
$T$ are distinct). Thus, the right hand side of
(\ref{pf.thm.spechtmod.vdm.a.2}) is the product of the $x_{m}^{i-1}$, where
$m$ ranges over all entries in the $i$-th row of $T$. But the left hand side
is the same product as well (since the entries in the $i$-th row of $T$ are
precisely the $T\left(  i,j\right)  $ for all integers $j$ satisfying $\left(
i,j\right)  \in D$). Thus, the two sides of (\ref{pf.thm.spechtmod.vdm.a.2})
are equal. Hence, (\ref{pf.thm.spechtmod.vdm.a.2}) is proved.

Forget that we fixed $i$. We thus have proved (\ref{pf.thm.spechtmod.vdm.a.2})
for each integer $i$. Using (\ref{pf.thm.spechtmod.vdm.a.2}), we can rewrite
(\ref{pf.thm.spechtmod.vdm.a.1}) as%
\[
\prod_{\left(  i,j\right)  \in D}x_{T\left(  i,j\right)  }^{i-1}=\prod
_{i\geq1}\ \ \prod_{m\in\operatorname*{Row}\left(  i,T\right)  }x_{m}^{i-1}.
\]
Thus, Claim 2 is proved.
\end{proof}

Claim 1 and Claim 2 yield that the map $\beta$ is well-defined. It remains to
prove the following two claims:

\begin{statement}
\textit{Claim 3:} The map $\beta$ is injective.
\end{statement}

\begin{statement}
\textit{Claim 4:} The map $\beta$ is left $\mathbf{k}\left[  S_{n}\right]  $-linear.
\end{statement}

\begin{proof}
[Proof of Claim 3.]Let $\mathfrak{M}$ denote the set of all monomials in the
$n$ indeterminates $x_{1},x_{2},\ldots,x_{n}$. Then, the $\mathbf{k}$-module
$\mathbf{k}\left[  x_{1},x_{2},\ldots,x_{n}\right]  $ is the free $\mathbf{k}%
$-module $\mathbf{k}^{\left(  \mathfrak{M}\right)  }$. Meanwhile,
$\mathcal{M}^{D}$ is the free $\mathbf{k}$-module $\mathbf{k}^{\left(
\left\{  n\text{-tabloids of shape }D\right\}  \right)  }$ (by its
definition). The $\mathbf{k}$-linear map $\beta:\mathcal{M}^{D}\rightarrow
\mathbf{k}\left[  x_{1},x_{2},\ldots,x_{n}\right]  $ sends each standard basis
vector $\overline{T}$ of $\mathcal{M}^{D}$ to the standard basis vector
$\prod_{\left(  i,j\right)  \in D}x_{T\left(  i,j\right)  }^{i-1}$ of
$\mathbf{k}\left[  x_{1},x_{2},\ldots,x_{n}\right]  $; thus, it is the
linearization\footnote{See Definition \ref{def.linearization.gen} for the
definition of \textquotedblleft linearization\textquotedblright.} $f_{\ast}$
of the map%
\begin{align*}
f:\left\{  n\text{-tabloids of shape }D\right\}   &  \rightarrow
\mathfrak{M},\\
\overline{T}  &  \mapsto\prod_{\left(  i,j\right)  \in D}x_{T\left(
i,j\right)  }^{i-1}.
\end{align*}
Consider this latter map $f$. If we can show that $f$ is injective, then it
will follow (by Proposition \ref{prop.linearization.ject}) that its
linearization $f_{\ast}$ is injective as well, and this will mean that $\beta$
is injective (since $\beta$ is the linearization $f_{\ast}$); thus, Claim 3
will be proven. Hence, it suffices to show that $f$ is injective.

So let us prove this. Let $\overline{T}$ and $\overline{S}$ be two
$n$-tabloids of shape $D$ such that $f\left(  \overline{T}\right)  =f\left(
\overline{S}\right)  $. We must show that $\overline{T}=\overline{S}$.

We know that $T$ is an $n$-tableau, i.e., an injective map from $D$ to
$\left\{  1,2,3,\ldots\right\}  $ whose image is $\left[  n\right]  $. Thus,
we can view $T$ as a bijection from $D$ to $\left[  n\right]  $. Similarly, we
can view $S$ as a bijection from $D$ to $\left[  n\right]  $. This is how we
shall view $T$ and $S$ henceforth.

For each $m\in\left[  n\right]  $, let us denote the cell $T^{-1}\left(
m\right)  \in D$ (that is, the cell of $T$ that contains $m$) by $\left(
a_{m},b_{m}\right)  $, and let us denote the cell $S^{-1}\left(  m\right)  \in
D$ (that is, the cell of $S$ that contains $m$) by $\left(  c_{m}%
,d_{m}\right)  $.

Now, the definition of $f$ easily yields that%
\begin{equation}
f\left(  \overline{T}\right)  =\prod_{m=1}^{n}x_{m}^{a_{m}-1}
\label{pf.thm.spechtmod.vdm.a.4}%
\end{equation}
\footnote{\textit{Proof.} Let $\left(  i,j\right)  \in D$. For each
$m\in\left[  n\right]  $, we have $\left(  a_{m},b_{m}\right)  =T^{-1}\left(
m\right)  $ (by the definition of $\left(  a_{m},b_{m}\right)  $). Applying
this to $m=T\left(  i,j\right)  $, we obtain $\left(  a_{T\left(  i,j\right)
},b_{T\left(  i,j\right)  }\right)  =T^{-1}\left(  T\left(  i,j\right)
\right)  =\left(  i,j\right)  $. In other words, $a_{T\left(  i,j\right)  }=i$
and $b_{T\left(  i,j\right)  }=j$.
\par
Forget that we fixed $\left(  i,j\right)  $. We thus have shown that all
$\left(  i,j\right)  \in D$ satisfy
\begin{equation}
a_{T\left(  i,j\right)  }=i\ \ \ \ \ \ \ \ \ \ \text{and}%
\ \ \ \ \ \ \ \ \ \ b_{T\left(  i,j\right)  }=j.
\label{pf.thm.spechtmod.vdm.a.c3.pf.fn1.1}%
\end{equation}
\par
But the definition of $f$ yields%
\[
f\left(  \overline{T}\right)  =\prod_{\left(  i,j\right)  \in D}%
\underbrace{x_{T\left(  i,j\right)  }^{i-1}}_{\substack{=x_{T\left(
i,j\right)  }^{a_{T\left(  i,j\right)  }-1}\\\text{(since
(\ref{pf.thm.spechtmod.vdm.a.c3.pf.fn1.1}) yields }a_{T\left(  i,j\right)
}=i\\\text{and thus }i=a_{T\left(  i,j\right)  }\text{)}}}=\prod_{\left(
i,j\right)  \in D}x_{T\left(  i,j\right)  }^{a_{T\left(  i,j\right)  }%
-1}=\prod_{m\in\left[  n\right]  }x_{m}^{a_{m}-1}%
\]
(here, we have substituted $m$ for $T\left(  i,j\right)  $ in the product,
since the map $T:D\rightarrow\left[  n\right]  $ is a bijection). Qed.}.
Similarly,%
\[
f\left(  \overline{S}\right)  =\prod_{m=1}^{n}x_{m}^{c_{m}-1}.
\]
The left hand sides of these two equalities are equal (since $f\left(
\overline{T}\right)  =f\left(  \overline{S}\right)  $). Thus, so are their
right hand sides. In other words, $\prod_{m=1}^{n}x_{m}^{a_{m}-1}=\prod
_{m=1}^{n}x_{m}^{c_{m}-1}$. Therefore,%
\[
a_{m}-1=c_{m}-1\ \ \ \ \ \ \ \ \ \ \text{for each }m\in\left[  n\right]
\]
(since equal monomials must have equal powers on each indeterminate). In other
words,%
\begin{equation}
a_{m}=c_{m}\ \ \ \ \ \ \ \ \ \ \text{for each }m\in\left[  n\right]  .
\label{pf.thm.spechtmod.vdm.a.6}%
\end{equation}

Now, fix any $i\in\mathbb{Z}$. Note that the entries in the $i$-th row of $T$
are distinct (since $T$ is injective). Likewise, so are the entries in the
$i$-th row of $S$.

For any $m\in\left[  n\right]  $, we have the logical equivalences%
\begin{align}
&  \ \left(  m\in\operatorname*{Row}\left(  i,T\right)  \right) \nonumber\\
&  \Longleftrightarrow\ \left(  m\text{ appears in the }i\text{-th row of
}T\right) \nonumber\\
&  \ \ \ \ \ \ \ \ \ \ \ \ \ \ \ \ \ \ \ \ \left(
\begin{array}
[c]{c}%
\text{since }\operatorname*{Row}\left(  i,T\right)  \text{ is defined to be
the}\\
\text{set of all entries in the }i\text{-th row of }T
\end{array}
\right) \nonumber\\
&  \Longleftrightarrow\ \left(  \text{the cell of }T\text{ that contains
}m\text{ lies in the }i\text{-th row}\right) \nonumber\\
&  \Longleftrightarrow\ \left(  \text{the cell }\left(  a_{m},b_{m}\right)
\text{ lies in the }i\text{-th row}\right) \nonumber\\
&  \ \ \ \ \ \ \ \ \ \ \ \ \ \ \ \ \ \ \ \ \left(  \text{since the cell of
}T\text{ that contains }m\text{ is }T^{-1}\left(  m\right)  =\left(
a_{m},b_{m}\right)  \right) \nonumber\\
&  \Longleftrightarrow\ \left(  a_{m}=i\right)
\label{pf.thm.spechtmod.vdm.a.7a}%
\end{align}
and similarly the equivalence%
\begin{equation}
\left(  m\in\operatorname*{Row}\left(  i,S\right)  \right)
\ \Longleftrightarrow\ \left(  c_{m}=i\right)  .
\label{pf.thm.spechtmod.vdm.a.7c}%
\end{equation}
Thus, for any $m\in\left[  n\right]  $, we have the logical equivalences%
\begin{align*}
\left(  m\in\operatorname*{Row}\left(  i,T\right)  \right)  \  &
\Longleftrightarrow\ \left(  a_{m}=i\right)  \ \ \ \ \ \ \ \ \ \ \left(
\text{by (\ref{pf.thm.spechtmod.vdm.a.7a})}\right) \\
&  \Longleftrightarrow\ \left(  c_{m}=i\right)  \ \ \ \ \ \ \ \ \ \ \left(
\text{since (\ref{pf.thm.spechtmod.vdm.a.6}) yields }a_{m}=c_{m}\right) \\
&  \Longleftrightarrow\ \left(  m\in\operatorname*{Row}\left(  i,S\right)
\right)  \ \ \ \ \ \ \ \ \ \ \left(  \text{by (\ref{pf.thm.spechtmod.vdm.a.7c}%
)}\right)  .
\end{align*}
Hence, $\operatorname*{Row}\left(  i,T\right)  =\operatorname*{Row}\left(
i,S\right)  $ (since $\operatorname*{Row}\left(  i,T\right)  $ and
$\operatorname*{Row}\left(  i,S\right)  $ are subsets of $\left[  n\right]
$). In other words, the set of entries in the $i$-th row of $T$ equals the set
of entries in the $i$-th row of $S$ (since these two sets are precisely
$\operatorname*{Row}\left(  i,T\right)  $ and $\operatorname*{Row}\left(
i,S\right)  $). We can furthermore replace the word \textquotedblleft
set\textquotedblright\ by \textquotedblleft multiset\textquotedblright\ here
(since the entries in the $i$-th row of $T$ are distinct, and so are the
entries in the $i$-th row of $S$, and therefore the multisets of these entries
do not have any nontrivial multiplicities). Thus, we conclude that the
multiset of entries in the $i$-th row of $T$ equals the multiset of entries in
the $i$-th row of $S$.

Forget that we fixed $i$. We thus have shown that for each $i\in\mathbb{Z}$,
the multiset of entries in the $i$-th row of $T$ equals the multiset of
entries in the $i$-th row of $S$. In other words, the $n$-tableaux $T$ and $S$
are row-equivalent. In other words, $\overline{T}=\overline{S}$ (since an
$n$-tabloid is an equivalence class with respect to row-equivalence).

Forget that we fixed $\overline{T}$ and $\overline{S}$. We thus have shown
that if $\overline{T}$ and $\overline{S}$ are two $n$-tabloids of shape $D$
such that $f\left(  \overline{T}\right)  =f\left(  \overline{S}\right)  $,
then $\overline{T}=\overline{S}$. In other words, $f$ is injective. As
explained above, this completes the proof of Claim 3.
\end{proof}

\begin{proof}
[Proof of Claim 4.]We shall first show that $\beta$ is $S_{n}$-equivariant.

For this purpose, we must prove that $\beta\left(  w\mathbf{a}\right)
=w\beta\left(  \mathbf{a}\right)  $ for all $w\in S_{n}$ and all
$\mathbf{a}\in\mathcal{M}^{D}$.

Fix $w\in S_{n}$ and $\mathbf{a}\in\mathcal{M}^{D}$. We must prove that
$\beta\left(  w\mathbf{a}\right)  =w\beta\left(  \mathbf{a}\right)  $. Both
sides of this equality are $\mathbf{k}$-linear in $\mathbf{a}$. Thus, by
linearity, we can WLOG assume that $\mathbf{a}$ is an $n$-tabloid of shape $D$
(since $\mathcal{M}^{D}$ is spanned by such $n$-tabloids). Assume this. Thus,
$\mathbf{a}=\overline{T}$ for some $n$-tabloid $\overline{T}$ of shape $D$.
Consider this $\overline{T}$. From $\mathbf{a}=\overline{T}$, we obtain
$\beta\left(  \mathbf{a}\right)  =\beta\left(  \overline{T}\right)
=\prod_{\left(  i,j\right)  \in D}x_{T\left(  i,j\right)  }^{i-1}$ (by the
definition of $\beta$), so that%
\begin{equation}
w\beta\left(  \mathbf{a}\right)  =w\prod_{\left(  i,j\right)  \in
D}x_{T\left(  i,j\right)  }^{i-1}=\prod_{\left(  i,j\right)  \in D}x_{w\left(
T\left(  i,j\right)  \right)  }^{i-1} \label{pf.thm.spechtmod.vdm.a.10}%
\end{equation}
(since the action of $w$ transforms each indeterminate $x_{k}$ into
$x_{w\left(  k\right)  }$). On the other hand, from $\mathbf{a}=\overline{T}$,
we obtain $w\mathbf{a}=w\overline{T}=\overline{wT}$, so that%
\begin{align*}
\beta\left(  w\mathbf{a}\right)   &  =\beta\left(  \overline{wT}\right)
=\prod_{\left(  i,j\right)  \in D}\underbrace{x_{\left(  wT\right)  \left(
i,j\right)  }^{i-1}}_{\substack{=x_{w\left(  T\left(  i,j\right)  \right)
}^{i-1}\\\text{(since }wT=w\circ T\text{ and}\\\text{thus }\left(  wT\right)
\left(  i,j\right)  =w\left(  T\left(  i,j\right)  \right)  \text{)}%
}}\ \ \ \ \ \ \ \ \ \ \left(  \text{by the definition of }\beta\right) \\
&  =\prod_{\left(  i,j\right)  \in D}x_{w\left(  T\left(  i,j\right)  \right)
}^{i-1}.
\end{align*}
Comparing this with (\ref{pf.thm.spechtmod.vdm.a.10}), we find $\beta\left(
w\mathbf{a}\right)  =w\beta\left(  \mathbf{a}\right)  $.

Forget that we fixed $w$ and $\mathbf{a}$. We thus have shown that
$\beta\left(  w\mathbf{a}\right)  =w\beta\left(  \mathbf{a}\right)  $ for all
$w\in S_{n}$ and $\mathbf{a}\in\mathcal{M}^{D}$. In other words, the map
$\beta$ is $S_{n}$-equivariant. Since $\beta$ is also $\mathbf{k}$-linear, we
thus conclude that $\beta$ is a morphism of $S_{n}$-representations. Therefore
(by Proposition \ref{prop.rep.G-rep.mor=mor}), the map $\beta$ is a morphism
of left $\mathbf{k}\left[  S_{n}\right]  $-modules, i.e., a left
$\mathbf{k}\left[  S_{n}\right]  $-linear map. This proves Claim 4.
\end{proof}

Claim 3 and Claim 4 show that $\beta$ is an injective left $\mathbf{k}\left[
S_{n}\right]  $-module morphism. Theorem \ref{thm.spechtmod.vdm} \textbf{(a)}
is thus proved. \medskip

\textbf{(b)} Let $T$ be an $n$-tableau of shape $D$. Then, the definition of
$\mathbf{e}_{T}$ yields%
\[
\mathbf{e}_{T}=\sum_{w\in\mathcal{C}\left(  T\right)  }\left(  -1\right)
^{w}\overline{wT}.
\]
Hence,%
\begin{align}
\beta\left(  \mathbf{e}_{T}\right)   &  =\beta\left(  \sum_{w\in
\mathcal{C}\left(  T\right)  }\left(  -1\right)  ^{w}\overline{wT}\right)
\nonumber\\
&  =\sum_{w\in\mathcal{C}\left(  T\right)  }\left(  -1\right)  ^{w}%
\beta\left(  \overline{wT}\right)  \ \ \ \ \ \ \ \ \ \ \left(  \text{since the
map }\beta\text{ is }\mathbf{k}\text{-linear}\right) \nonumber\\
&  =\sum_{w\in\mathcal{C}\left(  T\right)  }\left(  -1\right)  ^{w}%
\prod_{\left(  i,j\right)  \in D}x_{\left(  wT\right)  \left(  i,j\right)
}^{i-1} \label{pf.thm.spechtmod.vdm.b.1}%
\end{align}
(since the definition of $\beta$ yields $\beta\left(  \overline{wT}\right)
=\prod_{\left(  i,j\right)  \in D}x_{\left(  wT\right)  \left(  i,j\right)
}^{i-1}$ for every $w\in\mathcal{C}\left(  T\right)  $).

Now, let us WLOG assume that all cells $\left(  i,j\right)  \in D$ satisfy
$j>0$. (We can always achieve this by a horizontal translation of $D$ (since
$D$ is finite); this does not affect the validity of the theorem we are
proving.) Since $D$ is finite, there exists some $k\in\mathbb{N}$ such that
all cells of $D$ lie in columns $1,2,\ldots,k$. Pick such a $k$. Thus, each
cell $\left(  i,j\right)  \in D$ satisfies $j\in\left[  k\right]  $.
Therefore, we can rewrite the product sign $\prod_{\left(  i,j\right)  \in D}$
as $\prod_{j=1}^{k}\ \ \prod_{\substack{i\in\mathbb{Z};\\\left(  i,j\right)
\in D}}$. Furthermore, for any given $j\in\left[  k\right]  $, we can rewrite
the product sign $\prod_{\substack{i\in\mathbb{Z};\\\left(  i,j\right)  \in
D}}$ as $\prod_{i\in D\left\langle j\right\rangle }$ (since $D\left\langle
j\right\rangle $ is the set of all integers $i$ such that $\left(  i,j\right)
\in D$). Altogether, we thus find the following equality of product signs:%
\begin{equation}
\prod_{\left(  i,j\right)  \in D}=\prod_{j=1}^{k}\ \ \underbrace{\prod
_{\substack{i\in\mathbb{Z};\\\left(  i,j\right)  \in D}}}_{=\prod_{i\in
D\left\langle j\right\rangle }}=\prod_{j=1}^{k}\ \ \prod_{i\in D\left\langle
j\right\rangle }. \label{pf.thm.spechtmod.vdm.b.1.prodsigns}%
\end{equation}

Furthermore, for each integer $j\notin\left[  k\right]  $, we have
$D\left\langle j\right\rangle =\varnothing$ (since all cells of $D$ lie in
columns $1,2,\ldots,k$, and thus no cell of $D$ lies in column $j$), and
therefore%
\[
\det\underbrace{\left(  x_{T\left(  i,j\right)  }^{p-1}\right)  _{i,p\in
D\left\langle j\right\rangle }}_{\substack{\text{a }0\times0\text{-matrix}%
\\\text{(since }D\left\langle j\right\rangle =\varnothing\text{)}}}=1.
\]
In other words, in the infinite product $\prod_{j\in\mathbb{Z}}\det\left(
x_{T\left(  i,j\right)  }^{p-1}\right)  _{i,p\in D\left\langle j\right\rangle
}$, all factors with $j\notin\left[  k\right]  $ equal $1$. Hence, this
infinite product is well-defined, and equals the same product but restricted
only to the factors with $j\in\left[  k\right]  $. In other words,%
\begin{equation}
\prod_{j\in\mathbb{Z}}\det\left(  x_{T\left(  i,j\right)  }^{p-1}\right)
_{i,p\in D\left\langle j\right\rangle }=\prod_{j=1}^{k}\det\left(  x_{T\left(
i,j\right)  }^{p-1}\right)  _{i,p\in D\left\langle j\right\rangle }.
\label{pf.thm.spechtmod.vdm.b.1.get-rid}%
\end{equation}

For any $j\in\mathbb{Z}$, let $\operatorname*{Col}\left(  j,T\right)  $ denote
the set of all entries in the $j$-th column of $T$. Let us set $X_{j}%
:=\operatorname*{Col}\left(  j,T\right)  $ for each $j\in\left[  k\right]  $.
As in the proof of Proposition \ref{prop.symmetrizers.int} \textbf{(b)}, we
can see that $\operatorname*{Col}\left(  1,T\right)  ,\ \operatorname*{Col}%
\left(  2,T\right)  ,\ \ldots,\ \operatorname*{Col}\left(  k,T\right)  $ are
$k$ disjoint subsets of $\left[  n\right]  $ and satisfy $\operatorname*{Col}%
\left(  1,T\right)  \cup\operatorname*{Col}\left(  2,T\right)  \cup\cdots
\cup\operatorname*{Col}\left(  k,T\right)  =\left[  n\right]  $. In other
words, $X_{1},X_{2},\ldots,X_{k}$ are $k$ disjoint subsets of $\left[
n\right]  $ and satisfy $X_{1}\cup X_{2}\cup\cdots\cup X_{k}=\left[  n\right]
$ (since $X_{j}=\operatorname*{Col}\left(  j,T\right)  $ for each $j\in\left[
k\right]  $).

Define a subset $K$ of $S_{n}$ by%
\[
K:=\left\{  w\in S_{n}\ \mid\ w\left(  X_{i}\right)  \subseteq X_{i}\text{ for
all }i\in\left[  k\right]  \right\}  .
\]
Then,%
\begin{align}
K  &  =\left\{  w\in S_{n}\ \mid\ w\left(  X_{i}\right)  \subseteq X_{i}\text{
for all }i\in\left[  k\right]  \right\} \nonumber\\
&  =\left\{  w\in S_{n}\ \mid\ w\left(  \operatorname*{Col}\left(  i,T\right)
\right)  \subseteq\operatorname*{Col}\left(  i,T\right)  \text{ for all }%
i\in\left[  k\right]  \right\}  \label{pf.thm.spechtmod.vdm.b.3}%
\end{align}
(since $X_{i}=\operatorname*{Col}\left(  i,T\right)  $ for each $i\in\left[
k\right]  $).

However, for any permutation $w\in S_{n}$, we have the logical equivalence%
\begin{equation}
\left(  w\in\mathcal{C}\left(  T\right)  \right)  \ \Longleftrightarrow
\ \left(  w\left(  \operatorname*{Col}\left(  i,T\right)  \right)
\subseteq\operatorname*{Col}\left(  i,T\right)  \text{ for all }i\in\left[
k\right]  \right)  \label{pf.thm.spechtmod.vdm.b.4}%
\end{equation}
(this was already proved back in the proof of Proposition
\ref{prop.symmetrizers.int} \textbf{(b)}). From $\mathcal{C}\left(  T\right)
\subseteq S_{n}$, we obtain%
\begin{align*}
\mathcal{C}\left(  T\right)   &  =\left\{  w\in S_{n}\ \mid\ w\in
\mathcal{C}\left(  T\right)  \right\} \\
&  =\left\{  w\in S_{n}\ \mid\ w\left(  \operatorname*{Col}\left(  i,T\right)
\right)  \subseteq\operatorname*{Col}\left(  i,T\right)  \text{ for all }%
i\in\left[  k\right]  \right\} \\
&  \ \ \ \ \ \ \ \ \ \ \ \ \ \ \ \ \ \ \ \ \left(  \text{by the equivalence
(\ref{pf.thm.spechtmod.vdm.b.4})}\right) \\
&  =K\ \ \ \ \ \ \ \ \ \ \left(  \text{by (\ref{pf.thm.spechtmod.vdm.b.3}%
)}\right)  .
\end{align*}
Thus, we can rewrite (\ref{pf.thm.spechtmod.vdm.b.1}) as
\begin{equation}
\beta\left(  \mathbf{e}_{T}\right)  =\sum_{w\in K}\left(  -1\right)  ^{w}%
\prod_{\left(  i,j\right)  \in D}x_{\left(  wT\right)  \left(  i,j\right)
}^{i-1}. \label{pf.thm.spechtmod.vdm.b.5}%
\end{equation}

For each permutation $w\in K$ and each $j\in\left[  k\right]  $, we define a
map $w^{\left\langle j\right\rangle }:X_{j}\rightarrow X_{j}$ by setting%
\[
w^{\left\langle j\right\rangle }\left(  x\right)  =w\left(  x\right)
\ \ \ \ \ \ \ \ \ \ \text{for each }x\in X_{j}.
\]
As we know from Lemma \ref{lem.Xk.bij1} \textbf{(a)}, this map
$w^{\left\langle j\right\rangle }$ is well-defined and belongs to the
symmetric group $S_{X_{j}}$ for each $w\in K$ and each $j\in\left[  k\right]
$. Furthermore, Lemma \ref{lem.Xk.bij1} \textbf{(b)} tells us that the map%
\begin{align}
K  &  \rightarrow S_{X_{1}}\times S_{X_{2}}\times\cdots\times S_{X_{k}%
},\nonumber\\
w  &  \mapsto\left(  w^{\left\langle 1\right\rangle },w^{\left\langle
2\right\rangle },\ldots,w^{\left\langle k\right\rangle }\right)
\label{pf.thm.spechtmod.vdm.b.bij}%
\end{align}
is a bijection. It is easy to see that each $w\in K$ and each $\left(
i,j\right)  \in D$ satisfy%
\begin{equation}
\left(  wT\right)  \left(  i,j\right)  =w^{\left\langle j\right\rangle
}\left(  T\left(  i,j\right)  \right)  \label{pf.thm.spechtmod.vdm.b.6}%
\end{equation}
\footnote{\textit{Proof:} Let $w\in K$ and $\left(  i,j\right)  \in D$. Then,
$T\left(  i,j\right)  $ is an entry in the $j$-th column of $T$, and thus
belongs to $\operatorname*{Col}\left(  j,T\right)  $ (by the definition of
$\operatorname*{Col}\left(  j,T\right)  $). In other words, $T\left(
i,j\right)  \in\operatorname*{Col}\left(  j,T\right)  =X_{j}$ (since $X_{j}$
is defined to be $\operatorname*{Col}\left(  j,T\right)  $). Hence, the
definition of $w^{\left\langle j\right\rangle }$ yields $w^{\left\langle
j\right\rangle }\left(  T\left(  i,j\right)  \right)  =w\left(  T\left(
i,j\right)  \right)  $. But the definition of the action of $S_{n}$ on
$n$-tableaux shows that $wT=w\circ T$, so that $\left(  wT\right)  \left(
i,j\right)  =\left(  w\circ T\right)  \left(  i,j\right)  =w\left(  T\left(
i,j\right)  \right)  $. Comparing this with $w^{\left\langle j\right\rangle
}\left(  T\left(  i,j\right)  \right)  =w\left(  T\left(  i,j\right)  \right)
$, we obtain $\left(  wT\right)  \left(  i,j\right)  =w^{\left\langle
j\right\rangle }\left(  T\left(  i,j\right)  \right)  $. This proves
(\ref{pf.thm.spechtmod.vdm.b.6}).}.

Finally, each $w\in K$ satisfies%
\begin{align}
\left(  -1\right)  ^{w}  &  =\left(  -1\right)  ^{w^{\left\langle
1\right\rangle }}\left(  -1\right)  ^{w^{\left\langle 2\right\rangle }}%
\cdots\left(  -1\right)  ^{w^{\left\langle k\right\rangle }}%
\ \ \ \ \ \ \ \ \ \ \left(  \text{by Lemma \ref{lem.Xk.bij1} \textbf{(c)}%
}\right) \nonumber\\
&  =\prod_{j=1}^{k}\left(  -1\right)  ^{w^{\left\langle j\right\rangle }}.
\label{pf.thm.spechtmod.vdm.b.7}%
\end{align}

Hence, (\ref{pf.thm.spechtmod.vdm.b.5}) becomes%
\begin{align}
\beta\left(  \mathbf{e}_{T}\right)   &  =\sum_{w\in K}\underbrace{\left(
-1\right)  ^{w}}_{\substack{=\prod_{j=1}^{k}\left(  -1\right)
^{w^{\left\langle j\right\rangle }}\\\text{(by (\ref{pf.thm.spechtmod.vdm.b.7}%
))}}}\ \ \underbrace{\prod_{\left(  i,j\right)  \in D}}_{\substack{=\prod
_{j=1}^{k}\ \ \prod_{i\in D\left\langle j\right\rangle }\\\text{(by
(\ref{pf.thm.spechtmod.vdm.b.1.prodsigns}))}}}\ \ \underbrace{x_{\left(
wT\right)  \left(  i,j\right)  }^{i-1}}_{\substack{=x_{w^{\left\langle
j\right\rangle }\left(  T\left(  i,j\right)  \right)  }^{i-1}\\\text{(since
(\ref{pf.thm.spechtmod.vdm.b.6}) shows}\\\text{that }\left(  wT\right)
\left(  i,j\right)  =w^{\left\langle j\right\rangle }\left(  T\left(
i,j\right)  \right)  \text{)}}}\nonumber\\
&  =\sum_{w\in K}\underbrace{\left(  \prod_{j=1}^{k}\left(  -1\right)
^{w^{\left\langle j\right\rangle }}\right)  \prod_{j=1}^{k}\ \ \prod_{i\in
D\left\langle j\right\rangle }x_{w^{\left\langle j\right\rangle }\left(
T\left(  i,j\right)  \right)  }^{i-1}}_{=\prod_{j=1}^{k}\left(  \left(
-1\right)  ^{w^{\left\langle j\right\rangle }}\prod_{i\in D\left\langle
j\right\rangle }x_{w^{\left\langle j\right\rangle }\left(  T\left(
i,j\right)  \right)  }^{i-1}\right)  }\nonumber\\
&  =\sum_{w\in K}\ \ \prod_{j=1}^{k}\left(  \left(  -1\right)
^{w^{\left\langle j\right\rangle }}\prod_{i\in D\left\langle j\right\rangle
}x_{w^{\left\langle j\right\rangle }\left(  T\left(  i,j\right)  \right)
}^{i-1}\right) \nonumber\\
&  =\sum_{\left(  w_{1},w_{2},\ldots,w_{k}\right)  \in S_{X_{1}}\times
S_{X_{2}}\times\cdots\times S_{X_{k}}}\ \ \prod_{j=1}^{k}\left(  \left(
-1\right)  ^{w_{j}}\prod_{i\in D\left\langle j\right\rangle }x_{w_{j}\left(
T\left(  i,j\right)  \right)  }^{i-1}\right) \nonumber\\
&  \ \ \ \ \ \ \ \ \ \ \ \ \ \ \ \ \ \ \ \ \left(
\begin{array}
[c]{c}%
\text{here, we have substituted }\left(  w_{1},w_{2},\ldots,w_{k}\right) \\
\text{for }\left(  w^{\left\langle 1\right\rangle },w^{\left\langle
2\right\rangle },\ldots,w^{\left\langle k\right\rangle }\right)  \text{ in the
sum, since the}\\
\text{map (\ref{pf.thm.spechtmod.vdm.b.bij}) is a bijection}%
\end{array}
\right) \nonumber\\
&  =\prod_{j=1}^{k}\ \ \sum_{w\in S_{X_{j}}}\left(  -1\right)  ^{w}\prod_{i\in
D\left\langle j\right\rangle }x_{w\left(  T\left(  i,j\right)  \right)
}^{i-1} \label{pf.thm.spechtmod.vdm.b.as-prod}%
\end{align}
(by the product rule, applied backwards\footnote{i.e., since the product rule
shows that
\begin{align*}
&  \prod_{j=1}^{k}\ \ \sum_{w\in S_{X_{j}}}\left(  -1\right)  ^{w}\prod_{i\in
D\left\langle j\right\rangle }x_{w\left(  T\left(  i,j\right)  \right)
}^{i-1}\\
&  =\sum_{\left(  w_{1},w_{2},\ldots,w_{k}\right)  \in S_{X_{1}}\times
S_{X_{2}}\times\cdots\times S_{X_{k}}}\ \ \prod_{j=1}^{k}\left(  \left(
-1\right)  ^{w_{j}}\prod_{i\in D\left\langle j\right\rangle }x_{w_{j}\left(
T\left(  i,j\right)  \right)  }^{i-1}\right)
\end{align*}
}).

Now, we shall prove the following:

\begin{statement}
\textit{Claim 5.} Let $j\in\left[  k\right]  $. Then,%
\begin{equation}
\det\left(  x_{T\left(  i,j\right)  }^{p-1}\right)  _{i,p\in D\left\langle
j\right\rangle }=\sum_{w\in S_{X_{j}}}\left(  -1\right)  ^{w}\prod_{i\in
D\left\langle j\right\rangle }x_{w\left(  T\left(  i,j\right)  \right)
}^{i-1}. \label{pf.thm.spechtmod.vdm.b.det}%
\end{equation}

\end{statement}

\begin{proof}
[Proof of Claim 5.]We recall a basic property of determinants: The determinant
of a square matrix does not change if we relabel its rows and its columns (as
long as both rows and columns are relabelled simultaneously using the same
bijection). In other words, if $\alpha:P\rightarrow Q$ is a bijection between
two finite sets $P$ and $Q$, and if $\left(  a_{p,q}\right)  _{p,q\in Q}$ is a
$Q\times Q$-matrix, then%
\begin{equation}
\det\left(  a_{p,q}\right)  _{p,q\in Q}=\det\left(  a_{\alpha\left(  p\right)
,\alpha\left(  q\right)  }\right)  _{p,q\in P}.
\label{pf.thm.spechtmod.vdm.b.relabel}%
\end{equation}

Next, we recall that $T$ is an $n$-tableau, hence injective. In other words,
all entries $T\left(  i,j\right)  $ of $T$ are distinct.

Now, recall that $D\left\langle j\right\rangle $ is the set of all integers
$i$ such that $\left(  i,j\right)  \in D$. Thus, for each $i\in D\left\langle
j\right\rangle $, the entry $T\left(  i,j\right)  $ of $T$ is a well-defined
entry of the $j$-th column of $T$. In other words, for each $i\in
D\left\langle j\right\rangle $, we have $T\left(  i,j\right)  \in
\operatorname*{Col}\left(  j,T\right)  $ (since $\operatorname*{Col}\left(
j,T\right)  $ is the set of all entries of the $j$-th column of $T$). In other
words, for each $i\in D\left\langle j\right\rangle $, we have $T\left(
i,j\right)  \in X_{j}$ (since $X_{j}=\operatorname*{Col}\left(  j,T\right)
$). Hence, the map%
\begin{align*}
\alpha:D\left\langle j\right\rangle  &  \rightarrow X_{j},\\
i  &  \mapsto T\left(  i,j\right)
\end{align*}
is well-defined. This map $\alpha$ is furthermore injective (since all entries
$T\left(  i,j\right)  $ of $T$ are distinct) and surjective (since each $x\in
X_{j}$ is an entry of the $j$-th column of $T$ (because $x\in X_{j}%
=\operatorname*{Col}\left(  j,T\right)  $) and thus has the form $T\left(
i,j\right)  $ for some $i\in D\left\langle j\right\rangle $). Thus, it is
bijective. In other words, $\alpha:D\left\langle j\right\rangle \rightarrow
X_{j}$ is a bijection.

Hence, (\ref{pf.thm.spechtmod.vdm.b.relabel}) (applied to $P=D\left\langle
j\right\rangle $ and $Q=X_{j}$ and $a_{p,q}=x_{p}^{\alpha^{-1}\left(
q\right)  -1}$) yields%
\begin{align*}
&  \det\left(  x_{p}^{\alpha^{-1}\left(  q\right)  -1}\right)  _{p,q\in X_{j}%
}\\
&  =\det\left(  x_{\alpha\left(  p\right)  }^{\alpha^{-1}\left(  \alpha\left(
q\right)  \right)  -1}\right)  _{p,q\in D\left\langle j\right\rangle }\\
&  =\det\left(  x_{\alpha\left(  p\right)  }^{q-1}\right)  _{p,q\in
D\left\langle j\right\rangle }\ \ \ \ \ \ \ \ \ \ \left(  \text{since }%
\alpha^{-1}\left(  \alpha\left(  q\right)  \right)  =q\text{ for each }q\in
D\left\langle j\right\rangle \right) \\
&  =\det\left(  x_{\alpha\left(  i\right)  }^{p-1}\right)  _{i,p\in
D\left\langle j\right\rangle }\ \ \ \ \ \ \ \ \ \ \left(
\begin{array}
[c]{c}%
\text{here, we have renamed the indices }p\text{ and }q\\
\text{as }i\text{ and }p
\end{array}
\right) \\
&  =\det\left(  x_{T\left(  i,j\right)  }^{p-1}\right)  _{i,p\in D\left\langle
j\right\rangle }\ \ \ \ \ \ \ \ \ \ \left(
\begin{array}
[c]{c}%
\text{since }\alpha\left(  i\right)  =T\left(  i,j\right)  \text{ for each
}i\in D\left\langle j\right\rangle \\
\text{(by the definition of }\alpha\text{)}%
\end{array}
\right)  .
\end{align*}
Thus,%
\[
\det\left(  x_{T\left(  i,j\right)  }^{p-1}\right)  _{i,p\in D\left\langle
j\right\rangle }=\det\left(  x_{p}^{\alpha^{-1}\left(  q\right)  -1}\right)
_{p,q\in X_{j}}=\det\left(  x_{q}^{\alpha^{-1}\left(  p\right)  -1}\right)
_{p,q\in X_{j}}%
\]
(since the determinant of a matrix equals the determinant of its transpose).
Therefore,%
\begin{align*}
\det\left(  x_{T\left(  i,j\right)  }^{p-1}\right)  _{i,p\in D\left\langle
j\right\rangle }  &  =\det\left(  x_{q}^{\alpha^{-1}\left(  p\right)
-1}\right)  _{p,q\in X_{j}}\\
&  =\sum_{\sigma\in S_{X_{j}}}\left(  -1\right)  ^{\sigma}\prod_{p\in X_{j}%
}x_{\sigma\left(  p\right)  }^{\alpha^{-1}\left(  p\right)  -1}%
\ \ \ \ \ \ \ \ \ \ \left(
\begin{array}
[c]{c}%
\text{by the definition of}\\
\text{a determinant}%
\end{array}
\right) \\
&  =\sum_{\sigma\in S_{X_{j}}}\left(  -1\right)  ^{\sigma}\prod_{i\in
D\left\langle j\right\rangle }\underbrace{x_{\sigma\left(  \alpha\left(
i\right)  \right)  }^{\alpha^{-1}\left(  \alpha\left(  i\right)  \right)  -1}%
}_{\substack{=x_{\sigma\left(  T\left(  i,j\right)  \right)  }^{i-1}%
\\\text{(since }\alpha^{-1}\left(  \alpha\left(  i\right)  \right)
=i\\\text{and }\alpha\left(  i\right)  =T\left(  i,j\right)  \\\text{(by the
definition of }\alpha\text{))}}}\\
&  \ \ \ \ \ \ \ \ \ \ \ \ \ \ \ \ \ \ \ \ \left(
\begin{array}
[c]{c}%
\text{here, we have substituted }\alpha\left(  i\right)  \text{ for }p\text{
in}\\
\text{the product, since the map }\alpha:D\left\langle j\right\rangle
\rightarrow X_{j}\\
\text{is a bijection}%
\end{array}
\right) \\
&  =\sum_{\sigma\in S_{X_{j}}}\left(  -1\right)  ^{\sigma}\prod_{i\in
D\left\langle j\right\rangle }x_{\sigma\left(  T\left(  i,j\right)  \right)
}^{i-1}\\
&  =\sum_{w\in S_{X_{j}}}\left(  -1\right)  ^{w}\prod_{i\in D\left\langle
j\right\rangle }x_{w\left(  T\left(  i,j\right)  \right)  }^{i-1}%
\end{align*}
(here, we have renamed the summation index $\sigma$ as $w$). This proves Claim 5.
\end{proof}

Now, (\ref{pf.thm.spechtmod.vdm.b.as-prod}) becomes%
\begin{align*}
\beta\left(  \mathbf{e}_{T}\right)   &  =\prod_{j=1}^{k}\ \ \underbrace{\sum
_{w\in S_{X_{j}}}\left(  -1\right)  ^{w}\prod_{i\in D\left\langle
j\right\rangle }x_{w\left(  T\left(  i,j\right)  \right)  }^{i-1}%
}_{\substack{=\det\left(  x_{T\left(  i,j\right)  }^{p-1}\right)  _{i,p\in
D\left\langle j\right\rangle }\\\text{(by (\ref{pf.thm.spechtmod.vdm.b.det}%
))}}}\\
&  =\prod_{j=1}^{k}\det\left(  x_{T\left(  i,j\right)  }^{p-1}\right)
_{i,p\in D\left\langle j\right\rangle }\\
&  =\prod_{j\in\mathbb{Z}}\det\left(  x_{T\left(  i,j\right)  }^{p-1}\right)
_{i,p\in D\left\langle j\right\rangle }\ \ \ \ \ \ \ \ \ \ \left(  \text{by
(\ref{pf.thm.spechtmod.vdm.b.1.get-rid})}\right)  .
\end{align*}
In other words, $\beta$ sends the polytabloid $\mathbf{e}_{T}$ to the product
$\prod_{j\in\mathbb{Z}}\det\left(  x_{T\left(  i,j\right)  }^{p-1}\right)
_{i,p\in D\left\langle j\right\rangle }$. This proves the first claim of
Theorem \ref{thm.spechtmod.vdm} \textbf{(b)}.

It remains to prove the second claim -- which says that the map $\beta$ sends
the Specht module $\mathcal{S}^{D}$ to the span of the products shown in
(\ref{eq.thm.spechtmod.vdm.b.prod}). But this is easy: The map $\beta$ is
$\mathbf{k}$-linear, and sends each polytabloid $\mathbf{e}_{T}$ to the
product shown in (\ref{eq.thm.spechtmod.vdm.b.prod}) (as we have just shown).
Thus, it sends the span of the polytabloids $\mathbf{e}_{T}$ to the span of
the products shown in (\ref{eq.thm.spechtmod.vdm.b.prod}). In other words, it
sends the Specht module $\mathcal{S}^{D}$ to the span of the products shown in
(\ref{eq.thm.spechtmod.vdm.b.prod}) (since the Specht module $\mathcal{S}^{D}$
is precisely the span of the polytabloids $\mathbf{e}_{T}$). The proof of
Theorem \ref{thm.spechtmod.vdm} \textbf{(b)} is thus complete. \medskip

\textbf{(c)} Let $T$ be an $n$-tableau of shape $D$. Define the sets
$D\left\langle j\right\rangle $ as in Theorem \ref{thm.spechtmod.vdm}
\textbf{(b)}. Note that each integer $j\leq0$ satisfies $D\left\langle
j\right\rangle =\varnothing$ (since all cells of $D$ lie in columns
$1,2,3,\ldots$) and therefore%
\begin{equation}
\det\underbrace{\left(  x_{T\left(  i,j\right)  }^{p-1}\right)  _{i,p\in
D\left\langle j\right\rangle }}_{\substack{\text{a }0\times0\text{-matrix}%
\\\text{(since }D\left\langle j\right\rangle =\varnothing\text{)}}}=1.
\label{pf.thm.spechtmod.vdm.c.0}%
\end{equation}

But Theorem \ref{thm.spechtmod.vdm} \textbf{(b)} yields%
\begin{equation}
\beta\left(  \mathbf{e}_{T}\right)  =\prod_{j\in\mathbb{Z}}\det\left(
x_{T\left(  i,j\right)  }^{p-1}\right)  _{i,p\in D\left\langle j\right\rangle
}. \label{pf.thm.spechtmod.vdm.c.1a}%
\end{equation}
In this product, all the factors with $j\leq0$ equal $1$ (by
(\ref{pf.thm.spechtmod.vdm.c.0})). Hence, we can drop all these factors. Thus,
(\ref{pf.thm.spechtmod.vdm.c.1a}) simplifies to
\begin{equation}
\beta\left(  \mathbf{e}_{T}\right)  =\prod_{j\geq1}\det\left(  x_{T\left(
i,j\right)  }^{p-1}\right)  _{i,p\in D\left\langle j\right\rangle }.
\label{pf.thm.spechtmod.vdm.c.1}%
\end{equation}

Now, let $j$ be any positive integer. Then, the cells in the $j$-th column of
$Y\left(  \lambda\right)  $ are $\left(  1,j\right)  ,\ \left(  2,j\right)
,\ \ldots,\ \left(  \lambda_{j}^{t},j\right)  $ (by Theorem
\ref{thm.partitions.conj} \textbf{(e)}). In other words, the integers $i$ that
satisfy $\left(  i,j\right)  \in Y\left(  \lambda\right)  $ are $1,2,\ldots
,\lambda_{j}^{t}$. Since $D=Y\left(  \lambda\right)  $, we can rewrite this as
follows: The integers $i$ that satisfy $\left(  i,j\right)  \in D$ are
$1,2,\ldots,\lambda_{j}^{t}$. In other words, $D\left\langle j\right\rangle
=\left\{  1,2,\ldots,\lambda_{j}^{t}\right\}  $ (since $D\left\langle
j\right\rangle $ is the set of all integers $i$ that satisfy $\left(
i,j\right)  \in D$).

Set $q:=\lambda_{j}^{t}$. Thus,
\[
D\left\langle j\right\rangle =\left\{  1,2,\ldots,\lambda_{j}^{t}\right\}
=\left[  \lambda_{j}^{t}\right]  =\left[  q\right]
\ \ \ \ \ \ \ \ \ \ \left(  \text{since }\lambda_{j}^{t}=q\right)  .
\]

Hence, the determinant $\det\left(  x_{T\left(  i,j\right)  }^{p-1}\right)
_{i,p\in D\left\langle j\right\rangle }$ is an instance of the well-known
\emph{Vandermonde determinant}. Let us recall what this means: The
\emph{Vandermonde determinant formula} (\cite[Theorem 6.4.31 \textbf{(c)}%
]{21s}) says that if $k\in\mathbb{N}$ is arbitrary, and if $a_{1},a_{2}%
,\ldots,a_{k}$ are any $k$ elements of any commutative ring, then%
\begin{equation}
\det\left(
\begin{array}
[c]{cccc}%
a_{1}^{0} & a_{1}^{1} & \cdots & a_{1}^{k-1}\\
a_{2}^{0} & a_{2}^{1} & \cdots & a_{2}^{k-1}\\
\vdots & \vdots & \ddots & \vdots\\
a_{k}^{0} & a_{k}^{1} & \cdots & a_{k}^{k-1}%
\end{array}
\right)  =\prod_{1\leq i_{1}<i_{2}\leq k}\left(  a_{i_{2}}-a_{i_{1}}\right)  .
\label{pf.thm.spechtmod.vdm.c.vdm}%
\end{equation}
Thus,%
\begin{align*}
\det\left(  x_{T\left(  i,j\right)  }^{p-1}\right)  _{i,p\in\left[  q\right]
}  &  =\det\left(
\begin{array}
[c]{cccc}%
x_{T\left(  1,j\right)  }^{0} & x_{T\left(  1,j\right)  }^{1} & \cdots &
x_{T\left(  1,j\right)  }^{q-1}\\
x_{T\left(  2,j\right)  }^{0} & x_{T\left(  2,j\right)  }^{1} & \cdots &
x_{T\left(  2,j\right)  }^{q-1}\\
\vdots & \vdots & \ddots & \vdots\\
x_{T\left(  q,j\right)  }^{0} & x_{T\left(  q,j\right)  }^{1} & \cdots &
x_{T\left(  q,j\right)  }^{q-1}%
\end{array}
\right) \\
&  =\prod_{1\leq i_{1}<i_{2}\leq q}\left(  x_{T\left(  i_{2},j\right)
}-x_{T\left(  i_{1},j\right)  }\right)
\end{align*}
(by (\ref{pf.thm.spechtmod.vdm.c.vdm}), applied to $k=q$ and $a_{i}%
=x_{T\left(  i,j\right)  }$). In view of $\left[  q\right]  =D\left\langle
j\right\rangle $ and $q=\lambda_{j}^{t}$, we can rewrite this as%
\begin{equation}
\det\left(  x_{T\left(  i,j\right)  }^{p-1}\right)  _{i,p\in D\left\langle
j\right\rangle }=\prod_{1\leq i_{1}<i_{2}\leq\lambda_{j}^{t}}\left(
x_{T\left(  i_{2},j\right)  }-x_{T\left(  i_{1},j\right)  }\right)  .
\label{pf.thm.spechtmod.vdm.c.7}%
\end{equation}

Forget that we fixed $j$. We thus have proved (\ref{pf.thm.spechtmod.vdm.c.7})
for each positive integer $j$. Therefore, (\ref{pf.thm.spechtmod.vdm.c.1})
rewrites as%
\[
\beta\left(  \mathbf{e}_{T}\right)  =\prod_{j\geq1}\ \ \prod_{1\leq
i_{1}<i_{2}\leq\lambda_{j}^{t}}\left(  x_{T\left(  i_{2},j\right)
}-x_{T\left(  i_{1},j\right)  }\right)  .
\]
This proves Theorem \ref{thm.spechtmod.vdm} \textbf{(c)}.
\end{proof}

The Specht--Vandermonde avatar (Theorem \ref{thm.spechtmod.vdm} \textbf{(c)})
takes on a particularly simple form when $D$ is the Young diagram $Y\left(
\lambda\right)  $ of a hook partition $\lambda=\left(  i,1^{j}\right)  $. In
this case, Theorem \ref{thm.spechtmod.vdm} \textbf{(c)} leads to the following:

\begin{corollary}
\label{cor.peel.vdm-ava}Let $\lambda$ be a hook partition of $n$. (See
Definition \ref{def.partitions.hook} for the notion of a hook partition.) Let
$D$ be the Young diagram $Y\left(  \lambda\right)  $. Consider the injective
$\mathbf{k}\left[  S_{n}\right]  $-module morphism $\beta:\mathcal{M}%
^{D}\rightarrow\mathbf{k}\left[  x_{1},x_{2},\ldots,x_{n}\right]  $ from
Theorem \ref{thm.spechtmod.vdm}. Then: \medskip

\textbf{(a)} The map $\beta$ sends any polytabloid $\mathbf{e}_{T}$ to%
\[
\prod_{1\leq u<v\leq\ell\left(  \lambda\right)  }\left(  x_{T\left(
v,1\right)  }-x_{T\left(  u,1\right)  }\right)  .
\]

\textbf{(b)} The image of the Specht module $\mathcal{S}^{D}$ under the
morphism $\beta$ is the submodule%
\[
\operatorname*{span}\nolimits_{\mathbf{k}}\left\{  \prod_{1\leq u<v\leq
\ell\left(  \lambda\right)  }\left(  x_{\sigma\left(  v\right)  }%
-x_{\sigma\left(  u\right)  }\right)  \ \mid\ \sigma\in S_{n}\right\}
\]
of $\mathbf{k}\left[  x_{1},x_{2},\ldots,x_{n}\right]  $.
\end{corollary}

\begin{proof}
Let $\lambda^{t}$ be the conjugate of the partition $\lambda$ (as defined in
Theorem \ref{thm.partitions.conj}).

We know that $\lambda$ is a hook partition. In other words, $\lambda=\left(
p,1^{q}\right)  $ for some positive integer $p$ and some nonnegative integer
$q$. Thus, any column of $Y\left(  \lambda\right)  $ except for the first one
has at most one cell. In other words, for each integer $j\geq2$, we have%
\begin{equation}
\left(  \text{\# of cells in the }j\text{-th column of }Y\left(
\lambda\right)  \right)  \leq1. \label{pf.cor.peel.vdm-ava.1}%
\end{equation}
Hence, for each integer $j\geq2$, we have%
\begin{align}
\lambda_{j}^{t}  &  =\left(  \text{\# of cells in the }j\text{-th column of
}Y\left(  \lambda\right)  \right)  \ \ \ \ \ \ \ \ \ \ \left(  \text{by
Theorem \ref{thm.partitions.conj} \textbf{(b)}}\right) \nonumber\\
&  \leq1\ \ \ \ \ \ \ \ \ \ \left(  \text{by (\ref{pf.cor.peel.vdm-ava.1}%
)}\right)  . \label{pf.cor.peel.vdm-ava.2}%
\end{align}

Furthermore, Theorem \ref{thm.partitions.conj} \textbf{(c)} yields
$\lambda_{1}^{t}=\ell\left(  \lambda\right)  $. \medskip

\textbf{(a)} The map $\beta$ sends any polytabloid $\mathbf{e}_{T}$ to%
\begin{align*}
&  \prod_{j\geq1}\ \ \prod_{1\leq i_{1}<i_{2}\leq\lambda_{j}^{t}}\left(
x_{T\left(  i_{2},j\right)  }-x_{T\left(  i_{1},j\right)  }\right)
\ \ \ \ \ \ \ \ \ \ \left(  \text{by Theorem \ref{thm.spechtmod.vdm}
\textbf{(c)}}\right) \\
&  =\left(  \prod_{1\leq i_{1}<i_{2}\leq\lambda_{1}^{t}}\left(  x_{T\left(
i_{2},1\right)  }-x_{T\left(  i_{1},1\right)  }\right)  \right)  \cdot
\prod_{j\geq2}\ \ \underbrace{\prod_{1\leq i_{1}<i_{2}\leq\lambda_{j}^{t}%
}\left(  x_{T\left(  i_{2},j\right)  }-x_{T\left(  i_{1},j\right)  }\right)
}_{\substack{=\left(  \text{empty product}\right)  \\\text{(since
(\ref{pf.cor.peel.vdm-ava.2}) shows that }\lambda_{j}^{t}\leq1\text{, and thus
there are}\\\text{no two integers }i_{1}\text{ and }i_{2}\text{ satisfying
}1\leq i_{1}<i_{2}\leq\lambda_{j}^{t}\text{)}}}\\
&  \ \ \ \ \ \ \ \ \ \ \ \ \ \ \ \ \ \ \ \ \left(
\begin{array}
[c]{c}%
\text{here, we have split off the factor for }j=1\text{ from}\\
\text{the outer product}%
\end{array}
\right) \\
&  =\left(  \prod_{1\leq i_{1}<i_{2}\leq\lambda_{1}^{t}}\left(  x_{T\left(
i_{2},1\right)  }-x_{T\left(  i_{1},1\right)  }\right)  \right)  \cdot
\prod_{j\geq2}\underbrace{\left(  \text{empty product}\right)  }_{=1}\\
&  =\left(  \prod_{1\leq i_{1}<i_{2}\leq\lambda_{1}^{t}}\left(  x_{T\left(
i_{2},1\right)  }-x_{T\left(  i_{1},1\right)  }\right)  \right)
\cdot\underbrace{\prod_{j\geq2}1}_{=1}\\
&  =\prod_{1\leq i_{1}<i_{2}\leq\lambda_{1}^{t}}\left(  x_{T\left(
i_{2},1\right)  }-x_{T\left(  i_{1},1\right)  }\right) \\
&  =\prod_{1\leq u<v\leq\lambda_{1}^{t}}\left(  x_{T\left(  v,1\right)
}-x_{T\left(  u,1\right)  }\right)  \ \ \ \ \ \ \ \ \ \ \left(
\begin{array}
[c]{c}%
\text{here, we have renamed the}\\
\text{indices }i_{1}\text{ and }i_{2}\text{ as }u\text{ and }v
\end{array}
\right) \\
&  =\prod_{1\leq u<v\leq\ell\left(  \lambda\right)  }\left(  x_{T\left(
v,1\right)  }-x_{T\left(  u,1\right)  }\right)  \ \ \ \ \ \ \ \ \ \ \left(
\text{since }\lambda_{1}^{t}=\ell\left(  \lambda\right)  \right)  .
\end{align*}
This proves Corollary \ref{cor.peel.vdm-ava} \textbf{(a)}. \medskip

\textbf{(b)} We have $D=Y\left(  \lambda\right)  $. Thus, the cells in the
$1$-st column of $D=Y\left(  \lambda\right)  $ are $\left(  1,1\right)
,\ \left(  2,1\right)  ,\ \ldots,\ \left(  \lambda_{1}^{t},1\right)  $ (by
Theorem \ref{thm.partitions.conj} \textbf{(e)}, applied to $j=1$).

Let $T$ be the $n$-tableau of shape $D$ that is obtained by filling the
numbers $1,2,\ldots,\lambda_{1}^{t}$ into these cells $\left(  1,1\right)
,\ \left(  2,1\right)  ,\ \ldots,\ \left(  \lambda_{1}^{t},1\right)  $ and
filling the remaining numbers $\lambda_{1}^{t}+1,\ \lambda_{1}^{t}%
+2,\ \ldots,n$ into the remaining cells of $D$ (say, from left to right,
although the exact order is immaterial). Then,
\[
T\left(  j,1\right)  =j\ \ \ \ \ \ \ \ \ \ \text{for each }j\in\left[
\lambda_{1}^{t}\right]
\]
(since we filled the numbers $1,2,\ldots,\lambda_{1}^{t}$ into the cells
$\left(  1,1\right)  ,\ \left(  2,1\right)  ,\ \ldots,\ \left(  \lambda
_{1}^{t},1\right)  $). In other words,%
\begin{equation}
T\left(  j,1\right)  =j\ \ \ \ \ \ \ \ \ \ \text{for each }j\in\left[
\ell\left(  \lambda\right)  \right]  \label{pf.cor.peel.vdm-ava.b.1}%
\end{equation}
(since $\lambda_{1}^{t}=\ell\left(  \lambda\right)  $). Hence, for each
$u,v\in\left[  \ell\left(  \lambda\right)  \right]  $, we have%
\begin{equation}
\underbrace{x_{T\left(  v,1\right)  }}_{\substack{=x_{v}\\\text{(since
(\ref{pf.cor.peel.vdm-ava.b.1})}\\\text{yields }T\left(  v,1\right)
=v\text{)}}}-\underbrace{x_{T\left(  u,1\right)  }}_{\substack{=x_{u}%
\\\text{(since (\ref{pf.cor.peel.vdm-ava.b.1})}\\\text{yields }T\left(
u,1\right)  =u\text{)}}}=x_{v}-x_{u}. \label{pf.cor.peel.vdm-ava.b.2}%
\end{equation}

Now, the map $\beta$ sends the polytabloid $\mathbf{e}_{T}$ to%
\begin{align*}
&  \prod_{1\leq u<v\leq\ell\left(  \lambda\right)  }\underbrace{\left(
x_{T\left(  v,1\right)  }-x_{T\left(  u,1\right)  }\right)  }%
_{\substack{=x_{v}-x_{u}\\\text{(by (\ref{pf.cor.peel.vdm-ava.b.2}))}%
}}\ \ \ \ \ \ \ \ \ \ \left(  \text{by Corollary \ref{cor.peel.vdm-ava}
\textbf{(a)}}\right) \\
&  =\prod_{1\leq u<v\leq\ell\left(  \lambda\right)  }\left(  x_{v}%
-x_{u}\right)  .
\end{align*}
In other words,%
\begin{equation}
\beta\left(  \mathbf{e}_{T}\right)  =\prod_{1\leq u<v\leq\ell\left(
\lambda\right)  }\left(  x_{v}-x_{u}\right)  . \label{pf.cor.peel.vdm-ava.b.4}%
\end{equation}

But $\beta$ is a left $\mathbf{k}\left[  S_{n}\right]  $-module morphism.
Hence, for each $w\in S_{n}$, we have%
\begin{align}
\beta\left(  w\mathbf{e}_{T}\right)   &  =w\cdot\beta\left(  \mathbf{e}%
_{T}\right)  =w\cdot\prod_{1\leq u<v\leq\ell\left(  \lambda\right)  }\left(
x_{v}-x_{u}\right)  \ \ \ \ \ \ \ \ \ \ \left(  \text{by
(\ref{pf.cor.peel.vdm-ava.b.4})}\right) \nonumber\\
&  =\prod_{1\leq u<v\leq\ell\left(  \lambda\right)  }\left(  x_{w\left(
v\right)  }-x_{w\left(  u\right)  }\right)  \label{pf.cor.peel.vdm-ava.b.5}%
\end{align}
(since $S_{n}$ acts on the polynomial ring $\mathbf{k}\left[  x_{1}%
,x_{2},\ldots,x_{n}\right]  $ by permuting the indeterminates). However, for
each $w\in S_{n}$, we have $\mathbf{e}_{wT}=w\mathbf{e}_{T}$ (by Lemma
\ref{lem.spechtmod.submod} \textbf{(a)}) and thus
\begin{equation}
\beta\left(  \mathbf{e}_{wT}\right)  =\beta\left(  w\mathbf{e}_{T}\right)
=\prod_{1\leq u<v\leq\ell\left(  \lambda\right)  }\left(  x_{w\left(
v\right)  }-x_{w\left(  u\right)  }\right)  \label{pf.cor.peel.vdm-ava.b.6}%
\end{equation}
(by (\ref{pf.cor.peel.vdm-ava.b.5})).

Now, recall the equality (\ref{pf.thm.spechtmod.leftideal.b.2}), which we
proved above. It says that%
\[
\mathcal{S}^{D}=\operatorname*{span}\nolimits_{\mathbf{k}}\left\{
\mathbf{e}_{wT}\ \mid\ w\in S_{n}\right\}  .
\]
Applying the map $\beta$ to it, we obtain%
\begin{align*}
\beta\left(  \mathcal{S}^{D}\right)   &  =\beta\left(  \operatorname*{span}%
\nolimits_{\mathbf{k}}\left\{  \mathbf{e}_{wT}\ \mid\ w\in S_{n}\right\}
\right) \\
&  =\operatorname*{span}\nolimits_{\mathbf{k}}\left\{  \beta\left(
\mathbf{e}_{wT}\right)  \ \mid\ w\in S_{n}\right\}
\ \ \ \ \ \ \ \ \ \ \left(
\begin{array}
[c]{c}%
\text{since the map }\beta\text{ is }\mathbf{k}\text{-linear,}\\
\text{and thus respects spans}%
\end{array}
\right) \\
&  =\operatorname*{span}\nolimits_{\mathbf{k}}\left\{  \prod_{1\leq
u<v\leq\ell\left(  \lambda\right)  }\left(  x_{w\left(  v\right)
}-x_{w\left(  u\right)  }\right)  \ \mid\ w\in S_{n}\right\}
\ \ \ \ \ \ \ \ \ \ \left(  \text{by (\ref{pf.cor.peel.vdm-ava.b.6})}\right)
\\
&  =\operatorname*{span}\nolimits_{\mathbf{k}}\left\{  \prod_{1\leq
u<v\leq\ell\left(  \lambda\right)  }\left(  x_{\sigma\left(  v\right)
}-x_{\sigma\left(  u\right)  }\right)  \ \mid\ \sigma\in S_{n}\right\}
\end{align*}
(here, we have renamed the index $w$ as $\sigma$). Thus, we have proved
Corollary \ref{cor.peel.vdm-ava} \textbf{(b)}.
\end{proof}

Corollary \ref{cor.peel.vdm-ava} \textbf{(b)} has been used by Peel to study
the Specht modules $\mathcal{S}^{Y\left(  \lambda\right)  }$ for hook
partitions $\lambda$ in \cite{Peel71} (see \cite{peelMO} for an introduction).

\subsubsection{The letterplace avatar}

The Specht--Vandermonde avatar represents $n$-tabloids and polytabloids as
polynomials in $n$ indeterminates $x_{1},x_{2},\ldots,x_{n}$. There is another
avatar which is even simpler, in that it only uses multilinear polynomials
(i.e., polynomials of degree $\leq1$ in each variable). However, this
simplicity is obtained at the cost of requiring more indeterminates: Namely,
it needs $nk$ indeterminates $x_{1,1},x_{1,2},\ldots,x_{n,k}$ (the subscripts
range over all pairs $\left(  i,j\right)  $ with $i\in\left[  n\right]  $ and
$j\in\left[  k\right]  $), where $k$ is chosen such that all cells of $D$ lie
in rows $1,2,\ldots,k$. Other than this, it is very similar to the
Specht--Vandermonde avatar, mainly differing in that each $x_{T\left(
i,j\right)  }^{i-1}$ is replaced by $x_{T\left(  i,j\right)  ,i}$ (so to
speak, the exponent $i-1$ is \textquotedblleft lowered\textquotedblright\ to
become a second subscript):

\begin{theorem}
[letterplace avatar of Young and Specht modules]\label{thm.spechtmod.det}Let
$D$ be a diagram with $\left\vert D\right\vert =n$. Let $k\in\mathbb{N}$ be
such that all cells of $D$ lie in rows $1,2,\ldots,k$.

Consider the polynomial ring
\begin{align*}
&  \mathbf{k}\left[  x_{i,j}\ \mid\ i\in\left[  n\right]  \text{ and }%
j\in\left[  k\right]  \right] \\
&  =\mathbf{k}\left[  x_{1,1},\ x_{1,2},\ \ldots,\ x_{1,k},\ x_{2,1}%
,\ x_{2,2},\ \ldots,\ x_{2,k},\ \ldots,\ x_{n,k}\right]
\end{align*}
in $nk$ indeterminates that are named $x_{i,j}$ and indexed by the pairs
$\left(  i,j\right)  \in\left[  n\right]  \times\left[  k\right]  $. This ring
is known as the \emph{letterplace algebra} (and the indeterminates $x_{i,j}$
are called the \emph{letterplace variables}, with the index $i$ being called
the \emph{letter} and the index $j$ being called the \emph{place}). Let the
symmetric group $S_{n}$ act on this ring by permuting the variables according
to the rule%
\[
g\rightharpoonup x_{i,j}=x_{g\left(  i\right)  ,j}%
\ \ \ \ \ \ \ \ \ \ \text{for all }g\in S_{n}\text{ and }i\in\left[  n\right]
\text{ and }j\in\left[  k\right]  .
\]
(That is, a permutation $g\in S_{n}$ acts on a polynomial $p$ by replacing
each indeterminate $x_{i,j}$ with $x_{g\left(  i\right)  ,j}$. This is called
the \emph{action on letters}.)

We define a $\mathbf{k}$-linear map%
\begin{align*}
\gamma:\mathcal{M}^{D}  &  \rightarrow\mathbf{k}\left[  x_{i,j}\ \mid
\ i\in\left[  n\right]  \text{ and }j\in\left[  k\right]  \right]  ,\\
\overline{T}  &  \mapsto\prod_{\left(  i,j\right)  \in D}x_{T\left(
i,j\right)  ,i}=\prod_{i\geq1}\ \ \prod_{m\in\operatorname*{Row}\left(
i,T\right)  }x_{m,i},
\end{align*}
where $\operatorname*{Row}\left(  i,T\right)  $ denotes the set of all entries
in the $i$-th row of $T$. Then: \medskip

\textbf{(a)} This map $\gamma$ is an injective left $\mathbf{k}\left[
S_{n}\right]  $-module morphism. \medskip

\textbf{(b)} This map $\gamma$ sends any polytabloid $\mathbf{e}_{T}%
\in\mathcal{M}^{D}$ to the determinant product%
\begin{equation}
\prod_{j\in\mathbb{Z}}\det\left(  x_{T\left(  i,j\right)  ,p}\right)  _{i,p\in
D\left\langle j\right\rangle }, \label{eq.thm.spechtmod.det.b.prod}%
\end{equation}
where $D\left\langle j\right\rangle $ is the set of all integers $i$ that
satisfy $\left(  i,j\right)  \in D$. (Note that the matrix $\left(
x_{T\left(  i,j\right)  ,p}\right)  _{i,p\in D\left\langle j\right\rangle }$
is thus a square matrix, whose rows and columns are indexed by the elements of
$D\left\langle j\right\rangle $. Such square matrices have well-defined
determinants, since $D\left\langle j\right\rangle $ is a finite set. When
$D\left\langle j\right\rangle $ is empty, this matrix is a $0\times0$-matrix
and thus has determinant $1$, by definition. Hence, all but finitely many
factors of the product in (\ref{eq.thm.spechtmod.det.b.prod}) equal $1$.)

Thus, the map $\gamma$ sends the Specht module $\mathcal{S}^{D}$ to the span
of the products shown in (\ref{eq.thm.spechtmod.det.b.prod}).
\end{theorem}

\begin{example}
Let $n$, $D$ and $T$ be as in Example \ref{exa.spechtmod.vdm.8}. Let $k=3$.
Then, the map $\gamma$ from Theorem \ref{thm.spechtmod.det} sends the
$n$-tabloid $\overline{T}\in\mathcal{M}^{D}$ to the monomial
\begin{align*}
\prod_{\left(  i,j\right)  \in D}x_{T\left(  i,j\right)  ,i}  &  =x_{T\left(
1,1\right)  ,1}x_{T\left(  1,3\right)  ,1}x_{T\left(  1,4\right)
,1}x_{T\left(  2,1\right)  ,2}x_{T\left(  2,2\right)  ,2}x_{T\left(
3,1\right)  ,3}x_{T\left(  3,2\right)  ,3}x_{T\left(  3,4\right)  ,3}\\
&  =x_{4,1}x_{1,1}x_{3,1}x_{2,2}x_{8,2}x_{7,3}x_{6,3}x_{5,3},
\end{align*}
and (by Theorem \ref{thm.spechtmod.det} \textbf{(b)}) sends the polytabloid
$\mathbf{e}_{T}$ to the product%
\begin{align*}
&  \prod_{j\in\mathbb{Z}}\det\left(  x_{T\left(  i,j\right)  ,p}\right)
_{i,p\in D\left\langle j\right\rangle }\\
&  =\det\left(  x_{T\left(  i,1\right)  ,p}\right)  _{i,p\in D\left\langle
1\right\rangle }\cdot\det\left(  x_{T\left(  i,2\right)  ,p}\right)  _{i,p\in
D\left\langle 2\right\rangle }\\
&  \ \ \ \ \ \ \ \ \ \ \cdot\det\left(  x_{T\left(  i,3\right)  ,p}\right)
_{i,p\in D\left\langle 3\right\rangle }\cdot\det\left(  x_{T\left(
i,4\right)  ,p}\right)  _{i,p\in D\left\langle 4\right\rangle }\\
&  \ \ \ \ \ \ \ \ \ \ \ \ \ \ \ \ \ \ \ \ \left(
\begin{array}
[c]{c}%
\text{because for all }j\notin\left\{  1,2,3,4\right\}  \text{, the set
}D\left\langle j\right\rangle \text{ is empty,}\\
\text{and thus }\det\left(  x_{T\left(  i,j\right)  ,p}\right)  _{i,p\in
D\left\langle j\right\rangle }=1
\end{array}
\right) \\
&  =\det\left(
\begin{array}
[c]{ccc}%
x_{4,1} & x_{4,2} & x_{4,3}\\
x_{2,1} & x_{2,2} & x_{2,3}\\
x_{7,1} & x_{7,2} & x_{7,3}%
\end{array}
\right)  \cdot\det\left(
\begin{array}
[c]{cc}%
x_{8,2} & x_{8,3}\\
x_{6,2} & x_{6,3}%
\end{array}
\right)  \cdot\det\left(
\begin{array}
[c]{c}%
x_{1,1}%
\end{array}
\right)  \cdot\det\left(
\begin{array}
[c]{cc}%
x_{3,1} & x_{3,3}\\
x_{5,1} & x_{5,3}%
\end{array}
\right)  .
\end{align*}

\end{example}

\begin{proof}
[Proof of Theorem \ref{thm.spechtmod.det}.]Analogous to parts \textbf{(a)} and
\textbf{(b)} of Theorem \ref{thm.spechtmod.vdm}.
\end{proof}

The isomorphic avatar of Specht modules defined in Theorem
\ref{thm.spechtmod.det} (or, rather, an insubstantial generalization thereof)
has been used quite successfully to prove their properties, e.g., by Clausen
in \cite{Clause91} and by Krob in \cite[\S 2.2]{Krob95}.

\bigskip

We have now seen four avatars of Young and Specht modules (and the tableau
avatar depends on the choice of an $n$-tableau $T$, so it gives many different
isomorphic copies).

\bigskip

\subsection{The standard polytabloids are linearly independent}

\subsubsection{Introduction and main result}

Young modules are easy to understand, at least as $\mathbf{k}$-modules: For
any diagram $D$ of size $n$, the Young module $\mathcal{M}^{D}$ has a basis
\[
\left(  \overline{T}\right)  _{\overline{T}\text{ is an }n\text{-tabloid of
shape }D}%
\]
(as a $\mathbf{k}$-module), so it is a free $\mathbf{k}$-module of rank%
\[
\left(  \text{\# of }n\text{-tabloids of shape }D\right)  =\dfrac{n!}%
{r_{1}!r_{2}!\cdots r_{k}!},
\]
where $r_{1},r_{2},\ldots,r_{k}$ are the sizes of the nonempty rows of $D$
(why?). The symmetric group $S_{n}$ acts on $\mathcal{M}^{D}$ by permuting the
$n$-tabloids.

The Specht module $\mathcal{S}^{D}$ is much more mysterious. It is a submodule
of $\mathcal{M}^{D}$ spanned by the $n!$ polytabloids $\mathbf{e}_{T}$, but
these polytabloids are heavily linearly dependent, so it is not clear how many
of them are redundant, i.e., what the dimension (rank) of $\mathcal{S}^{D}$
really is. In fact, it is not even clear whether $\mathcal{S}^{D}$ has a basis
for all commutative rings $\mathbf{k}$.

\begin{question}
\label{quest.spechtmod.bad.free}For a general diagram $D$ and a general
commutative ring $\mathbf{k}$, is the $\mathbf{k}$-module $\mathcal{S}^{D}$
free? Does it have a basis independent of $\mathbf{k}$ ? Is its rank
independent of $\mathbf{k}$ ? Is it a direct addend of $\mathcal{M}^{D}$ as a
$\mathbf{k}$-module?
\end{question}

These questions have been open since the 1970s, and work has been done for
various classes of diagrams $D$.

\begin{noncompile}
In particular, a positive answer is known whenever $D$ is a skew Young
diagram, but also for some wider classes (\textquotedblleft\%-avoiding
shapes\textquotedblright\ for example, done by Reiner and Shimozono in the 1990s).
\end{noncompile}

We shall give a complete (and positive) answer to this question when $D$ is a
skew Young diagram: We will show that the so-called \emph{standard
polytabloids} form a basis of $\mathcal{S}^{D}$ in this case. We can define
the standard polytabloids for any $D$ (we just cannot expect them to be a
basis in this generality):

\begin{definition}
\label{def.spechtmod.stdpt}Let $D$ be any diagram with $\left\vert
D\right\vert =n$. The \emph{standard polytabloids} of shape $D$ will mean the
polytabloids $\mathbf{e}_{T}\in\mathcal{S}^{D}$, where $T$ is a standard
tableau of shape $D$. (Recall from Proposition \ref{prop.tableau.std-n} that
any standard tableau $T$ of shape $D$ is an $n$-tableau, which is why
$\mathbf{e}_{T}\in\mathcal{S}^{D}$ is well-defined.)
\end{definition}

\begin{example}
The standard polytabloids of shape $Y\left(  2,2\right)  $ are $\mathbf{e}%
_{12\backslash\backslash34}$ and $\mathbf{e}_{13\backslash\backslash24}$,
since the standard tableaux of this shape are $12\backslash\backslash34$ and
$13\backslash\backslash24$.
\end{example}

We shall now make a first major step towards understanding the Specht module
$\mathcal{S}^{D}$ by proving the following theorem:

\begin{theorem}
\label{thm.spechtmod.linind}Let $D$ be any diagram. Then, the standard
polytabloids of shape $D$ are $\mathbf{k}$-linearly independent.
\end{theorem}

\subsubsection{The Young last letter order}

The proof of Theorem \ref{thm.spechtmod.linind} is not hard. It relies on a
\textquotedblleft triangularity-like argument\textquotedblright\ (or
\textquotedblleft leading term argument\textquotedblright). For this, we need
a total order on the set of all $n$-tabloids of shape $D$:

\begin{definition}
\label{def.tabloid.llo}\textbf{(a)} If $\overline{T}$ is any $n$-tabloid (of
any shape), and if $i\in\left[  n\right]  $, then we let $r_{\overline{T}%
}\left(  i\right)  $ be the number of the row of $T$ that contains $i$.
(Clearly, this depends only on $\overline{T}$, not on $T$, since we could just
as well have written \textquotedblleft the row of $\overline{T}$ that contains
$i$\textquotedblright\ instead of \textquotedblleft the row of $T$ that
contains $i$\textquotedblright.) \medskip

\textbf{(b)} Let $D$ be a diagram. Define a binary relation $<$ on the set
$\left\{  n\text{-tabloids of shape }D\right\}  $ as follows:

Let $\overline{T}$ and $\overline{S}$ be two $n$-tabloids of shape $D$. Then,
we declare that $\overline{T}<\overline{S}$ if and only if

\begin{itemize}
\item there exists at least one $i\in\left[  n\right]  $ such that
$r_{\overline{T}}\left(  i\right)  \neq r_{\overline{S}}\left(  i\right)  $, and

\item the \textbf{largest} such $i$ satisfies $r_{\overline{T}}\left(
i\right)  <r_{\overline{S}}\left(  i\right)  $.
\end{itemize}

In other words, we say that $\overline{T}<\overline{S}$ if the
\textbf{largest} number that appears in different rows in $\overline{T}$ and
$\overline{S}$ appears further north in $\overline{T}$ than in $\overline{S}$.
\end{definition}

\begin{example}
Let $n=7$ and $D=Y\left(  3,2,2\right)  $. Consider the two $n$-tableaux%
\[
T=\ytableaushort{237,46,15}\ \ \ \ \ \ \ \ \ \ \text{and}%
\ \ \ \ \ \ \ \ \ \ S=\ytableaushort{357,14,26}.
\]
Then, $r_{\overline{T}}\left(  2\right)  =1$ and $r_{\overline{S}}\left(
2\right)  =3$ and so on. We have $r_{\overline{T}}\left(  7\right)
=r_{\overline{S}}\left(  7\right)  =1$ but $r_{\overline{T}}\left(  6\right)
=2<3=r_{\overline{S}}\left(  6\right)  $. Thus, the \textbf{largest} number
that appears in different rows in $\overline{T}$ and $\overline{S}$ is $6$ and
therefore appears further north in $\overline{T}$ than in $\overline{S}$. In
other words, $\overline{T}<\overline{S}$.
\end{example}

It is easy to see that the relation $<$ that we have just defined is a total order:

\begin{proposition}
\label{prop.spechtmod.row.order}Let $D$ be a diagram. Then, the relation $<$
introduced in Definition \ref{def.tabloid.llo} \textbf{(b)} is the smaller
relation of a total order on the set $\left\{  n\text{-tabloids of shape
}D\right\}  $.
\end{proposition}

This total order $<$ is called the \emph{Young last letter order}.

The proof of Proposition \ref{prop.spechtmod.row.order} is fairly
straightforward, and is built upon an even more obvious lemma:

\begin{lemma}
\label{lem.spechtmod.row.equal}Let $P$ and $Q$ be two $n$-tableaux. Assume
that $r_{\overline{P}}\left(  i\right)  =r_{\overline{Q}}\left(  i\right)  $
for each $i\in\left[  n\right]  $. Then, $P$ and $Q$ are row-equivalent.
\end{lemma}

\begin{fineprint}
\begin{proof}
[Proof of Lemma \ref{lem.spechtmod.row.equal}.]Let $k\in\mathbb{Z}$. Recall
that $P$ is an $n$-tableau, and hence injective. In other words, all entries
of $P$ are distinct. In particular, all entries in the $k$-th row of $P$ are
distinct. In other words, the multiset of entries in the $k$-th row of $P$ is
a set (i.e., contains no nontrivial multiplicities). Similarly, the multiset
of entries in the $k$-th row of $Q$ is a set.

Let $u$ be an entry of the $k$-th row of $P$. Thus, the row of $P$ that
contains $u$ is the $k$-th row. In other words, $r_{\overline{P}}\left(
u\right)  =k$ (by the definition of $r_{\overline{P}}\left(  u\right)  $).
However, we know that $r_{\overline{P}}\left(  i\right)  =r_{\overline{Q}%
}\left(  i\right)  $ for each $i\in\left[  n\right]  $. Applying this to
$i=u$, we obtain $r_{\overline{P}}\left(  u\right)  =r_{\overline{Q}}\left(
u\right)  $. Hence, $r_{\overline{Q}}\left(  u\right)  =r_{\overline{P}%
}\left(  u\right)  =k$. In other words, the row of $Q$ that contains $u$ is
the $k$-th row (by the definition of $r_{\overline{Q}}\left(  u\right)  $).
Hence, $u$ is an entry of the $k$-th row of $Q$.

Now, forget that we fixed $u$. We thus have shown that if $u$ is an entry of
the $k$-th row of $P$, then $u$ is an entry of the $k$-th row of $Q$. In other
words, each entry of the $k$-th row of $P$ is an entry of the $k$-th row of
$Q$. In other words, the set of entries in the $k$-th row of $P$ is a subset
of the set of entries in the $k$-th row of $Q$. But the same argument (with
the roles of $P$ and $Q$ interchanged) shows that the set of entries in the
$k$-th row of $Q$ is a subset of the set of entries in the $k$-th row of $P$.
Combining these two containment relations, we conclude that the set of entries
in the $k$-th row of $P$ is the set of entries in the $k$-th row of $Q$. We
can furthermore replace the word \textquotedblleft set\textquotedblright\ by
\textquotedblleft multiset\textquotedblright\ in this sentence (since the
multiset of entries in the $k$-th row of $P$ is a set, and the multiset of
entries in the $k$-th row of $Q$ is a set). Thus, we have shown that the
multiset of entries in the $k$-th row of $P$ is the multiset of entries in the
$k$-th row of $Q$.

Now, forget that we fixed $k$. We have shown that for every $k\in\mathbb{Z}$,
the multiset of entries in the $k$-th row of $P$ is the multiset of entries in
the $k$-th row of $Q$. In other words, the tableaux $P$ and $Q$ have the same
entries in corresponding rows. In other words, $P$ and $Q$ are row-equivalent.
This proves Lemma \ref{lem.spechtmod.row.equal}.
\end{proof}
\end{fineprint}

\begin{fineprint}
\begin{proof}
[Proof of Proposition \ref{prop.spechtmod.row.order}.]In a nutshell: This
total order is actually a disguised lexicographic order: We can encode any
$n$-tabloid $\overline{T}$ of shape $D$ by the $n$-tuple \newline$\left(
r_{\overline{T}}\left(  n\right)  ,\ r_{\overline{T}}\left(  n-1\right)
,\ \ldots,\ r_{\overline{T}}\left(  1\right)  \right)  \in\mathbb{Z}^{n}$ (it
is easy to see that this encoding is injective, i.e., we can reconstruct
$\overline{T}$ from $\left(  r_{\overline{T}}\left(  n\right)  ,\ r_{\overline
{T}}\left(  n-1\right)  ,\ \ldots,\ r_{\overline{T}}\left(  1\right)  \right)
$), and then our relation $<$ is simply given by%
\begin{align*}
&  \ \left(  \overline{T}<\overline{S}\right) \\
&  \Longleftrightarrow\ \left(  \text{the encoding of }\overline{T}\text{ is
lexicographically smaller than the encoding of }\overline{S}\right)  .
\end{align*}
Thus, $<$ is a total order, since the lexicographic order on $\mathbb{Z}^{n}$
is a total order.

But let me also give a more direct and detailed proof:

First, we shall show that the relation $<$ is the smaller relation of a
partial order on $\left\{  n\text{-tabloids of shape }D\right\}  $. In order
to prove this, we need to verify three claims:

\begin{statement}
\textit{Claim 1:} The relation $<$ is irreflexive (i.e., there exists no
$n$-tabloid $\overline{T}$ such that $\overline{T}<\overline{T}$).
\end{statement}

\begin{statement}
\textit{Claim 2:} The relation $<$ is transitive (i.e., if three $n$-tabloids
$\overline{P}$, $\overline{Q}$ and $\overline{R}$ satisfy $\overline
{P}<\overline{Q}$ and $\overline{Q}<\overline{R}$, then $\overline
{P}<\overline{R}$).
\end{statement}

\begin{statement}
\textit{Claim 3:} The relation $<$ is asymmetric (i.e., no two $n$-tabloids
$\overline{P}$ and $\overline{Q}$ satisfy $\overline{P}<\overline{Q}$ and
$\overline{Q}<\overline{P}$ at the same time).
\end{statement}

\begin{proof}
[Proof of Claim 1.]Let $\overline{T}$ be an $n$-tabloid such that
$\overline{T}<\overline{T}$. We shall derive a contradiction.

Indeed, we have $\overline{T}<\overline{T}$. By the definition of the relation
$<$, this means that

\begin{itemize}
\item there exists at least one $i\in\left[  n\right]  $ such that
$r_{\overline{T}}\left(  i\right)  \neq r_{\overline{T}}\left(  i\right)  $, and

\item the \textbf{largest} such $i$ satisfies $r_{\overline{T}}\left(
i\right)  <r_{\overline{T}}\left(  i\right)  $.
\end{itemize}

\noindent But the first of these two bullet points is clearly absurd (since
$r_{\overline{T}}\left(  i\right)  \neq r_{\overline{T}}\left(  i\right)  $
can never happen). Thus, we have found a contradiction.

Forget that we fixed $\overline{T}$. We thus have obtained a contradiction for
each $n$-tabloid $\overline{T}$ such that $\overline{T}<\overline{T}$. Hence,
there exists no such $\overline{T}$. In other words, the relation $<$ is
irreflexive. This proves Claim 1.
\end{proof}

\begin{proof}
[Proof of Claim 2.]Let $\overline{P}$, $\overline{Q}$ and $\overline{R}$ be
three $n$-tabloids that satisfy $\overline{P}<\overline{Q}$ and $\overline
{Q}<\overline{R}$. We must prove that $\overline{P}<\overline{R}$.

Indeed, we have $\overline{P}<\overline{Q}$. By the definition of the relation
$<$, this means that

\begin{itemize}
\item there exists at least one $i\in\left[  n\right]  $ such that
$r_{\overline{P}}\left(  i\right)  \neq r_{\overline{Q}}\left(  i\right)  $, and

\item the \textbf{largest} such $i$ satisfies $r_{\overline{P}}\left(
i\right)  <r_{\overline{Q}}\left(  i\right)  $.
\end{itemize}

Let us denote this largest $i$ by $u$. Thus, we have $r_{\overline{P}}\left(
u\right)  <r_{\overline{Q}}\left(  u\right)  $, whereas%
\begin{equation}
\text{each }i>u\text{ satisfies }r_{\overline{P}}\left(  i\right)
=r_{\overline{Q}}\left(  i\right)  \label{pf.prop.spechtmod.row.order.c2.pf.1}%
\end{equation}
(since $u$ is the \textbf{largest} $i\in\left[  n\right]  $ that satisfies
$r_{\overline{P}}\left(  i\right)  \neq r_{\overline{Q}}\left(  i\right)  $).

Furthermore, we have $\overline{Q}<\overline{R}$. By the definition of the
relation $<$, this means that

\begin{itemize}
\item there exists at least one $i\in\left[  n\right]  $ such that
$r_{\overline{Q}}\left(  i\right)  \neq r_{\overline{R}}\left(  i\right)  $, and

\item the \textbf{largest} such $i$ satisfies $r_{\overline{Q}}\left(
i\right)  <r_{\overline{R}}\left(  i\right)  $.
\end{itemize}

Let us denote this largest $i$ by $v$. Thus, we have $r_{\overline{Q}}\left(
v\right)  <r_{\overline{R}}\left(  v\right)  $, whereas%
\begin{equation}
\text{each }i>v\text{ satisfies }r_{\overline{Q}}\left(  i\right)
=r_{\overline{R}}\left(  i\right)  \label{pf.prop.spechtmod.row.order.c2.pf.3}%
\end{equation}
(since $v$ is the \textbf{largest} $i\in\left[  n\right]  $ that satisfies
$r_{\overline{Q}}\left(  i\right)  \neq r_{\overline{R}}\left(  i\right)  $).

Now, let $w:=\max\left\{  u,v\right\}  $. Then, $w\geq u$ and $w\geq v$.
Hence,%
\begin{equation}
\text{each }i>w\text{ satisfies }r_{\overline{P}}\left(  i\right)
=r_{\overline{R}}\left(  i\right)  \label{pf.prop.spechtmod.row.order.c2.pf.5}%
\end{equation}
\footnote{\textit{Proof.} Let $i\in\left[  n\right]  $ be such that $i>w$.
Then, $i>w\geq u$, so that $r_{\overline{P}}\left(  i\right)  =r_{\overline
{Q}}\left(  i\right)  $ (by (\ref{pf.prop.spechtmod.row.order.c2.pf.1})).
Also, $i>w\geq v$, so that $r_{\overline{Q}}\left(  i\right)  =r_{\overline
{R}}\left(  i\right)  $ (by (\ref{pf.prop.spechtmod.row.order.c2.pf.3})).
Hence, $r_{\overline{P}}\left(  i\right)  =r_{\overline{Q}}\left(  i\right)
=r_{\overline{R}}\left(  i\right)  $. This proves
(\ref{pf.prop.spechtmod.row.order.c2.pf.5}).}. Furthermore,
\begin{equation}
r_{\overline{P}}\left(  w\right)  <r_{\overline{R}}\left(  w\right)
\label{pf.prop.spechtmod.row.order.c2.pf.7}%
\end{equation}
\footnote{\textit{Proof of (\ref{pf.prop.spechtmod.row.order.c2.pf.7}):} We
are in one of the following three cases:
\par
\textit{Case 1:} We have $u<v$.
\par
\textit{Case 2:} We have $u=v$.
\par
\textit{Case 3:} We have $u>v$.
\par
Let us first consider Case 1. In this case, we have $u<v$. Hence,
$w=\max\left\{  u,v\right\}  =v$ (since $u<v$). But we know that
$r_{\overline{Q}}\left(  v\right)  <r_{\overline{R}}\left(  v\right)  $. Since
$w=v$, we can rewrite this as $r_{\overline{Q}}\left(  w\right)
<r_{\overline{R}}\left(  w\right)  $. Furthermore, $w=v>u$ (since $u<v$) and
thus $r_{\overline{P}}\left(  w\right)  =r_{\overline{Q}}\left(  w\right)  $
(by (\ref{pf.prop.spechtmod.row.order.c2.pf.1}), applied to $i=w$). Thus,
$r_{\overline{P}}\left(  w\right)  =r_{\overline{Q}}\left(  w\right)
<r_{\overline{R}}\left(  w\right)  $. Hence,
(\ref{pf.prop.spechtmod.row.order.c2.pf.7}) is proved in Case 1.
\par
Let us next consider Case 2. In this case, we have $u=v$. Hence,
$w=\max\left\{  u,v\right\}  =u$ (since $u=v$). But we know that
$r_{\overline{Q}}\left(  v\right)  <r_{\overline{R}}\left(  v\right)  $. Since
$w=u=v$, we can rewrite this as $r_{\overline{Q}}\left(  w\right)
<r_{\overline{R}}\left(  w\right)  $. Furthermore, we know that $r_{\overline
{P}}\left(  u\right)  <r_{\overline{Q}}\left(  u\right)  $. Since $w=u$, we
can rewrite this as $r_{\overline{P}}\left(  w\right)  <r_{\overline{Q}%
}\left(  w\right)  $. Thus, $r_{\overline{P}}\left(  w\right)  <r_{\overline
{Q}}\left(  w\right)  <r_{\overline{R}}\left(  w\right)  $. Hence,
(\ref{pf.prop.spechtmod.row.order.c2.pf.7}) is proved in Case 2.
\par
Let us finally consider Case 3. In this case, we have $u>v$. Hence,
$w=\max\left\{  u,v\right\}  =u$ (since $u>v$). But we know that
$r_{\overline{P}}\left(  u\right)  <r_{\overline{Q}}\left(  u\right)  $. Since
$w=u$, we can rewrite this as $r_{\overline{P}}\left(  w\right)
<r_{\overline{Q}}\left(  w\right)  $. Furthermore, $w=u>v$ and thus
$r_{\overline{Q}}\left(  w\right)  =r_{\overline{R}}\left(  w\right)  $ (by
(\ref{pf.prop.spechtmod.row.order.c2.pf.3}), applied to $i=w$). Thus,
$r_{\overline{P}}\left(  w\right)  <r_{\overline{Q}}\left(  w\right)
=r_{\overline{R}}\left(  w\right)  $. Hence,
(\ref{pf.prop.spechtmod.row.order.c2.pf.7}) is proved in Case 3.
\par
Thus, (\ref{pf.prop.spechtmod.row.order.c2.pf.7}) is proved in all three Cases
1, 2 and 3.} and therefore $r_{\overline{P}}\left(  w\right)  \neq
r_{\overline{R}}\left(  w\right)  $. Hence, $w$ is an $i\in\left[  n\right]  $
that satisfies $r_{\overline{P}}\left(  i\right)  \neq r_{\overline{R}}\left(
i\right)  $. Moreover, $w$ is the \textbf{largest} such $i$ (since
(\ref{pf.prop.spechtmod.row.order.c2.pf.5}) shows that all larger $i$'s
satisfy $r_{\overline{P}}\left(  i\right)  =r_{\overline{R}}\left(  i\right)
$). Hence, we have now shown that

\begin{itemize}
\item there exists at least one $i\in\left[  n\right]  $ such that
$r_{\overline{P}}\left(  i\right)  \neq r_{\overline{R}}\left(  i\right)  $
(namely, $i=w$), and

\item the \textbf{largest} such $i$ satisfies $r_{\overline{P}}\left(
i\right)  <r_{\overline{R}}\left(  i\right)  $ (because the largest such $i$
is $w$, and we know that $w$ satisfies $r_{\overline{P}}\left(  w\right)
<r_{\overline{R}}\left(  w\right)  $).
\end{itemize}

In other words, $\overline{P}<\overline{R}$ (by the definition of the relation
$<$). This proves Claim 2.
\end{proof}

\begin{proof}
[Proof of Claim 3.]Any binary relation that is irreflexive and transitive is
automatically asymmetric (since $\overline{P}<\overline{Q}$ and $\overline
{Q}<\overline{P}$ would entail $\overline{P}<\overline{P}$ by transitivity,
but this would contradict irreflexivity). Since the relation $<$ is
irreflexive (by Claim 1) and transitive (by Claim 2), we thus conclude that it
is also asymmetric. This proves Claim 3.
\end{proof}

Having proved Claims 1, 2 and 3, we thus conclude that the relation $<$ is the
smaller relation of a partial order. It remains to show that this partial
order is a total order. For this purpose, we need to prove the following last claim:

\begin{statement}
\textit{Claim 4:} Let $\overline{P}$ and $\overline{Q}$ be two distinct
$n$-tabloids of shape $D$. Then, we have $\overline{P}<\overline{Q}$ or
$\overline{Q}<\overline{P}$.
\end{statement}

\begin{proof}
[Proof of Claim 4.]We shall first show that there exists some $i\in\left[
n\right]  $ such that $r_{\overline{P}}\left(  i\right)  \neq r_{\overline{Q}%
}\left(  i\right)  $.

Indeed, assume the contrary. Thus, $r_{\overline{P}}\left(  i\right)
=r_{\overline{Q}}\left(  i\right)  $ for each $i\in\left[  n\right]  $. Hence,
Lemma \ref{lem.spechtmod.row.equal} yields that $P$ and $Q$ are
row-equivalent. In other words, $\overline{P}=\overline{Q}$ (since an
$n$-tabloid is an equivalence class of $n$-tableaux with respect to
row-equivalence). But this contradicts the fact that $\overline{P}$ and
$\overline{Q}$ are distinct.

This contradiction shows that our assumption was false. Hence, we have shown
that there exists some $i\in\left[  n\right]  $ such that $r_{\overline{P}%
}\left(  i\right)  \neq r_{\overline{Q}}\left(  i\right)  $.

Consider the \textbf{largest} such $i$. Then, this $i$ satisfies
$r_{\overline{P}}\left(  i\right)  \neq r_{\overline{Q}}\left(  i\right)  $.
Hence, either $r_{\overline{P}}\left(  i\right)  <r_{\overline{Q}}\left(
i\right)  $ or $r_{\overline{P}}\left(  i\right)  >r_{\overline{Q}}\left(
i\right)  $. In the former case, we have $\overline{P}<\overline{Q}$ (by the
definition of the relation $<$). In the latter case, we have $\overline
{Q}<\overline{P}$ (by the definition of the relation $<$
again\footnote{\textit{Proof.} Assume that $r_{\overline{P}}\left(  i\right)
>r_{\overline{Q}}\left(  i\right)  $. Thus, $r_{\overline{Q}}\left(  i\right)
<r_{\overline{P}}\left(  i\right)  $. But recall that our $i$ is the
\textbf{largest} $i\in\left[  n\right]  $ such that $r_{\overline{P}}\left(
i\right)  \neq r_{\overline{Q}}\left(  i\right)  $. In other words, our $i$ is
the \textbf{largest} $i\in\left[  n\right]  $ such that $r_{\overline{Q}%
}\left(  i\right)  \neq r_{\overline{P}}\left(  i\right)  $. Hence, there
exists an $i\in\left[  n\right]  $ such that $r_{\overline{Q}}\left(
i\right)  \neq r_{\overline{P}}\left(  i\right)  $, and the \textbf{largest}
such $i$ satisfies $r_{\overline{Q}}\left(  i\right)  <r_{\overline{P}}\left(
i\right)  $. In other words, we have $\overline{Q}<\overline{P}$ (by the
definition of the relation $<$), qed.}). Thus, we have found that
$\overline{P}<\overline{Q}$ or $\overline{Q}<\overline{P}$. This proves Claim 4.
\end{proof}

Claim 4 shows that the partial order whose smaller relation is $<$ is actually
a total order. The proof of Proposition \ref{prop.spechtmod.row.order} is thus complete.
\end{proof}
\end{fineprint}

Now we claim the following:

\begin{lemma}
\label{lem.spechtmod.lead-term}Let $D$ be a diagram. Let $T$ be a
column-standard $n$-tableau of shape $D$. Let $<$ be the relation introduced
in Definition \ref{def.tabloid.llo}. Then: \medskip

\textbf{(a)} Every $w\in\mathcal{C}\left(  T\right)  $ satisfying
$w\neq\operatorname*{id}$ satisfies $\overline{wT}<\overline{T}$. \medskip

\textbf{(b)} We have%
\[
\mathbf{e}_{T}=\overline{T}+\left(  \text{a linear combination of
}n\text{-tabloids }\overline{S}\text{ with }\overline{S}<\overline{T}\right)
.
\]

\end{lemma}

\begin{proof}
\textbf{(a)} Let $w\in\mathcal{C}\left(  T\right)  $ be such that
$w\neq\operatorname*{id}$. Then, there exists some $k\in\left[  n\right]  $
such that $w\left(  k\right)  \neq k$ (since $w\neq\operatorname*{id}$). Pick
the \textbf{largest} such $k$. Then,
\begin{equation}
w\left(  i\right)  =i\ \ \ \ \ \ \ \ \ \ \text{for all }i>k.
\label{pf.lem.spechtmod.lead-term.a.1}%
\end{equation}

Here is the rest of the proof in a nutshell: The permutation $w$ is vertical
for $T$, and thus must permute the set of entries of the column of $T$ that
contains the entry $k$. All entries south of this entry $k$ stay unchanged
under this permutation (by (\ref{pf.lem.spechtmod.lead-term.a.1})), but the
entry $k$ does not (since $w\left(  k\right)  \neq k$). Hence, this entry $k$
must move further north. In other words, the entry $k$ lies further north in
$wT$ than it does in $T$. But all entries larger than $k$ lie in the same row
(and, in fact, in the same cell) in $wT$ as they do in $T$ (again by
(\ref{pf.lem.spechtmod.lead-term.a.1})). Thus, the \textbf{largest} number
that appears in different rows in $\overline{wT}$ and $\overline{T}$ is the
number $k$, and this number appears further north in $\overline{wT}$ than in
$\overline{T}$. By the definition of the relation $<$, this means that
$\overline{wT}<\overline{T}$, and thus Lemma \ref{lem.spechtmod.lead-term}
\textbf{(a)} is proved. \medskip

\begin{fineprint}
Here is a rigorous way to formulate this proof:

Let $\left(  u,v\right)  $ be the cell of $T$ that contains the entry $k$.
Thus, $T\left(  u,v\right)  =k$. Hence, the number $k$ lies in the $u$-th row
and in the $v$-th column of $T$. In particular, $r_{\overline{T}}\left(
k\right)  =u$.

Let $\left(  u^{\prime},v^{\prime}\right)  $ be the cell of the $n$-tableau
$wT$ that contains the entry $k$. Thus, $\left(  wT\right)  \left(  u^{\prime
},v^{\prime}\right)  =k$. Hence, the number $k$ lies in the $u^{\prime}$-th
row and in the $v^{\prime}$-th column of $wT$. In particular, $r_{\overline
{wT}}\left(  k\right)  =u^{\prime}$.

The definition of the action of $S_{n}$ on the $n$-tableaux yields $wT=w\circ
T$, so that $\left(  wT\right)  \left(  u^{\prime},v^{\prime}\right)
=w\left(  T\left(  u^{\prime},v^{\prime}\right)  \right)  $. Hence, $w\left(
T\left(  u^{\prime},v^{\prime}\right)  \right)  =\left(  wT\right)  \left(
u^{\prime},v^{\prime}\right)  =k$. Therefore, $T\left(  u^{\prime},v^{\prime
}\right)  =w^{-1}\left(  k\right)  $. This shows that the number
$w^{-1}\left(  k\right)  $ lies in the $v^{\prime}$-th column of $T$.

Note that $w^{-1}\left(  k\right)  \neq k$ (since $w^{-1}\left(  k\right)  =k$
would yield $k=w\left(  k\right)  \neq k$, which is absurd). Hence, $T\left(
u^{\prime},v^{\prime}\right)  =w^{-1}\left(  k\right)  \neq k=T\left(
u,v\right)  $, and thus $\left(  u^{\prime},v^{\prime}\right)  \neq\left(
u,v\right)  $.

But the permutation $w$ is vertical for $T$ (since $w\in\mathcal{C}\left(
T\right)  $). Hence, for each $i\in\left[  n\right]  $, the number $w\left(
i\right)  $ lies in the same column of $T$ as $i$ does (by the definition of a
vertical permutation). Applying this to $i=w^{-1}\left(  k\right)  $, we
conclude that the number $w\left(  w^{-1}\left(  k\right)  \right)  $ lies in
the same column of $T$ as $w^{-1}\left(  k\right)  $ does. In other words, the
number $k$ lies in the same column of $T$ as $w^{-1}\left(  k\right)  $ does
(since $w\left(  w^{-1}\left(  k\right)  \right)  =k$). In other words, we
have $v^{\prime}=v$ (since the number $w^{-1}\left(  k\right)  $ lies in the
$v$-th column of $T$, whereas the number $k$ lies in the $v$-th column of
$T$). Thus, $v=v^{\prime}$. Hence, $\left(  u^{\prime},v\right)  =\left(
u^{\prime},v^{\prime}\right)  \in D$, and $T\left(  u^{\prime},v\right)
=T\left(  u^{\prime},v^{\prime}\right)  =w^{-1}\left(  k\right)  $.

If we had $u^{\prime}=u$, then we would have $\left(  u^{\prime},v^{\prime
}\right)  =\left(  u,v\right)  $ (since $v=v^{\prime}$), which would
contradict $\left(  u^{\prime},v^{\prime}\right)  \neq\left(  u,v\right)  $.
Hence, we cannot have $u^{\prime}=u$. Thus, $u^{\prime}\neq u$.

We now claim that $u^{\prime}<u$. Indeed, assume the contrary. Thus,
$u^{\prime}\geq u$, so that $u^{\prime}>u$ (since $u^{\prime}\neq u$). Hence,
the cell $\left(  u^{\prime},v\right)  $ lies strictly south of the cell
$\left(  u,v\right)  $ (and we know that both cells belong to $D$). Since the
entries of $T$ increase top-to-bottom down each column (because $T$ is
column-standard), we thus conclude that $T\left(  u^{\prime},v\right)
>T\left(  u,v\right)  =k$. Hence, (\ref{pf.lem.spechtmod.lead-term.a.1})
(applied to $i=T\left(  u^{\prime},v\right)  $) yields $w\left(  T\left(
u^{\prime},v\right)  \right)  =T\left(  u^{\prime},v\right)  $. In view of
$T\left(  u^{\prime},v\right)  =w^{-1}\left(  k\right)  $, we can rewrite this
as $w\left(  w^{-1}\left(  k\right)  \right)  =w^{-1}\left(  k\right)  $.
Hence, $w^{-1}\left(  k\right)  =w\left(  w^{-1}\left(  k\right)  \right)
=k$. But this contradicts $w^{-1}\left(  k\right)  \neq k$.

This contradiction shows that our assumption was false. Hence, $u^{\prime}<u$.
In other words,
\[
r_{\overline{wT}}\left(  k\right)  <r_{\overline{T}}\left(  k\right)
\]
(since $r_{\overline{T}}\left(  k\right)  =u$ and $r_{\overline{wT}}\left(
k\right)  =u^{\prime}$). Thus, $r_{\overline{wT}}\left(  k\right)  \neq
r_{\overline{T}}\left(  k\right)  $.

Moreover, using (\ref{pf.lem.spechtmod.lead-term.a.1}), we can easily see that
each $i\in\left[  n\right]  $ satisfying $i>k$ satisfies $r_{\overline{wT}%
}\left(  i\right)  =r_{\overline{T}}\left(  i\right)  $%
\ \ \ \ \footnote{\textit{Proof.} Let $i\in\left[  n\right]  $ be such that
$i>k$. We must prove that $r_{\overline{wT}}\left(  i\right)  =r_{\overline
{T}}\left(  i\right)  $.
\par
Indeed, let $\left(  x,y\right)  $ be the cell of $T$ that contains the entry
$i$. Thus, $T\left(  x,y\right)  =i$ and $r_{\overline{T}}\left(  i\right)
=x$. But $wT=w\circ T$, so that $\left(  wT\right)  \left(  x,y\right)
=w\left(  \underbrace{T\left(  x,y\right)  }_{=i}\right)  =w\left(  i\right)
=i$ (by (\ref{pf.lem.spechtmod.lead-term.a.1})). Thus, $\left(  x,y\right)  $
is the cell of $wT$ that contains $i$. Hence, $r_{\overline{wT}}\left(
i\right)  =x$. Comparing this with $r_{\overline{T}}\left(  i\right)  =x$, we
obtain $r_{\overline{wT}}\left(  i\right)  =r_{\overline{T}}\left(  i\right)
$, qed.}.

Thus, we conclude that $k$ is a number $i\in\left[  n\right]  $ such that
$r_{\overline{wT}}\left(  i\right)  \neq r_{\overline{T}}\left(  i\right)  $
(since $r_{\overline{wT}}\left(  k\right)  \neq r_{\overline{T}}\left(
k\right)  $), but there exists no number $i\in\left[  n\right]  $ higher than
$k$ that satisfies $r_{\overline{wT}}\left(  i\right)  \neq r_{\overline{T}%
}\left(  i\right)  $ (since each $i\in\left[  n\right]  $ satisfying $i>k$
satisfies $r_{\overline{wT}}\left(  i\right)  =r_{\overline{T}}\left(
i\right)  $). Hence, $k$ is the \textbf{largest} $i\in\left[  n\right]  $ such
that $r_{\overline{wT}}\left(  i\right)  \neq r_{\overline{T}}\left(
i\right)  $. Thus, we have shown that

\begin{itemize}
\item there exists at least one $i\in\left[  n\right]  $ such that
$r_{\overline{wT}}\left(  i\right)  \neq r_{\overline{T}}\left(  i\right)  $
(namely, $k$), and

\item the \textbf{largest} such $i$ satisfies $r_{\overline{wT}}\left(
i\right)  <r_{\overline{T}}\left(  i\right)  $ (since this largest $i$ is $k$,
and we know that $r_{\overline{wT}}\left(  k\right)  <r_{\overline{T}}\left(
k\right)  $).
\end{itemize}

In other words, $\overline{wT}<\overline{T}$ (by the definition of the
relation $<$). Lemma \ref{lem.spechtmod.lead-term} \textbf{(a)} is thus
proved. \medskip
\end{fineprint}

\textbf{(b)} The definition of $\mathbf{e}_{T}$ yields%
\begin{align*}
\mathbf{e}_{T}  &  =\sum_{w\in\mathcal{C}\left(  T\right)  }\left(  -1\right)
^{w}\overline{wT}\\
&  =\underbrace{\left(  -1\right)  ^{\operatorname*{id}}}_{=1}%
\underbrace{\overline{\operatorname*{id}T}}_{=\overline{T}}+\sum
_{\substack{w\in\mathcal{C}\left(  T\right)  ;\\w\neq\operatorname*{id}%
}}\left(  -1\right)  ^{w}\overline{wT}\\
&  \ \ \ \ \ \ \ \ \ \ \ \ \ \ \ \ \ \ \ \ \left(
\begin{array}
[c]{c}%
\text{here, we have split off the addend for }w=\operatorname*{id}\\
\text{from the sum, since }\operatorname*{id}\in\mathcal{C}\left(  T\right)
\end{array}
\right) \\
&  =\overline{T}+\underbrace{\sum_{\substack{w\in\mathcal{C}\left(  T\right)
;\\w\neq\operatorname*{id}}}\left(  -1\right)  ^{w}\overline{wT}%
}_{\substack{=\left(  \text{a linear combination of }n\text{-tabloids
}\overline{S}\text{ with }\overline{S}<\overline{T}\right)  \\\text{(since
Lemma \ref{lem.spechtmod.lead-term} \textbf{(a)} yields that }\overline
{wT}<\overline{T}\\\text{for each }w\in\mathcal{C}\left(  T\right)  \text{
satisfying }w\neq\operatorname*{id}\text{)}}}\\
&  =\overline{T}+\left(  \text{a linear combination of }n\text{-tabloids
}\overline{S}\text{ with }\overline{S}<\overline{T}\right)  .
\end{align*}
The proof of Lemma \ref{lem.spechtmod.lead-term} \textbf{(b)} is thus complete.
\end{proof}

\subsubsection{Proof of the linear independence}

We can now prove Theorem \ref{thm.spechtmod.linind}:

\begin{proof}
[Proof of Theorem \ref{thm.spechtmod.linind}.]We must prove that the standard
polytabloids $\mathbf{e}_{T}$ are $\mathbf{k}$-linearly independent.

Assume the contrary. Thus, they are $\mathbf{k}$-linearly dependent. That is,
there exists a nontrivial relation%
\begin{equation}
\sum_{T\in U}\alpha_{T}\mathbf{e}_{T}=0 \label{pf.thm.spechtmod.linind.2}%
\end{equation}
for some nonempty set $U$ of standard $n$-tableaux of shape $D$ and some
nonzero coefficients $\alpha_{T}\in\mathbf{k}$. Consider this $U$ and these
coefficients $\alpha_{T}$.

Proposition \ref{prop.tabloid.row-st} says that there is a bijection%
\begin{align*}
&  \text{from }\left\{  \text{row-standard }n\text{-tableaux of shape
}D\right\} \\
&  \text{to }\left\{  n\text{-tabloids of shape }D\right\}
\end{align*}
that sends each row-standard $n$-tableau $T$ to its tabloid $\overline{T}$. In
particular, this bijection is injective. In other words, distinct row-standard
$n$-tableaux $T$ of shape $D$ produce distinct $n$-tabloids $\overline{T}$.
Hence, distinct standard $n$-tableaux $T$ produce distinct $n$-tabloids
$\overline{T}$ (since any standard $n$-tableau is row-standard). Therefore,
among the standard $n$-tableaux $T\in U$, there is a unique one for which the
$n$-tabloid $\overline{T}$ is maximum with respect to the Young last letter
order (since this order is a total order). Let $T_{\max}$ be this $n$-tableau.
Thus,%
\begin{equation}
\overline{T}<\overline{T_{\max}}\ \ \ \ \ \ \ \ \ \ \text{for each }T\in
U\text{ satisfying }T\neq T_{\max}. \label{pf.thm.spechtmod.linind.max}%
\end{equation}

Now, recall that the family $\left(  \overline{T}\right)  _{\overline{T}\text{
is an }n\text{-tabloid of shape }D}$ is a basis of the $\mathbf{k}$-module
$\mathcal{M}^{D}$. Thus, each element $\mathbf{a}\in\mathcal{M}^{D}$ can be
uniquely expressed as a $\mathbf{k}$-linear combination $\sum_{\overline
{S}\text{ is an }n\text{-tabloid of shape }D}a_{\overline{S}}\overline{S}$ of
the $n$-tabloids of shape $D$. We shall denote the coefficients $a_{\overline
{S}}$ appearing in this $\mathbf{k}$-linear combination by $\left[
\overline{S}\right]  \mathbf{a}$. In other words, we let $\left[  \overline
{S}\right]  \mathbf{a}$ denote the $\overline{S}$-coordinate of an element
$\mathbf{a}\in\mathcal{M}^{D}$ with respect to the basis $\left(  \overline
{T}\right)  _{\overline{T}\text{ is an }n\text{-tabloid of shape }D}$.

Now, we claim the following:

\begin{statement}
\textit{Claim 1:} Let $T\in U$ satisfy $T\neq T_{\max}$. Then, $\left[
\overline{T_{\max}}\right]  \left(  \mathbf{e}_{T}\right)  =0$.
\end{statement}

\begin{statement}
\textit{Claim 2:} We have $\left[  \overline{T_{\max}}\right]  \left(
\mathbf{e}_{T_{\max}}\right)  =1$.
\end{statement}

Both of these claims follow easily from Lemma \ref{lem.spechtmod.lead-term}
\textbf{(b)} and (\ref{pf.thm.spechtmod.linind.max}). Namely:

\begin{fineprint}
\begin{proof}
[Proof of Claim 1.]From (\ref{pf.thm.spechtmod.linind.max}), we obtain
$\overline{T}<\overline{T_{\max}}$, so that $\overline{T}\neq\overline
{T_{\max}}$ and therefore $\left[  \overline{T_{\max}}\right]  \overline{T}%
=0$. Moreover, $T$ is a standard $n$-tableau (since $T\in U$), thus a
column-standard $n$-tableau. Hence, Lemma \ref{lem.spechtmod.lead-term}
\textbf{(b)} yields%
\[
\mathbf{e}_{T}=\overline{T}+\left(  \text{a linear combination of
}n\text{-tabloids }\overline{S}\text{ with }\overline{S}<\overline{T}\right)
.
\]
In other words,
\begin{equation}
\mathbf{e}_{T}=\overline{T}+\sum_{\substack{\overline{S}\text{ is an
}n\text{-tabloid}\\\text{with }\overline{S}<\overline{T}}}\gamma_{\overline
{S}}\overline{S} \label{pf.thm.spechtmod.linind.c1.pf.1}%
\end{equation}
for some coefficients $\gamma_{\overline{S}}\in\mathbf{k}$. Using these
coefficients $\gamma_{\overline{S}}$, we now obtain%
\begin{align*}
\left[  \overline{T_{\max}}\right]  \left(  \mathbf{e}_{T}\right)   &
=\left[  \overline{T_{\max}}\right]  \left(  \overline{T}+\sum
_{\substack{\overline{S}\text{ is an }n\text{-tabloid}\\\text{with }%
\overline{S}<\overline{T}}}\gamma_{\overline{S}}\overline{S}\right)
\ \ \ \ \ \ \ \ \ \ \left(  \text{by (\ref{pf.thm.spechtmod.linind.c1.pf.1}%
)}\right) \\
&  =\underbrace{\left[  \overline{T_{\max}}\right]  \overline{T}}_{=0}%
+\sum_{\substack{\overline{S}\text{ is an }n\text{-tabloid}\\\text{with
}\overline{S}<\overline{T}}}\gamma_{\overline{S}}\underbrace{\left[
\overline{T_{\max}}\right]  \overline{S}}_{\substack{=0\\\text{(since
}\overline{S}<\overline{T}<\overline{T_{\max}}\\\text{and thus }\overline
{S}\neq\overline{T_{\max}}\text{)}}}=0+\underbrace{\sum_{\substack{\overline
{S}\text{ is an }n\text{-tabloid}\\\text{with }\overline{S}<\overline{T}%
}}\gamma_{\overline{S}}0}_{=0}=0.
\end{align*}
This proves Claim 1.
\end{proof}

\begin{proof}
[Proof of Claim 2.]We know that $T_{\max}$ is a standard $n$-tableau (since
$T_{\max}\in U$), thus a column-standard $n$-tableau. Hence, Lemma
\ref{lem.spechtmod.lead-term} \textbf{(b)} (applied to $T=T_{\max}$) yields%
\[
\mathbf{e}_{T_{\max}}=\overline{T_{\max}}+\left(  \text{a linear combination
of }n\text{-tabloids }\overline{S}\text{ with }\overline{S}<\overline{T_{\max
}}\right)  .
\]
In other words,
\begin{equation}
\mathbf{e}_{T_{\max}}=\overline{T_{\max}}+\sum_{\substack{\overline{S}\text{
is an }n\text{-tabloid}\\\text{with }\overline{S}<\overline{T_{\max}}}%
}\gamma_{\overline{S}}\overline{S} \label{pf.thm.spechtmod.linind.c2.pf.1}%
\end{equation}
for some coefficients $\gamma_{\overline{S}}\in\mathbf{k}$. Using these
coefficients $\gamma_{\overline{S}}$, we now obtain%
\begin{align*}
\left[  \overline{T_{\max}}\right]  \left(  \mathbf{e}_{T_{\max}}\right)   &
=\left[  \overline{T_{\max}}\right]  \left(  \overline{T_{\max}}%
+\sum_{\substack{\overline{S}\text{ is an }n\text{-tabloid}\\\text{with
}\overline{S}<\overline{T_{\max}}}}\gamma_{\overline{S}}\overline{S}\right)
\ \ \ \ \ \ \ \ \ \ \left(  \text{by (\ref{pf.thm.spechtmod.linind.c2.pf.1}%
)}\right) \\
&  =\underbrace{\left[  \overline{T_{\max}}\right]  \overline{T_{\max}}}%
_{=1}+\sum_{\substack{\overline{S}\text{ is an }n\text{-tabloid}\\\text{with
}\overline{S}<\overline{T_{\max}}}}\gamma_{\overline{S}}\underbrace{\left[
\overline{T_{\max}}\right]  \overline{S}}_{\substack{=0\\\text{(since
}\overline{S}<\overline{T_{\max}}\\\text{and thus }\overline{S}\neq
\overline{T_{\max}}\text{)}}}=1+\underbrace{\sum_{\substack{\overline{S}\text{
is an }n\text{-tabloid}\\\text{with }\overline{S}<\overline{T_{\max}}}%
}\gamma_{\overline{S}}0}_{=0}=1.
\end{align*}
This proves Claim 2.
\end{proof}
\end{fineprint}

Now,%
\begin{align*}
\left[  \overline{T_{\max}}\right]  \left(  \sum_{T\in U}\alpha_{T}%
\mathbf{e}_{T}\right)   &  =\sum_{T\in U}\alpha_{T}\left[  \overline{T_{\max}%
}\right]  \left(  \mathbf{e}_{T}\right) \\
&  =\alpha_{T_{\max}}\underbrace{\left[  \overline{T_{\max}}\right]  \left(
\mathbf{e}_{T_{\max}}\right)  }_{\substack{=1\\\text{(by Claim 2)}}%
}+\sum_{\substack{T\in U;\\T\neq T_{\max}}}\alpha_{T}\underbrace{\left[
\overline{T_{\max}}\right]  \left(  \mathbf{e}_{T}\right)  }%
_{\substack{=0\\\text{(by Claim 1)}}}\\
&  \ \ \ \ \ \ \ \ \ \ \ \ \ \ \ \ \ \ \ \ \left(
\begin{array}
[c]{c}%
\text{here, we have split off the}\\
\text{addend for }T=T_{\max}\text{ from the sum}%
\end{array}
\right) \\
&  =\alpha_{T_{\max}}+\underbrace{\sum_{\substack{T\in U;\\T\neq T_{\max}%
}}\alpha_{T}0}_{=0}=\alpha_{T_{\max}}\neq0
\end{align*}
(since all our coefficients $\alpha_{T}$ are nonzero). But this contradicts%
\[
\left[  \overline{T_{\max}}\right]  \underbrace{\left(  \sum_{T\in U}%
\alpha_{T}\mathbf{e}_{T}\right)  }_{\substack{=0\\\text{(by
(\ref{pf.thm.spechtmod.linind.2}))}}}=\left[  \overline{T_{\max}}\right]
0=0.
\]
This contradiction shows that our assumption was false. Hence, the standard
polytabloids $\mathbf{e}_{T}$ are $\mathbf{k}$-linearly independent. This
proves Theorem \ref{thm.spechtmod.linind}.
\end{proof}

\begin{corollary}
\label{cor.spechtmod.nonzero}Let $D$ be any diagram with $\left\vert
D\right\vert =n$. Then, $\mathcal{S}^{D}\neq0$ (unless $\mathbf{k}$ is trivial).
\end{corollary}

\begin{proof}
First, we shall show that there exists at least one standard $n$-tableau $T$
of shape $D$. Indeed, we can construct such a $T$ explicitly:

The \emph{depth} of a cell $\left(  i,j\right)  \in D$ shall mean the integer
$i+j$. It is clear that the depth of a cell increases (strictly) as we move
east or south.

Now, number the $n$ cells of $D$ by $1,2,\ldots,n$ in the order of increasing
depth\footnote{There might be ties (i.e., distinct cells with equal depths),
but these ties can be broken arbitrarily. For instance, you can choose to
number cells of equal depth in the order of increasing x-coordinate (i.e.,
from northernmost to southernmost).}. Define an $n$-tableau $T$ by filling
each cell $c\in D$ with its number (i.e., we set $T\left(  c\right)  $ to be
the number that we have just assigned to the cell $c$). This $n$-tableau $T$
is row-standard (because if a cell $\left(  i,j\right)  $ lies strictly west
of a cell $\left(  i,j^{\prime}\right)  $, then $j<j^{\prime}$, so that the
cell $\left(  i,j\right)  $ has smaller depth than $\left(  i,j^{\prime
}\right)  $, and thus the cell $\left(  i,j\right)  $ has a lower number than
$\left(  i,j^{\prime}\right)  $; but this means that $T\left(  i,j\right)
<T\left(  i,j^{\prime}\right)  $) and column-standard (for similar reasons).
Thus, it is standard. Hence, we have found a standard $n$-tableau $T$ of shape
$D$.

Now consider the corresponding standard polytabloid $\mathbf{e}_{T}$. By
Theorem \ref{thm.spechtmod.linind}, this $\mathbf{e}_{T}$ is part of a
$\mathbf{k}$-linearly independent family of vectors. Hence, unless
$\mathbf{k}$ is trivial, it follows that $\mathbf{e}_{T}\neq0$, so that
$\mathcal{S}^{D}\neq0$ (since $\mathbf{e}_{T}$ belongs to $\mathcal{S}^{D}$).
This proves Corollary \ref{cor.spechtmod.nonzero}.
\end{proof}

\begin{remark}
\label{rmk.spechtmod.not-bas}In general, the standard polytabloids
$\mathbf{e}_{T}$ do \textbf{not} form a basis of $\mathcal{S}^{D}$. We will
soon see (Theorem \ref{thm.spechtmod.basis}) that they do whenever $D$ is a
skew Young diagram $Y\left(  \lambda/\mu\right)  $. But in general, there are
not enough of them. For example, if $n=3$ and
\[
D=\left\{  \left(  1,1\right)  ,\ \left(  1,2\right)  ,\ \left(  2,2\right)
\right\}  =%
%TCIMACRO{\TeXButton{tikz tromino}{\begin{tikzpicture}[scale=0.7]
%\draw[fill=red!50] (1, 0) rectangle (2, 1);
%\draw[fill=red!50] (0, 1) rectangle (1, 2);
%\draw[fill=red!50] (1, 1) rectangle (2, 2);
%\end{tikzpicture}}}%
%BeginExpansion
\begin{tikzpicture}[scale=0.7]
\draw[fill=red!50] (1, 0) rectangle (2, 1);
\draw[fill=red!50] (0, 1) rectangle (1, 2);
\draw[fill=red!50] (1, 1) rectangle (2, 2);
\end{tikzpicture}%
%EndExpansion
\ \ ,
\]
then there is only one standard tableau of shape $D$, but the Specht module
$\mathcal{S}^{D}$ is a free $\mathbf{k}$-module of rank $2$, so that a single
polytabloid $\mathbf{e}_{T}$ is not enough to span it.
\end{remark}

\subsection{The Garnir relations}

As we showed above, the standard polytabloids $\mathbf{e}_{T}$ in a Specht
module $\mathcal{S}^{D}$ are linearly independent, but don't always span
$\mathcal{S}^{D}$. Our goal is to show that they span $\mathcal{S}^{D}$
whenever $D$ is a skew Young diagram $Y\left(  \lambda/\mu\right)  $. To do
this, we need a useful construction called the \emph{Garnir relations}. This
construction underlies much of the modern study of $S_{n}$-representations;
see \cite[\S 3.4]{Welsh92}, \cite[\S 6.2]{Wildon18}, \cite[\S 7]{James78},
\cite[\S 7.2]{JamKer81}, \cite[\S 2.6]{Sagan01}, \cite[\S 11.4]{Howe22},
\cite[\S 4.5]{LakBro18} for various alternative treatments (some only handling
straight Young diagrams), and \cite[\S 3.2]{Liu16}, \cite{Welsh92},
\cite{Murphy81}, \cite[\S 3]{DipJam85} for generalizations, variants and applications.

Before we define the Garnir relations, we need to introduce some simple but
useful features of arbitrary diagrams.

\subsubsection{Notations and some harmless lemmas}

\begin{definition}
\label{def.diagram.Dj}Let $D$ be a diagram. Let $j\in\mathbb{Z}$. Then,
$D\left\langle j\right\rangle $ shall mean the set of all integers $i$ that
satisfy $\left(  i,j\right)  \in D$.
\end{definition}

For example, if
\[
D=\left\{  \left(  1,2\right)  ,\ \left(  1,3\right)  ,\ \left(  2,1\right)
,\ \left(  2,2\right)  ,\ \left(  3,3\right)  ,\ \left(  4,1\right)  \right\}
=%
%TCIMACRO{\TeXButton{tikz weird diagram}{\begin{tikzpicture}[scale=0.7]
%\draw[fill=red!50] (0, 0) rectangle (1, 1);
%\draw[fill=red!50] (2, 1) rectangle (3, 2);
%\draw[fill=red!50] (0, 2) rectangle (1, 3);
%\draw[fill=red!50] (1, 2) rectangle (2, 3);
%\draw[fill=red!50] (1, 3) rectangle (2, 4);
%\draw[fill=red!50] (2, 3) rectangle (3, 4);
%\end{tikzpicture}}}%
%BeginExpansion
\begin{tikzpicture}[scale=0.7]
\draw[fill=red!50] (0, 0) rectangle (1, 1);
\draw[fill=red!50] (2, 1) rectangle (3, 2);
\draw[fill=red!50] (0, 2) rectangle (1, 3);
\draw[fill=red!50] (1, 2) rectangle (2, 3);
\draw[fill=red!50] (1, 3) rectangle (2, 4);
\draw[fill=red!50] (2, 3) rectangle (3, 4);
\end{tikzpicture}%
%EndExpansion
\ \ ,
\]
then $D\left\langle 1\right\rangle =\left\{  2,4\right\}  $ and $D\left\langle
2\right\rangle =\left\{  1,2\right\}  $ and $D\left\langle 3\right\rangle
=\left\{  1,3\right\}  $ and $D\left\langle j\right\rangle =\varnothing$ for
all $j\notin\left[  3\right]  $. If $D$ is a skew Young diagram $Y\left(
\lambda/\mu\right)  $, then each set $D\left\langle j\right\rangle $ is an
integer interval (equalling $\left\{  \mu_{j}^{t}+1,\ \mu_{j}^{t}%
+2,\ \ldots,\ \lambda_{j}^{t}\right\}  $ when $j>0$ and being empty
otherwise), but of course this does not generally hold for bad-shape diagrams.

\begin{definition}
\label{def.tableau.ColTj}Let $T$ be an $n$-tableau (of any shape). Let
$j\in\mathbb{Z}$. Then, $\operatorname*{Col}\left(  j,T\right)  $ shall denote
the set of all entries in the $j$-th column of $T$.
\end{definition}

For example, if $T=\ytableaushort{214,65}$\ , then $\operatorname*{Col}\left(
1,T\right)  =\left\{  2,6\right\}  $ and $\operatorname*{Col}\left(
2,T\right)  =\left\{  1,5\right\}  $ and $\operatorname*{Col}\left(
3,T\right)  =\left\{  4\right\}  $ and $\operatorname*{Col}\left(  j,T\right)
=\varnothing$ for all $j\notin\left[  3\right]  $.

The following should come as no surprise:

\begin{proposition}
\label{prop.tableau.ColTj}Let $D$ be a diagram. Let $T$ be an $n$-tableau of
shape $D$. Let $j\in\mathbb{Z}$. Then, $\left\vert \operatorname*{Col}\left(
j,T\right)  \right\vert =\left\vert D\left\langle j\right\rangle \right\vert $.
\end{proposition}

\begin{fineprint}
\begin{proof}
We know that $T$ is an $n$-tableau, thus injective. Hence, all entries of $T$
are distinct. In particular, all entries in the $j$-th column of $T$ are distinct.

But $D\left\langle j\right\rangle $ is defined as the set of all integers $i$
that satisfy $\left(  i,j\right)  \in D$. Thus, the cells in the $j$-th column
of $D$ are precisely the cells $\left(  i,j\right)  $ with $i\in D\left\langle
j\right\rangle $. Hence, for any $i\in D\left\langle j\right\rangle $, the
cell $\left(  i,j\right)  $ belongs to the $j$-th column of $D$, and therefore
the corresponding entry $T\left(  i,j\right)  $ of $T$ is an entry of the
$j$-th column of $T$ and thus belongs to $\operatorname*{Col}\left(
j,T\right)  $ (by the definition of $\operatorname*{Col}\left(  j,T\right)
$). Therefore, the map%
\begin{align*}
D\left\langle j\right\rangle  &  \rightarrow\operatorname*{Col}\left(
j,T\right)  ,\\
i  &  \mapsto T\left(  i,j\right)
\end{align*}
is well-defined. This map is furthermore injective (since all entries in the
$j$-th column of $T$ are distinct) and surjective (since each element of
$\operatorname*{Col}\left(  j,T\right)  $ is an entry in the $j$-th column of
$T$, and thus has the form $T\left(  i,j\right)  $ for some $i\in
D\left\langle j\right\rangle $). Thus, this map is bijective, i.e., a
bijection. Hence, we have found a bijection from $D\left\langle j\right\rangle
$ to $\operatorname*{Col}\left(  j,T\right)  $. We thus conclude (by the
bijection principle) that $\left\vert \operatorname*{Col}\left(  j,T\right)
\right\vert =\left\vert D\left\langle j\right\rangle \right\vert $. This
proves Proposition \ref{prop.tableau.ColTj}.
\end{proof}
\end{fineprint}

Somewhat (but not much) less obvious is the following equality in a Young module:

\begin{lemma}
\label{lem.garnir.S}Let $D$ be a diagram. Let $T$ be an $n$-tableau of shape
$D$. Let $j$ and $k$ be two integers. Let $Z$ be a subset of
$\operatorname*{Col}\left(  j,T\right)  \cup\operatorname*{Col}\left(
k,T\right)  $ such that $\left\vert Z\right\vert >\left\vert D\left\langle
j\right\rangle \cup D\left\langle k\right\rangle \right\vert $. Then, in the
Young module $\mathcal{M}^{D}$, we have%
\[
\nabla_{Z}^{-}\overline{T}=0.
\]

\end{lemma}

\begin{example}
\label{exa.garnir.S.1}Let $n=8$ and%
\begin{align*}
D  &  =\left\{  \left(  1,1\right)  ,\ \left(  1,3\right)  ,\ \left(
2,1\right)  ,\ \left(  2,2\right)  ,\ \left(  2,3\right)  ,\ \left(
3,1\right)  ,\ \left(  3,3\right)  ,\ \left(  4,1\right)  \right\} \\
&  =%
%TCIMACRO{\TeXButton{tikz H diagram}{\begin{tikzpicture}[scale=0.7]
%\draw[fill=red!50] (0, -1) rectangle (1, 0);
%\draw[fill=red!50] (0, 0) rectangle (1, 1);
%\draw[fill=red!50] (2, 0) rectangle (3, 1);
%\draw[fill=red!50] (0, 1) rectangle (1, 2);
%\draw[fill=red!50] (1, 1) rectangle (2, 2);
%\draw[fill=red!50] (2, 1) rectangle (3, 2);
%\draw[fill=red!50] (0, 2) rectangle (1, 3);
%\draw[fill=red!50] (2, 2) rectangle (3, 3);
%\end{tikzpicture}}}%
%BeginExpansion
\begin{tikzpicture}[scale=0.7]
\draw[fill=red!50] (0, -1) rectangle (1, 0);
\draw[fill=red!50] (0, 0) rectangle (1, 1);
\draw[fill=red!50] (2, 0) rectangle (3, 1);
\draw[fill=red!50] (0, 1) rectangle (1, 2);
\draw[fill=red!50] (1, 1) rectangle (2, 2);
\draw[fill=red!50] (2, 1) rectangle (3, 2);
\draw[fill=red!50] (0, 2) rectangle (1, 3);
\draw[fill=red!50] (2, 2) rectangle (3, 3);
\end{tikzpicture}%
%EndExpansion
\ \ .
\end{align*}
Let $T$ be the $n$-tableau $\ytableaushort{3\none6,154,7\none2,8}$ of shape
$D$. Let $Z=\left\{  1,2,6,7,8\right\}  $. Then, $Z$ is a subset of
$\operatorname*{Col}\left(  1,T\right)  \cup\operatorname*{Col}\left(
3,T\right)  $ and satisfies $\left\vert Z\right\vert >\left\vert D\left\langle
1\right\rangle \cup D\left\langle 3\right\rangle \right\vert $ (since
$\left\vert Z\right\vert =5$ and $\left\vert D\left\langle 1\right\rangle \cup
D\left\langle 3\right\rangle \right\vert =4$). Hence, Lemma \ref{lem.garnir.S}
yields $\nabla_{Z}^{-}\overline{T}=0$. If we color the cells of $T$ that
contain the entries in $Z$ green, then $T$ looks as follows:%
\[
T=\ytableaushort{3\none{*(green)6},{*(green)1}54,{*(green)7}\none{*(green)2},{*(green)8}}\ \ .
\]

\end{example}

\begin{proof}
[Proof of Lemma \ref{lem.garnir.S}.]We note that $\left\vert D\right\vert =n$
(since $T$ is an $n$-tableau of shape $D$).

For each $m\in\left[  n\right]  $, we let $r_{T}\left(  m\right)  $ be the
number of the row of $T$ that contains the entry $m$.

We have $Z\subseteq\operatorname*{Col}\left(  j,T\right)  \cup
\operatorname*{Col}\left(  k,T\right)  $ (by the definition of $Z$). Thus,
each $m\in Z$ satisfies $r_{T}\left(  m\right)  \in D\left\langle
j\right\rangle \cup D\left\langle k\right\rangle $ (this is easy to
see\footnote{\textit{Proof.} Let $m\in Z$. We must show that $r_{T}\left(
m\right)  \in D\left\langle j\right\rangle \cup D\left\langle k\right\rangle
$.
\par
We have $m\in Z\subseteq\operatorname*{Col}\left(  j,T\right)  \cup
\operatorname*{Col}\left(  k,T\right)  $. Thus, we have $m\in
\operatorname*{Col}\left(  j,T\right)  $ or $m\in\operatorname*{Col}\left(
k,T\right)  $. We WLOG assume that $m\in\operatorname*{Col}\left(  j,T\right)
$ (since the other case is analogous). Thus, $m$ is an entry of the $j$-th
column of $T$ (by the definition of $\operatorname*{Col}\left(  j,T\right)
$). In other words, $m=T\left(  i,j\right)  $ for some $i\in\mathbb{Z}$
satisfying $\left(  i,j\right)  \in D$. Consider this $i$. From $\left(
i,j\right)  \in D$, we obtain $i\in D\left\langle j\right\rangle $ (by the
definition of $D\left\langle j\right\rangle $). But $m=T\left(  i,j\right)  $
shows that the entry $m$ appears in the cell $\left(  i,j\right)  $ of $T$.
Hence, the entry $m$ appears in the $i$-th row of $T$. In other words,
$r_{T}\left(  m\right)  =i$ (by the definition of $r_{T}\left(  m\right)  $).
Thus, $r_{T}\left(  m\right)  =i\in D\left\langle j\right\rangle \subseteq
D\left\langle j\right\rangle \cup D\left\langle k\right\rangle $, qed.}).
Hence, the map%
\begin{align*}
Z  &  \rightarrow D\left\langle j\right\rangle \cup D\left\langle
k\right\rangle ,\\
m  &  \mapsto r_{T}\left(  m\right)
\end{align*}
is well-defined. By the pigeonhole principle, this map cannot be injective
(since its domain $Z$ is larger than its target $D\left\langle j\right\rangle
\cup D\left\langle k\right\rangle $\ \ \ \ \footnote{In more detail: One form
of the pigeonhole principle says that if $X$ and $Y$ are two finite sets
satisfying $\left\vert X\right\vert >\left\vert Y\right\vert $, then there
exists no injective map from $X$ to $Y$. Applying this to $X=Z$ and
$Y=D\left\langle j\right\rangle \cup D\left\langle k\right\rangle $, we see
that there exists no injective map from $Z$ to $D\left\langle j\right\rangle
\cup D\left\langle k\right\rangle $ (since $\left\vert Z\right\vert
>\left\vert D\left\langle j\right\rangle \cup D\left\langle k\right\rangle
\right\vert $). Hence, in particular, the map%
\begin{align*}
Z  &  \rightarrow D\left\langle j\right\rangle \cup D\left\langle
k\right\rangle ,\\
m  &  \mapsto r_{T}\left(  m\right)
\end{align*}
cannot be injective.}). Thus, there exist two distinct elements $u,v\in Z$
satisfying $r_{T}\left(  u\right)  =r_{T}\left(  v\right)  $. Consider these
$u$ and $v$. Thus, the entries $u$ and $v$ lie in the same row of the tableau
$T$\ \ \ \ \footnote{\textit{Proof.} Recall that $r_{T}\left(  u\right)  $ is
defined as the number of the row of $T$ that contains the entry $u$, whereas
$r_{T}\left(  v\right)  $ is defined as the number of the row of $T$ that
contains the entry $v$. Hence, the equality $r_{T}\left(  u\right)
=r_{T}\left(  v\right)  $ (which we know to be true) is saying that the
entries $u$ and $v$ lie in the same row of the tableau $T$.}. Hence, the
permutation $t_{u,v}$ is horizontal for $T$\ \ \ \ \footnote{\textit{Proof.}
We need to show that for each $i\in\left[  n\right]  $, the number
$t_{u,v}\left(  i\right)  $ lies in the same row of $T$ as $i$ does. But this
is true for $i=u$ (since $t_{u,v}\left(  u\right)  =v$ lies in the same row of
$T$ as $u$ does) and true for $i=v$ (for similar reasons), and also true for
all remaining values of $i$ (since $t_{u,v}\left(  i\right)  =i$ whenever $i$
is neither $u$ nor $v$). Thus, we have shown that the permutation $t_{u,v}$ is
horizontal for $T$.}. In other words, $t_{u,v}\in\mathcal{R}\left(  T\right)
$. Hence, Proposition \ref{prop.tableau.Sn-act.0} \textbf{(a)} (applied to
$w=t_{u,v}$) shows that the tableau $t_{u,v}\rightharpoonup T$ is
row-equivalent to $T$. In other words, $\overline{t_{u,v}\rightharpoonup
T}=\overline{T}$ (since the $n$-tabloids are the equivalence classes of the
$n$-tableaux with respect to row-equivalence). Thus, in $\mathcal{M}^{D}$, we
have%
\[
\left(  1-t_{u,v}\right)  \overline{T}=\overline{T}-\underbrace{t_{u,v}%
\overline{T}}_{=\overline{t_{u,v}T}=\overline{t_{u,v}\rightharpoonup
T}=\overline{T}}=\overline{T}-\overline{T}=0.
\]

Let us set%
\[
\nabla_{Z}^{\operatorname*{even}}:=\sum_{\substack{w\in S_{n,Z};\\\left(
-1\right)  ^{w}=1}}w
\]
(where $S_{n,Z}$ is defined according to Proposition \ref{prop.intX.basics}).

However, $u$ and $v$ are two distinct elements of $Z$. Thus, the equality
(\ref{prop.intX.basics.d.2}) from Proposition \ref{prop.intX.basics}
\textbf{(d)} (applied to $Z$, $u$ and $v$ instead of $X$, $i$ and $j$) yields
that%
\[
\nabla_{Z}^{-}=\nabla_{Z}^{\operatorname*{even}}\cdot\left(  1-t_{u,v}\right)
=\left(  1-t_{u,v}\right)  \cdot\nabla_{Z}^{\operatorname*{even}}.
\]
Hence,%
\[
\underbrace{\nabla_{Z}^{-}}_{=\nabla_{Z}^{\operatorname*{even}}\cdot\left(
1-t_{u,v}\right)  }\overline{T}=\nabla_{Z}^{\operatorname*{even}}%
\cdot\underbrace{\left(  1-t_{u,v}\right)  \overline{T}}_{=0}=0.
\]
This proves Lemma \ref{lem.garnir.S}.
\end{proof}

Lemma \ref{lem.garnir.S} easily yields the following:

\begin{lemma}
\label{lem.garnir.ST}Let $D$ be a diagram. Let $S$ and $T$ be two
column-equivalent $n$-tableaux of shape $D$. Let $j$ and $k$ be two integers.
Let $Z$ be a subset of $\operatorname*{Col}\left(  j,T\right)  \cup
\operatorname*{Col}\left(  k,T\right)  $ such that $\left\vert Z\right\vert
>\left\vert D\left\langle j\right\rangle \cup D\left\langle k\right\rangle
\right\vert $. Then, in the Young module $\mathcal{M}^{D}$, we have%
\[
\nabla_{Z}^{-}\overline{S}=0.
\]

\end{lemma}

\begin{proof}
The tableau $T$ is column-equivalent to $S$. Thus, it contains the same
entries in its $j$-th column as $S$ does. In other words, $\operatorname*{Col}%
\left(  j,T\right)  =\operatorname*{Col}\left(  j,S\right)  $. Similarly,
$\operatorname*{Col}\left(  k,T\right)  =\operatorname*{Col}\left(
k,S\right)  $.

But we assumed that $Z$ is a subset of $\operatorname*{Col}\left(  j,T\right)
\cup\operatorname*{Col}\left(  k,T\right)  $. In other words, $Z$ is a subset
of $\operatorname*{Col}\left(  j,S\right)  \cup\operatorname*{Col}\left(
k,S\right)  $ (since $\operatorname*{Col}\left(  j,T\right)
=\operatorname*{Col}\left(  j,S\right)  $ and $\operatorname*{Col}\left(
k,T\right)  =\operatorname*{Col}\left(  k,S\right)  $). Hence, Lemma
\ref{lem.garnir.S} (applied to $S$ instead of $T$) yields $\nabla_{Z}%
^{-}\overline{S}=0$. This proves Lemma \ref{lem.garnir.ST}.
\end{proof}

In order to apply Lemma \ref{lem.garnir.ST}, we need to find a subset $Z$ of
$\operatorname*{Col}\left(  j,T\right)  \cup\operatorname*{Col}\left(
k,T\right)  $ that satisfies the condition $\left\vert Z\right\vert
>\left\vert D\left\langle j\right\rangle \cup D\left\langle k\right\rangle
\right\vert $. It is easy to see that such a subset $Z$ can only exist if the
sets $D\left\langle j\right\rangle $ and $D\left\langle k\right\rangle $ have
a nonempty intersection -- i.e., if $D$ has a row that has both a cell in the
$j$-th column and a cell in the $k$-th column.

The following lemma will be particularly useful for us as we will apply Lemma
\ref{lem.garnir.ST} to skew Young tableaux:

\begin{lemma}
\label{lem.garnir.skew}Let $D$ be the skew Young diagram $Y\left(  \lambda
/\mu\right)  $ of a skew partition $\lambda/\mu$. Let $T$ be an $n$-tableau of
shape $D$. Let $j$ and $k$ be two integers with $j<k$. Let $i\in\mathbb{Z}$ be
such that $\left(  i,j\right)  \in D$ and $\left(  i,k\right)  \in D$. Define
two sets%
\begin{align*}
X  &  =\left\{  T\left(  i^{\prime},j\right)  \ \mid\ i^{\prime}\in
\mathbb{Z}\text{ with }i^{\prime}\geq i\text{ and }\left(  i^{\prime
},j\right)  \in D\right\}  \ \ \ \ \ \ \ \ \ \ \text{and}\\
Y  &  =\left\{  T\left(  i^{\prime},k\right)  \ \mid\ i^{\prime}\in
\mathbb{Z}\text{ with }i^{\prime}\leq i\text{ and }\left(  i^{\prime
},k\right)  \in D\right\}  .
\end{align*}
(That is, $X$ is the set of all entries in the $j$-th column of $T$ at the
cell $\left(  i,j\right)  $ and further south, whereas $Y$ is the set of all
entries in the $k$-th column of $T$ at the cell $\left(  i,k\right)  $ and
further north.) Then: \medskip

\textbf{(a)} We have $X\subseteq\operatorname*{Col}\left(  j,T\right)  $ and
$Y\subseteq\operatorname*{Col}\left(  k,T\right)  $ and $\left\vert X\cup
Y\right\vert >\left\vert D\left\langle j\right\rangle \cup D\left\langle
k\right\rangle \right\vert $. \medskip

\textbf{(b)} If the tableau $T$ is column-standard and satisfies $T\left(
i,j\right)  >T\left(  i,k\right)  $, then any element of $X$ is larger than
any element of $Y$.
\end{lemma}

\Needspace{20pc}

\begin{example}
Let $n=10$ and
\[
D=Y\left(  \left(  3,3,3,3,1\right)  /\left(  2,1\right)  \right)
=\ydiagram{2+1,1+2,0+3,0+3,0+1}\ \ ,
\]
and let $T$ be the $n$-tableau
$\ytableaushort{\none\none4,\none16,258,379,{10}}$\ . Let $j=2$ and $k=3$ and
$i=3$. Then, the conditions of Lemma \ref{lem.garnir.skew} are satisfied, and
the sets $X$ and $Y$ in Lemma \ref{lem.garnir.skew} are $X=\left\{
5,7\right\}  $ and $Y=\left\{  4,6,8\right\}  $. If we color the cells of $T$
that contain elements in $X$ green, and color the cells of $T$ that contain
elements of $Y$ red, then $T$ looks as follows:%
\[
T=\ytableaushort{\none\none{*(red)4},\none1{*(red)6},2{*(green)5}{*(red)8},3{*(green)7}9,{10}}\ \ .
\]

However, Lemma \ref{lem.garnir.skew} \textbf{(b)} cannot be applied here,
since the condition $T\left(  i,j\right)  >T\left(  i,k\right)  $ is not satisfied.
\end{example}

\begin{fineprint}
\begin{proof}
[Proof of Lemma \ref{lem.garnir.skew}.]\textbf{(a)} The definition of $X$
makes it clear that each element of $X$ is an entry $T\left(  i^{\prime
},j\right)  $ in the $j$-th column of $T$, thus an element of
$\operatorname*{Col}\left(  j,T\right)  $ (by the definition of
$\operatorname*{Col}\left(  j,T\right)  $). In other words, $X\subseteq
\operatorname*{Col}\left(  j,T\right)  $. Similarly, $Y\subseteq
\operatorname*{Col}\left(  k,T\right)  $.

Set $Z:=X\cup Y$. Then, $Z$ is a subset of $\operatorname*{Col}\left(
j,T\right)  \cup\operatorname*{Col}\left(  k,T\right)  $ (since we have
$Z=\underbrace{X}_{\subseteq\operatorname*{Col}\left(  j,T\right)  }%
\cup\underbrace{Y}_{\subseteq\operatorname*{Col}\left(  k,T\right)  }%
\subseteq\operatorname*{Col}\left(  j,T\right)  \cup\operatorname*{Col}\left(
k,T\right)  $) and satisfies $X\subseteq X\cup Y=Z$ and $Y\subseteq X\cup Y=Z$.

The definition of $X$ yields $T\left(  i,j\right)  \in X$ (since
$i\in\mathbb{Z}$ and $i\geq i$ and $\left(  i,j\right)  \in D$). Similarly,
$T\left(  i,k\right)  \in Y$. Thus, $T\left(  i,j\right)  $ and $T\left(
i,k\right)  $ are two elements of $Z$ (since $T\left(  i,j\right)  \in
X\subseteq Z$ and $T\left(  i,k\right)  \in Y\subseteq Z$).

We know that $T$ is an $n$-tableau. Thus, $T$ is injective, i.e., all entries
of $T$ are distinct. Thus, in particular, $T\left(  i,j\right)  \neq T\left(
i,k\right)  $ (since $\left(  i,j\right)  \neq\left(  i,k\right)  $ (because
$j<k$)). Consequently, $T\left(  i,j\right)  $ and $T\left(  i,k\right)  $ are
two distinct elements of $Z$.

For each $m\in\left[  n\right]  $, we let $r_{T}\left(  m\right)  $ be the
number of the row of $T$ that contains the entry $m$. In particular,
$r_{T}\left(  T\left(  i,j\right)  \right)  =i$ (since the entry $T\left(
i,j\right)  $ obviously lies in the $i$-th row of $T$) and similarly
$r_{T}\left(  T\left(  i,k\right)  \right)  =i$. Hence, $r_{T}\left(  T\left(
i,j\right)  \right)  =i=r_{T}\left(  T\left(  i,k\right)  \right)  $.

Now we claim the following:

\begin{statement}
\textit{Claim 1:} For each $i^{\prime}\in D\left\langle j\right\rangle \cup
D\left\langle k\right\rangle $, there exists some $m\in Z$ such that
$r_{T}\left(  m\right)  =i^{\prime}$.
\end{statement}

\begin{proof}
[Proof of Claim 1.]Let $i^{\prime}\in D\left\langle j\right\rangle \cup
D\left\langle k\right\rangle $. We must prove that there exists some $m\in Z$
such that $r_{T}\left(  m\right)  =i^{\prime}$.

We are in one of the following two cases:

\textit{Case 1:} We have $i^{\prime}\geq i$.

\textit{Case 2:} We have $i^{\prime}\leq i$.

Let us consider Case 1. In this case, we have $i^{\prime}\geq i$. Thus, $i\leq
i^{\prime}$. Furthermore, we have $i^{\prime}\in D\left\langle j\right\rangle
\cup D\left\langle k\right\rangle $. This easily yields $i^{\prime}\in
D\left\langle j\right\rangle $ because $D$ is a skew Young
diagram\footnote{\textit{Proof.} We must prove that $i^{\prime}\in
D\left\langle j\right\rangle $.
\par
But we have $i^{\prime}\in D\left\langle j\right\rangle \cup D\left\langle
k\right\rangle $. Therefore, $i^{\prime}\in D\left\langle j\right\rangle $ or
$i^{\prime}\in D\left\langle k\right\rangle $. In the former case, we are
immediately done. Thus, we can WLOG assume that we are in the latter case. In
other words, $i^{\prime}\in D\left\langle k\right\rangle $. In other words,
$\left(  i^{\prime},k\right)  \in D$ (by the definition of $D\left\langle
k\right\rangle $).
\par
Now, recall the relation $\leq$ on cells introduced in Definition
\ref{def.diagrams.diagrams} \textbf{(c)}. The three cells $\left(  i,j\right)
$, $\left(  i^{\prime},j\right)  $ and $\left(  i^{\prime},k\right)  $ all
belong to $\left\{  1,2,3,\ldots\right\}  ^{2}$ and satisfy $\left(
i,j\right)  \leq\left(  i^{\prime},j\right)  $ (since $i\leq i^{\prime}$) and
$\left(  i^{\prime},j\right)  \leq\left(  i^{\prime},k\right)  $ (since $j\leq
k$) and $\left(  i,j\right)  \in D=Y\left(  \lambda/\mu\right)  $ and $\left(
i^{\prime},k\right)  \in D=Y\left(  \lambda/\mu\right)  $. Hence, Proposition
\ref{prop.young.convexity} (applied to $c=\left(  i,j\right)  $ and $d=\left(
i^{\prime},j\right)  $ and $e=\left(  i^{\prime},k\right)  $) yields $\left(
i^{\prime},j\right)  \in Y\left(  \lambda/\mu\right)  =D$. Thus, $i^{\prime
}\in D\left\langle j\right\rangle $ (by the definition of $D\left\langle
j\right\rangle $). This completes our proof.}. In other words, $\left(
i^{\prime},j\right)  \in D$ (by the definition of $D\left\langle
j\right\rangle $). Since $i^{\prime}\geq i$, this entails that $T\left(
i^{\prime},j\right)  \in X$ (by the definition of $X$) and therefore $T\left(
i^{\prime},j\right)  \in X\subseteq Z$. Of course, the entry $T\left(
i^{\prime},j\right)  $ of $T$ appears in the $i^{\prime}$-th row. In other
words, $r_{T}\left(  T\left(  i^{\prime},j\right)  \right)  =i^{\prime}$.
Hence, there exists some $m\in Z$ such that $r_{T}\left(  m\right)
=i^{\prime}$ (namely, $m=T\left(  i^{\prime},j\right)  $). This proves Claim 1
in Case 1.

Now, consider Case 2. In this case, we have $i^{\prime}\leq i$. Furthermore,
we have $i^{\prime}\in D\left\langle j\right\rangle \cup D\left\langle
k\right\rangle $. This easily yields $i^{\prime}\in D\left\langle
k\right\rangle $ because $D$ is a skew Young diagram\footnote{\textit{Proof.}
We must prove that $i^{\prime}\in D\left\langle k\right\rangle $.
\par
But we have $i^{\prime}\in D\left\langle j\right\rangle \cup D\left\langle
k\right\rangle $. Therefore, $i^{\prime}\in D\left\langle j\right\rangle $ or
$i^{\prime}\in D\left\langle k\right\rangle $. In the latter case, we are
immediately done. Thus, we can WLOG assume that we are in the former case. In
other words, $i^{\prime}\in D\left\langle j\right\rangle $. In other words,
$\left(  i^{\prime},j\right)  \in D$ (by the definition of $D\left\langle
j\right\rangle $).
\par
Now, recall the relation $\leq$ on cells introduced in Definition
\ref{def.diagrams.diagrams} \textbf{(c)}. The three cells $\left(  i^{\prime
},j\right)  $, $\left(  i^{\prime},k\right)  $ and $\left(  i,k\right)  $ all
belong to $\left\{  1,2,3,\ldots\right\}  ^{2}$ and satisfy $\left(
i^{\prime},j\right)  \leq\left(  i^{\prime},k\right)  $ (since $j\leq k$) and
$\left(  i^{\prime},k\right)  \leq\left(  i,k\right)  $ (since $i^{\prime}\leq
i$) and $\left(  i^{\prime},j\right)  \in D=Y\left(  \lambda/\mu\right)  $ and
$\left(  i,k\right)  \in D=Y\left(  \lambda/\mu\right)  $. Hence, Proposition
\ref{prop.young.convexity} (applied to $c=\left(  i^{\prime},j\right)  $ and
$d=\left(  i^{\prime},k\right)  $ and $e=\left(  i,k\right)  $) yields
$\left(  i^{\prime},k\right)  \in Y\left(  \lambda/\mu\right)  =D$. Thus,
$i^{\prime}\in D\left\langle k\right\rangle $ (by the definition of
$D\left\langle k\right\rangle $). This completes our proof.}. In other words,
$\left(  i^{\prime},k\right)  \in D$ (by the definition of $D\left\langle
k\right\rangle $). Since $i^{\prime}\leq i$, this entails that $T\left(
i^{\prime},k\right)  \in Y$ (by the definition of $Y$) and therefore $T\left(
i^{\prime},k\right)  \in Y\subseteq Z$. Of course, the entry $T\left(
i^{\prime},k\right)  $ of $T$ appears in the $i^{\prime}$-th row. In other
words, $r_{T}\left(  T\left(  i^{\prime},k\right)  \right)  =i^{\prime}$.
Hence, there exists some $m\in Z$ such that $r_{T}\left(  m\right)
=i^{\prime}$ (namely, $m=T\left(  i^{\prime},k\right)  $). This proves Claim 1
in Case 2.

We have now proved Claim 1 in both Cases 1 and 2. Thus, Claim 1 always holds.
\end{proof}

Now, recall that we have already shown that $X\subseteq\operatorname*{Col}%
\left(  j,T\right)  $ and $Y\subseteq\operatorname*{Col}\left(  k,T\right)  $.
It remains to prove that $\left\vert X\cup Y\right\vert >\left\vert
D\left\langle j\right\rangle \cup D\left\langle k\right\rangle \right\vert $.

Assume the contrary. Thus, $\left\vert X\cup Y\right\vert \leq\left\vert
D\left\langle j\right\rangle \cup D\left\langle k\right\rangle \right\vert $.
In other words, $\left\vert Z\right\vert \leq\left\vert D\left\langle
j\right\rangle \cup D\left\langle k\right\rangle \right\vert $ (since $Z=X\cup
Y$). Therefore, one form of the pigeonhole principle shows that any surjective
map from $Z$ to $D\left\langle j\right\rangle \cup D\left\langle
k\right\rangle $ is bijective\footnote{\textit{Proof.} One version of the
pigeonhole principle says that if $P$ and $Q$ are two finite sets such that
$\left\vert P\right\vert \leq\left\vert Q\right\vert $, then any surjective
map from $P$ to $Q$ is bijective. Applying this to $P=Z$ and $Q=D\left\langle
j\right\rangle \cup D\left\langle k\right\rangle $, we see that any surjective
map from $Z$ to $D\left\langle j\right\rangle \cup D\left\langle
k\right\rangle $ is bijective (since $\left\vert Z\right\vert \leq\left\vert
D\left\langle j\right\rangle \cup D\left\langle k\right\rangle \right\vert
$).}.

However, $Z$ is a subset of $\operatorname*{Col}\left(  j,T\right)
\cup\operatorname*{Col}\left(  k,T\right)  $. Thus, the map%
\begin{align*}
Z  &  \rightarrow D\left\langle j\right\rangle \cup D\left\langle
k\right\rangle ,\\
m  &  \mapsto r_{T}\left(  m\right)
\end{align*}
is well-defined (this has already been shown in the proof of Lemma
\ref{lem.garnir.S}). Claim 1 shows that this map is surjective. Hence, it is
bijective (since any surjective map from $Z$ to $D\left\langle j\right\rangle
\cup D\left\langle k\right\rangle $ is bijective), thus injective. Therefore,
it must send the two distinct elements $T\left(  i,j\right)  $ and $T\left(
i,k\right)  $ of $Z$ to two distinct values. In other words, $r_{T}\left(
T\left(  i,j\right)  \right)  \neq r_{T}\left(  T\left(  i,k\right)  \right)
$. But this contradicts $r_{T}\left(  T\left(  i,j\right)  \right)
=r_{T}\left(  T\left(  i,k\right)  \right)  $. This contradiction shows that
our assumption was false. Hence, $\left\vert X\cup Y\right\vert >\left\vert
D\left\langle j\right\rangle \cup D\left\langle k\right\rangle \right\vert $
is proved. This completes the proof of Lemma \ref{lem.garnir.skew}
\textbf{(a)}. \medskip

\textbf{(b)} Assume that the tableau $T$ is column-standard and satisfies
$T\left(  i,j\right)  >T\left(  i,k\right)  $. We must prove that any element
of $X$ is larger than any element of $Y$. In other words, we must prove that
$x>y$ for each $x\in X$ and each $y\in Y$.

So let $x\in X$ and $y\in Y$. We must prove that $x>y$.

We have $x\in X$. In other words, $x=T\left(  i^{\prime},j\right)  $ for some
$i^{\prime}\in\mathbb{Z}$ satisfying $i^{\prime}\geq i$ and $\left(
i^{\prime},j\right)  \in D$ (by the definition of $X$). Consider this
$i^{\prime}$. The cell $\left(  i^{\prime},j\right)  $ lies weakly south of
the cell $\left(  i,j\right)  $ in the $j$-th column (since $i^{\prime}\geq
i$). Since the entries of $T$ strictly increase top-to-bottom in each column
(because $T$ is column-standard), we thus conclude that $T\left(  i^{\prime
},j\right)  \geq T\left(  i,j\right)  $. Thus, $x=T\left(  i^{\prime
},j\right)  \geq T\left(  i,j\right)  $.

A similar argument can be used to show that $y\leq T\left(  i,k\right)  $
(since the cell of $T$ that contains $y$ lies weakly north of the cell
$\left(  i,k\right)  $ in the $k$-th column). Thus, $T\left(  i,k\right)  \geq
y$. Altogether, we now have $x\geq T\left(  i,j\right)  >T\left(  i,k\right)
\geq y$. Thus, $x>y$ is proved. As we said, this completes the proof of Lemma
\ref{lem.garnir.skew} \textbf{(b)}.
\end{proof}
\end{fineprint}

\subsubsection{Transversals (systems of coset representatives) in groups}

We will use two basic notions from group theory: that of a \emph{left
transversal}, and that of a \emph{right transversal}. We briefly mentioned the
latter in Remark \ref{rmk.intX.birep}, but let us now give a better definition:

\begin{definition}
\label{def.G/H-transversal}Let $G$ be a group, and let $H$ be a subgroup of
$G$. Then: \medskip

\textbf{(a)} A subset $L$ of $G$ is said to be a \emph{left transversal of
}$H$ \emph{in }$G$ (aka a \emph{system of representatives for the left cosets
of }$H$\emph{ in }$G$) if it has the property that each left coset $uH$ of $H$
in $G$ contains exactly one element of $L$ (that is, each $u\in G$ satisfies
$\left\vert uH\cap L\right\vert =1$). \medskip

\textbf{(b)} A subset $R$ of $G$ is said to be a \emph{right transversal of
}$H$ \emph{in }$G$ (aka a \emph{system of representatives for the right cosets
of }$H$\emph{ in }$G$) if it has the property that each right coset $Hu$ of
$H$ in $G$ contains exactly one element of $R$ (that is, each $u\in G$
satisfies $\left\vert Hu\cap R\right\vert =1$).
\end{definition}

\begin{example}
Let $G$ be the symmetric group $S_{3}$, and let $H$ be its subgroup $\left\{
1,s_{1}\right\}  $. Then:

\begin{itemize}
\item The left cosets of $H$ in $G$ are $\left\{  1,s_{1}\right\}  $ and
$\left\{  s_{2},\operatorname*{cyc}\nolimits_{1,3,2}\right\}  $ and $\left\{
t_{1,3},\operatorname*{cyc}\nolimits_{1,2,3}\right\}  $.

\item The right cosets of $H$ in $G$ are $\left\{  1,s_{1}\right\}  $ and
$\left\{  s_{2},\operatorname*{cyc}\nolimits_{1,2,3}\right\}  $ and $\left\{
t_{1,3},\operatorname*{cyc}\nolimits_{1,3,2}\right\}  $.
\end{itemize}

Thus:

\begin{itemize}
\item The set $\left\{  1,s_{2},\operatorname*{cyc}\nolimits_{1,2,3}\right\}
$ is a left transversal of $H$ in $G$, but not a right transversal of $H$ in
$G$ (since it contains two elements of the right coset $\left\{
s_{2},\operatorname*{cyc}\nolimits_{1,2,3}\right\}  $ but no elements of the
right coset $\left\{  t_{1,3},\operatorname*{cyc}\nolimits_{1,3,2}\right\}  $).

\item The set $\left\{  1,s_{2},\operatorname*{cyc}\nolimits_{1,3,2}\right\}
$ is a right transversal of $H$ in $G$, but not a left transversal.

\item The set $\left\{  s_{1},s_{2},t_{1,3}\right\}  $ is both a left
transversal of $H$ in $G$ and a right transversal of $H$ in $G$.
\end{itemize}
\end{example}

Clearly, the notions of a left transversal and of a right transversal are
mutually analogous, and in fact their roles get swapped if we reverse the
order of multiplication in $G$. Thus, any result we can prove about one notion
can be translated into a corresponding result about the other.

For instance, here is a nearly obvious (but important) fact:

\begin{proposition}
\label{prop.G/H-transversal.existR}Let $G$ be a group, and let $H$ be a
subgroup of $G$. Then: \medskip

\textbf{(a)} There exists at least one left transversal of $H$ in $G$.
\medskip

\textbf{(b)} There exists at least one right transversal of $H$ in $G$.
\end{proposition}

\begin{fineprint}
\begin{proof}
[Proof sketch.]\textbf{(a)} The left cosets of $H$ in $G$ are nonempty and
disjoint. Thus, we can obtain a left transversal simply by choosing one
element from each left coset.\footnote{If $G$ is infinite, then this uses the
axiom of choice.} \medskip

\textbf{(b)} Analogous to part \textbf{(a)}.
\end{proof}
\end{fineprint}

What is much less obvious is that (as we have already mentioned in Remark
\ref{rmk.intX.birep}) there exists a subset of $G$ that is simultaneously a
left transversal and a right transversal of $H$ in $G$, provided that $G$ is
finite (a theorem of Miller); but we will not need it.

We will, however, use the following simple fact:

\begin{proposition}
\label{prop.G/H-transversal.bijL}Let $G$ be a group, and let $H$ be a subgroup
of $G$. Let $L$ be a left transversal of $H$ in $G$. Then, the map%
\begin{align*}
H\times L  &  \rightarrow G,\\
\left(  h,x\right)   &  \mapsto xh
\end{align*}
is a bijection.
\end{proposition}

\begin{fineprint}
\begin{proof}
[Proof sketch.]The map is clearly well-defined. Furthermore:

\begin{itemize}
\item This map is injective.

[\textit{Proof:} Let $\left(  h,x\right)  \in H\times L$ and $\left(
h^{\prime},x^{\prime}\right)  \in H\times L$ be two pairs that this map sends
to the same value. Thus, $xh=x^{\prime}h^{\prime}$. Now,%
\[
x\underbrace{H}_{\substack{=hH\\\text{(since }h\in H\text{)}}}=\underbrace{xh}%
_{=x^{\prime}h^{\prime}}H=x^{\prime}\underbrace{h^{\prime}H}%
_{\substack{=H\\\text{(since }h^{\prime}\in H\text{)}}}=x^{\prime}H.
\]
But $L$ is a left transversal of $H$; thus, the left coset $xH$ contains
exactly one element of $L$. Therefore, any two elements of $L$ contained in
$xH$ must be equal. Since $x$ and $x^{\prime}$ are two elements of $L$
contained in $xH$ (because $x\in xH$ and $x^{\prime}\in x^{\prime}H=xH$), we
thus conclude that $x$ and $x^{\prime}$ are equal. In other words,
$x=x^{\prime}$. Now, $xh=\underbrace{x^{\prime}}_{=x}h^{\prime}=xh^{\prime}$,
so that $h=h^{\prime}$ (by cancelling $x$). Combined with $x=x^{\prime}$, this
yields $\left(  h,x\right)  =\left(  h^{\prime},x^{\prime}\right)  $.

Forget that we fixed $\left(  h,x\right)  $ and $\left(  h^{\prime},x^{\prime
}\right)  $. We thus have shown that if $\left(  h,x\right)  \in H\times L$
and $\left(  h^{\prime},x^{\prime}\right)  \in H\times L$ are two pairs that
our map sends to the same value, then $\left(  h,x\right)  =\left(  h^{\prime
},x^{\prime}\right)  $. In other words, our map is injective.]

\item This map is surjective.

[\textit{Proof:} Let $g\in G$. Then, $gH$ is a left coset of $H$ in $G$. But
$L$ is a left transversal of $H$; thus, the left coset $gH$ contains exactly
one element of $L$. Let $x$ be this element. Thus, $x\in L$ and $x\in gH$, so
that $x=gh$ for some $h\in H$. Consider this $h$. Since $H$ is a subgroup of
$G$, we obtain $h^{-1}\in H$ (since $h\in H$), so that $\left(  h^{-1}%
,x\right)  \in H\times L$ (since $x\in L$). But $x=gh$ leads to $g=xh^{-1}$.
Thus, our map sends the pair $\left(  h^{-1},x\right)  \in H\times L$ to
$xh^{-1}=g$. In particular, our map takes $g$ as a value.

Forget that we fixed $g$. We thus have shown that our map takes each $g\in G$
as a value. In other words, this map is surjective.]
\end{itemize}

Hence, this map is both injective and surjective. Thus, it is bijective, i.e.,
a bijection. This proves Proposition \ref{prop.G/H-transversal.bijL}.
\end{proof}
\end{fineprint}

Obviously, Proposition \ref{prop.G/H-transversal.bijL} has an analogue for
right transversals (with an analogous proof):

\begin{proposition}
\label{prop.G/H-transversal.bijR}Let $G$ be a group, and let $H$ be a subgroup
of $G$. Let $R$ be a right transversal of $H$ in $G$. Then, the map%
\begin{align*}
H\times R  &  \rightarrow G,\\
\left(  h,x\right)   &  \mapsto hx
\end{align*}
is a bijection.
\end{proposition}

When $G$ is a subgroup of a symmetric group, these propositions have some
useful consequences in the group algebra:

\begin{corollary}
\label{cor.garnir.left-transver}Let $G$ be a subgroup of $S_{n}$. Let $H$ be a
subgroup of $G$. Let $L$ be a left transversal of $H$ in $G$. Then, in
$\mathbf{k}\left[  S_{n}\right]  $, we have%
\begin{equation}
\sum_{g\in G}g=\left(  \sum_{x\in L}x\right)  \left(  \sum_{h\in H}h\right)
\label{eq.cor.garnir.left-transver.+}%
\end{equation}
and%
\begin{equation}
\sum_{g\in G}\left(  -1\right)  ^{g}g=\left(  \sum_{x\in L}\left(  -1\right)
^{x}x\right)  \left(  \sum_{h\in H}\left(  -1\right)  ^{h}h\right)  .
\label{eq.cor.garnir.left-transver.-}%
\end{equation}

\end{corollary}

\begin{proof}
Proposition \ref{prop.G/H-transversal.bijL} says that the map%
\begin{align*}
H\times L  &  \rightarrow G,\\
\left(  h,x\right)   &  \mapsto xh
\end{align*}
is a bijection. Hence, we can substitute $xh$ for $g$ in the sum $\sum_{g\in
G}\left(  -1\right)  ^{g}g$. Thus we obtain%
\begin{align*}
\sum_{g\in G}\left(  -1\right)  ^{g}g  &  =\underbrace{\sum_{\left(
h,x\right)  \in H\times L}}_{=\sum_{h\in H}\ \ \sum_{x\in L}}%
\underbrace{\left(  -1\right)  ^{xh}}_{\substack{=\left(  -1\right)
^{x}\left(  -1\right)  ^{h}\\\text{(by the multiplicativity}\\\text{of the
sign)}}}xh=\sum_{h\in H}\ \ \sum_{x\in L}\left(  -1\right)  ^{x}\left(
-1\right)  ^{h}xh\\
&  =\left(  \sum_{x\in L}\left(  -1\right)  ^{x}x\right)  \left(  \sum_{h\in
H}\left(  -1\right)  ^{h}h\right)  .
\end{align*}
This proves (\ref{eq.cor.garnir.left-transver.-}). An analogous argument (but
without the signs) proves (\ref{eq.cor.garnir.left-transver.+}). Thus,
Corollary \ref{cor.garnir.left-transver} is verified.
\end{proof}

Of course, the analogue of Corollary \ref{cor.garnir.left-transver} for right
transversals is just as easy to state and to prove:

\begin{corollary}
\label{cor.garnir.right-transver}Let $G$ be a subgroup of $S_{n}$. Let $H$ be
a subgroup of $G$. Let $R$ be a right transversal of $H$ in $G$. Then, in
$\mathbf{k}\left[  S_{n}\right]  $, we have%
\begin{equation}
\sum_{g\in G}g=\left(  \sum_{h\in H}h\right)  \left(  \sum_{x\in R}x\right)
\label{eq.cor.garnir.right-transver.+}%
\end{equation}
and%
\begin{equation}
\sum_{g\in G}\left(  -1\right)  ^{g}g=\left(  \sum_{h\in H}\left(  -1\right)
^{h}h\right)  \left(  \sum_{x\in R}\left(  -1\right)  ^{x}x\right)  .
\label{eq.cor.garnir.right-transver.-}%
\end{equation}

\end{corollary}

\begin{exercise}
\fbox{2} Explain how Lemma \ref{lem.intX.rec3} is a particular case of
(\ref{eq.cor.garnir.right-transver.+}).
\end{exercise}

The following simple property of left transversals will be used in what follows:

\begin{lemma}
\label{lem.G/H.transversal.diff}Let $G$ be a group. Let $H$ and $K$ be two
subgroups of $G$. Let $L$ be a left transversal of $H\cap K$ in $H$. (This
makes sense, since $H\cap K$ is a subgroup of $H$.) Then: \medskip

\textbf{(a)} There is exactly one $u\in L$ that belongs to $K$. \medskip

\textbf{(b)} This $u$ satisfies $L\setminus\left\{  u\right\}  =L\setminus K$.
\end{lemma}

\begin{proof}
We have $L\subseteq H$ (since $L$ is a left transversal of $H\cap K$ in $H$),
thus $H\cap L=L$. \medskip

\textbf{(a)} Consider the trivial left coset $1\left(  H\cap K\right)  $ of
the subgroup $H\cap K$ in $H$. This coset contains exactly one element of $L$
(since $L$ is a left transversal of $H\cap K$ in $H$, and therefore each left
coset of $H\cap K$ in $H$ contains exactly one element of $L$). In other
words, there is exactly one $u\in1\left(  H\cap K\right)  $ that belongs to
$L$. In other words, there is exactly one $u\in1\left(  H\cap K\right)  \cap
L$. Since $\underbrace{1\left(  H\cap K\right)  }_{=H\cap K}\cap\,L=H\cap
K\cap L=\underbrace{H\cap L}_{=L}\cap\,K=L\cap K$, we can rewrite this as
follows: There is exactly one $u\in L\cap K$. In other words, there is exactly
one $u\in L$ that belongs to $K$. This proves Lemma
\ref{lem.G/H.transversal.diff} \textbf{(a)}. \medskip

\textbf{(b)} Let $u$ be the unique $u\in L$ that belongs to $K$. Thus, $u$ is
the only element of $L$ that belongs to $K$. In other words, $u$ is the only
element of $L\cap K$. Hence, $L\cap K=\left\{  u\right\}  $.

But $L\setminus\left(  L\cap K\right)  =L\setminus K$ (since any two sets $X$
and $Y$ satisfy $X\setminus\left(  X\cap Y\right)  =X\setminus Y$). In view of
$L\cap K=\left\{  u\right\}  $, we can rewrite this as $L\setminus\left\{
u\right\}  =L\setminus K$. This proves Lemma \ref{lem.G/H.transversal.diff}
\textbf{(b)}.
\end{proof}

\subsubsection{The Garnir relations}

We are now ready to state the \emph{Garnir relations}:

\begin{theorem}
[Garnir relations]\label{thm.garnir.grel}Let $D$ be a diagram. Let $T$ be an
$n$-tableau of shape $D$. Let $j$ and $k$ be two distinct integers. Let $Z$ be
a subset of $\operatorname*{Col}\left(  j,T\right)  \cup\operatorname*{Col}%
\left(  k,T\right)  $ such that $\left\vert Z\right\vert >\left\vert
D\left\langle j\right\rangle \cup D\left\langle k\right\rangle \right\vert $.

Define $S_{n,Z}$ as in Proposition \ref{prop.intX.basics}. This set $S_{n,Z}$
is a subgroup of $S_{n}$ (by Proposition \ref{prop.intX.basics} \textbf{(a)},
applied to $X=Z$). Thus, both $\mathcal{C}\left(  T\right)  $ and $S_{n,Z}$
are subgroups of $S_{n}$. Hence, their intersection $S_{n,Z}\cap
\mathcal{C}\left(  T\right)  $ is a subgroup of $S_{n,Z}$. Let $L$ be a left
transversal of $S_{n,Z}\cap\mathcal{C}\left(  T\right)  $ in $S_{n,Z}$. Then,
in the Specht module $\mathcal{S}^{D}$, we have%
\begin{equation}
\sum_{x\in L}\left(  -1\right)  ^{x}\mathbf{e}_{xT}=0
\label{eq.thm.garnir.grel.gar}%
\end{equation}
and%
\begin{equation}
\mathbf{e}_{T}=-\sum_{x\in L\setminus\mathcal{C}\left(  T\right)  }\left(
-1\right)  ^{x}\mathbf{e}_{xT}. \label{eq.thm.garnir.grel.gar2}%
\end{equation}

\end{theorem}

The equality (\ref{eq.thm.garnir.grel.gar}) is known as a \emph{Garnir
relation} (after H. G. Garnir, who discovered it in 1950 \cite[\S 11,
Th\'{e}or\`{e}me III]{Garnir50}\footnote{To be more precise, Garnir only
proved this theorem in the case when $D$ is a Young diagram $Y\left(
\lambda\right)  $ (and he stated it in terms of the Specht--Vandermonde
avatar, i.e., as a determinantal identity). I don't know who generalized it to
general skew Young diagrams.}). The equality (\ref{eq.thm.garnir.grel.gar2})
is (as we will soon see) just a restatement of (\ref{eq.thm.garnir.grel.gar}),
in which all but one of the addends have been moved to the right hand side.

Before we prove Theorem \ref{thm.garnir.grel}, let us illustrate it on an example.

\begin{example}
\label{exa.garnir.grel.1}Let $n=8$ and $D=Y\left(  3,3,2\right)  $ and
$T=\ytableaushort{1{*(green)2}4,{*(green)6}{*(green)3}7,{*(green)8}5}$\ \ .
(This particular $T$ is straight and column-standard, although neither of
these properties is needed for Theorem \ref{thm.garnir.grel} to hold.) Set
$j=1$ and $k=2$ and $Z=\left\{  2,3,6,8\right\}  $ (the set of the entries
highlighted in green). This satisfies all requirements of Theorem
\ref{thm.garnir.grel} (since $Z\subseteq\operatorname*{Col}\left(  1,T\right)
\cup\operatorname*{Col}\left(  2,T\right)  $ and $\left\vert Z\right\vert =4$
and $\left\vert D\left\langle 1\right\rangle \cup D\left\langle 2\right\rangle
\right\vert =3$).

In order to apply Theorem \ref{thm.garnir.grel}, we need to choose a left
transversal $L$ of $S_{n,Z}\cap\mathcal{C}\left(  T\right)  $ in $S_{n,Z}$.
Let us find a good way to choose such an $L$. Note that%
\begin{align*}
S_{n,Z}  &  =\left\{  w\in S_{8}\ \mid\ w\text{ fixes }1,4,5,7\right\}
\ \ \ \ \ \ \ \ \ \ \text{and}\\
S_{n,Z}\cap\mathcal{C}\left(  T\right)   &  =\left\{  w\in S_{8}%
\ \mid\ w\text{ fixes }1,4,5,7\text{ and permutes }\left\{  6,8\right\}
\text{ and }\left\{  2,3\right\}  \right\} \\
&  =\left\{  w\in S_{8}\ \mid\ w\text{ fixes }1,4,5,7\text{ and permutes
}\left\{  6,8\right\}  \right\}
\end{align*}
(why?). Thus, two permutations $u$ and $v$ in $S_{n,Z}$ belong to the same
left coset of $S_{n,Z}\cap\mathcal{C}\left(  T\right)  $ if and only if we can
write $u$ as $u=vh$ for some $h\in S_{8}$ that fixes $1,4,5,7$ and permutes
$\left\{  6,8\right\}  $. In other words, two permutations $u$ and $v$ in
$S_{n,Z}$ belong to the same left coset of $S_{n,Z}\cap\mathcal{C}\left(
T\right)  $ if and only if $u\left(  \left\{  6,8\right\}  \right)  =v\left(
\left\{  6,8\right\}  \right)  $ (why?). Thus, each left coset of $S_{n,Z}%
\cap\mathcal{C}\left(  T\right)  $ in $S_{n,Z}$ consists of all permutations
$w\in S_{n,Z}$ that send the set $\left\{  6,8\right\}  $ to some specific
$2$-element subset of $Z$. Therefore, in order to choose a left transversal
$L$ of $S_{n,Z}\cap\mathcal{C}\left(  T\right)  $ in $S_{n,Z}$, we can proceed
as follows: For each $2$-element subset $U$ of $Z$, we choose some permutation
$x_{U}\in S_{n,Z}$ that sends $\left\{  6,8\right\}  $ to $U$. The set of all
these chosen permutations $x_{U}$ will then be our left transversal $L$.

There are many choices for these permutations $x_{U}$, but here is perhaps the
simplest one (we are also listing the corresponding tableaux $x_{U}T$):%
\[%
\begin{tabular}
[c]{|c||c|c|c|c|c|c|}\hline
$U$ & $\left\{  6,8\right\}  $ & $\left\{  2,8\right\}  $ & $\left\{
2,6\right\}  $ & $\left\{  3,8\right\}  $ & $\left\{  3,6\right\}  $ &
$\left\{  2,3\right\}  $\\\hline
$x_{U}$ & $\operatorname*{id}$ & $t_{2,6}$ & $t_{2,8}$ & $t_{3,6}$ & $t_{3,8}$
& $t_{2,6}t_{3,8}$\\\hline
$x_{U}T$ &
$\ytableaushort{1{*(green)2}4,{*(green)6}{*(green)3}7,{*(green)8}5}$ &
$\ytableaushort{1{*(green)6}4,{*(green)2}{*(green)3}7,{*(green)8}5}$ &
$\ytableaushort{1{*(green)8}4,{*(green)6}{*(green)3}7,{*(green)2}5}$ &
$\ytableaushort{1{*(green)2}4,{*(green)3}{*(green)6}7,{*(green)8}5}$ &
$\ytableaushort{1{*(green)2}4,{*(green)6}{*(green)8}7,{*(green)3}5}$ &
$\ytableaushort{1{*(green)6}4,{*(green)2}{*(green)8}7,{*(green)3}5}$\\\hline
\end{tabular}
\ \ \ .
\]

So we obtain the left transversal $L=\left\{  \operatorname*{id}%
,\ t_{2,6},\ t_{2,8},\ t_{3,6},\ t_{3,8},\ t_{2,6}t_{3,8}\right\}  $. The
Garnir relation (\ref{eq.thm.garnir.grel.gar}) thus becomes%
\[
\mathbf{e}_{\operatorname*{id}T}-\mathbf{e}_{t_{2,6}T}-\mathbf{e}_{t_{2,8}%
T}-\mathbf{e}_{t_{3,6}T}-\mathbf{e}_{t_{3,8}T}+\mathbf{e}_{t_{2,6}t_{3,8}%
T}=0.
\]
Solving this for $\mathbf{e}_{\operatorname*{id}T}=\mathbf{e}_{T}$, we obtain%
\begin{equation}
\mathbf{e}_{T}=\mathbf{e}_{t_{2,6}T}+\mathbf{e}_{t_{2,8}T}+\mathbf{e}%
_{t_{3,6}T}+\mathbf{e}_{t_{3,8}T}-\mathbf{e}_{t_{2,6}t_{3,8}T}.
\label{eq.exa.garnir.grel.1.eT=}%
\end{equation}
This is also what we would have obtained from (\ref{eq.thm.garnir.grel.gar2}).
\end{example}

\begin{remark}
\label{rmk.garnir.canonL}Let $D$, $T$, $j$, $k$ and $Z$ be as in Theorem
\ref{thm.garnir.grel}. Set $X=Z\cap\operatorname*{Col}\left(  j,T\right)  $
and $Y=Z\cap\operatorname*{Col}\left(  k,T\right)  $. There are many possible
choices for the left transversal $L$, but are there any \textquotedblleft
natural\textquotedblright\ ones? Yes, and in fact there are several:

\begin{enumerate}
\item We can choose $L$ to be the set of all permutations $w\in S_{n,Z}$ that
are increasing on $X$ and increasing on $Y$. (These permutations are known as
\emph{riffle shuffles}.)

\item For any $\left\vert X\right\vert $-element subset $U$ of $Z$, we define
an involution $w_{U}\in S_{n,Z}$ by $w_{U}=t_{u_{1},x_{1}}t_{u_{2},x_{2}%
}\cdots t_{u_{p},x_{p}}$, where $u_{1},u_{2},\ldots,u_{p}$ are the elements of
$U\setminus X$ listed in increasing order, and where $x_{1},x_{2},\ldots
,x_{p}$ are the elements of $X\setminus U$ listed in increasing order. (Thus,
$w_{U}$ exchanges the elements of $U\setminus X$ with those of $X\setminus U$
while leaving all other elements undisturbed.) We then can choose $L$ to be
the set of all these involutions $w_{U}$ for all possible $U$. (This is how
the $L$ in Example \ref{exa.garnir.grel.1} was chosen.)

\item We can obtain variants of these choices by replacing \textquotedblleft
increasing\textquotedblright\ by \textquotedblleft
decreasing\textquotedblright.
\end{enumerate}

All these choices lead to left transversals $L$, which in general are
different. Of course, they all satisfy $\left\vert L\right\vert =\dbinom
{\left\vert X\right\vert +\left\vert Y\right\vert }{\left\vert X\right\vert }$.
\end{remark}

\begin{proof}
[Proof of Theorem \ref{thm.garnir.grel}.]We note that $\left\vert D\right\vert
=n$ (since $T$ is an $n$-tableau of shape $D$).

In the Young module $\mathcal{M}^{D}$, we have%
\begin{align}
\sum_{x\in L}\left(  -1\right)  ^{x}\underbrace{\mathbf{e}_{xT}}%
_{\substack{=x\mathbf{e}_{T}\\\text{(by Lemma \ref{lem.spechtmod.submod}
\textbf{(a)})}}}  &  =\sum_{x\in L}\left(  -1\right)  ^{x}x\mathbf{e}%
_{T}=\left(  \sum_{x\in L}\left(  -1\right)  ^{x}x\right)
\underbrace{\mathbf{e}_{T}}_{\substack{=\nabla_{\operatorname*{Col}T}%
^{-}\overline{T}\\\text{(by Proposition \ref{prop.symmetrizers.eT})}%
}}\nonumber\\
&  =\left(  \sum_{x\in L}\left(  -1\right)  ^{x}x\right)  \nabla
_{\operatorname*{Col}T}^{-}\overline{T}. \label{pf.thm.garnir.grel.1}%
\end{align}

Now, recall that both $\mathcal{C}\left(  T\right)  $ and $S_{n,Z}$ are
subgroups of $S_{n}$. Hence, their intersection $S_{n,Z}\cap\mathcal{C}\left(
T\right)  $ is a subgroup of $\mathcal{C}\left(  T\right)  $. Proposition
\ref{prop.G/H-transversal.existR} \textbf{(b)} (applied to $G=\mathcal{C}%
\left(  T\right)  $ and $H=S_{n,Z}\cap\mathcal{C}\left(  T\right)  $) thus
shows that there exists at least one right transversal of $S_{n,Z}%
\cap\mathcal{C}\left(  T\right)  $ in $\mathcal{C}\left(  T\right)  $. Pick
such a right transversal and call it $R$. Thus,
(\ref{eq.cor.garnir.right-transver.-}) (applied to $G=\mathcal{C}\left(
T\right)  $ and $H=S_{n,Z}\cap\mathcal{C}\left(  T\right)  $) yields%
\[
\sum_{g\in\mathcal{C}\left(  T\right)  }\left(  -1\right)  ^{g}g=\left(
\sum_{h\in S_{n,Z}\cap\mathcal{C}\left(  T\right)  }\left(  -1\right)
^{h}h\right)  \left(  \sum_{x\in R}\left(  -1\right)  ^{x}x\right)  .
\]
But the definition of $\nabla_{\operatorname*{Col}T}^{-}$ yields%
\begin{align}
\nabla_{\operatorname*{Col}T}^{-}  &  =\sum_{w\in\mathcal{C}\left(  T\right)
}\left(  -1\right)  ^{w}w=\sum_{g\in\mathcal{C}\left(  T\right)  }\left(
-1\right)  ^{g}g\nonumber\\
&  =\left(  \sum_{h\in S_{n,Z}\cap\mathcal{C}\left(  T\right)  }\left(
-1\right)  ^{h}h\right)  \left(  \sum_{x\in R}\left(  -1\right)  ^{x}x\right)
. \label{pf.thm.garnir.grel.2}%
\end{align}

On the other hand, Proposition \ref{prop.intX.basics} \textbf{(b)} (applied to
$X=Z$) yields%
\[
\nabla_{Z}=\sum\limits_{w\in S_{n,Z}}w\ \ \ \ \ \ \ \ \ \ \text{and}%
\ \ \ \ \ \ \ \ \ \ \nabla_{Z}^{-}=\sum\limits_{w\in S_{n,Z}}\left(
-1\right)  ^{w}w.
\]
Furthermore, (\ref{eq.cor.garnir.left-transver.-}) (applied to $G=S_{n,Z}$ and
$H=S_{n,Z}\cap\mathcal{C}\left(  T\right)  $) yields that
\[
\sum_{g\in S_{n,Z}}\left(  -1\right)  ^{g}g=\left(  \sum_{x\in L}\left(
-1\right)  ^{x}x\right)  \left(  \sum_{h\in S_{n,Z}\cap\mathcal{C}\left(
T\right)  }\left(  -1\right)  ^{h}h\right)
\]
(since $L$ is a left transversal of $S_{n,Z}\cap\mathcal{C}\left(  T\right)  $
in $S_{n,Z}$). Thus,%
\begin{align}
\nabla_{Z}^{-}  &  =\sum\limits_{w\in S_{n,Z}}\left(  -1\right)  ^{w}%
w=\sum_{g\in S_{n,Z}}\left(  -1\right)  ^{g}g\nonumber\\
&  =\left(  \sum_{x\in L}\left(  -1\right)  ^{x}x\right)  \left(  \sum_{h\in
S_{n,Z}\cap\mathcal{C}\left(  T\right)  }\left(  -1\right)  ^{h}h\right)  .
\label{pf.thm.garnir.grel.Z}%
\end{align}

Now, (\ref{pf.thm.garnir.grel.1}) becomes%
\begin{align}
\sum_{x\in L}\left(  -1\right)  ^{x}\mathbf{e}_{xT}  &  =\left(  \sum_{x\in
L}\left(  -1\right)  ^{x}x\right)  \nabla_{\operatorname*{Col}T}^{-}%
\overline{T}\nonumber\\
&  =\underbrace{\left(  \sum_{x\in L}\left(  -1\right)  ^{x}x\right)  \left(
\sum_{h\in S_{n,Z}\cap\mathcal{C}\left(  T\right)  }\left(  -1\right)
^{h}h\right)  }_{\substack{=\nabla_{Z}^{-}\\\text{(by
(\ref{pf.thm.garnir.grel.Z}))}}}\underbrace{\left(  \sum_{x\in R}\left(
-1\right)  ^{x}x\right)  \overline{T}}_{=\sum_{x\in R}\left(  -1\right)
^{x}x\overline{T}}\nonumber\\
&  \ \ \ \ \ \ \ \ \ \ \ \ \ \ \ \ \ \ \ \ \left(  \text{by
(\ref{pf.thm.garnir.grel.2})}\right) \nonumber\\
&  =\nabla_{Z}^{-}\sum_{x\in R}\left(  -1\right)  ^{x}x\overline{T}\nonumber\\
&  =\sum_{x\in R}\left(  -1\right)  ^{x}\nabla_{Z}^{-}x\overline{T}.
\label{pf.thm.garnir.grel.5}%
\end{align}

Now, let $x\in R$. Then, $x\in R\subseteq\mathcal{C}\left(  T\right)  $ (by
the definition of $R$). Hence, Proposition \ref{prop.tableau.Sn-act.0}
\textbf{(b)} (applied to $w=x$) shows that the tableau $x\rightharpoonup T$ is
column-equivalent to $T$. Thus, Lemma \ref{lem.garnir.ST} (applied to
$S=x\rightharpoonup T$) shows that $\nabla_{Z}^{-}\overline{x\rightharpoonup
T}=0$. In other words, $\nabla_{Z}^{-}x\overline{T}=0$ (since $\overline
{x\rightharpoonup T}=\overline{xT}=x\overline{T}$).

Forget that we fixed $x$. We thus have shown that $\nabla_{Z}^{-}x\overline
{T}=0$ for each $x\in R$. Thus, (\ref{pf.thm.garnir.grel.5}) rewrites as%
\[
\sum_{x\in L}\left(  -1\right)  ^{x}\mathbf{e}_{xT}=\sum_{x\in R}\left(
-1\right)  ^{x}\underbrace{\nabla_{Z}^{-}x\overline{T}}_{=0}=0.
\]
This proves (\ref{eq.thm.garnir.grel.gar}). It remains to prove
(\ref{eq.thm.garnir.grel.gar2}).

Lemma \ref{lem.G/H.transversal.diff} \textbf{(a)} (applied to $G=S_{n}$ and
$H=S_{n,Z}$ and $K=\mathcal{C}\left(  T\right)  $) shows that there is exactly
one $u\in L$ that belongs to $\mathcal{C}\left(  T\right)  $ (since $L$ is a
left transversal of $S_{n,Z}\cap\mathcal{C}\left(  T\right)  $ in $S_{n,Z}$).
Consider this $u$. Thus, $u\in\mathcal{C}\left(  T\right)  $. Hence, Lemma
\ref{lem.spechtmod.submod} \textbf{(c)} yields $\mathbf{e}_{uT}=\left(
-1\right)  ^{u}\mathbf{e}_{T}$. Thus,
\[
\left(  -1\right)  ^{u}\mathbf{e}_{uT}=\underbrace{\left(  -1\right)
^{u}\left(  -1\right)  ^{u}}_{=\left(  \left(  -1\right)  ^{2}\right)  ^{u}%
=1}\mathbf{e}_{T}=\mathbf{e}_{T}.
\]

Moreover, Lemma \ref{lem.G/H.transversal.diff} \textbf{(b)} (applied to
$G=S_{n}$ and $H=S_{n,Z}$ and $K=\mathcal{C}\left(  T\right)  $) shows that
$L\setminus\left\{  u\right\}  =L\setminus\mathcal{C}\left(  T\right)  $.

Now, (\ref{eq.thm.garnir.grel.gar}) yields%
\begin{align*}
0  &  =\sum_{x\in L}\left(  -1\right)  ^{x}\mathbf{e}_{xT}=\underbrace{\left(
-1\right)  ^{u}\mathbf{e}_{uT}}_{=\mathbf{e}_{T}}+\underbrace{\sum
_{\substack{x\in L;\\x\neq u}}}_{=\sum_{x\in L\setminus\left\{  u\right\}  }%
}\left(  -1\right)  ^{x}\mathbf{e}_{xT}\\
&  \ \ \ \ \ \ \ \ \ \ \ \ \ \ \ \ \ \ \ \ \left(  \text{here, we have split
off the addend for }x=u\text{ from the sum}\right) \\
&  =\mathbf{e}_{T}+\sum_{x\in L\setminus\left\{  u\right\}  }\left(
-1\right)  ^{x}\mathbf{e}_{xT}.
\end{align*}
Solving this for $\mathbf{e}_{T}$, we find%
\[
\mathbf{e}_{T}=-\sum_{x\in L\setminus\left\{  u\right\}  }\left(  -1\right)
^{x}\mathbf{e}_{xT}=-\sum_{x\in L\setminus\mathcal{C}\left(  T\right)
}\left(  -1\right)  ^{x}\mathbf{e}_{xT}%
\]
(since $L\setminus\left\{  u\right\}  =L\setminus\mathcal{C}\left(  T\right)
$). This proves (\ref{eq.thm.garnir.grel.gar2}). Thus, the proof of Theorem
\ref{thm.garnir.grel} is complete.
\end{proof}

Theorem \ref{thm.garnir.grel} gives a linear dependence relation
(\ref{eq.thm.garnir.grel.gar}) between several polytabloids $\mathbf{e}_{xT}$;
this relation allows us to express any of these polytabloids through the
others. In particular, we can use it (in the form
(\ref{eq.thm.garnir.grel.gar2})) to express $\mathbf{e}_{T}$ as a $\mathbf{k}%
$-linear combination of other polytabloids. In the next section, we shall use
this to express the non-standard polytabloids through the standard ones when
$D$ is a skew Young diagram.

\begin{remark}
In Theorem \ref{thm.garnir.grel}, the left transversal $L$ can be chosen at
will. However, the resulting Garnir relation will not significantly depend on
the choice, because of the following:
\end{remark}

\begin{exercise}
\fbox{1} Prove that any addend $\left(  -1\right)  ^{x}\mathbf{e}_{xT}$ in
(\ref{eq.thm.garnir.grel.gar}) depends only on the left coset $x\left(
S_{n,Z}\cap\mathcal{C}\left(  T\right)  \right)  $, but not on $x$ itself.

[\textbf{Hint:} Lemma \ref{lem.spechtmod.submod} \textbf{(c)}.]
\end{exercise}

\begin{remark}
Theorem \ref{thm.garnir.grel} can be easily generalized to a \textquotedblleft
multi-column Garnir relation\textquotedblright, in which the two integers $j$
and $k$ are replaced by several integers. We will have no need for this
generality, however.
\end{remark}

\subsection{The standard basis theorem}

\subsubsection{The theorem}

We can now describe the linear structure of $\mathcal{S}^{D}$ completely when
$D$ is a skew Young diagram:

\begin{theorem}
[standard basis theorem]\label{thm.spechtmod.basis}Let $D$ be a skew Young
diagram. Then, the standard polytabloids $\mathbf{e}_{T}$ (that is, the
$\mathbf{e}_{T}$ where $T$ ranges over all standard tableaux of shape $D$)
form a basis of the $\mathbf{k}$-module $\mathcal{S}^{D}$.
\end{theorem}

This theorem was first proved by Garnir in 1950 for straight Young diagrams
$D=Y\left(  \lambda\right)  $ (see \cite[\S 12]{Garnir50}).

Before we properly prove Theorem \ref{thm.spechtmod.basis}, let me explain the
rough idea of the proof: We already know (from Theorem
\ref{thm.spechtmod.linind}) that the standard polytabloids $\mathbf{e}_{T}$
are $\mathbf{k}$-linearly independent. It remains to show that they span
$\mathcal{S}^{D}$. It suffices to show that any polytabloid $\mathbf{e}_{T}$
can be written as a $\mathbf{k}$-linear combination of standard polytabloids.
Furthermore, it suffices to prove this in the case when $T$ is
column-standard. Indeed, if an $n$-tableau $T$ is not column-standard, then we
can find a column-standard $n$-tableau $S$ that is column-equivalent to $T$
(just by sorting each column of $T$ in increasing order), and then Lemma
\ref{lem.spechtmod.submod} \textbf{(d)} yields $\mathbf{e}_{T}=\pm
\mathbf{e}_{S}$, so that any set of vectors that spans $\mathbf{e}_{S}$ will
also span $\mathbf{e}_{T}$.

Thus, we only need to show that any column-standard polytabloid $\mathbf{e}%
_{T}$ can be written as a $\mathbf{k}$-linear combination of standard polytabloids.

We try to leverage the Garnir relations (\ref{eq.thm.garnir.grel.gar2}) to
show this. Here is an example:

\begin{example}
\label{exa.spechtmod.basis.straighten1}Let $n=8$ and $D=Y\left(  3,3,2\right)
$ and $T=\ytableaushort{124,637,85}$\ \ . This tableau $T$ is column-standard,
but not standard, so we need to expand $\mathbf{e}_{T}$ as a $\mathbf{k}%
$-linear combination of standard polytabloids. To this purpose, we consider
the two entries $T\left(  2,1\right)  =6$ and $T\left(  2,2\right)  =3$ of
$T$, which witness the non-standardness of $T$ (since $6>3$). We denote the
respective two cells $\left(  2,1\right)  $ and $\left(  2,2\right)  $ by
$\left(  i,j\right)  $ and $\left(  i,k\right)  $ (so we set $i=2$ and $j=1$
and $k=2$), and we apply the Garnir relation (\ref{eq.thm.garnir.grel.gar2})
to the set $Z=X\cup Y$, where%
\begin{align*}
X  &  :=\left\{  T\left(  i^{\prime},j\right)  \ \mid\ i^{\prime}\in
\mathbb{Z}\text{ with }i^{\prime}\geq i\text{ and }\left(  i^{\prime
},j\right)  \in D\right\}  =\left\{  6,8\right\}
\ \ \ \ \ \ \ \ \ \ \text{and}\\
Y  &  :=\left\{  T\left(  i^{\prime},k\right)  \ \mid\ i^{\prime}\in
\mathbb{Z}\text{ with }i^{\prime}\leq i\text{ and }\left(  i^{\prime
},k\right)  \in D\right\}  =\left\{  2,3\right\}  .
\end{align*}
The result is the equality (\ref{eq.exa.garnir.grel.1.eT=}) (since our set $Z$
is precisely the set $Z=\left\{  2,3,6,8\right\}  $ we used in Example
\ref{exa.garnir.grel.1}). Thus, we have expressed our non-standard polytabloid
$\mathbf{e}_{T}$ as a $\mathbf{k}$-linear combination (and, in fact, a
$\mathbb{Z}$-linear combination) of the polytabloids $\mathbf{e}_{t_{2,6}%
T},\ \mathbf{e}_{t_{2,8}T},\ \mathbf{e}_{t_{3,6}T},\ \mathbf{e}_{t_{3,8}%
T},\ \mathbf{e}_{t_{2,6}t_{3,8}T}$.

Alas, these latter polytabloids are not all standard. For example, the tableau
$t_{2,6}T$ is neither row-standard ($6>4$ and $8>5$) nor column-standard
($6>3$). Thus, we have not yet expressed $\mathbf{e}_{T}$ as a $\mathbf{k}%
$-linear combination of standard polytabloids.

But have we at least moved closer to this goal?
\end{example}

As we will soon see, the answer is \textquotedblleft yes\textquotedblright. To
make this precise, we need to define an appropriate notion of
\textquotedblleft closer\textquotedblright, i.e., a partial order on the set
of all $n$-tableaux such that all our new tableaux $t_{2,6}T,\ t_{2,8}%
T,\ \ldots$ (more generally: all the tableaux appearing on the right hand side
of (\ref{eq.thm.garnir.grel.gar2})) are larger than $T$. More precisely, we
will introduce a total order on the set of all \textquotedblleft%
$n$-column-tabloids\textquotedblright\ (an analogue of $n$-tabloids relying on
column-equivalence instead of row-equivalence).

\subsubsection{Column-tabloids}

First, we define $n$-column-tabloids. The definition is completely analogous
to Definition \ref{def.tabloid.tabloid}:

\begin{definition}
\label{def.coltabloid.coltabloid}\textbf{(a)} Let $D$ be a diagram. An
$n$\emph{-column-tabloid} of shape $D$ is an equivalence class of $n$-tableaux
of shape $D$ with respect to column-equivalence. For example, for $D=Y\left(
2,1\right)  $ and $n=3$, the set $\left\{  12\backslash\backslash
3,\ 32\backslash\backslash1\right\}  $ (where we are using the shorthand
notation from Convention \ref{conv.tableau.poetic}) is an $n$-column-tabloid
of shape $D$. \medskip

\textbf{(b)} We use the notation $\widetilde{T}$ for the $n$-column-tabloid
corresponding to an $n$-tableau $T$ (that is, the equivalence class containing
$T$). For instance, the $3$-column-tabloid $\left\{  12\backslash
\backslash3,\ 32\backslash\backslash1\right\}  $ can thus be written as
$\widetilde{12\backslash\backslash3}$ or (equivalently) as
$\widetilde{32\backslash\backslash1}$.

To draw an $n$-column-tabloid, we just draw one of the $n$-tableaux in it
(often we pick the column-standard one), and we drop the horizontal
subdividers between the entries of a column. For instance, the $3$%
-column-tabloid $\widetilde{12\backslash\backslash3}$ is thus drawn as $%
\begin{tabular}
[c]{|c|c}%
$1$ & \multicolumn{1}{|c|}{$2$}\\
$3$ &
\end{tabular}
\ \ $ . \medskip

\textbf{(c)} If $\widetilde{T}$ is an $n$-column-tabloid, and if
$j\in\mathbb{Z}$, then the $j$\emph{-th column of }$\widetilde{T}$ shall mean
the set of the entries of the $j$-th column of $T$. This is well-defined,
since two column-equivalent tableaux will have the same entries in their
$j$-th column. For example, the $2$-nd column of the $5$-column-tabloid
\begin{tabular}
[c]{|c|c|c}%
$2$ & $4$ & \multicolumn{1}{|c|}{$3$}\\
$1$ & $5$ &
\end{tabular}
is the set $\left\{  4,5\right\}  $.
\end{definition}

\begin{example}
There are six $4$-column-tabloids of shape $Y\left(  2,2\right)  $, namely%
\[%
\begin{tabular}
[c]{|c|c|}%
$1$ & $3$\\
$2$ & $4$%
\end{tabular}
\ ,\qquad%
\begin{tabular}
[c]{|c|c|}%
$1$ & $2$\\
$3$ & $4$%
\end{tabular}
\ ,\qquad%
\begin{tabular}
[c]{|c|c|}%
$1$ & $2$\\
$4$ & $3$%
\end{tabular}
\ ,\qquad%
\begin{tabular}
[c]{|c|c|}%
$2$ & $1$\\
$3$ & $4$%
\end{tabular}
\ ,\qquad%
\begin{tabular}
[c]{|c|c|}%
$2$ & $1$\\
$4$ & $3$%
\end{tabular}
\ ,\qquad%
\begin{tabular}
[c]{|c|c|}%
$3$ & $1$\\
$4$ & $2$%
\end{tabular}
\ .
\]

\end{example}

While $n$-column-tabloids are equivalence classes, they also can be
represented in a canonical way, namely by column-standard $n$-tableaux:

\begin{proposition}
\label{prop.coltabloid.col-st}Let $D$ be any diagram. Then, there is a
bijection%
\begin{align*}
&  \text{from }\left\{  \text{column-standard }n\text{-tableaux of shape
}D\right\} \\
&  \text{to }\left\{  n\text{-column-tabloids of shape }D\right\}
\end{align*}
that sends each column-standard $n$-tableau $T$ to its column-tabloid
$\widetilde{T}$.
\end{proposition}

\begin{proof}
Analogous to Proposition \ref{prop.tabloid.row-st}.
\end{proof}

We could furthermore define an $S_{n}$-action on the set $\left\{
n\text{-column-tabloids of shape }D\right\}  $ (for a given diagram $D$),
similarly to the analogous action for tabloids (see Definition
\ref{def.tabloid.Sn-act}). However, we will have no need for this action until
much later, so we postpone this to Definition \ref{def.coltabloid.Sn-act}.

\subsubsection{The column last letter order}

Let us now define the column last letter order: an analogue of the Young last
letter order (see Definition \ref{def.tabloid.llo} and the paragraphs
afterwards) for $n$-column-tabloids instead of $n$-tabloids.

\begin{definition}
\label{def.spechtmod.garnir-order}\textbf{(a)} If $\widetilde{T}$ is any
$n$-column-tabloid (of any shape), and if $i\in\left[  n\right]  $, then we
let $c_{\widetilde{T}}\left(  i\right)  $ be the number of the column of $T$
that contains $i$. (Clearly, this depends only on $\widetilde{T}$, not on $T$,
since we could just as well have written \textquotedblleft the column of
$\widetilde{T}$ that contains $i$\textquotedblright\ instead of
\textquotedblleft the column of $T$ that contains $i$\textquotedblright.)
\medskip

\textbf{(b)} Let $D$ be a diagram. Define a binary relation $<$ on the set
$\left\{  n\text{-column-tabloids of shape }D\right\}  $ as follows:

Let $\widetilde{T}$ and $\widetilde{S}$ be two $n$-column-tabloids of shape
$D$. Then, we declare that $\widetilde{T}<\widetilde{S}$ if and only if

\begin{itemize}
\item there exists at least one $i\in\left[  n\right]  $ such that
$c_{\widetilde{T}}\left(  i\right)  \neq c_{\widetilde{S}}\left(  i\right)  $, and

\item the \textbf{largest} such $i$ satisfies $c_{\widetilde{T}}\left(
i\right)  <c_{\widetilde{S}}\left(  i\right)  $.
\end{itemize}

In other words, we say that $\widetilde{T}<\widetilde{S}$ if the
\textbf{largest} number that appears in different columns in $\widetilde{T}$
and $\widetilde{S}$ appears further left in $\widetilde{T}$ than in
$\widetilde{S}$.
\end{definition}

\begin{example}
If $\widetilde{T}=\ $%
\begin{tabular}
[c]{|c|c}%
$5$ & \multicolumn{1}{|c|}{$4$}\\
$2$ & \multicolumn{1}{|c|}{$1$}\\
$3$ &
\end{tabular}
and $\widetilde{S}=\
\begin{tabular}
[c]{|c|c}%
$1$ & \multicolumn{1}{|c|}{$2$}\\
$5$ & \multicolumn{1}{|c|}{$4$}\\
$3$ &
\end{tabular}
$, then $\widetilde{T}<\widetilde{S}$, since the largest number that appears
in different columns in $\widetilde{T}$ and $\widetilde{S}$ is $2$, and this
number appears further left in $\widetilde{T}$ than in $\widetilde{S}$.
\end{example}

The following proposition is an analogue to Proposition
\ref{prop.spechtmod.row.order}, and its proof is virtually identical to the
proof of the latter:

\begin{proposition}
\label{prop.spechtmod.garnir-order.order}Let $D$ be a diagram. Then, the
relation $<$ introduced in Definition \ref{def.spechtmod.garnir-order}
\textbf{(b)} is the smaller relation of a total order on the set $\left\{
n\text{-column-tabloids of shape }D\right\}  $.
\end{proposition}

This total order $<$ is called the \emph{column last letter order}. We define
the relation $>$ on the set $\left\{  n\text{-column-tabloids of shape
}D\right\}  $ accordingly (by setting $\widetilde{T}>\widetilde{S}$ if and
only if $\widetilde{S}<\widetilde{T}$).

\begin{exercise}
\fbox{1} Let $D$ be a diagram. Let $T$ and $S$ be two $n$-tableaux of shape
$D$. Consider also the two $n$-tableaux $T\mathbf{r}$ and $S\mathbf{r}$ of
shape $\mathbf{r}\left(  D\right)  $ (see Definition \ref{def.tableaux.r}).
Prove that $\widetilde{T}<\widetilde{S}$ (with respect to the column last
letter order) if and only if $\overline{T\mathbf{r}}<\overline{S\mathbf{r}}$
(with respect to the Young last letter order, as defined in Definition
\ref{def.tabloid.llo}).
\end{exercise}

An important (to us) combinatorial property of the column last letter order
will be the following lemma, which will help us show that the process
illustrated in Example \ref{exa.spechtmod.basis.straighten1} is really
bringing us closer to standardness:

\begin{lemma}
\label{lem.spechtmod.garnir-larger}Let $T$ be an $n$-tableau (of any shape).
Let $j$ and $k$ be two integers such that $j<k$. Let $X\subseteq
\operatorname*{Col}\left(  j,T\right)  $ and $Y\subseteq\operatorname*{Col}%
\left(  k,T\right)  $ be two subsets such that each element of $X$ is larger
than each element of $Y$. Let $w\in S_{n,X\cup Y}\setminus\mathcal{C}\left(
T\right)  $. Then, $\widetilde{wT}>\widetilde{T}$ (with respect to the column
last letter order).
\end{lemma}

\begin{example}
Let $n=10$ and
$T=\ytableaushort{\none\none{*(red)1},\none6{*(red)4},8{*(green)5}{*(red)2},3{*(green)7}9,{10}}$%
\ . Let $j=2$ and $k=3$ and $X=\left\{  5,7\right\}  $ and $Y=\left\{
1,2,4\right\}  $. (We are highlighting the entries of $X$ in green and the
entries of $Y$ in red.) Then, $j<k$ and $X\subseteq\operatorname*{Col}\left(
j,T\right)  $ and $Y\subseteq\operatorname*{Col}\left(  k,T\right)  $, and
each element of $X$ is larger than each element of $Y$. Hence, Lemma
\ref{lem.spechtmod.garnir-larger} yields that each $w\in S_{n,X\cup
Y}\setminus\mathcal{C}\left(  T\right)  $ satisfies $\widetilde{wT}%
>\widetilde{T}$. For instance, for $w=t_{4,5}$, we have
$wT=\ytableaushort{\none\none{*(red)1},\none6{*(green)5},8{*(red)4}{*(red)2},3{*(green)7}9,{10}}$
and thus visibly $\widetilde{wT}>\widetilde{T}$.
\end{example}

Lemma \ref{lem.spechtmod.garnir-larger} is responsible for the fact that all
the tableaux $t_{2,6}T$, $t_{2,8}T$, $t_{3,6}T$, $t_{3,8}T$ and $t_{2,6}%
t_{3,8}T$ in Example \ref{exa.spechtmod.basis.straighten1} are larger than $T$
(more precisely: that their $n$-column-tabloids are larger than $\widetilde{T}%
$). Hence, the formula (\ref{eq.exa.garnir.grel.1.eT=}), which expresses
$\mathbf{e}_{T}$ as a linear combination of the polytabloids corresponding to
these tableaux, is a step in the right direction.

\begin{proof}
[Proof of Lemma \ref{lem.spechtmod.garnir-larger} (sketched).]Let us first
give a slightly informal but short proof. A more formal but also more
intricate argument will be given afterwards (and in fact will prove a more
general version of the lemma).

From $w\in S_{n,X\cup Y}\setminus\mathcal{C}\left(  T\right)  $, we obtain
$w\in S_{n,X\cup Y}$ and $w\notin\mathcal{C}\left(  T\right)  $. In
particular, $w\notin\mathcal{C}\left(  T\right)  $ entails that the
$n$-tableau $w\rightharpoonup T$ is not column-equivalent to $T$ (by
Proposition \ref{prop.tableau.Sn-act.0} \textbf{(b)}), so that
$\widetilde{w\rightharpoonup T}\neq\widetilde{T}$. In other words,
$\widetilde{wT}\neq\widetilde{T}$ (since $wT=w\rightharpoonup T$).

From $w\in S_{n,X\cup Y}$, we see that the permutation $w$ fixes all elements
of $\left[  n\right]  $ that don't belong to $X\cup Y$, while permuting the
elements of $X\cup Y$ among themselves.

Thus, the $n$-column-tabloid $\widetilde{wT}$ can be obtained from
$\widetilde{T}$ by moving some elements of $X$ from column $j$ to column $k$
and moving some elements of $Y$ from column $k$ to column $j$ (since
$X\subseteq\operatorname*{Col}\left(  j,T\right)  $ and $Y\subseteq
\operatorname*{Col}\left(  k,T\right)  $).\ \ \ \ \footnote{Some elements of
$X$ and of $Y$ also move vertically within their columns, but this does not
affect the $n$-column-tabl\textbf{oid}, only the $n$-tableau.} The \# of
entries moving east (i.e., from column $j$ to column $k$) must equal the \# of
entries moving west (i.e., from column $k$ to column $j$) (since the columns
retain their sizes), and furthermore must be nonzero (since $\widetilde{wT}%
\neq\widetilde{T}$, so that at least one entry must change columns from
$\widetilde{T}$ to $\widetilde{wT}$). But any entry moving east is an element
of $X$, whereas any entry moving west is an element of $Y$. Hence, any entry
moving east is larger than any entry moving west (because each element of $X$
is larger than each element of $Y$). Consequently, the largest entry to change
columns from $\widetilde{T}$ to $\widetilde{wT}$ must be an entry moving east.
In other words, the largest number that appears in different columns in
$\widetilde{T}$ and $\widetilde{wT}$ appears further left in $\widetilde{T}$
than in $\widetilde{wT}$. In other words, $\widetilde{T}<\widetilde{wT}$. This
proves Lemma \ref{lem.spechtmod.garnir-larger}.
\end{proof}

\begin{fineprint}
As promised, let me also give a more formal proof of Lemma
\ref{lem.spechtmod.garnir-larger}. In fact, I will prove a more general result:

\begin{lemma}
\label{lem.spechtmod.garnir-larger.gen}Let $T$ be an $n$-tableau (of any
shape). Let $Z$ be a subset of $\left[  n\right]  $. Assume that $Z$ has the
following property:%
\begin{align}
&  \text{If }a,b\in Z\text{ are two numbers satisfying }c_{\widetilde{T}%
}\left(  a\right)  <c_{\widetilde{T}}\left(  b\right)  \text{,}\nonumber\\
&  \text{then }a>b\text{.} \label{eq.lem.spechtmod.garnir-larger.gen.ass}%
\end{align}
Let $w\in S_{n,Z}\setminus\mathcal{C}\left(  T\right)  $. Then,
$\widetilde{wT}>\widetilde{T}$ (with respect to the column last letter order).
\end{lemma}

\begin{example}
Let $n=9$ and
$T=\ytableaushort{\none\none\none{*(green)1},\none9{*(green)3}6,\none{*(green)5}{*(green)2},4{*(green)7},{*(green)8}}$
and $Z=\left\{  1,2,3,5,7,8\right\}  $. (The elements of $Z$ have been
highlighted in green.) Then, $Z$ has the property
(\ref{eq.lem.spechtmod.garnir-larger.gen.ass}). Thus, Lemma
\ref{lem.spechtmod.garnir-larger.gen} says that $\widetilde{wT}>\widetilde{T}$
for each $w\in S_{n,Z}\setminus\mathcal{C}\left(  T\right)  $. For example,
picking $w=t_{3,5}t_{7,8}$, we find
$wT=\ytableaushort{\none\none\none{*(green)1},\none9{*(green)5}6,\none{*(green)3}{*(green)2},4{*(green)8},{*(green)7}}$
and thus easily $\widetilde{wT}>\widetilde{T}$.
\end{example}

\begin{proof}
[Proof of Lemma \ref{lem.spechtmod.garnir-larger.gen}.]From $w\in
S_{n,Z}\setminus\mathcal{C}\left(  T\right)  $, we obtain $w\in S_{n,Z}$ and
$w\notin\mathcal{C}\left(  T\right)  $.

From $w\in S_{n,Z}$, we easily obtain $w^{-1}\left(  Z\right)  \subseteq
Z$\ \ \ \ \footnote{\textit{Proof.} Let $k\in w^{-1}\left(  Z\right)  $. Thus,
$w\left(  k\right)  \in Z$. We shall show that $k\in Z$.
\par
Indeed, assume the contrary. Thus, $k\notin Z$.
\par
However, $w\in S_{n,Z}$. This means that $w\left(  i\right)  =i$ for all
$i\notin Z$ (by the definition of $S_{n,Z}$). Applying this to $i=k$, we
obtain $w\left(  k\right)  =k$ (since $k\notin Z$). Hence, $w\left(  k\right)
=k\notin Z$, which contradicts $w\left(  k\right)  \in Z$. This contradiction
shows that our assumption was false. Hence, $k\in Z$ is proved.
\par
Forget that we fixed $k$. We thus have proved that $k\in Z$ for each $k\in
w^{-1}\left(  Z\right)  $. In other words, $w^{-1}\left(  Z\right)  \subseteq
Z$.}.

If the $n$-tableau $w\rightharpoonup T$ was column-equivalent to $T$, then we
would have $w\in\mathcal{C}\left(  T\right)  $ (by Proposition
\ref{prop.tableau.Sn-act.0} \textbf{(b)}), which would contradict
$w\notin\mathcal{C}\left(  T\right)  $. Hence, $w\rightharpoonup T$ is not
column-equivalent to $T$. In other words, $\widetilde{w\rightharpoonup T}%
\neq\widetilde{T}$ (since an $n$-column-tabloid is an equivalence class with
respect to column-equivalence). In other words, $\widetilde{wT}\neq
\widetilde{T}$ (since $w\rightharpoonup T=wT$).

We must prove that $\widetilde{wT}>\widetilde{T}$. Assume the contrary. Thus,
$\widetilde{wT}\leq\widetilde{T}$ (since the column last letter order is a
total order), so that $\widetilde{wT}<\widetilde{T}$ (since $\widetilde{wT}%
\neq\widetilde{T}$). By the definition of the relation $<$, this means that

\begin{itemize}
\item there exists at least one $i\in\left[  n\right]  $ such that
$c_{\widetilde{wT}}\left(  i\right)  \neq c_{\widetilde{T}}\left(  i\right)
$, and

\item the \textbf{largest} such $i$ satisfies $c_{\widetilde{wT}}\left(
i\right)  <c_{\widetilde{T}}\left(  i\right)  $.
\end{itemize}

\noindent Let $k$ be this largest $i$. Thus,%
\begin{equation}
c_{\widetilde{wT}}\left(  i\right)  =c_{\widetilde{T}}\left(  i\right)
\ \ \ \ \ \ \ \ \ \ \text{for all }i>k
\label{pf.lem.spechtmod.garnir-larger.gen.1}%
\end{equation}
(since $k$ is the \textbf{largest} $i\in\left[  n\right]  $ such that
$c_{\widetilde{wT}}\left(  i\right)  \neq c_{\widetilde{T}}\left(  i\right)
$) and%
\[
c_{\widetilde{wT}}\left(  k\right)  <c_{\widetilde{T}}\left(  k\right)
\]
(since this largest $i$ satisfies $c_{\widetilde{wT}}\left(  i\right)
<c_{\widetilde{T}}\left(  i\right)  $).

It is easy to see that%
\begin{equation}
c_{\widetilde{wT}}\left(  w\left(  m\right)  \right)  =c_{\widetilde{T}%
}\left(  m\right)  \ \ \ \ \ \ \ \ \ \ \text{for each }m\in\left[  n\right]
\label{pf.lem.spechtmod.garnir-larger.gen.wm}%
\end{equation}
\footnote{\textit{Proof:} Let $m\in\left[  n\right]  $. Let $\left(
i,j\right)  $ be the cell of $T$ that contains the entry $m$. Thus, $T\left(
i,j\right)  =m$. Hence, the number $m$ appears in the $j$-th column of $T$
(namely, in the cell $\left(  i,j\right)  $). Hence, $c_{\widetilde{T}}\left(
m\right)  =j$ (by the definition of $c_{\widetilde{T}}\left(  m\right)  $).
\par
However, the definition of $wT$ yields $wT=w\circ T$, so that $\left(
wT\right)  \left(  i,j\right)  =\left(  w\circ T\right)  \left(  i,j\right)
=w\left(  \underbrace{T\left(  i,j\right)  }_{=m}\right)  =w\left(  m\right)
$. Thus, the number $w\left(  m\right)  $ appears in the $j$-th column of $wT$
(namely, in the cell $\left(  i,j\right)  $). Hence, $c_{\widetilde{wT}%
}\left(  w\left(  m\right)  \right)  =j$ (by the definition of
$c_{\widetilde{wT}}\left(  w\left(  m\right)  \right)  $). In other words,
$c_{\widetilde{wT}}\left(  w\left(  m\right)  \right)  =c_{\widetilde{T}%
}\left(  m\right)  $ (since $c_{\widetilde{T}}\left(  m\right)  =j$). This
proves (\ref{pf.lem.spechtmod.garnir-larger.gen.wm}).}. Applying this to
$m=w^{-1}\left(  k\right)  $, we find $c_{\widetilde{wT}}\left(  w\left(
w^{-1}\left(  k\right)  \right)  \right)  =c_{\widetilde{T}}\left(
w^{-1}\left(  k\right)  \right)  $. Since $w\left(  w^{-1}\left(  k\right)
\right)  =k$, we can rewrite this as $c_{\widetilde{wT}}\left(  k\right)
=c_{\widetilde{T}}\left(  w^{-1}\left(  k\right)  \right)  $. Hence,
$c_{\widetilde{T}}\left(  w^{-1}\left(  k\right)  \right)  =c_{\widetilde{wT}%
}\left(  k\right)  <c_{\widetilde{T}}\left(  k\right)  $. Hence,
$c_{\widetilde{T}}\left(  w^{-1}\left(  k\right)  \right)  \neq
c_{\widetilde{T}}\left(  k\right)  $, so that $w^{-1}\left(  k\right)  \neq
k$. Hence, $w\left(  k\right)  \neq k$.

Recall that $w\in S_{n,Z}$. This means that the equality $w\left(  i\right)
=i$ holds for all $i\notin Z$ (by the definition of $S_{n,Z}$). If we had
$k\notin Z$, then we could apply this equality to $i=k$ and obtain $w\left(
k\right)  =k$, which would contradict $w\left(  k\right)  \neq k$. Hence, we
cannot have $k\notin Z$. Thus, $k\in Z$. Therefore, $w^{-1}\left(  k\right)
\in w^{-1}\left(  Z\right)  \subseteq Z$.

Hence, (\ref{eq.lem.spechtmod.garnir-larger.gen.ass}) (applied to
$a=w^{-1}\left(  k\right)  $ and $b=k$) yields $w^{-1}\left(  k\right)  >k$
(since $c_{\widetilde{T}}\left(  w^{-1}\left(  k\right)  \right)
<c_{\widetilde{T}}\left(  k\right)  $).

Now we claim the following:

\begin{statement}
\textit{Claim 1:} For each positive integer $i$, we have $w^{-i}\left(
k\right)  >k$ and $w^{-i}\left(  k\right)  \in Z$ and $c_{\widetilde{T}%
}\left(  w^{-i}\left(  k\right)  \right)  <c_{\widetilde{T}}\left(  k\right)
$.
\end{statement}

\begin{proof}
[Proof of Claim 1.]We induct on $i$:

\textit{Base case:} We have shown that $w^{-1}\left(  k\right)  >k$ and
$w^{-1}\left(  k\right)  \in Z$ and $c_{\widetilde{T}}\left(  w^{-1}\left(
k\right)  \right)  <c_{\widetilde{T}}\left(  k\right)  $. In other words,
Claim 1 holds for $i=1$.

\textit{Induction step:} Let $j>1$ be an integer. Assume (as the induction
hypothesis) that Claim 1 holds for $i=j-1$. In other words, we have
$w^{-\left(  j-1\right)  }\left(  k\right)  >k$ and $w^{-\left(  j-1\right)
}\left(  k\right)  \in Z$ and $c_{\widetilde{T}}\left(  w^{-\left(
j-1\right)  }\left(  k\right)  \right)  <c_{\widetilde{T}}\left(  k\right)  $.
We must show that Claim 1 holds for $i=j$.

Now,
\begin{align*}
w\left(  w^{-j}\left(  k\right)  \right)   &  =w^{1+\left(  -j\right)
}\left(  k\right)  =w^{-\left(  j-1\right)  }\left(  k\right)
\ \ \ \ \ \ \ \ \ \ \left(  \text{since }1+\left(  -j\right)  =-\left(
j-1\right)  \right) \\
&  \in Z,
\end{align*}
so that
\[
w^{-j}\left(  k\right)  \in w^{-1}\left(  Z\right)  \subseteq Z.
\]
Furthermore, from $w^{-\left(  j-1\right)  }\left(  k\right)  >k$, we obtain
\begin{align*}
c_{\widetilde{wT}}\left(  w^{-\left(  j-1\right)  }\left(  k\right)  \right)
&  =c_{\widetilde{T}}\left(  w^{-\left(  j-1\right)  }\left(  k\right)
\right)  \ \ \ \ \ \ \ \ \ \ \left(  \text{by
(\ref{pf.lem.spechtmod.garnir-larger.gen.1}), applied to }i=w^{-\left(
j-1\right)  }\left(  k\right)  \right) \\
&  <c_{\widetilde{T}}\left(  k\right)  .
\end{align*}
However, (\ref{pf.lem.spechtmod.garnir-larger.gen.wm}) (applied to
$m=w^{-j}\left(  k\right)  $) yields
\[
c_{\widetilde{wT}}\left(  w\left(  w^{-j}\left(  k\right)  \right)  \right)
=c_{\widetilde{T}}\left(  w^{-j}\left(  k\right)  \right)  .
\]
Therefore,%
\[
c_{\widetilde{T}}\left(  w^{-j}\left(  k\right)  \right)  =c_{\widetilde{wT}%
}\left(  \underbrace{w\left(  w^{-j}\left(  k\right)  \right)  }_{=w^{-\left(
j-1\right)  }\left(  k\right)  }\right)  =c_{\widetilde{wT}}\left(
w^{-\left(  j-1\right)  }\left(  k\right)  \right)  <c_{\widetilde{T}}\left(
k\right)  .
\]

Thus, (\ref{eq.lem.spechtmod.garnir-larger.gen.ass}) (applied to
$a=w^{-j}\left(  k\right)  $ and $b=k$) yields $w^{-j}\left(  k\right)  >k$
(since $w^{-j}\left(  k\right)  \in Z$ and $k\in Z$).

Altogether, we have now proved that $w^{-j}\left(  k\right)  >k$ and
$w^{-j}\left(  k\right)  \in Z$ and $c_{\widetilde{T}}\left(  w^{-j}\left(
k\right)  \right)  <c_{\widetilde{T}}\left(  k\right)  $. In other words,
Claim 1 holds for $i=j$. This completes the induction step. Thus, Claim 1 is proven.
\end{proof}

Now, recall that $w$ is an element of the finite group $S_{n}$, and thus has a
finite order (by Lagrange's theorem). In other words, there exists a positive
integer $i$ such that $w^{i}=1$. Consider this $i$. We find $w^{-i}=\left(
\underbrace{w^{i}}_{=1}\right)  ^{-1}=1^{-1}=1=\operatorname*{id}$ and thus
$w^{-i}\left(  k\right)  =k$. But Claim 1 yields $w^{-i}\left(  k\right)  >k$.
This contradicts $w^{-i}\left(  k\right)  =k$. This contradiction shows that
our assumption was wrong.

Thus, $\widetilde{wT}>\widetilde{T}$ is proved. This completes the proof of
Lemma \ref{lem.spechtmod.garnir-larger.gen}.
\end{proof}

\begin{proof}
[Proof of Lemma \ref{lem.spechtmod.garnir-larger}.]Set $Z=X\cup Y$. Then, the
set $Z$ has the property (\ref{eq.lem.spechtmod.garnir-larger.gen.ass}%
)\footnote{\textit{Proof.} Let $a,b\in Z$ be two numbers satisfying
$c_{\widetilde{T}}\left(  a\right)  <c_{\widetilde{T}}\left(  b\right)  $. We
must prove that $a>b$.
\par
We have $a\in Z=X\cup Y$. In other words, we have $a\in X$ or $a\in Y$. If
$a\in X$, then the number $a$ lies in the $j$-th column of $T$ (since $a\in
X\subseteq\operatorname*{Col}\left(  j,T\right)  $) and thus satisfies
$c_{\widetilde{T}}\left(  a\right)  =j$. If $a\in Y$, then we instead obtain
$c_{\widetilde{T}}\left(  a\right)  =k$ (by a similar argument using
$Y\subseteq\operatorname*{Col}\left(  k,T\right)  $). Thus, we always have
$c_{\widetilde{T}}\left(  a\right)  =j$ or $c_{\widetilde{T}}\left(  a\right)
=k$ (since we have $a\in X$ or $a\in Y$). In other words, $c_{\widetilde{T}%
}\left(  a\right)  \in\left\{  j,k\right\}  $. Similarly, $c_{\widetilde{T}%
}\left(  b\right)  \in\left\{  j,k\right\}  $.
\par
From $c_{\widetilde{T}}\left(  b\right)  \in\left\{  j,k\right\}  $, we obtain
$c_{\widetilde{T}}\left(  b\right)  \leq\max\left\{  j,k\right\}  =k$ (since
$j<k$). From $c_{\widetilde{T}}\left(  a\right)  \in\left\{  j,k\right\}  $,
we obtain $c_{\widetilde{T}}\left(  a\right)  \geq\min\left\{  j,k\right\}
=j$ (since $j<k$). However, we assumed that $c_{\widetilde{T}}\left(
a\right)  <c_{\widetilde{T}}\left(  b\right)  $. Thus, $c_{\widetilde{T}%
}\left(  a\right)  <c_{\widetilde{T}}\left(  b\right)  \leq k$, so that
$c_{\widetilde{T}}\left(  a\right)  \neq k$. In other words, the number $a$
does not lie in the $k$-th column of $T$. In other words, $a\notin%
\operatorname*{Col}\left(  k,T\right)  $. Hence, $a\notin Y$ (since $a\in Y$
would imply $a\in Y\subseteq\operatorname*{Col}\left(  k,T\right)  $, which
would contradict $a\notin\operatorname*{Col}\left(  k,T\right)  $). Combining
$a\in X\cup Y$ with $a\notin Y$, we obtain $a\in\left(  X\cup Y\right)
\setminus Y\subseteq X$.
\par
Furthermore, from $c_{\widetilde{T}}\left(  a\right)  <c_{\widetilde{T}%
}\left(  b\right)  $, we obtain $c_{\widetilde{T}}\left(  b\right)
>c_{\widetilde{T}}\left(  a\right)  \geq j$, so that $c_{\widetilde{T}}\left(
b\right)  \neq j$. In other words, the number $b$ does not lie in the $j$-th
column of $T$. In other words, $b\notin\operatorname*{Col}\left(  j,T\right)
$. Hence, $b\notin X$ (since $b\in X$ would imply $b\in X\subseteq
\operatorname*{Col}\left(  j,T\right)  $, which would contradict
$b\notin\operatorname*{Col}\left(  j,T\right)  $). Combining $b\in Z=X\cup Y$
with $b\notin X$, we obtain $b\in\left(  X\cup Y\right)  \setminus X\subseteq
Y$.
\par
Now, recall that each element of $X$ is larger than each element of $Y$. In
other words, each $x\in X$ and each $y\in Y$ satisfy $x>y$. Applying this to
$x=a$ and $y=b$, we obtain $a>b$ (since $a\in X$ and $b\in Y$). Thus, we have
shown that $Z$ satisfies (\ref{eq.lem.spechtmod.garnir-larger.gen.ass}).}.
Moreover, we have $w\in S_{n,X\cup Y}\setminus\mathcal{C}\left(  T\right)
=S_{n,Z}\setminus\mathcal{C}\left(  T\right)  $ (since $X\cup Y=Z$).
Therefore, Lemma \ref{lem.spechtmod.garnir-larger.gen} yields $\widetilde{wT}%
>\widetilde{T}$. This proves Lemma \ref{lem.spechtmod.garnir-larger}.
\end{proof}
\end{fineprint}

\subsubsection{Proof of the standard basis theorem}

Now we can actually use the Garnir relations to prove the standard basis
theorem. The main step is the following lemma:

\begin{lemma}
\label{lem.spechtmod.straighten1}Let $D$ be a skew Young diagram. Let $T$ be
an $n$-tableau of shape $D$ that is column-standard but not row-standard.
Then, $\mathbf{e}_{T}$ can be written as a $\mathbb{Z}$-linear combination of
polytabloids $\mathbf{e}_{P}$, where the subscripts $P$ are $n$-tableaux of
shape $D$ satisfying $\widetilde{P}>\widetilde{T}$ (with respect to the column
last letter order).
\end{lemma}

\begin{proof}
The tableau $T$ is not row-standard. Hence, $T$ has two entries in the same
row that are in the wrong order (i.e., the entry further left is $\geq$ to the
entry further right). In other words, there exist two entries $T\left(
i,j\right)  $ and $T\left(  i,k\right)  $ of $T$ with $j<k$ and $T\left(
i,j\right)  \geq T\left(  i,k\right)  $. Consider these $i$, $j$ and $k$.
Define the two subsets%
\begin{align*}
X  &  :=\left\{  T\left(  i^{\prime},j\right)  \ \mid\ i^{\prime}\in
\mathbb{Z}\text{ with }i^{\prime}\geq i\text{ and }\left(  i^{\prime
},j\right)  \in D\right\}  \ \ \ \ \ \ \ \ \ \ \text{and}\\
Y  &  :=\left\{  T\left(  i^{\prime},k\right)  \ \mid\ i^{\prime}\in
\mathbb{Z}\text{ with }i^{\prime}\leq i\text{ and }\left(  i^{\prime
},k\right)  \in D\right\}
\end{align*}
of $\left[  n\right]  $. Thus, $X$ is the set of entries of $T$ in the $j$-th
column weakly below the cell $\left(  i,j\right)  $, while $Y$ is the set of
entries of $T$ in the $k$-th column weakly above the cell $\left(  i,k\right)
$. In particular, $X\subseteq\operatorname*{Col}\left(  j,T\right)  $ and
$Y\subseteq\operatorname*{Col}\left(  k,T\right)  $.

Note that $j\neq k$ (since $j<k$), and thus $\left(  i,j\right)  \neq\left(
i,k\right)  $. Since $T$ is an injective map (because $T$ is an $n$-tableau),
we thus have $T\left(  i,j\right)  \neq T\left(  i,k\right)  $ as well.
Combining this with $T\left(  i,j\right)  \geq T\left(  i,k\right)  $, we
obtain $T\left(  i,j\right)  >T\left(  i,k\right)  $. Since $T$ is
column-standard, we thus readily find that each element of $X$ is larger than
each element of $Y$\ \ \ \ \footnote{\textit{Proof.} Let $x\in X$ and $y\in
Y$. We must prove that $x>y$.
\par
We have $x\in X=\left\{  T\left(  i^{\prime},j\right)  \ \mid\ i^{\prime}%
\in\mathbb{Z}\text{ with }i^{\prime}\geq i\text{ and }\left(  i^{\prime
},j\right)  \in D\right\}  $. In other words, $x$ is an entry of $T$ in the
$j$-th column weakly below the cell $\left(  i,j\right)  $. But $T$ is
column-standard; thus, the entries in the $j$-th column of $T$ increase from
top to bottom. Hence, $x\geq T\left(  i,j\right)  $ (since $x$ is an entry of
$T$ in the $j$-th column weakly below the cell $\left(  i,j\right)  $).
\par
We have $y\in Y=\left\{  T\left(  i^{\prime},k\right)  \ \mid\ i^{\prime}%
\in\mathbb{Z}\text{ with }i^{\prime}\leq i\text{ and }\left(  i^{\prime
},k\right)  \in D\right\}  $. In other words, $y$ is an entry of $T$ in the
$k$-th column weakly above the cell $\left(  i,k\right)  $. But $T$ is
column-standard; thus, the entries in the $k$-th column of $T$ increase from
top to bottom. Hence, $y\leq T\left(  i,k\right)  $ (since $y$ is an entry of
$T$ in the $k$-th column weakly above the cell $\left(  i,k\right)  $). In
other words, $T\left(  i,k\right)  \geq y$.
\par
Now, $x\geq T\left(  i,j\right)  >T\left(  i,k\right)  \geq y$, qed.}.
Moreover, $D$ is a skew Young diagram, and thus Lemma \ref{lem.garnir.skew}
\textbf{(a)} entails%
\begin{equation}
\left\vert X\cup Y\right\vert >\left\vert D\left\langle j\right\rangle \cup
D\left\langle k\right\rangle \right\vert
\label{pf.lem.spechtmod.straighten1.ineq}%
\end{equation}
\footnote{Even better: $\left\vert X\cup Y\right\vert =\left\vert
D\left\langle j\right\rangle \cup D\left\langle k\right\rangle \right\vert
+1$. Do you see why? (But we don't need this.)}. Let us set $Z=X\cup Y$. Thus,
we can rewrite (\ref{pf.lem.spechtmod.straighten1.ineq}) as $\left\vert
Z\right\vert >\left\vert D\left\langle j\right\rangle \cup D\left\langle
k\right\rangle \right\vert $. Furthermore, we have $Z=\underbrace{X}%
_{\subseteq\operatorname*{Col}\left(  j,T\right)  }\cup\underbrace{Y}%
_{\subseteq\operatorname*{Col}\left(  k,T\right)  }\subseteq
\operatorname*{Col}\left(  j,T\right)  \cup\operatorname*{Col}\left(
k,T\right)  $. Therefore, the Garnir relations (Theorem \ref{thm.garnir.grel})
can be applied, once we pick a left transversal $L$ of $S_{n,Z}\cap
\mathcal{C}\left(  T\right)  $ in $S_{n,Z}$. However, we can easily do this:
We can easily see (as in Theorem \ref{thm.garnir.grel}) that $S_{n,Z}%
\cap\mathcal{C}\left(  T\right)  $ is a subgroup of $S_{n,Z}$. Thus,
Proposition \ref{prop.G/H-transversal.existR} \textbf{(a)} (applied to
$G=S_{n,Z}$ and $H=S_{n,Z}\cap\mathcal{C}\left(  T\right)  $) shows that there
exists at least one left transversal $L$ of $S_{n,Z}\cap\mathcal{C}\left(
T\right)  $ in $S_{n,Z}$. Consider this $L$. Then,
(\ref{eq.thm.garnir.grel.gar2}) yields%
\begin{equation}
\mathbf{e}_{T}=-\sum\limits_{x\in L\setminus\mathcal{C}\left(  T\right)
}\left(  -1\right)  ^{x}\mathbf{e}_{xT}.
\label{pf.lem.spechtmod.straighten1.garn}%
\end{equation}

However, for each $x\in L\setminus\mathcal{C}\left(  T\right)  $, we have
$x\in\underbrace{L}_{\substack{\subseteq S_{n,Z}=S_{n,X\cup Y}\\\text{(since
}Z=X\cup Y\text{)}}}\setminus\,\mathcal{C}\left(  T\right)  \subseteq
S_{n,X\cup Y}\setminus\mathcal{C}\left(  T\right)  $ and therefore
$\widetilde{xT}>\widetilde{T}$ in the column last letter order (by Lemma
\ref{lem.spechtmod.garnir-larger} (applied to $w=x$), since each element of
$X$ is larger than each element of $Y$). In other words, each $x$ on the right
hand side of (\ref{pf.lem.spechtmod.straighten1.garn}) satisfies
$\widetilde{xT}>\widetilde{T}$. Hence, the equality
(\ref{pf.lem.spechtmod.straighten1.garn}) expresses $\mathbf{e}_{T}$ as a
$\mathbb{Z}$-linear combination of polytabloids $\mathbf{e}_{P}$, where the
subscripts $P$ are $n$-tableaux of shape $D$ satisfying $\widetilde{P}%
>\widetilde{T}$ (with respect to the column last letter order). This proves
Lemma \ref{lem.spechtmod.straighten1}.
\end{proof}

\begin{proposition}
[straightening law for standard polytabloids]\label{prop.spechtmod.straighten}%
Let $D$ be a skew Young diagram. Let $T$ be an $n$-tableau of shape $D$. Then,
$\mathbf{e}_{T}$ is a $\mathbb{Z}$-linear combination of standard polytabloids
$\mathbf{e}_{Q}$, where $Q$ ranges over all standard tableaux of shape $D$.
\end{proposition}

\begin{proof}
Forget that we fixed $T$. We must show that each polytabloid $\mathbf{e}_{T}$
is a $\mathbb{Z}$-linear combination of standard polytabloids.

In this proof, we will not consider any diagrams other than $D$. Thus, the
word \textquotedblleft$n$-tableau\textquotedblright\ will always mean
\textquotedblleft$n$-tableau of shape $D$\textquotedblright. Likewise for the
word \textquotedblleft$n$-column-tabloid\textquotedblright.

Define a $\mathbb{Z}$-submodule $\overline{\mathcal{S}}^{D}$ of $\mathcal{S}%
^{D}$ by
\[
\overline{\mathcal{S}}^{D}:=\operatorname*{span}\nolimits_{\mathbb{Z}}\left\{
\mathbf{e}_{Q}\ \mid\ Q\text{ is a standard }n\text{-tableau}\right\}  .
\]
Thus, $\overline{\mathcal{S}}^{D}$ is the $\mathbb{Z}$-linear span of all
standard polytabloids. We shall show that
\begin{equation}
\mathbf{e}_{T}\in\overline{\mathcal{S}}^{D}\ \ \ \ \ \ \ \ \ \ \text{for all
}n\text{-tableaux }T. \label{pf.prop.spechtmod.straighten.goal}%
\end{equation}

Indeed, let us define the \emph{depth} $\operatorname*{dep}\widetilde{T}$ of
an $n$-column-tabloid $\widetilde{T}$ to be the \# of all $n$-column-tabloids
$\widetilde{P}$ that satisfy $\widetilde{P}>\widetilde{T}$ (with respect to
the column last letter order). This is clearly a well-defined nonnegative
integer (since there are only finitely many $n$-column-tabloids $\widetilde{P}%
$). Moreover, if two $n$-column-tabloids $\widetilde{S}$ and $\widetilde{T}$
satisfy $\widetilde{S}>\widetilde{T}$ (with respect to the column last letter
order), then%
\begin{equation}
\operatorname*{dep}\widetilde{S}<\operatorname*{dep}\widetilde{T}
\label{pf.prop.spechtmod.straighten.depless}%
\end{equation}
\footnote{\textit{Proof:} Let $\widetilde{S}$ and $\widetilde{T}$ be two
$n$-column-tabloids that satisfy $\widetilde{S}>\widetilde{T}$ (with respect
to the column last letter order). Define the two sets%
\begin{align*}
X  &  :=\left\{  n\text{-column-tabloids }\widetilde{P}\text{ that satisfy
}\widetilde{P}>\widetilde{S}\right\}  \ \ \ \ \ \ \ \ \ \ \text{and}\\
Y  &  :=\left\{  n\text{-column-tabloids }\widetilde{P}\text{ that satisfy
}\widetilde{P}>\widetilde{T}\right\}  .
\end{align*}
Then, $\operatorname*{dep}\widetilde{T}=\left\vert Y\right\vert $ (since
$\operatorname*{dep}\widetilde{T}$ is defined as the \# of all $n$%
-column-tabloids $\widetilde{P}$ that satisfy $\widetilde{P}>\widetilde{T}$)
and $\operatorname*{dep}\widetilde{S}=\left\vert X\right\vert $ (similarly).
However, each element of $X$ is an $n$-column-tabloid $\widetilde{P}$ that
satisfies $\widetilde{P}>\widetilde{S}$ (by the definition of $X$), and thus
also satisfies $\widetilde{P}>\widetilde{T}$ (since $\widetilde{P}%
>\widetilde{S}>\widetilde{T}$), and therefore belongs to $Y$ as well (by the
definition of $Y$). Thus, $X\subseteq Y$. Furthermore, $\widetilde{S}$ belongs
to $Y$ (since $\widetilde{S}>\widetilde{T}$) but not to $X$ (since we don't
have $\widetilde{S}>\widetilde{S}$). Hence, the sets $X$ and $Y$ differ in the
element $\widetilde{S}$. Thus, $X\neq Y$. Combining this with $X\subseteq Y$,
we see that $X$ is a proper subset of $Y$. Since the set $Y$ is finite, this
entails $\left\vert X\right\vert <\left\vert Y\right\vert $. In other words,
$\operatorname*{dep}\widetilde{S}<\operatorname*{dep}\widetilde{T}$ (since
$\operatorname*{dep}\widetilde{S}=\left\vert X\right\vert $ and
$\operatorname*{dep}\widetilde{T}=\left\vert Y\right\vert $). This proves
(\ref{pf.prop.spechtmod.straighten.depless}).}.

Now, we shall prove (\ref{pf.prop.spechtmod.straighten.goal}) by strong
induction on $\operatorname*{dep}\widetilde{T}$:

\textit{Induction step:} Fix an $n$-tableau $U$, and assume (as the induction
hypothesis) that (\ref{pf.prop.spechtmod.straighten.goal}) is already proved
for all $n$-tableaux $T$ satisfying $\operatorname*{dep}\widetilde{T}%
<\operatorname*{dep}\widetilde{U}$. Our goal is now to prove
(\ref{pf.prop.spechtmod.straighten.goal}) for $T=U$. In other words, our goal
is to prove that $\mathbf{e}_{U}\in\overline{\mathcal{S}}^{D}$.

The tableau $U$ may or may not be column-standard. However, it is easy to find
a column-standard $n$-tableau $S$ that is column-equivalent to $U$ (just let
$S$ be the result of sorting the entries in each column of $U$ into increasing
order). Consider this $S$. Then, $\widetilde{S}=\widetilde{U}$ (since $S$ is
column-equivalent to $U$) and $\mathbf{e}_{U}=\pm\mathbf{e}_{S}$ (by Lemma
\ref{lem.spechtmod.submod} \textbf{(d)}, applied to $T=U$).

If the $n$-tableau $S$ is standard, then we have $\mathbf{e}_{S}\in
\overline{\mathcal{S}}^{D}$ (by the definition of $\overline{\mathcal{S}}^{D}%
$, since $\mathbf{e}_{S}$ is a standard polytabloid), and thus $\mathbf{e}%
_{U}=\pm\underbrace{\mathbf{e}_{S}}_{\in\overline{\mathcal{S}}^{D}}%
\in\overline{\mathcal{S}}^{D}$ (since $\overline{\mathcal{S}}^{D}$ is a
$\mathbb{Z}$-module), which immediately completes our induction step (since
our goal is to prove that $\mathbf{e}_{U}\in\overline{\mathcal{S}}^{D}$).

Hence, for the rest of the induction step, we WLOG assume that $S$ is not
standard. Since $S$ is column-standard, we thus conclude that $S$ is not
row-standard (since otherwise, $S$ would be both row-standard and
column-standard and therefore standard).

Hence, Lemma \ref{lem.spechtmod.straighten1} (applied to $S$ instead of $T$)
shows that $\mathbf{e}_{S}$ can be written as a $\mathbb{Z}$-linear
combination of polytabloids $\mathbf{e}_{P}$, where the subscripts $P$ are
$n$-tableaux of shape $D$ satisfying $\widetilde{P}>\widetilde{S}$ (with
respect to the column last letter order). In other words,
\begin{equation}
\mathbf{e}_{S}\in\operatorname*{span}\nolimits_{\mathbb{Z}}\left\{
\mathbf{e}_{P}\ \mid\ P\text{ is an }n\text{-tableau satisfying }%
\widetilde{P}>\widetilde{S}\right\}  . \label{pf.prop.spechtmod.straighten.5}%
\end{equation}

Now, let $P$ be an $n$-tableau satisfying $\widetilde{P}>\widetilde{S}$. Then,
$\widetilde{P}>\widetilde{S}=\widetilde{U}$ and therefore $\operatorname*{dep}%
\widetilde{P}<\operatorname*{dep}\widetilde{U}$ (by
(\ref{pf.prop.spechtmod.straighten.depless}), applied to $P$ and $U$ instead
of $S$ and $T$). But our induction hypothesis says that
(\ref{pf.prop.spechtmod.straighten.goal}) is already proved for all
$n$-tableaux $T$ satisfying $\operatorname*{dep}\widetilde{T}%
<\operatorname*{dep}\widetilde{U}$. Applying this to $T=P$, we conclude that
(\ref{pf.prop.spechtmod.straighten.goal}) is already proved for $T=P$ (since
$\operatorname*{dep}\widetilde{P}<\operatorname*{dep}\widetilde{U}$). Hence,
$\mathbf{e}_{P}\in\overline{\mathcal{S}}^{D}$.

Forget that we fixed $P$. We thus have shown that $\mathbf{e}_{P}\in
\overline{\mathcal{S}}^{D}$ whenever $P$ is an $n$-tableau satisfying
$\widetilde{P}>\widetilde{S}$. In other words,%
\[
\left\{  \mathbf{e}_{P}\ \mid\ P\text{ is an }n\text{-tableau satisfying
}\widetilde{P}>\widetilde{S}\right\}  \subseteq\overline{\mathcal{S}}^{D}.
\]

Now, (\ref{pf.prop.spechtmod.straighten.5}) becomes%
\begin{align*}
\mathbf{e}_{S}  &  \in\operatorname*{span}\nolimits_{\mathbb{Z}}%
\underbrace{\left\{  \mathbf{e}_{P}\ \mid\ P\text{ is an }n\text{-tableau
satisfying }\widetilde{P}>\widetilde{S}\right\}  }_{\subseteq\overline
{\mathcal{S}}^{D}}\\
&  \subseteq\operatorname*{span}\nolimits_{\mathbb{Z}}\left(  \overline
{\mathcal{S}}^{D}\right)  =\overline{\mathcal{S}}^{D}%
\ \ \ \ \ \ \ \ \ \ \left(  \text{since }\overline{\mathcal{S}}^{D}\text{ is a
}\mathbb{Z}\text{-module}\right)  .
\end{align*}
Hence, $\mathbf{e}_{U}=\pm\underbrace{\mathbf{e}_{S}}_{\in\overline
{\mathcal{S}}^{D}}\in\overline{\mathcal{S}}^{D}$ (since $\overline
{\mathcal{S}}^{D}$ is a $\mathbb{Z}$-module). But this is precisely what we
intended to prove. Thus, we have proved
(\ref{pf.prop.spechtmod.straighten.goal}) for $T=U$. This completes the
induction step. Thus, (\ref{pf.prop.spechtmod.straighten.goal}) is proved for
all $T$.

In other words, we have proved that each polytabloid $\mathbf{e}_{T}$ belongs
to $\overline{\mathcal{S}}^{D}$. In other words, each polytabloid
$\mathbf{e}_{T}$ belongs to the $\mathbb{Z}$-linear span of all standard
polytabloids (since $\overline{\mathcal{S}}^{D}$ is the $\mathbb{Z}$-linear
span of all standard polytabloids). In other words, each polytabloid
$\mathbf{e}_{T}$ is a $\mathbb{Z}$-linear combination of standard
polytabloids. This proves Proposition \ref{prop.spechtmod.straighten}.
\end{proof}

The above proof of Proposition \ref{prop.spechtmod.straighten} is fully
constructive: It provides a recursive algorithm for expressing any polytabloid
as a $\mathbb{Z}$-linear combination of standard polytabloids. This expansion
is called \emph{straightening}. Note that the word \textquotedblleft
straightening\textquotedblright\ here is an archaism; a better word would be
\textquotedblleft standardization\textquotedblright\ (since it makes the
tableau standard, not straight-shaped), but this word has a different meaning already.

We can now easily prove the standard basis theorem (Theorem
\ref{thm.spechtmod.basis}):

\begin{proof}
[Proof of Theorem \ref{thm.spechtmod.basis}.]In this proof, we will not
consider any diagrams other than $D$. Thus, the word \textquotedblleft
polytabloid\textquotedblright\ will always mean \textquotedblleft polytabloid
of shape $D$\textquotedblright.

The definition of the Specht module $\mathcal{S}^{D}$ yields%
\[
\mathcal{S}^{D}=\operatorname*{span}\nolimits_{\mathbf{k}}\left\{  \text{all
polytabloids }\mathbf{e}_{T}\right\}  .
\]

However, each polytabloid $\mathbf{e}_{T}$ is a $\mathbb{Z}$-linear
combination of standard polytabloids (by Proposition
\ref{prop.spechtmod.straighten}), and thus is a $\mathbf{k}$-linear
combination of standard polytabloids as well (since any $\mathbb{Z}$-linear
combination in a $\mathbf{k}$-module is automatically a $\mathbf{k}$-linear
combination). In other words, each polytabloid $\mathbf{e}_{T}$ belongs to
$\operatorname*{span}\nolimits_{\mathbf{k}}\left\{  \text{standard
polytabloids}\right\}  $. In other words,%
\[
\left\{  \text{all polytabloids }\mathbf{e}_{T}\right\}  \subseteq
\operatorname*{span}\nolimits_{\mathbf{k}}\left\{  \text{standard
polytabloids}\right\}  .
\]
Hence,%
\begin{align*}
\mathcal{S}^{D}  &  =\operatorname*{span}\nolimits_{\mathbf{k}}%
\underbrace{\left\{  \text{all polytabloids }\mathbf{e}_{T}\right\}
}_{\subseteq\operatorname*{span}\nolimits_{\mathbf{k}}\left\{  \text{standard
polytabloids}\right\}  }\\
&  \subseteq\operatorname*{span}\nolimits_{\mathbf{k}}\left(
\operatorname*{span}\nolimits_{\mathbf{k}}\left\{  \text{standard
polytabloids}\right\}  \right) \\
&  =\operatorname*{span}\nolimits_{\mathbf{k}}\left\{  \text{standard
polytabloids}\right\}
\end{align*}
(since $\operatorname*{span}\nolimits_{\mathbf{k}}\left\{  \text{standard
polytabloids}\right\}  $ is a $\mathbf{k}$-module). In other words, the
standard polytabloids span the $\mathbf{k}$-module $\mathcal{S}^{D}$. Since
these standard polytabloids are furthermore $\mathbf{k}$-linearly independent
(by Theorem \ref{thm.spechtmod.linind}), we thus conclude that they form a
basis of $\mathcal{S}^{D}$. Thus, Theorem \ref{thm.spechtmod.basis} is proved
at last.
\end{proof}

The specific coefficients that appear when a given polytabloid $\mathbf{e}%
_{T}$ is expanded in the standard polytabloid basis are not easily described.
(A recent paper by Hodges \cite{Hodges23} describes them in the case of a
straight Young diagram $D=Y\left(  \lambda\right)  $, but even that
description is still essentially recursive, as it involves a sum over a set of
sequences.\footnote{Hodges works in the somewhat more general setting of
$\operatorname*{GL}\left(  m\right)  $-representations and semistandard
tableaux instead of $S_{n}$-representations and standard tableaux; but we can
recover our setting by restricting ourselves to the \textquotedblleft
multilinear component\textquotedblright, i.e., to the span of the $n$-tabloids
that contain each of the entries $1,2,\ldots,n$ exactly once.} There are
various other algorithms for computing these coefficients -- such as the
\textquotedblleft turbo-straightening\textquotedblright\ of \cite{CaLaLe92}
and \cite[Exercice 2.2.19]{Krob95} -- that compute them faster than a naive
approach using Garnir relations, but are still a far cry from an explicit formula.)

\begin{remark}
Proposition \ref{prop.spechtmod.straighten} can be slightly strengthened:

Let $D$ be a skew Young diagram. Let $T$ be an $n$-tableau of shape $D$. Then,
$\mathbf{e}_{T}$ is a $\mathbb{Z}$-linear combination of standard polytabloids
$\mathbf{e}_{Q}$, where $Q$ ranges over all standard tableaux of shape $D$
that satisfy $\widetilde{Q}\geq\widetilde{T}$ in the column last letter order.

This can be proved using the same induction argument that we used to prove
Proposition \ref{prop.spechtmod.straighten}, with just a little more bookkeeping.
\end{remark}

\begin{exercise}
\label{exe.spechtmod.diradd}\fbox{2} Let $D$ be a skew Young diagram. Prove
that the Specht module $\mathcal{S}^{D}$ is a direct addend of $\mathcal{M}%
^{D}$ as a $\mathbf{k}$-module (although not necessarily as a representation
of $S_{n}$).
\end{exercise}

\subsubsection{Remarks on bad shapes}

What goes wrong if $D$ is a bad shape (i.e., not a skew Young diagram)?
Theorem \ref{thm.spechtmod.basis} does not hold in this case, as we already
found in Remark \ref{rmk.spechtmod.not-bas}; but let us see what goes wrong if
we try to straighten a polytabloid using the Garnir relations:

\begin{example}
Let $n$ and $D$ be as in Remark \ref{rmk.spechtmod.not-bas}. Then, the
$n$-tableau $\ytableaushort{21,\none3}=21\backslash\backslash\circ3$ of shape
$D$ is not standard. Can we simplify the polytabloid $\mathbf{e}_{T}$ using a
Garnir relation? The only Garnir relation that can be applied is one that
involves all three entries $1,2,3$ (that is, Theorem \ref{thm.garnir.grel} for
$j=1$ and $k=2$ and $Z=\left\{  1,2,3\right\}  $), since otherwise the
inequality $\left\vert Z\right\vert >\left\vert D\left\langle j\right\rangle
\cup D\left\langle k\right\rangle \right\vert $ would not be satisfied. This
Garnir relation (specifically, (\ref{eq.thm.garnir.grel.gar2})) gives%
\[
\mathbf{e}_{21\backslash\backslash\circ3}=\mathbf{e}_{12\backslash
\backslash\circ3}+\mathbf{e}_{31\backslash\backslash\circ2}.
\]
But the tableau $31\backslash\backslash\circ2$ is still not standard, and
applying a further Garnir relation to this tableau yields%
\[
\mathbf{e}_{31\backslash\backslash\circ2}=\mathbf{e}_{13\backslash
\backslash\circ2}+\mathbf{e}_{21\backslash\backslash\circ3}.
\]
Note that the $\mathbf{e}_{21\backslash\backslash\circ3}$ reappears again
here. Thus, we have expressed $\mathbf{e}_{21\backslash\backslash\circ3}$
through $\mathbf{e}_{31\backslash\backslash\circ2}$, but then expressed
$\mathbf{e}_{31\backslash\backslash\circ2}$ back through $\mathbf{e}%
_{21\backslash\backslash\circ3}$. In other words, we are running in circles.
Hence, we cannot use Garnir relations to express $\mathbf{e}_{21\backslash
\backslash\circ3}$ entirely in terms of standard polytabloids.
\end{example}

Of course, this example does not doom the Specht module $\mathcal{S}^{D}$ to
be an inscrutable mystery. Indeed, when $D=%
%TCIMACRO{\TeXButton{tikz tromino}{\begin{tikzpicture}[scale=0.7]
%\draw[fill=red!50] (1, 0) rectangle (2, 1);
%\draw[fill=red!50] (0, 1) rectangle (1, 2);
%\draw[fill=red!50] (1, 1) rectangle (2, 2);
%\end{tikzpicture}}}%
%BeginExpansion
\begin{tikzpicture}[scale=0.7]
\draw[fill=red!50] (1, 0) rectangle (2, 1);
\draw[fill=red!50] (0, 1) rectangle (1, 2);
\draw[fill=red!50] (1, 1) rectangle (2, 2);
\end{tikzpicture}%
%EndExpansion
$, we can find a basis of $\mathcal{S}^{D}$ by reflecting $D$ vertically (or
horizontally); the result will be a Young diagram (or a skew Young diagram),
to which the standard basis theorem (Theorem \ref{thm.spechtmod.basis})
applies. (The reflection replaces $\mathcal{S}^{D}$ by an isomorphic
representation, since $\mathcal{S}^{D}$ is defined in terms of all
$n$-tableaux, without any standardness requirement.)

More generally, if you permute the rows or the columns of a diagram $D$, then
its Young module $\mathcal{M}^{D}$ and its Specht module $\mathcal{S}^{D}$
remain unchanged (up to isomorphism). Thus, if a bad shape $D$ can be
transformed into a skew Young diagram using such permutations, then we can
still find a basis for $\mathcal{S}^{D}$. Many bad shapes (in particular, all
diagrams $D$ with at most $5$ cells) can be handled in this way.

However, there are some bad shapes $D$ that \textbf{cannot} be turned into
skew Young diagrams using such permutations. For example:

\begin{exercise}
\label{exe.spechtmod.bad.6}\fbox{2} Prove that
\[%
%TCIMACRO{\TeXButton{tikz 6-diagram}{\begin{tikzpicture}[scale=0.7]
%\draw[fill=red!50] (1, 0) rectangle (2, 1);
%\draw[fill=red!50] (0, 1) rectangle (1, 2);
%\draw[fill=red!50] (1, 1) rectangle (2, 2);
%\draw[fill=red!50] (0, 1) rectangle (-1, 0);
%\draw[fill=red!50] (1, 0) rectangle (0, -1);
%\draw[fill=red!50] (0, 0) rectangle (-1, -1);
%\end{tikzpicture}}}%
%BeginExpansion
\begin{tikzpicture}[scale=0.7]
\draw[fill=red!50] (1, 0) rectangle (2, 1);
\draw[fill=red!50] (0, 1) rectangle (1, 2);
\draw[fill=red!50] (1, 1) rectangle (2, 2);
\draw[fill=red!50] (0, 1) rectangle (-1, 0);
\draw[fill=red!50] (1, 0) rectangle (0, -1);
\draw[fill=red!50] (0, 0) rectangle (-1, -1);
\end{tikzpicture}%
%EndExpansion
=\left\{  \left(  i,j\right)  \in\left[  3\right]  ^{2}\ \mid\ i\neq
j\right\}
\]
is such a diagram.
\end{exercise}

There have been works by James, Peel, Taylor, Reiner, Shimozono, Liu and
others (\cite[Theorem 2.4]{JamPee79}, \cite{ReiShi95a}, \cite{ReiShi95b},
\cite{Taylor01}, \cite{Liu10}, \cite{Liu16}, \cite[\S 3]{BilPaw13}, etc.)
dealing with these kinds of bad-shape Specht modules, but the theory is still
mysterious and intricate and we don't know a basis in general. For instance,
explicit matrix computations using SageMath show that when $D$ is the diagram
from Exercise \ref{exe.spechtmod.bad.6}, the Specht module $\mathcal{S}^{D}$
is a free $\mathbf{k}$-module of rank $42$ (no matter what $\mathbf{k}$
is)\footnote{The method to compute this is explained in Exercise
\ref{exe.spechtmod.bad.bf} below.}, but this brute-force approach quickly
surpasses the abilities of modern computers as $D$ becomes larger. Thus,
Question \ref{quest.spechtmod.bad.free} remains unanswered in general. Some
diagrams $D$ lend themselves to easy answers:

\begin{exercise}
\textbf{(a)} \fbox{1} Let%
\[
D=\left\{  \left(  1,2\right)  ,\ \left(  2,1\right)  ,\ \left(  2,2\right)
,\ \left(  3,1\right)  \right\}  =%
%TCIMACRO{\TeXButton{tikz 6-diagram}{\begin{tikzpicture}[scale=0.7]
%\draw[fill=red!50] (0, 0) rectangle (1, 1);
%\draw[fill=red!50] (0, 1) rectangle (1, 2);
%\draw[fill=red!50] (1, 1) rectangle (2, 2);
%\draw[fill=red!50] (1, 2) rectangle (2, 3);
%\end{tikzpicture}}}%
%BeginExpansion
\begin{tikzpicture}[scale=0.7]
\draw[fill=red!50] (0, 0) rectangle (1, 1);
\draw[fill=red!50] (0, 1) rectangle (1, 2);
\draw[fill=red!50] (1, 1) rectangle (2, 2);
\draw[fill=red!50] (1, 2) rectangle (2, 3);
\end{tikzpicture}%
%EndExpansion
\ \ .
\]
Prove that $\mathcal{S}^{D}$ is a free $\mathbf{k}$-module, and find its rank.
\medskip

\textbf{(b)} \fbox{3} Let
\[
D=\left\{  \left(  1,2\right)  ,\ \left(  1,4\right)  ,\ \left(  1,6\right)
,\ \left(  2,1\right)  ,\ \left(  2,3\right)  ,\ \left(  2,6\right)  \right\}
=%
%TCIMACRO{\TeXButton{tikz turtle-and-wall diagram}{\begin{tikzpicture}%
%[scale=0.7]
%\draw[fill=red!50] (0, 0) rectangle (1, 1);
%\draw[fill=red!50] (1, 1) rectangle (2, 2);
%\draw[fill=red!50] (2, 0) rectangle (3, 1);
%\draw[fill=red!50] (3, 1) rectangle (4, 2);
%\draw[fill=red!50] (5, 0) rectangle (6, 1);
%\draw[fill=red!50] (5, 1) rectangle (6, 2);
%\end{tikzpicture}}}%
%BeginExpansion
\begin{tikzpicture}[scale=0.7]
\draw[fill=red!50] (0, 0) rectangle (1, 1);
\draw[fill=red!50] (1, 1) rectangle (2, 2);
\draw[fill=red!50] (2, 0) rectangle (3, 1);
\draw[fill=red!50] (3, 1) rectangle (4, 2);
\draw[fill=red!50] (5, 0) rectangle (6, 1);
\draw[fill=red!50] (5, 1) rectangle (6, 2);
\end{tikzpicture}%
%EndExpansion
\ \ .
\]
Prove that $\mathcal{S}^{D}$ is a free $\mathbf{k}$-module, and find its rank.
\end{exercise}

\begin{exercise}
\label{exe.spechtmod.bad.bf}Let $D$ be a diagram with $\left\vert D\right\vert
=n$. The following is a brute-force method for determining whether
$\mathcal{S}^{D}$ is \textquotedblleft well-behaved\textquotedblright.

Let $k$ be the \# of $n$-tabloids of shape $D$. Let $M\in\mathbb{Z}^{n!\times
k}$ be the $n!\times k$-matrix whose rows are indexed by the $n$-tableaux $T$
of shape $D$, and whose columns are indexed by the $n$-tabloids $\overline{S}$
of shape $D$, and whose $\left(  T,\overline{S}\right)  $-th entry is the
coefficient of $\overline{S}$ in the polytabloid $\mathbf{e}_{T}\in
\mathcal{M}^{D}$ (expanded in the basis consisting of all $n$-tabloids). Note
that all entries of $M$ are $0$'s, $1$'s and $-1$'s.

Let $M^{\operatorname*{sm}}$ be the Smith normal form of the matrix $M$. Prove
the following: \medskip

\textbf{(a)} \fbox{1} If all diagonal entries of $M^{\operatorname*{sm}}$ are
$0$'s and $1$'s, then $\mathcal{S}^{D}$ is a direct addend of $\mathcal{M}%
^{D}$ as a $\mathbf{k}$-module for any $\mathbf{k}$. \medskip

\textbf{(b)} \fbox{1} If all diagonal entries of $M^{\operatorname*{sm}}$ are
$0$'s and $1$'s, then $\mathcal{S}^{D}$ is a free $\mathbf{k}$-module of rank
$r$ for any $\mathbf{k}$, where $r$ is the \# of $1$'s on the diagonal of
$M^{\operatorname*{sm}}$. \medskip

\textbf{(c)} \fbox{1} If \textbf{not} all diagonal entries of
$M^{\operatorname*{sm}}$ are $0$'s and $1$'s, then there exists a prime number
$p$ such that the dimension of $\mathcal{S}^{D}$ for $\mathbf{k}%
=\mathbb{F}_{p}$ is smaller than the dimension of $\mathcal{S}^{D}$ for
$\mathbf{k}=\mathbb{Q}$. (Namely, $p$ can be taken to be a prime divisor of
some nonzero diagonal entry of $M^{\operatorname*{sm}}$.) \medskip

\textbf{(d)} \fbox{3} If \textbf{not} all diagonal entries of
$M^{\operatorname*{sm}}$ are $0$'s and $1$'s, then there exists a prime power
$p^{\ell}$ such that $\mathcal{S}^{D}$ is not a free $\mathbf{k}$-module for
$\mathbf{k}=\mathbb{Z}/p^{\ell}$. \medskip

\textbf{(e)} \fbox{1} All the above claims remain true if we remove from the
matrix $M$ all its rows except for those indexed by column-standard
$n$-tableaux $T$. (This produces a much smaller matrix.) \medskip

[\textbf{Hint:} See \cite[\S 43--\S 44]{Elman22} for the Smith normal form and
its uses.]
\end{exercise}

\subsubsection{The rank of $\mathcal{S}^{\lambda}$}

For future use, let us record a simple consequence of the standard basis theorem:

\begin{lemma}
\label{lem.specht.Slam-flam}Let $\lambda$ be a partition of $n$. Set
$\mathcal{S}^{\lambda}:=\mathcal{S}^{Y\left(  \lambda\right)  }$. Let
$f^{\lambda}$ denote the \# of standard tableaux of shape $Y\left(
\lambda\right)  $. Then, the Specht module $\mathcal{S}^{\lambda}$ is a free
$\mathbf{k}$-module of rank $f^{\lambda}$.
\end{lemma}

\begin{proof}
The Young diagram $Y\left(  \lambda\right)  $ is a skew Young diagram (since
$Y\left(  \lambda\right)  =Y\left(  \lambda/\varnothing\right)  $). Thus, the
standard basis theorem (Theorem \ref{thm.spechtmod.basis}, applied to
$D=Y\left(  \lambda\right)  $) shows that the standard polytabloids
$\mathbf{e}_{T}$ (that is, the $\mathbf{e}_{T}$ where $T$ ranges over all
standard tableaux of shape $Y\left(  \lambda\right)  $) form a basis of the
$\mathbf{k}$-module $\mathcal{S}^{Y\left(  \lambda\right)  }$. Hence,
$\mathcal{S}^{Y\left(  \lambda\right)  }$ is a free $\mathbf{k}$-module, and
its rank is equal to the \# of all standard tableaux of shape $Y\left(
\lambda\right)  $. Thus, its rank is $f^{\lambda}$ (since $f^{\lambda}$ is the
\# of standard tableaux of shape $Y\left(  \lambda\right)  $).

Thus we have shown that $\mathcal{S}^{Y\left(  \lambda\right)  }$ is a free
$\mathbf{k}$-module of rank $f^{\lambda}$. In other words, $\mathcal{S}%
^{\lambda}$ is a free $\mathbf{k}$-module of rank $f^{\lambda}$ (since
$\mathcal{S}^{\lambda}=\mathcal{S}^{Y\left(  \lambda\right)  }$). This proves
Lemma \ref{lem.specht.Slam-flam}.
\end{proof}

\subsection{The Young alternative}

Let us next study the case when $D$ is a straight Young diagram $Y\left(
\lambda\right)  $ of a partition $\lambda$. In this case, the Specht module
$\mathcal{S}^{Y\left(  \lambda\right)  }$ is also known as $\mathcal{S}%
^{\lambda}$, and such straight-shaped Specht modules $\mathcal{S}^{\lambda}$
play an important role: They are the irreducible representations of $S_{n}$
when $\mathbf{k}$ is a field of characteristic $0$. The proof of this fact
(Theorem \ref{thm.spechtmod.irred}) will take us some preparation.

\subsubsection{A bit of combinatorics: the Young alternative}

First, we need to lay some combinatorial groundwork. We begin with some
properties of straight-shaped Young tableaux.\setcounter{theo}{-2}

\begin{definition}
Let $T$ be an $n$-tableau of any shape. \medskip

\textbf{(a)} We say that some numbers $i_{1},i_{2},\ldots,i_{k}\in\left[
n\right]  $ are \emph{horizontally aligned} in $T$ if they all lie in the same
row of $T$. \medskip

\textbf{(b)} We say that some numbers $i_{1},i_{2},\ldots,i_{k}\in\left[
n\right]  $ are \emph{vertically aligned} in $T$ if they all lie in the same
column of $T$.
\end{definition}

\begin{example}
Let $T$ be the $8$-tableau $\ytableaushort{615,42,8,3,7}$ . Then, the numbers
$5$ and $6$ are horizontally aligned in $T$, while the numbers $1$ and $2$ are
vertically aligned in $T$.
\end{example}

\begin{theorem}
[Young alternative]\label{thm.youngtab.alt}Let $\lambda$ and $\mu$ be two
partitions of $n$.

Let $S$ be an $n$-tableau of shape $Y\left(  \lambda\right)  $. Let $T$ be an
$n$-tableau of shape $Y\left(  \mu\right)  $.

Assume that there are no two distinct integers that are horizontally aligned
in $S$ and simultaneously vertically aligned in $T$. Then: \medskip

\textbf{(a)} There exists a permutation $c\in\mathcal{C}\left(  T\right)  $
such that each $i\geq1$ satisfies the following property: Each entry of the
$i$-th row of $S$ lies in one of the rows $1,2,\ldots,i$ of $cT$. \medskip

\textbf{(b)} We have $\lambda_{1}+\lambda_{2}+\cdots+\lambda_{i}\leq\mu
_{1}+\mu_{2}+\cdots+\mu_{i}$ for each $i\geq1$. \medskip

\textbf{(c)} Now assume that $\lambda=\mu$. Then, there exist permutations
$r\in\mathcal{R}\left(  S\right)  $ and $c\in\mathcal{C}\left(  T\right)  $
such that $rS=cT$.
\end{theorem}

\Needspace{10pc}

\begin{example}
\label{exa.youngtab.alt.1}\textbf{(a)} Let $n=8$ and $\lambda=\left(
3,2,1,1,1\right)  $ and $\mu=\left(  4,3,1\right)  $ and
\[
S=\ytableaushort{615,42,8,3,7}\ \ \ \ \ \ \ \ \ \ \text{and}%
\ \ \ \ \ \ \ \ \ \ T=\ytableaushort{3786,152,4}\ \ .
\]
Then, it is easy to check that there are no two distinct integers that are
horizontally aligned in $S$ and simultaneously vertically aligned in $T$. (For
example, the entries $6,1,5$ of the first row of $S$ lie in distinct columns
of $T$, and so do the entries $4$ and $2$ of the second row of $S$.) Thus,
Theorem \ref{thm.youngtab.alt} \textbf{(a)} yields that there exists a
permutation $c\in\mathcal{C}\left(  T\right)  $ such that each $i\geq1$
satisfies the following property: Each entry of the $i$-th row of $S$ lies in
one of the rows $1,2,\ldots,i$ of $cT$. What is this permutation $c$ ? For
example, we can take $c=\operatorname*{cyc}\nolimits_{4,3,1}t_{5,7}$, which
produces the $n$-tableau $cT=\ytableaushort{1586,472,3}$\ \ . Another valid
choice is $c=\operatorname*{cyc}\nolimits_{4,3,1}t_{5,7}t_{2,8}$. \medskip

\textbf{(b)} Let $n=6$ and $\lambda=\left(  3,2,1\right)  $ and $\mu=\left(
2,2,2\right)  $ and%
\[
S=\ytableaushort{315,42,6}\ \ \ \ \ \ \ \ \ \ \text{and}%
\ \ \ \ \ \ \ \ \ \ T=\ytableaushort{14,25,63}\ \ .
\]
Then, the claim of Theorem \ref{thm.youngtab.alt} \textbf{(b)} is not
satisfied, since for $i=2$ we have $\lambda_{1}+\lambda_{2}=3+2>2+2=\mu
_{1}+\mu_{2}$. Thus, the contrapositive of Theorem \ref{thm.youngtab.alt}
\textbf{(b)} yields that there are two distinct integers that are horizontally
aligned in $S$ and simultaneously vertically aligned in $T$. And indeed, $3$
and $5$ are two such integers. \medskip

\textbf{(c)} Let $n=6$ and $\lambda=\left(  3,2,1\right)  $ and $\mu=\left(
3,2,1\right)  $ and
\[
S=\ytableaushort{162,43,5}\ \ \ \ \ \ \ \ \ \ \text{and}%
\ \ \ \ \ \ \ \ \ \ T=\ytableaushort{536,42,1}\ \ .
\]
Then, $\lambda=\mu$. Moreover, it is again easy to check that there are no two
distinct integers that are horizontally aligned in $S$ and simultaneously
vertically aligned in $T$. Hence, Theorem \ref{thm.youngtab.alt} \textbf{(c)}
yields that there exist permutations $r\in\mathcal{R}\left(  S\right)  $ and
$c\in\mathcal{C}\left(  T\right)  $ such that $rS=cT$. And indeed, $r=t_{2,6}$
and $c=t_{1,5}t_{2,3}$ are such permutations, since they satisfy
$rS=\ytableaushort{126,43,5}\ \ =cT$.
\end{example}

We will often use Theorem \ref{thm.youngtab.alt} in the equivalent form
\textquotedblleft\textbf{either} there are two distinct integers that are
horizontally aligned in $S$ and simultaneously vertically aligned in $T$,
\textbf{or} the claims of Theorem \ref{thm.youngtab.alt} \textbf{(a)},
\textbf{(b)} and \textbf{(c)} hold\footnote{where part \textbf{(c)} is still
subject to the condition that $\lambda=\mu$}\textquotedblright. (This is a
non-exclusive \textquotedblleft or\textquotedblright, at least as far as parts
\textbf{(a)} and \textbf{(b)} are concerned.) Thus I call it the \emph{Young
alternative}.

\begin{proof}
[Proof of Theorem \ref{thm.youngtab.alt}.]In Definition \ref{def.tabloid.llo}
\textbf{(a)}, we introduced the notation $r_{\overline{T}}\left(  i\right)  $
for the number of the row in which a given $n$-tableau $T$ contains a given
number $i\in\left[  n\right]  $. This notation will again be useful in our
current proof. (Of course, we will use it not only for our specific
$n$-tableau $T$.)

Furthermore, if $i\in\left[  n\right]  $ is arbitrary, then the integer
$r_{\overline{S}}\left(  i\right)  $ (which is simply the number of the row in
which $S$ contains $i$) will be called the $S$\emph{-depth} of $i$. For
instance, if $S$ is as in Example \ref{exa.youngtab.alt.1} \textbf{(a)}, then
$r_{\overline{S}}\left(  4\right)  =2$ and $r_{\overline{S}}\left(  7\right)
=5$. We make the following trivial observation:

\begin{statement}
\textit{Claim 0:} For each $j\in\left[  n\right]  $, we have $r_{\overline{S}%
}\left(  j\right)  \geq1$.
\end{statement}

\begin{proof}
[Proof of Claim 0.]The tableau $S$ has shape $Y\left(  \lambda\right)  $,
which is a subset of $\left\{  1,2,3,\ldots\right\}  ^{2}$. Thus, each cell of
$S$ belongs to one of the rows $1,2,3,\ldots$. Therefore, any of the numbers
$1,2,\ldots,n$ appears in $S$ in one of the rows $1,2,3,\ldots$. In other
words, $r_{\overline{S}}\left(  j\right)  \in\left\{  1,2,3,\ldots\right\}  $
for each $j\in\left[  n\right]  $ (since $r_{\overline{S}}\left(  j\right)  $
is defined to be the number of the row of $S$ that contains the entry $j$).
Hence, for each $j\in\left[  n\right]  $, we have $r_{\overline{S}}\left(
j\right)  \in\left\{  1,2,3,\ldots\right\}  $ and therefore $r_{\overline{S}%
}\left(  j\right)  \geq1$. This proves Claim 0.
\end{proof}

Let us now permute the entries of the $n$-tableau $T$ as follows: For each
column of $T$, we sort the entries of this column in the order of increasing
$S$-depth (so that entries with larger $S$-depth are moved south of entries
with smaller $S$-depth that lie in the same column). The resulting $n$-tableau
shall be called $R$. Thus, $R$ is an $n$-tableau of shape $Y\left(
\mu\right)  $ (since $T$ is an $n$-tableau of shape $Y\left(  \mu\right)  $),
and is column-equivalent to $T$ (since it is obtained from $T$ by sorting the
entries in each column).\footnote{For instance, if $S$ and $T$ are as in
Example \ref{exa.youngtab.alt.1} \textbf{(a)}, then%
\[
R=\ytableaushort{1526,478,3}\ \ .
\]
} Moreover, the very construction of $R$ shows that the entries in each column
of $R$ are sorted in the order of increasing $S$-depth (from top to bottom).
In other words, if $u$ and $v$ are two entries that lie in the same column of
$R$, with $u$ lying further north than $v$, then%
\begin{equation}
r_{\overline{S}}\left(  u\right)  \leq r_{\overline{S}}\left(  v\right)  .
\label{pf.thm.youngtab.alt.inc-order}%
\end{equation}

Now we claim the following:

\begin{statement}
\textit{Claim 1:} For each $k\in\left[  n\right]  $, we have $r_{\overline{R}%
}\left(  k\right)  \leq r_{\overline{S}}\left(  k\right)  $.
\end{statement}

\begin{proof}
[Proof of Claim 1.]Let $k\in\left[  n\right]  $. Let $\left(  i,j\right)  $ be
the cell of $R$ that contains the entry $k$. Thus, $\left(  i,j\right)  \in
Y\left(  \mu\right)  $ and $R\left(  i,j\right)  =k$, so that $r_{\overline
{R}}\left(  k\right)  =i$ (since $R\left(  i,j\right)  =k$ shows that the
number $k$ is contained in the $i$-th row of $R$).

We have $\left(  i,j\right)  \in Y\left(  \mu\right)  $. Hence, it can be
easily shown that all the $i$ cells $\left(  1,j\right)  ,\ \left(
2,j\right)  ,\ \ldots,\ \left(  i,j\right)  $ (that is, all the $i$ cells
$\left(  \ell,j\right)  $ for $\ell\in\left[  i\right]  $) belong to $Y\left(
\mu\right)  $ as well\footnote{\textit{Proof.} Let $\ell\in\left[  i\right]
$. Then, $\ell\leq i$. The cell $\left(  \ell,j\right)  $ belongs to $\left\{
1,2,3,\ldots\right\}  ^{2}$ and lies weakly northwest of the cell $\left(
i,j\right)  $ (since $\ell\leq i$ and $j\leq j$). Therefore, Proposition
\ref{prop.young.NE} (applied to $\mu$, $\left(  \ell,j\right)  $ and $\left(
i,j\right)  $ instead of $\lambda$, $c$ and $d$) shows that $\left(
\ell,j\right)  \in Y\left(  \mu\right)  $ (since $\left(  i,j\right)  \in
Y\left(  \mu\right)  $).
\par
Forget that we fixed $\ell$. We thus have shown that $\left(  \ell,j\right)
\in Y\left(  \mu\right)  $ for each $\ell\in\left[  i\right]  $. In other
words, all the $i$ cells $\left(  1,j\right)  ,\ \left(  2,j\right)
,\ \ldots,\ \left(  i,j\right)  $ belong to $Y\left(  \mu\right)  $.}. Hence,
their entries $R\left(  1,j\right)  ,\ R\left(  2,j\right)  ,\ \ldots
,\ R\left(  i,j\right)  $ are well-defined (since $R$ is an $n$-tableau of
shape $Y\left(  \mu\right)  $). In other words, for each $\ell\in\left[
i\right]  $, the entry $R\left(  \ell,j\right)  $ is well-defined and belongs
to $\left[  n\right]  $. Hence, let us set%
\[
r_{\ell}:=r_{\overline{S}}\left(  R\left(  \ell,j\right)  \right)
\ \ \ \ \ \ \ \ \ \ \text{for each }\ell\in\left[  i\right]  .
\]

In particular, $r_{i}=r_{\overline{S}}\left(  R\left(  i,j\right)  \right)
=r_{\overline{S}}\left(  k\right)  $ (since $R\left(  i,j\right)  =k$). Thus,
$r_{\overline{S}}\left(  k\right)  =r_{i}$.

Now, we shall show that $r_{\ell}-r_{\ell-1}\geq1$ for each $\ell\in\left\{
2,3,\ldots,i\right\}  $.

Indeed, let $\ell\in\left\{  2,3,\ldots,i\right\}  $. Thus, $r_{\ell
}=r_{\overline{S}}\left(  R\left(  \ell,j\right)  \right)  $ and $r_{\ell
-1}=r_{\overline{S}}\left(  R\left(  \ell-1,j\right)  \right)  $ (by the
definition of $r_{\ell-1}$). But $R\left(  \ell-1,j\right)  $ and $R\left(
\ell,j\right)  $ are two entries that lie in the same column of $R$ (namely,
in the $j$-th column), with $R\left(  \ell-1,j\right)  $ lying further north
than $R\left(  \ell,j\right)  $ (since the cell $\left(  \ell-1,j\right)  $
lies further north than $\left(  \ell,j\right)  $). Hence,
(\ref{pf.thm.youngtab.alt.inc-order}) (applied to $u=R\left(  \ell-1,j\right)
$ and $v=R\left(  \ell,j\right)  $) shows that $r_{\overline{S}}\left(
R\left(  \ell-1,j\right)  \right)  \leq r_{\overline{S}}\left(  R\left(
\ell,j\right)  \right)  $. In other words, $r_{\ell-1}\leq r_{\ell}$ (since
$r_{\ell}=r_{\overline{S}}\left(  R\left(  \ell,j\right)  \right)  $ and
$r_{\ell-1}=r_{\overline{S}}\left(  R\left(  \ell-1,j\right)  \right)  $).

But $r_{\overline{S}}\left(  R\left(  \ell,j\right)  \right)  $ is defined to
be the number of the row of $S$ that contains the entry $R\left(
\ell,j\right)  $. Hence, the number $R\left(  \ell,j\right)  $ lies in the
$r_{\overline{S}}\left(  R\left(  \ell,j\right)  \right)  $-th row of $S$. In
other words, the number $R\left(  \ell,j\right)  $ lies in the $r_{\ell}$-th
row of $S$ (since $r_{\ell}=r_{\overline{S}}\left(  R\left(  \ell,j\right)
\right)  $). The same argument (applied to $\ell-1$ instead of $\ell$) shows
that the number $R\left(  \ell-1,j\right)  $ lies in the $r_{\ell-1}$-th row
of $S$.

Moreover, $R$ is an $n$-tableau, and thus injective. In other words, all
entries of $R$ are distinct. Hence, $R\left(  \ell-1,j\right)  \neq R\left(
\ell,j\right)  $ (since $\left(  \ell-1,j\right)  \neq\left(  \ell,j\right)
$). Thus, $R\left(  \ell-1,j\right)  $ and $R\left(  \ell,j\right)  $ are two
distinct integers.

As we know, the tableau $R$ is column-equivalent to $T$. Thus, the $j$-th
column of $R$ contains the same entries as the $j$-th column of $T$.

The two integers $R\left(  \ell-1,j\right)  $ and $R\left(  \ell,j\right)  $
both lie in the $j$-th column of $R$ (obviously), and thus lie in the $j$-th
column of $T$ (since the $j$-th column of $R$ contains the same entries as the
$j$-th column of $T$). Thus, these two integers $R\left(  \ell-1,j\right)  $
and $R\left(  \ell,j\right)  $ lie in the same column of $T$. In other words,
they are vertically aligned in $T$.

Now, recall that there are no two distinct integers that are horizontally
aligned in $S$ and simultaneously vertically aligned in $T$. Hence, if two
distinct integers are vertically aligned in $T$, then they cannot be
horizontally aligned in $S$. Applying this to the two integers $R\left(
\ell-1,j\right)  $ and $R\left(  \ell,j\right)  $ (which are distinct and are
vertically aligned in $T$), we thus conclude that these two integers $R\left(
\ell-1,j\right)  $ and $R\left(  \ell,j\right)  $ cannot be horizontally
aligned in $S$. In other words, $R\left(  \ell-1,j\right)  $ and $R\left(
\ell,j\right)  $ cannot lie in the same row of $S$. In other words,
$r_{\ell-1}\neq r_{\ell}$ (since $R\left(  \ell-1,j\right)  $ lies in the
$r_{\ell-1}$-th row of $S$, and since $R\left(  \ell,j\right)  $ lies in the
$r_{\ell}$-th row of $S$).

Combining this with $r_{\ell-1}\leq r_{\ell}$, we obtain $r_{\ell-1}<r_{\ell}%
$. Hence, $r_{\ell-1}\leq r_{\ell}-1$ (since $r_{\ell-1}$ and $r_{\ell}$ are
integers). Thus, $r_{\ell}-r_{\ell-1}\geq1$.

Now, forget that we fixed $\ell$. We thus have proved that $r_{\ell}%
-r_{\ell-1}\geq1$ for each $\ell\in\left\{  2,3,\ldots,i\right\}  $. Adding
these inequalities for all $\ell\in\left\{  2,3,\ldots,i\right\}  $ together,
we obtain $\sum_{\ell=2}^{i}\left(  r_{\ell}-r_{\ell-1}\right)  \geq\sum
_{\ell=2}^{i}1=i-1$. Thus,%
\[
i-1\leq\sum_{\ell=2}^{i}\left(  r_{\ell}-r_{\ell-1}\right)  =r_{i}%
-r_{1}\ \ \ \ \ \ \ \ \ \ \left(  \text{by the telescope principle}\right)  .
\]

Next, recall that $r_{1}=r_{\overline{S}}\left(  R\left(  1,j\right)  \right)
$ (by the definition of $r_{1}$). But Claim 0 (applied to $R\left(
1,j\right)  $ instead of $j$) shows that $r_{\overline{S}}\left(  R\left(
1,j\right)  \right)  \geq1$. In other words, $r_{1}\geq1$ (since
$r_{1}=r_{\overline{S}}\left(  R\left(  1,j\right)  \right)  $). Now,%
\[
i-1\leq r_{i}-\underbrace{r_{1}}_{\geq1}\leq r_{i}-1,
\]
so that $i\leq r_{i}$. In other words, $r_{\overline{R}}\left(  k\right)  \leq
r_{\overline{S}}\left(  k\right)  $ (since $r_{\overline{R}}\left(  k\right)
=i$ and $r_{\overline{S}}\left(  k\right)  =r_{i}$). This proves Claim 1.
\end{proof}

Claim 1 easily yields the following:

\begin{statement}
\textit{Claim 2:} Let $i\geq1$. Then, each entry of the $i$-th row of $S$ lies
in one of the rows $1,2,\ldots,i$ of $R$.
\end{statement}

\begin{proof}
[Proof of Claim 2.]Let $k$ be an entry of the $i$-th row of $S$. We must prove
that $k$ lies in one of the rows $1,2,\ldots,i$ of $R$.

Indeed, the number $k$ lies in the $i$-th row of $S$ (by its definition). In
other words, $r_{\overline{S}}\left(  k\right)  =i$ (by the definition of
$r_{\overline{S}}\left(  k\right)  $).

Claim 0 (applied to $j=k$) yields $r_{\overline{S}}\left(  k\right)  \geq1$.
The same argument (but applied to $\mu$ and $R$ instead of $\lambda$ and $S$)
yields $r_{\overline{R}}\left(  k\right)  \geq1$.

However, Claim 1 yields $r_{\overline{R}}\left(  k\right)  \leq r_{\overline
{S}}\left(  k\right)  =i$. Combining this with $r_{\overline{R}}\left(
k\right)  \geq1$, we obtain $r_{\overline{R}}\left(  k\right)  \in\left\{
1,2,\ldots,i\right\}  $ (since $r_{\overline{R}}\left(  k\right)  $ is an
integer). In other words, $k$ lies in one of the rows $1,2,\ldots,i$ of $R$
(since $r_{\overline{R}}\left(  k\right)  $ is the number of the row of $R$ in
which $k$ lies). This completes our proof of Claim 2.
\end{proof}

Now, recall that the tableau $R$ is column-equivalent to $T$. In other words,
the $n$-tableaux $T$ and $R$ are column-equivalent. Hence, Proposition
\ref{prop.tableau.req-R} \textbf{(b)} (applied to $Y\left(  \mu\right)  $ and
$R$ instead of $D$ and $S$) shows that there exists some $w\in\mathcal{C}%
\left(  T\right)  $ such that $R=w\rightharpoonup T$. Consider this $w$. Thus,
$R=w\rightharpoonup T=wT$.

Claim 2 shows that each $i\geq1$ has the property that each entry of the
$i$-th row of $S$ lies in one of the rows $1,2,\ldots,i$ of $R$. In other
words, each $i\geq1$ has the property that each entry of the $i$-th row of $S$
lies in one of the rows $1,2,\ldots,i$ of $wT$ (since $R=wT$). Thus, there
exists a permutation $c\in\mathcal{C}\left(  T\right)  $ such that each
$i\geq1$ has the property that each entry of the $i$-th row of $S$ lies in one
of the rows $1,2,\ldots,i$ of $cT$ (namely, $c=w$). This proves Theorem
\ref{thm.youngtab.alt} \textbf{(a)}. \medskip

\textbf{(b)} Let $i\geq1$. Recall that $S$ is an $n$-tableau, thus injective.
Hence, all entries of $S$ are distinct. Moreover, $S$ has shape $Y\left(
\lambda\right)  $, so that the $k$-th row of $S$ has $\lambda_{k}$ cells for
each $k\geq1$. Hence, rows $1,2,\ldots,i$ of $S$ have $\lambda_{1}+\lambda
_{2}+\cdots+\lambda_{i}$ cells in total. Obviously, these cells contain
altogether $\lambda_{1}+\lambda_{2}+\cdots+\lambda_{i}$ many entries, which
are all distinct (since all entries of $S$ are distinct). Thus, we have shown
that rows $1,2,\ldots,i$ of $S$ have altogether $\lambda_{1}+\lambda
_{2}+\cdots+\lambda_{i}$ many distinct entries. In other words,%
\begin{equation}
\left\vert \left\{  \text{entries of rows }1,2,\ldots,i\text{ of }S\right\}
\right\vert =\lambda_{1}+\lambda_{2}+\cdots+\lambda_{i}.
\label{pf.thm.youngtab.alt.b.l}%
\end{equation}

The same argument (applied to $R$ and $\mu$ instead of $S$ and $\lambda$)
shows that%
\begin{equation}
\left\vert \left\{  \text{entries of rows }1,2,\ldots,i\text{ of }R\right\}
\right\vert =\mu_{1}+\mu_{2}+\cdots+\mu_{i} \label{pf.thm.youngtab.alt.b.m}%
\end{equation}
(since $R$ is an $n$-tableau of shape $Y\left(  \mu\right)  $).

However, Claim 1 easily shows that%
\[
\left\{  \text{entries of rows }1,2,\ldots,i\text{ of }S\right\}
\subseteq\left\{  \text{entries of rows }1,2,\ldots,i\text{ of }R\right\}
\]
\footnote{\textit{Proof.} Let $k\in\left\{  \text{entries of rows }%
1,2,\ldots,i\text{ of }S\right\}  $. Thus, $k$ is an entry of one of the rows
$1,2,\ldots,i$ of $S$. Thus, the number of the row of $S$ that contains $k$ is
$\in\left\{  1,2,\ldots,i\right\}  $. In other words, $r_{\overline{S}}\left(
k\right)  \in\left\{  1,2,\ldots,i\right\}  $ (since $r_{\overline{S}}\left(
k\right)  $ is defined as the number of the row of $S$ that contains $k$). But
Claim 1 yields $r_{\overline{R}}\left(  k\right)  \leq r_{\overline{S}}\left(
k\right)  \leq i$ (since $r_{\overline{S}}\left(  k\right)  \in\left\{
1,2,\ldots,i\right\}  $). Since we also have $r_{\overline{R}}\left(
k\right)  \geq1$ (this is proved just like in the proof of Claim 2 above), we
thus obtain $r_{\overline{R}}\left(  k\right)  \in\left\{  1,2,\ldots
,i\right\}  $. In other words, the number of the row of $R$ that contains $k$
is $\in\left\{  1,2,\ldots,i\right\}  $ (since $r_{\overline{R}}\left(
k\right)  $ is defined as the number of the row of $R$ that contains $k$). In
other words, $k$ is an entry of one of the rows $1,2,\ldots,i$ of $R$. In
other words, $k\in\left\{  \text{entries of rows }1,2,\ldots,i\text{ of
}R\right\}  $.
\par
Forget that we fixed $k$. We thus have shown that $k\in\left\{  \text{entries
of rows }1,2,\ldots,i\text{ of }R\right\}  $ for each $k\in\left\{
\text{entries of rows }1,2,\ldots,i\text{ of }S\right\}  $. In other words,
$\left\{  \text{entries of rows }1,2,\ldots,i\text{ of }S\right\}
\subseteq\left\{  \text{entries of rows }1,2,\ldots,i\text{ of }R\right\}  $%
.}. Hence,%
\[
\left\vert \left\{  \text{entries of rows }1,2,\ldots,i\text{ of }S\right\}
\right\vert \leq\left\vert \left\{  \text{entries of rows }1,2,\ldots,i\text{
of }R\right\}  \right\vert .
\]
In view of (\ref{pf.thm.youngtab.alt.b.l}) and (\ref{pf.thm.youngtab.alt.b.m}%
), we can rewrite this as
\[
\lambda_{1}+\lambda_{2}+\cdots+\lambda_{i}\leq\mu_{1}+\mu_{2}+\cdots+\mu_{i}.
\]
So Theorem \ref{thm.youngtab.alt} \textbf{(b)} is proved. \medskip

\textbf{(c)} We know that $S$ is an $n$-tableau of shape $Y\left(
\lambda\right)  $, thus an $n$-tableau of shape $Y\left(  \mu\right)  $ (since
$\lambda=\mu$). Thus, both $S$ and $R$ are $n$-tableaux of shape $Y\left(
\mu\right)  $.

We shall prove the following stronger version of Claim 1:

\begin{statement}
\textit{Claim 3:} For each $k\in\left[  n\right]  $, we have $r_{\overline{R}%
}\left(  k\right)  =r_{\overline{S}}\left(  k\right)  $.
\end{statement}

\begin{proof}
[Proof of Claim 3.]Let $k\in\left[  n\right]  $. We must prove that
$r_{\overline{R}}\left(  k\right)  =r_{\overline{S}}\left(  k\right)  $.

Assume the contrary. Thus, $r_{\overline{R}}\left(  k\right)  \neq
r_{\overline{S}}\left(  k\right)  $. But Claim 1 yields $r_{\overline{R}%
}\left(  k\right)  \leq r_{\overline{S}}\left(  k\right)  $. Hence,
$r_{\overline{R}}\left(  k\right)  <r_{\overline{S}}\left(  k\right)  $ (since
$r_{\overline{R}}\left(  k\right)  \neq r_{\overline{S}}\left(  k\right)  $).

Set $i:=r_{\overline{R}}\left(  k\right)  $. Thus, $i=r_{\overline{R}}\left(
k\right)  <r_{\overline{S}}\left(  k\right)  $, so that $r_{\overline{S}%
}\left(  k\right)  >i$. Also, obviously, $i=r_{\overline{R}}\left(  k\right)
\geq1$ (this is proved just as in the proof of Claim 2 above).

But $k$ is an entry of the $r_{\overline{R}}\left(  k\right)  $-th row of $R$
(by the definition of $r_{\overline{R}}\left(  k\right)  $). In other words,
$k$ is an entry of the $i$-th row of $R$ (since $i=r_{\overline{R}}\left(
k\right)  $). Hence, $k$ is an entry of one of rows $1,2,\ldots,i$ of $R$
(namely, of the $i$-th row). In other words,%
\begin{equation}
k\in\left\{  \text{entries of rows }1,2,\ldots,i\text{ of }R\right\}  .
\label{pf.thm.youngtab.alt.c.1R}%
\end{equation}

On the other hand, if we had $k\in\left\{  \text{entries of rows }%
1,2,\ldots,i\text{ of }S\right\}  $, then $k$ would be an entry of one of rows
$1,2,\ldots,i$ of $S$; but this would mean that $r_{\overline{S}}\left(
k\right)  \in\left\{  1,2,\ldots,i\right\}  $ (since $r_{\overline{S}}\left(
k\right)  $ is defined as the number of the row of $S$ that contains $k$),
which would contradict $r_{\overline{S}}\left(  k\right)  >i$. Hence, we
cannot have $k\in\left\{  \text{entries of rows }1,2,\ldots,i\text{ of
}S\right\}  $. Thus, we have%
\begin{equation}
k\notin\left\{  \text{entries of rows }1,2,\ldots,i\text{ of }S\right\}  .
\label{pf.thm.youngtab.alt.c.1S}%
\end{equation}

From (\ref{pf.thm.youngtab.alt.c.1R}) and (\ref{pf.thm.youngtab.alt.c.1S}), we
see that the sets $\left\{  \text{entries of rows }1,2,\ldots,i\text{ of
}R\right\}  $ and $\left\{  \text{entries of rows }1,2,\ldots,i\text{ of
}S\right\}  $ differ (indeed, the first one contains $k$, whereas the second
does not). However, in our above proof of Theorem \ref{thm.youngtab.alt}
\textbf{(b)}, we have shown that
\[
\left\{  \text{entries of rows }1,2,\ldots,i\text{ of }S\right\}
\subseteq\left\{  \text{entries of rows }1,2,\ldots,i\text{ of }R\right\}  .
\]
Thus, the set $\left\{  \text{entries of rows }1,2,\ldots,i\text{ of
}S\right\}  $ is a \textbf{proper} subset of \newline$\left\{  \text{entries
of rows }1,2,\ldots,i\text{ of }R\right\}  $ (it is proper because the sets
\newline$\left\{  \text{entries of rows }1,2,\ldots,i\text{ of }R\right\}  $
and $\left\{  \text{entries of rows }1,2,\ldots,i\text{ of }S\right\}  $
differ). Therefore,%
\[
\left\vert \left\{  \text{entries of rows }1,2,\ldots,i\text{ of }S\right\}
\right\vert <\left\vert \left\{  \text{entries of rows }1,2,\ldots,i\text{ of
}R\right\}  \right\vert
\]
(because if $X$ is a proper subset of a finite set $Y$, then $\left\vert
X\right\vert <\left\vert Y\right\vert $). In view of
(\ref{pf.thm.youngtab.alt.b.l}) and (\ref{pf.thm.youngtab.alt.b.m}), we can
rewrite this as
\[
\lambda_{1}+\lambda_{2}+\cdots+\lambda_{i}<\mu_{1}+\mu_{2}+\cdots+\mu_{i}.
\]
However, this contradicts the equality%
\[
\lambda_{1}+\lambda_{2}+\cdots+\lambda_{i}=\mu_{1}+\mu_{2}+\cdots+\mu_{i}%
\]
(which holds because $\lambda=\mu$). This contradiction shows that our
assumption was false. Thus, $r_{\overline{R}}\left(  k\right)  =r_{\overline
{S}}\left(  k\right)  $ is proved. This proves Claim 3.
\end{proof}

We have $r_{\overline{R}}\left(  i\right)  =r_{\overline{S}}\left(  i\right)
$ for each $i\in\left[  n\right]  $ (by Claim 3, applied to $k=i$). Hence,
Lemma \ref{lem.spechtmod.row.equal} (applied to $P=R$ and $Q=S$) shows that
$R$ and $S$ are row-equivalent. In other words, $S$ and $R$ are row-equivalent.

However, Proposition \ref{prop.tableau.req-R} \textbf{(a)} (applied to
$Y\left(  \mu\right)  $, $S$ and $R$ instead of $D$, $T$ and $S$) says that
the $n$-tableaux $S$ and $R$ are row-equivalent if and only if there exists
some $w^{\prime}\in\mathcal{R}\left(  S\right)  $ such that $R=w^{\prime
}\rightharpoonup S$ (note that we are using the letter $w^{\prime}$ here for
what was called $w$ in Proposition \ref{prop.tableau.req-R} \textbf{(a)},
since $w$ already has a different meaning in our context). Thus, there exists
some $w^{\prime}\in\mathcal{R}\left(  S\right)  $ such that $R=w^{\prime
}\rightharpoonup S$ (since $S$ and $R$ are row-equivalent). Consider this
$w^{\prime}$.

We have $R=w^{\prime}\rightharpoonup S=w^{\prime}S$. Comparing this with
$R=wT$, we obtain $w^{\prime}S=wT$. Thus, there exist permutations
$r\in\mathcal{R}\left(  S\right)  $ and $c\in\mathcal{C}\left(  T\right)  $
such that $rS=cT$ (namely, $r=w^{\prime}$ and $c=w$). This proves Theorem
\ref{thm.youngtab.alt} \textbf{(c)}.
\end{proof}

\begin{noncompile}
THE FOLLOWING WAS AN UNFINISHED DIFFERENT APPROACH TO PROVING Theorem
\ref{thm.youngtab.alt}. Its main advantage is that it generalizes more easily
to $m \neq n$ and generally to subsets of the entries. But it doesn't seem to
be useful and would take a while to finish.

The proof of Theorem \ref{thm.youngtab.alt} rests on a few simple lemmas. The
first is nearly trivial:

\begin{lemma}
\label{lem.youngtab.alt.cd}Let $T$ be an $n$-tableau (of any shape). Let
$d\in\mathcal{C}\left(  T\right)  $. Let $c\in\mathcal{C}\left(  dT\right)  $.
Then, $cd\in\mathcal{C}\left(  T\right)  $.
\end{lemma}

\begin{proof}
We have $d\in\mathcal{C}\left(  T\right)  $. Hence, the tableau
$d\rightharpoonup T$ is column-equivalent to $T$ (by Proposition
\ref{prop.tableau.Sn-act.0} \textbf{(b)}, applied to $w=d$). In other words,
the tableau $dT$ is column-equivalent to $T$ (since $d\rightharpoonup T=dT$).
In other words, the two $n$-tableaux $dT$ and $T$ are column-equivalent.
Hence, Proposition \ref{prop.tableau.roweq.H=H} \textbf{(b)} (applied to
$S=dT$) shows that $\mathcal{C}\left(  dT\right)  =\mathcal{C}\left(
T\right)  $. Thus, $c\in\mathcal{C}\left(  dT\right)  =\mathcal{C}\left(
T\right)  $. Combining this with $d\in\mathcal{C}\left(  T\right)  $, we
obtain $cd\in\mathcal{C}\left(  T\right)  $ (since $\mathcal{C}\left(
T\right)  $ is a group and thus is closed under multiplication). This proves
Lemma \ref{lem.youngtab.alt.cd}.
\end{proof}

In Definition \ref{def.tabloid.llo} \textbf{(a)}, we introduced the notation
$r_{\overline{T}}\left(  i\right)  $ for the number of the row in which a
given $n$-tabloid $\overline{T}$ contains a given number $i\in\left[
n\right]  $. Using this notation, we can restate the statement
\textquotedblleft each entry of the $i$-th row of $S$ lies in one of the rows
$1,2,\ldots,i$ of $cT$\textquotedblright\ in Theorem \ref{thm.youngtab.alt}
\textbf{(a)} as \textquotedblleft for each $k\in\left[  n\right]  $ satisfying
$r_{\overline{S}}\left(  k\right)  =i$, we have $r_{\overline{cT}}\left(
k\right)  \leq i$\textquotedblright. To claim that this statement holds for
all $i\geq1$ is thus tantamount to asserting that each $k\in\left[  n\right]
$ satisfies $r_{\overline{cT}}\left(  k\right)  \leq r_{\overline{S}}\left(
k\right)  $. Thus, the claim of Theorem \ref{thm.youngtab.alt} \textbf{(a)}
can be restated in the more convenient form \textquotedblleft there exists a
permutation $c\in\mathcal{C}\left(  T\right)  $ such that each $k\in\left[
n\right]  $ satisfies $r_{\overline{cT}}\left(  k\right)  \leq r_{\overline
{S}}\left(  k\right)  $\textquotedblright.

We can furthermore generalize this restated form of Theorem
\ref{thm.youngtab.alt} \textbf{(a)} to focus on a specific subset $K$ of
$\left[  n\right]  $, restricting both the assumption and the claim of the
theorem to entries coming from this set. The generalization thus claims the following:

\begin{proposition}
\label{prop.youngtab.alt.gen}Let $\lambda$ and $\mu$ be two partitions of $n$.

Let $S$ be an $n$-tableau of shape $Y\left(  \lambda\right)  $. Let $T$ be an
$n$-tableau of shape $Y\left(  \mu\right)  $.

Let $K$ be a subset of $\left[  n\right]  $. Assume that there are no two
distinct elements of $K$ that lie in the same row of $S$ and simultaneously
lie in the same column of $T$. Then, there exists a permutation $c\in
\mathcal{C}\left(  T\right)  $ such that each $k\in K$ satisfies
$r_{\overline{cT}}\left(  k\right)  \leq r_{\overline{S}}\left(  k\right)  $.
\end{proposition}

\begin{example}
Let $n=6$ and $\lambda=\left(  3,3\right)  $ and $\mu=\left(  3,2,1\right)  $.
Let $S$ and $T$ be the two $n$-tableaux given by%
\[
S=\ytableaushort{{*(green)5}{*(green)1}{*(green)6},2{*(green)4}{*(green)3}}\ \ \ \ \ \ \ \ \ \ \text{and}%
\ \ \ \ \ \ \ \ \ \ T=\ytableaushort{{*(green)1}{*(green)3}{*(green)6},{*(green)4}{*(green)5},2}\ \ ,
\]
of respective shapes $Y\left(  \lambda\right)  $ and $Y\left(  \mu\right)  $.
Let $K=\left\{  1,3,4,5,6\right\}  $ (this is the set of the entries
highlighted in green above). Then, there are no two distinct elements of $K$
that lie in the same row of $S$ and simultaneously lie in the same column of
$T$. Hence, Proposition \ref{prop.youngtab.alt.gen} says that there exists a
permutation $c\in\mathcal{C}\left(  T\right)  $ such that each $k\in K$
satisfies $r_{\overline{cT}}\left(  k\right)  \leq r_{\overline{S}}\left(
k\right)  $. And indeed, $t_{3,5}$ is such a permutation.
\end{example}

\begin{proof}
[Proof of Proposition \ref{prop.youngtab.alt.gen}.]We induct on $\left\vert
K\right\vert $:

\textit{Base case:} If $\left\vert K\right\vert =0$, then the set $K$ is
empty, so that the statement \textquotedblleft each $k\in K$ satisfies
$r_{\overline{cT}}\left(  k\right)  \leq r_{\overline{S}}\left(  k\right)
$\textquotedblright\ is vacuously true for $c=\operatorname*{id}$ (or, for
that matter, for any permutation $c\in\mathcal{C}\left(  T\right)  $). Thus,
Proposition \ref{prop.youngtab.alt.gen} is proved for $\left\vert K\right\vert
=0$. This completes the induction base.

\textit{Induction step:} Let $\kappa\in\mathbb{N}$. Assume (as the induction
hypothesis) that Proposition \ref{prop.youngtab.alt.gen} holds for $\left\vert
K\right\vert =\kappa$. We must now prove Proposition
\ref{prop.youngtab.alt.gen} for $\left\vert K\right\vert =\kappa+1$ as well.

So let $K$ be a subset of $\left[  n\right]  $ such that $\left\vert
K\right\vert =\kappa+1$. Assume that there are no two distinct elements of $K$
that lie in the same row of $S$ and simultaneously lie in the same column of
$T$.

From $\left\vert K\right\vert =\kappa+1>\kappa\geq0$, we see that $K$ is
nonempty. Hence, there exists some $u\in K$. Fix such a $u$ with
\textbf{maximum} $r_{\overline{S}}\left(  u\right)  $. (That is, fix such a
$u$ that lies as far south as possible.) Then,%
\[
r_{\overline{S}}\left(  u\right)  \geq r_{\overline{S}}\left(  \ell\right)
\ \ \ \ \ \ \ \ \ \ \text{for each }\ell\in K
\]
(since $u$ was chosen to have maximum $r_{\overline{S}}\left(  u\right)  $).

We have $u\in K$ and thus $\left\vert K\setminus\left\{  u\right\}
\right\vert =\left\vert K\right\vert -1=\kappa$ (since $\left\vert
K\right\vert =\kappa+1$). Moreover, the set $K$ has the property that there
are no two distinct elements of $K$ that lie in the same row of $S$ and
simultaneously lie in the same column of $T$ (by assumption). Hence, its
subset $K\setminus\left\{  u\right\}  $ has the analogous property (viz., that
there are no two distinct elements of $K\setminus\left\{  u\right\}  $ that
lie in the same row of $S$ and simultaneously lie in the same column of $T$)
as well. Therefore, by our induction hypothesis, we can apply Proposition
\ref{prop.youngtab.alt.gen} to $K\setminus\left\{  u\right\}  $ instead of
$K$. We thus conclude that there exists a permutation $c\in\mathcal{C}\left(
T\right)  $ such that each $k\in K\setminus\left\{  u\right\}  $ satisfies
$r_{\overline{cT}}\left(  k\right)  \leq r_{\overline{S}}\left(  k\right)  $.
Let us denote this permutation $c$ by $d$. Thus, $d\in\mathcal{C}\left(
T\right)  $ is a permutation, and each $k\in K\setminus\left\{  u\right\}  $
satisfies
\begin{equation}
r_{\overline{dT}}\left(  k\right)  \leq r_{\overline{S}}\left(  k\right)  .
\label{pf.prop.youngtab.alt.gen.IH}%
\end{equation}

We must prove that Proposition \ref{prop.youngtab.alt.gen} holds for our
subset $K$. In other words, we must prove that there exists a permutation
$c\in\mathcal{C}\left(  T\right)  $ such that each $k\in K$ satisfies
$r_{\overline{cT}}\left(  k\right)  \leq r_{\overline{S}}\left(  k\right)  $.
If every $k\in K$ satisfies $r_{\overline{dT}}\left(  k\right)  \leq
r_{\overline{S}}\left(  k\right)  $, then this is obviously true (we can just
take $c=d$ in this case). Thus, we WLOG assume that not every $k\in K$
satisfies $r_{\overline{dT}}\left(  k\right)  \leq r_{\overline{S}}\left(
k\right)  $. In other words, there exists some $k\in K$ that does not satisfy
$r_{\overline{dT}}\left(  k\right)  \leq r_{\overline{S}}\left(  k\right)  $.
This $k$ cannot belong to $K\setminus\left\{  u\right\}  $ (since otherwise,
we would have $k\in K\setminus\left\{  u\right\}  $ and therefore
$r_{\overline{dT}}\left(  k\right)  \leq r_{\overline{S}}\left(  k\right)  $
by (\ref{pf.prop.youngtab.alt.gen.IH})), and thus must be $u$ (since $k$
belongs to $K$ but cannot belong to $K\setminus\left\{  u\right\}  $, so that
$k$ belongs to $K\setminus\left(  K\setminus\left\{  u\right\}  \right)
\subseteq\left\{  u\right\}  $, and therefore $k$ equals $u$). Thus, we
conclude that $k=u$ does not satisfy $r_{\overline{dT}}\left(  k\right)  \leq
r_{\overline{S}}\left(  k\right)  $. In other words, we do not have
$r_{\overline{dT}}\left(  u\right)  \leq r_{\overline{S}}\left(  u\right)  $.
In other words, we have%
\begin{equation}
r_{\overline{dT}}\left(  u\right)  >r_{\overline{S}}\left(  u\right)  .
\label{pf.prop.youngtab.alt.gen.but-g}%
\end{equation}

Let $\left(  i,j\right)  $ be the cell of the $n$-tableau $S$ that contains
the entry $u$. Thus, $\left(  i,j\right)  \in Y\left(  \lambda\right)  $ and
$S\left(  i,j\right)  =u$. Hence, $r_{\overline{S}}\left(  u\right)  =i$
(since $S\left(  i,j\right)  =u$ shows that the number $u$ is contained in the
$i$-th row of the $n$-tableau $S$, thus in the $i$-th row of the $n$-tabloid
$\overline{S}$).

Let $T^{\prime}:=dT$. This is an $n$-tableau of shape $Y\left(  \mu\right)  $
(since $T$ is an $n$-tableau of shape $Y\left(  \mu\right)  $).

Let $\left(  i^{\prime},j^{\prime}\right)  $ be the cell of the $n$-tableau
$T^{\prime}$ that contains the entry $u$. Thus, $\left(  i^{\prime},j^{\prime
}\right)  \in Y\left(  \mu\right)  $ and $T^{\prime}\left(  i^{\prime
},j^{\prime}\right)  =u$. Hence, $r_{\overline{T^{\prime}}}\left(  u\right)
=i^{\prime}$ (since $T^{\prime}\left(  i^{\prime},j^{\prime}\right)  =u$ shows
that the number $u$ is contained in the $i^{\prime}$-th row of the $n$-tableau
$T^{\prime}$, thus in the $i^{\prime}$-th row of the $n$-tabloid
$\overline{T^{\prime}}$). But $T^{\prime}=dT$. Hence, $r_{\overline{T^{\prime
}}}\left(  u\right)  =r_{\overline{dT}}\left(  u\right)  >r_{\overline{S}%
}\left(  u\right)  $ (by (\ref{pf.prop.youngtab.alt.gen.but-g})), so that
$r_{\overline{S}}\left(  u\right)  <r_{\overline{T^{\prime}}}\left(  u\right)
=i^{\prime}$.

....

We now claim the following:

\begin{statement}
\textit{Claim 1:} There exists some $\ell\in\left[  i\right]  $ such that
$T^{\prime}\left(  \ell,j^{\prime}\right)  \notin K$ or $r_{\overline{S}%
}\left(  T^{\prime}\left(  \ell,j^{\prime}\right)  \right)  \geq i^{\prime}$.
\end{statement}

....

Now, we shall tweak this permutation $d$ a little bit to make the inequality
(\ref{pf.prop.youngtab.alt.gen.IH}) true not only for each $k\in
K\setminus\left\{  g\right\}  $ but also for $k=g$ (and thus for each $k\in
K$). To this purpose, we will find a transposition $t_{u,v}\in\mathcal{C}%
\left(  dT\right)  $ such that $r_{\overline{t_{u,v}dT}}\left(  g\right)  \leq
r_{\overline{S}}\left(  g\right)  $ but all $k\in K\setminus\left\{
g\right\}  $ satisfy $r_{\overline{t_{u,v}dT}}\left(  k\right)  =r_{\overline
{dT}}\left(  k\right)  $.

Here is how we find this transposition:....
\end{proof}
\end{noncompile}

\begin{fineprint}
\begin{remark}
Theorem \ref{thm.youngtab.alt} can also be proved in several other ways:

\begin{itemize}
\item The most widely known proof of Theorem \ref{thm.youngtab.alt}
\textbf{(a)} constructs the required permutation $c\in\mathcal{C}\left(
T\right)  $ iteratively, by successively raising those entries of $T$ that lie
further south in $T$ than they do in $S$ so that they no longer do so. This
proof appears, e.g., in \cite[proof of Lemma 5.13.1]{EGHetc11}, \cite[\S 7.1,
proof of Lemma 1]{Fulton97}, \cite[proof of Lemma 1.16]{BrMaPe16},
\cite[\S 9.2.3, proof of Lemma]{Proces07}, \cite[proof of Proposition
IV.2.7]{Morel19}, \cite[Lemma (4.2.A)]{Weyl53} and \cite[proof of Lemma
1.5.7]{JamKer81} (usually implicitly, since these sources only prove Theorem
\ref{thm.youngtab.alt} \textbf{(c)}, but the same argument can be used for
part \textbf{(a)}). Several variants of this argument are possible, depending
on how one chooses the order in which the entries are raised (most sources
work in the order of increasing $S$-depth, but actually any order can be
used). In a sense, our construction of $w$ in our above proof is doing the
same thing, but it is doing so in one swoop rather than iteratively.

\item Garsia proves parts \textbf{(b)} and \textbf{(c)} of Theorem
\ref{thm.youngtab.alt} rather creatively in \cite[\S 1.2]{GarEge20}, using an
extra $n$-tableau $S\wedge T$ (which he calls $T_{1}\wedge T_{2}$) that is
defined as follows: Its shape is the set $D$ of all pairs $\left(  i,j\right)
\in\left\{  1,2,3,\ldots\right\}  ^{2}$ such that the $i$-th row of $S$ has an
entry in common with the $j$-th column of $T$. For any such cell $\left(
i,j\right)  \in D$, the entry of $S\wedge T$ in this cell is defined to be the
unique common entry of the $i$-th row of $S$ and the $j$-th column of $T$ (the
uniqueness here is guaranteed by the \textquotedblleft there are no two
distinct integers...\textquotedblright\ assumption of Theorem
\ref{thm.youngtab.alt}). Note that $S\wedge T$ is a bad-shaped tableau in
general, but its very existence already yields Theorem \ref{thm.youngtab.alt}
\textbf{(b)} quite easily, and with some more work yields Theorem
\ref{thm.youngtab.alt} \textbf{(c)} as well.

\item Another approach to Theorem \ref{thm.youngtab.alt} \textbf{(a)} is to
construct the desired permutation $c\in\mathcal{C}\left(  T\right)  $ as a
solution to a maximization problem: Namely, we define a number%
\[
\omega\left(  c\right)  :=\sum_{k=1}^{n}r_{\overline{S}}\left(  k\right)
r_{\overline{cT}}\left(  k\right)  \in\mathbb{N}\ \ \ \ \ \ \ \ \ \ \left(
\text{see Definition \ref{def.tabloid.llo} \textbf{(a)} for the notation}%
\right)
\]
for each $c\in\mathcal{C}\left(  T\right)  $, and we pick a permutation
$c\in\mathcal{C}\left(  T\right)  $ that maximizes this number. We can then
check that this $c$ has the property desired in Theorem \ref{thm.youngtab.alt}
\textbf{(a)} (since otherwise, we can find two distinct entries $u$ and $v$
lying in the same column of $T$ such that $\omega\left(  t_{u,v}c\right)
>\omega\left(  c\right)  $).
\end{itemize}
\end{remark}
\end{fineprint}

\begin{exercise}
\fbox{1} Prove that Theorem \ref{thm.youngtab.alt} \textbf{(c)} has a
converse: If there exist permutations $r\in\mathcal{R}\left(  S\right)  $ and
$c\in\mathcal{C}\left(  T\right)  $ such that $rS=cT$, then there are no two
distinct integers that are horizontally aligned in $S$ and simultaneously
vertically aligned in $T$.
\end{exercise}

\begin{exercise}
\fbox{2} \textbf{(a)} Prove that the two permutations $r$ and $c$ in Theorem
\ref{thm.youngtab.alt} \textbf{(c)} are unique. \medskip

\fbox{?} \textbf{(b)} Under what conditions is the $c$ in Theorem
\ref{thm.youngtab.alt} \textbf{(a)} unique? Are there any good criteria that
guarantee uniqueness?
\end{exercise}

\begin{exercise}
\fbox{1} In Theorem \ref{thm.youngtab.alt} \textbf{(a)}, can we replace
$Y\left(  \lambda\right)  $ by a skew Young diagram? by an arbitrary diagram
contained in $\left\{  1,2,3,\ldots\right\}  ^{2}$ ?

Can we do this with $Y\left(  \mu\right)  $ ?
\end{exercise}

\begin{fineprint}
\begin{remark}
\label{rmk.youngtab.alt.n-m}Parts of Theorem \ref{thm.youngtab.alt} can be
generalized to partitions of different sizes. Indeed, let $m\in\mathbb{N}$ be
a further nonnegative integer, and let $\mu$ be a partition of $m$ (instead of
being a partition of $n$). Then, Theorem \ref{thm.youngtab.alt} remains valid
if we additionally require $n\leq m$. Moreover, Theorem \ref{thm.youngtab.alt}
\textbf{(a)} remains valid even if $n>m$ as long as we replace
\textquotedblleft Each entry of the $i$-th row of $S$\textquotedblright\ by
\textquotedblleft Each entry of the $i$-th row of $S$ that belongs to $\left[
m\right]  $\textquotedblright. This all is not too hard to show by adapting
our above proof of Theorem \ref{thm.youngtab.alt}. (For instance, we have to
define $r_{\overline{S}}\left(  k\right)  $ to be some sufficiently large
number when $k>n$, and likewise for $r_{\overline{T}}\left(  k\right)  $ when
$k>m$.) We leave the details to the reader.
\end{remark}
\end{fineprint}

\subsubsection{Algebraic consequences: the sandwich lemmas}

The Young alternative has some useful corollaries, which we (following Adriano
Garsia) call the \emph{von Neumann sandwich lemmas}. Some of them seem to
appear already in Young's original works, but John von Neumann realized how
they can be used to understand Specht modules.\footnote{John von Neumann
apparently never published on this topic, but his contributions appear in
\cite[volume 2, \S 14.7]{Waerde91} and \cite[\S 12]{Ruther48}.}

The simplest result of this type is a generalization of Theorem
\ref{thm.int.+-=0}:

\begin{lemma}
[symmetrizer cancellation lemma]\label{lem.specht.sandw0-ur}Let $S$ and $T$ be
two $n$-tableaux of arbitrary shapes. Assume that there are two distinct
integers that are horizontally aligned in $S$ and simultaneously vertically
aligned in $T$. Then,
\[
\nabla_{\operatorname*{Col}T}^{-}\nabla_{\operatorname*{Row}S}%
=0\ \ \ \ \ \ \ \ \ \ \text{and}\ \ \ \ \ \ \ \ \ \ \nabla
_{\operatorname*{Row}S}\nabla_{\operatorname*{Col}T}^{-}=0.
\]

\end{lemma}

\begin{proof}
We have assumed that there exist two distinct integers that are horizontally
aligned in $S$ and simultaneously vertically aligned in $T$. Let $i$ and $j$
be these integers.

In particular, $i$ and $j$ are vertically aligned in $T$. In other words, $i$
and $j$ lie in the same column of $T$. Hence, Proposition
\ref{prop.symmetrizers.factor-out-col} \textbf{(b)} shows that%
\[
\nabla_{\operatorname*{Col}T}^{-}=\nabla_{\operatorname*{Col}T}%
^{\operatorname*{even}}\cdot\left(  1-t_{i,j}\right)  =\left(  1-t_{i,j}%
\right)  \cdot\nabla_{\operatorname*{Col}T}^{\operatorname*{even}},
\]
where $\nabla_{\operatorname*{Col}T}^{\operatorname*{even}}$ is defined as in
Proposition \ref{prop.symmetrizers.factor-out-col} \textbf{(b)}.

But $i$ and $j$ are also horizontally aligned in $S$. In other words, $i$ and
$j$ lie in the same row of $S$. Hence, Proposition
\ref{prop.symmetrizers.factor-out-row} \textbf{(b)} (applied to $S$ instead of
$T$) shows that%
\[
\nabla_{\operatorname*{Row}S}=\nabla_{\operatorname*{Row}S}%
^{\operatorname*{even}}\cdot\left(  1+t_{i,j}\right)  =\left(  1+t_{i,j}%
\right)  \cdot\nabla_{\operatorname*{Row}S}^{\operatorname*{even}},
\]
where $\nabla_{\operatorname*{Row}S}^{\operatorname*{even}}$ is defined like
the $\nabla_{\operatorname*{Row}T}^{\operatorname*{even}}$ in Proposition
\ref{prop.symmetrizers.factor-out-row} \textbf{(b)} (but using $S$ instead of
$T$). Thus,%
\[
\underbrace{\nabla_{\operatorname*{Col}T}^{-}}_{=\nabla_{\operatorname*{Col}%
T}^{\operatorname*{even}}\cdot\left(  1-t_{i,j}\right)  }%
\ \ \underbrace{\nabla_{\operatorname*{Row}S}}_{=\left(  1+t_{i,j}\right)
\cdot\nabla_{\operatorname*{Row}S}^{\operatorname*{even}}}=\nabla
_{\operatorname*{Col}T}^{\operatorname*{even}}\cdot\underbrace{\left(
1-t_{i,j}\right)  \left(  1+t_{i,j}\right)  }_{\substack{=1-t_{i,j}%
^{2}=0\\\text{(since }t_{i,j}^{2}=\operatorname*{id}=1\text{)}}}\cdot
\,\nabla_{\operatorname*{Row}S}^{\operatorname*{even}}=0
\]
and%
\[
\underbrace{\nabla_{\operatorname*{Row}S}}_{=\nabla_{\operatorname*{Row}%
S}^{\operatorname*{even}}\cdot\left(  1+t_{i,j}\right)  }%
\ \ \underbrace{\nabla_{\operatorname*{Col}T}^{-}}_{=\left(  1-t_{i,j}\right)
\cdot\nabla_{\operatorname*{Col}T}^{\operatorname*{even}}}=\nabla
_{\operatorname*{Row}S}^{\operatorname*{even}}\cdot\underbrace{\left(
1+t_{i,j}\right)  \left(  1-t_{i,j}\right)  }_{\substack{=1-t_{i,j}%
^{2}=0\\\text{(since }t_{i,j}^{2}=\operatorname*{id}=1\text{)}}}\cdot
\,\nabla_{\operatorname*{Col}T}^{\operatorname*{even}}=0.
\]
Thus, Lemma \ref{lem.specht.sandw0-ur} is proved.
\end{proof}

Now, we can state the first real \textquotedblleft sandwich
lemma\textquotedblright:

\begin{theorem}
[von Neumann sandwich lemma, version 1]\label{thm.specht.sandw1}Let $\lambda$
be a partition of $n$. Let $T$ be an $n$-tableau of shape $Y\left(
\lambda\right)  $. Let $w\in S_{n}$. Then: \medskip

\textbf{(a)} If there exist two distinct integers that are horizontally
aligned in $wT$ and simultaneously vertically aligned in $T$, then
\[
\nabla_{\operatorname*{Col}T}^{-}w\nabla_{\operatorname*{Row}T}=0.
\]

\textbf{(b)} If not, then there exist two permutations $c\in\mathcal{C}\left(
T\right)  $ and $r\in\mathcal{R}\left(  T\right)  $ such that%
\[
w=cr\ \ \ \ \ \ \ \ \ \ \text{and}\ \ \ \ \ \ \ \ \ \ \nabla
_{\operatorname*{Col}T}^{-}w\nabla_{\operatorname*{Row}T}=\left(  -1\right)
^{c}\nabla_{\operatorname*{Col}T}^{-}\nabla_{\operatorname*{Row}T}.
\]

\end{theorem}

The name \textquotedblleft sandwich lemma\textquotedblright\ comes from
viewing the product $\nabla_{\operatorname*{Col}T}^{-}w\nabla
_{\operatorname*{Row}T}$ as \textquotedblleft sandwiching\textquotedblright%
\ the $w$ between the two bread-pieces $\nabla_{\operatorname*{Col}T}^{-}$ and
$\nabla_{\operatorname*{Row}T}$. Despite this name, we have chosen to declare
it a theorem due to its significance.

\begin{proof}
[Proof of Theorem \ref{thm.specht.sandw1}.]\textbf{(a)} Assume that there
exist two distinct integers that are horizontally aligned in $wT$ and
simultaneously vertically aligned in $T$. Then, Lemma
\ref{lem.specht.sandw0-ur} (applied to $S=wT$) yields
\begin{equation}
\nabla_{\operatorname*{Col}T}^{-}\nabla_{\operatorname*{Row}\left(  wT\right)
}=0. \label{pf.thm.specht.sandw1.a.1}%
\end{equation}

But Proposition \ref{prop.symmetrizers.conj} yields $\nabla
_{\operatorname*{Row}\left(  w\rightharpoonup T\right)  }=w\nabla
_{\operatorname*{Row}T}w^{-1}$. In other words, $\nabla_{\operatorname*{Row}%
\left(  wT\right)  }=w\nabla_{\operatorname*{Row}T}w^{-1}$ (since
$w\rightharpoonup T=wT$). Hence, we can rewrite
(\ref{pf.thm.specht.sandw1.a.1}) as%
\[
\nabla_{\operatorname*{Col}T}^{-}w\nabla_{\operatorname*{Row}T}w^{-1}=0.
\]
Multiplying this equality by $w$ on the right, we obtain $\nabla
_{\operatorname*{Col}T}^{-}w\nabla_{\operatorname*{Row}T}=0$. Thus, Theorem
\ref{thm.specht.sandw1} \textbf{(a)} is proved. \medskip

\textbf{(b)} Assume that there exist no two distinct integers that are
horizontally aligned in $wT$ and simultaneously vertically aligned in $T$.
Then, Theorem \ref{thm.youngtab.alt} \textbf{(c)} (applied to $S=wT$ and
$\mu=\lambda$) thus yields that there exist permutations $r\in\mathcal{R}%
\left(  wT\right)  $ and $c\in\mathcal{C}\left(  T\right)  $ such that
$rwT=cT$. Consider these $r$ and $c$, and let us denote them by $\widetilde{r}%
$ and $c$. Thus, $\widetilde{r}\in\mathcal{R}\left(  wT\right)  $ and
$c\in\mathcal{C}\left(  T\right)  $ and $\widetilde{r}wT=cT$. (We are using
the symbol $\widetilde{r}$ instead of $r$ because it is not the $r$ that we
are looking for.)

It is easy to see that any two permutations $u,v\in S_{n}$ satisfying $uT=vT$
must necessarily be equal\footnote{\textit{Proof.} Let $u,v\in S_{n}$ be two
permutations satisfying $uT=vT$. We must show that $u$ and $v$ are equal.
\par
Recall that $T$ is an $n$-tableau. Thus, we can view $T$ as a bijection from
$Y\left(  \lambda\right)  $ to $\left[  n\right]  $. Hence, its inverse
$T^{-1}$ is well-defined.
\par
The definition of the action on entries (Definition \ref{def.tableau.Sn-act})
yields $uT=u\circ T$ and $vT=v\circ T$. Hence, $u\circ T=uT=vT=v\circ T$.
Thus, $\underbrace{u\circ T}_{=v\circ T}\circ\,T^{-1}=v\circ\underbrace{T\circ
T^{-1}}_{=\operatorname*{id}}=v$. In view of $u\circ\underbrace{T\circ T^{-1}%
}_{=\operatorname*{id}}=u$, we can rewrite this as $u=v$. In other words, $u$
and $v$ are equal, qed.}. Applying this to $u=\widetilde{r}w$ and $v=c$, we
obtain $\widetilde{r}w=c$ (since $\widetilde{r}wT=cT$).

Proposition \ref{prop.tableau.Sn-act.0} \textbf{(c)} yields $\mathcal{R}%
\left(  w\rightharpoonup T\right)  =w\mathcal{R}\left(  T\right)  w^{-1}$. In
other words, $\mathcal{R}\left(  wT\right)  =w\mathcal{R}\left(  T\right)
w^{-1}$ (since $w\rightharpoonup T=wT$).

We have $\widetilde{r}\in\mathcal{R}\left(  wT\right)  =w\mathcal{R}\left(
T\right)  w^{-1}$. In other words, $\widetilde{r}=wr^{\prime}w^{-1}$ for some
$r^{\prime}\in\mathcal{R}\left(  T\right)  $. Consider this $r^{\prime}$.
Then, $\left(  r^{\prime}\right)  ^{-1}\in\mathcal{R}\left(  T\right)  $
(since $\mathcal{R}\left(  T\right)  $ is a group). Set $r:=\left(  r^{\prime
}\right)  ^{-1}$. Then, $r=\left(  r^{\prime}\right)  ^{-1}\in\mathcal{R}%
\left(  T\right)  $.

Hence, Proposition \ref{prop.symmetrizers.fix} \textbf{(a)} (applied to $r$
instead of $w$) yields $r\nabla_{\operatorname*{Row}T}=\nabla
_{\operatorname*{Row}T}r=\nabla_{\operatorname*{Row}T}$. Also, Proposition
\ref{prop.symmetrizers.fix} \textbf{(b)} (applied to $c$ instead of $w$)
yields $c\nabla_{\operatorname*{Col}T}^{-}=\nabla_{\operatorname*{Col}T}%
^{-}c=\left(  -1\right)  ^{c}\nabla_{\operatorname*{Col}T}^{-}$ (since
$c\in\mathcal{C}\left(  T\right)  $).

We have $\widetilde{r}w=wr^{\prime}$ (since $\widetilde{r}=wr^{\prime}w^{-1}%
$). Comparing this with $\widetilde{r}w=c$, we obtain $wr^{\prime}=c$, so that
$w=c\underbrace{\left(  r^{\prime}\right)  ^{-1}}_{=r}=cr$. Thus,%
\[
\nabla_{\operatorname*{Col}T}^{-}w\nabla_{\operatorname*{Row}T}%
=\underbrace{\nabla_{\operatorname*{Col}T}^{-}c}_{=\left(  -1\right)
^{c}\nabla_{\operatorname*{Col}T}^{-}}\underbrace{r\nabla_{\operatorname*{Row}%
T}}_{=\nabla_{\operatorname*{Row}T}}=\left(  -1\right)  ^{c}\nabla
_{\operatorname*{Col}T}^{-}\nabla_{\operatorname*{Row}T}.
\]
Thus, we have found two permutations $c\in\mathcal{C}\left(  T\right)  $ and
$r\in\mathcal{R}\left(  T\right)  $ such that $w=cr$ and $\nabla
_{\operatorname*{Col}T}^{-}w\nabla_{\operatorname*{Row}T}=\left(  -1\right)
^{c}\nabla_{\operatorname*{Col}T}^{-}\nabla_{\operatorname*{Row}T}$. This
proves Theorem \ref{thm.specht.sandw1} \textbf{(b)}.
\end{proof}

\begin{corollary}
[von Neumann sandwich lemma, version 2]\label{cor.specht.sandw2}Let $\lambda$
be a partition of $n$. Let $T$ be an $n$-tableau of shape $Y\left(
\lambda\right)  $. Let $\mathbf{a}\in\mathbf{k}\left[  S_{n}\right]  $. Then,
\[
\nabla_{\operatorname*{Col}T}^{-}\mathbf{a}\nabla_{\operatorname*{Row}%
T}=\kappa\nabla_{\operatorname*{Col}T}^{-}\nabla_{\operatorname*{Row}%
T}\ \ \ \ \ \ \ \ \ \ \text{for some scalar }\kappa\in\mathbf{k}.
\]
(Of course, $\kappa$ depends on $\mathbf{a}$.)
\end{corollary}

\begin{proof}
Theorem \ref{thm.specht.sandw1} has the following consequence: For each
permutation $w\in S_{n}$, there exists an integer $\kappa_{w}\in\mathbb{Z}$
such that%
\begin{equation}
\nabla_{\operatorname*{Col}T}^{-}w\nabla_{\operatorname*{Row}T}=\kappa
_{w}\nabla_{\operatorname*{Col}T}^{-}\nabla_{\operatorname*{Row}T}
\label{pf.thm.specht.sandw1.kappaw}%
\end{equation}
\footnote{\textit{Proof.} If there exist two distinct integers that are
horizontally aligned in $wT$ and simultaneously vertically aligned in $T$,
then we can choose $\kappa_{w}=0$, because Theorem \ref{thm.specht.sandw1}
\textbf{(a)} shows that $\nabla_{\operatorname*{Col}T}^{-}w\nabla
_{\operatorname*{Row}T}=0$ in this case. In the contrary case, Theorem
\ref{thm.specht.sandw1} \textbf{(b)} yields that $\nabla_{\operatorname*{Col}%
T}^{-}w\nabla_{\operatorname*{Row}T}=\left(  -1\right)  ^{c}\nabla
_{\operatorname*{Col}T}^{-}\nabla_{\operatorname*{Row}T}$ for some
$c\in\mathcal{C}\left(  T\right)  $ and some $r\in\mathcal{R}\left(  T\right)
$, and thus we can choose $\kappa_{w}=\left(  -1\right)  ^{c}$. In either
case, we have thus found an integer $\kappa_{w}\in\mathbb{Z}$ such that
$\nabla_{\operatorname*{Col}T}^{-}w\nabla_{\operatorname*{Row}T}=\kappa
_{w}\nabla_{\operatorname*{Col}T}^{-}\nabla_{\operatorname*{Row}T}$.}.
Consider this integer $\kappa_{w}$.

Now, let us expand $\mathbf{a}$ in the standard basis of $\mathbf{k}\left[
S_{n}\right]  $. That is, let us write $\mathbf{a}$ in the form $\mathbf{a}%
=\sum_{w\in S_{n}}\alpha_{w}w$ for some scalars $\alpha_{w}\in\mathbf{k}$.
Then,%
\begin{align*}
\nabla_{\operatorname*{Col}T}^{-}\mathbf{a}\nabla_{\operatorname*{Row}T}  &
=\nabla_{\operatorname*{Col}T}^{-}\left(  \sum_{w\in S_{n}}\alpha_{w}w\right)
\nabla_{\operatorname*{Row}T}=\sum_{w\in S_{n}}\alpha_{w}\underbrace{\nabla
_{\operatorname*{Col}T}^{-}w\nabla_{\operatorname*{Row}T}}_{\substack{=\kappa
_{w}\nabla_{\operatorname*{Col}T}^{-}\nabla_{\operatorname*{Row}T}\\\text{(by
(\ref{pf.thm.specht.sandw1.kappaw}))}}}\\
&  =\sum_{w\in S_{n}}\alpha_{w}\kappa_{w}\nabla_{\operatorname*{Col}T}%
^{-}\nabla_{\operatorname*{Row}T}=\left(  \sum_{w\in S_{n}}\alpha_{w}%
\kappa_{w}\right)  \nabla_{\operatorname*{Col}T}^{-}\nabla
_{\operatorname*{Row}T}.
\end{align*}
Thus, $\nabla_{\operatorname*{Col}T}^{-}\mathbf{a}\nabla_{\operatorname*{Row}%
T}=\kappa\nabla_{\operatorname*{Col}T}^{-}\nabla_{\operatorname*{Row}T}$ for
some scalar $\kappa\in\mathbf{k}$ (namely, for $\kappa=\sum_{w\in S_{n}}%
\alpha_{w}\kappa_{w}$). This proves Corollary \ref{cor.specht.sandw2}.
\end{proof}

\begin{exercise}
\fbox{2} Generalize Corollary \ref{cor.specht.sandw2} to $\nabla
_{\operatorname*{Col}T}^{-}\mathbf{a}\nabla_{\operatorname*{Row}S}$, where $T$
and $S$ are two $n$-tableaux of the same shape $Y\left(  \lambda\right)  $.
\medskip

[\textbf{Hint:} The claim will have the form $\nabla_{\operatorname*{Col}%
T}^{-}\mathbf{a}\nabla_{\operatorname*{Row}S}=\kappa\mathbf{b}$ for some
scalar $\kappa\in\mathbf{k}$ and some fixed $\mathbf{b}\in\mathbf{k}\left[
S_{n}\right]  $ that is independent of $\mathbf{a}$. But $\mathbf{b}$ will not
be $\nabla_{\operatorname*{Col}T}^{-}\nabla_{\operatorname*{Row}S}$ this time.]
\end{exercise}

\begin{exercise}
Prove that the scalar $\kappa$ in Corollary \ref{cor.specht.sandw2} can be
computed by the following two explicit formulas: \medskip

\textbf{(a)} \fbox{2} If $\mathbf{a}=\sum_{w\in S_{n}}\alpha_{w}w$ for some
scalars $\alpha_{w}\in\mathbf{k}$, then
\[
\kappa=\sum_{c\in\mathcal{C}\left(  T\right)  }\ \ \sum_{r\in\mathcal{R}%
\left(  T\right)  }\left(  -1\right)  ^{c}\alpha_{cr}.
\]

\textbf{(b)} \fbox{2} If $\mathbf{a}\nabla_{\operatorname*{Row}T}%
\nabla_{\operatorname*{Col}T}^{-}=\sum_{w\in S_{n}}\beta_{w}w$ for some
scalars $\beta_{w}\in\mathbf{k}$, then $\kappa=\beta_{\operatorname*{id}}$.
\end{exercise}

We furthermore will use an analogue of the von Neumann sandwich lemma which
involves two different shapes:

\begin{theorem}
[von Neumann zero-sandwich lemma]\label{thm.specht.sandw0}Let $\lambda$ and
$\mu$ be two partitions of $n$. Assume that $\lambda_{1}+\lambda_{2}%
+\cdots+\lambda_{i}>\mu_{1}+\mu_{2}+\cdots+\mu_{i}$ for some $i\geq1$.

Let $S$ be an $n$-tableau of shape $Y\left(  \lambda\right)  $. Let $T$ be an
$n$-tableau of shape $Y\left(  \mu\right)  $. Then: \medskip

\textbf{(a)} We have $\nabla_{\operatorname*{Col}T}^{-}w\nabla
_{\operatorname*{Row}S}=0$ for each $w\in S_{n}$. \medskip

\textbf{(b)} We have $\nabla_{\operatorname*{Col}T}^{-}\mathbf{a}%
\nabla_{\operatorname*{Row}S}=0$ for each $\mathbf{a}\in\mathbf{k}\left[
S_{n}\right]  $. \medskip

\textbf{(c)} We have $\nabla_{\operatorname*{Row}S}w\nabla
_{\operatorname*{Col}T}^{-}=0$ for each $w\in S_{n}$. \medskip

\textbf{(d)} We have $\nabla_{\operatorname*{Row}S}\mathbf{a}\nabla
_{\operatorname*{Col}T}^{-}=0$ for each $\mathbf{a}\in\mathbf{k}\left[
S_{n}\right]  $.
\end{theorem}

\begin{proof}
\textbf{(a)} Let $w\in S_{n}$. Then, $wS$ is an $n$-tableau of shape $Y\left(
\lambda\right)  $ (since $S$ is such an $n$-tableau). Hence, by Theorem
\ref{thm.youngtab.alt} \textbf{(b)} (or, rather, by its contrapositive), we
can see that there are two distinct integers that are horizontally aligned in
$wS$ and simultaneously vertically aligned in $T$%
\ \ \ \ \footnote{\textit{Proof.} Assume the contrary. Thus, there are no two
distinct integers that are horizontally aligned in $wS$ and simultaneously
vertically aligned in $T$. Hence, Theorem \ref{thm.youngtab.alt} \textbf{(b)}
(applied to $wS$ instead of $S$) shows that we have $\lambda_{1}+\lambda
_{2}+\cdots+\lambda_{i}\leq\mu_{1}+\mu_{2}+\cdots+\mu_{i}$ for each $i\geq1$.
But this contradicts the fact that $\lambda_{1}+\lambda_{2}+\cdots+\lambda
_{i}>\mu_{1}+\mu_{2}+\cdots+\mu_{i}$ for some $i\geq1$. This contradiction
shows that our assumption was false, qed.}. Hence, Lemma
\ref{lem.specht.sandw0-ur} (applied to $wS$ instead of $S$) yields%
\[
\nabla_{\operatorname*{Col}T}^{-}\nabla_{\operatorname*{Row}\left(  wS\right)
}=0\ \ \ \ \ \ \ \ \ \ \text{and}\ \ \ \ \ \ \ \ \ \ \nabla
_{\operatorname*{Row}\left(  wS\right)  }\nabla_{\operatorname*{Col}T}^{-}=0.
\]
But Proposition \ref{prop.symmetrizers.conj} (applied to $S$ instead of $T$)
yields $\nabla_{\operatorname*{Row}\left(  w\rightharpoonup S\right)
}=w\nabla_{\operatorname*{Row}S}w^{-1}$. In other words, $\nabla
_{\operatorname*{Row}\left(  wS\right)  }=w\nabla_{\operatorname*{Row}S}%
w^{-1}$ (since $w\rightharpoonup S=wS$). Hence, $\nabla_{\operatorname*{Row}%
\left(  wS\right)  }w=w\nabla_{\operatorname*{Row}S}$ and therefore%
\[
\nabla_{\operatorname*{Col}T}^{-}\underbrace{w\nabla_{\operatorname*{Row}S}%
}_{=\nabla_{\operatorname*{Row}\left(  wS\right)  }w}=\underbrace{\nabla
_{\operatorname*{Col}T}^{-}\nabla_{\operatorname*{Row}\left(  wS\right)  }%
}_{=0}w=0.
\]

This proves Theorem \ref{thm.specht.sandw0} \textbf{(a)}. \medskip

\textbf{(b)} This follows from part \textbf{(a)} by linearity.\footnote{In
more detail: We must prove the equality $\nabla_{\operatorname*{Col}T}%
^{-}\mathbf{a}\nabla_{\operatorname*{Row}S}=0$ for each $\mathbf{a}%
\in\mathbf{k}\left[  S_{n}\right]  $. Since both sides of this equality are
$\mathbf{k}$-linear in $\mathbf{a}$, it suffices to prove it when $\mathbf{a}$
is a standard basis vector of $\mathbf{k}\left[  S_{n}\right]  $ (since the
standard basis vectors span the $\mathbf{k}$-module $\mathbf{k}\left[
S_{n}\right]  $). In other words, it suffices to prove it when $\mathbf{a}=w$
for some permutation $w\in S_{n}$. In other words, it suffices to prove that
$\nabla_{\operatorname*{Col}T}^{-}w\nabla_{\operatorname*{Row}S}=0$ for each
$w\in S_{n}$. But this follows immediately from Theorem
\ref{thm.specht.sandw0} \textbf{(a)}. Thus, Theorem \ref{thm.specht.sandw0}
\textbf{(b)} is proved.} \medskip

\textbf{(c)} This can be done similarly to part \textbf{(a)}, but using the
$\nabla_{\operatorname*{Row}\left(  wS\right)  }\nabla_{\operatorname*{Col}%
T}^{-}=0$ equality instead of $\nabla_{\operatorname*{Col}T}^{-}%
\nabla_{\operatorname*{Row}\left(  wS\right)  }=0$. (More precisely, this way
we can obtain $\nabla_{\operatorname*{Row}S}w^{-1}\nabla_{\operatorname*{Col}%
T}^{-}=0$, but since $w$ is arbitrary, we can apply this to $w^{-1}$ instead
of $w$ to recover the equality we are looking for).

Alternatively, we can derive Theorem \ref{thm.specht.sandw0} \textbf{(c)} from
part \textbf{(a)} using the antipode. This technique is useful in various
situations, so let me show this proof in some detail:

Let us rename the $n$-tableau $S$ as $R$ (in order to avoid confusing it with
the antipode $S$ of $\mathbf{k}\left[  S_{n}\right]  $). Thus, we must prove
that $\nabla_{\operatorname*{Row}R}w\nabla_{\operatorname*{Col}T}^{-}=0$ for
each $w\in S_{n}$.

Fix $w\in S_{n}$. Then, Theorem \ref{thm.specht.sandw0} \textbf{(a)} (applied
to $w^{-1}$ instead of $w$) yields $\nabla_{\operatorname*{Col}T}^{-}%
w^{-1}\nabla_{\operatorname*{Row}R}=0$ (since we have renamed the $n$-tableau
$S$ as $R$). Applying the antipode $S$ to both sides of this equality, we
obtain $S\left(  \nabla_{\operatorname*{Col}T}^{-}w^{-1}\nabla
_{\operatorname*{Row}R}\right)  =S\left(  0\right)  =0$. Hence,%
\begin{align*}
0  &  =S\left(  \nabla_{\operatorname*{Col}T}^{-}w^{-1}\nabla
_{\operatorname*{Row}R}\right) \\
&  =\underbrace{S\left(  \nabla_{\operatorname*{Row}R}\right)  }%
_{\substack{=\nabla_{\operatorname*{Row}R}\\\text{(by Proposition
\ref{prop.symmetrizers.antipode},}\\\text{applied to }R\text{ instead of
}T\text{)}}}\ \ \underbrace{S\left(  w^{-1}\right)  }_{\substack{=\left(
w^{-1}\right)  ^{-1}\\\text{(by the definition of }S\text{)}}%
}\ \ \underbrace{S\left(  \nabla_{\operatorname*{Col}T}^{-}\right)
}_{\substack{=\nabla_{\operatorname*{Col}T}^{-}\\\text{(by Proposition
\ref{prop.symmetrizers.antipode})}}}\\
&  \ \ \ \ \ \ \ \ \ \ \ \ \ \ \ \ \ \ \ \ \left(  \text{since }S\text{ is a
}\mathbf{k}\text{-algebra anti-automorphism (by Theorem \ref{thm.S.auto}
\textbf{(a)})}\right) \\
&  =\nabla_{\operatorname*{Row}R}\underbrace{\left(  w^{-1}\right)  ^{-1}%
}_{=w}\nabla_{\operatorname*{Col}T}^{-}=\nabla_{\operatorname*{Row}R}%
w\nabla_{\operatorname*{Col}T}^{-}.
\end{align*}
Thus, $\nabla_{\operatorname*{Row}R}w\nabla_{\operatorname*{Col}T}^{-}=0$.
This is precisely what we wanted to prove (since we renamed the $n$-tableau
$S$ as $R$). Thus, Theorem \ref{thm.specht.sandw0} \textbf{(c)} is proved.
\medskip

\textbf{(d)} This follows from part \textbf{(c)} just like part \textbf{(b)}
follows from part \textbf{(a)}.
\end{proof}

Let us observe yet another easy consequence of Theorem
\ref{lem.specht.sandw0-ur}, connecting the row-symmetrizers and the
column-antisymmetrizers with the Young last letter order from Definition
\ref{def.tabloid.llo}:

\begin{lemma}
\label{lem.specht.sandw0-llo}Let $\lambda$ be a partition of $n$. Let $S$ and
$T$ be two $n$-tableaux of shape $Y\left(  \lambda\right)  $ such that $T$ is
column-standard. Assume that $\overline{T}<\overline{S}$ with respect to the
Young last letter order (i.e., the binary relation $<$ defined in Definition
\ref{def.tabloid.llo} \textbf{(b)}). Then: \medskip

\textbf{(a)} There exist no permutations $r\in\mathcal{R}\left(  S\right)  $
and $c\in\mathcal{C}\left(  T\right)  $ such that $rS=cT$. \medskip

\textbf{(b)} There are two distinct integers that are horizontally aligned in
$S$ and simultaneously vertically aligned in $T$. \medskip

\textbf{(c)} We have%
\[
\nabla_{\operatorname*{Col}T}^{-}\nabla_{\operatorname*{Row}S}%
=0\ \ \ \ \ \ \ \ \ \ \text{and}\ \ \ \ \ \ \ \ \ \ \nabla
_{\operatorname*{Row}S}\nabla_{\operatorname*{Col}T}^{-}=0.
\]

\end{lemma}

\begin{proof}
\textbf{(a)} Assume the contrary. Thus, there exist permutations
$r\in\mathcal{R}\left(  S\right)  $ and $c\in\mathcal{C}\left(  T\right)  $
such that $rS=cT$. Consider such $r$ and $c$. Consider the $n$-tabloids
$\overline{S}$ and $\overline{T}$. We have $r\in\mathcal{R}\left(  S\right)
=\left\{  w\in S_{n}\ \mid\ \overline{w\rightharpoonup S}=\overline
{S}\right\}  $ (by Proposition \ref{prop.tableau.Sn-act.1} \textbf{(a)},
applied to $S$ instead of $T$). In other words, $r\in S_{n}$ and
$\overline{r\rightharpoonup S}=\overline{S}$. But $rS$ is shorthand for
$r\rightharpoonup S$. Thus, $rS=r\rightharpoonup S$, so that $\overline
{rS}=\overline{r\rightharpoonup S}=\overline{S}$. In other words,
$\overline{cT}=\overline{S}$ (since $rS=cT$).

But we have $c\in\mathcal{C}\left(  T\right)  $; therefore, we obtain
$\overline{cT}\leq\overline{T}$ as an easy consequence of Lemma
\ref{lem.spechtmod.lead-term} \textbf{(a)}\footnote{\textit{Proof:} We must
show that $\overline{cT}\leq\overline{T}$. If $c\neq\operatorname*{id}$, then
Lemma \ref{lem.spechtmod.lead-term} \textbf{(a)} (applied to $w=c$) yields
that $\overline{cT}<\overline{T}$, and thus $\overline{cT}\leq\overline{T}$
follows. If $c=\operatorname*{id}$, then $cT=\operatorname*{id}T=T$ and thus
$\overline{cT}=\overline{T}$, so that $\overline{cT}\leq\overline{T}$. Thus,
in either case, we have shown that $\overline{cT}\leq\overline{T}$.}. In view
of $\overline{cT}=\overline{S}$, we can rewrite this as $\overline{S}%
\leq\overline{T}$. But this contradicts $\overline{T}<\overline{S}$. This
contradiction shows that our assumption was false. Thus, the proof of Lemma
\ref{lem.specht.sandw0-llo} \textbf{(a)} is complete. \medskip

\textbf{(b)} Assume the contrary. Thus, there are no two distinct integers
that are horizontally aligned in $S$ and simultaneously vertically aligned in
$T$. Hence, Theorem \ref{thm.youngtab.alt} \textbf{(c)} (applied to
$\mu=\lambda$) yields that there exist permutations $r\in\mathcal{R}\left(
S\right)  $ and $c\in\mathcal{C}\left(  T\right)  $ such that $rS=cT$. But
this contradicts Lemma \ref{lem.specht.sandw0-llo} \textbf{(a)}. This
contradiction shows that our assumption was false. Thus, the proof of Lemma
\ref{lem.specht.sandw0-llo} \textbf{(b)} is complete. \medskip

\textbf{(c)} Lemma \ref{lem.specht.sandw0-llo} \textbf{(b)} shows that there
are two distinct integers that are horizontally aligned in $S$ and
simultaneously vertically aligned in $T$. Hence, Lemma
\ref{lem.specht.sandw0-ur} yields
\[
\nabla_{\operatorname*{Col}T}^{-}\nabla_{\operatorname*{Row}S}%
=0\ \ \ \ \ \ \ \ \ \ \text{and}\ \ \ \ \ \ \ \ \ \ \nabla
_{\operatorname*{Row}S}\nabla_{\operatorname*{Col}T}^{-}=0.
\]
This proves Lemma \ref{lem.specht.sandw0-llo} \textbf{(c)}.
\end{proof}

\subsection{\label{sec.specht.ysym}Young symmetrizers}

\subsubsection{Definitions}

The sandwich lemmas will help us understand the Specht modules $\mathcal{S}%
^{Y\left(  \lambda\right)  }$ better. First, we introduce a shorthand for them:

\begin{definition}
\label{def.specht.ET.defs}\textbf{(a)} Given any partition $\lambda$ of size
$n$, we set $\mathcal{S}^{\lambda}:=\mathcal{S}^{Y\left(  \lambda\right)  }$
(the Specht module for shape $Y\left(  \lambda\right)  $). \medskip

\textbf{(b)} For any $n$-tableau $T$ (of any shape), we set%
\[
\mathbf{E}_{T}:=\nabla_{\operatorname*{Col}T}^{-}\nabla_{\operatorname*{Row}%
T}\in\mathbf{k}\left[  S_{n}\right]  .
\]
The element $\mathbf{E}_{T}$ is commonly called a \emph{Young symmetrizer}.
\medskip

\textbf{(c)} Furthermore, we set $\mathcal{A}:=\mathbf{k}\left[  S_{n}\right]
$ for the rest of this chapter.
\end{definition}

The straight-shaped Specht modules $\mathcal{S}^{\lambda}=\mathcal{S}%
^{Y\left(  \lambda\right)  }$ are a particular case of the skew-shaped Specht
modules $\mathcal{S}^{Y\left(  \lambda/\mu\right)  }$ (since $Y\left(
\lambda\right)  =Y\left(  \lambda/\varnothing\right)  $ for any partition
$\lambda$). As such, they are fairly well-understood (in particular, Theorem
\ref{thm.spechtmod.basis} provides a basis for each of them). However, they
have a special role to play: If $\mathbf{k}$ is a field of characteristic $0$,
then they are precisely the irreps (= irreducible representations) of $S_{n}$
over $\mathbf{k}$.

To prove this will take us some work. First, we recall the left ideal avatar
of a Specht module (Theorem \ref{thm.spechtmod.leftideal} \textbf{(b)}): If
$D$ is any diagram, and if $T$ is an $n$-tableau of shape $D$, then
\begin{align}
\mathcal{S}^{D}  &  \cong\underbrace{\mathbf{k}\left[  S_{n}\right]
}_{\substack{=\mathcal{A}\\\text{(by Definition \ref{def.specht.ET.defs}
\textbf{(c)})}}}\cdot\underbrace{\nabla_{\operatorname*{Col}T}^{-}%
\nabla_{\operatorname*{Row}T}}_{\substack{=\mathbf{E}_{T}\\\text{(by
Definition \ref{def.specht.ET.defs} \textbf{(b)})}}%
}\ \ \ \ \ \ \ \ \ \ \left(  \text{by Theorem \ref{thm.spechtmod.leftideal}
\textbf{(b)}}\right) \nonumber\\
&  =\mathcal{A}\mathbf{E}_{T} \label{eq.def.specht.ET.defs.SD=AET}%
\end{align}
as left $\mathbf{k}\left[  S_{n}\right]  $-modules (i.e., as $S_{n}%
$-representations). Thus, we can hope to understand $\mathcal{S}^{D}$ by
understanding the Young symmetrizer $\mathbf{E}_{T}$. This is true for any
diagram $D$, but will be particularly useful in the case when $D$ is a
straight Young diagram $Y\left(  \lambda\right)  $.

\begin{example}
\label{exa.specht.ET.n=3}Let $n=3$. Then, each $n$-tableau of straight shape
has one of the three forms%
\[
\ytableaushort{ijk}\qquad\text{and}\qquad\ytableaushort{i,j,k}\qquad
\text{and}\qquad\ytableaushort{ij,k}\ \ ,
\]
where $i,j,k$ are the three elements of $\left[  n\right]  $ in some order.
Using the shorthand notation from Convention \ref{conv.tableau.poetic}, we can
rewrite these three forms as $ijk$ and $i\backslash\backslash j\backslash
\backslash k$ and $ij\backslash\backslash k$, respectively. Thus, the Young
symmetrizers for straight-shaped $n$-tableaux are the elements%
\begin{align*}
\mathbf{E}_{ijk}  &  =\underbrace{\nabla_{\operatorname*{Col}\left(
ijk\right)  }^{-}}_{=1}\underbrace{\nabla_{\operatorname*{Row}\left(
ijk\right)  }}_{=\nabla}=\nabla=1+s_{1}+s_{2}+t_{1,3}+\operatorname*{cyc}%
\nolimits_{1,2,3}+\operatorname*{cyc}\nolimits_{1,3,2};\\
\mathbf{E}_{i\backslash\backslash j\backslash\backslash k}  &
=\underbrace{\nabla_{\operatorname*{Col}\left(  i\backslash\backslash
j\backslash\backslash k\right)  }^{-}}_{=\nabla^{-}}\underbrace{\nabla
_{\operatorname*{Row}\left(  i\backslash\backslash j\backslash\backslash
k\right)  }}_{=1}=\nabla^{-}=1-s_{1}-s_{2}-t_{1,3}+\operatorname*{cyc}%
\nolimits_{1,2,3}+\operatorname*{cyc}\nolimits_{1,3,2};\\
\mathbf{E}_{ij\backslash\backslash k}  &  =\underbrace{\nabla
_{\operatorname*{Col}\left(  ij\backslash\backslash k\right)  }^{-}%
}_{=1-t_{i,k}}\underbrace{\nabla_{\operatorname*{Row}\left(  ij\backslash
\backslash k\right)  }}_{=1+t_{i,j}}=\left(  1-t_{i,k}\right)  \left(
1+t_{i,j}\right)  =1-t_{i,k}+t_{i,j}-\operatorname*{cyc}\nolimits_{i,j,k}%
\end{align*}
for all such orderings $i,j,k$. (You are reading right: Both $\mathbf{E}%
_{ijk}$ and $\mathbf{E}_{i\backslash\backslash j\backslash\backslash k}$ are
independent of $i,j,k$.)
\end{example}

\begin{exercise}
\textbf{(a)} \fbox{1} Prove that $\mathbf{E}_{123\backslash\backslash
456}=\mathbf{E}_{465\backslash\backslash132}$ for $n=6$, without fully
expanding either side of the equality. \medskip

\textbf{(b)} \fbox{4} Let $P$ and $Q$ be two $n$-tableaux of any shapes (not
necessarily of the same shape!). Assume that $2\neq0$ in $\mathbf{k}$. Prove
that $\mathbf{E}_{P}=\mathbf{E}_{Q}$ holds if and only if the equalities
$\mathcal{R}\left(  P\right)  =\mathcal{R}\left(  Q\right)  $ and
$\mathcal{C}\left(  P\right)  =\mathcal{C}\left(  Q\right)  $ hold. \medskip

[\textbf{Hint:} In part \textbf{(b)}, consider the coefficients with which
transpositions occur in $\mathbf{E}_{P}$ and $\mathbf{E}_{Q}$.]
\end{exercise}

\subsubsection{\label{subsec.specht.ysym.qi}Quasi-idempotency}

When $D$ is a Young diagram, the Young symmetrizer $\mathbf{E}_{T}$ is
particularly nice -- it is what is called a \emph{quasi-idempotent}:

\begin{theorem}
[Young symmetrizer theorem]\label{thm.specht.ETidp}Let $\lambda$ be a
partition of $n$. Let $T$ be an $n$-tableau of shape $Y\left(  \lambda\right)
$. Let $f^{\lambda}$ be the \# of standard tableaux of shape $Y\left(
\lambda\right)  $. Then, $\dfrac{n!}{f^{\lambda}}$ is a positive integer, and
we have%
\[
\mathbf{E}_{T}^{2}=\dfrac{n!}{f^{\lambda}}\mathbf{E}_{T}.
\]
(Recall that $\mathbf{E}_{T}=\nabla_{\operatorname*{Col}T}^{-}\nabla
_{\operatorname*{Row}T}$, as we said in Definition \ref{def.specht.ET.defs}
\textbf{(b)}.)
\end{theorem}

\begin{example}
Let $n=3$ and $\lambda=\left(  2,1\right)  $ and $T=\ytableaushort{12,3}$\ .
Then,%
\[
\mathbf{E}_{T}=\underbrace{\nabla_{\operatorname*{Col}T}^{-}}_{=1-t_{1,3}%
}\ \ \underbrace{\nabla_{\operatorname*{Row}T}}_{=1+t_{1,2}}=\left(
1-t_{1,3}\right)  \left(  1+t_{1,2}\right)  =1-t_{1,3}+t_{1,2}%
-\operatorname*{cyc}\nolimits_{1,2,3}.
\]
So $\mathbf{E}_{T}^{2}$ is a sum of $16$ terms. Theorem \ref{thm.specht.ETidp}
tells us that this sum simplifies to%
\[
\dfrac{n!}{f^{\lambda}}\mathbf{E}_{T}=\dfrac{3!}{2}\mathbf{E}_{T}%
=3\mathbf{E}_{T}=3\left(  1-t_{1,3}+t_{1,2}-\operatorname*{cyc}%
\nolimits_{1,2,3}\right)  .
\]

\end{example}

Note that Theorem \ref{thm.specht.ETidp} is not (generally) true when
$Y\left(  \lambda\right)  $ is replaced by a skew Young diagram; in this case,
$\mathbf{E}_{T}^{2}$ is usually not a scalar multiple of $\mathbf{E}_{T}%
$.\ \ \ \ \footnote{The simplest counterexample is obtained when
$\mathbf{k}=\mathbb{Q}$ and $n=4$ and $T$ is a $4$-tableau of shape $Y\left(
\left(  3,2\right)  /\left(  1\right)  \right)  $. In this case,
$\mathbf{E}_{T}^{2}$ is not a scalar multiple of $\mathbf{E}_{T}$, and there
is no nontrivial quadratic polynomial that vanishes at $\mathbf{E}_{T}$.
Rather, the minimal polynomial of $\mathbf{E}_{T}$ (that is, the monic
polynomial $P\in\mathbb{Q}\left[  x\right]  $ of smallest degree that
satisfies $P\left(  \mathbf{E}_{T}\right)  =0$) is $\left(  x-6\right)
\left(  x-4\right)  x$ in this case.
\par
This example might suggest that in general, the minimal polynomial of
$\mathbf{E}_{T}$ might at least be split over $\mathbb{Z}$ (that is, all its
roots are integers). Alas, this pattern is not sustained either. For example,
if $\mathbf{k}=\mathbb{Q}$ and $n=7$ and if $T$ is a $7$-tableau of shape
$Y\left(  \left(  5,3,2\right)  /\left(  2,1\right)  \right)  $, then the
minimal polynomial of $\mathbf{E}_{T}$ is%
\[
\underbrace{\left(  x^{2}-62x+864\right)  }_{\text{irreducible over
}\mathbb{Q}}\left(  x-48\right)  \left(  x-38\right)  \left(  x-32\right)
\left(  x-14\right)  x.
\]
} \medskip

The proof of Theorem \ref{thm.specht.ETidp} is rather surprising: It is very
easy (using the von Neumann sandwich lemma) to show that $\mathbf{E}_{T}%
^{2}=\kappa\mathbf{E}_{T}$ for some $\kappa\in\mathbf{k}$, but it is much
harder to identify this $\kappa$ as $\dfrac{n!}{f^{\lambda}}$. We will thus
approach this proof slowly, using several lemmas that seem unrelated yet will
all find their use. First, we handle the sandwich part, in fact proving a more
general version thereof:

\begin{proposition}
[double sandwich lemma]\label{prop.specht.ET.ETaET}Let $\lambda$ be a
partition of $n$. Let $T$ be an $n$-tableau of shape $Y\left(  \lambda\right)
$. Let $\mathbf{a}\in\mathbf{k}\left[  S_{n}\right]  $. Then, $\mathbf{E}%
_{T}\mathbf{aE}_{T}=\kappa\mathbf{E}_{T}$ for some scalar $\kappa\in
\mathbf{k}$.
\end{proposition}

\begin{proof}
We have $\mathbf{E}_{T}=\nabla_{\operatorname*{Col}T}^{-}\nabla
_{\operatorname*{Row}T}$. Hence,%
\begin{align*}
\mathbf{E}_{T}\mathbf{aE}_{T}  &  =\nabla_{\operatorname*{Col}T}^{-}%
\nabla_{\operatorname*{Row}T}\mathbf{a}\nabla_{\operatorname*{Col}T}^{-}%
\nabla_{\operatorname*{Row}T}=\nabla_{\operatorname*{Col}T}^{-}\left(
\nabla_{\operatorname*{Row}T}\mathbf{a}\nabla_{\operatorname*{Col}T}%
^{-}\right)  \nabla_{\operatorname*{Row}T}\\
&  =\kappa\nabla_{\operatorname*{Col}T}^{-}\nabla_{\operatorname*{Row}%
T}\ \ \ \ \ \ \ \ \ \ \text{for some scalar }\kappa\in\mathbf{k}%
\end{align*}
(by Corollary \ref{cor.specht.sandw2}, applied to $\nabla_{\operatorname*{Row}%
T}\mathbf{a}\nabla_{\operatorname*{Col}T}^{-}$ instead of $\mathbf{a}$). In
other words, $\mathbf{E}_{T}\mathbf{aE}_{T}=\kappa\mathbf{E}_{T}$ for some
scalar $\kappa\in\mathbf{k}$ (since $\nabla_{\operatorname*{Col}T}^{-}%
\nabla_{\operatorname*{Row}T}=\mathbf{E}_{T}$). This proves Proposition
\ref{prop.specht.ET.ETaET}.
\end{proof}

Clearly, applying Proposition \ref{prop.specht.ET.ETaET} to $\mathbf{a}=1$
yields $\mathbf{E}_{T}1\mathbf{E}_{T}=\kappa\mathbf{E}_{T}$, that is,
$\mathbf{E}_{T}^{2}=\kappa\mathbf{E}_{T}$. Now, we \textquotedblleft
only\textquotedblright\ need to show that this $\kappa$ is $\dfrac
{n!}{f^{\lambda}}$. As we said, this will require a significant amount of
preparation. We begin by introducing a notation:

\begin{definition}
\label{def.specht.ET.wa}Let $G$ be a group. Let $w\in G$ and $\mathbf{a}%
\in\mathbf{k}\left[  G\right]  $. Then, $\left[  w\right]  \mathbf{a}$ denotes
the coefficient of $w$ in $\mathbf{a}$ (when $\mathbf{a}$ is expanded as a
linear combination of the standard basis vectors of $\mathbf{k}\left[
G\right]  $). That is, if $\mathbf{a}=\sum\limits_{g\in G}\alpha_{g}g$ for
some scalars $\alpha_{g}\in\mathbf{k}$, then $\left[  w\right]  \mathbf{a}%
:=\alpha_{w}$. The scalar $\left[  w\right]  \mathbf{a}$ is also known as the
$w$\emph{-coordinate} of $\mathbf{a}$.
\end{definition}

For instance, the third jucys--murphy $\mathbf{m}_{3}\in\mathbf{k}\left[
S_{n}\right]  $ equals $t_{1,3}+s_{2}$, so that we have%
\[
\left[  s_{2}\right]  \mathbf{m}_{3}=1\ \ \ \ \ \ \ \ \ \ \text{but}%
\ \ \ \ \ \ \ \ \ \ \left[  s_{1}\right]  \mathbf{m}_{3}=0.
\]
Also, for all $w\in S_{n}$, we have $\left[  w\right]  \nabla=1$ and $\left[
w\right]  \nabla^{-}=\left(  -1\right)  ^{w}$.

We shall furthermore use some elementary linear algebra, in particular the
notion of a \emph{trace}. The trace is defined both for a square matrix and
for an endomorphism of a free $\mathbf{k}$-module of finite rank. Let us
recall the definitions:

\begin{itemize}
\item The \emph{trace} of a square matrix $M$ is the sum of the diagonal
entries of $M$. It is denoted by $\operatorname*{Tr}M$. For example,
$\operatorname*{Tr}\left(
\begin{array}
[c]{cc}%
a & b\\
c & d
\end{array}
\right)  =a+d$.

The trace has several nice properties. In particular, any two conjugate
matrices have the same trace.

\item If $V$ is a free $\mathbf{k}$-module of finite rank (i.e., with a finite
basis), and if $f$ is any endomorphism of $V$ (that is, any $\mathbf{k}%
$-linear map from $V$ to $V$), then the \emph{trace} of $f$ is defined as the
trace of the matrix that represents $f$ with respect to any chosen basis of
$V$. (The choice of the basis does not affect the trace\footnote{Indeed, any
two matrices that each represent $f$ with respect to some basis of $V$ will
necessarily be conjugate to each other, and thus will have the same trace
(since any two conjugate matrices have the same trace).}.) Again, the trace of
$f$ is denoted by $\operatorname*{Tr}f$. Note that it belongs to $\mathbf{k}$.
\end{itemize}

Certain maps have particularly simple traces:

\begin{proposition}
\label{prop.groupalg.trace}Let $G$ be a finite group. Let $\mathbf{a}%
\in\mathbf{k}\left[  G\right]  $. Let $f_{\mathbf{a}}:\mathbf{k}\left[
G\right]  \rightarrow\mathbf{k}\left[  G\right]  $ be the $\mathbf{k}$-linear
map that sends each $\mathbf{b}\in\mathbf{k}\left[  G\right]  $ to
$\mathbf{ba}$. (For obvious reasons, this map $f_{\mathbf{a}}$ is called
\textquotedblleft right multiplication by $\mathbf{a}$\textquotedblright).
Then,
\[
\operatorname*{Tr}\left(  f_{\mathbf{a}}\right)  =\left\vert G\right\vert
\cdot\left[  1\right]  \mathbf{a}%
\]
(where $1$ stands for the neutral element of $G$).
\end{proposition}

\begin{proof}
First, we show a general formula for traces of $\mathbf{k}$-module
endomorphisms of $\mathbf{k}\left[  G\right]  $:

\begin{statement}
\textit{Claim 1:} Let $f$ be any $\mathbf{k}$-module endomorphism of
$\mathbf{k}\left[  G\right]  $. Then,%
\[
\operatorname*{Tr}f=\sum_{g\in G}\left[  g\right]  \left(  f\left(  g\right)
\right)  .
\]

\end{statement}

\begin{proof}
[Proof of Claim 1.]Let $g_{1},g_{2},\ldots,g_{m}$ be the elements of $G$,
listed without repetition (so that $m=\left\vert G\right\vert $). Then, the
$\mathbf{k}$-module $\mathbf{k}\left[  G\right]  $ is free, with basis
$\left(  g_{1},g_{2},\ldots,g_{m}\right)  $. The matrix that represents the
$\mathbf{k}$-linear map $f$ with respect to this basis is the $m\times
m$-matrix whose $\left(  i,j\right)  $-th entry is $\left[  g_{i}\right]
\left(  f\left(  g_{j}\right)  \right)  $ for all $i\in\left[  m\right]  $ and
$j\in\left[  m\right]  $. Hence, in particular, the diagonal entries of this
matrix are the elements $\left[  g_{i}\right]  \left(  f\left(  g_{i}\right)
\right)  $ for all $i\in\left[  m\right]  $. Therefore, the trace of this
matrix equals $\sum_{i=1}^{m}\left[  g_{i}\right]  \left(  f\left(
g_{i}\right)  \right)  $ (since the trace of a matrix is the sum of the
diagonal entries of this matrix). In other words,
\begin{equation}
\operatorname*{Tr}f=\sum_{i=1}^{m}\left[  g_{i}\right]  \left(  f\left(
g_{i}\right)  \right)  \label{pf.prop.groupalg.trace.0}%
\end{equation}
(since $\operatorname*{Tr}f$ is defined to be the trace of this matrix).

But $g_{1},g_{2},\ldots,g_{m}$ are the elements of $G$, listed without
repetition. Thus, $\sum_{g\in G}\left[  g\right]  \left(  f\left(  g\right)
\right)  =\sum_{i=1}^{m}\left[  g_{i}\right]  \left(  f\left(  g_{i}\right)
\right)  $. Comparing this with (\ref{pf.prop.groupalg.trace.0}), we obtain
$\operatorname*{Tr}f=\sum_{g\in G}\left[  g\right]  \left(  f\left(  g\right)
\right)  $. This proves Claim 1.
\end{proof}

Now, let us consider the $\mathbf{k}$-module endomorphism $f_{\mathbf{a}}$ in
particular. Write the element $\mathbf{a}\in\mathbf{k}\left[  G\right]  $ as
$\mathbf{a}=\sum\limits_{g\in G}\alpha_{g}g$ for some scalars $\alpha_{g}%
\in\mathbf{k}$. Then, $\left[  1\right]  \mathbf{a}=\alpha_{1}$ (by Definition
\ref{def.specht.ET.wa}). For each $h\in G$, we have%
\begin{align}
f_{\mathbf{a}}\left(  h\right)   &  =h\mathbf{a}\ \ \ \ \ \ \ \ \ \ \left(
\text{by the definition of }f_{\mathbf{a}}\right) \nonumber\\
&  =h\sum\limits_{g\in G}\alpha_{g}g\ \ \ \ \ \ \ \ \ \ \left(  \text{since
}\mathbf{a}=\sum\limits_{g\in G}\alpha_{g}g\right) \nonumber\\
&  =\sum\limits_{g\in G}\alpha_{g}hg=\sum\limits_{g\in G}\alpha_{h^{-1}%
g}\underbrace{hh^{-1}}_{=1}g\nonumber\\
&  \ \ \ \ \ \ \ \ \ \ \ \ \ \ \ \ \ \ \ \ \left(
\begin{array}
[c]{c}%
\text{here, we have substituted }h^{-1}g\text{ for }g\text{ in the sum,}\\
\text{since the map }G\rightarrow G,\ g\mapsto h^{-1}g\text{ is a bijection}\\
\text{(because }G\text{ is a group)}%
\end{array}
\right) \nonumber\\
&  =\sum\limits_{g\in G}\alpha_{h^{-1}g}g, \label{pf.prop.groupalg.trace.1}%
\end{align}
and thus%
\begin{equation}
\left[  g\right]  \left(  f_{\mathbf{a}}\left(  h\right)  \right)
=\alpha_{h^{-1}g}\ \ \ \ \ \ \ \ \ \ \text{for each }g\in G
\label{pf.prop.groupalg.trace.2}%
\end{equation}
(by Definition \ref{def.specht.ET.wa} again).

Now, Claim 1 (applied to $f=f_{\mathbf{a}}$) yields%
\[
\operatorname*{Tr}\left(  f_{\mathbf{a}}\right)  =\sum_{g\in G}%
\underbrace{\left[  g\right]  \left(  f_{\mathbf{a}}\left(  g\right)  \right)
}_{\substack{=\alpha_{g^{-1}g}\\\text{(by (\ref{pf.prop.groupalg.trace.2}),
applied to }h=g\text{)}}}=\sum_{g\in G}\underbrace{\alpha_{g^{-1}g}%
}_{\substack{=\alpha_{1}\\\text{(since }g^{-1}g=1\text{)}}}=\sum_{g\in
G}\alpha_{1}=\left\vert G\right\vert \cdot\alpha_{1}.
\]

We can rewrite this as $\operatorname*{Tr}\left(  f_{\mathbf{a}}\right)
=\left\vert G\right\vert \cdot\left[  1\right]  \mathbf{a}$ (since $\left[
1\right]  \mathbf{a}=\alpha_{1}$). This proves Proposition
\ref{prop.groupalg.trace}.
\end{proof}

We next need some more general properties of traces. The first one is a basic
and well-known property of traces of matrices (\cite[8.49]{Axler24}):

\begin{lemma}
\label{lem.linalg.TrAB}Let $k,d\in\mathbb{N}$. Let $A\in\mathbf{k}^{d\times
k}$ and $B\in\mathbf{k}^{k\times d}$ be two matrices. Then,
$\operatorname*{Tr}\left(  AB\right)  =\operatorname*{Tr}\left(  BA\right)  $.
\end{lemma}

Note that the two matrices $AB$ and $BA$ are both square, but of different sizes.

We can easily translate Lemma \ref{lem.linalg.TrAB} into the language of
linear maps:

\begin{lemma}
\label{lem.linalg.Trfg}Let $V$ and $W$ be two free $\mathbf{k}$-modules of
finite ranks (not necessarily of equal ranks!). Let $f:V\rightarrow W$ and
$g:W\rightarrow V$ be two $\mathbf{k}$-linear maps. Then, $\operatorname*{Tr}%
\left(  f\circ g\right)  =\operatorname*{Tr}\left(  g\circ f\right)  $.
\end{lemma}

\begin{proof}
Pick a basis $\left(  v_{1},v_{2},\ldots,v_{k}\right)  $ of $V$ and a basis
$\left(  w_{1},w_{2},\ldots,w_{d}\right)  $ of $W$. Let $A\in\mathbf{k}%
^{d\times k}$ be the matrix that represents the map $f:V\rightarrow W$ with
respect to these two bases, and let $B\in\mathbf{k}^{k\times d}$ be the matrix
that represents the map $g:W\rightarrow V$ with respect to these two bases.
Then, $AB$ is the matrix that represents the map $f\circ g:W\rightarrow W$
with respect to the basis $\left(  w_{1},w_{2},\ldots,w_{d}\right)  $ of $W$
(since matrix multiplication is \textquotedblleft the same
as\textquotedblright\ composition of linear maps). Thus, $\operatorname*{Tr}%
\left(  f\circ g\right)  =\operatorname*{Tr}\left(  AB\right)  $ (by the
definition of the trace of a linear map). Similarly, $\operatorname*{Tr}%
\left(  g\circ f\right)  =\operatorname*{Tr}\left(  BA\right)  $. However,
Lemma \ref{lem.linalg.TrAB} says that $\operatorname*{Tr}\left(  AB\right)
=\operatorname*{Tr}\left(  BA\right)  $. In other words, $\operatorname*{Tr}%
\left(  f\circ g\right)  =\operatorname*{Tr}\left(  g\circ f\right)  $ (since
$\operatorname*{Tr}\left(  f\circ g\right)  =\operatorname*{Tr}\left(
AB\right)  $ and $\operatorname*{Tr}\left(  g\circ f\right)
=\operatorname*{Tr}\left(  BA\right)  $). This proves Lemma
\ref{lem.linalg.Trfg}.
\end{proof}

Our next property of traces is even more obvious:

\begin{lemma}
\label{lem.linalg.Trkid}Let $d\in\mathbb{N}$. Let $W$ be a free $\mathbf{k}%
$-module of rank $d$. Let $\kappa\in\mathbf{k}$ be a scalar. Then, the
$\mathbf{k}$-module endomorphism $\kappa\operatorname*{id}\nolimits_{W}$ of
$W$ has trace $\operatorname*{Tr}\left(  \kappa\operatorname*{id}%
\nolimits_{W}\right)  =d\kappa$.
\end{lemma}

\begin{proof}
The $\mathbf{k}$-module $W$ is free of rank $d$. Thus, it has a basis $\left(
w_{1},w_{2},\ldots,w_{d}\right)  $ consisting of $d$ vectors. Consider this basis.

Now, consider the $\mathbf{k}$-module endomorphism $\operatorname*{id}%
\nolimits_{W}$ of $W$. The matrix that represents this endomorphism with
respect to the basis $\left(  w_{1},w_{2},\ldots,w_{d}\right)  $ is the
$d\times d$ identity matrix (since the matrix that represents an identity map
with respect to a basis is always the respective identity matrix). Thus, the
diagonal entries of this matrix are $\underbrace{1_{\mathbf{k}},1_{\mathbf{k}%
},\ldots,1_{\mathbf{k}}}_{d\text{ times}}$ (where $1_{\mathbf{k}}$ denotes the
unity of $\mathbf{k}$), and therefore its trace is $\underbrace{1_{\mathbf{k}%
}+1_{\mathbf{k}}+\cdots+1_{\mathbf{k}}}_{d\text{ times}}=d\cdot1_{\mathbf{k}}%
$. Hence, the trace of the $\mathbf{k}$-module endomorphism
$\operatorname*{id}\nolimits_{W}$ is $d\cdot1_{\mathbf{k}}$ as well (since the
trace of a $\mathbf{k}$-module endomorphism is defined as the trace of a
matrix that represents it with respect to a basis). In other words,
$\operatorname*{Tr}\left(  \operatorname*{id}\nolimits_{W}\right)
=d\cdot1_{\mathbf{k}}$.

Recall that the trace of a $\mathbf{k}$-module endomorphism depends
$\mathbf{k}$-linearly on this endomorphism. Hence, $\operatorname*{Tr}\left(
\omega j\right)  =\omega\operatorname*{Tr}j$ for any scalar $\omega
\in\mathbf{k}$ and any $\mathbf{k}$-module endomorphism $j$. Thus,
$\operatorname*{Tr}\left(  \kappa\operatorname*{id}\nolimits_{W}\right)
=\kappa\underbrace{\operatorname*{Tr}\left(  \operatorname*{id}\nolimits_{W}%
\right)  }_{=d\cdot1_{\mathbf{k}}}=\kappa\cdot d\cdot1_{\mathbf{k}}=d\kappa$.
This proves Lemma \ref{lem.linalg.Trkid}.
\end{proof}

The next lemma is technical but crucial:

\begin{lemma}
\label{lem.linalg.ladisch}Let $V$ be a free $\mathbf{k}$-module of finite
rank. Let $f:V\rightarrow V$ be an endomorphism of $V$. Assume that the image
$f\left(  V\right)  $ of $f$ is a free $\mathbf{k}$-module of rank $d$ for
some $d\in\mathbb{N}$. Assume furthermore that $f^{2}=\kappa f$ for some
scalar $\kappa\in\mathbf{k}$ (where, as usual, $f^{2}$ means $f\circ f$).
Then, $\operatorname*{Tr}f=d\kappa$.
\end{lemma}

Before we come to the proof of this lemma, let me say a few words about its
significance. We shall apply it to the situation where $V$ is $\mathcal{A}%
=\mathbf{k}\left[  S_{n}\right]  $ and where $f$ is right multiplication by
$\mathbf{E}_{T}$ for some $n$-tableau $T$ of shape $Y\left(  \lambda\right)
$. In this case, the image $f\left(  V\right)  $ is $\mathcal{A}\mathbf{E}%
_{T}$, which is isomorphic to the Specht module $\mathcal{S}^{Y\left(
\lambda\right)  }$ (by (\ref{eq.def.specht.ET.defs.SD=AET})) and therefore a
free $\mathbf{k}$-module of rank $f^{\lambda}$ (by Lemma
\ref{lem.specht.Slam-flam}). The equality $f^{2}=\kappa f$ will follow easily
from the equality $\mathbf{E}_{T}^{2}=\kappa\mathbf{E}_{T}$ that we have
obtained above as a consequence of Proposition \ref{prop.specht.ET.ETaET}.
Thus, Lemma \ref{lem.linalg.ladisch} will be applicable and allow us to
compute $\kappa$. Indeed, Lemma \ref{lem.linalg.ladisch} is tailor-made for
this application, but has also found uses in other representation-theoretical
questions (see the \textquotedblleft Motivation\textquotedblright\ in
\cite{Ladisc11}).

First, however, we must prove Lemma \ref{lem.linalg.ladisch}:

\begin{proof}
[Proof of Lemma \ref{lem.linalg.ladisch}.]The following beautiful proof is due
to Frieder Ladisch (\cite{Ladisc11}).

Clearly, $f\left(  v\right)  \in f\left(  V\right)  $ for each $v\in V$. Thus,
we can define a map $h:V\rightarrow f\left(  V\right)  $ that sends each $v\in
V$ to $f\left(  v\right)  $. Consider this map $h$. (This map $h$ is just the
map $f$, except that we have replaced its target by $f\left(  V\right)  $.)
Note that the map $h$ is $\mathbf{k}$-linear (since $f$ is $\mathbf{k}$-linear).

Let $g:f\left(  V\right)  \rightarrow V$ be the inclusion map of the
$\mathbf{k}$-submodule $f\left(  V\right)  $ into the $\mathbf{k}$-module $V$.
This is the map sending each $w\in f\left(  V\right)  $ to $w$ itself.
Clearly, this map $g$ is $\mathbf{k}$-linear.

From $f^{2}=\kappa f$, we easily obtain $h\circ g=\kappa\operatorname*{id}%
\nolimits_{f\left(  V\right)  }$.

[\textit{Proof:} Let $w\in f\left(  V\right)  $. Thus, $w=f\left(  v\right)  $
for some $v\in V$. Consider this $v$. But $g$ is the inclusion map, so that
$g\left(  w\right)  =w=f\left(  v\right)  $. Now,
\begin{align*}
\left(  h\circ g\right)  \left(  w\right)   &  =h\left(  \underbrace{g\left(
w\right)  }_{=f\left(  v\right)  }\right)  =h\left(  f\left(  v\right)
\right)  =f\left(  f\left(  v\right)  \right)  \ \ \ \ \ \ \ \ \ \ \left(
\text{by the definition of }h\right) \\
&  =\underbrace{f^{2}}_{=\kappa f}\left(  v\right)  =\left(  \kappa f\right)
\left(  v\right)  =\kappa\underbrace{f\left(  v\right)  }_{=w}=\kappa w.
\end{align*}

Forget that we fixed $w$. We thus have proved that $\left(  h\circ g\right)
\left(  w\right)  =\kappa w$ for each $w\in f\left(  V\right)  $. In other
words, $h\circ g=\kappa\operatorname*{id}\nolimits_{f\left(  V\right)  }$.]
\medskip

On the other hand, we have $g\circ h=f$.

[\textit{Proof:} For each $v\in V$, we have
\begin{align*}
\left(  g\circ h\right)  \left(  v\right)   &  =g\left(  h\left(  v\right)
\right)  =h\left(  v\right)  \ \ \ \ \ \ \ \ \ \ \left(  \text{since }g\text{
is the inclusion map}\right) \\
&  =f\left(  v\right)  \ \ \ \ \ \ \ \ \ \ \left(  \text{by the definition of
}h\right)  .
\end{align*}
In other words, $g\circ h=f$.] \medskip

We assumed that the $\mathbf{k}$-module $f\left(  V\right)  $ is free of rank
$d$. Thus, Lemma \ref{lem.linalg.Trkid} (applied to $W=f\left(  V\right)  $)
yields $\operatorname*{Tr}\left(  \kappa\operatorname*{id}\nolimits_{f\left(
V\right)  }\right)  =d\kappa$. Since $h\circ g=\kappa\operatorname*{id}%
\nolimits_{f\left(  V\right)  }$, we can rewrite this as%
\begin{equation}
\operatorname*{Tr}\left(  h\circ g\right)  =d\kappa.
\label{pf.lem.linalg.ladisch.Trhg}%
\end{equation}

But Lemma \ref{lem.linalg.Trfg} (applied to $f\left(  V\right)  $ and $h$
instead of $W$ and $f$) yields $\operatorname*{Tr}\left(  h\circ g\right)
=\operatorname*{Tr}\left(  g\circ h\right)  =\operatorname*{Tr}f$ (since
$g\circ h=f$). Comparing this with (\ref{pf.lem.linalg.ladisch.Trhg}), we
obtain $\operatorname*{Tr}f=d\kappa$. This proves Lemma
\ref{lem.linalg.ladisch}.
\end{proof}

\begin{noncompile}
Before I had written down Lemma \ref{lem.linalg.ladisch}, its role was played
by part \textbf{(c)} of the following simple linear-algebraic lemma:

\begin{lemma}
\label{lem.linalg.quasi-idp0}Assume that $\mathbf{k}$ is a field. Let $V$ be a
finite-dimensional $\mathbf{k}$-vector space. Let $f:V\rightarrow V$ be an
endomorphism of $V$. Then: \medskip

\textbf{(a)} The subspace $f\left(  V\right)  $ of $V$ is preserved under $f$,
so that the restriction $f\mid_{f\left(  V\right)  }$ is a $\mathbf{k}$-linear
map from $f\left(  V\right)  $ to $f\left(  V\right)  $. \medskip

\textbf{(b)} We have $\operatorname*{Tr}f=\operatorname*{Tr}\left(
f\mid_{f\left(  V\right)  }\right)  $. \medskip

\textbf{(c)} If $f^{2}=\kappa f$ for some scalar $\kappa\in\mathbf{k}$, then
$\operatorname*{Tr}\left(  f\right)  =\dim\left(  f\left(  V\right)  \right)
\cdot\kappa$.
\end{lemma}

\begin{proof}
[Proof sketch.]\textbf{(a)} We just need to show that $f\left(  v\right)  \in
f\left(  V\right)  $ for each $v\in f\left(  V\right)  $. But this is obvious.
\medskip

\textbf{(b)} Choose a basis $\left(  u_{1},u_{2},\ldots,u_{k}\right)  $ of the
vector space $f\left(  V\right)  $, and extend it to a basis $\left(
u_{1},u_{2},\ldots,u_{n}\right)  $ of the whole space $V$. (This can be done,
since $V$ is a finite-dimensional $\mathbf{k}$-vector space.) Now, let $A$ be
the $n\times n$-matrix that represents the linear map $f:V\rightarrow V$ with
respect to the basis $\left(  u_{1},u_{2},\ldots,u_{n}\right)  $ of $V$, and
let $B$ be the $n\times n$-matrix that represents the restriction $\left.
f\mid_{f\left(  V\right)  }\right.  :f\left(  V\right)  \rightarrow f\left(
V\right)  $ with respect to the basis $\left(  u_{1},u_{2},\ldots
,u_{k}\right)  $ of $f\left(  V\right)  $. Then, the definition of a trace
yields $\operatorname*{Tr}f=\operatorname*{Tr}A$ and $\operatorname*{Tr}%
\left(  f\mid_{f\left(  V\right)  }\right)  =\operatorname*{Tr}B$.

However, the columns of the matrix $A$ are (representations of) the vectors
$f\left(  u_{1}\right)  ,f\left(  u_{2}\right)  ,\ldots,f\left(  u_{n}\right)
$, which belong to $f\left(  V\right)  $ and thus can be written as
$\mathbf{k}$-linear combinations of the vectors $u_{1},u_{2},\ldots,u_{k}$
(without any need for $u_{k+1},u_{k+2},\ldots,u_{n}$). Hence, the last $n-k$
rows of the matrix $A$ are zero. Hence, the matrix $A$ has a block-matrix
structure%
\[
A=\left(
\begin{array}
[c]{cc}%
B & C\\
0_{\left(  n-k\right)  \times k} & 0_{\left(  n-k\right)  \times\left(
n-k\right)  }%
\end{array}
\right)
\]
for some $k\times\left(  n-k\right)  $-matrix $C$ (here, the top-left block is
$B$ because $B$ represents the restriction of $A$ to the span of $u_{1}%
,u_{2},\ldots,u_{k}$). Therefore, the diagonal entries of $A$ are just the
diagonal entries of $B$ followed by $n-k$ zeroes. In particular, the sum of
the diagonal entries of $A$ equals the sum of the diagonal entries of $B$. In
other words, $\operatorname*{Tr}A=\operatorname*{Tr}B$ (since the trace of a
matrix is the sum of its diagonal entries). In other words,
$\operatorname*{Tr}f=\operatorname*{Tr}\left(  f\mid_{f\left(  V\right)
}\right)  $ (since $\operatorname*{Tr}f=\operatorname*{Tr}A$ and
$\operatorname*{Tr}\left(  f\mid_{f\left(  V\right)  }\right)
=\operatorname*{Tr}B$). This proves Lemma \ref{lem.linalg.quasi-idp0}
\textbf{(b)}. \medskip

\textbf{(c)} Assume that $f^{2}=\kappa f$ for some scalar $\kappa\in
\mathbf{k}$. Consider this $\kappa$. Then, each vector $w\in f\left(
V\right)  $ satisfies $f\left(  w\right)  =\kappa w$ (because each vector
$w\in f\left(  V\right)  $ can be written as $f\left(  v\right)  $ for some
$v\in V$, and thus satisfies $f\left(  \underbrace{w}_{=f\left(  v\right)
}\right)  =f\left(  f\left(  v\right)  \right)  =\underbrace{f^{2}}_{=\kappa
f}\left(  v\right)  =\kappa\underbrace{f\left(  v\right)  }_{=w}=\kappa w$).
In other words, on the set $f\left(  V\right)  $, the map $f$ is just
multiplication by $\kappa$. In other words, $f\mid_{f\left(  V\right)
}=\kappa\operatorname*{id}\nolimits_{f\left(  V\right)  }$. Hence,%
\[
\operatorname*{Tr}\left(  f\mid_{f\left(  V\right)  }\right)
=\operatorname*{Tr}\left(  \kappa\operatorname*{id}\nolimits_{f\left(
V\right)  }\right)  =\kappa\operatorname*{Tr}\left(  \operatorname*{id}%
\nolimits_{f\left(  V\right)  }\right)
\]
(since the trace of a map depends $\mathbf{k}$-linearly on the map). But any
finite-dimensional vector space $W$ satisfies $\operatorname*{Tr}\left(
\operatorname*{id}\nolimits_{W}\right)  =\dim W$ in $\mathbf{k}$ (since the
identity endomorphism $\operatorname*{id}\nolimits_{W}$ is represented by the
identity matrix of size $\dim W$, and the trace of the latter matrix is $\dim
W$). Hence, $\operatorname*{Tr}\left(  \operatorname*{id}\nolimits_{f\left(
V\right)  }\right)  =\dim\left(  f\left(  V\right)  \right)  $.

Now, Lemma \ref{lem.linalg.quasi-idp0} \textbf{(b)} yields
\[
\operatorname*{Tr}\left(  f\right)  =\operatorname*{Tr}\left(  f\mid_{f\left(
V\right)  }\right)  =\kappa\underbrace{\operatorname*{Tr}\left(
\operatorname*{id}\nolimits_{f\left(  V\right)  }\right)  }_{=\dim\left(
f\left(  V\right)  \right)  }=\kappa\dim\left(  f\left(  V\right)  \right)
=\dim\left(  f\left(  V\right)  \right)  \cdot\kappa.
\]
This proves Lemma \ref{lem.linalg.quasi-idp0} \textbf{(c)}.

[\textbf{Note:} There is a better proof of Lemma \ref{lem.linalg.quasi-idp0}
\textbf{(c)} in
\url{https://mathoverflow.net/questions/58776/pseudo-idempotent-matrix-generating-a-free-module}
, which works for any commutative ring $\mathbf{k}$.]
\end{proof}
\end{noncompile}

Combining this lemma with Proposition \ref{prop.groupalg.trace}, we obtain the
following consequence for group algebras:

\begin{lemma}
\label{lem.linalg.quasi-idp1}Let $G$ be a finite group. Let $\mathcal{A}%
=\mathbf{k}\left[  G\right]  $. (In this, we break with Definition
\ref{def.specht.ET.defs} \textbf{(c)} in this lemma.)

Let $\mathbf{a}\in\mathcal{A}$ and $\kappa\in\mathbf{k}$ be two elements such
that $\mathbf{a}^{2}=\kappa\mathbf{a}$. Assume furthermore that the left ideal
$\mathcal{A}\mathbf{a}$ of $\mathcal{A}$ is a free $\mathbf{k}$-module of rank
$d$ for some $d\in\mathbb{N}$. Then, in $\mathbf{k}$, we have%
\[
d\kappa=\left\vert G\right\vert \cdot\left[  1\right]  \mathbf{a}%
\]
(where $1$ means the neutral element of $G$).
\end{lemma}

\begin{proof}
Let $f_{\mathbf{a}}:\mathcal{A}\rightarrow\mathcal{A}$ be the $\mathbf{k}%
$-linear map that sends each $\mathbf{b}\in\mathcal{A}$ to $\mathbf{ba}$.
Thus, $f_{\mathbf{a}}:\mathbf{k}\left[  G\right]  \rightarrow\mathbf{k}\left[
G\right]  $ is the $\mathbf{k}$-linear map that sends each $\mathbf{b}%
\in\mathbf{k}\left[  G\right]  $ to $\mathbf{ba}$ (since $\mathcal{A}%
=\mathbf{k}\left[  G\right]  $). Hence, Proposition \ref{prop.groupalg.trace}
yields
\begin{equation}
\operatorname*{Tr}\left(  f_{\mathbf{a}}\right)  =\left\vert G\right\vert
\cdot\left[  1\right]  \mathbf{a}. \label{pf.lem.linalg.quasi-idp1.1}%
\end{equation}

However, the $\mathbf{k}$-linear map $f_{\mathbf{a}}$ satisfies $f_{\mathbf{a}%
}^{2}=\kappa f_{\mathbf{a}}$ (since each $\mathbf{b}\in\mathcal{A}$ satisfies
\begin{align*}
f_{\mathbf{a}}^{2}\left(  \mathbf{b}\right)   &  =f_{\mathbf{a}}\left(
f_{\mathbf{a}}\left(  \mathbf{b}\right)  \right)  =\underbrace{f_{\mathbf{a}%
}\left(  \mathbf{b}\right)  }_{\substack{=\mathbf{ba}\\\text{(by the
definition}\\\text{of }f_{\mathbf{a}}\text{)}}}\cdot\,\mathbf{a}%
\ \ \ \ \ \ \ \ \ \ \left(  \text{by the definition of }f_{\mathbf{a}}\right)
\\
&  =\mathbf{b}\underbrace{\mathbf{a}\cdot\mathbf{a}}_{=\mathbf{a}^{2}%
=\kappa\mathbf{a}}=\mathbf{b}\kappa\mathbf{a}=\kappa\underbrace{\mathbf{ba}%
}_{=f_{\mathbf{a}}\left(  \mathbf{b}\right)  }=\kappa f_{\mathbf{a}}\left(
\mathbf{b}\right)
\end{align*}
). Moreover, its image is
\[
f_{\mathbf{a}}\left(  \mathcal{A}\right)  =\left\{  \underbrace{f_{\mathbf{a}%
}\left(  \mathbf{b}\right)  }_{=\mathbf{ba}}\ \mid\ \mathbf{b}\in
\mathcal{A}\right\}  =\left\{  \mathbf{ba}\ \mid\ \mathbf{b}\in\mathcal{A}%
\right\}  =\mathcal{A}\mathbf{a},
\]
and thus is a free $\mathbf{k}$-module of rank $d$ (since we assumed that
$\mathcal{A}\mathbf{a}$ is a free $\mathbf{k}$-module of rank $d$). Hence,
Lemma \ref{lem.linalg.ladisch} (applied to $V=\mathcal{A}$ and
$f=f_{\mathbf{a}}$) yields
\[
\operatorname*{Tr}\left(  f_{\mathbf{a}}\right)  =d\kappa.
\]
Comparing this with (\ref{pf.lem.linalg.quasi-idp1.1}), we find%
\[
d\kappa=\left\vert G\right\vert \cdot\left[  1\right]  \mathbf{a}.
\]
This proves Lemma \ref{lem.linalg.quasi-idp1}.
\end{proof}

\begin{lemma}
\label{lem.specht.ET.1-coord}Let $T$ be an $n$-tableau of any shape. Then,
$\left[  1\right]  \mathbf{E}_{T}=1$.
\end{lemma}

\begin{proof}
Recall that%
\begin{align}
\mathbf{E}_{T}  &  =\nabla_{\operatorname*{Col}T}^{-}\nabla
_{\operatorname*{Row}T}=\left(  \sum_{c\in\mathcal{C}\left(  T\right)
}\left(  -1\right)  ^{c}c\right)  \left(  \sum_{r\in\mathcal{R}\left(
T\right)  }r\right) \nonumber\\
&  \ \ \ \ \ \ \ \ \ \ \ \ \ \ \ \ \ \ \ \ \left(
\begin{array}
[c]{c}%
\text{since Definition \ref{def.symmetrizers.symmetrizers} yields }%
\nabla_{\operatorname*{Row}T}=\sum_{w\in\mathcal{R}\left(  T\right)  }%
w=\sum_{r\in\mathcal{R}\left(  T\right)  }r\\
\text{and }\nabla_{\operatorname*{Col}T}^{-}=\sum_{w\in\mathcal{C}\left(
T\right)  }\left(  -1\right)  ^{w}w=\sum_{c\in\mathcal{C}\left(  T\right)
}\left(  -1\right)  ^{c}c
\end{array}
\right) \nonumber\\
&  =\sum_{c\in\mathcal{C}\left(  T\right)  }\ \ \sum_{r\in\mathcal{R}\left(
T\right)  }\left(  -1\right)  ^{c}cr=\sum_{\left(  c,r\right)  \in
\mathcal{C}\left(  T\right)  \times\mathcal{R}\left(  T\right)  }\left(
-1\right)  ^{c}cr. \label{pf.lem.specht.ET.1-coord.1}%
\end{align}
Taking the coefficient of $\operatorname*{id}$ on both sides of this equality,
we find%
\[
\left[  \operatorname*{id}\right]  \mathbf{E}_{T}=\left[  \operatorname*{id}%
\right]  \left(  \sum_{\left(  c,r\right)  \in\mathcal{C}\left(  T\right)
\times\mathcal{R}\left(  T\right)  }\left(  -1\right)  ^{c}cr\right)
=\sum_{\substack{\left(  c,r\right)  \in\mathcal{C}\left(  T\right)
\times\mathcal{R}\left(  T\right)  ;\\cr=\operatorname*{id}}}\left(
-1\right)  ^{c}=\left(  -1\right)  ^{\operatorname*{id}},
\]
because Proposition \ref{prop.tableau.RnC} \textbf{(b)} shows that the only
pair $\left(  c,r\right)  \in\mathcal{C}\left(  T\right)  \times
\mathcal{R}\left(  T\right)  $ satisfying $cr=\operatorname*{id}$ is the pair
$\left(  \operatorname*{id},\operatorname*{id}\right)  $. Thus, $\left[
\operatorname*{id}\right]  \mathbf{E}_{T}=\left(  -1\right)
^{\operatorname*{id}}=1$. Since we denote $\operatorname*{id}$ by $1$, we can
rewrite this as $\left[  1\right]  \mathbf{E}_{T}=1$. This proves Lemma
\ref{lem.specht.ET.1-coord}.
\end{proof}

\begin{proof}
[Proof of Theorem \ref{thm.specht.ETidp}.]Let us first prove the theorem for
$\mathbf{k}=\mathbb{Z}$.

So let us assume that $\mathbf{k}=\mathbb{Z}$. Proposition
\ref{prop.specht.ET.ETaET} (applied to $\mathbf{a}=1$) yields $\mathbf{E}%
_{T}1\mathbf{E}_{T}=\kappa\mathbf{E}_{T}$ for some scalar $\kappa\in
\mathbf{k}$. Consider this $\kappa$. Thus, $\kappa\in\mathbf{k}=\mathbb{Z}$.
Of course, $\mathbf{E}_{T}1\mathbf{E}_{T}=\mathbf{E}_{T}\mathbf{E}%
_{T}=\mathbf{E}_{T}^{2}$, so that%
\begin{equation}
\mathbf{E}_{T}^{2}=\mathbf{E}_{T}1\mathbf{E}_{T}=\kappa\mathbf{E}_{T}.
\label{pf.thm.specht.ETidp.1}%
\end{equation}

It remains to show that $\kappa=\dfrac{n!}{f^{\lambda}}$.

Let $\mathcal{A}=\mathbf{k}\left[  S_{n}\right]  $. Lemma
\ref{lem.specht.ET.1-coord} yields $\left[  1\right]  \mathbf{E}_{T}=1\neq0$.
Hence, the element $\mathbf{E}_{T}$ is nonzero. Now, the $\mathbf{k}$-module
$\mathcal{A}\mathbf{E}_{T}$ contains the element $1\mathbf{E}_{T}%
=\mathbf{E}_{T}$, which is nonzero (as we have just shown). Hence,
$\mathcal{A}\mathbf{E}_{T}$ itself is nonzero.

\begin{noncompile}
Let $D=Y\left(  \lambda\right)  $. Then, $T$ is an $n$-tableau of shape $D$
(since $T$ is an $n$-tableau of shape $Y\left(  \lambda\right)  $). Hence,
(\ref{eq.def.specht.ET.defs.SD=AET}) yields $\mathcal{S}^{D}\cong%
\mathcal{A}\mathbf{E}_{T}$. In other words, $\mathcal{S}^{Y\left(
\lambda\right)  }\cong\mathcal{A}\mathbf{E}_{T}$ (since $D=Y\left(
\lambda\right)  $).
\end{noncompile}

But $T$ is an $n$-tableau of shape $Y\left(  \lambda\right)  $. Hence,
(\ref{eq.def.specht.ET.defs.SD=AET}) (applied to $D=Y\left(  \lambda\right)
$) yields $\mathcal{S}^{Y\left(  \lambda\right)  }\cong\mathcal{A}%
\mathbf{E}_{T}$.

\begin{noncompile}
Furthermore, $f^{\lambda}$ is the \# of standard tableaux of shape $Y\left(
\lambda\right)  $. In other words, $f^{\lambda}$ is the \# of standard
tableaux of shape $D$ (since $D=Y\left(  \lambda\right)  $).

Moreover, $D$ is a skew Young diagram (since $D=Y\left(  \lambda\right)
=Y\left(  \lambda/\varnothing\right)  $). Thus, the standard basis theorem
(Theorem \ref{thm.spechtmod.basis}) shows that the standard polytabloids
$\mathbf{e}_{T}$ (that is, the $\mathbf{e}_{T}$ where $T$ ranges over all
standard tableaux of shape $D$) form a basis of the $\mathbf{k}$-module
$\mathcal{S}^{D}$. Hence, $\mathcal{S}^{D}$ is a free $\mathbf{k}$-module, and
its rank is equal to the \# of all standard tableaux of shape $D$. Thus, its
rank is $f^{\lambda}$ (since $f^{\lambda}$ is the \# of standard tableaux of
shape $D$).

So we have shown that $\mathcal{S}^{D}$ is a free $\mathbf{k}$-module of rank
$f^{\lambda}$. Therefore, $\mathcal{A}\mathbf{E}_{T}$ is a free $\mathbf{k}%
$-module of rank $f^{\lambda}$ (since $\mathcal{S}^{D}\cong\mathcal{A}%
\mathbf{E}_{T}$). Since the $\mathbf{k}$-module $\mathcal{A}\mathbf{E}_{T}$ is
nonzero (as we showed above), this entails that its rank $f^{\lambda}$ is nonzero.
\end{noncompile}

Let $\mathcal{S}^{\lambda}$ be the Specht module $\mathcal{S}^{Y\left(
\lambda\right)  }$. Then, $\mathcal{S}^{\lambda}=\mathcal{S}^{Y\left(
\lambda\right)  }\cong\mathcal{A}\mathbf{E}_{T}$. But Lemma
\ref{lem.specht.Slam-flam} shows that $\mathcal{S}^{\lambda}$ is a free
$\mathbf{k}$-module of rank $f^{\lambda}$. Therefore, $\mathcal{A}%
\mathbf{E}_{T}$ is a free $\mathbf{k}$-module of rank $f^{\lambda}$ (since
$\mathcal{S}^{\lambda}\cong\mathcal{A}\mathbf{E}_{T}$). Since the $\mathbf{k}%
$-module $\mathcal{A}\mathbf{E}_{T}$ is nonzero (as we showed above), this
entails that its rank $f^{\lambda}$ is nonzero.

Now, recall that $\mathbf{E}_{T}^{2}=\kappa\mathbf{E}_{T}$. Hence, Lemma
\ref{lem.linalg.quasi-idp1} (applied to $G=S_{n}$ and $\mathbf{a}%
=\mathbf{E}_{T}$ and $d=f^{\lambda}$) yields%
\[
f^{\lambda}\kappa=\underbrace{\left\vert S_{n}\right\vert }_{=n!}%
\cdot\underbrace{\left[  1\right]  \mathbf{E}_{T}}_{=1}=n!.
\]
Hence,
\[
\kappa=\dfrac{n!}{f^{\lambda}}\ \ \ \ \ \ \ \ \ \ \left(  \text{since
}f^{\lambda}\text{ is nonzero}\right)  .
\]

Thus, $\dfrac{n!}{f^{\lambda}}\in\mathbb{Z}$ (since $\kappa\in\mathbb{Z}$).
Hence, $\dfrac{n!}{f^{\lambda}}$ is a positive integer (since its positivity
is obvious). Moreover, in view of $\kappa=\dfrac{n!}{f^{\lambda}}$, we can
rewrite (\ref{pf.thm.specht.ETidp.1}) as%
\begin{equation}
\mathbf{E}_{T}^{2}=\dfrac{n!}{f^{\lambda}}\mathbf{E}_{T}.
\label{pf.thm.specht.ETidp.5}%
\end{equation}

Thus, we have proved Theorem \ref{thm.specht.ETidp} in the case when
$\mathbf{k}=\mathbb{Z}$. It remains to prove the theorem in the general case,
i.e., for an arbitrary commutative ring $\mathbf{k}$.

But this is just an easy application of base change: Drop our assumption that
$\mathbf{k}=\mathbb{Z}$, and consider an arbitrary commutative ring
$\mathbf{k}$ instead. Then, there is a canonical ring morphism $f:\mathbb{Z}%
\rightarrow\mathbf{k}$ (which sends each integer $m$ to $m\cdot1_{\mathbf{k}}%
$). This ring morphism $f$ induces a base change morphism $f_{\ast}%
:\mathbb{Z}\left[  S_{n}\right]  \rightarrow\mathbf{k}\left[  S_{n}\right]  $
(as defined in Definition \ref{def.monalg.base-change}), and this latter
morphism $f_{\ast}$ sends the Young symmetrizer $\mathbf{E}_{T}\in
\mathbb{Z}\left[  S_{n}\right]  $ to the corresponding Young symmetrizer
$\mathbf{E}_{T}\in\mathbf{k}\left[  S_{n}\right]  $ (since the definition of
$\mathbf{E}_{T}$ is independent of the base ring\footnote{To make this more
concrete: As we saw in the proof of Lemma \ref{lem.specht.ET.1-coord}, we
have
\[
\mathbf{E}_{T}=\left(  \sum_{c\in\mathcal{C}\left(  T\right)  }\left(
-1\right)  ^{c}c\right)  \left(  \sum_{r\in\mathcal{R}\left(  T\right)
}r\right)  .
\]
This equality holds both in $\mathbb{Z}\left[  S_{n}\right]  $ and in
$\mathbf{k}\left[  S_{n}\right]  $. But the base change morphism $f_{\ast
}:\mathbb{Z}\left[  S_{n}\right]  \rightarrow\mathbf{k}\left[  S_{n}\right]  $
sends the element $\left(  \sum_{c\in\mathcal{C}\left(  T\right)  }\left(
-1\right)  ^{c}c\right)  \left(  \sum_{r\in\mathcal{R}\left(  T\right)
}r\right)  $ of $\mathbb{Z}\left[  S_{n}\right]  $ to the analogous element of
$\mathbf{k}\left[  S_{n}\right]  $. Thus, it sends the Young symmetrizer
$\mathbf{E}_{T}\in\mathbb{Z}\left[  S_{n}\right]  $ to the corresponding Young
symmetrizer $\mathbf{E}_{T}\in\mathbf{k}\left[  S_{n}\right]  $.}). Hence,
applying this morphism $f_{\ast}$ to both sides of the equality
(\ref{pf.thm.specht.ETidp.5}), we obtain the analogous equality%
\[
\mathbf{E}_{T}^{2}=\dfrac{n!}{f^{\lambda}}\mathbf{E}_{T}%
\ \ \ \ \ \ \ \ \ \ \text{in }\mathbf{k}\left[  S_{n}\right]
\]
(since $f_{\ast}$ is a ring morphism and thus respects squares and scaling).
This completes the proof of Theorem \ref{thm.specht.ETidp} in the general case
(since we have already shown that $\dfrac{n!}{f^{\lambda}}$ is a positive integer).
\end{proof}

We observe one more easy property of Young symmetrizers:

\begin{lemma}
\label{lem.specht.ET.wT}Let $T$ be an $n$-tableau. Let $w\in S_{n}$. Then,%
\[
\mathbf{E}_{wT}=w\mathbf{E}_{T}w^{-1}.
\]

\end{lemma}

\begin{proof}
We have $wT=w\rightharpoonup T$. Thus,
\begin{align*}
\mathbf{E}_{wT}  &  =\mathbf{E}_{w\rightharpoonup T}\\
&  =\underbrace{\nabla_{\operatorname*{Col}\left(  w\rightharpoonup T\right)
}^{-}}_{\substack{=w\nabla_{\operatorname*{Col}T}^{-}w^{-1}\\\text{(by
Proposition \ref{prop.symmetrizers.conj})}}}\ \ \underbrace{\nabla
_{\operatorname*{Row}\left(  w\rightharpoonup T\right)  }}_{\substack{=w\nabla
_{\operatorname*{Row}T}w^{-1}\\\text{(by Proposition
\ref{prop.symmetrizers.conj})}}}\ \ \ \ \ \ \ \ \ \ \left(  \text{by the
definition of }\mathbf{E}_{w\rightharpoonup T}\right) \\
&  =w\nabla_{\operatorname*{Col}T}^{-}\underbrace{w^{-1}w}_{=1}\nabla
_{\operatorname*{Row}T}w^{-1}=w\underbrace{\nabla_{\operatorname*{Col}T}%
^{-}\nabla_{\operatorname*{Row}T}}_{\substack{=\mathbf{E}_{T}\\\text{(by the
definition of }\mathbf{E}_{T}\text{)}}}w^{-1}=w\mathbf{E}_{T}w^{-1}.
\end{align*}
Thus, Lemma \ref{lem.specht.ET.wT} is proved.
\end{proof}

Lemma \ref{lem.specht.ET.wT} shows that if you know the Young symmetrizer
$\mathbf{E}_{T}$ for an $n$-tableau $T$ of a given shape, then you can easily
obtain the $\mathbf{E}_{S}$'s for all other $n$-tableaux $S$ of the same shape.

\subsubsection{Orthogonality}

Now let us see what happens when we multiply two Young symmetrizers
corresponding to different tableaux. First, we deal with the easy case where
the tableaux have different shapes:

\begin{proposition}
[double zero-sandwich lemma]\label{prop.specht.ET.ESaET}Let $\lambda$ and
$\mu$ be two distinct partitions of $n$. Let $S$ be an $n$-tableau of shape
$Y\left(  \lambda\right)  $. Let $T$ be an $n$-tableau of shape $Y\left(
\mu\right)  $. Let $\mathbf{a}\in\mathbf{k}\left[  S_{n}\right]  $. Then,%
\[
\mathbf{E}_{S}\mathbf{aE}_{T}=0.
\]

\end{proposition}

\begin{proof}
By the definition of a Young symmetrizer, we have $\mathbf{E}_{S}%
=\nabla_{\operatorname*{Col}S}^{-}\nabla_{\operatorname*{Row}S}$ and
$\mathbf{E}_{T}=\nabla_{\operatorname*{Col}T}^{-}\nabla_{\operatorname*{Row}%
T}$. Thus,
\begin{equation}
\mathbf{E}_{S}\mathbf{aE}_{T}=\nabla_{\operatorname*{Col}S}^{-}\nabla
_{\operatorname*{Row}S}\mathbf{a}\nabla_{\operatorname*{Col}T}^{-}%
\nabla_{\operatorname*{Row}T}. \label{pf.prop.specht.ET.ESaET.0}%
\end{equation}
But we have $\lambda\neq\mu$. Hence, there exists some $i\geq1$ such that
$\lambda_{i}\neq\mu_{i}$. Consider the \textbf{smallest} such $i$. Then,
$\lambda_{k}=\mu_{k}$ for all $k\in\left[  i-1\right]  $ (since $i$ is the
smallest). Summing up all these $i-1$ equalities, we obtain%
\begin{equation}
\lambda_{1}+\lambda_{2}+\cdots+\lambda_{i-1}=\mu_{1}+\mu_{2}+\cdots+\mu_{i-1}.
\label{pf.prop.specht.ET.ESaET.1}%
\end{equation}

However, $\lambda_{i}\neq\mu_{i}$. Thus, we are in one of the following two cases:

\textit{Case 1:} We have $\lambda_{i}>\mu_{i}$.

\textit{Case 2:} We have $\lambda_{i}<\mu_{i}$.

Let us consider Case 1 first. In this case, we have $\lambda_{i}>\mu_{i}$.
Adding this inequality to the equality (\ref{pf.prop.specht.ET.ESaET.1}), we
obtain%
\[
\left(  \lambda_{1}+\lambda_{2}+\cdots+\lambda_{i-1}\right)  +\lambda
_{i}>\left(  \mu_{1}+\mu_{2}+\cdots+\mu_{i-1}\right)  +\mu_{i}.
\]
In other words,%
\[
\lambda_{1}+\lambda_{2}+\cdots+\lambda_{i}>\mu_{1}+\mu_{2}+\cdots+\mu_{i}.
\]
Hence, Theorem \ref{thm.specht.sandw0} \textbf{(d)} yields $\nabla
_{\operatorname*{Row}S}\mathbf{a}\nabla_{\operatorname*{Col}T}^{-}=0$. And
thus, (\ref{pf.prop.specht.ET.ESaET.0}) becomes%
\[
\mathbf{E}_{S}\mathbf{aE}_{T}=\nabla_{\operatorname*{Col}S}^{-}%
\underbrace{\nabla_{\operatorname*{Row}S}\mathbf{a}\nabla_{\operatorname*{Col}%
T}^{-}}_{=0}\nabla_{\operatorname*{Row}T}=0.
\]
This proves Proposition \ref{prop.specht.ET.ESaET} in Case 1.

Let us now consider Case 2. In this case, we have $\lambda_{i}<\mu_{i}$.
Adding this inequality to the equality (\ref{pf.prop.specht.ET.ESaET.1}), we
obtain%
\[
\left(  \lambda_{1}+\lambda_{2}+\cdots+\lambda_{i-1}\right)  +\lambda
_{i}<\left(  \mu_{1}+\mu_{2}+\cdots+\mu_{i-1}\right)  +\mu_{i}.
\]
In other words,%
\[
\lambda_{1}+\lambda_{2}+\cdots+\lambda_{i}<\mu_{1}+\mu_{2}+\cdots+\mu_{i}.
\]
In other words,%
\[
\mu_{1}+\mu_{2}+\cdots+\mu_{i}>\lambda_{1}+\lambda_{2}+\cdots+\lambda_{i}.
\]
Thus, Theorem \ref{thm.specht.sandw0} \textbf{(b)} (applied to $\mu$,
$\lambda$, $T$, $S$ and $\nabla_{\operatorname*{Row}S}\mathbf{a}%
\nabla_{\operatorname*{Col}T}^{-}$ instead of $\lambda$, $\mu$, $S$, $T$ and
$\mathbf{a}$) yields
\[
\nabla_{\operatorname*{Col}S}^{-}\nabla_{\operatorname*{Row}S}\mathbf{a}%
\nabla_{\operatorname*{Col}T}^{-}\nabla_{\operatorname*{Row}T}=0.
\]
In view of (\ref{pf.prop.specht.ET.ESaET.0}), we can rewrite this as
$\mathbf{E}_{S}\mathbf{aE}_{T}=0$. This proves Proposition
\ref{prop.specht.ET.ESaET} in Case 2.

Hence, Proposition \ref{prop.specht.ET.ESaET} is proved in both Cases 1 and 2.
\end{proof}

What happens when you multiply two Young symmetrizers $\mathbf{E}_{S}$ and
$\mathbf{E}_{T}$ of the same shape (but possibly different tableaux)? The
following proposition gives an answer:

\begin{proposition}
\label{prop.specht.ET.ESET}Let $\lambda$ be a partition of $n$. Let $S$ and
$T$ be two $n$-tableaux of shape $Y\left(  \lambda\right)  $. Let $f^{\lambda
}$ be the \# of standard tableaux of shape $Y\left(  \lambda\right)  $. Then:
\medskip

\textbf{(a)} If there exist no $r\in\mathcal{R}\left(  S\right)  $ and
$c\in\mathcal{C}\left(  T\right)  $ such that $rS=cT$, then $\mathbf{E}%
_{S}\mathbf{E}_{T}=0$. \medskip

\textbf{(b)} If there exist $r\in\mathcal{R}\left(  S\right)  $ and
$c\in\mathcal{C}\left(  T\right)  $ such that $rS=cT$, then these $r$ and $c$
satisfy $\mathbf{E}_{S}\mathbf{E}_{T}=h^{\lambda}\left(  -1\right)  ^{c}%
r^{-1}\mathbf{E}_{U}c$, where we set $U:=rS=cT$ and $h^{\lambda}:=\dfrac
{n!}{f^{\lambda}}$.
\end{proposition}

\begin{proof}
\textbf{(a)} Assume that there exist no $r\in\mathcal{R}\left(  S\right)  $
and $c\in\mathcal{C}\left(  T\right)  $ such that $rS=cT$.

If there were no two distinct integers that are horizontally aligned in $S$
and simultaneously vertically aligned in $T$, then Theorem
\ref{thm.youngtab.alt} \textbf{(c)} (applied to $\mu=\lambda$) would show that
there exist permutations $r\in\mathcal{R}\left(  S\right)  $ and
$c\in\mathcal{C}\left(  T\right)  $ such that $rS=cT$; but this would
contradict the preceding sentence. Thus, there are two distinct integers that
are horizontally aligned in $S$ and simultaneously vertically aligned in $T$.
Hence, Lemma \ref{lem.specht.sandw0-ur} yields that $\nabla
_{\operatorname*{Row}S}\nabla_{\operatorname*{Col}T}^{-}=0$. But%
\[
\underbrace{\mathbf{E}_{S}}_{\substack{=\nabla_{\operatorname*{Col}S}%
^{-}\nabla_{\operatorname*{Row}S}\\\text{(by the definition}\\\text{of
}\mathbf{E}_{S}\text{)}}}\ \ \underbrace{\mathbf{E}_{T}}_{\substack{=\nabla
_{\operatorname*{Col}T}^{-}\nabla_{\operatorname*{Row}T}\\\text{(by the
definition}\\\text{of }\mathbf{E}_{T}\text{)}}}=\nabla_{\operatorname*{Col}%
S}^{-}\underbrace{\nabla_{\operatorname*{Row}S}\nabla_{\operatorname*{Col}%
T}^{-}}_{=0}\nabla_{\operatorname*{Row}T}=0.
\]
This proves Proposition \ref{prop.specht.ET.ESET} \textbf{(a)}. \medskip

\textbf{(b)} Let $r\in\mathcal{R}\left(  S\right)  $ and $c\in\mathcal{C}%
\left(  T\right)  $ be such that $rS=cT$. Set $U:=rS=cT$ and $h^{\lambda
}:=\dfrac{n!}{f^{\lambda}}$. We must prove that $\mathbf{E}_{S}\mathbf{E}%
_{T}=h^{\lambda}\left(  -1\right)  ^{c}r^{-1}\mathbf{E}_{U}c$.

Theorem \ref{thm.specht.ETidp} (applied to $U$ instead of $T$) yields
\begin{equation}
\mathbf{E}_{U}^{2}=\underbrace{\dfrac{n!}{f^{\lambda}}}_{=h^{\lambda}%
}\mathbf{E}_{U}=h^{\lambda}\mathbf{E}_{U}.
\label{pf.prop.specht.ET.ESET.b.EUEU}%
\end{equation}

Proposition \ref{prop.symmetrizers.fix} \textbf{(a)} (applied to $S$ and $r$
instead of $T$ and $w$) yields $r\nabla_{\operatorname*{Row}S}=\nabla
_{\operatorname*{Row}S}r=\nabla_{\operatorname*{Row}S}$ (since $r\in
\mathcal{R}\left(  S\right)  $). Hence,
\[
\underbrace{\mathbf{E}_{S}}_{\substack{=\nabla_{\operatorname*{Col}S}%
^{-}\nabla_{\operatorname*{Row}S}\\\text{(by the definition}\\\text{of
}\mathbf{E}_{S}\text{)}}}r=\nabla_{\operatorname*{Col}S}^{-}\underbrace{\nabla
_{\operatorname*{Row}S}r}_{=\nabla_{\operatorname*{Row}S}}=\nabla
_{\operatorname*{Col}S}^{-}\nabla_{\operatorname*{Row}S}=\mathbf{E}_{S}.
\]

Proposition \ref{prop.symmetrizers.fix} \textbf{(b)} (applied to $w=c$) yields
$c\nabla_{\operatorname*{Col}T}^{-}=\nabla_{\operatorname*{Col}T}^{-}c=\left(
-1\right)  ^{c}\nabla_{\operatorname*{Col}T}^{-}$ (since $c\in\mathcal{C}%
\left(  T\right)  $). Hence,
\[
c\underbrace{\mathbf{E}_{T}}_{\substack{=\nabla_{\operatorname*{Col}T}%
^{-}\nabla_{\operatorname*{Row}T}\\\text{(by the definition}\\\text{of
}\mathbf{E}_{T}\text{)}}}=\underbrace{c\nabla_{\operatorname*{Col}T}^{-}%
}_{=\left(  -1\right)  ^{c}\nabla_{\operatorname*{Col}T}^{-}}\nabla
_{\operatorname*{Row}T}=\left(  -1\right)  ^{c}\underbrace{\nabla
_{\operatorname*{Col}T}^{-}\nabla_{\operatorname*{Row}T}}_{=\mathbf{E}_{T}%
}=\left(  -1\right)  ^{c}\mathbf{E}_{T}.
\]

But $U=rS$ and thus
\begin{align*}
\mathbf{E}_{U}  &  =\mathbf{E}_{rS}=r\underbrace{\mathbf{E}_{S}r^{-1}%
}_{\substack{=\mathbf{E}_{S}\\\text{(since }\mathbf{E}_{S}r=\mathbf{E}%
_{S}\text{)}}}\ \ \ \ \ \ \ \ \ \ \left(  \text{by Lemma
\ref{lem.specht.ET.wT}}\right) \\
&  =r\mathbf{E}_{S},
\end{align*}
so that $\mathbf{E}_{S}=r^{-1}\mathbf{E}_{U}$. Furthermore, $U=cT$ and thus%
\begin{align*}
\mathbf{E}_{U}  &  =\mathbf{E}_{cT}=\underbrace{c\mathbf{E}_{T}}_{=\left(
-1\right)  ^{c}\mathbf{E}_{T}}c^{-1}\ \ \ \ \ \ \ \ \ \ \left(  \text{by Lemma
\ref{lem.specht.ET.wT}}\right) \\
&  =\left(  -1\right)  ^{c}\mathbf{E}_{T}c^{-1},
\end{align*}
so that $\mathbf{E}_{T}=\underbrace{\dfrac{1}{\left(  -1\right)  ^{c}}%
}_{=\left(  -1\right)  ^{c}}\mathbf{E}_{U}c=\left(  -1\right)  ^{c}%
\mathbf{E}_{U}c$. Thus,%
\[
\underbrace{\mathbf{E}_{S}}_{=r^{-1}\mathbf{E}_{U}}\ \ \underbrace{\mathbf{E}%
_{T}}_{=\left(  -1\right)  ^{c}\mathbf{E}_{U}c}=r^{-1}\mathbf{E}_{U}\left(
-1\right)  ^{c}\mathbf{E}_{U}c=\left(  -1\right)  ^{c}r^{-1}%
\underbrace{\mathbf{E}_{U}^{2}}_{\substack{=h^{\lambda}\mathbf{E}%
_{U}\\\text{(by (\ref{pf.prop.specht.ET.ESET.b.EUEU}))}}}c=h^{\lambda}\left(
-1\right)  ^{c}r^{-1}\mathbf{E}_{U}c.
\]
This proves Proposition \ref{prop.specht.ET.ESET} \textbf{(b)}.
\end{proof}

A formula for $\mathbf{E}_{S}\mathbf{aE}_{T}$ (for arbitrary $\mathbf{a}%
\in\mathbf{k}\left[  S_{n}\right]  $) can also be found, but it is less
pleasant and no stronger than Proposition \ref{prop.specht.ET.ESET} (as it
follows from the latter proposition using $\mathbf{E}_{wT}=w\mathbf{E}%
_{T}w^{-1}$).

\begin{example}
\label{exa.specht.ET.ESET.1}\textbf{(a)} Let us take $n=4$ and $\lambda
=\left(  2,2\right)  $ and $S=\ytableaushort{13,24}$ and
$T=\ytableaushort{12,34}$\ \ . Then, there exist no $r\in\mathcal{R}\left(
S\right)  $ and $c\in\mathcal{C}\left(  T\right)  $ such that $rS=cT$. Hence,
Proposition \ref{prop.specht.ET.ESET} \textbf{(a)} yields $\mathbf{E}%
_{S}\mathbf{E}_{T}=0$. A similar argument shows that $\mathbf{E}_{T}%
\mathbf{E}_{S}=0$. \medskip

\textbf{(b)} Let us instead take $n=5$ and $\lambda=\left(  3,2\right)  $ and
$S=\ytableaushort{135,24}$ and $T=\ytableaushort{123,45}$\ \ . Then, there
exist $r\in\mathcal{R}\left(  S\right)  $ and $c\in\mathcal{C}\left(
T\right)  $ such that $rS=cT$ (namely, $r=t_{3,5}t_{2,4}$ and $c=t_{2,5}$).
Hence, Proposition \ref{prop.specht.ET.ESET} \textbf{(b)} yields
$\mathbf{E}_{S}\mathbf{E}_{T}=h^{\lambda}\left(  -1\right)  ^{c}%
r^{-1}\mathbf{E}_{U}c$, where $U:=rS=cT=\ytableaushort{153,42}$ and
$h^{\lambda}:=\dfrac{n!}{f^{\lambda}}=\dfrac{5!}{5}=24$.
\end{example}

\begin{exercise}
\label{exe.specht.ET.ESET=0-when}\textbf{(a)} \fbox{4} Let $\lambda$ be a hook
partition (see Definition \ref{def.partitions.hook}). Let $S$ and $T$ be two
distinct standard $n$-tableaux of shape $Y\left(  \lambda\right)  $. Prove
that $\mathbf{E}_{S}\mathbf{E}_{T}=0$. \medskip

\textbf{(b)} \fbox{2} One might wonder what other partitions $\lambda$ have
the property stated in part \textbf{(a)}. For example, $\lambda=\left(
2,2\right)  $ does (as we saw in Example \ref{exa.specht.ET.ESET.1}
\textbf{(a)}), but $\lambda=\left(  3,2\right)  $ does not (at least not when
$\mathbf{k} = \mathbb{Q}$, as we saw in Example \ref{exa.specht.ET.ESET.1}
\textbf{(b)}).

Show that if $\mathbf{k}$ is a field of characteristic $0$, then the only
partitions $\lambda$ that have the property stated in part \textbf{(a)} are
the hook partitions and the partitions $\left(  {}\right)  $ and $\left(
2,2\right)  $.
\end{exercise}

Proposition \ref{prop.specht.ET.ESET} is nice (note how $S$ and $T$ appear in
symmetric roles in it), but we will eventually need a less symmetric (and more
\textquotedblleft$T$-biased\textquotedblright) form as well:

\begin{proposition}
\label{prop.specht.ET.ESET2}Let $\lambda$ be a partition of $n$. Let $S$ and
$T$ be two $n$-tableaux of shape $Y\left(  \lambda\right)  $. Let $f^{\lambda
}$ be the \# of standard tableaux of shape $Y\left(  \lambda\right)  $. Then:
\medskip

\textbf{(a)} If there exist no $\widetilde{r}\in\mathcal{R}\left(  T\right)  $
and $\widetilde{c}\in\mathcal{C}\left(  T\right)  $ such that $S=\widetilde{c}%
\widetilde{r}T$, then $\mathbf{E}_{S}\mathbf{E}_{T}=0$. \medskip

\textbf{(b)} If there exist $\widetilde{r}\in\mathcal{R}\left(  T\right)  $
and $\widetilde{c}\in\mathcal{C}\left(  T\right)  $ such that $S=\widetilde{c}%
\widetilde{r}T$, then these $\widetilde{r}$ and $\widetilde{c}$ satisfy
$\mathbf{E}_{S}\mathbf{E}_{T}=h^{\lambda}\left(  -1\right)  ^{\widetilde{c}%
}\widetilde{c}\widetilde{r}\mathbf{E}_{T}$, where we set $h^{\lambda}%
:=\dfrac{n!}{f^{\lambda}}$.
\end{proposition}

\begin{proof}
We shall derive this from Proposition \ref{prop.specht.ET.ESET}.

First, we observe that each permutation $c\in S_{n}$ satisfies
$cT=c\rightharpoonup T$ and thus%
\begin{equation}
\mathcal{R}\left(  cT\right)  =\mathcal{R}\left(  c\rightharpoonup T\right)
=c\mathcal{R}\left(  T\right)  c^{-1} \label{pf.prop.specht.ET.ESET2.RcT=}%
\end{equation}
(by Proposition \ref{prop.tableau.Sn-act.0} \textbf{(c)}, applied to $w=c$).
\medskip

\textbf{(a)} Assume that there exist no $\widetilde{r}\in\mathcal{R}\left(
T\right)  $ and $\widetilde{c}\in\mathcal{C}\left(  T\right)  $ such that
$S=\widetilde{c}\widetilde{r}T$. Then, it is not hard to see that there exist
no $r\in\mathcal{R}\left(  S\right)  $ and $c\in\mathcal{C}\left(  T\right)  $
such that $rS=cT$.

[\textit{Proof:} Assume the contrary. Thus, there exist some $r\in
\mathcal{R}\left(  S\right)  $ and $c\in\mathcal{C}\left(  T\right)  $ such
that $rS=cT$. Consider these $r$ and $c$. The $n$-tableau
$cT=rS=r\rightharpoonup S$ is row-equivalent to $S$ (by Proposition
\ref{prop.tableau.Sn-act.0} \textbf{(a)}, since $r\in\mathcal{R}\left(
S\right)  $), and thus we have $\mathcal{R}\left(  cT\right)  =\mathcal{R}%
\left(  S\right)  $ (by Proposition \ref{prop.tableau.roweq.H=H}
\textbf{(a)}). Hence, $\mathcal{R}\left(  S\right)  =\mathcal{R}\left(
cT\right)  =c\mathcal{R}\left(  T\right)  c^{-1}$ (by
(\ref{pf.prop.specht.ET.ESET2.RcT=})). Thus, $r\in\mathcal{R}\left(  S\right)
=c\mathcal{R}\left(  T\right)  c^{-1}$. In other words, $r=cr^{\prime}c^{-1}$
for some $r^{\prime}\in\mathcal{R}\left(  T\right)  $. Consider this
$r^{\prime}$. From $r=cr^{\prime}c^{-1}$, we obtain $r^{-1}=\left(
cr^{\prime}c^{-1}\right)  ^{-1}=c\left(  r^{\prime}\right)  ^{-1}c^{-1}$ and
therefore $c^{-1}r^{-1}c=\left(  r^{\prime}\right)  ^{-1}\in\mathcal{R}\left(
T\right)  $ (since $r^{\prime}\in\mathcal{R}\left(  T\right)  $, but
$\mathcal{R}\left(  T\right)  $ is a group).

However $rS=cT$ entails $S=r^{-1}cT=c\underbrace{c^{-1}r^{-1}c}_{=\left(
r^{\prime}\right)  ^{-1}}T=c\left(  r^{\prime}\right)  ^{-1}T$. Thus, there
exist $\widetilde{r}\in\mathcal{R}\left(  T\right)  $ and $\widetilde{c}%
\in\mathcal{C}\left(  T\right)  $ such that $S=\widetilde{c}\widetilde{r}T$
(namely, $\widetilde{c}=c$ and $\widetilde{r}=\left(  r^{\prime}\right)
^{-1}$). But we assumed that there exist no such $\widetilde{r}$ and
$\widetilde{c}$. Contradiction! This completes our proof.]

So we have shown that there exist no $r\in\mathcal{R}\left(  S\right)  $ and
$c\in\mathcal{C}\left(  T\right)  $ such that $rS=cT$. Thus, Proposition
\ref{prop.specht.ET.ESET} \textbf{(a)} yields $\mathbf{E}_{S}\mathbf{E}_{T}%
=0$. This proves Proposition \ref{prop.specht.ET.ESET2} \textbf{(a)}. \medskip

\textbf{(b)} Assume that there exist $\widetilde{r}\in\mathcal{R}\left(
T\right)  $ and $\widetilde{c}\in\mathcal{C}\left(  T\right)  $ such that
$S=\widetilde{c}\widetilde{r}T$. Consider these $\widetilde{r}$ and
$\widetilde{c}$. Set $h^{\lambda}:=\dfrac{n!}{f^{\lambda}}$. We must show that
$\mathbf{E}_{S}\mathbf{E}_{T}=h^{\lambda}\left(  -1\right)  ^{\widetilde{c}%
}\widetilde{c}\widetilde{r}\mathbf{E}_{T}$.

Set $c:=\widetilde{c}$ and $r:=\widetilde{c}\widetilde{r}^{-1}\widetilde{c}%
^{-1}$. We have $\widetilde{r}\in\mathcal{R}\left(  T\right)  $ and therefore
$\widetilde{r}^{-1}\in\mathcal{R}\left(  T\right)  $ (since $\mathcal{R}%
\left(  T\right)  $ is a group). Furthermore,%
\[
r=\underbrace{\widetilde{c}}_{=c}\widetilde{r}^{-1}\underbrace{\widetilde{c}%
^{-1}}_{\substack{=c^{-1}\\\text{(since }\widetilde{c}=c\text{)}%
}}=c\underbrace{\widetilde{r}^{-1}}_{\in\mathcal{R}\left(  T\right)  }%
c^{-1}\in c\mathcal{R}\left(  T\right)  c^{-1}=\mathcal{R}\left(  cT\right)
\ \ \ \ \ \ \ \ \ \ \left(  \text{by (\ref{pf.prop.specht.ET.ESET2.RcT=}%
)}\right)  .
\]
Thus, $r^{-1}\in\mathcal{R}\left(  cT\right)  $ (since $\mathcal{R}\left(
cT\right)  $ is a group). Also,%
\[
\underbrace{r}_{=\widetilde{c}\widetilde{r}^{-1}\widetilde{c}^{-1}%
}S=\underbrace{\widetilde{c}}_{=c}\widetilde{r}^{-1}\underbrace{\widetilde{c}%
^{-1}S}_{\substack{=\widetilde{r}T\\\text{(since }S=\widetilde{c}%
\widetilde{r}T\text{)}}}=c\underbrace{\widetilde{r}^{-1}\widetilde{r}}%
_{=1}T=cT.
\]

Hence, $S=r^{-1}cT=r^{-1}\rightharpoonup cT$. But the $n$-tableau
$r^{-1}\rightharpoonup cT$ is row-equivalent to $cT$ (by Proposition
\ref{prop.tableau.Sn-act.0} \textbf{(a)}, since $r^{-1}\in\mathcal{R}\left(
cT\right)  $). In other words, $S$ is row-equivalent to $cT$ (since
$S=r^{-1}\rightharpoonup cT$). Thus, $\mathcal{R}\left(  S\right)
=\mathcal{R}\left(  cT\right)  $ (by Proposition \ref{prop.tableau.roweq.H=H}
\textbf{(a)}). Hence, $r\in\mathcal{R}\left(  S\right)  $ (since
$r\in\mathcal{R}\left(  cT\right)  $).

Altogether, we have now shown that $r\in\mathcal{R}\left(  S\right)  $ and
$c\in\mathcal{C}\left(  T\right)  $ (since $c=\widetilde{c}\in\mathcal{C}%
\left(  T\right)  $) and $rS=cT$. Hence, Proposition \ref{prop.specht.ET.ESET}
\textbf{(b)} shows that $\mathbf{E}_{S}\mathbf{E}_{T}=h^{\lambda}\left(
-1\right)  ^{c}r^{-1}\mathbf{E}_{U}c$, where we set $U:=rS=cT$.

From $U=cT$, we obtain $\mathbf{E}_{U}=\mathbf{E}_{cT}=c\mathbf{E}_{T}c^{-1}$
(by Lemma \ref{lem.specht.ET.wT}), so that $\mathbf{E}_{U}c=c\mathbf{E}_{T}$.
Hence,%
\begin{align*}
\mathbf{E}_{S}\mathbf{E}_{T}  &  =h^{\lambda}\left(  -1\right)  ^{c}%
\underbrace{r^{-1}}_{\substack{=\left(  c\widetilde{r}^{-1}c^{-1}\right)
^{-1}\\\text{(since }r=c\widetilde{r}^{-1}c^{-1}\text{)}}%
}\underbrace{\mathbf{E}_{U}c}_{=c\mathbf{E}_{T}}=h^{\lambda}\left(  -1\right)
^{c}\underbrace{\left(  c\widetilde{r}^{-1}c^{-1}\right)  ^{-1}}%
_{=c\widetilde{r}c^{-1}}c\mathbf{E}_{T}\\
&  =h^{\lambda}\left(  -1\right)  ^{c}c\widetilde{r}\underbrace{c^{-1}c}%
_{=1}\mathbf{E}_{T}=h^{\lambda}\left(  -1\right)  ^{c}c\widetilde{r}%
\mathbf{E}_{T}=h^{\lambda}\left(  -1\right)  ^{\widetilde{c}}\widetilde{c}%
\widetilde{r}\mathbf{E}_{T}\ \ \ \ \ \ \ \ \ \ \left(  \text{since
}c=\widetilde{c}\right)  .
\end{align*}
This proves Proposition \ref{prop.specht.ET.ESET2} \textbf{(b)}.
\end{proof}

The formulas shown in this subsection answer the question \textquotedblleft
what is the product of two Young symmetrizers\textquotedblright. Roughly
speaking, the answer is \textquotedblleft$0$ in most cases, or some modified
version of a Young symmetrizer in the remaining cases\textquotedblright.

\subsubsection{The opposite Young symmetrizer}

The Young symmetrizer $\mathbf{E}_{T}$ was defined to be the product
$\nabla_{\operatorname*{Col}T}^{-}\nabla_{\operatorname*{Row}T}$. What can be
said about the opposite product, $\nabla_{\operatorname*{Row}T}\nabla
_{\operatorname*{Col}T}^{-}$ ? A few things, as it turns out:

\begin{proposition}
\label{prop.specht.FT.basics}Let $\lambda$ be a partition of $n$. Let $T$ be
an $n$-tableau of shape $\lambda$. Let%
\[
\mathbf{F}_{T}:=\nabla_{\operatorname*{Row}T}\nabla_{\operatorname*{Col}T}%
^{-}\in\mathcal{A}.
\]
Let $f^{\lambda}$ denote the \# of standard tableaux of shape $Y\left(
\lambda\right)  $. Then: \medskip

\textbf{(a)} The antipode $S$ of $\mathbf{k}\left[  S_{n}\right]  $ (which was
defined in Definition \ref{def.S.S}) satisfies $S\left(  \mathbf{E}%
_{T}\right)  =\mathbf{F}_{T}$. \medskip

\textbf{(b)} We have $\mathbf{F}_{T}^{2}=\dfrac{n!}{f^{\lambda}}\mathbf{F}%
_{T}$. \medskip

\textbf{(c)} Let $h^{\lambda}:=\dfrac{n!}{f^{\lambda}}$. If $h^{\lambda}$ is
invertible in $\mathbf{k}$, then $\mathcal{S}^{\lambda}\cong\mathcal{A}%
\mathbf{E}_{T}\cong\mathcal{A}\mathbf{F}_{T}$ as left $\mathcal{A}$-modules.
\end{proposition}

\begin{proof}
\textbf{(a)} The antipode $S$ is a $\mathbf{k}$-algebra anti-automorphism (by
Theorem \ref{thm.S.auto} \textbf{(a)}). Thus,%
\[
S\left(  \nabla_{\operatorname*{Col}T}^{-}\nabla_{\operatorname*{Row}%
T}\right)  =\underbrace{S\left(  \nabla_{\operatorname*{Row}T}\right)
}_{\substack{=\nabla_{\operatorname*{Row}T}\\\text{(by Proposition
\ref{prop.symmetrizers.antipode})}}}\ \ \underbrace{S\left(  \nabla
_{\operatorname*{Col}T}^{-}\right)  }_{\substack{=\nabla_{\operatorname*{Col}%
T}^{-}\\\text{(by Proposition \ref{prop.symmetrizers.antipode})}}%
}=\nabla_{\operatorname*{Row}T}\nabla_{\operatorname*{Col}T}^{-}%
=\mathbf{F}_{T}%
\]
(by the definition of $\mathbf{F}_{T}$). In view of $\mathbf{E}_{T}%
=\nabla_{\operatorname*{Col}T}^{-}\nabla_{\operatorname*{Row}T}$, we can
rewrite this as $S\left(  \mathbf{E}_{T}\right)  =\mathbf{F}_{T}$. This proves
Proposition \ref{prop.specht.FT.basics} \textbf{(a)}. \medskip

\textbf{(b)} The antipode $S$ is a $\mathbf{k}$-algebra anti-automorphism (by
Theorem \ref{thm.S.auto} \textbf{(a)}). Thus, $S\left(  \mathbf{E}%
_{T}\mathbf{E}_{T}\right)  =S\left(  \mathbf{E}_{T}\right)  S\left(
\mathbf{E}_{T}\right)  =\left(  S\left(  \mathbf{E}_{T}\right)  \right)  ^{2}%
$. In other words, $S\left(  \mathbf{E}_{T}^{2}\right)  =\left(  S\left(
\mathbf{E}_{T}\right)  \right)  ^{2}$ (since $\mathbf{E}_{T}\mathbf{E}%
_{T}=\mathbf{E}_{T}^{2}$).

But Theorem \ref{thm.specht.ETidp} yields $\mathbf{E}_{T}^{2}=\dfrac
{n!}{f^{\lambda}}\mathbf{E}_{T}$. Applying the map $S$ to both sides of this
equality, we obtain%
\[
S\left(  \mathbf{E}_{T}^{2}\right)  =S\left(  \dfrac{n!}{f^{\lambda}%
}\mathbf{E}_{T}\right)  =\dfrac{n!}{f^{\lambda}}S\left(  \mathbf{E}%
_{T}\right)  \ \ \ \ \ \ \ \ \ \ \left(  \text{since the map }S\text{ is
}\mathbf{k}\text{-linear}\right)  .
\]
Comparing this with $S\left(  \mathbf{E}_{T}^{2}\right)  =\left(  S\left(
\mathbf{E}_{T}\right)  \right)  ^{2}$, we find%
\[
\left(  S\left(  \mathbf{E}_{T}\right)  \right)  ^{2}=\dfrac{n!}{f^{\lambda}%
}S\left(  \mathbf{E}_{T}\right)  .
\]
Since $S\left(  \mathbf{E}_{T}\right)  =\mathbf{F}_{T}$ (by Proposition
\ref{prop.specht.FT.basics} \textbf{(a)}), we can rewrite this as
$\mathbf{F}_{T}^{2}=\dfrac{n!}{f^{\lambda}}\mathbf{F}_{T}$. This proves
Proposition \ref{prop.specht.FT.basics} \textbf{(b)}. \medskip

\textbf{(c)} Assume that $h^{\lambda}$ is invertible in $\mathbf{k}$. The
following is easy:

\begin{statement}
\textit{Claim 1:} For each $\mathbf{b}\in\mathcal{A}\mathbf{E}_{T}$, we have
$\mathbf{b}\nabla_{\operatorname*{Col}T}^{-}\in\mathcal{A}\mathbf{F}_{T}$.
\end{statement}

\begin{proof}
[Proof of Claim 1.]Let $\mathbf{b}\in\mathcal{A}\mathbf{E}_{T}$. Thus,
$\mathbf{b}=\mathbf{aE}_{T}$ for some $\mathbf{a}\in\mathcal{A}$. Consider
this $\mathbf{a}$. Then, $\mathbf{b}=\mathbf{aE}_{T}=\mathbf{a}\nabla
_{\operatorname*{Col}T}^{-}\nabla_{\operatorname*{Row}T}$ (since
$\mathbf{E}_{T}=\nabla_{\operatorname*{Col}T}^{-}\nabla_{\operatorname*{Row}%
T}$), so that%
\[
\mathbf{b}\nabla_{\operatorname*{Col}T}^{-}=\mathbf{a}\nabla
_{\operatorname*{Col}T}^{-}\underbrace{\nabla_{\operatorname*{Row}T}%
\nabla_{\operatorname*{Col}T}^{-}}_{=\mathbf{F}_{T}}=\underbrace{\mathbf{a}%
\nabla_{\operatorname*{Col}T}^{-}}_{\in\mathcal{A}}\mathbf{F}_{T}%
\in\mathcal{A}\mathbf{F}_{T}.
\]
This proves Claim 1.
\end{proof}

\begin{statement}
\textit{Claim 2:} For each $\mathbf{b}\in\mathcal{A}\mathbf{F}_{T}$, we have
$\mathbf{b}\nabla_{\operatorname*{Row}T}\in\mathcal{A}\mathbf{E}_{T}$.
\end{statement}

\begin{proof}
[Proof of Claim 2.]This is analogous to Claim 1, but now the roles of
$\nabla_{\operatorname*{Col}T}^{-}$ and $\nabla_{\operatorname*{Row}T}$ are
interchanged (and so are the roles of $\mathbf{E}_{T}$ and $\mathbf{F}_{T}$).
\end{proof}

\begin{statement}
\textit{Claim 3:} For each $\mathbf{b}\in\mathcal{A}\mathbf{E}_{T}$, we have
$\mathbf{bE}_{T}=h^{\lambda}\cdot\mathbf{b}$.
\end{statement}

\begin{proof}
[Proof of Claim 3.]Let $\mathbf{b}\in\mathcal{A}\mathbf{E}_{T}$. Thus,
$\mathbf{b}=\mathbf{aE}_{T}$ for some $\mathbf{a}\in\mathcal{A}$. Consider
this $\mathbf{a}$.

Theorem \ref{thm.specht.ETidp} yields $\mathbf{E}_{T}^{2}=\dfrac
{n!}{f^{\lambda}}\mathbf{E}_{T}$. In other words, $\mathbf{E}_{T}%
^{2}=h^{\lambda}\mathbf{E}_{T}$ (since $\dfrac{n!}{f^{\lambda}}=h^{\lambda}$).
But
\[
\underbrace{\mathbf{b}}_{=\mathbf{aE}_{T}}\mathbf{E}_{T}=\mathbf{a}%
\underbrace{\mathbf{E}_{T}\mathbf{E}_{T}}_{=\mathbf{E}_{T}^{2}=h^{\lambda
}\mathbf{E}_{T}}=h^{\lambda}\cdot\underbrace{\mathbf{aE}_{T}}_{=\mathbf{b}%
}=h^{\lambda}\cdot\mathbf{b}.
\]
This proves Claim 3.
\end{proof}

\begin{statement}
\textit{Claim 4:} For each $\mathbf{b}\in\mathcal{A}\mathbf{F}_{T}$, we have
$\mathbf{bF}_{T}=h^{\lambda}\cdot\mathbf{b}$.
\end{statement}

\begin{proof}
[Proof of Claim 4.]This is analogous to Claim 3, but now the roles of
$\nabla_{\operatorname*{Col}T}^{-}$ and $\nabla_{\operatorname*{Row}T}$ are
interchanged (and so are the roles of $\mathbf{E}_{T}$ and $\mathbf{F}_{T}$).
This time, we have to use Proposition \ref{prop.specht.FT.basics} \textbf{(b)}
instead of Theorem \ref{thm.specht.ETidp}.
\end{proof}

Each $\mathbf{b}\in\mathcal{A}\mathbf{E}_{T}$ satisfies $\mathbf{b}%
\nabla_{\operatorname*{Col}T}^{-}\in\mathcal{A}\mathbf{F}_{T}$ (by Claim 1).
Thus, we can define the map%
\begin{align*}
f:\mathcal{A}\mathbf{E}_{T}  &  \rightarrow\mathcal{A}\mathbf{F}_{T},\\
\mathbf{b}  &  \mapsto\mathbf{b}\nabla_{\operatorname*{Col}T}^{-}.
\end{align*}
Furthermore, each $\mathbf{b}\in\mathcal{A}\mathbf{F}_{T}$ satisfies
$\mathbf{b}\nabla_{\operatorname*{Row}T}\in\mathcal{A}\mathbf{E}_{T}$ (by
Claim 2) and therefore $\mathbf{b}\cdot\dfrac{1}{h^{\lambda}}\nabla
_{\operatorname*{Row}T}=\dfrac{1}{h^{\lambda}}\underbrace{\mathbf{b}%
\nabla_{\operatorname*{Row}T}}_{\in\mathcal{A}\mathbf{E}_{T}}\in
\mathcal{A}\mathbf{E}_{T}$ (since $\mathcal{A}\mathbf{E}_{T}$ is a
$\mathbf{k}$-module). Thus, we can define the map%
\begin{align*}
g:\mathcal{A}\mathbf{F}_{T}  &  \rightarrow\mathcal{A}\mathbf{E}_{T},\\
\mathbf{b}  &  \mapsto\mathbf{b}\cdot\dfrac{1}{h^{\lambda}}\nabla
_{\operatorname*{Row}T}.
\end{align*}
Consider these two maps $f$ and $g$. Both maps $f$ and $g$ are left
$\mathcal{A}$-linear (this is clear from their definitions), i.e., are left
$\mathcal{A}$-module morphisms. Moreover, we have $f\circ g=\operatorname*{id}%
$, since each $\mathbf{b}\in\mathcal{A}\mathbf{F}_{T}$ satisfies%
\begin{align*}
\left(  f\circ g\right)  \left(  \mathbf{b}\right)   &  =f\left(  g\left(
\mathbf{b}\right)  \right)  =\underbrace{g\left(  \mathbf{b}\right)
}_{\substack{=\mathbf{b}\cdot\dfrac{1}{h^{\lambda}}\nabla_{\operatorname*{Row}%
T}\\\text{(by the definition of }g\text{)}}}\cdot\,\nabla_{\operatorname*{Col}%
T}^{-}\ \ \ \ \ \ \ \ \ \ \left(  \text{by the definition of }f\right) \\
&  =\mathbf{b}\cdot\dfrac{1}{h^{\lambda}}\nabla_{\operatorname*{Row}T}%
\nabla_{\operatorname*{Col}T}^{-}=\dfrac{1}{h^{\lambda}}\mathbf{b}%
\underbrace{\nabla_{\operatorname*{Row}T}\nabla_{\operatorname*{Col}T}^{-}%
}_{=\mathbf{F}_{T}}=\dfrac{1}{h^{\lambda}}\underbrace{\mathbf{bF}_{T}%
}_{\substack{=h^{\lambda}\cdot\mathbf{b}\\\text{(by Claim 4)}}}\\
&  =\dfrac{1}{h^{\lambda}}h^{\lambda}\cdot\mathbf{b}=\mathbf{b}%
=\operatorname*{id}\left(  \mathbf{b}\right)  .
\end{align*}
Furthermore, we have $g\circ f=\operatorname*{id}$, since each $\mathbf{b}%
\in\mathcal{A}\mathbf{E}_{T}$ satisfies%
\begin{align*}
\left(  g\circ f\right)  \left(  \mathbf{b}\right)   &  =g\left(  f\left(
\mathbf{b}\right)  \right)  =\underbrace{f\left(  \mathbf{b}\right)
}_{\substack{=\mathbf{b}\nabla_{\operatorname*{Col}T}^{-}\\\text{(by the
definition of }f\text{)}}}\cdot\,\dfrac{1}{h^{\lambda}}\nabla
_{\operatorname*{Row}T}\ \ \ \ \ \ \ \ \ \ \left(  \text{by the definition of
}g\right) \\
&  =\mathbf{b}\nabla_{\operatorname*{Col}T}^{-}\cdot\dfrac{1}{h^{\lambda}%
}\nabla_{\operatorname*{Row}T}=\dfrac{1}{h^{\lambda}}\mathbf{b}%
\underbrace{\nabla_{\operatorname*{Col}T}^{-}\nabla_{\operatorname*{Row}T}%
}_{=\mathbf{E}_{T}}=\dfrac{1}{h^{\lambda}}\underbrace{\mathbf{bE}_{T}%
}_{\substack{=h^{\lambda}\cdot\mathbf{b}\\\text{(by Claim 3)}}}\\
&  =\dfrac{1}{h^{\lambda}}h^{\lambda}\cdot\mathbf{b}=\mathbf{b}%
=\operatorname*{id}\left(  \mathbf{b}\right)  .
\end{align*}

The maps $f$ and $g$ are mutually inverse (since $f\circ g=\operatorname*{id}$
and $g\circ f=\operatorname*{id}$), and thus are left $\mathcal{A}$-module
isomorphisms (since they are left $\mathcal{A}$-module morphisms). Thus,
$\mathcal{A}\mathbf{E}_{T}\cong\mathcal{A}\mathbf{F}_{T}$ as left
$\mathcal{A}$-modules.

But $T$ is an $n$-tableau of shape $\lambda$, that is, an $n$-tableau of shape
$Y\left(  \lambda\right)  $. Hence, (\ref{eq.def.specht.ET.defs.SD=AET})
(applied to $D=Y\left(  \lambda\right)  $) shows that $\mathcal{S}^{Y\left(
\lambda\right)  }\cong\mathcal{A}\mathbf{E}_{T}$ as left $\mathcal{A}%
$-modules. In other words, $\mathcal{S}^{\lambda}\cong\mathcal{A}%
\mathbf{E}_{T}$ as left $\mathcal{A}$-modules (since $\mathcal{S}^{\lambda}$
is a shorthand for $\mathcal{S}^{Y\left(  \lambda\right)  }$). Thus,
$\mathcal{S}^{\lambda}\cong\mathcal{A}\mathbf{E}_{T}\cong\mathcal{A}%
\mathbf{F}_{T}$ as left $\mathcal{A}$-modules. Proposition
\ref{prop.specht.FT.basics} \textbf{(c)} is thus proved.
\end{proof}

Proposition \ref{prop.specht.FT.basics} \textbf{(c)} would not be true if the
requirement that $h^{\lambda}$ be invertible was dropped. An example is given
in the following exercise:

\begin{exercise}
\label{exe.specht.FT.not-mod-p}Let $n\geq2$ and $\lambda=\left(  n-1,1\right)
$. Pick any $n$-tableau $T$ of shape $\lambda$. Let $\mathbf{F}_{T}$ be
defined as in Proposition \ref{prop.specht.FT.basics}. Show that: \medskip

\textbf{(a)} \fbox{1} The $\mathcal{A}$-module $\mathcal{A}\mathbf{E}_{T}$ is
isomorphic to the zero-sum representation $R\left(  \mathbf{k}^{n}\right)  $
(see Example \ref{exa.spechtmod.Yn-1-1}). \medskip

\textbf{(b)} \fbox{4} The $\mathcal{A}$-module $\mathcal{A}\mathbf{F}_{T}$ is
isomorphic to the quotient representation $\mathbf{k}^{n}/D\left(
\mathbf{k}^{n}\right)  $ (see Theorem \ref{thm.rep.G-rep.Sn-nat.quots}).
\medskip

\textbf{(c)} \fbox{1} Assume that $n\geq3$, and that $\mathbf{k}$ is a field
of characteristic $p$ with $p\mid n$. Then, we don't have $\mathcal{A}%
\mathbf{E}_{T}\cong\mathcal{A}\mathbf{F}_{T}$ as left $\mathcal{A}$-modules.
\end{exercise}

However, $\mathcal{A}\mathbf{E}_{T}$ and $\mathcal{A}\mathbf{F}_{T}$ are
always isomorphic as $\mathbf{k}$-modules (see Exercise
\ref{exe.specht.FT.iso-as-k-mod} below).

\subsection{Centralized Young symmetrizers}

Having studied the Young symmetrizers $\mathbf{E}_{T}$ corresponding to
straight-shaped $n$-tableaux $T$, we shall now explore certain
\textquotedblleft centralized\textquotedblright\ versions thereof, which are
the sums of the Young symmetrizers $\mathbf{E}_{T}$ over all $n$-tableaux $T$
of a given shape $Y\left(  \lambda\right)  $. We could just as well call them
\textquotedblleft symmetrized\textquotedblright\ or \textquotedblleft
averaged\textquotedblright, since they make up for the asymmetric placement of
the numbers in a given $n$-tableau by summing over \textbf{all} $n$-tableaux
of that shape. However, we will soon see why \textquotedblleft
centralized\textquotedblright\ is the most adequate word. \medskip

We will keep using Definition \ref{def.specht.ET.defs} throughout this
section. Also, recall that if $\lambda$ is a partition, then a tableau of
shape $\lambda$ means a tableau of shape $Y\left(  \lambda\right)  $.

\subsubsection{Definition and centrality}

\begin{definition}
\label{def.spechtmod.Elam.Elam}For any partition $\lambda$ of $n$, we define
the \emph{centralized Young symmetrizer}%
\begin{equation}
\mathbf{E}_{\lambda}:=\sum_{T\text{ is an }n\text{-tableau of shape }\lambda
}\mathbf{E}_{T}. \label{eq.def.spechtmod.Elam.Elam.def}%
\end{equation}
Note that this sum has $n!$ addends, since it ranges over all the $n!$ many
$n$-tableaux of shape $\lambda$ (not just the standard ones).
\end{definition}

\begin{example}
Let $n=3$. Then, using the formulas from Example \ref{exa.specht.ET.n=3}, we
can compute the centralized Young symmetrizers $\mathbf{E}_{\lambda}$
corresponding to all three partitions $\left(  3\right)  $, $\left(
1,1,1\right)  $ and $\left(  2,1\right)  $ of $3$. Namely, (where each $\sum$
sign means a sum over all the six orderings $\left(  i,j,k\right)  $ of
$\left[  3\right]  $) we have%
\begin{align*}
\mathbf{E}_{\left(  3\right)  }  &  =\sum\mathbf{E}_{ijk}=\sum\nabla
=6\nabla;\\
\mathbf{E}_{\left(  1,1,1\right)  }  &  =\sum\mathbf{E}_{i\backslash\backslash
j\backslash\backslash k}=\sum\nabla^{-}=6\nabla^{-};\\
\mathbf{E}_{\left(  2,1\right)  }  &  =\sum\mathbf{E}_{ij\backslash\backslash
k}=\sum\left(  1-t_{i,k}+t_{i,j}-\operatorname*{cyc}\nolimits_{i,j,k}\right)
\\
&  =6-\sum\operatorname*{cyc}\nolimits_{i,j,k}=6-3\left(  \operatorname*{cyc}%
\nolimits_{1,2,3}+\operatorname*{cyc}\nolimits_{1,3,2}\right)  .
\end{align*}

\end{example}

The word \textquotedblleft centralized\textquotedblright\ in the name of the
$\mathbf{E}_{\lambda}$ is well-deserved:

\begin{proposition}
\label{prop.spechtmod.Elam.incent}Let $\lambda$ be a partition of $n$. Then,
$\mathbf{E}_{\lambda}$ belongs to the center $Z\left(  \mathbf{k}\left[
S_{n}\right]  \right)  $ of $\mathbf{k}\left[  S_{n}\right]  $.
\end{proposition}

\begin{proof}
Let $w\in S_{n}$. Then, the action of $w\in S_{n}$ on the set $\left\{
n\text{-tableaux of shape }\lambda\right\}  $ is a permutation of this set.
Thus, the map%
\begin{align}
\left\{  n\text{-tableaux of shape }\lambda\right\}   &  \rightarrow\left\{
n\text{-tableaux of shape }\lambda\right\}  ,\nonumber\\
T  &  \mapsto wT \label{pf.prop.spechtmod.Elam.incent.bij}%
\end{align}
is a bijection (and its inverse is given by $T\mapsto w^{-1}T$). Now, we have%
\begin{align*}
w\mathbf{E}_{\lambda}w^{-1}  &  =w\left(  \sum_{T\text{ is an }n\text{-tableau
of shape }\lambda}\mathbf{E}_{T}\right)  w^{-1}\ \ \ \ \ \ \ \ \ \ \left(
\text{by (\ref{eq.def.spechtmod.Elam.Elam.def})}\right) \\
&  =\sum_{T\text{ is an }n\text{-tableau of shape }\lambda}%
\underbrace{w\mathbf{E}_{T}w^{-1}}_{\substack{=\mathbf{E}_{wT}\\\text{(by
Lemma \ref{lem.specht.ET.wT})}}}\\
&  =\sum_{T\text{ is an }n\text{-tableau of shape }\lambda}\mathbf{E}%
_{wT}=\sum_{T\text{ is an }n\text{-tableau of shape }\lambda}\mathbf{E}_{T}\\
&  \ \ \ \ \ \ \ \ \ \ \ \ \ \ \ \ \ \ \ \ \left(
\begin{array}
[c]{c}%
\text{here, we have substituted }T\text{ for }wT\text{ in the sum,}\\
\text{since the map (\ref{pf.prop.spechtmod.Elam.incent.bij}) is a bijection}%
\end{array}
\right) \\
&  =\mathbf{E}_{\lambda}\ \ \ \ \ \ \ \ \ \ \left(  \text{by
(\ref{eq.def.spechtmod.Elam.Elam.def})}\right)  .
\end{align*}

Forget that we fixed $w$. We thus have shown that each $w\in S_{n}$ satisfies
$w\mathbf{E}_{\lambda}w^{-1}=\mathbf{E}_{\lambda}$. Thus, Lemma
\ref{lem.center.waw} (applied to $G=S_{n}$ and $\mathbf{a}=\mathbf{E}%
_{\lambda}$) shows that $\mathbf{E}_{\lambda}\in Z\left(  \mathbf{k}\left[
S_{n}\right]  \right)  $. This proves Proposition
\ref{prop.spechtmod.Elam.incent}.
\end{proof}

If this centrality was the only property of the $\mathbf{E}_{\lambda}$, then
it would be worth an exercise at best. But in fact, the $\mathbf{E}_{\lambda}$
play an important role in the center of $\mathbf{k}\left[  S_{n}\right]  $:

\begin{theorem}
\label{thm.spechtmod.Elam.center}Assume that $\mathbf{k}$ is a field of
characteristic $0$. Then: \medskip

\textbf{(a)} The family%
\[
\left(  \mathbf{E}_{\lambda}\right)  _{\lambda\text{ is a partition of }n}%
\]
is a basis of the $\mathbf{k}$-vector space $Z\left(  \mathbf{k}\left[
S_{n}\right]  \right)  $. \medskip

\textbf{(b)} Each $\mathbf{E}_{\lambda}$ is quasi-idempotent, with
$\mathbf{E}_{\lambda}^{2}=c_{\lambda}\mathbf{E}_{\lambda}$ for a certain
nonzero $c_{\lambda}\in\mathbb{Q}$. \medskip

\textbf{(c)} The $\mathbf{E}_{\lambda}$ are orthogonal, in the sense that
$\mathbf{E}_{\lambda}\mathbf{E}_{\mu}=0$ whenever $\lambda\neq\mu$.
\end{theorem}

We will prove this theorem soon, along with some slightly stronger claims.
Once it is proved, it will follow easily that $Z\left(  \mathbf{k}\left[
S_{n}\right]  \right)  $ is isomorphic to $\mathbf{k}\times\mathbf{k}%
\times\cdots\times\mathbf{k}$ as a $\mathbf{k}$-algebra when $\mathbf{k}$ is a
field of characteristic $0$. From here, it won't be far to Theorem
\ref{thm.AWS.demo} \textbf{(b)}.

First, however, let us rewrite a centralized Young symmetrizer $\mathbf{E}%
_{\lambda}$:

\begin{proposition}
\label{prop.spechtmod.Elam.as-sum-w}Let $\lambda$ be a partition of $n$. Let
$P$ be an $n$-tableau of shape $Y\left(  \lambda\right)  $. Then,%
\[
\mathbf{E}_{\lambda}=\sum_{w\in S_{n}}w\mathbf{E}_{P}w^{-1}.
\]

\end{proposition}

\begin{proof}
Proposition \ref{prop.tableau.Sn-act.1} \textbf{(d)} (applied to $D=Y\left(
\lambda\right)  $ and $T=P$) shows that the map
\begin{align*}
f:S_{n}  &  \rightarrow\left\{  n\text{-tableaux of shape }Y\left(
\lambda\right)  \right\}  ,\\
w  &  \mapsto w\rightharpoonup P
\end{align*}
is a bijection. Consider this bijection $f$.

But (\ref{eq.def.spechtmod.Elam.Elam.def}) shows that%
\begin{align*}
\mathbf{E}_{\lambda}  &  =\sum_{T\text{ is an }n\text{-tableau of shape
}\lambda}\mathbf{E}_{T}\\
&  =\sum_{T\text{ is an }n\text{-tableau of shape }Y\left(  \lambda\right)
}\mathbf{E}_{T}\ \ \ \ \ \ \ \ \ \ \left(
\begin{array}
[c]{c}%
\text{since an }n\text{-tableau of shape }\lambda\\
\text{is the same thing as}\\
\text{an }n\text{-tableau of shape }Y\left(  \lambda\right)
\end{array}
\right) \\
&  =\sum_{T\in\left\{  n\text{-tableaux of shape }Y\left(  \lambda\right)
\right\}  }\mathbf{E}_{T}\\
&  =\sum_{w\in S_{n}}\underbrace{\mathbf{E}_{f\left(  w\right)  }%
}_{\substack{=\mathbf{E}_{wP}\\\text{(since the definition}\\\text{of }f\text{
yields }f\left(  w\right)  =w\rightharpoonup P=wP\\\text{(since }wP\text{ is
shorthand for }w\rightharpoonup P\text{))}}}\\
&  \ \ \ \ \ \ \ \ \ \ \ \ \ \ \ \ \ \ \ \ \left(
\begin{array}
[c]{c}%
\text{here, we have substituted }f\left(  w\right)  \text{ for }T\text{ in the
sum,}\\
\text{since the map }f:S_{n}\rightarrow\left\{  n\text{-tableaux of shape
}Y\left(  \lambda\right)  \right\} \\
\text{is a bijection}%
\end{array}
\right) \\
&  =\sum_{w\in S_{n}}\underbrace{\mathbf{E}_{wP}}_{\substack{=w\mathbf{E}%
_{P}w^{-1}\\\text{(by Lemma \ref{lem.specht.ET.wT},}\\\text{applied to
}T=P\text{)}}}=\sum_{w\in S_{n}}w\mathbf{E}_{P}w^{-1}.
\end{align*}
This proves Proposition \ref{prop.spechtmod.Elam.as-sum-w}.
\end{proof}

\subsubsection{Proof of the basic claims}

We shall prove Theorem \ref{thm.spechtmod.Elam.center} step by step. We begin
with the orthogonality part:

\begin{proposition}
\label{prop.spechtmod.Elam.orth}Let $\lambda$ and $\mu$ be two distinct
partitions of $n$. Then, $\mathbf{E}_{\lambda}\mathbf{E}_{\mu}=0$.
\end{proposition}

\begin{proof}
The double zero-sandwich lemma (Proposition \ref{prop.specht.ET.ESaET}) says
that
\begin{equation}
\mathbf{E}_{S}\mathbf{aE}_{T}=0 \label{pf.prop.spechtmod.Elam.orth.1}%
\end{equation}
for all $n$-tableaux $S$ of shape $Y\left(  \lambda\right)  $ and all
$n$-tableaux $T$ of shape $Y\left(  \mu\right)  $ and all $\mathbf{a}%
\in\mathbf{k}\left[  S_{n}\right]  $. But the definition of $\mathbf{E}%
_{\lambda}$ yields%
\begin{align*}
\mathbf{E}_{\lambda}  &  =\sum_{T\text{ is an }n\text{-tableau of shape
}\lambda}\mathbf{E}_{T}=\sum_{S\text{ is an }n\text{-tableau of shape }%
\lambda}\mathbf{E}_{S}\\
&  =\sum_{S\text{ is an }n\text{-tableau of shape }Y\left(  \lambda\right)
}\mathbf{E}_{S}%
\end{align*}
(since \textquotedblleft shape $\lambda$\textquotedblright\ means
\textquotedblleft shape $Y\left(  \lambda\right)  $\textquotedblright), and
similarly we find%
\[
\mathbf{E}_{\mu}=\sum_{T\text{ is an }n\text{-tableau of shape }Y\left(
\mu\right)  }\mathbf{E}_{T}.
\]
Multiplying these two equalities, we obtain%
\begin{align*}
\mathbf{E}_{\lambda}\mathbf{E}_{\mu}  &  =\left(  \sum_{S\text{ is an
}n\text{-tableau of shape }Y\left(  \lambda\right)  }\mathbf{E}_{S}\right)
\left(  \sum_{T\text{ is an }n\text{-tableau of shape }Y\left(  \mu\right)
}\mathbf{E}_{T}\right) \\
&  =\sum_{S\text{ is an }n\text{-tableau of shape }Y\left(  \lambda\right)
}\ \ \sum_{T\text{ is an }n\text{-tableau of shape }Y\left(  \mu\right)
}\underbrace{\mathbf{E}_{S}\mathbf{E}_{T}}_{\substack{=\mathbf{E}%
_{S}1\mathbf{E}_{T}\\=0\\\text{(by (\ref{pf.prop.spechtmod.Elam.orth.1}%
),}\\\text{applied to }\mathbf{a}=1\text{)}}}=0.
\end{align*}
This proves Proposition \ref{prop.spechtmod.Elam.orth}.
\end{proof}

Thus, Theorem \ref{thm.spechtmod.Elam.center} \textbf{(c)} is proven. \medskip

The following lemma is an easy consequence of Lemma
\ref{lem.specht.ET.1-coord}:

\begin{lemma}
\label{lem.specht.Elam.1-coord}Let $\lambda$ be a partition of $n$. Then,
$\left[  1\right]  \mathbf{E}_{\lambda}=n!$.
\end{lemma}

\begin{proof}
The $n$-tableaux of shape $\lambda$ are just the bijections from $Y\left(
\lambda\right)  $ to $\left[  n\right]  $. Since both $Y\left(  \lambda
\right)  $ and $\left[  n\right]  $ are $n$-element sets (because $\lambda$ is
a partition of $n$), there are exactly $n!$ such bijections. In other words,
there are $n!$ many $n$-tableaux of shape $\lambda$.

From (\ref{eq.def.spechtmod.Elam.Elam.def}), we obtain
\begin{align*}
\left[  1\right]  \mathbf{E}_{\lambda}  &  =\left[  1\right]  \sum_{T\text{ is
an }n\text{-tableau of shape }\lambda}\mathbf{E}_{T}\\
&  =\sum_{T\text{ is an }n\text{-tableau of shape }\lambda}\underbrace{\left[
1\right]  \mathbf{E}_{T}}_{\substack{=1\\\text{(by Lemma
\ref{lem.specht.ET.1-coord})}}}\\
&  =\sum_{T\text{ is an }n\text{-tableau of shape }\lambda}1=n!
\end{align*}
(since there are $n!$ many $n$-tableaux of shape $\lambda$). This proves Lemma
\ref{lem.specht.Elam.1-coord}.
\end{proof}

The next lemma looks hardly like progress, but it will serve a crucial role in
proving the rest of the theorem:

\begin{lemma}
\label{lem.spechtmod.Elam.sqn0}Let $\lambda$ be a partition of $n$. Then,
$\mathbf{E}_{\lambda}^{2}\neq0$ if $\mathbf{k}$ is a field of characteristic
$0$.
\end{lemma}

\begin{proof}
Assume that $\mathbf{k}$ is a field of characteristic $0$. Then, the left
ideal $\mathcal{A}\mathbf{E}_{\lambda}$ of $\mathcal{A}$ is a
finite-dimensional $\mathbf{k}$-vector space (since it is a subspace of the
finite-dimensional $\mathbf{k}$-vector space $\mathcal{A}$), and thus is a
free $\mathbf{k}$-module of rank $d$ for some $d\in\mathbb{N}$ (namely, for
$d=\dim\left(  \mathcal{A}\mathbf{E}_{\lambda}\right)  $). Consider this $d$.

We must prove that $\mathbf{E}_{\lambda}^{2}\neq0$. Assume the contrary. Thus,
$\mathbf{E}_{\lambda}^{2}=0=0\mathbf{E}_{\lambda}$. Hence, Lemma
\ref{lem.linalg.quasi-idp1} (applied to $G=S_{n}$ and $\mathbf{a}%
=\mathbf{E}_{\lambda}$ and $\kappa=0$) yields
\[
d\cdot0=\underbrace{\left\vert S_{n}\right\vert }_{=n!}\cdot
\underbrace{\left[  1\right]  \mathbf{E}_{\lambda}}_{\substack{=n!\\\text{(by
Lemma \ref{lem.specht.Elam.1-coord})}}}=n!\cdot n!=n!^{2}\neq
0\ \ \ \ \ \ \ \ \ \ \text{in }\mathbf{k}%
\]
(since $\mathbf{k}$ has characteristic $0$). This contradicts $d\cdot0=0$.
This contradiction shows that our assumption was false, so that $\mathbf{E}%
_{\lambda}^{2}\neq0$ is proved. Thus, Lemma \ref{lem.spechtmod.Elam.sqn0} is proved.
\end{proof}

We can now prove Theorem \ref{thm.spechtmod.Elam.center} \textbf{(a)}:

\begin{theorem}
\label{thm.spechtmod.Elam.basis}Assume that $\mathbf{k}$ is a field of
characteristic $0$. Then, the family $\left(  \mathbf{E}_{\lambda}\right)
_{\lambda\text{ is a partition of }n}$ is a basis of the $\mathbf{k}$-vector
space $Z\left(  \mathbf{k}\left[  S_{n}\right]  \right)  $.
\end{theorem}

\begin{proof}
Let $P$ be the set of all partitions of $n$. Thus, we must show that the
family $\left(  \mathbf{E}_{\lambda}\right)  _{\lambda\in P}$ is a basis of
the $\mathbf{k}$-vector space $Z\left(  \mathbf{k}\left[  S_{n}\right]
\right)  $. First of all, we observe that this family really consists of
vectors in $Z\left(  \mathbf{k}\left[  S_{n}\right]  \right)  $, since
Proposition \ref{prop.spechtmod.Elam.incent} shows that each $\lambda\in P$
satisfies $\mathbf{E}_{\lambda}\in Z\left(  \mathbf{k}\left[  S_{n}\right]
\right)  $.

Note that $P$ is a finite set, and its size is $\left\vert P\right\vert
=\left(  \text{\# of partitions of }n\right)  $ (since $P$ is the set of all
partitions of $n$).

Theorem \ref{thm.center.conjsums} (applied to $G=S_{n}$) shows that the
conjugacy class sums of $S_{n}$ form a basis of the $\mathbf{k}$-module
$Z\left(  \mathbf{k}\left[  S_{n}\right]  \right)  $. Hence, the dimension of
$Z\left(  \mathbf{k}\left[  S_{n}\right]  \right)  $ as a $\mathbf{k}$-vector
space is%
\begin{align*}
&  \left(  \text{\# of finite conjugacy classes of }S_{n}\right) \\
&  \ \ \ \ \ \ \ \ \ \ \ \ \ \ \ \ \ \ \ \ \left(
\begin{array}
[c]{c}%
\text{since the conjugacy class sums are}\\
\text{indexed by the finite conjugacy classes}%
\end{array}
\right) \\
&  =\left(  \text{\# of conjugacy classes of }S_{n}\right) \\
&  \ \ \ \ \ \ \ \ \ \ \ \ \ \ \ \ \ \ \ \ \left(  \text{since all conjugacy
classes of }S_{n}\text{ are finite}\right) \\
&  =\left(  \text{\# of partitions of }\left\vert \left[  n\right]
\right\vert \right) \\
&  \ \ \ \ \ \ \ \ \ \ \ \ \ \ \ \ \ \ \ \ \left(  \text{by Corollary
\ref{cor.type.conjclasses} \textbf{(c)}, applied to }X=\left[  n\right]
\right) \\
&  =\left(  \text{\# of partitions of }n\right)  \ \ \ \ \ \ \ \ \ \ \left(
\text{since }\left\vert \left[  n\right]  \right\vert =n\right) \\
&  =\left\vert P\right\vert \ \ \ \ \ \ \ \ \ \ \left(  \text{since the set of
all partitions of }n\text{ is }P\right)  .
\end{align*}
In other words, the $\mathbf{k}$-vector space $Z\left(  \mathbf{k}\left[
S_{n}\right]  \right)  $ is $\left\vert P\right\vert $-dimensional. Thus, our
family $\left(  \mathbf{E}_{\lambda}\right)  _{\lambda\in P}$ has the right
size to be a basis of $Z\left(  \mathbf{k}\left[  S_{n}\right]  \right)  $.
Therefore, in order to show that it is a basis of $Z\left(  \mathbf{k}\left[
S_{n}\right]  \right)  $, it suffices to prove that it is $\mathbf{k}%
$-linearly independent (by a basic fact from linear algebra, since we are
working over a field\footnote{We will explain this argument in more detail at
the end of this proof.}). Let us do this:

\begin{statement}
\textit{Claim 1:} The family $\left(  \mathbf{E}_{\lambda}\right)
_{\lambda\in P}$ is $\mathbf{k}$-linearly independent.
\end{statement}

\begin{proof}
[Proof of Claim 1.]Let $\left(  \alpha_{\mu}\right)  _{\mu\in P}\in
\mathbf{k}^{P}$ be a family of scalars such that%
\begin{equation}
\sum_{\mu\in P}\alpha_{\mu}\mathbf{E}_{\mu}=0.
\label{pf.thm.spechtmod.Elam.basis.sum=0}%
\end{equation}
We shall show that all the coefficients $\alpha_{\mu}$ are zero. Indeed, fix
any $\lambda\in P$. Then, $\lambda$ is a partition of $n$. However, from
(\ref{pf.thm.spechtmod.Elam.basis.sum=0}), we obtain $\mathbf{E}_{\lambda}%
\sum_{\mu\in P}\alpha_{\mu}\mathbf{E}_{\mu}=\mathbf{E}_{\lambda}\cdot0=0$.
Hence,%
\begin{align*}
0  &  =\mathbf{E}_{\lambda}\sum_{\mu\in P}\alpha_{\mu}\mathbf{E}_{\mu}%
=\sum_{\mu\in P}\alpha_{\mu}\mathbf{E}_{\lambda}\mathbf{E}_{\mu}\\
&  =\alpha_{\lambda}\underbrace{\mathbf{E}_{\lambda}\mathbf{E}_{\lambda}%
}_{=\mathbf{E}_{\lambda}^{2}}+\sum_{\substack{\mu\in P;\\\mu\neq\lambda
}}\alpha_{\mu}\underbrace{\mathbf{E}_{\lambda}\mathbf{E}_{\mu}}%
_{\substack{=0\\\text{(by Proposition \ref{prop.spechtmod.Elam.orth}%
,}\\\text{since }\lambda\text{ and }\mu\text{ are two}\\\text{distinct
partitions of }n\text{)}}}\\
&  \ \ \ \ \ \ \ \ \ \ \ \ \ \ \ \ \ \ \ \ \left(
\begin{array}
[c]{c}%
\text{here, we have split off the}\\
\text{addend for }\mu=\lambda\text{ from the sum}%
\end{array}
\right) \\
&  =\alpha_{\lambda}\mathbf{E}_{\lambda}^{2}+\underbrace{\sum_{\substack{\mu
\in P;\\\mu\neq\lambda}}\alpha_{\mu}0}_{=0}=\alpha_{\lambda}\mathbf{E}%
_{\lambda}^{2}.
\end{align*}

But Lemma \ref{lem.spechtmod.Elam.sqn0} yields $\mathbf{E}_{\lambda}^{2}\neq
0$. If $\alpha_{\lambda}$ was nonzero, then we could divide the equality
$0=\alpha_{\lambda}\mathbf{E}_{\lambda}^{2}$ by $\alpha_{\lambda}$ (since
$\mathbf{k}$ is a field!) and obtain $0=\mathbf{E}_{\lambda}^{2}$, which would
contradict $\mathbf{E}_{\lambda}^{2}\neq0$. Thus, $\alpha_{\lambda}$ cannot be
nonzero. Hence, $\alpha_{\lambda}=0$.

Forget that we fixed $\lambda$. We thus have shown that each $\lambda\in P$
satisfies $\alpha_{\lambda}=0$. In other words, each $\mu\in P$ satisfies
$\alpha_{\mu}=0$.

Forget that we fixed $\left(  \alpha_{\mu}\right)  _{\mu\in P}$. We thus have
shown that if $\left(  \alpha_{\mu}\right)  _{\mu\in P}\in\mathbf{k}^{P}$ is a
family of scalars such that $\sum_{\mu\in P}\alpha_{\mu}\mathbf{E}_{\mu}=0$,
then each $\mu\in P$ satisfies $\alpha_{\mu}=0$. In other words, the family
$\left(  \mathbf{E}_{\mu}\right)  _{\mu\in P}$ is $\mathbf{k}$-linearly
independent. In other words, the family $\left(  \mathbf{E}_{\lambda}\right)
_{\lambda\in P}$ is $\mathbf{k}$-linearly independent. This proves Claim 1.
\end{proof}

However, $\mathbf{k}$ is a field. Thus, a classical fact from linear algebra
(see, e.g., \cite[2.38]{Axler24} or \cite[\textquotedblleft The dimension of a
vector space\textquotedblright, Theorem 2.10, last sentence]{Conrad}) says
that if $V$ is a $d$-dimensional $\mathbf{k}$-vector space, then any
$\mathbf{k}$-linearly independent family of $d$ vectors in $V$ must be a basis
of $V$. Applying this to $V=Z\left(  \mathbf{k}\left[  S_{n}\right]  \right)
$ and $d=\left\vert P\right\vert $ and to the $\mathbf{k}$-linearly
independent family $\left(  \mathbf{E}_{\lambda}\right)  _{\lambda\in P}$, we
thus conclude that the family $\left(  \mathbf{E}_{\lambda}\right)
_{\lambda\in P}$ must be a basis of $Z\left(  \mathbf{k}\left[  S_{n}\right]
\right)  $. This proves Theorem \ref{thm.spechtmod.Elam.basis}.
\end{proof}

Next let us prove Theorem \ref{thm.spechtmod.Elam.center} \textbf{(b)}:

\begin{proposition}
\label{prop.spechtmod.Elam.idp-weak}Assume that $\mathbf{k}$ is a commutative
$\mathbb{Q}$-algebra. Let $\lambda$ be a partition of $n$. Then,
$\mathbf{E}_{\lambda}^{2}=c_{\lambda}\mathbf{E}_{\lambda}$ for a certain
nonzero $c_{\lambda}\in\mathbb{Q}$.
\end{proposition}

\begin{proof}
If we can prove Proposition \ref{prop.spechtmod.Elam.idp-weak} for
$\mathbf{k}=\mathbb{Q}$, then it will automatically follow that this
proposition holds for any commutative $\mathbb{Q}$-algebra $\mathbf{k}$ (since
there is a canonical ring morphism $\mathbb{Q}\rightarrow\mathbf{k}$, which
induces a base change morphism $\mathbb{Q}\left[  S_{n}\right]  \rightarrow
\mathbf{k}\left[  S_{n}\right]  $, and we can apply the latter base change
morphism to the equality $\mathbf{E}_{\lambda}^{2}=c_{\lambda}\mathbf{E}%
_{\lambda}$ in $\mathbb{Q}\left[  S_{n}\right]  $ to obtain the same equality
in $\mathbf{k}\left[  S_{n}\right]  $\ \ \ \ \footnote{This is analogous to
the base change argument that we did at the end of our proof of Theorem
\ref{thm.specht.ETidp}.}).

Thus, we can WLOG assume that $\mathbf{k}=\mathbb{Q}$. Assume this. Hence,
$\mathbf{k}$ is a field of characteristic $0$.

Theorem \ref{thm.spechtmod.Elam.basis} tells us that the family $\left(
\mathbf{E}_{\mu}\right)  _{\mu\text{ is a partition of }n}$ is a basis of the
$\mathbf{k}$-vector space $Z\left(  \mathbf{k}\left[  S_{n}\right]  \right)
$. Thus, in particular, this family spans $Z\left(  \mathbf{k}\left[
S_{n}\right]  \right)  $. Hence, the vector $1$ (which clearly belongs to
$Z\left(  \mathbf{k}\left[  S_{n}\right]  \right)  $) is a linear combination
of this family. In other words, we can write $1$ in the form%
\[
1=\sum_{\mu\text{ is a partition of }n}\alpha_{\mu}\mathbf{E}_{\mu
}\ \ \ \ \ \ \ \ \ \ \text{for some scalars }\alpha_{\mu}\in\mathbf{k}.
\]
Consider these scalars $\alpha_{\mu}$. Multiplying the equality $1=\sum
_{\mu\text{ is a partition of }n}\alpha_{\mu}\mathbf{E}_{\mu}$ by
$\mathbf{E}_{\lambda}$ from the left, we get%
\begin{align*}
\mathbf{E}_{\lambda}\cdot1  &  =\mathbf{E}_{\lambda}\sum_{\mu\text{ is a
partition of }n}\alpha_{\mu}\mathbf{E}_{\mu}=\sum_{\mu\text{ is a partition of
}n}\alpha_{\mu}\mathbf{E}_{\lambda}\mathbf{E}_{\mu}\\
&  =\alpha_{\lambda}\underbrace{\mathbf{E}_{\lambda}\mathbf{E}_{\lambda}%
}_{=\mathbf{E}_{\lambda}^{2}}+\sum_{\substack{\mu\text{ is a partition of
}n;\\\mu\neq\lambda}}\alpha_{\mu}\underbrace{\mathbf{E}_{\lambda}%
\mathbf{E}_{\mu}}_{\substack{=0\\\text{(by Proposition
\ref{prop.spechtmod.Elam.orth})}}}\\
&  \ \ \ \ \ \ \ \ \ \ \ \ \ \ \ \ \ \ \ \ \left(
\begin{array}
[c]{c}%
\text{here, we have split off the}\\
\text{addend for }\mu=\lambda\text{ from the sum}%
\end{array}
\right) \\
&  =\alpha_{\lambda}\mathbf{E}_{\lambda}^{2}+\underbrace{\sum_{\substack{\mu
\text{ is a partition of }n;\\\mu\neq\lambda}}\alpha_{\mu}0}_{=0}%
=\alpha_{\lambda}\mathbf{E}_{\lambda}^{2}.
\end{align*}
Hence, $\alpha_{\lambda}\mathbf{E}_{\lambda}^{2}=\mathbf{E}_{\lambda}%
\cdot1=\mathbf{E}_{\lambda}$.

But Lemma \ref{lem.spechtmod.Elam.sqn0} yields $\mathbf{E}_{\lambda}^{2}\neq
0$, and thus $\mathbf{E}_{\lambda}\neq0$. Hence, $\alpha_{\lambda}%
\mathbf{E}_{\lambda}^{2}=\mathbf{E}_{\lambda}\neq0$, so that $\alpha_{\lambda
}\neq0$. Since $\alpha_{\lambda}\in\mathbf{k}=\mathbb{Q}$, the inverse
$\dfrac{1}{\alpha_{\lambda}}$ is a well-defined nonzero element of
$\mathbb{Q}$, and we can divide the equality $\alpha_{\lambda}\mathbf{E}%
_{\lambda}^{2}=\mathbf{E}_{\lambda}$ by $\alpha_{\lambda}$ and obtain
$\mathbf{E}_{\lambda}^{2}=\dfrac{1}{\alpha_{\lambda}}\mathbf{E}_{\lambda}$.
Hence, $\mathbf{E}_{\lambda}^{2}=c_{\lambda}\mathbf{E}_{\lambda}$ for a
certain nonzero $c_{\lambda}\in\mathbb{Q}$ (namely, for $c_{\lambda}=\dfrac
{1}{\alpha_{\lambda}}$). This proves Proposition
\ref{prop.spechtmod.Elam.idp-weak}.
\end{proof}

\subsubsection{A stronger quasi-idempotence}

We have thus proved Theorem \ref{thm.spechtmod.Elam.center} completely. This
is a good beginning, but we can improve the theorem in several directions, and
prove some further properties of centralized Young symmetrizers. In doing so,
we will use the following notations:

\begin{definition}
\label{def.spechtmod.flam}Let $\lambda$ be a partition of $n$. Then: \medskip

\textbf{(a)} We let $f^{\lambda}$ be the \# of standard tableaux of shape
$Y\left(  \lambda\right)  $. \medskip

\textbf{(b)} We set $h^{\lambda}:=\dfrac{n!}{f^{\lambda}}$. This is a positive
integer (by Theorem \ref{thm.specht.ETidp}).
\end{definition}

These notations have been already used above, but we shall now use them
without explicit mention throughout the rest of this chapter.

\begin{example}
\label{exa.spechtmod.flam-hlam-6}Here are lists of all partitions of $n$ for
$n\leq5$ with the corresponding numbers $f^{\lambda}$ and $h^{\lambda}$:%
\[%
\begin{tabular}
[c]{|c||c||c||c|c||c|c|c||}\hline
$\lambda$ & $\left(  {}\right)  $ & $\left(  1\right)  $ & $\left(  2\right)
$ & $\left(  1,1\right)  $ & $\left(  3\right)  $ & $\left(  2,1\right)  $ &
$\left(  1,1,1\right)  $\\\hline
$f^{\lambda}$ & $1$ & $1$ & $1$ & $1$ & $1$ & $2$ & $1$\\\hline
$h^{\lambda}$ & $1$ & $1$ & $2$ & $2$ & $6$ & $3$ & $6$\\\hline
\end{tabular}
\]%
\[%
\begin{tabular}
[c]{|c||c|c|c|c|c||}\hline
$\lambda$ & $\left(  4\right)  $ & $\left(  3,1\right)  $ & $\left(
2,2\right)  $ & $\left(  2,1,1\right)  $ & $\left(  1,1,1,1\right)  $\\\hline
$f^{\lambda}$ & $1$ & $3$ & $2$ & $3$ & $1$\\\hline
$h^{\lambda}$ & $24$ & $8$ & $12$ & $8$ & $24$\\\hline
\end{tabular}
\]%
\[%
\begin{tabular}
[c]{|c||c|c|c|c|c|c|c||}\hline
$\lambda$ & $\left(  5\right)  $ & $\left(  4,1\right)  $ & $\left(
3,2\right)  $ & $\left(  3,1,1\right)  $ & $\left(  2,2,1\right)  $ & $\left(
2,1,1,1\right)  $ & $\left(  1,1,1,1,1\right)  $\\\hline
$f^{\lambda}$ & $1$ & $4$ & $5$ & $6$ & $5$ & $4$ & $1$\\\hline
$h^{\lambda}$ & $120$ & $30$ & $24$ & $20$ & $24$ & $30$ & $120$\\\hline
\end{tabular}
\]

\end{example}

\begin{remark}
Let $\lambda$ be a partition of $n$. Then, the hook length formula (Theorem
\ref{thm.tableau.hlf}) shows that the number $h^{\lambda}=\dfrac
{n!}{f^{\lambda}}$ equals the product $%
%TCIMACRO{\dprod \limits_{c\in Y\left(  \lambda\right)  }}%
%BeginExpansion
{\displaystyle\prod\limits_{c\in Y\left(  \lambda\right)  }}
%EndExpansion
\left\vert H_{\lambda}\left(  c\right)  \right\vert $ of the sizes of the
hooks of all cells $c\in Y\left(  \lambda\right)  $. We shall not need this
fact, but we are mentioning it here since it gives the easiest way to compute
$h^{\lambda}$.
\end{remark}

Our first improvement on Theorem \ref{thm.spechtmod.Elam.center} is an
explicit formula for the number $c_{\lambda}$ in Theorem
\ref{thm.spechtmod.Elam.center} \textbf{(b)} (or in Proposition
\ref{prop.spechtmod.Elam.idp-weak}). Namely, we claim that $c_{\lambda
}=\left(  h^{\lambda}\right)  ^{2}$. In other words, we claim the following
(using the notations of Definition \ref{def.spechtmod.flam}):

\begin{theorem}
\label{thm.spechtmod.Elam.idp-strong}Let $\lambda$ be a partition of $n$.
Then, $\mathbf{E}_{\lambda}^{2}=\left(  h^{\lambda}\right)  ^{2}%
\mathbf{E}_{\lambda}$.
\end{theorem}

To prove this, we need another fact (which is also useful in its own right):

\begin{proposition}
\label{prop.spechtmod.Elam.ETEl}Let $\lambda$ be a partition of $n$. Let $T$
be an $n$-tableau of shape $Y\left(  \lambda\right)  $. Then, $\mathbf{E}%
_{\lambda}\mathbf{E}_{T}=\left(  h^{\lambda}\right)  ^{2}\mathbf{E}_{T}$.
\end{proposition}

\begin{proof}
[Proof of Proposition \ref{prop.spechtmod.Elam.ETEl}.]We have $\mathbf{E}%
_{\lambda}=\sum_{S\text{ is an }n\text{-tableau of shape }\lambda}%
\mathbf{E}_{S}$ (by (\ref{eq.def.spechtmod.Elam.Elam.def}), with the index $T$
renamed as $S$). Thus,%
\begin{align*}
\mathbf{E}_{\lambda}\mathbf{E}_{T}  &  =\left(  \sum_{S\text{ is an
}n\text{-tableau of shape }\lambda}\mathbf{E}_{S}\right)  \mathbf{E}_{T}%
=\sum_{S\text{ is an }n\text{-tableau of shape }\lambda}\mathbf{E}%
_{S}\mathbf{E}_{T}\\
&  =\sum_{\substack{S\text{ is an }n\text{-tableau of shape }\lambda
;\\\text{there exist }\widetilde{r}\in\mathcal{R}\left(  T\right)  \text{ and
}\widetilde{c}\in\mathcal{C}\left(  T\right)  \\\text{such that }%
S=\widetilde{c}\widetilde{r}T}}\mathbf{E}_{S}\mathbf{E}_{T}+\sum
_{\substack{S\text{ is an }n\text{-tableau of shape }\lambda;\\\text{there
exist no }\widetilde{r}\in\mathcal{R}\left(  T\right)  \text{ and
}\widetilde{c}\in\mathcal{C}\left(  T\right)  \\\text{such that }%
S=\widetilde{c}\widetilde{r}T}}\mathbf{E}_{S}\mathbf{E}_{T}%
\end{align*}
(here, we have split the sum according to the existence or non-existence of
$\widetilde{r}\in\mathcal{R}\left(  T\right)  $ and $\widetilde{c}%
\in\mathcal{C}\left(  T\right)  $ such that $S=\widetilde{c}\widetilde{r}T$).
Now, let us consider the two sums on the right hand side here. The second sum
is easy to simplify:%
\begin{align}
&  \sum_{\substack{S\text{ is an }n\text{-tableau of shape }\lambda
;\\\text{there exist no }\widetilde{r}\in\mathcal{R}\left(  T\right)  \text{
and }\widetilde{c}\in\mathcal{C}\left(  T\right)  \\\text{such that
}S=\widetilde{c}\widetilde{r}T}}\ \ \underbrace{\mathbf{E}_{S}\mathbf{E}_{T}%
}_{\substack{=0\\\text{(by Proposition \ref{prop.specht.ET.ESET2}
\textbf{(a)})}}}\nonumber\\
&  =\sum_{\substack{S\text{ is an }n\text{-tableau of shape }\lambda
;\\\text{there exist no }\widetilde{r}\in\mathcal{R}\left(  T\right)  \text{
and }\widetilde{c}\in\mathcal{C}\left(  T\right)  \\\text{such that
}S=\widetilde{c}\widetilde{r}T}}0=0. \label{pf.prop.spechtmod.Elam.ETEl.1}%
\end{align}

It remains to simplify the first sum. Here it helps to prove the following
simple claim:

\begin{statement}
\textit{Claim 1:} Let $S$ be an $n$-tableau of shape $\lambda$. If there
exists a pair $\left(  \widetilde{r},\widetilde{c}\right)  \in\mathcal{R}%
\left(  T\right)  \times\mathcal{C}\left(  T\right)  $ such that
$S=\widetilde{c}\widetilde{r}T$, then this pair $\left(  \widetilde{r}%
,\widetilde{c}\right)  $ is unique.
\end{statement}

\begin{proof}
[Proof of Claim 1.]Assume that $\left(  \widetilde{r}_{1},\widetilde{c}%
_{1}\right)  $ and $\left(  \widetilde{r}_{2},\widetilde{c}_{2}\right)  $ are
two such pairs. We must show that these two pairs are equal -- i.e., that
$\left(  \widetilde{r}_{1},\widetilde{c}_{1}\right)  =\left(  \widetilde{r}%
_{2},\widetilde{c}_{2}\right)  $.

We have assumed that $\left(  \widetilde{r}_{1},\widetilde{c}_{1}\right)  $ is
a pair $\left(  \widetilde{r},\widetilde{c}\right)  \in\mathcal{R}\left(
T\right)  \times\mathcal{C}\left(  T\right)  $ such that $S=\widetilde{c}%
\widetilde{r}T$. In other words, $\left(  \widetilde{r}_{1},\widetilde{c}%
_{1}\right)  \in\mathcal{R}\left(  T\right)  \times\mathcal{C}\left(
T\right)  $ and $S=\widetilde{c}_{1}\widetilde{r}_{1}T$. Similarly, $\left(
\widetilde{r}_{2},\widetilde{c}_{2}\right)  \in\mathcal{R}\left(  T\right)
\times\mathcal{C}\left(  T\right)  $ and $S=\widetilde{c}_{2}\widetilde{r}%
_{2}T$. Comparing $S=\widetilde{c}_{1}\widetilde{r}_{1}T$ with
$S=\widetilde{c}_{2}\widetilde{r}_{2}T$, we obtain $\widetilde{c}%
_{1}\widetilde{r}_{1}T=\widetilde{c}_{2}\widetilde{r}_{2}T$.

Moreover, from $\left(  \widetilde{r}_{1},\widetilde{c}_{1}\right)
\in\mathcal{R}\left(  T\right)  \times\mathcal{C}\left(  T\right)  $, we
obtain $\widetilde{r}_{1}\in\mathcal{R}\left(  T\right)  $ and $\widetilde{c}%
_{1}\in\mathcal{C}\left(  T\right)  $. Similarly, $\widetilde{r}_{2}%
\in\mathcal{R}\left(  T\right)  $ and $\widetilde{c}_{2}\in\mathcal{C}\left(
T\right)  $. From $\widetilde{c}_{1}\in\mathcal{C}\left(  T\right)  $ and
$\widetilde{c}_{2}\in\mathcal{C}\left(  T\right)  $, we obtain $\widetilde{c}%
_{1}^{-1}\widetilde{c}_{2}\in\mathcal{C}\left(  T\right)  $ (since
$\mathcal{C}\left(  T\right)  $ is a group). Similarly, $\widetilde{r}%
_{1}\widetilde{r}_{2}^{-1}\in\mathcal{R}\left(  T\right)  $.

However, it is easy to see that any two permutations $u,v\in S_{n}$ satisfying
$uT=vT$ must necessarily be equal\footnote{This was shown during the proof of
Theorem \ref{thm.specht.sandw1} \textbf{(b)}.}. Applying this to
$u=\widetilde{c}_{1}\widetilde{r}_{1}$ and $v=\widetilde{c}_{2}\widetilde{r}%
_{2}$, we see that the two permutations $\widetilde{c}_{1}\widetilde{r}_{1}$
and $\widetilde{c}_{2}\widetilde{r}_{2}$ must be equal (since $\widetilde{c}%
_{1}\widetilde{r}_{1}T=\widetilde{c}_{2}\widetilde{r}_{2}T$). In other words,
$\widetilde{c}_{1}\widetilde{r}_{1}=\widetilde{c}_{2}\widetilde{r}_{2}$.
Hence, $\widetilde{r}_{1}=\widetilde{c}_{1}^{-1}\widetilde{c}_{2}%
\widetilde{r}_{2}$, so that $\widetilde{r}_{1}\widetilde{r}_{2}^{-1}%
=\widetilde{c}_{1}^{-1}\widetilde{c}_{2}\in\mathcal{C}\left(  T\right)  $.
Combining this with $\widetilde{r}_{1}\widetilde{r}_{2}^{-1}\in\mathcal{R}%
\left(  T\right)  $, we obtain $\widetilde{r}_{1}\widetilde{r}_{2}^{-1}%
\in\mathcal{R}\left(  T\right)  \cap\mathcal{C}\left(  T\right)  =\left\{
\operatorname*{id}\right\}  $ (by Proposition \ref{prop.tableau.RnC}
\textbf{(a)}). In other words, $\widetilde{r}_{1}\widetilde{r}_{2}%
^{-1}=\operatorname*{id}$. Hence, $\widetilde{r}_{1}=\widetilde{r}_{2}$.
Furthermore, comparing $\widetilde{r}_{1}\widetilde{r}_{2}^{-1}%
=\operatorname*{id}$ with $\widetilde{r}_{1}\widetilde{r}_{2}^{-1}%
=\widetilde{c}_{1}^{-1}\widetilde{c}_{2}$, we find $\widetilde{c}_{1}%
^{-1}\widetilde{c}_{2}=\operatorname*{id}$, and thus $\widetilde{c}%
_{1}=\widetilde{c}_{2}$. Combining $\widetilde{r}_{1}=\widetilde{r}_{2}$ with
$\widetilde{c}_{1}=\widetilde{c}_{2}$, we obtain $\left(  \widetilde{r}%
_{1},\widetilde{c}_{1}\right)  =\left(  \widetilde{r}_{2},\widetilde{c}%
_{2}\right)  $. This proves Claim 1.
\end{proof}

Now, the map%
\begin{align*}
&  \mathcal{R}\left(  T\right)  \times\mathcal{C}\left(  T\right) \\
&  \rightarrow\left\{  n\text{-tableaux }S\text{ of shape }\lambda
\text{\ }\mid\ \text{there exist }\widetilde{r}\in\mathcal{R}\left(  T\right)
\text{ and }\widetilde{c}\in\mathcal{C}\left(  T\right)  \right. \\
&
\ \ \ \ \ \ \ \ \ \ \ \ \ \ \ \ \ \ \ \ \ \ \ \ \ \ \ \ \ \ \ \ \ \ \ \ \ \ \ \ \ \ \ \ \ \ \ \ \ \ \left.
\text{such that }S=\widetilde{c}\widetilde{r}T\right\}
\end{align*}
that sends each pair $\left(  \widetilde{r},\widetilde{c}\right)
\in\mathcal{R}\left(  T\right)  \times\mathcal{C}\left(  T\right)  $ to the
$n$-tableau $\widetilde{c}\widetilde{r}T$ is clearly well-defined and
surjective (indeed, its target has been handpicked to contain precisely all
its values), and furthermore (by Claim 1) is injective. Hence, this map is a
bijection. Thus, we can substitute $\widetilde{c}\widetilde{r}T$ for $S$ in
the sum
\[
\sum_{\substack{S\text{ is an }n\text{-tableau of shape }\lambda;\\\text{there
exist }\widetilde{r}\in\mathcal{R}\left(  T\right)  \text{ and }%
\widetilde{c}\in\mathcal{C}\left(  T\right)  \\\text{such that }%
S=\widetilde{c}\widetilde{r}T}}\mathbf{E}_{S}\mathbf{E}_{T}.
\]
We thus obtain%
\begin{align*}
\sum_{\substack{S\text{ is an }n\text{-tableau of shape }\lambda;\\\text{there
exist }\widetilde{r}\in\mathcal{R}\left(  T\right)  \text{ and }%
\widetilde{c}\in\mathcal{C}\left(  T\right)  \\\text{such that }%
S=\widetilde{c}\widetilde{r}T}}\mathbf{E}_{S}\mathbf{E}_{T}  &
=\underbrace{\sum_{\left(  \widetilde{r},\widetilde{c}\right)  \in
\mathcal{R}\left(  T\right)  \times\mathcal{C}\left(  T\right)  }}%
_{=\sum_{\widetilde{r}\in\mathcal{R}\left(  T\right)  }\ \ \sum_{\widetilde{c}%
\in\mathcal{C}\left(  T\right)  }}\mathbf{E}_{\widetilde{c}\widetilde{r}%
T}\mathbf{E}_{T}\\
&  =\sum_{\widetilde{r}\in\mathcal{R}\left(  T\right)  }\ \ \sum
_{\widetilde{c}\in\mathcal{C}\left(  T\right)  }\underbrace{\mathbf{E}%
_{\widetilde{c}\widetilde{r}T}\mathbf{E}_{T}}_{\substack{=h^{\lambda}\left(
-1\right)  ^{\widetilde{c}}\widetilde{c}\widetilde{r}\mathbf{E}_{T}\\\text{(by
Proposition \ref{prop.specht.ET.ESET2} \textbf{(b)},}\\\text{applied to
}S=\widetilde{c}\widetilde{r}T\text{)}}}\\
&  =\sum_{\widetilde{r}\in\mathcal{R}\left(  T\right)  }\ \ \sum
_{\widetilde{c}\in\mathcal{C}\left(  T\right)  }h^{\lambda}\left(  -1\right)
^{\widetilde{c}}\widetilde{c}\widetilde{r}\mathbf{E}_{T}\\
&  =h^{\lambda}\underbrace{\left(  \sum_{\widetilde{c}\in\mathcal{C}\left(
T\right)  }\left(  -1\right)  ^{\widetilde{c}}\widetilde{c}\right)
}_{\substack{=\sum_{w\in\mathcal{C}\left(  T\right)  }\left(  -1\right)
^{w}w\\=\nabla_{\operatorname*{Col}T}^{-}\\\text{(by the definition}\\\text{of
}\nabla_{\operatorname*{Col}T}^{-}\text{)}}}\ \ \underbrace{\left(
\sum_{\widetilde{r}\in\mathcal{R}\left(  T\right)  }\widetilde{r}\right)
}_{\substack{=\sum_{w\in\mathcal{R}\left(  T\right)  }w\\=\nabla
_{\operatorname*{Row}T}\\\text{(by the definition}\\\text{of }\nabla
_{\operatorname*{Row}T}\text{)}}}\mathbf{E}_{T}\\
&  =h^{\lambda}\underbrace{\nabla_{\operatorname*{Col}T}^{-}\nabla
_{\operatorname*{Row}T}}_{\substack{=\mathbf{E}_{T}\\\text{(by the definition
of }\mathbf{E}_{T}\text{)}}}\mathbf{E}_{T}=h^{\lambda}\underbrace{\mathbf{E}%
_{T}\mathbf{E}_{T}}_{=\mathbf{E}_{T}^{2}}=h^{\lambda}\mathbf{E}_{T}^{2}.
\end{align*}
Moreover, Theorem \ref{thm.specht.ETidp} yields $\mathbf{E}_{T}^{2}=\dfrac
{n!}{f^{\lambda}}\mathbf{E}_{T}=h^{\lambda}\mathbf{E}_{T}$ (since $\dfrac
{n!}{f^{\lambda}}=h^{\lambda}$). Thus,%
\begin{align}
\sum_{\substack{S\text{ is an }n\text{-tableau of shape }\lambda;\\\text{there
exist }\widetilde{r}\in\mathcal{R}\left(  T\right)  \text{ and }%
\widetilde{c}\in\mathcal{C}\left(  T\right)  \\\text{such that }%
S=\widetilde{c}\widetilde{r}T}}\mathbf{E}_{S}\mathbf{E}_{T}  &  =h^{\lambda
}\underbrace{\mathbf{E}_{T}^{2}}_{=h^{\lambda}\mathbf{E}_{T}}=h^{\lambda
}h^{\lambda}\mathbf{E}_{T}\nonumber\\
&  =\left(  h^{\lambda}\right)  ^{2}\mathbf{E}_{T}.
\label{pf.prop.spechtmod.Elam.ETEl.3}%
\end{align}

Now, we can finish the calculation that we started at the beginning of this
proof:%
\begin{align*}
\mathbf{E}_{\lambda}\mathbf{E}_{T}  &  =\underbrace{\sum_{\substack{S\text{ is
an }n\text{-tableau of shape }\lambda;\\\text{there exist }\widetilde{r}%
\in\mathcal{R}\left(  T\right)  \text{ and }\widetilde{c}\in\mathcal{C}\left(
T\right)  \\\text{such that }S=\widetilde{c}\widetilde{r}T}}\mathbf{E}%
_{S}\mathbf{E}_{T}}_{\substack{=\left(  h^{\lambda}\right)  ^{2}\mathbf{E}%
_{T}\\\text{(by (\ref{pf.prop.spechtmod.Elam.ETEl.3}))}}}+\underbrace{\sum
_{\substack{S\text{ is an }n\text{-tableau of shape }\lambda;\\\text{there
exist no }\widetilde{r}\in\mathcal{R}\left(  T\right)  \text{ and
}\widetilde{c}\in\mathcal{C}\left(  T\right)  \\\text{such that }%
S=\widetilde{c}\widetilde{r}T}}\mathbf{E}_{S}\mathbf{E}_{T}}%
_{\substack{=0\\\text{(by (\ref{pf.prop.spechtmod.Elam.ETEl.1}))}}}\\
&  =\left(  h^{\lambda}\right)  ^{2}\mathbf{E}_{T}+0=\left(  h^{\lambda
}\right)  ^{2}\mathbf{E}_{T}.
\end{align*}
This proves Proposition \ref{prop.spechtmod.Elam.ETEl}.
\end{proof}

Note that the product $\mathbf{E}_{\lambda}\mathbf{E}_{T}$ in Proposition
\ref{prop.spechtmod.Elam.ETEl} can also be rewritten as $\mathbf{E}%
_{T}\mathbf{E}_{\lambda}$ (since Proposition \ref{prop.spechtmod.Elam.incent}
shows that $\mathbf{E}_{\lambda}\in Z\left(  \mathbf{k}\left[  S_{n}\right]
\right)  $ and thus $\mathbf{E}_{\lambda}\mathbf{E}_{T}=\mathbf{E}%
_{T}\mathbf{E}_{\lambda}$).

\begin{proof}
[Proof of Theorem \ref{thm.spechtmod.Elam.idp-strong}.]We have%
\begin{align*}
\mathbf{E}_{\lambda}^{2}  &  =\mathbf{E}_{\lambda}\mathbf{E}_{\lambda
}=\mathbf{E}_{\lambda}\sum_{T\text{ is an }n\text{-tableau of shape }\lambda
}\mathbf{E}_{T}\ \ \ \ \ \ \ \ \ \ \left(  \text{by
(\ref{eq.def.spechtmod.Elam.Elam.def})}\right) \\
&  =\sum_{T\text{ is an }n\text{-tableau of shape }\lambda}%
\underbrace{\mathbf{E}_{\lambda}\mathbf{E}_{T}}_{\substack{=\left(
h^{\lambda}\right)  ^{2}\mathbf{E}_{T}\\\text{(by Proposition
\ref{prop.spechtmod.Elam.ETEl})}}}\\
&  =\sum_{T\text{ is an }n\text{-tableau of shape }\lambda}\left(  h^{\lambda
}\right)  ^{2}\mathbf{E}_{T}=\left(  h^{\lambda}\right)  ^{2}\underbrace{\sum
_{T\text{ is an }n\text{-tableau of shape }\lambda}\mathbf{E}_{T}%
}_{\substack{=\mathbf{E}_{\lambda}\\\text{(by
(\ref{eq.def.spechtmod.Elam.Elam.def}))}}}\\
&  =\left(  h^{\lambda}\right)  ^{2}\mathbf{E}_{\lambda}.
\end{align*}
This proves Theorem \ref{thm.spechtmod.Elam.idp-strong}.
\end{proof}

Theorem \ref{thm.spechtmod.Elam.idp-strong} has an easy but surprising
consequence. As we recall, $\mathbf{k}$-algebras are also known as
\emph{unital }$\mathbf{k}$\emph{-algebras}, since they have unities. In
contrast, a \emph{nonunital }$\mathbf{k}$\emph{-algebra} is defined just like
a $\mathbf{k}$-algebra, except that it is not required to have a unity.
Obviously, each $\mathbf{k}$-algebra becomes a nonunital $\mathbf{k}$-algebra
if we forget its unity. A \emph{nonunital }$\mathbf{k}$\emph{-subalgebra} of a
nonunital $\mathbf{k}$-algebra $A$ is defined as a $\mathbf{k}$-submodule $B$
of $A$ that is closed under multiplication. Somewhat surprisingly, it can
happen that $A$ and $B$ are two unital $\mathbf{k}$-algebras, and $B$ is a
nonunital $\mathbf{k}$-subalgebra of $A$ but not a unital $\mathbf{k}%
$-subalgebra of $A$ (despite having a unity!). This is the case whenever the
unities of $A$ and $B$ are distinct, and this does indeed happen quite
often\footnote{For example: If $R$ and $S$ are two nontrivial $\mathbf{k}%
$-algebras, then $R\times\left\{  0_{S}\right\}  $ is a nonunital $\mathbf{k}%
$-subalgebra of $R\times S$, but is itself a unital $\mathbf{k}$-algebra (with
unity $\left(  1_{R},0_{S}\right)  $), despite not being a unital $\mathbf{k}%
$-subalgebra of $R\times S$ (since the unity of $R\times S$ is $\left(
1_{R},1_{S}\right)  $).}. The following corollary gives another example of
this paradoxical phenomenon:

\begin{corollary}
\label{cor.spechtmod.Elam.subalg}Recall that $\mathcal{A}=\mathbf{k}\left[
S_{n}\right]  $. Let $\lambda$ be a partition of $n$. \medskip

\textbf{(a)} The subset $\mathcal{A}\mathbf{E}_{\lambda}$ of $\mathcal{A}$ is
a $\mathbf{k}$-submodule and is closed under multiplication. That is, it is a
nonunital subalgebra of $\mathcal{A}$. \medskip

\textbf{(b)} Now assume that $h^{\lambda}$ is invertible in $\mathbf{k}$.
Then, the nonunital algebra $\mathcal{A}\mathbf{E}_{\lambda}$ from Corollary
\ref{cor.spechtmod.Elam.subalg} \textbf{(a)} has a unity, namely $\dfrac
{1}{\left(  h^{\lambda}\right)  ^{2}}\mathbf{E}_{\lambda}$. (But this does not
make $\mathcal{A}\mathbf{E}_{\lambda}$ a unital subalgebra of $\mathcal{A}$.)
\end{corollary}

\begin{proof}
\textbf{(a)} This is obvious. In fact, if $R$ is any $\mathbf{k}$-algebra, and
if $r\in R$ is any element, then the subset $Rr$ of $R$ is a $\mathbf{k}%
$-submodule and is closed under multiplication (since the product of two
elements of $Rr$ has the form $xr\cdot yr=\underbrace{\left(  xry\right)
}_{\in R}r$ and thus belongs to $Rr$ again), and thus is a nonunital
subalgebra of $R$. Applying this fact to $R=\mathcal{A}$ and $r=\mathbf{E}%
_{\lambda}$, we obtain the claim of Corollary \ref{cor.spechtmod.Elam.subalg}
\textbf{(a)}. \medskip

\textbf{(b)} Proposition \ref{prop.spechtmod.Elam.incent} yields
$\mathbf{E}_{\lambda}\in Z\left(  \mathbf{k}\left[  S_{n}\right]  \right)  $.
In other words, $\mathbf{E}_{\lambda}\in Z\left(  \mathcal{A}\right)  $ (since
$\mathcal{A}=\mathbf{k}\left[  S_{n}\right]  $).

Let $\mathbf{r}=\dfrac{1}{\left(  h^{\lambda}\right)  ^{2}}\mathbf{E}%
_{\lambda}$. Thus, $\mathbf{r}=\dfrac{1}{\left(  h^{\lambda}\right)  ^{2}%
}\mathbf{E}_{\lambda}=\underbrace{\dfrac{1}{\left(  h^{\lambda}\right)  ^{2}%
}\cdot1}_{\in\mathcal{A}}\mathbf{E}_{\lambda}\in\mathcal{A}\mathbf{E}%
_{\lambda}$ and%
\[
\mathbf{E}_{\lambda}\mathbf{r}=\mathbf{E}_{\lambda}\cdot\dfrac{1}{\left(
h^{\lambda}\right)  ^{2}}\mathbf{E}_{\lambda}=\dfrac{1}{\left(  h^{\lambda
}\right)  ^{2}}\underbrace{\mathbf{E}_{\lambda}^{2}}_{\substack{=\left(
h^{\lambda}\right)  ^{2}\mathbf{E}_{\lambda}\\\text{(by Theorem
\ref{thm.spechtmod.Elam.idp-strong})}}}=\dfrac{1}{\left(  h^{\lambda}\right)
^{2}}\left(  h^{\lambda}\right)  ^{2}\mathbf{E}_{\lambda}=\mathbf{E}_{\lambda
}.
\]

Now, we claim that each $\mathbf{a}\in\mathcal{A}\mathbf{E}_{\lambda}$
satisfies $\mathbf{ar}=\mathbf{ra}=\mathbf{a}$.

[\textit{Proof:} Let $\mathbf{a}\in\mathcal{A}\mathbf{E}_{\lambda}$. Thus,
$\mathbf{a}=\mathbf{bE}_{\lambda}$ for some $\mathbf{b}\in\mathcal{A}$.
Consider this $\mathbf{b}$. Now, from $\mathbf{a}=\mathbf{bE}_{\lambda}$, we
obtain%
\[
\mathbf{ar}=\mathbf{b}\underbrace{\mathbf{E}_{\lambda}\mathbf{r}}%
_{=\mathbf{E}_{\lambda}}=\mathbf{bE}_{\lambda}=\mathbf{a}.
\]
Furthermore, from $\mathbf{r}=\dfrac{1}{\left(  h^{\lambda}\right)  ^{2}%
}\mathbf{E}_{\lambda}$, we obtain%
\[
\mathbf{ra}=\dfrac{1}{\left(  h^{\lambda}\right)  ^{2}}\underbrace{\mathbf{E}%
_{\lambda}\mathbf{a}}_{\substack{=\mathbf{aE}_{\lambda}\\\text{(since
}\mathbf{E}_{\lambda}\in Z\left(  \mathcal{A}\right)  \text{)}}}=\dfrac
{1}{\left(  h^{\lambda}\right)  ^{2}}\mathbf{aE}_{\lambda}=\mathbf{a}%
\cdot\underbrace{\dfrac{1}{\left(  h^{\lambda}\right)  ^{2}}\mathbf{E}%
_{\lambda}}_{=\mathbf{r}}=\mathbf{ar},
\]
so that $\mathbf{ar}=\mathbf{ra}=\mathbf{a}$ (since $\mathbf{ra}%
=\mathbf{ar}=\mathbf{a}$). Qed.]

We thus have shown that $\mathbf{r}\in\mathcal{A}\mathbf{E}_{\lambda}$ and
that each $\mathbf{a}\in\mathcal{A}\mathbf{E}_{\lambda}$ satisfies
$\mathbf{ar}=\mathbf{ra}=\mathbf{a}$. In other words, $\mathbf{r}$ is a unity
of the nonunital algebra $\mathcal{A}\mathbf{E}_{\lambda}$. In other words,
$\dfrac{1}{\left(  h^{\lambda}\right)  ^{2}}\mathbf{E}_{\lambda}$ is a unity
of the nonunital algebra $\mathcal{A}\mathbf{E}_{\lambda}$ (since
$\mathbf{r}=\dfrac{1}{\left(  h^{\lambda}\right)  ^{2}}\mathbf{E}_{\lambda}$).
This proves Corollary \ref{cor.spechtmod.Elam.subalg} \textbf{(b)}.
\end{proof}

\subsubsection{A partition of unity into centralized Young symmetrizers}

One of our next goals is to make Theorem \ref{thm.spechtmod.Elam.center}
\textbf{(a)} more concrete by providing explicit expansions for elements of
$Z\left(  \mathbf{k}\left[  S_{n}\right]  \right)  $ as $\mathbf{k}$-linear
combinations of the $\mathbf{E}_{\lambda}$. We begin with the element $1$:

\begin{theorem}
\label{thm.spechtmod.Elam.pou}Let $\mathbf{k}$ be a commutative $\mathbb{Q}%
$-algebra. Then, in $\mathbf{k}\left[  S_{n}\right]  $, we have%
\[
1=\sum_{\lambda\text{ is a partition of }n}\dfrac{\mathbf{E}_{\lambda}%
}{\left(  h^{\lambda}\right)  ^{2}}.
\]

\end{theorem}

\begin{proof}
As in the proof of Proposition \ref{prop.spechtmod.Elam.idp-weak}, we can WLOG
assume that $\mathbf{k}=\mathbb{Q}$ (since base change allows us to extend the
result from $\mathbf{k}=\mathbb{Q}$ to any commutative $\mathbb{Q}$-algebra
$\mathbf{k}$). Assume this. Thus, $\mathbf{k}$ is a field.

As in the proof of Proposition \ref{prop.spechtmod.Elam.idp-weak}, we can
write $1$ in the form%
\[
1=\sum_{\mu\text{ is a partition of }n}\alpha_{\mu}\mathbf{E}_{\mu
}\ \ \ \ \ \ \ \ \ \ \text{for some scalars }\alpha_{\mu}\in\mathbf{k}.
\]
Consider these scalars $\alpha_{\mu}$.

Now, let $\lambda$ be a partition of $n$. As in the proof of Proposition
\ref{prop.spechtmod.Elam.idp-weak}, we can show that $\mathbf{E}_{\lambda}%
\neq0$ and
\begin{equation}
\alpha_{\lambda}\mathbf{E}_{\lambda}^{2}=\mathbf{E}_{\lambda}.
\label{pf.thm.spechtmod.Elam.pou.4}%
\end{equation}
But Theorem \ref{thm.spechtmod.Elam.idp-strong} shows that $\mathbf{E}%
_{\lambda}^{2}=\left(  h^{\lambda}\right)  ^{2}\mathbf{E}_{\lambda}$. Thus, we
can rewrite the equality (\ref{pf.thm.spechtmod.Elam.pou.4}) as
\[
\alpha_{\lambda}\left(  h^{\lambda}\right)  ^{2}\mathbf{E}_{\lambda
}=\mathbf{E}_{\lambda}.
\]

Thus, $\left(  \alpha_{\lambda}\left(  h^{\lambda}\right)  ^{2}-1\right)
\mathbf{E}_{\lambda}=\underbrace{\alpha_{\lambda}\left(  h^{\lambda}\right)
^{2}\mathbf{E}_{\lambda}}_{=\mathbf{E}_{\lambda}}-\,\mathbf{E}_{\lambda
}=\mathbf{E}_{\lambda}-\mathbf{E}_{\lambda}=0$.

But $\mathbf{k}$ is a field. Thus, if the scalar $\alpha_{\lambda}\left(
h^{\lambda}\right)  ^{2}-1$ was nonzero, then this scalar would be invertible,
and thus we could divide the equality $\left(  \alpha_{\lambda}\left(
h^{\lambda}\right)  ^{2}-1\right)  \mathbf{E}_{\lambda}=0$ by $\alpha
_{\lambda}\left(  h^{\lambda}\right)  ^{2}-1$, and obtain $\mathbf{E}%
_{\lambda}=0$, which would contradict $\mathbf{E}_{\lambda}\neq0$. Hence, the
scalar $\alpha_{\lambda}\left(  h^{\lambda}\right)  ^{2}-1$ cannot be nonzero.
In other words, $\alpha_{\lambda}\left(  h^{\lambda}\right)  ^{2}-1=0$.
Equivalently, $\alpha_{\lambda}\left(  h^{\lambda}\right)  ^{2}=1$. Thus,
$\alpha_{\lambda}=\dfrac{1}{\left(  h^{\lambda}\right)  ^{2}}$.

Forget that we fixed $\lambda$. We thus have shown that each partition
$\lambda$ of $n$ satisfies%
\begin{equation}
\alpha_{\lambda}=\dfrac{1}{\left(  h^{\lambda}\right)  ^{2}}.
\label{pf.thm.spechtmod.Elam.pou.7}%
\end{equation}

Now, recall that%
\begin{align*}
1  &  =\sum_{\mu\text{ is a partition of }n}\alpha_{\mu}\mathbf{E}_{\mu}%
=\sum_{\lambda\text{ is a partition of }n}\underbrace{\alpha_{\lambda}%
}_{\substack{=\dfrac{1}{\left(  h^{\lambda}\right)  ^{2}}\\\text{(by
(\ref{pf.thm.spechtmod.Elam.pou.7}))}}}\mathbf{E}_{\lambda}\\
&  =\sum_{\lambda\text{ is a partition of }n}\dfrac{1}{\left(  h^{\lambda
}\right)  ^{2}}\mathbf{E}_{\lambda}=\sum_{\lambda\text{ is a partition of }%
n}\dfrac{\mathbf{E}_{\lambda}}{\left(  h^{\lambda}\right)  ^{2}}.
\end{align*}
This proves Theorem \ref{thm.spechtmod.Elam.pou}.
\end{proof}

\subsubsection{A combinatorial consequence: the sum of all $\left(
f^{\lambda}\right)  ^{2}$}

Theorem \ref{thm.spechtmod.Elam.pou} leads to an important purely
combinatorial identity:

\begin{corollary}
\label{cor.spechtmod.sumflam2}We have%
\[
\sum_{\lambda\text{ is a partition of }n}\left(  f^{\lambda}\right)  ^{2}=n!.
\]

\end{corollary}

\begin{example}
Let $n=4$. Then, the partitions of $n$ are $\left(  4\right)  $, $\left(
3,1\right)  $, $\left(  2,2\right)  $, $\left(  2,1,1\right)  $ and $\left(
1,1,1,1\right)  $, and the corresponding numbers $f^{\lambda}$ are%
\[
f^{\left(  4\right)  }=1,\ \ \ \ \ \ \ \ \ \ f^{\left(  3,1\right)
}=3,\ \ \ \ \ \ \ \ \ \ f^{\left(  2,2\right)  }%
=2,\ \ \ \ \ \ \ \ \ \ f^{\left(  2,1,1\right)  }%
=3,\ \ \ \ \ \ \ \ \ \ f^{\left(  1,1,1,1\right)  }=1.
\]
Hence, Corollary \ref{cor.spechtmod.sumflam2} says that $1^{2}+3^{2}%
+2^{2}+3^{2}+1^{2}=4!$. And this is indeed true.
\end{example}

\begin{proof}
[Proof of Corollary \ref{cor.spechtmod.sumflam2}.]Let $\mathbf{k}=\mathbb{Q}$.
Then, Theorem \ref{thm.spechtmod.Elam.pou} yields%
\begin{equation}
1=\sum_{\lambda\text{ is a partition of }n}\dfrac{\mathbf{E}_{\lambda}%
}{\left(  h^{\lambda}\right)  ^{2}}. \label{pf.cor.spechtmod.sumflam2.1}%
\end{equation}
Now, $\left[  1\right]  1=1$, so that%
\begin{align}
1  &  =\left[  1\right]  1=\left[  1\right]  \sum_{\lambda\text{ is a
partition of }n}\dfrac{\mathbf{E}_{\lambda}}{\left(  h^{\lambda}\right)  ^{2}%
}\ \ \ \ \ \ \ \ \ \ \left(  \text{by (\ref{pf.cor.spechtmod.sumflam2.1}%
)}\right) \nonumber\\
&  =\sum_{\lambda\text{ is a partition of }n}\dfrac{\left[  1\right]
\mathbf{E}_{\lambda}}{\left(  h^{\lambda}\right)  ^{2}}.
\label{pf.cor.spechtmod.sumflam2.2}%
\end{align}

However, each partition $\lambda$ of $n$ satisfies $\left[  1\right]
\mathbf{E}_{\lambda}=n!$ (by Lemma \ref{lem.specht.Elam.1-coord}) and thus%
\begin{align*}
\dfrac{\left[  1\right]  \mathbf{E}_{\lambda}}{\left(  h^{\lambda}\right)
^{2}}  &  =\dfrac{n!}{\left(  h^{\lambda}\right)  ^{2}}=\dfrac{n!}{\left(
n!/f^{\lambda}\right)  ^{2}}\ \ \ \ \ \ \ \ \ \ \left(  \text{since
}h^{\lambda}=\dfrac{n!}{f^{\lambda}}=n!/f^{\lambda}\right) \\
&  =\dfrac{\left(  f^{\lambda}\right)  ^{2}}{n!}.
\end{align*}
Thus, we can rewrite (\ref{pf.cor.spechtmod.sumflam2.2}) as%
\[
1=\sum_{\lambda\text{ is a partition of }n}\dfrac{\left(  f^{\lambda}\right)
^{2}}{n!}.
\]
Multiplying both sides of this equality by $n!$, we find%
\[
n!=\sum_{\lambda\text{ is a partition of }n}\left(  f^{\lambda}\right)  ^{2}.
\]
This proves Corollary \ref{cor.spechtmod.sumflam2}.
\end{proof}

Corollary \ref{cor.spechtmod.sumflam2} is a crucial result in the
combinatorics of permutations, and can be proved without any algebra. Such
proofs can be found in \cite[Proposition 1.3.3]{vL-RSK}, \cite[page 26, first
formula]{Ruther48}, \cite[Corollary 8.5]{Stanley-AC}, \cite[Exercise
5.7]{MenRem15}, \cite[\S 4.3, (4)]{Fulton97}, \cite[Theorem 8.26 a.]%
{Aigner07}, \cite[Theorem 2.6.5 part 3.]{Sagan01}, \cite[(3.12)]{Prasad-rep},
\cite[\S 5.1.4, Theorem A]{Knuth-TAoCP3} and many other sources. (See also
\cite[Corollary 7.12.6]{Stanley-EC2} for a proof using symmetric functions.)

\subsubsection{Explicit expansions in centralized Young symmetrizers}

In Theorem \ref{thm.spechtmod.Elam.center} \textbf{(a)}, we required that
$\mathbf{k}$ be a field of characteristic $0$. However, this requirement is
stronger than necessary; it suffices to assume that $n!$ be invertible in
$\mathbf{k}$. In order to show this, we need some more lemmas.

We try to expand an arbitrary element $\mathbf{a}$ of $Z\left(  \mathbf{k}%
\left[  S_{n}\right]  \right)  $ as a $\mathbf{k}$-linear combination of the
$\mathbf{E}_{\lambda}$ (when $n!$ is invertible in $\mathbf{k}$). We first do
this in the case when $\mathbf{k}$ is a field of characteristic $0$:

\begin{lemma}
\label{lem.spechtmod.Elam.a-field}Let $\mathbf{k}$ be a field of
characteristic $0$. Let $\mathbf{a}\in Z\left(  \mathbf{k}\left[
S_{n}\right]  \right)  $. Then,%
\[
\mathbf{a}=\sum_{\lambda\text{ is a partition of }n}\dfrac{\left[  1\right]
\left(  \mathbf{E}_{\lambda}\mathbf{a}\right)  }{\left(  h^{\lambda}\right)
^{2}\cdot n!}\mathbf{E}_{\lambda}.
\]

\end{lemma}

\begin{proof}
Theorem \ref{thm.spechtmod.Elam.basis} tells us that the family $\left(
\mathbf{E}_{\mu}\right)  _{\mu\text{ is a partition of }n}$ is a basis of the
$\mathbf{k}$-vector space $Z\left(  \mathbf{k}\left[  S_{n}\right]  \right)
$. Thus, in particular, this family spans $Z\left(  \mathbf{k}\left[
S_{n}\right]  \right)  $. Hence, the vector $\mathbf{a}$ (which belongs to
$Z\left(  \mathbf{k}\left[  S_{n}\right]  \right)  $) is a linear combination
of this family. In other words, we can write $\mathbf{a}$ in the form%
\[
\mathbf{a}=\sum_{\mu\text{ is a partition of }n}\alpha_{\mu}\mathbf{E}_{\mu
}\ \ \ \ \ \ \ \ \ \ \text{for some scalars }\alpha_{\mu}\in\mathbf{k}.
\]
Consider these scalars $\alpha_{\mu}$.

Let $\lambda$ be a partition of $n$. Multiplying the equality $\mathbf{a}%
=\sum_{\mu\text{ is a partition of }n}\alpha_{\mu}\mathbf{E}_{\mu}$ by
$\mathbf{E}_{\lambda}$ from the left, we get%
\[
\mathbf{E}_{\lambda}\mathbf{a}=\mathbf{E}_{\lambda}\sum_{\mu\text{ is a
partition of }n}\alpha_{\mu}\mathbf{E}_{\mu}=\alpha_{\lambda}\mathbf{E}%
_{\lambda}^{2}%
\]
(by the same computation that we used to prove $\mathbf{E}_{\lambda}\sum
_{\mu\text{ is a partition of }n}\alpha_{\mu}\mathbf{E}_{\mu}=\alpha_{\lambda
}\mathbf{E}_{\lambda}^{2}$ during our above proof of Proposition
\ref{prop.spechtmod.Elam.idp-weak}). However, Theorem
\ref{thm.spechtmod.Elam.idp-strong} yields $\mathbf{E}_{\lambda}^{2}=\left(
h^{\lambda}\right)  ^{2}\mathbf{E}_{\lambda}$. Altogether, we thus find%
\[
\mathbf{E}_{\lambda}\mathbf{a}=\alpha_{\lambda}\underbrace{\mathbf{E}%
_{\lambda}^{2}}_{=\left(  h^{\lambda}\right)  ^{2}\mathbf{E}_{\lambda}}%
=\alpha_{\lambda}\left(  h^{\lambda}\right)  ^{2}\mathbf{E}_{\lambda}.
\]
Hence,%
\[
\left[  1\right]  \left(  \mathbf{E}_{\lambda}\mathbf{a}\right)  =\left[
1\right]  \left(  \alpha_{\lambda}\left(  h^{\lambda}\right)  ^{2}%
\mathbf{E}_{\lambda}\right)  =\alpha_{\lambda}\left(  h^{\lambda}\right)
^{2}\cdot\underbrace{\left[  1\right]  \mathbf{E}_{\lambda}}%
_{\substack{=n!\\\text{(by Lemma \ref{lem.specht.Elam.1-coord})}}%
}=\alpha_{\lambda}\left(  h^{\lambda}\right)  ^{2}\cdot n!.
\]
Solving this equality for $\alpha_{\lambda}$ (since $\left(  h^{\lambda
}\right)  ^{2}\cdot n!$ is a nonzero element of the field $\mathbf{k}$), we
obtain%
\begin{equation}
\alpha_{\lambda}=\dfrac{\left[  1\right]  \left(  \mathbf{E}_{\lambda
}\mathbf{a}\right)  }{\left(  h^{\lambda}\right)  ^{2}\cdot n!}.
\label{pf.lem.spechtmod.Elam.a-field.5}%
\end{equation}

Now, forget that we fixed $\lambda$. We thus have proved
(\ref{pf.lem.spechtmod.Elam.a-field.5}) for each partition $\lambda$ of $n$.

Now, recall that%
\begin{align*}
\mathbf{a}  &  =\sum_{\mu\text{ is a partition of }n}\alpha_{\mu}%
\mathbf{E}_{\mu}=\sum_{\lambda\text{ is a partition of }n}\underbrace{\alpha
_{\lambda}}_{\substack{=\dfrac{\left[  1\right]  \left(  \mathbf{E}_{\lambda
}\mathbf{a}\right)  }{\left(  h^{\lambda}\right)  ^{2}\cdot n!}\\\text{(by
(\ref{pf.lem.spechtmod.Elam.a-field.5}))}}}\mathbf{E}_{\lambda}\\
&  =\sum_{\lambda\text{ is a partition of }n}\dfrac{\left[  1\right]  \left(
\mathbf{E}_{\lambda}\mathbf{a}\right)  }{\left(  h^{\lambda}\right)  ^{2}\cdot
n!}\mathbf{E}_{\lambda}.
\end{align*}
This proves Lemma \ref{lem.spechtmod.Elam.a-field}.
\end{proof}

It stands to reason that Lemma \ref{lem.spechtmod.Elam.a-field} should hold
more generally for any commutative ring $\mathbf{k}$ in which $n!$ is
invertible (since all denominators are easily seen to divide $n!^{3}$). This
is indeed the case. To show this, we first consider conjugacy class sums over
an arbitrary commutative ring $\mathbf{k}$:

\begin{lemma}
\label{lem.spechtmod.Elam.expand1}Let $C$ be a conjugacy class of $S_{n}$. Let
$\mathbf{z}_{C}:=\sum_{c\in C}c\in\mathbf{k}\left[  S_{n}\right]  $. Then,%
\[
n!^{3}\cdot\mathbf{z}_{C}=\sum_{\lambda\text{ is a partition of }n}\left(
f^{\lambda}\right)  ^{2}\cdot\left[  1\right]  \left(  \mathbf{E}_{\lambda
}\mathbf{z}_{C}\right)  \cdot\mathbf{E}_{\lambda}.
\]

\end{lemma}

\begin{proof}
\textit{Step 1:} In order to prove Lemma \ref{lem.spechtmod.Elam.expand1} for
all commutative rings $\mathbf{k}$, it suffices to prove it for $\mathbf{k}%
=\mathbb{Z}$.

Indeed, the same base change argument that we used in the proof of Theorem
\ref{thm.specht.ETidp} can be used here to reduce the general case to the case
of $\mathbf{k}=\mathbb{Z}$. Namely, we consider the canonical ring morphism
$f:\mathbb{Z}\rightarrow\mathbf{k}$ and the base change morphism $f_{\ast
}:\mathbb{Z}\left[  S_{n}\right]  \rightarrow\mathbf{k}\left[  S_{n}\right]  $
corresponding to $f$. Then, $f_{\ast}$ sends $\mathbf{z}_{C}\in\mathbb{Z}%
\left[  S_{n}\right]  $ to $\mathbf{z}_{C}\in\mathbf{k}\left[  S_{n}\right]  $
and sends $\mathbf{E}_{\lambda}\in\mathbb{Z}\left[  S_{n}\right]  $ to
$\mathbf{E}_{\lambda}\in\mathbf{k}\left[  S_{n}\right]  $. Hence, $f_{\ast}$
sends the product $\mathbf{E}_{\lambda}\mathbf{z}_{C}\in\mathbb{Z}\left[
S_{n}\right]  $ to the corresponding product $\mathbf{E}_{\lambda}%
\mathbf{z}_{C}\in\mathbf{k}\left[  S_{n}\right]  $. Thus, the coefficient
$\left[  1\right]  \left(  \mathbf{E}_{\lambda}\mathbf{z}_{C}\right)  $
computed using $\mathbb{Z}\left[  S_{n}\right]  $ is transformed into the
coefficient $\left[  1\right]  \left(  \mathbf{E}_{\lambda}\mathbf{z}%
_{C}\right)  $ computed using $\mathbf{k}\left[  S_{n}\right]  $ by the
morphism $f$. Thus, if the equality%
\[
n!^{3}\cdot\mathbf{z}_{C}=\sum_{\lambda\text{ is a partition of }n}\left(
f^{\lambda}\right)  ^{2}\cdot\left[  1\right]  \left(  \mathbf{E}_{\lambda
}\mathbf{z}_{C}\right)  \cdot\mathbf{E}_{\lambda}\ \ \ \ \ \ \ \ \ \ \text{in
}\mathbb{Z}\left[  S_{n}\right]
\]
has been proved, then we can apply the base change morphism $f_{\ast}$ to this
equality, and obtain the analogous equality in $\mathbf{k}\left[
S_{n}\right]  $. In other words, if Lemma \ref{lem.spechtmod.Elam.expand1} has
been proved for $\mathbf{k}=\mathbb{Z}$, then it holds for all $\mathbf{k}$.
\medskip

\textit{Step 2:} In order to prove Lemma \ref{lem.spechtmod.Elam.expand1} for
$\mathbf{k}=\mathbb{Z}$, it suffices to prove it for $\mathbf{k}=\mathbb{Q}$.

Indeed, this too follows from a base change argument, except that we apply the
base change morphism backwards now. Namely, the canonical ring morphism
$f:\mathbb{Z}\rightarrow\mathbb{Q}$ is injective (indeed, it is just the
inclusion map), and thus the base change morphism $f_{\ast}:\mathbb{Z}\left[
S_{n}\right]  \rightarrow\mathbb{Q}\left[  S_{n}\right]  $ corresponding to
$f$ is injective as well. However, as we already saw in Step 1, the equality%
\[
n!^{3}\cdot\mathbf{z}_{C}=\sum_{\lambda\text{ is a partition of }n}\left(
f^{\lambda}\right)  ^{2}\cdot\left[  1\right]  \left(  \mathbf{E}_{\lambda
}\mathbf{z}_{C}\right)  \cdot\mathbf{E}_{\lambda}\ \ \ \ \ \ \ \ \ \ \text{in
}\mathbb{Q}\left[  S_{n}\right]
\]
can be obtained from the equality%
\[
n!^{3}\cdot\mathbf{z}_{C}=\sum_{\lambda\text{ is a partition of }n}\left(
f^{\lambda}\right)  ^{2}\cdot\left[  1\right]  \left(  \mathbf{E}_{\lambda
}\mathbf{z}_{C}\right)  \cdot\mathbf{E}_{\lambda}\ \ \ \ \ \ \ \ \ \ \text{in
}\mathbb{Z}\left[  S_{n}\right]
\]
by applying the base change morphism $f_{\ast}:\mathbb{Z}\left[  S_{n}\right]
\rightarrow\mathbb{Q}\left[  S_{n}\right]  $. Thus, conversely, the latter
equality can be obtained from the former equality by \textquotedblleft
un-applying\textquotedblright\ this morphism $f_{\ast}$ (which is allowed,
since $f_{\ast}$ is injective). Hence, if the former equality holds, then so
does the latter equality. In other words, if Lemma
\ref{lem.spechtmod.Elam.expand1} has been proved for $\mathbf{k}=\mathbb{Q}$,
then it holds for $\mathbf{k}=\mathbb{Z}$ and therefore for all $\mathbf{k}$
(as we saw in Step 1). \medskip

\textit{Step 3:} It thus remains to prove Lemma
\ref{lem.spechtmod.Elam.expand1} for $\mathbf{k}=\mathbb{Q}$.

So let us do this. Let $\mathbf{k}=\mathbb{Q}$. Then, $\mathbf{k}$ is a field
of characteristic $0$.

Theorem \ref{thm.center.conjsums} (applied to $G=S_{n}$) shows that the center
$Z\left(  \mathbf{k}\left[  S_{n}\right]  \right)  $ of the group algebra
$\mathbf{k}\left[  S_{n}\right]  $ is the $\mathbf{k}$-linear span of the
conjugacy class sums. Hence, in particular, each conjugacy class sum of
$S_{n}$ belongs to $Z\left(  \mathbf{k}\left[  S_{n}\right]  \right)  $.

The element $\mathbf{z}_{C}=\sum_{c\in C}c$ is a conjugacy class sum (since it
is the sum of all permutations in the conjugacy class $C$). Thus, it belongs
to $Z\left(  \mathbf{k}\left[  S_{n}\right]  \right)  $ (since each conjugacy
class sum of $S_{n}$ belongs to $Z\left(  \mathbf{k}\left[  S_{n}\right]
\right)  $). Hence, Lemma \ref{lem.spechtmod.Elam.a-field} (applied to
$\mathbf{a}=\mathbf{z}_{C}$) yields%
\[
\mathbf{z}_{C}=\sum_{\lambda\text{ is a partition of }n}\dfrac{\left[
1\right]  \left(  \mathbf{E}_{\lambda}\mathbf{z}_{C}\right)  }{\left(
h^{\lambda}\right)  ^{2}\cdot n!}\mathbf{E}_{\lambda}.
\]
Multiplying both sides of this equality by $n!^{3}$, we obtain%
\begin{align*}
n!^{3}\cdot\mathbf{z}_{C}  &  =n!^{3}\cdot\sum_{\lambda\text{ is a partition
of }n}\dfrac{\left[  1\right]  \left(  \mathbf{E}_{\lambda}\mathbf{z}%
_{C}\right)  }{\left(  h^{\lambda}\right)  ^{2}\cdot n!}\mathbf{E}_{\lambda}\\
&  =\sum_{\lambda\text{ is a partition of }n}\underbrace{n!^{3}\cdot
\dfrac{\left[  1\right]  \left(  \mathbf{E}_{\lambda}\mathbf{z}_{C}\right)
}{\left(  h^{\lambda}\right)  ^{2}\cdot n!}}_{\substack{=\left(  \dfrac
{n!}{h^{\lambda}}\right)  ^{2}\cdot\left[  1\right]  \left(  \mathbf{E}%
_{\lambda}\mathbf{z}_{C}\right)  \\=\left(  f^{\lambda}\right)  ^{2}%
\cdot\left[  1\right]  \left(  \mathbf{E}_{\lambda}\mathbf{z}_{C}\right)
\\\text{(since }\dfrac{n!}{h^{\lambda}}=f^{\lambda}\text{ (because }%
h^{\lambda}=\dfrac{n!}{f^{\lambda}}\text{))}}}\mathbf{E}_{\lambda}\\
&  =\sum_{\lambda\text{ is a partition of }n}\left(  f^{\lambda}\right)
^{2}\cdot\left[  1\right]  \left(  \mathbf{E}_{\lambda}\mathbf{z}_{C}\right)
\cdot\mathbf{E}_{\lambda}.
\end{align*}
Thus, Lemma \ref{lem.spechtmod.Elam.expand1} is proved for $\mathbf{k}%
=\mathbb{Q}$. As we explained above, this completes the proof of Lemma
\ref{lem.spechtmod.Elam.expand1} for all $\mathbf{k}$.
\end{proof}

We can now extend Lemma \ref{lem.spechtmod.Elam.a-field} to all commutative
rings $\mathbf{k}$ (upon multiplication by $n!^{3}$ in order to clear all denominators):

\begin{proposition}
\label{prop.spechtmod.Elam.a}Let $\mathbf{a}\in Z\left(  \mathbf{k}\left[
S_{n}\right]  \right)  $. Then,%
\[
n!^{3}\cdot\mathbf{a}=\sum_{\lambda\text{ is a partition of }n}\left(
f^{\lambda}\right)  ^{2}\cdot\left[  1\right]  \left(  \mathbf{E}_{\lambda
}\mathbf{a}\right)  \cdot\mathbf{E}_{\lambda}.
\]

\end{proposition}

\begin{proof}
Theorem \ref{thm.center.conjsums} (applied to $G=S_{n}$) shows that the center
$Z\left(  \mathbf{k}\left[  S_{n}\right]  \right)  $ of the group algebra
$\mathbf{k}\left[  S_{n}\right]  $ is the $\mathbf{k}$-linear span of the
conjugacy class sums. In other words, the conjugacy class sums span $Z\left(
\mathbf{k}\left[  S_{n}\right]  \right)  $.

We need to prove the equality%
\[
n!^{3}\cdot\mathbf{a}=\sum_{\lambda\text{ is a partition of }n}\left(
f^{\lambda}\right)  ^{2}\cdot\left[  1\right]  \left(  \mathbf{E}_{\lambda
}\mathbf{a}\right)  \cdot\mathbf{E}_{\lambda}.
\]
Both sides of this equality are $\mathbf{k}$-linear in $\mathbf{a}$. Hence, in
proving this equality, we can WLOG assume that $\mathbf{a}$ is a conjugacy
class sum (since the conjugacy class sums span $Z\left(  \mathbf{k}\left[
S_{n}\right]  \right)  $). Assume this. Thus, $\mathbf{a}=\mathbf{z}_{C}$ for
some conjugacy class $C$ of $S_{n}$ (where $\mathbf{z}_{C}$ is as defined in
Definition \ref{def.groups.conj} \textbf{(c)}). Consider this $C$. Thus,
$\mathbf{z}_{C}=\sum_{c\in C}c$ (by the definition of $\mathbf{z}_{C}$).
Hence, Lemma \ref{lem.spechtmod.Elam.expand1} yields%
\[
n!^{3}\cdot\mathbf{z}_{C}=\sum_{\lambda\text{ is a partition of }n}\left(
f^{\lambda}\right)  ^{2}\cdot\left[  1\right]  \left(  \mathbf{E}_{\lambda
}\mathbf{z}_{C}\right)  \cdot\mathbf{E}_{\lambda}.
\]
Since $\mathbf{a}=\mathbf{z}_{C}$, we can rewrite this as%
\[
n!^{3}\cdot\mathbf{a}=\sum_{\lambda\text{ is a partition of }n}\left(
f^{\lambda}\right)  ^{2}\cdot\left[  1\right]  \left(  \mathbf{E}_{\lambda
}\mathbf{a}\right)  \cdot\mathbf{E}_{\lambda}.
\]
Thus, Proposition \ref{prop.spechtmod.Elam.a} is proved.
\end{proof}

As a consequence, we can now generalize Theorem
\ref{thm.spechtmod.Elam.center} \textbf{(a)} to all commutative rings
$\mathbf{k}$ in which $n!$ is invertible:

\begin{theorem}
\label{thm.spechtmod.Elam.center-gen}Assume that $n!$ is invertible in
$\mathbf{k}$. Then, the family%
\[
\left(  \mathbf{E}_{\lambda}\right)  _{\lambda\text{ is a partition of }n}%
\]
is a basis of the $\mathbf{k}$-module $Z\left(  \mathbf{k}\left[
S_{n}\right]  \right)  $.
\end{theorem}

\begin{proof}
Let $P$ be the set of all partitions of $n$. Thus, we must show that the
family $\left(  \mathbf{E}_{\lambda}\right)  _{\lambda\in P}$ is a basis of
the $\mathbf{k}$-module $Z\left(  \mathbf{k}\left[  S_{n}\right]  \right)  $.

Clearly, this family really consists of vectors in $Z\left(  \mathbf{k}\left[
S_{n}\right]  \right)  $, since Proposition \ref{prop.spechtmod.Elam.incent}
shows that each $\lambda\in P$ satisfies $\mathbf{E}_{\lambda}\in Z\left(
\mathbf{k}\left[  S_{n}\right]  \right)  $. Thus, it remains to prove the
following two claims:

\begin{statement}
\textit{Claim 1:} The family $\left(  \mathbf{E}_{\lambda}\right)
_{\lambda\in P}$ is $\mathbf{k}$-linearly independent.
\end{statement}

\begin{statement}
\textit{Claim 2:} The family $\left(  \mathbf{E}_{\lambda}\right)
_{\lambda\in P}$ spans the $\mathbf{k}$-module $Z\left(  \mathbf{k}\left[
S_{n}\right]  \right)  $.
\end{statement}

\begin{proof}
[Proof of Claim 1.]This proof is mostly analogous to the proof of Claim 1
during the above proof of Theorem \ref{thm.spechtmod.Elam.basis}. The only
part that needs to be changed is the argument that we use to derive
$\alpha_{\lambda}=0$ from $0=\alpha_{\lambda}\mathbf{E}_{\lambda}^{2}$ (since
$\mathbf{k}$ is no longer assumed to be a field). This part can now be done as
follows: We have $f^{\lambda}h^{\lambda}=n!$ (since $h^{\lambda}=\dfrac
{n!}{f^{\lambda}}$). Thus, the number $h^{\lambda}$ is invertible in
$\mathbf{k}$ (since $f^{\lambda}h^{\lambda}=n!$ is invertible in $\mathbf{k}%
$). Hence, its square $\left(  h^{\lambda}\right)  ^{2}$ is invertible in
$\mathbf{k}$ as well. But we also have%
\[
0=\alpha_{\lambda}\underbrace{\mathbf{E}_{\lambda}^{2}}_{\substack{=\left(
h^{\lambda}\right)  ^{2}\mathbf{E}_{\lambda}\\\text{(by Theorem
\ref{thm.spechtmod.Elam.idp-strong})}}}=\alpha_{\lambda}\left(  h^{\lambda
}\right)  ^{2}\mathbf{E}_{\lambda}.
\]
We can divide this equality by $\left(  h^{\lambda}\right)  ^{2}$ (since
$\left(  h^{\lambda}\right)  ^{2}$ is invertible in $\mathbf{k}$), and thus
obtain $0=\alpha_{\lambda}\mathbf{E}_{\lambda}$. Hence,
\[
\left[  1\right]  0=\left[  1\right]  \left(  \alpha_{\lambda}\mathbf{E}%
_{\lambda}\right)  =\alpha_{\lambda}\cdot\underbrace{\left[  1\right]
\mathbf{E}_{\lambda}}_{\substack{=n!\\\text{(by Lemma
\ref{lem.specht.Elam.1-coord})}}}=\alpha_{\lambda}\cdot n!.
\]
In other words, $\alpha_{\lambda}\cdot n!=\left[  1\right]  0=0$. We can
divide this equality by $n!$ (since $n!$ is invertible in $\mathbf{k}$), and
obtain $\alpha_{\lambda}=0$. Thus, $\alpha_{\lambda}=0$ has been derived from
$0=\alpha_{\lambda}\mathbf{E}_{\lambda}^{2}$. As we said above, the rest of
the proof of Claim 1 proceeds precisely as in the proof of Theorem
\ref{thm.spechtmod.Elam.basis}.
\end{proof}

\begin{proof}
[Proof of Claim 2.]Let $\mathbf{a}\in Z\left(  \mathbf{k}\left[  S_{n}\right]
\right)  $. Then, Proposition \ref{prop.spechtmod.Elam.a} yields%
\[
n!^{3}\cdot\mathbf{a}=\sum_{\lambda\text{ is a partition of }n}\left(
f^{\lambda}\right)  ^{2}\cdot\left[  1\right]  \left(  \mathbf{E}_{\lambda
}\mathbf{a}\right)  \cdot\mathbf{E}_{\lambda}.
\]
We can divide this equality by $n!^{3}$ (since $n!$ and therefore also
$n!^{3}$ is invertible in $\mathbf{k}$), and obtain%
\begin{align*}
\mathbf{a}  &  =\dfrac{1}{n!^{3}}\sum_{\lambda\text{ is a partition of }%
n}\left(  f^{\lambda}\right)  ^{2}\cdot\left[  1\right]  \left(
\mathbf{E}_{\lambda}\mathbf{a}\right)  \cdot\mathbf{E}_{\lambda}\\
&  =\sum_{\lambda\text{ is a partition of }n}\dfrac{\left(  f^{\lambda
}\right)  ^{2}\cdot\left[  1\right]  \left(  \mathbf{E}_{\lambda}%
\mathbf{a}\right)  }{n!^{3}}\cdot\mathbf{E}_{\lambda}\\
&  =\sum_{\lambda\in P}\dfrac{\left(  f^{\lambda}\right)  ^{2}\cdot\left[
1\right]  \left(  \mathbf{E}_{\lambda}\mathbf{a}\right)  }{n!^{3}}%
\cdot\mathbf{E}_{\lambda}%
\end{align*}
(since the set of all partitions of $n$ is $P$). This shows that $\mathbf{a}$
is a $\mathbf{k}$-linear combination of the family $\left(  \mathbf{E}%
_{\lambda}\right)  _{\lambda\in P}$.

Forget that we fixed $\mathbf{a}$. We thus have shown that each $\mathbf{a}\in
Z\left(  \mathbf{k}\left[  S_{n}\right]  \right)  $ is a $\mathbf{k}$-linear
combination of the family $\left(  \mathbf{E}_{\lambda}\right)  _{\lambda\in
P}$. In other words, the family $\left(  \mathbf{E}_{\lambda}\right)
_{\lambda\in P}$ spans the $\mathbf{k}$-module $Z\left(  \mathbf{k}\left[
S_{n}\right]  \right)  $. This proves Claim 2.
\end{proof}

Combining Claim 1 with Claim 2, we conclude that the family $\left(
\mathbf{E}_{\lambda}\right)  _{\lambda\in P}$ is a basis of $Z\left(
\mathbf{k}\left[  S_{n}\right]  \right)  $. As we said, this proves Theorem
\ref{thm.spechtmod.Elam.center-gen}.
\end{proof}

Proposition \ref{prop.spechtmod.Elam.a} also allows us to extend Theorem
\ref{thm.spechtmod.Elam.pou} to a larger class of rings $\mathbf{k}$:

\begin{corollary}
\label{cor.spechtmod.Elam.pou}Assume that $n!$ is invertible in $\mathbf{k}$.
Then: \medskip

\textbf{(a)} For any partition $\lambda$ of $n$, the number $\left(
h^{\lambda}\right)  ^{2}$ is invertible in $\mathbf{k}$. \medskip

\textbf{(b)} In $\mathbf{k}\left[  S_{n}\right]  $, we have%
\[
1=\sum_{\lambda\text{ is a partition of }n}\dfrac{\mathbf{E}_{\lambda}%
}{\left(  h^{\lambda}\right)  ^{2}}.
\]

\end{corollary}

\begin{proof}
\textbf{(a)} Let $\lambda$ be a partition of $n$. Then, the definition of
$h^{\lambda}$ yields $h^{\lambda}=\dfrac{n!}{f^{\lambda}}\mid n!$, so that
$h^{\lambda}$ is invertible in $\mathbf{k}$ (since $n!$ is invertible in
$\mathbf{k}$). Therefore, $\left(  h^{\lambda}\right)  ^{2}$ is also
invertible in $\mathbf{k}$. This proves Corollary \ref{cor.spechtmod.Elam.pou}
\textbf{(a)}. \medskip

\textbf{(b)} Corollary \ref{cor.spechtmod.Elam.pou} \textbf{(a)} shows that
$\left(  h^{\lambda}\right)  ^{2}$ is invertible in $\mathbf{k}$ for any
partition $\lambda$ of $n$. Thus, the sum $\sum_{\lambda\text{ is a partition
of }n}\dfrac{\mathbf{E}_{\lambda}}{\left(  h^{\lambda}\right)  ^{2}}$ is well-defined.

We have $1\in Z\left(  \mathbf{k}\left[  S_{n}\right]  \right)  $. Thus,
Proposition \ref{prop.spechtmod.Elam.a} (applied to $\mathbf{a}=1$) yields%
\begin{align*}
n!^{3}\cdot1  &  =\sum_{\lambda\text{ is a partition of }n}\left(  f^{\lambda
}\right)  ^{2}\cdot\underbrace{\left[  1\right]  \left(  \mathbf{E}_{\lambda
}1\right)  }_{\substack{=\left[  1\right]  \mathbf{E}_{\lambda}=n!\\\text{(by
Lemma \ref{lem.specht.Elam.1-coord})}}}\cdot\,\mathbf{E}_{\lambda}\\
&  =\sum_{\lambda\text{ is a partition of }n}\left(  f^{\lambda}\right)
^{2}\cdot n!\cdot\mathbf{E}_{\lambda}.
\end{align*}
Dividing this equality by $n!^{3}$, we find%
\begin{align*}
1  &  =\dfrac{1}{n!^{3}}\sum_{\lambda\text{ is a partition of }n}\left(
f^{\lambda}\right)  ^{2}\cdot n!\cdot\mathbf{E}_{\lambda}\\
&  =\sum_{\lambda\text{ is a partition of }n}\underbrace{\dfrac{1}{n!^{3}%
}\left(  f^{\lambda}\right)  ^{2}\cdot n!}_{\substack{=\left(  f^{\lambda
}/n!\right)  ^{2}=1/\left(  h^{\lambda}\right)  ^{2}\\\text{(since }%
h^{\lambda}=n!/f^{\lambda}\text{ and}\\\text{thus }1/\left(  h^{\lambda
}\right)  ^{2}=1/\left(  n!/f^{\lambda}\right)  ^{2}=\left(  f^{\lambda
}/n!\right)  ^{2}\text{)}}}\cdot\,\mathbf{E}_{\lambda}\\
&  =\sum_{\lambda\text{ is a partition of }n}\left(  1/\left(  h^{\lambda
}\right)  ^{2}\right)  \cdot\mathbf{E}_{\lambda}=\sum_{\lambda\text{ is a
partition of }n}\dfrac{\mathbf{E}_{\lambda}}{\left(  h^{\lambda}\right)  ^{2}%
}.
\end{align*}
This proves Corollary \ref{cor.spechtmod.Elam.pou} \textbf{(b)}.
\end{proof}

\begin{exercise}
\fbox{2} Prove that $n!^{2}=\sum_{\lambda\text{ is a partition of }n}\left(
f^{\lambda}\right)  ^{2}\mathbf{E}_{\lambda}$ in $\mathbf{k}\left[
S_{n}\right]  $ (for arbitrary $\mathbf{k}$).

[Be careful not to divide by $n!$ without justification, since $n!$ might not
be invertible.]
\end{exercise}

Combining Theorem \ref{thm.spechtmod.Elam.idp-strong}, Corollary
\ref{cor.spechtmod.Elam.pou} and Proposition \ref{prop.spechtmod.Elam.orth},
we obtain a decomposition of $\mathcal{A}=\mathbf{k}\left[  S_{n}\right]  $ as
a left $\mathcal{A}$-module:

\begin{corollary}
\label{cor.spechtmod.Elam.A=sum1}Assume that $n!$ is invertible in
$\mathbf{k}$. Then,%
\[
\mathcal{A}=\bigoplus_{\lambda\text{ is a partition of }n}\mathcal{A}%
\mathbf{E}_{\lambda}\ \ \ \ \ \ \ \ \ \ \left(  \text{an internal direct
sum}\right)
\]
as left $\mathcal{A}$-modules.
\end{corollary}

\begin{proof}
Let $P$ be the set of all partitions of $n$. We must then show that
$\mathcal{A}=\bigoplus\limits_{\lambda\in P}\mathcal{A}\mathbf{E}_{\lambda}$
(an internal direct sum).

First, we shall show that the sum $\sum_{\lambda\in P}\mathcal{A}%
\mathbf{E}_{\lambda}$ is direct.

Indeed, let $\left(  \mathbf{a}_{\lambda}\right)  _{\lambda\in P}\in
\prod\limits_{\lambda\in P}\mathcal{A}\mathbf{E}_{\lambda}$ be a family of
elements such that $\sum_{\lambda\in P}\mathbf{a}_{\lambda}=0$. We shall show
that $\mathbf{a}_{\lambda}=0$ for all $\lambda\in P$.

Indeed, we have $\left(  \mathbf{a}_{\lambda}\right)  _{\lambda\in P}\in
\prod_{\lambda\in P}\mathcal{A}\mathbf{E}_{\lambda}$. In other words,%
\begin{equation}
\mathbf{a}_{\lambda}\in\mathcal{A}\mathbf{E}_{\lambda}%
\ \ \ \ \ \ \ \ \ \ \text{for each }\lambda\in P.
\label{pf.cor.spechtmod.Elam.A=sum1.2}%
\end{equation}

Now, let $\mu\in P$ be arbitrary. Then, $\mu$ is a partition of $n$ (since $P$
is the set of all partitions of $n$). Thus, Corollary
\ref{cor.spechtmod.Elam.pou} \textbf{(a)} (applied to $\lambda=\mu$) shows
that the number $\left(  h^{\mu}\right)  ^{2}$ is invertible in $\mathbf{k}$.
We have%
\begin{equation}
\sum_{\lambda\in P}\mathbf{a}_{\lambda}\mathbf{E}_{\mu}=\mathbf{a}_{\mu
}\mathbf{E}_{\mu}+\sum_{\substack{\lambda\in P;\\\lambda\neq\mu}%
}\mathbf{a}_{\lambda}\mathbf{E}_{\mu} \label{pf.cor.spechtmod.Elam.A=sum1.3}%
\end{equation}
(here, we have split off the addend for $\lambda=\mu$ from the sum).

We have $\mathbf{a}_{\mu}\in\mathcal{A}\mathbf{E}_{\mu}$ (by
(\ref{pf.cor.spechtmod.Elam.A=sum1.2}), applied to $\lambda=\mu$). In other
words, $\mathbf{a}_{\mu}=\mathbf{cE}_{\mu}$ for some $\mathbf{c}\in
\mathcal{A}$. Consider this $\mathbf{c}$. Then,%
\begin{align}
\underbrace{\mathbf{a}_{\mu}}_{=\mathbf{cE}_{\mu}}\mathbf{E}_{\mu}  &
=\mathbf{c}\underbrace{\mathbf{E}_{\mu}\mathbf{E}_{\mu}}%
_{\substack{=\mathbf{E}_{\mu}^{2}=\left(  h^{\mu}\right)  ^{2}\mathbf{E}_{\mu
}\\\text{(by Theorem \ref{thm.spechtmod.Elam.idp-strong},}\\\text{applied to
}\lambda=\mu\text{)}}}=\mathbf{c}\left(  h^{\mu}\right)  ^{2}\mathbf{E}_{\mu
}=\left(  h^{\mu}\right)  ^{2}\underbrace{\mathbf{cE}_{\mu}}_{=\mathbf{a}%
_{\mu}}\nonumber\\
&  =\left(  h^{\mu}\right)  ^{2}\mathbf{a}_{\mu}.
\label{pf.cor.spechtmod.Elam.A=sum1.4}%
\end{align}

For any $\lambda\in P$ satisfying $\lambda\neq\mu$, we have $\mathbf{a}%
_{\lambda}\mathbf{E}_{\mu}=0$ (this follows easily from Proposition
\ref{prop.spechtmod.Elam.orth})\footnote{\textit{Proof.} Let $\lambda\in P$ be
such that $\lambda\neq\mu$. Then, $\lambda,\mu\in P$. In other words,
$\lambda$ and $\mu$ are two partitions of $n$ (since $P$ is the set of all
partitions of $n$). Moreover, $\lambda$ and $\mu$ are distinct (since
$\lambda\neq\mu$). Thus, Proposition \ref{prop.spechtmod.Elam.orth} yields
$\mathbf{E}_{\lambda}\mathbf{E}_{\mu}=0$. However,
(\ref{pf.cor.spechtmod.Elam.A=sum1.2}) yields $\mathbf{a}_{\lambda}%
\in\mathcal{A}\mathbf{E}_{\lambda}$. In other words, $\mathbf{a}_{\lambda
}=\mathbf{bE}_{\lambda}$ for some $\mathbf{b}\in\mathcal{A}$. Consider this
$\mathbf{b}$. Then, $\underbrace{\mathbf{a}_{\lambda}}_{=\mathbf{bE}_{\lambda
}}\mathbf{E}_{\mu}=\mathbf{b}\underbrace{\mathbf{E}_{\lambda}\mathbf{E}_{\mu}%
}_{=0}=0$, qed.}. Hence, $\sum_{\substack{\lambda\in P;\\\lambda\neq\mu
}}\underbrace{\mathbf{a}_{\lambda}\mathbf{E}_{\mu}}_{=0}=\sum
_{\substack{\lambda\in P;\\\lambda\neq\mu}}0=0$. Thus,
(\ref{pf.cor.spechtmod.Elam.A=sum1.3}) becomes%
\[
\sum_{\lambda\in P}\mathbf{a}_{\lambda}\mathbf{E}_{\mu}=\underbrace{\mathbf{a}%
_{\mu}\mathbf{E}_{\mu}}_{\substack{=\left(  h^{\mu}\right)  ^{2}%
\mathbf{a}_{\mu}\\\text{(by (\ref{pf.cor.spechtmod.Elam.A=sum1.4}))}%
}}+\underbrace{\sum_{\substack{\lambda\in P;\\\lambda\neq\mu}}\mathbf{a}%
_{\lambda}\mathbf{E}_{\mu}}_{=0}=\left(  h^{\mu}\right)  ^{2}\mathbf{a}_{\mu
}.
\]
Hence,%
\[
\left(  h^{\mu}\right)  ^{2}\mathbf{a}_{\mu}=\sum_{\lambda\in P}%
\mathbf{a}_{\lambda}\mathbf{E}_{\mu}=\underbrace{\left(  \sum_{\lambda\in
P}\mathbf{a}_{\lambda}\right)  }_{=0}\mathbf{E}_{\mu}=0.
\]
We can divide this equality by $\left(  h^{\mu}\right)  ^{2}$ (since $\left(
h^{\mu}\right)  ^{2}$ is invertible in $\mathbf{k}$), and thus obtain
$\mathbf{a}_{\mu}=0$.

Forget that we fixed $\mu$. We thus have shown that $\mathbf{a}_{\mu}=0$ for
each $\mu\in P$. In other words, $\mathbf{a}_{\lambda}=0$ for each $\lambda\in
P$.

Forget that we fixed $\left(  \mathbf{a}_{\lambda}\right)  _{\lambda\in P}$.
We thus have shown that if $\left(  \mathbf{a}_{\lambda}\right)  _{\lambda\in
P}\in\prod_{\lambda\in P}\mathcal{A}\mathbf{E}_{\lambda}$ is a family of
elements such that $\sum_{\lambda\in P}\mathbf{a}_{\lambda}=0$, then we have
$\mathbf{a}_{\lambda}=0$ for all $\lambda\in P$. In other words, the sum
$\sum_{\lambda\in P}\mathcal{A}\mathbf{E}_{\lambda}$ is direct. Thus, we can
write it as $\bigoplus\limits_{\lambda\in P}\mathcal{A}\mathbf{E}_{\lambda}$
(an internal direct sum).

We shall now show that this sum is $\mathcal{A}$. Indeed, Corollary
\ref{cor.spechtmod.Elam.pou} \textbf{(b)} yields%
\begin{equation}
1=\underbrace{\sum_{\lambda\text{ is a partition of }n}}_{\substack{=\sum
_{\lambda\in P}\\\text{(by the definition of }P\text{)}}}\dfrac{\mathbf{E}%
_{\lambda}}{\left(  h^{\lambda}\right)  ^{2}}=\sum_{\lambda\in P}%
\dfrac{\mathbf{E}_{\lambda}}{\left(  h^{\lambda}\right)  ^{2}}.
\label{pf.cor.spechtmod.Elam.A=sum1.9}%
\end{equation}
Hence, for each $\mathbf{a}\in\mathcal{A}$, we have%
\begin{align*}
\mathbf{a}  &  =\mathbf{a}\cdot1=\mathbf{a}\cdot\sum_{\lambda\in P}%
\dfrac{\mathbf{E}_{\lambda}}{\left(  h^{\lambda}\right)  ^{2}}%
\ \ \ \ \ \ \ \ \ \ \left(  \text{by (\ref{pf.cor.spechtmod.Elam.A=sum1.9}%
)}\right) \\
&  =\sum_{\lambda\in P}\underbrace{\dfrac{\mathbf{a}\cdot\mathbf{E}_{\lambda}%
}{\left(  h^{\lambda}\right)  ^{2}}}_{=\dfrac{\mathbf{a}}{\left(  h^{\lambda
}\right)  ^{2}}\mathbf{E}_{\lambda}}=\sum_{\lambda\in P}\underbrace{\dfrac
{\mathbf{a}}{\left(  h^{\lambda}\right)  ^{2}}}_{\in\mathcal{A}}%
\mathbf{E}_{\lambda}\in\sum_{\lambda\in P}\mathcal{A}\mathbf{E}_{\lambda}.
\end{align*}
In other words, $\mathcal{A}\subseteq\sum_{\lambda\in P}\mathcal{A}%
\mathbf{E}_{\lambda}$. Since $\sum_{\lambda\in P}\mathcal{A}\mathbf{E}%
_{\lambda}$ is clearly a subset of $\mathcal{A}$, we thus obtain
\begin{align*}
\mathcal{A}  &  =\sum_{\lambda\in P}\mathcal{A}\mathbf{E}_{\lambda}%
=\bigoplus\limits_{\lambda\in P}\mathcal{A}\mathbf{E}_{\lambda}%
\ \ \ \ \ \ \ \ \ \ \left(  \text{as we saw above}\right) \\
&  =\bigoplus\limits_{\lambda\text{ is a partition of }n}\mathcal{A}%
\mathbf{E}_{\lambda}\ \ \ \ \ \ \ \ \ \ \left(  \text{since }P\text{ is the
set of all partitions of }n\right)  .
\end{align*}
This proves Corollary \ref{cor.spechtmod.Elam.A=sum1}.
\end{proof}

\subsubsection{The structure of $Z\left(  \mathbf{k}\left[  S_{n}\right]
\right)  $}

Now we are ready to describe $Z\left(  \mathbf{k}\left[  S_{n}\right]
\right)  $ as a $\mathbf{k}$-algebra when $n!$ is invertible in $\mathbf{k}$
(and thus prove Theorem \ref{thm.AWS.demo} \textbf{(b)}):

\begin{theorem}
\label{thm.spechtmod.Elam.iso}Assume that $n!$ is invertible in $\mathbf{k}$.
Then, there exists a $\mathbf{k}$-algebra isomorphism%
\begin{align*}
\mathbf{k}^{\left\{  \text{partitions of }n\right\}  }  &  \rightarrow
Z\left(  \mathbf{k}\left[  S_{n}\right]  \right)  ,\\
\left(  \beta_{\lambda}\right)  _{\lambda\text{ is a partition of }n}  &
\mapsto\sum_{\lambda\text{ is a partition of }n}\beta_{\lambda}\dfrac
{\mathbf{E}_{\lambda}}{\left(  h^{\lambda}\right)  ^{2}}.
\end{align*}
(Here, $\mathbf{k}^{\left\{  \text{partitions of }n\right\}  }$ is a direct
product of copies of $\mathbf{k}$, one copy for each partition of $n$. The
addition and the multiplication on this direct product are entrywise, as usual.)
\end{theorem}

\begin{proof}
We begin with some simple observations. Theorem
\ref{thm.spechtmod.Elam.center-gen} shows that the family%
\[
\left(  \mathbf{E}_{\lambda}\right)  _{\lambda\text{ is a partition of }n}%
\]
is a basis of the $\mathbf{k}$-module $Z\left(  \mathbf{k}\left[
S_{n}\right]  \right)  $. Hence, this family is $\mathbf{k}$-linearly
independent and spans $Z\left(  \mathbf{k}\left[  S_{n}\right]  \right)  $.

Next, we claim the following:

\begin{statement}
\textit{Claim 1:} For any family $\left(  \beta_{\lambda}\right)
_{\lambda\text{ is a partition of }n}\in\mathbf{k}^{\left\{  \text{partitions
of }n\right\}  }$, the sum $\sum_{\lambda\text{ is a partition of }n}%
\beta_{\lambda}\dfrac{\mathbf{E}_{\lambda}}{\left(  h^{\lambda}\right)  ^{2}}$
is a well-defined element of $Z\left(  \mathbf{k}\left[  S_{n}\right]
\right)  $.
\end{statement}

\begin{proof}
[Proof of Claim 1.]It clearly suffices to prove that $\dfrac{\mathbf{E}%
_{\lambda}}{\left(  h^{\lambda}\right)  ^{2}}$ is a well-defined element of
$Z\left(  \mathbf{k}\left[  S_{n}\right]  \right)  $ whenever $\lambda$ is a
partition of $n$.

So let us do this. Let $\lambda$ be a partition of $n$. Then, Proposition
\ref{prop.spechtmod.Elam.incent} shows that $\mathbf{E}_{\lambda}$ belongs to
$Z\left(  \mathbf{k}\left[  S_{n}\right]  \right)  $. Thus, $\dfrac
{\mathbf{E}_{\lambda}}{\left(  h^{\lambda}\right)  ^{2}}$ is a well-defined
element of $Z\left(  \mathbf{k}\left[  S_{n}\right]  \right)  $ (since
Corollary \ref{cor.spechtmod.Elam.pou} \textbf{(a)} shows that $\left(
h^{\lambda}\right)  ^{2}$ is invertible in $\mathbf{k}$). This completes our
proof of Claim 1.
\end{proof}

Claim 1 allows us to define a map%
\begin{align*}
\Phi:\mathbf{k}^{\left\{  \text{partitions of }n\right\}  }  &  \rightarrow
Z\left(  \mathbf{k}\left[  S_{n}\right]  \right)  ,\\
\left(  \beta_{\lambda}\right)  _{\lambda\text{ is a partition of }n}  &
\mapsto\sum_{\lambda\text{ is a partition of }n}\beta_{\lambda}\dfrac
{\mathbf{E}_{\lambda}}{\left(  h^{\lambda}\right)  ^{2}}.
\end{align*}
Consider this map $\Phi$. It is clearly $\mathbf{k}$-linear. Our goal is now
to show that $\Phi$ is a $\mathbf{k}$-algebra isomorphism (since this is
precisely what Theorem \ref{thm.spechtmod.Elam.iso} claims). We shall achieve
this by proving the following four claims:

\begin{statement}
\textit{Claim 2:} The map $\Phi$ is injective.
\end{statement}

\begin{statement}
\textit{Claim 3:} The map $\Phi$ is surjective.
\end{statement}

\begin{statement}
\textit{Claim 4:} The map $\Phi$ respects multiplication (i.e., we have
$\Phi\left(  \beta\gamma\right)  =\Phi\left(  \beta\right)  \Phi\left(
\gamma\right)  $ for all $\beta,\gamma\in\mathbf{k}^{\left\{  \text{partitions
of }n\right\}  }$).
\end{statement}

\begin{statement}
\textit{Claim 5:} The map $\Phi$ respects the unity (i.e., we have
$\Phi\left(  1_{\mathbf{k}^{\left\{  \text{partitions of }n\right\}  }%
}\right)  =1$).
\end{statement}

\begin{proof}
[Proof of Claim 2.]Let $\beta\in\operatorname*{Ker}\Phi$. We shall show that
$\beta=0$.

We have $\beta\in\operatorname*{Ker}\Phi\subseteq\mathbf{k}^{\left\{
\text{partitions of }n\right\}  }$. Therefore, write $\beta$ as $\beta=\left(
\beta_{\lambda}\right)  _{\lambda\text{ is a partition of }n}$ for some
scalars $\beta_{\lambda}\in\mathbf{k}$. Then, the definition of $\Phi$ yields%
\[
\Phi\left(  \beta\right)  =\sum_{\lambda\text{ is a partition of }%
n}\underbrace{\beta_{\lambda}\dfrac{\mathbf{E}_{\lambda}}{\left(  h^{\lambda
}\right)  ^{2}}}_{=\dfrac{\beta_{\lambda}}{\left(  h^{\lambda}\right)  ^{2}%
}\mathbf{E}_{\lambda}}=\sum_{\lambda\text{ is a partition of }n}\dfrac
{\beta_{\lambda}}{\left(  h^{\lambda}\right)  ^{2}}\mathbf{E}_{\lambda}.
\]
Hence,%
\[
\sum_{\lambda\text{ is a partition of }n}\dfrac{\beta_{\lambda}}{\left(
h^{\lambda}\right)  ^{2}}\mathbf{E}_{\lambda}=\Phi\left(  \beta\right)
=0\ \ \ \ \ \ \ \ \ \ \left(  \text{since }\beta\in\operatorname*{Ker}%
\Phi\right)  .
\]
From this equality, we conclude that%
\[
\dfrac{\beta_{\lambda}}{\left(  h^{\lambda}\right)  ^{2}}%
=0\ \ \ \ \ \ \ \ \ \ \text{for all partitions }\lambda\text{ of }n
\]
(since the family $\left(  \mathbf{E}_{\lambda}\right)  _{\lambda\text{ is a
partition of }n}$ is $\mathbf{k}$-linearly independent). Multiplying this
equality by $\left(  h^{\lambda}\right)  ^{2}$, we find that%
\[
\beta_{\lambda}=0\ \ \ \ \ \ \ \ \ \ \text{for all partitions }\lambda\text{
of }n.
\]
In other words, $\left(  \beta_{\lambda}\right)  _{\lambda\text{ is a
partition of }n}=\left(  0\right)  _{\lambda\text{ is a partition of }n}=0$.
Hence, $\beta=\left(  \beta_{\lambda}\right)  _{\lambda\text{ is a partition
of }n}=0$.

Forget that we fixed $\beta$. We thus have shown that $\beta=0$ for each
$\beta\in\operatorname*{Ker}\Phi$. In other words, $\operatorname*{Ker}\Phi
=0$. Since $\Phi$ is a $\mathbf{k}$-linear map, this entails that $\Phi$ is
injective. Thus, Claim 2 is proved.
\end{proof}

\begin{proof}
[Proof of Claim 3.]Let $\mathbf{a}\in Z\left(  \mathbf{k}\left[  S_{n}\right]
\right)  $. Then, $\mathbf{a}$ is a $\mathbf{k}$-linear combination of the
family $\left(  \mathbf{E}_{\lambda}\right)  _{\lambda\text{ is a partition of
}n}$ (since this family spans $Z\left(  \mathbf{k}\left[  S_{n}\right]
\right)  $). In other words, we can write $\mathbf{a}$ as $\mathbf{a}%
=\sum_{\lambda\text{ is a partition of }n}\gamma_{\lambda}\mathbf{E}_{\lambda
}$ for some scalars $\gamma_{\lambda}\in\mathbf{k}$. Consider these scalars
$\gamma_{\lambda}$. Now, the definition of $\Phi$ yields%
\begin{align*}
\Phi\left(  \left(  \left(  h^{\lambda}\right)  ^{2}\gamma_{\lambda}\right)
_{\lambda\text{ is a partition of }n}\right)   &  =\sum_{\lambda\text{ is a
partition of }n}\underbrace{\left(  h^{\lambda}\right)  ^{2}\gamma_{\lambda
}\dfrac{\mathbf{E}_{\lambda}}{\left(  h^{\lambda}\right)  ^{2}}}%
_{=\gamma_{\lambda}\mathbf{E}_{\lambda}}\\
&  =\sum_{\lambda\text{ is a partition of }n}\gamma_{\lambda}\mathbf{E}%
_{\lambda}=\mathbf{a},
\end{align*}
so that $\mathbf{a}=\Phi\left(  \left(  \left(  h^{\lambda}\right)  ^{2}%
\gamma_{\lambda}\right)  _{\lambda\text{ is a partition of }n}\right)
\in\operatorname{Im}\Phi$.

Forget that we fixed $\mathbf{a}$. We thus have shown that each $\mathbf{a}\in
Z\left(  \mathbf{k}\left[  S_{n}\right]  \right)  $ satisfies $\mathbf{a}%
\in\operatorname{Im}\Phi$. In other words, the map $\Phi$ is surjective. This
proves Claim 3.
\end{proof}

\begin{proof}
[Proof of Claim 4.]Let $\beta,\gamma\in\mathbf{k}^{\left\{  \text{partitions
of }n\right\}  }$. We must show that $\Phi\left(  \beta\gamma\right)
=\Phi\left(  \beta\right)  \Phi\left(  \gamma\right)  $.

Write $\beta$ as $\beta=\left(  \beta_{\lambda}\right)  _{\lambda\text{ is a
partition of }n}$ for some scalars $\beta_{\lambda}\in\mathbf{k}$. Likewise,
write $\gamma$ as $\gamma=\left(  \gamma_{\lambda}\right)  _{\lambda\text{ is
a partition of }n}$ for some scalars $\gamma_{\lambda}\in\mathbf{k}$. Thus,
\[
\beta\gamma=\left(  \beta_{\lambda}\right)  _{\lambda\text{ is a partition of
}n}\left(  \gamma_{\lambda}\right)  _{\lambda\text{ is a partition of }%
n}=\left(  \beta_{\lambda}\gamma_{\lambda}\right)  _{\lambda\text{ is a
partition of }n}.
\]
Hence, the definition of $\Phi$ yields%
\begin{equation}
\Phi\left(  \beta\gamma\right)  =\sum_{\lambda\text{ is a partition of }%
n}\beta_{\lambda}\gamma_{\lambda}\dfrac{\mathbf{E}_{\lambda}}{\left(
h^{\lambda}\right)  ^{2}}. \label{pf.thm.spechtmod.Elam.iso.c4.1}%
\end{equation}

However, the definition of $\Phi$ also yields%
\[
\Phi\left(  \beta\right)  =\sum_{\lambda\text{ is a partition of }n}%
\beta_{\lambda}\dfrac{\mathbf{E}_{\lambda}}{\left(  h^{\lambda}\right)  ^{2}%
}\ \ \ \ \ \ \ \ \ \ \left(  \text{since }\beta=\left(  \beta_{\lambda
}\right)  _{\lambda\text{ is a partition of }n}\right)
\]
and%
\begin{align*}
\Phi\left(  \gamma\right)   &  =\sum_{\lambda\text{ is a partition of }%
n}\gamma_{\lambda}\dfrac{\mathbf{E}_{\lambda}}{\left(  h^{\lambda}\right)
^{2}}\ \ \ \ \ \ \ \ \ \ \left(  \text{since }\gamma=\left(  \gamma_{\lambda
}\right)  _{\lambda\text{ is a partition of }n}\right) \\
&  =\sum_{\mu\text{ is a partition of }n}\gamma_{\mu}\dfrac{\mathbf{E}_{\mu}%
}{\left(  h^{\mu}\right)  ^{2}}.
\end{align*}
Multiplying these two equalities, we obtain%
\begin{align*}
&  \Phi\left(  \beta\right)  \Phi\left(  \gamma\right) \\
&  =\left(  \sum_{\lambda\text{ is a partition of }n}\beta_{\lambda}%
\dfrac{\mathbf{E}_{\lambda}}{\left(  h^{\lambda}\right)  ^{2}}\right)  \left(
\sum_{\mu\text{ is a partition of }n}\gamma_{\mu}\dfrac{\mathbf{E}_{\mu}%
}{\left(  h^{\mu}\right)  ^{2}}\right) \\
&  =\underbrace{\sum_{\lambda\text{ is a partition of }n}\ \ \sum_{\mu\text{
is a partition of }n}}_{=\sum_{\substack{\lambda\text{ and }\mu\text{ are
two}\\\text{partitions of }n}}}\underbrace{\beta_{\lambda}\dfrac
{\mathbf{E}_{\lambda}}{\left(  h^{\lambda}\right)  ^{2}}\gamma_{\mu}%
\dfrac{\mathbf{E}_{\mu}}{\left(  h^{\mu}\right)  ^{2}}}_{=\dfrac
{\beta_{\lambda}\gamma_{\mu}}{\left(  h^{\lambda}\right)  ^{2}\left(  h^{\mu
}\right)  ^{2}}\mathbf{E}_{\lambda}\mathbf{E}_{\mu}}\\
&  =\sum_{\substack{\lambda\text{ and }\mu\text{ are two}\\\text{partitions of
}n}}\dfrac{\beta_{\lambda}\gamma_{\mu}}{\left(  h^{\lambda}\right)
^{2}\left(  h^{\mu}\right)  ^{2}}\mathbf{E}_{\lambda}\mathbf{E}_{\mu}\\
&  =\sum_{\substack{\lambda\text{ and }\mu\text{ are two}\\\text{partitions of
}n;\\\lambda=\mu}}\dfrac{\beta_{\lambda}\gamma_{\mu}}{\left(  h^{\lambda
}\right)  ^{2}\left(  h^{\mu}\right)  ^{2}}\mathbf{E}_{\lambda}\mathbf{E}%
_{\mu}+\sum_{\substack{\lambda\text{ and }\mu\text{ are two}\\\text{partitions
of }n;\\\lambda\neq\mu}}\dfrac{\beta_{\lambda}\gamma_{\mu}}{\left(
h^{\lambda}\right)  ^{2}\left(  h^{\mu}\right)  ^{2}}\underbrace{\mathbf{E}%
_{\lambda}\mathbf{E}_{\mu}}_{\substack{=0\\\text{(by Proposition
\ref{prop.spechtmod.Elam.orth})}}}\\
&  =\sum_{\substack{\lambda\text{ and }\mu\text{ are two}\\\text{partitions of
}n;\\\lambda=\mu}}\dfrac{\beta_{\lambda}\gamma_{\mu}}{\left(  h^{\lambda
}\right)  ^{2}\left(  h^{\mu}\right)  ^{2}}\mathbf{E}_{\lambda}\mathbf{E}%
_{\mu}+\underbrace{\sum_{\substack{\lambda\text{ and }\mu\text{ are
two}\\\text{partitions of }n;\\\lambda\neq\mu}}\dfrac{\beta_{\lambda}%
\gamma_{\mu}}{\left(  h^{\lambda}\right)  ^{2}\left(  h^{\mu}\right)  ^{2}}%
0}_{=0}\\
&  =\sum_{\substack{\lambda\text{ and }\mu\text{ are two}\\\text{partitions of
}n;\\\lambda=\mu}}\dfrac{\beta_{\lambda}\gamma_{\mu}}{\left(  h^{\lambda
}\right)  ^{2}\left(  h^{\mu}\right)  ^{2}}\mathbf{E}_{\lambda}\mathbf{E}%
_{\mu}\\
&  =\sum_{\lambda\text{ is a partition of }n}\dfrac{\beta_{\lambda}%
\gamma_{\lambda}}{\left(  h^{\lambda}\right)  ^{2}\left(  h^{\lambda}\right)
^{2}}\underbrace{\mathbf{E}_{\lambda}\mathbf{E}_{\lambda}}%
_{\substack{=\mathbf{E}_{\lambda}^{2}=\left(  h^{\lambda}\right)
^{2}\mathbf{E}_{\lambda}\\\text{(by Theorem
\ref{thm.spechtmod.Elam.idp-strong})}}}\\
&  \ \ \ \ \ \ \ \ \ \ \ \ \ \ \ \ \ \ \ \ \left(
\begin{array}
[c]{c}%
\text{here, we substituted }\lambda\text{ for }\mu\text{ in the sum,}\\
\text{because of the requirement }\lambda=\mu
\end{array}
\right) \\
&  =\sum_{\lambda\text{ is a partition of }n}\underbrace{\dfrac{\beta
_{\lambda}\gamma_{\lambda}}{\left(  h^{\lambda}\right)  ^{2}\left(
h^{\lambda}\right)  ^{2}}\left(  h^{\lambda}\right)  ^{2}\mathbf{E}_{\lambda}%
}_{=\beta_{\lambda}\gamma_{\lambda}\dfrac{\mathbf{E}_{\lambda}}{\left(
h^{\lambda}\right)  ^{2}}}=\sum_{\lambda\text{ is a partition of }n}%
\beta_{\lambda}\gamma_{\lambda}\dfrac{\mathbf{E}_{\lambda}}{\left(
h^{\lambda}\right)  ^{2}}.
\end{align*}
Comparing this with (\ref{pf.thm.spechtmod.Elam.iso.c4.1}), we find
$\Phi\left(  \beta\gamma\right)  =\Phi\left(  \beta\right)  \Phi\left(
\gamma\right)  $. Thus, Claim 4 is proved.
\end{proof}

\begin{proof}
[Proof of Claim 5.]We have $1_{\mathbf{k}^{\left\{  \text{partitions of
}n\right\}  }}=\left(  1\right)  _{\lambda\text{ is a partition of }n}$. Thus,
the definition of $\Phi$ yields%
\[
\Phi\left(  1_{\mathbf{k}^{\left\{  \text{partitions of }n\right\}  }}\right)
=\sum_{\lambda\text{ is a partition of }n}1\dfrac{\mathbf{E}_{\lambda}%
}{\left(  h^{\lambda}\right)  ^{2}}=\sum_{\lambda\text{ is a partition of }%
n}\dfrac{\mathbf{E}_{\lambda}}{\left(  h^{\lambda}\right)  ^{2}}=1
\]
(by Corollary \ref{cor.spechtmod.Elam.pou} \textbf{(b)}). In other words,
$\Phi$ respects the unity. This proves Claim 5.
\end{proof}

Now, it is time to combine what we have shown. The map $\Phi$ is $\mathbf{k}%
$-linear and respects multiplication (by Claim 4) and respects the unity (by
Claim 5). Hence, it is a $\mathbf{k}$-algebra morphism. Moreover, the map
$\Phi$ is injective (by Claim 2) and surjective (by Claim 3), hence bijective,
thus invertible. Hence, it is a $\mathbf{k}$-algebra isomorphism (since it is
a $\mathbf{k}$-algebra morphism). This completes the proof of Theorem
\ref{thm.spechtmod.Elam.iso}.
\end{proof}

Thus, Theorem \ref{thm.AWS.demo} \textbf{(b)} is proved as well. \medskip

We have now reached a full understanding of the $\mathbf{k}$-algebra $Z\left(
\mathbf{k}\left[  S_{n}\right]  \right)  $ when $n!$ is invertible in
$\mathbf{k}$. This is more or less the most general case in which $Z\left(
\mathbf{k}\left[  S_{n}\right]  \right)  $ is a direct product of $\mathbf{k}%
$'s. The structure of $Z\left(  \mathbf{k}\left[  S_{n}\right]  \right)  $
when $\mathbf{k}$ has \textquotedblleft small characteristic\textquotedblright%
\ (i.e., $n!$ is not invertible in $\mathbf{k}$) is much less clear and might
not be easily described.\footnote{For example, if $\mathbf{k}$ is a field with
$0<\operatorname*{char}\mathbf{k}\leq n$, then the center $Z\left(
\mathbf{k}\left[  S_{n}\right]  \right)  $ contains the nonzero nilpotent
element $\nabla$ (it is nilpotent since $\nabla^{2}=n!\cdot\nabla=0$, because
$n!=0$ in $\mathbf{k}$), whereas a direct product of $\mathbf{k}$'s cannot
contain any nonzero nilpotent elements (since $\mathbf{k}$ is a field). Thus,
$Z\left(  \mathbf{k}\left[  S_{n}\right]  \right)  $ cannot be isomorphic to a
direct product of $\mathbf{k}$'s in this case.}

\begin{question}
What can we say about the structure of $Z\left(  \mathbf{k}\left[
S_{n}\right]  \right)  $ when $n$ is prime and $\mathbf{k}$ is a field of
characteristic $n$ ? (This case appears to be the most tractable, as it comes
the closest to the well-behaved \textquotedblleft$n!$ invertible in
$\mathbf{k}$\textquotedblright\ case.)
\end{question}

\subsection{Irreducibility, non-isomorphism and completeness}

Our next goal is the classification of irreps (i.e., irreducible
representations) of the symmetric group $S_{n}$ over a field of characteristic
$0$. The goal here is to prove the following:

\begin{theorem}
\label{thm.spechtmod.irred}Assume that $\mathbf{k}$ is a field of
characteristic $0$. Then: \medskip

\textbf{(a)} Each Specht module $\mathcal{S}^{\lambda}$ (where $\lambda$ is a
partition of $n$) is irreducible as a left $\mathbf{k}\left[  S_{n}\right]
$-module (i.e., as a representation of $S_{n}$). \medskip

\textbf{(b)} The Specht modules $\mathcal{S}^{\lambda}$ for different
partitions $\lambda$ of $n$ are mutually non-isomorphic. \medskip

\textbf{(c)} Each irreducible representation of $S_{n}$ is isomorphic to a
Specht module $\mathcal{S}^{\lambda}$ for some partition $\lambda$ of $n$.
\end{theorem}

This is a famous result (see, e.g., \cite[(28.15)]{CurRei62}, \cite[Theorem
5.12.2]{EGHetc11}, \cite[Theorem 3.3.2]{Prasad-rep}, \cite[Corollary 3.5,
Corollary 4.4, and the paragraph below Example 4.5]{Wildon18}, \cite[\S 14.7]%
{Waerde91} for proofs at various levels of generality\footnote{In
\cite[Theorem 5.12.2]{EGHetc11}, the result is only stated for $\mathbf{k}%
=\mathbb{C}$; in \cite[\S 14.7]{Waerde91}, only for algebraically closed
fields $\mathbf{k}$; in \cite[(28.15)]{CurRei62}, only for $\mathbf{k}%
=\mathbb{Q}$.}), and is commonly proved using some representation theory of
finite groups\footnote{In particular, there is a known fact (see, e.g.,
\cite[(27.25)]{CurRei62}) that if $G$ is an arbitrary finite group and
$\mathbf{k}$ is a field of characteristic $0$, then
\begin{align}
&  \left(  \text{\# of irreducible representations of }G\text{ (counted up to
isomorphism)}\right) \nonumber\\
&  \leq\left(  \text{\# of conjugacy classes of }G\right)  .
\label{eq.rmk.thm.spechtmod.irred.fn1.1}%
\end{align}
Using this fact, we can easily derive part \textbf{(c)} of Theorem
\ref{thm.spechtmod.irred} from parts \textbf{(a)} and \textbf{(b)} (since the
\# of conjugacy classes of $S_{n}$ is the \# of partitions of $n$). It is
worth noting (and this is much better known) that the inequality
(\ref{eq.rmk.thm.spechtmod.irred.fn1.1}) becomes an equality whenever
$\mathbf{k}$ is algebraically closed (see, e.g., \cite[Corollary
4.2.2]{EGHetc11} or \cite[\S 2.5, Theorem 7]{Serre77}).}. However, we shall
give a proof that requires no new representation theory apart from a few
lemmas that we will prove. It relies instead on the properties of Young
symmetrizers and centralized Young symmetrizers we have derived above.

\subsubsection{Irreducibility}

We begin by preparing for the proof of Theorem \ref{thm.spechtmod.irred}
\textbf{(a)}. To prove this, we need a lemma:

\begin{lemma}
\label{lem.spechtmod.irred.R}Let $\mathbf{k}$ be a field. Let $R$ be a
$\mathbf{k}$-algebra. Let $e\in R$ be an element. Assume that $e$ has the
following two properties:

\begin{enumerate}
\item Quasi-idempotency: There is a nonzero scalar $\kappa\in\mathbf{k}$ such
that $e^{2}=\kappa e$.

\item \textquotedblleft Sandwich property\textquotedblright: For each $a\in
R$, the product $eae$ is a scalar multiple of $e$ (that is, we have $eae=\tau
e$ for some $\tau\in\mathbf{k}$).
\end{enumerate}

Then: \medskip

\textbf{(a)} Any left $R$-module morphism from $Re$ to $Re$ is a scalar
multiple of the identity map $\operatorname*{id}$. \medskip

\textbf{(b)} In particular, any left $R$-module morphism from $Re$ to $Re$ is
either zero or an isomorphism. \medskip

\textbf{(c)} The left $R$-module $Re$ cannot be written as a direct sum of two
nonzero $R$-submodules. \medskip

\textbf{(d)} Assume that $e\neq0$. Furthermore, assume that $R=\mathbf{k}%
\left[  G\right]  $ for a finite group $G$, and that $\operatorname*{char}%
\mathbf{k}=0$. Then, the left $R$-module $Re$ is irreducible.
\end{lemma}

\begin{proof}
By the quasi-idempotency property, there is a nonzero scalar $\kappa
\in\mathbf{k}$ such that $e^{2}=\kappa e$. Consider this $\kappa$. Note that
$\kappa$ is nonzero and thus invertible (since $\mathbf{k}$ is a field).
\medskip

\textbf{(a)} Let $f:Re\rightarrow Re$ be a left $R$-module morphism. We must
show that $f$ is a scalar multiple of $\operatorname*{id}$.

We have $e\in Re$ (since $e=1e$); thus, $f\left(  e\right)  $ is well-defined.
Moreover, $f\left(  e\right)  \in Re$, so that $f\left(  e\right)  =ae$ for
some $a\in R$. Consider this $a$. The \textquotedblleft sandwich
property\textquotedblright\ shows that $eae$ is a scalar multiple of $e$. In
other words, we have $eae=\tau e$ for some $\tau\in\mathbf{k}$. Consider this
$\tau$.

We have $e^{2}=\kappa e$, so that $e=\dfrac{1}{\kappa}e^{2}=\dfrac{1}{\kappa
}ee$. Hence,%
\begin{align*}
f\left(  e\right)   &  =f\left(  \dfrac{1}{\kappa}ee\right)  =\dfrac{1}%
{\kappa}e\underbrace{f\left(  e\right)  }_{=ae}\ \ \ \ \ \ \ \ \ \ \left(
\text{since }f\text{ is left }R\text{-linear}\right) \\
&  =\dfrac{1}{\kappa}\underbrace{eae}_{=\tau e}=\dfrac{1}{\kappa}\tau
e=\dfrac{\tau}{\kappa}e.
\end{align*}
Setting $\lambda:=\dfrac{\tau}{\kappa}$, we can rewrite this as $f\left(
e\right)  =\lambda e$. Thus, for every $r\in R$, we have%
\begin{align*}
f\left(  re\right)   &  =r\underbrace{f\left(  e\right)  }_{=\lambda
e}\ \ \ \ \ \ \ \ \ \ \left(  \text{since }f\text{ is left }R\text{-linear}%
\right) \\
&  =r\cdot\lambda e=\lambda re.
\end{align*}
In other words, $f\left(  b\right)  =\lambda b$ for each $b\in Re$ (since each
$b\in Re$ has the form $re$ for some $r\in R$). In other words, $f=\lambda
\operatorname*{id}$. Thus, $f$ is a scalar multiple of $\operatorname*{id}$.
This proves Lemma \ref{lem.spechtmod.irred.R} \textbf{(a)}. \medskip

\textbf{(b)} Let $f:Re\rightarrow Re$ be a left $R$-module morphism. We must
show that $f$ is either zero or an isomorphism.

Lemma \ref{lem.spechtmod.irred.R} \textbf{(a)} shows that $f$ is a scalar
multiple of $\operatorname*{id}$. In other words, $f=\lambda\operatorname*{id}%
$ for some $\lambda\in\mathbf{k}$. If $\lambda=0$, then this entails that $f$
is zero. If $\lambda\neq0$, then this shows that $f$ is an isomorphism (with
inverse $\lambda^{-1}\operatorname*{id}$, which is well-defined because
$\mathbf{k}$ is a field). Thus, we conclude that $f$ is either zero or an
isomorphism. This proves Lemma \ref{lem.spechtmod.irred.R} \textbf{(b)}.
\medskip

\textbf{(c)} Assume the contrary. Thus, $Re=W\oplus U$ for some nonzero left
$R$-submodules $W$ and $U$. Consider these $W$ and $U$. Thus, each element of
$Re$ can be written uniquely as a sum $w+u$ for some $w\in W$ and some $u\in
U$. This allows us to define a map
\begin{align*}
f:Re  &  \rightarrow Re,\\
w+u  &  \mapsto w\ \ \ \ \ \ \ \ \ \ \text{for each }w\in W\text{ and }u\in U.
\end{align*}
This map $f$ is easily seen to be a left $R$-module morphism (since $W$ and
$U$ are $R$-submodules). Hence, by Lemma \ref{lem.spechtmod.irred.R}
\textbf{(b)}, it is either zero or an isomorphism. But it cannot be zero
(since its image is $f\left(  Re\right)  =W\neq0$), and it cannot be an
isomorphism either (since its kernel is $\operatorname*{Ker}f=U\neq0$). This
is a contradiction. Thus, our assumption was false, and Lemma
\ref{lem.spechtmod.irred.R} \textbf{(c)} is proved. \medskip

\textbf{(d)} We recall that $G$-representations are the same as left
$\mathbf{k}\left[  G\right]  $-modules, i.e., as left $R$-modules (since
$\mathbf{k}\left[  G\right]  =R$).

From $e\neq0$, we obtain $Re\neq0$ (since $e=1e\in Re$).

Next, we shall show that the $G$-representation $Re$ has no subrepresentation
$W$ that is distinct from $Re$ and $\left\{  0\right\}  $.

Indeed, assume the contrary. Thus, $Re$ has such a subrepresentation $W$.
Consider this $W$.

Now, $W$ is a $G$-subrepresentation, i.e., a left $\mathbf{k}\left[  G\right]
$-submodule of $Re$. Hence, it is a $\mathbf{k}$-vector subspace of $Re$, and
thus a direct addend of $Re$ as a $\mathbf{k}$-vector space (by Proposition
\ref{prop.mod.diradd.field}, since $\mathbf{k}$ is a field). Moreover,
$\left\vert G\right\vert $ is invertible in $\mathbf{k}$ (since
$\operatorname*{char}\mathbf{k}=0$). Thus, Maschke's theorem (Theorem
\ref{thm.rep.G-rep.maschke}, applied to $V=Re$) shows that the left
$\mathbf{k}\left[  G\right]  $-submodule $W$ is also a direct addend of $Re$
as a left $\mathbf{k}\left[  G\right]  $-module. In other words, $Re=W\oplus
U$ for some further left $\mathbf{k}\left[  G\right]  $-submodule $U$. Now,
this latter submodule $U$ must be nonzero (since otherwise, $Re=W\oplus U$
would become $Re=W$, contradicting the fact that $W$ is distinct from $Re$).
Hence, $Re=W\oplus U$ shows that $Re$ is a direct sum of two nonzero
$R$-submodules. But this contradicts Lemma \ref{lem.spechtmod.irred.R}
\textbf{(c)}.

This contradiction shows that our assumption was wrong. Thus we have proved
that the $G$-representation $Re$ has no subrepresentation $W$ that is distinct
from $Re$ and $\left\{  0\right\}  $. Since $Re\neq0$, this shows that $Re$ is
irreducible. This proves Lemma \ref{lem.spechtmod.irred.R} \textbf{(d)}.
\end{proof}

\begin{corollary}
\label{cor.spechtmod.irred}Assume that $\mathbf{k}$ is a field of
characteristic $0$. Let $\lambda$ be a partition of $n$. Then, the Specht
module $\mathcal{S}^{\lambda}$ is irreducible as a left $\mathbf{k}\left[
S_{n}\right]  $-module (i.e., as a representation of $S_{n}$).
\end{corollary}

\begin{proof}
Recall that $\mathcal{A}=\mathbf{k}\left[  S_{n}\right]  $. Pick any
$n$-tableau $T$ of shape $\lambda$ (clearly, such an $n$-tableau exists, since
$\left\vert Y\left(  \lambda\right)  \right\vert =\left\vert \lambda
\right\vert =n$). Then, (\ref{eq.def.specht.ET.defs.SD=AET}) (applied to
$D=Y\left(  \lambda\right)  $) yields $\mathcal{S}^{Y\left(  \lambda\right)
}\cong\mathcal{A}\mathbf{E}_{T}$ (as left $\mathcal{A}$-modules). In other
words, $\mathcal{S}^{\lambda}\cong\mathcal{A}\mathbf{E}_{T}$ (since
$\mathcal{S}^{\lambda}=\mathcal{S}^{Y\left(  \lambda\right)  }$). We observe
the following:

\begin{itemize}
\item The element $\mathbf{E}_{T}$ is quasi-idempotent: i.e., there is a
nonzero scalar $\kappa\in\mathbf{k}$ such that $\mathbf{E}_{T}^{2}%
=\kappa\mathbf{E}_{T}$. (Indeed, $\kappa=\dfrac{n!}{f^{\lambda}}$ works,
because of Theorem \ref{thm.specht.ETidp}.)

\item The element $\mathbf{E}_{T}$ has the \textquotedblleft sandwich
property\textquotedblright: For each $\mathbf{a}\in\mathcal{A}$, the product
$\mathbf{E}_{T}\mathbf{aE}_{T}$ is a scalar multiple of $\mathbf{E}_{T}$.
(This is just the claim of Proposition \ref{prop.specht.ET.ETaET}.)

\item We have $\mathbf{E}_{T}\neq0$ (since Lemma \ref{lem.specht.ET.1-coord}
yields $\left[  1\right]  \mathbf{E}_{T}=1\neq0$).

\item We have $\mathcal{A}=\mathbf{k}\left[  G\right]  $ for some finite group
$G$ (namely, for $G=S_{n}$).

\item We have $\operatorname*{char}\mathbf{k}=0$.
\end{itemize}

Thus, Lemma \ref{lem.spechtmod.irred.R} \textbf{(d)} (applied to
$R=\mathcal{A}$ and $e=\mathbf{E}_{T}$) shows that the left $\mathcal{A}%
$-module $\mathcal{A}\mathbf{E}_{T}$ is irreducible. In other words, the left
$\mathbf{k}\left[  S_{n}\right]  $-module $\mathcal{A}\mathbf{E}_{T}$ is
irreducible (since $\mathcal{A}=\mathbf{k}\left[  S_{n}\right]  $). Hence, the
left $\mathbf{k}\left[  S_{n}\right]  $-module $\mathcal{S}^{\lambda}$ is
irreducible (since $\mathcal{S}^{\lambda}\cong\mathcal{A}\mathbf{E}_{T}$).
This proves Corollary \ref{cor.spechtmod.irred}.
\end{proof}

This proves Theorem \ref{thm.spechtmod.irred} \textbf{(a)}. \medskip

It is easy to generalize Corollary \ref{cor.spechtmod.irred} to the case when
$\mathbf{k}$ is a field of characteristic larger than $n$. However, we cannot
allow $\mathbf{k}$ to be just a commutative $\mathbb{Q}$-algebra, since
irreducibility becomes a far less well-behaved notion when $\mathbf{k}$ is not
a field (for example, if $\mathbf{k}$ is the polynomial ring $\mathbb{Q}%
\left[  x\right]  $, then any Specht module $\mathcal{S}^{\lambda}$ has
subrepresentations $x\mathcal{S}^{\lambda},\ x^{2}\mathcal{S}^{\lambda
},\ x^{3}\mathcal{S}^{\lambda},\ \ldots$, which are all distinct and nonzero).

A much more interesting question is what happens when $\mathbf{k}$ is a field
of positive characteristic $\leq n$. The Specht modules $\mathcal{S}^{\lambda
}$ are no longer always irreducible in this case, and indeed there is a
blatant counterexample in characteristic $2$ (\cite[Example 4.5]{Wildon18}):

\begin{proposition}
If $2=0$ in $\mathbf{k}$, then%
\[
\mathcal{S}^{\left(  5,1,1\right)  }\cong\mathcal{S}^{\left(  5,2\right)
}\oplus\mathcal{S}^{\left(  7\right)  }\ \ \ \ \ \ \ \ \ \ \text{as left
}\mathbf{k}\left[  S_{7}\right]  \text{-modules.}%
\]

\end{proposition}

\begin{noncompile}
Wildon only proves this over a field of characteristic $2$, but this is true
more generally whenever $2=0$, because we can use base-change from
$\mathbb{F}_{2}$ to any commutative ring $\mathbf{k}$ with $2=0$.
\end{noncompile}

Thus, $\mathcal{S}^{\left(  5,1,1\right)  }$ is not only not irreducible over
a characteristic-$2$ field, but even decomposes as a direct sum of two other
Specht modules! Non-irreducible Specht modules $\mathcal{S}^{\lambda}$ also
exist in other positive characteristics, although nontrivial direct sum
decompositions only exist in characteristic $2$ (see \cite[Corollary
13.18]{James78} for a proof).

Some Specht modules do remain irreducible in positive characteristic. For
example, here is the situation for $n=4$:%
\[%
\begin{tabular}
[c]{c|ccc}
& $\operatorname*{char}\mathbf{k}=2$ & $\operatorname*{char}\mathbf{k}=3$ &
$\operatorname*{char}\mathbf{k}\notin\left\{  2,3\right\}  $\\\hline
$\mathcal{S}^{\left(  4\right)  }$ & irreducible & irreducible & irreducible\\
$\mathcal{S}^{\left(  3,1\right)  }$ & reducible & irreducible & irreducible\\
$\mathcal{S}^{\left(  2,2\right)  }$ & irreducible & reducible & irreducible\\
$\mathcal{S}^{\left(  2,1,1\right)  }$ & reducible & irreducible &
irreducible\\
$\mathcal{S}^{\left(  1,1,1,1\right)  }$ & irreducible & irreducible &
irreducible
\end{tabular}
\ \ .
\]

\begin{exercise}
\fbox{3+3+1} Verify all the claims in the above table. (The claims in the
first and the last rows are trivial, since $1$-dimensional $\mathbf{k}%
$-modules are always irreducible. Each remaining claim in any of the first two
columns is worth 1 point; the whole third column is worth 1 point.)
\end{exercise}

See \cite{JamMat98} for a classification of the Specht modules $\mathcal{S}%
^{\lambda}$ that are irreducible over a characteristic-$2$ field $\mathbf{k}$.

Arguably, for fields of positive characteristic, the proper candidates for the
irreducible $S_{n}$-representations are not the Specht modules $\mathcal{S}%
^{\lambda}$ themselves, but rather their quotients. See \cite[\S 5]{Wildon18}
and \cite[\S 11]{James78} for the details of their construction and
\cite[Chapter 6]{JamKer81} for more about the modular (i.e.,
positive-characteristic) representation theory of $S_{n}$. In a nutshell, the
irreducible $S_{n}$-representations can still be nicely characterized, but
their fine structure (bases, dimensions, etc.) are mostly unknown (e.g., a
\textquotedblleft standard basis theorem\textquotedblright\ is missing).

\subsubsection{Non-isomorphism}

Next we aim at proving Theorem \ref{thm.spechtmod.irred} \textbf{(b)}. Here,
we shall use the following lemma to distinguish $\mathcal{S}^{\lambda}$ from
$\mathcal{S}^{\mu}$:

\begin{lemma}
\label{lem.spechtmod.ET.on-S-mu}Let $\lambda$ be a partition of $n$. Let $T$
be an $n$-tableau of shape $\lambda$. As usual in algebra, we shall write
$\mathbf{a}M$ for the set $\left\{  \mathbf{a}m\ \mid\ m\in M\right\}  $ when
$\mathbf{a}$ is an element of $\mathbf{k}\left[  S_{n}\right]  $ and when $M$
is a left $\mathbf{k}\left[  S_{n}\right]  $-module. Then: \medskip

\textbf{(a)} If $h^{\lambda}\neq0$ in $\mathbf{k}$, then $\mathbf{E}%
_{T}\mathcal{S}^{\lambda}\neq0$. \medskip

\textbf{(b)} If $\mu$ is any partition of $n$ that is distinct from $\lambda$,
then $\mathbf{E}_{T}\mathcal{S}^{\mu}=0$.
\end{lemma}

\begin{proof}
\textbf{(a)} Assume that $h^{\lambda}\neq0$ in $\mathbf{k}$. Then,
\begin{align*}
\mathbf{E}_{T}1\mathbf{E}_{T}  &  =\mathbf{E}_{T}^{2}=\dfrac{n!}{f^{\lambda}%
}\mathbf{E}_{T}\ \ \ \ \ \ \ \ \ \ \left(  \text{by Theorem
\ref{thm.specht.ETidp}}\right) \\
&  =h^{\lambda}\mathbf{E}_{T}\ \ \ \ \ \ \ \ \ \ \left(  \text{since }%
\dfrac{n!}{f^{\lambda}}=h^{\lambda}\right)  ,
\end{align*}
so that $\left[  1\right]  \left(  \mathbf{E}_{T}1\mathbf{E}_{T}\right)
=\left[  1\right]  \left(  h^{\lambda}\mathbf{E}_{T}\right)  =h^{\lambda}%
\cdot\underbrace{\left[  1\right]  \mathbf{E}_{T}}_{\substack{=1\\\text{(by
Lemma \ref{lem.specht.ET.1-coord})}}}=h^{\lambda}\neq0$. Therefore,
$\mathbf{E}_{T}1\mathbf{E}_{T}\neq0$. Hence, $\mathbf{E}_{T}\mathcal{A}%
\mathbf{E}_{T}\neq0$ (since $\mathbf{E}_{T}1\mathbf{E}_{T}$ is clearly an
element of $\mathbf{E}_{T}\mathcal{A}\mathbf{E}_{T}$).

But (\ref{eq.def.specht.ET.defs.SD=AET}) (applied to $D=Y\left(
\lambda\right)  $) yields $\mathcal{S}^{Y\left(  \lambda\right)  }%
\cong\mathcal{A}\mathbf{E}_{T}$ (as left $\mathcal{A}$-modules). In other
words, $\mathcal{S}^{\lambda}\cong\mathcal{A}\mathbf{E}_{T}$ (since
$\mathcal{S}^{\lambda}=\mathcal{S}^{Y\left(  \lambda\right)  }$). Hence,
$\mathbf{E}_{T}\mathcal{S}^{\lambda}\cong\mathbf{E}_{T}\mathcal{A}%
\mathbf{E}_{T}$ (since the action of $\mathbf{E}_{T}$ on isomorphic left
$\mathbf{k}\left[  S_{n}\right]  $-modules is isomorphic). Therefore,
$\mathbf{E}_{T}\mathcal{S}^{\lambda}\neq0$ (since $\mathbf{E}_{T}%
\mathcal{A}\mathbf{E}_{T}\neq0$). This proves Lemma
\ref{lem.spechtmod.ET.on-S-mu} \textbf{(a)}. \medskip

\textbf{(b)} Let $\mu$ be any partition of $n$ that is distinct from $\lambda
$. Pick any $n$-tableau $S$ of shape $\mu$ (such an $n$-tableau clearly
exists). Then, $S$ is an $n$-tableau of shape $Y\left(  \mu\right)  $, whereas
$T$ is an $n$-tableau of shape $Y\left(  \lambda\right)  $. Therefore,
Proposition \ref{prop.specht.ET.ESaET} (applied to $T$ and $S$ instead of $S$
and $T$) yields that
\[
\mathbf{E}_{T}\mathbf{aE}_{S}=0\ \ \ \ \ \ \ \ \ \ \text{for each }%
\mathbf{a}\in\mathbf{k}\left[  S_{n}\right]  .
\]
In other words, $\mathbf{E}_{T}\mathcal{A}\mathbf{E}_{S}=0$.

However, in our above proof of Lemma \ref{lem.spechtmod.ET.on-S-mu}
\textbf{(a)}, we have shown that $\mathcal{S}^{\lambda}\cong\mathcal{A}%
\mathbf{E}_{T}$ as left $\mathcal{A}$-modules. Similarly, we can see that
$\mathcal{S}^{\mu}\cong\mathcal{A}\mathbf{E}_{S}$ as left $\mathcal{A}%
$-modules (since $S$ is an $n$-tableau of shape $\mu$). Therefore,
$\mathbf{E}_{T}\mathcal{S}^{\mu}\cong\mathbf{E}_{T}\mathcal{A}\mathbf{E}%
_{S}=0$, so that $\mathbf{E}_{T}\mathcal{S}^{\mu}=0$. This proves Lemma
\ref{lem.spechtmod.ET.on-S-mu} \textbf{(b).}
\end{proof}

\begin{corollary}
\label{cor.spechtmod.noniso}Assume that $\mathbf{k}$ is a field of
characteristic $0$. Let $\lambda$ and $\mu$ be two distinct partitions of $n$.
Then, $\mathcal{S}^{\lambda}\ncong\mathcal{S}^{\mu}$ as $S_{n}$-representations.
\end{corollary}

\begin{proof}
[First proof of Corollary \ref{cor.spechtmod.noniso}.]Assume the contrary.
Thus, we have $\mathcal{S}^{\lambda}\cong\mathcal{S}^{\mu}$ as $S_{n}$-representations.

Note that $h^{\lambda}\neq0$ in $\mathbf{k}$ (since $\mathbf{k}$ is a field of
characteristic $0$).

Pick any $n$-tableau $T$ of shape $\lambda$ (such a $T$ clearly exists). Then,
Lemma \ref{lem.spechtmod.ET.on-S-mu} \textbf{(a)} yields $\mathbf{E}%
_{T}\mathcal{S}^{\lambda}\neq0$. But Lemma \ref{lem.spechtmod.ET.on-S-mu}
\textbf{(b)} yields $\mathbf{E}_{T}\mathcal{S}^{\mu}=0$. However, from
$\mathcal{S}^{\lambda}\cong\mathcal{S}^{\mu}$, we obtain $\mathbf{E}%
_{T}\mathcal{S}^{\lambda}\cong\mathbf{E}_{T}\mathcal{S}^{\mu}=0$,
contradicting $\mathbf{E}_{T}\mathcal{S}^{\lambda}\neq0$. This contradiction
shows that our assumption was wrong. Thus, Corollary
\ref{cor.spechtmod.noniso} is proved.
\end{proof}

\begin{proof}
[Second proof of Corollary \ref{cor.spechtmod.noniso} (sketched).]Let me
sketch a different proof of Corollary \ref{cor.spechtmod.noniso} (coming from
\cite[Corollary 13.17]{James78}), which establishes a more general result:
Namely, it shows that Corollary \ref{cor.spechtmod.noniso} still holds
whenever $\mathbf{k}$ is a field of characteristic $\neq2$ (rather than
characteristic $0$). Even better, it works whenever $\mathbf{k}$ is any
nontrivial commutative ring in which $2$ is invertible.

So let $\mathbf{k}$ be a nontrivial commutative ring in which $2$ is
invertible. Let $\lambda$ and $\mu$ be two distinct partitions of $n$.

Consider the \textbf{smallest} $i\geq1$ satisfying $\lambda_{i}\neq\mu_{i}$.
Assume WLOG that $\lambda_{i}<\mu_{i}$ (otherwise, we can swap $\lambda$ with
$\mu$). Then,
\begin{equation}
\lambda_{1}+\lambda_{2}+\cdots+\lambda_{i}<\mu_{1}+\mu_{2}+\cdots+\mu_{i}
\label{pf.cor.spechtmod.noniso.ineq}%
\end{equation}
(as we saw in the proof of Proposition \ref{prop.specht.ET.ESaET}, except that
we have the $<$ sign instead of the $>$ sign).

We must show that $\mathcal{S}^{\lambda}\ncong\mathcal{S}^{\mu}$ as $S_{n}$-representations.

Indeed, assume the contrary. Thus, there is a left $\mathbf{k}\left[
S_{n}\right]  $-module isomorphism $f:\mathcal{S}^{\lambda}\rightarrow
\mathcal{S}^{\mu}$. Consider this $f$.

Pick any $n$-tableau $T$ of shape $\lambda$. Write the image $f\left(
\mathbf{e}_{T}\right)  \in\mathcal{S}^{\mu}=\mathcal{S}^{Y\left(  \mu\right)
}\subseteq\mathcal{M}^{Y\left(  \mu\right)  }$ as
\begin{equation}
f\left(  \mathbf{e}_{T}\right)  =\sum_{\overline{S}}a_{\overline{S}}%
\overline{S}, \label{pf.cor.spechtmod.noniso.feT=}%
\end{equation}
where the sum ranges over all $n$-tabloids $\overline{S}$ of shape $Y\left(
\mu\right)  $, and where the coefficients $a_{\overline{S}}$ are scalars in
$\mathbf{k}$. For any $w\in S_{n}$, we have%
\begin{align}
f\left(  \mathbf{e}_{wT}\right)   &  =f\left(  w\mathbf{e}_{T}\right)
\ \ \ \ \ \ \ \ \ \ \left(  \text{since Lemma \ref{lem.spechtmod.submod}
\textbf{(a)} yields }\mathbf{e}_{wT}=w\mathbf{e}_{T}\right) \nonumber\\
&  =wf\left(  \mathbf{e}_{T}\right)  \ \ \ \ \ \ \ \ \ \ \left(  \text{since
}f\text{ is a left }\mathbf{k}\left[  S_{n}\right]  \text{-module
morphism}\right) \nonumber\\
&  =w\sum_{\overline{S}}a_{\overline{S}}\overline{S}%
\ \ \ \ \ \ \ \ \ \ \left(  \text{by (\ref{pf.cor.spechtmod.noniso.feT=}%
)}\right) \nonumber\\
&  =\sum_{\overline{S}}a_{\overline{S}}w\overline{S}=\sum_{\overline{S}%
}a_{w^{-1}\overline{S}}\overline{S} \label{pf.cor.spechtmod.noniso.2nd.1}%
\end{align}
(here, we have substituted $w^{-1}\overline{S}$ for $\overline{S}$ in the sum).

Now, let $w\in S_{n}$ be a transposition $t_{u,v}$ swapping two distinct
numbers $u,v\in\left[  n\right]  $ that lie in the same column of $T$. Then,
$w\in\mathcal{C}\left(  T\right)  $, so that Lemma \ref{lem.spechtmod.submod}
\textbf{(c)} yields $\mathbf{e}_{wT}=\left(  -1\right)  ^{w}\mathbf{e}%
_{T}=-\mathbf{e}_{T}$ (since the transposition $w$ has sign $\left(
-1\right)  ^{w}=-1$), and thus we have%
\[
f\left(  \mathbf{e}_{wT}\right)  =f\left(  -\mathbf{e}_{T}\right)  =-f\left(
\mathbf{e}_{T}\right)  =-\sum_{\overline{S}}a_{\overline{S}}\overline
{S}\ \ \ \ \ \ \ \ \ \ \left(  \text{by (\ref{pf.cor.spechtmod.noniso.feT=}%
)}\right)  .
\]
Comparing this with (\ref{pf.cor.spechtmod.noniso.2nd.1}), we obtain%
\[
\sum_{\overline{S}}a_{w^{-1}\overline{S}}\overline{S}=-\sum_{\overline{S}%
}a_{\overline{S}}\overline{S}=\sum_{\overline{S}}\left(  -a_{\overline{S}%
}\right)  \overline{S}.
\]
Since the $n$-tabloids $\overline{S}$ form a basis of $\mathcal{M}^{Y\left(
\mu\right)  }$, we thus obtain%
\begin{equation}
a_{w^{-1}\overline{S}}=-a_{\overline{S}}\ \ \ \ \ \ \ \ \ \ \text{for each
}n\text{-tabloid }\overline{S}\text{ of shape }Y\left(  \mu\right)  .
\label{pf.cor.spechtmod.noniso.2nd.2}%
\end{equation}

Forget that we fixed $w$. We thus have shown that if $w\in S_{n}$ is a
transposition $t_{u,v}$ swapping two distinct numbers $u,v\in\left[  n\right]
$ that lie in the same column of $T$, then
(\ref{pf.cor.spechtmod.noniso.2nd.2}) holds.

Now, let $S$ be any $n$-tableau of shape $Y\left(  \mu\right)  $. Then, by
Theorem \ref{thm.youngtab.alt} \textbf{(b)} (or, rather, by its
contrapositive), we can see that there are two distinct integers $u$ and $v$
that are horizontally aligned in $S$ and simultaneously vertically aligned in
$T$ (since otherwise, (\ref{pf.cor.spechtmod.noniso.ineq}) would be violated).
Consider these $u$ and $v$. Thus, $u$ and $v$ lie in the same row of $S$;
moreover, $u$ and $v$ lie in the same column of $T$.

Let $w$ be the transposition $t_{u,v}\in S_{n}$. Then,
(\ref{pf.cor.spechtmod.noniso.2nd.2}) shows that $a_{w^{-1}\overline{S}%
}=-a_{\overline{S}}$. But $w=t_{u,v}\in\mathcal{R}\left(  S\right)  $ (since
$u$ and $v$ lie in the same row of $S$), so that $wS$ is row-equivalent to
$S$. Therefore, $\overline{wS}=\overline{S}$. Thus, $w\overline{S}%
=\overline{wS}=\overline{S}$, so that $w^{-1}\overline{S}=\overline{S}$ as
well. Hence, the equality $a_{w^{-1}\overline{S}}=-a_{\overline{S}}$ can be
rewritten as $a_{\overline{S}}=-a_{\overline{S}}$. In other words,
$2a_{\overline{S}}=0$. Thus, $a_{\overline{S}}=0$ (since $2$ is invertible in
$\mathbf{k}$).

Forget that we fixed $S$. We thus have shown that $a_{\overline{S}}=0$ for
each $n$-tableau $S$ of shape $Y\left(  \mu\right)  $. Hence,
(\ref{pf.cor.spechtmod.noniso.feT=}) rewrites as $f\left(  \mathbf{e}%
_{T}\right)  =\sum_{\overline{S}}0\overline{S}=0$.

Forget that we fixed $T$. We have now shown that $f\left(  \mathbf{e}%
_{T}\right)  =0$ for each $n$-tableau $T$ of shape $\lambda$. In other words,
the $\mathbf{k}$-linear map $f$ sends each polytabloid $\mathbf{e}_{T}$ to
$0$. Since the Specht module $\mathcal{S}^{\lambda}$ is spanned by the
polytabloids $\mathbf{e}_{T}$, this entails that $f$ sends every vector to
$0$. In other words, $f=0$. Hence, $f\left(  \mathcal{S}^{\lambda}\right)
=0$. But $f$ is an isomorphism, and thus $f\left(  \mathcal{S}^{\lambda
}\right)  =\mathcal{S}^{\mu}$. Thus, $\mathcal{S}^{\mu}=f\left(
\mathcal{S}^{\lambda}\right)  =0$. But Corollary \ref{cor.spechtmod.nonzero}
shows that $\mathcal{S}^{\mu}\neq0$. This clearly contradicts $\mathcal{S}%
^{\mu}=0$. This contradiction shows that our assumption was false, and thus
Corollary \ref{cor.spechtmod.noniso} is proved again.
\end{proof}

Note that Corollary \ref{cor.spechtmod.noniso} fails in characteristic $2$,
since we have%
\[
\mathcal{S}^{\left(  2\right)  }\cong\mathcal{S}^{\left(  1,1\right)
}\ \ \ \ \ \ \ \ \ \ \text{when }\mathbf{k}\text{ is a field of characteristic
}2
\]
(for fairly obvious reasons).

Of course, Theorem \ref{thm.spechtmod.irred} \textbf{(b)} follows from
Corollary \ref{cor.spechtmod.noniso} (in fact, it is just a restatement of the latter).

\begin{exercise}
\fbox{2} Tweak the above second proof of Corollary \ref{cor.spechtmod.noniso}
a bit further to work under the even weaker assumption \textquotedblleft%
$\mathbf{k}$ is a commutative ring such that $2\neq0$ in $\mathbf{k}%
$\textquotedblright\ instead of \textquotedblleft$\mathbf{k}$ is a field of
characteristic $0$\textquotedblright.
\end{exercise}

\subsubsection{Completeness}

To prove Theorem \ref{thm.spechtmod.irred} \textbf{(c)}, we need the following
general fact from representation theory:

\begin{proposition}
\label{prop.spechtmod.complete.sum}Let $\mathbf{k}$ be a field. Let $R$ be a
$\mathbf{k}$-algebra that is finite-dimensional as a $\mathbf{k}$-vector
space. Let $J_{1},J_{2},\ldots,J_{m}$ be $m$ irreducible left $R$-submodules
of the left regular $R$-module $R$ satisfying $R=J_{1}+J_{2}+\cdots+J_{m}$.
Let $K$ be a further irreducible left $R$-module. Then, $K\cong J_{i}$ for
some $i\in\left[  m\right]  $.
\end{proposition}

\begin{proof}
[Proof sketch.]If some $i\in\left[  m\right]  $ satisfies $J_{1}+J_{2}%
+\cdots+J_{i}=J_{1}+J_{2}+\cdots+J_{i-1}$, then we have%
\begin{align*}
R  &  =J_{1}+J_{2}+\cdots+J_{m}\\
&  =\underbrace{\left(  J_{1}+J_{2}+\cdots+J_{i}\right)  }_{=J_{1}%
+J_{2}+\cdots+J_{i-1}}+\left(  J_{i+1}+J_{i+2}+\cdots+J_{m}\right) \\
&  =\left(  J_{1}+J_{2}+\cdots+J_{i-1}\right)  +\left(  J_{i+1}+J_{i+2}%
+\cdots+J_{m}\right)  ,
\end{align*}
and therefore the equality $R=J_{1}+J_{2}+\cdots+J_{m}$ remains valid if we
remove the addend $J_{i}$ from its right hand side. Thus, in this case, we can
reduce the claim of Proposition \ref{prop.spechtmod.complete.sum} to an
analogous claim with a smaller value of $m$ (by removing $J_{i}$ from the list
$\left(  J_{1},J_{2},\ldots,J_{m}\right)  $ and replacing $m$ by $m-1$).

We can perform such reductions until there is no longer any $i\in\left[
m\right]  $ remaining that satisfies $J_{1}+J_{2}+\cdots+J_{i}=J_{1}%
+J_{2}+\cdots+J_{i-1}$. Hence, we can WLOG assume that no $i\in\left[
m\right]  $ satisfies $J_{1}+J_{2}+\cdots+J_{i}=J_{1}+J_{2}+\cdots+J_{i-1}$.
Assume this. Thus,%
\begin{equation}
J_{1}+J_{2}+\cdots+J_{i}\neq J_{1}+J_{2}+\cdots+J_{i-1}
\label{pf.prop.spechtmod.complete.sum.WLOG}%
\end{equation}
for each $i\in\left[  m\right]  $.

Now, for each $i\in\left\{  0,1,\ldots,m\right\}  $, let us define the left
$R$-submodule
\[
F_{i}:=J_{1}+J_{2}+\cdots+J_{i}\text{ of }R.
\]
Thus, $F_{0}=\left\{  0\right\}  $ (since an empty sum of submodules is
$\left\{  0\right\}  $ by definition) and $F_{m}=J_{1}+J_{2}+\cdots+J_{m}=R$.
Moreover, each $i\in\left[  m\right]  $ satisfies%
\begin{align}
F_{i}  &  =J_{1}+J_{2}+\cdots+J_{i}=\underbrace{\left(  J_{1}+J_{2}%
+\cdots+J_{i-1}\right)  }_{\substack{=F_{i-1}\\\text{(by the definition of
}F_{i-1}\text{)}}}+\,J_{i}\nonumber\\
&  =F_{i-1}+J_{i} \label{pf.prop.spechtmod.complete.sum.Fi=+}%
\end{align}
and therefore $F_{i-1}\subseteq F_{i}$. Hence, $\left\{  0\right\}
=F_{0}\subseteq F_{1}\subseteq F_{2}\subseteq\cdots\subseteq F_{m}=R$. This
shows that $\left(  F_{0},F_{1},\ldots,F_{m}\right)  $ is a filtration of the
left $R$-module $R$.

Next, let us show that the subquotients $F_{i}/F_{i-1}$ of this filtration are
isomorphic to the respective $J_{i}$:

\begin{statement}
\textit{Claim 1:} Let $i\in\left[  m\right]  $. Then, $F_{i}/F_{i-1}\cong
J_{i}$ as left $R$-modules.
\end{statement}

\begin{proof}
[Proof of Claim 1.]First of all, we have $F_{i-1}=J_{1}+J_{2}+\cdots+J_{i-1}$
(by the definition of $F_{i-1}$). But
(\ref{pf.prop.spechtmod.complete.sum.WLOG}) says that $J_{1}+J_{2}%
+\cdots+J_{i}\neq J_{1}+J_{2}+\cdots+J_{i-1}$. In view of $F_{i}=J_{1}%
+J_{2}+\cdots+J_{i}$ and $F_{i-1}=J_{1}+J_{2}+\cdots+J_{i-1}$, we can rewrite
this as $F_{i}\neq F_{i-1}$. Thus, $F_{i}/F_{i-1}\neq\left\{  0\right\}  $.

From (\ref{pf.prop.spechtmod.complete.sum.Fi=+}), we have $F_{i}=F_{i-1}%
+J_{i}$. Hence, $J_{i}\subseteq F_{i}$. Thus, there is a left $R$-linear map%
\begin{align*}
g:J_{i}  &  \rightarrow F_{i}/F_{i-1},\\
x  &  \mapsto\overline{x}%
\end{align*}
(where $\overline{x}$ denotes the residue class of $x\in F_{i}$ modulo
$F_{i-1}$). Consider this map $g$.

Let us show that $g$ is surjective. Indeed, let $v\in F_{i}/F_{i-1}$ be
arbitrary. Then, we can write $v$ as $v=\overline{y}$ for some $y\in F_{i}$
(obviously). Consider this $y$. Then, $y\in F_{i}=F_{i-1}+J_{i}$, so that we
can write $y$ as $y=z+w$ for some $z\in F_{i-1}$ and $w\in J_{i}$. Consider
these $z$ and $w$. From $y=z+w$, we obtain $y-w=z\in F_{i-1}$ and therefore
$\overline{y}=\overline{w}$. Hence, $v=\overline{y}=\overline{w}=g\left(
w\right)  $ (since the definition of $g$ yields $g\left(  w\right)
=\overline{w}$). Therefore, $v=g\left(  w\right)  \in\operatorname{Im}g$.

Forget that we fixed $v$. We thus have shown that $v\in\operatorname{Im}g$ for
each $v\in F_{i}/F_{i-1}$. In other words, the map $g$ is surjective. That is,
$g\left(  J_{i}\right)  =F_{i}/F_{i-1}\neq\left\{  0\right\}  $.

The kernel $\operatorname*{Ker}g$ of the map $g:J_{i}\rightarrow F_{i}%
/F_{i-1}$ is a left $R$-submodule of $J_{i}$ (since $g$ is left $R$-linear),
and thus must be either $J_{i}$ or $\left\{  0\right\}  $ (since the
$R$-module $J_{i}$ is irreducible, and thus its only $R$-submodules are
$J_{i}$ and $\left\{  0\right\}  $). If $\operatorname*{Ker}g$ was $J_{i}$,
then $g\left(  J_{i}\right)  $ would be $\left\{  0\right\}  $, which would
contradict $g\left(  J_{i}\right)  \neq\left\{  0\right\}  $. Hence,
$\operatorname*{Ker}g$ must be $\left\{  0\right\}  $ (since
$\operatorname*{Ker}g$ is either $J_{i}$ or $\left\{  0\right\}  $). In other
words, $g$ is injective (since $g$ is a left $R$-linear map). Since $g$ is
also surjective, we thus conclude that $g$ is bijective, i.e., invertible.
Hence, $g$ is a left $R$-module isomorphism (since $g$ is left $R$-linear).
Thus, $F_{i}/F_{i-1}\cong J_{i}$ as left $R$-modules. This proves Claim 1.
\end{proof}

Claim 1 shows that each subquotient $F_{i}/F_{i-1}$ of our filtration $\left(
F_{0},F_{1},\ldots,F_{m}\right)  $ is isomorphic to the corresponding $J_{i}$,
and thus is irreducible (since $J_{i}$ is irreducible). Hence, the filtration
$\left(  F_{0},F_{1},\ldots,F_{m}\right)  $ is a composition series of $R$.
Thus, Theorem \ref{thm.mod.filt.irrep-in-filt} (applied to $m$ and $\left(
F_{0},F_{1},\ldots,F_{m}\right)  $ instead of $k$ and $\left(  W_{0}%
,W_{1},\ldots,W_{k}\right)  $) shows that each irreducible left $R$-module is
isomorphic to $F_{i}/F_{i-1}$ for some $i\in\left[  m\right]  $. In
particular, $K$ is isomorphic to $F_{i}/F_{i-1}$ for some $i\in\left[
m\right]  $ (since $K$ is an irreducible left $R$-module). For this $i$, we
then have $K\cong F_{i}/F_{i-1}\cong J_{i}$ (by Claim 1). This proves
Proposition \ref{prop.spechtmod.complete.sum}.
\end{proof}

\begin{corollary}
\label{cor.spechtmod.complete}Assume that $\mathbf{k}$ is a field of
characteristic $0$. Let $K$ be any irreducible representation of $S_{n}$ over
$\mathbf{k}$. Then, $K\cong\mathcal{S}^{\lambda}$ for some partition $\lambda$
of $n$.
\end{corollary}

\begin{proof}
Recall that we denote the ring $\mathbf{k}\left[  S_{n}\right]  $ by
$\mathcal{A}$. Thus, $K$ is an irreducible left $\mathcal{A}$-module (since
$K$ is an irreducible representation of $S_{n}$ over $\mathbf{k}$).

Let $\left(  T_{1},T_{2},\ldots,T_{m}\right)  $ be a list of \textbf{all} the
$n$-tableaux of shape $\lambda$ for \textbf{all} the partitions $\lambda$ of
$n$. (Thus, the length $m$ of this list is $n!\cdot p\left(  n\right)  $,
where $p\left(  n\right)  $ is the \# of all partitions of $n$.)

Theorem \ref{thm.spechtmod.Elam.pou} yields%
\begin{align}
1  &  =\sum_{\lambda\text{ is a partition of }n}\dfrac{\mathbf{E}_{\lambda}%
}{\left(  h^{\lambda}\right)  ^{2}}\nonumber\\
&  =\sum_{\lambda\text{ is a partition of }n}\dfrac{\sum_{T\text{ is an
}n\text{-tableau of shape }\lambda}\mathbf{E}_{T}}{\left(  h^{\lambda}\right)
^{2}}\ \ \ \ \ \ \ \ \ \ \left(  \text{by
(\ref{eq.def.spechtmod.Elam.Elam.def})}\right) \nonumber\\
&  =\sum_{\lambda\text{ is a partition of }n}\ \ \sum_{T\text{ is an
}n\text{-tableau of shape }\lambda}\dfrac{\mathbf{E}_{T}}{\left(  h^{\lambda
}\right)  ^{2}}. \label{pf.cor.spechtmod.complete.1}%
\end{align}
Thus, each $\mathbf{a}\in\mathcal{A}$ satisfies%
\begin{align*}
\mathbf{a}  &  =\mathbf{a}\cdot1=\mathbf{a}\cdot\sum_{\lambda\text{ is a
partition of }n}\ \ \sum_{T\text{ is an }n\text{-tableau of shape }\lambda
}\dfrac{\mathbf{E}_{T}}{\left(  h^{\lambda}\right)  ^{2}}%
\ \ \ \ \ \ \ \ \ \ \left(  \text{by (\ref{pf.cor.spechtmod.complete.1}%
)}\right) \\
&  =\sum_{\lambda\text{ is a partition of }n}\ \ \sum_{T\text{ is an
}n\text{-tableau of shape }\lambda}\underbrace{\dfrac{\mathbf{aE}_{T}}{\left(
h^{\lambda}\right)  ^{2}}}_{\in\mathcal{A}\mathbf{E}_{T}}\\
&  \in\sum_{\lambda\text{ is a partition of }n}\ \ \sum_{T\text{ is an
}n\text{-tableau of shape }\lambda}\mathcal{A}\mathbf{E}_{T}\\
&  =\mathcal{A}\mathbf{E}_{T_{1}}+\mathcal{A}\mathbf{E}_{T_{2}}+\cdots
+\mathcal{A}\mathbf{E}_{T_{m}}%
\end{align*}
(since $T_{1},T_{2},\ldots,T_{m}$ are all the $n$-tableaux of shape $\lambda$
for all the partitions $\lambda$ of $n$). Hence,%
\[
\mathcal{A}\subseteq\mathcal{A}\mathbf{E}_{T_{1}}+\mathcal{A}\mathbf{E}%
_{T_{2}}+\cdots+\mathcal{A}\mathbf{E}_{T_{m}}.
\]
This entails%
\[
\mathcal{A}=\mathcal{A}\mathbf{E}_{T_{1}}+\mathcal{A}\mathbf{E}_{T_{2}}%
+\cdots+\mathcal{A}\mathbf{E}_{T_{m}}%
\]
(since we obviously have $\mathcal{A}\mathbf{E}_{T_{1}}+\mathcal{A}%
\mathbf{E}_{T_{2}}+\cdots+\mathcal{A}\mathbf{E}_{T_{m}}\subseteq\mathcal{A}$
as well). Moreover, the left $\mathcal{A}$-submodules $\mathcal{A}%
\mathbf{E}_{T_{i}}$ of $\mathcal{A}$ are isomorphic to Specht modules of the
form $\mathcal{S}^{\lambda}$ (by (\ref{eq.def.specht.ET.defs.SD=AET}), since
$T_{i}$ is an $n$-tableau of shape $\lambda$ for some partition $\lambda$ of
$n$), and thus irreducible (since Corollary \ref{cor.spechtmod.irred} says
that any Specht module of the form $\mathcal{S}^{\lambda}$ is irreducible).
Hence, Proposition \ref{prop.spechtmod.complete.sum} (applied to
$R=\mathcal{A}$ and $J_{i}=\mathcal{A}\mathbf{E}_{T_{i}}$) yields that
$K\cong\mathcal{A}\mathbf{E}_{T_{i}}$ for some $i\in\left[  m\right]  $.
Therefore, $K\cong\mathcal{S}^{\lambda}$ for some partition $\lambda$ of $n$
(since, as we just said, each $\mathcal{A}\mathbf{E}_{T_{i}}$ is isomorphic to
a Specht module of the form $\mathcal{S}^{\lambda}$). This proves Corollary
\ref{cor.spechtmod.complete}.
\end{proof}

Thus, Theorem \ref{thm.spechtmod.irred} \textbf{(c)} follows (being a
restatement of Corollary \ref{cor.spechtmod.complete}). \medskip

Furthermore, we obtain the following corollary, which says that each
finite-dimensional representation of $S_{n}$ over a field of characteristic
$0$ can be decomposed into a direct sum of Specht modules $\mathcal{S}%
^{\lambda}$:

\begin{corollary}
\label{cor.spechtmod.complete2}Assume that $\mathbf{k}$ is a field of
characteristic $0$. Let $P$ be the set of all partitions of $n$. Let $V$ be
any finite-dimensional representation of $S_{n}$ over $\mathbf{k}$. Then,%
\[
V\cong\bigoplus_{\lambda\in P}\left(  \mathcal{S}^{\lambda}\right)
^{k_{\lambda}}\ \ \ \ \ \ \ \ \ \ \text{for some integers }k_{\lambda}%
\in\mathbb{N}.
\]
Here, $M^{k}$ means the external direct sum $\underbrace{M\oplus M\oplus
\cdots\oplus M}_{k\text{ times}}$ whenever $M$ is an $S_{n}$-representation
and $k\in\mathbb{N}$.
\end{corollary}

\begin{proof}
Theorem \ref{thm.rep.G-rep.JH} \textbf{(a)} (applied to $R=\mathbf{k}\left[
S_{n}\right]  $) yields that there exists a composition series $\left(
W_{0},W_{1},\ldots,W_{k}\right)  $ of $V$. Consider this composition series.
For each $i\in\left[  k\right]  $, the quotient $W_{i}/W_{i-1}$ is an
irreducible left $\mathbf{k}\left[  S_{n}\right]  $-module (by the definition
of a composition series), i.e., an irreducible representation of $S_{n}$, and
thus satisfies $W_{i}/W_{i-1}\cong\mathcal{S}^{\lambda}$ for some partition
$\lambda$ of $n$ (by Corollary \ref{cor.spechtmod.complete}, applied to
$K=W_{i}/W_{i-1}$). Let us denote this partition $\lambda$ by $\lambda_{i}$.
Thus,%
\begin{equation}
W_{i}/W_{i-1}\cong\mathcal{S}^{\lambda_{i}}\ \ \ \ \ \ \ \ \ \ \text{for each
}i\in\left[  k\right]  . \label{pf.cor.spechtmod.complete2.quot}%
\end{equation}
Clearly, $\lambda_{i}\in P$ for each $i\in\left[  k\right]  $ (since
$\lambda_{i}$ is a partition of $n$).

But $\left(  W_{0},W_{1},\ldots,W_{k}\right)  $ is a composition series of
$V$, thus a filtration of $V$. Hence, Corollary
\ref{cor.rep.G-rep.maschke-filt} (applied to $G=S_{n}$) yields
\begin{align*}
V  &  \cong\left(  W_{1}/W_{0}\right)  \oplus\left(  W_{2}/W_{1}\right)
\oplus\cdots\oplus\left(  W_{k}/W_{k-1}\right) \\
&  =\bigoplus_{i\in\left[  k\right]  }\underbrace{\left(  W_{i}/W_{i-1}%
\right)  }_{\substack{\cong\mathcal{S}^{\lambda_{i}}\\\text{(by
(\ref{pf.cor.spechtmod.complete2.quot}))}}}\cong\bigoplus_{i\in\left[
k\right]  }\mathcal{S}^{\lambda_{i}}\cong\bigoplus_{\lambda\in P}\left(
\mathcal{S}^{\lambda}\right)  ^{\left(  \text{\# of }i\in\left[  k\right]
\text{ satisfying }\lambda_{i}=\lambda\right)  }%
\end{align*}
(here, we have rewritten our direct sum by collecting equal addends together
-- for instance, $\mathcal{S}^{\left(  4\right)  }\oplus\mathcal{S}^{\left(
3,1\right)  }\oplus\mathcal{S}^{\left(  4\right)  }\oplus\mathcal{S}^{\left(
4\right)  }=\left(  \mathcal{S}^{\left(  4\right)  }\right)  ^{3}\oplus\left(
\mathcal{S}^{\left(  3,1\right)  }\right)  ^{1}$). Hence,%
\[
V\cong\bigoplus_{\lambda\in P}\left(  \mathcal{S}^{\lambda}\right)
^{k_{\lambda}}\ \ \ \ \ \ \ \ \ \ \text{for some integers }k_{\lambda}%
\in\mathbb{N}%
\]
(namely, for $k_{\lambda}=\left(  \text{\# of }i\in\left[  k\right]  \text{
satisfying }\lambda_{i}=\lambda\right)  $). This proves Corollary
\ref{cor.spechtmod.complete2}.
\end{proof}

\subsection{The Artin--Wedderburn theorem}

\subsubsection{The theorem}

Another crucial fact that we can now prove is the \emph{Artin--Wedderburn
theorem for symmetric groups}, stated here in its coordinate-free version
(i.e., using endomorphism rings instead of matrix rings):

\begin{theorem}
\label{thm.specht.AW}Assume that $n!$ is invertible in $\mathbf{k}$. For each
partition $\lambda$ of $n$, we let $\rho_{\lambda}:\mathbf{k}\left[
S_{n}\right]  \rightarrow\operatorname*{End}\nolimits_{\mathbf{k}}\left(
\mathcal{S}^{\lambda}\right)  $ be the curried form of the left $\mathbf{k}%
\left[  S_{n}\right]  $-action on the Specht module $\mathcal{S}^{\lambda}$.
(See Definition \ref{def.mod.curry-into-EndV} for the meaning of
\textquotedblleft curried form\textquotedblright. Explicitly, $\rho_{\lambda}$
is given by $\left(  \rho_{\lambda}\left(  \mathbf{a}\right)  \right)  \left(
\mathbf{v}\right)  =\mathbf{av}$ for each $\mathbf{a}\in\mathbf{k}\left[
S_{n}\right]  $ and each $\mathbf{v}\in\mathcal{S}^{\lambda}$.)

Then, the map%
\begin{align*}
\mathbf{k}\left[  S_{n}\right]   &  \rightarrow\underbrace{\prod
_{\lambda\text{ is a partition of }n}\operatorname*{End}\nolimits_{\mathbf{k}%
}\left(  \mathcal{S}^{\lambda}\right)  }_{\text{a direct product of
}\mathbf{k}\text{-algebras}},\\
\mathbf{a}  &  \mapsto\left(  \rho_{\lambda}\left(  \mathbf{a}\right)
\right)  _{\lambda\text{ is a partition of }n}%
\end{align*}
is a $\mathbf{k}$-algebra isomorphism.
\end{theorem}

\subsubsection{Linear algebra lemmas}

To prove this theorem, we recall a few basic results from linear algebra over
commutative rings. The first one allows us to (sometimes) simplify our life
when proving that two linear maps are mutually inverse:

\begin{lemma}
\label{lem.linalg.AB=1-BA=1}Let $r\in\mathbb{N}$. Let $V$ and $W$ be two free
$\mathbf{k}$-modules of rank $r$. Let $f:V\rightarrow W$ and $g:W\rightarrow
V$ be two $\mathbf{k}$-linear maps such that $f\circ g=\operatorname*{id}%
\nolimits_{W}$. Then, $g\circ f=\operatorname*{id}\nolimits_{V}$.
\end{lemma}

This fact is well-known in the case when $\mathbf{k}$ is a field (see, e.g.,
\cite[3.68]{Axler24} for this case), but actually holds for an arbitrary
commutative ring $\mathbf{k}$. The general case can be reduced to matrix
algebra in a standard way:

\begin{proof}
[Proof of Lemma \ref{lem.linalg.AB=1-BA=1}.]Choose a basis $\left(
v_{1},v_{2},\ldots,v_{r}\right)  $ of $V$ and a basis $\left(  w_{1}%
,w_{2},\ldots,w_{r}\right)  $ of $W$. Let $A\in\mathbf{k}^{r\times r}$ and
$B\in\mathbf{k}^{r\times r}$ be the matrices that represent the $\mathbf{k}%
$-linear maps $f$ and $g$ with respect to these two bases. Then, the matrices
$AB$ and $BA$ represent the compositions $f\circ g$ and $g\circ f$ with
respect to these two bases (since matrix multiplication corresponds to
composition of linear maps). Hence, from $f\circ g=\operatorname*{id}%
\nolimits_{W}$, we obtain $AB=I_{r}$ (where $I_{r}$ denotes the $r\times r$
identity matrix). But this entails $BA=I_{r}$ by a known result about square
matrices (see, e.g., \cite[Corollary 6.112 \textbf{(a)}]{detnotes}). This, in
turn, entails $g\circ f=\operatorname*{id}\nolimits_{V}$ (since the matrix
$BA$ represents the linear map $g\circ f$). This proves Lemma
\ref{lem.linalg.AB=1-BA=1}.
\end{proof}

\begin{remark}
Lemma \ref{lem.linalg.AB=1-BA=1} can be generalized: Instead of requiring $V$
and $W$ to be free $\mathbf{k}$-modules of the same rank $r\in\mathbb{N}$, we
can get by just requiring $V$ and $W$ to be two isomorphic finitely generated
$\mathbf{k}$-modules. For a proof of this version, see \cite[\textquotedblleft
Universal identities\textquotedblright, Corollary 5.6]{Conrad} (why does this
generalize Lemma \ref{lem.linalg.AB=1-BA=1}?).
\end{remark}

Next, we recall two simple properties of ranks of free modules. The first one
is the straightforward generalization (to arbitrary rings) of the classical
fact that the dimension of a direct product is the sum of the dimensions of
the factors:

\begin{lemma}
\label{lem.linalg.dim-dirprod}Let $M_{1},M_{2},\ldots,M_{k}$ be $k$ free
$\mathbf{k}$-modules of ranks $r_{1},r_{2},\ldots,r_{k}$, respectively. Then,
their direct product $M_{1}\times M_{2}\times\cdots\times M_{k}$ is a free
$\mathbf{k}$-module of rank $r_{1}+r_{2}+\cdots+r_{k}$.
\end{lemma}

\begin{proof}
[Proof idea.]For each $i\in\left[  k\right]  $, choose a basis $\left(
m_{i,1},m_{i,2},\ldots,m_{i,r_{i}}\right)  $ of the free $\mathbf{k}$-module
$M_{i}$. Let $m_{i,j}^{\prime}$ be the image of each basis vector $m_{i,j}\in
M_{i}$ under the canonical embedding $M_{i}\rightarrow M_{1}\times M_{2}%
\times\cdots\times M_{k}$. Then, the list
\begin{align*}
&  (m_{1,1}^{\prime},m_{1,2}^{\prime},\ldots,m_{1,r_{1}}^{\prime},\\
&  \ \text{\ }m_{2,1}^{\prime},m_{2,2}^{\prime},\ldots,m_{2,r_{2}}^{\prime},\\
&  \ \ \ldots,\\
&  \ \ m_{k,1}^{\prime},m_{k,2}^{\prime},\ldots,m_{k,r_{k}}^{\prime})
\end{align*}
is a basis of $M_{1}\times M_{2}\times\cdots\times M_{k}$. Since this basis
consists of $r_{1}+r_{2}+\cdots+r_{k}$ vectors, it follows that $M_{1}\times
M_{2}\times\cdots\times M_{k}$ is a free $\mathbf{k}$-module of rank
$r_{1}+r_{2}+\cdots+r_{k}$.
\end{proof}

The next property we need is a formula for the rank (i.e., dimension) of the
endomorphism ring of a free $\mathbf{k}$-module:

\begin{lemma}
\label{lem.linalg.dim-End}Let $M$ be a free $\mathbf{k}$-module of rank $r$.
Then, $\operatorname*{End}\nolimits_{\mathbf{k}}M$ is a free $\mathbf{k}%
$-module of rank $r^{2}$.
\end{lemma}

\begin{proof}
[Proof idea.]Choose a basis $\left(  m_{1},m_{2},\ldots,m_{r}\right)  $ of
$M$. Then, each endomorphism $f\in\operatorname*{End}\nolimits_{\mathbf{k}}M$
is represented by a unique $r\times r$-matrix $A_{f}\in\mathbf{k}^{r\times r}$
with respect to this basis. This gives a one-to-one correspondence between
endomorphisms $f\in\operatorname*{End}\nolimits_{\mathbf{k}}M$ and $r\times
r$-matrices $A\in\mathbf{k}^{r\times r}$. That is, we obtain a bijection%
\begin{align*}
\operatorname*{End}\nolimits_{\mathbf{k}}M  &  \rightarrow\mathbf{k}^{r\times
r},\\
f  &  \mapsto A_{f}.
\end{align*}
Moreover, this bijection is $\mathbf{k}$-linear (since addition of matrices
corresponds to addition of linear maps, and the same holds for scaling), and
thus is a $\mathbf{k}$-module isomorphism. Hence, the $\mathbf{k}$-module
$\operatorname*{End}\nolimits_{\mathbf{k}}M$ is isomorphic to $\mathbf{k}%
^{r\times r}$.

But the $\mathbf{k}$-module $\mathbf{k}^{r\times r}$ has a basis consisting of
the $r^{2}$ elementary matrices (each of which has a $1$ in some cell and
$0$'s in all remaining cells). Hence, the $\mathbf{k}$-module $\mathbf{k}%
^{r\times r}$ is free of rank $r^{2}$. Thus, the $\mathbf{k}$-module
$\operatorname*{End}\nolimits_{\mathbf{k}}M$ is free of rank $r^{2}$ as well
(since the $\mathbf{k}$-module $\operatorname*{End}\nolimits_{\mathbf{k}}M$ is
isomorphic to $\mathbf{k}^{r\times r}$).
\end{proof}

We also recall that the endomorphism ring of a free $\mathbf{k}$-module is
isomorphic to a matrix ring:

\begin{lemma}
\label{lem.linalg.End-mtring}Let $M$ be a free $\mathbf{k}$-module of rank
$r$. Then, $\operatorname*{End}\nolimits_{\mathbf{k}}M\cong\mathbf{k}^{r\times
r}$ as $\mathbf{k}$-algebras.
\end{lemma}

\begin{proof}
[Proof idea.]Argue as in the proof of Lemma \ref{lem.linalg.dim-End}, but
observe additionally that the bijection $\operatorname*{End}%
\nolimits_{\mathbf{k}}M\rightarrow\mathbf{k}^{r\times r}$ is a $\mathbf{k}%
$-algebra morphism (since multiplication of matrices corresponds to
composition of linear maps) and thus is a $\mathbf{k}$-algebra isomorphism.
\end{proof}

\subsubsection{The proof}

We can now prove Theorem \ref{thm.specht.AW}:

\begin{proof}
[Proof of Theorem \ref{thm.specht.AW}.]\textit{Step 1:} In the following, for
any given partition $\lambda$ of $n$, we shall abbreviate the summation sign
\textquotedblleft$\sum_{T\text{ is an }n\text{-tableau of shape }\lambda}%
$\textquotedblright\ as \textquotedblleft$\sum_{T\text{ of shape }\lambda}%
$\textquotedblright. Using this abbreviation, we can rewrite the equality
(\ref{eq.def.spechtmod.Elam.Elam.def}) as%
\begin{equation}
\mathbf{E}_{\lambda}=\sum_{T\text{ of shape }\lambda}\mathbf{E}_{T}
\label{pf.thm.specht.AW.Elam=}%
\end{equation}
(for any partition $\lambda$ of $n$). \medskip

\textit{Step 2:} Define a map%
\begin{align*}
\Phi:\mathbf{k}\left[  S_{n}\right]   &  \rightarrow\prod_{\lambda\text{ is a
partition of }n}\operatorname*{End}\nolimits_{\mathbf{k}}\left(
\mathcal{S}^{\lambda}\right)  ,\\
\mathbf{a}  &  \mapsto\left(  \rho_{\lambda}\left(  \mathbf{a}\right)
\right)  _{\lambda\text{ is a partition of }n}.
\end{align*}
This map $\Phi$ is precisely the map%
\begin{align*}
\mathbf{k}\left[  S_{n}\right]   &  \rightarrow\prod_{\lambda\text{ is a
partition of }n}\operatorname*{End}\nolimits_{\mathbf{k}}\left(
\mathcal{S}^{\lambda}\right)  ,\\
\mathbf{a}  &  \mapsto\left(  \rho_{\lambda}\left(  \mathbf{a}\right)
\right)  _{\lambda\text{ is a partition of }n}%
\end{align*}
in Theorem \ref{thm.specht.AW}. Thus, our goal is to prove that this map
$\Phi$ is a $\mathbf{k}$-algebra isomorphism.

We recall that if $R$ is any $\mathbf{k}$-algebra, and if $V$ is any left
$R$-module, then the curried form $\rho:R\rightarrow\operatorname*{End}%
\nolimits_{\mathbf{k}}V$ of the left $R$-action on $V$ is a $\mathbf{k}%
$-algebra morphism (by Theorem \ref{thm.mod.curry-into-EndV}). Hence, for each
partition $\lambda$ of $n$, the curried form $\rho_{\lambda}:\mathbf{k}\left[
S_{n}\right]  \rightarrow\operatorname*{End}\nolimits_{\mathbf{k}}\left(
\mathcal{S}^{\lambda}\right)  $ of the left $\mathbf{k}\left[  S_{n}\right]
$-action on the Specht module $\mathcal{S}^{\lambda}$ is a $\mathbf{k}%
$-algebra morphism. Thus, the map $\Phi$ (which is defined by applying all
these morphisms $\rho_{\lambda}$ to its input and collecting the resulting
outputs in a family) is a $\mathbf{k}$-algebra morphism as well\footnote{since
the $\mathbf{k}$-algebra structure on $\prod_{\lambda\text{ is a partition of
}n}\operatorname*{End}\nolimits_{\mathbf{k}}\left(  \mathcal{S}^{\lambda
}\right)  $ is entrywise}.

We shall now focus on proving that $\Phi$ is invertible. \medskip

\begin{noncompile}
First, let us recall a trivial fact about $\mathbf{k}$-modules: If $M$ and $N$
are two isomorphic $\mathbf{k}$-modules, then any $\mathbf{k}$-module
endomorphism of $M$ can be converted into a $\mathbf{k}$-module isomorphism of
$N$. More concretely: If $\alpha:M\rightarrow N$ is a $\mathbf{k}$-module
isomorphism between two $\mathbf{k}$-modules $M$ and $N$, and if $\zeta
\in\operatorname*{End}\nolimits_{\mathbf{k}}M$ is a $\mathbf{k}$-module
endomorphism of $M$, then $\alpha\circ\zeta\circ\alpha^{-1}\in
\operatorname*{End}\nolimits_{\mathbf{k}}N$ is a $\mathbf{k}$-module
endomorphism of $N$.

We shall use this to convert endomorphisms of the Specht modules
$\mathcal{S}^{\lambda}$ into endomorphisms of their left ideal avatars
$\mathcal{A}\mathbf{E}_{T}$. Namely:
\end{noncompile}

\textit{Step 3:} Let $\lambda$ be a partition of $n$. Let $T$ be an
$n$-tableau of shape $\lambda$. Then, $T$ is an $n$-tableau of shape $Y\left(
\lambda\right)  $. Hence, (\ref{eq.def.specht.ET.defs.SD=AET}) (applied to
$D=Y\left(  \lambda\right)  $) yields $\mathcal{S}^{Y\left(  \lambda\right)
}\cong\mathcal{A}\mathbf{E}_{T}$ (as left $\mathbf{k}\left[  S_{n}\right]
$-modules, i.e., as left $\mathcal{A}$-modules). Hence, $\mathcal{S}^{\lambda
}=\mathcal{S}^{Y\left(  \lambda\right)  }\cong\mathcal{A}\mathbf{E}_{T}$. In
other words, there exists an $\mathcal{A}$-module isomorphism $\alpha
_{\lambda,T}:\mathcal{S}^{\lambda}\rightarrow\mathcal{A}\mathbf{E}_{T}$. Fix
such an isomorphism $\alpha_{\lambda,T}$ once and for all.\footnote{A
canonical choice of $\alpha_{\lambda,T}$ can be obtained from Theorem
\ref{thm.spechtmod.leftideal} \textbf{(b)}, but we don't need it to be
canonical here.}

Forget that we fixed $\lambda$ and $T$. Thus, for any partition $\lambda$ of
$n$ and any $n$-tableau $T$ of shape $\lambda$, we have defined an
$\mathcal{A}$-module isomorphism $\alpha_{\lambda,T}:\mathcal{S}^{\lambda
}\rightarrow\mathcal{A}\mathbf{E}_{T}$. In particular, $\alpha_{\lambda,T}$ is
thus a $\mathbf{k}$-module isomorphism. Therefore, if $\zeta\in
\operatorname*{End}\nolimits_{\mathbf{k}}\left(  \mathcal{S}^{\lambda}\right)
$ is a $\mathbf{k}$-module endomorphism of $\mathcal{S}^{\lambda}$, then
$\alpha_{\lambda,T}\circ\zeta\circ\alpha_{\lambda,T}^{-1}\in
\operatorname*{End}\nolimits_{\mathbf{k}}\left(  \mathcal{A}\mathbf{E}%
_{T}\right)  $ is a $\mathbf{k}$-module endomorphism of $\mathcal{A}%
\mathbf{E}_{T}$. \medskip

\textit{Step 4:} Now, define a map%
\[
\Psi:\prod_{\lambda\text{ is a partition of }n}\operatorname*{End}%
\nolimits_{\mathbf{k}}\left(  \mathcal{S}^{\lambda}\right)  \rightarrow
\mathbf{k}\left[  S_{n}\right]  ,
\]
which sends each family $\left(  \zeta_{\lambda}\right)  _{\lambda\text{ is a
partition of }n}\in\prod_{\lambda\text{ is a partition of }n}%
\operatorname*{End}\nolimits_{\mathbf{k}}\left(  \mathcal{S}^{\lambda}\right)
$ to%
\[
\sum_{\lambda\text{ is a partition of }n}\dfrac{1}{\left(  h^{\lambda}\right)
^{2}}\sum_{T\text{ of shape }\lambda}\underbrace{\left(  \alpha_{\lambda
,T}\circ\zeta_{\lambda}\circ\alpha_{\lambda,T}^{-1}\right)  \left(
\mathbf{E}_{T}\right)  }_{\substack{\text{(this is well-defined,}\\\text{since
}\mathbf{E}_{T}=1\mathbf{E}_{T}\in\mathcal{A}\mathbf{E}_{T}\\\text{and }%
\alpha_{\lambda,T}\circ\zeta_{\lambda}\circ\alpha_{\lambda,T}^{-1}%
\in\operatorname*{End}\nolimits_{\mathbf{k}}\left(  \mathcal{A}\mathbf{E}%
_{T}\right)  \text{)}}}\in\mathcal{A}=\mathbf{k}\left[  S_{n}\right]  .
\]
This map $\Psi$ is clearly $\mathbf{k}$-linear (since each $\left(
\alpha_{\lambda,T}\circ\zeta_{\lambda}\circ\alpha_{\lambda,T}^{-1}\right)
\left(  \mathbf{E}_{T}\right)  $ depends $\mathbf{k}$-linearly on
$\zeta_{\lambda}$ and thus on the input of $\Psi$). \medskip

\textit{Step 5:} We shall next show that $\Psi\circ\Phi=\operatorname*{id}$.

Indeed, let $\mathbf{a}\in\mathbf{k}\left[  S_{n}\right]  $. Then, the
definition of $\Phi$ yields%
\[
\Phi\left(  \mathbf{a}\right)  =\left(  \rho_{\lambda}\left(  \mathbf{a}%
\right)  \right)  _{\lambda\text{ is a partition of }n}.
\]
Hence,%
\begin{align*}
\Psi\left(  \Phi\left(  \mathbf{a}\right)  \right)   &  =\Psi\left(  \left(
\rho_{\lambda}\left(  \mathbf{a}\right)  \right)  _{\lambda\text{ is a
partition of }n}\right) \\
&  =\sum_{\lambda\text{ is a partition of }n}\dfrac{1}{\left(  h^{\lambda
}\right)  ^{2}}\sum_{T\text{ of shape }\lambda}\left(  \alpha_{\lambda,T}%
\circ\left(  \rho_{\lambda}\left(  \mathbf{a}\right)  \right)  \circ
\alpha_{\lambda,T}^{-1}\right)  \left(  \mathbf{E}_{T}\right)
\end{align*}
(by the definition of $\Psi$).

However, we can simplify the addends in this sum significantly, due to the
following claim:

\begin{statement}
\textit{Claim 1:} Let $\lambda$ be a partition of $n$. Let $T$ be an
$n$-tableau of shape $\lambda$. Then,%
\[
\left(  \alpha_{\lambda,T}\circ\left(  \rho_{\lambda}\left(  \mathbf{a}%
\right)  \right)  \circ\alpha_{\lambda,T}^{-1}\right)  \left(  \mathbf{E}%
_{T}\right)  =\mathbf{aE}_{T}.
\]

\end{statement}

\begin{proof}
[Proof of Claim 1.]Set $\mathbf{u}=\alpha_{\lambda,T}^{-1}\left(
\mathbf{E}_{T}\right)  $. Then, $\alpha_{\lambda,T}\left(  \mathbf{u}\right)
=\mathbf{E}_{T}$. Furthermore, we have $\left(  \rho_{\lambda}\left(
\mathbf{a}\right)  \right)  \left(  \mathbf{u}\right)  =\mathbf{au}$ (since
$\rho_{\lambda}$ is just the curried form of the left $\mathbf{k}\left[
S_{n}\right]  $-action on $\mathcal{S}^{\lambda}$). Moreover, $\alpha
_{\lambda,T}\left(  \mathbf{au}\right)  =\mathbf{a}\alpha_{\lambda,T}\left(
\mathbf{u}\right)  $ (since $\alpha_{\lambda,T}$ is an $\mathcal{A}$-module
morphism). Now,%
\begin{align*}
\left(  \alpha_{\lambda,T}\circ\left(  \rho_{\lambda}\left(  \mathbf{a}%
\right)  \right)  \circ\alpha_{\lambda,T}^{-1}\right)  \left(  \mathbf{E}%
_{T}\right)   &  =\alpha_{\lambda,T}\left(  \left(  \rho_{\lambda}\left(
\mathbf{a}\right)  \right)  \left(  \underbrace{\alpha_{\lambda,T}^{-1}\left(
\mathbf{E}_{T}\right)  }_{=\mathbf{u}}\right)  \right) \\
&  =\alpha_{\lambda,T}\left(  \underbrace{\left(  \rho_{\lambda}\left(
\mathbf{a}\right)  \right)  \left(  \mathbf{u}\right)  }_{=\mathbf{au}%
}\right)  =\alpha_{\lambda,T}\left(  \mathbf{au}\right)  =\mathbf{a}%
\underbrace{\alpha_{\lambda,T}\left(  \mathbf{u}\right)  }_{=\mathbf{E}_{T}}\\
&  =\mathbf{aE}_{T}.
\end{align*}
This proves Claim 1.
\end{proof}

Now,
\begin{align*}
\left(  \Psi\circ\Phi\right)  \left(  \mathbf{a}\right)   &  =\Psi\left(
\Phi\left(  \mathbf{a}\right)  \right) \\
&  =\sum_{\lambda\text{ is a partition of }n}\dfrac{1}{\left(  h^{\lambda
}\right)  ^{2}}\sum_{T\text{ of shape }\lambda}\underbrace{\left(
\alpha_{\lambda,T}\circ\left(  \rho_{\lambda}\left(  \mathbf{a}\right)
\right)  \circ\alpha_{\lambda,T}^{-1}\right)  \left(  \mathbf{E}_{T}\right)
}_{\substack{=\mathbf{aE}_{T}\\\text{(by Claim 1)}}}\\
&  \ \ \ \ \ \ \ \ \ \ \ \ \ \ \ \ \ \ \ \ \left(  \text{by our above
computation of }\Psi\left(  \Phi\left(  \mathbf{a}\right)  \right)  \right) \\
&  =\sum_{\lambda\text{ is a partition of }n}\dfrac{1}{\left(  h^{\lambda
}\right)  ^{2}}\sum_{T\text{ of shape }\lambda}\mathbf{aE}_{T}\\
&  =\mathbf{a}\sum_{\lambda\text{ is a partition of }n}\dfrac{1}{\left(
h^{\lambda}\right)  ^{2}}\underbrace{\sum_{T\text{ of shape }\lambda
}\mathbf{E}_{T}}_{\substack{=\mathbf{E}_{\lambda}\\\text{(by
(\ref{pf.thm.specht.AW.Elam=}))}}}\\
&  =\mathbf{a}\sum_{\lambda\text{ is a partition of }n}\dfrac{1}{\left(
h^{\lambda}\right)  ^{2}}\mathbf{E}_{\lambda}=\mathbf{a}\underbrace{\sum
_{\lambda\text{ is a partition of }n}\dfrac{\mathbf{E}_{\lambda}}{\left(
h^{\lambda}\right)  ^{2}}}_{\substack{=1\\\text{(by Corollary
\ref{cor.spechtmod.Elam.pou} \textbf{(b)})}}}=\mathbf{a}.
\end{align*}

Forget that we fixed $\mathbf{a}$. We thus have shown that $\left(  \Psi
\circ\Phi\right)  \left(  \mathbf{a}\right)  =\mathbf{a}$ for each
$\mathbf{a}\in\mathbf{k}\left[  S_{n}\right]  $. In other words, $\Psi
\circ\Phi=\operatorname*{id}$. \medskip

\textit{Step 6:} Next we want to show that $\Phi\circ\Psi=\operatorname*{id}$
as well. We will derive this from $\Psi\circ\Phi=\operatorname*{id}$ using
Lemma \ref{lem.linalg.AB=1-BA=1}. In order to do so, we need to check that
$\mathbf{k}\left[  S_{n}\right]  $ and $\prod_{\lambda\text{ is a partition of
}n}\operatorname*{End}\nolimits_{\mathbf{k}}\left(  \mathcal{S}^{\lambda
}\right)  $ are two free $\mathbf{k}$-modules of the same rank (namely, $n!$).

Let us now do this.

Let $\lambda$ be a partition of $n$. Then, Lemma \ref{lem.specht.Slam-flam}
shows that $\mathcal{S}^{\lambda}$ is a free $\mathbf{k}$-module of rank
$f^{\lambda}$. Hence, Lemma \ref{lem.linalg.dim-End} (applied to
$M=\mathcal{S}^{\lambda}$ and $r=f^{\lambda}$) shows that $\operatorname*{End}%
\nolimits_{\mathbf{k}}\left(  \mathcal{S}^{\lambda}\right)  $ is a free
$\mathbf{k}$-module of rank $\left(  f^{\lambda}\right)  ^{2}$.

Forget that we fixed $\lambda$. We thus have shown that $\operatorname*{End}%
\nolimits_{\mathbf{k}}\left(  \mathcal{S}^{\lambda}\right)  $ is a free
$\mathbf{k}$-module of rank $\left(  f^{\lambda}\right)  ^{2}$ whenever
$\lambda$ is a partition of $n$. Hence, by Lemma \ref{lem.linalg.dim-dirprod},
we conclude that the direct product $\prod_{\lambda\text{ is a partition of
}n}\operatorname*{End}\nolimits_{\mathbf{k}}\left(  \mathcal{S}^{\lambda
}\right)  $ of these free $\mathbf{k}$-modules $\operatorname*{End}%
\nolimits_{\mathbf{k}}\left(  \mathcal{S}^{\lambda}\right)  $ is a free
$\mathbf{k}$-module of rank $\sum_{\lambda\text{ is a partition of }n}\left(
f^{\lambda}\right)  ^{2}$. In other words, the direct product $\prod
_{\lambda\text{ is a partition of }n}\operatorname*{End}\nolimits_{\mathbf{k}%
}\left(  \mathcal{S}^{\lambda}\right)  $ is a free $\mathbf{k}$-module of rank
$n!$ (since Corollary \ref{cor.spechtmod.sumflam2} yields $\sum_{\lambda\text{
is a partition of }n}\left(  f^{\lambda}\right)  ^{2}=n!$).

But the $\mathbf{k}$-module $\mathbf{k}\left[  S_{n}\right]  $ is also a free
$\mathbf{k}$-module of rank $n!$ (since it has a basis $\left(  w\right)
_{w\in S_{n}}$, which consists of $\left\vert S_{n}\right\vert =n!$ many vectors).

Thus, Lemma \ref{lem.linalg.AB=1-BA=1} (applied to $V=\prod_{\lambda\text{ is
a partition of }n}\operatorname*{End}\nolimits_{\mathbf{k}}\left(
\mathcal{S}^{\lambda}\right)  $ and $W=\mathbf{k}\left[  S_{n}\right]  $ and
$r=n!$ and $f=\Psi$ and $g=\Phi$) shows that $\Phi\circ\Psi=\operatorname*{id}%
$ (since $\Psi\circ\Phi=\operatorname*{id}$). \medskip

\textit{Step 7:} Combining $\Phi\circ\Psi=\operatorname*{id}$ with $\Psi
\circ\Phi=\operatorname*{id}$, we see that the maps $\Phi$ and $\Psi$ are
mutually inverse. Hence, the map $\Phi$ is invertible, and thus is a
$\mathbf{k}$-algebra isomorphism (since it is a $\mathbf{k}$-algebra
morphism). As explained above, this proves Theorem \ref{thm.specht.AW}.
\end{proof}

\begin{exercise}
Let $\mathbf{a}\in\mathbf{k}\left[  S_{n}\right]  $. Assume that
$\mathbf{a}\mathcal{S}^{\lambda}=0$ for each partition $\lambda$ of $n$. (As
usual, $\mathbf{a}\mathcal{S}^{\lambda}$ means $\left\{  \mathbf{av}%
\ \mid\ \mathbf{v}\in\mathcal{S}^{\lambda}\right\}  $.) \medskip

\textbf{(a)} \fbox{1} Prove that $\mathbf{a}=0$ if $n!$ is invertible in
$\mathbf{k}$. \medskip

\textbf{(b)} \fbox{3} Prove that $n!^{2}\cdot\mathbf{a}=0$ always.
\end{exercise}

\subsubsection{Proof of Theorem \ref{thm.AWS.demo} \textbf{(a)}}

We can now easily verify Theorem \ref{thm.AWS.demo} \textbf{(a)} as well
(which is just a less specific version of Theorem \ref{thm.specht.AW}):

\begin{proof}
[Proof of Theorem \ref{thm.AWS.demo} \textbf{(a)}.]Let $\lambda$ be a
partition of $n$. Then, $\mathcal{S}^{\lambda}$ is a free $\mathbf{k}$-module
of rank $f^{\lambda}$ (by Lemma \ref{lem.specht.Slam-flam}). Hence,
\begin{equation}
\operatorname*{End}\nolimits_{\mathbf{k}}\left(  \mathcal{S}^{\lambda}\right)
\cong\mathbf{k}^{f^{\lambda}\times f^{\lambda}}\ \ \ \ \ \ \ \ \ \ \text{as
}\mathbf{k}\text{-algebras} \label{pf.thm.AWS.demo.a.1}%
\end{equation}
(by Lemma \ref{lem.linalg.End-mtring}, applied to $M=\mathcal{S}^{\lambda}$
and $r=f^{\lambda}$).

Forget that we fixed $\lambda$. We thus have proved (\ref{pf.thm.AWS.demo.a.1}%
) for each partition $\lambda$ of $n$. Now, Theorem \ref{thm.specht.AW}
constructs a $\mathbf{k}$-algebra isomorphism%
\[
\mathbf{k}\left[  S_{n}\right]  \rightarrow\prod_{\lambda\text{ is a partition
of }n}\operatorname*{End}\nolimits_{\mathbf{k}}\left(  \mathcal{S}^{\lambda
}\right)  .
\]
Thus, we see that%
\[
\mathbf{k}\left[  S_{n}\right]  \cong\prod_{\lambda\text{ is a partition of
}n}\underbrace{\operatorname*{End}\nolimits_{\mathbf{k}}\left(  \mathcal{S}%
^{\lambda}\right)  }_{\substack{\cong\mathbf{k}^{f^{\lambda}\times f^{\lambda
}}\\\text{(by (\ref{pf.thm.AWS.demo.a.1}))}}}\cong\prod_{\lambda\text{ is a
partition of }n}\mathbf{k}^{f^{\lambda}\times f^{\lambda}}%
\]
as $\mathbf{k}$-algebras. In other words,
\[
\mathbf{k}\left[  S_{n}\right]  \cong\prod_{\lambda\in\Lambda}\mathbf{k}%
^{f_{\lambda}\times f_{\lambda}},
\]
where $\Lambda$ is the set of all partitions of $n$, and where each
$f_{\lambda}$ is what we call $f^{\lambda}$. This proves Theorem
\ref{thm.AWS.demo} \textbf{(a)}.
\end{proof}

We have now proved all of Theorem \ref{thm.AWS.demo} except for part
\textbf{(c)}. As we mentioned, the latter part requires the \emph{seminormal
basis}, which connects Young tableaux with the jucys--murphies and would
require an entire new chapter.

\subsubsection{The kernel and the surjectivity of $\rho_{\lambda}%
:\mathbf{k}\left[  S_{n}\right]  \rightarrow\operatorname*{End}%
\nolimits_{\mathbf{k}}\left(  \mathcal{S}^{\lambda}\right)  $}

Theorem \ref{thm.specht.AW} claims something like a \textquotedblleft joint
bijectivity\textquotedblright\ of the maps $\rho_{\lambda}:\mathcal{A}%
\rightarrow\operatorname*{End}\nolimits_{\mathbf{k}}\left(  \mathcal{S}%
^{\lambda}\right)  $ -- meaning that even if each of these maps is not
bijective per se, their combination (sending $\mathbf{a}\mapsto\left(
\rho_{\lambda}\left(  \mathbf{a}\right)  \right)  _{\lambda\text{ is a
partition of }n}$) is (as long as $n!$ is invertible in $\mathbf{k}$). One
might wonder what can be said about a single map $\rho_{\lambda}$ rather than
all of them in combination. And indeed, there are things that can be said
about it, and even under a slightly weaker assumption (requiring the
invertibility of $h^{\lambda}$ rather than the stronger invertibility of
$n!$). We begin with an easy description of its kernel:

\begin{proposition}
\label{prop.specht.AWlam.ker}Let $\lambda$ be a partition of $n$. Assume that
$h^{\lambda}$ is invertible in $\mathbf{k}$. Let $\rho_{\lambda}%
:\mathbf{k}\left[  S_{n}\right]  \rightarrow\operatorname*{End}%
\nolimits_{\mathbf{k}}\left(  \mathcal{S}^{\lambda}\right)  $ be the curried
form of the left $\mathbf{k}\left[  S_{n}\right]  $-action on the Specht
module $\mathcal{S}^{\lambda}$ (as in Theorem \ref{thm.specht.AW}). Let
$\mathbf{G}_{\lambda}:=1-\dfrac{\mathbf{E}_{\lambda}}{\left(  h^{\lambda
}\right)  ^{2}}\in\mathcal{A}$. Then,%
\[
\operatorname*{Ker}\left(  \rho_{\lambda}\right)  =\mathcal{A}\mathbf{G}%
_{\lambda}.
\]

\end{proposition}

This will follow easily from the following fact, which is interesting by itself:

\begin{proposition}
\label{prop.specht.Elam-Slam}Let $\lambda$ be a partition of $n$. Let
$\mathbf{v}\in\mathcal{S}^{\lambda}$. Then, $\mathbf{E}_{\lambda}%
\mathbf{v}=\left(  h^{\lambda}\right)  ^{2}\mathbf{v}$.
\end{proposition}

\begin{proof}
The diagram $Y\left(  \lambda\right)  $ has size $\left\vert Y\left(
\lambda\right)  \right\vert =\left\vert \lambda\right\vert =n$ (since
$\lambda$ is a partition of $n$). Pick any $n$-tableau $T$ of shape $Y\left(
\lambda\right)  $. (Such a $T$ clearly exists, since $\left\vert Y\left(
\lambda\right)  \right\vert =n$.)

Thus, (\ref{eq.def.specht.ET.defs.SD=AET}) (applied to $D=Y\left(
\lambda\right)  $) yields $\mathcal{S}^{Y\left(  \lambda\right)  }%
\cong\mathcal{A}\mathbf{E}_{T}$ as left $\mathcal{A}$-modules. In other words,
$\mathcal{S}^{\lambda}\cong\mathcal{A}\mathbf{E}_{T}$ as left $\mathcal{A}%
$-modules (since $\mathcal{S}^{\lambda}=\mathcal{S}^{Y\left(  \lambda\right)
}$). In other words, there exists a left $\mathcal{A}$-module isomorphism
$f:\mathcal{S}^{\lambda}\rightarrow\mathcal{A}\mathbf{E}_{T}$. Consider this
$f$. Thus, $f\left(  \mathbf{v}\right)  \in\mathcal{A}\mathbf{E}_{T}$. In
other words, $f\left(  \mathbf{v}\right)  =\mathbf{aE}_{T}$ for some
$\mathbf{a}\in\mathcal{A}$. Consider this $\mathbf{a}$.

Proposition \ref{prop.spechtmod.Elam.incent} yields that $\mathbf{E}_{\lambda
}$ belongs to the center $Z\left(  \mathbf{k}\left[  S_{n}\right]  \right)  $
of $\mathbf{k}\left[  S_{n}\right]  $. Hence, $\mathbf{E}_{\lambda}%
\mathbf{a}=\mathbf{aE}_{\lambda}$ (since $\mathbf{a}\in\mathcal{A}%
=\mathbf{k}\left[  S_{n}\right]  $).

Now, recall that $f$ is a left $\mathcal{A}$-module morphism. Hence,%
\begin{align*}
f\left(  \mathbf{E}_{\lambda}\mathbf{v}\right)   &  =\mathbf{E}_{\lambda
}\underbrace{f\left(  \mathbf{v}\right)  }_{=\mathbf{aE}_{T}}%
=\underbrace{\mathbf{E}_{\lambda}\mathbf{a}}_{=\mathbf{aE}_{\lambda}%
}\mathbf{E}_{T}=\mathbf{a}\underbrace{\mathbf{E}_{\lambda}\mathbf{E}_{T}%
}_{\substack{=\left(  h^{\lambda}\right)  ^{2}\mathbf{E}_{T}\\\text{(by
Proposition \ref{prop.spechtmod.Elam.ETEl})}}}\\
&  =\left(  h^{\lambda}\right)  ^{2}\underbrace{\mathbf{aE}_{T}}_{=f\left(
\mathbf{v}\right)  }=\left(  h^{\lambda}\right)  ^{2}f\left(  \mathbf{v}%
\right)  =f\left(  \left(  h^{\lambda}\right)  ^{2}\mathbf{v}\right)
\end{align*}
(since $f$ is $\mathbf{k}$-linear). But the map $f$ is injective (since $f$ is
an isomorphism). Thus, from $f\left(  \mathbf{E}_{\lambda}\mathbf{v}\right)
=f\left(  \left(  h^{\lambda}\right)  ^{2}\mathbf{v}\right)  $, we obtain
$\mathbf{E}_{\lambda}\mathbf{v}=\left(  h^{\lambda}\right)  ^{2}\mathbf{v}$.
This proves Proposition \ref{prop.specht.Elam-Slam}.
\end{proof}

\begin{proof}
[Proof of Proposition \ref{prop.specht.AWlam.ker}.]\textit{Step 1:} We shall
first show that $\mathcal{A}\mathbf{G}_{\lambda}\subseteq\operatorname*{Ker}%
\left(  \rho_{\lambda}\right)  $.

\textit{Proof.} Let $\mathbf{a}\in\mathcal{A}\mathbf{G}_{\lambda}$. Then,
$\mathbf{a}=\mathbf{bG}_{\lambda}$ for some $\mathbf{b}\in\mathcal{A}$.
Consider this $\mathbf{b}$.

Let $\mathbf{v}\in\mathcal{S}^{\lambda}$. Then, from $\mathbf{G}_{\lambda
}=1-\dfrac{\mathbf{E}_{\lambda}}{\left(  h^{\lambda}\right)  ^{2}}$, we obtain%
\begin{align*}
\mathbf{G}_{\lambda}\mathbf{v}  &  =\left(  1-\dfrac{\mathbf{E}_{\lambda}%
}{\left(  h^{\lambda}\right)  ^{2}}\right)  \mathbf{v}=\mathbf{v}%
-\dfrac{\mathbf{E}_{\lambda}}{\left(  h^{\lambda}\right)  ^{2}}\mathbf{v}%
=\mathbf{v}-\dfrac{1}{\left(  h^{\lambda}\right)  ^{2}}\underbrace{\mathbf{E}%
_{\lambda}\mathbf{v}}_{\substack{=\left(  h^{\lambda}\right)  ^{2}%
\mathbf{v}\\\text{(by Proposition \ref{prop.specht.Elam-Slam})}}}\\
&  =\mathbf{v}-\underbrace{\dfrac{1}{\left(  h^{\lambda}\right)  ^{2}}\left(
h^{\lambda}\right)  ^{2}}_{=1}\mathbf{v}=\mathbf{v}-\mathbf{v}=0.
\end{align*}
Now,%
\[
\underbrace{\mathbf{a}}_{=\mathbf{bG}_{\lambda}}\mathbf{v}=\mathbf{b}%
\underbrace{\mathbf{G}_{\lambda}\mathbf{v}}_{=0}=0.
\]

The definition of $\rho_{\lambda}$ yields $\left(  \rho_{\lambda}\left(
\mathbf{a}\right)  \right)  \left(  \mathbf{v}\right)  =\mathbf{av}=0$.

Forget that we fixed $\mathbf{v}$. We thus have shown that $\left(
\rho_{\lambda}\left(  \mathbf{a}\right)  \right)  \left(  \mathbf{v}\right)
=0$ for each $\mathbf{v}\in\mathcal{S}^{\lambda}$. In other words,
$\rho_{\lambda}\left(  \mathbf{a}\right)  =0$. In other words, $\mathbf{a}%
\in\operatorname*{Ker}\left(  \rho_{\lambda}\right)  $.

Forget that we fixed $\mathbf{a}$. We thus have proved that $\mathbf{a}%
\in\operatorname*{Ker}\left(  \rho_{\lambda}\right)  $ for each $\mathbf{a}%
\in\mathcal{A}\mathbf{G}_{\lambda}$. In other words, $\mathcal{A}%
\mathbf{G}_{\lambda}\subseteq\operatorname*{Ker}\left(  \rho_{\lambda}\right)
$. \medskip

\textit{Step 2:} We shall now prove that $\operatorname*{Ker}\left(
\rho_{\lambda}\right)  \subseteq\mathcal{A}\mathbf{G}_{\lambda}$.

\textit{Proof.} Let $\mathbf{d}\in\operatorname*{Ker}\left(  \rho_{\lambda
}\right)  $. Thus, $\mathbf{d}\in\mathbf{k}\left[  S_{n}\right]  =\mathcal{A}$
and $\rho_{\lambda}\left(  \mathbf{d}\right)  =0$.

Now we shall show that $\mathbf{dE}_{T}=0$ for any $n$-tableau $T$ of shape
$\lambda$.

Indeed, let $T$ be an $n$-tableau of shape $\lambda$. Thus, $T$ is an
$n$-tableau of shape $Y\left(  \lambda\right)  $. Hence,
(\ref{eq.def.specht.ET.defs.SD=AET}) (applied to $D=Y\left(  \lambda\right)
$) yields $\mathcal{S}^{Y\left(  \lambda\right)  }\cong\mathcal{A}%
\mathbf{E}_{T}$ as left $\mathcal{A}$-modules. In other words, $\mathcal{S}%
^{\lambda}\cong\mathcal{A}\mathbf{E}_{T}$ as left $\mathcal{A}$-modules (since
$\mathcal{S}^{\lambda}=\mathcal{S}^{Y\left(  \lambda\right)  }$). In other
words, there exists a left $\mathcal{A}$-module isomorphism $f:\mathcal{S}%
^{\lambda}\rightarrow\mathcal{A}\mathbf{E}_{T}$. Consider this $f$. We have
$\mathbf{E}_{T}=\underbrace{1}_{\in\mathcal{A}}\mathbf{E}_{T}\in
\mathcal{A}\mathbf{E}_{T}$. Hence, there exists some $\mathbf{v}\in
\mathcal{S}^{\lambda}$ such that $\mathbf{E}_{T}=f\left(  \mathbf{v}\right)  $
(since $f$ is an isomorphism and thus surjective). Consider this $\mathbf{v}$.

The definition of $\rho_{\lambda}$ yields $\left(  \rho_{\lambda}\left(
\mathbf{d}\right)  \right)  \left(  \mathbf{v}\right)  =\mathbf{dv}$, so that
$\mathbf{dv}=\underbrace{\left(  \rho_{\lambda}\left(  \mathbf{d}\right)
\right)  }_{=0}\left(  \mathbf{v}\right)  =0$. Hence, $f\left(  \mathbf{dv}%
\right)  =f\left(  0\right)  =0$ (since $f$ is $\mathbf{k}$-linear). But $f$
is a left $\mathcal{A}$-module morphism; thus, $f\left(  \mathbf{dv}\right)
=\mathbf{d}\underbrace{f\left(  \mathbf{v}\right)  }_{=\mathbf{E}_{T}%
}=\mathbf{dE}_{T}$. Hence, $\mathbf{dE}_{T}=f\left(  \mathbf{dv}\right)  =0$.

Forget that we fixed $T$. We thus have shown that
\begin{equation}
\mathbf{dE}_{T}=0\ \ \ \ \ \ \ \ \ \ \text{for any }n\text{-tableau }T\text{
of shape }\lambda. \label{pf.prop.specht.AWlam.ker.b.4}%
\end{equation}

Now,%
\begin{align*}
\mathbf{dE}_{\lambda}  &  =\mathbf{d}\sum_{T\text{ is an }n\text{-tableau of
shape }\lambda}\mathbf{E}_{T}\ \ \ \ \ \ \ \ \ \ \left(  \text{by
(\ref{eq.def.spechtmod.Elam.Elam.def})}\right) \\
&  =\sum_{T\text{ is an }n\text{-tableau of shape }\lambda}%
\underbrace{\mathbf{dE}_{T}}_{\substack{=0\\\text{(by
(\ref{pf.prop.specht.AWlam.ker.b.4}))}}}=0.
\end{align*}
Hence,%
\begin{align*}
\mathbf{dG}_{\lambda}  &  =\mathbf{d}\left(  1-\dfrac{\mathbf{E}_{\lambda}%
}{\left(  h^{\lambda}\right)  ^{2}}\right)  \ \ \ \ \ \ \ \ \ \ \left(
\text{since }\mathbf{G}_{\lambda}=1-\dfrac{\mathbf{E}_{\lambda}}{\left(
h^{\lambda}\right)  ^{2}}\right) \\
&  =\mathbf{d}-\dfrac{\mathbf{dE}_{\lambda}}{\left(  h^{\lambda}\right)  ^{2}%
}=\mathbf{d}-\dfrac{0}{\left(  h^{\lambda}\right)  ^{2}}%
\ \ \ \ \ \ \ \ \ \ \left(  \text{since }\mathbf{dE}_{\lambda}=0\right) \\
&  =\mathbf{d},
\end{align*}
so that $\mathbf{d}=\underbrace{\mathbf{d}}_{\in\mathcal{A}}\mathbf{G}%
_{\lambda}\in\mathcal{A}\mathbf{G}_{\lambda}$.

Forget that we fixed $\mathbf{d}$. We thus have shown that $\mathbf{d}%
\in\mathcal{A}\mathbf{G}_{\lambda}$ for each $\mathbf{d}\in\operatorname*{Ker}%
\left(  \rho_{\lambda}\right)  $. In other words, $\operatorname*{Ker}\left(
\rho_{\lambda}\right)  \subseteq\mathcal{A}\mathbf{G}_{\lambda}$. \medskip

\textit{Step 3:} Combining $\operatorname*{Ker}\left(  \rho_{\lambda}\right)
\subseteq\mathcal{A}\mathbf{G}_{\lambda}$ with $\mathcal{A}\mathbf{G}%
_{\lambda}\subseteq\operatorname*{Ker}\left(  \rho_{\lambda}\right)  $, we
obtain $\operatorname*{Ker}\left(  \rho_{\lambda}\right)  =\mathcal{A}%
\mathbf{G}_{\lambda}$. This proves Proposition \ref{prop.specht.AWlam.ker}.
\end{proof}

More subtle is the surjectivity of $\rho_{\lambda}$:

\begin{theorem}
\label{thm.specht.AWlam.sur}Let $\lambda$ be a partition of $n$. Assume that
$h^{\lambda}$ is invertible in $\mathbf{k}$. Let $\rho_{\lambda}%
:\mathbf{k}\left[  S_{n}\right]  \rightarrow\operatorname*{End}%
\nolimits_{\mathbf{k}}\left(  \mathcal{S}^{\lambda}\right)  $ be the curried
form of the left $\mathbf{k}\left[  S_{n}\right]  $-action on the Specht
module $\mathcal{S}^{\lambda}$ (as in Theorem \ref{thm.specht.AW}). Then:
\medskip

\textbf{(a)} The map $\rho_{\lambda}$ is surjective. \medskip

\textbf{(b)} The restriction of $\rho_{\lambda}$ to $\mathcal{A}%
\mathbf{E}_{\lambda}$ is a $\mathbf{k}$-algebra isomorphism from
$\mathcal{A}\mathbf{E}_{\lambda}$ to $\operatorname*{End}\nolimits_{\mathbf{k}%
}\left(  \mathcal{S}^{\lambda}\right)  $. (See Corollary
\ref{cor.spechtmod.Elam.subalg} \textbf{(b)} for why $\mathcal{A}%
\mathbf{E}_{\lambda}$ is a $\mathbf{k}$-algebra.)
\end{theorem}

We will prove Theorem \ref{thm.specht.AWlam.sur} in full in the next section
(Subsection \ref{subsec.specht.nat-basis.AWlam}), but it is already not hard
to prove it under the slightly more restrictive requirement that $n!$ be invertible:

\begin{exercise}
\fbox{3} Prove Theorem \ref{thm.specht.AWlam.sur} in the case when $n!$ is
invertible in $\mathbf{k}$. \medskip

[\textbf{Hint:} Show first that if $e$ is a central idempotent element of a
$\mathbf{k}$-algebra $R$, then $Re\cong R/R\left(  1-e\right)  $ as
$\mathbf{k}$-algebras.]
\end{exercise}

\subsection{\label{sec.specht.nat-basis}The Young symmetrizer basis of
$\mathbf{k}\left[  S_{n}\right]  $}

Corollary \ref{cor.spechtmod.sumflam2} shows that $n!$ equals the number of
\textbf{pairs} of standard tableaux of shape $Y\left(  \lambda\right)  $,
summed over all partitions $\lambda$ of $n$. (Indeed, for any given partition
$\lambda$ of $n$, the number of such pairs is $\left(  f^{\lambda}\right)
^{2}$.) The Artin--Wedderburn theorem (Theorem \ref{thm.AWS.demo}
\textbf{(a)}) provides an algebraic interpretation of this fact. Another
linear-algebraic consequence of this fact is that $\mathcal{A}$ (as a
$\mathbf{k}$-module) has a basis indexed by such pairs.

Of course, the mere existence of such a basis is obvious (since the number of
such pairs is $n!$, and clearly $\mathcal{A}$ has a basis of size $n!$). But
it is reasonable to wonder whether there is a \textquotedblleft
natural\textquotedblright\ basis of this kind, i.e., a basis of
combinatorially defined elements that are meaningfully related to pairs of
standard tableaux.

The answer is \textquotedblleft yes\textquotedblright, at least when $n!$ is
invertible in $\mathbf{k}$. In this section, we will construct such a basis,
which we call the \emph{Young symmetrizer basis} (Corollary
\ref{cor.specht.A.nat-basis}). It consists of Young symmetrizers
$\mathbf{E}_{P}$, slightly modified to depend on two $n$-tableaux $P$ and $Q$
rather than a single $n$-tableau $P$. These tweaked Young symmetrizers are
called $\mathbf{E}_{P,Q}$, and we will define them shortly.

Once the elements $\mathbf{E}_{P,Q}$ are defined, the basis will be easy to
find: It is the family $\left(  \mathbf{E}_{U,V}\right)  _{\lambda\text{ is a
partition of }n\text{, and }U,V\in\operatorname*{SYT}\left(  \lambda\right)
}$, where $\operatorname*{SYT}\left(  \lambda\right)  $ denotes the set of all
standard $n$-tableaux of shape $Y\left(  \lambda\right)  $. Proving that this
family really is a basis of $\mathcal{A}$ (when $n!$ is invertible) is not all
that easy. The journey is long, although none of the steps involved is
particularly difficult.\footnote{See \cite[\S 6]{Nazaro07} for an alternative
proof, although likely not very different from ours.}

The standard basis theorem for Specht modules will play a role in the proof,
but we will also need to understand the behavior of the elements
$\mathbf{E}_{P,Q}$ under transposing the tableaux $P$ and $Q$ (that is,
reflecting them across the northwest-to-southeast diagonal). This behavior
will be described in (\ref{prop.specht.EPQ.r}) below, which is just one of a
number of auxiliary results that we will use as stepping stones towards the
proof. Another remarkable result we will find along the way (Corollary
\ref{cor.specht.A.decomp}) is an explicit decomposition of $\mathcal{A}$ into
a direct sum of submodules of the form $\mathcal{A}\mathbf{E}_{T}$, each of
which is isomorphic to a Specht module $\mathcal{S}^{\lambda}$. Some further
comments about the Young symmetrizer basis will be made in Section
\ref{sec.bas.ysb}.

One of our auxiliary results (Corollary \ref{cor.specht.AElam.decomp1}) will
then be used to prove Theorem \ref{thm.specht.AWlam.sur} (as promised in the
previous section).

We shall use Definition \ref{def.spechtmod.flam} throughout this section.

\subsubsection{The $w_{P,Q}$ and the $\mathbf{E}_{P,Q}$}

We begin our journey with the definitions of the elements $w_{P,Q}$ and
$\mathbf{E}_{P,Q}$.

The following lemma is nearly obvious, and has been implicitly used a few
times already (see, e.g., Proposition \ref{prop.tableau.Sn-act.1} \textbf{(b)}):

\begin{lemma}
\label{lem.specht.wPQ.wd}Let $D$ be any diagram. Let $P$ and $Q$ be two
$n$-tableaux of shape $D$. Then, there is a unique permutation $w\in S_{n}$
that satisfies $wQ=P$.
\end{lemma}

\begin{fineprint}
\begin{proof}
We know that $P$ is an $n$-tableau of shape $D$. In other words, $P$ is an
injective map from $D$ to $\left\{  1,2,3,\ldots\right\}  $ whose image is
$\left[  n\right]  $. Thus, we can view $P$ as a bijection from $D$ to
$\left[  n\right]  $. Likewise, we can view $Q$ as a bijection from $D$ to
$\left[  n\right]  $. Hence, the inverse $Q^{-1}$ of $Q$ is well-defined and
is a bijection from $\left[  n\right]  $ to $D$. Therefore, $P\circ Q^{-1}$ is
a well-defined bijection from $\left[  n\right]  $ to $\left[  n\right]  $. In
other words, $P\circ Q^{-1}$ is a well-defined permutation in $S_{n}$.

But any permutation $w\in S_{n}$ satisfies $wQ=w\rightharpoonup Q=w\circ Q$
(by Definition \ref{def.tableau.Sn-act}). Thus, for any permutation $w\in
S_{n}$, we have the following chain of equivalences:%
\begin{align*}
\left(  wQ=P\right)  \  &  \Longleftrightarrow\ \left(  w\circ Q=P\right)
\ \ \ \ \ \ \ \ \ \ \left(  \text{since }wQ=w\circ Q\right) \\
&  \Longleftrightarrow\ \left(  w\circ Q=P\circ Q^{-1}\circ Q\right)
\ \ \ \ \ \ \ \ \ \ \left(  \text{since }P=P\circ Q^{-1}\circ Q\right) \\
&  \Longleftrightarrow\ \left(  w=P\circ Q^{-1}\right)
\end{align*}
(here, we have cancelled the $Q$ on both sides, since $Q$ is a bijection).
Hence, a permutation $w\in S_{n}$ satisfies $wQ=P$ if and only if it satisfies
$w=P\circ Q^{-1}$. In other words, there is exactly one permutation $w\in
S_{n}$ that satisfies $wQ=P$, namely the permutation $P\circ Q^{-1}$ (since we
already know that $P\circ Q^{-1}$ is a well-defined permutation in $S_{n}$).
This proves Lemma \ref{lem.specht.wPQ.wd}.
\end{proof}
\end{fineprint}

Let us give the unique permutation $w$ in Lemma \ref{lem.specht.wPQ.wd} a name:

\begin{definition}
\label{def.specht.wPQ}Let $D$ be any diagram. Let $P$ and $Q$ be two
$n$-tableaux of shape $D$. Lemma \ref{lem.specht.wPQ.wd} shows that there is a
unique permutation $w\in S_{n}$ that satisfies $wQ=P$. We shall denote this
permutation $w$ by $w_{P,Q}$.

Thus, $w_{P,Q}\in S_{n}$ is the permutation satisfying%
\begin{equation}
w_{P,Q}Q=P. \label{eq.def.specht.wPQ.eq}%
\end{equation}

\end{definition}

\begin{example}
Let $n=4$ and $D=Y\left(  2,2\right)  $ and $P=\ytableaushort{14,23}$ and
$Q=\ytableaushort{21,43}\ \ $. Then, $w_{P,Q}=\operatorname*{oln}\left(
4132\right)  $, since $\operatorname*{oln}\left(  4132\right)  Q=P$.
\end{example}

We now define a slight generalization of the Young symmetrizers $\mathbf{E}%
_{T}$:

\begin{definition}
\label{def.specht.EPQ}Let $D$ be any diagram. Let $P$ and $Q$ be two
$n$-tableaux of shape $D$. Then, we define an element $\mathbf{E}_{P,Q}%
\in\mathbf{k}\left[  S_{n}\right]  $ by%
\[
\mathbf{E}_{P,Q}:=\nabla_{\operatorname*{Col}P}^{-}w_{P,Q}\nabla
_{\operatorname*{Row}Q}.
\]
(See Definition \ref{def.symmetrizers.symmetrizers} and Definition
\ref{def.specht.wPQ} for the meaning of the terms used here.)
\end{definition}

\begin{example}
\label{exa.specht.EPQ.21}Let $n=3$ and $D=Y\left(  2,1\right)  $. Then, for
the two $n$-tableaux $31\backslash\backslash2$ and $12\backslash\backslash3$
of shape $D$, we have%
\begin{align*}
\mathbf{E}_{31\backslash\backslash2,\ 12\backslash\backslash3}  &
=\underbrace{\nabla_{\operatorname*{Col}\left(  31\backslash\backslash
2\right)  }^{-}}_{\substack{=1-t_{3,2}\\=1-s_{2}}}\underbrace{w_{31\backslash
\backslash2,\ 12\backslash\backslash3}}_{\substack{=\operatorname*{oln}\left(
312\right)  \\=\operatorname*{cyc}\nolimits_{1,3,2}\\=s_{2}s_{1}%
}}\underbrace{\nabla_{\operatorname*{Row}\left(  12\backslash\backslash
3\right)  }}_{\substack{=1+t_{1,2}\\=1+s_{1}}}\\
&  =\underbrace{\left(  1-s_{2}\right)  s_{2}}_{=s_{2}-1}\underbrace{s_{1}%
\left(  1+s_{1}\right)  }_{=s_{1}+1}=\left(  s_{2}-1\right)  \left(
s_{1}+1\right)  =s_{2}s_{1}-s_{1}+s_{2}-1.
\end{align*}
The two standard tableaux of shape $D$ are $12\backslash\backslash3$ and
$13\backslash\backslash2$. The corresponding $\mathbf{E}_{P,Q}$ are%
\begin{align*}
\mathbf{E}_{12\backslash\backslash3,\ 12\backslash\backslash3}  &
=1-t_{1,3}+t_{1,2}-\operatorname*{cyc}\nolimits_{1,2,3};\\
\mathbf{E}_{12\backslash\backslash3,\ 13\backslash\backslash2}  &
=t_{2,3}-\operatorname*{cyc}\nolimits_{1,3,2}+\operatorname*{cyc}%
\nolimits_{1,2,3}-\,t_{1,2};\\
\mathbf{E}_{13\backslash\backslash2,\ 12\backslash\backslash3}  &
=t_{2,3}-\operatorname*{cyc}\nolimits_{1,2,3}+\operatorname*{cyc}%
\nolimits_{1,3,2}-\,t_{1,3};\\
\mathbf{E}_{13\backslash\backslash2,\ 13\backslash\backslash2}  &
=1-t_{1,2}+t_{1,3}-\operatorname*{cyc}\nolimits_{1,3,2}.
\end{align*}

\end{example}

\begin{remark}
The $w_{P,Q}$ and the $\mathbf{E}_{P,Q}$ go all the way back to Young's 1928
paper \cite[\textit{On Quantitative Substitutional Analysis (Third Paper)},
\S 4]{Young77}. A more recent source is Rutherford's \cite[\S 8 and
\S 11]{Ruther48}, which studies them in the case when the shape $D$ is the
Young diagram $Y\left(  \alpha\right)  $ of a partition $\alpha$. In this
setting, Rutherford numbers all $n$-tableaux of shape $D$ by $S_{1}^{\alpha
},S_{2}^{\alpha},\ldots,S_{n!}^{\alpha}$ (in an arbitrary but fixed order),
and uses the notations $\sigma_{rs}^{\alpha}$ and $E_{rs}^{\alpha}$ for what
we call $w_{P,Q}$ and $\mathbf{E}_{P,Q}$ when $P=S_{r}^{\alpha}$ and
$Q=S_{s}^{\alpha}$. (He often omits the superscript $\alpha$.)
\end{remark}

The elements $\mathbf{E}_{P,Q}$ we introduced in Definition
\ref{def.specht.EPQ} are very close to the Young symmetrizers $\mathbf{E}_{P}$
(which is why we had no need to introduce them until now). Indeed, we can
describe the former through the latter in two easy ways:

\begin{proposition}
\label{prop.specht.EPQ.EPEQ}Let $D$ be any diagram. Let $P$ and $Q$ be two
$n$-tableaux of shape $D$. Then,%
\begin{align}
\mathbf{E}_{P,Q}  &  =w_{P,Q}\mathbf{E}_{Q} \label{eq.prop.specht.EPQ.EPEQ.wE}%
\\
&  =\mathbf{E}_{P}w_{P,Q}. \label{eq.prop.specht.EPQ.EPEQ.Ew}%
\end{align}

\end{proposition}

\begin{proof}
Let us denote the permutation $w_{P,Q}\in S_{n}$ by $w$. Then, $w=w_{P,Q}$, so
that $wQ=w_{P,Q}Q=P$ (by (\ref{eq.def.specht.wPQ.eq})). In other words,
$w\rightharpoonup Q=P$ (since $wQ$ is just a shorthand for $w\rightharpoonup
Q$). But Proposition \ref{prop.symmetrizers.conj} (applied to $T=Q$) yields
$\nabla_{\operatorname*{Row}\left(  w\rightharpoonup Q\right)  }%
=w\nabla_{\operatorname*{Row}Q}w^{-1}$ and $\nabla_{\operatorname*{Col}\left(
w\rightharpoonup Q\right)  }^{-}=w\nabla_{\operatorname*{Col}Q}^{-}w^{-1}$. In
view of $w\rightharpoonup Q=P$, we can rewrite these two equalities as%
\[
\nabla_{\operatorname*{Row}P}=w\nabla_{\operatorname*{Row}Q}w^{-1}%
\ \ \ \ \ \ \ \ \ \ \text{and}\ \ \ \ \ \ \ \ \ \ \nabla_{\operatorname*{Col}%
P}^{-}=w\nabla_{\operatorname*{Col}Q}^{-}w^{-1}.
\]

Now, the definition of $\mathbf{E}_{Q}$ yields $\mathbf{E}_{Q}=\nabla
_{\operatorname*{Col}Q}^{-}\nabla_{\operatorname*{Row}Q}$. Likewise, we obtain
$\mathbf{E}_{P}=\nabla_{\operatorname*{Col}P}^{-}\nabla_{\operatorname*{Row}%
P}$. But the definition of $\mathbf{E}_{P,Q}$ yields%
\begin{align*}
\mathbf{E}_{P,Q}  &  =\underbrace{\nabla_{\operatorname*{Col}P}^{-}}%
_{=w\nabla_{\operatorname*{Col}Q}^{-}w^{-1}}\underbrace{w_{P,Q}}_{=w}%
\nabla_{\operatorname*{Row}Q}\\
&  =w\nabla_{\operatorname*{Col}Q}^{-}\underbrace{w^{-1}w}_{=1}\nabla
_{\operatorname*{Row}Q}=w\underbrace{\nabla_{\operatorname*{Col}Q}^{-}%
\nabla_{\operatorname*{Row}Q}}_{=\mathbf{E}_{Q}}=\underbrace{w}_{=w_{P,Q}%
}\mathbf{E}_{Q}=w_{P,Q}\mathbf{E}_{Q}%
\end{align*}
and%
\begin{align*}
\mathbf{E}_{P,Q}  &  =\nabla_{\operatorname*{Col}P}^{-}\underbrace{w_{P,Q}%
}_{=w}\underbrace{\nabla_{\operatorname*{Row}Q}}_{\substack{=w^{-1}%
\nabla_{\operatorname*{Row}P}w\\\text{(since }\nabla_{\operatorname*{Row}%
P}=w\nabla_{\operatorname*{Row}Q}w^{-1}\text{)}}}\\
&  =\nabla_{\operatorname*{Col}P}^{-}\underbrace{ww^{-1}}_{=1}\nabla
_{\operatorname*{Row}P}w=\underbrace{\nabla_{\operatorname*{Col}P}^{-}%
\nabla_{\operatorname*{Row}P}}_{=\mathbf{E}_{P}}w=\mathbf{E}_{P}%
\underbrace{w}_{=w_{P,Q}}=\mathbf{E}_{P}w_{P,Q}.
\end{align*}
Together, these two equalities prove Proposition \ref{prop.specht.EPQ.EPEQ}.
\end{proof}

\begin{exercise}
\label{exe.specht.EPQ.uEUVv}\fbox{1} Let $D$ be any diagram. Let $P$ and $Q$
be two $n$-tableaux of shape $D$. Let $u,v\in S_{n}$. Prove that
\[
u\mathbf{E}_{P,Q}v=\mathbf{E}_{uP,v^{-1}Q}.
\]

\end{exercise}

\begin{exercise}
\label{exe.specht.EPQ.EPQEQR}\fbox{1} Let $\lambda$ be a partition of $n$. Let
$P$, $Q$ and $R$ be three $n$-tableaux of shape $\lambda$. Show that
$\mathbf{E}_{P,Q}\mathbf{E}_{Q,R}=h^{\lambda}\mathbf{E}_{P,R}$.
\end{exercise}

The Young symmetrizers $\mathbf{E}_{P}$ are particular cases of the
newly-defined elements $\mathbf{E}_{P,P}$:

\begin{proposition}
\label{prop.specht.EPQ.EP}Let $D$ be any diagram. Let $P$ be an $n$-tableau of
shape $D$. Then, $w_{P,P}=\operatorname*{id}$ and $\mathbf{E}_{P,P}%
=\mathbf{E}_{P}$.
\end{proposition}

\begin{proof}
Definition \ref{def.specht.wPQ} tells us that $w_{P,P}$ is defined as the
unique permutation $w\in S_{n}$ that satisfies $wP=P$. Hence, if some
permutation $w\in S_{n}$ satisfies $wP=P$, then $w=w_{P,P}$ (by the uniqueness
claim in the preceding sentence). Applying this to $w=\operatorname*{id}$, we
obtain $\operatorname*{id}=w_{P,P}$ (since $\operatorname*{id}P=P$). Thus, we
conclude that $w_{P,P}=\operatorname*{id}$.

Furthermore, Definition \ref{def.specht.EPQ} yields $\mathbf{E}_{P,P}%
=\nabla_{\operatorname*{Col}P}^{-}\underbrace{w_{P,P}}_{=\operatorname*{id}%
}\nabla_{\operatorname*{Row}P}=\nabla_{\operatorname*{Col}P}^{-}%
\nabla_{\operatorname*{Row}P}=\mathbf{E}_{P}$ (since $\mathbf{E}_{P}$ was
defined to be $\nabla_{\operatorname*{Col}P}^{-}\nabla_{\operatorname*{Row}P}%
$). The proof of Proposition \ref{prop.specht.EPQ.EP} is thus complete.
\end{proof}

One reason we introduced the $\mathbf{E}_{P,Q}$'s is that they form a basis of
the left ideal avatar $\mathcal{A}\mathbf{E}_{T}$ of a skew Specht module
$\mathcal{S}^{Y\left(  \lambda/\mu\right)  }$:

\begin{proposition}
\label{prop.specht.EPQ.basis-AET}Let $D$ be a skew Young diagram. Let $T$ be
an $n$-tableau of shape $D$. Let $\operatorname*{SYT}\left(  D\right)  $ be
the set of all standard $n$-tableaux of shape $D$. Then, the family $\left(
\mathbf{E}_{P,T}\right)  _{P\in\operatorname*{SYT}\left(  D\right)  }$ is a
basis of the $\mathbf{k}$-module $\mathcal{A}\mathbf{E}_{T}$.
\end{proposition}

\begin{proof}
This is essentially just the standard basis theorem (Theorem
\ref{thm.spechtmod.basis}), after applying the isomorphism $\mathcal{S}%
^{D}\rightarrow\mathcal{A}\mathbf{E}_{T}$. In more detail:

We know that $T$ is an $n$-tableau of shape $D$, thus a bijection from $D$ to
$\left[  n\right]  $. Hence, $\left\vert D\right\vert =\left\vert \left[
n\right]  \right\vert =n$.

Theorem \ref{thm.spechtmod.basis} says that the standard polytabloids
$\mathbf{e}_{P}$ (that is, the $\mathbf{e}_{P}$ where $P$ ranges over all
standard tableaux of shape $D$) form a basis of the $\mathbf{k}$-module
$\mathcal{S}^{D}$. In other words, the $\mathbf{e}_{P}$ for $P\in
\operatorname*{SYT}\left(  D\right)  $ form a basis of the $\mathbf{k}$-module
$\mathcal{S}^{D}$ (since $\operatorname*{SYT}\left(  D\right)  $ is the set of
all standard $n$-tableaux of shape $D$). In other words, the family $\left(
\mathbf{e}_{P}\right)  _{P\in\operatorname*{SYT}\left(  D\right)  }$ is a
basis of the $\mathbf{k}$-module $\mathcal{S}^{D}$.

Corollary \ref{cor.spechtmod.leftideal.S} shows that there is a left
$\mathbf{k}\left[  S_{n}\right]  $-module isomorphism%
\[
\overline{\alpha}:\mathbf{k}\left[  S_{n}\right]  \cdot\nabla
_{\operatorname*{Col}T}^{-}\nabla_{\operatorname*{Row}T}\rightarrow
\mathcal{S}^{D}%
\]
that sends $\nabla_{\operatorname*{Col}T}^{-}\nabla_{\operatorname*{Row}T}$ to
$\mathbf{e}_{T}$. In other words, there is a left $\mathcal{A}$-module
isomorphism%
\[
\overline{\alpha}:\mathcal{A}\mathbf{E}_{T}\rightarrow\mathcal{S}^{D}%
\]
that sends $\mathbf{E}_{T}$ to $\mathbf{e}_{T}$ (since $\mathbf{k}\left[
S_{n}\right]  =\mathcal{A}$ and $\nabla_{\operatorname*{Col}T}^{-}%
\nabla_{\operatorname*{Row}T}=\mathbf{E}_{T}$). Consider this isomorphism
$\overline{\alpha}$. Thus, $\overline{\alpha}\left(  \mathbf{E}_{T}\right)
=\mathbf{e}_{T}$ (since $\overline{\alpha}$ sends $\mathbf{E}_{T}$ to
$\mathbf{e}_{T}$).

Next, let $P\in\operatorname*{SYT}\left(  D\right)  $. Then,
(\ref{eq.prop.specht.EPQ.EPEQ.wE}) (applied to $Q=T$) yields $\mathbf{E}%
_{P,T}=w_{P,T}\mathbf{E}_{T}$. Hence,%
\begin{align}
\overline{\alpha}\left(  \mathbf{E}_{P,T}\right)   &  =\overline{\alpha
}\left(  w_{P,T}\mathbf{E}_{T}\right)  =w_{P,T}\underbrace{\overline{\alpha
}\left(  \mathbf{E}_{T}\right)  }_{=\mathbf{e}_{T}}\ \ \ \ \ \ \ \ \ \ \left(
\text{since }\overline{\alpha}\text{ is a left }\mathcal{A}\text{-module
morphism}\right) \nonumber\\
&  =w_{P,T}\mathbf{e}_{T}=\mathbf{e}_{w_{P,T}T}
\label{pf.prop.specht.EPQ.basis-AET.6}%
\end{align}
(since Lemma \ref{lem.spechtmod.submod} \textbf{(a)} (applied to $u=w_{P,T}$)
yields $\mathbf{e}_{w_{P,T}T}=w_{P,T}\mathbf{e}_{T}$). But
(\ref{eq.def.specht.wPQ.eq}) (applied to $Q=T$) yields $w_{P,T}T=P$. Hence, we
can rewrite (\ref{pf.prop.specht.EPQ.basis-AET.6}) as $\overline{\alpha
}\left(  \mathbf{E}_{P,T}\right)  =\mathbf{e}_{P}$. Therefore, $\overline
{\alpha}^{-1}\left(  \mathbf{e}_{P}\right)  =\mathbf{E}_{P,T}$ (since
$\overline{\alpha}$ is an isomorphism and thus invertible).

Forget that we fixed $P$. We thus have shown that
\begin{equation}
\overline{\alpha}^{-1}\left(  \mathbf{e}_{P}\right)  =\mathbf{E}%
_{P,T}\ \ \ \ \ \ \ \ \ \ \text{for each }P\in\operatorname*{SYT}\left(
D\right)  . \label{pf.prop.specht.EPQ.basis-AET.8pre}%
\end{equation}
In other words,
\begin{equation}
\left(  \overline{\alpha}^{-1}\left(  \mathbf{e}_{P}\right)  \right)
_{P\in\operatorname*{SYT}\left(  D\right)  }=\left(  \mathbf{E}_{P,T}\right)
_{P\in\operatorname*{SYT}\left(  D\right)  }.
\label{pf.prop.specht.EPQ.basis-AET.8}%
\end{equation}

But the map $\overline{\alpha}:\mathcal{A}\mathbf{E}_{T}\rightarrow
\mathcal{S}^{D}$ is a $\mathbf{k}$-module isomorphism (since it is a left
$\mathcal{A}$-module isomorphism). Hence, it has an inverse $\overline{\alpha
}^{-1}:\mathcal{S}^{D}\rightarrow\mathcal{A}\mathbf{E}_{T}$, which is also a
$\mathbf{k}$-module isomorphism. It is well-known that a $\mathbf{k}$-module
isomorphism sends any basis of its domain to a basis of its target. Applying
this to the $\mathbf{k}$-module isomorphism $\overline{\alpha}^{-1}$, we see
that $\overline{\alpha}^{-1}$ sends any basis of $\mathcal{S}^{D}$ to a basis
of $\mathcal{A}\mathbf{E}_{T}$. In other words, if $\left(  v_{i}\right)
_{i\in I}$ is a basis of $\mathcal{S}^{D}$, then $\left(  \overline{\alpha
}^{-1}\left(  v_{i}\right)  \right)  _{i\in I}$ is a basis of $\mathcal{A}%
\mathbf{E}_{T}$.

Applying this to $\left(  v_{i}\right)  _{i\in I}=\left(  \mathbf{e}%
_{P}\right)  _{P\in\operatorname*{SYT}\left(  D\right)  }$, we conclude that
the family \newline$\left(  \overline{\alpha}^{-1}\left(  \mathbf{e}%
_{P}\right)  \right)  _{P\in\operatorname*{SYT}\left(  D\right)  }$ is a basis
of $\mathcal{A}\mathbf{E}_{T}$ (since $\left(  \mathbf{e}_{P}\right)
_{P\in\operatorname*{SYT}\left(  D\right)  }$ is a basis of the $\mathbf{k}%
$-module $\mathcal{S}^{D}$). In view of (\ref{pf.prop.specht.EPQ.basis-AET.8}%
), we can rewrite this as follows: The family $\left(  \mathbf{E}%
_{P,T}\right)  _{P\in\operatorname*{SYT}\left(  D\right)  }$ is a basis of
$\mathcal{A}\mathbf{E}_{T}$. This proves Proposition
\ref{prop.specht.EPQ.basis-AET}.
\end{proof}

Thus, we found an explicit basis of $\mathcal{A}\mathbf{E}_{T}$. This is not
very surprising, as $\mathcal{A}\mathbf{E}_{T}$ is isomorphic to the Specht
module $\mathcal{S}^{D}$. What is more interesting is to find an explicit
basis for $\mathcal{A}\mathbf{E}_{\lambda}$ when $\lambda$ is a partition of
$n$. This will take us the next few pages (and require some invertibility
conditions on $\mathbf{k}$). As an intermediate step, let us find a basis of
$\mathbf{E}_{T}\mathcal{A}$. This is quite similar to $\mathcal{A}%
\mathbf{E}_{T}$, and indeed we will derive a basis of $\mathbf{E}%
_{T}\mathcal{A}$ from our above basis of $\mathcal{A}\mathbf{E}_{T}$. To do
so, we shall now briefly study transposition of Young tableaux, which will
help us relate the Young symmetrizers $\mathbf{E}_{T}$ to the antipode of
$\mathcal{A}$.

\subsubsection{Transposes of Young tableaux}

In Theorem \ref{thm.partitions.conj}, we defined the map $\mathbf{r}%
:\mathbb{Z}^{2}\rightarrow\mathbb{Z}^{2}$, which sends each cell $\left(
i,j\right)  \in\mathbb{Z}^{2}$ to $\left(  j,i\right)  $. Geometrically, this
map is the reflection across the northwest-to-southeast diagonal. We observed
that the image of a Young diagram $Y\left(  \lambda\right)  $ under this map
$\mathbf{r}$ is a Young diagram again, namely the Young diagram of the
partition $\lambda^{t}$ known as the conjugate of $\lambda$ (see Theorem
\ref{thm.partitions.conj} \textbf{(a)}). It is also clear that the map
$\mathbf{r}$ is an involution (i.e., it is inverse to itself, i.e., it
satisfies $\mathbf{r}\circ\mathbf{r}=\operatorname*{id}$), and thus is a bijection.

We shall now use the map $\mathbf{r}$ to reflect a tableau across the
northwest-to-southeast diagonal:

\begin{definition}
\label{def.tableaux.r}Let $D$ be any diagram. Consider the map $\mathbf{r}%
:\mathbb{Z}^{2}\rightarrow\mathbb{Z}^{2}$ from Theorem
\ref{thm.partitions.conj}. Then, $\mathbf{r}\left(  \mathbf{r}\left(
D\right)  \right)  =D$ (since $\mathbf{r}$ is an involution). Hence, we can
restrict the bijection $\mathbf{r}$ to obtain a bijection $\left.
\mathbf{r}\mid_{\mathbf{r}\left(  D\right)  }\right.  :\mathbf{r}\left(
D\right)  \rightarrow D$.

Now, let $T$ be an $n$-tableau of shape $D$. Thus, $T$ is a bijection from $D$
to $\left[  n\right]  $. Composing this bijection $T$ with the bijection
$\left.  \mathbf{r}\mid_{\mathbf{r}\left(  D\right)  }\right.  :\mathbf{r}%
\left(  D\right)  \rightarrow D$, we obtain a bijection $T\circ\left(
\mathbf{r}\mid_{\mathbf{r}\left(  D\right)  }\right)  :\mathbf{r}\left(
D\right)  \rightarrow\left[  n\right]  $, which we shall call $T\mathbf{r}$.
Thus, $T\mathbf{r}$ is an $n$-tableau of shape $\mathbf{r}\left(  D\right)  $.
Explicitly, the entries of this $n$-tableau $T\mathbf{r}$ are given by%
\begin{align*}
\left(  T\mathbf{r}\right)  \left(  i,j\right)   &  =T\left(  \left(
\mathbf{r}\mid_{\mathbf{r}\left(  D\right)  }\right)  \left(  i,j\right)
\right)  \ \ \ \ \ \ \ \ \ \ \left(  \text{since }T\mathbf{r}=T\circ\left(
\mathbf{r}\mid_{\mathbf{r}\left(  D\right)  }\right)  \right) \\
&  =T\left(  j,i\right)  \ \ \ \ \ \ \ \ \ \ \left(
\begin{array}
[c]{c}%
\text{since }\left(  \mathbf{r}\mid_{\mathbf{r}\left(  D\right)  }\right)
\left(  i,j\right)  =\mathbf{r}\left(  i,j\right)  =\left(  j,i\right) \\
\text{(by the definition of }\mathbf{r}\text{)}%
\end{array}
\right) \\
&  \ \ \ \ \ \ \ \ \ \ \ \ \ \ \ \ \ \ \ \ \text{for each }\left(  i,j\right)
\in\mathbf{r}\left(  D\right)  .
\end{align*}
Thus, informally speaking, the $n$-tableau $T\mathbf{r}$ is obtained by
reflecting $T$ across the northwest-to-southeast diagonal.\footnotemark\ We
call $T\mathbf{r}$ the \emph{transpose} of the $n$-tableau $T$.
\end{definition}

\footnotetext{No, you don't replace $8$'s by $\infty$'s.}

\begin{example}
Let $D=Y\left(  \left(  4,4,2\right)  /\left(  2\right)  \right)  $ and
$T=\ytableaushort{\none\none36,1457,28}\ \ $. Then, $T\mathbf{r}%
=\ytableaushort{\none12,\none48,35,67}\ \ $.
\end{example}

\begin{proposition}
\label{prop.tableaux.r}Let $D$ be any diagram. Then: \medskip

\textbf{(a)} If $T$ is a standard $n$-tableau of shape $D$, then $T\mathbf{r}$
is a standard $n$-tableau of shape $\mathbf{r}\left(  D\right)  $. \medskip

\textbf{(b)} If $T$ is any $n$-tableau of shape $D$, then $\left(
T\mathbf{r}\right)  \mathbf{r}=T$. \medskip

\textbf{(c)} For any diagram $E$, let $\operatorname*{SYT}\left(  E\right)  $
be the set of all standard $n$-tableaux of shape $E$. Then, the map%
\begin{align*}
\operatorname*{SYT}\left(  D\right)   &  \rightarrow\operatorname*{SYT}\left(
\mathbf{r}\left(  D\right)  \right)  ,\\
T  &  \mapsto T\mathbf{r}%
\end{align*}
is a bijection. \medskip

\textbf{(d)} If $T$ is any $n$-tableau of shape $D$, and if $w\in S_{n}$ is
arbitrary, then $\left(  wT\right)  \mathbf{r}=w\left(  T\mathbf{r}\right)  $.
\medskip

\textbf{(e)} If $T$ is any $n$-tableau of shape $D$, then $\mathcal{R}\left(
T\right)  =\mathcal{C}\left(  T\mathbf{r}\right)  $ and $\mathcal{C}\left(
T\right)  =\mathcal{R}\left(  T\mathbf{r}\right)  $.
\end{proposition}

The proof of Proposition \ref{prop.tableaux.r} is essentially trivial,
exploiting obvious symmetries of the definitions and the fact that
$\mathbf{r}$ is an involution. Thus we are banishing it to the appendix
(Section \ref{sec.details.specht.nat-basis}).

\begin{corollary}
\label{cor.tableaux.flamt}Let $\lambda$ be a partition of $n$. Then,
$f^{\lambda}=f^{\lambda^{t}}$. (See Theorem \ref{thm.partitions.conj}
\textbf{(a)} for the meaning of $\lambda^{t}$.)
\end{corollary}

\begin{fineprint}
\begin{proof}
We have $\left\vert Y\left(  \lambda\right)  \right\vert =\left\vert
\lambda\right\vert =n$ (since $\lambda$ is a partition of $n$). Hence,
Proposition \ref{prop.tableau.std-n} (applied to $D=Y\left(  \lambda\right)
$) shows that the standard tableaux of shape $Y\left(  \lambda\right)  $ are
precisely the standard $n$-tableaux of shape $Y\left(  \lambda\right)  $.

We shall use the notation $\operatorname*{SYT}\left(  E\right)  $ from
Proposition \ref{prop.tableaux.r} \textbf{(c)}. Applying Proposition
\ref{prop.tableaux.r} \textbf{(c)} to $D=Y\left(  \lambda\right)  $, we obtain
that the map%
\begin{align*}
\operatorname*{SYT}\left(  Y\left(  \lambda\right)  \right)   &
\rightarrow\operatorname*{SYT}\left(  \mathbf{r}\left(  Y\left(
\lambda\right)  \right)  \right)  ,\\
T  &  \mapsto T\mathbf{r}%
\end{align*}
is a bijection. Hence, $\left\vert \operatorname*{SYT}\left(  Y\left(
\lambda\right)  \right)  \right\vert =\left\vert \operatorname*{SYT}\left(
\mathbf{r}\left(  Y\left(  \lambda\right)  \right)  \right)  \right\vert $. In
other words, $\left\vert \operatorname*{SYT}\left(  Y\left(  \lambda\right)
\right)  \right\vert =\left\vert \operatorname*{SYT}\left(  Y\left(
\lambda^{t}\right)  \right)  \right\vert $ (since the definition of
$\lambda^{t}$ yields $Y\left(  \lambda^{t}\right)  =\mathbf{r}\left(  Y\left(
\lambda\right)  \right)  $).

However, $f^{\lambda}$ is defined as the \# of standard tableaux of shape
$Y\left(  \lambda\right)  $. In other words, $f^{\lambda}$ is the \# of
standard $n$-tableaux of shape $Y\left(  \lambda\right)  $ (since the standard
tableaux of shape $Y\left(  \lambda\right)  $ are precisely the standard
$n$-tableaux of shape $Y\left(  \lambda\right)  $). In other words,
$f^{\lambda}=\left\vert \operatorname*{SYT}\left(  Y\left(  \lambda\right)
\right)  \right\vert $ (since $\operatorname*{SYT}\left(  Y\left(
\lambda\right)  \right)  $ is defined as the set of all standard $n$-tableaux
of shape $Y\left(  \lambda\right)  $). The same argument (applied to
$\lambda^{t}$ instead of $\lambda$) yields $f^{\lambda^{t}}=\left\vert
\operatorname*{SYT}\left(  Y\left(  \lambda^{t}\right)  \right)  \right\vert $
(since Theorem \ref{thm.partitions.conj} \textbf{(c)} yields $\left\vert
\lambda^{t}\right\vert =\left\vert \lambda\right\vert =n$). Hence,
\[
f^{\lambda}=\left\vert \operatorname*{SYT}\left(  Y\left(  \lambda\right)
\right)  \right\vert =\left\vert \operatorname*{SYT}\left(  Y\left(
\lambda^{t}\right)  \right)  \right\vert =f^{\lambda^{t}}.
\]
This proves Corollary \ref{cor.tableaux.flamt}.
\end{proof}
\end{fineprint}

Now, recall two $\mathbf{k}$-linear maps from $\mathbf{k}\left[  S_{n}\right]
$ to $\mathbf{k}\left[  S_{n}\right]  $: the sign-twist
$T_{\operatorname*{sign}}$ defined in Definition \ref{def.Tsign.Tsign}, and
the antipode $S$ defined in Definition \ref{def.S.S}. Using these maps, we can
relate the Young symmetrizer of a given $n$-tableau $P$ with that of its
transpose $P\mathbf{r}$:

\begin{proposition}
\label{prop.specht.ET.r}Let $D$ be any diagram. Let $P$ be any $n$-tableau of
shape $D$. Then,%
\begin{align}
\nabla_{\operatorname*{Row}\left(  P\mathbf{r}\right)  }  &
=T_{\operatorname*{sign}}\left(  \nabla_{\operatorname*{Col}P}^{-}\right)
;\label{eq.prop.specht.ET.r.1}\\
\nabla_{\operatorname*{Col}\left(  P\mathbf{r}\right)  }^{-}  &
=T_{\operatorname*{sign}}\left(  \nabla_{\operatorname*{Row}P}\right)
;\label{eq.prop.specht.ET.r.2}\\
\mathbf{E}_{P\mathbf{r}}  &  =T_{\operatorname*{sign}}\left(  S\left(
\mathbf{E}_{P}\right)  \right)  . \label{eq.prop.specht.ET.r.3}%
\end{align}

\end{proposition}

\begin{proof}
Proposition \ref{prop.tableaux.r} \textbf{(b)} (applied to $T=P$) yields
$\left(  P\mathbf{r}\right)  \mathbf{r}=P$. Proposition \ref{prop.tableaux.r}
\textbf{(e)} (applied to $T=P$) yields $\mathcal{R}\left(  P\right)
=\mathcal{C}\left(  P\mathbf{r}\right)  $ and $\mathcal{C}\left(  P\right)
=\mathcal{R}\left(  P\mathbf{r}\right)  $.

Now, the definition of $\nabla_{\operatorname*{Row}\left(  P\mathbf{r}\right)
}$ yields%
\begin{equation}
\nabla_{\operatorname*{Row}\left(  P\mathbf{r}\right)  }=\sum_{w\in
\mathcal{R}\left(  P\mathbf{r}\right)  }w=\sum_{w\in\mathcal{C}\left(
P\right)  }w \label{pf.prop.specht.ET.r.1}%
\end{equation}
(since $\mathcal{R}\left(  P\mathbf{r}\right)  =\mathcal{C}\left(  P\right)
$). On the other hand, the definition of $\nabla_{\operatorname*{Col}P}^{-}$
yields $\nabla_{\operatorname*{Col}P}^{-}=\sum_{w\in\mathcal{C}\left(
P\right)  }\left(  -1\right)  ^{w}w$. Applying the map
$T_{\operatorname*{sign}}$ to this equality, we obtain%
\begin{align*}
T_{\operatorname*{sign}}\left(  \nabla_{\operatorname*{Col}P}^{-}\right)   &
=T_{\operatorname*{sign}}\left(  \sum_{w\in\mathcal{C}\left(  P\right)
}\left(  -1\right)  ^{w}w\right) \\
&  =\sum_{w\in\mathcal{C}\left(  P\right)  }\left(  -1\right)  ^{w}%
\underbrace{T_{\operatorname*{sign}}\left(  w\right)  }_{\substack{=\left(
-1\right)  ^{w}w\\\text{(by the definition}\\\text{of }T_{\operatorname*{sign}%
}\text{)}}}\ \ \ \ \ \ \ \ \ \ \left(  \text{since the map }%
T_{\operatorname*{sign}}\text{ is }\mathbf{k}\text{-linear}\right) \\
&  =\sum_{w\in\mathcal{C}\left(  P\right)  }\underbrace{\left(  -1\right)
^{w}\left(  -1\right)  ^{w}}_{\substack{=\left(  \left(  -1\right)
^{w}\right)  ^{2}=1\\\text{(since }\left(  -1\right)  ^{w}\in\left\{
1,-1\right\}  \text{)}}}w=\sum_{w\in\mathcal{C}\left(  P\right)  }w.
\end{align*}
Comparing this with (\ref{pf.prop.specht.ET.r.1}), we obtain
\[
\nabla_{\operatorname*{Row}\left(  P\mathbf{r}\right)  }%
=T_{\operatorname*{sign}}\left(  \nabla_{\operatorname*{Col}P}^{-}\right)  .
\]
Thus, (\ref{eq.prop.specht.ET.r.1}) is proved.

Now, observe that $P\mathbf{r}$ is an $n$-tableau of shape $\mathbf{r}\left(
D\right)  $. Hence, applying (\ref{eq.prop.specht.ET.r.1}) to $\mathbf{r}%
\left(  D\right)  $ and $P\mathbf{r}$ instead of $D$ and $P$, we find%
\[
\nabla_{\operatorname*{Row}\left(  \left(  P\mathbf{r}\right)  \mathbf{r}%
\right)  }=T_{\operatorname*{sign}}\left(  \nabla_{\operatorname*{Col}\left(
P\mathbf{r}\right)  }^{-}\right)  .
\]
In view of $\left(  P\mathbf{r}\right)  \mathbf{r}=P$, we can rewrite this as%
\[
\nabla_{\operatorname*{Row}P}=T_{\operatorname*{sign}}\left(  \nabla
_{\operatorname*{Col}\left(  P\mathbf{r}\right)  }^{-}\right)  .
\]
Applying the map $T_{\operatorname*{sign}}$ to both sides of this equality, we
obtain%
\begin{align*}
T_{\operatorname*{sign}}\left(  \nabla_{\operatorname*{Row}P}\right)   &
=T_{\operatorname*{sign}}\left(  T_{\operatorname*{sign}}\left(
\nabla_{\operatorname*{Col}\left(  P\mathbf{r}\right)  }^{-}\right)  \right)
=\underbrace{\left(  T_{\operatorname*{sign}}\circ T_{\operatorname*{sign}%
}\right)  }_{\substack{=\operatorname*{id}\\\text{(by Theorem
\ref{thm.Tsign.auto} \textbf{(b)})}}}\left(  \nabla_{\operatorname*{Col}%
\left(  P\mathbf{r}\right)  }^{-}\right) \\
&  =\operatorname*{id}\left(  \nabla_{\operatorname*{Col}\left(
P\mathbf{r}\right)  }^{-}\right)  =\nabla_{\operatorname*{Col}\left(
P\mathbf{r}\right)  }^{-}.
\end{align*}
This proves (\ref{eq.prop.specht.ET.r.2}).

Finally, Proposition \ref{prop.symmetrizers.antipode} (applied to $T=P$)
yields that the antipode $S$ satisfies%
\[
S\left(  \nabla_{\operatorname*{Row}P}\right)  =\nabla_{\operatorname*{Row}%
P}\ \ \ \ \ \ \ \ \ \ \text{and}\ \ \ \ \ \ \ \ \ \ S\left(  \nabla
_{\operatorname*{Col}P}^{-}\right)  =\nabla_{\operatorname*{Col}P}^{-}.
\]
However, we have $\mathbf{E}_{P}=\nabla_{\operatorname*{Col}P}^{-}%
\nabla_{\operatorname*{Row}P}$ (by the definition of $\mathbf{E}_{P}$).
Applying the map $S$ to this equality, we obtain%
\[
S\left(  \mathbf{E}_{P}\right)  =S\left(  \nabla_{\operatorname*{Col}P}%
^{-}\nabla_{\operatorname*{Row}P}\right)  =S\left(  \nabla
_{\operatorname*{Row}P}\right)  S\left(  \nabla_{\operatorname*{Col}P}%
^{-}\right)
\]
(since Theorem \ref{thm.S.auto} \textbf{(a)} shows that $S$ is a $\mathbf{k}%
$-algebra anti-morphism). Thus,%
\[
S\left(  \mathbf{E}_{P}\right)  =\underbrace{S\left(  \nabla
_{\operatorname*{Row}P}\right)  }_{\substack{=\nabla_{\operatorname*{Row}P}%
}}\underbrace{S\left(  \nabla_{\operatorname*{Col}P}^{-}\right)
}_{\substack{=\nabla_{\operatorname*{Col}P}^{-}}}=\nabla_{\operatorname*{Row}%
P}\nabla_{\operatorname*{Col}P}^{-}.
\]
Applying the map $T_{\operatorname*{sign}}$ to this equality, we obtain%
\begin{align*}
T_{\operatorname*{sign}}\left(  S\left(  \mathbf{E}_{P}\right)  \right)   &
=T_{\operatorname*{sign}}\left(  \nabla_{\operatorname*{Row}P}\nabla
_{\operatorname*{Col}P}^{-}\right) \\
&  =\underbrace{T_{\operatorname*{sign}}\left(  \nabla_{\operatorname*{Row}%
P}\right)  }_{\substack{=\nabla_{\operatorname*{Col}\left(  P\mathbf{r}%
\right)  }^{-}\\\text{(by (\ref{eq.prop.specht.ET.r.2}))}}%
}\underbrace{T_{\operatorname*{sign}}\left(  \nabla_{\operatorname*{Col}P}%
^{-}\right)  }_{\substack{=\nabla_{\operatorname*{Row}\left(  P\mathbf{r}%
\right)  }\\\text{(by (\ref{eq.prop.specht.ET.r.1}))}}}\\
&  \ \ \ \ \ \ \ \ \ \ \ \ \ \ \ \ \ \ \ \ \left(
\begin{array}
[c]{c}%
\text{since Theorem \ref{thm.Tsign.auto} \textbf{(a)} shows}\\
\text{that }T_{\operatorname*{sign}}\text{ is a }\mathbf{k}\text{-algebra
morphism}%
\end{array}
\right) \\
&  =\nabla_{\operatorname*{Col}\left(  P\mathbf{r}\right)  }^{-}%
\nabla_{\operatorname*{Row}\left(  P\mathbf{r}\right)  }=\mathbf{E}%
_{P\mathbf{r}}%
\end{align*}
(since $\mathbf{E}_{P\mathbf{r}}$ is defined to be $\nabla
_{\operatorname*{Col}\left(  P\mathbf{r}\right)  }^{-}\nabla
_{\operatorname*{Row}\left(  P\mathbf{r}\right)  }$). This proves
(\ref{eq.prop.specht.ET.r.3}), thus concluding the proof of Proposition
\ref{prop.specht.ET.r}.
\end{proof}

\begin{exercise}
\fbox{2} Let $\lambda$ be a partition of $n$. Prove that
$T_{\operatorname*{sign}}\left(  \mathbf{E}_{\lambda}\right)  =\mathbf{E}%
_{\lambda^{t}}$.
\end{exercise}

\begin{exercise}
\label{exe.specht.FT.iso-as-k-mod}\fbox{2} Let $D$ be a skew Young diagram.
Let $T$ be an $n$-tableau of shape $D$. Let $\mathbf{F}_{T}:=\nabla
_{\operatorname*{Row}T}\nabla_{\operatorname*{Col}T}^{-}\in\mathcal{A}$. Prove
that $\mathcal{A}\mathbf{E}_{T}\cong\mathcal{A}\mathbf{E}_{T\mathbf{r}}%
\cong\mathcal{A}\mathbf{F}_{T}$ as $\mathbf{k}$-modules. (But not generally as
$\mathcal{A}$-modules, as we saw in Exercise \ref{exe.specht.FT.not-mod-p}
\textbf{(c)}.)
\end{exercise}

Next, we shall extend (\ref{eq.prop.specht.ET.r.3}) to the generalization
$\mathbf{E}_{P,Q}$ of the Young symmetrizers. First, however, we prove the
following simple properties of the permutations $w_{P,Q}$ from Definition
\ref{def.specht.wPQ}:

\begin{proposition}
\label{prop.specht.wPQ.wQP}Let $D$ be any diagram. Then: \medskip

\textbf{(a)} If $P$, $Q$ and $R$ are three $n$-tableaux of shape $D$, then
$w_{P,Q}w_{Q,R}=w_{P,R}$. \medskip

\textbf{(b)} If $P$ and $Q$ are two $n$-tableaux of shape $D$, then
$w_{Q,P}=w_{P,Q}^{-1}$. \medskip

\textbf{(c)} If $P$ and $Q$ are two $n$-tableaux of shape $D$, then
$w_{Q,P}=w_{Q\mathbf{r},P\mathbf{r}}$. \medskip

\textbf{(d)} If $T$ is an $n$-tableau of shape $D$, and if $u\in S_{n}$, then
$w_{uT,T}=u$.
\end{proposition}

\begin{proof}
\textbf{(a)} Let $P$, $Q$ and $R$ be three $n$-tableaux of shape $D$. Then,
\[
w_{P,Q}\underbrace{w_{Q,R}R}_{\substack{=Q\\\text{(by
(\ref{eq.def.specht.wPQ.eq}), applied to }Q\text{ and }R\\\text{instead of
}P\text{ and }Q\text{)}}}=w_{P,Q}Q=P\ \ \ \ \ \ \ \ \ \ \left(  \text{by
(\ref{eq.def.specht.wPQ.eq})}\right)  .
\]
But $w_{P,R}$ is defined as the unique permutation $w\in S_{n}$ that satisfies
$wR=P$. Hence, if a permutation $w\in S_{n}$ satisfies $wR=P$, then
$w=w_{P,R}$. Applying this to $w=w_{P,Q}w_{Q,R}$, we obtain $w_{P,Q}%
w_{Q,R}=w_{P,R}$ (since $w_{P,Q}w_{Q,R}R=P$). This proves Proposition
\ref{prop.specht.wPQ.wQP} \textbf{(a)}. \medskip

\textbf{(b)} Let $P$ and $Q$ be two $n$-tableaux of shape $D$. Then,
Proposition \ref{prop.specht.wPQ.wQP} \textbf{(a)} (applied to $R=P$) yields
$w_{P,Q}w_{Q,P}=w_{P,P}=\operatorname*{id}$ (by Proposition
\ref{prop.specht.EPQ.EP}). Hence, $w_{Q,P}=w_{P,Q}^{-1}$. This proves
Proposition \ref{prop.specht.wPQ.wQP} \textbf{(b)}. \medskip

\textbf{(c)} Let $P$ and $Q$ be two $n$-tableaux of shape $D$. Then,
(\ref{eq.def.specht.wPQ.eq}) (applied to $Q$ and $P$ instead of $P$ and $Q$)
yields $w_{Q,P}P=Q$. However, Proposition \ref{prop.tableaux.r} \textbf{(d)}
(applied to $P$ and $w_{Q,P}$ instead of $T$ and $w$) yields $\left(
w_{Q,P}P\right)  \mathbf{r}=w_{Q,P}\left(  P\mathbf{r}\right)  $. Hence,
$w_{Q,P}\left(  P\mathbf{r}\right)  =\underbrace{\left(  w_{Q,P}P\right)
}_{=Q}\mathbf{r}=Q\mathbf{r}$.

However, $w_{Q\mathbf{r},P\mathbf{r}}$ is defined as the unique permutation
$w\in S_{n}$ that satisfies $w\left(  P\mathbf{r}\right)  =Q\mathbf{r}$.
Hence, if a permutation $w\in S_{n}$ satisfies $w\left(  P\mathbf{r}\right)
=Q\mathbf{r}$, then $w=w_{Q\mathbf{r},P\mathbf{r}}$. Applying this to
$w=w_{Q,P}$, we obtain $w_{Q,P}=w_{Q\mathbf{r},P\mathbf{r}}$ (since
$w_{Q,P}\left(  P\mathbf{r}\right)  =Q\mathbf{r}$). Proposition
\ref{prop.specht.wPQ.wQP} \textbf{(c)} is now proved. \medskip

\textbf{(d)} Let $T$ be an $n$-tableau of shape $D$. Let $u\in S_{n}$. Then,
$w_{uT,T}$ is defined as the unique permutation $w\in S_{n}$ that satisfies
$wT=uT$. Hence, if a permutation $w\in S_{n}$ satisfies $wT=uT$, then
$w=w_{uT,T}$. Applying this to $w=u$, we obtain $u=w_{uT,T}$ (since $uT=uT$).
In other words, $w_{uT,T}=u$. This proves Proposition
\ref{prop.specht.wPQ.wQP} \textbf{(d)}.
\end{proof}

We are now ready to prove the generalized version of
(\ref{eq.prop.specht.ET.r.3}):

\begin{proposition}
\label{prop.specht.EPQ.r}Let $D$ be any diagram. Let $P$ and $Q$ be two
$n$-tableaux of shape $D$. Then,%
\[
\left(  -1\right)  ^{w_{Q,P}}\mathbf{E}_{Q\mathbf{r},P\mathbf{r}%
}=T_{\operatorname*{sign}}\left(  S\left(  \mathbf{E}_{P,Q}\right)  \right)
.
\]

\end{proposition}

\begin{proof}
The definition of $S$ yields $S\left(  w_{P,Q}\right)  =w_{P,Q}^{-1}$. But
Proposition \ref{prop.specht.wPQ.wQP} \textbf{(b)} yields $w_{Q,P}%
=w_{P,Q}^{-1}$. Comparing these two equalities, we find $S\left(
w_{P,Q}\right)  =w_{Q,P}$.

The definition of $\mathbf{E}_{P,Q}$ says that $\mathbf{E}_{P,Q}%
=\nabla_{\operatorname*{Col}P}^{-}w_{P,Q}\nabla_{\operatorname*{Row}Q}$.
Applying the map $S$ to this equality, we obtain%
\begin{align*}
S\left(  \mathbf{E}_{P,Q}\right)   &  =S\left(  \nabla_{\operatorname*{Col}%
P}^{-}w_{P,Q}\nabla_{\operatorname*{Row}Q}\right) \\
&  =\underbrace{S\left(  \nabla_{\operatorname*{Row}Q}\right)  }%
_{\substack{=\nabla_{\operatorname*{Row}Q}\\\text{(by Proposition
\ref{prop.symmetrizers.antipode},}\\\text{applied to }T=Q\text{)}%
}}\underbrace{S\left(  w_{P,Q}\right)  }_{=w_{Q,P}}\underbrace{S\left(
\nabla_{\operatorname*{Col}P}^{-}\right)  }_{\substack{=\nabla
_{\operatorname*{Col}P}^{-}\\\text{(by Proposition
\ref{prop.symmetrizers.antipode},}\\\text{applied to }T=P\text{)}}}\\
&  \ \ \ \ \ \ \ \ \ \ \ \ \ \ \ \ \ \ \ \ \left(
\begin{array}
[c]{c}%
\text{since Theorem \ref{thm.S.auto} \textbf{(a)} shows}\\
\text{that }S\text{ is a }\mathbf{k}\text{-algebra anti-morphism}%
\end{array}
\right) \\
&  =\nabla_{\operatorname*{Row}Q}w_{Q,P}\nabla_{\operatorname*{Col}P}^{-}.
\end{align*}

Applying the map $T_{\operatorname*{sign}}$ to this equality, we obtain%
\begin{align*}
&  T_{\operatorname*{sign}}\left(  S\left(  \mathbf{E}_{P,Q}\right)  \right)
\\
&  =T_{\operatorname*{sign}}\left(  \nabla_{\operatorname*{Row}Q}w_{Q,P}%
\nabla_{\operatorname*{Col}P}^{-}\right) \\
&  =\underbrace{T_{\operatorname*{sign}}\left(  \nabla_{\operatorname*{Row}%
Q}\right)  }_{\substack{=\nabla_{\operatorname*{Col}\left(  Q\mathbf{r}%
\right)  }^{-}\\\text{(since (\ref{eq.prop.specht.ET.r.2}) (applied to
}Q\text{ instead of }P\text{)}\\\text{yields }\nabla_{\operatorname*{Col}%
\left(  Q\mathbf{r}\right)  }^{-}=T_{\operatorname*{sign}}\left(
\nabla_{\operatorname*{Row}Q}\right)  \text{)}}%
}\ \ \underbrace{T_{\operatorname*{sign}}\left(  w_{Q,P}\right)
}_{\substack{=\left(  -1\right)  ^{w_{Q,P}}w_{Q,P}\\\text{(by the definition
of }T_{\operatorname*{sign}}\text{)}}}\underbrace{T_{\operatorname*{sign}%
}\left(  \nabla_{\operatorname*{Col}P}^{-}\right)  }_{\substack{=\nabla
_{\operatorname*{Row}\left(  P\mathbf{r}\right)  }\\\text{(by
(\ref{eq.prop.specht.ET.r.1}))}}}\\
&  \ \ \ \ \ \ \ \ \ \ \ \ \ \ \ \ \ \ \ \ \left(
\begin{array}
[c]{c}%
\text{since Theorem \ref{thm.Tsign.auto} \textbf{(a)} shows}\\
\text{that }T_{\operatorname*{sign}}\text{ is a }\mathbf{k}\text{-algebra
morphism}%
\end{array}
\right) \\
&  =\nabla_{\operatorname*{Col}\left(  Q\mathbf{r}\right)  }^{-}\left(
-1\right)  ^{w_{Q,P}}\underbrace{w_{Q,P}}_{\substack{=w_{Q\mathbf{r}%
,P\mathbf{r}}\\\text{(by Proposition \ref{prop.specht.wPQ.wQP} \textbf{(c)})}%
}}\nabla_{\operatorname*{Row}\left(  P\mathbf{r}\right)  }\\
&  =\nabla_{\operatorname*{Col}\left(  Q\mathbf{r}\right)  }^{-}\left(
-1\right)  ^{w_{Q,P}}w_{Q\mathbf{r},P\mathbf{r}}\nabla_{\operatorname*{Row}%
\left(  P\mathbf{r}\right)  }\\
&  =\left(  -1\right)  ^{w_{Q,P}}\cdot\underbrace{\nabla_{\operatorname*{Col}%
\left(  Q\mathbf{r}\right)  }^{-}w_{Q\mathbf{r},P\mathbf{r}}\nabla
_{\operatorname*{Row}\left(  P\mathbf{r}\right)  }}_{\substack{=\mathbf{E}%
_{Q\mathbf{r},P\mathbf{r}}\\\text{(since }\mathbf{E}_{Q\mathbf{r},P\mathbf{r}%
}\text{ is defined to}\\\text{be }\nabla_{\operatorname*{Col}\left(
Q\mathbf{r}\right)  }^{-}w_{Q\mathbf{r},P\mathbf{r}}\nabla
_{\operatorname*{Row}\left(  P\mathbf{r}\right)  }\text{)}}}=\left(
-1\right)  ^{w_{Q,P}}\mathbf{E}_{Q\mathbf{r},P\mathbf{r}}.
\end{align*}
This proves Proposition \ref{prop.specht.EPQ.r}.
\end{proof}

\subsubsection{Straightening left and right}

Combining the above properties of transposed tableaux with Proposition
\ref{prop.specht.EPQ.basis-AET}, we obtain the following:

\begin{lemma}
\label{lem.specht.EPQ.straighten}Let $D$ be a skew Young diagram. Let $Q$ be
an $n$-tableau of shape $D$. Let $\operatorname*{SYT}\left(  D\right)  $ be
the set of all standard $n$-tableaux of shape $D$. Then,
\begin{equation}
\mathcal{A}\mathbf{E}_{Q}=\operatorname*{span}\nolimits_{\mathbf{k}}\left\{
\mathbf{E}_{P,Q}\ \mid\ P\in\operatorname*{SYT}\left(  D\right)  \right\}
\label{eq.lem.specht.EPQ.straighten.AET}%
\end{equation}
and%
\begin{equation}
\mathbf{E}_{Q}\mathcal{A}=\operatorname*{span}\nolimits_{\mathbf{k}}\left\{
\mathbf{E}_{Q,P}\ \mid\ P\in\operatorname*{SYT}\left(  D\right)  \right\}  .
\label{eq.lem.specht.EPQ.straighten.ETA}%
\end{equation}

\end{lemma}

\begin{proof}
We know that $Q$ is an $n$-tableau of shape $D$, thus a bijection from $D$ to
$\left[  n\right]  $. Hence, $\left\vert D\right\vert =\left\vert \left[
n\right]  \right\vert =n$. But $\mathbf{r}$ is a bijection; thus, $\left\vert
\mathbf{r}\left(  D\right)  \right\vert =\left\vert D\right\vert =n$ as well.

Proposition \ref{prop.specht.EPQ.basis-AET} (applied to $T=Q$) yields that the
family $\left(  \mathbf{E}_{P,Q}\right)  _{P\in\operatorname*{SYT}\left(
D\right)  }$ is a basis of the $\mathbf{k}$-module $\mathcal{A}\mathbf{E}_{Q}%
$. Thus, this family spans $\mathcal{A}\mathbf{E}_{Q}$. In other words,
\newline$\mathcal{A}\mathbf{E}_{Q}=\operatorname*{span}\nolimits_{\mathbf{k}%
}\left\{  \mathbf{E}_{P,Q}\ \mid\ P\in\operatorname*{SYT}\left(  D\right)
\right\}  $. This proves (\ref{eq.lem.specht.EPQ.straighten.AET}).

It remains to prove (\ref{eq.lem.specht.EPQ.straighten.ETA}). To do so, we
will essentially apply (\ref{eq.lem.specht.EPQ.straighten.AET}) to
$\mathbf{r}\left(  D\right)  $ and $Q\mathbf{r}$ instead of $D$ and $Q$, and
then transform the result using the antipode $S$ (which is an anti-morphism
and thus turns $\mathcal{A}\mathbf{E}_{Q\mathbf{r}}$ into $S\left(
\mathbf{E}_{Q\mathbf{r}}\right)  \mathcal{A}$) and the sign-twist
$T_{\operatorname*{sign}}$. Here are the details:

For any diagram $E$, let $\operatorname*{SYT}\left(  E\right)  $ be the set of
all standard $n$-tableaux of shape $E$. (This generalizes the notation
$\operatorname*{SYT}\left(  D\right)  $ defined in the lemma.) Proposition
\ref{prop.tableaux.r} \textbf{(c)} then says that the map%
\begin{align*}
\operatorname*{SYT}\left(  D\right)   &  \rightarrow\operatorname*{SYT}\left(
\mathbf{r}\left(  D\right)  \right)  ,\\
T  &  \mapsto T\mathbf{r}%
\end{align*}
is a bijection.

But $D$ is a skew Young diagram. Hence, $\mathbf{r}\left(  D\right)  $ is a
skew Young diagram as well\footnote{\textit{Proof.} We assumed that $D$ is a
skew Young diagram. In other words, $D=Y\left(  \lambda/\mu\right)  $ for some
skew partition $\lambda/\mu$. Consider this $\lambda/\mu$. Now, Remark
\ref{rmk.partitions.conj.Ylm} shows that $\lambda^{t}/\mu^{t}$ is again a skew
partition, and that we have $Y\left(  \lambda^{t}/\mu^{t}\right)
=\mathbf{r}\left(  Y\left(  \lambda/\mu\right)  \right)  $. But $D=Y\left(
\lambda/\mu\right)  $ and therefore $\mathbf{r}\left(  D\right)
=\mathbf{r}\left(  Y\left(  \lambda/\mu\right)  \right)  =Y\left(  \lambda
^{t}/\mu^{t}\right)  $. This shows that $\mathbf{r}\left(  D\right)  $ is a
skew Young diagram.}. Therefore, (\ref{eq.lem.specht.EPQ.straighten.AET})
(applied to $\mathbf{r}\left(  D\right)  $ and $Q\mathbf{r}$ instead of $D$
and $Q$) yields%
\[
\mathcal{A}\mathbf{E}_{Q\mathbf{r}}=\operatorname*{span}\nolimits_{\mathbf{k}%
}\left\{  \mathbf{E}_{P,Q\mathbf{r}}\ \mid\ P\in\operatorname*{SYT}\left(
\mathbf{r}\left(  D\right)  \right)  \right\}
\]
(since $Q\mathbf{r}$ is an $n$-tableau of shape $\mathbf{r}\left(  D\right)
$). Applying the map $S$ to both sides of this equality, we obtain%
\begin{align*}
S\left(  \mathcal{A}\mathbf{E}_{Q\mathbf{r}}\right)   &  =S\left(
\operatorname*{span}\nolimits_{\mathbf{k}}\left\{  \mathbf{E}_{P,Q\mathbf{r}%
}\ \mid\ P\in\operatorname*{SYT}\left(  \mathbf{r}\left(  D\right)  \right)
\right\}  \right) \\
&  =\operatorname*{span}\nolimits_{\mathbf{k}}\left\{  S\left(  \mathbf{E}%
_{P,Q\mathbf{r}}\right)  \ \mid\ P\in\operatorname*{SYT}\left(  \mathbf{r}%
\left(  D\right)  \right)  \right\}  \ \ \ \ \ \ \ \ \ \ \left(
\begin{array}
[c]{c}%
\text{since the map }S\\
\text{is }\mathbf{k}\text{-linear}%
\end{array}
\right) \\
&  =\operatorname*{span}\nolimits_{\mathbf{k}}\left\{  S\left(  \mathbf{E}%
_{T\mathbf{r},Q\mathbf{r}}\right)  \ \mid\ T\in\operatorname*{SYT}\left(
D\right)  \right\} \\
&  \ \ \ \ \ \ \ \ \ \ \ \ \ \ \ \ \ \ \ \ \left(
\begin{array}
[c]{c}%
\text{here, we have substituted }T\mathbf{r}\text{ for }P\text{, since}\\
\text{the map }\operatorname*{SYT}\left(  D\right)  \rightarrow
\operatorname*{SYT}\left(  \mathbf{r}\left(  D\right)  \right)  ,\ T\mapsto
T\mathbf{r}\\
\text{is a bijection}%
\end{array}
\right) \\
&  =\operatorname*{span}\nolimits_{\mathbf{k}}\left\{  S\left(  \mathbf{E}%
_{P\mathbf{r},Q\mathbf{r}}\right)  \ \mid\ P\in\operatorname*{SYT}\left(
D\right)  \right\}
\end{align*}
(here, we have renamed the index $T$ as $P$). Applying the map
$T_{\operatorname*{sign}}$ to this equality, we furthermore obtain%
\begin{align}
&  T_{\operatorname*{sign}}\left(  S\left(  \mathcal{A}\mathbf{E}%
_{Q\mathbf{r}}\right)  \right) \nonumber\\
&  =T_{\operatorname*{sign}}\left(  \operatorname*{span}\nolimits_{\mathbf{k}%
}\left\{  S\left(  \mathbf{E}_{P\mathbf{r},Q\mathbf{r}}\right)  \ \mid
\ P\in\operatorname*{SYT}\left(  D\right)  \right\}  \right) \nonumber\\
&  =\operatorname*{span}\nolimits_{\mathbf{k}}\left\{  T_{\operatorname*{sign}%
}\left(  S\left(  \mathbf{E}_{P\mathbf{r},Q\mathbf{r}}\right)  \right)
\ \mid\ P\in\operatorname*{SYT}\left(  D\right)  \right\}
\label{pf.lem.specht.EPQ.straighten.4}%
\end{align}
(since the map $T_{\operatorname*{sign}}$ is $\mathbf{k}$-linear).

Now, let $P\in\operatorname*{SYT}\left(  D\right)  $. Then, $P\mathbf{r}$ and
$Q\mathbf{r}$ are two $n$-tableaux of shape $\mathbf{r}\left(  D\right)  $
(since $P$ and $Q$ are two $n$-tableaux of shape $D$). Thus, Proposition
\ref{prop.specht.EPQ.r} (applied to $\mathbf{r}\left(  D\right)  $,
$P\mathbf{r}$ and $Q\mathbf{r}$ instead of $D$, $P$ and $Q$) yields%
\begin{equation}
\left(  -1\right)  ^{w_{Q\mathbf{r},P\mathbf{r}}}\mathbf{E}_{\left(
Q\mathbf{r}\right)  \mathbf{r},\left(  P\mathbf{r}\right)  \mathbf{r}%
}=T_{\operatorname*{sign}}\left(  S\left(  \mathbf{E}_{P\mathbf{r}%
,Q\mathbf{r}}\right)  \right)  . \label{pf.lem.specht.EPQ.straighten.5}%
\end{equation}
But Proposition \ref{prop.tableaux.r} \textbf{(b)} yields $\left(
P\mathbf{r}\right)  \mathbf{r}=P$ and $\left(  Q\mathbf{r}\right)
\mathbf{r}=Q$. In view of this, we can rewrite
(\ref{pf.lem.specht.EPQ.straighten.5}) as%
\begin{equation}
\left(  -1\right)  ^{w_{Q\mathbf{r},P\mathbf{r}}}\mathbf{E}_{Q,P}%
=T_{\operatorname*{sign}}\left(  S\left(  \mathbf{E}_{P\mathbf{r},Q\mathbf{r}%
}\right)  \right)  . \label{pf.lem.specht.EPQ.straighten.6}%
\end{equation}

Forget that we fixed $P$. We thus have proved
(\ref{pf.lem.specht.EPQ.straighten.6}) for each $P\in\operatorname*{SYT}%
\left(  D\right)  $. Thus, (\ref{pf.lem.specht.EPQ.straighten.4}) becomes%
\begin{align}
T_{\operatorname*{sign}}\left(  S\left(  \mathcal{A}\mathbf{E}_{Q\mathbf{r}%
}\right)  \right)   &  =\operatorname*{span}\nolimits_{\mathbf{k}}\left\{
\underbrace{T_{\operatorname*{sign}}\left(  S\left(  \mathbf{E}_{P\mathbf{r}%
,Q\mathbf{r}}\right)  \right)  }_{\substack{=\left(  -1\right)
^{w_{Q\mathbf{r},P\mathbf{r}}}\mathbf{E}_{Q,P}\\\text{(by
(\ref{pf.lem.specht.EPQ.straighten.6}))}}}\ \mid\ P\in\operatorname*{SYT}%
\left(  D\right)  \right\} \nonumber\\
&  =\operatorname*{span}\nolimits_{\mathbf{k}}\left\{  \left(  -1\right)
^{w_{Q\mathbf{r},P\mathbf{r}}}\mathbf{E}_{Q,P}\ \mid\ P\in\operatorname*{SYT}%
\left(  D\right)  \right\} \nonumber\\
&  =\operatorname*{span}\nolimits_{\mathbf{k}}\left\{  \mathbf{E}_{Q,P}%
\ \mid\ P\in\operatorname*{SYT}\left(  D\right)  \right\}
\label{pf.lem.specht.EPQ.straighten.7}%
\end{align}
(here, we have removed the $\left(  -1\right)  ^{w_{Q\mathbf{r},P\mathbf{r}}}%
$\ factors inside the span, since the span of a set of vectors does not change
if we scale each of these vectors by a power of $-1$).

Proposition \ref{prop.tableaux.r} \textbf{(b)} yields $\left(  Q\mathbf{r}%
\right)  \mathbf{r}=Q$. But $Q\mathbf{r}$ is an $n$-tableau of shape
$\mathbf{r}\left(  D\right)  $ (since $Q$ is an $n$-tableau of shape $D$).
Thus, (\ref{eq.prop.specht.ET.r.3}) (applied to $\mathbf{r}\left(  D\right)  $
and $Q\mathbf{r}$ instead of $D$ and $P$) yields $\mathbf{E}_{\left(
Q\mathbf{r}\right)  \mathbf{r}}=T_{\operatorname*{sign}}\left(  S\left(
\mathbf{E}_{Q\mathbf{r}}\right)  \right)  $. In view of $\left(
Q\mathbf{r}\right)  \mathbf{r}=Q$, we can rewrite this as%
\begin{equation}
\mathbf{E}_{Q}=T_{\operatorname*{sign}}\left(  S\left(  \mathbf{E}%
_{Q\mathbf{r}}\right)  \right)  . \label{pf.lem.specht.EPQ.straighten.EQ=}%
\end{equation}

However, $S:\mathcal{A}\rightarrow\mathcal{A}$ is an involution (by Theorem
\ref{thm.S.auto} \textbf{(b)}), and thus a bijection. Hence, $S\left(
\mathcal{A}\right)  =\mathcal{A}$. Furthermore, $S$ is an algebra
anti-morphism (by Theorem \ref{thm.S.auto} \textbf{(a)}). Therefore,%
\begin{equation}
S\left(  \mathcal{A}\mathbf{E}_{Q\mathbf{r}}\right)  =S\left(  \mathbf{E}%
_{Q\mathbf{r}}\right)  \underbrace{S\left(  \mathcal{A}\right)  }%
_{=\mathcal{A}}=S\left(  \mathbf{E}_{Q\mathbf{r}}\right)  \mathcal{A}.
\label{pf.lem.specht.EPQ.straighten.8}%
\end{equation}

Furthermore, $T_{\operatorname*{sign}}:\mathcal{A}\rightarrow\mathcal{A}$ is
an involution (by Theorem \ref{thm.Tsign.auto} \textbf{(b)}), and thus a
bijection. Hence, $T_{\operatorname*{sign}}\left(  \mathcal{A}\right)
=\mathcal{A}$. Furthermore, $T_{\operatorname*{sign}}$ is an algebra morphism
(by Theorem \ref{thm.Tsign.auto} \textbf{(a)}). Thus,%
\begin{align}
T_{\operatorname*{sign}}\left(  S\left(  \mathbf{E}_{Q\mathbf{r}}\right)
\mathcal{A}\right)   &  =\underbrace{T_{\operatorname*{sign}}\left(  S\left(
\mathbf{E}_{Q\mathbf{r}}\right)  \right)  }_{\substack{=\mathbf{E}%
_{Q}\\\text{(by (\ref{pf.lem.specht.EPQ.straighten.EQ=}))}}%
}\underbrace{T_{\operatorname*{sign}}\left(  \mathcal{A}\right)
}_{=\mathcal{A}}\nonumber\\
&  =\mathbf{E}_{Q}\mathcal{A}. \label{pf.lem.specht.EPQ.straighten.9}%
\end{align}
Applying the map $T_{\operatorname*{sign}}$ to both sides of
(\ref{pf.lem.specht.EPQ.straighten.8}), we obtain%
\[
T_{\operatorname*{sign}}\left(  S\left(  \mathcal{A}\mathbf{E}_{Q\mathbf{r}%
}\right)  \right)  =T_{\operatorname*{sign}}\left(  S\left(  \mathbf{E}%
_{Q\mathbf{r}}\right)  \mathcal{A}\right)  =\mathbf{E}_{Q}\mathcal{A}%
\ \ \ \ \ \ \ \ \ \ \left(  \text{by (\ref{pf.lem.specht.EPQ.straighten.9}%
)}\right)  .
\]
Comparing this with (\ref{pf.lem.specht.EPQ.straighten.7}), we obtain%
\[
\mathbf{E}_{Q}\mathcal{A}=\operatorname*{span}\nolimits_{\mathbf{k}}\left\{
\mathbf{E}_{Q,P}\ \mid\ P\in\operatorname*{SYT}\left(  D\right)  \right\}  .
\]
Thus, (\ref{eq.lem.specht.EPQ.straighten.ETA}) is proved, and the proof of
Lemma \ref{lem.specht.EPQ.straighten} is complete.
\end{proof}

\begin{lemma}
\label{lem.specht.EPQ.in-span}Let $D$ be a skew Young diagram. Let
$\operatorname*{SYT}\left(  D\right)  $ be the set of all standard
$n$-tableaux of shape $D$. Let $P$ and $Q$ be two $n$-tableaux of shape $D$.
Then,%
\[
\mathbf{E}_{P,Q}\in\operatorname*{span}\nolimits_{\mathbf{k}}\left\{
\mathbf{E}_{U,V}\ \mid\ U,V\in\operatorname*{SYT}\left(  D\right)  \right\}
.
\]

\end{lemma}

\begin{proof}
Let
\begin{equation}
\Omega:=\operatorname*{span}\nolimits_{\mathbf{k}}\left\{  \mathbf{E}%
_{U,V}\ \mid\ U,V\in\operatorname*{SYT}\left(  D\right)  \right\}  .
\label{pf.lem.specht.EPQ.in-span.Omega=}%
\end{equation}
Thus, $\Omega$ is a $\mathbf{k}$-submodule of $\mathcal{A}$ and satisfies%
\begin{equation}
\mathbf{E}_{U,V}\in\Omega\ \ \ \ \ \ \ \ \ \ \text{for each }U,V\in
\operatorname*{SYT}\left(  D\right)  .
\label{pf.lem.specht.EPQ.in-span.in-Omega}%
\end{equation}

Forget that we fixed $P$ and $Q$. We must prove that \newline$\mathbf{E}%
_{P,Q}\in\operatorname*{span}\nolimits_{\mathbf{k}}\left\{  \mathbf{E}%
_{U,V}\ \mid\ U,V\in\operatorname*{SYT}\left(  D\right)  \right\}  $ for any
two $n$-tableaux $P,Q$ of shape $D$. Using
(\ref{pf.lem.specht.EPQ.in-span.in-Omega}), we can restate this as follows: We
must prove that $\mathbf{E}_{P,Q}\in\Omega$ for any two $n$-tableaux $P,Q$ of
shape $D$. Let us rename the indices $P$ and $Q$ as $X$ and $Y$ here. Thus, we
must prove that $\mathbf{E}_{X,Y}\in\Omega$ for any two $n$-tableaux $X,Y$ of
shape $D$.

Let us now prove this. Let $X$ and $Y$ be two $n$-tableaux of shape $D$. Then,
(\ref{eq.prop.specht.EPQ.EPEQ.wE}) (applied to $P=X$ and $Q=Y$) yields
\begin{align}
\mathbf{E}_{X,Y}  &  =\underbrace{w_{X,Y}}_{\in\mathcal{A}}\mathbf{E}%
_{Y}\nonumber\\
&  \in\mathcal{A}\mathbf{E}_{Y}\nonumber\\
&  =\operatorname*{span}\nolimits_{\mathbf{k}}\left\{  \mathbf{E}_{P,Y}%
\ \mid\ P\in\operatorname*{SYT}\left(  D\right)  \right\}
\ \ \ \ \ \ \ \ \ \ \left(  \text{by (\ref{eq.lem.specht.EPQ.straighten.AET}),
applied to }Q=Y\right) \nonumber\\
&  =\operatorname*{span}\nolimits_{\mathbf{k}}\left\{  \mathbf{E}_{U,Y}%
\ \mid\ U\in\operatorname*{SYT}\left(  D\right)  \right\}
\label{pf.lem.specht.EPQ.in-span.3}%
\end{align}
(here, we have renamed the index $P$ as $U$). Furthermore, for each
$U\in\operatorname*{SYT}\left(  D\right)  $, we have%
\begin{align*}
\mathbf{E}_{U,Y}  &  =\mathbf{E}_{U}\underbrace{w_{U,Y}}_{\in\mathcal{A}%
}\ \ \ \ \ \ \ \ \ \ \left(  \text{by (\ref{eq.prop.specht.EPQ.EPEQ.Ew}),
applied to }P=U\text{ and }Q=Y\right) \\
&  \in\mathbf{E}_{U}\mathcal{A}\\
&  =\operatorname*{span}\nolimits_{\mathbf{k}}\underbrace{\left\{
\mathbf{E}_{U,P}\ \mid\ P\in\operatorname*{SYT}\left(  D\right)  \right\}
}_{\substack{=\left\{  \mathbf{E}_{U,V}\ \mid\ V\in\operatorname*{SYT}\left(
D\right)  \right\}  \\\subseteq\Omega\\\text{(since
(\ref{pf.lem.specht.EPQ.in-span.in-Omega}) shows}\\\text{that }\mathbf{E}%
_{U,V}\in\Omega\text{ for each }V\in\operatorname*{SYT}\left(  D\right)
\text{)}}}\ \ \ \ \ \ \ \ \ \ \left(  \text{by
(\ref{eq.lem.specht.EPQ.straighten.ETA}), applied to }Q=U\right) \\
&  \subseteq\operatorname*{span}\nolimits_{\mathbf{k}}\Omega\subseteq
\Omega\ \ \ \ \ \ \ \ \ \ \left(  \text{since }\Omega\text{ is a }%
\mathbf{k}\text{-module}\right)  .
\end{align*}
In other words,%
\[
\left\{  \mathbf{E}_{U,Y}\ \mid\ U\in\operatorname*{SYT}\left(  D\right)
\right\}  \subseteq\Omega.
\]
Thus, (\ref{pf.lem.specht.EPQ.in-span.3}) becomes%
\[
\mathbf{E}_{X,Y}\in\operatorname*{span}\nolimits_{\mathbf{k}}%
\underbrace{\left\{  \mathbf{E}_{U,Y}\ \mid\ U\in\operatorname*{SYT}\left(
D\right)  \right\}  }_{\subseteq\Omega}\subseteq\operatorname*{span}%
\nolimits_{\mathbf{k}}\Omega\subseteq\Omega
\]
(since $\Omega$ is a $\mathbf{k}$-module).

Forget that we fixed $X$ and $Y$. We thus have shown that $\mathbf{E}_{X,Y}%
\in\Omega$ for any two $n$-tableaux $X,Y$ of shape $D$. As we said above, this
completes the proof of Lemma \ref{lem.specht.EPQ.in-span}.
\end{proof}

\begin{proposition}
\label{prop.specht.AElam.span}Let $\lambda$ be a partition of $n$. Let
$D=Y\left(  \lambda\right)  $. Let $\operatorname*{SYT}\left(  D\right)  $ be
the set of all standard $n$-tableaux of shape $D$. Then: \medskip

\textbf{(a)} If $T$ is any $n$-tableau of shape $D$, then%
\[
\mathcal{A}\mathbf{E}_{T}\subseteq\operatorname*{span}\nolimits_{\mathbf{k}%
}\left\{  \mathbf{E}_{U,V}\ \mid\ U,V\in\operatorname*{SYT}\left(  D\right)
\right\}  .
\]

\textbf{(b)} We have%
\[
\mathcal{A}\mathbf{E}_{\lambda}\subseteq\operatorname*{span}%
\nolimits_{\mathbf{k}}\left\{  \mathbf{E}_{U,V}\ \mid\ U,V\in
\operatorname*{SYT}\left(  D\right)  \right\}  .
\]

\end{proposition}

\begin{proof}
Note that $D=Y\left(  \lambda\right)  =Y\left(  \lambda/\varnothing\right)  $
is a skew Young diagram.

Let $\Omega:=\operatorname*{span}\nolimits_{\mathbf{k}}\left\{  \mathbf{E}%
_{U,V}\ \mid\ U,V\in\operatorname*{SYT}\left(  D\right)  \right\}  $. Then,
$\Omega$ is a $\mathbf{k}$-submodule of $\mathcal{A}$. Hence,
$\operatorname*{span}\nolimits_{\mathbf{k}}\Omega\subseteq\Omega$. \medskip

\textbf{(a)} Let $T$ be any $n$-tableau of shape $D$. Let $u\in S_{n}$. We
shall first show that $u\mathbf{E}_{T}\in\Omega$.

Indeed, Proposition \ref{prop.specht.wPQ.wQP} \textbf{(d)} yields $w_{uT,T}%
=u$. But (\ref{eq.prop.specht.EPQ.EPEQ.wE}) (applied to $P=uT$ and $Q=T$)
yields $\mathbf{E}_{uT,T}=w_{uT,T}\mathbf{E}_{T}=u\mathbf{E}_{T}$ (since
$w_{uT,T}=u$). Hence,%
\begin{align*}
u\mathbf{E}_{T}  &  =\mathbf{E}_{uT,T}\in\operatorname*{span}%
\nolimits_{\mathbf{k}}\left\{  \mathbf{E}_{U,V}\ \mid\ U,V\in
\operatorname*{SYT}\left(  D\right)  \right\} \\
&  \ \ \ \ \ \ \ \ \ \ \ \ \ \ \ \ \ \ \ \ \left(  \text{by Lemma
\ref{lem.specht.EPQ.in-span}, applied to }P=uT\text{ and }Q=T\right) \\
&  =\Omega.
\end{align*}

Forget that we fixed $u$. We thus have shown that $u\mathbf{E}_{T}\in\Omega$
for each $u\in S_{n}$. In other words, $\left\{  u\mathbf{E}_{T}\ \mid\ u\in
S_{n}\right\}  \subseteq\Omega$. However, $\mathcal{A}=\mathbf{k}\left[
S_{n}\right]  =\operatorname*{span}\nolimits_{\mathbf{k}}\left\{
u\ \mid\ u\in S_{n}\right\}  $ (since the permutations $u\in S_{n}$ form a
basis of $\mathbf{k}\left[  S_{n}\right]  $), and thus%
\begin{align*}
\mathcal{A}\mathbf{E}_{T}  &  =\operatorname*{span}\nolimits_{\mathbf{k}%
}\left\{  u\ \mid\ u\in S_{n}\right\}  \cdot\mathbf{E}_{T}\\
&  =\operatorname*{span}\nolimits_{\mathbf{k}}\underbrace{\left\{
u\mathbf{E}_{T}\ \mid\ u\in S_{n}\right\}  }_{\subseteq\Omega}\\
&  \subseteq\operatorname*{span}\nolimits_{\mathbf{k}}\Omega\subseteq
\Omega=\operatorname*{span}\nolimits_{\mathbf{k}}\left\{  \mathbf{E}%
_{U,V}\ \mid\ U,V\in\operatorname*{SYT}\left(  D\right)  \right\}  .
\end{align*}
This proves Proposition \ref{prop.specht.AElam.span} \textbf{(a)}. \medskip

\textbf{(b)} We have%
\begin{align*}
\mathcal{A}\mathbf{E}_{\lambda}  &  =\mathcal{A}\left(  \sum_{T\text{ is an
}n\text{-tableau of shape }\lambda}\mathbf{E}_{T}\right)
\ \ \ \ \ \ \ \ \ \ \left(  \text{by (\ref{eq.def.spechtmod.Elam.Elam.def}%
)}\right) \\
&  \subseteq\sum_{T\text{ is an }n\text{-tableau of shape }\lambda
}\underbrace{\mathcal{A}\mathbf{E}_{T}}_{\substack{\subseteq
\operatorname*{span}\nolimits_{\mathbf{k}}\left\{  \mathbf{E}_{U,V}%
\ \mid\ U,V\in\operatorname*{SYT}\left(  D\right)  \right\}  \\\text{(by
Proposition \ref{prop.specht.AElam.span} \textbf{(a)})}}%
}\ \ \ \ \ \ \ \ \ \ \left(  \text{by distributivity}\right) \\
&  \subseteq\sum_{T\text{ is an }n\text{-tableau of shape }\lambda
}\operatorname*{span}\nolimits_{\mathbf{k}}\left\{  \mathbf{E}_{U,V}%
\ \mid\ U,V\in\operatorname*{SYT}\left(  D\right)  \right\} \\
&  \subseteq\operatorname*{span}\nolimits_{\mathbf{k}}\left\{  \mathbf{E}%
_{U,V}\ \mid\ U,V\in\operatorname*{SYT}\left(  D\right)  \right\}
\end{align*}
(since $\operatorname*{span}\nolimits_{\mathbf{k}}\left\{  \mathbf{E}%
_{U,V}\ \mid\ U,V\in\operatorname*{SYT}\left(  D\right)  \right\}  $ is a
$\mathbf{k}$-module). This proves Proposition \ref{prop.specht.AElam.span}
\textbf{(b)}.
\end{proof}

\begin{exercise}
\fbox{2} Let $D$ be a skew Young diagram. Let $T$ be an $n$-tableau of shape
$D$. Let $\operatorname*{SYT}\left(  D\right)  $ be the set of all standard
$n$-tableaux of shape $D$. Prove that the family $\left(  \mathbf{E}%
_{T,P}\right)  _{P\in\operatorname*{SYT}\left(  D\right)  }$ is a basis of the
$\mathbf{k}$-module $\mathbf{E}_{T}\mathcal{A}$.
\end{exercise}

\subsubsection{Decomposing $\mathcal{A}$ and $\mathcal{A}\mathbf{E}_{\lambda}$
as left $\mathcal{A}$-modules}

We now take aim at the following crucial result:

\begin{theorem}
\label{thm.specht.AElam.decomp}Let $\lambda$ be a partition of $n$. Let
$D=Y\left(  \lambda\right)  $. Let $\operatorname*{SYT}\left(  D\right)  $ be
the set of all standard $n$-tableaux of shape $D$. Assume that $h^{\lambda}$
is invertible in $\mathbf{k}$. Then,%
\[
\mathcal{A}\mathbf{E}_{\lambda}=\bigoplus_{T\in\operatorname*{SYT}\left(
D\right)  }\mathcal{A}\mathbf{E}_{T}\ \ \ \ \ \ \ \ \ \ \left(  \text{an
internal direct sum}\right)  .
\]

\end{theorem}

The proof will rely on the following three lemmas:

\begin{lemma}
\label{lem.specht.AElam.triang}Let $\lambda$ be a partition of $n$. Let $S$
and $T$ be two standard $n$-tableaux of shape $Y\left(  \lambda\right)  $.
Assume that $\overline{S}>\overline{T}$ with respect to the Young last letter
order (this is the total order constructed in Proposition
\ref{prop.spechtmod.row.order}). Then, $\mathbf{E}_{S}\mathbf{E}_{T}=0$.
\end{lemma}

\begin{proof}
We have $\overline{S}>\overline{T}$, so that $\overline{T}<\overline{S}$.
Moreover, the $n$-tableau $T$ is standard, hence column-standard. Thus, Lemma
\ref{lem.specht.sandw0-llo} \textbf{(a)} shows that there exist no
permutations $r\in\mathcal{R}\left(  S\right)  $ and $c\in\mathcal{C}\left(
T\right)  $ such that $rS=cT$. Hence, Proposition \ref{prop.specht.ET.ESET}
\textbf{(a)} yields $\mathbf{E}_{S}\mathbf{E}_{T}=0$, so that Lemma
\ref{lem.specht.AElam.triang} is proven.
\end{proof}

\begin{lemma}
\label{lem.specht.AElam.direct}Let $\lambda$ be a partition of $n$. Let
$D=Y\left(  \lambda\right)  $. Let $\operatorname*{SYT}\left(  D\right)  $ be
the set of all standard $n$-tableaux of shape $D$. Assume that $h^{\lambda}$
is invertible in $\mathbf{k}$. Then, the sum $\sum_{T\in\operatorname*{SYT}%
\left(  D\right)  }\mathcal{A}\mathbf{E}_{T}$ is direct (i.e., the submodules
$\mathcal{A}\mathbf{E}_{T}$ for all $T\in\operatorname*{SYT}\left(  D\right)
$ are linearly disjoint).
\end{lemma}

\begin{proof}
Let us first recall Proposition \ref{prop.tabloid.row-st}, which says that
there is a bijection%
\begin{align*}
&  \text{from }\left\{  \text{row-standard }n\text{-tableaux of shape
}D\right\} \\
&  \text{to }\left\{  n\text{-tabloids of shape }D\right\}
\end{align*}
that sends each row-standard $n$-tableau $T$ to its tabloid $\overline{T}$.
Thus, in particular, this bijection is injective. In other words, if two
elements $T$ and $S$ of $\left\{  \text{row-standard }n\text{-tableaux of
shape }D\right\}  $ are distinct, then their images under this bijection are
also distinct, i.e., we have $\overline{T}\neq\overline{S}$. In other words,
the following holds:

\begin{statement}
\textit{Observation 1:} If $T$ and $S$ are two distinct row-standard
$n$-tableaux of shape $D$, then $\overline{T}\neq\overline{S}$.
\end{statement}

Next, we recall a fundamental criterion for when a sum of submodules is direct:

\begin{statement}
\textit{Observation 2:} Let $M$ be a $\mathbf{k}$-module, and let $\left(
N_{T}\right)  _{T\in I}$ be a finite family of $\mathbf{k}$-submodules of $M$.
Then, the sum $\sum_{T\in I}N_{T}$ is direct if and only if the only family
$\left(  \mathbf{a}_{T}\right)  _{T\in I}\in\prod_{T\in I}N_{T}$ that
satisfies $\sum_{T\in I}\mathbf{a}_{T}=0$ is the family $\left(  0\right)
_{T\in I}$.
\end{statement}

This observation is just part of Proposition \ref{prop.mod.dirsum=uniq},
restated with a different indexing. \medskip

Now, Lemma \ref{lem.specht.AElam.direct} is claiming that the sum $\sum
_{T\in\operatorname*{SYT}\left(  D\right)  }\mathcal{A}\mathbf{E}_{T}$ is
direct. This is what we need to prove. By Observation 2 (applied to
$M=\mathcal{A}$ and $I=\operatorname*{SYT}\left(  D\right)  $ and
$N_{T}=\mathcal{A}\mathbf{E}_{T}$), it suffices to show that the only family
$\left(  \mathbf{a}_{T}\right)  _{T\in\operatorname*{SYT}\left(  D\right)
}\in\prod_{T\in\operatorname*{SYT}\left(  D\right)  }\mathcal{A}\mathbf{E}%
_{T}$ that satisfies $\sum_{T\in\operatorname*{SYT}\left(  D\right)
}\mathbf{a}_{T}=0$ is the family $\left(  0\right)  _{T\in\operatorname*{SYT}%
\left(  D\right)  }$. So let us prove this.

Let $\left(  \mathbf{a}_{T}\right)  _{T\in\operatorname*{SYT}\left(  D\right)
}\in\prod_{T\in\operatorname*{SYT}\left(  D\right)  }\mathcal{A}\mathbf{E}%
_{T}$ be a family that satisfies $\sum_{T\in\operatorname*{SYT}\left(
D\right)  }\mathbf{a}_{T}=0$. We must prove that $\left(  \mathbf{a}%
_{T}\right)  _{T\in\operatorname*{SYT}\left(  D\right)  }=\left(  0\right)
_{T\in\operatorname*{SYT}\left(  D\right)  }$.

Assume the contrary. Thus, $\left(  \mathbf{a}_{T}\right)  _{T\in
\operatorname*{SYT}\left(  D\right)  }\neq\left(  0\right)  _{T\in
\operatorname*{SYT}\left(  D\right)  }$. Renaming the index $T$ as $S$, we can
rewrite this as follows: $\left(  \mathbf{a}_{S}\right)  _{S\in
\operatorname*{SYT}\left(  D\right)  }\neq\left(  0\right)  _{S\in
\operatorname*{SYT}\left(  D\right)  }$. In other words, there exists some
$n$-tableau $S\in\operatorname*{SYT}\left(  D\right)  $ such that
$\mathbf{a}_{S}\neq0$. Among all such $n$-tableaux $S$, let us choose one for
which the $n$-tabloid $\overline{S}$ is smallest with respect to the Young
last letter order (this is the total order constructed in Proposition
\ref{prop.spechtmod.row.order}). Thus, $S\in\operatorname*{SYT}\left(
D\right)  $ is an $n$-tableau such that $\mathbf{a}_{S}\neq0$, and moreover,
it has the smallest $n$-tabloid $\overline{S}$ among all such $n$-tableaux.
The \textquotedblleft moreover\textquotedblright\ part of the preceding
sentence can be restated as follows:%
\begin{align}
&  \text{Each }n\text{-tableau }T\in\operatorname*{SYT}\left(  D\right)
\text{ that satisfies }\mathbf{a}_{T}\neq0\nonumber\\
&  \text{must satisfy }\overline{T}\geq\overline{S}.
\label{pf.lem.specht.AElam.direct.max1}%
\end{align}

We recall that $\left(  \mathbf{a}_{T}\right)  _{T\in\operatorname*{SYT}%
\left(  D\right)  }\in\prod_{T\in\operatorname*{SYT}\left(  D\right)
}\mathcal{A}\mathbf{E}_{T}$. In other words,%
\begin{equation}
\mathbf{a}_{T}\in\mathcal{A}\mathbf{E}_{T}\ \ \ \ \ \ \ \ \ \ \text{for each
}T\in\operatorname*{SYT}\left(  D\right)  .
\label{pf.lem.specht.AElam.direct.aT-in-AET}%
\end{equation}

Now, we shall prove that each $n$-tableau $T\in\operatorname*{SYT}\left(
D\right)  $ that is distinct from $S$ must satisfy%
\begin{equation}
\mathbf{a}_{T}\mathbf{E}_{S}=0. \label{pf.lem.specht.AElam.direct.aTES=0}%
\end{equation}

\begin{proof}
[Proof of (\ref{pf.lem.specht.AElam.direct.aTES=0}):]Let $T\in
\operatorname*{SYT}\left(  D\right)  $ be an $n$-tableau that is distinct from
$S$. We must prove that $\mathbf{a}_{T}\mathbf{E}_{S}=0$.

Assume the contrary. Thus, $\mathbf{a}_{T}\mathbf{E}_{S}\neq0$, so that
$\mathbf{a}_{T}\neq0$. Hence, (\ref{pf.lem.specht.AElam.direct.max1}) yields
$\overline{T}\geq\overline{S}$. However, both $n$-tableaux $T$ and $S$ are
standard (since they belong to $\operatorname*{SYT}\left(  D\right)  $) and
therefore row-standard. Since $T$ and $S$ are furthermore distinct, we thus
conclude (by Observation 1) that $\overline{T}\neq\overline{S}$. Combining
this with $\overline{T}\geq\overline{S}$, we obtain $\overline{T}>\overline
{S}$. Moreover, $T$ and $S$ are two standard $n$-tableaux of shape $D=Y\left(
\lambda\right)  $. Thus, Lemma \ref{lem.specht.AElam.triang} (applied to $T$
and $S$ instead of $S$ and $T$) yields $\mathbf{E}_{T}\mathbf{E}_{S}=0$ (since
$\overline{T}>\overline{S}$).

However, $\mathbf{a}_{T}\in\mathcal{A}\mathbf{E}_{T}$ (by
(\ref{pf.lem.specht.AElam.direct.aT-in-AET})). In other words, $\mathbf{a}%
_{T}=\mathbf{bE}_{T}$ for some $\mathbf{b}\in\mathcal{A}$. Consider this
$\mathbf{b}$. Then, $\underbrace{\mathbf{a}_{T}}_{=\mathbf{bE}_{T}}%
\mathbf{E}_{S}=\mathbf{b}\underbrace{\mathbf{E}_{T}\mathbf{E}_{S}}_{=0}=0$.
This contradicts $\mathbf{a}_{T}\mathbf{E}_{S}\neq0$.

This contradiction shows that our assumption was false. Hence, $\mathbf{a}%
_{T}\mathbf{E}_{S}=0$ is proved. In other words,
(\ref{pf.lem.specht.AElam.direct.aTES=0}) is proved.
\end{proof}

Furthermore, we claim that%
\begin{equation}
\mathbf{a}_{S}\mathbf{E}_{S}=h^{\lambda}\mathbf{a}_{S}.
\label{pf.lem.specht.AElam.direct.aSES=}%
\end{equation}

\begin{proof}
[Proof of (\ref{pf.lem.specht.AElam.direct.aSES=}):]We know that $S$ is an
$n$-tableau of shape $D=Y\left(  \lambda\right)  $. Thus, Theorem
\ref{thm.specht.ETidp} (applied to $T=S$) yields $\mathbf{E}_{S}%
^{2}=\underbrace{\dfrac{n!}{f^{\lambda}}}_{\substack{=h^{\lambda}\\\text{(by
the definition of }h^{\lambda}\text{)}}}\mathbf{E}_{S}=h^{\lambda}%
\mathbf{E}_{S}$.

However, (\ref{pf.lem.specht.AElam.direct.aT-in-AET}) (applied to $T=S$) shows
that $\mathbf{a}_{S}\in\mathcal{A}\mathbf{E}_{S}$. In other words,
$\mathbf{a}_{S}=\mathbf{bE}_{S}$ for some $\mathbf{b}\in\mathcal{A}$. Consider
this $\mathbf{b}$. Thus,
\[
\underbrace{\mathbf{a}_{S}}_{=\mathbf{bE}_{S}}\mathbf{E}_{S}=\mathbf{b}%
\underbrace{\mathbf{E}_{S}\mathbf{E}_{S}}_{=\mathbf{E}_{S}^{2}=h^{\lambda
}\mathbf{E}_{S}}=\mathbf{b}h^{\lambda}\mathbf{E}_{S}=h^{\lambda}%
\underbrace{\mathbf{bE}_{S}}_{=\mathbf{a}_{S}}=h^{\lambda}\mathbf{a}_{S}.
\]
This proves (\ref{pf.lem.specht.AElam.direct.aSES=}).
\end{proof}

Now, recall that $\sum_{T\in\operatorname*{SYT}\left(  D\right)  }%
\mathbf{a}_{T}=0$. Multiplying this equality by $\mathbf{E}_{S}$ from the
right, we obtain $\left(  \sum_{T\in\operatorname*{SYT}\left(  D\right)
}\mathbf{a}_{T}\right)  \mathbf{E}_{S}=0\mathbf{E}_{S}=0$. Thus,%
\begin{align*}
0  &  =\left(  \sum_{T\in\operatorname*{SYT}\left(  D\right)  }\mathbf{a}%
_{T}\right)  \mathbf{E}_{S}=\sum_{T\in\operatorname*{SYT}\left(  D\right)
}\mathbf{a}_{T}\mathbf{E}_{S}\\
&  =\underbrace{\mathbf{a}_{S}\mathbf{E}_{S}}_{\substack{=h^{\lambda
}\mathbf{a}_{S}\\\text{(by (\ref{pf.lem.specht.AElam.direct.aSES=}))}}%
}+\sum_{\substack{T\in\operatorname*{SYT}\left(  D\right)  ;\\T\neq
S}}\underbrace{\mathbf{a}_{T}\mathbf{E}_{S}}_{\substack{=0\\\text{(by
(\ref{pf.lem.specht.AElam.direct.aTES=0}))}}}\ \ \ \ \ \ \ \ \ \ \left(
\begin{array}
[c]{c}%
\text{here, we have split off the}\\
\text{addend for }T=S\text{ from the sum}%
\end{array}
\right) \\
&  =h^{\lambda}\mathbf{a}_{S}+\underbrace{\sum_{\substack{T\in
\operatorname*{SYT}\left(  D\right)  ;\\T\neq S}}0}_{=0}=h^{\lambda}%
\mathbf{a}_{S}.
\end{align*}
We can divide this equality by $h^{\lambda}$ (since $h^{\lambda}$ is
invertible in $\mathbf{k}$), and thus obtain $0=\mathbf{a}_{S}$. But this
contradicts $\mathbf{a}_{S}\neq0$.

This contradiction shows that our assumption was false. Hence, we have proved
that $\left(  \mathbf{a}_{T}\right)  _{T\in\operatorname*{SYT}\left(
D\right)  }=\left(  0\right)  _{T\in\operatorname*{SYT}\left(  D\right)  }$.
As explained above, this completes our proof of Lemma
\ref{lem.specht.AElam.direct}.
\end{proof}

\begin{lemma}
\label{lem.specht.AElam.cts-AET}Let $\lambda$ be a partition of $n$. Assume
that $h^{\lambda}$ is invertible in $\mathbf{k}$. Let $T$ be an $n$-tableau of
shape $Y\left(  \lambda\right)  $. Then,%
\[
\mathcal{A}\mathbf{E}_{T}\subseteq\mathcal{A}\mathbf{E}_{\lambda}.
\]

\end{lemma}

\begin{proof}
Proposition \ref{prop.spechtmod.Elam.ETEl} yields $\mathbf{E}_{\lambda
}\mathbf{E}_{T}=\left(  h^{\lambda}\right)  ^{2}\mathbf{E}_{T}$. But
Proposition \ref{prop.spechtmod.Elam.incent} yields that $\mathbf{E}_{\lambda
}$ belongs to the center $Z\left(  \mathbf{k}\left[  S_{n}\right]  \right)  $
of $\mathbf{k}\left[  S_{n}\right]  $. Hence, $\mathbf{E}_{\lambda}%
\mathbf{E}_{T}=\mathbf{E}_{T}\mathbf{E}_{\lambda}$. Comparing this with
$\mathbf{E}_{\lambda}\mathbf{E}_{T}=\left(  h^{\lambda}\right)  ^{2}%
\mathbf{E}_{T}$, we find $\left(  h^{\lambda}\right)  ^{2}\mathbf{E}%
_{T}=\mathbf{E}_{T}\mathbf{E}_{\lambda}$. We can divide this equality by
$\left(  h^{\lambda}\right)  ^{2}$ (since the number $h^{\lambda}$ and thus
also its square $\left(  h^{\lambda}\right)  ^{2}$ is invertible in
$\mathbf{k}$), and thus obtain%
\[
\mathbf{E}_{T}=\underbrace{\dfrac{1}{\left(  h^{\lambda}\right)  ^{2}%
}\mathbf{E}_{T}}_{\in\mathcal{A}}\mathbf{E}_{\lambda}\in\mathcal{A}%
\mathbf{E}_{\lambda}.
\]
Hence, $\mathcal{A}\underbrace{\mathbf{E}_{T}}_{\in\mathcal{A}\mathbf{E}%
_{\lambda}}\subseteq\underbrace{\mathcal{A}\mathcal{A}}_{\subseteq\mathcal{A}%
}\mathbf{E}_{\lambda}\subseteq\mathcal{A}\mathbf{E}_{\lambda}$. This proves
Lemma \ref{lem.specht.AElam.cts-AET}.
\end{proof}

\begin{proof}
[Proof of Theorem \ref{thm.specht.AElam.decomp}.]\textit{Step 1:} Lemma
\ref{lem.specht.AElam.direct} shows that the sum $\sum_{T\in
\operatorname*{SYT}\left(  D\right)  }\mathcal{A}\mathbf{E}_{T}$ is direct.
Thus, the internal direct sum $\bigoplus\limits_{T\in\operatorname*{SYT}%
\left(  D\right)  }\mathcal{A}\mathbf{E}_{T}$ is well-defined. Let us denote
this internal direct sum by $\Omega$. Thus, $\Omega=\bigoplus\limits_{T\in
\operatorname*{SYT}\left(  D\right)  }\mathcal{A}\mathbf{E}_{T}$ is a
$\mathbf{k}$-submodule of $\mathcal{A}$. \medskip

\textit{Step 2:} We shall now show that
\begin{equation}
\Omega\subseteq\mathcal{A}\mathbf{E}_{\lambda}.
\label{pf.thm.specht.AElam.decomp.3}%
\end{equation}

Indeed, we have%
\[
\Omega=\bigoplus_{T\in\operatorname*{SYT}\left(  D\right)  }\mathcal{A}%
\mathbf{E}_{T}=\sum_{T\in\operatorname*{SYT}\left(  D\right)  }%
\underbrace{\mathcal{A}\mathbf{E}_{T}}_{\substack{\subseteq\mathcal{A}%
\mathbf{E}_{\lambda}\\\text{(by Lemma \ref{lem.specht.AElam.cts-AET}%
,}\\\text{since }T\in\operatorname*{SYT}\left(  D\right)  \text{
shows}\\\text{that }T\text{ is a standard }n\text{-tableau}\\\text{of shape
}D=Y\left(  \lambda\right)  \text{)}}}\subseteq\sum_{T\in\operatorname*{SYT}%
\left(  D\right)  }\mathcal{A}\mathbf{E}_{\lambda}\subseteq\mathcal{A}%
\mathbf{E}_{\lambda}%
\]
(since $\mathcal{A}\mathbf{E}_{\lambda}$ is a $\mathbf{k}$-module). This
proves (\ref{pf.thm.specht.AElam.decomp.3}). \medskip

\textit{Step 3:} Let us now show that
\begin{equation}
\mathcal{A}\mathbf{E}_{\lambda}\subseteq\Omega.
\label{pf.thm.specht.AElam.decomp.4}%
\end{equation}
This, too, is quite easy: Each $U,V\in\operatorname*{SYT}\left(  D\right)  $
satisfy%
\begin{align*}
\mathbf{E}_{U,V}  &  =\underbrace{w_{U,V}}_{\in\mathcal{A}}\mathbf{E}%
_{V}\ \ \ \ \ \ \ \ \ \ \left(  \text{by (\ref{eq.prop.specht.EPQ.EPEQ.wE}),
applied to }P=U\text{ and }Q=V\right) \\
&  \in\mathcal{A}\mathbf{E}_{V}\subseteq\bigoplus_{T\in\operatorname*{SYT}%
\left(  D\right)  }\mathcal{A}\mathbf{E}_{T}%
\end{align*}
(since $\mathcal{A}\mathbf{E}_{V}$ is one of the addends of the internal
direct sum $\bigoplus\limits_{T\in\operatorname*{SYT}\left(  D\right)
}\mathcal{A}\mathbf{E}_{T}$). In other words,%
\[
\left\{  \mathbf{E}_{U,V}\ \mid\ U,V\in\operatorname*{SYT}\left(  D\right)
\right\}  \subseteq\bigoplus\limits_{T\in\operatorname*{SYT}\left(  D\right)
}\mathcal{A}\mathbf{E}_{T}=\Omega.
\]
But Proposition \ref{prop.specht.AElam.span} \textbf{(b)} yields%
\[
\mathcal{A}\mathbf{E}_{\lambda}\subseteq\operatorname*{span}%
\nolimits_{\mathbf{k}}\underbrace{\left\{  \mathbf{E}_{U,V}\ \mid
\ U,V\in\operatorname*{SYT}\left(  D\right)  \right\}  }_{\subseteq\Omega
}\subseteq\operatorname*{span}\nolimits_{\mathbf{k}}\Omega\subseteq\Omega
\]
(since $\Omega$ is a $\mathbf{k}$-module). This proves
(\ref{pf.thm.specht.AElam.decomp.4}). \medskip

\textit{Step 4:} Combining (\ref{pf.thm.specht.AElam.decomp.4}) with
(\ref{pf.thm.specht.AElam.decomp.3}), we find%
\[
\mathcal{A}\mathbf{E}_{\lambda}=\Omega=\bigoplus\limits_{T\in
\operatorname*{SYT}\left(  D\right)  }\mathcal{A}\mathbf{E}_{T}.
\]
This proves Theorem \ref{thm.specht.AElam.decomp}.
\end{proof}

Theorem \ref{thm.specht.AElam.decomp} allows us to describe $\mathcal{A}%
\mathbf{E}_{\lambda}$ as a direct sum of Specht modules:

\begin{corollary}
\label{cor.specht.AElam.decomp1}Let $\lambda$ be a partition of $n$. Assume
that $h^{\lambda}$ is invertible in $\mathbf{k}$. Then,%
\[
\mathcal{A}\mathbf{E}_{\lambda}\cong\left(  \mathcal{S}^{\lambda}\right)
^{f^{\lambda}}\ \ \ \ \ \ \ \ \ \ \text{as an }S_{n}\text{-representation}.
\]
Here, $M^{k}$ just means the external direct sum $\underbrace{M\oplus
M\oplus\cdots\oplus M}_{k\text{ times}}$ whenever $M$ is an $S_{n}%
$-representation and $k\in\mathbb{N}$.
\end{corollary}

\begin{proof}
Let $D=Y\left(  \lambda\right)  $. Let $\operatorname*{SYT}\left(  D\right)  $
be the set of all standard $n$-tableaux of shape $D$.

From $D=Y\left(  \lambda\right)  $, we obtain $\left\vert \operatorname*{SYT}%
\left(  D\right)  \right\vert =\left\vert \operatorname*{SYT}\left(  Y\left(
\lambda\right)  \right)  \right\vert =f^{\lambda}$ (since we established the
equality $f^{\lambda}=\left\vert \operatorname*{SYT}\left(  Y\left(
\lambda\right)  \right)  \right\vert $ during our proof of Corollary
\ref{cor.tableaux.flamt}). From $D=Y\left(  \lambda\right)  $, we also obtain
$\mathcal{S}^{D}=\mathcal{S}^{Y\left(  \lambda\right)  }=\mathcal{S}^{\lambda
}$.

Let $T\in\operatorname*{SYT}\left(  D\right)  $. Then, $T$ is a standard
$n$-tableau of shape $D=Y\left(  \lambda\right)  $. Hence,
(\ref{eq.def.specht.ET.defs.SD=AET}) yields that $\mathcal{S}^{D}%
\cong\mathcal{A}\mathbf{E}_{T}$ as left $\mathbf{k}\left[  S_{n}\right]
$-modules. Hence, $\mathcal{A}\mathbf{E}_{T}\cong\mathcal{S}^{D}%
=\mathcal{S}^{\lambda}$ as left $\mathbf{k}\left[  S_{n}\right]  $-modules.

Forget that we fixed $T$. We thus have shown that
\begin{equation}
\mathcal{A}\mathbf{E}_{T}\cong\mathcal{S}^{\lambda}\text{ as left }%
\mathbf{k}\left[  S_{n}\right]  \text{-modules}
\label{pf.cor.specht.AElam.decomp1.1}%
\end{equation}
for each $T\in\operatorname*{SYT}\left(  D\right)  $.

But Theorem \ref{thm.specht.AElam.decomp} yields%
\begin{align*}
\mathcal{A}\mathbf{E}_{\lambda}  &  =\bigoplus\limits_{T\in\operatorname*{SYT}%
\left(  D\right)  }\underbrace{\mathcal{A}\mathbf{E}_{T}}_{\substack{\cong%
\mathcal{S}^{\lambda}\\\text{(by (\ref{pf.cor.specht.AElam.decomp1.1}))}}}\\
&  \cong\bigoplus\limits_{T\in\operatorname*{SYT}\left(  D\right)
}\mathcal{S}^{\lambda}\ \ \ \ \ \ \ \ \ \ \left(  \text{an external direct
sum}\right) \\
&  =\left(  \mathcal{S}^{\lambda}\right)  ^{\left\vert \operatorname*{SYT}%
\left(  D\right)  \right\vert }=\left(  \mathcal{S}^{\lambda}\right)
^{f^{\lambda}}\ \ \ \ \ \ \ \ \ \ \left(  \text{since }\left\vert
\operatorname*{SYT}\left(  D\right)  \right\vert =f^{\lambda}\right)
\end{align*}
as left $\mathbf{k}\left[  S_{n}\right]  $-modules. Hence, $\mathcal{A}%
\mathbf{E}_{\lambda}\cong\left(  \mathcal{S}^{\lambda}\right)  ^{f^{\lambda}}$
as $S_{n}$-representations (since Corollary \ref{cor.rep.G-rep.iso=iso} easily
yields that isomorphic left $\mathbf{k}\left[  S_{n}\right]  $-modules are
isomorphic $S_{n}$-representations). This proves Corollary
\ref{cor.specht.AElam.decomp1}.
\end{proof}

We also obtain a decomposition of the left regular representation
$\mathcal{A}=\mathbf{k}\left[  S_{n}\right]  $ of $S_{n}$:

\begin{corollary}
\label{cor.specht.A.decomp}For any partition $\lambda$ of $n$, let
$\operatorname*{SYT}\left(  \lambda\right)  $ be the set of all standard
$n$-tableaux of shape $Y\left(  \lambda\right)  $. Assume that $n!$ is
invertible in $\mathbf{k}$. Then,%
\begin{align*}
\mathcal{A}  &  =\bigoplus_{\lambda\text{ is a partition of }n}\ \ \bigoplus
_{T\in\operatorname*{SYT}\left(  \lambda\right)  }\mathcal{A}\mathbf{E}%
_{T}\ \ \ \ \ \ \ \ \ \ \left(  \text{an internal direct sum}\right) \\
&  \cong\bigoplus_{\lambda\text{ is a partition of }n}\left(  \mathcal{S}%
^{\lambda}\right)  ^{f^{\lambda}}\ \ \ \ \ \ \ \ \ \ \text{as an }%
S_{n}\text{-representation}.
\end{align*}
(See Corollary \ref{cor.specht.AElam.decomp1} for the meaning of $\left(
\mathcal{S}^{\lambda}\right)  ^{f^{\lambda}}$.)
\end{corollary}

\begin{proof}
Let $\lambda$ be any partition of $n$. Then, the number $h^{\lambda}$ is
invertible in $\mathbf{k}$ (as we saw in the proof of Corollary
\ref{cor.spechtmod.Elam.pou} \textbf{(a)}).

Let $D=Y\left(  \lambda\right)  $. Then, $\operatorname*{SYT}\left(
\lambda\right)  $ is the set of all standard $n$-tableaux of shape $D$ (since
$\operatorname*{SYT}\left(  \lambda\right)  $ is the set of all standard
$n$-tableaux of shape $Y\left(  \lambda\right)  $). In other words,
$\operatorname*{SYT}\left(  \lambda\right)  $ is what we called
$\operatorname*{SYT}\left(  D\right)  $ in Theorem
\ref{thm.specht.AElam.decomp}. Thus, Theorem \ref{thm.specht.AElam.decomp}
(with $\operatorname*{SYT}\left(  D\right)  $ renamed as $\operatorname*{SYT}%
\left(  \lambda\right)  $) yields that%
\begin{equation}
\mathcal{A}\mathbf{E}_{\lambda}=\bigoplus_{T\in\operatorname*{SYT}\left(
\lambda\right)  }\mathcal{A}\mathbf{E}_{T}\ \ \ \ \ \ \ \ \ \ \left(  \text{an
internal direct sum}\right)  . \label{pf.cor.specht.A.decomp.1}%
\end{equation}
Moreover, Corollary \ref{cor.specht.AElam.decomp1} shows that%
\begin{equation}
\mathcal{A}\mathbf{E}_{\lambda}\cong\left(  \mathcal{S}^{\lambda}\right)
^{f^{\lambda}}\ \ \ \ \ \ \ \ \ \ \text{as an }S_{n}\text{-representation}.
\label{pf.cor.specht.A.decomp.2}%
\end{equation}

Forget that we fixed $\lambda$. We thus have proved
(\ref{pf.cor.specht.A.decomp.1}) and (\ref{pf.cor.specht.A.decomp.2}) for each
partition $\lambda$ of $n$. Now, Corollary \ref{cor.spechtmod.Elam.A=sum1}
yields
\begin{align*}
\mathcal{A}  &  =\bigoplus_{\lambda\text{ is a partition of }n}%
\ \ \underbrace{\mathcal{A}\mathbf{E}_{\lambda}}_{\substack{=\bigoplus
\limits_{T\in\operatorname*{SYT}\left(  \lambda\right)  }\mathcal{A}%
\mathbf{E}_{T}\\\text{(by (\ref{pf.cor.specht.A.decomp.1}))}}%
}\ \ \ \ \ \ \ \ \ \ \left(  \text{an internal direct sum}\right) \\
&  =\bigoplus_{\lambda\text{ is a partition of }n}\ \ \bigoplus_{T\in
\operatorname*{SYT}\left(  \lambda\right)  }\mathcal{A}\mathbf{E}%
_{T}\ \ \ \ \ \ \ \ \ \ \left(  \text{an internal direct sum}\right)  .
\end{align*}
Moreover,%
\begin{align*}
\mathcal{A}  &  =\bigoplus_{\lambda\text{ is a partition of }n}%
\underbrace{\mathcal{A}\mathbf{E}_{\lambda}}_{\substack{\cong\left(
\mathcal{S}^{\lambda}\right)  ^{f^{\lambda}}\\\text{(by
(\ref{pf.cor.specht.A.decomp.2}))}}}\\
&  \cong\bigoplus_{\lambda\text{ is a partition of }n}\left(  \mathcal{S}%
^{\lambda}\right)  ^{f^{\lambda}}\ \ \ \ \ \ \ \ \ \ \left(  \text{an external
direct sum}\right)
\end{align*}
as left $\mathcal{A}$-modules, i.e., as left $\mathbf{k}\left[  S_{n}\right]
$-modules (since $\mathcal{A}=\mathbf{k}\left[  S_{n}\right]  $), i.e., as
$S_{n}$-representations (since Corollary \ref{cor.rep.G-rep.iso=iso} easily
yields that isomorphic left $\mathbf{k}\left[  S_{n}\right]  $-modules are
isomorphic $S_{n}$-representations). Thus, our proof of Corollary
\ref{cor.specht.A.decomp} is complete.
\end{proof}

If $\mathbf{k}$ is a field of characteristic $0$, then Corollary
\ref{cor.specht.A.decomp} gives a decomposition of the left regular
representation $\mathcal{A}=\mathbf{k}\left[  S_{n}\right]  $ of $S_{n}$ into
a direct sum of irreducible representations (since Corollary
\ref{cor.spechtmod.irred} shows that the Specht modules $\mathcal{S}^{\lambda
}$ are irreducible).

\subsubsection{The Young symmetrizer bases of $\mathcal{A}\mathbf{E}_{\lambda
}$ and $\mathcal{A}$}

We can now find a basis of the $\mathbf{k}$-module $\mathcal{A}\mathbf{E}%
_{\lambda}$ for any partition $\lambda$ for which $h^{\lambda}$ is invertible:

\begin{theorem}
\label{thm.specht.AElam.nat-basis}Let $\lambda$ be a partition of $n$. Let
$D=Y\left(  \lambda\right)  $. Let $\operatorname*{SYT}\left(  D\right)  $ be
the set of all standard $n$-tableaux of shape $D$. Assume that $h^{\lambda}$
is invertible in $\mathbf{k}$. Then, the family $\left(  \mathbf{E}%
_{U,V}\right)  _{U,V\in\operatorname*{SYT}\left(  D\right)  }$ is a basis of
the $\mathbf{k}$-module $\mathcal{A}\mathbf{E}_{\lambda}$.
\end{theorem}

\begin{remark}
This is not generally true if $h^{\lambda}$ is not invertible in $\mathbf{k}$.
In fact, in this general case, the elements $\mathbf{E}_{U,V}$ might not even
lie in the $\mathbf{k}$-module $\mathcal{A}\mathbf{E}_{\lambda}$ in the first
place. Moreover, they may fail to be $\mathbf{k}$-linearly independent. For an
example of the latter situation, let us consider the case $n=3$ and
$\lambda=\left(  2,1\right)  $. The formulas from Example
\ref{exa.specht.EPQ.21} easily entail%
\begin{align*}
&  \mathbf{E}_{12\backslash\backslash3,\ 12\backslash\backslash3}%
-\mathbf{E}_{12\backslash\backslash3,\ 13\backslash\backslash2}+\mathbf{E}%
_{13\backslash\backslash2,\ 12\backslash\backslash3}-\mathbf{E}_{13\backslash
\backslash2,\ 13\backslash\backslash2}\\
&  =-3t_{1,3}+3t_{1,2}-3\operatorname*{cyc}\nolimits_{1,2,3}%
+\,3\operatorname*{cyc}\nolimits_{1,3,2}\\
&  =0\ \ \ \ \ \ \ \ \ \ \text{if }3=0\text{ in }\mathbf{k}\text{.}%
\end{align*}
Thus, the family $\left(  \mathbf{E}_{U,V}\right)  _{U,V\in\operatorname*{SYT}%
\left(  D\right)  }$ is not linearly independent if $\mathbf{k}$ is a field of
characteristic $3$.
\end{remark}

\begin{proof}
[Proof of Theorem \ref{thm.specht.AElam.nat-basis}.]We first recall a basic
property of $\mathbf{k}$-modules:

\begin{statement}
\textit{Observation 1:} Let $M$ be a $\mathbf{k}$-module, and let $\left(
N_{T}\right)  _{T\in I}$ be a finite family of $\mathbf{k}$-submodules of $M$.
Assume that the sum $\sum_{T\in I}N_{T}$ is direct. Assume furthermore that
each $\mathbf{k}$-module $N_{T}$ has a basis $\left(  v_{T,j}\right)  _{j\in
J_{T}}$. Then, the internal direct sum $\bigoplus\limits_{T\in I}N_{T}$ has a
basis $\left(  v_{T,j}\right)  _{T\in I\text{, and }j\in J_{T}}$.
\end{statement}

Observation 1 is a basic fact in linear algebra and easy to prove. (For
example, in the case when $I=\left\{  1,2\right\}  $, this fact takes the
following simple form: If $N_{1}$ and $N_{2}$ are two $\mathbf{k}$-submodules
of $M$ such that the sum $N_{1}+N_{2}$ is direct, and if $\left(
v_{1,1},v_{1,2},\ldots,v_{1,i}\right)  $ and $\left(  v_{2,1},v_{2,2}%
,\ldots,v_{2,j}\right)  $ are bases of $N_{1}$ and $N_{2}$, respectively, then
the internal direct sum $N_{1}\oplus N_{2}$ has a basis $\left(
v_{1,1},v_{1,2},\ldots,v_{1,i},v_{2,1},v_{2,2},\ldots,v_{2,j}\right)  $. The
general case is no harder to prove.) \medskip

Let us note that $D=Y\left(  \lambda\right)  =Y\left(  \lambda/\varnothing
\right)  $, so that $D$ is a skew Young diagram.

Now, Theorem \ref{thm.specht.AElam.decomp} yields%
\[
\mathcal{A}\mathbf{E}_{\lambda}=\bigoplus\limits_{T\in\operatorname*{SYT}%
\left(  D\right)  }\mathcal{A}\mathbf{E}_{T}\ \ \ \ \ \ \ \ \ \ \left(
\text{an internal direct sum}\right)  .
\]
Thus, the sum $\sum_{T\in\operatorname*{SYT}\left(  D\right)  }\mathcal{A}%
\mathbf{E}_{T}$ is direct. Moreover, each $\mathbf{k}$-module $\mathcal{A}%
\mathbf{E}_{T}$ in this direct sum has a basis $\left(  \mathbf{E}%
_{P,T}\right)  _{P\in\operatorname*{SYT}\left(  D\right)  }$ (by Proposition
\ref{prop.specht.EPQ.basis-AET}). Thus, Observation 1 (applied to
$M=\mathcal{A}$ and $I=\operatorname*{SYT}\left(  D\right)  $ and
$N_{T}=\mathcal{A}\mathbf{E}_{T}$ and $\left(  v_{T,j}\right)  _{j\in J_{T}%
}=\left(  \mathbf{E}_{P,T}\right)  _{P\in\operatorname*{SYT}\left(  D\right)
}$) shows that the internal direct sum $\bigoplus\limits_{T\in
\operatorname*{SYT}\left(  D\right)  }\mathcal{A}\mathbf{E}_{T}$ has a basis%
\[
\left(  \mathbf{E}_{P,T}\right)  _{T\in\operatorname*{SYT}\left(  D\right)
,\text{ and }P\in\operatorname*{SYT}\left(  D\right)  }.
\]
Since $\mathcal{A}\mathbf{E}_{\lambda}=\bigoplus\limits_{T\in
\operatorname*{SYT}\left(  D\right)  }\mathcal{A}\mathbf{E}_{T}$, we can
rewrite this as follows: The $\mathbf{k}$-module $\mathcal{A}\mathbf{E}%
_{\lambda}$ has a basis
\[
\left(  \mathbf{E}_{P,T}\right)  _{T\in\operatorname*{SYT}\left(  D\right)
,\text{ and }P\in\operatorname*{SYT}\left(  D\right)  }=\left(  \mathbf{E}%
_{P,T}\right)  _{P,T\in\operatorname*{SYT}\left(  D\right)  }=\left(
\mathbf{E}_{U,V}\right)  _{U,V\in\operatorname*{SYT}\left(  D\right)  }%
\]
(here, we have renamed the indices $P$ and $T$ as $U$ and $V$). This proves
Theorem \ref{thm.specht.AElam.nat-basis}.
\end{proof}

\begin{corollary}
\label{cor.specht.A.nat-basis}For any partition $\lambda$ of $n$, let
$\operatorname*{SYT}\left(  \lambda\right)  $ be the set of all standard
$n$-tableaux of shape $Y\left(  \lambda\right)  $. Assume that $n!$ is
invertible in $\mathbf{k}$. Then, the family $\left(  \mathbf{E}_{U,V}\right)
_{\lambda\text{ is a partition of }n\text{, and }U,V\in\operatorname*{SYT}%
\left(  \lambda\right)  }$ is a basis of the $\mathbf{k}$-module $\mathcal{A}$.
\end{corollary}

\begin{proof}
Corollary \ref{cor.spechtmod.Elam.A=sum1} yields
\[
\mathcal{A}=\bigoplus\limits_{\lambda\text{ is a partition of }n}%
\mathcal{A}\mathbf{E}_{\lambda}\ \ \ \ \ \ \ \ \ \ \left(  \text{an internal
direct sum}\right)  .
\]
Thus, the sum $\sum_{\lambda\text{ is a partition of }n}\mathcal{A}%
\mathbf{E}_{\lambda}$ is direct.

Let $\lambda$ be any partition of $n$. Then, the number $h^{\lambda}$ is
invertible in $\mathbf{k}$ (as we saw in the proof of Corollary
\ref{cor.spechtmod.Elam.pou} \textbf{(a)}).

Let $D=Y\left(  \lambda\right)  $. Then, $\operatorname*{SYT}\left(
\lambda\right)  $ is the set of all standard $n$-tableaux of shape $D$ (since
$\operatorname*{SYT}\left(  \lambda\right)  $ is the set of all standard
$n$-tableaux of shape $Y\left(  \lambda\right)  $). In other words,
$\operatorname*{SYT}\left(  \lambda\right)  $ is what we called
$\operatorname*{SYT}\left(  D\right)  $ in Theorem
\ref{thm.specht.AElam.nat-basis}. Thus, Theorem
\ref{thm.specht.AElam.nat-basis} (with $\operatorname*{SYT}\left(  D\right)  $
renamed as $\operatorname*{SYT}\left(  \lambda\right)  $) yields that the
family $\left(  \mathbf{E}_{U,V}\right)  _{U,V\in\operatorname*{SYT}\left(
\lambda\right)  }$ is a basis of the $\mathbf{k}$-module $\mathcal{A}%
\mathbf{E}_{\lambda}$.

Forget that we fixed $\lambda$. We thus have shown that the family $\left(
\mathbf{E}_{U,V}\right)  _{U,V\in\operatorname*{SYT}\left(  \lambda\right)  }$
is a basis of the $\mathbf{k}$-module $\mathcal{A}\mathbf{E}_{\lambda}$
whenever $\lambda$ is a partition of $n$.

Now, recall Observation 1 from our above proof of Theorem
\ref{thm.specht.AElam.nat-basis}. Let us restate this observation with the
index $T$ renamed as $\lambda$:

\begin{statement}
\textit{Observation 1:} Let $M$ be a $\mathbf{k}$-module, and let $\left(
N_{\lambda}\right)  _{\lambda\in I}$ be a finite family of $\mathbf{k}%
$-submodules of $M$. Assume that the sum $\sum_{\lambda\in I}N_{\lambda}$ is
direct. Assume furthermore that each $\mathbf{k}$-module $N_{\lambda}$ has a
basis $\left(  v_{\lambda,j}\right)  _{j\in J_{\lambda}}$. Then, the internal
direct sum $\bigoplus\limits_{\lambda\in I}N_{\lambda}$ has a basis $\left(
v_{\lambda,j}\right)  _{\lambda\in I\text{, and }j\in J_{\lambda}}$.
\end{statement}

Now, we can apply this observation to $M=\mathcal{A}$ and $I=\left\{
\text{partitions of }n\right\}  $ and $N_{\lambda}=\mathcal{A}\mathbf{E}%
_{\lambda}$ and $\left(  v_{\lambda,j}\right)  _{j\in J_{\lambda}}=\left(
\mathbf{E}_{U,V}\right)  _{U,V\in\operatorname*{SYT}\left(  \lambda\right)  }$
(since the sum $\sum_{\lambda\text{ is a partition of }n}\mathcal{A}%
\mathbf{E}_{\lambda}$ is direct, and since each $\mathbf{k}$-module
$\mathcal{A}\mathbf{E}_{\lambda}$ has a basis $\left(  \mathbf{E}%
_{U,V}\right)  _{U,V\in\operatorname*{SYT}\left(  \lambda\right)  }$). Thus,
we conclude that the internal direct sum $\bigoplus\limits_{\lambda\text{ is a
partition of }n}\mathcal{A}\mathbf{E}_{\lambda}$ has a basis \newline$\left(
\mathbf{E}_{U,V}\right)  _{\lambda\text{ is a partition of }n\text{, and
}U,V\in\operatorname*{SYT}\left(  \lambda\right)  }$. In other words, the
$\mathbf{k}$-module $\mathcal{A}$ has a basis $\left(  \mathbf{E}%
_{U,V}\right)  _{\lambda\text{ is a partition of }n\text{, and }%
U,V\in\operatorname*{SYT}\left(  \lambda\right)  }$ (since $\mathcal{A}%
=\bigoplus\limits_{\lambda\text{ is a partition of }n}\mathcal{A}%
\mathbf{E}_{\lambda}$). This proves Corollary \ref{cor.specht.A.nat-basis}.
\end{proof}

The basis $\left(  \mathbf{E}_{U,V}\right)  _{\lambda\text{ is a partition of
}n\text{, and }U,V\in\operatorname*{SYT}\left(  \lambda\right)  }$ in
Corollary \ref{cor.specht.A.nat-basis} is called (by myself at least) the
\emph{Young symmetrizer basis} of $\mathbf{k}\left[  S_{n}\right]  $.
%(due to its
%relation to what is called the \emph{Young symmetrizer basis} of a Specht module
%$\mathcal{S}^{\lambda}$, which is just the standard basis consisting of the
%standard polytabloids $\mathbf{e}_{T}$).

One neat application of the Young symmetrizer basis is the following purely
enumerative identity:

\begin{exercise}
\label{exe.specht.sum-flam-selfconj}\fbox{3} Prove that%
\[
\sum_{w\in S_{n}\text{ is an involution}}\left(  -1\right)  ^{w}%
=\sum_{\substack{\lambda\text{ is a partition of }n;\\\lambda^{t}=\lambda
}}\left(  -1\right)  ^{\operatorname*{upper}\left(  \lambda\right)
}f^{\lambda},
\]
where $\operatorname*{upper}\left(  \lambda\right)  $ denotes the number of
all cells $\left(  i,j\right)  \in Y\left(  \lambda\right)  $ satisfying
$i<j$. \medskip

[\textbf{Hint:} Consider the endomorphism $T_{\operatorname*{sign}}\circ S$ of
the $\mathbf{k}$-module $\mathcal{A}$. Compute the trace of this endomorphism
in two ways: once using the standard basis $\left(  w\right)  _{w\in S_{n}}$
and once using the Young symmetrizer basis $\left(  \mathbf{E}_{U,V}\right)
_{\lambda\text{ is a partition of }n\text{, and }U,V\in\operatorname*{SYT}%
\left(  \lambda\right)  }$.]
\end{exercise}

The identity in Exercise \ref{exe.specht.sum-flam-selfconj} has a
simpler-looking (and better-known) partner, which says that%
\[
\left(  \text{\# of involutions }w\in S_{n}\right)  =\sum_{\lambda\text{ is a
partition of }n}f^{\lambda}.
\]
This, too, can be proved by computing the trace of an endomorphism of
$\mathcal{A}$ using two different bases. For this, see Corollary
\ref{cor.spechtmod.sumflam} below. \medskip

For our proof of Theorem \ref{thm.specht.AWlam.sur}, we need one last lemma:

\begin{lemma}
\label{lem.specht.dim-AElam}Let $\lambda$ be a partition of $n$. Assume that
$h^{\lambda}$ is invertible in $\mathbf{k}$. Then, $\mathcal{A}\mathbf{E}%
_{\lambda}$ is a free $\mathbf{k}$-module of rank $\left(  f^{\lambda}\right)
^{2}$.
\end{lemma}

\begin{proof}
From Lemma \ref{lem.specht.Slam-flam}, we know that $\mathcal{S}^{\lambda}$ is
a free $\mathbf{k}$-module of rank $f^{\lambda}$.

Moreover, if $M$ is a free $\mathbf{k}$-module of some rank $r\in\mathbb{N}$,
and if $k\in\mathbb{N}$ is an arbitrary nonnegative integer, then the $k$-th
power $M^{k}=\underbrace{M\oplus M\oplus\cdots\oplus M}_{k\text{ times}}$ (an
external direct sum) is a free $\mathbf{k}$-module of rank $kr$ (this follows
easily from Lemma \ref{lem.linalg.dim-dirprod}, since a finite direct sum is
the same as a finite direct product). Applying this to $M=\mathcal{S}%
^{\lambda}$ and $r=f^{\lambda}$ and $k=f^{\lambda}$, we conclude that $\left(
\mathcal{S}^{\lambda}\right)  ^{f^{\lambda}}$ is a free $\mathbf{k}$-module of
rank $f^{\lambda}f^{\lambda}$ (since $\mathcal{S}^{\lambda}$ is a free
$\mathbf{k}$-module of rank $f^{\lambda}$). In other words, $\left(
\mathcal{S}^{\lambda}\right)  ^{f^{\lambda}}$ is a free $\mathbf{k}$-module of
rank $\left(  f^{\lambda}\right)  ^{2}$ (since $f^{\lambda}f^{\lambda}=\left(
f^{\lambda}\right)  ^{2}$).

Furthermore, Corollary \ref{cor.specht.AElam.decomp1} yields $\mathcal{A}%
\mathbf{E}_{\lambda}\cong\left(  \mathcal{S}^{\lambda}\right)  ^{f^{\lambda}}$
as $S_{n}$-representations and thus as $\mathbf{k}$-modules as well. Hence,
$\mathcal{A}\mathbf{E}_{\lambda}$ is a free $\mathbf{k}$-module of rank
$\left(  f^{\lambda}\right)  ^{2}$ (since $\left(  \mathcal{S}^{\lambda
}\right)  ^{f^{\lambda}}$ is a free $\mathbf{k}$-module of rank $\left(
f^{\lambda}\right)  ^{2}$). This proves Lemma \ref{lem.specht.dim-AElam}.
\end{proof}

\subsubsection{\label{subsec.specht.nat-basis.AWlam}Proof of Theorem
\ref{thm.specht.AWlam.sur}}

We shall now finally prove Theorem \ref{thm.specht.AWlam.sur}:

\begin{proof}
[Proof of Theorem \ref{thm.specht.AWlam.sur}.]This is rather similar to our
proof of Theorem \ref{thm.specht.AW}. The main difference is that we no longer
take sums over all partitions $\lambda$ of $n$, but instead work with our
chosen partition $\lambda$. Instead of the two maps $\Phi:\mathbf{k}\left[
S_{n}\right]  \rightarrow\prod_{\lambda\text{ is a partition of }%
n}\operatorname*{End}\nolimits_{\mathbf{k}}\left(  \mathcal{S}^{\lambda
}\right)  $ and $\Psi:\prod_{\lambda\text{ is a partition of }n}%
\operatorname*{End}\nolimits_{\mathbf{k}}\left(  \mathcal{S}^{\lambda}\right)
\rightarrow\mathbf{k}\left[  S_{n}\right]  $ (the second of which might fail
to exist if $n!$ is not invertible), we consider two maps $\Phi_{\lambda
}:\mathcal{A}\mathbf{E}_{\lambda}\rightarrow\operatorname*{End}%
\nolimits_{\mathbf{k}}\left(  \mathcal{S}^{\lambda}\right)  $ and
$\Psi_{\lambda}:\operatorname*{End}\nolimits_{\mathbf{k}}\left(
\mathcal{S}^{\lambda}\right)  \rightarrow\mathcal{A}\mathbf{E}_{\lambda}$
which are (in a sense) the \textquotedblleft$\lambda$-parts\textquotedblright%
\ of $\Phi$ and $\Psi$. Here is the proof in detail: \medskip

\textit{Step 1:} In the following, we shall abbreviate the summation sign
\newline\textquotedblleft$\sum_{T\text{ is an }n\text{-tableau of shape
}\lambda}$\textquotedblright\ as \textquotedblleft$\sum_{T\text{ of shape
}\lambda}$\textquotedblright. Using this abbreviation, we can rewrite the
equality (\ref{eq.def.spechtmod.Elam.Elam.def}) as%
\begin{equation}
\mathbf{E}_{\lambda}=\sum_{T\text{ of shape }\lambda}\mathbf{E}_{T}.
\label{pf.thm.specht.AWlam.sur.Elam=}%
\end{equation}

\textit{Step 2:} Corollary \ref{cor.spechtmod.Elam.subalg} \textbf{(a)} shows
that $\mathcal{A}\mathbf{E}_{\lambda}$ is a nonunital subalgebra of
$\mathcal{A}$. Furthermore, Corollary \ref{cor.spechtmod.Elam.subalg}
\textbf{(b)} shows that this nonunital algebra $\mathcal{A}\mathbf{E}%
_{\lambda}$ has a unity, namely $\dfrac{1}{\left(  h^{\lambda}\right)  ^{2}%
}\mathbf{E}_{\lambda}$, and thus is a unital $\mathbf{k}$-algebra. Let us
denote this unity by $1_{\lambda}$. Thus, $1_{\lambda}=\dfrac{1}{\left(
h^{\lambda}\right)  ^{2}}\mathbf{E}_{\lambda}$ is the unity of the
$\mathbf{k}$-algebra $\mathcal{A}\mathbf{E}_{\lambda}$.

Define a map%
\begin{align*}
\Phi_{\lambda}:\mathcal{A}\mathbf{E}_{\lambda}  &  \rightarrow
\operatorname*{End}\nolimits_{\mathbf{k}}\left(  \mathcal{S}^{\lambda}\right)
,\\
\mathbf{a}  &  \mapsto\rho_{\lambda}\left(  \mathbf{a}\right)  .
\end{align*}
In other words, $\Phi_{\lambda}$ is the restriction of the map $\rho_{\lambda
}$ to the $\mathbf{k}$-submodule $\mathcal{A}\mathbf{E}_{\lambda}$. One of our
goals is to prove that this map $\Phi_{\lambda}$ is a $\mathbf{k}$-algebra isomorphism.

As in Step 2 of the proof of Theorem \ref{thm.specht.AW}, we can see that
$\rho_{\lambda}:\mathbf{k}\left[  S_{n}\right]  \rightarrow\operatorname*{End}%
\nolimits_{\mathbf{k}}\left(  \mathcal{S}^{\lambda}\right)  $ is a
$\mathbf{k}$-algebra morphism.

Now, we shall show that $\Phi_{\lambda}$ is a $\mathbf{k}$-algebra morphism.
Indeed, the map $\Phi_{\lambda}$ is clearly $\mathbf{k}$-linear. Furthermore,
for any $\mathbf{a},\mathbf{b}\in\mathcal{A}\mathbf{E}_{\lambda}$, we have%
\begin{align*}
\Phi_{\lambda}\left(  \mathbf{ab}\right)   &  =\rho_{\lambda}\left(
\mathbf{ab}\right)  \ \ \ \ \ \ \ \ \ \ \left(  \text{by the definition of
}\Phi_{\lambda}\right) \\
&  =\underbrace{\rho_{\lambda}\left(  \mathbf{a}\right)  }_{\substack{=\Phi
_{\lambda}\left(  \mathbf{a}\right)  \\\text{(since the}\\\text{definition of
}\Phi_{\lambda}\\\text{yields }\Phi_{\lambda}\left(  \mathbf{a}\right)
=\rho_{\lambda}\left(  \mathbf{a}\right)  \text{)}}}\underbrace{\rho_{\lambda
}\left(  \mathbf{b}\right)  }_{\substack{=\Phi_{\lambda}\left(  \mathbf{b}%
\right)  \\\text{(since the}\\\text{definition of }\Phi_{\lambda
}\\\text{yields }\Phi_{\lambda}\left(  \mathbf{b}\right)  =\rho_{\lambda
}\left(  \mathbf{b}\right)  \text{)}}}\ \ \ \ \ \ \ \ \ \ \left(
\begin{array}
[c]{c}%
\text{since }\rho_{\lambda}\text{ is a }\mathbf{k}\text{-algebra}\\
\text{morphism}%
\end{array}
\right) \\
&  =\Phi_{\lambda}\left(  \mathbf{a}\right)  \Phi_{\lambda}\left(
\mathbf{b}\right)  .
\end{align*}
In other words, the map $\Phi_{\lambda}$ respects multiplication. Now, let us
compute the endomorphism $\Phi_{\lambda}\left(  1_{\lambda}\right)
\in\operatorname*{End}\nolimits_{\mathbf{k}}\left(  \mathcal{S}^{\lambda
}\right)  $. Indeed, for each $\mathbf{v}\in\mathcal{S}^{\lambda}$, we have%
\begin{align*}
\underbrace{\left(  \Phi_{\lambda}\left(  1_{\lambda}\right)  \right)
}_{\substack{=\rho_{\lambda}\left(  1_{\lambda}\right)  \\\text{(by the
definition of }\Phi_{\lambda}\text{)}}}\left(  \mathbf{v}\right)   &  =\left(
\rho_{\lambda}\left(  1_{\lambda}\right)  \right)  \left(  \mathbf{v}\right)
=1_{\lambda}\mathbf{v}\ \ \ \ \ \ \ \ \ \ \left(
\begin{array}
[c]{c}%
\text{by the definition of}\\
\text{the curried form }\rho_{\lambda}%
\end{array}
\right) \\
&  =\dfrac{1}{\left(  h^{\lambda}\right)  ^{2}}\underbrace{\mathbf{E}%
_{\lambda}\mathbf{v}}_{\substack{=\left(  h^{\lambda}\right)  ^{2}%
\mathbf{v}\\\text{(by Proposition \ref{prop.specht.Elam-Slam})}}%
}\ \ \ \ \ \ \ \ \ \ \left(  \text{since }1_{\lambda}=\dfrac{1}{\left(
h^{\lambda}\right)  ^{2}}\mathbf{E}_{\lambda}\right) \\
&  =\underbrace{\dfrac{1}{\left(  h^{\lambda}\right)  ^{2}}\left(  h^{\lambda
}\right)  ^{2}}_{=1}\mathbf{v}=\mathbf{v}.
\end{align*}
In other words, $\Phi_{\lambda}\left(  1_{\lambda}\right)  =\operatorname*{id}%
$. In other words, $\Phi_{\lambda}$ sends the unity of the $\mathbf{k}%
$-algebra $\mathcal{A}\mathbf{E}_{\lambda}$ to the unity of the $\mathbf{k}%
$-algebra $\operatorname*{End}\nolimits_{\mathbf{k}}\left(  \mathcal{S}%
^{\lambda}\right)  $ (since $1_{\lambda}$ is the unity of the $\mathbf{k}%
$-algebra $\mathcal{A}\mathbf{E}_{\lambda}$, whereas $\operatorname*{id}$ is
the unity of the $\mathbf{k}$-algebra $\operatorname*{End}%
\nolimits_{\mathbf{k}}\left(  \mathcal{S}^{\lambda}\right)  $). In other
words, the map $\Phi_{\lambda}$ respects the unity. Since $\Phi_{\lambda}$
also respects multiplication and is $\mathbf{k}$-linear, we thus conclude that
the map $\Phi_{\lambda}$ is a $\mathbf{k}$-algebra morphism.

We shall now focus on proving that $\Phi_{\lambda}$ is invertible. \medskip

\textit{Step 3:} Let $T$ be an $n$-tableau of shape $\lambda$. Then, there
exists an $\mathcal{A}$-module isomorphism $\alpha_{\lambda,T}:\mathcal{S}%
^{\lambda}\rightarrow\mathcal{A}\mathbf{E}_{T}$ (as we have shown in Step 3 of
the proof of Theorem \ref{thm.specht.AW}). Fix such an isomorphism
$\alpha_{\lambda,T}$ once and for all.

Forget that we fixed $T$. Thus, for any $n$-tableau $T$ of shape $\lambda$, we
have defined an $\mathcal{A}$-module isomorphism $\alpha_{\lambda
,T}:\mathcal{S}^{\lambda}\rightarrow\mathcal{A}\mathbf{E}_{T}$. In particular,
$\alpha_{\lambda,T}$ is thus a $\mathbf{k}$-module isomorphism. Therefore, if
$\zeta\in\operatorname*{End}\nolimits_{\mathbf{k}}\left(  \mathcal{S}%
^{\lambda}\right)  $ is a $\mathbf{k}$-module endomorphism of $\mathcal{S}%
^{\lambda}$, then $\alpha_{\lambda,T}\circ\zeta\circ\alpha_{\lambda,T}^{-1}%
\in\operatorname*{End}\nolimits_{\mathbf{k}}\left(  \mathcal{A}\mathbf{E}%
_{T}\right)  $ is a $\mathbf{k}$-module endomorphism of $\mathcal{A}%
\mathbf{E}_{T}$. \medskip

\textit{Step 4:} Now, define a map%
\[
\Psi_{\lambda}:\operatorname*{End}\nolimits_{\mathbf{k}}\left(  \mathcal{S}%
^{\lambda}\right)  \rightarrow\mathcal{A}\mathbf{E}_{\lambda},
\]
which sends each endomorphism $\zeta\in\operatorname*{End}%
\nolimits_{\mathbf{k}}\left(  \mathcal{S}^{\lambda}\right)  $ to%
\[
\dfrac{1}{\left(  h^{\lambda}\right)  ^{2}}\sum_{T\text{ of shape }\lambda
}\underbrace{\left(  \alpha_{\lambda,T}\circ\zeta\circ\alpha_{\lambda,T}%
^{-1}\right)  \left(  \mathbf{E}_{T}\right)  }_{\substack{\text{(this is
well-defined,}\\\text{since }\mathbf{E}_{T}=1\mathbf{E}_{T}\in\mathcal{A}%
\mathbf{E}_{T}\\\text{and }\alpha_{\lambda,T}\circ\zeta\circ\alpha_{\lambda
,T}^{-1}\in\operatorname*{End}\nolimits_{\mathbf{k}}\left(  \mathcal{A}%
\mathbf{E}_{T}\right)  \text{)}}}.
\]
This map $\Psi_{\lambda}$ is well-defined (since each $\zeta\in
\operatorname*{End}\nolimits_{\mathbf{k}}\left(  \mathcal{S}^{\lambda}\right)
$ satisfies%
\begin{align*}
\dfrac{1}{\left(  h^{\lambda}\right)  ^{2}}\sum_{T\text{ of shape }\lambda
}\underbrace{\left(  \alpha_{\lambda,T}\circ\zeta\circ\alpha_{\lambda,T}%
^{-1}\right)  \left(  \mathbf{E}_{T}\right)  }_{\substack{\in\mathcal{A}%
\mathbf{E}_{T}\\\text{(since }\alpha_{\lambda,T}\circ\zeta\circ\alpha
_{\lambda,T}^{-1}\in\operatorname*{End}\nolimits_{\mathbf{k}}\left(
\mathcal{A}\mathbf{E}_{T}\right)  \text{)}}}  &  \in\dfrac{1}{\left(
h^{\lambda}\right)  ^{2}}\sum_{T\text{ of shape }\lambda}%
\underbrace{\mathcal{A}\mathbf{E}_{T}}_{\substack{\subseteq\mathcal{A}%
\mathbf{E}_{\lambda}\\\text{(by Lemma \ref{lem.specht.AElam.cts-AET})}}}\\
&  \subseteq\dfrac{1}{\left(  h^{\lambda}\right)  ^{2}}\sum_{T\text{ of shape
}\lambda}\mathcal{A}\mathbf{E}_{\lambda}\subseteq\mathcal{A}\mathbf{E}%
_{\lambda}%
\end{align*}
(since $\mathcal{A}\mathbf{E}_{\lambda}$ is a $\mathbf{k}$-module)) and
$\mathbf{k}$-linear (since each $\left(  \alpha_{\lambda,T}\circ\zeta
\circ\alpha_{\lambda,T}^{-1}\right)  \left(  \mathbf{E}_{T}\right)  $ depends
$\mathbf{k}$-linearly on $\zeta$). \medskip

\textit{Step 5:} We shall next show that $\Psi_{\lambda}\circ\Phi_{\lambda
}=\operatorname*{id}$.

Indeed, let $\mathbf{a}\in\mathcal{A}\mathbf{E}_{\lambda}$. Then, the
definition of $\Phi_{\lambda}$ yields%
\[
\Phi_{\lambda}\left(  \mathbf{a}\right)  =\rho_{\lambda}\left(  \mathbf{a}%
\right)  .
\]
Hence,%
\[
\Psi_{\lambda}\left(  \Phi_{\lambda}\left(  \mathbf{a}\right)  \right)
=\Psi_{\lambda}\left(  \rho_{\lambda}\left(  \mathbf{a}\right)  \right)
=\dfrac{1}{\left(  h^{\lambda}\right)  ^{2}}\sum_{T\text{ of shape }\lambda
}\left(  \alpha_{\lambda,T}\circ\left(  \rho_{\lambda}\left(  \mathbf{a}%
\right)  \right)  \circ\alpha_{\lambda,T}^{-1}\right)  \left(  \mathbf{E}%
_{T}\right)
\]
(by the definition of $\Psi_{\lambda}$).

However, we can simplify the addends in this sum significantly, due to the
following claim:

\begin{statement}
\textit{Claim 1:} Let $T$ be an $n$-tableau of shape $\lambda$. Then,%
\[
\left(  \alpha_{\lambda,T}\circ\left(  \rho_{\lambda}\left(  \mathbf{a}%
\right)  \right)  \circ\alpha_{\lambda,T}^{-1}\right)  \left(  \mathbf{E}%
_{T}\right)  =\mathbf{aE}_{T}.
\]

\end{statement}

\begin{proof}
[Proof of Claim 1.]This is the same Claim 1 that we proved in Step 5 of our
proof of Theorem \ref{thm.specht.AW}. Thus, we need not prove it again.
\end{proof}

Now,
\begin{align*}
\left(  \Psi_{\lambda}\circ\Phi_{\lambda}\right)  \left(  \mathbf{a}\right)
&  =\Psi_{\lambda}\left(  \Phi_{\lambda}\left(  \mathbf{a}\right)  \right) \\
&  =\dfrac{1}{\left(  h^{\lambda}\right)  ^{2}}\sum_{T\text{ of shape }%
\lambda}\underbrace{\left(  \alpha_{\lambda,T}\circ\left(  \rho_{\lambda
}\left(  \mathbf{a}\right)  \right)  \circ\alpha_{\lambda,T}^{-1}\right)
\left(  \mathbf{E}_{T}\right)  }_{\substack{=\mathbf{aE}_{T}\\\text{(by Claim
1)}}}\\
&  \ \ \ \ \ \ \ \ \ \ \ \ \ \ \ \ \ \ \ \ \left(  \text{by our above
computation of }\Psi_{\lambda}\left(  \Phi_{\lambda}\left(  \mathbf{a}\right)
\right)  \right) \\
&  =\dfrac{1}{\left(  h^{\lambda}\right)  ^{2}}\sum_{T\text{ of shape }%
\lambda}\mathbf{aE}_{T}=\mathbf{a}\dfrac{1}{\left(  h^{\lambda}\right)  ^{2}%
}\underbrace{\sum_{T\text{ of shape }\lambda}\mathbf{E}_{T}}%
_{\substack{=\mathbf{E}_{\lambda}\\\text{(by
(\ref{pf.thm.specht.AWlam.sur.Elam=}))}}}\\
&  =\mathbf{a}\underbrace{\dfrac{1}{\left(  h^{\lambda}\right)  ^{2}%
}\mathbf{E}_{\lambda}}_{\substack{=1_{\lambda}\\\text{(by the definition of
}1_{\lambda}\text{)}}}=\mathbf{a}1_{\lambda}=\mathbf{a}%
\end{align*}
(since $\mathbf{a}\in\mathcal{A}\mathbf{E}_{\lambda}$, whereas $1_{\lambda}$
is the unity of the $\mathbf{k}$-algebra $\mathcal{A}\mathbf{E}_{\lambda}$).

Forget that we fixed $\mathbf{a}$. We thus have shown that $\left(
\Psi_{\lambda}\circ\Phi_{\lambda}\right)  \left(  \mathbf{a}\right)
=\mathbf{a}$ for each $\mathbf{a}\in\mathcal{A}\mathbf{E}_{\lambda}$. In other
words, $\Psi_{\lambda}\circ\Phi_{\lambda}=\operatorname*{id}$. \medskip

\textit{Step 6:} Next we want to show that $\Phi_{\lambda}\circ\Psi_{\lambda
}=\operatorname*{id}$ as well. We will derive this from $\Psi_{\lambda}%
\circ\Phi_{\lambda}=\operatorname*{id}$ using Lemma \ref{lem.linalg.AB=1-BA=1}%
. In order to do so, we need to check that $\mathcal{A}\mathbf{E}_{\lambda}$
and $\operatorname*{End}\nolimits_{\mathbf{k}}\left(  \mathcal{S}^{\lambda
}\right)  $ are two free $\mathbf{k}$-modules of the same rank (namely,
$\left(  f^{\lambda}\right)  ^{2}$).

Let us now do this.

Lemma \ref{lem.specht.Slam-flam} tells us that $\mathcal{S}^{\lambda}$ is a
free $\mathbf{k}$-module of rank $f^{\lambda}$. Hence, Lemma
\ref{lem.linalg.dim-End} (applied to $M=\mathcal{S}^{\lambda}$ and
$r=f^{\lambda}$) shows that $\operatorname*{End}\nolimits_{\mathbf{k}}\left(
\mathcal{S}^{\lambda}\right)  $ is a free $\mathbf{k}$-module of rank $\left(
f^{\lambda}\right)  ^{2}$. Moreover, Lemma \ref{lem.specht.dim-AElam} tells us
that $\mathcal{A}\mathbf{E}_{\lambda}$ is a free $\mathbf{k}$-module of rank
$\left(  f^{\lambda}\right)  ^{2}$.

Thus, Lemma \ref{lem.linalg.AB=1-BA=1} (applied to $V=\operatorname*{End}%
\nolimits_{\mathbf{k}}\left(  \mathcal{S}^{\lambda}\right)  $ and
$W=\mathcal{A}\mathbf{E}_{\lambda}$ and $r=\left(  f^{\lambda}\right)  ^{2}$
and $f=\Psi_{\lambda}$ and $g=\Phi_{\lambda}$) shows that $\Phi_{\lambda}%
\circ\Psi_{\lambda}=\operatorname*{id}$ (since $\Psi_{\lambda}\circ
\Phi_{\lambda}=\operatorname*{id}$). \medskip

\textit{Step 7:} Combining $\Phi_{\lambda}\circ\Psi_{\lambda}%
=\operatorname*{id}$ with $\Psi_{\lambda}\circ\Phi_{\lambda}%
=\operatorname*{id}$, we see that the maps $\Phi_{\lambda}$ and $\Psi
_{\lambda}$ are mutually inverse. Hence, the map $\Phi_{\lambda}$ is
invertible, and thus is a $\mathbf{k}$-algebra isomorphism (since it is a
$\mathbf{k}$-algebra morphism).

In other words, the restriction of $\rho_{\lambda}$ to $\mathcal{A}%
\mathbf{E}_{\lambda}$ is a $\mathbf{k}$-algebra isomorphism (since
$\Phi_{\lambda}$ is the restriction of $\rho_{\lambda}$ to $\mathcal{A}%
\mathbf{E}_{\lambda}$). This proves Theorem \ref{thm.specht.AWlam.sur}
\textbf{(b)}. \medskip

\textit{Step 8:} It remains to prove Theorem \ref{thm.specht.AWlam.sur}
\textbf{(a)}.

But this is easy: We know that the map $\Phi_{\lambda}$ is invertible, thus
bijective, thus surjective. In other words, $\Phi_{\lambda}\left(
\mathcal{A}\mathbf{E}_{\lambda}\right)  =\operatorname*{End}%
\nolimits_{\mathbf{k}}\left(  \mathcal{S}^{\lambda}\right)  $. However,
$\Phi_{\lambda}$ is the restriction of $\rho_{\lambda}$ to $\mathcal{A}%
\mathbf{E}_{\lambda}$. Thus, the map $\Phi_{\lambda}$ is identical with
$\rho_{\lambda}$ on the submodule $\mathcal{A}\mathbf{E}_{\lambda}$. Hence,
$\Phi_{\lambda}\left(  \mathcal{A}\mathbf{E}_{\lambda}\right)  =\rho_{\lambda
}\left(  \underbrace{\mathcal{A}\mathbf{E}_{\lambda}}_{\subseteq\mathcal{A}%
}\right)  \subseteq\rho_{\lambda}\left(  \mathcal{A}\right)  $. Therefore,
$\rho_{\lambda}\left(  \mathcal{A}\right)  \supseteq\Phi_{\lambda}\left(
\mathcal{A}\mathbf{E}_{\lambda}\right)  =\operatorname*{End}%
\nolimits_{\mathbf{k}}\left(  \mathcal{S}^{\lambda}\right)  $. This entails
$\rho_{\lambda}\left(  \mathcal{A}\right)  =\operatorname*{End}%
\nolimits_{\mathbf{k}}\left(  \mathcal{S}^{\lambda}\right)  $ (since we
obviously have $\rho_{\lambda}\left(  \mathcal{A}\right)  \subseteq
\operatorname*{End}\nolimits_{\mathbf{k}}\left(  \mathcal{S}^{\lambda}\right)
$). In other words, the map $\rho_{\lambda}$ is surjective. Thus, Theorem
\ref{thm.specht.AWlam.sur} \textbf{(a)} is proved.
\end{proof}

We will learn more about the Young symmetrizer basis of $\mathcal{A}$ in
Section \ref{sec.bas.ysb}.

\subsection{\label{sec.specht.conj-class-action}The action of conjugacy class
sums on Specht modules}

\subsubsection{A general formula}

Let us next discuss how a conjugacy class sum of $S_{n}$ (see Definition
\ref{def.groups.conj} \textbf{(c)}) acts on a Specht module of the form
$\mathcal{S}^{\lambda}$. The answer turns out to be pretty nice: it just
scales every vector in $\mathcal{S}^{\lambda}$ by a certain number. Here is
the precise answer:

\begin{theorem}
\label{thm.spechtmod.zC}Let $\lambda$ be a partition of $n$. Let $T$ be any
$n$-tableau of shape $\lambda$. Let $C$ be a conjugacy class of $S_{n}$. Let
$\mathbf{z}_{C}$ denote the corresponding conjugacy class sum $\sum_{c\in
C}c\in\mathcal{A}$ (as in Definition \ref{def.groups.conj} \textbf{(c)}).
Then,%
\[
\mathbf{z}_{C}\mathbf{v}=\left(  \sum_{\substack{\left(  c,r\right)
\in\mathcal{C}\left(  T\right)  \times\mathcal{R}\left(  T\right)  ;\\cr\in
C}}\left(  -1\right)  ^{c}\right)  \mathbf{v}\ \ \ \ \ \ \ \ \ \ \text{for
each }\mathbf{v}\in\mathcal{S}^{\lambda}.
\]

\end{theorem}

\begin{proof}
Let $\mathbf{v}\in\mathcal{S}^{\lambda}$.

Theorem \ref{thm.center.conjsums} (applied to $G=S_{n}$) shows that the center
$Z\left(  \mathbf{k}\left[  S_{n}\right]  \right)  $ of the group algebra
$\mathbf{k}\left[  S_{n}\right]  $ is the $\mathbf{k}$-linear span of the
conjugacy class sums. Hence, in particular, the conjugacy class sums belong to
$Z\left(  \mathbf{k}\left[  S_{n}\right]  \right)  $. Thus, $\mathbf{z}_{C}$
belongs to $Z\left(  \mathbf{k}\left[  S_{n}\right]  \right)  $ (since
$\mathbf{z}_{C}$ is a conjugacy class sum of $S_{n}$).

Recall that $T$ is an $n$-tableau of shape $\lambda$, that is, an $n$-tableau
of shape $Y\left(  \lambda\right)  $. Thus,
(\ref{eq.def.specht.ET.defs.SD=AET}) (applied to $D=Y\left(  \lambda\right)
$) yields $\mathcal{S}^{Y\left(  \lambda\right)  }\cong\mathcal{A}%
\mathbf{E}_{T}$ as left $\mathcal{A}$-modules. In other words, $\mathcal{S}%
^{\lambda}\cong\mathcal{A}\mathbf{E}_{T}$ as left $\mathcal{A}$-modules (since
$\mathcal{S}^{\lambda}=\mathcal{S}^{Y\left(  \lambda\right)  }$). In other
words, there exists a left $\mathcal{A}$-module isomorphism $f:\mathcal{S}%
^{\lambda}\rightarrow\mathcal{A}\mathbf{E}_{T}$. Consider this $f$. Thus,
$f\left(  \mathbf{v}\right)  \in\mathcal{A}\mathbf{E}_{T}$. In other words,
\[
f\left(  \mathbf{v}\right)  =\mathbf{aE}_{T}\ \ \ \ \ \ \ \ \ \ \text{for some
}\mathbf{a}\in\mathcal{A}.
\]
Consider this $\mathbf{a}$.

Now, $f$ is a left $\mathcal{A}$-module morphism. Hence,%
\begin{align}
f\left(  \mathbf{z}_{C}\mathbf{v}\right)   &  =\mathbf{z}_{C}%
\underbrace{f\left(  \mathbf{v}\right)  }_{=\mathbf{aE}_{T}}=\mathbf{z}%
_{C}\mathbf{a}\underbrace{\mathbf{E}_{T}}_{\substack{=\nabla
_{\operatorname*{Col}T}^{-}\nabla_{\operatorname*{Row}T}\\\text{(by the
definition of }\mathbf{E}_{T}\text{)}}}=\mathbf{z}_{C}\mathbf{a}%
\nabla_{\operatorname*{Col}T}^{-}\nabla_{\operatorname*{Row}T}\nonumber\\
&  =\underbrace{\mathbf{z}_{C}\left(  \mathbf{a}\nabla_{\operatorname*{Col}%
T}^{-}\right)  }_{\substack{=\left(  \mathbf{a}\nabla_{\operatorname*{Col}%
T}^{-}\right)  \mathbf{z}_{C}\\\text{(since }\mathbf{z}_{C}\text{ belongs
to}\\\text{the center }Z\left(  \mathbf{k}\left[  S_{n}\right]  \right)
\text{)}}}\nabla_{\operatorname*{Row}T}=\left(  \mathbf{a}\nabla
_{\operatorname*{Col}T}^{-}\right)  \mathbf{z}_{C}\nabla_{\operatorname*{Row}%
T}\nonumber\\
&  =\mathbf{a}\nabla_{\operatorname*{Col}T}^{-}\underbrace{\mathbf{z}_{C}%
}_{\substack{=\sum_{c\in C}c\\=\sum_{w\in C}w}}\nabla_{\operatorname*{Row}%
T}=\mathbf{a}\nabla_{\operatorname*{Col}T}^{-}\left(  \sum_{w\in C}w\right)
\nabla_{\operatorname*{Row}T}\nonumber\\
&  =\sum_{w\in C}\mathbf{a}\nabla_{\operatorname*{Col}T}^{-}w\nabla
_{\operatorname*{Row}T}. \label{pf.thm.spechtmod.zC.2}%
\end{align}

Define a subset $U$ of $S_{n}$ by
\[
U:=\left\{  cr\ \mid\ \left(  c,r\right)  \in\mathcal{C}\left(  T\right)
\times\mathcal{R}\left(  T\right)  \right\}  .
\]
Then, obviously, $cr\in U$ for each $\left(  c,r\right)  \in\mathcal{C}\left(
T\right)  \times\mathcal{R}\left(  T\right)  $. Hence, the map%
\begin{align*}
\phi:\mathcal{C}\left(  T\right)  \times\mathcal{R}\left(  T\right)   &
\rightarrow U,\\
\left(  c,r\right)   &  \mapsto cr
\end{align*}
is well-defined. Consider this map $\phi$. We claim the following:

\begin{statement}
\textit{Claim 1:} The map $\phi$ is a bijection.
\end{statement}

\begin{proof}
[Proof of Claim 1.]The image of $\phi$ is%
\begin{align*}
\phi\left(  \mathcal{C}\left(  T\right)  \times\mathcal{R}\left(  T\right)
\right)   &  =\left\{  \underbrace{\phi\left(  c,r\right)  }%
_{\substack{=cr\\\text{(by the definition of }\phi\text{)}}}\ \mid\ \left(
c,r\right)  \in\mathcal{C}\left(  T\right)  \times\mathcal{R}\left(  T\right)
\right\} \\
&  =\left\{  cr\ \mid\ \left(  c,r\right)  \in\mathcal{C}\left(  T\right)
\times\mathcal{R}\left(  T\right)  \right\}  =U.
\end{align*}
Hence, the map $\phi$ is surjective.

Let us next show that $\phi$ is injective. Indeed, let $\left(  c_{1}%
,r_{1}\right)  $ and $\left(  c_{2},r_{2}\right)  $ be two elements of
$\mathcal{C}\left(  T\right)  \times\mathcal{R}\left(  T\right)  $ satisfying
$\phi\left(  c_{1},r_{1}\right)  =\phi\left(  c_{2},r_{2}\right)  $. We shall
show that $\left(  c_{1},r_{1}\right)  =\left(  c_{2},r_{2}\right)  $.

The definition of $\phi$ shows that $\phi\left(  c_{1},r_{1}\right)
=c_{1}r_{1}$ and $\phi\left(  c_{2},r_{2}\right)  =c_{2}r_{2}$, so that
$c_{1}r_{1}=\phi\left(  c_{1},r_{1}\right)  =\phi\left(  c_{2},r_{2}\right)
=c_{2}r_{2}$. Hence, $c_{2}^{-1}\underbrace{c_{1}r_{1}}_{=c_{2}r_{2}}%
r_{2}^{-1}=\underbrace{c_{2}^{-1}c_{2}}_{=\operatorname*{id}}\underbrace{r_{2}%
r_{2}^{-1}}_{=\operatorname*{id}}=\operatorname*{id}$. Moreover, $c_{1}%
\in\mathcal{C}\left(  T\right)  $ (since $\left(  c_{1},r_{1}\right)
\in\mathcal{C}\left(  T\right)  \times\mathcal{R}\left(  T\right)  $) and
$c_{2}\in\mathcal{C}\left(  T\right)  $ (similarly), so that $c_{2}^{-1}%
c_{1}\in\mathcal{C}\left(  T\right)  $ (since $\mathcal{C}\left(  T\right)  $
is a group). Similarly, $r_{1}r_{2}^{-1}\in\mathcal{R}\left(  T\right)  $.
Combining these two facts, we obtain $\left(  c_{2}^{-1}c_{1},r_{1}r_{2}%
^{-1}\right)  \in\mathcal{C}\left(  T\right)  \times\mathcal{R}\left(
T\right)  $.

But Proposition \ref{prop.tableau.RnC} \textbf{(b)} (applied to $D=Y\left(
\lambda\right)  $) yields that the only pair $\left(  c,r\right)
\in\mathcal{C}\left(  T\right)  \times\mathcal{R}\left(  T\right)  $
satisfying $cr=\operatorname*{id}$ is the pair $\left(  \operatorname*{id}%
,\operatorname*{id}\right)  $. Hence, any pair $\left(  c,r\right)
\in\mathcal{C}\left(  T\right)  \times\mathcal{R}\left(  T\right)  $
satisfying $cr=\operatorname*{id}$ must equal $\left(  \operatorname*{id}%
,\operatorname*{id}\right)  $. Applying this to $\left(  c,r\right)  =\left(
c_{2}^{-1}c_{1},r_{1}r_{2}^{-1}\right)  $, we obtain that $\left(  c_{2}%
^{-1}c_{1},r_{1}r_{2}^{-1}\right)  $ equals $\left(  \operatorname*{id}%
,\operatorname*{id}\right)  $ (since $\left(  c_{2}^{-1}c_{1},r_{1}r_{2}%
^{-1}\right)  \in\mathcal{C}\left(  T\right)  \times\mathcal{R}\left(
T\right)  $ and $c_{2}^{-1}c_{1}r_{1}r_{2}^{-1}=\operatorname*{id}$). In other
words, $c_{2}^{-1}c_{1}=\operatorname*{id}$ and $r_{1}r_{2}^{-1}%
=\operatorname*{id}$. Hence, $c_{1}=c_{2}$ (since $c_{2}^{-1}c_{1}%
=\operatorname*{id}$) and $r_{1}=r_{2}$ (since $r_{1}r_{2}^{-1}%
=\operatorname*{id}$). Therefore, $\left(  c_{1},r_{1}\right)  =\left(
c_{2},r_{2}\right)  $.

Forget that we fixed $\left(  c_{1},r_{1}\right)  $ and $\left(  c_{2}%
,r_{2}\right)  $. We have thus proved that if $\left(  c_{1},r_{1}\right)  $
and $\left(  c_{2},r_{2}\right)  $ are two elements of $\mathcal{C}\left(
T\right)  \times\mathcal{R}\left(  T\right)  $ satisfying $\phi\left(
c_{1},r_{1}\right)  =\phi\left(  c_{2},r_{2}\right)  $, then $\left(
c_{1},r_{1}\right)  =\left(  c_{2},r_{2}\right)  $. In other words, the map
$\phi$ is injective.

Now, we know that the map $\phi$ is both injective and surjective. Hence,
$\phi$ is bijective, i.e., is a bijection. This proves Claim 1.
\end{proof}

\begin{statement}
\textit{Claim 2:} Let $w\in S_{n}$ be such that $w\notin U$. Then,
$\nabla_{\operatorname*{Col}T}^{-}w\nabla_{\operatorname*{Row}T}=0$.
\end{statement}

\begin{proof}
[Proof of Claim 2.]We have%
\[
w\notin U=\left\{  cr\ \mid\ \left(  c,r\right)  \in\mathcal{C}\left(
T\right)  \times\mathcal{R}\left(  T\right)  \right\}  =\left\{
cr\ \mid\ c\in\mathcal{C}\left(  T\right)  \text{ and }r\in\mathcal{R}\left(
T\right)  \right\}  .
\]
In other words, there exist no two permutations $c\in\mathcal{C}\left(
T\right)  $ and $r\in\mathcal{R}\left(  T\right)  $ such that $w=cr$.

Using Theorem \ref{thm.specht.sandw1} \textbf{(b)}, we can thus see that there
exist two distinct integers that are horizontally aligned in $wT$ and
simultaneously vertically aligned in $T$ (since otherwise, Theorem
\ref{thm.specht.sandw1} \textbf{(b)} would yield that there exist two
permutations $c\in\mathcal{C}\left(  T\right)  $ and $r\in\mathcal{R}\left(
T\right)  $ such that $w=cr$ and $\nabla_{\operatorname*{Col}T}^{-}%
w\nabla_{\operatorname*{Row}T}=\left(  -1\right)  ^{c}\nabla
_{\operatorname*{Col}T}^{-}\nabla_{\operatorname*{Row}T}$; but this would
contradict the fact that there exist no two permutations $c\in\mathcal{C}%
\left(  T\right)  $ and $r\in\mathcal{R}\left(  T\right)  $ such that $w=cr$).
Hence, Theorem \ref{thm.specht.sandw1} \textbf{(a)} yields $\nabla
_{\operatorname*{Col}T}^{-}w\nabla_{\operatorname*{Row}T}=0$. This proves
Claim 2.
\end{proof}

\begin{statement}
\textit{Claim 3:} Let $\left(  c,r\right)  \in\mathcal{C}\left(  T\right)
\times\mathcal{R}\left(  T\right)  $. Then, $\nabla_{\operatorname*{Col}T}%
^{-}cr\nabla_{\operatorname*{Row}T}=\left(  -1\right)  ^{c}\mathbf{E}_{T}$.
\end{statement}

\begin{proof}
[Proof of Claim 3.]From $\left(  c,r\right)  \in\mathcal{C}\left(  T\right)
\times\mathcal{R}\left(  T\right)  $, we obtain $c\in\mathcal{C}\left(
T\right)  $ and $r\in\mathcal{R}\left(  T\right)  $. Thus, Proposition
\ref{prop.symmetrizers.fix} \textbf{(a)} (applied to $r$ instead of $w$)
yields $r\nabla_{\operatorname*{Row}T}=\nabla_{\operatorname*{Row}T}%
r=\nabla_{\operatorname*{Row}T}$. Also, Proposition
\ref{prop.symmetrizers.fix} \textbf{(b)} (applied to $c$ instead of $w$)
yields $c\nabla_{\operatorname*{Col}T}^{-}=\nabla_{\operatorname*{Col}T}%
^{-}c=\left(  -1\right)  ^{c}\nabla_{\operatorname*{Col}T}^{-}$ (since
$c\in\mathcal{C}\left(  T\right)  $). Using these facts, we find%
\[
\underbrace{\nabla_{\operatorname*{Col}T}^{-}c}_{=\left(  -1\right)
^{c}\nabla_{\operatorname*{Col}T}^{-}}\ \ \underbrace{r\nabla
_{\operatorname*{Row}T}}_{=\nabla_{\operatorname*{Row}T}}=\left(  -1\right)
^{c}\underbrace{\nabla_{\operatorname*{Col}T}^{-}\nabla_{\operatorname*{Row}%
T}}_{\substack{=\mathbf{E}_{T}\\\text{(by the definition of }\mathbf{E}%
_{T}\text{)}}}=\left(  -1\right)  ^{c}\mathbf{E}_{T}.
\]
This proves Claim 3.
\end{proof}

Now, each element $w\in C$ satisfies either $w\in U$ or $w\notin U$ (but not
both). Hence,%
\begin{align*}
&  \sum_{w\in C}\mathbf{a}\nabla_{\operatorname*{Col}T}^{-}w\nabla
_{\operatorname*{Row}T}\\
&  =\underbrace{\sum_{\substack{w\in C;\\w\in U}}}_{=\sum
\limits_{\substack{w\in U;\\w\in C}}}\mathbf{a}\nabla_{\operatorname*{Col}%
T}^{-}w\nabla_{\operatorname*{Row}T}+\sum_{\substack{w\in C;\\w\notin
U}}\mathbf{a}\underbrace{\nabla_{\operatorname*{Col}T}^{-}w\nabla
_{\operatorname*{Row}T}}_{\substack{=0\\\text{(by Claim 2)}}}\\
&  =\sum\limits_{\substack{w\in U;\\w\in C}}\mathbf{a}\nabla
_{\operatorname*{Col}T}^{-}w\nabla_{\operatorname*{Row}T}+\underbrace{\sum
_{\substack{w\in C;\\w\notin U}}\mathbf{a}0}_{=0}=\sum\limits_{\substack{w\in
U;\\w\in C}}\mathbf{a}\nabla_{\operatorname*{Col}T}^{-}w\nabla
_{\operatorname*{Row}T}\\
&  =\sum_{\substack{\left(  c,r\right)  \in\mathcal{C}\left(  T\right)
\times\mathcal{R}\left(  T\right)  ;\\\phi\left(  c,r\right)  \in
C}}\mathbf{a}\nabla_{\operatorname*{Col}T}^{-}\left(  \phi\left(  c,r\right)
\right)  \nabla_{\operatorname*{Row}T}\\
&  \ \ \ \ \ \ \ \ \ \ \ \ \ \ \ \ \ \ \ \ \left(
\begin{array}
[c]{c}%
\text{here, we have substituted }\phi\left(  c,r\right)  \text{ for }w\text{
in the sum,}\\
\text{since the map }\phi:\mathcal{C}\left(  T\right)  \times\mathcal{R}%
\left(  T\right)  \rightarrow U\text{ is a bijection}%
\end{array}
\right) \\
&  =\sum_{\substack{\left(  c,r\right)  \in\mathcal{C}\left(  T\right)
\times\mathcal{R}\left(  T\right)  ;\\cr\in C}}\mathbf{a}\underbrace{\nabla
_{\operatorname*{Col}T}^{-}cr\nabla_{\operatorname*{Row}T}}%
_{\substack{=\left(  -1\right)  ^{c}\mathbf{E}_{T}\\\text{(by Claim 3)}%
}}\ \ \ \ \ \ \ \ \ \ \left(
\begin{array}
[c]{c}%
\text{since }\phi\left(  c,r\right)  =cr\\
\text{for each }\left(  c,r\right)  \in\mathcal{C}\left(  T\right)
\times\mathcal{R}\left(  T\right) \\
\text{(by the definition of }\phi\text{)}%
\end{array}
\right) \\
&  =\sum_{\substack{\left(  c,r\right)  \in\mathcal{C}\left(  T\right)
\times\mathcal{R}\left(  T\right)  ;\\cr\in C}}\mathbf{a}\left(  -1\right)
^{c}\mathbf{E}_{T}=\left(  \sum_{\substack{\left(  c,r\right)  \in
\mathcal{C}\left(  T\right)  \times\mathcal{R}\left(  T\right)  ;\\cr\in
C}}\left(  -1\right)  ^{c}\right)  \underbrace{\mathbf{aE}_{T}}%
_{\substack{=f\left(  \mathbf{v}\right)  \\\text{(since }f\left(
\mathbf{v}\right)  =\mathbf{aE}_{T}\text{)}}}\\
&  =\left(  \sum_{\substack{\left(  c,r\right)  \in\mathcal{C}\left(
T\right)  \times\mathcal{R}\left(  T\right)  ;\\cr\in C}}\left(  -1\right)
^{c}\right)  f\left(  \mathbf{v}\right)  =f\left(  \left(  \sum
_{\substack{\left(  c,r\right)  \in\mathcal{C}\left(  T\right)  \times
\mathcal{R}\left(  T\right)  ;\\cr\in C}}\left(  -1\right)  ^{c}\right)
\mathbf{v}\right)
\end{align*}
(since the map $f$ is $\mathbf{k}$-linear). Thus, (\ref{pf.thm.spechtmod.zC.2}%
) becomes%
\[
f\left(  \mathbf{z}_{C}\mathbf{v}\right)  =\sum_{w\in C}\mathbf{a}%
\nabla_{\operatorname*{Col}T}^{-}w\nabla_{\operatorname*{Row}T}=f\left(
\left(  \sum_{\substack{\left(  c,r\right)  \in\mathcal{C}\left(  T\right)
\times\mathcal{R}\left(  T\right)  ;\\cr\in C}}\left(  -1\right)  ^{c}\right)
\mathbf{v}\right)  .
\]
Since the map $f$ is injective (because $f$ is an isomorphism), we can
\textquotedblleft un-apply\textquotedblright\ $f$ to this equality, and thus
obtain%
\[
\mathbf{z}_{C}\mathbf{v}=\left(  \sum_{\substack{\left(  c,r\right)
\in\mathcal{C}\left(  T\right)  \times\mathcal{R}\left(  T\right)  ;\\cr\in
C}}\left(  -1\right)  ^{c}\right)  \mathbf{v}.
\]
This proves Theorem \ref{thm.spechtmod.zC}.
\end{proof}

\subsubsection{The sum of all transpositions}

Recall that all transpositions in $S_{n}$ form a single conjugacy class. Thus,
the sum of all transpositions in $S_{n}$ is a conjugacy class sum. We can
therefore apply Theorem \ref{thm.spechtmod.zC} to see how this sum acts on a
Specht module $\mathcal{S}^{\lambda}$. What is interesting is that the result
can be simplified, to yield a much nicer formula for the result.

To state this formula, we introduce a notation:

\begin{definition}
\label{def.partitions.nlam}Let $\lambda=\left(  \lambda_{1},\lambda_{2}%
,\ldots,\lambda_{k}\right)  $ be a partition. Then, we set
\[
\operatorname*{n}\left(  \lambda\right)  :=\sum_{i=1}^{k}\underbrace{\left(
i-1\right)  }_{\in\mathbb{N}}\lambda_{i}\in\mathbb{N}.
\]

\end{definition}

For example, if $\lambda=\left(  4,2,2\right)  $, then $\operatorname*{n}%
\left(  \lambda\right)  =\left(  1-1\right)  \cdot4+\left(  2-1\right)
\cdot2+\left(  3-1\right)  \cdot2=6$.

Note that many authors (e.g., Macdonald in \cite[\S I.1, (1.5)]{Macdon95})
write $n\left(  \lambda\right)  $ for $\operatorname*{n}\left(  \lambda
\right)  $, but I am using the Roman-style $\operatorname*{n}$ to avoid
confusion with the number $n$.

The following three propositions give alternative descriptions of
$\operatorname*{n}\left(  \lambda\right)  $:

\begin{proposition}
\label{prop.partitions.nlam-sumij}Let $\lambda$ be a partition. Then,%
\[
\operatorname*{n}\left(  \lambda\right)  =\sum_{\left(  i,j\right)  \in
Y\left(  \lambda\right)  }\left(  i-1\right)  .
\]

\end{proposition}

\begin{proposition}
\label{prop.partitions.nlamt}Let $\lambda$ be a partition. Write the conjugate
$\lambda^{t}$ of $\lambda$ as $\lambda^{t}=\left(  \lambda_{1}^{t},\lambda
_{2}^{t},\ldots,\lambda_{m}^{t}\right)  $. Then,%
\[
\operatorname*{n}\left(  \lambda\right)  =\sum_{i=1}^{m}\dbinom{\lambda
_{i}^{t}}{2}.
\]

\end{proposition}

\begin{proposition}
\label{prop.partitions.nlam-counts-t}Let $\lambda$ be a partition of $n$. Let
$T$ be an $n$-tableau of shape $\lambda$. Then, $\operatorname*{n}\left(
\lambda\right)  $ is the number of transpositions in the subgroup
$\mathcal{C}\left(  T\right)  $ of $S_{n}$.
\end{proposition}

These three propositions are all elementary and quite easy to prove; detailed
proofs can be found in the appendix (Section
\ref{sec.details.specht.conj-class-action}).

We are now ready to state how the sum of all transpositions acts on a Specht
module $\mathcal{S}^{\lambda}$:

\begin{theorem}
\label{thm.spechtmod.sum-transp}Let $\lambda$ be a partition of $n$. Let
$\mathbf{t}\in\mathbf{k}\left[  S_{n}\right]  $ be the sum of all
transpositions in $S_{n}$. Then,%
\[
\mathbf{tv}=\left(  \operatorname*{n}\left(  \lambda^{t}\right)
-\operatorname*{n}\left(  \lambda\right)  \right)  \mathbf{v}%
\ \ \ \ \ \ \ \ \ \ \text{for each }\mathbf{v}\in\mathcal{S}^{\lambda}.
\]

\end{theorem}

Before we can prove this, we must establish three combinatorial lemmas:

\begin{lemma}
\label{lem.spechtmod.sum-transp.1}Let $T$ be an $n$-tableau (of any shape
$D$). Let $c\in\mathcal{C}\left(  T\right)  $ and $r\in\mathcal{R}\left(
T\right)  $. Let $k\in\left[  n\right]  $. Then: \medskip

\textbf{(a)} If $c\left(  k\right)  =r\left(  k\right)  $, then $c\left(
k\right)  =r\left(  k\right)  =k$. \medskip

\textbf{(b)} If $\left(  cr\right)  \left(  k\right)  =k$, then $c\left(
k\right)  =r\left(  k\right)  =k$.
\end{lemma}

\begin{proof}
\textbf{(a)} Assume that $c\left(  k\right)  =r\left(  k\right)  $. The
tableau $T$ is an $n$-tableau, and thus contains each of the integers
$1,2,\ldots,n$ exactly once. Hence, there is a unique cell $p$ of $T$ that
contains the number $k$, and there is a unique cell $q$ of $T$ that contains
the number $c\left(  k\right)  =r\left(  k\right)  $. Consider these two cells
$p$ and $q$. Thus, $T\left(  p\right)  =k$ (since the cell $p$ of $T$ contains
the number $k$) and similarly $T\left(  q\right)  =c\left(  k\right)
=r\left(  k\right)  $.

From $r\in\mathcal{R}\left(  T\right)  $, we see that the permutation $r$ is
horizontal for $T$ (by the definition of $\mathcal{R}\left(  T\right)  $).
Hence, the number $r\left(  k\right)  $ lies in the same row of $T$ as $k$
does (by the definition of \textquotedblleft horizontal\textquotedblright). In
other words, the cell of $T$ that contains the number $k$ lies in the same row
as the cell of $T$ that contains the number $r\left(  k\right)  $. Since the
former cell is $p$, and since the latter cell is $q$, we can rewrite this as
follows: The cell $p$ lies in the same row as the cell $q$.

Similarly, from $c\in\mathcal{C}\left(  T\right)  $, we conclude that the cell
$p$ lies in the same column as the cell $q$.

However, if two cells lie in the same row and in the same column, then they
are simply the same cell. Thus, $p$ and $q$ are the same cell (since $p$ lies
in the same row and in the same column as $q$). In other words, $p=q$. Hence,
$T\left(  p\right)  =T\left(  q\right)  =r\left(  k\right)  $, so that
$r\left(  k\right)  =T\left(  p\right)  =k$. Combining this with $c\left(
k\right)  =r\left(  k\right)  $, we obtain $c\left(  k\right)  =r\left(
k\right)  =k$. Thus, Lemma \ref{lem.spechtmod.sum-transp.1} \textbf{(a)} is
proved. \medskip

\textbf{(b)} Assume that $\left(  cr\right)  \left(  k\right)  =k$. Thus,
$c\left(  r\left(  k\right)  \right)  =\left(  cr\right)  \left(  k\right)
=k$. Therefore, $r\left(  k\right)  =c^{-1}\left(  k\right)  $, so that
$c^{-1}\left(  k\right)  =r\left(  k\right)  $. But $c\in\mathcal{C}\left(
T\right)  $ and thus $c^{-1}\in\mathcal{C}\left(  T\right)  $ (since
$\mathcal{C}\left(  T\right)  $ is a group). Hence, Lemma
\ref{lem.spechtmod.sum-transp.1} \textbf{(a)} (applied to $c^{-1}$ instead of
$c$) shows that $c^{-1}\left(  k\right)  =r\left(  k\right)  =k$. Hence,
$k=c\left(  k\right)  $, so that $c\left(  k\right)  =k=r\left(  k\right)
=k$. This proves Lemma \ref{lem.spechtmod.sum-transp.1} \textbf{(b)}.
\end{proof}

\begin{lemma}
\label{lem.spechtmod.sum-transp.2}Let $C$ be the set of all transpositions in
$S_{n}$. Let $T$ be an $n$-tableau (of any shape $D$). Let $\left(
c,r\right)  \in\mathcal{C}\left(  T\right)  \times\mathcal{R}\left(  T\right)
$ be such that $cr\in C$. Then: \medskip

\textbf{(a)} If $\left(  -1\right)  ^{c}=1$, then $c=\operatorname*{id}$ and
$r\in C$. \medskip

\textbf{(b)} If $\left(  -1\right)  ^{c}=-1$, then $r=\operatorname*{id}$ and
$c\in C$.
\end{lemma}

\begin{proof}
From $\left(  c,r\right)  \in\mathcal{C}\left(  T\right)  \times
\mathcal{R}\left(  T\right)  $, we obtain $c\in\mathcal{C}\left(  T\right)  $
and $r\in\mathcal{R}\left(  T\right)  $.

We have $cr\in C$ (by assumption). In other words, $cr$ is a transposition
(since $C$ is the set of all transpositions in $S_{n}$). In other words,
$cr=t_{u,v}$ for some two distinct numbers $u,v\in\left[  n\right]  $.
Consider these $u,v$. Define a subset $X$ of $\left[  n\right]  $ by
$X:=\left\{  u,v\right\}  $. This subset $X$ has size $2$ (since $u$ and $v$
are distinct). Thus, its symmetric group $S_{X}$ consists of two permutations
(namely, $\operatorname*{id}\nolimits_{X}$ and $t_{u,v}$). In other words,
$\left\vert S_{X}\right\vert =2$.

Define the set
\[
S_{n,X}:=\left\{  w\in S_{n}\ \mid\ w\left(  i\right)  =i\text{ for all
}i\notin X\right\}
\]
(where \textquotedblleft for all $i\notin X$\textquotedblright\ is shorthand
for \textquotedblleft for all $i\in\left[  n\right]  \setminus X$%
\textquotedblright). Proposition \ref{prop.intX.basics} \textbf{(a)} shows
that this set $S_{n,X}$ is a subgroup of $S_{n}$, and is isomorphic to $S_{X}%
$. Hence, $\left\vert S_{n,X}\right\vert =\left\vert S_{X}\right\vert =2$. But
the two permutations $\operatorname*{id}\in S_{n}$ and $t_{u,v}\in S_{n}$
clearly belong to $S_{n,X}$ (since $u,v\in X$), and are distinct (since
$t_{u,v}\left(  u\right)  =v\neq u=\operatorname*{id}\left(  u\right)  $ and
thus $t_{u,v}\neq\operatorname*{id}$). This easily yields that $S_{n,X}%
=\left\{  \operatorname*{id},t_{u,v}\right\}  $%
\ \ \ \ \footnote{\textit{Proof.} We know that $\operatorname*{id}$ and
$t_{u,v}$ are two distinct elements of $S_{n,X}$ (since $\operatorname*{id}$
and $t_{u,v}$ belong to $S_{n,X}$). Moreover, the set $S_{n,X}$ has no further
elements beyond $\operatorname*{id}$ and $t_{u,v}$ (since $S_{n,X}$ has only
two elements (because $\left\vert S_{n,X}\right\vert =2$)). Thus, the set
$S_{n,X}$ has the elements $\operatorname*{id}$ and $t_{u,v}$ and no others.
In other words, $S_{n,X}=\left\{  \operatorname*{id},t_{u,v}\right\}  $.}.

Note that $\left(  -1\right)  ^{t_{u,v}}=-1$ (since any transposition has sign
$-1$).

We now claim the following:

\begin{statement}
\textit{Claim 1:} We have $c\in S_{n,X}$.
\end{statement}

\begin{proof}
[Proof of Claim 1.]Let $i\in\left[  n\right]  \setminus X$. Then, $i\notin
X=\left\{  u,v\right\}  $, so that $t_{u,v}\left(  i\right)  =i$. But
$cr=t_{u,v}$, so that $\left(  cr\right)  \left(  i\right)  =t_{u,v}\left(
i\right)  =i$. Hence, Lemma \ref{lem.spechtmod.sum-transp.1} \textbf{(b)}
(applied to $k=i$) yields $c\left(  i\right)  =r\left(  i\right)  =i$. In
particular, $c\left(  i\right)  =i$.

Forget that we fixed $i$. We thus have shown that $c\left(  i\right)  =i$ for
all $i\in\left[  n\right]  \setminus X$. In other words, $c\left(  i\right)
=i$ for all $i\notin X$. In other words, $c\in S_{n,X}$ (by the definition of
$S_{n,X}$). This proves Claim 1.
\end{proof}

Now, we can easily prove both parts of the lemma: \medskip

\textbf{(a)} Assume that $\left(  -1\right)  ^{c}=1$. Hence, $c\neq t_{u,v}$
(since $\left(  -1\right)  ^{c}=1\neq-1=\left(  -1\right)  ^{t_{u,v}}$). But
Claim 1 yields $c\in S_{n,X}=\left\{  \operatorname*{id},t_{u,v}\right\}  $.
Combining this with $c\neq t_{u,v}$, we find $c\in\left\{  \operatorname*{id}%
,t_{u,v}\right\}  \setminus\left\{  t_{u,v}\right\}  \subseteq\left\{
\operatorname*{id}\right\}  $, so that $c=\operatorname*{id}$. Thus,
$cr=\operatorname*{id}r=r$, so that $r=cr\in C$. Thus, Lemma
\ref{lem.spechtmod.sum-transp.2} \textbf{(a)} is proved. \medskip

\textbf{(b)} Assume that $\left(  -1\right)  ^{c}=-1$. Hence, $c\neq
\operatorname*{id}$ (since $\left(  -1\right)  ^{c}=-1\neq1=\left(  -1\right)
^{\operatorname*{id}}$). But Claim 1 yields $c\in S_{n,X}=\left\{
\operatorname*{id},t_{u,v}\right\}  $. Combining this with $c\neq
\operatorname*{id}$, we find $c\in\left\{  \operatorname*{id},t_{u,v}\right\}
\setminus\left\{  \operatorname*{id}\right\}  \subseteq\left\{  t_{u,v}%
\right\}  $, so that $c=t_{u,v}=cr$. Cancelling $c$ from this equality, we
find $\operatorname*{id}=r$ (since $c$ is invertible in $S_{n}$), so that
$r=\operatorname*{id}$. Moreover, $c=cr\in C$. Thus, Lemma
\ref{lem.spechtmod.sum-transp.2} \textbf{(b)} is proved.
\end{proof}

The above lemma lets us show the following:

\begin{lemma}
\label{lem.spechtmod.sum-transp.3}Let $C$ be the set of all transpositions in
$S_{n}$. Let $\lambda$ be a partition of $n$. Let $T$ be any $n$-tableau of
shape $\lambda$. Then,%
\[
\sum_{\substack{\left(  c,r\right)  \in\mathcal{C}\left(  T\right)
\times\mathcal{R}\left(  T\right)  ;\\cr\in C}}\left(  -1\right)
^{c}=\operatorname*{n}\left(  \lambda^{t}\right)  -\operatorname*{n}\left(
\lambda\right)  .
\]

\end{lemma}

\begin{proof}
Define the two sets
\[
U:=\left\{  \left(  c,r\right)  \in\mathcal{C}\left(  T\right)  \times
\mathcal{R}\left(  T\right)  \ \mid\ cr\in C\text{ and }\left(  -1\right)
^{c}=1\right\}
\]
and%
\[
V:=\left\{  \left(  c,r\right)  \in\mathcal{C}\left(  T\right)  \times
\mathcal{R}\left(  T\right)  \ \mid\ cr\in C\text{ and }\left(  -1\right)
^{c}=-1\right\}  .
\]

Then, it is easy to see that $\left(  c,\operatorname*{id}\right)  \in V$ for
each $c\in C\cap\mathcal{C}\left(  T\right)  $%
\ \ \ \ \footnote{\textit{Proof.} Let $c\in C\cap\mathcal{C}\left(  T\right)
$. We must show that $\left(  c,\operatorname*{id}\right)  \in V$. According
to the definition of $V$, this means that we must prove that $\left(
c,\operatorname*{id}\right)  \in\mathcal{C}\left(  T\right)  \times
\mathcal{R}\left(  T\right)  $ and $c\operatorname*{id}\in C$ and $\left(
-1\right)  ^{c}=-1$.
\par
But this is pretty easy: We have $\operatorname*{id}\in\mathcal{R}\left(
T\right)  $ (since $\mathcal{R}\left(  T\right)  $ is a subgroup of $S_{n}$).
Combining $c\in C\cap\mathcal{C}\left(  T\right)  \subseteq\mathcal{C}\left(
T\right)  $ with $\operatorname*{id}\in\mathcal{R}\left(  T\right)  $, we
obtain $\left(  c,\operatorname*{id}\right)  \in\mathcal{C}\left(  T\right)
\times\mathcal{R}\left(  T\right)  $. Furthermore, $c\operatorname*{id}=c\in
C\cap\mathcal{C}\left(  T\right)  \subseteq C$. Finally, $c\in C$, so that $c$
is a transposition (since $C$ is the set of all transpositions in $S_{n}$).
Hence, $\left(  -1\right)  ^{c}=-1$ (since any transposition has sign $-1$).
\par
Altogether, we have shown that $\left(  c,\operatorname*{id}\right)
\in\mathcal{C}\left(  T\right)  \times\mathcal{R}\left(  T\right)  $ and
$c\operatorname*{id}\in C$ and $\left(  -1\right)  ^{c}=-1$. In other words,
$\left(  c,\operatorname*{id}\right)  \in V$ (by the definition of $V$).
Qed.}. Hence, the map%
\begin{align}
C\cap\mathcal{C}\left(  T\right)   &  \rightarrow V,\nonumber\\
c  &  \mapsto\left(  c,\operatorname*{id}\right)
\label{pf.lem.spechtmod.sum-transp.3.map1}%
\end{align}
is well-defined. This map is furthermore injective\footnote{since any $c\in
C\cap\mathcal{C}\left(  T\right)  $ can be reconstructed from the pair
$\left(  c,\operatorname*{id}\right)  $ (simply by taking the first entry of
the pair)} and surjective\footnote{\textit{Proof.} Let $v\in V$. We must show
that there exists some $c\in C\cap\mathcal{C}\left(  T\right)  $ such that
$v=\left(  c,\operatorname*{id}\right)  $.
\par
We have $v\in V$. In other words, $v$ is a pair $\left(  c,r\right)
\in\mathcal{C}\left(  T\right)  \times\mathcal{R}\left(  T\right)  $
satisfying $cr\in C$ and $\left(  -1\right)  ^{c}=-1$ (by the definition of
$V$). Consider this pair $\left(  c,r\right)  $. Then, Lemma
\ref{lem.spechtmod.sum-transp.2} \textbf{(b)} (applied to $D=Y\left(
\lambda\right)  $) yields that $r=\operatorname*{id}$ and $c\in C$. Also,
$c\in\mathcal{C}\left(  T\right)  $ (since $\left(  c,r\right)  \in
\mathcal{C}\left(  T\right)  \times\mathcal{R}\left(  T\right)  $). Combining
$c\in C$ with $c\in\mathcal{C}\left(  T\right)  $, we obtain $c\in
C\cap\mathcal{C}\left(  T\right)  $. Moreover, by the definition of $\left(
c,r\right)  $, we have $v=\left(  c,r\right)  =\left(  c,\operatorname*{id}%
\right)  $ (since $r=\operatorname*{id}$).
\par
Forget that we defined $c$. We thus have found a $c\in C\cap\mathcal{C}\left(
T\right)  $ such that $v=\left(  c,\operatorname*{id}\right)  $. Thus, there
exists some $c\in C\cap\mathcal{C}\left(  T\right)  $ such that $v=\left(
c,\operatorname*{id}\right)  $.
\par
Forget that we fixed $v$. We thus have shown that for each $v\in V$, there
exists some $c\in C\cap\mathcal{C}\left(  T\right)  $ such that $v=\left(
c,\operatorname*{id}\right)  $. In other words, the map
(\ref{pf.lem.spechtmod.sum-transp.3.map1}) is surjective.}. Hence, it is
bijective, i.e., a bijection. Thus,
\begin{equation}
\left\vert V\right\vert =\left\vert C\cap\mathcal{C}\left(  T\right)
\right\vert . \label{pf.lem.spechtmod.sum-transp.3.map1eq}%
\end{equation}

Furthermore, it is easy to see that $\left(  \operatorname*{id},r\right)  \in
U$ for each $r\in C\cap\mathcal{R}\left(  T\right)  $%
\ \ \ \ \footnote{\textit{Proof.} Let $r\in C\cap\mathcal{R}\left(  T\right)
$. We must show that $\left(  \operatorname*{id},r\right)  \in U$. According
to the definition of $U$, this means that we must prove that $\left(
\operatorname*{id},r\right)  \in\mathcal{C}\left(  T\right)  \times
\mathcal{R}\left(  T\right)  $ and $\operatorname*{id}r\in C$ and $\left(
-1\right)  ^{\operatorname*{id}}=1$.
\par
But this is pretty easy: We have $\operatorname*{id}\in\mathcal{C}\left(
T\right)  $ (since $\mathcal{C}\left(  T\right)  $ is a subgroup of $S_{n}$).
Combining this with $r\in C\cap\mathcal{R}\left(  T\right)  \subseteq
\mathcal{R}\left(  T\right)  $, we obtain $\left(  \operatorname*{id}%
,r\right)  \in\mathcal{C}\left(  T\right)  \times\mathcal{R}\left(  T\right)
$. Furthermore, $\operatorname*{id}r=r\in C\cap\mathcal{R}\left(  T\right)
\subseteq C$. Finally, $\left(  -1\right)  ^{\operatorname*{id}}=1$.
\par
Altogether, we have shown that $\left(  \operatorname*{id},r\right)
\in\mathcal{C}\left(  T\right)  \times\mathcal{R}\left(  T\right)  $ and
$\operatorname*{id}r\in C$ and $\left(  -1\right)  ^{\operatorname*{id}}=1$.
In other words, $\left(  \operatorname*{id},r\right)  \in U$ (by the
definition of $U$). Qed.}. Hence, the map%
\begin{align}
C\cap\mathcal{R}\left(  T\right)   &  \rightarrow U,\nonumber\\
r  &  \mapsto\left(  \operatorname*{id},r\right)
\label{pf.lem.spechtmod.sum-transp.3.map2}%
\end{align}
is well-defined. This map is furthermore injective\footnote{since any $r\in
C\cap\mathcal{R}\left(  T\right)  $ can be reconstructed from the pair
$\left(  \operatorname*{id},r\right)  $ (simply by taking the second entry of
the pair)} and surjective\footnote{\textit{Proof.} Let $u\in U$. We must show
that there exists some $r\in C\cap\mathcal{R}\left(  T\right)  $ such that
$u=\left(  \operatorname*{id},r\right)  $.
\par
We have $u\in U$. In other words, $u$ is a pair $\left(  c,r\right)
\in\mathcal{C}\left(  T\right)  \times\mathcal{R}\left(  T\right)  $
satisfying $cr\in C$ and $\left(  -1\right)  ^{c}=1$ (by the definition of
$U$). Consider this pair $\left(  c,r\right)  $. Then, Lemma
\ref{lem.spechtmod.sum-transp.2} \textbf{(a)} (applied to $D=Y\left(
\lambda\right)  $) yields that $c=\operatorname*{id}$ and $r\in C$. Also,
$r\in\mathcal{R}\left(  T\right)  $ (since $\left(  c,r\right)  \in
\mathcal{C}\left(  T\right)  \times\mathcal{R}\left(  T\right)  $). Combining
$r\in C$ with $r\in\mathcal{R}\left(  T\right)  $, we obtain $r\in
C\cap\mathcal{R}\left(  T\right)  $. Moreover, by the definition of $\left(
c,r\right)  $, we have $u=\left(  c,r\right)  =\left(  \operatorname*{id}%
,r\right)  $ (since $c=\operatorname*{id}$).
\par
Forget that we defined $r$. We thus have found an $r\in C\cap\mathcal{R}%
\left(  T\right)  $ such that $u=\left(  \operatorname*{id},r\right)  $. Thus,
there exists some $r\in C\cap\mathcal{R}\left(  T\right)  $ such that
$u=\left(  \operatorname*{id},r\right)  $.
\par
Forget that we fixed $u$. We thus have shown that for each $u\in U$, there
exists some $r\in C\cap\mathcal{R}\left(  T\right)  $ such that $u=\left(
\operatorname*{id},r\right)  $. In other words, the map
(\ref{pf.lem.spechtmod.sum-transp.3.map2}) is surjective.}. Hence, it is
bijective, i.e., a bijection. Thus,
\begin{equation}
\left\vert U\right\vert =\left\vert C\cap\mathcal{R}\left(  T\right)
\right\vert . \label{pf.lem.spechtmod.sum-transp.3.map2eq}%
\end{equation}

Furthermore, we have%
\begin{equation}
\left\vert C\cap\mathcal{C}\left(  T\right)  \right\vert =\operatorname*{n}%
\left(  \lambda\right)  \label{pf.lem.spechtmod.sum-transp.3.nlam}%
\end{equation}
\footnote{\textit{Proof:} Proposition \ref{prop.partitions.nlam-counts-t}
shows that $\operatorname*{n}\left(  \lambda\right)  $ is the number of
transpositions in $\mathcal{C}\left(  T\right)  $. In other words,%
\begin{equation}
\operatorname*{n}\left(  \lambda\right)  =\left(  \text{\# of all
transpositions in }\mathcal{C}\left(  T\right)  \right)  .
\label{pf.lem.spechtmod.sum-transp.3.nlam.pf.1}%
\end{equation}
\par
Recall that $C$ is the set of all transpositions in $S_{n}$. In other words,
$C=\left\{  \text{all transpositions in }S_{n}\right\}  $. Hence,%
\begin{align*}
C\cap\mathcal{C}\left(  T\right)   &  =\left\{  \text{all transpositions in
}S_{n}\right\}  \cap\mathcal{C}\left(  T\right) \\
&  =\left\{  \text{all transpositions in }S_{n}\text{ that belong to
}\mathcal{C}\left(  T\right)  \right\} \\
&  =\left\{  \text{all transpositions in }\mathcal{C}\left(  T\right)
\right\}  \ \ \ \ \ \ \ \ \ \ \left(  \text{since }\mathcal{C}\left(
T\right)  \subseteq S_{n}\right)  .
\end{align*}
Thus,%
\begin{align*}
\left\vert C\cap\mathcal{C}\left(  T\right)  \right\vert  &  =\left\vert
\left\{  \text{all transpositions in }\mathcal{C}\left(  T\right)  \right\}
\right\vert \\
&  =\left(  \text{\# of all transpositions in }\mathcal{C}\left(  T\right)
\right)  =\operatorname*{n}\left(  \lambda\right)
\end{align*}
(by (\ref{pf.lem.spechtmod.sum-transp.3.nlam.pf.1})).} and
\begin{equation}
\left\vert C\cap\mathcal{R}\left(  T\right)  \right\vert =\operatorname*{n}%
\left(  \lambda^{t}\right)  \label{pf.lem.spechtmod.sum-transp.3.nlamt}%
\end{equation}
\footnote{\textit{Proof:} Consider the $n$-tableau $T\mathbf{r}$ defined in
Definition \ref{def.tableaux.r}. This is an $n$-tableau of shape
$\mathbf{r}\left(  Y\left(  \lambda\right)  \right)  $ (since $T$ is an
$n$-tableau of shape $Y\left(  \lambda\right)  $). But the definition of
$\lambda^{t}$ shows that $Y\left(  \lambda^{t}\right)  =\mathbf{r}\left(
Y\left(  \lambda\right)  \right)  $. Hence, $T\mathbf{r}$ is an $n$-tableau of
shape $Y\left(  \lambda^{t}\right)  $ (since $T\mathbf{r}$ is an $n$-tableau
of shape $\mathbf{r}\left(  Y\left(  \lambda\right)  \right)  $).
\par
Proposition \ref{prop.partitions.nlam-counts-t} shows that $\operatorname*{n}%
\left(  \lambda\right)  $ is the number of transpositions in $\mathcal{C}%
\left(  T\right)  $. The same reasoning (applied to $\lambda^{t}$ and
$T\mathbf{r}$ instead of $\lambda$ and $T$) shows that $\operatorname*{n}%
\left(  \lambda^{t}\right)  $ is the number of transpositions in
$\mathcal{C}\left(  T\mathbf{r}\right)  $ (since $T\mathbf{r}$ is an
$n$-tableau of shape $Y\left(  \lambda^{t}\right)  $). In other words,%
\begin{equation}
\operatorname*{n}\left(  \lambda^{t}\right)  =\left(  \text{\# of all
transpositions in }\mathcal{C}\left(  T\mathbf{r}\right)  \right)  .
\label{pf.lem.spechtmod.sum-transp.3.nlamt.pf.1}%
\end{equation}
But Proposition \ref{prop.tableaux.r} \textbf{(e)} (applied to $D=Y\left(
\lambda\right)  $) yields $\mathcal{R}\left(  T\right)  =\mathcal{C}\left(
T\mathbf{r}\right)  $. Hence, we can rewrite
(\ref{pf.lem.spechtmod.sum-transp.3.nlamt.pf.1}) as%
\begin{equation}
\operatorname*{n}\left(  \lambda^{t}\right)  =\left(  \text{\# of all
transpositions in }\mathcal{R}\left(  T\right)  \right)  .
\label{pf.lem.spechtmod.sum-transp.3.nlamt.pf.2}%
\end{equation}
\par
Recall that $C$ is the set of all transpositions in $S_{n}$. In other words,
$C=\left\{  \text{all transpositions in }S_{n}\right\}  $. Hence,%
\begin{align*}
C\cap\mathcal{R}\left(  T\right)   &  =\left\{  \text{all transpositions in
}S_{n}\right\}  \cap\mathcal{R}\left(  T\right) \\
&  =\left\{  \text{all transpositions in }S_{n}\text{ that belong to
}\mathcal{R}\left(  T\right)  \right\} \\
&  =\left\{  \text{all transpositions in }\mathcal{R}\left(  T\right)
\right\}  \ \ \ \ \ \ \ \ \ \ \left(  \text{since }\mathcal{R}\left(
T\right)  \subseteq S_{n}\right)  .
\end{align*}
Thus,%
\begin{align*}
\left\vert C\cap\mathcal{R}\left(  T\right)  \right\vert  &  =\left\vert
\left\{  \text{all transpositions in }\mathcal{R}\left(  T\right)  \right\}
\right\vert \\
&  =\left(  \text{\# of all transpositions in }\mathcal{R}\left(  T\right)
\right)  =\operatorname*{n}\left(  \lambda^{t}\right)
\end{align*}
(by (\ref{pf.lem.spechtmod.sum-transp.3.nlamt.pf.2})).}.

Now, for each pair $\left(  c,r\right)  \in\mathcal{C}\left(  T\right)
\times\mathcal{R}\left(  T\right)  $, the sign $\left(  -1\right)  ^{c}$ is
either $1$ or $-1$ (since the sign of any permutation is $1$ or $-1$). Hence,
we can split the sum $\sum_{\substack{\left(  c,r\right)  \in\mathcal{C}%
\left(  T\right)  \times\mathcal{R}\left(  T\right)  ;\\cr\in C}}\left(
-1\right)  ^{c}$ as follows:%
\begin{align*}
\sum_{\substack{\left(  c,r\right)  \in\mathcal{C}\left(  T\right)
\times\mathcal{R}\left(  T\right)  ;\\cr\in C}}\left(  -1\right)  ^{c}  &
=\sum_{\substack{\left(  c,r\right)  \in\mathcal{C}\left(  T\right)
\times\mathcal{R}\left(  T\right)  ;\\cr\in C;\\\left(  -1\right)  ^{c}%
=1}}\underbrace{\left(  -1\right)  ^{c}}_{=1}+\sum_{\substack{\left(
c,r\right)  \in\mathcal{C}\left(  T\right)  \times\mathcal{R}\left(  T\right)
;\\cr\in C;\\\left(  -1\right)  ^{c}=-1}}\underbrace{\left(  -1\right)  ^{c}%
}_{=-1}\\
&  =\underbrace{\sum_{\substack{\left(  c,r\right)  \in\mathcal{C}\left(
T\right)  \times\mathcal{R}\left(  T\right)  ;\\cr\in C;\\\left(  -1\right)
^{c}=1}}}_{\substack{=\sum_{\left(  c,r\right)  \in U}\\\text{(by the
definition of }U\text{)}}}1+\underbrace{\sum_{\substack{\left(  c,r\right)
\in\mathcal{C}\left(  T\right)  \times\mathcal{R}\left(  T\right)  ;\\cr\in
C;\\\left(  -1\right)  ^{c}=-1}}}_{\substack{=\sum_{\left(  c,r\right)  \in
V}\\\text{(by the definition of }V\text{)}}}\left(  -1\right) \\
&  =\underbrace{\sum_{\left(  c,r\right)  \in U}1}_{\substack{=\left\vert
U\right\vert \cdot1=\left\vert U\right\vert \\=\left\vert C\cap\mathcal{R}%
\left(  T\right)  \right\vert \\\text{(by
(\ref{pf.lem.spechtmod.sum-transp.3.map2eq}))}}}+\underbrace{\sum_{\left(
c,r\right)  \in V}\left(  -1\right)  }_{\substack{=\left\vert V\right\vert
\cdot\left(  -1\right)  =-\left\vert V\right\vert \\=-\left\vert
C\cap\mathcal{C}\left(  T\right)  \right\vert \\\text{(by
(\ref{pf.lem.spechtmod.sum-transp.3.map1eq}))}}}\\
&  =\left\vert C\cap\mathcal{R}\left(  T\right)  \right\vert +\left(
-\left\vert C\cap\mathcal{C}\left(  T\right)  \right\vert \right) \\
&  =\underbrace{\left\vert C\cap\mathcal{R}\left(  T\right)  \right\vert
}_{\substack{=\operatorname*{n}\left(  \lambda^{t}\right)  \\\text{(by
(\ref{pf.lem.spechtmod.sum-transp.3.nlamt}))}}}-\underbrace{\left\vert
C\cap\mathcal{C}\left(  T\right)  \right\vert }_{\substack{=\operatorname*{n}%
\left(  \lambda\right)  \\\text{(by (\ref{pf.lem.spechtmod.sum-transp.3.nlam}%
))}}}=\operatorname*{n}\left(  \lambda^{t}\right)  -\operatorname*{n}\left(
\lambda\right)  .
\end{align*}
This proves Lemma \ref{lem.spechtmod.sum-transp.3}.
\end{proof}

Now we are ready to prove Theorem \ref{thm.spechtmod.sum-transp}:

\begin{proof}
[Proof of Theorem \ref{thm.spechtmod.sum-transp}.]Let $C$ be the set of all
transpositions in $S_{n}$. Then, $C$ is a conjugacy class of $S_{n}$ (since
all transpositions in $S_{n}$ form a single conjugacy class)\footnote{Strictly
speaking, this is only true for $n\geq2$ (since otherwise, $C$ is empty). We
leave it to the reader to make the corresponding adjustments for the case
$n<2$.}. Let $\mathbf{z}_{C}$ denote the corresponding conjugacy class sum
$\sum_{c\in C}c\in\mathcal{A}$ (as in Definition \ref{def.groups.conj}
\textbf{(c)}).

Recall that $\mathbf{t}$ is the sum of all transpositions in $S_{n}$. In other
words, $\mathbf{t}$ is the sum of all elements of $C$ (since the
transpositions in $S_{n}$ are precisely the elements of $C$ (since $C$ is the
set of all transpositions in $S_{n}$)). In other words, $\mathbf{t}=\sum_{c\in
C}c=\mathbf{z}_{C}$ (by the definition of $\mathbf{z}_{C}$).

Let $v\in\mathcal{S}^{\lambda}$. Pick any $n$-tableau $T$ of shape $Y\left(
\lambda\right)  $ (such a $T$ exists, since $\left\vert Y\left(
\lambda\right)  \right\vert =\left\vert \lambda\right\vert =n$). Thus, $T$ is
an $n$-tableau of shape $\lambda$. Hence, Theorem \ref{thm.spechtmod.zC}
yields%
\[
\mathbf{z}_{C}\mathbf{v}=\underbrace{\left(  \sum_{\substack{\left(
c,r\right)  \in\mathcal{C}\left(  T\right)  \times\mathcal{R}\left(  T\right)
;\\cr\in C}}\left(  -1\right)  ^{c}\right)  }_{\substack{=\operatorname*{n}%
\left(  \lambda^{t}\right)  -\operatorname*{n}\left(  \lambda\right)
\\\text{(by Lemma \ref{lem.spechtmod.sum-transp.3})}}}\mathbf{v}=\left(
\operatorname*{n}\left(  \lambda^{t}\right)  -\operatorname*{n}\left(
\lambda\right)  \right)  \mathbf{v}.
\]
In other words, $\mathbf{tv}=\left(  \operatorname*{n}\left(  \lambda
^{t}\right)  -\operatorname*{n}\left(  \lambda\right)  \right)  \mathbf{v}$
(since $\mathbf{t}=\mathbf{z}_{C}$). This proves Theorem
\ref{thm.spechtmod.sum-transp}.
\end{proof}

\begin{exercise}
\label{exe.spechtmod.sum-transp.gen}Assume that $\mathbf{k}$ is a field.
\medskip

\textbf{(a)} \fbox{3} Let $m\in\mathbb{N}$. Let $\alpha=\left(  \alpha
_{1},\alpha_{2},\ldots,\alpha_{m}\right)  \in\mathbf{k}^{m}$. Prove that the
element $\alpha$ generates the $\mathbf{k}$-algebra $\mathbf{k}^{m}$ (that is,
we have $\mathbf{k}^{m}=\operatorname*{AlgGen}\left\{  \alpha\right\}  $ in
the notation of Definition \ref{def.AlgGen.AlgGen}) if and only if its entries
$\alpha_{1},\alpha_{2},\ldots,\alpha_{m}$ are all distinct. \medskip

\textbf{(b)} \fbox{2} Now assume that $\mathbf{k}$ is a field of
characteristic $0$. Prove that $\mathbf{t}$ generates the center $Z\left(
\mathbf{k}\left[  S_{n}\right]  \right)  $ if and only if the differences
$\operatorname*{n}\left(  \lambda^{t}\right)  -\operatorname*{n}\left(
\lambda\right)  $ for all partitions $\lambda$ of $n$ are distinct.

(This distinctness condition is easily seen to hold for $n\leq5$ and for
$n=7$, but it does not hold for any other $n$.)
\end{exercise}

\subsection{Characters}

The \emph{character} of a representation is one of the key tools in classical
representation theory: It is a much simpler object than the representation
itself, yet encodes many of its properties and -- under reasonable assumptions
-- uniquely determines the representation (whence probably its name). It thus
can be viewed as a \textquotedblleft fingerprint\textquotedblright\ of a
representation. The general theory of characters is explained in many places
(e.g., \cite[Chapter 10]{Artin}, \cite[Chapters 3 and 4]{EGHetc11},
\cite[Chapter 1]{Sagan01}, \cite{Serre77}, \cite[Chapter 14]{Waerde91}), and
the specific case of symmetric groups can be found in many places as well
(e.g., \cite[\S 4.6--\S 4.10]{Sagan01}, \cite{CSScTo10}, \cite[\S 5.4 and
\S 5.11]{Prasad-rep}, \cite[\S 2.4]{JamKer81}, \cite{BleSch05}). We shall here
focus on the basics (once again trading some depth for generality) and on
connecting the characters of the Specht modules of $S_{n}$ to the centralized
Young symmetrizers $\mathbf{E}_{\lambda}$.

\subsubsection{Definition of a character}

As we know from Corollary \ref{cor.rep.G-rep.mod}, the representations of a
group $G$ are just the left $\mathbf{k}\left[  G\right]  $-modules. Thus, we
define characters not just for group representations but for left $R$-modules
for any $\mathbf{k}$-algebra $R$, subject to a freeness
property:\footnote{Recall the notion of the trace of a $\mathbf{k}$-module
endomorphism; we defined this in Subsection \ref{subsec.specht.ysym.qi}.}

\begin{definition}
\label{def.char.char}Let $R$ be a $\mathbf{k}$-algebra, and let $V$ be a left
$R$-module. Assume that $V$ is a free $\mathbf{k}$-module of finite rank
(i.e., with a finite basis).\footnotemark\ For each $r\in R$, we let $\tau
_{r}:V\rightarrow V$ be the $\mathbf{k}$-linear map $v\mapsto rv$ (which we
already used in Theorem \ref{thm.mod.curry-into-EndV}). Then, the
\emph{character} of the left $R$-module $V$ is defined as the map%
\begin{align*}
\chi_{V}:R  &  \rightarrow\mathbf{k},\\
r  &  \mapsto\operatorname*{Tr}\left(  \tau_{r}\right)  .
\end{align*}
(The trace here is well-defined since $V$ is a free $\mathbf{k}$-module of
finite rank.)
\end{definition}

\footnotetext{We could weaken the freeness requirement in Definition
\ref{def.char.char}: Instead of requiring $V$ to be a free $\mathbf{k}$-module
of finite rank, it would suffice to assume that $V$ is \textquotedblleft
finite projective\textquotedblright\ (i.e., a direct addend of a free
$\mathbf{k}$-module of finite rank), since the notion of a trace is still
defined in this case (see, e.g.,
\url{https://math.stackexchange.com/questions/4394614} ). But the above
definition already suffices for the modules we are planning to consider.}The
first, very simple, property of characters is their linearity:

\begin{proposition}
\label{prop.char.klin}Let $R$ be a $\mathbf{k}$-algebra, and let $V$ be a left
$R$-module. Assume that $V$ is a free $\mathbf{k}$-module of finite rank
(i.e., with a finite basis). Then, the character $\chi_{V}$ of $V$ is a
$\mathbf{k}$-linear map.
\end{proposition}

\begin{fineprint}
\begin{proof}
We must prove that $\chi_{V}\left(  \lambda a+\mu b\right)  =\lambda\chi
_{V}\left(  a\right)  +\mu\chi_{V}\left(  b\right)  $ for any $\lambda,\mu
\in\mathbf{k}$ and any $a,b\in R$.

Let $\lambda,\mu\in\mathbf{k}$ and $a,b\in R$. For each $r\in R$, we let
$\tau_{r}:V\rightarrow V$ be the $\mathbf{k}$-linear map $v\mapsto rv$. The
definition of the character $\chi_{V}$ then shows that $\chi_{V}\left(
r\right)  =\operatorname*{Tr}\left(  \tau_{r}\right)  $ for each $r\in R$.
Thus, $\chi_{V}\left(  \lambda a+\mu b\right)  =\operatorname*{Tr}\left(
\tau_{\lambda a+\mu b}\right)  $ and $\chi_{V}\left(  a\right)
=\operatorname*{Tr}\left(  \tau_{a}\right)  $ and $\chi_{V}\left(  b\right)
=\operatorname*{Tr}\left(  \tau_{b}\right)  $.

Recall two simple properties of traces: We have $\operatorname*{Tr}\left(
f+g\right)  =\operatorname*{Tr}f+\operatorname*{Tr}g$ for any two $\mathbf{k}%
$-linear maps $f,g:V\rightarrow V$, and we have $\operatorname*{Tr}\left(
\lambda f\right)  =\lambda\operatorname*{Tr}f$ for any $\lambda\in\mathbf{k}$
and any $\mathbf{k}$-linear map $f:V\rightarrow V$. (Both of these properties
can be easily derived from the analogous properties of traces of matrices.)
These two properties show that the map $\operatorname*{Tr}:\operatorname*{End}%
\nolimits_{\mathbf{k}}V\rightarrow\mathbf{k}$ is $\mathbf{k}$-linear.

Now, let $v\in V$. Then, the definition of $\tau_{a}$ yields that $\tau
_{a}\left(  v\right)  =av$. Similarly, $\tau_{b}\left(  v\right)  =bv$ and
$\tau_{\lambda a+\mu b}\left(  v\right)  =\left(  \lambda a+\mu b\right)  v$.
Thus,%
\[
\tau_{\lambda a+\mu b}\left(  v\right)  =\left(  \lambda a+\mu b\right)
v=\lambda\underbrace{av}_{=\tau_{a}\left(  v\right)  }+\,\mu\underbrace{bv}%
_{=\tau_{b}\left(  v\right)  }=\lambda\tau_{a}\left(  v\right)  +\mu\tau
_{b}\left(  v\right)  =\left(  \lambda\tau_{a}+\mu\tau_{b}\right)  \left(
v\right)  .
\]

Forget that we fixed $v$. We thus have shown that $\tau_{\lambda a+\mu
b}\left(  v\right)  =\left(  \lambda\tau_{a}+\mu\tau_{b}\right)  \left(
v\right)  $ for each $v\in V$. In other words, $\tau_{\lambda a+\mu b}%
=\lambda\tau_{a}+\mu\tau_{b}$.

Now,
\[
\chi_{V}\left(  \lambda a+\mu b\right)  =\operatorname*{Tr}\left(
\underbrace{\tau_{\lambda a+\mu b}}_{=\lambda\tau_{a}+\mu\tau_{b}}\right)
=\operatorname*{Tr}\left(  \lambda\tau_{a}+\mu\tau_{b}\right)  =\lambda
\operatorname*{Tr}\left(  \tau_{a}\right)  +\mu\operatorname*{Tr}\left(
\tau_{b}\right)
\]
(since the map $\operatorname*{Tr}:\operatorname*{End}\nolimits_{\mathbf{k}%
}V\rightarrow\mathbf{k}$ is $\mathbf{k}$-linear). In view of $\chi_{V}\left(
a\right)  =\operatorname*{Tr}\left(  \tau_{a}\right)  $ and $\chi_{V}\left(
b\right)  =\operatorname*{Tr}\left(  \tau_{b}\right)  $, we can rewrite this
as $\chi_{V}\left(  \lambda a+\mu b\right)  =\lambda\chi_{V}\left(  a\right)
+\mu\chi_{V}\left(  b\right)  $. This completes the proof of Proposition
\ref{prop.char.klin}.
\end{proof}
\end{fineprint}

A slightly less obvious yet crucial property of characters is their
\textquotedblleft pseudo-commutativity\textquotedblright:

\begin{proposition}
\label{prop.char.abba}Let $R$ be a $\mathbf{k}$-algebra, and let $V$ be a left
$R$-module. Assume that $V$ is a free $\mathbf{k}$-module of finite rank
(i.e., with a finite basis). Let $a,b\in R$. Then, $\chi_{V}\left(  ab\right)
=\chi_{V}\left(  ba\right)  $.
\end{proposition}

\begin{fineprint}
\begin{proof}
For each $r\in R$, we let $\tau_{r}:V\rightarrow V$ be the $\mathbf{k}$-linear
map $v\mapsto rv$. The definition of the character $\chi_{V}$ then shows that
$\chi_{V}\left(  r\right)  =\operatorname*{Tr}\left(  \tau_{r}\right)  $ for
each $r\in R$. Thus, $\chi_{V}\left(  ab\right)  =\operatorname*{Tr}\left(
\tau_{ab}\right)  $. On the other hand, it is easy to see the following:

\begin{statement}
\textit{Claim 1:} We have $\tau_{ab}=\tau_{a}\circ\tau_{b}$.
\end{statement}

\begin{proof}
[Proof of Claim 1.]Let $v\in V$. Then, the definition of $\tau_{b}$ yields
that $\tau_{b}\left(  v\right)  =bv$. Similarly, $\tau_{ab}\left(  v\right)
=\left(  ab\right)  v$. Furthermore, the definition of $\tau_{a}$ yields
$\tau_{a}\left(  bv\right)  =a\left(  bv\right)  $. Thus, $\tau_{ab}\left(
v\right)  =\left(  ab\right)  v=abv=a\left(  bv\right)  =\tau_{a}\left(
bv\right)  =\tau_{a}\left(  \tau_{b}\left(  v\right)  \right)  $ (since
$bv=\tau_{b}\left(  v\right)  $), so that $\tau_{ab}\left(  v\right)
=\tau_{a}\left(  \tau_{b}\left(  v\right)  \right)  =\left(  \tau_{a}\circ
\tau_{b}\right)  \left(  v\right)  $.

Forget that we fixed $v$. We thus have shown that $\tau_{ab}\left(  v\right)
=\left(  \tau_{a}\circ\tau_{b}\right)  \left(  v\right)  $ for each $v\in V$.
In other words, $\tau_{ab}=\tau_{a}\circ\tau_{b}$. Thus, Claim 1 is proved.
\end{proof}

Now, recall that $\chi_{V}\left(  ab\right)  =\operatorname*{Tr}\left(
\tau_{ab}\right)  =\operatorname*{Tr}\left(  \tau_{a}\circ\tau_{b}\right)  $
(since Claim 1 shows that $\tau_{ab}=\tau_{a}\circ\tau_{b}$). The same
argument (with the roles of $a$ and $b$ interchanged) shows that $\chi
_{V}\left(  ba\right)  =\operatorname*{Tr}\left(  \tau_{b}\circ\tau
_{a}\right)  $. However, Lemma \ref{lem.linalg.Trfg} (applied to $W=V$ and
$f=\tau_{a}$ and $g=\tau_{b}$) shows that $\operatorname*{Tr}\left(  \tau
_{a}\circ\tau_{b}\right)  =\operatorname*{Tr}\left(  \tau_{b}\circ\tau
_{a}\right)  $. Altogether, we thus find%
\[
\chi_{V}\left(  ab\right)  =\operatorname*{Tr}\left(  \tau_{a}\circ\tau
_{b}\right)  =\operatorname*{Tr}\left(  \tau_{b}\circ\tau_{a}\right)
=\chi_{V}\left(  ba\right)  \ \ \ \ \ \ \ \ \ \ \left(  \text{since }\chi
_{V}\left(  ba\right)  =\operatorname*{Tr}\left(  \tau_{b}\circ\tau
_{a}\right)  \right)  .
\]
This proves Proposition \ref{prop.char.abba}.
\end{proof}
\end{fineprint}

\begin{warning}
Proposition \ref{prop.char.abba} does \textbf{not} mean that products inside
characters can be arbitrarily reordered. For instance, when $a,b,c\in R$, the
values $\chi_{V}\left(  abc\right)  $ and $\chi_{V}\left(  bac\right)  $ are
usually not equal. (But the values $\chi_{V}\left(  abc\right)  $ and
$\chi_{V}\left(  bca\right)  $ are equal, since we can apply Proposition
\ref{prop.char.abba} to $bc$ instead of $b$.)
\end{warning}

Our main interest lies in representations of groups, which (as we know) are
particular cases of modules over rings. Thus we observe the following:

\begin{definition}
\label{def.char.G}Let $G$ be a group, and let $V$ be a representation of $G$
over $\mathbf{k}$. Assume that $V$ is a free $\mathbf{k}$-module of finite
rank (i.e., with a finite basis). Then, $V$ is a left $\mathbf{k}\left[
G\right]  $-module, and thus has a character $\chi_{V}$ (as defined in
Definition \ref{def.char.char}). This character $\chi_{V}$ is $\mathbf{k}%
$-linear (by Proposition \ref{prop.char.klin}), and thus is uniquely
determined by its values on the standard basis vectors $g\in G$, that is, by
its restriction $\chi_{V}\mid_{G}$ to the group $G$. This restriction
$\chi_{V}\mid_{G}$ is called the \emph{group character} of the $G$%
-representation $V$.

We will often blur the distinction between the character $\chi_{V}%
:\mathbf{k}\left[  G\right]  \rightarrow\mathbf{k}$ and its restriction
$\left.  \chi_{V}\mid_{G}\right.  :G\rightarrow\mathbf{k}$, since we can
easily compute one from the other.
\end{definition}

A simple property of group characters is that they assign equal values to
conjugate elements:

\begin{corollary}
\label{cor.char.conj=eq}Let $G$ be a group, and let $V$ be a representation of
$G$ over $\mathbf{k}$. Assume that $V$ is a free $\mathbf{k}$-module of finite
rank (i.e., with a finite basis). Let $w,g\in G$. Then, the character
$\chi_{V}$ of $V$ satisfies%
\begin{equation}
\chi_{V}\left(  wgw^{-1}\right)  =\chi_{V}\left(  g\right)  .
\label{eq.cor.char.conj=eq.eq}%
\end{equation}

\end{corollary}

\begin{fineprint}
\begin{proof}
Proposition \ref{prop.char.abba} (applied to $a=wg$ and $b=w^{-1}$) yields
$\chi_{V}\left(  wgw^{-1}\right)  =\chi_{V}\left(  \underbrace{w^{-1}w}%
_{=1}g\right)  =\chi_{V}\left(  g\right)  $. This proves Corollary
\ref{cor.char.conj=eq}.
\end{proof}
\end{fineprint}

\subsubsection{Examples}

As a warm-up, let us compute some group characters:

\begin{example}
\label{exa.char.triv}Let $G$ be any group. Let $\mathbf{k}%
_{\operatorname*{triv}}$ be the trivial $G$-representation on the $\mathbf{k}%
$-module $\mathbf{k}$. (This is the $\mathbf{k}$-module $\mathbf{k}$, equipped
with the trivial $G$-action, which is given by $g\rightharpoonup v=v$ for each
$g\in G$ and $v\in\mathbf{k}$.) Let us denote the character $\chi
_{\mathbf{k}_{\operatorname*{triv}}}$ of $\mathbf{k}_{\operatorname*{triv}}$
by $\chi_{\operatorname*{triv}}$. Then: \medskip

\textbf{(a)} We have
\[
\chi_{\operatorname*{triv}}\left(  g\right)  =1\ \ \ \ \ \ \ \ \ \ \text{for
each }g\in G.
\]

\textbf{(b)} The character $\chi_{\operatorname*{triv}}$ is precisely the
$\mathbf{k}$-algebra morphism $\varepsilon:\mathbf{k}\left[  G\right]
\rightarrow\mathbf{k}$ that sends each $g\in G$ to $1$ (as suggested in Remark
\ref{rmk.eps.all-G}).
\end{example}

\begin{fineprint}
\begin{proof}
\textbf{(a)} Let $g\in G$. We must prove that $\chi_{\operatorname*{triv}%
}\left(  g\right)  =1$.

Recall that $\chi_{\operatorname*{triv}}=\chi_{\mathbf{k}%
_{\operatorname*{triv}}}$ (by the definition of $\chi_{\operatorname*{triv}}%
$). Thus, $\chi_{\operatorname*{triv}}\left(  g\right)  =\chi_{\mathbf{k}%
_{\operatorname*{triv}}}\left(  g\right)  $.

For each $r\in\mathbf{k}\left[  G\right]  $, we let $\tau_{r}:\mathbf{k}%
_{\operatorname*{triv}}\rightarrow\mathbf{k}_{\operatorname*{triv}}$ be the
$\mathbf{k}$-linear map $v\mapsto rv$. Thus, for each $v\in\mathbf{k}%
_{\operatorname*{triv}}$, we have%
\begin{align*}
\tau_{g}\left(  v\right)   &  =gv=g\rightharpoonup
v=v\ \ \ \ \ \ \ \ \ \ \left(  \text{since the }G\text{-action on }%
\mathbf{k}_{\operatorname*{triv}}\text{ is trivial}\right) \\
&  =\operatorname*{id}\nolimits_{\mathbf{k}}\left(  v\right)  .
\end{align*}
In other words, $\tau_{g}=\operatorname*{id}\nolimits_{\mathbf{k}}$. But the
matrix that represents the $\mathbf{k}$-linear map $\operatorname*{id}%
\nolimits_{\mathbf{k}}:\mathbf{k}\rightarrow\mathbf{k}$ with respect to the
basis $\left(  1\right)  $ of $\mathbf{k}$ is the $1\times1$-matrix $\left(
\begin{array}
[c]{c}%
1
\end{array}
\right)  $, whose trace is clearly $1$. Thus, $\operatorname*{Tr}\left(
\operatorname*{id}\nolimits_{\mathbf{k}}\right)  =1$. Altogether, we now have%
\begin{align*}
\chi_{\operatorname*{triv}}\left(  g\right)   &  =\chi_{\mathbf{k}%
_{\operatorname*{triv}}}\left(  g\right)  =\operatorname*{Tr}\left(  \tau
_{g}\right)  \ \ \ \ \ \ \ \ \ \ \left(  \text{by the definition of the
character }\chi_{\mathbf{k}_{\operatorname*{triv}}}\right) \\
&  =\operatorname*{Tr}\left(  \operatorname*{id}\nolimits_{\mathbf{k}}\right)
\ \ \ \ \ \ \ \ \ \ \left(  \text{since }\tau_{g}=\operatorname*{id}%
\nolimits_{\mathbf{k}}\right) \\
&  =1.
\end{align*}
This proves Example \ref{exa.char.triv} \textbf{(a)}. \medskip

\textbf{(b)} Consider the $\mathbf{k}$-algebra morphism $\varepsilon
:\mathbf{k}\left[  G\right]  \rightarrow\mathbf{k}$ that sends each $g\in G$
to $1$ (as suggested in Remark \ref{rmk.eps.all-G}). We must show that
$\chi_{\operatorname*{triv}}=\varepsilon$. Since both maps $\chi
_{\operatorname*{triv}}$ and $\varepsilon$ are $\mathbf{k}$-linear, it
suffices to show that these two maps $\chi_{\operatorname*{triv}}$ and
$\varepsilon$ agree on all standard basis vectors $g$ of $\mathbf{k}\left[
G\right]  $ (since these standard basis vectors span $\mathbf{k}\left[
G\right]  $). In other words, it suffices to show that $\chi
_{\operatorname*{triv}}\left(  g\right)  =\varepsilon\left(  g\right)  $ for
all $g\in G$. But this is easy: For each $g\in G$, we have $\chi
_{\operatorname*{triv}}\left(  g\right)  =1$ (by Example \ref{exa.char.triv}
\textbf{(a)}) and $\varepsilon\left(  g\right)  =1$ (by the definition of
$\varepsilon$), so that $\chi_{\operatorname*{triv}}\left(  g\right)
=1=\varepsilon\left(  g\right)  $. Thus, Example \ref{exa.char.triv}
\textbf{(b)} is proved.
\end{proof}
\end{fineprint}

For symmetric groups, there is an analogue of Example \ref{exa.char.triv}
\textbf{(a)}:

\begin{example}
\label{exa.char.sign}Let $\mathbf{k}_{\operatorname*{sign}}$ be the sign
representation of the symmetric group $S_{n}$. (This is the $\mathbf{k}%
$-module $\mathbf{k}$, equipped with the left $S_{n}$-action given by
$g\rightharpoonup v=\left(  -1\right)  ^{g}v$ for each $g\in S_{n}$ and
$v\in\mathbf{k}$. See Example \ref{exa.rep.Sn-rep.sign} for details.) Let us
denote the character $\chi_{\mathbf{k}_{\operatorname*{sign}}}$ of
$\mathbf{k}_{\operatorname*{sign}}$ by $\chi_{\operatorname*{sign}}$. Then, we
have
\[
\chi_{\operatorname*{sign}}\left(  g\right)  =\left(  -1\right)
^{g}\ \ \ \ \ \ \ \ \ \ \text{for each }g\in S_{n}.
\]

\end{example}

\begin{fineprint}
\begin{proof}
[Proof sketch.]Similar to Example \ref{exa.char.triv} \textbf{(a)}, but now
involving a $\left(  -1\right)  ^{g}$ factor.
\end{proof}
\end{fineprint}

Example \ref{exa.char.triv} \textbf{(a)} can be generalized:

\begin{exercise}
\fbox{1} Let $G$ be any group. Let $V$ be a representation of $G$ over
$\mathbf{k}$. Assume that $V$ is a free $\mathbf{k}$-module of rank $d$, where
$d\in\mathbb{N}$. Let $g\in G$ be an element that acts on $V$ as scaling by a
scalar $\lambda\in\mathbf{k}$ (that is, satisfies $gv=\lambda v$ for all $v\in
V$). Prove that $\chi_{V}\left(  g\right)  =d\lambda$.
\end{exercise}

Slightly more interesting examples of characters come from the natural and
left regular representations of $S_{n}$:

\begin{example}
\label{exa.char.nat}Let $\mathbf{k}^{n}$ be the natural representation of the
symmetric group $S_{n}$. (See Example \ref{exa.rep.Sn-rep.nat} for its
definition.) Let us denote the character $\chi_{\mathbf{k}^{n}}$ of this
representation by $\chi_{\operatorname*{nat}}$. Then, we have
\[
\chi_{\operatorname*{nat}}\left(  g\right)  =\left(  \text{\# of fixed points
of }g\right)  \ \ \ \ \ \ \ \ \ \ \text{for each }g\in S_{n}.
\]

\end{example}

\begin{fineprint}
\begin{proof}
[Proof sketch.]For each $r\in\mathbf{k}\left[  S_{n}\right]  $, we let
$\tau_{r}:\mathbf{k}^{n}\rightarrow\mathbf{k}^{n}$ be the $\mathbf{k}$-linear
map $v\mapsto rv$. Thus, the definition of the character $\chi_{\mathbf{k}%
^{n}}$ shows that $\chi_{\mathbf{k}^{n}}\left(  r\right)  =\operatorname*{Tr}%
\left(  \tau_{r}\right)  $ for each $r\in\mathbf{k}\left[  S_{n}\right]  $.

Now, let $g\in S_{n}$. Then, $g\in S_{n}\subseteq\mathbf{k}\left[
S_{n}\right]  $, so that $\chi_{\mathbf{k}^{n}}\left(  g\right)
=\operatorname*{Tr}\left(  \tau_{g}\right)  $ (since $\chi_{\mathbf{k}^{n}%
}\left(  r\right)  =\operatorname*{Tr}\left(  \tau_{r}\right)  $ for each
$r\in\mathbf{k}\left[  S_{n}\right]  $). In other words,
\begin{equation}
\chi_{\operatorname*{nat}}\left(  g\right)  =\operatorname*{Tr}\left(
\tau_{g}\right)  \label{pf.exa.char.nat.1}%
\end{equation}
(since we denote the character $\chi_{\mathbf{k}^{n}}$ by $\chi
_{\operatorname*{nat}}$).

Let $e_{1},e_{2},\ldots,e_{n}$ be the standard basis vectors of $\mathbf{k}%
^{n}$ (so that each $e_{i}$ is the $n$-tuple $\left(  0,0,\ldots
,0,1,0,0,\ldots,0\right)  $ that has the $1$ in its $i$-th position). The
definition of the place permutation action of $S_{n}$ on $\mathbf{k}^{n}$
easily shows that%
\[
g\rightharpoonup e_{j}=e_{g\left(  j\right)  }\ \ \ \ \ \ \ \ \ \ \text{for
each }j\in\left[  n\right]  .
\]
Thus, for each $j\in\left[  n\right]  $, we have%
\begin{align*}
\tau_{g}\left(  e_{j}\right)   &  =ge_{j}\ \ \ \ \ \ \ \ \ \ \left(  \text{by
the definition of }\tau_{g}\right) \\
&  =g\rightharpoonup e_{j}=e_{g\left(  j\right)  }.
\end{align*}
Hence, the matrix that represents the $\mathbf{k}$-linear map $\tau
_{g}:\mathbf{k}^{n}\rightarrow\mathbf{k}^{n}$ with respect to the basis
$\left(  e_{1},e_{2},\ldots,e_{n}\right)  $ of $\mathbf{k}^{n}$ is the
$n\times n$-matrix $P$ whose $j$-th column is $e_{g\left(  j\right)  }$ for
each $j\in\left[  n\right]  $. Explicitly, the $\left(  i,j\right)  $-th entry
of this matrix $P$ is%
\[%
\begin{cases}
1, & \text{if }i=g\left(  j\right)  ;\\
0, & \text{if }i\neq g\left(  j\right)
\end{cases}
\]
for all $\left(  i,j\right)  \in\left[  n\right]  \times\left[  n\right]  $.
Hence, the diagonal entries of $P$ are the entries%
\[%
\begin{cases}
1, & \text{if }i=g\left(  i\right)  ;\\
0, & \text{if }i\neq g\left(  i\right)
\end{cases}
\ \ \ \ \ \ \ \ \ \ \text{for all }i\in\left[  n\right]  .
\]
Thus, the trace of $P$ is%
\begin{align*}
\operatorname*{Tr}P  &  =\sum_{i\in\left[  n\right]  }%
\begin{cases}
1, & \text{if }i=g\left(  i\right)  ;\\
0, & \text{if }i\neq g\left(  i\right)
\end{cases}
\\
&  =\sum_{\substack{i\in\left[  n\right]  ;\\i=g\left(  i\right)
}}1+\underbrace{\sum_{\substack{i\in\left[  n\right]  ;\\i\neq g\left(
i\right)  }}0}_{=0}=\sum_{\substack{i\in\left[  n\right]  ;\\i=g\left(
i\right)  }}1=\left(  \text{\# of }i\in\left[  n\right]  \text{ such that
}i=g\left(  i\right)  \right)  \cdot1\\
&  =\left(  \text{\# of }i\in\left[  n\right]  \text{ such that }i=g\left(
i\right)  \right) \\
&  =\left(  \text{\# of fixed points of }g\right)  .
\end{align*}

But we know that $P$ is the matrix that represents the $\mathbf{k}$-linear map
$\tau_{g}:\mathbf{k}^{n}\rightarrow\mathbf{k}^{n}$ with respect to the basis
$\left(  e_{1},e_{2},\ldots,e_{n}\right)  $ of $\mathbf{k}^{n}$. Hence,
$\operatorname*{Tr}\left(  \tau_{g}\right)  =\operatorname*{Tr}P=\left(
\text{\# of fixed points of }g\right)  $. Thus, we can rewrite
(\ref{pf.exa.char.nat.1}) as%
\[
\chi_{\operatorname*{nat}}\left(  g\right)  =\left(  \text{\# of fixed points
of }g\right)  .
\]
Thus, Example \ref{exa.char.nat} is proved.
\end{proof}
\end{fineprint}

\begin{example}
\label{exa.char.reg}Let $G$ be a finite group. Consider the left regular
representation $\mathbf{k}\left[  G\right]  $ of $G$. (Its definition was
given in Definition \ref{def.rep.G-rep.lreg}.) Let us denote the character
$\chi_{\mathbf{k}\left[  G\right]  }$ of this representation by $\chi
_{\operatorname*{reg}}$. Then, we have
\[
\chi_{\operatorname*{reg}}\left(  g\right)  =%
\begin{cases}
\left\vert G\right\vert , & \text{if }g=1;\\
0, & \text{if }g\neq1
\end{cases}
\ \ \ \ \ \ \ \ \ \ \text{for each }g\in G.
\]

\end{example}

\begin{fineprint}
\begin{proof}
Let us recall Definition \ref{def.specht.ET.wa}. Thus, for any two elements
$h$ and $h^{\prime}$ of $G$, we have%
\begin{equation}
\left[  h\right]  \left(  h^{\prime}\right)  =%
\begin{cases}
1, & \text{if }h^{\prime}=h;\\
0, & \text{if }h^{\prime}\neq h.
\end{cases}
\label{pf.exa.char.reg.1}%
\end{equation}

For each $r\in\mathbf{k}\left[  G\right]  $, we let $\tau_{r}:\mathbf{k}%
\left[  G\right]  \rightarrow\mathbf{k}\left[  G\right]  $ be the $\mathbf{k}%
$-linear map $v\mapsto rv$. Then, the definition of the character
$\chi_{\mathbf{k}\left[  G\right]  }$ shows that $\chi_{\mathbf{k}\left[
G\right]  }\left(  r\right)  =\operatorname*{Tr}\left(  \tau_{r}\right)  $ for
each $r\in\mathbf{k}\left[  G\right]  $.

Let $g\in G$. Then, $g\in G\subseteq\mathbf{k}\left[  G\right]  $, and thus
$\chi_{\mathbf{k}\left[  G\right]  }\left(  g\right)  =\operatorname*{Tr}%
\left(  \tau_{g}\right)  $ (since $\chi_{\mathbf{k}\left[  G\right]  }\left(
r\right)  =\operatorname*{Tr}\left(  \tau_{r}\right)  $ for each
$r\in\mathbf{k}\left[  G\right]  $). In other words, $\chi
_{\operatorname*{reg}}\left(  g\right)  =\operatorname*{Tr}\left(  \tau
_{g}\right)  $ (since we denote the character $\chi_{\mathbf{k}\left[
G\right]  }$ by $\chi_{\operatorname*{reg}}$).

Now let us recall the following fact:

\begin{statement}
\textit{Claim 1:} Let $f$ be any $\mathbf{k}$-module endomorphism of
$\mathbf{k}\left[  G\right]  $. Then,%
\[
\operatorname*{Tr}f=\sum_{h\in G}\left[  h\right]  \left(  f\left(  h\right)
\right)  .
\]

\end{statement}

\begin{proof}
[Proof of Claim 1.]This is just Claim 1 in the proof of Proposition
\ref{prop.groupalg.trace}, except that we have renamed the summation index $g$
as $h$.
\end{proof}

Applying Claim 1 to $f=\tau_{g}$, we obtain%
\begin{align*}
\operatorname*{Tr}\left(  \tau_{g}\right)   &  =\sum_{h\in G}\left[  h\right]
\left(  \underbrace{\tau_{g}\left(  h\right)  }_{\substack{=gh\\\text{(by the
definition of }\tau_{g}\text{)}}}\right)  =\sum_{h\in G}\underbrace{\left[
h\right]  \left(  gh\right)  }_{\substack{=%
\begin{cases}
1, & \text{if }gh=h;\\
0, & \text{if }gh\neq h
\end{cases}
\\\text{(by (\ref{pf.exa.char.reg.1}), applied to }h^{\prime}=gh\text{)}}}\\
&  =\sum_{h\in G}\underbrace{%
\begin{cases}
1, & \text{if }gh=h;\\
0, & \text{if }gh\neq h
\end{cases}
}_{\substack{=%
\begin{cases}
1, & \text{if }g=1;\\
0, & \text{if }g\neq1
\end{cases}
\\\text{(since the equation }gh=h\text{ is}\\\text{equivalent to }g=1\text{
(because }G\\\text{is a group))}}}=\sum_{h\in G}%
\begin{cases}
1, & \text{if }g=1;\\
0, & \text{if }g\neq1
\end{cases}
\\
&  =\left\vert G\right\vert \cdot%
\begin{cases}
1, & \text{if }g=1;\\
0, & \text{if }g\neq1
\end{cases}
\ \ =%
\begin{cases}
\left\vert G\right\vert \cdot1, & \text{if }g=1;\\
\left\vert G\right\vert \cdot0, & \text{if }g\neq1
\end{cases}
\ \ =%
\begin{cases}
\left\vert G\right\vert , & \text{if }g=1;\\
0, & \text{if }g\neq1
\end{cases}
\end{align*}
(since $\left\vert G\right\vert \cdot1=\left\vert G\right\vert $ and
$\left\vert G\right\vert \cdot0=0$). In view of $\chi_{\operatorname*{reg}%
}\left(  g\right)  =\operatorname*{Tr}\left(  \tau_{g}\right)  $, we can
rewrite this as follows:%
\[
\chi_{\operatorname*{reg}}\left(  g\right)  =%
\begin{cases}
\left\vert G\right\vert , & \text{if }g=1;\\
0, & \text{if }g\neq1.
\end{cases}
\]
This proves Example \ref{exa.char.reg}.
\end{proof}
\end{fineprint}

Soon (in Example \ref{exa.char.Rkn}), we will compute the character of the
zero-sum representation $R\left(  \mathbf{k}^{n}\right)  $ of $S_{n}$ (defined
in Subsection \ref{subsec.rep.G-rep.subreps}) as well. But it is a nice
exercise to compute it by hand for $n=3$:

\begin{example}
\label{exa.char.Rk3}Let $n=3$. Consider the zero-sum representation $R\left(
\mathbf{k}^{n}\right)  $ of $S_{n}$ (defined in Subsection
\ref{subsec.rep.G-rep.subreps}). Let us compute its character $\chi_{R\left(
\mathbf{k}^{n}\right)  }$.

For each $r\in\mathbf{k}\left[  S_{n}\right]  $, we let $\tau_{r}:R\left(
\mathbf{k}^{n}\right)  \rightarrow R\left(  \mathbf{k}^{n}\right)  $ be the
$\mathbf{k}$-linear map $v\mapsto rv$. Then,
\begin{equation}
\chi_{R\left(  \mathbf{k}^{n}\right)  }\left(  r\right)  =\operatorname*{Tr}%
\left(  \tau_{r}\right)  \label{eq.exa.char.Rk3.tr=}%
\end{equation}
for each $r\in\mathbf{k}\left[  S_{n}\right]  $ (by the definition of the character).

We know that the $\mathbf{k}$-module $R\left(  \mathbf{k}^{n}\right)  $ has a
basis
\[
\left(  e_{1}-e_{n},\ e_{2}-e_{n},\ \ldots,\ e_{n-1}-e_{n}\right)  =\left(
e_{1}-e_{3},\ e_{2}-e_{3}\right)  \ \ \ \ \ \ \ \ \ \ \left(  \text{since
}n=3\right)  .
\]
Let us denote this basis by $\left(  f_{1},f_{2}\right)  $, so that
$f_{1}=e_{1}-e_{3}$ and $f_{2}=e_{2}-e_{3}$.

The action of the simple transposition $s_{1}\in S_{3}$ on this basis $\left(
f_{1},f_{2}\right)  $ is given by%
\begin{align*}
\tau_{s_{1}}\left(  f_{1}\right)   &  =s_{1}\underbrace{f_{1}}_{=e_{1}-e_{3}%
}=s_{1}\left(  e_{1}-e_{3}\right)  =e_{2}-e_{3}=f_{2};\\
\tau_{s_{1}}\left(  f_{2}\right)   &  =s_{1}\underbrace{f_{2}}_{=e_{2}-e_{3}%
}=s_{1}\left(  e_{2}-e_{3}\right)  =e_{1}-e_{3}=f_{1}.
\end{align*}
Thus, the matrix that represents the map $\tau_{s_{1}}:R\left(  \mathbf{k}%
^{n}\right)  \rightarrow R\left(  \mathbf{k}^{n}\right)  $ with respect to
this basis is $\left(
\begin{array}
[c]{cc}%
0 & 1\\
1 & 0
\end{array}
\right)  $. Its trace is $0+0=0$, so that we have $\operatorname*{Tr}\left(
\tau_{s_{1}}\right)  =0$. Hence, (\ref{eq.exa.char.Rk3.tr=}) (applied to
$r=s_{1}$) yields $\chi_{R\left(  \mathbf{k}^{n}\right)  }\left(
s_{1}\right)  =\operatorname*{Tr}\left(  \tau_{s_{1}}\right)  =0$.

Let us next compute $\chi_{R\left(  \mathbf{k}^{n}\right)  }\left(
s_{2}\right)  $. Indeed, the action of the simple transposition $s_{2}\in
S_{3}$ on our basis $\left(  f_{1},f_{2}\right)  $ is given by%
\begin{align*}
\tau_{s_{2}}\left(  f_{1}\right)   &  =s_{2}f_{1}=s_{2}\left(  e_{1}%
-e_{3}\right)  =e_{1}-e_{2}\\
&  =f_{1}-f_{2}\ \ \ \ \ \ \ \ \ \ \left(  \text{by a simple calculation}%
\right)  ;\\
\tau_{s_{2}}\left(  f_{2}\right)   &  =s_{2}f_{2}=s_{2}\left(  e_{2}%
-e_{3}\right)  =e_{3}-e_{2}=-f_{2}.
\end{align*}
Thus, the matrix that represents the map $\tau_{s_{2}}:R\left(  \mathbf{k}%
^{n}\right)  \rightarrow R\left(  \mathbf{k}^{n}\right)  $ with respect to
this basis is $\left(
\begin{array}
[c]{cc}%
1 & 0\\
-1 & -1
\end{array}
\right)  $. Its trace is $1+\left(  -1\right)  =0$, so that we have
$\operatorname*{Tr}\left(  \tau_{s_{2}}\right)  =0$. Hence,
(\ref{eq.exa.char.Rk3.tr=}) (applied to $r=s_{2}$) yields $\chi_{R\left(
\mathbf{k}^{n}\right)  }\left(  s_{2}\right)  =\operatorname*{Tr}\left(
\tau_{s_{2}}\right)  =0$. (Alternatively, we could have derived this from
(\ref{eq.cor.char.conj=eq.eq}) and $\chi_{R\left(  \mathbf{k}^{n}\right)
}\left(  s_{1}\right)  =0$, since $s_{2}=t_{1,3}s_{1}t_{1,3}^{-1}$.)

Similarly, we can compute $\chi_{R\left(  \mathbf{k}^{n}\right)  }\left(
g\right)  $ for the other four elements $g\in S_{3}$. A full list of character
values is
\[%
\begin{tabular}
[c]{|c||c|c|c|c|c|c|}\hline
$g$ & $1$ & $s_{1}$ & $s_{2}$ & $t_{1,3}$ & $\operatorname*{cyc}%
\nolimits_{1,2,3}$ & $\operatorname*{cyc}\nolimits_{1,3,2}$\\\hline\hline
$\chi_{R\left(  \mathbf{k}^{n}\right)  }\left(  g\right)  $ & $2$ & $0$ & $0$
& $0$ & $-1$ & $-1$\\\hline
$\chi_{\operatorname*{triv}}\left(  g\right)  $ & $1$ & $1$ & $1$ & $1$ & $1$
& $1$\\\hline
$\chi_{\operatorname*{sign}}\left(  g\right)  $ & $1$ & $-1$ & $-1$ & $-1$ &
$1$ & $1$\\\hline
$\chi_{\operatorname*{nat}}\left(  g\right)  $ & $3$ & $1$ & $1$ & $1$ & $0$ &
$0$\\\hline
$\chi_{\operatorname*{reg}}\left(  g\right)  $ & $6$ & $0$ & $0$ & $0$ & $0$ &
$0$\\\hline
\end{tabular}
\ \ \ .
\]
(For comparison, we have also listed the values of the characters
$\chi_{\operatorname*{triv}},\chi_{\operatorname*{sign}},\chi
_{\operatorname*{nat}},\chi_{\operatorname*{reg}}$ computed in Examples
\ref{exa.char.triv}, \ref{exa.char.sign}, \ref{exa.char.nat},
\ref{exa.char.reg}.)
\end{example}

The following exercise generalizes both Example \ref{exa.char.nat} and Example
\ref{exa.char.reg} (why?):

\begin{exercise}
\fbox{1} Let $G$ be a group, and let $X$ be a finite left $G$-set. Let
$\mathbf{k}^{\left(  X\right)  }$ be the permutation module corresponding to
the left $G$-set $X$. (See Example \ref{exa.mod.kG-on-kX} for its definition.)
Then, $\mathbf{k}^{\left(  X\right)  }$ is a left $\mathbf{k}\left[  G\right]
$-module, i.e., a $G$-representation. Let us denote its character
$\chi_{\mathbf{k}^{\left(  X\right)  }}$ by $\chi_{X}$. Prove that each $g\in
G$ satisfies%
\[
\chi_{X}\left(  g\right)  =\left(  \text{\# of }x\in X\text{ that satisfy
}g\rightharpoonup x=x\right)  .
\]

\end{exercise}

Another example concerns an arbitrary $R$-module, but a rather specific
element of $R$:

\begin{example}
\label{exa.char.1}Let $R$ be a $\mathbf{k}$-algebra, and let $V$ be a left
$R$-module. Assume that $V$ is a free $\mathbf{k}$-module of rank $d$, where
$d\in\mathbb{N}$. Then, $\chi_{V}\left(  1_{R}\right)  =d\cdot1_{\mathbf{k}}$.
\end{example}

\begin{fineprint}
\begin{proof}
For each $r\in R$, we let $\tau_{r}:V\rightarrow V$ be the $\mathbf{k}$-linear
map $v\mapsto rv$. The definition of the character $\chi_{V}$ then shows that
$\chi_{V}\left(  r\right)  =\operatorname*{Tr}\left(  \tau_{r}\right)  $ for
each $r\in R$. Hence, $\chi_{V}\left(  1_{R}\right)  =\operatorname*{Tr}%
\left(  \tau_{1_{R}}\right)  $.

But each $v\in V$ satisfies%
\begin{align*}
\tau_{1_{R}}\left(  v\right)   &  =1_{R}v\ \ \ \ \ \ \ \ \ \ \left(  \text{by
the definition of }\tau_{1_{R}}\right) \\
&  =v\ \ \ \ \ \ \ \ \ \ \left(  \text{by one of the module axioms}\right) \\
&  =\operatorname*{id}\nolimits_{V}\left(  v\right)  .
\end{align*}
Hence, $\tau_{1_{R}}=\operatorname*{id}\nolimits_{V}$. Thus,
$\operatorname*{Tr}\left(  \tau_{1_{R}}\right)  =\operatorname*{Tr}\left(
\operatorname*{id}\nolimits_{V}\right)  $.

But $V$ is a free $\mathbf{k}$-module of rank $d$, and thus has a basis
consisting of $d$ vectors. The matrix that represents the $\mathbf{k}$-linear
map $\operatorname*{id}\nolimits_{V}:V\rightarrow V$ with respect to this
basis is the $d\times d$-identity matrix, and thus has trace $d\cdot
1_{\mathbf{k}}$. Therefore, the map $\operatorname*{id}\nolimits_{V}$ has
trace $d\cdot1_{\mathbf{k}}$. In other words, $\operatorname*{Tr}\left(
\operatorname*{id}\nolimits_{V}\right)  =d\cdot1_{\mathbf{k}}$. Combining all
we have found, we obtain%
\[
\chi_{V}\left(  1_{R}\right)  =\operatorname*{Tr}\left(  \tau_{1_{R}}\right)
=\operatorname*{Tr}\left(  \operatorname*{id}\nolimits_{V}\right)
=d\cdot1_{\mathbf{k}}.
\]
This proves Example \ref{exa.char.1}.
\end{proof}
\end{fineprint}

\subsubsection{Characters and isomorphism}

Isomorphic $R$-modules have the same character:

\begin{proposition}
\label{prop.char.iso}Let $R$ be a $\mathbf{k}$-algebra, and let $V$ and $W$ be
two left $R$-modules. Assume that $V$ and $W$ are free $\mathbf{k}$-modules of
finite ranks. If $V\cong W$ as left $R$-modules, then $\chi_{V}=\chi_{W}$.
\end{proposition}

\begin{fineprint}
\begin{proof}
Assume that $V\cong W$ as left $R$-modules. Thus, there exists a left
$R$-module isomorphism $f:V\rightarrow W$. Consider this $f$.

Let $r\in R$.

Let $\tau_{r}:V\rightarrow V$ be the $\mathbf{k}$-linear map $v\mapsto rv$.
The definition of the character $\chi_{V}$ then shows that $\chi_{V}\left(
r\right)  =\operatorname*{Tr}\left(  \tau_{r}\right)  $.

Let $\sigma_{r}:W\rightarrow W$ be the $\mathbf{k}$-linear map $v\mapsto rv$.
Thus, $\sigma_{r}:W\rightarrow W$ is the $\mathbf{k}$-linear map that is
defined just like $\tau_{r}$ but using $W$ instead of $V$. Hence, the
definition of the character $\chi_{W}$ shows that $\chi_{W}\left(  r\right)
=\operatorname*{Tr}\left(  \sigma_{r}\right)  $.

The map $f:V\rightarrow W$ is a left $R$-module isomorphism, thus left
$R$-linear. Hence, $f\left(  rv\right)  =rf\left(  v\right)  $ for each $v\in
V$. Thus, we easily obtain $f\circ\tau_{r}=\sigma_{r}\circ f$%
\ \ \ \ \footnote{\textit{Proof.} Let $v\in V$. Then, $\tau_{r}\left(
v\right)  =rv$ (by the definition of $\tau_{r}$) and $\sigma_{r}\left(
f\left(  v\right)  \right)  =rf\left(  v\right)  $ (by the definition of
$\sigma_{r}$). But $f\left(  rv\right)  =rf\left(  v\right)  $ (since the map
$f$ is left $R$-linear). Hence,
\[
\left(  f\circ\tau_{r}\right)  \left(  v\right)  =f\left(  \underbrace{\tau
_{r}\left(  v\right)  }_{=rv}\right)  =f\left(  rv\right)  =rf\left(
v\right)  =\sigma_{r}\left(  f\left(  v\right)  \right)  =\left(  \sigma
_{r}\circ f\right)  \left(  v\right)  .
\]
\par
Forget that we fixed $v$. We thus have shown that $\left(  f\circ\tau
_{r}\right)  \left(  v\right)  =\left(  \sigma_{r}\circ f\right)  \left(
v\right)  $ for each $v\in V$. In other words, $f\circ\tau_{r}=\sigma_{r}\circ
f$.}. Hence, $\sigma_{r}\circ f=f\circ\tau_{r}$.

But Lemma \ref{lem.linalg.Trfg} (applied to $g=f^{-1}\circ\sigma_{r}$) yields
$\operatorname*{Tr}\left(  f\circ f^{-1}\circ\sigma_{r}\right)
=\operatorname*{Tr}\left(  f^{-1}\circ\sigma_{r}\circ f\right)  $. In view of
$\underbrace{f\circ f^{-1}}_{=\operatorname*{id}}\circ\,\sigma_{r}=\sigma_{r}$
and $f^{-1}\circ\underbrace{\sigma_{r}\circ f}_{=f\circ\tau_{r}}%
=\underbrace{f^{-1}\circ f}_{=\operatorname*{id}}\circ\,\tau_{r}=\tau_{r}$, we
can rewrite this as $\operatorname*{Tr}\left(  \sigma_{r}\right)
=\operatorname*{Tr}\left(  \tau_{r}\right)  $. Hence, $\chi_{W}\left(
r\right)  =\operatorname*{Tr}\left(  \sigma_{r}\right)  =\operatorname*{Tr}%
\left(  \tau_{r}\right)  =\chi_{V}\left(  r\right)  $, so that $\chi
_{V}\left(  r\right)  =\chi_{W}\left(  r\right)  $.

Forget that we fixed $r$. We thus have shown that $\chi_{V}\left(  r\right)
=\chi_{W}\left(  r\right)  $ for each $r\in R$. In other words, $\chi_{V}%
=\chi_{W}$. This proves Proposition \ref{prop.char.iso}.
\end{proof}
\end{fineprint}

What is more surprising is that under certain \textquotedblleft
fair-weather\textquotedblright\ conditions -- specifically, in the classical
representation theory of a finite group -- the converse of Proposition
\ref{prop.char.iso} holds as well:

\begin{theorem}
\label{thm.char.iso=char}Assume that $\mathbf{k}$ is a field of characteristic
$0$. Let $G$ be a finite group. Let $V$ and $W$ be two finite-dimensional
representations of $G$ over $\mathbf{k}$. Then, $V\cong W$ if and only if
$\chi_{V}=\chi_{W}$.
\end{theorem}

We will not use this theorem here, so we will not prove it either (but we will
sketch a proof of the particular case $G=S_{n}$ further below: see Proposition
\ref{prop.char.iso-Sn}). Proofs of Theorem \ref{thm.char.iso=char} under
various additional assumptions can be easily found in the literature (e.g., in
\cite[Corollary 4.2.4]{EGHetc11} when $\mathbf{k}$ is algebraically closed, or
in \cite[Corollary 3.3.3]{Webb16} or \cite[\S 2.3, Corollary 2]{Serre77} when
$\mathbf{k}=\mathbb{C}$). A fully general proof appears in \cite[Theorem
1.45]{Lorenz18}.

\subsubsection{The central element corresponding to a character}

Characters of a finite group $G$ give rise to certain elements of
$\mathbf{k}\left[  G\right]  $:

\begin{definition}
\label{def.char.XV}Let $G$ be a finite group. Let $V$ be a representation of
$G$. Assume that $V$ is a free $\mathbf{k}$-module of finite rank. Then, we
define an element $\mathbf{X}_{V}\in\mathbf{k}\left[  G\right]  $ by
\[
\mathbf{X}_{V}:=\sum_{g\in G}\chi_{V}\left(  g\right)  g.
\]
(Here, of course, $\chi_{V}$ means the character of $V$.)
\end{definition}

The element $\mathbf{X}_{V}$ clearly encodes the character $\chi_{V}$, since
each value $\chi_{V}\left(  g\right)  $ (and thus, by linearity, all values of
$\chi_{V}$) can be reconstructed from it by taking the coefficient of $g$.

\begin{example}
\label{exa.char.XV.triv}Let $G$ be any finite group. Let $\mathbf{k}%
_{\operatorname*{triv}}$ be the trivial $G$-representation on the $\mathbf{k}%
$-module $\mathbf{k}$ (as in Example \ref{exa.char.triv}). Then, the element
$\mathbf{X}_{\mathbf{k}_{\operatorname*{triv}}}$ (as defined in Definition
\ref{def.char.XV}) is%
\[
\mathbf{X}_{\mathbf{k}_{\operatorname*{triv}}}=\sum_{g\in G}\underbrace{\chi
_{\mathbf{k}_{\operatorname*{triv}}}\left(  g\right)  }%
_{\substack{=1\\\text{(by Example \ref{exa.char.triv} \textbf{(a)})}}%
}g=\sum_{g\in G}g.
\]
This is precisely the element $\nabla_{G}$ from Remark \ref{rmk.integral.G}.
\end{example}

\begin{proposition}
\label{prop.char.XV.central}Let $G$ be a finite group. Let $V$ be a
representation of $G$. Assume that $V$ is a free $\mathbf{k}$-module of finite
rank. Then, $\mathbf{X}_{V}$ belongs to the center $Z\left(  \mathbf{k}\left[
G\right]  \right)  $ of the ring $\mathbf{k}\left[  G\right]  $.
\end{proposition}

\begin{proof}
Let $\mathbf{a}:=\mathbf{X}_{V}$. Now, let $w\in G$. But $G$ is a group.
Hence, the map%
\begin{align*}
G  &  \rightarrow G,\\
g  &  \mapsto wgw^{-1}%
\end{align*}
is a bijection (and its inverse is the map $G\rightarrow G,\ g\mapsto
w^{-1}gw$).

Now, recall that $\mathbf{a}=\mathbf{X}_{V}=\sum_{g\in G}\chi_{V}\left(
g\right)  g$ (by Definition \ref{def.char.XV}). Thus,
\begin{align*}
w\mathbf{a}w^{-1}  &  =w\left(  \sum_{g\in G}\chi_{V}\left(  g\right)
g\right)  w^{-1}=\sum_{g\in G}\underbrace{\chi_{V}\left(  g\right)
}_{\substack{=\chi_{V}\left(  wgw^{-1}\right)  \\\text{(by
(\ref{eq.cor.char.conj=eq.eq}))}}}wgw^{-1}\\
&  =\sum_{g\in G}\chi_{V}\left(  wgw^{-1}\right)  wgw^{-1}\\
&  =\sum_{g\in G}\chi_{V}\left(  g\right)  g\ \ \ \ \ \ \ \ \ \ \left(
\begin{array}
[c]{c}%
\text{here, we have substituted }g\text{ for }wgw^{-1}\text{ in}\\
\text{the sum, since the map }G\rightarrow G,\ g\mapsto wgw^{-1}\\
\text{is a bijection}%
\end{array}
\right) \\
&  =\mathbf{a}.
\end{align*}

Forget that we fixed $w$. We thus have shown that each $w\in G$ satisfies
$w\mathbf{a}w^{-1}=\mathbf{a}$. Hence, Lemma \ref{lem.center.waw} reveals that
$\mathbf{a}\in Z\left(  \mathbf{k}\left[  G\right]  \right)  $. In other
words, $\mathbf{X}_{V}\in Z\left(  \mathbf{k}\left[  G\right]  \right)  $
(since $\mathbf{a}=\mathbf{X}_{V}$). In other words, $\mathbf{X}_{V}$ belongs
to the center $Z\left(  \mathbf{k}\left[  G\right]  \right)  $ of the ring
$\mathbf{k}\left[  G\right]  $. This proves Proposition
\ref{prop.char.XV.central}.
\end{proof}

\begin{corollary}
\label{cor.char.XV.iso}Let $G$ be a finite group. Let $V$ and $W$ be two
representations of $G$. Assume that $V$ and $W$ are free $\mathbf{k}$-modules
of finite ranks. If $V\cong W$, then $\mathbf{X}_{V}=\mathbf{X}_{W}$.
\end{corollary}

\begin{proof}
Assume that $V\cong W$. Then, Proposition \ref{prop.char.iso} (applied to
$R=\mathbf{k}\left[  G\right]  $) yields $\chi_{V}=\chi_{W}$. But
$\mathbf{X}_{V}$ is defined by the equality $\mathbf{X}_{V}=\sum_{g\in G}%
\chi_{V}\left(  g\right)  g$, whereas $\mathbf{X}_{W}$ is defined by the
equality $\mathbf{X}_{W}=\sum_{g\in G}\chi_{W}\left(  g\right)  g$. The right
hand sides of these two equalities are equal (since $\chi_{V}=\chi_{W}$).
Hence, so are the left hand sides. In other words, $\mathbf{X}_{V}%
=\mathbf{X}_{W}$. This proves Corollary \ref{cor.char.XV.iso}.
\end{proof}

The above facts are the tip of the iceberg, but they are the few things that
can be said about the elements $\mathbf{X}_{V}$ in full generality. In the
classical setup, more is known:

\begin{theorem}
\label{thm.char.XV.basis}Let $\mathbf{k}$ be a field of characteristic $0$.
Let $G$ be a finite group. Let $J_{1},J_{2},\ldots,J_{m}$ be all the
irreducible representations of $G$ over $\mathbf{k}$, each listed exactly once
(up to isomorphism). Then: \medskip

\textbf{(a)} The elements $\mathbf{X}_{J_{1}},\mathbf{X}_{J_{2}}%
,\ldots,\mathbf{X}_{J_{m}}$ of $Z\left(  \mathbf{k}\left[  G\right]  \right)
$ are linearly independent. \medskip

\textbf{(b)} If $\mathbf{k}$ is algebraically closed, then $\left(
\mathbf{X}_{J_{1}},\mathbf{X}_{J_{2}},\ldots,\mathbf{X}_{J_{m}}\right)  $ is a
basis of $Z\left(  \mathbf{k}\left[  G\right]  \right)  $.
\end{theorem}

Theorem \ref{thm.char.XV.basis} \textbf{(b)} is proved in most texts on
representation theory -- e.g., in \cite[Theorem 4.2.1]{EGHetc11}. A proof of
the whole Theorem \ref{thm.char.XV.basis} (and, in fact, of a more general
fact) can be found in \cite[Theorem 1.44]{Lorenz18}. (Note that both of these
sources work not with the elements $\mathbf{X}_{J_{k}}$ but rather with the
characters $\chi_{J_{k}}$ themselves, which belong to the vector space of
class functions\footnote{In \cite{Lorenz18}, class functions are called
\textquotedblleft trace forms\textquotedblright.} on $\mathbf{k}\left[
G\right]  $. This has the advantage of working in a more general setting than
just group algebras. But this can easily be translated into our language when
dealing with a group algebra.)

In the case when $G=S_{n}$, Theorem \ref{thm.char.XV.basis} \textbf{(b)} holds
even without the assumption that $\mathbf{k}$ be algebraically closed; this
follows from Corollary \ref{cor.char.Xlam-prop} \textbf{(d)} (which we will
prove further below).

\subsubsection{Quotients and direct sums}

The next theorem is a highly useful tool in computing characters:

\begin{theorem}
\label{thm.char.M/N}Let $R$ be a $\mathbf{k}$-algebra. Let $M$ be a left
$R$-module. Let $N$ be a left $R$-submodule of $M$. Assume that, as
$\mathbf{k}$-modules, both $N$ and $M/N$ are free of finite ranks (not
necessarily of the same rank). Then: \medskip

\textbf{(a)} The $\mathbf{k}$-module $M$ is free of finite rank as well.
\medskip

\textbf{(b)} We have $\chi_{M/N}=\chi_{M}-\chi_{N}$. (Here, the subtraction is
pointwise, as usual for $\mathbf{k}$-linear maps from $R$ to $\mathbf{k}$.)
\end{theorem}

\begin{proof}
For any vector $m\in M$, we let $\overline{m}$ denote its projection onto the
quotient $M/N$ (that is, the coset $m+N$ of $N$ in $M$).

We assumed that both $\mathbf{k}$-modules $N$ and $M/N$ are free of finite
ranks. Thus, they have bases $\left(  n_{1},n_{2},\ldots,n_{p}\right)  $ and
$\left(  c_{1},c_{2},\ldots,c_{q}\right)  $, respectively. Consider these
bases. For each $i\in\left[  q\right]  $, we pick some vector $m_{i}\in M$
satisfying
\begin{equation}
\overline{m_{i}}=c_{i} \label{pf.thm.char.M/N.mi=ci}%
\end{equation}
(such a vector exists, since $c_{i}\in M/N$). Thus, we have defined $q$
vectors $m_{1},m_{2},\ldots,m_{q}$ in $M$.

Now, we shall show that the list $\left(  n_{1},n_{2},\ldots,n_{p},m_{1}%
,m_{2},\ldots,m_{q}\right)  $ is a basis of the $\mathbf{k}$-module $M$. This
is a classical fact in linear algebra, but let us nevertheless outline the
proof. We must prove the following two claims:

\begin{statement}
\textit{Claim 1:} The list $\left(  n_{1},n_{2},\ldots,n_{p},m_{1}%
,m_{2},\ldots,m_{q}\right)  $ is $\mathbf{k}$-linearly independent.
\end{statement}

\begin{statement}
\textit{Claim 2:} The list $\left(  n_{1},n_{2},\ldots,n_{p},m_{1}%
,m_{2},\ldots,m_{q}\right)  $ spans the $\mathbf{k}$-module $M$.
\end{statement}

\begin{proof}
[Proof of Claim 1.]Let $\nu_{1},\nu_{2},\ldots,\nu_{p},\mu_{1},\mu_{2}%
,\ldots,\mu_{q}\in\mathbf{k}$ be scalars satisfying%
\begin{equation}
\nu_{1}n_{1}+\nu_{2}n_{2}+\cdots+\nu_{p}n_{p}+\mu_{1}m_{1}+\mu_{2}m_{2}%
+\cdots+\mu_{q}m_{q}=0. \label{pf.thm.char.M/N.c1.pf.1}%
\end{equation}
We shall show that all these scalars $\nu_{1},\nu_{2},\ldots,\nu_{p},\mu
_{1},\mu_{2},\ldots,\mu_{q}$ are $0$.

Indeed, (\ref{pf.thm.char.M/N.c1.pf.1}) shows that%
\[
\mu_{1}m_{1}+\mu_{2}m_{2}+\cdots+\mu_{q}m_{q}=-\left(  \nu_{1}n_{1}+\nu
_{2}n_{2}+\cdots+\nu_{p}n_{p}\right)  \in N
\]
(since all the vectors $n_{1},n_{2},\ldots,n_{p}$ belong to the $\mathbf{k}%
$-module $N$). In other words,%
\[
\overline{\mu_{1}m_{1}+\mu_{2}m_{2}+\cdots+\mu_{q}m_{q}}=0
\]
in $M/N$. In view of
\begin{align*}
\overline{\mu_{1}m_{1}+\mu_{2}m_{2}+\cdots+\mu_{q}m_{q}}  &  =\mu
_{1}\underbrace{\overline{m_{1}}}_{\substack{=c_{1}\\\text{(by
(\ref{pf.thm.char.M/N.mi=ci}))}}}+\,\mu_{2}\underbrace{\overline{m_{2}}%
}_{\substack{=c_{2}\\\text{(by (\ref{pf.thm.char.M/N.mi=ci}))}}}+\cdots
+\mu_{q}\underbrace{\overline{m_{q}}}_{\substack{=c_{q}\\\text{(by
(\ref{pf.thm.char.M/N.mi=ci}))}}}\\
&  =\mu_{1}c_{1}+\mu_{2}c_{2}+\cdots+\mu_{q}c_{q},
\end{align*}
we can rewrite this as $\mu_{1}c_{1}+\mu_{2}c_{2}+\cdots+\mu_{q}c_{q}=0$.
Since the list $\left(  c_{1},c_{2},\ldots,c_{q}\right)  $ is $\mathbf{k}%
$-linearly independent (because it is a basis of $M/N$), this shows that%
\begin{equation}
\mu_{1}=\mu_{2}=\cdots=\mu_{q}=0. \label{pf.thm.char.M/N.c1.pf.3}%
\end{equation}
However, (\ref{pf.thm.char.M/N.c1.pf.1}) also shows that%
\begin{align*}
\nu_{1}n_{1}+\nu_{2}n_{2}+\cdots+\nu_{p}n_{p}  &  =-\left(  \mu_{1}m_{1}%
+\mu_{2}m_{2}+\cdots+\mu_{q}m_{q}\right) \\
&  =-\underbrace{\mu_{1}}_{\substack{=0\\\text{(by
(\ref{pf.thm.char.M/N.c1.pf.3}))}}}m_{1}-\underbrace{\mu_{2}}%
_{\substack{=0\\\text{(by (\ref{pf.thm.char.M/N.c1.pf.3}))}}}m_{2}%
-\cdots-\underbrace{\mu_{q}}_{\substack{=0\\\text{(by
(\ref{pf.thm.char.M/N.c1.pf.3}))}}}m_{q}\\
&  =-0m_{1}-0m_{2}-\cdots-0m_{q}=0.
\end{align*}
Since the list $\left(  n_{1},n_{2},\ldots,n_{p}\right)  $ is $\mathbf{k}%
$-linearly independent (because it is a basis of $N$), this shows that%
\[
\nu_{1}=\nu_{2}=\cdots=\nu_{p}=0.
\]
Combining this with (\ref{pf.thm.char.M/N.c1.pf.3}), we conclude that all the
scalars $\nu_{1},\nu_{2},\ldots,\nu_{p},\mu_{1},\mu_{2},\ldots,\mu_{q}$ are
$0$.

Forget that we fixed $\nu_{1},\nu_{2},\ldots,\nu_{p},\mu_{1},\mu_{2}%
,\ldots,\mu_{q}$. We thus have shown that if $\nu_{1},\nu_{2},\ldots,\nu
_{p},\mu_{1},\mu_{2},\ldots,\mu_{q}\in\mathbf{k}$ are scalars satisfying%
\[
\nu_{1}n_{1}+\nu_{2}n_{2}+\cdots+\nu_{p}n_{p}+\mu_{1}m_{1}+\mu_{2}m_{2}%
+\cdots+\mu_{q}m_{q}=0,
\]
then all the scalars $\nu_{1},\nu_{2},\ldots,\nu_{p},\mu_{1},\mu_{2}%
,\ldots,\mu_{q}$ are $0$. In other words, the list $\left(  n_{1},n_{2}%
,\ldots,n_{p},m_{1},m_{2},\ldots,m_{q}\right)  $ is $\mathbf{k}$-linearly
independent. This proves Claim 1.
\end{proof}

\begin{proof}
[Proof of Claim 2.]Let $v\in M$ be arbitrary. Then, the vector $\overline
{v}\in M/N$ can be written as a $\mathbf{k}$-linear combination of
$c_{1},c_{2},\ldots,c_{q}$ (since $\left(  c_{1},c_{2},\ldots,c_{q}\right)  $
is a basis of $M/N$). In other words, $\overline{v}$ can be written as%
\begin{equation}
\overline{v}=\alpha_{1}c_{1}+\alpha_{2}c_{2}+\cdots+\alpha_{q}c_{q}
\label{pf.thm.char.M/N.c2.pf.1}%
\end{equation}
for some scalars $\alpha_{1},\alpha_{2},\ldots,\alpha_{q}\in\mathbf{k}$.
Consider these scalars $\alpha_{1},\alpha_{2},\ldots,\alpha_{q}$.

Let $w:=\alpha_{1}m_{1}+\alpha_{2}m_{2}+\cdots+\alpha_{q}m_{q}$. Thus,%
\begin{align*}
\overline{w}  &  =\overline{\alpha_{1}m_{1}+\alpha_{2}m_{2}+\cdots+\alpha
_{q}m_{q}}=\alpha_{1}\underbrace{\overline{m_{1}}}_{\substack{=c_{1}%
\\\text{(by (\ref{pf.thm.char.M/N.mi=ci}))}}}+\,\alpha_{2}%
\underbrace{\overline{m_{2}}}_{\substack{=c_{2}\\\text{(by
(\ref{pf.thm.char.M/N.mi=ci}))}}}+\cdots+\alpha_{q}\underbrace{\overline
{m_{q}}}_{\substack{=c_{q}\\\text{(by (\ref{pf.thm.char.M/N.mi=ci}))}}}\\
&  =\alpha_{1}c_{1}+\alpha_{2}c_{2}+\cdots+\alpha_{q}c_{q}=\overline
{v}\ \ \ \ \ \ \ \ \ \ \left(  \text{by (\ref{pf.thm.char.M/N.c2.pf.1}%
)}\right)  .
\end{align*}
Thus, $\overline{v}=\overline{w}$. In other words, $v-w\in N$. Hence, the
vector $v-w$ can be written as a $\mathbf{k}$-linear combination of
$n_{1},n_{2},\ldots,n_{p}$ (since $\left(  n_{1},n_{2},\ldots,n_{p}\right)  $
is a basis of $N$). In other words, $v-w$ can be written as%
\begin{equation}
v-w=\beta_{1}n_{1}+\beta_{2}n_{2}+\cdots+\beta_{p}n_{p}
\label{pf.thm.char.M/N.c2.pf.2}%
\end{equation}
for some scalars $\beta_{1},\beta_{2},\ldots,\beta_{p}\in\mathbf{k}$. Consider
these scalars $\beta_{1},\beta_{2},\ldots,\beta_{p}$.

From (\ref{pf.thm.char.M/N.c2.pf.2}), we obtain%
\begin{align*}
v  &  =\beta_{1}n_{1}+\beta_{2}n_{2}+\cdots+\beta_{p}n_{p}+\underbrace{w}%
_{=\alpha_{1}m_{1}+\alpha_{2}m_{2}+\cdots+\alpha_{q}m_{q}}\\
&  =\beta_{1}n_{1}+\beta_{2}n_{2}+\cdots+\beta_{p}n_{p}+\alpha_{1}m_{1}%
+\alpha_{2}m_{2}+\cdots+\alpha_{q}m_{q}.
\end{align*}
Hence, $v$ is a $\mathbf{k}$-linear combination of the vectors $n_{1}%
,n_{2},\ldots,n_{p},m_{1},m_{2},\ldots,m_{q}$.

Forget that we fixed $v$. We thus have shown that each $v\in M$ is a
$\mathbf{k}$-linear combination of the vectors $n_{1},n_{2},\ldots,n_{p}%
,m_{1},m_{2},\ldots,m_{q}$. In other words, the list $\left(  n_{1}%
,n_{2},\ldots,n_{p},m_{1},m_{2},\ldots,m_{q}\right)  $ spans the $\mathbf{k}%
$-module $M$. This proves Claim 2.
\end{proof}

Combining Claim 1 with Claim 2, we see that the list $\left(  n_{1}%
,n_{2},\ldots,n_{p},m_{1},m_{2},\ldots,m_{q}\right)  $ is a basis of the
$\mathbf{k}$-module $M$. \medskip

\textbf{(a)} We have shown above that the list $\left(  n_{1},n_{2}%
,\ldots,n_{p},m_{1},m_{2},\ldots,m_{q}\right)  $ (which has size $p+q$) is a
basis of $M$. Hence, the $\mathbf{k}$-module $M$ has a finite basis (of size
$p+q$). In other words, the $\mathbf{k}$-module $M$ is free of finite rank.
This proves Theorem \ref{thm.char.M/N} \textbf{(a)}. \medskip

\textbf{(b)} Recall that the $p$ vectors $n_{1},n_{2},\ldots,n_{p}$ belong to
$N$ and thus to $M$ (since $N\subseteq M$). Let us define $q$ further vectors
$n_{p+1},n_{p+2},\ldots,n_{p+q}$ in $M$ by setting%
\begin{equation}
n_{i}:=m_{i-p}\ \ \ \ \ \ \ \ \ \ \text{for each }i\in\left\{  p+1,p+2,\ldots
,p+q\right\}  . \label{pf.thm.char.M/N.b.ni=mi-p}%
\end{equation}
Thus, we have defined altogether $p+q$ vectors $n_{1},n_{2},\ldots,n_{p+q}$ in
$M$. Note that the $q$ new vectors $n_{p+1},n_{p+2},\ldots,n_{p+q}$ do not
belong to $N$ any more.

For each $j\in\left[  q\right]  $, we have $p+j\in\left\{  p+1,p+2,\ldots
,p+q\right\}  $ and thus%
\begin{align}
n_{p+j}  &  =m_{p+j-p}\ \ \ \ \ \ \ \ \ \ \left(  \text{by
(\ref{pf.thm.char.M/N.b.ni=mi-p}), applied to }i=p+j\right) \nonumber\\
&  =m_{j}\ \ \ \ \ \ \ \ \ \ \left(  \text{since }p+j-p=j\right)  .
\label{pf.thm.char.M/N.b.np+j=}%
\end{align}

In other words,%
\[
\left(  n_{p+1},n_{p+2},\ldots,n_{p+q}\right)  =\left(  m_{1},m_{2}%
,\ldots,m_{q}\right)  .
\]
Moreover,%
\begin{align*}
\left(  n_{1},n_{2},\ldots,n_{p+q}\right)   &  =\left(  n_{1},n_{2}%
,\ldots,n_{p},\underbrace{n_{p+1}}_{\substack{=m_{1}\\\text{(by
(\ref{pf.thm.char.M/N.b.np+j=}))}}},\underbrace{n_{p+2}}_{\substack{=m_{2}%
\\\text{(by (\ref{pf.thm.char.M/N.b.np+j=}))}}},\ldots,\underbrace{n_{p+q}%
}_{\substack{=m_{q}\\\text{(by (\ref{pf.thm.char.M/N.b.np+j=}))}}}\right) \\
&  =\left(  n_{1},n_{2},\ldots,n_{p},m_{1},m_{2},\ldots,m_{q}\right)  .
\end{align*}
Hence, the list $\left(  n_{1},n_{2},\ldots,n_{p+q}\right)  $ is a basis of
the $\mathbf{k}$-module $M$ (since the list $\left(  n_{1},n_{2},\ldots
,n_{p},m_{1},m_{2},\ldots,m_{q}\right)  $ is a basis of the $\mathbf{k}%
$-module $M$). \medskip

Fix $r\in R$. For any left $R$-module $V$, we let $\tau_{r,V}:V\rightarrow V$
be the $\mathbf{k}$-linear map $v\mapsto rv$. (This is the map denoted by
$\tau_{r}$ in Definition \ref{def.char.char}.) The definition of the character
$\chi_{V}$ yields that%
\begin{equation}
\chi_{V}\left(  r\right)  =\operatorname*{Tr}\left(  \tau_{r,V}\right)
\label{pf.thm.char.M/N.b.chi=}%
\end{equation}
whenever $V$ is a left $R$-module that is free of finite rank as a
$\mathbf{k}$-module.

We shall use the notation $\left(  a_{i,j}\right)  _{i,j\in\left[  k\right]
}$ for the $k\times k$-matrix whose $\left(  i,j\right)  $-th entry is
$a_{i,j}$ for all $\left(  i,j\right)  \in\left[  k\right]  \times\left[
k\right]  $. Clearly, the trace of this matrix is $\sum_{i=1}^{k}a_{i,i}$.

Consider the $\mathbf{k}$-linear map $\tau_{r,M}:M\rightarrow M$. Let $\left(
a_{i,j}\right)  _{i,j\in\left[  p+q\right]  }$ be the matrix that represents
this map $\tau_{r,M}$ with respect to the basis $\left(  n_{1},n_{2}%
,\ldots,n_{p+q}\right)  $ of $M$. Thus, we have%
\begin{equation}
\tau_{r,M}\left(  n_{j}\right)  =\sum_{i=1}^{p+q}a_{i,j}n_{i}
\label{pf.thm.char.M/N.b.tauM}%
\end{equation}
for each $j\in\left[  p+q\right]  $. Hence, for each $j\in\left[  p+q\right]
$, we have%
\begin{align}
rn_{j}  &  =\tau_{r,M}\left(  n_{j}\right)  \ \ \ \ \ \ \ \ \ \ \left(
\begin{array}
[c]{c}%
\text{since the definition of }\tau_{r,M}\\
\text{shows that }\tau_{r,M}\left(  n_{j}\right)  =rn_{j}%
\end{array}
\right) \nonumber\\
&  =\sum_{i=1}^{p+q}a_{i,j}n_{i}\ \ \ \ \ \ \ \ \ \ \left(  \text{by
(\ref{pf.thm.char.M/N.b.tauM})}\right) \nonumber\\
&  =\sum_{i=1}^{p}a_{i,j}n_{i}+\sum_{i=p+1}^{p+q}a_{i,j}\underbrace{n_{i}%
}_{\substack{=m_{i-p}\\\text{(by (\ref{pf.thm.char.M/N.b.ni=mi-p}))}%
}}\nonumber\\
&  =\sum_{i=1}^{p}a_{i,j}n_{i}+\sum_{i=p+1}^{p+q}a_{i,j}m_{i-p}\nonumber\\
&  =\sum_{i=1}^{p}a_{i,j}n_{i}+\sum_{i=1}^{q}a_{i+p,j}m_{i}
\label{pf.thm.char.M/N.b.rnj=s+s}%
\end{align}
(here, we have substituted $i+p$ for $i$ in the second sum).

Now we claim the following:

\begin{statement}
\textit{Claim 3:} We have $\tau_{r,N}\left(  n_{j}\right)  =\sum_{i=1}%
^{p}a_{i,j}n_{i}$ for each $j\in\left[  p\right]  $.
\end{statement}

\begin{proof}
[Proof of Claim 3.]Let $j\in\left[  p\right]  $. Then, $n_{j}\in N$ (since the
$p$ vectors $n_{1},n_{2},\ldots,n_{p}$ belong to $N$). Hence, $\tau
_{r,N}\left(  n_{j}\right)  =rn_{j}$ (by the definition of $\tau_{r,N}$).
Thus,%
\begin{equation}
\tau_{r,N}\left(  n_{j}\right)  =rn_{j}=\sum_{i=1}^{p}a_{i,j}n_{i}+\sum
_{i=1}^{q}a_{i+p,j}m_{i} \label{pf.thm.char.M/N.b.c3.pf.1}%
\end{equation}
(by (\ref{pf.thm.char.M/N.b.rnj=s+s})). Solving this for $\sum_{i=1}%
^{q}a_{i+p,j}m_{i}$, we obtain%
\[
\sum_{i=1}^{q}a_{i+p,j}m_{i}=\tau_{r,N}\left(  n_{j}\right)  -\sum_{i=1}%
^{p}a_{i,j}n_{i}\in N
\]
(since the vector $\tau_{r,N}\left(  n_{j}\right)  $ as well as the $p$
vectors $n_{1},n_{2},\ldots,n_{p}$ belong to $N$). Hence, in $M/N$, we have%
\[
\overline{\sum_{i=1}^{q}a_{i+p,j}m_{i}}=0.
\]
In view of%
\[
\overline{\sum_{i=1}^{q}a_{i+p,j}m_{i}}=\sum_{i=1}^{q}a_{i+p,j}%
\underbrace{\overline{m_{i}}}_{\substack{=c_{i}\\\text{(by
(\ref{pf.thm.char.M/N.mi=ci}))}}}=\sum_{i=1}^{q}a_{i+p,j}c_{i}=a_{1+p,j}%
c_{1}+a_{2+p,j}c_{2}+\cdots+a_{q+p,j}c_{q},
\]
we can rewrite this as $a_{1+p,j}c_{1}+a_{2+p,j}c_{2}+\cdots+a_{q+p,j}c_{q}%
=0$. Since the list $\left(  c_{1},c_{2},\ldots,c_{q}\right)  $ is
$\mathbf{k}$-linearly independent (because it is a basis of $M/N$), this shows
that all the coefficients $a_{1+p,j},a_{2+p,j},\ldots,a_{q+p,j}$ are $0$. In
other words,
\begin{equation}
a_{i+p,j}=0\ \ \ \ \ \ \ \ \ \ \text{for each }i\in\left[  q\right]  .
\label{pf.thm.char.M/N.b.c3.pf.5}%
\end{equation}
Hence, (\ref{pf.thm.char.M/N.b.c3.pf.1}) becomes%
\[
\tau_{r,N}\left(  n_{j}\right)  =\sum_{i=1}^{p}a_{i,j}n_{i}+\sum_{i=1}%
^{q}\underbrace{a_{i+p,j}}_{\substack{=0\\\text{(by
(\ref{pf.thm.char.M/N.b.c3.pf.5}))}}}m_{i}=\sum_{i=1}^{p}a_{i,j}%
n_{i}+\underbrace{\sum_{i=1}^{q}0m_{i}}_{=0}=\sum_{i=1}^{p}a_{i,j}n_{i}.
\]
This proves Claim 3.
\end{proof}

\begin{statement}
\textit{Claim 4:} We have $\tau_{r,M/N}\left(  c_{j}\right)  =\sum_{i=1}%
^{q}a_{i+p,j+p}c_{i}$ for each $j\in\left[  q\right]  $.
\end{statement}

\begin{proof}
[Proof of Claim 4.]Let $j\in\left[  q\right]  $. Then,
(\ref{pf.thm.char.M/N.b.np+j=}) yields $n_{p+j}=m_{j}$. In other words,
$n_{j+p}=m_{j}$ (since $p+j=j+p$). Hence, $rn_{j+p}=rm_{j}$, so that%
\begin{equation}
rm_{j}=rn_{j+p}=\sum_{i=1}^{p}a_{i,j+p}n_{i}+\sum_{i=1}^{q}a_{i+p,j+p}m_{i}
\label{pf.thm.char.M/N.b.c4.pf.2}%
\end{equation}
(by (\ref{pf.thm.char.M/N.b.rnj=s+s}), applied to $j+p$ instead of $j$).

Furthermore, $\overline{m_{j}}=c_{j}$ (by (\ref{pf.thm.char.M/N.mi=ci}),
applied to $i=j$). Hence, $r\overline{m_{j}}=rc_{j}$, so that
\begin{align*}
rc_{j}  &  =r\overline{m_{j}}=\overline{rm_{j}}=\overline{\sum_{i=1}%
^{p}a_{i,j+p}n_{i}+\sum_{i=1}^{q}a_{i+p,j+p}m_{i}}\ \ \ \ \ \ \ \ \ \ \left(
\text{by (\ref{pf.thm.char.M/N.b.c4.pf.2})}\right) \\
&  =\sum_{i=1}^{p}a_{i,j+p}\underbrace{\overline{n_{i}}}%
_{\substack{=0\\\text{(since }n_{i}\in N\\\text{(because the }p\\\text{vectors
}n_{1},n_{2},\ldots,n_{p}\\\text{belong to }N\text{))}}}+\sum_{i=1}%
^{q}a_{i+p,j+p}\underbrace{\overline{m_{i}}}_{\substack{=c_{i}\\\text{(by
(\ref{pf.thm.char.M/N.mi=ci}))}}}\\
&  =\underbrace{\sum_{i=1}^{p}a_{i,j+p}0}_{=0}+\sum_{i=1}^{q}a_{i+p,j+p}%
c_{i}=\sum_{i=1}^{q}a_{i+p,j+p}c_{i}.
\end{align*}
But the definition of $\tau_{r,M/N}$ yields $\tau_{r,M/N}\left(  c_{j}\right)
=rc_{j}=\sum_{i=1}^{q}a_{i+p,j+p}c_{i}$. This proves Claim 4.
\end{proof}

We are now ready to compute the character values $\chi_{M}\left(  r\right)  $,
$\chi_{N}\left(  r\right)  $ and $\chi_{M/N}\left(  r\right)  $:

\begin{itemize}
\item Applying (\ref{pf.thm.char.M/N.b.chi=}) to $V=M$, we obtain $\chi
_{M}\left(  r\right)  =\operatorname*{Tr}\left(  \tau_{r,M}\right)  $. But the
matrix that represents the $\mathbf{k}$-linear map $\tau_{r,M}:M\rightarrow M$
with respect to the basis $\left(  n_{1},n_{2},\ldots,n_{p+q}\right)  $ of $M$
is the matrix $\left(  a_{i,j}\right)  _{i,j\in\left[  p+q\right]  }$. Hence,%
\[
\operatorname*{Tr}\left(  \tau_{r,M}\right)  =\operatorname*{Tr}\left(
\left(  a_{i,j}\right)  _{i,j\in\left[  p+q\right]  }\right)  =\sum
_{i=1}^{p+q}a_{i,i}=\sum_{i=1}^{p}a_{i,i}+\sum_{i=p+1}^{p+q}a_{i,i}.
\]
In view of $\chi_{M}\left(  r\right)  =\operatorname*{Tr}\left(  \tau
_{r,M}\right)  $, we can rewrite this as%
\begin{equation}
\chi_{M}\left(  r\right)  =\sum_{i=1}^{p}a_{i,i}+\sum_{i=p+1}^{p+q}a_{i,i}.
\label{pf.thm.char.M/N.b.chiMr=}%
\end{equation}

\item Applying (\ref{pf.thm.char.M/N.b.chi=}) to $V=N$, we obtain $\chi
_{N}\left(  r\right)  =\operatorname*{Tr}\left(  \tau_{r,N}\right)  $. But the
matrix that represents the $\mathbf{k}$-linear map $\tau_{r,N}:N\rightarrow N$
with respect to the basis $\left(  n_{1},n_{2},\ldots,n_{p}\right)  $ of $N$
is the matrix $\left(  a_{i,j}\right)  _{i,j\in\left[  p\right]  }$, since
each $j\in\left[  p\right]  $ satisfies%
\[
\tau_{r,N}\left(  n_{j}\right)  =\sum_{i=1}^{p}a_{i,j}n_{i}%
\ \ \ \ \ \ \ \ \ \ \left(  \text{by Claim 3}\right)  .
\]
Hence,%
\[
\operatorname*{Tr}\left(  \tau_{r,N}\right)  =\operatorname*{Tr}\left(
\left(  a_{i,j}\right)  _{i,j\in\left[  p\right]  }\right)  =\sum_{i=1}%
^{p}a_{i,i}.
\]
In view of $\chi_{N}\left(  r\right)  =\operatorname*{Tr}\left(  \tau
_{r,N}\right)  $, we can rewrite this as%
\begin{equation}
\chi_{N}\left(  r\right)  =\sum_{i=1}^{p}a_{i,i}.
\label{pf.thm.char.M/N.b.chiNr=}%
\end{equation}

\item Applying (\ref{pf.thm.char.M/N.b.chi=}) to $V=M/N$, we obtain
$\chi_{M/N}\left(  r\right)  =\operatorname*{Tr}\left(  \tau_{r,M/N}\right)
$. But the matrix that represents the $\mathbf{k}$-linear map $\tau
_{r,M/N}:M/N\rightarrow M/N$ with respect to the basis $\left(  c_{1}%
,c_{2},\ldots,c_{q}\right)  $ of $M/N$ is the matrix $\left(  a_{i+p,j+p}%
\right)  _{i,j\in\left[  q\right]  }$, since each $j\in\left[  q\right]  $
satisfies%
\[
\tau_{r,M/N}\left(  c_{j}\right)  =\sum_{i=1}^{q}a_{i+p,j+p}c_{i}%
\ \ \ \ \ \ \ \ \ \ \left(  \text{by Claim 4}\right)  .
\]
Hence,%
\[
\operatorname*{Tr}\left(  \tau_{r,M/N}\right)  =\operatorname*{Tr}\left(
\left(  a_{i+p,j+p}\right)  _{i,j\in\left[  q\right]  }\right)  =\sum
_{i=1}^{q}a_{i+p,i+p}=\sum_{i=p+1}^{p+q}a_{i,i}%
\]
(here, we have substituted $i$ for $i+p$ in the sum). In view of $\chi
_{M/N}\left(  r\right)  =\operatorname*{Tr}\left(  \tau_{r,M/N}\right)  $, we
can rewrite this as%
\begin{equation}
\chi_{M/N}\left(  r\right)  =\sum_{i=p+1}^{p+q}a_{i,i}.
\label{pf.thm.char.M/N.b.chiM/Nr=}%
\end{equation}

\end{itemize}

Now,
\begin{align*}
\left(  \chi_{M}-\chi_{N}\right)  \left(  r\right)   &  =\chi_{M}\left(
r\right)  -\chi_{N}\left(  r\right) \\
&  =\left(  \sum_{i=1}^{p}a_{i,i}+\sum_{i=p+1}^{p+q}a_{i,i}\right)  -\left(
\sum_{i=1}^{p}a_{i,i}\right)  \ \ \ \ \ \ \ \ \ \ \left(  \text{by
(\ref{pf.thm.char.M/N.b.chiMr=}) and (\ref{pf.thm.char.M/N.b.chiNr=})}\right)
\\
&  =\sum_{i=p+1}^{p+q}a_{i,i}=\chi_{M/N}\left(  r\right)
\ \ \ \ \ \ \ \ \ \ \left(  \text{by (\ref{pf.thm.char.M/N.b.chiM/Nr=}%
)}\right)  .
\end{align*}

Forget that we fixed $r$. We thus have shown that $\left(  \chi_{M}-\chi
_{N}\right)  \left(  r\right)  =\chi_{M/N}\left(  r\right)  $ for each $r\in
R$. In other words, $\chi_{M}-\chi_{N}=\chi_{M/N}$. This proves Theorem
\ref{thm.char.M/N} \textbf{(b)}.
\end{proof}

As an example for the use of Theorem \ref{thm.char.M/N}, we can now easily
compute the character of the zero-sum representation $R\left(  \mathbf{k}%
^{n}\right)  $ of $S_{n}$ defined in Subsection \ref{subsec.rep.G-rep.subreps}:

\begin{example}
\label{exa.char.Rkn}Let $n\geq1$. Consider the natural representation
$\mathbf{k}^{n}$ of the symmetric group $S_{n}$ (defined in Example
\ref{exa.rep.Sn-rep.nat}), and its two subrepresentations $D\left(
\mathbf{k}^{n}\right)  $ and $R\left(  \mathbf{k}^{n}\right)  $ (defined in
Subsection \ref{subsec.rep.G-rep.subreps}). Then, for each $g\in S_{n}$, we
have%
\[
\chi_{D\left(  \mathbf{k}^{n}\right)  }\left(  g\right)
=1\ \ \ \ \ \ \ \ \ \ \text{and}\ \ \ \ \ \ \ \ \ \ \chi_{R\left(
\mathbf{k}^{n}\right)  }\left(  g\right)  =\left(  \text{\# of fixed points of
}g\right)  -1.
\]

\end{example}

\begin{proof}
Consider the trivial representation $\mathbf{k}_{\operatorname*{triv}}$ of
$S_{n}$ (defined in Example \ref{exa.rep.Sn-rep.triv}). Let us denote its
character $\chi_{\mathbf{k}_{\operatorname*{triv}}}$ by $\chi
_{\operatorname*{triv}}$. Thus, $\chi_{\mathbf{k}_{\operatorname*{triv}}}%
=\chi_{\operatorname*{triv}}$.

Let us furthermore denote the character $\chi_{\mathbf{k}^{n}}$ of the natural
representation $\mathbf{k}^{n}$ of $S_{n}$ by $\chi_{\operatorname*{nat}}$.
Thus, $\chi_{\mathbf{k}^{n}}=\chi_{\operatorname*{nat}}$.

Let $g\in S_{n}$. We already know (from Subsection
\ref{subsec.rep.G-rep.subreps}) that as an $S_{n}$-representation, $D\left(
\mathbf{k}^{n}\right)  $ is isomorphic to the trivial representation
$\mathbf{k}_{\operatorname*{triv}}$. Thus, Proposition \ref{prop.char.iso}
(applied to $R=\mathbf{k}\left[  S_{n}\right]  $ and $V=D\left(
\mathbf{k}^{n}\right)  $ and $W=\mathbf{k}_{\operatorname*{triv}}$) shows that
$\chi_{D\left(  \mathbf{k}^{n}\right)  }=\chi_{\mathbf{k}%
_{\operatorname*{triv}}}=\chi_{\operatorname*{triv}}$. But Example
\ref{exa.char.triv} \textbf{(a)} (applied to $G=S_{n}$) yields that
$\chi_{\operatorname*{triv}}\left(  g\right)  =1$. In view of $\chi_{D\left(
\mathbf{k}^{n}\right)  }=\chi_{\operatorname*{triv}}$, we can rewrite this as
$\chi_{D\left(  \mathbf{k}^{n}\right)  }\left(  g\right)  =1$.

It remains to show that $\chi_{R\left(  \mathbf{k}^{n}\right)  }\left(
g\right)  =\left(  \text{\# of fixed points of }g\right)  -1$. For this, we
recall that Example \ref{exa.char.nat} yields
\[
\chi_{\operatorname*{nat}}\left(  g\right)  =\left(  \text{\# of fixed points
of }g\right)  .
\]
But Theorem \ref{thm.rep.G-rep.Sn-nat.quots} \textbf{(d)} yields that
$\mathbf{k}^{n}/R\left(  \mathbf{k}^{n}\right)  \cong D\left(  \mathbf{k}%
^{n}\right)  $ as $S_{n}$-representations. Hence, $\mathbf{k}^{n}/R\left(
\mathbf{k}^{n}\right)  $ is a free $\mathbf{k}$-module of finite rank (since
$D\left(  \mathbf{k}^{n}\right)  $ is a free $\mathbf{k}$-module of finite
rank). Thus, Proposition \ref{prop.char.iso} (applied to $R=\mathbf{k}\left[
S_{n}\right]  $ and $V=\mathbf{k}^{n}/R\left(  \mathbf{k}^{n}\right)  $ and
$W=D\left(  \mathbf{k}^{n}\right)  $) shows that
\begin{equation}
\chi_{\mathbf{k}^{n}/R\left(  \mathbf{k}^{n}\right)  }=\chi_{D\left(
\mathbf{k}^{n}\right)  } \label{pf.exa.char.Rkn.5}%
\end{equation}
(since $\mathbf{k}^{n}/R\left(  \mathbf{k}^{n}\right)  \cong D\left(
\mathbf{k}^{n}\right)  $ as $S_{n}$-representations).

But Theorem \ref{thm.char.M/N} \textbf{(b)} (applied to $R=\mathbf{k}\left[
S_{n}\right]  $ and $M=\mathbf{k}^{n}$ and $N=R\left(  \mathbf{k}^{n}\right)
$) yields $\chi_{\mathbf{k}^{n}/R\left(  \mathbf{k}^{n}\right)  }%
=\chi_{\mathbf{k}^{n}}-\chi_{R\left(  \mathbf{k}^{n}\right)  }$. In view of
(\ref{pf.exa.char.Rkn.5}), we can rewrite this as $\chi_{D\left(
\mathbf{k}^{n}\right)  }=\chi_{\mathbf{k}^{n}}-\chi_{R\left(  \mathbf{k}%
^{n}\right)  }$. Solving this for $\chi_{R\left(  \mathbf{k}^{n}\right)  }$,
we obtain
\[
\chi_{R\left(  \mathbf{k}^{n}\right)  }=\underbrace{\chi_{\mathbf{k}^{n}}%
}_{=\chi_{\operatorname*{nat}}}-\underbrace{\chi_{D\left(  \mathbf{k}%
^{n}\right)  }}_{=\chi_{\operatorname*{triv}}}=\chi_{\operatorname*{nat}}%
-\chi_{\operatorname*{triv}}.
\]
Thus,%
\begin{align*}
\chi_{R\left(  \mathbf{k}^{n}\right)  }\left(  g\right)   &  =\left(
\chi_{\operatorname*{nat}}-\chi_{\operatorname*{triv}}\right)  \left(
g\right)  =\underbrace{\chi_{\operatorname*{nat}}\left(  g\right)  }_{=\left(
\text{\# of fixed points of }g\right)  }-\underbrace{\chi
_{\operatorname*{triv}}\left(  g\right)  }_{=1}\\
&  =\left(  \text{\# of fixed points of }g\right)  -1.
\end{align*}
This completes the proof of Example \ref{exa.char.Rkn}.
\end{proof}

Another consequence of Theorem \ref{thm.char.M/N} is the following formula for
characters of direct sums:

\begin{corollary}
\label{cor.char.V+W}Let $R$ be a $\mathbf{k}$-algebra. Let $V$ and $W$ be two
left $R$-modules. Assume that, as $\mathbf{k}$-modules, both $V$ and $W$ are
free of finite ranks (not necessarily of the same rank). Then,
\[
\chi_{V\oplus W}=\chi_{V}+\chi_{W}.
\]
(Here, the addition is pointwise, as usual for $\mathbf{k}$-linear maps from
$R$ to $\mathbf{k}$.)
\end{corollary}

\begin{fineprint}
\begin{proof}
It is well-known that the $\mathbf{k}$-module $V\oplus W$ is free of finite
rank (since $V$ and $W$ are).

Consider the left $R$-submodule%
\[
V\oplus0:=\left\{  \left(  v,0\right)  \ \mid\ v\in V\right\}
\]
of $V\oplus W$. It is well-known that this $R$-submodule $V\oplus0$ is
isomorphic to $V$ (via the isomorphism $V\oplus0\rightarrow V$ that sends each
$\left(  v,0\right)  $ to $v$). Thus, the $\mathbf{k}$-module $V\oplus0$ is
free of finite rank (since $V$ is). Hence, Proposition \ref{prop.char.iso}
(applied to $V\oplus0$ and $V$ instead of $V$ and $W$) yields $\chi_{V\oplus
0}=\chi_{V}$ (since $V\oplus0\cong V$ as left $R$-modules).

But it is also well-known that the quotient $R$-module $\left(  V\oplus
W\right)  /\left(  V\oplus0\right)  $ is isomorphic to the $R$-module $W$ (via
the isomorphism $\left(  V\oplus W\right)  /\left(  V\oplus0\right)
\rightarrow W$ that sends each $\overline{\left(  v,w\right)  }\in\left(
V\oplus W\right)  /\left(  V\oplus0\right)  $ to $w$). Thus, the $\mathbf{k}%
$-module $\left(  V\oplus W\right)  /\left(  V\oplus0\right)  $ is free of
finite rank (since $W$ is). Hence, Proposition \ref{prop.char.iso} (applied to
$\left(  V\oplus W\right)  /\left(  V\oplus0\right)  $ instead of $V$) yields
$\chi_{\left(  V\oplus W\right)  /\left(  V\oplus0\right)  }=\chi_{W}$ (since
$\left(  V\oplus W\right)  /\left(  V\oplus0\right)  \cong W$ as left $R$-modules).

Now, Theorem \ref{thm.char.M/N} \textbf{(b)} (applied to $M=V\oplus W$ and
$N=V\oplus0$) yields $\chi_{\left(  V\oplus W\right)  /\left(  V\oplus
0\right)  }=\chi_{V\oplus W}-\chi_{V\oplus0}$. In view of $\chi_{\left(
V\oplus W\right)  /\left(  V\oplus0\right)  }=\chi_{W}$ and $\chi_{V\oplus
0}=\chi_{V}$, we can rewrite this as $\chi_{W}=\chi_{V\oplus W}-\chi_{V}$.
Thus, $\chi_{V\oplus W}=\chi_{V}+\chi_{W}$. This proves Corollary
\ref{cor.char.V+W}.
\end{proof}
\end{fineprint}

We can easily extend Corollary \ref{cor.char.V+W} to direct sums of multiple
$R$-modules:

\begin{corollary}
\label{cor.char.V+W.m}Let $R$ be a $\mathbf{k}$-algebra. Let $V_{1}%
,V_{2},\ldots,V_{m}$ be finitely many left $R$-modules. Assume that, as a
$\mathbf{k}$-module, each $V_{i}$ is free of finite rank. Then,%
\[
\chi_{V_{1}\oplus V_{2}\oplus\cdots\oplus V_{m}}=\chi_{V_{1}}+\chi_{V_{2}%
}+\cdots+\chi_{V_{m}}.
\]
(Here, the addition is pointwise, as usual for $\mathbf{k}$-linear maps from
$R$ to $\mathbf{k}$.)
\end{corollary}

\begin{fineprint}
\begin{proof}
[Proof sketch.]Induct on $m$. The \textit{base case} $m=0$ follows from the
easy fact that the character of the $R$-module $\left\{  0\right\}  $ is $0$.
The \textit{induction step} uses Corollary \ref{cor.char.V+W} and the
observation that $V_{1}\oplus V_{2}\oplus\cdots\oplus V_{m}\cong\left(
V_{1}\oplus V_{2}\oplus\cdots\oplus V_{m-1}\right)  \oplus V_{m}$.
\end{proof}
\end{fineprint}

\subsubsection{The quasiidempotent character formula}

We shall now state an important formula for the characters of certain
representations of a finite group $G$: specifically, for those representations
that can be written as $\mathcal{A}\mathbf{a}$ for some quasi-idempotent
element $\mathbf{a}\in\mathbf{k}\left[  G\right]  $. This class of
representations includes the Specht modules $\mathcal{S}^{\lambda}$ for
partitions $\lambda$ of $n$, since Theorem \ref{thm.specht.ETidp} shows that
the respective Young symmetrizers $\mathbf{E}_{T}$ are quasi-idempotent.

Before we can state our formula, we recall Definition \ref{def.specht.ET.wa}.
The notation $\left[  w\right]  \mathbf{a}$ introduced in this definition has
the following simple but important property:

\begin{proposition}
\label{prop.specht.ET.wa.uva}Let $G$ be a group. Let $u,v\in G$ and
$\mathbf{b}\in\mathbf{k}\left[  G\right]  $. Then,
\begin{align}
\left[  uv\right]  \mathbf{b}  &  =\left[  v\right]  \left(  u^{-1}%
\mathbf{b}\right) \label{eq.prop.specht.ET.wa.uva.1}\\
&  =\left[  u\right]  \left(  \mathbf{b}v^{-1}\right)  .
\label{eq.prop.specht.ET.wa.uva.2}%
\end{align}

\end{proposition}

\begin{proof}
Recall that $G$ is a group. Thus, the map $G\rightarrow G,\ g\mapsto ug$ is a
bijection (and its inverse is the map $G\rightarrow G,\ g\mapsto u^{-1}g$).
For similar reasons, the map $G\rightarrow G,\ g\mapsto gv$ is a bijection.

Write the element $\mathbf{b}\in\mathbf{k}\left[  G\right]  $ as a
$\mathbf{k}$-linear combination of the standard basis vectors $g\in G$ as
follows:%
\begin{equation}
\mathbf{b}=\sum_{g\in G}\beta_{g}g, \label{pf.thm.char.ladisch.c2.pf.1}%
\end{equation}
where the coefficients $\beta_{g}$ belong to $\mathbf{k}$. Then, $\left[
g\right]  \mathbf{b}=\beta_{g}$ for each $g\in G$ (by Definition
\ref{def.specht.ET.wa}). Applying this to $g=uv$, we obtain $\left[
uv\right]  \mathbf{b}=\beta_{uv}$.

On the other hand, multiplying the equality (\ref{pf.thm.char.ladisch.c2.pf.1}%
) by $u^{-1}$ from the left, we obtain
\begin{align*}
u^{-1}\mathbf{b}  &  =u^{-1}\sum_{g\in G}\beta_{g}g=\sum_{g\in G}\beta
_{g}u^{-1}g=\sum_{g\in G}\beta_{ug}\underbrace{u^{-1}u}_{=1}g\\
&  \ \ \ \ \ \ \ \ \ \ \ \ \ \ \ \ \ \ \ \ \left(
\begin{array}
[c]{c}%
\text{here, we have substituted }ug\text{ for }g\text{ in the sum, since}\\
\text{the map }G\rightarrow G,\ g\mapsto ug\text{ is a bijection}%
\end{array}
\right) \\
&  =\sum_{g\in G}\beta_{ug}g.
\end{align*}
Hence, $\left[  g\right]  \left(  u^{-1}\mathbf{b}\right)  =\beta_{ug}$ for
each $g\in G$ (by Definition \ref{def.specht.ET.wa}). Applying this to $g=v$,
we obtain $\left[  v\right]  \left(  u^{-1}\mathbf{b}\right)  =\beta_{uv}$.
Comparing this with $\left[  uv\right]  \mathbf{b}=\beta_{uv}$, we find
$\left[  uv\right]  \mathbf{b}=\left[  v\right]  \left(  u^{-1}\mathbf{b}%
\right)  $. This proves (\ref{eq.prop.specht.ET.wa.uva.1}).

Furthermore, multiplying the equality (\ref{pf.thm.char.ladisch.c2.pf.1}) by
$v^{-1}$ from the right, we obtain
\begin{align*}
\mathbf{b}v^{-1}  &  =\left(  \sum_{g\in G}\beta_{g}g\right)  v^{-1}%
=\sum_{g\in G}\beta_{g}gv^{-1}=\sum_{g\in G}\beta_{gv}g\underbrace{vv^{-1}%
}_{=1}\\
&  \ \ \ \ \ \ \ \ \ \ \ \ \ \ \ \ \ \ \ \ \left(
\begin{array}
[c]{c}%
\text{here, we have substituted }gv\text{ for }g\text{ in the sum, since}\\
\text{the map }G\rightarrow G,\ g\mapsto gv\text{ is a bijection}%
\end{array}
\right) \\
&  =\sum_{g\in G}\beta_{gv}g.
\end{align*}
Hence, $\left[  g\right]  \left(  \mathbf{b}v^{-1}\right)  =\beta_{gv}$ for
each $g\in G$ (by Definition \ref{def.specht.ET.wa}). Applying this to $g=u$,
we obtain $\left[  u\right]  \left(  \mathbf{b}v^{-1}\right)  =\beta_{uv}$.
Comparing this with $\left[  uv\right]  \mathbf{b}=\beta_{uv}$, we find
$\left[  uv\right]  \mathbf{b}=\left[  u\right]  \left(  \mathbf{b}%
v^{-1}\right)  $. This proves (\ref{eq.prop.specht.ET.wa.uva.2}).

Since we have proved both (\ref{eq.prop.specht.ET.wa.uva.1}) and
(\ref{eq.prop.specht.ET.wa.uva.2}), we thus know that Proposition
\ref{prop.specht.ET.wa.uva} holds.
\end{proof}

We are now ready to state our formula for characters:

\begin{theorem}
\label{thm.char.ladisch}Let $G$ be a finite group. Let $\mathcal{A}%
=\mathbf{k}\left[  G\right]  $. (We do \textbf{not} follow the convention that
$\mathcal{A}=\mathbf{k}\left[  S_{n}\right]  $ in this theorem.) Let
$\mathbf{a}\in\mathcal{A}$ and $\kappa\in\mathbf{k}$ be such that
$\mathbf{a}^{2}=\kappa\mathbf{a}$. Assume that the $\mathbf{k}$-module
$\mathcal{A}\mathbf{a}$ is free of finite rank. Then, the character
$\chi_{\mathcal{A}\mathbf{a}}$ of the left $\mathbf{k}\left[  G\right]
$-module $\mathcal{A}\mathbf{a}$ satisfies%
\[
\kappa\cdot\chi_{\mathcal{A}\mathbf{a}}\left(  h\right)  =\sum_{w\in G}\left[
h^{-1}\right]  \left(  w\mathbf{a}w^{-1}\right)  \ \ \ \ \ \ \ \ \ \ \text{for
each }h\in G.
\]

\end{theorem}

\begin{proof}
The following proof is somewhat similar to the proof of Lemma
\ref{lem.linalg.ladisch} above.

First, we observe that $\mathcal{A}\mathbf{a}$ is a left $\mathcal{A}%
$-submodule of $\mathcal{A}$, thus a left $\mathcal{A}$-module, i.e., a left
$\mathbf{k}\left[  G\right]  $-module (since $\mathcal{A}=\mathbf{k}\left[
G\right]  $). Thus, its character $\chi_{\mathcal{A}\mathbf{a}}$ is
well-defined (since the $\mathbf{k}$-module $\mathcal{A}\mathbf{a}$ is free of
finite rank).

We will need the following two facts:

\begin{statement}
\textit{Claim 1:} Let $f$ be any $\mathbf{k}$-module endomorphism of
$\mathbf{k}\left[  G\right]  $. Then,%
\[
\operatorname*{Tr}f=\sum_{w\in G}\left[  w\right]  \left(  f\left(  w\right)
\right)  .
\]

\end{statement}

\begin{proof}
[Proof of Claim 1.]This is just Claim 1 in the proof of Proposition
\ref{prop.groupalg.trace}, except that we have renamed the summation index $g$
as $w$.
\end{proof}

\begin{statement}
\textit{Claim 2:} Let $h,w\in G$. Then,
\[
\left[  w\right]  \left(  hw\mathbf{a}\right)  =\left[  h^{-1}\right]  \left(
w\mathbf{a}w^{-1}\right)  .
\]

\end{statement}

\begin{proof}
[Proof of Claim 2.]Applying (\ref{eq.prop.specht.ET.wa.uva.2}) to $u=1_{G}$
and $v=w$ and $\mathbf{b}=hw\mathbf{a}$, we obtain $\left[  1_{G}w\right]
\left(  hw\mathbf{a}\right)  =\left[  1_{G}\right]  \left(  hw\mathbf{a}%
w^{-1}\right)  $. In other words, $\left[  w\right]  \left(  hw\mathbf{a}%
\right)  =\left[  1_{G}\right]  \left(  hw\mathbf{a}w^{-1}\right)  $ (since
$1_{G}w=w$).

On the other hand, applying (\ref{eq.prop.specht.ET.wa.uva.1}) to $u=h$ and
$v=h^{-1}$ and $\mathbf{b}=hw\mathbf{a}w^{-1}$, we obtain $\left[
hh^{-1}\right]  \left(  hw\mathbf{a}w^{-1}\right)  =\left[  h^{-1}\right]
\left(  \underbrace{h^{-1}h}_{=1}w\mathbf{a}w^{-1}\right)  =\left[
h^{-1}\right]  \left(  w\mathbf{a}w^{-1}\right)  $. In view of $hh^{-1}=1_{G}%
$, we can rewrite this as $\left[  1_{G}\right]  \left(  hw\mathbf{a}%
w^{-1}\right)  =\left[  h^{-1}\right]  \left(  w\mathbf{a}w^{-1}\right)  $.

Altogether, we now have $\left[  w\right]  \left(  hw\mathbf{a}\right)
=\left[  1_{G}\right]  \left(  hw\mathbf{a}w^{-1}\right)  =\left[
h^{-1}\right]  \left(  w\mathbf{a}w^{-1}\right)  $. This proves Claim 2.
\end{proof}

Now, fix $h\in G$. Let $f:\mathcal{A}\rightarrow\mathcal{A}\mathbf{a}$ be the
map that sends each $\mathbf{b}\in\mathcal{A}$ to $h\mathbf{ba}$. (This is
well-defined, since each $\mathbf{b}\in\mathcal{A}$ satisfies
$\underbrace{h\mathbf{b}}_{\in\mathcal{A}}\mathbf{a}\in\mathcal{A}\mathbf{a}%
$.) Let $g:\mathcal{A}\mathbf{a}\rightarrow\mathcal{A}$ be the inclusion map
(i.e., the map that sends each $\mathbf{b}\in\mathcal{A}\mathbf{a}$ to
$\mathbf{b}$ itself). Both maps $f$ and $g$ are clearly $\mathbf{k}$-linear.
Hence, Lemma \ref{lem.linalg.Trfg} (applied to $V=\mathcal{A}$ and
$W=\mathcal{A}\mathbf{a}$) yields
\begin{equation}
\operatorname*{Tr}\left(  f\circ g\right)  =\operatorname*{Tr}\left(  g\circ
f\right)  . \label{pf.thm.char.ladisch.2}%
\end{equation}

However, the map $g\circ f:\mathcal{A}\rightarrow\mathcal{A}$ is $\mathbf{k}%
$-linear (since $g$ and $f$ are), and thus is a $\mathbf{k}$-module
endomorphism of $\mathcal{A}=\mathbf{k}\left[  G\right]  $. Hence, Claim 1
(applied to $g\circ f$ instead of $f$) yields%
\begin{equation}
\operatorname*{Tr}\left(  g\circ f\right)  =\sum_{w\in G}\left[  w\right]
\left(  \left(  g\circ f\right)  \left(  w\right)  \right)  .
\label{pf.thm.char.ladisch.3}%
\end{equation}
However, each $w\in G$ satisfies%
\begin{align*}
\left(  g\circ f\right)  \left(  w\right)   &  =g\left(  f\left(  w\right)
\right)  =f\left(  w\right)  \ \ \ \ \ \ \ \ \ \ \left(  \text{by the
definition of }g\right) \\
&  =hw\mathbf{a}\ \ \ \ \ \ \ \ \ \ \left(  \text{by the definition of
}f\right)
\end{align*}
and therefore%
\begin{equation}
\left[  w\right]  \left(  \underbrace{\left(  g\circ f\right)  \left(
w\right)  }_{=hw\mathbf{a}}\right)  =\left[  w\right]  \left(  hw\mathbf{a}%
\right)  =\left[  h^{-1}\right]  \left(  w\mathbf{a}w^{-1}\right)
\label{pf.thm.char.ladisch.4}%
\end{equation}
(by Claim 2). Hence, (\ref{pf.thm.char.ladisch.3}) becomes%
\begin{align}
\operatorname*{Tr}\left(  g\circ f\right)   &  =\sum_{w\in G}%
\underbrace{\left[  w\right]  \left(  \left(  g\circ f\right)  \left(
w\right)  \right)  }_{\substack{=\left[  h^{-1}\right]  \left(  w\mathbf{a}%
w^{-1}\right)  \\\text{(by (\ref{pf.thm.char.ladisch.4}))}}}\nonumber\\
&  =\sum_{w\in G}\left[  h^{-1}\right]  \left(  w\mathbf{a}w^{-1}\right)  .
\label{pf.thm.char.ladisch.5}%
\end{align}

On the other hand, for each $r\in\mathcal{A}$, we let $\tau_{r}:\mathcal{A}%
\mathbf{a}\rightarrow\mathcal{A}\mathbf{a}$ be the $\mathbf{k}$-linear map
$v\mapsto rv$. Then, the definition of the character $\chi_{\mathcal{A}%
\mathbf{a}}$ yields%
\[
\chi_{\mathcal{A}\mathbf{a}}\left(  r\right)  =\operatorname*{Tr}\left(
\tau_{r}\right)  \ \ \ \ \ \ \ \ \ \ \text{for each }r\in\mathcal{A}.
\]
Applying this to $r=\kappa h$, we obtain%
\begin{equation}
\chi_{\mathcal{A}\mathbf{a}}\left(  \kappa h\right)  =\operatorname*{Tr}%
\left(  \tau_{\kappa h}\right)  . \label{pf.thm.char.ladisch.6b}%
\end{equation}

Next, we claim that%
\begin{equation}
\tau_{\kappa h}=f\circ g. \label{pf.thm.char.ladisch.7}%
\end{equation}

\begin{proof}
[Proof of (\ref{pf.thm.char.ladisch.7}):]Let $\mathbf{b}\in\mathcal{A}%
\mathbf{a}$. Thus, $\mathbf{b}=\mathbf{ca}$ for some $\mathbf{c}\in
\mathcal{A}$. Consider this $\mathbf{c}$. Then, the definition of
$\tau_{\kappa h}$ yields $\tau_{\kappa h}\left(  \mathbf{b}\right)  =\kappa
h\underbrace{\mathbf{b}}_{=\mathbf{ca}}=\kappa h\mathbf{ca}$. On the other
hand, the definition of $g$ yields $g\left(  \mathbf{b}\right)  =\mathbf{b}$.
Applying the map $f$ to this equality, we find $f\left(  g\left(
\mathbf{b}\right)  \right)  =f\left(  \mathbf{b}\right)  =h\mathbf{ba}$ (by
the definition of $f$). Hence,%
\[
\left(  f\circ g\right)  \left(  \mathbf{b}\right)  =f\left(  g\left(
\mathbf{b}\right)  \right)  =h\underbrace{\mathbf{b}}_{=\mathbf{ca}}%
\mathbf{a}=h\mathbf{c}\underbrace{\mathbf{aa}}_{=\mathbf{a}^{2}=\kappa
\mathbf{a}}=h\mathbf{c}\cdot\kappa\mathbf{a}=\kappa h\mathbf{ca}.
\]
Comparing this with $\tau_{\kappa h}\left(  \mathbf{b}\right)  =\kappa
h\mathbf{ca}$, we obtain $\tau_{\kappa h}\left(  \mathbf{b}\right)  =\left(
f\circ g\right)  \left(  \mathbf{b}\right)  $.

Forget that we fixed $\mathbf{b}$. We thus have shown that $\tau_{\kappa
h}\left(  \mathbf{b}\right)  =\left(  f\circ g\right)  \left(  \mathbf{b}%
\right)  $ for each $\mathbf{b}\in\mathcal{A}\mathbf{a}$. In other words,
$\tau_{\kappa h}=f\circ g$. This proves (\ref{pf.thm.char.ladisch.7}).
\end{proof}

Now, (\ref{pf.thm.char.ladisch.2}) yields%
\begin{align*}
\operatorname*{Tr}\left(  g\circ f\right)   &  =\operatorname*{Tr}%
\underbrace{\left(  f\circ g\right)  }_{\substack{=\tau_{\kappa h}\\\text{(by
(\ref{pf.thm.char.ladisch.7}))}}}=\operatorname*{Tr}\left(  \tau_{\kappa
h}\right)  =\chi_{\mathcal{A}\mathbf{a}}\left(  \kappa h\right)
\ \ \ \ \ \ \ \ \ \ \left(  \text{by (\ref{pf.thm.char.ladisch.6b})}\right) \\
&  =\kappa\cdot\chi_{\mathcal{A}\mathbf{a}}\left(  h\right)
\end{align*}
(since Proposition \ref{prop.char.klin} shows that $\chi_{\mathcal{A}%
\mathbf{a}}$ is a $\mathbf{k}$-linear map). Comparing this with
(\ref{pf.thm.char.ladisch.5}), we obtain
\[
\kappa\cdot\chi_{\mathcal{A}\mathbf{a}}\left(  h\right)  =\sum_{w\in G}\left[
h^{-1}\right]  \left(  w\mathbf{a}w^{-1}\right)  .
\]
This proves Theorem \ref{thm.char.ladisch}.
\end{proof}

\begin{corollary}
\label{cor.char.ladisch-X}Let $G$ be a finite group. Define the $\mathbf{k}%
$-linear map%
\begin{align*}
S:\mathbf{k}\left[  G\right]   &  \rightarrow\mathbf{k}\left[  G\right]  ,\\
w  &  \mapsto w^{-1}\ \ \ \ \ \ \ \ \ \ \text{for all }w\in G.
\end{align*}
This is understood as in Definition \ref{def.Tsign.Tsign}. (Note that this map
$S$ generalizes the antipode $S$ of $\mathbf{k}\left[  S_{n}\right]  $, and
was already mentioned in Remark \ref{rmk.S.all-G}.)

Let $\mathcal{A}=\mathbf{k}\left[  G\right]  $. (We do \textbf{not} follow the
convention that $\mathcal{A}=\mathbf{k}\left[  S_{n}\right]  $ in this
theorem.) Let $\mathbf{a}\in\mathcal{A}$ and $\kappa\in\mathbf{k}$ be such
that $\mathbf{a}^{2}=\kappa\mathbf{a}$. Assume that the $\mathbf{k}$-module
$\mathcal{A}\mathbf{a}$ is free of finite rank. Then,%
\[
\kappa\cdot S\left(  \mathbf{X}_{\mathcal{A}\mathbf{a}}\right)  =\sum_{w\in
G}w\mathbf{a}w^{-1}.
\]
(See Definition \ref{def.char.XV} for the meaning of $\mathbf{X}%
_{\mathcal{A}\mathbf{a}}$.)
\end{corollary}

\begin{proof}
First, we observe that $\mathcal{A}\mathbf{a}$ is a left $\mathcal{A}%
$-submodule of $\mathcal{A}$, thus a left $\mathcal{A}$-module, i.e., a left
$\mathbf{k}\left[  G\right]  $-module (since $\mathcal{A}=\mathbf{k}\left[
G\right]  $). Thus, $\mathbf{X}_{\mathcal{A}\mathbf{a}}$ is well-defined
(since the $\mathbf{k}$-module $\mathcal{A}\mathbf{a}$ is free of finite rank).

Note that the map $G\rightarrow G,\ g\mapsto g^{-1}$ is a bijection (since it
is inverse to itself).

The definition of $\mathbf{X}_{\mathcal{A}\mathbf{a}}$ yields
\[
\mathbf{X}_{\mathcal{A}\mathbf{a}}=\sum_{g\in G}\chi_{\mathcal{A}\mathbf{a}%
}\left(  g\right)  g.
\]
Applying the map $S$ to this equality, we obtain%
\begin{align*}
S\left(  \mathbf{X}_{\mathcal{A}\mathbf{a}}\right)   &  =S\left(  \sum_{g\in
G}\chi_{\mathcal{A}\mathbf{a}}\left(  g\right)  g\right) \\
&  =\sum_{g\in G}\chi_{\mathcal{A}\mathbf{a}}\left(  g\right)
\underbrace{S\left(  g\right)  }_{\substack{=g^{-1}\\\text{(by the definition
of }S\text{)}}}\ \ \ \ \ \ \ \ \ \ \left(  \text{since the map }S\text{ is
}\mathbf{k}\text{-linear}\right) \\
&  =\sum_{g\in G}\chi_{\mathcal{A}\mathbf{a}}\left(  g\right)  g^{-1}.
\end{align*}
Multiplying both sides of this equality by $\kappa$, we find%
\begin{align*}
\kappa\cdot S\left(  \mathbf{X}_{\mathcal{A}\mathbf{a}}\right)   &
=\kappa\cdot\sum_{g\in G}\chi_{\mathcal{A}\mathbf{a}}\left(  g\right)
g^{-1}=\sum_{g\in G}\underbrace{\kappa\cdot\chi_{\mathcal{A}\mathbf{a}}\left(
g\right)  }_{\substack{=\sum_{w\in G}\left[  g^{-1}\right]  \left(
w\mathbf{a}w^{-1}\right)  \\\text{(by Theorem \ref{thm.char.ladisch}%
,}\\\text{applied to }h=g\text{)}}}g^{-1}\\
&  =\sum_{g\in G}\left(  \sum_{w\in G}\left[  g^{-1}\right]  \left(
w\mathbf{a}w^{-1}\right)  \right)  g^{-1}=\sum_{g\in G}\ \ \sum_{w\in
G}\left(  \left[  g^{-1}\right]  \left(  w\mathbf{a}w^{-1}\right)  \right)
\cdot g^{-1}\\
&  =\sum_{g\in G}\ \ \sum_{w\in G}\left(  \left[  g\right]  \left(
w\mathbf{a}w^{-1}\right)  \right)  \cdot g
\end{align*}
(here, we have substituted $g$ for $g^{-1}$ in the outer sum, since the map
$G\rightarrow G,\ g\mapsto g^{-1}$ is a bijection). Thus,%
\begin{align}
\kappa\cdot S\left(  \mathbf{X}_{\mathcal{A}\mathbf{a}}\right)   &
=\underbrace{\sum_{g\in G}\ \ \sum_{w\in G}}_{=\sum_{w\in G}\ \ \sum_{g\in G}%
}\left(  \left[  g\right]  \left(  w\mathbf{a}w^{-1}\right)  \right)  \cdot
g\nonumber\\
&  =\sum_{w\in G}\ \ \sum_{g\in G}\left(  \left[  g\right]  \left(
w\mathbf{a}w^{-1}\right)  \right)  \cdot g. \label{pf.thm.char.ladisch-X.4}%
\end{align}

However, it is easy to see that each $\mathbf{b}\in\mathcal{A}$ satisfies
\begin{equation}
\sum_{g\in G}\left(  \left[  g\right]  \mathbf{b}\right)  \cdot g=\mathbf{b}
\label{pf.thm.char.ladisch-X.5}%
\end{equation}
\footnote{\textit{Proof:} Let $\mathbf{b}\in\mathcal{A}$. Write the element
$\mathbf{b}\in\mathcal{A}=\mathbf{k}\left[  G\right]  $ as a $\mathbf{k}%
$-linear combination of the standard basis vectors $g\in G$ as follows:%
\[
\mathbf{b}=\sum_{g\in G}\beta_{g}g,
\]
where the coefficients $\beta_{g}$ belong to $\mathbf{k}$. Then, $\left[
g\right]  \mathbf{b}=\beta_{g}$ for each $g\in G$ (by Definition
\ref{def.specht.ET.wa}). Hence, $\sum_{g\in G}\underbrace{\left(  \left[
g\right]  \mathbf{b}\right)  }_{=\beta_{g}}\cdot\,g=\sum_{g\in G}\beta
_{g}g=\mathbf{b}$. This proves (\ref{pf.thm.char.ladisch-X.5}).}. Hence,
(\ref{pf.thm.char.ladisch-X.4}) becomes%
\[
\kappa\cdot S\left(  \mathbf{X}_{\mathcal{A}\mathbf{a}}\right)  =\sum_{w\in
G}\ \ \underbrace{\sum_{g\in G}\left(  \left[  g\right]  \left(
w\mathbf{a}w^{-1}\right)  \right)  \cdot g}_{\substack{=w\mathbf{a}%
w^{-1}\\\text{(by (\ref{pf.thm.char.ladisch-X.5}), applied to }\mathbf{b}%
=w\mathbf{a}w^{-1}\text{)}}}=\sum_{w\in G}w\mathbf{a}w^{-1}.
\]
This proves Corollary \ref{cor.char.ladisch-X}.
\end{proof}

When $G$ is the symmetric group $S_{n}$, we can simplify Corollary
\ref{cor.char.ladisch-X} by removing the map $S$:

\begin{corollary}
\label{cor.char.ladisch-XSn}Let $\mathcal{A}=\mathbf{k}\left[  S_{n}\right]
$. Let $\mathbf{a}\in\mathcal{A}$ and $\kappa\in\mathbf{k}$ be such that
$\mathbf{a}^{2}=\kappa\mathbf{a}$. Assume that the $\mathbf{k}$-module
$\mathcal{A}\mathbf{a}$ is free of finite rank. Then,%
\[
\kappa\cdot\mathbf{X}_{\mathcal{A}\mathbf{a}}=\sum_{w\in S_{n}}w\mathbf{a}%
w^{-1}.
\]
(See Definition \ref{def.char.XV} for the meaning of $\mathbf{X}%
_{\mathcal{A}\mathbf{a}}$.)
\end{corollary}

\begin{proof}
Let $G=S_{n}$. Then, the map $S$ defined in Corollary \ref{cor.char.ladisch-X}
is precisely the antipode $S$ of $\mathbf{k}\left[  S_{n}\right]  $. Hence,
Corollary \ref{cor.char.ladisch-X} yields%
\begin{equation}
\kappa\cdot S\left(  \mathbf{X}_{\mathcal{A}\mathbf{a}}\right)  =\sum_{w\in
G}w\mathbf{a}w^{-1}. \label{pf.cor.char.ladisch-XSn.1}%
\end{equation}

But Proposition \ref{prop.char.XV.central} (applied to $V=\mathcal{A}%
\mathbf{a}$) yields that $\mathbf{X}_{\mathcal{A}\mathbf{a}}$ belongs to the
center $Z\left(  \mathbf{k}\left[  S_{n}\right]  \right)  $ of the ring
$\mathbf{k}\left[  S_{n}\right]  $. In other words, $\mathbf{X}_{\mathcal{A}%
\mathbf{a}}\in Z\left(  \mathbf{k}\left[  S_{n}\right]  \right)  $. Hence,
Proposition \ref{prop.ZkSn.S=id} (applied to $\mathbf{X}_{\mathcal{A}%
\mathbf{a}}$ instead of $\mathbf{a}$) yields $S\left(  \mathbf{X}%
_{\mathcal{A}\mathbf{a}}\right)  =\mathbf{X}_{\mathcal{A}\mathbf{a}}$. Hence,
we can rewrite (\ref{pf.cor.char.ladisch-XSn.1}) as
\[
\kappa\cdot\mathbf{X}_{\mathcal{A}\mathbf{a}}=\sum_{w\in G}w\mathbf{a}%
w^{-1}=\sum_{w\in S_{n}}w\mathbf{a}w^{-1}\ \ \ \ \ \ \ \ \ \ \left(
\text{since }G=S_{n}\right)  .
\]
This proves Corollary \ref{cor.char.ladisch-XSn}.
\end{proof}

\subsubsection{Characters and centralized Young symmetrizers}

We are now ready to state the main theorem about characters of Specht modules:

\begin{theorem}
\label{thm.char.Elam-Xlam}Let $\lambda$ be a partition of $n$. Consider the
Specht module $\mathcal{S}^{\lambda}$, and let $\mathbf{X}_{\lambda}$ denote
the element $\mathbf{X}_{\mathcal{S}^{\lambda}}$ corresponding to it (defined
as in Definition \ref{def.char.XV}). Then,%
\[
\mathbf{E}_{\lambda}=h^{\lambda}\mathbf{X}_{\lambda}.
\]
(See Definition \ref{def.spechtmod.flam} \textbf{(b)} for the meaning of
$h^{\lambda}$, and Definition \ref{def.spechtmod.Elam.Elam} for the meaning of
$\mathbf{E}_{\lambda}$.)
\end{theorem}

\begin{proof}
Recall that $\mathcal{A}=\mathbf{k}\left[  S_{n}\right]  $. Pick any
$n$-tableau $T$ of shape $Y\left(  \lambda\right)  $ (clearly, such an
$n$-tableau exists, since $\left\vert Y\left(  \lambda\right)  \right\vert
=\left\vert \lambda\right\vert =n$). Then, (\ref{eq.def.specht.ET.defs.SD=AET}%
) (applied to $D=Y\left(  \lambda\right)  $) yields $\mathcal{S}^{Y\left(
\lambda\right)  }\cong\mathcal{A}\mathbf{E}_{T}$ (as left $\mathcal{A}%
$-modules). This is an isomorphism of left $\mathcal{A}$-modules, i.e., of
left $\mathbf{k}\left[  S_{n}\right]  $-modules (since $\mathcal{A}%
=\mathbf{k}\left[  S_{n}\right]  $), thus an isomorphism of representations of
$S_{n}$, and therefore also an isomorphism of $\mathbf{k}$-modules.

Moreover, $Y\left(  \lambda\right)  =Y\left(  \lambda/\varnothing\right)  $ is
a skew Young diagram. Hence, Theorem \ref{thm.spechtmod.basis} (applied to
$D=Y\left(  \lambda\right)  $) shows that the standard polytabloids form a
basis of the $\mathbf{k}$-module $\mathcal{S}^{Y\left(  \lambda\right)  }$.
Hence, $\mathcal{S}^{Y\left(  \lambda\right)  }$ is a free $\mathbf{k}$-module
of finite rank. Thus, $\mathcal{A}\mathbf{E}_{T}$ is a free $\mathbf{k}%
$-module of finite rank as well (since $\mathcal{S}^{Y\left(  \lambda\right)
}\cong\mathcal{A}\mathbf{E}_{T}$).

Hence, Corollary \ref{cor.char.XV.iso} (applied to $G=S_{n}$, $V=\mathcal{S}%
^{Y\left(  \lambda\right)  }$ and $W=\mathcal{A}\mathbf{E}_{T}$) yields
$\mathbf{X}_{\mathcal{S}^{Y\left(  \lambda\right)  }}=\mathbf{X}%
_{\mathcal{A}\mathbf{E}_{T}}$ (since $\mathcal{S}^{Y\left(  \lambda\right)
}\cong\mathcal{A}\mathbf{E}_{T}$ as representations of $S_{n}$). However, the
definition of $\mathbf{X}_{\lambda}$ shows that $\mathbf{X}_{\lambda
}=\mathbf{X}_{\mathcal{S}^{\lambda}}=\mathbf{X}_{\mathcal{S}^{Y\left(
\lambda\right)  }}$ (since $\mathcal{S}^{\lambda}=\mathcal{S}^{Y\left(
\lambda\right)  }$). Thus,%
\begin{equation}
\mathbf{X}_{\lambda}=\mathbf{X}_{\mathcal{S}^{Y\left(  \lambda\right)  }%
}=\mathbf{X}_{\mathcal{A}\mathbf{E}_{T}}. \label{pf.thm.char.Elam-Xlam.2}%
\end{equation}

On the other hand, Theorem \ref{thm.specht.ETidp} yields $\mathbf{E}_{T}%
^{2}=\dfrac{n!}{f^{\lambda}}\mathbf{E}_{T}$. In other words, $\mathbf{E}%
_{T}^{2}=h^{\lambda}\mathbf{E}_{T}$ (since Definition \ref{def.spechtmod.flam}
\textbf{(b)} yields $\dfrac{n!}{f^{\lambda}}=h^{\lambda}$). Hence, Corollary
\ref{cor.char.ladisch-XSn} (applied to $\mathbf{a}=\mathbf{E}_{T}$ and
$\kappa=h^{\lambda}$) yields that
\[
h^{\lambda}\cdot\mathbf{X}_{\mathcal{A}\mathbf{E}_{T}}=\sum_{w\in S_{n}%
}w\mathbf{E}_{T}w^{-1}.
\]
On the other hand, Proposition \ref{prop.spechtmod.Elam.as-sum-w} (applied to
$P=T$) shows that%
\[
\mathbf{E}_{\lambda}=\sum_{w\in S_{n}}w\mathbf{E}_{T}w^{-1}.
\]
Comparing these two equalities, we obtain $\mathbf{E}_{\lambda}=h^{\lambda
}\cdot\underbrace{\mathbf{X}_{\mathcal{A}\mathbf{E}_{T}}}%
_{\substack{=\mathbf{X}_{\lambda}\\\text{(by (\ref{pf.thm.char.Elam-Xlam.2}%
))}}}=h^{\lambda}\mathbf{X}_{\lambda}$. This proves Theorem
\ref{thm.char.Elam-Xlam}.
\end{proof}

\begin{corollary}
\label{cor.char.Elam-mulhlam}Let $\lambda$ be a partition of $n$. Then, all
coefficients of the element $\mathbf{E}_{\lambda}\in\mathbf{k}\left[
S_{n}\right]  $ are multiples of $h^{\lambda}$ (inside $\mathbf{k}$).
\end{corollary}

\begin{proof}
Theorem \ref{thm.char.Elam-Xlam} yields that $\mathbf{E}_{\lambda}=h^{\lambda
}\mathbf{X}_{\lambda}$ (where $\mathbf{X}_{\lambda}$ is as defined in Theorem
\ref{thm.char.Elam-Xlam}). Thus, each coefficient of $\mathbf{E}_{\lambda}$ is
$h^{\lambda}$ times the corresponding coefficient of $\mathbf{X}_{\lambda}$.
Since the latter coefficient belongs to $\mathbf{k}$, we thus conclude that
the former coefficient is a multiple of $h^{\lambda}$. This proves Corollary
\ref{cor.char.Elam-mulhlam}.
\end{proof}

\begin{remark}
Let $\lambda$ be a partition of $n$. Then, \cite[Chapter 7, supplementary
problem 101]{Stanley-EC2} shows that there exists a permutation $w\in S_{n}$
such that $\chi_{\mathcal{S}^{\lambda}}\left(  w\right)  =\pm1$. Consider this
$w$. Thus, the element $\mathbf{X}_{\lambda}$ of $\mathbf{k}\left[
S_{n}\right]  $ (defined in Theorem \ref{thm.char.Elam-Xlam}) has the
coefficient $\pm1$ in front of $w$. Consequently, the element $\mathbf{E}%
_{\lambda}$ of $\mathbf{k}\left[  S_{n}\right]  $ has the coefficient $\pm
h^{\lambda}$ in front of $w$ (because Theorem \ref{thm.char.Elam-Xlam} yields
that $\mathbf{E}_{\lambda}=h^{\lambda}\mathbf{X}_{\lambda}$). Combined with
Corollary \ref{cor.char.Elam-mulhlam}, this easily yields that $h^{\lambda}$
is the greatest common divisor of all coefficients of $\mathbf{E}_{\lambda}$.
(The notion of \textquotedblleft greatest common divisor\textquotedblright%
\ here should be understood appropriately -- after all, $\mathbf{k}$ is not
necessarily a unique factorization domain; but we can take it literally when
$\mathbf{k}=\mathbb{Z}$.)
\end{remark}

If $h^{\lambda}$ is invertible in $\mathbf{k}$, then the equality
$\mathbf{E}_{\lambda}=h^{\lambda}\mathbf{X}_{\lambda}$ in Theorem
\ref{thm.char.Elam-Xlam} can be solved for $\mathbf{X}_{\lambda}$, yielding
$\mathbf{X}_{\lambda}=\dfrac{1}{h^{\lambda}}\mathbf{E}_{\lambda}$. This
provides a good way of computing $\mathbf{X}_{\lambda}$ and thus the character
$\chi_{\mathcal{S}^{\lambda}}$.

Using Theorem \ref{thm.char.Elam-Xlam}, we can furthermore derive some
properties of the elements $\mathbf{X}_{\lambda}$ from analogous properties of
the $\mathbf{E}_{\lambda}$:

\begin{corollary}
\label{cor.char.Xlam-prop}For each partition $\lambda$ of $n$, we let
$\mathbf{X}_{\lambda}$ denote the element $\mathbf{X}_{\mathcal{S}^{\lambda}}$
corresponding to the Specht module $\mathcal{S}^{\lambda}$. Then: \medskip

\textbf{(a)} We have $\mathbf{X}_{\lambda}^{2}=h^{\lambda}\mathbf{X}_{\lambda
}$ for each partition $\lambda$ of $n$. \medskip

\textbf{(b)} We have $\mathbf{X}_{\lambda}\mathbf{X}_{\mu}=0$ whenever
$\lambda$ and $\mu$ are two distinct partitions of $n$. \medskip

\textbf{(c)} For each $\mathbf{a}\in Z\left(  \mathbf{k}\left[  S_{n}\right]
\right)  $, we have
\[
n!\cdot\mathbf{a}=\sum_{\lambda\text{ is a partition of }n}\left[  1\right]
\left(  \mathbf{X}_{\lambda}\mathbf{a}\right)  \cdot\mathbf{X}_{\lambda}.
\]

\textbf{(d)} If $n!$ is invertible in $\mathbf{k}$, then the family $\left(
\mathbf{X}_{\lambda}\right)  _{\lambda\text{ is a partition of }n}$ is a basis
of the $\mathbf{k}$-module $Z\left(  \mathbf{k}\left[  S_{n}\right]  \right)
$.
\end{corollary}

This corollary is easy to prove in the case when $n!$ is invertible in
$\mathbf{k}$, because in this case we can use the above-mentioned equality
$\mathbf{X}_{\lambda}=\dfrac{1}{h^{\lambda}}\mathbf{E}_{\lambda}$ to reduce
everything to the properties of centralized Young symmetrizers (Theorem
\ref{thm.spechtmod.Elam.idp-strong}, Proposition
\ref{prop.spechtmod.Elam.orth}, Proposition \ref{prop.spechtmod.Elam.a} and
Theorem \ref{thm.spechtmod.Elam.center-gen}). Fortunately, the general case
can be reduced to this easy case using a base change argument similar to the
one we used back in the proof of Proposition
\ref{prop.spechtmod.Elam.idp-weak}. This reduction relies on the following two lemmas:

\begin{lemma}
\label{lem.char.Xlam-prop.reduce-chi}Let $\lambda$ be a partition of $n$. Let
us denote the Specht module $\mathcal{S}^{\lambda}$ (defined using the base
ring $\mathbf{k}$) by $\mathcal{S}_{\mathbf{k}}^{\lambda}$ in order to stress
its dependence on $\mathbf{k}$. Let $g\in S_{n}$. Then,%
\[
\chi_{\mathcal{S}_{\mathbf{k}}^{\lambda}}\left(  g\right)  =\chi
_{\mathcal{S}_{\mathbb{Z}}^{\lambda}}\left(  g\right)  \cdot1_{\mathbf{k}}.
\]

\end{lemma}

\begin{proof}
Before I give a rigorous proof, let me explain why this lemma is obvious (even
though the formal explanation of this obviousness will take us a couple boring
pages). In a nutshell, the argument goes as follows: The Specht module
$\mathcal{S}^{\lambda}$ has the same basis for any base ring $\mathbf{k}$
(namely, the basis described in Theorem \ref{thm.spechtmod.basis}). The matrix
that represents the action of $g\in S_{n}$ with respect to this basis is
independent of $\mathbf{k}$ as well (more precisely, the entries of this
matrix are the same integers for every $\mathbf{k}$, at least if you equate
each integer $m\in\mathbb{Z}$ with the corresponding element $m\cdot
1_{\mathbf{k}}$ of $\mathbf{k}$). Thus, the trace of this matrix does not
depend on $\mathbf{k}$ either. But this trace is precisely $\chi
_{\mathcal{S}_{\mathbf{k}}^{\lambda}}\left(  g\right)  $ (by the definition of
a character). Thus, $\chi_{\mathcal{S}_{\mathbf{k}}^{\lambda}}\left(
g\right)  $ is independent of $\mathbf{k}$. In particular, $\chi
_{\mathcal{S}_{\mathbf{k}}^{\lambda}}\left(  g\right)  $ is the same as
$\chi_{\mathcal{S}_{\mathbb{Z}}^{\lambda}}\left(  g\right)  $, or, to be more
precise, we have $\chi_{\mathcal{S}_{\mathbf{k}}^{\lambda}}\left(  g\right)
=\chi_{\mathcal{S}_{\mathbb{Z}}^{\lambda}}\left(  g\right)  \cdot
1_{\mathbf{k}}$. Thus, Lemma \ref{lem.char.Xlam-prop.reduce-chi} is proved, as
long as you believe the above less-than-rigorous argument. \medskip

\begin{fineprint}
Let me now give a rigorous version of this proof.

Let $D$ be the Young diagram $Y\left(  \lambda\right)  $. Let
$\operatorname*{SYT}\left(  D\right)  $ be the set of all standard tableaux of
shape $D$.

Let us denote the polytabloids $\mathbf{e}_{T}$ introduced in Definition
\ref{def.spechtmod.spechtmod} \textbf{(a)} by $\mathbf{e}_{T,\mathbf{k}}$ in
order to stress their dependence on $\mathbf{k}$. Thus, $\mathbf{e}%
_{T}=\mathbf{e}_{T,\mathbf{k}}$ for any $n$-tableau $T$ of shape $D$.

From $D=Y\left(  \lambda\right)  $, we obtain $\mathcal{S}^{D}=\mathcal{S}%
^{Y\left(  \lambda\right)  }=\mathcal{S}^{\lambda}=\mathcal{S}_{\mathbf{k}%
}^{\lambda}$ (since we are denoting $\mathcal{S}^{\lambda}$ by $\mathcal{S}%
_{\mathbf{k}}^{\lambda}$). Moreover, $D$ is a skew Young diagram (since
$D=Y\left(  \lambda\right)  =Y\left(  \lambda/\varnothing\right)  $). Thus,
Theorem \ref{thm.spechtmod.basis} shows that the standard polytabloids
$\mathbf{e}_{T}$ form a basis of the $\mathbf{k}$-module $\mathcal{S}^{D}$. In
other words, the family $\left(  \mathbf{e}_{T}\right)  _{T\in
\operatorname*{SYT}\left(  D\right)  }$ is a basis of the $\mathbf{k}$-module
$\mathcal{S}^{D}$ (since the standard polytabloids $\mathbf{e}_{T}$ are
precisely the $\mathbf{e}_{T}$ for $T\in\operatorname*{SYT}\left(  D\right)
$). In other words, the family $\left(  \mathbf{e}_{T,\mathbf{k}}\right)
_{T\in\operatorname*{SYT}\left(  D\right)  }$ is a basis of the $\mathbf{k}%
$-module $\mathcal{S}_{\mathbf{k}}^{\lambda}$ (since $\mathcal{S}%
^{D}=\mathcal{S}_{\mathbf{k}}^{\lambda}$, and since $\mathbf{e}_{T}%
=\mathbf{e}_{T,\mathbf{k}}$ for any $n$-tableau $T$ of shape $D$). The same
argument (applied to $\mathbb{Z}$ instead of $\mathbf{k}$) shows that the
family $\left(  \mathbf{e}_{T,\mathbb{Z}}\right)  _{T\in\operatorname*{SYT}%
\left(  D\right)  }$ is a basis of the $\mathbb{Z}$-module $\mathcal{S}%
_{\mathbb{Z}}^{\lambda}$.

In particular, this family spans this $\mathbb{Z}$-module $\mathcal{S}%
_{\mathbb{Z}}^{\lambda}$. Hence, for each $P\in\operatorname*{SYT}\left(
D\right)  $, the element $g\mathbf{e}_{P,\mathbb{Z}}$ of $\mathcal{S}%
_{\mathbb{Z}}^{\lambda}$ can be written as a $\mathbb{Z}$-linear combination
of this family:%
\begin{equation}
g\mathbf{e}_{P,\mathbb{Z}}=\sum_{T\in\operatorname*{SYT}\left(  D\right)
}a_{T,P}\mathbf{e}_{T,\mathbb{Z}} \label{pf.lem.char.Xlam-prop.reduce-chi.a1}%
\end{equation}
for some scalars $a_{T,P}\in\mathbb{Z}$. Consider these scalars $a_{T,P}$. We
now claim the following:

\begin{statement}
\textit{Claim 1:} We have
\[
\chi_{\mathcal{S}_{\mathbb{Z}}^{\lambda}}\left(  g\right)  =\sum
_{T\in\operatorname*{SYT}\left(  D\right)  }a_{T,T}.
\]

\end{statement}

\begin{proof}
[Proof of Claim 1.]Let $\tau_{g,\mathbb{Z}}:\mathcal{S}_{\mathbb{Z}}^{\lambda
}\rightarrow\mathcal{S}_{\mathbb{Z}}^{\lambda}$ be the $\mathbb{Z}$-linear map
$v\mapsto gv$. Then, $\chi_{\mathcal{S}_{\mathbb{Z}}^{\lambda}}\left(
g\right)  =\operatorname*{Tr}\left(  \tau_{g,\mathbb{Z}}\right)  $ (by the
definition of a character). However, the definition of $\tau_{g,\mathbb{Z}}$
shows that each $P\in\operatorname*{SYT}\left(  D\right)  $ satisfies%
\[
\tau_{g,\mathbb{Z}}\left(  \mathbf{e}_{P,\mathbb{Z}}\right)  =g\mathbf{e}%
_{P,\mathbb{Z}}=\sum_{T\in\operatorname*{SYT}\left(  D\right)  }%
a_{T,P}\mathbf{e}_{T,\mathbb{Z}}\ \ \ \ \ \ \ \ \ \ \left(  \text{by
(\ref{pf.lem.char.Xlam-prop.reduce-chi.a1})}\right)  .
\]
Thus, the matrix that represents the $\mathbb{Z}$-linear map $\tau
_{g,\mathbb{Z}}:\mathcal{S}_{\mathbb{Z}}^{\lambda}\rightarrow\mathcal{S}%
_{\mathbb{Z}}^{\lambda}$ with respect to the basis $\left(  \mathbf{e}%
_{T,\mathbb{Z}}\right)  _{T\in\operatorname*{SYT}\left(  D\right)  }$ of the
$\mathbb{Z}$-module $\mathcal{S}_{\mathbb{Z}}^{\lambda}$ is the
matrix\footnote{This is a $\operatorname*{SYT}\left(  D\right)  \times
\operatorname*{SYT}\left(  D\right)  $-matrix, i.e., a matrix whose rows and
columns are indexed by the elements of $\operatorname*{SYT}\left(  D\right)  $
rather than by integers. The trace of such a matrix can be defined in the same
way as the trace of a \textquotedblleft usual\textquotedblright\ $m\times
m$-matrix: It is the sum of its diagonal entries, i.e., of all its entries
whose row index and column index are equal.} $\left(  a_{T,P}\right)
_{T,P\in\operatorname*{SYT}\left(  D\right)  }$. Hence,%
\[
\operatorname*{Tr}\left(  \tau_{g,\mathbb{Z}}\right)  =\operatorname*{Tr}%
\left(  \left(  a_{T,P}\right)  _{T,P\in\operatorname*{SYT}\left(  D\right)
}\right)  =\sum_{T\in\operatorname*{SYT}\left(  D\right)  }a_{T,T}%
\]
(by the definition of the trace of a matrix). In view of $\chi_{\mathcal{S}%
_{\mathbb{Z}}^{\lambda}}\left(  g\right)  =\operatorname*{Tr}\left(
\tau_{g,\mathbb{Z}}\right)  $, we can rewrite this as $\chi_{\mathcal{S}%
_{\mathbb{Z}}^{\lambda}}\left(  g\right)  =\sum_{T\in\operatorname*{SYT}%
\left(  D\right)  }a_{T,T}$. This proves Claim 1.
\end{proof}

Let us denote the Young module $\mathcal{M}^{D}$ (defined using the base ring
$\mathbf{k}$) by $\mathcal{M}_{\mathbf{k}}^{D}$ in order to stress its
dependence on $\mathbf{k}$. Note that $\mathcal{S}_{\mathbf{k}}^{\lambda
}=\mathcal{S}^{D}\subseteq\mathcal{M}^{D}=\mathcal{M}_{\mathbf{k}}^{D}$ and
similarly $\mathcal{S}_{\mathbb{Z}}^{\lambda}\subseteq\mathcal{M}_{\mathbb{Z}%
}^{D}$.

Consider the canonical ring morphism $f:\mathbb{Z}\rightarrow\mathbf{k}$,
which sends each $m\in\mathbb{Z}$ to $m\cdot1_{\mathbf{k}}\in\mathbf{k}$. This
morphism canonically induces a $\mathbb{Z}$-linear map $f_{D}:\mathcal{M}%
_{\mathbb{Z}}^{D}\rightarrow\mathcal{M}_{\mathbf{k}}^{D}$ which sends each
element $\sum_{\substack{\overline{T}\text{ is an }n\text{-tabloid}\\\text{of
shape }D}}c_{\overline{T}}\overline{T}\in\mathcal{M}_{\mathbb{Z}}^{D}$ (where
the coefficients $c_{\overline{T}}$ belong to $\mathbb{Z}$) to $\sum
_{\substack{\overline{T}\text{ is an }n\text{-tabloid}\\\text{of shape }%
D}}f\left(  c_{\overline{T}}\right)  \overline{T}\in\mathcal{M}_{\mathbf{k}%
}^{D}$. It is clear that every $n$-tableau $T$ of shape $D$ satisfies
\begin{equation}
f_{D}\left(  \mathbf{e}_{T,\mathbb{Z}}\right)  =\mathbf{e}_{T,\mathbf{k}}
\label{pf.lem.char.Xlam-prop.reduce-chi.fe}%
\end{equation}
(since the definition of a polytabloid says that $\mathbf{e}_{T}=\sum
_{w\in\mathcal{C}\left(  T\right)  }\left(  -1\right)  ^{w}\overline{wT}$,
independently of $\mathbf{k}$). It is also clear that the map $f_{D}$ is
$S_{n}$-equivariant (since $S_{n}$ acts on $\mathcal{M}_{\mathbb{Z}}^{D}$ and
on $\mathcal{M}_{\mathbf{k}}^{D}$ in the same way, namely by permuting the
$n$-tabloids). Using these two observations, we can easily send the equality
(\ref{pf.lem.char.Xlam-prop.reduce-chi.a1}) \textquotedblleft over into
$\mathcal{S}_{\mathbf{k}}^{\lambda}$\textquotedblright:

\begin{statement}
\textit{Claim 2:} For each $P\in\operatorname*{SYT}\left(  D\right)  $, we
have%
\[
g\mathbf{e}_{P,\mathbf{k}}=\sum_{T\in\operatorname*{SYT}\left(  D\right)
}f\left(  a_{T,P}\right)  \mathbf{e}_{T,\mathbf{k}}.
\]

\end{statement}

\begin{proof}
[Proof of Claim 2.]Let $P\in\operatorname*{SYT}\left(  D\right)  $. Since the
map $f_{D}$ is $S_{n}$-equivariant, we have
\[
f_{D}\left(  g\mathbf{e}_{P,\mathbb{Z}}\right)  =g\cdot\underbrace{f_{D}%
\left(  \mathbf{e}_{P,\mathbb{Z}}\right)  }_{\substack{=\mathbf{e}%
_{P,\mathbf{k}}\\\text{(by (\ref{pf.lem.char.Xlam-prop.reduce-chi.fe}),
applied to }T=P\text{)}}}=g\mathbf{e}_{P,\mathbf{k}}.
\]
Hence,%
\begin{align*}
g\mathbf{e}_{P,\mathbf{k}}  &  =f_{D}\left(  g\mathbf{e}_{P,\mathbb{Z}%
}\right)  =f_{D}\left(  \sum_{T\in\operatorname*{SYT}\left(  D\right)
}a_{T,P}\mathbf{e}_{T,\mathbb{Z}}\right)  \ \ \ \ \ \ \ \ \ \ \left(  \text{by
(\ref{pf.lem.char.Xlam-prop.reduce-chi.a1})}\right) \\
&  =\sum_{T\in\operatorname*{SYT}\left(  D\right)  }a_{T,P}\underbrace{f_{D}%
\left(  \mathbf{e}_{T,\mathbb{Z}}\right)  }_{\substack{=\mathbf{e}%
_{T,\mathbf{k}}\\\text{(by (\ref{pf.lem.char.Xlam-prop.reduce-chi.fe}))}%
}}\ \ \ \ \ \ \ \ \ \ \left(  \text{since the map }f_{D}\text{ is }%
\mathbb{Z}\text{-linear}\right) \\
&  =\sum_{T\in\operatorname*{SYT}\left(  D\right)  }\underbrace{a_{T,P}%
\mathbf{e}_{T,\mathbf{k}}}_{=\left(  a_{T,P}\cdot1_{\mathbf{k}}\right)
\mathbf{e}_{T,\mathbf{k}}}=\sum_{T\in\operatorname*{SYT}\left(  D\right)
}\underbrace{\left(  a_{T,P}\cdot1_{\mathbf{k}}\right)  }_{\substack{=f\left(
a_{T,P}\right)  \\\text{(since the definition of }f\\\text{yields }f\left(
a_{T,P}\right)  =a_{T,P}\cdot1_{\mathbf{k}}\text{)}}}\mathbf{e}_{T,\mathbf{k}%
}\\
&  =\sum_{T\in\operatorname*{SYT}\left(  D\right)  }f\left(  a_{T,P}\right)
\mathbf{e}_{T,\mathbf{k}}.
\end{align*}
This proves Claim 2.
\end{proof}

As a consequence of Claim 2, we can see the following analogue of Claim 1:

\begin{statement}
\textit{Claim 3:} We have $\chi_{\mathcal{S}_{\mathbf{k}}^{\lambda}}\left(
g\right)  =\sum_{T\in\operatorname*{SYT}\left(  D\right)  }f\left(
a_{T,T}\right)  $.
\end{statement}

\begin{proof}
[Proof of Claim 3.]Let $\tau_{g,\mathbf{k}}:\mathcal{S}_{\mathbf{k}}^{\lambda
}\rightarrow\mathcal{S}_{\mathbf{k}}^{\lambda}$ be the $\mathbf{k}$-linear map
$v\mapsto gv$. Then, $\chi_{\mathcal{S}_{\mathbf{k}}^{\lambda}}\left(
g\right)  =\operatorname*{Tr}\left(  \tau_{g,\mathbf{k}}\right)  $ (by the
definition of a character). However, the definition of $\tau_{g,\mathbf{k}}$
shows that each $P\in\operatorname*{SYT}\left(  D\right)  $ satisfies%
\[
\tau_{g,\mathbf{k}}\left(  \mathbf{e}_{P,\mathbf{k}}\right)  =g\mathbf{e}%
_{P,\mathbf{k}}=\sum_{T\in\operatorname*{SYT}\left(  D\right)  }f\left(
a_{T,P}\right)  \mathbf{e}_{T,\mathbf{k}}\ \ \ \ \ \ \ \ \ \ \left(  \text{by
Claim 2}\right)  .
\]
Thus, the matrix that represents the $\mathbf{k}$-linear map $\tau
_{g,\mathbf{k}}:\mathcal{S}_{\mathbf{k}}^{\lambda}\rightarrow\mathcal{S}%
_{\mathbf{k}}^{\lambda}$ with respect to the basis $\left(  \mathbf{e}%
_{T,\mathbf{k}}\right)  _{T\in\operatorname*{SYT}\left(  D\right)  }$ of the
$\mathbf{k}$-module $\mathcal{S}_{\mathbf{k}}^{\lambda}$ is the matrix
$\left(  f\left(  a_{T,P}\right)  \right)  _{T,P\in\operatorname*{SYT}\left(
D\right)  }$. Hence,%
\[
\operatorname*{Tr}\left(  \tau_{g,\mathbf{k}}\right)  =\operatorname*{Tr}%
\left(  \left(  f\left(  a_{T,P}\right)  \right)  _{T,P\in\operatorname*{SYT}%
\left(  D\right)  }\right)  =\sum_{T\in\operatorname*{SYT}\left(  D\right)
}f\left(  a_{T,T}\right)
\]
(by the definition of the trace of a matrix). In view of $\chi_{\mathcal{S}%
_{\mathbf{k}}^{\lambda}}\left(  g\right)  =\operatorname*{Tr}\left(
\tau_{g,\mathbf{k}}\right)  $, we can rewrite this as $\chi_{\mathcal{S}%
_{\mathbf{k}}^{\lambda}}\left(  g\right)  =\sum_{T\in\operatorname*{SYT}%
\left(  D\right)  }f\left(  a_{T,T}\right)  $. This proves Claim 3.
\end{proof}

Now, Claim 3 yields%
\begin{align*}
\chi_{\mathcal{S}_{\mathbf{k}}^{\lambda}}\left(  g\right)   &  =\sum
_{T\in\operatorname*{SYT}\left(  D\right)  }\underbrace{f\left(
a_{T,T}\right)  }_{\substack{=a_{T,T}\cdot1_{\mathbf{k}}\\\text{(by the
definition of }f\text{)}}}=\sum_{T\in\operatorname*{SYT}\left(  D\right)
}a_{T,T}\cdot1_{\mathbf{k}}=\underbrace{\left(  \sum_{T\in\operatorname*{SYT}%
\left(  D\right)  }a_{T,T}\right)  }_{\substack{=\chi_{\mathcal{S}%
_{\mathbb{Z}}^{\lambda}}\left(  g\right)  \\\text{(by Claim 1)}}%
}\cdot\,1_{\mathbf{k}}\\
&  =\chi_{\mathcal{S}_{\mathbb{Z}}^{\lambda}}\left(  g\right)  \cdot
1_{\mathbf{k}}.
\end{align*}
This proves Lemma \ref{lem.char.Xlam-prop.reduce-chi} (rigorously this time).
\end{fineprint}
\end{proof}

\begin{lemma}
\label{lem.char.Xlam-prop.reduce}Let $\lambda$ be a partition of $n$. We let
$\mathbf{X}_{\lambda,\mathbf{k}}$ denote the element $\mathbf{X}%
_{\mathcal{S}^{\lambda}}$ corresponding to the Specht module $\mathcal{S}%
^{\lambda}$.

Consider the canonical ring morphism $f:\mathbb{Z}\rightarrow\mathbf{k}$ and
the base change morphism $f_{\ast}:\mathbb{Z}\left[  S_{n}\right]
\rightarrow\mathbf{k}\left[  S_{n}\right]  $ corresponding to $f$. Then,%
\[
f_{\ast}\left(  \mathbf{X}_{\lambda,\mathbb{Z}}\right)  =\mathbf{X}%
_{\lambda,\mathbf{k}}.
\]

\end{lemma}

\begin{proof}
Let us denote the Specht module $\mathcal{S}^{\lambda}$ (defined using the
base ring $\mathbf{k}$) by $\mathcal{S}_{\mathbf{k}}^{\lambda}$ in order to
stress its dependence on $\mathbf{k}$. Thus, $\mathcal{S}^{\lambda
}=\mathcal{S}_{\mathbf{k}}^{\lambda}$.

Recall that $\mathbf{X}_{\lambda,\mathbf{k}}$ is a notation for $\mathbf{X}%
_{\mathcal{S}^{\lambda}}$. Thus,%
\begin{align}
\mathbf{X}_{\lambda,\mathbf{k}}  &  =\mathbf{X}_{\mathcal{S}^{\lambda}}%
=\sum_{g\in S_{n}}\chi_{\mathcal{S}^{\lambda}}\left(  g\right)
g\ \ \ \ \ \ \ \ \ \ \left(  \text{by Definition \ref{def.char.XV}}\right)
\nonumber\\
&  =\sum_{g\in S_{n}}\chi_{\mathcal{S}_{\mathbf{k}}^{\lambda}}\left(
g\right)  g\ \ \ \ \ \ \ \ \ \ \left(  \text{since }\mathcal{S}^{\lambda
}=\mathcal{S}_{\mathbf{k}}^{\lambda}\right)  .
\label{pf.cor.char.Xlam-prop.a.1kk}%
\end{align}
The same argument (applied to $\mathbb{Z}$ instead of $\mathbf{k}$) shows that%
\[
\mathbf{X}_{\lambda,\mathbb{Z}}=\sum_{g\in S_{n}}\chi_{\mathcal{S}%
_{\mathbb{Z}}^{\lambda}}\left(  g\right)  g.
\]
Applying the map $f_{\ast}$ to both sides of this equality, we obtain%
\begin{align*}
f_{\ast}\left(  \mathbf{X}_{\lambda,\mathbb{Z}}\right)   &  =f_{\ast}\left(
\sum_{g\in S_{n}}\chi_{\mathcal{S}_{\mathbb{Z}}^{\lambda}}\left(  g\right)
g\right) \\
&  =\sum_{g\in S_{n}}\underbrace{f\left(  \chi_{\mathcal{S}_{\mathbb{Z}%
}^{\lambda}}\left(  g\right)  \right)  }_{\substack{=\chi_{\mathcal{S}%
_{\mathbb{Z}}^{\lambda}}\left(  g\right)  \cdot1_{\mathbf{k}}\\\text{(by the
definition of }f\text{)}}}g\ \ \ \ \ \ \ \ \ \ \left(  \text{by the definition
of }f_{\ast}\right) \\
&  =\sum_{g\in S_{n}}\underbrace{\left(  \chi_{\mathcal{S}_{\mathbb{Z}%
}^{\lambda}}\left(  g\right)  \cdot1_{\mathbf{k}}\right)  }_{\substack{=\chi
_{\mathcal{S}_{\mathbf{k}}^{\lambda}}\left(  g\right)  \\\text{(by Lemma
\ref{lem.char.Xlam-prop.reduce-chi})}}}g=\sum_{g\in S_{n}}\chi_{\mathcal{S}%
_{\mathbf{k}}^{\lambda}}\left(  g\right)  g=\mathbf{X}_{\lambda,\mathbf{k}}%
\end{align*}
(by (\ref{pf.cor.char.Xlam-prop.a.1kk})). This proves Lemma
\ref{lem.char.Xlam-prop.reduce}.
\end{proof}

\begin{proof}
[Proof of Corollary \ref{cor.char.Xlam-prop}.]We shall prove part \textbf{(a)}
in detail, and leave the rest as an exercise (Exercise
\ref{exe.char.Xlam-prop}). \medskip

\textbf{(a)} Let $\lambda$ be a partition of $n$. We must prove that
$\mathbf{X}_{\lambda}^{2}=h^{\lambda}\mathbf{X}_{\lambda}$.

We shall first prove this for $\mathbf{k}=\mathbb{Q}$, and then extend it to
the general case using a base change argument similar to the one we used back
in the proof of Proposition \ref{prop.spechtmod.Elam.idp-weak} above. \medskip

\textit{Proof in the case when }$\mathbf{k}=\mathbb{Q}$\textit{:} Assume that
$\mathbf{k}=\mathbb{Q}$. Theorem \ref{thm.spechtmod.Elam.idp-strong} says that
$\mathbf{E}_{\lambda}^{2}=\left(  h^{\lambda}\right)  ^{2}\mathbf{E}_{\lambda
}$. But Theorem \ref{thm.char.Elam-Xlam} yields that $\mathbf{E}_{\lambda
}=h^{\lambda}\mathbf{X}_{\lambda}$. In view of this, we can rewrite
$\mathbf{E}_{\lambda}^{2}=\left(  h^{\lambda}\right)  ^{2}\mathbf{E}_{\lambda
}$ as%
\[
\left(  h^{\lambda}\mathbf{X}_{\lambda}\right)  ^{2}=\left(  h^{\lambda
}\right)  ^{2}h^{\lambda}\mathbf{X}_{\lambda}.
\]
Thus,%
\[
\left(  h^{\lambda}\right)  ^{2}\mathbf{X}_{\lambda}^{2}=\left(  h^{\lambda
}\mathbf{X}_{\lambda}\right)  ^{2}=\left(  h^{\lambda}\right)  ^{2}h^{\lambda
}\mathbf{X}_{\lambda}.
\]
Dividing both sides of this equality by $\left(  h^{\lambda}\right)  ^{2}$
(this is allowed, since $h^{\lambda}$ is a positive integer and since
$\mathbf{k}=\mathbb{Q}$), we obtain
\[
\mathbf{X}_{\lambda}^{2}=h^{\lambda}\mathbf{X}_{\lambda}.
\]

Thus, we have proved $\mathbf{X}_{\lambda}^{2}=h^{\lambda}\mathbf{X}_{\lambda
}$ in the case when $\mathbf{k}=\mathbb{Q}$. \medskip

\textit{General discussion:} Now, forget our assumption that $\mathbf{k}%
=\mathbb{Q}$. Instead, let $\mathbf{k}$ be an arbitrary commutative ring again.

Let us denote the element $\mathbf{X}_{\lambda}$ of $\mathbf{k}\left[
S_{n}\right]  $ by $\mathbf{X}_{\lambda,\mathbf{k}}$ in order to stress its
dependence on $\mathbf{k}$. Thus, $\mathbf{X}_{\lambda,\mathbf{k}}$ is the
element $\mathbf{X}_{\mathcal{S}^{\lambda}}$ corresponding to the Specht
module $\mathcal{S}^{\lambda}$.

Above, we have proved $\mathbf{X}_{\lambda}^{2}=h^{\lambda}\mathbf{X}%
_{\lambda}$ in the case when $\mathbf{k}=\mathbb{Q}$. In other words, we have
proved that%
\begin{equation}
\mathbf{X}_{\lambda,\mathbb{Q}}^{2}=h^{\lambda}\mathbf{X}_{\lambda,\mathbb{Q}%
}. \label{pf.cor.char.Xlam-prop.a.Q}%
\end{equation}

\textit{Proof in the case when }$\mathbf{k}=\mathbb{Z}$\textit{:} Next, let us
prove that $\mathbf{X}_{\lambda,\mathbb{Z}}^{2}=h^{\lambda}\mathbf{X}%
_{\lambda,\mathbb{Z}}$.

Indeed, consider the canonical ring morphism $f:\mathbb{Z}\rightarrow
\mathbb{Q}$ and the base change morphism $f_{\ast}:\mathbb{Z}\left[
S_{n}\right]  \rightarrow\mathbb{Q}\left[  S_{n}\right]  $ corresponding to
$f$. Then, Lemma \ref{lem.char.Xlam-prop.reduce} (applied to $\mathbb{Q}$
instead of $\mathbf{k}$) says that%
\[
f_{\ast}\left(  \mathbf{X}_{\lambda,\mathbb{Z}}\right)  =\mathbf{X}%
_{\lambda,\mathbb{Q}}.
\]
The morphism $f$ is injective (indeed, it is simply the embedding of
$\mathbb{Z}$ into $\mathbb{Q}$). Thus, the morphism $f_{\ast}$ is injective as
well. But $f_{\ast}$ is furthermore a ring morphism; thus,%
\begin{align*}
f_{\ast}\left(  \mathbf{X}_{\lambda,\mathbb{Z}}^{2}\right)   &  =\left(
\underbrace{f_{\ast}\left(  \mathbf{X}_{\lambda,\mathbb{Z}}\right)
}_{=\mathbf{X}_{\lambda,\mathbb{Q}}}\right)  ^{2}=\mathbf{X}_{\lambda
,\mathbb{Q}}^{2}=h^{\lambda}\underbrace{\mathbf{X}_{\lambda,\mathbb{Q}}%
}_{=f_{\ast}\left(  \mathbf{X}_{\lambda,\mathbb{Z}}\right)  }%
\ \ \ \ \ \ \ \ \ \ \left(  \text{by (\ref{pf.cor.char.Xlam-prop.a.Q})}\right)
\\
&  =h^{\lambda}f_{\ast}\left(  \mathbf{X}_{\lambda,\mathbb{Z}}\right)
=f_{\ast}\left(  h^{\lambda}\mathbf{X}_{\lambda,\mathbb{Z}}\right)
\ \ \ \ \ \ \ \ \ \ \left(  \text{since }f_{\ast}\text{ is a ring
morphism}\right)  .
\end{align*}
Since $f_{\ast}$ is injective, we can \textquotedblleft
un-apply\textquotedblright\ $f_{\ast}$ to this equality, and obtain%
\begin{equation}
\mathbf{X}_{\lambda,\mathbb{Z}}^{2}=h^{\lambda}\mathbf{X}_{\lambda,\mathbb{Z}%
}. \label{pf.cor.char.Xlam-prop.a.Z}%
\end{equation}
Thus, we have proved $\mathbf{X}_{\lambda}^{2}=h^{\lambda}\mathbf{X}_{\lambda
}$ in the case when $\mathbf{k}=\mathbb{Z}$. \medskip

\textit{Proof in the general case:} We are now ready to prove $\mathbf{X}%
_{\lambda}^{2}=h^{\lambda}\mathbf{X}_{\lambda}$ for an arbitrary base ring
$\mathbf{k}$.

Indeed, forget the morphisms $f$ and $f_{\ast}$ constructed above. Instead,
let us consider the canonical ring morphism $f:\mathbb{Z}\rightarrow
\mathbf{k}$ and the base change morphism $f_{\ast}:\mathbb{Z}\left[
S_{n}\right]  \rightarrow\mathbf{k}\left[  S_{n}\right]  $ corresponding to
$f$. Then, Lemma \ref{lem.char.Xlam-prop.reduce} says that%
\[
f_{\ast}\left(  \mathbf{X}_{\lambda,\mathbb{Z}}\right)  =\mathbf{X}%
_{\lambda,\mathbf{k}}.
\]
But $f_{\ast}$ is a ring morphism. Thus,
\[
f_{\ast}\left(  \mathbf{X}_{\lambda,\mathbb{Z}}^{2}\right)  =\left(
\underbrace{f_{\ast}\left(  \mathbf{X}_{\lambda,\mathbb{Z}}\right)
}_{=\mathbf{X}_{\lambda,\mathbf{k}}}\right)  ^{2}=\mathbf{X}_{\lambda
,\mathbf{k}}^{2}.
\]
Hence,%
\begin{align*}
\mathbf{X}_{\lambda,\mathbf{k}}^{2}  &  =f_{\ast}\left(
\underbrace{\mathbf{X}_{\lambda,\mathbb{Z}}^{2}}_{\substack{=h^{\lambda
}\mathbf{X}_{\lambda,\mathbb{Z}}\\\text{(by (\ref{pf.cor.char.Xlam-prop.a.Z}%
))}}}\right)  =f_{\ast}\left(  h^{\lambda}\mathbf{X}_{\lambda,\mathbb{Z}%
}\right)  =h^{\lambda}\underbrace{f_{\ast}\left(  \mathbf{X}_{\lambda
,\mathbb{Z}}\right)  }_{=\mathbf{X}_{\lambda,\mathbf{k}}}%
\ \ \ \ \ \ \ \ \ \ \left(
\begin{array}
[c]{c}%
\text{since }f_{\ast}\text{ is a}\\
\text{ring morphism}%
\end{array}
\right) \\
&  =h^{\lambda}\mathbf{X}_{\lambda,\mathbf{k}}.
\end{align*}
Since $\mathbf{X}_{\lambda,\mathbf{k}}=\mathbf{X}_{\lambda}$, we can rewrite
this as $\mathbf{X}_{\lambda}^{2}=h^{\lambda}\mathbf{X}_{\lambda}$. Thus, we
have proved $\mathbf{X}_{\lambda}^{2}=h^{\lambda}\mathbf{X}_{\lambda}$ in the
general case (i.e., for all commutative rings $\mathbf{k}$). This proves
Corollary \ref{cor.char.Xlam-prop} \textbf{(a)}.
\end{proof}

\begin{exercise}
\label{exe.char.Xlam-prop}\fbox{1+2+1} Prove parts \textbf{(b)}, \textbf{(c)}
and \textbf{(d)} of Corollary \ref{cor.char.Xlam-prop}.
\end{exercise}

\subsubsection{Characters of small symmetric groups}

Theorem \ref{thm.char.Elam-Xlam} provides an easy way to compute the elements
$\mathbf{X}_{\lambda}$ when $\mathbf{k}=\mathbb{Q}$. Using Lemma
\ref{lem.char.Xlam-prop.reduce}, we can then extend these formulas for
$\mathbf{X}_{\lambda}$ to $\mathbf{k}=\mathbb{Z}$ and thence to arbitrary
rings $\mathbf{k}$. Then we can compute all character values $\chi
_{\mathcal{S}^{\lambda}}\left(  g\right)  $ by taking appropriate coefficients
of $\mathbf{X}_{\lambda}$.

We shall now give tables of these character values for the first five
symmetric groups $S_{0},S_{1},S_{2},S_{3},S_{4}$. The tables should be read as follows:

\begin{itemize}
\item The rows of our tables correspond to the Specht modules $\mathcal{S}%
^{\lambda}$ indexed by the partitions $\lambda$ of $n$. We denote the
character $\chi_{\mathcal{S}^{\lambda}}$ of a Specht module $\mathcal{S}%
^{\lambda}$ by $\chi_{\lambda}$. Note that $\chi_{\left(  n\right)  }%
=\chi_{\operatorname*{triv}}$ (see Example \ref{exa.char.triv}) and
$\chi_{\left(  1^{n}\right)  }=\chi_{\operatorname*{sign}}$ (see Example
\ref{exa.char.sign}).

\item If two permutations $u,v\in S_{n}$ have the same cycle type, then they
are conjugate (by Theorem \ref{thm.type.conj=type}) and thus satisfy
$\chi_{\lambda}\left(  u\right)  =\chi_{\lambda}\left(  v\right)  $ for each
partition $\lambda$ of $n$ (since Corollary \ref{cor.char.conj=eq} shows that
a character of a representation sends conjugate elements to the same value).
Hence, the value $\chi_{\lambda}\left(  g\right)  $ for a permutation $g\in
S_{n}$ depends only on the cycle type $\operatorname*{type}g$ of $g$, but not
on $g$ itself. Thus, the columns of our tables correspond to the possible
cycle types that a permutation $g\in S_{n}$ can have; we don't need a separate
column for each permutation. Note that the possible cycle types are the
partitions of $n$. Thus, our tables are of square shapes, with both rows and
columns indexed by the partitions of $n$.
\end{itemize}

Here come the tables themselves:%
\[
\underbrace{%
\begin{tabular}
[c]{|c||c|}\hline
$\operatorname*{type}g$ & $\left(  {}\right)  $\\\hline\hline
$\chi_{\left(  {}\right)  }\left(  g\right)  $ & $1$\\\hline
\end{tabular}
}_{n=0}\ \ \ \ \ \ \ \ \ \ \ \ \ \ \ \ \ \ \ \ \underbrace{%
\begin{tabular}
[c]{|c||c|}\hline
$\operatorname*{type}g$ & $\left(  1\right)  $\\\hline\hline
$\chi_{\left(  1\right)  }\left(  g\right)  $ & $1$\\\hline
\end{tabular}
}_{n=1}\ \ \ \ \ \ \ \ \ \ \ \ \ \ \ \ \ \ \ \ \underbrace{%
\begin{tabular}
[c]{|c||c|c|}\hline
$\operatorname*{type}g$ & $\left(  2\right)  $ & $\left(  1,1\right)
$\\\hline\hline
$\chi_{\left(  2\right)  }\left(  g\right)  $ & $1$ & $1$\\\hline
$\chi_{\left(  1,1\right)  }\left(  g\right)  $ & $-1$ & $1$\\\hline
\end{tabular}
}_{n=2}%
\]%
\[
\underbrace{%
\begin{tabular}
[c]{|c||c|c|c|}\hline
$\operatorname*{type}g$ & $\left(  3\right)  $ & $\left(  2,1\right)  $ &
$\left(  1,1,1\right)  $\\\hline\hline
$\chi_{\left(  3\right)  }\left(  g\right)  $ & $1$ & $1$ & $1$\\\hline
$\chi_{\left(  2,1\right)  }\left(  g\right)  $ & $-1$ & $0$ & $2$\\\hline
$\chi_{\left(  1,1,1\right)  }\left(  g\right)  $ & $1$ & $-1$ & $1$\\\hline
\end{tabular}
}_{n=3}%
\]%
\[
\underbrace{%
\begin{tabular}
[c]{|c||c|c|c|c|c|}\hline
$\operatorname*{type}g$ & $\left(  4\right)  $ & $\left(  3,1\right)  $ &
$\left(  2,2\right)  $ & $\left(  2,1,1\right)  $ & $\left(  1,1,1,1\right)
$\\\hline\hline
$\chi_{\left(  4\right)  }\left(  g\right)  $ & $1$ & $1$ & $1$ & $1$ &
$1$\\\hline
$\chi_{\left(  3,1\right)  }\left(  g\right)  $ & $-1$ & $0$ & $-1$ & $1$ &
$3$\\\hline
$\chi_{\left(  2,2\right)  }\left(  g\right)  $ & $0$ & $-1$ & $2$ & $0$ &
$2$\\\hline
$\chi_{\left(  2,1,1\right)  }\left(  g\right)  $ & $1$ & $0$ & $-1$ & $-1$ &
$3$\\\hline
$\chi_{\left(  1,1,1,1\right)  }\left(  g\right)  $ & $-1$ & $1$ & $1$ & $-1$
& $1$\\\hline
\end{tabular}
}_{n=4}\ \ .
\]
Character tables of $S_{5},S_{6},\ldots,S_{10}$ can be found in \cite[\S I.A]%
{JamKer81}. Further character tables can be calculated on a computer, e.g.,
using SageMath\footnote{For instance,
\texttt{SymmetricGroup(8).character\_table()} returns the character table of
$S_{8}$, and similarly for any $S_{n}$. Sadly, I don't know in which order the
rows and the columns are printed in this table (I think it is in reverse
lexicographic order of the corresponding partitions; this would be the reverse
of the order of my tables above).
\par
The easiest way to obtain a given value $\chi_{\lambda}\left(  g\right)  $ in
SageMath is probably using symmetric functions: The code
\par
\qquad\texttt{Sym = SymmetricFunctions(QQ); p = Sym.p(); s = Sym.s()}
\par
\qquad\texttt{s[Partition([3,2])].scalar(p[Partition([2,1,1,1])])}
\par
returns the value $\chi_{\left(  3,2\right)  }\left(  g\right)  $ for a
permutation $g\in S_{5}$ with cycle type $\operatorname*{type}g=\left(
2,1,1,1\right)  $. Similarly for other partitions and other permutations. The
theory that underlies this method is that of the Frobenius characteristic map
(see, e.g., \cite[\S 4.7]{Sagan01} or \cite[\S 7.18]{Stanley-EC2}).}.

\begin{exercise}
\label{exe.char.Slam.sums}Let $\lambda$ be a partition of $n$. Let us denote
the character $\chi_{\mathcal{S}^{\lambda}}$ of the Specht module
$\mathcal{S}^{\lambda}$ by $\chi_{\lambda}$. \medskip

\textbf{(a)} \fbox{2} Prove that $\sum_{g\in S_{n}}\chi_{\lambda}\left(
g\right)  =%
\begin{cases}
n!, & \text{if }\lambda=\left(  n\right)  ;\\
0, & \text{if }\lambda\neq\left(  n\right)  .
\end{cases}
$ \medskip

\textbf{(b)} \fbox{3} Prove that $\sum_{g\in S_{n}}\left(  \chi_{\lambda
}\left(  g\right)  \right)  ^{2}=n!$. \medskip

\textbf{(c)} \fbox{?} Prove that $\sum_{g\in S_{n}}\chi_{\lambda}\left(
g^{2}\right)  =n!$. \medskip

[\textbf{Hint:} For part \textbf{(a)}, show that $\nabla\mathbf{X}_{\lambda
}=0$ for $\lambda\neq\left(  n\right)  $ first. What does this have to do with
part \textbf{(a)}? Likewise, what does part \textbf{(b)} have to do with
$\mathbf{X}_{\lambda}^{2}=h^{\lambda}\mathbf{X}_{\lambda}$ ?

I don't know a proof of part \textbf{(c)} that does not use any deeper
representation theory (at least multilinear algebra). I would welcome any suggestions!]
\end{exercise}

\subsubsection{An explicit formula for the character of $\mathcal{S}^{\lambda
}$}

The following combinatorial formula gives a faster way to compute the
character $\chi_{\mathcal{S}^{\lambda}}$ of a straight-shaped Specht module
$\mathcal{S}^{\lambda}$ than Theorem \ref{thm.char.Elam-Xlam}.

\begin{theorem}
\label{thm.char.chilam-sum}Let $\lambda$ be a partition of $n$. Consider the
character $\chi_{\mathcal{S}^{\lambda}}$ of the corresponding Specht module
$\mathcal{S}^{\lambda}$. Let $g\in S_{n}$ be any permutation. Let $C$ be the
conjugacy class of $S_{n}$ that contains $g$. Let $T$ be any $n$-tableau of
shape $\lambda$. Then,%
\[
\chi_{\mathcal{S}^{\lambda}}\left(  g\right)  =\underbrace{\dfrac{f^{\lambda}%
}{\left\vert C\right\vert }\left(  \sum_{\substack{\left(  c,r\right)
\in\mathcal{C}\left(  T\right)  \times\mathcal{R}\left(  T\right)  ;\\cr\in
C}}\left(  -1\right)  ^{c}\right)  }_{\text{This is an integer.}}%
\cdot\,1_{\mathbf{k}}.
\]
(See Definition \ref{def.spechtmod.flam} \textbf{(a)} for the meaning of
$f^{\lambda}$.)
\end{theorem}

\begin{proof}
Define the integer%
\begin{equation}
k:=\sum_{\substack{\left(  c,r\right)  \in\mathcal{C}\left(  T\right)
\times\mathcal{R}\left(  T\right)  ;\\cr\in C}}\left(  -1\right)  ^{c}.
\label{pf.thm.char.chilam-sum.k=}%
\end{equation}

We start with the following observation:

\begin{statement}
\textit{Claim 1:} We have $\chi_{\mathcal{S}^{\lambda}}\left(  c\right)
=\chi_{\mathcal{S}^{\lambda}}\left(  g\right)  $ for each $c\in C$.
\end{statement}

\begin{proof}
[Proof of Claim 1.]Let $c\in C$. Then, $c$ and $g$ belong to the same
conjugacy class of $S_{n}$ (namely, to the class $C$). In other words, the
elements $c$ and $g$ are conjugate in $S_{n}$. In other words, $c=wgw^{-1}$
for some $w\in S_{n}$. Consider this $w$. From $c=wgw^{-1}$, we obtain
$\chi_{\mathcal{S}^{\lambda}}\left(  c\right)  =\chi_{\mathcal{S}^{\lambda}%
}\left(  wgw^{-1}\right)  =\chi_{\mathcal{S}^{\lambda}}\left(  g\right)  $ (by
Corollary \ref{cor.char.conj=eq}, applied to $S_{n}$ and $\mathcal{S}%
^{\lambda}$ instead of $G$ and $V$). This proves Claim 1.
\end{proof}

Now, let $\mathbf{z}_{C}$ denote the conjugacy class sum $\sum_{c\in C}%
c\in\mathbf{k}\left[  S_{n}\right]  $ (as in Definition \ref{def.groups.conj}
\textbf{(c)}). Then, each $\mathbf{v}\in\mathcal{S}^{\lambda}$ satisfies%
\begin{align}
\mathbf{z}_{C}\mathbf{v}  &  =\underbrace{\left(  \sum_{\substack{\left(
c,r\right)  \in\mathcal{C}\left(  T\right)  \times\mathcal{R}\left(  T\right)
;\\cr\in C}}\left(  -1\right)  ^{c}\right)  }_{\substack{=k\\\text{(by
(\ref{pf.thm.char.chilam-sum.k=}))}}}\mathbf{v}\ \ \ \ \ \ \ \ \ \ \left(
\text{by Theorem \ref{thm.spechtmod.zC}}\right) \nonumber\\
&  =k\mathbf{v}. \label{pf.thm.char.chilam-sum.k1}%
\end{align}
Thus, we can easily see the following:

\begin{statement}
\textit{Claim 2:} We have $\chi_{\mathcal{S}^{\lambda}}\left(  \mathbf{z}%
_{C}\right)  =f^{\lambda}k\cdot1_{\mathbf{k}}$.
\end{statement}

\begin{proof}
[Proof of Claim 2.]The $\mathbf{k}$-module $\mathcal{S}^{\lambda}$ is free of
rank $f^{\lambda}$ (by Lemma \ref{lem.specht.Slam-flam}).

For each $r\in\mathbf{k}\left[  S_{n}\right]  $, we let $\tau_{r}%
:\mathcal{S}^{\lambda}\rightarrow\mathcal{S}^{\lambda}$ be the $\mathbf{k}%
$-linear map $v\mapsto rv$. The definition of the character $\chi
_{\mathcal{S}^{\lambda}}$ then shows that $\chi_{\mathcal{S}^{\lambda}}\left(
r\right)  =\operatorname*{Tr}\left(  \tau_{r}\right)  $ for each
$r\in\mathbf{k}\left[  S_{n}\right]  $. Thus, in particular, $\chi
_{\mathcal{S}^{\lambda}}\left(  \mathbf{z}_{C}\right)  =\operatorname*{Tr}%
\left(  \tau_{\mathbf{z}_{C}}\right)  $.

However, from (\ref{pf.thm.char.chilam-sum.k1}), we easily obtain
$\tau_{\mathbf{z}_{C}}=\left(  k\cdot1_{\mathbf{k}}\right)  \operatorname*{id}%
\nolimits_{\mathcal{S}^{\lambda}}$\ \ \ \ \footnote{\textit{Proof.} For each
$\mathbf{v}\in\mathcal{S}^{\lambda}$, we have%
\begin{align*}
\tau_{\mathbf{z}_{C}}\left(  \mathbf{v}\right)   &  =\mathbf{z}_{C}%
\mathbf{v}\ \ \ \ \ \ \ \ \ \ \left(  \text{by the definition of }%
\tau_{\mathbf{z}_{C}}\right) \\
&  =k\mathbf{v}\ \ \ \ \ \ \ \ \ \ \left(  \text{by
(\ref{pf.thm.char.chilam-sum.k1})}\right) \\
&  =\left(  k\cdot1_{\mathbf{k}}\right)  \mathbf{v}=\left(  \left(
k\cdot1_{\mathbf{k}}\right)  \operatorname*{id}\nolimits_{\mathcal{S}%
^{\lambda}}\right)  \left(  \mathbf{v}\right)  .
\end{align*}
In other words, $\tau_{\mathbf{z}_{C}}=\left(  k\cdot1_{\mathbf{k}}\right)
\operatorname*{id}\nolimits_{\mathcal{S}^{\lambda}}$.}. Hence,
\begin{align*}
\operatorname*{Tr}\left(  \tau_{\mathbf{z}_{C}}\right)   &
=\operatorname*{Tr}\left(  \left(  k\cdot1_{\mathbf{k}}\right)
\operatorname*{id}\nolimits_{\mathcal{S}^{\lambda}}\right) \\
&  =f^{\lambda}\left(  k\cdot1_{\mathbf{k}}\right)
\ \ \ \ \ \ \ \ \ \ \left(
\begin{array}
[c]{c}%
\text{by Lemma \ref{lem.linalg.Trkid},}\\
\text{applied to }W=\mathcal{S}^{\lambda}\text{ and }d=f^{\lambda}\text{ and
}\kappa=k\cdot1_{\mathbf{k}}\\
\text{(since }\mathcal{S}^{\lambda}\text{ is a free }\mathbf{k}\text{-module
of rank }f^{\lambda}\text{)}%
\end{array}
\right) \\
&  =f^{\lambda}k\cdot1_{\mathbf{k}}.
\end{align*}
In view of $\chi_{\mathcal{S}^{\lambda}}\left(  \mathbf{z}_{C}\right)
=\operatorname*{Tr}\left(  \tau_{\mathbf{z}_{C}}\right)  $, we can rewrite
this as $\chi_{\mathcal{S}^{\lambda}}\left(  \mathbf{z}_{C}\right)
=f^{\lambda}k\cdot1_{\mathbf{k}}$. Thus, Claim 2 is proved.
\end{proof}

Combining Claim 1 with Claim 2, we obtain:

\begin{statement}
\textit{Claim 3:} We have $\left\vert C\right\vert \cdot\chi_{\mathcal{S}%
^{\lambda}}\left(  g\right)  =f^{\lambda}k\cdot1_{\mathbf{k}}$.
\end{statement}

\begin{proof}
[Proof of Claim 3.]The character $\chi_{\mathcal{S}^{\lambda}}:\mathbf{k}%
\left[  S_{n}\right]  \rightarrow\mathbf{k}$ is a $\mathbf{k}$-linear map (by
Proposition \ref{prop.char.klin}). Now, from $\mathbf{z}_{C}=\sum_{c\in C}c$,
we obtain%
\begin{align*}
\chi_{\mathcal{S}^{\lambda}}\left(  \mathbf{z}_{C}\right)   &  =\chi
_{\mathcal{S}^{\lambda}}\left(  \sum_{c\in C}c\right)  =\sum_{c\in
C}\underbrace{\chi_{\mathcal{S}^{\lambda}}\left(  c\right)  }_{\substack{=\chi
_{\mathcal{S}^{\lambda}}\left(  g\right)  \\\text{(by Claim 1)}}%
}\ \ \ \ \ \ \ \ \ \ \left(  \text{since }\chi_{\mathcal{S}^{\lambda}}\text{
is a }\mathbf{k}\text{-linear map}\right) \\
&  =\sum_{c\in C}\chi_{\mathcal{S}^{\lambda}}\left(  g\right)  =\left\vert
C\right\vert \cdot\chi_{\mathcal{S}^{\lambda}}\left(  g\right)  .
\end{align*}
Thus, $\left\vert C\right\vert \cdot\chi_{\mathcal{S}^{\lambda}}\left(
g\right)  =\chi_{\mathcal{S}^{\lambda}}\left(  \mathbf{z}_{C}\right)
=f^{\lambda}k\cdot1_{\mathbf{k}}$ (by Claim 2). This proves Claim 3.
\end{proof}

From Claim 3, we aim to recover $\chi_{\mathcal{S}^{\lambda}}\left(  g\right)
$. Unfortunately, we cannot divide by $\left\vert C\right\vert $ when
$\mathbf{k}$ is just an arbitrary commutative ring. Thus, we make a short
detour via the ring $\mathbb{Z}$ and a base change argument.

Let us denote the Specht module $\mathcal{S}^{\lambda}$ (defined using the
base ring $\mathbf{k}$) by $\mathcal{S}_{\mathbf{k}}^{\lambda}$ in order to
stress its dependence on $\mathbf{k}$. Now, we claim the following:

\begin{statement}
\textit{Claim 4:} We have $\dfrac{f^{\lambda}k}{\left\vert C\right\vert }%
\in\mathbb{Z}$ and $\chi_{\mathcal{S}_{\mathbb{Z}}^{\lambda}}\left(  g\right)
=\dfrac{f^{\lambda}k}{\left\vert C\right\vert }$.
\end{statement}

\begin{proof}
[Proof of Claim 4.]Applying Claim 3 to $\mathbb{Z}$ instead of $\mathbf{k}$,
we find $\left\vert C\right\vert \cdot\chi_{\mathcal{S}_{\mathbb{Z}}^{\lambda
}}\left(  g\right)  =f^{\lambda}k\cdot1_{\mathbb{Z}}=f^{\lambda}k$. Since
$\left\vert C\right\vert $ is a positive integer (because the conjugacy class
$C$ is clearly a nonempty finite set), we can divide this equality by
$\left\vert C\right\vert $, and find $\chi_{\mathcal{S}_{\mathbb{Z}}^{\lambda
}}\left(  g\right)  =\dfrac{f^{\lambda}k}{\left\vert C\right\vert }$ in
$\mathbb{Q}$. Thus, $\dfrac{f^{\lambda}k}{\left\vert C\right\vert }%
=\chi_{\mathcal{S}_{\mathbb{Z}}^{\lambda}}\left(  g\right)  \in\mathbb{Z}$
(since $\chi_{\mathcal{S}_{\mathbb{Z}}^{\lambda}}$ is a map from
$\mathbb{Z}\left[  S_{n}\right]  $ to $\mathbb{Z}$). Thus, Claim 4 is proved.
\end{proof}

We now return to $\mathbf{k}$:

\begin{statement}
\textit{Claim 5:} We have $\chi_{\mathcal{S}^{\lambda}}\left(  g\right)
=\dfrac{f^{\lambda}k}{\left\vert C\right\vert }\cdot1_{\mathbf{k}}$.
\end{statement}

\begin{proof}
[Proof of Claim 5.]Lemma \ref{lem.char.Xlam-prop.reduce-chi} yields
\[
\chi_{\mathcal{S}_{\mathbf{k}}^{\lambda}}\left(  g\right)  =\underbrace{\chi
_{\mathcal{S}_{\mathbb{Z}}^{\lambda}}\left(  g\right)  }_{\substack{=\dfrac
{f^{\lambda}k}{\left\vert C\right\vert }\\\text{(by Claim 4)}}}\cdot
\,1_{\mathbf{k}}=\dfrac{f^{\lambda}k}{\left\vert C\right\vert }\cdot
1_{\mathbf{k}}.
\]
Since $\mathcal{S}_{\mathbf{k}}^{\lambda}=\mathcal{S}^{\lambda}$, we can
rewrite this as $\chi_{\mathcal{S}^{\lambda}}\left(  g\right)  =\dfrac
{f^{\lambda}k}{\left\vert C\right\vert }\cdot1_{\mathbf{k}}$. Thus, Claim 5 is proved.
\end{proof}

Now, Claim 5 yields%
\[
\chi_{\mathcal{S}^{\lambda}}\left(  g\right)  =\dfrac{f^{\lambda}k}{\left\vert
C\right\vert }\cdot1_{\mathbf{k}}=\dfrac{f^{\lambda}}{\left\vert C\right\vert
}k\cdot1_{\mathbf{k}}=\dfrac{f^{\lambda}}{\left\vert C\right\vert }\left(
\sum_{\substack{\left(  c,r\right)  \in\mathcal{C}\left(  T\right)
\times\mathcal{R}\left(  T\right)  ;\\cr\in C}}\left(  -1\right)  ^{c}\right)
\cdot1_{\mathbf{k}}%
\]
(by (\ref{pf.thm.char.chilam-sum.k=})). Moreover, $\dfrac{f^{\lambda}%
}{\left\vert C\right\vert }\underbrace{\left(  \sum_{\substack{\left(
c,r\right)  \in\mathcal{C}\left(  T\right)  \times\mathcal{R}\left(  T\right)
;\\cr\in C}}\left(  -1\right)  ^{c}\right)  }_{=k}=\dfrac{f^{\lambda}%
}{\left\vert C\right\vert }k=\dfrac{f^{\lambda}k}{\left\vert C\right\vert }%
\in\mathbb{Z}$ (by Claim 4). In other words, $\dfrac{f^{\lambda}}{\left\vert
C\right\vert }\left(  \sum_{\substack{\left(  c,r\right)  \in\mathcal{C}%
\left(  T\right)  \times\mathcal{R}\left(  T\right)  ;\\cr\in C}}\left(
-1\right)  ^{c}\right)  $ is an integer. Thus, Theorem
\ref{thm.char.chilam-sum} is proved.
\end{proof}

Specifically for transpositions, we obtain the following simple formula:

\begin{corollary}
\label{cor.char.chilam-transp}Let $\lambda$ be a partition of $n$. Consider
the character $\chi_{\mathcal{S}^{\lambda}}$ of the corresponding Specht
module $\mathcal{S}^{\lambda}$. Let $g\in S_{n}$ be any transposition. Then,%
\[
\chi_{\mathcal{S}^{\lambda}}\left(  g\right)  =\underbrace{\dfrac{f^{\lambda}%
}{\dbinom{n}{2}}\left(  \operatorname*{n}\left(  \lambda^{t}\right)
-\operatorname*{n}\left(  \lambda\right)  \right)  }_{\text{This is an
integer.}}\cdot\,1_{\mathbf{k}}.
\]
(See Definition \ref{def.spechtmod.flam} \textbf{(a)} for the meaning of
$f^{\lambda}$, and Definition \ref{def.partitions.nlam} for the meaning of
$\operatorname*{n}\left(  \lambda\right)  $.)
\end{corollary}

\begin{proof}
Pick any $n$-tableau $T$ of shape $Y\left(  \lambda\right)  $ (such a $T$
exists, since $\left\vert Y\left(  \lambda\right)  \right\vert =\left\vert
\lambda\right\vert =n$). Thus, $T$ is an $n$-tableau of shape $\lambda$.

Let $C$ be the set of all transpositions in $S_{n}$. Then, $C$ is a conjugacy
class of $S_{n}$ (since we know that all transpositions in $S_{n}$ form a
single conjugacy class), and contains $g$ (since $g$ is a transposition).
Thus, $C$ is the conjugacy class of $S_{n}$ that contains $g$. Theorem
\ref{thm.char.chilam-sum} thus yields%
\begin{equation}
\chi_{\mathcal{S}^{\lambda}}\left(  g\right)  =\underbrace{\dfrac{f^{\lambda}%
}{\left\vert C\right\vert }\left(  \sum_{\substack{\left(  c,r\right)
\in\mathcal{C}\left(  T\right)  \times\mathcal{R}\left(  T\right)  ;\\cr\in
C}}\left(  -1\right)  ^{c}\right)  }_{\text{This is an integer.}}%
\cdot\,1_{\mathbf{k}}. \label{pf.cor.char.chilam-transp.1}%
\end{equation}
However, we have $\left\vert C\right\vert =\dbinom{n}{2}$%
\ \ \ \ \footnote{\textit{Proof.} The transpositions in $S_{n}$ are the
permutations of the form $t_{i,j}$, where $i$ and $j$ are two distinct
elements of $\left[  n\right]  $. Moreover, the order of $i$ and $j$ does not
matter (i.e., we have $t_{i,j}=t_{j,i}$). There are $\dbinom{n}{2}$ ways to
choose two distinct elements $i$ and $j$ of $\left[  n\right]  $ if we
disregard the order; and each of these $\dbinom{n}{2}$ choices produces a
different transposition $t_{i,j}$. Thus, there are exactly $\dbinom{n}{2}$
transpositions in $S_{n}$. In other words, $\left\vert C\right\vert
=\dbinom{n}{2}$ (since $C$ is the set of all transpositions in $S_{n}$).} and
\[
\sum_{\substack{\left(  c,r\right)  \in\mathcal{C}\left(  T\right)
\times\mathcal{R}\left(  T\right)  ;\\cr\in C}}\left(  -1\right)
^{c}=\operatorname*{n}\left(  \lambda^{t}\right)  -\operatorname*{n}\left(
\lambda\right)  \ \ \ \ \ \ \ \ \ \ \left(  \text{by Lemma
\ref{lem.spechtmod.sum-transp.3}}\right)  .
\]
In view of this, we can rewrite (\ref{pf.cor.char.chilam-transp.1}) as%
\[
\chi_{\mathcal{S}^{\lambda}}\left(  g\right)  =\underbrace{\dfrac{f^{\lambda}%
}{\dbinom{n}{2}}\left(  \operatorname*{n}\left(  \lambda^{t}\right)
-\operatorname*{n}\left(  \lambda\right)  \right)  }_{\text{This is an
integer.}}\cdot\,1_{\mathbf{k}}.
\]
Thus, Corollary \ref{cor.char.chilam-transp} is proved.
\end{proof}

\subsubsection{An outlook on the Murnaghan--Nakayama rule}

Let us now briefly note a different -- combinatorial -- approach to computing
the characters of Specht modules. This is the famous \emph{Murnaghan--Nakayama
rule}. It relies on the concept of a \emph{rim hook}:

\begin{definition}
\label{def.rimhook.rimhook}A \emph{rim hook} means a nonempty diagram $D$ that
can be written in the form%
\[
D=\left\{  c_{1},c_{2},\ldots,c_{k}\right\}  ,
\]
where $c_{1},c_{2},\ldots,c_{k}\in\mathbb{Z}^{2}$ are cells such that for each
$i\in\left\{  2,3,\ldots,k\right\}  $, the cell $c_{i}$ is either the northern
neighbor or the eastern neighbor of $c_{i-1}$.
\end{definition}

For example, the following diagrams are rim hooks:%
\[
\begin{ytableau}
\none & \none & \none & {c_6} \\
\none & {c_3} & {c_4} & {c_5} \\
{c_1} & {c_2}
\end{ytableau}\qquad\text{and}\qquad\begin{ytableau}
\none & c_4 \\
c_2 & c_3 \\
c_1
\end{ytableau}\qquad\text{and}\qquad\begin{ytableau}
c_2 & c_3 & c_4 \\
c_1
\end{ytableau}
\]
(where we have labelled the cells by $c_{1},c_{2},\ldots,c_{k}$ in a way that
makes the requirement in Definition \ref{def.rimhook.rimhook} true). In
contrast, none of the diagrams%
\[
\ydiagram{1+3,0+3}\qquad\text{and}\qquad\ydiagram{2+1,0+2}\qquad
\text{and}\qquad\ydiagram{0+1,0+2}
\]
is a rim hook. Rim hooks are also known as \emph{ribbons} or \emph{border
strips} or \emph{skew hooks}.

\begin{exercise}
\fbox{2} Prove that every rim hook contained in $\left\{  1,2,3,\ldots
\right\}  ^{2}$ is a skew Young diagram.
\end{exercise}

\begin{exercise}
\fbox{1} Let $\lambda$ be a partition. Prove that the diagram $Y\left(
\lambda\right)  $ is a rim hook if and only if $\lambda$ is a hook partition.
\end{exercise}

\begin{exercise}
\fbox{1} Recall the map $\mathbf{r}:\mathbb{Z}^{2}\rightarrow\mathbb{Z}^{2}$
defined in Theorem \ref{thm.partitions.conj}. Let $D$ be a diagram. Prove that
$D$ is a rim hook if and only if $\mathbf{r}\left(  D\right)  $ is a rim hook.
\end{exercise}

The following exercises give several equivalent characterizations of rim hooks:

\begin{exercise}
\fbox{2} Define a map $h:\mathbb{Z}^{2}\rightarrow\mathbb{Z}$ by setting%
\[
h\left(  i,j\right)  =j-i\ \ \ \ \ \ \ \ \ \ \text{for each }\left(
i,j\right)  \in\mathbb{Z}^{2}.
\]

Let $\lambda/\mu$ be a skew partition. Prove that the diagram $Y\left(
\lambda/\mu\right)  $ is a rim hook if and only if the restriction of the map
$h$ to $Y\left(  \lambda/\mu\right)  $ is injective and its image $h\left(
Y\left(  \lambda/\mu\right)  \right)  $ is the integer interval $\left\{
u,u+1,\ldots,v\right\}  $ for some integers $u\leq v$.
\end{exercise}

\begin{exercise}
\fbox{2} Let $\lambda/\mu$ be a skew partition. Let $G$ be the undirected
graph whose vertices are the cells in $Y\left(  \lambda/\mu\right)  $, and
whose edges are defined as follows: Two vertices $c$ and $d$ are adjacent if
and only if the cells $c$ and $d$ are neighbors (i.e., the Euclidean distance
between $c$ and $d$ is $1$).

Prove that the diagram $Y\left(  \lambda/\mu\right)  $ is a rim hook if and
only if the graph $G$ is connected and the diagram $Y\left(  \lambda
/\mu\right)  $ contains no $2\times2$-square (i.e., there exists no $\left(
i,j\right)  \in Y\left(  \lambda/\mu\right)  $ such that $\left(
i,j+1\right)  $ and $\left(  i+1,j\right)  $ and $\left(  i+1,j+1\right)  $
all belong to $Y\left(  \lambda/\mu\right)  $ as well).
\end{exercise}

\begin{exercise}
\fbox{4} Let $\lambda/\mu$ be a skew partition. Prove that the diagram
$Y\left(  \lambda/\mu\right)  $ is a rim hook if and only if there exist two
positive integers $u\leq v$ such that
\begin{align*}
&  \lambda_{u}>\mu_{u},\\
&  \left(  \lambda_{i}=\mu_{i-1}+1\text{ for each }i\in\left\{  u+1,u+2,\ldots
,v\right\}  \right)  \ \ \ \ \ \ \ \ \ \ \text{and}\\
&  \left(  \lambda_{i}=\mu_{i}\text{ for each }i\notin\left\{  u,u+1,\ldots
,v\right\}  \right)  .
\end{align*}

\end{exercise}

We next define the \emph{leg length} of a rim hook:

\begin{definition}
Let $D$ be a rim hook. Then, the \emph{leg length} of $D$ is defined to be
$\left\vert \operatorname*{Rows}D\right\vert -1$, where $\operatorname*{Rows}%
D:=\left\{  i\ \mid\ \left(  i,j\right)  \in D\right\}  $.
\end{definition}

For instance, the rim hooks%
\[
\begin{ytableau}
\none & \none & \none & {c_6} \\
\none & {c_3} & {c_4} & {c_5} \\
{c_1} & {c_2}
\end{ytableau}\qquad\text{and}\qquad\begin{ytableau}
\none & c_4 \\
c_2 & c_3 \\
c_1
\end{ytableau}\qquad\text{and}\qquad\begin{ytableau}
c_2 & c_3 & c_4 \\
c_1
\end{ytableau}
\]
have leg lengths $2$ and $2$ and $1$, respectively. \medskip

Next, we define \emph{rim hook tableaux}:

\begin{definition}
Let $D$ be a diagram. Let $\mu=\left(  \mu_{1},\mu_{2},\ldots,\mu_{k}\right)
$ be a tuple of nonnegative integers. \medskip

\textbf{(a)} A \emph{rim hook tableau} of shape $D$ and content $\mu$ means a
tableau $T$ of shape $D$ with the following properties:

\begin{itemize}
\item All entries of $T$ belong to the set $\left[  k\right]  $.

\item The tableau $T$ is row-semistandard and column-semistandard.

\item For each $i\in\left[  k\right]  $, the set of all cells of $T$ that
contain the entry $i$ has size $\mu_{i}$ and is either a rim hook or empty. We
shall denote this set by $R_{i}\left(  T\right)  $.
\end{itemize}

\textbf{(b)} If $T$ is such a rim hook tableau, then the \emph{leg length} of
$T$ means the sum of the leg lengths of the rim hooks $R_{1}\left(  T\right)
,\ R_{2}\left(  T\right)  ,\ \ldots,\ R_{k}\left(  T\right)  $ (where the
empty set is considered to have leg length $0$).
\end{definition}

\begin{example}
Let $D$ be the skew Young diagram $Y\left(  \left(  5,3,2\right)  /\left(
1,1\right)  \right)  $. Then, the only rim hook tableaux of shape $D$ and
content $\left(  4,4\right)  $ are%
\[
\ytableaushort{\none1222,\none12,11}\ \ \ \ \ \ \ \ \ \ \text{and}%
\ \ \ \ \ \ \ \ \ \ \ytableaushort{\none1111,\none22,22}\ \ .
\]
Their respective leg lengths are $2+1=3$ and $0+1=1$.
\end{example}

Note that a rim hook tableau of shape $D$ and content $\left(  1,1,\ldots
,1\right)  $ is the same thing as a standard tableau of shape $D$ (assuming
that the number of $1$'s in the content is the size of $D$).

We can finally state the \emph{Murnaghan--Nakayama rule}:

\begin{theorem}
[Murnaghan--Nakayama rule]\label{thm.char.MNrule}Let $D$ be a skew Young
diagram of size $\left\vert D\right\vert =n$. Let $\chi_{D}$ be the character
of the corresponding Specht module $\mathcal{S}^{D}$. Let $g\in S_{n}$ be a
permutation, and let $\mu=\operatorname*{type}g$ be its cycle type (or, more
generally, any tuple of nonnegative integers such that the nonzero entries of
$\mu$ are the lengths of the cycles of $g$, not necessarily in decreasing
order). Then,%
\[
\chi_{D}\left(  g\right)  =\sum_{\substack{T\text{ is a rim hook
tableau}\\\text{of shape }D\text{ and content }\mu}}\left(  -1\right)
^{\left(  \text{leg length of }T\right)  }.
\]

\end{theorem}

We will not prove this theorem here, as it would require some further
character theory that we have not introduced. A proof can be found in
\cite[Remark 3.5.18]{CSScTo10} or in \cite[Corollary 2.3.7]{Kleshc05}. (In the
case of a straight Young diagram, Theorem \ref{thm.char.MNrule} also appears,
e.g., in \cite[Theorem 7.18.5]{Stanley-EC2}, \cite[Corollary 4.10.6]{Sagan01},
\cite[Theorem 3.10]{Meliot17}. Various other texts contain recursive rules
that are easily seen to be equivalent.)

We observe a simple particular case of Theorem \ref{thm.char.MNrule}:

\begin{corollary}
\label{cor.char.MN-1cycle}Let $D$ be a skew Young diagram of size $\left\vert
D\right\vert =n>0$. Let $\chi_{D}$ be the character of the corresponding
Specht module $\mathcal{S}^{D}$. Let $g\in S_{n}$ be an $n$-cycle. Then,%
\[
\chi_{D}\left(  g\right)  =%
\begin{cases}
\left(  -1\right)  ^{\left(  \text{leg length of }D\right)  }, & \text{if
}D\text{ is a rim hook};\\
0, & \text{otherwise}.
\end{cases}
\]

\end{corollary}

\begin{exercise}
\textbf{(a)} \fbox{1} Derive Corollary \ref{cor.char.MN-1cycle} from Theorem
\ref{thm.char.MNrule}. \medskip

\textbf{(b)} \fbox{?} Prove Corollary \ref{cor.char.MN-1cycle} independently.
\end{exercise}

\subsubsection{Characters uniquely determine $S_{n}$-representations}

Let us next prove another useful fact, which is often applied to prove that
two $S_{n}$-representations are isomorphic without constructing an explicit isomorphism:

\begin{proposition}
\label{prop.char.iso-Sn}Assume that $\mathbf{k}$ is a field of characteristic
$0$. Let $V$ and $W$ be two finite-dimensional representations of $S_{n}$.
Then, $V\cong W$ if and only if $\chi_{V}=\chi_{W}$.
\end{proposition}

This proposition is a particular case of Theorem \ref{thm.char.iso=char}, but
let us outline an independent proof:

\begin{proof}
[Proof of Proposition \ref{prop.char.iso-Sn} (sketched).]The \textquotedblleft
only if\textquotedblright\ part of Proposition \ref{prop.char.iso-Sn} follows
from Proposition \ref{prop.char.iso}. Thus, it remains to prove the
\textquotedblleft if\textquotedblright\ part. So let us assume that $\chi
_{V}=\chi_{W}$. We must prove that $V\cong W$.

From $\chi_{V}=\chi_{W}$, we obtain $\mathbf{X}_{V}=\mathbf{X}_{W}$ (just as
in the proof of Corollary \ref{cor.char.XV.iso}).

For each partition $\lambda$ of $n$, we let $\mathbf{X}_{\lambda}$ denote the
element $\mathbf{X}_{\mathcal{S}^{\lambda}}$ corresponding to the Specht
module $\mathcal{S}^{\lambda}$.

Let $P$ be the set of all partitions of $n$. Then, the family $\left(
\mathbf{X}_{\lambda}\right)  _{\lambda\in P}=\left(  \mathbf{X}_{\lambda
}\right)  _{\lambda\text{ is a partition of }n}$ is a basis of the
$\mathbf{k}$-module $Z\left(  \mathbf{k}\left[  S_{n}\right]  \right)  $ (by
Corollary \ref{cor.char.Xlam-prop} \textbf{(d)}). Hence, this family is
$\mathbf{k}$-linearly independent.

Corollary \ref{cor.spechtmod.complete2} shows that%
\[
V\cong\bigoplus_{\lambda\in P}\left(  \mathcal{S}^{\lambda}\right)
^{i_{\lambda}}\ \ \ \ \ \ \ \ \ \ \text{for some integers }i_{\lambda}%
\in\mathbb{N}.
\]
Likewise, Corollary \ref{cor.spechtmod.complete2} (applied to $W$ instead of
$V$) shows that%
\[
W\cong\bigoplus_{\lambda\in P}\left(  \mathcal{S}^{\lambda}\right)
^{j_{\lambda}}\ \ \ \ \ \ \ \ \ \ \text{for some integers }j_{\lambda}%
\in\mathbb{N}.
\]
Consider these integers $i_{\lambda}$ and $j_{\lambda}$.

Corollary \ref{cor.char.V+W.m} shows that the character of a direct sum of
several $S_{n}$-representations is the sum of the characters of these
representations. Thus, from $V\cong\bigoplus\limits_{\lambda\in P}\left(
\mathcal{S}^{\lambda}\right)  ^{i_{\lambda}}$, we obtain $\chi_{V}%
=\sum\limits_{\lambda\in P}i_{\lambda}\chi_{\mathcal{S}^{\lambda}}$ (since
isomorphic $S_{n}$-representations have equal characters). Therefore, for each
$g\in S_{n}$, we have%
\begin{equation}
\chi_{V}\left(  g\right)  =\left(  \sum\limits_{\lambda\in P}i_{\lambda}%
\chi_{\mathcal{S}^{\lambda}}\right)  \left(  g\right)  =\sum\limits_{\lambda
\in P}i_{\lambda}\chi_{\mathcal{S}^{\lambda}}\left(  g\right)  .
\label{pf.prop.char.iso-Sn.4}%
\end{equation}
Now, Definition \ref{def.char.XV} yields%
\begin{align*}
\mathbf{X}_{V}  &  =\sum_{g\in S_{n}}\underbrace{\chi_{V}\left(  g\right)
}_{\substack{=\sum\limits_{\lambda\in P}i_{\lambda}\chi_{\mathcal{S}^{\lambda
}}\left(  g\right)  \\\text{(by (\ref{pf.prop.char.iso-Sn.4}))}}}g=\sum_{g\in
S_{n}}\left(  \sum\limits_{\lambda\in P}i_{\lambda}\chi_{\mathcal{S}^{\lambda
}}\left(  g\right)  \right)  g\\
&  =\sum\limits_{\lambda\in P}i_{\lambda}\underbrace{\sum_{g\in S_{n}}%
\chi_{\mathcal{S}^{\lambda}}\left(  g\right)  g}_{\substack{=\mathbf{X}%
_{\mathcal{S}^{\lambda}}\\\text{(by Definition \ref{def.char.XV})}}%
}=\sum\limits_{\lambda\in P}i_{\lambda}\underbrace{\mathbf{X}_{\mathcal{S}%
^{\lambda}}}_{=\mathbf{X}_{\lambda}}=\sum\limits_{\lambda\in P}i_{\lambda
}\mathbf{X}_{\lambda}.
\end{align*}
Likewise,
\[
\mathbf{X}_{W}=\sum\limits_{\lambda\in P}j_{\lambda}\mathbf{X}_{\lambda}.
\]
The left hand sides of these two equalities are equal (since $\mathbf{X}%
_{V}=\mathbf{X}_{W}$). Thus, so are the right hand sides. In other words,%
\[
\sum\limits_{\lambda\in P}i_{\lambda}\mathbf{X}_{\lambda}=\sum\limits_{\lambda
\in P}j_{\lambda}\mathbf{X}_{\lambda}.
\]
Since the family $\left(  \mathbf{X}_{\lambda}\right)  _{\lambda\in P}$ is
$\mathbf{k}$-linearly independent, this equality entails that all $\lambda\in
P$ satisfy $i_{\lambda}=j_{\lambda}$ in $\mathbf{k}$. Since $\mathbb{Z}$ is a
subring of $\mathbf{k}$ (because $\mathbf{k}$ is a field of characteristic
$0$), this entails that all $\lambda\in P$ satisfy $i_{\lambda}=j_{\lambda}$
in $\mathbb{Z}$ as well. Now, recall that
\begin{align*}
V  &  \cong\bigoplus_{\lambda\in P}\left(  \mathcal{S}^{\lambda}\right)
^{i_{\lambda}}=\bigoplus_{\lambda\in P}\left(  \mathcal{S}^{\lambda}\right)
^{j_{\lambda}}\ \ \ \ \ \ \ \ \ \ \left(  \text{since all }\lambda\in P\text{
satisfy }i_{\lambda}=j_{\lambda}\right) \\
&  \cong W\ \ \ \ \ \ \ \ \ \ \left(  \text{since }W\cong\bigoplus_{\lambda\in
P}\left(  \mathcal{S}^{\lambda}\right)  ^{j_{\lambda}}\right)  .
\end{align*}
This proves the \textquotedblleft if\textquotedblright\ part of Proposition
\ref{prop.char.iso-Sn}. The proof of this proposition is thus complete.
\end{proof}

\begin{remark}
Note that Proposition \ref{prop.char.iso-Sn} really requires $\mathbf{k}$ to
be a field of characteristic $0$. Merely requiring $\mathbf{k}$ to be a field
in which $n!$ is invertible would not suffice, since (e.g.) the
representations $\mathbf{k}_{\operatorname*{triv}}^{m}=\underbrace{\mathbf{k}%
_{\operatorname*{triv}}\oplus\mathbf{k}_{\operatorname*{triv}}\oplus
\cdots\oplus\mathbf{k}_{\operatorname*{triv}}}_{m\text{ times}}$ and $0$ have
the same character whenever $m=0$ in $\mathbf{k}$.
\end{remark}

\subsection{\label{sec.rep.sign-twist}The sign-twist of an $S_{n}%
$-representation}

In this short section, we will explore a simple way to transform $S_{n}%
$-representations into new ones by \textquotedblleft flipping a
sign\textquotedblright, and describe how it acts on Specht modules.

\subsubsection{The sign-twist in general}

In this section, we will oftentimes be dealing with two different
$\mathbf{k}\left[  S_{n}\right]  $-module structures on one and the same set.
Thus, it behooves us to have a notation that allows us to tell them apart.
Here is the notation that we will use:

\begin{convention}
\label{conv.sign-twist.dot1}Let $R$ be a ring, and let $M$ be an $R$-module.
Then, the action of $R$ on $M$ (that is, the map $R\times M\rightarrow M$ that
sends each pair $\left(  r,m\right)  $ to $rm$) will be denoted
$\overset{M}{\cdot}$. Thus, for any $r\in R$ and $m\in M$, we will write
$r\overset{M}{\cdot}m$ for the image of $\left(  r,m\right)  $ under this
action. Of course, this image has been hitherto denoted by $rm$, but the
notation \textquotedblleft$rm$\textquotedblright\ can be ambiguous when there
are two different actions on the same set $M$, whereas the notation
\textquotedblleft$r\overset{M}{\cdot}m$\textquotedblright\ is unambiguous in
this case, since it explicitly mentions the $R$-module $M$.
\end{convention}

We could also introduce a similar notation for the addition on $M$, but we
have no need for this, since we won't have to deal with differing additions on
the same set $M$ here. \medskip

We shall now give two equivalent definitions of the sign-twist: one for
representations of $S_{n}$ and one for left $\mathbf{k}\left[  S_{n}\right]
$-modules. As we know, a representation of $S_{n}$ (over $\mathbf{k}$) is
\textquotedblleft the same as\textquotedblright\ a left $\mathbf{k}\left[
S_{n}\right]  $-module (by Corollary \ref{cor.rep.G-rep.mod}, applied to
$G=S_{n}$), so that both definitions apply to the same objects; as we will
soon see (Proposition \ref{prop.sign-twist.exist} \textbf{(c)}), they also
produce the same objects.

Here are the two definitions:

\begin{definition}
\label{def.sign-twist.sign-twist.Sn}Let $V$ be a representation of $S_{n}$
over $\mathbf{k}$. Then, the \emph{sign-twist} of $V$ is a new representation
of $S_{n}$ over $\mathbf{k}$, which is called $V^{\operatorname*{sign}}$ and
is defined as follows:

\begin{itemize}
\item As a $\mathbf{k}$-module, $V^{\operatorname*{sign}}$ is just $V$.

\item The $S_{n}$-action on $V^{\operatorname*{sign}}$ is the map
$\overset{\operatorname*{sign}}{\rightharpoonup}\ :S_{n}\times V\rightarrow V$
defined by the formula%
\[
g\overset{\operatorname*{sign}}{\rightharpoonup}v:=\left(  -1\right)
^{g}g\rightharpoonup v\ \ \ \ \ \ \ \ \ \ \text{for all }g\in S_{n}\text{ and
}v\in V.
\]
Here, the \textquotedblleft$\rightharpoonup$\textquotedblright\ symbol on the
right hand side denotes the $S_{n}$-action on the original $S_{n}%
$-representation $V$.
\end{itemize}

This representation $V^{\operatorname*{sign}}$ really exists (by Proposition
\ref{prop.sign-twist.exist} \textbf{(a)} below), so the sign-twist of $V$ is well-defined.
\end{definition}

\begin{definition}
\label{def.sign-twist.sign-twist.kSn}Let $V$ be a left $\mathbf{k}\left[
S_{n}\right]  $-module. Then, the \emph{sign-twist} of $V$ is a new left
$\mathbf{k}\left[  S_{n}\right]  $-module, which is called
$V^{\operatorname*{sign}}$ and is defined as follows:

\begin{itemize}
\item The addition and the zero of $V^{\operatorname*{sign}}$ are the same as
those of $V$. (In other words, $V^{\operatorname*{sign}}=V$ as an additive group.)

\item The action of $\mathbf{k}\left[  S_{n}\right]  $ on
$V^{\operatorname*{sign}}$ is given by\footnotemark%
\[
\mathbf{a}\overset{V^{\operatorname*{sign}}}{\cdot}v:=T_{\operatorname*{sign}%
}\left(  \mathbf{a}\right)  \overset{V}{\cdot}v\ \ \ \ \ \ \ \ \ \ \text{for
all }\mathbf{a}\in\mathbf{k}\left[  S_{n}\right]  \text{ and }v\in
V^{\operatorname*{sign}}.
\]
In other words, the action of $\mathbf{k}\left[  S_{n}\right]  $ on
$V^{\operatorname*{sign}}$ is the map%
\begin{align*}
\mathbf{k}\left[  S_{n}\right]  \times V  &  \rightarrow V,\\
\left(  \mathbf{a},v\right)   &  \mapsto T_{\operatorname*{sign}}\left(
\mathbf{a}\right)  \overset{V}{\cdot}v.
\end{align*}

\end{itemize}

This left $\mathbf{k}\left[  S_{n}\right]  $-module $V^{\operatorname*{sign}}$
is well-defined (by Proposition \ref{prop.sign-twist.exist} \textbf{(b)} below).
\end{definition}

\footnotetext{We are using Convention \ref{conv.sign-twist.dot1} and
Definition \ref{def.Tsign.Tsign} here.}Soon (in Proposition
\ref{prop.sign-twist.exist}) we will show that these two definitions are
legitimate and equivalent. But first, let us give a simple example:

\begin{example}
\label{exa.sign-twist.ksign}Consider the trivial representation $\mathbf{k}%
_{\operatorname*{triv}}$ of $S_{n}$ constructed in Example
\ref{exa.rep.Sn-rep.triv}. Let us find its sign-twist $\left(  \mathbf{k}%
_{\operatorname*{triv}}\right)  ^{\operatorname*{sign}}$ according to
Definition \ref{def.sign-twist.sign-twist.Sn}. Indeed, as a $\mathbf{k}%
$-module, $\left(  \mathbf{k}_{\operatorname*{triv}}\right)
^{\operatorname*{sign}}=\mathbf{k}_{\operatorname*{triv}}=\mathbf{k}$, whereas
the $S_{n}$-action on $\left(  \mathbf{k}_{\operatorname*{triv}}\right)
^{\operatorname*{sign}}$ is the map $\overset{\operatorname*{sign}%
}{\rightharpoonup}\ :S_{n}\times\mathbf{k}\rightarrow\mathbf{k}$ defined by
the formula%
\[
g\overset{\operatorname*{sign}}{\rightharpoonup}v:=\left(  -1\right)
^{g}g\rightharpoonup v\ \ \ \ \ \ \ \ \ \ \text{for all }g\in S_{n}\text{ and
}v\in\mathbf{k}_{\operatorname*{triv}}.
\]
This formula can be simplified to
\[
g\overset{\operatorname*{sign}}{\rightharpoonup}v:=\left(  -1\right)
^{g}v\ \ \ \ \ \ \ \ \ \ \text{for all }g\in S_{n}\text{ and }v\in
\mathbf{k}_{\operatorname*{triv}}%
\]
(since $g\rightharpoonup v=v$ (because the $S_{n}$-action on $\mathbf{k}%
_{\operatorname*{triv}}$ is trivial)). But this is precisely the left $S_{n}%
$-action on the sign representation $\mathbf{k}_{\operatorname*{sign}}$ (as
defined in Example \ref{exa.rep.Sn-rep.sign}). Hence, we conclude that
$\left(  \mathbf{k}_{\operatorname*{triv}}\right)  ^{\operatorname*{sign}%
}=\mathbf{k}_{\operatorname*{sign}}$.

More generally, the same arguments show that if $V$ is any trivial
representation of $S_{n}$ (see Example \ref{exa.rep.Sn-rep.triv}), then its
sign-twist $V^{\operatorname*{sign}}$ is a sign representation of $S_{n}$ (see
Example \ref{exa.rep.Sn-rep.sign}).
\end{example}

It is not much harder to see that $\left(  \mathbf{k}_{\operatorname*{sign}%
}\right)  ^{\operatorname*{sign}}=\mathbf{k}_{\operatorname*{triv}}$. (This
also follows from $\left(  \mathbf{k}_{\operatorname*{triv}}\right)
^{\operatorname*{sign}}=\mathbf{k}_{\operatorname*{sign}}$ using Proposition
\ref{prop.sign-twist.twice-back}.) We can easily find the sign-twist of the
regular representation:

\begin{exercise}
\fbox{1} Consider the left regular representation $\mathbf{k}\left[
S_{n}\right]  $ of $S_{n}$. Prove that its sign-twist $\left(  \mathbf{k}%
\left[  S_{n}\right]  \right)  ^{\operatorname*{sign}}$ is isomorphic to
$\mathbf{k}\left[  S_{n}\right]  $ itself, and that the map
$T_{\operatorname*{sign}}:\mathbf{k}\left[  S_{n}\right]  \rightarrow\left(
\mathbf{k}\left[  S_{n}\right]  \right)  ^{\operatorname*{sign}}$ is an
isomorphism of $S_{n}$-representations.
\end{exercise}

We will see other examples of sign-twists below. \medskip

We still owe a proof of the well-definedness of $V^{\operatorname*{sign}}$, so
let us give it now:

\begin{proposition}
\label{prop.sign-twist.exist}Let $V$ be a left $\mathbf{k}\left[
S_{n}\right]  $-module, i.e., a representation of $S_{n}$ over $\mathbf{k}$.
Then: \medskip

\textbf{(a)} The representation $V^{\operatorname*{sign}}$ in Definition
\ref{def.sign-twist.sign-twist.Sn} is well-defined. \medskip

\textbf{(b)} The left $\mathbf{k}\left[  S_{n}\right]  $-module
$V^{\operatorname*{sign}}$ in Definition \ref{def.sign-twist.sign-twist.kSn}
is well-defined. \medskip

\textbf{(c)} The representation $V^{\operatorname*{sign}}$ in Definition
\ref{def.sign-twist.sign-twist.Sn} and the left $\mathbf{k}\left[
S_{n}\right]  $-module $V^{\operatorname*{sign}}$ in Definition
\ref{def.sign-twist.sign-twist.kSn} are the same (i.e., converting the former
into a left $\mathbf{k}\left[  S_{n}\right]  $-module using Theorem
\ref{thm.rep.G-rep.mod} yields the latter).
\end{proposition}

\begin{fineprint}
\begin{proof}
We begin with part \textbf{(b)}, since this is the simplest part, and since
establishing it will render the other two parts easier. \medskip

\textbf{(b)} Consider the map%
\begin{align*}
\mathbf{k}\left[  S_{n}\right]  \times V  &  \rightarrow V,\\
\left(  \mathbf{a},v\right)   &  \mapsto T_{\operatorname*{sign}}\left(
\mathbf{a}\right)  \overset{V}{\cdot}v.
\end{align*}
This is the map that is intended to serve as the action of $\mathbf{k}\left[
S_{n}\right]  $ on $V^{\operatorname*{sign}}$ (according to Definition
\ref{def.sign-twist.sign-twist.kSn}). Hence, let us denote this map by
$\overset{V^{\operatorname*{sign}}}{\cdot}$ (that is, let us write
$\mathbf{a}\overset{V^{\operatorname*{sign}}}{\cdot}v$ for the image of a pair
$\left(  \mathbf{a},v\right)  \in\mathbf{k}\left[  S_{n}\right]  \times V$
under this map). Thus,%
\begin{equation}
\mathbf{a}\overset{V^{\operatorname*{sign}}}{\cdot}v=T_{\operatorname*{sign}%
}\left(  \mathbf{a}\right)  \overset{V}{\cdot}v
\label{pf.prop.sign-twist.exist.b.1}%
\end{equation}
for all $\mathbf{a}\in\mathbf{k}\left[  S_{n}\right]  $ and $v\in V$. The map
$\overset{V^{\operatorname*{sign}}}{\cdot}$ is clearly $\mathbf{k}$-bilinear
(since the map $T_{\operatorname*{sign}}$ is $\mathbf{k}$-linear, and since
the action $\overset{V}{\cdot}$ of $\mathbf{k}\left[  S_{n}\right]  $ on $V$
is $\mathbf{k}$-bilinear).

In Definition \ref{def.sign-twist.sign-twist.kSn}, we defined the left
$\mathbf{k}\left[  S_{n}\right]  $-module $V^{\operatorname*{sign}}$ to be the
additive group $V$, equipped with the map $\overset{V^{\operatorname*{sign}%
}}{\cdot}$ as the action of $\mathbf{k}\left[  S_{n}\right]  $ on
$V^{\operatorname*{sign}}$. In order to prove that this left $\mathbf{k}%
\left[  S_{n}\right]  $-module $V^{\operatorname*{sign}}$ is well-defined, we
must show that this structure satisfies all the module axioms (see Definition
\ref{def.mod.leftmod}). Most of these axioms are straightforward to check
(e.g., both right and left distributivity laws follow from the $\mathbf{k}%
$-bilinearity of the map $\overset{V^{\operatorname*{sign}}}{\cdot}$), so we
only focus on two axioms: the associativity law and the \textquotedblleft%
$1m=m$\textquotedblright\ axiom. Let us check them both:

\begin{itemize}
\item \textit{The associativity law:} We must prove the associativity law for
the left $\mathbf{k}\left[  S_{n}\right]  $-module $V^{\operatorname*{sign}}$.
In other words, we must prove that $\left(  rs\right)
\overset{V^{\operatorname*{sign}}}{\cdot}m=r\overset{V^{\operatorname*{sign}%
}}{\cdot}\left(  s\overset{V^{\operatorname*{sign}}}{\cdot}m\right)  $ for all
$r,s\in\mathbf{k}\left[  S_{n}\right]  $ and $m\in V$. So let $r,s\in
\mathbf{k}\left[  S_{n}\right]  $ and $m\in V$ be arbitrary. Then, the
definition of $\overset{V^{\operatorname*{sign}}}{\cdot}$ yields
$s\overset{V^{\operatorname*{sign}}}{\cdot}m=T_{\operatorname*{sign}}\left(
s\right)  \overset{V}{\cdot}m$ and $r\overset{V^{\operatorname*{sign}}}{\cdot
}\left(  s\overset{V^{\operatorname*{sign}}}{\cdot}m\right)
=T_{\operatorname*{sign}}\left(  r\right)  \overset{V}{\cdot}\left(
s\overset{V^{\operatorname*{sign}}}{\cdot}m\right)  $ and $\left(  rs\right)
\overset{V^{\operatorname*{sign}}}{\cdot}m=T_{\operatorname*{sign}}\left(
rs\right)  \overset{V}{\cdot}m$. But $T_{\operatorname*{sign}}$ is a
$\mathbf{k}$-algebra morphism (by Theorem \ref{thm.Tsign.auto} \textbf{(a)}),
and thus we have $T_{\operatorname*{sign}}\left(  rs\right)
=T_{\operatorname*{sign}}\left(  r\right)  \cdot T_{\operatorname*{sign}%
}\left(  s\right)  $. Hence,%
\begin{align*}
\left(  rs\right)  \overset{V^{\operatorname*{sign}}}{\cdot}m  &
=\underbrace{T_{\operatorname*{sign}}\left(  rs\right)  }%
_{=T_{\operatorname*{sign}}\left(  r\right)  \cdot T_{\operatorname*{sign}%
}\left(  s\right)  }\overset{V}{\cdot}\,m=\left(  T_{\operatorname*{sign}%
}\left(  r\right)  \cdot T_{\operatorname*{sign}}\left(  s\right)  \right)
\overset{V}{\cdot}m\\
&  =T_{\operatorname*{sign}}\left(  r\right)  \overset{V}{\cdot}\left(
T_{\operatorname*{sign}}\left(  s\right)  \overset{V}{\cdot}m\right)
\end{align*}
(since the original action $\overset{V}{\cdot}$ of $\mathbf{k}\left[
S_{n}\right]  $ on $V$ satisfies the associativity law). Comparing this with
\[
r\overset{V^{\operatorname*{sign}}}{\cdot}\left(
s\overset{V^{\operatorname*{sign}}}{\cdot}m\right)  =T_{\operatorname*{sign}%
}\left(  r\right)  \overset{V}{\cdot}\left(
\underbrace{s\overset{V^{\operatorname*{sign}}}{\cdot}m}%
_{=T_{\operatorname*{sign}}\left(  s\right)  \overset{V}{\cdot}m}\right)
=T_{\operatorname*{sign}}\left(  r\right)  \overset{V}{\cdot}\left(
T_{\operatorname*{sign}}\left(  s\right)  \overset{V}{\cdot}m\right)  ,
\]
we obtain $\left(  rs\right)  \overset{V^{\operatorname*{sign}}}{\cdot
}m=r\overset{V^{\operatorname*{sign}}}{\cdot}\left(
s\overset{V^{\operatorname*{sign}}}{\cdot}m\right)  $. Thus, the associativity
law for the $\mathbf{k}\left[  S_{n}\right]  $-module $V^{\operatorname*{sign}%
}$ is proved.

\item \textit{The }$1m=m$ \textit{axiom:} We must prove the \textquotedblleft%
$1m=m$\textquotedblright\ axiom for the left $\mathbf{k}\left[  S_{n}\right]
$-module $V^{\operatorname*{sign}}$. In other words, we must prove that
$1\overset{V^{\operatorname*{sign}}}{\cdot}m=m$ for all $m\in V$ (where $1$ is
the unity of $\mathbf{k}\left[  S_{n}\right]  $). So let $m\in V$. Then, the
definition of $\overset{V^{\operatorname*{sign}}}{\cdot}$ yields
$1\overset{V^{\operatorname*{sign}}}{\cdot}m=T_{\operatorname*{sign}}\left(
1\right)  \overset{V}{\cdot}m$. But $T_{\operatorname*{sign}}$ is a
$\mathbf{k}$-algebra morphism (by Theorem \ref{thm.Tsign.auto} \textbf{(a)}),
and thus we have $T_{\operatorname*{sign}}\left(  1\right)  =1$. Hence,
$1\overset{V^{\operatorname*{sign}}}{\cdot}%
m=\underbrace{T_{\operatorname*{sign}}\left(  1\right)  }_{=1}%
\overset{V}{\cdot}\,m=1\overset{V}{\cdot}m=m$ (since the $\mathbf{k}\left[
S_{n}\right]  $-module $V$ satisfies the module axioms). Thus, we have proved
the \textquotedblleft$1m=m$\textquotedblright\ axiom for the left
$\mathbf{k}\left[  S_{n}\right]  $-module $V^{\operatorname*{sign}}$.
\end{itemize}

As we said, the other axioms are even easier to check. Thus, the left
$\mathbf{k}\left[  S_{n}\right]  $-module $V^{\operatorname*{sign}}$ is
well-defined. This proves Proposition \ref{prop.sign-twist.exist}
\textbf{(b)}. \medskip

\textbf{(a)} Instead of proving the claim directly, we shall derive it from
part \textbf{(b)} (which we have already proved).

First, we introduce some notations (exclusively for this proof):

\begin{itemize}
\item We shall write $V_{1}^{\operatorname*{sign}}$ for the representation
$V^{\operatorname*{sign}}$ defined in Definition
\ref{def.sign-twist.sign-twist.Sn} (even though we have not yet shown that it
is well-defined). Our goal is to prove that this representation $V_{1}%
^{\operatorname*{sign}}$ is well-defined.

\item We shall write $V_{2}^{\operatorname*{sign}}$ for the left
$\mathbf{k}\left[  S_{n}\right]  $-module $V^{\operatorname*{sign}}$ defined
in Definition \ref{def.sign-twist.sign-twist.kSn}. Proposition
\ref{prop.sign-twist.exist} \textbf{(b)} (which we have already proved) shows
that the latter left $\mathbf{k}\left[  S_{n}\right]  $-module $V_{2}%
^{\operatorname*{sign}}$ is well-defined. We reinterpret this left
$\mathbf{k}\left[  S_{n}\right]  $-module $V_{2}^{\operatorname*{sign}}$ as a
representation of $S_{n}$ (using Corollary \ref{cor.rep.G-rep.mod}, as usual).
\end{itemize}

We shall now show that this representation $V_{2}^{\operatorname*{sign}}$ is
precisely the (purported) representation $V_{1}^{\operatorname*{sign}}$. This
will immediately entail that the latter representation is well-defined, and
thus Proposition \ref{prop.sign-twist.exist} \textbf{(a)} will follow (since
$V_{1}^{\operatorname*{sign}}$ is precisely the representation
$V^{\operatorname*{sign}}$ defined in Definition
\ref{def.sign-twist.sign-twist.Sn}).

First, we recall that $T_{\operatorname*{sign}}$ is a $\mathbf{k}$-algebra
morphism (by Theorem \ref{thm.Tsign.auto} \textbf{(a)}), and thus we have
$T_{\operatorname*{sign}}\left(  1_{\mathbf{k}\left[  S_{n}\right]  }\right)
=1_{\mathbf{k}\left[  S_{n}\right]  }$. But $T_{\operatorname*{sign}}$ is
$\mathbf{k}$-linear. Hence, for each $\lambda\in\mathbf{k}$, we have%
\begin{equation}
T_{\operatorname*{sign}}\left(  \lambda1_{\mathbf{k}\left[  S_{n}\right]
}\right)  =\lambda\cdot\underbrace{T_{\operatorname*{sign}}\left(
1_{\mathbf{k}\left[  S_{n}\right]  }\right)  }_{=1_{\mathbf{k}\left[
S_{n}\right]  }}=\lambda1_{\mathbf{k}\left[  S_{n}\right]  }.
\label{pf.prop.sign-twist.exist.a.0}%
\end{equation}

Next, we recall that $V_{2}^{\operatorname*{sign}}$ is the left $\mathbf{k}%
\left[  S_{n}\right]  $-module $V^{\operatorname*{sign}}$ from Definition
\ref{def.sign-twist.sign-twist.kSn}. Hence, the action of $\mathbf{k}\left[
S_{n}\right]  $ on $V_{2}^{\operatorname*{sign}}$ is given by%
\begin{equation}
\mathbf{a}\overset{V_{2}^{\operatorname*{sign}}}{\cdot}%
v=T_{\operatorname*{sign}}\left(  \mathbf{a}\right)  \overset{V}{\cdot}v
\label{pf.prop.sign-twist.exist.a.1}%
\end{equation}
for all $\mathbf{a}\in\mathbf{k}\left[  S_{n}\right]  $ and $v\in V$ (by
Definition \ref{def.sign-twist.sign-twist.kSn}).

Now, we shall show that our $S_{n}$-representation $V_{2}%
^{\operatorname*{sign}}$ (which we obtained from the left $\mathbf{k}\left[
S_{n}\right]  $-module $V_{2}^{\operatorname*{sign}}$ defined in Definition
\ref{def.sign-twist.sign-twist.kSn}) is the (purported) representation
$V_{1}^{\operatorname*{sign}}$ defined in Definition
\ref{def.sign-twist.sign-twist.Sn}. Indeed:

\begin{itemize}
\item The underlying $\mathbf{k}$-module of our $S_{n}$-representation
$V_{2}^{\operatorname*{sign}}$ is the original $\mathbf{k}$-module $V$.

[\textit{Proof:} We must show that $V_{2}^{\operatorname*{sign}}=V$ as
$\mathbf{k}$-modules. By its definition, our $S_{n}$-representation
$V_{2}^{\operatorname*{sign}}$ has the same addition and zero as $V$. Thus, it
remains to show that $V_{2}^{\operatorname*{sign}}$ also has the same scaling
(by elements of $\mathbf{k}$) as $V$. In other words, it remains to prove that
every $\lambda\in\mathbf{k}$ and $v\in V$ satisfy $\lambda\overset{V_{2}%
^{\operatorname*{sign}}}{\cdot}v=\lambda\overset{V}{\cdot}v$. So let us prove
this. For any $\lambda\in\mathbf{k}$ and $v\in V$, we have%
\begin{align*}
\lambda\overset{V_{2}^{\operatorname*{sign}}}{\cdot}v  &  =\left(
\lambda1_{\mathbf{k}\left[  S_{n}\right]  }\right)  \overset{V_{2}%
^{\operatorname*{sign}}}{\cdot}v\ \ \ \ \ \ \ \ \ \ \left(
\begin{array}
[c]{c}%
\text{since the }\mathbf{k}\text{-module structure on }V_{2}%
^{\operatorname*{sign}}\text{ is}\\
\text{obtained from the }\mathbf{k}\left[  S_{n}\right]  \text{-module
structure}\\
\text{by restriction of scalars}%
\end{array}
\right) \\
&  =\underbrace{T_{\operatorname*{sign}}\left(  \lambda1_{\mathbf{k}\left[
S_{n}\right]  }\right)  }_{\substack{=\lambda1_{\mathbf{k}\left[
S_{n}\right]  }\\\text{(by (\ref{pf.prop.sign-twist.exist.a.0}))}%
}}\,\overset{V}{\cdot}\,v\ \ \ \ \ \ \ \ \ \ \left(  \text{by
(\ref{pf.prop.sign-twist.exist.a.1}), applied to }\mathbf{a}=\lambda
1_{\mathbf{k}\left[  S_{n}\right]  }\right) \\
&  =\left(  \lambda1_{\mathbf{k}\left[  S_{n}\right]  }\right)
\overset{V}{\cdot}v\\
&  =\lambda\overset{V}{\cdot}\underbrace{\left(  1_{\mathbf{k}\left[
S_{n}\right]  }\overset{V}{\cdot}v\right)  }_{\substack{=v\\\text{(by the
module axioms)}}}\ \ \ \ \ \ \ \ \ \ \left(  \text{since the action
}\overset{V}{\cdot}\text{ is }\mathbf{k}\text{-bilinear}\right) \\
&  =\lambda\overset{V}{\cdot}v.
\end{align*}
Thus, the $\mathbf{k}$-module $V_{2}^{\operatorname*{sign}}$ has the same
scaling as $V$. Since $V_{2}^{\operatorname*{sign}}$ also has the same
addition and zero as $V$, we thus conclude that the $\mathbf{k}$-module
$V_{2}^{\operatorname*{sign}}$ is identical with the $\mathbf{k}$-module $V$.]

\item The action of $S_{n}$ on our $S_{n}$-representation $V_{2}%
^{\operatorname*{sign}}$ is exactly the map $\overset{\operatorname*{sign}%
}{\rightharpoonup}\ :S_{n}\times V\rightarrow V$ defined in Definition
\ref{def.sign-twist.sign-twist.Sn}.

[\textit{Proof:} Let $\underset{2}{\overset{\operatorname*{sign}%
}{\rightharpoonup}}\ :S_{n}\times V_{2}^{\operatorname*{sign}}\rightarrow
V_{2}^{\operatorname*{sign}}$ be the action of $S_{n}$ on the representation
$V_{2}^{\operatorname*{sign}}$ of $S_{n}$. Then, for any $g\in S_{n}$ and
$v\in V_{2}^{\operatorname*{sign}}$, we have%
\begin{align*}
g\underset{2}{\overset{\operatorname*{sign}}{\rightharpoonup}}v  &
=g\overset{V_{2}^{\operatorname*{sign}}}{\cdot}%
v=\underbrace{T_{\operatorname*{sign}}\left(  g\right)  }_{\substack{=\left(
-1\right)  ^{g}g\\\text{(by the definition of }T_{\operatorname*{sign}%
}\text{)}}}\overset{V}{\cdot}\,v\ \ \ \ \ \ \ \ \ \ \left(  \text{by
(\ref{pf.prop.sign-twist.exist.a.1}), applied to }\mathbf{a}=g\right) \\
&  =\left(  \left(  -1\right)  ^{g}g\right)  \overset{V}{\cdot}v=\left(
-1\right)  ^{g}\underbrace{g\overset{V}{\cdot}v}_{\substack{=g\rightharpoonup
v\\\text{(where }\rightharpoonup\text{ denotes the action of }S_{n}\\\text{on
the original representation }V\text{ of }S_{n}\text{)}}}=\left(  -1\right)
^{g}g\rightharpoonup v.
\end{align*}
But this is precisely the formula for the map $\overset{\operatorname*{sign}%
}{\rightharpoonup}\ :S_{n}\times V\rightarrow V$ defined in Definition
\ref{def.sign-twist.sign-twist.Sn}. Thus, the map
$\underset{2}{\overset{\operatorname*{sign}}{\rightharpoonup}}\ :S_{n}\times
V_{2}^{\operatorname*{sign}}\rightarrow V_{2}^{\operatorname*{sign}}$ is
precisely the map $\overset{\operatorname*{sign}}{\rightharpoonup}%
\ :S_{n}\times V\rightarrow V$ defined in Definition
\ref{def.sign-twist.sign-twist.Sn}. In other words, the action of $S_{n}$ on
our $S_{n}$-representation $V_{2}^{\operatorname*{sign}}$ is exactly the map
$\overset{\operatorname*{sign}}{\rightharpoonup}\ :S_{n}\times V\rightarrow V$
defined in Definition \ref{def.sign-twist.sign-twist.Sn}.]
\end{itemize}

Thus, we have shown that our representation $V_{2}^{\operatorname*{sign}}$ of
$S_{n}$ has the same underlying $\mathbf{k}$-module as $V$, and its action of
$S_{n}$ is precisely the map $\overset{\operatorname*{sign}}{\rightharpoonup
}\ :S_{n}\times V\rightarrow V$ from Definition
\ref{def.sign-twist.sign-twist.Sn}. Consequently, this representation
$V_{2}^{\operatorname*{sign}}$ is exactly the representation
$V^{\operatorname*{sign}}$ defined in Definition
\ref{def.sign-twist.sign-twist.Sn} (i.e., the representation that we now call
$V_{1}^{\operatorname*{sign}}$). Therefore, the latter representation
$V^{\operatorname*{sign}}$ is well-defined (since the former representation is
well-defined). This proves Proposition \ref{prop.sign-twist.exist}
\textbf{(a)}. \medskip

\textbf{(c)} We must prove that the representation $V^{\operatorname*{sign}}$
in Definition \ref{def.sign-twist.sign-twist.Sn} and the left $\mathbf{k}%
\left[  S_{n}\right]  $-module $V^{\operatorname*{sign}}$ in Definition
\ref{def.sign-twist.sign-twist.kSn} are the same. In our above proof of part
\textbf{(a)}, we have denoted the former representation by $V_{1}%
^{\operatorname*{sign}}$ and the latter $\mathbf{k}\left[  S_{n}\right]
$-module by $V_{2}^{\operatorname*{sign}}$; thus, we must prove that
$V_{1}^{\operatorname*{sign}}$ and $V_{2}^{\operatorname*{sign}}$ are the same.

But we have essentially done this already. Indeed, in our above proof of part
\textbf{(a)}, we have shown that the representation $V_{2}%
^{\operatorname*{sign}}$ is exactly the representation
$V^{\operatorname*{sign}}$ defined in Definition
\ref{def.sign-twist.sign-twist.Sn}; but the latter representation is precisely
the one that we called $V_{1}^{\operatorname*{sign}}$. Hence, we have shown
that $V_{1}^{\operatorname*{sign}}$ and $V_{2}^{\operatorname*{sign}}$ are the
same. This proves Proposition \ref{prop.sign-twist.exist} \textbf{(c)}.
\end{proof}
\end{fineprint}

Proposition \ref{prop.sign-twist.exist} \textbf{(a)} and \textbf{(b)} shows
that our two definitions of $V^{\operatorname*{sign}}$ are well-defined, and
Proposition \ref{prop.sign-twist.exist} \textbf{(c)} shows that they are
equivalent. We now turn to some simple properties of sign-twists:

\begin{proposition}
\label{prop.sign-twist.twice-back}Let $V$ be a representation of $S_{n}$.
Then, $\left(  V^{\operatorname*{sign}}\right)  ^{\operatorname*{sign}}=V$.
(This is saying that the $S_{n}$-representation $\left(
V^{\operatorname*{sign}}\right)  ^{\operatorname*{sign}}$ is literally
identical to $V$ -- not just isomorphic.)
\end{proposition}

\begin{fineprint}
\begin{proof}
We shall use Definition \ref{def.sign-twist.sign-twist.kSn} as a definition of
$V^{\operatorname*{sign}}$ and also of $\left(  V^{\operatorname*{sign}%
}\right)  ^{\operatorname*{sign}}$. The three left $\mathbf{k}\left[
S_{n}\right]  $-modules $V$ and $V^{\operatorname*{sign}}$ and $\left(
V^{\operatorname*{sign}}\right)  ^{\operatorname*{sign}}$ all have the same
underlying additive group (by Definition \ref{def.sign-twist.sign-twist.kSn}),
but usually don't all have the same action of $\mathbf{k}\left[  S_{n}\right]
$. However, we claim that the first and the third of them (i.e., the left
$\mathbf{k}\left[  S_{n}\right]  $-modules $V$ and $\left(
V^{\operatorname*{sign}}\right)  ^{\operatorname*{sign}}$) do have the same
action of $\mathbf{k}\left[  S_{n}\right]  $.

Indeed, the definition of the action of $\mathbf{k}\left[  S_{n}\right]  $ on
$V^{\operatorname*{sign}}$ (see Definition \ref{def.sign-twist.sign-twist.kSn}%
) shows that%
\begin{equation}
\mathbf{a}\overset{V^{\operatorname*{sign}}}{\cdot}v=T_{\operatorname*{sign}%
}\left(  \mathbf{a}\right)  \overset{V}{\cdot}v
\label{pf.prop.sign-twist.twice-back.1}%
\end{equation}
for all $\mathbf{a}\in\mathbf{k}\left[  S_{n}\right]  $ and $v\in
V^{\operatorname*{sign}}$. Applying the same reasoning to
$V^{\operatorname*{sign}}$ instead of $V$, we obtain
\begin{equation}
\mathbf{a}\overset{\left(  V^{\operatorname*{sign}}\right)
^{\operatorname*{sign}}}{\cdot}v=T_{\operatorname*{sign}}\left(
\mathbf{a}\right)  \overset{V^{\operatorname*{sign}}}{\cdot}v
\label{pf.prop.sign-twist.twice-back.2}%
\end{equation}
for all $\mathbf{a}\in\mathbf{k}\left[  S_{n}\right]  $ and $v\in\left(
V^{\operatorname*{sign}}\right)  ^{\operatorname*{sign}}$. Now, for any
$\mathbf{a}\in\mathbf{k}\left[  S_{n}\right]  $ and $v\in V$, we have $v\in
V=V^{\operatorname*{sign}}=\left(  V^{\operatorname*{sign}}\right)
^{\operatorname*{sign}}$ and therefore%
\begin{align*}
\mathbf{a}\overset{\left(  V^{\operatorname*{sign}}\right)
^{\operatorname*{sign}}}{\cdot}v  &  =T_{\operatorname*{sign}}\left(
\mathbf{a}\right)  \overset{V^{\operatorname*{sign}}}{\cdot}%
v\ \ \ \ \ \ \ \ \ \ \left(  \text{by (\ref{pf.prop.sign-twist.twice-back.2}%
)}\right) \\
&  =\underbrace{T_{\operatorname*{sign}}\left(  T_{\operatorname*{sign}%
}\left(  \mathbf{a}\right)  \right)  }_{\substack{=\left(
T_{\operatorname*{sign}}\circ T_{\operatorname*{sign}}\right)  \left(
\mathbf{a}\right)  \\=\mathbf{a}\\\text{(since Theorem \ref{thm.Tsign.auto}
\textbf{(b)}}\\\text{says that }T_{\operatorname*{sign}}\circ
T_{\operatorname*{sign}}=\operatorname*{id}\text{)}}}\overset{V}{\cdot
}\,v\ \ \ \ \ \ \ \ \ \ \left(
\begin{array}
[c]{c}%
\text{by (\ref{pf.prop.sign-twist.twice-back.1}), applied to }%
T_{\operatorname*{sign}}\left(  \mathbf{a}\right) \\
\text{instead of }\mathbf{a}%
\end{array}
\right) \\
&  =\mathbf{a}\overset{V}{\cdot}v.
\end{align*}
This shows that the maps $\overset{\left(  V^{\operatorname*{sign}}\right)
^{\operatorname*{sign}}}{\cdot}$ and $\overset{V}{\cdot}$ are equal. In other
words, the actions of $\mathbf{k}\left[  S_{n}\right]  $ on the left
$\mathbf{k}\left[  S_{n}\right]  $-modules $\left(  V^{\operatorname*{sign}%
}\right)  ^{\operatorname*{sign}}$ and $V$ are equal (since these actions are
$\overset{\left(  V^{\operatorname*{sign}}\right)  ^{\operatorname*{sign}%
}}{\cdot}$ and $\overset{V}{\cdot}$, respectively).

Now, the left $\mathbf{k}\left[  S_{n}\right]  $-modules $\left(
V^{\operatorname*{sign}}\right)  ^{\operatorname*{sign}}$ and $V$ have the
same underlying additive group (as we saw above), but also have the same
action of $\mathbf{k}\left[  S_{n}\right]  $ (since we just showed that the
actions of $\mathbf{k}\left[  S_{n}\right]  $ on the left $\mathbf{k}\left[
S_{n}\right]  $-modules $\left(  V^{\operatorname*{sign}}\right)
^{\operatorname*{sign}}$ and $V$ are equal). Hence, these two left
$\mathbf{k}\left[  S_{n}\right]  $-modules are completely identical. In other
words, $\left(  V^{\operatorname*{sign}}\right)  ^{\operatorname*{sign}}=V$.
This proves Proposition \ref{prop.sign-twist.twice-back}.
\end{proof}
\end{fineprint}

We give two more basic properties of sign-twists without proof:

\begin{proposition}
\label{prop.sign-twist.dirsums}Let $V$ and $W$ be two representations of
$S_{n}$. Then, $\left(  V\oplus W\right)  ^{\operatorname*{sign}%
}=V^{\operatorname*{sign}}\oplus W^{\operatorname*{sign}}$.
\end{proposition}

\begin{proposition}
\label{prop.sign-twist.mors}Let $V$ and $W$ be two representations of $S_{n}$.
Let $f:V\rightarrow W$ be a map. Then: \medskip

\textbf{(a)} The map $f$ is a morphism of representations from $V$ to $W$ if
and only if $f$ is a morphism of representations from $V^{\operatorname*{sign}%
}$ to $W^{\operatorname*{sign}}$. \medskip

\textbf{(b)} The map $f$ is an isomorphism of representations from $V$ to $W$
if and only if $f$ is an isomorphism of representations from
$V^{\operatorname*{sign}}$ to $W^{\operatorname*{sign}}$.
\end{proposition}

\begin{corollary}
\label{cor.sign-twist.iso}Let $V$ and $W$ be two representations of $S_{n}$.
If $V\cong W$, then $V^{\operatorname*{sign}}\cong W^{\operatorname*{sign}}$.
\end{corollary}

\begin{exercise}
\textbf{(a)} \fbox{1} Prove Proposition \ref{prop.sign-twist.dirsums}.
\medskip

\textbf{(b)} \fbox{2} Prove Proposition \ref{prop.sign-twist.mors} and
Corollary \ref{cor.sign-twist.iso}.
\end{exercise}

\begin{remark}
\label{rmk.sign-twist.functor}Proposition \ref{prop.sign-twist.mors}
\textbf{(a)} shows that we can define a \textquotedblleft sign-twist
functor\textquotedblright\ on the category of $S_{n}$-representations. This is
a (covariant) functor from the category of $S_{n}$-representations to itself.
On objects, it acts by sending each $S_{n}$-representation $V$ to
$V^{\operatorname*{sign}}$. On morphisms, it does nothing (i.e., it sends any
morphism $f$ to $f$ itself).
\end{remark}

\begin{remark}
The definition of the sign-twist (Definition
\ref{def.sign-twist.sign-twist.kSn}) can be generalized: Let $R$ and
$R^{\prime}$ be two rings, and let $\varphi:R\rightarrow R^{\prime}$ be any
ring morphism. Then, any left $R^{\prime}$-module $V$ can be turned into a
left $R$-module by setting%
\[
r\cdot v:=\varphi\left(  r\right)  \cdot v\ \ \ \ \ \ \ \ \ \ \text{for all
}r\in R\text{ and }v\in V
\]
(where the \textquotedblleft$\cdot$\textquotedblright\ on the left hand side
means the action of $R$ on $V$ that we are defining, whereas the
\textquotedblleft$\cdot$\textquotedblright\ on the right hand side means the
existing action of $R^{\prime}$ on $V$). The addition on this new $R$-module
$V$ is the same as on the original $R^{\prime}$-module $V$. We denote this new
$R$-module $V$ by $V^{\varphi}$.

This way of transforming $R^{\prime}$-modules $V$ into $R$-modules
$V^{\varphi}$ (with the same ground set and the same addition) is known as
\emph{restriction of scalars}, and has already been introduced in Definition
\ref{def.mod.R-to-k} in a slightly less general context (with the letters
$\mathbf{k}$ and $R$ used instead of $R$ and $R^{\prime}$). Nothing prevents
us from taking $R=\mathbf{k}\left[  S_{n}\right]  $ and $R^{\prime}%
=\mathbf{k}\left[  S_{n}\right]  $ here, and letting $\varphi$ be the ring
morphism $T_{\operatorname*{sign}}:\mathbf{k}\left[  S_{n}\right]
\rightarrow\mathbf{k}\left[  S_{n}\right]  $. If we do so, then $V^{\varphi
}=V^{T_{\operatorname*{sign}}}$ is precisely the sign-twist
$V^{\operatorname*{sign}}$ defined in Definition
\ref{def.sign-twist.sign-twist.kSn}. Many properties of sign-twists are
particular cases of general facts about restriction of scalars. For example,
if $\varphi:R\rightarrow R^{\prime}$ and $\psi:R^{\prime}\rightarrow
R^{\prime\prime}$ are two ring morphisms, and if $V$ is a left $R^{\prime
\prime}$-module, then $\left(  V^{\psi}\right)  ^{\varphi}=V^{\psi\circ
\varphi}$. Applying this to $R=\mathbf{k}\left[  S_{n}\right]  $ and
$R^{\prime}=\mathbf{k}\left[  S_{n}\right]  $ and $R^{\prime\prime}%
=\mathbf{k}\left[  S_{n}\right]  $ and $\varphi=T_{\operatorname*{sign}}$ and
$\psi=T_{\operatorname*{sign}}$, we obtain $\left(  V^{T_{\operatorname*{sign}%
}}\right)  ^{T_{\operatorname*{sign}}}=V^{T_{\operatorname*{sign}}\circ
T_{\operatorname*{sign}}}=V^{\operatorname*{id}}$ (since
$T_{\operatorname*{sign}}\circ T_{\operatorname*{sign}}=\operatorname*{id}$).
Since $V^{\operatorname*{id}}$ is easily seen to be just $V$, whereas
$V^{T_{\operatorname*{sign}}}$ is the sign-twist $V^{\operatorname*{sign}}$,
this is precisely the claim of Proposition \ref{prop.sign-twist.twice-back}.
\end{remark}

\begin{fineprint}
\begin{remark}
There is yet another way of defining the sign-twist $V^{\operatorname*{sign}}$
of an $S_{n}$-representation $V$: Namely, it can be viewed as an instance of
the \emph{internal tensor product}, namely the tensor product $\mathbf{k}%
_{\operatorname*{sign}}\otimes V$. Personally, I find this definition
unnecessarily convoluted (internal tensor products in general are rather
complicated), but it is worth knowing about its existence, so here is an outline:

The \emph{internal tensor product} of two $S_{n}$-representations $U$ and $V$
is defined to be the $S_{n}$-representation which (as a $\mathbf{k}$-module)
is the tensor product $U\otimes V$ and whose $S_{n}$-action is given by the
rule%
\[
g\rightharpoonup\left(  u\otimes v\right)  =\left(  g\rightharpoonup u\right)
\otimes\left(  g\rightharpoonup v\right)  \ \ \ \ \ \ \ \ \ \ \text{for all
}g\in S_{n}\text{ and }u\in U\text{ and }v\in V.
\]
(This rule determines the $S_{n}$-action on all pure tensors, and thus -- by
linearity -- on all tensors.) It is not hard to see that this really is a
well-defined $S_{n}$-action (indeed, this works just as well for any group $G$
instead of $S_{n}$), and that each $S_{n}$-representation $V$ satisfies
$\mathbf{k}_{\operatorname*{triv}}\otimes V\cong V$ (by the isomorphism
$\lambda\otimes v\mapsto\lambda v$) and $\mathbf{k}_{\operatorname*{sign}%
}\otimes V\cong V^{\operatorname*{sign}}$ (by the isomorphism $\lambda\otimes
v\mapsto\lambda v$ again). The latter isomorphism allows us to work with
$\mathbf{k}_{\operatorname*{sign}}\otimes V$ instead of
$V^{\operatorname*{sign}}$, thus sidestepping the need for defining
$V^{\operatorname*{sign}}$.
\end{remark}
\end{fineprint}

\subsubsection{The sign-twist of a Specht module}

We have seen some examples of sign-twists in Example
\ref{exa.sign-twist.ksign} and shortly thereafter. Let us now state the main
result of this section:

\begin{theorem}
\label{thm.sign-twist.Slam}Let $\lambda$ be any partition of $n$. Set
$\mathcal{S}^{\lambda}:=\mathcal{S}^{Y\left(  \lambda\right)  }$ (the Specht
module for shape $Y\left(  \lambda\right)  $). Let $\lambda^{t}$ be the
conjugate partition of $\lambda$ (defined in Theorem \ref{thm.partitions.conj}%
). Let $h^{\lambda}:=\dfrac{n!}{f^{\lambda}}$.

Assume that $h^{\lambda}$ is invertible in $\mathbf{k}$. Then, the sign-twist
$\left(  \mathcal{S}^{\lambda}\right)  ^{\operatorname*{sign}}$ of the Specht
module $\mathcal{S}^{\lambda}$ is isomorphic to the Specht module
$\mathcal{S}^{\lambda^{t}}$.
\end{theorem}

For example, for $n=4$ and $\lambda=\left(  3,1\right)  $, this is saying that
$\left(  \mathcal{S}^{\left(  3,1\right)  }\right)  ^{\operatorname*{sign}%
}\cong\mathcal{S}^{\left(  3,1\right)  ^{t}}=\mathcal{S}^{\left(
2,1,1\right)  }$ as long as $h^{\left(  3,1\right)  }=\dfrac{4!}{f^{\left(
3,1\right)  }}=8$ is invertible in $\mathbf{k}$.

To prove Theorem \ref{thm.sign-twist.Slam}, we recall that we denote the
$\mathbf{k}$-algebra $\mathbf{k}\left[  S_{n}\right]  $ by $\mathcal{A}$. Now,
any left ideal $\mathcal{A}\mathbf{a}$ of $\mathcal{A}$ (with $\mathbf{a}%
\in\mathcal{A}$) is a left $\mathcal{A}$-module, that is, a left
$\mathbf{k}\left[  S_{n}\right]  $-module, i.e., a representation of $S_{n}$.
Hence, it has a well-defined sign-twist $\left(  \mathcal{A}\mathbf{a}\right)
^{\operatorname*{sign}}$. We claim the following general formula for this sign-twist:

\begin{lemma}
\label{lem.sign-twist.a}Let $\mathbf{a}\in\mathcal{A}$. Then,
\[
\left(  \mathcal{A}\mathbf{a}\right)  ^{\operatorname*{sign}}\cong%
\mathcal{A}\cdot T_{\operatorname*{sign}}\left(  \mathbf{a}\right)
\ \ \ \ \ \ \ \ \ \ \text{as }S_{n}\text{-representations.}%
\]

\end{lemma}

\begin{proof}
We will use an easy general fact about $\mathbf{k}$-algebras:

\begin{statement}
\textit{Claim 1:} Let $R$ and $R^{\prime}$ be two $\mathbf{k}$-algebras. Let
$\varphi:R\rightarrow R^{\prime}$ be a $\mathbf{k}$-algebra isomorphism. Let
$r\in R$. Then, the restriction of $\varphi$ to the left ideal $Rr$ is a
$\mathbf{k}$-module isomorphism from $Rr$ to $R^{\prime}\cdot\varphi\left(
r\right)  $.
\end{statement}

The map $T_{\operatorname*{sign}}:\mathcal{A}\rightarrow\mathcal{A}$ is a
$\mathbf{k}$-algebra isomorphism (by Theorem \ref{thm.Tsign.auto}
\textbf{(a)}, since $\mathcal{A}=\mathbf{k}\left[  S_{n}\right]  $). Hence,
Claim 1 (applied to $R=\mathcal{A}$ and $R^{\prime}=\mathcal{A}$ and
$\varphi=T_{\operatorname*{sign}}$ and $r=\mathbf{a}$) shows that the
restriction of $T_{\operatorname*{sign}}$ to the left ideal $\mathcal{A}%
\mathbf{a}$ is a $\mathbf{k}$-module isomorphism from $\mathcal{A}\mathbf{a}$
to $\mathcal{A}\cdot T_{\operatorname*{sign}}\left(  \mathbf{a}\right)  $. Let
us denote this restriction by $f$. Thus, $f:\mathcal{A}\mathbf{a}%
\rightarrow\mathcal{A}\cdot T_{\operatorname*{sign}}\left(  \mathbf{a}\right)
$ is a $\mathbf{k}$-module isomorphism from $\mathcal{A}\mathbf{a}$ to
$\mathcal{A}\cdot T_{\operatorname*{sign}}\left(  \mathbf{a}\right)  $, and
hence is invertible. Moreover,
\begin{equation}
f\left(  \mathbf{u}\right)  =T_{\operatorname*{sign}}\left(  \mathbf{u}%
\right)  \ \ \ \ \ \ \ \ \ \ \text{for each }\mathbf{u}\in\mathcal{A}%
\mathbf{a} \label{pf.lem.sign-twist.a.fu=}%
\end{equation}
(since $f$ is a restriction of $T_{\operatorname*{sign}}$).

The map $f$ is $\mathbf{k}$-linear, but not left $\mathcal{A}$-linear -- at
least not when regarded as a map from $\mathcal{A}\mathbf{a}$ to
$\mathcal{A}\cdot T_{\operatorname*{sign}}\left(  \mathbf{a}\right)  $.
However, we can also regard $f$ as a map from the sign-twist $\left(
\mathcal{A}\mathbf{a}\right)  ^{\operatorname*{sign}}$ to $\mathcal{A}\cdot
T_{\operatorname*{sign}}\left(  \mathbf{a}\right)  $ (since $\left(
\mathcal{A}\mathbf{a}\right)  ^{\operatorname*{sign}}=\mathcal{A}\mathbf{a}$
as sets), and then the situation improves:

\begin{statement}
\textit{Claim 2:} The map $f$ is a left $\mathcal{A}$-linear map from the
sign-twist $\left(  \mathcal{A}\mathbf{a}\right)  ^{\operatorname*{sign}}$ to
$\mathcal{A}\cdot T_{\operatorname*{sign}}\left(  \mathbf{a}\right)  $.
\end{statement}

\begin{proof}
[Proof of Claim 2.]Let us denote the left $\mathcal{A}$-module $\mathcal{A}%
\mathbf{a}$ by $V$. Thus, $\mathcal{A}\mathbf{a}=V$. Definition
\ref{def.sign-twist.sign-twist.kSn} yields that
\begin{align}
\mathbf{b}\overset{V^{\operatorname*{sign}}}{\cdot}v  &
=T_{\operatorname*{sign}}\left(  \mathbf{b}\right)  \overset{V}{\cdot
}v\label{pf.lem.sign-twist.a.fu=.C2.pf.1}\\
&  \ \ \ \ \ \ \ \ \ \ \text{for all }\mathbf{b}\in\mathbf{k}\left[
S_{n}\right]  \text{ and }v\in V^{\operatorname*{sign}}\nonumber
\end{align}
(using Convention \ref{conv.sign-twist.dot1}).

Let us denote the left $\mathcal{A}$-module $\mathcal{A}\cdot
T_{\operatorname*{sign}}\left(  \mathbf{a}\right)  $ by $W$. Thus,
$\mathcal{A}\cdot T_{\operatorname*{sign}}\left(  \mathbf{a}\right)  =W$.

We must show that $f$ is a left $\mathcal{A}$-linear map from $\left(
\mathcal{A}\mathbf{a}\right)  ^{\operatorname*{sign}}$ to $\mathcal{A}\cdot
T_{\operatorname*{sign}}\left(  \mathbf{a}\right)  $. In other words, we must
show that $f$ is a left $\mathcal{A}$-linear map from $V^{\operatorname*{sign}%
}$ to $W$ (since $\mathcal{A}\mathbf{a}=V$ and $\mathcal{A}\cdot
T_{\operatorname*{sign}}\left(  \mathbf{a}\right)  =W$). Since we already know
that $f$ respects addition and zero (because $f$ is $\mathbf{k}$-linear), we
only need to show that%
\[
f\left(  \mathbf{b}\overset{V^{\operatorname*{sign}}}{\cdot}v\right)
=\mathbf{b}\overset{W}{\cdot}f\left(  v\right)  \ \ \ \ \ \ \ \ \ \ \text{for
all }\mathbf{b}\in\mathcal{A}\text{ and }v\in V^{\operatorname*{sign}}.
\]
So let us show this. Let $\mathbf{b}\in\mathcal{A}$ and $v\in
V^{\operatorname*{sign}}$. Thus, $\mathbf{b}\in\mathcal{A}=\mathbf{k}\left[
S_{n}\right]  $, so that (\ref{pf.lem.sign-twist.a.fu=.C2.pf.1}) yields
$\mathbf{b}\overset{V^{\operatorname*{sign}}}{\cdot}v=T_{\operatorname*{sign}%
}\left(  \mathbf{b}\right)  \overset{V}{\cdot}v$. Furthermore, $v\in
V^{\operatorname*{sign}}=V$ (we are talking about equality of sets here, not
of $\mathcal{A}$-modules), so that $v\in V=\mathcal{A}\mathbf{a}$.

But $V$ is a left ideal of $\mathcal{A}$ (since $V=\mathcal{A}\mathbf{a}$);
thus, the action of $\mathcal{A}$ on $V$ is just multiplication inside
$\mathcal{A}$. In other words, $\mathbf{c}\overset{V}{\cdot}w=\mathbf{c}w$ for
all $\mathbf{c}\in\mathcal{A}$ and $w\in V$. Hence, in particular,
$T_{\operatorname*{sign}}\left(  \mathbf{b}\right)  \overset{V}{\cdot
}v=T_{\operatorname*{sign}}\left(  \mathbf{b}\right)  \cdot v$. Similarly, we
find $\mathbf{b}\overset{W}{\cdot}f\left(  v\right)  =\mathbf{b}\cdot f\left(
v\right)  $ (since $W$ is also a left ideal of $\mathcal{A}$). Now,%
\[
\mathbf{b}\overset{V^{\operatorname*{sign}}}{\cdot}v=T_{\operatorname*{sign}%
}\left(  \mathbf{b}\right)  \overset{V}{\cdot}v=T_{\operatorname*{sign}%
}\left(  \mathbf{b}\right)  \cdot v.
\]
But $\mathbf{b}\overset{V^{\operatorname*{sign}}}{\cdot}v\in
V^{\operatorname*{sign}}=V=\mathcal{A}\mathbf{a}$. Hence,
(\ref{pf.lem.sign-twist.a.fu=}) (applied to $\mathbf{u}=\mathbf{b}%
\overset{V^{\operatorname*{sign}}}{\cdot}v$) yields
\begin{align*}
f\left(  \mathbf{b}\overset{V^{\operatorname*{sign}}}{\cdot}v\right)   &
=T_{\operatorname*{sign}}\left(  \underbrace{\mathbf{b}%
\overset{V^{\operatorname*{sign}}}{\cdot}v}_{=T_{\operatorname*{sign}}\left(
\mathbf{b}\right)  \cdot v}\right)  =T_{\operatorname*{sign}}\left(
T_{\operatorname*{sign}}\left(  \mathbf{b}\right)  \cdot v\right) \\
&  =\underbrace{T_{\operatorname*{sign}}\left(  T_{\operatorname*{sign}%
}\left(  \mathbf{b}\right)  \right)  }_{\substack{=\left(
T_{\operatorname*{sign}}\circ T_{\operatorname*{sign}}\right)  \left(
\mathbf{b}\right)  \\=\mathbf{b}\\\text{(since Theorem \ref{thm.Tsign.auto}
\textbf{(b)}}\\\text{says that }T_{\operatorname*{sign}}\circ
T_{\operatorname*{sign}}=\operatorname*{id}\text{)}}}\cdot
\,T_{\operatorname*{sign}}\left(  v\right)  \ \ \ \ \ \ \ \ \ \ \left(
\begin{array}
[c]{c}%
\text{since }T_{\operatorname*{sign}}\text{ is a }\mathbf{k}\text{-algebra}\\
\text{morphism}%
\end{array}
\right) \\
&  =\mathbf{b}\cdot T_{\operatorname*{sign}}\left(  v\right)  .
\end{align*}
Comparing this with%
\[
\mathbf{b}\overset{W}{\cdot}f\left(  v\right)  =\mathbf{b}\cdot
\underbrace{f\left(  v\right)  }_{\substack{=T_{\operatorname*{sign}}\left(
v\right)  \\\text{(by (\ref{pf.lem.sign-twist.a.fu=}),}\\\text{applied to
}\mathbf{u}=v\text{)}}}=\mathbf{b}\cdot T_{\operatorname*{sign}}\left(
v\right)  ,
\]
we obtain $f\left(  \mathbf{b}\overset{V^{\operatorname*{sign}}}{\cdot
}v\right)  =\mathbf{b}\overset{W}{\cdot}f\left(  v\right)  $.

Forget that we fixed $\mathbf{b}$ and $v$. We thus have proved that%
\[
f\left(  \mathbf{b}\overset{V^{\operatorname*{sign}}}{\cdot}v\right)
=\mathbf{b}\overset{W}{\cdot}f\left(  v\right)  \ \ \ \ \ \ \ \ \ \ \text{for
all }\mathbf{b}\in\mathcal{A}\text{ and }v\in V^{\operatorname*{sign}}.
\]
As explained above, this completes the proof of Claim 2.
\end{proof}

Now, Claim 2 shows that the map $f$ is a left $\mathcal{A}$-linear map from
the sign-twist $\left(  \mathcal{A}\mathbf{a}\right)  ^{\operatorname*{sign}}$
to $\mathcal{A}\cdot T_{\operatorname*{sign}}\left(  \mathbf{a}\right)  $.
Since $f$ is invertible, we thus conclude that $f$ is a left $\mathcal{A}%
$-module isomorphism (since any invertible left $\mathcal{A}$-linear map is a
left $\mathcal{A}$-module isomorphism). Hence, $\left(  \mathcal{A}%
\mathbf{a}\right)  ^{\operatorname*{sign}}\cong\mathcal{A}\cdot
T_{\operatorname*{sign}}\left(  \mathbf{a}\right)  $ as left $\mathcal{A}%
$-modules, i.e., as left $\mathbf{k}\left[  S_{n}\right]  $-modules, i.e., as
$S_{n}$-representations. This proves Lemma \ref{lem.sign-twist.a}.
\end{proof}

Our next lemma is a general \textquotedblleft left ideal
avatar\textquotedblright\ for the sign-twist of a Specht module:

\begin{lemma}
\label{lem.sign-twist.AFT}Let $D$ be a diagram with $\left\vert D\right\vert
=n$. Let $T$ be an $n$-tableau of shape $D$. Let $\mathbf{F}_{T}%
:=\nabla_{\operatorname*{Row}T}\nabla_{\operatorname*{Col}T}^{-}\in
\mathcal{A}$. Then,%
\[
\mathcal{S}^{D}\cong\mathcal{A}\mathbf{E}_{T}\ \ \ \ \ \ \ \ \ \ \text{and}%
\ \ \ \ \ \ \ \ \ \ \left(  \mathcal{S}^{\mathbf{r}\left(  D\right)  }\right)
^{\operatorname*{sign}}\cong\mathcal{A}\mathbf{F}_{T}%
\]
as $S_{n}$-representations. (See Theorem \ref{thm.partitions.conj} for the
meaning of $\mathbf{r}$.)
\end{lemma}

\begin{proof}
From (\ref{eq.def.specht.ET.defs.SD=AET}), we know that $\mathcal{S}^{D}%
\cong\mathcal{A}\mathbf{E}_{T}$ as left $\mathcal{A}$-modules, i.e., as
$S_{n}$-representations (since $\mathcal{A}=\mathbf{k}\left[  S_{n}\right]
$). It remains to prove that $\left(  \mathcal{S}^{\mathbf{r}\left(  D\right)
}\right)  ^{\operatorname*{sign}}\cong\mathcal{A}\mathbf{F}_{T}$.

A straightforward generalization of Proposition \ref{prop.specht.FT.basics}
\textbf{(a)} (with the same proof) says that the antipode $S$ of
$\mathbf{k}\left[  S_{n}\right]  $ satisfies $S\left(  \mathbf{E}_{T}\right)
=\mathbf{F}_{T}$.

\begin{noncompile}
Applying the map $S$ to this equality, we obtain $S\left(  S\left(
\mathbf{E}_{T}\right)  \right)  =S\left(  \mathbf{F}_{T}\right)  $. But
$S\left(  S\left(  \mathbf{E}_{T}\right)  \right)  =\mathbf{E}_{T}$ (since
Theorem \ref{thm.S.auto} \textbf{(b)} says that $S\circ S=\operatorname*{id}%
$). Comparing these two equalities, we find $\mathbf{E}_{T}=S\left(
\mathbf{F}_{T}\right)  $.
\end{noncompile}

We know that $T$ is an $n$-tableau of shape $D$. By Definition
\ref{def.tableaux.r}, we conclude that $T\mathbf{r}$ is an $n$-tableau of
shape $\mathbf{r}\left(  D\right)  $. Hence,
(\ref{eq.def.specht.ET.defs.SD=AET}) (applied to $\mathbf{r}\left(  D\right)
$ and $T\mathbf{r}$ instead of $D$ and $T$) shows that $\mathcal{S}%
^{\mathbf{r}\left(  D\right)  }\cong\mathcal{A}\mathbf{E}_{T\mathbf{r}}$ as
left $\mathcal{A}$-modules.

Furthermore, (\ref{eq.prop.specht.ET.r.3}) (applied to $P=T$) yields
$\mathbf{E}_{T\mathbf{r}}=T_{\operatorname*{sign}}\left(  S\left(
\mathbf{E}_{T}\right)  \right)  $. In view of $S\left(  \mathbf{E}_{T}\right)
=\mathbf{F}_{T}$, we can rewrite this as $\mathbf{E}_{T\mathbf{r}%
}=T_{\operatorname*{sign}}\left(  \mathbf{F}_{T}\right)  $. Applying the map
$T_{\operatorname*{sign}}$ to this equality, we obtain
\[
T_{\operatorname*{sign}}\left(  \mathbf{E}_{T\mathbf{r}}\right)
=T_{\operatorname*{sign}}\left(  T_{\operatorname*{sign}}\left(
\mathbf{F}_{T}\right)  \right)  =\mathbf{F}_{T}%
\]
(since Theorem \ref{thm.Tsign.auto} \textbf{(b)} says that
$T_{\operatorname*{sign}}\circ T_{\operatorname*{sign}}=\operatorname*{id}$).

But we have $\mathcal{S}^{\mathbf{r}\left(  D\right)  }\cong\mathcal{A}%
\mathbf{E}_{T\mathbf{r}}$ as left $\mathcal{A}$-modules, i.e., as $S_{n}%
$-representations (since $\mathcal{A}=\mathbf{k}\left[  S_{n}\right]  $).
Hence, by Corollary \ref{cor.sign-twist.iso} (applied to $V=\mathcal{S}%
^{\mathbf{r}\left(  D\right)  }$ and $W=\mathcal{A}\mathbf{E}_{T\mathbf{r}}$),
we obtain $\left(  \mathcal{S}^{\mathbf{r}\left(  D\right)  }\right)
^{\operatorname*{sign}}\cong\left(  \mathcal{A}\mathbf{E}_{T\mathbf{r}%
}\right)  ^{\operatorname*{sign}}$. Furthermore, Lemma \ref{lem.sign-twist.a}
(applied to $\mathbf{a}=\mathbf{E}_{T\mathbf{r}}$) yields $\left(
\mathcal{A}\mathbf{E}_{T\mathbf{r}}\right)  ^{\operatorname*{sign}}%
\cong\mathcal{A}\cdot T_{\operatorname*{sign}}\left(  \mathbf{E}_{T\mathbf{r}%
}\right)  $ as $S_{n}$-representations. In view of $T_{\operatorname*{sign}%
}\left(  \mathbf{E}_{T\mathbf{r}}\right)  =\mathbf{F}_{T}$, we can rewrite
this as $\left(  \mathcal{A}\mathbf{E}_{T\mathbf{r}}\right)
^{\operatorname*{sign}}\cong\mathcal{A}\mathbf{F}_{T}$. Thus, $\left(
\mathcal{S}^{\mathbf{r}\left(  D\right)  }\right)  ^{\operatorname*{sign}%
}\cong\left(  \mathcal{A}\mathbf{E}_{T\mathbf{r}}\right)
^{\operatorname*{sign}}\cong\mathcal{A}\mathbf{F}_{T}$. This completes the
proof of Lemma \ref{lem.sign-twist.AFT}.
\end{proof}

We can now prove Theorem \ref{thm.sign-twist.Slam}:

\begin{proof}
[Proof of Theorem \ref{thm.sign-twist.Slam}.]Pick any $n$-tableau $T$ of shape
$Y\left(  \lambda\right)  $. (Such a $T$ exists, as we saw in the proof of
Corollary \ref{cor.spechtmod.nonzero}.) Thus, $T$ is an $n$-tableau of shape
$\lambda$. Note that $\left\vert Y\left(  \lambda\right)  \right\vert
=\left\vert \lambda\right\vert =n$ (since $\lambda$ is a partition of $n$).

Set $\mathbf{F}_{T}:=\nabla_{\operatorname*{Row}T}\nabla_{\operatorname*{Col}%
T}^{-}\in\mathcal{A}$. Thus, Lemma \ref{lem.sign-twist.AFT} (applied to
$D=Y\left(  \lambda\right)  $) yields that%
\[
\mathcal{S}^{Y\left(  \lambda\right)  }\cong\mathcal{A}\mathbf{E}%
_{T}\ \ \ \ \ \ \ \ \ \ \text{and}\ \ \ \ \ \ \ \ \ \ \left(  \mathcal{S}%
^{\mathbf{r}\left(  Y\left(  \lambda\right)  \right)  }\right)
^{\operatorname*{sign}}\cong\mathcal{A}\mathbf{F}_{T}%
\]
as $S_{n}$-representations. But Proposition \ref{prop.specht.FT.basics}
\textbf{(c)} says that $\mathcal{S}^{\lambda}\cong\mathcal{A}\mathbf{E}%
_{T}\cong\mathcal{A}\mathbf{F}_{T}$ as left $\mathcal{A}$-modules, i.e., as
$S_{n}$-representations.

However, $\mathbf{r}\left(  Y\left(  \lambda\right)  \right)  =Y\left(
\lambda^{t}\right)  $ (by the definition of $\lambda^{t}$). Thus,
$\mathcal{S}^{\mathbf{r}\left(  Y\left(  \lambda\right)  \right)
}=\mathcal{S}^{Y\left(  \lambda^{t}\right)  }=\mathcal{S}^{\lambda^{t}}$
(since $\mathcal{S}^{\lambda^{t}}$ is a shorthand for $\mathcal{S}^{Y\left(
\lambda^{t}\right)  }$). Hence, $\mathcal{S}^{\lambda^{t}}=\mathcal{S}%
^{\mathbf{r}\left(  Y\left(  \lambda\right)  \right)  }$, so that $\left(
\mathcal{S}^{\lambda^{t}}\right)  ^{\operatorname*{sign}}=\left(
\mathcal{S}^{\mathbf{r}\left(  Y\left(  \lambda\right)  \right)  }\right)
^{\operatorname*{sign}}\cong\mathcal{A}\mathbf{F}_{T}\cong\mathcal{S}%
^{\lambda}$ as $S_{n}$-representations (since $\mathcal{S}^{\lambda}%
\cong\mathcal{A}\mathbf{F}_{T}$ as $S_{n}$-representations). Hence, by
Corollary \ref{cor.sign-twist.iso} (applied to $V=\left(  \mathcal{S}%
^{\lambda^{t}}\right)  ^{\operatorname*{sign}}$ and $W=\mathcal{S}^{\lambda}%
$), we obtain
\begin{equation}
\left(  \left(  \mathcal{S}^{\lambda^{t}}\right)  ^{\operatorname*{sign}%
}\right)  ^{\operatorname*{sign}}\cong\left(  \mathcal{S}^{\lambda}\right)
^{\operatorname*{sign}} \label{pf.thm.sign-twist.Slam.4}%
\end{equation}
as $S_{n}$-representations.

Finally, Proposition \ref{prop.sign-twist.twice-back} (applied to
$V=\mathcal{S}^{\lambda^{t}}$) yields $\left(  \left(  \mathcal{S}%
^{\lambda^{t}}\right)  ^{\operatorname*{sign}}\right)  ^{\operatorname*{sign}%
}=\mathcal{S}^{\lambda^{t}}$. Hence, we can rewrite
(\ref{pf.thm.sign-twist.Slam.4}) as $\mathcal{S}^{\lambda^{t}}\cong\left(
\mathcal{S}^{\lambda}\right)  ^{\operatorname*{sign}}$. In other words,
$\left(  \mathcal{S}^{\lambda}\right)  ^{\operatorname*{sign}}\cong%
\mathcal{S}^{\lambda^{t}}$. This proves Theorem \ref{thm.sign-twist.Slam}.
\end{proof}

Note that Theorem \ref{thm.sign-twist.Slam} is not generally true without the
\textquotedblleft$h^{\lambda}$ is invertible\textquotedblright\ condition, as
the following exercise shows:

\begin{exercise}
\label{exe.sign-twist.Slam-charp}\fbox{2} As we know from Example
\ref{exa.spechtmod.Yn-1-1}, the Specht module $\mathcal{S}^{\left(
2,1\right)  }$ is isomorphic to the zero-sum subrepresentation $R\left(
\mathbf{k}^{3}\right)  $ of the natural $S_{3}$-representation $\mathbf{k}%
^{3}$.

Assume that $\mathbf{k}$ is a field of characteristic $3$. Prove that $\left(
\mathcal{S}^{\left(  2,1\right)  }\right)  ^{\operatorname*{sign}}$ is not
isomorphic to $\mathcal{S}^{\left(  2,1\right)  ^{t}}$. (Note that $\left(
2,1\right)  ^{t}=\left(  2,1\right)  $.)
\end{exercise}

In many situations, Theorem \ref{thm.sign-twist.Slam} can be generalized from
straight-shaped Young diagrams to arbitrary diagrams:

\begin{theorem}
\label{thm.sign-twist.D}Assume that $\mathbf{k}$ is a field of characteristic
$0$. Let $D$ be any diagram with $\left\vert D\right\vert =n$. Then, the
sign-twist $\left(  \mathcal{S}^{D}\right)  ^{\operatorname*{sign}}$ of the
Specht module $\mathcal{S}^{D}$ is isomorphic to the Specht module
$\mathcal{S}^{\mathbf{r}\left(  D\right)  }$ (where $\mathbf{r}:\mathbb{Z}%
^{2}\rightarrow\mathbb{Z}^{2}$ is the map from Theorem
\ref{thm.partitions.conj}).
\end{theorem}

However, we will need the tools of the next section to prove this theorem.

\begin{question}
Does Theorem \ref{thm.sign-twist.D} still hold if we merely assume
$\mathbf{k}$ to be a commutative ring in which $n!$ is invertible (as opposed
to a field of characteristic $0$)?
\end{question}

\begin{question}
Is there an explicit construction of an isomorphism in Theorem
\ref{thm.sign-twist.D}?
\end{question}

\subsection{The dual of an $S_{n}$-representation}

We now move to a different operation that transforms $S_{n}$-representations
(and, more generally, $G$-representations for any group $G$): dualization.

\subsubsection{Reminders on dual $\mathbf{k}$-modules}

This operation is based on the linear-algebraic notion of a dual $\mathbf{k}%
$-module, which, in turn, generalizes the classical concept of a dual
$\mathbf{k}$-vector space. We begin by recalling how these notions are
defined. While duals are generally better-behaved for $\mathbf{k}$-vector
spaces (i.e., when $\mathbf{k}$ is a field) than for $\mathbf{k}$-modules
(i.e., when $\mathbf{k}$ is just a commutative ring), their definition is the
same in both settings; thus we shall waste no time talking about the case of
$\mathbf{k}$-vector spaces except in those situations where we can genuinely
say more about it than about the general case.

We recall the definition of duals (\cite[3.110]{Axler24}):

\begin{definition}
\label{def.dual.dual-k-mod}\textbf{(a)} If $V$ and $W$ are two $\mathbf{k}%
$-modules, then $\operatorname*{Hom}\left(  V,W\right)  $ means the
$\mathbf{k}$-module of all $\mathbf{k}$-linear maps from $V$ to $W$. (Its
addition and scaling are entrywise: e.g., we have $\left(  f+g\right)  \left(
v\right)  =f\left(  v\right)  +g\left(  v\right)  $ for all $f,g\in
\operatorname*{Hom}\left(  V,W\right)  $ and $v\in V$.) \medskip

\textbf{(b)} If $V$ is any $\mathbf{k}$-module, then the $\mathbf{k}$-module
$\operatorname*{Hom}\left(  V,\mathbf{k}\right)  $ is also denoted by
$V^{\ast}$ and is called the \emph{dual} (or \emph{dual }$\mathbf{k}%
$\emph{-module}, or \emph{dual space}) of $V$.
\end{definition}

Thus, the dual $V^{\ast}$ of a $\mathbf{k}$-module $V$ consists of all
$\mathbf{k}$-linear maps from $V$ to $\mathbf{k}$; these are often called
\textquotedblleft linear functions\textquotedblright\ or \textquotedblleft
linear functionals\textquotedblright\ on $V$. Some authors also prefer the
notation $V^{\vee}$ for $V^{\ast}$.

\begin{warning}
The dual $V^{\ast}$ of a $\mathbf{k}$-module $V$ depends not just on $V$ but
also on the base ring $\mathbf{k}$. Thus, for example, a $\mathbb{Q}$-vector
space $V$ is always a $\mathbb{Z}$-module (by restriction of scalars), but the
dual $V^{\ast}$ of the $\mathbb{Q}$-vector space $V$ is not the same as the
dual $V^{\ast}$ of the $\mathbb{Z}$-module $V$, not even as a set. (Indeed, it
can be shown that the latter dual is always trivial, since the only
$\mathbb{Z}$-linear map $f:V\rightarrow\mathbb{Z}$ from a $\mathbb{Q}$-vector
space $V$ to $\mathbb{Z}$ is the constant-zero map.)

We can avoid this pitfall by using the unambiguous notation $V_{\mathbf{k}%
}^{\ast}$ (instead of $V^{\ast}$) for the dual of a $\mathbf{k}$-module $V$.
However, we shall keep such precautions restricted to situations in which they
are genuinely necessary (i.e., when we are considering one and the same module
as a module over different base rings).
\end{warning}

The following properties of duals are known from linear algebra (see
\cite[\textquotedblleft Dual modules\textquotedblright, \S 5]{Conrad}%
\footnote{In the case when $\mathbf{k}$ is a field, you can find these facts
and more in any good text on linear algebra, such as \cite[\S 6]{StoLui19}.}):

\begin{proposition}
\label{prop.dual.basis}Let $V$ be a free $\mathbf{k}$-module with a finite
basis $\left(  e_{i}\right)  _{i\in I}$. Then, its dual $V^{\ast}$ has a
finite basis $\left(  e_{i}^{\ast}\right)  _{i\in I}$, which is defined as
follows: For each $i\in I$, we let $e_{i}^{\ast}:V\rightarrow\mathbf{k}$ be
the $\mathbf{k}$-linear map that sends each basis vector $e_{j}$ of $V$ to the
Kronecker delta $\delta_{i,j}\in\mathbf{k}$.

This basis $\left(  e_{i}^{\ast}\right)  _{i\in I}$ is called the \emph{dual
basis} of the original basis $\left(  e_{i}\right)  _{i\in I}$. Its definition
works even if the basis $\left(  e_{i}\right)  _{i\in I}$ is not finite, but
in that case it is not a basis of $V^{\ast}$ but merely a $\mathbf{k}%
$-linearly independent family.
\end{proposition}

\begin{proposition}
\label{prop.dual.functorial}Let $V$ and $W$ be two $\mathbf{k}$-modules.
\medskip

\textbf{(a)} For each $\mathbf{k}$-linear map $f:V\rightarrow W$, there is a
$\mathbf{k}$-linear map $f^{\ast}:W^{\ast}\rightarrow V^{\ast}$ (called the
\emph{dual} or the \emph{adjoint} of $f$) defined by%
\[
f^{\ast}\left(  g\right)  :=g\circ f\ \ \ \ \ \ \ \ \ \ \text{for each }g\in
W^{\ast}.
\]

\textbf{(b)} If $U$ is a further $\mathbf{k}$-module, and if $f:V\rightarrow
W$ and $g:W\rightarrow U$ are two $\mathbf{k}$-linear maps, then $\left(
g\circ f\right)  ^{\ast}=f^{\ast}\circ g^{\ast}$. \medskip

\textbf{(c)} If $f:V\rightarrow W$ is a $\mathbf{k}$-module isomorphism, then
its dual $f^{\ast}:W^{\ast}\rightarrow V^{\ast}$ is also a $\mathbf{k}$-module
isomorphism. \medskip

\textbf{(d)} Assume that the $\mathbf{k}$-modules $V$ and $W$ are free, with
finite bases $\left(  v_{i}\right)  _{i\in I}$ and $\left(  w_{j}\right)
_{j\in J}$. Let $f:V\rightarrow W$ be a $\mathbf{k}$-linear map, and let
$A\in\mathbf{k}^{J\times I}$ be the $J\times I$-matrix that represents $f$
with respect to the bases $\left(  v_{i}\right)  _{i\in I}$ and $\left(
w_{j}\right)  _{j\in J}$. Let $\left(  v_{i}^{\ast}\right)  _{i\in I}$ and
$\left(  w_{j}^{\ast}\right)  _{j\in J}$ be the dual bases of these bases
$\left(  v_{i}\right)  _{i\in I}$ and $\left(  w_{j}\right)  _{j\in J}$. Then,
the matrix that represents the map $f^{\ast}:W^{\ast}\rightarrow V^{\ast}$
with respect to these dual bases $\left(  w_{j}^{\ast}\right)  _{j\in J}$ and
$\left(  v_{i}^{\ast}\right)  _{i\in I}$ is the transpose $A^{T}$ of the
matrix $A$.
\end{proposition}

Proposition \ref{prop.dual.functorial} \textbf{(d)} is the conceptual meaning
of transpose matrices (if you view matrices as models for linear maps).
Another cluster of well-known facts concerns double duals (i.e., duals of duals):

\begin{proposition}
\label{prop.dual.linalg1}Let $V$ be a $\mathbf{k}$-module. \medskip

\textbf{(a)} There is a canonical $\mathbf{k}$-linear map $\iota:V\rightarrow
V^{\ast\ast}$ (where $V^{\ast\ast}=\left(  V^{\ast}\right)  ^{\ast}$) that
sends each vector $v\in V$ to the $\mathbf{k}$-linear map%
\begin{align*}
\operatorname*{ev}\nolimits_{v}:V^{\ast}  &  \rightarrow\mathbf{k},\\
g  &  \mapsto g\left(  v\right)  .
\end{align*}

\textbf{(b)} If the $\mathbf{k}$-module $V$ has a finite basis, then this map
$\iota$ is an isomorphism. \medskip

\textbf{(c)} In particular, if $\mathbf{k}$ is a field and $V$ is a
finite-dimensional $\mathbf{k}$-vector space, then $\dim\left(  V^{\ast
}\right)  =\dim V$.
\end{proposition}

A proof of Proposition \ref{prop.dual.linalg1} can be found in
\cite[\textquotedblleft Dual modules\textquotedblright, \S 4]{Conrad}%
.\footnote{Those who believe in the axiom of choice (and classical logic) can
say a bit more: The map $\iota$ in Proposition \ref{prop.dual.linalg1}
\textbf{(a)} is injective whenever $\mathbf{k}$ is a field. See, e.g.,
\url{https://math.stackexchange.com/questions/4894044} for the proof of this
fact (which we will not need in these notes). Keep in mind that it is false
when $\mathbf{k}$ is not a field (not even the axiom of choice can make the
map $\iota:\mathbb{Q}\rightarrow\mathbb{Q}^{\ast\ast}$ injective for
$\mathbf{k}=\mathbb{Z}$, since $\mathbb{Q}^{\ast}=0$ entails $\mathbb{Q}%
^{\ast\ast}=0^{\ast}=0$ here).} Note that Proposition \ref{prop.dual.linalg1}
\textbf{(b)} requires a \textbf{finite basis}; neither a finite generating set
nor an infinite basis would be sufficient.

See \cite[\textquotedblleft Dual modules\textquotedblright]{Conrad} for more
about duals of $\mathbf{k}$-modules.

\subsubsection{Bilinear forms}

Bilinear forms are an important tool for understanding duals. To remind:

\begin{itemize}
\item A $\mathbf{k}$-bilinear map is a map $f:U\times V\rightarrow W$ (where
$U$, $V$ and $W$ are three $\mathbf{k}$-modules) such that if we fix either
argument, the map $f$ becomes a linear map in the other (see below for a
detailed definition).

\item A $\mathbf{k}$-bilinear form is just a $\mathbf{k}$-bilinear map to the
base ring $\mathbf{k}$ itself.
\end{itemize}

Let us state these definitions in more detail, as the relevant notations will
be useful later on. We begin with a notation for turning a two-argument map
into a one-argument map by fixing one of the arguments to a constant:

\begin{definition}
\label{def.dual.fu?}Let $U$, $V$ and $W$ be three sets, and $f:U\times
V\rightarrow W$ be any map. Then: \medskip

\textbf{(a)} For any $u\in U$, we let $f\left(  u,?\right)  $ be the map from
$V$ to $W$ that sends each $v\in V$ to $f\left(  u,v\right)  $. \medskip

\textbf{(b)} For any $v\in V$, we let $f\left(  ?,v\right)  $ be the map from
$U$ to $W$ that sends each $u\in U$ to $f\left(  u,v\right)  $.
\end{definition}

For instance, if $f:\mathbb{Z}\times\mathbb{Z}\rightarrow\mathbb{Z}$ is the
addition map (i.e., the map $\left(  a,b\right)  \mapsto a+b$), then $f\left(
3,?\right)  $ is the map that sends each integer $b$ to $3+b$. (This is the
same as the map $f\left(  ?,3\right)  $, since addition is commutative.)

Next, we define bilinearity in detail:

\begin{definition}
\label{def.dual.bilinear}Let $U$, $V$ and $W$ be three $\mathbf{k}$-modules,
and $f:U\times V\rightarrow W$ be any map. \medskip

\textbf{(a)} We say that the map $f$ is \emph{linear in the second argument}
if for each $u\in U$, the map $f\left(  u,?\right)  :V\rightarrow W$ is
$\mathbf{k}$-linear. In this case, we define the map%
\begin{align*}
f_{L}:U  &  \rightarrow\operatorname*{Hom}\left(  V,W\right)  ,\\
u  &  \mapsto f\left(  u,?\right)  .
\end{align*}
We call $f_{L}$ the \emph{left-curried form of }$f$. Explicitly, $f_{L}$ is
thus given by the equality%
\[
\left(  f_{L}\left(  u\right)  \right)  \left(  v\right)  =f\left(
u,v\right)  \ \ \ \ \ \ \ \ \ \ \text{for all }u\in U\text{ and }v\in V.
\]

\textbf{(b)} We say that the map $f$ is \emph{linear in the first argument} if
for each $v\in V$, the map $f\left(  ?,v\right)  :U\rightarrow W$ is
$\mathbf{k}$-linear. In this case, we define the map%
\begin{align*}
f_{R}:V  &  \rightarrow\operatorname*{Hom}\left(  U,W\right)  ,\\
v  &  \mapsto f\left(  ?,v\right)  .
\end{align*}
We call $f_{R}$ the \emph{right-curried form of }$f$. Explicitly, $f_{R}$ is
thus given by the equality%
\[
\left(  f_{R}\left(  v\right)  \right)  \left(  u\right)  =f\left(
u,v\right)  \ \ \ \ \ \ \ \ \ \ \text{for all }u\in U\text{ and }v\in V.
\]

\textbf{(c)} We say that the map $f$ is \emph{bilinear} (or $\mathbf{k}%
$\emph{-bilinear}, to be more precise) if $f$ is both linear in the first
argument and linear in the second argument. In this case, both maps
$f_{L}:U\rightarrow\operatorname*{Hom}\left(  V,W\right)  $ and $f_{R}%
:V\rightarrow\operatorname*{Hom}\left(  U,W\right)  $ are easily seen to be
$\mathbf{k}$-linear.
\end{definition}

We can also rewrite this definition of \textquotedblleft
bilinear\textquotedblright\ explicitly as a list of axioms:

\begin{remark}
\label{rmk.dual.bilin.axioms}Let $U$, $V$ and $W$ be three $\mathbf{k}%
$-modules. Let $f:U\times V\rightarrow W$ be a map. Then, $f$ is $\mathbf{k}%
$-bilinear if and only if it satisfies the six axioms
\begin{align*}
f\left(  u,\ v_{1}+v_{2}\right)   &  =f\left(  u,v_{1}\right)  +f\left(
u,v_{2}\right)  \ \ \ \ \ \ \ \ \ \ \text{for all }u\in U\text{ and }%
v_{1},v_{2}\in V;\\
f\left(  u,\ \lambda v\right)   &  =\lambda f\left(  u,v\right)
\ \ \ \ \ \ \ \ \ \ \text{for all }u\in U\text{ and }v\in V\text{ and }%
\lambda\in\mathbf{k};\\
f\left(  u,0\right)   &  =0\ \ \ \ \ \ \ \ \ \ \text{for all }u\in U;\\
f\left(  u_{1}+u_{2},\ v\right)   &  =f\left(  u_{1},v\right)  +f\left(
u_{2},v\right)  \ \ \ \ \ \ \ \ \ \ \text{for all }u_{1},u_{2}\in U\text{ and
}v\in V;\\
f\left(  \lambda u,\ v\right)   &  =\lambda f\left(  u,v\right)
\ \ \ \ \ \ \ \ \ \ \text{for all }u\in U\text{ and }v\in V\text{ and }%
\lambda\in\mathbf{k};\\
f\left(  0,v\right)   &  =0\ \ \ \ \ \ \ \ \ \ \text{for all }v\in V.
\end{align*}

\end{remark}

\begin{definition}
\label{def.dual.bilform}Let $U$ and $V$ be two $\mathbf{k}$-modules. A
\emph{bilinear form} on $U$ and $V$ means a bilinear map $f:U\times
V\rightarrow\mathbf{k}$.
\end{definition}

Thus, if $f:U\times V\rightarrow\mathbf{k}$ is a bilinear form, then its
left-curried form $f_{L}$ is a $\mathbf{k}$-linear map from $U$ to $V^{\ast}$,
whereas its right-curried form $f_{R}$ is a $\mathbf{k}$-linear map from $V$
to $U^{\ast}$.

\begin{noncompile}
Let $U$ and $V$ be two $\mathbf{k}$-modules. Let $f:U\times V\rightarrow
\mathbf{k}$ be a bilinear form. Then, its left-curried form $f_{L}$ is a
$\mathbf{k}$-linear map from $U$ to $V^{\ast}$, whereas its right-curried form
$f_{R}$ is a $\mathbf{k}$-linear map from $V$ to $U^{\ast}$. \medskip

\textbf{(a)} We say that the form $f$ is \emph{nondegenerate} if the map
$f_{L}$ is injective. \medskip

\textbf{(b)} We say that the form $f$ is \emph{perfect} if the map $f_{L}$ is bijective.

Note that some authors (e.g., Conrad in \cite[\textquotedblleft Dual
modules\textquotedblright, \S 4]{Conrad}) have a slightly stronger notion of
perfectness, which requires not just $f_{L}$ but also $f_{R}$ to be bijective.

\label{prop.bilform.perfect}Let $U$ and $V$ be two $\mathbf{k}$-modules. Let
$f:U\times V\rightarrow\mathbf{k}$ be a bilinear form. Then: \medskip

\textbf{(a)} If $f$ is perfect, then $f$ is nondegenerate. \medskip

\textbf{(b)} Assume that $\mathbf{k}$ is a field, and that $U$ and $V$ are
finite-dimensional vector spaces...
\end{noncompile}

Bilinear forms on free $\mathbf{k}$-modules can be described by the values
they take on the basis vectors:

\begin{proposition}
\label{prop.dual.bilf.on-basis}Let $U$ and $V$ be two free $\mathbf{k}%
$-modules with bases $\left(  u_{i}\right)  _{i\in I}$ and $\left(
v_{j}\right)  _{j\in J}$. Then, any bilinear form $f:U\times V\rightarrow
\mathbf{k}$ is uniquely determined by its values $f\left(  u_{i},v_{j}\right)
$ for all $i\in I$ and $j\in J$. Moreover, for any family $\left(
a_{i,j}\right)  _{\left(  i,j\right)  \in I\times J}\in\mathbf{k}^{I\times J}$
of scalars, there is a unique bilinear form $f:U\times V\rightarrow\mathbf{k}$
that satisfies%
\[
\left(  f\left(  u_{i},v_{j}\right)  =a_{i,j}\ \ \ \ \ \ \ \ \ \ \text{for all
}i\in I\text{ and }j\in J\right)  .
\]
Thus, a bilinear form $f:U\times V\rightarrow\mathbf{k}$ can be defined by
providing its values $f\left(  u_{i},v_{j}\right)  $ for all $i\in I$ and
$j\in J$.
\end{proposition}

\begin{proof}
See \cite[Theorem 3.9.2]{23wa} (which proves the same claim in the more
general setting where bilinear forms $f:U\times V\rightarrow\mathbf{k}$ are
replaced by bilinear maps $f:U\times V\rightarrow W$).
\end{proof}

\begin{proposition}
\label{prop.dual.bilf.rep-mat}Let $U$ and $V$ be two free $\mathbf{k}$-modules
with finite bases $\left(  u_{i}\right)  _{i\in I}$ and $\left(  v_{j}\right)
_{j\in J}$. Let $f:U\times V\rightarrow\mathbf{k}$ be a bilinear form. Then:
\medskip

\textbf{(a)} The $\mathbf{k}$-linear map $f_{L}:U\rightarrow V^{\ast}$ is
invertible if and only if the matrix $\left(  f\left(  u_{i},v_{j}\right)
\right)  _{\left(  i,j\right)  \in I\times J}\in\mathbf{k}^{I\times J}$ is
invertible. \medskip

\textbf{(b)} The $\mathbf{k}$-linear map $f_{R}:V\rightarrow U^{\ast}$ is
invertible if and only if the matrix $\left(  f\left(  u_{i},v_{j}\right)
\right)  _{\left(  i,j\right)  \in I\times J}\in\mathbf{k}^{I\times J}$ is invertible.
\end{proposition}

\begin{proof}
\textbf{(b)} Let $\left(  u_{i}^{\ast}\right)  _{i\in I}$ and $\left(
v_{j}^{\ast}\right)  _{j\in J}$ be the dual bases of the bases $\left(
u_{i}\right)  _{i\in I}$ and $\left(  v_{j}\right)  _{j\in J}$ of $U$ and $V$.
For each $j\in J$, we have%
\begin{equation}
f_{R}\left(  v_{j}\right)  =\sum_{i\in I}f\left(  u_{i},v_{j}\right)
u_{i}^{\ast}. \label{pf.prop.dual.bilf.rep-mat.a.1}%
\end{equation}

\begin{proof}
[Proof of (\ref{pf.prop.dual.bilf.rep-mat.a.1}).]Let $j\in J$. Then, both
$f_{R}\left(  v_{j}\right)  $ and $\sum_{i\in I}f\left(  u_{i},v_{j}\right)
u_{i}^{\ast}$ are $\mathbf{k}$-linear maps from $U$ to $\mathbf{k}$. Hence, in
order to prove that these two maps are equal, it suffices to show that they
agree on each of the basis vectors $u_{k}$ of $U$. But this is easy: For each
$k\in I$, we have%
\begin{align}
\left(  \sum_{i\in I}f\left(  u_{i},v_{j}\right)  u_{i}^{\ast}\right)  \left(
u_{k}\right)   &  =\sum_{i\in I}f\left(  u_{i},v_{j}\right)  \underbrace{u_{i}%
^{\ast}\left(  u_{k}\right)  }_{\substack{=\delta_{i,k}\\\text{(by the
definition of the}\\\text{dual basis }\left(  u_{i}^{\ast}\right)  _{i\in
I}\text{)}}}=\sum_{i\in I}f\left(  u_{i},v_{j}\right)  \delta_{i,k}\nonumber\\
&  =f\left(  u_{k},v_{j}\right)  \underbrace{\delta_{k,k}}_{=1}+\sum
_{\substack{i\in I;\\i\neq k}}f\left(  u_{i},v_{j}\right)  \underbrace{\delta
_{i,k}}_{\substack{=0\\\text{(since }i\neq k\text{)}}}\nonumber\\
&  =f\left(  u_{k},v_{j}\right)  +\underbrace{\sum_{\substack{i\in I;\\i\neq
k}}f\left(  u_{i},v_{j}\right)  0}_{=0}\nonumber\\
&  =f\left(  u_{k},v_{j}\right)  . \label{pf.prop.dual.bilf.rep-mat.a.2}%
\end{align}
On the other hand, for each $k\in I$, we have%
\begin{align*}
\left(  f_{R}\left(  v_{j}\right)  \right)  \left(  u_{k}\right)   &
=f\left(  u_{k},v_{j}\right)  \ \ \ \ \ \ \ \ \ \ \left(  \text{by the
definition of }f_{R}\right) \\
&  =\left(  \sum_{i\in I}f\left(  u_{i},v_{j}\right)  u_{i}^{\ast}\right)
\left(  u_{k}\right)  \ \ \ \ \ \ \ \ \ \ \left(  \text{by
(\ref{pf.prop.dual.bilf.rep-mat.a.2})}\right)  .
\end{align*}
In other words, the two $\mathbf{k}$-linear maps $f_{R}\left(  v_{j}\right)  $
and $\sum_{i\in I}f\left(  u_{i},v_{j}\right)  u_{i}^{\ast}$ agree on each of
the basis vectors $u_{k}$ of $U$. Because of their $\mathbf{k}$-linearity,
this entails that these two maps are equal. Thus,
(\ref{pf.prop.dual.bilf.rep-mat.a.1}) is proved.
\end{proof}

Now, the equality (\ref{pf.prop.dual.bilf.rep-mat.a.1}) shows that the matrix
$\left(  f\left(  u_{i},v_{j}\right)  \right)  _{\left(  i,j\right)  \in
I\times J}$ represents the $\mathbf{k}$-linear map $f_{R}$ with respect to the
bases $\left(  v_{j}\right)  _{j\in J}$ and $\left(  u_{i}^{\ast}\right)
_{i\in I}$ of $V$ and $U^{\ast}$. Hence, the latter map is invertible if and
only if the former matrix is invertible (since invertible matrices represent
invertible linear maps). This proves Proposition \ref{prop.dual.bilf.rep-mat}
\textbf{(b)}. \medskip

\textbf{(a)} In our proof of part \textbf{(b)}, we have shown that the
$\mathbf{k}$-linear map $f_{R}:V\rightarrow U^{\ast}$ is invertible if and
only if the matrix $\left(  f\left(  u_{i},v_{j}\right)  \right)  _{\left(
i,j\right)  \in I\times J}\in\mathbf{k}^{I\times J}$ is invertible. Similarly,
we can show that the $\mathbf{k}$-linear map $f_{L}:U\rightarrow V^{\ast}$ is
invertible if and only if the matrix $\left(  f\left(  u_{i},v_{j}\right)
\right)  _{\left(  j,i\right)  \in J\times I}\in\mathbf{k}^{J\times I}$ (note
that the two indices now appear in reverse order!) is invertible.\footnote{The
proof relies on the identity%
\[
f_{L}\left(  u_{i}\right)  =\sum_{j\in J}f\left(  u_{i},v_{j}\right)
v_{j}^{\ast}\ \ \ \ \ \ \ \ \ \ \text{for all }i\in I;
\]
this is an analogue of (\ref{pf.prop.dual.bilf.rep-mat.a.1}).} Hence, we have
the following chain of logical equivalences:%
\begin{align*}
&  \ \left(  \text{the }\mathbf{k}\text{-linear map }f_{L}:U\rightarrow
V^{\ast}\text{ is invertible}\right) \\
&  \Longleftrightarrow\ \left(  \text{the matrix }\left(  f\left(  u_{i}%
,v_{j}\right)  \right)  _{\left(  j,i\right)  \in J\times I}\in\mathbf{k}%
^{J\times I}\text{ is invertible}\right) \\
&  \Longleftrightarrow\ \left(  \text{the transpose of the matrix }\left(
f\left(  u_{i},v_{j}\right)  \right)  _{\left(  i,j\right)  \in I\times J}%
\in\mathbf{k}^{I\times J}\text{ is invertible}\right) \\
&  \ \ \ \ \ \ \ \ \ \ \ \ \ \ \ \ \ \ \ \ \left(
\begin{array}
[c]{c}%
\text{since the matrix }\left(  f\left(  u_{i},v_{j}\right)  \right)
_{\left(  j,i\right)  \in J\times I}\in\mathbf{k}^{J\times I}\text{ is the}\\
\text{transpose of the matrix }\left(  f\left(  u_{i},v_{j}\right)  \right)
_{\left(  i,j\right)  \in I\times J}\in\mathbf{k}^{I\times J}%
\end{array}
\right) \\
&  \Longleftrightarrow\ \left(  \text{the matrix }\left(  f\left(  u_{i}%
,v_{j}\right)  \right)  _{\left(  i,j\right)  \in I\times J}\in\mathbf{k}%
^{I\times J}\text{ is invertible}\right)
\end{align*}
(since the transpose of a matrix $A$ is invertible if and only if the matrix
$A$ is invertible). This proves Proposition \ref{prop.dual.bilf.rep-mat}
\textbf{(a)}.
\end{proof}

\begin{definition}
\label{def.dual.gram}The matrix $\left(  f\left(  u_{i},v_{j}\right)  \right)
_{\left(  i,j\right)  \in I\times J}\in\mathbf{k}^{I\times J}$ in Proposition
\ref{prop.dual.bilf.rep-mat} is called the \emph{Gram matrix} of the bilinear
form $f$ with respect to the bases $\left(  u_{i}\right)  _{i\in I}$ and
$\left(  v_{j}\right)  _{j\in J}$.
\end{definition}

Another useful fact from linear algebra is the following:

\begin{lemma}
\label{lem.dual.field-fLfRinj}Assume that $\mathbf{k}$ is a field. Let $U$ and
$V$ be two finite-dimensional $\mathbf{k}$-vector spaces. Let $\alpha
:U\rightarrow V^{\ast}$ and $\beta:V\rightarrow U^{\ast}$ be two injective
$\mathbf{k}$-linear maps. Then, these two maps $\alpha$ and $\beta$ are bijective.
\end{lemma}

\begin{proof}
Since $V$ is finite-dimensional, we have $\dim\left(  V^{\ast}\right)  =\dim
V$ (by Proposition \ref{prop.dual.linalg1} \textbf{(c)}). Likewise,
$\dim\left(  U^{\ast}\right)  =\dim U$.

Now, recall the following two classical facts from linear algebra (which are
easy consequences of the rank-nullity theorem):

\begin{statement}
\textit{Fact 1:} If $\varphi:X\rightarrow Y$ is an injective $\mathbf{k}%
$-linear map between two finite-dimensional $\mathbf{k}$-vector spaces $X$ and
$Y$, then $\dim X\leq\dim Y$.
\end{statement}

\begin{statement}
\textit{Fact 2:} If $\varphi:X\rightarrow Y$ is an injective $\mathbf{k}%
$-linear map between two finite-dimensional $\mathbf{k}$-vector spaces $X$ and
$Y$ that satisfy $\dim X\geq\dim Y$, then $\varphi$ is bijective.
\end{statement}

We have assumed that the $\mathbf{k}$-linear map $\alpha:U\rightarrow V^{\ast
}$ is injective. Hence, Fact 1 (applied to $X=U$ and $Y=V^{\ast}$ and
$\varphi=\alpha$) yields $\dim U\leq\dim\left(  V^{\ast}\right)  =\dim V$.
Hence, $\dim V\geq\dim U=\dim\left(  U^{\ast}\right)  $ (since $\dim\left(
U^{\ast}\right)  =\dim U$).

But we have also assumed that the $\mathbf{k}$-linear map $\beta:V\rightarrow
U^{\ast}$ is injective. Hence, Fact 2 (applied to $X=V$ and $Y=U^{\ast}$ and
$\varphi=\beta$) yields that $\beta$ is bijective (since $\dim V\geq
\dim\left(  U^{\ast}\right)  $). Similarly, we can show that $\alpha$ is
bijective. This proves Lemma \ref{lem.dual.field-fLfRinj}.
\end{proof}

\subsubsection{The dual of a group representation}

We will now define the \emph{dual} of a representation $V$ of a group $G$.
Just as for sign-twists (in Section \ref{sec.rep.sign-twist}), we will give
two definitions: one for representations of $G$ and one for left
$\mathbf{k}\left[  G\right]  $-modules. As we know, a representation of $G$
(over $\mathbf{k}$) is \textquotedblleft the same as\textquotedblright\ a left
$\mathbf{k}\left[  G\right]  $-module (by Corollary \ref{cor.rep.G-rep.mod}),
so that both definitions apply to the same objects; as we will soon see
(Proposition \ref{prop.dual.exist} \textbf{(c)}), they also produce the same objects.

We begin with the first definition:

\begin{definition}
\label{def.dual.dual.G}Let $G$ be a group. Let $V$ be a representation of $G$
over $\mathbf{k}$. For each $g\in G$, we consider the map%
\begin{align*}
\tau_{g}:V  &  \rightarrow V,\\
x  &  \mapsto g\rightharpoonup x.
\end{align*}
This map $\tau_{g}$ is $\mathbf{k}$-linear (by Definition
\ref{def.rep.G-rep.G-rep}).

Now, the \emph{dual} of $V$ is a new representation of $G$ over $\mathbf{k}$,
which is called $V^{\ast}$ and is defined as follows:

\begin{itemize}
\item As a $\mathbf{k}$-module, $V^{\ast}$ is just the dual $V^{\ast}$ of the
$\mathbf{k}$-module $V$ (as defined in Definition \ref{def.dual.dual-k-mod}
\textbf{(b)}).

\item The $G$-action on $V^{\ast}$ is the map $\overset{\ast}{\rightharpoonup
}\ :G\times V^{\ast}\rightarrow V^{\ast}$ defined by the formula%
\[
g\overset{\ast}{\rightharpoonup}f:=f\circ\tau_{g^{-1}}%
\ \ \ \ \ \ \ \ \ \ \text{for all }g\in G\text{ and }f\in V^{\ast}.
\]

\end{itemize}

This representation $V^{\ast}$ really exists (by Proposition
\ref{prop.dual.exist} \textbf{(a)} below), so the dual of $V$ is well-defined.
\end{definition}

\begin{remark}
\label{rmk.dual.dual.G.on-v}The definition of the $G$-action on $V^{\ast}$ in
Definition \ref{def.dual.dual.G} can be restated as follows: The $G$-action on
$V^{\ast}$ is the map $\overset{\ast}{\rightharpoonup}\ :G\times V^{\ast
}\rightarrow V^{\ast}$ that satisfies
\begin{align}
&  \left(  g\overset{\ast}{\rightharpoonup}f\right)  \left(  v\right)
=f\left(  g^{-1}\rightharpoonup v\right) \label{eq.def.dual.dual.G.on-v}\\
&  \ \ \ \ \ \ \ \ \ \ \text{for all }g\in G\text{ and }f\in V^{\ast}\text{
and }v\in V\nonumber
\end{align}
(where the \textquotedblleft$\rightharpoonup$\textquotedblright\ on the right
hand side means the action of $G$ on $V$). Indeed, the formula%
\[
g\overset{\ast}{\rightharpoonup}f:=f\circ\tau_{g^{-1}}%
\ \ \ \ \ \ \ \ \ \ \text{for all }g\in G\text{ and }f\in V^{\ast}%
\]
(which defines $\overset{\ast}{\rightharpoonup}$ in Definition
\ref{def.dual.dual.G}) is clearly equivalent to the formula%
\begin{align*}
&  \left(  g\overset{\ast}{\rightharpoonup}f\right)  \left(  v\right)
=\left(  f\circ\tau_{g^{-1}}\right)  \left(  v\right) \\
&  \ \ \ \ \ \ \ \ \ \ \text{for all }g\in G\text{ and }f\in V^{\ast}\text{
and }v\in V
\end{align*}
(since two maps from $V$ to $\mathbf{k}$ are equal if and only if they give
equal values on every $v\in V$). But the latter formula is, in turn,
equivalent to (\ref{eq.def.dual.dual.G.on-v}), since all $g\in G$ and $f\in
V^{\ast}$ and $v\in V$ satisfy%
\[
\left(  f\circ\tau_{g^{-1}}\right)  \left(  v\right)  =f\left(
\underbrace{\tau_{g^{-1}}\left(  v\right)  }_{\substack{=g^{-1}\rightharpoonup
v\\\text{(by the definition of }\tau_{g^{-1}}\text{)}}}\right)  =f\left(
g^{-1}\rightharpoonup v\right)  .
\]

\end{remark}

For the second definition, we need a bit of preparation.:

\begin{definition}
\label{def.dual.S}Let $G$ be a group. Consider the $\mathbf{k}$-linear map%
\begin{align*}
S:\mathbf{k}\left[  G\right]   &  \rightarrow\mathbf{k}\left[  G\right]  ,\\
w  &  \mapsto w^{-1}\ \ \ \ \ \ \ \ \ \ \text{for all }w\in G.
\end{align*}
This map $S$ is called the \emph{antipode} of $\mathbf{k}\left[  G\right]  $;
it was already defined in Corollary \ref{cor.char.ladisch-X} (and, for
$G=S_{n}$, it was defined even earlier, in Definition \ref{def.S.S}). Its main
properties are the following: \medskip

\textbf{(a)} The map $S:\mathbf{k}\left[  G\right]  \rightarrow\mathbf{k}%
\left[  G\right]  $ is a $\mathbf{k}$-algebra anti-automorphism. \medskip

\textbf{(b)} It is furthermore an involution (i.e., it satisfies $S\circ
S=\operatorname*{id}$).\medskip

Indeed, both of these properties were proved for $G=S_{n}$ in Theorem
\ref{thm.S.auto}; the proofs in the general case are analogous.
\end{definition}

We are now ready to define duals for left $\mathbf{k}\left[  G\right]  $-modules:

\begin{definition}
\label{def.dual.dual.kG}Let $G$ be a group. Let $V$ be a left $\mathbf{k}%
\left[  G\right]  $-module. Then, the \emph{dual} of $V$ is a new left
$\mathbf{k}\left[  G\right]  $-module, which is called $V^{\ast}$ and is
defined as follows:

\begin{itemize}
\item As an additive group, $V^{\ast}$ is just the dual $V^{\ast}$ of the
$\mathbf{k}$-module $V$ (as defined in Definition \ref{def.dual.dual-k-mod}
\textbf{(b)}).

\item The action of $\mathbf{k}\left[  G\right]  $ on $V^{\ast}$ is given by%
\[
\mathbf{a}\cdot f:=f\circ\tau_{S\left(  \mathbf{a}\right)  }%
\ \ \ \ \ \ \ \ \ \ \text{for all }\mathbf{a}\in\mathbf{k}\left[  G\right]
\text{ and }f\in V^{\ast},
\]
where $\tau_{S\left(  \mathbf{a}\right)  }$ denotes the $\mathbf{k}$-linear
map%
\begin{align*}
V  &  \rightarrow V,\\
v  &  \mapsto S\left(  \mathbf{a}\right)  \cdot v.
\end{align*}
In other words, the action of $\mathbf{k}\left[  G\right]  $ on $V^{\ast}$ is
the map%
\begin{align*}
\mathbf{k}\left[  G\right]  \times V^{\ast}  &  \rightarrow V^{\ast},\\
\left(  \mathbf{a},f\right)   &  \mapsto f\circ\tau_{S\left(  \mathbf{a}%
\right)  }.
\end{align*}

\end{itemize}

This left $\mathbf{k}\left[  G\right]  $-module $V^{\ast}$ is well-defined (by
Proposition \ref{prop.dual.exist} \textbf{(b)} below).
\end{definition}

\begin{remark}
\label{rmk.dual.dual.kG.on-v}The definition of the action of $\mathbf{k}%
\left[  G\right]  $ on $V^{\ast}$ in Definition \ref{def.dual.dual.kG} can be
restated as follows: The action of $\mathbf{k}\left[  G\right]  $ on $V^{\ast
}$ is given by the rule%
\begin{align}
&  \left(  \mathbf{a}\cdot f\right)  \left(  v\right)  =f\left(  S\left(
\mathbf{a}\right)  \cdot v\right) \label{eq.def.dual.dual.kG.on-v}\\
&  \ \ \ \ \ \ \ \ \ \ \text{for all }\mathbf{a}\in\mathbf{k}\left[  G\right]
\text{ and }f\in V^{\ast}\text{ and }v\in V\nonumber
\end{align}
(where the \textquotedblleft$\cdot$\textquotedblright\ on the right hand side
means the action of $\mathbf{k}\left[  G\right]  $ on $V$). Indeed, the
formula%
\[
\mathbf{a}\cdot f:=f\circ\tau_{S\left(  \mathbf{a}\right)  }%
\ \ \ \ \ \ \ \ \ \ \text{for all }\mathbf{a}\in\mathbf{k}\left[  G\right]
\text{ and }f\in V^{\ast}%
\]
(which defines the action of $\mathbf{k}\left[  G\right]  $ on $V^{\ast}$ in
Definition \ref{def.dual.dual.kG}) is clearly equivalent to the formula%
\begin{align*}
&  \left(  \mathbf{a}\cdot f\right)  \left(  v\right)  =\left(  f\circ
\tau_{S\left(  \mathbf{a}\right)  }\right)  \left(  v\right) \\
&  \ \ \ \ \ \ \ \ \ \ \text{for all }\mathbf{a}\in\mathbf{k}\left[  G\right]
\text{ and }f\in V^{\ast}\text{ and }v\in V
\end{align*}
(since two maps from $V$ to $\mathbf{k}$ are equal if and only if they give
equal values on every $v\in V$). But the latter formula is, in turn,
equivalent to (\ref{eq.def.dual.dual.kG.on-v}), since all $\mathbf{a}%
\in\mathbf{k}\left[  G\right]  $ and $f\in V^{\ast}$ and $v\in V$ satisfy%
\[
\left(  f\circ\tau_{S\left(  \mathbf{a}\right)  }\right)  \left(  v\right)
=f\left(  \underbrace{\tau_{S\left(  \mathbf{a}\right)  }\left(  v\right)
}_{\substack{=S\left(  \mathbf{a}\right)  \cdot v\\\text{(by the definition of
}\tau_{S\left(  \mathbf{a}\right)  }\text{)}}}\right)  =f\left(  S\left(
\mathbf{a}\right)  \cdot v\right)  .
\]

\end{remark}

Soon (in Proposition \ref{prop.dual.exist}) we will show that the above two
definitions of duals are legitimate and equivalent. But first, let us give a
simple example:

\begin{example}
\label{exa.dual.ksign}Let $G$ be a group. Consider the trivial representation
$\mathbf{k}_{\operatorname*{triv}}$ of $G$, which is the $\mathbf{k}$-module
$\mathbf{k}$ equipped with the trivial $G$-action given by%
\[
g\rightharpoonup v=v\ \ \ \ \ \ \ \ \ \ \text{for all }g\in G\text{ and }%
v\in\mathbf{k}.
\]
(This generalizes the $\mathbf{k}_{\operatorname*{triv}}$ from Example
\ref{exa.rep.Sn-rep.triv} to arbitrary groups.)

Let us find its dual $\left(  \mathbf{k}_{\operatorname*{triv}}\right)
^{\ast}$ according to Definition \ref{def.dual.dual.G}. Indeed, as a
$\mathbf{k}$-module, $\left(  \mathbf{k}_{\operatorname*{triv}}\right)
^{\ast}=\mathbf{k}^{\ast}\cong\mathbf{k}$, whereas the $G$-action on $\left(
\mathbf{k}_{\operatorname*{triv}}\right)  ^{\ast}$ is the map $\overset{\ast
}{\rightharpoonup}\ :G\times\mathbf{k}^{\ast}\rightarrow\mathbf{k}^{\ast}$
defined by the formula%
\[
g\overset{\ast}{\rightharpoonup}f:=f\circ\tau_{g^{-1}}%
\ \ \ \ \ \ \ \ \ \ \text{for all }g\in G\text{ and }f\in\mathbf{k}^{\ast}.
\]
This formula can be simplified to
\[
g\overset{\ast}{\rightharpoonup}f:=f\ \ \ \ \ \ \ \ \ \ \text{for all }g\in
G\text{ and }f\in\mathbf{k}^{\ast}%
\]
(since the triviality of the $G$-action on $\mathbf{k}_{\operatorname*{triv}}$
yields that $\tau_{g}=\operatorname*{id}\nolimits_{V}$ for each $g\in G$, and
thus $f\circ\underbrace{\tau_{g^{-1}}}_{=\operatorname*{id}\nolimits_{V}%
}=f\circ\operatorname*{id}\nolimits_{V}=f$ for each $g\in G$ and
$f\in\mathbf{k}^{\ast}$). This means that the $G$-action on $\left(
\mathbf{k}_{\operatorname*{triv}}\right)  ^{\ast}$ is again trivial. Hence, we
conclude that $\left(  \mathbf{k}_{\operatorname*{triv}}\right)  ^{\ast}%
\cong\mathbf{k}_{\operatorname*{triv}}$ (since $\left(  \mathbf{k}%
_{\operatorname*{triv}}\right)  ^{\ast}\cong\mathbf{k}$ as a $\mathbf{k}$-module).

More generally, the same arguments show that if $V$ is any trivial
representation of $G$ (that is, any representation of $G$ that satisfies
$g\rightharpoonup v=v$ for all $g\in G$ and $v\in V$), then its dual $V^{\ast
}$ is again a trivial representation of $G$.
\end{example}

It is not much harder to see that the sign representation $\mathbf{k}%
_{\operatorname*{sign}}$ of $S_{n}$ satisfies $\left(  \mathbf{k}%
_{\operatorname*{sign}}\right)  ^{\ast}\cong\mathbf{k}_{\operatorname*{sign}}$.

We can easily find the dual of the regular representation of a finite group
$G$ and of the natural representation of a symmetric group $S_{n}$:

\begin{exercise}
\label{exe.dual.reg.1}\fbox{2} Let $G$ be a finite group. Consider the left
regular representation $\mathbf{k}\left[  G\right]  $ of $G$. Prove that its
dual $\left(  \mathbf{k}\left[  G\right]  \right)  ^{\ast}$ is isomorphic to
$\mathbf{k}\left[  G\right]  $ itself, and that the $\mathbf{k}$-linear map
\begin{align*}
\mathbf{k}\left[  G\right]   &  \rightarrow\left(  \mathbf{k}\left[  G\right]
\right)  ^{\ast},\\
g  &  \mapsto g^{\ast}%
\end{align*}
(where $\left(  h^{\ast}\right)  _{h\in G}$ denotes the dual basis of the
standard basis $\left(  h\right)  _{h\in G}$ of $\mathbf{k}\left[  G\right]
$) is an isomorphism of $G$-representations.
\end{exercise}

\begin{exercise}
\label{exe.dual.nat.1}\fbox{1} Consider the natural representation
$\mathbf{k}^{n}$ of the symmetric group $S_{n}$ (defined in Example
\ref{exa.rep.Sn-rep.nat}). Prove that its dual $\left(  \mathbf{k}^{n}\right)
^{\ast}$ is isomorphic to $\mathbf{k}^{n}$ itself, and that the $\mathbf{k}%
$-linear map%
\begin{align*}
\mathbf{k}^{n}  &  \rightarrow\left(  \mathbf{k}^{n}\right)  ^{\ast},\\
e_{i}  &  \mapsto e_{i}^{\ast}%
\end{align*}
(where $\left(  e_{1}^{\ast},e_{2}^{\ast},\ldots,e_{n}^{\ast}\right)  $
denotes the dual basis of the standard basis $\left(  e_{1},e_{2},\ldots
,e_{n}\right)  $ of $\mathbf{k}^{n}$) is an isomorphism of $S_{n}$-representations.
\end{exercise}

These examples might suggest that each representation $V$ of a finite group
$G$ satisfies $V^{\ast}\cong V$; but this is not always the case (for example,
a nontrivial irrep of the cyclic group $C_{3}$ over $\mathbb{C}$ is not
isomorphic to its dual\footnote{In more detail: Let $C_{3}=\left\{
1,g,g^{2}\right\}  $ be the cyclic group of size $3$. Let $\zeta$ be one of
the two nontrivial $3$-rd roots of unity in $\mathbb{C}$ (these are $e^{2\pi
i/3}$ and $e^{-2\pi i/3}$). Let $\mathbb{C}_{\zeta}$ be the $C_{3}%
$-representation on the vector space $\mathbb{C}$ in which $g$ acts as
multiplication by $\zeta$, and let $\mathbb{C}_{\zeta^{2}}$ be the $C_{3}%
$-representation on the vector space $\mathbb{C}$ in which $g$ acts as
multiplication by $\zeta^{2}=\zeta^{-1}$. Then, $\left(  \mathbb{C}_{\zeta
}\right)  ^{\ast}$ is isomorphic to $\mathbb{C}_{\zeta^{2}}$, not to
$\mathbb{C}_{\zeta}$.}). Nevertheless, this is the case when $G$ is a
symmetric group and $\mathbf{k}$ is a field of characteristic $0$ and $V$ is
finite-dimensional; this is a nontrivial and rather useful result, which we
will prove below. We will also describe the duals of Specht modules in most
cases. Note that Specht modules are not always isomorphic to their own duals
in positive characteristic:

\begin{exercise}
\fbox{2} Let $\mathbf{k}$ be a field of characteristic $3$. Prove that
$\left(  \mathcal{S}^{\left(  2,1\right)  }\right)  ^{\ast}\ncong
\mathcal{S}^{\left(  2,1\right)  }$ as $S_{3}$-representations.
\end{exercise}

We still owe a proof of the well-definedness of $V^{\ast}$, so let us give it now:

\begin{proposition}
\label{prop.dual.exist}Let $G$ be a group. Let $V$ be a left $\mathbf{k}%
\left[  G\right]  $-module, i.e., a representation of $G$ over $\mathbf{k}$.
Then: \medskip

\textbf{(a)} The representation $V^{\ast}$ in Definition \ref{def.dual.dual.G}
is well-defined. \medskip

\textbf{(b)} The left $\mathbf{k}\left[  G\right]  $-module $V^{\ast}$ in
Definition \ref{def.dual.dual.kG} is well-defined. \medskip

\textbf{(c)} The representation $V^{\ast}$ in Definition \ref{def.dual.dual.G}
and the left $\mathbf{k}\left[  G\right]  $-module $V^{\ast}$ in Definition
\ref{def.dual.dual.kG} are the same (i.e., converting the former into a left
$\mathbf{k}\left[  G\right]  $-module using Theorem \ref{thm.rep.G-rep.mod}
yields the latter).
\end{proposition}

\begin{fineprint}
\begin{proof}
Structurally, this is similar to the proof of Proposition
\ref{prop.sign-twist.exist}, although the details of the computations differ.
We again shall first prove part \textbf{(b)}, then use it to prove the other
two parts.

However, before we start with the proof of part \textbf{(b)}, let us state
some features of $V$.

For each $\mathbf{b}\in\mathbf{k}\left[  G\right]  $, we let $\tau
_{\mathbf{b}}$ denote the map
\begin{align*}
V  &  \rightarrow V,\\
v  &  \mapsto\mathbf{b}\cdot v.
\end{align*}
This generalizes the map $\tau_{S\left(  \mathbf{a}\right)  }$ in Definition
\ref{def.dual.dual.kG}.

Recall that the action of $\mathbf{k}\left[  G\right]  $ on $V$ (viewed as the
map $\mathbf{k}\left[  G\right]  \times V\rightarrow V,\ \left(
\mathbf{b},v\right)  \mapsto\mathbf{b}\cdot v$) is $\mathbf{k}$-bilinear.
Thus, for each $\mathbf{b}\in\mathbf{k}\left[  G\right]  $, the map
$\tau_{\mathbf{b}}$ is $\mathbf{k}$-linear, and depends $\mathbf{k}$-linearly
on $\mathbf{b}$ (meaning that $\tau_{\lambda\mathbf{b}+\mu\mathbf{c}}%
=\lambda\tau_{\mathbf{b}}+\mu\tau_{\mathbf{c}}$ for any $\lambda,\mu
\in\mathbf{k}$ and any $\mathbf{b},\mathbf{c}\in\mathbf{k}\left[  G\right]  $).

Furthermore, for any $\mathbf{b},\mathbf{c}\in\mathbf{k}\left[  G\right]  $,
we have%
\begin{equation}
\tau_{\mathbf{bc}}=\tau_{\mathbf{b}}\circ\tau_{\mathbf{c}}
\label{pf.prop.dual.exist.taubc}%
\end{equation}
\footnote{\textit{Proof.} Let $\mathbf{b},\mathbf{c}\in\mathbf{k}\left[
G\right]  $. Let $v\in V$. Then, $\tau_{\mathbf{bc}}\left(  v\right)
=\mathbf{bc}\cdot v$ (by the definition of $\tau_{\mathbf{bc}}$) and
$\tau_{\mathbf{c}}\left(  v\right)  =\mathbf{c}\cdot v$ (similarly). Now,%
\begin{align*}
\left(  \tau_{\mathbf{b}}\circ\tau_{\mathbf{c}}\right)  \left(  v\right)   &
=\tau_{\mathbf{b}}\left(  \tau_{\mathbf{c}}\left(  v\right)  \right)
=\mathbf{b}\cdot\underbrace{\tau_{\mathbf{c}}\left(  v\right)  }%
_{=\mathbf{c}\cdot v}\ \ \ \ \ \ \ \ \ \ \left(  \text{by the definition of
}\tau_{\mathbf{b}}\right) \\
&  =\mathbf{b}\cdot\left(  \mathbf{c}\cdot v\right)  =\mathbf{bc}\cdot
v\ \ \ \ \ \ \ \ \ \ \left(  \text{since }V\text{ is a left }\mathbf{k}\left[
G\right]  \text{-module}\right)  .
\end{align*}
Comparing this with $\tau_{\mathbf{bc}}\left(  v\right)  =\mathbf{bc}\cdot v$,
we obtain $\tau_{\mathbf{bc}}\left(  v\right)  =\left(  \tau_{\mathbf{b}}%
\circ\tau_{\mathbf{c}}\right)  \left(  v\right)  $.
\par
Forget that we fixed $v$. We thus have proved that $\tau_{\mathbf{bc}}\left(
v\right)  =\left(  \tau_{\mathbf{b}}\circ\tau_{\mathbf{c}}\right)  \left(
v\right)  $ for each $v\in V$. In other words, $\tau_{\mathbf{bc}}%
=\tau_{\mathbf{b}}\circ\tau_{\mathbf{c}}$. This proves
(\ref{pf.prop.dual.exist.taubc}).}. Moreover, the unity $1$ of $\mathbf{k}%
\left[  G\right]  $ satisfies $\tau_{1}=\operatorname*{id}$ (since each $v\in
V$ satisfies $\tau_{1}\left(  v\right)  =1\cdot v=v=\operatorname*{id}\left(
v\right)  $). Hence, for each $\lambda\in\mathbf{k}$, we have%
\begin{align}
\tau_{\lambda\cdot1}  &  =\lambda\underbrace{\tau_{1}}_{=\operatorname*{id}%
}\ \ \ \ \ \ \ \ \ \ \left(  \text{since }\tau_{\mathbf{b}}\text{ depends
}\mathbf{k}\text{-linearly on }\mathbf{b}\right) \nonumber\\
&  =\lambda\operatorname*{id}. \label{pf.prop.dual.exist.taulam}%
\end{align}

Now, we step to the proofs of parts \textbf{(b)}, \textbf{(a)} and
\textbf{(c)}. \medskip

\textbf{(b)} Consider the map%
\begin{align*}
\mathbf{k}\left[  G\right]  \times V^{\ast}  &  \rightarrow V^{\ast},\\
\left(  \mathbf{a},f\right)   &  \mapsto f\circ\tau_{S\left(  \mathbf{a}%
\right)  }.
\end{align*}
This is the map that is intended to serve as the action of $\mathbf{k}\left[
G\right]  $ on $V^{\ast}$ (according to Definition \ref{def.dual.dual.kG}).
Hence, let us denote this map by $\overset{V^{\ast}}{\cdot}$ (that is, let us
write $\mathbf{a}\overset{V^{\ast}}{\cdot}f$ for the image of a pair $\left(
\mathbf{a},f\right)  \in\mathbf{k}\left[  G\right]  \times V^{\ast}$ under
this map). Thus,%
\begin{equation}
\mathbf{a}\overset{V^{\ast}}{\cdot}f=f\circ\tau_{S\left(  \mathbf{a}\right)  }
\label{pf.prop.dual.exist.b.1}%
\end{equation}
for all $\mathbf{a}\in\mathbf{k}\left[  G\right]  $ and $f\in V^{\ast}$.

The map $\overset{V^{\ast}}{\cdot}$ is easily seen to be $\mathbf{k}$-bilinear
(since composition of $\mathbf{k}$-linear maps is $\mathbf{k}$-bilinear; since
the map $S$ is $\mathbf{k}$-linear; and since the linear map $\tau
_{\mathbf{b}}:V\rightarrow V$ depends $\mathbf{k}$-linearly on $\mathbf{b}$).

In Definition \ref{def.dual.dual.kG}, we defined the left $\mathbf{k}\left[
G\right]  $-module $V^{\ast}$ to be the additive group $V^{\ast}$ (the dual of
the $\mathbf{k}$-module $V$), equipped with the map $\overset{V^{\ast}}{\cdot
}$ as the action of $\mathbf{k}\left[  G\right]  $ on $V^{\ast}$. In order to
prove that this left $\mathbf{k}\left[  G\right]  $-module $V^{\ast}$ is
well-defined, we must show that this structure satisfies all the module axioms
(see Definition \ref{def.mod.leftmod}). Most of these axioms are
straightforward to check (e.g., both right and left distributivity laws follow
from the $\mathbf{k}$-bilinearity of the map $\overset{V^{\ast}}{\cdot}$), so
we only focus on two axioms: the associativity law and the \textquotedblleft%
$1m=m$\textquotedblright\ axiom. Let us check them both:

\begin{itemize}
\item \textit{The associativity law:} We must prove the associativity law for
the left $\mathbf{k}\left[  G\right]  $-module $V^{\ast}$. In other words, we
must prove that $\left(  rs\right)  \overset{V^{\ast}}{\cdot}%
m=r\overset{V^{\ast}}{\cdot}\left(  s\overset{V^{\ast}}{\cdot}m\right)  $ for
all $r,s\in\mathbf{k}\left[  G\right]  $ and $m\in V^{\ast}$. So let
$r,s\in\mathbf{k}\left[  G\right]  $ and $m\in V^{\ast}$ be arbitrary. Then,
the definition of $\overset{V^{\ast}}{\cdot}$ yields $s\overset{V^{\ast
}}{\cdot}m=m\circ\tau_{S\left(  s\right)  }$ and $r\overset{V^{\ast}}{\cdot
}\left(  s\overset{V^{\ast}}{\cdot}m\right)  =\left(  s\overset{V^{\ast
}}{\cdot}m\right)  \circ\tau_{S\left(  r\right)  }$ and $\left(  rs\right)
\overset{V^{\ast}}{\cdot}m=m\circ\tau_{S\left(  rs\right)  }$. But $S$ is a
$\mathbf{k}$-algebra anti-morphism (by Definition \ref{def.dual.S}
\textbf{(a)}), and thus we have $S\left(  rs\right)  =S\left(  s\right)  \cdot
S\left(  r\right)  $. Hence, $\tau_{S\left(  rs\right)  }=\tau_{S\left(
s\right)  \cdot S\left(  r\right)  }=\tau_{S\left(  s\right)  }\circ
\tau_{S\left(  r\right)  }$ (by (\ref{pf.prop.dual.exist.taubc})). Thus,%
\[
\left(  rs\right)  \overset{V^{\ast}}{\cdot}m=m\circ\underbrace{\tau_{S\left(
rs\right)  }}_{=\tau_{S\left(  s\right)  }\circ\tau_{S\left(  r\right)  }%
}=m\circ\tau_{S\left(  s\right)  }\circ\tau_{S\left(  r\right)  }.
\]
Comparing this with
\[
r\overset{V^{\ast}}{\cdot}\left(  s\overset{V^{\ast}}{\cdot}m\right)
=\underbrace{\left(  s\overset{V^{\ast}}{\cdot}m\right)  }_{=m\circ
\tau_{S\left(  s\right)  }}\circ\,\tau_{S\left(  r\right)  }=m\circ
\tau_{S\left(  s\right)  }\circ\tau_{S\left(  r\right)  },
\]
we obtain $\left(  rs\right)  \overset{V^{\ast}}{\cdot}m=r\overset{V^{\ast
}}{\cdot}\left(  s\overset{V^{\ast}}{\cdot}m\right)  $. Thus, the
associativity law for the $\mathbf{k}\left[  G\right]  $-module $V^{\ast}$ is proved.

\item \textit{The }$1m=m$ \textit{axiom:} We must prove the \textquotedblleft%
$1m=m$\textquotedblright\ axiom for the left $\mathbf{k}\left[  G\right]
$-module $V^{\ast}$. In other words, we must prove that $1\overset{V^{\ast
}}{\cdot}m=m$ for all $m\in V^{\ast}$ (where $1$ is the unity of
$\mathbf{k}\left[  G\right]  $). So let $m\in V^{\ast}$. Then, the definition
of $\overset{V^{\ast}}{\cdot}$ yields $1\overset{V^{\ast}}{\cdot}m=m\circ
\tau_{S\left(  1\right)  }$. But $S$ is a $\mathbf{k}$-algebra anti-morphism
(by Definition \ref{def.dual.S} \textbf{(a)}), and thus we have $S\left(
1\right)  =1$. Hence, $\tau_{S\left(  1\right)  }=\tau_{1}=\operatorname*{id}%
$. Thus, $1\overset{V^{\ast}}{\cdot}m=m\circ\underbrace{\tau_{S\left(
1\right)  }}_{=\operatorname*{id}}=m$. Thus, we have proved the
\textquotedblleft$1m=m$\textquotedblright\ axiom for the left $\mathbf{k}%
\left[  G\right]  $-module $V^{\ast}$.
\end{itemize}

As we said, the other axioms are even easier to check. Thus, the left
$\mathbf{k}\left[  G\right]  $-module $V^{\ast}$ is well-defined. This proves
Proposition \ref{prop.dual.exist} \textbf{(b)}. \medskip

\textbf{(a)} Instead of proving the claim directly, we shall derive it from
part \textbf{(b)} (which we have already proved).

First, we introduce some notations (exclusively for this proof):

\begin{itemize}
\item We shall write $V_{1}^{\ast}$ for the representation $V^{\ast}$ defined
in Definition \ref{def.dual.dual.G} (even though we have not yet shown that it
is well-defined). Our goal is to prove that this representation $V_{1}^{\ast}$
is well-defined.

\item We shall write $V_{2}^{\ast}$ for the left $\mathbf{k}\left[  G\right]
$-module $V^{\ast}$ defined in Definition \ref{def.dual.dual.kG}. Proposition
\ref{prop.dual.exist} \textbf{(b)} (which we have already proved) shows that
the latter left $\mathbf{k}\left[  G\right]  $-module $V_{2}^{\ast}$ is
well-defined. We reinterpret this left $\mathbf{k}\left[  G\right]  $-module
$V_{2}^{\ast}$ as a representation of $G$ (using Corollary
\ref{cor.rep.G-rep.mod}, as usual).
\end{itemize}

We shall now show that this representation $V_{2}^{\ast}$ is precisely the
(purported) representation $V_{1}^{\ast}$. This will immediately entail that
the latter representation is well-defined, and thus Proposition
\ref{prop.dual.exist} \textbf{(a)} will follow (since $V_{1}^{\ast}$ is
precisely the representation $V^{\ast}$ defined in Definition
\ref{def.dual.dual.G}).

First, we recall that $S$ is a $\mathbf{k}$-algebra anti-morphism (by
Definition \ref{def.dual.S} \textbf{(a)}), and thus we have $S\left(
1_{\mathbf{k}\left[  G\right]  }\right)  =1_{\mathbf{k}\left[  G\right]  }$.
But $S$ is also $\mathbf{k}$-linear. Hence, for each $\lambda\in\mathbf{k}$,
we have%
\begin{equation}
S\left(  \lambda1_{\mathbf{k}\left[  G\right]  }\right)  =\lambda
\cdot\underbrace{S\left(  1_{\mathbf{k}\left[  G\right]  }\right)
}_{=1_{\mathbf{k}\left[  G\right]  }}=\lambda1_{\mathbf{k}\left[  G\right]  }.
\label{pf.prop.dual.exist.a.0}%
\end{equation}

Next, we recall that $V_{2}^{\ast}$ is the left $\mathbf{k}\left[  G\right]
$-module $V^{\ast}$ from Definition \ref{def.dual.dual.kG}. Hence, the action
of $\mathbf{k}\left[  G\right]  $ on $V_{2}^{\ast}$ is given by%
\begin{equation}
\mathbf{a}\overset{V_{2}^{\ast}}{\cdot}f:=f\circ\tau_{S\left(  \mathbf{a}%
\right)  } \label{pf.prop.dual.exist.a.1}%
\end{equation}
for all $\mathbf{a}\in\mathbf{k}\left[  G\right]  $ and $f\in V^{\ast}$ (by
Definition \ref{def.dual.dual.kG}).

Now, we shall show that our $G$-representation $V_{2}^{\ast}$ (which we
obtained from the left $\mathbf{k}\left[  G\right]  $-module $V_{2}^{\ast}$
defined in Definition \ref{def.dual.dual.kG}) is the (purported)
representation $V_{1}^{\ast}$ defined in Definition \ref{def.dual.dual.G}. Indeed:

\begin{itemize}
\item The underlying $\mathbf{k}$-module of our $G$-representation
$V_{2}^{\ast}$ is the dual $V^{\ast}$ of the $\mathbf{k}$-module $V$.

[\textit{Proof:} We must show that $V_{2}^{\ast}=V^{\ast}$ as $\mathbf{k}%
$-modules. As an additive group, $V_{2}^{\ast}$ is just the dual $V^{\ast}$ of
the $\mathbf{k}$-module $V$ (by its definition). In other words, our
$G$-representation $V_{2}^{\ast}$ has the same addition and zero as $V^{\ast}%
$. Thus, it remains to show that $V_{2}^{\ast}$ also has the same scaling (by
elements of $\mathbf{k}$) as $V^{\ast}$. In other words, it remains to prove
that every $\lambda\in\mathbf{k}$ and $f\in V^{\ast}$ satisfy $\lambda
\overset{V_{2}^{\ast}}{\cdot}f=\lambda f$ (where the right hand side is
computed in the dual $V^{\ast}$ of the $\mathbf{k}$-module $V$). So let us
prove this. For any $\lambda\in\mathbf{k}$ and $f\in V^{\ast}$, we have%
\begin{align*}
\lambda\overset{V_{2}^{\ast}}{\cdot}f  &  =\left(  \lambda1_{\mathbf{k}\left[
G\right]  }\right)  \overset{V_{2}^{\ast}}{\cdot}f\ \ \ \ \ \ \ \ \ \ \left(
\begin{array}
[c]{c}%
\text{since the }\mathbf{k}\text{-module structure on }V_{2}^{\ast}\text{
is}\\
\text{obtained from the }\mathbf{k}\left[  G\right]  \text{-module
structure}\\
\text{by restriction of scalars}%
\end{array}
\right) \\
&  =f\circ\tau_{S\left(  \lambda1_{\mathbf{k}\left[  G\right]  }\right)
}\ \ \ \ \ \ \ \ \ \ \left(  \text{by (\ref{pf.prop.dual.exist.a.1}), applied
to }\mathbf{a}=\lambda1_{\mathbf{k}\left[  G\right]  }\right) \\
&  =f\circ\underbrace{\tau_{\lambda1_{\mathbf{k}\left[  G\right]  }}%
}_{\substack{=\tau_{\lambda1}\\=\lambda\operatorname*{id}\\\text{(by
(\ref{pf.prop.dual.exist.taulam}))}}}\ \ \ \ \ \ \ \ \ \ \left(  \text{since
(\ref{pf.prop.dual.exist.a.0}) yields }S\left(  \lambda1_{\mathbf{k}\left[
G\right]  }\right)  =\lambda1_{\mathbf{k}\left[  G\right]  }\right) \\
&  =f\circ\left(  \lambda\operatorname*{id}\right)  =\lambda f
\end{align*}
(since each $v\in V$ satisfies $\left(  f\circ\left(  \lambda
\operatorname*{id}\right)  \right)  \left(  v\right)  =f\left(
\underbrace{\left(  \lambda\operatorname*{id}\right)  \left(  v\right)
}_{=\lambda v}\right)  =f\left(  \lambda v\right)  =\lambda f\left(  v\right)
$ due to the $\mathbf{k}$-linearity of $f$). Thus, we have shown that every
$\lambda\in\mathbf{k}$ and $f\in V^{\ast}$ satisfy $\lambda\overset{V_{2}%
^{\ast}}{\cdot}f=\lambda f$. In other words, the $\mathbf{k}$-module
$V_{2}^{\ast}$ has the same scaling as $V^{\ast}$. Since $V_{2}^{\ast}$ also
has the same addition and zero as $V^{\ast}$, we thus conclude that the
$\mathbf{k}$-module $V_{2}^{\ast}$ is identical with the $\mathbf{k}$-module
$V^{\ast}$.]

\item The action of $G$ on our $G$-representation $V_{2}^{\ast}$ is exactly
the map $\overset{\ast}{\rightharpoonup}\ :G\times V^{\ast}\rightarrow
V^{\ast}$ defined in Definition \ref{def.dual.dual.G}.

[\textit{Proof:} Let $\underset{2}{\overset{\ast}{\rightharpoonup}}\ :G\times
V_{2}^{\ast}\rightarrow V_{2}^{\ast}$ be the action of $G$ on the
representation $V_{2}^{\ast}$ of $G$. Then, for any $g\in G$ and $f\in
V_{2}^{\ast}$, we have%
\begin{align*}
g\underset{2}{\overset{\ast}{\rightharpoonup}}f  &  =g\overset{V_{2}^{\ast
}}{\cdot}f=f\circ\tau_{S\left(  g\right)  }\ \ \ \ \ \ \ \ \ \ \left(
\text{by (\ref{pf.prop.dual.exist.a.1}), applied to }\mathbf{a}=g\right) \\
&  =f\circ\tau_{g^{-1}}\ \ \ \ \ \ \ \ \ \ \left(
\begin{array}
[c]{c}%
\text{since the definition of }S\text{ yields }S\left(  g\right)  =g^{-1}\\
\text{(because }g\in G\text{)}%
\end{array}
\right)  .
\end{align*}
But this is precisely the formula for the map $\overset{\ast}{\rightharpoonup
}\ :G\times V^{\ast}\rightarrow V^{\ast}$ defined in Definition
\ref{def.dual.dual.G}. Thus, the map $\underset{2}{\overset{\ast
}{\rightharpoonup}}\ :G\times V_{2}^{\ast}\rightarrow V_{2}^{\ast}$ is
precisely the map $\overset{\ast}{\rightharpoonup}\ :G\times V^{\ast
}\rightarrow V^{\ast}$ defined in Definition \ref{def.dual.dual.G}. In other
words, the action of $G$ on our $G$-representation $V_{2}^{\ast}$ is exactly
the map $\overset{\ast}{\rightharpoonup}\ :G\times V^{\ast}\rightarrow
V^{\ast}$ defined in Definition \ref{def.dual.dual.G}.]
\end{itemize}

Thus, we have shown that our representation $V_{2}^{\ast}$ of $G$ has the same
underlying $\mathbf{k}$-module as $V^{\ast}$, and its action of $G$ is
precisely the map $\overset{\ast}{\rightharpoonup}\ :G\times V^{\ast
}\rightarrow V^{\ast}$ from Definition \ref{def.dual.dual.G}. Consequently,
this representation $V_{2}^{\ast}$ is exactly the representation $V^{\ast}$
defined in Definition \ref{def.dual.dual.G} (i.e., the representation that we
now call $V_{1}^{\ast}$). Therefore, the latter representation $V^{\ast}$ is
well-defined (since the former representation is well-defined). This proves
Proposition \ref{prop.dual.exist} \textbf{(a)}. \medskip

\textbf{(c)} We must prove that the representation $V^{\ast}$ in Definition
\ref{def.dual.dual.G} and the left $\mathbf{k}\left[  G\right]  $-module
$V^{\ast}$ in Definition \ref{def.dual.dual.kG} are the same. In our above
proof of part \textbf{(a)}, we have denoted the former representation by
$V_{1}^{\ast}$ and the latter $\mathbf{k}\left[  G\right]  $-module by
$V_{2}^{\ast}$; thus, we must prove that $V_{1}^{\ast}$ and $V_{2}^{\ast}$ are
the same.

But we have essentially done this already. Indeed, in our above proof of part
\textbf{(a)}, we have shown that the representation $V_{2}^{\ast}$ is exactly
the representation $V^{\ast}$ defined in Definition \ref{def.dual.dual.G}; but
the latter representation is precisely the one that we called $V_{1}^{\ast}$.
Hence, we have shown that $V_{1}^{\ast}$ and $V_{2}^{\ast}$ are the same. This
proves Proposition \ref{prop.dual.exist} \textbf{(c)}.
\end{proof}
\end{fineprint}

Parts \textbf{(a)} and \textbf{(b)} of Proposition \ref{prop.dual.exist} show
that our two definitions of $V^{\ast}$ are well-defined, and Proposition
\ref{prop.dual.exist} \textbf{(c)} shows that they are equivalent. We now turn
to some simple properties of duals. For the following proposition, we recall
the notion of the dual $f^{\ast}$ of a $\mathbf{k}$-linear map $f$, as defined
in Proposition \ref{prop.dual.functorial} \textbf{(a)}:

\begin{proposition}
\label{prop.dual.dual-map}Let $G$ be a group. Let $V$ be a representation of
$G$. Then: \medskip

\textbf{(a)} For each $g\in G$, we have%
\[
\left(  \tau_{g}\text{ on }V^{\ast}\right)  =\left(  \tau_{g^{-1}}\text{ on
}V\right)  ^{\ast}.
\]
Here, \textquotedblleft$\tau_{g}$ on $V^{\ast}$\textquotedblright\ means the
$\mathbf{k}$-linear map $V^{\ast}\rightarrow V^{\ast},\ f\mapsto
g\rightharpoonup f$, whereas \textquotedblleft$\tau_{g^{-1}}$ on
$V$\textquotedblright\ means the $\mathbf{k}$-linear map $V\rightarrow
V,\ v\mapsto g^{-1}\rightharpoonup v$. \medskip

\textbf{(b)} More generally, for each $\mathbf{a}\in\mathbf{k}\left[
G\right]  $, we have%
\[
\left(  \tau_{\mathbf{a}}\text{ on }V^{\ast}\right)  =\left(  \tau_{S\left(
\mathbf{a}\right)  }\text{ on }V\right)  ^{\ast}.
\]
Here, \textquotedblleft$\tau_{\mathbf{a}}$ on $V^{\ast}$\textquotedblright%
\ means the $\mathbf{k}$-linear map $V^{\ast}\rightarrow V^{\ast}%
,\ f\mapsto\mathbf{a}\cdot f$, whereas \textquotedblleft$\tau_{S\left(
\mathbf{a}\right)  }$ on $V$\textquotedblright\ means the $\mathbf{k}$-linear
map $V\rightarrow V,\ v\mapsto S\left(  \mathbf{a}\right)  \cdot v$.
\end{proposition}

\begin{fineprint}
\begin{proof}
These are essentially just the definitions of the action of $G$ (resp.
$\mathbf{k}\left[  G\right]  $) on $V^{\ast}$, restated in terms of the dual
of a map. In more detail: \medskip

\textbf{(b)} Let $\mathbf{a}\in\mathbf{k}\left[  G\right]  $. Then, for each
$f\in V^{\ast}$, we have%
\begin{align*}
\left(  \tau_{\mathbf{a}}\text{ on }V^{\ast}\right)  \left(  f\right)   &
=\mathbf{a}\cdot f\ \ \ \ \ \ \ \ \ \ \left(  \text{by the definition of
}\left(  \tau_{\mathbf{a}}\text{ on }V^{\ast}\right)  \right) \\
&  =f\circ\tau_{S\left(  \mathbf{a}\right)  }\ \ \ \ \ \ \ \ \ \ \left(
\begin{array}
[c]{c}%
\text{by the definition of the }\mathbf{k}\left[  G\right]  \text{-action on
}V^{\ast}\\
\text{(Definition \ref{def.dual.dual.kG})}%
\end{array}
\right) \\
&  =f\circ\left(  \tau_{S\left(  \mathbf{a}\right)  }\text{ on }V\right)
\ \ \ \ \ \ \ \ \ \ \left(
\begin{array}
[c]{c}%
\text{since the }\tau_{S\left(  \mathbf{a}\right)  }\text{ in Definition
\ref{def.dual.dual.kG}}\\
\text{is what we call }\left(  \tau_{S\left(  \mathbf{a}\right)  }\text{ on
}V\right)  \text{ here}%
\end{array}
\right) \\
&  =\left(  \tau_{S\left(  \mathbf{a}\right)  }\text{ on }V\right)  ^{\ast
}\left(  f\right)
\end{align*}
(since the definition of the dual $\left(  \tau_{S\left(  \mathbf{a}\right)
}\text{ on }V\right)  ^{\ast}$ yields $\left(  \tau_{S\left(  \mathbf{a}%
\right)  }\text{ on }V\right)  ^{\ast}\left(  f\right)  =f\circ\left(
\tau_{S\left(  \mathbf{a}\right)  }\text{ on }V\right)  $). In other words,
$\left(  \tau_{\mathbf{a}}\text{ on }V^{\ast}\right)  =\left(  \tau_{S\left(
\mathbf{a}\right)  }\text{ on }V\right)  ^{\ast}$. This proves Proposition
\ref{prop.dual.dual-map} \textbf{(b)}. \medskip

\textbf{(a)} This is analogous to part \textbf{(b)}, except that we are using
Definition \ref{def.dual.dual.G} instead of Definition \ref{def.dual.dual.kG}
as the definition of $V^{\ast}$ here. Alternatively, because Definition
\ref{def.dual.dual.G} and Definition \ref{def.dual.dual.kG} are equivalent, we
can derive this from Proposition \ref{prop.dual.dual-map} \textbf{(b)}
(applied to $\mathbf{a}=g$).
\end{proof}
\end{fineprint}

\begin{proposition}
\label{prop.dual.twice-back}Let $G$ be a group. Let $V$ be a representation of
$G$. Then: \medskip

\textbf{(a)} The map $\iota:V\rightarrow V^{\ast\ast}$ from Proposition
\ref{prop.dual.linalg1} \textbf{(a)} is a morphism of $G$-representations.
\medskip

\textbf{(b)} Thus, if $V$ has a finite basis, then $V^{\ast\ast}\cong V$ as
$G$-representations.
\end{proposition}

\begin{fineprint}
\begin{proof}
\textbf{(a)} For any $G$-representation $W$ and any $\mathbf{b}\in
\mathbf{k}\left[  G\right]  $, we let $\tau_{\mathbf{b}}^{W}$ denote the
$\mathbf{k}$-linear map $W\rightarrow W,\ v\mapsto\mathbf{b}v$. (This map
would have been previously denoted by $\tau_{\mathbf{b}}$, but we want to make
the role of $W$ explicit.) Thus, the action of $\mathbf{k}\left[  G\right]  $
on $V^{\ast}$ is given by%
\begin{equation}
\mathbf{a}\cdot f=f\circ\tau_{S\left(  \mathbf{a}\right)  }^{V}%
\ \ \ \ \ \ \ \ \ \ \text{for all }\mathbf{a}\in\mathbf{k}\left[  G\right]
\text{ and }f\in V^{\ast} \label{pf.prop.dual.twice-back.af1}%
\end{equation}
(indeed, this is precisely the formula $\mathbf{a}\cdot f=f\circ\tau_{S\left(
\mathbf{a}\right)  }$ from Definition \ref{def.dual.dual.kG}, since the map
that was called $\tau_{S\left(  \mathbf{a}\right)  }$ in Definition
\ref{def.dual.dual.kG} is now called $\tau_{S\left(  \mathbf{a}\right)  }^{V}%
$). Applying the same reasoning to $V^{\ast}$ instead of $V$, we see that the
action of $\mathbf{k}\left[  G\right]  $ on $V^{\ast\ast}$ is given by%
\begin{equation}
\mathbf{a}\cdot f=f\circ\tau_{S\left(  \mathbf{a}\right)  }^{V^{\ast}%
}\ \ \ \ \ \ \ \ \ \ \text{for all }\mathbf{a}\in\mathbf{k}\left[  G\right]
\text{ and }f\in V^{\ast\ast}. \label{pf.prop.dual.twice-back.af2}%
\end{equation}

The definition of $\iota$ yields that%
\begin{equation}
\left(  \iota\left(  v\right)  \right)  \left(  g\right)  =g\left(  v\right)
\ \ \ \ \ \ \ \ \ \ \text{for each }g\in V^{\ast}\text{ and each }v\in V.
\label{pf.prop.dual.twice-back.a.iv}%
\end{equation}

We must prove that the map $\iota:V\rightarrow V^{\ast\ast}$ is a morphism of
$G$-representations. Equivalently, we must prove that this map $\iota$ is left
$\mathbf{k}\left[  G\right]  $-linear. But we already know (from Proposition
\ref{prop.dual.linalg1} \textbf{(a)}) that it is $\mathbf{k}$-linear. Thus,
all we need to show is that it satisfies $\iota\left(  \mathbf{a}\cdot
v\right)  =\mathbf{a}\cdot\iota\left(  v\right)  $ for all $\mathbf{a}%
\in\mathbf{k}\left[  G\right]  $ and all $v\in V$.

Let us show this.

Let $\mathbf{a}\in\mathbf{k}\left[  G\right]  $ and $v\in V$. We must show
that $\iota\left(  \mathbf{a}\cdot v\right)  =\mathbf{a}\cdot\iota\left(
v\right)  $.

Definition \ref{def.dual.S} \textbf{(b)} says that $S\circ
S=\operatorname*{id}$. Thus, $S\left(  S\left(  \mathbf{a}\right)  \right)
=\mathbf{a}$.

Applying (\ref{pf.prop.dual.twice-back.af2}) to $f=\iota\left(  v\right)  $,
we find $\mathbf{a}\cdot\iota\left(  v\right)  =\iota\left(  v\right)
\circ\tau_{S\left(  \mathbf{a}\right)  }^{V^{\ast}}$. Hence, for each $f\in
V^{\ast}$, we have%
\begin{align*}
\underbrace{\left(  \mathbf{a}\cdot\iota\left(  v\right)  \right)  }%
_{=\iota\left(  v\right)  \circ\tau_{S\left(  \mathbf{a}\right)  }^{V^{\ast}}%
}\left(  f\right)   &  =\left(  \iota\left(  v\right)  \circ\tau_{S\left(
\mathbf{a}\right)  }^{V^{\ast}}\right)  \left(  f\right)  =\left(
\iota\left(  v\right)  \right)  \underbrace{\left(  \tau_{S\left(
\mathbf{a}\right)  }^{V^{\ast}}\left(  f\right)  \right)  }%
_{\substack{=S\left(  \mathbf{a}\right)  \cdot f\\\text{(by the definition of
}\tau_{S\left(  \mathbf{a}\right)  }^{V^{\ast}}\text{)}}}=\left(  \iota\left(
v\right)  \right)  \left(  S\left(  \mathbf{a}\right)  \cdot f\right) \\
&  =\underbrace{\left(  S\left(  \mathbf{a}\right)  \cdot f\right)
}_{\substack{=f\circ\tau_{S\left(  S\left(  \mathbf{a}\right)  \right)  }%
^{V}\\\text{(by (\ref{pf.prop.dual.twice-back.af1}), applied to }S\left(
\mathbf{a}\right)  \\\text{instead of }\mathbf{a}\text{)}}}\left(  v\right)
\ \ \ \ \ \ \ \ \ \ \left(  \text{by (\ref{pf.prop.dual.twice-back.a.iv}),
applied to }g=S\left(  \mathbf{a}\right)  \cdot f\right) \\
&  =\left(  f\circ\tau_{S\left(  S\left(  \mathbf{a}\right)  \right)  }%
^{V}\right)  \left(  v\right)  =\left(  f\circ\tau_{\mathbf{a}}^{V}\right)
\left(  v\right)  \ \ \ \ \ \ \ \ \ \ \left(  \text{since }S\left(  S\left(
\mathbf{a}\right)  \right)  =\mathbf{a}\right) \\
&  =f\left(  \underbrace{\tau_{\mathbf{a}}^{V}\left(  v\right)  }%
_{\substack{=\mathbf{a}\cdot v\\\text{(by the definition of }\tau_{\mathbf{a}%
}^{V}\text{)}}}\right)  =f\left(  \mathbf{a}\cdot v\right)  =\left(
\iota\left(  \mathbf{a}\cdot v\right)  \right)  \left(  f\right)
\end{align*}
(since (\ref{pf.prop.dual.twice-back.a.iv}) (applied to $\mathbf{a}\cdot v$
and $f$ instead of $v$ and $g$) yields $\left(  \iota\left(  \mathbf{a}\cdot
v\right)  \right)  \left(  f\right)  =f\left(  \mathbf{a}\cdot v\right)  $).
In other words, $\mathbf{a}\cdot\iota\left(  v\right)  =\iota\left(
\mathbf{a}\cdot v\right)  $. In other words, $\iota\left(  \mathbf{a}\cdot
v\right)  =\mathbf{a}\cdot\iota\left(  v\right)  $. As we explained above,
this completes the proof of Proposition \ref{prop.dual.twice-back}
\textbf{(a)}. \medskip

\textbf{(b)} Assume that $V$ has a finite basis. Then, Proposition
\ref{prop.dual.linalg1} \textbf{(b)} shows that the map $\iota:V\rightarrow
V^{\ast\ast}$ is an isomorphism of $\mathbf{k}$-modules, hence is invertible.
But Proposition \ref{prop.dual.twice-back} \textbf{(a)} shows that this map
$\iota$ is a morphism of $G$-representations. Thus, $\iota$ is an invertible
morphism of $G$-representations, hence an isomorphism of $G$-representations.
Therefore, $V^{\ast\ast}\cong V$. This proves Proposition
\ref{prop.dual.twice-back} \textbf{(b)}.
\end{proof}
\end{fineprint}

We give two more basic properties of duals without proof:

\begin{proposition}
\label{prop.dual.dirsums}Let $G$ be a group. Let $V$ and $W$ be two
representations of $G$. Then, $\left(  V\oplus W\right)  ^{\ast}\cong V^{\ast
}\oplus W^{\ast}$ canonically. (To be more specific: The standard $\mathbf{k}%
$-module isomorphism $\left(  V\oplus W\right)  ^{\ast}\rightarrow V^{\ast
}\oplus W^{\ast}$ is an isomorphism of $G$-representations.)
\end{proposition}

\begin{proposition}
\label{prop.dual.mors}Let $G$ be a group. Let $V$ and $W$ be two
representations of $G$. Let $f:V\rightarrow W$ be a map. Then: \medskip

\textbf{(a)} If $f$ is a morphism of representations from $V$ to $W$, then its
dual $f^{\ast}$ (as defined in Proposition \ref{prop.dual.functorial}
\textbf{(a)}) is a morphism of representations from $W^{\ast}$ to $V^{\ast}$.
\medskip

\textbf{(b)} If $f$ is an isomorphism of representations from $V$ to $W$, then
its dual $f^{\ast}$ is an isomorphism of representations from $W^{\ast}$ to
$V^{\ast}$.
\end{proposition}

\begin{corollary}
\label{cor.dual.iso}Let $G$ be a group. Let $V$ and $W$ be two representations
of $G$. If $V\cong W$, then $V^{\ast}\cong W^{\ast}$.
\end{corollary}

\begin{exercise}
\textbf{(a)} \fbox{1} Prove Proposition \ref{prop.dual.dirsums}. \medskip

\textbf{(b)} \fbox{2} Prove Proposition \ref{prop.dual.mors} and Corollary
\ref{cor.dual.iso}.
\end{exercise}

\begin{remark}
\label{rmk.dual.functor}Proposition \ref{prop.dual.mors} \textbf{(a)}
(combined with the facts that $\left(  f\circ g\right)  ^{\ast}=g^{\ast}\circ
f^{\ast}$ and $\operatorname*{id}\nolimits^{\ast}=\operatorname*{id}$) shows
that we can define a \textquotedblleft dualization functor\textquotedblright%
\ on the category of $G$-representations for any group $G$. This is a
contravariant functor from the category of $G$-representations to itself. On
objects, it acts by sending each $G$-representation $V$ to $V^{\ast}$. On
morphisms, it acts by sending each morphism $f$ to $f^{\ast}$.
\end{remark}

\subsubsection{$G$-invariant bilinear forms and duals}

Bilinear forms are often the tool of choice to understand duals of
$\mathbf{k}$-modules. To understand duals of $G$-representations, we need a
slightly restrictive kind of bilinear forms, namely the $G$\emph{-invariant
bilinear forms}:

\begin{definition}
\label{def.dual.bilf.G-inv}Let $G$ be a group. Let $U$ and $V$ be two
representations of $G$ over $\mathbf{k}$. Let $f:U\times V\rightarrow
\mathbf{k}$ be a bilinear form. Then, $f$ is said to be $G$\emph{-invariant}
if and only if every $g\in G$ and $u\in U$ and $v\in V$ satisfy the equality%
\[
f\left(  gu,\ gv\right)  =f\left(  u,\ v\right)  .
\]
(You can rewrite this equality as $f\left(  g\rightharpoonup
u,\ g\rightharpoonup v\right)  =f\left(  u,\ v\right)  $ if you wish.)
\end{definition}

\begin{example}
\label{exa.dual.bilf.G-inv.nat}Consider the natural representation
$\mathbf{k}^{n}$ of the symmetric group $S_{n}$ (defined in Example
\ref{exa.rep.Sn-rep.nat}). Then: \medskip

\textbf{(a)} Define the \emph{dot product} $\left\langle u,v\right\rangle $ of
two vectors $u=\left(  u_{1},u_{2},\ldots,u_{n}\right)  \in\mathbf{k}^{n}$ and
$v=\left(  v_{1},v_{2},\ldots,v_{n}\right)  \in\mathbf{k}^{n}$ to be the
scalar $\sum_{i=1}^{n}u_{i}v_{i}$. (This can also be rewritten as $u^{T}v$, if
you identify the $1\times1$-matrix $u^{T}v$ with its only entry.) The map%
\begin{align*}
\mathbf{k}^{n}\times\mathbf{k}^{n}  &  \rightarrow\mathbf{k},\\
\left(  u,v\right)   &  \mapsto\left\langle u,v\right\rangle
\end{align*}
(which sends a pair of vectors to their dot product) is a bilinear form, known
as the \emph{dot product} or the \emph{standard scalar product} on
$\mathbf{k}^{n}$. This form is $S_{n}$-invariant, since each $g\in S_{n}$ and
any two vectors $u\in\mathbf{k}^{n}$ and $v\in\mathbf{k}^{n}$ satisfy the
equality $\left\langle gu,gv\right\rangle =\left\langle u,v\right\rangle $.
(Indeed, if we write $u$ and $v$ as $u=\left(  u_{1},u_{2},\ldots
,u_{n}\right)  $ and $v=\left(  v_{1},v_{2},\ldots,v_{n}\right)  $, then this
equality rewrites as $\sum_{i=1}^{n}u_{g^{-1}\left(  i\right)  }%
v_{g^{-1}\left(  i\right)  }=\sum_{i=1}^{n}u_{i}v_{i}$, which is true because
$g^{-1}$ is a permutation of $\left[  n\right]  $.) \medskip

\textbf{(b)} Another $S_{n}$-invariant bilinear form on the natural
representation $\mathbf{k}^{n}$ is the map%
\begin{align*}
\mathbf{k}^{n}\times\mathbf{k}^{n}  &  \rightarrow\mathbf{k},\\
\left(  u,v\right)   &  \mapsto\left(  \operatorname*{sum}u\right)  \left(
\operatorname*{sum}v\right)  ,
\end{align*}
where $\operatorname*{sum}w$ denotes the sum of all entries of any vector
$w\in\mathbf{k}^{n}$ (that is, $\operatorname*{sum}\left(  w_{1},w_{2}%
,\ldots,w_{n}\right)  :=w_{1}+w_{2}+\cdots+w_{n}$). This follows easily from
the fact that $\operatorname*{sum}\left(  gw\right)  =\operatorname*{sum}w$
for any $w\in\mathbf{k}^{n}$ and any $g\in S_{n}$. \medskip

\textbf{(c)} Of course, the zero form $\mathbf{k}^{n}\times\mathbf{k}%
^{n}\rightarrow\mathbf{k}$ (which sends every $\left(  u,v\right)  $ to $0$)
is also $S_{n}$-invariant. \medskip

\textbf{(d)} For any $n\geq1$, an example of a bilinear form that is not
$S_{n}$-invariant is%
\begin{align*}
\mathbf{k}^{n}\times\mathbf{k}^{n}  &  \rightarrow\mathbf{k},\\
\left(  u,v\right)   &  \mapsto u_{1}v_{1}\ \ \ \ \ \ \ \ \ \ \left(
\text{where }u=\left(  u_{1},u_{2},\ldots,u_{n}\right)  \text{ and }v=\left(
v_{1},v_{2},\ldots,v_{n}\right)  \right)
\end{align*}
(unless $n=1$, in which case this form is the dot product form and thus is
$S_{n}$-invariant).
\end{example}

\begin{exercise}
\textbf{(a)} \fbox{1} Let $G$ be a group. Let $U$ and $V$ be two trivial
$G$-representations. Prove that any bilinear form $f:U\times V\rightarrow
\mathbf{k}$ is $G$-invariant. \medskip

\textbf{(b)} \fbox{1} Consider the trivial representation $\mathbf{k}%
_{\operatorname*{triv}}$ and the sign representation $\mathbf{k}%
_{\operatorname*{sign}}$ of $S_{n}$ (see Example \ref{exa.rep.Sn-rep.triv} and
Example \ref{exa.rep.Sn-rep.sign}). Assume that $2$ is invertible in
$\mathbf{k}$, and that $n\geq2$. Prove that the only $S_{n}$-invariant
bilinear form $f:\mathbf{k}_{\operatorname*{triv}}\times\mathbf{k}%
_{\operatorname*{sign}}\rightarrow\mathbf{k}$ is the zero form (sending
everything to $0$).
\end{exercise}

\begin{exercise}
\fbox{1} Let $U$ and $V$ be two $\mathbf{k}$-modules. \medskip

\textbf{(a)} Prove that the set of all bilinear forms $f:U\times
V\rightarrow\mathbf{k}$ is a $\mathbf{k}$-module (with respect to pointwise
addition and scaling). \medskip

\textbf{(b)} Now assume that $U$ and $V$ are representations of a group $G$.
Prove that the set of all $G$-invariant bilinear forms $f:U\times
V\rightarrow\mathbf{k}$ is a $\mathbf{k}$-submodule of the $\mathbf{k}$-module
of all bilinear forms $f:U\times V\rightarrow\mathbf{k}$.
\end{exercise}

\begin{exercise}
\fbox{1} Let $G$ be a group. Define a map $S:\mathbf{k}\left[  G\right]
\rightarrow\mathbf{k}\left[  G\right]  $ as in Definition \ref{def.dual.S}.

Let $U$ and $V$ be two $G$-representations. Let $f:U\times V\rightarrow
\mathbf{k}$ be a bilinear form. Prove that $f$ is $G$-invariant if and only if
every $\mathbf{a}\in\mathbf{k}\left[  G\right]  $ and every $u\in U$ and $v\in
V$ satisfy%
\[
f\left(  \mathbf{a}\cdot u,\ v\right)  =f\left(  u,\ S\left(  \mathbf{a}%
\right)  \cdot v\right)  .
\]

\end{exercise}

As we said, $G$-invariant forms provide information about dual
representations. This is due to the following easy fact:

\begin{proposition}
\label{prop.dual.bilf.mor}Let $G$ be a group. Let $U$ and $V$ be two
$G$-representations. Let $f:U\times V\rightarrow\mathbf{k}$ be a $G$-invariant
bilinear form. Then: \medskip

\textbf{(a)} The maps $f_{L}:U\rightarrow V^{\ast}$ and $f_{R}:V\rightarrow
U^{\ast}$ (defined in Definition \ref{def.dual.bilinear} for $W=\mathbf{k}$)
are morphisms of $G$-representations. \medskip

\textbf{(b)} If $f_{L}$ is bijective, then $U\cong V^{\ast}$ as $G$-representations.
\end{proposition}

\begin{fineprint}
\begin{proof}
The form $f$ is $G$-invariant. In other words, every $g\in G$ and $u\in U$ and
$v\in V$ satisfy the equality%
\begin{equation}
f\left(  gu,\ gv\right)  =f\left(  u,\ v\right)  .
\label{pf.prop.dual.bilf.mor.inv}%
\end{equation}

\textbf{(a)} Let us prove the claim about $f_{L}$ only (the claim about
$f_{R}$ is analogous).

So we must show that $f_{L}$ is a morphism of $G$-representations. Since we
already know that $f_{L}$ is $\mathbf{k}$-linear, we must only show that
$f_{L}$ is $G$-equivariant. In other words, we must show that $f_{L}\left(
g\rightharpoonup u\right)  =g\overset{\ast}{\rightharpoonup}f_{L}\left(
u\right)  $ for all $g\in G$ and $u\in U$, where $\overset{\ast
}{\rightharpoonup}$ denotes the $G$-action on $V^{\ast}$ (as defined in
Definition \ref{def.dual.dual.G}).

So let us fix $g\in G$ and $u\in U$. Then, for any $v\in V$, we have%
\begin{align*}
\left(  f_{L}\left(  g\rightharpoonup u\right)  \right)  \left(  v\right)   &
=f\left(  \underbrace{g\rightharpoonup u}_{=gu},\ \underbrace{v}_{=gg^{-1}%
v}\right)  \ \ \ \ \ \ \ \ \ \ \left(  \text{by the definition of }%
f_{L}\right) \\
&  =f\left(  gu,\ gg^{-1}v\right) \\
&  =f\left(  u,\ g^{-1}v\right)  \ \ \ \ \ \ \ \ \ \ \left(  \text{by
(\ref{pf.prop.dual.bilf.mor.inv}), applied to }g^{-1}v\text{ instead of
}v\right) \\
&  =\left(  f_{L}\left(  u\right)  \right)  \left(  g^{-1}v\right)
\ \ \ \ \ \ \ \ \ \ \left(
\begin{array}
[c]{c}%
\text{since the definition of }f_{L}\\
\text{yields }\left(  f_{L}\left(  u\right)  \right)  \left(  g^{-1}v\right)
=f\left(  u,\ g^{-1}v\right)
\end{array}
\right)
\end{align*}
and%
\begin{align*}
\left(  g\overset{\ast}{\rightharpoonup}f_{L}\left(  u\right)  \right)
\left(  v\right)   &  =\left(  f_{L}\left(  u\right)  \right)  \left(
g^{-1}\rightharpoonup v\right)  \ \ \ \ \ \ \ \ \ \ \left(  \text{by
(\ref{eq.def.dual.dual.G.on-v}), applied to }f_{L}\left(  u\right)  \text{
instead of }f\right) \\
&  =\left(  f_{L}\left(  u\right)  \right)  \left(  g^{-1}v\right)
\ \ \ \ \ \ \ \ \ \ \left(  \text{since }g^{-1}\rightharpoonup v=g^{-1}%
v\right)  .
\end{align*}
Comparing these two equalities, we obtain%
\[
\left(  f_{L}\left(  g\rightharpoonup u\right)  \right)  \left(  v\right)
=\left(  g\overset{\ast}{\rightharpoonup}f_{L}\left(  u\right)  \right)
\left(  v\right)  \ \ \ \ \ \ \ \ \ \ \text{for each }v\in V.
\]
In other words, $f_{L}\left(  g\rightharpoonup u\right)  =g\overset{\ast
}{\rightharpoonup}f_{L}\left(  u\right)  $.

Forget that we fixed $g$ and $u$. We thus have shown that $f_{L}\left(
g\rightharpoonup u\right)  =g\overset{\ast}{\rightharpoonup}f_{L}\left(
u\right)  $ for all $g\in G$ and $u\in U$. In other words, the map $f_{L}$ is
$G$-equivariant. As we said, this completes our proof of Proposition
\ref{prop.dual.bilf.mor} \textbf{(a)}. \medskip

\textbf{(b)} Assume that $f_{L}$ is bijective. Thus, $f_{L}$ is invertible.
But Proposition \ref{prop.dual.bilf.mor} \textbf{(a)} shows that
$f_{L}:U\rightarrow V^{\ast}$ is a morphism of $G$-representations. Hence,
$f_{L}:U\rightarrow V^{\ast}$ is an invertible morphism of $G$%
-representations, thus an isomorphism (by Proposition
\ref{prop.rep.G-rep.iso=bij}). Hence, $U\cong V^{\ast}$. This proves
Proposition \ref{prop.dual.bilf.mor} \textbf{(b)}.
\end{proof}
\end{fineprint}

As an example of how Proposition \ref{prop.dual.bilf.mor} can be used, let us
show that any finite permutation module of any group $G$ is isomorphic to its
own dual:

\begin{theorem}
\label{thm.dual.perm-rep}Let $G$ be a group, and let $X$ be a finite left
$G$-set. Recall the permutation module $\mathbf{k}^{\left(  X\right)  }$
defined in Example \ref{exa.mod.kG-on-kX}; this is the left $\mathbf{k}\left[
G\right]  $-module on which $G$ acts by permuting the standard basis vectors
(see the equality (\ref{eq.exa.mod.kG-on-kX.2}) for details).

Then, $\left(  \mathbf{k}^{\left(  X\right)  }\right)  ^{\ast}\cong%
\mathbf{k}^{\left(  X\right)  }$ as $G$-representations.
\end{theorem}

\begin{proof}
As in Example \ref{exa.mod.kG-on-kX}, we denote the standard basis vector
$e_{x}$ of $\mathbf{k}^{\left(  X\right)  }$ as $x$ for each $x\in X$. Then,
$\left(  x\right)  _{x\in X}$ is the standard basis of $\mathbf{k}^{\left(
X\right)  }$. Thus, by Proposition \ref{prop.dual.bilf.on-basis}, we can
define a bilinear form $f:\mathbf{k}^{\left(  X\right)  }\times\mathbf{k}%
^{\left(  X\right)  }\rightarrow\mathbf{k}$ by setting%
\[
f\left(  x,y\right)  =\delta_{x,y}\ \ \ \ \ \ \ \ \ \ \text{for all }x\in
X\text{ and }y\in X
\]
(where $\delta_{x,y}$ denotes the Kronecker delta). Consider this bilinear
form $f$.

We shall now show that the bilinear form $f$ is $G$-invariant. Indeed, by
Definition \ref{def.dual.bilf.G-inv}, this amounts to showing that every $g\in
G$ and $u\in\mathbf{k}^{\left(  X\right)  }$ and $v\in\mathbf{k}^{\left(
X\right)  }$ satisfy the equality $f\left(  gu,\ gv\right)  =f\left(
u,\ v\right)  $. So let us show this. Fix $g\in G$ and $u\in\mathbf{k}%
^{\left(  X\right)  }$ and $v\in\mathbf{k}^{\left(  X\right)  }$. We must
prove the equality $f\left(  gu,\ gv\right)  =f\left(  u,\ v\right)  $. Both
sides of this equality depend $\mathbf{k}$-linearly on each of $u$ and $v$
(since $f$ is $\mathbf{k}$-bilinear). Hence, in proving it, we can WLOG assume
that both $u$ and $v$ belong to the standard basis $\left(  x\right)  _{x\in
X}$ of $\mathbf{k}^{\left(  X\right)  }$. Assume this. Thus, $u=x$ and $v=y$
for some $x,y\in X$. Consider these $x,y$.

Now, the definition of $f$ yields%
\begin{align*}
f\left(  gx,\ gy\right)   &  =\delta_{gx,gy}\\
&  =%
\begin{cases}
1, & \text{if }gx=gy;\\
0, & \text{if }gx\neq gy
\end{cases}
\ \ \ \ \ \ \ \ \ \ \left(  \text{by the definition of the Kronecker
delta}\right) \\
&  =%
\begin{cases}
1, & \text{if }x=y;\\
0, & \text{if }x\neq y
\end{cases}
\ \ \ \ \ \ \ \ \ \ \left(
\begin{array}
[c]{c}%
\text{since }gx=gy\text{ if and only if }x=y\\
\text{(because }g\text{ is invertible)}%
\end{array}
\right) \\
&  =\delta_{x,y}\ \ \ \ \ \ \ \ \ \ \left(  \text{by the definition of the
Kronecker delta}\right) \\
&  =f\left(  x,\ y\right)  \ \ \ \ \ \ \ \ \ \ \left(  \text{by the definition
of }f\right)  .
\end{align*}
In other words, $f\left(  gu,\ gv\right)  =f\left(  u,\ v\right)  $ (since
$u=x$ and $v=y$).

Forget that we fixed $u$ and $v$. We thus have shown that every $g\in G$ and
$u\in\mathbf{k}^{\left(  X\right)  }$ and $v\in\mathbf{k}^{\left(  X\right)
}$ satisfy the equality $f\left(  gu,\ gv\right)  =f\left(  u,\ v\right)  $.
In other words, the bilinear form $f$ is $G$-invariant.

Now, recall that $\mathbf{k}^{\left(  X\right)  }$ is a free $\mathbf{k}%
$-module with a finite basis $\left(  x\right)  _{x\in X}$ (it is finite since
$X$ is finite). Let $\left(  x^{\ast}\right)  _{x\in X}$ be the dual basis of
this basis $\left(  x\right)  _{x\in X}$. Then, $\left(  x^{\ast}\right)
_{x\in X}$ is a basis of the dual $\mathbf{k}$-module $\left(  \mathbf{k}%
^{\left(  X\right)  }\right)  ^{\ast}$ (since the set $X$ is finite).

However, the map $f_{L}:\mathbf{k}^{\left(  X\right)  }\rightarrow\left(
\mathbf{k}^{\left(  X\right)  }\right)  ^{\ast}$ (defined as in Definition
\ref{def.dual.bilinear}) is easily seen to satisfy%
\[
f_{L}\left(  x\right)  =x^{\ast}\ \ \ \ \ \ \ \ \ \ \text{for each }x\in X
\]
\footnote{\textit{Proof.} Let $x\in X$. Then, for each $y\in X$, we have%
\begin{align*}
\left(  f_{L}\left(  x\right)  \right)  \left(  y\right)   &  =f\left(
x,y\right)  \ \ \ \ \ \ \ \ \ \ \left(  \text{by the definition of }%
f_{L}\right) \\
&  =\delta_{x,y}\ \ \ \ \ \ \ \ \ \ \left(  \text{by the definition of
}f\right) \\
&  =x^{\ast}\left(  y\right)
\end{align*}
(since the definition of a dual basis yields $x^{\ast}\left(  y\right)
=\delta_{x,y}$). In other words, the two maps $f_{L}\left(  x\right)  $ and
$x^{\ast}$ (from $\mathbf{k}^{\left(  X\right)  }$ to $\mathbf{k}$) agree on
each basis vector $y$ of the standard basis of $\mathbf{k}^{\left(  X\right)
}$. Thus, these two maps must be identical (because they are $\mathbf{k}%
$-linear, but any two $\mathbf{k}$-linear maps that agree on a basis of their
domain must be identical). In other words, $f_{L}\left(  x\right)  =x^{\ast}$,
qed.}. In other words, the map $f_{L}$ sends the entries $x$ of the basis
$\left(  x\right)  _{x\in X}$ of $\mathbf{k}^{\left(  X\right)  }$ to the
corresponding entries $x^{\ast}$ of the basis $\left(  x^{\ast}\right)  _{x\in
X}$ of $\left(  \mathbf{k}^{\left(  X\right)  }\right)  ^{\ast}$. Hence, this
map $f_{L}$ is a $\mathbf{k}$-module isomorphism (since it is $\mathbf{k}%
$-linear and sends a basis to a basis), thus bijective. Hence, Proposition
\ref{prop.dual.bilf.mor} \textbf{(b)} (applied to $U=\mathbf{k}^{\left(
X\right)  }$ and $V=\mathbf{k}^{\left(  X\right)  }$) shows that
$\mathbf{k}^{\left(  X\right)  }\cong\left(  \mathbf{k}^{\left(  X\right)
}\right)  ^{\ast}$ as $G$-representations. Thus, Theorem
\ref{thm.dual.perm-rep} is proved.
\end{proof}

Two particular cases of Theorem \ref{thm.dual.perm-rep} are worth mentioning.
The first has already been stated in Exercise \ref{exe.dual.nat.1}, but we
shall now prove it anew by interpreting it as a special case of Theorem
\ref{thm.dual.perm-rep}:

\begin{corollary}
\label{cor.dual.nat}Consider the natural representation $\mathbf{k}^{n}$ of
the symmetric group $S_{n}$ (defined in Example \ref{exa.rep.Sn-rep.nat}).
Then, $\left(  \mathbf{k}^{n}\right)  ^{\ast}\cong\mathbf{k}^{n}$ as $S_{n}$-representations.
\end{corollary}

\begin{proof}
[Proof sketch.]I claim that the natural representation $\mathbf{k}^{n}$ of
$S_{n}$ is a permutation module. Indeed, consider the natural action of
$S_{n}$ on the set $\left[  n\right]  =\left\{  1,2,\ldots,n\right\}  $ given
by the rule%
\[
g\rightharpoonup i:=g\left(  i\right)  \ \ \ \ \ \ \ \ \ \ \text{for all }g\in
S_{n}\text{ and }i\in\left[  n\right]
\]
(this is the action from Example \ref{exa.rep.G-sets.nat}, applied to
$A=\left[  n\right]  $). Thus, $\left[  n\right]  $ is a left $S_{n}$-set. The
corresponding permutation module $\mathbf{k}^{\left(  \left[  n\right]
\right)  }$ (as defined in Example \ref{exa.mod.kG-on-kX}) is the free
$\mathbf{k}$-module with basis $\left(  e_{i}\right)  _{i\in\left[  n\right]
}=\left(  e_{1},e_{2},\ldots,e_{n}\right)  $ and with a linear $S_{n}$-action
given by
\[
g\rightharpoonup e_{i}=e_{g\rightharpoonup i}=e_{g\left(  i\right)
}\ \ \ \ \ \ \ \ \ \ \text{for all }g\in S_{n}\text{ and }i\in\left[
n\right]  .
\]
But this is precisely the natural representation $\mathbf{k}^{n}$ of $S_{n}$
with its standard basis.

Thus, this natural representation $\mathbf{k}^{n}$ is the permutation module
$\mathbf{k}^{\left(  \left[  n\right]  \right)  }$ corresponding to the left
$S_{n}$-set $\left[  n\right]  $. But Theorem \ref{thm.dual.perm-rep} (applied
to $G=S_{n}$ and $X=\left[  n\right]  $) shows that $\left(  \mathbf{k}%
^{\left(  \left[  n\right]  \right)  }\right)  ^{\ast}\cong\mathbf{k}^{\left(
\left[  n\right]  \right)  }$ as $S_{n}$-representations. Since the natural
representation $\mathbf{k}^{n}$ is the permutation module $\mathbf{k}^{\left(
\left[  n\right]  \right)  }$, we can rewrite this as $\left(  \mathbf{k}%
^{n}\right)  ^{\ast}\cong\mathbf{k}^{n}$ (as $S_{n}$-representations). This
proves Corollary \ref{cor.dual.nat}.
\end{proof}

The next particular case of Theorem \ref{thm.dual.perm-rep} is the case of
Young modules:

\begin{corollary}
\label{cor.dual.youngmod}Let $D$ be a diagram with $\left\vert D\right\vert
=n$. Then, the Young module $\mathcal{M}^{D}$ satisfies $\left(
\mathcal{M}^{D}\right)  ^{\ast}\cong\mathcal{M}^{D}$ as $S_{n}$-representations.
\end{corollary}

\begin{proof}
The Young module $\mathcal{M}^{D}$ is a permutation module. Indeed, recall
(from Definition \ref{def.tabloid.Sn-act}) that the set $\left\{
n\text{-tabloids of shape }D\right\}  $ is a left $S_{n}$-set. The permutation
module corresponding to this $S_{n}$-set (as defined in Example
\ref{exa.mod.kG-on-kX}) is precisely $\mathcal{M}^{D}$ (by Definition
\ref{def.youngmod.youngmod}). In other words,%
\begin{equation}
\mathcal{M}^{D}=\mathbf{k}^{\left(  \left\{  n\text{-tabloids of shape
}D\right\}  \right)  } \label{pf.cor.dual.youngmod.M=K}%
\end{equation}
as $S_{n}$-representations. But Theorem \ref{thm.dual.perm-rep} (applied to
$G=S_{n}$ and $X=\left\{  n\text{-tabloids of shape }D\right\}  $) shows that%
\[
\left(  \mathbf{k}^{\left(  \left\{  n\text{-tabloids of shape }D\right\}
\right)  }\right)  ^{\ast}\cong\mathbf{k}^{\left(  \left\{  n\text{-tabloids
of shape }D\right\}  \right)  }\text{ as }S_{n}\text{-representations.}%
\]
In view of (\ref{pf.cor.dual.youngmod.M=K}), we can rewrite this as follows:
\[
\left(  \mathcal{M}^{D}\right)  ^{\ast}\cong\mathcal{M}^{D}\text{ as }%
S_{n}\text{-representations.}%
\]
This proves Corollary \ref{cor.dual.youngmod}.
\end{proof}

\subsubsection{The dual of a Specht module: The case of a field}

So the dual of a Young module is the Young module itself (up to a rather
obvious isomorphism). We shall now study the dual of a Specht module, which is
far less obvious. We begin with the case when $\mathbf{k}$ is a field.

Recall the map $\mathbf{r}:\mathbb{Z}^{2}\rightarrow\mathbb{Z}^{2}$ (defined
in Theorem \ref{thm.partitions.conj}), which sends each cell $\left(
i,j\right)  \in\mathbb{Z}^{2}$ to $\left(  j,i\right)  $. Recall also the
sign-twist $V^{\operatorname*{sign}}$ of an $S_{n}$-representation $V$,
defined in Definition \ref{def.sign-twist.sign-twist.Sn} or in Definition
\ref{def.sign-twist.sign-twist.kSn}. Now we claim (generalizing \cite[Theorems
6.7 and 8.15]{James78}):

\begin{theorem}
\label{thm.dual.specht-field}Let $D$ be a diagram with $\left\vert
D\right\vert =n$. Assume that $\mathbf{k}$ is a field. Then, the Specht module
$\mathcal{S}^{D}$ satisfies
\[
\left(  \mathcal{S}^{D}\right)  ^{\ast}\cong\left(  \mathcal{S}^{\mathbf{r}%
\left(  D\right)  }\right)  ^{\operatorname*{sign}}\text{ as }S_{n}%
\text{-representations.}%
\]

\end{theorem}

I don't know if the \textquotedblleft$\mathbf{k}$ is a field\textquotedblright%
\ assumption is actually necessary here; the proof below makes use of it, but
I don't know of any counterexample when $\mathbf{k}$ is not a field. I will
say more about this below, after proving Theorem \ref{thm.dual.specht-field}
and Theorem \ref{thm.dual.specht-skew}.

\begin{proof}
[Proof of Theorem \ref{thm.dual.specht-field}.]Let $\mathcal{A}=\mathbf{k}%
\left[  S_{n}\right]  $. Recall the antipode $S$ of $\mathbf{k}\left[
S_{n}\right]  $.

Furthermore, recall Definition \ref{def.specht.ET.wa}. Thus, if an element
$\mathbf{a}\in\mathcal{A}$ is written as $\mathbf{a}=\sum\limits_{g\in S_{n}%
}\alpha_{g}g$ for some scalars $\alpha_{g}\in\mathbf{k}$, then $\left[
w\right]  \mathbf{a}:=\alpha_{w}$ for any $w\in S_{n}$. We begin by showing a
few properties of coefficients of elements of $\mathcal{A}$:

\begin{statement}
\textit{Claim 1:} Let $\mathbf{a}\in\mathbf{k}\left[  S_{n}\right]  $ and
$g\in S_{n}$. Then,%
\begin{align}
\left[  1\right]  \left(  g\mathbf{a}g^{-1}\right)   &  =\left[  1\right]
\mathbf{a};\label{pf.thm.dual.specht-field.c1.1}\\
\left[  1\right]  \left(  \mathbf{a}g^{-1}\right)   &  =\left[  g\right]
\mathbf{a};\label{pf.thm.dual.specht-field.c1.2}\\
\left[  g\right]  \mathbf{a}  &  =\left[  1\right]  \left(  gS\left(
\mathbf{a}\right)  \right)  . \label{pf.thm.dual.specht-field.c1.3}%
\end{align}

\end{statement}

\begin{proof}
[Proof of Claim 1.]All three equalities (\ref{pf.thm.dual.specht-field.c1.1}),
(\ref{pf.thm.dual.specht-field.c1.2}) and (\ref{pf.thm.dual.specht-field.c1.3}%
) depend $\mathbf{k}$-linearly on $\mathbf{a}$ (since $S$ is a $\mathbf{k}%
$-linear map). Thus, by linearity, we can WLOG assume that $\mathbf{a}$ is a
standard basis vector of $\mathbf{k}\left[  S_{n}\right]  $. Assume this.
Thus, $\mathbf{a}=w$ for some $w\in S_{n}$. Consider this $w$. Hence,
$S\left(  \mathbf{a}\right)  =S\left(  w\right)  =w^{-1}$ (by the definition
of the antipode $S$).

Definition \ref{def.specht.ET.wa} yields that any two elements $u,v\in S_{n}$
satisfy $\left[  u\right]  v=\delta_{u,v}$ (Kronecker delta). Hence, $\left[
g\right]  w=\delta_{g,w}$ and $\left[  1\right]  \left(  gw^{-1}\right)
=\delta_{1,gw^{-1}}=\delta_{g,w}$ (since the equation $1=gw^{-1}$ is
equivalent to $g=w$). Comparing these two equalities, we obtain $\left[
g\right]  w=\left[  1\right]  \left(  gw^{-1}\right)  $. In other words,
$\left[  g\right]  \mathbf{a}=\left[  1\right]  \left(  gS\left(
\mathbf{a}\right)  \right)  $ (since $\mathbf{a}=w$ and $S\left(
\mathbf{a}\right)  =w^{-1}$). Thus, (\ref{pf.thm.dual.specht-field.c1.3}) is proved.

Recall again the fact that any two elements $u,v\in S_{n}$ satisfy $\left[
u\right]  v=\delta_{u,v}$. Hence, $\left[  1\right]  \left(  wg^{-1}\right)
=\delta_{1,wg^{-1}}=\delta_{g,w}$ (since the equation $1=wg^{-1}$ holds if and
only if $g=w$) and $\left[  g\right]  w=\delta_{g,w}$. Comparing these two
equalities, we obtain $\left[  1\right]  \left(  wg^{-1}\right)  =\left[
g\right]  w$. In other words, $\left[  1\right]  \left(  \mathbf{a}%
g^{-1}\right)  =\left[  g\right]  \mathbf{a}$ (since $\mathbf{a}=w$). Thus,
(\ref{pf.thm.dual.specht-field.c1.2}) is proved.

It remains to prove (\ref{pf.thm.dual.specht-field.c1.1}). For this, we
observe that the equation $1=gwg^{-1}$ is equivalent to $1=w$ (because
multiplying the equation $1=gwg^{-1}$ by $g$ on the right transforms it into
$g=gw$, which is clearly equivalent to $1=w$). Now, recall once again the fact
that any two elements $u,v\in S_{n}$ satisfy $\left[  u\right]  v=\delta
_{u,v}$. Hence, $\left[  1\right]  \left(  gwg^{-1}\right)  =\delta
_{1,gwg^{-1}}=\delta_{1,w}$ (since the equation $1=gwg^{-1}$ is equivalent to
$1=w$) and $\left[  1\right]  w=\delta_{1,w}$. Comparing these two equalities,
we obtain $\left[  1\right]  \left(  gwg^{-1}\right)  =\left[  1\right]  w$.
In other words, $\left[  1\right]  \left(  g\mathbf{a}g^{-1}\right)  =\left[
1\right]  \mathbf{a}$ (since $\mathbf{a}=w$). Thus,
(\ref{pf.thm.dual.specht-field.c1.1}) is proved.

Now, we have proved all three equalities (\ref{pf.thm.dual.specht-field.c1.1}%
), (\ref{pf.thm.dual.specht-field.c1.2}) and
(\ref{pf.thm.dual.specht-field.c1.3}). Thus, Claim 1 is proved.
\end{proof}

Fix an $n$-tableau $T$ of shape $D$ (such an $n$-tableau clearly exists, since
$\left\vert D\right\vert =n$). Consider the elements
\[
\mathbf{E}_{T}:=\nabla_{\operatorname*{Col}T}^{-}\nabla_{\operatorname*{Row}%
T}\ \ \ \ \ \ \ \ \ \ \text{and}\ \ \ \ \ \ \ \ \ \ \mathbf{F}_{T}%
:=\nabla_{\operatorname*{Row}T}\nabla_{\operatorname*{Col}T}^{-}%
\]
of $\mathcal{A}$. Proposition \ref{prop.specht.FT.basics} \textbf{(a)} (or,
rather, its obvious generalization to tableaux $T$ of arbitrary shape) says
that the antipode $S$ of $\mathbf{k}\left[  S_{n}\right]  $ satisfies
$S\left(  \mathbf{E}_{T}\right)  =\mathbf{F}_{T}$. Applying the map $S$ to
this equality, we obtain $S\left(  S\left(  \mathbf{E}_{T}\right)  \right)
=S\left(  \mathbf{F}_{T}\right)  $. But $S\left(  S\left(  \mathbf{E}%
_{T}\right)  \right)  =\mathbf{E}_{T}$ (since Definition \ref{def.dual.S}
\textbf{(b)} says that $S\circ S=\operatorname*{id}$). Comparing these two
equalities, we find $S\left(  \mathbf{F}_{T}\right)  =\mathbf{E}_{T}$.
Moreover, recall (from Definition \ref{def.dual.S} \textbf{(a)}) that $S$ is
an algebra anti-morphism. Lemma \ref{lem.sign-twist.AFT} yields that%
\begin{equation}
\mathcal{S}^{D}\cong\mathcal{A}\mathbf{E}_{T}\ \ \ \ \ \ \ \ \ \ \text{and}%
\ \ \ \ \ \ \ \ \ \ \left(  \mathcal{S}^{\mathbf{r}\left(  D\right)  }\right)
^{\operatorname*{sign}}\cong\mathcal{A}\mathbf{F}_{T}
\label{pf.thm.dual.specht-field.avas}%
\end{equation}
as $S_{n}$-representations. From $\mathcal{S}^{D}\cong\mathcal{A}%
\mathbf{E}_{T}$, we obtain $\left(  \mathcal{S}^{D}\right)  ^{\ast}%
\cong\left(  \mathcal{A}\mathbf{E}_{T}\right)  ^{\ast}$ (by Corollary
\ref{cor.dual.iso}).

We must prove that $\left(  \mathcal{S}^{D}\right)  ^{\ast}\cong\left(
\mathcal{S}^{\mathbf{r}\left(  D\right)  }\right)  ^{\operatorname*{sign}}$ as
$S_{n}$-representations. In view of (\ref{pf.thm.dual.specht-field.avas}), it
suffices to show that $\left(  \mathcal{A}\mathbf{E}_{T}\right)  ^{\ast}%
\cong\mathcal{A}\mathbf{F}_{T}$ as $S_{n}$-representations (because once this
is proved, then we will immediately obtain $\left(  \mathcal{S}^{D}\right)
^{\ast}\cong\left(  \mathcal{A}\mathbf{E}_{T}\right)  ^{\ast}\cong%
\mathcal{A}\mathbf{F}_{T}\cong\left(  \mathcal{S}^{\mathbf{r}\left(  D\right)
}\right)  ^{\operatorname*{sign}}$ by (\ref{pf.thm.dual.specht-field.avas})).
In order to show this, we will construct an $S_{n}$-invariant bilinear form
$f:\mathcal{A}\mathbf{E}_{T}\times\mathcal{A}\mathbf{F}_{T}\rightarrow
\mathbf{k}$ such that the map $f_{R}:\mathcal{A}\mathbf{F}_{T}\rightarrow
\left(  \mathcal{A}\mathbf{E}_{T}\right)  ^{\ast}$ is an isomorphism of
$S_{n}$-representations.

We define this bilinear form $f:\mathcal{A}\mathbf{E}_{T}\times\mathcal{A}%
\mathbf{F}_{T}\rightarrow\mathbf{k}$ by the equality%
\begin{equation}
f\left(  u\mathbf{E}_{T},\ v\mathbf{F}_{T}\right)  :=\left[  1\right]  \left(
u\mathbf{E}_{T}S\left(  v\right)  \right)  \ \ \ \ \ \ \ \ \ \ \text{for all
}u,v\in\mathcal{A}. \label{pf.thm.dual.specht-field.f=}%
\end{equation}
But first, we must show that this form $f$ is well-defined, i.e., that the
equality (\ref{pf.thm.dual.specht-field.f=}) really uniquely determines each
value of $f$ and that the resulting map $f:\mathcal{A}\mathbf{E}_{T}%
\times\mathcal{A}\mathbf{F}_{T}\rightarrow\mathbf{k}$ is a bilinear form.

In other words, we must prove the following:

\begin{statement}
\textit{Claim 2:} The bilinear form $f:\mathcal{A}\mathbf{E}_{T}%
\times\mathcal{A}\mathbf{F}_{T}\rightarrow\mathbf{k}$ is well-defined.
\end{statement}

\begin{proof}
[Proof of Claim 2.]First, we shall show that the equality
(\ref{pf.thm.dual.specht-field.f=}) defines a valid \textbf{map}
$f:\mathcal{A}\mathbf{E}_{T}\times\mathcal{A}\mathbf{F}_{T}\rightarrow
\mathbf{k}$. (We will prove its bilinearity later.)

Each element of $\mathcal{A}\mathbf{E}_{T}$ can be written as $u\mathbf{E}%
_{T}$ for some $u\in\mathcal{A}$. Each element of $\mathcal{A}\mathbf{F}_{T}$
can be written as $v\mathbf{F}_{T}$ for some $v\in\mathcal{A}$. Thus, if we
want to compute $f\left(  a,b\right)  $ for some $a\in\mathcal{A}%
\mathbf{E}_{T}$ and some $b\in\mathcal{A}\mathbf{F}_{T}$, then we can write
$a$ as $u\mathbf{E}_{T}$ and write $b$ as $v\mathbf{F}_{T}$, and then compute
the value $f\left(  a,b\right)  $ as follows:%
\[
f\left(  \underbrace{a}_{=u\mathbf{E}_{T}},\ \underbrace{b}_{=v\mathbf{F}_{T}%
}\right)  =f\left(  u\mathbf{E}_{T},\ v\mathbf{F}_{T}\right)  =\left[
1\right]  \left(  u\mathbf{E}_{T}S\left(  v\right)  \right)
\ \ \ \ \ \ \ \ \ \ \left(  \text{by (\ref{pf.thm.dual.specht-field.f=}%
)}\right)  .
\]
Thus, the equality (\ref{pf.thm.dual.specht-field.f=}) allows us to compute
each value $f\left(  a,b\right)  $ of $f$. However, this does not yet mean
that $f$ is well-defined: Indeed, when we choose the elements $u$ and $v$
satisfying $a=u\mathbf{E}_{T}$ and $b=v\mathbf{F}_{T}$, we have some freedom
(they are not uniquely determined by $a$ and $b$, since (e.g.) one single $a$
can be written as $u\mathbf{E}_{T}$ for several values of $u$), and it is
conceivable that the resulting value $\left[  1\right]  \left(  u\mathbf{E}%
_{T}S\left(  v\right)  \right)  $ could depend on these choices; but then we
would get different values for the same $f\left(  a,b\right)  $ depending on
our choices of $u$ and $v$, and this would render $f$ undefined.

So we need to show that this cannot happen. In other words, we need to show
that the element $\left[  1\right]  \left(  u\mathbf{E}_{T}S\left(  v\right)
\right)  $ on the right hand side of (\ref{pf.thm.dual.specht-field.f=})
depends only on $u\mathbf{E}_{T}$ and $v\mathbf{F}_{T}$ but not on the
specific values of $u$ and $v$. In other words, we need to show that if
$u,v\in\mathcal{A}$ and $u^{\prime},v^{\prime}\in\mathcal{A}$ satisfy
$u\mathbf{E}_{T}=u^{\prime}\mathbf{E}_{T}$ and $v\mathbf{F}_{T}=v^{\prime
}\mathbf{F}_{T}$, then $\left[  1\right]  \left(  u\mathbf{E}_{T}S\left(
v\right)  \right)  =\left[  1\right]  \left(  u^{\prime}\mathbf{E}_{T}S\left(
v^{\prime}\right)  \right)  $.

Let us show this. Indeed, let $u,v\in\mathcal{A}$ and $u^{\prime},v^{\prime
}\in\mathcal{A}$ be such that $u\mathbf{E}_{T}=u^{\prime}\mathbf{E}_{T}$ and
$v\mathbf{F}_{T}=v^{\prime}\mathbf{F}_{T}$. Recall that $S$ is an algebra
anti-morphism. Thus, $S\left(  v\mathbf{F}_{T}\right)  =\underbrace{S\left(
\mathbf{F}_{T}\right)  }_{=\mathbf{E}_{T}}S\left(  v\right)  =\mathbf{E}%
_{T}S\left(  v\right)  $. Likewise, $S\left(  v^{\prime}\mathbf{F}_{T}\right)
=\mathbf{E}_{T}S\left(  v^{\prime}\right)  $. From $S\left(  v\mathbf{F}%
_{T}\right)  =\mathbf{E}_{T}S\left(  v\right)  $, we obtain%
\[
\mathbf{E}_{T}S\left(  v\right)  =S\left(  \underbrace{v\mathbf{F}_{T}%
}_{=v^{\prime}\mathbf{F}_{T}}\right)  =S\left(  v^{\prime}\mathbf{F}%
_{T}\right)  =\mathbf{E}_{T}S\left(  v^{\prime}\right)  .
\]

But now,%
\[
\left[  1\right]  \left(  \underbrace{u\mathbf{E}_{T}}_{=u^{\prime}%
\mathbf{E}_{T}}S\left(  v\right)  \right)  =\left[  1\right]  \left(
u^{\prime}\underbrace{\mathbf{E}_{T}S\left(  v\right)  }_{=\mathbf{E}%
_{T}S\left(  v^{\prime}\right)  }\right)  =\left[  1\right]  \left(
u^{\prime}\mathbf{E}_{T}S\left(  v^{\prime}\right)  \right)  .
\]

Forget that we fixed $u,v$ and $u^{\prime},v^{\prime}$. Thus, we have shown
that if $u,v\in\mathcal{A}$ and $u^{\prime},v^{\prime}\in\mathcal{A}$ satisfy
$u\mathbf{E}_{T}=u^{\prime}\mathbf{E}_{T}$ and $v\mathbf{F}_{T}=v^{\prime
}\mathbf{F}_{T}$, then $\left[  1\right]  \left(  u\mathbf{E}_{T}S\left(
v\right)  \right)  =\left[  1\right]  \left(  u^{\prime}\mathbf{E}_{T}S\left(
v^{\prime}\right)  \right)  $. In other words, the right hand side of
(\ref{pf.thm.dual.specht-field.f=}) depends only on $u\mathbf{E}_{T}$ and on
$v\mathbf{F}_{T}$ but not on the specific choices of $u$ and $v$. Therefore,
the map $f:\mathcal{A}\mathbf{E}_{T}\times\mathcal{A}\mathbf{F}_{T}%
\rightarrow\mathbf{k}$ is well-defined (in the sense that its output value
does not depend on how exactly we write its two inputs as $u\mathbf{E}_{T}$
and $v\mathbf{F}_{T}$).

It remains to prove that this map is a bilinear form. In other words, it
remains to prove that it is $\mathbf{k}$-bilinear. But this is
straightforward\footnote{Here is the proof in some detail:
\par
We show that the map $f$ is $\mathbf{k}$-bilinear by verifying that it
satisfies all six axioms from Remark \ref{rmk.dual.bilin.axioms}. Let us only
verify the first of these six axioms (the other five are analogous).
\par
Thus, we must show that $f\left(  a,\ b_{1}+b_{2}\right)  =f\left(
a,b_{1}\right)  +f\left(  a,b_{2}\right)  $ for all $a\in\mathcal{A}%
\mathbf{E}_{T}$ and $b_{1},b_{2}\in\mathcal{A}\mathbf{F}_{T}$.
\par
Indeed, let $a\in\mathcal{A}\mathbf{E}_{T}$ and $b_{1},b_{2}\in\mathcal{A}%
\mathbf{F}_{T}$. Then, we can write $a$ as $a=u\mathbf{E}_{T}$ for some
$u\in\mathcal{A}$ (since $a\in\mathcal{A}\mathbf{E}_{T}$). Consider this $u$.
Furthermore, we can write $b_{1}$ as $b_{1}=v_{1}\mathbf{F}_{T}$ for some
$v_{1}\in\mathcal{A}$ (since $b_{1}\in\mathcal{A}\mathbf{F}_{T}$), and we can
write $b_{2}$ as $b_{2}=v_{2}\mathbf{F}_{T}$ for some $v_{2}\in\mathcal{A}$
(since $b_{2}\in\mathcal{A}\mathbf{F}_{T}$). Consider these $v_{1}$ and
$v_{2}$.
\par
From $a=u\mathbf{E}_{T}$ and $b_{1}=v_{1}\mathbf{F}_{T}$, we obtain%
\[
f\left(  a,b_{1}\right)  =f\left(  u\mathbf{E}_{T},\ v_{1}\mathbf{F}%
_{T}\right)  =\left[  1\right]  \left(  u\mathbf{E}_{T}S\left(  v_{1}\right)
\right)  \ \ \ \ \ \ \ \ \ \ \left(  \text{by
(\ref{pf.thm.dual.specht-field.f=}), applied to }v=v_{1}\right)  .
\]
In other words, $\left[  1\right]  \left(  u\mathbf{E}_{T}S\left(
v_{1}\right)  \right)  =f\left(  a,b_{1}\right)  $. Likewise, $\left[
1\right]  \left(  u\mathbf{E}_{T}S\left(  v_{2}\right)  \right)  =f\left(
a,b_{2}\right)  $. Furthermore, adding the equalities $b_{1}=v_{1}%
\mathbf{F}_{T}$ and $b_{2}=v_{2}\mathbf{F}_{T}$ together, we obtain
$b_{1}+b_{2}=v_{1}\mathbf{F}_{T}+v_{2}\mathbf{F}_{T}=\left(  v_{1}%
+v_{2}\right)  \mathbf{F}_{T}$. From this and from $a=u\mathbf{E}_{T}$, we
obtain%
\begin{align*}
f\left(  a,\ b_{1}+b_{2}\right)   &  =f\left(  u\mathbf{E}_{T},\ \left(
v_{1}+v_{2}\right)  \mathbf{F}_{T}\right) \\
&  =\left[  1\right]  \left(  u\mathbf{E}_{T}\underbrace{S\left(  v_{1}%
+v_{2}\right)  }_{\substack{=S\left(  v_{1}\right)  +S\left(  v_{2}\right)
\\\text{(since the map }S\text{ is }\mathbf{k}\text{-linear)}}}\right)
\ \ \ \ \ \ \ \ \ \ \left(  \text{by (\ref{pf.thm.dual.specht-field.f=}),
applied to }v=v_{1}+v_{2}\right) \\
&  =\left[  1\right]  \underbrace{\left(  u\mathbf{E}_{T}\left(  S\left(
v_{1}\right)  +S\left(  v_{2}\right)  \right)  \right)  }_{=u\mathbf{E}%
_{T}S\left(  v_{1}\right)  +u\mathbf{E}_{T}S\left(  v_{2}\right)  }=\left[
1\right]  \left(  u\mathbf{E}_{T}S\left(  v_{1}\right)  +u\mathbf{E}%
_{T}S\left(  v_{2}\right)  \right) \\
&  =\underbrace{\left[  1\right]  \left(  u\mathbf{E}_{T}S\left(
v_{1}\right)  \right)  }_{=f\left(  a,b_{1}\right)  }+\underbrace{\left[
1\right]  \left(  u\mathbf{E}_{T}S\left(  v_{2}\right)  \right)  }_{=f\left(
a,b_{2}\right)  }=f\left(  a,b_{1}\right)  +f\left(  a,b_{2}\right)  .
\end{align*}
\par
Forget that we fixed $a,b_{1},b_{2}$. We thus have shown that $f\left(
a,\ b_{1}+b_{2}\right)  =f\left(  a,b_{1}\right)  +f\left(  a,b_{2}\right)  $
for all $a\in\mathcal{A}\mathbf{E}_{T}$ and $b_{1},b_{2}\in\mathcal{A}%
\mathbf{F}_{T}$. This is the first of the six axioms that must hold in order
for $f$ to be $\mathbf{k}$-bilinear. The remaining five axioms are just as
easy to check using the same kind of reasoning. Thus, we have shown that $f$
is $\mathbf{k}$-bilinear.}. Thus, Claim 2 is proved.
\end{proof}

\begin{statement}
\textit{Claim 3:} The bilinear form $f:\mathcal{A}\mathbf{E}_{T}%
\times\mathcal{A}\mathbf{F}_{T}\rightarrow\mathbf{k}$ is $S_{n}$-invariant.
\end{statement}

\begin{proof}
[Proof of Claim 3.]We must show that every $g\in S_{n}$ and $a\in
\mathcal{A}\mathbf{E}_{T}$ and $b\in\mathcal{A}\mathbf{F}_{T}$ satisfy the
equality $f\left(  ga,\ gb\right)  =f\left(  a,\ b\right)  $ (since this is
what it means for $f$ to be $S_{n}$-invariant, according to Definition
\ref{def.dual.bilf.G-inv}).

So let us prove this. Fix $g\in S_{n}$ and $a\in\mathcal{A}\mathbf{E}_{T}$ and
$b\in\mathcal{A}\mathbf{F}_{T}$. We must show that $f\left(  ga,\ gb\right)
=f\left(  a,\ b\right)  $.

Write $a$ as $a=u\mathbf{E}_{T}$ for some $u\in\mathcal{A}$ (we can do this,
since $a\in\mathcal{A}\mathbf{E}_{T}$). Write $b$ as $b=v\mathbf{F}_{T}$ for
some $v\in\mathcal{A}$ (we can do this, since $b\in\mathcal{A}\mathbf{F}_{T}%
$). Since $S$ is an algebra anti-morphism, we have $S\left(  gv\right)
=S\left(  v\right)  S\left(  g\right)  =S\left(  v\right)  g^{-1}$ (since the
definition of $S$ yields $S\left(  g\right)  =g^{-1}$).

However, (\ref{pf.thm.dual.specht-field.f=}) (applied to $gu$ and $gv$ instead
of $u$ and $v$) yields%
\begin{align*}
f\left(  gu\mathbf{E}_{T},\ gv\mathbf{F}_{T}\right)   &  =\left[  1\right]
\left(  gu\mathbf{E}_{T}\underbrace{S\left(  gv\right)  }_{=S\left(  v\right)
g^{-1}}\right)  =\left[  1\right]  \left(  gu\mathbf{E}_{T}S\left(  v\right)
g^{-1}\right) \\
&  =\left[  1\right]  \left(  u\mathbf{E}_{T}S\left(  v\right)  \right)
\ \ \ \ \ \ \ \ \ \ \left(  \text{by (\ref{pf.thm.dual.specht-field.c1.1}),
applied to }\mathbf{a}=u\mathbf{E}_{T}S\left(  v\right)  \right) \\
&  =f\left(  u\mathbf{E}_{T},\ v\mathbf{F}_{T}\right)
\ \ \ \ \ \ \ \ \ \ \left(  \text{by (\ref{pf.thm.dual.specht-field.f=}%
)}\right)  .
\end{align*}
In view of $a=u\mathbf{E}_{T}$ and $b=v\mathbf{F}_{T}$, we can rewrite this as
$f\left(  ga,\ gb\right)  =f\left(  a,\ b\right)  $. But this is precisely
what we desired to show. Thus, Claim 3 is proved.
\end{proof}

Now, recall that our bilinear form $f:\mathcal{A}\mathbf{E}_{T}\times
\mathcal{A}\mathbf{F}_{T}\rightarrow\mathbf{k}$ gives rise to two $\mathbf{k}%
$-linear maps $f_{L}:\mathcal{A}\mathbf{E}_{T}\rightarrow\operatorname*{Hom}%
\left(  \mathcal{A}\mathbf{F}_{T},\mathbf{k}\right)  =\left(  \mathcal{A}%
\mathbf{F}_{T}\right)  ^{\ast}$ and $f_{R}:\mathcal{A}\mathbf{F}%
_{T}\rightarrow\operatorname*{Hom}\left(  \mathcal{A}\mathbf{E}_{T}%
,\mathbf{k}\right)  =\left(  \mathcal{A}\mathbf{E}_{T}\right)  ^{\ast}$, as
explained in Definition \ref{def.dual.bilinear}. We shall show that these two
maps $f_{L}$ and $f_{R}$ are bijective. Indeed, we first claim that they are injective:

\begin{statement}
\textit{Claim 4:} The $\mathbf{k}$-linear map $f_{L}:\mathcal{A}\mathbf{E}%
_{T}\rightarrow\left(  \mathcal{A}\mathbf{F}_{T}\right)  ^{\ast}$ is injective.
\end{statement}

\begin{proof}
[Proof of Claim 4.]The map $f_{L}$ is $\mathbf{k}$-linear. Thus, in order to
prove its injectivity, it suffices to show that its kernel
$\operatorname*{Ker}\left(  f_{L}\right)  $ is $0$. In other words, it
suffices to show that $\mathbf{a}=0$ for each $\mathbf{a}\in
\operatorname*{Ker}\left(  f_{L}\right)  $. This is what we shall now show.

Fix $\mathbf{a}\in\operatorname*{Ker}\left(  f_{L}\right)  $. Thus,
$\mathbf{a}\in\mathcal{A}\mathbf{E}_{T}$ and $f_{L}\left(  \mathbf{a}\right)
=0$. We must show that $\mathbf{a}=0$.

Write $\mathbf{a}$ as $\mathbf{a}=u\mathbf{E}_{T}$ for some $u\in\mathcal{A}$
(this can be done, since $\mathbf{a}\in\mathcal{A}\mathbf{E}_{T}$).

Let $w\in S_{n}$ be any permutation. We shall show that $\left[  w\right]
\mathbf{a}=0$. Indeed, we have $\left[  1\right]  \left(  \mathbf{a}%
w^{-1}\right)  =\left[  w\right]  \mathbf{a}$ (by
(\ref{pf.thm.dual.specht-field.c1.2}), applied to $g=w$). Furthermore, the
definition of the antipode $S$ yields $S\left(  w\right)  =w^{-1}$. We have
$w\mathbf{F}_{T}\in\mathcal{A}\mathbf{F}_{T}$ (since $w\in\mathcal{A}$). But
the definition of $f_{L}$ yields%
\begin{align*}
\left(  f_{L}\left(  \mathbf{a}\right)  \right)  \left(  w\mathbf{F}%
_{T}\right)   &  =f\left(  \underbrace{\mathbf{a}}_{=u\mathbf{E}_{T}%
},\ w\mathbf{F}_{T}\right)  =f\left(  u\mathbf{E}_{T},\ w\mathbf{F}_{T}\right)
\\
&  =\left[  1\right]  \left(  \underbrace{u\mathbf{E}_{T}}_{=\mathbf{a}%
}\underbrace{S\left(  w\right)  }_{=w^{-1}}\right)
\ \ \ \ \ \ \ \ \ \ \left(  \text{by (\ref{pf.thm.dual.specht-field.f=}),
applied to }v=w\right) \\
&  =\left[  1\right]  \left(  \mathbf{a}w^{-1}\right)  =\left[  w\right]
\mathbf{a},
\end{align*}
so that%
\[
\left[  w\right]  \mathbf{a}=\underbrace{\left(  f_{L}\left(  \mathbf{a}%
\right)  \right)  }_{=0}\left(  w\mathbf{F}_{T}\right)  =0\left(
w\mathbf{F}_{T}\right)  =0.
\]

Forget that we fixed $w$. We thus have shown that $\left[  w\right]
\mathbf{a}=0$ for all $w\in S_{n}$. In other words, all coordinates of
$\mathbf{a}$ are $0$ (when $\mathbf{a}$ is expanded in the standard basis of
$\mathbf{k}\left[  S_{n}\right]  $). Hence, $\mathbf{a}=0$. As explained
above, this completes the proof of Claim 4.
\end{proof}

\begin{statement}
\textit{Claim 5:} The $\mathbf{k}$-linear map $f_{R}:\mathcal{A}\mathbf{F}%
_{T}\rightarrow\left(  \mathcal{A}\mathbf{E}_{T}\right)  ^{\ast}$ is injective.
\end{statement}

\begin{proof}
[Proof of Claim 5.]The map $f_{R}$ is $\mathbf{k}$-linear. Thus, in order to
prove its injectivity, it suffices to show that its kernel
$\operatorname*{Ker}\left(  f_{R}\right)  $ is $0$. In other words, it
suffices to show that $\mathbf{a}=0$ for each $\mathbf{a}\in
\operatorname*{Ker}\left(  f_{R}\right)  $. This is what we shall now show.

Fix $\mathbf{a}\in\operatorname*{Ker}\left(  f_{R}\right)  $. Thus,
$\mathbf{a}\in\mathcal{A}\mathbf{F}_{T}$ and $f_{R}\left(  \mathbf{a}\right)
=0$. We must show that $\mathbf{a}=0$.

Write $\mathbf{a}$ as $\mathbf{a}=v\mathbf{F}_{T}$ for some $v\in\mathcal{A}$
(this can be done, since $\mathbf{a}\in\mathcal{A}\mathbf{F}_{T}$). Thus,
\begin{align}
S\left(  \mathbf{a}\right)   &  =S\left(  v\mathbf{F}_{T}\right)
=\underbrace{S\left(  \mathbf{F}_{T}\right)  }_{=\mathbf{E}_{T}}S\left(
v\right)  \ \ \ \ \ \ \ \ \ \ \left(  \text{since }S\text{ is an algebra
anti-morphism}\right) \nonumber\\
&  =\mathbf{E}_{T}S\left(  v\right)  . \label{pf.thm.dual.specht-field.c5.2}%
\end{align}

Let $w\in S_{n}$ be any permutation. We shall show that $\left[  w\right]
\mathbf{a}=0$. Indeed, we have $\left[  w\right]  \mathbf{a}=\left[  1\right]
\left(  wS\left(  \mathbf{a}\right)  \right)  $ (by
(\ref{pf.thm.dual.specht-field.c1.3}), applied to $g=w$). We have
$w\mathbf{E}_{T}\in\mathcal{A}\mathbf{E}_{T}$ (since $w\in\mathcal{A}$). But
the definition of $f_{R}$ yields%
\begin{align*}
\left(  f_{R}\left(  \mathbf{a}\right)  \right)  \left(  w\mathbf{E}%
_{T}\right)   &  =f\left(  w\mathbf{E}_{T},\ \underbrace{\mathbf{a}%
}_{=v\mathbf{F}_{T}}\right)  =f\left(  w\mathbf{E}_{T},\ v\mathbf{F}%
_{T}\right) \\
&  =\left[  1\right]  \left(  w\underbrace{\mathbf{E}_{T}S\left(  v\right)
}_{\substack{=S\left(  \mathbf{a}\right)  \\\text{(by
(\ref{pf.thm.dual.specht-field.c5.2}))}}}\right)  \ \ \ \ \ \ \ \ \ \ \left(
\text{by (\ref{pf.thm.dual.specht-field.f=}), applied to }u=w\right) \\
&  =\left[  1\right]  \left(  wS\left(  \mathbf{a}\right)  \right)  =\left[
w\right]  \mathbf{a}\ \ \ \ \ \ \ \ \ \ \left(  \text{since }\left[  w\right]
\mathbf{a}=\left[  1\right]  \left(  wS\left(  \mathbf{a}\right)  \right)
\right)  .
\end{align*}
Hence,
\[
\left[  w\right]  \mathbf{a}=\underbrace{\left(  f_{R}\left(  \mathbf{a}%
\right)  \right)  }_{=0}\left(  w\mathbf{E}_{T}\right)  =0\left(
w\mathbf{E}_{T}\right)  =0.
\]

Forget that we fixed $w$. We thus have shown that $\left[  w\right]
\mathbf{a}=0$ for all $w\in S_{n}$. In other words, all coordinates of
$\mathbf{a}$ are $0$ (when $\mathbf{a}$ is expanded in the standard basis of
$\mathbf{k}\left[  S_{n}\right]  $). Hence, $\mathbf{a}=0$. As we explained
above, this completes the proof of Claim 5.
\end{proof}

Now, recall that $\mathbf{k}$ is a field. Hence, $\mathcal{A}=\mathbf{k}%
\left[  S_{n}\right]  $ is a finite-dimensional $\mathbf{k}$-vector space (of
dimension $\left\vert S_{n}\right\vert =n!$). Its subspaces $\mathcal{A}%
\mathbf{E}_{T}$ and $\mathcal{A}\mathbf{F}_{T}$ are therefore
finite-dimensional as well. Moreover, the $\mathbf{k}$-linear maps
$f_{L}:\mathcal{A}\mathbf{E}_{T}\rightarrow\left(  \mathcal{A}\mathbf{F}%
_{T}\right)  ^{\ast}$ and $f_{R}:\mathcal{A}\mathbf{F}_{T}\rightarrow\left(
\mathcal{A}\mathbf{E}_{T}\right)  ^{\ast}$ are injective (by Claim 4 and Claim
5). Hence, Lemma \ref{lem.dual.field-fLfRinj} (applied to $U=\mathcal{A}%
\mathbf{E}_{T}$ and $V=\mathcal{A}\mathbf{F}_{T}$ and $\alpha=f_{L}$ and
$\beta=f_{R}$) yields that these two maps $f_{L}$ and $f_{R}$ are bijective.
In other words, the maps $f_{L}$ and $f_{R}$ are invertible.

But the bilinear form $f$ is $S_{n}$-invariant (by Claim 3), and thus the map
$f_{R}$ is a morphism of $S_{n}$-representations (by Proposition
\ref{prop.dual.bilf.mor} \textbf{(a)}, applied to $G=S_{n}$ and $U=\mathcal{A}%
\mathbf{E}_{T}$ and $V=\mathcal{A}\mathbf{F}_{T}$). Since $f_{R}$ is
invertible, we thus conclude that $f_{R}$ is an isomorphism of $S_{n}%
$-representations (since an invertible morphism of representations is an
isomorphism). Hence, $\left(  \mathcal{A}\mathbf{E}_{T}\right)  ^{\ast}%
\cong\mathcal{A}\mathbf{F}_{T}$ as $S_{n}$-representations. As we explained
above, this proves Theorem \ref{thm.dual.specht-field}.

(The above proof is somewhat inspired by Fayers's work \cite[\S 4]{Fayers03},
though it is more direct.)
\end{proof}

\subsubsection{The dual of a Specht module: The case of a skew Young diagram}

Can we lift the \textquotedblleft$\mathbf{k}$ is a field\textquotedblright%
\ requirement in Theorem \ref{thm.dual.specht-field}? As I said, I don't know
-- in general. At least, we can lift it when $D$ is a skew Young diagram:

\begin{theorem}
\label{thm.dual.specht-skew}Let $D$ be a skew Young diagram with $\left\vert
D\right\vert =n$. Then, the Specht module $\mathcal{S}^{D}$ satisfies
\[
\left(  \mathcal{S}^{D}\right)  ^{\ast}\cong\left(  \mathcal{S}^{\mathbf{r}%
\left(  D\right)  }\right)  ^{\operatorname*{sign}}\text{ as }S_{n}%
\text{-representations.}%
\]

\end{theorem}

\begin{proof}
[Proof of Theorem \ref{thm.dual.specht-skew}.]We proceed as in the above proof
of Theorem \ref{thm.dual.specht-field}, up until the point where we used the
assumption \textquotedblleft$\mathbf{k}$ is a field\textquotedblright. Since
we no longer have this assumption, we must find a new reason why $f_{R}$ is bijective.

For any diagram $E$, we let $\operatorname*{SYT}\left(  E\right)  $ be the set
of all standard $n$-tableaux of shape $E$. We shall first show that both
$\mathbf{k}$-modules $\mathcal{A}\mathbf{E}_{T}$ and $\mathcal{A}%
\mathbf{F}_{T}$ have finite bases indexed by this set $\operatorname*{SYT}%
\left(  D\right)  $.

Indeed, recall Definition \ref{def.specht.wPQ} and Definition
\ref{def.specht.EPQ}. Recall also that if $P$ is any $n$-tableau of shape $D$,
then $P\mathbf{r}$ is an $n$-tableau of shape $\mathbf{r}\left(  D\right)  $
(see Definition \ref{def.tableaux.r}). Thus, in particular, $T\mathbf{r}$ is
an $n$-tableau of shape $\mathbf{r}\left(  D\right)  $. Now, for each
$P\in\operatorname*{SYT}\left(  D\right)  $, we define the two elements%
\begin{align}
\mathbf{U}_{P}  &  :=\mathbf{E}_{P,T}\ \ \ \ \ \ \ \ \ \ \text{and}%
\label{pf.thm.dual.specht-skew.UP=}\\
\mathbf{V}_{P}  &  :=T_{\operatorname*{sign}}\left(  \mathbf{E}_{P\mathbf{r}%
,T\mathbf{r}}\right)  \label{pf.thm.dual.specht-skew.VP=}%
\end{align}
of $\mathcal{A}$. Recall that $T_{\operatorname*{sign}}$ is a $\mathbf{k}%
$-algebra automorphism of $\mathbf{k}\left[  S_{n}\right]  $ (by Theorem
\ref{thm.Tsign.auto} \textbf{(a)}) and satisfies
\[
T_{\operatorname*{sign}}\left(  \mathbf{E}_{T\mathbf{r}}\right)
=\mathbf{F}_{T}%
\]
(as we saw in the proof of Lemma \ref{lem.sign-twist.AFT}).

Recall that $D$ is a skew Young diagram. Thus, $\mathbf{r}\left(  D\right)  $
is a skew Young diagram as well (by Remark \ref{rmk.partitions.conj.Ylm}%
)\footnote{See the proof of Lemma \ref{lem.specht.EPQ.straighten} for the
details of this argument.}.

Proposition \ref{prop.tableaux.r} \textbf{(c)} (with the variable $T$ renamed
as $Q$) says that the map%
\begin{align}
\operatorname*{SYT}\left(  D\right)   &  \rightarrow\operatorname*{SYT}\left(
\mathbf{r}\left(  D\right)  \right)  ,\nonumber\\
Q  &  \mapsto Q\mathbf{r} \label{pf.thm.dual.specht-skew.bij}%
\end{align}
is a bijection.

Now, we claim the following:\footnote{We are numbering the claims starting
from 6 since we are continuing the proof of Theorem
\ref{thm.dual.specht-field}.}

\begin{statement}
\textit{Claim 6:} The family $\left(  \mathbf{U}_{P}\right)  _{P\in
\operatorname*{SYT}\left(  D\right)  }$ is a basis of the $\mathbf{k}$-module
$\mathcal{A}\mathbf{E}_{T}$.
\end{statement}

\begin{proof}
[Proof of Claim 6.]From (\ref{pf.thm.dual.specht-skew.UP=}), we obtain
$\left(  \mathbf{U}_{P}\right)  _{P\in\operatorname*{SYT}\left(  D\right)
}=\left(  \mathbf{E}_{P,T}\right)  _{P\in\operatorname*{SYT}\left(  D\right)
}$.

But Proposition \ref{prop.specht.EPQ.basis-AET} shows that the family $\left(
\mathbf{E}_{P,T}\right)  _{P\in\operatorname*{SYT}\left(  D\right)  }$ is a
basis of the $\mathbf{k}$-module $\mathcal{A}\mathbf{E}_{T}$ (since $D$ is a
skew Young diagram). In view of $\left(  \mathbf{U}_{P}\right)  _{P\in
\operatorname*{SYT}\left(  D\right)  }=\left(  \mathbf{E}_{P,T}\right)
_{P\in\operatorname*{SYT}\left(  D\right)  }$, we can rewrite this as follows:
The family $\left(  \mathbf{U}_{P}\right)  _{P\in\operatorname*{SYT}\left(
D\right)  }$ is a basis of the $\mathbf{k}$-module $\mathcal{A}\mathbf{E}_{T}%
$. This proves Claim 6.
\end{proof}

\begin{statement}
\textit{Claim 7:} The family $\left(  \mathbf{V}_{Q}\right)  _{Q\in
\operatorname*{SYT}\left(  D\right)  }$ is a basis of the $\mathbf{k}$-module
$\mathcal{A}\mathbf{F}_{T}$.
\end{statement}

\begin{proof}
[Proof of Claim 7.]For each $Q\in\operatorname*{SYT}\left(  D\right)  $, we
have $\mathbf{V}_{Q}=T_{\operatorname*{sign}}\left(  \mathbf{E}_{Q\mathbf{r}%
,T\mathbf{r}}\right)  $ (by (\ref{pf.thm.dual.specht-skew.VP=}), applied to
$P=Q$). In other words, $\left(  \mathbf{V}_{Q}\right)  _{Q\in
\operatorname*{SYT}\left(  D\right)  }=\left(  T_{\operatorname*{sign}}\left(
\mathbf{E}_{Q\mathbf{r},T\mathbf{r}}\right)  \right)  _{Q\in
\operatorname*{SYT}\left(  D\right)  }$.

Furthermore, $T_{\operatorname*{sign}}\left(  \mathbf{k}\left[  S_{n}\right]
\right)  =\mathbf{k}\left[  S_{n}\right]  $ (since $T_{\operatorname*{sign}}$
is an automorphism). In other words, $T_{\operatorname*{sign}}\left(
\mathcal{A}\right)  =\mathcal{A}$ (since $\mathcal{A}=\mathbf{k}\left[
S_{n}\right]  $). Now, since $T_{\operatorname*{sign}}$ is a $\mathbf{k}%
$-algebra morphism, we have%
\[
T_{\operatorname*{sign}}\left(  \mathcal{A}\mathbf{E}_{T\mathbf{r}}\right)
=\underbrace{T_{\operatorname*{sign}}\left(  \mathcal{A}\right)
}_{=\mathcal{A}}\cdot\underbrace{T_{\operatorname*{sign}}\left(
\mathbf{E}_{T\mathbf{r}}\right)  }_{=\mathbf{F}_{T}}=\mathcal{A}\mathbf{F}%
_{T}.
\]

But $T\mathbf{r}$ is an $n$-tableau of shape $\mathbf{r}\left(  D\right)  $
(since $T$ is an $n$-tableau of shape $D$). Hence, Proposition
\ref{prop.specht.EPQ.basis-AET} (applied to $\mathbf{r}\left(  D\right)  $ and
$T\mathbf{r}$ instead of $D$ and $T$) shows that the family $\left(
\mathbf{E}_{P,T\mathbf{r}}\right)  _{P\in\operatorname*{SYT}\left(
\mathbf{r}\left(  D\right)  \right)  }$ is a basis of the $\mathbf{k}$-module
$\mathcal{A}\mathbf{E}_{T\mathbf{r}}$ (since $\mathbf{r}\left(  D\right)  $ is
a skew Young diagram).

Using the bijection (\ref{pf.thm.dual.specht-skew.bij}), we can reindex this
family as $\left(  \mathbf{E}_{Q\mathbf{r},T\mathbf{r}}\right)  _{Q\in
\operatorname*{SYT}\left(  D\right)  }$ (by substituting $Q\mathbf{r}$ for $P$
in this family). Thus, we conclude that the latter family $\left(
\mathbf{E}_{Q\mathbf{r},T\mathbf{r}}\right)  _{Q\in\operatorname*{SYT}\left(
D\right)  }$ is a basis of the $\mathbf{k}$-module $\mathcal{A}\mathbf{E}%
_{T\mathbf{r}}$ as well (since reindexing a basis using a bijection yields a
basis again).

But $T_{\operatorname*{sign}}$ is a $\mathbf{k}$-module isomorphism (being a
$\mathbf{k}$-algebra automorphism), and sends the $\mathbf{k}$-module
$\mathcal{A}\mathbf{E}_{T\mathbf{r}}$ to $\mathcal{A}\mathbf{F}_{T}$ (since
$T_{\operatorname*{sign}}\left(  \mathcal{A}\mathbf{E}_{T\mathbf{r}}\right)
=\mathcal{A}\mathbf{F}_{T}$). Hence, it sends any basis of $\mathcal{A}%
\mathbf{E}_{T\mathbf{r}}$ to a basis of $\mathcal{A}\mathbf{F}_{T}$ (since a
$\mathbf{k}$-module isomorphism sends any basis to a basis). Therefore, since
$\left(  \mathbf{E}_{Q\mathbf{r},T\mathbf{r}}\right)  _{Q\in
\operatorname*{SYT}\left(  D\right)  }$ is a basis of the $\mathbf{k}$-module
$\mathcal{A}\mathbf{E}_{T\mathbf{r}}$, we conclude that the family $\left(
T_{\operatorname*{sign}}\left(  \mathbf{E}_{Q\mathbf{r},T\mathbf{r}}\right)
\right)  _{Q\in\operatorname*{SYT}\left(  D\right)  }$ is a basis of the
$\mathbf{k}$-module $\mathcal{A}\mathbf{F}_{T}$ (since this family $\left(
T_{\operatorname*{sign}}\left(  \mathbf{E}_{Q\mathbf{r},T\mathbf{r}}\right)
\right)  _{Q\in\operatorname*{SYT}\left(  D\right)  }$ is the image of
$\left(  \mathbf{E}_{Q\mathbf{r},T\mathbf{r}}\right)  _{Q\in
\operatorname*{SYT}\left(  D\right)  }$ under $T_{\operatorname*{sign}}$). In
view of $\left(  \mathbf{V}_{Q}\right)  _{Q\in\operatorname*{SYT}\left(
D\right)  }=\left(  T_{\operatorname*{sign}}\left(  \mathbf{E}_{Q\mathbf{r}%
,T\mathbf{r}}\right)  \right)  _{Q\in\operatorname*{SYT}\left(  D\right)  }$,
we can rewrite this as follows: The family $\left(  \mathbf{V}_{Q}\right)
_{Q\in\operatorname*{SYT}\left(  D\right)  }$ is a basis of the $\mathbf{k}%
$-module $\mathcal{A}\mathbf{F}_{T}$. This proves Claim 7.
\end{proof}

We have now found finite bases $\left(  \mathbf{U}_{P}\right)  _{P\in
\operatorname*{SYT}\left(  D\right)  }$ and $\left(  \mathbf{V}_{Q}\right)
_{Q\in\operatorname*{SYT}\left(  D\right)  }$ of the two $\mathbf{k}$-modules
$\mathcal{A}\mathbf{E}_{T}$ and $\mathcal{A}\mathbf{F}_{T}$, respectively (in
Claim 6 and Claim 7). Thus, we can construct the Gram matrix\footnote{See
Definition \ref{def.dual.gram} for the meaning of \textquotedblleft Gram
matrix\textquotedblright.}%
\[
\left(  f\left(  \mathbf{U}_{P},\ \mathbf{V}_{Q}\right)  \right)  _{\left(
P,Q\right)  \in\operatorname*{SYT}\left(  D\right)  \times\operatorname*{SYT}%
\left(  D\right)  }%
\]
of the bilinear form $f:\mathcal{A}\mathbf{E}_{T}\times\mathcal{A}%
\mathbf{F}_{T}\rightarrow\mathbf{k}$ with respect to these two bases. Let us
denote this matrix by $M_{\mathbf{k}}$, in order to make its dependence on
$\mathbf{k}$ explicit. Note that this Gram matrix $M_{\mathbf{k}}$ is a square
matrix, since both its rows and its columns are indexed by
$\operatorname*{SYT}\left(  D\right)  $. Thus, it has a determinant
$\det\left(  M_{\mathbf{k}}\right)  $. This determinant determines whether
$f_{R}$ is invertible:

\begin{statement}
\textit{Claim 8:} The map $f_{R}:\mathcal{A}\mathbf{F}_{T}\rightarrow\left(
\mathcal{A}\mathbf{E}_{T}\right)  ^{\ast}$ is invertible if and only if the
element $\det\left(  M_{\mathbf{k}}\right)  $ of $\mathbf{k}$ is invertible.
\end{statement}

\begin{proof}
[Proof of Claim 8.]Proposition \ref{prop.dual.bilf.rep-mat} \textbf{(b)}
(applied to $U=\mathcal{A}\mathbf{E}_{T}$ and $V=\mathcal{A}\mathbf{F}_{T}$
and $I=\operatorname*{SYT}\left(  D\right)  $ and $J=\operatorname*{SYT}%
\left(  D\right)  $ and $\left(  u_{i}\right)  _{i\in I}=\left(
\mathbf{U}_{P}\right)  _{P\in\operatorname*{SYT}\left(  D\right)  }$ and
$\left(  v_{j}\right)  _{j\in J}=\left(  \mathbf{V}_{Q}\right)  _{Q\in
\operatorname*{SYT}\left(  D\right)  }$) shows that the $\mathbf{k}$-linear
map $f_{R}:\mathcal{A}\mathbf{F}_{T}\rightarrow\left(  \mathcal{A}%
\mathbf{E}_{T}\right)  ^{\ast}$ is invertible if and only if the matrix
$\left(  f\left(  \mathbf{U}_{P},\ \mathbf{V}_{Q}\right)  \right)  _{\left(
P,Q\right)  \in\operatorname*{SYT}\left(  D\right)  \times\operatorname*{SYT}%
\left(  D\right)  }\in\mathbf{k}^{\operatorname*{SYT}\left(  D\right)
\times\operatorname*{SYT}\left(  D\right)  }$ is invertible. Hence, we have
the following chain of logical equivalences:%
\begin{align*}
&  \ \left(  \text{the map }f_{R}:\mathcal{A}\mathbf{F}_{T}\rightarrow\left(
\mathcal{A}\mathbf{E}_{T}\right)  ^{\ast}\text{ is invertible}\right) \\
&  \Longleftrightarrow\ \left(  \text{the matrix }\left(  f\left(
\mathbf{U}_{P},\ \mathbf{V}_{Q}\right)  \right)  _{\left(  P,Q\right)
\in\operatorname*{SYT}\left(  D\right)  \times\operatorname*{SYT}\left(
D\right)  }\text{ is invertible}\right) \\
&  \Longleftrightarrow\ \left(  \text{the matrix }M_{\mathbf{k}}\text{ is
invertible}\right) \\
&  \ \ \ \ \ \ \ \ \ \ \ \ \ \ \ \ \ \ \ \ \left(  \text{since }\left(
f\left(  \mathbf{U}_{P},\ \mathbf{V}_{Q}\right)  \right)  _{\left(
P,Q\right)  \in\operatorname*{SYT}\left(  D\right)  \times\operatorname*{SYT}%
\left(  D\right)  }=M_{\mathbf{k}}\right) \\
&  \Longleftrightarrow\ \left(  \text{the element }\det\left(  M_{\mathbf{k}%
}\right)  \text{ of }\mathbf{k}\text{ is invertible}\right)
\end{align*}
(since a well-known fact (e.g., \cite[Theorem 6.110 \textbf{(a)}]{detnotes} or
\cite[Corollary 9.51 and Remark 9.52]{Loehr-BC}) says that a square matrix
over a commutative ring is invertible if and only if its determinant is
invertible). This proves Claim 8.
\end{proof}

Now, let us take a closer look at how the Gram matrix $M_{\mathbf{k}}$ depends
on the ring $\mathbf{k}$. We claim that it essentially does not: i.e., the
entries of $M_{\mathbf{k}}$ are the same integers for all $\mathbf{k}$, just
mapped into $\mathbf{k}$ using the canonical ring morphisms from $\mathbb{Z}$
to $\mathbf{k}$. In fact, we can describe these entries explicitly:

\begin{statement}
\textit{Claim 9:} Let $P\in\operatorname*{SYT}\left(  D\right)  $ and
$Q\in\operatorname*{SYT}\left(  D\right)  $. Then,%
\[
f\left(  \mathbf{U}_{P},\ \mathbf{V}_{Q}\right)  =\sum_{\substack{\left(
c,r\right)  \in\mathcal{C}\left(  T\right)  \times\mathcal{R}\left(  T\right)
;\\w_{Q,T}=w_{P,T}cr}}\left(  -1\right)  ^{c}\left(  -1\right)  ^{w_{Q,T}}.
\]

\end{statement}

\begin{proof}
[Proof of Claim 9.]From (\ref{pf.thm.dual.specht-skew.UP=}), we obtain%
\[
\mathbf{U}_{P}=\mathbf{E}_{P,T}=w_{P,T}\mathbf{E}_{T}%
\]
(by (\ref{eq.prop.specht.EPQ.EPEQ.wE}), applied to $T$ instead of $Q$).
Furthermore, Proposition \ref{prop.specht.wPQ.wQP} \textbf{(c)} (applied to
$T$ instead of $P$) yields $w_{Q,T}=w_{Q\mathbf{r},T\mathbf{r}}$, so that
$w_{Q\mathbf{r},T\mathbf{r}}=w_{Q,T}$. However,
(\ref{eq.prop.specht.EPQ.EPEQ.wE}) (applied to $\mathbf{r}\left(  D\right)  $,
$Q\mathbf{r}$ and $T\mathbf{r}$ instead of $D$, $P$ and $Q$) yields%
\[
\mathbf{E}_{Q\mathbf{r},T\mathbf{r}}=\underbrace{w_{Q\mathbf{r},T\mathbf{r}}%
}_{\substack{=w_{Q,T}}}\mathbf{E}_{T\mathbf{r}}=w_{Q,T}\mathbf{E}%
_{T\mathbf{r}}.
\]
Now, from (\ref{pf.thm.dual.specht-skew.VP=}) (applied to $Q$ instead of $P$),
we obtain%
\begin{align*}
\mathbf{V}_{Q}  &  =T_{\operatorname*{sign}}\left(  \underbrace{\mathbf{E}%
_{Q\mathbf{r},T\mathbf{r}}}_{=w_{Q,T}\mathbf{E}_{T\mathbf{r}}}\right)
=T_{\operatorname*{sign}}\left(  w_{Q,T}\mathbf{E}_{T\mathbf{r}}\right) \\
&  =\underbrace{T_{\operatorname*{sign}}\left(  w_{Q,T}\right)  }%
_{\substack{=\left(  -1\right)  ^{w_{Q,T}}w_{Q,T}\\\text{(by the
definition}\\\text{of }T_{\operatorname*{sign}}\text{, since }w_{Q,T}\in
S_{n}\text{)}}}\cdot\underbrace{T_{\operatorname*{sign}}\left(  \mathbf{E}%
_{T\mathbf{r}}\right)  }_{=\mathbf{F}_{T}}\ \ \ \ \ \ \ \ \ \ \left(
\text{since }T_{\operatorname*{sign}}\text{ is a }\mathbf{k}\text{-algebra
morphism}\right) \\
&  =\left(  -1\right)  ^{w_{Q,T}}w_{Q,T}\mathbf{F}_{T}.
\end{align*}
Now,
\begin{align*}
f\left(  \underbrace{\mathbf{U}_{P}}_{=w_{P,T}\mathbf{E}_{T}}%
,\ \underbrace{\mathbf{V}_{Q}}_{=\left(  -1\right)  ^{w_{Q,T}}w_{Q,T}%
\mathbf{F}_{T}}\right)   &  =f\left(  w_{P,T}\mathbf{E}_{T},\ \left(
-1\right)  ^{w_{Q,T}}w_{Q,T}\mathbf{F}_{T}\right) \\
&  =\left[  1\right]  \left(  w_{P,T}\mathbf{E}_{T}S\left(  \left(  -1\right)
^{w_{Q,T}}w_{Q,T}\right)  \right)
\end{align*}
(by (\ref{pf.thm.dual.specht-field.f=}), applied to $u=w_{P,T}$ and $v=\left(
-1\right)  ^{w_{Q,T}}w_{Q,T}$). In view of%
\begin{align*}
&  w_{P,T}\mathbf{E}_{T}\underbrace{S\left(  \left(  -1\right)  ^{w_{Q,T}%
}w_{Q,T}\right)  }_{\substack{=\left(  -1\right)  ^{w_{Q,T}}S\left(
w_{Q,T}\right)  \\\text{(since the map }S\text{ is }\mathbf{k}\text{-linear)}%
}}\\
&  =w_{P,T}\underbrace{\mathbf{E}_{T}}_{\substack{=\sum_{\left(  c,r\right)
\in\mathcal{C}\left(  T\right)  \times\mathcal{R}\left(  T\right)  }\left(
-1\right)  ^{c}cr\\\text{(by (\ref{pf.lem.specht.ET.1-coord.1}))}}}\left(
-1\right)  ^{w_{Q,T}}\underbrace{S\left(  w_{Q,T}\right)  }%
_{\substack{=w_{Q,T}^{-1}\\\text{(by the definition of }S\text{,}\\\text{since
}w_{Q,T}\in S_{n}\text{)}}}\\
&  =w_{P,T}\left(  \sum_{\left(  c,r\right)  \in\mathcal{C}\left(  T\right)
\times\mathcal{R}\left(  T\right)  }\left(  -1\right)  ^{c}cr\right)  \left(
-1\right)  ^{w_{Q,T}}w_{Q,T}^{-1}\\
&  =\sum_{\left(  c,r\right)  \in\mathcal{C}\left(  T\right)  \times
\mathcal{R}\left(  T\right)  }\underbrace{w_{P,T}\left(  -1\right)
^{c}cr\left(  -1\right)  ^{w_{Q,T}}w_{Q,T}^{-1}}_{=\left(  -1\right)
^{c}\left(  -1\right)  ^{w_{Q,T}}w_{P,T}crw_{Q,T}^{-1}}\\
&  =\sum_{\left(  c,r\right)  \in\mathcal{C}\left(  T\right)  \times
\mathcal{R}\left(  T\right)  }\left(  -1\right)  ^{c}\left(  -1\right)
^{w_{Q,T}}w_{P,T}crw_{Q,T}^{-1},
\end{align*}
we can rewrite this as%
\begin{align*}
f\left(  \mathbf{U}_{P},\ \mathbf{V}_{Q}\right)   &  =\left[  1\right]
\left(  \sum_{\left(  c,r\right)  \in\mathcal{C}\left(  T\right)
\times\mathcal{R}\left(  T\right)  }\left(  -1\right)  ^{c}\left(  -1\right)
^{w_{Q,T}}w_{P,T}crw_{Q,T}^{-1}\right) \\
&  =\sum_{\left(  c,r\right)  \in\mathcal{C}\left(  T\right)  \times
\mathcal{R}\left(  T\right)  }\left(  -1\right)  ^{c}\left(  -1\right)
^{w_{Q,T}}\cdot\underbrace{\left[  1\right]  \left(  w_{P,T}crw_{Q,T}%
^{-1}\right)  }_{\substack{=\left[  w_{Q,T}\right]  \left(  w_{P,T}cr\right)
\\\text{(by (\ref{pf.thm.dual.specht-field.c1.2}), applied to }g=w_{Q,T}%
\\\text{and }\mathbf{a}=w_{P,T}cr\text{)}}}\\
&  =\sum_{\left(  c,r\right)  \in\mathcal{C}\left(  T\right)  \times
\mathcal{R}\left(  T\right)  }\left(  -1\right)  ^{c}\left(  -1\right)
^{w_{Q,T}}\cdot\underbrace{\left[  w_{Q,T}\right]  \left(  w_{P,T}cr\right)
}_{\substack{=\delta_{w_{Q,T},w_{P,T}cr}\\\text{(since }\left[  u\right]
v=\delta_{u,v}\text{ for any }u,v\in S_{n}\text{)}}}\\
&  \ \ \ \ \ \ \ \ \ \ \ \ \ \ \ \ \ \ \ \ \left(
\begin{array}
[c]{c}%
\text{since addition and scaling in }\mathbf{k}\left[  S_{n}\right] \\
\text{are coordinatewise}%
\end{array}
\right) \\
&  =\sum_{\left(  c,r\right)  \in\mathcal{C}\left(  T\right)  \times
\mathcal{R}\left(  T\right)  }\left(  -1\right)  ^{c}\left(  -1\right)
^{w_{Q,T}}\delta_{w_{Q,T},w_{P,T}cr}\\
&  =\sum_{\substack{\left(  c,r\right)  \in\mathcal{C}\left(  T\right)
\times\mathcal{R}\left(  T\right)  ;\\w_{Q,T}=w_{P,T}cr}}\left(  -1\right)
^{c}\left(  -1\right)  ^{w_{Q,T}}\underbrace{\delta_{w_{Q,T},w_{P,T}cr}%
}_{\substack{=1\\\text{(since }1=w_{P,T}crw_{Q,T}^{-1}\text{)}}}\\
&  \ \ \ \ \ \ \ \ \ \ +\sum_{\substack{\left(  c,r\right)  \in\mathcal{C}%
\left(  T\right)  \times\mathcal{R}\left(  T\right)  ;\\w_{Q,T}\neq w_{P,T}%
cr}}\left(  -1\right)  ^{c}\left(  -1\right)  ^{w_{Q,T}}\underbrace{\delta
_{w_{Q,T},w_{P,T}cr}}_{\substack{=0\\\text{(since }1\neq w_{P,T}crw_{Q,T}%
^{-1}\text{)}}}\\
&  \ \ \ \ \ \ \ \ \ \ \ \ \ \ \ \ \ \ \ \ \left(
\begin{array}
[c]{c}%
\text{since each }\left(  c,r\right)  \in\mathcal{C}\left(  T\right)
\times\mathcal{R}\left(  T\right)  \text{ satisfies}\\
\text{either }w_{Q,T}=w_{P,T}cr\text{ or }w_{Q,T}\neq w_{P,T}cr\\
\text{(but not both)}%
\end{array}
\right) \\
&  =\sum_{\substack{\left(  c,r\right)  \in\mathcal{C}\left(  T\right)
\times\mathcal{R}\left(  T\right)  ;\\w_{Q,T}=w_{P,T}cr}}\left(  -1\right)
^{c}\left(  -1\right)  ^{w_{Q,T}}+\underbrace{\sum_{\substack{\left(
c,r\right)  \in\mathcal{C}\left(  T\right)  \times\mathcal{R}\left(  T\right)
;\\w_{Q,T}\neq w_{P,T}cr}}\left(  -1\right)  ^{c}\left(  -1\right)  ^{w_{Q,T}%
}0}_{=0}\\
&  =\sum_{\substack{\left(  c,r\right)  \in\mathcal{C}\left(  T\right)
\times\mathcal{R}\left(  T\right)  ;\\w_{Q,T}=w_{P,T}cr}}\left(  -1\right)
^{c}\left(  -1\right)  ^{w_{Q,T}}.
\end{align*}
This proves Claim 9.
\end{proof}

We shall now use a \textquotedblleft base change\textquotedblright\ argument
somewhat similar to the one in the proof of Proposition
\ref{prop.spechtmod.Elam.idp-weak} but even stranger, in that we will use
finite fields $\mathbb{Z}/p$ for prime numbers $p$.

We let $\varphi_{\mathbf{k}}$ denote the canonical ring morphism from
$\mathbb{Z}$ to $\mathbf{k}$ (sending each integer $m$ to $m\cdot
1_{\mathbf{k}}$). (This was denoted $f$ in some prior proofs, but the letter
$f$ now means something different.) We can now relate the determinant of the
Gram matrix $M_{\mathbf{k}}$ defined over the ring $\mathbf{k}$ with the
determinant of the Gram matrix $M_{\mathbb{Z}}$ defined over the ring
$\mathbb{Z}$:

\begin{statement}
\textit{Claim 10:} We have%
\[
\det\left(  M_{\mathbf{k}}\right)  =\varphi_{\mathbf{k}}\left(  \det\left(
M_{\mathbb{Z}}\right)  \right)  .
\]

\end{statement}

\begin{proof}
[Proof of Claim 10.]We shall compare respective entries of the matrices
$M_{\mathbb{Z}}$ and $M_{\mathbf{k}}$. Indeed, fix $P\in\operatorname*{SYT}%
\left(  D\right)  $ and $Q\in\operatorname*{SYT}\left(  D\right)  $. Then, by
the definition of $M_{\mathbf{k}}$, we know that the $\left(  P,Q\right)  $-th
entry of the matrix $M_{\mathbf{k}}$ is%
\begin{equation}
f\left(  \mathbf{U}_{P},\ \mathbf{V}_{Q}\right)  =\sum_{\substack{\left(
c,r\right)  \in\mathcal{C}\left(  T\right)  \times\mathcal{R}\left(  T\right)
;\\w_{Q,T}=w_{P,T}cr}}\left(  -1\right)  ^{c}\left(  -1\right)  ^{w_{Q,T}}
\label{pf.thm.dual.specht-skew.C10.1}%
\end{equation}
(by Claim 9). But the same argument (applied to the ring $\mathbb{Z}$ instead
of $\mathbf{k}$) shows that the $\left(  P,Q\right)  $-th entry of the matrix
$M_{\mathbb{Z}}$ is%
\begin{equation}
\sum_{\substack{\left(  c,r\right)  \in\mathcal{C}\left(  T\right)
\times\mathcal{R}\left(  T\right)  ;\\w_{Q,T}=w_{P,T}cr}}\left(  -1\right)
^{c}\left(  -1\right)  ^{w_{Q,T}} \label{pf.thm.dual.specht-skew.C10.2}%
\end{equation}
as well, but now interpreted as an actual integer rather than as an element of
$\mathbf{k}$. Clearly, the ring morphism $\varphi_{\mathbf{k}}:\mathbb{Z}%
\rightarrow\mathbf{k}$ sends the sum in (\ref{pf.thm.dual.specht-skew.C10.2})
to the sum in (\ref{pf.thm.dual.specht-skew.C10.1}) (since this morphism
$\varphi_{\mathbf{k}}$ sends each integer $m$ to the corresponding element
$m\cdot1_{\mathbf{k}}$ of $\mathbf{k}$). In other words, the morphism
$\varphi_{\mathbf{k}}$ sends the $\left(  P,Q\right)  $-th entry of the matrix
$M_{\mathbb{Z}}$ to the $\left(  P,Q\right)  $-th entry of the matrix
$M_{\mathbf{k}}$ (since the former entry is the sum in
(\ref{pf.thm.dual.specht-skew.C10.2}), whereas the latter entry is the sum in
(\ref{pf.thm.dual.specht-skew.C10.1})).

Forget that we fixed $P$ and $Q$. We thus have shown that for each
$P\in\operatorname*{SYT}\left(  D\right)  $ and $Q\in\operatorname*{SYT}%
\left(  D\right)  $, the morphism $\varphi_{\mathbf{k}}$ sends the $\left(
P,Q\right)  $-th entry of the matrix $M_{\mathbb{Z}}$ to the $\left(
P,Q\right)  $-th entry of the matrix $M_{\mathbf{k}}$. In other words, the
morphism $\varphi_{\mathbf{k}}$ sends each entry of the matrix $M_{\mathbb{Z}%
}$ to the corresponding entry of the matrix $M_{\mathbf{k}}$. Since
$\varphi_{\mathbf{k}}$ is a ring morphism, this entails that $\varphi
_{\mathbf{k}}$ sends the determinant of $M_{\mathbb{Z}}$ to the determinant of
$M_{\mathbf{k}}$ (since ring morphisms respect determinants). In other words,
$\varphi_{\mathbf{k}}\left(  \det\left(  M_{\mathbb{Z}}\right)  \right)
=\det\left(  M_{\mathbf{k}}\right)  $. This proves Claim 10.
\end{proof}

Now comes the actually surprising part of the proof:

\begin{statement}
\textit{Claim 11:} The integer $\det\left(  M_{\mathbb{Z}}\right)  $ is
invertible in $\mathbb{Z}$.
\end{statement}

\begin{proof}
[Proof of Claim 11.]Assume the contrary. Thus, $\det\left(  M_{\mathbb{Z}%
}\right)  $ is not invertible in $\mathbb{Z}$. Hence, $\det\left(
M_{\mathbb{Z}}\right)  $ is neither $1$ nor $-1$ (since the integers $1$ and
$-1$ are invertible in $\mathbb{Z}$). Therefore, $\det\left(  M_{\mathbb{Z}%
}\right)  $ is divisible by some prime number $p$ (since any integer other
than $1$ and $-1$ is divisible by some prime). Consider this $p$. The
canonical ring morphism $\varphi_{\mathbb{Z}/p}:\mathbb{Z}\rightarrow
\mathbb{Z}/p$ sends each integer that is divisible by $p$ to $0$. Hence,
$\varphi_{\mathbb{Z}/p}$ sends $\det\left(  M_{\mathbb{Z}}\right)  $ to $0$
(since $\det\left(  M_{\mathbb{Z}}\right)  $ is divisible by $p$). However,
applying Claim 10 to the field $\mathbb{Z}/p$ instead of $\mathbf{k}$, we
obtain%
\[
\det\left(  M_{\mathbb{Z}/p}\right)  =\varphi_{\mathbb{Z}/p}\left(
\det\left(  M_{\mathbb{Z}}\right)  \right)  =0
\]
(since $\varphi_{\mathbb{Z}/p}$ sends $\det\left(  M_{\mathbb{Z}}\right)  $ to
$0$).

However, in the above proof of Theorem \ref{thm.dual.specht-field}, we have
shown that the map $f_{R}$ is invertible \textbf{if} $\mathbf{k}$ is a field.
We cannot apply this fact to our ring $\mathbf{k}$, since $\mathbf{k}$ might
not be a field. However, we can apply it to $\mathbb{Z}/p$ instead of
$\mathbf{k}$ (since $\mathbb{Z}/p$ is a field). Thus, we find that the map
$f_{R}$ defined over the field $\mathbb{Z}/p$ instead of $\mathbf{k}$ is invertible.

But we can also apply Claim 8 to $\mathbb{Z}/p$ instead of $\mathbf{k}$. As a
result, we see that the map $f_{R}$ defined over the field $\mathbb{Z}/p$
instead of $\mathbf{k}$ is invertible if and only if the element $\det\left(
M_{\mathbb{Z}/p}\right)  $ of $\mathbb{Z}/p$ is invertible. Hence, the element
$\det\left(  M_{\mathbb{Z}/p}\right)  $ of $\mathbb{Z}/p$ is invertible (since
the map $f_{R}$ defined over the field $\mathbb{Z}/p$ instead of $\mathbf{k}$
is invertible). In other words, the element $0$ of $\mathbb{Z}/p$ is
invertible (since $\det\left(  M_{\mathbb{Z}/p}\right)  =0$). This contradicts
the fact that the element $0$ of $\mathbb{Z}/p$ is clearly not invertible
(since $0$ is not invertible in any field).

This contradiction shows that our assumption was false. Hence, Claim 11 is proved.
\end{proof}

Now, we are almost done. Claim 11 shows that the integer $\det\left(
M_{\mathbb{Z}}\right)  $ is invertible in $\mathbb{Z}$. Hence, its image
$\varphi_{\mathbf{k}}\left(  \det\left(  M_{\mathbb{Z}}\right)  \right)  $ is
invertible in $\mathbf{k}$ (since $\varphi_{\mathbf{k}}$ is a ring morphism
and thus sends invertible elements to invertible elements). In other words,
$\det\left(  M_{\mathbf{k}}\right)  $ is invertible in $\mathbf{k}$ (since
Claim 10 says that $\det\left(  M_{\mathbf{k}}\right)  =\varphi_{\mathbf{k}%
}\left(  \det\left(  M_{\mathbb{Z}}\right)  \right)  $).

In other words, the element $\det\left(  M_{\mathbf{k}}\right)  $ of
$\mathbf{k}$ is invertible. By Claim 8, this entails that the map
$f_{R}:\mathcal{A}\mathbf{F}_{T}\rightarrow\left(  \mathcal{A}\mathbf{E}%
_{T}\right)  ^{\ast}$ is invertible.

From here, we can finish the proof of Theorem \ref{thm.dual.specht-skew} in
the exact same way as we finished the proof of Theorem
\ref{thm.dual.specht-field} after we showed that the map $f_{R}$ is invertible.
\end{proof}

The above proof of Theorem \ref{thm.dual.specht-skew} works not only when $D$
is a skew Young diagram, but also more generally if we merely assume that the
$\mathbf{k}$-modules $\mathcal{A}\mathbf{E}_{T}$ and $\mathcal{A}%
\mathbf{F}_{T}$ (or, equivalently, $\mathcal{S}^{D}$ and $\mathcal{S}%
^{\mathbf{r}\left(  D\right)  }$) have bases independent of $\mathbf{k}$ (that
is, bases indexed by combinatorial objects that are \textquotedblleft
respected\textquotedblright\ by the canonical ring morphisms $\varphi
_{\mathbf{k}}:\mathbb{Z}\rightarrow\mathbf{k}$). This holds for various
diagrams that are not skew Young diagrams, such as the diagram from Exercise
\ref{exe.spechtmod.bad.6}. Whether this holds for arbitrary $D$ is an open
question. Thus, we do not know how far Theorem \ref{thm.dual.specht-skew} can
actually be generalized.

\begin{exercise}
There is an alternative proof of Theorem \ref{thm.dual.specht-skew}, which
avoids the use of base change but instead argues by computing $\det\left(
M_{\mathbf{k}}\right)  $. Namely, using the notations from our above proof of
Theorem \ref{thm.dual.specht-skew}: \medskip

\textbf{(a)} \fbox{3} Show that every $P,Q\in\operatorname*{SYT}\left(
D\right)  $ satisfy
\[
f\left(  \mathbf{U}_{P},\ \mathbf{V}_{Q}\right)  =\sum_{\substack{\left(
c,r\right)  \in\mathcal{C}\left(  P\right)  \times\mathcal{R}\left(  P\right)
;\\Q=crP}}\left(  -1\right)  ^{c}\left(  -1\right)  ^{w_{Q,T}}.
\]

\textbf{(b)} \fbox{2} Show that $f\left(  \mathbf{U}_{P},\ \mathbf{V}%
_{Q}\right)  =0$ for every $P,Q\in\operatorname*{SYT}\left(  D\right)  $ that
satisfy $\overline{P}<\overline{Q}$ in the Young last letter order (recall
Definition \ref{def.tabloid.llo}). \medskip

\textbf{(c)} \fbox{1} Show that $f\left(  \mathbf{U}_{P},\ \mathbf{V}%
_{P}\right)  =\left(  -1\right)  ^{w_{P,T}}$ for all $P\in\operatorname*{SYT}%
\left(  D\right)  $. \medskip

\textbf{(d)} \fbox{1} Conclude that the matrix $M_{\mathbf{k}}$ is triangular
and has determinant $\pm1$, and thus is invertible.
\end{exercise}

\subsubsection{Consequences in characteristic $0$}

Under certain special conditions, the dual of a Specht module $\mathcal{S}%
^{D}$ has an even simpler description than the one given by Theorem
\ref{thm.dual.specht-skew}. Namely, if $D$ is the Young diagram $Y\left(
\lambda\right)  $ of a partition $\lambda$, and if the number $h^{\lambda}$
from Definition \ref{def.spechtmod.flam} is invertible in $\mathbf{k}$, then
the dual of $\mathcal{S}^{D}=\mathcal{S}^{\lambda}$ is just $\mathcal{S}%
^{\lambda}$ itself:

\begin{theorem}
\label{thm.dual.specht-Slam-self}Let $\lambda$ be a partition of $n$. Assume
that the integer $h^{\lambda}$ (defined in Definition \ref{def.spechtmod.flam}%
) is invertible in $\mathbf{k}$. Then, $\left(  \mathcal{S}^{\lambda}\right)
^{\ast}\cong\mathcal{S}^{\lambda}$ as $S_{n}$-representations.
\end{theorem}

\begin{proof}
Recall that $\mathcal{S}^{\lambda}=\mathcal{S}^{Y\left(  \lambda\right)  }$.
Likewise, $\mathcal{S}^{\lambda^{t}}=\mathcal{S}^{Y\left(  \lambda^{t}\right)
}=\mathcal{S}^{\mathbf{r}\left(  Y\left(  \lambda\right)  \right)  }$ (since
the definition of the conjugate partition $\lambda^{t}$ yields $Y\left(
\lambda^{t}\right)  =\mathbf{r}\left(  Y\left(  \lambda\right)  \right)  $).

But $Y\left(  \lambda\right)  =Y\left(  \lambda/\varnothing\right)  $ is a
skew Young diagram of size $\left\vert Y\left(  \lambda\right)  \right\vert
=\left\vert \lambda\right\vert =n$ (since $\lambda$ is a partition of $n$).
Hence, Theorem \ref{thm.dual.specht-skew} (applied to $D=Y\left(
\lambda\right)  $) yields%
\[
\left(  \mathcal{S}^{Y\left(  \lambda\right)  }\right)  ^{\ast}\cong\left(
\mathcal{S}^{\mathbf{r}\left(  Y\left(  \lambda\right)  \right)  }\right)
^{\operatorname*{sign}}\text{ as }S_{n}\text{-representations.}%
\]
In other words,
\[
\left(  \mathcal{S}^{\lambda}\right)  ^{\ast}\cong\left(  \mathcal{S}%
^{\lambda^{t}}\right)  ^{\operatorname*{sign}}\text{ as }S_{n}%
\text{-representations}%
\]
(since $\mathcal{S}^{\lambda}=\mathcal{S}^{Y\left(  \lambda\right)  }$ and
$\mathcal{S}^{\lambda^{t}}=\mathcal{S}^{\mathbf{r}\left(  Y\left(
\lambda\right)  \right)  }$). But Theorem \ref{thm.sign-twist.Slam} yields
$\left(  \mathcal{S}^{\lambda}\right)  ^{\operatorname*{sign}}\cong%
\mathcal{S}^{\lambda^{t}}$. Thus, $\mathcal{S}^{\lambda^{t}}\cong\left(
\mathcal{S}^{\lambda}\right)  ^{\operatorname*{sign}}$.

By Corollary \ref{cor.sign-twist.iso}, this entails $\left(  \mathcal{S}%
^{\lambda^{t}}\right)  ^{\operatorname*{sign}}\cong\left(  \left(
\mathcal{S}^{\lambda}\right)  ^{\operatorname*{sign}}\right)
^{\operatorname*{sign}}=\mathcal{S}^{\lambda}$ (by Proposition
\ref{prop.sign-twist.twice-back}). Altogether, we thus find $\left(
\mathcal{S}^{\lambda}\right)  ^{\ast}\cong\left(  \mathcal{S}^{\lambda^{t}%
}\right)  ^{\operatorname*{sign}}\cong\mathcal{S}^{\lambda}$. This proves
Theorem \ref{thm.dual.specht-Slam-self}.
\end{proof}

\begin{corollary}
[self-duality theorem for $S_{n}$-representations]\label{cor.dual.char0-self}%
Let $\mathbf{k}$ be a field of characteristic $0$. Let $V$ be any $S_{n}%
$-representation that is finite-dimensional as a $\mathbf{k}$-vector space.
Then, $V^{\ast}\cong V$ as $S_{n}$-representations.
\end{corollary}

\begin{proof}
[Proof sketch.]We shall say that an $S_{n}$-representation $M$ is
\emph{self-dual} if it satisfies $M^{\ast}\cong M$. Thus, our goal is to prove
that $V$ is self-dual.

Clearly, if $M$ is a self-dual $S_{n}$-representation, then any $S_{n}%
$-representation $N$ that is isomorphic to $M$ must also be self-dual (since
$N\cong M$ entails $N^{\ast}\cong M^{\ast}\cong M\cong N$). In other words,
self-duality is preserved under isomorphism.

If $\lambda$ is any partition of $n$, then the positive integer $h^{\lambda}$
(defined in Definition \ref{def.spechtmod.flam}) is invertible in $\mathbf{k}$
(since $\mathbf{k}$ has characteristic $0$), and thus Theorem
\ref{thm.dual.specht-Slam-self} shows that $\left(  \mathcal{S}^{\lambda
}\right)  ^{\ast}\cong\mathcal{S}^{\lambda}$; but this is just saying that the
Specht module $\mathcal{S}^{\lambda}$ is self-dual (as an $S_{n}%
$-representation). In other words, we have shown the following claim:

\begin{statement}
\textit{Claim 1:} If $\lambda$ is any partition of $n$, then the Specht module
$\mathcal{S}^{\lambda}$ is self-dual.
\end{statement}

On the other hand, if $P$ and $Q$ are two self-dual $S_{n}$-representations,
then Proposition \ref{prop.dual.dirsums} yields%
\[
\left(  P\oplus Q\right)  ^{\ast}\cong\underbrace{P^{\ast}}_{\substack{\cong
P\\\text{(since }P\text{ is self-dual)}}}\oplus\underbrace{Q^{\ast}%
}_{\substack{\cong Q\\\text{(since }Q\text{ is self-dual)}}}\cong P\oplus Q,
\]
which shows that $P\oplus Q$ is again self-dual. In other words, any direct
sum of two self-dual $S_{n}$-representations is again self-dual. Hence, by
induction, we conclude that any direct sum of finitely many self-dual $S_{n}%
$-representations is again self-dual.

But Corollary \ref{cor.spechtmod.complete2} shows that $V$ is isomorphic to a
direct sum of finitely many copies of Specht modules of the form
$\mathcal{S}^{\lambda}$ with $\lambda$ being partitions of $n$. Each of these
Specht modules $\mathcal{S}^{\lambda}$ is self-dual (by Claim 1). Hence, their
direct sum is self-dual again (since any direct sum of finitely many self-dual
$S_{n}$-representations is again self-dual). In other words, $V$ is self-dual
(since $V$ is isomorphic to this direct sum, but self-duality is preserved
under isomorphism). This proves Corollary \ref{cor.dual.char0-self}.
\end{proof}

\begin{remark}
Corollary \ref{cor.dual.char0-self} would fail if $V$ was allowed to be
infinite-dimensional (here, $V^{\ast}\cong V$ will often fail even on the
level of $\mathbf{k}$-vector spaces). It can also fail if $\mathbf{k}$ is not
a field but merely a commutative $\mathbb{Q}$-algebra.
\end{remark}

\begin{remark}
Corollary \ref{cor.dual.char0-self} can be generalized: Let $G$ be a finite
group such that each $g\in G$ is conjugate to its inverse $g^{-1}$. (Such
groups $G$ are said to be \emph{ambivalent}.) Let $\mathbf{k}$ be a field of
characteristic $0$. Let $V$ be any $G$-representation that is
finite-dimensional as a $\mathbf{k}$-vector space. Then, $V^{\ast}\cong V$.

This result can be applied to $G=S_{n}$ because Theorem
\ref{thm.antip-conj.perm} shows that each $g\in S_{n}$ is conjugate to
$g^{-1}$. A proof of this generalization is implicit in \cite[proof of Lemma
3.24]{Lorenz18}.
\end{remark}

\begin{exercise}
\fbox{3} A bilinear form $f:V\times V\rightarrow\mathbf{k}$ for a $\mathbf{k}%
$-module $V$ is said to be \emph{symmetric} if all $u,v\in V$ satisfy
$f\left(  u,v\right)  =f\left(  v,u\right)  $.

Let $\lambda$ be a partition of $n$. Assume that $h^{\lambda}$ is invertible
in $\mathbf{k}$. Prove that there exists a symmetric $S_{n}$-invariant
bilinear form $f:\mathcal{S}^{\lambda}\times\mathcal{S}^{\lambda}%
\rightarrow\mathbf{k}$ such that the map $f_{R}:\mathcal{S}^{\lambda
}\rightarrow\left(  \mathcal{S}^{\lambda}\right)  ^{\ast}$ is bijective.
\end{exercise}

\subsubsection{Characters of duals}

Recall the notion of a character (Definition \ref{def.char.char}). We can
easily describe the character of the dual of a $G$-representation in general:

\begin{proposition}
\label{prop.dual.char-gen}Let $G$ be a group. Let $V$ be a representation of
$G$ over $\mathbf{k}$. Assume that $V$ is free of finite rank (i.e., has a
finite basis) as a $\mathbf{k}$-module. Then, the characters $\chi_{V}$ and
$\chi_{V^{\ast}}$ of $V$ and $V^{\ast}$ are related by the following formula:%
\[
\chi_{V^{\ast}}\left(  g\right)  =\chi_{V}\left(  g^{-1}\right)
\ \ \ \ \ \ \ \ \ \ \text{for each }g\in G\text{.}%
\]

\end{proposition}

For $G=S_{n}$, we can simplify this by replacing $g^{-1}$ by $g$:

\begin{corollary}
\label{cor.dual.char-Sn}Let $V$ be a representation of $S_{n}$ over
$\mathbf{k}$. Assume that $V$ is free of finite rank (i.e., has a finite
basis) as a $\mathbf{k}$-module. Then, the characters $\chi_{V}$ and
$\chi_{V^{\ast}}$ of $V$ and $V^{\ast}$ are equal: i.e., we have $\chi
_{V}=\chi_{V^{\ast}}$.
\end{corollary}

\begin{exercise}
\fbox{2} Prove Proposition \ref{prop.dual.char-gen} and Corollary
\ref{cor.dual.char-Sn}.
\end{exercise}

\begin{exercise}
\fbox{3} Many texts on representation theory restrict themselves to the case
when $\mathbf{k}=\mathbb{C}$ and the group $G$ is finite. In such texts, a
version of Proposition \ref{prop.dual.char-gen} (for $\mathbf{k}=\mathbb{C}$
and for $G$ finite) is often found with the equality $\chi_{V^{\ast}}\left(
g\right)  =\chi_{V}\left(  g^{-1}\right)  $ replaced by $\chi_{V^{\ast}%
}\left(  g\right)  =\overline{\chi_{V}\left(  g\right)  }$, where
$\overline{z}$ denotes the complex conjugate of a complex number $z$.

Derive this version from Proposition \ref{prop.dual.char-gen}. \medskip

[\textbf{Hint:} Argue that the endomorphism $\tau_{g}$ of $V$ has finite order
in the group $\operatorname*{GL}\left(  V\right)  $, and thus all its
eigenvalues lie on the unit circle.]
\end{exercise}

\subsubsection{Proving Theorem \ref{thm.sign-twist.D}}

We can now prove Theorem \ref{thm.sign-twist.D}, which we left unproved in the
previous section:

\begin{proof}
[Proof of Theorem \ref{thm.sign-twist.D}.]Theorem \ref{thm.dual.specht-field}
shows that $\left(  \mathcal{S}^{D}\right)  ^{\ast}\cong\left(  \mathcal{S}%
^{\mathbf{r}\left(  D\right)  }\right)  ^{\operatorname*{sign}}$ as $S_{n}%
$-representations. But Corollary \ref{cor.dual.char0-self} (applied to
$V=\mathcal{S}^{D}$) yields $\left(  \mathcal{S}^{D}\right)  ^{\ast}%
\cong\mathcal{S}^{D}$. Hence, $\mathcal{S}^{D}\cong\left(  \mathcal{S}%
^{D}\right)  ^{\ast}\cong\left(  \mathcal{S}^{\mathbf{r}\left(  D\right)
}\right)  ^{\operatorname*{sign}}$. By Corollary \ref{cor.sign-twist.iso},
this entails%
\[
\left(  \mathcal{S}^{D}\right)  ^{\operatorname*{sign}}\cong\left(  \left(
\mathcal{S}^{\mathbf{r}\left(  D\right)  }\right)  ^{\operatorname*{sign}%
}\right)  ^{\operatorname*{sign}}=\mathcal{S}^{\mathbf{r}\left(  D\right)
}\ \ \ \ \ \ \ \ \ \ \left(  \text{by Proposition
\ref{prop.sign-twist.twice-back}}\right)  .
\]
This proves Theorem \ref{thm.sign-twist.D}.
\end{proof}

\subsection{\label{sec.rep.antip-conj}Application: Proof of antipodal
conjugacy}

With the Artin--Wedderburn theorem (Theorem \ref{thm.specht.AW}) and the
self-duality theorem (Corollary \ref{cor.dual.char0-self}), we now have most
of the tools needed to prove Theorem \ref{thm.antip-conj.linear}. Thus we will
give the proof in the present section, albeit somewhat tersely.

\subsubsection{Conjugacy and the Taussky--Zassenhaus theorem}

The proof requires a few auxiliary facts. One is a linear-algebraic result by
Taussky and Zassenhaus from 1959 \cite[Theorem 1]{TauZas59}:\footnote{See
\url{https://math.stackexchange.com/questions/62497/matrix-is-conjugate-to-its-own-transpose}
or \cite[Corollary 45.14]{Elman22} for a short proof using the Smith normal
form over the polynomial ring $\mathbf{k}\left[  x\right]  $.}

\begin{theorem}
[Taussky--Zassenhaus theorem]\label{thm.tau-zas}Assume that $\mathbf{k}$ is a
field. Let $A\in\mathbf{k}^{m\times m}$ be any square matrix over $\mathbf{k}%
$. Then, the matrix $A$ is similar (i.e., conjugate in the matrix ring
$\mathbf{k}^{m\times m}$) to its transpose $A^{T}$.
\end{theorem}

Note that this theorem would fail for $\mathbf{k}=\mathbb{Z}$ (see
\url{https://math.stackexchange.com/questions/1732276/is-a-matrix-over-a-pid-similar-to-its-transpose}
for a counterexample), so the \textquotedblleft$\mathbf{k}$ is a
field\textquotedblright\ assumption is important.

We can give the Taussky--Zassenhaus theorem a more conceptual formulation by
replacing matrices by linear maps. Similarity is then replaced by a slightly
generalized notion of conjugacy, defined as follows:

\begin{definition}
\label{def.conjugacy-of-ends}Let $U$ and $V$ be two $\mathbf{k}$-modules.
Then, two endomorphisms $f\in\operatorname*{End}U$ and $g\in
\operatorname*{End}V$ will be called \emph{conjugate} (to each other) if there
exists some $\mathbf{k}$-module isomorphism $\phi:U\rightarrow V$ satisfying
$g=\phi\circ f\circ\phi^{-1}$ (that is, $g\circ\phi=\phi\circ f$). We will
write \textquotedblleft$f\sim g$\textquotedblright\ for this.
\end{definition}

Note that this notion of conjugacy generalizes the usual concept of conjugacy
in an endomorphism ring. Indeed, if $V$ is any $\mathbf{k}$-module, then two
endomorphisms $f,g\in\operatorname*{End}V$ satisfy $f\sim g$ (in the sense of
Definition \ref{def.conjugacy-of-ends}) if and only if $f$ and $g$ are
conjugate in the ring $\operatorname*{End}V$ (because a $\mathbf{k}$-module
isomorphism $\phi:V\rightarrow V$ is the same thing as an invertible element
$\phi$ of $\operatorname*{End}V$).

It is easy to see that conjugacy of endomorphisms is an equivalence relation.
In particular, it is symmetric (i.e., if $\alpha\sim\beta$, then $\beta
\sim\alpha$) and transitive (i.e., if $\alpha\sim\beta$ and $\beta\sim\gamma$,
then $\alpha\sim\gamma$).

Conjugacy of endomorphisms is the coordinate-free version of similarity of matrices:

\begin{proposition}
\label{prop.conjugacy-of-ends.by-matrices}Let $U$ and $V$ be two free
$\mathbf{k}$-modules with finite bases $\left(  u_{1},u_{2},\ldots
,u_{m}\right)  $ and $\left(  v_{1},v_{2},\ldots,v_{m}\right)  $,
respectively. Let $f\in\operatorname*{End}U$ and $g\in\operatorname*{End}V$ be
two endomorphisms. Let $A\in\mathbf{k}^{m\times m}$ be the matrix representing
the endomorphism $f$ with respect to the basis $\left(  u_{1},u_{2}%
,\ldots,u_{m}\right)  $ of $U$. Let $B\in\mathbf{k}^{m\times m}$ be the matrix
representing the endomorphism $g$ with respect to the basis $\left(
v_{1},v_{2},\ldots,v_{m}\right)  $ of $V$. Then, the conjugacy $f\sim g$ holds
if and only if the matrices $A$ and $B$ are similar (i.e., conjugate in the
matrix ring $\mathbf{k}^{m\times m}$).
\end{proposition}

\begin{fineprint}
\begin{proof}
[Proof sketch.]This is basic linear algebra. For the \textquotedblleft%
$\Longrightarrow$\textquotedblright\ direction, assume that $f\sim g$ holds.
Thus, there exists some $\mathbf{k}$-module isomorphism $\phi:U\rightarrow V$
satisfying $g=\phi\circ f\circ\phi^{-1}$. Consider this $\phi$, and let
$P\in\mathbf{k}^{m\times m}$ be the matrix representing it with respect to the
bases $\left(  u_{1},u_{2},\ldots,u_{m}\right)  $ and $\left(  v_{1}%
,v_{2},\ldots,v_{m}\right)  $ of $U$ and $V$. Then, this matrix $P$ is
invertible (since it represents the invertible map $\phi$) and satisfies
$B=PAP^{-1}$ (since $g=\phi\circ f\circ\phi^{-1}$). Hence, the matrices $A$
and $B$ are similar. This proves the \textquotedblleft$\Longrightarrow
$\textquotedblright\ direction of Proposition
\ref{prop.conjugacy-of-ends.by-matrices}. To prove the \textquotedblleft%
$\Longleftarrow$\textquotedblright\ direction, reverse the argument (recalling
that any $m\times m$-matrix represents some linear map $\phi:U\rightarrow V$
with respect to the bases $\left(  u_{1},u_{2},\ldots,u_{m}\right)  $ and
$\left(  v_{1},v_{2},\ldots,v_{m}\right)  $).
\end{proof}
\end{fineprint}

We can now restate Theorem \ref{thm.tau-zas} as follows:

\begin{corollary}
[Taussky--Zassenhaus theorem in coordinate-free form]\label{cor.tau-zas.end}%
Assume that $\mathbf{k}$ is a field. Let $V$ be a finite-dimensional
$\mathbf{k}$-vector space. Let $f\in\operatorname*{End}V$ be any endomorphism
of $V$. Then, the endomorphism $f\in\operatorname*{End}V$ is conjugate to its
dual $f^{\ast}\in\operatorname*{End}\left(  V^{\ast}\right)  $.
\end{corollary}

\begin{proof}
[Proof sketch.]Fix a basis $\left(  v_{1},v_{2},\ldots,v_{m}\right)  $ of $V$.
Let $\left(  v_{1}^{\ast},v_{2}^{\ast},\ldots,v_{m}^{\ast}\right)  $ be its
dual basis, which of course is a basis of $V^{\ast}$. Let $A$ be the matrix
representing the endomorphism $f$ with respect to the basis $\left(
v_{1},v_{2},\ldots,v_{m}\right)  $ of $V$. Then, the matrix representing its
dual $f^{\ast}$ with respect to the dual basis $\left(  v_{1}^{\ast}%
,v_{2}^{\ast},\ldots,v_{m}^{\ast}\right)  $ is the transpose $A^{T}$ of $A$
(by Proposition \ref{prop.dual.functorial} \textbf{(d)}). Hence, the conjugacy
$f\sim f^{\ast}$ holds if and only if the matrices $A$ and $A^{T}$ are similar
(by Proposition \ref{prop.conjugacy-of-ends.by-matrices}). But Theorem
\ref{thm.tau-zas} shows that the matrices $A$ and $A^{T}$ are indeed similar.
Thus, $f\sim f^{\ast}$ holds. This proves Corollary \ref{cor.tau-zas.end}.
\end{proof}

Another source of conjugate endomorphisms are isomorphic representations of a
group or, more generally, isomorphic modules over a $\mathbf{k}$-algebra:

\begin{proposition}
\label{prop.conjugacy-of-ends.iso-reps}Let $R$ be a $\mathbf{k}$-algebra. Let
$U$ and $V$ be two isomorphic left $R$-modules. Let $\mathbf{a}\in R$. Then,%
\[
\left(  \tau_{\mathbf{a}}\text{ on }U\right)  \sim\left(  \tau_{\mathbf{a}%
}\text{ on }V\right)  .
\]
Here, \textquotedblleft$\tau_{\mathbf{a}}$ on $U$\textquotedblright\ means the
$\mathbf{k}$-linear map $U\rightarrow U,\ u\mapsto\mathbf{a}\cdot u$ (viewed
as an endomorphism of the $\mathbf{k}$-module $U$), whereas \textquotedblleft%
$\tau_{\mathbf{a}}$ on $V$\textquotedblright\ means the $\mathbf{k}$-linear
map $V\rightarrow V,\ v\mapsto\mathbf{a}\cdot v$ (viewed as an endomorphism of
the $\mathbf{k}$-module $V$).
\end{proposition}

\begin{fineprint}
\begin{proof}
There exists a left $R$-module isomorphism $\phi:U\rightarrow V$ (since $U$
and $V$ are isomorphic). Consider this $\phi$. Thus,%
\begin{equation}
\phi\left(  \mathbf{a}\cdot u\right)  =\mathbf{a}\cdot\phi\left(  u\right)
\ \ \ \ \ \ \ \ \ \ \text{for each }u\in U
\label{pf.prop.conjugacy-of-ends.iso-reps.1}%
\end{equation}
(since $\phi$ is a left $R$-module morphism). In other words,%
\[
\phi\circ\left(  \tau_{\mathbf{a}}\text{ on }U\right)  =\left(  \tau
_{\mathbf{a}}\text{ on }V\right)  \circ\phi
\]
\footnote{\textit{Proof.} Let $u\in U$. Then,%
\[
\left(  \phi\circ\left(  \tau_{\mathbf{a}}\text{ on }U\right)  \right)
\left(  u\right)  =\phi\left(  \left(  \tau_{\mathbf{a}}\text{ on }U\right)
\left(  u\right)  \right)  =\phi\left(  \mathbf{a}\cdot u\right)
\]
(since the definition of \textquotedblleft$\tau_{\mathbf{a}}$ on
$U$\textquotedblright\ yields $\left(  \tau_{\mathbf{a}}\text{ on }U\right)
\left(  u\right)  =\mathbf{a}\cdot u$) and%
\[
\left(  \left(  \tau_{\mathbf{a}}\text{ on }V\right)  \circ\phi\right)
\left(  u\right)  =\left(  \tau_{\mathbf{a}}\text{ on }V\right)  \left(
\phi\left(  u\right)  \right)  =\mathbf{a}\cdot\phi\left(  u\right)
\]
(by the definition of \textquotedblleft$\tau_{\mathbf{a}}$ on $V$%
\textquotedblright). The right hand sides of these two equalities are equal
(by (\ref{pf.prop.conjugacy-of-ends.iso-reps.1})). Hence, so are their left
hand sides. In other words, $\left(  \phi\circ\left(  \tau_{\mathbf{a}}\text{
on }U\right)  \right)  \left(  u\right)  =\left(  \left(  \tau_{\mathbf{a}%
}\text{ on }V\right)  \circ\phi\right)  \left(  u\right)  $.
\par
Having proved this equality for all $u\in U$, we thus conclude that $\phi
\circ\left(  \tau_{\mathbf{a}}\text{ on }U\right)  =\left(  \tau_{\mathbf{a}%
}\text{ on }V\right)  \circ\phi$.}. Since $\phi$ is invertible, we can solve
this equality for $\left(  \tau_{\mathbf{a}}\text{ on }V\right)  $, obtaining%
\[
\left(  \tau_{\mathbf{a}}\text{ on }V\right)  =\phi\circ\left(  \tau
_{\mathbf{a}}\text{ on }U\right)  \circ\phi^{-1}.
\]
Since $\phi$ is a $\mathbf{k}$-module isomorphism (being an $R$-module
isomorphism), we thus conclude that $\left(  \tau_{\mathbf{a}}\text{ on
}U\right)  \sim\left(  \tau_{\mathbf{a}}\text{ on }V\right)  $ (by Definition
\ref{def.conjugacy-of-ends}). This proves Proposition
\ref{prop.conjugacy-of-ends.iso-reps}.
\end{proof}
\end{fineprint}

\subsubsection{Using Artin--Wedderburn}

Our next lemma is a simple consequence of the Artin--Wedderburn theorem:

\begin{lemma}
\label{lem.antip-conj.linear.2}Assume that $n!$ is invertible in $\mathbf{k}$.
Let $\mathbf{a},\mathbf{b}\in\mathbf{k}\left[  S_{n}\right]  $. Assume that
\begin{equation}
\rho_{\lambda}\left(  \mathbf{a}\right)  \text{ is conjugate to }\rho
_{\lambda}\left(  \mathbf{b}\right)  \text{ in }\operatorname*{End}\left(
\mathcal{S}^{\lambda}\right)  \label{eq.lem.antip-conj.linear.2.ass}%
\end{equation}
for each partition $\lambda$ of $n$, where the morphism $\rho_{\lambda}$ is
defined in Theorem \ref{thm.specht.AW}. Then, $\mathbf{a}$ is conjugate to
$\mathbf{b}$ in $\mathbf{k}\left[  S_{n}\right]  $.
\end{lemma}

\begin{proof}
Let $P$ be the set of all partitions of $n$. Consider the map%
\begin{align*}
\operatorname*{AH}:\mathbf{k}\left[  S_{n}\right]   &  \rightarrow
\underbrace{\prod_{\lambda\in P}\operatorname*{End}\nolimits_{\mathbf{k}%
}\left(  \mathcal{S}^{\lambda}\right)  }_{\text{a direct product of
}\mathbf{k}\text{-algebras}},\\
\mathbf{c}  &  \mapsto\left(  \rho_{\lambda}\left(  \mathbf{c}\right)
\right)  _{\lambda\in P}.
\end{align*}
By Theorem \ref{thm.specht.AW}, this map $\operatorname*{AH}$ is a
$\mathbf{k}$-algebra isomorphism\footnote{The notation in Theorem
\ref{thm.specht.AW} differs slightly from the one we are using here, in that
we wrote \textquotedblleft$\lambda$ is a partition of $n$\textquotedblright%
\ instead of the short notation \textquotedblleft$\lambda\in P$%
\textquotedblright, and in that we used the letter $\mathbf{a}$ (rather than
$\mathbf{c}$) for the variable.}. Hence, in particular, this map
$\operatorname*{AH}$ is injective, and its inverse $\operatorname*{AH}%
\nolimits^{-1}$ is again a $\mathbf{k}$-algebra isomorphism.

Now, let $\lambda\in P$. Thus, $\lambda$ is a partition of $n$ (by the
definition of $P$). Hence, $\rho_{\lambda}\left(  \mathbf{a}\right)  $ is
conjugate to $\rho_{\lambda}\left(  \mathbf{b}\right)  $ in
$\operatorname*{End}\left(  \mathcal{S}^{\lambda}\right)  $ (by our assumption
(\ref{eq.lem.antip-conj.linear.2.ass})). In other words, there exists some
invertible element $\psi_{\lambda}\in\operatorname*{End}\left(  \mathcal{S}%
^{\lambda}\right)  $ such that%
\begin{equation}
\rho_{\lambda}\left(  \mathbf{a}\right)  =\psi_{\lambda}\cdot\rho_{\lambda
}\left(  \mathbf{b}\right)  \cdot\psi_{\lambda}^{-1}
\label{pf.lem.antip-conj.linear.2.fn1.1}%
\end{equation}
(by the definition of \textquotedblleft conjugate\textquotedblright\ in a
ring\footnote{Note that the multiplication $\cdot$ on the right hand side of
(\ref{pf.lem.antip-conj.linear.2.fn1.1}) is composition of endomorphisms, so
we could just as well write $\circ$ for it.}). Consider this $\psi_{\lambda}$.

Forget that we fixed $\lambda$. We thus have constructed an invertible element
$\psi_{\lambda}\in\operatorname*{End}\left(  \mathcal{S}^{\lambda}\right)  $
satisfying (\ref{pf.lem.antip-conj.linear.2.fn1.1}) for each $\lambda\in P$.

Next, define the family%
\[
\psi:=\left(  \psi_{\lambda}\right)  _{\lambda\in P}\in\prod_{\lambda\in
P}\operatorname*{End}\nolimits_{\mathbf{k}}\left(  \mathcal{S}^{\lambda
}\right)  .
\]
Then, $\psi$ is a family of invertible elements (since $\psi_{\lambda}$ is
invertible for each $\lambda\in P$), and thus is itself invertible (since the
multiplication in the direct product $\prod_{\lambda\in P}\operatorname*{End}%
\nolimits_{\mathbf{k}}\left(  \mathcal{S}^{\lambda}\right)  $ is entrywise),
with inverse%
\[
\psi^{-1}=\left(  \psi_{\lambda}^{-1}\right)  _{\lambda\in P}%
\]
(since $\psi=\left(  \psi_{\lambda}\right)  _{\lambda\in P}$ and again since
the multiplication is entrywise).

Define an element $\mathbf{c}\in\mathbf{k}\left[  S_{n}\right]  $ by
$\mathbf{c}:=\operatorname*{AH}\nolimits^{-1}\left(  \psi\right)  $. Then,
$\mathbf{c}$ is the image of the invertible element $\psi$ under the
$\mathbf{k}$-algebra isomorphism $\operatorname*{AH}\nolimits^{-1}$;
therefore, $\mathbf{c}$ itself is invertible (since $\mathbf{k}$-algebra
isomorphisms send invertible elements to invertible elements). From
$\mathbf{c}=\operatorname*{AH}\nolimits^{-1}\left(  \psi\right)  $, we obtain%
\[
\operatorname*{AH}\left(  \mathbf{c}\right)  =\psi=\left(  \psi_{\lambda
}\right)  _{\lambda\in P}.
\]

Moreover, since $\operatorname*{AH}$ is a $\mathbf{k}$-algebra morphism, we
have%
\[
\operatorname*{AH}\left(  \mathbf{c}^{-1}\right)  =\left(
\underbrace{\operatorname*{AH}\left(  \mathbf{c}\right)  }_{=\psi}\right)
^{-1}=\psi^{-1}=\left(  \psi_{\lambda}^{-1}\right)  _{\lambda\in P}.
\]
Since $\operatorname*{AH}$ is a $\mathbf{k}$-algebra morphism, we furthermore
have%
\begin{align*}
&  \operatorname*{AH}\left(  \mathbf{cbc}^{-1}\right) \\
&  =\operatorname*{AH}\left(  \mathbf{c}\right)  \cdot\operatorname*{AH}%
\left(  \mathbf{b}\right)  \cdot\operatorname*{AH}\left(  \mathbf{c}%
^{-1}\right) \\
&  =\left(  \psi_{\lambda}\right)  _{\lambda\in P}\cdot\left(  \rho_{\lambda
}\left(  \mathbf{b}\right)  \right)  _{\lambda\in P}\cdot\left(  \psi
_{\lambda}^{-1}\right)  _{\lambda\in P}\\
&  \ \ \ \ \ \ \ \ \ \ \ \ \ \ \ \ \ \ \ \ \left(
\begin{array}
[c]{c}%
\text{since }\operatorname*{AH}\left(  \mathbf{c}\right)  =\left(
\psi_{\lambda}\right)  _{\lambda\in P}\text{ and }\operatorname*{AH}\left(
\mathbf{b}\right)  =\left(  \rho_{\lambda}\left(  \mathbf{b}\right)  \right)
_{\lambda\in P}\\
\text{(by the definition of }\operatorname*{AH}\text{) and }\operatorname*{AH}%
\left(  \mathbf{c}^{-1}\right)  =\left(  \psi_{\lambda}^{-1}\right)
_{\lambda\in P}%
\end{array}
\right) \\
&  =\left(  \psi_{\lambda}\cdot\rho_{\lambda}\left(  \mathbf{b}\right)
\cdot\psi_{\lambda}^{-1}\right)  _{\lambda\in P}\\
&  \ \ \ \ \ \ \ \ \ \ \ \ \ \ \ \ \ \ \ \ \left(  \text{since the
multiplication in }\prod_{\lambda\in P}\operatorname*{End}%
\nolimits_{\mathbf{k}}\left(  \mathcal{S}^{\lambda}\right)  \text{ is
entrywise}\right) \\
&  =\left(  \rho_{\lambda}\left(  \mathbf{a}\right)  \right)  _{\lambda\in
P}\ \ \ \ \ \ \ \ \ \ \left(
\begin{array}
[c]{c}%
\text{since (\ref{pf.lem.antip-conj.linear.2.fn1.1}) yields }\psi_{\lambda
}\cdot\rho_{\lambda}\left(  \mathbf{b}\right)  \cdot\psi_{\lambda}^{-1}%
=\rho_{\lambda}\left(  \mathbf{a}\right) \\
\text{for each }\lambda\in P
\end{array}
\right) \\
&  =\operatorname*{AH}\left(  \mathbf{a}\right)  \ \ \ \ \ \ \ \ \ \ \left(
\text{by the definition of }\operatorname*{AH}\right)  .
\end{align*}
Since $\operatorname*{AH}$ is injective, we thus conclude that $\mathbf{cbc}%
^{-1}=\mathbf{a}$. In other words, $\mathbf{a}=\mathbf{cbc}^{-1}$. Hence, the
elements $\mathbf{a}$ and $\mathbf{b}$ of $\mathbf{k}\left[  S_{n}\right]  $
are conjugate. Lemma \ref{lem.antip-conj.linear.2} is thus proved.
\end{proof}

\subsubsection{The proof}

We can now finally prove Theorem \ref{thm.antip-conj.linear}:

\begin{proof}
[Proof of Theorem \ref{thm.antip-conj.linear} (sketched).]We must show that
$\mathbf{a}$ is conjugate to $S\left(  \mathbf{a}\right)  $ in $\mathbf{k}%
\left[  S_{n}\right]  $. We know that $n!$ is invertible in $\mathbf{k}$
(since $\mathbf{k}$ is a field of characteristic $0$). Hence, by Lemma
\ref{lem.antip-conj.linear.2} (applied to $\mathbf{b}=S\left(  \mathbf{a}%
\right)  $), it suffices to show that $\rho_{\lambda}\left(  \mathbf{a}%
\right)  $ is conjugate to $\rho_{\lambda}\left(  S\left(  \mathbf{a}\right)
\right)  $ in $\operatorname*{End}\left(  \mathcal{S}^{\lambda}\right)  $ for
each partition $\lambda$ of $n$.

So let $\lambda$ be any partition of $n$. We must prove that $\rho_{\lambda
}\left(  \mathbf{a}\right)  $ is conjugate to $\rho_{\lambda}\left(  S\left(
\mathbf{a}\right)  \right)  $ in $\operatorname*{End}\left(  \mathcal{S}%
^{\lambda}\right)  $.

If $G$ is any group, if $\mathbf{b}\in\mathbf{k}\left[  G\right]  $ is any
element, and if $U$ is any left $\mathbf{k}\left[  G\right]  $-module, then we
let \textquotedblleft$\tau_{\mathbf{b}}$ on $U$\textquotedblright\ denote the
$\mathbf{k}$-linear map $U\rightarrow U,\ u\mapsto\mathbf{b}\cdot u$ (viewed
as an endomorphism of the $\mathbf{k}$-module $U$). This generalizes the
notations \textquotedblleft$\tau_{\mathbf{a}}$ on $V^{\ast}$\textquotedblright%
\ and \textquotedblleft$\tau_{S\left(  \mathbf{a}\right)  }$ on $V$%
\textquotedblright\ from Proposition \ref{prop.dual.dual-map} \textbf{(b)}.

Note that the endomorphism $\rho_{\lambda}\left(  \mathbf{a}\right)  $ of
$\mathcal{S}^{\lambda}$ is precisely the endomorphism \textquotedblleft%
$\tau_{\mathbf{a}}$ on $\mathcal{S}^{\lambda}$\textquotedblright\ (since both
endomorphisms send each $v\in\mathcal{S}^{\lambda}$ to $\mathbf{a}\cdot v$).
In other words,
\begin{equation}
\rho_{\lambda}\left(  \mathbf{a}\right)  =\left(  \tau_{\mathbf{a}}\text{ on
}\mathcal{S}^{\lambda}\right)  . \label{pf.thm.antip-conj.linear.rho1}%
\end{equation}
The same reasoning (applied to $S\left(  \mathbf{a}\right)  $ instead of
$\mathbf{a}$) shows that%
\begin{equation}
\rho_{\lambda}\left(  S\left(  \mathbf{a}\right)  \right)  =\left(
\tau_{S\left(  \mathbf{a}\right)  }\text{ on }\mathcal{S}^{\lambda}\right)  .
\label{pf.thm.antip-conj.linear.rho2}%
\end{equation}

Corollary \ref{cor.tau-zas.end} (applied to $V=\mathcal{S}^{\lambda}$ and
$f=\rho_{\lambda}\left(  S\left(  \mathbf{a}\right)  \right)  $) shows that
the endomorphism $\rho_{\lambda}\left(  S\left(  \mathbf{a}\right)  \right)  $
of $\mathcal{S}^{\lambda}$ is conjugate to its dual $\left(  \rho_{\lambda
}\left(  S\left(  \mathbf{a}\right)  \right)  \right)  ^{\ast}$. In other
words,%
\begin{equation}
\rho_{\lambda}\left(  S\left(  \mathbf{a}\right)  \right)  \sim\left(
\rho_{\lambda}\left(  S\left(  \mathbf{a}\right)  \right)  \right)  ^{\ast}.
\label{pf.thm.antip-conj.linear.rho3}%
\end{equation}

But Proposition \ref{prop.dual.dual-map} \textbf{(b)} (applied to $G=S_{n}$
and $V=\mathcal{S}^{\lambda}$) yields that
\[
\left(  \tau_{\mathbf{a}}\text{ on }\left(  \mathcal{S}^{\lambda}\right)
^{\ast}\right)  =\left(  \tau_{S\left(  \mathbf{a}\right)  }\text{ on
}\mathcal{S}^{\lambda}\right)  ^{\ast}=\left(  \rho_{\lambda}\left(  S\left(
\mathbf{a}\right)  \right)  \right)  ^{\ast}%
\]
(since (\ref{pf.thm.antip-conj.linear.rho2}) yields $\left(  \tau_{S\left(
\mathbf{a}\right)  }\text{ on }\mathcal{S}^{\lambda}\right)  =\rho_{\lambda
}\left(  S\left(  \mathbf{a}\right)  \right)  $).

However, the integer $h^{\lambda}$ (defined in Definition
\ref{def.spechtmod.flam}) is a positive integer, and thus is invertible in
$\mathbf{k}$ (again since $\mathbf{k}$ is a field of characteristic $0$).
Hence, Theorem \ref{thm.dual.specht-Slam-self} shows that the $S_{n}%
$-representations $\mathcal{S}^{\lambda}$ and $\left(  \mathcal{S}^{\lambda
}\right)  ^{\ast}$ are isomorphic. In other words, the left $\mathbf{k}\left[
S_{n}\right]  $-modules $\mathcal{S}^{\lambda}$ and $\left(  \mathcal{S}%
^{\lambda}\right)  ^{\ast}$ are isomorphic. Hence, Proposition
\ref{prop.conjugacy-of-ends.iso-reps} (applied to $R=\mathbf{k}\left[
S_{n}\right]  $ and $U=\mathcal{S}^{\lambda}$ and $V=\left(  \mathcal{S}%
^{\lambda}\right)  ^{\ast}$) yields%
\[
\left(  \tau_{\mathbf{a}}\text{ on }\mathcal{S}^{\lambda}\right)  \sim\left(
\tau_{\mathbf{a}}\text{ on }\left(  \mathcal{S}^{\lambda}\right)  ^{\ast
}\right)  .
\]

Now, we have the following chain of equivalences:%
\[
\rho_{\lambda}\left(  \mathbf{a}\right)  =\left(  \tau_{\mathbf{a}}\text{ on
}\mathcal{S}^{\lambda}\right)  \sim\left(  \tau_{\mathbf{a}}\text{ on }\left(
\mathcal{S}^{\lambda}\right)  ^{\ast}\right)  =\left(  \rho_{\lambda}\left(
S\left(  \mathbf{a}\right)  \right)  \right)  ^{\ast}\sim\rho_{\lambda}\left(
S\left(  \mathbf{a}\right)  \right)
\]
(by (\ref{pf.thm.antip-conj.linear.rho3})). Hence, $\rho_{\lambda}\left(
\mathbf{a}\right)  \sim\rho_{\lambda}\left(  S\left(  \mathbf{a}\right)
\right)  $ (since conjugacy of endomorphisms is an equivalence relation). In
other words, $\rho_{\lambda}\left(  \mathbf{a}\right)  $ is conjugate to
$\rho_{\lambda}\left(  S\left(  \mathbf{a}\right)  \right)  $ in
$\operatorname*{End}\left(  \mathcal{S}^{\lambda}\right)  $. As we explained,
this proves Theorem \ref{thm.antip-conj.linear}.
\end{proof}

Another proof of Theorem \ref{thm.antip-conj.linear} can be given using the
seminormal basis of $\mathbf{k}\left[  S_{n}\right]  $.

\subsection{The $\mathbf{F}_{U}w_{U,V}\mathbf{E}_{V}$ basis}

Throughout this section, we will use the following conventions:

\begin{convention}
\label{conv.specht.FwE.conv}\textbf{(a)} We let $\mathcal{A}=\mathbf{k}\left[
S_{n}\right]  $ as usual. \medskip

\textbf{(b)} If $P$ and $Q$ are two $n$-tableaux of the same shape $D$, then
the permutation $w_{P,Q}\in S_{n}$ is defined as in Definition
\ref{def.specht.wPQ}. (That is, it is the unique permutation $w\in S_{n}$ that
satisfies $wQ=P$.) \medskip

\textbf{(c)} If $T$ is an $n$-tableau of any shape, then we set%
\[
\mathbf{E}_{T}:=\nabla_{\operatorname*{Col}T}^{-}\nabla_{\operatorname*{Row}%
T}\in\mathcal{A}\ \ \ \ \ \ \ \ \ \ \text{and}\ \ \ \ \ \ \ \ \ \ \mathbf{F}%
_{T}:=\nabla_{\operatorname*{Row}T}\nabla_{\operatorname*{Col}T}^{-}%
\in\mathcal{A}.
\]
(See Proposition \ref{prop.specht.FT.basics} for the properties of
$\mathbf{F}_{T}$ in the case when $T$ is straight-shaped.)
\end{convention}

\subsubsection{The basis}

We shall next construct another basis of $\mathcal{A}=\mathbf{k}\left[
S_{n}\right]  $ when $n!$ is invertible in $\mathbf{k}$. While similar to
Young's symmetrizer basis $\left(  \mathbf{E}_{U,V}\right)  _{\lambda\text{ is
a partition of }n\text{, and }U,V\in\operatorname*{SYT}\left(  \lambda\right)
}$ from Corollary \ref{cor.specht.A.nat-basis}, it will interact with the
antipode $S$ in an even simpler way:

\begin{theorem}
\label{thm.specht.FwE.basis}For any partition $\lambda$ of $n$, let
$\operatorname*{SYT}\left(  \lambda\right)  $ be the set of all standard
$n$-tableaux of shape $Y\left(  \lambda\right)  $. Assume that $n!$ is
invertible in $\mathbf{k}$. Then, the family $\left(  \mathbf{F}_{U}%
w_{U,V}\mathbf{E}_{V}\right)  _{\lambda\text{ is a partition of }n\text{, and
}U,V\in\operatorname*{SYT}\left(  \lambda\right)  }$ is a basis of the
$\mathbf{k}$-module $\mathcal{A}$.
\end{theorem}

\begin{example}
\label{exa.specht.FwE.n=3}Let $n=3$, and consider the two $n$-tableaux
$U:=12\backslash\backslash3$ and $V:=13\backslash\backslash2$ of shape
$\left(  2,1\right)  $. Then,%
\begin{align*}
\underbrace{\mathbf{F}_{U}}_{=\nabla_{\operatorname*{Row}U}\nabla
_{\operatorname*{Col}U}^{-}}w_{U,V}\underbrace{\mathbf{E}_{V}}_{=\nabla
_{\operatorname*{Col}V}^{-}\nabla_{\operatorname*{Row}V}}  &
=\underbrace{\nabla_{\operatorname*{Row}U}}_{=1+t_{1,2}}\ \ \underbrace{\nabla
_{\operatorname*{Col}U}^{-}}_{=1-t_{1,3}}\ \ \underbrace{w_{U,V}}_{=t_{2,3}%
}\ \ \underbrace{\nabla_{\operatorname*{Col}V}^{-}}_{=1-t_{1,2}}%
\ \ \underbrace{\nabla_{\operatorname*{Row}V}}_{=1+t_{1,3}}\\
&  =-2-2s_{1}+4s_{2}-2t_{1,3}+4\operatorname*{cyc}\nolimits_{1,2,3}%
-\,2\operatorname*{cyc}\nolimits_{1,3,2}.
\end{align*}
In total, denoting the four standard $3$-tableaux of partition shapes by
\[
P:=123,\ \ \ \ \ \ \ \ \ \ U:=12\backslash\backslash
3,\ \ \ \ \ \ \ \ \ \ V:=13\backslash\backslash
2,\ \ \ \ \ \ \ \ \ \ N:=1\backslash\backslash2\backslash\backslash3,
\]
the family $\left(  \mathbf{F}_{U}w_{U,V}\mathbf{E}_{V}\right)  _{\lambda
\text{ is a partition of }n\text{, and }U,V\in\operatorname*{SYT}\left(
\lambda\right)  }$ looks as follows:%
\begin{align*}
\mathbf{F}_{P}w_{P,P}\mathbf{E}_{P}  &  =6\nabla=6+6s_{1}+6s_{2}%
+6t_{1,3}+6\operatorname*{cyc}\nolimits_{1,2,3}+\,6\operatorname*{cyc}%
\nolimits_{1,3,2};\\
\mathbf{F}_{U}w_{U,U}\mathbf{E}_{U}  &  =4+4s_{1}-2s_{2}-2t_{1,3}%
-2\operatorname*{cyc}\nolimits_{1,2,3}-\,2\operatorname*{cyc}\nolimits_{1,3,2}%
;\\
\mathbf{F}_{U}w_{U,V}\mathbf{E}_{V}  &  =-2-2s_{1}+4s_{2}-2t_{1,3}%
+4\operatorname*{cyc}\nolimits_{1,2,3}-\,2\operatorname*{cyc}\nolimits_{1,3,2}%
;\\
\mathbf{F}_{V}w_{V,U}\mathbf{E}_{U}  &  =-2-2s_{1}+4s_{2}-2t_{1,3}%
-2\operatorname*{cyc}\nolimits_{1,2,3}+\,4\operatorname*{cyc}\nolimits_{1,3,2}%
;\\
\mathbf{F}_{V}w_{V,V}\mathbf{E}_{V}  &  =4-2s_{1}-2s_{2}+4t_{1,3}%
-2\operatorname*{cyc}\nolimits_{1,2,3}-\,2\operatorname*{cyc}\nolimits_{1,3,2}%
;\\
\mathbf{F}_{N}w_{N,N}\mathbf{E}_{N}  &  =6\nabla^{-}=6-6s_{1}-6s_{2}%
-6t_{1,3}+6\operatorname*{cyc}\nolimits_{1,2,3}+\,6\operatorname*{cyc}%
\nolimits_{1,3,2}.
\end{align*}

\end{example}

\begin{noncompile}
The change-of-basis matrix to the standard basis $\left(  w\right)  _{w\in
S_{n}}$ is thus%
\[
\left(
\begin{array}
[c]{cccccc}%
6 & 6 & 6 & 6 & 6 & 6\\
4 & 4 & -2 & -2 & -2 & -2\\
-2 & -2 & 4 & -2 & 4 & -2\\
-2 & -2 & 4 & -2 & -2 & 4\\
4 & -2 & -2 & 4 & -2 & -2\\
6 & -6 & -6 & -6 & 6 & 6
\end{array}
\right)  .
\]
Its determinant is $-6^{7}$.
\end{noncompile}

The proof of Theorem \ref{thm.specht.FwE.basis} will require some work. First,
however, let us reveal how the antipode $S$ acts on this basis:

\begin{proposition}
\label{prop.specht.FwE.S}Let $\lambda$ be a partition of $n$. Let $U$ and $V$
be any two $n$-tableaux of shape $Y\left(  \lambda\right)  $. Then,%
\[
S\left(  \mathbf{F}_{U}w_{U,V}\mathbf{E}_{V}\right)  =\mathbf{F}_{V}%
w_{V,U}\mathbf{E}_{U}.
\]

\end{proposition}

\begin{proof}
Proposition \ref{prop.specht.FT.basics} \textbf{(a)} (applied to $T=V$) shows
that $S\left(  \mathbf{E}_{V}\right)  =\mathbf{F}_{V}$. Applying the map $S$
to this equality, we obtain $S\left(  S\left(  \mathbf{E}_{V}\right)  \right)
=S\left(  \mathbf{F}_{V}\right)  $. But $S\left(  S\left(  \mathbf{E}%
_{V}\right)  \right)  =\mathbf{E}_{V}$ (since Theorem \ref{thm.S.auto}
\textbf{(b)} says that $S\circ S=\operatorname*{id}$). Comparing these two
equalities, we find $S\left(  \mathbf{F}_{V}\right)  =\mathbf{E}_{V}$.
Applying the same argument to $U$ instead of $V$, we obtain $S\left(
\mathbf{F}_{U}\right)  =\mathbf{E}_{U}$.

Proposition \ref{prop.specht.wPQ.wQP} \textbf{(b)} (applied to $D=Y\left(
\lambda\right)  $ and $P=U$ and $Q=V$) yields $w_{V,U}=w_{U,V}^{-1}$. However,
the definition of the antipode $S$ yields $S\left(  w_{U,V}\right)
=w_{U,V}^{-1}$. Comparing these two equalities, we find $S\left(
w_{U,V}\right)  =w_{V,U}$.

We know from Theorem \ref{thm.S.auto} \textbf{(a)} that the antipode $S$ is an
algebra anti-morphism. Thus,%
\[
S\left(  \mathbf{F}_{U}w_{U,V}\mathbf{E}_{V}\right)  =\underbrace{S\left(
\mathbf{E}_{V}\right)  }_{=\mathbf{F}_{V}}\cdot\underbrace{S\left(
w_{U,V}\right)  }_{=w_{V,U}}\cdot\underbrace{S\left(  \mathbf{F}_{U}\right)
}_{=\mathbf{E}_{U}}=\mathbf{F}_{V}w_{V,U}\mathbf{E}_{U}.
\]
This proves Proposition \ref{prop.specht.FwE.S}.
\end{proof}

The following exercise can be used to speed up the computation of the family
$\left(  \mathbf{F}_{U}w_{U,V}\mathbf{E}_{V}\right)  _{\lambda\text{ is a
partition of }n\text{, and }U,V\in\operatorname*{SYT}\left(  \lambda\right)
}$, as it mostly reduces this computation to picking an arbitrary $n$-tableau
$T$ of shape $Y\left(  \lambda\right)  $ for each partition $\lambda$ of $n$
and computing the products $\nabla_{\operatorname*{Row}T}\nabla
_{\operatorname*{Col}T}^{-}\nabla_{\operatorname*{Row}T}$ for these chosen $T$'s:

\begin{exercise}
\fbox{2} Let $\lambda$ be a partition of $n$. Let $T$, $U$ and $V$ be any
three $n$-tableaux of shape $Y\left(  \lambda\right)  $. Prove that%
\[
\mathbf{F}_{U}w_{U,V}\mathbf{E}_{V}=\left\vert \mathcal{C}\left(  T\right)
\right\vert \cdot w_{U,T}\nabla_{\operatorname*{Row}T}\nabla
_{\operatorname*{Col}T}^{-}\nabla_{\operatorname*{Row}T}w_{T,V}.
\]
(See Definition \ref{def.tableau.cell-cts} for the meaning of $\mathcal{C}%
\left(  T\right)  $.)
\end{exercise}

\subsubsection{Proof of the basis property}

We now take aim at the proof of Theorem \ref{thm.specht.FwE.basis}. A few
simple lemmas pave the way:

\begin{lemma}
\label{lem.specht.FwE.1}Let $\lambda$ be a partition of $n$. Let $T$ be an
$n$-tableau of shape $Y\left(  \lambda\right)  $. Then,%
\[
\mathbf{F}_{T}\mathbf{E}_{T}\mathbf{F}_{T}=\left\vert \mathcal{C}\left(
T\right)  \right\vert \cdot\left\vert \mathcal{R}\left(  T\right)  \right\vert
\cdot h^{\lambda}\cdot\mathbf{F}_{T}.
\]
(See Definition \ref{def.spechtmod.flam} for the meaning of $h^{\lambda}$, and
Definition \ref{def.tableau.cell-cts} for the meanings of $\mathcal{R}\left(
T\right)  $ and $\mathcal{C}\left(  T\right)  $.)
\end{lemma}

\begin{proof}
We know that $T$ is an $n$-tableau of shape $Y\left(  \lambda\right)  $, that
is, an $n$-tableau of shape $\lambda$. Hence, Proposition
\ref{prop.specht.FT.basics} \textbf{(b)} yields%
\begin{equation}
\mathbf{F}_{T}^{2}=\underbrace{\dfrac{n!}{f^{\lambda}}}_{\substack{=h^{\lambda
}\\\text{(by the definition of }h^{\lambda}\text{)}}}\mathbf{F}_{T}%
=h^{\lambda}\cdot\mathbf{F}_{T}. \label{pf.lem.specht.FwE.1.1}%
\end{equation}

On the other hand, $\mathbf{F}_{T}=\nabla_{\operatorname*{Row}T}%
\nabla_{\operatorname*{Col}T}^{-}$ (by the definition of $\mathbf{F}_{T}$) and
$\mathbf{E}_{T}=\nabla_{\operatorname*{Col}T}^{-}\nabla_{\operatorname*{Row}%
T}$ (by the definition of $\mathbf{E}_{T}$). Hence,%
\begin{align*}
\mathbf{F}_{T}\mathbf{E}_{T}\mathbf{F}_{T}  &  =\left(  \nabla
_{\operatorname*{Row}T}\nabla_{\operatorname*{Col}T}^{-}\right)  \left(
\nabla_{\operatorname*{Col}T}^{-}\nabla_{\operatorname*{Row}T}\right)  \left(
\nabla_{\operatorname*{Row}T}\nabla_{\operatorname*{Col}T}^{-}\right) \\
&  =\nabla_{\operatorname*{Row}T}\underbrace{\left(  \nabla
_{\operatorname*{Col}T}^{-}\right)  ^{2}}_{\substack{=\left\vert
\mathcal{C}\left(  T\right)  \right\vert \cdot\nabla_{\operatorname*{Col}%
T}^{-}\\\text{(by Corollary \ref{cor.symmetrizers.square} \textbf{(b)})}%
}}\ \ \underbrace{\left(  \nabla_{\operatorname*{Row}T}\right)  ^{2}%
}_{\substack{=\left\vert \mathcal{R}\left(  T\right)  \right\vert \cdot
\nabla_{\operatorname*{Row}T}\\\text{(by Corollary
\ref{cor.symmetrizers.square} \textbf{(a)})}}}\nabla_{\operatorname*{Col}%
T}^{-}\\
&  =\nabla_{\operatorname*{Row}T}\left\vert \mathcal{C}\left(  T\right)
\right\vert \cdot\nabla_{\operatorname*{Col}T}^{-}\cdot\left\vert
\mathcal{R}\left(  T\right)  \right\vert \cdot\nabla_{\operatorname*{Row}%
T}\nabla_{\operatorname*{Col}T}^{-}\\
&  =\left\vert \mathcal{C}\left(  T\right)  \right\vert \cdot\left\vert
\mathcal{R}\left(  T\right)  \right\vert \cdot\underbrace{\nabla
_{\operatorname*{Row}T}\nabla_{\operatorname*{Col}T}^{-}}_{=\mathbf{F}_{T}%
}\cdot\underbrace{\nabla_{\operatorname*{Row}T}\nabla_{\operatorname*{Col}%
T}^{-}}_{=\mathbf{F}_{T}}\\
&  =\left\vert \mathcal{C}\left(  T\right)  \right\vert \cdot\left\vert
\mathcal{R}\left(  T\right)  \right\vert \cdot\underbrace{\mathbf{F}%
_{T}\mathbf{F}_{T}}_{\substack{=\mathbf{F}_{T}^{2}=h^{\lambda}\cdot
\mathbf{F}_{T}\\\text{(by (\ref{pf.lem.specht.FwE.1.1}))}}}=\left\vert
\mathcal{C}\left(  T\right)  \right\vert \cdot\left\vert \mathcal{R}\left(
T\right)  \right\vert \cdot h^{\lambda}\cdot\mathbf{F}_{T}.
\end{align*}
This proves Lemma \ref{lem.specht.FwE.1}.
\end{proof}

\begin{lemma}
\label{lem.specht.FwE.2}Let $\lambda$ be a partition of $n$. Let $T$ be an
$n$-tableau of shape $Y\left(  \lambda\right)  $. Assume that $n!$ is
invertible in $\mathbf{k}$. Then,%
\[
\mathbf{F}_{T}\mathcal{A}\mathbf{E}_{T}=\mathbf{k}\cdot\mathbf{F}%
_{T}\mathbf{E}_{T}.
\]
(Here, as usual, for any $\mathbf{a}\in\mathcal{A}$, we let $\mathbf{k}%
\cdot\mathbf{a}$ denote the $\mathbf{k}$-submodule $\left\{  u\mathbf{a}%
\ \mid\ u\in\mathbf{k}\right\}  $ of $\mathcal{A}$, that is, the $\mathbf{k}%
$-linear span of $\left\{  \mathbf{a}\right\}  $.)
\end{lemma}

\begin{proof}
We have $\mathbf{F}_{T}\mathbf{E}_{T}=\mathbf{F}_{T}\underbrace{1_{\mathcal{A}%
}}_{\in\mathcal{A}}\mathbf{E}_{T}\in\mathbf{F}_{T}\mathcal{A}\mathbf{E}_{T}$
and thus
\begin{equation}
\mathbf{k}\cdot\underbrace{\mathbf{F}_{T}\mathbf{E}_{T}}_{\in\mathbf{F}%
_{T}\mathcal{A}\mathbf{E}_{T}}\subseteq\mathbf{k}\cdot\mathbf{F}%
_{T}\mathcal{A}\mathbf{E}_{T}\subseteq\mathbf{F}_{T}\mathcal{A}\mathbf{E}_{T}
\label{pf.lem.specht.FwE.2.easy}%
\end{equation}
(since $\mathbf{F}_{T}\mathcal{A}\mathbf{E}_{T}$ is a $\mathbf{k}$-module).

\begin{noncompile}
On the other hand, we claim that%
\begin{equation}
\mathbf{E}_{T}\mathcal{A}\mathbf{E}_{T}\subseteq\mathbf{kE}_{T}.
\label{pf.lem.specht.FwE.2.EAE}%
\end{equation}

\begin{proof}
[Proof of (\ref{pf.lem.specht.FwE.2.EAE}).]Let $\mathbf{a}\in\mathcal{A}$.
Then, $\mathbf{a}\in\mathcal{A}=\mathbf{k}\left[  S_{n}\right]  $. Hence,
Proposition \ref{prop.specht.ET.ETaET} shows that $\mathbf{E}_{T}%
\mathbf{aE}_{T}=\kappa\mathbf{E}_{T}$ for some scalar $\kappa\in\mathbf{k}$.
Thus, for this $\kappa$, we have $\mathbf{E}_{T}\mathbf{aE}_{T}=\kappa
\mathbf{E}_{T}\in\mathbf{kE}_{T}$.

Forget that we fixed $\mathbf{a}$. We thus have shown that $\mathbf{E}%
_{T}\mathbf{aE}_{T}\in\mathbf{kE}_{T}$ for each $\mathbf{a}\in\mathcal{A}$. In
other words,%
\[
\left\{  \mathbf{E}_{T}\mathbf{aE}_{T}\ \mid\ \mathbf{a}\in\mathcal{A}%
\right\}  \subseteq\mathbf{kE}_{T}.
\]
In other words, $\mathbf{E}_{T}\mathcal{A}\mathbf{E}_{T}\subseteq
\mathbf{kE}_{T}$ (since $\mathbf{E}_{T}\mathcal{A}\mathbf{E}_{T}=\left\{
\mathbf{E}_{T}\mathbf{aE}_{T}\ \mid\ \mathbf{a}\in\mathcal{A}\right\}  $).
This proves (\ref{pf.lem.specht.FwE.2.EAE}).
\end{proof}
\end{noncompile}

On the other hand, Lemma \ref{lem.specht.FwE.1} yields
\begin{equation}
\mathbf{F}_{T}\mathbf{E}_{T}\mathbf{F}_{T}=\left\vert \mathcal{C}\left(
T\right)  \right\vert \cdot\left\vert \mathcal{R}\left(  T\right)  \right\vert
\cdot h^{\lambda}\cdot\mathbf{F}_{T}, \label{pf.lem.specht.FwE.2.1}%
\end{equation}
where $h^{\lambda}$ and $\mathcal{R}\left(  T\right)  $ and $\mathcal{C}%
\left(  T\right)  $ are defined as in Lemma \ref{lem.specht.FwE.1}.

But $h^{\lambda}=\dfrac{n!}{f^{\lambda}}$ (by the definition of $h^{\lambda}%
$). Thus, $h^{\lambda}$ is a divisor of $n!$. Hence, $h^{\lambda}$ is
invertible in $\mathbf{k}$ (since $n!$ is invertible in $\mathbf{k}$).
Furthermore, $\mathcal{R}\left(  T\right)  $ is a subgroup of $S_{n}$ (by its
definition). Hence, Lagrange's theorem shows that $\left\vert \mathcal{R}%
\left(  T\right)  \right\vert $ is a divisor of $\left\vert S_{n}\right\vert
=n!$. Thus, $\left\vert \mathcal{R}\left(  T\right)  \right\vert $ is
invertible in $\mathbf{k}$ (since $n!$ is invertible in $\mathbf{k}$).
Similarly, we can see that $\left\vert \mathcal{C}\left(  T\right)
\right\vert $ is invertible in $\mathbf{k}$. Now, we can divide the equality
(\ref{pf.lem.specht.FwE.2.1}) by the three scalars $\left\vert \mathcal{C}%
\left(  T\right)  \right\vert $, $\left\vert \mathcal{R}\left(  T\right)
\right\vert $ and $h^{\lambda}$ (since we have just shown that these scalars
$\left\vert \mathcal{C}\left(  T\right)  \right\vert $, $\left\vert
\mathcal{R}\left(  T\right)  \right\vert $ and $h^{\lambda}$ are invertible in
$\mathbf{k}$). As a result, we obtain
\begin{equation}
\dfrac{1}{\left\vert \mathcal{C}\left(  T\right)  \right\vert \cdot\left\vert
\mathcal{R}\left(  T\right)  \right\vert \cdot h^{\lambda}}\cdot\mathbf{F}%
_{T}\mathbf{E}_{T}\mathbf{F}_{T}=\mathbf{F}_{T}. \label{pf.lem.specht.FwE.2.2}%
\end{equation}

Now, we shall show that $\mathbf{F}_{T}\mathcal{A}\mathbf{E}_{T}%
\subseteq\mathbf{k}\cdot\mathbf{F}_{T}\mathbf{E}_{T}$. Indeed, fix
$\mathbf{b}\in\mathcal{A}$. Then, $\mathbf{F}_{T}\mathbf{b}\in\mathcal{A}%
=\mathbf{k}\left[  S_{n}\right]  $. Hence, Proposition
\ref{prop.specht.ET.ETaET} (applied to $\mathbf{a}=\mathbf{F}_{T}\mathbf{b}$)
shows that $\mathbf{E}_{T}\mathbf{F}_{T}\mathbf{bE}_{T}=\kappa\mathbf{E}_{T}$
for some scalar $\kappa\in\mathbf{k}$. Consider this $\kappa$. We have%
\begin{align*}
\underbrace{\mathbf{F}_{T}}_{\substack{=\dfrac{1}{\left\vert \mathcal{C}%
\left(  T\right)  \right\vert \cdot\left\vert \mathcal{R}\left(  T\right)
\right\vert \cdot h^{\lambda}}\cdot\mathbf{F}_{T}\mathbf{E}_{T}\mathbf{F}%
_{T}\\\text{(by (\ref{pf.lem.specht.FwE.2.2}))}}}\mathbf{bE}_{T}  &
=\dfrac{1}{\left\vert \mathcal{C}\left(  T\right)  \right\vert \cdot\left\vert
\mathcal{R}\left(  T\right)  \right\vert \cdot h^{\lambda}}\cdot\mathbf{F}%
_{T}\underbrace{\mathbf{E}_{T}\mathbf{F}_{T}\mathbf{bE}_{T}}_{=\kappa
\mathbf{E}_{T}}\\
&  =\dfrac{1}{\left\vert \mathcal{C}\left(  T\right)  \right\vert
\cdot\left\vert \mathcal{R}\left(  T\right)  \right\vert \cdot h^{\lambda}%
}\cdot\mathbf{F}_{T}\cdot\kappa\mathbf{E}_{T}\\
&  =\underbrace{\dfrac{\kappa}{\left\vert \mathcal{C}\left(  T\right)
\right\vert \cdot\left\vert \mathcal{R}\left(  T\right)  \right\vert \cdot
h^{\lambda}}}_{\in\mathbf{k}}\cdot\,\mathbf{F}_{T}\mathbf{E}_{T}\\
&  \in\mathbf{k}\cdot\mathbf{F}_{T}\mathbf{E}_{T}.
\end{align*}

Forget that we fixed $\mathbf{b}$. We thus have shown that $\mathbf{F}%
_{T}\mathbf{bE}_{T}\in\mathbf{k}\cdot\mathbf{F}_{T}\mathbf{E}_{T}$ for each
$\mathbf{b}\in\mathcal{A}$. In other words,%
\[
\left\{  \mathbf{F}_{T}\mathbf{bE}_{T}\ \mid\ \mathbf{b}\in\mathcal{A}%
\right\}  \subseteq\mathbf{k}\cdot\mathbf{F}_{T}\mathbf{E}_{T}.
\]
In other words, $\mathbf{F}_{T}\mathcal{A}\mathbf{E}_{T}\subseteq
\mathbf{k}\cdot\mathbf{F}_{T}\mathbf{E}_{T}$ (since $\mathbf{F}_{T}%
\mathcal{A}\mathbf{E}_{T}=\left\{  \mathbf{F}_{T}\mathbf{bE}_{T}%
\ \mid\ \mathbf{b}\in\mathcal{A}\right\}  $). Combining this with
(\ref{pf.lem.specht.FwE.2.easy}), we obtain $\mathbf{F}_{T}\mathcal{A}%
\mathbf{E}_{T}=\mathbf{k}\cdot\mathbf{F}_{T}\mathbf{E}_{T}$. Thus, Lemma
\ref{lem.specht.FwE.2} is proved.
\end{proof}

\begin{lemma}
\label{lem.specht.FwE.3}Let $\lambda$ be a partition of $n$. Let $U$ and $V$
be two $n$-tableaux of shape $Y\left(  \lambda\right)  $. Assume that $n!$ is
invertible in $\mathbf{k}$. Then,%
\[
\mathbf{F}_{U}\mathcal{A}\mathbf{E}_{V}=\mathbf{k}\cdot\mathbf{F}_{U}%
w_{U,V}\mathbf{E}_{V}.
\]

\end{lemma}

\begin{proof}
Proposition \ref{prop.specht.EPQ.EPEQ} (applied to $D=Y\left(  \lambda\right)
$ and $P=U$ and $Q=V$) yields%
\[
\mathbf{E}_{U,V}=w_{U,V}\mathbf{E}_{V}=\mathbf{E}_{U}w_{U,V}.
\]
Hence, in particular, $w_{U,V}\mathbf{E}_{V}=\mathbf{E}_{U}w_{U,V}$.

Moreover, $w_{U,V}$ is a permutation in $S_{n}$, thus an invertible element of
$\mathcal{A}$ (since any element of a group $G$ is invertible in the group
algebra $\mathbf{k}\left[  G\right]  $). Thus, $\mathcal{A}w_{U,V}%
=\mathcal{A}$ (since any invertible element $\mathbf{p}$ of $\mathcal{A}$
satisfies $\mathcal{A}\mathbf{p}=\mathcal{A}$). Hence, $\mathcal{A}%
=\mathcal{A}w_{U,V}$, so that%
\begin{align*}
\mathbf{F}_{U}\underbrace{\mathcal{A}}_{=\mathcal{A}w_{U,V}}\mathbf{E}_{V}  &
=\mathbf{F}_{U}\mathcal{A}\underbrace{w_{U,V}\mathbf{E}_{V}}_{=\mathbf{E}%
_{U}w_{U,V}}=\underbrace{\mathbf{F}_{U}\mathcal{A}\mathbf{E}_{U}%
}_{\substack{=\mathbf{k}\cdot\mathbf{F}_{U}\mathbf{E}_{U}\\\text{(by Lemma
\ref{lem.specht.FwE.2},}\\\text{applied to }T=U\text{)}}}w_{U,V}%
=\mathbf{k}\cdot\mathbf{F}_{U}\underbrace{\mathbf{E}_{U}w_{U,V}}%
_{=w_{U,V}\mathbf{E}_{V}}\\
&  =\mathbf{k}\cdot\mathbf{F}_{U}w_{U,V}\mathbf{E}_{V}.
\end{align*}
This proves Lemma \ref{lem.specht.FwE.3}.
\end{proof}

Our next two lemmas are general linear-algebraic facts about free modules of
finite rank. When $\mathbf{k}$ is a field, they are both well-known from any
good course on linear algebra, but they are in fact true for any commutative
ring $\mathbf{k}$.

\begin{lemma}
\label{lem.free-mod.sur=bij}Let $r\in\mathbb{N}$. Let $V$ and $W$ be two free
$\mathbf{k}$-modules of rank $r$. Let $f:V\rightarrow W$ be a surjective
$\mathbf{k}$-linear map. Then, $f$ is a $\mathbf{k}$-module isomorphism.
\end{lemma}

\begin{proof}
Both $\mathbf{k}$-modules $V$ and $W$ are free of rank $r$. Hence, they have
bases $\left(  v_{1},v_{2},\ldots,v_{r}\right)  $ and $\left(  w_{1}%
,w_{2},\ldots,w_{r}\right)  $ of size $r$. Consider these bases.

The map $f$ is surjective. In other words, $W=f\left(  V\right)  $. Hence, for
each $i\in\left\{  1,2,\ldots,r\right\}  $, there exists some $x_{i}\in V$
such that $w_{i}=f\left(  x_{i}\right)  $ (since $w_{i}\in W=f\left(
V\right)  $). Consider this $x_{i}$. Thus, we have constructed $r$ vectors
$x_{1},x_{2},\ldots,x_{r}$ in $V$.

Recall that $\left(  w_{1},w_{2},\ldots,w_{r}\right)  $ is a basis of the
$\mathbf{k}$-module $W$. Hence, we can define a $\mathbf{k}$-linear map
$g:W\rightarrow V$ by setting%
\begin{equation}
\left(  g\left(  w_{i}\right)  =x_{i}\ \ \ \ \ \ \ \ \ \ \text{for each }%
i\in\left\{  1,2,\ldots,r\right\}  \right)  \label{pf.lem.free-mod.sur=bij.1}%
\end{equation}
(since a $\mathbf{k}$-linear map from a free $\mathbf{k}$-module can be
defined by choosing its values on the elements of a given basis). Consider
this map $g$.

Now, for each $i\in\left\{  1,2,\ldots,r\right\}  $, we have%
\begin{align*}
\left(  f\circ g\right)  \left(  w_{i}\right)   &  =f\left(  g\left(
w_{i}\right)  \right) \\
&  =f\left(  x_{i}\right)  \ \ \ \ \ \ \ \ \ \ \left(  \text{since
(\ref{pf.lem.free-mod.sur=bij.1}) yields }g\left(  w_{i}\right)  =x_{i}\right)
\\
&  =w_{i}\ \ \ \ \ \ \ \ \ \ \left(  \text{since }x_{i}\text{ was defined to
satisfy }w_{i}=f\left(  x_{i}\right)  \right) \\
&  =\operatorname*{id}\nolimits_{W}\left(  w_{i}\right)  .
\end{align*}
Thus, the two $\mathbf{k}$-linear maps $f\circ g:W\rightarrow W$ and
$\operatorname*{id}\nolimits_{W}:W\rightarrow W$ agree on each of the $r$
vectors $w_{1},w_{2},\ldots,w_{r}$. In other words, these two $\mathbf{k}%
$-linear maps agree on the basis $\left(  w_{1},w_{2},\ldots,w_{r}\right)  $
of $W$. Hence, these two maps must be identical (because if two $\mathbf{k}%
$-linear maps agree on a basis of their domain, then they must be identical).
In other words, $f\circ g=\operatorname*{id}\nolimits_{W}$. Hence, Lemma
\ref{lem.linalg.AB=1-BA=1} yields that $g\circ f=\operatorname*{id}%
\nolimits_{V}$. Combining this with $f\circ g=\operatorname*{id}\nolimits_{W}%
$, we conclude that the maps $f$ and $g$ are mutually inverse. Hence, the map
$f$ is invertible, and hence a $\mathbf{k}$-module isomorphism (since it is
$\mathbf{k}$-linear). This proves Lemma \ref{lem.free-mod.sur=bij}.
\end{proof}

\begin{lemma}
\label{lem.free-mod.span=basis}Let $r\in\mathbb{N}$. Let $M$ be a free
$\mathbf{k}$-module of rank $r$. Let $\left(  m_{i}\right)  _{i\in I}$ be a
family of vectors in $M$ with $\left\vert I\right\vert =r$. Assume that this
family $\left(  m_{i}\right)  _{i\in I}$ spans $M$. Then, $\left(
m_{i}\right)  _{i\in I}$ is a basis of $M$.
\end{lemma}

\begin{proof}
The set $I$ has $r$ elements (since $\left\vert I\right\vert =r$), and is only
used for indexing our family $\left(  m_{i}\right)  _{i\in I}$. Thus, we can
WLOG assume that $I=\left\{  1,2,\ldots,r\right\}  $ (since we can otherwise
rename the $r$ elements of $I$ as $1,2,\ldots,r$ and accordingly rename the
entries of our family $\left(  m_{i}\right)  _{i\in I}$, without changing the
essence of anything). Assume this. Then, the family $\left(  m_{i}\right)
_{i\in I}=\left(  m_{i}\right)  _{i\in\left\{  1,2,\ldots,r\right\}  }$ is
just an $r$-tuple $\left(  m_{1},m_{2},\ldots,m_{r}\right)  $. We have assumed
that this family $\left(  m_{i}\right)  _{i\in I}$ spans $M$. In other words,
the $r$-tuple $\left(  m_{1},m_{2},\ldots,m_{r}\right)  $ spans $M$.

But $M$ is a free $\mathbf{k}$-module of rank $r$. Hence, $M$ has a basis
$\left(  b_{1},b_{2},\ldots,b_{r}\right)  $ consisting of $r$ vectors.
Consider this basis.

Thus, $\left(  b_{1},b_{2},\ldots,b_{r}\right)  $ is a basis of $M$. Hence, we
can define a $\mathbf{k}$-linear map $f:M\rightarrow M$ by setting%
\[
\left(  f\left(  b_{i}\right)  =m_{i}\ \ \ \ \ \ \ \ \ \ \text{for each }%
i\in\left\{  1,2,\ldots,r\right\}  \right)
\]
(since a $\mathbf{k}$-linear map from a free $\mathbf{k}$-module can be
defined by choosing its values on the elements of a given basis). Consider
this map $f$.

The definition of $f$ shows that $f\left(  b_{i}\right)  =m_{i}$ for each
$i\in\left\{  1,2,\ldots,r\right\}  $. Hence,%
\begin{equation}
\left(  f\left(  b_{1}\right)  ,f\left(  b_{2}\right)  ,\ldots,f\left(
b_{r}\right)  \right)  =\left(  m_{1},m_{2},\ldots,m_{r}\right)  .
\label{pf.lem.free-mod.span=basis.3}%
\end{equation}

Since $\left(  b_{1},b_{2},\ldots,b_{r}\right)  $ is a basis of $M$, we have
$M=\operatorname*{span}\nolimits_{\mathbf{k}}\left\{  b_{1},b_{2},\ldots
,b_{r}\right\}  $. Applying the map $f$ to this equality, we obtain%
\begin{align*}
f\left(  M\right)   &  =f\left(  \operatorname*{span}\nolimits_{\mathbf{k}%
}\left\{  b_{1},b_{2},\ldots,b_{r}\right\}  \right) \\
&  =\operatorname*{span}\nolimits_{\mathbf{k}}\underbrace{\left\{  f\left(
b_{1}\right)  ,f\left(  b_{2}\right)  ,\ldots,f\left(  b_{r}\right)  \right\}
}_{\substack{=\left\{  m_{1},m_{2},\ldots,m_{r}\right\}  \\\text{(by
(\ref{pf.lem.free-mod.span=basis.3}))}}}\\
&  \ \ \ \ \ \ \ \ \ \ \ \ \ \ \ \ \ \ \ \ \left(
\begin{array}
[c]{c}%
\text{since the map }f\text{ is }\mathbf{k}\text{-linear, and thus}\\
\text{respects }\mathbf{k}\text{-linear combinations}%
\end{array}
\right) \\
&  =\operatorname*{span}\nolimits_{\mathbf{k}}\left\{  m_{1},m_{2}%
,\ldots,m_{r}\right\}  =M
\end{align*}
(since the $r$-tuple $\left(  m_{1},m_{2},\ldots,m_{r}\right)  $ spans $M$).
In other words, the map $f$ is surjective. Hence, Lemma
\ref{lem.free-mod.sur=bij} (applied to $V=M$ and $W=M$) shows that this map
$f$ is a $\mathbf{k}$-module isomorphism. Therefore, $f$ must send any basis
of $M$ to a basis of $M$ (since a $\mathbf{k}$-module isomorphism always sends
a basis of its domain to a basis of its target). Since $\left(  b_{1}%
,b_{2},\ldots,b_{r}\right)  $ is a basis of $M$, we thus conclude that
$\left(  f\left(  b_{1}\right)  ,f\left(  b_{2}\right)  ,\ldots,f\left(
b_{r}\right)  \right)  $ is a basis of $M$ as well (since $f$ sends the basis
$\left(  b_{1},b_{2},\ldots,b_{r}\right)  $ to $\left(  f\left(  b_{1}\right)
,f\left(  b_{2}\right)  ,\ldots,f\left(  b_{r}\right)  \right)  $). In other
words, $\left(  m_{1},m_{2},\ldots,m_{r}\right)  $ is a basis of $M$ (since
(\ref{pf.lem.free-mod.span=basis.3}) shows that $\left(  m_{1},m_{2}%
,\ldots,m_{r}\right)  =\left(  f\left(  b_{1}\right)  ,f\left(  b_{2}\right)
,\ldots,f\left(  b_{r}\right)  \right)  $). In other words, $\left(
m_{i}\right)  _{i\in I}$ is a basis of $M$ (since the family $\left(
m_{i}\right)  _{i\in I}$ is just the $r$-tuple $\left(  m_{1},m_{2}%
,\ldots,m_{r}\right)  $). This proves Lemma \ref{lem.free-mod.span=basis}.
\end{proof}

Next, we state a version of Theorem \ref{thm.specht.FwE.basis} for a specific partition:

\begin{proposition}
\label{prop.specht.FwE.basis-lam}Fix a partition $\lambda$ of $n$. Let
$\operatorname*{SYT}\left(  \lambda\right)  $ be the set of all standard
$n$-tableaux of shape $Y\left(  \lambda\right)  $. Assume that $n!$ is
invertible in $\mathbf{k}$. Then, the family $\left(  \mathbf{F}_{U}%
w_{U,V}\mathbf{E}_{V}\right)  _{U,V\in\operatorname*{SYT}\left(
\lambda\right)  }$ is a basis of the $\mathbf{k}$-module $\mathcal{A}%
\mathbf{E}_{\lambda}$.
\end{proposition}

\begin{proof}
Recall Definition \ref{def.spechtmod.flam}. In particular, $h^{\lambda}%
=\dfrac{n!}{f^{\lambda}}$ (by the definition of $h^{\lambda}$). Thus,
$h^{\lambda}$ is a divisor of $n!$. Hence, $h^{\lambda}$ is invertible in
$\mathbf{k}$ (since $n!$ is invertible in $\mathbf{k}$). Furthermore, note
that $\left\vert \operatorname*{SYT}\left(  \lambda\right)  \right\vert
=f^{\lambda}$\ \ \ \ \footnote{\textit{Proof.} The set of all standard
$n$-tableaux of shape $Y\left(  \lambda\right)  $ was denoted by
$\operatorname*{SYT}\left(  Y\left(  \lambda\right)  \right)  $ in the proof
of Corollary \ref{cor.tableaux.flamt}, but we are now denoting it by
$\operatorname*{SYT}\left(  \lambda\right)  $. Hence, $\operatorname*{SYT}%
\left(  Y\left(  \lambda\right)  \right)  =\operatorname*{SYT}\left(
\lambda\right)  $.
\par
However, in our proof of Corollary \ref{cor.tableaux.flamt}, we saw that
$f^{\lambda}=\left\vert \operatorname*{SYT}\left(  Y\left(  \lambda\right)
\right)  \right\vert $. In view of $\operatorname*{SYT}\left(  Y\left(
\lambda\right)  \right)  =\operatorname*{SYT}\left(  \lambda\right)  $, we can
rewrite this as $f^{\lambda}=\left\vert \operatorname*{SYT}\left(
\lambda\right)  \right\vert $. In other words, $\left\vert \operatorname*{SYT}%
\left(  \lambda\right)  \right\vert =f^{\lambda}$.}.

Let $D=Y\left(  \lambda\right)  $. Thus, the standard $n$-tableaux of shape
$D$ are just the standard $n$-tableaux of shape $Y\left(  \lambda\right)  $.
Hence, the set $\operatorname*{SYT}\left(  D\right)  $ in Theorem
\ref{thm.specht.AElam.decomp} (defined as the set of all standard $n$-tableaux
of shape $D$) is precisely our set $\operatorname*{SYT}\left(  \lambda\right)
$ (defined as the set of all standard $n$-tableaux of shape $Y\left(
\lambda\right)  $). Hence, Theorem \ref{thm.specht.AElam.decomp} (with
$\operatorname*{SYT}\left(  D\right)  $ renamed as $\operatorname*{SYT}\left(
\lambda\right)  $) shows that%
\begin{align}
\mathcal{A}\mathbf{E}_{\lambda}  &  =\bigoplus_{T\in\operatorname*{SYT}\left(
\lambda\right)  }\mathcal{A}\mathbf{E}_{T}\ \ \ \ \ \ \ \ \ \ \left(  \text{an
internal direct sum}\right) \nonumber\\
&  =\sum_{T\in\operatorname*{SYT}\left(  \lambda\right)  }\mathcal{A}%
\mathbf{E}_{T} \label{pf.prop.specht.FwE.basis-lam.2}%
\end{align}
(since direct sums are sums).

On the other hand, Proposition \ref{prop.spechtmod.Elam.incent} yields
$\mathbf{E}_{\lambda}\in Z\left(  \mathbf{k}\left[  S_{n}\right]  \right)  $.
Hence, Proposition \ref{prop.ZkSn.S=id} (applied to $\mathbf{a}=\mathbf{E}%
_{\lambda}$) yields $S\left(  \mathbf{E}_{\lambda}\right)  =\mathbf{E}%
_{\lambda}$. Furthermore, $S:\mathcal{A}\rightarrow\mathcal{A}$ is an
involution (by Theorem \ref{thm.S.auto} \textbf{(b)}), and thus a bijection.
Hence, $S\left(  \mathcal{A}\right)  =\mathcal{A}$. Moreover, $S$ is an
algebra anti-morphism (by Theorem \ref{thm.S.auto} \textbf{(a)}). Therefore,%
\[
S\left(  \mathcal{A}\mathbf{E}_{\lambda}\right)  =\underbrace{S\left(
\mathbf{E}_{\lambda}\right)  }_{=\mathbf{E}_{\lambda}}\underbrace{S\left(
\mathcal{A}\right)  }_{=\mathcal{A}}=\mathbf{E}_{\lambda}\mathcal{A}.
\]
Hence,%
\begin{align*}
\mathbf{E}_{\lambda}\mathcal{A}  &  =S\left(  \mathcal{A}\mathbf{E}_{\lambda
}\right)  =S\left(  \sum_{T\in\operatorname*{SYT}\left(  \lambda\right)
}\mathcal{A}\mathbf{E}_{T}\right)  \ \ \ \ \ \ \ \ \ \ \left(  \text{by
(\ref{pf.prop.specht.FwE.basis-lam.2})}\right) \\
&  =\sum_{T\in\operatorname*{SYT}\left(  \lambda\right)  }\underbrace{S\left(
\mathbf{E}_{T}\right)  }_{\substack{=\mathbf{F}_{T}\\\text{(by Proposition
\ref{prop.specht.FT.basics} \textbf{(a)})}}}\underbrace{S\left(
\mathcal{A}\right)  }_{=\mathcal{A}}\ \ \ \ \ \ \ \ \ \ \left(
\begin{array}
[c]{c}%
\text{since }S\text{ is an algebra}\\
\text{anti-morphism}%
\end{array}
\right) \\
&  =\sum_{T\in\operatorname*{SYT}\left(  \lambda\right)  }\mathbf{F}%
_{T}\mathcal{A}.
\end{align*}
Multiplying this equality with (\ref{pf.prop.specht.FwE.basis-lam.2}), we
obtain%
\begin{align*}
\left(  \mathbf{E}_{\lambda}\mathcal{A}\right)  \left(  \mathcal{A}%
\mathbf{E}_{\lambda}\right)   &  =\left(  \sum_{T\in\operatorname*{SYT}\left(
\lambda\right)  }\mathbf{F}_{T}\mathcal{A}\right)  \left(  \sum_{T\in
\operatorname*{SYT}\left(  \lambda\right)  }\mathcal{A}\mathbf{E}_{T}\right)
\\
&  =\left(  \sum_{U\in\operatorname*{SYT}\left(  \lambda\right)  }%
\mathbf{F}_{U}\mathcal{A}\right)  \left(  \sum_{V\in\operatorname*{SYT}\left(
\lambda\right)  }\mathcal{A}\mathbf{E}_{V}\right) \\
&  \ \ \ \ \ \ \ \ \ \ \ \ \ \ \ \ \ \ \ \ \left(
\begin{array}
[c]{c}%
\text{here, we have renamed the summation}\\
\text{indices }T\text{ and }T\text{ as }U\text{ and }V
\end{array}
\right) \\
&  =\underbrace{\sum_{U\in\operatorname*{SYT}\left(  \lambda\right)  }%
\ \ \sum_{V\in\operatorname*{SYT}\left(  \lambda\right)  }}_{=\sum_{\left(
U,V\right)  \in\operatorname*{SYT}\left(  \lambda\right)  \times
\operatorname*{SYT}\left(  \lambda\right)  }}\underbrace{\left(
\mathbf{F}_{U}\mathcal{A}\right)  \left(  \mathcal{A}\mathbf{E}_{V}\right)
}_{=\mathbf{F}_{U}\mathcal{AA}\mathbf{E}_{V}}\\
&  =\sum_{\left(  U,V\right)  \in\operatorname*{SYT}\left(  \lambda\right)
\times\operatorname*{SYT}\left(  \lambda\right)  }\mathbf{F}_{U}%
\underbrace{\mathcal{AA}}_{=\mathcal{A}}\mathbf{E}_{V}\\
&  =\sum_{\left(  U,V\right)  \in\operatorname*{SYT}\left(  \lambda\right)
\times\operatorname*{SYT}\left(  \lambda\right)  }\underbrace{\mathbf{F}%
_{U}\mathcal{A}\mathbf{E}_{V}}_{\substack{=\mathbf{k}\cdot\mathbf{F}%
_{U}w_{U,V}\mathbf{E}_{V}\\\text{(by Lemma \ref{lem.specht.FwE.3})}}}\\
&  =\sum_{\left(  U,V\right)  \in\operatorname*{SYT}\left(  \lambda\right)
\times\operatorname*{SYT}\left(  \lambda\right)  }\underbrace{\mathbf{k}%
\cdot\mathbf{F}_{U}w_{U,V}\mathbf{E}_{V}}_{\substack{=\operatorname*{span}%
\nolimits_{\mathbf{k}}\left\{  \mathbf{F}_{U}w_{U,V}\mathbf{E}_{V}\right\}
\\\text{(since }\mathbf{k}\cdot v=\operatorname*{span}\nolimits_{\mathbf{k}%
}\left\{  v\right\}  \text{ for any vector }v\text{)}}}\\
&  =\sum_{\left(  U,V\right)  \in\operatorname*{SYT}\left(  \lambda\right)
\times\operatorname*{SYT}\left(  \lambda\right)  }\operatorname*{span}%
\nolimits_{\mathbf{k}}\left\{  \mathbf{F}_{U}w_{U,V}\mathbf{E}_{V}\right\} \\
&  =\operatorname*{span}\nolimits_{\mathbf{k}}\left\{  \mathbf{F}_{U}%
w_{U,V}\mathbf{E}_{V}\ \mid\ \left(  U,V\right)  \in\operatorname*{SYT}\left(
\lambda\right)  \times\operatorname*{SYT}\left(  \lambda\right)  \right\}  .
\end{align*}
In view of%
\begin{align*}
\left(  \mathbf{E}_{\lambda}\mathcal{A}\right)  \left(  \mathcal{A}%
\mathbf{E}_{\lambda}\right)   &  =\mathbf{E}_{\lambda}\underbrace{\mathcal{AA}%
}_{=\mathcal{A}}\mathbf{E}_{\lambda}=\underbrace{\mathbf{E}_{\lambda
}\mathcal{A}}_{\substack{=\mathcal{A}\mathbf{E}_{\lambda}\\\text{(since
}\mathbf{E}_{\lambda}\in Z\left(  \mathbf{k}\left[  S_{n}\right]  \right)
\text{)}}}\mathbf{E}_{\lambda}\\
&  =\mathcal{A}\underbrace{\mathbf{E}_{\lambda}\mathbf{E}_{\lambda}%
}_{\substack{=\mathbf{E}_{\lambda}^{2}=\left(  h^{\lambda}\right)
^{2}\mathbf{E}_{\lambda}\\\text{(by Theorem
\ref{thm.spechtmod.Elam.idp-strong})}}}=\underbrace{\mathcal{A}\left(
h^{\lambda}\right)  ^{2}}_{\substack{=\mathcal{A}\\\text{(since }\left(
h^{\lambda}\right)  ^{2}\text{ is invertible in }\mathbf{k}\\\text{(because
}h^{\lambda}\text{ is invertible in }\mathbf{k}\text{))}}}\mathbf{E}_{\lambda
}=\mathcal{A}\mathbf{E}_{\lambda},
\end{align*}
we can rewrite this as%
\[
\mathcal{A}\mathbf{E}_{\lambda}=\operatorname*{span}\nolimits_{\mathbf{k}%
}\left\{  \mathbf{F}_{U}w_{U,V}\mathbf{E}_{V}\ \mid\ \left(  U,V\right)
\in\operatorname*{SYT}\left(  \lambda\right)  \times\operatorname*{SYT}\left(
\lambda\right)  \right\}  .
\]
In other words, the family $\left(  \mathbf{F}_{U}w_{U,V}\mathbf{E}%
_{V}\right)  _{\left(  U,V\right)  \in\operatorname*{SYT}\left(
\lambda\right)  \times\operatorname*{SYT}\left(  \lambda\right)  }$ spans the
$\mathbf{k}$-module $\mathcal{A}\mathbf{E}_{\lambda}$. This family has size%
\[
\left\vert \operatorname*{SYT}\left(  \lambda\right)  \times
\operatorname*{SYT}\left(  \lambda\right)  \right\vert =\underbrace{\left\vert
\operatorname*{SYT}\left(  \lambda\right)  \right\vert }_{=f^{\lambda}}%
\cdot\underbrace{\left\vert \operatorname*{SYT}\left(  \lambda\right)
\right\vert }_{=f^{\lambda}}=f^{\lambda}\cdot f^{\lambda}=\left(  f^{\lambda
}\right)  ^{2}.
\]

But Lemma \ref{lem.specht.dim-AElam} shows that $\mathcal{A}\mathbf{E}%
_{\lambda}$ is a free $\mathbf{k}$-module of rank $\left(  f^{\lambda}\right)
^{2}$.

Hence, Lemma \ref{lem.free-mod.span=basis} (applied to $r=\left(  f^{\lambda
}\right)  ^{2}$ and $M=\mathcal{A}\mathbf{E}_{\lambda}$ and
$I=\operatorname*{SYT}\left(  \lambda\right)  \times\operatorname*{SYT}\left(
\lambda\right)  $ and $\left(  m_{i}\right)  _{i\in I}=\left(  \mathbf{F}%
_{U}w_{U,V}\mathbf{E}_{V}\right)  _{\left(  U,V\right)  \in\operatorname*{SYT}%
\left(  \lambda\right)  \times\operatorname*{SYT}\left(  \lambda\right)  }$)
shows that the family $\left(  \mathbf{F}_{U}w_{U,V}\mathbf{E}_{V}\right)
_{\left(  U,V\right)  \in\operatorname*{SYT}\left(  \lambda\right)
\times\operatorname*{SYT}\left(  \lambda\right)  }$ is a basis of the
$\mathbf{k}$-module $\mathcal{A}\mathbf{E}_{\lambda}$ (since we have shown
that this family spans $\mathcal{A}\mathbf{E}_{\lambda}$, and since its size
is $\left\vert \operatorname*{SYT}\left(  \lambda\right)  \times
\operatorname*{SYT}\left(  \lambda\right)  \right\vert =\left(  f^{\lambda
}\right)  ^{2}$). In other words, the family $\left(  \mathbf{F}_{U}%
w_{U,V}\mathbf{E}_{V}\right)  _{U,V\in\operatorname*{SYT}\left(
\lambda\right)  }$ is a basis of the $\mathbf{k}$-module $\mathcal{A}%
\mathbf{E}_{\lambda}$ (since \textquotedblleft$U,V\in\operatorname*{SYT}%
\left(  \lambda\right)  $\textquotedblright\ is just shorthand for
\textquotedblleft$\left(  U,V\right)  \in\operatorname*{SYT}\left(
\lambda\right)  \times\operatorname*{SYT}\left(  \lambda\right)
$\textquotedblright). This proves Proposition \ref{prop.specht.FwE.basis-lam}.
\end{proof}

We are now ready for the proof of Theorem \ref{thm.specht.FwE.basis}:

\begin{proof}
[Proof of Theorem \ref{thm.specht.FwE.basis}.]This is very similar to the
proof of Corollary \ref{cor.specht.A.nat-basis}.

Corollary \ref{cor.spechtmod.Elam.A=sum1} shows that
\[
\mathcal{A}=\bigoplus\limits_{\lambda\text{ is a partition of }n}%
\mathcal{A}\mathbf{E}_{\lambda}\ \ \ \ \ \ \ \ \ \ \left(  \text{an internal
direct sum}\right)  .
\]
Thus, the sum $\sum_{\lambda\text{ is a partition of }n}\mathcal{A}%
\mathbf{E}_{\lambda}$ is direct.

From Proposition \ref{prop.specht.FwE.basis-lam}, we know that the family
$\left(  \mathbf{F}_{U}w_{U,V}\mathbf{E}_{V}\right)  _{U,V\in
\operatorname*{SYT}\left(  \lambda\right)  }$ is a basis of the $\mathbf{k}%
$-module $\mathcal{A}\mathbf{E}_{\lambda}$ whenever $\lambda$ is a partition
of $n$. Thus, Observation 1 from our above proof of Corollary
\ref{cor.specht.A.nat-basis} can be applied to $M=\mathcal{A}$ and $I=\left\{
\text{partitions of }n\right\}  $ and $N_{\lambda}=\mathcal{A}\mathbf{E}%
_{\lambda}$ and $\left(  v_{\lambda,j}\right)  _{j\in J_{\lambda}}=\left(
\mathbf{F}_{U}w_{U,V}\mathbf{E}_{V}\right)  _{U,V\in\operatorname*{SYT}\left(
\lambda\right)  }$ (since the sum $\sum_{\lambda\text{ is a partition of }%
n}\mathcal{A}\mathbf{E}_{\lambda}$ is direct, and since each $\mathbf{k}%
$-module $\mathcal{A}\mathbf{E}_{\lambda}$ has a basis $\left(  \mathbf{F}%
_{U}w_{U,V}\mathbf{E}_{V}\right)  _{U,V\in\operatorname*{SYT}\left(
\lambda\right)  }$). Thus, we conclude that the internal direct sum
$\bigoplus\limits_{\lambda\text{ is a partition of }n}\mathcal{A}%
\mathbf{E}_{\lambda}$ has a basis $\left(  \mathbf{F}_{U}w_{U,V}\mathbf{E}%
_{V}\right)  _{\lambda\text{ is a partition of }n\text{, and }U,V\in
\operatorname*{SYT}\left(  \lambda\right)  }$. In other words, the
$\mathbf{k}$-module $\mathcal{A}$ has a basis $\left(  \mathbf{F}_{U}%
w_{U,V}\mathbf{E}_{V}\right)  _{\lambda\text{ is a partition of }n\text{, and
}U,V\in\operatorname*{SYT}\left(  \lambda\right)  }$ (since $\mathcal{A}%
=\bigoplus\limits_{\lambda\text{ is a partition of }n}\mathcal{A}%
\mathbf{E}_{\lambda}$). This proves Theorem \ref{thm.specht.FwE.basis}.
\end{proof}

The basis $\left(  \mathbf{F}_{U}w_{U,V}\mathbf{E}_{V}\right)  _{\lambda\text{
is a partition of }n\text{, and }U,V\in\operatorname*{SYT}\left(
\lambda\right)  }$ of $\mathcal{A}$ (when $n!$ is invertible in $\mathbf{k}$)
will be discussed further in Section \ref{sec.bas.FwE}.

\subsubsection{A combinatorial consequence: the sum of all $f^{\lambda}$}

Theorem \ref{thm.specht.FwE.basis} leads to a combinatorial identity in the
vein of Corollary \ref{cor.spechtmod.sumflam2}:

\begin{corollary}
\label{cor.spechtmod.sumflam}We have%
\[
\sum_{\lambda\text{ is a partition of }n}f^{\lambda}=\left(  \text{\# of
involutions }w\in S_{n}\right)  .
\]
(We recall that an \emph{involution} means a permutation $w$ satisfying
$w^{2}=\operatorname*{id}$, or, equivalently, $w^{-1}=w$. Also recall
Definition \ref{def.spechtmod.flam}.)
\end{corollary}

\begin{example}
Let $n=4$. Then, the partitions of $n$ are $\left(  4\right)  $, $\left(
3,1\right)  $, $\left(  2,2\right)  $, $\left(  2,1,1\right)  $ and $\left(
1,1,1,1\right)  $, and the corresponding numbers $f^{\lambda}$ are%
\[
f^{\left(  4\right)  }=1,\ \ \ \ \ \ \ \ \ \ f^{\left(  3,1\right)
}=3,\ \ \ \ \ \ \ \ \ \ f^{\left(  2,2\right)  }%
=2,\ \ \ \ \ \ \ \ \ \ f^{\left(  2,1,1\right)  }%
=3,\ \ \ \ \ \ \ \ \ \ f^{\left(  1,1,1,1\right)  }=1.
\]
Meanwhile, the involutions $w\in S_{4}$ are $\operatorname*{id}$, the six
transpositions $t_{i,j}$ and the three products $t_{1,2}t_{3,4}$ and
$t_{1,3}t_{2,4}$ and $t_{1,4}t_{2,3}$, which make altogether $10$ involutions.
Hence, Corollary \ref{cor.spechtmod.sumflam} says that $1+3+2+3+1=10$. And
this is indeed true.
\end{example}

To prove Corollary \ref{cor.spechtmod.sumflam}, we need a basic lemma about
traces of matrices:

\begin{lemma}
\label{lem.permatrix.trace}Let $M$ be a free $\mathbf{k}$-module with basis
$\left(  m_{i}\right)  _{i\in I}$, where $I$ is a finite set. Let
$f:I\rightarrow I$ be any map. Let $\phi:M\rightarrow M$ be a $\mathbf{k}%
$-linear map that satisfies%
\begin{equation}
\left(  \phi\left(  m_{i}\right)  =m_{f\left(  i\right)  }%
\ \ \ \ \ \ \ \ \ \ \text{for all }i\in I\right)  .
\label{eq.lem.permatrix.trace.1}%
\end{equation}
(This condition actually determines $\phi$ uniquely, but this is immaterial
here.) Then, the trace of this map $\phi$ is%
\[
\operatorname*{Tr}\phi=\left(  \text{\# of all }i\in I\text{ satisfying
}f\left(  i\right)  =i\right)  \cdot1_{\mathbf{k}}.
\]

\end{lemma}

\begin{example}
Let $M=\mathbf{k}^{5}$; this is a free $\mathbf{k}$-module with basis $\left(
e_{1},e_{2},e_{3},e_{4},e_{5}\right)  =\left(  e_{i}\right)  _{i\in\left[
5\right]  }$. Let $f:\left[  5\right]  \rightarrow\left[  5\right]  $ be the
map that sends the numbers $1,2,3,4,5$ to $2,4,3,3,5$, respectively. Let
$\phi:M\rightarrow M$ be the $\mathbf{k}$-linear map that sends each basis
vector $e_{i}$ to $e_{f\left(  i\right)  }$ (that is, that satisfies
$\phi\left(  e_{i}\right)  =e_{f\left(  i\right)  }$ for all $i\in\left[
5\right]  $). Then, Lemma \ref{lem.permatrix.trace} (applied to $I=\left[
5\right]  $ and $m_{i}=e_{i}$) says that%
\[
\operatorname*{Tr}\phi=\underbrace{\left(  \text{\# of all }i\in\left[
5\right]  \text{ satisfying }f\left(  i\right)  =i\right)  }%
_{\substack{=2\\\text{(since these }i\in\left[  5\right]  \text{ are }3\text{
and }5\text{)}}}\cdot\,1_{\mathbf{k}}=2\cdot1_{\mathbf{k}}.
\]
And indeed, this can be easily checked: The matrix that represents the
$\mathbf{k}$-linear map $\phi$ with respect to the basis $\left(  e_{1}%
,e_{2},e_{3},e_{4},e_{5}\right)  $ is%
\[
\left(
\begin{array}
[c]{ccccc}%
0 & 0 & 0 & 0 & 0\\
1_{\mathbf{k}} & 0 & 0 & 0 & 0\\
0 & 0 & 1_{\mathbf{k}} & 1_{\mathbf{k}} & 0\\
0 & 1_{\mathbf{k}} & 0 & 0 & 0\\
0 & 0 & 0 & 0 & 1_{\mathbf{k}}%
\end{array}
\right)  ,
\]
and thus its trace is $0+0+1_{\mathbf{k}}+0+1_{\mathbf{k}}=2\cdot
1_{\mathbf{k}}$. This example illustrates why Lemma \ref{lem.permatrix.trace}
holds in general: The matrix that represents the $\mathbf{k}$-linear map
$\phi$ with respect to the basis $\left(  e_{1},e_{2},e_{3},e_{4}%
,e_{5}\right)  $ (or $\left(  m_{i}\right)  _{i\in I}$ in the general case)
consists only of $0$'s and $1_{\mathbf{k}}$'s, and the $1_{\mathbf{k}}$'s lie
precisely in the $\left(  f\left(  i\right)  ,i\right)  $-th cells of the
matrix (for all $i\in I$). Thus, the diagonal entries of this matrix are $0$'s
and $1_{\mathbf{k}}$'s, with a $1_{\mathbf{k}}$ for each $i\in I$ that
satisfies $f\left(  i\right)  =i$. Hence, the sum of these diagonal entries
(i.e., the trace of the matrix) is a sum of $1_{\mathbf{k}}$'s for each $i\in
I$ that satisfies $f\left(  i\right)  =i$. Of course, this sum is precisely
$\left(  \text{\# of all }i\in I\text{ satisfying }f\left(  i\right)
=i\right)  \cdot1_{\mathbf{k}}$.
\end{example}

\begin{proof}
[Proof of Lemma \ref{lem.permatrix.trace}.]Let us first assume that $I=\left[
r\right]  $ for some $r\in\mathbb{N}$. (We will discuss the general case
further below.)

We have $I=\left[  r\right]  $. Thus, the family $\left(  m_{i}\right)  _{i\in
I}=\left(  m_{i}\right)  _{i\in\left[  r\right]  }$ is just an $r$-tuple
$\left(  m_{1},m_{2},\ldots,m_{r}\right)  $. Hence, $\left(  m_{1}%
,m_{2},\ldots,m_{r}\right)  $ is a basis of $M$ (since $\left(  m_{i}\right)
_{i\in I}$ is a basis of $M$). Let $A$ be the matrix that represents the
$\mathbf{k}$-linear map $\phi:M\rightarrow M$ with respect to this basis
$\left(  m_{1},m_{2},\ldots,m_{r}\right)  $. Then, $\operatorname*{Tr}%
\phi=\operatorname*{Tr}A$ (by the definition of the trace of $\phi$).

For any vector $u\in M$ and any $i\in\left[  r\right]  $, we let $\left[
m_{i}\right]  u$ denote the $m_{i}$-coordinate of $u$ with respect to the
basis $\left(  m_{1},m_{2},\ldots,m_{r}\right)  $ (that is, the coefficient of
$m_{i}$ when $u$ is expanded in the basis $\left(  m_{1},m_{2},\ldots
,m_{r}\right)  $ of $M$). Then, the matrix $A$ can be written explicitly as
\[
A=\left(  \left[  m_{i}\right]  \left(  \phi\left(  m_{j}\right)  \right)
\right)  _{i,j\in\left[  r\right]  }%
\]
(since $A$ is the matrix that represents the $\mathbf{k}$-linear map
$\phi:M\rightarrow M$ with respect to the basis $\left(  m_{1},m_{2}%
,\ldots,m_{r}\right)  $). Hence,%
\[
\operatorname*{Tr}A=\operatorname*{Tr}\left(  \left(  \left[  m_{i}\right]
\left(  \phi\left(  m_{j}\right)  \right)  \right)  _{i,j\in\left[  r\right]
}\right)  =\sum_{i=1}^{r}\left[  m_{i}\right]  \left(  \phi\left(
m_{i}\right)  \right)
\]
(since the trace of a matrix is defined as the sum of its diagonal entries).
In view of $\operatorname*{Tr}\phi=\operatorname*{Tr}A$, we can rewrite this
as%
\begin{align*}
\operatorname*{Tr}\phi &  =\sum_{i=1}^{r}\left[  m_{i}\right]  \left(
\underbrace{\phi\left(  m_{i}\right)  }_{\substack{=m_{f\left(  i\right)
}\\\text{(by (\ref{eq.lem.permatrix.trace.1}))}}}\right)  =\sum_{i=1}%
^{r}\left[  m_{i}\right]  \left(  m_{f\left(  i\right)  }\right)  =\sum
_{i\in\left[  r\right]  }\left[  m_{i}\right]  \left(  m_{f\left(  i\right)
}\right) \\
&  =\sum_{\substack{i\in\left[  r\right]  ;\\f\left(  i\right)  =i}%
}\underbrace{\left[  m_{i}\right]  \left(  m_{f\left(  i\right)  }\right)
}_{\substack{=1_{\mathbf{k}}\\\text{(since }f\left(  i\right)  =i\text{, and
thus}\\\text{the }m_{i}\text{-coordinate of }m_{f\left(  i\right)  }\text{ is
}1_{\mathbf{k}}\text{)}}}+\sum_{\substack{i\in\left[  r\right]  ;\\f\left(
i\right)  \neq i}}\underbrace{\left[  m_{i}\right]  \left(  m_{f\left(
i\right)  }\right)  }_{\substack{=0\\\text{(since }f\left(  i\right)  \neq
i\text{, and thus}\\\text{the }m_{i}\text{-coordinate of }m_{f\left(
i\right)  }\text{ is }0\text{)}}}\\
&  \ \ \ \ \ \ \ \ \ \ \ \ \ \ \ \ \ \ \ \ \left(
\begin{array}
[c]{c}%
\text{since each }i\in\left[  r\right]  \text{ satisfies either }f\left(
i\right)  =i\\
\text{or }f\left(  i\right)  \neq i\text{, but not both}%
\end{array}
\right) \\
&  =\sum_{\substack{i\in\left[  r\right]  ;\\f\left(  i\right)  =i}%
}1_{\mathbf{k}}+\underbrace{\sum_{\substack{i\in\left[  r\right]  ;\\f\left(
i\right)  \neq i}}0}_{=0}=\sum_{\substack{i\in\left[  r\right]  ;\\f\left(
i\right)  =i}}1_{\mathbf{k}}\\
&  =\left(  \text{\# of all }i\in\left[  r\right]  \text{ satisfying }f\left(
i\right)  =i\right)  \cdot1_{\mathbf{k}}\\
&  =\left(  \text{\# of all }i\in I\text{ satisfying }f\left(  i\right)
=i\right)  \cdot1_{\mathbf{k}}\ \ \ \ \ \ \ \ \ \ \left(  \text{since }\left[
r\right]  =I\right)  .
\end{align*}
Thus, Lemma \ref{lem.permatrix.trace} is proved in the particular case when
$I=\left[  r\right]  $ for some $r\in\mathbb{N}$.

The general case can easily be reduced to this particular case: We just rename
the elements of the set $I$ as $1,2,\ldots,r$ (where $r=\left\vert
I\right\vert $), and accordingly change the map $f:I\rightarrow I$ (by
replacing both its inputs and its output values by $1,2,\ldots,r$%
).\ \ \ \ \footnote{To make this rigorous: We pick any bijection
$\beta:I\rightarrow\left[  r\right]  $, and we define a new map $f^{\prime
}:\left[  r\right]  \rightarrow\left[  r\right]  $ by $f^{\prime}:=\beta\circ
f\circ\beta^{-1}$. Then, we apply Lemma \ref{lem.permatrix.trace} to the set
$\left[  r\right]  $, the basis $\left(  m_{\beta^{-1}\left(  i\right)
}\right)  _{i\in\left[  r\right]  }$ and the map $f^{\prime}:\left[  r\right]
\rightarrow\left[  r\right]  $ instead of the set $I$, the basis $\left(
m_{i}\right)  _{i\in I}$ and the map $f:I\rightarrow I$. This is allowed,
since we have $\phi\left(  m_{\beta^{-1}\left(  i\right)  }\right)
=m_{\beta^{-1}\left(  f^{\prime}\left(  i\right)  \right)  }$ for all
$i\in\left[  r\right]  $ (an easy consequence of
(\ref{eq.lem.permatrix.trace.1}) and $f^{\prime}=\beta\circ f\circ\beta^{-1}%
$).} These operations do not change the number%
\[
\left(  \text{\# of all }i\in I\text{ satisfying }f\left(  i\right)
=i\right)  ,
\]
since the map $f$ is modified along with the set $I$.\ \ \ \ \footnote{To make
this rigorous: We have%
\[
\left(  \text{\# of all }i\in I\text{ satisfying }f\left(  i\right)
=i\right)  =\left(  \text{\# of all }i\in\left[  r\right]  \text{ satisfying
}f^{\prime}\left(  i\right)  =i\right)  ,
\]
where $f^{\prime}$ is defined as in the previous footnote.} Hence, the claim
of Lemma \ref{lem.permatrix.trace} remains unaffected under this operation.
But as a result of this operation, we obtain the situation when $I=\left[
r\right]  $, and we have already proved Lemma \ref{lem.permatrix.trace} in
this situation. Hence, Lemma \ref{lem.permatrix.trace} must hold in the
original (fully general) situation as well.
\end{proof}

\begin{proof}
[Proof of Corollary \ref{cor.spechtmod.sumflam}.]Let $\mathbf{k}$ be the field
$\mathbb{Q}$. Then, $n!$ is invertible in $\mathbf{k}$. Hence, Theorem
\ref{thm.specht.FwE.basis} shows that the family $\left(  \mathbf{F}%
_{U}w_{U,V}\mathbf{E}_{V}\right)  _{\lambda\text{ is a partition of }n\text{,
and }U,V\in\operatorname*{SYT}\left(  \lambda\right)  }$ is a basis of the
$\mathbf{k}$-module $\mathcal{A}$.

The antipode $S:\mathbf{k}\left[  S_{n}\right]  \rightarrow\mathbf{k}\left[
S_{n}\right]  $ is a $\mathbf{k}$-linear map from $\mathcal{A}$ to
$\mathcal{A}$. We shall compute its trace $\operatorname*{Tr}S$ in two
different ways. Both times, we will apply Lemma \ref{lem.permatrix.trace} to
$M=\mathcal{A}$ and $\phi=S$, but we will use two different bases $\left(
m_{i}\right)  _{i\in I}$ and two different maps $f$.

Here are our two computations of $\operatorname*{Tr}S$:

\begin{statement}
\textit{Claim 1:} We have
\[
\operatorname*{Tr}S=\left(  \text{\# of involutions }w\in S_{n}\right)  .
\]

\end{statement}

\begin{proof}
[Proof of Claim 1.]Define a map $f:S_{n}\rightarrow S_{n}$ by setting%
\[
f\left(  w\right)  =w^{-1}\ \ \ \ \ \ \ \ \ \ \text{for all }w\in S_{n}.
\]
The $\mathbf{k}$-module $\mathcal{A}=\mathbf{k}\left[  S_{n}\right]  $ is free
with basis $\left(  w\right)  _{w\in S_{n}}$. The map $f:S_{n}\rightarrow
S_{n}$ and the $\mathbf{k}$-linear map $S:\mathcal{A}\rightarrow\mathcal{A}$
satisfy%
\begin{equation}
S\left(  w\right)  =f\left(  w\right)  \ \ \ \ \ \ \ \ \ \ \text{for all }w\in
S_{n} \label{pf.cor.spechtmod.sumflam.11}%
\end{equation}
(since both $S\left(  w\right)  $ and $f\left(  w\right)  $ are defined to be
$w^{-1}$). In other words,%
\[
S\left(  i\right)  =f\left(  i\right)  \ \ \ \ \ \ \ \ \ \ \text{for all }i\in
S_{n}%
\]
(here, we have renamed the variable $w$ as $i$ in
(\ref{pf.cor.spechtmod.sumflam.11})). Hence, Lemma \ref{lem.permatrix.trace}
(applied to $M=\mathcal{A}$ and $I=S_{n}$ and $\left(  m_{i}\right)  _{i\in
I}=\left(  w\right)  _{w\in S_{n}}$ and $\phi=S$) yields%
\begin{align*}
\operatorname*{Tr}S  &  =\left(  \text{\# of all }i\in S_{n}\text{ satisfying
}f\left(  i\right)  =i\right)  \cdot\underbrace{1_{\mathbf{k}}}%
_{\substack{=1_{\mathbb{Q}}\\\text{(since }\mathbf{k}=\mathbb{Q}\text{)}}}\\
&  =\left(  \text{\# of all }i\in S_{n}\text{ satisfying }f\left(  i\right)
=i\right) \\
&  =\left(  \text{\# of all }w\in S_{n}\text{ satisfying }\underbrace{f\left(
w\right)  }_{\substack{=w^{-1}\\\text{(by the definition of }f\text{)}%
}}=w\right) \\
&  \ \ \ \ \ \ \ \ \ \ \ \ \ \ \ \ \ \ \ \ \left(  \text{here, we have renamed
the index }i\text{ as }w\right) \\
&  =\left(  \text{\# of all }w\in S_{n}\text{ satisfying }w^{-1}=w\right) \\
&  =\left(  \text{\# of involutions }w\in S_{n}\right)
\end{align*}
(since the permutations $w\in S_{n}$ satisfying $w^{-1}=w$ are precisely the
involutions $w\in S_{n}$). This proves Claim 1.
\end{proof}

\begin{statement}
\textit{Claim 2:} We have
\[
\operatorname*{Tr}S=\sum_{\lambda\text{ is a partition of }n}f^{\lambda}.
\]

\end{statement}

\begin{proof}
[Proof of Claim 2.]For each partition $\lambda$ of $n$, we let
$\operatorname*{SYT}\left(  \lambda\right)  $ denote the set of all standard
$n$-tableaux of shape $Y\left(  \lambda\right)  $. Then, for each partition
$\lambda$ of $n$, we have%
\[
\left\vert \operatorname*{SYT}\left(  \lambda\right)  \right\vert =f^{\lambda}%
\]
(as we saw in the proof of Proposition \ref{prop.specht.FwE.basis-lam}).

Let
\[
I=\bigsqcup\limits_{\lambda\text{ is a partition of }n}\left(
\operatorname*{SYT}\left(  \lambda\right)  \times\operatorname*{SYT}\left(
\lambda\right)  \right)  .
\]
In other words, let $I$ be the set of all pairs $\left(  U,V\right)  $, where
$\lambda$ is a partition of $n$ and where $U,V\in\operatorname*{SYT}\left(
\lambda\right)  $.

If $\left(  U,V\right)  \in I$ is any pair, then the tableaux $U$ and $V$ have
the same shape (since $\left(  U,V\right)  \in I$ means that $U,V\in
\operatorname*{SYT}\left(  \lambda\right)  $ for some partition $\lambda$ of
$n$, and therefore both tableaux $U$ and $V$ have shape $Y\left(
\lambda\right)  $), and thus the permutation $w_{U,V}$ is well-defined. Thus,
for each $\left(  U,V\right)  \in I$, we can define a vector%
\[
m_{\left(  U,V\right)  }:=\mathbf{F}_{U}w_{U,V}\mathbf{E}_{V}\in\mathcal{A}.
\]
Thus, we have defined a vector $m_{i}\in\mathcal{A}$ for each $i\in I$ (since
each $i\in I$ can be written as a pair $\left(  U,V\right)  $). We therefore
have defined a family $\left(  m_{i}\right)  _{i\in I}$ of vectors in
$\mathcal{A}$. This family%
\begin{align*}
\left(  m_{i}\right)  _{i\in I}  &  =\left(  m_{\left(  U,V\right)  }\right)
_{\left(  U,V\right)  \in I}\ \ \ \ \ \ \ \ \ \ \left(
\begin{array}
[c]{c}%
\text{here, we have renamed the}\\
\text{index }i\text{ as }\left(  U,V\right)
\end{array}
\right) \\
&  =\left(  \mathbf{F}_{U}w_{U,V}\mathbf{E}_{V}\right)  _{\left(  U,V\right)
\in I}\ \ \ \ \ \ \ \ \ \ \left(
\begin{array}
[c]{c}%
\text{since }m_{\left(  U,V\right)  }=\mathbf{F}_{U}w_{U,V}\mathbf{E}_{V}\\
\text{for each }\left(  U,V\right)  \in I
\end{array}
\right) \\
&  =\left(  \mathbf{F}_{U}w_{U,V}\mathbf{E}_{V}\right)  _{\lambda\text{ is a
partition of }n\text{, and }U,V\in\operatorname*{SYT}\left(  \lambda\right)
}\\
&  \ \ \ \ \ \ \ \ \ \ \ \ \ \ \ \ \ \ \ \ \left(
\begin{array}
[c]{c}%
\text{since our definition of }I\text{ allows us to rewrite}\\
\text{the subscript \textquotedblleft}\left(  U,V\right)  \in
I\text{\textquotedblright}\\
\text{as \textquotedblleft}\lambda\text{ is a partition of }n\text{, and
}U,V\in\operatorname*{SYT}\left(  \lambda\right)  \text{\textquotedblright}%
\end{array}
\right)
\end{align*}
is a basis of the $\mathbf{k}$-module $\mathcal{A}$ (since we saw above that
the family \newline$\left(  \mathbf{F}_{U}w_{U,V}\mathbf{E}_{V}\right)
_{\lambda\text{ is a partition of }n\text{, and }U,V\in\operatorname*{SYT}%
\left(  \lambda\right)  }$ is a basis of the $\mathbf{k}$-module $\mathcal{A}$).

Define the map $f:I\rightarrow I$ by setting%
\[
f\left(  U,V\right)  =\left(  V,U\right)  \ \ \ \ \ \ \ \ \ \ \text{for all
}\left(  U,V\right)  \in I.
\]
(This is well-defined, since $U,V\in\operatorname*{SYT}\left(  \lambda\right)
$ clearly entails $V,U\in\operatorname*{SYT}\left(  \lambda\right)  $.) Then,
for each $\left(  U,V\right)  \in I$, we have%
\begin{align*}
S\left(  m_{\left(  U,V\right)  }\right)   &  =S\left(  \mathbf{F}_{U}%
w_{U,V}\mathbf{E}_{V}\right)  \ \ \ \ \ \ \ \ \ \ \left(  \text{since
}m_{\left(  U,V\right)  }=\mathbf{F}_{U}w_{U,V}\mathbf{E}_{V}\right) \\
&  =\mathbf{F}_{V}w_{V,U}\mathbf{E}_{U}\ \ \ \ \ \ \ \ \ \ \left(
\begin{array}
[c]{c}%
\text{by Proposition \ref{prop.specht.FwE.S}, since }U,V\in\operatorname*{SYT}%
\left(  \lambda\right) \\
\text{for some partition }\lambda\text{ of }n
\end{array}
\right) \\
&  =m_{\left(  V,U\right)  }\ \ \ \ \ \ \ \ \ \ \left(
\begin{array}
[c]{c}%
\text{since the definition of }m_{\left(  V,U\right)  }\\
\text{yields }m_{\left(  V,U\right)  }=\mathbf{F}_{V}w_{V,U}\mathbf{E}_{U}%
\end{array}
\right) \\
&  =m_{f\left(  U,V\right)  }\ \ \ \ \ \ \ \ \ \ \left(  \text{since }\left(
V,U\right)  =f\left(  U,V\right)  \text{ (because }f\left(  U,V\right)
=\left(  V,U\right)  \text{)}\right)  .
\end{align*}
Renaming the variable $\left(  U,V\right)  $ as $i$ in this sentence, we
obtain the following: For each $i\in I$, we have $S\left(  m_{i}\right)
=m_{f\left(  i\right)  }$.

Hence, we can apply Lemma \ref{lem.permatrix.trace} to $M=\mathcal{A}$ and
$\phi=S$. Thus, we obtain%
\begin{align*}
\operatorname*{Tr}S  &  =\left(  \text{\# of all }i\in I\text{ satisfying
}f\left(  i\right)  =i\right)  \cdot\underbrace{1_{\mathbf{k}}}%
_{\substack{=1_{\mathbb{Q}}\\\text{(since }\mathbf{k}=\mathbb{Q}\text{)}}}\\
&  =\left(  \text{\# of all }i\in I\text{ satisfying }f\left(  i\right)
=i\right) \\
&  =\left(  \text{\# of all }\left(  U,V\right)  \in I\text{ satisfying
}\underbrace{f\left(  U,V\right)  }_{\substack{=\left(  V,U\right)
\\\text{(by the definition of }f\text{)}}}=\left(  U,V\right)  \right) \\
&  \ \ \ \ \ \ \ \ \ \ \ \ \ \ \ \ \ \ \ \ \left(
\begin{array}
[c]{c}%
\text{here, we have renamed the index }i\text{ as }\left(  U,V\right)
\text{,}\\
\text{since each }i\in I\text{ is a pair of the form }\left(  U,V\right)
\end{array}
\right) \\
&  =\left(  \text{\# of all }\left(  U,V\right)  \in I\text{ satisfying
}\left(  V,U\right)  =\left(  U,V\right)  \right) \\
&  =\left(  \text{\# of all }\left(  U,V\right)  \in I\text{ satisfying
}U=V\right) \\
&  \ \ \ \ \ \ \ \ \ \ \ \ \ \ \ \ \ \ \ \ \left(  \text{since the condition
\textquotedblleft}\left(  V,U\right)  =\left(  U,V\right)
\text{\textquotedblright\ is equivalent to \textquotedblleft}%
U=V\text{\textquotedblright}\right) \\
&  =\left(  \text{\# of all }\left(  U,V\right)  \in\bigsqcup\limits_{\lambda
\text{ is a partition of }n}\left(  \operatorname*{SYT}\left(  \lambda\right)
\times\operatorname*{SYT}\left(  \lambda\right)  \right)  \text{ satisfying
}U=V\right) \\
&  \ \ \ \ \ \ \ \ \ \ \ \ \ \ \ \ \ \ \ \ \left(  \text{since }%
I=\bigsqcup\limits_{\lambda\text{ is a partition of }n}\left(
\operatorname*{SYT}\left(  \lambda\right)  \times\operatorname*{SYT}\left(
\lambda\right)  \right)  \right) \\
&  =\sum_{\lambda\text{ is a partition of }n}\underbrace{\left(  \text{\# of
all }\left(  U,V\right)  \in\operatorname*{SYT}\left(  \lambda\right)
\times\operatorname*{SYT}\left(  \lambda\right)  \text{ satisfying
}U=V\right)  }_{\substack{=\left(  \text{\# of all }U\in\operatorname*{SYT}%
\left(  \lambda\right)  \right)  \\\text{(since a pair }\left(  U,V\right)
\text{ satisfying }U=V\\\text{is uniquely determined by its first entry
}U\text{)}}}\\
&  =\sum_{\lambda\text{ is a partition of }n}\underbrace{\left(  \text{\# of
all }U\in\operatorname*{SYT}\left(  \lambda\right)  \right)  }_{=\left\vert
\operatorname*{SYT}\left(  \lambda\right)  \right\vert =f^{\lambda}}%
=\sum_{\lambda\text{ is a partition of }n}f^{\lambda}.
\end{align*}
This proves Claim 2.
\end{proof}

Now, Claim 2 yields%
\[
\sum_{\lambda\text{ is a partition of }n}f^{\lambda}=\operatorname*{Tr}%
S=\left(  \text{\# of involutions }w\in S_{n}\right)
\]
(by Claim 1). Thus, Corollary \ref{cor.spechtmod.sumflam} is proved.
\end{proof}

Corollary \ref{cor.spechtmod.sumflam} is not as fundamental as Corollary
\ref{cor.spechtmod.sumflam2}, but still has uses in the combinatorics of
permutations. It has combinatorial proofs as well -- see, e.g.,
\cite[Proposition 1.3.2 and the paragraph after its proof]{vL-RSK},
\cite[Exercise 5.7]{MenRem15}, \cite[\S 4.3, Exercise 6 and the paragraph
above it]{Fulton97}, \cite[Theorem 8.26 b.]{Aigner07}, \cite[Exercise 3.12.7
(c)]{Sagan01}, \cite[Exercise 3.2.12]{Prasad-rep}, \cite[\S 5.1.4, Corollary
B]{Knuth-TAoCP3}, \cite[Corollary 7.13.9]{Stanley-EC2}, \cite{Beissi87}.

The following counting exercise gives some equivalent expressions for the
right hand side of Corollary \ref{cor.spechtmod.sumflam}:

\begin{exercise}
\textbf{(a)} \fbox{1} Prove that
\begin{align*}
\left(  \text{\# of involutions }w\in S_{n}\right)   &  =\sum_{k=0}^{n}%
\dbinom{n}{2k}\cdot\underbrace{1\cdot3\cdot5\cdot\cdots\cdot\left(
2k-1\right)  }_{\substack{\text{This is the product of all}\\\text{odd
integers from }1\text{ to }2k-1}}\\
&  =\sum_{k=0}^{n}\dfrac{\dbinom{n}{2k}\dbinom{2k}{k}k!}{2^{k}}.
\end{align*}

\textbf{(b)} \fbox{1} Let $t\left(  m\right)  :=\left(  \text{\# of
involutions }w\in S_{m}\right)  $ for each $m\in\mathbb{N}$. Prove that
$t\left(  n\right)  =t\left(  n-1\right)  +\left(  n-1\right)  t\left(
n-2\right)  $ if $n\geq2$.
\end{exercise}

The numbers $t\left(  n\right)  =\left(  \text{\# of involutions }w\in
S_{n}\right)  $ appear in the OEIS as \href{https://oeis.org/A000085}{Sequence
A000085}. See Remark \ref{rmk.AWS.OEIS} for a table of the first few.

\subsection{\label{sec.specht.TD/Kerpi}The quotient presentation of a Specht
module}

\subsubsection{\label{subsec.specht.TD/Kerpi.thms}Definitions and theorems}

We have already seen several avatars of a Specht module $\mathcal{S}^{D}$: the
original definition (Definition \ref{def.spechtmod.spechtmod}), the left ideal
avatar (Theorem \ref{thm.spechtmod.leftideal} \textbf{(b)}), the
Specht--Vandermonde avatar (Theorem \ref{thm.spechtmod.vdm} \textbf{(b)}) and
the letterplace avatar (Theorem \ref{thm.spechtmod.det} \textbf{(b)}). All
these avatars have one thing in common: They present $\mathcal{S}^{D}$ as a
submodule of some larger module (the Young module or the symmetric group
algebra $\mathbf{k}\left[  S_{n}\right]  $ itself or a polynomial ring). What
we have not seen so far is a presentation of $\mathcal{S}^{D}$ as a
\textbf{quotient} module. Ideally, we would like to write $\mathcal{S}^{D}$ as
a quotient of a free module (over $\mathbf{k}$ or, better, over $\mathcal{A}$)
by an explicitly generated submodule; this would be a presentation of
$\mathcal{S}^{D}$ by generators and relations.

This is what we shall do now. In part, this is very easy: The Specht module
$\mathcal{S}^{D}$ is already defined by a set of generators (viz., the
polytabloids $\mathbf{e}_{T}$). Can we describe the relations between these
generators? I don't know how to do this in general (the Garnir relations from
Theorem \ref{thm.garnir.grel} always hold, but they are not always sufficient
to derive all relations between the $\mathbf{e}_{T}$), but at least it can be
done when $D$ is a skew Young diagram.

In this section, we will see how this is done. We begin with some generalities
that make sense for any diagram $D$ of size $\left\vert D\right\vert =n$.

Let $\operatorname*{Tab}\left(  D\right)  $ be the set of all $n$-tableaux of
shape $D$. As we said, the Specht module $\mathcal{S}^{D}$ is spanned (as a
$\mathbf{k}$-module) by the polytabloids $\mathbf{e}_{T}$ with $T\in
\operatorname*{Tab}\left(  D\right)  $. There are $n!$ of these polytabloids
(since there are $n!$ many $n$-tableaux $T\in\operatorname*{Tab}\left(
D\right)  $). Thus, we can write $\mathcal{S}^{D}$ as a quotient of a free
$\mathbf{k}$-module of rank $n!$, which has a basis $\left(  e_{T}\right)
_{T\in\operatorname*{Tab}\left(  D\right)  }$ (do not mistake the $e_{T}$ for
the $\mathbf{e}_{T}$). We shall call this free $\mathbf{k}$-module
$\mathcal{T}^{D}$. The best way to define it (which gives it a left
$\mathbf{k}\left[  S_{n}\right]  $-module structure \textquotedblleft for
free\textquotedblright) is as a permutation module:

\begin{definition}
\label{def.specht.TD}Let $D$ be any diagram with $\left\vert D\right\vert =n$.

As we know from Definition \ref{def.tableau.Sn-act}, the group $S_{n}$ acts on
the set $\left\{  n\text{-tableaux of shape }D\right\}  $. Let us denote this
set by $\operatorname*{Tab}\left(  D\right)  $. Thus, $\operatorname*{Tab}%
\left(  D\right)  =\left\{  n\text{-tableaux of shape }D\right\}  $ is a left
$S_{n}$-set.

Let $\mathcal{T}^{D}$ be the permutation module\footnotemark\ $\mathbf{k}%
^{\left(  \operatorname*{Tab}\left(  D\right)  \right)  }$ corresponding to
this left $S_{n}$-set $\operatorname*{Tab}\left(  D\right)  $. Thus,
$\mathcal{T}^{D}$ is the free $\mathbf{k}$-module $\mathbf{k}^{\left(
\operatorname*{Tab}\left(  D\right)  \right)  }$ with its standard basis
$\left(  e_{T}\right)  _{T\in\operatorname*{Tab}\left(  D\right)  }$ (do not
mistake these basis vectors $e_{T}$ for the polytabloids $\mathbf{e}_{T}$),
and it is furthermore a left $\mathbf{k}\left[  S_{n}\right]  $-module whose
action satisfies%
\begin{equation}
we_{T}=e_{wT}\ \ \ \ \ \ \ \ \ \ \text{for each }w\in S_{n}\text{ and }%
T\in\operatorname*{Tab}\left(  D\right)  \label{eq.def.specht.TD.weT}%
\end{equation}
(by a general property of permutation modules -- see
(\ref{pf.exa.mod.kG-on-kX.ekey})). Hence, it is a representation of $S_{n}$.
\end{definition}

\footnotetext{See Example \ref{exa.mod.kG-on-kX} for the definition of a
permutation module.}The basis vectors $e_{T}$ of $\mathcal{T}^{D}$ are not the
polytabloids $\mathbf{e}_{T}$ from the Specht module $\mathcal{S}^{D}$, even
though the rule for the $S_{n}$-action is the same for both types of vectors:
The polytabloids $\mathbf{e}_{T}\in\mathcal{S}^{D}$ also satisfy
\begin{equation}
w\mathbf{e}_{T}=\mathbf{e}_{wT}\ \ \ \ \ \ \ \ \ \ \text{for each }w\in
S_{n}\text{ and }T\in\operatorname*{Tab}\left(  D\right)
\label{eq.def.specht.TD.weTSD}%
\end{equation}
(since Lemma \ref{lem.spechtmod.submod} \textbf{(a)} (applied to $u=w$) yields
$\mathbf{e}_{wT}=w\mathbf{e}_{T}$). The equalities (\ref{eq.def.specht.TD.weT}%
) and (\ref{eq.def.specht.TD.weTSD}) look alike, and this will be used later on.

As a left $\mathbf{k}\left[  S_{n}\right]  $-module (i.e., as an $S_{n}%
$-representation), $\mathcal{T}^{D}$ is not particularly interesting:

\begin{proposition}
\label{prop.specht.TD.isoA}Let $D$ be any diagram with $\left\vert
D\right\vert =n$. Consider the left regular representation $\mathbf{k}\left[
S_{n}\right]  $ of $S_{n}$. Then,
\[
\mathcal{T}^{D}\cong\mathbf{k}\left[  S_{n}\right]
\ \ \ \ \ \ \ \ \ \ \text{as }S_{n}\text{-representations.}%
\]
Specifically, we can obtain an isomorphism as follows: Fix any $n$-tableau
$P\in\operatorname*{Tab}\left(  D\right)  $. Let $\phi_{P}:\mathbf{k}\left[
S_{n}\right]  \rightarrow\mathcal{T}^{D}$ be the $\mathbf{k}$-linear map that
sends each standard basis vector $w\in S_{n}$ to $e_{wP}$. Then, $\phi_{P}$ is
an isomorphism of $S_{n}$-representations.
\end{proposition}

\begin{proof}
First, we shall show that $\phi_{P}$ is an isomorphism of $\mathbf{k}%
$-modules. Indeed, $P$ is an $n$-tableau of shape $D$ (since $P\in
\operatorname*{Tab}\left(  D\right)  $). Hence, the map%
\begin{align*}
f:S_{n}  &  \rightarrow\left\{  n\text{-tableaux of shape }D\right\}  ,\\
w  &  \mapsto w\rightharpoonup P
\end{align*}
is a bijection (by Proposition \ref{prop.tableau.Sn-act.1} \textbf{(d)},
applied to $T=P$). In other words, the map%
\begin{align*}
f:S_{n}  &  \rightarrow\operatorname*{Tab}\left(  D\right)  ,\\
w  &  \mapsto wP
\end{align*}
is a bijection (since $\left\{  n\text{-tableaux of shape }D\right\}
=\operatorname*{Tab}\left(  D\right)  $ and since $w\rightharpoonup P=wP$ for
each $w\in S_{n}$). Therefore, the family $\left(  e_{wP}\right)  _{w\in
S_{n}}$ is a reindexing of the family $\left(  e_{T}\right)  _{T\in
\operatorname*{Tab}\left(  D\right)  }$. Hence, the family $\left(
e_{wP}\right)  _{w\in S_{n}}$ is a basis of the $\mathbf{k}$-module
$\mathcal{T}^{D}$ (since the family $\left(  e_{T}\right)  _{T\in
\operatorname*{Tab}\left(  D\right)  }$ is a basis of $\mathcal{T}^{D}$).

Now, the $\mathbf{k}$-linear map $\phi_{P}:\mathbf{k}\left[  S_{n}\right]
\rightarrow\mathcal{T}^{D}$ sends the basis $\left(  w\right)  _{w\in S_{n}}$
of $\mathbf{k}\left[  S_{n}\right]  $ to the basis $\left(  e_{wP}\right)
_{w\in S_{n}}$ of $\mathcal{T}^{D}$ (since it sends each $w$ to $e_{wP}$ by
its definition). Hence, this map is an isomorphism of $\mathbf{k}$-modules
(because any $\mathbf{k}$-linear map that sends a basis of its domain to a
basis of its target must be an isomorphism), and thus is invertible.

Moreover, it is easy to show that the map $\phi_{P}$ is $S_{n}$-equivariant
(this essentially follows from (\ref{eq.def.specht.TD.weT})\footnote{In more
detail: We must show that $\phi_{P}$ is $S_{n}$-equivariant. In other words,
we must prove that all $g\in S_{n}$ and all $x\in\mathbf{k}\left[
S_{n}\right]  $ satisfy $\phi_{P}\left(  g\rightharpoonup x\right)
=g\rightharpoonup\phi_{P}\left(  x\right)  $. So let us prove this. Let $g\in
S_{n}$ and $x\in\mathbf{k}\left[  S_{n}\right]  $. We must prove the equality
$\phi_{P}\left(  g\rightharpoonup x\right)  =g\rightharpoonup\phi_{P}\left(
x\right)  $. Since both sides of this equality are $\mathbf{k}$-linear in $x$,
we can WLOG assume that $x$ is an element of the standard basis $\left(
w\right)  _{w\in S_{n}}$ of $\mathbf{k}\left[  S_{n}\right]  $. Assume this.
Thus, $x=w$ for some $w\in S_{n}$. Consider this $w$.
\par
From $x=w$, we obtain $\phi_{P}\left(  x\right)  =\phi_{P}\left(  w\right)
=e_{wP}$ (by the definition of $\phi_{P}$) and $g\rightharpoonup
x=g\rightharpoonup w=gw$. From $g\rightharpoonup x=gw$, we obtain $\phi
_{P}\left(  g\rightharpoonup x\right)  =\phi_{P}\left(  gw\right)  =e_{gwP}$
(by the definition of $\phi_{P}$, since $gw\in S_{n}$). But
(\ref{eq.def.specht.TD.weT}) (applied to $g$ and $wP$ instead of $w$ and $T$)
yields $ge_{wP}=e_{gwP}$. Comparing this with $\phi_{P}\left(
g\rightharpoonup x\right)  =e_{gwP}$, we obtain $\phi_{P}\left(
g\rightharpoonup x\right)  =ge_{wP}=g\rightharpoonup e_{wP}$. In view of
$\phi_{P}\left(  x\right)  =e_{wP}$, we can rewrite this as $\phi_{P}\left(
g\rightharpoonup x\right)  =g\rightharpoonup\phi_{P}\left(  x\right)  $. This
completes our proof that $\phi_{P}$ is $S_{n}$-equivariant.}). Hence, this map
$\phi_{P}$ is a morphism of $S_{n}$-representations (since it is $\mathbf{k}%
$-linear). In other words, the map $\phi_{P}$ is a morphism of left
$\mathbf{k}\left[  S_{n}\right]  $-modules (by Proposition
\ref{prop.rep.G-rep.mor=mor}). Since $\phi_{P}$ is invertible, we thus
conclude that $\phi_{P}$ is a left $\mathbf{k}\left[  S_{n}\right]  $-module
isomorphism. In other words, $\phi_{P}$ is an isomorphism of $S_{n}%
$-representations. Hence, $\mathcal{T}^{D}\cong\mathbf{k}\left[  S_{n}\right]
$ as $S_{n}$-representations. Proposition \ref{prop.specht.TD.isoA} is thus proven.

(Alternatively, we could have argued that the above map $f:S_{n}%
\rightarrow\operatorname*{Tab}\left(  D\right)  $ is an isomorphism of left
$S_{n}$-sets, and thus (by the functoriality of the \textquotedblleft
permutation module\textquotedblright\ functor) induces an isomorphism between
the corresponding permutation modules, which are $\mathbf{k}\left[
S_{n}\right]  $ and $\mathcal{T}^{D}$; this latter isomorphism is precisely
our map $\phi_{P}$.)
\end{proof}

Another easy observation is the existence of a surjective morphism of $S_{n}%
$-representations from $\mathcal{T}^{D}$ to the Specht module $\mathcal{S}%
^{D}$:

\begin{proposition}
\label{prop.specht.TD.pi}Let $D$ be any diagram with $\left\vert D\right\vert
=n$. Let $\pi_{D}:\mathcal{T}^{D}\rightarrow\mathcal{S}^{D}$ be the
$\mathbf{k}$-linear map that sends each standard basis vector $e_{T}$ of
$\mathcal{T}^{D}$ to the polytabloid $\mathbf{e}_{T}\in\mathcal{S}^{D}$. Then,
$\pi_{D}$ is a surjective morphism of $S_{n}$-representations.
\end{proposition}

\begin{proof}
Recall that $\mathcal{S}^{D}$ is the span of the polytabloids $\mathbf{e}_{T}$
for all $T\in\operatorname*{Tab}\left(  D\right)  $ (by the definition of
$\mathcal{S}^{D}$). This yields that the map $\pi_{D}$ is
surjective\footnote{In more detail: The $\mathbf{k}$-module $\mathcal{T}^{D}$
is free with basis $\left(  e_{T}\right)  _{T\in\operatorname*{Tab}\left(
D\right)  }$. Hence, $\mathcal{T}^{D}=\operatorname*{span}%
\nolimits_{\mathbf{k}}\left\{  e_{T}\ \mid\ T\in\operatorname*{Tab}\left(
D\right)  \right\}  $. Thus,%
\begin{align*}
\pi_{D}\left(  \mathcal{T}^{D}\right)   &  =\pi_{D}\left(
\operatorname*{span}\nolimits_{\mathbf{k}}\left\{  e_{T}\ \mid\ T\in
\operatorname*{Tab}\left(  D\right)  \right\}  \right) \\
&  =\operatorname*{span}\nolimits_{\mathbf{k}}\left\{  \underbrace{\pi
_{D}\left(  e_{T}\right)  }_{\substack{=\mathbf{e}_{T}\\\text{(by the
definition of }\pi_{D}\text{)}}}\ \mid\ T\in\operatorname*{Tab}\left(
D\right)  \right\} \\
&  \ \ \ \ \ \ \ \ \ \ \ \ \ \ \ \ \ \ \ \ \left(
\begin{array}
[c]{c}%
\text{since the map }\pi_{D}\text{ is }\mathbf{k}\text{-linear, and thus}\\
\text{respects }\mathbf{k}\text{-linear combinations}%
\end{array}
\right) \\
&  =\operatorname*{span}\nolimits_{\mathbf{k}}\left\{  \mathbf{e}_{T}%
\ \mid\ T\in\operatorname*{Tab}\left(  D\right)  \right\} \\
&  =\operatorname*{span}\nolimits_{\mathbf{k}}\left\{  \mathbf{e}_{T}%
\ \mid\ T\text{ is an }n\text{-tableau of shape }D\right\} \\
&  \ \ \ \ \ \ \ \ \ \ \ \ \ \ \ \ \ \ \ \ \left(  \text{since }%
\operatorname*{Tab}\left(  D\right)  \text{ is the set of all }%
n\text{-tableaux of shape }D\right) \\
&  =\mathcal{S}^{D}\ \ \ \ \ \ \ \ \ \ \left(  \text{by the definition of
}\mathcal{S}^{D}\right)  .
\end{align*}
In other words, the map $\pi_{D}$ is surjective.}. It remains to prove that
$\pi_{D}$ is a morphism of $S_{n}$-representations. It clearly suffices to
show that $\pi_{D}$ is $S_{n}$-equivariant (since $\pi_{D}$ is $\mathbf{k}%
$-linear). But this is an easy consequence of comparing
(\ref{eq.def.specht.TD.weT}) with (\ref{eq.def.specht.TD.weTSD}).\footnote{In
more detail: We must prove that $\pi_{D}$ is $S_{n}$-equivariant. In other
words, we must prove that all $g\in S_{n}$ and all $x\in\mathcal{T}^{D}$
satisfy $\pi_{D}\left(  g\rightharpoonup x\right)  =g\rightharpoonup\pi
_{D}\left(  x\right)  $. So let us prove this. Let $g\in S_{n}$ and
$x\in\mathcal{T}^{D}$. We must prove the equality $\pi_{D}\left(
g\rightharpoonup x\right)  =g\rightharpoonup\pi_{D}\left(  x\right)  $. Since
both sides of this equality are $\mathbf{k}$-linear in $x$, we can WLOG assume
that $x$ is an element of the standard basis $\left(  e_{T}\right)
_{T\in\operatorname*{Tab}\left(  D\right)  }$ of $\mathcal{T}^{D}$. Assume
this. Thus, $x=e_{T}$ for some $T\in\operatorname*{Tab}\left(  D\right)  $.
Consider this $T$.
\par
From $x=e_{T}$, we obtain $\pi_{D}\left(  x\right)  =\pi_{D}\left(
e_{T}\right)  =\mathbf{e}_{T}$ (by the definition of $\pi_{D}$) and thus
$g\rightharpoonup\pi_{D}\left(  x\right)  =g\rightharpoonup\mathbf{e}%
_{T}=g\mathbf{e}_{T}=\mathbf{e}_{gT}$ (by (\ref{eq.def.specht.TD.weTSD}),
applied to $w=g$). On the other hand, from $x=e_{T}$, we obtain
$g\rightharpoonup x=g\rightharpoonup e_{T}=ge_{T}=e_{gT}$ (by
(\ref{eq.def.specht.TD.weT}), applied to $w=g$) and therefore $\pi_{D}\left(
g\rightharpoonup x\right)  =\pi_{D}\left(  e_{gT}\right)  =\mathbf{e}_{gT}$
(by the definition of $\pi_{D}$). Comparing this with $g\rightharpoonup\pi
_{D}\left(  x\right)  =\mathbf{e}_{gT}$, we find $\pi_{D}\left(
g\rightharpoonup x\right)  =g\rightharpoonup\pi_{D}\left(  x\right)  $. This
completes our proof that $\pi_{D}$ is $S_{n}$-equivariant.} Thus, the proof of
Proposition \ref{prop.specht.TD.pi} is complete.
\end{proof}

The morphism $\pi_{D}:\mathcal{T}^{D}\rightarrow\mathcal{S}^{D}$ from
Proposition \ref{prop.specht.TD.pi} is surjective (by said proposition), and
thus we have%
\begin{equation}
\mathcal{S}^{D}=\operatorname{Im}\left(  \pi_{D}\right)  =\pi_{D}\left(
\mathcal{T}^{D}\right)  \cong\mathcal{T}^{D}/\operatorname*{Ker}\left(
\pi_{D}\right)  \label{eq.prop.specht.TD.pi.quot}%
\end{equation}
(by the first isomorphism theorem), and this is an isomorphism of $S_{n}%
$-representations (since $\pi_{D}$ is a morphism of $S_{n}$-representations).
This can be viewed as a new avatar of $\mathcal{S}^{D}$. Moreover, since
$\mathcal{T}^{D}$ is a free $\mathbf{k}$-module, this yields a presentation of
the $\mathbf{k}$-module $\mathcal{S}^{D}$ by generators and relations (the
generators being the projections of the standard basis vectors $e_{T}%
\in\mathcal{T}^{D}$ onto $\mathcal{T}^{D}/\operatorname*{Ker}\left(  \pi
_{D}\right)  $, while the relations are given by the elements of
$\operatorname*{Ker}\left(  \pi_{D}\right)  $), provided that we can
explicitly express $\operatorname*{Ker}\left(  \pi_{D}\right)  $ as a span of
some concrete vectors (called \textquotedblleft relators\textquotedblright,
since they determine the relations in this presentation). We do not know how
to do this in general, but we shall do this in the case when $D$ is a skew
Young diagram. The relevant vectors spanning $\operatorname*{Ker}\left(
\pi_{D}\right)  $ will be the so-called \emph{column relators} and
\emph{Garnir relators}. These relators can be defined for arbitrary $D$ and
always belong to $\operatorname*{Ker}\left(  \pi_{D}\right)  $, but they might
not always suffice to span $\operatorname*{Ker}\left(  \pi_{D}\right)  $ when
$D$ is not a skew Young diagram.

We begin with their definition:

\begin{definition}
\label{def.specht.TD.relators}Let $D$ be any diagram with $\left\vert
D\right\vert =n$. We define the following two types of elements of
$\mathcal{T}^{D}$: \medskip

\textbf{(a)} A \emph{column relator} (in $\mathcal{T}^{D}$) shall mean an
element of the form%
\[
e_{T}-\left(  -1\right)  ^{u}e_{uT}\in\mathcal{T}^{D},
\]
where $T\in\operatorname*{Tab}\left(  D\right)  $ and $u\in\mathcal{C}\left(
T\right)  $. \medskip

\textbf{(b)} A \emph{Garnir pentuple} (of shape $D$) shall mean a $5$-tuple
$\left(  T,j,k,Z,L\right)  $ whose entries $T,j,k,Z,L$ satisfy the assumptions
of Theorem \ref{thm.garnir.grel}. Concretely, this means the following:

\begin{itemize}
\item The first entry $T$ of the $5$-tuple is an $n$-tableau of shape $D$
(that is, we have $T\in\operatorname*{Tab}\left(  D\right)  $).

\item The next two entries $j$ and $k$ are two distinct integers.

\item The fourth entry $Z$ is a subset of $\operatorname*{Col}\left(
j,T\right)  \cup\operatorname*{Col}\left(  k,T\right)  $ such that $\left\vert
Z\right\vert >\left\vert D\left\langle j\right\rangle \cup D\left\langle
k\right\rangle \right\vert $. (See Definition \ref{def.diagram.Dj} for the
meaning of $D\left\langle j\right\rangle $.)

\item The final entry $L$ is a left transversal of $S_{n,Z}\cap\mathcal{C}%
\left(  T\right)  $ in $S_{n,Z}$. (See Proposition \ref{prop.intX.basics} for
the definition of $S_{n,Z}$.)
\end{itemize}

\textbf{(c)} A \emph{Garnir relator} (in $\mathcal{T}^{D}$) shall mean an
element of the form%
\[
\sum_{x\in L}\left(  -1\right)  ^{x}e_{xT}\in\mathcal{T}^{D},
\]
where $\left(  T,j,k,Z,L\right)  $ is a Garnir pentuple (of shape $D$).
\end{definition}

\begin{example}
Let $n=8$ and $D=Y\left(  3,3,2\right)  $ for this example. Then: \medskip

\textbf{(a)} There are $8!$ distinct $n$-tableaux $T$ of shape $D$, and for
each such $T$ there are $\left\vert \mathcal{C}\left(  T\right)  \right\vert
=6\cdot6\cdot2=72$ possible permutations $u\in\mathcal{C}\left(  T\right)  $.
Thus, there are altogether $8!\cdot72=2\,903\,040$ many column relators in
$\mathcal{T}^{D}$ (or, more precisely, this is how many pairs $\left(
T,u\right)  $ there are). Here are two of them (using the notation from
Convention \ref{conv.tableau.poetic}):%
\[
e_{123\backslash\backslash456\backslash\backslash78}+e_{423\backslash
\backslash156\backslash\backslash78}%
\]
(corresponding to $T=123\backslash\backslash456\backslash\backslash78$ and
$u=t_{1,4}\in\mathcal{C}\left(  T\right)  $) and%
\[
e_{123\backslash\backslash456\backslash\backslash78}-e_{453\backslash
\backslash126\backslash\backslash78}%
\]
(corresponding to $T=123\backslash\backslash456\backslash\backslash78$ and
$u=t_{1,4}t_{2,5}\in\mathcal{C}\left(  T\right)  $). Note that each column
relator $e_{T}-\left(  -1\right)  ^{u}e_{uT}$ with $u=\operatorname*{id}$ is
simply $0$; the other column relators also have many linear dependencies
(usually). \medskip

\textbf{(b)} There are $8!$ distinct $n$-tableaux $T$ of shape $D$, and for
each such $T$ there are many possible Garnir pentuples $\left(
T,j,k,Z,L\right)  $ that start with $T$. One such pentuple has been given in
Example \ref{exa.garnir.grel.1}:%
\begin{align*}
T  &
=\ytableaushort{1{*(green)2}4,{*(green)6}{*(green)3}7,{*(green)8}5}\ \ \ \ \ \ \ \ \ \ \text{and}%
\ \ \ \ \ \ \ \ \ \ j=1\ \ \ \ \ \ \ \ \ \ \text{and}%
\ \ \ \ \ \ \ \ \ \ k=2\ \ \ \ \ \ \ \ \ \ \text{and}\\
Z  &  =\left\{  2,3,6,8\right\}  \ \ \ \ \ \ \ \ \ \ \text{and}%
\ \ \ \ \ \ \ \ \ \ L=\left\{  \operatorname*{id},\ t_{2,6},\ t_{2,8}%
,\ t_{3,6},\ t_{3,8},\ t_{2,6}t_{3,8}\right\}  .
\end{align*}
Hence,
\[
\sum_{x\in L}\left(  -1\right)  ^{x}e_{xT}=e_{\operatorname*{id}T}%
-e_{t_{2,6}T}-e_{t_{2,8}T}-e_{t_{3,6}T}-e_{t_{3,8}T}+e_{t_{2,6}t_{3,8}T}%
\]
is a Garnir relator in $\mathcal{T}^{D}$.
\end{example}

As we already mentioned, both column and Garnir relators belong to
$\operatorname*{Ker}\left(  \pi_{D}\right)  $:

\begin{proposition}
\label{prop.specht.TD.relators-in-ker}Let $D$ be any diagram with $\left\vert
D\right\vert =n$. Consider the morphism $\pi_{D}:\mathcal{T}^{D}%
\rightarrow\mathcal{S}^{D}$ from Proposition \ref{prop.specht.TD.pi}. Then:
\medskip

\textbf{(a)} Each column relator in $\mathcal{T}^{D}$ belongs to
$\operatorname*{Ker}\left(  \pi_{D}\right)  $. \medskip

\textbf{(b)} Each Garnir relator in $\mathcal{T}^{D}$ belongs to
$\operatorname*{Ker}\left(  \pi_{D}\right)  $.
\end{proposition}

\begin{proof}
\textbf{(a)} This is a consequence of Lemma \ref{lem.spechtmod.submod}
\textbf{(c)} (rewritten in the equivalent form $\mathbf{e}_{T}-\left(
-1\right)  ^{u}\mathbf{e}_{uT}=0$). In more detail:

\begin{fineprint}
Let $\mathfrak{c}$ be a column relator in $\mathcal{T}^{D}$. We must prove
that $\mathfrak{c}\in\operatorname*{Ker}\left(  \pi_{D}\right)  $.

We know that $\mathfrak{c}$ is a column relator. In other words,
$\mathfrak{c}=e_{T}-\left(  -1\right)  ^{u}e_{uT}$ for some $T\in
\operatorname*{Tab}\left(  D\right)  $ and some $u\in\mathcal{C}\left(
T\right)  $ (by the definition of a column relator). Consider these $T$ and
$u$. Then, Lemma \ref{lem.spechtmod.submod} \textbf{(c)} yields $\mathbf{e}%
_{uT}=\left(  -1\right)  ^{u}\mathbf{e}_{T}$. But applying the map $\pi_{D}$
to the equality $\mathfrak{c}=e_{T}-\left(  -1\right)  ^{u}e_{uT}$, we obtain%
\begin{align*}
\pi_{D}\left(  \mathfrak{c}\right)   &  =\pi_{D}\left(  e_{T}-\left(
-1\right)  ^{u}e_{uT}\right) \\
&  =\underbrace{\pi_{D}\left(  e_{T}\right)  }_{\substack{=\mathbf{e}%
_{T}\\\text{(by the definition of }\pi_{D}\text{)}}}-\left(  -1\right)
^{u}\underbrace{\pi_{D}\left(  e_{uT}\right)  }_{\substack{=\mathbf{e}%
_{uT}\\\text{(by the definition of }\pi_{D}\text{)}}%
}\ \ \ \ \ \ \ \ \ \ \left(  \text{since }\pi_{D}\text{ is }\mathbf{k}%
\text{-linear}\right) \\
&  =\mathbf{e}_{T}-\left(  -1\right)  ^{u}\underbrace{\mathbf{e}_{uT}%
}_{=\left(  -1\right)  ^{u}\mathbf{e}_{T}}=\mathbf{e}_{T}-\underbrace{\left(
-1\right)  ^{u}\left(  -1\right)  ^{u}}_{\substack{=\left(  -1\right)
^{2u}=1\\\text{(since }2u\text{ is even)}}}\mathbf{e}_{T}\\
&  =\mathbf{e}_{T}-\mathbf{e}_{T}=0.
\end{align*}
In other words, $\mathfrak{c}\in\operatorname*{Ker}\left(  \pi_{D}\right)  $.
Thus, Proposition \ref{prop.specht.TD.relators-in-ker} \textbf{(a)} is proved.
\end{fineprint}

\textbf{(b)} This is a consequence of the Garnir relation
(\ref{eq.thm.garnir.grel.gar}). In more detail:

\begin{fineprint}
Let $\mathfrak{g}$ be a Garnir relator in $\mathcal{T}^{D}$. We must prove
that $\mathfrak{g}\in\operatorname*{Ker}\left(  \pi_{D}\right)  $.

We know that $\mathfrak{g}$ is a Garnir relator. Thus, $\mathfrak{g}%
=\sum_{x\in L}\left(  -1\right)  ^{x}e_{xT}$ for some Garnir pentuple $\left(
T,j,k,Z,L\right)  $. Consider this $\left(  T,j,k,Z,L\right)  $. We know that
$\left(  T,j,k,Z,L\right)  $ is a Garnir pentuple, so that its entries
$T,j,k,Z,L$ satisfy the assumptions of Theorem \ref{thm.garnir.grel}. Hence,
(\ref{eq.thm.garnir.grel.gar}) says that%
\[
\sum_{x\in L}\left(  -1\right)  ^{x}\mathbf{e}_{xT}=0.
\]

However, applying the map $\pi_{D}$ to the equality $\mathfrak{g}=\sum_{x\in
L}\left(  -1\right)  ^{x}e_{xT}$, we obtain%
\begin{align*}
\pi_{D}\left(  \mathfrak{g}\right)   &  =\pi_{D}\left(  \sum_{x\in L}\left(
-1\right)  ^{x}e_{xT}\right) \\
&  =\sum_{x\in L}\left(  -1\right)  ^{x}\underbrace{\pi_{D}\left(
e_{xT}\right)  }_{\substack{=\mathbf{e}_{xT}\\\text{(by the definition of }%
\pi_{D}\text{)}}}\ \ \ \ \ \ \ \ \ \ \left(  \text{since }\pi_{D}\text{ is
}\mathbf{k}\text{-linear}\right) \\
&  =\sum_{x\in L}\left(  -1\right)  ^{x}\mathbf{e}_{xT}=0.
\end{align*}
That is, $\mathfrak{g}\in\operatorname*{Ker}\left(  \pi_{D}\right)  $. Thus,
Proposition \ref{prop.specht.TD.relators-in-ker} \textbf{(b)} is proved.
\end{fineprint}
\end{proof}

The main result of this section is the following:

\begin{theorem}
\label{thm.specht.TD.garnir-all-span}Let $D$ be a skew Young diagram with
$\left\vert D\right\vert =n$. Consider the morphism $\pi_{D}:\mathcal{T}%
^{D}\rightarrow\mathcal{S}^{D}$ from Proposition \ref{prop.specht.TD.pi}.
Then,%
\begin{align*}
&  \operatorname*{Ker}\left(  \pi_{D}\right) \\
&  =\operatorname*{span}\nolimits_{\mathbf{k}}\left(  \left\{  \text{all
column relators in }\mathcal{T}^{D}\right\}  \cup\left\{  \text{all Garnir
relators in }\mathcal{T}^{D}\right\}  \right)  .
\end{align*}

\end{theorem}

We will prove this theorem in Subsection \ref{subsec.specht.TD/Kerpi.pfs}. But
first, we shall generalize it by replacing the two (finite but rather large)
sets $\left\{  \text{all column relators in }\mathcal{T}^{D}\right\}  $ and
$\left\{  \text{all Garnir relators in }\mathcal{T}^{D}\right\}  $ by certain
subsets $\mathfrak{C}$ and $\mathfrak{G}$ that can be much smaller (we shall
not uniquely specify these subsets but only impose conditions on them). The
following types of subsets do the trick:

\begin{definition}
\label{def.specht.TD.garnir-triasub}Let $D$ be a skew Young diagram with
$\left\vert D\right\vert =n$. \medskip

\textbf{(a)} A set $\mathfrak{C}$ of column relators (in $\mathcal{T}^{D}$) is
said to be \emph{triangulating} if it has the following property: For any
$n$-tableau $T\in\operatorname*{Tab}\left(  D\right)  $ that is not
column-standard, there exists a $u\in\mathcal{C}\left(  T\right)  $ such that
the $n$-tableau $uT$ is column-standard and such that $e_{T}-\left(
-1\right)  ^{u}e_{uT}\in\mathfrak{C}$. \medskip

\textbf{(b)} A Garnir pentuple $\left(  T,j,k,Z,L\right)  $ is said to be
\emph{triangulating} if each $x\in L\setminus\mathcal{C}\left(  T\right)  $
satisfies $\widetilde{xT}>\widetilde{T}$ in the column last letter order. (For
the definition of the column last letter order, see Definition
\ref{def.spechtmod.garnir-order} and the discussion that follows.) \medskip

\textbf{(c)} A set $\mathfrak{G}$ of Garnir relators (in $\mathcal{T}^{D}$) is
said to be \emph{triangulating} if it has the following property: For any
column-standard $n$-tableau $T\in\operatorname*{Tab}\left(  D\right)  $ that
is not standard, there exists a triangulating Garnir pentuple $\left(
T,j,k,Z,L\right)  $ with first entry $T$ such that $\sum_{x\in L}\left(
-1\right)  ^{x}e_{xT}\in\mathfrak{G}$.
\end{definition}

Triangulating sets of column relators and of Garnir relators can be much
smaller than the sets of all such relators: It is enough to have just one
column relator for each tableau $T\in\operatorname*{Tab}\left(  D\right)  $
that is not column-standard, and only one Garnir relator for each
column-standard tableau $T$ that is not standard. Altogether, fewer than $n!$
relators are needed.\Needspace{50pc}

\begin{example}
For this example, let $n=5$ and $D=Y\left(  2,2,1\right)  $. Then, there are
$10$ column-standard $n$-tableaux of shape $D$. Five of these are standard;
the other five are (using the notation from Convention
\ref{conv.tableau.poetic})%
\[
\ytableaushort{12,{*(yellow)4}{*(yellow)3},5}\ , \quad
\ytableaushort{{*(yellow)2}{*(yellow)1},35,4}\ , \quad
\ytableaushort{{*(yellow)2}{*(yellow)1},34,5}\ , \quad
\ytableaushort{{*(yellow)2}{*(yellow)1},{*(yellow)4}{*(yellow)3},5}\ ,
\quad\ytableaushort{{*(yellow)3}{*(yellow)1},{*(yellow)4}{*(yellow)2},5}
\]
(where we used yellow background color to mark the rows that violate
standardness). For a set $\mathfrak{G}$ of Garnir relators (in $\mathcal{T}%
^{D}$) to be triangulating, it must contain -- for each of these five
non-standard tableaux $T$ -- a Garnir relator $\sum_{x\in L}\left(  -1\right)
^{x}e_{xT}$ coming from a triangulating Garnir pentuple $\left(
T,j,k,Z,L\right)  $ with first entry $T$. One possible choice of such Garnir
relators is given as follows (note that $j$ and $k$ are always $1$ and $2$
here, and the set $Z$ is always the set of the entries highlighted in green,
whereas the set $L$ is easily reconstructed from the resulting tableaux $xT$):%
\[%
\begin{tabular}
[c]{|c|c|}\hline
$T$ & Garnir relator $\sum_{x\in L}\left(  -1\right)  ^{x}e_{xT}%
$\\\hline\hline
$\ytableaushort{1{*(green)2},{*(green)4}{*(green)3},{*(green)5}}$ &
$e_{12\backslash\backslash43\backslash\backslash5}-e_{12\backslash
\backslash34\backslash\backslash5}-e_{14\backslash\backslash23\backslash
\backslash5}-e_{15\backslash\backslash43\backslash\backslash2}-e_{12\backslash
\backslash45\backslash\backslash3}+e_{14\backslash\backslash25\backslash
\backslash3}$\\\hline
$\ytableaushort{{*(green)2}{*(green)1},{*(green)3}5,{*(green)4}}$ &
$e_{21\backslash\backslash35\backslash\backslash4}-e_{12\backslash
\backslash35\backslash\backslash4}-e_{23\backslash\backslash15\backslash
\backslash4}-e_{24\backslash\backslash35\backslash\backslash1}$\\\hline
$\ytableaushort{{*(green)2}{*(green)1},{*(green)3}4,{*(green)5}}$ &
$e_{21\backslash\backslash34\backslash\backslash5}-e_{12\backslash
\backslash34\backslash\backslash5}-e_{23\backslash\backslash14\backslash
\backslash5}-e_{25\backslash\backslash34\backslash\backslash1}$\\\hline
$\ytableaushort{{*(green)2}{*(green)1},{*(green)4}3,{*(green)5}}$ &
$e_{21\backslash\backslash43\backslash\backslash5}-e_{12\backslash
\backslash43\backslash\backslash5}-e_{24\backslash\backslash13\backslash
\backslash5}-e_{25\backslash\backslash43\backslash\backslash1}$\\\hline
$\ytableaushort{{*(green)3}{*(green)1},{*(green)4}2,{*(green)5}}$ &
$e_{31\backslash\backslash42\backslash\backslash5}-e_{13\backslash
\backslash42\backslash\backslash5}-e_{34\backslash\backslash12\backslash
\backslash5}-e_{35\backslash\backslash42\backslash\backslash1}$\\\hline
\end{tabular}
\]
(check that these Garnir pentuples are all triangulating!). Thus, if we take
$\mathfrak{G}$ to be the set of these five Garnir relators, then
$\mathfrak{G}$ is triangulating. There are several alternative choices (in
particular, for the last tableau $T$ in our list, we could have picked
$Z=\left\{  1,2,3,4,5\right\}  $, thus obtaining a Garnir relator with $10$ addends).
\end{example}

We now claim that triangulating sets are sufficient to span
$\operatorname*{Ker}\left(  \pi_{D}\right)  $:

\begin{theorem}
\label{thm.specht.TD.garnir-tria-span}Let $D$ be a skew Young diagram with
$\left\vert D\right\vert =n$. Consider the morphism $\pi_{D}:\mathcal{T}%
^{D}\rightarrow\mathcal{S}^{D}$ from Proposition \ref{prop.specht.TD.pi}. Let
$\mathfrak{C}$ be a triangulating set of column relators (in $\mathcal{T}^{D}%
$). Let $\mathfrak{G}$ be a triangulating set of Garnir relators (in
$\mathcal{T}^{D}$). Then,%
\[
\operatorname*{Ker}\left(  \pi_{D}\right)  =\operatorname*{span}%
\nolimits_{\mathbf{k}}\left(  \mathfrak{C}\cup\mathfrak{G}\right)  .
\]

\end{theorem}

This will be proved in Subsection \ref{subsec.specht.TD/Kerpi.pfs}.

As we said, Theorem \ref{thm.specht.TD.garnir-tria-span} and Theorem
\ref{thm.specht.TD.garnir-all-span} can be combined with
(\ref{eq.prop.specht.TD.pi.quot}) to yield explicit presentations of a
skew-shaped Specht module $\mathcal{S}^{D}$ in terms of generators and
relations. Let us make this explicit (for Theorem
\ref{thm.specht.TD.garnir-tria-span}):

\begin{corollary}
\label{cor.specht.TD.SD=TD/(CuG)-tria}Let $D$ be a skew Young diagram with
$\left\vert D\right\vert =n$. Consider the morphism $\pi_{D}:\mathcal{T}%
^{D}\rightarrow\mathcal{S}^{D}$ from Proposition \ref{prop.specht.TD.pi}. Let
$\mathfrak{C}$ be a triangulating set of column relators (in $\mathcal{T}^{D}%
$). Let $\mathfrak{G}$ be a triangulating set of Garnir relators (in
$\mathcal{T}^{D}$). Then, $\operatorname*{span}\nolimits_{\mathbf{k}}\left(
\mathfrak{C}\cup\mathfrak{G}\right)  $ is an $S_{n}$-subrepresentation of
$\mathcal{T}^{D}$, and we have%
\[
\mathcal{S}^{D}\cong\mathcal{T}^{D}/\operatorname*{span}\nolimits_{\mathbf{k}%
}\left(  \mathfrak{C}\cup\mathfrak{G}\right)  \ \ \ \ \ \ \ \ \ \ \text{as
}S_{n}\text{-representations.}%
\]

\end{corollary}

\begin{proof}
The map $\pi_{D}:\mathcal{T}^{D}\rightarrow\mathcal{S}^{D}$ from Proposition
\ref{prop.specht.TD.pi} is a morphism of $S_{n}$-representations, i.e., a
morphism of left $\mathbf{k}\left[  S_{n}\right]  $-modules. Thus, its kernel
$\operatorname*{Ker}\left(  \pi_{D}\right)  $ is a left $\mathbf{k}\left[
S_{n}\right]  $-submodule of $\mathcal{T}^{D}$, that is, an $S_{n}%
$-subrepresentation of $\mathcal{T}^{D}$. But Theorem
\ref{thm.specht.TD.garnir-tria-span} yields%
\[
\operatorname*{Ker}\left(  \pi_{D}\right)  =\operatorname*{span}%
\nolimits_{\mathbf{k}}\left(  \mathfrak{C}\cup\mathfrak{G}\right)  .
\]
Thus, $\operatorname*{span}\nolimits_{\mathbf{k}}\left(  \mathfrak{C}%
\cup\mathfrak{G}\right)  $ is an $S_{n}$-subrepresentation of $\mathcal{T}%
^{D}$ (since $\operatorname*{Ker}\left(  \pi_{D}\right)  $ is an $S_{n}%
$-subrepresentation of $\mathcal{T}^{D}$). Now,
(\ref{eq.prop.specht.TD.pi.quot}) becomes%
\[
\mathcal{S}^{D}\cong\mathcal{T}^{D}/\underbrace{\operatorname*{Ker}\left(
\pi_{D}\right)  }_{=\operatorname*{span}\nolimits_{\mathbf{k}}\left(
\mathfrak{C}\cup\mathfrak{G}\right)  }=\mathcal{T}^{D}/\operatorname*{span}%
\nolimits_{\mathbf{k}}\left(  \mathfrak{C}\cup\mathfrak{G}\right)
\ \ \ \ \ \ \ \ \ \ \text{as }S_{n}\text{-representations.}%
\]
Thus, Corollary \ref{cor.specht.TD.SD=TD/(CuG)-tria} is proved.
\end{proof}

\subsubsection{\label{subsec.specht.TD/Kerpi.pfs}Proofs}

\begin{fineprint}
We shall now provide the missing proofs of Theorem
\ref{thm.specht.TD.garnir-tria-span} and Theorem
\ref{thm.specht.TD.garnir-all-span} that we promised in Subsection
\ref{subsec.specht.TD/Kerpi.thms}.

To prove Theorem \ref{thm.specht.TD.garnir-tria-span}, we shall study the
quotient module $\mathcal{T}^{D}/\operatorname*{span}\nolimits_{\mathbf{k}%
}\left(  \mathfrak{C}\cup\mathfrak{G}\right)  $, showing in particular that it
has a basis $\left(  \overline{e_{T}}\right)  _{T\in\operatorname*{SYT}\left(
D\right)  }$ analogous to the standard basis $\left(  \mathbf{e}_{T}\right)
_{T\in\operatorname*{SYT}\left(  D\right)  }$ of the Specht module
$\mathcal{S}^{D}$. More precisely, we will only need to show that the family
$\left(  \overline{e_{T}}\right)  _{T\in\operatorname*{SYT}\left(  D\right)
}$ spans the quotient $\mathcal{T}^{D}/\operatorname*{span}%
\nolimits_{\mathbf{k}}\left(  \mathfrak{C}\cup\mathfrak{G}\right)  $ (its
linear independence, while true, won't be needed). This will give us an
injective morphism from $\mathcal{T}^{D}/\operatorname*{span}%
\nolimits_{\mathbf{k}}\left(  \mathfrak{C}\cup\mathfrak{G}\right)  $ to
$\mathcal{S}^{D}$ compatible with the projection $\pi_{D}:\mathcal{T}%
^{D}\rightarrow\mathcal{S}^{D}$, and this will easily yield Theorem
\ref{thm.specht.TD.garnir-tria-span}. This proof will be somewhat long, but
almost all of it will be a calque of our above proof of the standard basis
theorem (Theorem \ref{thm.spechtmod.basis}), and the rest is easy enough that
it could be given as a \fbox{2}-point exercise. Yet, the quest for
completeness is forcing me to present this proof in full, for what it is worth.

\begin{convention}
Throughout Subsection \ref{subsec.specht.TD/Kerpi.pfs}, we shall use the
following notations:

\begin{itemize}
\item We fix a skew Young diagram $D$ with $\left\vert D\right\vert =n$.

\item We let $\mathfrak{C}$ be a triangulating set of column relators (in
$\mathcal{T}^{D}$), and we let $\mathfrak{G}$ be a triangulating set of Garnir
relators (in $\mathcal{T}^{D}$).

\item We let $\operatorname*{SYT}\left(  D\right)  $ denote the set of all
standard tableaux of shape $D$. In other words, $\operatorname*{SYT}\left(
D\right)  $ is the set of all standard $n$-tableaux of shape $D$ (by
Proposition \ref{prop.tableau.std-n}). Thus, $\operatorname*{SYT}\left(
D\right)  \subseteq\operatorname*{Tab}\left(  D\right)  $.

\item We consider the quotient $\mathbf{k}$-module $\mathcal{T}^{D}%
/\operatorname*{span}\nolimits_{\mathbf{k}}\left(  \mathfrak{C}\cup
\mathfrak{G}\right)  $. For each $v\in\mathcal{T}^{D}$, we let $\overline{v}$
be the projection of $v$ on this quotient module. In particular, the
projections $\overline{e_{T}}$ of the standard basis vectors $e_{T}%
\in\mathcal{T}^{D}$ will be called the \emph{quasi-polytabloids} (to stress
their analogy with the polytabloids $\mathbf{e}_{T}\in\mathcal{S}^{D}$). The
quasi-polytabloids $\overline{e_{T}}$ with $T\in\operatorname*{SYT}\left(
D\right)  $ will be called the \emph{standard quasi-polytabloids}.
\end{itemize}
\end{convention}

The first step on our way to Theorem \ref{thm.specht.TD.garnir-tria-span} is a
restatement of Proposition \ref{prop.coltabloid.col-st}:

\begin{lemma}
\label{lem.coltabloid.col-st-unique}Each $n$-column-tabloid of shape $D$
contains a unique column-standard $n$-tableau\footnotemark.
\end{lemma}

\footnotetext{The word \textquotedblleft contains\textquotedblright\ here
should be understood literally, as a containment of an element in a set (since
an $n$-column-tabloid is an equivalence class of $n$-tableaux of shape $D$
with respect to column-equivalence, and thus is a set of $n$-tableaux).}

\begin{proof}
Proposition \ref{prop.coltabloid.col-st} shows that there is a bijection%
\begin{align*}
&  \text{from }\left\{  \text{column-standard }n\text{-tableaux of shape
}D\right\} \\
&  \text{to }\left\{  n\text{-column-tabloids of shape }D\right\}
\end{align*}
that sends each column-standard $n$-tableau $T$ to its column-tabloid
$\widetilde{T}$. Hence, each $n$-column-tabloid of shape $D$ has a unique
preimage under this bijection (because this is what being a bijection means).
In other words, for each $n$-column-tabloid $C$ of shape $D$, there exists a
unique column-standard $n$-tableau $T$ of shape $D$ such that $\widetilde{T}%
=C$. In other words, for each $n$-column-tabloid $C$ of shape $D$, there
exists a unique column-standard $n$-tableau $T$ of shape $D$ such that $T\in
C$ (because the statement \textquotedblleft$\widetilde{T}=C$\textquotedblright%
\ is equivalent to \textquotedblleft$T\in C$\textquotedblright). In other
words, each $n$-column-tabloid $C$ of shape $D$ contains a unique
column-standard $n$-tableau $T$. This proves Lemma
\ref{lem.coltabloid.col-st-unique}.
\end{proof}

Next, we show a lemma saying that the quasi-polytabloids $\overline{e_{T}}%
\in\mathcal{T}^{D}/\operatorname*{span}\nolimits_{\mathbf{k}}\left(
\mathfrak{C}\cup\mathfrak{G}\right)  $ corresponding to column-equivalent
$n$-tableaux $T$ are equal up to sign:

\begin{lemma}
\label{lem.specht.TD.garnir-tria-span.coleq}Let $P$ and $Q$ be two
column-equivalent $n$-tableaux of shape $D$. Then, $\overline{e_{P}}%
=\pm\overline{e_{Q}}$ in $\mathcal{T}^{D}/\operatorname*{span}%
\nolimits_{\mathbf{k}}\left(  \mathfrak{C}\cup\mathfrak{G}\right)  $.
\end{lemma}

\begin{proof}
Recall that an $n$-column-tabloid is an equivalence class of $n$-tableaux of
shape $D$ with respect to column-equivalence. Hence, the $n$-tableaux $P$ and
$Q$ belong to the same $n$-column-tabloid (since they are column-equivalent).
In other words, $\widetilde{P}=\widetilde{Q}$.

However, Lemma \ref{lem.coltabloid.col-st-unique} says that each
$n$-column-tabloid of shape $D$ contains a unique column-standard $n$-tableau.
Thus, in particular, the $n$-column-tabloid $\widetilde{P}=\widetilde{Q}$
contains a unique column-standard $n$-tableau. Let us denote this unique
column-standard $n$-tableau by $K$. Hence, $K\in\widetilde{P}=\widetilde{Q}$,
so that $\widetilde{K}=\widetilde{Q}=\widetilde{P}$. Furthermore, $K$ is
defined as the \textbf{unique} column-standard $n$-tableau contained in the
$n$-column-tabloid $\widetilde{P}$. Hence, if $W\in\widetilde{P}$ is any
column-standard $n$-tableau, then%
\begin{equation}
W=K. \label{pf.lem.specht.TD.garnir-tria-span.coleq.W=T}%
\end{equation}

We will now prove that%
\begin{equation}
\overline{e_{P}}=\pm\overline{e_{K}}.
\label{pf.lem.specht.TD.garnir-tria-span.coleq.PT}%
\end{equation}

\begin{proof}
[Proof of (\ref{pf.lem.specht.TD.garnir-tria-span.coleq.PT}):]We must show
that $\overline{e_{P}}=\pm\overline{e_{K}}$.

The $n$-tableau $P$ is clearly contained in the $n$-column-tabloid
$\widetilde{P}$. If this $n$-tableau $P$ is furthermore column-standard, then
(\ref{pf.lem.specht.TD.garnir-tria-span.coleq.W=T}) (applied to $W=P$) shows
that $P=K$ and therefore $\overline{e_{P}}=\overline{e_{K}}=\pm\overline
{e_{K}}$, so that (\ref{pf.lem.specht.TD.garnir-tria-span.coleq.PT}) is
proved. Hence, for the rest of this proof, we WLOG assume that $P$ is
\textbf{not} column-standard.

Now recall that the set $\mathfrak{C}$ of column relators is triangulating.
Hence, (by Definition \ref{def.specht.TD.garnir-triasub} \textbf{(a)}) we
conclude that for any $n$-tableau $T\in\operatorname*{Tab}\left(  D\right)  $
that is not column-standard, there exists a $u\in\mathcal{C}\left(  T\right)
$ such that the $n$-tableau $uT$ is column-standard and such that
$e_{T}-\left(  -1\right)  ^{u}e_{uT}\in\mathfrak{C}$. Applying this to $T=P$,
we conclude that there exists a $u\in\mathcal{C}\left(  P\right)  $ such that
the $n$-tableau $uP$ is column-standard and such that $e_{P}-\left(
-1\right)  ^{u}e_{uP}\in\mathfrak{C}$ (since $P\in\operatorname*{Tab}\left(
D\right)  $ is not column-standard). Consider this $u$. Note that
\[
e_{P}-\left(  -1\right)  ^{u}e_{uP}\in\mathfrak{C}\subseteq\mathfrak{C}%
\cup\mathfrak{G}\subseteq\operatorname*{span}\nolimits_{\mathbf{k}}\left(
\mathfrak{C}\cup\mathfrak{G}\right)  .
\]
In other words, $\overline{e_{P}-\left(  -1\right)  ^{u}e_{uP}}=0$ in
$\mathcal{T}^{D}/\operatorname*{span}\nolimits_{\mathbf{k}}\left(
\mathfrak{C}\cup\mathfrak{G}\right)  $. Hence, $\overline{e_{P}}-\left(
-1\right)  ^{u}\overline{e_{uP}}=\overline{e_{P}-\left(  -1\right)  ^{u}%
e_{uP}}=0$, so that%
\begin{equation}
\overline{e_{P}}=\underbrace{\left(  -1\right)  ^{u}}_{=\pm1}\overline{e_{uP}%
}=\pm\overline{e_{uP}}. \label{pf.lem.specht.TD.garnir-tria-span.coleq.5}%
\end{equation}

But $u\in\mathcal{C}\left(  P\right)  $. Hence, the $n$-tableau
$uP=u\rightharpoonup P$ is column-equivalent to $P$ (by Proposition
\ref{prop.tableau.Sn-act.0} \textbf{(b)}, applied to $w=u$ and $T=P$), and
thus belongs to the same $n$-column-tabloid as $P$. In other words,
$uP\in\widetilde{P}$. Since $uP$ is furthermore column-standard, we thus
conclude that $uP=K$ (by (\ref{pf.lem.specht.TD.garnir-tria-span.coleq.W=T}),
applied to $W=uP$). In view of this, we can rewrite
(\ref{pf.lem.specht.TD.garnir-tria-span.coleq.5}) as $\overline{e_{P}}%
=\pm\overline{e_{K}}$. Thus, (\ref{pf.lem.specht.TD.garnir-tria-span.coleq.PT}%
) is proved.
\end{proof}

We have now shown that $\overline{e_{P}}=\pm\overline{e_{K}}$. Since $P$ and
$Q$ play analogous roles in our situation, we can likewise show that
$\overline{e_{Q}}=\pm\overline{e_{K}}$ (the $\pm$ sign here is not necessarily
the same as in $\overline{e_{P}}=\pm\overline{e_{K}}$). Hence, $\overline
{e_{K}}=\pm\overline{e_{Q}}$. Combining our results, we obtain%
\[
\overline{e_{P}}=\pm\underbrace{\overline{e_{K}}}_{=\pm\overline{e_{Q}}}%
=\pm\left(  \pm\overline{e_{Q}}\right)  =\pm\overline{e_{Q}}.
\]
Thus, Lemma \ref{lem.specht.TD.garnir-tria-span.coleq} is proved.
\end{proof}

Next we show an analogue of Lemma \ref{lem.spechtmod.straighten1}:

\begin{lemma}
\label{lem.specht.TD.garnir-tria-span.straighten1}Let $T$ be an $n$-tableau of
shape $D$ that is column-standard but not standard. Then, $\overline{e_{T}}$
can be written as a $\mathbb{Z}$-linear combination of quasi-polytabloids
$\overline{e_{P}}$, where the subscripts $P$ are $n$-tableaux of shape $D$
satisfying $\widetilde{P}>\widetilde{T}$ (with respect to the column last
letter order).
\end{lemma}

\begin{proof}
The set $\mathfrak{G}$ of Garnir relators is triangulating. Therefore, by
Definition \ref{def.specht.TD.garnir-triasub} \textbf{(c)}, we conclude that
there exists a triangulating Garnir pentuple $\left(  T,j,k,Z,L\right)  $ with
first entry $T$ such that $\sum_{x\in L}\left(  -1\right)  ^{x}e_{xT}%
\in\mathfrak{G}$ (since $T\in\operatorname*{Tab}\left(  D\right)  $ is
column-standard but not standard). Consider this Garnir pentuple $\left(
T,j,k,Z,L\right)  $. By the definition of a Garnir pentuple, it thus satisfies
the following:

\begin{itemize}
\item The entries $j$ and $k$ are two distinct integers.

\item The entry $Z$ is a subset of $\operatorname*{Col}\left(  j,T\right)
\cup\operatorname*{Col}\left(  k,T\right)  $ such that $\left\vert
Z\right\vert >\left\vert D\left\langle j\right\rangle \cup D\left\langle
k\right\rangle \right\vert $.

\item The entry $L$ is a left transversal of $S_{n,Z}\cap\mathcal{C}\left(
T\right)  $ in $S_{n,Z}$.
\end{itemize}

Moreover, the Garnir pentuple $\left(  T,j,k,Z,L\right)  $ is triangulating.
In other words, each $x\in L\setminus\mathcal{C}\left(  T\right)  $ satisfies
\begin{equation}
\widetilde{xT}>\widetilde{T}\ \ \ \ \ \ \ \ \ \ \text{in the column last
letter order} \label{pf.lem.specht.TD.garnir-tria-span.straighten1.xT>T}%
\end{equation}
(by Definition \ref{def.specht.TD.garnir-triasub} \textbf{(b)}).

However,
\[
\sum_{x\in L}\left(  -1\right)  ^{x}e_{xT}\in\mathfrak{G}\subseteq
\mathfrak{C}\cup\mathfrak{G}\subseteq\operatorname*{span}\nolimits_{\mathbf{k}%
}\left(  \mathfrak{C}\cup\mathfrak{G}\right)  .
\]
In other words, $\overline{\sum_{x\in L}\left(  -1\right)  ^{x}e_{xT}}=0$ in
$\mathcal{T}^{D}/\operatorname*{span}\nolimits_{\mathbf{k}}\left(
\mathfrak{C}\cup\mathfrak{G}\right)  $. Hence,
\begin{equation}
0=\overline{\sum_{x\in L}\left(  -1\right)  ^{x}e_{xT}}=\sum_{x\in L}\left(
-1\right)  ^{x}\overline{e_{xT}}.
\label{pf.lem.specht.TD.garnir-tria-span.straighten1.4}%
\end{equation}

Now, recall that $L$ is a left transversal of $S_{n,Z}\cap\mathcal{C}\left(
T\right)  $ in $S_{n,Z}$. Thus, each left coset of the subgroup $S_{n,Z}%
\cap\mathcal{C}\left(  T\right)  $ in $S_{n,Z}$ contains exactly one element
of $L$ (by Definition \ref{def.G/H-transversal} \textbf{(a)}). Hence, in
particular, the trivial left coset $1\left(  S_{n,Z}\cap\mathcal{C}\left(
T\right)  \right)  $ of $S_{n,Z}\cap\mathcal{C}\left(  T\right)  $ in
$S_{n,Z}$ contains exactly one element of $L$. Let $v\in L$ be this element.
Thus, $v\in1\left(  S_{n,Z}\cap\mathcal{C}\left(  T\right)  \right)
=S_{n,Z}\cap\mathcal{C}\left(  T\right)  \subseteq\mathcal{C}\left(  T\right)
$. Therefore, the $n$-tableau $vT=v\rightharpoonup T$ is column-equivalent to
$T$ (by Proposition \ref{prop.tableau.Sn-act.0} \textbf{(b)}, applied to
$w=v$). Thus, Lemma \ref{lem.specht.TD.garnir-tria-span.coleq} (applied to
$P=vT$ and $Q=T$) yields%
\[
\overline{e_{vT}}=\pm\overline{e_{T}}\ \ \ \ \ \ \ \ \ \ \text{in }%
\mathcal{T}^{D}/\operatorname*{span}\nolimits_{\mathbf{k}}\left(
\mathfrak{C}\cup\mathfrak{G}\right)  .
\]

Now, recall that $v$ is the only element of $L$ that is contained in $1\left(
S_{n,Z}\cap\mathcal{C}\left(  T\right)  \right)  $ (by the definition of $v$).
Thus, each $x\in L\setminus\left\{  v\right\}  $ satisfies $x\notin1\left(
S_{n,Z}\cap\mathcal{C}\left(  T\right)  \right)  $ and therefore
\begin{equation}
\widetilde{xT}>\widetilde{T}\ \ \ \ \ \ \ \ \ \ \text{in the column last
letter order} \label{pf.lem.specht.TD.garnir-tria-span.straighten1.xT>T2}%
\end{equation}
(because combining $x\in L\setminus\left\{  v\right\}  \subseteq L\subseteq
S_{n,Z}$ with $x\notin1\left(  S_{n,Z}\cap\mathcal{C}\left(  T\right)
\right)  =S_{n,Z}\cap\mathcal{C}\left(  T\right)  $, we obtain $x\in
S_{n,Z}\setminus\left(  S_{n,Z}\cap\mathcal{C}\left(  T\right)  \right)
=S_{n,Z}\setminus\mathcal{C}\left(  T\right)  $ and therefore $x\notin%
\mathcal{C}\left(  T\right)  $, so that $x\in L\setminus\mathcal{C}\left(
T\right)  $ (since $x\in L$ and $x\notin\mathcal{C}\left(  T\right)  $), and
therefore (\ref{pf.lem.specht.TD.garnir-tria-span.straighten1.xT>T}) shows
that $\widetilde{xT}>\widetilde{T}$ in the column last letter order). But
(\ref{pf.lem.specht.TD.garnir-tria-span.straighten1.4}) becomes%
\begin{align*}
0  &  =\sum_{x\in L}\left(  -1\right)  ^{x}\overline{e_{xT}}%
=\underbrace{\left(  -1\right)  ^{v}}_{=\pm1}\underbrace{\overline{e_{vT}}%
}_{=\pm\overline{e_{T}}}+\sum_{x\in L\setminus\left\{  v\right\}  }\left(
-1\right)  ^{x}\overline{e_{xT}}\\
&  \ \ \ \ \ \ \ \ \ \ \ \ \ \ \ \ \ \ \ \ \left(
\begin{array}
[c]{c}%
\text{here, we have split off the addend for }x=v\\
\text{from the sum, since }v\in L
\end{array}
\right) \\
&  =\underbrace{\left(  \pm1\right)  \left(  \pm\overline{e_{T}}\right)
}_{=\pm\overline{e_{T}}}+\sum_{x\in L\setminus\left\{  v\right\}  }\left(
-1\right)  ^{x}\overline{e_{xT}}=\pm\overline{e_{T}}+\sum_{x\in L\setminus
\left\{  v\right\}  }\left(  -1\right)  ^{x}\overline{e_{xT}}.
\end{align*}
Moving the $\pm\overline{e_{T}}$ addend onto the left hand side here, we
obtain
\begin{equation}
\overline{e_{T}}=\pm\sum_{x\in L\setminus\left\{  v\right\}  }\left(
-1\right)  ^{x}\overline{e_{xT}}.
\label{pf.lem.specht.TD.garnir-tria-span.straighten1.at}%
\end{equation}

But (\ref{pf.lem.specht.TD.garnir-tria-span.straighten1.xT>T2}) tells us that
each $x\in L\setminus\left\{  v\right\}  $ satisfies $\widetilde{xT}%
>\widetilde{T}$ in the column last letter order. Hence, each of the vectors
$\overline{e_{xT}}$ on the right hand side of
(\ref{pf.lem.specht.TD.garnir-tria-span.straighten1.at}) has the form
$\overline{e_{P}}$ for some $n$-tableau $P$ of shape $D$ satisfying
$\widetilde{P}>\widetilde{T}$ in the column last letter order (namely, for
$P=xT$). Thus, the equality
(\ref{pf.lem.specht.TD.garnir-tria-span.straighten1.at}) shows that
$\overline{e_{T}}$ can be written as a $\mathbb{Z}$-linear combination of
quasi-polytabloids $\overline{e_{P}}$, where the subscripts $P$ are
$n$-tableaux of shape $D$ satisfying $\widetilde{P}>\widetilde{T}$ (with
respect to the column last letter order). This proves Lemma
\ref{lem.specht.TD.garnir-tria-span.straighten1}.
\end{proof}

The following lemma is an analogue of Proposition
\ref{prop.spechtmod.straighten}:

\begin{lemma}
\label{lem.specht.TD.garnir-tria-span.straighten2}Let $T$ be an $n$-tableau of
shape $D$. Then, $\overline{e_{T}}$ is a $\mathbb{Z}$-linear combination of
standard quasi-polytabloids $\overline{e_{Q}}$, where $Q$ ranges over all
standard tableaux of shape $D$.
\end{lemma}

\begin{proof}
Forget that we fixed $T$. We must show that each quasi-polytabloid
$\overline{e_{T}}$ is a $\mathbb{Z}$-linear combination of standard quasi-polytabloids.

In this proof, we will not consider any diagrams other than $D$. Thus, the
word \textquotedblleft$n$-tableau\textquotedblright\ will always mean
\textquotedblleft$n$-tableau of shape $D$\textquotedblright. Likewise for the
word \textquotedblleft$n$-column-tabloid\textquotedblright.

Define a $\mathbb{Z}$-submodule $\overline{\mathcal{Q}}^{D}$ of $\mathcal{T}%
^{D}/\operatorname*{span}\nolimits_{\mathbf{k}}\left(  \mathfrak{C}%
\cup\mathfrak{G}\right)  $ by
\[
\overline{\mathcal{Q}}^{D}:=\operatorname*{span}\nolimits_{\mathbb{Z}}\left\{
\overline{e_{Q}}\ \mid\ Q\text{ is a standard }n\text{-tableau}\right\}  .
\]
Thus, $\overline{\mathcal{Q}}^{D}$ is the $\mathbb{Z}$-linear span of all
standard quasi-polytabloids. We shall show that
\begin{equation}
\overline{e_{T}}\in\overline{\mathcal{Q}}^{D}\ \ \ \ \ \ \ \ \ \ \text{for all
}n\text{-tableaux }T.
\label{pf.lem.specht.TD.garnir-tria-span.straighten2.goal}%
\end{equation}

Indeed, let us define the \emph{depth} $\operatorname*{dep}\widetilde{T}$ of
an $n$-column-tabloid $\widetilde{T}$ to be the \# of all $n$-column-tabloids
$\widetilde{P}$ that satisfy $\widetilde{P}>\widetilde{T}$ (with respect to
the column last letter order). This is clearly a well-defined nonnegative
integer (since there are only finitely many $n$-column-tabloids $\widetilde{P}%
$). Moreover, if two $n$-column-tabloids $\widetilde{S}$ and $\widetilde{T}$
satisfy $\widetilde{S}>\widetilde{T}$ (with respect to the column last letter
order), then%
\begin{equation}
\operatorname*{dep}\widetilde{S}<\operatorname*{dep}\widetilde{T}
\label{pf.lem.specht.TD.garnir-tria-span.straighten2.depless}%
\end{equation}
(indeed, this was already proved during the proof of Proposition
\ref{prop.spechtmod.straighten}).

Now, we shall prove (\ref{pf.lem.specht.TD.garnir-tria-span.straighten2.goal})
by strong induction on $\operatorname*{dep}\widetilde{T}$:

\textit{Induction step:} Fix an $n$-tableau $U$, and assume (as the induction
hypothesis) that (\ref{pf.lem.specht.TD.garnir-tria-span.straighten2.goal}) is
already proved for all $n$-tableaux $T$ satisfying $\operatorname*{dep}%
\widetilde{T}<\operatorname*{dep}\widetilde{U}$. Our goal is now to prove
(\ref{pf.lem.specht.TD.garnir-tria-span.straighten2.goal}) for $T=U$. In other
words, our goal is to prove that $\overline{e_{U}}\in\overline{\mathcal{Q}%
}^{D}$.

The tableau $U$ may or may not be column-standard. However, it is easy to find
a column-standard $n$-tableau $S$ that is column-equivalent to $U$ (just let
$S$ be the result of sorting the entries in each column of $U$ into increasing
order). Consider this $S$. Then, $\widetilde{S}=\widetilde{U}$ (since $S$ is
column-equivalent to $U$) and $\overline{e_{U}}=\pm\overline{e_{S}}$ (by Lemma
\ref{lem.specht.TD.garnir-tria-span.coleq}, applied to $P=U$ and $Q=S$).

If the $n$-tableau $S$ is standard, then we have $\overline{e_{S}}\in
\overline{\mathcal{Q}}^{D}$ (by the definition of $\overline{\mathcal{Q}}^{D}%
$, since $\overline{e_{S}}$ is a standard quasi-polytabloid), and thus
$\overline{e_{U}}=\pm\underbrace{\overline{e_{S}}}_{\in\overline{\mathcal{Q}%
}^{D}}\in\overline{\mathcal{Q}}^{D}$ (since $\overline{\mathcal{Q}}^{D}$ is a
$\mathbb{Z}$-module), which immediately completes our induction step (since
our goal is to prove that $\overline{e_{U}}\in\overline{\mathcal{Q}}^{D}$).

Hence, for the rest of the induction step, we WLOG assume that $S$ is not standard.

Thus, Lemma \ref{lem.specht.TD.garnir-tria-span.straighten1} (applied to $S$
instead of $T$) shows that $\overline{e_{S}}$ can be written as a $\mathbb{Z}%
$-linear combination of quasi-polytabloids $\overline{e_{P}}$, where the
subscripts $P$ are $n$-tableaux of shape $D$ satisfying $\widetilde{P}%
>\widetilde{S}$ (with respect to the column last letter order). In other
words,
\begin{equation}
\overline{e_{S}}\in\operatorname*{span}\nolimits_{\mathbb{Z}}\left\{
\overline{e_{P}}\ \mid\ P\text{ is an }n\text{-tableau satisfying
}\widetilde{P}>\widetilde{S}\right\}  .
\label{pf.lem.specht.TD.garnir-tria-span.straighten2.5}%
\end{equation}

Now, let $P$ be an $n$-tableau satisfying $\widetilde{P}>\widetilde{S}$. Then,
$\widetilde{P}>\widetilde{S}=\widetilde{U}$ and therefore $\operatorname*{dep}%
\widetilde{P}<\operatorname*{dep}\widetilde{U}$ (by
(\ref{pf.lem.specht.TD.garnir-tria-span.straighten2.depless}), applied to $P$
and $U$ instead of $S$ and $T$). But our induction hypothesis says that
(\ref{pf.lem.specht.TD.garnir-tria-span.straighten2.goal}) is already proved
for all $n$-tableaux $T$ satisfying $\operatorname*{dep}\widetilde{T}%
<\operatorname*{dep}\widetilde{U}$. Applying this to $T=P$, we conclude that
(\ref{pf.lem.specht.TD.garnir-tria-span.straighten2.goal}) is already proved
for $T=P$ (since $\operatorname*{dep}\widetilde{P}<\operatorname*{dep}%
\widetilde{U}$). Hence, $\overline{e_{P}}\in\overline{\mathcal{Q}}^{D}$.

Forget that we fixed $P$. We thus have shown that $\overline{e_{P}}%
\in\overline{\mathcal{Q}}^{D}$ whenever $P$ is an $n$-tableau satisfying
$\widetilde{P}>\widetilde{S}$. In other words,%
\[
\left\{  \overline{e_{P}}\ \mid\ P\text{ is an }n\text{-tableau satisfying
}\widetilde{P}>\widetilde{S}\right\}  \subseteq\overline{\mathcal{Q}}^{D}.
\]

Now, (\ref{pf.lem.specht.TD.garnir-tria-span.straighten2.5}) becomes%
\[
\overline{e_{S}}\in\operatorname*{span}\nolimits_{\mathbb{Z}}%
\underbrace{\left\{  \overline{e_{P}}\ \mid\ P\text{ is an }n\text{-tableau
satisfying }\widetilde{P}>\widetilde{S}\right\}  }_{\subseteq\overline
{\mathcal{Q}}^{D}}\subseteq\operatorname*{span}\nolimits_{\mathbb{Z}}\left(
\overline{\mathcal{Q}}^{D}\right)  =\overline{\mathcal{Q}}^{D}%
\]
(since $\overline{\mathcal{Q}}^{D}$ is a $\mathbb{Z}$-module). Hence,
$\overline{e_{U}}=\pm\underbrace{\overline{e_{S}}}_{\in\overline{\mathcal{Q}%
}^{D}}\in\overline{\mathcal{Q}}^{D}$ (since $\overline{\mathcal{Q}}^{D}$ is a
$\mathbb{Z}$-module). But this is precisely what we intended to prove. Thus,
we have proved (\ref{pf.lem.specht.TD.garnir-tria-span.straighten2.goal}) for
$T=U$. This completes the induction step. Thus,
(\ref{pf.lem.specht.TD.garnir-tria-span.straighten2.goal}) is proved for all
$T$.

In other words, we have proved that each quasi-polytabloid $\overline{e_{T}}$
belongs to $\overline{\mathcal{Q}}^{D}$. In other words, each
quasi-polytabloid $\overline{e_{T}}$ belongs to the $\mathbb{Z}$-linear span
of all standard quasi-polytabloids (since $\overline{\mathcal{Q}}^{D}$ is the
$\mathbb{Z}$-linear span of all standard quasi-polytabloids). In other words,
each quasi-polytabloid $\overline{e_{T}}$ is a $\mathbb{Z}$-linear combination
of standard quasi-polytabloids. This proves Lemma
\ref{lem.specht.TD.garnir-tria-span.straighten2}.
\end{proof}

Lemma \ref{lem.specht.TD.garnir-tria-span.straighten2} easily yields the following:

\begin{lemma}
\label{lem.specht.TD.garnir-tria-span.straighten3}The $\mathbf{k}$-module
$\mathcal{T}^{D}/\operatorname*{span}\nolimits_{\mathbf{k}}\left(
\mathfrak{C}\cup\mathfrak{G}\right)  $ is spanned by the family $\left(
\overline{e_{T}}\right)  _{T\in\operatorname*{SYT}\left(  D\right)  }$ of all
standard quasi-polytabloids.
\end{lemma}

\begin{proof}
The $\mathbf{k}$-module $\mathcal{T}^{D}$ is spanned by the family $\left(
e_{T}\right)  _{T\in\operatorname*{Tab}\left(  D\right)  }$ (since it is free
with basis $\left(  e_{T}\right)  _{T\in\operatorname*{Tab}\left(  D\right)
}$). Hence, its quotient module $\mathcal{T}^{D}/\operatorname*{span}%
\nolimits_{\mathbf{k}}\left(  \mathfrak{C}\cup\mathfrak{G}\right)  $ is
spanned by the family $\left(  \overline{e_{T}}\right)  _{T\in
\operatorname*{Tab}\left(  D\right)  }$. In other words,%
\begin{equation}
\mathcal{T}^{D}/\operatorname*{span}\nolimits_{\mathbf{k}}\left(
\mathfrak{C}\cup\mathfrak{G}\right)  =\operatorname*{span}%
\nolimits_{\mathbf{k}}\left\{  \overline{e_{T}}\ \mid\ T\in\operatorname*{Tab}%
\left(  D\right)  \right\}  .
\label{pf.lem.specht.TD.garnir-tria-span.straighten3.1}%
\end{equation}

On the other hand, each $T\in\operatorname*{Tab}\left(  D\right)  $ satisfies%
\begin{align*}
\overline{e_{T}}  &  \in\operatorname*{span}\nolimits_{\mathbb{Z}}\left\{
\overline{e_{Q}}\ \mid\ Q\text{ is a standard tableau of shape }D\right\} \\
&  \ \ \ \ \ \ \ \ \ \ \ \ \ \ \ \ \ \ \ \ \left(
\begin{array}
[c]{c}%
\text{since Lemma \ref{lem.specht.TD.garnir-tria-span.straighten2} shows that
}\overline{e_{T}}\text{ is a }\mathbb{Z}\text{-linear combination of}\\
\text{the }\overline{e_{Q}}\text{, where }Q\text{ ranges over all standard
tableaux of shape }D
\end{array}
\right) \\
&  =\operatorname*{span}\nolimits_{\mathbb{Z}}\left\{  \overline{e_{Q}}%
\ \mid\ Q\in\operatorname*{SYT}\left(  D\right)  \right\} \\
&  \ \ \ \ \ \ \ \ \ \ \ \ \ \ \ \ \ \ \ \ \left(  \text{since }%
\operatorname*{SYT}\left(  D\right)  \text{ is the set of all standard
tableaux of shape }D\right) \\
&  \subseteq\operatorname*{span}\nolimits_{\mathbf{k}}\left\{  \overline
{e_{Q}}\ \mid\ Q\in\operatorname*{SYT}\left(  D\right)  \right\} \\
&  \ \ \ \ \ \ \ \ \ \ \ \ \ \ \ \ \ \ \ \ \left(  \text{since }%
\operatorname*{span}\nolimits_{\mathbb{Z}}A\subseteq\operatorname*{span}%
\nolimits_{\mathbf{k}}A\text{ for any subset }A\text{ of any }\mathbf{k}%
\text{-module}\right)  .
\end{align*}
In other words,%
\[
\left\{  \overline{e_{T}}\ \mid\ T\in\operatorname*{Tab}\left(  D\right)
\right\}  \subseteq\operatorname*{span}\nolimits_{\mathbf{k}}\left\{
\overline{e_{Q}}\ \mid\ Q\in\operatorname*{SYT}\left(  D\right)  \right\}  .
\]
Thus,%
\begin{align*}
&  \operatorname*{span}\nolimits_{\mathbf{k}}\underbrace{\left\{
\overline{e_{T}}\ \mid\ T\in\operatorname*{Tab}\left(  D\right)  \right\}
}_{\subseteq\operatorname*{span}\nolimits_{\mathbf{k}}\left\{  \overline
{e_{Q}}\ \mid\ Q\in\operatorname*{SYT}\left(  D\right)  \right\}  }\\
&  \subseteq\operatorname*{span}\nolimits_{\mathbf{k}}\left(
\operatorname*{span}\nolimits_{\mathbf{k}}\left\{  \overline{e_{Q}}%
\ \mid\ Q\in\operatorname*{SYT}\left(  D\right)  \right\}  \right) \\
&  \subseteq\operatorname*{span}\nolimits_{\mathbf{k}}\left\{  \overline
{e_{Q}}\ \mid\ Q\in\operatorname*{SYT}\left(  D\right)  \right\} \\
&  \ \ \ \ \ \ \ \ \ \ \ \ \ \ \ \ \ \ \ \ \left(  \text{since }%
\operatorname*{span}\nolimits_{\mathbf{k}}\left\{  \overline{e_{Q}}%
\ \mid\ Q\in\operatorname*{SYT}\left(  D\right)  \right\}  \text{ is a
}\mathbf{k}\text{-module}\right) \\
&  =\operatorname*{span}\nolimits_{\mathbf{k}}\left\{  \overline{e_{T}}%
\ \mid\ T\in\operatorname*{SYT}\left(  D\right)  \right\}  .
\end{align*}
Using (\ref{pf.lem.specht.TD.garnir-tria-span.straighten3.1}), we can rewrite
this as%
\[
\mathcal{T}^{D}/\operatorname*{span}\nolimits_{\mathbf{k}}\left(
\mathfrak{C}\cup\mathfrak{G}\right)  \subseteq\operatorname*{span}%
\nolimits_{\mathbf{k}}\left\{  \overline{e_{T}}\ \mid\ T\in\operatorname*{SYT}%
\left(  D\right)  \right\}  .
\]
Combining this with the obvious inclusion%
\[
\operatorname*{span}\nolimits_{\mathbf{k}}\left\{  \overline{e_{T}}%
\ \mid\ T\in\operatorname*{SYT}\left(  D\right)  \right\}  \subseteq
\mathcal{T}^{D}/\operatorname*{span}\nolimits_{\mathbf{k}}\left(
\mathfrak{C}\cup\mathfrak{G}\right)  ,
\]
we conclude that%
\[
\mathcal{T}^{D}/\operatorname*{span}\nolimits_{\mathbf{k}}\left(
\mathfrak{C}\cup\mathfrak{G}\right)  =\operatorname*{span}%
\nolimits_{\mathbf{k}}\left\{  \overline{e_{T}}\ \mid\ T\in\operatorname*{SYT}%
\left(  D\right)  \right\}  .
\]
In other words, the $\mathbf{k}$-module $\mathcal{T}^{D}/\operatorname*{span}%
\nolimits_{\mathbf{k}}\left(  \mathfrak{C}\cup\mathfrak{G}\right)  $ is
spanned by the family $\left(  \overline{e_{T}}\right)  _{T\in
\operatorname*{SYT}\left(  D\right)  }$. This proves Lemma
\ref{lem.specht.TD.garnir-tria-span.straighten3}.
\end{proof}

We can now prove Theorem \ref{thm.specht.TD.garnir-tria-span} at last:

\begin{proof}
[Proof of Theorem \ref{thm.specht.TD.garnir-tria-span}.]Consider the morphism
$\pi_{D}:\mathcal{T}^{D}\rightarrow\mathcal{S}^{D}$ from Proposition
\ref{prop.specht.TD.pi}.

Each element of $\mathfrak{C}$ is a column relator (since $\mathfrak{C}$ is a
set of column relators) and thus belongs to $\operatorname*{Ker}\left(
\pi_{D}\right)  $ (by Proposition \ref{prop.specht.TD.relators-in-ker}
\textbf{(a)}). In other words, $\mathfrak{C}\subseteq\operatorname*{Ker}%
\left(  \pi_{D}\right)  $.

Each element of $\mathfrak{G}$ is a Garnir relator (since $\mathfrak{G}$ is a
set of Garnir relators) and thus belongs to $\operatorname*{Ker}\left(
\pi_{D}\right)  $ (by Proposition \ref{prop.specht.TD.relators-in-ker}
\textbf{(b)}). In other words, $\mathfrak{G}\subseteq\operatorname*{Ker}%
\left(  \pi_{D}\right)  $.

Combining $\mathfrak{C}\subseteq\operatorname*{Ker}\left(  \pi_{D}\right)  $
and $\mathfrak{G}\subseteq\operatorname*{Ker}\left(  \pi_{D}\right)  $, we
obtain $\mathfrak{C}\cup\mathfrak{G}\subseteq\operatorname*{Ker}\left(
\pi_{D}\right)  $. Hence,
\[
\operatorname*{span}\nolimits_{\mathbf{k}}\underbrace{\left(  \mathfrak{C}%
\cup\mathfrak{G}\right)  }_{\subseteq\operatorname*{Ker}\left(  \pi
_{D}\right)  }\subseteq\operatorname*{span}\nolimits_{\mathbf{k}}\left(
\operatorname*{Ker}\left(  \pi_{D}\right)  \right)  \subseteq
\operatorname*{Ker}\left(  \pi_{D}\right)  \ \ \ \ \ \ \ \ \ \ \left(
\text{since }\operatorname*{Ker}\left(  \pi_{D}\right)  \text{ is a
}\mathbf{k}\text{-module}\right)  .
\]

Now, let us show the reverse inclusion. Indeed, let $v\in\operatorname*{Ker}%
\left(  \pi_{D}\right)  $. We shall prove that $v\in\operatorname*{span}%
\nolimits_{\mathbf{k}}\left(  \mathfrak{C}\cup\mathfrak{G}\right)  $. For this
purpose, we shall show that $\overline{v}=0$ in the quotient $\mathbf{k}%
$-module $\mathcal{T}^{D}/\operatorname*{span}\nolimits_{\mathbf{k}}\left(
\mathfrak{C}\cup\mathfrak{G}\right)  $.

Lemma \ref{lem.specht.TD.garnir-tria-span.straighten3} shows that the
$\mathbf{k}$-module $\mathcal{T}^{D}/\operatorname*{span}\nolimits_{\mathbf{k}%
}\left(  \mathfrak{C}\cup\mathfrak{G}\right)  $ is spanned by the family
$\left(  \overline{e_{T}}\right)  _{T\in\operatorname*{SYT}\left(  D\right)
}$ of all standard quasi-polytabloids. Hence, the vector $\overline{v}%
\in\mathcal{T}^{D}/\operatorname*{span}\nolimits_{\mathbf{k}}\left(
\mathfrak{C}\cup\mathfrak{G}\right)  $ is a $\mathbf{k}$-linear combination of
these standard quasi-polytabloids. In other words, we can write $\overline{v}$
as
\begin{equation}
\overline{v}=\sum_{T\in\operatorname*{SYT}\left(  D\right)  }\alpha
_{T}\overline{e_{T}} \label{pf.thm.specht.TD.garnir-tria-span.vbar=}%
\end{equation}
for some scalars $\alpha_{T}\in\mathbf{k}$. Consider these $\alpha_{T}$.
Hence,%
\[
\overline{\sum_{T\in\operatorname*{SYT}\left(  D\right)  }\alpha_{T}e_{T}%
}=\sum_{T\in\operatorname*{SYT}\left(  D\right)  }\alpha_{T}\overline{e_{T}%
}=\overline{v}\ \ \ \ \ \ \ \ \ \ \left(  \text{by
(\ref{pf.thm.specht.TD.garnir-tria-span.vbar=})}\right)  .
\]
In other words,
\[
\sum_{T\in\operatorname*{SYT}\left(  D\right)  }\alpha_{T}e_{T}-v\in
\operatorname*{span}\nolimits_{\mathbf{k}}\left(  \mathfrak{C}\cup
\mathfrak{G}\right)  \subseteq\operatorname*{Ker}\left(  \pi_{D}\right)  .
\]
Hence,
\begin{equation}
\pi_{D}\left(  \sum_{T\in\operatorname*{SYT}\left(  D\right)  }\alpha_{T}%
e_{T}-v\right)  =0. \label{pf.thm.specht.TD.garnir-tria-span.4}%
\end{equation}

But the map $\pi_{D}$ is $\mathbf{k}$-linear. Hence,%
\begin{align*}
\pi_{D}\left(  \sum_{T\in\operatorname*{SYT}\left(  D\right)  }\alpha_{T}%
e_{T}-v\right)   &  =\sum_{T\in\operatorname*{SYT}\left(  D\right)  }%
\alpha_{T}\underbrace{\pi_{D}\left(  e_{T}\right)  }_{\substack{=\mathbf{e}%
_{T}\\\text{(by the definition of }\pi_{D}\text{)}}}-\underbrace{\pi
_{D}\left(  v\right)  }_{\substack{=0\\\text{(since }v\in\operatorname*{Ker}%
\left(  \pi_{D}\right)  \text{)}}}\\
&  =\sum_{T\in\operatorname*{SYT}\left(  D\right)  }\alpha_{T}\mathbf{e}_{T}.
\end{align*}
Comparing this with (\ref{pf.thm.specht.TD.garnir-tria-span.4}), we obtain%
\begin{equation}
\sum_{T\in\operatorname*{SYT}\left(  D\right)  }\alpha_{T}\mathbf{e}_{T}=0.
\label{pf.thm.specht.TD.garnir-tria-span.5}%
\end{equation}

However, Theorem \ref{thm.spechtmod.linind} shows that the standard
polytabloids of shape $D$ are $\mathbf{k}$-linearly independent. In other
words, the vectors $\mathbf{e}_{T}$, where $T$ ranges over all standard
tableaux of shape $D$, are $\mathbf{k}$-linearly independent (since these
vectors are the standard polytabloids of shape $D$). In other words, the
vectors $\mathbf{e}_{T}$ with $T\in\operatorname*{SYT}\left(  D\right)  $ are
$\mathbf{k}$-linearly independent (since $\operatorname*{SYT}\left(  D\right)
$ is the set of all standard tableaux of shape $D$). Hence, from
(\ref{pf.thm.specht.TD.garnir-tria-span.5}), we conclude that all the
coefficients $\alpha_{T}$ must be $0$. In other words, $\alpha_{T}=0$ for each
$T\in\operatorname*{SYT}\left(  D\right)  $. Hence,
(\ref{pf.thm.specht.TD.garnir-tria-span.vbar=}) becomes%
\[
\overline{v}=\sum_{T\in\operatorname*{SYT}\left(  D\right)  }%
\underbrace{\alpha_{T}}_{=0}\overline{e_{T}}=0.
\]
In other words, $v\in\operatorname*{span}\nolimits_{\mathbf{k}}\left(
\mathfrak{C}\cup\mathfrak{G}\right)  $.

Forget that we fixed $v$. We thus have proved that $v\in\operatorname*{span}%
\nolimits_{\mathbf{k}}\left(  \mathfrak{C}\cup\mathfrak{G}\right)  $ for each
$v\in\operatorname*{Ker}\left(  \pi_{D}\right)  $. In other words,
$\operatorname*{Ker}\left(  \pi_{D}\right)  \subseteq\operatorname*{span}%
\nolimits_{\mathbf{k}}\left(  \mathfrak{C}\cup\mathfrak{G}\right)  $.
Combining this with $\operatorname*{span}\nolimits_{\mathbf{k}}\left(
\mathfrak{C}\cup\mathfrak{G}\right)  \subseteq\operatorname*{Ker}\left(
\pi_{D}\right)  $, we conclude that $\operatorname*{Ker}\left(  \pi
_{D}\right)  =\operatorname*{span}\nolimits_{\mathbf{k}}\left(  \mathfrak{C}%
\cup\mathfrak{G}\right)  $. This proves Theorem
\ref{thm.specht.TD.garnir-tria-span}.
\end{proof}

\begin{proof}
[Proof of Theorem \ref{thm.specht.TD.garnir-all-span}.]Define the two sets%
\begin{align}
\mathfrak{C}  &  =\left\{  \text{all column relators in }\mathcal{T}%
^{D}\right\}  \ \ \ \ \ \ \ \ \ \ \text{and}%
\label{pf.thm.specht.TD.garnir-all-span.C=}\\
\mathfrak{G}  &  =\left\{  \text{all Garnir relators in }\mathcal{T}%
^{D}\right\}  . \label{pf.thm.specht.TD.garnir-all-span.G=}%
\end{align}
Then, we note the following:

\begin{itemize}
\item The set $\mathfrak{C}$ is a triangulating set of column relators.

[\textit{Proof:} Clearly, $\mathfrak{C}$ is a set of column relators. We must
prove that this set $\mathfrak{C}$ is triangulating. In other words, we must
prove that for any $n$-tableau $T\in\operatorname*{Tab}\left(  D\right)  $
that is not column-standard, there exists a $u\in\mathcal{C}\left(  T\right)
$ such that the $n$-tableau $uT$ is column-standard and such that
$e_{T}-\left(  -1\right)  ^{u}e_{uT}\in\mathfrak{C}$ (because this is what
\textquotedblleft triangulating\textquotedblright\ means, by Definition
\ref{def.specht.TD.garnir-triasub} \textbf{(a)}).

So let us do this. Let $T\in\operatorname*{Tab}\left(  D\right)  $ be an
$n$-tableau that is not column-standard. We must show that there exists a
$u\in\mathcal{C}\left(  T\right)  $ such that the $n$-tableau $uT$ is
column-standard and such that $e_{T}-\left(  -1\right)  ^{u}e_{uT}%
\in\mathfrak{C}$.

The $n$-column-tabloid $\widetilde{T}$ contains a unique column-standard
$n$-tableau (since Lemma \ref{lem.coltabloid.col-st-unique} says that each
$n$-column-tabloid of shape $D$ contains a unique column-standard
$n$-tableau). Let $P$ be the latter column-standard $n$-tableau. Then, $P$ is
contained in the $n$-column-tabloid $\widetilde{T}$. In other words, $P$ is
column-equivalent to $T$. In other words, the $n$-tableaux $T$ and $P$ are column-equivalent.

Now, Proposition \ref{prop.tableau.req-R} \textbf{(b)} (applied to $S=P$)
shows that the $n$-tableaux $T$ and $P$ are column-equivalent if and only if
there exists some $w\in\mathcal{C}\left(  T\right)  $ such that
$P=w\rightharpoonup T$. Hence, there exists some $w\in\mathcal{C}\left(
T\right)  $ such that $P=w\rightharpoonup T$ (since the $n$-tableaux $T$ and
$P$ are column-equivalent). Consider this $w$.

Now, the definition of a column relator (Definition
\ref{def.specht.TD.relators} \textbf{(a)}, applied to $u=w$) shows that
$e_{T}-\left(  -1\right)  ^{w}e_{wT}$ is a column relator (since
$T\in\operatorname*{Tab}\left(  D\right)  $ and $w\in\mathcal{C}\left(
T\right)  $). In other words, $e_{T}-\left(  -1\right)  ^{w}e_{wT}%
\in\mathfrak{C}$ (since $\mathfrak{C}$ is the set of all column relators).
Moreover, the tableau $wT=w\rightharpoonup T=P$ is column-standard. Hence,
there exists a $u\in\mathcal{C}\left(  T\right)  $ such that the $n$-tableau
$uT$ is column-standard and such that $e_{T}-\left(  -1\right)  ^{u}e_{uT}%
\in\mathfrak{C}$ (namely, $u=w$). This is precisely what we needed to prove.
Thus, we have shown that $\mathfrak{C}$ is a triangulating set of column relators.]

\item The set $\mathfrak{G}$ is a triangulating set of Garnir relators.

[\textit{Proof:} Clearly, $\mathfrak{G}$ is a set of Garnir relators. We must
prove that this set $\mathfrak{G}$ is triangulating. In other words, we must
prove that for any column-standard $n$-tableau $T\in\operatorname*{Tab}\left(
D\right)  $ that is not standard, there exists a triangulating Garnir pentuple
$\left(  T,j,k,Z,L\right)  $ with first entry $T$ such that $\sum_{x\in
L}\left(  -1\right)  ^{x}e_{xT}\in\mathfrak{G}$ (because this is what
\textquotedblleft triangulating\textquotedblright\ means, by Definition
\ref{def.specht.TD.garnir-triasub} \textbf{(c)}).

So let us do this. Let $T\in\operatorname*{Tab}\left(  D\right)  $ be a
column-standard $n$-tableau that is not standard. We must construct a
triangulating Garnir pentuple $\left(  T,j,k,Z,L\right)  $ with first entry
$T$ such that $\sum_{x\in L}\left(  -1\right)  ^{x}e_{xT}\in\mathfrak{G}$.

The tableau $T$ is column-standard but not standard. Hence, $T$ is not
row-standard (since otherwise, $T$ would be both column-standard and
row-standard, and therefore standard). Therefore, $T$ has two entries in the
same row that are in the wrong order (i.e., the entry further left is $\geq$
to the entry further right). In other words, there exist two entries $T\left(
i,j\right)  $ and $T\left(  i,k\right)  $ of $T$ with $j<k$ and $T\left(
i,j\right)  \geq T\left(  i,k\right)  $. Consider these $i$, $j$ and $k$.
Define two sets $Z$ and $L$ as in the above proof of Lemma
\ref{lem.spechtmod.straighten1}. In that proof, we have shown that

\begin{itemize}
\item the $j$ and $k$ are two distinct integers;

\item the set $Z$ is a subset of $\operatorname*{Col}\left(  j,T\right)
\cup\operatorname*{Col}\left(  k,T\right)  $ such that $\left\vert
Z\right\vert >\left\vert D\left\langle j\right\rangle \cup D\left\langle
k\right\rangle \right\vert $;

\item the set $L$ is a left transversal of $S_{n,Z}\cap\mathcal{C}\left(
T\right)  $ in $S_{n,Z}$.
\end{itemize}

\medskip Therefore, $\left(  T,j,k,Z,L\right)  $ is a Garnir pentuple (by
Definition \ref{def.specht.TD.relators} \textbf{(b)}). Moreover, each $x\in
L\setminus\mathcal{C}\left(  T\right)  $ satisfies $\widetilde{xT}%
>\widetilde{T}$ in the column last letter order (as we showed during our above
proof of Lemma \ref{lem.spechtmod.straighten1}). In other words, the Garnir
pentuple $\left(  T,j,k,Z,L\right)  $ is triangulating (by Definition
\ref{def.specht.TD.garnir-triasub} \textbf{(b)}). Finally, $\sum_{x\in
L}\left(  -1\right)  ^{x}e_{xT}$ is a Garnir relator (by Definition
\ref{def.specht.TD.relators} \textbf{(c)}), and thus belongs to $\mathfrak{G}$
(since $\mathfrak{G}$ is the set of all Garnir relators). Thus, we have
constructed a triangulating Garnir pentuple $\left(  T,j,k,Z,L\right)  $ with
first entry $T$ such that $\sum_{x\in L}\left(  -1\right)  ^{x}e_{xT}%
\in\mathfrak{G}$. As we said above, this completes our proof that
$\mathfrak{G}$ is a triangulating set of Garnir relators.]
\end{itemize}

Hence, we can apply Theorem \ref{thm.specht.TD.garnir-tria-span} and conclude
that%
\begin{align*}
\operatorname*{Ker}\left(  \pi_{D}\right)   &  =\operatorname*{span}%
\nolimits_{\mathbf{k}}\left(  \mathfrak{C}\cup\mathfrak{G}\right) \\
&  =\operatorname*{span}\nolimits_{\mathbf{k}}\left(  \left\{  \text{all
column relators in }\mathcal{T}^{D}\right\}  \cup\left\{  \text{all Garnir
relators in }\mathcal{T}^{D}\right\}  \right)
\end{align*}
(by (\ref{pf.thm.specht.TD.garnir-all-span.C=}) and
(\ref{pf.thm.specht.TD.garnir-all-span.G=})). Hence, Theorem
\ref{thm.specht.TD.garnir-all-span} is proved.
\end{proof}
\end{fineprint}

\subsection{\label{sec.specht.cogarnir}Linear equations defining the Specht
module}

The Specht module $\mathcal{S}^{D}$ corresponding to a diagram $D$ was defined
as a span of certain elements of the (usually much larger) Young module
$\mathcal{M}^{D}$ (see Definition \ref{def.spechtmod.spechtmod} \textbf{(b)}).
A natural question to ask is whether it can also be described as the solution
set of a system of linear equations -- i.e., as a kernel of a linear map. In
this section, we shall present such a description, at least in the case when
$D$ is a skew Young diagram. The linear equations in this description will be
called the \emph{coGarnir equations}, due to the fact that we derive them from
the Garnir relations. They hold for all vectors in $\mathcal{S}^{D}$ for any
diagram $D$, but they do not determine $\mathcal{S}^{D}$ in general; however,
they do when $D$ is a skew Young diagram.

\subsubsection{\label{subsec.specht.cogarnir.cog-eqs}The coGarnir equations}

To state the coGarnir equations, we need some preparation. First, we give
analogues of Definition \ref{def.diagram.Dj} and Definition
\ref{def.tableau.ColTj}:

\begin{definition}
\label{def.diagram.Djr}Let $D$ be a diagram. Let $i\in\mathbb{Z}$. Then,
$D\left\lfloor i\right\rfloor $ shall mean the set of all integers $j$ that
satisfy $\left(  i,j\right)  \in D$.
\end{definition}

For example, if
\[
D=\left\{  \left(  1,2\right)  ,\ \left(  1,3\right)  ,\ \left(  2,1\right)
,\ \left(  2,2\right)  ,\ \left(  3,3\right)  ,\ \left(  4,1\right)  \right\}
=%
%TCIMACRO{\TeXButton{tikz weird diagram}{\begin{tikzpicture}[scale=0.7]
%\draw[fill=red!50] (0, 0) rectangle (1, 1);
%\draw[fill=red!50] (2, 1) rectangle (3, 2);
%\draw[fill=red!50] (0, 2) rectangle (1, 3);
%\draw[fill=red!50] (1, 2) rectangle (2, 3);
%\draw[fill=red!50] (1, 3) rectangle (2, 4);
%\draw[fill=red!50] (2, 3) rectangle (3, 4);
%\end{tikzpicture}}}%
%BeginExpansion
\begin{tikzpicture}[scale=0.7]
\draw[fill=red!50] (0, 0) rectangle (1, 1);
\draw[fill=red!50] (2, 1) rectangle (3, 2);
\draw[fill=red!50] (0, 2) rectangle (1, 3);
\draw[fill=red!50] (1, 2) rectangle (2, 3);
\draw[fill=red!50] (1, 3) rectangle (2, 4);
\draw[fill=red!50] (2, 3) rectangle (3, 4);
\end{tikzpicture}%
%EndExpansion
\ \ ,
\]
then $D\left\lfloor 1\right\rfloor =\left\{  2,3\right\}  $ and $D\left\lfloor
2\right\rfloor =\left\{  1,2\right\}  $ and $D\left\lfloor 3\right\rfloor
=\left\{  3\right\}  $ and $D\left\lfloor 4\right\rfloor =\left\{  1\right\}
$ and $D\left\lfloor i\right\rfloor =\varnothing$ for all $i\notin\left[
4\right]  $. If $D$ is a skew Young diagram $Y\left(  \lambda/\mu\right)  $,
then each set $D\left\lfloor i\right\rfloor $ is an integer interval
(equalling $\left\{  \mu_{i}+1,\ \mu_{i}+2,\ \ldots,\ \lambda_{i}\right\}  $
when $i>0$ and being empty otherwise), but of course this does not generally
hold for bad-shape diagrams.

\begin{definition}
\label{def.tableau.RowTi}Let $T$ be an $n$-tableau (of any shape). Let
$i\in\mathbb{Z}$. Then, $\operatorname*{Row}\left(  i,T\right)  $ shall denote
the set of all entries in the $i$-th row of $T$.
\end{definition}

The following (somewhat ugly) notation will shorten some of our formulas:

\begin{definition}
Let $D$ be a diagram with $\left\vert D\right\vert =n$. Then, $\overline
{\operatorname{Tab}}\left(  D\right)  $ shall denote the set $\left\{
n\text{-tabloids of shape }D\right\}  $. As we know from Definition
\ref{def.tabloid.Sn-act}, this set is a left $S_{n}$-set.
\end{definition}

As we recall (see Definition \ref{def.youngmod.youngmod}), the Young module
$\mathcal{M}^{D}$ is the permutation module corresponding to this left $S_{n}%
$-set $\overline{\operatorname{Tab}}\left(  D\right)  $. Hence, as a
$\mathbf{k}$-module, it is the free $\mathbf{k}$-module $\mathbf{k}^{\left(
\overline{\operatorname{Tab}}\left(  D\right)  \right)  }$ and has a basis
$\left(  \overline{P}\right)  _{\overline{P}\in\overline{\operatorname{Tab}%
}\left(  D\right)  }$. An element $\mathbf{a}$ of $\mathcal{M}^{D}$ can thus
be written as a family $\left(  a_{\overline{P}}\right)  _{\overline{P}%
\in\overline{\operatorname{Tab}}\left(  D\right)  }$ of scalars $a_{\overline
{P}}\in\mathbf{k}$, or, equivalently, as a formal $\mathbf{k}$-linear
combination $\sum_{\overline{P}\in\overline{\operatorname{Tab}}\left(
D\right)  }a_{\overline{P}}\overline{P}$ with coefficients $a_{\overline{P}%
}\in\mathbf{k}$.

We now define the coGarnir equations. As the name suggests, these are somewhat
similar (in their construction) to the Garnir relations (Theorem
\ref{thm.garnir.grel}); in particular, they too rely on left transversals of
subgroups (Definition \ref{def.G/H-transversal}). However, they differ in a
few other aspects (rows are used instead of columns, and permutations appear
without sign). Here is the precise definition:

\begin{definition}
\label{def.cogarnir.cogarnirTL}Let $D$ be a diagram with $\left\vert
D\right\vert =n$. Let $T$ be an $n$-tableau of shape $D$. Let $j$ and $k$ be
two distinct integers. Let $Z$ be a subset of $\operatorname*{Row}\left(
j,T\right)  \cup\operatorname*{Row}\left(  k,T\right)  $ such that $\left\vert
Z\right\vert >\left\vert D\left\lfloor j\right\rfloor \cup D\left\lfloor
k\right\rfloor \right\vert $. (Here we are using the notations introduced in
Definition \ref{def.diagram.Djr} and Definition \ref{def.tableau.RowTi}.)

Define $S_{n,Z}$ as in Proposition \ref{prop.intX.basics}. This set $S_{n,Z}$
is a subgroup of $S_{n}$ (by Proposition \ref{prop.intX.basics} \textbf{(a)},
applied to $X=Z$). Thus, both $\mathcal{R}\left(  T\right)  $ and $S_{n,Z}$
are subgroups of $S_{n}$. Hence, their intersection $S_{n,Z}\cap
\mathcal{R}\left(  T\right)  $ is a subgroup of $S_{n,Z}$. Let $L$ be a left
transversal of $S_{n,Z}\cap\mathcal{R}\left(  T\right)  $ in $S_{n,Z}$. Then:
\medskip

\textbf{(a)} The set $L$ is called a \emph{row-Garnir set} for $T$. \medskip

\textbf{(b)} Now let $\mathbf{a}=\left(  a_{\overline{P}}\right)
_{\overline{P}\in\overline{\operatorname{Tab}}\left(  D\right)  }%
\in\mathcal{M}^{D}$ be an element of the Young module $\mathcal{M}^{D}$. Then,
we say that $\mathbf{a}$ satisfies the $\left(  T,L\right)  $\emph{-coGarnir
equation} if we have%
\begin{equation}
\sum_{x\in L}a_{\overline{xT}}=0. \label{eq.def.cogarnir.cogarnirTL.b.eq}%
\end{equation}

\end{definition}

\begin{definition}
\label{def.cogarnir.cogarnirall}Let $D$ be a diagram with $\left\vert
D\right\vert =n$. Let $\mathbf{a}\in\mathcal{M}^{D}$ be an element of the
Young module $\mathcal{M}^{D}$. Then, we say that $\mathbf{a}$ satisfies
\emph{all coGarnir equations} if $\mathbf{a}$ satisfies the $\left(
T,L\right)  $-coGarnir equation for each $n$-tableau $T$ of shape $D$ and each
row-Garnir set $L$ for $T$.
\end{definition}

\begin{example}
\label{exa.cogarnir.cogarnirTL.1}Let $n=7$ and $D=Y\left(  3,3,1\right)  $
and
\[
T=\ytableaushort{1{*(green)2}{*(green)3},{*(green)4}{*(green)5}6,7}=123\backslash
\backslash456\backslash\backslash7
\]
(using the shorthand notation from Convention \ref{conv.tableau.poetic}). Set
$j=1$ and $k=2$ and $Z=\left\{  2,3,4,5\right\}  $ (the set of entries
highlighted in green in the tableau shown above). These $j$, $k$ and $Z$
satisfy all requirements of Definition \ref{def.cogarnir.cogarnirTL}.

Now, let us choose a left transversal $L$ of $S_{n,Z}\cap\mathcal{R}\left(
T\right)  $ in $S_{n,Z}$, for example%
\[
L=\left\{  \operatorname*{id},\ t_{2,4},\ t_{2,5},\ t_{3,4},\ t_{3,5}%
,\ t_{2,4}t_{3,5}\right\}  .
\]
Then, $L$ is a row-Garnir set for $T$. Thus, an element $\mathbf{a}=\left(
a_{\overline{P}}\right)  _{\overline{P}\in\overline{\operatorname{Tab}}\left(
D\right)  }\in\mathcal{M}^{D}$ satisfies the $\left(  T,L\right)  $-coGarnir
equation if and only if%
\[
a_{\overline{\operatorname*{id}T}}+a_{\overline{t_{2,4}T}}+a_{\overline
{t_{2,5}T}}+a_{\overline{t_{3,4}T}}+a_{\overline{t_{3,5}T}}+a_{\overline
{t_{2,4}t_{3,5}T}}=0,
\]
that is, if and only if%
\begin{align*}
&  a_{\overline{123\backslash\backslash456\backslash\backslash7}}%
+a_{\overline{143\backslash\backslash256\backslash\backslash7}}+a_{\overline
{153\backslash\backslash426\backslash\backslash7}}+a_{\overline{124\backslash
\backslash356\backslash\backslash7}}\\
&  \ \ \ \ \ \ \ \ \ \ +a_{\overline{125\backslash\backslash436\backslash
\backslash7}}+a_{\overline{145\backslash\backslash236\backslash\backslash7}%
}=0.
\end{align*}
Note that there are no minus signs here!

Of course, for $\mathbf{a}$ to satisfy all coGarnir equations, several more
equalities need to hold.
\end{example}

\begin{exercise}
\fbox{2} Prove that the $\left(  T,L\right)  $-coGarnir equation (as defined
in Definition \ref{def.cogarnir.cogarnirTL}) can be equivalently rewritten as
\[
\sum_{\overline{P}\in\Omega}a_{\overline{P}}=0,
\]
where $\Omega$ is the orbit of $\overline{T}\in\overline{\operatorname{Tab}%
}\left(  D\right)  $ under the action of the subgroup $S_{n,Z}$ of $S_{n}$.
Thus, this equation depends only on $T$ and $Z$, not on $L$.
\end{exercise}

We now state our main results:

\begin{theorem}
\label{thm.cogarnir.cogarnir-sats}Let $D$ be a diagram with $\left\vert
D\right\vert =n$. Let $\mathbf{a}\in\mathcal{S}^{D}$ be an element of the
Specht module $\mathcal{S}^{D}$. Then, $\mathbf{a}$ (viewed as an element of
the Young module $\mathcal{M}^{D}$) satisfies all coGarnir equations.
\end{theorem}

\begin{theorem}
\label{thm.cogarnir.cogarnir-skew}Let $D$ be a skew Young diagram with
$\left\vert D\right\vert =n$. Then, the Specht module $\mathcal{S}^{D}$ can be
described as%
\[
\mathcal{S}^{D}=\left\{  \mathbf{a}\in\mathcal{M}^{D}\ \mid\ \mathbf{a}\text{
satisfies all coGarnir equations}\right\}  .
\]
(This is an identity, not just an isomorphism!)
\end{theorem}

We shall prove these two theorems in the three subsections that follow.

\subsubsection{\label{subsec.specht.cogarnir.garnir-prep}Transposing the
Garnir relations}

In preparation for proving Theorem \ref{thm.cogarnir.cogarnir-sats}, we shall
transform the (original) Garnir relations (\ref{eq.thm.garnir.grel.gar}).
First, we translate (\ref{eq.thm.garnir.grel.gar}) from the Specht module
$\mathcal{S}^{D}$ to its left ideal avatar $\mathbf{k}\left[  S_{n}\right]
\cdot\mathbf{E}_{T}$:

\begin{corollary}
\label{cor.hlf-young.grel2}Let $D$ be a diagram. Let $T$ be an $n$-tableau of
shape $D$. Let $j$ and $k$ be two distinct integers. Let $Z$ be a subset of
$\operatorname*{Col}\left(  j,T\right)  \cup\operatorname*{Col}\left(
k,T\right)  $ such that $\left\vert Z\right\vert >\left\vert D\left\langle
j\right\rangle \cup D\left\langle k\right\rangle \right\vert $. (Here we are
using the notations introduced in Definition \ref{def.diagram.Dj} and
Definition \ref{def.tableau.ColTj}.)

Define $S_{n,Z}$ as in Proposition \ref{prop.intX.basics}. This set $S_{n,Z}$
is a subgroup of $S_{n}$ (by Proposition \ref{prop.intX.basics} \textbf{(a)},
applied to $X=Z$). Thus, both $\mathcal{C}\left(  T\right)  $ and $S_{n,Z}$
are subgroups of $S_{n}$. Hence, their intersection $S_{n,Z}\cap
\mathcal{C}\left(  T\right)  $ is a subgroup of $S_{n,Z}$. Let $L$ be a left
transversal of $S_{n,Z}\cap\mathcal{C}\left(  T\right)  $ in $S_{n,Z}$. Then,
in $\mathbf{k}\left[  S_{n}\right]  $, we have%
\[
\sum_{x\in L}\left(  -1\right)  ^{x}x\nabla_{\operatorname*{Col}T}^{-}%
\nabla_{\operatorname*{Row}T}=0.
\]

\end{corollary}

\begin{proof}
We have $\left\vert D\right\vert =n$ (since $T$ is an $n$-tableau of shape
$D$). Corollary \ref{cor.spechtmod.leftideal.S} thus shows that there is a
left $\mathbf{k}\left[  S_{n}\right]  $-module isomorphism%
\[
\overline{\alpha}:\mathbf{k}\left[  S_{n}\right]  \cdot\nabla
_{\operatorname*{Col}T}^{-}\nabla_{\operatorname*{Row}T}\rightarrow
\mathcal{S}^{D}%
\]
that sends $\nabla_{\operatorname*{Col}T}^{-}\nabla_{\operatorname*{Row}T}$ to
$\mathbf{e}_{T}$. Consider this $\overline{\alpha}$.

Now, (\ref{eq.thm.garnir.grel.gar}) yields%
\begin{equation}
\sum_{x\in L}\left(  -1\right)  ^{x}\mathbf{e}_{xT}=0.
\label{pf.cor.hlf-young.grel2.2}%
\end{equation}
But $\underbrace{\sum_{x\in L}\left(  -1\right)  ^{x}x}_{\in\mathbf{k}\left[
S_{n}\right]  }\nabla_{\operatorname*{Col}T}^{-}\nabla_{\operatorname*{Row}%
T}\in\mathbf{k}\left[  S_{n}\right]  \cdot\nabla_{\operatorname*{Col}T}%
^{-}\nabla_{\operatorname*{Row}T}$, so that we can apply the map
$\overline{\alpha}$ to $\sum_{x\in L}\left(  -1\right)  ^{x}x\nabla
_{\operatorname*{Col}T}^{-}\nabla_{\operatorname*{Row}T}$. We thus obtain%
\begin{align*}
&  \overline{\alpha}\left(  \sum_{x\in L}\left(  -1\right)  ^{x}%
x\nabla_{\operatorname*{Col}T}^{-}\nabla_{\operatorname*{Row}T}\right) \\
&  =\sum_{x\in L}\left(  -1\right)  ^{x}x\underbrace{\overline{\alpha}\left(
\nabla_{\operatorname*{Col}T}^{-}\nabla_{\operatorname*{Row}T}\right)
}_{\substack{=\mathbf{e}_{T}\\\text{(since }\overline{\alpha}\text{ sends
}\nabla_{\operatorname*{Col}T}^{-}\nabla_{\operatorname*{Row}T}\text{ to
}\mathbf{e}_{T}\text{)}}}\ \ \ \ \ \ \ \ \ \ \left(
\begin{array}
[c]{c}%
\text{since }\overline{\alpha}\text{ is a left }\mathbf{k}\left[
S_{n}\right]  \text{-module}\\
\text{morphism}%
\end{array}
\right) \\
&  =\sum_{x\in L}\left(  -1\right)  ^{x}\underbrace{x\mathbf{e}_{T}%
}_{\substack{=\mathbf{e}_{xT}\\\text{(since Lemma \ref{lem.spechtmod.submod}
\textbf{(a)} (applied to }u=x\text{)}\\\text{yields }\mathbf{e}_{xT}%
=x\mathbf{e}_{T}\text{)}}}=\sum_{x\in L}\left(  -1\right)  ^{x}\mathbf{e}%
_{xT}=0
\end{align*}
(by (\ref{pf.cor.hlf-young.grel2.2})). In other words, $\sum_{x\in L}\left(
-1\right)  ^{x}x\nabla_{\operatorname*{Col}T}^{-}\nabla_{\operatorname*{Row}%
T}\in\operatorname*{Ker}\overline{\alpha}=0$ (since $\overline{\alpha}$ is an
isomorphism and thus injective). Hence, $\sum_{x\in L}\left(  -1\right)
^{x}x\nabla_{\operatorname*{Col}T}^{-}\nabla_{\operatorname*{Row}T}=0$. This
proves Corollary \ref{cor.hlf-young.grel2}.
\end{proof}

Next, we transform Corollary \ref{cor.hlf-young.grel2} by reflecting the
diagram $D$ and the tableau $T$ across the diagonal:

\begin{corollary}
\label{cor.hlf-young.grel3}Let $D$ be a diagram. Let $T$ be an $n$-tableau of
shape $D$. Let $j$ and $k$ be two distinct integers. Let $Z$ be a subset of
$\operatorname*{Row}\left(  j,T\right)  \cup\operatorname*{Row}\left(
k,T\right)  $ such that $\left\vert Z\right\vert >\left\vert D\left\lfloor
j\right\rfloor \cup D\left\lfloor k\right\rfloor \right\vert $. (Here we are
using the notations introduced in Definition \ref{def.diagram.Djr} and
Definition \ref{def.tableau.RowTi}.)

Define $S_{n,Z}$ as in Proposition \ref{prop.intX.basics}. This set $S_{n,Z}$
is a subgroup of $S_{n}$ (by Proposition \ref{prop.intX.basics} \textbf{(a)},
applied to $X=Z$). Thus, both $\mathcal{R}\left(  T\right)  $ and $S_{n,Z}$
are subgroups of $S_{n}$. Hence, their intersection $S_{n,Z}\cap
\mathcal{R}\left(  T\right)  $ is a subgroup of $S_{n,Z}$. Let $L$ be a left
transversal of $S_{n,Z}\cap\mathcal{R}\left(  T\right)  $ in $S_{n,Z}$. Then,
in $\mathbf{k}\left[  S_{n}\right]  $, we have%
\[
\sum_{x\in L}x\nabla_{\operatorname*{Row}T}\nabla_{\operatorname*{Col}T}%
^{-}=0.
\]

\end{corollary}

\begin{proof}
[Proof of Corollary \ref{cor.hlf-young.grel3}.]We shall use the notations
introduced in Definition \ref{def.diagram.Dj} and Definition
\ref{def.tableau.ColTj} as well.

Consider the map $\mathbf{r}:\mathbb{Z}^{2}\rightarrow\mathbb{Z}^{2}$ from
Theorem \ref{thm.partitions.conj}. This map $\mathbf{r}$ is an involution, so
that it is invertible and its inverse is $\mathbf{r}^{-1}=\mathbf{r}$. For any
cell $\left(  i,j\right)  \in\mathbb{Z}^{2}$, we have the following chain of
logical equivalences:%
\begin{align}
\left(  \left(  i,j\right)  \in\mathbf{r}\left(  D\right)  \right)  \  &
\Longleftrightarrow\ \left(  \mathbf{r}^{-1}\left(  i,j\right)  \in D\right)
\nonumber\\
&  \Longleftrightarrow\ \left(  \left(  j,i\right)  \in D\right)
\label{pf.cor.hlf-young.grel3.equiv1}%
\end{align}
(since $\underbrace{\mathbf{r}^{-1}}_{=\mathbf{r}}\left(  i,j\right)
=\mathbf{r}\left(  i,j\right)  =\left(  j,i\right)  $ (by the definition of
$\mathbf{r}$)).

In Definition \ref{def.tableaux.r}, we introduced an $n$-tableau $T\mathbf{r}$
of shape $\mathbf{r}\left(  D\right)  $. This $n$-tableau $T\mathbf{r}$ is
obtained by reflecting $T$ across the northwest-to-southeast diagonal. This
reflection turns the rows of $T$ into the columns of $T\mathbf{r}$. More
specifically, it turns the $i$-th row of $T$ into the $i$-th column of
$T\mathbf{r}$ for each $i\in\mathbb{Z}$. Thus, for each $i\in\mathbb{Z}$, the
entries in the $i$-th row of $T$ are precisely the entries in the $i$-th
column of $T\mathbf{r}$. Hence, for each $i\in\mathbb{Z}$, we have%
\begin{equation}
\operatorname*{Row}\left(  i,T\right)  =\operatorname*{Col}\left(
i,T\mathbf{r}\right)  \label{pf.cor.hlf-young.grel3.2}%
\end{equation}
(since $\operatorname*{Row}\left(  i,T\right)  $ is the set of all entries in
the $i$-th row of $T$, while $\operatorname*{Col}\left(  i,T\mathbf{r}\right)
$ is the set of all entries in the $i$-th column of $T\mathbf{r}$). Thus, in
particular,%
\[
\underbrace{\operatorname*{Row}\left(  j,T\right)  }%
_{\substack{=\operatorname*{Col}\left(  j,T\mathbf{r}\right)  \\\text{(by
(\ref{pf.cor.hlf-young.grel3.2}))}}}\cup\underbrace{\operatorname*{Row}\left(
k,T\right)  }_{\substack{=\operatorname*{Col}\left(  k,T\mathbf{r}\right)
\\\text{(by (\ref{pf.cor.hlf-young.grel3.2}))}}}=\operatorname*{Col}\left(
j,T\mathbf{r}\right)  \cup\operatorname*{Col}\left(  k,T\mathbf{r}\right)  .
\]
Hence, $Z$ is a subset of $\operatorname*{Col}\left(  j,T\mathbf{r}\right)
\cup\operatorname*{Col}\left(  k,T\mathbf{r}\right)  $ (since $Z$ is a subset
of $\operatorname*{Row}\left(  j,T\right)  \cup\operatorname*{Row}\left(
k,T\right)  $).

Moreover, each $p\in\mathbb{Z}$ satisfies
\begin{equation}
D\left\lfloor p\right\rfloor =\left(  \mathbf{r}\left(  D\right)  \right)
\left\langle p\right\rangle \label{pf.cor.hlf-young.grel3.3}%
\end{equation}
\footnote{\textit{Proof.} Let $p\in\mathbb{Z}$. Then, $D\left\lfloor
p\right\rfloor $ is defined as the set of all integers $j$ that satisfy
$\left(  p,j\right)  \in D$. In other words,%
\begin{equation}
D\left\lfloor p\right\rfloor =\left\{  j\in\mathbb{Z}\ \mid\ \left(
p,j\right)  \in D\right\}  =\left\{  i\in\mathbb{Z}\ \mid\ \left(  p,i\right)
\in D\right\}  \label{pf.cor.hlf-young.grel3.3.pf.1}%
\end{equation}
(here, we have renamed the index $j$ as $i$). On the other hand, $\left(
\mathbf{r}\left(  D\right)  \right)  \left\langle p\right\rangle $ is defined
as the set of all integers $i$ that satisfy $\left(  i,p\right)  \in
\mathbf{r}\left(  D\right)  $ (by Definition \ref{def.diagram.Dj}). In other
words,%
\begin{equation}
\left(  \mathbf{r}\left(  D\right)  \right)  \left\langle p\right\rangle
=\left\{  i\in\mathbb{Z}\ \mid\ \left(  i,p\right)  \in\mathbf{r}\left(
D\right)  \right\}  . \label{pf.cor.hlf-young.grel3.3.pf.2}%
\end{equation}
\par
However, if $i\in\mathbb{Z}$ is any integer, then the logical equivalence
$\left(  \left(  i,p\right)  \in\mathbf{r}\left(  D\right)  \right)
\ \Longleftrightarrow\ \left(  \left(  p,i\right)  \in D\right)  $ holds (by
(\ref{pf.cor.hlf-young.grel3.equiv1}), applied to $\left(  i,p\right)  $
instead of $\left(  i,j\right)  $). Thus, the right hand sides of the
equalities (\ref{pf.cor.hlf-young.grel3.3.pf.2}) and
(\ref{pf.cor.hlf-young.grel3.3.pf.1}) are equal (since the condition $\left(
i,p\right)  \in\mathbf{r}\left(  D\right)  $ in the former can be replaced by
the equivalent condition $\left(  p,i\right)  \in D$ in the latter). Thus, the
left hand sides of these two equalities are also equal. In other words,
\[
\left(  \mathbf{r}\left(  D\right)  \right)  \left\langle p\right\rangle
=D\left\lfloor p\right\rfloor .
\]
Thus, (\ref{pf.cor.hlf-young.grel3.3}) is proved.}. Hence, in particular,
$D\left\lfloor j\right\rfloor =\left(  \mathbf{r}\left(  D\right)  \right)
\left\langle j\right\rangle $ and $D\left\lfloor k\right\rfloor =\left(
\mathbf{r}\left(  D\right)  \right)  \left\langle k\right\rangle $.

We assumed that $\left\vert Z\right\vert >\left\vert D\left\lfloor
j\right\rfloor \cup D\left\lfloor k\right\rfloor \right\vert $. In view of
$D\left\lfloor j\right\rfloor =\left(  \mathbf{r}\left(  D\right)  \right)
\left\langle j\right\rangle $ and $D\left\lfloor k\right\rfloor =\left(
\mathbf{r}\left(  D\right)  \right)  \left\langle k\right\rangle $, we can
rewrite this as%
\[
\left\vert Z\right\vert >\left\vert \left(  \mathbf{r}\left(  D\right)
\right)  \left\langle j\right\rangle \cup\left(  \mathbf{r}\left(  D\right)
\right)  \left\langle k\right\rangle \right\vert .
\]

From Proposition \ref{prop.tableaux.r} \textbf{(e)}, we know that
$\mathcal{R}\left(  T\right)  =\mathcal{C}\left(  T\mathbf{r}\right)  $ and
$\mathcal{C}\left(  T\right)  =\mathcal{R}\left(  T\mathbf{r}\right)  $.

Recall that $L$ is a left transversal of $S_{n,Z}\cap\mathcal{R}\left(
T\right)  $ in $S_{n,Z}$. In other words, $L$ is a left transversal of
$S_{n,Z}\cap\mathcal{C}\left(  T\mathbf{r}\right)  $ in $S_{n,Z}$ (since
$\mathcal{R}\left(  T\right)  =\mathcal{C}\left(  T\mathbf{r}\right)  $).
Hence, we can apply Corollary \ref{cor.hlf-young.grel2} to $\mathbf{r}\left(
D\right)  $ and $T\mathbf{r}$ instead of $D$ and $T$. As a result, we obtain%
\begin{equation}
\sum_{x\in L}\left(  -1\right)  ^{x}x\nabla_{\operatorname*{Col}\left(
T\mathbf{r}\right)  }^{-}\nabla_{\operatorname*{Row}\left(  T\mathbf{r}%
\right)  }=0. \label{pf.def.tableau.RowTi.6}%
\end{equation}

Now, recall the map $T_{\operatorname*{sign}}:\mathbf{k}\left[  S_{n}\right]
\rightarrow\mathbf{k}\left[  S_{n}\right]  $ from Definition
\ref{def.Tsign.Tsign}. (This map has nothing to do with the tableau $T$,
despite the \textquotedblleft$T$\textquotedblright\ in its name.) Theorem
\ref{thm.Tsign.auto} \textbf{(a)} says that this map $T_{\operatorname*{sign}%
}:\mathbf{k}\left[  S_{n}\right]  \rightarrow\mathbf{k}\left[  S_{n}\right]  $
is a $\mathbf{k}$-algebra automorphism. Hence, $T_{\operatorname*{sign}}$ is
invertible, thus injective.

Furthermore, $T_{\operatorname*{sign}}$ is a $\mathbf{k}$-algebra morphism.
Thus,%
\begin{align*}
&  T_{\operatorname*{sign}}\left(  \sum_{x\in L}x\nabla_{\operatorname*{Row}%
T}\nabla_{\operatorname*{Col}T}^{-}\right) \\
&  =\sum_{x\in L}\underbrace{T_{\operatorname*{sign}}\left(  x\right)
}_{\substack{=\left(  -1\right)  ^{x}x\\\text{(by the definition of
}T_{\operatorname*{sign}}\text{,}\\\text{since }x\in L\subseteq S_{n}\text{)}%
}}\cdot\underbrace{T_{\operatorname*{sign}}\left(  \nabla_{\operatorname*{Row}%
T}\right)  }_{\substack{=\nabla_{\operatorname*{Col}\left(  T\mathbf{r}%
\right)  }^{-}\\\text{(since (\ref{eq.prop.specht.ET.r.2}) (applied to
}P=T\text{)}\\\text{yields }\nabla_{\operatorname*{Col}\left(  T\mathbf{r}%
\right)  }^{-}=T_{\operatorname*{sign}}\left(  \nabla_{\operatorname*{Row}%
T}\right)  \text{)}}}\cdot\underbrace{T_{\operatorname*{sign}}\left(
\nabla_{\operatorname*{Col}T}^{-}\right)  }_{\substack{=\nabla
_{\operatorname*{Row}\left(  T\mathbf{r}\right)  }\\\text{(since
(\ref{eq.prop.specht.ET.r.1}) (applied to }P=T\text{)}\\\text{yields }%
\nabla_{\operatorname*{Row}\left(  T\mathbf{r}\right)  }%
=T_{\operatorname*{sign}}\left(  \nabla_{\operatorname*{Col}T}^{-}\right)
\text{)}}}\\
&  =\sum_{x\in L}\left(  -1\right)  ^{x}x\nabla_{\operatorname*{Col}\left(
T\mathbf{r}\right)  }^{-}\nabla_{\operatorname*{Row}\left(  T\mathbf{r}%
\right)  }=0\ \ \ \ \ \ \ \ \ \ \left(  \text{by (\ref{pf.def.tableau.RowTi.6}%
)}\right)  .
\end{align*}
Thus, $\sum_{x\in L}x\nabla_{\operatorname*{Row}T}\nabla_{\operatorname*{Col}%
T}^{-}\in\operatorname*{Ker}\left(  T_{\operatorname*{sign}}\right)  =0$
(since $T_{\operatorname*{sign}}$ is injective), so that $\sum_{x\in L}%
x\nabla_{\operatorname*{Row}T}\nabla_{\operatorname*{Col}T}^{-}=0$. This
proves Corollary \ref{cor.hlf-young.grel3}.
\end{proof}

We can restate Corollary \ref{cor.hlf-young.grel3} in a form particularly
convenient for our needs:

\begin{corollary}
\label{cor.hlf-young.grel3rg}Let $D$ be a diagram. Let $T$ be an $n$-tableau
of shape $D$. Let $L$ be a row-Garnir set for $T$. Define the Young
symmetrizer $\mathbf{E}_{T}$ as in Definition \ref{def.specht.ET.defs}
\textbf{(b)}. Then, in $\mathbf{k}\left[  S_{n}\right]  $, we have%
\[
\mathbf{E}_{T}\sum_{x\in L}x^{-1}=0.
\]

\end{corollary}

\begin{proof}
Recall that $L$ is a row-Garnir set for $T$. According to Definition
\ref{def.cogarnir.cogarnirTL} \textbf{(a)}, this means that there exist two
distinct integers $j$ and $k$ as well as a subset $Z$ of $\operatorname*{Row}%
\left(  j,T\right)  \cup\operatorname*{Row}\left(  k,T\right)  $ such that
$\left\vert Z\right\vert >\left\vert D\left\lfloor j\right\rfloor \cup
D\left\lfloor k\right\rfloor \right\vert $ and such that $L$ is a left
transversal of $S_{n,Z}\cap\mathcal{R}\left(  T\right)  $ in $S_{n,Z}$.
Consider these $j$, $k$ and $Z$. Then, Corollary \ref{cor.hlf-young.grel3}
shows that%
\begin{equation}
\sum_{x\in L}x\nabla_{\operatorname*{Row}T}\nabla_{\operatorname*{Col}T}%
^{-}=0. \label{pf.cor.hlf-young.grel3rg.1}%
\end{equation}

However, Theorem \ref{thm.S.auto} \textbf{(a)} says that the map
$S:\mathbf{k}\left[  S_{n}\right]  \rightarrow\mathbf{k}\left[  S_{n}\right]
$ is a $\mathbf{k}$-algebra anti-automorphism. Applying this map $S$ to the
equality (\ref{pf.cor.hlf-young.grel3rg.1}), we obtain%
\[
S\left(  \sum_{x\in L}x\nabla_{\operatorname*{Row}T}\nabla
_{\operatorname*{Col}T}^{-}\right)  =S\left(  0\right)  =0
\]
(since $S$ is $\mathbf{k}$-linear). Hence,%
\begin{align*}
0  &  =S\left(  \sum_{x\in L}x\nabla_{\operatorname*{Row}T}\nabla
_{\operatorname*{Col}T}^{-}\right) \\
&  =\sum_{x\in L}\underbrace{S\left(  \nabla_{\operatorname*{Col}T}%
^{-}\right)  }_{\substack{=\nabla_{\operatorname*{Col}T}^{-}\\\text{(by
Proposition \ref{prop.symmetrizers.antipode})}}}\ \ \underbrace{S\left(
\nabla_{\operatorname*{Row}T}\right)  }_{\substack{=\nabla
_{\operatorname*{Row}T}\\\text{(by Proposition
\ref{prop.symmetrizers.antipode})}}}\ \ \underbrace{S\left(  x\right)
}_{\substack{=x^{-1}\\\text{(by the definition of }S\text{,}\\\text{since
}x\in L\subseteq S_{n,Z}\subseteq S_{n}\text{)}}}\\
&  \ \ \ \ \ \ \ \ \ \ \ \ \ \ \ \ \ \ \ \ \left(  \text{since }S\text{ is a
}\mathbf{k}\text{-algebra anti-morphism}\right) \\
&  =\sum_{x\in L}\underbrace{\nabla_{\operatorname*{Col}T}^{-}\nabla
_{\operatorname*{Row}T}}_{\substack{=\mathbf{E}_{T}\\\text{(by the definition
of }\mathbf{E}_{T}\text{)}}}x^{-1}=\sum_{x\in L}\mathbf{E}_{T}x^{-1}%
=\mathbf{E}_{T}\sum_{x\in L}x^{-1}.
\end{align*}
This proves Corollary \ref{cor.hlf-young.grel3rg}.
\end{proof}

\subsubsection{\label{subsec.specht.cogarnir.corr}Proof of the coGarnir
equations: Correctness}

We are now ready to prove Theorem \ref{thm.cogarnir.cogarnir-sats} (and, with
it, the \textquotedblleft$\subseteq$\textquotedblright\ part of Theorem
\ref{thm.cogarnir.cogarnir-skew}):

\begin{proof}
[Proof of Theorem \ref{thm.cogarnir.cogarnir-sats}.]We have $\mathbf{a}%
\in\mathcal{S}^{D}\subseteq\mathcal{M}^{D}=\mathbf{k}^{\left(  \overline
{\operatorname{Tab}}\left(  D\right)  \right)  }$. Thus, $\mathbf{a}$ is a
family $\left(  a_{\overline{P}}\right)  _{\overline{P}\in\overline
{\operatorname{Tab}}\left(  D\right)  }$ of scalars $a_{\overline{P}}%
\in\mathbf{k}$. Consider these scalars $a_{\overline{P}}$. Thus,
$\mathbf{a}=\left(  a_{\overline{P}}\right)  _{\overline{P}\in\overline
{\operatorname{Tab}}\left(  D\right)  }$.

We must show that $\mathbf{a}$ satisfies all coGarnir equations. In other
words, we must show that $\mathbf{a}$ satisfies the $\left(  T,L\right)
$-coGarnir equation for each $n$-tableau $T$ of shape $D$ and each row-Garnir
set $L$ for $T$ (because this is what \textquotedblleft$\mathbf{a}$ satisfies
all coGarnir equations\textquotedblright\ means, according to Definition
\ref{def.cogarnir.cogarnirall}).

So let $T$ be an $n$-tableau of shape $D$, and let $L$ be a row-Garnir set for
$T$. Then, we must show that $\mathbf{a}$ satisfies the $\left(  T,L\right)
$-coGarnir equation. In other words, we must prove that%
\begin{equation}
\sum_{x\in L}a_{\overline{xT}}=0 \label{pf.thm.cogarnir.cogarnir-sats.goal}%
\end{equation}
(because this is what \textquotedblleft$\mathbf{a}$ satisfies the $\left(
T,L\right)  $-coGarnir equation\textquotedblright\ means, according to
Definition \ref{def.cogarnir.cogarnirTL} \textbf{(b)}, since $\mathbf{a}%
=\left(  a_{\overline{P}}\right)  _{\overline{P}\in\overline
{\operatorname{Tab}}\left(  D\right)  }$).

Recall that $L$ is a row-Garnir set for $T$. According to Definition
\ref{def.cogarnir.cogarnirTL} \textbf{(a)}, this means that there exist two
distinct integers $j$ and $k$ as well as a subset $Z$ of $\operatorname*{Row}%
\left(  j,T\right)  \cup\operatorname*{Row}\left(  k,T\right)  $ such that
$\left\vert Z\right\vert >\left\vert D\left\lfloor j\right\rfloor \cup
D\left\lfloor k\right\rfloor \right\vert $ and such that $L$ is a left
transversal of $S_{n,Z}\cap\mathcal{R}\left(  T\right)  $ in $S_{n,Z}$.
Consider these $j$, $k$ and $Z$.

Consider the left $\mathbf{k}\left[  S_{n}\right]  $-module isomorphism
\[
\alpha:\mathbf{k}\left[  S_{n}\right]  \cdot\nabla_{\operatorname*{Row}%
T}\rightarrow\mathcal{M}^{D}%
\]
from Theorem \ref{thm.spechtmod.leftideal} \textbf{(a)}. From Proposition
\ref{prop.spechtmod.leftideal.al-1} \textbf{(c)}, we obtain%
\begin{equation}
\alpha^{-1}\left(  \sum_{\substack{\overline{S}\text{ is an }n\text{-tabloid}%
\\\text{of shape }D}}a_{\overline{S}}\overline{S}\right)  =\sum_{u\in S_{n}%
}a_{u\overline{T}}u. \label{pf.thm.cogarnir.cogarnir-sats.al-1}%
\end{equation}
We can rewrite the summation sign \textquotedblleft$\sum_{\substack{\overline
{S}\text{ is an }n\text{-tabloid}\\\text{of shape }D}}$\textquotedblright\ as
\textquotedblleft$\sum_{\overline{S}\in\overline{\operatorname{Tab}}\left(
D\right)  }$\textquotedblright\ (since we have $\overline{\operatorname{Tab}%
}\left(  D\right)  =\left\{  n\text{-tabloids of shape }D\right\}  $). Hence,
the equality (\ref{pf.thm.cogarnir.cogarnir-sats.al-1}) rewrites as%
\begin{equation}
\alpha^{-1}\left(  \sum_{\overline{S}\in\overline{\operatorname{Tab}}\left(
D\right)  }a_{\overline{S}}\overline{S}\right)  =\sum_{u\in S_{n}%
}a_{u\overline{T}}u. \label{pf.thm.cogarnir.cogarnir-sats.al-1b}%
\end{equation}

But%
\begin{align}
\mathbf{a}  &  =\left(  a_{\overline{P}}\right)  _{\overline{P}\in
\overline{\operatorname{Tab}}\left(  D\right)  }=\left(  a_{\overline{S}%
}\right)  _{\overline{S}\in\overline{\operatorname{Tab}}\left(  D\right)
}\nonumber\\
&  =\sum_{\overline{S}\in\overline{\operatorname{Tab}}\left(  D\right)
}a_{\overline{S}}\overline{S} \label{pf.thm.cogarnir.cogarnir-sats.2}%
\end{align}
(since $\mathcal{M}^{D}$ is the free $\mathbf{k}$-module $\mathbf{k}^{\left(
\overline{\operatorname{Tab}}\left(  D\right)  \right)  }$, and thus the
family $\left(  a_{\overline{S}}\right)  _{\overline{S}\in\overline
{\operatorname{Tab}}\left(  D\right)  }$ can be rewritten as the $\mathbf{k}%
$-linear combination $\sum_{\overline{S}\in\overline{\operatorname{Tab}%
}\left(  D\right)  }a_{\overline{S}}\overline{S}$ of its standard basis
vectors $\overline{S}$). In view of this, we can rewrite
(\ref{pf.thm.cogarnir.cogarnir-sats.al-1b}) as
\begin{equation}
\alpha^{-1}\left(  \mathbf{a}\right)  =\sum_{u\in S_{n}}a_{u\overline{T}}u.
\label{pf.thm.cogarnir.cogarnir-sats.3}%
\end{equation}

Now, define the Young symmetrizer $\mathbf{E}_{T}$ as in Definition
\ref{def.specht.ET.defs} \textbf{(b)}. Theorem \ref{thm.spechtmod.leftideal}
\textbf{(b)} shows that the isomorphism $\alpha$ sends $\mathbf{k}\left[
S_{n}\right]  \cdot\nabla_{\operatorname*{Col}T}^{-}\nabla
_{\operatorname*{Row}T}$ to $\mathcal{S}^{D}$. Hence, its inverse $\alpha
^{-1}$ sends $\mathcal{S}^{D}$ back to $\mathbf{k}\left[  S_{n}\right]
\cdot\nabla_{\operatorname*{Col}T}^{-}\nabla_{\operatorname*{Row}T}$. In other
words,
\[
\alpha^{-1}\left(  \mathcal{S}^{D}\right)  =\mathbf{k}\left[  S_{n}\right]
\cdot\underbrace{\nabla_{\operatorname*{Col}T}^{-}\nabla_{\operatorname*{Row}%
T}}_{\substack{=\mathbf{E}_{T}\\\text{(by the definition of }\mathbf{E}%
_{T}\text{)}}}=\mathbf{k}\left[  S_{n}\right]  \cdot\mathbf{E}_{T}.
\]
Hence, from $\mathbf{a}\in\mathcal{S}^{D}$, we obtain $\alpha^{-1}\left(
\mathbf{a}\right)  \in\alpha^{-1}\left(  \mathcal{S}^{D}\right)
=\mathbf{k}\left[  S_{n}\right]  \cdot\mathbf{E}_{T}$. In other words,
$\alpha^{-1}\left(  \mathbf{a}\right)  =\mathbf{bE}_{T}$ for some
$\mathbf{b}\in\mathbf{k}\left[  S_{n}\right]  $. Consider this $\mathbf{b}$.
Thus,%
\[
\underbrace{\alpha^{-1}\left(  \mathbf{a}\right)  }_{=\mathbf{bE}_{T}}%
\cdot\sum_{x\in L}x^{-1}=\mathbf{b}\underbrace{\mathbf{E}_{T}\sum_{x\in
L}x^{-1}}_{\substack{=0\\\text{(by Corollary \ref{cor.hlf-young.grel3rg})}%
}}=0.
\]
Comparing this with%
\begin{align*}
&  \underbrace{\alpha^{-1}\left(  \mathbf{a}\right)  }_{\substack{=\sum_{u\in
S_{n}}a_{u\overline{T}}u\\\text{(by (\ref{pf.thm.cogarnir.cogarnir-sats.3}))}%
}}\cdot\sum_{x\in L}x^{-1}\\
&  =\sum_{u\in S_{n}}a_{u\overline{T}}u\cdot\sum_{x\in L}x^{-1}=\sum_{x\in
L}\ \ \sum_{u\in S_{n}}a_{u\overline{T}}ux^{-1}\\
&  =\underbrace{\sum_{x\in L}\ \ \sum_{w\in S_{n}}}_{=\sum_{w\in S_{n}%
}\ \ \sum_{x\in L}}a_{wx\overline{T}}\underbrace{wxx^{-1}}_{=w}%
\ \ \ \ \ \ \ \ \ \ \left(
\begin{array}
[c]{c}%
\text{here, we substituted }wx\text{ for }u\\
\text{in the inner sum, since}\\
\text{the map }S_{n}\rightarrow S_{n},\ w\mapsto wx\\
\text{(for a fixed }x\in L\text{) is a bijection}%
\end{array}
\right) \\
&  =\sum_{w\in S_{n}}\ \ \sum_{x\in L}a_{wx\overline{T}}w=\sum_{w\in S_{n}%
}\left(  \sum_{x\in L}a_{wx\overline{T}}\right)  w,
\end{align*}
we obtain%
\[
\sum_{w\in S_{n}}\left(  \sum_{x\in L}a_{wx\overline{T}}\right)  w=0.
\]
Since the standard basis $\left(  w\right)  _{w\in S_{n}}$ of $\mathbf{k}%
\left[  S_{n}\right]  $ is $\mathbf{k}$-linearly independent (because it is a
basis), this equality entails that all the coefficients $\sum_{x\in
L}a_{wx\overline{T}}$ must be $0$. In other words,%
\[
\sum_{x\in L}a_{wx\overline{T}}=0\ \ \ \ \ \ \ \ \ \ \text{for each }w\in
S_{n}.
\]
Applying this to $w=1$ (where $1$ means the identity permutation
$\operatorname*{id}\in S_{n}$), we obtain
\begin{equation}
\sum_{x\in L}a_{1x\overline{T}}=0. \label{pf.thm.cogarnir.cogarnir-sats.at}%
\end{equation}
But each $x\in L$ satisfies $1x\overline{T}=x\overline{T}=\overline{xT}$.
Hence, we can rewrite (\ref{pf.thm.cogarnir.cogarnir-sats.at}) as $\sum_{x\in
L}a_{\overline{xT}}=0$. This proves (\ref{pf.thm.cogarnir.cogarnir-sats.goal}%
). As we explained, our proof of Theorem \ref{thm.cogarnir.cogarnir-sats} is
thus complete.
\end{proof}

\subsubsection{\label{subsec.specht.cogarnir.suff}Proof of the coGarnir
equations: Sufficiency}

We now approach the proof of Theorem \ref{thm.cogarnir.cogarnir-skew}. First,
we define two subsets of $\mathcal{M}^{D}$:

\begin{definition}
\label{def.cogarnir.C-and-G}Let $D$ be a diagram with $\left\vert D\right\vert
=n$. Let $\operatorname*{Tab}\left(  D\right)  $ denote the set $\left\{
n\text{-tableaux of shape }D\right\}  $. Define two subsets $\mathcal{C}^{D}$
and $\mathcal{G}^{D}$ of $\mathcal{M}^{D}$ by%
\[
\mathcal{C}^{D}:=\operatorname*{span}\nolimits_{\mathbf{k}}\left\{
\overline{T}\ \mid\ T\in\operatorname*{Tab}\left(  D\right)  \text{ is
row-standard but not standard}\right\}
\]
and%
\[
\mathcal{G}^{D}:=\left\{  \mathbf{a}\in\mathcal{M}^{D}\ \mid\ \mathbf{a}\text{
satisfies all coGarnir equations}\right\}  .
\]

\end{definition}

In a moment, we will see that both $\mathcal{C}^{D}$ and $\mathcal{G}^{D}$ are
$\mathbf{k}$-submodules of $\mathcal{M}^{D}$. In general, $\mathcal{C}^{D}$ is
not a left $\mathbf{k}\left[  S_{n}\right]  $-submodule, however (i.e., not a
subrepresentation). As for $\mathcal{G}^{D}$, we will eventually see that
$\mathcal{G}^{D}=\mathcal{S}^{D}$ when $D$ is a skew Young diagram; thus, in
this case, $\mathcal{G}^{D}$ is a left $\mathbf{k}\left[  S_{n}\right]
$-submodule of $\mathcal{M}^{D}$.

We begin with the easy parts:

\begin{proposition}
\label{prop.cogarnir.C-and-G-sub}Let $D$ be a diagram with $\left\vert
D\right\vert =n$. Then, both $\mathcal{C}^{D}$ and $\mathcal{G}^{D}$ are
$\mathbf{k}$-submodules of $\mathcal{M}^{D}$.
\end{proposition}

\begin{proof}
Clearly, the set $\mathcal{C}^{D}$ is a $\mathbf{k}$-submodule of
$\mathcal{M}^{D}$ (since it is defined as a $\mathbf{k}$-linear span). It thus
remains to prove that $\mathcal{G}^{D}$ is a $\mathbf{k}$-submodule of
$\mathcal{M}^{D}$.

Recall that $\mathcal{M}^{D}$ is the free $\mathbf{k}$-module $\mathbf{k}%
^{\left(  \overline{\operatorname{Tab}}\left(  D\right)  \right)  }$, so that
each vector $\mathbf{a}\in\mathcal{M}^{D}$ can be written as a family $\left(
a_{\overline{P}}\right)  _{\overline{P}\in\overline{\operatorname{Tab}}\left(
D\right)  }$, and the entries $a_{\overline{P}}$ of this family are the
coordinates of $\mathbf{a}$ (with respect to the standard basis).

By its definition, $\mathcal{G}^{D}$ is the set of all vectors $\mathbf{a}%
\in\mathcal{M}^{D}$ that satisfy all coGarnir equations. However, for any
given vector $\mathbf{a}=\left(  a_{\overline{P}}\right)  _{\overline{P}%
\in\overline{\operatorname{Tab}}\left(  D\right)  }\in\mathcal{M}^{D}$, we
have the following chain of equivalences:%
\begin{align*}
&  \ \left(  \mathbf{a}\text{ satisfies all coGarnir equations}\right) \\
&  \Longleftrightarrow\ \left(  \mathbf{a}\text{ satisfies the }\left(
T,L\right)  \text{-coGarnir equation for each }n\text{-tableau }T\text{ of
shape }D\right. \\
&  \ \ \ \ \ \ \ \ \ \ \ \ \ \ \ \left.  \text{and each row-Garnir set
}L\text{ for }T\right)  \ \ \ \ \ \ \ \ \ \ \left(  \text{by Definition
\ref{def.cogarnir.cogarnirall}}\right) \\
&  \Longleftrightarrow\ \left(  \text{the equation
(\ref{eq.def.cogarnir.cogarnirTL.b.eq}) holds for each }n\text{-tableau
}T\text{ of shape }D\right. \\
&  \ \ \ \ \ \ \ \ \ \ \ \ \ \ \ \left.  \text{and each row-Garnir set
}L\text{ for }T\right)  \ \ \ \ \ \ \ \ \ \ \left(  \text{by Definition
\ref{def.cogarnir.cogarnirTL} \textbf{(b)}}\right)  .
\end{align*}
Thus, the statement \textquotedblleft$\mathbf{a}$ satisfies all coGarnir
equations\textquotedblright\ (for a given vector $\mathbf{a}\in\mathcal{M}%
^{D}$) is tantamount to a set of linear equations on the coordinates of
$\mathbf{a}$ (namely, the equations (\ref{eq.def.cogarnir.cogarnirTL.b.eq})
for all $n$-tableaux $T$ of shape $D$ and all row-Garnir sets $L$ for $T$).
Therefore, $\mathcal{G}^{D}$ is the set of all vectors $\mathbf{a}%
\in\mathcal{M}^{D}$ that satisfy these linear equations (since $\mathcal{G}%
^{D}$ is the set of all vectors $\mathbf{a}\in\mathcal{M}^{D}$ that satisfy
all coGarnir equations).

Thus we have shown that $\mathcal{G}^{D}$ is the set of solutions of a system
of linear equations. Hence, $\mathcal{G}^{D}$ is a $\mathbf{k}$-submodule of
$\mathcal{M}^{D}$ (since the set of solutions of a system of linear equations
is always a $\mathbf{k}$-submodule). This completes the proof of Proposition
\ref{prop.cogarnir.C-and-G-sub}.
\end{proof}

\begin{lemma}
\label{lem.cogarnir.SsG}Let $D$ be a diagram with $\left\vert D\right\vert
=n$. Then, $\mathcal{S}^{D}\subseteq\mathcal{G}^{D}$.
\end{lemma}

\begin{proof}
Each element $\mathbf{b}\in\mathcal{S}^{D}$ satisfies all coGarnir equations
(by Theorem \ref{thm.cogarnir.cogarnir-sats}, applied to $\mathbf{a}%
=\mathbf{b}$), and thus satisfies
\[
\mathbf{b}\in\left\{  \mathbf{a}\in\mathcal{M}^{D}\ \mid\ \mathbf{a}\text{
satisfies all coGarnir equations}\right\}  =\mathcal{G}^{D}%
\]
(by the definition of $\mathcal{G}^{D}$). In other words, $\mathcal{S}%
^{D}\subseteq\mathcal{G}^{D}$. This proves Lemma \ref{lem.cogarnir.SsG}.
\end{proof}

\begin{lemma}
\label{lem.cogarnir.sum}Let $D$ be a diagram with $\left\vert D\right\vert
=n$. Then, $\mathcal{M}^{D}=\mathcal{C}^{D}+\mathcal{S}^{D}$.
\end{lemma}

\begin{proof}
Note that $\mathcal{C}^{D}+\mathcal{S}^{D}$ is a $\mathbf{k}$-submodule of
$\mathcal{M}^{D}$ (since $\mathcal{C}^{D}$ and $\mathcal{S}^{D}$ are
$\mathbf{k}$-submodules), and thus is closed under $\mathbf{k}$-linear
combinations. In other words, $\operatorname*{span}\nolimits_{\mathbf{k}%
}\left(  \mathcal{C}^{D}+\mathcal{S}^{D}\right)  \subseteq\mathcal{C}%
^{D}+\mathcal{S}^{D}$.

We shall now show that%
\begin{equation}
\overline{T}\in\mathcal{C}^{D}+\mathcal{S}^{D}\ \ \ \ \ \ \ \ \ \ \text{for
each }\overline{T}\in\overline{\operatorname{Tab}}\left(  D\right)  .
\label{pf.lem.cogarnir.sum.1}%
\end{equation}

\begin{proof}
[Proof of (\ref{pf.lem.cogarnir.sum.1}):]Recall that the set $\overline
{\operatorname{Tab}}\left(  D\right)  =\left\{  n\text{-tabloids of shape
}D\right\}  $ is equipped with a total order -- namely, the Young last letter
order (which we defined in Proposition \ref{prop.spechtmod.row.order}). We
shall prove (\ref{pf.lem.cogarnir.sum.1}) by strong induction on $\overline
{T}$ with respect to this order:

\textit{Induction step:} Fix some $\overline{P}\in\overline{\operatorname{Tab}%
}\left(  D\right)  $. Assume (as the induction hypothesis) that
(\ref{pf.lem.cogarnir.sum.1}) holds for each $\overline{T}\in\overline
{\operatorname{Tab}}\left(  D\right)  $ satisfying $\overline{T}<\overline{P}$
(where $<$ is the smaller relation of the Young last letter order, i.e., the
relation introduced in Definition \ref{def.tabloid.llo}). Our goal is to show
that (\ref{pf.lem.cogarnir.sum.1}) holds for $\overline{T}=\overline{P}$. In
other words, our goal is to prove that $\overline{P}\in\mathcal{C}%
^{D}+\mathcal{S}^{D}$.

Clearly, $P$ is an $n$-tableau of shape $D$. In other words, $P\in
\operatorname*{Tab}\left(  D\right)  $. By Proposition
\ref{prop.tabloid.row-st}, we can WLOG assume that $P$ is
row-standard\footnote{Indeed, Proposition \ref{prop.tabloid.row-st} shows that
there is a bijection%
\begin{align*}
&  \text{from }\left\{  \text{row-standard }n\text{-tableaux of shape
}D\right\} \\
&  \text{to }\left\{  n\text{-tabloids of shape }D\right\}
\end{align*}
that sends each row-standard $n$-tableau $T$ to its tabloid $\overline{T}$.
Hence, the preimage of the $n$-tabloid $\overline{P}\in\overline
{\operatorname{Tab}}\left(  D\right)  =\left\{  n\text{-tabloids of shape
}D\right\}  $ under this bijection is a row-standard $n$-tableau $T$ of shape
$D$ satisfying $\overline{T}=\overline{P}$ (since it is a preimage of
$\overline{P}$). By replacing $P$ with this row-standard $n$-tableau $T$, we
can thus make $P$ row-standard without changing the $n$-tabloid $\overline{P}%
$.}. Assume this. We are in one of the following two cases:

\textit{Case 1:} The tableau $P$ is not standard.

\textit{Case 2:} The tableau $P$ is standard.

Let us first consider Case 1. In this case, the tableau $P$ is not standard.
Thus, $P\in\operatorname*{Tab}\left(  D\right)  $ is row-standard but not
standard. Hence, the standard basis vector $\overline{P}$ of $\mathcal{M}^{D}$
can be written as $\overline{T}$ for some $n$-tableau $T\in\operatorname*{Tab}%
\left(  D\right)  $ that is row-standard but not standard (namely, for $T=P$).
In other words,%
\begin{align*}
\overline{P}  &  \in\left\{  \overline{T}\ \mid\ T\in\operatorname*{Tab}%
\left(  D\right)  \text{ is row-standard but not standard}\right\} \\
&  \subseteq\operatorname*{span}\nolimits_{\mathbf{k}}\left\{  \overline
{T}\ \mid\ T\in\operatorname*{Tab}\left(  D\right)  \text{ is row-standard but
not standard}\right\} \\
&  =\mathcal{C}^{D}\ \ \ \ \ \ \ \ \ \ \left(  \text{by the definition of
}\mathcal{C}^{D}\right) \\
&  \subseteq\mathcal{C}^{D}+\mathcal{S}^{D}.
\end{align*}
Thus, $\overline{P}\in\mathcal{C}^{D}+\mathcal{S}^{D}$ is proved in Case 1.

Let us now consider Case 2. In this case, the tableau $P$ is standard. Hence,
Lemma \ref{lem.spechtmod.lead-term} \textbf{(b)} (applied to $T=P$) shows
that
\[
\mathbf{e}_{P}=\overline{P}+\left(  \text{a linear combination of
}n\text{-tabloids }\overline{S}\text{ with }\overline{S}<\overline{P}\right)
.
\]
Hence,%
\begin{align*}
\mathbf{e}_{P}-\overline{P}  &  =\left(  \text{a linear combination of
}n\text{-tabloids }\overline{S}\text{ with }\overline{S}<\overline{P}\right)
\\
&  =\left(  \text{a linear combination of }\overline{S}\text{ for }%
\overline{S}\in\overline{\operatorname*{Tab}}\left(  D\right)  \text{
satisfying }\overline{S}<\overline{P}\right)  .
\end{align*}
In other words,%
\begin{align*}
\mathbf{e}_{P}-\overline{P}  &  \in\operatorname*{span}\nolimits_{\mathbf{k}%
}\left\{  \overline{S}\ \mid\ \overline{S}\in\overline{\operatorname*{Tab}%
}\left(  D\right)  \text{ satisfying }\overline{S}<\overline{P}\right\} \\
&  =\operatorname*{span}\nolimits_{\mathbf{k}}\left\{  \overline{T}%
\ \mid\ \overline{T}\in\overline{\operatorname*{Tab}}\left(  D\right)  \text{
satisfying }\overline{T}<\overline{P}\right\}
\end{align*}
(here, we have renamed the index $\overline{S}$ as $\overline{T}$). But our
induction hypothesis says that (\ref{pf.lem.cogarnir.sum.1}) holds for each
$\overline{T}\in\overline{\operatorname{Tab}}\left(  D\right)  $ satisfying
$\overline{T}<\overline{P}$. In other words, each $\overline{T}\in
\overline{\operatorname{Tab}}\left(  D\right)  $ satisfying $\overline
{T}<\overline{P}$ satisfies $\overline{T}\in\mathcal{C}^{D}+\mathcal{S}^{D}$.
In other words,%
\[
\left\{  \overline{T}\ \mid\ \overline{T}\in\overline{\operatorname*{Tab}%
}\left(  D\right)  \text{ satisfying }\overline{T}<\overline{P}\right\}
\subseteq\mathcal{C}^{D}+\mathcal{S}^{D}.
\]
Hence, our above computation becomes%
\begin{align*}
\mathbf{e}_{P}-\overline{P}  &  \in\operatorname*{span}\nolimits_{\mathbf{k}%
}\underbrace{\left\{  \overline{T}\ \mid\ \overline{T}\in\overline
{\operatorname*{Tab}}\left(  D\right)  \text{ satisfying }\overline
{T}<\overline{P}\right\}  }_{\subseteq\mathcal{C}^{D}+\mathcal{S}^{D}}\\
&  \subseteq\operatorname*{span}\nolimits_{\mathbf{k}}\left(  \mathcal{C}%
^{D}+\mathcal{S}^{D}\right)  \subseteq\mathcal{C}^{D}+\mathcal{S}^{D}.
\end{align*}
In other words, $\mathbf{e}_{P}-\overline{P}=\mathbf{c}+\mathbf{s}$ for some
$\mathbf{c}\in\mathcal{C}^{D}$ and some $\mathbf{s}\in\mathcal{S}^{D}$.
Consider these $\mathbf{c}$ and $\mathbf{s}$. Observe that $\mathbf{e}_{P}$ is
a polytabloid and thus belongs to $\mathcal{S}^{D}$ (since $\mathcal{S}^{D}$
is defined as the span of all polytabloids $\mathbf{e}_{T}$ with
$T\in\operatorname*{Tab}\left(  D\right)  $). That is, $\mathbf{e}_{P}%
\in\mathcal{S}^{D}$. Now, solving the equation $\mathbf{e}_{P}-\overline
{P}=\mathbf{c}+\mathbf{s}$ for $\overline{P}$, we find
\begin{align*}
\overline{P}  &  =\mathbf{e}_{P}-\left(  \mathbf{c}+\mathbf{s}\right)
=\underbrace{\mathbf{e}_{P}}_{\in\mathcal{S}^{D}}-\underbrace{\mathbf{s}}%
_{\in\mathcal{S}^{D}}-\underbrace{\mathbf{c}}_{\in\mathcal{C}^{D}}%
\in\underbrace{\mathcal{S}^{D}-\mathcal{S}^{D}}_{\substack{\subseteq
\mathcal{S}^{D}\\\text{(since }\mathcal{S}^{D}\text{ is a }\mathbf{k}%
\text{-module)}}}-\,\mathcal{C}^{D}\\
&  \subseteq\mathcal{S}^{D}-\mathcal{C}^{D}=\underbrace{-\mathcal{C}^{D}%
}_{\substack{=\mathcal{C}^{D}\\\text{(since }\mathcal{C}^{D}\text{ is a
}\mathbf{k}\text{-module)}}}+\,\mathcal{S}^{D}=\mathcal{C}^{D}+\mathcal{S}%
^{D}.
\end{align*}
Thus, $\overline{P}\in\mathcal{C}^{D}+\mathcal{S}^{D}$ is proved in Case 2.

We have now proved $\overline{P}\in\mathcal{C}^{D}+\mathcal{S}^{D}$ in both
Cases 1 and 2. Hence, $\overline{P}\in\mathcal{C}^{D}+\mathcal{S}^{D}$ always
holds. This completes the induction step. Thus, (\ref{pf.lem.cogarnir.sum.1})
is proved.
\end{proof}

Now, recall that the $\mathbf{k}$-module $\mathcal{M}^{D}$ has a basis
$\left(  \overline{T}\right)  _{\overline{T}\in\overline{\operatorname*{Tab}%
}\left(  D\right)  }$. Hence,%
\[
\mathcal{M}^{D}=\operatorname*{span}\nolimits_{\mathbf{k}}\underbrace{\left\{
\overline{T}\ \mid\ \overline{T}\in\overline{\operatorname*{Tab}}\left(
D\right)  \right\}  }_{\substack{\subseteq\mathcal{C}^{D}+\mathcal{S}%
^{D}\\\text{(by (\ref{pf.lem.cogarnir.sum.1}))}}}\subseteq\operatorname*{span}%
\nolimits_{\mathbf{k}}\left(  \mathcal{C}^{D}+\mathcal{S}^{D}\right)
\subseteq\mathcal{C}^{D}+\mathcal{S}^{D}.
\]
Combining this with the obvious inclusion $\mathcal{C}^{D}+\mathcal{S}%
^{D}\subseteq\mathcal{M}^{D}$, we obtain $\mathcal{M}^{D}=\mathcal{C}%
^{D}+\mathcal{S}^{D}$. This proves Lemma \ref{lem.cogarnir.sum}.
\end{proof}

We now come to some combinatorial lemmas. Our next lemma is an analogue of
Lemma \ref{lem.garnir.skew}:

\begin{lemma}
\label{lem.cogarnir.skew}Let $D$ be the skew Young diagram $Y\left(
\lambda/\mu\right)  $ of a skew partition $\lambda/\mu$. Let $T$ be an
$n$-tableau of shape $D$. Let $j$ and $k$ be two integers with $j<k$. Let
$i\in\mathbb{Z}$ be such that $\left(  j,i\right)  \in D$ and $\left(
k,i\right)  \in D$. Define two sets%
\begin{align*}
X  &  =\left\{  T\left(  j,i^{\prime}\right)  \ \mid\ i^{\prime}\in
\mathbb{Z}\text{ with }i^{\prime}\geq i\text{ and }\left(  j,i^{\prime
}\right)  \in D\right\}  \ \ \ \ \ \ \ \ \ \ \text{and}\\
Y  &  =\left\{  T\left(  k,i^{\prime}\right)  \ \mid\ i^{\prime}\in
\mathbb{Z}\text{ with }i^{\prime}\leq i\text{ and }\left(  k,i^{\prime
}\right)  \in D\right\}  .
\end{align*}
(That is, $X$ is the set of all entries in the $j$-th row of $T$ at the cell
$\left(  j,i\right)  $ and further east, whereas $Y$ is the set of all entries
in the $k$-th row of $T$ at the cell $\left(  k,i\right)  $ and further west.)
Then: \medskip

\textbf{(a)} We have $X\subseteq\operatorname*{Row}\left(  j,T\right)  $ and
$Y\subseteq\operatorname*{Row}\left(  k,T\right)  $ and $\left\vert X\cup
Y\right\vert >\left\vert D\left\lfloor j\right\rfloor \cup D\left\lfloor
k\right\rfloor \right\vert $. \medskip

\textbf{(b)} If the tableau $T$ is row-standard and satisfies $T\left(
j,i\right)  >T\left(  k,i\right)  $, then any element of $X$ is larger than
any element of $Y$.
\end{lemma}

\begin{proof}
This is entirely analogous to Lemma \ref{lem.garnir.skew}; the only difference
is that the horizontal and vertical directions are interchanged (i.e., rows
become columns, north becomes west, each cell $\left(  i,j\right)  $ becomes
$\left(  j,i\right)  $, etc.).
\end{proof}

Our next lemma is an analogue of Lemma \ref{lem.spechtmod.garnir-larger}:

\begin{lemma}
\label{lem.cogarnir.garnir-larger}Let $T$ be an $n$-tableau (of any shape).
Let $j$ and $k$ be two integers such that $j<k$. Let $X\subseteq
\operatorname*{Row}\left(  j,T\right)  $ and $Y\subseteq\operatorname*{Row}%
\left(  k,T\right)  $ be two subsets such that each element of $X$ is larger
than each element of $Y$. Let $w\in S_{n,X\cup Y}\setminus\mathcal{R}\left(
T\right)  $. Then, $\overline{wT}>\overline{T}$ (with respect to the Young
last letter order).
\end{lemma}

\begin{proof}
This is entirely analogous to Lemma \ref{lem.spechtmod.garnir-larger}; the
only difference is that the horizontal and vertical directions are
interchanged (i.e., rows become columns, $\mathcal{C}\left(  T\right)  $
becomes $\mathcal{R}\left(  T\right)  $, each cell $\left(  i,j\right)  $
becomes $\left(  j,i\right)  $, etc.). (Of course, the same change needs to be
made in Lemma \ref{lem.spechtmod.garnir-larger.gen} and its proof.)
\end{proof}

\begin{lemma}
\label{lem.cogarnir.triv0}Let $D$ be a diagram with $\left\vert D\right\vert
=n$. Let $\mathbf{a}\in\mathcal{M}^{D}$. Write $\mathbf{a}$ in the form
$\mathbf{a}=\left(  a_{\overline{P}}\right)  _{\overline{P}\in\overline
{\operatorname{Tab}}\left(  D\right)  }$ (this can be done, since each element
of $\mathcal{M}^{D}$ can be written in this form). Assume that $\mathbf{a}%
\in\mathcal{C}^{D}$. Then,%
\[
a_{\overline{Q}}=0\ \ \ \ \ \ \ \ \ \ \text{for all standard tableaux }%
Q\in\operatorname*{Tab}\left(  D\right)  .
\]

\end{lemma}

\begin{proof}
Recall that an element of $\mathcal{M}^{D}=\mathbf{k}^{\left(  \overline
{\operatorname*{Tab}}\left(  D\right)  \right)  }$ is just a family of
scalars, indexed by the $n$-tabloids of shape $D$. These scalars are called
the coordinates of this element. If $\mathbf{b}=\left(  b_{\overline{P}%
}\right)  _{\overline{P}\in\overline{\operatorname{Tab}}\left(  D\right)  }$
is any element of $\mathcal{M}^{D}$, and if $\overline{Q}\in\overline
{\operatorname*{Tab}}\left(  D\right)  $ is any $n$-tabloid of shape $D$, then
$\left[  \overline{Q}\right]  \mathbf{b}$ shall denote the $\overline{Q}%
$-coordinate of $\mathbf{b}$ (that is, the scalar $b_{\overline{Q}}$).
Clearly,%
\begin{equation}
\left[  \overline{Q}\right]  \overline{P}=0\ \ \ \ \ \ \ \ \ \ \text{for any
two distinct }\overline{P},\overline{Q}\in\overline{\operatorname*{Tab}%
}\left(  D\right)  . \label{pf.lem.cogarnir.triv0.1}%
\end{equation}

But $\mathbf{a}\in\mathcal{C}^{D}=\operatorname*{span}\nolimits_{\mathbf{k}%
}\left\{  \overline{T}\ \mid\ T\in\operatorname*{Tab}\left(  D\right)  \text{
is row-standard but not standard}\right\}  $. Hence, we can write $\mathbf{a}$
as a $\mathbf{k}$-linear combination%
\begin{equation}
\mathbf{a}=\sum_{\substack{T\in\operatorname*{Tab}\left(  D\right)  \text{ is
row-standard}\\\text{but not standard}}}c_{\overline{T}}\overline{T}
\label{pf.lem.cogarnir.triv0.2}%
\end{equation}
for some scalars $c_{\overline{T}}\in\mathbf{k}$. Consider these scalars
$c_{\overline{T}}$.

Now, fix a standard tableau $Q\in\operatorname*{Tab}\left(  D\right)  $. Then,
$Q$ is clearly row-standard. Hence, if some $n$-tableau $T\in
\operatorname*{Tab}\left(  D\right)  $ is row-standard but not standard, then
$\overline{T}\neq\overline{Q}$\ \ \ \ \footnote{\textit{Proof.} Proposition
\ref{prop.tabloid.row-st} shows that there is a bijection%
\begin{align*}
&  \text{from }\left\{  \text{row-standard }n\text{-tableaux of shape
}D\right\} \\
&  \text{to }\left\{  n\text{-tabloids of shape }D\right\}
\end{align*}
that sends each row-standard $n$-tableau $T$ to its tabloid $\overline{T}$.
Thus, in particular, this bijection is injective. In other words, if $T_{1}$
and $T_{2}$ are two distinct row-standard $n$-tableaux of shape $D$, then
$\overline{T_{1}}\neq\overline{T_{2}}$.
\par
Now, let $T\in\operatorname*{Tab}\left(  D\right)  $ be an $n$-tableau that is
row-standard but not standard. Then, $T\neq Q$ (since $Q$ is standard but $T$
is not). Therefore, $T$ and $Q$ are two distinct row-standard $n$-tableaux of
shape $D$. Therefore, $\overline{T}\neq\overline{Q}$ (because if $T_{1}$ and
$T_{2}$ are two distinct row-standard $n$-tableaux of shape $D$, then
$\overline{T_{1}}\neq\overline{T_{2}}$).} and therefore
\begin{equation}
\left[  \overline{Q}\right]  \overline{T}=0 \label{pf.lem.cogarnir.triv0.4}%
\end{equation}
(by (\ref{pf.lem.cogarnir.triv0.1}) (applied to $P=T$), since $\overline
{T}\neq\overline{Q}$ shows that $\overline{T}$ and $\overline{Q}$ are distinct).

Now, from (\ref{pf.lem.cogarnir.triv0.2}), we obtain%
\[
\left[  \overline{Q}\right]  \mathbf{a}=\left[  \overline{Q}\right]  \left(
\sum_{\substack{T\in\operatorname*{Tab}\left(  D\right)  \text{ is
row-standard}\\\text{but not standard}}}c_{\overline{T}}\overline{T}\right)
=\sum_{\substack{T\in\operatorname*{Tab}\left(  D\right)  \text{ is
row-standard}\\\text{but not standard}}}c_{\overline{T}}\underbrace{\left[
\overline{Q}\right]  \overline{T}}_{\substack{=0\\\text{(by
(\ref{pf.lem.cogarnir.triv0.4}))}}}=0.
\]

But $\mathbf{a}=\left(  a_{\overline{P}}\right)  _{\overline{P}\in
\overline{\operatorname{Tab}}\left(  D\right)  }$. Hence, the definition of
$\left[  \overline{Q}\right]  \mathbf{a}$ shows that $\left[  \overline
{Q}\right]  \mathbf{a}=a_{\overline{Q}}$. Therefore, $a_{\overline{Q}}=\left[
\overline{Q}\right]  \mathbf{a}=0$. This proves Lemma \ref{lem.cogarnir.triv0}.
\end{proof}

\begin{lemma}
\label{lem.cogarnir.cut0}Let $D$ be a skew Young diagram with $\left\vert
D\right\vert =n$. Then, $\mathcal{C}^{D}\cap\mathcal{G}^{D}=0$.
\end{lemma}

\begin{proof}
Clearly, $\mathcal{C}^{D}\cap\mathcal{G}^{D}$ is a $\mathbf{k}$-submodule of
$\mathcal{M}^{D}$ (since $\mathcal{C}^{D}$ and $\mathcal{G}^{D}$ are
$\mathbf{k}$-submodules of $\mathcal{M}^{D}$, by Proposition
\ref{prop.cogarnir.C-and-G-sub}). Thus, in order to prove that $\mathcal{C}%
^{D}\cap\mathcal{G}^{D}=0$, it suffices to show that each $\mathbf{a}%
\in\mathcal{C}^{D}\cap\mathcal{G}^{D}$ satisfies $\mathbf{a}=0$. This is what
we will now do.

Let $\mathbf{a}\in\mathcal{C}^{D}\cap\mathcal{G}^{D}$. We must show that
$\mathbf{a}=0$.

Write $\mathbf{a}\in\mathcal{C}^{D}\cap\mathcal{G}^{D}\subseteq\mathcal{M}%
^{D}$ in the form $\mathbf{a}=\left(  a_{\overline{P}}\right)  _{\overline
{P}\in\overline{\operatorname{Tab}}\left(  D\right)  }$ (this can be done,
since each element of $\mathcal{M}^{D}$ can be written in this form).

We have $\mathbf{a}\in\mathcal{C}^{D}\cap\mathcal{G}^{D}\subseteq
\mathcal{G}^{D}$. In other words, $\mathbf{a}$ satisfies all coGarnir
equations (by the definition of $\mathcal{G}^{D}$). Moreover, $\mathbf{a}%
\in\mathcal{C}^{D}\cap\mathcal{G}^{D}\subseteq\mathcal{C}^{D}$.

We must show that $\mathbf{a}=0$.

Assume the contrary. Thus, $\mathbf{a}\neq0$. In other words, $\left(
a_{\overline{P}}\right)  _{\overline{P}\in\overline{\operatorname{Tab}}\left(
D\right)  }\neq0$ (since $\mathbf{a}=\left(  a_{\overline{P}}\right)
_{\overline{P}\in\overline{\operatorname{Tab}}\left(  D\right)  }$). In other
words, there exists some $\overline{P}\in\overline{\operatorname{Tab}}\left(
D\right)  $ such that $a_{\overline{P}}\neq0$. Consider the \textbf{largest}
such $\overline{P}$ with respect to the Young last letter order (i.e., to the
total order defined in Proposition \ref{prop.spechtmod.row.order}). Thus,
$a_{\overline{P}}\neq0$, but every $n$-tabloid $\overline{Q}\in\overline
{\operatorname*{Tab}}\left(  D\right)  $ satisfying $\overline{Q}>\overline
{P}$ satisfies%
\begin{equation}
a_{\overline{Q}}=0 \label{pf.lem.cogarnir.cut0.aQ=0}%
\end{equation}
(since $\overline{P}$ is the \textbf{largest} $n$-tabloid such that
$a_{\overline{P}}\neq0$).

Clearly, $P$ is an $n$-tableau of shape $D$. In other words, $P\in
\operatorname*{Tab}\left(  D\right)  $. By Proposition
\ref{prop.tabloid.row-st}, we can WLOG assume that $P$ is
row-standard\footnote{Indeed, Proposition \ref{prop.tabloid.row-st} shows that
there is a bijection%
\begin{align*}
&  \text{from }\left\{  \text{row-standard }n\text{-tableaux of shape
}D\right\} \\
&  \text{to }\left\{  n\text{-tabloids of shape }D\right\}
\end{align*}
that sends each row-standard $n$-tableau $T$ to its tabloid $\overline{T}$.
Hence, the preimage of the $n$-tabloid $\overline{P}\in\overline
{\operatorname{Tab}}\left(  D\right)  =\left\{  n\text{-tabloids of shape
}D\right\}  $ under this bijection is a row-standard $n$-tableau $T$ of shape
$D$ satisfying $\overline{T}=\overline{P}$ (since it is a preimage of
$\overline{P}$). By replacing $P$ with this row-standard $n$-tableau $T$, we
can thus make $P$ row-standard without changing the $n$-tabloid $\overline{P}%
$.}. Assume this.

If the tableau $P$ was standard, then Lemma \ref{lem.cogarnir.triv0} (applied
to $Q=P$) would yield $a_{\overline{P}}=0$ (since $\mathbf{a}\in
\mathcal{C}^{D}$), which would contradict $a_{\overline{P}}\neq0$. Thus, $P$
cannot be standard. Hence, $P$ cannot be column-standard (since $P$ is row-standard).

In other words, there exist two adjacent cells $\left(  j,i\right)  \in D$ and
$\left(  j+1,i\right)  \in D$ such that $P\left(  j,i\right)  \geq P\left(
j+1,i\right)  $\ \ \ \ \footnote{\textit{Proof.} If every $\left(  i,j\right)
\in D$ satisfying $\left(  i+1,j\right)  \in D$ satisfied $P\left(
i,j\right)  <P\left(  i+1,j\right)  $, then $P$ would be column-standard (by
Proposition \ref{prop.tableau.std-loc} \textbf{(b)}, applied to $T=P$), which
would contradict the fact that $P$ is not column-standard. Hence, not every
$\left(  i,j\right)  \in D$ satisfying $\left(  i+1,j\right)  \in D$ satisfies
$P\left(  i,j\right)  <P\left(  i+1,j\right)  $. In other words, there exists
some $\left(  i,j\right)  \in D$ satisfying $\left(  i+1,j\right)  \in D$ but
not $P\left(  i,j\right)  <P\left(  i+1,j\right)  $. In other words, there
exists some $\left(  i,j\right)  \in D$ satisfying $\left(  i+1,j\right)  \in
D$ and $P\left(  i,j\right)  \geq P\left(  i+1,j\right)  $. Renaming the cell
$\left(  i,j\right)  $ as $\left(  j,i\right)  $ in this result, we can
rewrite it as follows: There exists some $\left(  j,i\right)  \in D$
satisfying $\left(  j+1,i\right)  \in D$ and $P\left(  j,i\right)  \geq
P\left(  j+1,i\right)  $. In other words, there exist two adjacent cells
$\left(  j,i\right)  \in D$ and $\left(  j+1,i\right)  \in D$ such that
$P\left(  j,i\right)  \geq P\left(  j+1,i\right)  $.}. Consider these two
cells. Since all entries of $P$ are distinct (because $P$ is an $n$-tableau),
we have $P\left(  j,i\right)  \neq P\left(  j+1,i\right)  $. Hence, from
$P\left(  j,i\right)  \geq P\left(  j+1,i\right)  $, we obtain $P\left(
j,i\right)  >P\left(  j+1,i\right)  $.

Set $k:=j+1$. Thus, $k=j+1>j$, whence $j<k$. Thus, the integers $j$ and $k$
are distinct. Also, from $k=j+1$, we obtain $\left(  k,i\right)  =\left(
j+1,i\right)  \in D$. Finally, $P\left(  j,i\right)  >P\left(  j+1,i\right)
=P\left(  k,i\right)  $ (since $j+1=k$).

Now, define the two sets%
\begin{align*}
X  &  =\left\{  P\left(  j,i^{\prime}\right)  \ \mid\ i^{\prime}\in
\mathbb{Z}\text{ with }i^{\prime}\geq i\text{ and }\left(  j,i^{\prime
}\right)  \in D\right\}  \ \ \ \ \ \ \ \ \ \ \text{and}\\
Y  &  =\left\{  P\left(  k,i^{\prime}\right)  \ \mid\ i^{\prime}\in
\mathbb{Z}\text{ with }i^{\prime}\leq i\text{ and }\left(  k,i^{\prime
}\right)  \in D\right\}  .
\end{align*}
(That is, $X$ is the set of all entries in the $j$-th row of $P$ at the cell
$\left(  j,i\right)  $ and further east, whereas $Y$ is the set of all entries
in the $k$-th row of $P$ at the cell $\left(  k,i\right)  $ and further west.)

Now, Lemma \ref{lem.cogarnir.skew} \textbf{(a)} (applied to $T=P$) yields
$X\subseteq\operatorname*{Row}\left(  j,P\right)  $ and $Y\subseteq
\operatorname*{Row}\left(  k,P\right)  $ and $\left\vert X\cup Y\right\vert
>\left\vert D\left\lfloor j\right\rfloor \cup D\left\lfloor k\right\rfloor
\right\vert $ (since $j<k$). Furthermore, Lemma \ref{lem.cogarnir.skew}
\textbf{(b)} (applied to $T=P$) yields that any element of $X$ is larger than
any element of $Y$ (since $P$ is row-standard and satisfies $P\left(
j,i\right)  >P\left(  k,i\right)  $). Hence, Lemma
\ref{lem.cogarnir.garnir-larger} (applied to $T=P$) shows that
\begin{equation}
\text{each }w\in S_{n,X\cup Y}\setminus\mathcal{R}\left(  P\right)  \text{
satisfies }\overline{wP}>\overline{P} \label{pf.lem.cogarnir.cut0.ineq}%
\end{equation}
(with respect to the Young last letter order).

Now, set $Z:=X\cup Y$. Then, $Z$ is a subset of $\operatorname*{Row}\left(
j,P\right)  \cup\operatorname*{Row}\left(  k,P\right)  $ (since
$Z=\underbrace{X}_{\subseteq\operatorname*{Row}\left(  j,P\right)  }%
\cup\underbrace{Y}_{\subseteq\operatorname*{Row}\left(  k,P\right)  }%
\subseteq\operatorname*{Row}\left(  j,P\right)  \cup\operatorname*{Row}\left(
k,P\right)  $) and has size $\left\vert Z\right\vert =\left\vert X\cup
Y\right\vert >\left\vert D\left\lfloor j\right\rfloor \cup D\left\lfloor
k\right\rfloor \right\vert $.

Moreover, the set $S_{n,Z}$ (defined as in Proposition \ref{prop.intX.basics})
is a subgroup of $S_{n}$ (by Proposition \ref{prop.intX.basics} \textbf{(a)},
applied to $Z$ instead of $X$). Thus, both $\mathcal{R}\left(  P\right)  $ and
$S_{n,Z}$ are subgroups of $S_{n}$. Hence, their intersection $S_{n,Z}%
\cap\mathcal{R}\left(  P\right)  $ is a subgroup of $S_{n,Z}$. Pick any left
transversal $L$ of $S_{n,Z}\cap\mathcal{R}\left(  P\right)  $ in $S_{n,Z}$.
Then, $L$ is a row-Garnir set for $P$ (by Definition
\ref{def.cogarnir.cogarnirTL} \textbf{(a)}, applied to $T=P$). Hence,
$\mathbf{a}$ satisfies the $\left(  P,L\right)  $-coGarnir equation (by
Definition \ref{def.cogarnir.cogarnirall}, since $\mathbf{a}$ satisfies all
coGarnir equations). In other words, we have
\begin{equation}
\sum_{x\in L}a_{\overline{xP}}=0 \label{pf.lem.cogarnir.cut0.sum=0}%
\end{equation}
(by Definition \ref{def.cogarnir.cogarnirTL} \textbf{(b)}).

Recall that $L$ is a left transversal of $S_{n,Z}\cap\mathcal{R}\left(
P\right)  $ in $S_{n,Z}$. Hence, $L\subseteq S_{n,Z}=S_{n,X\cup Y}$ (since
$Z=X\cup Y$). Next, we claim that
\begin{equation}
\text{each }w\in L\setminus\mathcal{R}\left(  P\right)  \text{ satisfies
}\overline{wP}>\overline{P} \label{pf.lem.cogarnir.cut0.ine2}%
\end{equation}
(with respect to the Young last letter order). Indeed, each $w\in
L\setminus\mathcal{R}\left(  P\right)  $ satisfies $w\in\underbrace{L}%
_{\subseteq S_{n,X\cup Y}}\setminus\,\mathcal{R}\left(  P\right)  \subseteq
S_{n,X\cup Y}\setminus\mathcal{R}\left(  P\right)  $ and thus $\overline
{wP}>\overline{P}$ (by (\ref{pf.lem.cogarnir.cut0.ineq})); thus,
(\ref{pf.lem.cogarnir.cut0.ine2}) is proved.

Recall again that $L$ is a left transversal of $S_{n,Z}\cap\mathcal{R}\left(
P\right)  $ in $S_{n,Z}$. Hence, Lemma \ref{lem.G/H.transversal.diff}
\textbf{(a)} (applied to $G=S_{n}$ and $H=S_{n,Z}$ and $K=\mathcal{R}\left(
P\right)  $) shows that there is exactly one $u\in L$ that belongs to
$\mathcal{R}\left(  P\right)  $ (since $L$ is a left transversal of
$S_{n,Z}\cap\mathcal{R}\left(  P\right)  $ in $S_{n,Z}$). Consider this $u$.
Thus, $u\in\mathcal{R}\left(  P\right)  $. Moreover, Lemma
\ref{lem.G/H.transversal.diff} \textbf{(b)} (applied to $G=S_{n}$ and
$H=S_{n,Z}$ and $K=\mathcal{R}\left(  P\right)  $) shows that
\begin{equation}
L\setminus\left\{  u\right\}  =L\setminus\mathcal{R}\left(  P\right)  .
\label{pf.lem.cogarnir.cut0.LL}%
\end{equation}

Proposition \ref{prop.tableau.Sn-act.0} \textbf{(a)} (applied to $w=u$ and
$T=P$) shows that the tableau $u\rightharpoonup P$ is row-equivalent to $P$
(since $u\in\mathcal{R}\left(  P\right)  $). In other words, $\overline
{u\rightharpoonup P}=\overline{P}$ (since row-equivalent $n$-tableaux belong
to the same $n$-tabloid). In other words, $\overline{uP}=\overline{P}$ (since
$uP=u\rightharpoonup P$).

Now, from (\ref{pf.lem.cogarnir.cut0.sum=0}), we obtain%
\begin{align}
0  &  =\sum_{x\in L}a_{\overline{xP}}=\underbrace{a_{\overline{uP}}%
}_{\substack{=a_{\overline{P}}\\\text{(since }\overline{uP}=\overline
{P}\text{)}}}+\underbrace{\sum_{\substack{x\in L;\\x\neq u}}}_{\substack{=\sum
_{x\in L\setminus\left\{  u\right\}  }}}a_{\overline{xP}}%
\ \ \ \ \ \ \ \ \ \ \left(
\begin{array}
[c]{c}%
\text{here, we have split off the}\\
\text{addend for }x=u\text{ from the sum}%
\end{array}
\right) \nonumber\\
&  =a_{\overline{P}}+\sum_{x\in L\setminus\left\{  u\right\}  }a_{\overline
{xP}}. \label{pf.lem.cogarnir.cut0.0=sum}%
\end{align}

Now, let $x\in L\setminus\left\{  u\right\}  $ be arbitrary. Thus, $x\in
L\setminus\left\{  u\right\}  =L\setminus\mathcal{R}\left(  P\right)  $ (by
(\ref{pf.lem.cogarnir.cut0.LL})), so that (\ref{pf.lem.cogarnir.cut0.ineq})
shows that $\overline{xP}>\overline{P}$. Thus,
(\ref{pf.lem.cogarnir.cut0.aQ=0}) (applied to $\overline{Q}=\overline{xP}$)
yields $a_{\overline{xP}}=0$.

Forget that we fixed $x$. We thus have shown that $a_{\overline{xP}}=0$ for
each $x\in L\setminus\left\{  u\right\}  $. Hence, $\sum_{x\in L\setminus
\left\{  u\right\}  }\underbrace{a_{\overline{xP}}}_{=0}=0$. Thus,
(\ref{pf.lem.cogarnir.cut0.0=sum}) becomes%
\[
0=a_{\overline{P}}+\underbrace{\sum_{x\in L\setminus\left\{  u\right\}
}a_{\overline{xP}}}_{=0}=a_{\overline{P}}\neq0.
\]
This is clearly a contradiction. Thus, our assumption was wrong, and
$\mathbf{a}=0$ is proved. This completes our proof of Lemma
\ref{lem.cogarnir.cut0}.
\end{proof}

\begin{theorem}
\label{thm.cogarnir.CSG}Let $D$ be a skew Young diagram with $\left\vert
D\right\vert =n$. Then: \medskip

\textbf{(a)} We have $\mathcal{G}^{D}=\mathcal{S}^{D}$. \medskip

\textbf{(b)} We have $\mathcal{M}^{D}=\mathcal{S}^{D}\oplus\mathcal{C}^{D}$
(an internal direct sum) as $\mathbf{k}$-modules.
\end{theorem}

\begin{proof}
Proposition \ref{prop.cogarnir.C-and-G-sub} shows that $\mathcal{C}^{D}$ and
$\mathcal{G}^{D}$ are $\mathbf{k}$-submodules of $\mathcal{M}^{D}$. \medskip

\textbf{(a)} We have $\mathcal{S}^{D}\subseteq\mathcal{G}^{D}$ by Lemma
\ref{lem.cogarnir.SsG}.

Now, let $\mathbf{a}\in\mathcal{G}^{D}$. Then, $\mathbf{a}\in\mathcal{G}%
^{D}\subseteq\mathcal{M}^{D}=\mathcal{C}^{D}+\mathcal{S}^{D}$ by Lemma
\ref{lem.cogarnir.sum}. Hence, $\mathbf{a}=\mathbf{c}+\mathbf{s}$ for some
$\mathbf{c}\in\mathcal{C}^{D}$ and some $\mathbf{s}\in\mathcal{S}^{D}$.
Consider these $\mathbf{c}$ and $\mathbf{s}$. Now, from $\mathbf{a}%
=\mathbf{c}+\mathbf{s}$, we obtain
\[
\mathbf{c}=\underbrace{\mathbf{a}}_{\in\mathcal{G}^{D}}-\underbrace{\mathbf{s}%
}_{\in\mathcal{S}^{D}\subseteq\mathcal{G}^{D}}\in\mathcal{G}^{D}%
-\mathcal{G}^{D}\subseteq\mathcal{G}^{D}\ \ \ \ \ \ \ \ \ \ \left(
\text{since }\mathcal{G}^{D}\text{ is a }\mathbf{k}\text{-module}\right)  .
\]
Combining this with $\mathbf{c}\in\mathcal{C}^{D}$, we obtain $\mathbf{c}%
\in\mathcal{C}^{D}\cap\mathcal{G}^{D}=0$ by Lemma \ref{lem.cogarnir.cut0}.
Hence, $\mathbf{c}=0$. Thus, $\mathbf{a}=\underbrace{\mathbf{c}}%
_{=0}+\,\mathbf{s}=\mathbf{s}\in\mathcal{S}^{D}$.

Forget that we fixed $\mathbf{a}$. We thus have shown that $\mathbf{a}%
\in\mathcal{S}^{D}$ for each $\mathbf{a}\in\mathcal{G}^{D}$. In other words,
$\mathcal{G}^{D}\subseteq\mathcal{S}^{D}$. Combining this with $\mathcal{S}%
^{D}\subseteq\mathcal{G}^{D}$, we obtain $\mathcal{G}^{D}=\mathcal{S}^{D}$.
This proves Theorem \ref{thm.cogarnir.CSG} \textbf{(a)}. \medskip

\textbf{(b)} Lemma \ref{lem.cogarnir.cut0} yields $\mathcal{C}^{D}%
\cap\mathcal{G}^{D}=0$. Since $\mathcal{G}^{D}=\mathcal{S}^{D}$ (by Theorem
\ref{thm.cogarnir.CSG} \textbf{(a)}), we can rewrite this as $\mathcal{C}%
^{D}\cap\mathcal{S}^{D}=0$. Hence, the sum $\mathcal{C}^{D}+\mathcal{S}^{D}$
is a direct sum -- i.e., we have $\mathcal{C}^{D}+\mathcal{S}^{D}%
=\mathcal{C}^{D}\oplus\mathcal{S}^{D}$.

But Lemma \ref{lem.cogarnir.sum} yields $\mathcal{M}^{D}=\mathcal{C}%
^{D}+\mathcal{S}^{D}=\mathcal{C}^{D}\oplus\mathcal{S}^{D}=\mathcal{S}%
^{D}\oplus\mathcal{C}^{D}$ as $\mathbf{k}$-modules. This proves Theorem
\ref{thm.cogarnir.CSG} \textbf{(b)}.
\end{proof}

\begin{proof}
[Proof of Theorem \ref{thm.cogarnir.cogarnir-skew}.]Theorem
\ref{thm.cogarnir.CSG} \textbf{(a)} yields $\mathcal{G}^{D}=\mathcal{S}^{D}$.
Hence,%
\[
\mathcal{S}^{D}=\mathcal{G}^{D}=\left\{  \mathbf{a}\in\mathcal{M}^{D}%
\ \mid\ \mathbf{a}\text{ satisfies all coGarnir equations}\right\}  .
\]
This proves Theorem \ref{thm.cogarnir.cogarnir-skew}.
\end{proof}

\subsubsection{\label{subsec.specht.cogarnir.triang}Triangulating row-Garnir
sets}

Theorem \ref{thm.cogarnir.cogarnir-skew} expresses the Specht module
$\mathcal{S}^{D}$ for a skew Young diagram $D$ as the solution set of a system
of linear equations -- namely, of the coGarnir equations
(\ref{eq.def.cogarnir.cogarnirTL.b.eq}). However, this system is vastly
redundant, in the sense that it only takes a (much smaller) subset to
determine the same solution set. A class of such sufficient subsets shall now
be revealed (in Theorem \ref{thm.cogarnir.cogarnir-essG-skew}). These come
from \emph{triangulating} sets, which are defined as follows (note the
similarity to Definition \ref{def.specht.TD.garnir-triasub} \textbf{(b)} and
\textbf{(c)}):

\begin{definition}
\label{def.cogarnir.essTL}Let $D$ be a diagram with $\left\vert D\right\vert
=n$. \medskip

\textbf{(a)} We let $\operatorname*{RowGar}\left(  D\right)  $ be the set of
all pairs $\left(  T,L\right)  $, where $T$ is an $n$-tableau of shape $D$ and
where $L$ is a row-Garnir set for $T$. \medskip

\textbf{(b)} A subset $\mathbf{G}$ of $\operatorname*{RowGar}\left(  D\right)
$ will be called \emph{triangulating} if it has the following property: For
each $n$-tableau $T$ of shape $D$ that is row-standard but not standard, there
exists some row-Garnir set $L$ for $T$ such that $\left(  T,L\right)
\in\mathbf{G}$ and such that
\begin{equation}
\text{each }w\in L\setminus\mathcal{R}\left(  T\right)  \text{ satisfies
}\overline{wT}>\overline{T} \label{eq.def.cogarnir.essTL.b.gr}%
\end{equation}
(with respect to the Young last letter order).
\end{definition}

\begin{example}
\label{exa.cogarnir.essTL.2}Let $n=4$ and $D=Y\left(  2,2\right)  $. There are
$4!=24$ many $n$-tableaux of shape $D$, but only $6$ of them are row-standard,
and only $4$ of these $6$ fail to be standard. Namely, the latter $4$ tableaux
are%
\[
\ytableaushort{14,23}\ \ ,\quad\ytableaushort{23,14}\ \ ,\quad
\ytableaushort{24,13}\ \ ,\quad\ytableaushort{34,12}\ \ .
\]
Thus, a subset $\mathbf{G}$ of $\operatorname*{RowGar}\left(  D\right)  $ is
said to be triangulating if it has the property that for each of these $4$
tableaux $T$, there exists some row-Garnir set $L$ for $T$ such that $\left(
T,L\right)  \in\mathbf{G}$ and such that (\ref{eq.def.cogarnir.essTL.b.gr})
holds. This property can be ensured by picking one such $L$ for each $T$; here
is one possibility (where the middle column shows the $Z$ from Definition
\ref{def.cogarnir.cogarnirTL}):%
\[%
\begin{tabular}
[c]{|c|c|c|}\hline
$T$ & $Z$ & $L$\\\hline\hline
$\ytableaushort{1{*(green)4},{*(green)2}{*(green)3}}$ & $\left\{
2,3,4\right\}  $ & $\left\{  \operatorname*{id},t_{2,4},t_{3,4}\right\}
$\\\hline
$\ytableaushort{{*(green)2}{*(green)3},{*(green)1}4}$ & $\left\{
1,2,3\right\}  $ & $\left\{  \operatorname*{id},t_{2,1},t_{3,1}\right\}
$\\\hline
$\ytableaushort{{*(green)2}{*(green)4},{*(green)1}3}$ & $\left\{
1,2,4\right\}  $ & $\left\{  \operatorname*{id},t_{2,1},t_{4,1}\right\}
$\\\hline
$\ytableaushort{{*(green)3}{*(green)4},{*(green)1}2}$ & $\left\{
1,3,4\right\}  $ & $\left\{  \operatorname*{id},t_{3,1},t_{4,1}\right\}
$\\\hline
\end{tabular}
\ \ .
\]
We have constructed these $L$'s as follows:

\begin{itemize}
\item We have picked some $\left(  j,i\right)  \in D$ such that $\left(
j+1,i\right)  \in D$ and $T\left(  j,i\right)  >T\left(  j+1,i\right)  $.
(Such a $\left(  j,i\right)  $ exists, since $T$ is not column-standard; we
have seen this in our above proof of Lemma \ref{lem.cogarnir.cut0}, where the
tableau $P$ played the role that $T$ is now playing. Sometimes there are two
choices for this $\left(  j,i\right)  $; in such cases, either choice is fine,
but we chose the lexicographically smaller one.)

\item We then defined the two sets%
\begin{align*}
X  &  =\left\{  T\left(  j,i^{\prime}\right)  \ \mid\ i^{\prime}\in
\mathbb{Z}\text{ with }i^{\prime}\geq i\text{ and }\left(  j,i^{\prime
}\right)  \in D\right\}  \ \ \ \ \ \ \ \ \ \ \text{and}\\
Y  &  =\left\{  T\left(  j+1,i^{\prime}\right)  \ \mid\ i^{\prime}%
\in\mathbb{Z}\text{ with }i^{\prime}\leq i\text{ and }\left(  j+1,i^{\prime
}\right)  \in D\right\}
\end{align*}
(just as we did in the proof of Lemma \ref{lem.cogarnir.cut0}, taking $k=j+1$).

\item We set $Z=X\cup Y$ (and we highlighted the elements of $Z$ in green in
the above table).

\item We chose a left transversal $L$ of $S_{n,Z}\cap\mathcal{R}\left(
T\right)  $ in $S_{n,Z}$. (Here again, there are several choices, but we
picked the simplest one, consisting of the identity and transpositions.)
\end{itemize}

The $L$'s defined in this way always satisfy (\ref{eq.def.cogarnir.essTL.b.gr}%
), as we have shown in the proof of Lemma \ref{lem.cogarnir.cut0} (see the
equality (\ref{pf.lem.cogarnir.cut0.ine2}), and keep in mind that $T$ is now
playing the role of $P$). Thus, the subset%
\begin{align*}
\mathbf{G}  &  :=\left\{  \left(
\ytableaushort{1{*(green)4},{*(green)2}{*(green)3}}\ \ ,\ \ \left\{
\operatorname*{id},t_{2,4},t_{3,4}\right\}  \right)  ,\qquad\left(
\ytableaushort{{*(green)2}{*(green)3},{*(green)1}4}\ \ ,\ \ \left\{
\operatorname*{id},t_{2,1},t_{3,1}\right\}  \right)  ,\right. \\
&  \ \ \ \ \ \ \ \ \ \ \left.  \left(
\ytableaushort{{*(green)2}{*(green)4},{*(green)1}3}\ \ ,\ \ \left\{
\operatorname*{id},t_{2,1},t_{4,1}\right\}  \right)  ,\qquad\left(
\ytableaushort{{*(green)3}{*(green)4},{*(green)1}2}\ \ ,\ \ \left\{
\operatorname*{id},t_{3,1},t_{4,1}\right\}  \right)  \right\}
\end{align*}
of $\operatorname*{RowGar}\left(  D\right)  $ is triangulating.
\end{example}

A safe but wasteful way to construct a triangulating subset of
$\operatorname*{RowGar}\left(  D\right)  $ (when $D$ is a skew Young diagram)
is by taking the whole $\operatorname*{RowGar}\left(  D\right)  $:

\begin{proposition}
\label{prop.cogarnir.rg-tri}Let $D$ be a skew Young diagram with $\left\vert
D\right\vert =n$. Then, the set $\operatorname*{RowGar}\left(  D\right)  $ is triangulating.
\end{proposition}

\begin{proof}
[Proof idea.]We must show that for each $n$-tableau $T$ of shape $D$ that is
row-standard but not standard, there exists some row-Garnir set $L$ for $T$
such that $\left(  T,L\right)  \in\operatorname*{RowGar}\left(  D\right)  $
and such that
\[
\text{each }w\in L\setminus\mathcal{R}\left(  T\right)  \text{ satisfies
}\overline{wT}>\overline{T}%
\]
(with respect to the Young last letter order). But this was (implicitly) shown
in our proof of Lemma \ref{lem.cogarnir.cut0} (except that $T$ was called $P$
in that proof). Thus, Proposition \ref{prop.cogarnir.rg-tri} is proved.
\end{proof}

\begin{definition}
\label{def.cogarnir.essallG}Let $D$ be a diagram with $\left\vert D\right\vert
=n$. Let $\mathbf{G}$ be a subset of $\operatorname*{RowGar}\left(  D\right)
$.

Let $\mathbf{a}\in\mathcal{M}^{D}$ be an element of the Young module
$\mathcal{M}^{D}$. Then, we say that $\mathbf{a}$ satisfies \emph{all
}$\mathbf{G}$\emph{-coGarnir equations} if $\mathbf{a}$ satisfies the $\left(
T,L\right)  $-coGarnir equation for each $\left(  T,L\right)  \in\mathbf{G}$.
\end{definition}

\begin{example}
\label{exa.cogarnir.essallG.2}Let $n=4$ and $D=Y\left(  2,2\right)  $. Let
$\mathbf{G}$ be the triangulating subset of $\operatorname*{RowGar}\left(
D\right)  $ constructed in Example \ref{exa.cogarnir.essTL.2}. Let
$\mathbf{a}=\left(  a_{\overline{P}}\right)  _{\overline{P}\in\overline
{\operatorname{Tab}}\left(  D\right)  }\in\mathcal{M}^{D}$ be an element of
$\mathcal{M}^{D}$. Then, $\mathbf{a}$ satisfies all $\mathbf{G}$-coGarnir
equations if and only if it satisfies the four equations%
\begin{align*}
a_{\overline{14\backslash\backslash23}}+a_{\overline{12\backslash\backslash
43}}+a_{\overline{13\backslash\backslash24}}  &  =0,\\
a_{\overline{23\backslash\backslash14}}+a_{\overline{13\backslash\backslash
24}}+a_{\overline{21\backslash\backslash34}}  &  =0,\\
a_{\overline{24\backslash\backslash13}}+a_{\overline{14\backslash\backslash
23}}+a_{\overline{21\backslash\backslash43}}  &  =0,\\
a_{\overline{34\backslash\backslash12}}+a_{\overline{14\backslash\backslash
32}}+a_{\overline{31\backslash\backslash42}}  &  =0.
\end{align*}
It should be kept in mind that the subscripts here are $n$-tabloids, not
$n$-tableaux, so several of them are equal (e.g., we have $\overline
{31\backslash\backslash42}=\overline{13\backslash\backslash24}$). Thus, these
are four linear equations on the altogether six coordinates of $\mathbf{a}$.
The solution set of these four equations turns out to be precisely the Specht
module $\mathcal{S}^{D}$. This is no coincidence, as the following theorem shows:
\end{example}

\begin{theorem}
\label{thm.cogarnir.cogarnir-essG-skew}Let $D$ be a skew Young diagram with
$\left\vert D\right\vert =n$. Let $\mathbf{G}$ be a triangulating subset of
$\operatorname*{RowGar}\left(  D\right)  $. Then, the Specht module
$\mathcal{S}^{D}$ can be described as%
\[
\mathcal{S}^{D}=\left\{  \mathbf{a}\in\mathcal{M}^{D}\ \mid\ \mathbf{a}\text{
satisfies all }\mathbf{G}\text{-coGarnir equations}\right\}  .
\]
(This is an identity, not just an isomorphism!)
\end{theorem}

\begin{proof}
This can be proved in the same way as Theorem \ref{thm.cogarnir.cogarnir-skew}%
, once we make a few changes to Subsection \ref{subsec.specht.cogarnir.suff}:

\begin{itemize}
\item The definition of $\mathcal{G}^{D}$ (in Definition
\ref{def.cogarnir.C-and-G}) must be replaced by
\[
\mathcal{G}^{D}:=\left\{  \mathbf{a}\in\mathcal{M}^{D}\ \mid\ \mathbf{a}\text{
satisfies all }\mathbf{G}\text{-coGarnir equations}\right\}  .
\]
(That is, instead of requiring $\mathbf{a}$ to satisfy all coGarnir equations,
we only require $\mathbf{a}$ to satisfy all $\mathbf{G}$-coGarnir equations.)

\item In the proof of Proposition \ref{prop.cogarnir.C-and-G-sub}, we must
replace \textquotedblleft all coGarnir equations\textquotedblright\ by
\textquotedblleft all $\mathbf{G}$-coGarnir equations\textquotedblright, and
we must replace \textquotedblleft for each $n$-tableau $T$ of shape $D$ and
each row-Garnir set $L$ for $T$\textquotedblright\ by \textquotedblleft for
each $\left(  T,L\right)  \in\mathbf{G}$\textquotedblright.

\item A similar change needs to be done in the proof of Lemma
\ref{lem.cogarnir.SsG}.

\item In the proof of Lemma \ref{lem.cogarnir.cut0}, we need to pick the
row-Garnir set $L$ strategically in order to ensure that $\left(  P,L\right)
\in\mathbf{G}$ (because this is needed to guarantee that $\mathbf{a}$
satisfies the $\left(  P,L\right)  $-coGarnir equation). Fortunately, the
triangulating property of $\mathbf{G}$ makes this easy: Indeed, since
$\mathbf{G}$ is triangulating, we know (from Definition
\ref{def.cogarnir.essTL} \textbf{(b)}) that for each $n$-tableau $T$ of shape
$D$ that is row-standard but not standard, there exists some row-Garnir set
$L$ for $T$ such that $\left(  T,L\right)  \in\mathbf{G}$ and such that
\[
\text{each }w\in L\setminus\mathcal{R}\left(  T\right)  \text{ satisfies
}\overline{wT}>\overline{T}%
\]
(with respect to the Young last letter order). We can apply this fact to $T=P$
(since $P$ is an $n$-tableau of shape $D$ that is row-standard but not
standard), and thus conclude that there exists some row-Garnir set $L$ for $P$
such that $\left(  P,L\right)  \in\mathbf{G}$ and such that
\[
\text{each }w\in L\setminus\mathcal{R}\left(  P\right)  \text{ satisfies
}\overline{wP}>\overline{P}%
\]
(with respect to the Young last letter order). Consider this $L$. Since
$\mathbf{a}$ satisfies all $\mathbf{G}$-coGarnir equations (because
$\mathbf{a}\in\mathcal{G}^{D}$), we thus see that $\mathbf{a}$ satisfies the
$\left(  P,L\right)  $-coGarnir equation (since $\left(  P,L\right)
\in\mathbf{G}$). In other words,
\[
\sum_{x\in L}a_{\overline{xP}}=0\ \ \ \ \ \ \ \ \ \ \left(  \text{by
Definition \ref{def.cogarnir.cogarnirTL} \textbf{(b)}}\right)  .
\]

Now, recall that $L$ is a left transversal of $S_{n,Z}\cap\mathcal{R}\left(
P\right)  $ in $S_{n,Z}$ for some subset $Z$ of $\left[  n\right]  $ (because
$L$ is a row-Garnir set for $P$). Consider this $Z$. Now, proceed as in the
original proof of Lemma \ref{lem.cogarnir.cut0}, starting with the sentence
\textquotedblleft Recall again that $L$ is a left transversal of $S_{n,Z}%
\cap\mathcal{R}\left(  P\right)  $ in $S_{n,Z}$\textquotedblright. This
results in a contradiction, just as in the original proof, and so the proof is
again complete.
\end{itemize}

Altogether, this modified version of Subsection
\ref{subsec.specht.cogarnir.suff} culminates in a proof of Theorem
\ref{thm.cogarnir.cogarnir-skew} but with
\[
\left\{  \mathbf{a}\in\mathcal{M}^{D}\ \mid\ \mathbf{a}\text{ satisfies all
coGarnir equations}\right\}
\]
replaced by
\[
\left\{  \mathbf{a}\in\mathcal{M}^{D}\ \mid\ \mathbf{a}\text{ satisfies all
}\mathbf{G}\text{-coGarnir equations}\right\}  .
\]
But this is precisely the claim of Theorem
\ref{thm.cogarnir.cogarnir-essG-skew}.
\end{proof}

\subsection{\label{sec.specht.hlf-young}Appendix: Young's proof of the hook
length formula}

In this section, we shall prove Theorem \ref{thm.tableau.hlf} (the hook length
formula). The proof I will present is neither the simplest nor the most
elementary; but it has the distinction of being the historically first proof,
as it comes from a 1902 paper by Young himself \cite[\textit{On Quantitative
Substitutional Analysis (Second Paper)}, \S 4--\S 5]{Young77}. Arguably, Young
did not know Theorem \ref{thm.tableau.hlf} in the form we stated it in, but
rather was proving something that only later turned out to be equivalent to
Theorem \ref{thm.tableau.hlf-young}, which is itself an equivalent version of
Theorem \ref{thm.tableau.hlf}. Thus, the proof I shall give in this section is
a somewhat liberal modification of Young's, although all the main ideas have
been preserved. The proof is long and rather intricate, and bears more
semblance to a proof in a modern research paper than to a proof of a classical
result in a standard textbook.\footnote{Garsia refers to it as a
\textquotedblleft gruesome brute force proof\textquotedblright\ in
\cite[Remark 1.5]{GarEge20}. This is somewhat unfair, as it falsely suggests a
straightforward approach that merely requires a lot of work. But the
gruesomeness is hard to argue against!} In my presentation, the proof gains
additional length, as numerous steps that Young considered obvious enough to
omit (often for lack of rigorous language to explain them in) have been
expanded into lemmas and auxiliary claims with detailed proofs. In some cases,
I was unable to understand Young's original arguments and had to replace them
by longer but clearer justification (Proposition \ref{prop.hlf-young.grel5},
which loosely corresponds to \cite[\textit{On Quantitative Substitutional
Analysis (Second Paper)}, \S 4]{Young77}, is one such case).\footnote{Another
exposition of Young's proof appears in \cite[Chapter II]{RotGui71}, which
stays closer to the original argument but avoids formalizing it completely as
\textquotedblleft the notation would get hopelessly
complicated\textquotedblright.}

\subsubsection{A consequence of the Garnir relations}

At the core of our proof will be an identity that relates a row symmetrizer
$\nabla_{\operatorname*{Row}T}$, a column antisymmetrizer $\nabla
_{\operatorname*{Col}T}^{-}$ and several transpositions $t_{u,v}$:

\begin{proposition}
\label{prop.hlf-young.grel5}Let $\lambda$ be a partition of $n$. Let $T$ be an
$n$-tableau of shape $\lambda$. Let $j$ and $k$ be two positive integers with
$j>k$. Let $q\in\left[  \lambda_{j}\right]  $. Then,%
\[
\left(  1+\sum_{s=1}^{\lambda_{k}}t_{T\left(  k,s\right)  ,\ T\left(
j,q\right)  }\right)  \nabla_{\operatorname*{Row}T}\nabla_{\operatorname*{Col}%
T}^{-}=0.
\]

\end{proposition}

\Needspace{8pc}

\begin{example}
\label{exa.hlf-young.grel5.1}Let $n=14$ and $\lambda=\left(  5,4,3,2\right)
$. Let $j=3$ and $k=2$ and $q=2$. Let $T$ be the following $n$-tableau of
shape $\lambda$:%
\[
T=\ytableaushort{{10}{11}{12}{13}{14},{*(yellow)6}{*(yellow)7}{*(yellow)8}{*(yellow)9},3{*(green)4}5,12}\ \ .
\]
Here, we have highlighted the cell $\left(  j,q\right)  =\left(  3,2\right)  $
green and the cells of the form $\left(  k,s\right)  $ (for $s\in\left[
\lambda_{k}\right]  $) yellow. Proposition \ref{prop.hlf-young.grel5} yields
that
\[
\left(  1+t_{6,4}+t_{7,4}+t_{8,4}+t_{9,4}\right)  \nabla_{\operatorname*{Row}%
T}\nabla_{\operatorname*{Col}T}^{-}=0.
\]

\end{example}

There are several ways to prove Proposition \ref{prop.hlf-young.grel5}. Rota
(in \cite[Chapter II, Prop 2]{RotGui71}) derives it by an artful application
of Theorem \ref{thm.specht.sandw0}.\footnote{This might be the argument Young
intended to make in \cite[\textit{On Quantitative Substitutional Analysis
(Second Paper)}, \S 4]{Young77}.} We shall instead prove Proposition
\ref{prop.hlf-young.grel5} by identifying it as a particular case of Corollary
\ref{cor.hlf-young.grel3} (which, as we recall, is a consequence of the Garnir
relations), obtained by choosing a specific left transversal $L$ (consisting
of the identity and some transpositions). This is made possible by the
following simple lemma:

\begin{lemma}
\label{lem.hlf-young.grel4}Let $X$ be a subset of $\left[  n\right]  $. Let
$x\in X$. Consider the set $S_{n,X}$ defined in Proposition
\ref{prop.intX.basics}, and the set $S_{n,X\setminus\left\{  x\right\}  }$
defined similarly. Then: \medskip

\textbf{(a)} We have%
\begin{equation}
\left\{  w\in S_{n,X}\ \mid\ w\left(  x\right)  =x\right\}  =S_{n,X\setminus
\left\{  x\right\}  }. \label{pf.lem.hlf-young.grel4.sub}%
\end{equation}

\textbf{(b)} The set $S_{n,X}$ is a subgroup of $S_{n}$, whereas the set
$S_{n,X\setminus\left\{  x\right\}  }$ is a subgroup of $S_{n,X}$.\medskip

\textbf{(c)} For each $u\in S_{n,X}$, we have $u\left(  x\right)  \in X$.
\medskip

\textbf{(d)} The set $\left\{  \operatorname*{id}\right\}  \cup\left\{
t_{y,x}\ \mid\ y\in X\setminus\left\{  x\right\}  \right\}  $ is a left
transversal of $S_{n,X\setminus\left\{  x\right\}  }$ in $S_{n,X}$. \medskip

\textbf{(e)} For each $y\in X$, fix a permutation $\sigma_{y}\in S_{n,X}$ that
sends $x$ to $y$. Then, the set $\left\{  \sigma_{y}\ \mid\ y\in X\right\}  $
is a left transversal of $S_{n,X\setminus\left\{  x\right\}  }$ in $S_{n,X}$.
\end{lemma}

\begin{fineprint}
\begin{proof}
This was already implicit in the proof of Lemma \ref{lem.intX.rec3} (except
that our $\sigma_{y}$ now correspond to the $\sigma_{y}^{-1}$ in Lemma
\ref{lem.intX.rec3}), but let us give a detailed proof nevertheless. \medskip

We know that $S_{n,X}$ is a subgroup of $S_{n}$ (by Proposition
\ref{prop.intX.basics} \textbf{(a)}). The same argument (applied to
$X\setminus\left\{  x\right\}  $ instead of $X$) shows that $S_{n,X\setminus
\left\{  x\right\}  }$ is a subgroup of $S_{n}$. \medskip

\textbf{(a)} The equality (\ref{pf.lem.hlf-young.grel4.sub}) is precisely the
equality (\ref{pf.lem.intX.rec3.sub}) that we showed in our proof of Lemma
\ref{lem.intX.rec3} above. Thus, Lemma \ref{lem.hlf-young.grel4} \textbf{(a)}
is proved. \medskip

\textbf{(b)} We have already shown that $S_{n,X}$ is a subgroup of $S_{n}$.
From (\ref{pf.lem.hlf-young.grel4.sub}), we obtain%
\[
S_{n,X\setminus\left\{  x\right\}  }=\left\{  w\in S_{n,X}\ \mid\ w\left(
x\right)  =x\right\}  \subseteq S_{n,X}.
\]
Thus, $S_{n,X\setminus\left\{  x\right\}  }$ is a subgroup of $S_{n,X}$ (since
both $S_{n,X\setminus\left\{  x\right\}  }$ and $S_{n,X}$ are subgroups of
$S_{n}$). This completes the proof of Lemma \ref{lem.hlf-young.grel4}
\textbf{(b)}. \medskip

\textbf{(c)} Let $u\in S_{n,X}$. Thus, $u\in S_{n}$ and $u\left(  i\right)
=i$ for all $i\notin X$ (by the definition of $S_{n,X}$).

We must show that $u\left(  x\right)  \in X$. Assume the contrary. Thus,
$u\left(  x\right)  \notin X$. Set $z:=u\left(  x\right)  $. Thus, $z=u\left(
x\right)  \notin X$. But recall that $u\left(  i\right)  =i$ for all $i\notin
X$. Applying this to $i=z$, we find $u\left(  z\right)  =z$ (since $z\notin
X$). Hence, $u\left(  z\right)  =z=u\left(  x\right)  $. Hence, $z=x$ (since
$u$ is a permutation and thus injective), so that $z=x\in X$. This contradicts
$z\notin X$. This contradiction shows that our assumption was false. Hence,
$u\left(  x\right)  \in X$ must hold. This proves Lemma
\ref{lem.hlf-young.grel4} \textbf{(c)}. \medskip

\textbf{(e)} Let $L:=\left\{  \sigma_{y}\ \mid\ y\in X\right\}  $. This is
clearly a subset of $S_{n,X}$.

We know (from Lemma \ref{lem.hlf-young.grel4} \textbf{(b)}) that
$S_{n,X\setminus\left\{  x\right\}  }$ is a subgroup of $S_{n,X}$. Let
$uS_{n,X\setminus\left\{  x\right\}  }$ be any left coset of $S_{n,X\setminus
\left\{  x\right\}  }$ in $S_{n,X}$. We shall show that $uS_{n,X\setminus
\left\{  x\right\}  }$ contains exactly one element of $L$.

Indeed, $u\in S_{n,X}$.

Let $z:=u\left(  x\right)  $. Then, $z=u\left(  x\right)  \in X$ (by Lemma
\ref{lem.hlf-young.grel4} \textbf{(c)}). Hence, the permutation $\sigma_{z}\in
S_{n,X}$ is well-defined, and belongs to $L$ (since $L=\left\{  \sigma
_{y}\ \mid\ y\in X\right\}  $). Thus, $\sigma_{z}\in L$. Moreover, the
definition of $\sigma_{z}$ shows that $\sigma_{z}$ sends $x$ to $z$. In other
words, $\sigma_{z}\left(  x\right)  =z$. From $\sigma_{z}\in S_{n,X}$ and
$u\in S_{n,X}$, we obtain $u^{-1}\sigma_{z}\in S_{n,X}$ (since $S_{n,X}$ is a
group). Moreover,%
\[
\left(  u^{-1}\sigma_{z}\right)  \left(  x\right)  =u^{-1}\left(
\underbrace{\sigma_{z}\left(  x\right)  }_{=z=u\left(  x\right)  }\right)
=u^{-1}\left(  u\left(  x\right)  \right)  =x.
\]
Thus, $u^{-1}\sigma_{z}$ is a $w\in S_{n,X}$ that satisfies $w\left(
x\right)  =x$ (since $u^{-1}\sigma_{z}\in S_{n,X}$). In other words,%
\[
u^{-1}\sigma_{z}\in\left\{  w\in S_{n,X}\ \mid\ w\left(  x\right)  =x\right\}
=S_{n,X\setminus\left\{  x\right\}  }\ \ \ \ \ \ \ \ \ \ \left(  \text{by
(\ref{pf.lem.hlf-young.grel4.sub})}\right)  .
\]
In other words, $\sigma_{z}\in uS_{n,X\setminus\left\{  x\right\}  }$.
Combining this with $\sigma_{z}\in L$, we obtain $\sigma_{z}\in\left(
uS_{n,X\setminus\left\{  x\right\}  }\right)  \cap L$. Thus,%
\begin{equation}
\left\{  \sigma_{z}\right\}  \subseteq\left(  uS_{n,X\setminus\left\{
x\right\}  }\right)  \cap L. \label{pf.lem.hlf-young.grel4.c.4}%
\end{equation}

Let us now prove that $\left(  uS_{n,X\setminus\left\{  x\right\}  }\right)
\cap L=\left\{  \sigma_{z}\right\}  $. Indeed, let $g\in\left(
uS_{n,X\setminus\left\{  x\right\}  }\right)  \cap L$. Thus, $g\in
uS_{n,X\setminus\left\{  x\right\}  }$ and $g\in L$. From $g\in L=\left\{
\sigma_{y}\ \mid\ y\in X\right\}  $, we see that $g=\sigma_{y}$ for some $y\in
X$. Consider this $y$. From $g\in uS_{n,X\setminus\left\{  x\right\}  }$, we
see that $g=uv$ for some $v\in S_{n,X\setminus\left\{  x\right\}  }$. Consider
this $v$. We have $v\in S_{n,X\setminus\left\{  x\right\}  }=\left\{  w\in
S_{n,X}\ \mid\ w\left(  x\right)  =x\right\}  $ (by
(\ref{pf.lem.hlf-young.grel4.sub})), so that $v\left(  x\right)  =x$. Now,
$\underbrace{g}_{=uv}\left(  x\right)  =\left(  uv\right)  \left(  x\right)
=u\left(  \underbrace{v\left(  x\right)  }_{=x}\right)  =u\left(  x\right)
=z$. However, from $g=\sigma_{y}$, we obtain $g\left(  x\right)  =\sigma
_{y}\left(  x\right)  =y$ (since $\sigma_{y}$ sends $x$ to $y$). Comparing
this with $g\left(  x\right)  =z$, we find $y=z$. Now, $g=\sigma_{y}%
=\sigma_{z}$ (since $y=z$), so that $g\in\left\{  \sigma_{z}\right\}  $.

Forget that we fixed $g$. We thus have shown that $g\in\left\{  \sigma
_{z}\right\}  $ for each $g\in\left(  uS_{n,X\setminus\left\{  x\right\}
}\right)  \cap L$. In other words, $\left(  uS_{n,X\setminus\left\{
x\right\}  }\right)  \cap L\subseteq\left\{  \sigma_{z}\right\}  $. Combining
this with (\ref{pf.lem.hlf-young.grel4.c.4}), we obtain%
\[
\left(  uS_{n,X\setminus\left\{  x\right\}  }\right)  \cap L=\left\{
\sigma_{z}\right\}  .
\]
Hence, $\left\vert \left(  uS_{n,X\setminus\left\{  x\right\}  }\right)  \cap
L\right\vert =\left\vert \left\{  \sigma_{z}\right\}  \right\vert =1$. In
other words, the intersection $\left(  uS_{n,X\setminus\left\{  x\right\}
}\right)  \cap L$ contains exactly one element. In other words, the set
$uS_{n,X\setminus\left\{  x\right\}  }$ contains exactly one element of $L$.

Forget that we fixed $uS_{n,X\setminus\left\{  x\right\}  }$. We thus have
shown that each left coset $uS_{n,X\setminus\left\{  x\right\}  }$ contains
exactly one element of $L$. In other words, $L$ is a left transversal of
$S_{n,X\setminus\left\{  x\right\}  }$ in $S_{n,X}$ (by the definition of a
left transversal). In other words, $\left\{  \sigma_{y}\ \mid\ y\in X\right\}
$ is a left transversal of $S_{n,X\setminus\left\{  x\right\}  }$ in $S_{n,X}$
(since $L=\left\{  \sigma_{y}\ \mid\ y\in X\right\}  $). This proves Lemma
\ref{lem.hlf-young.grel4} \textbf{(e)}. \medskip

\textbf{(d)} The transposition $t_{y,x}$ is well-defined for each $y\in
X\setminus\left\{  x\right\}  $, and belongs to $S_{n,X}$ (since both $y$ and
$x$ belong to $X$). Let us furthermore define a permutation $t_{x,x}\in
S_{n,X}$ by $t_{x,x}:=\operatorname*{id}$. Then, for each $y\in X$, the
permutation $t_{y,x}$ is well-defined and belongs to $S_{n,X}$. Moreover, this
permutation $t_{y,x}$ sends $x$ to $y$ (by its definition when $y\neq x$, and
for obvious reasons when $y=x$). Hence, Lemma \ref{lem.hlf-young.grel4}
\textbf{(e)} (applied to $\sigma_{y}=t_{y,x}$) shows that the set $\left\{
t_{y,x}\ \mid\ y\in X\right\}  $ is a left transversal of $S_{n,X\setminus
\left\{  x\right\}  }$ in $S_{n,X}$. Since%
\begin{align*}
\left\{  t_{y,x}\ \mid\ y\in X\right\}   &  =\left\{  \underbrace{t_{x,x}%
}_{=\operatorname*{id}}\right\}  \cup\left\{  t_{y,x}\ \mid\ y\in
X\setminus\left\{  x\right\}  \right\} \\
&  \ \ \ \ \ \ \ \ \ \ \ \ \ \ \ \ \ \ \ \ \left(
\begin{array}
[c]{c}%
\text{here, we have split off the element for }y=x\text{,}\\
\text{since }x\in X
\end{array}
\right) \\
&  =\left\{  \operatorname*{id}\right\}  \cup\left\{  t_{y,x}\ \mid\ y\in
X\setminus\left\{  x\right\}  \right\}  ,
\end{align*}
we can rewrite this as follows: The set $\left\{  \operatorname*{id}\right\}
\cup\left\{  t_{y,x}\ \mid\ y\in X\setminus\left\{  x\right\}  \right\}  $ is
a left transversal of $S_{n,X\setminus\left\{  x\right\}  }$ in $S_{n,X}$.
This proves Lemma \ref{lem.hlf-young.grel4} \textbf{(d)}.
\end{proof}
\end{fineprint}

We will also use two nearly trivial facts about the rows of a Young diagram
$Y\left(  \lambda\right)  $. Recall Definition \ref{def.diagram.Djr} and
Definition \ref{def.tableau.RowTi}.

\begin{lemma}
\label{lem.hlf-young.grel-lamd}Let $\lambda$ be a partition of $n$. Let
$D=Y\left(  \lambda\right)  $. Let $i$ be a positive integer. Then,
$D\left\lfloor i\right\rfloor =\left[  \lambda_{i}\right]  $.
\end{lemma}

\begin{fineprint}
\begin{proof}
We have $D=Y\left(  \lambda\right)  $. Thus, the cells in the $i$-th row of
$D$ are $\left(  i,1\right)  ,\ \left(  i,2\right)  ,\ \ldots,\ \left(
i,\lambda_{i}\right)  $ (by the definition of $Y\left(  \lambda\right)  $). In
other words, the integers $j$ that satisfy $\left(  i,j\right)  \in D$ are
$1,2,\ldots,\lambda_{i}$.

But $D\left\lfloor i\right\rfloor $ was defined as the set of all integers $j$
that satisfy $\left(  i,j\right)  \in D$. Hence, $D\left\lfloor i\right\rfloor
=\left\{  1,2,\ldots,\lambda_{i}\right\}  $ (since the integers $j$ that
satisfy $\left(  i,j\right)  \in D$ are $1,2,\ldots,\lambda_{i}$). Thus,
$D\left\lfloor i\right\rfloor =\left\{  1,2,\ldots,\lambda_{i}\right\}
=\left[  \lambda_{i}\right]  $. This proves Lemma
\ref{lem.hlf-young.grel-lamd}.
\end{proof}
\end{fineprint}

\begin{proposition}
\label{prop.tableau.RowTj}Let $D$ be a diagram. Let $T$ be an $n$-tableau of
shape $D$. Let $j\in\mathbb{Z}$. Then, $\left\vert \operatorname*{Row}\left(
j,T\right)  \right\vert =\left\vert D\left\lfloor j\right\rfloor \right\vert $.
\end{proposition}

\begin{fineprint}
\begin{proof}
This is the analogue of Proposition \ref{prop.tableau.ColTj} for rows instead
of columns. The proof is analogous.
\end{proof}
\end{fineprint}

Two more lemmas will be useful. Both lemmas are easy, and I recommend proving
them yourself:

\begin{lemma}
\label{lem.hlf-young.grel6}Let $D$ be a diagram. Let $T$ be an $n$-tableau of
shape $D$. Let $k\in\mathbb{Z}$. Then, $S_{n,\operatorname*{Row}\left(
k,T\right)  }\subseteq\mathcal{R}\left(  T\right)  $. (For the definition of
$S_{n,\operatorname*{Row}\left(  k,T\right)  }$, see Proposition
\ref{prop.intX.basics}.)
\end{lemma}

\begin{fineprint}
\begin{proof}
Let $u\in S_{n,\operatorname*{Row}\left(  k,T\right)  }$. Thus,%
\[
u\in S_{n,\operatorname*{Row}\left(  k,T\right)  }=\left\{  w\in S_{n}%
\ \mid\ w\left(  i\right)  =i\text{ for all }i\notin\operatorname*{Row}\left(
k,T\right)  \right\}
\]
(by the definition of $S_{n,\operatorname*{Row}\left(  k,T\right)  }$). In
other words, we have $u\in S_{n}$ and%
\begin{equation}
u\left(  i\right)  =i\text{ for all }i\notin\operatorname*{Row}\left(
k,T\right)  . \label{pf.lem.hlf-young.grel6.1}%
\end{equation}

Now, let $i\in\left[  n\right]  $. We shall show that the number $u\left(
i\right)  $ lies in the same row of $T$ as $i$ does. If $u\left(  i\right)
=i$, then this is obviously true. Thus, for the rest of this proof, we WLOG
assume that $u\left(  i\right)  \neq i$. Thus, $i\in\operatorname*{Row}\left(
k,T\right)  $ (because otherwise, we would have $i\notin\operatorname*{Row}%
\left(  k,T\right)  $ and thus $u\left(  i\right)  =i$ (by
(\ref{pf.lem.hlf-young.grel6.1})), which would contradict $u\left(  i\right)
\neq i$). In other words, $i$ lies in the $k$-th row of $T$. Furthermore,
Lemma \ref{lem.hlf-young.grel4} \textbf{(c)} (applied to
$X=\operatorname*{Row}\left(  k,T\right)  $ and $x=i$) yields $u\left(
i\right)  \in\operatorname*{Row}\left(  k,T\right)  $ (since $i\in
\operatorname*{Row}\left(  k,T\right)  $ and $u\in S_{n,\operatorname*{Row}%
\left(  k,T\right)  }$). In other words, $u\left(  i\right)  $ lies in the
$k$-th row of $T$. Thus, we have shown that both $u\left(  i\right)  $ and $i$
lie in the $k$-th row of $T$. Therefore, the number $u\left(  i\right)  $ lies
in the same row of $T$ as $i$ does.

Forget that we fixed $i$. We thus have shown that for each $i\in\left[
n\right]  $, the number $u\left(  i\right)  $ lies in the same row of $T$ as
$i$ does. In other words, the permutation $u$ is horizontal for $T$ (by the
definition of \textquotedblleft horizontal\textquotedblright). In other words,
$u\in\mathcal{R}\left(  T\right)  $.

Forget that we fixed $u$. We thus have proved that $u\in\mathcal{R}\left(
T\right)  $ for each $u\in S_{n,\operatorname*{Row}\left(  k,T\right)  }$. In
other words, $S_{n,\operatorname*{Row}\left(  k,T\right)  }\subseteq
\mathcal{R}\left(  T\right)  $. This proves Lemma \ref{lem.hlf-young.grel6}.
\end{proof}
\end{fineprint}

\begin{lemma}
\label{lem.hlf-young.grel7}Let $D$ be a diagram. Let $T$ be an $n$-tableau of
shape $D$. Let $j$ and $k$ be two distinct integers. Let $z\in
\operatorname*{Row}\left(  j,T\right)  $ be arbitrary. Let $Z$ be the set
$\left\{  z\right\}  \cup\operatorname*{Row}\left(  k,T\right)  $. Then,
$\left\vert Z\right\vert =\left\vert D\left\lfloor k\right\rfloor \right\vert
+1$ and $Z\setminus\left\{  z\right\}  =\operatorname*{Row}\left(  k,T\right)
$ and $S_{n,Z}\cap\mathcal{R}\left(  T\right)  =S_{n,Z\setminus\left\{
z\right\}  }$.
\end{lemma}

\begin{fineprint}
\begin{proof}
Proposition \ref{prop.tableau.RowTj} (applied to $k$ instead of $j$) yields
$\left\vert \operatorname*{Row}\left(  k,T\right)  \right\vert =\left\vert
D\left\lfloor k\right\rfloor \right\vert $.

We have $z\in\operatorname*{Row}\left(  j,T\right)  $. In other words, $z$ is
an entry in the $j$-th row of $T$ (since $\operatorname*{Row}\left(
j,T\right)  $ is defined as the set of all entries in the $j$-th row of $T$).

But $T$ is an $n$-tableau, and thus is injective. Hence, all entries of $T$
are distinct. In particular, an entry in the $j$-th row of $T$ cannot also be
found in the $k$-th row of $T$ (since $j$ and $k$ are distinct). In other
words, an entry in the $j$-th row of $T$ cannot belong to $\operatorname*{Row}%
\left(  k,T\right)  $ (since $\operatorname*{Row}\left(  k,T\right)  $ is
defined as the set of all entries in the $k$-th row of $T$). Hence, in
particular, $z$ cannot belong to $\operatorname*{Row}\left(  k,T\right)  $
(since $z$ is an entry in the $j$-th row of $T$). In other words,
$z\notin\operatorname*{Row}\left(  k,T\right)  $. Hence, $\left\vert
\operatorname*{Row}\left(  k,T\right)  \cup\left\{  z\right\}  \right\vert
=\left\vert \operatorname*{Row}\left(  k,T\right)  \right\vert +1$.

The definition of $Z$ shows that
\[
Z=\left\{  z\right\}  \cup\operatorname*{Row}\left(  k,T\right)
=\operatorname*{Row}\left(  k,T\right)  \cup\left\{  z\right\}  .
\]
Hence,%
\[
\left\vert Z\right\vert =\left\vert \operatorname*{Row}\left(  k,T\right)
\cup\left\{  z\right\}  \right\vert =\underbrace{\left\vert
\operatorname*{Row}\left(  k,T\right)  \right\vert }_{=\left\vert
D\left\lfloor k\right\rfloor \right\vert }+\,1=\left\vert D\left\lfloor
k\right\rfloor \right\vert +1.
\]

From $Z=\operatorname*{Row}\left(  k,T\right)  \cup\left\{  z\right\}  $, we
obtain
\[
Z\setminus\left\{  z\right\}  =\left(  \operatorname*{Row}\left(  k,T\right)
\cup\left\{  z\right\}  \right)  \setminus\left\{  z\right\}
=\operatorname*{Row}\left(  k,T\right)
\]
(since $z\notin\operatorname*{Row}\left(  k,T\right)  $) and%
\[
Z\setminus\operatorname*{Row}\left(  k,T\right)  =\left(  \operatorname*{Row}%
\left(  k,T\right)  \cup\left\{  z\right\}  \right)  \setminus
\operatorname*{Row}\left(  k,T\right)  =\left\{  z\right\}
\]
(again since $z\notin\operatorname*{Row}\left(  k,T\right)  $). Also,
$z\in\left\{  z\right\}  \subseteq\left\{  z\right\}  \cup\operatorname*{Row}%
\left(  k,T\right)  =Z$. Hence, Lemma \ref{lem.hlf-young.grel4} \textbf{(a)}
(applied to $X=Z$ and $x=z$) shows that
\begin{equation}
\left\{  w\in S_{n,Z}\ \mid\ w\left(  z\right)  =z\right\}  =S_{n,Z\setminus
\left\{  z\right\}  }. \label{pf.lem.hlf-young.grel7.w=S}%
\end{equation}
Hence, $S_{n,Z\setminus\left\{  z\right\}  }=\left\{  w\in S_{n,Z}%
\ \mid\ w\left(  z\right)  =z\right\}  \subseteq S_{n,Z}$.

It remains to show that $S_{n,Z}\cap\mathcal{R}\left(  T\right)
=S_{n,Z\setminus\left\{  z\right\}  }$. For this purpose, we shall prove that
$S_{n,Z}\cap\mathcal{R}\left(  T\right)  \subseteq S_{n,Z\setminus\left\{
z\right\}  }$ and $S_{n,Z\setminus\left\{  z\right\}  }\subseteq S_{n,Z}%
\cap\mathcal{R}\left(  T\right)  $.

Let $u\in S_{n,Z}\cap\mathcal{R}\left(  T\right)  $ be arbitrary. Thus, $u\in
S_{n,Z}$ and $u\in\mathcal{R}\left(  T\right)  $. We have $z\in Z$. Thus, from
$u\in S_{n,Z}$, we obtain $u\left(  z\right)  \in Z$ (by Lemma
\ref{lem.hlf-young.grel4} \textbf{(c)}, applied to $X=Z$ and $x=z$). However,
from $u\in\mathcal{R}\left(  T\right)  $, we see that the permutation $u$ is
horizontal for $T$. Thus, for each $i\in\left[  n\right]  $, the number
$u\left(  i\right)  $ lies in the same row of $T$ as $i$ does (by the
definition of a horizontal permutation). Applying this to $i=z$, we see that
the number $u\left(  z\right)  $ lies in the same row of $T$ as $z$ does.
Since $z$ lies in the $j$-th row of $T$, we thus conclude that $u\left(
z\right)  $ lies in the $j$-th row of $T$ as well. In other words, $u\left(
z\right)  $ is an entry in the $j$-th row of $T$. Hence, $u\left(  z\right)  $
cannot belong to $\operatorname*{Row}\left(  k,T\right)  $ (since an entry in
the $j$-th row of $T$ cannot belong to $\operatorname*{Row}\left(  k,T\right)
$). In other words, $u\left(  z\right)  \notin\operatorname*{Row}\left(
k,T\right)  $.

Combining $u\left(  z\right)  \in Z$ with $u\left(  z\right)  \notin%
\operatorname*{Row}\left(  k,T\right)  $, we find%
\[
u\left(  z\right)  \in Z\setminus\operatorname*{Row}\left(  k,T\right)
=\left\{  z\right\}  .
\]
In other words, $u\left(  z\right)  =z$. Hence, $u\in S_{n,Z}$ and $u\left(
z\right)  =z$. In other words,%
\[
u\in\left\{  w\in S_{n,Z}\ \mid\ w\left(  z\right)  =z\right\}
=S_{n,Z\setminus\left\{  z\right\}  }\ \ \ \ \ \ \ \ \ \ \left(  \text{by
(\ref{pf.lem.hlf-young.grel7.w=S})}\right)  .
\]

Forget that we fixed $u$. We thus have shown that $u\in S_{n,Z\setminus
\left\{  z\right\}  }$ for each $u\in S_{n,Z}\cap\mathcal{R}\left(  T\right)
$. In other words,%
\begin{equation}
S_{n,Z}\cap\mathcal{R}\left(  T\right)  \subseteq S_{n,Z\setminus\left\{
z\right\}  }. \label{pf.lem.hlf-young.grel7.dir1}%
\end{equation}

Let us now prove the reverse inclusion. We have
\[
S_{n,Z\setminus\left\{  z\right\}  }=\underbrace{S_{n,Z\setminus\left\{
z\right\}  }}_{\subseteq S_{n,Z}}\cap\underbrace{S_{n,Z\setminus\left\{
z\right\}  }}_{\substack{=S_{n,\operatorname*{Row}\left(  k,T\right)
}\\\text{(since }Z\setminus\left\{  z\right\}  =\operatorname*{Row}\left(
k,T\right)  \text{)}}}\subseteq S_{n,Z}\cap
\underbrace{S_{n,\operatorname*{Row}\left(  k,T\right)  }}%
_{\substack{\subseteq\mathcal{R}\left(  T\right)  \\\text{(by Lemma
\ref{lem.hlf-young.grel6})}}}\subseteq S_{n,Z}\cap\mathcal{R}\left(  T\right)
.
\]
Combining this with (\ref{pf.lem.hlf-young.grel7.dir1}), we obtain
$S_{n,Z}\cap\mathcal{R}\left(  T\right)  =S_{n,Z\setminus\left\{  z\right\}
}$. This completes the proof of Lemma \ref{lem.hlf-young.grel7} (since
$\left\vert Z\right\vert =\left\vert D\left\lfloor k\right\rfloor \right\vert
+1$ and $Z\setminus\left\{  z\right\}  =\operatorname*{Row}\left(  k,T\right)
$ have already been proved).
\end{proof}
\end{fineprint}

We are now ready to prove Proposition \ref{prop.hlf-young.grel5}:

\begin{proof}
[Proof of Proposition \ref{prop.hlf-young.grel5}.]Let $D:=Y\left(
\lambda\right)  $.

We know that $T$ is an $n$-tableau of shape $\lambda$. In other words, $T$ is
an $n$-tableau of shape $Y\left(  \lambda\right)  $. In other words, $T$ is an
$n$-tableau of shape $D$ (since $D=Y\left(  \lambda\right)  $). Hence, $T$ is
injective (since any $n$-tableau is injective). In other words, all entries of
$T$ are distinct.

Let $z:=T\left(  j,q\right)  $. Thus, $z$ is an entry in the $j$-th row of $T$
(since $\left(  j,q\right)  $ is a cell in the $j$-th row). In other words,
$z\in\operatorname*{Row}\left(  j,T\right)  $. Thus, $\left\{  z\right\}
\subseteq\operatorname*{Row}\left(  j,T\right)  $.

Let $Z$ be the set $\left\{  z\right\}  \cup\operatorname*{Row}\left(
k,T\right)  $. Thus,%
\[
Z=\underbrace{\left\{  z\right\}  }_{\subseteq\operatorname*{Row}\left(
j,T\right)  }\cup\operatorname*{Row}\left(  k,T\right)  \subseteq
\operatorname*{Row}\left(  j,T\right)  \cup\operatorname*{Row}\left(
k,T\right)  .
\]
Thus, $Z$ is a subset of $\operatorname*{Row}\left(  j,T\right)
\cup\operatorname*{Row}\left(  k,T\right)  $.

The numbers $j$ and $k$ are distinct (since $j>k$). Hence, Lemma
\ref{lem.hlf-young.grel7} yields that%
\begin{align*}
\left\vert Z\right\vert  &  =\left\vert D\left\lfloor k\right\rfloor
\right\vert +1\ \ \ \ \ \ \ \ \ \ \text{and}\ \ \ \ \ \ \ \ \ \ Z\setminus
\left\{  z\right\}  =\operatorname*{Row}\left(  k,T\right) \\
&  \text{and}\ \ \ \ \ \ \ \ \ \ S_{n,Z}\cap\mathcal{R}\left(  T\right)
=S_{n,Z\setminus\left\{  z\right\}  }.
\end{align*}

Also, obviously, $z\in\left\{  z\right\}  \subseteq\left\{  z\right\}
\cup\operatorname*{Row}\left(  k,T\right)  =Z$.

Lemma \ref{lem.hlf-young.grel-lamd} yields $D\left\lfloor j\right\rfloor
=\left[  \lambda_{j}\right]  $ and $D\left\lfloor k\right\rfloor =\left[
\lambda_{k}\right]  $. But $\lambda$ is a partition, thus a weakly decreasing
sequence. Hence, $\lambda_{j}\leq\lambda_{k}$ (since $j>k$). Hence, $\left[
\lambda_{j}\right]  \subseteq\left[  \lambda_{k}\right]  $. In other words,
$D\left\lfloor j\right\rfloor \subseteq D\left\lfloor k\right\rfloor $ (since
$D\left\lfloor j\right\rfloor =\left[  \lambda_{j}\right]  $ and
$D\left\lfloor k\right\rfloor =\left[  \lambda_{k}\right]  $), so that%
\[
D\left\lfloor j\right\rfloor \cup D\left\lfloor k\right\rfloor =D\left\lfloor
k\right\rfloor .
\]
However,
\[
\left\vert Z\right\vert =\left\vert D\left\lfloor k\right\rfloor \right\vert
+1>\left\vert \underbrace{D\left\lfloor k\right\rfloor }_{=D\left\lfloor
j\right\rfloor \cup D\left\lfloor k\right\rfloor }\right\vert =\left\vert
D\left\lfloor j\right\rfloor \cup D\left\lfloor k\right\rfloor \right\vert .
\]

Lemma \ref{lem.hlf-young.grel4} \textbf{(d)} (applied to $X=Z$ and $x=z$)
shows that the set $\left\{  \operatorname*{id}\right\}  \cup\left\{
t_{y,z}\ \mid\ y\in Z\setminus\left\{  z\right\}  \right\}  $ is a left
transversal of $S_{n,Z\setminus\left\{  z\right\}  }$ in $S_{n,Z}$. In other
words, the set $\left\{  \operatorname*{id}\right\}  \cup\left\{
t_{y,z}\ \mid\ y\in Z\setminus\left\{  z\right\}  \right\}  $ is a left
transversal of $S_{n,Z}\cap\mathcal{R}\left(  T\right)  $ in $S_{n,Z}$ (since
$S_{n,Z}\cap\mathcal{R}\left(  T\right)  =S_{n,Z\setminus\left\{  z\right\}
}$). Hence, Corollary \ref{cor.hlf-young.grel3} (applied to $L=\left\{
\operatorname*{id}\right\}  \cup\left\{  t_{y,z}\ \mid\ y\in Z\setminus
\left\{  z\right\}  \right\}  $) shows that%
\begin{equation}
\sum_{x\in\left\{  \operatorname*{id}\right\}  \cup\left\{  t_{y,z}%
\ \mid\ y\in Z\setminus\left\{  z\right\}  \right\}  }x\nabla
_{\operatorname*{Row}T}\nabla_{\operatorname*{Col}T}^{-}=0.
\label{pf.prop.hlf-young.grel5.at}%
\end{equation}

But $D=Y\left(  \lambda\right)  $. Hence, the cells in the $k$-th row of $D$
are $\left(  k,1\right)  ,\ \left(  k,2\right)  ,\ \ldots,\ \left(
k,\lambda_{k}\right)  $. Thus, the entries in the $k$-th row of $T$ are
$T\left(  k,1\right)  ,\ T\left(  k,2\right)  ,\ \ldots,\ T\left(
k,\lambda_{k}\right)  $. Now, recall that
\[
Z\setminus\left\{  z\right\}  =\operatorname*{Row}\left(  k,T\right)
=\left\{  T\left(  k,1\right)  ,\ T\left(  k,2\right)  ,\ \ldots,\ T\left(
k,\lambda_{k}\right)  \right\}
\]
(since the entries in the $k$-th row of $T$ are $T\left(  k,1\right)
,\ T\left(  k,2\right)  ,\ \ldots,\ T\left(  k,\lambda_{k}\right)  $). Hence,%
\begin{align*}
\left\{  t_{y,z}\ \mid\ y\in Z\setminus\left\{  z\right\}  \right\}   &
=\left\{  t_{y,z}\ \mid\ y\in\left\{  T\left(  k,1\right)  ,\ T\left(
k,2\right)  ,\ \ldots,\ T\left(  k,\lambda_{k}\right)  \right\}  \right\} \\
&  =\left\{  t_{T\left(  k,1\right)  ,\ z},\ t_{T\left(  k,2\right)
,\ z},\ \ldots,\ t_{T\left(  k,\lambda_{k}\right)  ,\ z}\right\} \\
&  =\left\{  t_{T\left(  k,1\right)  ,\ T\left(  j,q\right)  },\ t_{T\left(
k,2\right)  ,\ T\left(  j,q\right)  },\ \ldots,\ t_{T\left(  k,\lambda
_{k}\right)  ,\ T\left(  j,q\right)  }\right\}
\end{align*}
(since $z=T\left(  j,q\right)  $). Thus,%
\begin{align}
&  \left\{  \operatorname*{id}\right\}  \cup\left\{  t_{y,z}\ \mid\ y\in
Z\setminus\left\{  z\right\}  \right\} \nonumber\\
&  =\left\{  \operatorname*{id}\right\}  \cup\left\{  t_{T\left(  k,1\right)
,\ T\left(  j,q\right)  },\ t_{T\left(  k,2\right)  ,\ T\left(  j,q\right)
},\ \ldots,\ t_{T\left(  k,\lambda_{k}\right)  ,\ T\left(  j,q\right)
}\right\} \nonumber\\
&  =\left\{  \operatorname*{id},\ t_{T\left(  k,1\right)  ,\ T\left(
j,q\right)  },\ t_{T\left(  k,2\right)  ,\ T\left(  j,q\right)  }%
,\ \ldots,\ t_{T\left(  k,\lambda_{k}\right)  ,\ T\left(  j,q\right)
}\right\}  . \label{pf.prop.hlf-young.grel5.7}%
\end{align}
Moreover, all the $\lambda_{k}+1$ permutations $\operatorname*{id}%
,\ t_{T\left(  k,1\right)  ,\ T\left(  j,q\right)  },\ t_{T\left(  k,2\right)
,\ T\left(  j,q\right)  },\ \ldots,\ t_{T\left(  k,\lambda_{k}\right)
,\ T\left(  j,q\right)  }$ are distinct\footnote{\textit{Proof.} The
$\lambda_{k}$ entries $T\left(  k,1\right)  ,\ T\left(  k,2\right)
,\ \ldots,\ T\left(  k,\lambda_{k}\right)  $ of $T$ are distinct (since all
entries of $T$ are distinct). Hence, the $\lambda_{k}$ transpositions
$t_{T\left(  k,1\right)  ,\ T\left(  j,q\right)  },\ t_{T\left(  k,2\right)
,\ T\left(  j,q\right)  },\ \ldots,\ t_{T\left(  k,\lambda_{k}\right)
,\ T\left(  j,q\right)  }$ are distinct (since they send $T\left(  j,q\right)
$ to the distinct entries $T\left(  k,1\right)  ,\ T\left(  k,2\right)
,\ \ldots,\ T\left(  k,\lambda_{k}\right)  $). These $\lambda_{k}$
transpositions are furthermore all distinct from the identity permutation
$\operatorname*{id}$ (since any transposition is distinct from
$\operatorname*{id}$). Hence, altogether, all the $\lambda_{k}+1$ permutations
$\operatorname*{id},\ t_{T\left(  k,1\right)  ,\ T\left(  j,q\right)
},\ t_{T\left(  k,2\right)  ,\ T\left(  j,q\right)  },\ \ldots,\ t_{T\left(
k,\lambda_{k}\right)  ,\ T\left(  j,q\right)  }$ are distinct.}. Combining
this with (\ref{pf.prop.hlf-young.grel5.7}), we see that the set $\left\{
\operatorname*{id}\right\}  \cup\left\{  t_{y,z}\ \mid\ y\in Z\setminus
\left\{  z\right\}  \right\}  $ consists of the $\lambda_{k}+1$ distinct
elements $\operatorname*{id},\ t_{T\left(  k,1\right)  ,\ T\left(  j,q\right)
},\ t_{T\left(  k,2\right)  ,\ T\left(  j,q\right)  },\ \ldots,\ t_{T\left(
k,\lambda_{k}\right)  ,\ T\left(  j,q\right)  }$. Hence, the sum of the
elements of this set is%
\begin{align*}
\sum_{x\in\left\{  \operatorname*{id}\right\}  \cup\left\{  t_{y,z}%
\ \mid\ y\in Z\setminus\left\{  z\right\}  \right\}  }x  &
=\operatorname*{id}+\underbrace{t_{T\left(  k,1\right)  ,\ T\left(
j,q\right)  }+t_{T\left(  k,2\right)  ,\ T\left(  j,q\right)  }+\cdots
+t_{T\left(  k,\lambda_{k}\right)  ,\ T\left(  j,q\right)  }}_{=\sum
_{s=1}^{\lambda_{k}}t_{T\left(  k,s\right)  ,\ T\left(  j,q\right)  }}\\
&  =\operatorname*{id}+\sum_{s=1}^{\lambda_{k}}t_{T\left(  k,s\right)
,\ T\left(  j,q\right)  }.
\end{align*}
In other words,%
\[
\operatorname*{id}+\sum_{s=1}^{\lambda_{k}}t_{T\left(  k,s\right)  ,\ T\left(
j,q\right)  }=\sum_{x\in\left\{  \operatorname*{id}\right\}  \cup\left\{
t_{y,z}\ \mid\ y\in Z\setminus\left\{  z\right\}  \right\}  }x.
\]
Hence,%
\begin{align*}
\left(  1+\sum_{s=1}^{\lambda_{k}}t_{T\left(  k,s\right)  ,\ T\left(
j,q\right)  }\right)  \nabla_{\operatorname*{Row}T}\nabla_{\operatorname*{Col}%
T}^{-}  &  =\left(  \sum_{x\in\left\{  \operatorname*{id}\right\}
\cup\left\{  t_{y,z}\ \mid\ y\in Z\setminus\left\{  z\right\}  \right\}
}x\right)  \nabla_{\operatorname*{Row}T}\nabla_{\operatorname*{Col}T}^{-}\\
&  =\sum_{x\in\left\{  \operatorname*{id}\right\}  \cup\left\{  t_{y,z}%
\ \mid\ y\in Z\setminus\left\{  z\right\}  \right\}  }x\nabla
_{\operatorname*{Row}T}\nabla_{\operatorname*{Col}T}^{-}=0
\end{align*}
(by (\ref{pf.prop.hlf-young.grel5.at})). This proves Proposition
\ref{prop.hlf-young.grel5}.
\end{proof}

\subsubsection{Arms, legs and hooks}

Next, we shall take a closer look at the hooks $H_{\lambda}\left(  c\right)  $
introduced in Theorem \ref{thm.tableau.hlf}. We will subdivide them into a
\textquotedblleft long arm\textquotedblright\ and a \textquotedblleft short
leg\textquotedblright:

\begin{proposition}
\label{prop.hlf-young.arm-leg.1}Let $\lambda$ be a partition of $n$. Let
$c=\left(  i,j\right)  $ be a cell in $Y\left(  \lambda\right)  $. Define the
\emph{long arm} of $c$ in $\lambda$ to be the set%
\[
\operatorname*{Arm}\nolimits_{\lambda}^{+}\left(  c\right)  :=\left\{
c\right\}  \cup\left\{  \text{all cells in }Y\left(  \lambda\right)  \text{
that lie east of }c\text{ (in the same row as }c\text{)}\right\}  .
\]
Define the \emph{short leg} of $c$ in $\lambda$ to be the set%
\[
\operatorname*{Leg}\nolimits_{\lambda}^{-}\left(  c\right)  :=\left\{
\text{all cells in }Y\left(  \lambda\right)  \text{ that lie south of }c\text{
(in the same column as }c\text{)}\right\}  .
\]
Also, define the \emph{hook} $H_{\lambda}\left(  c\right)  $ of $c$ in
$\lambda$ as in Theorem \ref{thm.tableau.hlf}. Then: \medskip

\textbf{(a)} We have $\left\vert \operatorname*{Arm}\nolimits_{\lambda}%
^{+}\left(  c\right)  \right\vert =\lambda_{i}-j+1$. \medskip

\textbf{(b)} We have $\left\vert \operatorname*{Leg}\nolimits_{\lambda}%
^{-}\left(  c\right)  \right\vert =\lambda_{j}^{t}-i$. \medskip

\textbf{(c)} We have $\left\vert \operatorname*{Arm}\nolimits_{\lambda}%
^{+}\left(  c\right)  \right\vert +\left\vert \operatorname*{Leg}%
\nolimits_{\lambda}^{-}\left(  c\right)  \right\vert =\left\vert H_{\lambda
}\left(  c\right)  \right\vert $.
\end{proposition}

\begin{example}
Let $n=13$ and $\lambda=\left(  4,4,3,2\right)  $. Let $c$ be the cell
$\left(  2,1\right)  \in Y\left(  \lambda\right)  $. Then, the long arm of $c$
consists of the green cells in the following picture, whereas the short leg
consists of the orange cells:%
\[%
%TCIMACRO{\TeXButton{ydiagram}{\ydiagram[*(white)]{0+4, 0+0, 1+2, 1+1}%
%*[*(green)]{0+0, 0+4, 0+0,0+0}*[*(orange)]{0+0, 0+0, 0+1, 0+1}}}%
%BeginExpansion
\ydiagram[*(white)]{0+4, 0+0, 1+2, 1+1}*[*(green)]{0+0, 0+4, 0+0,0+0}%
*[*(orange)]{0+0, 0+0, 0+1, 0+1}%
%EndExpansion
\ \ .
\]
In other words,%
\begin{align*}
\operatorname*{Arm}\nolimits_{\lambda}^{+}\left(  c\right)   &  =\left\{
\left(  2,1\right)  ,\ \left(  2,2\right)  ,\ \left(  2,3\right)  ,\ \left(
2,4\right)  \right\}  \ \ \ \ \ \ \ \ \ \ \text{and}\\
\operatorname*{Leg}\nolimits_{\lambda}^{-}\left(  c\right)   &  =\left\{
\left(  3,1\right)  ,\ \left(  4,1\right)  \right\}  .
\end{align*}

\end{example}

The picture in the above example should explain the words \textquotedblleft
arm\textquotedblright\ and \textquotedblleft leg\textquotedblright. The words
\textquotedblleft long\textquotedblright\ and \textquotedblleft
short\textquotedblright\ refer to the respective presence and absence of $c$.
Of course, we could just as well define the \textquotedblleft short
arm\textquotedblright\ and the \textquotedblleft long leg\textquotedblright,
but we will not need these.

The proof of Proposition \ref{prop.hlf-young.arm-leg.1} presents only
typographical difficulties (due to the long equations):

\begin{fineprint}
\begin{proof}
[Proof of Proposition \ref{prop.hlf-young.arm-leg.1}.]\textbf{(a)} The cells
in the $i$-th row of $Y\left(  \lambda\right)  $ are $\left(  i,1\right)
,\ \left(  i,2\right)  ,\ \ldots,\ \left(  i,\lambda_{i}\right)  $ (by the
definition of $Y\left(  \lambda\right)  $). The definition of
$\operatorname*{Arm}\nolimits_{\lambda}^{+}\left(  c\right)  $ yields
\begin{align*}
\operatorname*{Arm}\nolimits_{\lambda}^{+}\left(  c\right)   &  =\left\{
c\right\}  \cup\left\{  \text{all cells in }Y\left(  \lambda\right)  \text{
that lie east of }c\text{ (in the same row as }c\text{)}\right\} \\
&  =\left\{  \left(  i,j\right)  \right\}  \cup\underbrace{\left\{  \text{all
cells in }Y\left(  \lambda\right)  \text{ that lie east of }\left(
i,j\right)  \text{ (in the same row as }\left(  i,j\right)  \text{)}\right\}
}_{\substack{=\left\{  \text{all cells in }Y\left(  \lambda\right)  \text{
that lie east of }\left(  i,j\right)  \text{ (in the }i\text{-th
row)}\right\}  \\=\left\{  \text{all cells in the }i\text{-th row of }Y\left(
\lambda\right)  \text{ that lie east of }\left(  i,j\right)  \right\}
\\=\left\{  \left(  i,j+1\right)  ,\ \left(  i,j+2\right)  ,\ \ldots,\ \left(
i,\lambda_{i}\right)  \right\}  \\\text{(since the cells in the }i\text{-th
row of }Y\left(  \lambda\right)  \text{ are }\left(  i,1\right)  ,\ \left(
i,2\right)  ,\ \ldots,\ \left(  i,\lambda_{i}\right)  \text{)}}}\\
&  \ \ \ \ \ \ \ \ \ \ \ \ \ \ \ \ \ \ \ \ \left(  \text{since }c=\left(
i,j\right)  \right) \\
&  =\left\{  \left(  i,j\right)  \right\}  \cup\left\{  \left(  i,j+1\right)
,\ \left(  i,j+2\right)  ,\ \ldots,\ \left(  i,\lambda_{i}\right)  \right\} \\
&  =\left\{  \left(  i,j\right)  ,\ \left(  i,j+1\right)  ,\ \left(
i,j+2\right)  ,\ \ldots,\ \left(  i,\lambda_{i}\right)  \right\}  .
\end{align*}
Hence,%
\[
\left\vert \operatorname*{Arm}\nolimits_{\lambda}^{+}\left(  c\right)
\right\vert =\left\vert \left\{  \left(  i,j\right)  ,\ \left(  i,j+1\right)
,\ \left(  i,j+2\right)  ,\ \ldots,\ \left(  i,\lambda_{i}\right)  \right\}
\right\vert =\lambda_{i}-j+1
\]
(since the $\lambda_{i}-j+1$ cells $\left(  i,j\right)  ,\ \left(
i,j+1\right)  ,\ \left(  i,j+2\right)  ,\ \ldots,\ \left(  i,\lambda
_{i}\right)  $ are obviously distinct). This proves Proposition
\ref{prop.hlf-young.arm-leg.1} \textbf{(a)}. \medskip

\textbf{(b)} The cells in the $j$-th column of $Y\left(  \lambda\right)  $ are
$\left(  1,j\right)  ,\ \left(  2,j\right)  ,\ \ldots,\ \left(  \lambda
_{j}^{t},j\right)  $ (by Theorem \ref{thm.partitions.conj} \textbf{(e)}). The
definition of $\operatorname*{Leg}\nolimits_{\lambda}^{-}\left(  c\right)  $
yields%
\begin{align*}
\operatorname*{Leg}\nolimits_{\lambda}^{-}\left(  c\right)   &  =\left\{
\text{all cells in }Y\left(  \lambda\right)  \text{ that lie south of }c\text{
(in the same column as }c\text{)}\right\} \\
&  =\left\{  \text{all cells in }Y\left(  \lambda\right)  \text{ that lie
south of }\left(  i,j\right)  \text{ (in the same column as }\left(
i,j\right)  \text{)}\right\} \\
&  \ \ \ \ \ \ \ \ \ \ \ \ \ \ \ \ \ \ \ \ \left(  \text{since }c=\left(
i,j\right)  \right) \\
&  =\left\{  \text{all cells in }Y\left(  \lambda\right)  \text{ that lie
south of }\left(  i,j\right)  \text{ (in the }j\text{-th column)}\right\} \\
&  =\left\{  \text{all cells in the }j\text{-th column of }Y\left(
\lambda\right)  \text{ that lie south of }\left(  i,j\right)  \right\} \\
&  =\left\{  \left(  i+1,j\right)  ,\ \left(  i+2,j\right)  ,\ \ldots
,\ \left(  \lambda_{j}^{t},j\right)  \right\}
\end{align*}
(since the cells in the $j$-th column of $Y\left(  \lambda\right)  $ are
$\left(  1,j\right)  ,\ \left(  2,j\right)  ,\ \ldots,\ \left(  \lambda
_{j}^{t},j\right)  $). Hence,%
\[
\left\vert \operatorname*{Leg}\nolimits_{\lambda}^{-}\left(  c\right)
\right\vert =\left\vert \left\{  \left(  i+1,j\right)  ,\ \left(
i+2,j\right)  ,\ \ldots,\ \left(  \lambda_{j}^{t},j\right)  \right\}
\right\vert =\lambda_{j}^{t}-i
\]
(since the $\lambda_{j}^{t}-i$ cells $\left(  i+1,j\right)  ,\ \left(
i+2,j\right)  ,\ \ldots,\ \left(  \lambda_{j}^{t},j\right)  $ are obviously
distinct). This proves Proposition \ref{prop.hlf-young.arm-leg.1}
\textbf{(b)}. \medskip

\textbf{(c)} The set $\operatorname*{Arm}\nolimits_{\lambda}^{+}\left(
c\right)  $ consists entirely of cells that lie in the same row as $c$ (by its
definition), whereas the set $\operatorname*{Leg}\nolimits_{\lambda}%
^{-}\left(  c\right)  $ consists entirely of cells that lie in rows further
south than $c$ (again by its definition). Thus, there is no cell that belongs
to both of these two sets. In other words, the sets $\operatorname*{Arm}%
\nolimits_{\lambda}^{+}\left(  c\right)  $ and $\operatorname*{Leg}%
\nolimits_{\lambda}^{-}\left(  c\right)  $ are disjoint.

However, the definition of $H_{\lambda}\left(  c\right)  $ in Theorem
\ref{thm.tableau.hlf} says that%
\begin{align*}
H_{\lambda}\left(  c\right)   &  =\underbrace{\left\{  c\right\}  \cup\left\{
\text{all cells in }Y\left(  \lambda\right)  \text{ that lie east of }c\text{
(in the same row as }c\text{)}\right\}  }_{\substack{=\operatorname*{Arm}%
\nolimits_{\lambda}^{+}\left(  c\right)  \\\text{(by the definition of
}\operatorname*{Arm}\nolimits_{\lambda}^{+}\left(  c\right)  \text{)}}}\\
&  \ \ \ \ \ \ \ \ \ \ \cup\underbrace{\left\{  \text{all cells in }Y\left(
\lambda\right)  \text{ that lie south of }c\text{ (in the same column as
}c\text{)}\right\}  }_{\substack{=\operatorname*{Leg}\nolimits_{\lambda}%
^{-}\left(  c\right)  \\\text{(by the definition of }\operatorname*{Leg}%
\nolimits_{\lambda}^{-}\left(  c\right)  \text{)}}}\\
&  =\operatorname*{Arm}\nolimits_{\lambda}^{+}\left(  c\right)  \cup
\operatorname*{Leg}\nolimits_{\lambda}^{-}\left(  c\right)  .
\end{align*}
Thus,%
\[
\left\vert H_{\lambda}\left(  c\right)  \right\vert =\left\vert
\operatorname*{Arm}\nolimits_{\lambda}^{+}\left(  c\right)  \cup
\operatorname*{Leg}\nolimits_{\lambda}^{-}\left(  c\right)  \right\vert
=\left\vert \operatorname*{Arm}\nolimits_{\lambda}^{+}\left(  c\right)
\right\vert +\left\vert \operatorname*{Leg}\nolimits_{\lambda}^{-}\left(
c\right)  \right\vert
\]
(since the sets $\operatorname*{Arm}\nolimits_{\lambda}^{+}\left(  c\right)  $
and $\operatorname*{Leg}\nolimits_{\lambda}^{-}\left(  c\right)  $ are
disjoint). This proves Proposition \ref{prop.hlf-young.arm-leg.1} \textbf{(c)}.
\end{proof}
\end{fineprint}

The following property of long arms will also be useful:

\begin{proposition}
\label{prop.hlf-young.arm-prod}Let $\lambda$ be a partition of $n$. For any
cell $c\in Y\left(  \lambda\right)  $, define the long arm
$\operatorname*{Arm}\nolimits_{\lambda}^{+}\left(  c\right)  $ as in
Proposition \ref{prop.hlf-young.arm-leg.1}. Then: \medskip

\textbf{(a)} If $\lambda=\left(  \lambda_{1},\lambda_{2},\ldots,\lambda_{\ell
}\right)  $ for some $\ell\in\mathbb{N}$, then%
\[
\prod_{c\in Y\left(  \lambda\right)  }\left\vert \operatorname*{Arm}%
\nolimits_{\lambda}^{+}\left(  c\right)  \right\vert =\prod_{i=1}^{\ell
}\lambda_{i}!.
\]

\textbf{(b)} If $T$ is any $n$-tableau of shape $\lambda$, then%
\[
\prod_{c\in Y\left(  \lambda\right)  }\left\vert \operatorname*{Arm}%
\nolimits_{\lambda}^{+}\left(  c\right)  \right\vert =\left\vert
\mathcal{R}\left(  T\right)  \right\vert .
\]

\end{proposition}

\begin{fineprint}
\begin{proof}
\textbf{(a)} Assume that $\lambda=\left(  \lambda_{1},\lambda_{2}%
,\ldots,\lambda_{\ell}\right)  $ for some $\ell\in\mathbb{N}$. Then, the
elements of the Young diagram $Y\left(  \lambda\right)  $ are precisely the
pairs $\left(  i,j\right)  $ with $i\in\left[  \ell\right]  $ and $j\in\left[
\lambda_{i}\right]  $. Hence, the product sign $\prod_{\left(  i,j\right)  \in
Y\left(  \lambda\right)  }$ can be rewritten as $\prod_{i\in\left[
\ell\right]  }\ \ \prod_{j\in\left[  \lambda_{i}\right]  }$.

However,
\begin{align}
\prod_{c\in Y\left(  \lambda\right)  }\left\vert \operatorname*{Arm}%
\nolimits_{\lambda}^{+}\left(  c\right)  \right\vert  &  =\prod_{\left(
i,j\right)  \in Y\left(  \lambda\right)  }\underbrace{\left\vert
\operatorname*{Arm}\nolimits_{\lambda}^{+}\left(  i,j\right)  \right\vert
}_{\substack{=\lambda_{i}-j+1\\\text{(by Proposition
\ref{prop.hlf-young.arm-leg.1} \textbf{(a)},}\\\text{applied to }c=\left(
i,j\right)  \text{)}}}\ \ \ \ \ \ \ \ \ \ \left(
\begin{array}
[c]{c}%
\text{here, we have renamed}\\
\text{the index }c\text{ as }\left(  i,j\right)
\end{array}
\right) \nonumber\\
&  =\prod_{\left(  i,j\right)  \in Y\left(  \lambda\right)  }\left(
\lambda_{i}-j+1\right) \nonumber\\
&  =\prod_{i\in\left[  \ell\right]  }\ \ \prod_{j\in\left[  \lambda
_{i}\right]  }\left(  \lambda_{i}-j+1\right)
\label{pf.prop.hlf-young.arm-prod.a.1}%
\end{align}
(since the product sign $\prod_{\left(  i,j\right)  \in Y\left(
\lambda\right)  }$ can be rewritten as $\prod_{i\in\left[  \ell\right]
}\ \ \prod_{j\in\left[  \lambda_{i}\right]  }$). However, for each
$i\in\left[  \ell\right]  $, we have%
\begin{align*}
\underbrace{\prod_{j\in\left[  \lambda_{i}\right]  }}_{=\prod_{j=1}%
^{\lambda_{i}}}\underbrace{\left(  \lambda_{i}-j+1\right)  }_{=\lambda
_{i}+1-j}  &  =\prod_{j=1}^{\lambda_{i}}\left(  \lambda_{i}+1-j\right)
=\prod_{k=1}^{\lambda_{i}}k\ \ \ \ \ \ \ \ \ \ \left(
\begin{array}
[c]{c}%
\text{here, we have substituted }k\\
\text{for }\lambda_{i}+1-j\text{ in the product}%
\end{array}
\right) \\
&  =1\cdot2\cdot\cdots\cdot\lambda_{i}=\lambda_{i}!.
\end{align*}
Thus, we can rewrite (\ref{pf.prop.hlf-young.arm-prod.a.1}) as
\[
\prod_{c\in Y\left(  \lambda\right)  }\left\vert \operatorname*{Arm}%
\nolimits_{\lambda}^{+}\left(  c\right)  \right\vert =\underbrace{\prod
_{i\in\left[  \ell\right]  }}_{=\prod_{i=1}^{\ell}}\lambda_{i}!=\prod
_{i=1}^{\ell}\lambda_{i}!.
\]
This proves Proposition \ref{prop.hlf-young.arm-prod} \textbf{(a)}. \medskip

\textbf{(b)} Let $D$ be the diagram $Y\left(  \lambda\right)  $. Let $T$ be
any $n$-tableau of shape $\lambda$. Thus, $T$ is an $n$-tableau of shape
$Y\left(  \lambda\right)  =D$.

For any $i\in\mathbb{Z}$, let $\operatorname*{Row}\left(  i,T\right)  $ denote
the set of all entries in the $i$-th row of $T$. Then, for each positive
integer $i$, we have%
\begin{align}
\left\vert \operatorname*{Row}\left(  i,T\right)  \right\vert  &  =\left\vert
D\left\lfloor i\right\rfloor \right\vert \ \ \ \ \ \ \ \ \ \ \left(  \text{by
Proposition \ref{prop.tableau.RowTj}, applied to }j=i\right) \nonumber\\
&  =\left\vert \left[  \lambda_{i}\right]  \right\vert
\ \ \ \ \ \ \ \ \ \ \left(  \text{since Lemma \ref{lem.hlf-young.grel-lamd}
yields }D\left\lfloor i\right\rfloor =\left[  \lambda_{i}\right]  \right)
\nonumber\\
&  =\lambda_{i}. \label{pf.prop.hlf-young.arm-prod.b.1}%
\end{align}

Let $\ell$ be the length of $\lambda$. Thus, $\lambda=\left(  \lambda
_{1},\lambda_{2},\ldots,\lambda_{\ell}\right)  $. Hence, all cells of the
Young diagram $Y\left(  \lambda\right)  $ lie in rows $1,2,\ldots,\ell$.
Therefore, Proposition \ref{prop.tableau.RTsize} \textbf{(a)} (applied to
$D=Y\left(  \lambda\right)  $) yields%
\[
\left\vert \mathcal{R}\left(  T\right)  \right\vert =\prod_{i=1}^{\ell
}\underbrace{\left\vert \operatorname*{Row}\left(  i,T\right)  \right\vert
}_{\substack{=\lambda_{i}\\\text{(by (\ref{pf.prop.hlf-young.arm-prod.b.1}))}%
}}!=\prod_{i=1}^{\ell}\lambda_{i}!=\prod_{c\in Y\left(  \lambda\right)
}\left\vert \operatorname*{Arm}\nolimits_{\lambda}^{+}\left(  c\right)
\right\vert
\]
(by Proposition \ref{prop.hlf-young.arm-prod} \textbf{(a)}). This proves
Proposition \ref{prop.hlf-young.arm-prod} \textbf{(b)}.
\end{proof}
\end{fineprint}

\subsubsection{A few last lemmas}

To simplify the proof of Theorem \ref{thm.tableau.hlf}, we outsource a few
more arguments to lemmas. The following two lemmas are simple properties of
row symmetrizers:

\begin{lemma}
\label{lem.hlf-young.row1}Let $T$ be an $n$-tableau of any shape $D$. Let
$p,q,r$ be three distinct elements of $\left[  n\right]  $. Assume that $p$
and $r$ lie in the same row of $T$. Then,%
\[
\left(  1-t_{p,q}\right)  t_{q,r}\nabla_{\operatorname*{Row}T}=0.
\]

\end{lemma}

\begin{proof}
We know that $p$ and $r$ lie in the same row of $T$. Thus, Proposition
\ref{prop.symmetrizers.factor-out-row} \textbf{(c)} (applied to $i=p$ and
$j=r$) shows that%
\[
t_{p,r}\nabla_{\operatorname*{Row}T}=\nabla_{\operatorname*{Row}T}%
t_{p,r}=\nabla_{\operatorname*{Row}T}.
\]

On the other hand, we have $p\notin\left\{  q,r\right\}  $ (since $p,q,r$ are
distinct) and thus $t_{q,r}\left(  p\right)  =p$. Applying
(\ref{eq.intro.perms.cycs.sts-1}) to $\sigma=t_{q,r}$ and $i=p$ and $j=r$, we
find
\[
t_{q,r}\circ t_{p,r}\circ t_{q,r}^{-1}=t_{t_{q,r}\left(  p\right)
,t_{q,r}\left(  r\right)  }=t_{p,q}\ \ \ \ \ \ \ \ \ \ \left(  \text{since
}t_{q,r}\left(  p\right)  =p\text{ and }t_{q,r}\left(  r\right)  =q\right)  .
\]
Thus, $t_{p,q}=t_{q,r}\circ t_{p,r}\circ t_{q,r}^{-1}=t_{q,r}t_{p,r}%
t_{q,r}^{-1}$. Multiplying this equality by $t_{q,r}$ from the right, we
obtain
\[
t_{p,q}t_{q,r}=t_{q,r}t_{p,r}\underbrace{t_{q,r}^{-1}t_{q,r}}_{=1}%
=t_{q,r}t_{p,r}.
\]
Now,%
\[
\left(  1-t_{p,q}\right)  t_{q,r}=t_{q,r}-\underbrace{t_{p,q}t_{q,r}%
}_{=t_{q,r}t_{p,r}}=t_{q,r}-t_{q,r}t_{p,r}=t_{q,r}\left(  1-t_{p,r}\right)  .
\]
Hence,%
\[
\underbrace{\left(  1-t_{p,q}\right)  t_{q,r}}_{=t_{q,r}\left(  1-t_{p,r}%
\right)  }\nabla_{\operatorname*{Row}T}=t_{q,r}\underbrace{\left(
1-t_{p,r}\right)  \nabla_{\operatorname*{Row}T}}_{\substack{=\nabla
_{\operatorname*{Row}T}-t_{p,r}\nabla_{\operatorname*{Row}T}\\=0\\\text{(since
}t_{p,r}\nabla_{\operatorname*{Row}T}=\nabla_{\operatorname*{Row}T}\text{)}%
}}=0.
\]
This proves Lemma \ref{lem.hlf-young.row1}.
\end{proof}

\begin{lemma}
\label{lem.hlf-young.row2}Let $T$ be an $n$-tableau of any shape $D$. Let
$p,q,r$ be three distinct elements of $\left[  n\right]  $. Assume that $p$
and $r$ lie in the same row of $T$. Let $\mathbf{a}\in\mathcal{A}$ be an
element that commutes with $t_{p,r}$. Then,%
\[
\nabla_{\operatorname*{Row}T}\mathbf{a}t_{r,q}\nabla_{\operatorname*{Row}%
T}=\nabla_{\operatorname*{Row}T}\mathbf{a}t_{p,q}\nabla_{\operatorname*{Row}%
T}.
\]

\end{lemma}

\begin{proof}
As in the proof of Lemma \ref{lem.hlf-young.row1}, we can see that
\begin{equation}
t_{p,r}\nabla_{\operatorname*{Row}T}=\nabla_{\operatorname*{Row}T}%
t_{p,r}=\nabla_{\operatorname*{Row}T}. \label{pf.lem.hlf-young.row2.1}%
\end{equation}

Moreover, (\ref{eq.intro.perms.cycs.c=ttt}) (applied to $k=3$ and $\left(
i_{1},i_{2},\ldots,i_{k}\right)  =\left(  p,r,q\right)  $) yields
$\operatorname*{cyc}\nolimits_{p,r,q}=t_{p,r}\circ t_{r,q}$. Meanwhile,
(\ref{eq.intro.perms.cycs.c=ttt}) (applied to $k=3$ and $\left(  i_{1}%
,i_{2},\ldots,i_{k}\right)  =\left(  q,p,r\right)  $) yields
$\operatorname*{cyc}\nolimits_{q,p,r}=t_{q,p}\circ t_{p,r}$. But we have
$\operatorname*{cyc}\nolimits_{p,r,q}=\operatorname*{cyc}\nolimits_{q,p,r}$.
Thus,
\[
t_{p,r}\circ t_{r,q}=\operatorname*{cyc}\nolimits_{p,r,q}=\operatorname*{cyc}%
\nolimits_{q,p,r}=t_{q,p}\circ t_{p,r}.
\]

Now,%
\begin{align*}
\underbrace{\nabla_{\operatorname*{Row}T}}_{\substack{=\nabla
_{\operatorname*{Row}T}t_{p,r}\\\text{(by (\ref{pf.lem.hlf-young.row2.1}))}%
}}\mathbf{a}t_{r,q}\nabla_{\operatorname*{Row}T}  &  =\nabla
_{\operatorname*{Row}T}\underbrace{t_{p,r}\mathbf{a}}_{\substack{=\mathbf{a}%
t_{p,r}\\\text{(since }\mathbf{a}\text{ commutes}\\\text{with }t_{p,r}%
\text{)}}}t_{r,q}\nabla_{\operatorname*{Row}T}\\
&  =\nabla_{\operatorname*{Row}T}\mathbf{a}\underbrace{t_{p,r}t_{r,q}%
}_{\substack{=t_{p,r}\circ t_{r,q}\\=t_{q,p}\circ t_{p,r}\\=t_{q,p}t_{p,r}%
}}\nabla_{\operatorname*{Row}T}=\nabla_{\operatorname*{Row}T}\mathbf{a}%
\underbrace{t_{q,p}}_{=t_{p,q}}\underbrace{t_{p,r}\nabla_{\operatorname*{Row}%
T}}_{\substack{=\nabla_{\operatorname*{Row}T}\\\text{(by
(\ref{pf.lem.hlf-young.row2.1}))}}}\\
&  =\nabla_{\operatorname*{Row}T}\mathbf{a}t_{p,q}\nabla_{\operatorname*{Row}%
T}.
\end{align*}
This proves Lemma \ref{lem.hlf-young.row2}.
\end{proof}

The following lemma is again an easy consequence of things proved before. Its
main purpose is to isolate a simple argument for the proof of the hook length formula:

\begin{lemma}
\label{lem.hlf-young.last}Let $\mathbf{k}=\mathbb{Q}$. Let $\lambda$ be a
partition of $n$. Let $T$ be an $n$-tableau of shape $Y\left(  \lambda\right)
$. Set $P:=\nabla_{\operatorname*{Row}T}$ and $N:=\nabla_{\operatorname*{Col}%
T}^{-}$. Let $\kappa\in\mathbb{Q}$ be such that $\left(  PN\right)
^{2}=\kappa\cdot PN$. Then, the \# of standard tableaux of shape $Y\left(
\lambda\right)  $ is $\dfrac{n!}{\kappa}$.
\end{lemma}

\begin{proof}
Multiplying the equalities $N=\nabla_{\operatorname*{Col}T}^{-}$ and
$P=\nabla_{\operatorname*{Row}T}$, we find%
\[
NP=\nabla_{\operatorname*{Col}T}^{-}\nabla_{\operatorname*{Row}T}%
=\mathbf{E}_{T}\ \ \ \ \ \ \ \ \ \ \left(  \text{by Definition
\ref{def.specht.ET.defs} \textbf{(b)}}\right)  .
\]

Lemma \ref{lem.specht.ET.1-coord} yields $\left[  1\right]  \mathbf{E}_{T}=1$
(see Definition \ref{def.specht.ET.wa} for the meaning of $\left[  1\right]
\mathbf{E}_{T}$).

Let $f^{\lambda}$ be the \# of standard tableaux of shape $Y\left(
\lambda\right)  $. Then, Theorem \ref{thm.specht.ETidp} yields
\begin{equation}
\mathbf{E}_{T}^{2}=\dfrac{n!}{f^{\lambda}}\mathbf{E}_{T}.
\label{pf.lem.hlf-young.last.1ET1}%
\end{equation}

Recall the antipode $S$ of $\mathbf{k}\left[  S_{n}\right]  $ (as defined in
Definition \ref{def.S.S}). Then, Proposition \ref{prop.symmetrizers.antipode}
says that $S\left(  \nabla_{\operatorname*{Row}T}\right)  =\nabla
_{\operatorname*{Row}T}$ and $S\left(  \nabla_{\operatorname*{Col}T}%
^{-}\right)  =\nabla_{\operatorname*{Col}T}^{-}$. In view of $P=\nabla
_{\operatorname*{Row}T}$ and $N=\nabla_{\operatorname*{Col}T}^{-}$, we can
rewrite these two equalities as $S\left(  P\right)  =P$ and $S\left(
N\right)  =N$.

We have $PNPN=\left(  PN\right)  ^{2}=\kappa\cdot PN$. Applying the map $S$ to
both sides of this equality, we find%
\begin{equation}
S\left(  PNPN\right)  =S\left(  \kappa\cdot PN\right)  =\kappa\cdot S\left(
PN\right)  \label{pf.lem.hlf-young.last.1}%
\end{equation}
(since the map $S$ is $\mathbf{k}$-linear).

But $S$ is a $\mathbf{k}$-algebra anti-automorphism (by Theorem
\ref{thm.S.auto} \textbf{(a)}). Thus,
\[
S\left(  PN\right)  =\underbrace{S\left(  N\right)  }_{=N}\cdot
\underbrace{S\left(  P\right)  }_{=P}=NP=\mathbf{E}_{T}%
\]
and%
\[
S\left(  PNPN\right)  =\underbrace{S\left(  N\right)  }_{=N}\cdot
\underbrace{S\left(  P\right)  }_{=P}\cdot\underbrace{S\left(  N\right)
}_{=N}\cdot\underbrace{S\left(  P\right)  }_{=P}=\underbrace{NP}%
_{=\mathbf{E}_{T}}\underbrace{NP}_{=\mathbf{E}_{T}}=\mathbf{E}_{T}%
\mathbf{E}_{T}=\mathbf{E}_{T}^{2}=\dfrac{n!}{f^{\lambda}}\mathbf{E}_{T}%
\]
(by (\ref{pf.lem.hlf-young.last.1ET1})). Using these two equalities, we can
rewrite (\ref{pf.lem.hlf-young.last.1}) as%
\[
\dfrac{n!}{f^{\lambda}}\mathbf{E}_{T}=\kappa\cdot\mathbf{E}_{T}.
\]
Thus,
\[
\left[  1\right]  \left(  \dfrac{n!}{f^{\lambda}}\mathbf{E}_{T}\right)
=\left[  1\right]  \left(  \kappa\cdot\mathbf{E}_{T}\right)  =\kappa
\cdot\underbrace{\left[  1\right]  \mathbf{E}_{T}}_{=1}=\kappa.
\]
Hence,%
\[
\kappa=\left[  1\right]  \left(  \dfrac{n!}{f^{\lambda}}\mathbf{E}_{T}\right)
=\dfrac{n!}{f^{\lambda}}\underbrace{\left[  1\right]  \mathbf{E}_{T}}%
_{=1}=\dfrac{n!}{f^{\lambda}}.
\]
Hence, $f^{\lambda}=\dfrac{n!}{\kappa}$. In other words, the \# of standard
tableaux of shape $Y\left(  \lambda\right)  $ is $\dfrac{n!}{\kappa}$ (since
$f^{\lambda}$ is the \# of standard tableaux of shape $Y\left(  \lambda
\right)  $). This proves Lemma \ref{lem.hlf-young.last}.
\end{proof}

\subsubsection{Proof of the hook length formula}

We now have all the tools in place to prove the hook length formula:

\begin{proof}
[Proof of Theorem \ref{thm.tableau.hlf}.]Let $\mathbf{k}=\mathbb{Q}$.

We have $\left\vert Y\left(  \lambda\right)  \right\vert =\left\vert
\lambda\right\vert =n$. In other words, the diagram $Y\left(  \lambda\right)
$ consists of $n$ cells. Let us list these $n$ cells row by row, from the
bottommost row to the topmost row, where the cells of each row are listed from
left to right:%
\begin{align}
&  \left(  \ell,1\right)  ,\ \left(  \ell,2\right)  ,\ \ldots,\ \left(
\ell,\lambda_{\ell}\right)  ,\nonumber\\
&  \left(  \ell-1,1\right)  ,\ \left(  \ell-1,2\right)  ,\ \ldots,\ \left(
\ell-1,\lambda_{\ell-1}\right)  ,\nonumber\\
&  \ldots,\nonumber\\
&  \left(  1,1\right)  ,\ \left(  1,2\right)  ,\ \ldots,\ \left(
1,\lambda_{1}\right)  \label{pf.thm.tableau.hlf.list}%
\end{align}
(where $\ell$ is the length of $\lambda$). Now, we define an $n$-tableau $T$
of shape $Y\left(  \lambda\right)  $ by filling these $n$ cells (in the order
in which we just listed them) with the entries $1,2,\ldots,n$ (respectively).

\begin{example}
\label{exa.pf.thm.tableau.hlf.T}If $n=14$ and $\lambda=\left(  5,4,3,2\right)
$, then%
\[
T=\ytableaushort{{10}{11}{12}{13}{14},6789,345,12}\ \ .
\]

\end{example}

Thus we have defined an $n$-tableau $T$ of shape $Y\left(  \lambda\right)  $,
that is, an $n$-tableau $T$ of shape $\lambda$. In particular, all entries of
$T$ are distinct (since $T$ is an $n$-tableau, thus injective).

Let $P=\nabla_{\operatorname*{Row}T}$ and $N=\nabla_{\operatorname*{Col}T}%
^{-}$. Then, Corollary \ref{cor.symmetrizers.square} \textbf{(a)} says that
$\left(  \nabla_{\operatorname*{Row}T}\right)  ^{2}=\left\vert \mathcal{R}%
\left(  T\right)  \right\vert \cdot\nabla_{\operatorname*{Row}T}$. In other
words, $P^{2}=\left\vert \mathcal{R}\left(  T\right)  \right\vert \cdot P$
(since $P=\nabla_{\operatorname*{Row}T}$).

For any $j\in\mathbb{Z}$, let $\operatorname*{Col}\left(  j,T\right)  $ denote
the set of all entries in the $j$-th column of $T$. Choose some $k\in
\mathbb{N}$ such that all cells of $Y\left(  \lambda\right)  $ lie in columns
$1,2,\ldots,k$. (Such a $k$ clearly exists, since $Y\left(  \lambda\right)  $
is a Young diagram; in fact we can take $k=\lambda_{1}$.)

Proposition \ref{prop.symmetrizers.int} \textbf{(b)} (applied to $D=Y\left(
\lambda\right)  $) shows that
\begin{equation}
\nabla_{\operatorname*{Col}T}^{-}=\prod_{j=1}^{k}\nabla_{\operatorname*{Col}%
\left(  j,T\right)  }^{-}. \label{pf.thm.tableau.hlf.Nform}%
\end{equation}
Here, we are not specifying the order of the factors in the product, since all
the factors commute.

Now, for each $m\in\left\{  0,1,\ldots,n\right\}  $, we define an element
$N_{m}\in\mathcal{A}$ by%
\begin{equation}
N_{m}:=\prod_{j=1}^{k}\nabla_{\left[  m\right]  \cap\operatorname*{Col}\left(
j,T\right)  }^{-}. \label{pf.thm.tableau.hlf.Nm=}%
\end{equation}
Here, again, we are not specifying the order of the factors in the product,
since all the factors commute\footnote{\textit{Proof.} Let $m\in\left\{
0,1,\ldots,n\right\}  $. We must prove that all factors in the product on the
right hand side of (\ref{pf.thm.tableau.hlf.Nm=}) commute. In other words, we
must prove that the elements $\nabla_{\left[  m\right]  \cap
\operatorname*{Col}\left(  j,T\right)  }^{-}$ for all $j\in\left[  k\right]  $
commute. In other words, we must prove that $\nabla_{\left[  m\right]
\cap\operatorname*{Col}\left(  u,T\right)  }^{-}$ commutes with $\nabla
_{\left[  m\right]  \cap\operatorname*{Col}\left(  v,T\right)  }^{-}$ for any
two elements $u,v\in\left[  k\right]  $. Let us prove this.
\par
Let $u,v\in\left[  k\right]  $ be two elements. We must prove that
$\nabla_{\left[  m\right]  \cap\operatorname*{Col}\left(  u,T\right)  }^{-}$
commutes with $\nabla_{\left[  m\right]  \cap\operatorname*{Col}\left(
v,T\right)  }^{-}$. If $u=v$, then this is obvious. Thus, we WLOG assume that
$u\neq v$.
\par
Recall that all entries of $T$ are distinct. Thus, in particular, an entry in
the $u$-th column of $T$ cannot also appear in the $v$-th column of $T$ (since
$u\neq v$). In other words, a number cannot belong to $\operatorname*{Col}%
\left(  u,T\right)  $ and to $\operatorname*{Col}\left(  v,T\right)  $ at the
same time. In other words, $\operatorname*{Col}\left(  u,T\right)
\cap\operatorname*{Col}\left(  v,T\right)  =\varnothing$. Hence,
\begin{align*}
\left(  \left[  m\right]  \cap\operatorname*{Col}\left(  u,T\right)  \right)
\cap\left(  \left[  m\right]  \cap\operatorname*{Col}\left(  v,T\right)
\right)   &  =\left[  m\right]  \cap\left[  m\right]  \cap
\underbrace{\operatorname*{Col}\left(  u,T\right)  \cap\operatorname*{Col}%
\left(  v,T\right)  }_{=\varnothing}\\
&  =\left[  m\right]  \cap\left[  m\right]  \cap\varnothing=\varnothing.
\end{align*}
Thus, the two subsets $\left[  m\right]  \cap\operatorname*{Col}\left(
u,T\right)  $ and $\left[  m\right]  \cap\operatorname*{Col}\left(
v,T\right)  $ of $\left[  n\right]  $ are disjoint. Hence, Proposition
\ref{prop.int.commute} \textbf{(b)} (applied to $X=\left[  m\right]
\cap\operatorname*{Col}\left(  u,T\right)  $ and $Y=\left[  m\right]
\cap\operatorname*{Col}\left(  v,T\right)  $) yields that
\begin{align*}
\nabla_{\left[  m\right]  \cap\operatorname*{Col}\left(  u,T\right)  }%
\nabla_{\left[  m\right]  \cap\operatorname*{Col}\left(  v,T\right)  }  &
=\nabla_{\left[  m\right]  \cap\operatorname*{Col}\left(  v,T\right)  }%
\nabla_{\left[  m\right]  \cap\operatorname*{Col}\left(  u,T\right)  };\\
\nabla_{\left[  m\right]  \cap\operatorname*{Col}\left(  u,T\right)  }%
\nabla_{\left[  m\right]  \cap\operatorname*{Col}\left(  v,T\right)  }^{-}  &
=\nabla_{\left[  m\right]  \cap\operatorname*{Col}\left(  v,T\right)  }%
^{-}\nabla_{\left[  m\right]  \cap\operatorname*{Col}\left(  u,T\right)  };\\
\nabla_{\left[  m\right]  \cap\operatorname*{Col}\left(  u,T\right)  }%
^{-}\nabla_{\left[  m\right]  \cap\operatorname*{Col}\left(  v,T\right)  }  &
=\nabla_{\left[  m\right]  \cap\operatorname*{Col}\left(  v,T\right)  }%
\nabla_{\left[  m\right]  \cap\operatorname*{Col}\left(  u,T\right)  }^{-};\\
\nabla_{\left[  m\right]  \cap\operatorname*{Col}\left(  u,T\right)  }%
^{-}\nabla_{\left[  m\right]  \cap\operatorname*{Col}\left(  v,T\right)
}^{-}  &  =\nabla_{\left[  m\right]  \cap\operatorname*{Col}\left(
v,T\right)  }^{-}\nabla_{\left[  m\right]  \cap\operatorname*{Col}\left(
u,T\right)  }^{-}.
\end{align*}
The last of these four equalities shows that $\nabla_{\left[  m\right]
\cap\operatorname*{Col}\left(  u,T\right)  }^{-}$ commutes with $\nabla
_{\left[  m\right]  \cap\operatorname*{Col}\left(  v,T\right)  }^{-}$. This
completes our proof.}.

\begin{example}
If $n$, $\lambda$ and $T$ are as in Example \ref{exa.pf.thm.tableau.hlf.T},
then we can pick $k=5$ (since all cells of $Y\left(  \lambda\right)  $ lie in
columns $1,2,3,4,5$) and obtain%
\begin{align*}
N_{7}  &  =\prod_{j=1}^{k}\nabla_{\left[  7\right]  \cap\operatorname*{Col}%
\left(  j,T\right)  }^{-}\\
&  =\underbrace{\nabla_{\left[  7\right]  \cap\operatorname*{Col}\left(
1,T\right)  }^{-}}_{\substack{=\nabla_{\left[  7\right]  \cap\left\{
10,6,3,1\right\}  }^{-}\\=\nabla_{\left\{  6,3,1\right\}  }^{-}}%
}\underbrace{\nabla_{\left[  7\right]  \cap\operatorname*{Col}\left(
2,T\right)  }^{-}}_{\substack{=\nabla_{\left[  7\right]  \cap\left\{
11,7,4,2\right\}  }^{-}\\=\nabla_{\left\{  7,4,2\right\}  }^{-}}%
}\underbrace{\nabla_{\left[  7\right]  \cap\operatorname*{Col}\left(
3,T\right)  }^{-}}_{\substack{=\nabla_{\left[  7\right]  \cap\left\{
12,8,5\right\}  }^{-}\\=\nabla_{\left\{  5\right\}  }^{-}=1}%
}\underbrace{\nabla_{\left[  7\right]  \cap\operatorname*{Col}\left(
4,T\right)  }^{-}}_{\substack{=\nabla_{\left[  7\right]  \cap\left\{
13,9\right\}  }^{-}\\=\nabla_{\varnothing}^{-}=1}}\underbrace{\nabla_{\left[
7\right]  \cap\operatorname*{Col}\left(  5,T\right)  }^{-}}_{\substack{=\nabla
_{\left[  7\right]  \cap\left\{  14\right\}  }^{-}\\=\nabla_{\varnothing}%
^{-}=1}}\\
&  =\nabla_{\left\{  6,3,1\right\}  }^{-}\nabla_{\left\{  7,4,2\right\}  }%
^{-}.
\end{align*}

\end{example}

Note that each $j\in\left[  k\right]  $ satisfies $\left[  n\right]
\cap\operatorname*{Col}\left(  j,T\right)  =\operatorname*{Col}\left(
j,T\right)  $ (since $\operatorname*{Col}\left(  j,T\right)  \subseteq\left[
n\right]  $) and thus%
\begin{equation}
\nabla_{\left[  n\right]  \cap\operatorname*{Col}\left(  j,T\right)  }%
^{-}=\nabla_{\operatorname*{Col}\left(  j,T\right)  }^{-}.
\label{pf.thm.tableau.hlf.Nn.pf.1}%
\end{equation}
The definition of $N_{n}$ yields%
\begin{align}
N_{n}  &  =\prod_{j=1}^{k}\underbrace{\nabla_{\left[  n\right]  \cap
\operatorname*{Col}\left(  j,T\right)  }^{-}}_{\substack{=\nabla
_{\operatorname*{Col}\left(  j,T\right)  }^{-}\\\text{(by
(\ref{pf.thm.tableau.hlf.Nn.pf.1}))}}}=\prod_{j=1}^{k}\nabla
_{\operatorname*{Col}\left(  j,T\right)  }^{-}=\nabla_{\operatorname*{Col}%
T}^{-}\ \ \ \ \ \ \ \ \ \ \left(  \text{by (\ref{pf.thm.tableau.hlf.Nform}%
)}\right) \nonumber\\
&  =N\ \ \ \ \ \ \ \ \ \ \left(  \text{since }N=\nabla_{\operatorname*{Col}%
T}^{-}\right)  . \label{pf.thm.tableau.hlf.Nn=}%
\end{align}

Furthermore, each $j\in\left[  k\right]  $ satisfies $\underbrace{\left[
0\right]  }_{=\varnothing}\cap\operatorname*{Col}\left(  j,T\right)
=\varnothing\cap\operatorname*{Col}\left(  j,T\right)  =\varnothing$ and thus%
\begin{equation}
\nabla_{\left[  0\right]  \cap\operatorname*{Col}\left(  j,T\right)  }%
^{-}=\nabla_{\varnothing}^{-}=1 \label{pf.thm.tableau.hlf.N0.pf.1}%
\end{equation}
(by Example \ref{exa.intX.X=0123} \textbf{(a)}). The definition of $N_{0}$
yields%
\begin{equation}
N_{0}=\prod_{j=1}^{k}\underbrace{\nabla_{\left[  0\right]  \cap
\operatorname*{Col}\left(  j,T\right)  }^{-}}_{\substack{=1\\\text{(by
(\ref{pf.thm.tableau.hlf.N0.pf.1}))}}}=\prod_{j=1}^{k}1=1.
\label{pf.thm.tableau.hlf.N0=}%
\end{equation}

\Needspace{8pc}

\begin{remark}
\label{exa.pf.thm.tableau.hlf.Tk}Here is another way to view $N_{m}$: Let
$m\in\left\{  0,1,\ldots,n\right\}  $. Let $T_{m}$ be the $m$-tableau that is
obtained from $T$ by removing the cells filled with $m+1,m+2,\ldots,n$ (thus
leaving only the cells filled with $1,2,\ldots,m$ in the tableau). For
instance, if $n$, $\lambda$ and $T$ are as in Example
\ref{exa.pf.thm.tableau.hlf.T}, then%
\[
T_{7}=\ytableaushort{\none,67,345,12}\ \ .
\]
Now, $N_{m}$ is the column antisymmetrizer $\nabla_{\operatorname*{Col}\left(
T_{m}\right)  }^{-}\in\mathbf{k}\left[  S_{m}\right]  $, mapped into
$\mathbf{k}\left[  S_{n}\right]  $ using the default embedding $\left(
I_{m\rightarrow n}\right)  _{\ast}:\mathbf{k}\left[  S_{m}\right]
\rightarrow\mathbf{k}\left[  S_{n}\right]  $ (see Corollary
\ref{cor.default-ex.injmor}).
\end{remark}

We shall use the notations $\operatorname*{Arm}\nolimits_{\lambda}^{+}\left(
c\right)  $ and $\operatorname*{Leg}\nolimits_{\lambda}^{-}\left(  c\right)  $
as defined in Proposition \ref{prop.hlf-young.arm-leg.1}.

Furthermore, for each $m\in\left\{  0,1,\ldots,n\right\}  $, we define a
rational number\footnote{As usual, $T^{-1}\left(  \left[  m\right]  \right)  $
denotes the set of all $c\in Y\left(  \lambda\right)  $ that satisfy $T\left(
c\right)  \in\left[  m\right]  $.}%
\begin{equation}
a_{m}:=\prod_{c\in T^{-1}\left(  \left[  m\right]  \right)  }\dfrac{\left\vert
H_{\lambda}\left(  c\right)  \right\vert }{\left\vert \operatorname*{Arm}%
\nolimits_{\lambda}^{+}\left(  c\right)  \right\vert }.
\label{pf.thm.tableau.hlf.am=}%
\end{equation}
Thus,
\begin{align*}
a_{0}  &  =\prod_{c\in T^{-1}\left(  \left[  0\right]  \right)  }%
\dfrac{\left\vert H_{\lambda}\left(  c\right)  \right\vert }{\left\vert
\operatorname*{Arm}\nolimits_{\lambda}^{+}\left(  c\right)  \right\vert
}=\left(  \text{empty product}\right)  \ \ \ \ \ \ \ \ \ \ \left(  \text{since
}T^{-1}\left(  \left[  0\right]  \right)  =\varnothing\right) \\
&  =1
\end{align*}
and%
\begin{align*}
a_{n}  &  =\prod_{c\in T^{-1}\left(  \left[  n\right]  \right)  }%
\dfrac{\left\vert H_{\lambda}\left(  c\right)  \right\vert }{\left\vert
\operatorname*{Arm}\nolimits_{\lambda}^{+}\left(  c\right)  \right\vert
}=\prod_{c\in Y\left(  \lambda\right)  }\dfrac{\left\vert H_{\lambda}\left(
c\right)  \right\vert }{\left\vert \operatorname*{Arm}\nolimits_{\lambda}%
^{+}\left(  c\right)  \right\vert }\ \ \ \ \ \ \ \ \ \ \left(  \text{since
}T^{-1}\left(  \left[  n\right]  \right)  =Y\left(  \lambda\right)  \right) \\
&  =\dfrac{\prod_{c\in Y\left(  \lambda\right)  }\left\vert H_{\lambda}\left(
c\right)  \right\vert }{\prod_{c\in Y\left(  \lambda\right)  }\left\vert
\operatorname*{Arm}\nolimits_{\lambda}^{+}\left(  c\right)  \right\vert
}=\dfrac{\prod_{c\in Y\left(  \lambda\right)  }\left\vert H_{\lambda}\left(
c\right)  \right\vert }{\left\vert \mathcal{R}\left(  T\right)  \right\vert }%
\end{align*}
(since Proposition \ref{prop.hlf-young.arm-prod} \textbf{(b)} yields
$\prod_{c\in Y\left(  \lambda\right)  }\left\vert \operatorname*{Arm}%
\nolimits_{\lambda}^{+}\left(  c\right)  \right\vert =\left\vert
\mathcal{R}\left(  T\right)  \right\vert $).

Now, I claim that each $m\in\left\{  0,1,\ldots,n\right\}  $ satisfies
\begin{equation}
PN_{m}PN=a_{m}P^{2}N. \label{pf.thm.tableau.hlf.main}%
\end{equation}

Before we come to the proof of this claim, let me explain why it yields
Theorem \ref{thm.tableau.hlf}. Indeed, assume (for a moment) that
(\ref{pf.thm.tableau.hlf.main}) is proved. Then, we can apply
(\ref{pf.thm.tableau.hlf.main}) to $m=n$. As a result, we obtain
$PN_{n}PN=a_{n}P^{2}N$. Comparing this with $P\underbrace{N_{n}}%
_{=N}PN=PNPN=\left(  PN\right)  ^{2}$, we obtain%
\begin{align*}
\left(  PN\right)  ^{2}  &  =\underbrace{a_{n}}_{=\dfrac{\prod_{c\in Y\left(
\lambda\right)  }\left\vert H_{\lambda}\left(  c\right)  \right\vert
}{\left\vert \mathcal{R}\left(  T\right)  \right\vert }}\ \ \underbrace{P^{2}%
}_{=\left\vert \mathcal{R}\left(  T\right)  \right\vert \cdot P}%
N=\underbrace{\dfrac{\prod_{c\in Y\left(  \lambda\right)  }\left\vert
H_{\lambda}\left(  c\right)  \right\vert }{\left\vert \mathcal{R}\left(
T\right)  \right\vert }\cdot\left\vert \mathcal{R}\left(  T\right)
\right\vert }_{=\prod_{c\in Y\left(  \lambda\right)  }\left\vert H_{\lambda
}\left(  c\right)  \right\vert }\cdot\,PN\\
&  =\left(  \prod_{c\in Y\left(  \lambda\right)  }\left\vert H_{\lambda
}\left(  c\right)  \right\vert \right)  \cdot PN.
\end{align*}
Hence, Lemma \ref{lem.hlf-young.last} (applied to $\kappa=\prod_{c\in Y\left(
\lambda\right)  }\left\vert H_{\lambda}\left(  c\right)  \right\vert $) yields
that the \# of standard tableaux of shape $Y\left(  \lambda\right)  $ is
$\dfrac{n!}{\prod_{c\in Y\left(  \lambda\right)  }\left\vert H_{\lambda
}\left(  c\right)  \right\vert }$. Thus, Theorem \ref{thm.tableau.hlf} is
proved, provided that we can show (\ref{pf.thm.tableau.hlf.main}).

It thus remains to prove (\ref{pf.thm.tableau.hlf.main}). We shall do so by
induction on $m$:

\textit{Base case:} Recall that $a_{0}=1$ and $N_{0}=1$. Thus,
$P\underbrace{N_{0}}_{=1}PN=PPN=P^{2}N$ and $\underbrace{a_{0}}_{=1}%
P^{2}N=P^{2}N$. Comparing these two equalities, we obtain $PN_{0}PN=a_{0}%
P^{2}N$. In other words, (\ref{pf.thm.tableau.hlf.main}) holds for $m=0$. This
completes the base case.

\textit{Induction step:} Let $m\in\left[  n\right]  $. Assume (as the
induction hypothesis) that (\ref{pf.thm.tableau.hlf.main}) holds for $m-1$
instead of $m$. We must prove that (\ref{pf.thm.tableau.hlf.main}) also holds
for $m$.

Let $c=\left(  i,j\right)  $ be the cell of $T$ that contains the entry $m$.
Thus, $T\left(  c\right)  =m$. In other words, $T\left(  i,j\right)  =m$
(since $c=\left(  i,j\right)  $).

\begin{example}
\label{exa.pf.thm.tableau.hlf.Tk2}For instance, let $n$, $\lambda$ and $T$ be
as in Example \ref{exa.pf.thm.tableau.hlf.T}, and let $m=7$. Then, we have%
\[
T=\ytableaushort{{10}{11}{12}{13}{14},{*(green)6}{*(red)7}89,{*(green)3}{*(green)4}{*(green)5},{*(green)1}{*(green)2}}\ \ ,
\]
where the cells with entries $1,2,\ldots,m-1$ are colored green, while the
cell with entry $m$ is colored red. The element $N_{m-1}$ is entirely
determined by the entries in green cells (see Remark
\ref{exa.pf.thm.tableau.hlf.Tk}, which writes $N_{m-1}$ as the column
antisymmetrizer $\nabla_{\operatorname*{Col}\left(  T_{m-1}\right)  }^{-}$ of
the green subtableau $T_{m-1}$), whereas the element $N_{m}$ is determined by
the entries in green and red cells.
\end{example}

Theorem \ref{thm.partitions.conj} \textbf{(e)} shows that the cells in the
$j$-th column of $Y\left(  \lambda\right)  $ are%
\[
\left(  1,j\right)  ,\ \left(  2,j\right)  ,\ \ldots,\ \left(  \lambda_{j}%
^{t},j\right)  .
\]
Hence, the entries in the $j$-th column of $T$ are $T\left(  1,j\right)
,\ T\left(  2,j\right)  ,\ \ldots,\ T\left(  \lambda_{j}^{t},j\right)  $. In
other words,%
\begin{equation}
\operatorname*{Col}\left(  j,T\right)  =\left\{  T\left(  1,j\right)
,\ T\left(  2,j\right)  ,\ \ldots,\ T\left(  \lambda_{j}^{t},j\right)
\right\}  \label{pf.thm.tableau.hlf.ColjT=}%
\end{equation}
(since $\operatorname*{Col}\left(  j,T\right)  $ is defined as the set of all
entries in the $j$-th column of $T$). Thus, $T\left(  \ell,j\right)  $ is
well-defined for each $\ell\in\left[  \lambda_{j}^{t}\right]  $.

We note that $j\in\left[  k\right]  $\ \ \ \ \footnote{\textit{Proof.} We know
that $\left(  i,j\right)  $ is a cell of $Y\left(  \lambda\right)  $ (by its
definition). Thus, $\left(  i,j\right)  $ lies in one of the columns
$1,2,\ldots,k$ (since all cells of $Y\left(  \lambda\right)  $ lie in columns
$1,2,\ldots,k$). In other words, column $j$ is one of the columns
$1,2,\ldots,k$ (since $\left(  i,j\right)  $ lies in column $j$). In other
words, $j\in\left\{  1,2,\ldots,k\right\}  =\left[  k\right]  $.}. Moreover,
$m$ is an entry in the $j$-th column of $T$ (since $m=T\left(  i,j\right)  $).
In other words, $m\in\operatorname*{Col}\left(  j,T\right)  $.

Now we shall show a sequence of claims. The first three claims are obvious
from a quick look at the construction of $T$ and the fact that $T\left(
i,j\right)  =m$, but we give detailed proofs just for the sake of completeness.

\begin{statement}
\textit{Claim 1:} Let $v$ be a positive integer such that $v\neq j$. Then,%
\[
\left[  m\right]  \cap\operatorname*{Col}\left(  v,T\right)  =\left[
m-1\right]  \cap\operatorname*{Col}\left(  v,T\right)  .
\]

\end{statement}

\begin{fineprint}
\begin{proof}
[Proof of Claim 1.]Recall that all entries of $T$ are distinct. Hence, an
entry in the $j$-th column of $T$ cannot also be found in the $v$-th column of
$T$ (since $v\neq j$). In particular, $m$ cannot also be found in the $v$-th
column of $T$ (since $m$ is an entry in the $j$-th column of $T$). In other
words, $m\notin\operatorname*{Col}\left(  v,T\right)  $. But $\left[
m\right]  =\left[  m-1\right]  \cup\left\{  m\right\}  $, and thus%
\begin{align*}
\underbrace{\left[  m\right]  }_{=\left[  m-1\right]  \cup\left\{  m\right\}
}\cap\operatorname*{Col}\left(  v,T\right)   &  =\left(  \left[  m-1\right]
\cup\left\{  m\right\}  \right)  \cap\operatorname*{Col}\left(  v,T\right) \\
&  =\left(  \left[  m-1\right]  \cap\operatorname*{Col}\left(  v,T\right)
\right)  \cup\underbrace{\left(  \left\{  m\right\}  \cap\operatorname*{Col}%
\left(  v,T\right)  \right)  }_{\substack{=\varnothing\\\text{(since }%
m\notin\operatorname*{Col}\left(  v,T\right)  \text{)}}}\\
&  =\left(  \left[  m-1\right]  \cap\operatorname*{Col}\left(  v,T\right)
\right)  \cup\varnothing=\left[  m-1\right]  \cap\operatorname*{Col}\left(
v,T\right)  .
\end{align*}
This proves Claim 1.
\end{proof}
\end{fineprint}

\begin{statement}
\textit{Claim 2:} We have $m\in\left[  m\right]  \cap\operatorname*{Col}%
\left(  j,T\right)  $ and%
\[
\left[  m-1\right]  \cap\operatorname*{Col}\left(  j,T\right)  =\left(
\left[  m\right]  \cap\operatorname*{Col}\left(  j,T\right)  \right)
\setminus\left\{  m\right\}  .
\]

\end{statement}

\begin{fineprint}
\begin{proof}
[Proof of Claim 2.]Combining $m\in\left[  m\right]  $ with $m\in
\operatorname*{Col}\left(  j,T\right)  $, we find $m\in\left[  m\right]
\cap\operatorname*{Col}\left(  j,T\right)  $. Furthermore, $\left[
m-1\right]  =\left[  m\right]  \setminus\left\{  m\right\}  $, so that
\[
\underbrace{\left[  m-1\right]  }_{=\left[  m\right]  \setminus\left\{
m\right\}  }\cap\operatorname*{Col}\left(  j,T\right)  =\left(  \left[
m\right]  \setminus\left\{  m\right\}  \right)  \cap\operatorname*{Col}\left(
j,T\right)  =\left(  \left[  m\right]  \cap\operatorname*{Col}\left(
j,T\right)  \right)  \setminus\left\{  m\right\}  .
\]
Thus, Claim 2 is proved.
\end{proof}
\end{fineprint}

\begin{statement}
\textit{Claim 3:} We have%
\[
\left[  m-1\right]  \cap\operatorname*{Col}\left(  j,T\right)  =\left\{
T\left(  i+1,j\right)  ,\ T\left(  i+2,j\right)  ,\ \ldots,\ T\left(
\lambda_{j}^{t},j\right)  \right\}  .
\]

\end{statement}

\begin{fineprint}
\begin{proof}
[Proof of Claim 3.]In our construction of the tableau $T$, we filled the cells
of $Y\left(  \lambda\right)  $ with the entries $1,2,\ldots,n$ in a specific
order -- namely in the order in which these cells appeared in the list
(\ref{pf.thm.tableau.hlf.list}). But the list (\ref{pf.thm.tableau.hlf.list})
was listing the cells of $Y\left(  \lambda\right)  $ row by row, from bottom
to top. Thus, if a cell $c\in Y\left(  \lambda\right)  $ lies further south
than another cell $d\in Y\left(  \lambda\right)  $, then $c$ appears before
$d$ in the list (\ref{pf.thm.tableau.hlf.list}), and therefore the cell $c$ is
filled with a smaller entry in the tableau $T$ than the cell $d$ (since the
cells are filled with $1,2,\ldots,n$ in the order in which they appear in the
list); in other words, we have $T\left(  c\right)  <T\left(  d\right)  $ in
this case. Hence, in particular, the entries of $T$ strictly decrease down
each column. In particular, we thus have%
\[
T\left(  1,j\right)  >T\left(  2,j\right)  >\cdots>T\left(  \lambda_{j}%
^{t},j\right)  .
\]
In other words, if $p$ and $q$ are two numbers in $\left[  \lambda_{j}%
^{t}\right]  $ satisfying $p<q$, then%
\begin{equation}
T\left(  p,j\right)  >T\left(  q,j\right)  .
\label{pf.thm.tableau.hlf.c3.pf.1}%
\end{equation}
Thus,%
\[
\left\{  T\left(  i+1,j\right)  ,\ T\left(  i+2,j\right)  ,\ \ldots,\ T\left(
\lambda_{j}^{t},j\right)  \right\}  \subseteq\left[  m-1\right]
\]
\footnote{\textit{Proof.} Let $g\in\left\{  T\left(  i+1,j\right)  ,\ T\left(
i+2,j\right)  ,\ \ldots,\ T\left(  \lambda_{j}^{t},j\right)  \right\}  $.
Thus, $g=T\left(  s,j\right)  $ for some $s\in\left\{  i+1,i+2,\ldots
,\lambda_{j}^{t}\right\}  $. Consider this $s$.
\par
From $s\in\left\{  i+1,i+2,\ldots,\lambda_{j}^{t}\right\}  $, we obtain $s\geq
i+1>i$, so that $i<s$. Therefore, $T\left(  i,j\right)  >T\left(  s,j\right)
$ (by (\ref{pf.thm.tableau.hlf.c3.pf.1}), applied to $p=i$ and $q=s$). Hence,
$T\left(  s,j\right)  <T\left(  i,j\right)  =m$. Hence, $g=T\left(
s,j\right)  <m$, so that $g\leq m-1$ (since $g$ and $m$ are integers).
Therefore, $g\in\left[  m-1\right]  $.
\par
Forget that we fixed $g$. We thus have shown that $g\in\left[  m-1\right]  $
for each $g\in\left\{  T\left(  i+1,j\right)  ,\ T\left(  i+2,j\right)
,\ \ldots,\ T\left(  \lambda_{j}^{t},j\right)  \right\}  $. In other words,
$\left\{  T\left(  i+1,j\right)  ,\ T\left(  i+2,j\right)  ,\ \ldots
,\ T\left(  \lambda_{j}^{t},j\right)  \right\}  \subseteq\left[  m-1\right]
$.} and%
\[
\left\{  T\left(  1,j\right)  ,\ T\left(  2,j\right)  ,\ \ldots,\ T\left(
i,j\right)  \right\}  \cap\left[  m-1\right]  =\varnothing
\]
\footnote{\textit{Proof.} Let $g\in\left\{  T\left(  1,j\right)  ,\ T\left(
2,j\right)  ,\ \ldots,\ T\left(  i,j\right)  \right\}  \cap\left[  m-1\right]
$. We shall derive a contradiction.
\par
We have $g\in\left\{  T\left(  1,j\right)  ,\ T\left(  2,j\right)
,\ \ldots,\ T\left(  i,j\right)  \right\}  \cap\left[  m-1\right]
\subseteq\left\{  T\left(  1,j\right)  ,\ T\left(  2,j\right)  ,\ \ldots
,\ T\left(  i,j\right)  \right\}  $. In other words, $g=T\left(  s,j\right)  $
for some $s\in\left\{  1,2,\ldots,i\right\}  $. Consider this $s$. Thus,
$s\leq i$ (since $s\in\left\{  1,2,\ldots,i\right\}  $). On the other hand,
$g\in\left\{  T\left(  1,j\right)  ,\ T\left(  2,j\right)  ,\ \ldots
,\ T\left(  i,j\right)  \right\}  \cap\left[  m-1\right]  \subseteq\left[
m-1\right]  $ and therefore $g\leq m-1<m$. Hence, $T\left(  s,j\right)
=g<m=T\left(  i,j\right)  $. Therefore, $T\left(  s,j\right)  \neq T\left(
i,j\right)  $, so that $s\neq i$. Combining this with $s\leq i$, we obtain
$s<i$. Therefore, $T\left(  s,j\right)  >T\left(  i,j\right)  $ (by
(\ref{pf.thm.tableau.hlf.c3.pf.1}), applied to $p=s$ and $q=i$). But this
contradicts $T\left(  s,j\right)  <T\left(  i,j\right)  $.
\par
Forget that we fixed $g$. We thus have found a contradiction for each
$g\in\left\{  T\left(  1,j\right)  ,\ T\left(  2,j\right)  ,\ \ldots
,\ T\left(  i,j\right)  \right\}  \cap\left[  m-1\right]  $. Hence, there
exists no such $g$. In other words, $\left\{  T\left(  1,j\right)  ,\ T\left(
2,j\right)  ,\ \ldots,\ T\left(  i,j\right)  \right\}  \cap\left[  m-1\right]
=\varnothing$.}. However, (\ref{pf.thm.tableau.hlf.ColjT=}) shows that%
\begin{align*}
\operatorname*{Col}\left(  j,T\right)   &  =\left\{  T\left(  1,j\right)
,\ T\left(  2,j\right)  ,\ \ldots,\ T\left(  \lambda_{j}^{t},j\right)
\right\} \\
&  =\left\{  T\left(  1,j\right)  ,\ T\left(  2,j\right)  ,\ \ldots,\ T\left(
i,j\right)  ,\ T\left(  i+1,j\right)  ,\ T\left(  i+2,j\right)  ,\ \ldots
,\ T\left(  \lambda_{j}^{t},j\right)  \right\} \\
&  =\left\{  T\left(  1,j\right)  ,\ T\left(  2,j\right)  ,\ \ldots,\ T\left(
i,j\right)  \right\}  \cup\left\{  T\left(  i+1,j\right)  ,\ T\left(
i+2,j\right)  ,\ \ldots,\ T\left(  \lambda_{j}^{t},j\right)  \right\}  .
\end{align*}
Therefore,%
\begin{align*}
&  \left[  m-1\right]  \cap\underbrace{\operatorname*{Col}\left(  j,T\right)
}_{=\left\{  T\left(  1,j\right)  ,\ T\left(  2,j\right)  ,\ \ldots,\ T\left(
i,j\right)  \right\}  \cup\left\{  T\left(  i+1,j\right)  ,\ T\left(
i+2,j\right)  ,\ \ldots,\ T\left(  \lambda_{j}^{t},j\right)  \right\}  }\\
&  =\left[  m-1\right]  \cap\left(  \left\{  T\left(  1,j\right)  ,\ T\left(
2,j\right)  ,\ \ldots,\ T\left(  i,j\right)  \right\}  \cup\left\{  T\left(
i+1,j\right)  ,\ T\left(  i+2,j\right)  ,\ \ldots,\ T\left(  \lambda_{j}%
^{t},j\right)  \right\}  \right) \\
&  =\underbrace{\left(  \left[  m-1\right]  \cap\left\{  T\left(  1,j\right)
,\ T\left(  2,j\right)  ,\ \ldots,\ T\left(  i,j\right)  \right\}  \right)
}_{\substack{=\left\{  T\left(  1,j\right)  ,\ T\left(  2,j\right)
,\ \ldots,\ T\left(  i,j\right)  \right\}  \cap\left[  m-1\right]
\\=\varnothing}}\\
&  \ \ \ \ \ \ \ \ \ \ \cup\underbrace{\left(  \left[  m-1\right]
\cap\left\{  T\left(  i+1,j\right)  ,\ T\left(  i+2,j\right)  ,\ \ldots
,\ T\left(  \lambda_{j}^{t},j\right)  \right\}  \right)  }%
_{\substack{=\left\{  T\left(  i+1,j\right)  ,\ T\left(  i+2,j\right)
,\ \ldots,\ T\left(  \lambda_{j}^{t},j\right)  \right\}  \\\text{(since
}\left\{  T\left(  i+1,j\right)  ,\ T\left(  i+2,j\right)  ,\ \ldots
,\ T\left(  \lambda_{j}^{t},j\right)  \right\}  \subseteq\left[  m-1\right]
\text{)}}}\\
&  =\varnothing\cup\left\{  T\left(  i+1,j\right)  ,\ T\left(  i+2,j\right)
,\ \ldots,\ T\left(  \lambda_{j}^{t},j\right)  \right\} \\
&  =\left\{  T\left(  i+1,j\right)  ,\ T\left(  i+2,j\right)  ,\ \ldots
,\ T\left(  \lambda_{j}^{t},j\right)  \right\}  .
\end{align*}
This proves Claim 3.
\end{proof}
\end{fineprint}

\begin{statement}
\textit{Claim 4:} We have%
\begin{equation}
N_{m-1}\sum_{\ell=i+1}^{\lambda_{j}^{t}}t_{m,\ T\left(  \ell,j\right)
}=N_{m-1}-N_{m}. \label{pf.thm.tableau.hlf.c4}%
\end{equation}

\end{statement}

\begin{proof}
[Proof of Claim 4.]The equality (\ref{pf.thm.tableau.hlf.Nm=}) yields%
\[
N_{m}=\prod_{v=1}^{k}\nabla_{\left[  m\right]  \cap\operatorname*{Col}\left(
v,T\right)  }^{-}%
\]
(here, we have renamed the index $j$ in (\ref{pf.thm.tableau.hlf.Nm=}) as $v$,
since the letter $j$ already has another meaning in our context). As we know,
the order of the factors in this product is immaterial, since they all
commute. Thus,%
\begin{align}
N_{m}  &  =\prod_{v=1}^{k}\nabla_{\left[  m\right]  \cap\operatorname*{Col}%
\left(  v,T\right)  }^{-}=\prod_{v\in\left[  k\right]  }\nabla_{\left[
m\right]  \cap\operatorname*{Col}\left(  v,T\right)  }^{-}%
\label{pf.thm.tableau.hlf.c4.pf.0}\\
&  =\left(  \prod_{\substack{v\in\left[  k\right]  ;\\v\neq j}}\nabla_{\left[
m\right]  \cap\operatorname*{Col}\left(  v,T\right)  }^{-}\right)
\nabla_{\left[  m\right]  \cap\operatorname*{Col}\left(  j,T\right)  }^{-}
\label{pf.thm.tableau.hlf.c4.pf.1}%
\end{align}
(here, we have split off the factor for $v=j$ from our product, since we know
that $j\in\left[  k\right]  $). The same argument can be applied to $m-1$
instead of $m$, and thus we find%
\begin{equation}
N_{m-1}=\left(  \prod_{\substack{v\in\left[  k\right]  ;\\v\neq j}%
}\nabla_{\left[  m-1\right]  \cap\operatorname*{Col}\left(  v,T\right)  }%
^{-}\right)  \nabla_{\left[  m-1\right]  \cap\operatorname*{Col}\left(
j,T\right)  }^{-}. \label{pf.thm.tableau.hlf.c4.pf.2}%
\end{equation}
However, for each $v\in\left[  k\right]  $ satisfying $v\neq j$, we have
$\left[  m\right]  \cap\operatorname*{Col}\left(  v,T\right)  =\left[
m-1\right]  \cap\operatorname*{Col}\left(  v,T\right)  $ (by Claim 1) and
therefore%
\begin{equation}
\nabla_{\left[  m\right]  \cap\operatorname*{Col}\left(  v,T\right)  }%
^{-}=\nabla_{\left[  m-1\right]  \cap\operatorname*{Col}\left(  v,T\right)
}^{-}. \label{pf.thm.tableau.hlf.c4.pf.3}%
\end{equation}
Moreover, Claim 2 shows that
\[
\left[  m-1\right]  \cap\operatorname*{Col}\left(  j,T\right)  =\left(
\left[  m\right]  \cap\operatorname*{Col}\left(  j,T\right)  \right)
\setminus\left\{  m\right\}  .
\]
But Claim 2 also shows that $m\in\left[  m\right]  \cap\operatorname*{Col}%
\left(  j,T\right)  $. Thus, Lemma \ref{lem.intX.rec2sign} (applied to
$X=\left[  m\right]  \cap\operatorname*{Col}\left(  j,T\right)  $ and $x=m$)
yields%
\[
\nabla_{\left[  m\right]  \cap\operatorname*{Col}\left(  j,T\right)  }%
^{-}=\nabla_{\left(  \left[  m\right]  \cap\operatorname*{Col}\left(
j,T\right)  \right)  \setminus\left\{  m\right\}  }^{-}\left(  1-\sum
_{y\in\left(  \left[  m\right]  \cap\operatorname*{Col}\left(  j,T\right)
\right)  \setminus\left\{  m\right\}  }t_{y,m}\right)  .
\]
In view of $\left[  m-1\right]  \cap\operatorname*{Col}\left(  j,T\right)
=\left(  \left[  m\right]  \cap\operatorname*{Col}\left(  j,T\right)  \right)
\setminus\left\{  m\right\}  $, we can rewrite this as
\begin{equation}
\nabla_{\left[  m\right]  \cap\operatorname*{Col}\left(  j,T\right)  }%
^{-}=\nabla_{\left[  m-1\right]  \cap\operatorname*{Col}\left(  j,T\right)
}^{-}\left(  1-\sum_{y\in\left[  m-1\right]  \cap\operatorname*{Col}\left(
j,T\right)  }t_{y,m}\right)  . \label{pf.thm.tableau.hlf.c4.pf.5}%
\end{equation}

Claim 3 shows that $\left[  m-1\right]  \cap\operatorname*{Col}\left(
j,T\right)  =\left\{  T\left(  i+1,j\right)  ,\ T\left(  i+2,j\right)
,\ \ldots,\ T\left(  \lambda_{j}^{t},j\right)  \right\}  $. Hence,
\begin{align*}
&  \sum_{y\in\left[  m-1\right]  \cap\operatorname*{Col}\left(  j,T\right)
}t_{y,m}\\
&  =\sum_{y\in\left\{  T\left(  i+1,j\right)  ,\ T\left(  i+2,j\right)
,\ \ldots,\ T\left(  \lambda_{j}^{t},j\right)  \right\}  }t_{y,m}\\
&  =t_{T\left(  i+1,j\right)  ,\ m}+t_{T\left(  i+2,j\right)  ,\ m}%
+\cdots+t_{T\left(  \lambda_{j}^{t},j\right)  ,\ m}\\
&  \ \ \ \ \ \ \ \ \ \ \ \ \ \ \ \ \ \ \ \ \left(
\begin{array}
[c]{c}%
\text{since the elements }T\left(  i+1,j\right)  ,\ T\left(  i+2,j\right)
,\ \ldots,\ T\left(  \lambda_{j}^{t},j\right) \\
\text{are distinct (because all entries of }T\text{ are distinct)}%
\end{array}
\right) \\
&  =\sum_{\ell=i+1}^{\lambda_{j}^{t}}\underbrace{t_{T\left(  \ell,j\right)
,\ m}}_{\substack{=t_{m,\ T\left(  \ell,j\right)  }\\\text{(since }%
t_{p,q}=t_{q,p}\text{ for any }p\neq q\text{)}}}=\sum_{\ell=i+1}^{\lambda
_{j}^{t}}t_{m,\ T\left(  \ell,j\right)  }.
\end{align*}
Thus, we can rewrite (\ref{pf.thm.tableau.hlf.c4.pf.5}) as
\begin{equation}
\nabla_{\left[  m\right]  \cap\operatorname*{Col}\left(  j,T\right)  }%
^{-}=\nabla_{\left[  m-1\right]  \cap\operatorname*{Col}\left(  j,T\right)
}^{-}\left(  1-\sum_{\ell=i+1}^{\lambda_{j}^{t}}t_{m,\ T\left(  \ell,j\right)
}\right)  . \label{pf.thm.tableau.hlf.c4.pf.6}%
\end{equation}

Now, (\ref{pf.thm.tableau.hlf.c4.pf.1}) becomes%
\begin{align*}
N_{m}  &  =\left(  \prod_{\substack{v\in\left[  k\right]  ;\\v\neq
j}}\underbrace{\nabla_{\left[  m\right]  \cap\operatorname*{Col}\left(
v,T\right)  }^{-}}_{\substack{=\nabla_{\left[  m-1\right]  \cap
\operatorname*{Col}\left(  v,T\right)  }^{-}\\\text{(by
(\ref{pf.thm.tableau.hlf.c4.pf.3}))}}}\right)  \underbrace{\nabla_{\left[
m\right]  \cap\operatorname*{Col}\left(  j,T\right)  }^{-}}_{\substack{=\nabla
_{\left[  m-1\right]  \cap\operatorname*{Col}\left(  j,T\right)  }^{-}\left(
1-\sum_{\ell=i+1}^{\lambda_{j}^{t}}t_{m,\ T\left(  \ell,j\right)  }\right)
\\\text{(by (\ref{pf.thm.tableau.hlf.c4.pf.6}))}}}\\
&  =\underbrace{\left(  \prod_{\substack{v\in\left[  k\right]  ;\\v\neq
j}}\nabla_{\left[  m-1\right]  \cap\operatorname*{Col}\left(  v,T\right)
}^{-}\right)  \nabla_{\left[  m-1\right]  \cap\operatorname*{Col}\left(
j,T\right)  }^{-}}_{\substack{=N_{m-1}\\\text{(by
(\ref{pf.thm.tableau.hlf.c4.pf.2}))}}}\left(  1-\sum_{\ell=i+1}^{\lambda
_{j}^{t}}t_{m,\ T\left(  \ell,j\right)  }\right) \\
&  =N_{m-1}\left(  1-\sum_{\ell=i+1}^{\lambda_{j}^{t}}t_{m,\ T\left(
\ell,j\right)  }\right)  =N_{m-1}-N_{m-1}\sum_{\ell=i+1}^{\lambda_{j}^{t}%
}t_{m,\ T\left(  \ell,j\right)  }.
\end{align*}
Subtracting $N_{m-1}$ from both sides of this equality, we find%
\[
N_{m}-N_{m-1}=-N_{m-1}\sum_{\ell=i+1}^{\lambda_{j}^{t}}t_{m,\ T\left(
\ell,j\right)  }.
\]
Thus,%
\[
N_{m-1}\sum_{\ell=i+1}^{\lambda_{j}^{t}}t_{m,\ T\left(  \ell,j\right)
}=-\left(  N_{m}-N_{m-1}\right)  =N_{m-1}-N_{m}.
\]
This proves Claim 4.
\end{proof}

\begin{statement}
\textit{Claim 5:} Let $\ell\in\left\{  i+1,i+2,\ldots,\lambda_{j}^{t}\right\}
$. Then,%
\[
\left(  1+\sum_{s=1}^{\lambda_{i}}t_{T\left(  i,s\right)  ,\ T\left(
\ell,j\right)  }\right)  PN=0.
\]

\end{statement}

\begin{proof}
[Proof of Claim 5.]We have $\ell\in\left\{  i+1,i+2,\ldots,\lambda_{j}%
^{t}\right\}  $, so that $\ell\geq i+1>i$ and $\ell\leq\lambda_{j}^{t}$. The
cell $\left(  \ell,j\right)  $ is one of the cells $\left(  1,j\right)
,\ \left(  2,j\right)  ,\ \ldots,\ \left(  \lambda_{j}^{t},j\right)  $ (since
$\ell\geq i+1\geq1$ and $\ell\leq\lambda_{j}^{t}$), and therefore is a cell in
the $j$-th column of $Y\left(  \lambda\right)  $ (since the cells in the
$j$-th column of $Y\left(  \lambda\right)  $ are $\left(  1,j\right)
,\ \left(  2,j\right)  ,\ \ldots,\ \left(  \lambda_{j}^{t},j\right)  $). Thus,
in particular, $\left(  \ell,j\right)  \in Y\left(  \lambda\right)  $. In
other words, $j\in\left[  \lambda_{\ell}\right]  $ (by the definition of
$Y\left(  \lambda\right)  $). Thus, Proposition \ref{prop.hlf-young.grel5}
(applied to $\ell$, $i$ and $j$ instead of $j$, $k$ and $q$) yields
\[
\left(  1+\sum_{s=1}^{\lambda_{i}}t_{T\left(  i,s\right)  ,\ T\left(
\ell,j\right)  }\right)  \nabla_{\operatorname*{Row}T}\nabla
_{\operatorname*{Col}T}^{-}=0.
\]
This rewrites as%
\[
\left(  1+\sum_{s=1}^{\lambda_{i}}t_{T\left(  i,s\right)  ,\ T\left(
\ell,j\right)  }\right)  PN=0
\]
(since $P=\nabla_{\operatorname*{Row}T}$ and $N=\nabla_{\operatorname*{Col}%
T}^{-}$). Thus, Claim 5 is proved.
\end{proof}

\begin{statement}
\textit{Claim 6:} Let $p,q\in\left[  n\right]  $ be two distinct numbers that
lie in the same column of $T$ and satisfy $p<m$ and $q<m$. Then,
$N_{m-1}=N^{\prime}\cdot\left(  1-t_{p,q}\right)  $ for some $N^{\prime}%
\in\mathcal{A}$.
\end{statement}

\begin{proof}
[Proof of Claim 6.]We have assumed that $p$ and $q$ lie in the same column of
$T$. In other words, $p$ and $q$ lie in the $s$-th column of $T$ for some
$s\in\left[  k\right]  $ (since all cells of $Y\left(  \lambda\right)  $ lie
in columns $1,2,\ldots,k$). Consider this $s$. Set $X:=\left[  m-1\right]
\cap\operatorname*{Col}\left(  s,T\right)  $ and%
\[
\nabla_{X}^{\operatorname*{even}}:=\sum_{\substack{w\in S_{n,X};\\\left(
-1\right)  ^{w}=1}}w
\]
(where $S_{n,X}$ is defined according to Proposition \ref{prop.intX.basics}).

We have $p\in\operatorname*{Col}\left(  s,T\right)  $ (since $p$ lies in the
$s$-th column of $T$). Moreover, from $p<m$, we obtain $p\leq m-1$ (since $p$
and $m$ are integers), so that $p\in\left[  m-1\right]  $ (since $p\in\left[
n\right]  $ entails $p\geq1$). Combining this with $p\in\operatorname*{Col}%
\left(  s,T\right)  $, we obtain $p\in\left[  m-1\right]  \cap
\operatorname*{Col}\left(  s,T\right)  =X$. Similarly, $q\in X$. Hence, $p$
and $q$ are two distinct elements of $X$. Therefore, Proposition
\ref{prop.intX.basics} \textbf{(d)} (applied to $p$ and $q$ instead of $i$ and
$j$) yields%
\[
\nabla_{X}=\nabla_{X}^{\operatorname*{even}}\cdot\left(  1+t_{p,q}\right)
=\left(  1+t_{p,q}\right)  \cdot\nabla_{X}^{\operatorname*{even}}%
\]
and%
\[
\nabla_{X}^{-}=\nabla_{X}^{\operatorname*{even}}\cdot\left(  1-t_{p,q}\right)
=\left(  1-t_{p,q}\right)  \cdot\nabla_{X}^{\operatorname*{even}}.
\]
However, the equality (\ref{pf.thm.tableau.hlf.c4.pf.0}) (which we showed
during our proof of Claim 4) shows that%
\[
N_{m}=\prod_{v\in\left[  k\right]  }\nabla_{\left[  m\right]  \cap
\operatorname*{Col}\left(  v,T\right)  }^{-}.
\]
The same argument (applied to $m-1$ instead of $m$) yields%
\begin{align*}
N_{m-1}  &  =\prod_{v\in\left[  k\right]  }\nabla_{\left[  m-1\right]
\cap\operatorname*{Col}\left(  v,T\right)  }^{-}=\left(  \prod_{\substack{v\in
\left[  k\right]  ;\\v\neq s}}\nabla_{\left[  m-1\right]  \cap
\operatorname*{Col}\left(  v,T\right)  }^{-}\right)  \underbrace{\nabla
_{\left[  m-1\right]  \cap\operatorname*{Col}\left(  s,T\right)  }^{-}%
}_{\substack{=\nabla_{X}^{-}\\\text{(since }\left[  m-1\right]  \cap
\operatorname*{Col}\left(  s,T\right)  =X\text{)}}}\\
&  \ \ \ \ \ \ \ \ \ \ \ \ \ \ \ \ \ \ \ \ \left(
\begin{array}
[c]{c}%
\text{here, we split off the factor for }v=s\\
\text{from the product, since }s\in\left[  k\right]
\end{array}
\right) \\
&  =\left(  \prod_{\substack{v\in\left[  k\right]  ;\\v\neq s}}\nabla_{\left[
m-1\right]  \cap\operatorname*{Col}\left(  v,T\right)  }^{-}\right)
\underbrace{\nabla_{X}^{-}}_{=\nabla_{X}^{\operatorname*{even}}\cdot\left(
1-t_{p,q}\right)  }\\
&  =\left(  \prod_{\substack{v\in\left[  k\right]  ;\\v\neq s}}\nabla_{\left[
m-1\right]  \cap\operatorname*{Col}\left(  v,T\right)  }^{-}\right)
\nabla_{X}^{\operatorname*{even}}\cdot\left(  1-t_{p,q}\right)  .
\end{align*}
Therefore, $N_{m-1}=N^{\prime}\cdot\left(  1-t_{p,q}\right)  $ for some
$N^{\prime}\in\mathcal{A}$ (namely, for \newline$N^{\prime}=\left(
\prod_{\substack{v\in\left[  k\right]  ;\\v\neq s}}\nabla_{\left[  m-1\right]
\cap\operatorname*{Col}\left(  v,T\right)  }^{-}\right)  \nabla_{X}%
^{\operatorname*{even}}$). This proves Claim 6.
\end{proof}

\begin{statement}
\textit{Claim 7:} Let $p,q\in\left[  n\right]  $ be two distinct numbers that
lie in the same row of $T$. Then, $t_{p,q}P=Pt_{p,q}=P$.
\end{statement}

\begin{proof}
[Proof of Claim 7.]Proposition \ref{prop.symmetrizers.factor-out-row}
\textbf{(c)} (applied to $Y\left(  \lambda\right)  $, $p$ and $q$ instead of
$D$, $i$ and $j$) shows that $t_{p,q}\nabla_{\operatorname*{Row}T}%
=\nabla_{\operatorname*{Row}T}t_{p,q}=\nabla_{\operatorname*{Row}T}$. In view
of $P=\nabla_{\operatorname*{Row}T}$, we can rewrite this as $t_{p,q}%
P=Pt_{p,q}=P$. Thus, Claim 7 is proved.
\end{proof}

\begin{statement}
\textit{Claim 8:} Let $\ell\in\left\{  i+1,i+2,\ldots,\lambda_{j}^{t}\right\}
$. Let $s\in\left[  j-1\right]  $. Then,%
\[
N_{m-1}t_{T\left(  i,s\right)  ,\ T\left(  \ell,j\right)  }P=0.
\]

\end{statement}

\begin{proof}
[Proof of Claim 8.]We have $s\in\left[  j-1\right]  $, so that $s\leq
j-1<j\leq\lambda_{i}$ (since $\left(  i,j\right)  \in Y\left(  \lambda\right)
$). Thus, $\left(  i,s\right)  \in Y\left(  \lambda\right)  $. Therefore, the
entry $T\left(  i,s\right)  $ of $T$ is well-defined. Furthermore, $\left(
\ell,j\right)  \in Y\left(  \lambda\right)  $ (as we saw in the proof of Claim
5). Hence, the entry $T\left(  \ell,j\right)  $ of $T$ is well-defined.
Moreover, from $\left(  \ell,j\right)  \in Y\left(  \lambda\right)  $, we
obtain $j\leq\lambda_{\ell}$. Altogether, $s<j\leq\lambda_{\ell}$, so that
$\left(  \ell,s\right)  \in Y\left(  \lambda\right)  $. Thus, the entry
$T\left(  \ell,s\right)  $ of $T$ is well-defined.

We have $\ell\in\left\{  i+1,i+2,\ldots,\lambda_{j}^{t}\right\}  $, thus
$\ell\geq i+1>i$ and therefore $\ell\neq i$. Also, $s\neq j$ (since $s<j$).

Set $p:=T\left(  \ell,s\right)  $ and $q:=T\left(  i,s\right)  $ and
$r:=T\left(  \ell,j\right)  $. The two cells $\left(  \ell,s\right)  $ and
$\left(  \ell,j\right)  $ are distinct (since $s\neq j$), and are both
distinct from $\left(  i,s\right)  $ (since $\ell\neq i$). In other words, the
three cells $\left(  \ell,s\right)  $ and $\left(  i,s\right)  $ and $\left(
\ell,j\right)  $ are distinct. Hence, the corresponding entries $T\left(
\ell,s\right)  $ and $T\left(  i,s\right)  $ and $T\left(  \ell,j\right)  $
are also distinct (since all entries of $T$ are distinct). In other words, $p$
and $q$ and $r$ are distinct (since $p=T\left(  \ell,s\right)  $ and
$q=T\left(  i,s\right)  $ and $r=T\left(  \ell,j\right)  $).

The numbers $p$ and $q$ both lie in the $s$-th column of $T$ (since
$p=T\left(  \ell,s\right)  $ and $q=T\left(  i,s\right)  $). Thus, they lie in
the same column of $T$.

The numbers $p$ and $r$ both lie in the $\ell$-th row of $T$ (since
$p=T\left(  \ell,s\right)  $ and $r=T\left(  \ell,j\right)  $). Thus, they lie
in the same row of $T$. Hence, Lemma \ref{lem.hlf-young.row1} (applied to
$D=Y\left(  \lambda\right)  $) yields%
\[
\left(  1-t_{p,q}\right)  t_{q,r}\nabla_{\operatorname*{Row}T}=0.
\]

Moreover, the construction of $T$ easily shows that $p<m$%
\ \ \ \ \footnote{\textit{Proof.} We have $\ell>i$. Thus, the cell $\left(
\ell,s\right)  $ lies further south than the cell $\left(  i,j\right)  $.
\par
However, during our proof of Claim 3, we saw that if a cell $c\in Y\left(
\lambda\right)  $ lies further south than another cell $d\in Y\left(
\lambda\right)  $, then $T\left(  c\right)  <T\left(  d\right)  $. We can
apply this to $c=\left(  \ell,s\right)  $ and $d=\left(  i,j\right)  $, and
thus conclude that $T\left(  \ell,s\right)  <T\left(  i,j\right)  $ (since the
cell $\left(  \ell,s\right)  $ lies further south than the cell $\left(
i,j\right)  $). In other words, $p<m$ (since $p=T\left(  \ell,s\right)  $ and
$m=T\left(  i,j\right)  $).} and $q<m$\ \ \ \ \footnote{\textit{Proof.} The
cell $\left(  i,s\right)  $ lies in the same row as the cell $\left(
i,j\right)  $ but further west (since $s<j$).
\par
In our construction of the tableau $T$, we filled the cells of $Y\left(
\lambda\right)  $ with the entries $1,2,\ldots,n$ in a specific order --
namely in the order in which these cells appeared in the list
(\ref{pf.thm.tableau.hlf.list}). But the list (\ref{pf.thm.tableau.hlf.list})
was listing the cells of $Y\left(  \lambda\right)  $ row by row, where the
cells of each row were listed from left to right. Thus, if a cell $c\in
Y\left(  \lambda\right)  $ lies in the same row as another cell $d\in Y\left(
\lambda\right)  $ but further west, then $c$ appears before $d$ in the list
(\ref{pf.thm.tableau.hlf.list}), and therefore the cell $c$ is filled with a
smaller entry in the tableau $T$ than the cell $d$ (since the cells are filled
with $1,2,\ldots,n$ in the order in which they appear in the list); in other
words, we have $T\left(  c\right)  <T\left(  d\right)  $ in this case.
Applying this to $c=\left(  i,s\right)  $ and $d=\left(  i,j\right)  $, we
conclude that $T\left(  i,s\right)  <T\left(  i,j\right)  $ (since the cell
$\left(  i,s\right)  $ lies in the same row as the cell $\left(  i,j\right)  $
but further west). In other words, $q<m$ (since $q=T\left(  i,s\right)  $ and
$m=T\left(  i,j\right)  $).}. Moreover, the two numbers $p$ and $q$ are
distinct (as we saw above) and lie in the same column of $T$ (as we also saw
above). Hence, Claim 6 yields that $N_{m-1}=N^{\prime}\cdot\left(
1-t_{p,q}\right)  $ for some $N^{\prime}\in\mathcal{A}$. Consider this
$N^{\prime}$.

Now,%
\[
\underbrace{N_{m-1}}_{=N^{\prime}\cdot\left(  1-t_{p,q}\right)  }%
\underbrace{t_{T\left(  i,s\right)  ,\ T\left(  \ell,j\right)  }%
}_{\substack{=t_{q,r}\\\text{(since }T\left(  i,s\right)  =q\\\text{and
}T\left(  \ell,j\right)  =r\text{)}}}\underbrace{P}_{=\nabla
_{\operatorname*{Row}T}}=N^{\prime}\cdot\underbrace{\left(  1-t_{p,q}\right)
t_{q,r}\nabla_{\operatorname*{Row}T}}_{=0}=0.
\]
This proves Claim 8.
\end{proof}

\begin{statement}
\textit{Claim 9:} Let $p,q\in\left[  n\right]  $ be two distinct numbers that
satisfy $p\geq m$ and $q\geq m$. Then, $t_{p,q}N_{m-1}=N_{m-1}t_{p,q}$.
\end{statement}

\begin{proof}
[Proof of Claim 9.]The equality (\ref{pf.thm.tableau.hlf.Nm=}) yields%
\[
N_{m}=\prod_{v=1}^{k}\nabla_{\left[  m\right]  \cap\operatorname*{Col}\left(
v,T\right)  }^{-}%
\]
(here, we have renamed the index $j$ in (\ref{pf.thm.tableau.hlf.Nm=}) as $v$,
since the letter $j$ already has another meaning in our context). The same
argument can be applied to $m-1$ instead of $m$, and thus we find%
\begin{align*}
N_{m-1}  &  =\prod_{v=1}^{k}\nabla_{\left[  m-1\right]  \cap
\operatorname*{Col}\left(  v,T\right)  }^{-}\\
&  =\nabla_{\left[  m-1\right]  \cap\operatorname*{Col}\left(  1,T\right)
}^{-}\nabla_{\left[  m-1\right]  \cap\operatorname*{Col}\left(  2,T\right)
}^{-}\cdots\nabla_{\left[  m-1\right]  \cap\operatorname*{Col}\left(
k,T\right)  }^{-}.
\end{align*}

However, it is easy to see (using Proposition \ref{prop.int.commute}
\textbf{(d)}) that the transposition $t_{p,q}$ commutes with the element
$\nabla_{\left[  m-1\right]  \cap\operatorname*{Col}\left(  v,T\right)  }^{-}$
for each $v\in\left[  k\right]  $\ \ \ \ \footnote{\textit{Proof.} Let
$v\in\left[  k\right]  $. We must prove that $t_{p,q}$ commutes with
$\nabla_{\left[  m-1\right]  \cap\operatorname*{Col}\left(  v,T\right)  }^{-}%
$.
\par
Let $X:=\left[  m-1\right]  \cap\operatorname*{Col}\left(  v,T\right)  $ and
$Y:=\left\{  m,m+1,\ldots,n\right\}  $. Then, both $X$ and $Y$ are subsets of
$\left[  n\right]  $. Moreover,%
\[
\underbrace{X}_{\substack{=\left[  m-1\right]  \cap\operatorname*{Col}\left(
v,T\right)  \\\subseteq\left[  m-1\right]  =\left\{  1,2,\ldots,m-1\right\}
}}\cap\underbrace{Y}_{=\left\{  m,m+1,\ldots,n\right\}  }\subseteq\left\{
1,2,\ldots,m-1\right\}  \cap\left\{  m,m+1,\ldots,n\right\}  =\varnothing
\]
(since the sets $\left\{  1,2,\ldots,m-1\right\}  $ and $\left\{
m,m+1,\ldots,n\right\}  $ are obviously disjoint). Thus, $X\cap Y=\varnothing
$. In other words, the sets $X$ and $Y$ are disjoint. Furthermore, from $p\geq
m$ and $p\leq n$ (since $p\in\left[  n\right]  $), we obtain $p\in\left\{
m,m+1,\ldots,n\right\}  =Y$. Similarly, $q\in Y$. Thus, $p$ and $q$ are two
distinct elements of $Y$. Hence, Proposition \ref{prop.int.commute}
\textbf{(d)} shows that $\nabla_{X}^{-}t_{p,q}=t_{p,q}\nabla_{X}^{-}$. In
other words, $t_{p,q}$ commutes with $\nabla_{X}^{-}$. In other words,
$t_{p,q}$ commutes with $\nabla_{\left[  m-1\right]  \cap\operatorname*{Col}%
\left(  v,T\right)  }^{-}$ (since $X=\left[  m-1\right]  \cap
\operatorname*{Col}\left(  v,T\right)  $). This completes our proof.}. In
other words, $t_{p,q}$ commutes with each of the elements \newline%
$\nabla_{\left[  m-1\right]  \cap\operatorname*{Col}\left(  1,T\right)  }%
^{-},\nabla_{\left[  m-1\right]  \cap\operatorname*{Col}\left(  2,T\right)
}^{-},\ldots,\nabla_{\left[  m-1\right]  \cap\operatorname*{Col}\left(
k,T\right)  }^{-}$.

However, recall a well-known fact from algebra:

\begin{statement}
\textit{Fact:} Let $A$ be a ring. Let $a\in A$ and $b_{1},b_{2},\ldots
,b_{k}\in A$ be some elements. If $a$ commutes with each of the elements
$b_{1},b_{2},\ldots,b_{k}$, then $a$ must also commute with their product
$b_{1}b_{2}\cdots b_{k}$.\ \ \ \ \footnote{The easiest way to prove this is
just by observing that%
\begin{align*}
a\left(  b_{1}b_{2}\cdots b_{k}\right)   &  =\underbrace{ab_{1}}_{=b_{1}%
a}b_{2}b_{3}\cdots b_{k}=b_{1}\underbrace{ab_{2}}_{=b_{2}a}b_{3}\cdots
b_{k}=b_{1}b_{2}\underbrace{ab_{3}}_{=b_{3}a}\cdots b_{k}\\
&  =b_{1}b_{2}b_{3}a\cdots b_{k}=\cdots=b_{1}b_{2}b_{3}\cdots b_{k}a=\left(
b_{1}b_{2}\cdots b_{k}\right)  a.
\end{align*}
This can be easily formalized as an induction on $k$.}
\end{statement}

Applying this fact to $A=\mathcal{A}$ and $a=t_{p,q}$ and $b_{v}%
=\nabla_{\left[  m-1\right]  \cap\operatorname*{Col}\left(  v,T\right)  }^{-}%
$, we conclude that $t_{p,q}$ must commute with the product \newline%
$\nabla_{\left[  m-1\right]  \cap\operatorname*{Col}\left(  1,T\right)  }%
^{-}\nabla_{\left[  m-1\right]  \cap\operatorname*{Col}\left(  2,T\right)
}^{-}\cdots\nabla_{\left[  m-1\right]  \cap\operatorname*{Col}\left(
k,T\right)  }^{-}$ (since $t_{p,q}$ commutes with each of the elements
$\nabla_{\left[  m-1\right]  \cap\operatorname*{Col}\left(  1,T\right)  }%
^{-},\nabla_{\left[  m-1\right]  \cap\operatorname*{Col}\left(  2,T\right)
}^{-},\ldots,\nabla_{\left[  m-1\right]  \cap\operatorname*{Col}\left(
k,T\right)  }^{-}$). In other words, $t_{p,q}$ must commute with $N_{m-1}$
(since \newline$N_{m-1}=\nabla_{\left[  m-1\right]  \cap\operatorname*{Col}%
\left(  1,T\right)  }^{-}\nabla_{\left[  m-1\right]  \cap\operatorname*{Col}%
\left(  2,T\right)  }^{-}\cdots\nabla_{\left[  m-1\right]  \cap
\operatorname*{Col}\left(  k,T\right)  }^{-}$). In other words, we have
$t_{p,q}N_{m-1}=N_{m-1}t_{p,q}$. This proves Claim 9.
\end{proof}

\begin{statement}
\textit{Claim 10:} Let $\ell\in\left\{  i+1,i+2,\ldots,\lambda_{j}%
^{t}\right\}  $. Let $s\in\left\{  j,j+1,\ldots,\lambda_{i}\right\}  $. Then,%
\[
PN_{m-1}t_{T\left(  i,s\right)  ,\ T\left(  \ell,j\right)  }P=PN_{m-1}%
t_{m,\ T\left(  \ell,j\right)  }P.
\]

\end{statement}

\begin{proof}
[Proof of Claim 10.]We have $\left(  \ell,j\right)  \in Y\left(
\lambda\right)  $ (as we saw in the proof of Claim 5). Hence, the entry
$T\left(  \ell,j\right)  $ of $T$ is well-defined. From $\ell\in\left\{
i+1,i+2,\ldots,\lambda_{j}^{t}\right\}  $, we obtain $\ell\geq i+1>i$. From
$s\in\left\{  j,j+1,\ldots,\lambda_{i}\right\}  $, we obtain $s\leq\lambda
_{i}$. Hence, $\left(  i,s\right)  \in Y\left(  \lambda\right)  $. Thus, the
entry $T\left(  i,s\right)  $ of $T$ is well-defined.

Set $p:=T\left(  i,j\right)  =m$ and $q:=T\left(  \ell,j\right)  $ and
$r:=T\left(  i,s\right)  $. From $\ell>i$, we obtain $i\neq\ell$ and thus
$\left(  i,j\right)  \neq\left(  \ell,j\right)  $. Hence, $T\left(
i,j\right)  \neq T\left(  \ell,j\right)  $ (since all entries of $T$ are
distinct). In other words, $p\neq q$ (since $p=T\left(  i,j\right)  $ and
$q=T\left(  \ell,j\right)  $). Furthermore, from $\ell>i$, we obtain $\ell\neq
i$ and thus $\left(  \ell,j\right)  \neq\left(  i,s\right)  $. Hence,
$T\left(  \ell,j\right)  \neq T\left(  i,s\right)  $ (since all entries of $T$
are distinct). In other words, $q\neq r$ (since $q=T\left(  \ell,j\right)  $
and $r=T\left(  i,s\right)  $). Clearly, the numbers $p$ and $r$ both lie in
the $i$-th row of $T$ (since $p=T\left(  i,j\right)  $ and $r=T\left(
i,s\right)  $). Thus, the numbers $p$ and $r$ lie in the same row of $T$.

We must prove that $PN_{m-1}t_{T\left(  i,s\right)  ,\ T\left(  \ell,j\right)
}P=PN_{m-1}t_{m,\ T\left(  \ell,j\right)  }P$. If $T\left(  i,s\right)  =m$,
then this is obvious. Thus, we WLOG assume that $T\left(  i,s\right)  \neq m$.
Hence, $T\left(  i,s\right)  \neq m=T\left(  i,j\right)  $, so that $s\neq j$.
But $s\geq j$ (since $s\in\left\{  j,j+1,\ldots,\lambda_{i}\right\}  $).
Combining this with $s\neq j$, we find $s>j$, so that $j<s$.

Also, $r\neq p$ (since $r=T\left(  i,s\right)  \neq m=p$). Combining this with
$p\neq q$ and $q\neq r$, we see that the three numbers $p,q,r$ are distinct.

The numbers $p$ and $r$ are distinct (since $r\neq p$) and satisfy $p\geq m$
(since $p=m$) and $r\geq m$\ \ \ \ \footnote{\textit{Proof.} The cell $\left(
i,j\right)  $ lies in the same row as the cell $\left(  i,s\right)  $ but
further west (since $j<s$).
\par
In our construction of the tableau $T$, we filled the cells of $Y\left(
\lambda\right)  $ with the entries $1,2,\ldots,n$ in a specific order --
namely in the order in which these cells appeared in the list
(\ref{pf.thm.tableau.hlf.list}). But the list (\ref{pf.thm.tableau.hlf.list})
was listing the cells of $Y\left(  \lambda\right)  $ row by row, where the
cells of each row were listed from left to right. Thus, if a cell $c\in
Y\left(  \lambda\right)  $ lies in the same row as another cell $d\in Y\left(
\lambda\right)  $ but further west, then $c$ appears before $d$ in the list
(\ref{pf.thm.tableau.hlf.list}), and therefore the cell $c$ is filled with a
smaller entry in the tableau $T$ than the cell $d$ (since the cells are filled
with $1,2,\ldots,n$ in the order in which they appear in the list); in other
words, we have $T\left(  c\right)  <T\left(  d\right)  $ in this case.
Applying this to $c=\left(  i,j\right)  $ and $d=\left(  i,s\right)  $, we
conclude that $T\left(  i,j\right)  <T\left(  i,s\right)  $ (since the cell
$\left(  i,j\right)  $ lies in the same row as the cell $\left(  i,s\right)  $
but further west). In other words, $m<r$ (since $m=T\left(  i,j\right)  $ and
$r=T\left(  i,s\right)  $). Thus, $r>m$, so that $r\geq m$.}. Hence,
$t_{p,r}N_{m-1}=N_{m-1}t_{p,r}$ (by Claim 9, applied to $r$ instead of $q$).
In other words, the element $N_{m-1}$ of $\mathcal{A}$ commutes with $t_{p,r}%
$. Thus, Lemma \ref{lem.hlf-young.row2} (applied to $D=Y\left(  \lambda
\right)  $ and $\mathbf{a}=N_{m-1}$) yields
\[
\nabla_{\operatorname*{Row}T}N_{m-1}t_{r,q}\nabla_{\operatorname*{Row}%
T}=\nabla_{\operatorname*{Row}T}N_{m-1}t_{p,q}\nabla_{\operatorname*{Row}T}.
\]
In view of $\nabla_{\operatorname*{Row}T}=P$ and $r=T\left(  i,s\right)  $ and
$q=T\left(  \ell,j\right)  $ and $p=m$, we can rewrite this as%
\[
PN_{m-1}t_{T\left(  i,s\right)  ,\ T\left(  \ell,j\right)  }P=PN_{m-1}%
t_{m,\ T\left(  \ell,j\right)  }P.
\]
Thus, Claim 10 is proved.
\end{proof}

\begin{statement}
\textit{Claim 11:} Let $\ell\in\left\{  i+1,i+2,\ldots,\lambda_{j}%
^{t}\right\}  $. Then,%
\[
PN_{m-1}PN=-\left\vert \operatorname*{Arm}\nolimits_{\lambda}^{+}\left(
c\right)  \right\vert \cdot PN_{m-1}t_{m,\ T\left(  \ell,j\right)  }PN.
\]

\end{statement}

\begin{proof}
[Proof of Claim 11.]From $\left(  i,j\right)  \in Y\left(  \lambda\right)  $,
we obtain $j\in\left[  \lambda_{i}\right]  $, so that $1\leq j\leq\lambda_{i}%
$. We have%
\[
PN_{m-1}\underbrace{\left(  1+\sum_{s=1}^{\lambda_{i}}t_{T\left(  i,s\right)
,\ T\left(  \ell,j\right)  }\right)  PN}_{\substack{=0\\\text{(by Claim 5)}%
}}=0.
\]
Thus,%
\begin{align*}
0  &  =PN_{m-1}\left(  1+\sum_{s=1}^{\lambda_{i}}t_{T\left(  i,s\right)
,\ T\left(  \ell,j\right)  }\right)  PN\\
&  =PN_{m-1}PN+PN_{m-1}\left(  \sum_{s=1}^{\lambda_{i}}t_{T\left(  i,s\right)
,\ T\left(  \ell,j\right)  }\right)  PN.
\end{align*}
Subtracting $PN_{m-1}PN$ from this equality, we find%
\begin{align*}
-PN_{m-1}PN  &  =PN_{m-1}\left(  \sum_{s=1}^{\lambda_{i}}t_{T\left(
i,s\right)  ,\ T\left(  \ell,j\right)  }\right)  PN=\sum_{s=1}^{\lambda_{i}%
}PN_{m-1}t_{T\left(  i,s\right)  ,\ T\left(  \ell,j\right)  }PN\\
&  =\sum_{s=1}^{j-1}P\underbrace{N_{m-1}t_{T\left(  i,s\right)  ,\ T\left(
\ell,j\right)  }P}_{\substack{=0\\\text{(by Claim 8)}}}N+\sum_{s=j}%
^{\lambda_{i}}\underbrace{PN_{m-1}t_{T\left(  i,s\right)  ,\ T\left(
\ell,j\right)  }P}_{\substack{=PN_{m-1}t_{m,\ T\left(  \ell,j\right)
}P\\\text{(by Claim 10)}}}N\\
&  \ \ \ \ \ \ \ \ \ \ \ \ \ \ \ \ \ \ \ \ \left(  \text{since }1\leq
j\leq\lambda_{i}\right) \\
&  =\underbrace{\sum_{s=1}^{j-1}P0N}_{=0}+\underbrace{\sum_{s=j}^{\lambda_{i}%
}PN_{m-1}t_{m,\ T\left(  \ell,j\right)  }PN}_{\substack{=\left(  \lambda
_{i}-j+1\right)  PN_{m-1}t_{m,\ T\left(  \ell,j\right)  }PN\\\text{(since
}j\leq\lambda_{i}\text{)}}}\\
&  =\underbrace{\left(  \lambda_{i}-j+1\right)  }_{\substack{=\left\vert
\operatorname*{Arm}\nolimits_{\lambda}^{+}\left(  c\right)  \right\vert
\\\text{(by Proposition \ref{prop.hlf-young.arm-leg.1} \textbf{(a)})}%
}}PN_{m-1}t_{m,\ T\left(  \ell,j\right)  }PN\\
&  =\left\vert \operatorname*{Arm}\nolimits_{\lambda}^{+}\left(  c\right)
\right\vert \cdot PN_{m-1}t_{m,\ T\left(  \ell,j\right)  }PN.
\end{align*}
Multiplying this equality by $-1$, we obtain%
\[
PN_{m-1}PN=-\left\vert \operatorname*{Arm}\nolimits_{\lambda}^{+}\left(
c\right)  \right\vert \cdot PN_{m-1}t_{m,\ T\left(  \ell,j\right)  }PN.
\]
This proves Claim 11.
\end{proof}

\begin{statement}
\textit{Claim 12:} We have%
\[
\left\vert \operatorname*{Leg}\nolimits_{\lambda}^{-}\left(  c\right)
\right\vert \cdot PN_{m-1}PN=-\left\vert \operatorname*{Arm}\nolimits_{\lambda
}^{+}\left(  c\right)  \right\vert \cdot P\left(  N_{m-1}-N_{m}\right)  PN.
\]

\end{statement}

\begin{proof}
[Proof of Claim 12.]Proposition \ref{prop.hlf-young.arm-leg.1} \textbf{(b)}
yields $\left\vert \operatorname*{Leg}\nolimits_{\lambda}^{-}\left(  c\right)
\right\vert =\lambda_{j}^{t}-i$, so that $\lambda_{j}^{t}-i=\left\vert
\operatorname*{Leg}\nolimits_{\lambda}^{-}\left(  c\right)  \right\vert \geq0$
and therefore $i\leq\lambda_{j}^{t}$. Hence,%
\[
\sum_{\ell=i+1}^{\lambda_{j}^{t}}PN_{m-1}PN=\underbrace{\left(  \lambda
_{j}^{t}-i\right)  }_{=\left\vert \operatorname*{Leg}\nolimits_{\lambda}%
^{-}\left(  c\right)  \right\vert }\cdot\,PN_{m-1}PN=\left\vert
\operatorname*{Leg}\nolimits_{\lambda}^{-}\left(  c\right)  \right\vert \cdot
PN_{m-1}PN.
\]
Thus,%
\begin{align*}
\left\vert \operatorname*{Leg}\nolimits_{\lambda}^{-}\left(  c\right)
\right\vert \cdot PN_{m-1}PN  &  =\sum_{\ell=i+1}^{\lambda_{j}^{t}%
}\underbrace{PN_{m-1}PN}_{\substack{=-\left\vert \operatorname*{Arm}%
\nolimits_{\lambda}^{+}\left(  c\right)  \right\vert \cdot PN_{m-1}%
t_{m,\ T\left(  \ell,j\right)  }PN\\\text{(by Claim 11)}}}\\
&  =\sum_{\ell=i+1}^{\lambda_{j}^{t}}\left(  -\left\vert \operatorname*{Arm}%
\nolimits_{\lambda}^{+}\left(  c\right)  \right\vert \cdot PN_{m-1}%
t_{m,\ T\left(  \ell,j\right)  }PN\right) \\
&  =-\left\vert \operatorname*{Arm}\nolimits_{\lambda}^{+}\left(  c\right)
\right\vert \cdot P\underbrace{N_{m-1}\left(  \sum_{\ell=i+1}^{\lambda_{j}%
^{t}}t_{m,\ T\left(  \ell,j\right)  }\right)  }_{\substack{=N_{m-1}%
-N_{m}\\\text{(by Claim 4)}}}PN\\
&  =-\left\vert \operatorname*{Arm}\nolimits_{\lambda}^{+}\left(  c\right)
\right\vert \cdot P\left(  N_{m-1}-N_{m}\right)  PN.
\end{align*}
This proves Claim 12.
\end{proof}

\begin{statement}
\textit{Claim 13:} We have%
\[
\left\vert H_{\lambda}\left(  c\right)  \right\vert \cdot a_{m-1}=\left\vert
\operatorname*{Arm}\nolimits_{\lambda}^{+}\left(  c\right)  \right\vert \cdot
a_{m}.
\]

\end{statement}

\begin{proof}
[Proof of Claim 13.]The equality (\ref{pf.thm.tableau.hlf.am=}) (with the
index $c$ renamed as $d$) says that%
\begin{equation}
a_{m}=\prod_{d\in T^{-1}\left(  \left[  m\right]  \right)  }\dfrac{\left\vert
H_{\lambda}\left(  d\right)  \right\vert }{\left\vert \operatorname*{Arm}%
\nolimits_{\lambda}^{+}\left(  d\right)  \right\vert }.
\label{pf.thm.tableau.hlf.c13.pf.m}%
\end{equation}
The same argument (applied to $m-1$ instead of $m$) yields%
\begin{equation}
a_{m-1}=\prod_{d\in T^{-1}\left(  \left[  m-1\right]  \right)  }%
\dfrac{\left\vert H_{\lambda}\left(  d\right)  \right\vert }{\left\vert
\operatorname*{Arm}\nolimits_{\lambda}^{+}\left(  d\right)  \right\vert }.
\label{pf.thm.tableau.hlf.c13.pf.m-1}%
\end{equation}
But we can easily see that $c\in T^{-1}\left(  \left[  m\right]  \right)
$\ \ \ \ \footnote{\textit{Proof.} We have $T\left(  c\right)  =m\in\left[
m\right]  $. Thus, $c\in T^{-1}\left(  \left[  m\right]  \right)  $.} and
$T^{-1}\left(  \left[  m\right]  \right)  \setminus\left\{  c\right\}
=T^{-1}\left(  \left[  m-1\right]  \right)  $\ \ \ \ \footnote{\textit{Proof.}
We know that $T$ is an $n$-tableau of shape $Y\left(  \lambda\right)  $, thus
a bijection from $Y\left(  \lambda\right)  $ to $\left[  n\right]  $. Hence,
$T^{-1}\left(  X\setminus Y\right)  =T^{-1}\left(  X\right)  \setminus
T^{-1}\left(  Y\right)  $ for any two subsets $X$ and $Y$ of $\left[
n\right]  $ (indeed, this is true for any map, not just for a bijection).
Applying this to $X=\left[  m\right]  $ and $Y=\left\{  m\right\}  $, we
obtain%
\[
T^{-1}\left(  \left[  m\right]  \setminus\left\{  m\right\}  \right)
=T^{-1}\left(  \left[  m\right]  \right)  \setminus T^{-1}\left(  \left\{
m\right\}  \right)  .
\]
\par
However, $T$ is a bijection from $Y\left(  \lambda\right)  $ to $\left[
n\right]  $. Thus, from $T\left(  c\right)  =m$, we obtain $T^{-1}\left(
m\right)  =c$ and $T^{-1}\left(  \left\{  m\right\}  \right)  =\left\{
T^{-1}\left(  m\right)  \right\}  =\left\{  c\right\}  $ (since $T^{-1}\left(
m\right)  =c$). Thus,%
\[
T^{-1}\left(  \left[  m\right]  \setminus\left\{  m\right\}  \right)
=T^{-1}\left(  \left[  m\right]  \right)  \setminus\underbrace{T^{-1}\left(
\left\{  m\right\}  \right)  }_{=\left\{  c\right\}  }=T^{-1}\left(  \left[
m\right]  \right)  \setminus\left\{  c\right\}  .
\]
Therefore,%
\[
T^{-1}\left(  \left[  m\right]  \right)  \setminus\left\{  c\right\}
=T^{-1}\left(  \underbrace{\left[  m\right]  \setminus\left\{  m\right\}
}_{=\left[  m-1\right]  }\right)  =T^{-1}\left(  \left[  m-1\right]  \right)
.
\]
}.

Now, (\ref{pf.thm.tableau.hlf.c13.pf.m}) becomes%
\begin{align*}
a_{m}  &  =\prod_{d\in T^{-1}\left(  \left[  m\right]  \right)  }%
\dfrac{\left\vert H_{\lambda}\left(  d\right)  \right\vert }{\left\vert
\operatorname*{Arm}\nolimits_{\lambda}^{+}\left(  d\right)  \right\vert
}=\dfrac{\left\vert H_{\lambda}\left(  c\right)  \right\vert }{\left\vert
\operatorname*{Arm}\nolimits_{\lambda}^{+}\left(  c\right)  \right\vert }%
\cdot\prod_{d\in T^{-1}\left(  \left[  m\right]  \right)  \setminus\left\{
c\right\}  }\dfrac{\left\vert H_{\lambda}\left(  d\right)  \right\vert
}{\left\vert \operatorname*{Arm}\nolimits_{\lambda}^{+}\left(  d\right)
\right\vert }\\
&  \ \ \ \ \ \ \ \ \ \ \ \ \ \ \ \ \ \ \ \ \left(
\begin{array}
[c]{c}%
\text{here, we have split off the factor for }d=c\\
\text{from the product, since }c\in T^{-1}\left(  \left[  m\right]  \right)
\end{array}
\right) \\
&  =\dfrac{\left\vert H_{\lambda}\left(  c\right)  \right\vert }{\left\vert
\operatorname*{Arm}\nolimits_{\lambda}^{+}\left(  c\right)  \right\vert }%
\cdot\underbrace{\prod_{d\in T^{-1}\left(  \left[  m-1\right]  \right)
}\dfrac{\left\vert H_{\lambda}\left(  d\right)  \right\vert }{\left\vert
\operatorname*{Arm}\nolimits_{\lambda}^{+}\left(  d\right)  \right\vert }%
}_{\substack{=a_{m-1}\\\text{(by (\ref{pf.thm.tableau.hlf.c13.pf.m-1}))}}}\\
&  \ \ \ \ \ \ \ \ \ \ \ \ \ \ \ \ \ \ \ \ \left(  \text{since }T^{-1}\left(
\left[  m\right]  \right)  \setminus\left\{  c\right\}  =T^{-1}\left(  \left[
m-1\right]  \right)  \right) \\
&  =\dfrac{\left\vert H_{\lambda}\left(  c\right)  \right\vert }{\left\vert
\operatorname*{Arm}\nolimits_{\lambda}^{+}\left(  c\right)  \right\vert }\cdot
a_{m-1}.
\end{align*}
Multiplying this equality by $\left\vert \operatorname*{Arm}\nolimits_{\lambda
}^{+}\left(  c\right)  \right\vert $, we find
\[
\left\vert \operatorname*{Arm}\nolimits_{\lambda}^{+}\left(  c\right)
\right\vert \cdot a_{m}=\left\vert H_{\lambda}\left(  c\right)  \right\vert
\cdot a_{m-1}.
\]
This proves Claim 13.
\end{proof}

Now, Claim 12 yields%
\begin{align*}
\left\vert \operatorname*{Leg}\nolimits_{\lambda}^{-}\left(  c\right)
\right\vert \cdot PN_{m-1}PN  &  =-\left\vert \operatorname*{Arm}%
\nolimits_{\lambda}^{+}\left(  c\right)  \right\vert \cdot\underbrace{P\left(
N_{m-1}-N_{m}\right)  PN}_{=PN_{m-1}PN-PN_{m}PN}\\
&  =-\left\vert \operatorname*{Arm}\nolimits_{\lambda}^{+}\left(  c\right)
\right\vert \cdot\left(  PN_{m-1}PN-PN_{m}PN\right) \\
&  =-\left\vert \operatorname*{Arm}\nolimits_{\lambda}^{+}\left(  c\right)
\right\vert \cdot PN_{m-1}PN+\left\vert \operatorname*{Arm}\nolimits_{\lambda
}^{+}\left(  c\right)  \right\vert \cdot PN_{m}PN.
\end{align*}
Solving this equality for $\left\vert \operatorname*{Arm}\nolimits_{\lambda
}^{+}\left(  c\right)  \right\vert \cdot PN_{m}PN$, we obtain%
\begin{align}
\left\vert \operatorname*{Arm}\nolimits_{\lambda}^{+}\left(  c\right)
\right\vert \cdot PN_{m}PN  &  =\left\vert \operatorname*{Arm}%
\nolimits_{\lambda}^{+}\left(  c\right)  \right\vert \cdot PN_{m-1}%
PN+\left\vert \operatorname*{Leg}\nolimits_{\lambda}^{-}\left(  c\right)
\right\vert \cdot PN_{m-1}PN\nonumber\\
&  =\underbrace{\left(  \left\vert \operatorname*{Arm}\nolimits_{\lambda}%
^{+}\left(  c\right)  \right\vert +\left\vert \operatorname*{Leg}%
\nolimits_{\lambda}^{-}\left(  c\right)  \right\vert \right)  }%
_{\substack{=\left\vert H_{\lambda}\left(  c\right)  \right\vert \\\text{(by
Proposition \ref{prop.hlf-young.arm-leg.1} \textbf{(c)})}}}\cdot
\,PN_{m-1}PN\nonumber\\
&  =\left\vert H_{\lambda}\left(  c\right)  \right\vert \cdot
\underbrace{PN_{m-1}PN}_{\substack{=a_{m-1}P^{2}N\\\text{(since we assumed (as
our induction hypothesis)}\\\text{that (\ref{pf.thm.tableau.hlf.main}) holds
for }m-1\text{ instead of }m\text{)}}}\nonumber\\
&  =\underbrace{\left\vert H_{\lambda}\left(  c\right)  \right\vert \cdot
a_{m-1}}_{\substack{=\left\vert \operatorname*{Arm}\nolimits_{\lambda}%
^{+}\left(  c\right)  \right\vert \cdot a_{m}\\\text{(by Claim 13)}}%
}P^{2}N\nonumber\\
&  =\left\vert \operatorname*{Arm}\nolimits_{\lambda}^{+}\left(  c\right)
\right\vert \cdot a_{m}P^{2}N. \label{pf.thm.tableau.hlf.at}%
\end{align}
But the set $\operatorname*{Arm}\nolimits_{\lambda}^{+}\left(  c\right)  $ is
nonempty (since it contains $c$ by its definition). Hence, its size
$\left\vert \operatorname*{Arm}\nolimits_{\lambda}^{+}\left(  c\right)
\right\vert $ is a positive integer. Thus, we can divide the equality
(\ref{pf.thm.tableau.hlf.at}) by $\left\vert \operatorname*{Arm}%
\nolimits_{\lambda}^{+}\left(  c\right)  \right\vert $ (since our base ring
$\mathbf{k}$ is $\mathbb{Q}$ and thus contains the rational number $\dfrac
{1}{\left\vert \operatorname*{Arm}\nolimits_{\lambda}^{+}\left(  c\right)
\right\vert }$). We thus obtain%
\[
PN_{m}PN=a_{m}P^{2}N.
\]
Thus, we have shown that (\ref{pf.thm.tableau.hlf.main}) holds for $m$. This
completes the induction step. Hence, (\ref{pf.thm.tableau.hlf.main}) is
proved. As explained above, Theorem \ref{thm.tableau.hlf} now follows.
\end{proof}

\section{\label{chp.bas}Bases of $\mathbf{k}\left[  S_{n}\right]  $}

In this chapter, we will systematically study various bases of the
$\mathbf{k}$-module $\mathcal{A}=\mathbf{k}\left[  S_{n}\right]  $. This
includes revisiting the bases we have already encountered:

\begin{itemize}
\item the standard basis $\left(  w\right)  _{w\in S_{n}}$,

\item the Young symmetrizer basis $\left(  \mathbf{E}_{U,V}\right)
_{\lambda\text{ is a partition of }n\text{, and }U,V\in\operatorname*{SYT}%
\left(  \lambda\right)  }$ from Corollary \ref{cor.specht.A.nat-basis},

\item the basis $\left(  \mathbf{F}_{U}w_{U,V}\mathbf{E}_{V}\right)
_{\lambda\text{ is a partition of }n\text{, and }U,V\in\operatorname*{SYT}%
\left(  \lambda\right)  }$ from Theorem \ref{thm.specht.FwE.basis},
\end{itemize}

\noindent but also introducing further important bases, including the
\textit{Murphy bases}.

\subsection{\label{sec.bas.not}Notations}

We will use the following notations throughout this chapter:

\begin{definition}
\label{def.bas.not}\textbf{(a)} We let $\mathcal{A}$ be the $\mathbf{k}%
$-algebra $\mathbf{k}\left[  S_{n}\right]  $. \medskip

\textbf{(b)} We write \textquotedblleft$\lambda\vdash n$\textquotedblright%
\ for \textquotedblleft$\lambda$ is a partition of $n$\textquotedblright.
\medskip

\textbf{(c)} For any partition $\lambda$ of $n$, we let $\operatorname*{SYT}%
\left(  \lambda\right)  $ be the set of all standard $n$-tableaux of shape
$Y\left(  \lambda\right)  $ (that is, of shape $\lambda$). \medskip

\textbf{(d)} For any partition $\lambda$ of $n$, we let $f^{\lambda}$ denote
the \# of standard tableaux of shape $Y\left(  \lambda\right)  $. \medskip

\textbf{(e)} A \emph{standard }$n$\emph{-bitableau} shall mean a triple
$\left(  \lambda,U,V\right)  $, where $\lambda$ is a partition of $n$ and
where $U$ and $V$ are two standard $n$-tableaux of shape $\lambda$. (In
particular, $U$ and $V$ must be straight-shaped and of the same shape.)\footnotemark

We let $\operatorname*{SBT}\left(  n\right)  $ be the set of all standard
$n$-bitableaux. \medskip

\textbf{(f)} We recall the sign-twist $T_{\operatorname*{sign}}:\mathbf{k}%
\left[  S_{n}\right]  \rightarrow\mathbf{k}\left[  S_{n}\right]  $ introduced
in Definition \ref{def.Tsign.Tsign}. This is a $\mathbf{k}$-linear map from
$\mathcal{A}$ to $\mathcal{A}$ (since $\mathcal{A}=\mathbf{k}\left[
S_{n}\right]  $), and sends each $w\in S_{n}$ to $\left(  -1\right)  ^{w}w$
(by its definition). \medskip

\textbf{(g)} We recall the antipode $S:\mathbf{k}\left[  S_{n}\right]
\rightarrow\mathbf{k}\left[  S_{n}\right]  $ introduced in Definition
\ref{def.S.S}. This is a $\mathbf{k}$-linear map from $\mathcal{A}$ to
$\mathcal{A}$ (since $\mathcal{A}=\mathbf{k}\left[  S_{n}\right]  $), and
sends each $w\in S_{n}$ to $w^{-1}$ (by its definition). \medskip

\textbf{(h)} The symbol \textquotedblleft$\operatorname*{span}$%
\textquotedblright\ shall mean \textquotedblleft$\operatorname*{span}%
\nolimits_{\mathbf{k}}$\textquotedblright\ by default.
\end{definition}

\footnotetext{Note that the $\lambda$ in a standard $n$-bitableau $\left(
\lambda,U,V\right)  $ can be easily reconstructed from $U$ (indeed, each entry
$\lambda_{i}$ of $\lambda$ is simply the number of entries in row $i$ of $U$).
But we still include it in the $n$-bitableau for convenience.}The only one of
these notations that is new is the concept of a standard $n$-bitableau. Let us
give an example:

\begin{example}
\label{exa.sbt.sbt3}The standard $3$-bitableaux are the six triples
\begin{align*}
&  \left(  \left(  3\right)
,\ \ \ytableaushort{123}\ ,\ \ \ytableaushort{123}\right)  ,\\
&  \left(  \left(  2,1\right)
,\ \ \ytableaushort{12,3}\ ,\ \ \ytableaushort{12,3}\right)
,\ \ \ \ \ \ \ \ \ \ \left(  \left(  2,1\right)
,\ \ \ytableaushort{12,3}\ ,\ \ \ytableaushort{13,2}\right)  ,\\
&  \left(  \left(  2,1\right)
,\ \ \ytableaushort{13,2}\ ,\ \ \ytableaushort{12,3}\right)
,\ \ \ \ \ \ \ \ \ \ \left(  \left(  2,1\right)
,\ \ \ytableaushort{13,2}\ ,\ \ \ytableaushort{13,2}\right)  ,\\
&  \left(  \left(  1,1,1\right)
,\ \ \ytableaushort{1,2,3}\ ,\ \ \ytableaushort{1,2,3}\right)  .
\end{align*}
Using the shorthand notation from Convention \ref{conv.tableau.poetic}, we can
rewrite them as%
\begin{align*}
&  \left(  \left(  3\right)  ,\ 123,\ 123\right)  ,\\
&  \left(  \left(  2,1\right)  ,\ 12\backslash\backslash3,\ 12\backslash
\backslash3\right)  ,\ \ \ \ \ \ \ \ \ \ \left(  \left(  2,1\right)
,\ 12\backslash\backslash3,\ 13\backslash\backslash2\right)  ,\\
&  \left(  \left(  2,1\right)  ,\ 13\backslash\backslash2,\ 12\backslash
\backslash3\right)  ,\ \ \ \ \ \ \ \ \ \ \left(  \left(  2,1\right)
,\ 13\backslash\backslash2,\ 13\backslash\backslash2\right)  ,\\
&  \left(  \left(  1,1,1\right)  ,\ 1\backslash\backslash2\backslash
\backslash3,\ 1\backslash\backslash2\backslash\backslash3\right)  .
\end{align*}

\end{example}

We also recall Definition \ref{def.symmetrizers.symmetrizers}, which
introduces the row symmetrizer%
\[
\nabla_{\operatorname*{Row}T}:=\sum_{w\in\mathcal{R}\left(  T\right)  }%
w\in\mathbf{k}\left[  S_{n}\right]  =\mathcal{A}%
\]
and the column antisymmetrizer%
\[
\nabla_{\operatorname*{Col}T}^{-}:=\sum_{w\in\mathcal{C}\left(  T\right)
}\left(  -1\right)  ^{w}w\in\mathbf{k}\left[  S_{n}\right]  =\mathcal{A}%
\]
for each $n$-tableau $T$. Here, $\mathcal{R}\left(  T\right)  $ and
$\mathcal{C}\left(  T\right)  $ are two subgroups of $S_{n}$, defined in
Definition \ref{def.tableau.cell-cts}. (In a nutshell: $\mathcal{R}\left(
T\right)  $ consists of those $w\in S_{n}$ such that each $i\in\left[
n\right]  $ lies in the same row of $T$ as $w\left(  i\right)  $ does, whereas
$\mathcal{C}\left(  T\right)  $ consists of those $w\in S_{n}$ such that each
$i\in\left[  n\right]  $ lies in the same column of $T$ as $w\left(  i\right)
$ does.)

\subsection{\label{sec.bas.what}What to do with a basis}

\subsubsection{\label{subsec.bas.what.size}The size of a basis}

The $\mathbf{k}$-module $\mathcal{A}=\mathbf{k}\left[  S_{n}\right]  $ has a
basis $\left(  w\right)  _{w\in S_{n}}$, and thus is free of rank $\left\vert
S_{n}\right\vert =n!$. At least when $\mathbf{k}$ is a field, this entails
that \textbf{every} basis of $\mathcal{A}$ has size $n!$ (in fact, this is
even true for any nontrivial commutative ring $\mathbf{k}$; see e.g.
\url{https://math.stackexchange.com/questions/2651777/} for a proof). Thus it
should not come as a surprise that all bases of $\mathcal{A}$ we shall discuss
have size $n!$. However, not all of them are indexed by the set $S_{n}$.
Another possible index set is $\operatorname*{SBT}\left(  n\right)  $, because
of the following simple fact:

\begin{lemma}
\label{lem.sbt.sizen!}We have $\left\vert \operatorname*{SBT}\left(  n\right)
\right\vert =n!$.
\end{lemma}

\begin{proof}
This is just a restatement of Corollary \ref{cor.spechtmod.sumflam2}, since
the definition of standard $n$-bitableaux shows that $\left\vert
\operatorname*{SBT}\left(  n\right)  \right\vert =\sum_{\lambda\text{ is a
partition of }n}\left(  f^{\lambda}\right)  ^{2}$. For more details, see the
appendix (Section \ref{sec.details.bas.what}).
\end{proof}

All our bases of $\mathcal{A}$ will have indexing set $S_{n}$ or
$\operatorname*{SBT}\left(  n\right)  $.

\subsubsection{\label{subsec.bas.what.mtx}The matrix representing a linear
map}

Bases of $\mathbf{k}$-modules allow us to represent linear maps by matrices.
Let us recall how this is defined:\footnote{Recall the notion of a $P\times
Q$-matrix; this was defined in Section \ref{sec.spechtmod.detava}.}

\begin{definition}
\label{def.bas.repmat} Let $V$ and $W$ be two free $\mathbf{k}$-modules with
bases $\overrightarrow{v}=\left(  v_{q}\right)  _{q\in Q}$ and
$\overrightarrow{w}=\left(  w_{p}\right)  _{p\in P}$, respectively. Let
$f:V\rightarrow W$ be a $\mathbf{k}$-linear map. Then, $\left[  f\right]
_{\overrightarrow{v}\rightarrow\overrightarrow{w}}$ shall denote the matrix
representing this map $f$ with respect to the bases $\overrightarrow{v}$ and
$\overrightarrow{w}$. Explicitly, this is the $P\times Q$-matrix $\left(
a_{p,q}\right)  _{\left(  p,q\right)  \in P\times Q}$, whose entries $a_{p,q}$
are determined by the requirement that%
\[
f\left(  v_{q}\right)  =\sum_{p\in P}a_{p,q}w_{p}\ \ \ \ \ \ \ \ \ \ \text{for
all }q\in Q.
\]
Thus, each column of this matrix $\left[  f\right]  _{\overrightarrow{v}%
\rightarrow\overrightarrow{w}}$ corresponds to some index $q\in Q$, and its
entries are the coordinates of $f\left(  v_{q}\right)  $ with respect to the
basis $\overrightarrow{w}$ of $W$.
\end{definition}

We note that a $P\times Q$-matrix has its rows indexed by elements of $P$ and
its columns indexed by elements of $Q$. If the sets $P$ and $Q$ are equipped
with total orders, then we can arrange the rows from top to bottom (by
following the order on $P$) and arrange the columns from left to right (by
following the order on $Q$). This is how we will be writing down such matrices
in this chapter. For instance, an $\left\{  \alpha,\beta\right\}
\times\left\{  \varphi,\psi\right\}  $-matrix $\left(  a_{p,q}\right)
_{\left(  p,q\right)  \in\left\{  \alpha,\beta\right\}  \times\left\{
\varphi,\psi\right\}  }$ can be written down as $\left(
\begin{array}
[c]{cc}%
a_{\alpha,\varphi} & a_{\alpha,\psi}\\
a_{\beta,\varphi} & a_{\beta,\psi}%
\end{array}
\right)  $ if we order its indexing sets by $\alpha<\beta$ and $\varphi<\psi$.

Representing matrices have the following property: If $U,V,W$ are three free
$\mathbf{k}$-modules with respective bases $\overrightarrow{u}%
,\overrightarrow{v},\overrightarrow{w}$, and if $f:V\rightarrow W$ and
$g:U\rightarrow V$ are two $\mathbf{k}$-linear maps, then%
\begin{equation}
\left[  f\circ g\right]  _{\overrightarrow{u}\rightarrow\overrightarrow{w}%
}=\left[  f\right]  _{\overrightarrow{v}\rightarrow\overrightarrow{w}}%
\cdot\left[  g\right]  _{\overrightarrow{u}\rightarrow\overrightarrow{v}}.
\label{eq.bas.what.mtx.fg}%
\end{equation}

As a particular case of Definition~\ref{def.bas.repmat}, if
$\overrightarrow{\mathbf{b}}$ is any basis of $\mathcal{A}=\mathbf{k}\left[
S_{n}\right]  $, and if $f:\mathcal{A}\rightarrow\mathcal{A}$ is any
$\mathbf{k}$-linear map, then there exists a matrix $\left[  f\right]
_{\overrightarrow{\mathbf{b}}\rightarrow\overrightarrow{\mathbf{b}}}$
representing $f$ with respect to the basis $\overrightarrow{\mathbf{b}}$. In
particular, we can apply this to the $\mathbf{k}$-linear maps $f=S$ and
$f=T_{\operatorname*{sign}}$. The resulting matrices $\left[  f\right]
_{\overrightarrow{\mathbf{b}}\rightarrow\overrightarrow{\mathbf{b}}}$ will be
$S_{n}\times S_{n}$-matrices or $\operatorname*{SBT}\left(  n\right)
\times\operatorname*{SBT}\left(  n\right)  $-matrices depending on the
indexing set of our basis being $S_{n}$ or $\operatorname*{SBT}\left(
n\right)  $; in either case, they will have $n!$ rows and $n!$ columns.

\subsubsection{\label{subsec.bas.what.LR}The endomorphisms $L\left(
\mathbf{a}\right)  $ and $R\left(  \mathbf{a}\right)  $ and the matrices
$L_{\protect\overrightarrow{\mathbf{b}}}\left(  \mathbf{a}\right)  $ and
$R_{\protect\overrightarrow{\mathbf{b}}}\left(  \mathbf{a}\right)  $}

The following is a basic and general construction in algebra:

\begin{definition}
\label{def.LRmul.LR}Let $A$ be a $\mathbf{k}$-algebra. Let $\mathbf{a}\in A$.
Then, we define two $\mathbf{k}$-linear maps%
\begin{align*}
L\left(  \mathbf{a}\right)  :A  &  \rightarrow A,\\
\mathbf{x}  &  \mapsto\mathbf{ax}%
\end{align*}
and%
\begin{align*}
R\left(  \mathbf{a}\right)  :A  &  \rightarrow A,\\
\mathbf{x}  &  \mapsto\mathbf{xa}.
\end{align*}
These maps $L\left(  \mathbf{a}\right)  $ and $R\left(  \mathbf{a}\right)  $
are called \emph{left multiplication by }$\mathbf{a}$ and \emph{right
multiplication by }$\mathbf{a}$, respectively.
\end{definition}

\begin{proposition}
\label{prop.LRmul.alg}Let $A$ be a $\mathbf{k}$-algebra. Then: \medskip

\textbf{(a)} We have $L\left(  1\right)  =R\left(  1\right)
=\operatorname*{id}\nolimits_{A}$. \medskip

\textbf{(b)} For any $\mathbf{a},\mathbf{b}\in A$, we have $L\left(
\mathbf{a}+\mathbf{b}\right)  =L\left(  \mathbf{a}\right)  +L\left(
\mathbf{b}\right)  $ and $R\left(  \mathbf{a}+\mathbf{b}\right)  =R\left(
\mathbf{a}\right)  +R\left(  \mathbf{b}\right)  $. \medskip

\textbf{(c)} For any $\mathbf{a},\mathbf{b}\in A$, we have $L\left(
\mathbf{ab}\right)  =L\left(  \mathbf{a}\right)  \circ L\left(  \mathbf{b}%
\right)  $ and $R\left(  \mathbf{ab}\right)  =R\left(  \mathbf{b}\right)
\circ R\left(  \mathbf{a}\right)  $. \medskip

\textbf{(d)} For any $\lambda\in\mathbf{k}$ and $\mathbf{a}\in A$, we have
$L\left(  \lambda\mathbf{a}\right)  =\lambda L\left(  \mathbf{a}\right)  $ and
$R\left(  \lambda\mathbf{a}\right)  =\lambda R\left(  \mathbf{a}\right)  $.
\end{proposition}

\begin{proof}
Straightforward. For instance, in order to verify the equality $R\left(
\mathbf{ab}\right)  =R\left(  \mathbf{b}\right)  \circ R\left(  \mathbf{a}%
\right)  $ in part \textbf{(c)}, we need to show that $\left(  R\left(
\mathbf{ab}\right)  \right)  \left(  \mathbf{x}\right)  =\left(  R\left(
\mathbf{b}\right)  \circ R\left(  \mathbf{a}\right)  \right)  \left(
\mathbf{x}\right)  $ for each $\mathbf{x}\in A$. But this follows from
computing both sides:%
\[
\left(  R\left(  \mathbf{ab}\right)  \right)  \left(  \mathbf{x}\right)
=\mathbf{x}\left(  \mathbf{ab}\right)  \ \ \ \ \ \ \ \ \ \ \left(  \text{by
the definition of }R\left(  \mathbf{ab}\right)  \right)
\]
and%
\begin{align*}
\left(  R\left(  \mathbf{b}\right)  \circ R\left(  \mathbf{a}\right)  \right)
\left(  \mathbf{x}\right)   &  =\left(  R\left(  \mathbf{b}\right)  \right)
\left(  \left(  R\left(  \mathbf{a}\right)  \right)  \left(  \mathbf{x}%
\right)  \right) \\
&  =\left(  R\left(  \mathbf{b}\right)  \right)  \left(  \mathbf{xa}\right)
\ \ \ \ \ \ \ \ \ \ \left(
\begin{array}
[c]{c}%
\text{since the definition of }R\left(  \mathbf{a}\right) \\
\text{yields }\left(  R\left(  \mathbf{a}\right)  \right)  \left(
\mathbf{x}\right)  =\mathbf{xa}%
\end{array}
\right) \\
&  =\left(  \mathbf{xa}\right)  \mathbf{b}\ \ \ \ \ \ \ \ \ \ \left(  \text{by
the definition of }R\left(  \mathbf{b}\right)  \right) \\
&  =\mathbf{x}\left(  \mathbf{ab}\right)  \ \ \ \ \ \ \ \ \ \ \left(  \text{by
associativity}\right)  .
\end{align*}
Comparing these two equalities, we obtain $\left(  R\left(  \mathbf{ab}%
\right)  \right)  \left(  \mathbf{x}\right)  =\left(  R\left(  \mathbf{b}%
\right)  \circ R\left(  \mathbf{a}\right)  \right)  \left(  \mathbf{x}\right)
$, as desired.
\end{proof}

\begin{corollary}
\label{cor.LRmul.mor}Let $A$ be a $\mathbf{k}$-algebra. Then: \medskip

\textbf{(a)} The map%
\begin{align*}
L:A  &  \rightarrow\operatorname*{End}\nolimits_{\mathbf{k}}A,\\
\mathbf{a}  &  \mapsto L\left(  \mathbf{a}\right)
\end{align*}
is an injective $\mathbf{k}$-algebra morphism. \medskip

\textbf{(b)} The map%
\begin{align*}
R:A  &  \rightarrow\operatorname*{End}\nolimits_{\mathbf{k}}A,\\
\mathbf{a}  &  \mapsto R\left(  \mathbf{a}\right)
\end{align*}
is an injective $\mathbf{k}$-algebra anti-morphism.
\end{corollary}

\begin{proof}
Proposition \ref{prop.LRmul.alg} shows that $L$ is a $\mathbf{k}$-algebra
morphism and that $R$ is a $\mathbf{k}$-algebra anti-morphism. It remains to
prove that both $L$ and $R$ are injective. Observe that each $\mathbf{a}\in A$
satisfies%
\begin{align*}
\left(  L\left(  \mathbf{a}\right)  \right)  \left(  1_{A}\right)   &
=\mathbf{a}\cdot1_{A}\ \ \ \ \ \ \ \ \ \ \left(  \text{by the definition of
}L\left(  \mathbf{a}\right)  \right) \\
&  =\mathbf{a};
\end{align*}
this equality allows us to reconstruct $\mathbf{a}$ from $L\left(
\mathbf{a}\right)  $. Hence, the map $L$ is injective. Similarly, we can show
that $R$ is injective (by observing that each $\mathbf{a}\in A$ satisfies
$\left(  R\left(  \mathbf{a}\right)  \right)  \left(  1_{A}\right)
=\mathbf{a}$). Thus, the proof of Corollary \ref{cor.LRmul.mor} is complete.
\end{proof}

We can apply Definition \ref{def.LRmul.LR}, Proposition \ref{prop.LRmul.alg}
and Corollary \ref{cor.LRmul.mor} to the case $A=\mathcal{A}$; thus we obtain
an injective $\mathbf{k}$-algebra morphism $L:\mathcal{A}\rightarrow
\operatorname*{End}\nolimits_{\mathbf{k}}\mathcal{A}$ and an injective
$\mathbf{k}$-algebra anti-morphism $R:\mathcal{A}\rightarrow
\operatorname*{End}\nolimits_{\mathbf{k}}\mathcal{A}$. Given an element
$\mathbf{a}\in\mathcal{A}$, we thus obtain two $\mathbf{k}$-linear
endomorphisms $L\left(  \mathbf{a}\right)  $ and $R\left(  \mathbf{a}\right)
$ of $\mathcal{A}$. If we are furthermore given a basis
$\overrightarrow{\mathbf{b}}$ of $\mathcal{A}$, then these endomorphisms
$L\left(  \mathbf{a}\right)  $ and $R\left(  \mathbf{a}\right)  $ are
represented by matrices:

\begin{definition}
\label{def.LRmul.mats}Let $\overrightarrow{\mathbf{b}}=\left(  \mathbf{b}%
_{i}\right)  _{i\in I}$ be a basis of $\mathcal{A}=\mathbf{k}\left[
S_{n}\right]  $, and let $\mathbf{a}\in\mathcal{A}$. Then, we define the two
$I\times I$-matrices%
\[
L_{\overrightarrow{\mathbf{b}}}\left(  \mathbf{a}\right)  :=\left[  L\left(
\mathbf{a}\right)  \right]  _{\overrightarrow{\mathbf{b}}\rightarrow
\overrightarrow{\mathbf{b}}}\ \ \ \ \ \ \ \ \ \ \text{and}%
\ \ \ \ \ \ \ \ \ \ R_{\overrightarrow{\mathbf{b}}}\left(  \mathbf{a}\right)
:=\left[  R\left(  \mathbf{a}\right)  \right]  _{\overrightarrow{\mathbf{b}%
}\rightarrow\overrightarrow{\mathbf{b}}}.
\]
These are the matrices representing the endomorphisms $L\left(  \mathbf{a}%
\right)  $ and $R\left(  \mathbf{a}\right)  $ with respect to the basis
$\overrightarrow{\mathbf{b}}$.
\end{definition}

The main advantage of these matrices compared to the endomorphisms $L\left(
\mathbf{a}\right)  $ and $R\left(  \mathbf{a}\right)  $ themselves is, of
course, that they are concrete objects that we can just write down. Moreover,
we need not write them down for all $\mathbf{a}\in\mathcal{A}$; it suffices to
know them for $\mathbf{a}$ being the simple transpositions $s_{1},s_{2}%
,\ldots,s_{n-1}$ as defined in (\ref{eq.intro.perms.cycs.si=}). Here is why:

\begin{itemize}
\item The simple transpositions $s_{1},s_{2},\ldots,s_{n-1}$ generate the
symmetric group $S_{n}$ and thus also its group algebra $\mathcal{A}$. Thus,
each element of $\mathcal{A}$ is a linear combination of products of these
simple transpositions.

\item But the matrices $L_{\overrightarrow{\mathbf{b}}}\left(  \mathbf{a}%
\right)  $ and $R_{\overrightarrow{\mathbf{b}}}\left(  \mathbf{a}\right)  $
behave nicely under products and linear combinations: For all $\mathbf{a}%
_{1},\mathbf{a}_{2}\in\mathcal{A}$ and $\lambda_{1},\lambda_{2}\in\mathbf{k}$,
we have%
\begin{align*}
L_{\overrightarrow{\mathbf{b}}}\left(  \mathbf{a}_{1}\mathbf{a}_{2}\right)
&  =L_{\overrightarrow{\mathbf{b}}}\left(  \mathbf{a}_{1}\right)  \cdot
L_{\overrightarrow{\mathbf{b}}}\left(  \mathbf{a}_{2}\right)
\ \ \ \ \ \ \ \ \ \ \text{and}\\
R_{\overrightarrow{\mathbf{b}}}\left(  \mathbf{a}_{1}\mathbf{a}_{2}\right)
&  =R_{\overrightarrow{\mathbf{b}}}\left(  \mathbf{a}_{2}\right)  \cdot
R_{\overrightarrow{\mathbf{b}}}\left(  \mathbf{a}_{1}\right)
\ \ \ \ \ \ \ \ \ \ \text{and}\\
L_{\overrightarrow{\mathbf{b}}}\left(  \lambda_{1}\mathbf{a}_{1}+\lambda
_{2}\mathbf{a}_{2}\right)   &  =\lambda_{1}L_{\overrightarrow{\mathbf{b}}%
}\left(  \mathbf{a}_{1}\right)  +\lambda_{2}L_{\overrightarrow{\mathbf{b}}%
}\left(  \mathbf{a}_{2}\right)  \ \ \ \ \ \ \ \ \ \ \text{and}\\
R_{\overrightarrow{\mathbf{b}}}\left(  \lambda_{1}\mathbf{a}_{1}+\lambda
_{2}\mathbf{a}_{2}\right)   &  =\lambda_{1}R_{\overrightarrow{\mathbf{b}}%
}\left(  \mathbf{a}_{1}\right)  +\lambda_{2}R_{\overrightarrow{\mathbf{b}}%
}\left(  \mathbf{a}_{2}\right)
\end{align*}
(thanks to (\ref{eq.bas.what.mtx.fg}) and to Proposition \ref{prop.LRmul.alg}).
\end{itemize}

Hence, knowing these matrices $L_{\overrightarrow{\mathbf{b}}}\left(
\mathbf{a}\right)  $ and $R_{\overrightarrow{\mathbf{b}}}\left(
\mathbf{a}\right)  $ for all $\mathbf{a}\in\left\{  s_{1},s_{2},\ldots
,s_{n-1}\right\}  $ is the key to comfortably computing with a given basis
$\overrightarrow{\mathbf{b}}$ of $\mathcal{A}$.

\subsubsection{\label{subsec.bas.what.sub}Sub-spans}

What else can we do with a basis $\overrightarrow{\mathbf{b}}=\left(
\mathbf{b}_{i}\right)  _{i\in I}$ of $\mathcal{A}$ ? We can consider its
subfamilies $\left(  \mathbf{b}_{i}\right)  _{i\in J}$ for various subsets $J$
of $I$, and try to understand their spans $\operatorname*{span}\left\{
\mathbf{b}_{i}\ \mid\ i\in J\right\}  $. Are some of these spans $\mathbf{k}%
$-subalgebras of $\mathcal{A}$ ? Ideals? Left or right ideals?

In the rest of this chapter, we shall analyze various bases of $\mathcal{A}$
with regard to these and other questions.

\subsection{\label{sec.bas.std}The standard basis of $\mathbf{k}\left[
S_{n}\right]  $, revisited}

The simplest of all bases of $\mathcal{A}$ is the standard basis
$\overrightarrow{\operatorname*{std}}:=\left(  w\right)  _{w\in S_{n}}$, which
consists of all the $n!$ permutations $w\in S_{n}$ (or, if we are pedantic, of
the corresponding basis vectors of $\mathbf{k}\left[  S_{n}\right]  $). The
antipode $S$ and the sign-twist $T_{\operatorname*{sign}}$ act on this basis
by the simple formulas%
\[
S\left(  w\right)  =w^{-1}\ \ \ \ \ \ \ \ \ \ \text{and}%
\ \ \ \ \ \ \ \ \ \ T_{\operatorname*{sign}}\left(  w\right)  =\left(
-1\right)  ^{w}w\ \ \ \ \ \ \ \ \ \ \text{for each }w\in S_{n};
\]
thus, the matrix $\left[  S\right]  _{\overrightarrow{\operatorname*{std}%
}\rightarrow\overrightarrow{\operatorname*{std}}}$ is a permutation matrix,
whereas the matrix $\left[  T_{\operatorname*{sign}}\right]
_{\overrightarrow{\operatorname*{std}}\rightarrow
\overrightarrow{\operatorname*{std}}}$ is diagonal (with diagonal entries $1$
and $-1$). For each permutation $u\in S_{n}$, the matrices
$L_{\overrightarrow{\operatorname*{std}}}\left(  u\right)  $ and
$R_{\overrightarrow{\operatorname*{std}}}\left(  u\right)  $ (representing the
linear maps $L\left(  u\right)  $ and $R\left(  u\right)  $) are permutation
matrices as well, since both maps $L\left(  u\right)  $ and $R\left(
u\right)  $ just permute the standard basis
$\overrightarrow{\operatorname*{std}}$.

\begin{example}
For $n=3$, the entries of the standard basis
$\overrightarrow{\operatorname*{std}}=\left(  w\right)  _{w\in S_{n}}$ are the
elements of $S_{3}$, namely $\operatorname*{id},\ s_{1},\ s_{2},\ t_{1,3}%
,\ \operatorname*{cyc}\nolimits_{1,2,3}$ and $\operatorname*{cyc}%
\nolimits_{1,3,2}$. The matrices representing the $\mathbf{k}$-linear maps
$S$, $T_{\operatorname*{sign}}$, $L\left(  s_{1}\right)  $, $L\left(
s_{2}\right)  $, $R\left(  s_{1}\right)  $ and $R\left(  s_{2}\right)  $ are
as follows (where the rows and columns are indexed by the elements of $S_{3}$,
in the order just listed):
\[
\left[  S\right]  _{\overrightarrow{\operatorname*{std}}\rightarrow
\overrightarrow{\operatorname*{std}}}=\left(
\begin{array}
[c]{cccccc}%
1 & 0 & 0 & 0 & 0 & 0\\
0 & 1 & 0 & 0 & 0 & 0\\
0 & 0 & 1 & 0 & 0 & 0\\
0 & 0 & 0 & 1 & 0 & 0\\
0 & 0 & 0 & 0 & 0 & 1\\
0 & 0 & 0 & 0 & 1 & 0
\end{array}
\right)  ;
\]%
\[
\left[  T_{\operatorname*{sign}}\right]  _{\overrightarrow{\operatorname*{std}%
}\rightarrow\overrightarrow{\operatorname*{std}}}=\left(
\begin{array}
[c]{cccccc}%
1 & 0 & 0 & 0 & 0 & 0\\
0 & -1 & 0 & 0 & 0 & 0\\
0 & 0 & -1 & 0 & 0 & 0\\
0 & 0 & 0 & -1 & 0 & 0\\
0 & 0 & 0 & 0 & 1 & 0\\
0 & 0 & 0 & 0 & 0 & 1
\end{array}
\right)  ;
\]%
\[
L_{\overrightarrow{\operatorname*{std}}}\left(  s_{1}\right)  =\left[
L\left(  s_{1}\right)  \right]  _{\overrightarrow{\operatorname*{std}%
}\rightarrow\overrightarrow{\operatorname*{std}}}=\left(
\begin{array}
[c]{cccccc}%
0 & 1 & 0 & 0 & 0 & 0\\
1 & 0 & 0 & 0 & 0 & 0\\
0 & 0 & 0 & 0 & 1 & 0\\
0 & 0 & 0 & 0 & 0 & 1\\
0 & 0 & 1 & 0 & 0 & 0\\
0 & 0 & 0 & 1 & 0 & 0
\end{array}
\right)  ;
\]%
\[
L_{\overrightarrow{\operatorname*{std}}}\left(  s_{2}\right)  =\left[
L\left(  s_{2}\right)  \right]  _{\overrightarrow{\operatorname*{std}%
}\rightarrow\overrightarrow{\operatorname*{std}}}=\left(
\begin{array}
[c]{cccccc}%
0 & 0 & 1 & 0 & 0 & 0\\
0 & 0 & 0 & 0 & 0 & 1\\
1 & 0 & 0 & 0 & 0 & 0\\
0 & 0 & 0 & 0 & 1 & 0\\
0 & 0 & 0 & 1 & 0 & 0\\
0 & 1 & 0 & 0 & 0 & 0
\end{array}
\right)  ;
\]%
\[
R_{\overrightarrow{\operatorname*{std}}}\left(  s_{1}\right)  =\left[
R\left(  s_{1}\right)  \right]  _{\overrightarrow{\operatorname*{std}%
}\rightarrow\overrightarrow{\operatorname*{std}}}=\left(
\begin{array}
[c]{cccccc}%
0 & 1 & 0 & 0 & 0 & 0\\
1 & 0 & 0 & 0 & 0 & 0\\
0 & 0 & 0 & 0 & 0 & 1\\
0 & 0 & 0 & 0 & 1 & 0\\
0 & 0 & 0 & 1 & 0 & 0\\
0 & 0 & 1 & 0 & 0 & 0
\end{array}
\right)  ;
\]%
\[
R_{\overrightarrow{\operatorname*{std}}}\left(  s_{2}\right)  =\left[
R\left(  s_{2}\right)  \right]  _{\overrightarrow{\operatorname*{std}%
}\rightarrow\overrightarrow{\operatorname*{std}}}=\left(
\begin{array}
[c]{cccccc}%
0 & 0 & 1 & 0 & 0 & 0\\
0 & 0 & 0 & 0 & 1 & 0\\
1 & 0 & 0 & 0 & 0 & 0\\
0 & 0 & 0 & 0 & 0 & 1\\
0 & 1 & 0 & 0 & 0 & 0\\
0 & 0 & 0 & 1 & 0 & 0
\end{array}
\right)  .
\]

\end{example}

What about spans of subfamilies of the standard basis
$\overrightarrow{\operatorname*{std}}=\left(  w\right)  _{w\in S_{n}}$ ?
Several of these spans are subalgebras, due to the following simple exercise:

\begin{exercise}
\fbox{1} Let $H$ be a subgroup of a group $G$. Prove that
$\operatorname*{span}\left\{  w\ \mid\ w\in H\right\}  $ is a $\mathbf{k}%
$-subalgebra of $\mathbf{k}\left[  G\right]  $ that is isomorphic to the group
algebra $\mathbf{k}\left[  H\right]  $.
\end{exercise}

Hence, many subfamilies of the standard basis
$\overrightarrow{\operatorname*{std}}$ of $\mathcal{A}$ span $\mathbf{k}%
$-subalgebras of $\mathcal{A}$. But very few of them span ideals, or even left
or right ideals. Indeed, if a left ideal of $\mathcal{A}$ contains any of the
permutations $w\in S_{n}$, then it must also contain $w^{-1}w=1$, and thus
equals the whole algebra $\mathcal{A}$. The same applies to right ideals.

\subsection{\label{sec.bas.ysb}The Young symmetrizer basis of $\mathbf{k}%
\left[  S_{n}\right]  $, revisited}

For the present Section \ref{sec.bas.ysb}, we assume that $n!$ is invertible
in $\mathbf{k}$. Then, Corollary \ref{cor.specht.A.nat-basis} shows that
$\left(  \mathbf{E}_{U,V}\right)  _{\lambda\text{ is a partition of }n\text{,
and }U,V\in\operatorname*{SYT}\left(  \lambda\right)  }$ is a basis of the
$\mathbf{k}$-module $\mathcal{A}$ (where $\mathbf{E}_{U,V}$ is defined as in
Definition \ref{def.specht.EPQ}). We called this basis the \emph{Young
symmetrizer basis}. In view of Definition \ref{def.bas.not} \textbf{(e)}, we
can rewrite the condition \textquotedblleft$\lambda$ is a partition of $n$,
and $U,V\in\operatorname*{SYT}\left(  \lambda\right)  $\textquotedblright\ as
\textquotedblleft$\left(  \lambda,U,V\right)  $ is a standard $n$%
-bitableau\textquotedblright, or -- even shorter -- as \textquotedblleft%
$\left(  \lambda,U,V\right)  \in\operatorname*{SBT}\left(  n\right)
$\textquotedblright. Thus, the Young symmetrizer basis can be written as
\[
\left(  \mathbf{E}_{U,V}\right)  _{\left(  \lambda,U,V\right)  \in
\operatorname*{SBT}\left(  n\right)  }.
\]
Let us denote this basis by $\overrightarrow{\operatorname*{ysb}}$.

Describing how the operators $S$, $T_{\operatorname*{sign}}$, $L\left(
\mathbf{a}\right)  $ and $R\left(  \mathbf{a}\right)  $ (for $\mathbf{a}%
\in\mathbf{k}\left[  S_{n}\right]  $) act on this basis is not easy. At least
Proposition \ref{prop.specht.EPQ.r} shows that the composition
$T_{\operatorname*{sign}}\circ S$ acts in a rather simple way (namely, sending
each $\mathbf{E}_{P,Q}$ to $\left(  -1\right)  ^{w_{Q,P}}\mathbf{E}%
_{Q\mathbf{r},P\mathbf{r}}$). No such nice formula exists for the operators
$S$ and $T_{\operatorname*{sign}}$ in isolation. The left- and
right-multiplication maps $L\left(  \mathbf{a}\right)  $ and $R\left(
\mathbf{a}\right)  $ are somewhat better behaved, but still not
explicit.\footnote{Exercise \ref{exe.specht.EPQ.uEUVv} allows us to rewrite a
product $u\mathbf{E}_{P,Q}v$ with $u,v\in S_{n}$ as $\mathbf{E}_{uP,v^{-1}Q}$;
but the result is usually not a basis vector of the Young symmetrizer basis,
since the tableaux $uP$ and $v^{-1}Q$ may not be standard even if $P$ and $Q$
are.} Let us see an example:

\begin{example}
\label{exa.bas.ysb.n=3}Let $n=3$. Then, the Young symmetrizer basis
$\overrightarrow{\operatorname*{ysb}}=\left(  \mathbf{E}_{U,V}\right)
_{\left(  \lambda,U,V\right)  \in\operatorname*{SBT}\left(  n\right)  }$ of
$\mathcal{A}$ consists of the six vectors%
\begin{align*}
\mathbf{E}_{123,\ 123}  &  =\nabla=1+s_{1}+s_{2}+t_{1,3}+\operatorname*{cyc}%
\nolimits_{1,2,3}+\operatorname*{cyc}\nolimits_{1,3,2};\\
\mathbf{E}_{12\backslash\backslash3,\ 12\backslash\backslash3}  &
=1+s_{1}-t_{1,3}-\operatorname*{cyc}\nolimits_{1,2,3};\\
\mathbf{E}_{12\backslash\backslash3,\ 13\backslash\backslash2}  &
=-s_{1}+s_{2}+\operatorname*{cyc}\nolimits_{1,2,3}-\operatorname*{cyc}%
\nolimits_{1,3,2};\\
\mathbf{E}_{13\backslash\backslash2,\ 12\backslash\backslash3}  &
=s_{2}-t_{1,3}-\operatorname*{cyc}\nolimits_{1,2,3}+\operatorname*{cyc}%
\nolimits_{1,3,2};\\
\mathbf{E}_{13\backslash\backslash2,\ 13\backslash\backslash2}  &
=1-s_{1}+t_{1,3}-\operatorname*{cyc}\nolimits_{1,3,2};\\
\mathbf{E}_{1\backslash\backslash2\backslash\backslash3,\ 1\backslash
\backslash2\backslash\backslash3}  &  =\nabla^{-}=1-s_{1}-s_{2}-t_{1,3}%
+\operatorname*{cyc}\nolimits_{1,2,3}+\operatorname*{cyc}\nolimits_{1,3,2}.
\end{align*}
Let us agree to order the indexing set $\operatorname*{SBT}\left(  3\right)  $
of this basis as follows:%
\begin{align*}
&  \left(  \left(  3\right)  ,\ 123,\ 123\right)  ,\ \ \ \ \ \ \ \ \ \ \left(
\left(  2,1\right)  ,\ 12\backslash\backslash3,\ 12\backslash\backslash
3\right)  ,\ \ \ \ \ \ \ \ \ \ \left(  \left(  2,1\right)  ,\ 12\backslash
\backslash3,\ 13\backslash\backslash2\right)  ,\\
&  \left(  \left(  2,1\right)  ,\ 13\backslash\backslash2,\ 12\backslash
\backslash3\right)  ,\ \ \ \ \ \ \ \ \ \ \left(  \left(  2,1\right)
,\ 13\backslash\backslash2,\ 13\backslash\backslash2\right)  ,\\
&  \left(  \left(  1,1,1\right)  ,\ 1\backslash\backslash2\backslash
\backslash3,\ 1\backslash\backslash2\backslash\backslash3\right)  .
\end{align*}
Then, the matrices representing the $\mathbf{k}$-linear maps $S$,
$T_{\operatorname*{sign}}$, $L\left(  s_{1}\right)  $, $L\left(  s_{2}\right)
$, $R\left(  s_{1}\right)  $ and $R\left(  s_{2}\right)  $ are as follows:%
\[
\left[  S\right]  _{\overrightarrow{\operatorname*{ysb}}\rightarrow
\overrightarrow{\operatorname*{ysb}}}=\left(
\begin{array}
[c]{cccccc}%
1 & 0 & 0 & 0 & 0 & 0\\
0 & 4/3 & -2/3 & 2/3 & -1/3 & 0\\
0 & 2/3 & -1/3 & 4/3 & -2/3 & 0\\
0 & -2/3 & 4/3 & -1/3 & 2/3 & 0\\
0 & -1/3 & 2/3 & -2/3 & 4/3 & 0\\
0 & 0 & 0 & 0 & 0 & 1
\end{array}
\right)  ;
\]%
\[
\left[  T_{\operatorname*{sign}}\right]  _{\overrightarrow{\operatorname*{ysb}%
}\rightarrow\overrightarrow{\operatorname*{ysb}}}=\left(
\begin{array}
[c]{cccccc}%
0 & 0 & 0 & 0 & 0 & 1\\
0 & -1/3 & 2/3 & -2/3 & 4/3 & 0\\
0 & -2/3 & 1/3 & -4/3 & 2/3 & 0\\
0 & 2/3 & -4/3 & 1/3 & -2/3 & 0\\
0 & 4/3 & -2/3 & 2/3 & -1/3 & 0\\
1 & 0 & 0 & 0 & 0 & 0
\end{array}
\right)  ;
\]%
\[
L_{\overrightarrow{\operatorname*{ysb}}}\left(  s_{1}\right)  =\left[
L\left(  s_{1}\right)  \right]  _{\overrightarrow{\operatorname*{ysb}%
}\rightarrow\overrightarrow{\operatorname*{ysb}}}=\left(
\begin{array}
[c]{cccccc}%
1 & 0 & 0 & 0 & 0 & 0\\
0 & 1 & 0 & 0 & 0 & 0\\
0 & 0 & 1 & 0 & 0 & 0\\
0 & -1 & 0 & -1 & 0 & 0\\
0 & 0 & -1 & 0 & -1 & 0\\
0 & 0 & 0 & 0 & 0 & -1
\end{array}
\right)  ;
\]%
\[
L_{\overrightarrow{\operatorname*{ysb}}}\left(  s_{2}\right)  =\left[
L\left(  s_{2}\right)  \right]  _{\overrightarrow{\operatorname*{ysb}%
}\rightarrow\overrightarrow{\operatorname*{ysb}}}=\left(
\begin{array}
[c]{cccccc}%
1 & 0 & 0 & 0 & 0 & 0\\
0 & 0 & 0 & 1 & 0 & 0\\
0 & 0 & 0 & 0 & 1 & 0\\
0 & 1 & 0 & 0 & 0 & 0\\
0 & 0 & 1 & 0 & 0 & 0\\
0 & 0 & 0 & 0 & 0 & -1
\end{array}
\right)  ;
\]%
\[
R_{\overrightarrow{\operatorname*{ysb}}}\left(  s_{1}\right)  =\left[
R\left(  s_{1}\right)  \right]  _{\overrightarrow{\operatorname*{ysb}%
}\rightarrow\overrightarrow{\operatorname*{ysb}}}=\left(
\begin{array}
[c]{cccccc}%
1 & 0 & 0 & 0 & 0 & 0\\
0 & 1 & -1 & 0 & 0 & 0\\
0 & 0 & -1 & 0 & 0 & 0\\
0 & 0 & 0 & 1 & -1 & 0\\
0 & 0 & 0 & 0 & -1 & 0\\
0 & 0 & 0 & 0 & 0 & -1
\end{array}
\right)  ;
\]%
\[
R_{\overrightarrow{\operatorname*{ysb}}}\left(  s_{2}\right)  =\left[
R\left(  s_{2}\right)  \right]  _{\overrightarrow{\operatorname*{ysb}%
}\rightarrow\overrightarrow{\operatorname*{ysb}}}=\left(
\begin{array}
[c]{cccccc}%
1 & 0 & 0 & 0 & 0 & 0\\
0 & 0 & 1 & 0 & 0 & 0\\
0 & 1 & 0 & 0 & 0 & 0\\
0 & 0 & 0 & 0 & 1 & 0\\
0 & 0 & 0 & 1 & 0 & 0\\
0 & 0 & 0 & 0 & 0 & -1
\end{array}
\right)  .
\]

\end{example}

The simplicity of the last four matrices is somewhat misleading: For higher
values of $n$, the matrix $L_{\overrightarrow{\operatorname*{ysb}}}\left(
s_{2}\right)  $ (for example) can contain entries other than $0$, $1$ and
$-1$. (For instance, it contains a $2$ when $n=8$.) The same applies to
$R_{\overrightarrow{\operatorname*{ysb}}}\left(  s_{2}\right)  $. The same
applies to $L_{\overrightarrow{\operatorname*{ysb}}}\left(  s_{1}\right)  $,
although this requires a higher value of $n$\ \ \ \ \footnote{Indeed, if
$n=11$ and $\lambda=\left(  4,3,2,2\right)  $ and
$T=\ytableaushort{1236,457,8{10},9{11}}$\ , then the standard polytabloid
$\mathbf{e}_{T}$ in the Specht module $\mathcal{S}^{\lambda}$ has the property
that the expansion of $s_{1}\mathbf{e}_{T}$ as a linear combination of
standard polytabloids contains a coefficient of $-2$. (This has been checked
using both SageMath and Macaulay2.) By Proposition \ref{prop.bas.ysb.str}
further below, this entails that the matrix
$L_{\overrightarrow{\operatorname*{ysb}}}\left(  s_{1}\right)  $ has an entry
of $-2$.}.

However, some other patterns that can be spotted on the above matrices are
generalizable. For instance, the matrices
$R_{\overrightarrow{\operatorname*{ysb}}}\left(  s_{1}\right)  $ and
$R_{\overrightarrow{\operatorname*{ysb}}}\left(  s_{2}\right)  $ (for $n=3$)
are block-diagonal, where the diagonal blocks have sizes $1,2,2,1$. The same
holds for the matrices $L_{\overrightarrow{\operatorname*{ysb}}}\left(
s_{1}\right)  $ and $L_{\overrightarrow{\operatorname*{ysb}}}\left(
s_{2}\right)  $ after a permutation of the rows and the columns. These
observations are not a fluke; they are manifestations of Lemma
\ref{lem.specht.EPQ.straighten}. Namely:

\begin{itemize}
\item The equality (\ref{eq.lem.specht.EPQ.straighten.AET}) shows that for
each partition $\lambda$ of $n$ and each standard tableau $Q\in
\operatorname*{SYT}\left(  \lambda\right)  $, the subfamily $\left(
\mathbf{E}_{P,Q}\right)  _{P\in\operatorname*{SYT}\left(  \lambda\right)  }$
of the Young symmetrizer basis spans the left ideal $\mathcal{A}\mathbf{E}%
_{Q}$. In particular, this means that its span is preserved by left
multiplication with each $s_{i}$. Thus, the matrices
$L_{\overrightarrow{\operatorname*{ysb}}}\left(  s_{i}\right)  $ for all
$i\in\left[  n-1\right]  $ (and, more generally, the matrices
$L_{\overrightarrow{\operatorname*{ysb}}}\left(  \mathbf{a}\right)  $ for all
$\mathbf{a}\in\mathcal{A}$) become block-diagonal if the rows and the columns
are ordered in such a way that the standard $n$-bitableaux $\left(
\lambda,P,Q\right)  $ with equal $Q$'s appear as a contiguous block.

\item The equality (\ref{eq.lem.specht.EPQ.straighten.ETA}) shows that for
each partition $\lambda$ of $n$ and each standard tableau $Q\in
\operatorname*{SYT}\left(  \lambda\right)  $, the subfamily $\left(
\mathbf{E}_{Q,P}\right)  _{P\in\operatorname*{SYT}\left(  \lambda\right)  }$
of the Young symmetrizer basis spans the right ideal $\mathbf{E}%
_{Q}\mathcal{A}$. In particular, this means that its span is preserved by
right multiplication with each $s_{i}$. Thus, the matrices
$R_{\overrightarrow{\operatorname*{ysb}}}\left(  s_{i}\right)  $ for all
$i\in\left[  n-1\right]  $ (and, more generally, the matrices
$R_{\overrightarrow{\operatorname*{ysb}}}\left(  \mathbf{a}\right)  $ for all
$\mathbf{a}\in\mathcal{A}$) become block-diagonal if the rows and the columns
are ordered in such a way that the standard $n$-bitableaux $\left(
\lambda,Q,P\right)  $ with equal $Q$'s appear as a contiguous block.

\item Finally, Theorem \ref{thm.specht.AElam.nat-basis} shows that for each
partition $\lambda$ of $n$, the subfamily $\left(  \mathbf{E}_{U,V}\right)
_{U,V\in\operatorname*{SYT}\left(  \lambda\right)  }$ of the Young symmetrizer
basis spans the $\mathbf{k}$-module $\mathcal{A}\mathbf{E}_{\lambda}$, which
is a two-sided ideal of $\mathcal{A}$ (an easy consequence of $\mathbf{E}%
_{\lambda}\in Z\left(  \mathcal{A}\right)  $).
\end{itemize}

Even better, the entries of the matrices
$L_{\overrightarrow{\operatorname*{ysb}}}\left(  \mathbf{a}\right)  $ can be
described in terms of the straightening coefficients (the coefficients of the
$\mathbb{Z}$-linear combination in Proposition \ref{prop.spechtmod.straighten}):

\begin{proposition}
\label{prop.bas.ysb.str}Let $\left(  \lambda,U,V\right)  \in
\operatorname*{SBT}\left(  n\right)  $. Let $\mathbf{a}\in\mathcal{A}$. Write
the element $\mathbf{ae}_{U}$ in the Specht module $\mathcal{S}^{\lambda}$ as
a $\mathbf{k}$-linear combination%
\begin{equation}
\mathbf{ae}_{U}=\sum_{Q\in\operatorname*{SYT}\left(  \lambda\right)  }%
\alpha_{U,Q}\mathbf{e}_{Q} \label{eq.prop.bas.ysb.str.ass}%
\end{equation}
of standard polytabloids $\mathbf{e}_{Q}$. (This can be done, by Proposition
\ref{prop.spechtmod.straighten}.) Then,%
\begin{equation}
\mathbf{aE}_{U,V}=\sum_{Q\in\operatorname*{SYT}\left(  \lambda\right)  }%
\alpha_{U,Q}\mathbf{E}_{Q,V}. \label{eq.prop.bas.ysb.str.clm}%
\end{equation}
This holds for arbitrary commutative rings $\mathbf{k}$, even if $\left(
\mathbf{E}_{U,V}\right)  _{\left(  \lambda,U,V\right)  \in\operatorname*{SBT}%
\left(  n\right)  }$ is not a basis of $\mathcal{A}$.
\end{proposition}

\begin{proof}
Rename the $n$-tableau $V$ as $T$. Set $D=Y\left(  \lambda\right)  $, so that
$\operatorname*{SYT}\left(  \lambda\right)  =\operatorname*{SYT}\left(
D\right)  $; hence, in particular, $U$ and $V$ belong to $\operatorname*{SYT}%
\left(  D\right)  $. Now, recall the proof of Corollary
\ref{prop.specht.EPQ.basis-AET}, in which we constructed a left $\mathcal{A}%
$-module isomorphism $\overline{\alpha}:\mathcal{A}\mathbf{E}_{T}%
\rightarrow\mathcal{S}^{D}$. The inverse $\overline{\alpha}^{-1}$ of this
isomorphism is thus a left $\mathcal{A}$-module isomorphism $\overline{\alpha
}^{-1}:\mathcal{S}^{D}\rightarrow\mathcal{A}\mathbf{E}_{T}$. Applying this
isomorphism $\overline{\alpha}^{-1}$ to the equality
(\ref{eq.prop.bas.ysb.str.ass}), we obtain%
\begin{equation}
\mathbf{a}\cdot\overline{\alpha}^{-1}\left(  \mathbf{e}_{U}\right)
=\sum_{Q\in\operatorname*{SYT}\left(  \lambda\right)  }\alpha_{U,Q}%
\overline{\alpha}^{-1}\left(  \mathbf{e}_{Q}\right)
\label{pf.prop.bas.ysb.str.1}%
\end{equation}
(since $\overline{\alpha}^{-1}$ is left $\mathcal{A}$-linear). But the
equality (\ref{pf.prop.specht.EPQ.basis-AET.8pre}) (which we have shown in the
above-mentioned proof) says that $\overline{\alpha}^{-1}\left(  \mathbf{e}%
_{P}\right)  =\mathbf{E}_{P,T}$ for each $P\in\operatorname*{SYT}\left(
D\right)  $. In other words, $\overline{\alpha}^{-1}\left(  \mathbf{e}%
_{P}\right)  =\mathbf{E}_{P,T}$ for each $P\in\operatorname*{SYT}\left(
\lambda\right)  $ (since $\operatorname*{SYT}\left(  D\right)
=\operatorname*{SYT}\left(  \lambda\right)  $). Thus, we can rewrite the
equality (\ref{pf.prop.bas.ysb.str.1}) as%
\[
\mathbf{aE}_{U,T}=\sum_{Q\in\operatorname*{SYT}\left(  \lambda\right)  }%
\alpha_{U,Q}\mathbf{E}_{Q,T}.
\]
This proves Proposition \ref{prop.bas.ysb.str} (since $T=V$).
\end{proof}

The equality (\ref{eq.prop.bas.ysb.str.clm}) can be rewritten as%
\begin{align}
\mathbf{aE}_{U,V}  &  =\sum_{\left(  \mu,Q,P\right)  \in\operatorname*{SBT}%
\left(  n\right)  }\alpha_{\left(  \lambda,U,V\right)  ,\left(  \mu
,Q,P\right)  }\mathbf{E}_{Q,P},\label{eq.prop.bas.ysb.str.clm2}\\
&  \ \ \ \ \ \ \ \ \ \ \text{where }\alpha_{\left(  \lambda,U,V\right)
,\left(  \mu,Q,P\right)  }=%
\begin{cases}
\alpha_{U,Q}, & \text{if }\lambda=\mu\text{ and }V=P;\\
0, & \text{else}%
\end{cases}
\nonumber
\end{align}
(since all addends on the right hand side of (\ref{eq.prop.bas.ysb.str.clm2})
are $0$ except for those that satisfy $\lambda=\mu$ and $V=P$; but the latter
are just the addends on the right hand side of (\ref{eq.prop.bas.ysb.str.clm})).

The coefficients $\alpha_{\left(  \lambda,U,V\right)  ,\left(  \mu,Q,P\right)
}$ on the right hand side of (\ref{eq.prop.bas.ysb.str.clm2}) are thus the
entries of the matrix $L_{\overrightarrow{\operatorname*{ysb}}}\left(
\mathbf{a}\right)  $. This shows once again that the matrix
$L_{\overrightarrow{\operatorname*{ysb}}}\left(  \mathbf{a}\right)  $ is
block-diagonal after permuting its rows and columns, and moreover explains why
its entries are integers when the coefficients of $\mathbf{a}$ are integers
(since the $\alpha_{U,Q}$ are integers in this case). Furthermore, while the
diagonal blocks of $L_{\overrightarrow{\operatorname*{ysb}}}\left(
\mathbf{a}\right)  $ are indexed by the pairs $\left(  \lambda,V\right)  $ of
a partition $\lambda$ and a standard tableau $V$, the actual entries of these
blocks depend only on $\lambda$ (not on $V$), since they are the scalars
$\alpha_{U,Q}$.

The matrices $R_{\overrightarrow{\operatorname*{ysb}}}\left(  \mathbf{a}%
\right)  $ can be explained similarly using the \textquotedblleft right Specht
modules\textquotedblright\ (the right $\mathcal{A}$-module analogues of the
Specht modules). Alternatively, they can be obtained from the
$L_{\overrightarrow{\operatorname*{ysb}}}\left(  T_{\operatorname*{sign}%
}\left(  S\left(  \mathbf{a}\right)  \right)  \right)  $ by conjugating with
the $\mathbf{k}$-algebra anti-automorphism $T_{\operatorname*{sign}}\circ S$
(which is represented by a rather simple matrix in the basis
$\overrightarrow{\operatorname*{ysb}}$, thanks to Proposition
\ref{prop.specht.EPQ.r}\footnote{specifically, a permutation matrix with some
of the $1$'s replaced by $-1$'s}). Indeed, it is easy to see that each
$\mathbf{a}\in\mathcal{A}$ satisfies the equality of endomorphisms
\begin{align*}
R\left(  \mathbf{a}\right)   &  =S^{-1}\circ L\left(  S\left(  \mathbf{a}%
\right)  \right)  \circ S\\
&  =\left(  T_{\operatorname*{sign}}\circ S\right)  ^{-1}\circ L\left(
T_{\operatorname*{sign}}\left(  S\left(  \mathbf{a}\right)  \right)  \right)
\circ\left(  T_{\operatorname*{sign}}\circ S\right)
\end{align*}
and thus the equality of matrices%
\[
R_{\overrightarrow{\operatorname*{ysb}}}\left(  \mathbf{a}\right)
=\Omega^{-1}\cdot L_{\overrightarrow{\operatorname*{ysb}}}\left(
T_{\operatorname*{sign}}\left(  S\left(  \mathbf{a}\right)  \right)  \right)
\cdot\Omega,\ \ \ \ \ \ \ \ \ \ \text{where }\Omega:=\left[
T_{\operatorname*{sign}}\circ S\right]  _{\overrightarrow{\operatorname*{ysb}%
}\rightarrow\overrightarrow{\operatorname*{ysb}}}.
\]
When $\mathbf{a}=s_{i}$ for some $i\in\left[  n-1\right]  $, then $S\left(
\mathbf{a}\right)  =S\left(  s_{i}\right)  =s_{i}^{-1}=s_{i}$ and thus
$T_{\operatorname*{sign}}\left(  S\left(  \mathbf{a}\right)  \right)
=T_{\operatorname*{sign}}\left(  s_{i}\right)  =-s_{i}$, so this equality
becomes even simpler:%
\[
R_{\overrightarrow{\operatorname*{ysb}}}\left(  s_{i}\right)  =-\Omega
^{-1}\cdot L_{\overrightarrow{\operatorname*{ysb}}}\left(  s_{i}\right)
\cdot\Omega.
\]
Since $\Omega$ is a permutation matrix \textquotedblleft up to
sign\textquotedblright, this ensures that the entries of
$R_{\overrightarrow{\operatorname*{ysb}}}\left(  s_{i}\right)  $ agree with
those of $L_{\overrightarrow{\operatorname*{ysb}}}\left(  s_{i}\right)  $ up
to sign and position. It is worth noting that $\Omega^{-1}=\Omega$, since
$T_{\operatorname*{sign}}\circ S$ is an involution (by Exercise
\ref{exe.S.STsign} \textbf{(b)}).

\subsection{\label{sec.bas.FwE}The $\mathbf{F}_{U}w_{U,V}\mathbf{E}_{V}$ basis
of $\mathbf{k}\left[  S_{n}\right]  $, revisited}

For the present Section \ref{sec.bas.FwE}, we again assume that $n!$ is
invertible in $\mathbf{k}$. Then, Theorem \ref{thm.specht.FwE.basis} shows
that $\left(  \mathbf{F}_{U}w_{U,V}\mathbf{E}_{V}\right)  _{\lambda\text{ is a
partition of }n\text{, and }U,V\in\operatorname*{SYT}\left(  \lambda\right)
}$ is a basis of the $\mathbf{k}$-module $\mathcal{A}$ (where the notations we
use are defined as in Convention \ref{conv.specht.FwE.conv}). Let us call this
basis the \emph{FwE basis} of $\mathcal{A}$. In view of Definition
\ref{def.bas.not} \textbf{(e)}, we can rewrite the condition \textquotedblleft%
$\lambda$ is a partition of $n$, and $U,V\in\operatorname*{SYT}\left(
\lambda\right)  $\textquotedblright\ as \textquotedblleft$\left(
\lambda,U,V\right)  $ is a standard $n$-bitableau\textquotedblright, or --
even shorter -- as \textquotedblleft$\left(  \lambda,U,V\right)
\in\operatorname*{SBT}\left(  n\right)  $\textquotedblright. Thus, the FwE
basis can be written as
\[
\left(  \mathbf{F}_{U}w_{U,V}\mathbf{E}_{V}\right)  _{\left(  \lambda
,U,V\right)  \in\operatorname*{SBT}\left(  n\right)  }.
\]
Let us denote this basis by $\overrightarrow{\operatorname*{fwe}}$.

Proposition \ref{prop.specht.FwE.S} shows that the antipode $S$ of
$\mathcal{A}$ permutes this basis, sending its $\left(  \lambda,U,V\right)
$-entry to its $\left(  \lambda,V,U\right)  $-entry. The actions of the
operators $T_{\operatorname*{sign}}$, $L\left(  \mathbf{a}\right)  $ and
$R\left(  \mathbf{a}\right)  $ (for $\mathbf{a}\in\mathbf{k}\left[
S_{n}\right]  $) on this basis are far less simple.

\begin{example}
Let $n=3$. Let us agree to order the indexing set $\operatorname*{SBT}\left(
3\right)  $ of the FwE basis as follows:%
\begin{align*}
&  \left(  \left(  3\right)  ,\ 123,\ 123\right)  ,\ \ \ \ \ \ \ \ \ \ \left(
\left(  2,1\right)  ,\ 12\backslash\backslash3,\ 12\backslash\backslash
3\right)  ,\ \ \ \ \ \ \ \ \ \ \left(  \left(  2,1\right)  ,\ 12\backslash
\backslash3,\ 13\backslash\backslash2\right)  ,\\
&  \left(  \left(  2,1\right)  ,\ 13\backslash\backslash2,\ 12\backslash
\backslash3\right)  ,\ \ \ \ \ \ \ \ \ \ \left(  \left(  2,1\right)
,\ 13\backslash\backslash2,\ 13\backslash\backslash2\right)  ,\\
&  \left(  \left(  1,1,1\right)  ,\ 1\backslash\backslash2\backslash
\backslash3,\ 1\backslash\backslash2\backslash\backslash3\right)  .
\end{align*}
We have already listed the six entries of this basis in Example
\ref{exa.specht.FwE.n=3}. The matrices representing the $\mathbf{k}$-linear
maps $S$, $T_{\operatorname*{sign}}$, $L\left(  s_{1}\right)  $, $L\left(
s_{2}\right)  $, $R\left(  s_{1}\right)  $ and $R\left(  s_{2}\right)  $ are
as follows:%
\[
\left[  S\right]  _{\overrightarrow{\operatorname*{fwe}}\rightarrow
\overrightarrow{\operatorname*{fwe}}}=\left(
\begin{array}
[c]{cccccc}%
1 & 0 & 0 & 0 & 0 & 0\\
0 & 1 & 0 & 0 & 0 & 0\\
0 & 0 & 0 & 1 & 0 & 0\\
0 & 0 & 1 & 0 & 0 & 0\\
0 & 0 & 0 & 0 & 1 & 0\\
0 & 0 & 0 & 0 & 0 & 1
\end{array}
\right)  ;
\]%
\[
\left[  T_{\operatorname*{sign}}\right]  _{\overrightarrow{\operatorname*{fwe}%
}\rightarrow\overrightarrow{\operatorname*{fwe}}}=\left(
\begin{array}
[c]{cccccc}%
0 & 0 & 0 & 0 & 0 & 1\\
0 & 1/3 & -2/3 & -2/3 & 4/3 & 0\\
0 & 2/3 & -1/3 & -4/3 & 2/3 & 0\\
0 & 2/3 & -4/3 & -1/3 & 2/3 & 0\\
0 & 4/3 & -2/3 & -2/3 & 1/3 & 0\\
1 & 0 & 0 & 0 & 0 & 0
\end{array}
\right)  ;
\]%
\[
L_{\overrightarrow{\operatorname*{fwe}}}\left(  s_{1}\right)  =\left[
L\left(  s_{1}\right)  \right]  _{\overrightarrow{\operatorname*{fwe}%
}\rightarrow\overrightarrow{\operatorname*{fwe}}}=\left(
\begin{array}
[c]{cccccc}%
1 & 0 & 0 & 0 & 0 & 0\\
0 & 1 & 0 & -1 & 0 & 0\\
0 & 0 & 1 & 0 & -1 & 0\\
0 & 0 & 0 & -1 & 0 & 0\\
0 & 0 & 0 & 0 & -1 & 0\\
0 & 0 & 0 & 0 & 0 & -1
\end{array}
\right)  ;
\]%
\[
L_{\overrightarrow{\operatorname*{fwe}}}\left(  s_{2}\right)  =\left[
L\left(  s_{2}\right)  \right]  _{\overrightarrow{\operatorname*{fwe}%
}\rightarrow\overrightarrow{\operatorname*{fwe}}}=\left(
\begin{array}
[c]{cccccc}%
1 & 0 & 0 & 0 & 0 & 0\\
0 & 0 & 0 & 1 & 0 & 0\\
0 & 0 & 0 & 0 & 1 & 0\\
0 & 1 & 0 & 0 & 0 & 0\\
0 & 0 & 1 & 0 & 0 & 0\\
0 & 0 & 0 & 0 & 0 & -1
\end{array}
\right)  ;
\]%
\[
R_{\overrightarrow{\operatorname*{fwe}}}\left(  s_{1}\right)  =\left[
R\left(  s_{1}\right)  \right]  _{\overrightarrow{\operatorname*{fwe}%
}\rightarrow\overrightarrow{\operatorname*{fwe}}}=\left(
\begin{array}
[c]{cccccc}%
1 & 0 & 0 & 0 & 0 & 0\\
0 & 1 & -1 & 0 & 0 & 0\\
0 & 0 & -1 & 0 & 0 & 0\\
0 & 0 & 0 & 1 & -1 & 0\\
0 & 0 & 0 & 0 & -1 & 0\\
0 & 0 & 0 & 0 & 0 & -1
\end{array}
\right)  ;
\]%
\[
R_{\overrightarrow{\operatorname*{fwe}}}\left(  s_{2}\right)  =\left[
R\left(  s_{2}\right)  \right]  _{\overrightarrow{\operatorname*{fwe}%
}\rightarrow\overrightarrow{\operatorname*{fwe}}}=\left(
\begin{array}
[c]{cccccc}%
1 & 0 & 0 & 0 & 0 & 0\\
0 & 0 & 1 & 0 & 0 & 0\\
0 & 1 & 0 & 0 & 0 & 0\\
0 & 0 & 0 & 0 & 1 & 0\\
0 & 0 & 0 & 1 & 0 & 0\\
0 & 0 & 0 & 0 & 0 & -1
\end{array}
\right)  .
\]
Note that the last three of these six matrices agree with those for the Young
symmetrizer basis (which we saw in Example \ref{exa.bas.ysb.n=3})! The
equality $L_{\overrightarrow{\operatorname*{fwe}}}\left(  s_{2}\right)
=L_{\overrightarrow{\operatorname*{ysb}}}\left(  s_{2}\right)  $ is a fluke
(it fails for $n=5$), but the equalities
$R_{\overrightarrow{\operatorname*{fwe}}}\left(  s_{i}\right)
=R_{\overrightarrow{\operatorname*{ysb}}}\left(  s_{i}\right)  $ do indeed
hold for all $n\in\mathbb{N}$ and $i\in\left[  n-1\right]  $. This follows
from the following exercise:
\end{example}

\begin{exercise}
\label{exa.bas.fwe.ysb}\fbox{3} Show that
$R_{\overrightarrow{\operatorname*{fwe}}}\left(  \mathbf{a}\right)
=R_{\overrightarrow{\operatorname*{ysb}}}\left(  \mathbf{a}\right)  $ for all
$n\in\mathbb{N}$ and $\mathbf{a}\in\mathcal{A}$. \medskip

[\textbf{Hint:} For each partition $\lambda$ of $n$ and each standard tableau
$U\in\operatorname*{SYT}\left(  \lambda\right)  $, construct a right
$\mathcal{A}$-module morphism%
\begin{align*}
\Xi_{U}:\mathbf{E}_{U}\mathcal{A}  &  \rightarrow\mathbf{F}_{U}\mathcal{A},\\
\mathbf{x}  &  \mapsto\mathbf{F}_{U}\mathbf{x}.
\end{align*}
Show that this morphism sends the basis vectors $\mathbf{E}_{U,V}$ of
$\overrightarrow{\operatorname*{ysb}}$ to the respective basis vectors
$\mathbf{F}_{U}w_{U,V}\mathbf{E}_{V}$ of $\overrightarrow{\operatorname*{fwe}%
}$. Now apply it to an expansion of the form $\mathbf{E}_{U,V}\mathbf{a}%
=\sum_{Q\in\operatorname*{SYT}\left(  \lambda\right)  }\beta_{U,Q}%
\mathbf{E}_{U,Q}$, which exists thanks to
(\ref{eq.lem.specht.EPQ.straighten.ETA}).]
\end{exercise}

Again, the matrices $L_{\overrightarrow{\operatorname*{fwe}}}\left(
\mathbf{a}\right)  $ and $R_{\overrightarrow{\operatorname*{fwe}}}\left(
\mathbf{a}\right)  $ for all $\mathbf{a}\in\mathcal{A}$ are block-diagonal,
after appropriate permutations of rows and columns. We leave it to the reader
to explain this. The two types of matrices are connected by the equality%
\[
R_{\overrightarrow{\operatorname*{fwe}}}\left(  \mathbf{a}\right)
=\Omega_{\overrightarrow{\operatorname*{fwe}}}^{-1}\cdot
L_{\overrightarrow{\operatorname*{fwe}}}\left(  S\left(  \mathbf{a}\right)
\right)  \cdot\Omega_{\overrightarrow{\operatorname*{fwe}}}%
,\ \ \ \ \ \ \ \ \ \ \text{where }\Omega_{\overrightarrow{\operatorname*{fwe}%
}}:=\left[  S\right]  _{\overrightarrow{\operatorname*{fwe}}\rightarrow
\overrightarrow{\operatorname*{fwe}}}%
\]
(why?), which is particularly nice because $\Omega
_{\overrightarrow{\operatorname*{fwe}}}$ is just a permutation matrix (since
the antipode $S$ permutes the basis $\overrightarrow{\operatorname*{fwe}}$)
and satisfies $\Omega_{\overrightarrow{\operatorname*{fwe}}}^{-1}%
=\Omega_{\overrightarrow{\operatorname*{fwe}}}$ (since $S$ is an involution).

\subsection{\label{sec.bas.AWS}The Artin--Wedderburn--Specht bases}

One might believe that the block-diagonality properties make the bases
$\overrightarrow{\operatorname*{ysb}}$ and
$\overrightarrow{\operatorname*{fwe}}$ special; but they don't. There is a
general machine for constructing bases of $\mathcal{A}$ with such properties,
assuming that $n!$ is invertible in $\mathbf{k}$:

For the present Section \ref{sec.bas.AWS}, we again assume that $n!$ is
invertible in $\mathbf{k}$. The Artin--Wedderburn theorem (Theorem
\ref{thm.specht.AW}) yields that the map%
\begin{align*}
\rho:\mathcal{A}  &  \rightarrow\prod_{\lambda\text{ is a partition of }%
n}\operatorname*{End}\nolimits_{\mathbf{k}}\left(  \mathcal{S}^{\lambda
}\right)  ,\\
\mathbf{a}  &  \mapsto\left(  \rho_{\lambda}\left(  \mathbf{a}\right)
\right)  _{\lambda\text{ is a partition of }n}%
\end{align*}
is a $\mathbf{k}$-algebra isomorphism. But each Specht module $\mathcal{S}%
^{\lambda}$ is a free $\mathbf{k}$-module of rank $f^{\lambda}$ (by Lemma
\ref{lem.specht.Slam-flam}); hence, its endomorphism algebra
$\operatorname*{End}\nolimits_{\mathbf{k}}\left(  \mathcal{S}^{\lambda
}\right)  $ is a free $\mathbf{k}$-module of rank $\left(  f^{\lambda}\right)
^{2}$ (by Lemma \ref{lem.linalg.dim-End}). Thus, we can construct a basis of
$\mathcal{A}$ as follows:

\begin{enumerate}
\item[Step 1:] Pick any basis of $\operatorname*{End}\nolimits_{\mathbf{k}%
}\left(  \mathcal{S}^{\lambda}\right)  $ for each partition $\lambda$. The
simplest way to do so is to pick a basis $\left(  \mathbf{b}_{T}\right)
_{T\in\operatorname*{SYT}\left(  \lambda\right)  }$ of $\mathcal{S}^{\lambda}$
(indexed by standard tableaux because $\mathcal{S}^{\lambda}$ is free of rank
$f^{\lambda}=\left\vert \operatorname*{SYT}\left(  \lambda\right)  \right\vert
$; theoretically, other indexing sets are possible, but I have never seen a
such), and then construct the \textquotedblleft elementary matrix
basis\textquotedblright\ $\left(  \mathbf{M}_{U,V}\right)  _{\left(
U,V\right)  \in\operatorname*{SYT}\left(  \lambda\right)  \times
\operatorname*{SYT}\left(  \lambda\right)  }$ of $\operatorname*{End}%
\nolimits_{\mathbf{k}}\left(  \mathcal{S}^{\lambda}\right)  $ induced by it
(explicitly: $\mathbf{M}_{U,V}$ is the $\mathbf{k}$-module endomorphism of
$\mathcal{S}^{\lambda}$ that sends $\mathbf{b}_{V}$ to $\mathbf{b}_{U}$ and
sends all $\mathbf{b}_{Q}$ with $Q\neq V$ to $0$).

\item[Step 2:] Combine all these bases $\left(  \mathbf{M}_{U,V}\right)
_{\left(  U,V\right)  \in\operatorname*{SYT}\left(  \lambda\right)
\times\operatorname*{SYT}\left(  \lambda\right)  }$ (for all partitions
$\lambda$ of $n$) into a basis $\left(  \mathbf{M}_{\lambda,U,V}\right)
_{\left(  \lambda,U,V\right)  \in\operatorname*{SBT}\left(  n\right)  }$ of
$\prod_{\lambda\text{ is a partition of }n}\operatorname*{End}%
\nolimits_{\mathbf{k}}\left(  \mathcal{S}^{\lambda}\right)  $.

\item[Step 3:] Then, apply $\rho^{-1}$ to this basis -- that is, set
$\mathbf{M}_{\lambda,U,V}^{\prime}:=\rho^{-1}\left(  \mathbf{M}_{\lambda
,U,V}\right)  $ for each $\left(  \lambda,U,V\right)  \in\operatorname*{SBT}%
\left(  n\right)  $. This yields a basis $\left(  \mathbf{M}_{\lambda
,U,V}^{\prime}\right)  _{\left(  \lambda,U,V\right)  \in\operatorname*{SBT}%
\left(  n\right)  }$ of $\mathcal{A}$.
\end{enumerate}

Bases of $\mathcal{A}$ obtained in such a way have some properties in common.
In particular, they all rely on the invertibility of $n!$ (since $\rho$ is
usually not an isomorphism otherwise), and they are all indexed by the
standard $n$-bitableaux $\left(  \lambda,U,V\right)  $. Their subfamilies
often span ideals: For each given partition $\lambda$ of $n$, the span of
$\left\{  \mathbf{M}_{\lambda,U,V}^{\prime}\ \mid\ \left(  U,V\right)
\in\operatorname*{SYT}\left(  \lambda\right)  \times\operatorname*{SYT}\left(
\lambda\right)  \right\}  $ is the two-sided ideal $\mathcal{A}\mathbf{E}%
_{\lambda}$ of $\mathcal{A}$. Moreover, if the underlying bases $\left(
\mathbf{M}_{U,V}\right)  _{\left(  U,V\right)  \in\operatorname*{SYT}\left(
\lambda\right)  \times\operatorname*{SYT}\left(  \lambda\right)  }$ of
$\operatorname*{End}\nolimits_{\mathbf{k}}\left(  \mathcal{S}^{\lambda
}\right)  $ are \textquotedblleft elementary matrix bases\textquotedblright%
\ (see Step 1 above), then the span of $\left\{  \mathbf{M}_{\lambda
,U,V}^{\prime}\ \mid\ V\in\operatorname*{SYT}\left(  \lambda\right)  \right\}
$ for fixed $\lambda$ and $U$ is a right ideal of $\mathcal{A}$, while the
span of $\left\{  \mathbf{M}_{\lambda,U,V}^{\prime}\ \mid\ U\in
\operatorname*{SYT}\left(  \lambda\right)  \right\}  $ for fixed $\lambda$ and
$V$ is a left ideal of $\mathcal{A}$. All of this is rather easy to prove. I
refer to such bases as \emph{Artin--Wedderburn bases}.

Let us apply the above construction, picking the easiest option in Step 1: As
our basis $\left(  \mathbf{b}_{T}\right)  _{T\in\operatorname*{SYT}\left(
\lambda\right)  }$ of $\mathcal{S}^{\lambda}$, we choose the basis $\left(
\mathbf{e}_{T}\right)  _{T\in\operatorname*{SYT}\left(  \lambda\right)  }$
consisting of the standard polytabloids (see Theorem \ref{thm.spechtmod.basis}%
). We then construct the \textquotedblleft elementary matrix
basis\textquotedblright\ $\left(  \mathbf{M}_{U,V}\right)  _{\left(
U,V\right)  \in\operatorname*{SYT}\left(  \lambda\right)  \times
\operatorname*{SYT}\left(  \lambda\right)  }$, and proceed as explained above.
Let us denote the resulting basis $\left(  \mathbf{M}_{\lambda,U,V}^{\prime
}\right)  _{\left(  \lambda,U,V\right)  \in\operatorname*{SBT}\left(
n\right)  }$ of $\mathcal{A}$ by $\left(  \mathbf{W}_{U,V}\right)  _{\left(
\lambda,U,V\right)  \in\operatorname*{SBT}\left(  n\right)  }$ (we omit the
mention of $\lambda$ since $U$ uniquely determines $\lambda$) and by
$\overrightarrow{\operatorname*{aws}}$ (short for \textquotedblleft
Artin--Wedderburn--Specht\textquotedblright\ or \textquotedblleft
Artin--Wedderburn--standard\textquotedblright). Explicitly, $\mathbf{W}_{U,V}$
can be characterized as the unique element of $\mathcal{A}$ that satisfies
$\mathbf{W}_{U,V}\mathbf{e}_{V}=\mathbf{e}_{U}$ in the Specht module
$\mathcal{S}^{\lambda}$ (where $V\in\operatorname*{SYT}\left(  \lambda\right)
$) and satisfies $\mathbf{W}_{U,V}\mathbf{e}_{W}=0$ for all standard
$n$-tableaux $W\neq V$ of partition shapes.

\begin{example}
For $n=3$, this basis $\overrightarrow{\operatorname*{aws}}=\left(
\mathbf{W}_{U,V}\right)  _{\left(  \lambda,U,V\right)  \in\operatorname*{SBT}%
\left(  n\right)  }$ of $\mathcal{A}$ consists of the vectors
\begin{align*}
\mathbf{W}_{123,\ 123}  &  =\dfrac{1}{6}\nabla=\dfrac{1}{6}+\dfrac{1}{6}%
s_{1}+\dfrac{1}{6}s_{2}+\dfrac{1}{6}t_{1,3}+\dfrac{1}{6}\operatorname*{cyc}%
\nolimits_{1,2,3}+\,\dfrac{1}{6}\operatorname*{cyc}\nolimits_{1,3,2};\\
\mathbf{W}_{12\backslash\backslash3,\ 12\backslash\backslash3}  &  =\dfrac
{1}{3}+\dfrac{1}{3}s_{1}-\dfrac{1}{3}t_{1,3}-\dfrac{1}{3}\operatorname*{cyc}%
\nolimits_{1,2,3};\\
\mathbf{W}_{12\backslash\backslash3,\ 13\backslash\backslash2}  &  =-\dfrac
{1}{3}s_{1}+\dfrac{1}{3}s_{2}+\dfrac{1}{3}\operatorname*{cyc}\nolimits_{1,2,3}%
-\,\dfrac{1}{3}\operatorname*{cyc}\nolimits_{1,3,2};\\
\mathbf{W}_{13\backslash\backslash2,\ 12\backslash\backslash3}  &  =\dfrac
{1}{3}s_{2}-\dfrac{1}{3}t_{1,3}-\dfrac{1}{3}\operatorname*{cyc}%
\nolimits_{1,2,3}+\,\dfrac{1}{3}\operatorname*{cyc}\nolimits_{1,3,2};\\
\mathbf{W}_{13\backslash\backslash2,\ 13\backslash\backslash2}  &  =\dfrac
{1}{3}-\dfrac{1}{3}s_{1}+\dfrac{1}{3}t_{1,3}-\dfrac{1}{3}\operatorname*{cyc}%
\nolimits_{1,3,2};\\
\mathbf{W}_{1\backslash\backslash2\backslash\backslash3,\ 1\backslash
\backslash2\backslash\backslash3}  &  =\dfrac{1}{6}\nabla^{-}=\dfrac{1}%
{6}-\dfrac{1}{6}s_{1}-\dfrac{1}{6}s_{2}-\dfrac{1}{6}t_{1,3}+\dfrac{1}%
{6}\operatorname*{cyc}\nolimits_{1,2,3}+\,\dfrac{1}{6}\operatorname*{cyc}%
\nolimits_{1,3,2}.
\end{align*}

The attentive reader will have certainly noticed that these vectors are just
scalar multiples of the respective entries of the Young symmetrizer basis
$\overrightarrow{\operatorname*{ysb}}=\left(  \mathbf{E}_{U,V}\right)
_{\left(  \lambda,U,V\right)  \in\operatorname*{SBT}\left(  n\right)  }$ from
Example \ref{exa.bas.ysb.n=3}! The scalars are $\dfrac{1}{6}$ for the
$n$-bitableaux $\left(  \lambda,U,V\right)  $ with $\lambda=\left(  3\right)
$ or $\lambda=\left(  1,1,1\right)  $ and are $\dfrac{1}{3}$ for those with
$\lambda=\left(  2,1\right)  $. Thus, the matrices representing the
$\mathbf{k}$-linear maps $S$, $T_{\operatorname*{sign}}$, $L\left(
s_{1}\right)  $, $L\left(  s_{2}\right)  $, $R\left(  s_{1}\right)  $ and
$R\left(  s_{2}\right)  $ with respect to the basis
$\overrightarrow{\operatorname*{aws}}$ are exactly the same as for the basis
$\overrightarrow{\operatorname*{ysb}}$. (Normally, the fact that two different
scalars are involved would cause some entries to be rescaled; but this is
immaterial to us here, since these entries are all $0$.)
\end{example}

Alas, this pattern does not generalize. The vectors $\mathbf{W}_{U,V}$ are
scalar multiples of the $\mathbf{E}_{U,V}$ for all $n\leq4$, but not for
$n>4$. Here is an easy way to see the latter: For any $n\in\mathbb{N}$, the
basis $\overrightarrow{\operatorname*{aws}}=\left(  \mathbf{W}_{U,V}\right)
_{\left(  \lambda,U,V\right)  \in\operatorname*{SBT}\left(  n\right)  }$ has
the property that its entries multiply like elementary matrices:%
\begin{align*}
\mathbf{W}_{P,Q}\mathbf{W}_{Q,R}  &  =\mathbf{W}_{P,R};\\
\mathbf{W}_{P,Q}\mathbf{W}_{T,R}  &  =0\ \ \ \ \ \ \ \ \ \ \text{whenever
}Q\neq T.
\end{align*}
Up to a scalar factor, the Young symmetrizer basis
$\overrightarrow{\operatorname*{ysb}}=\left(  \mathbf{E}_{U,V}\right)
_{\left(  \lambda,U,V\right)  \in\operatorname*{SBT}\left(  n\right)  }$ also
satisfies the former of these two equalities (by Exercise
\ref{exe.specht.EPQ.EPQEQR}), but not the latter, unless $n\leq4$ (see Example
\ref{exa.specht.ET.ESET.1} and Exercise \ref{exe.specht.ET.ESET=0-when}, and
recall that $\mathbf{E}_{P}=\mathbf{E}_{P,P}$). Thus, the two bases cannot be
equal -- even up to scalar -- for $n\geq5$.

Nevertheless, the matrices representing the $\mathbf{k}$-linear maps $L\left(
\mathbf{a}\right)  $ for all $\mathbf{a}\in\mathcal{A}$ with respect to the
basis $\overrightarrow{\operatorname*{aws}}$ agree with those for the basis
$\overrightarrow{\operatorname*{ysb}}$ for general $n$:

\begin{exercise}
\label{exa.bas.aws.ysb}\fbox{3} Show that
$L_{\overrightarrow{\operatorname*{aws}}}\left(  \mathbf{a}\right)
=L_{\overrightarrow{\operatorname*{ysb}}}\left(  \mathbf{a}\right)  $ for all
$n\in\mathbb{N}$ and $\mathbf{a}\in\mathcal{A}$. \medskip

[\textbf{Hint:} By linearity, it suffices to prove this for $\mathbf{a}%
=\mathbf{W}_{P,Q}$. Show that $\mathbf{W}_{P,Q}\mathbf{W}_{U,V}=\delta
_{Q,U}\mathbf{W}_{P,V}$ and $\mathbf{W}_{P,Q}\mathbf{E}_{U,V}=\delta
_{Q,U}\mathbf{E}_{P,V}$.]
\end{exercise}

However, $R_{\overrightarrow{\operatorname*{aws}}}\left(  \mathbf{a}\right)  $
generally does not equal $R_{\overrightarrow{\operatorname*{ysb}}}\left(
\mathbf{a}\right)  $ for $n\in\mathbb{N}$ and $\mathbf{a}\in\mathcal{A}$. For
example, for $n=5$, we have $R_{\overrightarrow{\operatorname*{aws}}}\left(
s_{1}\right)  \neq R_{\overrightarrow{\operatorname*{ysb}}}\left(
s_{1}\right)  $ as well as $\left[  S\right]
_{\overrightarrow{\operatorname*{aws}}\rightarrow
\overrightarrow{\operatorname*{aws}}}\neq\left[  S\right]
_{\overrightarrow{\operatorname*{ysb}}\rightarrow
\overrightarrow{\operatorname*{ysb}}}$ and $\left[  T_{\operatorname*{sign}%
}\right]  _{\overrightarrow{\operatorname*{aws}}\rightarrow
\overrightarrow{\operatorname*{aws}}}\neq\left[  T_{\operatorname*{sign}%
}\right]  _{\overrightarrow{\operatorname*{ysb}}\rightarrow
\overrightarrow{\operatorname*{ysb}}}$.

\subsection{\label{sec.bas.tord}Combinatorial interlude: Lexicographic orders}

We now take a break from algebra and introduce some total orders on the sets
of partitions, of standard tableaux and of standard bitableaux.

\subsubsection{Partitions}

Recall that a partition is a weakly decreasing finite list of positive
integers, such as $\left(  3,2,2\right)  $. The length of a partition
$\lambda$ is denoted by $\ell\left(  \lambda\right)  $. If $\lambda$ is a
partition and $i$ is a positive integer, then $\lambda_{i}$ denotes the $i$-th
entry of $\lambda$; we understand this to mean $0$ if $i>\ell\left(
\lambda\right)  $ (see Definition \ref{def.partitions.partitions}
\textbf{(h)}). Thus, in a sense, we are pretending that a partition $\lambda$
has infinitely many entries (rather than just $\ell\left(  \lambda\right)  $
many), starting with its $\ell\left(  \lambda\right)  $ actual entries and
then continuing with infinitely many zeroes. This legal fiction will come in
handy in the following definition:

\begin{definition}
\label{def.bas.tord.pars}Define a binary relation $<$ on the set of all
partitions as follows:

Let $\lambda$ and $\mu$ be two partitions. Then, we declare that $\lambda<\mu$
if and only if

\begin{itemize}
\item there exists at least one $i\geq1$ such that $\lambda_{i}\neq\mu_{i}$;

\item the \textbf{smallest} such $i$ satisfies $\lambda_{i}<\mu_{i}$.
\end{itemize}

In other words, we say that $\lambda<\mu$ if the \textbf{leftmost} position in
which $\lambda$ and $\mu$ differ is occupied by a smaller entry in $\lambda$
than in $\mu$.
\end{definition}

For example:

\begin{itemize}
\item We have $\left(  5,3,3,1\right)  <\left(  5,4,1\right)  $, because the
leftmost position in which the partitions $\left(  5,3,3,1\right)  $ and
$\left(  5,4,1\right)  $ differ is the second position, and the former
partition has a smaller entry in this position than the latter partition does
(since $3<4$).

\item We have $\left(  5,2\right)  <\left(  5,2,1\right)  $, since the
leftmost position in which the partitions $\left(  5,2\right)  $ and $\left(
5,2,1\right)  $ differ is the third position (recall that we pretend that the
third entry of the partition $\left(  5,2\right)  $ is $0$), and the former
partition has a smaller entry in this position than the latter partition does
(since $0<1$).

\item The empty partition $\varnothing=\left(  {}\right)  $ satisfies
$\varnothing<\lambda$ for any partition $\lambda\neq\varnothing$, because all
entries of $\varnothing$ are $0$.
\end{itemize}

We shall prove the following fact shortly:

\begin{proposition}
\label{prop.bas.tord.pars.tord}The relation $<$ defined in Definition
\ref{def.bas.tord.pars} is the smaller relation of a total order on the set of
all partitions.
\end{proposition}

\begin{definition}
\label{def.bas.tord.pars.tord}This total order is called the
\emph{lexicographic order on partitions}. We shall use it whenever we need a
total order on the set of all partitions (unless we explicitly introduce a
different order). Thus, the symbol $<$ when applied to partitions will always
be understood as in Definition \ref{def.bas.tord.pars}. Moreover, the symbols
$\leq$, $>$ and $\geq$ will be understood accordingly (e.g., we let
\textquotedblleft$\lambda\leq\mu$\textquotedblright\ mean \textquotedblleft%
$\lambda=\mu$ or $\lambda<\mu$\textquotedblright, and we let \textquotedblleft%
$\lambda>\mu$\textquotedblright\ mean \textquotedblleft$\mu<\lambda
$\textquotedblright).
\end{definition}

\begin{example}
\textbf{(a)} Here are all seven partitions of $5$, ordered in lexicographic
order:%
\[
\left(  1,1,1,1,1\right)  <\left(  2,1,1,1\right)  <\left(  2,2,1\right)
<\left(  3,1,1\right)  <\left(  3,2\right)  <\left(  4,1\right)  <\left(
5\right)  .
\]

\textbf{(b)} Here is a chain of some partitions ordered in lexicographic
order:
\[
\left(  {}\right)  <\left(  1\right)  <\left(  1,1\right)  <\left(  2\right)
<\left(  2,1\right)  <\left(  2,2\right)  <\left(  2,2,1\right)  <\left(
3\right)  <\left(  3,1\right)  .
\]
It is easy to extend this chain to the right. We can also interject any number
of further partitions between $\left(  1\right)  $ and $\left(  2\right)  $:
namely,%
\[
\left(  1\right)  <\left(  1,1\right)  <\left(  1,1,1\right)  <\left(
1,1,1,1\right)  <\cdots<\left(  2\right)  .
\]

\end{example}

\begin{proof}
[Proof of Proposition \ref{prop.bas.tord.pars.tord} (sketched).]For a detailed
proof, see the Appendix (Section \ref{sec.details.bas.tord}). In a nutshell:
The relation $<$ is one of many instances of a lexicographic order (i.e., a
way of comparing tuples by first comparing the first entries, then resolving
ties using the second entries, then resolving any remaining ties using the
third entries, and so on). The only unusual feature of our current situation
is that different partitions can have different lengths; but we consider them
as having infinitely many entries (thanks to Definition
\ref{def.partitions.partitions} \textbf{(h)}).
\end{proof}

\begin{exercise}
\fbox{1} Let $\lambda$ be a partition of $n$. Let $k\in\left[  n\right]  $.
Prove that $\lambda_{1}\geq k$ holds if and only if $\lambda\geq\left(
k,\underbrace{1,1,\ldots,1}_{n-k\text{ times}}\right)  $ in the lexicographic
order on partitions.
\end{exercise}

\begin{exercise}
We shall use \emph{exponential notation}, meaning that a sequence
$\underbrace{p,p,\ldots,p}_{k\text{ times}}$ of equal numbers will be denoted
by $p^{k}$ (not to be mistaken for the $k$-th power of $p$). Thus, for
example, the partition $\left(  4,4,3,1,1,1\right)  $ can be written as
$\left(  4^{2},3,1^{3}\right)  $ or even as $\left(  4^{2},3^{1},1^{3}\right)
$.

Consider the lexicographic order on partitions, restricted to the set of all
partitions of $n$. \medskip

\textbf{(a)} \fbox{1} Prove that the smallest $\left\lfloor n/2\right\rfloor
+1$ partitions in this order are%
\[
\left(  1^{n}\right)  <\left(  2,1^{n-2}\right)  <\left(  2,2,1^{n-4}\right)
<\cdots<\left(  2^{\left\lfloor n/2\right\rfloor },1\right)  \text{ or
}\left(  2^{\left\lfloor n/2\right\rfloor }\right)
\]
(the last partition here is $\left(  2^{\left\lfloor n/2\right\rfloor
},1\right)  $ if $n$ is odd and $\left(  2^{\left\lfloor n/2\right\rfloor
}\right)  $ if $n$ is even). \medskip

\textbf{(b)} \fbox{1} What are the largest four partitions in this order
(assuming $n\geq4$) ?
\end{exercise}

The lexicographic order is just one of many orders that can be imposed on
partitions. Several other orders are known (both total and partial orders),
such as the colexicographic order, the dominance order (see, e.g.,
\cite[Definition 2.2.7]{GriRei14}, \cite[\S 1.4]{JamKer81} or \cite[Definition
4.1]{Wildon18}) and the refinement order (see, e.g., \cite[\S I.6, around
(6.10)]{Macdon95}). But we shall have no use for these other orders in the
present chapter.

\subsubsection{Standard tableaux of partition shapes}

Next, we move on from partitions to tableaux.

If $\kappa$ is a partition and if $U$ is an $n$-tableau of shape $Y\left(
\kappa\right)  $, then we can uniquely recover $\kappa$ from $U$ (in fact,
each entry $\kappa_{i}$ of $\kappa$ is the number of entries in the $i$-th row
of $U$). Hence, the sets $\operatorname*{SYT}\left(  \kappa\right)  $ for
different partitions $\kappa$ are disjoint. Therefore, the union
$\bigcup\limits_{\kappa\vdash n}\operatorname*{SYT}\left(  \kappa\right)  $ is
a disjoint union. (Recall that \textquotedblleft$\kappa\vdash n$%
\textquotedblright\ means \textquotedblleft$\kappa$ is a partition of
$n$\textquotedblright.) Explicitly, the set $\bigcup\limits_{\kappa\vdash
n}\operatorname*{SYT}\left(  \kappa\right)  $ consists of all standard
$n$-tableaux of partition shape (i.e., of shape $Y\left(  \kappa\right)  $ for
some partition $\kappa$ of $n$). This set $\bigcup\limits_{\kappa\vdash
n}\operatorname*{SYT}\left(  \kappa\right)  $ is finite (since all sets
$\operatorname*{SYT}\left(  \kappa\right)  $ are finite, and there are
finitely many partitions $\kappa$ of $n$).

For example, for $n=3$, the set $\bigcup\limits_{\kappa\vdash n}%
\operatorname*{SYT}\left(  \kappa\right)  $ consists of the four standard
tableaux $123,\ \ 12\backslash\backslash3,\ \ 13\backslash\backslash
2,\ \ 1\backslash\backslash2\backslash\backslash3$ (where we are using the
shorthand notation from Convention \ref{conv.tableau.poetic}).

Now we shall define a total order on the set $\bigcup\limits_{\kappa\vdash
n}\operatorname*{SYT}\left(  \kappa\right)  $ as follows:

\begin{definition}
\label{def.bas.tord.syt}We define a binary relation $<$ on the set
$\bigcup\limits_{\kappa\vdash n}\operatorname*{SYT}\left(  \kappa\right)  $ as follows:

Let $P,Q\in\bigcup\limits_{\kappa\vdash n}\operatorname*{SYT}\left(
\kappa\right)  $ be two $n$-tableaux. Thus, $P\in\operatorname*{SYT}\left(
\lambda\right)  $ and $Q\in\operatorname*{SYT}\left(  \mu\right)  $ for some
partitions $\lambda$ and $\mu$ of $n$. Consider these two partitions $\lambda$
and $\mu$. (They are uniquely determined by $P$ and $Q$, since the sets
$\operatorname*{SYT}\left(  \kappa\right)  $ for different partitions $\kappa$
are disjoint.) Consider also the $n$-tabloids $\overline{P}$ and $\overline
{Q}$ of shapes $Y\left(  \lambda\right)  $ and $Y\left(  \mu\right)  $. Now,
we declare that $P<Q$ if and only if%
\[
\left(  \lambda<\mu\text{ or }\left(  \lambda=\mu\text{ and }\overline
{P}<\overline{Q}\right)  \right)  ,
\]
where:

\begin{itemize}
\item the \textquotedblleft$<$\textquotedblright\ sign in \textquotedblleft%
$\lambda<\mu$\textquotedblright\ refers to the lexicographic order on
partitions (defined in Definition \ref{def.bas.tord.pars});

\item the \textquotedblleft$<$\textquotedblright\ sign in \textquotedblleft%
$\overline{P}<\overline{Q}$\textquotedblright\ refers to the Young last letter
order on the set $\left\{  n\text{-tabloids of shape }Y\left(  \lambda\right)
\right\}  $ (defined in Definition \ref{def.tabloid.llo} \textbf{(b)}), which
makes sense here because the preceding statement $\lambda=\mu$ ensures that
both $\overline{P}$ and $\overline{Q}$ are $n$-tabloids of shape $Y\left(
\lambda\right)  $ (in fact, $\overline{Q}$ is an $n$-tabloid of shape
$Y\left(  \mu\right)  $, but because of $\lambda=\mu$, this is the same as
being an $n$-tabloid of shape $Y\left(  \lambda\right)  $).
\end{itemize}
\end{definition}

For example, for $n=5$:

\begin{itemize}
\item We have $\ytableaushort{124,3,5} < \ytableaushort{123,45}\ $, because
the shapes $\left(  3,1,1\right)  $ and $\left(  3,2\right)  $ of these two
tableaux satisfy $\left(  3,1,1\right)  <\left(  3,2\right)  $.

\item We have $\ytableaushort{135,24} < \ytableaushort{123,45}\ $, because the
shapes $\left(  3,2\right)  $ and $\left(  3,2\right)  $ of these two tableaux
are equal whereas the corresponding $n$-tabloids satisfy
$\ytableausetup{tabloids}\ytableaushort{135,24} <
\ytableaushort{123,45}\ytableausetup{notabloids}$ .
\end{itemize}

We shall prove the following fact shortly:

\begin{proposition}
\label{prop.bas.tord.syt.tord}The relation $<$ defined in Definition
\ref{def.bas.tord.syt} is the smaller relation of a total order on the set
$\bigcup\limits_{\kappa\vdash n}\operatorname*{SYT}\left(  \kappa\right)  $.
\end{proposition}

\begin{definition}
\label{def.bas.tord.syt.tord}We shall use this total order whenever we need a
total order on the set $\bigcup\limits_{\kappa\vdash n}\operatorname*{SYT}%
\left(  \kappa\right)  $ (unless we explicitly introduce a different order).
Thus, the symbol $<$ when applied to standard $n$-tableaux will always be
understood as in Definition \ref{def.bas.tord.syt}. Moreover, the symbols
$\leq$, $>$ and $\geq$ will be understood accordingly (e.g., we let
\textquotedblleft$P\leq Q$\textquotedblright\ mean \textquotedblleft$P=Q$ or
$P<Q$\textquotedblright, and we let \textquotedblleft$P>Q$\textquotedblright%
\ mean \textquotedblleft$Q<P$\textquotedblright).
\end{definition}

\begin{example}
\label{exa.bas.tord.syt.n=3}For $n=3$, here are the elements of the set
$\bigcup\limits_{\kappa\vdash n}\operatorname*{SYT}\left(  \kappa\right)  $
listed in the order we just introduced:%
\[
\ytableaushort{1,2,3}<\ytableaushort{13,2}<\ytableaushort{12,3}<\ytableaushort{123}\ \ .
\]

\end{example}

\begin{proof}
[Proof of Proposition \ref{prop.bas.tord.syt.tord}.]See the appendix (Section
\ref{sec.details.bas.tord}) for a detailed proof. Here is a quick outline:
Again, the relation $<$ is a kind of lexicographic order, although not in the
most direct way, since it is defined not on tuples but rather on single
standard tableaux. To reveal its lexicographic nature, we need to re-encode
any standard $n$-tableau $T\in\bigcup\limits_{\kappa\vdash n}%
\operatorname*{SYT}\left(  \kappa\right)  $ as the pair $\left(
\lambda,\overline{T}\right)  $, where $\lambda$ is the partition of $n$
satisfying $T\in\operatorname*{SYT}\left(  \lambda\right)  $, and where
$\overline{T}$ is the $n$-tabloid of $T$. This encoding $T\mapsto\left(
\lambda,\overline{T}\right)  $ is injective, since each standard (or even just
row-standard) $n$-tableau $T$ can be uniquely reconstructed from its
$n$-tabloid $\overline{T}$ (see Proposition \ref{prop.tabloid.row-st}). Once
we replace the tableaux $T\in\bigcup\limits_{\kappa\vdash n}%
\operatorname*{SYT}\left(  \kappa\right)  $ by the corresponding pairs
$\left(  \lambda,\overline{T}\right)  $ (where $\lambda$ is the partition such
that $T\in\operatorname*{SYT}\left(  \lambda\right)  $), our definition of the
relation $<$ takes the form
\[
\left(  \left(  \lambda,\overline{P}\right)  <\left(  \mu,\overline{Q}\right)
\right)  \ \Longleftrightarrow\ \left(  \lambda<\mu\text{ or }\left(
\lambda=\mu\text{ and }\overline{P}<\overline{Q}\right)  \right)  ;
\]
but this is just a lexicographic order on pairs, and thus is well-known to be
a total order. So the original relation $<$ is a total order as well.
\end{proof}

The total order we constructed on $\bigcup\limits_{\kappa\vdash n}%
\operatorname*{SYT}\left(  \kappa\right)  $ has a rather convenient property:

\begin{lemma}
[Combined zero-sandwich lemma]\label{lem.bas.tord.syt.czs}Let $P$ and $Q$ be
two standard $n$-tableaux in $\bigcup\limits_{\kappa\vdash n}%
\operatorname*{SYT}\left(  \kappa\right)  $. Assume that $P<Q$ (where the
relation $<$ is the one defined in Definition \ref{def.bas.tord.syt}). Then,
\[
\nabla_{\operatorname*{Col}P}^{-}\nabla_{\operatorname*{Row}Q}%
=0\ \ \ \ \ \ \ \ \ \ \text{and}\ \ \ \ \ \ \ \ \ \ \nabla
_{\operatorname*{Row}Q}\nabla_{\operatorname*{Col}P}^{-}=0.
\]

\end{lemma}

\begin{proof}
We have $P\in\bigcup\limits_{\kappa\vdash n}\operatorname*{SYT}\left(
\kappa\right)  $. That is, we have $P\in\operatorname*{SYT}\left(
\lambda\right)  $ for some partition $\lambda$ of $n$. Similarly, we have
$Q\in\operatorname*{SYT}\left(  \mu\right)  $ for some partition $\mu$ of $n$.
Consider these two partitions $\lambda$ and $\mu$. Thus, $P$ is an $n$-tableau
of shape $Y\left(  \lambda\right)  $, whereas $Q$ is an $n$-tableau of shape
$Y\left(  \mu\right)  $.

Recall that $P<Q$. By Definition \ref{def.bas.tord.syt}, this means that
\[
\left(  \lambda<\mu\text{ or }\left(  \lambda=\mu\text{ and }\overline
{P}<\overline{Q}\right)  \right)
\]
(since $P\in\operatorname*{SYT}\left(  \lambda\right)  $ and $Q\in
\operatorname*{SYT}\left(  \mu\right)  $). Hence, we are in one of the
following two cases:

\textit{Case 1:} We have $\lambda<\mu$.

\textit{Case 2:} We have $\lambda=\mu$ and $\overline{P}<\overline{Q}$.

Let us first consider Case 1. In this case, we have $\lambda<\mu$. By
Definition \ref{def.bas.tord.pars}, this means that

\begin{itemize}
\item there exists at least one $i\geq1$ such that $\lambda_{i}\neq\mu_{i}$;

\item the \textbf{smallest} such $i$ satisfies $\lambda_{i}<\mu_{i}$.
\end{itemize}

Consider this smallest $i$. Thus, $i\geq1$ and $\lambda_{i}<\mu_{i}$, but each
positive integer $j<i$ satisfies
\begin{equation}
\lambda_{j}=\mu_{j} \label{pf.lem.bas.tord.syt.czs.c1.1}%
\end{equation}
(since $i$ is the \textbf{smallest }$i\geq1$ such that $\lambda_{i}\neq\mu
_{i}$). Adding the equalities (\ref{pf.lem.bas.tord.syt.czs.c1.1}) for all
$j\in\left\{  1,2,\ldots,i-1\right\}  $ together, we find
\[
\lambda_{1}+\lambda_{2}+\cdots+\lambda_{i-1}=\mu_{1}+\mu_{2}+\cdots+\mu
_{i-1}.
\]
Adding the inequality $\lambda_{i}<\mu_{i}$ to this equality, we find%
\[
\left(  \lambda_{1}+\lambda_{2}+\cdots+\lambda_{i-1}\right)  +\lambda
_{i}<\left(  \mu_{1}+\mu_{2}+\cdots+\mu_{i-1}\right)  +\mu_{i}.
\]
In other words,%
\[
\lambda_{1}+\lambda_{2}+\cdots+\lambda_{i}<\mu_{1}+\mu_{2}+\cdots+\mu_{i}.
\]
In other words, $\mu_{1}+\mu_{2}+\cdots+\mu_{i}>\lambda_{1}+\lambda_{2}%
+\cdots+\lambda_{i}$. Therefore:

\begin{itemize}
\item Theorem \ref{thm.specht.sandw0} \textbf{(b)} (applied to $\mu$,
$\lambda$, $Q$ and $P$ instead of $\lambda$, $\mu$, $S$ and $T$) yields
$\nabla_{\operatorname*{Col}P}^{-}\mathbf{a}\nabla_{\operatorname*{Row}Q}=0$
for each $\mathbf{a}\in\mathbf{k}\left[  S_{n}\right]  $. Applying this to
$\mathbf{a}=1$, we obtain $\nabla_{\operatorname*{Col}P}^{-}1\nabla
_{\operatorname*{Row}Q}=0$. In other words, $\nabla_{\operatorname*{Col}P}%
^{-}\nabla_{\operatorname*{Row}Q}=0$.

\item Theorem \ref{thm.specht.sandw0} \textbf{(d)} (applied to $\mu$,
$\lambda$, $Q$ and $P$ instead of $\lambda$, $\mu$, $S$ and $T$) yields
$\nabla_{\operatorname*{Row}Q}\mathbf{a}\nabla_{\operatorname*{Col}P}^{-}=0$
for each $\mathbf{a}\in\mathbf{k}\left[  S_{n}\right]  $. Applying this to
$\mathbf{a}=1$, we obtain $\nabla_{\operatorname*{Row}Q}1\nabla
_{\operatorname*{Col}P}^{-}=0$. In other words, $\nabla_{\operatorname*{Row}%
Q}\nabla_{\operatorname*{Col}P}^{-}=0$.
\end{itemize}

Hence, the proof of Lemma \ref{lem.bas.tord.syt.czs} is complete in Case 1.

Now, let us consider Case 2. In this case, we have $\lambda=\mu$ and
$\overline{P}<\overline{Q}$. From $\lambda=\mu$, we obtain
$\operatorname*{SYT}\left(  \lambda\right)  =\operatorname*{SYT}\left(
\mu\right)  $, so that $Q\in\operatorname*{SYT}\left(  \mu\right)
=\operatorname*{SYT}\left(  \lambda\right)  $. Hence, both tableaux $P$ and
$Q$ belong to $\operatorname*{SYT}\left(  \lambda\right)  $, that is, are
$n$-tableaux of shape $Y\left(  \lambda\right)  $. Moreover, they are
standard, hence column-standard. Now, recall that $\overline{P}<\overline{Q}$
in the Young last letter order. Hence, Lemma \ref{lem.specht.sandw0-llo}
\textbf{(c)} (applied to $S=Q$ and $T=P$) yields
\[
\nabla_{\operatorname*{Col}P}^{-}\nabla_{\operatorname*{Row}Q}%
=0\ \ \ \ \ \ \ \ \ \ \text{and}\ \ \ \ \ \ \ \ \ \ \nabla
_{\operatorname*{Row}Q}\nabla_{\operatorname*{Col}P}^{-}=0.
\]
Thus, Lemma \ref{lem.bas.tord.syt.czs} is proved in Case 2.

We have now proved Lemma \ref{lem.bas.tord.syt.czs} in both Cases 1 and 2.
Hence, the proof is complete.
\end{proof}

\subsubsection{Standard $n$-bitableaux}

Finally, we define a total order on the set $\operatorname*{SBT}\left(
n\right)  $ of all standard $n$-bitableaux. We recall that any standard
$n$-bitableau $\left(  \lambda,U,V\right)  \in\operatorname*{SBT}\left(
n\right)  $ satisfies $U,V\in\operatorname*{SYT}\left(  \lambda\right)
\subseteq\bigcup\limits_{\kappa\vdash n}\operatorname*{SYT}\left(
\kappa\right)  $.

\begin{definition}
\label{def.bas.tord.sbt}We define a binary relation $<$ on the set
$\operatorname*{SBT}\left(  n\right)  $ as follows:

For two standard $n$-bitableaux $\left(  \lambda,A,B\right)  $ and $\left(
\mu,C,D\right)  $ in $\operatorname*{SBT}\left(  n\right)  $, we declare that
$\left(  \lambda,A,B\right)  <\left(  \mu,C,D\right)  $ if and only if%
\[
\left(  A<C\text{ or }\left(  A=C\text{ and }B<D\right)  \right)  .
\]
Here, the symbol \textquotedblleft$<$\textquotedblright\ in \textquotedblleft%
$A<C$\textquotedblright\ and \textquotedblleft$B<D$\textquotedblright\ refers
to the binary relation $<$ introduced in Definition \ref{def.bas.tord.syt}.
\end{definition}

For example, for $n=5$:

\begin{itemize}
\item We have
\[
\left(  \left(  3,2\right)
,\ \ytableaushort{135,24},\ \ytableaushort{134,25}\right)  < \left(  \left(
3,2\right)  ,\ \ytableaushort{123,45},\ \ytableaushort{135,24}\right)  \ ,
\]
because the first tableaux in these two $n$-bitableaux satisfy
$\ytableaushort{135,24} < \ytableaushort{123,45}\ $.

\item We have
\[
\left(  \left(  3,2\right)
,\ \ytableaushort{134,25},\ \ytableaushort{135,24}\right)  < \left(  \left(
3,2\right)  ,\ \ytableaushort{134,25},\ \ytableaushort{123,45}\right)  \ ,
\]
because the first tableaux in these two $n$-bitableaux are equal whereas the
second tableaux satisfy $\ytableaushort{135,24} < \ytableaushort{123,45}\ $.
\end{itemize}

We shall prove the following fact shortly:

\begin{proposition}
\label{prop.bas.tord.sbt.tord}The relation $<$ defined in Definition
\ref{def.bas.tord.sbt} is the smaller relation of a total order on the set
$\operatorname*{SBT}\left(  n\right)  $.
\end{proposition}

\begin{definition}
\label{def.bas.tord.sbt.tord}We shall use this total order whenever we need a
total order on the set $\operatorname*{SBT}\left(  n\right)  $ (unless we
explicitly introduce a different order). Thus, the symbol $<$ when applied to
standard $n$-bitableaux will always be understood as in Definition
\ref{def.bas.tord.sbt}. Moreover, the symbols $\leq$, $>$ and $\geq$ will be
understood accordingly (e.g., we let \textquotedblleft$\mathfrak{P}%
\leq\mathfrak{Q}$\textquotedblright\ mean \textquotedblleft$\mathfrak{P}%
=\mathfrak{Q}$ or $\mathfrak{P}<\mathfrak{Q}$\textquotedblright).
\end{definition}

\begin{example}
\label{exa.bas.tord.sbt.n=3}For $n=3$, here are the six standard
$n$-bitableaux in $\operatorname*{SBT}\left(  n\right)  $ listed in the order
we just introduced:%
\begin{align*}
&  \left(  \left(  1,1,1\right)
,\ \ \ytableaushort{1,2,3}\ ,\ \ \ytableaushort{1,2,3}\right)  ,\\
&  \left(  \left(  2,1\right)
,\ \ \ytableaushort{13,2}\ ,\ \ \ytableaushort{13,2}\right)
,\ \ \ \ \ \ \ \ \ \ \left(  \left(  2,1\right)
,\ \ \ytableaushort{13,2}\ ,\ \ \ytableaushort{12,3}\right)  ,\\
&  \left(  \left(  2,1\right)
,\ \ \ytableaushort{12,3}\ ,\ \ \ytableaushort{13,2}\right)
,\ \ \ \ \ \ \ \ \ \ \left(  \left(  2,1\right)
,\ \ \ytableaushort{12,3}\ ,\ \ \ytableaushort{12,3}\right)  ,\\
&  \left(  \left(  3\right)
,\ \ \ytableaushort{123}\ ,\ \ \ytableaushort{123}\right)  .
\end{align*}

\end{example}

\begin{proof}
[Proof of Proposition \ref{prop.bas.tord.sbt.tord}.]See the appendix (Section
\ref{sec.details.bas.tord}) for a detailed proof. But again, this relation is
just a lexicographic order in disguise. Indeed, a standard $n$-bitableau
$\left(  \lambda,U,V\right)  \in\operatorname*{SBT}\left(  n\right)  $ can be
re-encoded as the pair $\left(  U,V\right)  \in\left(  \bigcup\limits_{\kappa
\vdash n}\operatorname*{SYT}\left(  \kappa\right)  \right)  ^{2}$, because the
partition $\lambda$ is uniquely determined by the tableau $U$. But the
condition that defines the binary relation $<$ on $\operatorname*{SBT}\left(
n\right)  $, namely%
\[
\left(  A<C\text{ or }\left(  A=C\text{ and }B<D\right)  \right)  ,
\]
is precisely the lexicographic inequality $\left(  A,B\right)  <\left(
C,D\right)  $ between the latter pairs.
\end{proof}

Let us state another simple property of the relation $<$ from Definition
\ref{def.bas.tord.sbt}:

\begin{lemma}
\label{lem.bas.tord.sbt.1}Let $\left(  \lambda,A,B\right)  $ and $\left(
\mu,C,D\right)  $ be two standard $n$-bitableaux in $\operatorname*{SBT}%
\left(  n\right)  $ satisfying $\left(  \lambda,A,B\right)  <\left(
\mu,C,D\right)  $ (with respect to the relation $<$ from Definition
\ref{def.bas.tord.sbt}). Then: \medskip

\textbf{(a)} We have $A\leq C$. (Here, the symbol \textquotedblleft$\leq
$\textquotedblright\ refers to the total order on $\bigcup\limits_{\kappa
\vdash n}\operatorname*{SYT}\left(  \kappa\right)  $ introduced in Definition
\ref{def.bas.tord.syt}.) \medskip

\textbf{(b)} We have $A<C$ or $B<D$. (Here, the symbol \textquotedblleft%
$<$\textquotedblright\ refers to the total order on $\bigcup\limits_{\kappa
\vdash n}\operatorname*{SYT}\left(  \kappa\right)  $ introduced in Definition
\ref{def.bas.tord.syt}.)
\end{lemma}

\begin{proof}
See the appendix (Section \ref{sec.details.bas.tord}) for the very simple proof.
\end{proof}

\subsection{\label{sec.bas.mur}The Murphy bases}

We now come to some new bases of $\mathcal{A}$. The \emph{Murphy bases} (also
known as the \emph{Murphy cellular bases}) are two closely related bases of
$\mathcal{A}$ that share a few of the nice properties of the Young symmetrizer
basis (Section \ref{sec.bas.ysb}), the FwE basis (Section \ref{sec.bas.FwE})
and the Artin--Wedderburn bases (Section \ref{sec.bas.AWS}), although in a
weakened form (the matrices representing the $L\left(  \mathbf{a}\right)  $
and $R\left(  \mathbf{a}\right)  $ maps shall only be block-triangular, not
block-diagonal). In exchange, the bases are much simpler to define, and exist
\textbf{for every commutative ring} $\mathbf{k}$. In a sense, these bases are
the best one can hope for without assuming that $n!$ be invertible in
$\mathbf{k}$.

\subsubsection{\label{subsec.bas.mur.n=3}Introductory example: $n=3$}

We have already seen one of the Murphy bases in Subsection
\ref{subsec.rep.maschke.kS3} for $n=3$. This was the basis
\[
\overrightarrow{\operatorname*{rmu}}_{3}:=\left(  \nabla,\ \nabla
_{3,3},\ \nabla_{2,3},\ \nabla_{3,2},\ \nabla_{2,2},\ 1\right)  ,
\]
where we set
\[
\nabla_{b,a}:=\sum_{\substack{w\in S_{n};\\w\left(  a\right)  =b}%
}w\ \ \ \ \ \ \ \ \ \ \text{for all }a,b\in\left[  n\right]  .
\]
We called this basis the \textquotedblleft row Murphy basis\textquotedblright%
\ of $\mathbf{k}\left[  S_{3}\right]  $. It consists of the six vectors%
\begin{align*}
\nabla &  =1+s_{1}+s_{2}+t_{1,3}+\operatorname*{cyc}\nolimits_{1,2,3}%
+\operatorname*{cyc}\nolimits_{1,3,2},\\
\nabla_{3,3}  &  =1+s_{1},\\
\nabla_{2,3}  &  =s_{2}+\operatorname*{cyc}\nolimits_{1,3,2},\\
\nabla_{3,2}  &  =s_{2}+\operatorname*{cyc}\nolimits_{1,2,3},\\
\nabla_{2,2}  &  =1+t_{1,3},\\
1  &  =1.
\end{align*}

The matrices representing the $\mathbf{k}$-linear maps $S$,
$T_{\operatorname*{sign}}$, $L\left(  s_{1}\right)  $, $L\left(  s_{2}\right)
$, $R\left(  s_{1}\right)  $ and $R\left(  s_{2}\right)  $ with respect to
this basis are as follows:%
\[
\left[  S\right]  _{\overrightarrow{\operatorname*{rmu}}\rightarrow
\overrightarrow{\operatorname*{rmu}}}=\left(
\begin{array}
[c]{cccccc}%
1 & 0 & 0 & 0 & 0 & 0\\
0 & 1 & 0 & 0 & 0 & 0\\
0 & 0 & 0 & 1 & 0 & 0\\
0 & 0 & 1 & 0 & 0 & 0\\
0 & 0 & 0 & 0 & 1 & 0\\
0 & 0 & 0 & 0 & 0 & 1
\end{array}
\right)  ;
\]%
\[
\left[  T_{\operatorname*{sign}}\right]  _{\overrightarrow{\operatorname*{rmu}%
}\rightarrow\overrightarrow{\operatorname*{rmu}}}=\left(
\begin{array}
[c]{cccccc}%
3 & 0 & 2 & 2 & 0 & 0\\
-4 & -1 & -2 & -2 & 0 & 0\\
-2 & 0 & -1 & -2 & 0 & 0\\
-2 & 0 & -2 & -1 & 0 & 0\\
-4 & 0 & -2 & -2 & -1 & 0\\
6 & 2 & 2 & 2 & 2 & 1
\end{array}
\right)  ;
\]%
\[
L_{\overrightarrow{\operatorname*{rmu}}}\left(  s_{1}\right)  =\left[
L\left(  s_{1}\right)  \right]  _{\overrightarrow{\operatorname*{rmu}%
}\rightarrow\overrightarrow{\operatorname*{rmu}}}=\left(
\begin{array}
[c]{cccccc}%
1 & 0 & 1 & 0 & 1 & 0\\
0 & 1 & -1 & 0 & 0 & 1\\
0 & 0 & -1 & 0 & 0 & 0\\
0 & 0 & 0 & 1 & -1 & 0\\
0 & 0 & 0 & 0 & -1 & 0\\
0 & 0 & 0 & 0 & 0 & -1
\end{array}
\right)  ;
\]%
\[
L_{\overrightarrow{\operatorname*{rmu}}}\left(  s_{2}\right)  =\left[
L\left(  s_{2}\right)  \right]  _{\overrightarrow{\operatorname*{rmu}%
}\rightarrow\overrightarrow{\operatorname*{rmu}}}=\left(
\begin{array}
[c]{cccccc}%
1 & 0 & 0 & 0 & 0 & -1\\
0 & 0 & 1 & 0 & 0 & 1\\
0 & 1 & 0 & 0 & 0 & 1\\
0 & 0 & 0 & 0 & 1 & 1\\
0 & 0 & 0 & 1 & 0 & 1\\
0 & 0 & 0 & 0 & 0 & -1
\end{array}
\right)  ;
\]%
\[
R_{\overrightarrow{\operatorname*{rmu}}}\left(  s_{1}\right)  =\left[
R\left(  s_{1}\right)  \right]  _{\overrightarrow{\operatorname*{rmu}%
}\rightarrow\overrightarrow{\operatorname*{rmu}}}= \left(
\begin{array}
[c]{cccccc}%
1 & 0 & 0 & 1 & 1 & 0\\
0 & 1 & 0 & -1 & 0 & 1\\
0 & 0 & 1 & 0 & -1 & 0\\
0 & 0 & 0 & -1 & 0 & 0\\
0 & 0 & 0 & 0 & -1 & 0\\
0 & 0 & 0 & 0 & 0 & -1
\end{array}
\right)  ;
\]%
\[
R_{\overrightarrow{\operatorname*{rmu}}}\left(  s_{2}\right)  =\left[
R\left(  s_{2}\right)  \right]  _{\overrightarrow{\operatorname*{rmu}%
}\rightarrow\overrightarrow{\operatorname*{rmu}}}= \left(
\begin{array}
[c]{cccccc}%
1 & 0 & 0 & 0 & 0 & -1\\
0 & 0 & 0 & 1 & 0 & 1\\
0 & 0 & 0 & 0 & 1 & 1\\
0 & 1 & 0 & 0 & 0 & 1\\
0 & 0 & 1 & 0 & 0 & 1\\
0 & 0 & 0 & 0 & 0 & -1
\end{array}
\right)  .
\]

We observe a few properties:

\begin{itemize}
\item All the matrices have integer entries. This is because the basis is
defined over any commutative ring $\mathbf{k}$, including $\mathbb{Z}$.

\item The matrix $\left[  S\right]  _{\overrightarrow{\operatorname*{rmu}%
}\rightarrow\overrightarrow{\operatorname*{rmu}}}$ is a permutation matrix. In
other words, the antipode $S$ permutes the row Murphy basis.

\item The matrix $\left[  T_{\operatorname*{sign}}\right]
_{\overrightarrow{\operatorname*{rmu}}\rightarrow
\overrightarrow{\operatorname*{rmu}}}$ is not so nice. Indeed, the sign-twist
$T_{\operatorname*{sign}}$ does not permute the row Murphy basis, but rather
sends it to (a permuted) column Murphy basis, which we will soon define as well.

\item The matrices $L_{\overrightarrow{\operatorname*{rmu}}}\left(
\mathbf{a}\right)  $ and $R_{\overrightarrow{\operatorname*{rmu}}}\left(
\mathbf{a}\right)  $ are not block-diagonal any more, but they are
block-upper-triangular (if the rows/columns are appropriately ordered).
\end{itemize}

We will now define a row Murphy basis for all $n\in\mathbb{N}$, and then show
that it satisfies all these properties.

\subsubsection{\label{subsec.bas.mur.def}Definitions}

We shall define the Murphy bases soon. First, we recall the notions of an
$n$-tabloid (Definition \ref{def.tabloid.tabloid}) and of an $n$%
-column-tabloid (Definition \ref{def.coltabloid.coltabloid}). For any diagram
$D$, we introduced a left $S_{n}$-action on the set $\left\{  n\text{-tabloids
of shape }D\right\}  $ in Definition \ref{def.tabloid.Sn-act}. We now make a
similar definition for $n$-column-tabloids:

\begin{definition}
\label{def.coltabloid.Sn-act}Let $D$ be any diagram. Then, we define a left
$S_{n}$-action on the set $\left\{  n\text{-column-tabloids of shape
}D\right\}  $ by setting%
\[
w\rightharpoonup\widetilde{T}:=\widetilde{w\rightharpoonup T}%
=\widetilde{w\circ T}\ \ \ \ \ \ \ \ \ \ \text{for all }w\in S_{n}\text{ and
all }n\text{-tableaux }T.
\]
This action is well-defined (we will justify this in a moment), and is again
called the \emph{action on entries} (for column-tabloids). This action makes
the set $\left\{  n\text{-column-tabloids of shape }D\right\}  $ into a left
$S_{n}$-set.
\end{definition}

The reason why this action is well-defined is the same as for the action in
Definition \ref{def.tabloid.Sn-act} (except that now, \textquotedblleft
row-equivalence\textquotedblright\ must be replaced by \textquotedblleft
column-equivalence\textquotedblright).

Next, we define a lot of elements of $\mathcal{A}$, some (not all) of which
will constitute the Murphy bases:

\begin{definition}
\label{def.bas.mur.mur}Let $P$ and $Q$ be two $n$-tableaux of the same shape.
Then: \medskip

\textbf{(a)} We define the \emph{row-to-row sum} $\nabla_{P,Q}%
^{\operatorname*{Row}}$ to be the element
\begin{equation}
\nabla_{P,Q}^{\operatorname*{Row}}:=\sum_{\substack{u\in S_{n};\\u\overline
{Q}=\overline{P}}}u\in\mathcal{A} \label{eq.def.bas.mur.mur.Row}%
\end{equation}
(recall that $\overline{T}$ denotes the $n$-tabloid corresponding to the
$n$-tableau $T$, as defined in Definition \ref{def.tabloid.tabloid}). \medskip

\textbf{(b)} We define the \emph{column-to-column signed sum} $\nabla
_{P,Q}^{-\operatorname*{Col}}$ by%
\begin{equation}
\nabla_{P,Q}^{-\operatorname*{Col}}=\sum_{\substack{u\in S_{n}%
;\\u\widetilde{Q}=\widetilde{P}}}\left(  -1\right)  ^{u}u\in\mathcal{A}
\label{eq.def.bas.mur.mur.Col}%
\end{equation}
(recall that $\widetilde{T}$ denotes the $n$-column-tabloid corresponding to
the $n$-tableau $T$, as defined in Definition \ref{def.coltabloid.coltabloid}).
\end{definition}

Using the notations $\operatorname*{Row}\left(  i,T\right)  $ and
$\operatorname*{Col}\left(  j,T\right)  $ from Proposition
\ref{prop.symmetrizers.int}, we can rewrite the definitions of $\nabla
_{P,Q}^{\operatorname*{Row}}$ and $\nabla_{P,Q}^{-\operatorname*{Col}}$ in the
following more concrete way:

\begin{definition}
\label{def.bas.mur.Row-Col}For any $n$-tableau $T$ and any $i\in\mathbb{Z}$,
we let $\operatorname*{Row}\left(  i,T\right)  $ denote the set of all entries
in the $i$-th row of $T$.

For any $n$-tableau $T$ and any $j\in\mathbb{Z}$, we let $\operatorname*{Col}%
\left(  j,T\right)  $ denote the set of all entries in the $j$-th column of
$T$.
\end{definition}

\begin{proposition}
\label{prop.bas.mur.equi-Row-Col}Let $P$ and $Q$ be two $n$-tableaux of the
same shape. Then,
\begin{equation}
\nabla_{P,Q}^{\operatorname*{Row}}=\sum_{\substack{u\in S_{n};\\u\left(
\operatorname*{Row}\left(  i,Q\right)  \right)  =\operatorname*{Row}\left(
i,P\right)  \text{ for all }i\in\mathbb{Z}}}u
\label{eq.rmk.bas.mur.equi-Row-Col.Row}%
\end{equation}
and%
\begin{equation}
\nabla_{P,Q}^{-\operatorname*{Col}}=\sum_{\substack{u\in S_{n};\\u\left(
\operatorname*{Col}\left(  j,Q\right)  \right)  =\operatorname*{Col}\left(
j,P\right)  \text{ for all }j\in\mathbb{Z}}}\left(  -1\right)  ^{u}u.
\label{eq.rmk.bas.mur.equi-Row-Col.Col}%
\end{equation}

\end{proposition}

\begin{example}
Let $P=\ytableaushort{12,34}$ and $Q=\ytableaushort{14,23}$ . Then, we can use
Proposition \ref{prop.bas.mur.equi-Row-Col} to compute $\nabla_{P,Q}%
^{\operatorname*{Row}}$ and $\nabla_{P,Q}^{-\operatorname*{Col}}$ as follows:
Let us write $\operatorname*{oln}\left(  i_{1}i_{2}\ldots i_{n}\right)  $ for
the permutation of $\left[  n\right]  $ that sends $1,2,\ldots,n$ to
$i_{1},i_{2},\ldots,i_{n}$, respectively. Then, $\operatorname*{Row}\left(
1,Q\right)  =\left\{  1,4\right\}  $ and $\operatorname*{Col}\left(
1,Q\right)  =\left\{  1,2\right\}  $, and so on. Thus,
(\ref{eq.rmk.bas.mur.equi-Row-Col.Row}) yields%
\[
\nabla_{P,Q}^{\operatorname*{Row}}=\sum_{\substack{u\in S_{n};\\u\left(
\left\{  1,4\right\}  \right)  =\left\{  1,2\right\}  ;\\u\left(  \left\{
2,3\right\}  \right)  =\left\{  3,4\right\}  }}u=\operatorname*{oln}\left(
1342\right)  +\operatorname*{oln}\left(  1432\right)  +\operatorname*{oln}%
\left(  2341\right)  +\operatorname*{oln}\left(  2431\right)  ,
\]
whereas (\ref{eq.rmk.bas.mur.equi-Row-Col.Col}) yields
\begin{align*}
\nabla_{P,Q}^{-\operatorname*{Col}}  &  =\sum_{\substack{u\in S_{n};\\u\left(
\left\{  1,2\right\}  \right)  =\left\{  1,3\right\}  ;\\u\left(  \left\{
4,3\right\}  \right)  =\left\{  2,4\right\}  }}\left(  -1\right)  ^{u}u\\
&  =-\operatorname*{oln}\left(  1324\right)  +\operatorname*{oln}\left(
1342\right)  +\operatorname*{oln}\left(  3124\right)  -\operatorname*{oln}%
\left(  3142\right)  .
\end{align*}

\end{example}

To prove Proposition \ref{prop.bas.mur.equi-Row-Col}, we need a simple
combinatorial lemma:

\begin{lemma}
\label{lem.bas.mur.row-equ-as-Row-Row}Let $U$ and $V$ be two $n$-tableaux of
the same shape. Then: \medskip

\textbf{(a)} We have the logical equivalence%
\[
\left(  \overline{U}=\overline{V}\right)  \ \Longleftrightarrow\ \left(
\operatorname*{Row}\left(  i,U\right)  =\operatorname*{Row}\left(  i,V\right)
\text{ for all }i\in\mathbb{Z}\right)  .
\]

\textbf{(b)} We have the logical equivalence%
\[
\left(  \widetilde{U}=\widetilde{V}\right)  \ \Longleftrightarrow\ \left(
\operatorname*{Col}\left(  j,U\right)  =\operatorname*{Col}\left(  j,V\right)
\text{ for all }j\in\mathbb{Z}\right)  .
\]

\end{lemma}

The proof of Lemma \ref{lem.bas.mur.row-equ-as-Row-Row} should be completely
routine by now; it can be found in the appendix (Section
\ref{sec.details.bas.mur}).

\begin{proof}
[Proof of Proposition \ref{prop.bas.mur.equi-Row-Col}.]The core of the proof
are the following two claims:

\begin{statement}
\textit{Claim 1:} Let $u\in S_{n}$. Then, we have the logical equivalence
\[
\left(  u\overline{Q}=\overline{P}\right)  \ \Longleftrightarrow\ \left(
u\left(  \operatorname*{Row}\left(  i,Q\right)  \right)  =\operatorname*{Row}%
\left(  i,P\right)  \text{ for all }i\in\mathbb{Z}\right)  .
\]

\end{statement}

\begin{statement}
\textit{Claim 2:} Let $u\in S_{n}$. Then, we have the logical equivalence
\[
\left(  u\widetilde{Q}=\widetilde{P}\right)  \ \Longleftrightarrow\ \left(
u\left(  \operatorname*{Col}\left(  j,Q\right)  \right)  =\operatorname*{Col}%
\left(  j,P\right)  \text{ for all }j\in\mathbb{Z}\right)  .
\]

\end{statement}

\begin{fineprint}
\begin{proof}
[Proof of Claim 1.]Definition \ref{def.tabloid.Sn-act} yields
$u\rightharpoonup\overline{Q}=\overline{u\rightharpoonup Q}$.

Let $i\in\mathbb{Z}$. By Definition \ref{def.tableau.Sn-act}, each entry of
the tableau $u\rightharpoonup Q$ is obtained from the corresponding entry of
$Q$ by applying the permutation $u$. Hence, in particular, the set of all
entries in the $i$-th row of $u\rightharpoonup Q$ is obtained from the set of
all entries in the $i$-th row of $Q$ by applying the permutation $u$ to each
element. In other words,
\begin{equation}
\operatorname*{Row}\left(  i,u\rightharpoonup Q\right)  =u\left(
\operatorname*{Row}\left(  i,Q\right)  \right)
\label{pf.prop.bas.mur.equi-Row-Col.a1}%
\end{equation}
(since the former set is $\operatorname*{Row}\left(  i,u\rightharpoonup
Q\right)  $ whereas the latter set is $\operatorname*{Row}\left(  i,Q\right)
$).

Forget that we fixed $i$. We thus have proved the equality
(\ref{pf.prop.bas.mur.equi-Row-Col.a1}) for each $i\in\mathbb{Z}$.

Now, we have the following chain of logical equivalences:%
\begin{align*}
\left(  u\overline{Q}=\overline{P}\right)  \  &  \Longleftrightarrow\ \left(
\overline{u\rightharpoonup Q}=\overline{P}\right)  \ \ \ \ \ \ \ \ \ \ \left(
\text{since }u\overline{Q}=u\rightharpoonup\overline{Q}=\overline
{u\rightharpoonup Q}\right) \\
&  \Longleftrightarrow\ \left(  \operatorname*{Row}\left(  i,u\rightharpoonup
Q\right)  =\operatorname*{Row}\left(  i,P\right)  \text{ for all }%
i\in\mathbb{Z}\right) \\
&  \ \ \ \ \ \ \ \ \ \ \ \ \ \ \ \ \ \ \ \ \left(  \text{by Lemma
\ref{lem.bas.mur.row-equ-as-Row-Row} \textbf{(a)}, applied to }%
U=u\rightharpoonup Q\text{ and }V=P\right) \\
&  \Longleftrightarrow\ \left(  u\left(  \operatorname*{Row}\left(
i,Q\right)  \right)  =\operatorname*{Row}\left(  i,P\right)  \text{ for all
}i\in\mathbb{Z}\right)
\end{align*}
(by (\ref{pf.prop.bas.mur.equi-Row-Col.a1})). This proves Claim 1.
\end{proof}

\begin{proof}
[Proof of Claim 2.]This is analogous to Claim 1 (using Lemma
\ref{lem.bas.mur.row-equ-as-Row-Row} \textbf{(b)} instead of Lemma
\ref{lem.bas.mur.row-equ-as-Row-Row} \textbf{(a)}).
\end{proof}
\end{fineprint}

Now, (\ref{eq.def.bas.mur.mur.Row}) becomes%
\[
\nabla_{P,Q}^{\operatorname*{Row}}=\sum_{\substack{u\in S_{n};\\u\overline
{Q}=\overline{P}}}u=\sum_{\substack{u\in S_{n};\\u\left(  \operatorname*{Row}%
\left(  i,Q\right)  \right)  =\operatorname*{Row}\left(  i,P\right)  \text{
for all }i\in\mathbb{Z}}}u
\]
(here, we have replaced the condition \textquotedblleft$u\overline
{Q}=\overline{P}$\textquotedblright\ under the summation sign by the condition
\textquotedblleft$u\left(  \operatorname*{Row}\left(  i,Q\right)  \right)
=\operatorname*{Row}\left(  i,P\right)  $ for all $i\in\mathbb{Z}%
$\textquotedblright, which is equivalent according to Claim 1). Thus,
(\ref{eq.rmk.bas.mur.equi-Row-Col.Row}) is proved. Likewise, we can obtain
(\ref{eq.rmk.bas.mur.equi-Row-Col.Col}) from (\ref{eq.def.bas.mur.mur.Col})
using Claim 2. This proves Proposition \ref{prop.bas.mur.equi-Row-Col}.
\end{proof}

\begin{remark}
\label{rmk.bas.mur.Nab=Nab}We have encountered the row-to-row sums
$\nabla_{P,Q}^{\operatorname*{Row}}$ before: Let $P$ and $Q$ be two
$n$-tableaux of the same shape. Then, in Definition
\ref{def.row-to-row.row-to-row}, we have defined $\nabla_{\overline
{P},\overline{Q}}=\sum_{\substack{u\in S_{n};\\u\overline{Q}=\overline{P}}}u$.
But in Definition \ref{def.bas.mur.mur}, we have defined $\nabla
_{P,Q}^{\operatorname*{Row}}=\sum_{\substack{u\in S_{n};\\u\overline
{Q}=\overline{P}}}u$. Comparing these two equalities, we obtain%
\begin{equation}
\nabla_{\overline{P},\overline{Q}}=\nabla_{P,Q}^{\operatorname*{Row}}.
\label{eq.rmk.bas.mur.Nab=Nab.main}%
\end{equation}

\end{remark}

\subsubsection{\label{subsec.bas.mur.basic}Basic identities}

We shall continue with a series of simple but useful properties of the
$\nabla_{P,Q}^{\operatorname*{Row}}$ and $\nabla_{P,Q}^{-\operatorname*{Col}}$
elements. We begin by identifying the row-symmetrizers $\nabla
_{\operatorname*{Row}T}$ and the column-antisymmetrizers $\nabla
_{\operatorname*{Col}T}^{-}$ as particular cases of row-to-row sums and signed
column-to-column sums, respectively:

\begin{proposition}
\label{prop.bas.mur.rsym}Let $T$ be an $n$-tableau of any shape. Then,%
\begin{align}
\nabla_{T,T}^{\operatorname*{Row}}  &  =\nabla_{\operatorname*{Row}%
T}\ \ \ \ \ \ \ \ \ \ \text{and}\label{eq.prop.bas.mur.rsym.Row}\\
\nabla_{T,T}^{-\operatorname*{Col}}  &  =\nabla_{\operatorname*{Col}T}^{-}.
\label{eq.prop.bas.mur.rsym.Col}%
\end{align}

\end{proposition}

\begin{proof}
The equality (\ref{eq.prop.bas.mur.rsym.Row}) is just Proposition
\ref{prop.row-to-row.rsym} (since (\ref{eq.rmk.bas.mur.Nab=Nab.main}) tells us
that $\nabla_{\overline{T},\overline{T}}$ is precisely our row-to-row sum
$\nabla_{T,T}^{\operatorname*{Row}}$). Thus, let us prove
(\ref{eq.prop.bas.mur.rsym.Col}).

Definition \ref{def.symmetrizers.symmetrizers} yields $\nabla
_{\operatorname*{Col}T}^{-}=\sum_{w\in\mathcal{C}\left(  T\right)  }\left(
-1\right)  ^{w}w$. Recall from Proposition \ref{prop.tableau.Sn-act.1}
\textbf{(a)} that $\mathcal{R}\left(  T\right)  =\left\{  w\in S_{n}%
\ \mid\ \overline{w\rightharpoonup T}=\overline{T}\right\}  $. The same
argument (but applied to columns instead of rows, and to $n$-column-tabloids
instead of $n$-tabloids) shows that%
\[
\mathcal{C}\left(  T\right)  =\left\{  w\in S_{n}\ \mid
\ \widetilde{w\rightharpoonup T}=\widetilde{T}\right\}  .
\]
In other words,%
\[
\mathcal{C}\left(  T\right)  =\left\{  w\in S_{n}\ \mid\ w\widetilde{T}%
=\widetilde{T}\right\}
\]
(since $\widetilde{w\rightharpoonup T}=w\rightharpoonup\widetilde{T}%
=w\widetilde{T}$ for each $w\in S_{n}$). Thus, the summation sign $\sum
_{w\in\mathcal{C}\left(  T\right)  }$ can be rewritten as $\sum
_{\substack{w\in S_{n};\\w\widetilde{T}=\widetilde{T}}}$. Hence,%
\[
\sum_{w\in\mathcal{C}\left(  T\right)  }\left(  -1\right)  ^{w}w=\sum
_{\substack{w\in S_{n};\\w\widetilde{T}=\widetilde{T}}}\left(  -1\right)
^{w}w=\sum_{\substack{u\in S_{n};\\u\widetilde{T}=\widetilde{T}}}\left(
-1\right)  ^{u}u
\]
(here, we have renamed the summation index $w$ as $u$). Therefore,%
\[
\nabla_{\operatorname*{Col}T}^{-}=\sum_{w\in\mathcal{C}\left(  T\right)
}\left(  -1\right)  ^{w}w=\sum_{\substack{u\in S_{n};\\u\widetilde{T}%
=\widetilde{T}}}\left(  -1\right)  ^{u}u=\nabla_{T,T}^{-\operatorname*{Col}}%
\]
(since $\nabla_{T,T}^{-\operatorname*{Col}}$ was defined as $\sum
_{\substack{u\in S_{n};\\u\widetilde{T}=\widetilde{T}}}\left(  -1\right)
^{u}u$). This proves (\ref{eq.prop.bas.mur.rsym.Col}). The proof of
(\ref{eq.prop.bas.mur.rsym.Row}) is analogous (but, as we said, it was already
done when we proved Proposition \ref{prop.row-to-row.rsym}). Thus, Proposition
\ref{prop.bas.mur.rsym} is proved.
\end{proof}

Next, we prove an equivariance property:

\begin{proposition}
\label{prop.bas.mur.Sn-act}Let $P$ and $Q$ be two $n$-tableaux of the same
shape. Let $g\in S_{n}$ and $h\in S_{n}$ be two permutations. Then,%
\begin{align}
\nabla_{gP,hQ}^{\operatorname*{Row}}  &  =g\nabla_{P,Q}^{\operatorname*{Row}%
}h^{-1}\ \ \ \ \ \ \ \ \ \ \text{and}\label{eq.prop.bas.mur.Sn-act.Row}\\
\nabla_{gP,hQ}^{-\operatorname*{Col}}  &  =\left(  -1\right)  ^{g}\left(
-1\right)  ^{h}g\nabla_{P,Q}^{-\operatorname*{Col}}h^{-1}.
\label{eq.prop.bas.mur.Sn-act.Col}%
\end{align}

\end{proposition}

\begin{proof}
Definition \ref{def.coltabloid.Sn-act} yields $h\rightharpoonup\widetilde{Q}%
=\widetilde{h\rightharpoonup Q}$. In other words, $h\widetilde{Q}%
=\widetilde{hQ}$. Thus, $\widetilde{hQ}=h\widetilde{Q}$. Similarly,
$\widetilde{gP}=g\widetilde{P}$.

For any permutation $u\in S_{n}$, we have the following chain of logical
equivalences:%
\begin{align}
&  \ \left(  \left(  guh^{-1}\right)  \widetilde{hQ}=\widetilde{gP}\right)
\nonumber\\
&  \Longleftrightarrow\ \left(  \left(  guh^{-1}\right)  h\widetilde{Q}%
=g\widetilde{P}\right)  \ \ \ \ \ \ \ \ \ \ \left(  \text{since }%
\widetilde{hQ}=h\widetilde{Q}\text{ and }\widetilde{gP}=g\widetilde{P}\right)
\nonumber\\
&  \Longleftrightarrow\ \left(  gu\widetilde{Q}=g\widetilde{P}\right)
\ \ \ \ \ \ \ \ \ \ \left(  \text{since }\left(  guh^{-1}\right)
h\widetilde{Q}=gu\underbrace{h^{-1}h}_{=1}\widetilde{Q}=gu\widetilde{Q}\right)
\nonumber\\
&  \Longleftrightarrow\ \left(  u\widetilde{Q}=\widetilde{P}\right)
\label{pf.prop.bas.mur.Sn-act.eq}%
\end{align}
(here, we have cancelled the \textquotedblleft factor\textquotedblright\ $g$,
since the action of $g$ on the set of all $n$-column-tabloids is a permutation
and thus injective). Now, the definition of $\nabla_{gP,hQ}%
^{-\operatorname*{Col}}$ yields%
\begin{align}
\nabla_{gP,hQ}^{-\operatorname*{Col}}  &  =\sum_{\substack{u\in S_{n}%
;\\u\widetilde{hQ}=\widetilde{gP}}}\left(  -1\right)  ^{u}u\nonumber\\
&  =\sum_{\substack{u\in S_{n};\\\left(  guh^{-1}\right)  \widetilde{hQ}%
=\widetilde{gP}}}\left(  -1\right)  ^{guh^{-1}}guh^{-1}
\label{pf.prop.bas.mur.Sn-act.1}%
\end{align}
(here, we have substituted $guh^{-1}$ for $u$ in the sum, since the map
$S_{n}\rightarrow S_{n},\ u\mapsto guh^{-1}$ is a bijection (its inverse is
the map $S_{n}\rightarrow S_{n},\ v\mapsto g^{-1}vh$)).

On the other hand, $\left(  -1\right)  ^{h}\in\left\{  1,-1\right\}  $, so
that $\left(  \left(  -1\right)  ^{h}\right)  ^{-1}=\left(  -1\right)  ^{h}$
(since each $x\in\left\{  1,-1\right\}  $ satisfies $x^{-1}=x$). Next, recall
that the sign homomorphism%
\begin{align*}
S_{n}  &  \rightarrow\left\{  1,-1\right\}  ,\\
\sigma &  \mapsto\left(  -1\right)  ^{\sigma}%
\end{align*}
is a group morphism. Thus, each $u\in S_{n}$ satisfies%
\[
\left(  -1\right)  ^{guh^{-1}}=\left(  -1\right)  ^{g}\left(  -1\right)
^{u}\underbrace{\left(  \left(  -1\right)  ^{h}\right)  ^{-1}}_{=\left(
-1\right)  ^{h}}=\left(  -1\right)  ^{g}\left(  -1\right)  ^{u}\left(
-1\right)  ^{h}.
\]
Hence, we can rewrite (\ref{pf.prop.bas.mur.Sn-act.1}) as%
\begin{align*}
\nabla_{gP,hQ}^{-\operatorname*{Col}}  &  =\sum_{\substack{u\in S_{n}%
;\\\left(  guh^{-1}\right)  \widetilde{hQ}=\widetilde{gP}}}\left(  -1\right)
^{g}\left(  -1\right)  ^{u}\left(  -1\right)  ^{h}guh^{-1}\\
&  =\sum_{\substack{u\in S_{n};\\u\widetilde{Q}=\widetilde{P}}}\left(
-1\right)  ^{g}\left(  -1\right)  ^{u}\left(  -1\right)  ^{h}guh^{-1}\\
&  \ \ \ \ \ \ \ \ \ \ \ \ \ \ \ \ \ \ \ \ \left(
\begin{array}
[c]{c}%
\text{here, we have replaced the condition \textquotedblleft}\left(
guh^{-1}\right)  \widetilde{hQ}=\widetilde{gP}\text{\textquotedblright}\\
\text{under the summation sign by \textquotedblleft}u\widetilde{Q}%
=\widetilde{P}\text{\textquotedblright,}\\
\text{which is equivalent because of (\ref{pf.prop.bas.mur.Sn-act.eq})}%
\end{array}
\right) \\
&  =\left(  -1\right)  ^{g}\left(  -1\right)  ^{h}g\underbrace{\left(
\sum_{\substack{u\in S_{n};\\u\widetilde{Q}=\widetilde{P}}}\left(  -1\right)
^{u}u\right)  }_{\substack{=\nabla_{P,Q}^{-\operatorname*{Col}}\\\text{(by
(\ref{eq.def.bas.mur.mur.Col}))}}}h^{-1}=\left(  -1\right)  ^{g}\left(
-1\right)  ^{h}g\nabla_{P,Q}^{-\operatorname*{Col}}h^{-1}.
\end{align*}
This proves (\ref{eq.prop.bas.mur.Sn-act.Col}). A similar argument (but easier
due to absence of signs) proves (\ref{eq.prop.bas.mur.Sn-act.Row}). Hence,
Proposition \ref{prop.bas.mur.Sn-act} is established.
\end{proof}

The row-to-row sums and the column-to-column signed sums are actually just
shifted versions of row symmetrizers and the column antisymmetrizers:

\begin{lemma}
[permutation lemma]\label{lem.bas.mur.perm}Let $P$ and $Q$ be two $n$-tableaux
of the same shape. Let $u\in S_{n}$ be a permutation such that $P=uQ$. Then,%
\begin{equation}
\nabla_{P,Q}^{\operatorname*{Row}}=\nabla_{\operatorname*{Row}P}%
u=u\nabla_{\operatorname*{Row}Q} \label{eq.lem.bas.mur.perm.R}%
\end{equation}
and%
\begin{equation}
\nabla_{P,Q}^{-\operatorname*{Col}}=\left(  -1\right)  ^{u}\nabla
_{\operatorname*{Col}P}^{-}u=\left(  -1\right)  ^{u}u\nabla
_{\operatorname*{Col}Q}^{-}. \label{eq.lem.bas.mur.perm.C}%
\end{equation}

\end{lemma}

\begin{proof}
We observe that $\left(  -1\right)  ^{u}\in\left\{  1,-1\right\}  $ and thus
$\left(  \left(  -1\right)  ^{u}\right)  ^{-1}=\left(  -1\right)  ^{u}$ (since
$x^{-1}=x$ for each $x\in\left\{  1,-1\right\}  $).

Applying (\ref{eq.prop.bas.mur.Sn-act.Col}) to $g=\operatorname*{id}$ and
$h=u$, we obtain%
\[
\nabla_{\operatorname*{id}P,uQ}^{-\operatorname*{Col}}=\underbrace{\left(
-1\right)  ^{\operatorname*{id}}}_{=1}\left(  -1\right)  ^{u}%
\underbrace{\operatorname*{id}}_{=1}\nabla_{P,Q}^{-\operatorname*{Col}}%
u^{-1}=\left(  -1\right)  ^{u}\nabla_{P,Q}^{-\operatorname*{Col}}u^{-1}.
\]
In view of $\operatorname*{id}P=P$ and $uQ=P$, we can rewrite this as
$\nabla_{P,P}^{-\operatorname*{Col}}=\left(  -1\right)  ^{u}\nabla
_{P,Q}^{-\operatorname*{Col}}u^{-1}$. Comparing this with $\nabla
_{P,P}^{-\operatorname*{Col}}=\nabla_{\operatorname*{Col}P}^{-}$ (which
follows from (\ref{eq.prop.bas.mur.rsym.Col}), applied to $T=P$), we obtain
$\left(  -1\right)  ^{u}\nabla_{P,Q}^{-\operatorname*{Col}}u^{-1}%
=\nabla_{\operatorname*{Col}P}^{-}$. Solving this for $\nabla_{P,Q}%
^{-\operatorname*{Col}}$, we find%
\begin{equation}
\nabla_{P,Q}^{-\operatorname*{Col}}=\underbrace{\dfrac{1}{\left(  -1\right)
^{u}}}_{=\left(  \left(  -1\right)  ^{u}\right)  ^{-1}=\left(  -1\right)
^{u}}\nabla_{\operatorname*{Col}P}^{-}u=\left(  -1\right)  ^{u}\nabla
_{\operatorname*{Col}P}^{-}u. \label{pf.lem.bas.mur.perm.1}%
\end{equation}
This proves the first equality sign in (\ref{eq.lem.bas.mur.perm.C}).

On the other hand, from $P=uQ$, we obtain $Q=u^{-1}P$. But
(\ref{eq.prop.bas.mur.Sn-act.Col}) (applied to $g=u^{-1}$ and
$h=\operatorname*{id}$) yields%
\[
\nabla_{u^{-1}P,\operatorname*{id}Q}^{-\operatorname*{Col}}%
=\underbrace{\left(  -1\right)  ^{u^{-1}}}_{=\left(  \left(  -1\right)
^{u}\right)  ^{-1}}\underbrace{\left(  -1\right)  ^{\operatorname*{id}}}%
_{=1}u^{-1}\nabla_{P,Q}^{-\operatorname*{Col}}\underbrace{\operatorname*{id}%
\nolimits^{-1}}_{=\operatorname*{id}=1}=\left(  \left(  -1\right)
^{u}\right)  ^{-1}u^{-1}\nabla_{P,Q}^{-\operatorname*{Col}}.
\]
In view of $u^{-1}P=Q$ and $\operatorname*{id}Q=Q$, we can rewrite this
equality as $\nabla_{Q,Q}^{-\operatorname*{Col}}=\left(  \left(  -1\right)
^{u}\right)  ^{-1}u^{-1}\nabla_{P,Q}^{-\operatorname*{Col}}$. Comparing this
with $\nabla_{Q,Q}^{-\operatorname*{Col}}=\nabla_{\operatorname*{Col}Q}^{-}$
(which follows from (\ref{eq.prop.bas.mur.rsym.Col}), applied to $T=Q$), we
obtain $\left(  \left(  -1\right)  ^{u}\right)  ^{-1}u^{-1}\nabla
_{P,Q}^{-\operatorname*{Col}}=\nabla_{\operatorname*{Col}Q}^{-}$. Solving this
for $\nabla_{P,Q}^{-\operatorname*{Col}}$, we find%
\begin{equation}
\nabla_{P,Q}^{-\operatorname*{Col}}=\left(  -1\right)  ^{u}u\nabla
_{\operatorname*{Col}Q}^{-}. \label{pf.lem.bas.mur.perm.2}%
\end{equation}
Combining this with (\ref{pf.lem.bas.mur.perm.1}), we obtain the equality
(\ref{eq.lem.bas.mur.perm.C}).

An analogous argument (but without any signs, and using
(\ref{eq.prop.bas.mur.rsym.Row}) and (\ref{eq.prop.bas.mur.Sn-act.Row}))
yields (\ref{eq.lem.bas.mur.perm.R}). Thus, Lemma \ref{lem.bas.mur.perm} is proved.
\end{proof}

We can easily describe the actions of $S$ and $T_{\operatorname*{sign}}$ on
the $\nabla_{P,Q}^{\operatorname*{Row}}$ and $\nabla_{P,Q}%
^{-\operatorname*{Col}}$ elements:

\begin{proposition}
\label{prop.bas.mur.S}Let $P$ and $Q$ be two $n$-tableaux of the same shape.
Then,%
\begin{align}
S\left(  \nabla_{P,Q}^{\operatorname*{Row}}\right)   &  =\nabla_{Q,P}%
^{\operatorname*{Row}}\ \ \ \ \ \ \ \ \ \ \text{and}%
\label{eq.prop.bas.mur.S.Row}\\
S\left(  \nabla_{P,Q}^{-\operatorname*{Col}}\right)   &  =\nabla
_{Q,P}^{-\operatorname*{Col}}. \label{eq.prop.bas.mur.S.Col}%
\end{align}

\end{proposition}

\begin{proof}
Let $D$ be the shape of the $n$-tableaux $P$ and $Q$. Let $g$ be the
permutation $w_{P,Q}\in S_{n}$ as in Definition \ref{def.specht.wPQ}. Thus,
(\ref{eq.def.specht.wPQ.eq}) says that $w_{P,Q}Q=P$. In other words, $gQ=P$
(since we defined $g$ to be $w_{P,Q}$), so that $P=gQ$. Note also that
$\left(  -1\right)  ^{g^{-1}}=\left(  -1\right)  ^{g}$ (since any permutation
$\sigma\in S_{n}$ satisfies $\left(  -1\right)  ^{\sigma^{-1}}=\left(
-1\right)  ^{\sigma}$).

Proposition \ref{prop.symmetrizers.antipode} shows that every $n$-tableau $T$
of shape $D$ satisfies $S\left(  \nabla_{\operatorname*{Col}T}^{-}\right)
=\nabla_{\operatorname*{Col}T}^{-}$. Applying this to $T=P$, we obtain
$S\left(  \nabla_{\operatorname*{Col}P}^{-}\right)  =\nabla
_{\operatorname*{Col}P}^{-}$.

Now, $P=gQ$; hence, (\ref{eq.lem.bas.mur.perm.C}) (applied to $u=g$) yields
\[
\nabla_{P,Q}^{-\operatorname*{Col}}=\left(  -1\right)  ^{g}\nabla
_{\operatorname*{Col}P}^{-}g=\left(  -1\right)  ^{g}g\nabla
_{\operatorname*{Col}Q}^{-}.
\]
Hence, in particular, $\nabla_{P,Q}^{-\operatorname*{Col}}=\left(  -1\right)
^{g}\nabla_{\operatorname*{Col}P}^{-}g$. Applying the map $S$ to this
equality, we find%
\[
S\left(  \nabla_{P,Q}^{-\operatorname*{Col}}\right)  =S\left(  \left(
-1\right)  ^{g}\nabla_{\operatorname*{Col}P}^{-}g\right)  =\left(  -1\right)
^{g}\cdot S\left(  g\right)  \cdot S\left(  \nabla_{\operatorname*{Col}P}%
^{-}\right)
\]
(since $S$ is a $\mathbf{k}$-algebra anti-automorphism (by Theorem
\ref{thm.S.auto} \textbf{(a)})). Thus,%
\begin{align}
S\left(  \nabla_{P,Q}^{-\operatorname*{Col}}\right)   &  =\left(  -1\right)
^{g}\cdot\underbrace{S\left(  g\right)  }_{\substack{=g^{-1}\\\text{(by the
definition of }S\text{)}}}\cdot\underbrace{S\left(  \nabla
_{\operatorname*{Col}P}^{-}\right)  }_{=\nabla_{\operatorname*{Col}P}^{-}%
}\nonumber\\
&  =\left(  -1\right)  ^{g}g^{-1}\nabla_{\operatorname*{Col}P}^{-}.
\label{pf.prop.bas.mur.S.4}%
\end{align}

On the other hand, $Q=g^{-1}P$ (since $P=gQ$); thus,
(\ref{eq.lem.bas.mur.perm.R}) (applied to $Q$, $P$ and $g^{-1}$ instead of
$P$, $Q$ and $u$) yields%
\[
\nabla_{Q,P}^{-\operatorname*{Col}}=\left(  -1\right)  ^{g^{-1}}%
\nabla_{\operatorname*{Col}Q}^{-}g^{-1}=\left(  -1\right)  ^{g^{-1}}%
g^{-1}\nabla_{\operatorname*{Col}P}^{-}.
\]
Hence,%
\[
\nabla_{Q,P}^{-\operatorname*{Col}}=\underbrace{\left(  -1\right)  ^{g^{-1}}%
}_{=\left(  -1\right)  ^{g}}g^{-1}\nabla_{\operatorname*{Col}P}^{-}=\left(
-1\right)  ^{g}g^{-1}\nabla_{\operatorname*{Col}P}^{-}.
\]
Comparing this with (\ref{pf.prop.bas.mur.S.4}), we obtain $S\left(
\nabla_{P,Q}^{-\operatorname*{Col}}\right)  =\nabla_{Q,P}%
^{-\operatorname*{Col}}$. This proves (\ref{eq.prop.bas.mur.S.Col}). A similar
argument (but without any signs) yields (\ref{eq.prop.bas.mur.S.Row}). Thus,
Proposition \ref{prop.bas.mur.S} is proved.
\end{proof}

\begin{proposition}
\label{prop.bas.mur.Tsign}Let $P$ and $Q$ be two $n$-tableaux of the same
shape. Then,%
\begin{align}
T_{\operatorname*{sign}}\left(  \nabla_{P,Q}^{\operatorname*{Row}}\right)   &
=\nabla_{P\mathbf{r},Q\mathbf{r}}^{-\operatorname*{Col}}%
\ \ \ \ \ \ \ \ \ \ \text{and}\label{eq.prop.bas.mur.Tsign.Row}\\
T_{\operatorname*{sign}}\left(  \nabla_{P,Q}^{-\operatorname*{Col}}\right)
&  =\nabla_{P\mathbf{r},Q\mathbf{r}}^{\operatorname*{Row}}.
\label{eq.prop.bas.mur.Tsign.Col}%
\end{align}
(Recall from Definition \ref{def.tableaux.r} that $T\mathbf{r}$ denotes the
reflection of a tableau $T$ across the main diagonal.)
\end{proposition}

\begin{proof}
Let $D$ be the shape of the $n$-tableaux $P$ and $Q$. Let $g$ be the
permutation $w_{P,Q}\in S_{n}$ as in Definition \ref{def.specht.wPQ}. Thus,
(\ref{eq.def.specht.wPQ.eq}) says that $w_{P,Q}Q=P$. In other words, $gQ=P$
(since we defined $g$ to be $w_{P,Q}$), so that $P=gQ$. Note also that
$\left(  \left(  -1\right)  ^{g}\right)  ^{2}=1$ (since any permutation
$\sigma\in S_{n}$ satisfies $\left(  \left(  -1\right)  ^{\sigma}\right)
^{2}=1$).

\begin{noncompile}
From Proposition \ref{prop.specht.ET.r}, we know that $\nabla
_{\operatorname*{Row}\left(  P\mathbf{r}\right)  }=T_{\operatorname*{sign}%
}\left(  \nabla_{\operatorname*{Col}P}^{-}\right)  $ and $\nabla
_{\operatorname*{Col}\left(  P\mathbf{r}\right)  }^{-}=T_{\operatorname*{sign}%
}\left(  \nabla_{\operatorname*{Row}P}\right)  $.
\end{noncompile}

Now, $P=gQ$; hence, (\ref{eq.lem.bas.mur.perm.C}) (applied to $u=g$) yields
\[
\nabla_{P,Q}^{-\operatorname*{Col}}=\left(  -1\right)  ^{g}\nabla
_{\operatorname*{Col}P}^{-}g=\left(  -1\right)  ^{g}g\nabla
_{\operatorname*{Col}Q}^{-}.
\]
Hence, in particular, $\nabla_{P,Q}^{-\operatorname*{Col}}=\left(  -1\right)
^{g}\nabla_{\operatorname*{Col}P}^{-}g$. Applying the map
$T_{\operatorname*{sign}}$ to this equality, we find%
\begin{align}
T_{\operatorname*{sign}}\left(  \nabla_{P,Q}^{-\operatorname*{Col}}\right)
&  =T_{\operatorname*{sign}}\left(  \left(  -1\right)  ^{g}\nabla
_{\operatorname*{Col}P}^{-}g\right) \nonumber\\
&  =\left(  -1\right)  ^{g}\cdot\underbrace{T_{\operatorname*{sign}}\left(
\nabla_{\operatorname*{Col}P}^{-}\right)  }_{\substack{=\nabla
_{\operatorname*{Row}\left(  P\mathbf{r}\right)  }\\\text{(by
(\ref{eq.prop.specht.ET.r.2}))}}}\cdot\underbrace{T_{\operatorname*{sign}%
}\left(  g\right)  }_{\substack{=\left(  -1\right)  ^{g}g\\\text{(by the
definition of }T_{\operatorname*{sign}}\text{)}}}\nonumber\\
&  \ \ \ \ \ \ \ \ \ \ \ \ \ \ \ \ \ \ \ \ \left(
\begin{array}
[c]{c}%
\text{since }T_{\operatorname*{sign}}\text{ is a }\mathbf{k}\text{-algebra
morphism}\\
\text{(by Theorem \ref{thm.Tsign.auto} \textbf{(a)})}%
\end{array}
\right) \nonumber\\
&  =\left(  -1\right)  ^{g}\cdot\nabla_{\operatorname*{Row}\left(
P\mathbf{r}\right)  }\cdot\left(  -1\right)  ^{g}g=\underbrace{\left(  \left(
-1\right)  ^{g}\right)  ^{2}}_{=1}\cdot\,\nabla_{\operatorname*{Row}\left(
P\mathbf{r}\right)  }g\nonumber\\
&  =\nabla_{\operatorname*{Row}\left(  P\mathbf{r}\right)  }g.
\label{pf.prop.bas.mur.Tsign.4}%
\end{align}

On the other hand, Proposition \ref{prop.tableaux.r} \textbf{(d)} (applied to
$Q$ and $g$ instead of $T$ and $w$) yields $\left(  gQ\right)  \mathbf{r}%
=g\left(  Q\mathbf{r}\right)  $. Since $P=gQ$, we can rewrite this as
$P\mathbf{r}=g\left(  Q\mathbf{r}\right)  $. Of course, this shows that the
tableaux $P\mathbf{r}$ and $Q\mathbf{r}$ have the same shape. Thus,
(\ref{eq.lem.bas.mur.perm.R}) (applied to $P\mathbf{r}$, $Q\mathbf{r}$ and $g$
instead of $P$, $Q$ and $u$) yields%
\[
\nabla_{P\mathbf{r},Q\mathbf{r}}^{\operatorname*{Row}}=\nabla
_{\operatorname*{Row}\left(  P\mathbf{r}\right)  }g=g\nabla
_{\operatorname*{Row}\left(  Q\mathbf{r}\right)  }.
\]
In particular, $\nabla_{P\mathbf{r},Q\mathbf{r}}^{\operatorname*{Row}}%
=\nabla_{\operatorname*{Row}\left(  P\mathbf{r}\right)  }g$. Comparing this
equality with (\ref{pf.prop.bas.mur.Tsign.4}), we obtain
$T_{\operatorname*{sign}}\left(  \nabla_{P,Q}^{-\operatorname*{Col}}\right)
=\nabla_{P\mathbf{r},Q\mathbf{r}}^{\operatorname*{Row}}$. This proves
(\ref{eq.prop.bas.mur.Tsign.Col}).

The equality (\ref{eq.prop.bas.mur.Tsign.Row}) can be proved similarly.
Alternatively, we can derive it from (\ref{eq.prop.bas.mur.Tsign.Col}) easily:
Apply (\ref{eq.prop.bas.mur.Tsign.Col}) to $P\mathbf{r}$ and $Q\mathbf{r}$
instead of $P$ and $Q$. Thus we find%
\[
T_{\operatorname*{sign}}\left(  \nabla_{P\mathbf{r},Q\mathbf{r}}%
^{-\operatorname*{Col}}\right)  =\nabla_{\left(  P\mathbf{r}\right)
\mathbf{r},\left(  Q\mathbf{r}\right)  \mathbf{r}}^{\operatorname*{Row}%
}=\nabla_{P,Q}^{\operatorname*{Row}}%
\]
(since Proposition \ref{prop.tableaux.r} \textbf{(b)} yields $\left(
P\mathbf{r}\right)  \mathbf{r}=P$ and $\left(  Q\mathbf{r}\right)
\mathbf{r}=Q$). Applying the map $T_{\operatorname*{sign}}$ to this equality,
we obtain%
\[
T_{\operatorname*{sign}}\left(  T_{\operatorname*{sign}}\left(  \nabla
_{P\mathbf{r},Q\mathbf{r}}^{-\operatorname*{Col}}\right)  \right)
=T_{\operatorname*{sign}}\left(  \nabla_{P,Q}^{\operatorname*{Row}}\right)  .
\]
Hence,%
\begin{align*}
T_{\operatorname*{sign}}\left(  \nabla_{P,Q}^{\operatorname*{Row}}\right)   &
=T_{\operatorname*{sign}}\left(  T_{\operatorname*{sign}}\left(
\nabla_{P\mathbf{r},Q\mathbf{r}}^{-\operatorname*{Col}}\right)  \right)
=\underbrace{\left(  T_{\operatorname*{sign}}\circ T_{\operatorname*{sign}%
}\right)  }_{\substack{=\operatorname*{id}\\\text{(by Theorem
\ref{thm.Tsign.auto} \textbf{(b)})}}}\left(  \nabla_{P\mathbf{r},Q\mathbf{r}%
}^{-\operatorname*{Col}}\right) \\
&  =\operatorname*{id}\left(  \nabla_{P\mathbf{r},Q\mathbf{r}}%
^{-\operatorname*{Col}}\right)  =\nabla_{P\mathbf{r},Q\mathbf{r}%
}^{-\operatorname*{Col}}.
\end{align*}
This proves (\ref{eq.prop.bas.mur.Tsign.Row}). Thus, the proof of Proposition
\ref{prop.bas.mur.Tsign} is complete.
\end{proof}

For future use, let us prove another identity:

\begin{lemma}
\label{lem.bas.mur.PQQP}Let $P$ and $Q$ be two $n$-tableaux of the same shape.
Then,%
\[
\nabla_{P,Q}^{-\operatorname*{Col}}\nabla_{Q,P}^{\operatorname*{Row}}%
=\pm\mathbf{E}_{P},
\]
where $\mathbf{E}_{P}$ is the Young symmetrizer (as defined in Definition
\ref{def.specht.ET.defs} \textbf{(b)}).
\end{lemma}

\begin{proof}
Let $D$ be the shape of the $n$-tableaux $P$ and $Q$. Then, Proposition
\ref{prop.tableau.Sn-act.1} \textbf{(b)} (applied to $T=Q$) shows that any
$n$-tableau of shape $D$ can be written as $w\rightharpoonup Q$ for some $w\in
S_{n}$. Thus, in particular, $P$ can be written in this way. In other words,
there exists some $w\in S_{n}$ such that $P=w\rightharpoonup Q$. Consider this
$w$.

Thus, $P=w\rightharpoonup Q=wQ$. Hence, (\ref{eq.lem.bas.mur.perm.C}) (applied
to $u=w$) yields%
\[
\nabla_{P,Q}^{-\operatorname*{Col}}=\left(  -1\right)  ^{w}\nabla
_{\operatorname*{Col}P}^{-}w=\left(  -1\right)  ^{w}w\nabla
_{\operatorname*{Col}Q}^{-}.
\]

Moreover, $Q=w^{-1}P$ (since $P=wQ$). Hence, (\ref{eq.lem.bas.mur.perm.R})
(applied to $Q$, $P$ and $w^{-1}$ instead of $P$, $Q$ and $u$) yields
\[
\nabla_{Q,P}^{\operatorname*{Row}}=\nabla_{\operatorname*{Row}Q}w^{-1}%
=w^{-1}\nabla_{\operatorname*{Row}P}.
\]

Hence,%
\begin{align*}
\underbrace{\nabla_{P,Q}^{-\operatorname*{Col}}}_{=\left(  -1\right)
^{w}\nabla_{\operatorname*{Col}P}^{-}w}\ \ \underbrace{\nabla_{Q,P}%
^{\operatorname*{Row}}}_{=w^{-1}\nabla_{\operatorname*{Row}P}}  &  =\left(
-1\right)  ^{w}\nabla_{\operatorname*{Col}P}^{-}\underbrace{ww^{-1}}%
_{=1}\nabla_{\operatorname*{Row}P}\\
&  =\left(  -1\right)  ^{w}\underbrace{\nabla_{\operatorname*{Col}P}^{-}%
\nabla_{\operatorname*{Row}P}}_{\substack{=\mathbf{E}_{P}\\\text{(by the
definition of }\mathbf{E}_{P}\text{)}}}=\underbrace{\left(  -1\right)  ^{w}%
}_{=\pm1}\mathbf{E}_{P}=\pm\mathbf{E}_{P}.
\end{align*}
This proves Lemma \ref{lem.bas.mur.PQQP}.
\end{proof}

\subsubsection{\label{subsec.bas.mur.bas}The basis theorem}

Recall from Definition \ref{def.bas.not} \textbf{(e)} that
$\operatorname*{SBT}\left(  n\right)  $ is the set of all standard
$n$-bitableaux, i.e., of all triples $\left(  \lambda,U,V\right)  $, where
$\lambda$ is a partition of $n$ and where $U$ and $V$ are two standard
$n$-tableaux of shape $\lambda$. We now state one of the main results in this
section (the proof will be given later):

\begin{theorem}
\label{thm.bas.mur.bas}Both families%
\[
\left(  \nabla_{V,U}^{\operatorname*{Row}}\right)  _{\left(  \lambda
,U,V\right)  \in\operatorname*{SBT}\left(  n\right)  }%
\ \ \ \ \ \ \ \ \ \ \text{and}\ \ \ \ \ \ \ \ \ \ \left(  \nabla
_{V,U}^{-\operatorname*{Col}}\right)  _{\left(  \lambda,U,V\right)
\in\operatorname*{SBT}\left(  n\right)  }%
\]
are bases of $\mathcal{A}$.
\end{theorem}

\begin{definition}
\label{def.bas.mur.bas}These two bases are called the \emph{Murphy bases} (or
the \emph{Murphy cellular bases}), specifically the \emph{row Murphy basis}
and the \emph{column Murphy basis}.
\end{definition}

For $n=3$, the row Murphy basis is the basis
$\overrightarrow{\operatorname*{rmu}}_{3}$ from Subsection
\ref{subsec.bas.mur.n=3}, since its basis vectors are%
\begin{align*}
\nabla_{123,\ 123}^{\operatorname*{Row}}  &  =\nabla;\\
\nabla_{12\backslash\backslash3,\ 12\backslash\backslash3}%
^{\operatorname*{Row}}  &  =\nabla_{3,3};\\
\nabla_{12\backslash\backslash3,\ 13\backslash\backslash2}%
^{\operatorname*{Row}}  &  =\nabla_{3,2};\\
\nabla_{13\backslash\backslash2,\ 12\backslash\backslash3}%
^{\operatorname*{Row}}  &  =\nabla_{2,3};\\
\nabla_{13\backslash\backslash2,\ 13\backslash\backslash2}%
^{\operatorname*{Row}}  &  =\nabla_{2,2};\\
\nabla_{1\backslash\backslash2\backslash\backslash3,\ 1\backslash
\backslash2\backslash\backslash3}^{\operatorname*{Row}}  &
=\operatorname*{id}=1.
\end{align*}

One of our main goals in this section is to prove Theorem
\ref{thm.bas.mur.bas}, and then to study the Murphy bases in more detail, in
particular their nice block-triangularity properties.

The Murphy bases were discovered by Canfield and Williamson in 1989
\cite[Theorems 3.12, 3.14 and 6.13]{CanWil89} (where the $\nabla
_{V,U}^{\operatorname*{Row}}$ are called \emph{matched sums}, and the
$\nabla_{V,U}^{-\operatorname*{Col}}$ are called \emph{signed sums}, although
both are indexed by pairs of ordered set partitions rather than of tableaux)
and -- in the more general context of a Hecke algebra -- by Murphy in 1992
\cite[\S 3]{Murphy92}, \cite[Theorem 4.17]{Murphy95} (whose $x_{st}$
generalize our $\nabla_{V,U}^{\operatorname*{Row}}$, although their definition
is somewhat subtler).

One more little remark: The unity $1$ of $\mathcal{A}$ is always an element of
each of the two Murphy bases. Namely:

\begin{exercise}
\fbox{1} \textbf{(a)} Prove that $\nabla_{P,P}^{\operatorname*{Row}}=1$, where
$P$ is any $n$-tableau of shape $\left(  1,1,\ldots,1\right)  $. \medskip

\textbf{(b)} Prove that $\nabla_{Q,Q}^{-\operatorname*{Col}}=1$, where $Q$ is
any $n$-tableau of shape $\left(  n\right)  $.
\end{exercise}

\subsubsection{\label{subsec.bas.mur.dot}The dot product on $\mathcal{A}$}

Recall the notion of a bilinear form (Definition \ref{def.dual.bilform}).

\begin{definition}
\label{def.bas.bilin}Recall that $\left(  w\right)  _{w\in S_{n}}$ is the
standard basis of $\mathbf{k}\left[  S_{n}\right]  =\mathcal{A}$. Thus, by
Proposition \ref{prop.dual.bilf.on-basis}, we can define a bilinear form
$f:\mathcal{A}\times\mathcal{A}\rightarrow\mathbf{k}$ by setting%
\[
f\left(  u,v\right)  =\delta_{u,v}\ \ \ \ \ \ \ \ \ \ \text{for all }u\in
S_{n}\text{ and }v\in S_{n}%
\]
(where $\delta_{u,v}$ denotes the Kronecker delta). We shall denote this form
$f$ by $\left\langle \cdot,\cdot\right\rangle $, which means that we shall
denote each value $f\left(  \mathbf{a},\mathbf{b}\right)  $ of $f$ as
$\left\langle \mathbf{a},\mathbf{b}\right\rangle $. Thus, for all $u\in S_{n}$
and $v\in S_{n}$, we have%
\begin{equation}
\left\langle u,v\right\rangle =f\left(  u,v\right)  =\delta_{u,v}
\label{eq.def.bas.bilin.onbas}%
\end{equation}
(by the definition of $f$).

This form $\left\langle \cdot,\cdot\right\rangle $ will be called the
\emph{dot product} on $\mathcal{A}$ (even though our notation for it is a pair
of angular brackets rather than a dot).
\end{definition}

As the name \textquotedblleft dot product\textquotedblright\ suggests, this
form $\left\langle \cdot,\cdot\right\rangle $ is the obvious analogue of the
dot product on $\mathbf{k}^{n}$ on the free $\mathbf{k}$-module $\mathbf{k}%
\left[  S_{n}\right]  $. Here is another way to say this:

\begin{exercise}
\label{exe.bas.bilin.explicit}\fbox{1} Let $\left(  a_{w}\right)  _{w\in
S_{n}}\in\mathbf{k}^{S_{n}}$ and $\left(  b_{w}\right)  _{w\in S_{n}}%
\in\mathbf{k}^{S_{n}}$ be two families of scalars. Prove that
\[
\left\langle \sum_{w\in S_{n}}a_{w}w,\ \sum_{w\in S_{n}}b_{w}w\right\rangle
=\sum_{w\in S_{n}}a_{w}b_{w}.
\]

\end{exercise}

We shall give several more explicit formulas for $\left\langle \mathbf{a}%
,\mathbf{b}\right\rangle $ in a moment. First, we recall from Definition
\ref{def.specht.ET.wa} that if $\mathbf{a}$ is any element of $\mathcal{A}%
=\mathbf{k}\left[  S_{n}\right]  $ and if $w\in S_{n}$ is any permutation,
then $\left[  w\right]  \mathbf{a}$ denotes the coefficient of $w$ in
$\mathbf{a}$ (when $\mathbf{a}$ is expanded as a linear combination of the
standard basis vectors of $\mathbf{k}\left[  S_{n}\right]  $). We now use this
notation to rewrite the bilinear form $\left\langle \cdot,\cdot\right\rangle $
as follows:

\begin{proposition}
\label{prop.bas.bilin.sym-etc}Let $\operatorname*{id}$ mean the identity
permutation $\operatorname*{id}\nolimits_{\left[  n\right]  }\in S_{n}$. Let
$\mathbf{a},\mathbf{b}\in\mathcal{A}$. Then,%
\begin{align}
\left\langle \mathbf{a},\mathbf{b}\right\rangle  &  =\left[
\operatorname*{id}\right]  \left(  S\left(  \mathbf{a}\right)  \mathbf{b}%
\right) \label{eq.prop.bas.bilin.sym-etc.1}\\
&  =\left[  \operatorname*{id}\right]  \left(  \mathbf{b}S\left(
\mathbf{a}\right)  \right) \label{eq.prop.bas.bilin.sym-etc.2}\\
&  =\left[  \operatorname*{id}\right]  \left(  S\left(  \mathbf{b}\right)
\mathbf{a}\right) \label{eq.prop.bas.bilin.sym-etc.3}\\
&  =\left[  \operatorname*{id}\right]  \left(  \mathbf{a}S\left(
\mathbf{b}\right)  \right) \label{eq.prop.bas.bilin.sym-etc.4}\\
&  =\left\langle \mathbf{b},\mathbf{a}\right\rangle .
\label{eq.prop.bas.bilin.sym-etc.sym}%
\end{align}

\end{proposition}

\begin{proof}
We must prove that the scalars $\left\langle \mathbf{a},\mathbf{b}%
\right\rangle $, $\left[  \operatorname*{id}\right]  \left(  S\left(
\mathbf{a}\right)  \mathbf{b}\right)  $, $\left[  \operatorname*{id}\right]
\left(  \mathbf{b}S\left(  \mathbf{a}\right)  \right)  $, $\left[
\operatorname*{id}\right]  \left(  S\left(  \mathbf{b}\right)  \mathbf{a}%
\right)  $, $\left[  \operatorname*{id}\right]  \left(  \mathbf{a}S\left(
\mathbf{b}\right)  \right)  $ and $\left\langle \mathbf{b},\mathbf{a}%
\right\rangle $ are equal. All these scalars depend $\mathbf{k}$-linearly on
each of $\mathbf{a}$ and $\mathbf{b}$ (since $\left\langle \cdot
,\cdot\right\rangle $ is a bilinear form and since $S$ is a $\mathbf{k}%
$-linear map). Hence, by linearity, we can WLOG assume that both $\mathbf{a}$
and $\mathbf{b}$ belong to the standard basis $\left(  w\right)  _{w\in S_{n}%
}$ of $\mathcal{A}$. Assume this. Thus, $\mathbf{a}=u$ and $\mathbf{b}=v$ for
some $u,v\in S_{n}$. Consider these $u,v$.

From $\mathbf{a}=u$ and $\mathbf{b}=v$, we obtain
\begin{align}
\left\langle \mathbf{a},\mathbf{b}\right\rangle  &  =\left\langle
u,v\right\rangle =\delta_{u,v}\ \ \ \ \ \ \ \ \ \ \left(  \text{by
(\ref{eq.def.bas.bilin.onbas})}\right) \nonumber\\
&  =%
\begin{cases}
1, & \text{if }u=v;\\
0, & \text{if }u\neq v
\end{cases}
\label{pf.prop.bas.bilin.sym-etc.1}%
\end{align}
and $S\left(  \mathbf{a}\right)  =S\left(  u\right)  =u^{-1}$ (by the
definition of $S$) and $S\left(  \mathbf{b}\right)  =S\left(  v\right)
=v^{-1}$ (again by the definition of $S$).

Next, we recall that each permutation $w\in S_{n}$ satisfies
\begin{equation}
\left[  \operatorname*{id}\right]  w=%
\begin{cases}
1, & \text{if }w=\operatorname*{id};\\
0, & \text{if }w\neq\operatorname*{id}%
\end{cases}
\label{pf.prop.bas.bilin.sym-etc.2}%
\end{equation}
(by the very definition of $\left[  \operatorname*{id}\right]  \mathbf{a}$ as
the coefficient of $\operatorname*{id}$ in $\mathbf{a}$). However, from
$S\left(  \mathbf{a}\right)  =u^{-1}$ and $\mathbf{b}=v$, we obtain $S\left(
\mathbf{a}\right)  \mathbf{b}=u^{-1}v$ and thus%
\begin{align*}
\left[  \operatorname*{id}\right]  \left(  S\left(  \mathbf{a}\right)
\mathbf{b}\right)   &  =\left[  \operatorname*{id}\right]  \left(
u^{-1}v\right) \\
&  =%
\begin{cases}
1, & \text{if }u^{-1}v=\operatorname*{id};\\
0, & \text{if }u^{-1}v\neq\operatorname*{id}%
\end{cases}
\ \ \ \ \ \ \ \ \ \ \left(  \text{by (\ref{pf.prop.bas.bilin.sym-etc.2}),
applied to }w=u^{-1}v\right) \\
&  =%
\begin{cases}
1, & \text{if }u=v;\\
0, & \text{if }u\neq v
\end{cases}
\end{align*}
(since the equation $u^{-1}v=\operatorname*{id}$ in $S_{n}$ is equivalent to
$u=v$). Comparing this with (\ref{pf.prop.bas.bilin.sym-etc.1}), we obtain
$\left\langle \mathbf{a},\mathbf{b}\right\rangle =\left[  \operatorname*{id}%
\right]  \left(  S\left(  \mathbf{a}\right)  \mathbf{b}\right)  $. Hence,
(\ref{eq.prop.bas.bilin.sym-etc.1}) is proved.

A similar argument can be used to prove (\ref{eq.prop.bas.bilin.sym-etc.2})
(here, we have to argue that $\mathbf{b}S\left(  \mathbf{a}\right)  =vu^{-1}$,
and that the equation $vu^{-1}=\operatorname*{id}$ in $S_{n}$ is equivalent to
$u=v$). Likewise, we can prove (\ref{eq.prop.bas.bilin.sym-etc.3}) (here, we
must argue that $S\left(  \mathbf{b}\right)  \mathbf{a}=v^{-1}u$, and that the
equation $v^{-1}u=\operatorname*{id}$ in $S_{n}$ is equivalent to $u=v$).
Similarly, we can prove (\ref{eq.prop.bas.bilin.sym-etc.4}) (here, we must
argue that $\mathbf{a}S\left(  \mathbf{b}\right)  =uv^{-1}$, and that the
equation $uv^{-1}=\operatorname*{id}$ in $S_{n}$ is equivalent to $u=v$).

It remains to prove (\ref{eq.prop.bas.bilin.sym-etc.sym}). For this purpose,
let us apply (\ref{eq.prop.bas.bilin.sym-etc.1}) to $\mathbf{b}$ and
$\mathbf{a}$ instead of $\mathbf{a}$ and $\mathbf{b}$ (we can do this, since
we have already proved (\ref{eq.prop.bas.bilin.sym-etc.1})). Thus we obtain
$\left\langle \mathbf{b},\mathbf{a}\right\rangle =\left[  \operatorname*{id}%
\right]  \left(  S\left(  \mathbf{b}\right)  \mathbf{a}\right)  $. Comparing
this with (\ref{eq.prop.bas.bilin.sym-etc.3}), we obtain $\left\langle
\mathbf{a},\mathbf{b}\right\rangle =\left\langle \mathbf{b},\mathbf{a}%
\right\rangle $. Thus, (\ref{eq.prop.bas.bilin.sym-etc.sym}) is proved, and so
the proof of Proposition \ref{prop.bas.bilin.sym-etc} is complete.
\end{proof}

Note that the dot product on $\mathcal{A}$ is precisely the $G$-invariant form
$f$ constructed in the proof of Theorem \ref{thm.dual.perm-rep}, where we take
$G$ to be the symmetric group $S_{n}$ and we take the $G$-set $X$ to be
$G=S_{n}$ itself (with the left regular $G$-action).

In the standard language of bilinear algebra, the equality
(\ref{eq.prop.bas.bilin.sym-etc.sym}) says that the dot product $\left\langle
\cdot,\cdot\right\rangle $ is symmetric, whereas the equality
(\ref{eq.def.bas.bilin.onbas}) says that the standard basis $\left(  w\right)
_{w\in S_{n}}$ of $\mathcal{A}$ is orthonormal with respect to the dot
product. Here are some further properties of the dot product:

\begin{definition}
\label{def.bas.bilin.curries}Recall the concepts of left-curried and
right-curried forms (Definition \ref{def.dual.bilinear}). A $\mathbf{k}%
$-bilinear form $f:U\times V\rightarrow\mathbf{k}$ on two $\mathbf{k}$-modules
$U$ and $V$ is said to be \emph{perfect} if its left-curried form
$f_{L}:U\rightarrow V^{\ast}$ and its right-curried form $f_{R}:V\rightarrow
U^{\ast}$ are $\mathbf{k}$-module isomorphisms. The dot product (i.e., the
bilinear form $\left\langle \cdot,\cdot\right\rangle $ from Definition
\ref{def.bas.bilin}) is perfect; this is easy to see from Exercise
\ref{exe.bas.bilin.explicit}.
\end{definition}

\begin{exercise}
\label{exe.bas.bilin.Sninv}\fbox{1} \textbf{(a)} Prove that the bilinear form
$\left\langle \cdot,\cdot\right\rangle $ is $S_{n}$-invariant, meaning that
all $w\in S_{n}$ and $\mathbf{a},\mathbf{b}\in\mathcal{A}$ satisfy%
\[
\left\langle \mathbf{a},\mathbf{b}\right\rangle =\left\langle w\mathbf{a}%
,w\mathbf{b}\right\rangle .
\]

\textbf{(b)} Prove that the bilinear form $\left\langle \cdot,\cdot
\right\rangle $ is invariant under $S$ and $T_{\operatorname*{sign}}$, meaning
that all $\mathbf{a},\mathbf{b}\in\mathcal{A}$ satisfy%
\[
\left\langle \mathbf{a},\mathbf{b}\right\rangle =\left\langle S\left(
\mathbf{a}\right)  ,S\left(  \mathbf{b}\right)  \right\rangle =\left\langle
T_{\operatorname*{sign}}\left(  \mathbf{a}\right)  ,T_{\operatorname*{sign}%
}\left(  \mathbf{b}\right)  \right\rangle .
\]

\end{exercise}

\subsubsection{Some lemmas from linear algebra}

Next we shall state a few facts from linear algebra that will be used later
on. We recall the notion of a $P\times Q$-matrix (as defined in Section
\ref{sec.spechtmod.detava}).

\begin{noncompile}
If $I$ is a finite set, then we write an $I\times I$-matrix $\left(
a_{i,j}\right)  _{\left(  i,j\right)  \in I\times I}$ as $\left(
a_{i,j}\right)  _{i,j\in I}$.
\end{noncompile}

\begin{lemma}
\label{lem.bas.mur.tria-inv}Let $I$ be a finite totally ordered set. Let
$A=\left(  a_{i,j}\right)  _{i,j\in I}\in\mathbf{k}^{I\times I}$ be an
$I\times I$-matrix. Assume that

\begin{itemize}
\item the scalar $a_{i,i}\in\mathbf{k}$ is invertible for each $i\in I$;

\item we have $a_{i,j}=0$ for all $i,j\in I$ that satisfy $i<j$.
\end{itemize}

\footnotemark\ Then, the matrix $A$ is invertible.
\end{lemma}

\footnotetext{These two assumptions can be reworded as \textquotedblleft the
matrix $A$ is lower-triangular (if its rows and columns are ordered according
to the total order on $I$) and its diagonal entries are
invertible\textquotedblright.}

\begin{proof}
For each $m\in\mathbb{N}$, the set $\left[  m\right]  =\left\{  1,2,\ldots
,m\right\}  $ is totally ordered (with the usual order: $1<2<\cdots<m$).

Any finite totally ordered set is isomorphic to the totally ordered set
$\left[  m\right]  $ for some $m\in\mathbb{N}$. In particular, this must hold
for our finite totally ordered set $I$. Thus, we can WLOG assume that this
totally ordered set $I$ \textbf{is} the totally ordered set $\left[  m\right]
$ for some $m\in\mathbb{N}$\ \ \ \ \footnote{In more detail: Let us label the
elements of $I$ as $t_{1},t_{2},\ldots,t_{m}$ in increasing order. (This can
be done, since $I$ is a finite totally ordered set.) Then, the map $\left[
m\right]  \rightarrow I,\ i\mapsto t_{i}$ is an isomorphism of totally ordered
sets. Hence, the $I\times I$-matrix $A=\left(  a_{i,j}\right)  _{i,j\in I}$
can be transformed into the $m\times m$-matrix $A^{\prime}:=\left(
a_{t_{i},t_{j}}\right)  _{i,j\in\left[  m\right]  }$ by merely changing the
indexing of its rows and columns (and this change does not affect their
ordering). Replacing the matrix $A$ by the matrix $A^{\prime}$ (and the
totally ordered set $I$ by the totally ordered set $\left[  m\right]  $) does
not change the claim of Lemma \ref{lem.bas.mur.tria-inv} (since the matrix $A$
is invertible if and only if its reindexed version $A^{\prime}$ is
invertible), nor does it change the assumptions (for example, the assumption
\textquotedblleft$a_{i,j}=0$ for all $i,j\in I$ that satisfy $i<j$%
\textquotedblright\ is equivalent to \textquotedblleft$a_{t_{i},t_{j}}=0$ for
all $i,j\in\left[  m\right]  $ that satisfy $i<j$\textquotedblright, because
the map $\left[  m\right]  \rightarrow I,\ i\mapsto t_{i}$ is an isomorphism
of totally ordered sets). Hence, we can WLOG assume that $I=\left[  m\right]
$ (otherwise, we just replace $I$ and $A$ by $\left[  m\right]  $ and
$A^{\prime}$).}. Assume this.

\begin{noncompile}
Let us label the elements of $I$ as $t_{1},t_{2},\ldots,t_{m}$ in increasing
order. (This can be done, since $I$ is a finite totally ordered set.) Then,
the map $\left[  m\right]  \rightarrow I,\ i\mapsto t_{i}$ is an isomorphism
of totally ordered sets. Hence, the $I\times I$-matrix $A=\left(
a_{i,j}\right)  _{i,j\in I}$ can be re-encoded as the $m\times m$-matrix
$A^{\prime}:=\left(  a_{t_{i},t_{j}}\right)  _{i,j\in\left[  m\right]  }$, in
the sense that the latter matrix differs from the former only in the indexing
of its rows and columns (and this reindexing does not change the ordering).

In particular, the matrix $A$ is invertible if and only if the matrix
$A^{\prime}$ is invertible. Hence, it suffices to show that the matrix
$A^{\prime}$ is invertible.

Now we claim that the matrix $A^{\prime}$ is lower-triangular.
[\textit{Proof:} Let $p,q\in\left[  m\right]  $ be such that $p<q$. Then,
$t_{p}<t_{q}$ (since we have labelled the elements of $I$ as $t_{1}%
,t_{2},\ldots,t_{m}$ in increasing order). But we assumed that we have
$a_{i,j}=0$ for all $i,j\in I$ that satisfy $i<j$. Applying this to $i=t_{p}$
and $j=t_{q}$, we obtain $a_{t_{p},t_{q}}=0$ (since $t_{p}<t_{q}$). Now forget
that we fixed $p,q$. We thus have shown that $a_{t_{p},t_{q}}=0$ for all
$p,q\in\left[  m\right]  $ satisfying $p<q$. In other words, the matrix
$\left(  a_{t_{i},t_{j}}\right)  _{i,j\in\left[  m\right]  }$ is
upper-triangular. In other words, the matrix $A^{\prime}$ is upper-triangular.]
\end{noncompile}

Now, recall our assumption that we have $a_{i,j}=0$ for all $i,j\in I$ that
satisfy $i<j$. In view of $I=\left[  m\right]  $, this is just saying that
$a_{i,j}=0$ for all $i,j\in\left[  m\right]  $ that satisfy $i<j$. In other
words, the matrix $A$ is lower-triangular. Hence, its determinant $\det A$ is
the product of its diagonal entries (since the determinant of a
lower-triangular matrix is the product of its diagonal entries). In other
words,%
\[
\det A=a_{1,1}a_{2,2}\cdots a_{m,m}=\prod_{i\in\left[  m\right]  }%
a_{i,i}=\prod_{i\in I}a_{i,i}\ \ \ \ \ \ \ \ \ \ \left(  \text{since }\left[
m\right]  =I\right)  .
\]
But we assumed that the scalar $a_{i,i}\in\mathbf{k}$ is invertible for each
$i\in I$. Hence, the product $\prod_{i\in I}a_{i,i}$ of these scalars is also
invertible (since any product of invertible scalars is invertible). In other
words, $\det A$ is invertible (since $\det A=\prod_{i\in I}a_{i,i}$). Thus,
the matrix $A$ is itself invertible (because a square matrix whose determinant
is invertible must itself be invertible\footnote{This is a classical fact
(see, e.g., \cite[Theorem 6.110 \textbf{(a)}]{detnotes}).}). This proves Lemma
\ref{lem.bas.mur.tria-inv}.
\end{proof}

\begin{lemma}
[Gram trick]\label{lem.bilin.gram}Let $n\in\mathbb{N}$. Let $V$ be a free
$\mathbf{k}$-module of rank $n$, and let $\left(  u_{1},u_{2},\ldots
,u_{n}\right)  $ and $\left(  v_{1},v_{2},\ldots,v_{n}\right)  $ be two
$n$-tuples of elements of $V$. Let $f:V\times V\rightarrow\mathbf{k}$ be a
bilinear form. Assume that the $n\times n$-matrix
\[
\left(  f\left(  u_{i},v_{j}\right)  \right)  _{i,j\in\left[  n\right]  }%
\in\mathbf{k}^{n\times n}%
\]
is invertible. (This matrix is called the \emph{Gram matrix} of $f$ with
respect to the $n$-tuples $\left(  u_{1},u_{2},\ldots,u_{n}\right)  $ and
$\left(  v_{1},v_{2},\ldots,v_{n}\right)  $.)

Then, both $n$-tuples $\left(  u_{1},u_{2},\ldots,u_{n}\right)  $ and $\left(
v_{1},v_{2},\ldots,v_{n}\right)  $ are bases of $V$.
\end{lemma}

\begin{proof}
We are given that the $\mathbf{k}$-module $V$ is free of rank $n$. In other
words, $V$ has a basis $\overrightarrow{b}=\left(  b_{1},b_{2},\ldots
,b_{n}\right)  $ consisting of $n$ vectors. Consider this basis.

For each $i\in\left[  n\right]  $, we let $x_{i,1},x_{i,2},\ldots,x_{i,n}$ be
the coordinates of the vector $u_{i}\in V$ with respect to the basis
$\overrightarrow{b}$; that is, we write $u_{i}$ as%
\begin{equation}
u_{i}=\sum_{k\in\left[  n\right]  }x_{i,k}b_{k}%
,\ \ \ \ \ \ \ \ \ \ \text{where }x_{i,k}\in\mathbf{k}.
\label{pf.lem.bilin.gram.ui=}%
\end{equation}

For each $j\in\left[  n\right]  $, we let $y_{j,1},y_{j,2},\ldots,y_{j,n}$ be
the coordinates of the vector $v_{j}\in V$ with respect to the basis
$\overrightarrow{b}$; that is, we write $v_{j}$ as%
\begin{equation}
v_{j}=\sum_{\ell\in\left[  n\right]  }y_{j,\ell}b_{\ell}%
,\ \ \ \ \ \ \ \ \ \ \text{where }y_{j,\ell}\in\mathbf{k}.
\label{pf.lem.bilin.gram.vj=}%
\end{equation}

We let $A_{i,j}$ denote the $\left(  i,j\right)  $-th entry of a matrix $A$.

Let $X$ be the $n\times n$-matrix $\left(  x_{i,k}\right)  _{i,k\in\left[
n\right]  }\in\mathbf{k}^{n\times n}$, and let $Y$ be the $n\times n$-matrix
$\left(  y_{j,\ell}\right)  _{j,\ell\in\left[  n\right]  }\in\mathbf{k}%
^{n\times n}$. Thus, we have%
\begin{equation}
X_{i,k}=x_{i,k}\ \ \ \ \ \ \ \ \ \ \text{for each }i,k\in\left[  n\right]
\label{pf.lem.bilin.gram.X}%
\end{equation}
and%
\begin{equation}
Y_{j,\ell}=y_{j,\ell}\ \ \ \ \ \ \ \ \ \ \text{for each }j,\ell\in\left[
n\right]  . \label{pf.lem.bilin.gram.Y}%
\end{equation}

Let $G$ be the matrix $\left(  f\left(  b_{k},b_{\ell}\right)  \right)
_{k,\ell\in\left[  n\right]  }\in\mathbf{k}^{n\times n}$. Thus,%
\begin{equation}
G_{k,\ell}=f\left(  b_{k},b_{\ell}\right)  \ \ \ \ \ \ \ \ \ \ \text{for each
}k,\ell\in\left[  n\right]  . \label{pf.lem.bilin.gram.G}%
\end{equation}

Now, we claim that%
\begin{equation}
\left(  f\left(  u_{i},v_{j}\right)  \right)  _{i,j\in\left[  n\right]
}=XGY^{T}. \label{pf.lem.bilin.gram.XGYT}%
\end{equation}

\begin{proof}
[Proof of (\ref{pf.lem.bilin.gram.XGYT}).]If $A$ and $B$ are two $n\times
n$-matrices, then the entries of their product $AB$ are given by the formula%
\[
\left(  AB\right)  _{i,j}=\sum_{k\in\left[  n\right]  }A_{i,k}B_{k,j}%
\ \ \ \ \ \ \ \ \ \ \text{for all }i,j\in\left[  n\right]
\]
(by the definition of the product of two matrices). Applying this formula
twice, we obtain the following variant of it: If $A$, $B$ and $C$ are three
$n\times n$-matrices, then the entries of their product $ABC$ are given by the
formula%
\[
\left(  ABC\right)  _{i,j}=\sum_{k\in\left[  n\right]  }\ \ \sum_{\ell
\in\left[  n\right]  }A_{i,k}B_{k,\ell}C_{\ell,j}\ \ \ \ \ \ \ \ \ \ \text{for
all }i,j\in\left[  n\right]  .
\]
Applying this formula to $A=X$ and $B=G$ and $C=Y^{T}$, we obtain
\begin{align*}
\left(  XGY^{T}\right)  _{i,j}  &  =\sum_{k\in\left[  n\right]  }%
\ \ \sum_{\ell\in\left[  n\right]  }X_{i,k}G_{k,\ell}\underbrace{\left(
Y^{T}\right)  _{\ell,j}}_{\substack{=Y_{j,\ell}\\\text{(by the definition
of}\\\text{the transpose of a matrix)}}}\\
&  =\sum_{k\in\left[  n\right]  }\ \ \sum_{\ell\in\left[  n\right]
}\underbrace{X_{i,k}}_{\substack{=x_{i,k}\\\text{(by
(\ref{pf.lem.bilin.gram.X}))}}}\ \ \underbrace{G_{k,\ell}}%
_{\substack{=f\left(  b_{k},b_{\ell}\right)  \\\text{(by
(\ref{pf.lem.bilin.gram.G}))}}}\ \ \underbrace{Y_{j,\ell}}%
_{\substack{=y_{j,\ell}\\\text{(by (\ref{pf.lem.bilin.gram.Y}))}}}\\
&  =\sum_{k\in\left[  n\right]  }\ \ \sum_{\ell\in\left[  n\right]  }%
x_{i,k}f\left(  b_{k},b_{\ell}\right)  y_{j,\ell}\ \ \ \ \ \ \ \ \ \ \text{for
all }i,j\in\left[  n\right]  .
\end{align*}
On the other hand, for all $i,j\in\left[  n\right]  $, we have%
\begin{align*}
f\left(  u_{i},v_{j}\right)   &  =f\left(  \sum_{k\in\left[  n\right]
}x_{i,k}b_{k},\ \ \sum_{\ell\in\left[  n\right]  }y_{j,\ell}b_{\ell}\right)
\ \ \ \ \ \ \ \ \ \ \left(  \text{by (\ref{pf.lem.bilin.gram.ui=}) and
(\ref{pf.lem.bilin.gram.vj=})}\right) \\
&  =\sum_{k\in\left[  n\right]  }\ \ \sum_{\ell\in\left[  n\right]  }%
x_{i,k}\underbrace{y_{j,\ell}f\left(  b_{k},b_{\ell}\right)  }_{=f\left(
b_{k},b_{\ell}\right)  y_{j,\ell}}\ \ \ \ \ \ \ \ \ \ \left(  \text{since
}f\text{ is }\mathbf{k}\text{-bilinear}\right) \\
&  =\sum_{k\in\left[  n\right]  }\ \ \sum_{\ell\in\left[  n\right]  }%
x_{i,k}f\left(  b_{k},b_{\ell}\right)  y_{j,\ell}.
\end{align*}
Comparing these two equalities, we obtain%
\[
\left(  XGY^{T}\right)  _{i,j}=f\left(  u_{i},v_{j}\right)
\ \ \ \ \ \ \ \ \ \ \text{for all }i,j\in\left[  n\right]  .
\]
Hence, the $n\times n$-matrix $XGY^{T}$ equals $\left(  f\left(  u_{i}%
,v_{j}\right)  \right)  _{i,j\in\left[  n\right]  }$. This proves
(\ref{pf.lem.bilin.gram.XGYT}).
\end{proof}

But we assumed that the $n\times n$-matrix $\left(  f\left(  u_{i}%
,v_{j}\right)  \right)  _{i,j\in\left[  n\right]  }$ is invertible. Because of
(\ref{pf.lem.bilin.gram.XGYT}), we can rewrite this as follows: The $n\times
n$-matrix $XGY^{T}$ is invertible. Hence, its determinant $\det\left(
XGY^{T}\right)  $ is invertible in $\mathbf{k}$ (because a square matrix is
invertible if and only if its determinant is invertible\footnote{This is a
classical fact (see, e.g., \cite[Theorem 6.110 \textbf{(a)}]{detnotes}).}).

But it is well-known that the determinant function is multiplicative (i.e.,
the determinant of a product of $n\times n$-matrices is the product of their
determinants). Since $X$, $G$ and $Y^{T}$ are square matrices, we thus obtain%
\[
\det\left(  XGY^{T}\right)  =\det X\cdot\det G\cdot\underbrace{\det\left(
Y^{T}\right)  }_{\substack{=\det Y\\\text{(since }\det\left(  A^{T}\right)
=\det A\\\text{for any square matrix }A\text{)}}}=\det X\cdot\det G\cdot\det
Y.
\]
Thus, the product $\det X\cdot\det G\cdot\det Y$ is invertible in $\mathbf{k}$
(since $\det\left(  XGY^{T}\right)  $ is invertible in $\mathbf{k}$).
Therefore, its three factors $\det X$ and $\det G$ and $\det Y$ are all
invertible in $\mathbf{k}$ (since a product of elements of a commutative ring
is invertible if and only if all its factors are invertible).

In particular, $\det X$ is invertible. Hence, the matrix $X$ is invertible
(since a square matrix is invertible if and only if its determinant is
invertible). Thus, its transpose $X^{T}$ is invertible as well (since the
transpose of an invertible matrix is invertible). Therefore, the $n$-tuple
$\left(  u_{1},u_{2},\ldots,u_{n}\right)  $ is a basis of $V$%
\ \ \ \ \footnote{\textit{Proof.} Let $\alpha:V\rightarrow V$ be the
$\mathbf{k}$-linear map that sends the basis vectors $b_{1},b_{2},\ldots
,b_{n}$ of the basis $\overrightarrow{b}$ to the vectors $u_{1},u_{2}%
,\ldots,u_{n}$, respectively (this map is well-defined, since a $\mathbf{k}%
$-linear map from a free $\mathbf{k}$-module $M$ can be defined by specifying
its values on a given basis of $M$). Then, for each $i\in\left[  n\right]  $,
we have%
\begin{align*}
\alpha\left(  b_{i}\right)   &  =u_{i}\ \ \ \ \ \ \ \ \ \ \left(  \text{since
the map }\alpha\text{ sends }b_{i}\text{ to }u_{i}\text{, by its
definition}\right) \\
&  =\sum_{k\in\left[  n\right]  }x_{i,k}b_{k}\ \ \ \ \ \ \ \ \ \ \left(
\text{by (\ref{pf.lem.bilin.gram.ui=})}\right)  .
\end{align*}
Thus, the matrix representing the $\mathbf{k}$-linear map $\alpha:V\rightarrow
V$ with respect to the basis $\overrightarrow{b}=\left(  b_{1},b_{2}%
,\ldots,b_{n}\right)  $ is the matrix $\left(  x_{i,k}\right)  _{k,i\in\left[
n\right]  }$ (by Definition \ref{def.bas.repmat}). But this matrix $\left(
x_{i,k}\right)  _{k,i\in\left[  n\right]  }$ is just $X^{T}$ (since $X=\left(
x_{i,k}\right)  _{i,k\in\left[  n\right]  }$ and thus $X^{T}=\left(
x_{i,k}\right)  _{k,i\in\left[  n\right]  }$), and thus is invertible (since
we have shown that $X^{T}$ is invertible). Hence, we have shown that the
matrix representing the $\mathbf{k}$-linear map $\alpha:V\rightarrow V$ with
respect to the basis $\overrightarrow{b}$ is invertible. Therefore, the
$\mathbf{k}$-linear map $\alpha:V\rightarrow V$ is invertible as well (since a
$\mathbf{k}$-linear map that is represented by an invertible matrix is itself
invertible). In other words, $\alpha:V\rightarrow V$ is a $\mathbf{k}$-module
isomorphism.
\par
This isomorphism $\alpha:V\rightarrow V$ sends the $n$ vectors $b_{1}%
,b_{2},\ldots,b_{n}$ to the $n$ vectors $u_{1},u_{2},\ldots,u_{n}$ (by its
definition). But a $\mathbf{k}$-module isomorphism always sends a basis of its
domain to a basis of its target. Since the isomorphism $\alpha$ sends
$b_{1},b_{2},\ldots,b_{n}$ to $u_{1},u_{2},\ldots,u_{n}$, we thus conclude
that $\left(  u_{1},u_{2},\ldots,u_{n}\right)  $ is a basis of $V$ (since
$\left(  b_{1},b_{2},\ldots,b_{n}\right)  $ is a basis of $V$).}. A similar
argument (using the invertibility of $\det Y$) shows that $\left(  v_{1}%
,v_{2},\ldots,v_{n}\right)  $ is a basis of $V$ as well. Thus, Lemma
\ref{lem.bilin.gram} is proved.
\end{proof}

\begin{exercise}
\fbox{2} Under the assumptions of Lemma \ref{lem.bilin.gram}, prove that the
bilinear form $f$ is perfect (i.e., its left-curried form $f_{L}:V\rightarrow
V^{\ast}$ and its right-curried form $f_{R}:V\rightarrow V^{\ast}$ are
$\mathbf{k}$-module isomorphisms).
\end{exercise}

\begin{exercise}
\fbox{2} Consider the setting of Lemma \ref{lem.bilin.gram} again, but drop
the assumption that $V$ be free of rank $n$. Prove that both $n$-tuples
$\left(  u_{1},u_{2},\ldots,u_{n}\right)  $ and $\left(  v_{1},v_{2}%
,\ldots,v_{n}\right)  $ are $\mathbf{k}$-linearly independent (although no
longer necessarily bases).
\end{exercise}

For convenience, let us state a version of Lemma \ref{lem.bilin.gram} in which
the two tuples $\left(  u_{1},u_{2},\ldots,u_{n}\right)  $ and $\left(
v_{1},v_{2},\ldots,v_{n}\right)  $ have been reindexed as families $\left(
u_{i}\right)  _{i\in I}$ and $\left(  v_{i}\right)  _{i\in I}$:

\begin{lemma}
[Gram trick, redux]\label{lem.bilin.gram-fam}Let $n\in\mathbb{N}$. Let $V$ be
a free $\mathbf{k}$-module of rank $n$, and let $\left(  u_{i}\right)  _{i\in
I}$ and $\left(  v_{i}\right)  _{i\in I}$ be two families of elements of $V$,
where $I$ is a finite set of size $\left\vert I\right\vert =n$. Let $f:V\times
V\rightarrow\mathbf{k}$ be a bilinear form. Assume that the $I\times
I$-matrix
\[
\left(  f\left(  u_{i},v_{j}\right)  \right)  _{i,j\in I}\in\mathbf{k}%
^{I\times I}%
\]
is invertible. (This matrix is called the \emph{Gram matrix} of $f$ with
respect to the families $\left(  u_{i}\right)  _{i\in I}$ and $\left(
v_{i}\right)  _{i\in I}$.)

Then, both families $\left(  u_{i}\right)  _{i\in I}$ and $\left(
v_{i}\right)  _{i\in I}$ are bases of $V$.
\end{lemma}

\begin{proof}
Rename the $n$ elements of $I$ as $1,2,\ldots,n$; then, this is just Lemma
\ref{lem.bilin.gram}.
\end{proof}

\subsubsection{\label{subsec.bas.mur.proof-bas}Proof of the basis property}

We shall now show two lemmas about the bilinear form $\left\langle \cdot
,\cdot\right\rangle $ on $\mathcal{A}$ and its interaction with the Murphy
bases. These lemmas will be crucial for proving Theorem \ref{thm.bas.mur.bas}:

\begin{lemma}
[Murphy diagonal lemma]\label{lem.bas.mur.dia}Let $\left(  \lambda,A,B\right)
\in\operatorname*{SBT}\left(  n\right)  $. Then,
\[
\left\langle \nabla_{B,A}^{-\operatorname*{Col}},\ \nabla_{B,A}%
^{\operatorname*{Row}}\right\rangle =\pm1.
\]

\end{lemma}

\begin{proof}
[Proof of Lemma \ref{lem.bas.mur.dia}.]Let $\mathbf{E}_{P}$ denote the Young
symmetrizer of an $n$-tableau $P$ (as defined in Definition
\ref{def.specht.ET.defs} \textbf{(b)}).

We have
\[
\underbrace{S\left(  \nabla_{B,A}^{-\operatorname*{Col}}\right)
}_{\substack{=\nabla_{A,B}^{-\operatorname*{Col}}\\\text{(by
(\ref{eq.prop.bas.mur.S.Col}))}}}\nabla_{B,A}^{\operatorname*{Row}}%
=\nabla_{A,B}^{-\operatorname*{Col}}\nabla_{B,A}^{\operatorname*{Row}}%
=\pm\mathbf{E}_{A}%
\]
(by Lemma \ref{lem.bas.mur.PQQP}, applied to $P=A$ and $Q=B$). Now, from
(\ref{eq.prop.bas.bilin.sym-etc.1}), we obtain%
\begin{align*}
\left\langle \nabla_{B,A}^{-\operatorname*{Col}},\ \nabla_{B,A}%
^{\operatorname*{Row}}\right\rangle  &  =\left[
\underbrace{\operatorname*{id}}_{=1}\right]  \left(  \underbrace{S\left(
\nabla_{B,A}^{-\operatorname*{Col}}\right)  \nabla_{B,A}^{\operatorname*{Row}%
}}_{=\pm\mathbf{E}_{A}}\right)  =\left[  1\right]  \left(  \pm\mathbf{E}%
_{A}\right) \\
&  =\pm\underbrace{\left[  1\right]  \left(  \mathbf{E}_{A}\right)
}_{\substack{=1\\\text{(by Lemma \ref{lem.specht.ET.1-coord})}}}=\pm1.
\end{align*}
This proves Lemma \ref{lem.bas.mur.dia}.
\end{proof}

\begin{lemma}
[Murphy orthogonality lemma]\label{lem.bas.mur.col}Let $\left(  \lambda
,A,B\right)  \in\operatorname*{SBT}\left(  n\right)  $ and $\left(
\mu,C,D\right)  \in\operatorname*{SBT}\left(  n\right)  $. Assume that $A<C$
or $B<D$, where the symbol \textquotedblleft$<$\textquotedblright\ refers to
the relation $<$ from Definition \ref{def.bas.tord.syt}. Then,%
\[
\left\langle \nabla_{B,A}^{-\operatorname*{Col}},\ \nabla_{D,C}%
^{\operatorname*{Row}}\right\rangle =0.
\]

\end{lemma}

\begin{proof}
[Proof of Lemma \ref{lem.bas.mur.col}.]Both $A$ and $B$ are standard
$n$-tableaux of shape $Y\left(  \lambda\right)  $ (since $\left(
\lambda,A,B\right)  \in\operatorname*{SBT}\left(  n\right)  $). Hence,
Proposition \ref{prop.tableau.Sn-act.1} \textbf{(b)} (applied to $Y\left(
\lambda\right)  $ and $B$ instead of $D$ and $T$) shows that any $n$-tableau
of shape $Y\left(  \lambda\right)  $ can be written as $w\rightharpoonup B$
for some $w\in S_{n}$. Hence, in particular, the $n$-tableau $A$ can be
written in this way. In other words, $A=w\rightharpoonup B$ for some $w\in
S_{n}$. Let us denote this $w$ by $u$. Thus, $A=u\rightharpoonup B=uB$.
Applying (\ref{eq.lem.bas.mur.perm.C}) to $P=A$ and $Q=B$, we thus obtain%
\[
\nabla_{A,B}^{-\operatorname*{Col}}=\left(  -1\right)  ^{u}\nabla
_{\operatorname*{Col}A}^{-}u=\left(  -1\right)  ^{u}u\nabla
_{\operatorname*{Col}B}^{-}.
\]

Both $C$ and $D$ are standard $n$-tableaux of shape $Y\left(  \mu\right)  $
(since $\left(  \mu,C,D\right)  \in\operatorname*{SBT}\left(  n\right)  $).
Hence, Proposition \ref{prop.tableau.Sn-act.1} \textbf{(b)} (applied to
$Y\left(  \mu\right)  $ and $C$ instead of $D$ and $T$) shows that any
$n$-tableau of shape $Y\left(  \mu\right)  $ can be written as
$w\rightharpoonup C$ for some $w\in S_{n}$. Hence, in particular, the
$n$-tableau $D$ can be written in this way. In other words,
$D=w\rightharpoonup C$ for some $w\in S_{n}$. Let us denote this $w$ by $v$.
Thus, $D=v\rightharpoonup C=vC$. Applying (\ref{eq.lem.bas.mur.perm.R}) to
$D$, $C$ and $v$ instead of $P$, $Q$ and $u$, we thus obtain%
\[
\nabla_{D,C}^{\operatorname*{Row}}=\nabla_{\operatorname*{Row}D}%
v=v\nabla_{\operatorname*{Row}C}.
\]

Furthermore, $S\left(  \nabla_{B,A}^{-\operatorname*{Col}}\right)
=\nabla_{A,B}^{-\operatorname*{Col}}$ (by (\ref{eq.prop.bas.mur.S.Col})).

We have assumed that $A<C$ or $B<D$. Thus, we are in one of the following two
cases (which may overlap):

\textit{Case 1:} We have $A<C$.

\textit{Case 2:} We have $B<D$.

Let us consider Case 1 first. In this case, we have $A<C$. Hence, Lemma
\ref{lem.bas.tord.syt.czs} (applied to $P=A$ and $Q=C$) yields
\[
\nabla_{\operatorname*{Col}A}^{-}\nabla_{\operatorname*{Row}C}%
=0\ \ \ \ \ \ \ \ \ \ \text{and}\ \ \ \ \ \ \ \ \ \ \nabla
_{\operatorname*{Row}C}\nabla_{\operatorname*{Col}A}^{-}=0.
\]
Now,%
\[
\underbrace{\nabla_{D,C}^{\operatorname*{Row}}}_{=v\nabla_{\operatorname*{Row}%
C}}\ \ \underbrace{S\left(  \nabla_{B,A}^{-\operatorname*{Col}}\right)
}_{\substack{=\nabla_{A,B}^{-\operatorname*{Col}}\\=\left(  -1\right)
^{u}\nabla_{\operatorname*{Col}A}^{-}u}}=\left(  -1\right)  ^{u}\cdot
v\underbrace{\nabla_{\operatorname*{Row}C}\nabla_{\operatorname*{Col}A}^{-}%
}_{=0}u=0.
\]
Hence, (\ref{eq.prop.bas.bilin.sym-etc.2}) yields $\left\langle \nabla
_{B,A}^{-\operatorname*{Col}},\ \nabla_{D,C}^{\operatorname*{Row}%
}\right\rangle =\left[  \operatorname*{id}\right]  \left(  \underbrace{\nabla
_{D,C}^{\operatorname*{Row}}S\left(  \nabla_{B,A}^{-\operatorname*{Col}%
}\right)  }_{=0}\right)  =\left[  \operatorname*{id}\right]  0=0$. Thus, Lemma
\ref{lem.bas.mur.col} is proved in Case 1.

Let us consider Case 2 next. In this case, we have $B<D$. Hence, Lemma
\ref{lem.bas.tord.syt.czs} (applied to $P=B$ and $Q=D$) yields
\[
\nabla_{\operatorname*{Col}B}^{-}\nabla_{\operatorname*{Row}D}%
=0\ \ \ \ \ \ \ \ \ \ \text{and}\ \ \ \ \ \ \ \ \ \ \nabla
_{\operatorname*{Row}D}\nabla_{\operatorname*{Col}B}^{-}=0.
\]
Now,
\[
\underbrace{S\left(  \nabla_{B,A}^{-\operatorname*{Col}}\right)
}_{\substack{=\nabla_{A,B}^{-\operatorname*{Col}}\\=\left(  -1\right)
^{u}u\nabla_{\operatorname*{Col}B}^{-}}}\ \ \underbrace{\nabla_{D,C}%
^{\operatorname*{Row}}}_{=\nabla_{\operatorname*{Row}D}v}=\left(  -1\right)
^{u}u\underbrace{\nabla_{\operatorname*{Col}B}^{-}\nabla_{\operatorname*{Row}%
D}}_{=0}v=0.
\]
Hence, (\ref{eq.prop.bas.bilin.sym-etc.1}) yields $\left\langle \nabla
_{B,A}^{-\operatorname*{Col}},\ \nabla_{D,C}^{\operatorname*{Row}%
}\right\rangle =\left[  \operatorname*{id}\right]  \left(
\underbrace{S\left(  \nabla_{B,A}^{-\operatorname*{Col}}\right)  \nabla
_{D,C}^{\operatorname*{Row}}}_{=0}\right)  =\left[  \operatorname*{id}\right]
0=0$. Thus, Lemma \ref{lem.bas.mur.col} is proved in Case 2.

So the proof of Lemma \ref{lem.bas.mur.col} is complete in both Cases 1 and 2.
\end{proof}

Now, we can prove Theorem \ref{thm.bas.mur.bas}:

\begin{proof}
[Proof of Theorem \ref{thm.bas.mur.bas}.]Recall the total order on
$\operatorname*{SBT}\left(  n\right)  $ defined in Definition
\ref{def.bas.tord.sbt.tord}. Thus, $\operatorname*{SBT}\left(  n\right)  $ is
a finite totally ordered set with size $\left\vert \operatorname*{SBT}\left(
n\right)  \right\vert =n!$ (by Lemma \ref{lem.sbt.sizen!}).

Now, let $M$ be the $\operatorname*{SBT}\left(  n\right)  \times
\operatorname*{SBT}\left(  n\right)  $-matrix%
\[
\left(  \left\langle \nabla_{B,A}^{-\operatorname*{Col}},\ \nabla
_{D,C}^{\operatorname*{Row}}\right\rangle \right)  _{\left(  \lambda
,A,B\right)  ,\left(  \mu,C,D\right)  \in\operatorname*{SBT}\left(  n\right)
}\in\mathbf{k}^{\operatorname*{SBT}\left(  n\right)  \times\operatorname*{SBT}%
\left(  n\right)  }.
\]
We claim that this matrix $M$ is \textquotedblleft lower-triangular with all
its diagonal entries being invertible\textquotedblright. Here is what this
means in concrete terms:

\begin{itemize}
\item The scalar $\left\langle \nabla_{B,A}^{-\operatorname*{Col}}%
,\ \nabla_{B,A}^{\operatorname*{Row}}\right\rangle \in\mathbf{k}$ is
invertible for each $\left(  \lambda,A,B\right)  \in\operatorname*{SBT}\left(
n\right)  $ (because Lemma \ref{lem.bas.mur.dia} shows that this scalar is
$\pm1$, but $\pm1$ is always invertible).

\item We have $\left\langle \nabla_{B,A}^{-\operatorname*{Col}},\ \nabla
_{D,C}^{\operatorname*{Row}}\right\rangle =0$ for all $\left(  \lambda
,A,B\right)  ,\left(  \mu,C,D\right)  \in\operatorname*{SBT}\left(  n\right)
$ that satisfy $\left(  \lambda,A,B\right)  <\left(  \mu,C,D\right)  $.
(Indeed, if $\left(  \lambda,A,B\right)  ,\left(  \mu,C,D\right)
\in\operatorname*{SBT}\left(  n\right)  $ satisfy $\left(  \lambda,A,B\right)
<\left(  \mu,C,D\right)  $, then we have $A<C$ or $B<D$ (by Lemma
\ref{lem.bas.tord.sbt.1} \textbf{(b)}), and therefore Lemma
\ref{lem.bas.mur.col} yields $\left\langle \nabla_{B,A}^{-\operatorname*{Col}%
},\ \nabla_{D,C}^{\operatorname*{Row}}\right\rangle =0$.)
\end{itemize}

Thus, Lemma \ref{lem.bas.mur.tria-inv} (applied to the set
$\operatorname*{SBT}\left(  n\right)  $ and the matrix \newline$M=\left(
\left\langle \nabla_{B,A}^{-\operatorname*{Col}},\ \nabla_{D,C}%
^{\operatorname*{Row}}\right\rangle \right)  _{\left(  \lambda,A,B\right)
,\left(  \mu,C,D\right)  \in\operatorname*{SBT}\left(  n\right)  }$ instead of
the set $I$ and the matrix $A=\left(  a_{i,j}\right)  _{i,j\in I}$) shows that
the matrix $M$ is invertible.

But the definition of $M$ shows that $M$ is actually the Gram matrix of the
form $\left\langle \cdot,\cdot\right\rangle :\mathcal{A}\times\mathcal{A}%
\rightarrow\mathbf{k}$ with respect to the two families%
\[
\left(  \nabla_{V,U}^{-\operatorname*{Col}}\right)  _{\left(  \lambda
,U,V\right)  \in\operatorname*{SBT}\left(  n\right)  }%
\ \ \ \ \ \ \ \ \ \ \text{and}\ \ \ \ \ \ \ \ \ \ \left(  \nabla
_{V,U}^{\operatorname*{Row}}\right)  _{\left(  \lambda,U,V\right)
\in\operatorname*{SBT}\left(  n\right)  }%
\]
(see Lemma \ref{lem.bilin.gram-fam} for the definition of \textquotedblleft
Gram matrix\textquotedblright). Hence, Lemma \ref{lem.bilin.gram-fam} (applied
to the number $n!$, the free $\mathbf{k}$-module $\mathcal{A}$, the two
families $\left(  \nabla_{V,U}^{-\operatorname*{Col}}\right)  _{\left(
\lambda,U,V\right)  \in\operatorname*{SBT}\left(  n\right)  }$ and $\left(
\nabla_{V,U}^{\operatorname*{Row}}\right)  _{\left(  \lambda,U,V\right)
\in\operatorname*{SBT}\left(  n\right)  }$ and the bilinear form $\left\langle
\cdot,\cdot\right\rangle $ instead of the number $n$, the free $\mathbf{k}%
$-module $V$, the two families $\left(  u_{i}\right)  _{i\in I}$ and $\left(
v_{i}\right)  _{i\in I}$ and the bilinear form $f$) shows that both families%
\[
\left(  \nabla_{V,U}^{-\operatorname*{Col}}\right)  _{\left(  \lambda
,U,V\right)  \in\operatorname*{SBT}\left(  n\right)  }%
\ \ \ \ \ \ \ \ \ \ \text{and}\ \ \ \ \ \ \ \ \ \ \left(  \nabla
_{V,U}^{\operatorname*{Row}}\right)  _{\left(  \lambda,U,V\right)
\in\operatorname*{SBT}\left(  n\right)  }%
\]
are bases of $\mathcal{A}$ (since $\mathcal{A}$ is a free $\mathbf{k}$-module
of rank $n!$, and since the Gram matrix $\left(  \left\langle \nabla
_{B,A}^{-\operatorname*{Col}},\ \nabla_{D,C}^{\operatorname*{Row}%
}\right\rangle \right)  _{\left(  \lambda,A,B\right)  ,\left(  \mu,C,D\right)
\in\operatorname*{SBT}\left(  n\right)  }=M$ is invertible). Hence, Theorem
\ref{thm.bas.mur.bas} is proved.
\end{proof}

In other words, we now know that the Murphy bases are bases indeed.

\begin{remark}
Other proofs of Theorem \ref{thm.bas.mur.bas} appear in \cite[Theorems 3.12,
3.14 and 6.13]{CanWil89} and \cite[Lemma 3.5 and Lemma 3.7]{Murphy92}. These
proofs have the advantage of giving a more efficient algorithm to expand a
standard basis vector $w$ in the Murphy bases. (Our proof gives a rather
horrible algorithm, since it relies on the proof of Lemma
\ref{lem.bilin.gram-fam} and thus involves inverting an $n!\times n!$-matrix.)
\end{remark}

\begin{exercise}
\fbox{1} Reprove Corollary \ref{cor.spechtmod.sumflam} using Theorem
\ref{thm.bas.mur.bas}.
\end{exercise}

\subsubsection{\label{subsec.bas.mur.triang}The Murphy spans}

Now we come to the study of spans of subfamilies of the Murphy bases. We will
show that some of these spans are left ideals (i.e., left $\mathcal{A}%
$-submodules) of $\mathcal{A}$, which explains the block-upper-triangularity
of the matrices $L_{\overrightarrow{\operatorname*{rmu}}}\left(
\mathbf{a}\right)  $ in Subsection \ref{subsec.bas.mur.n=3}. (The
block-upper-triangularity of the matrices
$R_{\overrightarrow{\operatorname*{rmu}}}\left(  \mathbf{a}\right)  $ can be
explained similarly, or derived from the former via the antipode.)

We recall that $\operatorname*{SBT}\left(  n\right)  $ is the set of all
standard $n$-bitableaux, i.e., of all triples $\left(  \lambda,U,V\right)  $,
where $\lambda$ is a partition of $n$ and where $U$ and $V$ are two standard
$n$-tableaux of shape $\lambda$. We shall now define a close relative of this
set: the set $\operatorname*{HSBT}\left(  n\right)  $ of all
\emph{half-standard }$n$\emph{-bitableaux}:

\begin{definition}
\label{def.bas.mur.HSBT}A \emph{half-standard }$n$\emph{-bitableau} shall mean
a triple $\left(  \lambda,U,V\right)  $, where $\lambda$ is a partition of $n$
and where $U$ and $V$ are two $n$-tableaux of shape $\lambda$ such that $U$ is
standard. Thus, the difference between a half-standard $n$-bitableau $\left(
\lambda,U,V\right)  $ and a standard one is that in the former, $V$ is not
required to be standard (but $U$ is), whereas both $U$ and $V$ are required to
be standard in the latter.

We let $\operatorname*{HSBT}\left(  n\right)  $ denote the set of all
half-standard $n$-bitableaux.
\end{definition}

For example, $\left(  \left(  2,1\right)
,\ \ \ytableaushort{12,3}\ ,\ \ \ytableaushort{23,1}\right)  $ is a
half-standard $3$-bitableau, thus belongs to $\operatorname*{HSBT}\left(
3\right)  $.

Obviously, $\operatorname*{SBT}\left(  n\right)  \subseteq\operatorname*{HSBT}%
\left(  n\right)  $, since any standard $n$-bitableau is a half-standard
$n$-bitableau as well.

\begin{definition}
\label{def.bas.mur.submods}Let $T\in\bigcup\limits_{\kappa\vdash
n}\operatorname*{SYT}\left(  \kappa\right)  $ be a standard $n$-tableau. Then,
we define the following eight $\mathbf{k}$-submodules of $\mathcal{A}$:

\begin{itemize}
\item the \emph{standard row-Murphy spans}%
\begin{align*}
\MurpF_{\operatorname*{std},\geq T}^{\operatorname*{Row}}  &
:=\operatorname*{span}\left\{  \nabla_{V,U}^{\operatorname*{Row}}%
\ \mid\ \left(  \lambda,U,V\right)  \in\operatorname*{SBT}\left(  n\right)
\text{ and }U\geq T\right\}  ;\\
\MurpF_{\operatorname*{std},>T}^{\operatorname*{Row}}  &
:=\operatorname*{span}\left\{  \nabla_{V,U}^{\operatorname*{Row}}%
\ \mid\ \left(  \lambda,U,V\right)  \in\operatorname*{SBT}\left(  n\right)
\text{ and }U>T\right\}  ;
\end{align*}

\item the \emph{standard column-Murphy spans}%
\begin{align*}
\MurpF_{\operatorname*{std},\leq T}^{-\operatorname*{Col}}  &
:=\operatorname*{span}\left\{  \nabla_{V,U}^{-\operatorname*{Col}}%
\ \mid\ \left(  \lambda,U,V\right)  \in\operatorname*{SBT}\left(  n\right)
\text{ and }U\leq T\right\}  ;\\
\MurpF_{\operatorname*{std},<T}^{-\operatorname*{Col}}  &
:=\operatorname*{span}\left\{  \nabla_{V,U}^{-\operatorname*{Col}}%
\ \mid\ \left(  \lambda,U,V\right)  \in\operatorname*{SBT}\left(  n\right)
\text{ and }U<T\right\}  ;
\end{align*}

\item the \emph{non-standard row-Murphy spans}%
\begin{align*}
\MurpF_{\operatorname*{all},\geq T}^{\operatorname*{Row}}  &
:=\operatorname*{span}\left\{  \nabla_{V,U}^{\operatorname*{Row}}%
\ \mid\ \left(  \lambda,U,V\right)  \in\operatorname*{HSBT}\left(  n\right)
\text{ and }U\geq T\right\}  ;\\
\MurpF_{\operatorname*{all},>T}^{\operatorname*{Row}}  &
:=\operatorname*{span}\left\{  \nabla_{V,U}^{\operatorname*{Row}}%
\ \mid\ \left(  \lambda,U,V\right)  \in\operatorname*{HSBT}\left(  n\right)
\text{ and }U>T\right\}  ;
\end{align*}

\item the \emph{non-standard column-Murphy spans}%
\begin{align*}
\MurpF_{\operatorname*{all},\leq T}^{-\operatorname*{Col}}  &
:=\operatorname*{span}\left\{  \nabla_{V,U}^{-\operatorname*{Col}}%
\ \mid\ \left(  \lambda,U,V\right)  \in\operatorname*{HSBT}\left(  n\right)
\text{ and }U\leq T\right\}  ;\\
\MurpF_{\operatorname*{all},<T}^{-\operatorname*{Col}}  &
:=\operatorname*{span}\left\{  \nabla_{V,U}^{-\operatorname*{Col}}%
\ \mid\ \left(  \lambda,U,V\right)  \in\operatorname*{HSBT}\left(  n\right)
\text{ and }U<T\right\}  ,
\end{align*}

\end{itemize}

\noindent where all the inequality signs ($\geq$, $>$, $\leq$ and $<$) are
understood with respect to the total order on $\bigcup\limits_{\kappa\vdash
n}\operatorname*{SYT}\left(  \kappa\right)  $ (see Definition
\ref{def.bas.tord.syt.tord}).

(The reader can easily see the pattern behind these naming conventions: The
\textquotedblleft$\operatorname*{std}$\textquotedblright\ and
\textquotedblleft$\operatorname*{all}$\textquotedblright\ subscripts refer to
the $n$-tableau $V$ being either standard or arbitrary (respectively), whereas
the \textquotedblleft$\geq T$\textquotedblright, \textquotedblleft%
$>T$\textquotedblright, \textquotedblleft$\leq T$\textquotedblright\ and
\textquotedblleft$<T$\textquotedblright\ subscripts refer to the inequalities
that $U$ is required to satisfy. Other combinations of requirements are
possible, but we will only use the above eight.)
\end{definition}

\begin{example}
For $n=3$ and $T=13\backslash\backslash2$, we have%
\begin{align*}
\MurpF_{\operatorname*{std},\geq T}^{\operatorname*{Row}}  &
=\operatorname*{span}\left\{  \nabla_{13\backslash\backslash2,\ 13\backslash
\backslash2}^{\operatorname*{Row}},\ \nabla_{13\backslash\backslash
2,\ 12\backslash\backslash3}^{\operatorname*{Row}},\ \nabla_{12\backslash
\backslash3,\ 13\backslash\backslash2}^{\operatorname*{Row}},\right. \\
&  \ \ \ \ \ \ \ \ \ \ \ \ \ \ \ \ \ \ \ \ \left.  \ \nabla_{12\backslash
\backslash3,\ 12\backslash\backslash3}^{\operatorname*{Row}},\ \nabla
_{123,\ 123}^{\operatorname*{Row}}\right\}  ;\\
\MurpF_{\operatorname*{std},>T}^{\operatorname*{Row}}  &
=\operatorname*{span}\left\{  \nabla_{13\backslash\backslash2,\ 12\backslash
\backslash3}^{\operatorname*{Row}},\ \nabla_{12\backslash\backslash
3,\ 12\backslash\backslash3}^{\operatorname*{Row}},\ \nabla_{123,\ 123}%
^{\operatorname*{Row}}\right\}  ;
\end{align*}%
\begin{align*}
\MurpF_{\operatorname*{std},\leq T}^{-\operatorname*{Col}}  &
=\operatorname*{span}\left\{  \nabla_{1\backslash\backslash2\backslash
\backslash3,\ 1\backslash\backslash2\backslash\backslash3}%
^{-\operatorname*{Col}},\ \nabla_{13\backslash\backslash2,\ 13\backslash
\backslash2}^{-\operatorname*{Col}},\ \nabla_{12\backslash\backslash
3,\ 13\backslash\backslash2}^{-\operatorname*{Col}}\right\}  ;\\
\MurpF_{\operatorname*{std},<T}^{-\operatorname*{Col}}  &
=\operatorname*{span}\left\{  \nabla_{1\backslash\backslash2\backslash
\backslash3,\ 1\backslash\backslash2\backslash\backslash3}%
^{-\operatorname*{Col}}\right\}  .
\end{align*}
The $\mathbf{k}$-modules $\MurpF_{\operatorname*{all},\geq T}%
^{\operatorname*{Row}},\ \MurpF_{\operatorname*{all},>T}^{\operatorname*{Row}%
},\ \MurpF_{\operatorname*{all},\leq T}^{-\operatorname*{Col}}%
,\ \MurpF_{\operatorname*{all},<T}^{-\operatorname*{Col}}$ are similar but
have additional generators $\nabla_{V,U}^{\operatorname*{Row}}$ or
$\nabla_{V,U}^{-\operatorname*{Col}}$ in which $V$ is a non-standard tableau;
for instance, one of the additional generators of $\MurpF_{\operatorname*{all}%
,>T}^{\operatorname*{Row}}$ is $\nabla_{23\backslash\backslash1,\ 12\backslash
\backslash3}^{\operatorname*{Row}}$.
\end{example}

Thus we have defined eight submodules of $\mathcal{A}$. Soon we will see
(Theorem \ref{thm.bas.mur.triang.orth}) that there are in fact only four of
them, since each standard span equals the corresponding non-standard one
(e.g., we have $\MurpF_{\operatorname*{std},\geq T}^{\operatorname*{Row}%
}=\MurpF_{\operatorname*{all},\geq T}^{\operatorname*{Row}}$). This is why
many of our results about these submodules will come in four variants.

To state the main result of this subsection, we need a classical notion from
the linear algebra of bilinear forms:

\begin{definition}
\label{def.bas.bilin.orthcomp}For each subset $V$ of $\mathcal{A}$, we define
a subset $V^{\perp}$ of $\mathcal{A}$ by
\begin{equation}
V^{\perp}:=\left\{  \mathbf{a}\in\mathcal{A}\ \mid\ \left\langle
\mathbf{a},\mathbf{v}\right\rangle =0\text{ for each }\mathbf{v}\in V\right\}
. \label{eq.def.bas.bilin.orthcomp.def}%
\end{equation}
This subset $V^{\perp}$ is easily seen to be a $\mathbf{k}$-submodule of
$\mathcal{A}$ (since the dot product $\left\langle \cdot,\cdot\right\rangle $
is bilinear, and thus $\left\langle \mathbf{a},\mathbf{v}\right\rangle $
depends $\mathbf{k}$-linearly on $\mathbf{a}$). We call $V^{\perp}$ the
\emph{orthogonal complement} of $V$.
\end{definition}

Because of the equality (\ref{eq.prop.bas.bilin.sym-etc.sym}), we can replace
the condition \textquotedblleft$\left\langle \mathbf{a},\mathbf{v}%
\right\rangle =0$\textquotedblright\ in (\ref{eq.def.bas.bilin.orthcomp.def})
by the equivalent condition \textquotedblleft$\left\langle \mathbf{v}%
,\mathbf{a}\right\rangle =0$\textquotedblright.

We can now state the main claim of this subsection:\footnote{We note that a
left $\mathcal{A}$-submodule of $\mathcal{A}$ is the same as a left ideal of
$\mathcal{A}$. Thus, Theorem \ref{thm.bas.mur.triang.orth} provides several
left ideals of $\mathcal{A}$.}

\begin{theorem}
\label{thm.bas.mur.triang.orth}Let $T\in\bigcup\limits_{\kappa\vdash
n}\operatorname*{SYT}\left(  \kappa\right)  $ be a standard $n$-tableau. Then:
\medskip

\textbf{(a)} The set $\MurpF_{\operatorname*{std},\geq T}^{\operatorname*{Row}%
}$ is a left $\mathcal{A}$-submodule of $\mathcal{A}$ and satisfies%
\[
\MurpF_{\operatorname*{std},\geq T}^{\operatorname*{Row}}%
=\MurpF_{\operatorname*{all},\geq T}^{\operatorname*{Row}}=\left(
\MurpF_{\operatorname*{std},<T}^{-\operatorname*{Col}}\right)  ^{\perp}.
\]

\textbf{(b)} The set $\MurpF_{\operatorname*{std},>T}^{\operatorname*{Row}}$
is a left $\mathcal{A}$-submodule of $\mathcal{A}$ and satisfies%
\[
\MurpF_{\operatorname*{std},>T}^{\operatorname*{Row}}%
=\MurpF_{\operatorname*{all},>T}^{\operatorname*{Row}}=\left(
\MurpF_{\operatorname*{std},\leq T}^{-\operatorname*{Col}}\right)  ^{\perp}.
\]

\textbf{(c)} The set $\MurpF_{\operatorname*{std},<T}^{-\operatorname*{Col}}$
is a left $\mathcal{A}$-submodule of $\mathcal{A}$ and satisfies
\[
\MurpF_{\operatorname*{std},<T}^{-\operatorname*{Col}}%
=\MurpF_{\operatorname*{all},<T}^{-\operatorname*{Col}}=\left(
\MurpF_{\operatorname*{std},\geq T}^{\operatorname*{Row}}\right)  ^{\perp}.
\]

\textbf{(d)} The set $\MurpF_{\operatorname*{std},\leq T}%
^{-\operatorname*{Col}}$ is a left $\mathcal{A}$-submodule of $\mathcal{A}$
and satisfies
\[
\MurpF_{\operatorname*{std},\leq T}^{-\operatorname*{Col}}%
=\MurpF_{\operatorname*{all},\leq T}^{-\operatorname*{Col}}=\left(
\MurpF_{\operatorname*{std},>T}^{\operatorname*{Row}}\right)  ^{\perp}.
\]

\end{theorem}

It will take us a while to prove this, not least because of the
similar-but-not-quite-the-same nature of the four parts. First, however, we
observe some near-obvious equalities:

\begin{proposition}
\label{prop.bas.mur.submods.eq1}Let $T_{1},T_{2},\ldots,T_{m}$ be all the
standard tableaux in the set $\bigcup\limits_{\kappa\vdash n}%
\operatorname*{SYT}\left(  \kappa\right)  $, listed in increasing order (with
respect to the total order defined in Definition \ref{def.bas.tord.syt.tord}).
Then: \medskip

\textbf{(a)} We have
\[
\MurpF_{\operatorname*{std},\geq T_{1}}^{\operatorname*{Row}}=\mathcal{A}%
\ \ \ \ \ \ \ \ \ \ \text{and}\ \ \ \ \ \ \ \ \ \ \MurpF_{\operatorname*{std}%
,\leq T_{m}}^{-\operatorname*{Col}}=\mathcal{A}.
\]

\textbf{(b)} We have
\[
\MurpF_{\operatorname*{std},>T_{m}}^{\operatorname*{Row}}%
=0\ \ \ \ \ \ \ \ \ \ \text{and}%
\ \ \ \ \ \ \ \ \ \ \MurpF_{\operatorname*{std},<T_{1}}^{-\operatorname*{Col}%
}=0.
\]

\textbf{(c)} Let $i\in\left[  m-1\right]  $. Then,
\begin{align*}
\MurpF_{\operatorname*{std},\geq T_{i+1}}^{\operatorname*{Row}}  &
=\MurpF_{\operatorname*{std},>T_{i}}^{\operatorname*{Row}}%
\ \ \ \ \ \ \ \ \ \ \text{and}\ \ \ \ \ \ \ \ \ \ \MurpF_{\operatorname*{std}%
,\leq T_{i}}^{-\operatorname*{Col}}=\MurpF_{\operatorname*{std},<T_{i+1}%
}^{-\operatorname*{Col}}\ \ \ \ \ \ \ \ \ \ \text{and}\\
\MurpF_{\operatorname*{all},\geq T_{i+1}}^{\operatorname*{Row}}  &
=\MurpF_{\operatorname*{all},>T_{i}}^{\operatorname*{Row}}%
\ \ \ \ \ \ \ \ \ \ \text{and}\ \ \ \ \ \ \ \ \ \ \MurpF_{\operatorname*{all}%
,\leq T_{i}}^{-\operatorname*{Col}}=\MurpF_{\operatorname*{all},<T_{i+1}%
}^{-\operatorname*{Col}}.
\end{align*}

\end{proposition}

\begin{proof}
We first observe that the set $\bigcup\limits_{\kappa\vdash n}%
\operatorname*{SYT}\left(  \kappa\right)  $ is nonempty\footnote{Indeed, the
$n$-tableau of shape $\left(  n\right)  $ with entries $1,2,\ldots,n$ arranged
from left to right in its only row is clearly an element of this set. Thus,
this set has at least one element.}. However, $T_{1},T_{2},\ldots,T_{m}$ were
defined to be all the standard tableaux in this set, listed without
repetition. Hence, $m=\left\vert \bigcup\limits_{\kappa\vdash n}%
\operatorname*{SYT}\left(  \kappa\right)  \right\vert >0$ (since the set
$\bigcup\limits_{\kappa\vdash n}\operatorname*{SYT}\left(  \kappa\right)  $ is
nonempty). This shows that $T_{1}$ and $T_{m}$ are well-defined.

Recall that $T_{1},T_{2},\ldots,T_{m}$ are all the standard tableaux in the
set $\bigcup\limits_{\kappa\vdash n}\operatorname*{SYT}\left(  \kappa\right)
$, listed in increasing order. Hence, $T_{1}$ is the smallest element of the
set $\bigcup\limits_{\kappa\vdash n}\operatorname*{SYT}\left(  \kappa\right)
$, and $T_{m}$ is the largest element of this set. \medskip

\textbf{(a)} Each tableau $U\in\bigcup\limits_{\kappa\vdash n}%
\operatorname*{SYT}\left(  \kappa\right)  $ satisfies $U\geq T_{1}$ (since
$T_{1}$ is the smallest element of the set $\bigcup\limits_{\kappa\vdash
n}\operatorname*{SYT}\left(  \kappa\right)  $). Thus,
\begin{equation}
\text{each }\left(  \lambda,U,V\right)  \in\operatorname*{SBT}\left(
n\right)  \text{ satisfies }U\geq T_{1}
\label{pf.prop.bas.mur.submods.eq1.UT1}%
\end{equation}
(because each $\left(  \lambda,U,V\right)  \in\operatorname*{SBT}\left(
n\right)  $ satisfies $U\in\bigcup\limits_{\kappa\vdash n}\operatorname*{SYT}%
\left(  \kappa\right)  $).

On the other hand, each tableau $U\in\bigcup\limits_{\kappa\vdash
n}\operatorname*{SYT}\left(  \kappa\right)  $ satisfies $U\leq T_{m}$ (since
$T_{m}$ is the largest element of the set $\bigcup\limits_{\kappa\vdash
n}\operatorname*{SYT}\left(  \kappa\right)  $). Thus,
\begin{equation}
\text{each }\left(  \lambda,U,V\right)  \in\operatorname*{SBT}\left(
n\right)  \text{ satisfies }U\leq T_{m}
\label{pf.prop.bas.mur.submods.eq1.UTm}%
\end{equation}
(because each $\left(  \lambda,U,V\right)  \in\operatorname*{SBT}\left(
n\right)  $ satisfies $U\in\bigcup\limits_{\kappa\vdash n}\operatorname*{SYT}%
\left(  \kappa\right)  $).

Theorem \ref{thm.bas.mur.bas} shows that both families%
\[
\left(  \nabla_{V,U}^{\operatorname*{Row}}\right)  _{\left(  \lambda
,U,V\right)  \in\operatorname*{SBT}\left(  n\right)  }%
\ \ \ \ \ \ \ \ \ \ \text{and}\ \ \ \ \ \ \ \ \ \ \left(  \nabla
_{V,U}^{-\operatorname*{Col}}\right)  _{\left(  \lambda,U,V\right)
\in\operatorname*{SBT}\left(  n\right)  }%
\]
are bases of $\mathcal{A}$. Thus, they span $\mathcal{A}$. In other words,%
\begin{align*}
\operatorname*{span}\left\{  \nabla_{V,U}^{\operatorname*{Row}}\ \mid\ \left(
\lambda,U,V\right)  \in\operatorname*{SBT}\left(  n\right)  \right\}   &
=\mathcal{A}\ \ \ \ \ \ \ \ \ \ \text{and}\\
\operatorname*{span}\left\{  \nabla_{V,U}^{-\operatorname*{Col}}%
\ \mid\ \left(  \lambda,U,V\right)  \in\operatorname*{SBT}\left(  n\right)
\right\}   &  =\mathcal{A}.
\end{align*}

Now, the definition of $\MurpF_{\operatorname*{std},\geq T_{1}}%
^{\operatorname*{Row}}$ shows that%
\begin{align*}
\MurpF_{\operatorname*{std},\geq T_{1}}^{\operatorname*{Row}}  &
=\operatorname*{span}\underbrace{\left\{  \nabla_{V,U}^{\operatorname*{Row}%
}\ \mid\ \left(  \lambda,U,V\right)  \in\operatorname*{SBT}\left(  n\right)
\text{ and }U\geq T_{1}\right\}  }_{\substack{=\left\{  \nabla_{V,U}%
^{\operatorname*{Row}}\ \mid\ \left(  \lambda,U,V\right)  \in
\operatorname*{SBT}\left(  n\right)  \right\}  \\\text{(here, we have removed
the condition \textquotedblleft}U\geq T_{1}\text{\textquotedblright%
,}\\\text{since (\ref{pf.prop.bas.mur.submods.eq1.UT1}) shows that this
condition}\\\text{is automatically satisfied for each }\left(  \lambda
,U,V\right)  \in\operatorname*{SBT}\left(  n\right)  \text{)}}}\\
&  =\operatorname*{span}\left\{  \nabla_{V,U}^{\operatorname*{Row}}%
\ \mid\ \left(  \lambda,U,V\right)  \in\operatorname*{SBT}\left(  n\right)
\right\}  =\mathcal{A}.
\end{align*}
Likewise, the definition of $\MurpF_{\operatorname*{std},\leq T_{m}%
}^{-\operatorname*{Col}}$ shows that%
\begin{align*}
\MurpF_{\operatorname*{std},\leq T_{m}}^{-\operatorname*{Col}}  &
=\operatorname*{span}\underbrace{\left\{  \nabla_{V,U}^{-\operatorname*{Col}%
}\ \mid\ \left(  \lambda,U,V\right)  \in\operatorname*{SBT}\left(  n\right)
\text{ and }U\leq T_{m}\right\}  }_{\substack{=\left\{  \nabla_{V,U}%
^{-\operatorname*{Col}}\ \mid\ \left(  \lambda,U,V\right)  \in
\operatorname*{SBT}\left(  n\right)  \right\}  \\\text{(here, we have removed
the condition \textquotedblleft}U\leq T_{m}\text{\textquotedblright%
,}\\\text{since (\ref{pf.prop.bas.mur.submods.eq1.UTm}) shows that this
condition}\\\text{is automatically satisfied for each }\left(  \lambda
,U,V\right)  \in\operatorname*{SBT}\left(  n\right)  \text{)}}}\\
&  =\operatorname*{span}\left\{  \nabla_{V,U}^{-\operatorname*{Col}}%
\ \mid\ \left(  \lambda,U,V\right)  \in\operatorname*{SBT}\left(  n\right)
\right\}  =\mathcal{A}.
\end{align*}
Hence, the proof of Proposition \ref{prop.bas.mur.submods.eq1} \textbf{(a)} is
complete. \medskip

\textbf{(b)} There exists no $\left(  \lambda,U,V\right)  \in
\operatorname*{SBT}\left(  n\right)  $ that satisfies $U<T_{1}$ (since every
$\left(  \lambda,U,V\right)  \in\operatorname*{SBT}\left(  n\right)  $
satisfies $U\geq T_{1}$ by (\ref{pf.prop.bas.mur.submods.eq1.UT1}), and thus
cannot satisfy $U<T_{1}$). Now, the definition of $\MurpF_{\operatorname*{std}%
,<T_{1}}^{-\operatorname*{Col}}$ shows that%
\[
\MurpF_{\operatorname*{std},<T_{1}}^{-\operatorname*{Col}}%
=\operatorname*{span}\underbrace{\left\{  \nabla_{V,U}^{-\operatorname*{Col}%
}\ \mid\ \left(  \lambda,U,V\right)  \in\operatorname*{SBT}\left(  n\right)
\text{ and }U<T_{1}\right\}  }_{\substack{=\varnothing\\\text{(since there
exists no }\left(  \lambda,U,V\right)  \in\operatorname*{SBT}\left(  n\right)
\text{ that satisfies }U<T_{1}\text{)}}}=\operatorname*{span}\varnothing=0.
\]

There exists no $\left(  \lambda,U,V\right)  \in\operatorname*{SBT}\left(
n\right)  $ that satisfies $U>T_{m}$ (since every $\left(  \lambda,U,V\right)
\in\operatorname*{SBT}\left(  n\right)  $ satisfies $U\leq T_{m}$ by
(\ref{pf.prop.bas.mur.submods.eq1.UTm}), and thus cannot satisfy $U>T_{m}$).
Now, the definition of $\MurpF_{\operatorname*{std},>T_{m}}%
^{\operatorname*{Row}}$ shows that%
\[
\MurpF_{\operatorname*{std},>T_{m}}^{\operatorname*{Row}}=\operatorname*{span}%
\underbrace{\left\{  \nabla_{V,U}^{\operatorname*{Row}}\ \mid\ \left(
\lambda,U,V\right)  \in\operatorname*{SBT}\left(  n\right)  \text{ and
}U>T_{m}\right\}  }_{\substack{=\varnothing\\\text{(since there exists no
}\left(  \lambda,U,V\right)  \in\operatorname*{SBT}\left(  n\right)  \text{
that satisfies }U>T_{m}\text{)}}}=\operatorname*{span}\varnothing=0.
\]
Thus, the proof of Proposition \ref{prop.bas.mur.submods.eq1} \textbf{(b)} is
complete. \medskip

\textbf{(c)} The tableaux $T_{i}$ and $T_{i+1}$ are consecutive entries of the
list $\left(  T_{1},T_{2},\ldots,T_{m}\right)  $, which contains each standard
tableau in $\bigcup\limits_{\kappa\vdash n}\operatorname*{SYT}\left(
\kappa\right)  $. Hence, for an arbitrary standard $n$-tableau $U\in
\bigcup\limits_{\kappa\vdash n}\operatorname*{SYT}\left(  \kappa\right)  $, we
have the logical equivalence%
\begin{equation}
\left(  U\geq T_{i+1}\right)  \ \Longleftrightarrow\ \left(  U>T_{i}\right)
\label{pf.prop.bas.mur.submods.eq1.1}%
\end{equation}
\footnote{\textit{Proof.} Let $U\in\bigcup\limits_{\kappa\vdash n}%
\operatorname*{SYT}\left(  \kappa\right)  $ be arbitrary. Then, $U=T_{j}$ for
some $j\in\left[  m\right]  $ (since $T_{1},T_{2},\ldots,T_{m}$ are all the
standard tableaux in the set $\bigcup\limits_{\kappa\vdash n}%
\operatorname*{SYT}\left(  \kappa\right)  $). Consider this $j$.
\par
But the list $\left(  T_{1},T_{2},\ldots,T_{m}\right)  $ is ordered in
increasing order (since this is how we defined $T_{1},T_{2},\ldots,T_{m}$).
Hence, we have $T_{j}>T_{i}$ if and only if $j>i$. For the same reason, we
have $T_{j}\geq T_{i+1}$ if and only if $j\geq i+1$. Thus, we have the
following chain of logical equivalences:%
\begin{align*}
\left(  T_{j}\geq T_{i+1}\right)  \  &  \Longleftrightarrow\ \left(  j\geq
i+1\right)  \ \Longleftrightarrow\ \left(  j>i\right)
\ \ \ \ \ \ \ \ \ \ \left(  \text{since }j\text{ and }i\text{ are
integers}\right) \\
&  \Longleftrightarrow\ \left(  T_{j}>T_{i}\right)
\ \ \ \ \ \ \ \ \ \ \left(  \text{since we have }T_{j}>T_{i}\text{ if and only
if }j>i\right)  .
\end{align*}
In other words,%
\[
\left(  U\geq T_{i+1}\right)  \ \Longleftrightarrow\ \left(  U>T_{i}\right)
\]
(since $U=T_{j}$). Qed.}. But the definition of $\MurpF_{\operatorname*{std}%
,\geq T_{i+1}}^{\operatorname*{Row}}$ shows that%
\begin{align*}
\MurpF_{\operatorname*{std},\geq T_{i+1}}^{\operatorname*{Row}}  &
=\operatorname*{span}\underbrace{\left\{  \nabla_{V,U}^{\operatorname*{Row}%
}\ \mid\ \left(  \lambda,U,V\right)  \in\operatorname*{SBT}\left(  n\right)
\text{ and }U\geq T_{i+1}\right\}  }_{\substack{=\left\{  \nabla
_{V,U}^{\operatorname*{Row}}\ \mid\ \left(  \lambda,U,V\right)  \in
\operatorname*{SBT}\left(  n\right)  \text{ and }U>T_{i}\right\}
\\\text{(because of the equivalence (\ref{pf.prop.bas.mur.submods.eq1.1}))}%
}}\\
&  =\operatorname*{span}\left\{  \nabla_{V,U}^{\operatorname*{Row}}%
\ \mid\ \left(  \lambda,U,V\right)  \in\operatorname*{SBT}\left(  n\right)
\text{ and }U>T_{i}\right\} \\
&  =\MurpF_{\operatorname*{std},>T_{i}}^{\operatorname*{Row}}%
\ \ \ \ \ \ \ \ \ \ \left(  \text{by the definition of }%
\MurpF_{\operatorname*{std},>T_{i}}^{\operatorname*{Row}}\right)  .
\end{align*}
Hence, we have proved that $\MurpF_{\operatorname*{std},\geq T_{i+1}%
}^{\operatorname*{Row}}=\MurpF_{\operatorname*{std},>T_{i}}%
^{\operatorname*{Row}}$. Similarly, we can show that
$\MurpF_{\operatorname*{all},\geq T_{i+1}}^{\operatorname*{Row}}%
=\MurpF_{\operatorname*{all},>T_{i}}^{\operatorname*{Row}}$. The remaining two
equalities $\MurpF_{\operatorname*{std},\leq T_{i}}^{-\operatorname*{Col}%
}=\MurpF_{\operatorname*{std},<T_{i+1}}^{-\operatorname*{Col}}$ and
$\MurpF_{\operatorname*{all},\leq T_{i}}^{-\operatorname*{Col}}%
=\MurpF_{\operatorname*{all},<T_{i+1}}^{-\operatorname*{Col}}$ can be proved
in an analogous way, where instead of the equivalence
(\ref{pf.prop.bas.mur.submods.eq1.1}) we must now use the equivalence%
\begin{equation}
\left(  U\leq T_{i}\right)  \ \Longleftrightarrow\ \left(  U<T_{i+1}\right)  .
\label{pf.prop.bas.mur.submods.eq1.2}%
\end{equation}
(The proof of the equivalence (\ref{pf.prop.bas.mur.submods.eq1.2}) is similar
to the proof of (\ref{pf.prop.bas.mur.submods.eq1.1}); alternatively, it can
be easily derived from (\ref{pf.prop.bas.mur.submods.eq1.1})\footnote{Here is
how: Let $U\in\bigcup\limits_{\kappa\vdash n}\operatorname*{SYT}\left(
\kappa\right)  $ be arbitrary. Then, we have the logical equivalence $\left(
U\leq T_{i}\right)  \ \Longleftrightarrow\ \left(  \text{not }U>T_{i}\right)
$ (since the order on $\bigcup\limits_{\kappa\vdash n}\operatorname*{SYT}%
\left(  \kappa\right)  $ is a total order) and the logical equivalence
$\left(  U<T_{i+1}\right)  \ \Longleftrightarrow\ \left(  \text{not }U\geq
T_{i+1}\right)  $ (similarly). Hence, we have the following chain of
equivalences:%
\begin{align*}
\left(  U<T_{i+1}\right)  \  &  \Longleftrightarrow\ \left(  \text{not }U\geq
T_{i+1}\right)  \ \Longleftrightarrow\ \left(  \text{not }U>T_{i}\right)
\ \ \ \ \ \ \ \ \ \ \left(  \text{by (\ref{pf.prop.bas.mur.submods.eq1.1}%
)}\right) \\
&  \Longleftrightarrow\ \left(  U\leq T_{i}\right)
\ \ \ \ \ \ \ \ \ \ \left(  \text{since }\left(  U\leq T_{i}\right)
\ \Longleftrightarrow\ \left(  \text{not }U>T_{i}\right)  \right)  .
\end{align*}
In other words, $\left(  U\leq T_{i}\right)  \ \Longleftrightarrow\ \left(
U<T_{i+1}\right)  $. This proves (\ref{pf.prop.bas.mur.submods.eq1.2}).}.)
Altogether, all claims of Proposition \ref{prop.bas.mur.submods.eq1}
\textbf{(c)} are now proved.
\end{proof}

The following two facts are also easy to see:

\begin{proposition}
\label{prop.bas.mur.Sn-act.all}Let $T\in\bigcup\limits_{\kappa\vdash
n}\operatorname*{SYT}\left(  \kappa\right)  $ be a standard $n$-tableau. Then,
the four sets $\MurpF_{\operatorname*{all},\geq T}^{\operatorname*{Row}}$,
$\MurpF_{\operatorname*{all},>T}^{\operatorname*{Row}}$,
$\MurpF_{\operatorname*{all},\leq T}^{-\operatorname*{Col}}$ and
$\MurpF_{\operatorname*{all},<T}^{-\operatorname*{Col}}$ are left
$\mathcal{A}$-submodules of $\mathcal{A}$ (that is, subrepresentations of the
left regular representation $\mathcal{A}$ of $S_{n}$).
\end{proposition}

\begin{proof}
Here is a sketch of the (easy) proof; details can be proved in the appendix
(Section \ref{sec.details.bas.mur}).

We shall only prove the claim about $\MurpF_{\operatorname*{all},\geq
T}^{\operatorname*{Row}}$, since the other three claims are analogous. We thus
need to show that $\MurpF_{\operatorname*{all},\geq T}^{\operatorname*{Row}}$
is a left $\mathcal{A}$-submodule of $\mathcal{A}$, that is, a
subrepresentation of the $S_{n}$-representation $\mathcal{A}$ (since
$\mathcal{A}=\mathbf{k}\left[  S_{n}\right]  $).

We defined the $\mathbf{k}$-module $\MurpF_{\operatorname*{all},\geq
T}^{\operatorname*{Row}}$ as the span of the vectors $\nabla_{V,U}%
^{\operatorname*{Row}}$ for $\left(  \lambda,U,V\right)  \in
\operatorname*{HSBT}\left(  n\right)  $ satisfying $U\geq T$; let us refer to
these vectors as the \textquotedblleft good nablas\textquotedblright. But
every $w\in S_{n}$ and every $\left(  \lambda,U,V\right)  \in
\operatorname*{HSBT}\left(  n\right)  $ satisfy%
\[
\nabla_{wV,U}^{\operatorname*{Row}}=w\nabla_{V,U}^{\operatorname*{Row}}%
\]
(by (\ref{eq.prop.bas.mur.Sn-act.Row}), applied to $P=V$ and $Q=U$ and $g=w$
and $h=\operatorname*{id}$). Hence, if $\nabla_{V,U}^{\operatorname*{Row}}$ is
a good nabla and $w\in S_{n}$ is a permutation, then $w\nabla_{V,U}%
^{\operatorname*{Row}}$ is again a good nabla (namely, $\nabla_{wV,U}%
^{\operatorname*{Row}}$). Thus, the span of all good nablas is an $S_{n}%
$-subset of $\mathcal{A}$, and hence a subrepresentation of the $S_{n}%
$-representation $\mathcal{A}$ (since it is a $\mathbf{k}$-submodule of
$\mathcal{A}$ as well). But this span is $\MurpF_{\operatorname*{all},\geq
T}^{\operatorname*{Row}}$. Hence, $\MurpF_{\operatorname*{all},\geq
T}^{\operatorname*{Row}}$ is a subrepresentation of the $S_{n}$-representation
$\mathcal{A}$, just as we wanted to show.
\end{proof}

\begin{proposition}
\label{prop.bas.mur.subset}Let $T\in\bigcup\limits_{\kappa\vdash
n}\operatorname*{SYT}\left(  \kappa\right)  $ be a standard $n$-tableau.
Then,
\[
\MurpF_{\operatorname*{std},>T}^{\operatorname*{Row}}\subseteq
\MurpF_{\operatorname*{std},\geq T}^{\operatorname*{Row}}%
\ \ \ \ \ \ \ \ \ \ \text{and}\ \ \ \ \ \ \ \ \ \ \MurpF_{\operatorname*{std}%
,<T}^{-\operatorname*{Col}}\subseteq\MurpF_{\operatorname*{std},\leq
T}^{-\operatorname*{Col}}.
\]

\end{proposition}

\begin{proof}
We have%
\begin{align*}
\MurpF_{\operatorname*{std},>T}^{\operatorname*{Row}}  &
=\operatorname*{span}\underbrace{\left\{  \nabla_{V,U}^{\operatorname*{Row}%
}\ \mid\ \left(  \lambda,U,V\right)  \in\operatorname*{SBT}\left(  n\right)
\text{ and }U>T\right\}  }_{\substack{\subseteq\left\{  \nabla_{V,U}%
^{\operatorname*{Row}}\ \mid\ \left(  \lambda,U,V\right)  \in
\operatorname*{SBT}\left(  n\right)  \text{ and }U\geq T\right\}
\\\text{(since the condition \textquotedblleft}U>T\text{\textquotedblright%
\ implies \textquotedblleft}U\geq T\text{\textquotedblright)}}}\\
&  \ \ \ \ \ \ \ \ \ \ \ \ \ \ \ \ \ \ \ \ \left(  \text{by the definition of
}\MurpF_{\operatorname*{std},>T}^{\operatorname*{Row}}\right) \\
&  \subseteq\operatorname*{span}\left\{  \nabla_{V,U}^{\operatorname*{Row}%
}\ \mid\ \left(  \lambda,U,V\right)  \in\operatorname*{SBT}\left(  n\right)
\text{ and }U\geq T\right\}  =\MurpF_{\operatorname*{std},\geq T}%
^{\operatorname*{Row}}%
\end{align*}
(by the definition of $\MurpF_{\operatorname*{std},\geq T}%
^{\operatorname*{Row}}$). Thus, $\MurpF_{\operatorname*{std},>T}%
^{\operatorname*{Row}}\subseteq\MurpF_{\operatorname*{std},\geq T}%
^{\operatorname*{Row}}$ is proved. A similar argument establishes
$\MurpF_{\operatorname*{std},<T}^{-\operatorname*{Col}}\subseteq
\MurpF_{\operatorname*{std},\leq T}^{-\operatorname*{Col}}$. This completes
the proof of Proposition \ref{prop.bas.mur.subset}.
\end{proof}

To proceed, we need a variant of Lemma \ref{lem.bas.mur.col}:

\begin{lemma}
[Murphy orthogonality lemma, take 2]\label{lem.bas.mur.col2}Let $\left(
\lambda,A,B\right)  \in\operatorname*{HSBT}\left(  n\right)  $ and $\left(
\mu,C,D\right)  \in\operatorname*{HSBT}\left(  n\right)  $. Assume that $A<C$,
where the symbol \textquotedblleft$<$\textquotedblright\ refers to the
relation $<$ from Definition \ref{def.bas.tord.syt}. Then,%
\[
\left\langle \nabla_{B,A}^{-\operatorname*{Col}},\ \nabla_{D,C}%
^{\operatorname*{Row}}\right\rangle =0.
\]

\end{lemma}

\begin{proof}
This lemma is very similar to Lemma \ref{lem.bas.mur.col} except for two
differences: We are no longer requiring that the tableaux $B$ and $D$ are
standard (since $\left(  \lambda,A,B\right)  $ and $\left(  \mu,C,D\right)  $
now belong to $\operatorname*{HSBT}\left(  n\right)  $ rather than to
$\operatorname*{SBT}\left(  n\right)  $), but in exchange, we are now
requiring that $A<C$ (rather than merely $A<C$ or $B<D$). Fortunately, the
proof of Lemma \ref{lem.bas.mur.col} still works, except that Case 2 no longer exists.
\end{proof}

An easy consequence of Lemma \ref{lem.bas.mur.col2} is the following:

\begin{lemma}
\label{lem.bas.mur.orth}Let $T\in\bigcup\limits_{\kappa\vdash n}%
\operatorname*{SYT}\left(  \kappa\right)  $ be a standard $n$-tableau. Then:
\medskip

\textbf{(a)} We have $\left\langle \mathbf{a},\mathbf{b}\right\rangle =0$ for
any $\mathbf{a}\in\MurpF_{\operatorname*{all},<T}^{-\operatorname*{Col}}$ and
$\mathbf{b}\in\MurpF_{\operatorname*{all},\geq T}^{\operatorname*{Row}}$.
\medskip

\textbf{(b)} We have $\left\langle \mathbf{a},\mathbf{b}\right\rangle =0$ for
any $\mathbf{a}\in\MurpF_{\operatorname*{all},\leq T}^{-\operatorname*{Col}}$
and $\mathbf{b}\in\MurpF_{\operatorname*{all},>T}^{\operatorname*{Row}}$.
\end{lemma}

\begin{proof}
\textbf{(a)} Let $\mathbf{a}\in\MurpF_{\operatorname*{all},<T}%
^{-\operatorname*{Col}}$ and $\mathbf{b}\in\MurpF_{\operatorname*{all},\geq
T}^{\operatorname*{Row}}$. We must prove that $\left\langle \mathbf{a}%
,\mathbf{b}\right\rangle =0$.

This claim depends $\mathbf{k}$-linearly on each of $\mathbf{a}$ and
$\mathbf{b}$ (because $\left\langle \cdot,\cdot\right\rangle $ is a bilinear
form). Hence, we can WLOG assume

\begin{itemize}
\item that $\mathbf{a}$ is one of the vectors $\nabla_{V,U}%
^{-\operatorname*{Col}}$ with $\left(  \lambda,U,V\right)  \in
\operatorname*{HSBT}\left(  n\right)  $ satisfying $U<T$ (since $\mathbf{a}$
belongs to $\MurpF_{\operatorname*{all},<T}^{-\operatorname*{Col}}$, which is
defined to be the span of these vectors), and

\item that $\mathbf{b}$ is one of the vectors $\nabla_{V,U}%
^{\operatorname*{Row}}$ with $\left(  \lambda,U,V\right)  \in
\operatorname*{HSBT}\left(  n\right)  $ satisfying $U\geq T$ (since
$\mathbf{b}$ belongs to $\MurpF_{\operatorname*{all},\geq T}%
^{\operatorname*{Row}}$, which is defined to be the span of these vectors).
\end{itemize}

\noindent Assume this. Thus,
\[
\mathbf{a}=\nabla_{B,A}^{-\operatorname*{Col}}\ \ \ \ \ \ \ \ \ \ \text{for
some }\left(  \lambda,A,B\right)  \in\operatorname*{HSBT}\left(  n\right)
\text{ satisfying }A<T,
\]
and%
\[
\mathbf{b}=\nabla_{D,C}^{\operatorname*{Row}}\ \ \ \ \ \ \ \ \ \ \text{for
some }\left(  \mu,C,D\right)  \in\operatorname*{HSBT}\left(  n\right)  \text{
satisfying }C\geq T.
\]
Consider these $\left(  \lambda,A,B\right)  $ and $\left(  \mu,C,D\right)  $.
Then, $A<T\leq C$ (since $C\geq T$), so that $A<C$ (since the relation $<$ is
a total order). Thus, $\left\langle \nabla_{B,A}^{-\operatorname*{Col}%
},\ \nabla_{D,C}^{\operatorname*{Row}}\right\rangle =0$ (by Lemma
\ref{lem.bas.mur.col2}). In other words, $\left\langle \mathbf{a}%
,\mathbf{b}\right\rangle =0$ (since $\mathbf{a}=\nabla_{B,A}%
^{-\operatorname*{Col}}$ and $\mathbf{b}=\nabla_{D,C}^{\operatorname*{Row}}$).
This proves Lemma \ref{lem.bas.mur.orth} \textbf{(a)}. \medskip

\textbf{(b)} This is analogous to part \textbf{(a)}. (The main difference is
that the chain of inequalities $A<T\leq C$ must now be replaced by $A\leq T<C$.)
\end{proof}

\begin{lemma}
\label{lem.bas.mur.triang.1}Let $T\in\bigcup\limits_{\kappa\vdash
n}\operatorname*{SYT}\left(  \kappa\right)  $ be a standard $n$-tableau. Then:
\medskip

\textbf{(a)} We have $\MurpF_{\operatorname*{all},\geq T}^{\operatorname*{Row}%
}\subseteq\left(  \MurpF_{\operatorname*{std},<T}^{-\operatorname*{Col}%
}\right)  ^{\perp}$. \medskip

\textbf{(b)} We have $\MurpF_{\operatorname*{all},>T}^{\operatorname*{Row}%
}\subseteq\left(  \MurpF_{\operatorname*{std},\leq T}^{-\operatorname*{Col}%
}\right)  ^{\perp}$. \medskip

\textbf{(c)} We have $\MurpF_{\operatorname*{all},<T}^{-\operatorname*{Col}%
}\subseteq\left(  \MurpF_{\operatorname*{std},\geq T}^{\operatorname*{Row}%
}\right)  ^{\perp}$. \medskip

\textbf{(d)} We have $\MurpF_{\operatorname*{all},\leq T}%
^{-\operatorname*{Col}}\subseteq\left(  \MurpF_{\operatorname*{std}%
,>T}^{\operatorname*{Row}}\right)  ^{\perp}$.
\end{lemma}

\begin{proof}
\textbf{(a)} We must prove that $\mathbf{b}\in\left(
\MurpF_{\operatorname*{std},<T}^{-\operatorname*{Col}}\right)  ^{\perp}$ for
each $\mathbf{b}\in\MurpF_{\operatorname*{all},\geq T}^{\operatorname*{Row}}$.

So let $\mathbf{b}\in\MurpF_{\operatorname*{all},\geq T}^{\operatorname*{Row}%
}$ be arbitrary. We must then prove that $\mathbf{b}\in\left(
\MurpF_{\operatorname*{std},<T}^{-\operatorname*{Col}}\right)  ^{\perp}$. In
other words, we must prove that $\left\langle \mathbf{b},\mathbf{v}%
\right\rangle =0$ for each $\mathbf{v}\in\MurpF_{\operatorname*{std}%
,<T}^{-\operatorname*{Col}}$ (because $\left(  \MurpF_{\operatorname*{std}%
,<T}^{-\operatorname*{Col}}\right)  ^{\perp}$ is defined as $\left\{
\mathbf{a}\in\mathcal{A}\ \mid\ \left\langle \mathbf{a},\mathbf{v}%
\right\rangle =0\text{ for each }\mathbf{v}\in\MurpF_{\operatorname*{std}%
,<T}^{-\operatorname*{Col}}\right\}  $).

So let $\mathbf{v}\in\MurpF_{\operatorname*{std},<T}^{-\operatorname*{Col}}$
be arbitrary. We must prove that $\left\langle \mathbf{b},\mathbf{v}%
\right\rangle =0$.

We have
\begin{align*}
\mathbf{v}  &  \in\MurpF_{\operatorname*{std},<T}^{-\operatorname*{Col}%
}=\operatorname*{span}\left\{  \nabla_{V,U}^{-\operatorname*{Col}}%
\ \mid\ \left(  \lambda,U,V\right)  \in\underbrace{\operatorname*{SBT}\left(
n\right)  }_{\subseteq\operatorname*{HSBT}\left(  n\right)  }\text{ and
}U<T\right\} \\
&  \ \ \ \ \ \ \ \ \ \ \ \ \ \ \ \ \ \ \ \ \left(  \text{by the definition of
}\MurpF_{\operatorname*{std},<T}^{-\operatorname*{Col}}\right) \\
&  \subseteq\operatorname*{span}\left\{  \nabla_{V,U}^{-\operatorname*{Col}%
}\ \mid\ \left(  \lambda,U,V\right)  \in\operatorname*{HSBT}\left(  n\right)
\text{ and }U<T\right\}  =\MurpF_{\operatorname*{all},<T}%
^{-\operatorname*{Col}}%
\end{align*}
(by the definition of $\MurpF_{\operatorname*{all},<T}^{-\operatorname*{Col}}%
$). Hence, Lemma \ref{lem.bas.mur.orth} \textbf{(a)} (applied to
$\mathbf{a}=\mathbf{v}$) yields $\left\langle \mathbf{v},\mathbf{b}%
\right\rangle =0$. Finally, (\ref{eq.prop.bas.bilin.sym-etc.sym}) shows that
$\left\langle \mathbf{v},\mathbf{b}\right\rangle =\left\langle \mathbf{b}%
,\mathbf{v}\right\rangle $, so that $\left\langle \mathbf{b},\mathbf{v}%
\right\rangle =\left\langle \mathbf{v},\mathbf{b}\right\rangle =0$. As we
explained above, this completes the proof of Lemma \ref{lem.bas.mur.triang.1}
\textbf{(a)}. \medskip

\textbf{(b)} This is analogous to part \textbf{(a)}, but using Lemma
\ref{lem.bas.mur.orth} \textbf{(b)} instead of Lemma \ref{lem.bas.mur.orth}
\textbf{(a)}. \medskip

\textbf{(c)} This proof, too, is similar to that of part \textbf{(a)}, but
there are some minor changes, so we give it in full.

We must prove that $\mathbf{b}\in\left(  \MurpF_{\operatorname*{std},\geq
T}^{\operatorname*{Row}}\right)  ^{\perp}$ for each $\mathbf{b}\in
\MurpF_{\operatorname*{all},<T}^{-\operatorname*{Col}}$.

So let $\mathbf{b}\in\MurpF_{\operatorname*{all},<T}^{-\operatorname*{Col}}$
be arbitrary. We must then prove that $\mathbf{b}\in\left(
\MurpF_{\operatorname*{std},\geq T}^{\operatorname*{Row}}\right)  ^{\perp}$.
In other words, we must prove that $\left\langle \mathbf{b},\mathbf{v}%
\right\rangle =0$ for each $\mathbf{v}\in\MurpF_{\operatorname*{std},\geq
T}^{\operatorname*{Row}}$ (because $\left(  \MurpF_{\operatorname*{std},\geq
T}^{\operatorname*{Row}}\right)  ^{\perp}$ is defined as $\left\{
\mathbf{a}\in\mathcal{A}\ \mid\ \left\langle \mathbf{a},\mathbf{v}%
\right\rangle =0\text{ for each }\mathbf{v}\in\MurpF_{\operatorname*{std},\geq
T}^{\operatorname*{Row}}\right\}  $).

So let $\mathbf{v}\in\MurpF_{\operatorname*{std},\geq T}^{\operatorname*{Row}%
}$ be arbitrary. We must prove that $\left\langle \mathbf{b},\mathbf{v}%
\right\rangle =0$.

We have
\begin{align*}
\mathbf{v}  &  \in\MurpF_{\operatorname*{std},\geq T}^{\operatorname*{Row}%
}=\operatorname*{span}\left\{  \nabla_{V,U}^{\operatorname*{Row}}%
\ \mid\ \left(  \lambda,U,V\right)  \in\underbrace{\operatorname*{SBT}\left(
n\right)  }_{\subseteq\operatorname*{HSBT}\left(  n\right)  }\text{ and }U\geq
T\right\} \\
&  \ \ \ \ \ \ \ \ \ \ \ \ \ \ \ \ \ \ \ \ \left(  \text{by the definition of
}\MurpF_{\operatorname*{std},\geq T}^{\operatorname*{Row}}\right) \\
&  \subseteq\operatorname*{span}\left\{  \nabla_{V,U}^{\operatorname*{Row}%
}\ \mid\ \left(  \lambda,U,V\right)  \in\operatorname*{HSBT}\left(  n\right)
\text{ and }U\geq T\right\}  =\MurpF_{\operatorname*{all},\geq T}%
^{\operatorname*{Row}}%
\end{align*}
(by the definition of $\MurpF_{\operatorname*{all},\geq T}%
^{\operatorname*{Row}}$). Hence, Lemma \ref{lem.bas.mur.orth} \textbf{(a)}
(applied to $\mathbf{b}$ and $\mathbf{v}$ instead of $\mathbf{a}$ and
$\mathbf{b}$) yields $\left\langle \mathbf{b},\mathbf{v}\right\rangle =0$. As
we explained above, this completes the proof of Lemma
\ref{lem.bas.mur.triang.1} \textbf{(c)}. \medskip

\textbf{(d)} This is analogous to part \textbf{(c)}, but using Lemma
\ref{lem.bas.mur.orth} \textbf{(b)} instead of Lemma \ref{lem.bas.mur.orth}
\textbf{(a)}.
\end{proof}

\begin{lemma}
\label{lem.bas.mur.triang.2}Let $T\in\bigcup\limits_{\kappa\vdash
n}\operatorname*{SYT}\left(  \kappa\right)  $ be a standard $n$-tableau. Then:
\medskip

\textbf{(a)} We have $\left(  \MurpF_{\operatorname*{std},<T}%
^{-\operatorname*{Col}}\right)  ^{\perp}\subseteq\MurpF_{\operatorname*{std}%
,\geq T}^{\operatorname*{Row}}$. \medskip

\textbf{(b)} We have $\left(  \MurpF_{\operatorname*{std},\leq T}%
^{-\operatorname*{Col}}\right)  ^{\perp}\subseteq\MurpF_{\operatorname*{std}%
,>T}^{\operatorname*{Row}}$. \medskip

\textbf{(c)} We have $\left(  \MurpF_{\operatorname*{std},\geq T}%
^{\operatorname*{Row}}\right)  ^{\perp}\subseteq\MurpF_{\operatorname*{std}%
,<T}^{-\operatorname*{Col}}$. \medskip

\textbf{(d)} We have $\left(  \MurpF_{\operatorname*{std},>T}%
^{\operatorname*{Row}}\right)  ^{\perp}\subseteq\MurpF_{\operatorname*{std}%
,\leq T}^{-\operatorname*{Col}}$.
\end{lemma}

\begin{proof}
\textbf{(a)} We must show that $\mathbf{a}\in\MurpF_{\operatorname*{std},\geq
T}^{\operatorname*{Row}}$ for each $\mathbf{a}\in\left(
\MurpF_{\operatorname*{std},<T}^{-\operatorname*{Col}}\right)  ^{\perp}$.

So let $\mathbf{a}\in\left(  \MurpF_{\operatorname*{std},<T}%
^{-\operatorname*{Col}}\right)  ^{\perp}$. We must show that $\mathbf{a}%
\in\MurpF_{\operatorname*{std},\geq T}^{\operatorname*{Row}}$.

We have $\mathbf{a}\in\left(  \MurpF_{\operatorname*{std},<T}%
^{-\operatorname*{Col}}\right)  ^{\perp}\subseteq\mathcal{A}$. Since the row
Murphy basis $\left(  \nabla_{V,U}^{\operatorname*{Row}}\right)  _{\left(
\lambda,U,V\right)  \in\operatorname*{SBT}\left(  n\right)  }$ is a basis of
$\mathcal{A}$ (by Theorem \ref{thm.bas.mur.bas}), we can thus write
$\mathbf{a}$ as a $\mathbf{k}$-linear combination of its elements. In other
words, we can write $\mathbf{a}$ as
\begin{equation}
\mathbf{a}=\sum_{\left(  \lambda,U,V\right)  \in\operatorname*{SBT}\left(
n\right)  }\omega_{\lambda,U,V}\nabla_{V,U}^{\operatorname*{Row}}
\label{pf.lem.bas.mur.triang.2.a=wsum}%
\end{equation}
for some scalars $\omega_{\lambda,U,V}\in\mathbf{k}$. Consider these
$\omega_{\lambda,U,V}$. We can rewrite the equality
(\ref{pf.lem.bas.mur.triang.2.a=wsum}) as
\begin{equation}
\mathbf{a}=\sum_{\left(  \mu,P,Q\right)  \in\operatorname*{SBT}\left(
n\right)  }\omega_{\mu,P,Q}\nabla_{Q,P}^{\operatorname*{Row}}
\label{pf.lem.bas.mur.triang.2.a=wsum2}%
\end{equation}
(by renaming the summation index $\left(  \lambda,U,V\right)  $ as $\left(
\mu,P,Q\right)  $).

Now, we shall show the following:

\begin{statement}
\textit{Claim 1:} We have%
\[
\omega_{\lambda,U,V}=0\ \ \ \ \ \ \ \ \ \ \text{for all }\left(
\lambda,U,V\right)  \in\operatorname*{SBT}\left(  n\right)  \text{ satisfying
}U<T.
\]

\end{statement}

\begin{proof}
[Proof of Claim 1.]We proceed by strong induction on $\left(  \lambda
,U,V\right)  $ (using the total order on $\operatorname*{SBT}\left(  n\right)
$ introduced in Definition \ref{def.bas.tord.sbt.tord}). So we fix some
$\left(  \lambda,U,V\right)  \in\operatorname*{SBT}\left(  n\right)  $
satisfying $U<T$. We assume (as the induction hypothesis) that
\begin{equation}
\omega_{\mu,P,Q}=0 \label{pf.lem.bas.mur.triang.2.IH}%
\end{equation}
for all $\left(  \mu,P,Q\right)  \in\operatorname*{SBT}\left(  n\right)  $
satisfying $P<T$ and $\left(  \mu,P,Q\right)  <\left(  \lambda,U,V\right)  $.
Our goal is to prove that $\omega_{\lambda,U,V}=0$.

Since $\left(  \lambda,U,V\right)  \in\operatorname*{SBT}\left(  n\right)  $
and $U<T$, we have $\nabla_{V,U}^{-\operatorname*{Col}}\in
\MurpF_{\operatorname*{std},<T}^{-\operatorname*{Col}}$ by the definition of
$\MurpF_{\operatorname*{std},<T}^{-\operatorname*{Col}}$ (in fact,
$\nabla_{V,U}^{-\operatorname*{Col}}$ is one of the spanning vectors in the
definition of $\MurpF_{\operatorname*{std},<T}^{-\operatorname*{Col}}$).

Since $\mathbf{a}\in\left(  \MurpF_{\operatorname*{std},<T}%
^{-\operatorname*{Col}}\right)  ^{\perp}$, we have $\left\langle
\mathbf{a},\mathbf{v}\right\rangle =0$ for each $\mathbf{v}\in
\MurpF_{\operatorname*{std},<T}^{-\operatorname*{Col}}$ (by the definition of
$\left(  \MurpF_{\operatorname*{std},<T}^{-\operatorname*{Col}}\right)
^{\perp}$). Applying this to $\mathbf{v}=\nabla_{V,U}^{-\operatorname*{Col}}$,
we obtain $\left\langle \mathbf{a},\ \nabla_{V,U}^{-\operatorname*{Col}%
}\right\rangle =0$ (since $\nabla_{V,U}^{-\operatorname*{Col}}\in
\MurpF_{\operatorname*{std},<T}^{-\operatorname*{Col}}$). But
(\ref{eq.prop.bas.bilin.sym-etc.sym}) shows that $\left\langle \mathbf{a}%
,\ \nabla_{V,U}^{-\operatorname*{Col}}\right\rangle =\left\langle \nabla
_{V,U}^{-\operatorname*{Col}},\ \mathbf{a}\right\rangle $, hence $\left\langle
\nabla_{V,U}^{-\operatorname*{Col}},\ \mathbf{a}\right\rangle =\left\langle
\mathbf{a},\ \nabla_{V,U}^{-\operatorname*{Col}}\right\rangle =0$. Thus,%
\begin{align}
0  &  =\left\langle \nabla_{V,U}^{-\operatorname*{Col}},\ \mathbf{a}%
\right\rangle =\left\langle \nabla_{V,U}^{-\operatorname*{Col}},\ \sum
_{\left(  \mu,P,Q\right)  \in\operatorname*{SBT}\left(  n\right)  }\omega
_{\mu,P,Q}\nabla_{Q,P}^{\operatorname*{Row}}\right\rangle
\ \ \ \ \ \ \ \ \ \ \left(  \text{by (\ref{pf.lem.bas.mur.triang.2.a=wsum2}%
)}\right) \nonumber\\
&  =\sum_{\left(  \mu,P,Q\right)  \in\operatorname*{SBT}\left(  n\right)
}\omega_{\mu,P,Q}\left\langle \nabla_{V,U}^{-\operatorname*{Col}}%
,\ \nabla_{Q,P}^{\operatorname*{Row}}\right\rangle
\label{pf.lem.bas.mur.triang.2.0}%
\end{align}
(since the form $\left\langle \cdot,\cdot\right\rangle $ is bilinear). Now, we
observe the following:

\begin{itemize}
\item For each $\left(  \mu,P,Q\right)  \in\operatorname*{SBT}\left(
n\right)  $ satisfying $\left(  \mu,P,Q\right)  <\left(  \lambda,U,V\right)
$, we have%
\begin{equation}
\omega_{\mu,P,Q}=0. \label{pf.lem.bas.mur.triang.2.1}%
\end{equation}

[\textit{Proof:} Let $\left(  \mu,P,Q\right)  \in\operatorname*{SBT}\left(
n\right)  $ be such that $\left(  \mu,P,Q\right)  <\left(  \lambda,U,V\right)
$. Then, Lemma \ref{lem.bas.tord.sbt.1} \textbf{(a)} (applied to $\left(
\mu,P,Q\right)  $ and $\left(  \lambda,U,V\right)  $ instead of $\left(
\lambda,A,B\right)  $ and $\left(  \mu,C,D\right)  $) shows that $P\leq U$.
Thus, $P\leq U<T$, so that $P<T$ (since the relation $<$ is a total order).
Thus, (\ref{pf.lem.bas.mur.triang.2.IH}) yields $\omega_{\mu,P,Q}=0$. This
proves (\ref{pf.lem.bas.mur.triang.2.1}).]

\item For each $\left(  \mu,P,Q\right)  \in\operatorname*{SBT}\left(
n\right)  $ satisfying $\left(  \mu,P,Q\right)  >\left(  \lambda,U,V\right)
$, we have%
\begin{equation}
\left\langle \nabla_{V,U}^{-\operatorname*{Col}},\ \nabla_{Q,P}%
^{\operatorname*{Row}}\right\rangle =0. \label{pf.lem.bas.mur.triang.2.2}%
\end{equation}

[\textit{Proof:} Let $\left(  \mu,P,Q\right)  \in\operatorname*{SBT}\left(
n\right)  $ be such that $\left(  \mu,P,Q\right)  >\left(  \lambda,U,V\right)
$. Thus, $\left(  \lambda,U,V\right)  <\left(  \mu,P,Q\right)  $. Hence, Lemma
\ref{lem.bas.tord.sbt.1} \textbf{(b)} (applied to $\left(  \lambda,U,V\right)
$ and $\left(  \mu,P,Q\right)  $ instead of $\left(  \lambda,A,B\right)  $ and
$\left(  \mu,C,D\right)  $) yields $U<P$ or $V<Q$. Thus, Lemma
\ref{lem.bas.mur.col} (applied to $A=U$, $B=V$, $C=P$ and $D=Q$) yields
$\left\langle \nabla_{V,U}^{-\operatorname*{Col}},\ \nabla_{Q,P}%
^{\operatorname*{Row}}\right\rangle =0$. This proves
(\ref{pf.lem.bas.mur.triang.2.2}).]
\end{itemize}

Now, recall that the relation $<$ on $\operatorname*{SBT}\left(  n\right)  $
is the smaller relation of a total order. Hence, each $\left(  \mu,P,Q\right)
\in\operatorname*{SBT}\left(  n\right)  $ satisfies exactly one of the three
statements $\left(  \mu,P,Q\right)  <\left(  \lambda,U,V\right)  $ and
$\left(  \mu,P,Q\right)  >\left(  \lambda,U,V\right)  $ and $\left(
\mu,P,Q\right)  =\left(  \lambda,U,V\right)  $.

Accordingly, the sum on the right hand side of
(\ref{pf.lem.bas.mur.triang.2.0}) can be split into three parts: the part
comprising all $\left(  \mu,P,Q\right)  \in\operatorname*{SBT}\left(
n\right)  $ that satisfy $\left(  \mu,P,Q\right)  <\left(  \lambda,U,V\right)
$; the part comprising all $\left(  \mu,P,Q\right)  \in\operatorname*{SBT}%
\left(  n\right)  $ that satisfy $\left(  \mu,P,Q\right)  >\left(
\lambda,U,V\right)  $; and a single addend corresponding to $\left(
\mu,P,Q\right)  =\left(  \lambda,U,V\right)  $. Splitting the sum in this way,
we rewrite (\ref{pf.lem.bas.mur.triang.2.0}) as follows:
\begin{align*}
0  &  =\sum_{\substack{\left(  \mu,P,Q\right)  \in\operatorname*{SBT}\left(
n\right)  ;\\\left(  \mu,P,Q\right)  <\left(  \lambda,U,V\right)
}}\ \ \underbrace{\omega_{\mu,P,Q}}_{\substack{=0\\\text{(by
(\ref{pf.lem.bas.mur.triang.2.1}))}}}\left\langle \nabla_{V,U}%
^{-\operatorname*{Col}},\ \nabla_{Q,P}^{\operatorname*{Row}}\right\rangle \\
&  \ \ \ \ \ \ \ \ \ \ +\sum_{\substack{\left(  \mu,P,Q\right)  \in
\operatorname*{SBT}\left(  n\right)  ;\\\left(  \mu,P,Q\right)  >\left(
\lambda,U,V\right)  }}\omega_{\mu,P,Q}\underbrace{\left\langle \nabla
_{V,U}^{-\operatorname*{Col}},\ \nabla_{Q,P}^{\operatorname*{Row}%
}\right\rangle }_{\substack{=0\\\text{(by (\ref{pf.lem.bas.mur.triang.2.2}))}%
}}+\,\omega_{\lambda,U,V}\left\langle \nabla_{V,U}^{-\operatorname*{Col}%
},\ \nabla_{V,U}^{\operatorname*{Row}}\right\rangle \\
&  =\underbrace{\sum_{\substack{\left(  \mu,P,Q\right)  \in\operatorname*{SBT}%
\left(  n\right)  ;\\\left(  \mu,P,Q\right)  <\left(  \lambda,U,V\right)
}}0\left\langle \nabla_{V,U}^{-\operatorname*{Col}},\ \nabla_{Q,P}%
^{\operatorname*{Row}}\right\rangle }_{=0}+\underbrace{\sum_{\substack{\left(
\mu,P,Q\right)  \in\operatorname*{SBT}\left(  n\right)  ;\\\left(
\mu,P,Q\right)  >\left(  \lambda,U,V\right)  }}\omega_{\mu,P,Q}0}%
_{=0}+\,\omega_{\lambda,U,V}\left\langle \nabla_{V,U}^{-\operatorname*{Col}%
},\ \nabla_{V,U}^{\operatorname*{Row}}\right\rangle \\
&  =\omega_{\lambda,U,V}\underbrace{\left\langle \nabla_{V,U}%
^{-\operatorname*{Col}},\ \nabla_{V,U}^{\operatorname*{Row}}\right\rangle
}_{\substack{=\pm1\\\text{(by Lemma \ref{lem.bas.mur.dia},}\\\text{applied to
}A=U\text{ and }B=V\text{)}}}=\pm\omega_{\lambda,U,V}.
\end{align*}
Hence, $\pm\omega_{\lambda,U,V}=0$. In other words, $\omega_{\lambda,U,V}=0$.
This completes the induction, and thus the proof of Claim 1.
\end{proof}

Now, (\ref{pf.lem.bas.mur.triang.2.a=wsum}) becomes%
\begin{align*}
\mathbf{a}  &  =\sum_{\left(  \lambda,U,V\right)  \in\operatorname*{SBT}%
\left(  n\right)  }\omega_{\lambda,U,V}\nabla_{V,U}^{\operatorname*{Row}}\\
&  =\sum_{\substack{\left(  \lambda,U,V\right)  \in\operatorname*{SBT}\left(
n\right)  ;\\U<T}}\underbrace{\omega_{\lambda,U,V}}_{\substack{=0\\\text{(by
Claim 1)}}}\nabla_{V,U}^{\operatorname*{Row}}+\sum_{\substack{\left(
\lambda,U,V\right)  \in\operatorname*{SBT}\left(  n\right)  ;\\U\geq T}%
}\omega_{\lambda,U,V}\nabla_{V,U}^{\operatorname*{Row}}\\
&  \ \ \ \ \ \ \ \ \ \ \ \ \ \ \ \ \ \ \ \ \left(
\begin{array}
[c]{c}%
\text{since each }\left(  \lambda,U,V\right)  \in\operatorname*{SBT}\left(
n\right)  \text{ satisfies}\\
\text{either }U<T\text{ or }U\geq T\text{ (but not both)}%
\end{array}
\right) \\
&  =\underbrace{\sum_{\substack{\left(  \lambda,U,V\right)  \in
\operatorname*{SBT}\left(  n\right)  ;\\U<T}}0\nabla_{V,U}%
^{\operatorname*{Row}}}_{=0}+\sum_{\substack{\left(  \lambda,U,V\right)
\in\operatorname*{SBT}\left(  n\right)  ;\\U\geq T}}\omega_{\lambda,U,V}%
\nabla_{V,U}^{\operatorname*{Row}}\\
&  =\sum_{\substack{\left(  \lambda,U,V\right)  \in\operatorname*{SBT}\left(
n\right)  ;\\U\geq T}}\omega_{\lambda,U,V}\nabla_{V,U}^{\operatorname*{Row}}\\
&  \in\operatorname*{span}\left\{  \nabla_{V,U}^{\operatorname*{Row}}%
\ \mid\ \left(  \lambda,U,V\right)  \in\operatorname*{SBT}\left(  n\right)
\text{ and }U\geq T\right\} \\
&  =\MurpF_{\operatorname*{std},\geq T}^{\operatorname*{Row}}%
\ \ \ \ \ \ \ \ \ \ \left(  \text{by the definition of }%
\MurpF_{\operatorname*{std},\geq T}^{\operatorname*{Row}}\right)  .
\end{align*}
This completes our proof of Lemma \ref{lem.bas.mur.triang.2} \textbf{(a)}.
\medskip

\textbf{(b)} This is analogous to the above proof of Lemma
\ref{lem.bas.mur.triang.2} \textbf{(a)}; the only difference is that some
inequality signs (between standard $n$-tableaux) change from strong to weak
and vice versa. (For instance, the \textquotedblleft$U<T$\textquotedblright%
\ in Claim 1 must be replaced by \textquotedblleft$U\leq T$\textquotedblright%
). We leave it to the reader to make the necessary changes. \medskip

\textbf{(c)} This is sufficiently similar to part \textbf{(a)} that I am
omitting the proof here, but still sufficiently different that I am
nevertheless giving the proof in the appendix (Section
\ref{sec.details.bas.mur}). \medskip

\textbf{(d)} Just as for part \textbf{(c)}, see the appendix (Section
\ref{sec.details.bas.mur}).
\end{proof}

We are now ready to prove Theorem \ref{thm.bas.mur.triang.orth}:

\begin{proof}
[Proof of Theorem \ref{thm.bas.mur.triang.orth}.]\textbf{(a)} We have
$\MurpF_{\operatorname*{std},\geq T}^{\operatorname*{Row}}\subseteq
\MurpF_{\operatorname*{all},\geq T}^{\operatorname*{Row}}$ (this follows from
the definitions of these two sets\footnote{\textit{Proof.} The definition of
$\MurpF_{\operatorname*{std},\geq T}^{\operatorname*{Row}}$ yields%
\begin{align*}
\MurpF_{\operatorname*{std},\geq T}^{\operatorname*{Row}}  &
=\operatorname*{span}\left\{  \nabla_{V,U}^{\operatorname*{Row}}%
\ \mid\ \left(  \lambda,U,V\right)  \in\underbrace{\operatorname*{SBT}\left(
n\right)  }_{\subseteq\operatorname*{HSBT}\left(  n\right)  }\text{ and }U\geq
T\right\} \\
&  \subseteq\operatorname*{span}\left\{  \nabla_{V,U}^{\operatorname*{Row}%
}\ \mid\ \left(  \lambda,U,V\right)  \in\operatorname*{HSBT}\left(  n\right)
\text{ and }U\geq T\right\}  =\MurpF_{\operatorname*{all},\geq T}%
^{\operatorname*{Row}}%
\end{align*}
(by the definition of $\MurpF_{\operatorname*{all},\geq T}%
^{\operatorname*{Row}}$).}). Hence,%
\begin{align*}
\MurpF_{\operatorname*{std},\geq T}^{\operatorname*{Row}}  &  \subseteq
\MurpF_{\operatorname*{all},\geq T}^{\operatorname*{Row}}\subseteq\left(
\MurpF_{\operatorname*{std},<T}^{-\operatorname*{Col}}\right)  ^{\perp
}\ \ \ \ \ \ \ \ \ \ \left(  \text{by Lemma \ref{lem.bas.mur.triang.1}
\textbf{(a)}}\right) \\
&  \subseteq\MurpF_{\operatorname*{std},\geq T}^{\operatorname*{Row}%
}\ \ \ \ \ \ \ \ \ \ \left(  \text{by Lemma \ref{lem.bas.mur.triang.2}
\textbf{(a)}}\right)  .
\end{align*}
This is a chain of inclusions that starts and ends with the same set (namely,
$\MurpF_{\operatorname*{std},\geq T}^{\operatorname*{Row}}$). Thus, all
inclusion signs in this chain must actually be equality signs (since any three
sets $X,Y,Z$ satisfying $X\subseteq Y\subseteq Z\subseteq X$ must satisfy
$X=Y=Z$). That is, we must have $\MurpF_{\operatorname*{std},\geq
T}^{\operatorname*{Row}}=\MurpF_{\operatorname*{all},\geq T}%
^{\operatorname*{Row}}=\left(  \MurpF_{\operatorname*{std},<T}%
^{-\operatorname*{Col}}\right)  ^{\perp}$.

It remains to show that $\MurpF_{\operatorname*{std},\geq T}%
^{\operatorname*{Row}}$ is a left $\mathcal{A}$-submodule of $\mathcal{A}$.

But Proposition \ref{prop.bas.mur.Sn-act.all} shows that
$\MurpF_{\operatorname*{all},\geq T}^{\operatorname*{Row}}$ is a left
$\mathcal{A}$-submodule of $\mathcal{A}$. In other words,
$\MurpF_{\operatorname*{std},\geq T}^{\operatorname*{Row}}$ is a left
$\mathcal{A}$-submodule of $\mathcal{A}$ (since $\MurpF_{\operatorname*{std}%
,\geq T}^{\operatorname*{Row}}=\MurpF_{\operatorname*{all},\geq T}%
^{\operatorname*{Row}}$). Hence, the proof of Theorem
\ref{thm.bas.mur.triang.orth} \textbf{(a)} is complete. \medskip

\textbf{(b)} The proof of part \textbf{(a)} applies to part \textbf{(b)} with
minor changes (the inequality signs need to be changed, and we must use Lemma
\ref{lem.bas.mur.triang.1} \textbf{(b)} instead of Lemma
\ref{lem.bas.mur.triang.1} \textbf{(a)} and use Lemma
\ref{lem.bas.mur.triang.2} \textbf{(b)} instead of Lemma
\ref{lem.bas.mur.triang.2} \textbf{(a)}). \medskip

\textbf{(c)} The proof of part \textbf{(a)} applies to part \textbf{(c)} with
minor changes (the inequality signs need to be changed, and we must use Lemma
\ref{lem.bas.mur.triang.1} \textbf{(c)} instead of Lemma
\ref{lem.bas.mur.triang.1} \textbf{(a)} and use Lemma
\ref{lem.bas.mur.triang.2} \textbf{(c)} instead of Lemma
\ref{lem.bas.mur.triang.2} \textbf{(a)}). \medskip

\textbf{(d)} The proof of part \textbf{(a)} applies to part \textbf{(d)} with
minor changes (the inequality signs need to be changed, and we must use Lemma
\ref{lem.bas.mur.triang.1} \textbf{(d)} instead of Lemma
\ref{lem.bas.mur.triang.1} \textbf{(a)} and use Lemma
\ref{lem.bas.mur.triang.2} \textbf{(d)} instead of Lemma
\ref{lem.bas.mur.triang.2} \textbf{(a)}).
\end{proof}

\subsubsection{\label{subsec.bas.mur.triang2}The Murphy span filtrations}

We shall now reformulate the results of the previous subsection (mainly
Theorem \ref{thm.bas.mur.triang.orth}) as two filtrations of the regular
representation of $S_{n}$ (that is, the left $\mathcal{A}$-module
$\mathcal{A}$) by left $\mathcal{A}$-submodules (the row- or column-Murphy
spans). We will then recognize the subquotients of these filtrations as the
Specht modules or their duals.

The two filtrations are defined in the following theorem:

\begin{theorem}
\label{thm.bas.mur.filt}Let $T_{1},T_{2},\ldots,T_{m}$ be all the standard
tableaux in the set $\bigcup\limits_{\kappa\vdash n}\operatorname*{SYT}\left(
\kappa\right)  $, listed in increasing order (with respect to the total order
defined in Definition \ref{def.bas.tord.syt.tord}). For each $i\in\left\{
0,1,\ldots,m\right\}  $, we define the $\mathbf{k}$-submodules%
\begin{align*}
\MurpF_{i}^{\operatorname*{Row}}  &  :=%
\begin{cases}
\MurpF_{\operatorname*{std},\geq T_{i+1}}^{\operatorname*{Row}}, & \text{if
}i<m;\\
0, & \text{if }i=m
\end{cases}
\ \ \ \ \ \ \ \ \ \ \text{and}\\
\MurpF_{i}^{-\operatorname*{Col}}  &  :=%
\begin{cases}
\MurpF_{\operatorname*{std},\leq T_{i}}^{-\operatorname*{Col}}, & \text{if
}i>0;\\
0, & \text{if }i=0
\end{cases}
\ \ \ \ \ \ \ \ \ \ \text{of }\mathcal{A}.
\end{align*}
Then: \medskip

\textbf{(a)} We have $\MurpF_{i}^{\operatorname*{Row}}=%
\begin{cases}
\MurpF_{\operatorname*{std},>T_{i}}^{\operatorname*{Row}}, & \text{if }i>0;\\
\mathcal{A}, & \text{if }i=0
\end{cases}
\ \ $ for each $i\in\left\{  0,1,\ldots,m\right\}  $. \medskip

\textbf{(b)} We have $\MurpF_{i}^{-\operatorname*{Col}}=%
\begin{cases}
\MurpF_{\operatorname*{std},<T_{i+1}}^{-\operatorname*{Col}}, & \text{if
}i<m;\\
\mathcal{A}, & \text{if }i=m
\end{cases}
\ \ $ for each $i\in\left\{  0,1,\ldots,m\right\}  $. \medskip

\textbf{(c)} For each $i\in\left\{  0,1,\ldots,m\right\}  $, the set
$\MurpF_{i}^{\operatorname*{Row}}$ is a left $\mathcal{A}$-submodule of
$\mathcal{A}$. \medskip

\textbf{(d)} For each $i\in\left\{  0,1,\ldots,m\right\}  $, the set
$\MurpF_{i}^{-\operatorname*{Col}}$ is a left $\mathcal{A}$-submodule of
$\mathcal{A}$. \medskip

\textbf{(e)} We have
\[
0=\MurpF_{m}^{\operatorname*{Row}}\subseteq\MurpF_{m-1}^{\operatorname*{Row}%
}\subseteq\cdots\subseteq\MurpF_{2}^{\operatorname*{Row}}\subseteq
\MurpF_{1}^{\operatorname*{Row}}\subseteq\MurpF_{0}^{\operatorname*{Row}%
}=\mathcal{A}.
\]
That is, $\left(  \MurpF_{m}^{\operatorname*{Row}},\MurpF_{m-1}%
^{\operatorname*{Row}},\ldots,\MurpF_{2}^{\operatorname*{Row}},\MurpF_{1}%
^{\operatorname*{Row}},\MurpF_{0}^{\operatorname*{Row}}\right)  $ is a
filtration of $\mathcal{A}$ by left $\mathcal{A}$-submodules. \medskip

\textbf{(f)} We have%
\[
0=\MurpF_{0}^{-\operatorname*{Col}}\subseteq\MurpF_{1}^{-\operatorname*{Col}%
}\subseteq\MurpF_{2}^{-\operatorname*{Col}}\subseteq\cdots\subseteq
\MurpF_{m-1}^{-\operatorname*{Col}}\subseteq\MurpF_{m}^{-\operatorname*{Col}%
}=\mathcal{A}.
\]
That is, $\left(  \MurpF_{0}^{-\operatorname*{Col}},\MurpF_{1}%
^{-\operatorname*{Col}},\MurpF_{2}^{-\operatorname*{Col}},\ldots
,\MurpF_{m-1}^{-\operatorname*{Col}},\MurpF_{m}^{-\operatorname*{Col}}\right)
$ is a filtration of $\mathcal{A}$ by left $\mathcal{A}$-submodules.
\end{theorem}

\begin{example}
\label{exa.bas.mur.filt.n=3}Let $n=3$. Recall that the set $\bigcup
\limits_{\kappa\vdash n}\operatorname*{SYT}\left(  \kappa\right)  $ consists
of just four standard tableaux, which are ordered as follows:%
\[
1\backslash\backslash2\backslash\backslash3<13\backslash\backslash
2<12\backslash\backslash3<123
\]
(see Example \ref{exa.bas.tord.syt.n=3}). Thus, following the notations of
Theorem \ref{thm.bas.mur.filt}, we have $m=4$ and $T_{1}=1\backslash
\backslash2\backslash\backslash3$ and $T_{2}=13\backslash\backslash2$ and
$T_{3}=12\backslash\backslash3$ and $T_{4}=123$.

Recall also the six standard $n$-bitableaux in $\operatorname*{SBT}\left(
n\right)  $ (see Example \ref{exa.bas.tord.sbt.n=3} for a list).

Now, the filtration $\left(  \MurpF_{m}^{\operatorname*{Row}},\MurpF_{m-1}%
^{\operatorname*{Row}},\ldots,\MurpF_{2}^{\operatorname*{Row}},\MurpF_{1}%
^{\operatorname*{Row}},\MurpF_{0}^{\operatorname*{Row}}\right)  $ in Theorem
\ref{thm.bas.mur.filt} \textbf{(e)} takes the form $\left(  \MurpF_{4}%
^{\operatorname*{Row}},\MurpF_{3}^{\operatorname*{Row}},\MurpF_{2}%
^{\operatorname*{Row}},\MurpF_{1}^{\operatorname*{Row}},\MurpF_{0}%
^{\operatorname*{Row}}\right)  $ with%
\begin{align*}
\MurpF_{0}^{\operatorname*{Row}}  &  =\MurpF_{\operatorname*{std},\geq T_{1}%
}^{\operatorname*{Row}}=\mathcal{A}\\
&  =\operatorname*{span}\left\{  \nabla_{1\backslash\backslash2\backslash
\backslash3,\ 1\backslash\backslash2\backslash\backslash3}%
^{\operatorname*{Row}},\ \nabla_{13\backslash\backslash2,\ 13\backslash
\backslash2}^{\operatorname*{Row}},\ \nabla_{13\backslash\backslash
2,\ 12\backslash\backslash3}^{\operatorname*{Row}},\ \nabla_{12\backslash
\backslash3,\ 13\backslash\backslash2}^{\operatorname*{Row}},\right. \\
&  \ \ \ \ \ \ \ \ \ \ \ \ \ \ \ \ \ \ \ \ \left.  \ \nabla_{12\backslash
\backslash3,\ 12\backslash\backslash3}^{\operatorname*{Row}},\ \nabla
_{123,\ 123}^{\operatorname*{Row}}\right\}  ;
\end{align*}%
\begin{align*}
\MurpF_{1}^{\operatorname*{Row}}  &  =\MurpF_{\operatorname*{std},\geq T_{2}%
}^{\operatorname*{Row}}=\MurpF_{\operatorname*{std},>T_{1}}%
^{\operatorname*{Row}}\\
&  =\operatorname*{span}\left\{  \nabla_{13\backslash\backslash
2,\ 13\backslash\backslash2}^{\operatorname*{Row}},\ \nabla_{13\backslash
\backslash2,\ 12\backslash\backslash3}^{\operatorname*{Row}},\ \nabla
_{12\backslash\backslash3,\ 13\backslash\backslash2}^{\operatorname*{Row}%
},\ \nabla_{12\backslash\backslash3,\ 12\backslash\backslash3}%
^{\operatorname*{Row}},\right. \\
&  \ \ \ \ \ \ \ \ \ \ \ \ \ \ \ \ \ \ \ \ \left.  \ \nabla_{123,\ 123}%
^{\operatorname*{Row}}\right\}  ;
\end{align*}%
\[
\MurpF_{2}^{\operatorname*{Row}}=\MurpF_{\operatorname*{std},\geq T_{3}%
}^{\operatorname*{Row}}=\MurpF_{\operatorname*{std},>T_{2}}%
^{\operatorname*{Row}}=\operatorname*{span}\left\{  \nabla_{13\backslash
\backslash2,\ 12\backslash\backslash3}^{\operatorname*{Row}},\ \nabla
_{12\backslash\backslash3,\ 12\backslash\backslash3}^{\operatorname*{Row}%
},\ \nabla_{123,\ 123}^{\operatorname*{Row}}\right\}  ;
\]%
\[
\MurpF_{3}^{\operatorname*{Row}}=\MurpF_{\operatorname*{std},\geq T_{4}%
}^{\operatorname*{Row}}=\MurpF_{\operatorname*{std},>T_{3}}%
^{\operatorname*{Row}}=\operatorname*{span}\left\{  \nabla_{123,\ 123}%
^{\operatorname*{Row}}\right\}  ;
\]%
\[
\MurpF_{4}^{\operatorname*{Row}}=0=\MurpF_{\operatorname*{std},>T_{4}%
}^{\operatorname*{Row}}=\operatorname*{span}\left\{  {}\right\}  .
\]
Thus, this filtration is precisely the filtration
(\ref{eq.subsec.rep.maschke.kS3.filt}). We described its subquotients back in
Subsection \ref{subsec.rep.maschke.kS3}.

Likewise, the filtration $\left(  \MurpF_{0}^{-\operatorname*{Col}}%
,\MurpF_{1}^{-\operatorname*{Col}},\MurpF_{2}^{-\operatorname*{Col}}%
,\ldots,\MurpF_{m-1}^{-\operatorname*{Col}},\MurpF_{m}^{-\operatorname*{Col}%
}\right)  $ in Theorem \ref{thm.bas.mur.filt} \textbf{(f)} takes the form
$\left(  \MurpF_{0}^{-\operatorname*{Col}},\MurpF_{1}^{-\operatorname*{Col}%
},\MurpF_{2}^{-\operatorname*{Col}},\MurpF_{3}^{-\operatorname*{Col}%
},\MurpF_{4}^{-\operatorname*{Col}}\right)  $ with%
\[
\MurpF_{0}^{-\operatorname*{Col}}=0=\MurpF_{\operatorname*{std},<T_{1}%
}^{-\operatorname*{Col}}=\operatorname*{span}\left\{  {}\right\}  ;
\]%
\[
\MurpF_{1}^{-\operatorname*{Col}}=\MurpF_{\operatorname*{std},\leq T_{1}%
}^{-\operatorname*{Col}}=\MurpF_{\operatorname*{std},<T_{2}}%
^{-\operatorname*{Col}}=\operatorname*{span}\left\{  \nabla_{1\backslash
\backslash2\backslash\backslash3,\ 1\backslash\backslash2\backslash
\backslash3}^{-\operatorname*{Col}}\right\}  ;
\]%
\begin{align*}
\MurpF_{2}^{-\operatorname*{Col}}  &  =\MurpF_{\operatorname*{std},\leq T_{2}%
}^{-\operatorname*{Col}}=\MurpF_{\operatorname*{std},<T_{3}}%
^{-\operatorname*{Col}}\\
&  =\operatorname*{span}\left\{  \nabla_{1\backslash\backslash2\backslash
\backslash3,\ 1\backslash\backslash2\backslash\backslash3}%
^{-\operatorname*{Col}},\ \nabla_{13\backslash\backslash2,\ 13\backslash
\backslash2}^{-\operatorname*{Col}},\ \nabla_{12\backslash\backslash
3,\ 13\backslash\backslash2}^{-\operatorname*{Col}}\right\}  ;
\end{align*}%
\begin{align*}
\MurpF_{3}^{-\operatorname*{Col}}  &  =\MurpF_{\operatorname*{std},\leq T_{3}%
}^{-\operatorname*{Col}}=\MurpF_{\operatorname*{std},<T_{4}}%
^{-\operatorname*{Col}}\\
&  =\operatorname*{span}\left\{  \nabla_{1\backslash\backslash2\backslash
\backslash3,\ 1\backslash\backslash2\backslash\backslash3}%
^{-\operatorname*{Col}},\ \nabla_{13\backslash\backslash2,\ 13\backslash
\backslash2}^{-\operatorname*{Col}},\ \nabla_{13\backslash\backslash
2,\ 12\backslash\backslash3}^{-\operatorname*{Col}},\ \nabla_{12\backslash
\backslash3,\ 13\backslash\backslash2}^{-\operatorname*{Col}},\right. \\
&  \ \ \ \ \ \ \ \ \ \ \ \ \ \ \ \ \ \ \ \ \ \left.  \nabla_{12\backslash
\backslash3,\ 12\backslash\backslash3}^{-\operatorname*{Col}}\right\}  ;
\end{align*}%
\begin{align*}
\MurpF_{4}^{-\operatorname*{Col}}  &  =\MurpF_{\operatorname*{std},\leq T_{4}%
}^{-\operatorname*{Col}}=\mathcal{A}\\
&  =\operatorname*{span}\left\{  \nabla_{1\backslash\backslash2\backslash
\backslash3,\ 1\backslash\backslash2\backslash\backslash3}%
^{-\operatorname*{Col}},\ \nabla_{13\backslash\backslash2,\ 13\backslash
\backslash2}^{-\operatorname*{Col}},\ \nabla_{13\backslash\backslash
2,\ 12\backslash\backslash3}^{-\operatorname*{Col}},\ \nabla_{12\backslash
\backslash3,\ 13\backslash\backslash2}^{-\operatorname*{Col}},\right. \\
&  \ \ \ \ \ \ \ \ \ \ \ \ \ \ \ \ \ \ \ \ \left.  \ \nabla_{12\backslash
\backslash3,\ 12\backslash\backslash3}^{-\operatorname*{Col}},\ \nabla
_{123,\ 123}^{-\operatorname*{Col}}\right\}  .
\end{align*}

Anticipating Theorem \ref{thm.bas.mur.filtC-subq}, we observe that the
subquotients of the latter filtration are isomorphic to Specht modules:%
\begin{align*}
\MurpF_{4}^{-\operatorname*{Col}}/\MurpF_{3}^{-\operatorname*{Col}}  &
\cong\mathcal{S}^{\left(  3\right)  };\\
\MurpF_{3}^{-\operatorname*{Col}}/\MurpF_{2}^{-\operatorname*{Col}}  &
\cong\mathcal{S}^{\left(  2,1\right)  };\\
\MurpF_{2}^{-\operatorname*{Col}}/\MurpF_{1}^{-\operatorname*{Col}}  &
\cong\mathcal{S}^{\left(  2,1\right)  };\\
\MurpF_{1}^{-\operatorname*{Col}}/\MurpF_{0}^{-\operatorname*{Col}}  &
\cong\mathcal{S}^{\left(  1,1,1\right)  }.
\end{align*}

\end{example}

\begin{proof}
[Proof of Theorem \ref{thm.bas.mur.filt}.]This is all just a matter of
combining previously proved facts (Proposition \ref{prop.bas.mur.submods.eq1},
Theorem \ref{thm.bas.mur.triang.orth} and Proposition
\ref{prop.bas.mur.subset}). See the appendix (Section
\ref{sec.details.bas.mur}) for a detailed proof.
\end{proof}

As we said, the filtration of $\mathcal{A}$ constructed in Theorem
\ref{thm.bas.mur.filt} \textbf{(e)} is a filtration by left $\mathcal{A}%
$-submodules; i.e., each part $\MurpF_{i}^{\operatorname*{Row}}$ of this
filtration is a left $\mathcal{A}$-submodule of $\mathcal{A}$. Thus, left
multiplication by any given $\mathbf{a}\in\mathcal{A}$ preserves all the parts
$\MurpF_{i}^{\operatorname*{Row}}$ of this filtration, and therefore is
represented by a block-upper-triangular matrix on the Murphy basis $\left(
\nabla_{V,U}^{\operatorname*{Row}}\right)  _{\left(  \lambda,U,V\right)
\in\operatorname*{SBT}\left(  n\right)  }$. The diagonal blocks have sizes
$f^{\lambda}$ for the various partitions $\lambda$ of $n$, but each
$f^{\lambda}$ appears $f^{\lambda}$ many times.

The filtration in Theorem \ref{thm.bas.mur.filt} \textbf{(f)} also has the
same block-triangularity property. But it also has the additional nice
property that its subquotients are the Specht modules $\mathcal{S}^{\lambda}$
corresponding to the partitions $\lambda$ satisfying $T\in\operatorname*{SYT}%
\left(  \lambda\right)  $. This follows from the following property of
column-Murphy spans:

\begin{theorem}
\label{thm.bas.mur.triang.iso1}Let $\lambda$ be a partition of $n$. Let $T$ be
any standard $n$-tableau of shape $\lambda$. Then,
\[
\MurpF_{\operatorname*{std},\leq T}^{-\operatorname*{Col}}%
/\MurpF_{\operatorname*{std},<T}^{-\operatorname*{Col}}\cong\mathcal{A}%
\mathbf{E}_{T}\cong\mathcal{S}^{\lambda}\ \ \ \ \ \ \ \ \ \ \text{as }%
S_{n}\text{-representations.}%
\]
More concretely, the map%
\begin{align*}
\MurpF_{\operatorname*{std},\leq T}^{-\operatorname*{Col}}%
/\MurpF_{\operatorname*{std},<T}^{-\operatorname*{Col}}  &  \rightarrow
\mathcal{A}\mathbf{E}_{T},\\
\overline{\mathbf{a}}  &  \mapsto\mathbf{a}\nabla_{\operatorname*{Row}T}%
\end{align*}
(where $\overline{\mathbf{a}}$ denotes the projection of a vector
$\mathbf{a}\in\MurpF_{\operatorname*{std},\leq T}^{-\operatorname*{Col}}$ on
the quotient $\MurpF_{\operatorname*{std},\leq T}^{-\operatorname*{Col}%
}/\MurpF_{\operatorname*{std},<T}^{-\operatorname*{Col}}$) is a left
$\mathcal{A}$-module isomorphism.
\end{theorem}

To prove this theorem, we need a lemma:

\begin{lemma}
\label{lem.bas.mur.triang.iso1lem}Let $\lambda$ be a partition of $n$. Let $T$
be any standard $n$-tableau of shape $\lambda$. Then: \medskip

\textbf{(a)} For any $\left(  \mu,U,V\right)  \in\operatorname*{SBT}\left(
n\right)  $ satisfying $U<T$, we have $\nabla_{V,U}^{-\operatorname*{Col}%
}\nabla_{\operatorname*{Row}T}=0$. \medskip

\textbf{(b)} For any $n$-tableau $V$ of shape $\lambda$, we have $\nabla
_{V,T}^{-\operatorname*{Col}}\nabla_{\operatorname*{Row}T}=\pm\mathbf{E}%
_{V,T}$. (See Definition \ref{def.specht.EPQ} for the meaning of
$\mathbf{E}_{V,T}$.) \medskip

\textbf{(c)} Each $\mathbf{a}\in\MurpF_{\operatorname*{std},<T}%
^{-\operatorname*{Col}}$ satisfies $\mathbf{a}\nabla_{\operatorname*{Row}T}%
=0$. \medskip\ 

\textbf{(d)} Each $\mathbf{a}\in\MurpF_{\operatorname*{std},\leq
T}^{-\operatorname*{Col}}$ satisfies $\mathbf{a}\nabla_{\operatorname*{Row}%
T}\in\mathcal{A}\mathbf{E}_{T}$. \medskip\ 

\textbf{(e)} We have $\MurpF_{\operatorname*{std},\leq T}%
^{-\operatorname*{Col}}=\MurpF_{\operatorname*{std},<T}^{-\operatorname*{Col}%
}+\operatorname*{span}\left\{  \nabla_{V,T}^{-\operatorname*{Col}}\ \mid
\ V\in\operatorname*{SYT}\left(  \lambda\right)  \right\}  $.
\end{lemma}

\begin{proof}
\textbf{(a)} Let $\left(  \mu,U,V\right)  \in\operatorname*{SBT}\left(
n\right)  $ be such that $U<T$. Then, $U,V\in\operatorname*{SYT}\left(
\mu\right)  $; thus, $V$ and $U$ are $n$-tableaux of the same shape (namely,
of shape $Y\left(  \mu\right)  $). Hence, Definition \ref{def.specht.wPQ}
defines a permutation $w_{V,U}\in S_{n}$ that satisfies $w_{V,U}U=V$. Consider
this $w_{V,U}$.

We have $U\in\operatorname*{SYT}\left(  \mu\right)  $ and $T\in
\operatorname*{SYT}\left(  \lambda\right)  $ (since $T$ is a standard
$n$-tableau of shape $\lambda$). Thus, $U,T\in\bigcup\limits_{\kappa\vdash
n}\operatorname*{SYT}\left(  \kappa\right)  $. Moreover, $U<T$. Thus, Lemma
\ref{lem.bas.tord.syt.czs} (applied to $P=U$ and $Q=T$) yields $\nabla
_{\operatorname*{Col}U}^{-}\nabla_{\operatorname*{Row}T}=0$ and $\nabla
_{\operatorname*{Row}T}\nabla_{\operatorname*{Col}U}^{-}=0$.

Set $u:=w_{V,U}$; thus, $uU=w_{V,U}U=V$, so that $V=uU$. Hence,
(\ref{eq.lem.bas.mur.perm.C}) (applied to $P=V$ and $Q=U$) yields
$\nabla_{V,U}^{-\operatorname*{Col}}=\left(  -1\right)  ^{u}\nabla
_{\operatorname*{Col}V}^{-}u=\left(  -1\right)  ^{u}u\nabla
_{\operatorname*{Col}U}^{-}$. Therefore,%
\begin{equation}
\underbrace{\nabla_{V,U}^{-\operatorname*{Col}}}_{=\left(  -1\right)
^{u}u\nabla_{\operatorname*{Col}U}^{-}}\nabla_{\operatorname*{Row}T}=\left(
-1\right)  ^{u}u\underbrace{\nabla_{\operatorname*{Col}U}^{-}\nabla
_{\operatorname*{Row}T}}_{=0}=0.\nonumber
\end{equation}
This proves Lemma \ref{lem.bas.mur.triang.iso1lem} \textbf{(a)}. \medskip

\textbf{(b)} Let $V$ be any $n$-tableau of shape $\lambda$. Then, $V$ and $T$
are $n$-tableaux of the same shape (namely, of shape $Y\left(  \lambda\right)
$). Hence, Definition \ref{def.specht.wPQ} defines a permutation $w_{V,T}\in
S_{n}$ that satisfies $w_{V,T}T=V$. Consider this $w_{V,T}$.

Set $u:=w_{V,T}$; thus, $uT=w_{V,T}T=V$, so that $V=uT$. Hence,
(\ref{eq.lem.bas.mur.perm.C}) (applied to $P=V$ and $Q=T$) yields
$\nabla_{V,T}^{-\operatorname*{Col}}=\left(  -1\right)  ^{u}\nabla
_{\operatorname*{Col}V}^{-}u=\left(  -1\right)  ^{u}u\nabla
_{\operatorname*{Col}T}^{-}$. Therefore,%
\begin{align}
\underbrace{\nabla_{V,T}^{-\operatorname*{Col}}}_{=\left(  -1\right)
^{u}u\nabla_{\operatorname*{Col}T}^{-}}\nabla_{\operatorname*{Row}T}  &
=\left(  -1\right)  ^{u}\underbrace{u}_{=w_{V,T}}\underbrace{\nabla
_{\operatorname*{Col}T}^{-}\nabla_{\operatorname*{Row}T}}%
_{\substack{=\mathbf{E}_{T}\\\text{(by the definition of }\mathbf{E}%
_{T}\text{)}}}\nonumber\\
&  =\left(  -1\right)  ^{u}\underbrace{w_{V,T}\mathbf{E}_{T}}%
_{\substack{=\mathbf{E}_{V,T}\\\text{(since (\ref{eq.prop.specht.EPQ.EPEQ.wE})
yields }\mathbf{E}_{V,T}=w_{V,T}\mathbf{E}_{T}\text{)}}%
}\label{pf.lem.bas.mur.triang.iso1lem.b.2}\\
&  =\underbrace{\left(  -1\right)  ^{u}}_{=\pm1}\mathbf{E}_{V,T}=\pm
\mathbf{E}_{V,T}.\nonumber
\end{align}
This proves Lemma \ref{lem.bas.mur.triang.iso1lem} \textbf{(b)}. \medskip

\textbf{(c)} The definition of $\MurpF_{\operatorname*{std},<T}%
^{-\operatorname*{Col}}$ says that%
\[
\MurpF_{\operatorname*{std},<T}^{-\operatorname*{Col}}=\operatorname*{span}%
\left\{  \nabla_{V,U}^{-\operatorname*{Col}}\ \mid\ \left(  \mu,U,V\right)
\in\operatorname*{SBT}\left(  n\right)  \text{ and }U<T\right\}
\]
(note that we used the notation \textquotedblleft$\left(  \lambda,U,V\right)
$\textquotedblright\ instead of \textquotedblleft$\left(  \mu,U,V\right)
$\textquotedblright\ in the original Definition \ref{def.bas.mur.submods}, but
we have now renamed it as \textquotedblleft$\left(  \mu,U,V\right)
$\textquotedblright\ because the letter $\lambda$ is already taken).

Now, let $\mathbf{a}\in\MurpF_{\operatorname*{std},<T}^{-\operatorname*{Col}}%
$. We must show that $\mathbf{a}\nabla_{\operatorname*{Row}T}=0$. This claim
depends $\mathbf{k}$-linearly on $\mathbf{a}$. Thus, by linearity, we can WLOG
assume that $\mathbf{a}$ is a vector of the form $\mathbf{a}=\nabla
_{V,U}^{-\operatorname*{Col}}$ for some $\left(  \mu,U,V\right)
\in\operatorname*{SBT}\left(  n\right)  $ satisfying $U<T$ (because
$\mathbf{a}\in\MurpF_{\operatorname*{std},<T}^{-\operatorname*{Col}%
}=\operatorname*{span}\left\{  \nabla_{V,U}^{-\operatorname*{Col}}%
\ \mid\ \left(  \mu,U,V\right)  \in\operatorname*{SBT}\left(  n\right)  \text{
and }U<T\right\}  $ shows that $\mathbf{a}$ is always a $\mathbf{k}$-linear
combination of vectors of this form). Assume this, and consider this $\left(
\mu,U,V\right)  $.

Now, from $\mathbf{a}=\nabla_{V,U}^{-\operatorname*{Col}}$, we obtain
$\mathbf{a}\nabla_{\operatorname*{Row}T}=\nabla_{V,U}^{-\operatorname*{Col}%
}\nabla_{\operatorname*{Row}T}=0$ (by Lemma \ref{lem.bas.mur.triang.iso1lem}
\textbf{(a)}). Thus, we have proved Lemma \ref{lem.bas.mur.triang.iso1lem}
\textbf{(c)}. \medskip

\textbf{(d)} The definition of $\MurpF_{\operatorname*{std},\leq
T}^{-\operatorname*{Col}}$ says that%
\[
\MurpF_{\operatorname*{std},\leq T}^{-\operatorname*{Col}}%
=\operatorname*{span}\left\{  \nabla_{V,U}^{-\operatorname*{Col}}%
\ \mid\ \left(  \mu,U,V\right)  \in\operatorname*{SBT}\left(  n\right)  \text{
and }U\leq T\right\}
\]
(note that we used the notation \textquotedblleft$\left(  \lambda,U,V\right)
$\textquotedblright\ instead of \textquotedblleft$\left(  \mu,U,V\right)
$\textquotedblright\ in the original Definition \ref{def.bas.mur.submods}, but
we have now renamed it as \textquotedblleft$\left(  \mu,U,V\right)
$\textquotedblright\ because the letter $\lambda$ is already taken).

Now, let $\mathbf{a}\in\MurpF_{\operatorname*{std},\leq T}%
^{-\operatorname*{Col}}$. We must show that $\mathbf{a}\nabla
_{\operatorname*{Row}T}\in\mathcal{A}\mathbf{E}_{T}$. This claim depends
$\mathbf{k}$-linearly on $\mathbf{a}$ (since $\mathcal{A}\mathbf{E}_{T}$ is a
$\mathbf{k}$-submodule of $\mathcal{A}$). Thus, by linearity, we can WLOG
assume that $\mathbf{a}$ is a vector of the form $\mathbf{a}=\nabla
_{V,U}^{-\operatorname*{Col}}$ for some $\left(  \mu,U,V\right)
\in\operatorname*{SBT}\left(  n\right)  $ satisfying $U\leq T$ (because
\newline$\mathbf{a}\in\MurpF_{\operatorname*{std},\leq T}%
^{-\operatorname*{Col}}=\operatorname*{span}\left\{  \nabla_{V,U}%
^{-\operatorname*{Col}}\ \mid\ \left(  \mu,U,V\right)  \in\operatorname*{SBT}%
\left(  n\right)  \text{ and }U\leq T\right\}  $ shows that $\mathbf{a}$ is
always a $\mathbf{k}$-linear combination of vectors of this form). Assume
this, and consider this $\left(  \mu,U,V\right)  $.

If $U<T$, then%
\begin{align*}
\underbrace{\mathbf{a}}_{=\nabla_{V,U}^{-\operatorname*{Col}}}\nabla
_{\operatorname*{Row}T}  &  =\nabla_{V,U}^{-\operatorname*{Col}}%
\nabla_{\operatorname*{Row}T}=0\ \ \ \ \ \ \ \ \ \ \left(  \text{by Lemma
\ref{lem.bas.mur.triang.iso1lem} \textbf{(a)}}\right) \\
&  \in\mathcal{A}\mathbf{E}_{T}\ \ \ \ \ \ \ \ \ \ \left(  \text{since
}\mathcal{A}\mathbf{E}_{T}\text{ is a }\mathbf{k}\text{-module}\right)  .
\end{align*}
Hence, $\mathbf{a}\nabla_{\operatorname*{Row}T}\in\mathcal{A}\mathbf{E}_{T}$
is proved in the case when $U<T$. Thus, for the rest of this proof, we WLOG
assume that we don't have $U<T$. Thus, we must have $U=T$ (since $U\leq T$).
Thus, $U=T\in\operatorname*{SYT}\left(  \lambda\right)  $ (since $T$ is a
standard $n$-tableau of shape $\lambda$). But $U\in\operatorname*{SYT}\left(
\mu\right)  $ (since $\left(  \mu,U,V\right)  \in\operatorname*{SBT}\left(
n\right)  $). Hence, the sets $\operatorname*{SYT}\left(  \lambda\right)  $
and $\operatorname*{SYT}\left(  \mu\right)  $ have at least one element in
common (namely, $U$). Thus, $\lambda=\mu$ (since the sets $\operatorname*{SYT}%
\left(  \kappa\right)  $ for different partitions $\kappa$ are disjoint).
Thus, $V$ is an $n$-tableau of shape $\lambda$ (since $V$ is an $n$-tableau of
shape $\mu$ (because $\left(  \mu,U,V\right)  \in\operatorname*{SBT}\left(
n\right)  $)).

Therefore, everything we wrote in our above proof of Lemma
\ref{lem.bas.mur.triang.iso1lem} \textbf{(b)} applies. In particular,
(\ref{pf.lem.bas.mur.triang.iso1lem.b.2}) holds, where $u\in S_{n}$ is defined
as in the proof just mentioned. Thus,%
\[
\nabla_{V,T}^{-\operatorname*{Col}}\nabla_{\operatorname*{Row}T}%
=\underbrace{\left(  -1\right)  ^{u}w_{V,T}}_{\in\mathcal{A}}\mathbf{E}_{T}%
\in\mathcal{A}\mathbf{E}_{T}.
\]
In other words, $\mathbf{a}\nabla_{\operatorname*{Row}T}\in\mathcal{A}%
\mathbf{E}_{T}$ (since $\mathbf{a}=\nabla_{V,U}^{-\operatorname*{Col}}$). This
completes the proof of Lemma \ref{lem.bas.mur.triang.iso1lem} \textbf{(d)}.
\medskip

\textbf{(e)} This follows from the definitions of $\MurpF_{\operatorname*{std}%
,\leq T}^{-\operatorname*{Col}}$ and $\MurpF_{\operatorname*{std}%
,<T}^{-\operatorname*{Col}}$ -- and from the fact that the only vectors
$\nabla_{V,U}^{-\operatorname*{Col}}$ in the definition of
$\MurpF_{\operatorname*{std},\leq T}^{-\operatorname*{Col}}$ that do not also
appear in the definition of $\MurpF_{\operatorname*{std},<T}%
^{-\operatorname*{Col}}$ are precisely the ones that satisfy $U=T$, that is,
the ones of the form $\nabla_{V,T}^{-\operatorname*{Col}}$ with $V\in
\operatorname*{SYT}\left(  \lambda\right)  $. See the appendix (Section
\ref{sec.details.bas.mur}) for details.
\end{proof}

\begin{proof}
[Proof of Theorem \ref{thm.bas.mur.triang.iso1}.]We know that $T$ is an
$n$-tableau of shape $Y\left(  \lambda\right)  $. Thus, $\mathcal{S}^{Y\left(
\lambda\right)  }\cong\mathcal{A}\mathbf{E}_{T}$ as left $\mathcal{A}$-modules
(by (\ref{eq.def.specht.ET.defs.SD=AET}), applied to $D=Y\left(
\lambda\right)  $). Hence, $\mathcal{A}\mathbf{E}_{T}\cong\mathcal{S}%
^{Y\left(  \lambda\right)  }=\mathcal{S}^{\lambda}$ (since $\mathcal{S}%
^{\lambda}$ is defined as $\mathcal{S}^{Y\left(  \lambda\right)  }$).

Theorem \ref{thm.bas.mur.triang.orth} \textbf{(c)} shows that
$\MurpF_{\operatorname*{std},<T}^{-\operatorname*{Col}}$ is a left
$\mathcal{A}$-submodule of $\mathcal{A}$. Theorem
\ref{thm.bas.mur.triang.orth} \textbf{(d)} shows that
$\MurpF_{\operatorname*{std},\leq T}^{-\operatorname*{Col}}$ is a left
$\mathcal{A}$-submodule of $\mathcal{A}$. Proposition
\ref{prop.bas.mur.subset} yields $\MurpF_{\operatorname*{std},<T}%
^{-\operatorname*{Col}}\subseteq\MurpF_{\operatorname*{std},\leq
T}^{-\operatorname*{Col}}$. Hence, the quotient $\MurpF_{\operatorname*{std}%
,\leq T}^{-\operatorname*{Col}}/\MurpF_{\operatorname*{std},<T}%
^{-\operatorname*{Col}}$ is well-defined, and is a left $\mathcal{A}$-module,
i.e., an $S_{n}$-representation.

Consider the $\mathbf{k}$-linear map%
\begin{align*}
\widehat{\Phi}:\MurpF_{\operatorname*{std},\leq T}^{-\operatorname*{Col}}  &
\rightarrow\mathcal{A}\mathbf{E}_{T},\\
\mathbf{a}  &  \mapsto\mathbf{a}\nabla_{\operatorname*{Row}T}%
\end{align*}
(this is well-defined, since Lemma \ref{lem.bas.mur.triang.iso1lem}
\textbf{(d)} shows that each $\mathbf{a}\in\MurpF_{\operatorname*{std},\leq
T}^{-\operatorname*{Col}}$ satisfies $\mathbf{a}\nabla_{\operatorname*{Row}%
T}\in\mathcal{A}\mathbf{E}_{T}$). This map $\widehat{\Phi}$ is easily seen to
be left $\mathcal{A}$-linear (because for any $\mathbf{b}\in\mathcal{A}$ and
any $\mathbf{a}\in\MurpF_{\operatorname*{std},\leq T}^{-\operatorname*{Col}}$,
we have $\widehat{\Phi}\left(  \mathbf{ba}\right)  =\mathbf{b}%
\underbrace{\mathbf{a}\nabla_{\operatorname*{Row}T}}_{=\widehat{\Phi}\left(
\mathbf{a}\right)  }=\mathbf{b}\widehat{\Phi}\left(  \mathbf{a}\right)  $). In
other words, $\widehat{\Phi}$ is a left $\mathcal{A}$-module morphism.
Moreover, this morphism $\widehat{\Phi}$ sends the left $\mathcal{A}%
$-submodule $\MurpF_{\operatorname*{std},<T}^{-\operatorname*{Col}}$ of
$\MurpF_{\operatorname*{std},\leq T}^{-\operatorname*{Col}}$ to $0$ (because
if $\mathbf{a}\in\MurpF_{\operatorname*{std},<T}^{-\operatorname*{Col}}$, then
$\widehat{\Phi}\left(  \mathbf{a}\right)  =\mathbf{a}\nabla
_{\operatorname*{Row}T}=0$ by Lemma \ref{lem.bas.mur.triang.iso1lem}
\textbf{(c)}). Hence, the universal property of quotient modules shows that
this morphism $\widehat{\Phi}$ factors through the quotient $\mathcal{A}%
$-module $\MurpF_{\operatorname*{std},\leq T}^{-\operatorname*{Col}%
}/\MurpF_{\operatorname*{std},<T}^{-\operatorname*{Col}}$. In other words,
$\widehat{\Phi}=\Phi\circ\pi$, where $\pi$ is the canonical projection from
$\MurpF_{\operatorname*{std},\leq T}^{-\operatorname*{Col}}$ to
$\MurpF_{\operatorname*{std},\leq T}^{-\operatorname*{Col}}%
/\MurpF_{\operatorname*{std},<T}^{-\operatorname*{Col}}$, and where
$\Phi:\MurpF_{\operatorname*{std},\leq T}^{-\operatorname*{Col}}%
/\MurpF_{\operatorname*{std},<T}^{-\operatorname*{Col}}\rightarrow
\mathcal{A}\mathbf{E}_{T}$ is a left $\mathcal{A}$-module morphism that is
explicitly given as follows:%
\begin{align}
\Phi:\MurpF_{\operatorname*{std},\leq T}^{-\operatorname*{Col}}%
/\MurpF_{\operatorname*{std},<T}^{-\operatorname*{Col}}  &  \rightarrow
\mathcal{A}\mathbf{E}_{T},\nonumber\\
\overline{\mathbf{a}}  &  \mapsto\mathbf{a}\nabla_{\operatorname*{Row}T}.
\label{pf.thm.bas.mur.triang.iso1.Phi}%
\end{align}
Consider this latter morphism $\Phi$. We shall now prove the following:

\begin{statement}
\textit{Claim 1:} The map $\Phi$ is surjective.
\end{statement}

\begin{statement}
\textit{Claim 2:} The map $\Phi$ is injective.
\end{statement}

\begin{proof}
[Proof of Claim 1.]Let $\mathbf{v}\in\mathcal{A}\mathbf{E}_{T}$. Thus,
$\mathbf{v}=\mathbf{bE}_{T}$ for some $\mathbf{b}\in\mathcal{A}$. Consider
this $\mathbf{b}$. We shall show that $\mathbf{v}$ lies in the image of $\Phi$.

We have $\nabla_{\operatorname*{Col}T}^{-}=\nabla_{T,T}^{-\operatorname*{Col}%
}$ (by (\ref{eq.prop.bas.mur.rsym.Col})). But $\nabla_{T,T}%
^{-\operatorname*{Col}}\in\MurpF_{\operatorname*{std},\leq T}%
^{-\operatorname*{Col}}$ (this follows from the definition of
$\MurpF_{\operatorname*{std},\leq T}^{-\operatorname*{Col}}$ pretty
easily\footnote{\textit{Proof.} The definition of $\MurpF_{\operatorname*{std}%
,\leq T}^{-\operatorname*{Col}}$ says that%
\[
\MurpF_{\operatorname*{std},\leq T}^{-\operatorname*{Col}}%
=\operatorname*{span}\left\{  \nabla_{V,U}^{-\operatorname*{Col}}%
\ \mid\ \left(  \mu,U,V\right)  \in\operatorname*{SBT}\left(  n\right)  \text{
and }U\leq T\right\}
\]
(note that we used the notation \textquotedblleft$\left(  \lambda,U,V\right)
$\textquotedblright\ instead of \textquotedblleft$\left(  \mu,U,V\right)
$\textquotedblright\ in the original Definition \ref{def.bas.mur.submods}, but
we have now renamed it as \textquotedblleft$\left(  \mu,U,V\right)
$\textquotedblright\ because the letter $\lambda$ is already taken). But
$\left(  \lambda,T,T\right)  \in\operatorname*{SBT}\left(  n\right)  $ (since
$T,T\in\operatorname*{SYT}\left(  \lambda\right)  $). Hence, $\left(
\lambda,T,T\right)  $ is a $\left(  \mu,U,V\right)  \in\operatorname*{SBT}%
\left(  n\right)  $ satisfying $U\leq T$ (since $T\leq T$). Therefore, the
element $\nabla_{T,T}^{-\operatorname*{Col}}$ has the form $\nabla
_{V,U}^{-\operatorname*{Col}}$ for some $\left(  \mu,U,V\right)
\in\operatorname*{SBT}\left(  n\right)  $ satisfying $U\leq T$ (namely, for
$\left(  \mu,U,V\right)  =\left(  \lambda,T,T\right)  $). That is,%
\begin{align*}
\nabla_{T,T}^{-\operatorname*{Col}}  &  \in\left\{  \nabla_{V,U}%
^{-\operatorname*{Col}}\ \mid\ \left(  \mu,U,V\right)  \in\operatorname*{SBT}%
\left(  n\right)  \text{ and }U\leq T\right\} \\
&  \subseteq\operatorname*{span}\left\{  \nabla_{V,U}^{-\operatorname*{Col}%
}\ \mid\ \left(  \mu,U,V\right)  \in\operatorname*{SBT}\left(  n\right)
\text{ and }U\leq T\right\}  =\MurpF_{\operatorname*{std},\leq T}%
^{-\operatorname*{Col}}.
\end{align*}
}). Therefore, $\nabla_{\operatorname*{Col}T}^{-}=\nabla_{T,T}%
^{-\operatorname*{Col}}\in\MurpF_{\operatorname*{std},\leq T}%
^{-\operatorname*{Col}}$, so that $\overline{\nabla_{\operatorname*{Col}T}%
^{-}}\in\MurpF_{\operatorname*{std},\leq T}^{-\operatorname*{Col}%
}/\MurpF_{\operatorname*{std},<T}^{-\operatorname*{Col}}$. Thus, $\Phi\left(
\overline{\nabla_{\operatorname*{Col}T}^{-}}\right)  $ is well-defined. The
explicit formula (\ref{pf.thm.bas.mur.triang.iso1.Phi}) for the map $\Phi$
shows that $\Phi\left(  \overline{\nabla_{\operatorname*{Col}T}^{-}}\right)
=\nabla_{\operatorname*{Col}T}^{-}\nabla_{\operatorname*{Row}T}=\mathbf{E}%
_{T}$ (since $\mathbf{E}_{T}$ is defined as $\nabla_{\operatorname*{Col}T}%
^{-}\nabla_{\operatorname*{Row}T}$). Now, since $\Phi$ is a left $\mathcal{A}%
$-module morphism, we have
\[
\Phi\left(  \mathbf{b}\overline{\nabla_{\operatorname*{Col}T}^{-}}\right)
=\mathbf{b}\underbrace{\Phi\left(  \overline{\nabla_{\operatorname*{Col}T}%
^{-}}\right)  }_{=\mathbf{E}_{T}}=\mathbf{bE}_{T}=\mathbf{v}%
\ \ \ \ \ \ \ \ \ \ \left(  \text{since }\mathbf{v}=\mathbf{bE}_{T}\right)  ,
\]
and thus $\mathbf{v}=\Phi\left(  \mathbf{b}\overline{\nabla
_{\operatorname*{Col}T}^{-}}\right)  $. This shows that $\mathbf{v}$ lies in
the image of $\Phi$.

Forget that we fixed $\mathbf{v}$. We thus have shown that each $\mathbf{v}%
\in\mathcal{A}\mathbf{E}_{T}$ lies in the image of $\Phi$. In other words,
$\Phi$ is surjective. This proves Claim 1.
\end{proof}

\begin{proof}
[Proof of Claim 2.]We must prove that $\Phi$ is injective. Since $\Phi$ is
$\mathbf{k}$-linear, it suffices to show that $\operatorname*{Ker}\Phi=0$
(because any $\mathbf{k}$-linear map whose kernel is $0$ is automatically injective).

So let $\mathbf{v}\in\operatorname*{Ker}\Phi$. We must show that
$\mathbf{v}=0$.

Clearly, $Y\left(  \lambda\right)  $ is a skew Young diagram, and we have
$\operatorname*{SYT}\left(  Y\left(  \lambda\right)  \right)
=\operatorname*{SYT}\left(  \lambda\right)  $ (since tableaux of shape
$Y\left(  \lambda\right)  $ are the same thing as tableaux of shape $\lambda
$). Thus, Proposition \ref{prop.specht.EPQ.basis-AET} (applied to $D=Y\left(
\lambda\right)  $) shows that the family $\left(  \mathbf{E}_{P,T}\right)
_{P\in\operatorname*{SYT}\left(  Y\left(  \lambda\right)  \right)  }$ is a
basis of the $\mathbf{k}$-module $\mathcal{A}\mathbf{E}_{T}$. Since
$\operatorname*{SYT}\left(  Y\left(  \lambda\right)  \right)
=\operatorname*{SYT}\left(  \lambda\right)  $, we can rewrite this as follows:
The family $\left(  \mathbf{E}_{P,T}\right)  _{P\in\operatorname*{SYT}\left(
\lambda\right)  }$ is a basis of the $\mathbf{k}$-module $\mathcal{A}%
\mathbf{E}_{T}$. In other words, the family $\left(  \mathbf{E}_{V,T}\right)
_{V\in\operatorname*{SYT}\left(  \lambda\right)  }$ is a basis of the
$\mathbf{k}$-module $\mathcal{A}\mathbf{E}_{T}$ (because this is the same
family as $\left(  \mathbf{E}_{P,T}\right)  _{P\in\operatorname*{SYT}\left(
\lambda\right)  }$; we just renamed the index $P$ as $V$). Hence, this family
is $\mathbf{k}$-linearly independent.

We have $\mathbf{v}\in\operatorname*{Ker}\Phi\subseteq
\MurpF_{\operatorname*{std},\leq T}^{-\operatorname*{Col}}%
/\MurpF_{\operatorname*{std},<T}^{-\operatorname*{Col}}$. Thus, $\mathbf{v}%
=\overline{\mathbf{c}}$ for some $\mathbf{c}\in\MurpF_{\operatorname*{std}%
,\leq T}^{-\operatorname*{Col}}$. Consider this $\mathbf{c}$. Then,%
\[
\mathbf{c}\in\MurpF_{\operatorname*{std},\leq T}^{-\operatorname*{Col}%
}=\MurpF_{\operatorname*{std},<T}^{-\operatorname*{Col}}+\operatorname*{span}%
\left\{  \nabla_{V,T}^{-\operatorname*{Col}}\ \mid\ V\in\operatorname*{SYT}%
\left(  \lambda\right)  \right\}
\]
(by Lemma \ref{lem.bas.mur.triang.iso1lem} \textbf{(e)}). In other words,
$\mathbf{c}$ can be written as $\mathbf{c}=\mathbf{a}+\mathbf{b}$ for some
$\mathbf{a}\in\MurpF_{\operatorname*{std},<T}^{-\operatorname*{Col}}$ and some
$\mathbf{b}\in\operatorname*{span}\left\{  \nabla_{V,T}^{-\operatorname*{Col}%
}\ \mid\ V\in\operatorname*{SYT}\left(  \lambda\right)  \right\}  $. Consider
these $\mathbf{a}$ and $\mathbf{b}$. Then, $\mathbf{c}=\mathbf{a}+\mathbf{b}$,
so that $\mathbf{c}-\mathbf{b}=\mathbf{a}\in\MurpF_{\operatorname*{std}%
,<T}^{-\operatorname*{Col}}$. Therefore, $\overline{\mathbf{c}}=\overline
{\mathbf{b}}$ (because two vectors in $\MurpF_{\operatorname*{std},\leq
T}^{-\operatorname*{Col}}$ whose difference belongs to
$\MurpF_{\operatorname*{std},<T}^{-\operatorname*{Col}}$ must have the same
projection onto the quotient $\MurpF_{\operatorname*{std},\leq T}%
^{-\operatorname*{Col}}/\MurpF_{\operatorname*{std},<T}^{-\operatorname*{Col}%
}$). Thus, $\mathbf{v}=\overline{\mathbf{c}}=\overline{\mathbf{b}}$, so that
$\overline{\mathbf{b}}=\mathbf{v}$.

But $\mathbf{b}\in\operatorname*{span}\left\{  \nabla_{V,T}%
^{-\operatorname*{Col}}\ \mid\ V\in\operatorname*{SYT}\left(  \lambda\right)
\right\}  $. Hence, we can write $\mathbf{b}$ as a $\mathbf{k}$-linear
combination%
\begin{equation}
\mathbf{b}=\sum_{V\in\operatorname*{SYT}\left(  \lambda\right)  }\omega
_{V}\nabla_{V,T}^{-\operatorname*{Col}}
\label{pf.thm.bas.mur.triang.iso1.c2.b=}%
\end{equation}
for some scalars $\omega_{V}\in\mathbf{k}$. Consider these $\omega_{V}$. But
$\overline{\mathbf{b}}=\mathbf{v}\in\operatorname*{Ker}\Phi$, so that
$\Phi\left(  \overline{\mathbf{b}}\right)  =0$. Hence,%
\begin{align*}
0  &  =\Phi\left(  \overline{\mathbf{b}}\right)  =\mathbf{b}\nabla
_{\operatorname*{Row}T}\ \ \ \ \ \ \ \ \ \ \left(  \text{by the explicit
formula (\ref{pf.thm.bas.mur.triang.iso1.Phi}) for the map }\Phi\right) \\
&  =\left(  \sum_{V\in\operatorname*{SYT}\left(  \lambda\right)  }\omega
_{V}\nabla_{V,T}^{-\operatorname*{Col}}\right)  \nabla_{\operatorname*{Row}%
T}\ \ \ \ \ \ \ \ \ \ \left(  \text{by (\ref{pf.thm.bas.mur.triang.iso1.c2.b=}%
)}\right) \\
&  =\sum_{V\in\operatorname*{SYT}\left(  \lambda\right)  }\omega
_{V}\underbrace{\nabla_{V,T}^{-\operatorname*{Col}}\nabla_{\operatorname*{Row}%
T}}_{\substack{=\pm\mathbf{E}_{V,T}\\\text{(by Lemma
\ref{lem.bas.mur.triang.iso1lem} \textbf{(b)})}}}=\sum_{V\in
\operatorname*{SYT}\left(  \lambda\right)  }\underbrace{\omega_{V}\left(
\pm\mathbf{E}_{V,T}\right)  }_{=\left(  \pm\omega_{V}\right)  \mathbf{E}%
_{V,T}}\\
&  =\sum_{V\in\operatorname*{SYT}\left(  \lambda\right)  }\left(  \pm
\omega_{V}\right)  \mathbf{E}_{V,T}.
\end{align*}
In other words,
\[
\sum_{V\in\operatorname*{SYT}\left(  \lambda\right)  }\left(  \pm\omega
_{V}\right)  \mathbf{E}_{V,T}=0\ \ \ \ \ \ \ \ \ \ \left(  \text{for an
appropriate choice of }\pm\text{ signs}\right)  .
\]

Since the family $\left(  \mathbf{E}_{V,T}\right)  _{V\in\operatorname*{SYT}%
\left(  \lambda\right)  }$ is $\mathbf{k}$-linearly independent, we thus
conclude that all the coefficients $\pm\omega_{V}$ in this sum are $0$. In
other words, each $V\in\operatorname*{SYT}\left(  \lambda\right)  $ satisfies
$\pm\omega_{V}=0$ and therefore $\omega_{V}=0$. Thus,
(\ref{pf.thm.bas.mur.triang.iso1.c2.b=}) simplifies to
\[
\mathbf{b}=\sum_{V\in\operatorname*{SYT}\left(  \lambda\right)  }%
\underbrace{\omega_{V}}_{=0}\nabla_{V,T}^{-\operatorname*{Col}}=\sum
_{V\in\operatorname*{SYT}\left(  \lambda\right)  }0\nabla_{V,T}%
^{-\operatorname*{Col}}=0.
\]
Now recall that $\mathbf{v}=\overline{\mathbf{b}}$. In view of $\mathbf{b}=0$,
we can rewrite this as $\mathbf{v}=\overline{0}=0$.

Thus we have shown that $\mathbf{v}=0$ for each $\mathbf{v}\in
\operatorname*{Ker}\Phi$. In other words, $\operatorname*{Ker}\Phi=0$. As we
said, this completes the proof of Claim 2.
\end{proof}

Now, the map $\Phi$ is injective (by Claim 2) and surjective (by Claim 1).
Hence, $\Phi$ is bijective, i.e., invertible. Thus, $\Phi$ is a left
$\mathcal{A}$-module isomorphism (since $\Phi$ is a left $\mathcal{A}$-module
morphism). In other words, the map%
\begin{align*}
\MurpF_{\operatorname*{std},\leq T}^{-\operatorname*{Col}}%
/\MurpF_{\operatorname*{std},<T}^{-\operatorname*{Col}}  &  \rightarrow
\mathcal{A}\mathbf{E}_{T},\\
\overline{\mathbf{a}}  &  \mapsto\mathbf{a}\nabla_{\operatorname*{Row}T}%
\end{align*}
is a left $\mathcal{A}$-module isomorphism (since this map is precisely $\Phi
$, because of (\ref{pf.thm.bas.mur.triang.iso1.Phi})). Therefore,
$\MurpF_{\operatorname*{std},\leq T}^{-\operatorname*{Col}}%
/\MurpF_{\operatorname*{std},<T}^{-\operatorname*{Col}}\cong\mathcal{A}%
\mathbf{E}_{T}$ as left $\mathcal{A}$-modules. Combining this with
$\mathcal{A}\mathbf{E}_{T}\cong\mathcal{S}^{\lambda}$, we obtain
$\MurpF_{\operatorname*{std},\leq T}^{-\operatorname*{Col}}%
/\MurpF_{\operatorname*{std},<T}^{-\operatorname*{Col}}\cong\mathcal{A}%
\mathbf{E}_{T}\cong\mathcal{S}^{\lambda}$ as left $\mathcal{A}$-modules, i.e.,
as $S_{n}$-representations. Thus, Theorem \ref{thm.bas.mur.triang.iso1} is proved.
\end{proof}

The column-Murphy sums have a property analogous to Theorem
\ref{thm.bas.mur.triang.iso1}, which we leave as an exercise:

\begin{exercise}
\fbox{3} Let $\lambda$ be a partition of $n$. Let $T$ be any standard
$n$-tableau of shape $\lambda$. Then, prove that
\[
\MurpF_{\operatorname*{std},\geq T}^{\operatorname*{Row}}%
/\MurpF_{\operatorname*{std},>T}^{\operatorname*{Row}}\cong\mathcal{A}%
\mathbf{F}_{T}\cong\left(  \mathcal{S}^{\lambda}\right)  ^{\ast}%
\ \ \ \ \ \ \ \ \ \ \text{as }S_{n}\text{-representations}%
\]
(where $\mathbf{F}_{T}$ is as defined in Proposition
\ref{prop.specht.FT.basics}, and where $\left(  \mathcal{S}^{\lambda}\right)
^{\ast}$ is the dual of the Specht module $\mathcal{S}^{\lambda}$). More
concretely, prove that the map%
\begin{align*}
\MurpF_{\operatorname*{std},\geq T}^{\operatorname*{Row}}%
/\MurpF_{\operatorname*{std},>T}^{\operatorname*{Row}}  &  \rightarrow
\mathcal{A}\mathbf{F}_{T},\\
\overline{\mathbf{a}}  &  \mapsto\mathbf{a}\nabla_{\operatorname*{Col}T}^{-}%
\end{align*}
(where $\overline{\mathbf{a}}$ denotes the projection of a vector
$\mathbf{a}\in\MurpF_{\operatorname*{std},\geq T}^{\operatorname*{Row}}$ on
the quotient $\MurpF_{\operatorname*{std},\geq T}^{\operatorname*{Row}%
}/\MurpF_{\operatorname*{std},>T}^{\operatorname*{Row}}$) is a left
$\mathcal{A}$-module isomorphism. \medskip

[\textbf{Hint:} From Lemma \ref{lem.sign-twist.AFT} and Theorem
\ref{thm.dual.specht-skew}, we have $\mathcal{A}\mathbf{F}_{T}\cong\left(
\mathcal{S}^{\mathbf{r}\left(  D\right)  }\right)  ^{\operatorname*{sign}%
}\cong\left(  \mathcal{S}^{D}\right)  ^{\ast}$ for $D=Y\left(  \lambda\right)
$. Also recall (\ref{eq.lem.specht.EPQ.straighten.ETA}).]
\end{exercise}

Theorem \ref{thm.bas.mur.triang.iso1} was stated in terms of a single tableau
$T$, but we can easily restate it in terms of the entire filtration:

\begin{theorem}
\label{thm.bas.mur.filtC-subq}Let $T_{1},T_{2},\ldots,T_{m}$ be all the
standard tableaux in the set $\bigcup\limits_{\kappa\vdash n}%
\operatorname*{SYT}\left(  \kappa\right)  $, listed in increasing order (with
respect to the total order defined in Definition \ref{def.bas.tord.syt.tord}).
Recall from Theorem \ref{thm.bas.mur.filt} \textbf{(f)} that $\mathcal{A}$ has
a filtration%
\begin{equation}
\left(  \MurpF_{0}^{-\operatorname*{Col}},\MurpF_{1}^{-\operatorname*{Col}%
},\MurpF_{2}^{-\operatorname*{Col}},\ldots,\MurpF_{m-1}^{-\operatorname*{Col}%
},\MurpF_{m}^{-\operatorname*{Col}}\right)
\label{eq.thm.bas.mur.filtC-subq.filt}%
\end{equation}
by left $\mathcal{A}$-submodules.

Let $i\in\left[  m\right]  $. Let $\lambda^{\left(  i\right)  }$ be the
partition of $n$ that satisfies $T_{i}\in\operatorname*{SYT}\left(
\lambda^{\left(  i\right)  }\right)  $. Then, the subquotient $\MurpF_{i}%
^{-\operatorname*{Col}}/\MurpF_{i-1}^{-\operatorname*{Col}}$ of the filtration
(\ref{eq.thm.bas.mur.filtC-subq.filt}) is isomorphic (as a left $\mathcal{A}%
$-module) to the Specht module $\mathcal{S}^{\lambda^{\left(  i\right)  }}$.
\end{theorem}

\begin{proof}
Here is a sketch (see Section \ref{sec.details.bas.mur} for more details). Set
$T:=T_{i}$ and $\lambda:=\lambda^{\left(  i\right)  }$. Thus, $T\in
\operatorname*{SYT}\left(  \lambda\right)  $. Hence, Theorem
\ref{thm.bas.mur.triang.iso1} yields that $\MurpF_{\operatorname*{std},\leq
T}^{-\operatorname*{Col}}/\MurpF_{\operatorname*{std},<T}%
^{-\operatorname*{Col}}\cong\mathcal{S}^{\lambda}$ as $S_{n}$-representations,
i.e., as left $\mathcal{A}$-modules.

But $\MurpF_{i}^{-\operatorname*{Col}}=\MurpF_{\operatorname*{std},\leq
T}^{-\operatorname*{Col}}$ (by the definition of $\MurpF_{i}%
^{-\operatorname*{Col}}$) and $\MurpF_{i-1}^{-\operatorname*{Col}%
}=\MurpF_{\operatorname*{std},<T}^{-\operatorname*{Col}}$ (by Theorem
\ref{thm.bas.mur.filt} \textbf{(b)}, applied to $i-1$ instead of $i$). Hence,%
\[
\MurpF_{i}^{-\operatorname*{Col}}/\MurpF_{i-1}^{-\operatorname*{Col}%
}=\MurpF_{\operatorname*{std},\leq T}^{-\operatorname*{Col}}%
/\MurpF_{\operatorname*{std},<T}^{-\operatorname*{Col}}\cong\mathcal{S}%
^{\lambda}=\mathcal{S}^{\lambda^{\left(  i\right)  }}%
\ \ \ \ \ \ \ \ \ \ \left(  \text{since }\lambda=\lambda^{\left(  i\right)
}\right)
\]
as left $\mathcal{A}$-modules. This proves Theorem
\ref{thm.bas.mur.filtC-subq}.
\end{proof}

Theorem \ref{thm.bas.mur.filtC-subq} shows that the filtration
(\ref{eq.thm.bas.mur.filtC-subq.filt}) is a so-called \emph{Specht filtration}
-- that is, a filtration whose each subquotient is isomorphic to a Specht
module of partition shape (i.e., to a Specht module $\mathcal{S}^{\lambda}$
with $\lambda\vdash n$). When $\mathbf{k}$ is a field of characteristic $0$,
these subquotients are therefore irreducible as left $\mathcal{A}$-modules
(since Corollary \ref{cor.spechtmod.irred} shows that the Specht modules
$\mathcal{S}^{\lambda}$ are irreducible), so that the filtration is a
composition series. Thus we have proved the following:

\begin{corollary}
\label{cor.bas.mur.filtC-comps}Assume that $\mathbf{k}$ is a field of
characteristic $0$. Then, the filtration (\ref{eq.thm.bas.mur.filtC-subq.filt}%
) is a composition series.
\end{corollary}

\begin{exercise}
\fbox{1} Prove an analogous property of the filtration $\left(  \MurpF_{m}%
^{\operatorname*{Row}},\MurpF_{m-1}^{\operatorname*{Row}},\ldots
,\MurpF_{2}^{\operatorname*{Row}},\MurpF_{1}^{\operatorname*{Row}}%
,\MurpF_{0}^{\operatorname*{Row}}\right)  $.
\end{exercise}

Corollary \ref{cor.bas.mur.filtC-comps} gives us a new proof of Corollary
\ref{cor.spechtmod.complete}:

\begin{proof}
[Another proof of Corollary \ref{cor.spechtmod.complete}.]As we just said
(Corollary \ref{cor.bas.mur.filtC-comps}), the filtration
(\ref{eq.thm.bas.mur.filtC-subq.filt}) of $\mathcal{A}$ is a composition
series. Hence, Theorem \ref{thm.mod.filt.irrep-in-filt} (applied to the ring
$\mathcal{A}$, the number $m$ and the composition series $\left(
\MurpF_{0}^{-\operatorname*{Col}},\MurpF_{1}^{-\operatorname*{Col}}%
,\ldots,\MurpF_{m}^{-\operatorname*{Col}}\right)  $ instead of the ring $R$,
the number $k$, and the composition series $\left(  W_{0},W_{1},\ldots
,W_{k}\right)  $) shows that each irreducible left $\mathcal{A}$-module is
isomorphic to $\MurpF_{i}^{-\operatorname*{Col}}/\MurpF_{i-1}%
^{-\operatorname*{Col}}$ for some $i\in\left[  m\right]  $. But Theorem
\ref{thm.bas.mur.filtC-subq} shows that each such $\MurpF_{i}%
^{-\operatorname*{Col}}/\MurpF_{i-1}^{-\operatorname*{Col}}$ is isomorphic to
$\mathcal{S}^{\lambda}$ for some partition $\lambda$ of $n$. Combining these
two isomorphisms, we conclude that each irreducible left $\mathcal{A}$-module
is isomorphic to $\mathcal{S}^{\lambda}$ for some partition $\lambda$ of $n$.
This proves Corollary \ref{cor.spechtmod.complete} again.
\end{proof}

\begin{remark}
The filtration (\ref{eq.thm.bas.mur.filtC-subq.filt}) is not the only one that
has the nice properties we proved (such as Theorem
\ref{thm.bas.mur.filtC-subq} and Corollary \ref{cor.bas.mur.filtC-comps}). One
easy way to construct other such filtrations is by replacing our total order
on $\bigcup\limits_{\kappa\vdash n}\operatorname*{SYT}\left(  \kappa\right)  $
by another that still satisfies Lemma \ref{lem.bas.tord.syt.czs}. Such
alternative orders are not unusual, although not every total order would work.
\end{remark}

\subsubsection{\label{subsec.bas.mur.twosid}Two-sided ideals from the Murphy
bases}

In Theorem \ref{thm.bas.mur.triang.orth}, we have seen how certain subfamilies
of the Murphy bases span left $\mathcal{A}$-submodules of $\mathcal{A}$ --
that is, left ideals of $\mathcal{A}$. We shall now discuss a similar way to
construct two-sided ideals of $\mathcal{A}$.

We will require the notion of an \textquotedblleft$n$%
-bitableau\textquotedblright\ (analogous to Definition \ref{def.bas.not}
\textbf{(e)} and Definition \ref{def.bas.mur.HSBT}, but without requiring any
of the tableaux to be standard):

\begin{definition}
\label{def.bas.mur.BT}An $n$\emph{-bitableau} shall mean a triple $\left(
\lambda,U,V\right)  $, where $\lambda$ is a partition of $n$ and where $U$ and
$V$ are two $n$-tableaux of shape $\lambda$.

We let $\operatorname*{BT}\left(  n\right)  $ be the set of all $n$-bitableaux.
\end{definition}

Obviously, $\operatorname*{SBT}\left(  n\right)  \subseteq\operatorname*{BT}%
\left(  n\right)  $, since any standard $n$-bitableau is an $n$-bitableau.

\begin{definition}
\label{def.bas.mur.submods-len}Let $k\in\mathbb{N}$. Then, we define the
following four $\mathbf{k}$-submodules of $\mathcal{A}$:%
\begin{align*}
\MurpF_{\operatorname*{std},\operatorname*{len}\leq k}^{\operatorname*{Row}}
&  :=\operatorname*{span}\left\{  \nabla_{V,U}^{\operatorname*{Row}}%
\ \mid\ \left(  \lambda,U,V\right)  \in\operatorname*{SBT}\left(  n\right)
\text{ and }\ell\left(  \lambda\right)  \leq k\right\}  ;\\
\MurpF_{\operatorname*{std},\operatorname*{len}>k}^{-\operatorname*{Col}}  &
:=\operatorname*{span}\left\{  \nabla_{V,U}^{-\operatorname*{Col}}%
\ \mid\ \left(  \lambda,U,V\right)  \in\operatorname*{SBT}\left(  n\right)
\text{ and }\ell\left(  \lambda\right)  >k\right\}  ;\\
\MurpF_{\operatorname*{all},\operatorname*{len}\leq k}^{\operatorname*{Row}}
&  :=\operatorname*{span}\left\{  \nabla_{V,U}^{\operatorname*{Row}}%
\ \mid\ \left(  \lambda,U,V\right)  \in\operatorname*{BT}\left(  n\right)
\text{ and }\ell\left(  \lambda\right)  \leq k\right\}  ;\\
\MurpF_{\operatorname*{all},\operatorname*{len}>k}^{-\operatorname*{Col}}  &
:=\operatorname*{span}\left\{  \nabla_{V,U}^{-\operatorname*{Col}}%
\ \mid\ \left(  \lambda,U,V\right)  \in\operatorname*{BT}\left(  n\right)
\text{ and }\ell\left(  \lambda\right)  >k\right\}  .
\end{align*}

\end{definition}

\begin{example}
For $n=3$, we have%
\begin{align*}
\MurpF_{\operatorname*{std},\operatorname*{len}\leq1}^{\operatorname*{Row}}
&  =\operatorname*{span}\left\{  \nabla_{123,\ 123}^{\operatorname*{Row}%
}\right\}  ;\\
\MurpF_{\operatorname*{std},\operatorname*{len}\leq2}^{\operatorname*{Row}}
&  =\operatorname*{span}\left\{  \nabla_{13\backslash\backslash
2,\ 13\backslash\backslash2}^{\operatorname*{Row}},\ \nabla_{13\backslash
\backslash2,\ 12\backslash\backslash3}^{\operatorname*{Row}},\ \nabla
_{12\backslash\backslash3,\ 13\backslash\backslash2}^{\operatorname*{Row}%
},\right. \\
&  \ \ \ \ \ \ \ \ \ \ \ \ \ \ \ \ \ \ \ \ \left.  \ \nabla_{12\backslash
\backslash3,\ 12\backslash\backslash3}^{\operatorname*{Row}},\ \nabla
_{123,\ 123}^{\operatorname*{Row}}\right\}  ;\\
\MurpF_{\operatorname*{std},\operatorname*{len}\leq k}^{\operatorname*{Row}}
&  =\operatorname*{span}\left\{  \nabla_{1\backslash\backslash2\backslash
\backslash3,\ 1\backslash\backslash2\backslash\backslash3}%
^{\operatorname*{Row}},\ \nabla_{13\backslash\backslash2,\ 13\backslash
\backslash2}^{\operatorname*{Row}},\ \nabla_{13\backslash\backslash
2,\ 12\backslash\backslash3}^{\operatorname*{Row}},\ \nabla_{12\backslash
\backslash3,\ 13\backslash\backslash2}^{\operatorname*{Row}},\right. \\
&  \ \ \ \ \ \ \ \ \ \ \ \ \ \ \ \ \ \ \ \ \left.  \ \nabla_{12\backslash
\backslash3,\ 12\backslash\backslash3}^{\operatorname*{Row}},\ \nabla
_{123,\ 123}^{\operatorname*{Row}}\right\}  \ \ \ \ \ \ \ \ \ \ \text{for all
}k\geq3
\end{align*}
(and, of course, $\MurpF_{\operatorname*{std},\operatorname*{len}\leq
0}^{\operatorname*{Row}}=\operatorname*{span}\varnothing=0$). Likewise, the
reader can compute $\MurpF_{\operatorname*{std},\operatorname*{len}%
>k}^{-\operatorname*{Col}}$ for $n=3$.
\end{example}

We now claim the following theorem, which resembles Theorem
\ref{thm.bas.mur.triang.orth}:

\begin{theorem}
\label{thm.bas.mur.twosid.orth-len}Let $k\in\mathbb{N}$. Then: \medskip

\textbf{(a)} The set $\MurpF_{\operatorname*{std},\operatorname*{len}\leq
k}^{\operatorname*{Row}}$ is a two-sided ideal of $\mathcal{A}$ and satisfies%
\[
\MurpF_{\operatorname*{std},\operatorname*{len}\leq k}^{\operatorname*{Row}%
}=\MurpF_{\operatorname*{all},\operatorname*{len}\leq k}^{\operatorname*{Row}%
}=\left(  \MurpF_{\operatorname*{std},\operatorname*{len}>k}%
^{-\operatorname*{Col}}\right)  ^{\perp}.
\]

\textbf{(b)} The set $\MurpF_{\operatorname*{std},\operatorname*{len}%
>k}^{-\operatorname*{Col}}$ is a two-sided ideal of $\mathcal{A}$ and
satisfies
\[
\MurpF_{\operatorname*{std},\operatorname*{len}>k}^{-\operatorname*{Col}%
}=\MurpF_{\operatorname*{all},\operatorname*{len}>k}^{-\operatorname*{Col}%
}=\left(  \MurpF_{\operatorname*{std},\operatorname*{len}\leq k}%
^{\operatorname*{Row}}\right)  ^{\perp}.
\]

\end{theorem}

Rather than prove this directly, we first generalize it by replacing the
conditions \textquotedblleft$\ell\left(  \lambda\right)  \leq k$%
\textquotedblright\ and \textquotedblleft$\ell\left(  \lambda\right)
>k$\textquotedblright\ by more flexible requirements (namely,
\textquotedblleft$\lambda\in X$\textquotedblright\ and \textquotedblleft%
$\lambda\in Y$\textquotedblright\ for a certain class of sets $X$ and $Y$):

\begin{theorem}
\label{thm.bas.mur.twosid.orth-XY}Let $\operatorname*{Par}\left(  n\right)  $
be the set of all partitions of $n$. Let $X$ and $Y$ be two disjoint subsets
of $\operatorname*{Par}\left(  n\right)  $ such that $X\cup
Y=\operatorname*{Par}\left(  n\right)  $. Assume that the following condition holds:

\begin{statement}
\textit{Non-dominance condition:} For every $\lambda\in X$ and every $\mu\in
Y$, there exists some $i\geq1$ satisfying $\lambda_{1}+\lambda_{2}%
+\cdots+\lambda_{i}>\mu_{1}+\mu_{2}+\cdots+\mu_{i}$.
\end{statement}

Define the following four $\mathbf{k}$-submodules of $\mathcal{A}$:%
\begin{align*}
\MurpF_{\operatorname*{std},X}^{\operatorname*{Row}}  &
:=\operatorname*{span}\left\{  \nabla_{V,U}^{\operatorname*{Row}}%
\ \mid\ \left(  \lambda,U,V\right)  \in\operatorname*{SBT}\left(  n\right)
\text{ and }\lambda\in X\right\}  ;\\
\MurpF_{\operatorname*{std},Y}^{-\operatorname*{Col}}  &
:=\operatorname*{span}\left\{  \nabla_{V,U}^{-\operatorname*{Col}}%
\ \mid\ \left(  \lambda,U,V\right)  \in\operatorname*{SBT}\left(  n\right)
\text{ and }\lambda\in Y\right\}  ;\\
\MurpF_{\operatorname*{all},X}^{\operatorname*{Row}}  &
:=\operatorname*{span}\left\{  \nabla_{V,U}^{\operatorname*{Row}}%
\ \mid\ \left(  \lambda,U,V\right)  \in\operatorname*{BT}\left(  n\right)
\text{ and }\lambda\in X\right\}  ;\\
\MurpF_{\operatorname*{all},Y}^{-\operatorname*{Col}}  &
:=\operatorname*{span}\left\{  \nabla_{V,U}^{-\operatorname*{Col}}%
\ \mid\ \left(  \lambda,U,V\right)  \in\operatorname*{BT}\left(  n\right)
\text{ and }\lambda\in Y\right\}  .
\end{align*}
Then: \medskip

\textbf{(a)} The set $\MurpF_{\operatorname*{std},X}^{\operatorname*{Row}}$ is
a two-sided ideal of $\mathcal{A}$ and satisfies%
\[
\MurpF_{\operatorname*{std},X}^{\operatorname*{Row}}%
=\MurpF_{\operatorname*{all},X}^{\operatorname*{Row}}=\left(
\MurpF_{\operatorname*{std},Y}^{-\operatorname*{Col}}\right)  ^{\perp}.
\]

\textbf{(b)} The set $\MurpF_{\operatorname*{std},Y}^{-\operatorname*{Col}}$
is a two-sided ideal of $\mathcal{A}$ and satisfies
\[
\MurpF_{\operatorname*{std},Y}^{-\operatorname*{Col}}%
=\MurpF_{\operatorname*{all},Y}^{-\operatorname*{Col}}=\left(
\MurpF_{\operatorname*{std},X}^{\operatorname*{Row}}\right)  ^{\perp}.
\]

\end{theorem}

Before we prove this theorem (in a manner fairly similar to our proof of
Theorem \ref{thm.bas.mur.triang.orth}), let us observe how Theorem
\ref{thm.bas.mur.twosid.orth-len} follows from it:

\begin{proof}
[Proof of Theorem \ref{thm.bas.mur.twosid.orth-len} using Theorem
\ref{thm.bas.mur.twosid.orth-XY}.]This is a sketch; a detailed version of this
proof can be found in the appendix (Section \ref{sec.details.bas.mur}).

Let $\operatorname*{Par}\left(  n\right)  $ be as in Theorem
\ref{thm.bas.mur.twosid.orth-XY}. Define the two subsets
\begin{align*}
X  &  :=\left\{  \lambda\in\operatorname*{Par}\left(  n\right)  \ \mid
\ \ell\left(  \lambda\right)  \leq k\right\}  \ \ \ \ \ \ \ \ \ \ \text{and}\\
Y  &  :=\left\{  \lambda\in\operatorname*{Par}\left(  n\right)  \ \mid
\ \ell\left(  \lambda\right)  >k\right\}
\end{align*}
of $\operatorname*{Par}\left(  n\right)  $. It is easy to see that for every
$\lambda\in X$ and every $\mu\in Y$, there exists some $i\geq1$ satisfying
$\lambda_{1}+\lambda_{2}+\cdots+\lambda_{i}>\mu_{1}+\mu_{2}+\cdots+\mu_{i}$
(namely, we can pick $i=\ell\left(  \lambda\right)  $). Thus, applying Theorem
\ref{thm.bas.mur.twosid.orth-XY}, we quickly obtain the claim of Theorem
\ref{thm.bas.mur.twosid.orth-len}.
\end{proof}

We now approach the proof of Theorem \ref{thm.bas.mur.twosid.orth-XY}. As we
said, this is similar to the proof of Theorem \ref{thm.bas.mur.triang.orth};
in particular, it relies on analogues of certain lemmas used in the latter
proof. The first one is an analogue to Lemma \ref{lem.bas.mur.orth}:

\begin{lemma}
\label{lem.bas.mur.twosid.orth}Let $\operatorname*{Par}\left(  n\right)  $,
$X$ and $Y$ be as in Theorem \ref{thm.bas.mur.twosid.orth-XY}. Then: \medskip

\textbf{(a)} If $\left(  \lambda,A,B\right)  \in\operatorname*{BT}\left(
n\right)  $ and $\left(  \mu,C,D\right)  \in\operatorname*{BT}\left(
n\right)  $ are such that $\lambda\in X$ and $\mu\in Y$, then $\left\langle
\nabla_{D,C}^{-\operatorname*{Col}},\ \nabla_{B,A}^{\operatorname*{Row}%
}\right\rangle =0$. \medskip

\textbf{(b)} We have $\left\langle \mathbf{a},\mathbf{b}\right\rangle =0$ for
any $\mathbf{a}\in\MurpF_{\operatorname*{all},Y}^{-\operatorname*{Col}}$ and
$\mathbf{b}\in\MurpF_{\operatorname*{all},X}^{\operatorname*{Row}}$.
\end{lemma}

\begin{proof}
\textbf{(a)} Let $\left(  \lambda,A,B\right)  \in\operatorname*{BT}\left(
n\right)  $ and $\left(  \mu,C,D\right)  \in\operatorname*{BT}\left(
n\right)  $ be such that $\lambda\in X$ and $\mu\in Y$.

Both $A$ and $B$ are $n$-tableaux of shape $Y\left(  \lambda\right)  $ (since
$\left(  \lambda,A,B\right)  \in\operatorname*{BT}\left(  n\right)  $). Hence,
Proposition \ref{prop.tableau.Sn-act.1} \textbf{(b)} (applied to $Y\left(
\lambda\right)  $ and $B$ instead of $D$ and $T$) shows that any $n$-tableau
of shape $Y\left(  \lambda\right)  $ can be written as $w\rightharpoonup B$
for some $w\in S_{n}$. Hence, in particular, the $n$-tableau $A$ can be
written in this way. In other words, $A=w\rightharpoonup B$ for some $w\in
S_{n}$. Let us denote this $w$ by $u$. Thus, $A=u\rightharpoonup B=uB$.
Applying (\ref{eq.lem.bas.mur.perm.R}) to $P=A$ and $Q=B$, we thus obtain%
\[
\nabla_{A,B}^{\operatorname*{Row}}=\nabla_{\operatorname*{Row}A}%
u=u\nabla_{\operatorname*{Row}B}.
\]

Both $C$ and $D$ are $n$-tableaux of shape $Y\left(  \mu\right)  $ (since
$\left(  \mu,C,D\right)  \in\operatorname*{BT}\left(  n\right)  $). Hence,
Proposition \ref{prop.tableau.Sn-act.1} \textbf{(b)} (applied to $Y\left(
\mu\right)  $ and $C$ instead of $D$ and $T$) shows that any $n$-tableau of
shape $Y\left(  \mu\right)  $ can be written as $w\rightharpoonup C$ for some
$w\in S_{n}$. Hence, in particular, the $n$-tableau $D$ can be written in this
way. In other words, $D=w\rightharpoonup C$ for some $w\in S_{n}$. Let us
denote this $w$ by $v$. Thus, $D=v\rightharpoonup C=vC$. Applying
(\ref{eq.lem.bas.mur.perm.C}) to $D$, $C$ and $v$ instead of $P$, $Q$ and $u$,
we thus obtain%
\[
\nabla_{D,C}^{-\operatorname*{Col}}=\left(  -1\right)  ^{v}\nabla
_{\operatorname*{Col}D}^{-}v=\left(  -1\right)  ^{v}v\nabla
_{\operatorname*{Col}C}^{-}.
\]

Furthermore, $S\left(  \nabla_{B,A}^{\operatorname*{Row}}\right)
=\nabla_{A,B}^{\operatorname*{Row}}$ (by (\ref{eq.prop.bas.mur.S.Row})).

Now, recall that $\lambda\in X$ and $\mu\in Y$. Hence, there exists some
$i\geq1$ satisfying $\lambda_{1}+\lambda_{2}+\cdots+\lambda_{i}>\mu_{1}%
+\mu_{2}+\cdots+\mu_{i}$ (by the non-dominance condition assumed in Theorem
\ref{thm.bas.mur.twosid.orth-XY}). In other words, $\lambda_{1}+\lambda
_{2}+\cdots+\lambda_{i}>\mu_{1}+\mu_{2}+\cdots+\mu_{i}$ for some $i\geq1$.
Hence, Theorem \ref{thm.specht.sandw0} \textbf{(a)} (applied to $S=A$ and
$T=C$ and $w=\operatorname*{id}$) yields $\nabla_{\operatorname*{Col}C}%
^{-}\operatorname*{id}\nabla_{\operatorname*{Row}A}=0$. In other words,
$\nabla_{\operatorname*{Col}C}^{-}\nabla_{\operatorname*{Row}A}=0$ (since
$\nabla_{\operatorname*{Col}C}^{-}\underbrace{\operatorname*{id}%
\nabla_{\operatorname*{Row}A}}_{=\nabla_{\operatorname*{Row}A}}=\nabla
_{\operatorname*{Col}C}^{-}\nabla_{\operatorname*{Row}A}$). Now,%
\[
\underbrace{\nabla_{D,C}^{-\operatorname*{Col}}}_{=\left(  -1\right)
^{v}v\nabla_{\operatorname*{Col}C}^{-}}\ \ \underbrace{S\left(  \nabla
_{B,A}^{\operatorname*{Row}}\right)  }_{\substack{=\nabla_{A,B}%
^{\operatorname*{Row}}\\=\nabla_{\operatorname*{Row}A}u}}=\left(  -1\right)
^{v}v\underbrace{\nabla_{\operatorname*{Col}C}^{-}\nabla_{\operatorname*{Row}%
A}}_{=0}u=0.
\]

Finally, (\ref{eq.prop.bas.bilin.sym-etc.4}) yields $\left\langle \nabla
_{D,C}^{-\operatorname*{Col}},\ \nabla_{B,A}^{\operatorname*{Row}%
}\right\rangle =\left[  \operatorname*{id}\right]  \left(  \underbrace{\nabla
_{D,C}^{-\operatorname*{Col}}S\left(  \nabla_{B,A}^{\operatorname*{Row}%
}\right)  }_{=0}\right)  =\left[  \operatorname*{id}\right]  0=0$. This proves
Lemma \ref{lem.bas.mur.twosid.orth} \textbf{(a)}. \medskip

\textbf{(b)} Let $\mathbf{a}\in\MurpF_{\operatorname*{all},Y}%
^{-\operatorname*{Col}}$ and $\mathbf{b}\in\MurpF_{\operatorname*{all}%
,X}^{\operatorname*{Row}}$. We must prove that $\left\langle \mathbf{a}%
,\mathbf{b}\right\rangle =0$.

This claim depends $\mathbf{k}$-linearly on each of $\mathbf{b}$ and
$\mathbf{a}$ (because $\left\langle \cdot,\cdot\right\rangle $ is a bilinear
form). Hence, we can WLOG assume

\begin{itemize}
\item that $\mathbf{b}$ is one of the vectors $\nabla_{V,U}%
^{\operatorname*{Row}}$ with $\left(  \lambda,U,V\right)  \in
\operatorname*{BT}\left(  n\right)  $ satisfying $\lambda\in X$ (since
$\mathbf{b}$ belongs to $\MurpF_{\operatorname*{all},X}^{\operatorname*{Row}}%
$, which is defined to be the span of these vectors), and

\item that $\mathbf{a}$ is one of the vectors $\nabla_{V,U}%
^{-\operatorname*{Col}}$ with $\left(  \lambda,U,V\right)  \in
\operatorname*{BT}\left(  n\right)  $ satisfying $\lambda\in Y$ (since
$\mathbf{a}$ belongs to $\MurpF_{\operatorname*{all},Y}^{-\operatorname*{Col}%
}$, which is defined to be the span of these vectors).
\end{itemize}

\noindent Assume this. Thus,
\[
\mathbf{b}=\nabla_{B,A}^{\operatorname*{Row}}\ \ \ \ \ \ \ \ \ \ \text{for
some }\left(  \lambda,A,B\right)  \in\operatorname*{BT}\left(  n\right)
\text{ satisfying }\lambda\in X,
\]
and%
\[
\mathbf{a}=\nabla_{D,C}^{-\operatorname*{Col}}\ \ \ \ \ \ \ \ \ \ \text{for
some }\left(  \mu,C,D\right)  \in\operatorname*{BT}\left(  n\right)  \text{
satisfying }\mu\in Y.
\]
Consider these $\left(  \lambda,A,B\right)  $ and $\left(  \mu,C,D\right)  $.
From $\mathbf{a}=\nabla_{D,C}^{-\operatorname*{Col}}$ and $\mathbf{b}%
=\nabla_{B,A}^{\operatorname*{Row}}$, we obtain $\left\langle \mathbf{a}%
,\mathbf{b}\right\rangle =\left\langle \nabla_{D,C}^{-\operatorname*{Col}%
},\ \nabla_{B,A}^{\operatorname*{Row}}\right\rangle =0$ (by Lemma
\ref{lem.bas.mur.twosid.orth} \textbf{(a)}). Thus, Lemma
\ref{lem.bas.mur.twosid.orth} \textbf{(b)} is proved.
\end{proof}

Now let us show an analogue of Proposition \ref{prop.bas.mur.Sn-act.all}:

\begin{proposition}
\label{prop.bas.mur.twosid.Sn-act.all}Let $\operatorname*{Par}\left(
n\right)  $, $X$ and $Y$ be as in Theorem \ref{thm.bas.mur.twosid.orth-XY}.
Then, the two sets $\MurpF_{\operatorname*{all},X}^{\operatorname*{Row}}$ and
$\MurpF_{\operatorname*{all},Y}^{-\operatorname*{Col}}$ are two-sided ideals
of $\mathcal{A}$.
\end{proposition}

\begin{proof}
Similar to Proposition \ref{prop.bas.mur.Sn-act.all}. Details can be found in
the appendix (Section \ref{sec.details.bas.mur}).
\end{proof}

Next we show an analogue of Lemma \ref{lem.bas.mur.triang.1}:

\begin{lemma}
\label{lem.bas.mur.twosid.triang.1}Let $\operatorname*{Par}\left(  n\right)
$, $X$ and $Y$ be as in Theorem \ref{thm.bas.mur.twosid.orth-XY}. Then:
\medskip

\textbf{(a)} We have $\MurpF_{\operatorname*{all},Y}^{-\operatorname*{Col}%
}\subseteq\left(  \MurpF_{\operatorname*{std},X}^{\operatorname*{Row}}\right)
^{\perp}$. \medskip

\textbf{(b)} We have $\MurpF_{\operatorname*{all},X}^{\operatorname*{Row}%
}\subseteq\left(  \MurpF_{\operatorname*{std},Y}^{-\operatorname*{Col}%
}\right)  ^{\perp}$.
\end{lemma}

\begin{proof}
Similar to the proof of Lemma \ref{lem.bas.mur.triang.1}. Details can be found
in the appendix (Section \ref{sec.details.bas.mur}).
\end{proof}

More interesting is the following analogue of Lemma \ref{lem.bas.mur.triang.2}:

\begin{lemma}
\label{lem.bas.mur.twosid.triang.2}Let $\operatorname*{Par}\left(  n\right)
$, $X$ and $Y$ be as in Theorem \ref{thm.bas.mur.twosid.orth-XY}. Then:
\medskip

\textbf{(a)} We have $\left(  \MurpF_{\operatorname*{std},X}%
^{\operatorname*{Row}}\right)  ^{\perp}\subseteq\MurpF_{\operatorname*{std}%
,Y}^{-\operatorname*{Col}}$. \medskip

\textbf{(b)} We have $\left(  \MurpF_{\operatorname*{std},Y}%
^{-\operatorname*{Col}}\right)  ^{\perp}\subseteq\MurpF_{\operatorname*{std}%
,X}^{\operatorname*{Row}}$.
\end{lemma}

\begin{proof}
Similar to Lemma \ref{lem.bas.mur.triang.2}, but using Lemma
\ref{lem.bas.mur.twosid.orth} \textbf{(a)} to get rid of certain addends of
the sum. The details can be found in the appendix (Section
\ref{sec.details.bas.mur}).
\end{proof}

We are now ready to prove Theorem \ref{thm.bas.mur.twosid.orth-XY}:

\begin{proof}
[Proof of Theorem \ref{thm.bas.mur.twosid.orth-XY}.]Similar to Theorem
\ref{thm.bas.mur.triang.orth}. The details can be found in the appendix
(Section \ref{sec.details.bas.mur}).
\end{proof}

We end with an alternative characterization of the two-sided ideal
$\MurpF_{\operatorname*{std},\operatorname*{len}>k}^{-\operatorname*{Col}}$
from Theorem \ref{thm.bas.mur.twosid.orth-len} \textbf{(b)}:

\begin{proposition}
\label{prop.bas.mur.twosid.Uk}Let $k\in\mathbb{N}$. Recall the notion of an
$X$-sign integral $\nabla_{X}^{-}$ (see Definition \ref{def.intX.intX}).
Define a $\mathbf{k}$-submodule $\mathcal{U}_{k}$ of $\mathcal{A}$ by%
\[
\mathcal{U}_{k}:=\operatorname*{span}\left\{  \nabla_{X}^{-}\ \mid\ X\text{ is
a subset of }\left[  n\right]  \text{ with size }\left\vert X\right\vert
>k\right\}  .
\]
Then, the two-sided ideal $\MurpF_{\operatorname*{std},\operatorname*{len}%
>k}^{-\operatorname*{Col}}$ from Theorem \ref{thm.bas.mur.twosid.orth-len}
\textbf{(b)} satisfies%
\[
\MurpF_{\operatorname*{std},\operatorname*{len}>k}^{-\operatorname*{Col}%
}=\mathcal{AU}_{k}=\mathcal{U}_{k}\mathcal{A}=\mathcal{AU}_{k}\mathcal{A}.
\]

\end{proposition}

\begin{proof}
We shall show several claims:

\begin{statement}
\textit{Claim 1:} We have $\MurpF_{\operatorname*{std},\operatorname*{len}%
>k}^{-\operatorname*{Col}}\subseteq\mathcal{A}\mathcal{U}_{k}$.
\end{statement}

\begin{proof}
[Proof of Claim 1.]Recall that $\MurpF_{\operatorname*{std}%
,\operatorname*{len}>k}^{-\operatorname*{Col}}$ is defined as the $\mathbf{k}%
$-linear span of the vectors $\nabla_{V,U}^{-\operatorname*{Col}}$ for all
$\left(  \lambda,U,V\right)  \in\operatorname*{SBT}\left(  n\right)  $
satisfying $\ell\left(  \lambda\right)  >k$. Hence, it will suffice to show
that all these vectors $\nabla_{V,U}^{-\operatorname*{Col}}$ belong to
$\mathcal{AU}_{k}$ (because then, by linearity, it will follow that their
whole span $\MurpF_{\operatorname*{std},\operatorname*{len}>k}%
^{-\operatorname*{Col}}$ is $\subseteq\mathcal{AU}_{k}$). In other words, we
must show that $\nabla_{V,U}^{-\operatorname*{Col}}\in\mathcal{AU}_{k}$ for
each $\left(  \lambda,U,V\right)  \in\operatorname*{SBT}\left(  n\right)  $
satisfying $\ell\left(  \lambda\right)  >k$.

So let $\left(  \lambda,U,V\right)  \in\operatorname*{SBT}\left(  n\right)  $
be such that $\ell\left(  \lambda\right)  >k$. We must prove that
$\nabla_{V,U}^{-\operatorname*{Col}}\in\mathcal{AU}_{k}$.

Both $U$ and $V$ are $n$-tableaux of shape $Y\left(  \lambda\right)  $ (since
$\left(  \lambda,U,V\right)  \in\operatorname*{SBT}\left(  n\right)  $).

Using (\ref{eq.lem.bas.mur.perm.C}), we can easily see that there exists some
$w\in S_{n}$ satisfying
\[
\nabla_{V,U}^{-\operatorname*{Col}}=\left(  -1\right)  ^{w}\nabla
_{\operatorname*{Col}V}^{-}w=\left(  -1\right)  ^{w}w\nabla
_{\operatorname*{Col}U}^{-}%
\]
\footnote{\textit{Proof.} Both $U$ and $V$ are $n$-tableaux of shape $Y\left(
\lambda\right)  $. Hence, Proposition \ref{prop.tableau.Sn-act.1} \textbf{(b)}
(applied to $D=Y\left(  \lambda\right)  $ and $T=U$) shows that any
$n$-tableau of shape $Y\left(  \lambda\right)  $ can be written as
$w\rightharpoonup U$ for some $w\in S_{n}$. Hence, in particular, the
$n$-tableau $V$ can be written in this way. In other words,
$V=w\rightharpoonup U$ for some $w\in S_{n}$. Let us denote this $w$ by $u$.
Thus, $V=u\rightharpoonup U=uU$. Applying (\ref{eq.lem.bas.mur.perm.C}) to
$P=V$ and $Q=U$, we thus obtain%
\[
\nabla_{V,U}^{-\operatorname*{Col}}=\left(  -1\right)  ^{u}\nabla
_{\operatorname*{Col}V}^{-}u=\left(  -1\right)  ^{u}u\nabla
_{\operatorname*{Col}U}^{-}.
\]
Thus, we have $\nabla_{V,U}^{-\operatorname*{Col}}=\left(  -1\right)
^{w}\nabla_{\operatorname*{Col}V}^{-}w=\left(  -1\right)  ^{w}w\nabla
_{\operatorname*{Col}U}^{-}$ for some $w\in S_{n}$ (namely, for $w=u$).}.
Consider this $w$.

Choose a positive integer $g$ such that all cells of $Y\left(  \lambda\right)
$ lie in columns $1,2,\ldots,g$ (this $g$ clearly exists, since all cells of
$Y\left(  \lambda\right)  $ lie in finitely many columns with positive
indices; for example, we can take $g=\lambda_{1}$ if $\lambda_{1}>0$ and
otherwise $g=1$). Define the notation $\operatorname*{Col}\left(  j,U\right)
$ for all $j\in\mathbb{Z}$ as in Proposition \ref{prop.symmetrizers.int}
\textbf{(b)}. Then, Proposition \ref{prop.symmetrizers.int} \textbf{(b)}
(applied to $Y\left(  \lambda\right)  $, $U$ and $g$ instead of $D$, $T$ and
$k$) yields $\nabla_{\operatorname*{Col}U}^{-}=\prod_{j=1}^{g}\nabla
_{\operatorname*{Col}\left(  j,U\right)  }^{-}$ (where the order of the
factors in the product is immaterial, since these factors all commute). Thus,%
\begin{align}
\nabla_{\operatorname*{Col}U}^{-}  &  =\prod_{j=1}^{g}\nabla
_{\operatorname*{Col}\left(  j,U\right)  }^{-}\nonumber\\
&  =\left(  \prod_{j=2}^{g}\nabla_{\operatorname*{Col}\left(  j,U\right)
}^{-}\right)  \cdot\nabla_{\operatorname*{Col}\left(  1,U\right)  }^{-}
\label{pf.prop.bas.mur.twosid.Uk.c1.3}%
\end{align}
(here, we have split off the $j=1$ factor from the product, using the fact
that all the factors $\nabla_{\operatorname*{Col}\left(  j,U\right)  }^{-}$ commute).

However, it is easy to see that $\left\vert \operatorname*{Col}\left(
1,U\right)  \right\vert =\ell\left(  \lambda\right)  $%
\ \ \ \ \footnote{\textit{Proof.} Recall that $\operatorname*{Col}\left(
1,U\right)  $ is the set of all entries in the $1$-st column of $U$ (by its
definition). Hence, $\left\vert \operatorname*{Col}\left(  1,U\right)
\right\vert $ is the \# of these entries. In other words,%
\[
\left\vert \operatorname*{Col}\left(  1,U\right)  \right\vert =\left(
\text{\# of entries in the }1\text{-st column of }U\right)  .
\]
But $U$ is an $n$-tableau, so all its entries are distinct. In particular, the
entries in the $1$-st column of $U$ are all distinct. Therefore, the \# of
these entries is the \# of the cells in the $1$-st column of $Y\left(
\lambda\right)  $. In other words,
\begin{align*}
&  \left(  \text{\# of entries in the }1\text{-st column of }U\right) \\
&  =\left(  \text{\# of cells in the }1\text{-st column of }Y\left(
\lambda\right)  \right)  .
\end{align*}
\par
But Theorem \ref{thm.partitions.conj} \textbf{(b)} (applied to $j=1$) shows
that%
\begin{align}
\lambda_{1}^{t}  &  =\left(  \text{\# of positive integers }i\text{ such that
}\lambda_{i}\geq1\right) \nonumber\\
&  =\left(  \text{\# of cells in the }1\text{-st column of }Y\left(
\lambda\right)  \right)  . \label{pf.prop.bas.mur.twosid.Uk.c1.fn.5}%
\end{align}
Furthermore, Theorem \ref{thm.partitions.conj} \textbf{(c)} yields
$\lambda_{1}^{t}=\ell\left(  \lambda\right)  $. Altogether, we now find%
\begin{align*}
\left\vert \operatorname*{Col}\left(  1,U\right)  \right\vert  &  =\left(
\text{\# of entries in the }1\text{-st column of }U\right) \\
&  =\left(  \text{\# of cells in the }1\text{-st column of }Y\left(
\lambda\right)  \right) \\
&  =\lambda_{1}^{t}\ \ \ \ \ \ \ \ \ \ \left(  \text{by
(\ref{pf.prop.bas.mur.twosid.Uk.c1.fn.5})}\right) \\
&  =\ell\left(  \lambda\right)  .
\end{align*}
}. Hence,
\[
\left\vert \operatorname*{Col}\left(  1,U\right)  \right\vert =\ell\left(
\lambda\right)  >k.
\]
Thus, $\operatorname*{Col}\left(  1,U\right)  $ is a subset of $\left[
n\right]  $ with size $\left\vert \operatorname*{Col}\left(  1,U\right)
\right\vert >k$. Thus,%
\begin{align*}
\nabla_{\operatorname*{Col}\left(  1,U\right)  }^{-}  &  \in\left\{
\nabla_{X}^{-}\ \mid\ X\text{ is a subset of }\left[  n\right]  \text{ with
size }\left\vert X\right\vert >k\right\} \\
&  \subseteq\operatorname*{span}\left\{  \nabla_{X}^{-}\ \mid\ X\text{ is a
subset of }\left[  n\right]  \text{ with size }\left\vert X\right\vert
>k\right\} \\
&  =\mathcal{U}_{k}\ \ \ \ \ \ \ \ \ \ \left(  \text{by the definition of
}\mathcal{U}_{k}\right)  .
\end{align*}
Combining the above results, we find%
\begin{align}
\nabla_{V,U}^{-\operatorname*{Col}}  &  =\left(  -1\right)  ^{w}%
w\underbrace{\nabla_{\operatorname*{Col}U}^{-}}_{=\left(  \prod_{j=2}%
^{g}\nabla_{\operatorname*{Col}\left(  j,U\right)  }^{-}\right)  \cdot
\nabla_{\operatorname*{Col}\left(  1,U\right)  }^{-}}=\underbrace{\left(
-1\right)  ^{w}w\left(  \prod_{j=2}^{g}\nabla_{\operatorname*{Col}\left(
j,U\right)  }^{-}\right)  }_{\in\mathcal{A}}\cdot\,\underbrace{\nabla
_{\operatorname*{Col}\left(  1,U\right)  }^{-}}_{\in\mathcal{U}_{k}%
}\nonumber\\
&  \in\mathcal{AU}_{k}. \label{pf.prop.bas.mur.twosid.Uk.c1.-1}%
\end{align}
As we said above, this completes the proof of Claim 1.
\end{proof}

\begin{statement}
\textit{Claim 2:} We have $\MurpF_{\operatorname*{std},\operatorname*{len}%
>k}^{-\operatorname*{Col}}\subseteq\mathcal{U}_{k}\mathcal{A}$.
\end{statement}

\begin{proof}
This is analogous to Claim 1, with just a few minor differences:

\begin{itemize}
\item Instead of the equality (\ref{pf.prop.bas.mur.twosid.Uk.c1.3}), we need
to show that%
\[
\nabla_{\operatorname*{Col}V}^{-}=\nabla_{\operatorname*{Col}\left(
1,V\right)  }^{-}\cdot\left(  \prod_{j=2}^{g}\nabla_{\operatorname*{Col}%
\left(  j,V\right)  }^{-}\right)  .
\]
This, too, is obtained by applying Proposition \ref{prop.symmetrizers.int}
\textbf{(b)} (this time to $D=Y\left(  \lambda\right)  $ and $T=V$) and then
splitting off the $j=1$ factor from the resulting product (but now putting it
on the left of the other factors).

\item Instead of showing that $\nabla_{\operatorname*{Col}\left(  1,U\right)
}^{-}\in\mathcal{U}_{k}$, we now need to show that $\nabla
_{\operatorname*{Col}\left(  1,V\right)  }^{-}\in\mathcal{U}_{k}$ (but the
proof is entirely analogous: just replace $U$ by $V$).

\item Instead of the computation (\ref{pf.prop.bas.mur.twosid.Uk.c1.-1}), we
now need to argue that
\begin{align*}
\nabla_{V,U}^{-\operatorname*{Col}}  &  =\left(  -1\right)  ^{w}%
\underbrace{\nabla_{\operatorname*{Col}V}^{-}}_{=\nabla_{\operatorname*{Col}%
\left(  1,V\right)  }^{-}\cdot\left(  \prod_{j=2}^{g}\nabla
_{\operatorname*{Col}\left(  j,V\right)  }^{-}\right)  }w=\left(  -1\right)
^{w}\nabla_{\operatorname*{Col}\left(  1,V\right)  }^{-}\cdot\left(
\prod_{j=2}^{g}\nabla_{\operatorname*{Col}\left(  j,V\right)  }^{-}\right)
w\\
&  =\underbrace{\nabla_{\operatorname*{Col}\left(  1,V\right)  }^{-}}%
_{\in\mathcal{U}_{k}}\cdot\underbrace{\left(  -1\right)  ^{w}\left(
\prod_{j=2}^{g}\nabla_{\operatorname*{Col}\left(  j,V\right)  }^{-}\right)
w}_{\in\mathcal{A}}\in\mathcal{U}_{k}\mathcal{A}.
\end{align*}

\end{itemize}
\end{proof}

\begin{statement}
\textit{Claim 3:} We have $\mathcal{AU}_{k}\mathcal{A}\subseteq
\MurpF_{\operatorname*{std},\operatorname*{len}>k}^{-\operatorname*{Col}}$.
\end{statement}

\begin{proof}
[Proof of Claim 3.]Since $\MurpF_{\operatorname*{std},\operatorname*{len}%
>k}^{-\operatorname*{Col}}$ is a two-sided ideal of $\mathcal{A}$ (by Theorem
\ref{thm.bas.mur.twosid.orth-len} \textbf{(b)}), we have $\mathcal{A}%
\MurpF_{\operatorname*{std},\operatorname*{len}>k}^{-\operatorname*{Col}%
}\mathcal{A}\subseteq\MurpF_{\operatorname*{std},\operatorname*{len}%
>k}^{-\operatorname*{Col}}$ (because if $J$ is a two-sided ideal of a
$\mathbf{k}$-algebra $A$, then $AJA\subseteq J$).

Recall that $\mathcal{U}_{k}$ is defined as $\operatorname*{span}\left\{
\nabla_{X}^{-}\ \mid\ X\text{ is a subset of }\left[  n\right]  \text{ with
size }\left\vert X\right\vert >k\right\}  $. Hence, if we can show that every
subset $X$ of $\left[  n\right]  $ with size $\left\vert X\right\vert >k$
satisfies $\nabla_{X}^{-}\in\MurpF_{\operatorname*{std},\operatorname*{len}%
>k}^{-\operatorname*{Col}}$, then it will follow that $\mathcal{U}_{k}$ is a
span of some elements of $\MurpF_{\operatorname*{std},\operatorname*{len}%
>k}^{-\operatorname*{Col}}$, so that we have $\mathcal{U}_{k}\subseteq
\MurpF_{\operatorname*{std},\operatorname*{len}>k}^{-\operatorname*{Col}}$
(since any span of elements of $\MurpF_{\operatorname*{std}%
,\operatorname*{len}>k}^{-\operatorname*{Col}}$ must be a subset of
$\MurpF_{\operatorname*{std},\operatorname*{len}>k}^{-\operatorname*{Col}}$)
and therefore%
\[
\mathcal{A}\underbrace{\mathcal{U}_{k}}_{\subseteq\MurpF_{\operatorname*{std}%
,\operatorname*{len}>k}^{-\operatorname*{Col}}}\mathcal{A}\subseteq
\mathcal{A}\MurpF_{\operatorname*{std},\operatorname*{len}>k}%
^{-\operatorname*{Col}}\mathcal{A}\subseteq\MurpF_{\operatorname*{std}%
,\operatorname*{len}>k}^{-\operatorname*{Col}},
\]
and thus Claim 3 will be proved. Hence, it suffices to show that every subset
$X$ of $\left[  n\right]  $ with size $\left\vert X\right\vert >k$ satisfies
$\nabla_{X}^{-}\in\MurpF_{\operatorname*{std},\operatorname*{len}%
>k}^{-\operatorname*{Col}}$.

Let us do this now. So let $X$ be a subset of $\left[  n\right]  $ with size
$\left\vert X\right\vert >k$. We must show that $\nabla_{X}^{-}\in
\MurpF_{\operatorname*{std},\operatorname*{len}>k}^{-\operatorname*{Col}}$.

Since $X$ is a subset of $\left[  n\right]  $, we have $\left\vert
X\right\vert \leq\left\vert \left[  n\right]  \right\vert =n$. Also, from
$\left\vert X\right\vert >k$, we obtain $\left\vert X\right\vert \geq k+1$.

Let $p:=\left\vert X\right\vert -1$. Then, $p=\left\vert X\right\vert -1\leq
n-1$ (since $\left\vert X\right\vert \leq n$) and $p=\left\vert X\right\vert
-1\geq k$ (since $\left\vert X\right\vert \geq k+1$). Thus, $p\geq k\geq0$.
Also, $p=\left\vert X\right\vert -1$, so that $\left\vert X\right\vert =p+1$.
Furthermore, $n-p\geq1$ (since $p\leq n-1$).

Let $\mu$ be the partition $\left(  n-p,\underbrace{1,1,\ldots,1}_{p\text{
times}}\right)  $; this is indeed a partition (since $p\geq0$ and $n-p\geq1$)
with length $\ell\left(  \mu\right)  =1+p=p+1>p\geq k$ and size $\left\vert
\mu\right\vert =\left(  n-p\right)  +\underbrace{1+1+\cdots+1}_{p\text{
times}}=\left(  n-p\right)  +p=n$. Its Young diagram $Y\left(  \mu\right)  $
has a $1$-st column with $p+1$ cells, whereas all its other columns have at
most one cell each. (Note that $\mu$ is a hook partition, in the terminology
of Definition \ref{def.partitions.hook}.)

Let $T$ be an $n$-tableau of shape $Y\left(  \mu\right)  $ constructed as follows:

\begin{itemize}
\item Fill the $1$-st column of $Y\left(  \mu\right)  $ with the $p+1$
elements of $X$ (in some arbitrarily chosen order). (This is possible, since
the $1$-st column of $Y\left(  \mu\right)  $ has exactly $p+1$ cells, and the
set $X$ has exactly $\left\vert X\right\vert =p+1$ many elements.)

\item Place the remaining elements of $\left[  n\right]  $ (that is, the
elements of $\left[  n\right]  \setminus X$) in the remaining columns of
$Y\left(  \mu\right)  $. (Again, the order does not matter.)
\end{itemize}

For instance, if $n=7$ and $X=\left\{  2,3,5\right\}  $, then one possible
choice for $T$ is $\ytableaushort{21467,3,5}$ .

Recall that $\operatorname*{Col}\left(  1,T\right)  $ is defined as the set of
all entries in the $1$-st column of $T$. But these entries are the elements of
$X$ (since the $1$-st column of $T$ was filled with the $p+1$ elements of
$X$). Thus, $\operatorname*{Col}\left(  1,T\right)  =X$.

Now, let $j>1$ be an integer. Then, the $j$-th column of $Y\left(  \mu\right)
$ contains at most $1$ cell (since all columns of $Y\left(  \mu\right)  $
except for the $1$-st one have at most one cell each). Hence, the $j$-th
column of $T$ contains at most $1$ entry. In other words, $\left\vert
\operatorname*{Col}\left(  j,T\right)  \right\vert \leq1$ (since
$\operatorname*{Col}\left(  j,T\right)  $ is the set of all entries in the
$j$-th column of $T$), and therefore%
\begin{equation}
\nabla_{\operatorname*{Col}\left(  j,T\right)  }^{-}=1
\label{pf.prop.bas.mur.twosid.Uk.c2.5}%
\end{equation}
(since Example \ref{exa.intX.X=0123} shows that $\nabla_{Y}^{-}=1$ for any
subset $Y$ of $\left[  n\right]  $ that has size $\left\vert Y\right\vert
\leq1$).

Forget that we fixed $j$. We thus have proved
(\ref{pf.prop.bas.mur.twosid.Uk.c2.5}) for every integer $j>1$.

Now, $\mu_{1}=n-p$; hence, all cells of $Y\left(  \mu\right)  $ lie in columns
$1,2,\ldots,n-p$. Hence, Proposition \ref{prop.symmetrizers.int} \textbf{(b)}
(applied to $Y\left(  \mu\right)  $ and $n-p$ instead of $D$ and $k$) yields
\begin{align*}
\nabla_{\operatorname*{Col}T}^{-}  &  =\prod_{j=1}^{n-p}\nabla
_{\operatorname*{Col}\left(  j,T\right)  }^{-}=\underbrace{\nabla
_{\operatorname*{Col}\left(  1,T\right)  }^{-}}_{\substack{=\nabla_{X}%
^{-}\\\text{(since }\operatorname*{Col}\left(  1,T\right)  =X\text{)}}%
}\cdot\prod_{j=2}^{n-p}\underbrace{\nabla_{\operatorname*{Col}\left(
j,T\right)  }^{-}}_{\substack{=1\\\text{(by
(\ref{pf.prop.bas.mur.twosid.Uk.c2.5}))}}}\\
&  \ \ \ \ \ \ \ \ \ \ \ \ \ \ \ \ \ \ \ \ \left(
\begin{array}
[c]{c}%
\text{here, we have split off the factor for }j=1\\
\text{from the product, since }n-p\geq1
\end{array}
\right) \\
&  =\nabla_{X}^{-}\cdot\prod_{j=2}^{n-p}1=\nabla_{X}^{-},
\end{align*}
so that
\[
\nabla_{X}^{-}=\nabla_{\operatorname*{Col}T}^{-}=\nabla_{T,T}%
^{-\operatorname*{Col}}\ \ \ \ \ \ \ \ \ \ \left(  \text{by
(\ref{eq.prop.bas.mur.rsym.Col})}\right)  .
\]
But $T$ is an $n$-tableau of shape $Y\left(  \mu\right)  $. Hence, $\left(
\mu,T,T\right)  $ is an $n$-bitableau; in other words, $\left(  \mu
,T,T\right)  \in\operatorname*{BT}\left(  n\right)  $. This fact, along with
$\ell\left(  \mu\right)  >k$, shows that%
\begin{align*}
\nabla_{T,T}^{-\operatorname*{Col}}  &  \in\left\{  \nabla_{V,U}%
^{-\operatorname*{Col}}\ \mid\ \left(  \lambda,U,V\right)  \in
\operatorname*{BT}\left(  n\right)  \text{ and }\ell\left(  \lambda\right)
>k\right\} \\
&  \subseteq\operatorname*{span}\left\{  \nabla_{V,U}^{-\operatorname*{Col}%
}\ \mid\ \left(  \lambda,U,V\right)  \in\operatorname*{BT}\left(  n\right)
\text{ and }\ell\left(  \lambda\right)  >k\right\} \\
&  =\MurpF_{\operatorname*{all},\operatorname*{len}>k}^{-\operatorname*{Col}%
}\ \ \ \ \ \ \ \ \ \ \left(  \text{by the definition of }%
\MurpF_{\operatorname*{all},\operatorname*{len}>k}^{-\operatorname*{Col}%
}\right)  .
\end{align*}
Hence,
\[
\nabla_{X}^{-}=\nabla_{T,T}^{-\operatorname*{Col}}\in
\MurpF_{\operatorname*{all},\operatorname*{len}>k}^{-\operatorname*{Col}%
}=\MurpF_{\operatorname*{std},\operatorname*{len}>k}^{-\operatorname*{Col}%
}\ \ \ \ \ \ \ \ \ \ \left(  \text{by Theorem
\ref{thm.bas.mur.twosid.orth-len} \textbf{(b)}}\right)  .
\]

As explained above, this completes the proof of Claim 3.
\end{proof}

Now we are almost there: Claim 1 says that $\MurpF_{\operatorname*{std}%
,\operatorname*{len}>k}^{-\operatorname*{Col}}\subseteq\mathcal{AU}_{k}$.
Combining this with%
\[
\mathcal{A}\underbrace{\mathcal{U}_{k}}_{=\mathcal{U}_{k}1}=\mathcal{AU}%
_{k}\underbrace{1}_{\in\mathcal{A}}\subseteq\mathcal{AU}_{k}\mathcal{A}%
\subseteq\MurpF_{\operatorname*{std},\operatorname*{len}>k}%
^{-\operatorname*{Col}}\ \ \ \ \ \ \ \ \ \ \left(  \text{by Claim 3}\right)
,
\]
we obtain
\begin{equation}
\MurpF_{\operatorname*{std},\operatorname*{len}>k}^{-\operatorname*{Col}%
}=\mathcal{AU}_{k}. \label{pf.prop.bas.mur.twosid.Uk.fin.1}%
\end{equation}

Claim 2 says that $\MurpF_{\operatorname*{std},\operatorname*{len}%
>k}^{-\operatorname*{Col}}\subseteq\mathcal{U}_{k}\mathcal{A}$. Combining this
with%
\[
\underbrace{\mathcal{U}_{k}}_{=1\mathcal{U}_{k}}\mathcal{A}=\underbrace{1}%
_{\in\mathcal{A}}\mathcal{U}_{k}\mathcal{A}\subseteq\mathcal{AU}%
_{k}\mathcal{A}\subseteq\MurpF_{\operatorname*{std},\operatorname*{len}%
>k}^{-\operatorname*{Col}}\ \ \ \ \ \ \ \ \ \ \left(  \text{by Claim
3}\right)  ,
\]
we obtain
\begin{equation}
\MurpF_{\operatorname*{std},\operatorname*{len}>k}^{-\operatorname*{Col}%
}=\mathcal{U}_{k}\mathcal{A}. \label{pf.prop.bas.mur.twosid.Uk.fin.2}%
\end{equation}

Finally, combining $\MurpF_{\operatorname*{std},\operatorname*{len}%
>k}^{-\operatorname*{Col}}=\mathcal{U}_{k}\mathcal{A}\subseteq\mathcal{AU}%
_{k}\mathcal{A}$ with $\mathcal{AU}_{k}\mathcal{A}\subseteq
\MurpF_{\operatorname*{std},\operatorname*{len}>k}^{-\operatorname*{Col}}$ (by
Claim 3), we obtain%
\begin{equation}
\MurpF_{\operatorname*{std},\operatorname*{len}>k}^{-\operatorname*{Col}%
}=\mathcal{AU}_{k}\mathcal{A}. \label{pf.prop.bas.mur.twosid.Uk.fin.3}%
\end{equation}
Combining the equalities (\ref{pf.prop.bas.mur.twosid.Uk.fin.1}),
(\ref{pf.prop.bas.mur.twosid.Uk.fin.2}) and
(\ref{pf.prop.bas.mur.twosid.Uk.fin.3}), we obtain%
\[
\MurpF_{\operatorname*{std},\operatorname*{len}>k}^{-\operatorname*{Col}%
}=\mathcal{AU}_{k}=\mathcal{U}_{k}\mathcal{A}=\mathcal{AU}_{k}\mathcal{A}.
\]
Thus, the proof of Proposition \ref{prop.bas.mur.twosid.Uk} is complete.
\end{proof}

\bigskip

\appendix

\section{\label{chp.details}Omitted details and proofs}

This chapter contains some proofs (and parts of proofs) that have been omitted
from the text above -- usually because they are technical arguments or of
tangential interest only.

Each of the proofs given below uses the notations and conventions of the
chapter and section in which the respective claim appears.

\subsection{\label{sec.details.specht.tabs}Tableaux}

\begin{fineprint}
Here we supply detailed proofs for some claims made in Section
\ref{sec.specht.tabs}.
\end{fineprint}

\begin{fineprint}
\begin{proof}
[Proof of Proposition \ref{prop.tableau.std-loc}.]\textbf{(a)} Assume that
every $\left(  i,j\right)  \in D$ satisfying $\left(  i,j+1\right)  \in D$
satisfies
\begin{equation}
T\left(  i,j\right)  <T\left(  i,j+1\right)  .
\label{pf.prop.tableau.std-loc.a.ass}%
\end{equation}
We must show that $T$ is row-standard.

In other words, we must show that the entries of $T$ strictly increase (left
to right) along each row. In other words, we must show that if $\left(
i,j\right)  $ and $\left(  i,k\right)  $ are two cells in the same row of $D$,
with $\left(  i,j\right)  $ lying west of $\left(  i,k\right)  $ (that is,
$j<k$), then $T\left(  i,j\right)  <T\left(  i,k\right)  $.

So let us show this. We fix the cell $\left(  i,j\right)  \in D$; then we must
prove the following claim:

\begin{statement}
\textit{Claim 1:} If $\left(  i,k\right)  $ is a cell of $D$ satisfying $j<k$,
then $T\left(  i,j\right)  <T\left(  i,k\right)  $.
\end{statement}

\begin{proof}
[Proof of Claim 1.]We shall prove Claim 1 by induction on $k$.

The \textit{base case} is the case $k=j+1$ (since the assumption $j<k$ renders
any smaller value of $k$ impossible). In this case, Claim 1 is saying that if
$\left(  i,j+1\right)  $ is a cell of $D$, then $T\left(  i,j\right)
<T\left(  i,j+1\right)  $. But this is precisely our assumption
(\ref{pf.prop.tableau.std-loc.a.ass}). Thus, Claim 1 holds for $k=j+1$.

Now we come to the \textit{induction step:} Fix some integer $p\geq j+1$.
Assume (as the induction hypothesis) that Claim 1 holds for $k=p$. We must
prove that Claim 1 holds for $k=p+1$. In other words, we must prove that if
$\left(  i,p+1\right)  $ is a cell of $D$, then $T\left(  i,j\right)
<T\left(  i,p+1\right)  $.

So let us assume that $\left(  i,p+1\right)  $ is a cell of $D$. Thus,
$\left(  i,p+1\right)  \in D$.

But $D$ is a skew Young diagram; that is, we have $D=Y\left(  \lambda
/\mu\right)  $ for some skew partition $\lambda/\mu$. Consider this
$\lambda/\mu$. We have $\left(  i,j\right)  \in D=Y\left(  \lambda/\mu\right)
$ and $\left(  i,p+1\right)  \in D=Y\left(  \lambda/\mu\right)  $.
Furthermore, $p\geq j+1>j$, so that $j<p$ and thus $j\leq p$, so that $\left(
i,j\right)  \leq\left(  i,p\right)  $ (since $i\leq i$ and $j\leq p$).
Moreover, $\left(  i,p\right)  \leq\left(  i,p+1\right)  $ (since $i\leq i$
and $p\leq p+1$). Therefore, Proposition \ref{prop.young.convexity} (applied
to $c=\left(  i,j\right)  $ and $d=\left(  i,p\right)  $ and $e=\left(
i,p+1\right)  $) yields $\left(  i,p\right)  \in Y\left(  \lambda/\mu\right)
$. In other words, $\left(  i,p\right)  \in D$ (since $D=Y\left(  \lambda
/\mu\right)  $). Hence, (\ref{pf.prop.tableau.std-loc.a.ass}) (applied to $p$
instead of $j$) yields $T\left(  i,p\right)  <T\left(  i,p+1\right)  $ (since
$\left(  i,p+1\right)  \in D$).

But we assumed that Claim 1 holds for $k=p$. In other words, if $\left(
i,p\right)  $ is a cell of $D$ satisfying $j<p$, then $T\left(  i,j\right)
<T\left(  i,p\right)  $. Hence, we conclude that $T\left(  i,j\right)
<T\left(  i,p\right)  $ (since $\left(  i,p\right)  \in D$ and $j<p$).
Altogether, we now have $T\left(  i,j\right)  <T\left(  i,p\right)  <T\left(
i,p+1\right)  $.

Thus, we have shown that $T\left(  i,j\right)  <T\left(  i,p+1\right)  $ under
the assumption that $\left(  i,p+1\right)  $ is a cell of $D$. In other words,
Claim 1 holds for $k=p+1$. This completes the induction step, so Claim 1 is proved.
\end{proof}

As we explained, this completes the proof of Proposition
\ref{prop.tableau.std-loc} \textbf{(a)}. \medskip

\textbf{(b)} This is analogous to part \textbf{(a)}, which we have just
proved. The only necessary change is to interchange the horizontal and
vertical directions (i.e., rows become columns, x-coordinates become
y-coordinates, etc.).
\end{proof}
\end{fineprint}

\subsection{\label{sec.details.specht.nat-basis}The Young symmetrizer basis of
$\mathbf{k}\left[  S_{n}\right]  $}

\begin{fineprint}
Here we supply detailed proofs for some claims made in Section
\ref{sec.specht.nat-basis}.

\begin{proof}
[Proof of Proposition \ref{prop.tableaux.r}.]\textbf{(a)} Let $T$ be a
standard $n$-tableau of shape $D$. Thus, $T$ is standard, i.e., both
row-standard and column-standard. In particular, $T$ is row-standard. In other
words, the entries of $T$ strictly increase (left to right) along each row.

However, the $n$-tableau $T\mathbf{r}$ is obtained by reflecting $T$ across
the northwest-to-southeast diagonal. This reflection clearly turns the rows of
$T$ into the columns of $T\mathbf{r}$. More precisely, if some row of $T$ has
entries $a_{1},a_{2},\ldots,a_{k}$ from left to right, then the corresponding
column of $T\mathbf{r}$ has entries $a_{1},a_{2},\ldots,a_{k}$ from top to
bottom. Thus, the entries of $T\mathbf{r}$ strictly increase (top to bottom)
along each column (since the entries of $T$ strictly increase (left to right)
along each row). In other words, the tableau $T\mathbf{r}$ is column-standard.
A similar argument (using the column-standardness of $T$) shows that the
tableau $T\mathbf{r}$ is row-standard.

Now, we know that the $n$-tableau $T\mathbf{r}$ is both row-standard and
column-standard. In other words, it is standard. This proves Proposition
\ref{prop.tableaux.r} \textbf{(a)}. \medskip

\textbf{(b)} Let $T$ be any $n$-tableau of shape $D$. We must prove that
$\left(  T\mathbf{r}\right)  \mathbf{r}=T$.

Visually, this is obvious (since $T\mathbf{r}$ is the reflection of $T$ across
the diagonal, and $\left(  T\mathbf{r}\right)  \mathbf{r}$ is the reflection
of this reflection across the same diagonal). The rigorous proof is not much
harder: Recall that the map $\mathbf{r}:\mathbb{Z}^{2}\rightarrow
\mathbb{Z}^{2}$ is an involution, and thus we have $\mathbf{r}\left(
\mathbf{r}\left(  D\right)  \right)  =D$. Moreover, the definition of
$T\mathbf{r}$ yields $T\mathbf{r}=T\circ\left(  \mathbf{r}\mid_{\mathbf{r}%
\left(  D\right)  }\right)  $, so that%
\begin{align*}
\left(  T\mathbf{r}\right)  \mathbf{r}  &  =\left(  T\circ\left(
\mathbf{r}\mid_{\mathbf{r}\left(  D\right)  }\right)  \right)  \mathbf{r}\\
&  =T\circ\underbrace{\left(  \mathbf{r}\mid_{\mathbf{r}\left(  D\right)
}\right)  \circ\left(  \mathbf{r}\mid_{\mathbf{r}\left(  \mathbf{r}\left(
D\right)  \right)  }\right)  }_{\substack{=\operatorname*{id}\nolimits_{D}%
\\\text{(since }\mathbf{r}\text{ is an involution)}}%
}\ \ \ \ \ \ \ \ \ \ \left(  \text{by the definition of }\left(  T\circ\left(
\mathbf{r}\mid_{\mathbf{r}\left(  D\right)  }\right)  \right)  \mathbf{r}%
\right) \\
&  =T\circ\operatorname*{id}\nolimits_{D}=T.
\end{align*}
Thus, Proposition \ref{prop.tableaux.r} \textbf{(b)} is proved. \medskip

\textbf{(c)} Recall that the map $\mathbf{r}:\mathbb{Z}^{2}\rightarrow
\mathbb{Z}^{2}$ is an involution, and thus we have $\mathbf{r}\left(
\mathbf{r}\left(  D\right)  \right)  =D$. Furthermore, Proposition
\ref{prop.tableaux.r} \textbf{(b)} yields that $\left(  T\mathbf{r}\right)
\mathbf{r}=T$ for any $T\in\operatorname*{SYT}\left(  D\right)  $. Likewise,
we have $\left(  T\mathbf{r}\right)  \mathbf{r}=T$ for any $T\in
\operatorname*{SYT}\left(  \mathbf{r}\left(  D\right)  \right)  $.

Proposition \ref{prop.tableaux.r} \textbf{(a)} is saying that if $T$ is a
standard $n$-tableau of shape $D$, then $T\mathbf{r}$ is a standard
$n$-tableau of shape $\mathbf{r}\left(  D\right)  $. In other words, if
$T\in\operatorname*{SYT}\left(  D\right)  $, then $T\mathbf{r}\in
\operatorname*{SYT}\left(  \mathbf{r}\left(  D\right)  \right)  $ (since
$T\in\operatorname*{SYT}\left(  D\right)  $ means that $T$ is a standard
$n$-tableau of shape $D$, whereas $T\mathbf{r}\in\operatorname*{SYT}\left(
\mathbf{r}\left(  D\right)  \right)  $ means that $T\mathbf{r}$ is a standard
$n$-tableau of shape $\mathbf{r}\left(  D\right)  $). Thus, the map
\begin{align}
\operatorname*{SYT}\left(  D\right)   &  \rightarrow\operatorname*{SYT}\left(
\mathbf{r}\left(  D\right)  \right)  ,\nonumber\\
T  &  \mapsto T\mathbf{r} \label{pf.prop.tableaux.r.c.1}%
\end{align}
is well-defined. The same argument (applied to $\mathbf{r}\left(  D\right)  $
instead of $D$) yields that the map%
\begin{align*}
\operatorname*{SYT}\left(  \mathbf{r}\left(  D\right)  \right)   &
\rightarrow\operatorname*{SYT}\left(  \mathbf{r}\left(  \mathbf{r}\left(
D\right)  \right)  \right)  ,\\
T  &  \mapsto T\mathbf{r}%
\end{align*}
is well-defined as well. In other words, the map%
\begin{align}
\operatorname*{SYT}\left(  \mathbf{r}\left(  D\right)  \right)   &
\rightarrow\operatorname*{SYT}\left(  D\right)  ,\nonumber\\
T  &  \mapsto T\mathbf{r} \label{pf.prop.tableaux.r.c.2}%
\end{align}
is well-defined (since $\mathbf{r}\left(  \mathbf{r}\left(  D\right)  \right)
=D$). These two maps (\ref{pf.prop.tableaux.r.c.1}) and
(\ref{pf.prop.tableaux.r.c.2}) are mutually inverse (since $\left(
T\mathbf{r}\right)  \mathbf{r}=T$ for any $T\in\operatorname*{SYT}\left(
D\right)  $, and since $\left(  T\mathbf{r}\right)  \mathbf{r}=T$ for any
$T\in\operatorname*{SYT}\left(  \mathbf{r}\left(  D\right)  \right)  $).
Hence, in particular, the former map is invertible, i.e., is a bijection. This
proves Proposition \ref{prop.tableaux.r} \textbf{(c)}. \medskip

\textbf{(d)} Let $T$ be an $n$-tableau of shape $D$. Let $w\in S_{n}$. Then,
$T\mathbf{r}=T\circ\left(  \mathbf{r}\mid_{\mathbf{r}\left(  D\right)
}\right)  $ (by the definition of $T\mathbf{r}$) and $\left(  w\circ T\right)
\mathbf{r}=w\circ T\circ\left(  \mathbf{r}\mid_{\mathbf{r}\left(  D\right)
}\right)  $ (likewise). But $wT=w\rightharpoonup T=w\circ T$ (by Definition
\ref{def.tableau.Sn-act}). Thus,
\begin{equation}
\left(  wT\right)  \mathbf{r}=\left(  w\circ T\right)  \mathbf{r}=w\circ
T\circ\left(  \mathbf{r}\mid_{\mathbf{r}\left(  D\right)  }\right)  .
\label{pf.prop.tableaux.r.d.1}%
\end{equation}
On the other hand, $w\left(  T\mathbf{r}\right)  =w\rightharpoonup\left(
T\mathbf{r}\right)  =w\circ\left(  T\mathbf{r}\right)  $ (again by Definition
\ref{def.tableau.Sn-act}), and thus%
\[
w\left(  T\mathbf{r}\right)  =w\circ\underbrace{\left(  T\mathbf{r}\right)
}_{=T\circ\left(  \mathbf{r}\mid_{\mathbf{r}\left(  D\right)  }\right)
}=w\circ T\circ\left(  \mathbf{r}\mid_{\mathbf{r}\left(  D\right)  }\right)
.
\]
Comparing this with (\ref{pf.prop.tableaux.r.d.1}), we obtain $\left(
wT\right)  \mathbf{r}=w\left(  T\mathbf{r}\right)  $. This proves Proposition
\ref{prop.tableaux.r} \textbf{(d)}. \medskip

\textbf{(e)} Let $T$ be an $n$-tableau of shape $D$. We must prove that
$\mathcal{R}\left(  T\right)  =\mathcal{C}\left(  T\mathbf{r}\right)  $ and
$\mathcal{C}\left(  T\right)  =\mathcal{R}\left(  T\mathbf{r}\right)  $.

Recall once again that the $n$-tableau $T\mathbf{r}$ is obtained by reflecting
$T$ across the northwest-to-southeast diagonal. This reflection turns the rows
of $T$ into the columns of $T\mathbf{r}$. Thus, two numbers $i,j\in\left[
n\right]  $ lie in the same row of $T$ if and only if they lie in the same
column of $T\mathbf{r}$. Hence, a permutation $w\in S_{n}$ is horizontal for
$T$ if and only if it is vertical for $T\mathbf{r}$ (because \textquotedblleft
horizontal for $T$\textquotedblright\ means that for each $i\in\left[
n\right]  $, the number $w\left(  i\right)  $ lies in the same row of $T$ as
$i$ does, whereas \textquotedblleft vertical for $T\mathbf{r}$%
\textquotedblright\ means that for each $i\in\left[  n\right]  $, the number
$w\left(  i\right)  $ lies in the same column of $T\mathbf{r}$ as $i$ does).

But $\mathcal{R}\left(  T\right)  $ is defined as the subgroup of $S_{n}$
consisting of all permutations $w\in S_{n}$ that are horizontal for $T$,
whereas $\mathcal{C}\left(  T\mathbf{r}\right)  $ is defined as the subgroup
of $S_{n}$ consisting of all permutations $w\in S_{n}$ that are vertical for
$T\mathbf{r}$. Thus, $\mathcal{R}\left(  T\right)  =\mathcal{C}\left(
T\mathbf{r}\right)  $ (since a permutation $w\in S_{n}$ is horizontal for $T$
if and only if it is vertical for $T\mathbf{r}$).

An analogous argument (but with the roles of rows and columns interchanged)
shows that $\mathcal{C}\left(  T\right)  =\mathcal{R}\left(  T\mathbf{r}%
\right)  $. (Alternatively, we can derive $\mathcal{C}\left(  T\right)
=\mathcal{R}\left(  T\mathbf{r}\right)  $ by applying the already proved
equality $\mathcal{R}\left(  T\right)  =\mathcal{C}\left(  T\mathbf{r}\right)
$ to $\mathbf{r}\left(  D\right)  $ and $T\mathbf{r}$ instead of $D$ and $T$,
and then simplifying the result using Proposition \ref{prop.tableaux.r}
\textbf{(b)}.)

The proof of Proposition \ref{prop.tableaux.r} \textbf{(e)} is now complete
(since we have proved both $\mathcal{R}\left(  T\right)  =\mathcal{C}\left(
T\mathbf{r}\right)  $ and $\mathcal{C}\left(  T\right)  =\mathcal{R}\left(
T\mathbf{r}\right)  $).
\end{proof}
\end{fineprint}

\subsection{\label{sec.details.specht.conj-class-action}The action of
conjugacy class sums on Specht modules}

\begin{fineprint}
Here we supply detailed proofs for some claims made in Section
\ref{sec.specht.conj-class-action}.

\begin{proof}
[Proof of Proposition \ref{prop.partitions.nlam-sumij}.]Let us write the
partition $\lambda$ as $\lambda=\left(  \lambda_{1},\lambda_{2},\ldots
,\lambda_{k}\right)  $. Then, the definition of $Y\left(  \lambda\right)  $
yields%
\[
Y\left(  \lambda\right)  =\left\{  \left(  i,j\right)  \in\left\{
1,2,\ldots,k\right\}  \times\left\{  1,2,3,\ldots\right\}  \ \mid
\ j\leq\lambda_{i}\right\}  .
\]
Thus, the summation sign $\sum_{\left(  i,j\right)  \in Y\left(
\lambda\right)  }$ can be rewritten as $\sum_{\substack{\left(  i,j\right)
\in\left\{  1,2,\ldots,k\right\}  \times\left\{  1,2,3,\ldots\right\}
;\\j\leq\lambda_{i}}}$. Hence,%
\begin{align*}
\sum_{\left(  i,j\right)  \in Y\left(  \lambda\right)  }\left(  i-1\right)
&  =\underbrace{\sum_{\substack{\left(  i,j\right)  \in\left\{  1,2,\ldots
,k\right\}  \times\left\{  1,2,3,\ldots\right\}  ;\\j\leq\lambda_{i}}}}%
_{=\sum_{i\in\left\{  1,2,\ldots,k\right\}  }\ \ \sum_{\substack{j\in\left\{
1,2,3,\ldots\right\}  ;\\j\leq\lambda_{i}}}}\left(  i-1\right)
=\underbrace{\sum_{i\in\left\{  1,2,\ldots,k\right\}  }}_{=\sum_{i=1}^{k}%
}\ \ \underbrace{\sum_{\substack{j\in\left\{  1,2,3,\ldots\right\}
;\\j\leq\lambda_{i}}}}_{=\sum_{j=1}^{\lambda_{i}}}\left(  i-1\right) \\
&  =\sum_{i=1}^{k}\ \ \underbrace{\sum_{j=1}^{\lambda_{i}}\left(  i-1\right)
}_{\substack{=\lambda_{i}\left(  i-1\right)  \\=\left(  i-1\right)
\lambda_{i}}}=\sum_{i=1}^{k}\left(  i-1\right)  \lambda_{i}=\operatorname*{n}%
\left(  \lambda\right)
\end{align*}
(by the definition of $\operatorname*{n}\left(  \lambda\right)  $). This
proves Proposition \ref{prop.partitions.nlam-sumij}.
\end{proof}

\begin{proof}
[Proof of Proposition \ref{prop.partitions.nlamt}.]Let us write the partition
$\lambda$ itself as $\lambda=\left(  \lambda_{1},\lambda_{2},\ldots
,\lambda_{k}\right)  $. We shall first show that%
\[
\operatorname*{n}\left(  \lambda^{t}\right)  =\sum_{i=1}^{k}\dbinom
{\lambda_{i}}{2}.
\]
Once this is proved, we will then apply this to $\lambda^{t}$ instead of
$\lambda$ and easily derive the proposition.

Theorem \ref{thm.partitions.conj} \textbf{(b)} shows that
\begin{align}
\lambda_{j}^{t}  &  =\left(  \text{\# of positive integers }i\text{ such that
}\lambda_{i}\geq j\right) \label{pf.prop.partitions.nlamt.1}\\
&  =\left(  \text{\# of cells in the }j\text{-th column of }Y\left(
\lambda\right)  \right)  \ \ \ \ \ \ \ \ \ \ \text{for all }j\geq1.\nonumber
\end{align}
Thus, for all $j\geq1$, we have%
\begin{equation}
\lambda_{j}^{t}=\sum_{\substack{i\in\left[  k\right]  ;\\j\leq\lambda_{i}}}1
\label{pf.prop.partitions.nlamt.2}%
\end{equation}
\footnote{\textit{Proof.} Let $j\geq1$. Recall that $\lambda=\left(
\lambda_{1},\lambda_{2},\ldots,\lambda_{k}\right)  $. Hence,
\begin{equation}
\lambda_{i}=0\ \ \ \ \ \ \ \ \ \ \text{for each }i>k.
\label{pf.prop.partitions.nlamt.2.pf.1}%
\end{equation}
Hence, if $i$ is a positive integer such that $\lambda_{i}\geq j$, then we
must have $i\leq k$ (since otherwise, we would have $i>k$ and thus
$\lambda_{i}=0$ (by (\ref{pf.prop.partitions.nlamt.2.pf.1})), which would
contradict $\lambda_{i}\geq j\geq1>0$) and therefore $i\in\left[  k\right]  $.
In other words, the positive integers $i$ that satisfy $\lambda_{i}\geq j$ all
belong to $\left[  k\right]  $. Hence, these positive integers are precisely
the $i\in\left[  k\right]  $ that satisfy $\lambda_{i}\geq j$. Therefore,%
\begin{align}
&  \left(  \text{\# of positive integers }i\text{ that satisfy }\lambda
_{i}\geq j\right) \nonumber\\
&  =\left(  \text{\# of }i\in\left[  k\right]  \text{ that satisfy }%
\lambda_{i}\geq j\right)  . \label{pf.prop.partitions.nlamt.2.pf.2}%
\end{align}
Now,%
\begin{align*}
\sum_{\substack{i\in\left[  k\right]  ;\\j\leq\lambda_{i}}}1  &  =\left(
\text{\# of }i\in\left[  k\right]  \text{ that satisfy }j\leq\lambda
_{i}\right)  \cdot1\\
&  =\left(  \text{\# of }i\in\left[  k\right]  \text{ that satisfy }%
j\leq\lambda_{i}\right) \\
&  =\left(  \text{\# of }i\in\left[  k\right]  \text{ that satisfy }%
\lambda_{i}\geq j\right)  \ \ \ \ \ \ \ \ \ \ \left(  \text{here, we rewrote
\textquotedblleft}j\leq\lambda_{i}\text{\textquotedblright\ as
\textquotedblleft}\lambda_{i}\geq j\text{\textquotedblright}\right) \\
&  =\left(  \text{\# of positive integers }i\text{ that satisfy }\lambda
_{i}\geq j\right)  \ \ \ \ \ \ \ \ \ \ \left(  \text{by
(\ref{pf.prop.partitions.nlamt.2.pf.2})}\right) \\
&  =\left(  \text{\# of positive integers }i\text{ such that }\lambda_{i}\geq
j\right)  =\lambda_{j}^{t}\ \ \ \ \ \ \ \ \ \ \left(  \text{by
(\ref{pf.prop.partitions.nlamt.1})}\right)  .
\end{align*}
This proves (\ref{pf.prop.partitions.nlamt.2}).}.

Moreover, for each $u\geq1$, we have $\lambda_{u}\leq m$%
\ \ \ \ \footnote{\textit{Proof.} Let $u\geq1$. We must prove that
$\lambda_{u}\leq m$.
\par
If $\lambda_{u}=0$, then this is obvious (since $0\leq m$). Thus, we WLOG
assume that $\lambda_{u}\neq0$. Hence, $\lambda_{u}\geq1$.
\par
Let $j=\lambda_{u}$. Thus, $\lambda_{u}\geq j$ (since $\lambda_{u}=j$). Hence,
there exists at least one positive integer $i$ such that $\lambda_{i}\geq j$
(namely, the integer $i=u$). In other words, $\left(  \text{\# of positive
integers }i\text{ such that }\lambda_{i}\geq j\right)  \geq1$. Using
(\ref{pf.prop.partitions.nlamt.1}), we can rewrite this as $\lambda_{j}%
^{t}\geq1$ (since $j=\lambda_{u}\geq1$).
\par
However, $\lambda^{t}=\left(  \lambda_{1}^{t},\lambda_{2}^{t},\ldots
,\lambda_{m}^{t}\right)  $, and thus $\lambda_{p}^{t}=0$ for all $p>m$. Hence,
if we had $j>m$, then we would have $\lambda_{j}^{t}=0$, which would
contradict $\lambda_{j}^{t}\geq1>0$. Thus, we cannot have $j>m$. Hence, we
have $j\leq m$, so that $\lambda_{u}=j\leq m$. Qed.} and thus $\left[
\lambda_{u}\right]  \subseteq\left[  m\right]  $, so that%
\begin{equation}
\left[  m\right]  \cap\left[  \lambda_{u}\right]  =\left[  \lambda_{u}\right]
. \label{pf.prop.partitions.nlamt.3}%
\end{equation}
On the other hand, each $p\in\mathbb{N}$ satisfies%
\begin{align}
\sum_{j\in\left[  p\right]  }\left(  j-1\right)   &  =\sum_{j=1}^{p}\left(
j-1\right)  =0+1+\cdots+\left(  p-1\right) \nonumber\\
&  =1+2+\cdots+\left(  p-1\right) \nonumber\\
&  =\dfrac{\left(  p-1\right)  \left(  \left(  p-1\right)  +1\right)  }%
{2}\ \ \ \ \ \ \ \ \ \ \left(  \text{by the Little Gauss formula}\right)
\nonumber\\
&  =\dfrac{p\left(  p-1\right)  }{2}=\dbinom{p}{2}.
\label{pf.prop.partitions.nlamt.LG}%
\end{align}
Hence, for each $i\geq1$, we have%
\begin{align}
\underbrace{\sum_{\substack{j\in\left[  m\right]  ;\\j\leq\lambda_{i}}%
}}_{\substack{=\sum_{\substack{j\in\left[  m\right]  ;\\j\in\left[
\lambda_{i}\right]  }}\\\text{(since the condition \textquotedblleft}%
j\leq\lambda_{i}\text{\textquotedblright}\\\text{on a positive integer
}j\text{ is}\\\text{equivalent to \textquotedblleft}j\in\left[  \lambda
_{i}\right]  \text{\textquotedblright)}}}\left(  j-1\right)   &
=\underbrace{\sum_{\substack{j\in\left[  m\right]  ;\\j\in\left[  \lambda
_{i}\right]  }}}_{\substack{=\sum_{j\in\left[  m\right]  \cap\left[
\lambda_{i}\right]  }=\sum_{j\in\left[  \lambda_{i}\right]  }\\\text{(since
(\ref{pf.prop.partitions.nlamt.3}) (applied to }u=i\text{)}\\\text{yields
}\left[  m\right]  \cap\left[  \lambda_{i}\right]  =\left[  \lambda
_{i}\right]  \text{)}}}\left(  j-1\right)  =\sum_{j\in\left[  \lambda
_{i}\right]  }\left(  j-1\right) \nonumber\\
&  =\dbinom{\lambda_{i}}{2} \label{pf.prop.partitions.nlamt.4}%
\end{align}
(by (\ref{pf.prop.partitions.nlamt.LG}), applied to $p=\lambda_{i}$).

Now, $\lambda^{t}=\left(  \lambda_{1}^{t},\lambda_{2}^{t},\ldots,\lambda
_{m}^{t}\right)  $. Thus, the definition of $\operatorname*{n}\left(
\lambda^{t}\right)  $ yields%
\begin{align*}
\operatorname*{n}\left(  \lambda^{t}\right)   &  =\sum_{i=1}^{m}\left(
i-1\right)  \lambda_{i}^{t}=\underbrace{\sum_{j=1}^{m}}_{=\sum_{j\in\left[
m\right]  }}\left(  j-1\right)  \underbrace{\lambda_{j}^{t}}_{\substack{=\sum
_{\substack{i\in\left[  k\right]  ;\\j\leq\lambda_{i}}}1\\\text{(by
(\ref{pf.prop.partitions.nlamt.2}))}}}=\sum_{j\in\left[  m\right]  }\left(
j-1\right)  \sum_{\substack{i\in\left[  k\right]  ;\\j\leq\lambda_{i}}}1\\
&  =\underbrace{\sum_{j\in\left[  m\right]  }\ \ \sum_{\substack{i\in\left[
k\right]  ;\\j\leq\lambda_{i}}}}_{=\sum_{i\in\left[  k\right]  }%
\ \ \sum_{\substack{j\in\left[  m\right]  ;\\j\leq\lambda_{i}}}}%
\underbrace{\left(  j-1\right)  \cdot1}_{=j-1}=\underbrace{\sum_{i\in\left[
k\right]  }}_{=\sum_{i=1}^{k}}\ \ \underbrace{\sum_{\substack{j\in\left[
m\right]  ;\\j\leq\lambda_{i}}}\left(  j-1\right)  }_{\substack{=\dbinom
{\lambda_{i}}{2}\\\text{(by (\ref{pf.prop.partitions.nlamt.4}))}}}=\sum
_{i=1}^{k}\dbinom{\lambda_{i}}{2}.
\end{align*}

The same reasoning (applied to $\lambda^{t}$, $k$ and $m$ instead of $\lambda
$, $m$ and $k$) shows that%
\begin{equation}
\operatorname*{n}\left(  \left(  \lambda^{t}\right)  ^{t}\right)  =\sum
_{i=1}^{m}\dbinom{\lambda_{i}^{t}}{2}. \label{pf.prop.partitions.nlamt.7}%
\end{equation}
However, Theorem \ref{thm.partitions.conj} \textbf{(d)} says that $\left(
\lambda^{t}\right)  ^{t}=\lambda$. Hence, we can rewrite
(\ref{pf.prop.partitions.nlamt.7}) as
\[
\operatorname*{n}\left(  \lambda\right)  =\sum_{i=1}^{m}\dbinom{\lambda
_{i}^{t}}{2}.
\]
This proves Proposition \ref{prop.partitions.nlamt}.
\end{proof}

\begin{proof}
[Proof of Proposition \ref{prop.partitions.nlam-counts-t}.]We know that $T$ is
an $n$-tableau of shape $\lambda$. That is, $T$ is an $n$-tableau of shape
$Y\left(  \lambda\right)  $. In particular, $T$ is injective (like any
$n$-tableau), i.e., all entries of $T$ are distinct.

Write the conjugate $\lambda^{t}$ of $\lambda$ as $\lambda^{t}=\left(
\lambda_{1}^{t},\lambda_{2}^{t},\ldots,\lambda_{m}^{t}\right)  $. Thus,%
\begin{equation}
\lambda_{k}^{t}=0\ \ \ \ \ \ \ \ \ \ \text{for all }k>m.
\label{pf.prop.partitions.nlam-counts-t.1}%
\end{equation}

Any transposition in $S_{n}$ has the form $t_{i,j}$ for some $2$-element
subset $\left\{  i,j\right\}  $ of $\left[  n\right]  $; moreover, the subset
$\left\{  i,j\right\}  $ is uniquely determined by this transposition.
Moreover, a transposition $t_{i,j}\in S_{n}$ belongs to $\mathcal{C}\left(
T\right)  $ if and only if the numbers $i$ and $j$ lie in the same column of
$T$\ \ \ \ \footnote{\textit{Proof.} Let $t_{i,j}\in S_{n}$ be a
transposition. If the numbers $i$ and $j$ lie in the same column of $T$, then
the transposition $t_{i,j}$ is vertical for $T$ and thus belongs to
$\mathcal{C}\left(  T\right)  $. On the other hand, if the numbers $i$ and $j$
do not lie in the same column of $T$, then the transposition $t_{i,j}$ is not
vertical for $T$ (since $t_{i,j}\left(  i\right)  =j$ does not lie in the same
column of $T$ as $i$ does), and thus does not belong to $\mathcal{C}\left(
T\right)  $. Combining the preceding two sentences, we conclude that $t_{i,j}$
belongs to $\mathcal{C}\left(  T\right)  $ if and only if the numbers $i$ and
$j$ lie in the same column of $T$. Qed.}. Hence, the transpositions in
$\mathcal{C}\left(  T\right)  $ are precisely the transpositions of the form
$t_{i,j}$, where $i$ and $j$ are two distinct numbers lying in the same column
of $T$. Therefore,%
\begin{align}
&  \left(  \text{\# of transpositions in }\mathcal{C}\left(  T\right)  \right)
\nonumber\\
&  =\left(  \text{\# of transpositions of the form }t_{i,j}\text{, where
}i\text{ and }j\text{ are}\right. \nonumber\\
&  \ \ \ \ \ \ \ \ \ \ \left.  \text{two distinct numbers lying in the same
column of }T\right) \nonumber\\
&  =\sum_{k\geq1}\left(  \text{\# of transpositions of the form }%
t_{i,j}\text{, where }i\text{ and }j\text{ are}\right. \nonumber\\
&  \ \ \ \ \ \ \ \ \ \ \left.  \text{two distinct numbers lying in the
}k\text{-th column of }T\right)  \label{pf.prop.partitions.nlam-counts-t.2}%
\end{align}
(since all entries of $T$ are distinct, and thus a given number $i$ cannot lie
in two different columns of $T$ at once).

However, we can prove the following:

\begin{statement}
\textit{Claim 1:} Each $k\geq1$ satisfies%
\begin{align}
&  \left(  \text{\# of transpositions of the form }t_{i,j}\text{, where
}i\text{ and }j\text{ are}\right. \nonumber\\
&  \ \ \ \ \ \ \ \ \ \ \left.  \text{two distinct numbers lying in the
}k\text{-th column of }T\right) \nonumber\\
&  =\dbinom{\lambda_{k}^{t}}{2}. \label{pf.prop.partitions.nlam-counts-t.3}%
\end{align}

\end{statement}

\begin{proof}
[Proof of Claim 1.]Let $k\geq1$. From Theorem \ref{thm.partitions.conj}
\textbf{(b)} (applied to $j=k$), we obtain%
\begin{align*}
\lambda_{k}^{t}  &  =\left(  \text{\# of positive integers }i\text{ such that
}\lambda_{i}\geq k\right) \\
&  =\left(  \text{\# of cells in the }k\text{-th column of }Y\left(
\lambda\right)  \right)  .
\end{align*}
Hence, the $k$-th column of $Y\left(  \lambda\right)  $ has $\lambda_{k}^{t}$
cells. Thus, the $k$-th column of $T$ has $\lambda_{k}^{t}$ entries (since $T$
is an $n$-tableau of shape $Y\left(  \lambda\right)  $). All these
$\lambda_{k}^{t}$ entries are distinct (since all entries of $T$ are
distinct). Hence, there are $\lambda_{k}^{t}$ distinct entries in the $k$-th
column of $T$. Thus, there are $\dbinom{\lambda_{k}^{t}}{2}$ ways to choose
two distinct entries $i$ and $j$ from the $k$-th column of $T$, provided that
the order does not matter (i.e., we do not distinguish between choosing $i$
and $j$ and choosing $j$ and $i$). Each of these $\dbinom{\lambda_{k}^{t}}{2}$
ways yields a transposition $t_{i,j}$, and all these $\dbinom{\lambda_{k}^{t}%
}{2}$ transpositions are distinct (since a transposition $t_{i,j}$ uniquely
determines the $2$-element set $\left\{  i,j\right\}  $). Hence,%
\begin{align*}
&  \left(  \text{\# of transpositions of the form }t_{i,j}\text{, where
}i\text{ and }j\text{ are}\right. \\
&  \ \ \ \ \ \ \ \ \ \ \left.  \text{two distinct numbers lying in the
}k\text{-th column of }T\right) \\
&  =\dbinom{\lambda_{k}^{t}}{2}.
\end{align*}
Thus, Claim 1 is proved.
\end{proof}

Now, using (\ref{pf.prop.partitions.nlam-counts-t.3}), we can rewrite
(\ref{pf.prop.partitions.nlam-counts-t.2}) as%
\begin{align*}
\left(  \text{\# of transpositions in }\mathcal{C}\left(  T\right)  \right)
&  =\sum_{k\geq1}\dbinom{\lambda_{k}^{t}}{2}=\sum_{k=1}^{m}\dbinom{\lambda
_{k}^{t}}{2}+\sum_{k=m+1}^{\infty}\underbrace{\dbinom{\lambda_{k}^{t}}{2}%
}_{\substack{=0\\\text{(since }k\geq m+1>m\text{,}\\\text{so that
(\ref{pf.prop.partitions.nlam-counts-t.1}) yields }\lambda_{k}^{t}%
=0\\\text{and thus }\dbinom{\lambda_{k}^{t}}{2}=\dbinom{0}{2}=0\text{)}}}\\
&  =\sum_{k=1}^{m}\dbinom{\lambda_{k}^{t}}{2}+\underbrace{\sum_{k=m+1}%
^{\infty}0}_{=0}=\sum_{k=1}^{m}\dbinom{\lambda_{k}^{t}}{2}=\sum_{i=1}%
^{m}\dbinom{\lambda_{i}^{t}}{2}\\
&  =\operatorname*{n}\left(  \lambda\right)  \ \ \ \ \ \ \ \ \ \ \left(
\text{by Proposition \ref{prop.partitions.nlamt}}\right)  .
\end{align*}
In other words, $\operatorname*{n}\left(  \lambda\right)  $ is the number of
transpositions in $\mathcal{C}\left(  T\right)  $. This proves Proposition
\ref{prop.partitions.nlam-counts-t}.
\end{proof}
\end{fineprint}

\subsection{\label{sec.details.bas.what}What to do with a basis}

\begin{fineprint}
Here we supply detailed proofs for some claims made in Section
\ref{sec.bas.what}.

\begin{proof}
[Detailed proof of Lemma \ref{lem.sbt.sizen!}.]For each partition $\lambda$ of
$n$, we have%
\begin{align}
&  \left(  \text{\# of standard }n\text{-tableaux of shape }\lambda\right)
\nonumber\\
&  =\left(  \text{\# of standard }n\text{-tableaux of shape }Y\left(
\lambda\right)  \right) \nonumber\\
&  =\left(  \text{\# of standard tableaux of shape }Y\left(  \lambda\right)
\right) \nonumber\\
&  \ \ \ \ \ \ \ \ \ \ \ \ \ \ \ \ \ \ \ \ \left(
\begin{array}
[c]{c}%
\text{since Proposition \ref{prop.tableau.std-n} (applied to }D=Y\left(
\lambda\right)  \text{)}\\
\text{shows that the standard tableaux of shape }Y\left(  \lambda\right) \\
\text{are precisely the standard }n\text{-tableaux of shape }Y\left(
\lambda\right) \\
\text{(because }\left\vert Y\left(  \lambda\right)  \right\vert =\left\vert
\lambda\right\vert =n\text{)}%
\end{array}
\right) \nonumber\\
&  =f^{\lambda}\ \ \ \ \ \ \ \ \ \ \left(  \text{by the definition of
}f^{\lambda}\right)  . \label{pf.lem.sbt.sizen!.1}%
\end{align}
Recall that $\operatorname*{SBT}\left(  n\right)  $ is the set of all standard
$n$-bitableaux, i.e., of all triples $\left(  \lambda,U,V\right)  $ consisting
of a partition $\lambda$ of $n$ and two standard $n$-tableaux $U,V$ of shape
$\lambda$. In other words,%
\[
\operatorname*{SBT}\left(  n\right)  =\bigcup_{\lambda\text{ is a partition of
}n}\left\{  \left(  \lambda,U,V\right)  \ \mid\ U\text{ and }V\text{ are
standard }n\text{-tableaux of shape }\lambda\right\}  .
\]
This union is a disjoint union (since $\lambda$ can be uniquely recovered from
the triple $\left(  \lambda,U,V\right)  $). Since the size of a disjoint union
equals the sum of the sizes of the sets involved, we thus obtain
\begin{align*}
&  \left\vert \operatorname*{SBT}\left(  n\right)  \right\vert \\
&  =\sum_{\lambda\text{ is a partition of }n}\underbrace{\left\vert \left\{
\left(  \lambda,U,V\right)  \ \mid\ U\text{ and }V\text{ are standard
}n\text{-tableaux of shape }\lambda\right\}  \right\vert }_{\substack{=\left(
\text{\# of pairs }\left(  U,V\right)  \text{ where }U\text{ and }V\text{ are
standard }n\text{-tableaux of shape }\lambda\right)  \\=\left(  \text{\# of
standard }n\text{-tableaux of shape }\lambda\right)  ^{2}\\=\left(
f^{\lambda}\right)  ^{2}\\\text{(by (\ref{pf.lem.sbt.sizen!.1}))}}}\\
&  =\sum_{\lambda\text{ is a partition of }n}\left(  f^{\lambda}\right)
^{2}=n!\ \ \ \ \ \ \ \ \ \ \left(  \text{by Corollary
\ref{cor.spechtmod.sumflam2}}\right)  .
\end{align*}
This proves Lemma \ref{lem.sbt.sizen!}.
\end{proof}
\end{fineprint}

\subsection{\label{sec.details.bas.tord}Combinatorial interlude: Lexicographic
orders}

\begin{fineprint}
Here we supply detailed proofs for some claims made in Section
\ref{sec.bas.tord}.

\begin{proof}
[Detailed proof of Proposition \ref{prop.bas.tord.pars.tord}.]First, we shall
show that the relation $<$ is the smaller relation of a partial order on the
set of all partitions. In order to prove this, we need to verify three claims:

\begin{statement}
\textit{Claim 1:} The relation $<$ is irreflexive (i.e., there exists no
partition $\lambda$ such that $\lambda<\lambda$).
\end{statement}

\begin{statement}
\textit{Claim 2:} The relation $<$ is transitive (i.e., if three partitions
$\lambda$, $\mu$ and $\nu$ satisfy $\lambda<\mu$ and $\mu<\nu$, then
$\lambda<\nu$).
\end{statement}

\begin{statement}
\textit{Claim 3:} The relation $<$ is asymmetric (i.e., no two partitions
$\lambda$ and $\mu$ satisfy $\lambda<\mu$ and $\mu<\lambda$ at the same time).
\end{statement}

\begin{proof}
[Proof of Claim 1.]Let $\lambda$ be a partition such that $\lambda<\lambda$.
We shall derive a contradiction.

Indeed, we have $\lambda<\lambda$. By the definition of the relation $<$, this
means that

\begin{itemize}
\item there exists at least one $i\geq1$ such that $\lambda_{i}\neq\lambda
_{i}$, and

\item the \textbf{smallest} such $i$ satisfies $\lambda_{i}<\lambda_{i}$.
\end{itemize}

\noindent But the first of these two bullet points is clearly absurd (since
$\lambda_{i}\neq\lambda_{i}$ can never happen). Thus, we have found a contradiction.

Forget that we fixed $\lambda$. We thus have obtained a contradiction for each
partition $\lambda$ such that $\lambda<\lambda$. Hence, there exists no such
$\lambda$. In other words, the relation $<$ is irreflexive. This proves Claim 1.
\end{proof}

\begin{proof}
[Proof of Claim 2.]Let $\lambda$, $\mu$ and $\nu$ be three partitions that
satisfy $\lambda<\mu$ and $\mu<\nu$. We must prove that $\lambda<\nu$.

Indeed, we have $\lambda<\mu$. By the definition of the relation $<$, this
means that

\begin{itemize}
\item there exists at least one $i\geq1$ such that $\lambda_{i}\neq\mu_{i}$, and

\item the \textbf{smallest} such $i$ satisfies $\lambda_{i}<\mu_{i}$.
\end{itemize}

Let us denote this smallest $i$ by $p$. Thus, we have $\lambda_{p}<\mu_{p}$,
whereas%
\begin{equation}
\text{each }i<p\text{ satisfies }\lambda_{i}=\mu_{i}
\label{pf.prop.bas.tord.pars.tord.c2.pf.1}%
\end{equation}
(since $p$ is the \textbf{smallest} $i\geq1$ that satisfies $\lambda_{i}%
\neq\mu_{i}$).

Furthermore, we have $\mu<\nu$. By the definition of the relation $<$, this
means that

\begin{itemize}
\item there exists at least one $i\geq1$ such that $\mu_{i}\neq\nu_{i}$, and

\item the \textbf{smallest} such $i$ satisfies $\mu_{i}<\nu_{i}$.
\end{itemize}

Let us denote this smallest $i$ by $q$. Thus, we have $\mu_{q}<\nu_{q}$,
whereas%
\begin{equation}
\text{each }i<q\text{ satisfies }\mu_{i}=\nu_{i}
\label{pf.prop.bas.tord.pars.tord.c2.pf.3}%
\end{equation}
(since $q$ is the \textbf{smallest} $i\geq1$ that satisfies $\mu_{i}\neq
\nu_{i}$). Hence, it is easy to see that%
\begin{equation}
\text{each }i\leq q\text{ satisfies }\mu_{i}\leq\nu_{i}
\label{pf.prop.bas.tord.pars.tord.c2.pf.4}%
\end{equation}
\footnote{\textit{Proof of (\ref{pf.prop.bas.tord.pars.tord.c2.pf.4}):} Let
$i$ be an integer such that $i\leq q$. We must show that $\mu_{i}\leq\nu_{i}$.
If $i<q$, then this follows from (\ref{pf.prop.bas.tord.pars.tord.c2.pf.3})
(since (\ref{pf.prop.bas.tord.pars.tord.c2.pf.3}) says $\mu_{i}=\nu_{i}$ in
this case, and hence $\mu_{i}\leq\nu_{i}$). Hence, we WLOG assume that we
don't have $i<q$. Thus, $i\geq q$. Combining this with $i\leq q$, we obtain
$i=q$. Hence, $\mu_{i}=\mu_{q}<\nu_{q}=\nu_{i}$ (since $q=i$). Therefore,
$\mu_{i}\leq\nu_{i}$. Thus, (\ref{pf.prop.bas.tord.pars.tord.c2.pf.4}) is
proved.}.

Now, let $w:=\min\left\{  p,q\right\}  $. Then, $w\leq p$ and $w\leq q$.
Hence,%
\begin{equation}
\text{each }i<w\text{ satisfies }\lambda_{i}=\nu_{i}
\label{pf.prop.bas.tord.pars.tord.c2.pf.5}%
\end{equation}
\footnote{\textit{Proof.} Let $i\geq1$ be such that $i<w$. Then,%
\begin{align*}
\lambda_{i}  &  =\mu_{i}\ \ \ \ \ \ \ \ \ \ \left(  \text{by
(\ref{pf.prop.bas.tord.pars.tord.c2.pf.1}), since }i<w\leq p\right) \\
&  =\nu_{i}\ \ \ \ \ \ \ \ \ \ \left(  \text{by
(\ref{pf.prop.bas.tord.pars.tord.c2.pf.3}), since }i<w\leq q\right)  .
\end{align*}
This proves (\ref{pf.prop.bas.tord.pars.tord.c2.pf.5}).}. Furthermore,
\begin{equation}
\lambda_{w}<\nu_{w} \label{pf.prop.bas.tord.pars.tord.c2.pf.7}%
\end{equation}
\footnote{\textit{Proof of (\ref{pf.prop.bas.tord.pars.tord.c2.pf.7}):} We are
in one of the following two cases:
\par
\textit{Case 1:} We have $p\leq q$.
\par
\textit{Case 2:} We have $p>q$.
\par
Let us first consider Case 1. In this case, we have $p\leq q$. Hence,
$w=\min\left\{  p,q\right\}  =p$ (since $p\leq q$). Applying
(\ref{pf.prop.bas.tord.pars.tord.c2.pf.4}) to $i=w$, we obtain $\mu_{w}\leq
\nu_{w}$ (since $w\leq q$). Recall that $\lambda_{p}<\mu_{p}$. Since $w=p$, we
can rewrite this as $\lambda_{w}<\mu_{w}$. Hence, $\lambda_{w}<\mu_{w}\leq
\nu_{w}$. Hence, (\ref{pf.prop.bas.tord.pars.tord.c2.pf.7}) is proved in Case
1.
\par
Let us next consider Case 2. In this case, we have $p>q$. Hence,
$w=\min\left\{  p,q\right\}  =q$ (since $p>q$). Thus, $w=q<p$ (since $p>q$)
and therefore $\lambda_{w}=\mu_{w}$ (by
(\ref{pf.prop.bas.tord.pars.tord.c2.pf.1}), applied to $i=w$). Recall that
$\mu_{q}<\nu_{q}$. Since $w=q$, we can rewrite this as $\mu_{w}<\nu_{w}$.
Hence, $\lambda_{w}=\mu_{w}<\nu_{w}$. Hence,
(\ref{pf.prop.bas.tord.pars.tord.c2.pf.7}) is proved in Case 2.
\par
Thus, (\ref{pf.prop.bas.tord.pars.tord.c2.pf.7}) is proved in both Cases 1 and
2.} and therefore $\lambda_{w}\neq\nu_{w}$. Hence, $w$ is an $i\geq1$ such
that $\lambda_{i}\neq\nu_{i}$. Moreover, $w$ is the \textbf{smallest} such $i$
(since (\ref{pf.prop.bas.tord.pars.tord.c2.pf.5}) shows that all smaller $i$'s
satisfy $\lambda_{i}=\nu_{i}$). Hence, we have now shown that

\begin{itemize}
\item there exists at least one $i\geq1$ such that $\lambda_{i}\neq\nu_{i}$
(namely, $i=w$), and

\item the \textbf{smallest} such $i$ satisfies $\lambda_{i}<\nu_{i}$ (because
the smallest such $i$ is $w$, and we know that $w$ satisfies $\lambda_{w}%
<\nu_{w}$).
\end{itemize}

In other words, $\lambda<\nu$ (by the definition of the relation $<$). This
proves Claim 2.
\end{proof}

\begin{proof}
[Proof of Claim 3.]Any binary relation that is irreflexive and transitive is
automatically asymmetric (since $\lambda<\mu$ and $\mu<\lambda$ would entail
$\lambda<\lambda$ by transitivity, but this would contradict irreflexivity).
Since the relation $<$ is irreflexive (by Claim 1) and transitive (by Claim
2), we thus conclude that it is also asymmetric. This proves Claim 3.
\end{proof}

Having proved Claims 1, 2 and 3, we thus conclude that the relation $<$ is the
smaller relation of a partial order. It remains to show that this partial
order is a total order. For this purpose, we need to prove the following last claim:

\begin{statement}
\textit{Claim 4:} Let $\lambda$ and $\mu$ be two distinct partitions. Then, we
have $\lambda<\mu$ or $\mu<\lambda$.
\end{statement}

\begin{proof}
[Proof of Claim 4.]We shall first show that there exists some $i\geq1$ such
that $\lambda_{i}\neq\mu_{i}$. Indeed, we must be in one of the three cases
$\ell\left(  \lambda\right)  >\ell\left(  \mu\right)  $ and $\ell\left(
\mu\right)  >\ell\left(  \lambda\right)  $ and $\ell\left(  \lambda\right)
=\ell\left(  \mu\right)  $, but we can easily find such an $i$ in each of
these three cases:

\begin{itemize}
\item If $\ell\left(  \lambda\right)  >\ell\left(  \mu\right)  $, then it is
easy to see that $\ell\left(  \lambda\right)  $ is such an $i$%
\ \ \ \ \footnote{\textit{Proof.} Assume that $\ell\left(  \lambda\right)
>\ell\left(  \mu\right)  $. We must show that $\ell\left(  \lambda\right)  $
is an $i\geq1$ such that $\lambda_{i}\neq\mu_{i}$. In other words, we must
prove that $\ell\left(  \lambda\right)  \geq1$ and $\lambda_{\ell\left(
\lambda\right)  }\neq\mu_{\ell\left(  \lambda\right)  }$.
\par
Indeed, we have $\ell\left(  \lambda\right)  >\ell\left(  \mu\right)  \geq0$
and thus $\ell\left(  \lambda\right)  \geq1$ (since $\ell\left(
\lambda\right)  $ is an integer). Moreover, from $\ell\left(  \lambda\right)
>\ell\left(  \mu\right)  $, we obtain $\mu_{\ell\left(  \lambda\right)  }=0$
(since each integer $i>\ell\left(  \mu\right)  $ satisfies $\mu_{i}=0$
according to Definition \ref{def.partitions.partitions} \textbf{(h)}). But
$\lambda_{\ell\left(  \lambda\right)  }>0$ (since $\lambda_{\ell\left(
\lambda\right)  }$ is the last entry of the partition $\lambda$, and thus is a
positive integer). Hence, $\lambda_{\ell\left(  \lambda\right)  }>0=\mu
_{\ell\left(  \lambda\right)  }$, so that $\lambda_{\ell\left(  \lambda
\right)  }\neq\mu_{\ell\left(  \lambda\right)  }$. This completes our proof.}.

\item If $\ell\left(  \mu\right)  >\ell\left(  \lambda\right)  $, then it is
easy to see that $\ell\left(  \mu\right)  $ is such an $i$%
\ \ \ \ \footnote{\textit{Proof.} The argument is the same as for the previous
sentence, but the roles of $\lambda$ and $\mu$ are switched.}.

\item If $\ell\left(  \lambda\right)  =\ell\left(  \mu\right)  $, then such an
$i$ also exists\footnote{\textit{Proof.} Assume that $\ell\left(
\lambda\right)  =\ell\left(  \mu\right)  $. Thus, $\ell\left(  \mu\right)
=\ell\left(  \lambda\right)  $. But $\lambda\neq\mu$ (since $\lambda$ and
$\mu$ are distinct). In light of $\lambda=\left(  \lambda_{1},\lambda
_{2},\ldots,\lambda_{\ell\left(  \lambda\right)  }\right)  $ and $\mu=\left(
\mu_{1},\mu_{2},\ldots,\mu_{\ell\left(  \mu\right)  }\right)  =\left(  \mu
_{1},\mu_{2},\ldots,\mu_{\ell\left(  \lambda\right)  }\right)  $ (since
$\ell\left(  \mu\right)  =\ell\left(  \lambda\right)  $), we can rewrite this
as%
\[
\left(  \lambda_{1},\lambda_{2},\ldots,\lambda_{\ell\left(  \lambda\right)
}\right)  \neq\left(  \mu_{1},\mu_{2},\ldots,\mu_{\ell\left(  \lambda\right)
}\right)  .
\]
Hence, there exists some $i\in\left[  \ell\left(  \lambda\right)  \right]  $
such that $\lambda_{i}\neq\mu_{i}$ (because two distinct $\ell\left(
\lambda\right)  $-tuples must differ in at least one position). Therefore,
there exists some $i\geq1$ such that $\lambda_{i}\neq\mu_{i}$ (because
$i\in\left[  \ell\left(  \lambda\right)  \right]  $ entails $i\geq1$).}.
\end{itemize}

\noindent Hence, we have shown that there exists some $i\geq1$ such that
$\lambda_{i}\neq\mu_{i}$.

Consider the \textbf{smallest} such $i$. Then, this $i$ satisfies $\lambda
_{i}\neq\mu_{i}$. Hence, either $\lambda_{i}<\mu_{i}$ or $\lambda_{i}>\mu_{i}%
$. In the former case, we have $\lambda<\mu$ (by the definition of the
relation $<$). In the latter case, we have $\mu<\lambda$ (by the definition of
the relation $<$ again\footnote{\textit{Proof.} Assume that $\lambda_{i}%
>\mu_{i}$. Thus, $\mu_{i}<\lambda_{i}$. But recall that our $i$ is the
\textbf{smallest} $i\geq1$ such that $\lambda_{i}\neq\mu_{i}$. In other words,
our $i$ is the \textbf{smallest} $i\geq1$ such that $\mu_{i}\neq\lambda_{i}$.
Hence, there exists an $i\geq1$ such that $\mu_{i}\neq\lambda_{i}$, and the
\textbf{smallest} such $i$ satisfies $\mu_{i}<\lambda_{i}$. In other words, we
have $\mu<\lambda$ (by the definition of the relation $<$), qed.}). Thus, we
have found that $\lambda<\mu$ or $\mu<\lambda$. This proves Claim 4.
\end{proof}

Claim 4 shows that the partial order whose smaller relation is $<$ is actually
a total order. The proof of Proposition \ref{prop.bas.tord.pars.tord} is thus complete.
\end{proof}

\begin{proof}
[Detailed proof of Proposition \ref{prop.bas.tord.syt.tord}.]First, we shall
show that the relation $<$ is the smaller relation of a partial order on
$\bigcup\limits_{\kappa\vdash n}\operatorname*{SYT}\left(  \kappa\right)  $.
In order to prove this, we need to verify three claims:

\begin{statement}
\textit{Claim 1:} The relation $<$ is irreflexive (i.e., there exists no
$P\in\bigcup\limits_{\kappa\vdash n}\operatorname*{SYT}\left(  \kappa\right)
$ such that $P<P$).
\end{statement}

\begin{statement}
\textit{Claim 2:} The relation $<$ is transitive (i.e., if three tableaux
$P,Q,R\in\bigcup\limits_{\kappa\vdash n}\operatorname*{SYT}\left(
\kappa\right)  $ satisfy $P<Q$ and $Q<R$, then $P<R$).
\end{statement}

\begin{statement}
\textit{Claim 3:} The relation $<$ is asymmetric (i.e., no two tableaux $P$
and $Q$ satisfy $P<Q$ and $Q<P$ at the same time).
\end{statement}

\begin{proof}
[Proof of Claim 1.]Let $P\in\bigcup\limits_{\kappa\vdash n}\operatorname*{SYT}%
\left(  \kappa\right)  $ be a tableau such that $P<P$. We shall derive a contradiction.

Indeed, we have $P\in\bigcup\limits_{\kappa\vdash n}\operatorname*{SYT}\left(
\kappa\right)  $. Hence, $P\in\operatorname*{SYT}\left(  \lambda\right)  $ for
some partition $\lambda$. Consider this $\lambda$. Now, we have assumed that
$P<P$. By the definition of the relation $<$, this means that%
\[
\left(  \lambda<\lambda\text{ or }\left(  \lambda=\lambda\text{ and }%
\overline{P}<\overline{P}\right)  \right)
\]
(since $P\in\operatorname*{SYT}\left(  \lambda\right)  $). But this is clearly
false (since neither $\lambda<\lambda$ nor $\overline{P}<\overline{P}$ can
hold (because both relations \textquotedblleft$<$\textquotedblright\ in these
two inequalities are total orders)). Thus, we have found a contradiction.

Forget that we fixed $P$. We thus have obtained a contradiction for each
tableau $P$ such that $P<P$. Hence, there exists no such $P$. In other words,
the relation $<$ is irreflexive. This proves Claim 1.
\end{proof}

\begin{proof}
[Proof of Claim 2.]Let $P,Q,R\in\bigcup\limits_{\kappa\vdash n}%
\operatorname*{SYT}\left(  \kappa\right)  $ be three tableaux that satisfy
$P<Q$ and $Q<R$. We must prove that $P<R$.

Indeed, we have $P\in\bigcup\limits_{\kappa\vdash n}\operatorname*{SYT}\left(
\kappa\right)  $. Hence, we have $P\in\operatorname*{SYT}\left(
\lambda\right)  $ for some partition $\lambda$. Similarly, we have
$Q\in\operatorname*{SYT}\left(  \mu\right)  $ for some partition $\mu$.
Likewise, we have $R\in\operatorname*{SYT}\left(  \nu\right)  $ for some
partition $\nu$. Consider these three partitions $\lambda,\mu,\nu$.

Now, we have assumed that $P<Q$. By the definition of the relation $<$ on
$\bigcup\limits_{\kappa\vdash n}\operatorname*{SYT}\left(  \kappa\right)  $,
this means that%
\[
\left(  \lambda<\mu\text{ or }\left(  \lambda=\mu\text{ and }\overline
{P}<\overline{Q}\right)  \right)
\]
(since $P\in\operatorname*{SYT}\left(  \lambda\right)  $ and $Q\in
\operatorname*{SYT}\left(  \mu\right)  $). Thus, we must have either
$\lambda<\mu$ or $\left(  \lambda=\mu\text{ and }\overline{P}<\overline
{Q}\right)  $. In both of these cases, we have $\lambda\leq\mu$ (since
$\lambda\leq\mu$ follows both from $\lambda<\mu$ and from $\lambda=\mu$).

So we have derived $\lambda\leq\mu$ from $P<Q$. Likewise, from $Q<R$, we can
derive $\mu\leq\nu$. But the lexicographic order on partitions is a total
order; thus, the relation $\leq$ on the set of all partitions is transitive.
Hence, from $\lambda\leq\mu$ and $\mu\leq\nu$, we obtain $\lambda\leq\nu$. In
other words, we have either $\lambda=\nu$ or $\lambda<\nu$. Hence, we are in
one of the following two cases:

\textit{Case 1:} We have $\lambda=\nu$.

\textit{Case 2:} We have $\lambda<\nu$.

Let us first consider Case 1. In this case, we have $\lambda=\nu$. Thus,
$\nu=\lambda\leq\mu$. Combining this with $\mu\leq\nu$, we obtain $\mu=\nu$
(since the lexicographic order on partitions is a total order, so that the
relation $\leq$ on the set of all partitions is antisymmetric). Hence,
$\mu=\nu=\lambda$, so that $\lambda=\mu$.

Now, recall that we have either $\lambda<\mu$ or $\left(  \lambda=\mu\text{
and }\overline{P}<\overline{Q}\right)  $. Since $\lambda<\mu$ is impossible
(because $\lambda=\mu$), we thus conclude that $\left(  \lambda=\mu\text{ and
}\overline{P}<\overline{Q}\right)  $ must hold. Hence, $\overline{P}%
<\overline{Q}$.

Moreover, we have assumed that $Q<R$. By the definition of the relation $<$ on
$\bigcup\limits_{\kappa\vdash n}\operatorname*{SYT}\left(  \kappa\right)  $,
this means that%
\[
\left(  \mu<\nu\text{ or }\left(  \mu=\nu\text{ and }\overline{Q}<\overline
{R}\right)  \right)
\]
(since $Q\in\operatorname*{SYT}\left(  \mu\right)  $ and $R\in
\operatorname*{SYT}\left(  \nu\right)  $). Since $\mu<\nu$ is impossible
(because $\mu=\nu$), we thus must have $\left(  \mu=\nu\text{ and }%
\overline{Q}<\overline{R}\right)  $. Hence, $\overline{Q}<\overline{R}$.

But the Young last letter order on the $n$-tabloids is a total order, thus
transitive. Hence, from $\overline{P}<\overline{Q}$ and $\overline
{Q}<\overline{R}$, we obtain $\overline{P}<\overline{R}$. Combining this with
$\lambda=\nu$, we obtain $\left(  \lambda=\nu\text{ and }\overline
{P}<\overline{R}\right)  $. Hence, we have%
\[
\left(  \lambda<\nu\text{ or }\left(  \lambda=\nu\text{ and }\overline
{P}<\overline{R}\right)  \right)  .
\]
But this is equivalent to saying that $P<R$ (by the definition of the relation
$<$ on $\bigcup\limits_{\kappa\vdash n}\operatorname*{SYT}\left(
\kappa\right)  $, because $P\in\operatorname*{SYT}\left(  \lambda\right)  $
and $R\in\operatorname*{SYT}\left(  \nu\right)  $). Thus, $P<R$ is proved in
Case 1.

Let us now consider Case 2. In this case, we have $\lambda<\nu$. Hence, we
have%
\[
\left(  \lambda<\nu\text{ or }\left(  \lambda=\nu\text{ and }\overline
{P}<\overline{R}\right)  \right)  .
\]
But this is equivalent to saying that $P<R$ (by the definition of the relation
$<$ on $\bigcup\limits_{\kappa\vdash n}\operatorname*{SYT}\left(
\kappa\right)  $, because $P\in\operatorname*{SYT}\left(  \lambda\right)  $
and $R\in\operatorname*{SYT}\left(  \nu\right)  $). Thus, $P<R$ is proved in
Case 2.

We have now proved $P<R$ in both Cases 1 and 2. Hence, $P<R$ always holds.
This proves Claim 2.
\end{proof}

\begin{proof}
[Proof of Claim 3.]Any binary relation that is irreflexive and transitive is
automatically asymmetric (since $P<Q$ and $Q<P$ would entail $P<P$ by
transitivity, but this would contradict irreflexivity). Since the relation $<$
is irreflexive (by Claim 1) and transitive (by Claim 2), we thus conclude that
it is also asymmetric. This proves Claim 3.
\end{proof}

Having proved Claims 1, 2 and 3, we thus conclude that the relation $<$ is the
smaller relation of a partial order. It remains to show that this partial
order is a total order. For this purpose, we need to prove the following last claim:

\begin{statement}
\textit{Claim 4:} Let $P$ and $Q$ be two distinct tableaux in $\bigcup
\limits_{\kappa\vdash n}\operatorname*{SYT}\left(  \kappa\right)  $. Then, we
have $P<Q$ or $Q<P$.
\end{statement}

\begin{proof}
[Proof of Claim 4.]We have $P\in\bigcup\limits_{\kappa\vdash n}%
\operatorname*{SYT}\left(  \kappa\right)  $. Hence, we have $P\in
\operatorname*{SYT}\left(  \lambda\right)  $ for some partition $\lambda$.
Similarly, we have $Q\in\operatorname*{SYT}\left(  \mu\right)  $ for some
partition $\mu$. Consider these two partitions $\lambda$ and $\mu$.

The lexicographic order on partitions is a total order. Hence, we have either
$\lambda<\mu$ or $\lambda=\mu$ or $\lambda>\mu$. In other words, we are in one
of the following three cases:

\textit{Case 1:} We have $\lambda<\mu$.

\textit{Case 2:} We have $\lambda=\mu$.

\textit{Case 3:} We have $\lambda>\mu$.

Let us consider Case 1 first. In this case, we have $\lambda<\mu$. Thus,%
\[
\left(  \lambda<\mu\text{ or }\left(  \lambda=\mu\text{ and }\overline
{P}<\overline{Q}\right)  \right)  .
\]
But this is equivalent to saying that $P<Q$ (by the definition of the relation
$<$ on $\bigcup\limits_{\kappa\vdash n}\operatorname*{SYT}\left(
\kappa\right)  $, since $P\in\operatorname*{SYT}\left(  \lambda\right)  $ and
$Q\in\operatorname*{SYT}\left(  \mu\right)  $). Thus, we have $P<Q$, and
therefore $\left(  P<Q\text{ or }Q<P\right)  $. Thus, Claim 4 is proved in
Case 1.

Now consider Case 2. In this case, we have $\lambda=\mu$. Hence,
$\operatorname*{SYT}\left(  \lambda\right)  =\operatorname*{SYT}\left(
\mu\right)  $. Thus, both $P$ and $Q$ are standard $n$-tableaux of shape
$Y\left(  \lambda\right)  $ (since $P\in\operatorname*{SYT}\left(
\lambda\right)  $ and $Q\in\operatorname*{SYT}\left(  \mu\right)
=\operatorname*{SYT}\left(  \lambda\right)  $). Thus, $P$ and $Q$ are
row-standard $n$-tableaux of shape $Y\left(  \lambda\right)  $ (since standard
tableaux are always row-standard). Proposition \ref{prop.tabloid.row-st}
(applied to $D=Y\left(  \lambda\right)  $) shows that there is a bijection%
\begin{align*}
&  \text{from }\left\{  \text{row-standard }n\text{-tableaux of shape
}Y\left(  \lambda\right)  \right\} \\
&  \text{to }\left\{  n\text{-tabloids of shape }Y\left(  \lambda\right)
\right\}
\end{align*}
that sends each row-standard $n$-tableau $T$ to its tabloid $\overline{T}$.
Hence, in particular, this bijection is injective. In other words, if $T_{1}$
and $T_{2}$ are two distinct row-standard $n$-tableaux of shape $Y\left(
\lambda\right)  $, then the $n$-tabloids $\overline{T_{1}}$ and $\overline
{T_{2}}$ are also distinct. Applying this to $T_{1}=P$ and $T_{2}=Q$, we
conclude that the $n$-tabloids $\overline{P}$ and $\overline{Q}$ are distinct
(since the $n$-tableaux $P$ and $Q$ are row-standard and distinct). Since the
Young last letter order on $\left\{  n\text{-tabloids of shape }Y\left(
\lambda\right)  \right\}  $ is a total order, we thus conclude that either
$\overline{P}<\overline{Q}$ or $\overline{Q}<\overline{P}$. We WLOG assume
that $\overline{P}<\overline{Q}$ holds, since in the other case ($\overline
{Q}<\overline{P}$) we can achieve this by swapping $P$ with $Q$. Thus, we have
$\lambda=\mu$ and $\overline{P}<\overline{Q}$. Thus,%
\[
\left(  \lambda<\mu\text{ or }\left(  \lambda=\mu\text{ and }\overline
{P}<\overline{Q}\right)  \right)  .
\]
But this is equivalent to saying that $P<Q$ (by the definition of the relation
$<$ on $\bigcup\limits_{\kappa\vdash n}\operatorname*{SYT}\left(
\kappa\right)  $, since $P\in\operatorname*{SYT}\left(  \lambda\right)  $ and
$Q\in\operatorname*{SYT}\left(  \mu\right)  $). Thus, we have $P<Q$, and
therefore $\left(  P<Q\text{ or }Q<P\right)  $. Thus, Claim 4 is proved in
Case 2.

It remains to handle Case 3. In this case, we have $\lambda>\mu$. This is
equivalent to $\mu<\lambda$. Thus, Case 3 is the same case as Case 1, except
that the roles of $P$ and $Q$ (and thus the roles of $\lambda$ and $\mu$) are
switched. Hence, the proof of Claim 4 in Case 3 is analogous to the proof in
Case 1 (since Claim 4 is symmetric in $P$ and $Q$).

We have now proved Claim 4 in all three Cases 1, 2 and 3; thus, its proof is complete.
\end{proof}

Claim 4 shows that the partial order whose smaller relation is $<$ is actually
a total order. The proof of Proposition \ref{prop.bas.tord.syt.tord} is thus complete.
\end{proof}

\begin{proof}
[Proof of Proposition \ref{prop.bas.tord.sbt.tord}.]First, we shall show that
the relation $<$ is the smaller relation of a partial order on
$\operatorname*{SBT}\left(  n\right)  $. In order to prove this, we need to
verify three claims:\footnote{We are using capital Fraktur letters like
$\mathfrak{P},\mathfrak{Q},\mathfrak{R}$ for standard $n$-bitableaux.}

\begin{statement}
\textit{Claim 1:} The relation $<$ is irreflexive (i.e., there exists no
$\mathfrak{P}\in\operatorname*{SBT}\left(  n\right)  $ such that
$\mathfrak{P}<\mathfrak{P}$).
\end{statement}

\begin{statement}
\textit{Claim 2:} The relation $<$ is transitive (i.e., if three
$n$-bitableaux $\mathfrak{P},\mathfrak{Q},\mathfrak{R}\in\operatorname*{SBT}%
\left(  n\right)  $ satisfy $\mathfrak{P}<\mathfrak{Q}$ and $\mathfrak{Q}%
<\mathfrak{R}$, then $\mathfrak{P}<\mathfrak{R}$).
\end{statement}

\begin{statement}
\textit{Claim 3:} The relation $<$ is asymmetric (i.e., no two $n$-bitableaux
$\mathfrak{P}$ and $\mathfrak{Q}$ satisfy $\mathfrak{P}<\mathfrak{Q}$ and
$\mathfrak{Q}<\mathfrak{P}$ at the same time).
\end{statement}

\begin{proof}
[Proof of Claim 1.]Let $\mathfrak{P}\in\operatorname*{SBT}\left(  n\right)  $
be an $n$-bitableau such that $\mathfrak{P}<\mathfrak{P}$. We shall derive a contradiction.

Indeed, let us write the $n$-bitableau $\mathfrak{P}$ as $\mathfrak{P}=\left(
\lambda,U,V\right)  $. Now, we have assumed that $\mathfrak{P}<\mathfrak{P}$.
In other words, $\left(  \lambda,U,V\right)  <\left(  \lambda,U,V\right)  $
(since $\mathfrak{P}=\left(  \lambda,U,V\right)  $). By the definition of the
relation $<$, this means that%
\[
\left(  U<U\text{ or }\left(  U=U\text{ and }V<V\right)  \right)  .
\]
But this is clearly false (since neither $U<U$ nor $V<V$ is true (because the
relation $<$ on $\bigcup\limits_{\kappa\vdash n}\operatorname*{SYT}\left(
\kappa\right)  $ is a total order)). Thus, we have found a contradiction.

Forget that we fixed $\mathfrak{P}$. We thus have obtained a contradiction for
each $n$-bitableau $\mathfrak{P}$ such that $\mathfrak{P}<\mathfrak{P}$.
Hence, there exists no such $\mathfrak{P}$. In other words, the relation $<$
is irreflexive. This proves Claim 1.
\end{proof}

\begin{proof}
[Proof of Claim 2.]Let $\mathfrak{P},\mathfrak{Q},\mathfrak{R}\in
\operatorname*{SBT}\left(  n\right)  $ be three $n$-bitableaux that satisfy
$\mathfrak{P}<\mathfrak{Q}$ and $\mathfrak{Q}<\mathfrak{R}$. We must prove
that $\mathfrak{P}<\mathfrak{R}$.

Indeed, let us write the $n$-bitableaux $\mathfrak{P},\mathfrak{Q}%
,\mathfrak{R}$ as
\[
\mathfrak{P}=\left(  \lambda,A,B\right)  ,\ \ \ \ \ \ \ \ \ \ \mathfrak{Q}%
=\left(  \mu,C,D\right)  ,\ \ \ \ \ \ \ \ \ \ \mathfrak{R}=\left(
\nu,E,F\right)  ,
\]
respectively. Now, recall that $\mathfrak{P}<\mathfrak{Q}$. In other words,
$\left(  \lambda,A,B\right)  <\left(  \mu,C,D\right)  $ (since $\mathfrak{P}%
=\left(  \lambda,A,B\right)  $ and $\mathfrak{Q}=\left(  \mu,C,D\right)  $).
By the definition of the relation $<$, this means that%
\[
\left(  A<C\text{ or }\left(  A=C\text{ and }B<D\right)  \right)  .
\]
So we must have $A<C$ or $\left(  A=C\text{ and }B<D\right)  $. Likewise, from
$\mathfrak{Q}<\mathfrak{R}$, we conclude that we must have $C<E$ or $\left(
C=E\text{ and }D<F\right)  $.

We just have shown that we have $A<C$ or $\left(  A=C\text{ and }B<D\right)
$. In both of these cases, we have $A\leq C$ (since $A\leq C$ follows both
from $A<C$ and from $A=C$).

So we have derived $A\leq C$ from $\mathfrak{P}<\mathfrak{Q}$. Likewise, from
$\mathfrak{Q}<\mathfrak{R}$, we can derive $C\leq E$. But the relation $\leq$
on the set $\bigcup\limits_{\kappa\vdash n}\operatorname*{SYT}\left(
\kappa\right)  $ is transitive (since it belongs to a total order). Hence,
from $A\leq C$ and $C\leq E$, we obtain $A\leq E$. In other words, we have
either $A=E$ or $A<E$. Hence, we are in one of the following two cases:

\textit{Case 1:} We have $A=E$.

\textit{Case 2:} We have $A<E$.

Let us first consider Case 1. In this case, we have $A=E$. Thus, $E=A\leq C$.
Combining this with $C\leq E$, we obtain $C=E$ (since the relation $\leq$ on
$\bigcup\limits_{\kappa\vdash n}\operatorname*{SYT}\left(  \kappa\right)  $
belongs to a total order and thus is antisymmetric). Hence, $C=E=A$, so that
$A=C$.

Now, recall that we have either $A<C$ or $\left(  A=C\text{ and }B<D\right)
$. Since $A<C$ is impossible (because $A=C$), we thus conclude that $\left(
A=C\text{ and }B<D\right)  $ must hold. Hence, $B<D$.

Recall that we showed that we must have $C<E$ or $\left(  C=E\text{ and
}D<F\right)  $. Since $C<E$ is impossible (because $C=E$), we thus must have
$\left(  C=E\text{ and }D<F\right)  $. Hence, $D<F$.

But the relation $<$ on the set $\bigcup\limits_{\kappa\vdash n}%
\operatorname*{SYT}\left(  \kappa\right)  $ is a total order, thus transitive.
Hence, from $B<D$ and $D<F$, we obtain $B<F$. Combining this with $A=E$, we
obtain $\left(  A=E\text{ and }B<F\right)  $. Hence, we have%
\[
\left(  A<E\text{ or }\left(  A=E\text{ and }B<F\right)  \right)  .
\]
But this is equivalent to saying that $\left(  \lambda,A,B\right)  <\left(
\nu,E,F\right)  $ (by the definition of the relation $<$ on
$\operatorname*{SBT}\left(  n\right)  $). In other words, $\mathfrak{P}%
<\mathfrak{R}$ (since $\mathfrak{P}=\left(  \lambda,A,B\right)  $ and
$\mathfrak{R}=\left(  \nu,E,F\right)  $). Thus, Claim 2 is proved in Case 1.

Let us now consider Case 2. In this case, we have $A<E$. Hence, we have%
\[
\left(  A<E\text{ or }\left(  A=E\text{ and }B<F\right)  \right)  .
\]
But this is equivalent to saying that $\left(  \lambda,A,B\right)  <\left(
\nu,E,F\right)  $ (by the definition of the relation $<$ on
$\operatorname*{SBT}\left(  n\right)  $). In other words, $\mathfrak{P}%
<\mathfrak{R}$ (since $\mathfrak{P}=\left(  \lambda,A,B\right)  $ and
$\mathfrak{R}=\left(  \nu,E,F\right)  $). Thus, Claim 2 is proved in Case 2.

We have now proved Claim 2 in both Cases 1 and 2. Hence, the proof of Claim 2
is complete.
\end{proof}

\begin{proof}
[Proof of Claim 3.]Any binary relation that is irreflexive and transitive is
automatically asymmetric (since $\mathfrak{P}<\mathfrak{Q}$ and $\mathfrak{Q}%
<\mathfrak{P}$ would entail $\mathfrak{P}<\mathfrak{P}$ by transitivity, but
this would contradict irreflexivity). Since the relation $<$ is irreflexive
(by Claim 1) and transitive (by Claim 2), we thus conclude that it is also
asymmetric. This proves Claim 3.
\end{proof}

Having proved Claims 1, 2 and 3, we thus conclude that the relation $<$ is the
smaller relation of a partial order. It remains to show that this partial
order is a total order. For this purpose, we need to prove the following last claim:

\begin{statement}
\textit{Claim 4:} Let $\mathfrak{P}$ and $\mathfrak{Q}$ be two distinct
$n$-bitableaux in $\operatorname*{SBT}\left(  n\right)  $. Then, we have
$\mathfrak{P}<\mathfrak{Q}$ or $\mathfrak{Q}<\mathfrak{P}$.
\end{statement}

\begin{proof}
[Proof of Claim 4.]Let us write the $n$-bitableaux $\mathfrak{P}$ and
$\mathfrak{Q}$ as
\[
\mathfrak{P}=\left(  \lambda,A,B\right)  \ \ \ \ \ \ \ \ \ \ \text{and}%
\ \ \ \ \ \ \ \ \ \ \mathfrak{Q}=\left(  \mu,C,D\right)  ,
\]
respectively. We are clearly in one of the following two cases:

\textit{Case 1:} We have $A=C$.

\textit{Case 2:} We have $A\neq C$.

Let us consider Case 1 first. In this case, we have $A=C$. But $A=C\in
\operatorname*{SYT}\left(  \mu\right)  $ (since $\left(  \mu,C,D\right)  $ is
a standard $n$-bitableau) and $A\in\operatorname*{SYT}\left(  \lambda\right)
$ (since $\left(  \lambda,A,B\right)  $ is a standard $n$-bitableau). Hence,
the sets $\operatorname*{SYT}\left(  \lambda\right)  $ and
$\operatorname*{SYT}\left(  \mu\right)  $ have at least one element in common
(namely, $A$). Since the sets $\operatorname*{SYT}\left(  \kappa\right)  $ for
different partitions $\kappa$ are disjoint, we thus conclude that the
partitions $\lambda$ and $\mu$ must be identical. In other words, $\lambda
=\mu$.

If we also had $B=D$, then we would have $\left(  \lambda,A,B\right)  =\left(
\mu,C,D\right)  $ (since $\lambda=\mu$ and $A=C$ and $B=D$), so that
$\mathfrak{P}=\left(  \lambda,A,B\right)  =\left(  \mu,C,D\right)
=\mathfrak{Q}$; but this would contradict the fact that $\mathfrak{P}$ and
$\mathfrak{Q}$ are distinct (by assumption). Hence, we cannot have $B=D$.
Thus, $B\neq D$.

But the relation $<$ on the set $\bigcup\limits_{\kappa\vdash n}%
\operatorname*{SYT}\left(  \kappa\right)  $ is a total order. Hence, from
$B\neq D$, we conclude that either $B<D$ or $D<B$. These two cases are
completely analogous (indeed, we can switch between them by swapping
$\mathfrak{P}$ with $\mathfrak{Q}$; note that this switch does not affect the
equality $A=C$), so we can WLOG assume that we are in the first of them --
i.e., that we have $B<D$. Now, we have%
\[
\left(  A<C\text{ or }\left(  A=C\text{ and }B<D\right)  \right)
\]
(since we have $A=C$ and $B<D$). But this is equivalent to saying that
$\left(  \lambda,A,B\right)  <\left(  \mu,C,D\right)  $ (by the definition of
the relation $<$ on $\operatorname*{SBT}\left(  n\right)  $). Thus, we have
$\left(  \lambda,A,B\right)  <\left(  \mu,C,D\right)  $. In other words,
$\mathfrak{P}<\mathfrak{Q}$ (since $\mathfrak{P}=\left(  \lambda,A,B\right)  $
and $\mathfrak{Q}=\left(  \mu,C,D\right)  $). Therefore, $\left(
\mathfrak{P}<\mathfrak{Q}\text{ or }\mathfrak{Q}<\mathfrak{P}\right)  $. Thus,
Claim 4 is proved in Case 1.

Now consider Case 2. In this case, we have $A\neq C$. But the relation $<$ on
the set $\bigcup\limits_{\kappa\vdash n}\operatorname*{SYT}\left(
\kappa\right)  $ is a total order. Hence, from $A\neq C$, we conclude that
either $A<C$ or $C<A$. These two cases are completely analogous (indeed, we
can switch between them by swapping $\mathfrak{P}$ with $\mathfrak{Q}$), so we
can WLOG assume that we are in the first of them -- i.e., that we have $A<C$.
Thus, we have%
\[
\left(  A<C\text{ or }\left(  A=C\text{ and }B<D\right)  \right)  .
\]
But this is equivalent to saying that $\left(  \lambda,A,B\right)  <\left(
\mu,C,D\right)  $ (by the definition of the relation $<$ on
$\operatorname*{SBT}\left(  n\right)  $). Thus, we have $\left(
\lambda,A,B\right)  <\left(  \mu,C,D\right)  $. In other words, $\mathfrak{P}%
<\mathfrak{Q}$ (since $\mathfrak{P}=\left(  \lambda,A,B\right)  $ and
$\mathfrak{Q}=\left(  \mu,C,D\right)  $). Therefore, $\left(  \mathfrak{P}%
<\mathfrak{Q}\text{ or }\mathfrak{Q}<\mathfrak{P}\right)  $. Thus, Claim 4 is
proved in Case 2.

We have now proved Claim 4 in both Cases 1 and 2; thus, its proof is complete.
\end{proof}

Claim 4 shows that the partial order whose smaller relation is $<$ is actually
a total order. The proof of Proposition \ref{prop.bas.tord.sbt.tord} is thus complete.
\end{proof}

\begin{proof}
[Proof of Lemma \ref{lem.bas.tord.sbt.1}.]We have $\left(  \lambda,A,B\right)
<\left(  \mu,C,D\right)  $. By the definition of the relation $<$, this means
that%
\[
\left(  A<C\text{ or }\left(  A=C\text{ and }B<D\right)  \right)  .
\]
So we must have $A<C$ or $\left(  A=C\text{ and }B<D\right)  $. Thus, we must
have $A<C$ or $A=C$ (since the statement \textquotedblleft$\left(  A=C\text{
and }B<D\right)  $\textquotedblright\ implies the statement \textquotedblleft%
$A=C$\textquotedblright). In other words, $A\leq C$. This proves Lemma
\ref{lem.bas.tord.sbt.1} \textbf{(a)}. \medskip

\textbf{(b)} Recall that we must have $A<C$ or $\left(  A=C\text{ and
}B<D\right)  $. Hence, we must have $A<C$ or $B<D$ (since the statement
\textquotedblleft$\left(  A=C\text{ and }B<D\right)  $\textquotedblright%
\ implies the statement \textquotedblleft$B<D$\textquotedblright). This proves
Lemma \ref{lem.bas.tord.sbt.1} \textbf{(b)}.
\end{proof}
\end{fineprint}

\subsection{\label{sec.details.bas.mur}The Murphy bases}

\begin{fineprint}
Here we supply detailed proofs for some claims made in Section
\ref{sec.bas.mur}.

\begin{proof}
[Proof of Lemma \ref{lem.bas.mur.row-equ-as-Row-Row}.]Both $U$ and $V$ are
$n$-tableaux. Thus, all entries of $U$ are distinct, and all entries of $V$
are distinct.{}\medskip

\textbf{(a)} Assume that we have $\left(  \operatorname*{Row}\left(
i,U\right)  =\operatorname*{Row}\left(  i,V\right)  \text{ for all }%
i\in\mathbb{Z}\right)  $.

Fix $i\in\mathbb{Z}$. Then, $\operatorname*{Row}\left(  i,U\right)
=\operatorname*{Row}\left(  i,V\right)  $ (by the assumption we just made). In
other words, the set of entries in the $i$-th row of $U$ equals the set of
entries in the $i$-th row of $V$ (since these two sets are precisely
$\operatorname*{Row}\left(  i,U\right)  $ and $\operatorname*{Row}\left(
i,V\right)  $). We can furthermore replace the word \textquotedblleft
set\textquotedblright\ by \textquotedblleft multiset\textquotedblright\ here
(since the entries in the $i$-th row of $U$ are distinct (because all entries
of $U$ are distinct), and so are the entries in the $i$-th row of $V$, and
therefore the multisets of these entries do not have any nontrivial
multiplicities). Thus, we conclude that the multiset of entries in the $i$-th
row of $U$ equals the multiset of entries in the $i$-th row of $V$.

Forget that we fixed $i$. We thus have shown that for each $i\in\mathbb{Z}$,
the multiset of entries in the $i$-th row of $U$ equals the multiset of
entries in the $i$-th row of $V$. In other words, the $n$-tableaux $U$ and $V$
are row-equivalent (by the definition of \textquotedblleft
row-equivalent\textquotedblright). In other words, $\overline{U}=\overline{V}$
(since an $n$-tabloid is an equivalence class with respect to row-equivalence).

Forget our assumption that $\left(  \operatorname*{Row}\left(  i,U\right)
=\operatorname*{Row}\left(  i,V\right)  \text{ for all }i\in\mathbb{Z}\right)
$. We thus have proved that $\overline{U}=\overline{V}$ if $\left(
\operatorname*{Row}\left(  i,U\right)  =\operatorname*{Row}\left(  i,V\right)
\text{ for all }i\in\mathbb{Z}\right)  $. In other words, we have proved the
implication%
\[
\left(  \operatorname*{Row}\left(  i,U\right)  =\operatorname*{Row}\left(
i,V\right)  \text{ for all }i\in\mathbb{Z}\right)  \ \Longrightarrow\ \left(
\overline{U}=\overline{V}\right)  .
\]
The same argument (but read in reverse) shows the converse implication%
\[
\left(  \overline{U}=\overline{V}\right)  \ \Longrightarrow\ \left(
\operatorname*{Row}\left(  i,U\right)  =\operatorname*{Row}\left(  i,V\right)
\text{ for all }i\in\mathbb{Z}\right)
\]
(indeed, this is even easier, since the equality of two multisets
automatically implies the equality of the corresponding sets). Combining these
two implications, we conclude the equivalence%
\[
\left(  \overline{U}=\overline{V}\right)  \ \Longleftrightarrow\ \left(
\operatorname*{Row}\left(  i,U\right)  =\operatorname*{Row}\left(  i,V\right)
\text{ for all }i\in\mathbb{Z}\right)  .
\]
This proves Lemma \ref{lem.bas.mur.row-equ-as-Row-Row} \textbf{(a)}. \medskip

\textbf{(b)} This is entirely analogous to part \textbf{(a)}.
\end{proof}

\begin{proof}
[Detailed proof of Proposition \ref{prop.bas.mur.Sn-act.all}.]Recall that we
defined $\MurpF_{\operatorname*{all},\geq T}^{\operatorname*{Row}}$ as
\[
\MurpF_{\operatorname*{all},\geq T}^{\operatorname*{Row}}=\operatorname*{span}%
\left\{  \nabla_{V,U}^{\operatorname*{Row}}\ \mid\ \left(  \lambda,U,V\right)
\in\operatorname*{HSBT}\left(  n\right)  \text{ and }U\geq T\right\}  .
\]
In other words, $\MurpF_{\operatorname*{all},\geq T}^{\operatorname*{Row}}$ is
the span of all the vectors of the form $\nabla_{V,U}^{\operatorname*{Row}}$
with $\left(  \lambda,U,V\right)  \in\operatorname*{HSBT}\left(  n\right)  $
satisfying $U\geq T$. We shall refer to such vectors as the \textquotedblleft
good nablas\textquotedblright. Thus, $\MurpF_{\operatorname*{all},\geq
T}^{\operatorname*{Row}}$ is the span of all the good nablas.

Let us now prove that $\MurpF_{\operatorname*{all},\geq T}%
^{\operatorname*{Row}}$ is a left $\mathcal{A}$-submodule of $\mathcal{A}$.
For this, it suffices to show that $\MurpF_{\operatorname*{all},\geq
T}^{\operatorname*{Row}}$ is an $S_{n}$-subset of $\mathcal{A}$ (since
$\MurpF_{\operatorname*{all},\geq T}^{\operatorname*{Row}}$ is clearly a
$\mathbf{k}$-submodule of $\mathcal{A}$). In other words, it suffices to show
that $w\mathbf{a}\in\MurpF_{\operatorname*{all},\geq T}^{\operatorname*{Row}}$
for each $w\in S_{n}$ and each $\mathbf{a}\in\MurpF_{\operatorname*{all},\geq
T}^{\operatorname*{Row}}$.

So let us show this. Let $w\in S_{n}$ and $\mathbf{a}\in
\MurpF_{\operatorname*{all},\geq T}^{\operatorname*{Row}}$. We must prove that
$w\mathbf{a}\in\MurpF_{\operatorname*{all},\geq T}^{\operatorname*{Row}}$.
Since this claim depends $\mathbf{k}$-linearly on $\mathbf{a}$ (because
$\MurpF_{\operatorname*{all},\geq T}^{\operatorname*{Row}}$ is a $\mathbf{k}%
$-submodule of $\mathcal{A}$), we can WLOG assume (by linearity) that
$\mathbf{a}$ is one of the good nablas (since $\MurpF_{\operatorname*{all}%
,\geq T}^{\operatorname*{Row}}$ is the span of all the good nablas). In other
words, we can WLOG assume that $\mathbf{a}=\nabla_{V,U}^{\operatorname*{Row}}$
for some $\left(  \lambda,U,V\right)  \in\operatorname*{HSBT}\left(  n\right)
$ satisfying $U\geq T$. Assume this, and consider this $\left(  \lambda
,U,V\right)  \in\operatorname*{HSBT}\left(  n\right)  $. Then, $\left(
\lambda,U,V\right)  \in\operatorname*{HSBT}\left(  n\right)  $ shows that
$\left(  \lambda,U,V\right)  $ is a half-standard $n$-bitableau. In other
words, $\lambda$ is a partition of $n$ and $U$ and $V$ are two $n$-tableaux of
shape $\lambda$ such that $U$ is standard. Hence, of course, $wV$ is an
$n$-tableau of shape $\lambda$ as well. Hence, $\left(  \lambda,U,wV\right)  $
is again a half-standard $n$-bitableau, i.e., an element of
$\operatorname*{HSBT}\left(  n\right)  $. Therefore, $\nabla_{wV,U}%
^{\operatorname*{Row}}$ is one of the good nablas as well (since $U\geq T$),
and hence belongs to $\MurpF_{\operatorname*{all},\geq T}^{\operatorname*{Row}%
}$ (since $\MurpF_{\operatorname*{all},\geq T}^{\operatorname*{Row}}$ is the
span of all the good nablas). In other words, $\nabla_{wV,U}%
^{\operatorname*{Row}}\in\MurpF_{\operatorname*{all},\geq T}%
^{\operatorname*{Row}}$.

Now, (\ref{eq.prop.bas.mur.Sn-act.Row}) (applied to $P=V$ and $Q=U$ and $g=w$
and $h=\operatorname*{id}$) yields $\nabla_{wV,\operatorname*{id}%
U}^{\operatorname*{Row}}=w\nabla_{V,U}^{\operatorname*{Row}}%
\underbrace{\operatorname*{id}\nolimits^{-1}}_{=\operatorname*{id}=1}%
=w\nabla_{V,U}^{\operatorname*{Row}}$. In other words, $\nabla_{wV,U}%
^{\operatorname*{Row}}=w\nabla_{V,U}^{\operatorname*{Row}}$ (since
$\operatorname*{id}U=U$). Hence, $w\nabla_{V,U}^{\operatorname*{Row}}%
=\nabla_{wV,U}^{\operatorname*{Row}}\in\MurpF_{\operatorname*{all},\geq
T}^{\operatorname*{Row}}$. In other words, $w\mathbf{a}\in
\MurpF_{\operatorname*{all},\geq T}^{\operatorname*{Row}}$ (since
$\mathbf{a}=\nabla_{V,U}^{\operatorname*{Row}}$).

Forget that we fixed $w$ and $\mathbf{a}$. We thus have shown that
$w\mathbf{a}\in\MurpF_{\operatorname*{all},\geq T}^{\operatorname*{Row}}$ for
each $w\in S_{n}$ and each $\mathbf{a}\in\MurpF_{\operatorname*{all},\geq
T}^{\operatorname*{Row}}$. In other words, $\MurpF_{\operatorname*{all},\geq
T}^{\operatorname*{Row}}$ is an $S_{n}$-subset of $\mathcal{A}$. Hence,
$\MurpF_{\operatorname*{all},\geq T}^{\operatorname*{Row}}$ is a
subrepresentation of the $S_{n}$-representation $\mathcal{A}$ (since
$\MurpF_{\operatorname*{all},\geq T}^{\operatorname*{Row}}$ is a $\mathbf{k}%
$-submodule of $\mathcal{A}$). In other words, $\MurpF_{\operatorname*{all}%
,\geq T}^{\operatorname*{Row}}$ is a left $\mathbf{k}\left[  S_{n}\right]
$-submodule of $\mathcal{A}$ (since Proposition \ref{prop.rep.G-rep.sub=sub}
shows that the subrepresentations of the $S_{n}$-representation $\mathcal{A}$
are precisely the left $\mathbf{k}\left[  S_{n}\right]  $-submodules of
$\mathcal{A}$). In other words, $\MurpF_{\operatorname*{all},\geq
T}^{\operatorname*{Row}}$ is a left $\mathcal{A}$-submodule of $\mathcal{A}$
(since $\mathbf{k}\left[  S_{n}\right]  =\mathcal{A}$). Likewise, we can show
that $\MurpF_{\operatorname*{all},>T}^{\operatorname*{Row}}$,
$\MurpF_{\operatorname*{all},\leq T}^{-\operatorname*{Col}}$ and
$\MurpF_{\operatorname*{all},<T}^{-\operatorname*{Col}}$ are left
$\mathcal{A}$-submodules of $\mathcal{A}$. (The proofs for
$\MurpF_{\operatorname*{all},\leq T}^{-\operatorname*{Col}}$ and
$\MurpF_{\operatorname*{all},<T}^{-\operatorname*{Col}}$ require application
of (\ref{eq.prop.bas.mur.Sn-act.Col}) instead of
(\ref{eq.prop.bas.mur.Sn-act.Row}), which is slightly more complicated:
Instead of $w\nabla_{V,U}^{\operatorname*{Row}}=\nabla_{wV,U}%
^{\operatorname*{Row}}$, we now obtain $\left(  -1\right)  ^{w}w\nabla
_{V,U}^{-\operatorname*{Col}}=\nabla_{wV,U}^{-\operatorname*{Col}}$. But this
leads to the same result, since we can divide the equality $\left(  -1\right)
^{w}w\nabla_{V,U}^{-\operatorname*{Col}}=\nabla_{wV,U}^{-\operatorname*{Col}}$
by the invertible scalar $\left(  -1\right)  ^{w}$ to obtain $w\nabla
_{V,U}^{-\operatorname*{Col}}=\dfrac{1}{\left(  -1\right)  ^{w}}\nabla
_{wV,U}^{-\operatorname*{Col}}$, and the $\left(  -1\right)  ^{w}$ factor does
not affect the property of $\nabla_{wV,U}^{-\operatorname*{Col}}$ to belong to
$\MurpF_{\operatorname*{all},\leq T}^{-\operatorname*{Col}}$ or
$\MurpF_{\operatorname*{all},<T}^{-\operatorname*{Col}}$.)

Altogether, Proposition \ref{prop.bas.mur.Sn-act.all} is thus proved.
\end{proof}

\begin{proof}
[Proof of Lemma \ref{lem.bas.mur.triang.2}, parts \textbf{(c)} and
\textbf{(d)}.]We shall now prove parts \textbf{(c)} and \textbf{(d)} of Lemma
\ref{lem.bas.mur.triang.2}. (Parts \textbf{(a)} and \textbf{(b)} were already
proved in Section \ref{sec.bas.mur}.) \medskip

\textbf{(c)} We must show that $\mathbf{a}\in\MurpF_{\operatorname*{std}%
,<T}^{-\operatorname*{Col}}$ for each $\mathbf{a}\in\left(
\MurpF_{\operatorname*{std},\geq T}^{\operatorname*{Row}}\right)  ^{\perp}$.

So let $\mathbf{a}\in\left(  \MurpF_{\operatorname*{std},\geq T}%
^{\operatorname*{Row}}\right)  ^{\perp}$. We must show that $\mathbf{a}%
\in\MurpF_{\operatorname*{std},<T}^{-\operatorname*{Col}}$.

We have $\mathbf{a}\in\left(  \MurpF_{\operatorname*{std},\geq T}%
^{\operatorname*{Row}}\right)  ^{\perp}\subseteq\mathcal{A}$. Since the column
Murphy basis $\left(  \nabla_{V,U}^{-\operatorname*{Col}}\right)  _{\left(
\lambda,U,V\right)  \in\operatorname*{SBT}\left(  n\right)  }$ is a basis of
$\mathcal{A}$ (by Theorem \ref{thm.bas.mur.bas}), we can thus write
$\mathbf{a}$ as a $\mathbf{k}$-linear combination of its elements. In other
words, we can write $\mathbf{a}$ as
\begin{equation}
\mathbf{a}=\sum_{\left(  \lambda,U,V\right)  \in\operatorname*{SBT}\left(
n\right)  }\omega_{\lambda,U,V}\nabla_{V,U}^{-\operatorname*{Col}}
\label{pf.lem.bas.mur.triang.2.c=wsum}%
\end{equation}
for some scalars $\omega_{\lambda,U,V}\in\mathbf{k}$. Consider these
$\omega_{\lambda,U,V}$. We can rewrite the equality
(\ref{pf.lem.bas.mur.triang.2.c=wsum}) as
\begin{equation}
\mathbf{a}=\sum_{\left(  \mu,P,Q\right)  \in\operatorname*{SBT}\left(
n\right)  }\omega_{\mu,P,Q}\nabla_{Q,P}^{-\operatorname*{Col}}
\label{pf.lem.bas.mur.triang.2.c=wsum2}%
\end{equation}
(by renaming the summation index $\left(  \lambda,U,V\right)  $ as $\left(
\mu,P,Q\right)  $).

Now, we shall show the following:

\begin{statement}
\textit{Claim 2:} We have%
\[
\omega_{\lambda,U,V}=0\ \ \ \ \ \ \ \ \ \ \text{for all }\left(
\lambda,U,V\right)  \in\operatorname*{SBT}\left(  n\right)  \text{ satisfying
}U\geq T.
\]

\end{statement}

\begin{proof}
[Proof of Claim 2.]We proceed by strong induction on $\left(  \lambda
,U,V\right)  $ using the \textbf{reverse} of the total order on
$\operatorname*{SBT}\left(  n\right)  $ introduced in Definition
\ref{def.bas.tord.sbt.tord}. So we fix some $\left(  \lambda,U,V\right)
\in\operatorname*{SBT}\left(  n\right)  $ satisfying $U\geq T$. We assume (as
the induction hypothesis) that
\begin{equation}
\omega_{\mu,P,Q}=0 \label{pf.lem.bas.mur.triang.2.cIH}%
\end{equation}
for all $\left(  \mu,P,Q\right)  \in\operatorname*{SBT}\left(  n\right)  $
satisfying $P\geq T$ and $\left(  \mu,P,Q\right)  >\left(  \lambda,U,V\right)
$. Our goal is to prove that $\omega_{\lambda,U,V}=0$.

Since $\left(  \lambda,U,V\right)  \in\operatorname*{SBT}\left(  n\right)  $
and $U\geq T$, we have $\nabla_{V,U}^{\operatorname*{Row}}\in
\MurpF_{\operatorname*{std},\geq T}^{\operatorname*{Row}}$ by the definition
of $\MurpF_{\operatorname*{std},\geq T}^{\operatorname*{Row}}$ (in fact,
$\nabla_{V,U}^{\operatorname*{Row}}$ is one of the spanning vectors in the
definition of $\MurpF_{\operatorname*{std},\geq T}^{\operatorname*{Row}}$).

Since $\mathbf{a}\in\left(  \MurpF_{\operatorname*{std},\geq T}%
^{\operatorname*{Row}}\right)  ^{\perp}$, we have $\left\langle \mathbf{a}%
,\mathbf{v}\right\rangle =0$ for each $\mathbf{v}\in
\MurpF_{\operatorname*{std},\geq T}^{\operatorname*{Row}}$ (by the definition
of $\left(  \MurpF_{\operatorname*{std},\geq T}^{\operatorname*{Row}}\right)
^{\perp}$). Applying this to $\mathbf{v}=\nabla_{V,U}^{\operatorname*{Row}}$,
we obtain $\left\langle \mathbf{a},\ \nabla_{V,U}^{\operatorname*{Row}%
}\right\rangle =0$ (since $\nabla_{V,U}^{\operatorname*{Row}}\in
\MurpF_{\operatorname*{std},\geq T}^{\operatorname*{Row}}$). Thus,%
\begin{align}
0  &  =\left\langle \mathbf{a},\ \nabla_{V,U}^{\operatorname*{Row}%
}\right\rangle =\left\langle \sum_{\left(  \mu,P,Q\right)  \in
\operatorname*{SBT}\left(  n\right)  }\omega_{\mu,P,Q}\nabla_{Q,P}%
^{-\operatorname*{Col}},\ \nabla_{V,U}^{\operatorname*{Row}}\right\rangle
\ \ \ \ \ \ \ \ \ \ \left(  \text{by (\ref{pf.lem.bas.mur.triang.2.c=wsum2}%
)}\right) \nonumber\\
&  =\sum_{\left(  \mu,P,Q\right)  \in\operatorname*{SBT}\left(  n\right)
}\omega_{\mu,P,Q}\left\langle \nabla_{Q,P}^{-\operatorname*{Col}}%
,\ \nabla_{V,U}^{\operatorname*{Row}}\right\rangle
\label{pf.lem.bas.mur.triang.2.c0}%
\end{align}
(since the form $\left\langle \cdot,\cdot\right\rangle $ is bilinear). Now, we
observe the following:

\begin{itemize}
\item For each $\left(  \mu,P,Q\right)  \in\operatorname*{SBT}\left(
n\right)  $ satisfying $\left(  \mu,P,Q\right)  >\left(  \lambda,U,V\right)
$, we have%
\begin{equation}
\omega_{\mu,P,Q}=0. \label{pf.lem.bas.mur.triang.2.c1}%
\end{equation}

[\textit{Proof:} Let $\left(  \mu,P,Q\right)  \in\operatorname*{SBT}\left(
n\right)  $ be such that $\left(  \mu,P,Q\right)  >\left(  \lambda,U,V\right)
$. Thus, $\left(  \lambda,U,V\right)  <\left(  \mu,P,Q\right)  $. Hence, Lemma
\ref{lem.bas.tord.sbt.1} \textbf{(a)} (applied to $\left(  \lambda,U,V\right)
$ and $\left(  \mu,P,Q\right)  $ instead of $\left(  \lambda,A,B\right)  $ and
$\left(  \mu,C,D\right)  $) shows that $U\leq P$. Thus, $P\geq U\geq T$, so
that $P\geq T$ (since the relation $<$ is a total order). Thus,
(\ref{pf.lem.bas.mur.triang.2.cIH}) yields $\omega_{\mu,P,Q}=0$. This proves
(\ref{pf.lem.bas.mur.triang.2.c1}).]

\item For each $\left(  \mu,P,Q\right)  \in\operatorname*{SBT}\left(
n\right)  $ satisfying $\left(  \mu,P,Q\right)  <\left(  \lambda,U,V\right)
$, we have%
\begin{equation}
\left\langle \nabla_{Q,P}^{-\operatorname*{Col}},\ \nabla_{V,U}%
^{\operatorname*{Row}}\right\rangle =0. \label{pf.lem.bas.mur.triang.2.c2}%
\end{equation}

[\textit{Proof:} Let $\left(  \mu,P,Q\right)  \in\operatorname*{SBT}\left(
n\right)  $ be such that $\left(  \mu,P,Q\right)  <\left(  \lambda,U,V\right)
$. Thus, Lemma \ref{lem.bas.tord.sbt.1} \textbf{(b)} (applied to $\left(
\mu,P,Q\right)  $ and $\left(  \lambda,U,V\right)  $ instead of $\left(
\lambda,A,B\right)  $ and $\left(  \mu,C,D\right)  $) yields $P<U$ or $Q<V$.
Thus, Lemma \ref{lem.bas.mur.col} (applied to $A=P$, $B=Q$, $C=U$ and $D=V$)
yields $\left\langle \nabla_{Q,P}^{-\operatorname*{Col}},\ \nabla
_{V,U}^{\operatorname*{Row}}\right\rangle =0$. This proves
(\ref{pf.lem.bas.mur.triang.2.c2}).]
\end{itemize}

Now, recall that the relation $<$ on $\operatorname*{SBT}\left(  n\right)  $
is the smaller relation of a total order. Hence, each $\left(  \mu,P,Q\right)
\in\operatorname*{SBT}\left(  n\right)  $ satisfies exactly one of the three
statements $\left(  \mu,P,Q\right)  <\left(  \lambda,U,V\right)  $ and
$\left(  \mu,P,Q\right)  >\left(  \lambda,U,V\right)  $ and $\left(
\mu,P,Q\right)  =\left(  \lambda,U,V\right)  $.

Accordingly, the sum on the right hand side of
(\ref{pf.lem.bas.mur.triang.2.c0}) can be split into three parts: the part
comprising all $\left(  \mu,P,Q\right)  \in\operatorname*{SBT}\left(
n\right)  $ that satisfy $\left(  \mu,P,Q\right)  <\left(  \lambda,U,V\right)
$; the part comprising all $\left(  \mu,P,Q\right)  \in\operatorname*{SBT}%
\left(  n\right)  $ that satisfy $\left(  \mu,P,Q\right)  >\left(
\lambda,U,V\right)  $; and a single addend corresponding to $\left(
\mu,P,Q\right)  =\left(  \lambda,U,V\right)  $. Splitting the sum in this way,
we rewrite (\ref{pf.lem.bas.mur.triang.2.c0}) as follows:
\begin{align*}
0  &  =\sum_{\substack{\left(  \mu,P,Q\right)  \in\operatorname*{SBT}\left(
n\right)  ;\\\left(  \mu,P,Q\right)  <\left(  \lambda,U,V\right)  }%
}\omega_{\mu,P,Q}\underbrace{\left\langle \nabla_{Q,P}^{-\operatorname*{Col}%
},\ \nabla_{V,U}^{\operatorname*{Row}}\right\rangle }_{\substack{=0\\\text{(by
(\ref{pf.lem.bas.mur.triang.2.c2}))}}}\\
&  \ \ \ \ \ \ \ \ \ \ +\sum_{\substack{\left(  \mu,P,Q\right)  \in
\operatorname*{SBT}\left(  n\right)  ;\\\left(  \mu,P,Q\right)  >\left(
\lambda,U,V\right)  }}\ \ \underbrace{\omega_{\mu,P,Q}}%
_{\substack{=0\\\text{(by (\ref{pf.lem.bas.mur.triang.2.c1}))}}}\left\langle
\nabla_{Q,P}^{-\operatorname*{Col}},\ \nabla_{V,U}^{\operatorname*{Row}%
}\right\rangle +\omega_{\lambda,U,V}\left\langle \nabla_{V,U}%
^{-\operatorname*{Col}},\ \nabla_{V,U}^{\operatorname*{Row}}\right\rangle \\
&  =\underbrace{\sum_{\substack{\left(  \mu,P,Q\right)  \in\operatorname*{SBT}%
\left(  n\right)  ;\\\left(  \mu,P,Q\right)  <\left(  \lambda,U,V\right)
}}\omega_{\mu,P,Q}0}_{=0}+\underbrace{\sum_{\substack{\left(  \mu,P,Q\right)
\in\operatorname*{SBT}\left(  n\right)  ;\\\left(  \mu,P,Q\right)  >\left(
\lambda,U,V\right)  }}0\left\langle \nabla_{Q,P}^{-\operatorname*{Col}%
},\ \nabla_{V,U}^{\operatorname*{Row}}\right\rangle }_{=0}+\,\omega
_{\lambda,U,V}\left\langle \nabla_{V,U}^{-\operatorname*{Col}},\ \nabla
_{V,U}^{\operatorname*{Row}}\right\rangle \\
&  =\omega_{\lambda,U,V}\underbrace{\left\langle \nabla_{V,U}%
^{-\operatorname*{Col}},\ \nabla_{V,U}^{\operatorname*{Row}}\right\rangle
}_{\substack{=\pm1\\\text{(by Lemma \ref{lem.bas.mur.dia},}\\\text{applied to
}A=U\text{ and }B=V\text{)}}}=\pm\omega_{\lambda,U,V}.
\end{align*}
Hence, $\pm\omega_{\lambda,U,V}=0$. In other words, $\omega_{\lambda,U,V}=0$.
This completes the induction, and thus the proof of Claim 2.
\end{proof}

Now, (\ref{pf.lem.bas.mur.triang.2.c=wsum}) becomes%
\begin{align*}
\mathbf{a}  &  =\sum_{\left(  \lambda,U,V\right)  \in\operatorname*{SBT}%
\left(  n\right)  }\omega_{\lambda,U,V}\nabla_{V,U}^{-\operatorname*{Col}}\\
&  =\sum_{\substack{\left(  \lambda,U,V\right)  \in\operatorname*{SBT}\left(
n\right)  ;\\U\geq T}}\underbrace{\omega_{\lambda,U,V}}%
_{\substack{=0\\\text{(by Claim 2)}}}\nabla_{V,U}^{-\operatorname*{Col}}%
+\sum_{\substack{\left(  \lambda,U,V\right)  \in\operatorname*{SBT}\left(
n\right)  ;\\U<T}}\omega_{\lambda,U,V}\nabla_{V,U}^{-\operatorname*{Col}}\\
&  \ \ \ \ \ \ \ \ \ \ \ \ \ \ \ \ \ \ \ \ \left(
\begin{array}
[c]{c}%
\text{since each }\left(  \lambda,U,V\right)  \in\operatorname*{SBT}\left(
n\right)  \text{ satisfies}\\
\text{either }U\geq T\text{ or }U<T\text{ (but not both)}%
\end{array}
\right) \\
&  =\underbrace{\sum_{\substack{\left(  \lambda,U,V\right)  \in
\operatorname*{SBT}\left(  n\right)  ;\\U\geq T}}0\nabla_{V,U}%
^{-\operatorname*{Col}}}_{=0}+\sum_{\substack{\left(  \lambda,U,V\right)
\in\operatorname*{SBT}\left(  n\right)  ;\\U<T}}\omega_{\lambda,U,V}%
\nabla_{V,U}^{-\operatorname*{Col}}\\
&  =\sum_{\substack{\left(  \lambda,U,V\right)  \in\operatorname*{SBT}\left(
n\right)  ;\\U<T}}\omega_{\lambda,U,V}\nabla_{V,U}^{-\operatorname*{Col}}\\
&  \in\operatorname*{span}\left\{  \nabla_{V,U}^{-\operatorname*{Col}}%
\ \mid\ \left(  \lambda,U,V\right)  \in\operatorname*{SBT}\left(  n\right)
\text{ and }U<T\right\} \\
&  =\MurpF_{\operatorname*{std},<T}^{-\operatorname*{Col}}%
\ \ \ \ \ \ \ \ \ \ \left(  \text{by the definition of }%
\MurpF_{\operatorname*{std},<T}^{-\operatorname*{Col}}\right)  .
\end{align*}
This completes our proof of Lemma \ref{lem.bas.mur.triang.2} \textbf{(c)}.
\medskip

\textbf{(d)} This is analogous to the above proof of Lemma
\ref{lem.bas.mur.triang.2} \textbf{(c)}; the only difference is that some
inequality signs (between standard $n$-tableaux) change from strong to weak
and vice versa. (For instance, the \textquotedblleft$U\geq T$%
\textquotedblright\ in Claim 2 must be replaced by \textquotedblleft%
$U>T$\textquotedblright). We leave it to the reader to make the necessary changes.
\end{proof}

\begin{proof}
[Detailed proof of Theorem \ref{thm.bas.mur.filt}.]As in the proof of
Proposition \ref{prop.bas.mur.submods.eq1}, we can see that $m>0$. Hence,
$0<m$. \medskip

\textbf{(a)} Let $i\in\left\{  0,1,\ldots,m\right\}  $. We must prove that%
\[
\MurpF_{i}^{\operatorname*{Row}}=%
\begin{cases}
\MurpF_{\operatorname*{std},>T_{i}}^{\operatorname*{Row}}, & \text{if }i>0;\\
\mathcal{A}, & \text{if }i=0.
\end{cases}
\]

We have $i\in\left\{  0,1,\ldots,m\right\}  $. Thus, we are in one of the
following three cases:

\textit{Case 1:} We have $i=0$.

\textit{Case 2:} We have $0<i<m$.

\textit{Case 3:} We have $i=m$.

First, let us consider Case 1. In this case, we have $i=0$. Hence, $i+1=1$ and
$i=0<m$. The definition of $\MurpF_{i}^{\operatorname*{Row}}$ yields
\begin{align*}
\MurpF_{i}^{\operatorname*{Row}}  &  =%
\begin{cases}
\MurpF_{\operatorname*{std},\geq T_{i+1}}^{\operatorname*{Row}}, & \text{if
}i<m;\\
0, & \text{if }i=m
\end{cases}
\ \ =\MurpF_{\operatorname*{std},\geq T_{i+1}}^{\operatorname*{Row}%
}\ \ \ \ \ \ \ \ \ \ \left(  \text{since }i<m\right) \\
&  =\MurpF_{\operatorname*{std},\geq T_{1}}^{\operatorname*{Row}%
}\ \ \ \ \ \ \ \ \ \ \left(  \text{since }i+1=1\right) \\
&  =\mathcal{A}\ \ \ \ \ \ \ \ \ \ \left(  \text{by Proposition
\ref{prop.bas.mur.submods.eq1} \textbf{(a)}}\right) \\
&  =%
\begin{cases}
\MurpF_{\operatorname*{std},>T_{i}}^{\operatorname*{Row}}, & \text{if }i>0;\\
\mathcal{A}, & \text{if }i=0
\end{cases}
\ \ \ \ \ \ \ \ \ \ \left(  \text{since }i=0\right)  .
\end{align*}
Thus, Theorem \ref{thm.bas.mur.filt} \textbf{(a)} is proved in Case 1.

Let us next consider Case 2. In this case, we have $0<i<m$. Thus, $i\in\left[
m-1\right]  $ and $i>0$. The definition of $\MurpF_{i}^{\operatorname*{Row}}$
yields
\begin{align*}
\MurpF_{i}^{\operatorname*{Row}}  &  =%
\begin{cases}
\MurpF_{\operatorname*{std},\geq T_{i+1}}^{\operatorname*{Row}}, & \text{if
}i<m;\\
0, & \text{if }i=m
\end{cases}
\ \ =\MurpF_{\operatorname*{std},\geq T_{i+1}}^{\operatorname*{Row}%
}\ \ \ \ \ \ \ \ \ \ \left(  \text{since }i<m\right) \\
&  =\MurpF_{\operatorname*{std},>T_{i}}^{\operatorname*{Row}}%
\ \ \ \ \ \ \ \ \ \ \left(  \text{by Proposition
\ref{prop.bas.mur.submods.eq1} \textbf{(c)}}\right) \\
&  =%
\begin{cases}
\MurpF_{\operatorname*{std},>T_{i}}^{\operatorname*{Row}}, & \text{if }i>0;\\
\mathcal{A}, & \text{if }i=0
\end{cases}
\ \ \ \ \ \ \ \ \ \ \left(  \text{since }i>0\right)  .
\end{align*}
Thus, Theorem \ref{thm.bas.mur.filt} \textbf{(a)} is proved in Case 2.

Finally, let us consider Case 3. In this case, we have $i=m$. Hence, $i=m>0$.
The definition of $\MurpF_{i}^{\operatorname*{Row}}$ yields
\[
\MurpF_{i}^{\operatorname*{Row}}=%
\begin{cases}
\MurpF_{\operatorname*{std},\geq T_{i+1}}^{\operatorname*{Row}}, & \text{if
}i<m;\\
0, & \text{if }i=m
\end{cases}
\ \ =0\ \ \ \ \ \ \ \ \ \ \left(  \text{since }i=m\right)  .
\]
Comparing this with%
\begin{align*}%
\begin{cases}
\MurpF_{\operatorname*{std},>T_{i}}^{\operatorname*{Row}}, & \text{if }i>0;\\
\mathcal{A}, & \text{if }i=0
\end{cases}
\ \  &  =\MurpF_{\operatorname*{std},>T_{i}}^{\operatorname*{Row}%
}\ \ \ \ \ \ \ \ \ \ \left(  \text{since }i>0\right) \\
&  =\MurpF_{\operatorname*{std},>T_{m}}^{\operatorname*{Row}}%
\ \ \ \ \ \ \ \ \ \ \left(  \text{since }i=m\right) \\
&  =0\ \ \ \ \ \ \ \ \ \ \left(  \text{by Proposition
\ref{prop.bas.mur.submods.eq1} \textbf{(b)}}\right)  ,
\end{align*}
we obtain%
\[
\MurpF_{i}^{\operatorname*{Row}}=%
\begin{cases}
\MurpF_{\operatorname*{std},>T_{i}}^{\operatorname*{Row}}, & \text{if }i>0;\\
\mathcal{A}, & \text{if }i=0.
\end{cases}
\]
Thus, Theorem \ref{thm.bas.mur.filt} \textbf{(a)} is proved in Case 3.

We have now proved Theorem \ref{thm.bas.mur.filt} \textbf{(a)} in each of the
three Cases 1, 2 and 3. Hence, Theorem \ref{thm.bas.mur.filt} \textbf{(a)}
always holds. \medskip

\textbf{(b)} Let $i\in\left\{  0,1,\ldots,m\right\}  $. We must prove that%
\[
\MurpF_{i}^{-\operatorname*{Col}}=%
\begin{cases}
\MurpF_{\operatorname*{std},<T_{i+1}}^{-\operatorname*{Col}}, & \text{if
}i<m;\\
\mathcal{A}, & \text{if }i=m.
\end{cases}
\]

We have $i\in\left\{  0,1,\ldots,m\right\}  $. Thus, we are in one of the
following three cases:

\textit{Case 1:} We have $i=0$.

\textit{Case 2:} We have $0<i<m$.

\textit{Case 3:} We have $i=m$.

First, let us consider Case 1. In this case, we have $i=0$. Hence, $i+1=1$ and
$i=0<m$. The definition of $\MurpF_{i}^{-\operatorname*{Col}}$ yields
\[
\MurpF_{i}^{-\operatorname*{Col}}=%
\begin{cases}
\MurpF_{\operatorname*{std},\leq T_{i}}^{-\operatorname*{Col}}, & \text{if
}i>0;\\
0, & \text{if }i=0
\end{cases}
\ \ =0\ \ \ \ \ \ \ \ \ \ \left(  \text{since }i=0\right)  .
\]
Comparing this with%
\begin{align*}%
\begin{cases}
\MurpF_{\operatorname*{std},<T_{i+1}}^{-\operatorname*{Col}}, & \text{if
}i<m;\\
\mathcal{A}, & \text{if }i=m
\end{cases}
\ \  &  =\MurpF_{\operatorname*{std},<T_{i+1}}^{-\operatorname*{Col}%
}\ \ \ \ \ \ \ \ \ \ \left(  \text{since }i<m\right) \\
&  =\MurpF_{\operatorname*{std},<T_{1}}^{-\operatorname*{Col}}%
\ \ \ \ \ \ \ \ \ \ \left(  \text{since }i+1=1\right) \\
&  =0\ \ \ \ \ \ \ \ \ \ \left(  \text{by Proposition
\ref{prop.bas.mur.submods.eq1} \textbf{(b)}}\right)  ,
\end{align*}
we obtain%
\[
\MurpF_{i}^{-\operatorname*{Col}}=%
\begin{cases}
\MurpF_{\operatorname*{std},<T_{i+1}}^{-\operatorname*{Col}}, & \text{if
}i<m;\\
\mathcal{A}, & \text{if }i=m.
\end{cases}
\]
Thus, Theorem \ref{thm.bas.mur.filt} \textbf{(b)} is proved in Case 1.

Let us next consider Case 2. In this case, we have $0<i<m$. Thus, $i\in\left[
m-1\right]  $ and $i>0$. The definition of $\MurpF_{i}^{-\operatorname*{Col}}$
yields
\begin{align*}
\MurpF_{i}^{-\operatorname*{Col}}  &  =%
\begin{cases}
\MurpF_{\operatorname*{std},\leq T_{i}}^{-\operatorname*{Col}}, & \text{if
}i>0;\\
0, & \text{if }i=0
\end{cases}
\ \ =\MurpF_{\operatorname*{std},\leq T_{i}}^{-\operatorname*{Col}%
}\ \ \ \ \ \ \ \ \ \ \left(  \text{since }i>0\right) \\
&  =\MurpF_{\operatorname*{std},<T_{i+1}}^{-\operatorname*{Col}}%
\ \ \ \ \ \ \ \ \ \ \left(  \text{by Proposition
\ref{prop.bas.mur.submods.eq1} \textbf{(c)}}\right) \\
&  =%
\begin{cases}
\MurpF_{\operatorname*{std},<T_{i+1}}^{-\operatorname*{Col}}, & \text{if
}i<m;\\
\mathcal{A}, & \text{if }i=m
\end{cases}
\ \ \ \ \ \ \ \ \ \ \left(  \text{since }i<m\right)  .
\end{align*}
Thus, Theorem \ref{thm.bas.mur.filt} \textbf{(b)} is proved in Case 2.

Finally, let us consider Case 3. In this case, we have $i=m$. Hence, $i=m>0$.
The definition of $\MurpF_{i}^{-\operatorname*{Col}}$ yields
\begin{align*}
\MurpF_{i}^{-\operatorname*{Col}}  &  =%
\begin{cases}
\MurpF_{\operatorname*{std},\leq T_{i}}^{-\operatorname*{Col}}, & \text{if
}i>0;\\
0, & \text{if }i=0
\end{cases}
\ \ =\MurpF_{\operatorname*{std},\leq T_{i}}^{-\operatorname*{Col}%
}\ \ \ \ \ \ \ \ \ \ \left(  \text{since }i>0\right) \\
&  =\MurpF_{\operatorname*{std},\leq T_{m}}^{-\operatorname*{Col}%
}\ \ \ \ \ \ \ \ \ \ \left(  \text{since }i=m\right) \\
&  =\mathcal{A}\ \ \ \ \ \ \ \ \ \ \left(  \text{by Proposition
\ref{prop.bas.mur.submods.eq1} \textbf{(a)}}\right) \\
&  =%
\begin{cases}
\MurpF_{\operatorname*{std},<T_{i+1}}^{-\operatorname*{Col}}, & \text{if
}i<m;\\
\mathcal{A}, & \text{if }i=m
\end{cases}
\ \ \ \ \ \ \ \ \ \ \left(  \text{since }i=m\right)  .
\end{align*}
Thus, Theorem \ref{thm.bas.mur.filt} \textbf{(b)} is proved in Case 3.

We have now proved Theorem \ref{thm.bas.mur.filt} \textbf{(b)} in each of the
three Cases 1, 2 and 3. Hence, Theorem \ref{thm.bas.mur.filt} \textbf{(b)}
always holds. \medskip

\textbf{(c)} Let $i\in\left\{  0,1,\ldots,m\right\}  $. We must prove that the
set $\MurpF_{i}^{\operatorname*{Row}}$ is a left $\mathcal{A}$-submodule of
$\mathcal{A}$.

If $i=m$, then this is obvious (because in this case, $\MurpF_{i}%
^{\operatorname*{Row}}$ is defined to be $0$, which is clearly a left
$\mathcal{A}$-submodule of $\mathcal{A}$). Thus, we WLOG assume that $i\neq
m$. Hence, $i<m$. Thus, the definition of $\MurpF_{i}^{\operatorname*{Row}}$
shows that $\MurpF_{i}^{\operatorname*{Row}}=\MurpF_{\operatorname*{std},\geq
T_{i+1}}^{\operatorname*{Row}}$. But Theorem \ref{thm.bas.mur.triang.orth}
\textbf{(a)} (applied to $T=T_{i+1}$) shows that the set
$\MurpF_{\operatorname*{std},\geq T_{i+1}}^{\operatorname*{Row}}$ is a left
$\mathcal{A}$-submodule of $\mathcal{A}$. In other words, $\MurpF_{i}%
^{\operatorname*{Row}}$ is a left $\mathcal{A}$-submodule of $\mathcal{A}$
(since $\MurpF_{i}^{\operatorname*{Row}}=\MurpF_{\operatorname*{std},\geq
T_{i+1}}^{\operatorname*{Row}}$). Thus, Theorem \ref{thm.bas.mur.filt}
\textbf{(c)} is proved. \medskip

\textbf{(d)} Let $i\in\left\{  0,1,\ldots,m\right\}  $. We must prove that the
set $\MurpF_{i}^{-\operatorname*{Col}}$ is a left $\mathcal{A}$-submodule of
$\mathcal{A}$.

If $i=0$, then this is obvious (because in this case, $\MurpF_{i}%
^{-\operatorname*{Col}}$ is defined to be $0$, which is clearly a left
$\mathcal{A}$-submodule of $\mathcal{A}$). Thus, we WLOG assume that $i\neq0$.
Hence, $i>0$. Thus, the definition of $\MurpF_{i}^{-\operatorname*{Col}}$
shows that $\MurpF_{i}^{-\operatorname*{Col}}=\MurpF_{\operatorname*{std},\leq
T_{i}}^{-\operatorname*{Col}}$. But Theorem \ref{thm.bas.mur.triang.orth}
\textbf{(d)} (applied to $T=T_{i}$) shows that the set
$\MurpF_{\operatorname*{std},\leq T_{i}}^{-\operatorname*{Col}}$ is a left
$\mathcal{A}$-submodule of $\mathcal{A}$. In other words, $\MurpF_{i}%
^{-\operatorname*{Col}}$ is a left $\mathcal{A}$-submodule of $\mathcal{A}$
(since $\MurpF_{i}^{-\operatorname*{Col}}=\MurpF_{\operatorname*{std},\leq
T_{i}}^{-\operatorname*{Col}}$). Thus, Theorem \ref{thm.bas.mur.filt}
\textbf{(d)} is proved. \medskip

\textbf{(e)} We already know that the sets $\MurpF_{i}^{\operatorname*{Row}}$
are left $\mathcal{A}$-submodules of $\mathcal{A}$ (by Theorem
\ref{thm.bas.mur.filt} \textbf{(c)}). Thus, it remains to prove that
\[
0=\MurpF_{m}^{\operatorname*{Row}}\subseteq\MurpF_{m-1}^{\operatorname*{Row}%
}\subseteq\cdots\subseteq\MurpF_{2}^{\operatorname*{Row}}\subseteq
\MurpF_{1}^{\operatorname*{Row}}\subseteq\MurpF_{0}^{\operatorname*{Row}%
}=\mathcal{A}.
\]
In other words, it remains to prove that $\MurpF_{m}^{\operatorname*{Row}}=0$
and $\MurpF_{0}^{\operatorname*{Row}}=\mathcal{A}$ and that each $i\in\left[
m\right]  $ satisfies $\MurpF_{i}^{\operatorname*{Row}}\subseteq
\MurpF_{i-1}^{\operatorname*{Row}}$.

But $\MurpF_{m}^{\operatorname*{Row}}=0$ is obvious (an immediate consequence
of the definition of $\MurpF_{m}^{\operatorname*{Row}}$), and $\MurpF_{0}%
^{\operatorname*{Row}}=\mathcal{A}$ is also easy to check (since Theorem
\ref{thm.bas.mur.filt} \textbf{(a)} (applied to $i=0$) shows that
$\MurpF_{0}^{\operatorname*{Row}}=%
\begin{cases}
\MurpF_{\operatorname*{std},>T_{0}}^{\operatorname*{Row}}, & \text{if }0>0;\\
\mathcal{A}, & \text{if }0=0
\end{cases}
\ \ =\mathcal{A}$ (since $0=0$)). Hence, it remains to prove that each
$i\in\left[  m\right]  $ satisfies $\MurpF_{i}^{\operatorname*{Row}}%
\subseteq\MurpF_{i-1}^{\operatorname*{Row}}$.

So let $i\in\left[  m\right]  $. We must show that $\MurpF_{i}%
^{\operatorname*{Row}}\subseteq\MurpF_{i-1}^{\operatorname*{Row}}$.

We have $i\in\left[  m\right]  $, so that $i>0$ and $i\leq m$ and therefore
$i-1<i\leq m$. The definition of $\MurpF_{i-1}^{\operatorname*{Row}}$ yields%
\begin{align}
\MurpF_{i-1}^{\operatorname*{Row}}  &  =%
\begin{cases}
\MurpF_{\operatorname*{std},\geq T_{\left(  i-1\right)  +1}}%
^{\operatorname*{Row}}, & \text{if }i-1<m;\\
0, & \text{if }i-1=m
\end{cases}
\ \ =\MurpF_{\operatorname*{std},\geq T_{\left(  i-1\right)  +1}%
}^{\operatorname*{Row}}\ \ \ \ \ \ \ \ \ \ \left(  \text{since }i-1<m\right)
\nonumber\\
&  =\MurpF_{\operatorname*{std},\geq T_{i}}^{\operatorname*{Row}%
}\ \ \ \ \ \ \ \ \ \ \left(  \text{since }\left(  i-1\right)  +1=i\right)
.\nonumber
\end{align}
On the other hand, Theorem \ref{thm.bas.mur.filt} \textbf{(a)} yields%
\[
\MurpF_{i}^{\operatorname*{Row}}=%
\begin{cases}
\MurpF_{\operatorname*{std},>T_{i}}^{\operatorname*{Row}}, & \text{if }i>0;\\
\mathcal{A}, & \text{if }i=0
\end{cases}
\ \ =\MurpF_{\operatorname*{std},>T_{i}}^{\operatorname*{Row}}%
\ \ \ \ \ \ \ \ \ \ \left(  \text{since }i>0\right)  .
\]

But Proposition \ref{prop.bas.mur.subset} (applied to $T=T_{i}$) yields
\[
\MurpF_{\operatorname*{std},>T_{i}}^{\operatorname*{Row}}\subseteq
\MurpF_{\operatorname*{std},\geq T_{i}}^{\operatorname*{Row}}%
\ \ \ \ \ \ \ \ \ \ \text{and}\ \ \ \ \ \ \ \ \ \ \MurpF_{\operatorname*{std}%
,<T_{i}}^{-\operatorname*{Col}}\subseteq\MurpF_{\operatorname*{std},\leq
T_{i}}^{-\operatorname*{Col}}.
\]
Hence, $\MurpF_{i}^{\operatorname*{Row}}=\MurpF_{\operatorname*{std},>T_{i}%
}^{\operatorname*{Row}}\subseteq\MurpF_{\operatorname*{std},\geq T_{i}%
}^{\operatorname*{Row}}=\MurpF_{i-1}^{\operatorname*{Row}}$ (since
$\MurpF_{i-1}^{\operatorname*{Row}}=\MurpF_{\operatorname*{std},\geq T_{i}%
}^{\operatorname*{Row}}$).

Thus we have shown that each $i\in\left[  m\right]  $ satisfies $\MurpF_{i}%
^{\operatorname*{Row}}\subseteq\MurpF_{i-1}^{\operatorname*{Row}}$. As we
said, this completes the proof of Theorem \ref{thm.bas.mur.filt} \textbf{(e)}.
\medskip

\textbf{(f)} We already know that the sets $\MurpF_{i}^{-\operatorname*{Col}}$
are left $\mathcal{A}$-submodules of $\mathcal{A}$ (by Theorem
\ref{thm.bas.mur.filt} \textbf{(d)}). Thus, it remains to prove that
\[
0=\MurpF_{0}^{-\operatorname*{Col}}\subseteq\MurpF_{1}^{-\operatorname*{Col}%
}\subseteq\MurpF_{2}^{-\operatorname*{Col}}\subseteq\cdots\subseteq
\MurpF_{m-1}^{-\operatorname*{Col}}\subseteq\MurpF_{m}^{-\operatorname*{Col}%
}=\mathcal{A}.
\]
In other words, it remains to prove that $\MurpF_{0}^{-\operatorname*{Col}}=0$
and $\MurpF_{m}^{-\operatorname*{Col}}=\mathcal{A}$ and that each $i\in\left[
m\right]  $ satisfies $\MurpF_{i-1}^{-\operatorname*{Col}}\subseteq
\MurpF_{i}^{-\operatorname*{Col}}$.

But $\MurpF_{0}^{-\operatorname*{Col}}=0$ is obvious (an immediate consequence
of the definition of $\MurpF_{0}^{-\operatorname*{Col}}$), and $\MurpF_{m}%
^{-\operatorname*{Col}}=\mathcal{A}$ is also easy to check (since Theorem
\ref{thm.bas.mur.filt} \textbf{(b)} (applied to $i=m$) shows that
$\MurpF_{m}^{-\operatorname*{Col}}=%
\begin{cases}
\MurpF_{\operatorname*{std},<T_{m+1}}^{-\operatorname*{Col}}, & \text{if
}m<m;\\
\mathcal{A}, & \text{if }m=m
\end{cases}
\ \ =\mathcal{A}$ (since $m=m$)). Hence, it remains to prove that each
$i\in\left[  m\right]  $ satisfies $\MurpF_{i-1}^{-\operatorname*{Col}%
}\subseteq\MurpF_{i}^{-\operatorname*{Col}}$.

So let $i\in\left[  m\right]  $. We must show that $\MurpF_{i-1}%
^{-\operatorname*{Col}}\subseteq\MurpF_{i}^{-\operatorname*{Col}}$.

We have $i\in\left[  m\right]  $, so that $i>0$ and $i\leq m$ and therefore
$i-1<i\leq m$. Theorem \ref{thm.bas.mur.filt} \textbf{(b)} (applied to $i-1$
instead of $i$) yields%
\begin{align}
\MurpF_{i-1}^{-\operatorname*{Col}}  &  =%
\begin{cases}
\MurpF_{\operatorname*{std},<T_{\left(  i-1\right)  +1}}^{-\operatorname*{Col}%
}, & \text{if }i-1<m;\\
\mathcal{A}, & \text{if }i-1=m
\end{cases}
\ \ =\MurpF_{\operatorname*{std},<T_{\left(  i-1\right)  +1}}%
^{-\operatorname*{Col}}\ \ \ \ \ \ \ \ \ \ \left(  \text{since }i-1<m\right)
\nonumber\\
&  =\MurpF_{\operatorname*{std},<T_{i}}^{-\operatorname*{Col}}%
\ \ \ \ \ \ \ \ \ \ \left(  \text{since }\left(  i-1\right)  +1=i\right)
.\nonumber
\end{align}
On the other hand, the definition of $\MurpF_{i}^{-\operatorname*{Col}}$
yields%
\[
\MurpF_{i}^{-\operatorname*{Col}}=%
\begin{cases}
\MurpF_{\operatorname*{std},\leq T_{i}}^{-\operatorname*{Col}}, & \text{if
}i>0;\\
0, & \text{if }i=0
\end{cases}
\ \ =\MurpF_{\operatorname*{std},\leq T_{i}}^{-\operatorname*{Col}%
}\ \ \ \ \ \ \ \ \ \ \left(  \text{since }i>0\right)  .
\]

But Proposition \ref{prop.bas.mur.subset} (applied to $T=T_{i}$) yields
\[
\MurpF_{\operatorname*{std},>T_{i}}^{\operatorname*{Row}}\subseteq
\MurpF_{\operatorname*{std},\geq T_{i}}^{\operatorname*{Row}}%
\ \ \ \ \ \ \ \ \ \ \text{and}\ \ \ \ \ \ \ \ \ \ \MurpF_{\operatorname*{std}%
,<T_{i}}^{-\operatorname*{Col}}\subseteq\MurpF_{\operatorname*{std},\leq
T_{i}}^{-\operatorname*{Col}}.
\]
Hence, $\MurpF_{i-1}^{-\operatorname*{Col}}=\MurpF_{\operatorname*{std}%
,<T_{i}}^{-\operatorname*{Col}}\subseteq\MurpF_{\operatorname*{std},\leq
T_{i}}^{-\operatorname*{Col}}=\MurpF_{i}^{-\operatorname*{Col}}$ (since
$\MurpF_{i}^{-\operatorname*{Col}}=\MurpF_{\operatorname*{std},\leq T_{i}%
}^{-\operatorname*{Col}}$).

Thus we have shown that each $i\in\left[  m\right]  $ satisfies $\MurpF_{i-1}%
^{-\operatorname*{Col}}\subseteq\MurpF_{i}^{-\operatorname*{Col}}$. As we
said, this completes the proof of Theorem \ref{thm.bas.mur.filt} \textbf{(f)}.
\end{proof}
\end{fineprint}

\begin{fineprint}
\begin{proof}
[Detailed proof of Lemma \ref{lem.bas.mur.triang.iso1lem} \textbf{(e)}.]We
observe that $T\in\operatorname*{SYT}\left(  \lambda\right)  $ (since $T$ is a
standard $n$-tableau of shape $\lambda$).

As we saw in the proof of Lemma \ref{lem.bas.mur.triang.iso1lem} \textbf{(c)},
we have%
\[
\MurpF_{\operatorname*{std},<T}^{-\operatorname*{Col}}=\operatorname*{span}%
\left\{  \nabla_{V,U}^{-\operatorname*{Col}}\ \mid\ \left(  \mu,U,V\right)
\in\operatorname*{SBT}\left(  n\right)  \text{ and }U<T\right\}  .
\]
As we saw in the proof of Lemma \ref{lem.bas.mur.triang.iso1lem} \textbf{(d)},
we have%
\begin{align*}
\MurpF_{\operatorname*{std},\leq T}^{-\operatorname*{Col}}  &
=\operatorname*{span}\underbrace{\left\{  \nabla_{V,U}^{-\operatorname*{Col}%
}\ \mid\ \left(  \mu,U,V\right)  \in\operatorname*{SBT}\left(  n\right)
\text{ and }U\leq T\right\}  }_{\substack{=\left\{  \nabla_{V,U}%
^{-\operatorname*{Col}}\ \mid\ \left(  \mu,U,V\right)  \in\operatorname*{SBT}%
\left(  n\right)  \text{ and }U<T\text{ or }U=T\right\}  \\\text{(since the
condition \textquotedblleft}U\leq T\text{\textquotedblright\ is equivalent to
\textquotedblleft}U<T\text{ or }U=T\text{\textquotedblright)}}}\\
&  =\operatorname*{span}\underbrace{\left\{  \nabla_{V,U}%
^{-\operatorname*{Col}}\ \mid\ \left(  \mu,U,V\right)  \in\operatorname*{SBT}%
\left(  n\right)  \text{ and }U<T\text{ or }U=T\right\}  }%
_{\substack{=\left\{  \nabla_{V,U}^{-\operatorname*{Col}}\ \mid\ \left(
\mu,U,V\right)  \in\operatorname*{SBT}\left(  n\right)  \text{ and
}U<T\right\}  \\\cup\left\{  \nabla_{V,U}^{-\operatorname*{Col}}%
\ \mid\ \left(  \mu,U,V\right)  \in\operatorname*{SBT}\left(  n\right)  \text{
and }U=T\right\}  }}\\
&  =\operatorname*{span}\left(  \left\{  \nabla_{V,U}^{-\operatorname*{Col}%
}\ \mid\ \left(  \mu,U,V\right)  \in\operatorname*{SBT}\left(  n\right)
\text{ and }U<T\right\}  \right. \\
&  \ \ \ \ \ \ \ \ \ \ \ \ \ \ \ \ \ \ \ \ \left.  \cup\left\{  \nabla
_{V,U}^{-\operatorname*{Col}}\ \mid\ \left(  \mu,U,V\right)  \in
\operatorname*{SBT}\left(  n\right)  \text{ and }U=T\right\}  \right) \\
&  =\operatorname*{span}\left\{  \nabla_{V,U}^{-\operatorname*{Col}}%
\ \mid\ \left(  \mu,U,V\right)  \in\operatorname*{SBT}\left(  n\right)  \text{
and }U<T\right\} \\
&  \ \ \ \ \ \ \ \ \ \ +\operatorname*{span}\left\{  \nabla_{V,U}%
^{-\operatorname*{Col}}\ \mid\ \left(  \mu,U,V\right)  \in\operatorname*{SBT}%
\left(  n\right)  \text{ and }U=T\right\}
\end{align*}
(since any two subsets $X$ and $Y$ of a $\mathbf{k}$-module satisfy
$\operatorname*{span}\left(  X\cup Y\right)  =\operatorname*{span}%
X+\operatorname*{span}Y$). In view of%
\[
\MurpF_{\operatorname*{std},<T}^{-\operatorname*{Col}}=\operatorname*{span}%
\left\{  \nabla_{V,U}^{-\operatorname*{Col}}\ \mid\ \left(  \mu,U,V\right)
\in\operatorname*{SBT}\left(  n\right)  \text{ and }U<T\right\}  ,
\]
we can rewrite this as%
\[
\MurpF_{\operatorname*{std},\leq T}^{-\operatorname*{Col}}%
=\MurpF_{\operatorname*{std},<T}^{-\operatorname*{Col}}+\operatorname*{span}%
\left\{  \nabla_{V,U}^{-\operatorname*{Col}}\ \mid\ \left(  \mu,U,V\right)
\in\operatorname*{SBT}\left(  n\right)  \text{ and }U=T\right\}  .
\]

However, we have the two inclusions%
\[
\left\{  \nabla_{V,U}^{-\operatorname*{Col}}\ \mid\ \left(  \mu,U,V\right)
\in\operatorname*{SBT}\left(  n\right)  \text{ and }U=T\right\}
\subseteq\left\{  \nabla_{W,T}^{-\operatorname*{Col}}\ \mid\ W\in
\operatorname*{SYT}\left(  \lambda\right)  \right\}
\]
\footnote{\textit{Proof.} This inclusion is just saying that each element of
the form $\nabla_{V,U}^{-\operatorname*{Col}}$ with $\left(  \mu,U,V\right)
\in\operatorname*{SBT}\left(  n\right)  $ satisfying $U=T$ can also be written
in the form $\nabla_{W,T}^{-\operatorname*{Col}}$ for some $W\in
\operatorname*{SYT}\left(  \lambda\right)  $. So let us prove this.
\par
Let $\left(  \mu,U,V\right)  \in\operatorname*{SBT}\left(  n\right)  $ be such
that $U=T$. We must show that $\nabla_{V,U}^{-\operatorname*{Col}}$ can be
written in the form $\nabla_{W,T}^{-\operatorname*{Col}}$ for some
$W\in\operatorname*{SYT}\left(  \lambda\right)  $.
\par
We have $U=T\in\operatorname*{SYT}\left(  \lambda\right)  $ and $U\in
\operatorname*{SYT}\left(  \mu\right)  $ (since $\left(  \mu,U,V\right)
\in\operatorname*{SBT}\left(  n\right)  $). Hence, the sets
$\operatorname*{SYT}\left(  \lambda\right)  $ and $\operatorname*{SYT}\left(
\mu\right)  $ have at least one element in common (namely, $U$). Thus,
$\lambda=\mu$ (since the sets $\operatorname*{SYT}\left(  \kappa\right)  $ for
different partitions $\kappa$ are disjoint), so that $\mu=\lambda$. However,
from $\left(  \mu,U,V\right)  \in\operatorname*{SBT}\left(  n\right)  $, we
also obtain $V\in\operatorname*{SYT}\left(  \mu\right)  =\operatorname*{SYT}%
\left(  \lambda\right)  $ (since $\mu=\lambda$).
\par
Now we know that $V\in\operatorname*{SYT}\left(  \lambda\right)  $ and
$\nabla_{V,U}^{-\operatorname*{Col}}=\nabla_{V,T}^{-\operatorname*{Col}}$
(since $U=T$). Therefore, $\nabla_{V,U}^{-\operatorname*{Col}}$ can be written
in the form $\nabla_{W,T}^{-\operatorname*{Col}}$ for some $W\in
\operatorname*{SYT}\left(  \lambda\right)  $ (namely, for $W=V$). This is
precisely what we wanted to prove.} and%
\[
\left\{  \nabla_{W,T}^{-\operatorname*{Col}}\ \mid\ W\in\operatorname*{SYT}%
\left(  \lambda\right)  \right\}  \subseteq\left\{  \nabla_{V,U}%
^{-\operatorname*{Col}}\ \mid\ \left(  \mu,U,V\right)  \in\operatorname*{SBT}%
\left(  n\right)  \text{ and }U=T\right\}
\]
\footnote{\textit{Proof.} This inclusion is just saying that each element of
the form $\nabla_{W,T}^{-\operatorname*{Col}}$ with $W\in\operatorname*{SYT}%
\left(  \lambda\right)  $ can also be written in the form $\nabla
_{V,U}^{-\operatorname*{Col}}$ for some $\left(  \mu,U,V\right)
\in\operatorname*{SBT}\left(  n\right)  $ satisfying $U=T$. But this is clear:
We just need to take $\left(  \mu,U,V\right)  =\left(  \lambda,T,W\right)  $
(this is an element of $\operatorname*{SBT}\left(  n\right)  $ because
$T\in\operatorname*{SYT}\left(  \lambda\right)  $ and $W\in\operatorname*{SYT}%
\left(  \lambda\right)  $; and it satisfies $U=T$ by its definition), and then
we obviously have $\nabla_{V,U}^{-\operatorname*{Col}}=\nabla_{W,T}%
^{-\operatorname*{Col}}$, so that we have indeed written the element
$\nabla_{W,T}^{-\operatorname*{Col}}$ in the form $\nabla_{V,U}%
^{-\operatorname*{Col}}$ for some $\left(  \mu,U,V\right)  \in
\operatorname*{SBT}\left(  n\right)  $ satisfying $U=T$.}. Combining these two
inclusions, we obtain the equality%
\begin{align*}
\left\{  \nabla_{V,U}^{-\operatorname*{Col}}\ \mid\ \left(  \mu,U,V\right)
\in\operatorname*{SBT}\left(  n\right)  \text{ and }U=T\right\}   &  =\left\{
\nabla_{W,T}^{-\operatorname*{Col}}\ \mid\ W\in\operatorname*{SYT}\left(
\lambda\right)  \right\} \\
&  =\left\{  \nabla_{V,T}^{-\operatorname*{Col}}\ \mid\ V\in
\operatorname*{SYT}\left(  \lambda\right)  \right\}
\end{align*}
(here, we have renamed the index $W$ as $V$).

Now, recall that%
\begin{align*}
\MurpF_{\operatorname*{std},\leq T}^{-\operatorname*{Col}}  &
=\MurpF_{\operatorname*{std},<T}^{-\operatorname*{Col}}+\operatorname*{span}%
\underbrace{\left\{  \nabla_{V,U}^{-\operatorname*{Col}}\ \mid\ \left(
\mu,U,V\right)  \in\operatorname*{SBT}\left(  n\right)  \text{ and
}U=T\right\}  }_{=\left\{  \nabla_{V,T}^{-\operatorname*{Col}}\ \mid
\ V\in\operatorname*{SYT}\left(  \lambda\right)  \right\}  }\\
&  =\MurpF_{\operatorname*{std},<T}^{-\operatorname*{Col}}%
+\operatorname*{span}\left\{  \nabla_{V,T}^{-\operatorname*{Col}}\ \mid
\ V\in\operatorname*{SYT}\left(  \lambda\right)  \right\}  .
\end{align*}
This proves Lemma \ref{lem.bas.mur.triang.iso1lem} \textbf{(e)}.
\end{proof}
\end{fineprint}

\begin{fineprint}
\begin{proof}
[Detailed proof of Theorem \ref{thm.bas.mur.filtC-subq}.]Set $T:=T_{i}$ and
$\lambda:=\lambda^{\left(  i\right)  }$. Thus, $T\in\operatorname*{SYT}\left(
\lambda\right)  $ (since $T_{i}\in\operatorname*{SYT}\left(  \lambda^{\left(
i\right)  }\right)  $). In other words, $T$ is a standard $n$-tableau of shape
$\lambda$. Hence, Theorem \ref{thm.bas.mur.triang.iso1} yields that
$\MurpF_{\operatorname*{std},\leq T}^{-\operatorname*{Col}}%
/\MurpF_{\operatorname*{std},<T}^{-\operatorname*{Col}}\cong\mathcal{A}%
\mathbf{E}_{T}\cong\mathcal{S}^{\lambda}$ as $S_{n}$-representations, i.e., as
left $\mathcal{A}$-modules.

But the definition of $\MurpF_{i}^{-\operatorname*{Col}}$ yields%
\begin{align*}
\MurpF_{i}^{-\operatorname*{Col}}  &  =%
\begin{cases}
\MurpF_{\operatorname*{std},\leq T_{i}}^{-\operatorname*{Col}}, & \text{if
}i>0;\\
0, & \text{if }i=0
\end{cases}
\ \ =\MurpF_{\operatorname*{std},\leq T_{i}}^{-\operatorname*{Col}%
}\ \ \ \ \ \ \ \ \ \ \left(  \text{since }i>0\right) \\
&  =\MurpF_{\operatorname*{std},\leq T}^{-\operatorname*{Col}}%
\ \ \ \ \ \ \ \ \ \ \left(  \text{since }T_{i}=T\right)  .
\end{align*}
Moreover, $i-1<i\leq m$. Now, Theorem \ref{thm.bas.mur.filt} \textbf{(b)}
(applied to $i-1$ instead of $i$) yields%
\begin{align*}
\MurpF_{i-1}^{-\operatorname*{Col}}  &  =%
\begin{cases}
\MurpF_{\operatorname*{std},<T_{\left(  i-1\right)  +1}}^{-\operatorname*{Col}%
}, & \text{if }i-1<m;\\
\mathcal{A}, & \text{if }i-1=m
\end{cases}
\ \ =\MurpF_{\operatorname*{std},<T_{\left(  i-1\right)  +1}}%
^{-\operatorname*{Col}}\ \ \ \ \ \ \ \ \ \ \left(  \text{since }i-1<m\right)
\\
&  =\MurpF_{\operatorname*{std},<T_{i}}^{-\operatorname*{Col}}%
\ \ \ \ \ \ \ \ \ \ \left(  \text{since }\left(  i-1\right)  +1=i\right) \\
&  =\MurpF_{\operatorname*{std},<T}^{-\operatorname*{Col}}%
\ \ \ \ \ \ \ \ \ \ \left(  \text{since }T_{i}=T\right)  .
\end{align*}
Hence,%
\begin{align*}
\underbrace{\MurpF_{i}^{-\operatorname*{Col}}}_{=\MurpF_{\operatorname*{std}%
,\leq T}^{-\operatorname*{Col}}}/\underbrace{\MurpF_{i-1}%
^{-\operatorname*{Col}}}_{=\MurpF_{\operatorname*{std},<T}%
^{-\operatorname*{Col}}}  &  =\MurpF_{\operatorname*{std},\leq T}%
^{-\operatorname*{Col}}/\MurpF_{\operatorname*{std},<T}^{-\operatorname*{Col}%
}\cong\mathcal{S}^{\lambda}\ \ \ \ \ \ \ \ \ \ \left(  \text{as we proved
above}\right) \\
&  =\mathcal{S}^{\lambda^{\left(  i\right)  }}\ \ \ \ \ \ \ \ \ \ \left(
\text{since }\lambda=\lambda^{\left(  i\right)  }\right)
\end{align*}
as left $\mathcal{A}$-modules. This proves Theorem
\ref{thm.bas.mur.filtC-subq}.
\end{proof}
\end{fineprint}

\begin{fineprint}
\begin{proof}
[Detailed proof of Theorem \ref{thm.bas.mur.twosid.orth-len} using Theorem
\ref{thm.bas.mur.twosid.orth-XY}.]Let $\operatorname*{Par}\left(  n\right)  $
be the set of all partitions of $n$. Define two subsets $X$ and $Y$ of
$\operatorname*{Par}\left(  n\right)  $ by%
\begin{align*}
X  &  :=\left\{  \lambda\in\operatorname*{Par}\left(  n\right)  \ \mid
\ \ell\left(  \lambda\right)  \leq k\right\}  ;\\
Y  &  :=\left\{  \lambda\in\operatorname*{Par}\left(  n\right)  \ \mid
\ \ell\left(  \lambda\right)  >k\right\}  .
\end{align*}
Clearly, each partition $\lambda\in\operatorname*{Par}\left(  n\right)  $
belongs to exactly one of the sets $X$ and $Y$ (since it satisfies exactly one
of the two conditions $\ell\left(  \lambda\right)  \leq k$ and $\ell\left(
\lambda\right)  >k$). In other words, the sets $X$ and $Y$ are disjoint and
their union is $X\cup Y=\operatorname*{Par}\left(  n\right)  $.

Next we shall show the following:

\begin{statement}
\textit{Claim 1:} For every $\lambda\in X$ and every $\mu\in Y$, there exists
some $i\geq1$ satisfying $\lambda_{1}+\lambda_{2}+\cdots+\lambda_{i}>\mu
_{1}+\mu_{2}+\cdots+\mu_{i}$.
\end{statement}

\begin{proof}
[Proof of Claim 1.]Let $\lambda\in X$ and $\mu\in Y$. We must show that there
exists some $i\geq1$ satisfying $\lambda_{1}+\lambda_{2}+\cdots+\lambda
_{i}>\mu_{1}+\mu_{2}+\cdots+\mu_{i}$.

We have $\lambda\in X$. By the definition of $X$, this means that $\lambda
\in\operatorname*{Par}\left(  n\right)  $ and $\ell\left(  \lambda\right)
\leq k$. From $\lambda\in\operatorname*{Par}\left(  n\right)  $, we see that
$\lambda$ is a partition of $n$. Thus, $\left\vert \lambda\right\vert =n$.

We have $\mu\in Y$. By the definition of $Y$, this means that $\mu
\in\operatorname*{Par}\left(  n\right)  $ and $\ell\left(  \mu\right)  >k$.
From $\mu\in\operatorname*{Par}\left(  n\right)  $, we see that $\mu$ is a
partition of $n$. Thus, $\left\vert \mu\right\vert =n$.

Write $\lambda$ as $\lambda=\left(  \lambda_{1},\lambda_{2},\ldots,\lambda
_{i}\right)  $, where $\lambda_{1},\lambda_{2},\ldots,\lambda_{i}$ are
positive integers. Thus, $i=\ell\left(  \lambda\right)  \leq k$.

Write $\mu$ as $\mu=\left(  \mu_{1},\mu_{2},\ldots,\mu_{j}\right)  $, where
$\mu_{1},\mu_{2},\ldots,\mu_{j}$ are positive integers. Thus, $j=\ell\left(
\mu\right)  >k\geq i$ (since $i\leq k$).

Recall that $\mu$ is a partition of $n$. Hence, if we had $n=0$, then we would
have $\mu=\left(  {}\right)  $ (since the only partition of $0$ is the empty
partition $\left(  {}\right)  $), which would entail $\ell\left(  \mu\right)
=\ell\left(  \left(  {}\right)  \right)  =0$, which would contradict
$\ell\left(  \mu\right)  >k\geq0$. Hence, we cannot have $n=0$. Thus, $n>0$.
Hence, $\left\vert \lambda\right\vert =n>0$, so that the partition $\lambda$
cannot be empty. In other words, $\ell\left(  \lambda\right)  \geq1$. Hence,
$i=\ell\left(  \lambda\right)  \geq1$.

Recall that $\mu_{1},\mu_{2},\ldots,\mu_{j}$ are positive integers. Hence,
$\mu_{i+1}+\mu_{i+2}+\cdots+\mu_{j}$ is a nonempty sum of positive integers
(nonempty because $j>i$), and therefore itself positive. In other words,
$\mu_{i+1}+\mu_{i+2}+\cdots+\mu_{j}>0$.

From $\lambda=\left(  \lambda_{1},\lambda_{2},\ldots,\lambda_{i}\right)  $, we
obtain $\left\vert \lambda\right\vert =\lambda_{1}+\lambda_{2}+\cdots
+\lambda_{i}$. Thus, $\lambda_{1}+\lambda_{2}+\cdots+\lambda_{i}=\left\vert
\lambda\right\vert =n$. Likewise, $\mu_{1}+\mu_{2}+\cdots+\mu_{j}=n$.
Comparing these two equalities, we find%
\begin{align*}
&  \lambda_{1}+\lambda_{2}+\cdots+\lambda_{i}\\
&  =\mu_{1}+\mu_{2}+\cdots+\mu_{j}\\
&  =\left(  \mu_{1}+\mu_{2}+\cdots+\mu_{i}\right)  +\underbrace{\left(
\mu_{i+1}+\mu_{i+2}+\cdots+\mu_{j}\right)  }_{>0}\ \ \ \ \ \ \ \ \ \ \left(
\text{since }j>i\right) \\
&  >\mu_{1}+\mu_{2}+\cdots+\mu_{i}.
\end{align*}
Hence, we have found an $i\geq1$ satisfying $\lambda_{1}+\lambda_{2}%
+\cdots+\lambda_{i}>\mu_{1}+\mu_{2}+\cdots+\mu_{i}$. Thus, such an $i$ exists.
This proves Claim 1.
\end{proof}

Now, our sets $X$ and $Y$ satisfy the assumptions of Theorem
\ref{thm.bas.mur.twosid.orth-XY} (since they are disjoint and have union
$X\cup Y=\operatorname*{Par}\left(  n\right)  $, and since Claim 1 shows that
the non-dominance condition holds). Thus, we can apply Theorem
\ref{thm.bas.mur.twosid.orth-XY} to our $X$ and $Y$. Let us do this.

Define the four $\mathbf{k}$-submodules $\MurpF_{\operatorname*{std}%
,X}^{\operatorname*{Row}}$, $\MurpF_{\operatorname*{std},Y}%
^{-\operatorname*{Col}}$, $\MurpF_{\operatorname*{all},X}^{\operatorname*{Row}%
}$ and $\MurpF_{\operatorname*{all},Y}^{-\operatorname*{Col}}$ of
$\mathcal{A}$ as in Theorem \ref{thm.bas.mur.twosid.orth-XY}. Recall that $X$
is defined as the set $\left\{  \lambda\in\operatorname*{Par}\left(  n\right)
\ \mid\ \ell\left(  \lambda\right)  \leq k\right\}  $. Thus, a partition
$\lambda$ of $n$ satisfies $\lambda\in X$ if and only if $\ell\left(
\lambda\right)  \leq k$. Hence, an $n$-bitableau $\left(  \lambda,U,V\right)
\in\operatorname*{SBT}\left(  n\right)  $ satisfies $\lambda\in X$ if and only
if $\ell\left(  \lambda\right)  \leq k$ (since $\lambda$ is a partition of $n$).

However, the definition of $\MurpF_{\operatorname*{std},X}%
^{\operatorname*{Row}}$ says that%
\begin{align*}
\MurpF_{\operatorname*{std},X}^{\operatorname*{Row}}  &  =\operatorname*{span}%
\left\{  \nabla_{V,U}^{\operatorname*{Row}}\ \mid\ \left(  \lambda,U,V\right)
\in\operatorname*{SBT}\left(  n\right)  \text{ and }\lambda\in X\right\} \\
&  =\operatorname*{span}\left\{  \nabla_{V,U}^{\operatorname*{Row}}%
\ \mid\ \left(  \lambda,U,V\right)  \in\operatorname*{SBT}\left(  n\right)
\text{ and }\ell\left(  \lambda\right)  \leq k\right\} \\
&  \ \ \ \ \ \ \ \ \ \ \ \ \ \ \ \ \ \ \ \ \left(
\begin{array}
[c]{c}%
\text{here, we have replaced the condition \textquotedblleft}\lambda\in
X\text{\textquotedblright}\\
\text{by the equivalent condition \textquotedblleft}\ell\left(  \lambda
\right)  \leq k\text{\textquotedblright, since}\\
\text{an }n\text{-bitableau }\left(  \lambda,U,V\right)  \in
\operatorname*{SBT}\left(  n\right)  \text{ satisfies }\lambda\in X\\
\text{if and only if }\ell\left(  \lambda\right)  \leq k
\end{array}
\right) \\
&  =\MurpF_{\operatorname*{std},\operatorname*{len}\leq k}%
^{\operatorname*{Row}}\ \ \ \ \ \ \ \ \ \ \left(  \text{by the definition of
}\MurpF_{\operatorname*{std},\operatorname*{len}\leq k}^{\operatorname*{Row}%
}\right)  .
\end{align*}
Thus, we have proved that%
\begin{equation}
\MurpF_{\operatorname*{std},X}^{\operatorname*{Row}}%
=\MurpF_{\operatorname*{std},\operatorname*{len}\leq k}^{\operatorname*{Row}}.
\label{pf.thm.bas.mur.twosid.orth-len.=1}%
\end{equation}
Similarly, we can find that
\begin{align}
\MurpF_{\operatorname*{std},Y}^{-\operatorname*{Col}}  &
=\MurpF_{\operatorname*{std},\operatorname*{len}>k}^{-\operatorname*{Col}%
}\ \ \ \ \ \ \ \ \ \ \text{and}\label{pf.thm.bas.mur.twosid.orth-len.=2}\\
\MurpF_{\operatorname*{all},X}^{\operatorname*{Row}}  &
=\MurpF_{\operatorname*{all},\operatorname*{len}\leq k}^{\operatorname*{Row}%
}\ \ \ \ \ \ \ \ \ \ \text{and}\label{pf.thm.bas.mur.twosid.orth-len.=3}\\
\MurpF_{\operatorname*{all},Y}^{-\operatorname*{Col}}  &
=\MurpF_{\operatorname*{all},\operatorname*{len}>k}^{-\operatorname*{Col}}.
\label{pf.thm.bas.mur.twosid.orth-len.=4}%
\end{align}

Now, we can easily prove both parts of Theorem
\ref{thm.bas.mur.twosid.orth-len}: \medskip

\textbf{(a)} Theorem \ref{thm.bas.mur.twosid.orth-XY} \textbf{(a)} says that
the set $\MurpF_{\operatorname*{std},X}^{\operatorname*{Row}}$ is a two-sided
ideal of $\mathcal{A}$ and satisfies%
\[
\MurpF_{\operatorname*{std},X}^{\operatorname*{Row}}%
=\MurpF_{\operatorname*{all},X}^{\operatorname*{Row}}=\left(
\MurpF_{\operatorname*{std},Y}^{-\operatorname*{Col}}\right)  ^{\perp}.
\]
In view of (\ref{pf.thm.bas.mur.twosid.orth-len.=1}),
(\ref{pf.thm.bas.mur.twosid.orth-len.=3}) and
(\ref{pf.thm.bas.mur.twosid.orth-len.=2}), we can rewrite this as follows: The
set $\MurpF_{\operatorname*{std},\operatorname*{len}\leq k}%
^{\operatorname*{Row}}$ is a two-sided ideal of $\mathcal{A}$ and satisfies%
\[
\MurpF_{\operatorname*{std},\operatorname*{len}\leq k}^{\operatorname*{Row}%
}=\MurpF_{\operatorname*{all},\operatorname*{len}\leq k}^{\operatorname*{Row}%
}=\left(  \MurpF_{\operatorname*{std},\operatorname*{len}>k}%
^{-\operatorname*{Col}}\right)  ^{\perp}.
\]
Thus, Theorem \ref{thm.bas.mur.twosid.orth-len} \textbf{(a)} is proved.
\medskip

\textbf{(b)} This follows from Theorem \ref{thm.bas.mur.twosid.orth-XY}
\textbf{(b)} in the same way as we derived Theorem
\ref{thm.bas.mur.twosid.orth-len} \textbf{(a)} from Theorem
\ref{thm.bas.mur.twosid.orth-XY} \textbf{(a)} (but of course, we need to use
(\ref{pf.thm.bas.mur.twosid.orth-len.=4}) now).
\end{proof}
\end{fineprint}

\begin{fineprint}
\begin{proof}
[Detailed proof of Proposition \ref{prop.bas.mur.twosid.Sn-act.all}.]Recall
that we defined $\MurpF_{\operatorname*{all},X}^{\operatorname*{Row}}$ as
\[
\MurpF_{\operatorname*{all},X}^{\operatorname*{Row}}=\operatorname*{span}%
\left\{  \nabla_{V,U}^{\operatorname*{Row}}\ \mid\ \left(  \lambda,U,V\right)
\in\operatorname*{BT}\left(  n\right)  \text{ and }\lambda\in X\right\}  .
\]
In other words, $\MurpF_{\operatorname*{all},X}^{\operatorname*{Row}}$ is the
span of all the vectors of the form $\nabla_{V,U}^{\operatorname*{Row}}$ with
$\left(  \lambda,U,V\right)  \in\operatorname*{BT}\left(  n\right)  $
satisfying $\lambda\in X$. We shall refer to such vectors as the
\textquotedblleft good nablas\textquotedblright. Thus,
$\MurpF_{\operatorname*{all},X}^{\operatorname*{Row}}$ is the span of all the
good nablas. Hence, $\MurpF_{\operatorname*{all},X}^{\operatorname*{Row}}$ is
a $\mathbf{k}$-submodule of $\mathcal{A}$.

We shall now show that $\MurpF_{\operatorname*{all},X}^{\operatorname*{Row}}$
is a (two-sided) ideal of $\mathcal{A}$. For this purpose, we will show the
following two claims:

\begin{statement}
\textit{Claim 1:} We have $\mathbf{ca}\in\MurpF_{\operatorname*{all}%
,X}^{\operatorname*{Row}}$ for all $\mathbf{c}\in\mathcal{A}$ and all
$\mathbf{a}\in\MurpF_{\operatorname*{all},X}^{\operatorname*{Row}}$.
\end{statement}

\begin{statement}
\textit{Claim 2:} We have $\mathbf{ac}\in\MurpF_{\operatorname*{all}%
,X}^{\operatorname*{Row}}$ for all $\mathbf{c}\in\mathcal{A}$ and all
$\mathbf{a}\in\MurpF_{\operatorname*{all},X}^{\operatorname*{Row}}$.
\end{statement}

\begin{proof}
[Proof of Claim 1.]Let $\mathbf{c}\in\mathcal{A}$ and $\mathbf{a}%
\in\MurpF_{\operatorname*{all},X}^{\operatorname*{Row}}$. We must prove that
$\mathbf{ca}\in\MurpF_{\operatorname*{all},X}^{\operatorname*{Row}}$. Since
this claim depends $\mathbf{k}$-linearly on each of $\mathbf{c}$ and
$\mathbf{a}$ (because $\MurpF_{\operatorname*{all},X}^{\operatorname*{Row}}$
is a $\mathbf{k}$-submodule of $\mathcal{A}$), we can WLOG assume (by
linearity) that $\mathbf{c}$ is one of the standard basis vectors $w\in S_{n}$
of $\mathcal{A}=\mathbf{k}\left[  S_{n}\right]  $ (since $\mathcal{A}$ is the
span of these basis vectors), and that $\mathbf{a}$ is one of the good nablas
(since $\MurpF_{\operatorname*{all},X}^{\operatorname*{Row}}$ is the span of
all the good nablas). Let us assume this. Thus, $\mathbf{c}=w$ for some $w\in
S_{n}$, whereas $\mathbf{a}=\nabla_{V,U}^{\operatorname*{Row}}$ for some
$\left(  \lambda,U,V\right)  \in\operatorname*{BT}\left(  n\right)  $
satisfying $\lambda\in X$ (because this is what a good nabla looks like).
Consider these $w$ and $\left(  \lambda,U,V\right)  $.

Now, $\left(  \lambda,U,V\right)  \in\operatorname*{BT}\left(  n\right)  $
shows that $\left(  \lambda,U,V\right)  $ is an $n$-bitableau. In other words,
$\lambda$ is a partition of $n$ and $U$ and $V$ are two $n$-tableaux of shape
$\lambda$. Hence, of course, $wV$ is an $n$-tableau of shape $\lambda$ as
well. Hence, $\left(  \lambda,U,wV\right)  $ is again an $n$-bitableau, i.e.,
an element of $\operatorname*{BT}\left(  n\right)  $. Therefore,
$\nabla_{wV,U}^{\operatorname*{Row}}$ is one of the good nablas as well (since
$\lambda\in X$), and hence belongs to $\MurpF_{\operatorname*{all}%
,X}^{\operatorname*{Row}}$ (since $\MurpF_{\operatorname*{all},X}%
^{\operatorname*{Row}}$ is the span of all the good nablas). In other words,
$\nabla_{wV,U}^{\operatorname*{Row}}\in\MurpF_{\operatorname*{all}%
,X}^{\operatorname*{Row}}$.

But (\ref{eq.prop.bas.mur.Sn-act.Row}) (applied to $P=V$ and $Q=U$ and $g=w$
and $h=\operatorname*{id}$) yields $\nabla_{wV,\operatorname*{id}%
U}^{\operatorname*{Row}}=w\nabla_{V,U}^{\operatorname*{Row}}%
\underbrace{\operatorname*{id}\nolimits^{-1}}_{=\operatorname*{id}=1}%
=w\nabla_{V,U}^{\operatorname*{Row}}$. In other words, $\nabla_{wV,U}%
^{\operatorname*{Row}}=w\nabla_{V,U}^{\operatorname*{Row}}$ (since
$\operatorname*{id}U=U$). Hence, $w\nabla_{V,U}^{\operatorname*{Row}}%
=\nabla_{wV,U}^{\operatorname*{Row}}\in\MurpF_{\operatorname*{all}%
,X}^{\operatorname*{Row}}$. In other words, $\mathbf{ca}\in
\MurpF_{\operatorname*{all},X}^{\operatorname*{Row}}$ (since $\mathbf{c}=w$
and $\mathbf{a}=\nabla_{V,U}^{\operatorname*{Row}}$). This proves Claim 1.
\end{proof}

\begin{proof}
[Proof of Claim 2.]Let $\mathbf{c}\in\mathcal{A}$ and $\mathbf{a}%
\in\MurpF_{\operatorname*{all},X}^{\operatorname*{Row}}$. We must prove that
$\mathbf{ac}\in\MurpF_{\operatorname*{all},X}^{\operatorname*{Row}}$. Since
this claim depends $\mathbf{k}$-linearly on each of $\mathbf{c}$ and
$\mathbf{a}$ (because $\MurpF_{\operatorname*{all},X}^{\operatorname*{Row}}$
is a $\mathbf{k}$-submodule of $\mathcal{A}$), we can WLOG assume (by
linearity) that $\mathbf{c}$ is one of the standard basis vectors $w\in S_{n}$
of $\mathcal{A}=\mathbf{k}\left[  S_{n}\right]  $ (since $\mathcal{A}$ is the
span of these basis vectors), and that $\mathbf{a}$ is one of the good nablas
(since $\MurpF_{\operatorname*{all},X}^{\operatorname*{Row}}$ is the span of
all the good nablas). Let us assume this. Thus, $\mathbf{c}=w$ for some $w\in
S_{n}$, whereas $\mathbf{a}=\nabla_{V,U}^{\operatorname*{Row}}$ for some
$\left(  \lambda,U,V\right)  \in\operatorname*{BT}\left(  n\right)  $
satisfying $\lambda\in X$ (because this is what a good nabla looks like).
Consider these $w$ and $\left(  \lambda,U,V\right)  $.

Now, $\left(  \lambda,U,V\right)  \in\operatorname*{BT}\left(  n\right)  $
shows that $\left(  \lambda,U,V\right)  $ is an $n$-bitableau. In other words,
$\lambda$ is a partition of $n$ and $U$ and $V$ are two $n$-tableaux of shape
$\lambda$. Hence, of course, $w^{-1}U$ is an $n$-tableau of shape $\lambda$ as
well. Hence, $\left(  \lambda,w^{-1}U,V\right)  $ is again an $n$-bitableau,
i.e., an element of $\operatorname*{BT}\left(  n\right)  $. Therefore,
$\nabla_{V,w^{-1}U}^{\operatorname*{Row}}$ is one of the good nablas as well
(since $\lambda\in X$), and hence belongs to $\MurpF_{\operatorname*{all}%
,X}^{\operatorname*{Row}}$ (since $\MurpF_{\operatorname*{all},X}%
^{\operatorname*{Row}}$ is the span of all the good nablas). In other words,
$\nabla_{V,w^{-1}U}^{\operatorname*{Row}}\in\MurpF_{\operatorname*{all}%
,X}^{\operatorname*{Row}}$.

But (\ref{eq.prop.bas.mur.Sn-act.Row}) (applied to $P=V$ and $Q=U$ and
$g=\operatorname*{id}$ and $h=w^{-1}$) yields $\nabla_{\operatorname*{id}%
V,w^{-1}U}^{\operatorname*{Row}}=\underbrace{\operatorname*{id}}_{=1}%
\nabla_{V,U}^{\operatorname*{Row}}\underbrace{\left(  w^{-1}\right)  ^{-1}%
}_{=w}=\nabla_{V,U}^{\operatorname*{Row}}w$. In other words, $\nabla
_{V,w^{-1}U}^{\operatorname*{Row}}=\nabla_{V,U}^{\operatorname*{Row}}w$ (since
$\operatorname*{id}V=V$). Hence, $\nabla_{V,U}^{\operatorname*{Row}}%
w=\nabla_{V,w^{-1}U}^{\operatorname*{Row}}\in\MurpF_{\operatorname*{all}%
,X}^{\operatorname*{Row}}$. In other words, $\mathbf{ac}\in
\MurpF_{\operatorname*{all},X}^{\operatorname*{Row}}$ (since $\mathbf{c}=w$
and $\mathbf{a}=\nabla_{V,U}^{\operatorname*{Row}}$). This proves Claim 2.
\end{proof}

Now, recall that $\MurpF_{\operatorname*{all},X}^{\operatorname*{Row}}$ is a
$\mathbf{k}$-submodule of $\mathcal{A}$. Combining this with Claim 1 and Claim
2, we conclude that $\MurpF_{\operatorname*{all},X}^{\operatorname*{Row}}$ is
a two-sided ideal of $\mathcal{A}$.

Likewise, we can show that $\MurpF_{\operatorname*{all},Y}%
^{-\operatorname*{Col}}$ is a two-sided ideal of $\mathcal{A}$. (This requires
us to use (\ref{eq.prop.bas.mur.Sn-act.Col}) instead of
(\ref{eq.prop.bas.mur.Sn-act.Row}), which is slightly more complicated because
the equality (\ref{eq.prop.bas.mur.Sn-act.Col}) involves the factors $\left(
-1\right)  ^{g}$ and $\left(  -1\right)  ^{h}$ that don't appear in
(\ref{eq.prop.bas.mur.Sn-act.Row}). However, these factors are invertible
scalars, so we can divide a vector by them without affecting its property of
belonging to $\MurpF_{\operatorname*{all},Y}^{-\operatorname*{Col}}$.)

Altogether, Proposition \ref{prop.bas.mur.twosid.Sn-act.all} is thus proved.
\end{proof}
\end{fineprint}

\begin{fineprint}
\begin{proof}
[Detailed proof of Lemma \ref{lem.bas.mur.twosid.triang.1}.]\textbf{(a)} We
must prove that $\mathbf{b}\in\left(  \MurpF_{\operatorname*{std}%
,X}^{\operatorname*{Row}}\right)  ^{\perp}$ for each $\mathbf{b}%
\in\MurpF_{\operatorname*{all},Y}^{-\operatorname*{Col}}$.

So let $\mathbf{b}\in\MurpF_{\operatorname*{all},Y}^{-\operatorname*{Col}}$ be
arbitrary. We must then prove that $\mathbf{b}\in\left(
\MurpF_{\operatorname*{std},X}^{\operatorname*{Row}}\right)  ^{\perp}$. In
other words, we must prove that $\left\langle \mathbf{b},\mathbf{v}%
\right\rangle =0$ for each $\mathbf{v}\in\MurpF_{\operatorname*{std}%
,X}^{\operatorname*{Row}}$ (because $\left(  \MurpF_{\operatorname*{std}%
,X}^{\operatorname*{Row}}\right)  ^{\perp}$ is defined as $\left\{
\mathbf{a}\in\mathcal{A}\ \mid\ \left\langle \mathbf{a},\mathbf{v}%
\right\rangle =0\text{ for each }\mathbf{v}\in\MurpF_{\operatorname*{std}%
,X}^{\operatorname*{Row}}\right\}  $).

So let $\mathbf{v}\in\MurpF_{\operatorname*{std},X}^{\operatorname*{Row}}$ be
arbitrary. We must prove that $\left\langle \mathbf{b},\mathbf{v}\right\rangle
=0$.

We have
\begin{align*}
\mathbf{v}  &  \in\MurpF_{\operatorname*{std},X}^{\operatorname*{Row}%
}=\operatorname*{span}\left\{  \nabla_{V,U}^{\operatorname*{Row}}%
\ \mid\ \left(  \lambda,U,V\right)  \in\underbrace{\operatorname*{SBT}\left(
n\right)  }_{\subseteq\operatorname*{BT}\left(  n\right)  }\text{ and }%
\lambda\in X\right\} \\
&  \ \ \ \ \ \ \ \ \ \ \ \ \ \ \ \ \ \ \ \ \left(  \text{by the definition of
}\MurpF_{\operatorname*{std},X}^{\operatorname*{Row}}\right) \\
&  \subseteq\operatorname*{span}\left\{  \nabla_{V,U}^{\operatorname*{Row}%
}\ \mid\ \left(  \lambda,U,V\right)  \in\operatorname*{BT}\left(  n\right)
\text{ and }\lambda\in X\right\}  =\MurpF_{\operatorname*{all},X}%
^{\operatorname*{Row}}%
\end{align*}
(by the definition of $\MurpF_{\operatorname*{all},X}^{\operatorname*{Row}}$).
Hence, Lemma \ref{lem.bas.mur.twosid.orth} \textbf{(b)} (applied to
$\mathbf{b}$ and $\mathbf{v}$ instead of $\mathbf{a}$ and $\mathbf{b}$) yields
$\left\langle \mathbf{b},\mathbf{v}\right\rangle =0$. As we explained above,
this completes the proof of Lemma \ref{lem.bas.mur.twosid.triang.1}
\textbf{(a)}. \medskip

\textbf{(b)} We must prove that $\mathbf{b}\in\left(
\MurpF_{\operatorname*{std},Y}^{-\operatorname*{Col}}\right)  ^{\perp}$ for
each $\mathbf{b}\in\MurpF_{\operatorname*{all},X}^{\operatorname*{Row}}$.

So let $\mathbf{b}\in\MurpF_{\operatorname*{all},X}^{\operatorname*{Row}}$ be
arbitrary. We must then prove that $\mathbf{b}\in\left(
\MurpF_{\operatorname*{std},Y}^{-\operatorname*{Col}}\right)  ^{\perp}$. In
other words, we must prove that $\left\langle \mathbf{b},\mathbf{v}%
\right\rangle =0$ for each $\mathbf{v}\in\MurpF_{\operatorname*{std}%
,Y}^{-\operatorname*{Col}}$ (because $\left(  \MurpF_{\operatorname*{std}%
,Y}^{-\operatorname*{Col}}\right)  ^{\perp}$ is defined as $\left\{
\mathbf{a}\in\mathcal{A}\ \mid\ \left\langle \mathbf{a},\mathbf{v}%
\right\rangle =0\text{ for each }\mathbf{v}\in\MurpF_{\operatorname*{std}%
,Y}^{-\operatorname*{Col}}\right\}  $).

So let $\mathbf{v}\in\MurpF_{\operatorname*{std},Y}^{-\operatorname*{Col}}$ be
arbitrary. We must prove that $\left\langle \mathbf{b},\mathbf{v}\right\rangle
=0$.

We have
\begin{align*}
\mathbf{v}  &  \in\MurpF_{\operatorname*{std},Y}^{-\operatorname*{Col}%
}=\operatorname*{span}\left\{  \nabla_{V,U}^{-\operatorname*{Col}}%
\ \mid\ \left(  \lambda,U,V\right)  \in\underbrace{\operatorname*{SBT}\left(
n\right)  }_{\subseteq\operatorname*{BT}\left(  n\right)  }\text{ and }%
\lambda\in Y\right\} \\
&  \ \ \ \ \ \ \ \ \ \ \ \ \ \ \ \ \ \ \ \ \left(  \text{by the definition of
}\MurpF_{\operatorname*{std},Y}^{-\operatorname*{Col}}\right) \\
&  \subseteq\operatorname*{span}\left\{  \nabla_{V,U}^{-\operatorname*{Col}%
}\ \mid\ \left(  \lambda,U,V\right)  \in\operatorname*{BT}\left(  n\right)
\text{ and }\lambda\in Y\right\}  =\MurpF_{\operatorname*{all},Y}%
^{-\operatorname*{Col}}%
\end{align*}
(by the definition of $\MurpF_{\operatorname*{all},Y}^{-\operatorname*{Col}}%
$). Hence, Lemma \ref{lem.bas.mur.twosid.orth} \textbf{(b)} (applied to
$\mathbf{a}=\mathbf{v}$) yields $\left\langle \mathbf{v},\mathbf{b}%
\right\rangle =0$. Finally, (\ref{eq.prop.bas.bilin.sym-etc.sym}) shows that
$\left\langle \mathbf{v},\mathbf{b}\right\rangle =\left\langle \mathbf{b}%
,\mathbf{v}\right\rangle $. Hence, $\left\langle \mathbf{b},\mathbf{v}%
\right\rangle =\left\langle \mathbf{v},\mathbf{b}\right\rangle =0$. As we
explained above, this completes the proof of Lemma
\ref{lem.bas.mur.twosid.triang.1} \textbf{(b)}.
\end{proof}
\end{fineprint}

\begin{fineprint}
\begin{proof}
[Detailed proof of Lemma \ref{lem.bas.mur.twosid.triang.2}.]According to the
assumptions of Theorem \ref{thm.bas.mur.twosid.orth-XY}, we know that the two
sets $X$ and $Y$ are disjoint and have union $X\cup Y=\operatorname*{Par}%
\left(  n\right)  $. Hence, each $\lambda\in\operatorname*{Par}\left(
n\right)  $ satisfies either $\lambda\in X$ or $\lambda\in Y$, but not both.
Therefore, each $\left(  \lambda,U,V\right)  \in\operatorname*{SBT}\left(
n\right)  $ satisfies either $\lambda\in X$ or $\lambda\in Y$, but not both
(because $\left(  \lambda,U,V\right)  \in\operatorname*{SBT}\left(  n\right)
$ entails $\lambda\in\operatorname*{Par}\left(  n\right)  $). \medskip

\textbf{(b)} We must show that $\mathbf{a}\in\MurpF_{\operatorname*{std}%
,X}^{\operatorname*{Row}}$ for each $\mathbf{a}\in\left(
\MurpF_{\operatorname*{std},Y}^{-\operatorname*{Col}}\right)  ^{\perp}$.

So let $\mathbf{a}\in\left(  \MurpF_{\operatorname*{std},Y}%
^{-\operatorname*{Col}}\right)  ^{\perp}$. We must show that $\mathbf{a}%
\in\MurpF_{\operatorname*{std},X}^{\operatorname*{Row}}$.

We have $\mathbf{a}\in\left(  \MurpF_{\operatorname*{std},Y}%
^{-\operatorname*{Col}}\right)  ^{\perp}\subseteq\mathcal{A}$. Since the row
Murphy basis $\left(  \nabla_{V,U}^{\operatorname*{Row}}\right)  _{\left(
\lambda,U,V\right)  \in\operatorname*{SBT}\left(  n\right)  }$ is a basis of
$\mathcal{A}$ (by Theorem \ref{thm.bas.mur.bas}), we can thus write
$\mathbf{a}$ as a $\mathbf{k}$-linear combination of its elements. In other
words, we can write $\mathbf{a}$ as
\begin{equation}
\mathbf{a}=\sum_{\left(  \lambda,U,V\right)  \in\operatorname*{SBT}\left(
n\right)  }\omega_{\lambda,U,V}\nabla_{V,U}^{\operatorname*{Row}}
\label{pf.lem.bas.mur.twosid.triang.2.a=wsum}%
\end{equation}
for some scalars $\omega_{\lambda,U,V}\in\mathbf{k}$. Consider these
$\omega_{\lambda,U,V}$. We can rewrite the equality
(\ref{pf.lem.bas.mur.twosid.triang.2.a=wsum}) as
\begin{equation}
\mathbf{a}=\sum_{\left(  \mu,P,Q\right)  \in\operatorname*{SBT}\left(
n\right)  }\omega_{\mu,P,Q}\nabla_{Q,P}^{\operatorname*{Row}}
\label{pf.lem.bas.mur.twosid.triang.2.a=wsum2}%
\end{equation}
(by renaming the summation index $\left(  \lambda,U,V\right)  $ as $\left(
\mu,P,Q\right)  $).

Now, we shall show the following:

\begin{statement}
\textit{Claim 1:} We have%
\[
\omega_{\lambda,U,V}=0\ \ \ \ \ \ \ \ \ \ \text{for all }\left(
\lambda,U,V\right)  \in\operatorname*{SBT}\left(  n\right)  \text{ satisfying
}\lambda\in Y.
\]

\end{statement}

\begin{proof}
[Proof of Claim 1.]We proceed by strong induction on $\left(  \lambda
,U,V\right)  $ (using the total order on $\operatorname*{SBT}\left(  n\right)
$ introduced in Definition \ref{def.bas.tord.sbt.tord}). So we fix some
$\left(  \lambda,U,V\right)  \in\operatorname*{SBT}\left(  n\right)  $
satisfying $\lambda\in Y$. We assume (as the induction hypothesis) that
\begin{equation}
\omega_{\mu,P,Q}=0 \label{pf.lem.bas.mur.twosid.triang.2.IH}%
\end{equation}
for all $\left(  \mu,P,Q\right)  \in\operatorname*{SBT}\left(  n\right)  $
satisfying $\mu\in Y$ and $\left(  \mu,P,Q\right)  <\left(  \lambda
,U,V\right)  $. Our goal is to prove that $\omega_{\lambda,U,V}=0$.

Since $\left(  \lambda,U,V\right)  \in\operatorname*{SBT}\left(  n\right)  $
and $\lambda\in Y$, we have $\nabla_{V,U}^{-\operatorname*{Col}}%
\in\MurpF_{\operatorname*{std},Y}^{-\operatorname*{Col}}$ by the definition of
$\MurpF_{\operatorname*{std},Y}^{-\operatorname*{Col}}$ (in fact,
$\nabla_{V,U}^{-\operatorname*{Col}}$ is one of the spanning vectors in the
definition of $\MurpF_{\operatorname*{std},Y}^{-\operatorname*{Col}}$).

Since $\mathbf{a}\in\left(  \MurpF_{\operatorname*{std},Y}%
^{-\operatorname*{Col}}\right)  ^{\perp}$, we have $\left\langle
\mathbf{a},\mathbf{v}\right\rangle =0$ for each $\mathbf{v}\in
\MurpF_{\operatorname*{std},Y}^{-\operatorname*{Col}}$ (by the definition of
$\left(  \MurpF_{\operatorname*{std},Y}^{-\operatorname*{Col}}\right)
^{\perp}$). Applying this to $\mathbf{v}=\nabla_{V,U}^{-\operatorname*{Col}}$,
we obtain $\left\langle \mathbf{a},\ \nabla_{V,U}^{-\operatorname*{Col}%
}\right\rangle =0$ (since $\nabla_{V,U}^{-\operatorname*{Col}}\in
\MurpF_{\operatorname*{std},Y}^{-\operatorname*{Col}}$). But
(\ref{eq.prop.bas.bilin.sym-etc.sym}) shows that $\left\langle \mathbf{a}%
,\ \nabla_{V,U}^{-\operatorname*{Col}}\right\rangle =\left\langle \nabla
_{V,U}^{-\operatorname*{Col}},\ \mathbf{a}\right\rangle $, hence $\left\langle
\nabla_{V,U}^{-\operatorname*{Col}},\ \mathbf{a}\right\rangle =\left\langle
\mathbf{a},\ \nabla_{V,U}^{-\operatorname*{Col}}\right\rangle =0$. Thus,%
\begin{align}
0  &  =\left\langle \nabla_{V,U}^{-\operatorname*{Col}},\ \mathbf{a}%
\right\rangle =\left\langle \nabla_{V,U}^{-\operatorname*{Col}},\ \sum
_{\left(  \mu,P,Q\right)  \in\operatorname*{SBT}\left(  n\right)  }\omega
_{\mu,P,Q}\nabla_{Q,P}^{\operatorname*{Row}}\right\rangle
\ \ \ \ \ \ \ \ \ \ \left(  \text{by
(\ref{pf.lem.bas.mur.twosid.triang.2.a=wsum2})}\right) \nonumber\\
&  =\sum_{\left(  \mu,P,Q\right)  \in\operatorname*{SBT}\left(  n\right)
}\omega_{\mu,P,Q}\left\langle \nabla_{V,U}^{-\operatorname*{Col}}%
,\ \nabla_{Q,P}^{\operatorname*{Row}}\right\rangle
\label{pf.lem.bas.mur.twosid.triang.2.0}%
\end{align}
(since the form $\left\langle \cdot,\cdot\right\rangle $ is bilinear). Now, we
observe the following:

\begin{itemize}
\item For each $\left(  \mu,P,Q\right)  \in\operatorname*{SBT}\left(
n\right)  $ satisfying $\mu\notin Y$, then%
\begin{equation}
\left\langle \nabla_{V,U}^{-\operatorname*{Col}},\ \nabla_{Q,P}%
^{\operatorname*{Row}}\right\rangle =0.
\label{pf.lem.bas.mur.twosid.triang.2.1b}%
\end{equation}

[\textit{Proof:} Let $\left(  \mu,P,Q\right)  \in\operatorname*{SBT}\left(
n\right)  $ be such that $\mu\notin Y$. Hence, $\mu\in\operatorname*{Par}%
\left(  n\right)  \setminus Y$ (since $\mu\in\operatorname*{Par}\left(
n\right)  $ and $\mu\notin Y$). Recall that $X\cup Y=\operatorname*{Par}%
\left(  n\right)  $, so that $\operatorname*{Par}\left(  n\right)  \setminus
Y\subseteq X$. Hence, $\mu\in\operatorname*{Par}\left(  n\right)  \setminus
Y\subseteq X$. Therefore, Lemma \ref{lem.bas.mur.twosid.orth} \textbf{(a)}
(applied to $\left(  \mu,P,Q\right)  $ and $\left(  \lambda,U,V\right)  $
instead of $\left(  \lambda,A,B\right)  $ and $\left(  \mu,C,D\right)  $)
shows that $\left\langle \nabla_{V,U}^{-\operatorname*{Col}},\ \nabla
_{Q,P}^{\operatorname*{Row}}\right\rangle =0$ (since $\left(  \mu,P,Q\right)
\in\operatorname*{SBT}\left(  n\right)  \subseteq\operatorname*{BT}\left(
n\right)  $ and $\left(  \lambda,U,V\right)  \in\operatorname*{SBT}\left(
n\right)  \subseteq\operatorname*{BT}\left(  n\right)  $ and $\mu\in X$ and
$\lambda\in Y$). This proves (\ref{pf.lem.bas.mur.twosid.triang.2.1b}).]

\item For each $\left(  \mu,P,Q\right)  \in\operatorname*{SBT}\left(
n\right)  $ satisfying $\left(  \mu,P,Q\right)  >\left(  \lambda,U,V\right)
$, we have%
\begin{equation}
\left\langle \nabla_{V,U}^{-\operatorname*{Col}},\ \nabla_{Q,P}%
^{\operatorname*{Row}}\right\rangle =0.
\label{pf.lem.bas.mur.twosid.triang.2.2}%
\end{equation}

[\textit{Proof:} This is just the equality (\ref{pf.lem.bas.mur.triang.2.2}),
which we showed during our proof of Lemma \ref{lem.bas.mur.triang.2}
\textbf{(a)}.]
\end{itemize}

Now, recall that the relation $<$ on $\operatorname*{SBT}\left(  n\right)  $
is the smaller relation of a total order. Hence, each $\left(  \mu,P,Q\right)
\in\operatorname*{SBT}\left(  n\right)  $ satisfies exactly one of the three
statements $\left(  \mu,P,Q\right)  <\left(  \lambda,U,V\right)  $ and
$\left(  \mu,P,Q\right)  >\left(  \lambda,U,V\right)  $ and $\left(
\mu,P,Q\right)  =\left(  \lambda,U,V\right)  $.

Accordingly, the sum on the right hand side of
(\ref{pf.lem.bas.mur.twosid.triang.2.0}) can be split into three parts: the
part comprising all $\left(  \mu,P,Q\right)  \in\operatorname*{SBT}\left(
n\right)  $ that satisfy $\left(  \mu,P,Q\right)  <\left(  \lambda,U,V\right)
$; the part comprising all $\left(  \mu,P,Q\right)  \in\operatorname*{SBT}%
\left(  n\right)  $ that satisfy $\left(  \mu,P,Q\right)  >\left(
\lambda,U,V\right)  $; and a single addend corresponding to $\left(
\mu,P,Q\right)  =\left(  \lambda,U,V\right)  $. Splitting the sum in this way,
we rewrite (\ref{pf.lem.bas.mur.twosid.triang.2.0}) as follows:
\begin{align*}
0  &  =\sum_{\substack{\left(  \mu,P,Q\right)  \in\operatorname*{SBT}\left(
n\right)  ;\\\left(  \mu,P,Q\right)  <\left(  \lambda,U,V\right)  }%
}\omega_{\mu,P,Q}\left\langle \nabla_{V,U}^{-\operatorname*{Col}}%
,\ \nabla_{Q,P}^{\operatorname*{Row}}\right\rangle \\
&  \ \ \ \ \ \ \ \ \ \ +\sum_{\substack{\left(  \mu,P,Q\right)  \in
\operatorname*{SBT}\left(  n\right)  ;\\\left(  \mu,P,Q\right)  >\left(
\lambda,U,V\right)  }}\omega_{\mu,P,Q}\underbrace{\left\langle \nabla
_{V,U}^{-\operatorname*{Col}},\ \nabla_{Q,P}^{\operatorname*{Row}%
}\right\rangle }_{\substack{=0\\\text{(by
(\ref{pf.lem.bas.mur.twosid.triang.2.2}))}}}+\,\omega_{\lambda,U,V}%
\left\langle \nabla_{V,U}^{-\operatorname*{Col}},\ \nabla_{V,U}%
^{\operatorname*{Row}}\right\rangle \\
&  =\sum_{\substack{\left(  \mu,P,Q\right)  \in\operatorname*{SBT}\left(
n\right)  ;\\\left(  \mu,P,Q\right)  <\left(  \lambda,U,V\right)  }%
}\omega_{\mu,P,Q}\left\langle \nabla_{V,U}^{-\operatorname*{Col}}%
,\ \nabla_{Q,P}^{\operatorname*{Row}}\right\rangle +\underbrace{\sum
_{\substack{\left(  \mu,P,Q\right)  \in\operatorname*{SBT}\left(  n\right)
;\\\left(  \mu,P,Q\right)  >\left(  \lambda,U,V\right)  }}\omega_{\mu,P,Q}%
0}_{=0}+\,\omega_{\lambda,U,V}\left\langle \nabla_{V,U}^{-\operatorname*{Col}%
},\ \nabla_{V,U}^{\operatorname*{Row}}\right\rangle \\
&  =\sum_{\substack{\left(  \mu,P,Q\right)  \in\operatorname*{SBT}\left(
n\right)  ;\\\left(  \mu,P,Q\right)  <\left(  \lambda,U,V\right)  }%
}\omega_{\mu,P,Q}\left\langle \nabla_{V,U}^{-\operatorname*{Col}}%
,\ \nabla_{Q,P}^{\operatorname*{Row}}\right\rangle +\omega_{\lambda
,U,V}\left\langle \nabla_{V,U}^{-\operatorname*{Col}},\ \nabla_{V,U}%
^{\operatorname*{Row}}\right\rangle \\
&  =\sum_{\substack{\left(  \mu,P,Q\right)  \in\operatorname*{SBT}\left(
n\right)  ;\\\left(  \mu,P,Q\right)  <\left(  \lambda,U,V\right)  ;\\\mu\in
Y}}\ \ \underbrace{\omega_{\mu,P,Q}}_{\substack{=0\\\text{(by
(\ref{pf.lem.bas.mur.twosid.triang.2.IH}))}}}\left\langle \nabla
_{V,U}^{-\operatorname*{Col}},\ \nabla_{Q,P}^{\operatorname*{Row}%
}\right\rangle +\sum_{\substack{\left(  \mu,P,Q\right)  \in\operatorname*{SBT}%
\left(  n\right)  ;\\\left(  \mu,P,Q\right)  <\left(  \lambda,U,V\right)
;\\\mu\notin Y}}\omega_{\mu,P,Q}\underbrace{\left\langle \nabla_{V,U}%
^{-\operatorname*{Col}},\ \nabla_{Q,P}^{\operatorname*{Row}}\right\rangle
}_{\substack{=0\\\text{(by (\ref{pf.lem.bas.mur.twosid.triang.2.1b}))}}}\\
&  \ \ \ \ \ \ \ \ \ \ +\omega_{\lambda,U,V}\left\langle \nabla_{V,U}%
^{-\operatorname*{Col}},\ \nabla_{V,U}^{\operatorname*{Row}}\right\rangle \\
&  \ \ \ \ \ \ \ \ \ \ \ \ \ \ \ \ \ \ \ \ \left(
\begin{array}
[c]{c}%
\text{here, we have split the sum into two subsums,}\\
\text{according to whether }\mu\in Y\text{ or }\mu\notin Y
\end{array}
\right) \\
&  =\underbrace{\sum_{\substack{\left(  \mu,P,Q\right)  \in\operatorname*{SBT}%
\left(  n\right)  ;\\\left(  \mu,P,Q\right)  <\left(  \lambda,U,V\right)
;\\\mu\in Y}}0\left\langle \nabla_{V,U}^{-\operatorname*{Col}},\ \nabla
_{Q,P}^{\operatorname*{Row}}\right\rangle }_{=0}+\underbrace{\sum
_{\substack{\left(  \mu,P,Q\right)  \in\operatorname*{SBT}\left(  n\right)
;\\\left(  \mu,P,Q\right)  <\left(  \lambda,U,V\right)  ;\\\mu\notin Y}%
}\omega_{\mu,P,Q}0}_{=0}+\,\omega_{\lambda,U,V}\left\langle \nabla
_{V,U}^{-\operatorname*{Col}},\ \nabla_{V,U}^{\operatorname*{Row}%
}\right\rangle \\
&  =\omega_{\lambda,U,V}\underbrace{\left\langle \nabla_{V,U}%
^{-\operatorname*{Col}},\ \nabla_{V,U}^{\operatorname*{Row}}\right\rangle
}_{\substack{=\pm1\\\text{(by Lemma \ref{lem.bas.mur.dia},}\\\text{applied to
}A=U\text{ and }B=V\text{)}}}=\pm\omega_{\lambda,U,V}.
\end{align*}
Hence, $\pm\omega_{\lambda,U,V}=0$. In other words, $\omega_{\lambda,U,V}=0$.
This completes the induction, and thus the proof of Claim 1.
\end{proof}

Now, (\ref{pf.lem.bas.mur.twosid.triang.2.a=wsum}) becomes%
\begin{align*}
\mathbf{a}  &  =\sum_{\left(  \lambda,U,V\right)  \in\operatorname*{SBT}%
\left(  n\right)  }\omega_{\lambda,U,V}\nabla_{V,U}^{\operatorname*{Row}}\\
&  =\sum_{\substack{\left(  \lambda,U,V\right)  \in\operatorname*{SBT}\left(
n\right)  ;\\\lambda\in X}}\omega_{\lambda,U,V}\nabla_{V,U}%
^{\operatorname*{Row}}+\sum_{\substack{\left(  \lambda,U,V\right)
\in\operatorname*{SBT}\left(  n\right)  ;\\\lambda\in Y}}\underbrace{\omega
_{\lambda,U,V}}_{\substack{=0\\\text{(by Claim 1)}}}\nabla_{V,U}%
^{\operatorname*{Row}}\\
&  \ \ \ \ \ \ \ \ \ \ \ \ \ \ \ \ \ \ \ \ \left(
\begin{array}
[c]{c}%
\text{since each }\left(  \lambda,U,V\right)  \in\operatorname*{SBT}\left(
n\right)  \text{ satisfies}\\
\text{either }\lambda\in X\text{ or }\lambda\in Y\text{ (but not both)}%
\end{array}
\right) \\
&  =\sum_{\substack{\left(  \lambda,U,V\right)  \in\operatorname*{SBT}\left(
n\right)  ;\\\lambda\in X}}\omega_{\lambda,U,V}\nabla_{V,U}%
^{\operatorname*{Row}}+\underbrace{\sum_{\substack{\left(  \lambda,U,V\right)
\in\operatorname*{SBT}\left(  n\right)  ;\\\lambda\in Y}}0\nabla
_{V,U}^{\operatorname*{Row}}}_{=0}\\
&  =\sum_{\substack{\left(  \lambda,U,V\right)  \in\operatorname*{SBT}\left(
n\right)  ;\\\lambda\in X}}\omega_{\lambda,U,V}\nabla_{V,U}%
^{\operatorname*{Row}}\\
&  \in\operatorname*{span}\left\{  \nabla_{V,U}^{\operatorname*{Row}}%
\ \mid\ \left(  \lambda,U,V\right)  \in\operatorname*{SBT}\left(  n\right)
\text{ and }\lambda\in X\right\} \\
&  =\MurpF_{\operatorname*{std},X}^{\operatorname*{Row}}%
\ \ \ \ \ \ \ \ \ \ \left(  \text{by the definition of }%
\MurpF_{\operatorname*{std},X}^{\operatorname*{Row}}\right)  .
\end{align*}
This completes our proof of Lemma \ref{lem.bas.mur.twosid.triang.2}
\textbf{(b)}. \medskip

\textbf{(a)} We must show that $\mathbf{a}\in\MurpF_{\operatorname*{std}%
,Y}^{-\operatorname*{Col}}$ for each $\mathbf{a}\in\left(
\MurpF_{\operatorname*{std},X}^{\operatorname*{Row}}\right)  ^{\perp}$.

So let $\mathbf{a}\in\left(  \MurpF_{\operatorname*{std},X}%
^{\operatorname*{Row}}\right)  ^{\perp}$. We must show that $\mathbf{a}%
\in\MurpF_{\operatorname*{std},Y}^{-\operatorname*{Col}}$.

We have $\mathbf{a}\in\left(  \MurpF_{\operatorname*{std},X}%
^{\operatorname*{Row}}\right)  ^{\perp}\subseteq\mathcal{A}$. Since the column
Murphy basis $\left(  \nabla_{V,U}^{-\operatorname*{Col}}\right)  _{\left(
\lambda,U,V\right)  \in\operatorname*{SBT}\left(  n\right)  }$ is a basis of
$\mathcal{A}$ (by Theorem \ref{thm.bas.mur.bas}), we can thus write
$\mathbf{a}$ as a $\mathbf{k}$-linear combination of its elements. In other
words, we can write $\mathbf{a}$ as
\begin{equation}
\mathbf{a}=\sum_{\left(  \lambda,U,V\right)  \in\operatorname*{SBT}\left(
n\right)  }\omega_{\lambda,U,V}\nabla_{V,U}^{-\operatorname*{Col}}
\label{pf.lem.bas.mur.twosid.triang.2.c=wsum}%
\end{equation}
for some scalars $\omega_{\lambda,U,V}\in\mathbf{k}$. Consider these
$\omega_{\lambda,U,V}$. We can rewrite the equality
(\ref{pf.lem.bas.mur.twosid.triang.2.c=wsum}) as
\begin{equation}
\mathbf{a}=\sum_{\left(  \mu,P,Q\right)  \in\operatorname*{SBT}\left(
n\right)  }\omega_{\mu,P,Q}\nabla_{Q,P}^{-\operatorname*{Col}}
\label{pf.lem.bas.mur.twosid.triang.2.c=wsum2}%
\end{equation}
(by renaming the summation index $\left(  \lambda,U,V\right)  $ as $\left(
\mu,P,Q\right)  $).

Now, we shall show the following:

\begin{statement}
\textit{Claim 2:} We have%
\[
\omega_{\lambda,U,V}=0\ \ \ \ \ \ \ \ \ \ \text{for all }\left(
\lambda,U,V\right)  \in\operatorname*{SBT}\left(  n\right)  \text{ satisfying
}\lambda\in X.
\]

\end{statement}

\begin{proof}
[Proof of Claim 2.]We proceed by strong induction on $\left(  \lambda
,U,V\right)  $ using the \textbf{reverse} of the total order on
$\operatorname*{SBT}\left(  n\right)  $ introduced in Definition
\ref{def.bas.tord.sbt.tord}. So we fix some $\left(  \lambda,U,V\right)
\in\operatorname*{SBT}\left(  n\right)  $ satisfying $\lambda\in X$. We assume
(as the induction hypothesis) that
\begin{equation}
\omega_{\mu,P,Q}=0 \label{pf.lem.bas.mur.twosid.triang.2.cIH}%
\end{equation}
for all $\left(  \mu,P,Q\right)  \in\operatorname*{SBT}\left(  n\right)  $
satisfying $\mu\in X$ and $\left(  \mu,P,Q\right)  >\left(  \lambda
,U,V\right)  $. Our goal is to prove that $\omega_{\lambda,U,V}=0$.

Since $\left(  \lambda,U,V\right)  \in\operatorname*{SBT}\left(  n\right)  $
and $\lambda\in X$, we have $\nabla_{V,U}^{\operatorname*{Row}}\in
\MurpF_{\operatorname*{std},X}^{\operatorname*{Row}}$ by the definition of
$\MurpF_{\operatorname*{std},X}^{\operatorname*{Row}}$ (in fact, $\nabla
_{V,U}^{\operatorname*{Row}}$ is one of the spanning vectors in the definition
of $\MurpF_{\operatorname*{std},X}^{\operatorname*{Row}}$).

Since $\mathbf{a}\in\left(  \MurpF_{\operatorname*{std},X}%
^{\operatorname*{Row}}\right)  ^{\perp}$, we have $\left\langle \mathbf{a}%
,\mathbf{v}\right\rangle =0$ for each $\mathbf{v}\in
\MurpF_{\operatorname*{std},X}^{\operatorname*{Row}}$ (by the definition of
$\left(  \MurpF_{\operatorname*{std},X}^{\operatorname*{Row}}\right)  ^{\perp
}$). Applying this to $\mathbf{v}=\nabla_{V,U}^{\operatorname*{Row}}$, we
obtain $\left\langle \mathbf{a},\ \nabla_{V,U}^{\operatorname*{Row}%
}\right\rangle =0$ (since $\nabla_{V,U}^{\operatorname*{Row}}\in
\MurpF_{\operatorname*{std},X}^{\operatorname*{Row}}$). Thus,%
\begin{align}
0  &  =\left\langle \mathbf{a},\ \nabla_{V,U}^{\operatorname*{Row}%
}\right\rangle =\left\langle \sum_{\left(  \mu,P,Q\right)  \in
\operatorname*{SBT}\left(  n\right)  }\omega_{\mu,P,Q}\nabla_{Q,P}%
^{-\operatorname*{Col}},\ \nabla_{V,U}^{\operatorname*{Row}}\right\rangle
\ \ \ \ \ \ \ \ \ \ \left(  \text{by
(\ref{pf.lem.bas.mur.twosid.triang.2.c=wsum2})}\right) \nonumber\\
&  =\sum_{\left(  \mu,P,Q\right)  \in\operatorname*{SBT}\left(  n\right)
}\omega_{\mu,P,Q}\left\langle \nabla_{Q,P}^{-\operatorname*{Col}}%
,\ \nabla_{V,U}^{\operatorname*{Row}}\right\rangle
\label{pf.lem.bas.mur.twosid.triang.2.c0}%
\end{align}
(since the form $\left\langle \cdot,\cdot\right\rangle $ is bilinear). Now, we
observe the following:

\begin{itemize}
\item For each $\left(  \mu,P,Q\right)  \in\operatorname*{SBT}\left(
n\right)  $ satisfying $\mu\notin X$, we have
\begin{equation}
\left\langle \nabla_{Q,P}^{-\operatorname*{Col}},\ \nabla_{V,U}%
^{\operatorname*{Row}}\right\rangle =0.
\label{pf.lem.bas.mur.twosid.triang.2.c1b}%
\end{equation}

[\textit{Proof:} Let $\left(  \mu,P,Q\right)  \in\operatorname*{SBT}\left(
n\right)  $ be such that $\mu\notin X$. Hence, $\mu\in\operatorname*{Par}%
\left(  n\right)  \setminus X$ (since $\mu\in\operatorname*{Par}\left(
n\right)  $ and $\mu\notin X$). Recall that $X\cup Y=\operatorname*{Par}%
\left(  n\right)  $, so that $\operatorname*{Par}\left(  n\right)  \setminus
X\subseteq Y$. Hence, $\mu\in\operatorname*{Par}\left(  n\right)  \setminus
X\subseteq Y$. Therefore, Lemma \ref{lem.bas.mur.twosid.orth} \textbf{(a)}
(applied to $\left(  \lambda,U,V\right)  $ and $\left(  \mu,P,Q\right)  $
instead of $\left(  \lambda,A,B\right)  $ and $\left(  \mu,C,D\right)  $)
shows that $\left\langle \nabla_{Q,P}^{-\operatorname*{Col}},\ \nabla
_{V,U}^{\operatorname*{Row}}\right\rangle =0$ (since $\left(  \lambda
,U,V\right)  \in\operatorname*{SBT}\left(  n\right)  \subseteq
\operatorname*{BT}\left(  n\right)  $ and $\left(  \mu,P,Q\right)
\in\operatorname*{SBT}\left(  n\right)  \subseteq\operatorname*{BT}\left(
n\right)  $ and $\lambda\in X$ and $\mu\in Y$). This proves
(\ref{pf.lem.bas.mur.twosid.triang.2.c1b}).]

\item For each $\left(  \mu,P,Q\right)  \in\operatorname*{SBT}\left(
n\right)  $ satisfying $\left(  \mu,P,Q\right)  <\left(  \lambda,U,V\right)
$, we have%
\begin{equation}
\left\langle \nabla_{Q,P}^{-\operatorname*{Col}},\ \nabla_{V,U}%
^{\operatorname*{Row}}\right\rangle =0.
\label{pf.lem.bas.mur.twosid.triang.2.c2}%
\end{equation}

[\textit{Proof:} This is just the equality (\ref{pf.lem.bas.mur.triang.2.c2}),
which we showed during our proof of Lemma \ref{lem.bas.mur.triang.2}
\textbf{(c)}.]
\end{itemize}

Now, recall that the relation $<$ on $\operatorname*{SBT}\left(  n\right)  $
is the smaller relation of a total order. Hence, each $\left(  \mu,P,Q\right)
\in\operatorname*{SBT}\left(  n\right)  $ satisfies exactly one of the three
statements $\left(  \mu,P,Q\right)  <\left(  \lambda,U,V\right)  $ and
$\left(  \mu,P,Q\right)  >\left(  \lambda,U,V\right)  $ and $\left(
\mu,P,Q\right)  =\left(  \lambda,U,V\right)  $.

Accordingly, the sum on the right hand side of
(\ref{pf.lem.bas.mur.twosid.triang.2.c0}) can be split into three parts: the
part comprising all $\left(  \mu,P,Q\right)  \in\operatorname*{SBT}\left(
n\right)  $ that satisfy $\left(  \mu,P,Q\right)  <\left(  \lambda,U,V\right)
$; the part comprising all $\left(  \mu,P,Q\right)  \in\operatorname*{SBT}%
\left(  n\right)  $ that satisfy $\left(  \mu,P,Q\right)  >\left(
\lambda,U,V\right)  $; and a single addend corresponding to $\left(
\mu,P,Q\right)  =\left(  \lambda,U,V\right)  $. Splitting the sum in this way,
we rewrite (\ref{pf.lem.bas.mur.twosid.triang.2.c0}) as follows:
\begin{align*}
0  &  =\sum_{\substack{\left(  \mu,P,Q\right)  \in\operatorname*{SBT}\left(
n\right)  ;\\\left(  \mu,P,Q\right)  <\left(  \lambda,U,V\right)  }%
}\omega_{\mu,P,Q}\underbrace{\left\langle \nabla_{Q,P}^{-\operatorname*{Col}%
},\ \nabla_{V,U}^{\operatorname*{Row}}\right\rangle }_{\substack{=0\\\text{(by
(\ref{pf.lem.bas.mur.twosid.triang.2.c2}))}}}\\
&  \ \ \ \ \ \ \ \ \ \ +\sum_{\substack{\left(  \mu,P,Q\right)  \in
\operatorname*{SBT}\left(  n\right)  ;\\\left(  \mu,P,Q\right)  >\left(
\lambda,U,V\right)  }}\omega_{\mu,P,Q}\left\langle \nabla_{Q,P}%
^{-\operatorname*{Col}},\ \nabla_{V,U}^{\operatorname*{Row}}\right\rangle
+\omega_{\lambda,U,V}\left\langle \nabla_{V,U}^{-\operatorname*{Col}}%
,\ \nabla_{V,U}^{\operatorname*{Row}}\right\rangle \\
&  =\underbrace{\sum_{\substack{\left(  \mu,P,Q\right)  \in\operatorname*{SBT}%
\left(  n\right)  ;\\\left(  \mu,P,Q\right)  <\left(  \lambda,U,V\right)
}}\omega_{\mu,P,Q}0}_{=0}+\sum_{\substack{\left(  \mu,P,Q\right)
\in\operatorname*{SBT}\left(  n\right)  ;\\\left(  \mu,P,Q\right)  >\left(
\lambda,U,V\right)  }}\omega_{\mu,P,Q}\left\langle \nabla_{Q,P}%
^{-\operatorname*{Col}},\ \nabla_{V,U}^{\operatorname*{Row}}\right\rangle
+\omega_{\lambda,U,V}\left\langle \nabla_{V,U}^{-\operatorname*{Col}}%
,\ \nabla_{V,U}^{\operatorname*{Row}}\right\rangle \\
&  =\sum_{\substack{\left(  \mu,P,Q\right)  \in\operatorname*{SBT}\left(
n\right)  ;\\\left(  \mu,P,Q\right)  >\left(  \lambda,U,V\right)  }%
}\omega_{\mu,P,Q}\left\langle \nabla_{Q,P}^{-\operatorname*{Col}}%
,\ \nabla_{V,U}^{\operatorname*{Row}}\right\rangle +\omega_{\lambda
,U,V}\left\langle \nabla_{V,U}^{-\operatorname*{Col}},\ \nabla_{V,U}%
^{\operatorname*{Row}}\right\rangle \\
&  =\sum_{\substack{\left(  \mu,P,Q\right)  \in\operatorname*{SBT}\left(
n\right)  ;\\\left(  \mu,P,Q\right)  >\left(  \lambda,U,V\right)  ;\\\mu\in
X}}\ \ \underbrace{\omega_{\mu,P,Q}}_{\substack{=0\\\text{(by
(\ref{pf.lem.bas.mur.twosid.triang.2.cIH}))}}}\left\langle \nabla
_{Q,P}^{-\operatorname*{Col}},\ \nabla_{V,U}^{\operatorname*{Row}%
}\right\rangle +\sum_{\substack{\left(  \mu,P,Q\right)  \in\operatorname*{SBT}%
\left(  n\right)  ;\\\left(  \mu,P,Q\right)  >\left(  \lambda,U,V\right)
;\\\mu\notin X}}\omega_{\mu,P,Q}\underbrace{\left\langle \nabla_{Q,P}%
^{-\operatorname*{Col}},\ \nabla_{V,U}^{\operatorname*{Row}}\right\rangle
}_{\substack{=0\\\text{(by (\ref{pf.lem.bas.mur.twosid.triang.2.c1b}))}}}\\
&  \ \ \ \ \ \ \ \ \ \ +\omega_{\lambda,U,V}\left\langle \nabla_{V,U}%
^{-\operatorname*{Col}},\ \nabla_{V,U}^{\operatorname*{Row}}\right\rangle \\
&  \ \ \ \ \ \ \ \ \ \ \ \ \ \ \ \ \ \ \ \ \left(
\begin{array}
[c]{c}%
\text{here, we have split the sum into two subsums,}\\
\text{according to whether }\mu\in X\text{ or }\mu\notin X
\end{array}
\right) \\
&  =\underbrace{\sum_{\substack{\left(  \mu,P,Q\right)  \in\operatorname*{SBT}%
\left(  n\right)  ;\\\left(  \mu,P,Q\right)  >\left(  \lambda,U,V\right)
;\\\mu\in X}}0\left\langle \nabla_{Q,P}^{-\operatorname*{Col}},\ \nabla
_{V,U}^{\operatorname*{Row}}\right\rangle }_{=0}+\underbrace{\sum
_{\substack{\left(  \mu,P,Q\right)  \in\operatorname*{SBT}\left(  n\right)
;\\\left(  \mu,P,Q\right)  >\left(  \lambda,U,V\right)  ;\\\mu\notin X}%
}\omega_{\mu,P,Q}0}_{=0}+\,\omega_{\lambda,U,V}\left\langle \nabla
_{V,U}^{-\operatorname*{Col}},\ \nabla_{V,U}^{\operatorname*{Row}%
}\right\rangle \\
&  =\omega_{\lambda,U,V}\underbrace{\left\langle \nabla_{V,U}%
^{-\operatorname*{Col}},\ \nabla_{V,U}^{\operatorname*{Row}}\right\rangle
}_{\substack{=\pm1\\\text{(by Lemma \ref{lem.bas.mur.dia},}\\\text{applied to
}A=U\text{ and }B=V\text{)}}}=\pm\omega_{\lambda,U,V}.
\end{align*}
Hence, $\pm\omega_{\lambda,U,V}=0$. In other words, $\omega_{\lambda,U,V}=0$.
This completes the induction, and thus the proof of Claim 2.
\end{proof}

Now, (\ref{pf.lem.bas.mur.twosid.triang.2.c=wsum}) becomes%
\begin{align*}
\mathbf{a}  &  =\sum_{\left(  \lambda,U,V\right)  \in\operatorname*{SBT}%
\left(  n\right)  }\omega_{\lambda,U,V}\nabla_{V,U}^{-\operatorname*{Col}}\\
&  =\sum_{\substack{\left(  \lambda,U,V\right)  \in\operatorname*{SBT}\left(
n\right)  ;\\\lambda\in X}}\underbrace{\omega_{\lambda,U,V}}%
_{\substack{=0\\\text{(by Claim 2)}}}\nabla_{V,U}^{-\operatorname*{Col}}%
+\sum_{\substack{\left(  \lambda,U,V\right)  \in\operatorname*{SBT}\left(
n\right)  ;\\\lambda\in Y}}\omega_{\lambda,U,V}\nabla_{V,U}%
^{-\operatorname*{Col}}\\
&  \ \ \ \ \ \ \ \ \ \ \ \ \ \ \ \ \ \ \ \ \left(
\begin{array}
[c]{c}%
\text{since each }\left(  \lambda,U,V\right)  \in\operatorname*{SBT}\left(
n\right)  \text{ satisfies}\\
\text{either }\lambda\in X\text{ or }\lambda\in Y\text{ (but not both)}%
\end{array}
\right) \\
&  =\underbrace{\sum_{\substack{\left(  \lambda,U,V\right)  \in
\operatorname*{SBT}\left(  n\right)  ;\\\lambda\in X}}0\nabla_{V,U}%
^{-\operatorname*{Col}}}_{=0}+\sum_{\substack{\left(  \lambda,U,V\right)
\in\operatorname*{SBT}\left(  n\right)  ;\\\lambda\in Y}}\omega_{\lambda
,U,V}\nabla_{V,U}^{-\operatorname*{Col}}\\
&  =\sum_{\substack{\left(  \lambda,U,V\right)  \in\operatorname*{SBT}\left(
n\right)  ;\\\lambda\in Y}}\omega_{\lambda,U,V}\nabla_{V,U}%
^{-\operatorname*{Col}}\\
&  \in\operatorname*{span}\left\{  \nabla_{V,U}^{-\operatorname*{Col}}%
\ \mid\ \left(  \lambda,U,V\right)  \in\operatorname*{SBT}\left(  n\right)
\text{ and }\lambda\in Y\right\} \\
&  =\MurpF_{\operatorname*{std},Y}^{-\operatorname*{Col}}%
\ \ \ \ \ \ \ \ \ \ \left(  \text{by the definition of }%
\MurpF_{\operatorname*{std},Y}^{-\operatorname*{Col}}\right)  .
\end{align*}
This completes our proof of Lemma \ref{lem.bas.mur.twosid.triang.2}
\textbf{(a)}.
\end{proof}
\end{fineprint}

\begin{fineprint}
\begin{proof}
[Detailed proof of Theorem \ref{thm.bas.mur.twosid.orth-XY}.]\textbf{(a)} We
have $\MurpF_{\operatorname*{std},X}^{\operatorname*{Row}}\subseteq
\MurpF_{\operatorname*{all},X}^{\operatorname*{Row}}$ (this follows from the
definitions of these two sets\footnote{\textit{Proof.} The definition of
$\MurpF_{\operatorname*{std},X}^{\operatorname*{Row}}$ yields%
\begin{align*}
\MurpF_{\operatorname*{std},X}^{\operatorname*{Row}}  &  =\operatorname*{span}%
\left\{  \nabla_{V,U}^{\operatorname*{Row}}\ \mid\ \left(  \lambda,U,V\right)
\in\underbrace{\operatorname*{SBT}\left(  n\right)  }_{\subseteq
\operatorname*{BT}\left(  n\right)  }\text{ and }\lambda\in X\right\} \\
&  \subseteq\operatorname*{span}\left\{  \nabla_{V,U}^{\operatorname*{Row}%
}\ \mid\ \left(  \lambda,U,V\right)  \in\operatorname*{BT}\left(  n\right)
\text{ and }\lambda\in X\right\}  =\MurpF_{\operatorname*{all},X}%
^{\operatorname*{Row}}%
\end{align*}
(by the definition of $\MurpF_{\operatorname*{all},X}^{\operatorname*{Row}}%
$).}). Hence,%
\begin{align*}
\MurpF_{\operatorname*{std},X}^{\operatorname*{Row}}  &  \subseteq
\MurpF_{\operatorname*{all},X}^{\operatorname*{Row}}\subseteq\left(
\MurpF_{\operatorname*{std},Y}^{-\operatorname*{Col}}\right)  ^{\perp
}\ \ \ \ \ \ \ \ \ \ \left(  \text{by Lemma \ref{lem.bas.mur.twosid.triang.1}
\textbf{(b)}}\right) \\
&  \subseteq\MurpF_{\operatorname*{std},X}^{\operatorname*{Row}}%
\ \ \ \ \ \ \ \ \ \ \left(  \text{by Lemma \ref{lem.bas.mur.twosid.triang.2}
\textbf{(b)}}\right)  .
\end{align*}
This is a chain of inclusions that starts and ends with the same set (namely,
$\MurpF_{\operatorname*{std},X}^{\operatorname*{Row}}$). Thus, all inclusion
signs in this chain must actually be equality signs (since any three sets
$X,Y,Z$ satisfying $X\subseteq Y\subseteq Z\subseteq X$ must satisfy $X=Y=Z$).
That is, we must have $\MurpF_{\operatorname*{std},X}^{\operatorname*{Row}%
}=\MurpF_{\operatorname*{all},X}^{\operatorname*{Row}}=\left(
\MurpF_{\operatorname*{std},Y}^{-\operatorname*{Col}}\right)  ^{\perp}$.

It remains to show that $\MurpF_{\operatorname*{std},X}^{\operatorname*{Row}}$
is a two-sided ideal of $\mathcal{A}$.

But Proposition \ref{prop.bas.mur.twosid.Sn-act.all} shows that
$\MurpF_{\operatorname*{all},X}^{\operatorname*{Row}}$ is a two-sided ideal of
$\mathcal{A}$. In other words, $\MurpF_{\operatorname*{std},X}%
^{\operatorname*{Row}}$ is a two-sided ideal of $\mathcal{A}$ (since
$\MurpF_{\operatorname*{std},X}^{\operatorname*{Row}}%
=\MurpF_{\operatorname*{all},X}^{\operatorname*{Row}}$). Hence, the proof of
Theorem \ref{thm.bas.mur.twosid.orth-XY} \textbf{(a)} is complete. \medskip

\textbf{(b)} The proof of part \textbf{(a)} applies to part \textbf{(b)} with
minor changes (the inequality signs need to be changed, and we must use Lemma
\ref{lem.bas.mur.twosid.triang.1} \textbf{(a)} instead of Lemma
\ref{lem.bas.mur.twosid.triang.1} \textbf{(b)} and use Lemma
\ref{lem.bas.mur.twosid.triang.2} \textbf{(a)} instead of Lemma
\ref{lem.bas.mur.twosid.triang.2} \textbf{(b)}).
\end{proof}
\end{fineprint}

\end{document}